%% file: knotoid_class.tex
\documentclass[reqno]{amsart}
\usepackage{graphicx}
\usepackage{hyperref}
\usepackage{array,amsmath}
\usepackage{multirow}
\usepackage{xcolor}
\newcolumntype{L}{>$l<$}
\allowdisplaybreaks

\let\mypdfximage\pdfximage
\def\pdfximage{\immediate\mypdfximage}

\title[A systematic classification of knotoids on the plane and on $S^2$]{A systematic classification of knotoids on the plane and on the sphere}

\author{Dimos Goundaroulis}
\address{Center for Integrative Genomics, University of Lausanne, Lausanne, Switzerland.\\
SIB Swiss Institute of Bioinformatics, Lausanne, Switzerland}
\email{dimoklis.gkountaroulis@unil.ch}
\author{Julien Dorier}
\address{Vital-IT, SIB Swiss Institute of Bioinformatics, Lausanne, Switzerland} 
\email{julien.dorier@sib.swiss}
\author{Andrzej Stasiak}
\address{Center for Integrative Genomics, University of Lausanne, Lausanne, Switzerland.\\
SIB Swiss Institute of Bioinformatics, Lausanne, Switzerland}
\email{andrzej.stasiak@unil.ch}
\begin{document}

\begin{abstract}
In this paper we generate and systematically classify all  prime planar knotoids with up to 5 crossings. We also extend the existing list of knotoids in $S^2$ and add all knotoids with 6 crossings.
\end{abstract}

\maketitle

\section{Introduction}

The theory of knotoids  forms a diagrammatic theory of open-ended oriented arcs that extends knot theory and it was introduced in 2012 by Turaev \cite{turaev2012}.  The notion of knotoids is defined through the different equivalence classes of open-ended diagrams up to a specific set of isotopy moves (Fig.~\ref{fig:Rmoves}).  Each equivalence class gives rise to a particular knotoid type. Several studies appeared recently that either extend concepts from knot theory to the case of knotoids or are dedicated studies on the theory of knotoids \cite{barbensi, chapman, gugumcu2017, gugumcu2017_2, gugumcu2018, kodokostas, korablev2017}.  

Even though the classification of knotoids is interesting on its own merit, our motivation for this work comes from biology. Knotoids can be considered as projections of open-ended embedded curves in 3-space \cite{gugumcu2017, przybyl}. For this reason they are used in the topological characterization of protein structures \cite{goundaroulis_sr2017, goundaroulis_pol2017, knoto-id, knotprot2}. Proteins are long, linear, open-ended biopolymers that are made out of amino acids connected by peptide bonds. These chains fold into conformations that in the case of some proteins reproducibly fold into open-ended knots. Knotted proteins retained their entangled structure throughout evolution despite the fact that their folding is less efficient and slower than the folding of unknotted proteins of similar size.  \cite{mallam2012,dabrowski2015,sulkowska_on2012}. Apparently, the knotted structure of some proteins, gives them advantages that cannot be achieved otherwise (see for example  \cite{virnau2006, sriramoju2018, dabrowski2016}). 
Having a topological characterization of proteins is needed for better understanding of the relation between structure and function of knotted proteins. Earlier methods for the topological characterization of proteins (e.g. \cite{sulkowska_2012}) required the artificial closure of the protein chain in order to form a knot, altering, thus, its geometry. Analyzing protein structures using the concept of knotoids not only preserves their geometry but 
also provides a more detail overview of a protein's topology \cite{goundaroulis_sr2017}, especially when the analysis uses planar knotoids \cite{goundaroulis_pol2017}. 

In this paper we provide a systematic classification of all planar knotoids with up to  5 crossings. We also extend the table of known knotoids in $S^2$ \cite{barth,korablev2018} and we add to it all knotoids with six crossings. The classification is up to reversion and up to three different types of symmetry-related involutions presented in Fig.~\ref{fig:invos}.  We propose a two-number notation of knotoids, where the first number indicates the minimal number of crossings and the second number indicates the rank of a given knotoid type among all knotoids types with the same minimal number of crossings. The rank is based on a total order of the corresponding realizable extended Gauss codes of a given type of knotoids. Our notation is consistent with both the planar and the $S^2$ case. It includes also all knot-type knotoids which are given the same notation as their corresponding knots in the Rolfsen table \cite{Adams}.

\section{Knotoids}
\begin{figure}[!htb]
\centering
\includegraphics{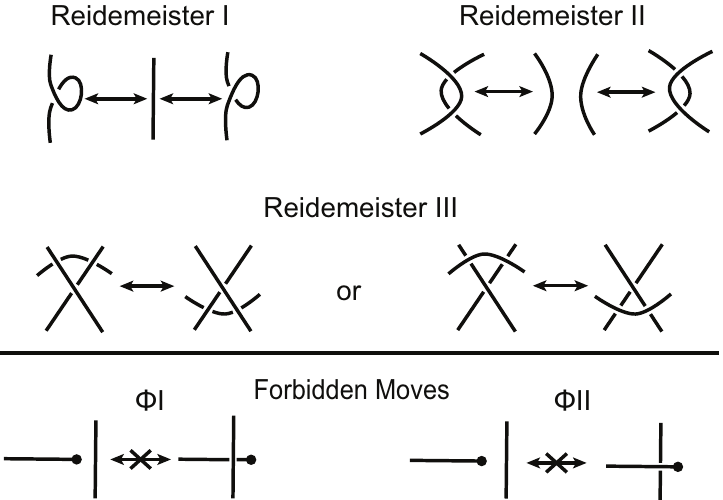}
\caption{The three Reidemeister moves and the two forbidden moves for knotoid diagrams.}\label{fig:Rmoves}
\end{figure}

Knotoids are oriented, open-ended, knot-like objects that are conveniently represented by knotoid diagrams. Knotoid diagrams are defined formally as a generic immersion of the interval $[0,\ 1]$ into $\mathbb{R}^2$ or the sphere $S^2$.  Only finitely-many double-points or self-intersections of the diagram  are allowed, while the ends of the interval are mapped to distinct points called the {\it tail} and the {\it head} of the diagram respectively. Additional information is included at the double points indicating which arc goes under and which goes over. A knotoid diagram is usually oriented from tail to head.  Furthermore, the endpoints are considered fixed in the local region of the diagram that they lie, meaning that they are not allowed to cross over or under any arc of the diagram. Two diagrams that can be deformed to one another using planar isotopy or isotopy in $S^2$, and a finite sequence of moves called the Reidemeister moves ($\Omega$-moves), are considered equivalent. The $\Omega$-moves are performed locally on a diagram and always away from the endpoints (See Fig.~\ref{fig:Rmoves}). Different knotoid types are represented by different equivalence classes of knotoid diagrams. The definition of knotoids can be generalized to any orientable surface. For example, virtual knotoids are introduced in \cite{gugumcu2017}.  

\begin{figure}[!htb]
\centering
\includegraphics{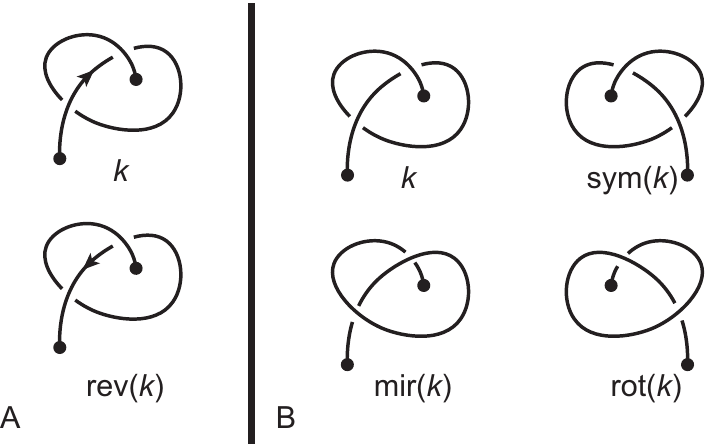}
\caption{Reversing the orientation of a knotoid $k$ with the reversion involution (A). The symmetry, mirror reflection and rotation involutions applied on the knotoid $k$ (B).}\label{fig:invos}
\end{figure}

There are four involutive operations that can be defined on knotoid diagrams \cite{turaev2012}. Reversion, ${\rm rev}(k)$, reverses the orientation of a knotoid diagram. Symmetry, ${\rm sym}(k)$, reflects a knotoid diagram with respect to the vertical line ${ 0 } \times \mathbb{R} \subset \mathbb{R}^2$ and can be extended to a self-homeomorphism of $S^2 \cong \mathbb{R}^2 \cup { \infty }$ by $\infty \mapsto \infty$. The next involution that can be defined on knotoids is the mirror reflection, ${\rm mir}(k)$, that changes undercrossings to overcrossings and vice versa (see Fig~\ref{fig:invos}). The composition of the mirror reflection and the symmetry is called {\it rotation}, ${\rm rot}(k)$, and it can be thought of as the rotation of a knotoid $k$ around the axis that passes through the endpoints  \cite{barbensi}.

\begin{figure}[!htb]
\centering
\includegraphics{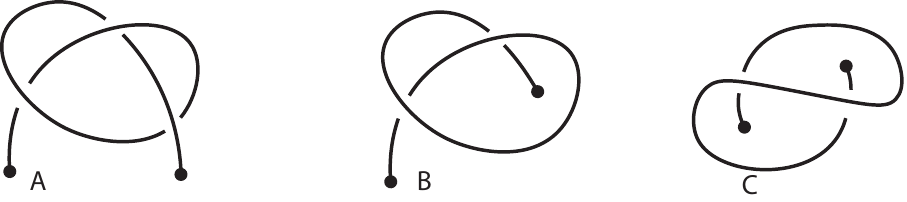}
\caption{(A) A knot-type knotoid, (B) a proper knotoid and (C) a planar knotoid that is neither knot-type nor proper.}\label{fig:simplify}
\end{figure}

There is a natural way to obtain a knot diagram from a knotoid diagram by connecting the endpoints with an arc that passes always under the rest of the diagram. A different knot may be obtained if the closure arc passes always over all other arcs. Planar knotoid diagrams with both endpoints in the outer region of the diagram are called {\it knot-type knotoids} (see Fig.~\ref{fig:simplify}A) and they are identified with their corresponding knot \cite{turaev2012}. In $S^2$ it is sufficient for a knot-type knotoid to have both endpoints in the same region. {\it Proper knotoids} are those knotoids that have exactly one endpoint in the outer region of the diagram and every knotoid in $S^2$ can be represented by a proper knotoid \cite{turaev2012}.

We call {\it rotatable} the knotoids that are isotopic to their rotation involution. Knot-type knotoids are known to be rotatable \cite{barbensi}. {\it Achiral} are the knoitoids that are isotopic to their mirror reflection while {\it strongly achiral} are those knotoids that are both rotatable and achiral. Achiral knots correspond to strongly achiral knot-type knotoids. 

Finally,  it is worth mentioning that non-equivalent planar knotoids may become equivalent when they are considered in $S^2$ \cite{turaev2012}. This is because we are free to move arcs around the surface of the sphere using isotopy. For example in Figure~\ref{fig:simplify} knotoids B and C are not equivalent as planar knotoids but are equivalent as $S^2$-knotoids.

 \section{Knotoid invariants}\label{sec:invariants}
 Our classification algorithm uses a number of knotoid invariants in order to distinguish non-isotopic knotoid diagrams.
For the case of knotoids in $S^2$ we use the arrow polynomial \cite{gugumcu2017}. The case of planar knotoids proved more challenging and so we used a combination of two different methods, the loop arrow polynomial and the double-branched covers of knotoids \cite{barbensi}. 

 \subsection{The arrow polynomial}
\begin{figure}[h] 
\centering
\includegraphics{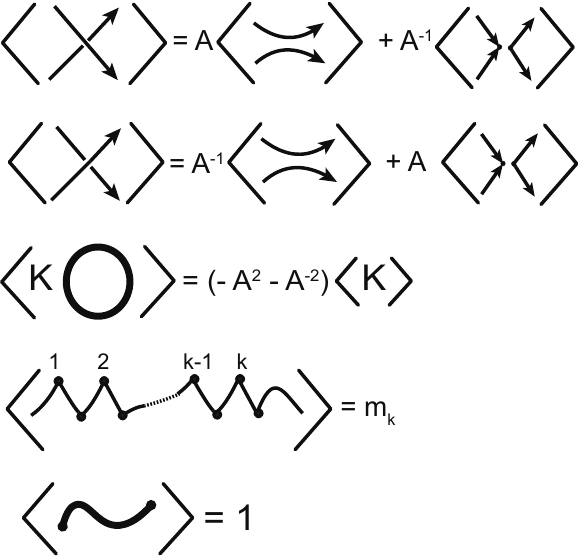}
\caption{\footnotesize The skein relation and the axioms of the arrow polynomial.}\label{fig:arrow}
\end{figure}

\begin{figure}[h] 
\centering
\includegraphics{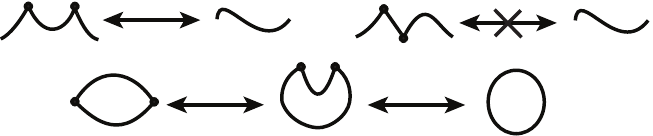}
\caption{\footnotesize The allowed and the forbidden cancellation rules.}\label{fig:cancel}
\end{figure}

 The arrow polynomial is based on the oriented state expansion of the bracket polynomial and it was initially  defined in \cite{kauffmandye2009} as an invariant for virtual knots. It is a Laurent polynomial that takes values in the ring $\mathbb{Z}[A, A^{-1}, m_1, m_2, \ldots ],$ where $m_i$ are an infinite set of independent commuting variables that also commute with the variable $A$.
 In \cite{gugumcu2017} it has been extended to the cases of classical knotoids in $S^2$ and virtual knotoids. The arrow polynomial can be defined recursively using the skein relation and the set of rules  shown in Fig.~\ref{fig:arrow}, that involve  smoothings with matching or conflicting orientations of arcs.
 Smoothing a crossing in a disorienting way results in a pair of cusps. If the acute angles of two consecutive cusps are in the same local region of the diagram, then they can be cancelled out  (see Fig.~\ref{fig:cancel}). Otherwise, if the acute angles of two consecutive cusps are in different local regions of the diagram then the cancelation is not possible. Each state includes a number of circular components and a long segment component, all of which may contain a number of consecutive cusps. The arrow polynomial assigns a new variable to each long segment of a state with a number of surviving cusps. In particular, two consecutive surviving cusps form a zigzag and a long segment with $2k$ surviving cusps is evaluated at $m_k$.

\subsection{The loop arrow polynomial}

The loop arrow polynomial is the extension of the arrow polynomial to the case of planar knotoids and it was first mentioned in \cite{goundaroulis_pol2017}. It is a Laurent polynomial that takes values the ring:
\[
\mathbb{Z}[A, A^{-1}, v, m_1, m_2, \ldots , w_1, w_2, \ldots, p_1, p_2, \ldots , q_1, q_2 \ldots ],
\]
 \begin{figure}[h] 
\centering
\includegraphics{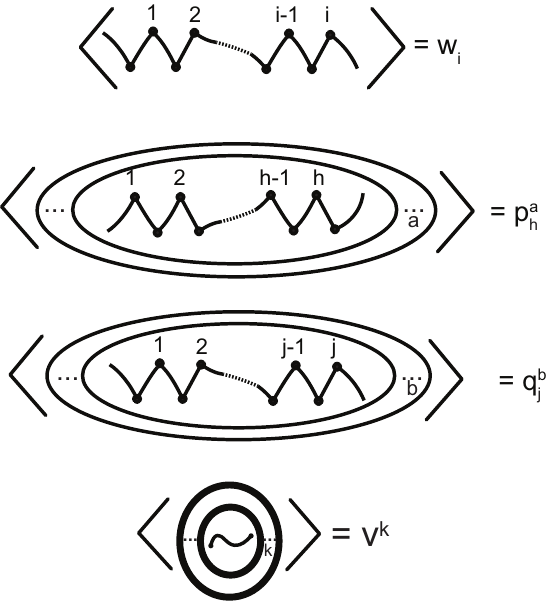}
\caption{\footnotesize The additional rules that the loop arrow polynomial has to satisfy.}\label{fig:looparrow}
\end{figure}
where the $m_i$, $w_j$, $p_k$, $q_\ell$ are all infinite sets of independent commuting variables that also commute with the each other but also with variables $A$ and $v$. The loop arrow distinguishes two different types of zigzags and it assigns one of the variables  $m_i$ or $w_i$, depending on the type of zigzag. Furthermore, the loop arrow polynomial assigns different variables to circular components that enclose the long segment with or without any of the two types of zigzags. The additional rules of Fig.~\ref{fig:looparrow} together with those of the arrow polynomial, define recursively the loop arrow polynomial.

\subsection{Double-branch covers of knotoids}\label{sec:branch} A planar knotoid diagram can be embedded in $\mathbb{R}^3$ in the following way  \cite{gugumcu2017}. We identify the plane where the planar knotoid diagram lies in with $\mathbb{R}^2 \times \{0\} \subset \mathbb{R}^3$ and we push the over-crossings in the upper half-space and the under-passes in the lower half-space. The endpoints of the diagram lie on and can move along two infinite lines that are perpendicular to the plane. Using the methods of \cite{barbensi} we consider a planar knotoid diagram  as an embedded arc inside the cylinder $D^2 \times I$. The double-branched cover of the cylinder over the two infinite lines is the solid torus $S^1 \times D^2$. Under this operation, a planar knotoid diagram is associated to a knot in the solid torus and its knot type is a knotoid invariant.  Therefore, if two planar knotoid diagrams correspond to non-equivalent knots in the solid torus then they are non-isotopic. 

\begin{figure}[!htb]
\centering
\includegraphics{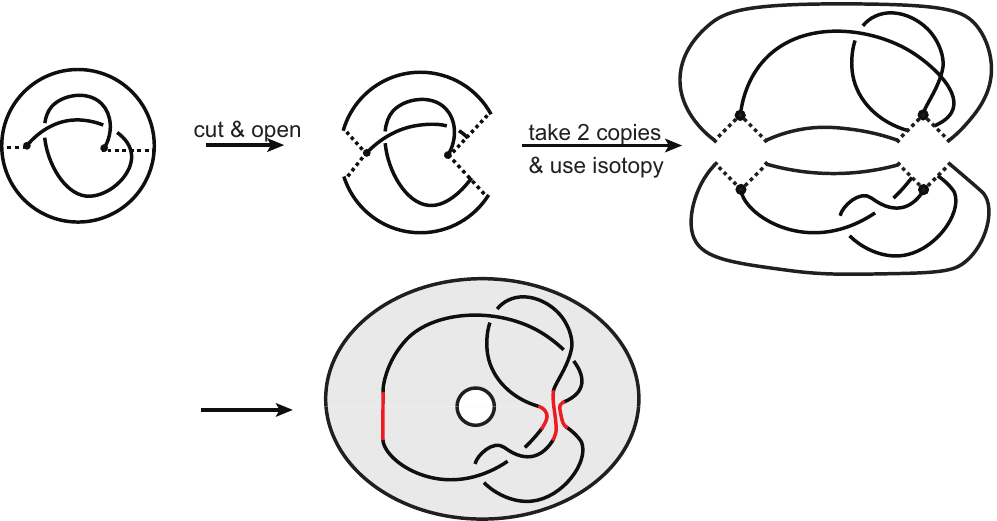}
\caption{The pre-image of a knotoid in $D^2$ in its double branch cover is a knot in the annulus which lifts to a knot in the solid torus.}\label{fig:branch}
\end{figure}

In brief, the process of finding the pre-image of  a planar knotoid diagram $k$ in its double branch cover is as follows. Consider  $k$ lying inside a disk and extend two lines, one from each endpoint towards the boundary of the disk. Cut and open the disk along those lines (see Fig.~\ref{fig:branch}). Next,  two copies of the cut and opened disk are considered, that we isotope them so that each cut opposes its copy and then we glue them using an appropriate homeomorphism.  The result is a knot in the annulus which lifts to a knot in the solid  torus. 

\section{Encoding diagrams} 
\subsection{Gauss Codes}\label{sec:gc} In order to classify all planar knotoids up to five crossings, we would like to encode their diagrams in a way that is easily handled by a computer. One of the standard notations for encoding knot diagrams as well as knotoid diagrams in $S^2$, that fulfills this criterion, is  the oriented Gauss code. 

The oriented Gauss code is a pair that consists in the {\it Gauss word} which is a sequence of labels that are assigned to a diagram's crossings as one, starting from the tail of a diagram, travels around the diagram and a sequence of the signs of each of the crossings of the diagram. Each crossing appears twice in the Gauss word since the crossings are encountered twice during this trip, once as an undercrossing and once as an overcrossing. To indicate in the Gauss word an undercrossing we add a ``-'' before the label and to indicate an overcrossing we add a ``+''. The length of a Gauss word is $2n$, where $n$ is the number of crossings of the diagram while the length of the signs' sequence is $n$.  The oriented Gauss code  together with the set of signs represents uniquely a knot or a knotoid diagram in $S^2$ up to isotopy. 

\begin{figure}[!htb]
\centering
\includegraphics{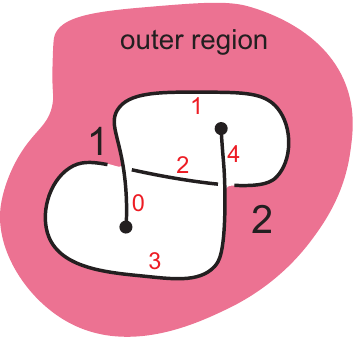}
\caption{A labeled knotoid diagram. The black numbers correspond to the crossings while the red ones to the arcs.}\label{fig:gc_example}
\end{figure}

For planar knotoids, we adapt to our case the extended oriented Gauss code \cite{gabrov, gabrov1}. More specifically, we attach to the oriented Gauss code a third piece of information, that is the list of labels of arcs that are adjacent to the outer or unbound region of the diagram.  The labels are assigned to the arcs by travelling around the diagram and labelling arcs as we meet them. Note that the labelling of the arcs starts from 0 and each time we pass through a crossing it increases by one. The extended oriented Gauss code allows the unique encoding of a planar knotoid diagram, up to isotopy. For example, the planar knotoid diagram in Figure~\ref{fig:gc_example} corresponds to the following extended Gauss code:

\begin{center}
 1 -2 -1 2 \ \ ++\ \ 1 3
\end{center}

\subsection{A total order on Gauss codes.}\label{sec:order} Following \cite{gabrov, gabrov1}, we impose a total ordering on all (extended) oriented Gauss codes. The codes are ordered by taking into consideration the following:
\begin{enumerate}
\item[i.] The length of the Gauss word. 
\item[ii.] The ordering $-1 \ < \ 1 < \ -2 \ < \ \ldots \ < \ -n \ <\ n $ of the crossings of a diagram.
\item[iii.] The ordering $ - \ < \ +$ of the signs of the crossings.
\item[iv.] The length of the third part of the code (outer region).
\item[v.] The ordering $ 0 \ < \ 1 < \ 2 \ < \ \ldots $ of the arcs touching the outer region of the diagram.
\end{enumerate}
Note that steps ${\rm iv}$ and ${\rm v}$ are required only for the case of planar knotoids. For example, the code 
\begin{center}
-1 1\ \ +\ \ 1
\end{center} comes before the code 
\begin{center}
1 -2 -1 2 \ \ -\ - \ \ 1 3
\end{center} which, in turn, precedes 
\begin{center}
1 -2 -1 2 \ \ -\ - \ \ 0\ 2\ 3. 
\end{center}

 Finally, an (extended) oriented Gauss code, $G$, is called {\it reducible} if there exists a finite sequence of $\Omega$-moves that transforms it to a code $G^\prime$ of smaller order. The resulting code is called a {\it reduction} of $G$ and we say that $G$ accepts a reduction  \cite{gabrov, gabrov1}.  Note that $G$ and $G^\prime$ represent the same knotoid.

\section{Reidemeister moves and oriented Gauss codes} 

Using the oriented Gauss code, one can encode the application of $\Omega$-moves on a knotoid diagram in $S^2$. Let $w_i$, $i \in I$ be a subword of the Gauss word and let $A, B, C, \ldots$ be individual crossings. Then, applying or removing an $\Omega_1$-move has the following effect on a Gauss word:
\[
\pm A \mp A \ w \longleftrightarrow w
\]

There are two cases for the application or removal an $\Omega_2$-move on a Gauss word, which depend on the orientation of the arcs that take part in the $\Omega_2$-move move:
\begin{align*}
& \pm A \pm B \ w_1 \mp A \mp B \ w_2 \longleftrightarrow w_1 \ w_2\\
& \pm A \pm B \ w_1 \mp B \mp A \ w_2 \longleftrightarrow w_1 \ w_2
\end{align*}
Finally, we have the following four cases for the application of an $\Omega_3$-move:
{\small
\begin{align*}
& \pm A \pm B \ w_1 \mp A \pm C \ w_2 \mp B \mp C\ w_3 \longleftrightarrow \pm B \pm A \ w_1 \pm C \mp A \ w_2 \mp C \mp B\ w_3\\
& \pm A \pm B \ w_1  \mp A \pm C \ w_2 \mp C \mp B\ w_3 \longleftrightarrow \pm B \pm A \ w_1 \pm C \mp A \ w_2 \mp B \mp C\ w_3\\
& \pm A \pm B \ w_1  \mp C \pm A \ w_2 \mp B \mp C\ w_3 \longleftrightarrow \pm B \pm A \ w_1 \pm A \mp C \ w_2 \mp C \mp B\ w_3\\
& \pm A \pm B \ w_1  \mp C \pm A \ w_2 \mp C \mp B\ w_3 \longleftrightarrow \pm B \pm A \ w_1 \pm A \mp C \ w_2 \mp B \mp C\ w_3\\
\end{align*}}
The set of signs is affected in the following way: A positive or negative $\Omega_1$-move one an arc, adds or removes a negative or positive crossing, a $\Omega_2$-move adds or removes two consecutive opposing signs while a $\Omega_3$-move just reorders the set of signs.

\section{$\Omega$-moves, planar diagrams and the outer region} 

As mentioned in Section~\ref{sec:gc}, for the case of planar knotoids we use the extended oriented Gauss code that contains also the information of which arcs are adjacent to the unbounded region of the diagram. Therefore, when performing $\Omega$-moves on an extended oriented Gauss code, one has to keep track also of how the third part of the code is affected by the applied $\Omega$-move. Before discussing how the $\Omega$-moves affect the third part of the extended oriented Gauss code of a knotoid diagram, we will make a digression on how to can determine the local regions of a knotoid diagram.

\subsection{Local regions of diagrams}
If the over/undercrossing information of a knotoid diagram $k$ is ignored, one obtains a planar graph where the vertices correspond to the crossings of $k$ and the edges to the arcs of $k$. There are two additional vertices that correspond to the endpoints $k$, making the number of vertices equal to $n+2$, where $n$ is the number of crossings of $k$. This graph is called the underlying graph of $k$ \cite{gugumcu2017} and it will be denoted by $G$. In what follows, it will be very helpful to know which arcs of the diagram bound each of the regions of $G$ and also which arcs touch the outer region of the diagram, that is the unbounded region of the plane.

\begin{figure}[!htb]
\centering
\includegraphics{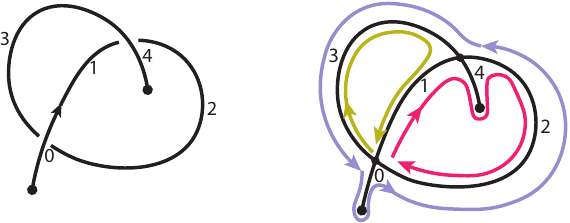}
\caption{Computing the arcs that bound the local regions of a knotoid diagram. The knotoid diagram is on the upper left corner of the figure while its underlying graph can be seen in the lower right.}\label{fig:regions}
\end{figure}

In order to determine the arcs of the local regions of a diagram, we consider $G$ and work as follows: Starting from an endpoint, travel around the graph and we label all arcs as we meet them. Next, pick a vertex that corresponds to a crossing and then pick an edge that is adjacent to that vertex. Move in a clockwise fashion and follow the closed path on the graph that loops back to the chosen edge. Each time a new edge is met, we note its label. If during the trip around the loop an endpoint is met, we go around it. Note that each edge is adjacent to two local regions, except the two edges that are directly connected to an endpoint. This process is repeated for all edges. For example, the regions of the diagram in Fig.~\ref{fig:regions} are:
\begin{align*}
&r_1: \ 0,\, 2,\, 3 \\
&r_2: \ 1,\, 4,\, 2 \\
&r_3:\ 3,\, 1
\end{align*}

Having knowledge of all local regions of a knotoid diagram, we can now start studying how the application of Reidemeister moves alters the third part of the extended oriented Gauss code. There is a number of cases to consider, especially for the second Reidemeister move and for this reason, we shall discuss this in a separate section.

If a Reidemeister move is applied to arcs that don't touch the outer region of the diagram then the third part of the extended oriented Gauss code (or, for simplicity, Gauss code) doesn't change. In what follows, we shall discuss how the application of a Reidemeister move on an arc that touches the outer region of the diagram affects the third part of a Gauss code. Recall that this part of a Gauss code contains, in increasing order, the labels of the arcs that touch the outer region of the diagram. Note also that for this, we consider the diagrams flat since the sign of a crossing doesn't contribute to the third part of a Gauss code.

\subsection{Reidemeister I} 
\begin{figure}[!htb]
\centering
\includegraphics{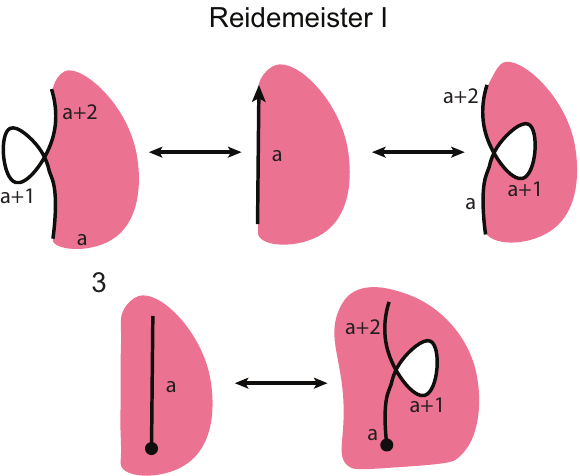}
\caption{Reidemeister I on outer regions. The shaded region corresponds to the outer region of the diagram.}\label{fig:RMI_outside}
\end{figure}

The $\Omega_1$-move creates or removes a crossing that divides an arc into three subarcs or merges them into a single one respectively. Consequently, the initial labels of the arcs  of the diagram following $a$ are shifted by 2 or -2 after the application of the move. There are two possible results for the third part of the Gauss code after the move. In the first case (Fig.~\ref{fig:RMI_outside}.I) the kink that is introduced by the move doesn't touch the outer region and so the third part, after renumbering the arcs, becomes:
\[ \ldots , a ,\ \ldots \longleftrightarrow  \ldots , a,\ a+2, \ldots 
\]
In  the second case the kink touches the outer region (Fig.~\ref{fig:RMI_outside}.II) or  an arc that contains an endpoint is involved (Fig.~\ref{fig:RMI_outside}.III)  and so we have that:
\[ \ldots ,  a ,\  \ldots \longleftrightarrow  \ldots , a,\ a+1,\ a+2, \ldots 
\]

\subsection{Reidemeister II} The application of the $\Omega_2$-move  divides each one of the two arcs into three subarcs. The initial labels of the arcs following arc $a$ shift by $\pm 2$ while the arcs following arc $b$ are shifted additionally by $\pm2$. For example, assume that the third part of a Gauss code contains the arcs $a \ b \ c $, with $a<b<c$, and that we apply an $\Omega_2$-move between arcs $a$ and $b$. The new arcs $a+1$ and $a+2$ that are introduced shift the labels of both of the arcs $b$ and $c$ by 2. The $\Omega_2$-move however splits, the already shifted, arc $b$ into three subarcs with labels $b+2$, $b+3$ and $b+4$. This means that the label of $c$ is shifted once more by 2 and therefore the third part of the Gauss code after the application of an $\Omega_2$-move becomes: $ a, \ \ a+1 ,\ \ a+2, \ \ b+2 ,\ \ b+3, \ \ b+4 ,\ \ c+4$.  

Depending on whether the arcs involved in the $\Omega_2$-move are parallel or anti-parallel, on whether they both touch the outer region or not, and also on whether they include endpoint(s) of the diagram or not, we distinguish four different cases for the $\Omega_2$-move which we discuss separately below.
 \begin{figure}[!htb]
\centering
\includegraphics{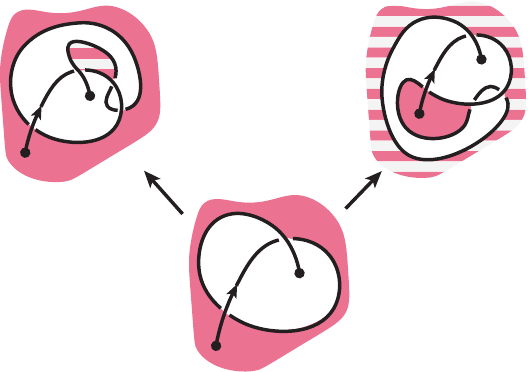}
\caption{The two different choices for the outside region of a knotoid diagram after an $\Omega_2$-move.}\label{fig:two_regions}
\end{figure}
Note that in all cases but one that are discussed, the application of a crossing increasing $\Omega_2$-move  creates two different regions, that both can be chosen as the outer region of the knotoid diagram. This, in turn, yields two different (non-minimal) diagrams respectively (see Fig.~\ref{fig:two_regions}). When  a crossing decreasing $\Omega_2$-move is applied, these two regions merge to a single region.
 \begin{figure}[!htb]
\centering
\includegraphics{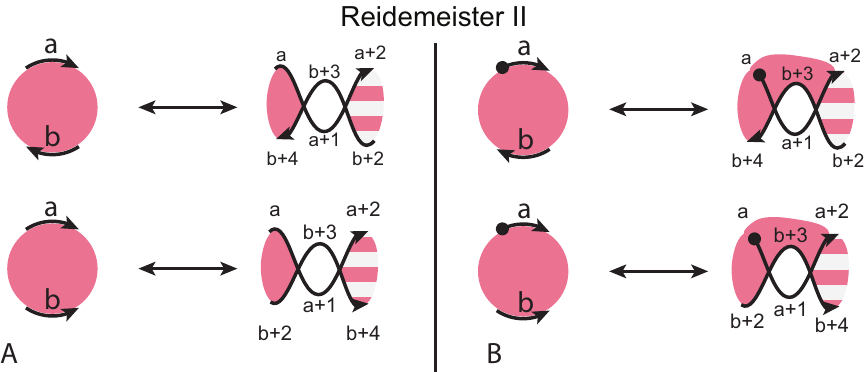}
\caption{The Reidemeister II move and how it is affecting the outer region of a knotoid diagram. The shaded region in the left-hand side indicates the outer region of the diagram. (A) The two arrows indicate the two arcs touching the outer region that participate in the RII move. In the right hand side the two different shadings indicate the two possible choices for the outer region for each case. (B) Same as before with the difference that the big black dot corresponds to an endpoint. Observe that a region bleeds over the endpoint.}\label{fig:RMII_both}
\end{figure}

\noindent {\bf Case 1: Both arcs of touch the outer region}. We distinguish two cases that depend on whether the arcs are parallel or antiparallel (Fig.\ref{fig:RMII_both}(A)). 
\begin{enumerate}
\item[i.] If the arcs are anti-parallel we have the following two choices:
\begin{align*}
& \ldots ,\ a,\ b,\ \ldots \longleftrightarrow \ldots, \ a, \  b+4, \ \ldots\\
& \ldots ,\ a,\ b,\ \ldots \longleftrightarrow \ldots, \ a+2, \  b+2, \ \ldots\\
\end{align*}

\item[ii.] If the arcs are parallel, then we have:
\begin{align*}
& \ldots ,\ a,\ b,\ \ldots \longleftrightarrow \ldots, \ a, \  b+2, \ \ldots\\
& \ldots ,\ a,\ b,\ \ldots \longleftrightarrow \ldots, \ a+2, \  b+4, \ \ldots\\
\end{align*}
\end{enumerate}

\noindent {\bf Case 2: Both arcs of touch the outer region and one of the arcs contains an endpoint}. The two candidate outer regions for this cases are described in the following two cases (see Fig.\ref{fig:RMII_both}(B)).
\begin{enumerate}
\item[i.] If the arcs are anti-parallel the two potential outer regions are:
\begin{align*}
& \ldots ,\ a,\ b,\ \ldots \longleftrightarrow \ldots, \ a, \ a+2 \ b+3,\  b+4, \ \ldots\\
& \ldots ,\ a,\ b,\ \ldots \longleftrightarrow \ldots, \ a+2, \  b+2, \ \ldots\\
\end{align*}
\item[ii.] If the arcs are parallel, then we have:
\begin{align*}
& \ldots ,\ a,\ b,\ \ldots \longleftrightarrow \ldots, \  a, \ a+2,\  b+2, \ b+3,\ \ldots\\
& \ldots ,\ a,\ b,\ \ldots \longleftrightarrow \ldots, \ a+2, \  b+4, \ \ldots\\
\end{align*}
\end{enumerate}
 \begin{figure}[!htb]
\centering
\includegraphics{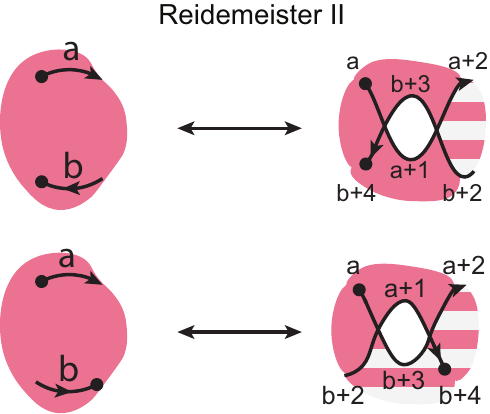}
\caption{Reidemeister II on outer regions (cont.). The two cases correspond to the possible orientations (parallel or anti-parallel) of the endpoints. }\label{fig:RMII_both_arc_end}
\end{figure}

\noindent {\bf Case 3: Both arcs of touch the outer region and both contain endpoints}. Once again we distinguish the cases where the two arcs are parallel or anti-parallel (see also Fig.\ref{fig:RMII_both_arc_end}).
\begin{enumerate}
\item[i.] If the arcs are anti-parallel the two potential outer regions are:
\begin{align*}
& \ldots ,\ a,\ b,\ \ldots \longleftrightarrow \ldots, \ a, \ a+1,\ a+2,\ b+2,\ b+3,\  b+4\ \ldots\\
& \ldots ,\ a,\ b,\ \ldots \longleftrightarrow \ldots, \ a+2, \  b+2, \ \ldots\\
\end{align*}
\item[ii.]
\begin{align*}
& \ldots ,\ a,\ b,\ \ldots \longleftrightarrow \ldots, \ a, \ a+1,\ a+2,\ b+2,\ b+3,\ \ldots\\
& \ldots ,\ a,\ b,\ \ldots \longleftrightarrow \ldots, \ a+1,\ a+2, \ b+2,\ b+4,\  \ldots\\
\end{align*}
\end{enumerate}

 \begin{figure}[!htb]
\centering
\includegraphics{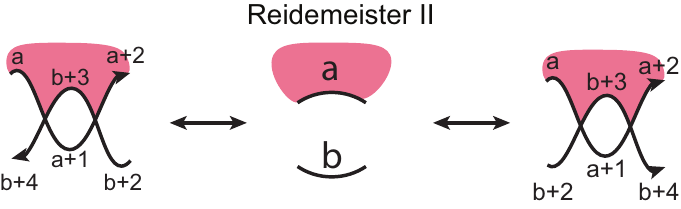}
\caption{Reidemeister II on outer regions (cont.). The RII move between an arc that touches the outer region and one that lies in the inner part of the diagram. The result is the same if the arcs are parallel or anti-parallel.}\label{fig:RMII_only_one_arc}
\end{figure}

\noindent {\bf Case 4: Only one arc touches the outer region} . 
In this last case there is only one possible outcome after applying an $\Omega_2$-move (see also Fig.\ref{fig:RMII_only_one_arc}):
\[
\ldots a, \ b, \ \ldots \longleftrightarrow \ldots a , \ a+2,  \ b+3, \  \ldots \\
\]

Finally, we discuss how the third Reidemeister move affects the outer region of a diagram. 

 \begin{figure}[!htb]
\centering
\includegraphics{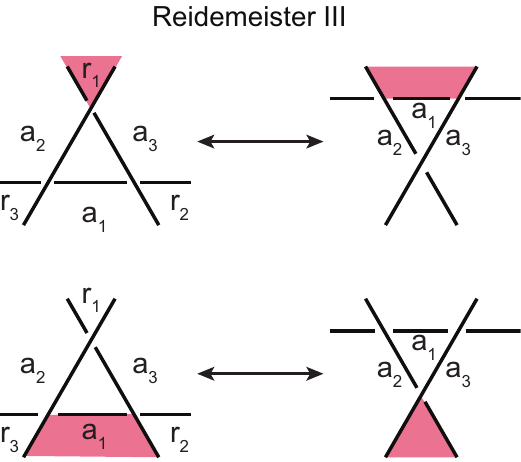}
\caption{Reidemeister III move and the outer region. In this example, the arc the moves during the RIII move is the arc $a_1$.}\label{fig:RMIII_outside}
\end{figure}

\subsection{Reidemeister III} Notice that the RIII move involves three arcs of the diagram that form a triangular region and let $a_1$, $a_2$, $a_3$ be those arcs. Opposite from an arc $a_i$ and adjacent to the intersection of the other two arcs, lies a local region of the diagram, $r_i$ , $i=1,2,3$ (see Fig.~\ref{fig:RMIII_outside}). If the region $r_i$ is a subset of the outer region of the diagram or not and if the arc $a_i$ is in the outer region or not, we have the following two cases:
\begin{align*}
& \mbox{if}\ r_i \subseteq \mbox{outer region} \longleftrightarrow \mbox{outer region} + { a_i } \\
& \mbox{if} \ r_i \not\subseteq \mbox{outer region and}\ a_i \in \mbox{outer region} \longleftrightarrow  \mbox{outer region} \setminus { a_i } \\
\end{align*}

\subsection{Gauss codes and involutions}
Applying each of the involutions on a knotoid diagram has the following effects on its Gauss code.

\noindent\textbf{Reversion:} The reversion involution reverses and renumbers the Gauss word and reverses the list of signs. Finally, it reverses and renumbers the list of arcs the are adjacent to the outer region:
\begin{center}
1 -2 2 3 -1 -3 \ \ +-\-- \ \ 2 \ \ $\longleftrightarrow$ \ \ -1 -2 1 3 -3 2 \ \ -\--+ \ \ 4
\end{center}
\textbf{Mirror reflection:} The mirror reflection involution changes the under-crossings to over-crossings and vice versa as well as the signs of the crossings in the second part of the Gauss code:
\begin{center}
1 -2 -1 2 \ \ ++ \ \ 0 2 3 \ \ $\longleftrightarrow$ \ \ -1 2 1 -2 \ \ -\-- \ \ 0 2 3
\end{center}
\textbf{Symmetry:} The symmetric involutions changes the signs of the crossings in the Gauss code.
\begin{center}
1 -2 -1 2 \ \ ++ \ \ 0 2 3 $\longleftrightarrow$ \ \ -1 2 1 -2 \ \ -\-- \ \ 0 2 3
\end{center}
\textbf{Rotation:} The rotation involution changes under-crossings to over-crossings and vice versa.
\begin{center}
1 -2 -1 2 \ \ ++ \ \ 0 2 3 $\longleftrightarrow$ \ \ -1 2 1 -2 \ \ ++\ \ 0 2 3
\end{center}
 
  \section{Classification algorithm}\label{sec:algorithm}
  Recall that knotoids are equivalence classes of knotoid diagrams, where the equivalence relation is generated by isotopy and $\Omega$-moves away from the endpoints. Without loss of generality, we can assume that the representative diagram of each class is  {\it minimal}, that is, it is a diagram with the minimal number of crossings among all diagrams within the same equivalence class.
  
A knotoid is called {\it composite}, if it can be written as a product of two other knotoids  \cite{turaev2012}. Knotoids that are not composite shall be called {\it primes}. Our goal is to classify, up to all involutions, all extended oriented Gauss codes that correspond to planar knotoid diagrams with up to 5 crossings and determine a unique representative diagram for each isotopy class. The representative is defined as the smallest diagram within a class in terms of the order of Section~\ref{sec:order}  . From this point on, abusing notation, we use the term prime for both the isotopy class and its representative diagram.  Our classification algorithm was implemented in python  $\mathtt{python\ 2.7}$. We also took the opportunity and applied the same algorithm to generate all prime knotoids in $S^2$ with up to 6 crossings. 
 
\subsection{Generate all extended oriented Gauss codes.} We start by generating all possible oriented Gauss codes, that is, all possible tuples whose first entry is a word of length $2n$ in the alphabet ${ \pm1, \ \pm2 , \ \ldots \ \pm n}$ and the second entry is a word of length $n$ in the alphabet ${ -,\ + }$. In total, there are $(2n)! \cdot 2^n$ possible oriented Gauss codes but not all of them represent a knotoid diagram. Such Gauss codes will be called {\it nonrealizable}. The realizability criterion that we followed is based on \cite{realize}. 
 
 For each of the realizable oriented Gauss codes, we determine the local regions of their corresponding knotoid diagram and each time we pick a different one as its outer region. In this way, one obtains all realizable extended oriented Gauss codes. Since a knotoid diagram with $n$ crossings has $n+1$ local regions, the total number of realizable Gauss codes is $(\#\ \mbox{realizable oriented}$ $\mbox{ Gauss codes})\ \times \ (n+1)$. This is summariazed in Table~\ref{table:planar}. For the case of planar knotoids with up to five crossings we have found in total 832904 realizable diagrams, while for the case of knotoids in $S^2$ with up to six crossings there are 2363766. Since a realizable Gauss code corresponds to a diagram, from now on we will use the terms ``(planar) diagram'' and ``(extended) oriented Gauss code'' interchangeably.

\begin{table}[]
\begin{tabular}{ccllcc}
\multicolumn{2}{c}{Planar}                                      &  &                       & \multicolumn{2}{c}{In $S^2$}                                   \\ \cline{1-2} \cline{5-6} 
\multicolumn{1}{|l|}{Crossings} & \multicolumn{1}{l|}{Diagrams} &  & \multicolumn{1}{l|}{} & \multicolumn{1}{l|}{Crossings} & \multicolumn{1}{l|}{Diagrams} \\ \cline{1-2} \cline{5-6} 
\multicolumn{1}{|c|}{1}         & \multicolumn{1}{r|}{8}        &  & \multicolumn{1}{l|}{} & \multicolumn{1}{c|}{1}         & \multicolumn{1}{r|}{4}        \\ \cline{1-2} \cline{5-6} 
\multicolumn{1}{|c|}{2}         & \multicolumn{1}{r|}{120}      &  & \multicolumn{1}{l|}{} & \multicolumn{1}{c|}{2}         & \multicolumn{1}{r|}{40}       \\ \cline{1-2} \cline{5-6} 
\multicolumn{1}{|c|}{3}         & \multicolumn{1}{r|}{2112}     &  & \multicolumn{1}{l|}{} & \multicolumn{1}{c|}{3}         & \multicolumn{1}{r|}{528}      \\ \cline{1-2} \cline{5-6} 
\multicolumn{1}{|c|}{4}         & \multicolumn{1}{r|}{39840}    &  & \multicolumn{1}{l|}{} & \multicolumn{1}{c|}{4}         & \multicolumn{1}{r|}{7968}     \\ \cline{1-2} \cline{5-6} 
\multicolumn{1}{|c|}{5}         & \multicolumn{1}{r|}{781824}   &  & \multicolumn{1}{l|}{} & \multicolumn{1}{c|}{5}         & \multicolumn{1}{r|}{130304}   \\ \cline{1-2} \cline{5-6} 
Total                                &    \multicolumn{1}{r}{832904}                          &  & \multicolumn{1}{l|}{} & \multicolumn{1}{c|}{6}         & \multicolumn{1}{r|}{2224922}  \\ \cline{5-6} 
                           &                         &  &                       & Total                          & \multicolumn{1}{r}{2363766}  
\end{tabular}

\caption{Realizable extended oriented Gauss codes.}\label{table:planar}
\end{table}

\subsection{Isotopy classes of knotoid diagrams}\label{sec:isotopy}
This step describes the algorithm that determines the different isotopy classes of knotoid diagrams with up to 5 crossings for the case of planar knotoids. In what follows, we first describe the building blocks of the algorithm before discussing the algorithm itself. The case of $S^2$-knotoids with up to 6 crossings is analogous, except where specifically mentioned.

\noindent \textbf{1. Partitioning.} For the case of planar knotoids, we partition the set of all diagrams using both the loop arrow polynomial and their double-branched cover, since this way we obtain a finer partition. In fact, even though the double-branched cover is in principle a stronger invariant than the loop arrow polynomial, there are pairs of knotoids that are distinguished only by the latter invariant. More details on this are provided in Section~\ref{sec:comp_inv}. Each subset of the partition contains all Gauss codes that correspond to all planar diagrams having the same loop arrow polynomial and the same double-branch cover. The case of knotoids in $S^2$ is simpler, as the arrow polynomial is sufficient to obtain a full classification. For knotoids in $S^2$, each subset of the partition contains all diagrams with the same arrow polynomial. The loop arrow polynomial and the arrow polynomial for each diagram have been computed using Knoto-ID \cite{knoto-id}. The computation of the double-branched cover for planar diagrams has been implemented in \texttt{python 2.7}. 

\noindent \textbf{2. Reachability Graph.} For each subset of the partition we create a {\it reachability graph}. The reachability graph is an undirected graph where each node corresponds to a diagram, while an edge between two nodes means that the corresponding diagrams are related either via a $\Omega$-move or via reversion. There is a one-to-one correspondence between connected components of the graph and isotopy classes, assuming that all possible sequences of $\Omega$-moves and reversions have been added to the reachability graph.

The graph is initially populated with nodes corresponding to all planar diagrams with up to 5 crossings. Edges corresponding to all possible crossing-decreasing $\Omega_1$-moves and $\Omega_2$-moves, all possible $\Omega_3$-moves as well as reversions are then added to the reachability graph. 

\noindent \textbf{3. Random walk.} After step 2, the reachability graph doesn't contain any edges that correspond to crossing-increasing $\Omega$-moves. However, to connect all diagrams in the same isotopy class, one has to consider also crossing-increasing moves. Since it is not possible to exhaustively enumerate all possible Reidemeister moves, the algorithm performs a random walk of predetermined length in the graph by applying a randomly chosen sequence of $\Omega$-moves, including crossing-increasing moves. The length of the random walk is empirically set to 10,000 steps. Note that throughout the random walk, diagrams with more than five crossings can be created. Whenever this happens, an edge is added only when the random walk reaches a diagram with up to five crossings, summarizing a sequence of $\Omega$-moves. Since the increase in the number of crossings happens through the application of Reidemeister moves, we are ensured that we never leave the isotopy class. Furthermore, consecutive applications of crossing-increasing moves will not contribute to finding a path connecting two components and so we bias our approach by giving higher probability of application to crossing-reducing $\Omega_1$ and $\Omega_2$ moves than to their crossing-increasing counterparts. The $\Omega_3$-move gets the same probability of application as the crossing reducing moves because applying an $\Omega_3$-move will never result in a diagram with increased number of crossings. The ratio of crossing reducing/preserving moves to crossing increasing moves is set empirically to 2:1.  Once the random walk reaches its predefined length, the algorithm attempts to reduce the resulting diagram by systematic application of crossings-reducing $\Omega$-moves. If the reduction of the diagram leads to a diagram already in the graph, then an edge between the corresponding nodes of the graph is added. Note that the algorithm does not add nodes to the graph for any intermediate step in the random walk that corresponds to a diagram with more than five crossings since this would be computationally too costly.

\noindent \textbf{4. Diagram status.} In this step, we go through all subsets of the partition and we determine the connected components of each corresponding graph. For each connected component, if the component contains a node marked as composite (see next step), we mark all nodes of the component as composite. Otherwise, we mark as candidate prime the smallest diagram in terms of the order of Section~\ref{sec:order} and all other diagrams within the connected component are ignored. 

If the corresponding reachability graph contains a unique connected component, it is marked as {\it terminated} and it is ignored for the rest of algorithm. If it doesn't contain a composite diagram, the smallest diagram in terms of the order of Section~\ref{sec:order} is marked as prime and all other diagrams are marked as not prime.

\begin{figure}[!htb]
\centering
\includegraphics{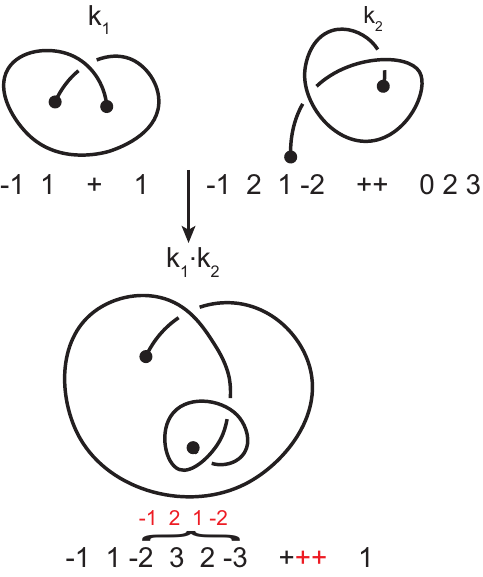}
\caption{The composition of two planar knotoids and the corresponding extended oriented Gauss codes. We note with red  the labels of the crossings of the second diagram, before renumbering, as well as it corresponding crossings. The outer region of the composite diagram is the same as the outer region of the first diagram.}\label{fig:composite}
\end{figure}

\noindent \textbf{5. Composite diagrams.} Since the double-branch cover detects the trivial planar knotoid \cite[Theorem~1.8]{barbensi}, we know that the corresponding reachability graph will have a single connected component corresponding to the set of all trivial planar diagrams with up to 5 crossings. For the case of knotoids in $S^2$, since there is no equivalent result for the arrow polynomial available, we repeatedly perform Step 3 on the subset of the partition containing the trivial knotoid diagram, until we obtain a single connected component corresponding to the set of all trivial diagrams in $S^2$ with up to 6 crossings. Next, ignoring the reachability graph corresponding to trivial knotoids, we consider all possible candidate primes throughout all partitions and all possible products of diagrams that result in a composite knotoid diagram with 5 or less crossings.  Not all combinations of planar knotoids yield a valid composite diagram. In fact, at least one diagram of the product must have one endpoint in the outside region of the diagram. Once we have a viable combination of planar diagrams, it  is straightforward  to produce a composite knotoid using the (extended) oriented Gauss codes since one has to concatenate the Gauss codes, relabel the crossings and arcs of the second diagram and recompute the outer region of the diagram (See Figure~\ref{fig:composite}).

 The case of $S^2$-knotoids is simpler since any knotoid in $S^2$ can be represented by a diagram that has one endpoint in the unbound region of the diagram \cite{turaev2012}. A pair of planar  knotoids yields up to four diagrams while a pair of $S^2$-knotoids yields exactly four diagrams. 

\noindent \textbf{6. The algorithm.} The algorithm for the determination of all isotopy classes of knotoid diagrams starts by partitioning the set of all diagrams using polynomial invariants (Step 1). Reachability graphs are then created for each subset of partition (Step 2). An initial round of random walks is performed for each subset of partition, starting one random walk for each node in the graph (Step 3). For each subset of partition, diagrams are then marked as either prime or candidate prime (Step 4). The goal of this initial round of random walks is to reduce the number of candidate primes, so as to minimize the computational cost for the next step, which is the evaluation of composite diagrams (Step 5). After this step, the algorithm repeatedly applies rounds of random walks (Step 3) for all reachability graph not marked as terminated, followed by a determination of the status of each diagram (Step 4). The algorithm terminates when each subset of the partition are marked as terminated in Step 4, i.e. all graphs have a unique connected component.

\subsection{Prime diagrams up to all involutions.} So far, the algorithm has focused on determining all isotopy classes of planar knotoid diagrams with up to five crossings. In this step we consider the isotopy class of a knotoid up to all of its involutions. 

We start by merging the graphs corresponding to subsets of the partition related by mirror reflection,  symmetry and mirror symmetric reflection. Note that the reversion involution has been already considered in Step 2 of Section~\ref{sec:isotopy}. For each node, edges corresponding to all involutions are then added to the reachability graph, thus connecting isotopy classes related by involutions. We finally apply Step 4 in order to obtain the status of each diagram from each reachability graph and then we  extract the list of primes.

\subsection{Knotoid naming conventions}
In this final step, we name the representatives of all distinct isotopy classes, up to involutions. All prime (planar) knotoids are labelled using the scheme $X_Y$, where $X$ is the minimal number of crossings of the knotoid and $Y$ is its relative position among all knotoids with the same number of crossings. Note that, since several different planar knotoids may correspond to the same knotoid in $S^2$, it is possible that diagrams may change label depending on whether they are considered on the plane or in $S^2$. Moreover, the order of Section~\ref{sec:order} doesn't guarantee that knot-type knotoids will get the label that their corresponding knots have in the  Rolfsen table notation \cite{Adams}. In order to avoid this, we impose the following additional order on the representative of the isotopy classes of planar knotoid diagrams:
 \[
\mbox{Knot-type knotoids} \ <\  \mbox{proper knotoids}  \ <\  \mbox{non-proper knotoids}.
\]
Since any $S^2$-knotoid that is not a knot-type can be represented by a proper knotoid \cite{turaev2012}, our naming scheme is consistent even if we restrict to the case of $S^2$ knotoids. Knot-type diagrams follow the labels of the Rolfsen table, while the diagrams within each of the other two categories follow the order of Section~\ref{sec:order}. 

%
The table of all planar knotoids with up to 5 crossings can be found in Appendix ~\ref{appendix:plane}, while the table of all knotoids in $S^2$ with up to 6 crossings can be found in Appendix ~\ref{appendix:s2}.  All knotoid diagrams were plotted using a modified version of the knot plotting function found in SageMath \cite{sage}.

 \section{Results}

 \subsection{Planar knotoids}

 \begin{table}[h]
\[
\begin{array}{l@{}  l}
\begin{array}{l@{{}:{}}@{\quad}L}
5_7& -1\ -2\ 3\ 4\ -3\ 2\ -5\ 1\ 5\ -4\ \ -\ -\ -\ +\ +\ \ 0\ 7\ 8 \\  5_{421}& -1\ -2\ 3\ 4\ -3\ 2\ -5\ 1\ 5\ -4\ \ -\ -\ -\ +\ +\ \ 3\ 4\ 10\\
\end{array}
\\ \vspace{0.2cm} \\
\begin{array}{l@{{}:{}}@{\quad}L}
 5_{9}& -1\ 2\ -3\ 1\ -4\ 5\ -2\ 3\ 4\ -5\ \ -\ -\ -\ +\ +\ \ 0\  3\  4\ 8\\ 5_{561}& -1\ 2\ -3\ 1\ -4\ 5\ -2\ 3\ 4\ -5\ \ -\ -\ -\ +\ +\ \ 2\ 5\ 6\ 8\ 10
\end{array}
\vspace{.2cm} \\ \\
\begin{array}{l@{{}:{}}@{\quad}L}
5_{12}& -1\ 2\ -3\ 1\ 4\ -5\ -2\ 3\ -4\ 5\ \ -\ -\ -\ -\ -\ \ 0\ 3\ 4\ 8 \\ 5_{593}& -1\ 2\ -3\ 1\ 4\ -5\ -2\ 3\ -4\ 5\ \ -\ -\ -\ -\ -\ \ 2\ 5\ 6\ 8\ 10
\end{array}\vspace{.2cm} \\ \\
\begin{array}{l@{{}:{}}@{\quad}L}
5_{19}& -1\ 2\ -3\ 4\ -5\ 1\ -2\ 3\ 5\ -4\ \ -\ -\ -\ +\ +\  \ 3\ 4\ 8\ 10 \\ 5_{796}& -1\ 2\ -3\ 4\ -5\ 1\ -2\ 3\ 5\ -4\ \ -\ -\ -\ +\ +\  0\ 2\ 5\ 6\ 8 
\end{array}
\\ \vspace{0.2cm} \\
\begin{array}{l@{{}:{}}@{\quad}L}
5_{21}& -1\ 2\ -3\ 4\ -5\ 1\ 5\ -2\ -4\ 3\	 \ -\ +\ -\ -\ +\ \ 0\ 5\ 6 \\ 5_{814}& -1\ 2\ -3\ 4\ -5\ 1\ 5\ -2\ -4\ 3\  \ -\ +\ -\ -\ +\ \ 2\ 3\ 8\ 10
\end{array}
\\ \vspace{0.2cm} \\
\begin{array}{l@{{}:{}}@{\quad}L}
5_{24}& -1\ 2\ -3\ 4\ 5\ -4\ -2\ 1\ 3\ -5\ \ -\ -\ -\ -\ +\ 4\ 5\ 10 \\ 5_{891}& -1\ 2\ -3\ 4\ 5\ -4\ -2\ 1\ 3\ -5\ \ -\ -\ -\ -\ +\ \ 0\ 2\ 7\ 8
\end{array}

 \end{array}
 \]
 \caption{The six unresolved pairs.}\label{table:unresolved}
\end{table}
 
 In this case our algorithm doesn't meet the termination criterion as there are 6 subsets within the partition, each containing a graph with two connected components. The unresolved cases can be seen in Table~\ref{table:unresolved}. There are a few interesting observations to be made here. 

We observe that in each pair, both diagrams represent the same knotoid in $S^2$, because they have the same  oriented Gauss code. However, they have different extended oriented Gauss codes, since each time a different region of the $S^2$-knotoid is marked as the outer region of the planar diagram.  

\begin{table}[h]
\begin{tabular}{crrrr}
\multicolumn{5}{c}{Planar Knotoids}                                                              \\ \hline
\multicolumn{1}{|l|}{Crossings} & \multicolumn{1}{l|}{Primes} & \multicolumn{1}{l|}{Rotatable} & \multicolumn{1}{l|}{Achiral} & \multicolumn{1}{l|}{Strongly Achiral}  \\ \hline
\multicolumn{1}{|c|}{1}         & \multicolumn{1}{r|}{1}      & \multicolumn{1}{r|}{1}&  \multicolumn{1}{r|}{0}     & \multicolumn{1}{r|}{0}     \\ \hline
\multicolumn{1}{|c|}{2}         & \multicolumn{1}{r|}{6}        & \multicolumn{1}{r|}{3} & \multicolumn{1}{r|}{1}     & \multicolumn{1}{r|}{0}   \\ \hline
\multicolumn{1}{|c|}{3}         & \multicolumn{1}{r|}{26}       & \multicolumn{1}{r|}{12} & \multicolumn{1}{r|}{0}     & \multicolumn{1}{r|}{0}     \\ \hline
\multicolumn{1}{|c|}{4}         & \multicolumn{1}{r|}{154}      & \multicolumn{1}{r|}{41} & \multicolumn{1}{r|}{5}    & \multicolumn{1}{r|}{1}      \\ \hline
\multicolumn{1}{|c|}{5}         & \multicolumn{1}{r|}{950}      & \multicolumn{1}{r|}{163} & \multicolumn{1}{r|}{0}  & \multicolumn{1}{r|}{0}     \\ \hline
Total                           & 1137       & 220 &               6 & 1                                   
\end{tabular}
\caption{Prime planar knotoids  with up to 5 crossings, assuming that all entries in Table~\ref{table:unresolved} are non-isotopic pairs.}\label{table:inplane}
\end{table}

Consecutive applications of longer random walks didn't reveal a connection between the two components in any of these cases.  For this reason, we hypothesize that each of these pairs represent a pair of non-isotopic planar knotoid diagrams. Under this hypothesis, the summary of the tabulation of all 5-crossing planar knotoids is shown in Table~\ref{table:inplane} and in Appendix ~\ref{appendix:plane}.

Strongly achiral is, as expected, the knot-type knotoid $\mathbf{4_1}$ while achiral are the knotoids $\mathbf{2_2}$, $\mathbf{4_{22}}$, $\mathbf{4_{84}}$, $\mathbf{4_{103}}$, $\mathbf{4_{148}}$. Table~\ref{table:rotatable} lists all rotatable planar knotoids per number of crossings.

\begin{table}[h]
\begin{tabular}{l|l|l|l}
Crossings & \multicolumn{3}{c}{Rotatable knotoids} \\ \hline
1         & \multicolumn{3}{l}{$1_1$} \\ \hline
2         & \multicolumn{3}{l}{$2_3$, $2_5$,  $2_6$} \\ \hline
3         & \multicolumn{3}{l}{\begin{tabular}[c]{@{}l@{}}$3_1$, $3_9$, $3_{11}$, $3_{12}$, $3_{15}$, $3_{16}$, $3_{21}$, $3_{22}$, $3_{23}$, $3_{24}$, $3_{25}$, $3_{26}$\end{tabular}}                                                                                                                                                                                                                                                                                                                                                                                                                                                                                                                                                                                                                                                                                                                                                                                                                                                                                                                                                                                                                                                                                                                                                                                                                                                                                                                                                                                                                                                                                                                                                                                                                                                                                                                                                                                                                                           \\ \hline
4         & \multicolumn{3}{l}{\begin{tabular}[c]{@{}l@{}}$4_{27}$, $4_{28}$, $4_{41}$, $4_{42}$, $4_{43}$, $4_{62}$, $4_{63}$, $4_{64}$, $4_{65}$, $4_{72}$, $4_{73}$, $4_{81}$, \\ $4_{92}$, $4_{93}$, $4_{94}$, $4_{95}$, $4_{96}$, $4_{101}$, $4_{109}$, $4_{113}$, $4_{114}$, $4_{115}$, $4_{116}$, \\ $4_{129}$, $4_{130}$, $4_{131}$, $4_{132}$, $4_{139}$ $4_{140}$, $4_{141}$, $4_{142}$, $4_{143}$, $4_{144}$, \\ $4_{145}$, $4_{146}$, $4_{149}$, $4_{150}$, $4_{151}$ $4_{152}$, $4_{153}$, $4_{154}$\end{tabular}}                                                                                                                                                                                                                                                                                                                                                                                                                                                                                                                                                                                                                                                                                                                                                                                                                                                                                                                                                                                                                                                                                                                                                                                                                                                                                                                                                                                                                                                                                                                                                      \\ \hline
 5         & \multicolumn{3}{l}{\begin{tabular}[c]{@{}l@{}}$5_1$, $5_2$, $5_{99}$, $5_{100}$, $5_{146}$, $5_{167}$, $5_{168}$, $5_{183}$, $5_{184}$, $5_{185}$, $5_{195}$, \\ $5_{196}$, $5_{197}$, $5_{198}$, $5_{202}$, $5_{203}$, $5_{229}$, $5_{268}$, $5_{270}$, $5_{328}$, $5_{329}$, \\ $5_{332}$, $5_{333}$, $5_{344}$, $5_{345}$, $5_{348}$, $5_{366}$, $5_{369}$, $5_{370}$, $5_{371}$, $5_{372}$, \\ $5_{383}$, $5_{384}$, $5_{385}$, $5_{388}$, $5_{393}$, $5_{394}$, $5_{395}$, $5_{396}$, $5_{397}$, $5_{398}$, \\ $5_{399}$, $5_{400}$, $5_{413}$, $5_{449}$, $5_{450}$, $5_{451}$, $5_{452}$, $5_{458}$, $5_{459}$, $5_{544}$, \\ $5_{545}$, $5_{552}$, $5_{553}$, $5_{554}$, $5_{563}$, $5_{565}$, $5_{575}$, $5_{576}$, $5_{577}$, $5_{578}$, \\ $5_{579}$, $5_{580}$, $5_{588}$, $5_{607}$, $5_{608}$, $5_{611}$, $5_{679}$, $5_{680}$, $5_{698}$, $5_{700}$, \\ $5_{702}$, $5_{706}$, $5_{707}$, $5_{708}$, $5_{709}$, $5_{721}$, $5_{722}$, $5_{725}$, $5_{726}$, $5_{733}$, \\ $5_{734}$, $5_{735}$, $5_{736}$, $5_{737}$, $5_{738}$, $5_{739}$, $5_{740}$, $5_{760}$, $5_{761}$, $5_{762}$, \\ $5_{763}$, $5_{768}$, $5_{769}$, $5_{770}$, $5_{771}$, $5_{772}$, $5_{773}$, $5_{774}$, $5_{777}$, $5_{790}$, \\ $5_{792}$, $5_{800}$, $5_{801}$, $5_{806}$, $5_{807}$, $5_{838}$, $5_{839}$, $5_{840}$, $5_{841}$, $5_{842}$, \\ $5_{843}$, $5_{855}$, $5_{856}$, $5_{857}$, $5_{858}$, $5_{871}$, $5_{872}$, $5_{873}$, $5_{874}$, $5_{875}$, \\ $5_{876}$, $5_{877}$, $5_{878}$, $5_{879}$, $5_{880}$, $5_{881}$, $5_{882}$, $5_{883}$, $5_{884}$, $5_{885}$, \\ $5_{886}$, $5_{898}$, $5_{899}$, $5_{900}$, $5_{907}$, $5_{908}$, $5_{909}$, $5_{910}$, $5_{911}$, $5_{912}$, \\ $5_{913}$, $5_{914}$, $5_{928}$, $5_{931}$, $5_{932}$, $5_{933}$, $5_{934}$, $5_{936}$, $5_{937}$, $5_{938}$, \\ $5_{939}$, $5_{940}$, $5_{941}$, $5_{942}$, $5_{943}$, $5_{944}$, $5_{945}$, $5_{946}$, $5_{947}$, $5_{948}$, \\ $5_{949}$, $5_{950}$  \end{tabular}}
\end{tabular}
\caption{All rotatable planar knotoids.}\label{table:rotatable}
\end{table}

\subsubsection{Comparing invariants}\label{sec:comp_inv}
In this section we take the opportunity to compare the loop arrow polynomial to the method of double-branched covers of knotoids. To do so, we count the number of isotopy classes that each of these methods independently distinguish.

\begin{table}[h]
\begin{tabular}{lll}
\hline
\multicolumn{1}{|l|}{}                                                                        & \multicolumn{1}{c|}{\begin{tabular}[c]{@{}c@{}}Loop arrow\\ polynomial\end{tabular}} & \multicolumn{1}{c|}{\begin{tabular}[c]{@{}c@{}}Double-branched\\ cover\end{tabular}} \\ \hline
\multicolumn{1}{|l|}{\# isotopy classes}                                                      & \multicolumn{1}{c|}{916}                                                             & \multicolumn{1}{c|}{1121}                                                            \\ \hline
\multicolumn{1}{|c|}{\begin{tabular}[c]{@{}c@{}}\# not distinguished\\ diagrams\end{tabular}} & \multicolumn{1}{c|}{221}                                                             & \multicolumn{1}{c|}{16}                                                              \\ \hline
                                                                                              &                                                                                      &                                                                                      \\
                                                                                              &                                                                                      &                                                                                      \\
                                                                                              &                                                                                      &                                                                                     
\end{tabular}
\caption{Comparing the strength of the loop arrow polynomial versus the double-branched cover.}\label{table:compar}
\end{table}

As shown in Table~\ref{table:compar}, the number of isotopy classes distinguished by the method of branched covers is equal to 1121, while those that are distinguished by the loop arrow polynomial is 916. The number of isotopy classes are distinguished by both invariants is 912. This suggests that the method of double branched covers is clearly stronger than the loop arrow polynomial. However, there are 4 diagrams that distinguished by the loop arrow but not by the double branched covers method.

In fact, the knotoids $\mathbf{2_5}$ and $\mathbf{5_{146}}$  lift to knots in the solid torus that have the same Jones polynomial:
\[
-A^{16}v + A^{12}v^3 - 2A^{8}v^3 + 2A^{8}v + A^{4}v^{3} 
\]
and so they are not distinguished by the double-branched cover. The loop arrow polynomial, however, manages to distinguish them:
{\small
\begin{align*}
\mathbf{2_5}: \ &  - A^{2}v + A^{6}v + A^{8}\\
\mathbf{5_{146}}: \ &  - A^{2}w_1 - A^{2}v - A^{2}m_1 - 2A^{4} - A^{4}q_1 - A^{4}p_1 - A^{6}w_1 - A^{6}v - A^{6}m_1 - A^{8}
\end{align*}
}

The other pair of knotoids that the double-branched cover doesn't distinguish is $\mathbf{3_7}$ and $\mathbf{4_{37}}$. In the case, the Jones polynomial for both lifts in the solid torus is:
\[
- A^{-12} v^3 + A^{-12}v^2 + A^{-8}v^3 - A^{-8}v 
\]
The loop arrow polynomial for each case is:
{\small
\begin{align*}
\mathbf{3_7}: \ &  - A^{6}v - A^{8} - A^{8} p_1 - 2*A^{10} m_1\\
\mathbf{4_{37}}: \ &  A^{4} + A^{4} m_2 + A^{6}*p_2 + 2A^{6} m_1 + A^{8}p_1 + A^{8} m_2
\end{align*}
}

We note here that, despite being a weaker invariant than the double branched cover method, the loop arrow polynomial is significantly faster in terms of computational speed. The double-branched cover method is computationally costly since the pre-image of a knotoid in the solid torus is a significantly larger knot diagram, in terms of number of crossings. This has an impact on the computational speed of the invariant that we choose to evaluate the resulting knot diagram.

\begin{table}[h]
\begin{tabular}{crrrr}
\multicolumn{5}{c}{$S^2$-knotoids}                                                              \\ \hline
\multicolumn{1}{|l|}{Crossings} & \multicolumn{1}{l|}{Primes} & \multicolumn{1}{l|}{Rotatable} & \multicolumn{1}{l|}{Achiral} & \multicolumn{1}{l|}{Strongly Achiral}  \\ \hline
\multicolumn{1}{|c|}{1}         & \multicolumn{1}{r|}{0}      & \multicolumn{1}{r|}{0}&  \multicolumn{1}{r|}{0}     & \multicolumn{1}{r|}{0}     \\ \hline
\multicolumn{1}{|c|}{2}         & \multicolumn{1}{r|}{1}        & \multicolumn{1}{r|}{0} & \multicolumn{1}{r|}{0}     & \multicolumn{1}{r|}{0}   \\ \hline
\multicolumn{1}{|c|}{3}         & \multicolumn{1}{r|}{2}       & \multicolumn{1}{r|}{1} & \multicolumn{1}{r|}{0}     & \multicolumn{1}{r|}{0}     \\ \hline
\multicolumn{1}{|c|}{4}         & \multicolumn{1}{r|}{8}      & \multicolumn{1}{r|}{0} & \multicolumn{1}{r|}{0}    & \multicolumn{1}{r|}{1}      \\ \hline
\multicolumn{1}{|c|}{5}         & \multicolumn{1}{r|}{24}      & \multicolumn{1}{r|}{5} & \multicolumn{1}{r|}{0}  & \multicolumn{1}{r|}{0}     \\ \hline
\multicolumn{1}{|c|}{6}         & \multicolumn{1}{r|}{121}      & \multicolumn{1}{r|}{7} & \multicolumn{1}{r|}{3}  & \multicolumn{1}{r|}{1}     \\ \hline
Total                           & 156 &  13       &        3 & 2                                       
\end{tabular}
\caption{Prime $S^2$-knotoids  with up to 6 crossings.}\label{table:ins2}
\end{table}

\subsection{Knotoids in $S^2$}
Using the algorithm in Section~\ref{sec:algorithm} and  the arrow polynomial, we systematically classified all knotoids in $S^2$ up to 6 crossings. Contrary to the case of planar knotoids, the algorithm terminates after distinguishing all isotopy classes of $S^2$-knotoids. The results are summarized in Table~\ref{table:ins2}. Compared to the table with up to five crossings \cite{barth} we found one missing entry, knotoid $\mathbf{5_{24}}$ in our notation, which confirms the findings of \cite{korablev2017}. Here we rearrange the knotoids with up to five crossings so that they follow the total order of Section~\ref{sec:order} and we extend the table so that it includes also all knotoids with six crossings.
\begin{table}[h]
\begin{tabular}{l|l|l|l}
Crossings & \multicolumn{3}{c}{Rotatable knotoids} \\ \hline
1         & \multicolumn{3}{c}{ $ - $} \\ \hline
2         & \multicolumn{3}{c}{$ - $} \\ \hline
3         & \multicolumn{3}{c}{\begin{tabular}[c]{@{}l@{}}$3_1$\end{tabular}}                                                                                                                                                                                                                                                                                                                                                                                                                                                                                                                                                                                                                                                                                                                                                                                                                                                                                                                                                                                                                                                                                                                                                                                                                                                                                                                                                                                                                                                                                                                                                                                                                                                                                                                                                                                                                                           \\ \hline
4         & \multicolumn{3}{c}{\begin{tabular}[c]{@{}l@{}}$ - $ \end{tabular}}                                                                                                                                                                                                                                                                                                                                                                                                                                                                                                                                                                                                                                                                                                                                                                                                                                                                                                                                                                                                                                                                                                                                                                                                                                                                                                                                                                                                                                                                                                                                                      \\ \hline
 5         & \multicolumn{3}{c}{\begin{tabular}[c]{@{}l@{}} $5_1$, $5_2$, $5_{15}$, $5_{17}$, $5_{18}$ \end{tabular}}
 \\ \hline
  6        & \multicolumn{3}{l}{\begin{tabular}[c]{@{}l@{}} $6_1$, $6_2$, $6_{62}$, $6_{80}$, $6_{85}$, $6_{87}$, $6_{96}$ \end{tabular}}
  \\ \hline
\end{tabular}
\caption{All rotatable knotoids in $S^2$.}\label{table:rotatable_sphere}
\end{table}

In the case of $S^2$-knotoids, there are three achiral knotoids, $\mathbf{6_{54}}$, $\mathbf{6_{86}}$, $\mathbf{6_{120}}$, and two strongly achiral knotoids,  the knot-type knotoid $\mathbf{4_1}$ and $\mathbf{6_3}$.  Finally, there are thirteen rotatable knotoids that are shown in Table~\ref{table:rotatable_sphere} and in Appendix~\ref{appendix:s2}.

\section*{Acknowledgements}
This work has been funded by Leverhulme Trust (Grant RP2013-K-017 to A.S.) and by the Swiss National Science Foundation (Grant 31003A 166684 to A.S.). The authors would like to thank Louis H. Kauffman, Neslihan G\"ug\"umc\"u, Dorothy Buck and Agnese Barbensi for several fruitful conversations.

 \newpage
 
  \appendix
 \section{Planar knotoids}\label{appendix:plane}
We distinguish the different type of knotoids using the following colour-code. Labels in red correspond to knot-type knotoids, labels in blue correspond to $S^2$-isotopy representative when it is seen as a planar knotoid and labels in black to planar knotoids.\\

\vspace{0.4cm}
 \input{input_plane.tex}

 \newpage
 \section{Knotoids in $S^2$}\label{appendix:s2}
 We distinguish the different type of knotoids using the following colour-code. Labels in red correspond to knot-type knotoids and labels in black to $S^2$-knotoids.\\

\vspace{0.4cm}

\input{input_sphere.tex}

 \end{document}

%% file: input_plane.tex
\begin{minipage}[t]{.25\linewidth}
\centering
\includegraphics[width=0.9\textwidth,height=3.5cm,keepaspectratio]{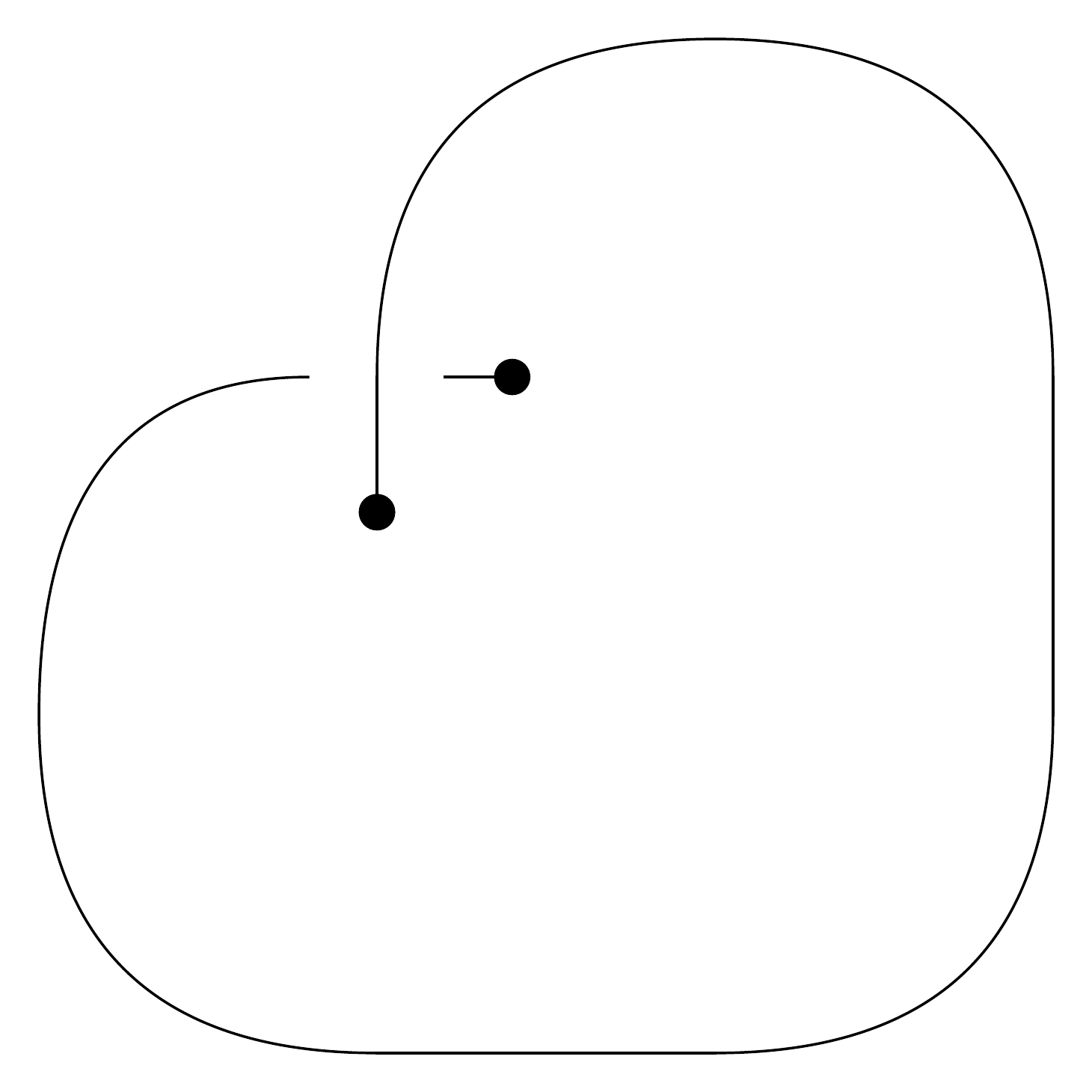}\\
\textcolor{black}{$1_{1}$}
\vspace{1cm}
\end{minipage}
\begin{minipage}[t]{.25\linewidth}
\centering
\includegraphics[width=0.9\textwidth,height=3.5cm,keepaspectratio]{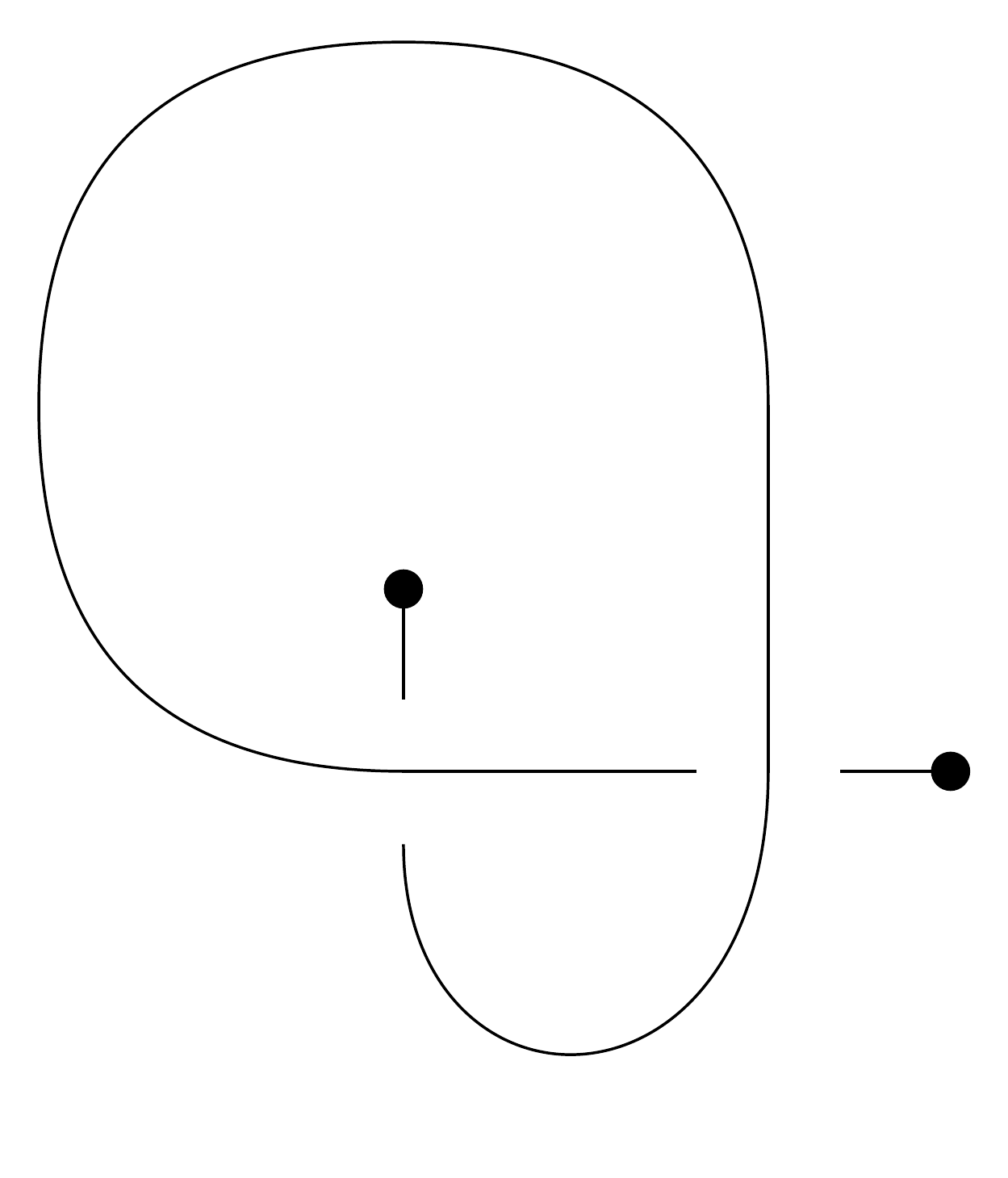}\\
\textcolor{blue}{$2_{1}$}
\vspace{1cm}
\end{minipage}
\begin{minipage}[t]{.25\linewidth}
\centering
\includegraphics[width=0.9\textwidth,height=3.5cm,keepaspectratio]{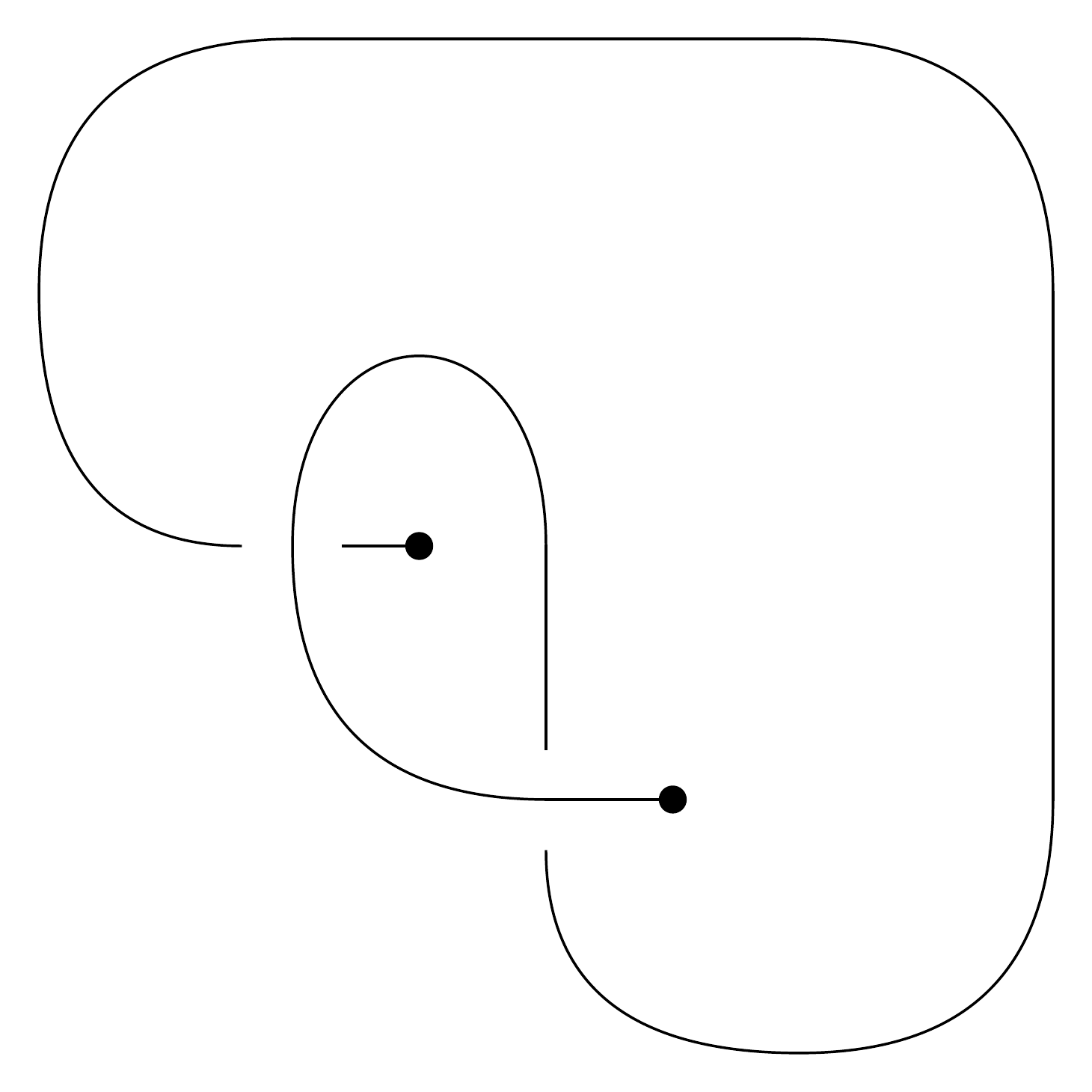}\\
\textcolor{black}{$2_{2}$}
\vspace{1cm}
\end{minipage}
\begin{minipage}[t]{.25\linewidth}
\centering
\includegraphics[width=0.9\textwidth,height=3.5cm,keepaspectratio]{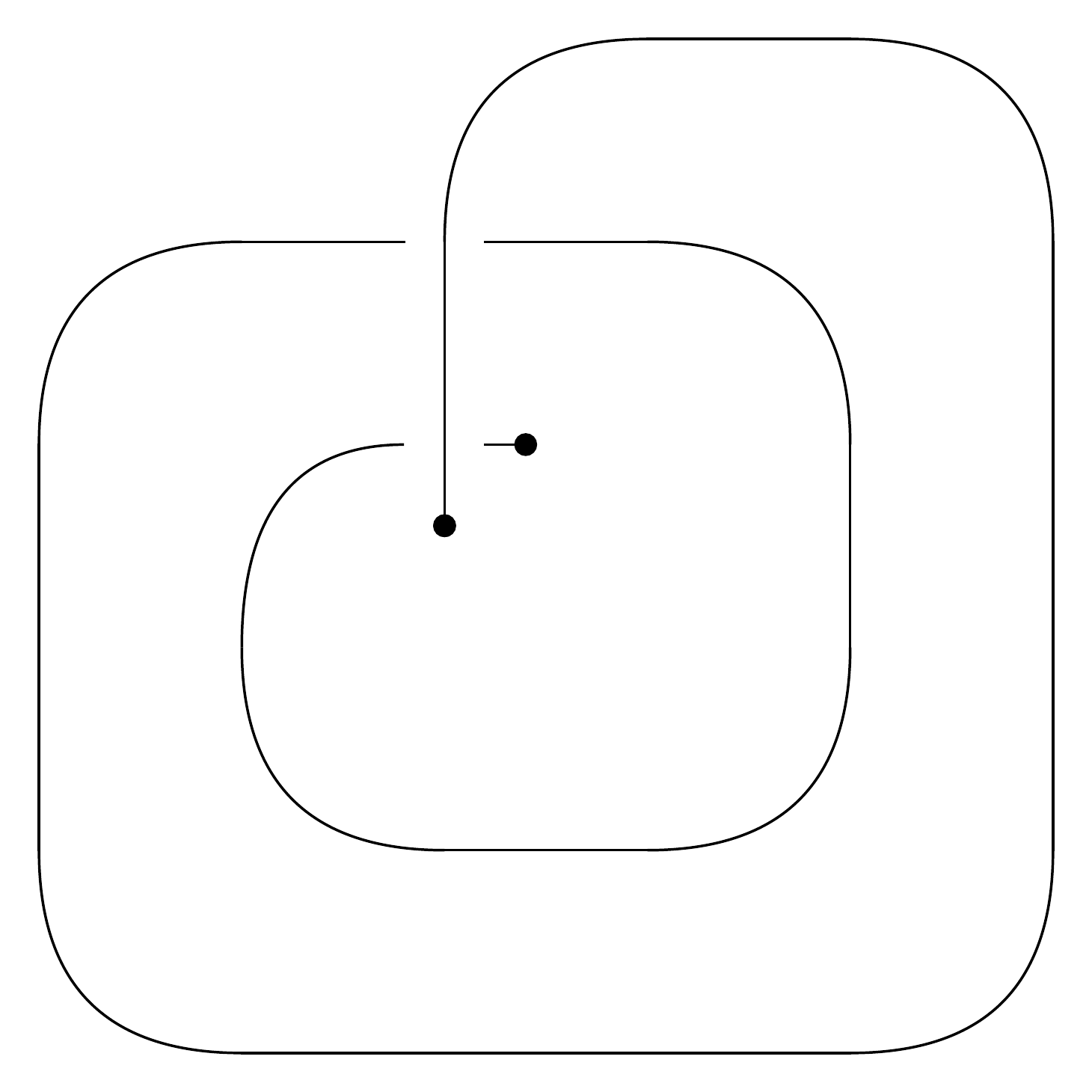}\\
\textcolor{black}{$2_{3}$}
\vspace{1cm}
\end{minipage}
\begin{minipage}[t]{.25\linewidth}
\centering
\includegraphics[width=0.9\textwidth,height=3.5cm,keepaspectratio]{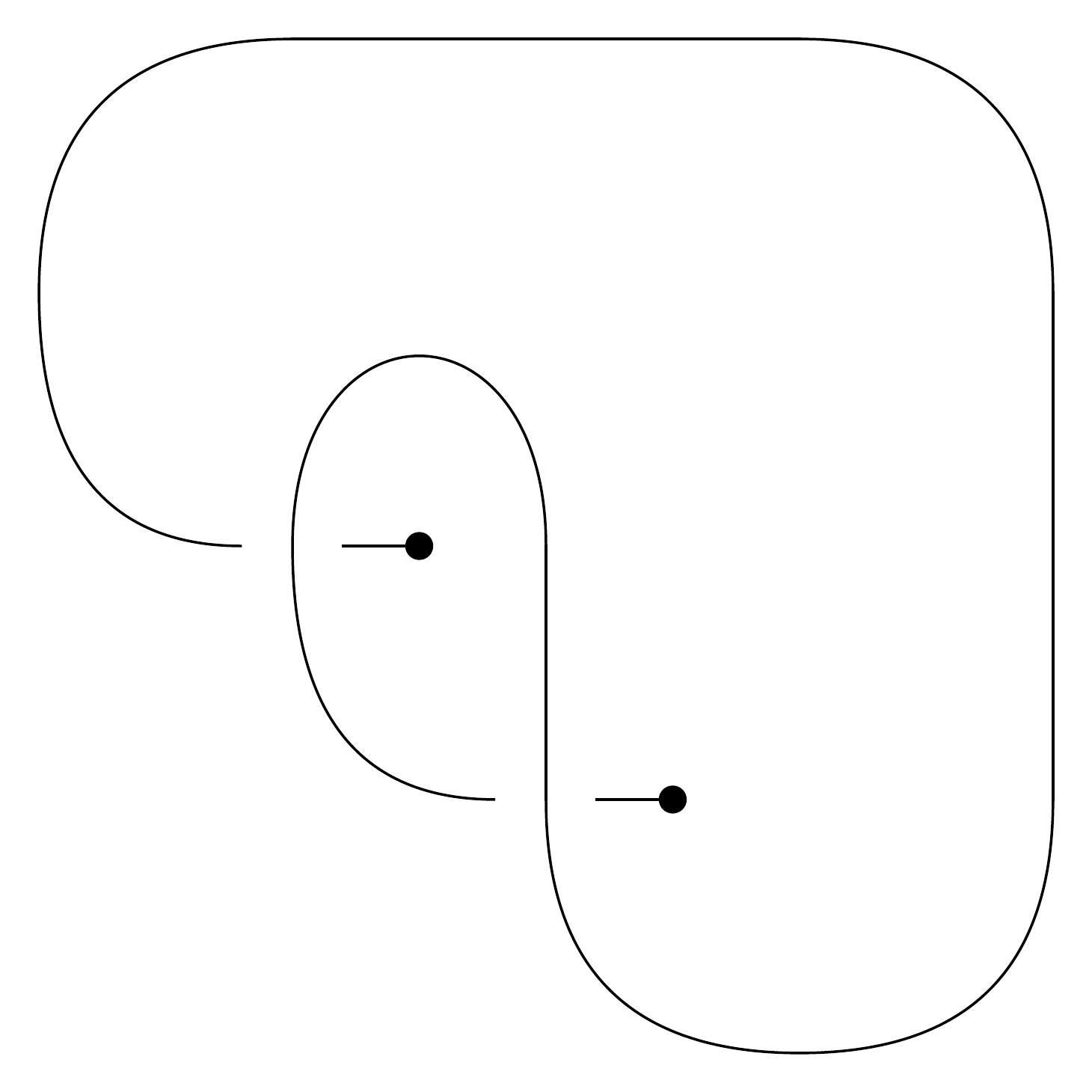}\\
\textcolor{black}{$2_{4}$}
\vspace{1cm}
\end{minipage}
\begin{minipage}[t]{.25\linewidth}
\centering
\includegraphics[width=0.9\textwidth,height=3.5cm,keepaspectratio]{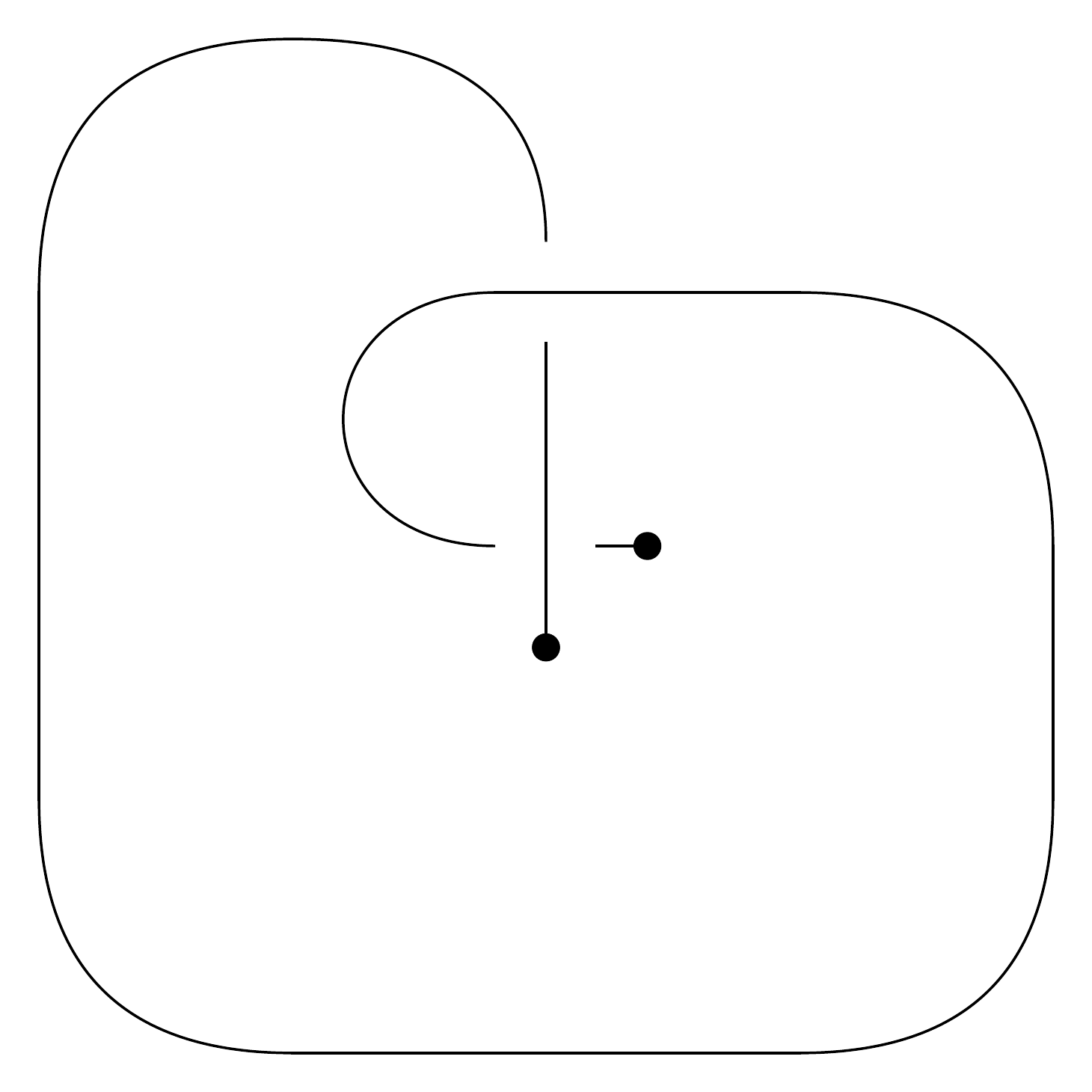}\\
\textcolor{black}{$2_{5}$}
\vspace{1cm}
\end{minipage}
\begin{minipage}[t]{.25\linewidth}
\centering
\includegraphics[width=0.9\textwidth,height=3.5cm,keepaspectratio]{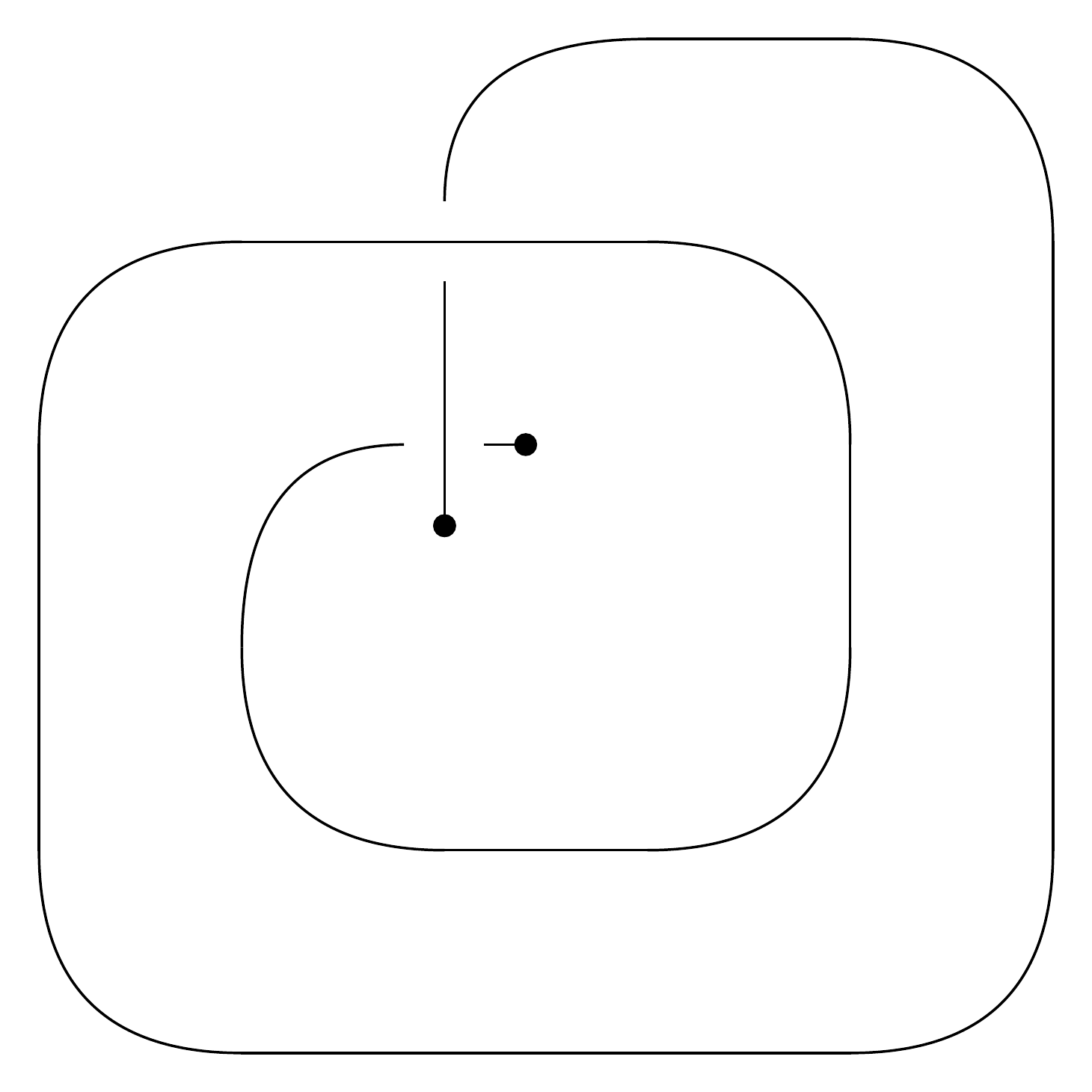}\\
\textcolor{black}{$2_{6}$}
\vspace{1cm}
\end{minipage}
\begin{minipage}[t]{.25\linewidth}
\centering
\includegraphics[width=0.9\textwidth,height=3.5cm,keepaspectratio]{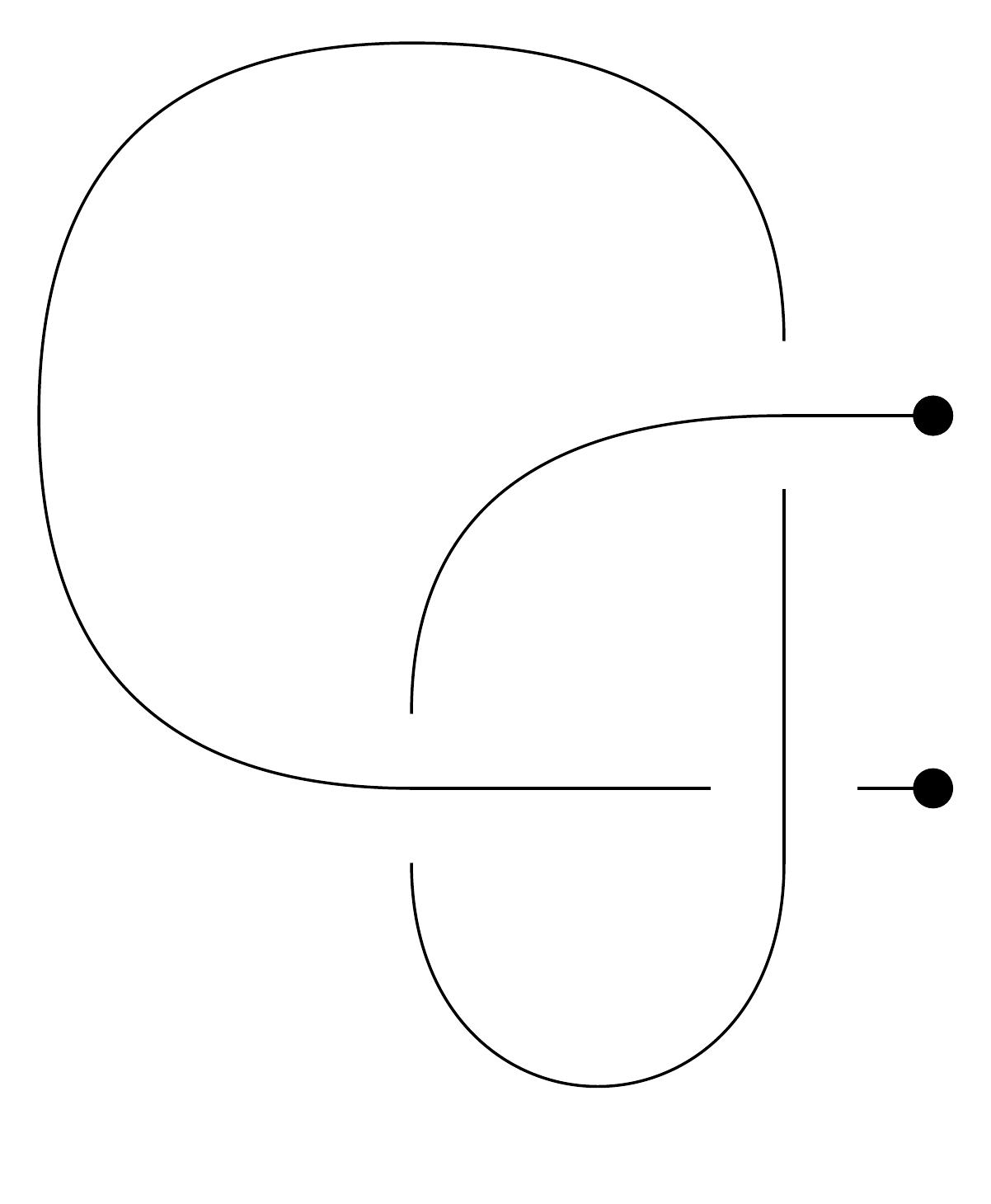}\\
\textcolor{red}{$3_{1}$}
\vspace{1cm}
\end{minipage}
\begin{minipage}[t]{.25\linewidth}
\centering
\includegraphics[width=0.9\textwidth,height=3.5cm,keepaspectratio]{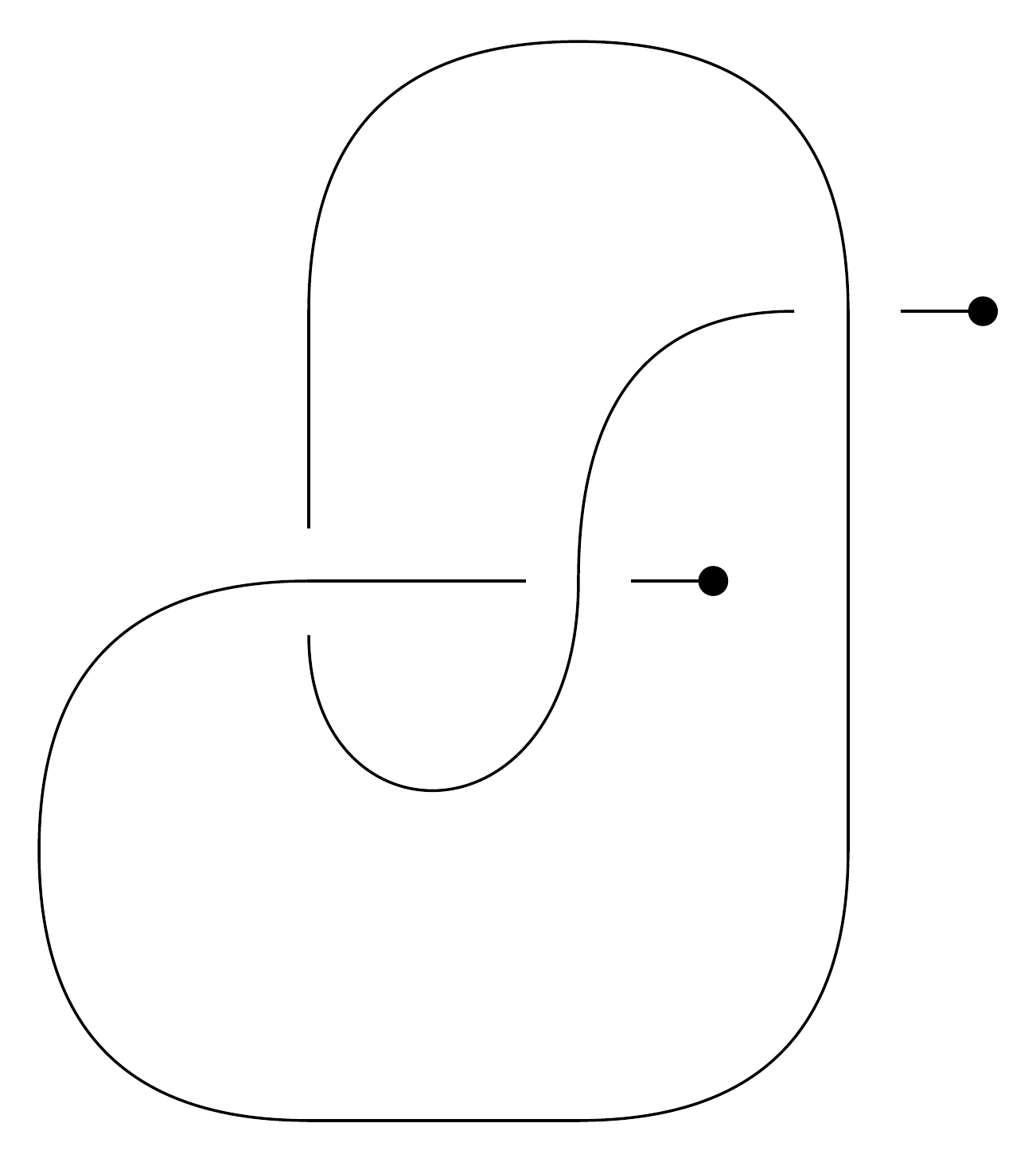}\\
\textcolor{blue}{$3_{2}$}
\vspace{1cm}
\end{minipage}
\begin{minipage}[t]{.25\linewidth}
\centering
\includegraphics[width=0.9\textwidth,height=3.5cm,keepaspectratio]{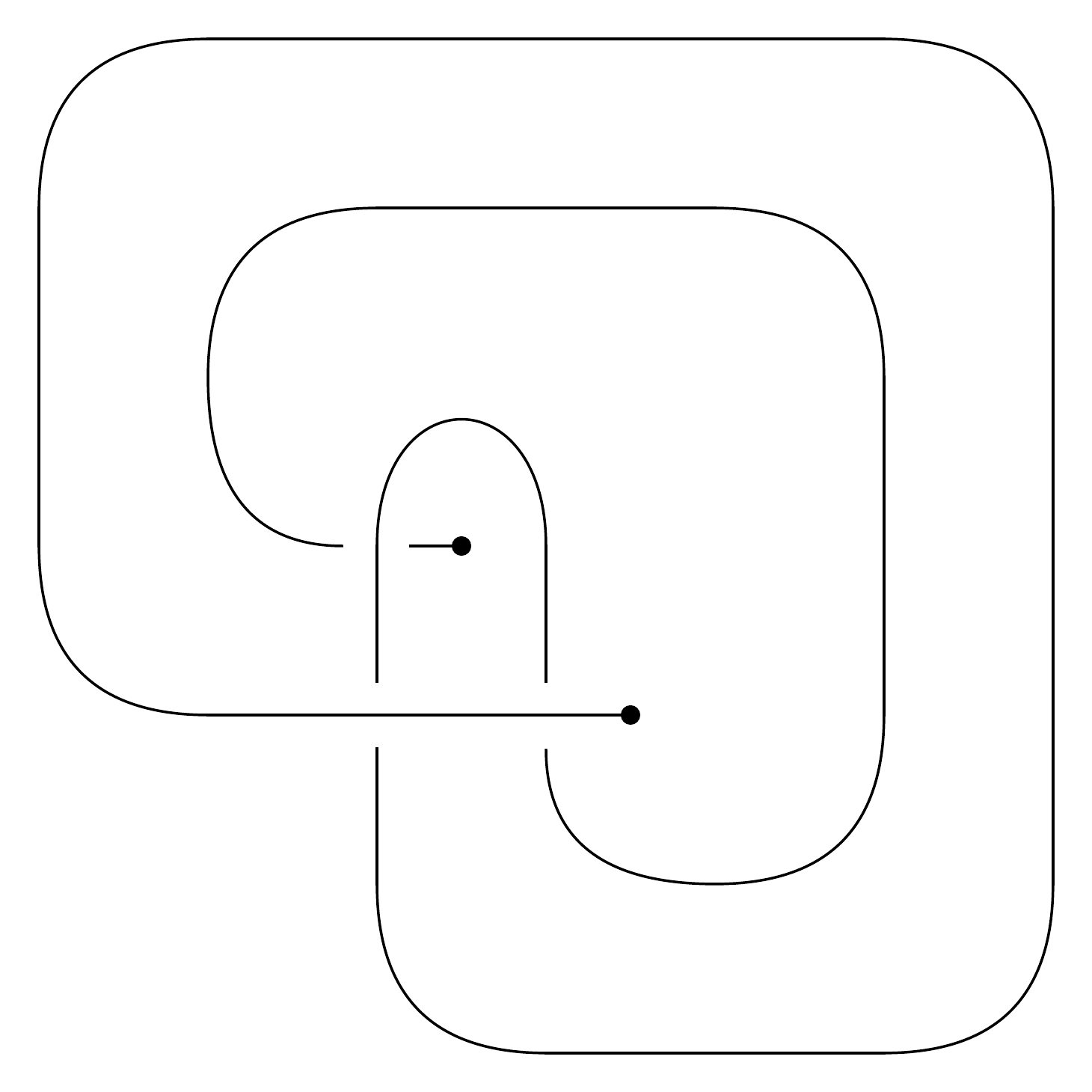}\\
\textcolor{black}{$3_{3}$}
\vspace{1cm}
\end{minipage}
\begin{minipage}[t]{.25\linewidth}
\centering
\includegraphics[width=0.9\textwidth,height=3.5cm,keepaspectratio]{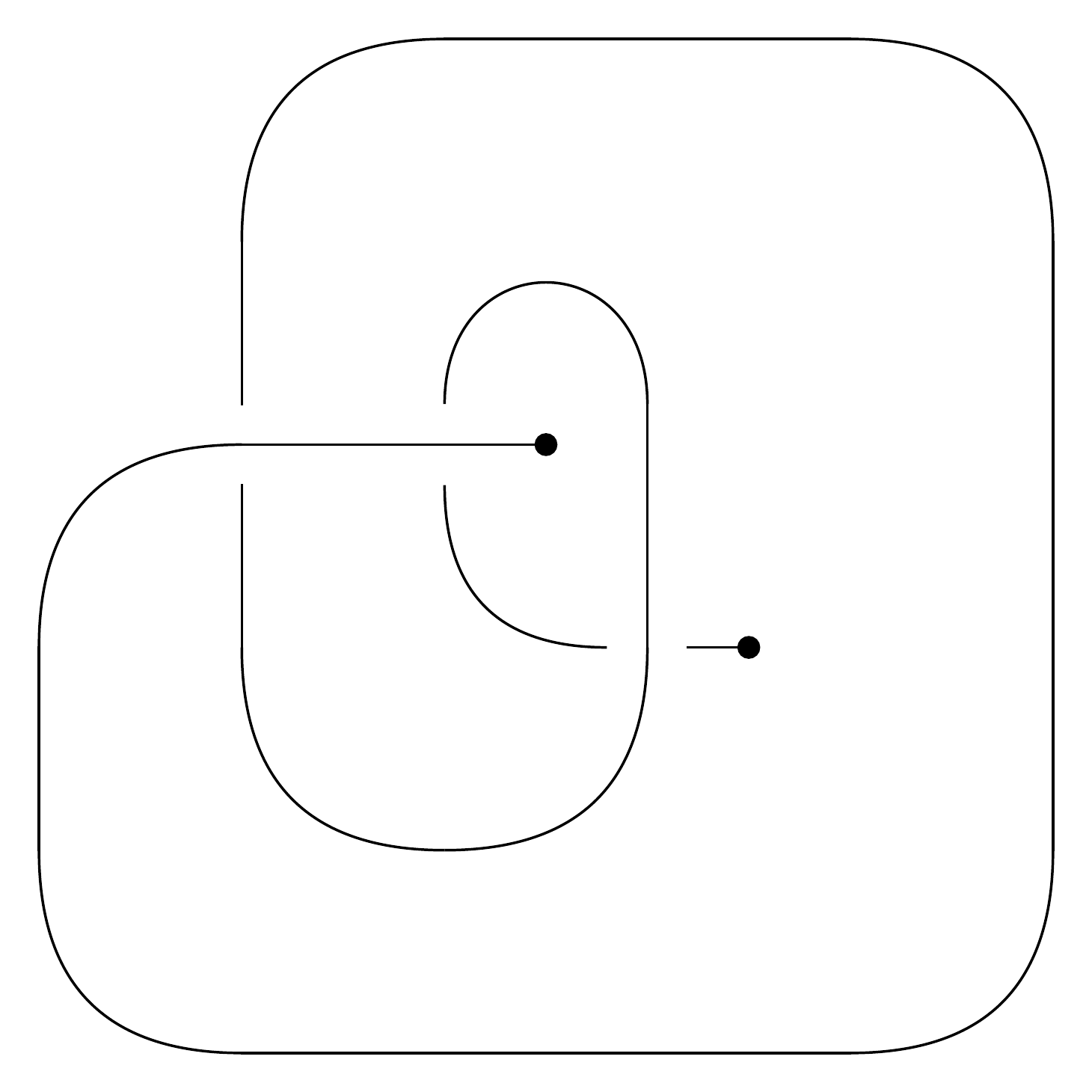}\\
\textcolor{black}{$3_{4}$}
\vspace{1cm}
\end{minipage}
\begin{minipage}[t]{.25\linewidth}
\centering
\includegraphics[width=0.9\textwidth,height=3.5cm,keepaspectratio]{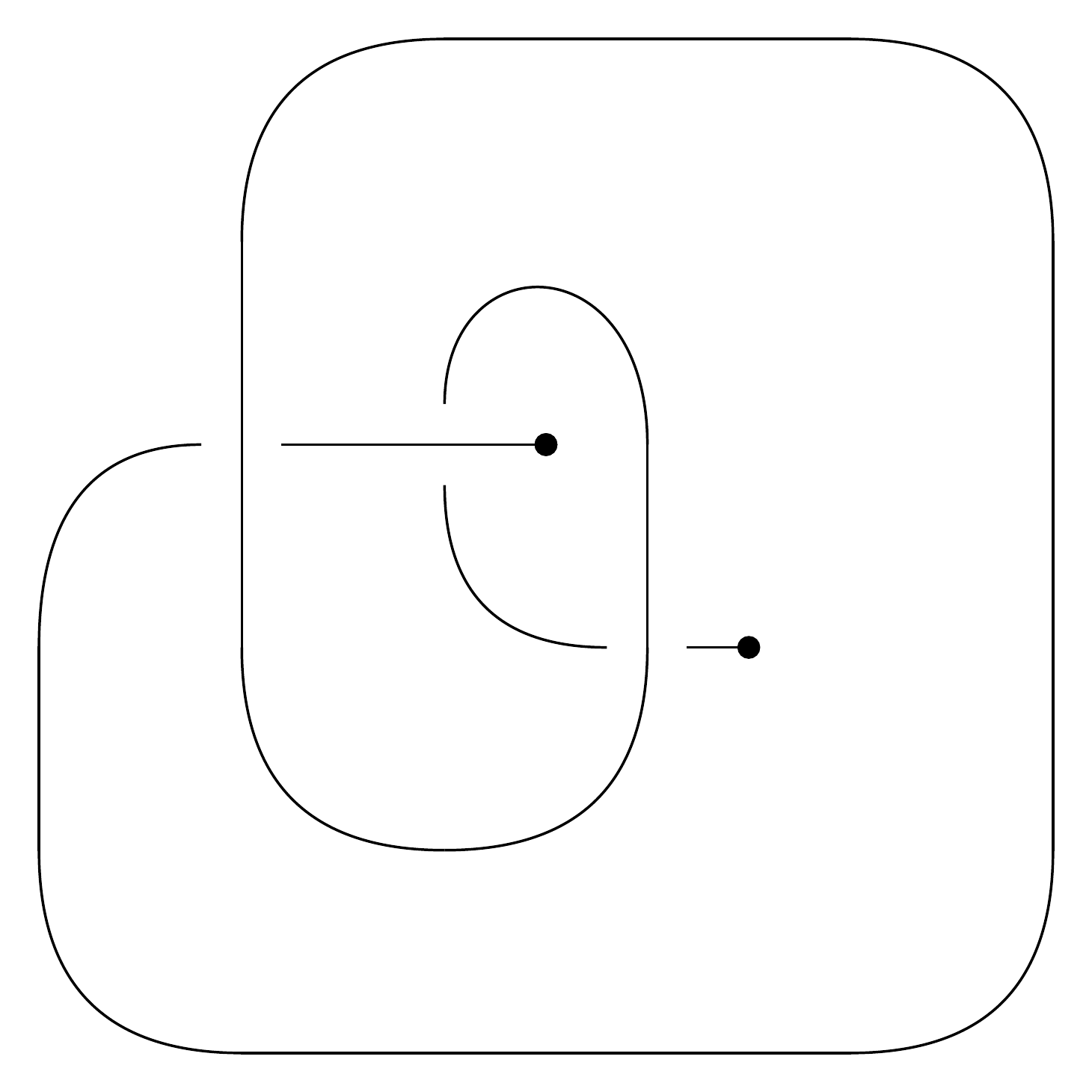}\\
\textcolor{black}{$3_{5}$}
\vspace{1cm}
\end{minipage}
\begin{minipage}[t]{.25\linewidth}
\centering
\includegraphics[width=0.9\textwidth,height=3.5cm,keepaspectratio]{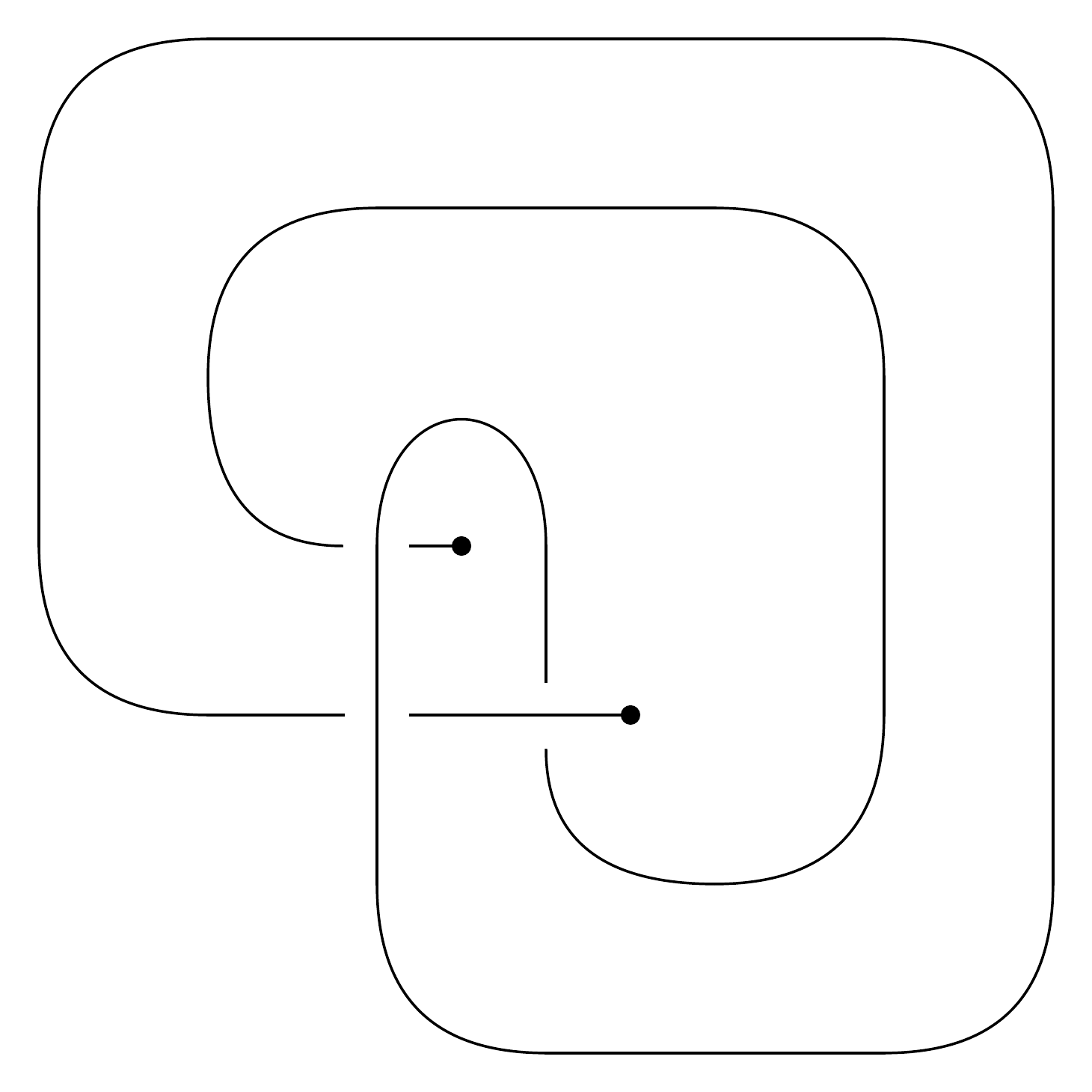}\\
\textcolor{black}{$3_{6}$}
\vspace{1cm}
\end{minipage}
\begin{minipage}[t]{.25\linewidth}
\centering
\includegraphics[width=0.9\textwidth,height=3.5cm,keepaspectratio]{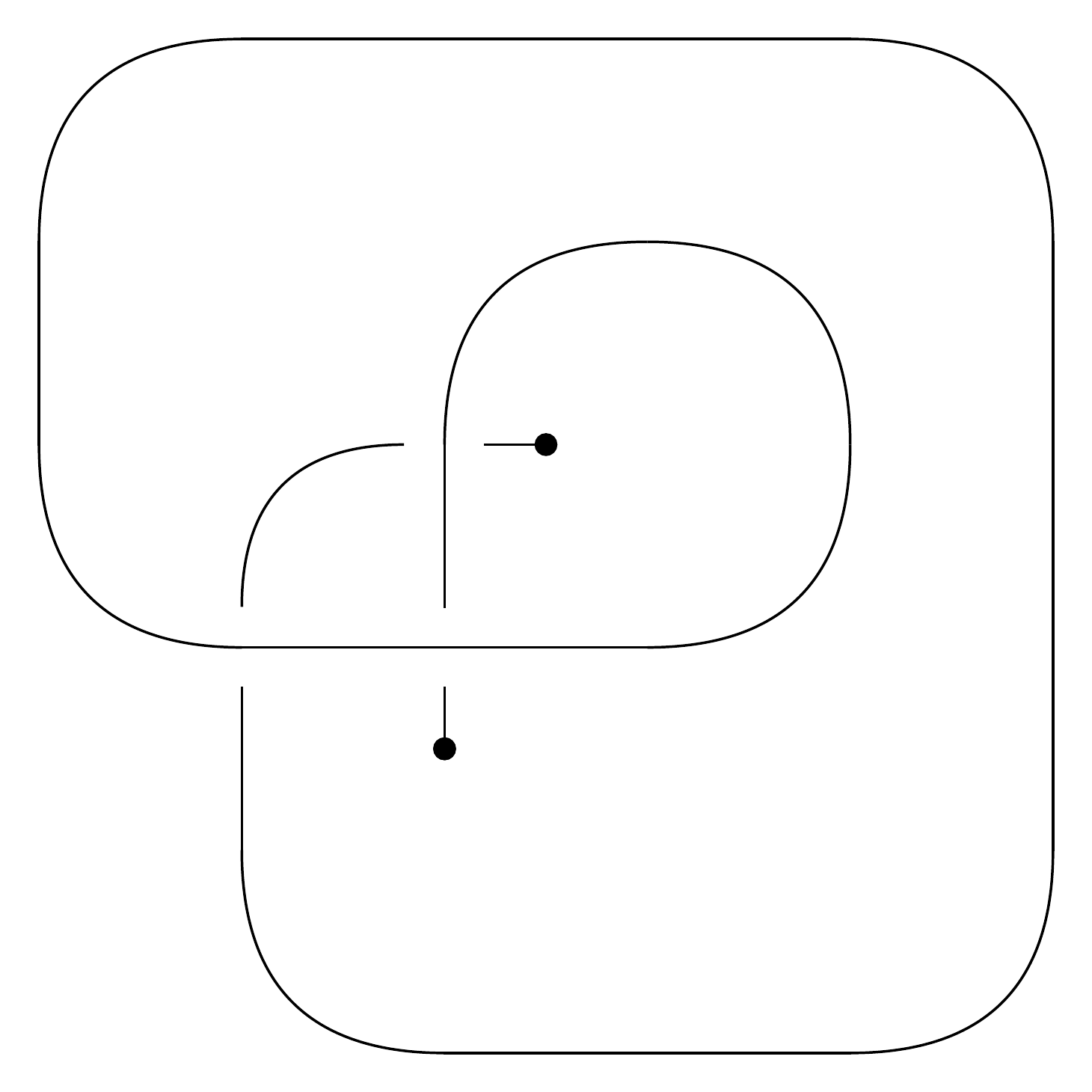}\\
\textcolor{black}{$3_{7}$}
\vspace{1cm}
\end{minipage}
\begin{minipage}[t]{.25\linewidth}
\centering
\includegraphics[width=0.9\textwidth,height=3.5cm,keepaspectratio]{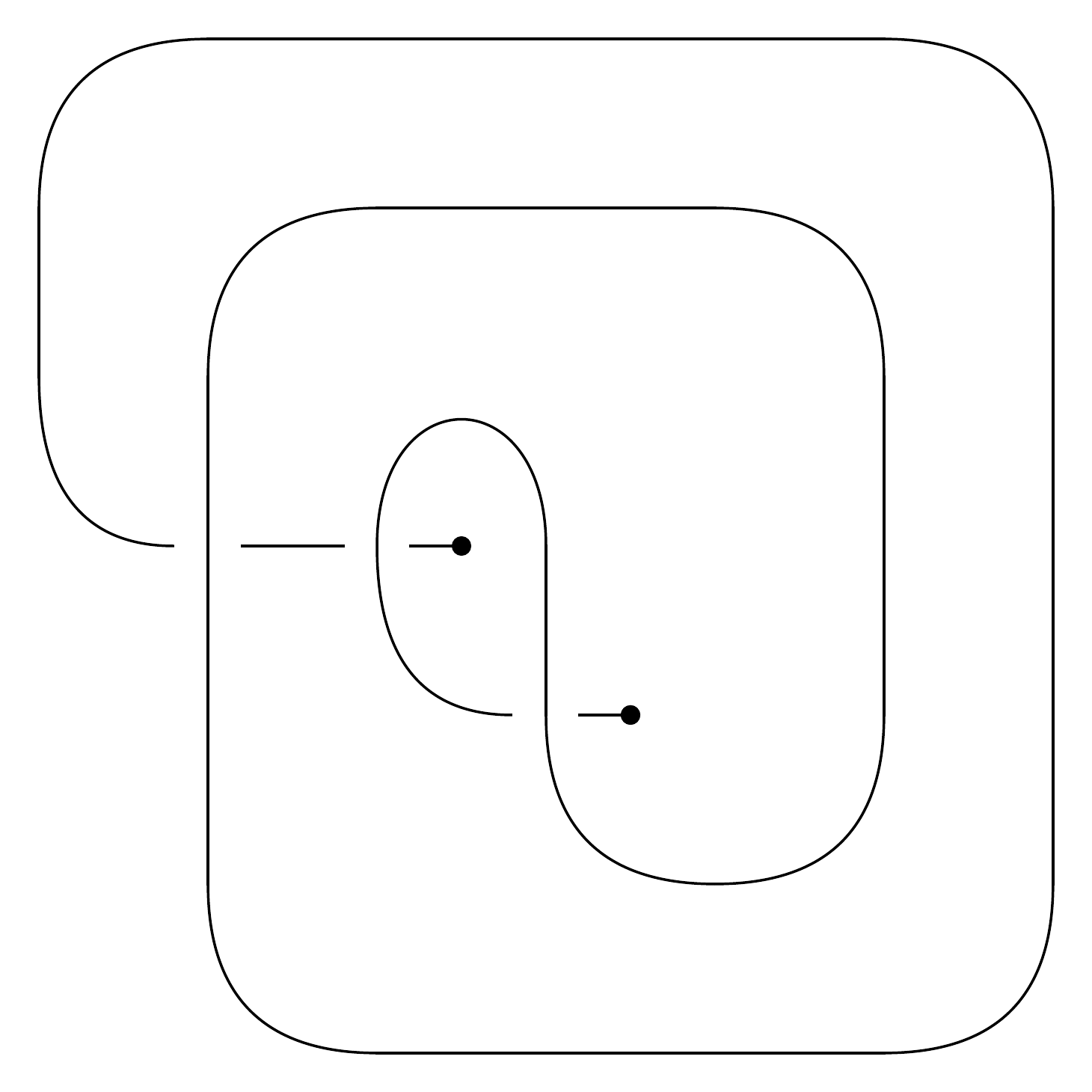}\\
\textcolor{black}{$3_{8}$}
\vspace{1cm}
\end{minipage}
\begin{minipage}[t]{.25\linewidth}
\centering
\includegraphics[width=0.9\textwidth,height=3.5cm,keepaspectratio]{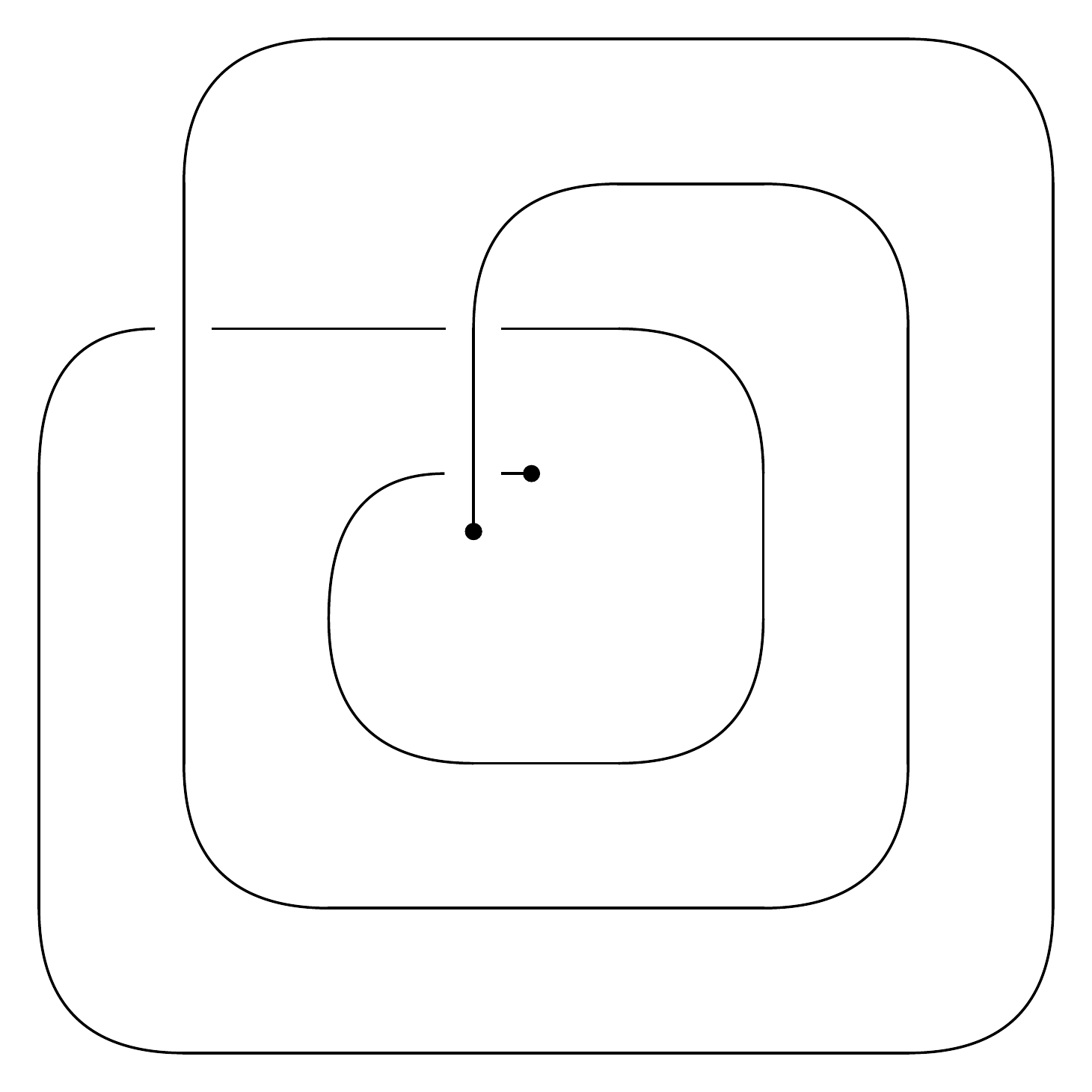}\\
\textcolor{black}{$3_{9}$}
\vspace{1cm}
\end{minipage}
\begin{minipage}[t]{.25\linewidth}
\centering
\includegraphics[width=0.9\textwidth,height=3.5cm,keepaspectratio]{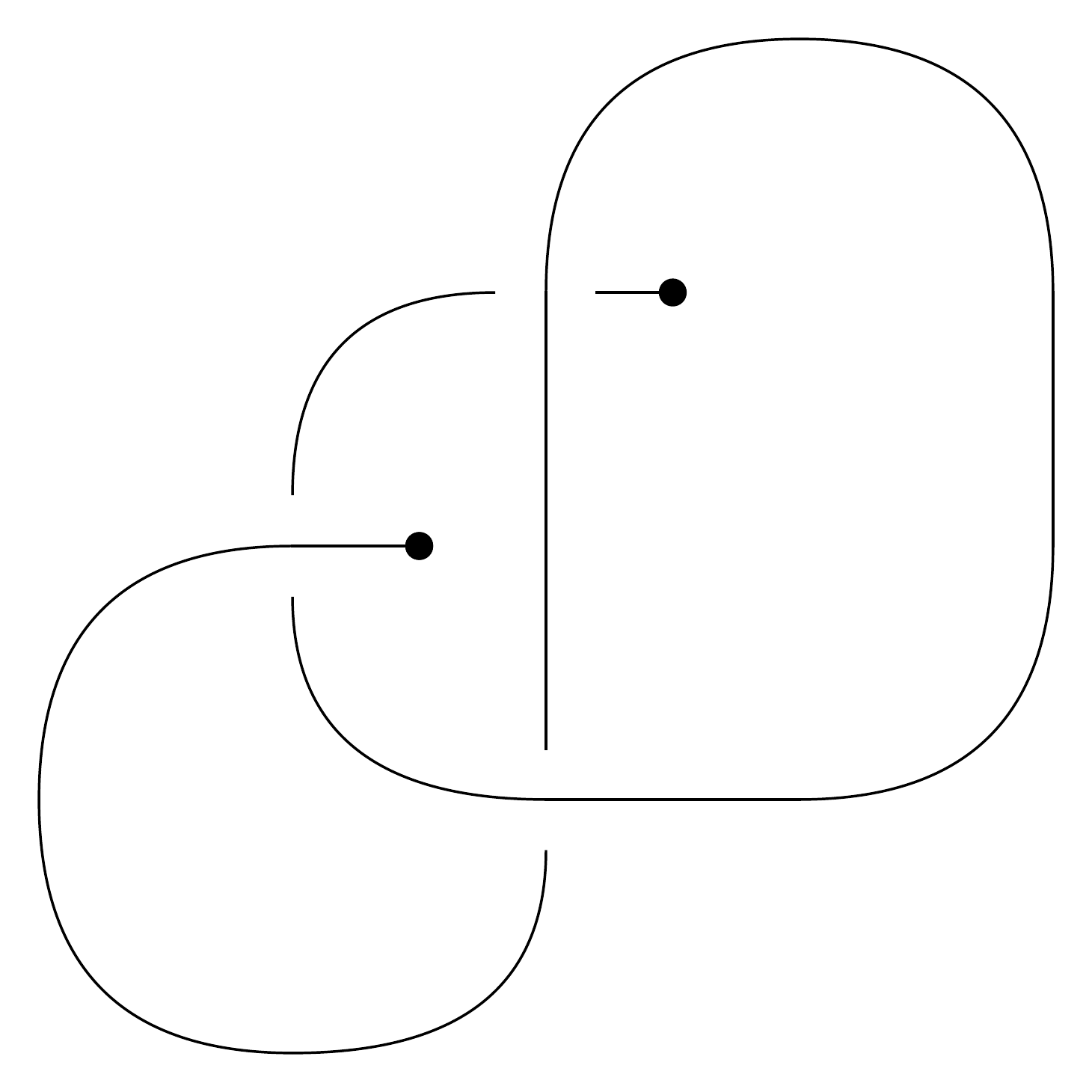}\\
\textcolor{black}{$3_{10}$}
\vspace{1cm}
\end{minipage}
\begin{minipage}[t]{.25\linewidth}
\centering
\includegraphics[width=0.9\textwidth,height=3.5cm,keepaspectratio]{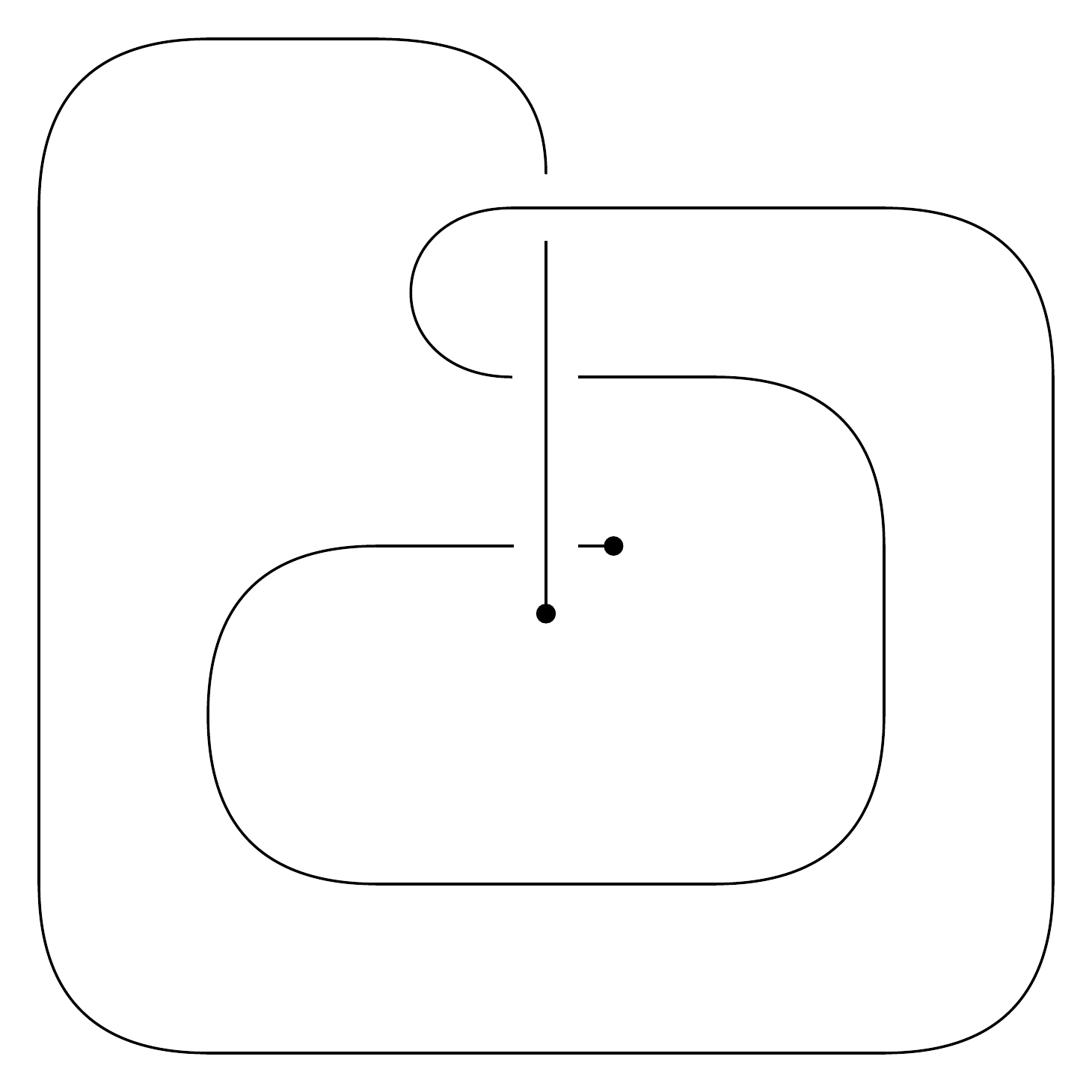}\\
\textcolor{black}{$3_{11}$}
\vspace{1cm}
\end{minipage}
\begin{minipage}[t]{.25\linewidth}
\centering
\includegraphics[width=0.9\textwidth,height=3.5cm,keepaspectratio]{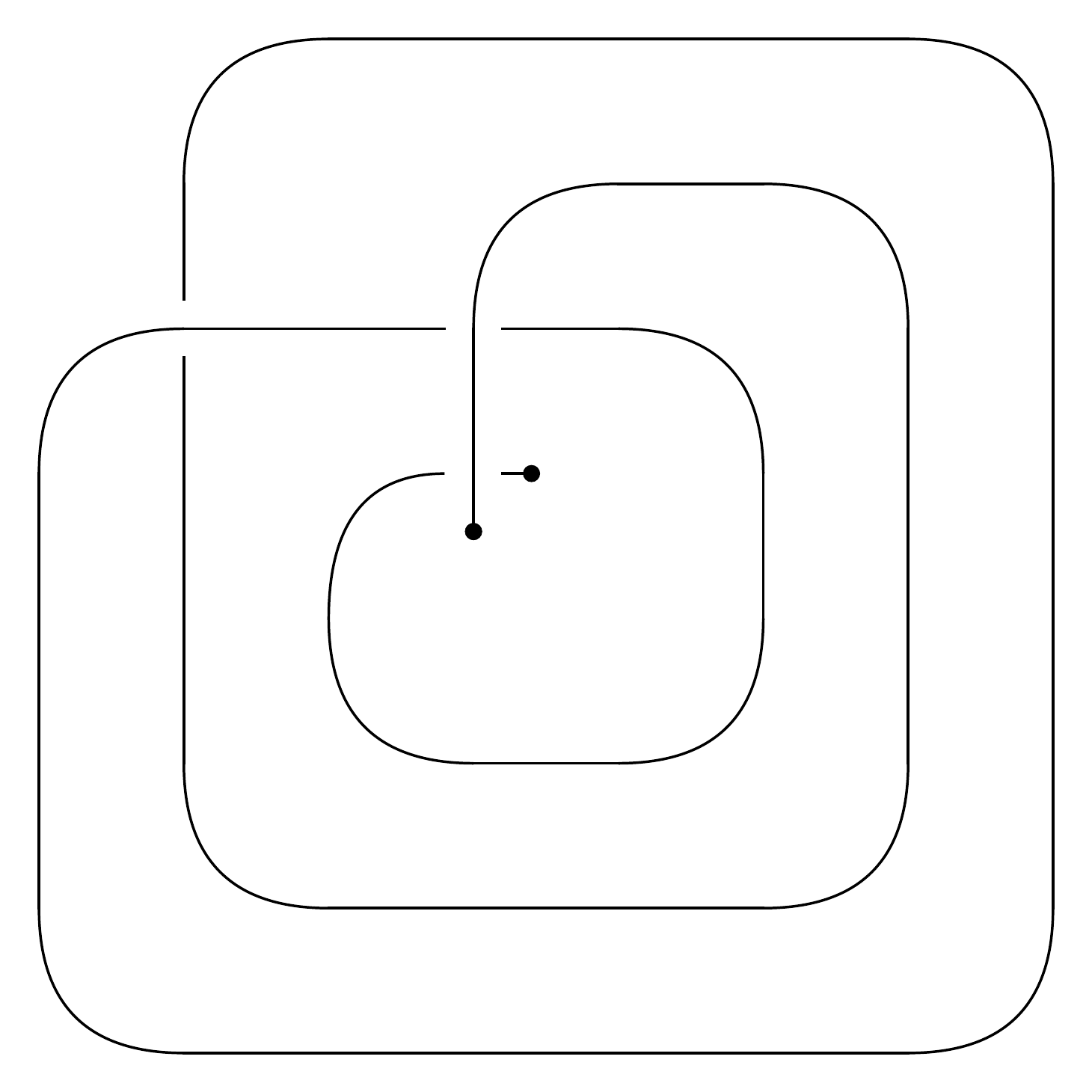}\\
\textcolor{black}{$3_{12}$}
\vspace{1cm}
\end{minipage}
\begin{minipage}[t]{.25\linewidth}
\centering
\includegraphics[width=0.9\textwidth,height=3.5cm,keepaspectratio]{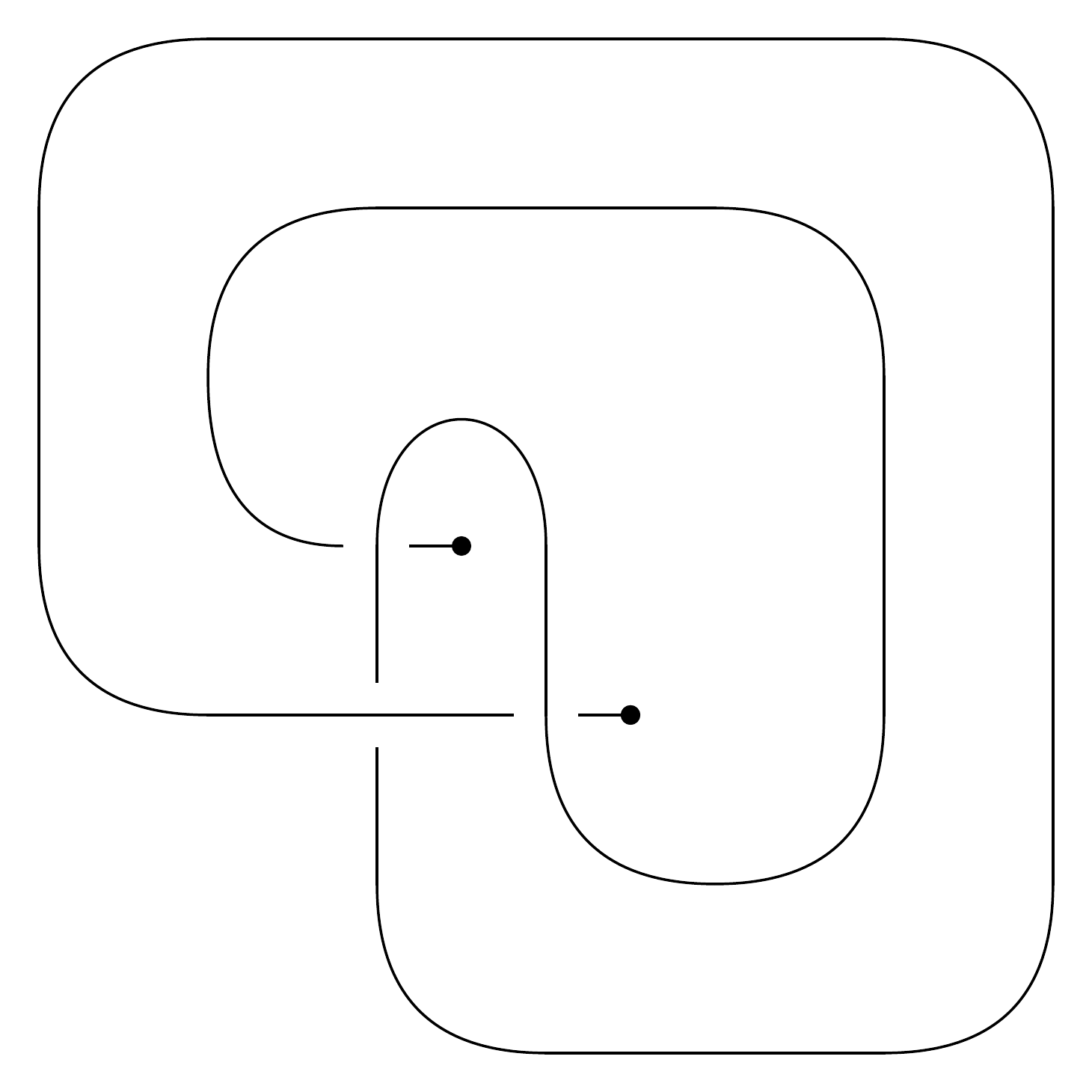}\\
\textcolor{black}{$3_{13}$}
\vspace{1cm}
\end{minipage}
\begin{minipage}[t]{.25\linewidth}
\centering
\includegraphics[width=0.9\textwidth,height=3.5cm,keepaspectratio]{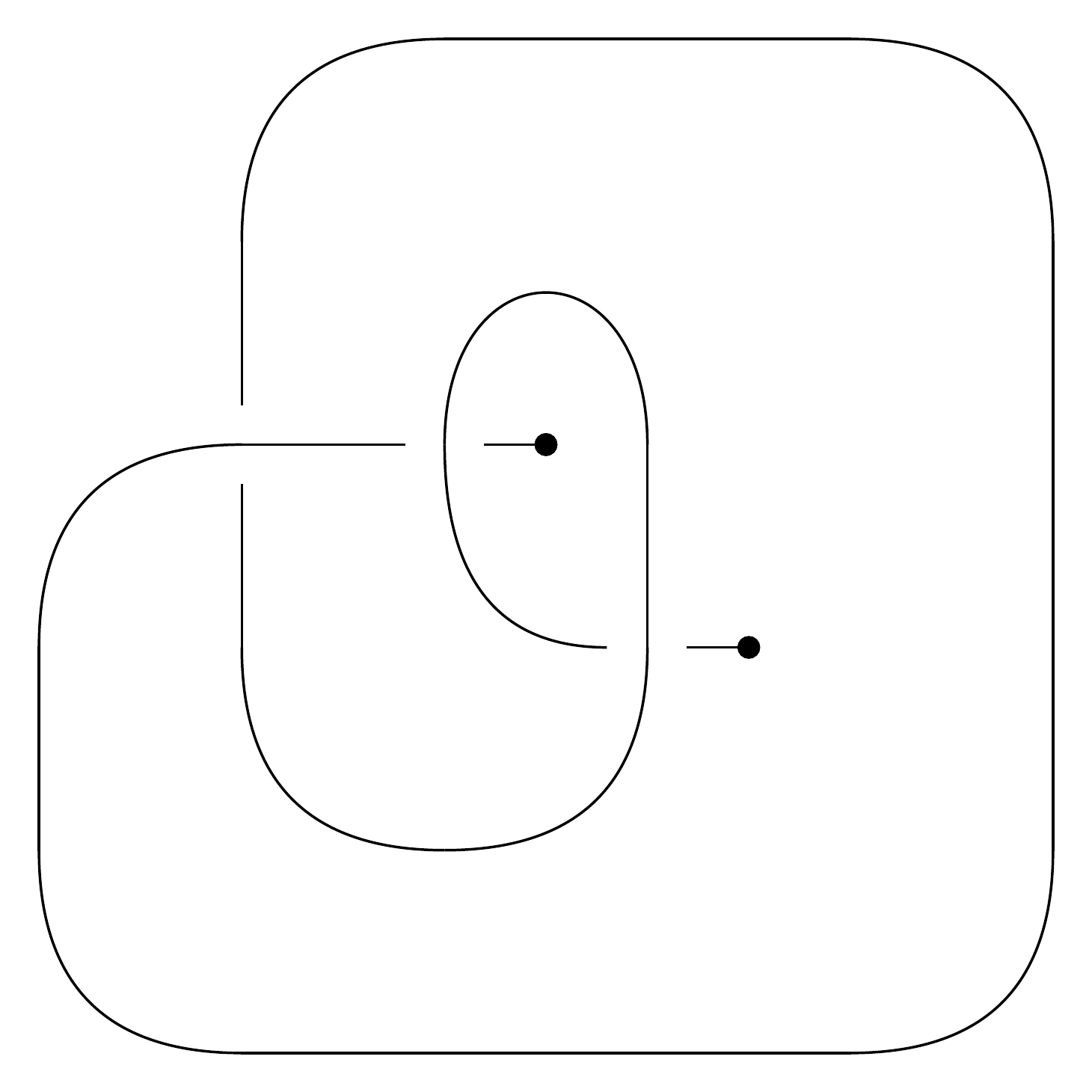}\\
\textcolor{black}{$3_{14}$}
\vspace{1cm}
\end{minipage}
\begin{minipage}[t]{.25\linewidth}
\centering
\includegraphics[width=0.9\textwidth,height=3.5cm,keepaspectratio]{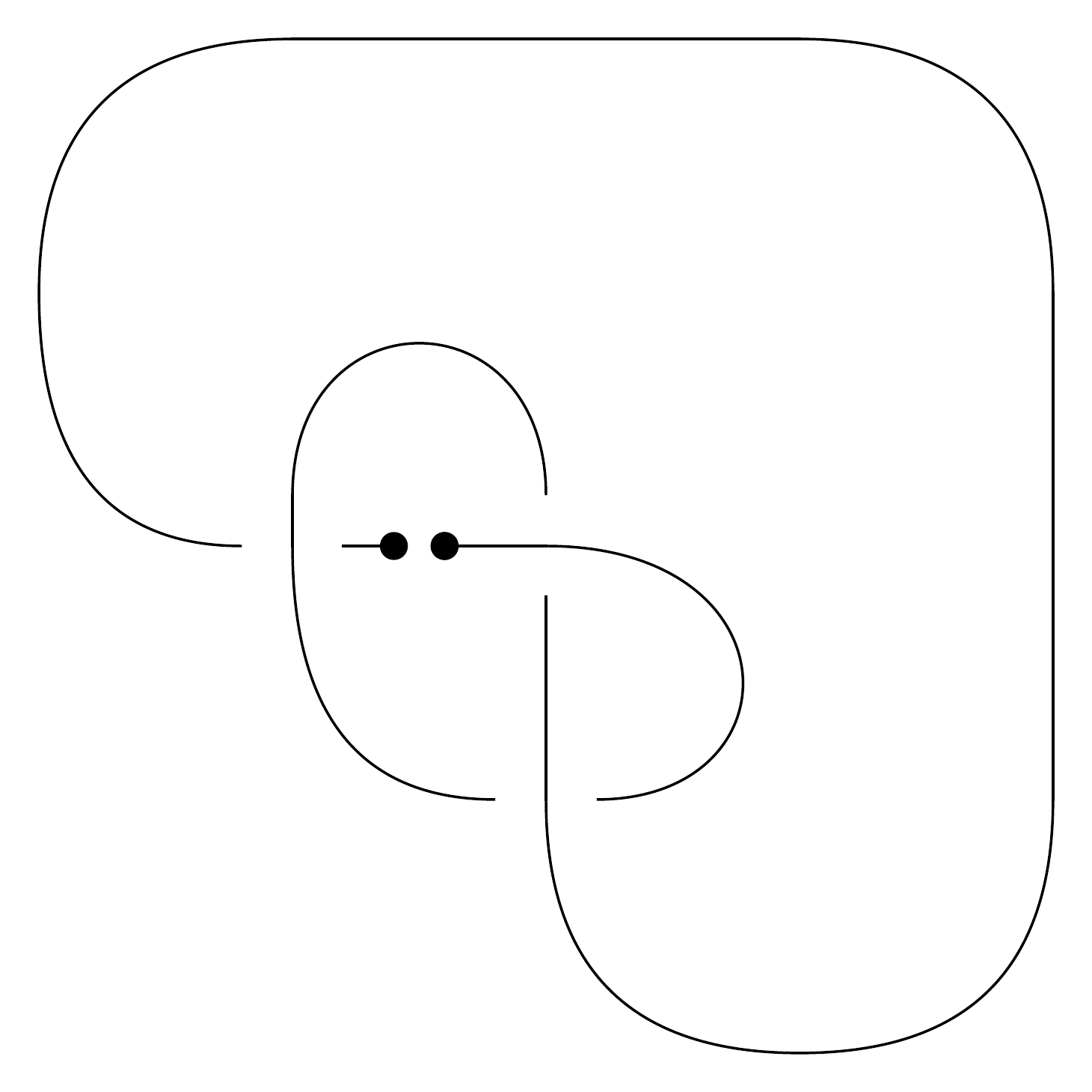}\\
\textcolor{black}{$3_{15}$}
\vspace{1cm}
\end{minipage}
\begin{minipage}[t]{.25\linewidth}
\centering
\includegraphics[width=0.9\textwidth,height=3.5cm,keepaspectratio]{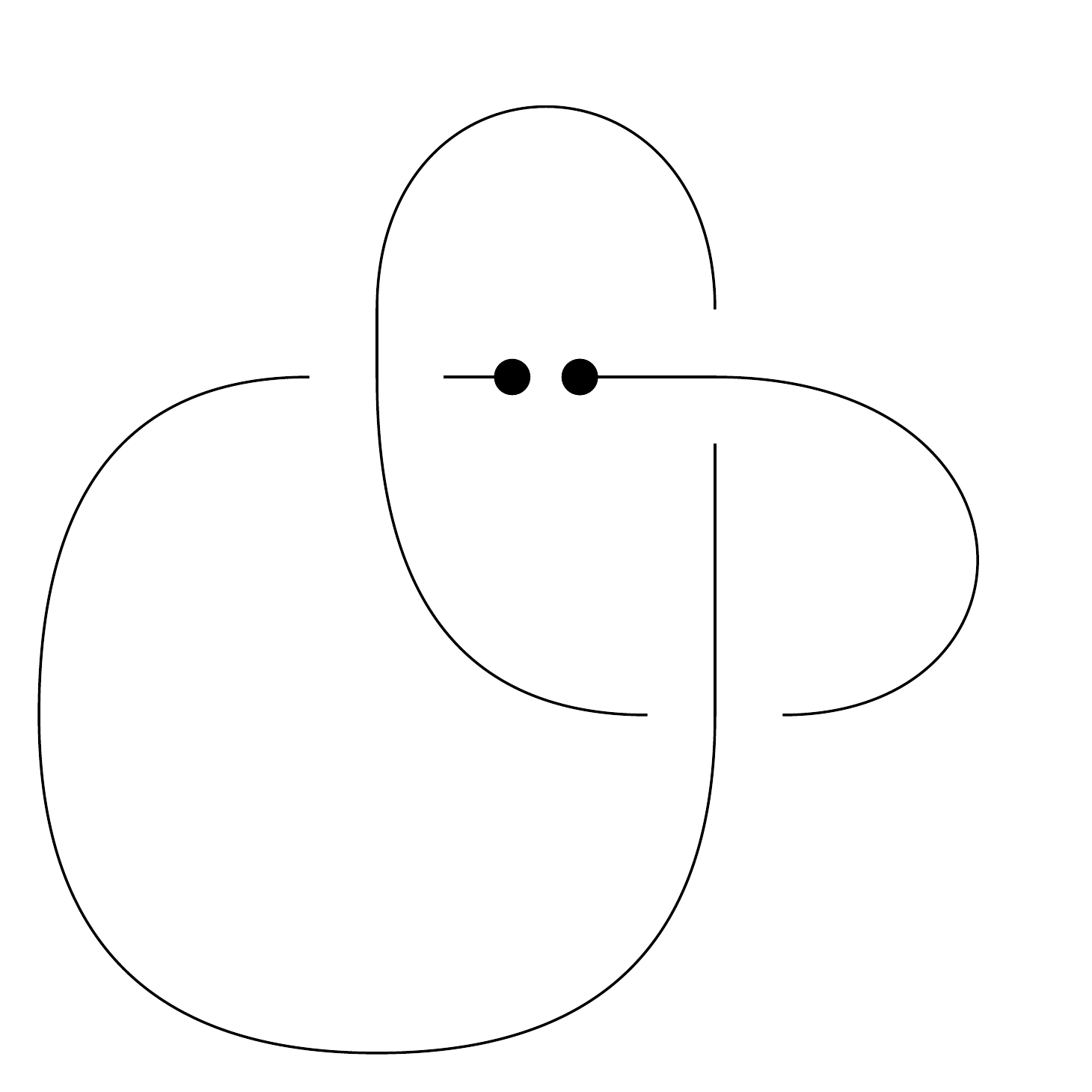}\\
\textcolor{black}{$3_{16}$}
\vspace{1cm}
\end{minipage}
\begin{minipage}[t]{.25\linewidth}
\centering
\includegraphics[width=0.9\textwidth,height=3.5cm,keepaspectratio]{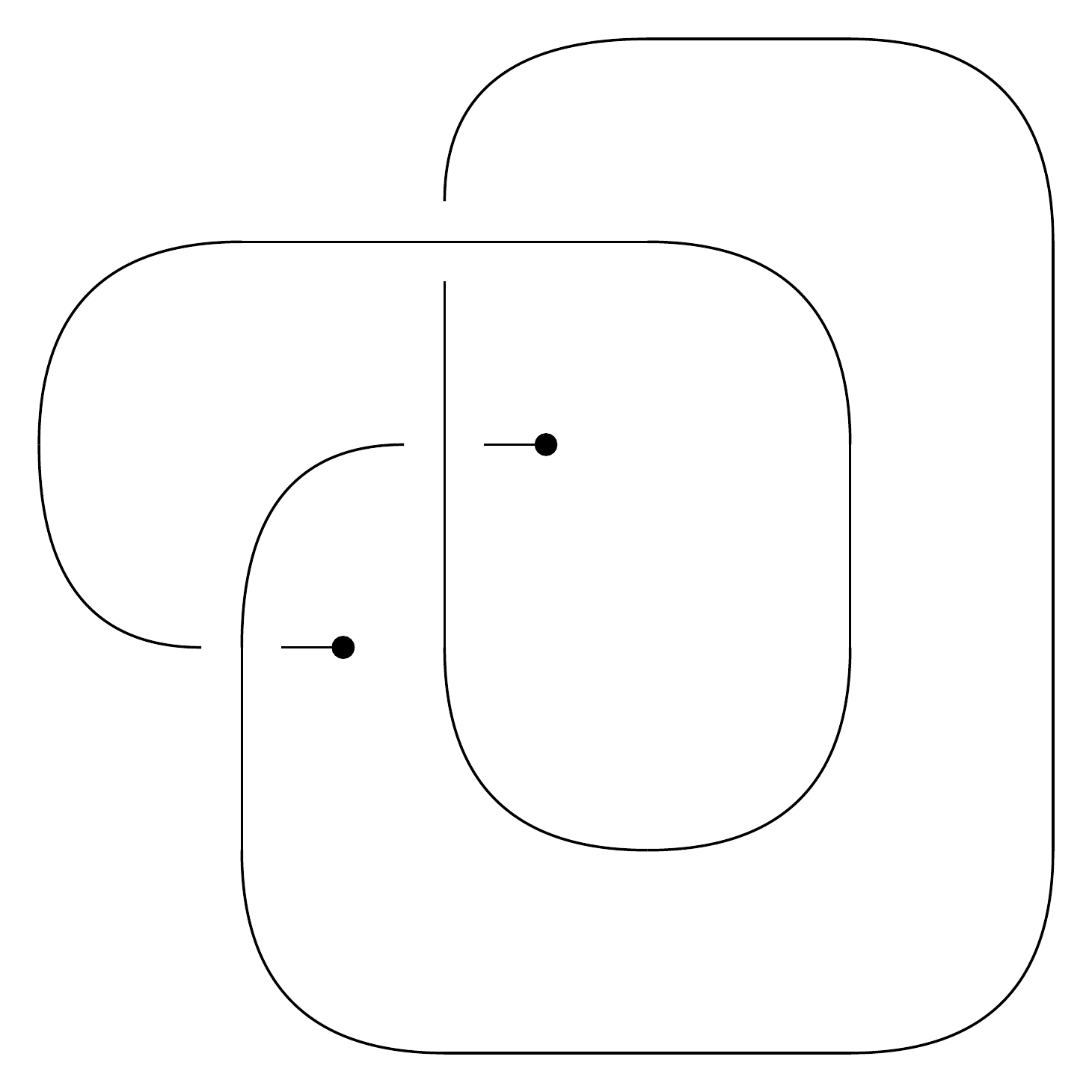}\\
\textcolor{black}{$3_{17}$}
\vspace{1cm}
\end{minipage}
\begin{minipage}[t]{.25\linewidth}
\centering
\includegraphics[width=0.9\textwidth,height=3.5cm,keepaspectratio]{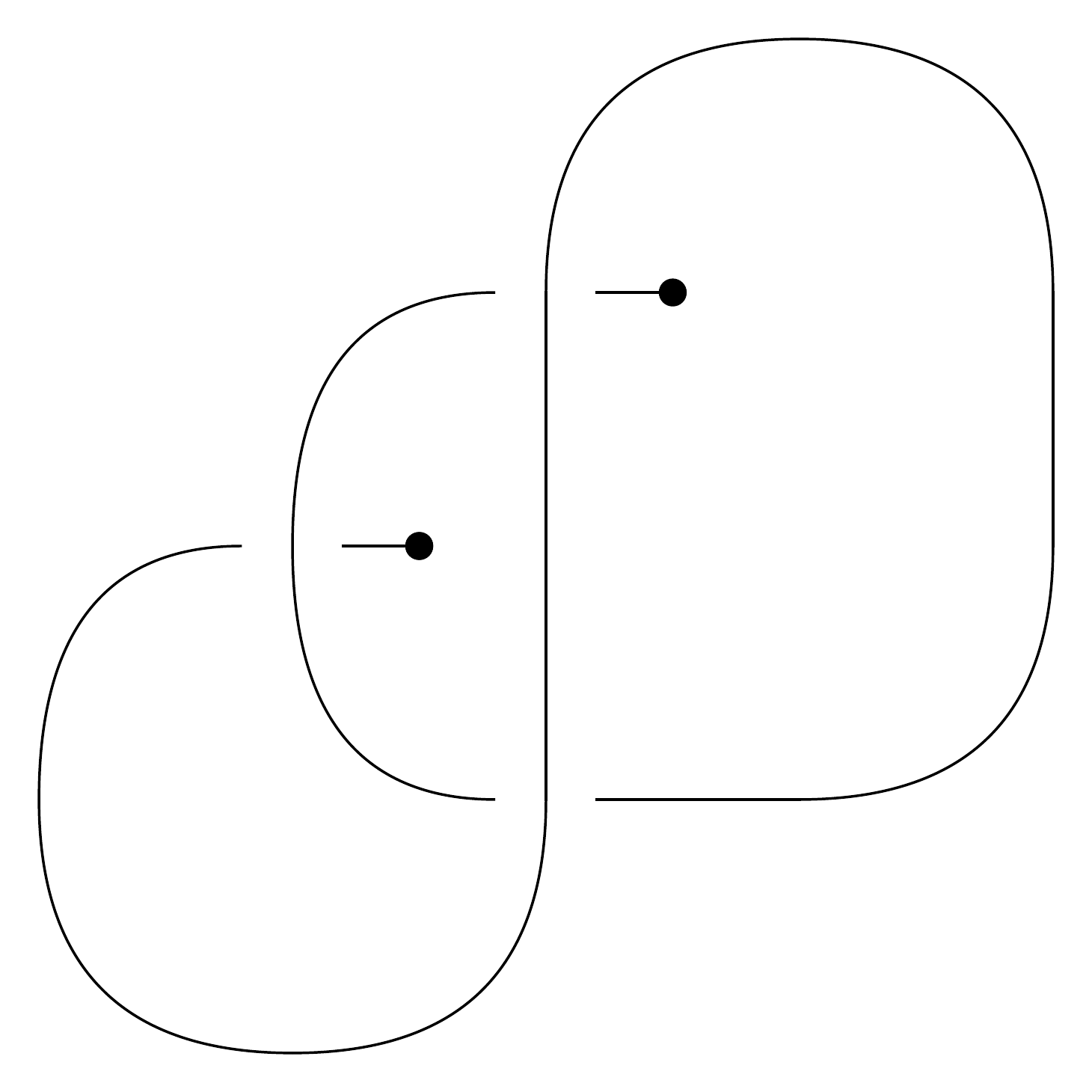}\\
\textcolor{black}{$3_{18}$}
\vspace{1cm}
\end{minipage}
\begin{minipage}[t]{.25\linewidth}
\centering
\includegraphics[width=0.9\textwidth,height=3.5cm,keepaspectratio]{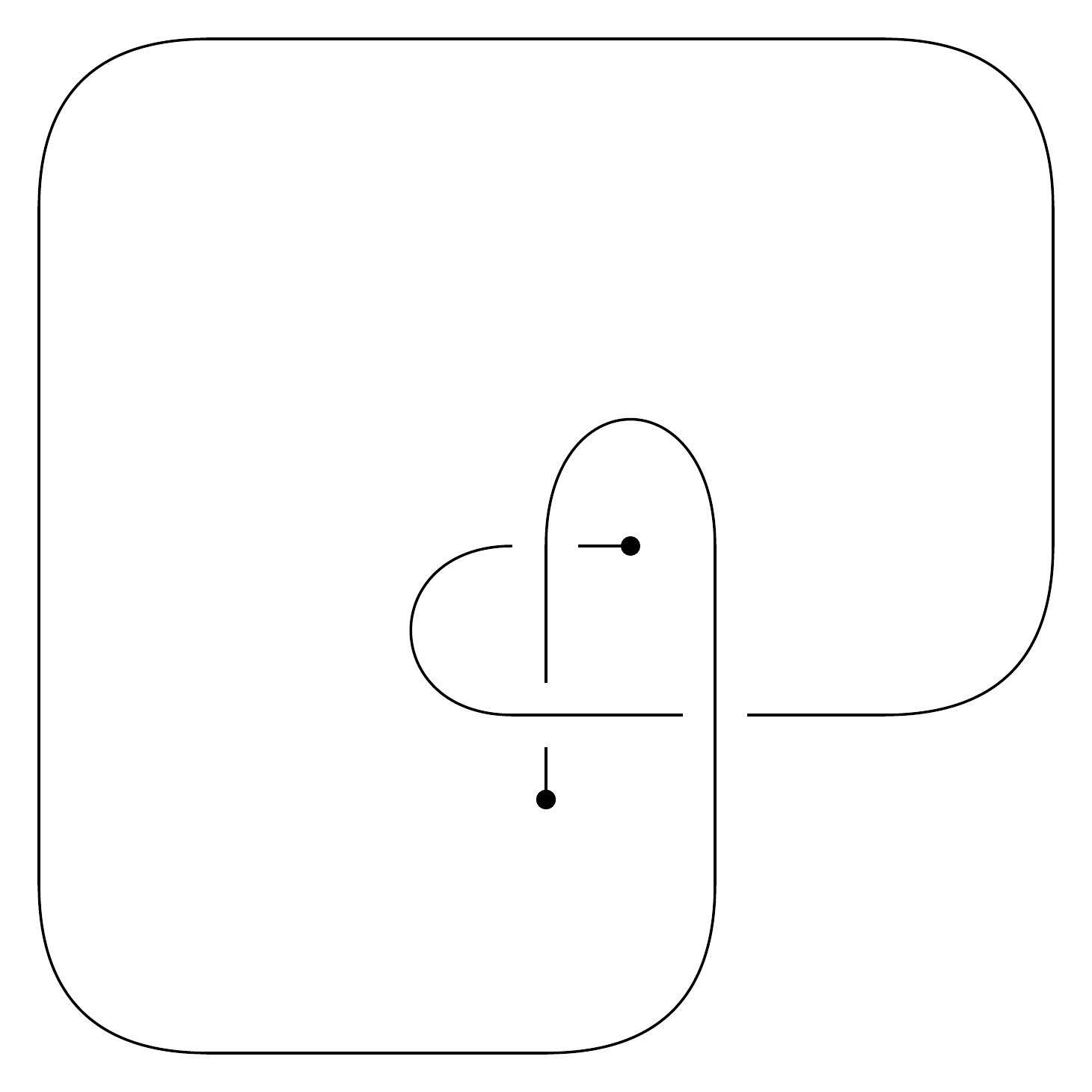}\\
\textcolor{black}{$3_{19}$}
\vspace{1cm}
\end{minipage}
\begin{minipage}[t]{.25\linewidth}
\centering
\includegraphics[width=0.9\textwidth,height=3.5cm,keepaspectratio]{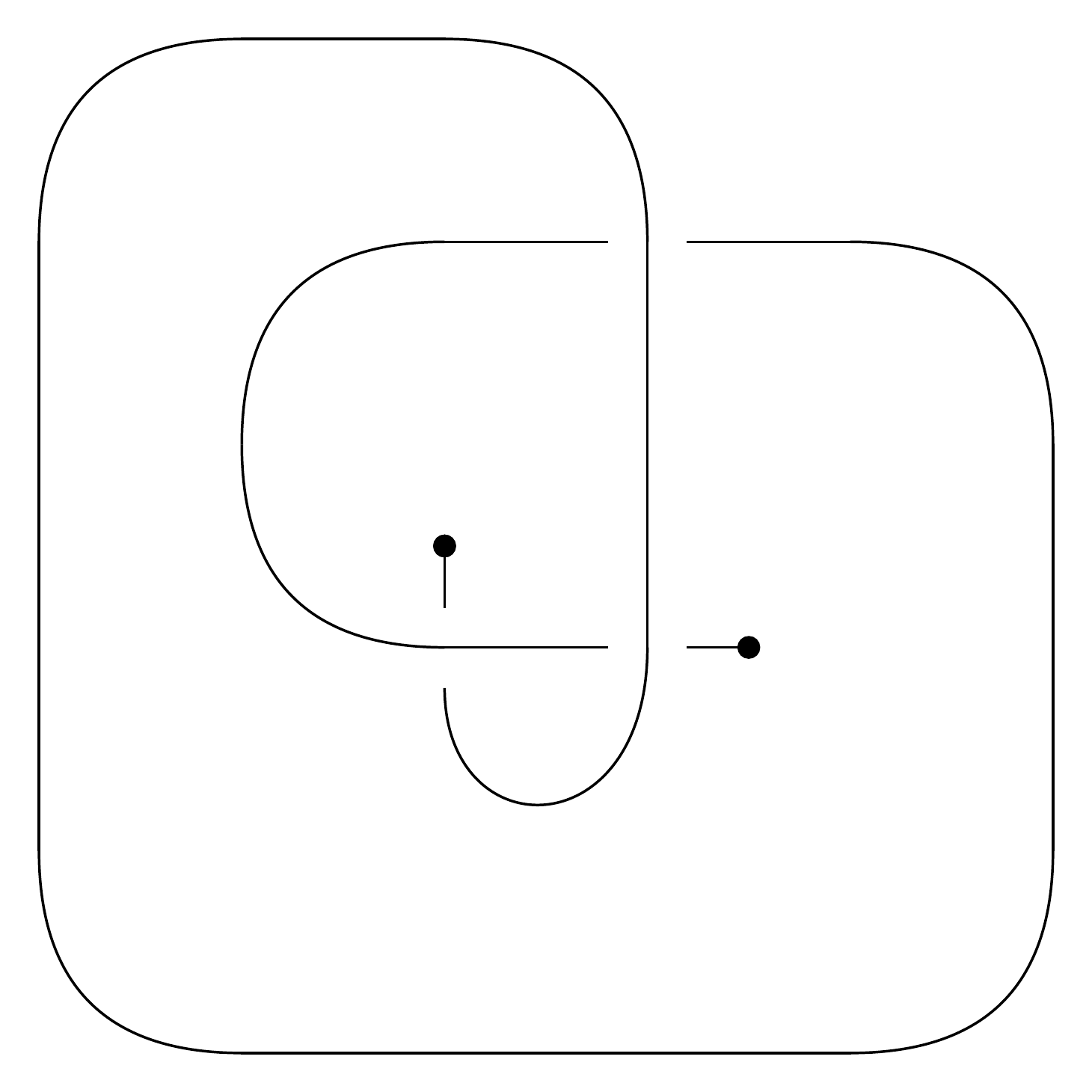}\\
\textcolor{black}{$3_{20}$}
\vspace{1cm}
\end{minipage}
\begin{minipage}[t]{.25\linewidth}
\centering
\includegraphics[width=0.9\textwidth,height=3.5cm,keepaspectratio]{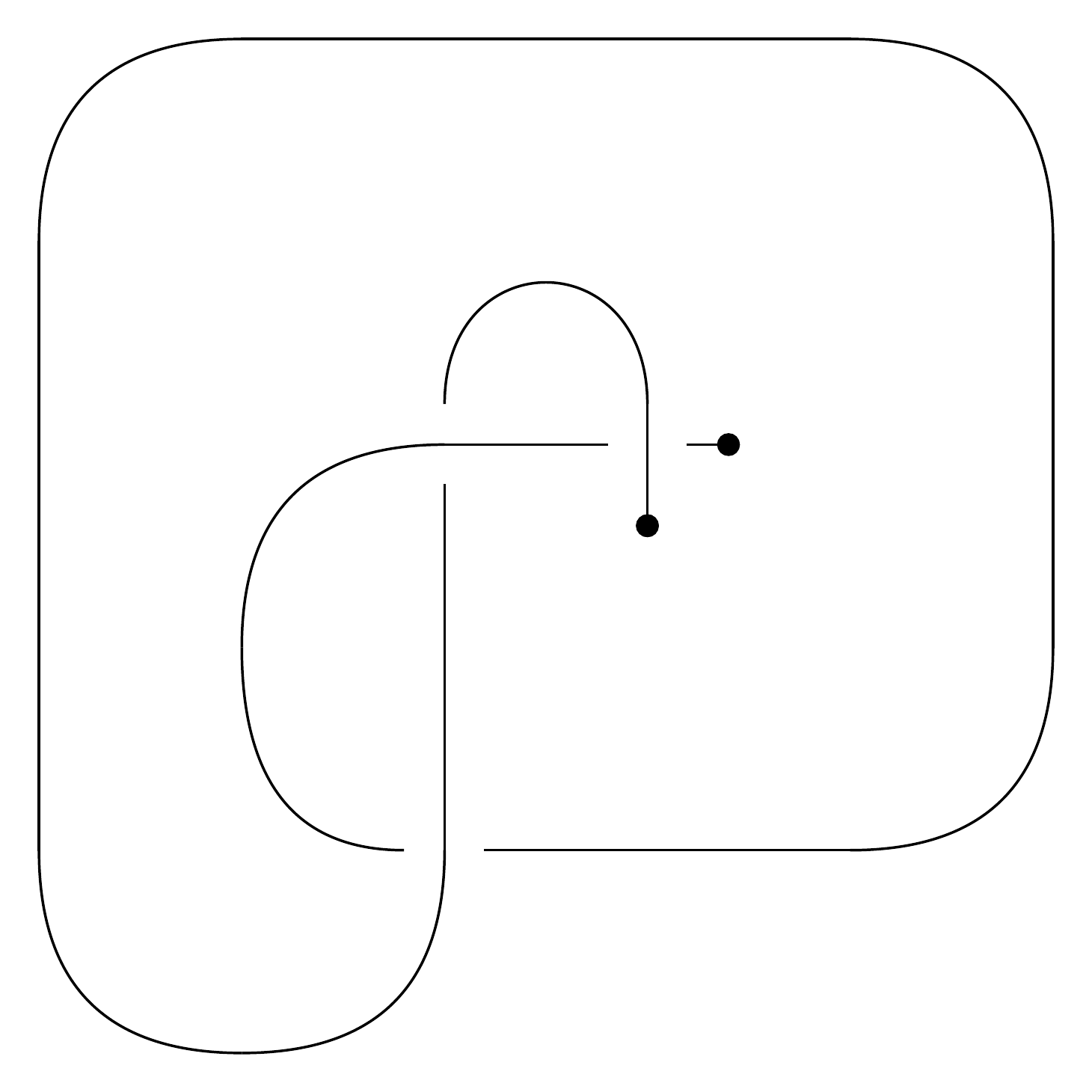}\\
\textcolor{black}{$3_{21}$}
\vspace{1cm}
\end{minipage}
\begin{minipage}[t]{.25\linewidth}
\centering
\includegraphics[width=0.9\textwidth,height=3.5cm,keepaspectratio]{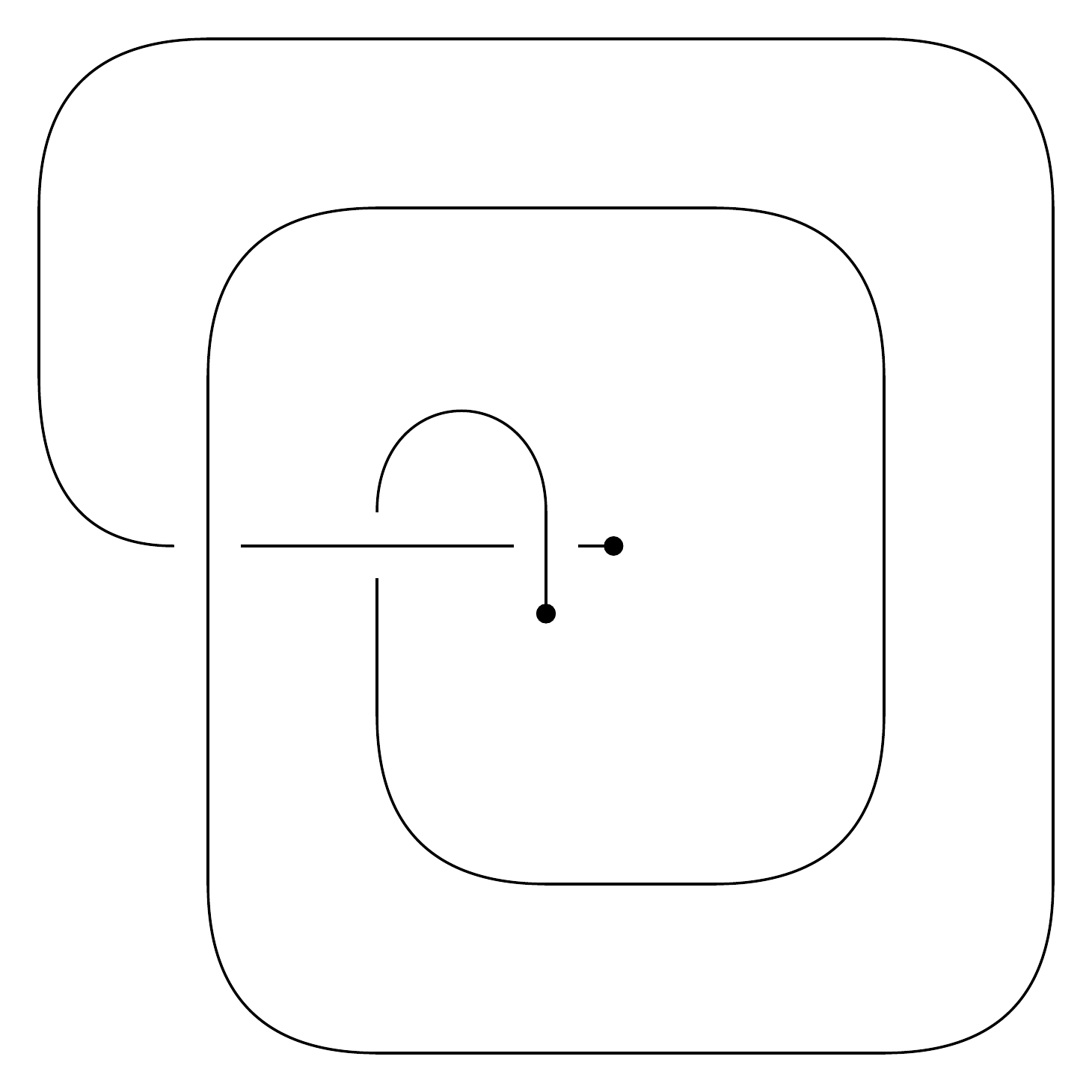}\\
\textcolor{black}{$3_{22}$}
\vspace{1cm}
\end{minipage}
\begin{minipage}[t]{.25\linewidth}
\centering
\includegraphics[width=0.9\textwidth,height=3.5cm,keepaspectratio]{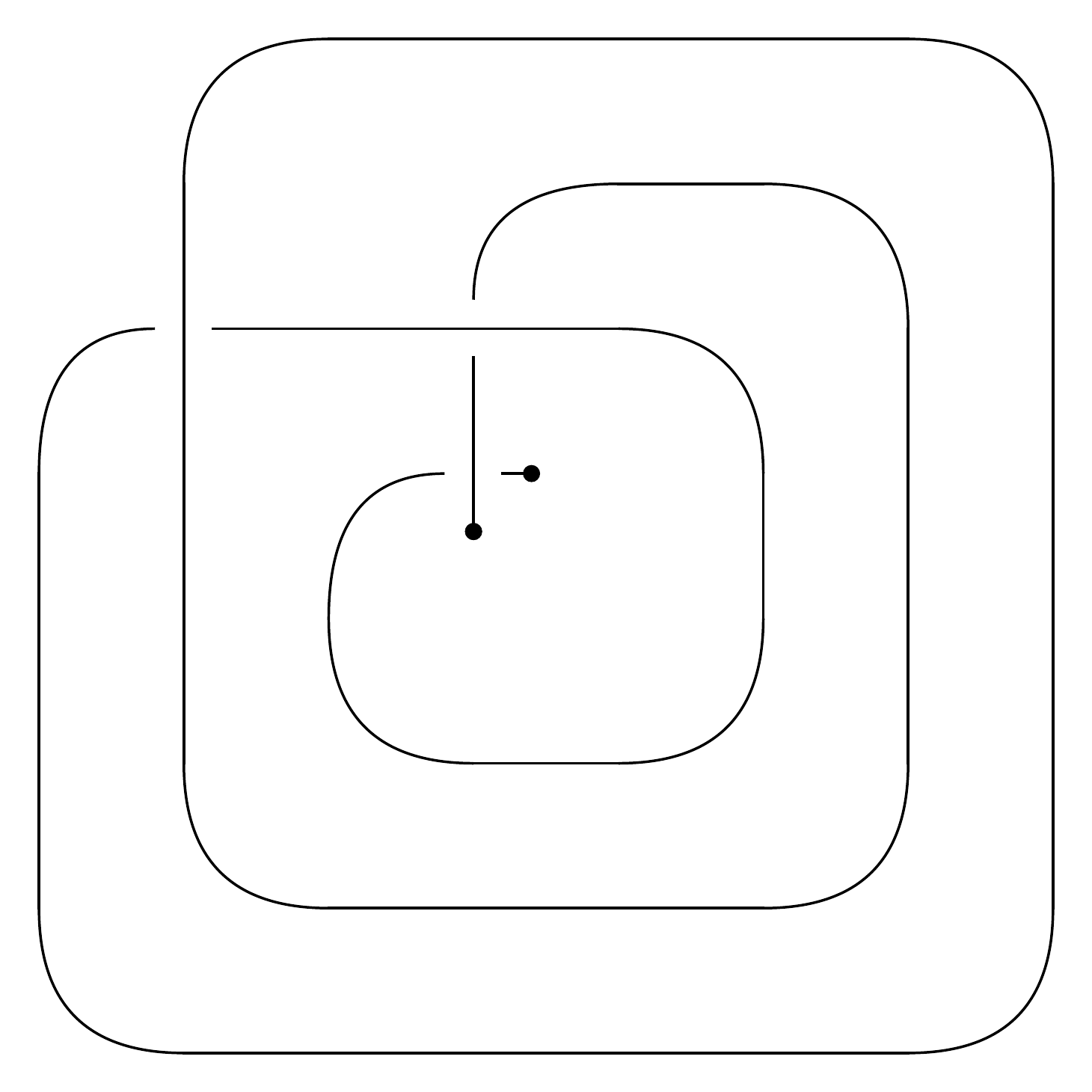}\\
\textcolor{black}{$3_{23}$}
\vspace{1cm}
\end{minipage}
\begin{minipage}[t]{.25\linewidth}
\centering
\includegraphics[width=0.9\textwidth,height=3.5cm,keepaspectratio]{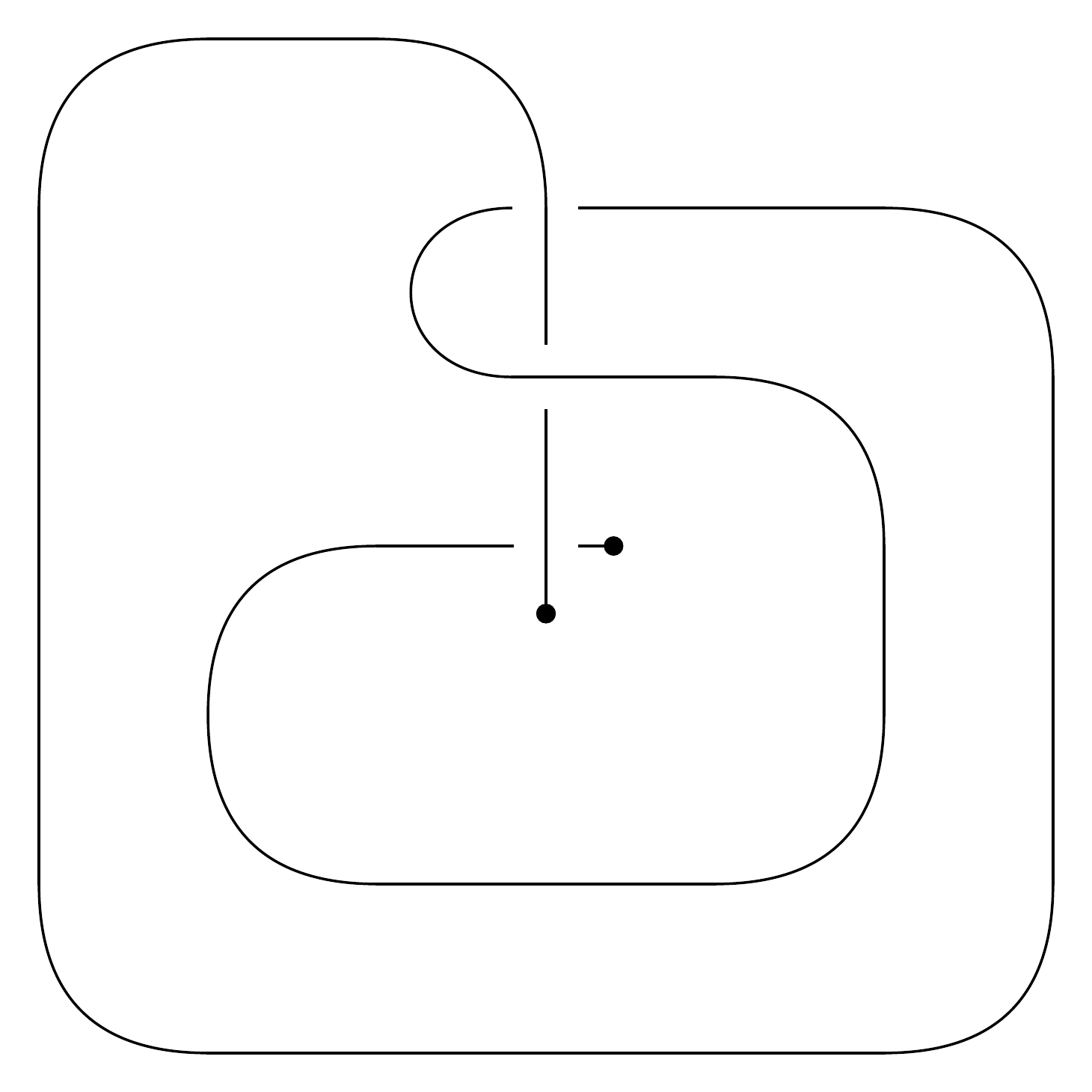}\\
\textcolor{black}{$3_{24}$}
\vspace{1cm}
\end{minipage}
\begin{minipage}[t]{.25\linewidth}
\centering
\includegraphics[width=0.9\textwidth,height=3.5cm,keepaspectratio]{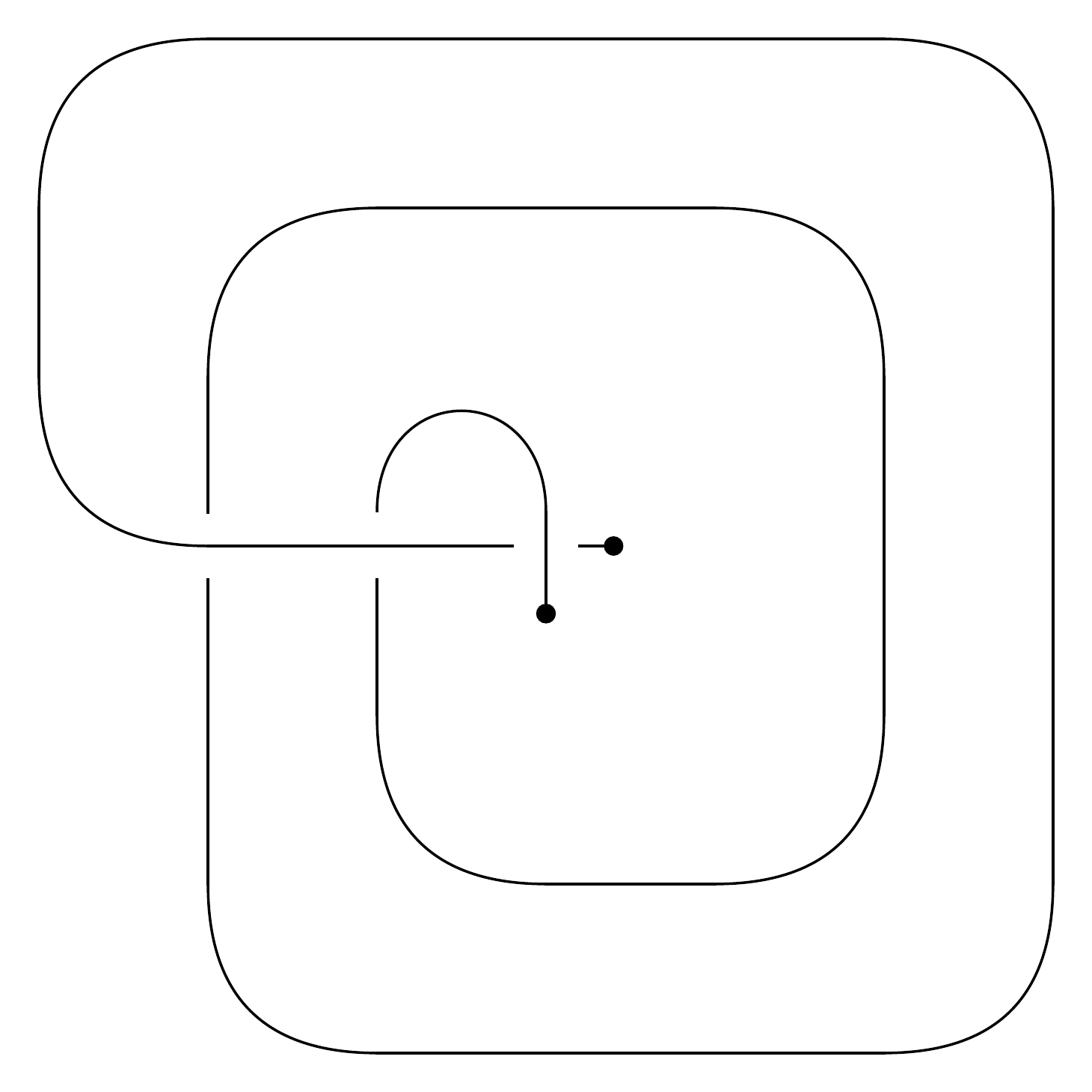}\\
\textcolor{black}{$3_{25}$}
\vspace{1cm}
\end{minipage}
\begin{minipage}[t]{.25\linewidth}
\centering
\includegraphics[width=0.9\textwidth,height=3.5cm,keepaspectratio]{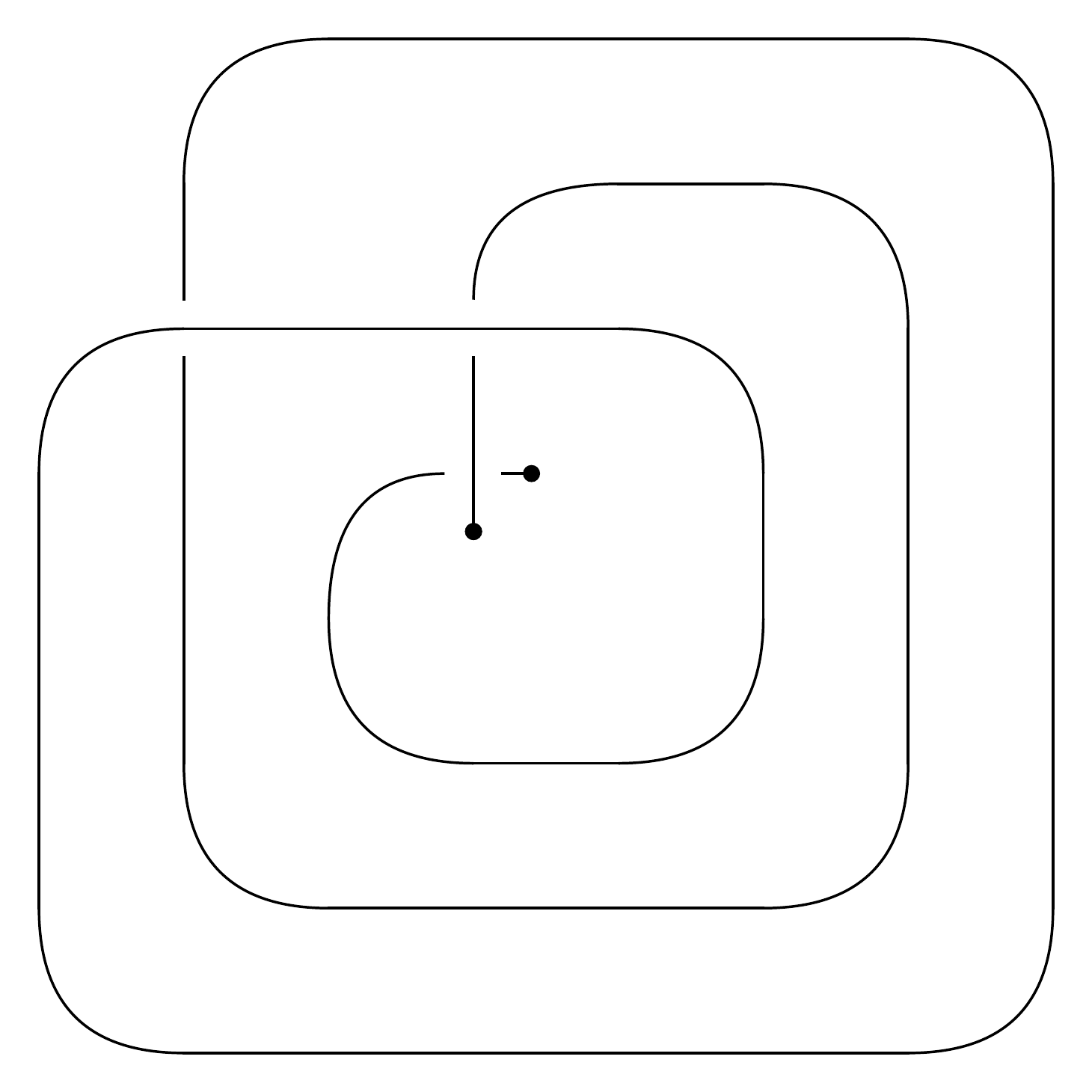}\\
\textcolor{black}{$3_{26}$}
\vspace{1cm}
\end{minipage}
\begin{minipage}[t]{.25\linewidth}
\centering
\includegraphics[width=0.9\textwidth,height=3.5cm,keepaspectratio]{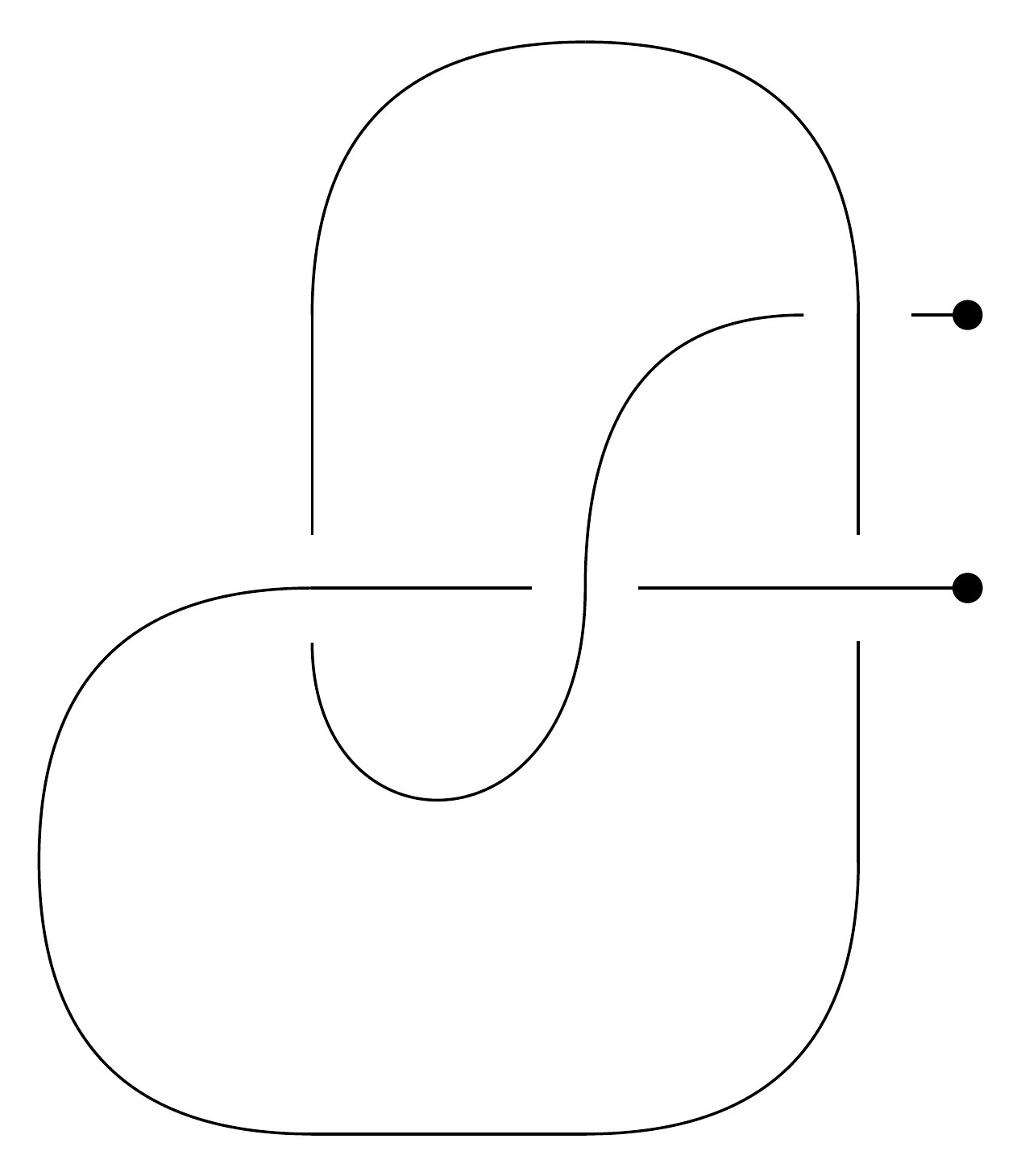}\\
\textcolor{red}{$4_{1}$}
\vspace{1cm}
\end{minipage}
\begin{minipage}[t]{.25\linewidth}
\centering
\includegraphics[width=0.9\textwidth,height=3.5cm,keepaspectratio]{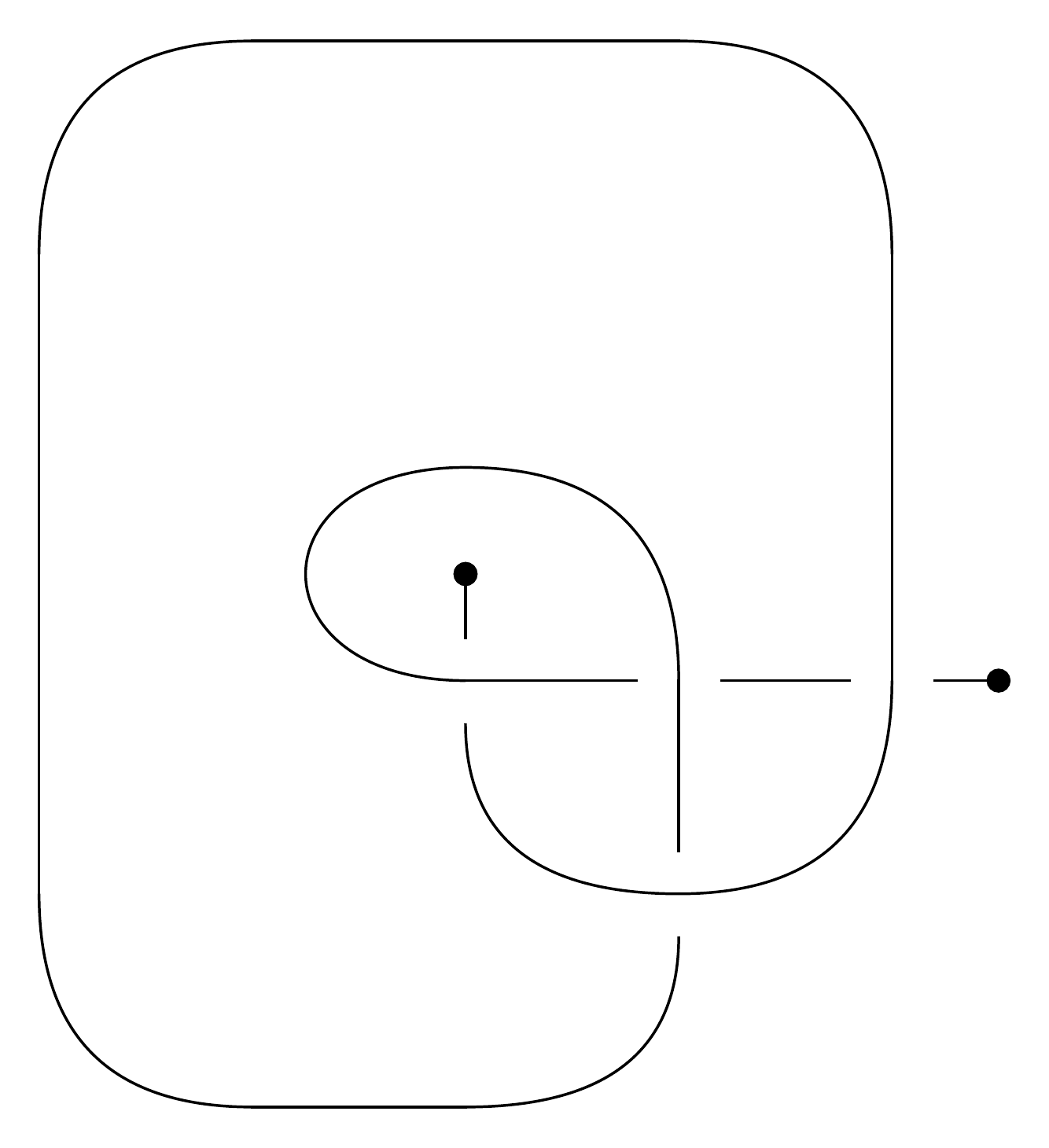}\\
\textcolor{blue}{$4_{2}$}
\vspace{1cm}
\end{minipage}
\begin{minipage}[t]{.25\linewidth}
\centering
\includegraphics[width=0.9\textwidth,height=3.5cm,keepaspectratio]{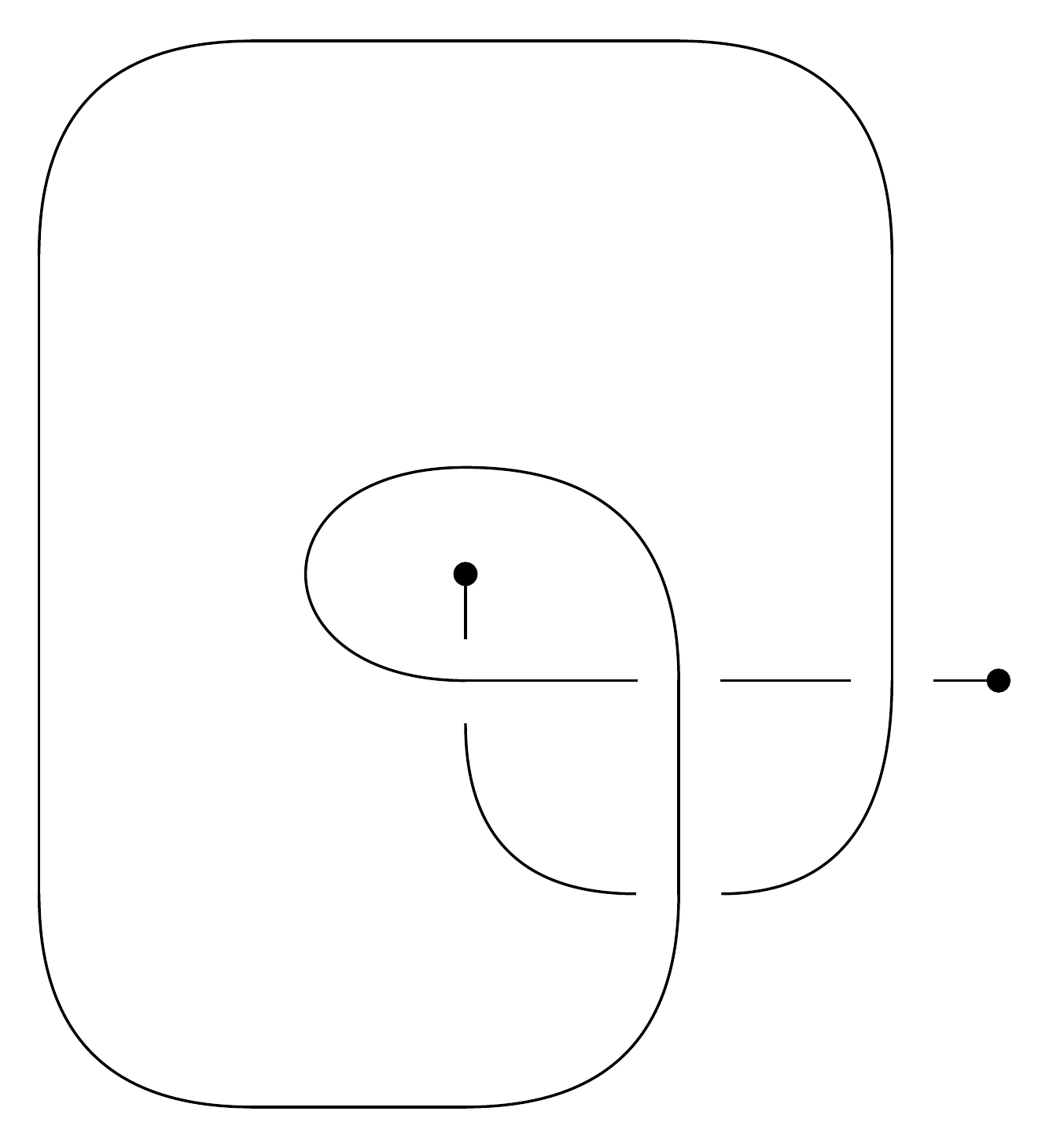}\\
\textcolor{blue}{$4_{3}$}
\vspace{1cm}
\end{minipage}
\begin{minipage}[t]{.25\linewidth}
\centering
\includegraphics[width=0.9\textwidth,height=3.5cm,keepaspectratio]{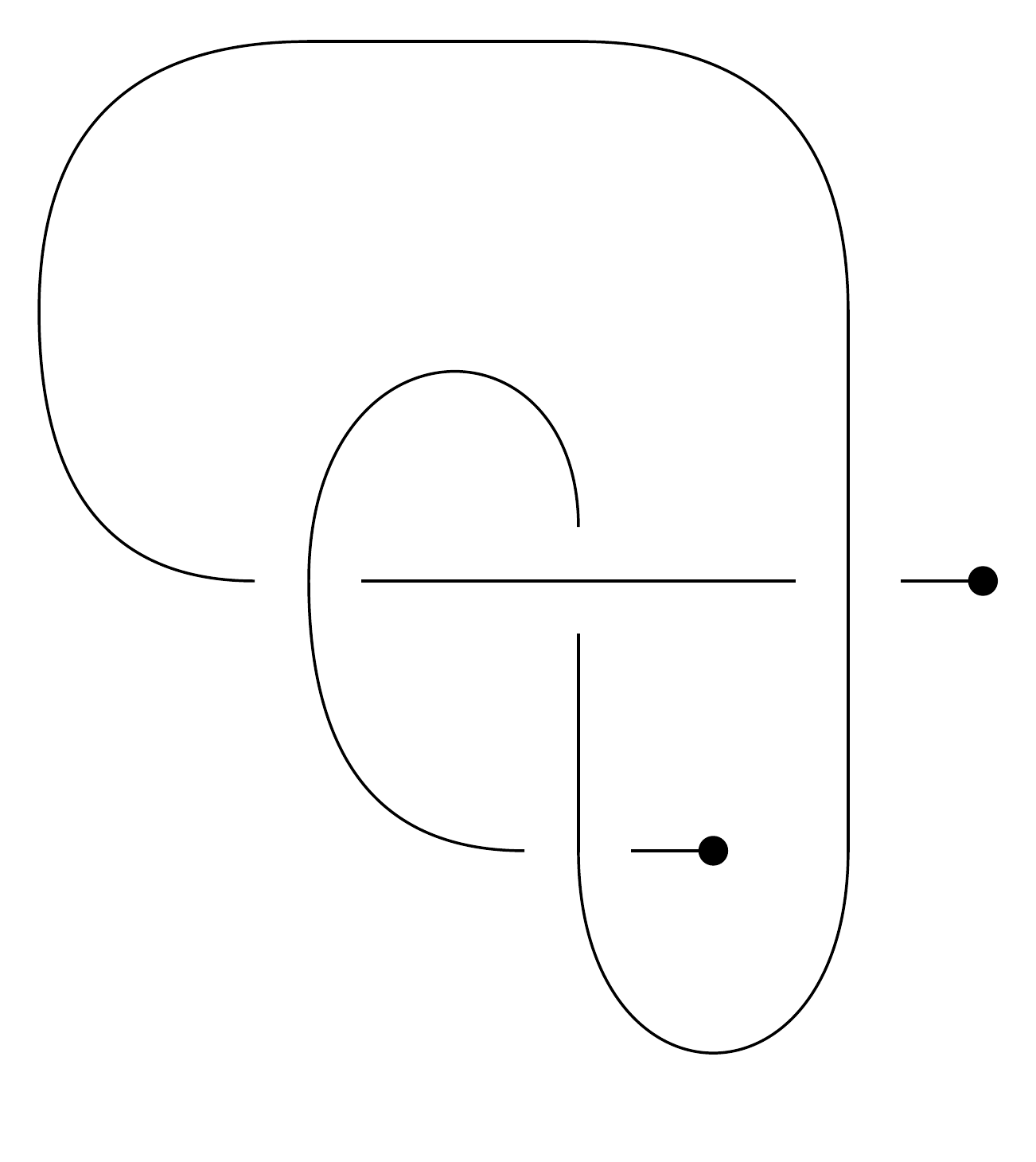}\\
\textcolor{blue}{$4_{4}$}
\vspace{1cm}
\end{minipage}
\begin{minipage}[t]{.25\linewidth}
\centering
\includegraphics[width=0.9\textwidth,height=3.5cm,keepaspectratio]{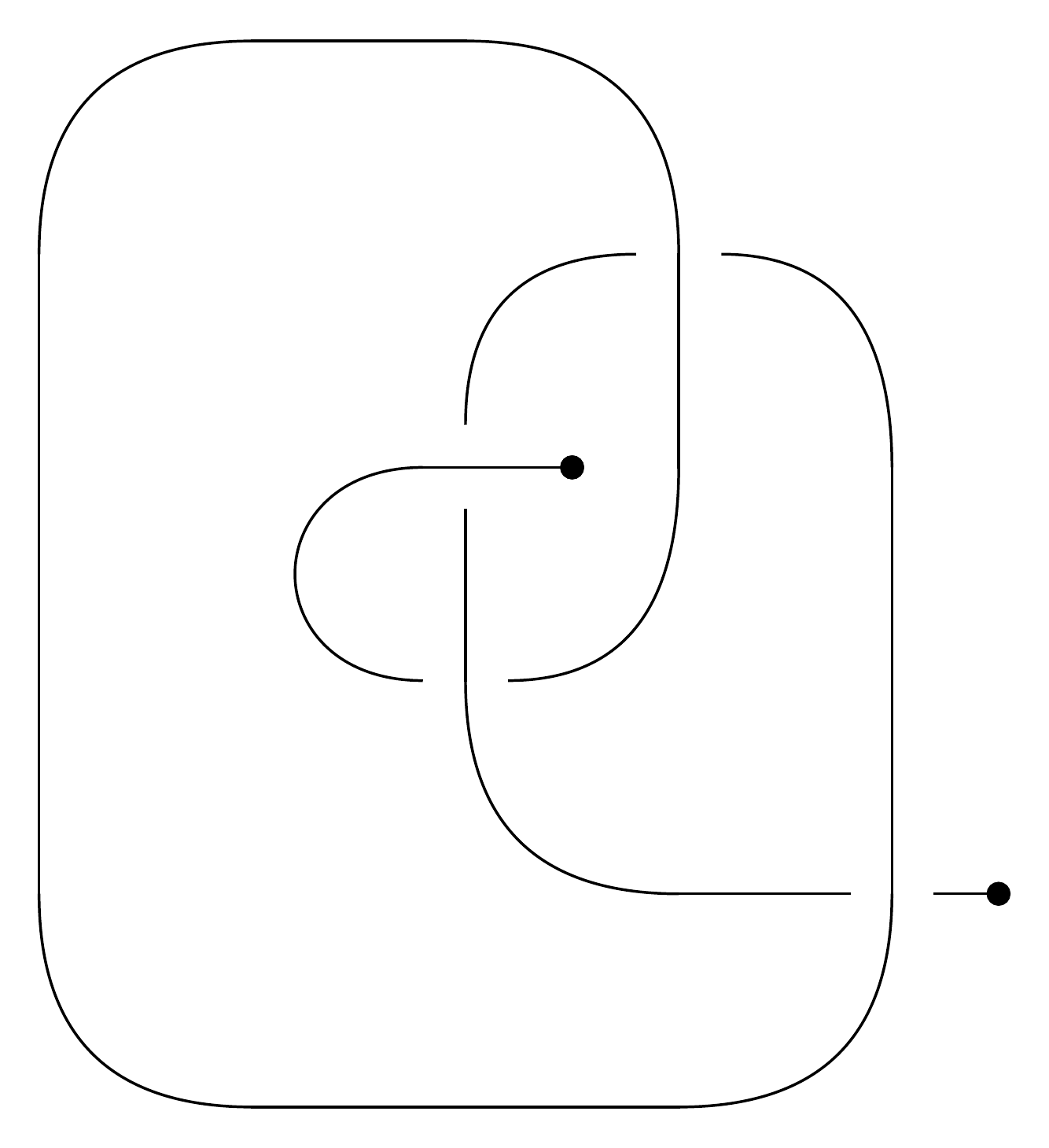}\\
\textcolor{blue}{$4_{5}$}
\vspace{1cm}
\end{minipage}
\begin{minipage}[t]{.25\linewidth}
\centering
\includegraphics[width=0.9\textwidth,height=3.5cm,keepaspectratio]{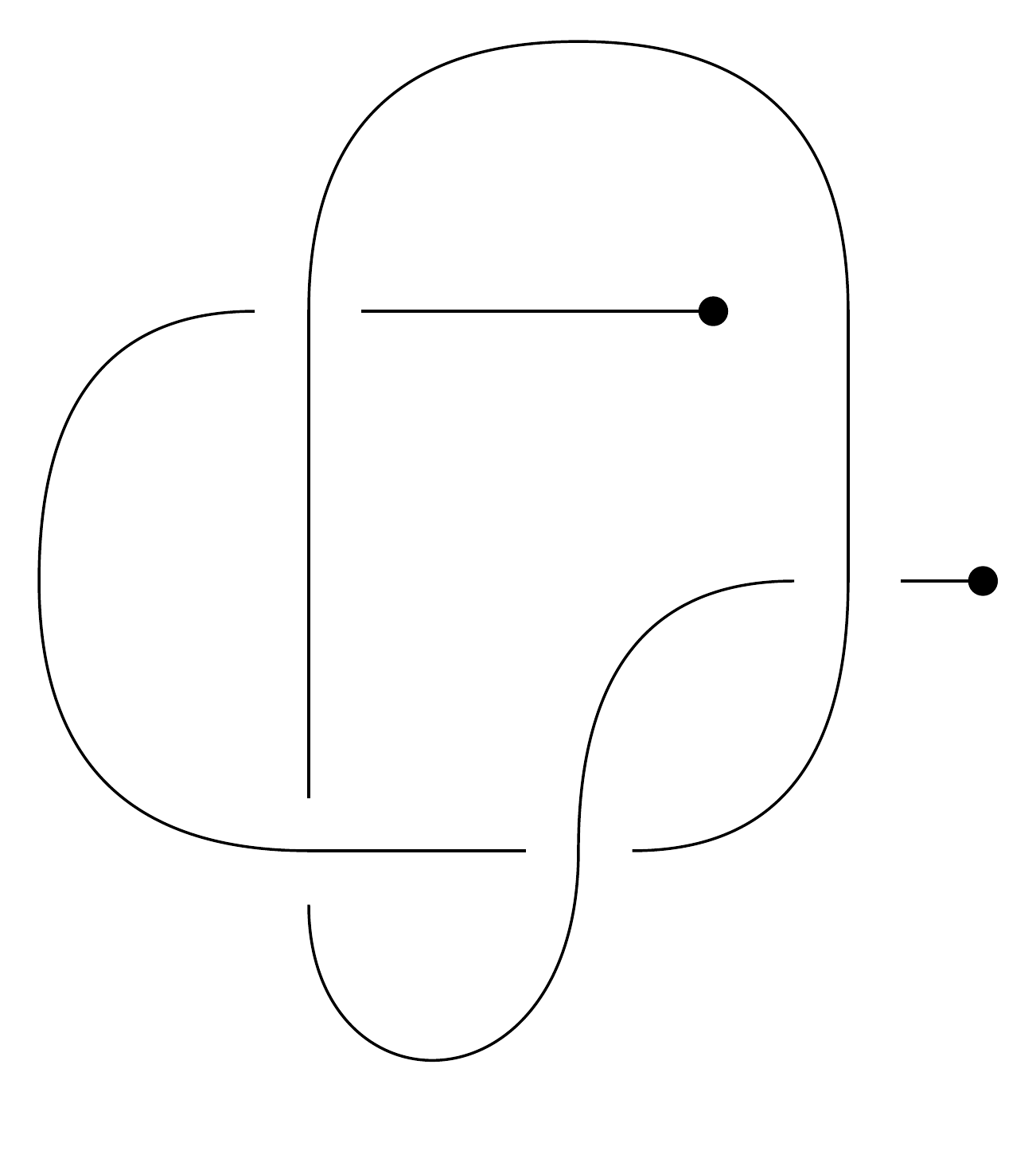}\\
\textcolor{blue}{$4_{6}$}
\vspace{1cm}
\end{minipage}
\begin{minipage}[t]{.25\linewidth}
\centering
\includegraphics[width=0.9\textwidth,height=3.5cm,keepaspectratio]{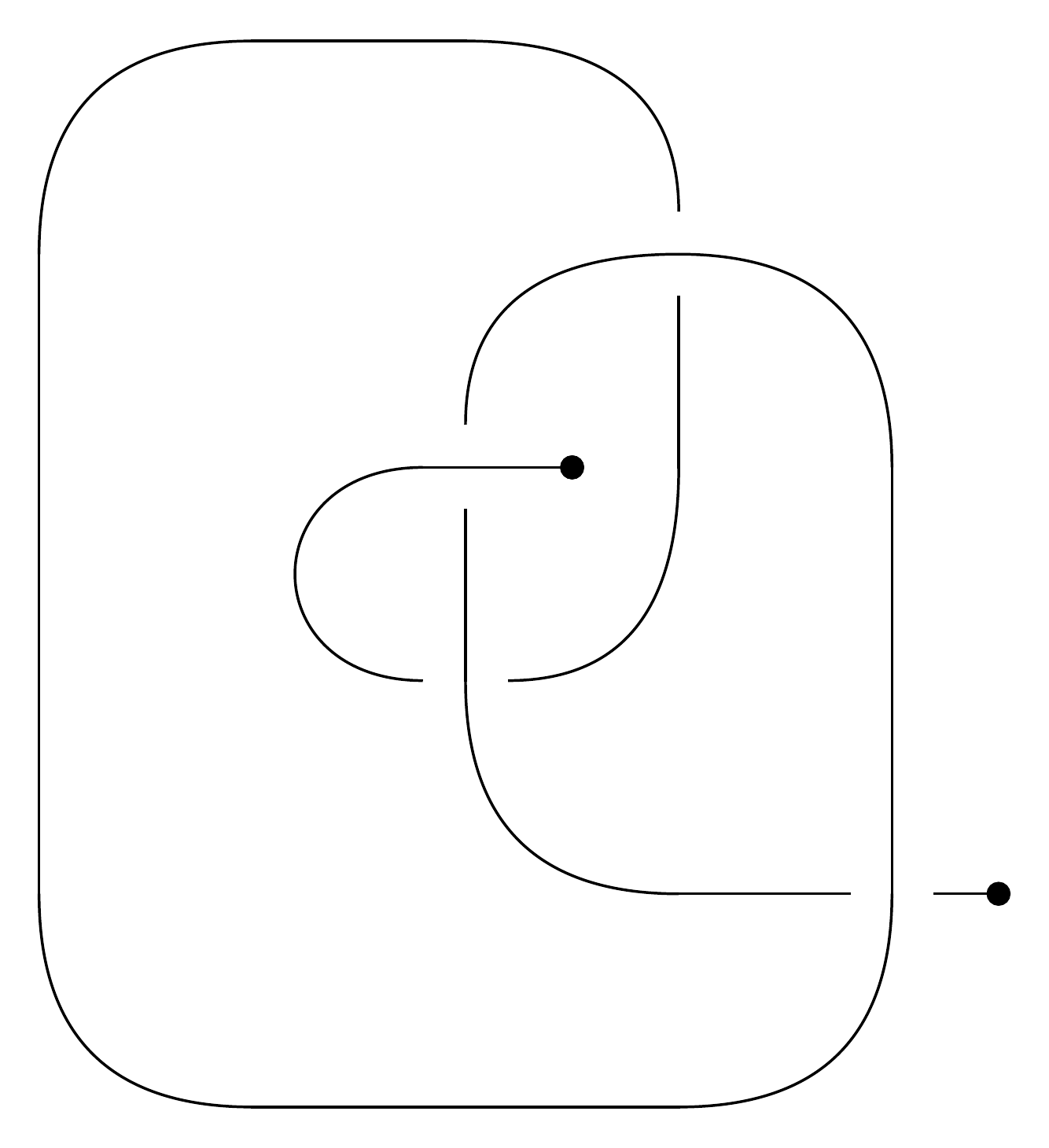}\\
\textcolor{blue}{$4_{7}$}
\vspace{1cm}
\end{minipage}
\begin{minipage}[t]{.25\linewidth}
\centering
\includegraphics[width=0.9\textwidth,height=3.5cm,keepaspectratio]{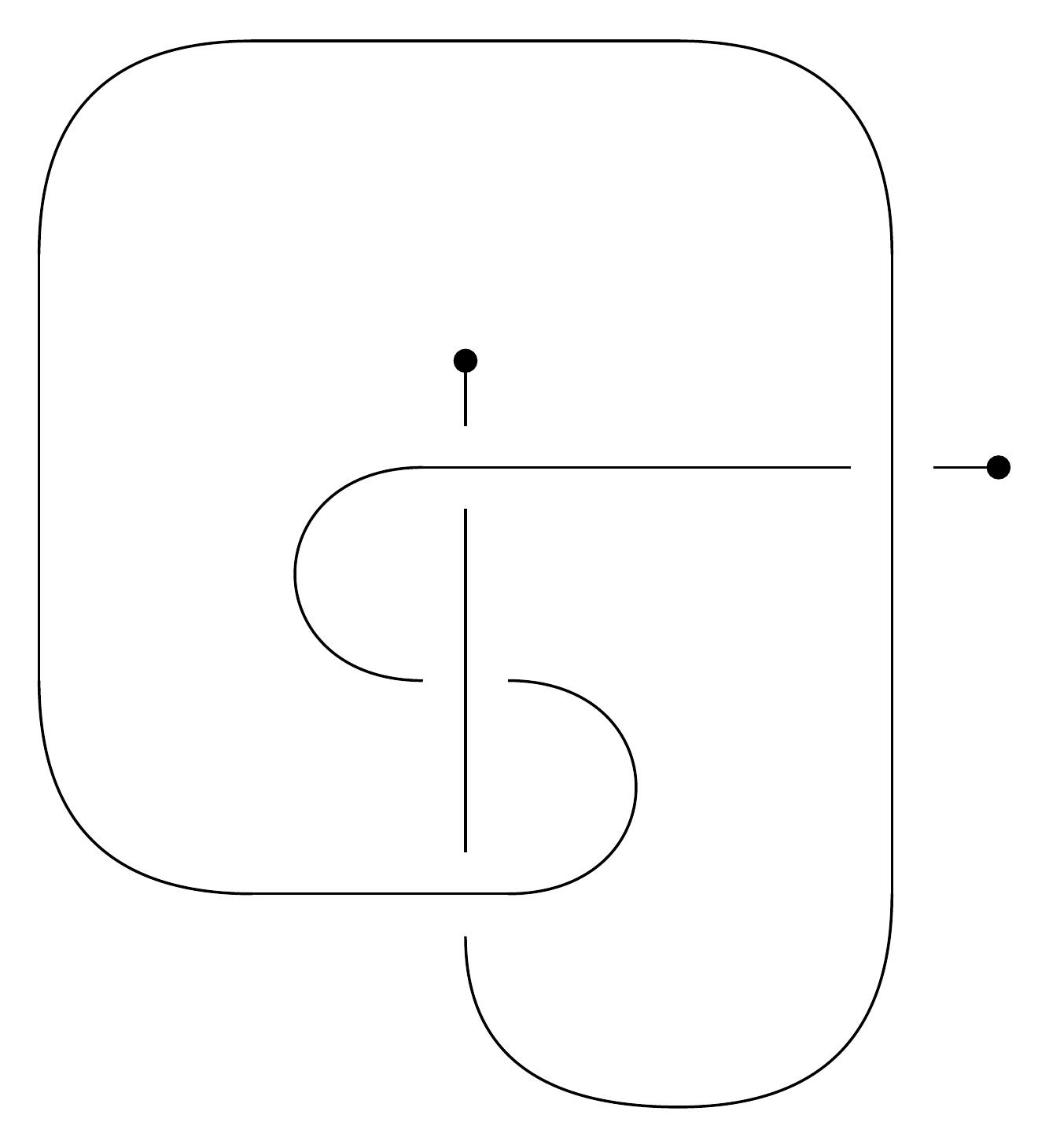}\\
\textcolor{blue}{$4_{8}$}
\vspace{1cm}
\end{minipage}
\begin{minipage}[t]{.25\linewidth}
\centering
\includegraphics[width=0.9\textwidth,height=3.5cm,keepaspectratio]{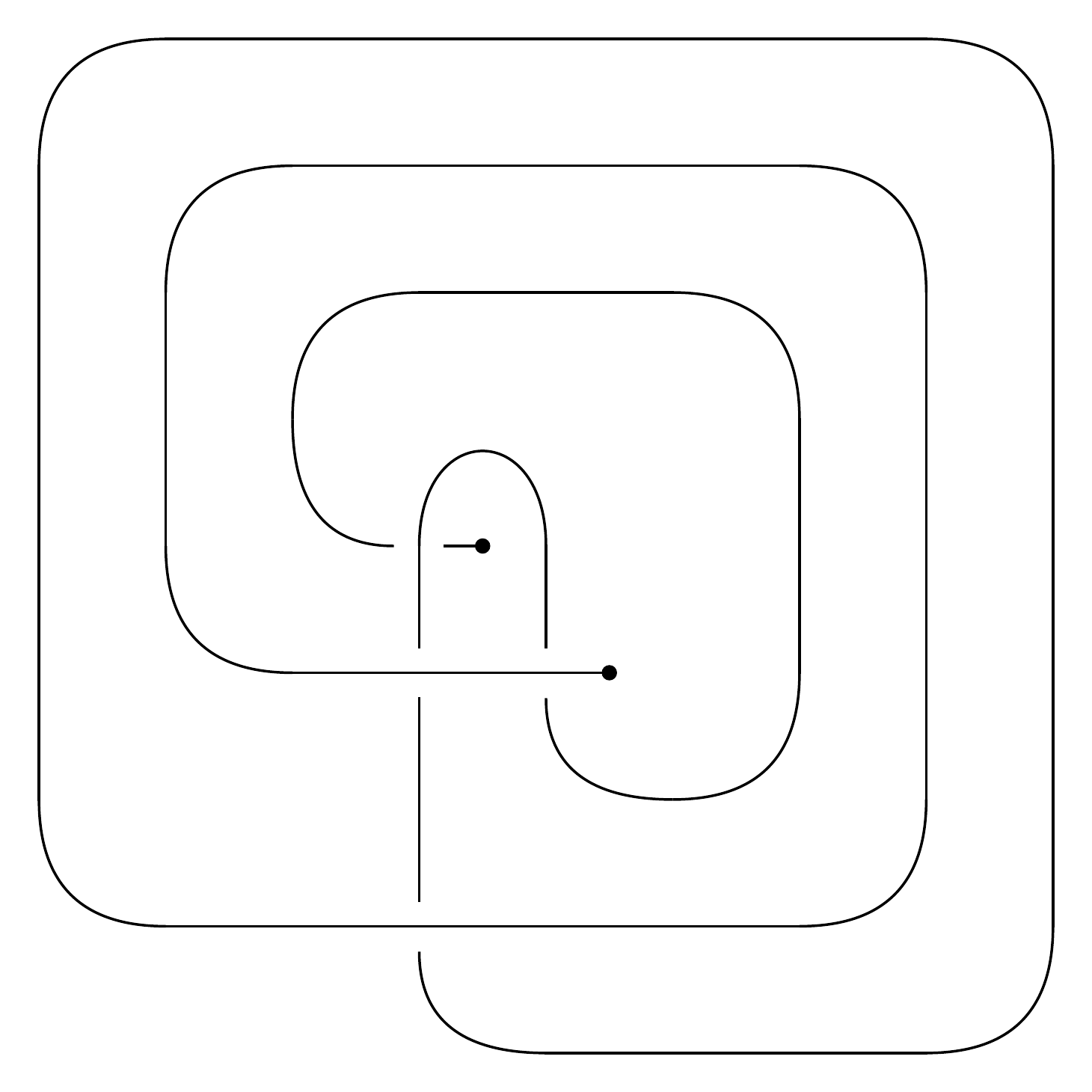}\\
\textcolor{black}{$4_{9}$}
\vspace{1cm}
\end{minipage}
\begin{minipage}[t]{.25\linewidth}
\centering
\includegraphics[width=0.9\textwidth,height=3.5cm,keepaspectratio]{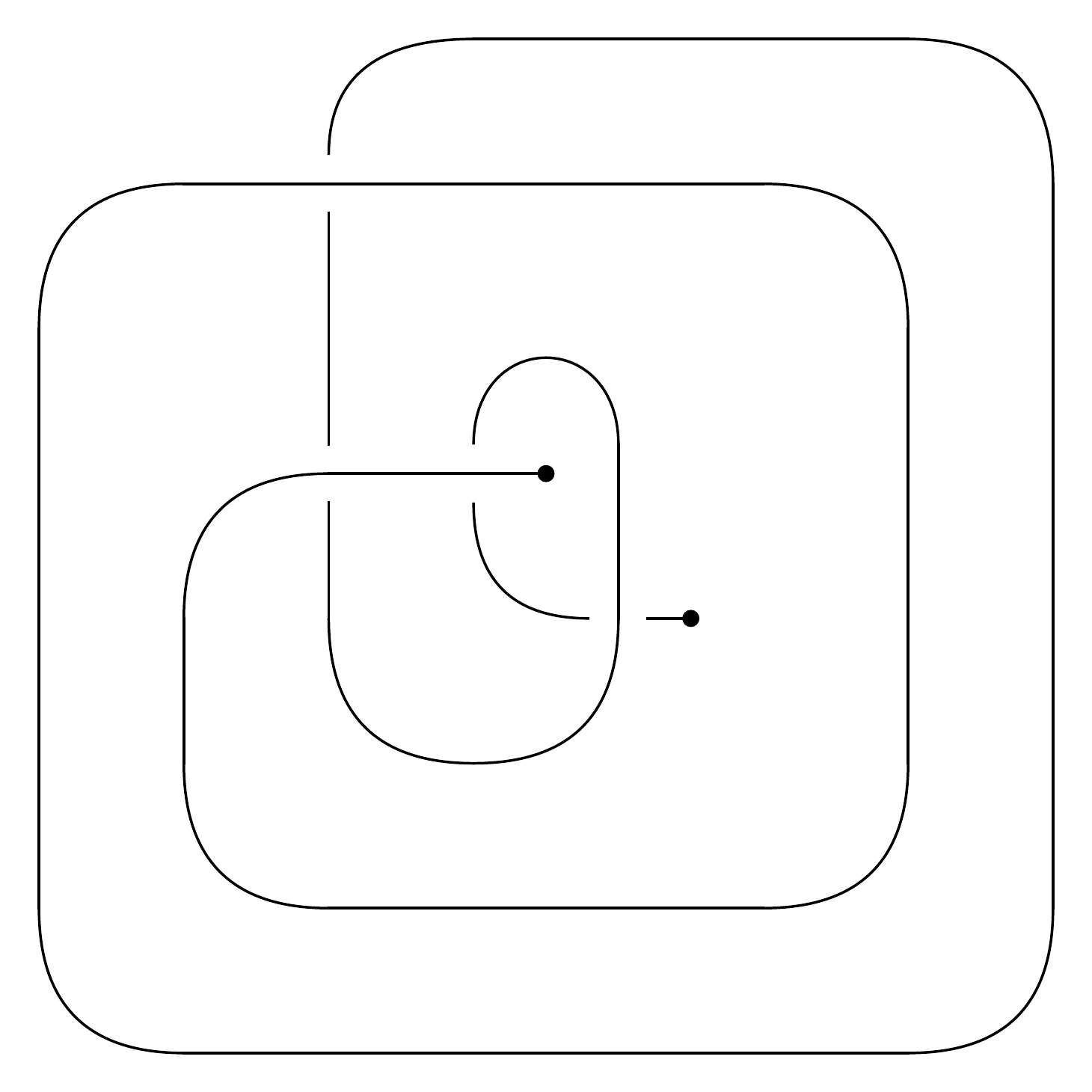}\\
\textcolor{black}{$4_{10}$}
\vspace{1cm}
\end{minipage}
\begin{minipage}[t]{.25\linewidth}
\centering
\includegraphics[width=0.9\textwidth,height=3.5cm,keepaspectratio]{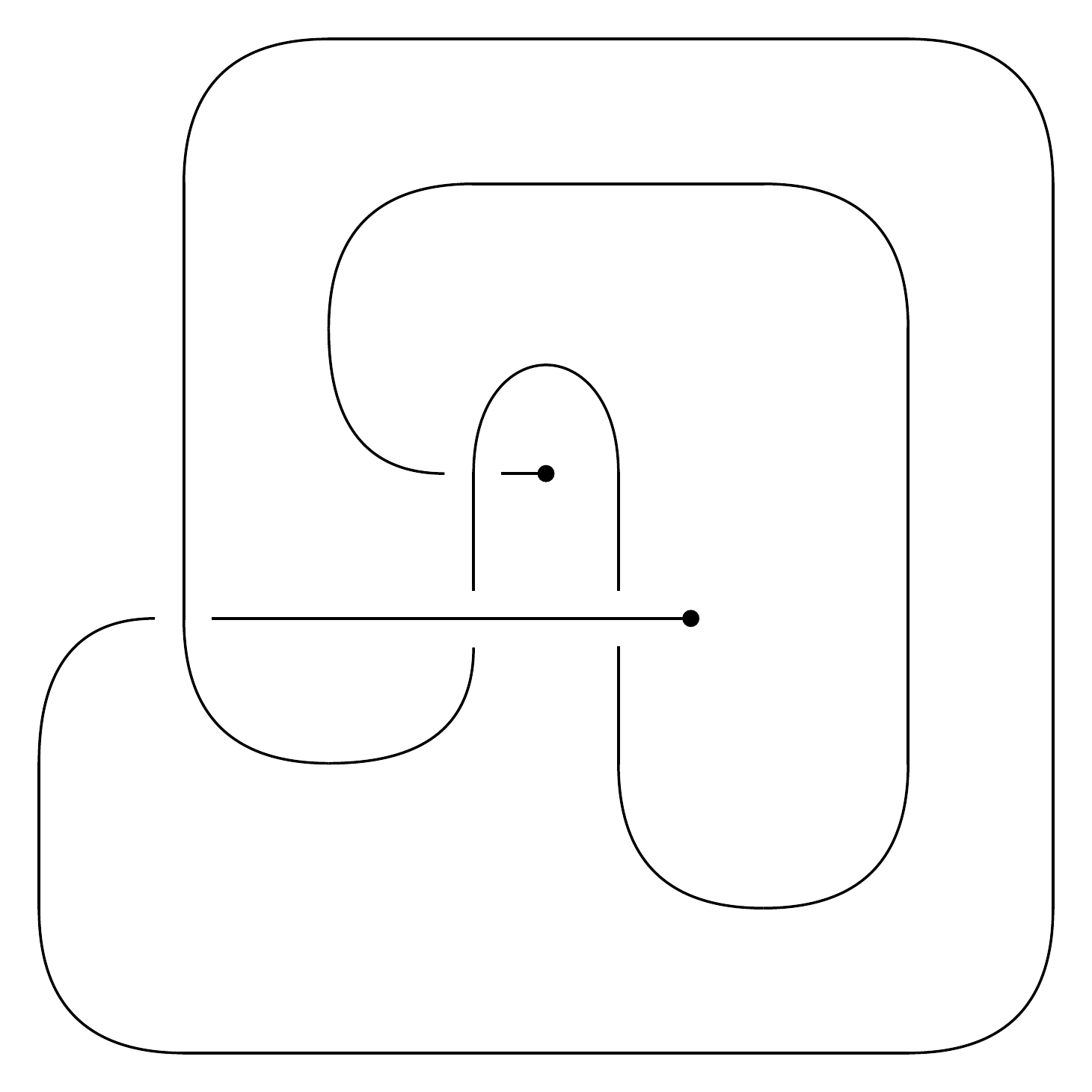}\\
\textcolor{black}{$4_{11}$}
\vspace{1cm}
\end{minipage}
\begin{minipage}[t]{.25\linewidth}
\centering
\includegraphics[width=0.9\textwidth,height=3.5cm,keepaspectratio]{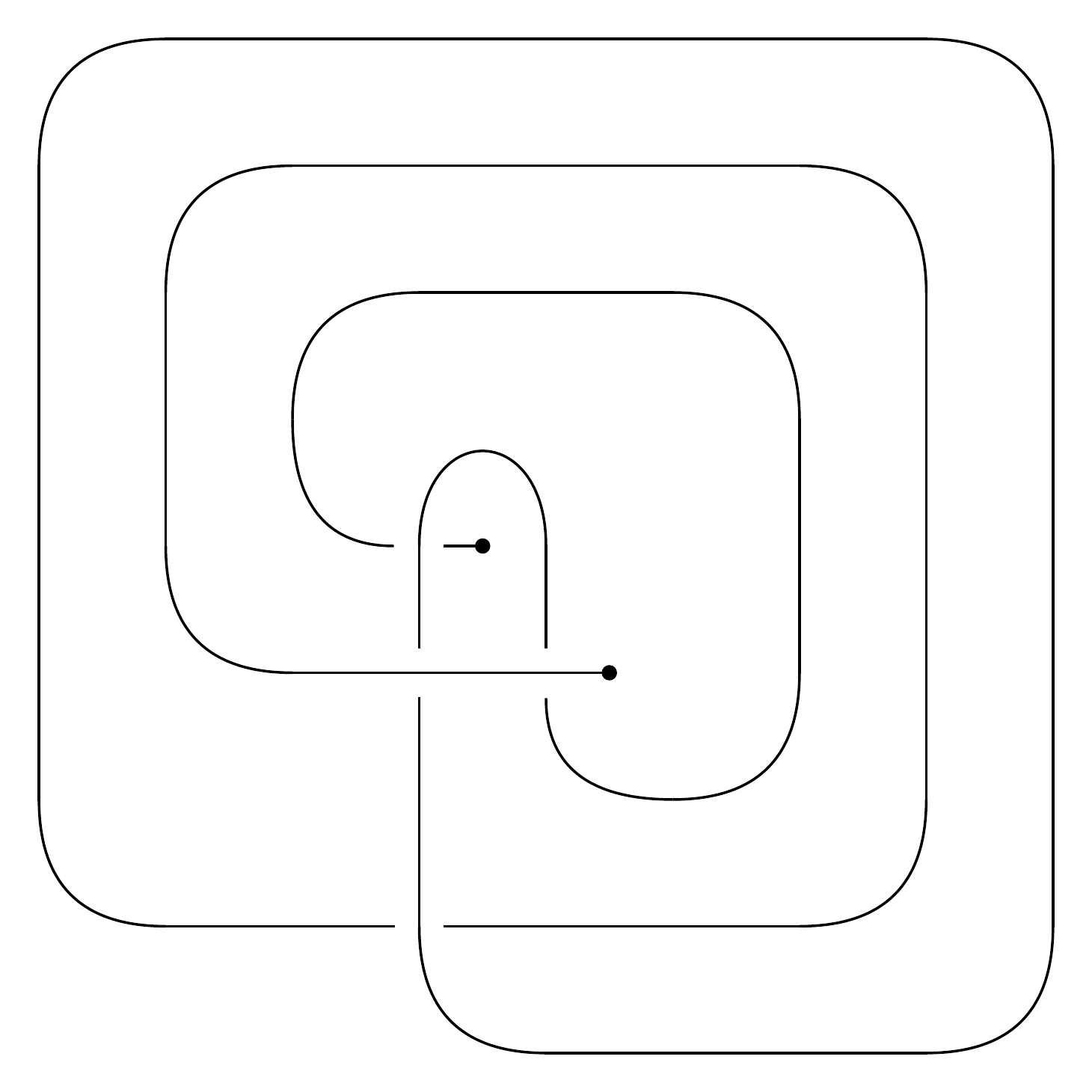}\\
\textcolor{black}{$4_{12}$}
\vspace{1cm}
\end{minipage}
\begin{minipage}[t]{.25\linewidth}
\centering
\includegraphics[width=0.9\textwidth,height=3.5cm,keepaspectratio]{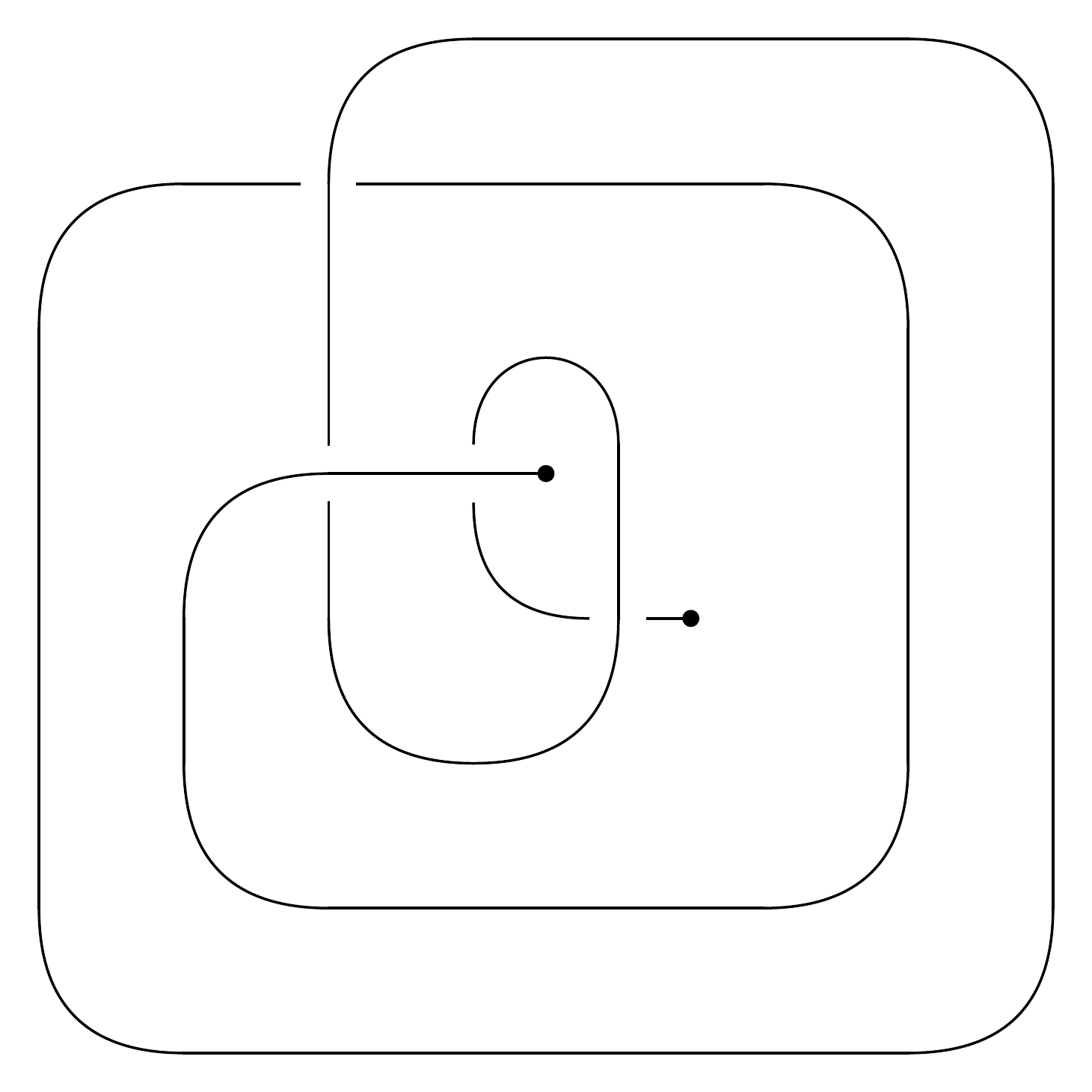}\\
\textcolor{black}{$4_{13}$}
\vspace{1cm}
\end{minipage}
\begin{minipage}[t]{.25\linewidth}
\centering
\includegraphics[width=0.9\textwidth,height=3.5cm,keepaspectratio]{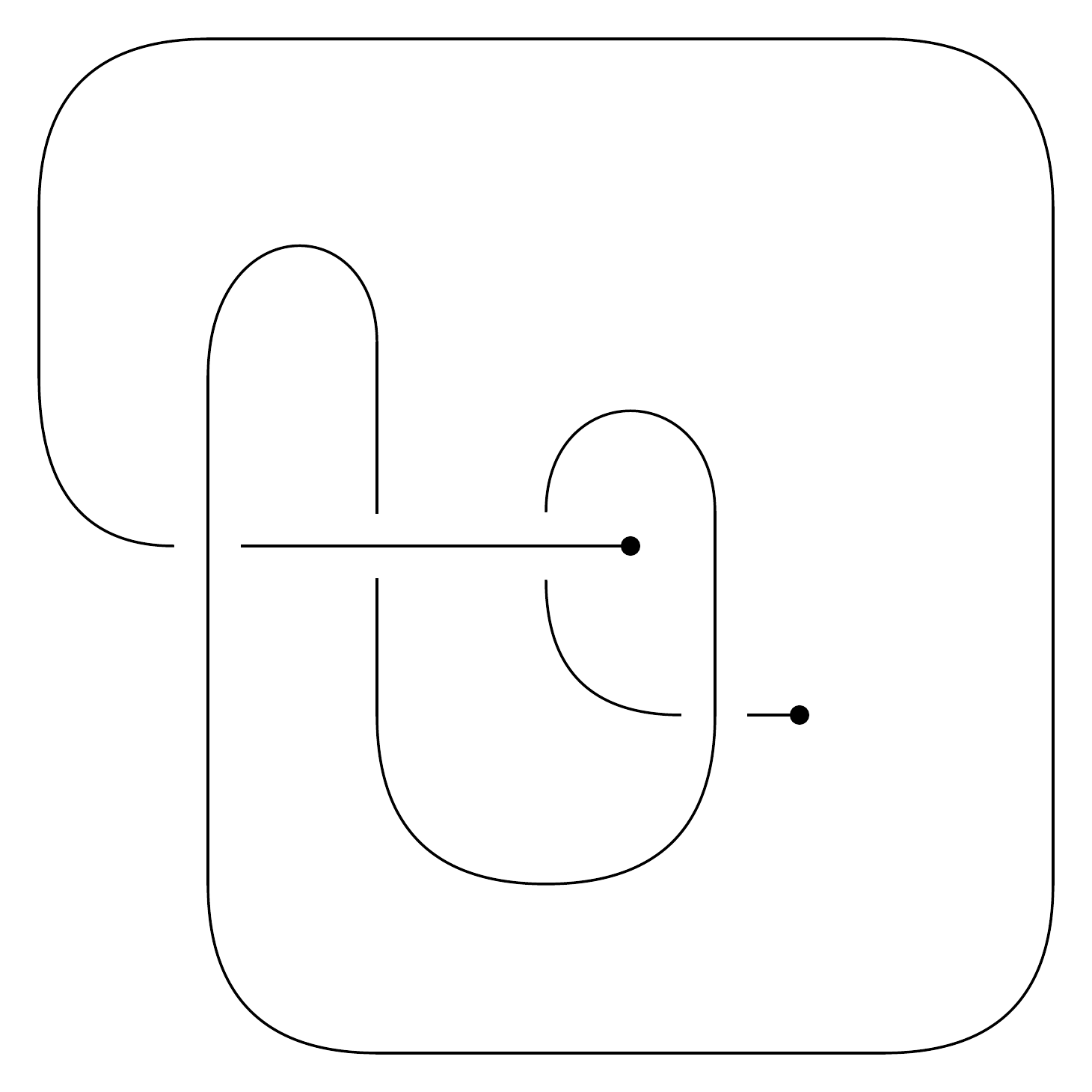}\\
\textcolor{black}{$4_{14}$}
\vspace{1cm}
\end{minipage}
\begin{minipage}[t]{.25\linewidth}
\centering
\includegraphics[width=0.9\textwidth,height=3.5cm,keepaspectratio]{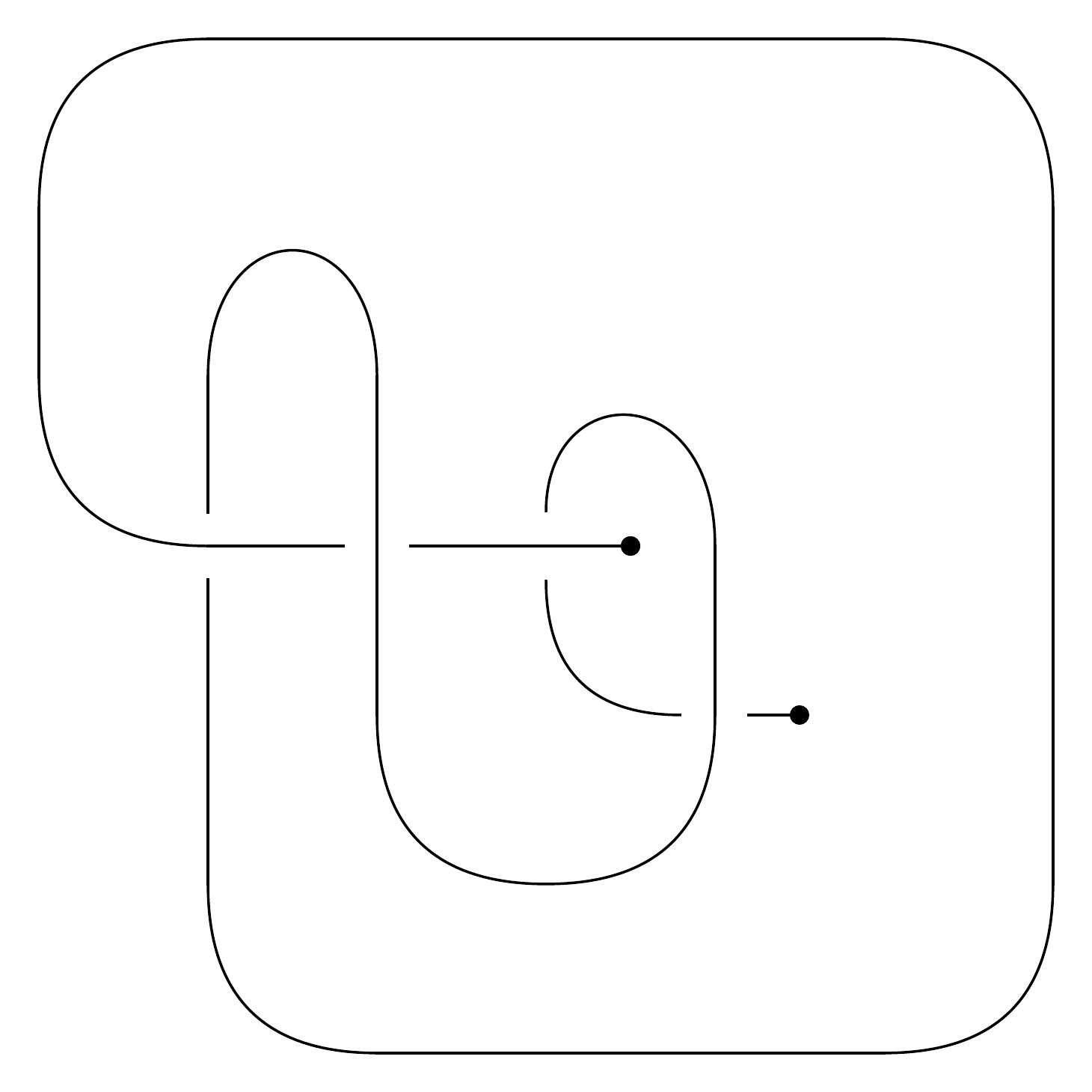}\\
\textcolor{black}{$4_{15}$}
\vspace{1cm}
\end{minipage}
\begin{minipage}[t]{.25\linewidth}
\centering
\includegraphics[width=0.9\textwidth,height=3.5cm,keepaspectratio]{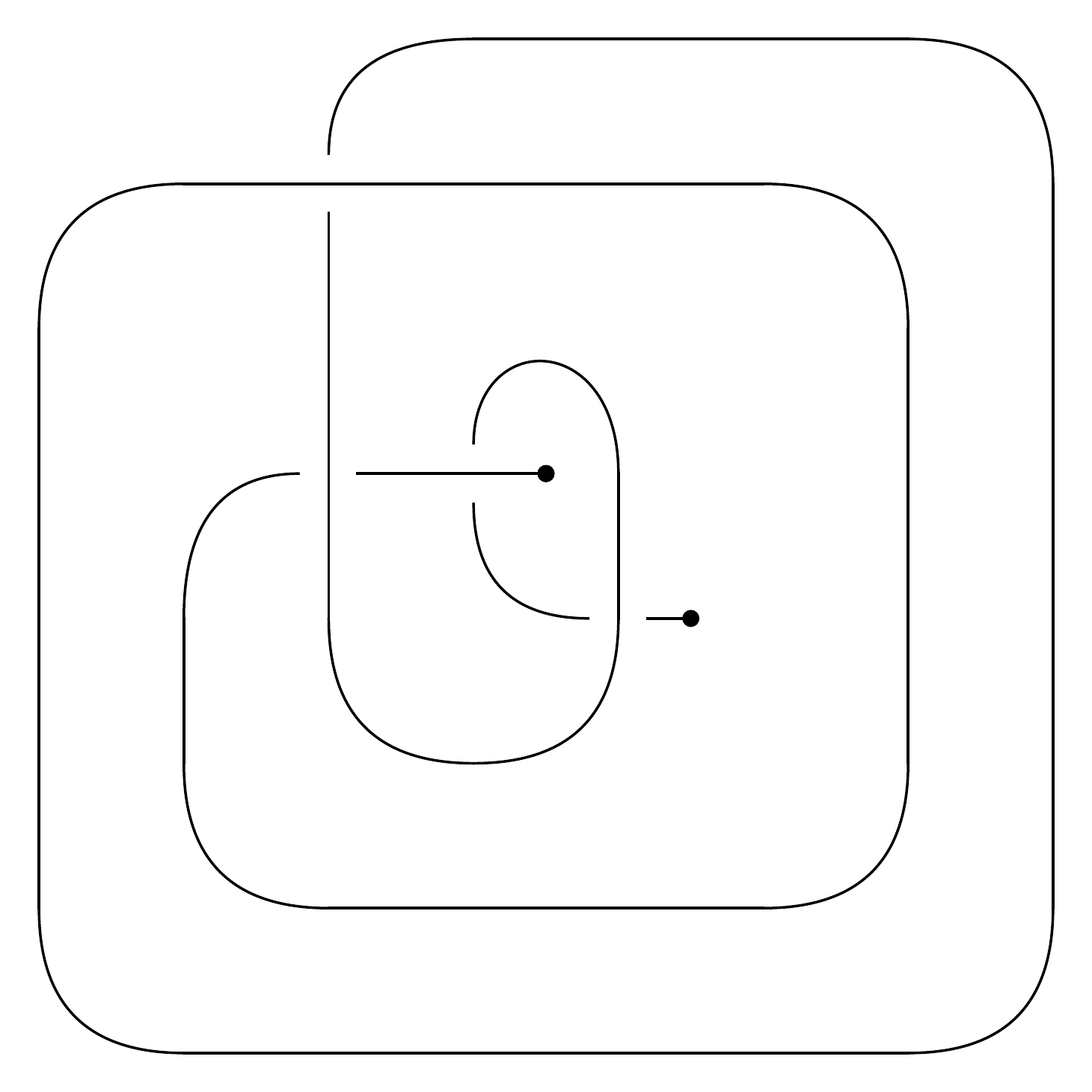}\\
\textcolor{black}{$4_{16}$}
\vspace{1cm}
\end{minipage}
\begin{minipage}[t]{.25\linewidth}
\centering
\includegraphics[width=0.9\textwidth,height=3.5cm,keepaspectratio]{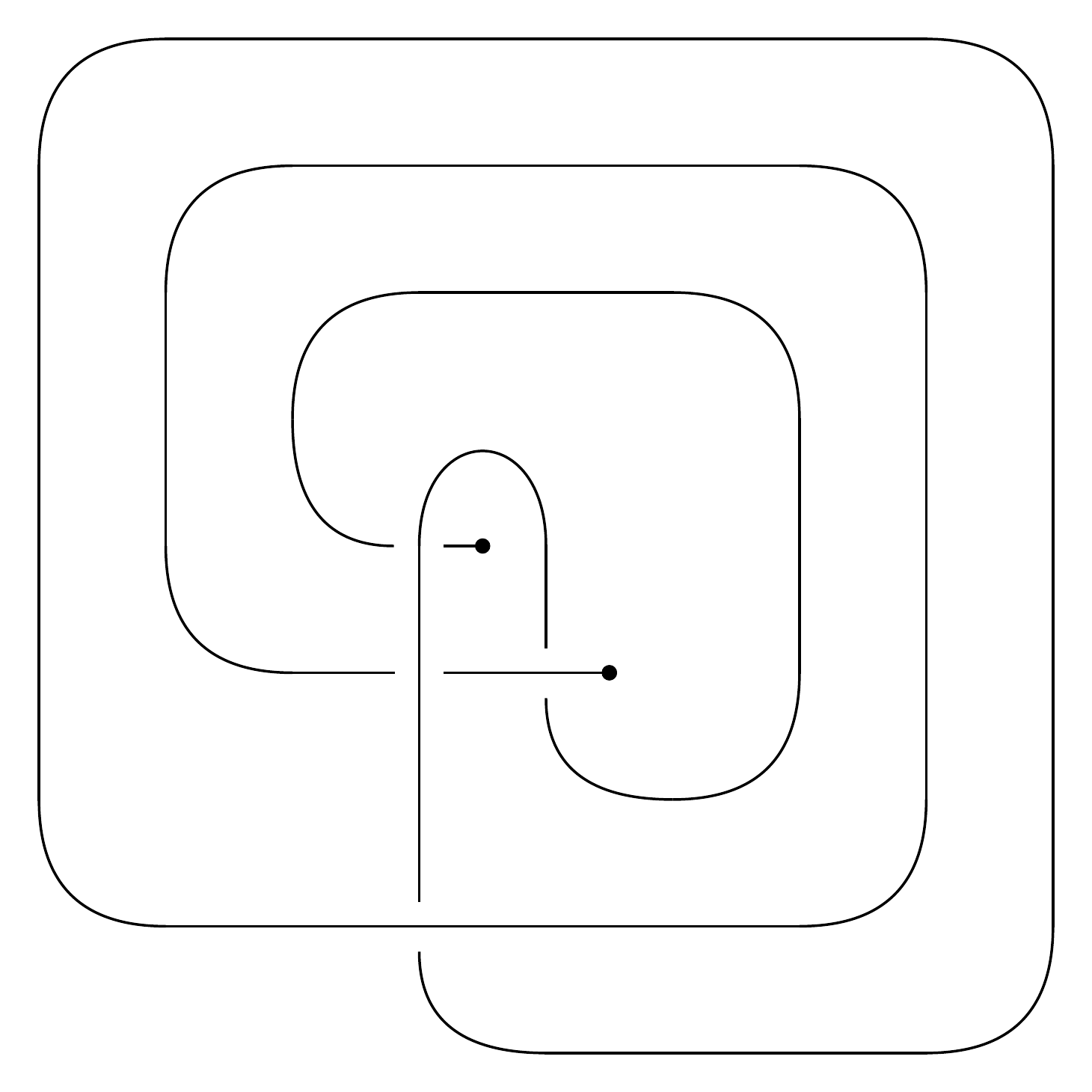}\\
\textcolor{black}{$4_{17}$}
\vspace{1cm}
\end{minipage}
\begin{minipage}[t]{.25\linewidth}
\centering
\includegraphics[width=0.9\textwidth,height=3.5cm,keepaspectratio]{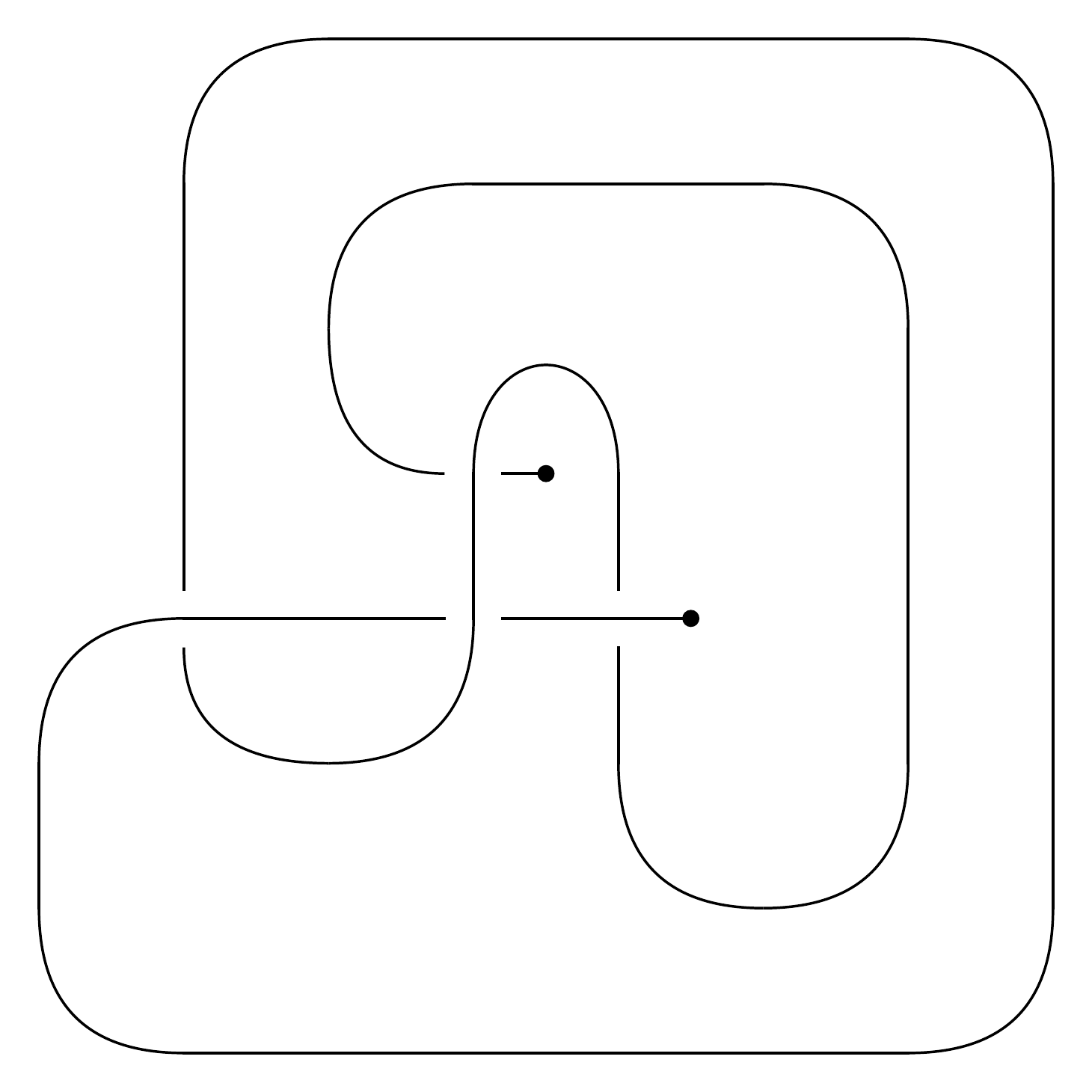}\\
\textcolor{black}{$4_{18}$}
\vspace{1cm}
\end{minipage}
\begin{minipage}[t]{.25\linewidth}
\centering
\includegraphics[width=0.9\textwidth,height=3.5cm,keepaspectratio]{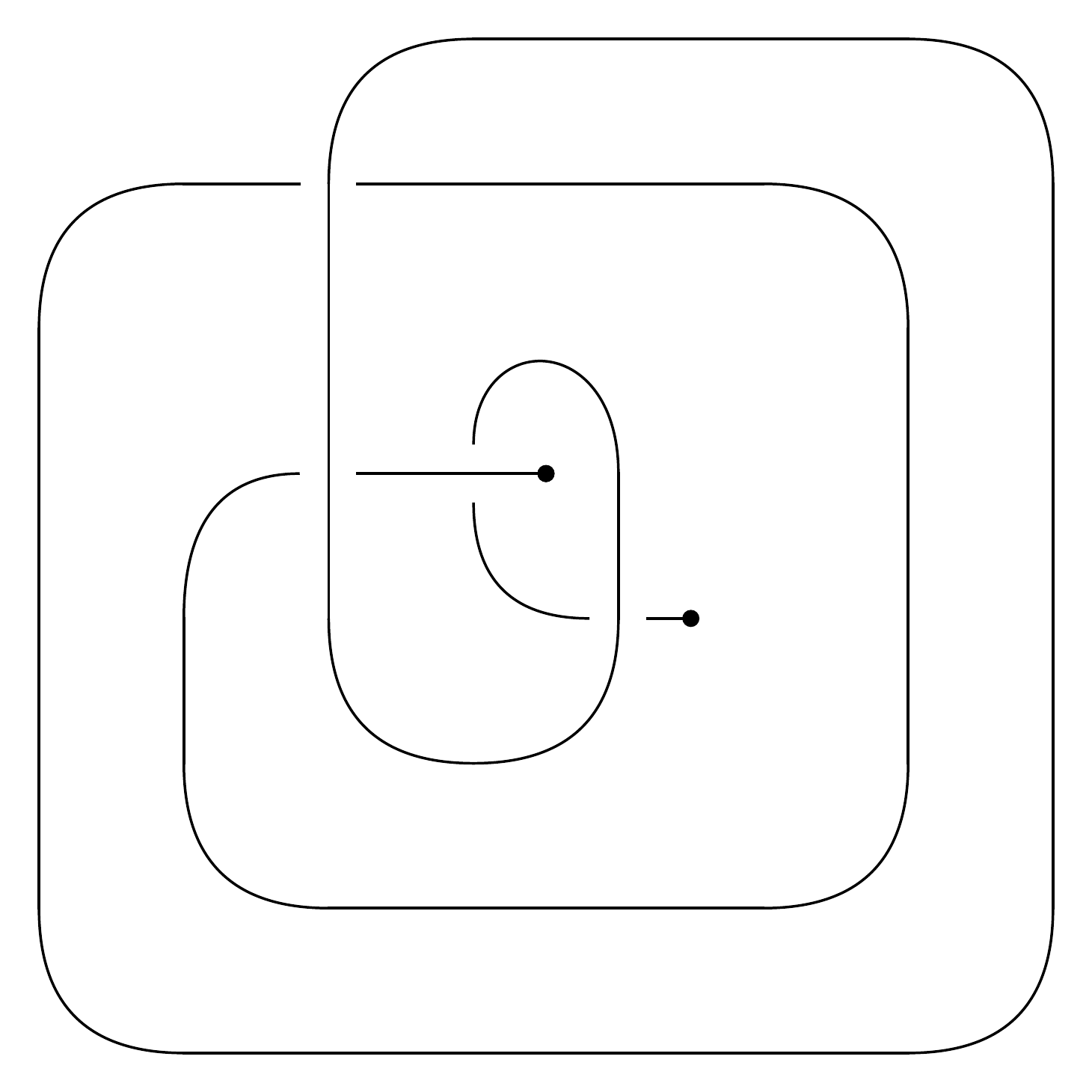}\\
\textcolor{black}{$4_{19}$}
\vspace{1cm}
\end{minipage}
\begin{minipage}[t]{.25\linewidth}
\centering
\includegraphics[width=0.9\textwidth,height=3.5cm,keepaspectratio]{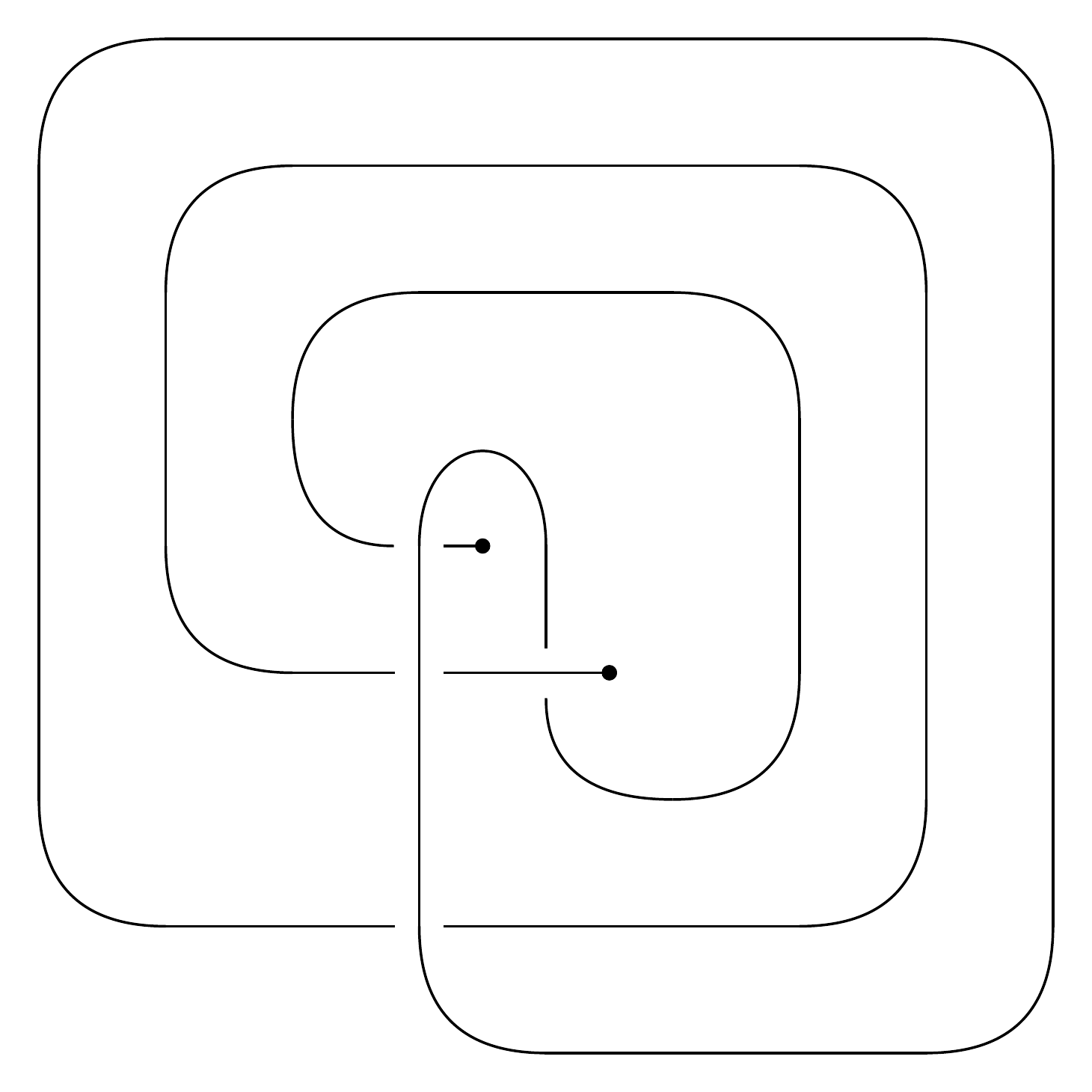}\\
\textcolor{black}{$4_{20}$}
\vspace{1cm}
\end{minipage}
\begin{minipage}[t]{.25\linewidth}
\centering
\includegraphics[width=0.9\textwidth,height=3.5cm,keepaspectratio]{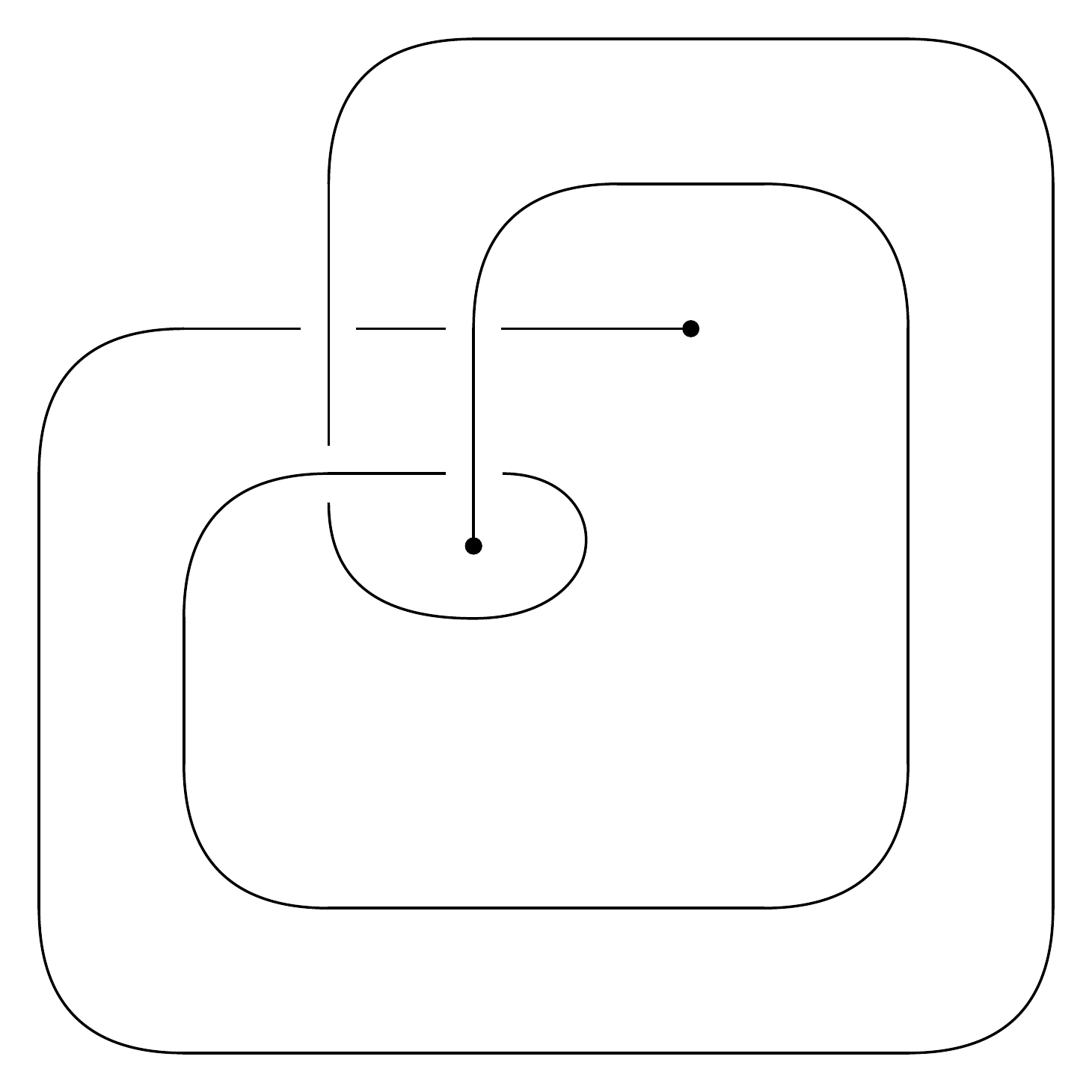}\\
\textcolor{black}{$4_{21}$}
\vspace{1cm}
\end{minipage}
\begin{minipage}[t]{.25\linewidth}
\centering
\includegraphics[width=0.9\textwidth,height=3.5cm,keepaspectratio]{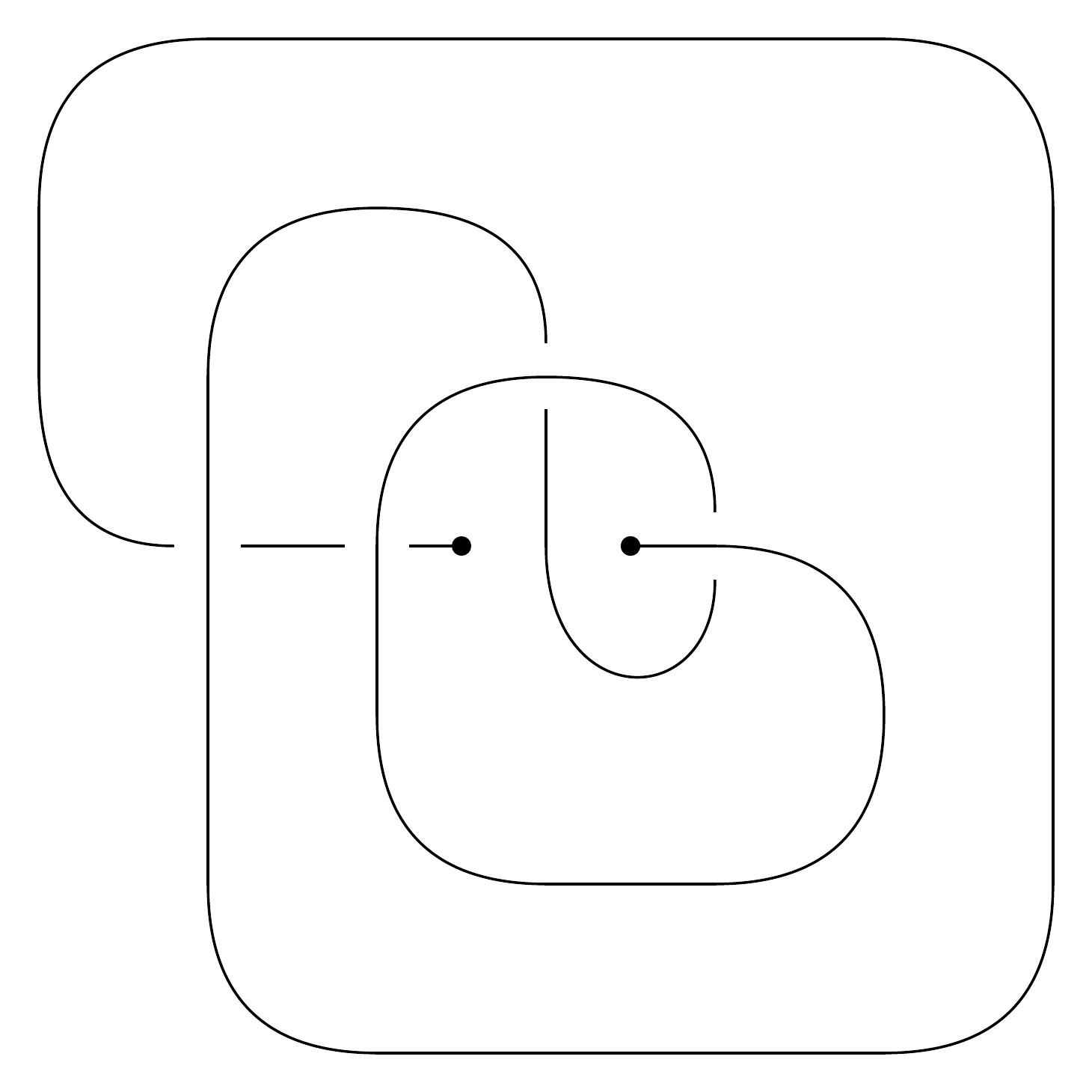}\\
\textcolor{black}{$4_{22}$}
\vspace{1cm}
\end{minipage}
\begin{minipage}[t]{.25\linewidth}
\centering
\includegraphics[width=0.9\textwidth,height=3.5cm,keepaspectratio]{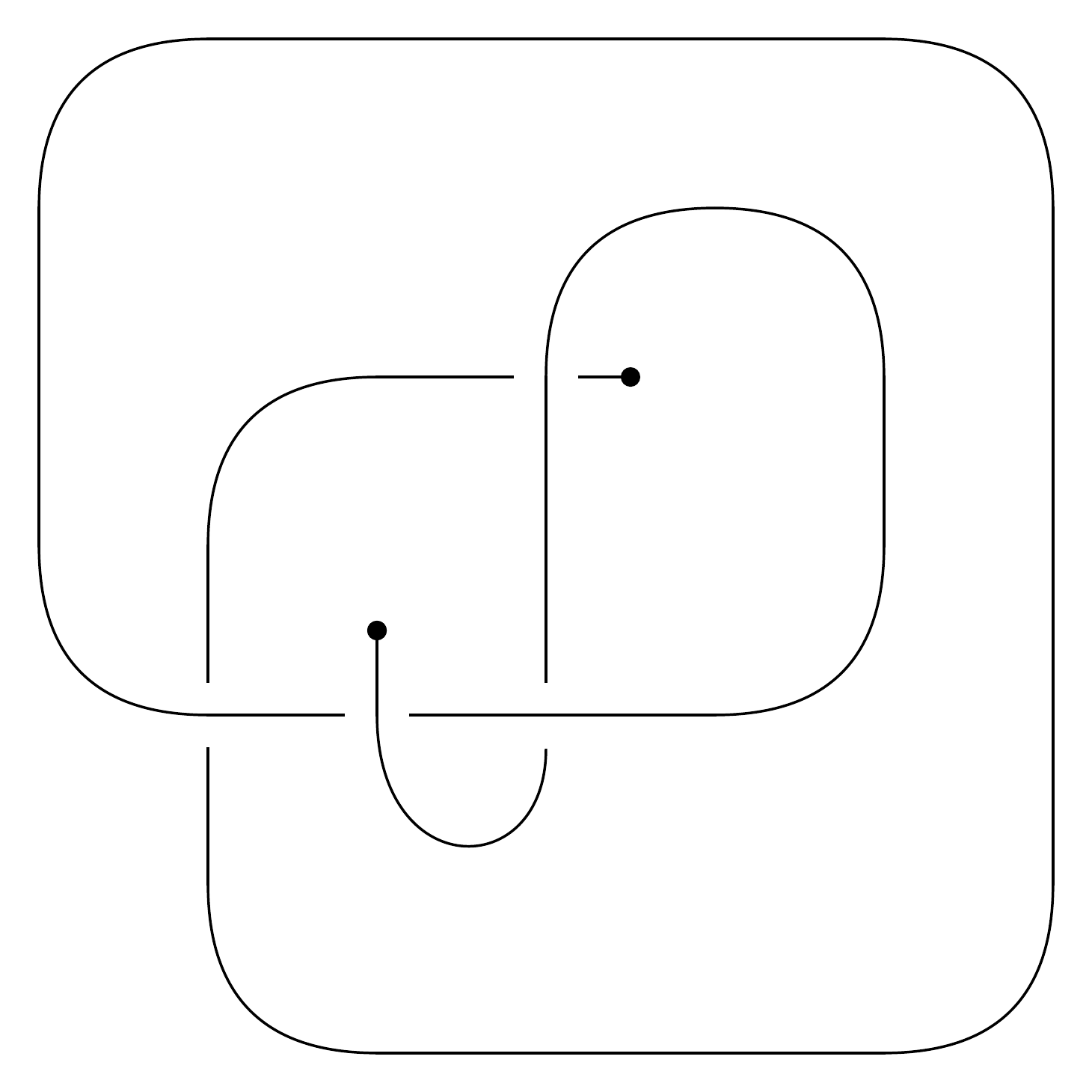}\\
\textcolor{black}{$4_{23}$}
\vspace{1cm}
\end{minipage}
\begin{minipage}[t]{.25\linewidth}
\centering
\includegraphics[width=0.9\textwidth,height=3.5cm,keepaspectratio]{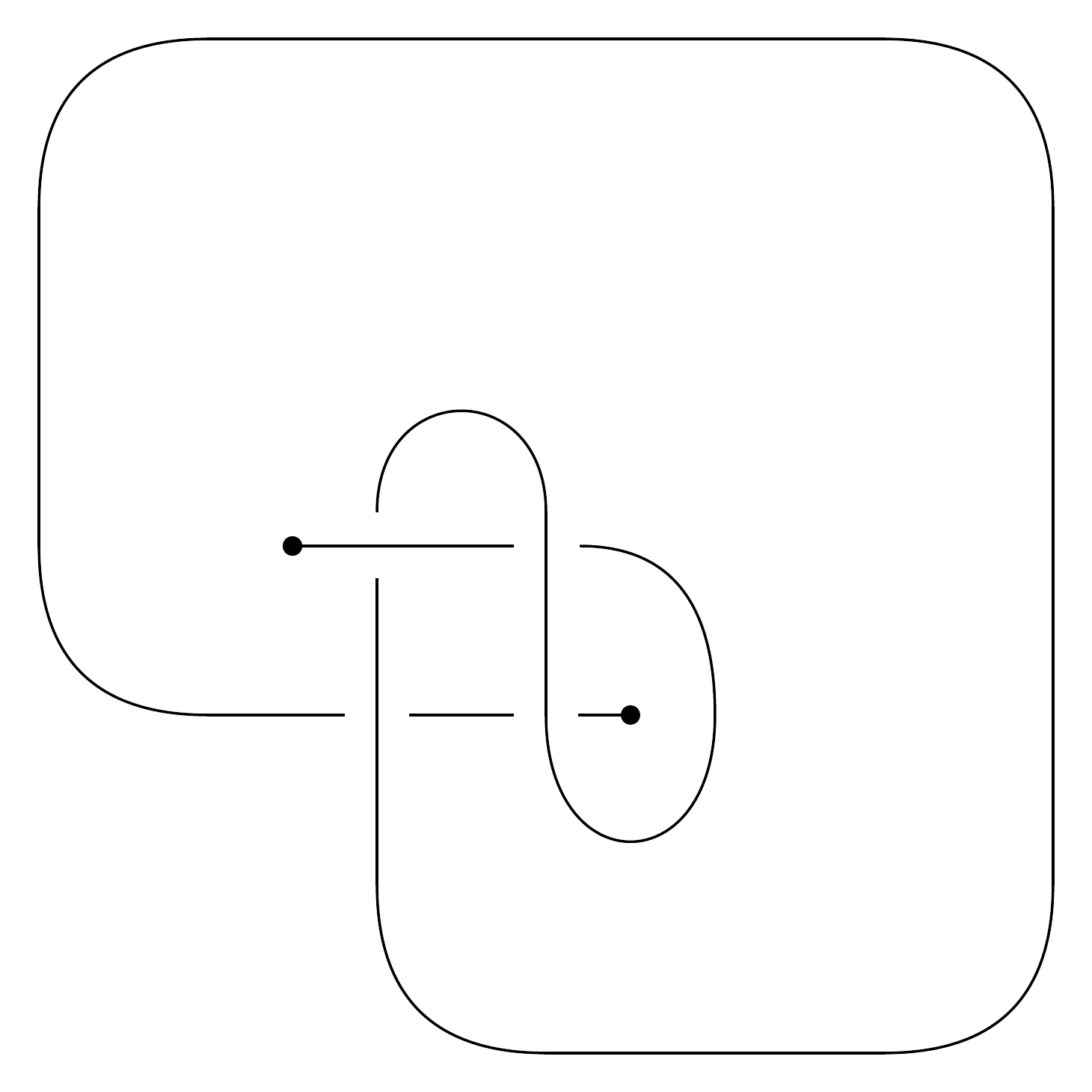}\\
\textcolor{black}{$4_{24}$}
\vspace{1cm}
\end{minipage}
\begin{minipage}[t]{.25\linewidth}
\centering
\includegraphics[width=0.9\textwidth,height=3.5cm,keepaspectratio]{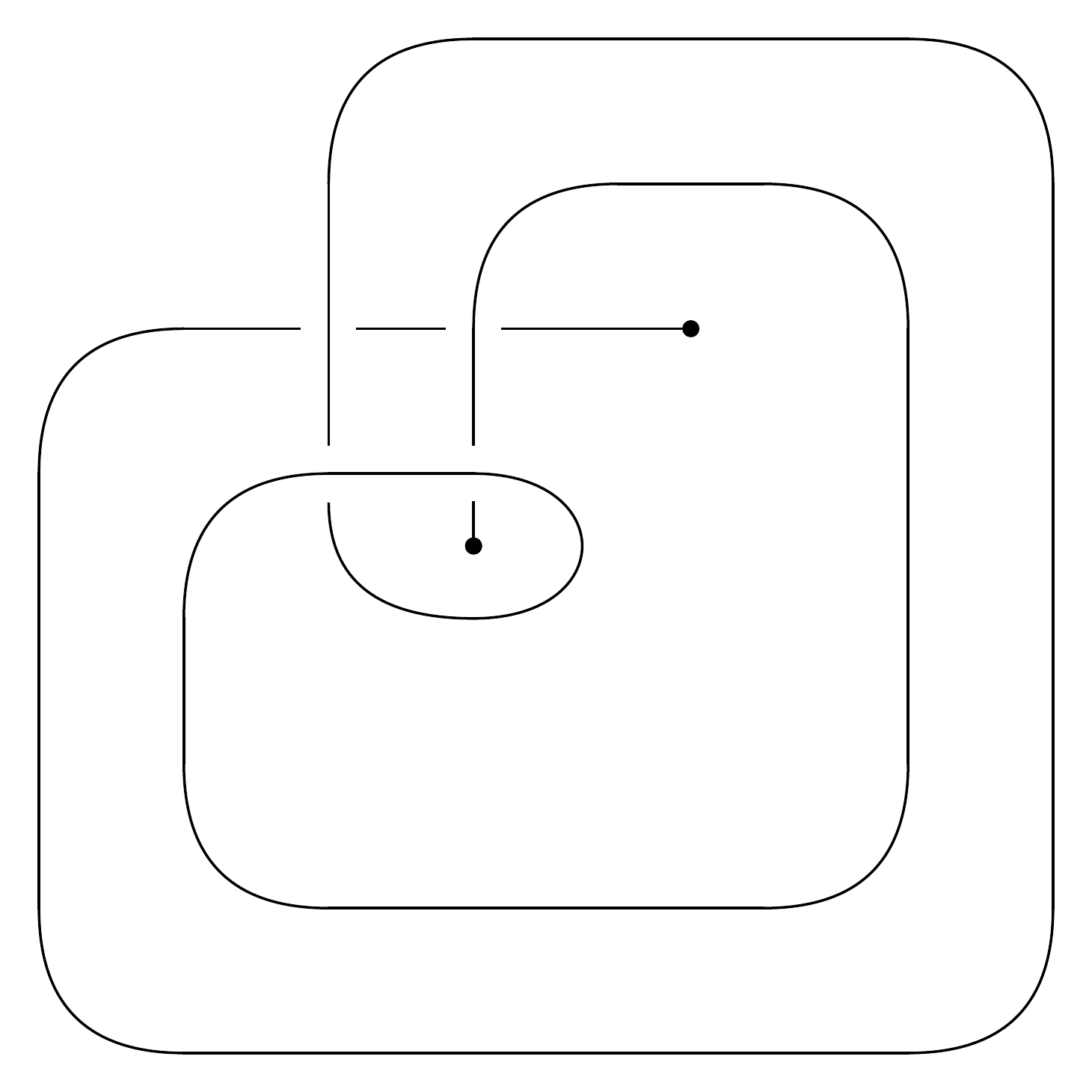}\\
\textcolor{black}{$4_{25}$}
\vspace{1cm}
\end{minipage}
\begin{minipage}[t]{.25\linewidth}
\centering
\includegraphics[width=0.9\textwidth,height=3.5cm,keepaspectratio]{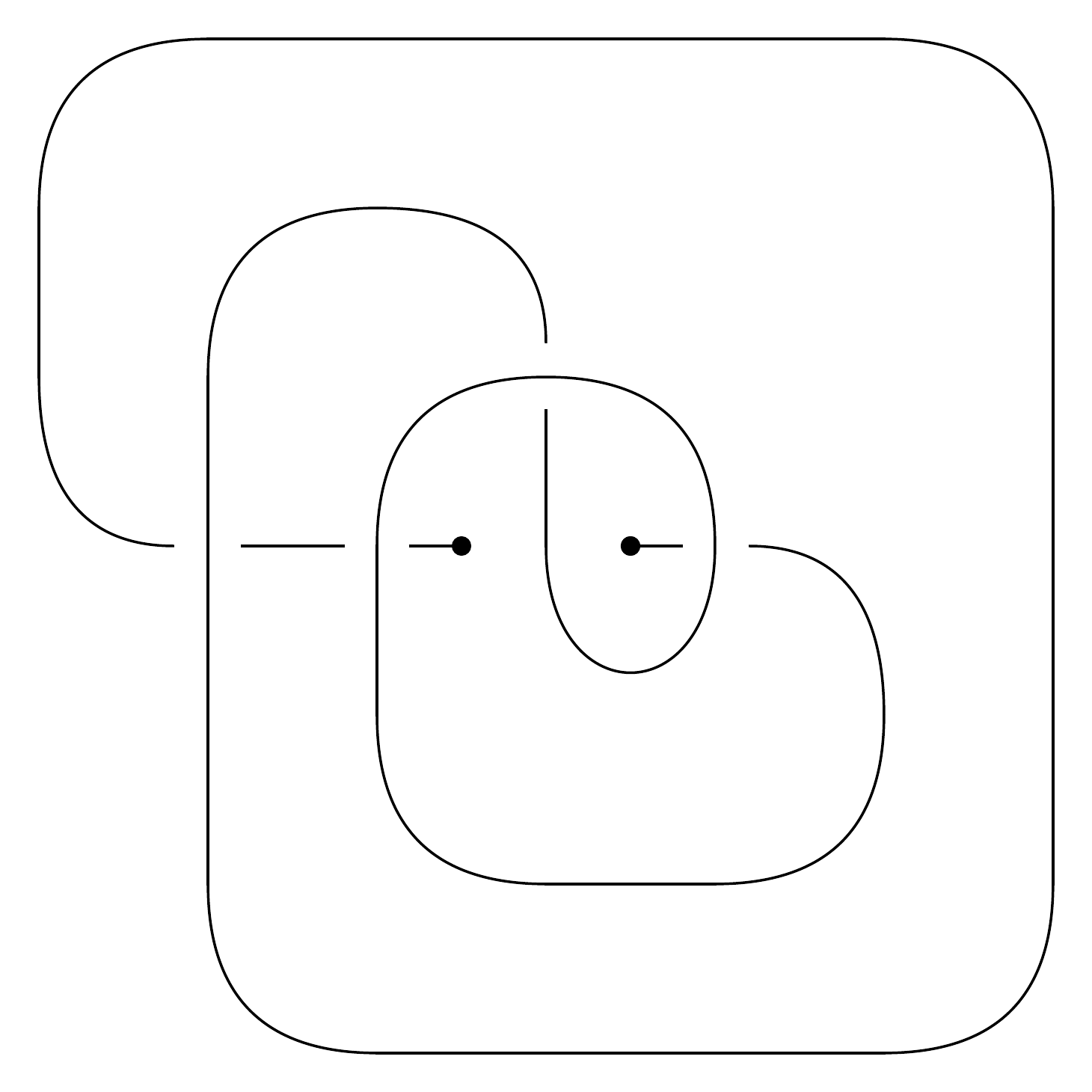}\\
\textcolor{black}{$4_{26}$}
\vspace{1cm}
\end{minipage}
\begin{minipage}[t]{.25\linewidth}
\centering
\includegraphics[width=0.9\textwidth,height=3.5cm,keepaspectratio]{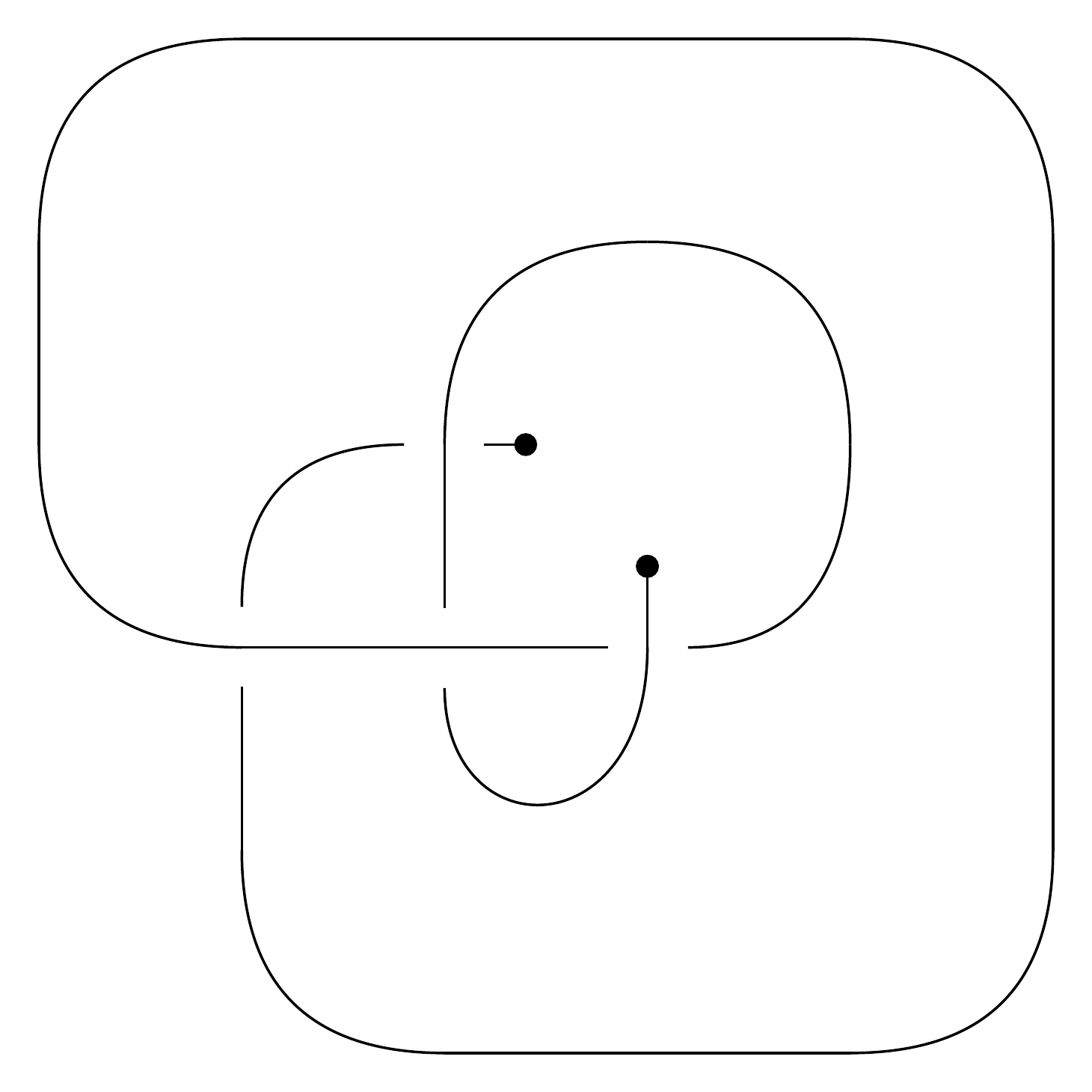}\\
\textcolor{black}{$4_{27}$}
\vspace{1cm}
\end{minipage}
\begin{minipage}[t]{.25\linewidth}
\centering
\includegraphics[width=0.9\textwidth,height=3.5cm,keepaspectratio]{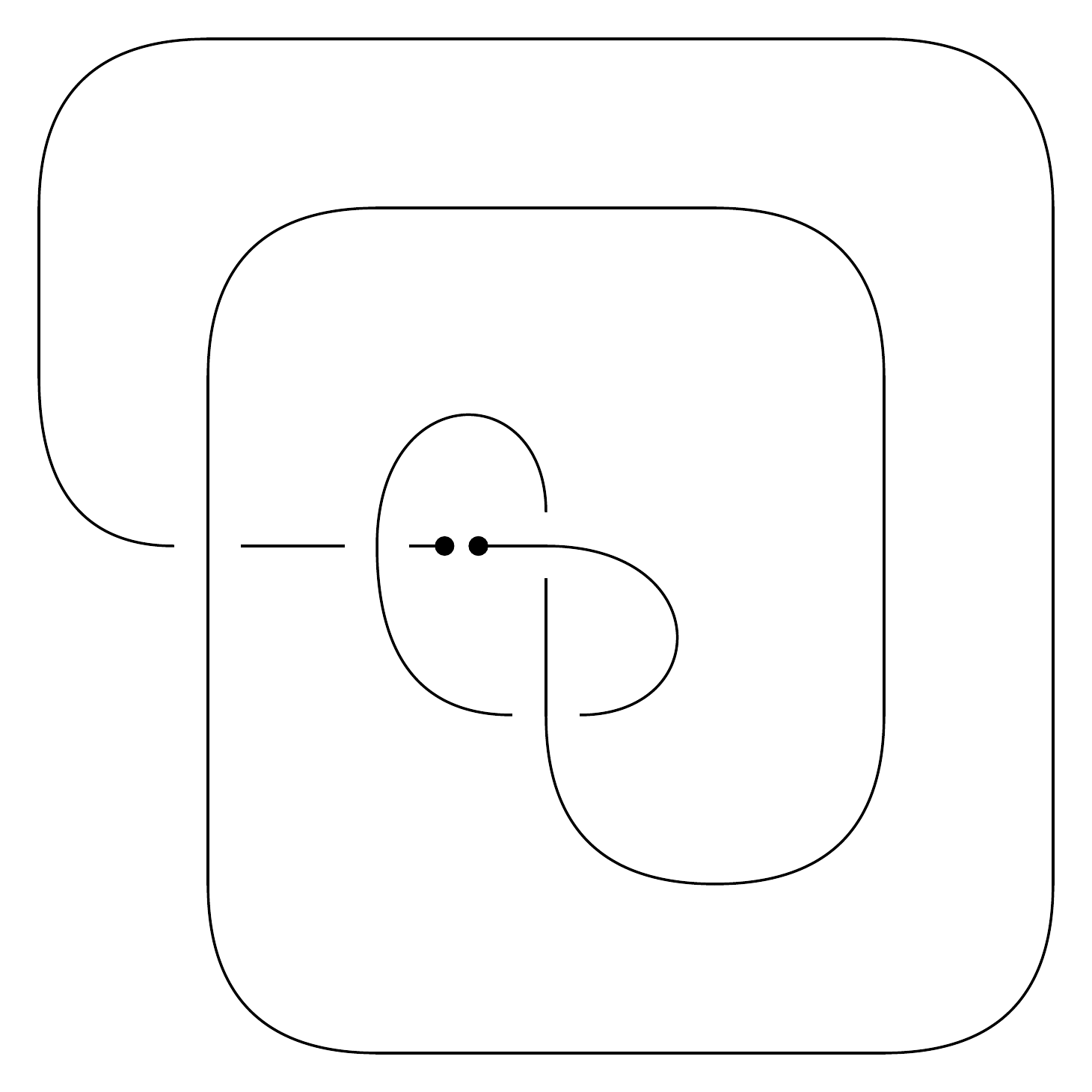}\\
\textcolor{black}{$4_{28}$}
\vspace{1cm}
\end{minipage}
\begin{minipage}[t]{.25\linewidth}
\centering
\includegraphics[width=0.9\textwidth,height=3.5cm,keepaspectratio]{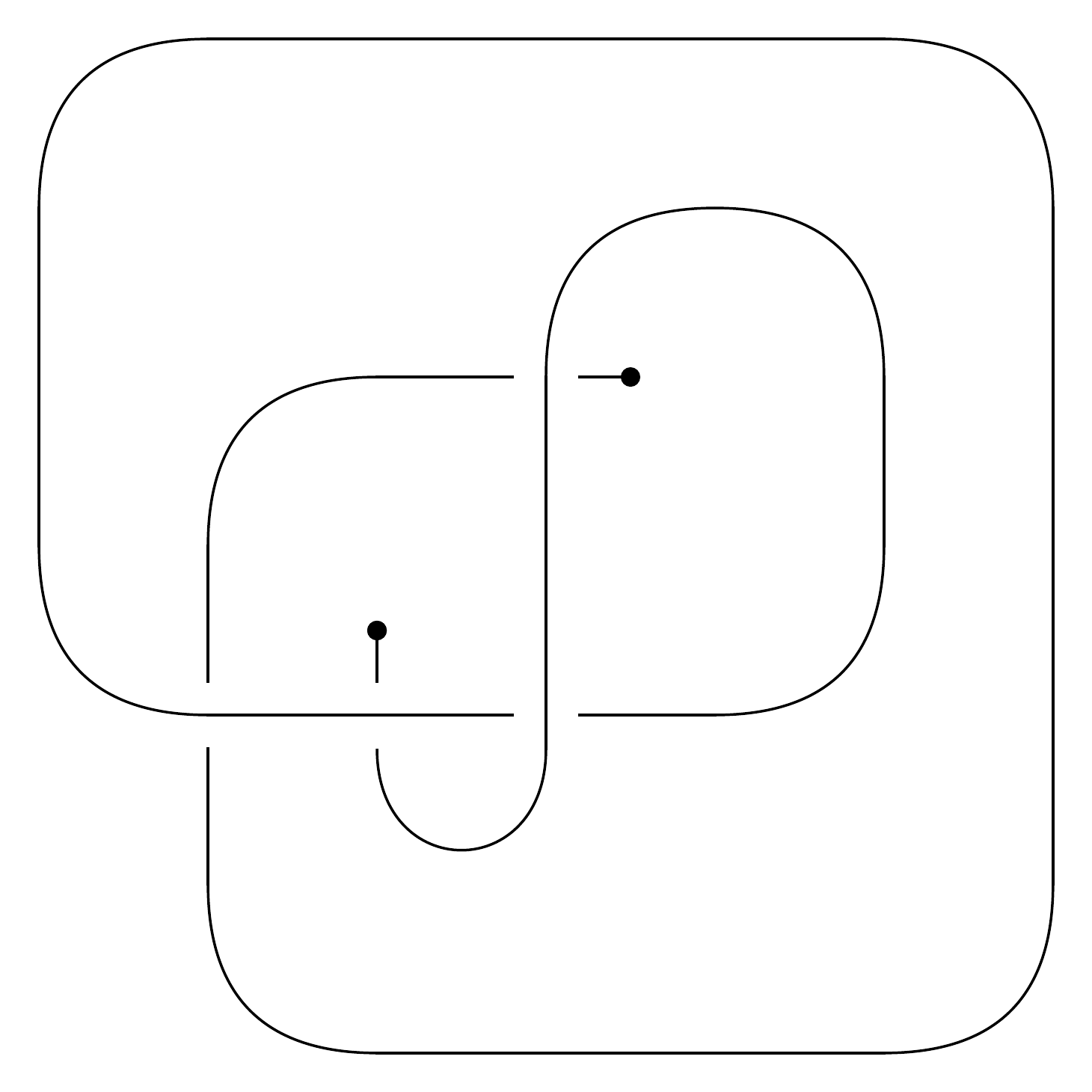}\\
\textcolor{black}{$4_{29}$}
\vspace{1cm}
\end{minipage}
\begin{minipage}[t]{.25\linewidth}
\centering
\includegraphics[width=0.9\textwidth,height=3.5cm,keepaspectratio]{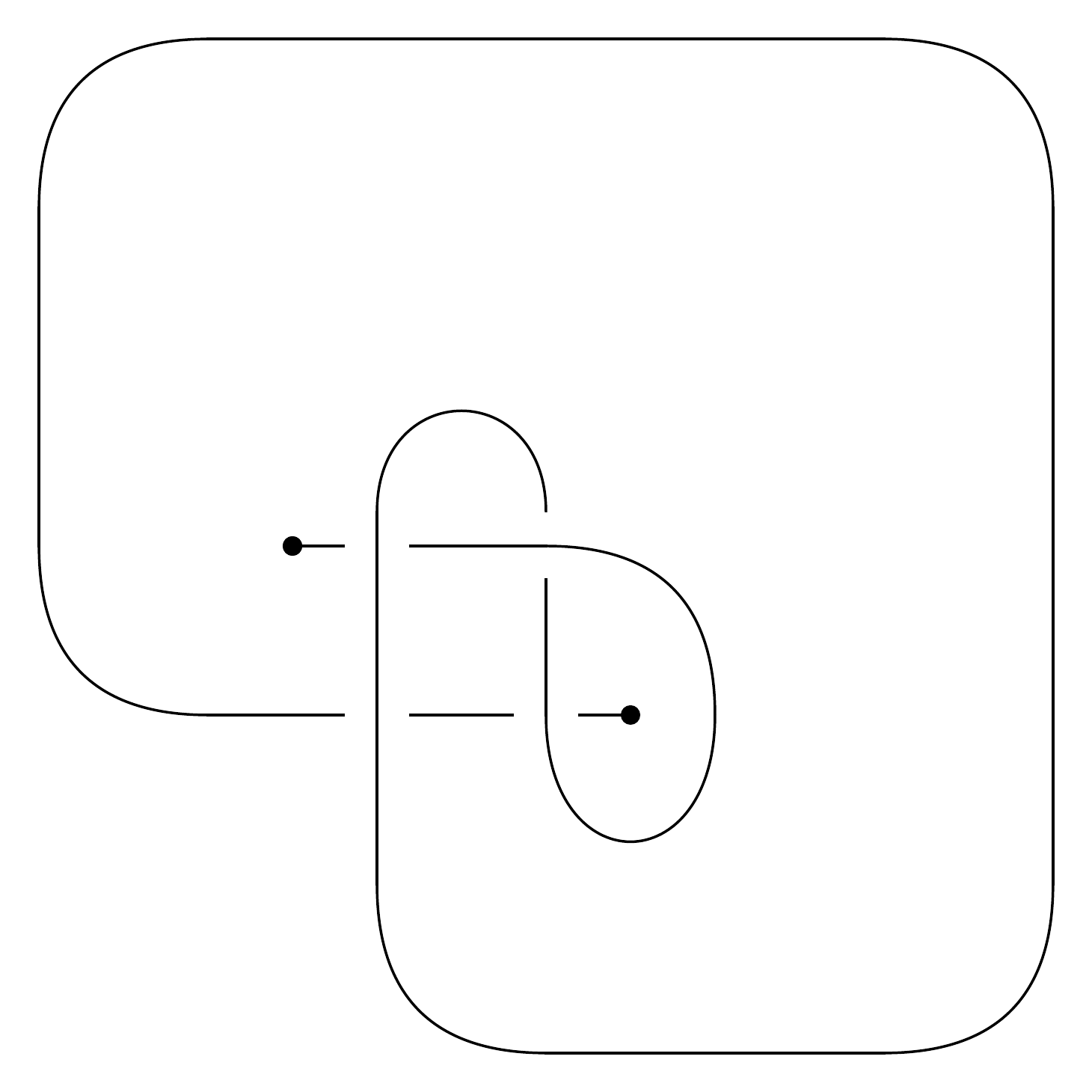}\\
\textcolor{black}{$4_{30}$}
\vspace{1cm}
\end{minipage}
\begin{minipage}[t]{.25\linewidth}
\centering
\includegraphics[width=0.9\textwidth,height=3.5cm,keepaspectratio]{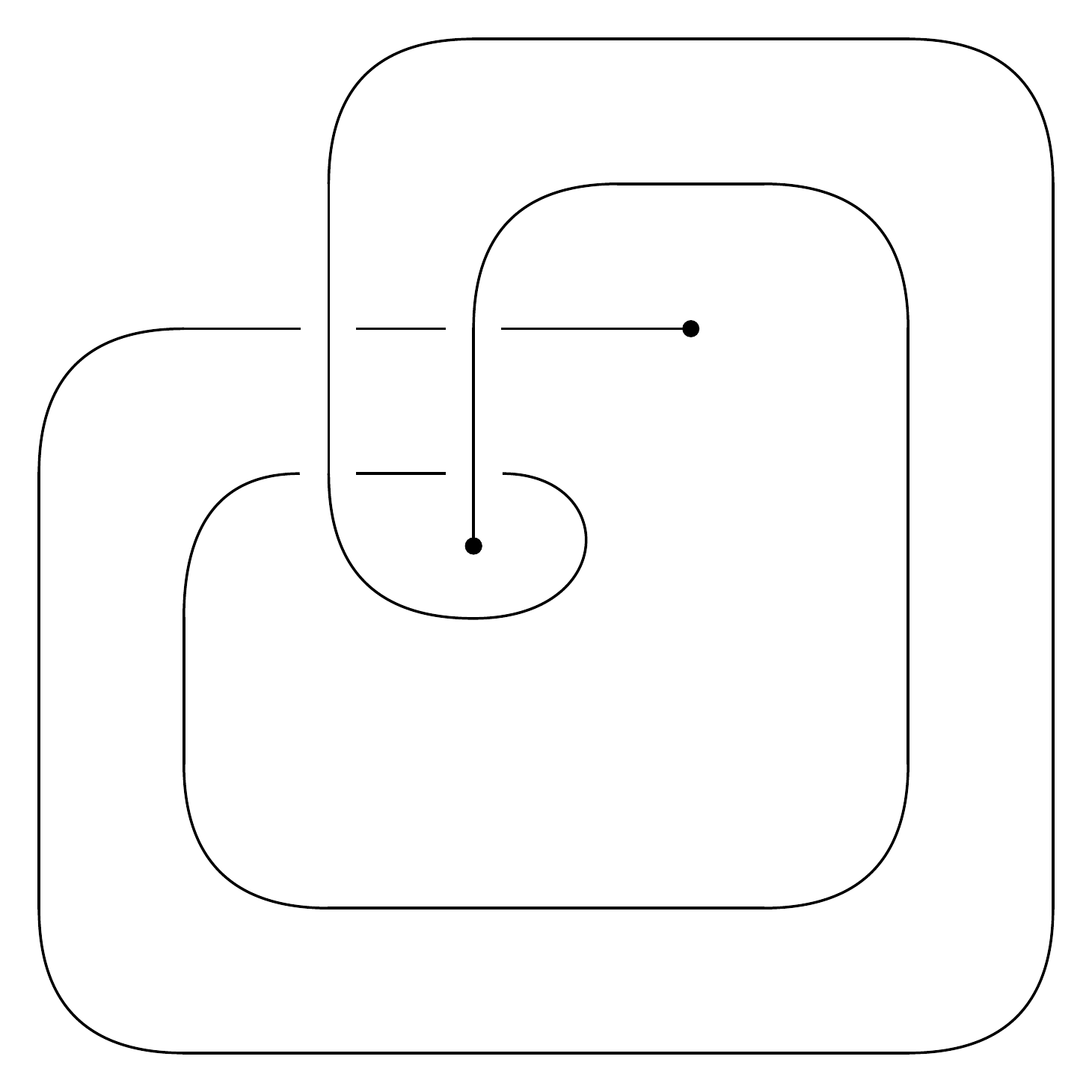}\\
\textcolor{black}{$4_{31}$}
\vspace{1cm}
\end{minipage}
\begin{minipage}[t]{.25\linewidth}
\centering
\includegraphics[width=0.9\textwidth,height=3.5cm,keepaspectratio]{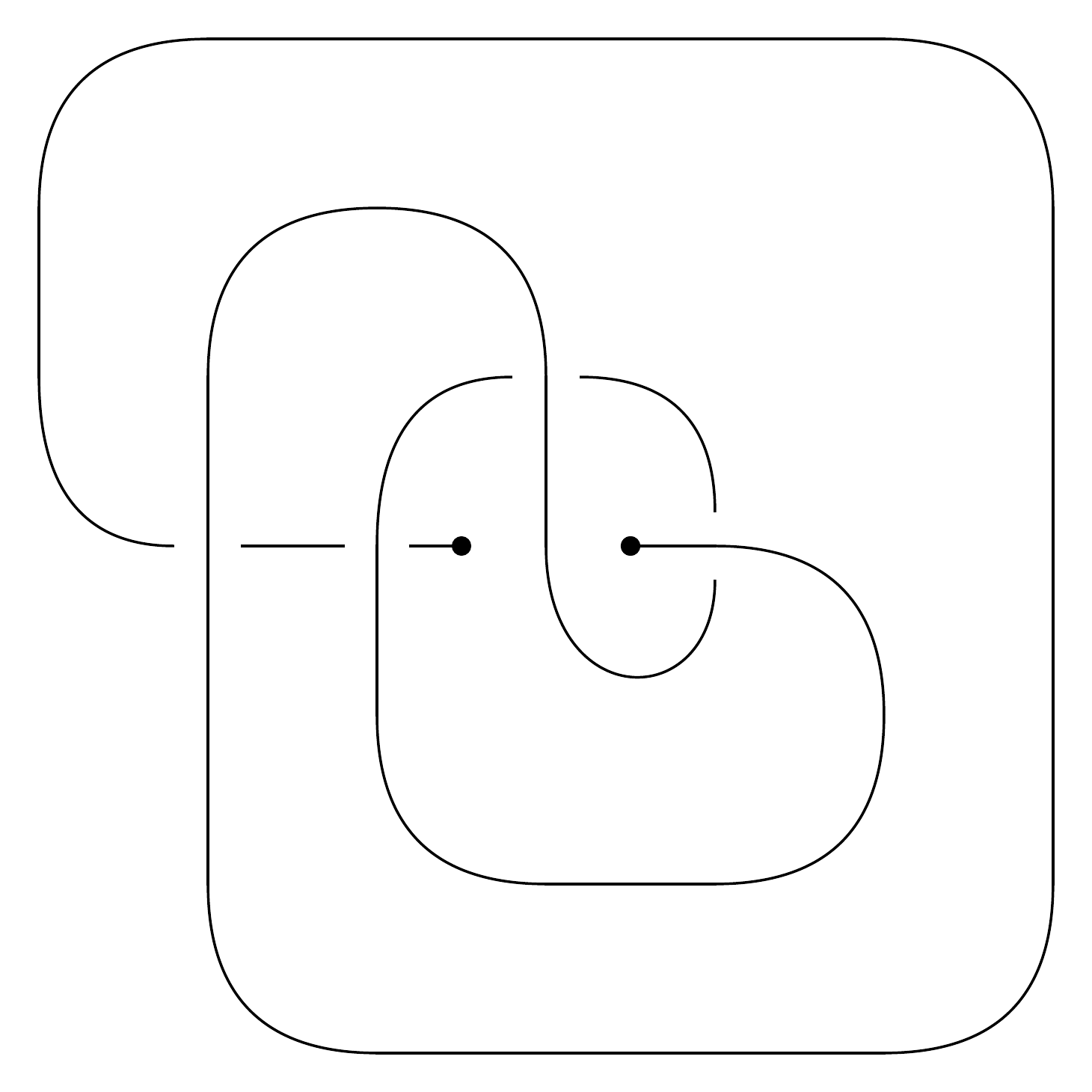}\\
\textcolor{black}{$4_{32}$}
\vspace{1cm}
\end{minipage}
\begin{minipage}[t]{.25\linewidth}
\centering
\includegraphics[width=0.9\textwidth,height=3.5cm,keepaspectratio]{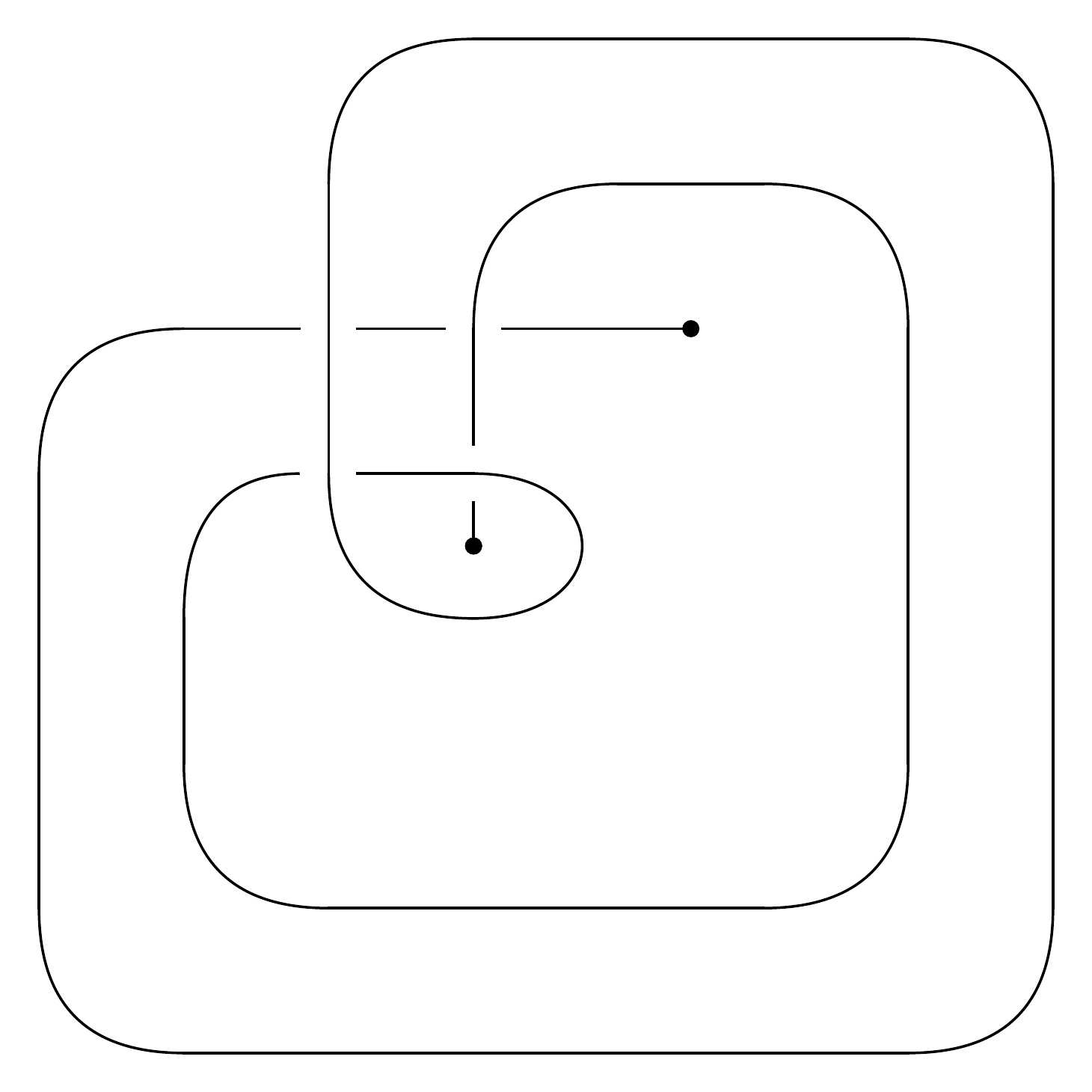}\\
\textcolor{black}{$4_{33}$}
\vspace{1cm}
\end{minipage}
\begin{minipage}[t]{.25\linewidth}
\centering
\includegraphics[width=0.9\textwidth,height=3.5cm,keepaspectratio]{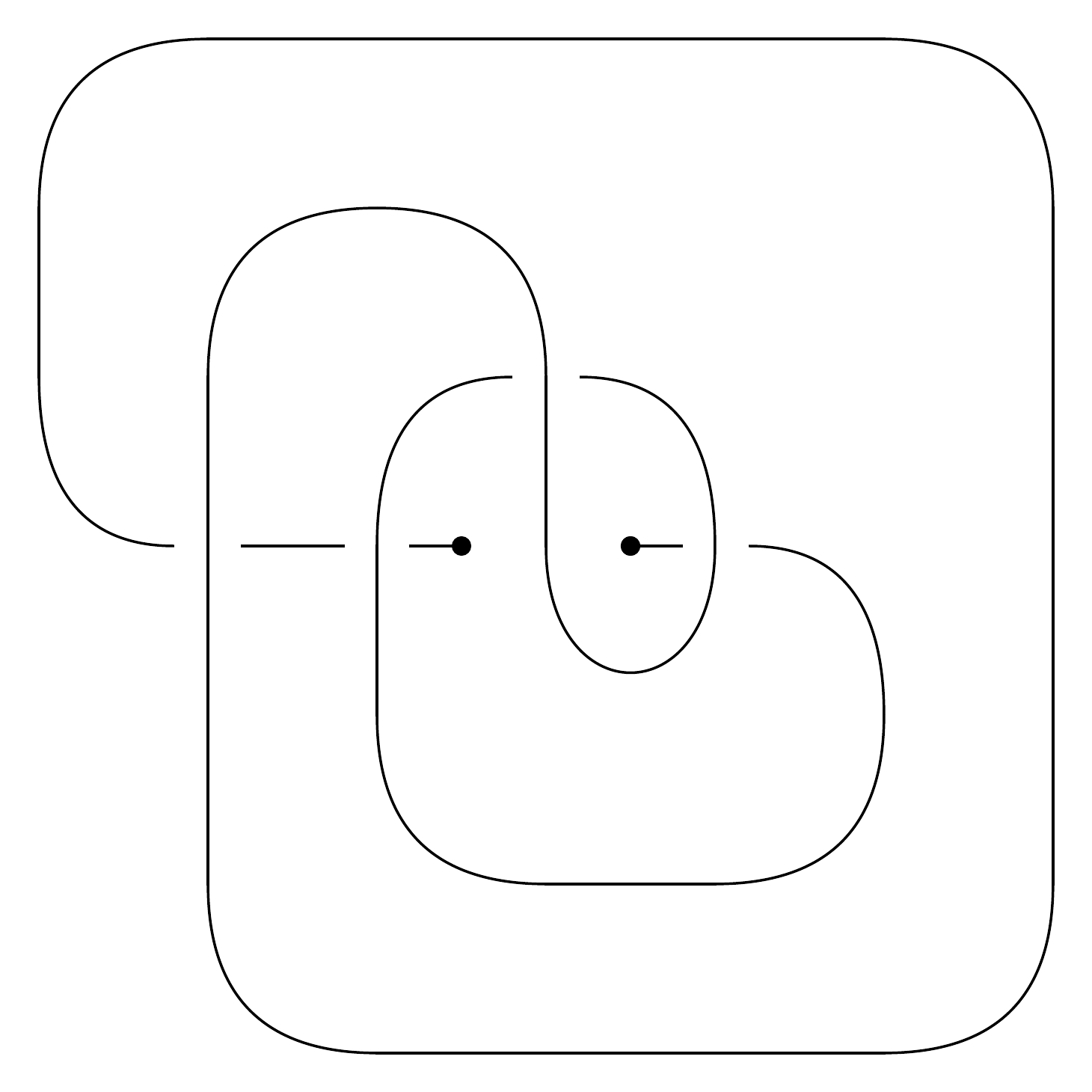}\\
\textcolor{black}{$4_{34}$}
\vspace{1cm}
\end{minipage}
\begin{minipage}[t]{.25\linewidth}
\centering
\includegraphics[width=0.9\textwidth,height=3.5cm,keepaspectratio]{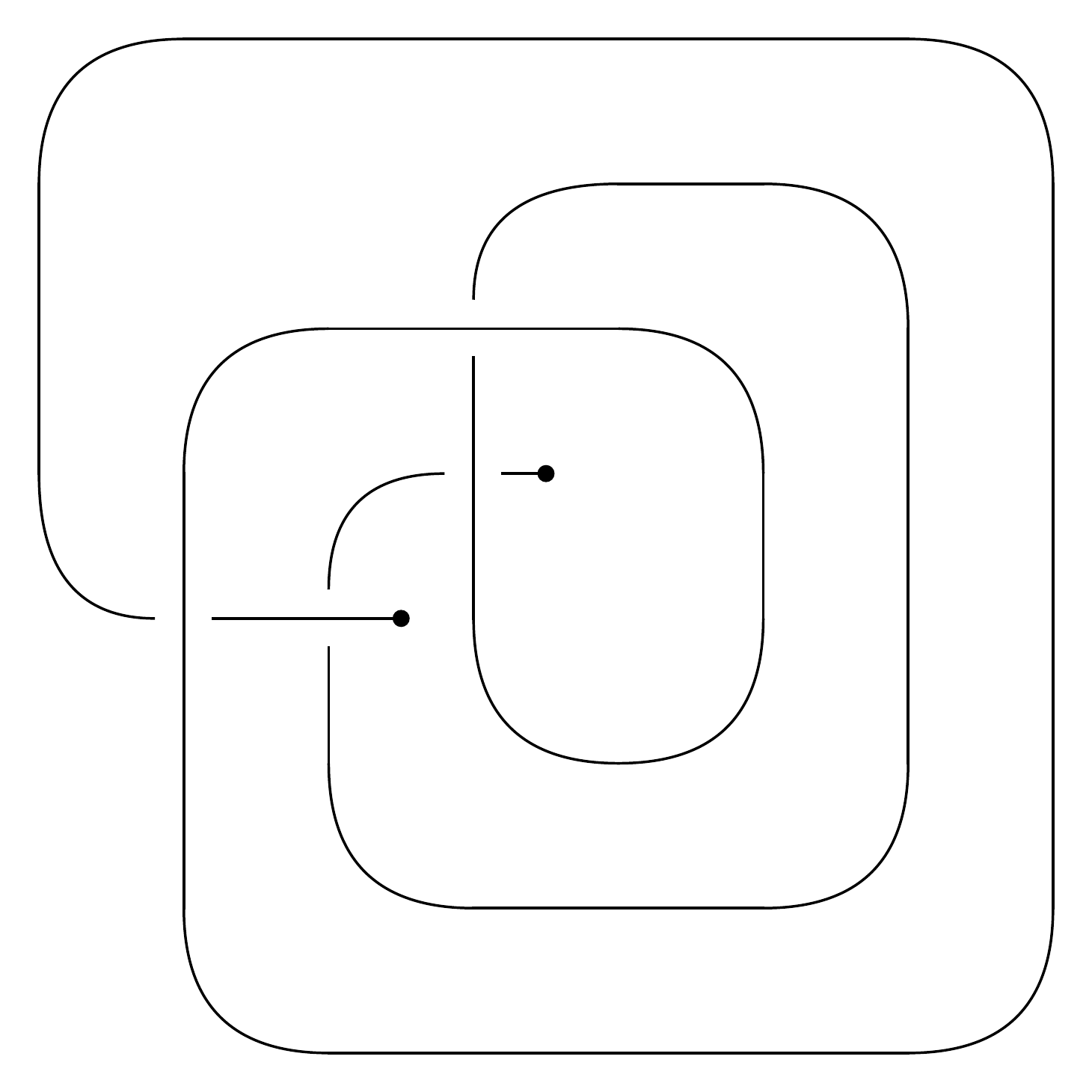}\\
\textcolor{black}{$4_{35}$}
\vspace{1cm}
\end{minipage}
\begin{minipage}[t]{.25\linewidth}
\centering
\includegraphics[width=0.9\textwidth,height=3.5cm,keepaspectratio]{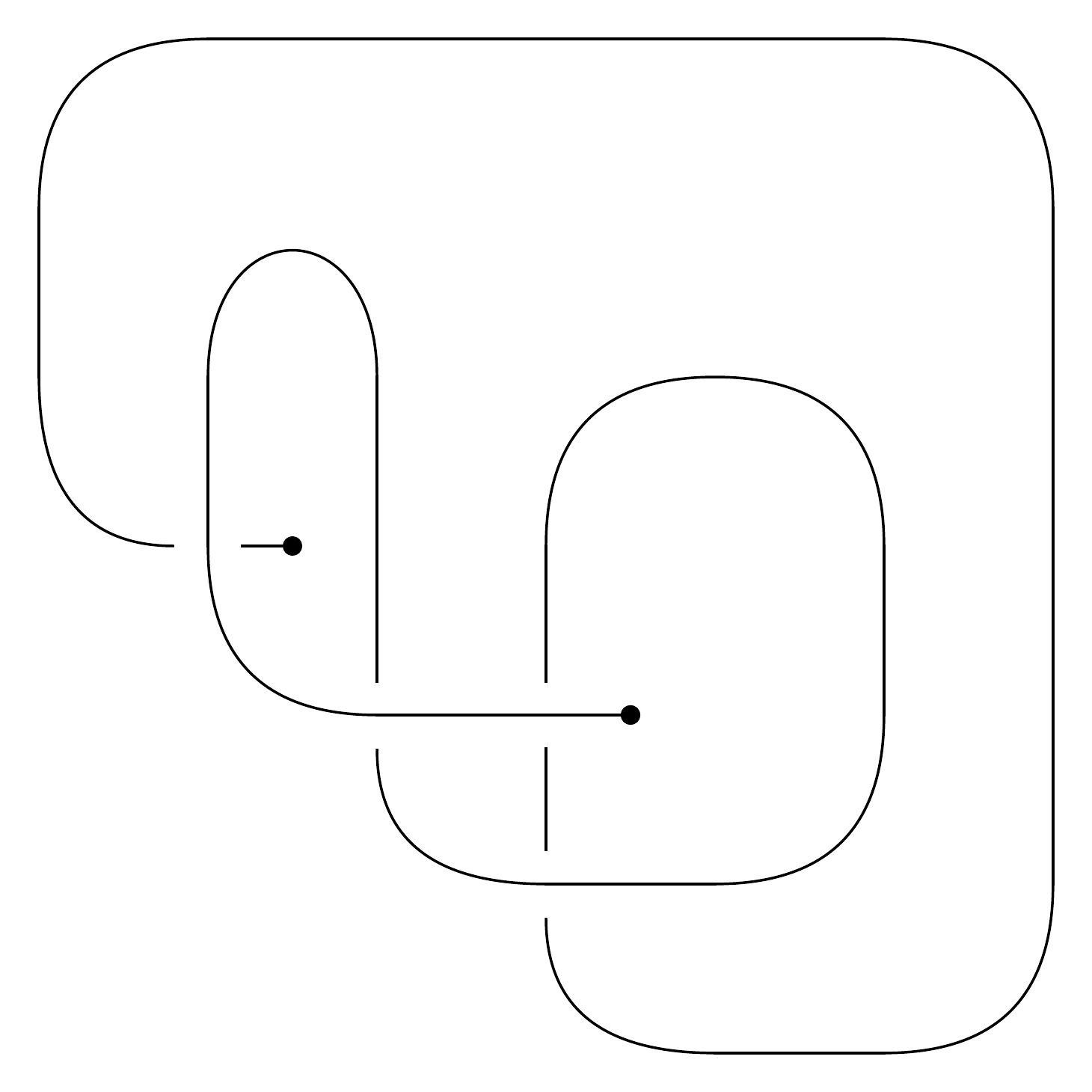}\\
\textcolor{black}{$4_{36}$}
\vspace{1cm}
\end{minipage}
\begin{minipage}[t]{.25\linewidth}
\centering
\includegraphics[width=0.9\textwidth,height=3.5cm,keepaspectratio]{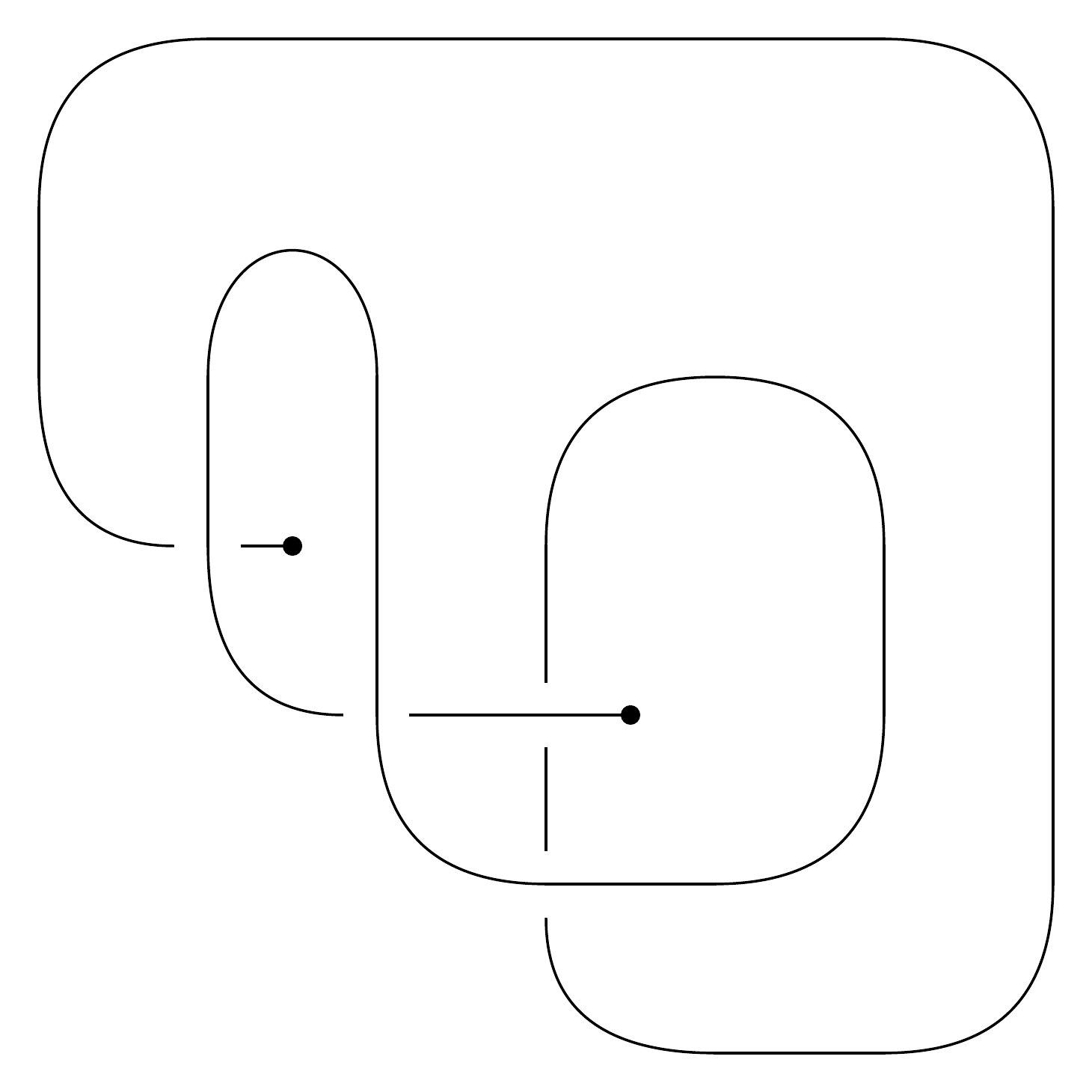}\\
\textcolor{black}{$4_{37}$}
\vspace{1cm}
\end{minipage}
\begin{minipage}[t]{.25\linewidth}
\centering
\includegraphics[width=0.9\textwidth,height=3.5cm,keepaspectratio]{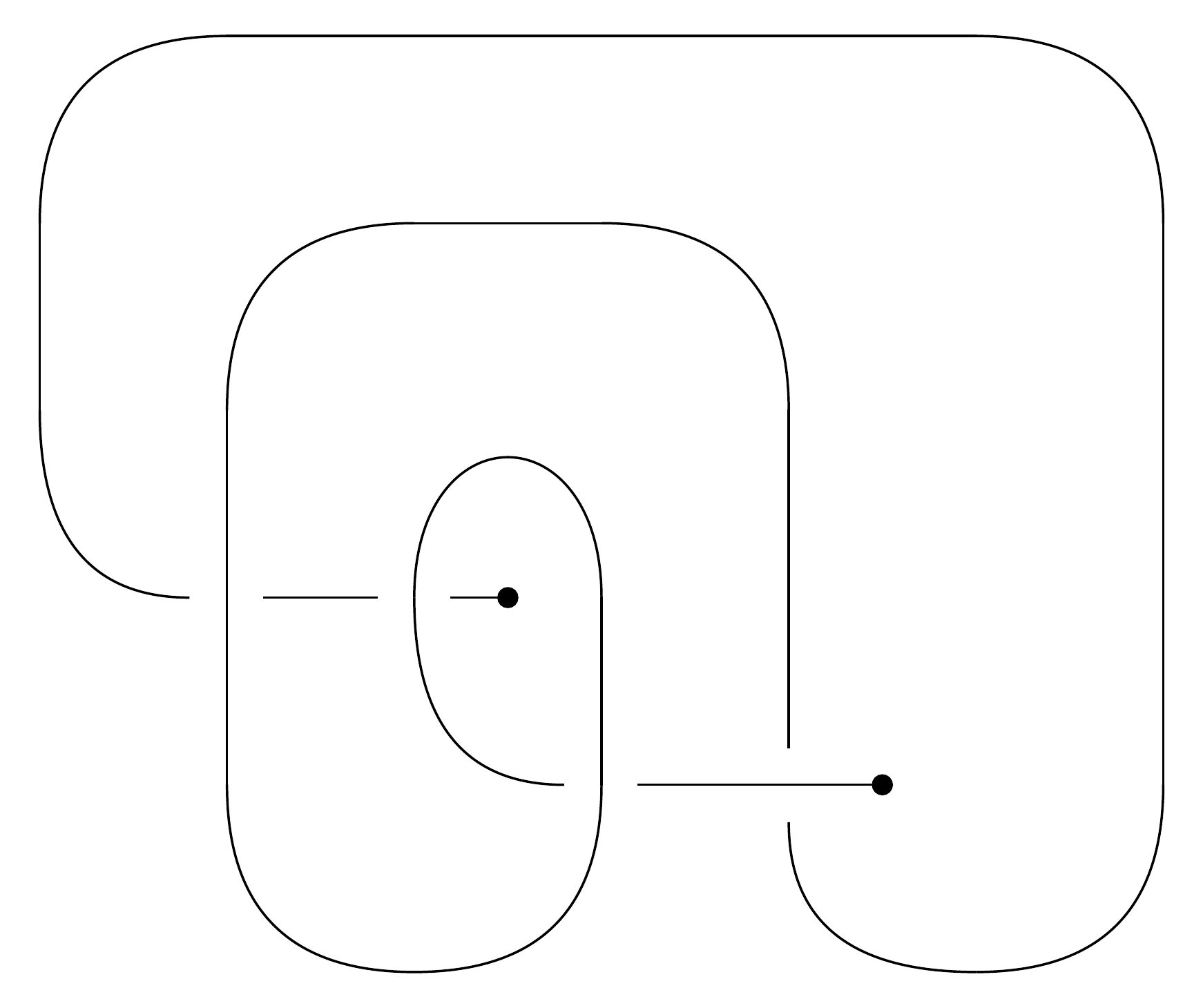}\\
\textcolor{black}{$4_{38}$}
\vspace{1cm}
\end{minipage}
\begin{minipage}[t]{.25\linewidth}
\centering
\includegraphics[width=0.9\textwidth,height=3.5cm,keepaspectratio]{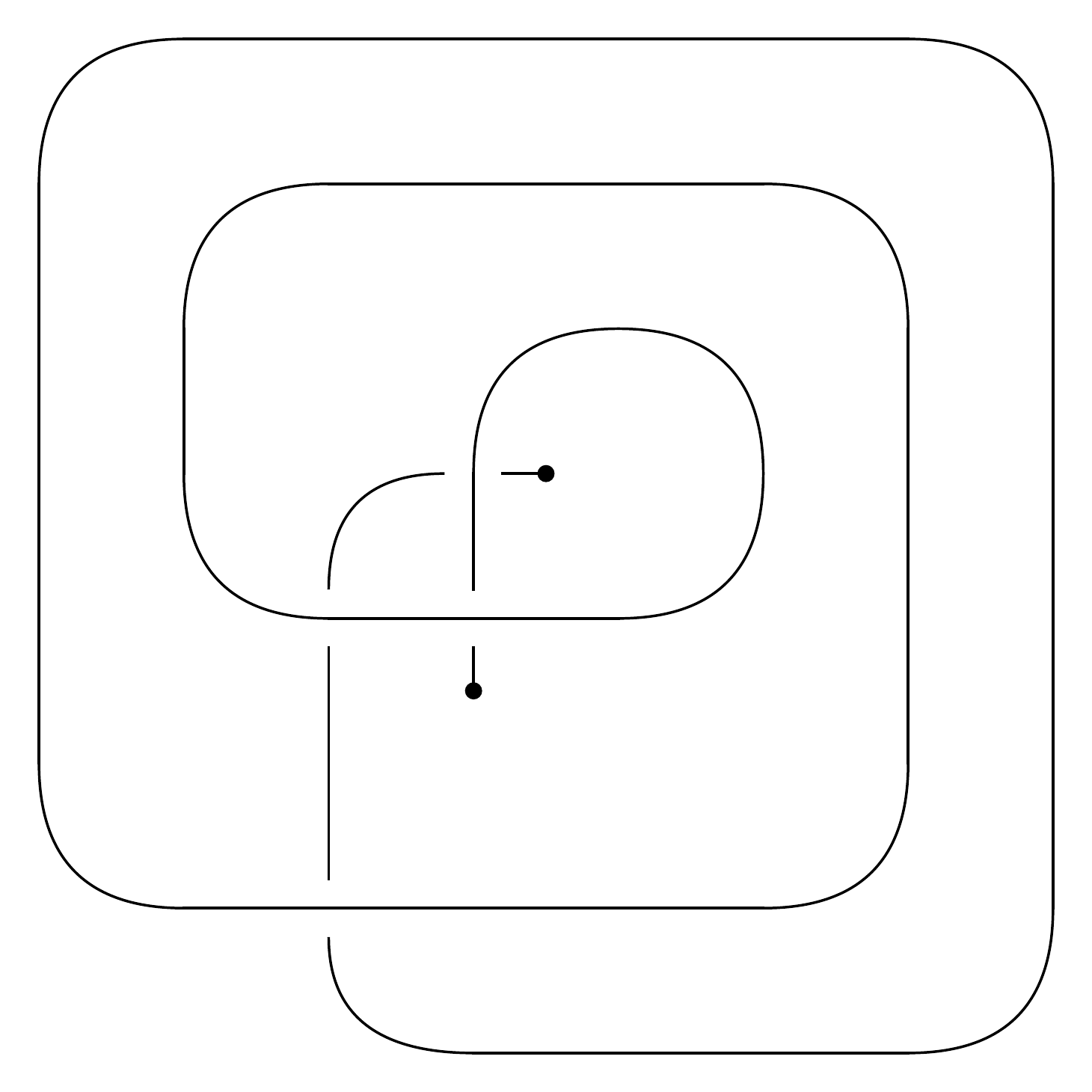}\\
\textcolor{black}{$4_{39}$}
\vspace{1cm}
\end{minipage}
\begin{minipage}[t]{.25\linewidth}
\centering
\includegraphics[width=0.9\textwidth,height=3.5cm,keepaspectratio]{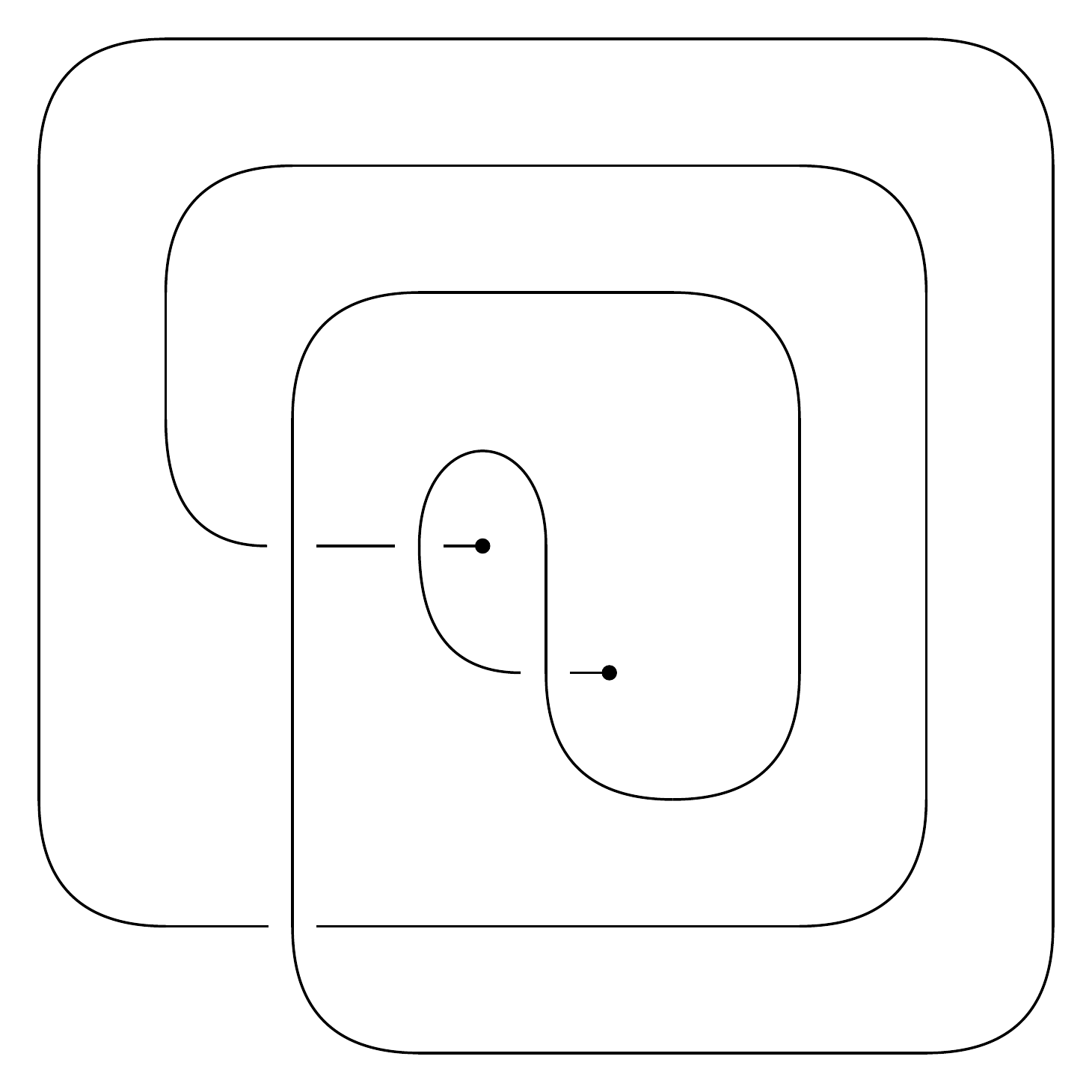}\\
\textcolor{black}{$4_{40}$}
\vspace{1cm}
\end{minipage}
\begin{minipage}[t]{.25\linewidth}
\centering
\includegraphics[width=0.9\textwidth,height=3.5cm,keepaspectratio]{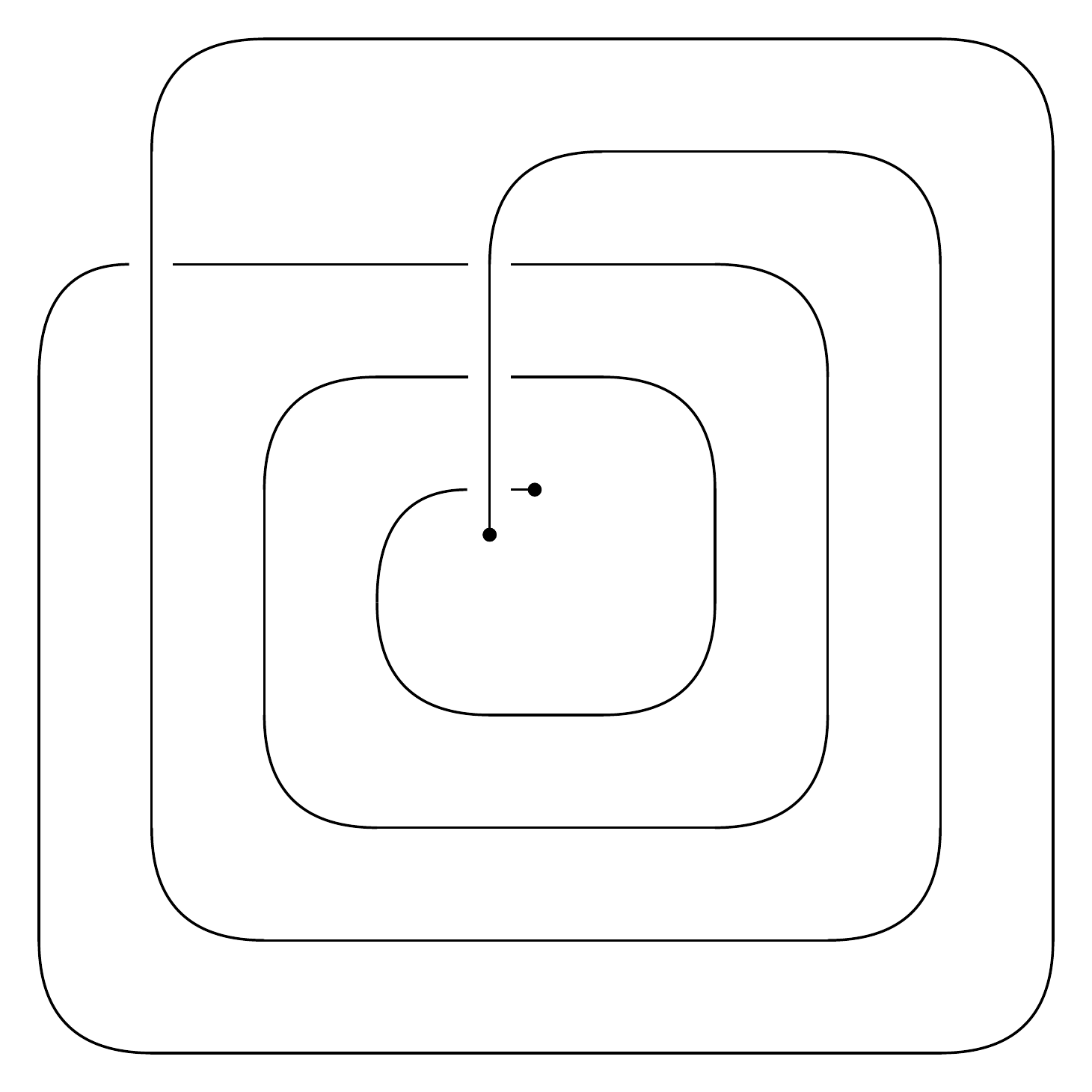}\\
\textcolor{black}{$4_{41}$}
\vspace{1cm}
\end{minipage}
\begin{minipage}[t]{.25\linewidth}
\centering
\includegraphics[width=0.9\textwidth,height=3.5cm,keepaspectratio]{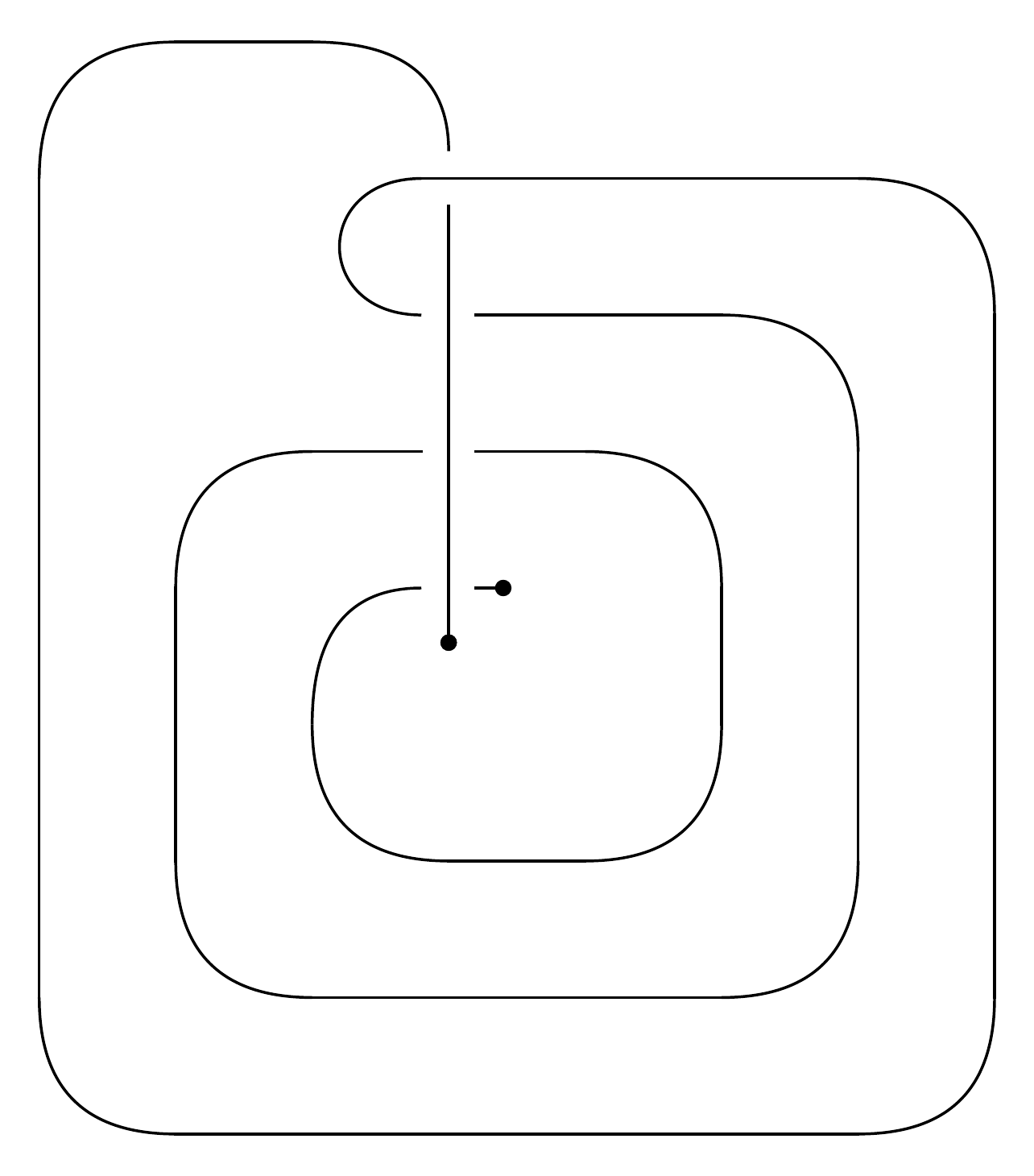}\\
\textcolor{black}{$4_{42}$}
\vspace{1cm}
\end{minipage}
\begin{minipage}[t]{.25\linewidth}
\centering
\includegraphics[width=0.9\textwidth,height=3.5cm,keepaspectratio]{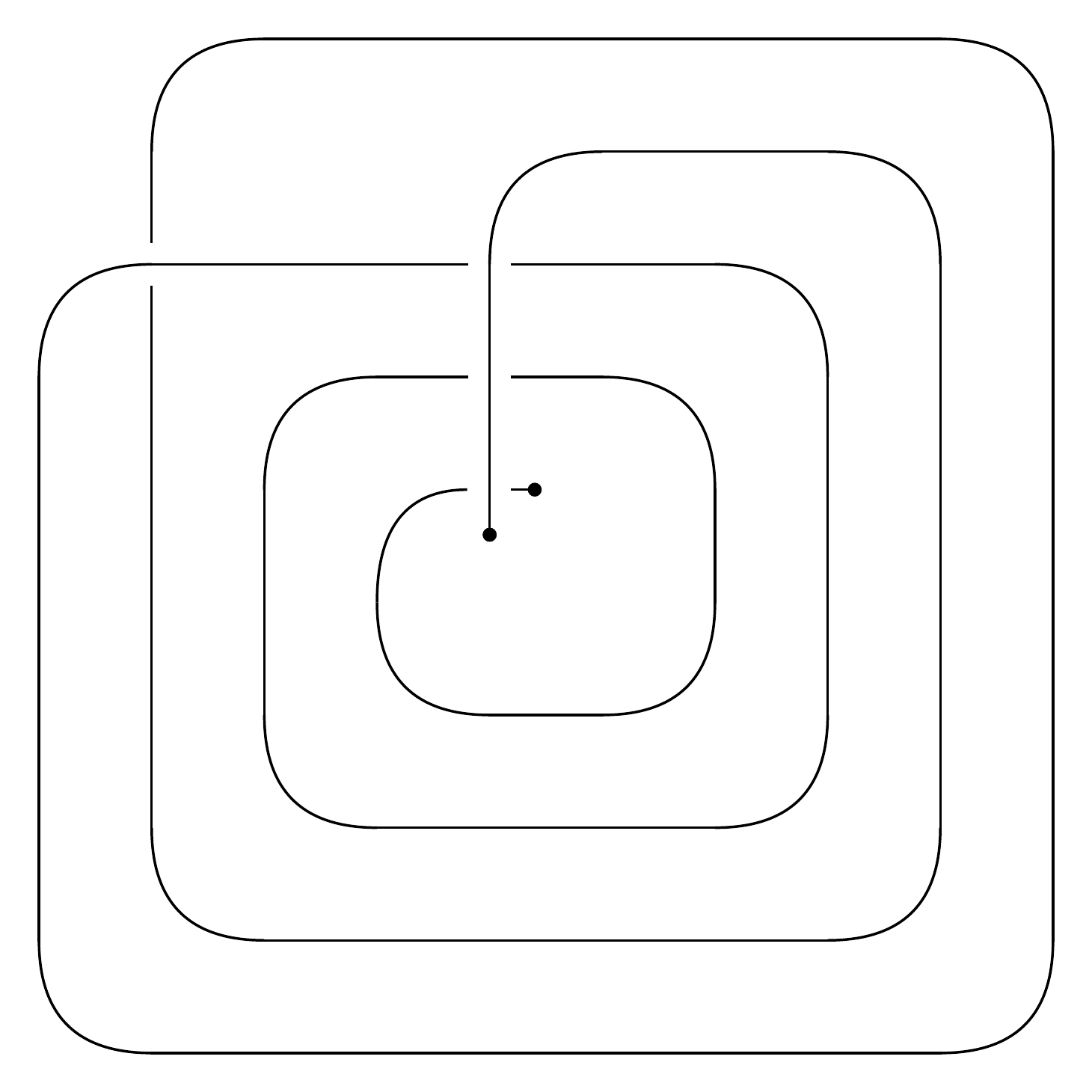}\\
\textcolor{black}{$4_{43}$}
\vspace{1cm}
\end{minipage}
\begin{minipage}[t]{.25\linewidth}
\centering
\includegraphics[width=0.9\textwidth,height=3.5cm,keepaspectratio]{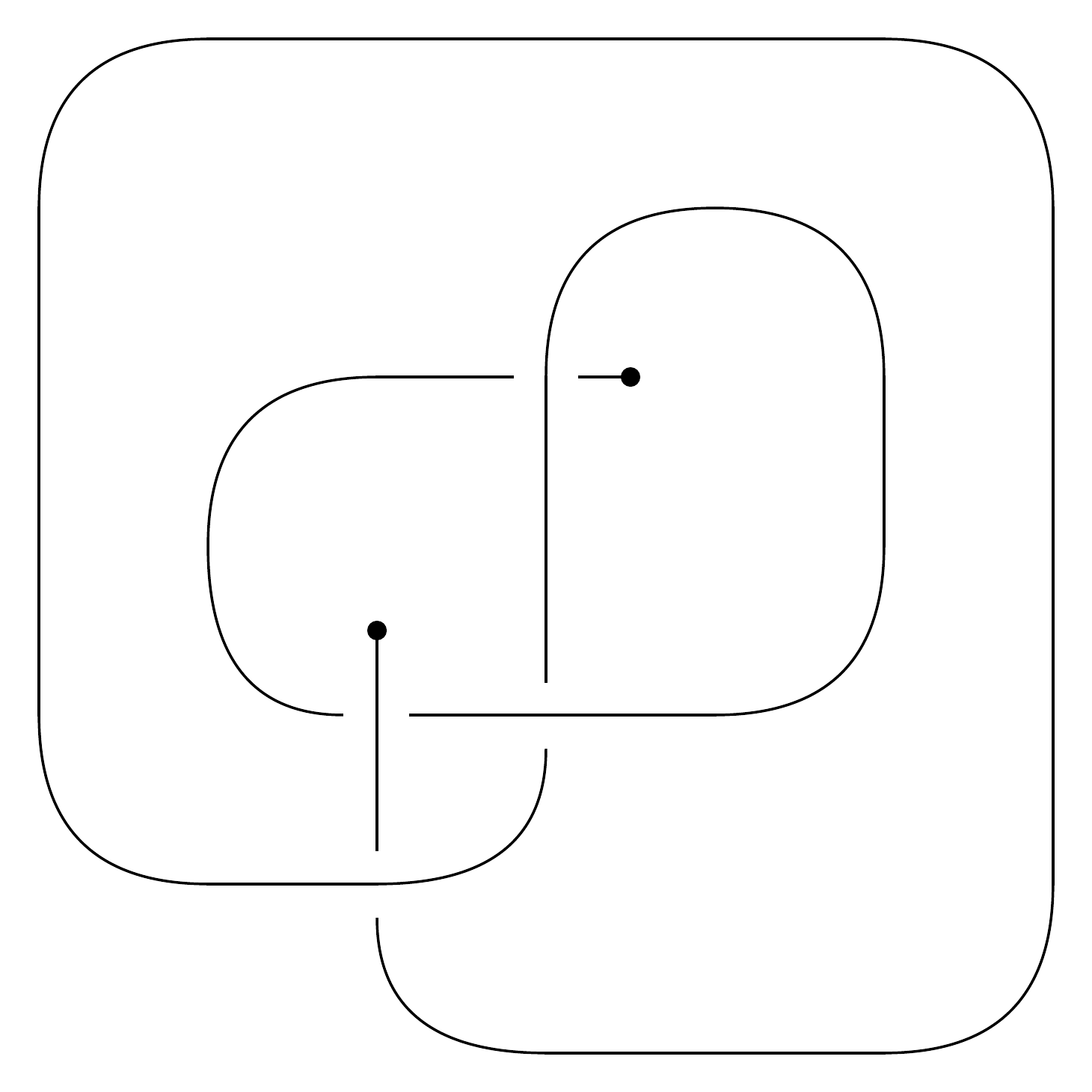}\\
\textcolor{black}{$4_{44}$}
\vspace{1cm}
\end{minipage}
\begin{minipage}[t]{.25\linewidth}
\centering
\includegraphics[width=0.9\textwidth,height=3.5cm,keepaspectratio]{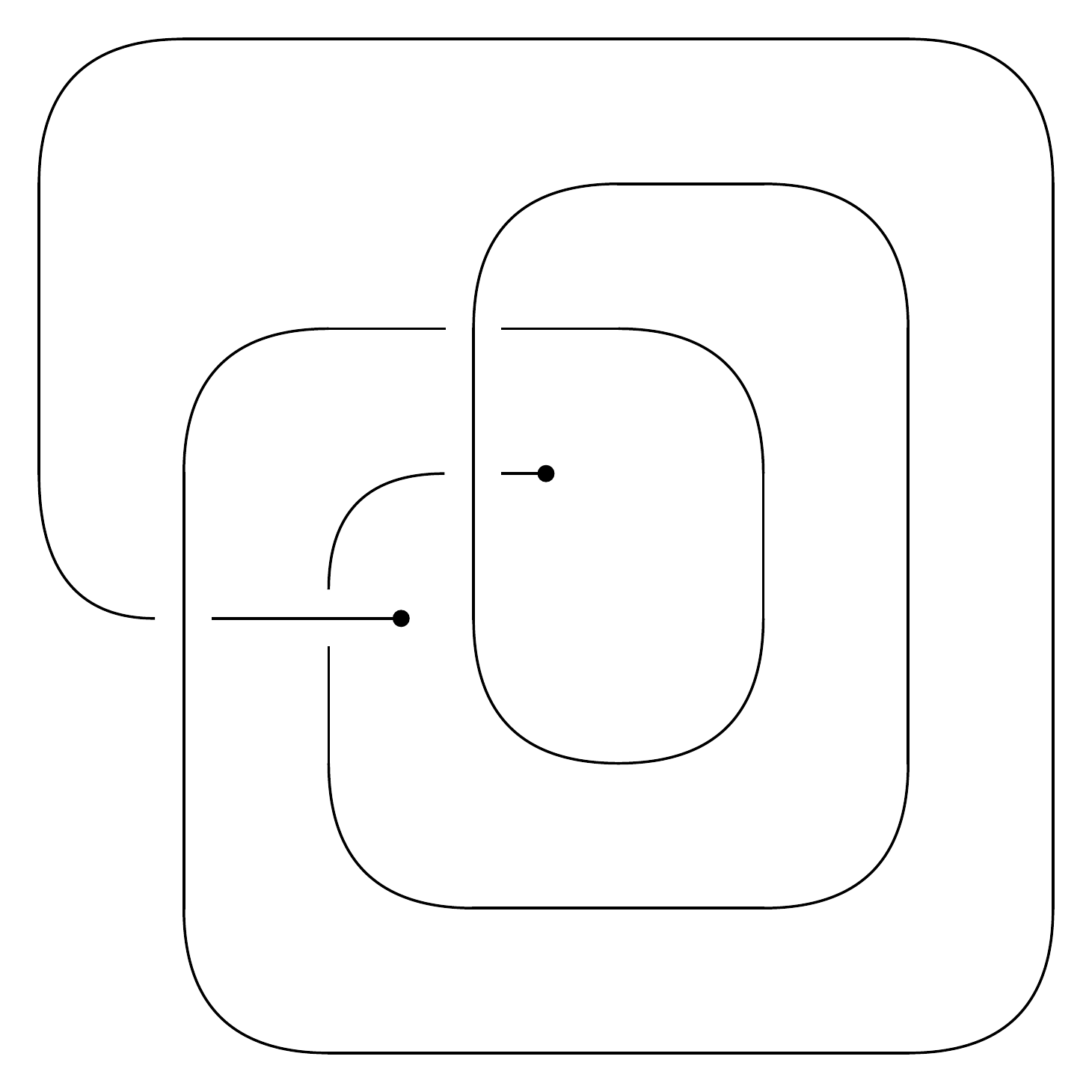}\\
\textcolor{black}{$4_{45}$}
\vspace{1cm}
\end{minipage}
\begin{minipage}[t]{.25\linewidth}
\centering
\includegraphics[width=0.9\textwidth,height=3.5cm,keepaspectratio]{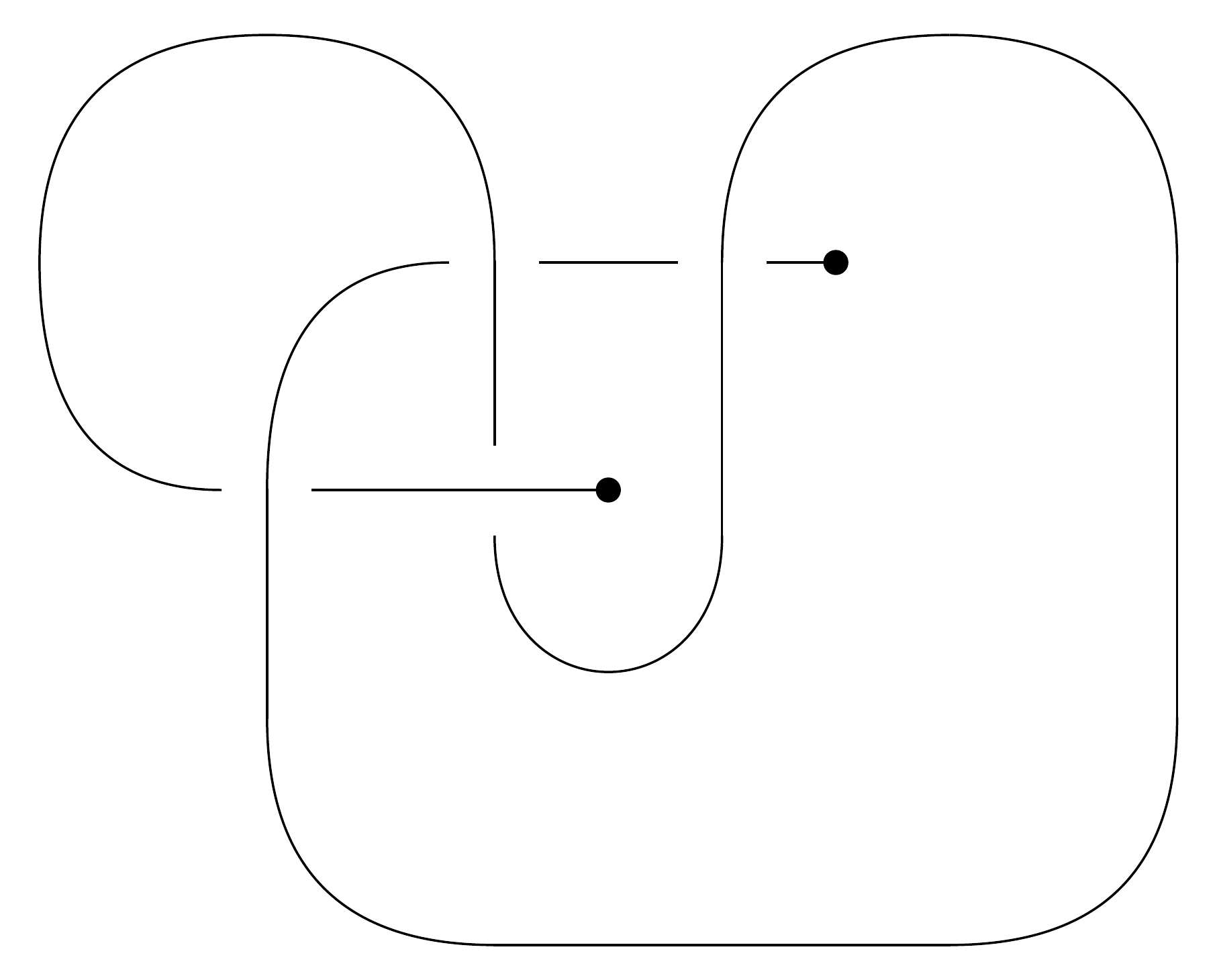}\\
\textcolor{black}{$4_{46}$}
\vspace{1cm}
\end{minipage}
\begin{minipage}[t]{.25\linewidth}
\centering
\includegraphics[width=0.9\textwidth,height=3.5cm,keepaspectratio]{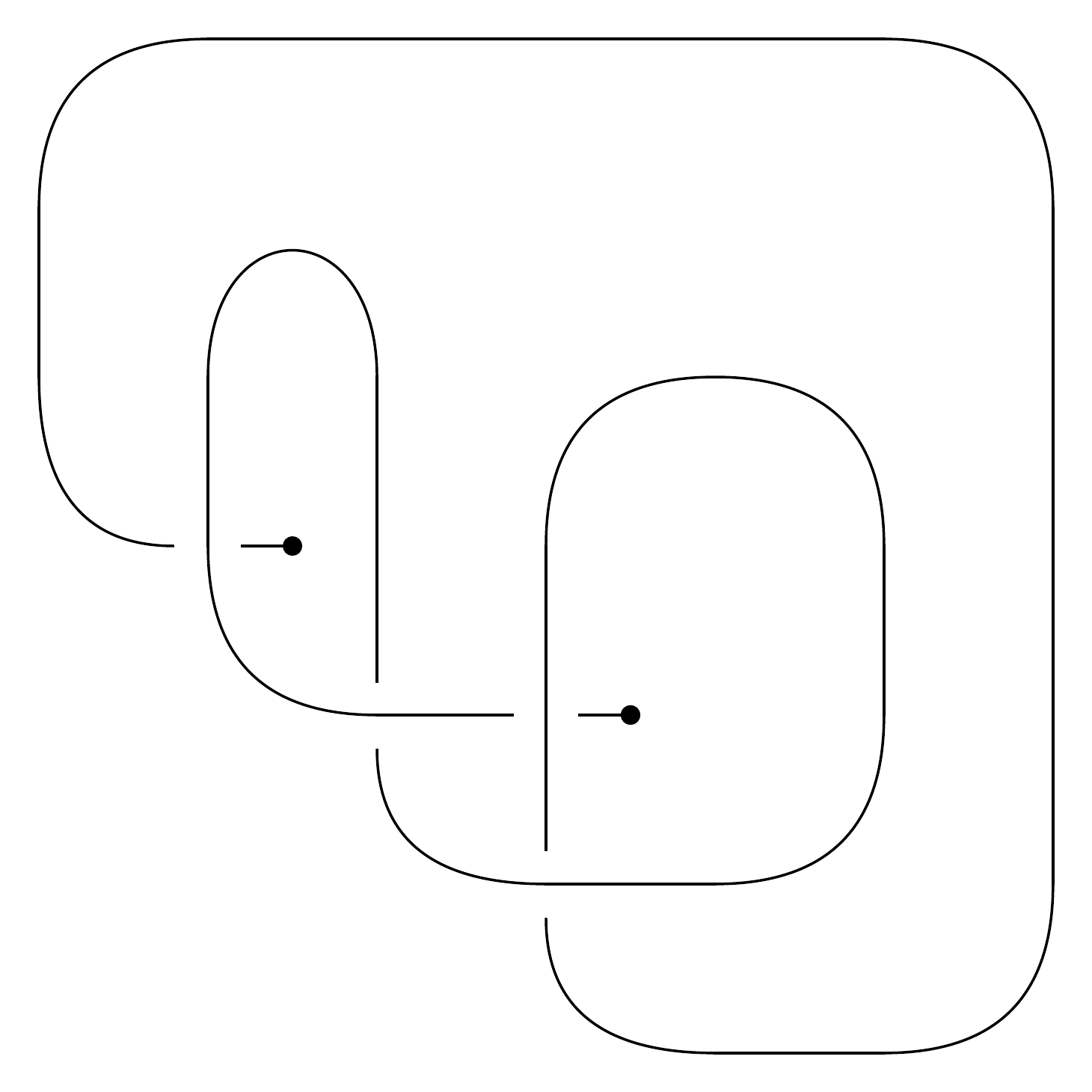}\\
\textcolor{black}{$4_{47}$}
\vspace{1cm}
\end{minipage}
\begin{minipage}[t]{.25\linewidth}
\centering
\includegraphics[width=0.9\textwidth,height=3.5cm,keepaspectratio]{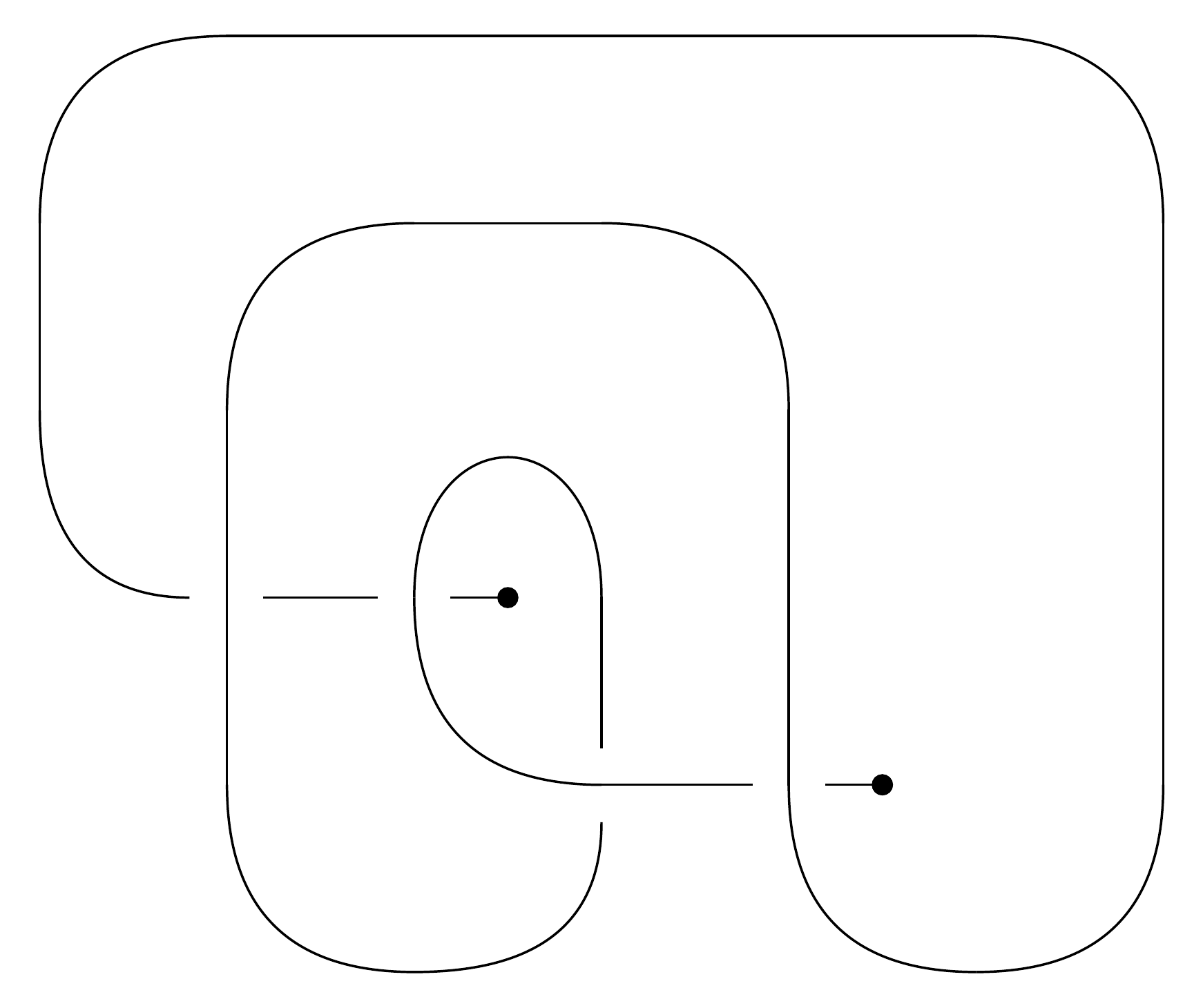}\\
\textcolor{black}{$4_{48}$}
\vspace{1cm}
\end{minipage}
\begin{minipage}[t]{.25\linewidth}
\centering
\includegraphics[width=0.9\textwidth,height=3.5cm,keepaspectratio]{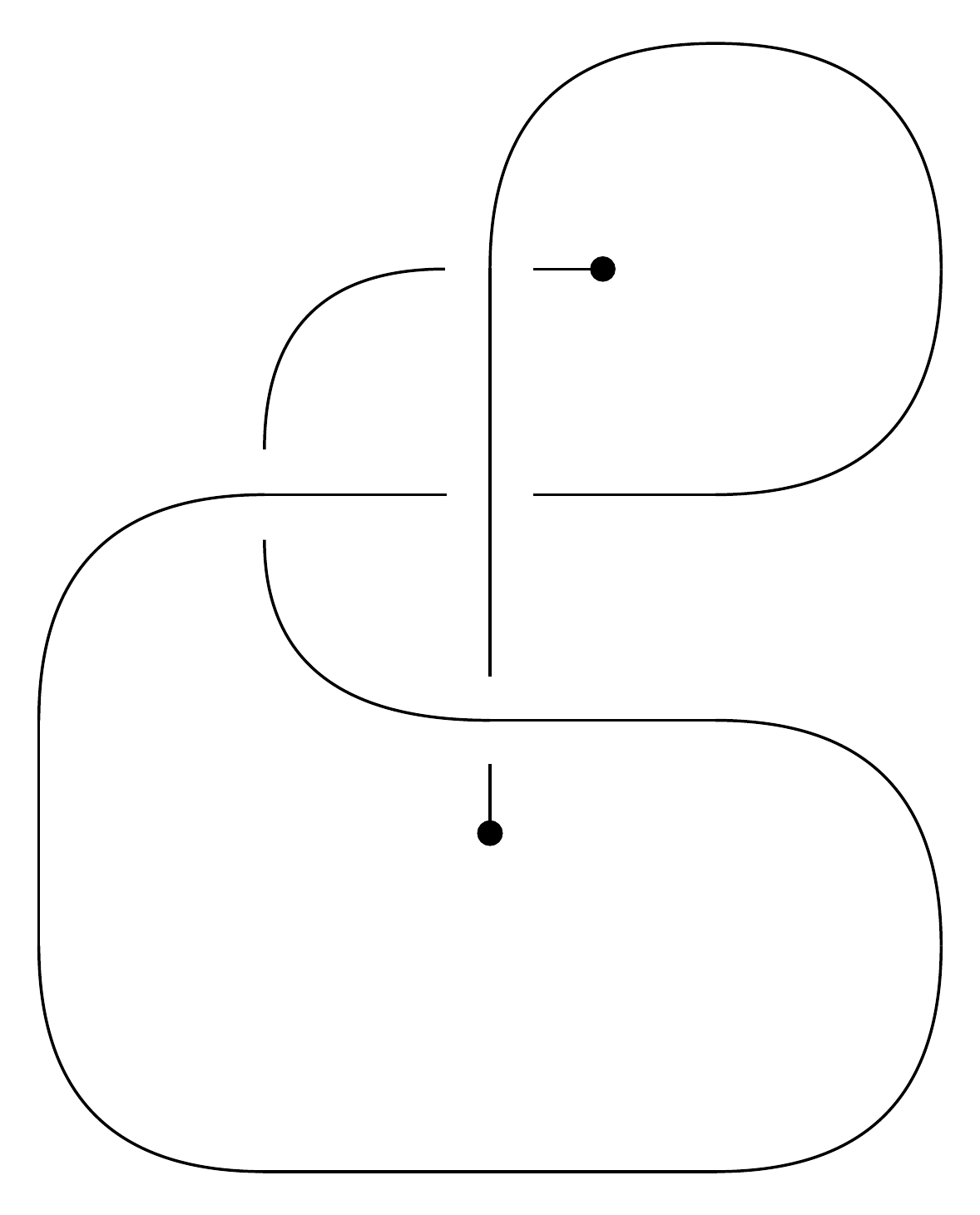}\\
\textcolor{black}{$4_{49}$}
\vspace{1cm}
\end{minipage}
\begin{minipage}[t]{.25\linewidth}
\centering
\includegraphics[width=0.9\textwidth,height=3.5cm,keepaspectratio]{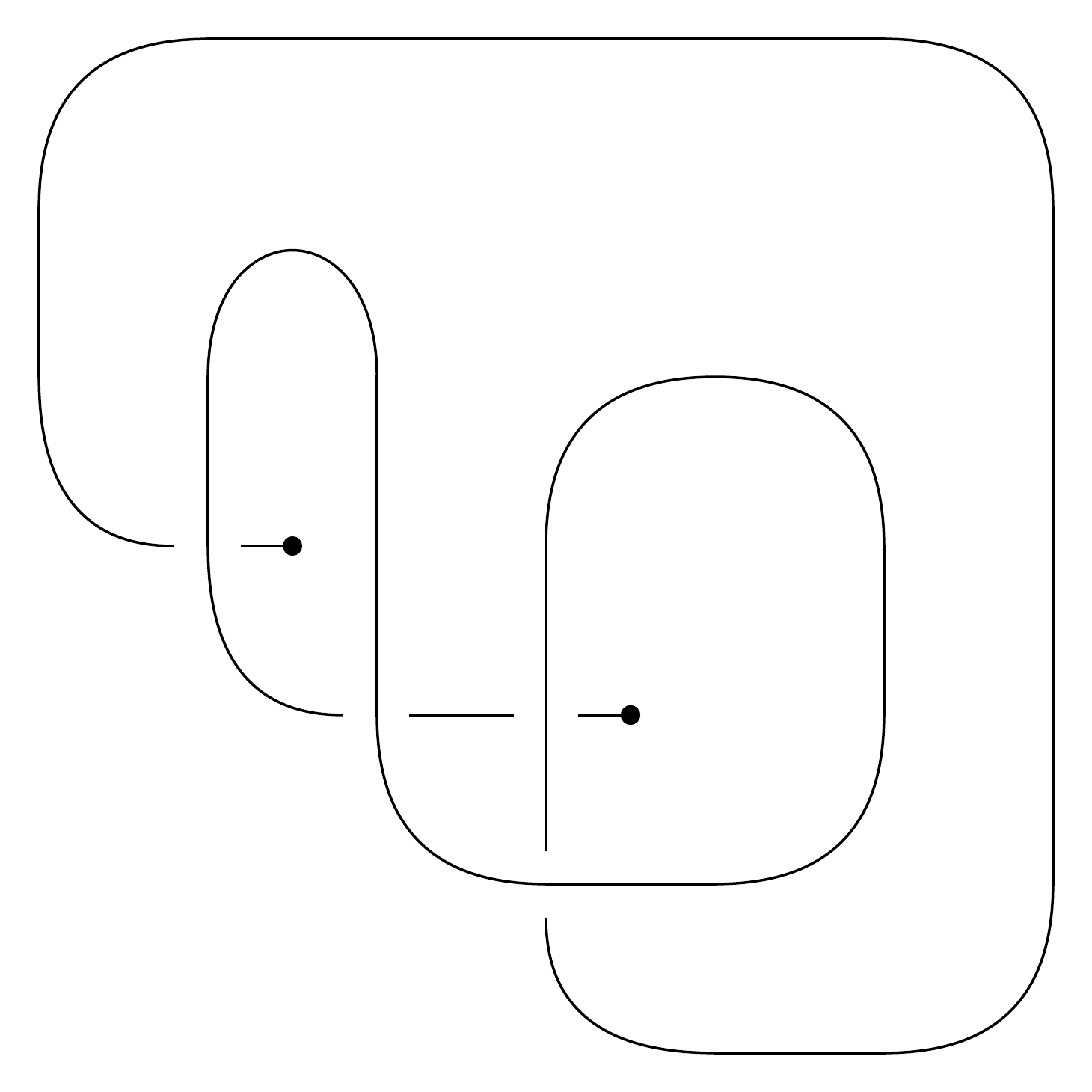}\\
\textcolor{black}{$4_{50}$}
\vspace{1cm}
\end{minipage}
\begin{minipage}[t]{.25\linewidth}
\centering
\includegraphics[width=0.9\textwidth,height=3.5cm,keepaspectratio]{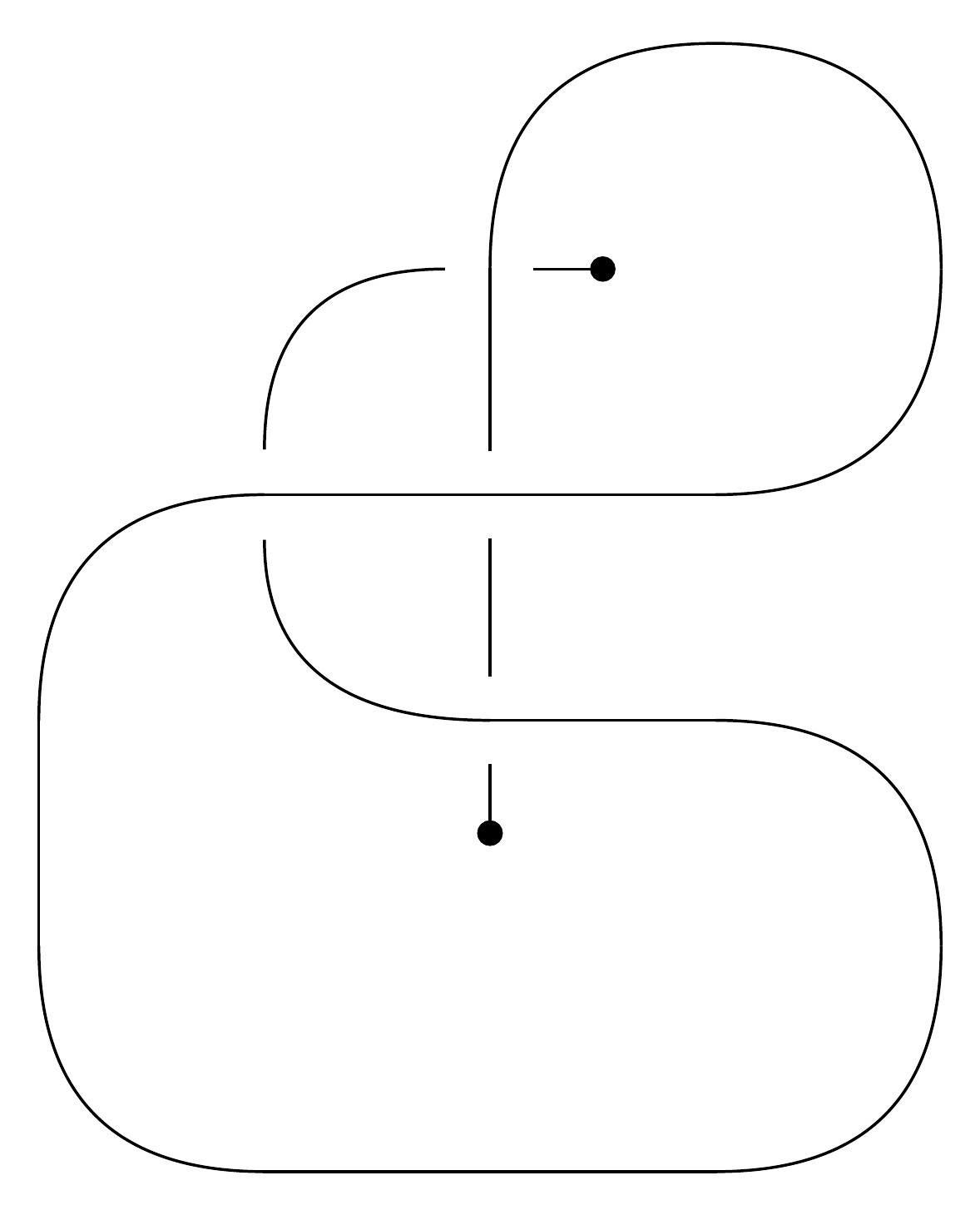}\\
\textcolor{black}{$4_{51}$}
\vspace{1cm}
\end{minipage}
\begin{minipage}[t]{.25\linewidth}
\centering
\includegraphics[width=0.9\textwidth,height=3.5cm,keepaspectratio]{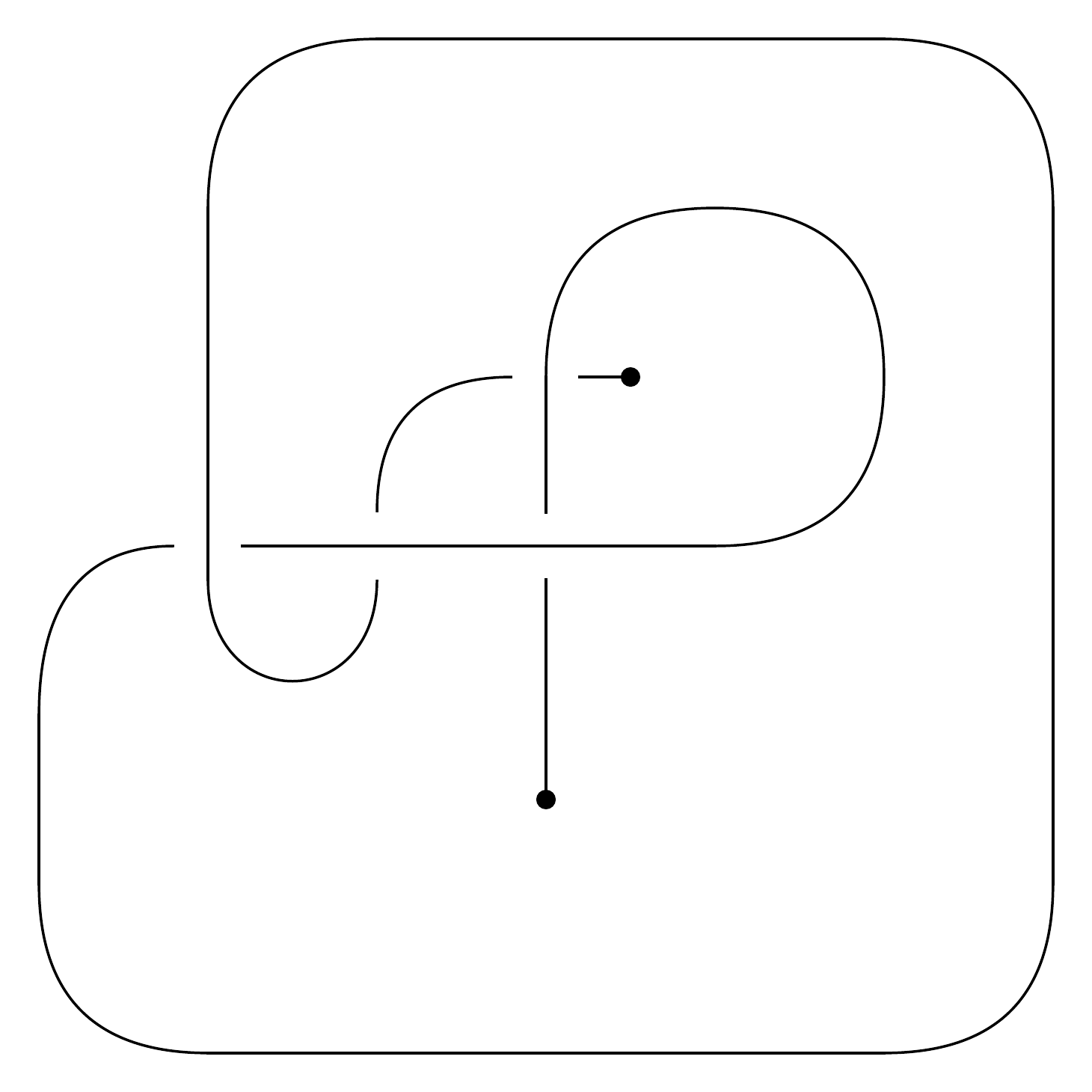}\\
\textcolor{black}{$4_{52}$}
\vspace{1cm}
\end{minipage}
\begin{minipage}[t]{.25\linewidth}
\centering
\includegraphics[width=0.9\textwidth,height=3.5cm,keepaspectratio]{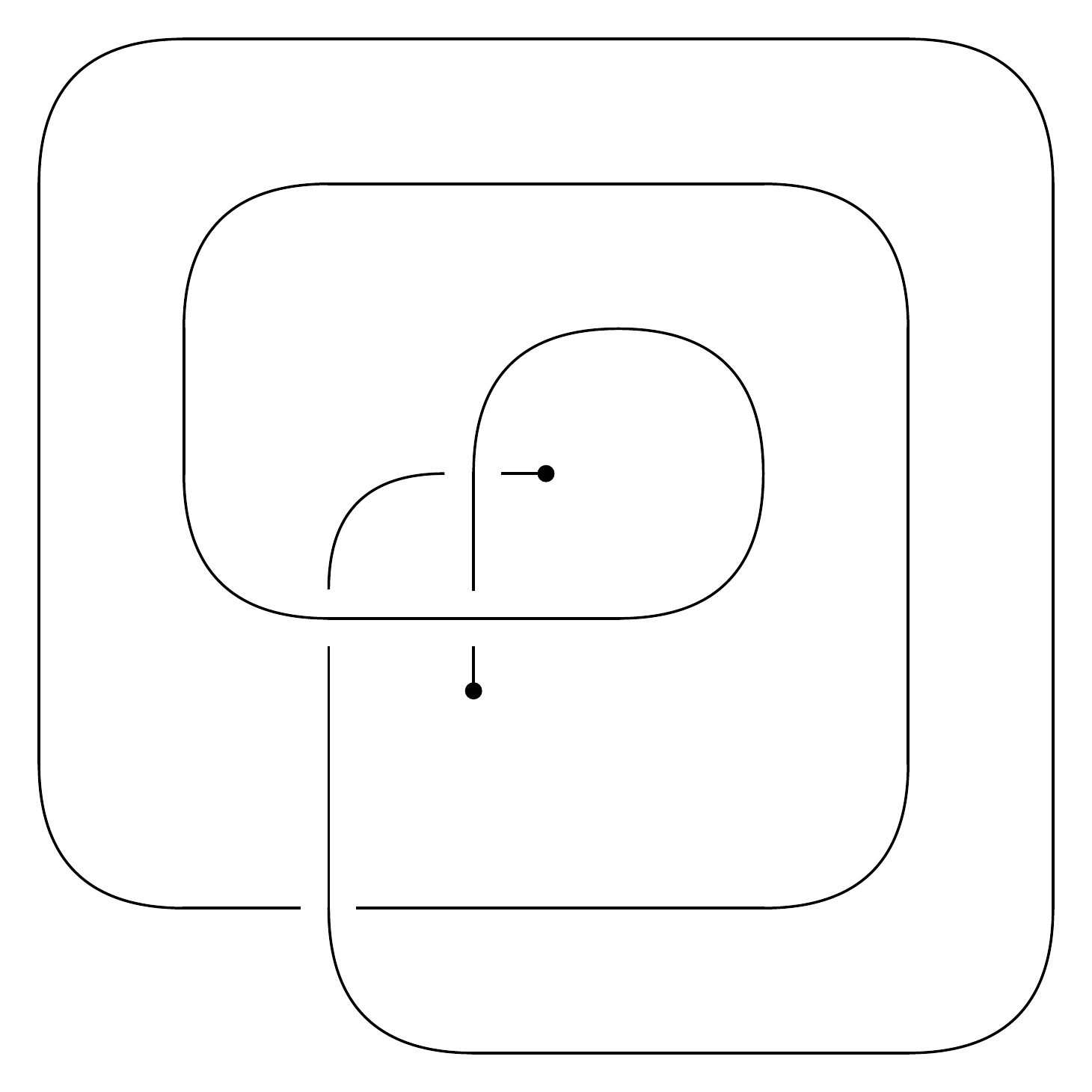}\\
\textcolor{black}{$4_{53}$}
\vspace{1cm}
\end{minipage}
\begin{minipage}[t]{.25\linewidth}
\centering
\includegraphics[width=0.9\textwidth,height=3.5cm,keepaspectratio]{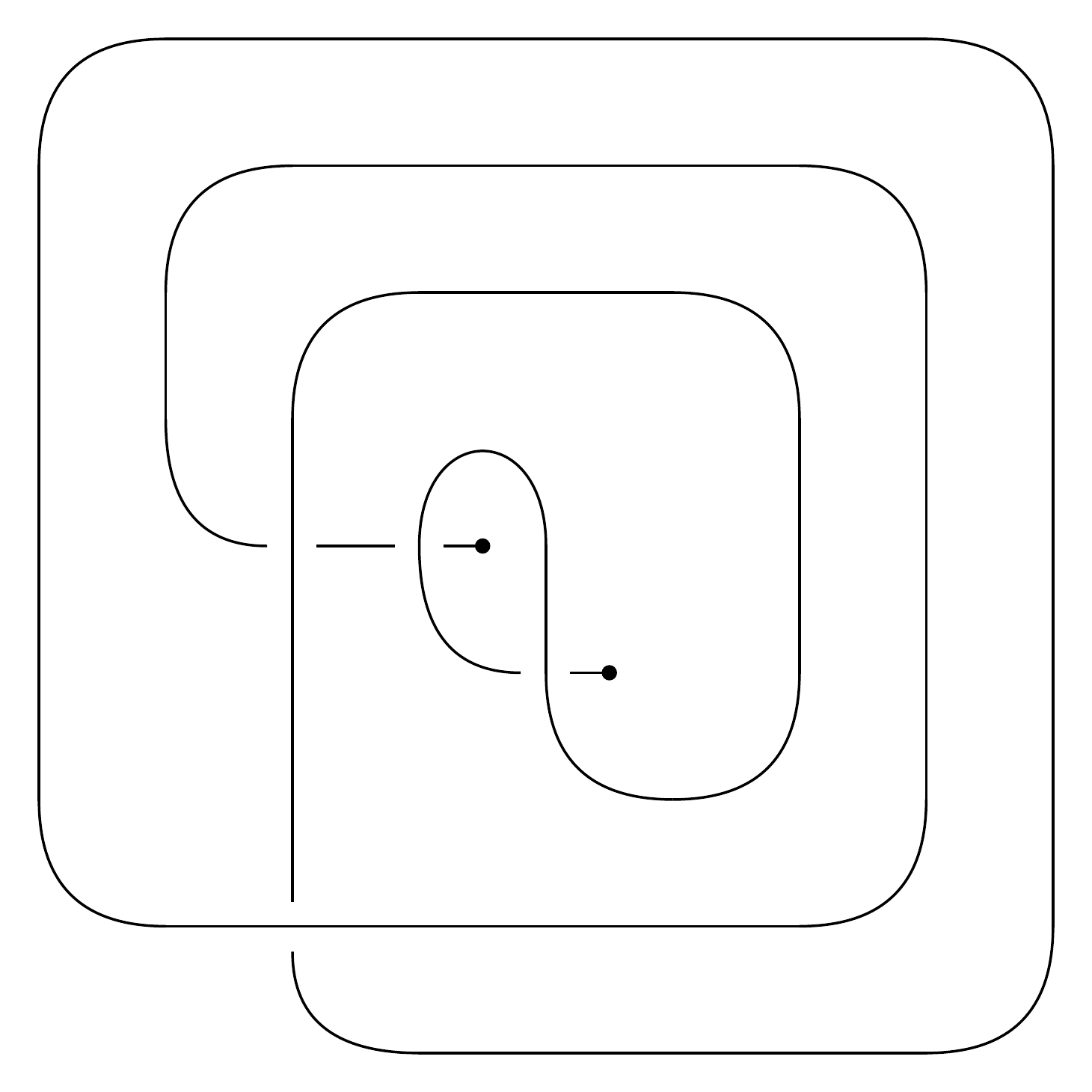}\\
\textcolor{black}{$4_{54}$}
\vspace{1cm}
\end{minipage}
\begin{minipage}[t]{.25\linewidth}
\centering
\includegraphics[width=0.9\textwidth,height=3.5cm,keepaspectratio]{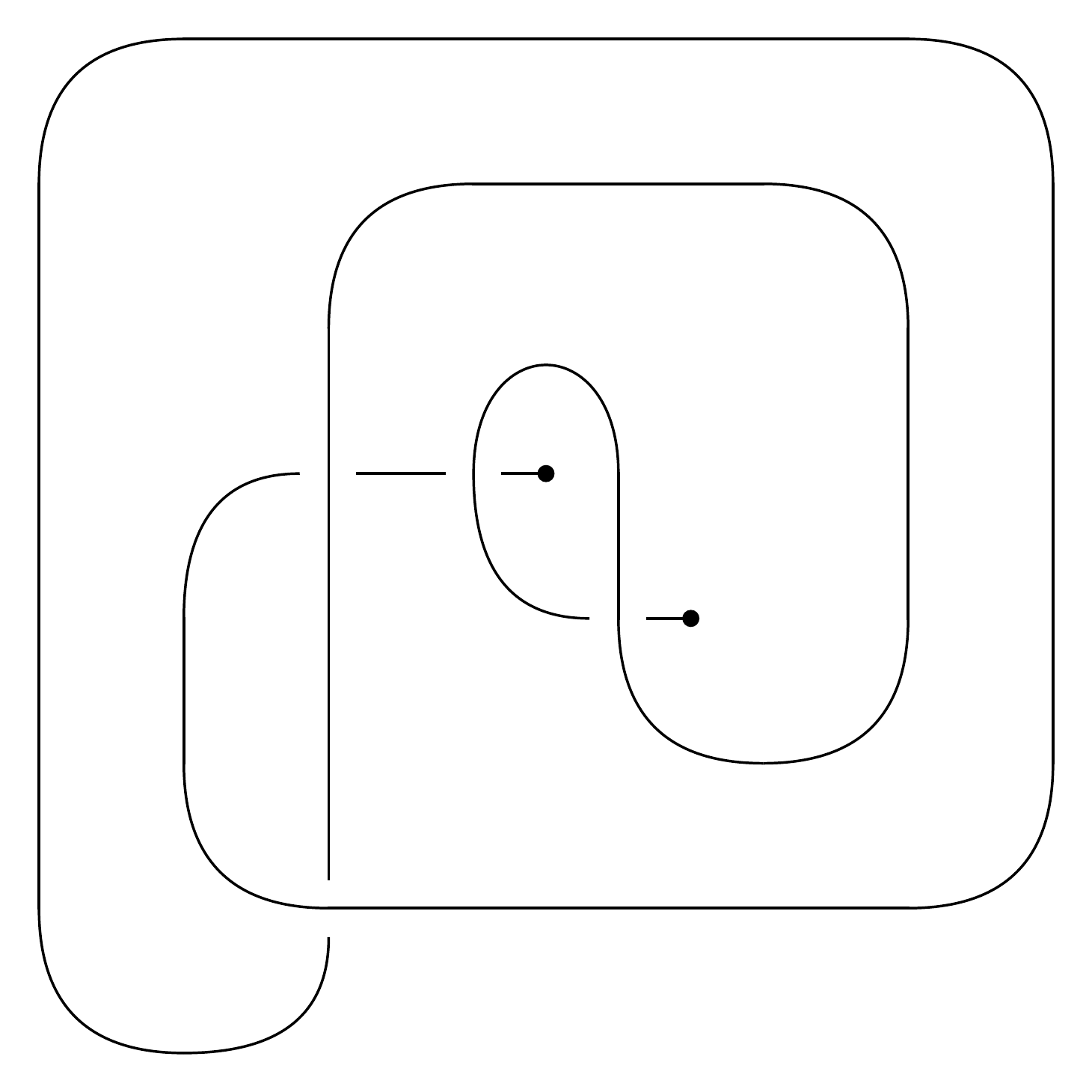}\\
\textcolor{black}{$4_{55}$}
\vspace{1cm}
\end{minipage}
\begin{minipage}[t]{.25\linewidth}
\centering
\includegraphics[width=0.9\textwidth,height=3.5cm,keepaspectratio]{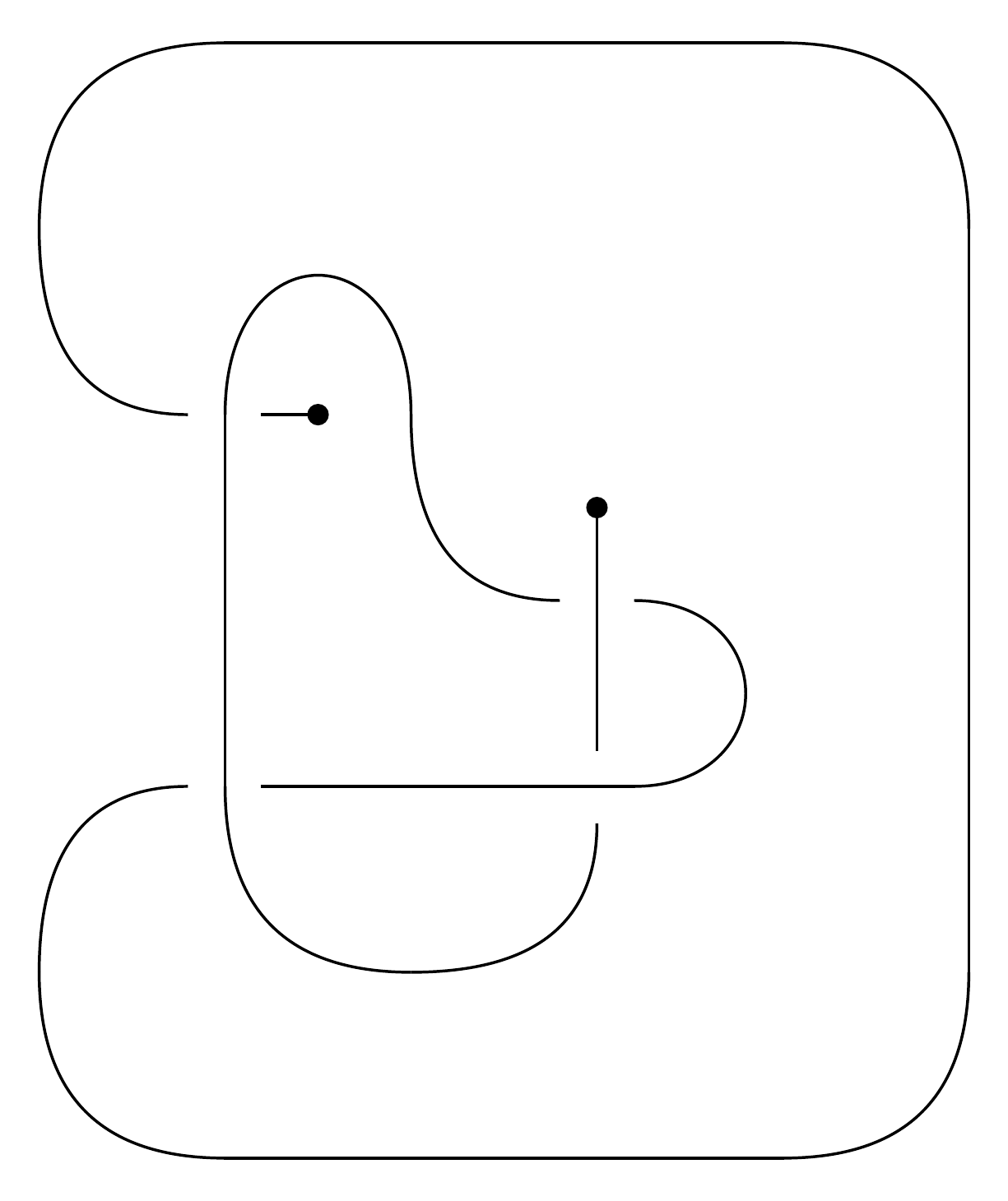}\\
\textcolor{black}{$4_{56}$}
\vspace{1cm}
\end{minipage}
\begin{minipage}[t]{.25\linewidth}
\centering
\includegraphics[width=0.9\textwidth,height=3.5cm,keepaspectratio]{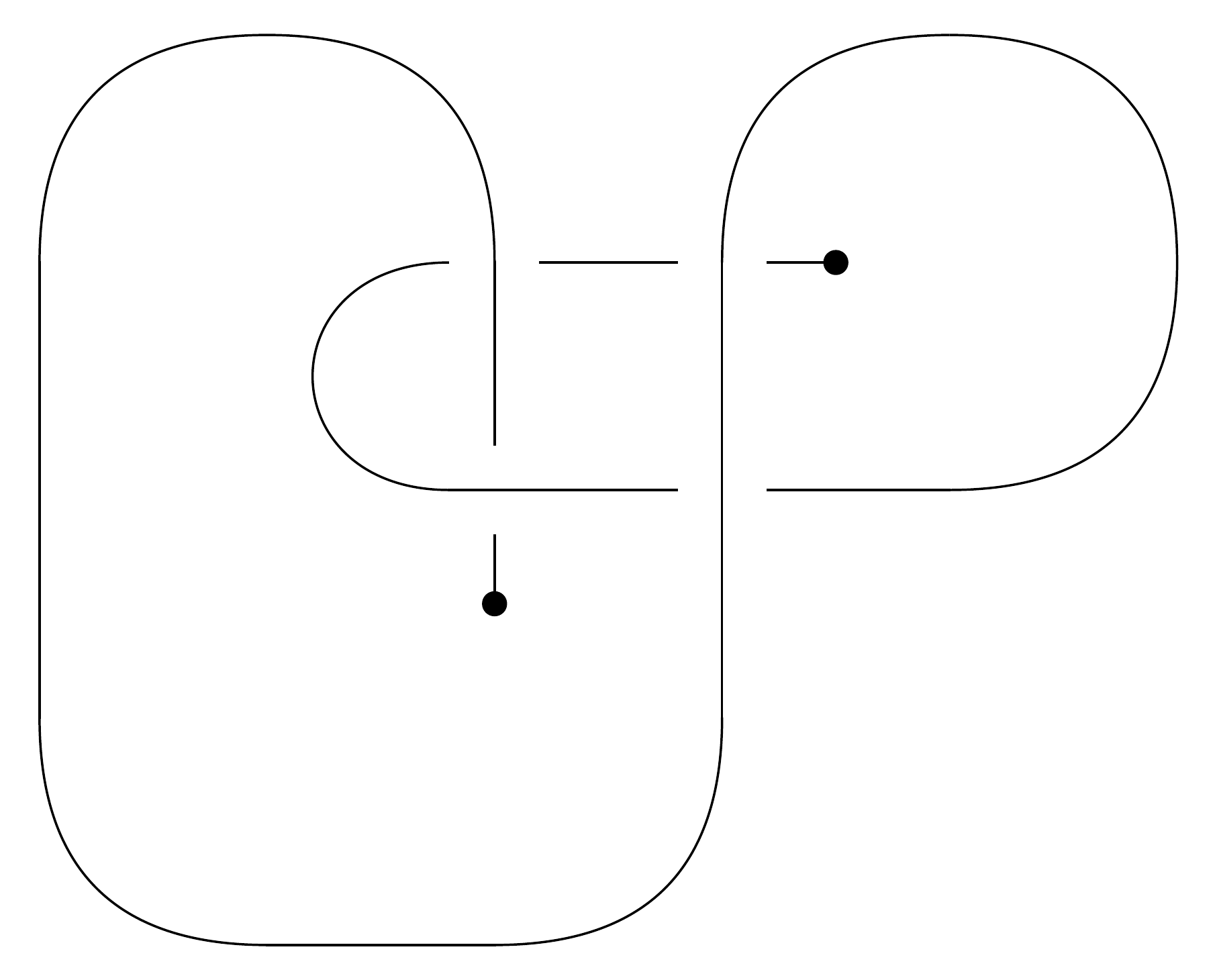}\\
\textcolor{black}{$4_{57}$}
\vspace{1cm}
\end{minipage}
\begin{minipage}[t]{.25\linewidth}
\centering
\includegraphics[width=0.9\textwidth,height=3.5cm,keepaspectratio]{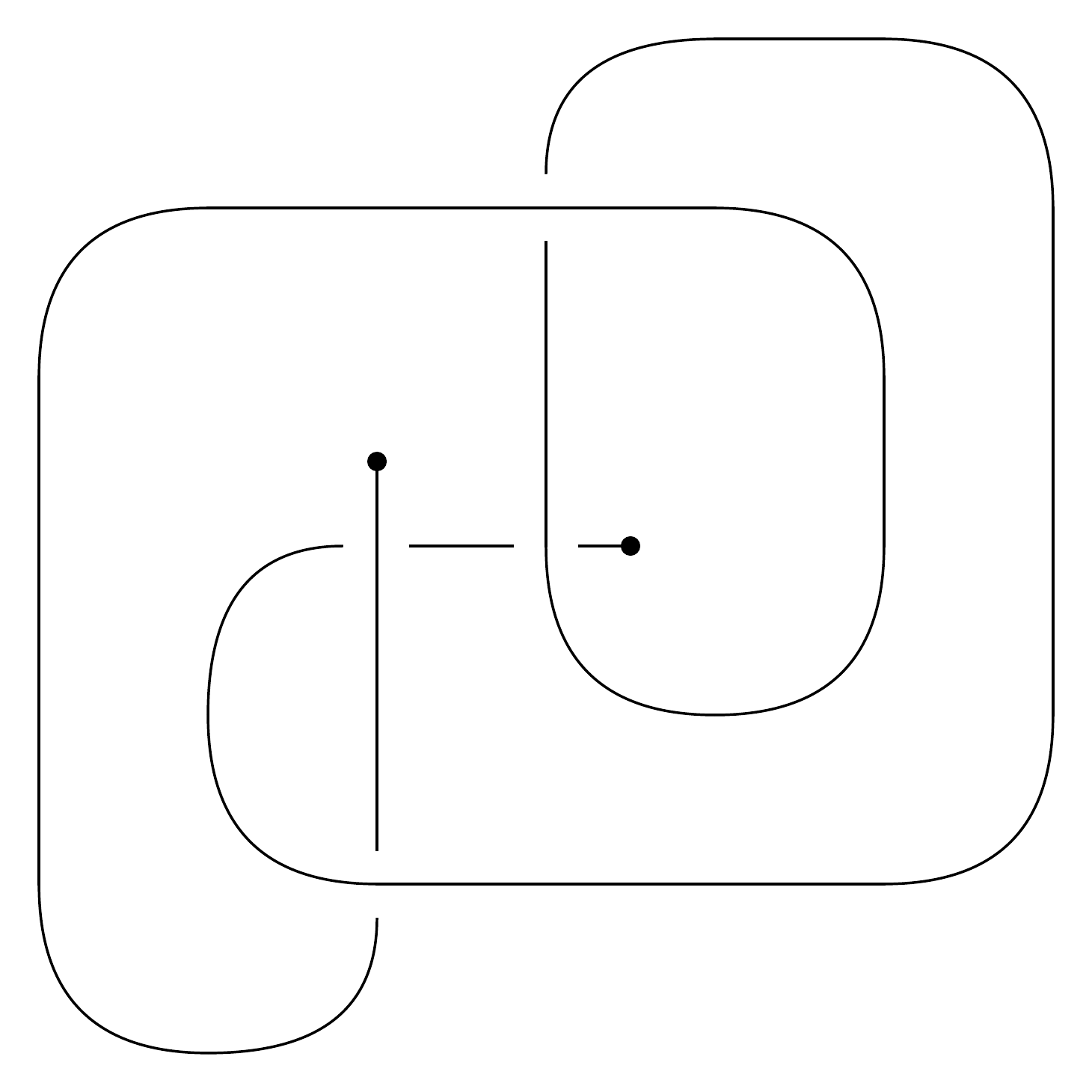}\\
\textcolor{black}{$4_{58}$}
\vspace{1cm}
\end{minipage}
\begin{minipage}[t]{.25\linewidth}
\centering
\includegraphics[width=0.9\textwidth,height=3.5cm,keepaspectratio]{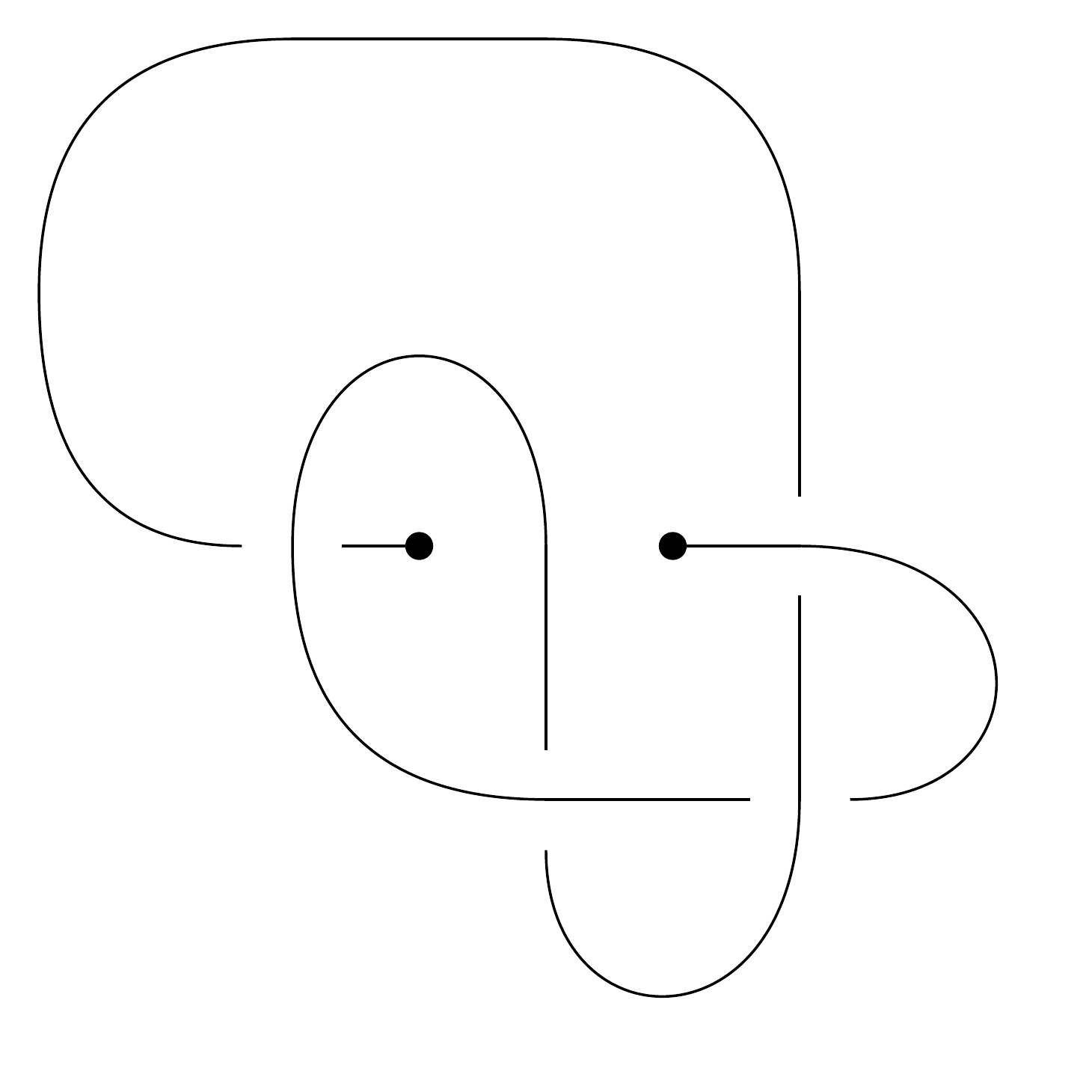}\\
\textcolor{black}{$4_{59}$}
\vspace{1cm}
\end{minipage}
\begin{minipage}[t]{.25\linewidth}
\centering
\includegraphics[width=0.9\textwidth,height=3.5cm,keepaspectratio]{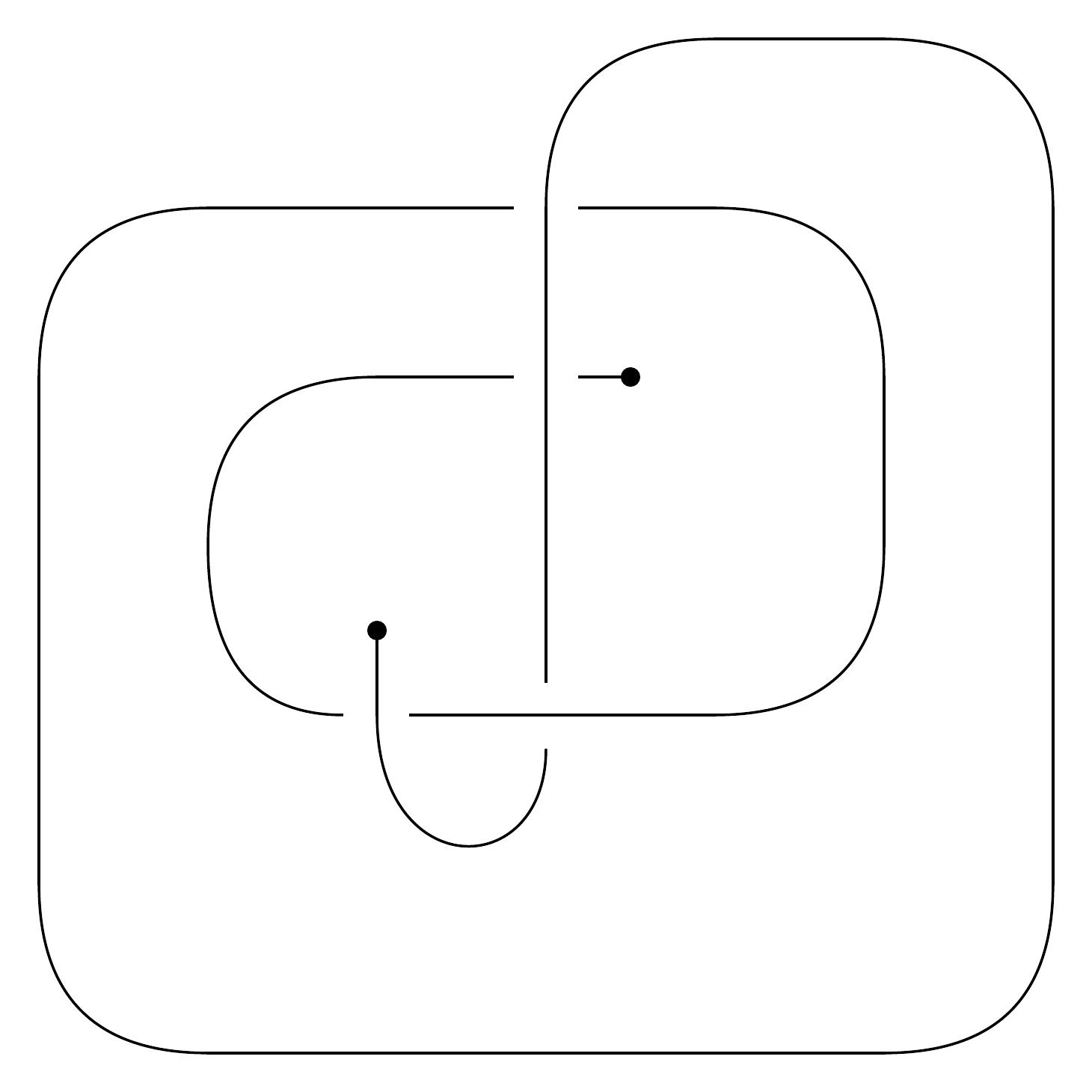}\\
\textcolor{black}{$4_{60}$}
\vspace{1cm}
\end{minipage}
\begin{minipage}[t]{.25\linewidth}
\centering
\includegraphics[width=0.9\textwidth,height=3.5cm,keepaspectratio]{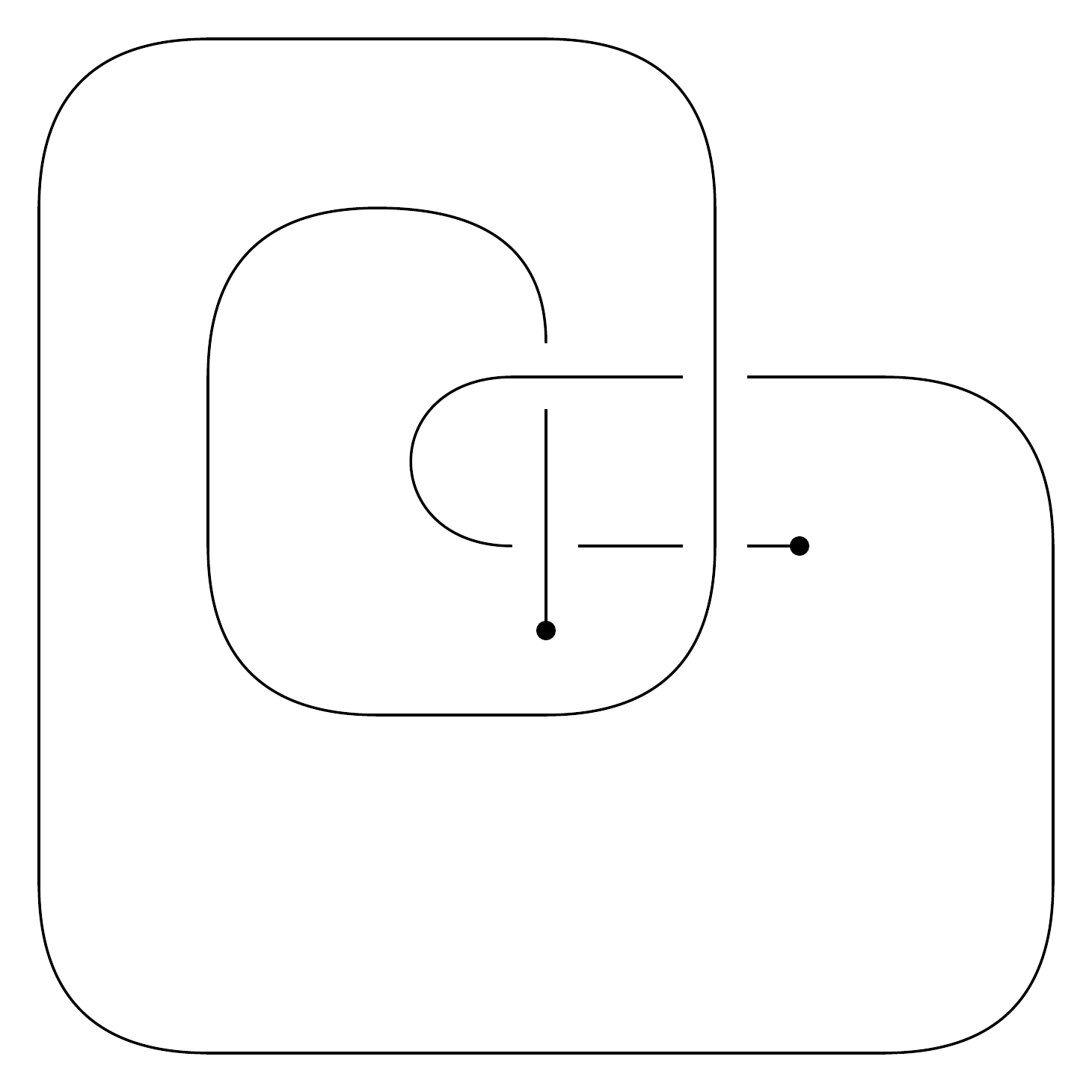}\\
\textcolor{black}{$4_{61}$}
\vspace{1cm}
\end{minipage}
\begin{minipage}[t]{.25\linewidth}
\centering
\includegraphics[width=0.9\textwidth,height=3.5cm,keepaspectratio]{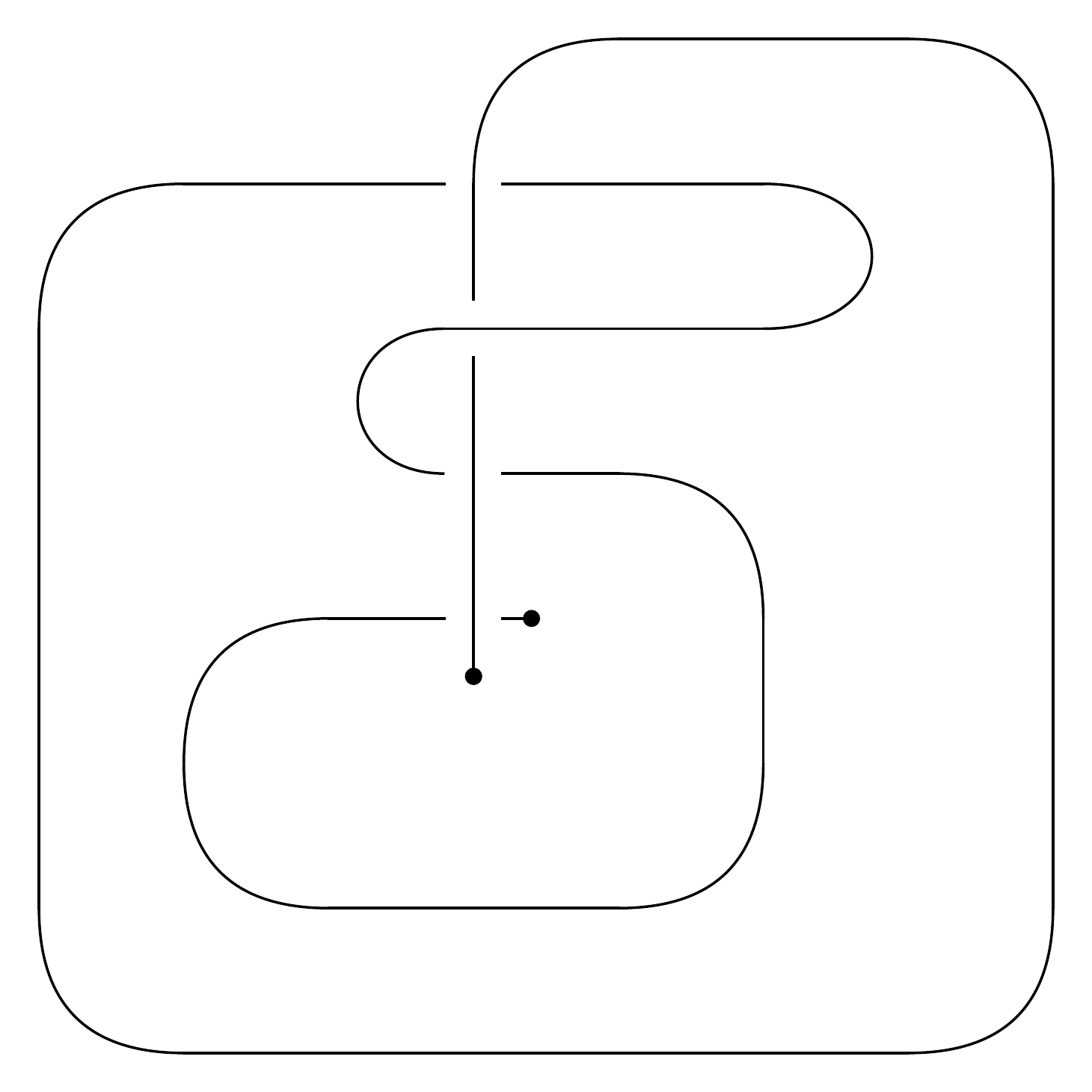}\\
\textcolor{black}{$4_{62}$}
\vspace{1cm}
\end{minipage}
\begin{minipage}[t]{.25\linewidth}
\centering
\includegraphics[width=0.9\textwidth,height=3.5cm,keepaspectratio]{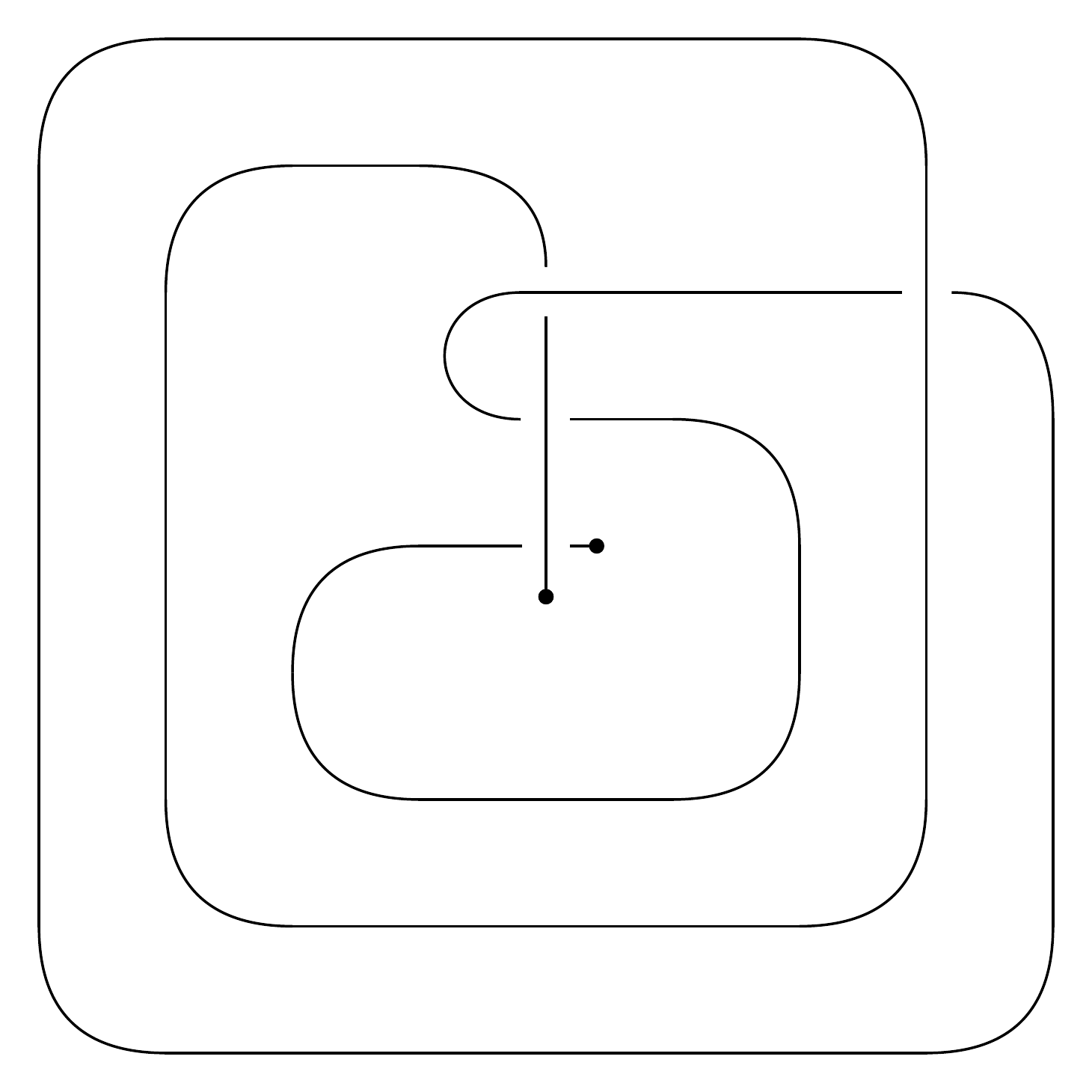}\\
\textcolor{black}{$4_{63}$}
\vspace{1cm}
\end{minipage}
\begin{minipage}[t]{.25\linewidth}
\centering
\includegraphics[width=0.9\textwidth,height=3.5cm,keepaspectratio]{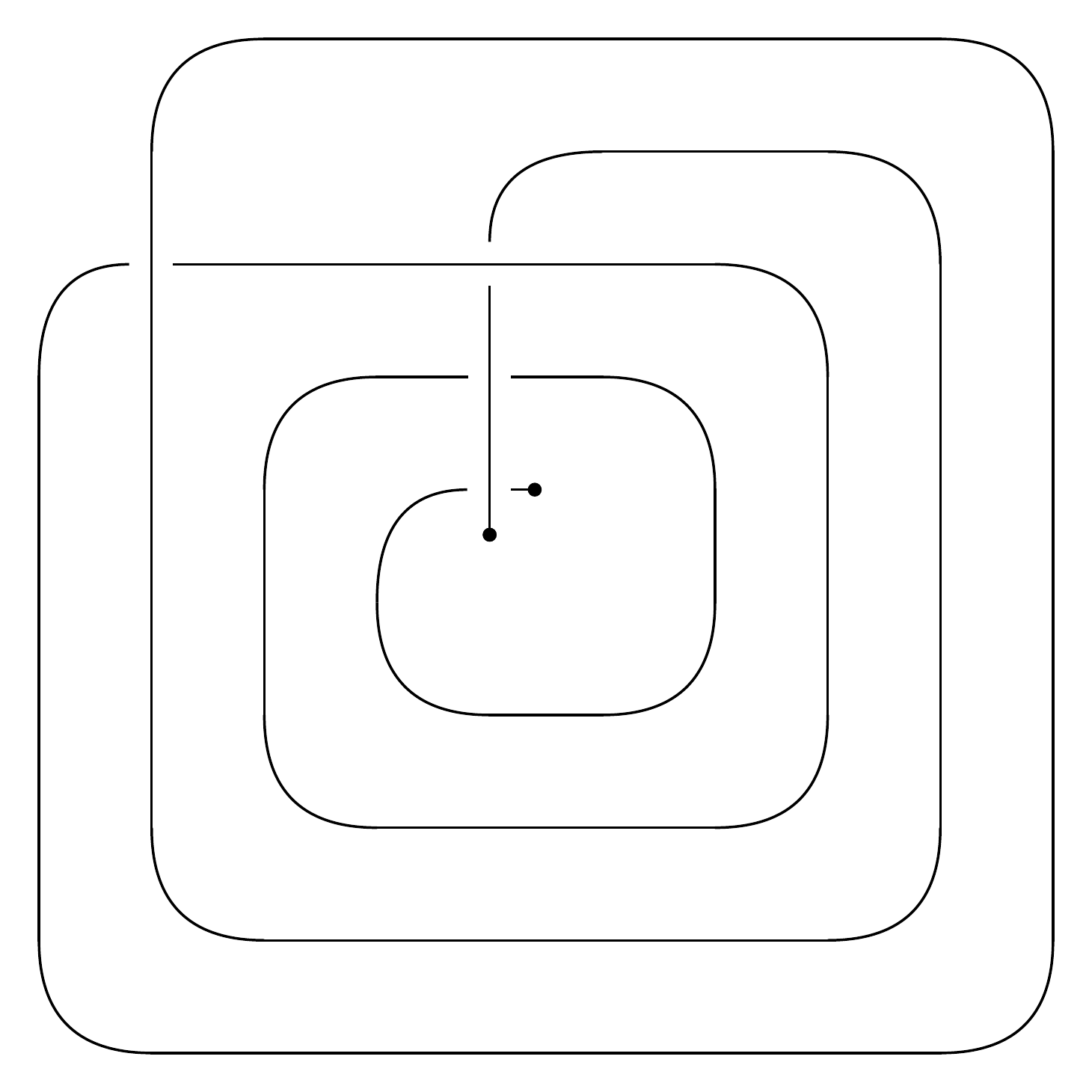}\\
\textcolor{black}{$4_{64}$}
\vspace{1cm}
\end{minipage}
\begin{minipage}[t]{.25\linewidth}
\centering
\includegraphics[width=0.9\textwidth,height=3.5cm,keepaspectratio]{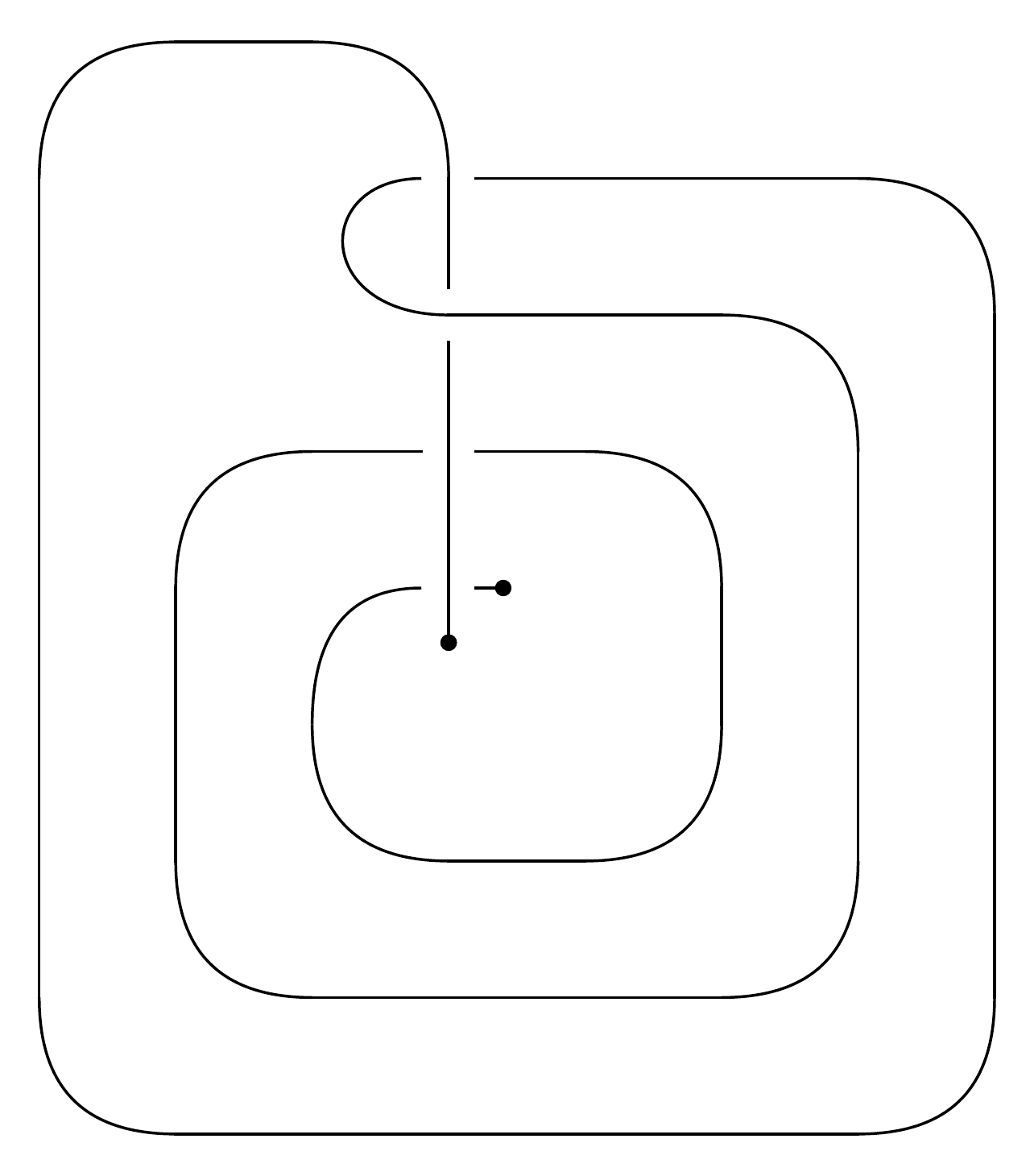}\\
\textcolor{black}{$4_{65}$}
\vspace{1cm}
\end{minipage}
\begin{minipage}[t]{.25\linewidth}
\centering
\includegraphics[width=0.9\textwidth,height=3.5cm,keepaspectratio]{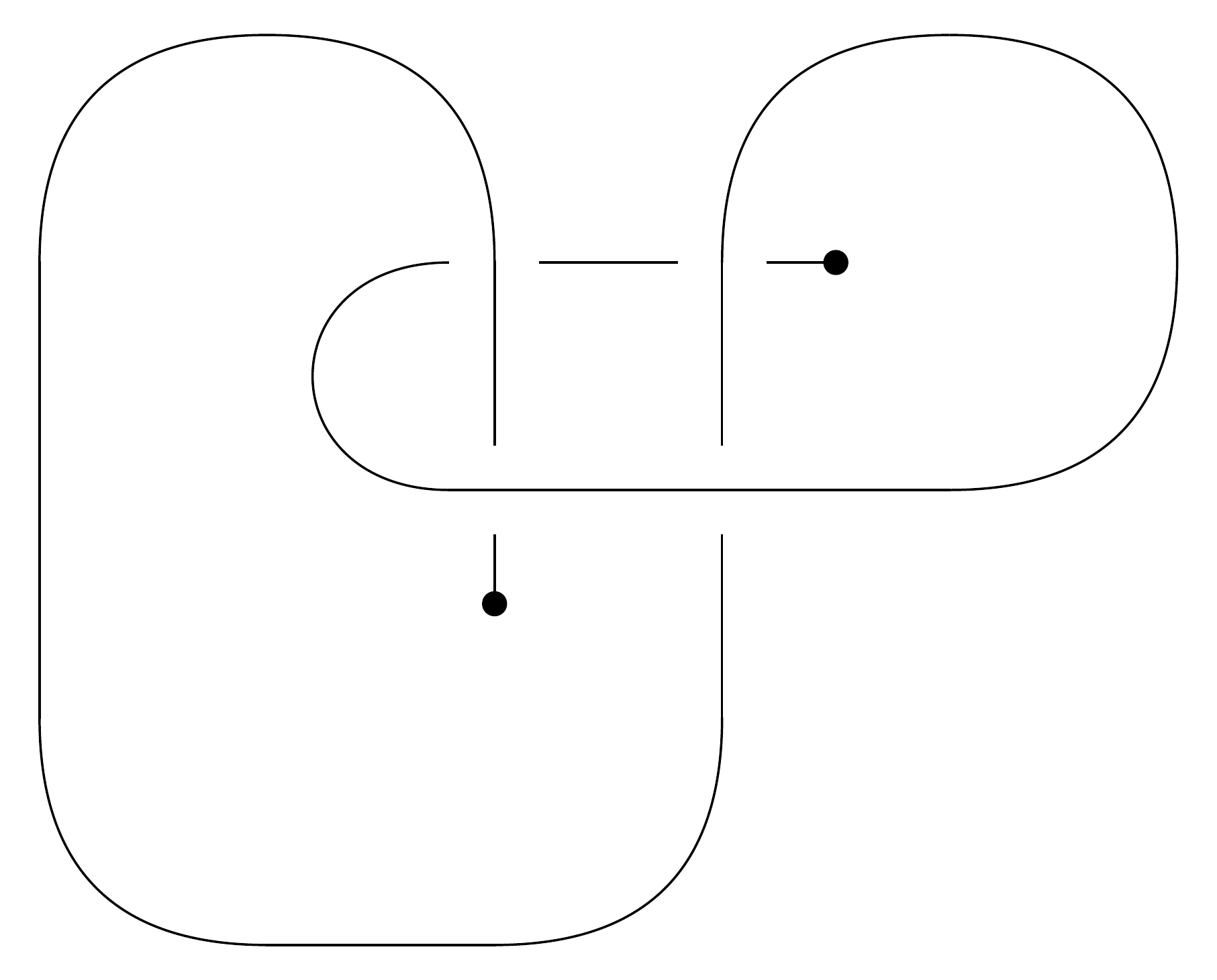}\\
\textcolor{black}{$4_{66}$}
\vspace{1cm}
\end{minipage}
\begin{minipage}[t]{.25\linewidth}
\centering
\includegraphics[width=0.9\textwidth,height=3.5cm,keepaspectratio]{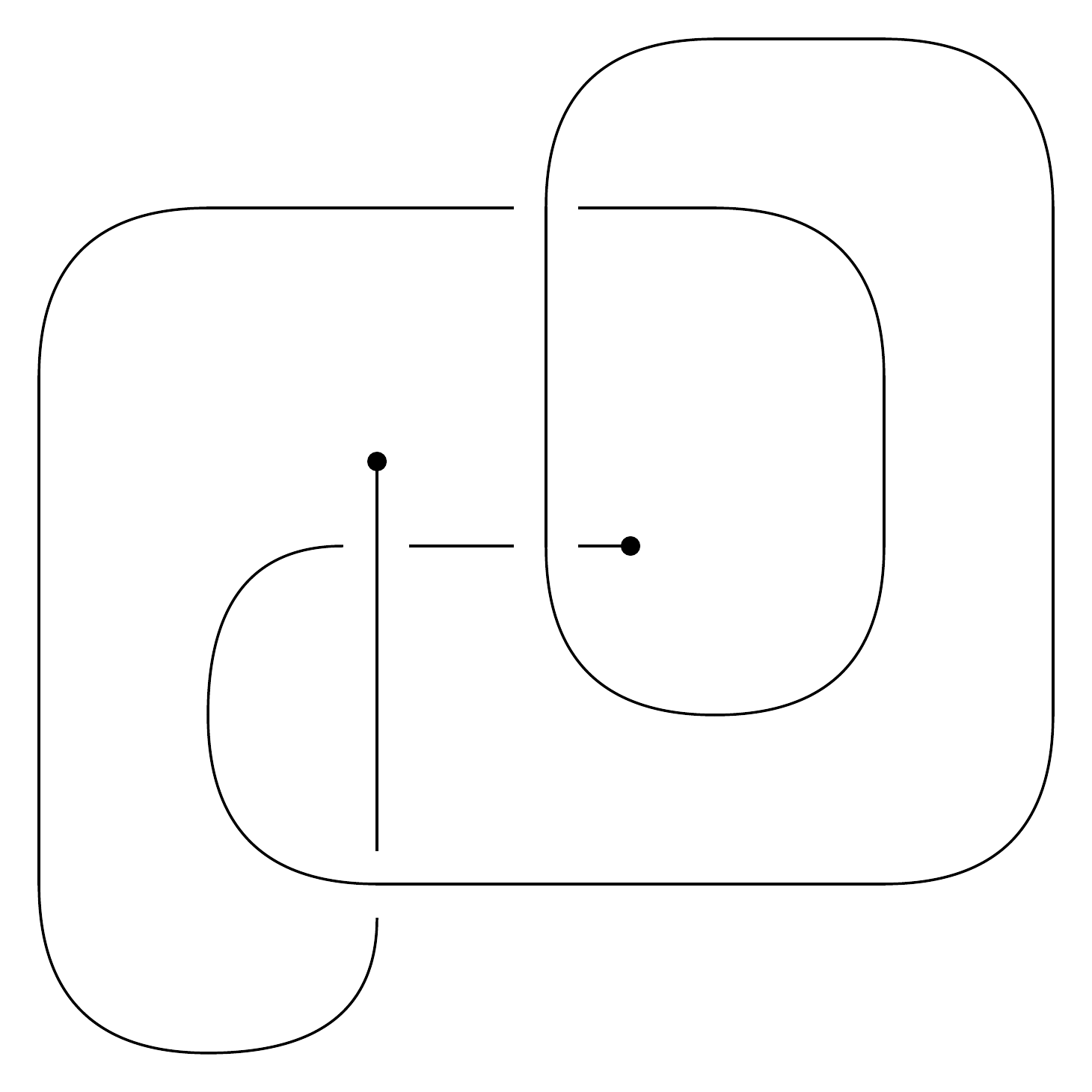}\\
\textcolor{black}{$4_{67}$}
\vspace{1cm}
\end{minipage}
\begin{minipage}[t]{.25\linewidth}
\centering
\includegraphics[width=0.9\textwidth,height=3.5cm,keepaspectratio]{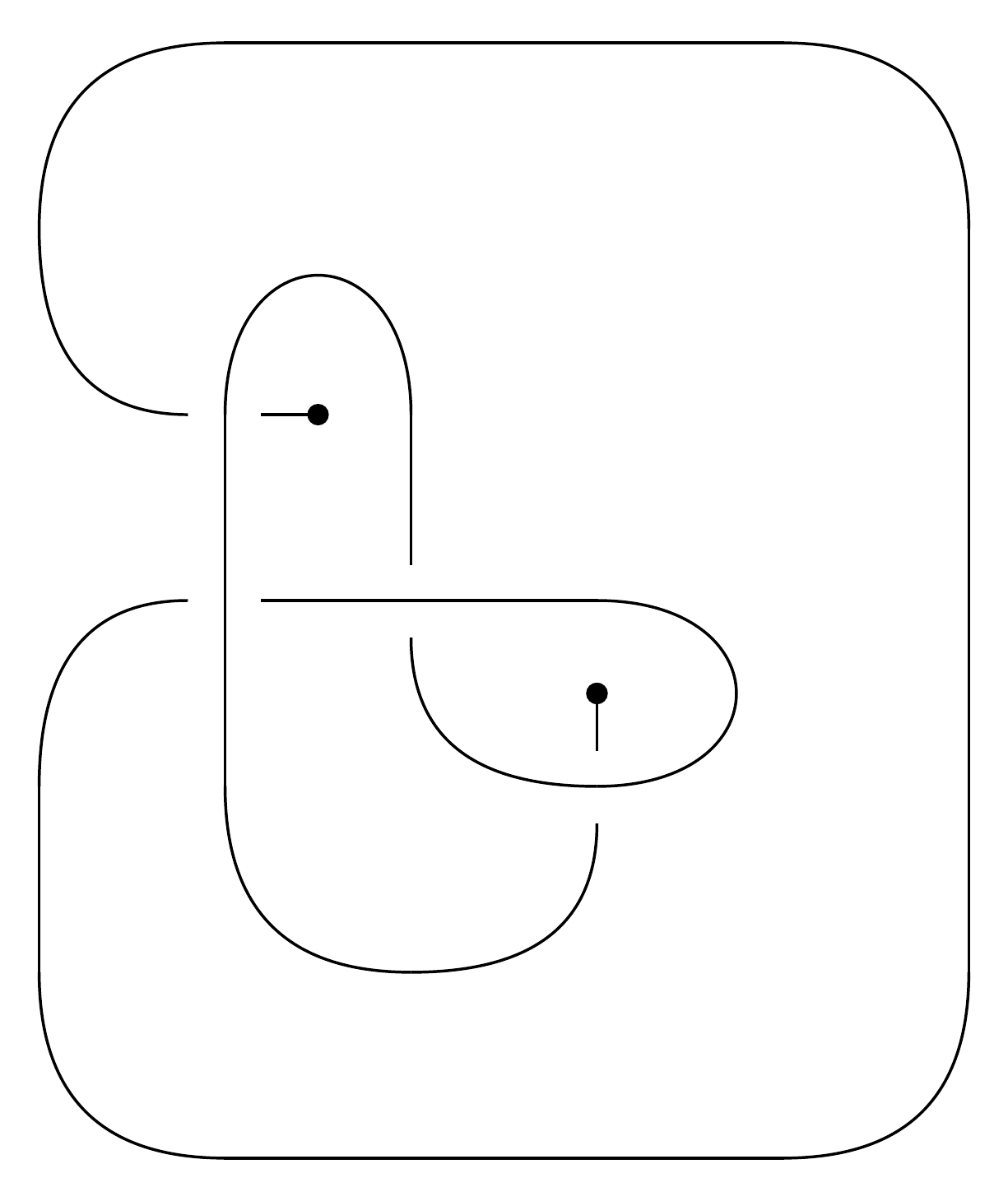}\\
\textcolor{black}{$4_{68}$}
\vspace{1cm}
\end{minipage}
\begin{minipage}[t]{.25\linewidth}
\centering
\includegraphics[width=0.9\textwidth,height=3.5cm,keepaspectratio]{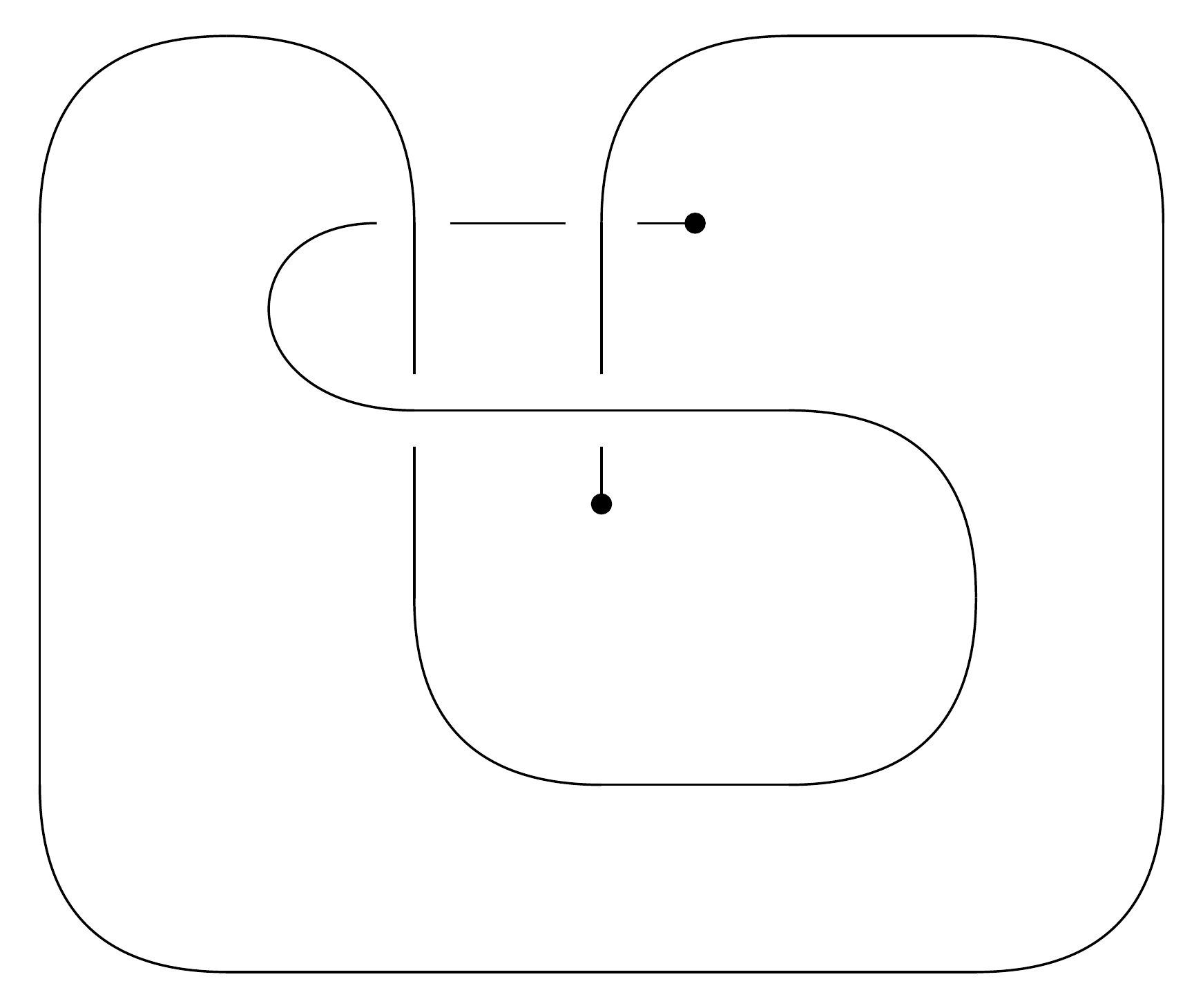}\\
\textcolor{black}{$4_{69}$}
\vspace{1cm}
\end{minipage}
\begin{minipage}[t]{.25\linewidth}
\centering
\includegraphics[width=0.9\textwidth,height=3.5cm,keepaspectratio]{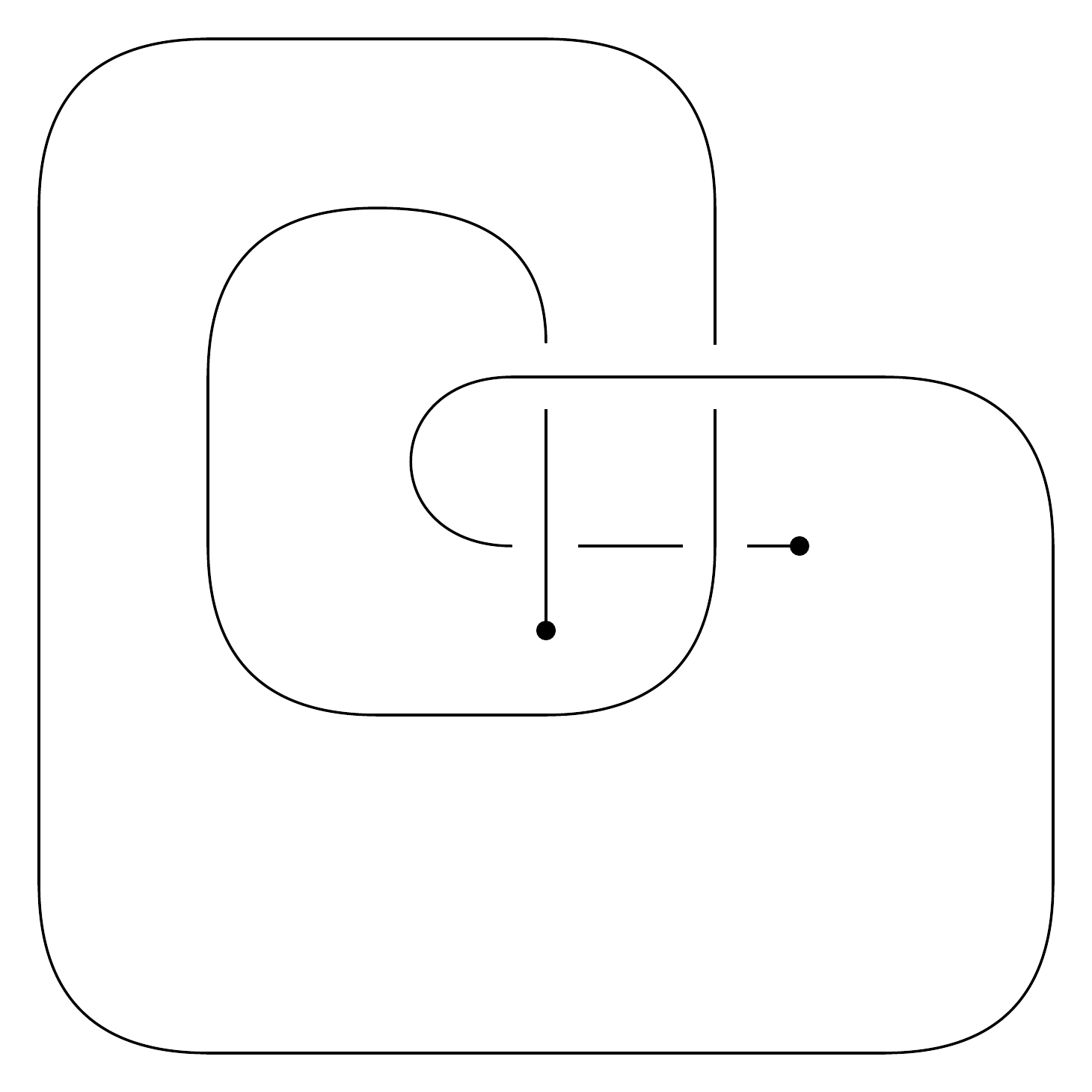}\\
\textcolor{black}{$4_{70}$}
\vspace{1cm}
\end{minipage}
\begin{minipage}[t]{.25\linewidth}
\centering
\includegraphics[width=0.9\textwidth,height=3.5cm,keepaspectratio]{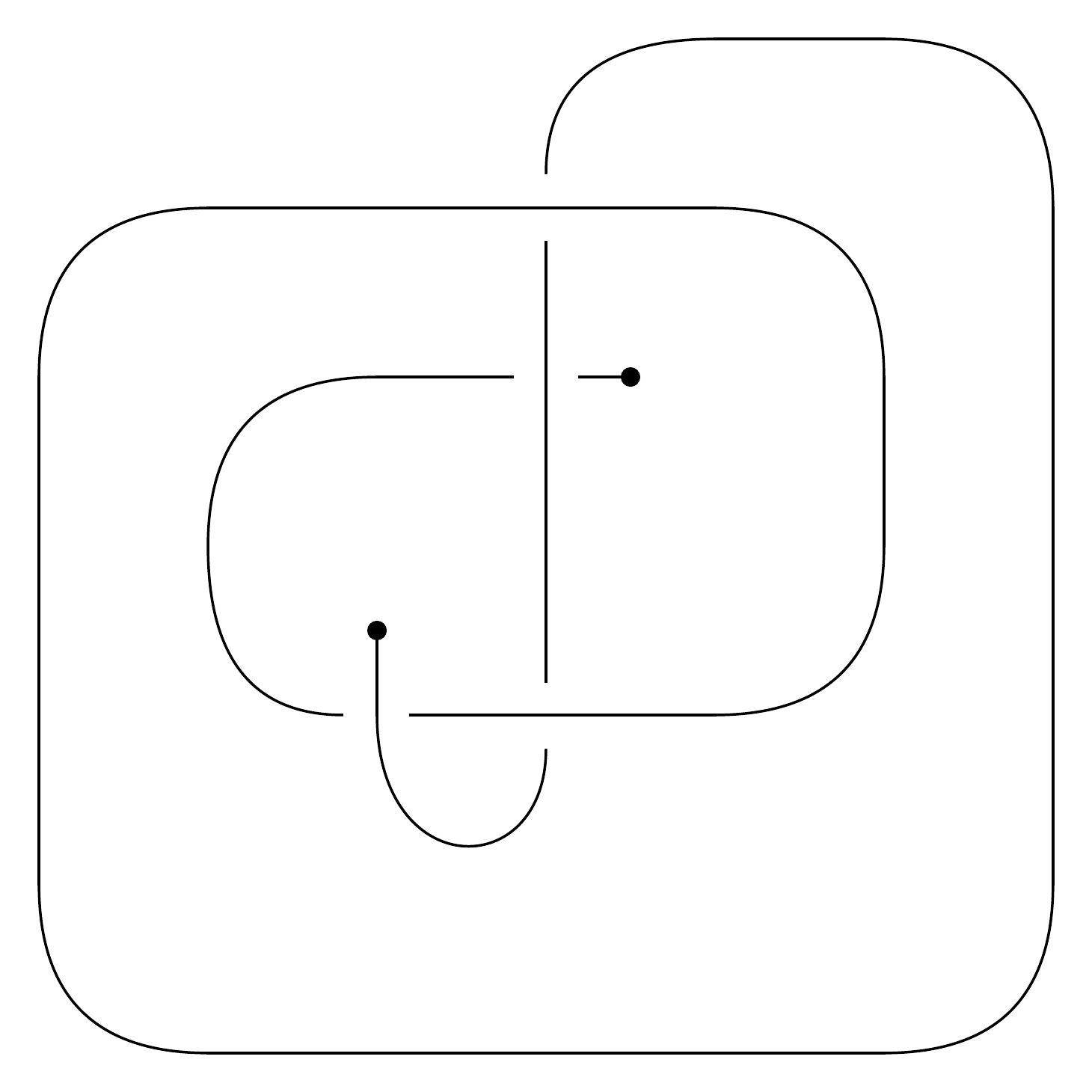}\\
\textcolor{black}{$4_{71}$}
\vspace{1cm}
\end{minipage}
\begin{minipage}[t]{.25\linewidth}
\centering
\includegraphics[width=0.9\textwidth,height=3.5cm,keepaspectratio]{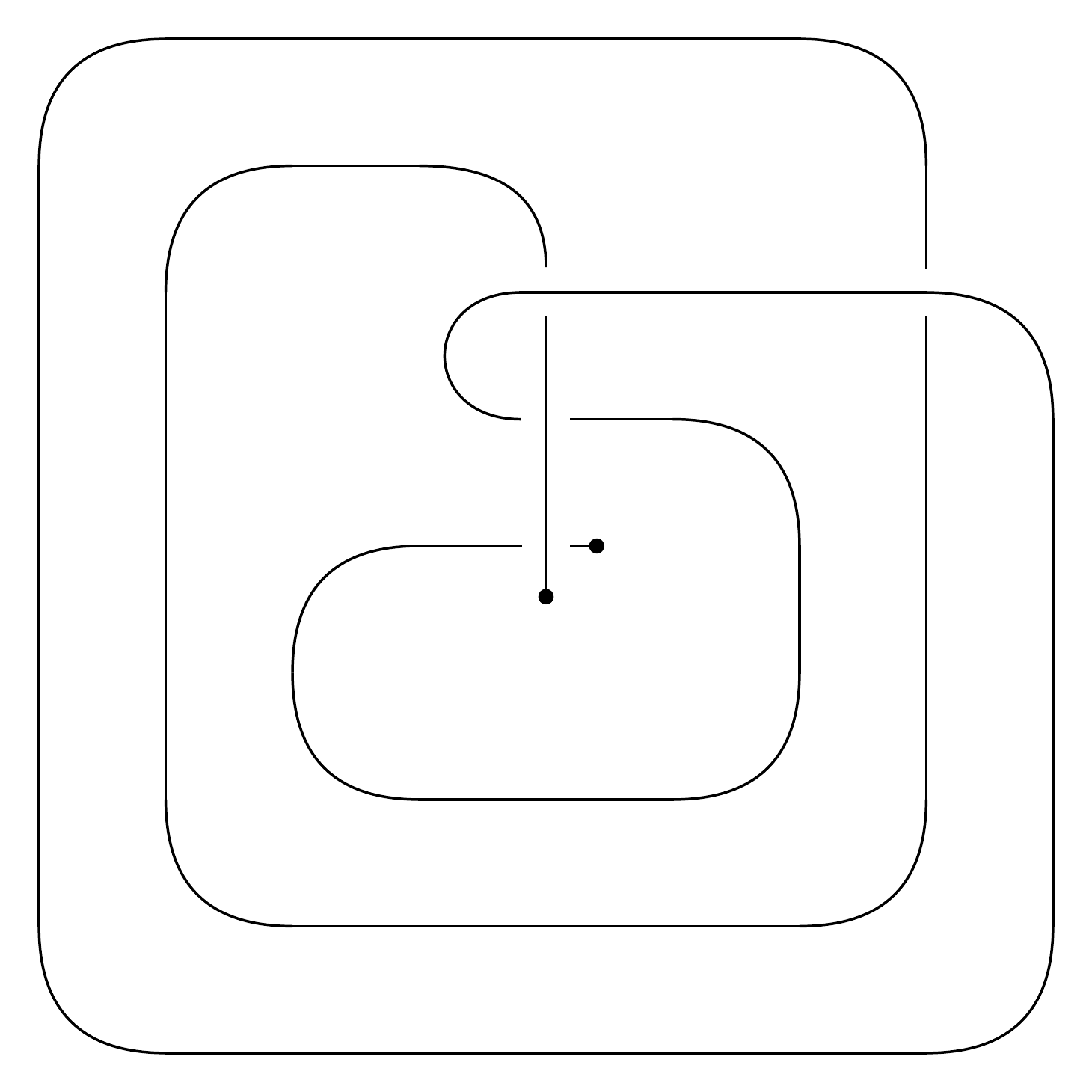}\\
\textcolor{black}{$4_{72}$}
\vspace{1cm}
\end{minipage}
\begin{minipage}[t]{.25\linewidth}
\centering
\includegraphics[width=0.9\textwidth,height=3.5cm,keepaspectratio]{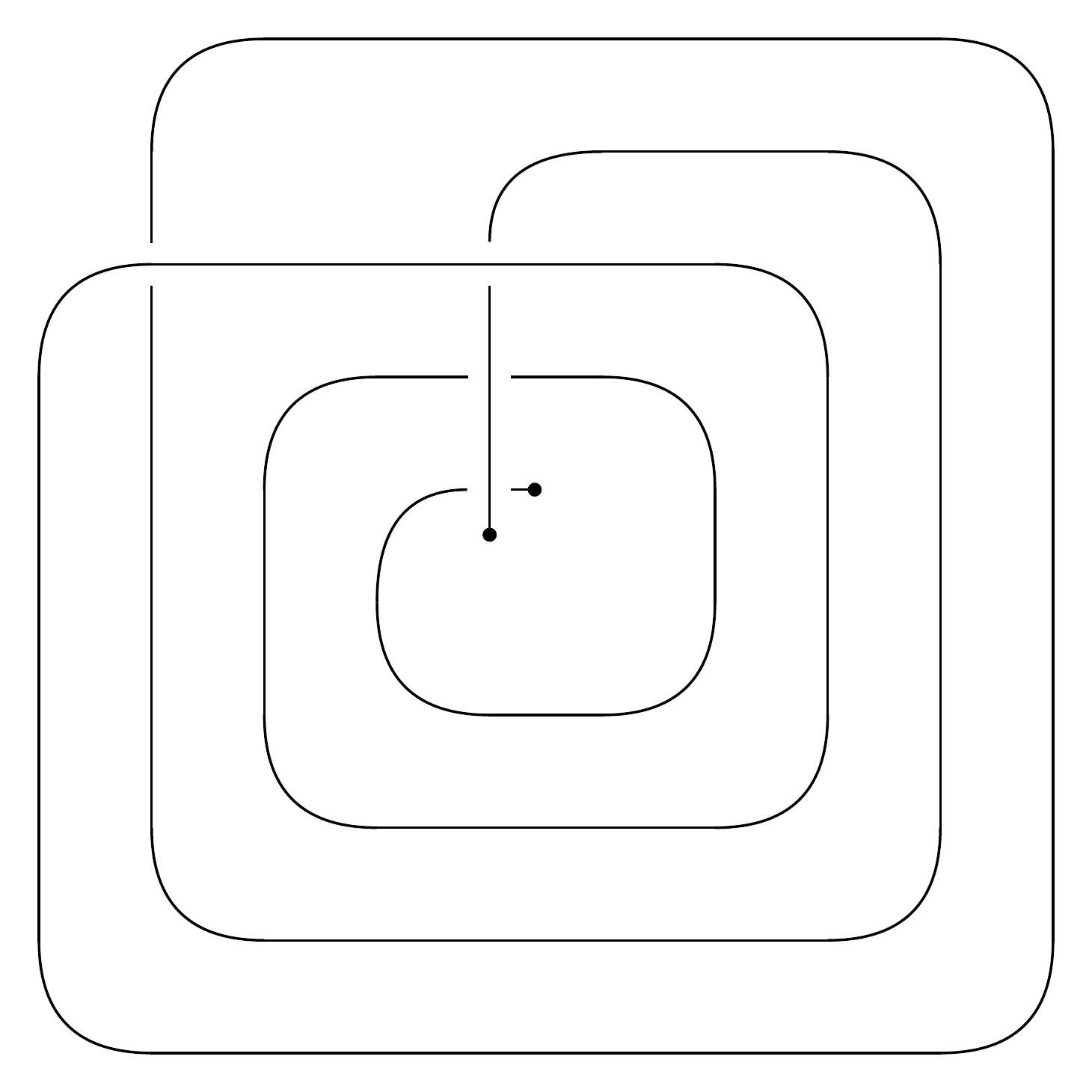}\\
\textcolor{black}{$4_{73}$}
\vspace{1cm}
\end{minipage}
\begin{minipage}[t]{.25\linewidth}
\centering
\includegraphics[width=0.9\textwidth,height=3.5cm,keepaspectratio]{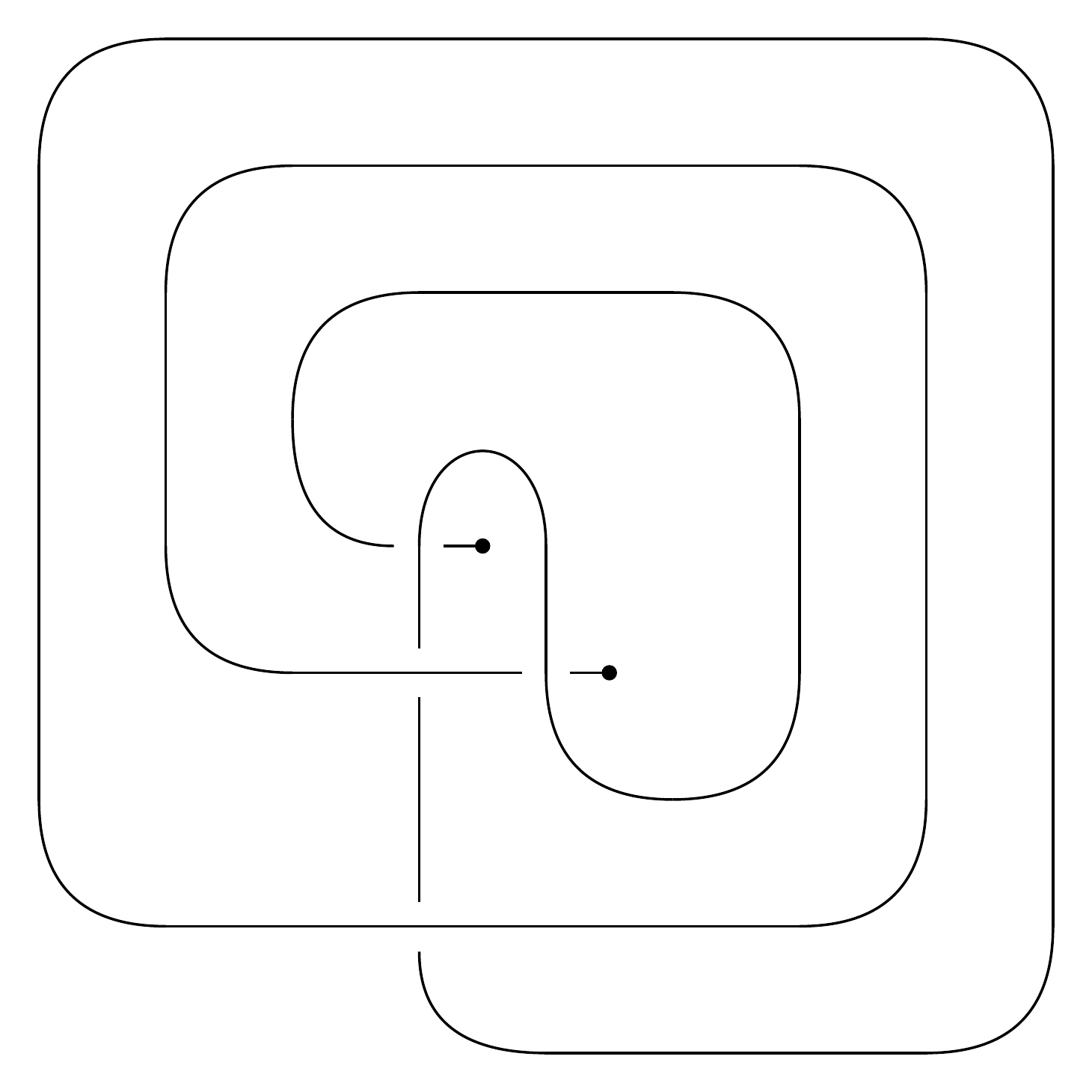}\\
\textcolor{black}{$4_{74}$}
\vspace{1cm}
\end{minipage}
\begin{minipage}[t]{.25\linewidth}
\centering
\includegraphics[width=0.9\textwidth,height=3.5cm,keepaspectratio]{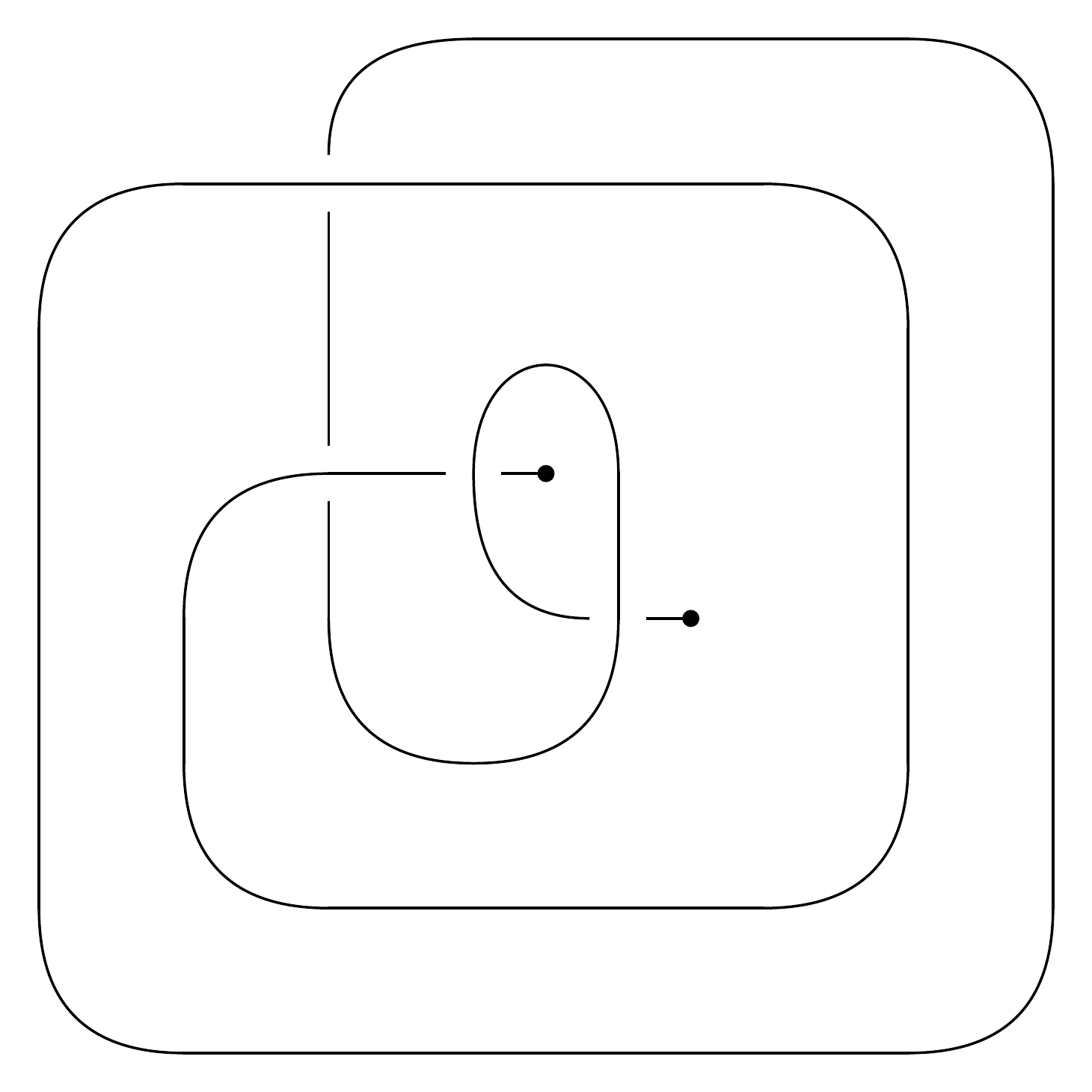}\\
\textcolor{black}{$4_{75}$}
\vspace{1cm}
\end{minipage}
\begin{minipage}[t]{.25\linewidth}
\centering
\includegraphics[width=0.9\textwidth,height=3.5cm,keepaspectratio]{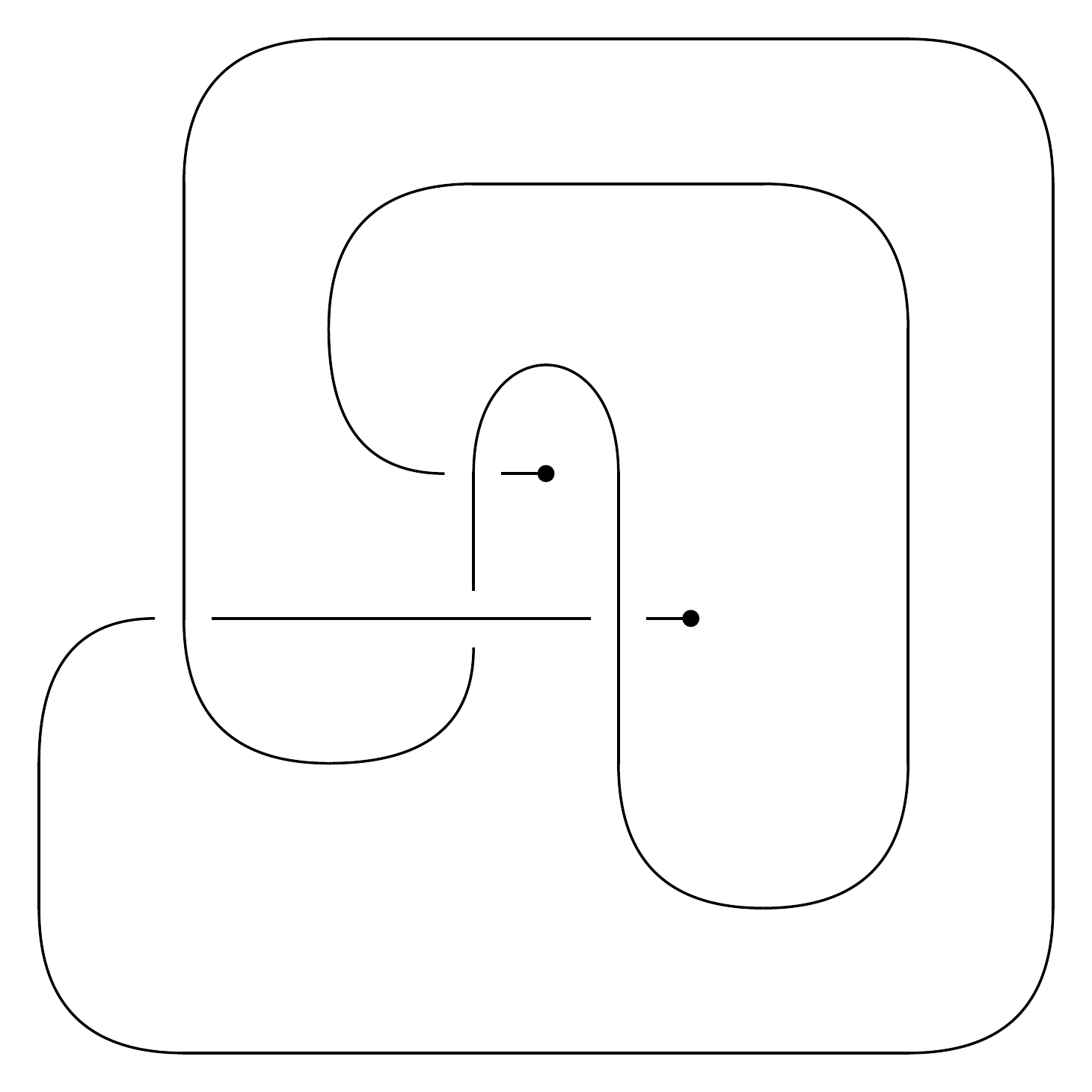}\\
\textcolor{black}{$4_{76}$}
\vspace{1cm}
\end{minipage}
\begin{minipage}[t]{.25\linewidth}
\centering
\includegraphics[width=0.9\textwidth,height=3.5cm,keepaspectratio]{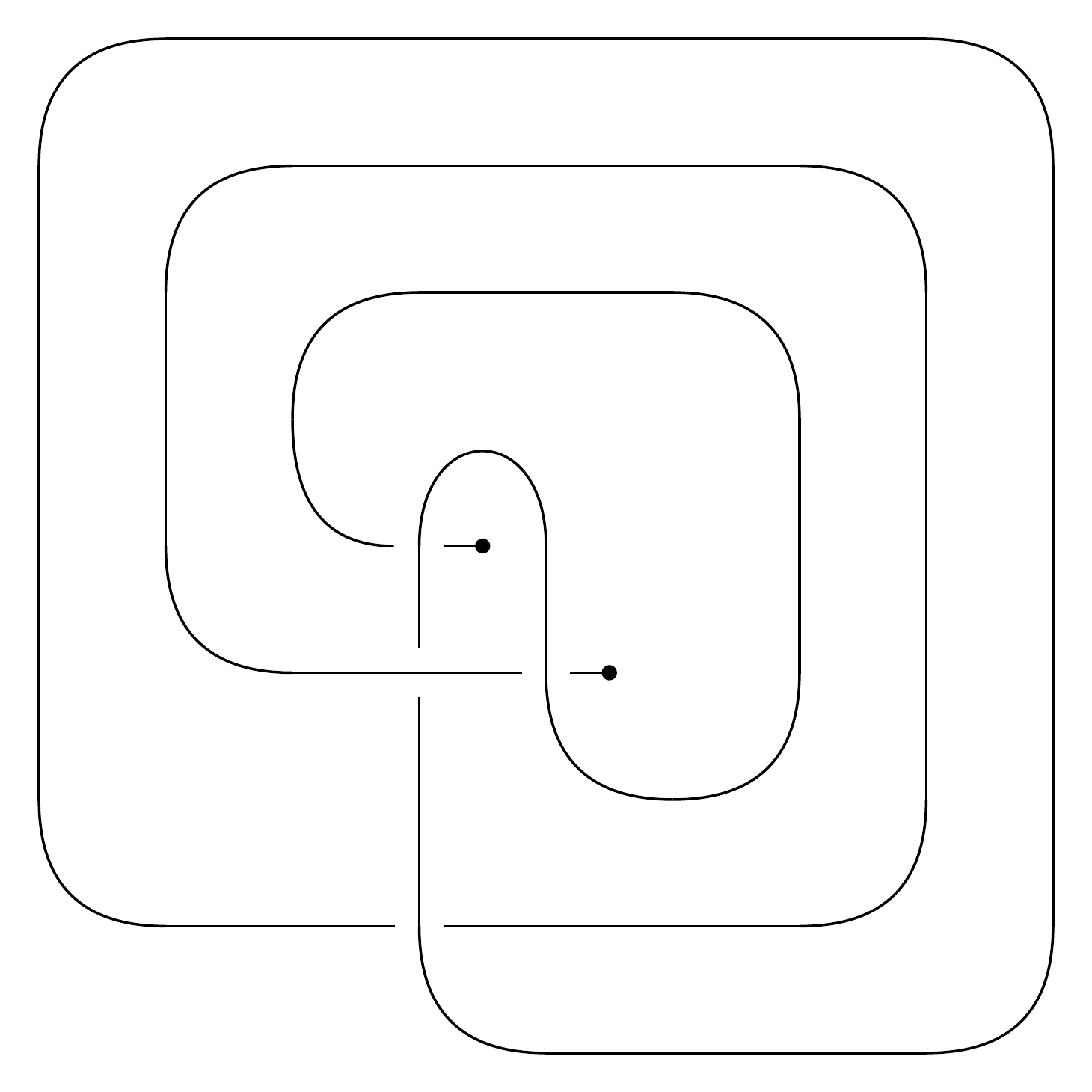}\\
\textcolor{black}{$4_{77}$}
\vspace{1cm}
\end{minipage}
\begin{minipage}[t]{.25\linewidth}
\centering
\includegraphics[width=0.9\textwidth,height=3.5cm,keepaspectratio]{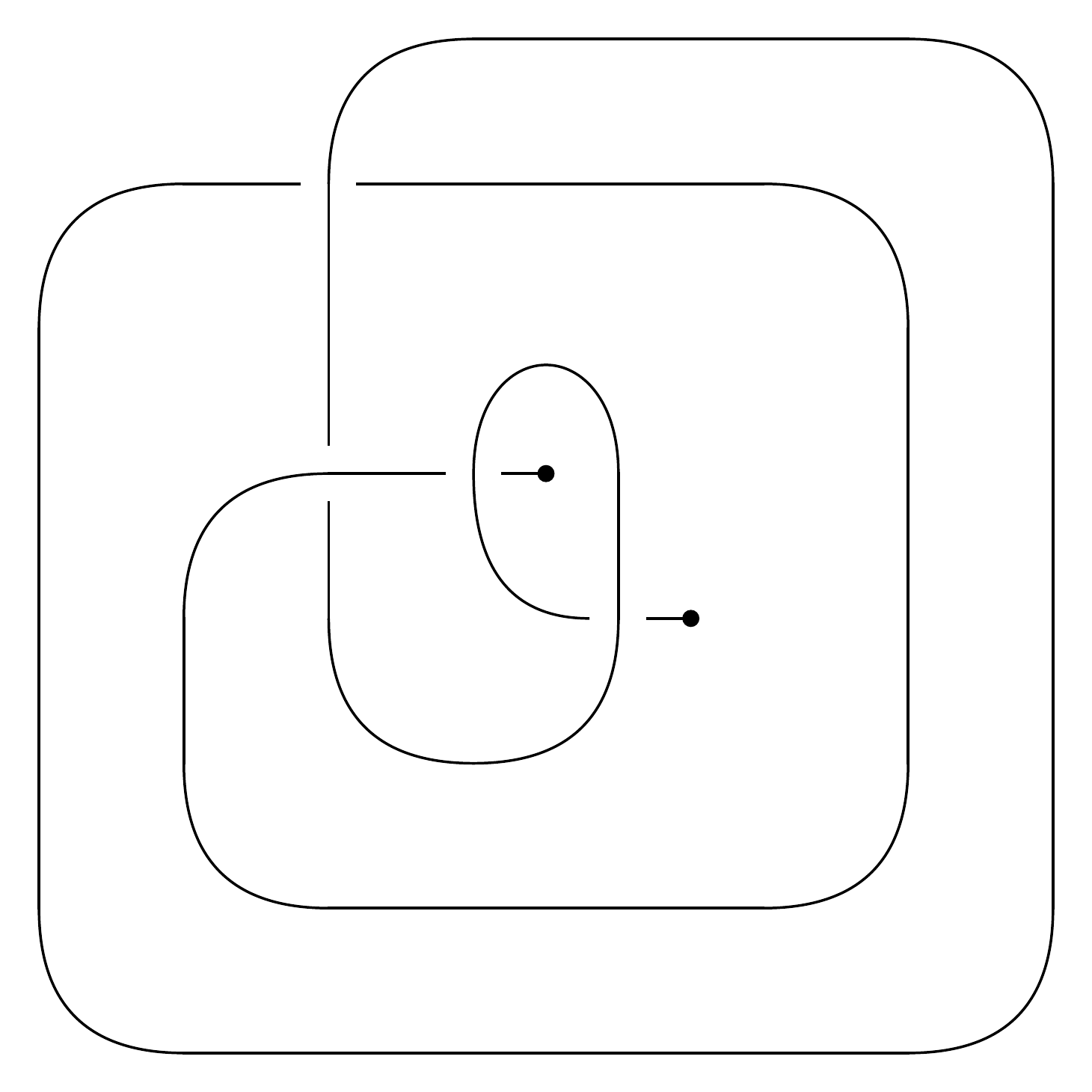}\\
\textcolor{black}{$4_{78}$}
\vspace{1cm}
\end{minipage}
\begin{minipage}[t]{.25\linewidth}
\centering
\includegraphics[width=0.9\textwidth,height=3.5cm,keepaspectratio]{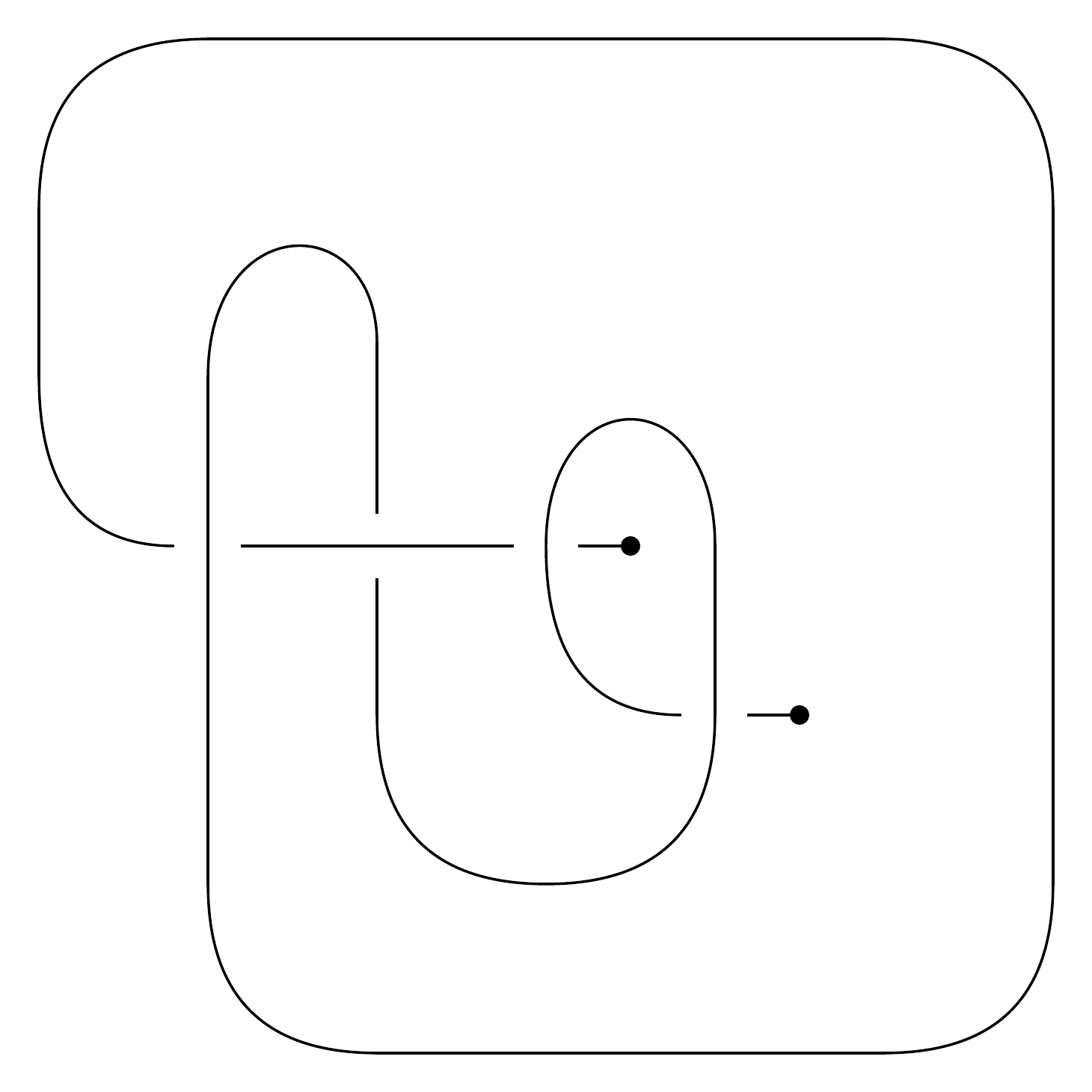}\\
\textcolor{black}{$4_{79}$}
\vspace{1cm}
\end{minipage}
\begin{minipage}[t]{.25\linewidth}
\centering
\includegraphics[width=0.9\textwidth,height=3.5cm,keepaspectratio]{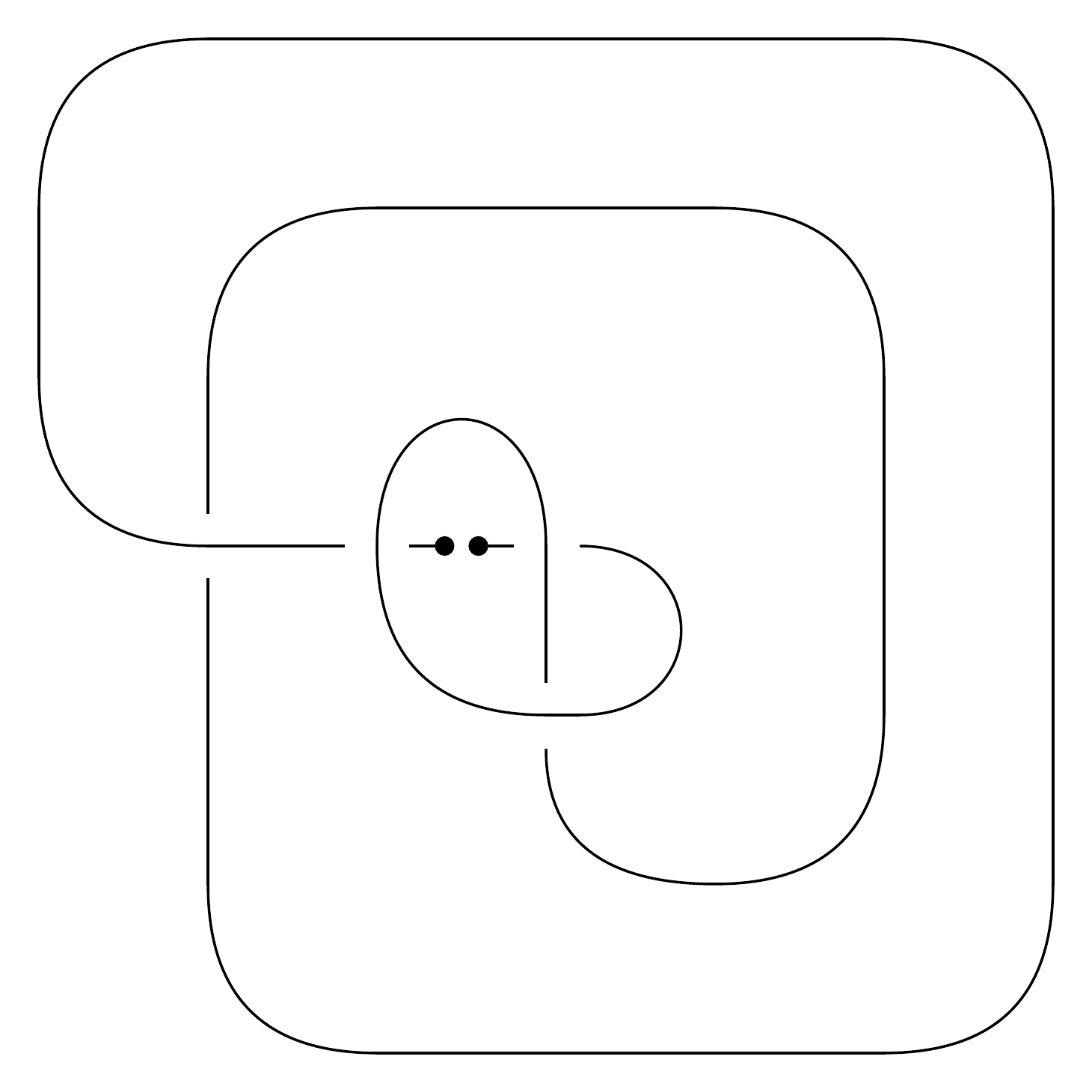}\\
\textcolor{black}{$4_{80}$}
\vspace{1cm}
\end{minipage}
\begin{minipage}[t]{.25\linewidth}
\centering
\includegraphics[width=0.9\textwidth,height=3.5cm,keepaspectratio]{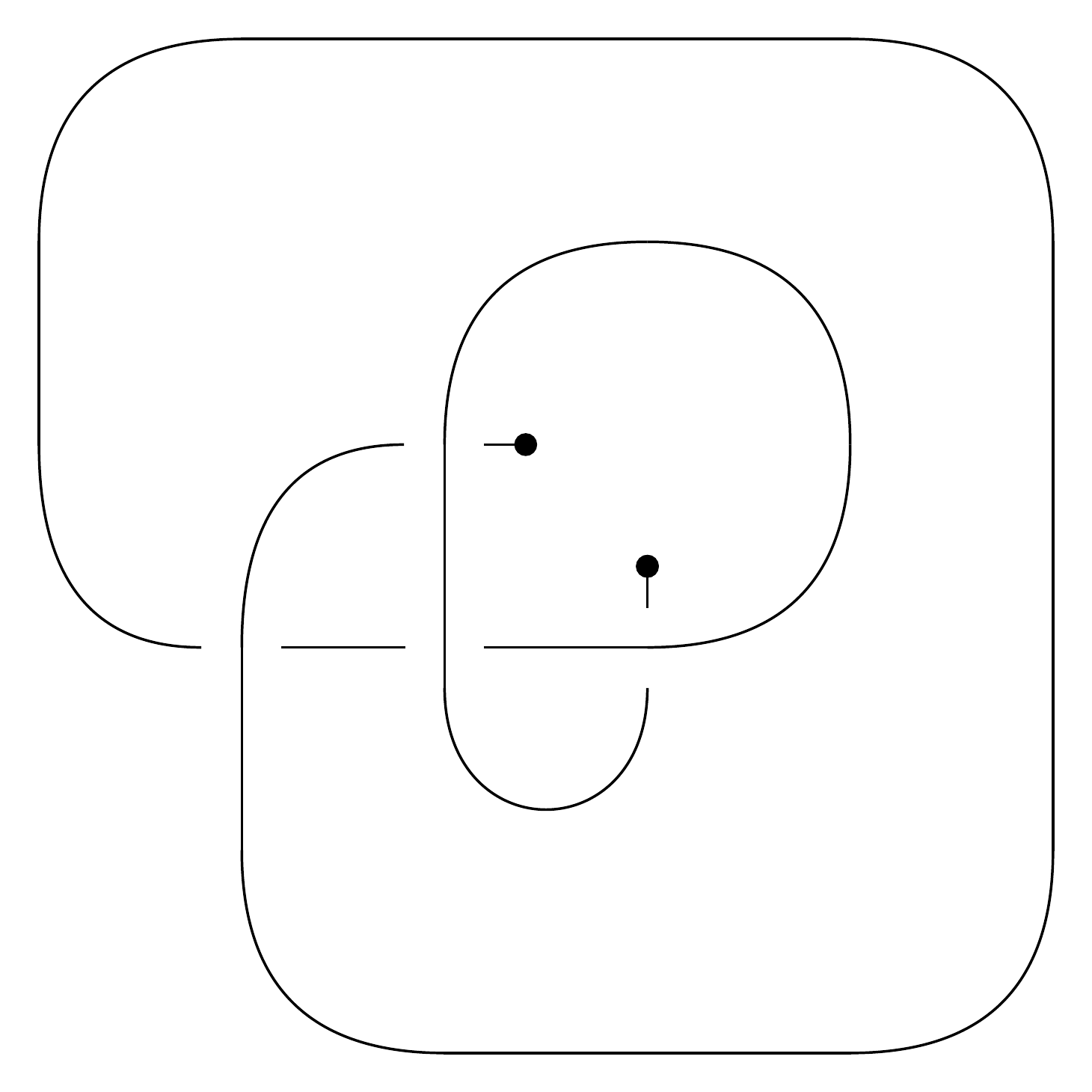}\\
\textcolor{black}{$4_{81}$}
\vspace{1cm}
\end{minipage}
\begin{minipage}[t]{.25\linewidth}
\centering
\includegraphics[width=0.9\textwidth,height=3.5cm,keepaspectratio]{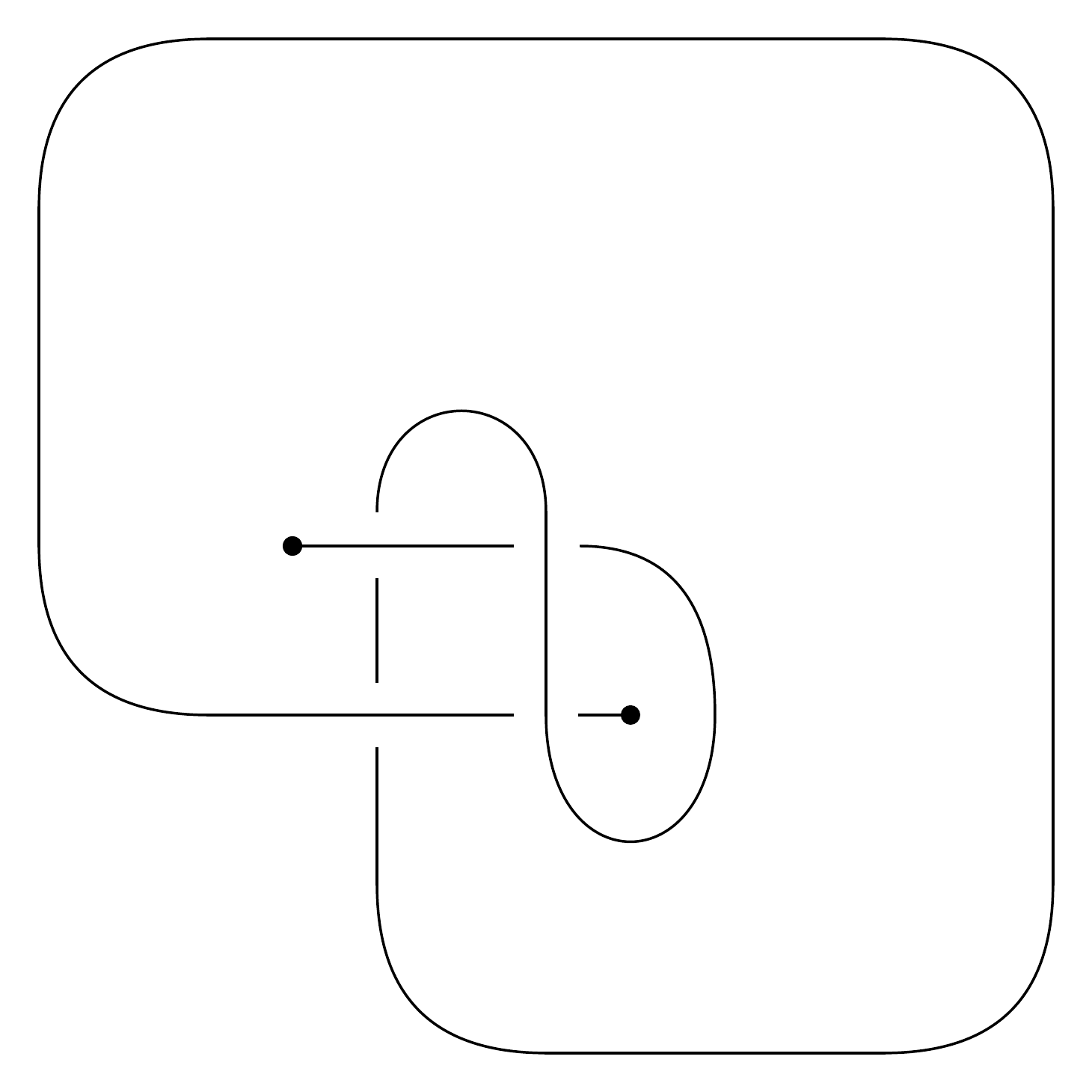}\\
\textcolor{black}{$4_{82}$}
\vspace{1cm}
\end{minipage}
\begin{minipage}[t]{.25\linewidth}
\centering
\includegraphics[width=0.9\textwidth,height=3.5cm,keepaspectratio]{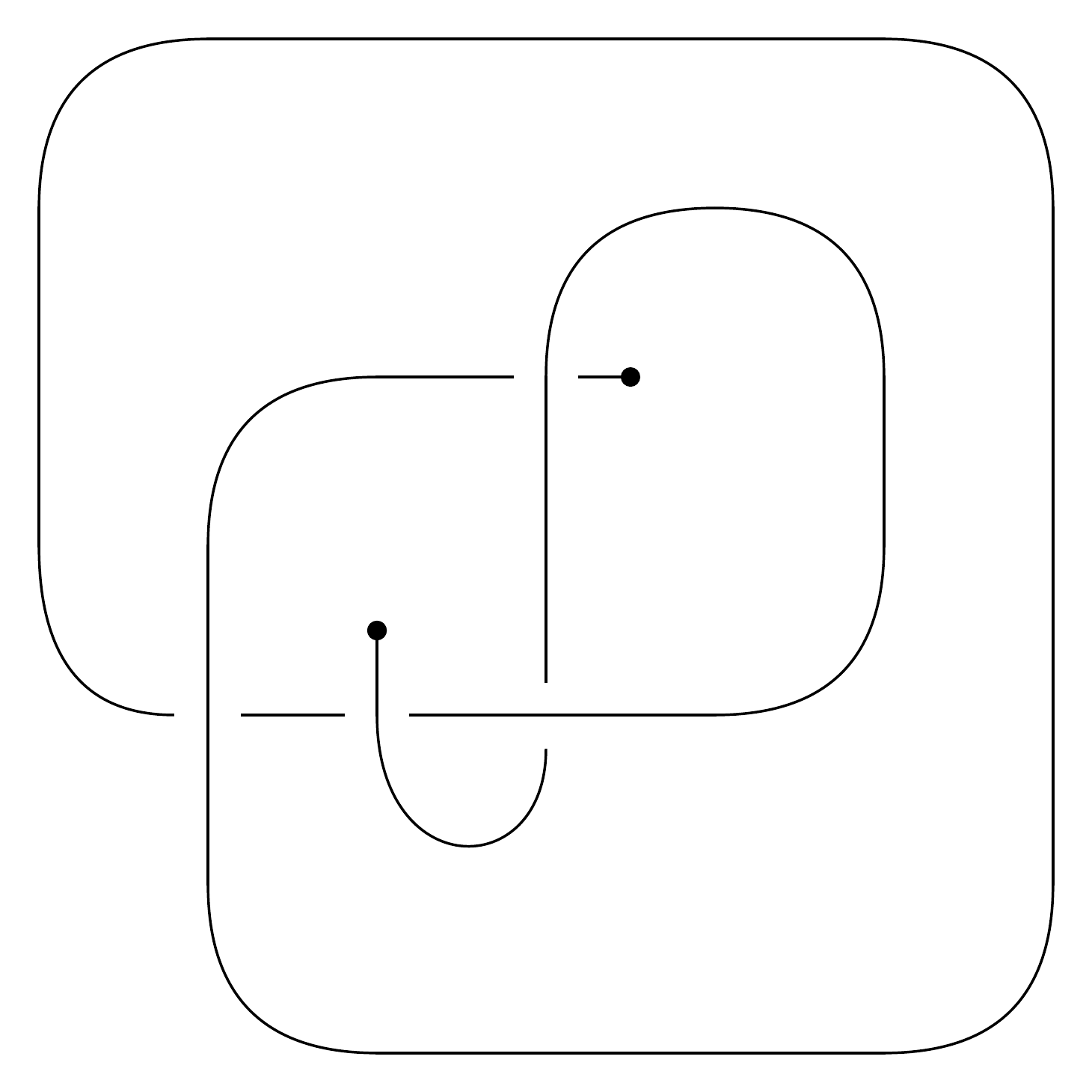}\\
\textcolor{black}{$4_{83}$}
\vspace{1cm}
\end{minipage}
\begin{minipage}[t]{.25\linewidth}
\centering
\includegraphics[width=0.9\textwidth,height=3.5cm,keepaspectratio]{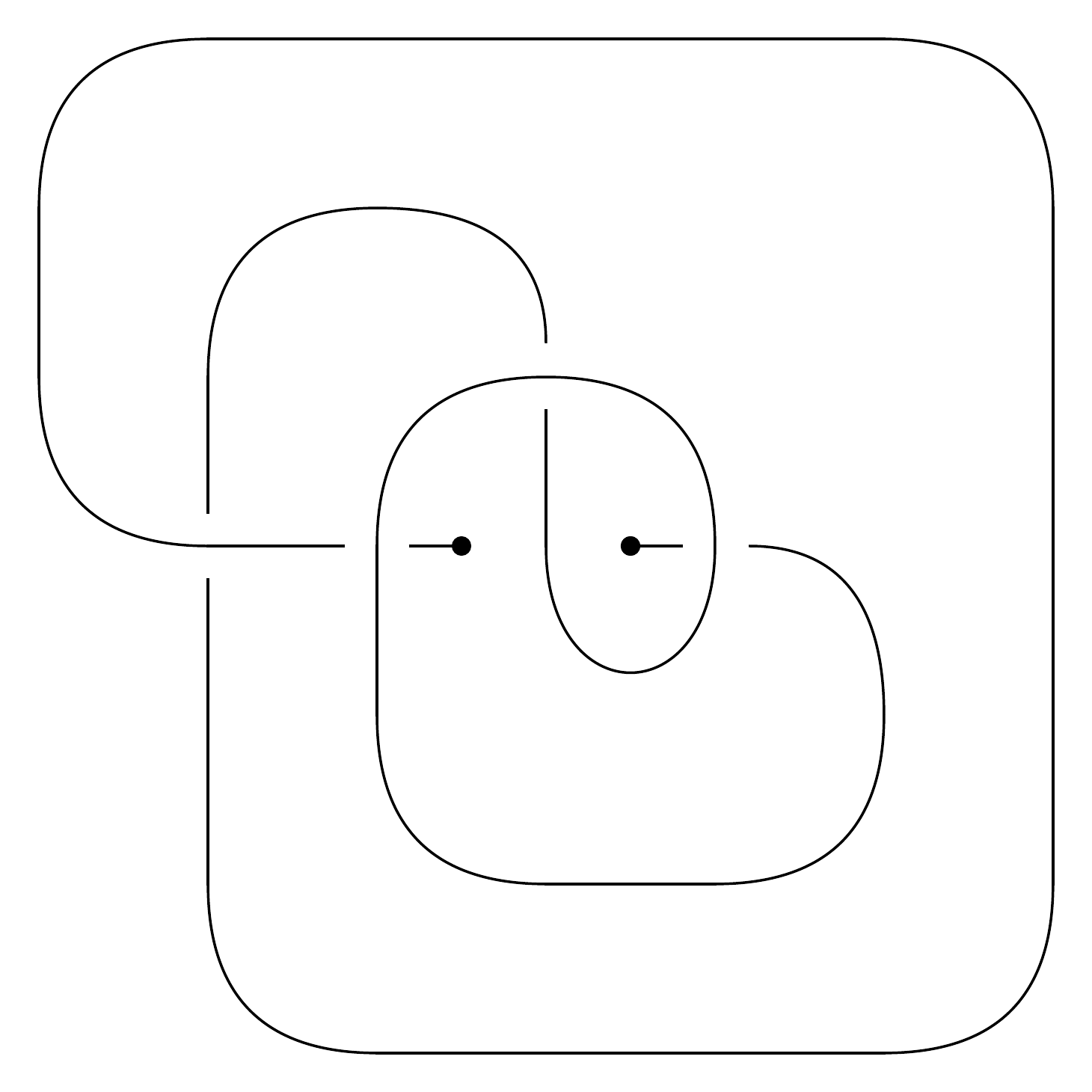}\\
\textcolor{black}{$4_{84}$}
\vspace{1cm}
\end{minipage}
\begin{minipage}[t]{.25\linewidth}
\centering
\includegraphics[width=0.9\textwidth,height=3.5cm,keepaspectratio]{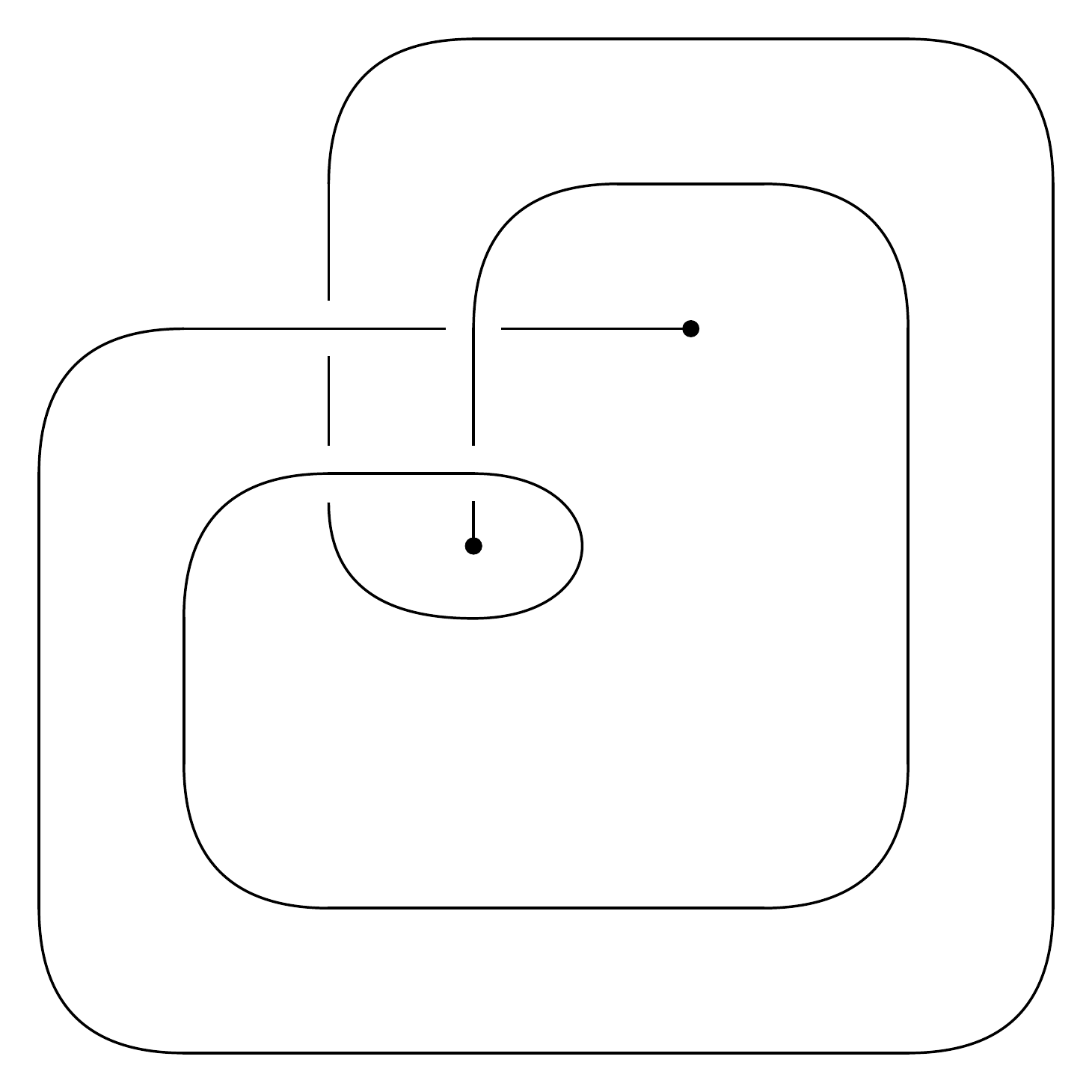}\\
\textcolor{black}{$4_{85}$}
\vspace{1cm}
\end{minipage}
\begin{minipage}[t]{.25\linewidth}
\centering
\includegraphics[width=0.9\textwidth,height=3.5cm,keepaspectratio]{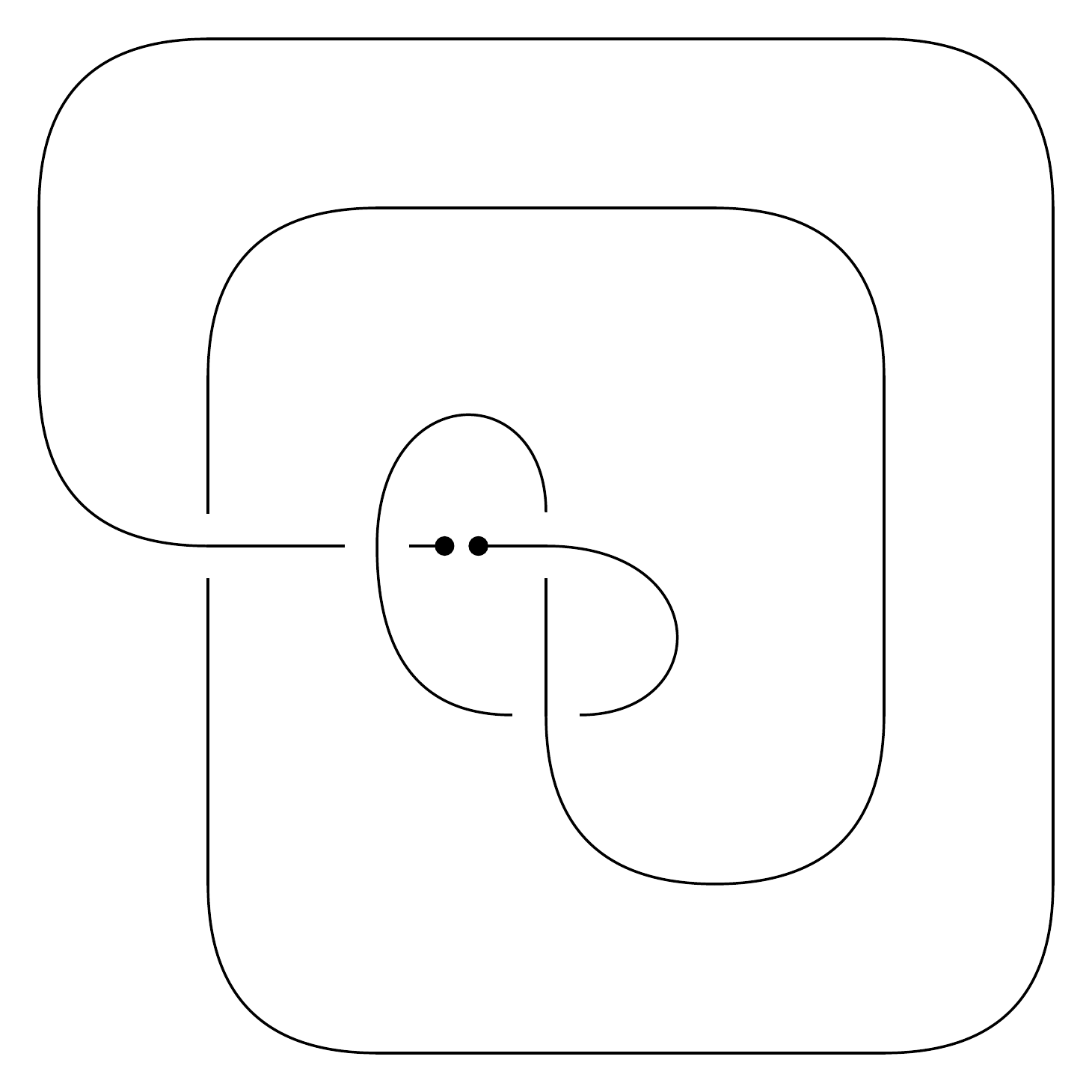}\\
\textcolor{black}{$4_{86}$}
\vspace{1cm}
\end{minipage}
\begin{minipage}[t]{.25\linewidth}
\centering
\includegraphics[width=0.9\textwidth,height=3.5cm,keepaspectratio]{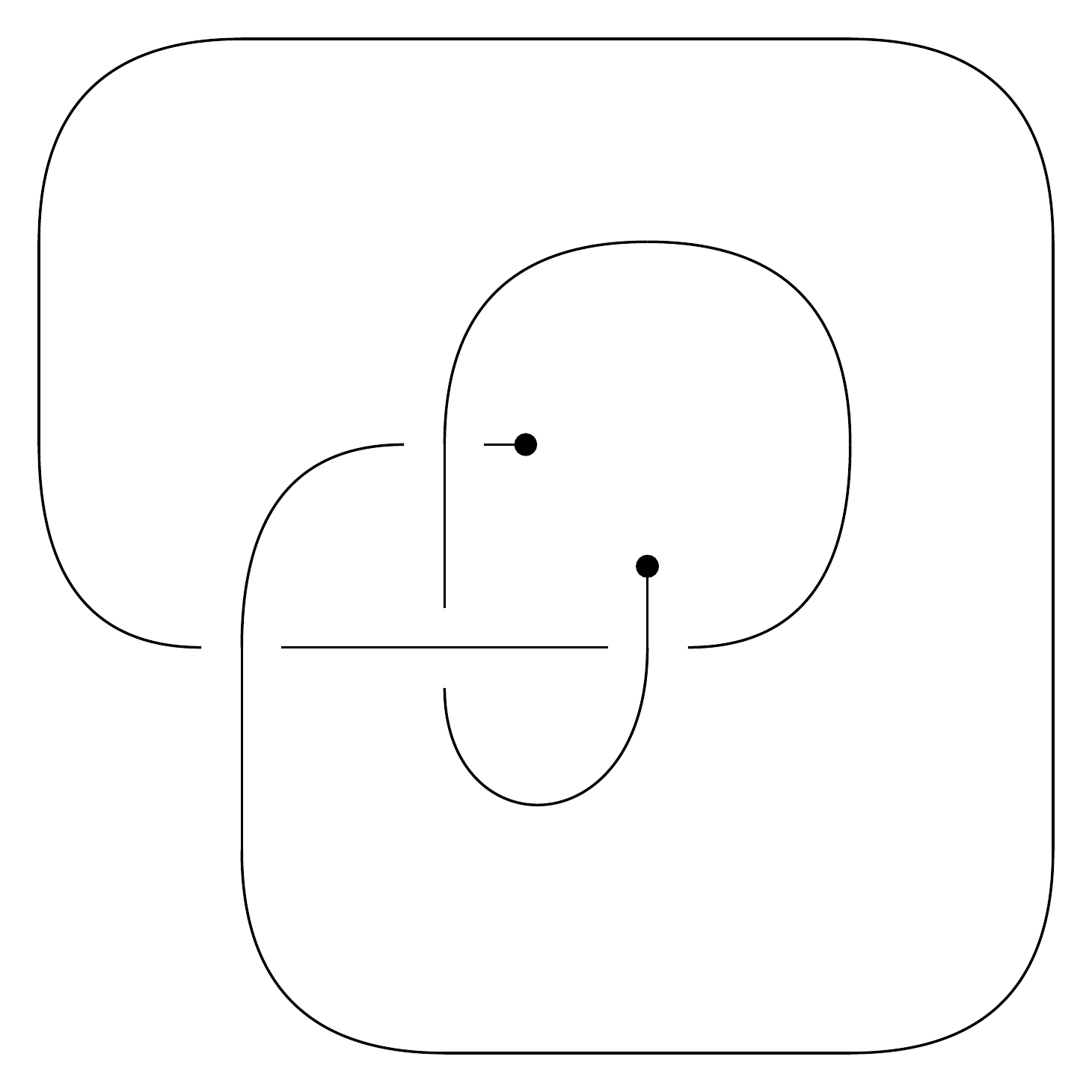}\\
\textcolor{black}{$4_{87}$}
\vspace{1cm}
\end{minipage}
\begin{minipage}[t]{.25\linewidth}
\centering
\includegraphics[width=0.9\textwidth,height=3.5cm,keepaspectratio]{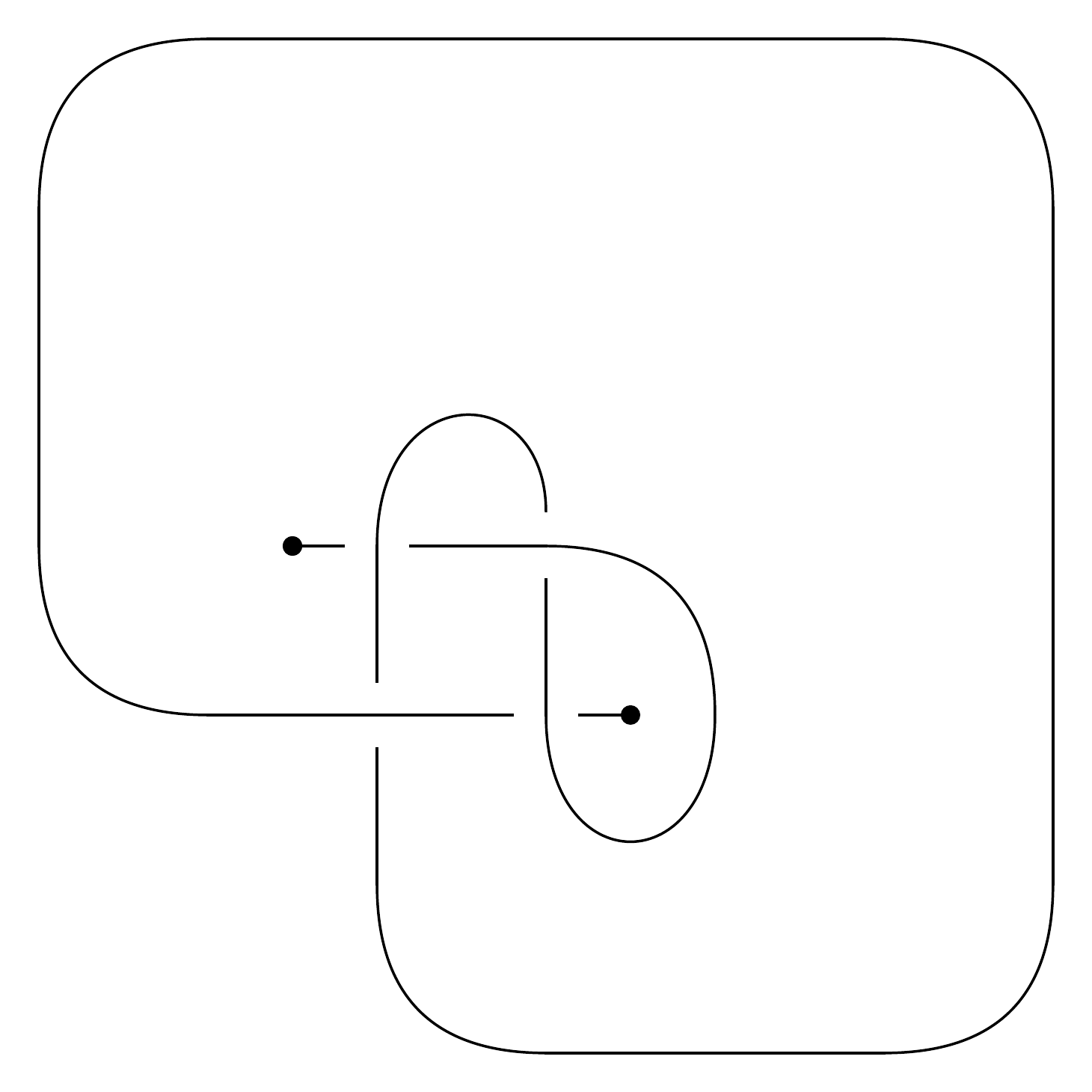}\\
\textcolor{black}{$4_{88}$}
\vspace{1cm}
\end{minipage}
\begin{minipage}[t]{.25\linewidth}
\centering
\includegraphics[width=0.9\textwidth,height=3.5cm,keepaspectratio]{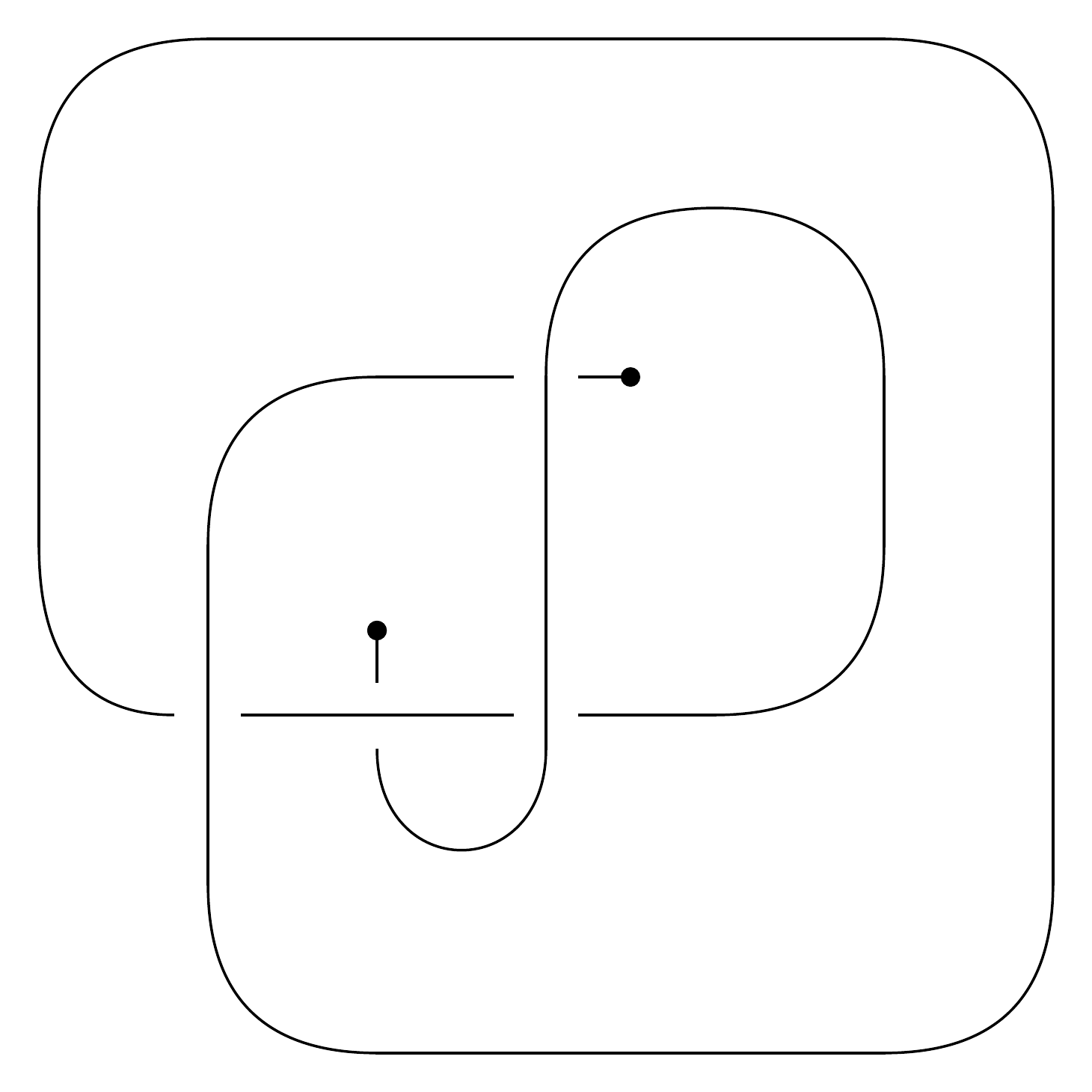}\\
\textcolor{black}{$4_{89}$}
\vspace{1cm}
\end{minipage}
\begin{minipage}[t]{.25\linewidth}
\centering
\includegraphics[width=0.9\textwidth,height=3.5cm,keepaspectratio]{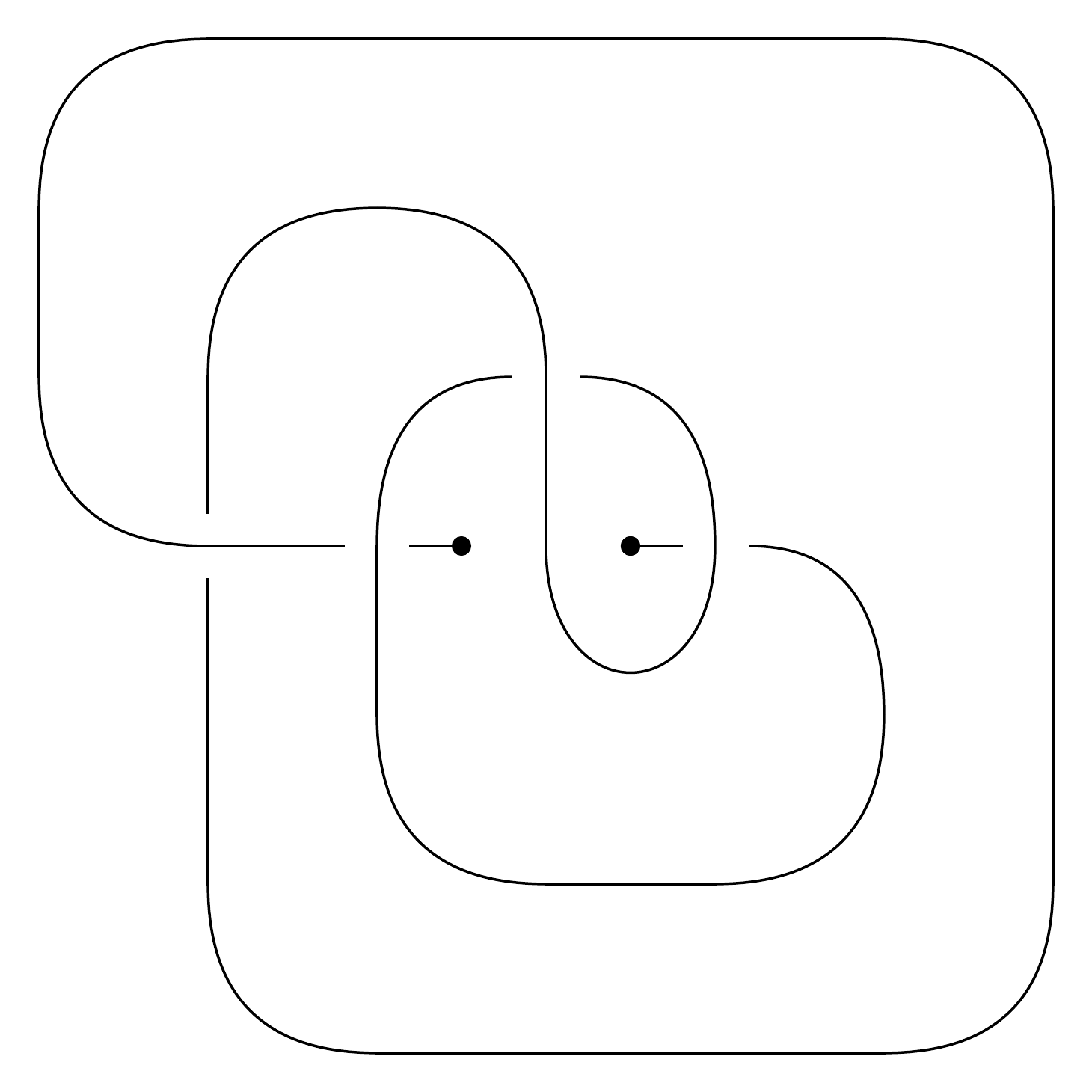}\\
\textcolor{black}{$4_{90}$}
\vspace{1cm}
\end{minipage}
\begin{minipage}[t]{.25\linewidth}
\centering
\includegraphics[width=0.9\textwidth,height=3.5cm,keepaspectratio]{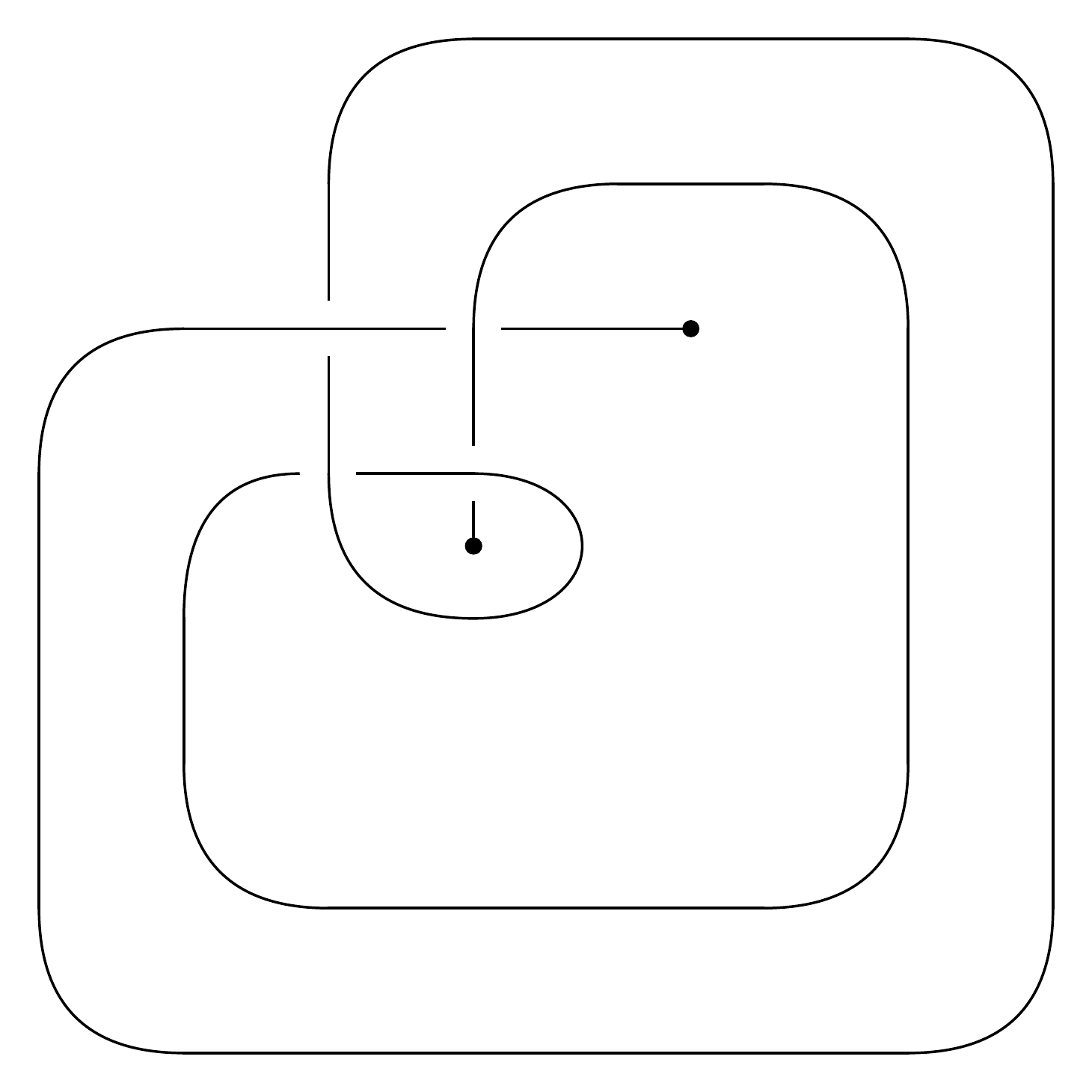}\\
\textcolor{black}{$4_{91}$}
\vspace{1cm}
\end{minipage}
\begin{minipage}[t]{.25\linewidth}
\centering
\includegraphics[width=0.9\textwidth,height=3.5cm,keepaspectratio]{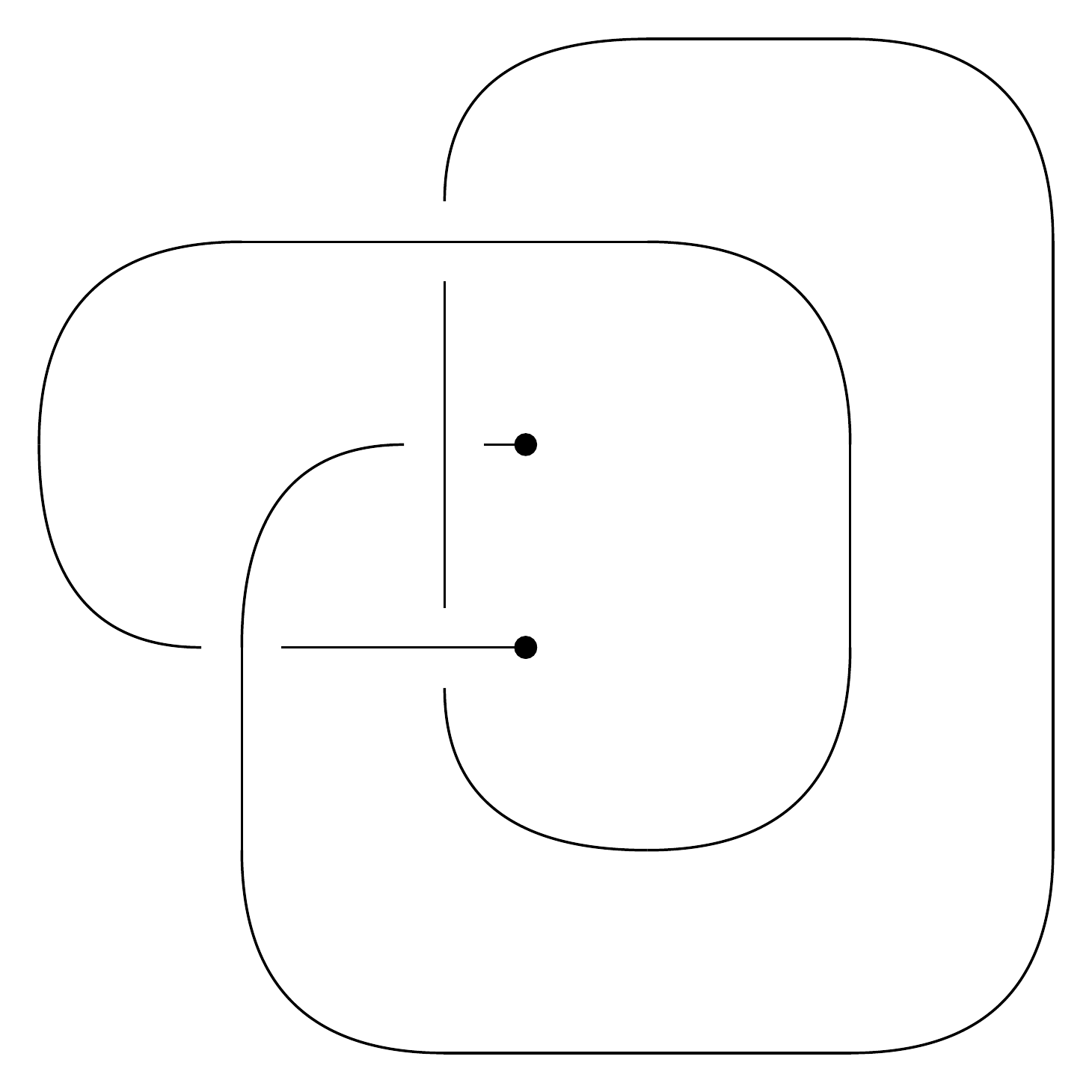}\\
\textcolor{black}{$4_{92}$}
\vspace{1cm}
\end{minipage}
\begin{minipage}[t]{.25\linewidth}
\centering
\includegraphics[width=0.9\textwidth,height=3.5cm,keepaspectratio]{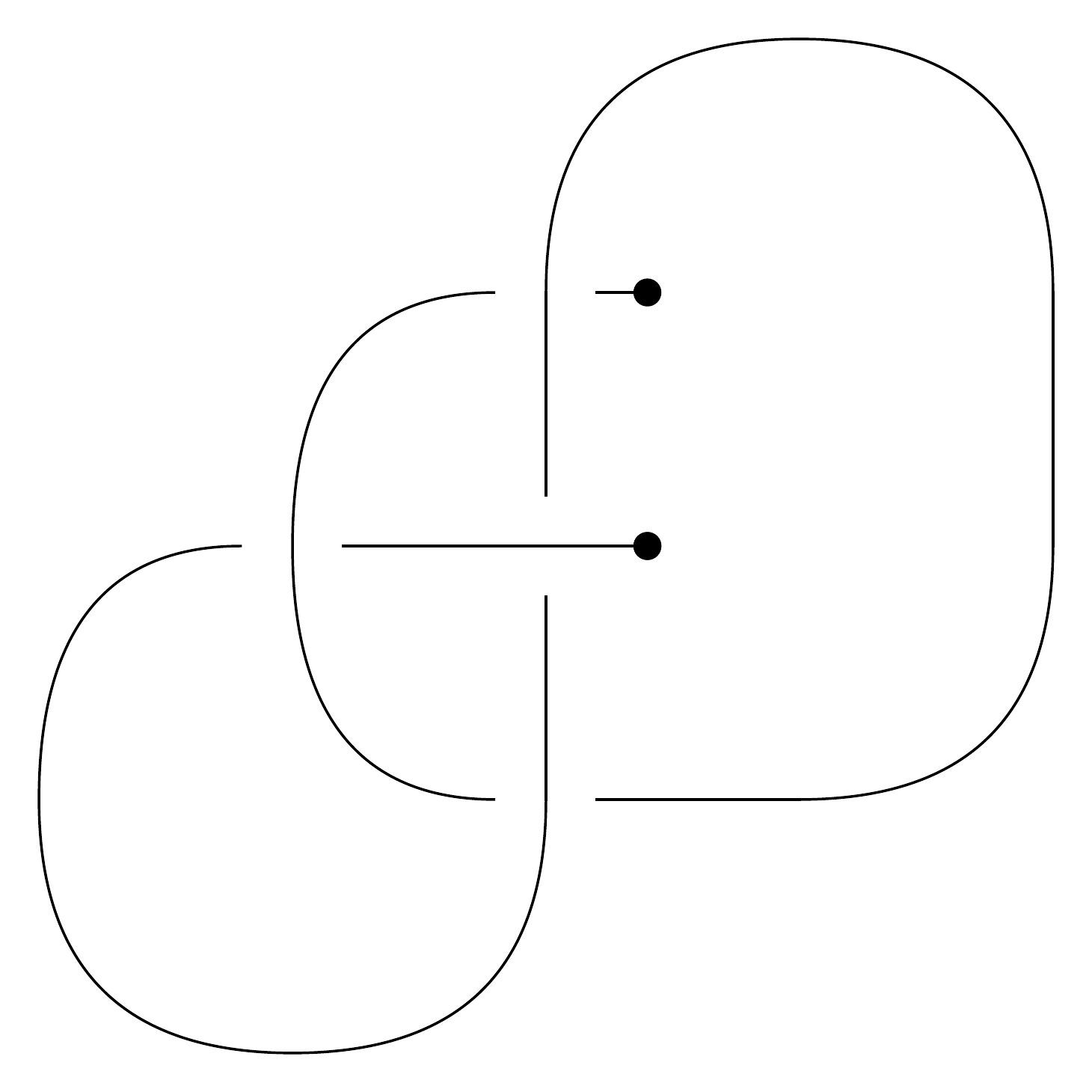}\\
\textcolor{black}{$4_{93}$}
\vspace{1cm}
\end{minipage}
\begin{minipage}[t]{.25\linewidth}
\centering
\includegraphics[width=0.9\textwidth,height=3.5cm,keepaspectratio]{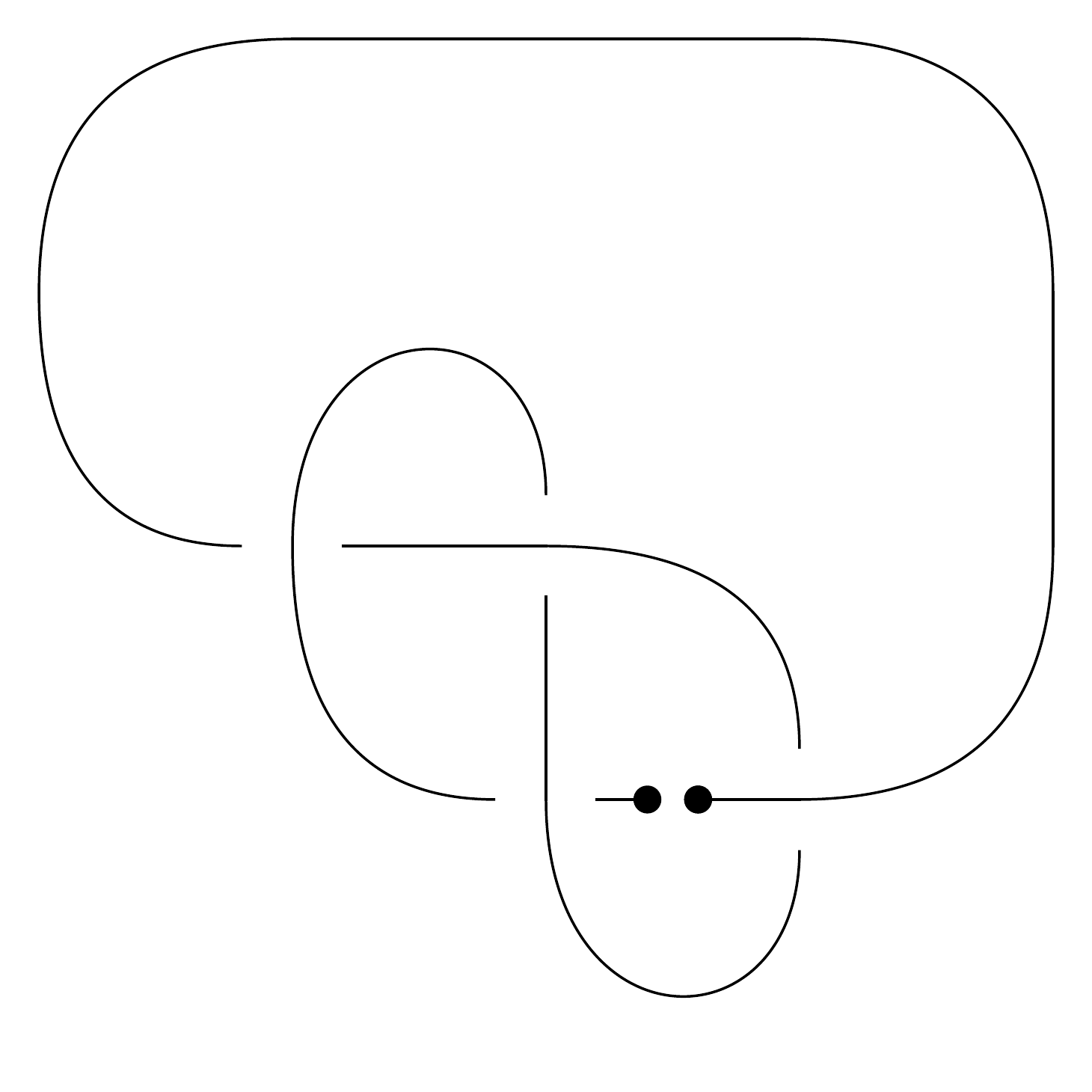}\\
\textcolor{black}{$4_{94}$}
\vspace{1cm}
\end{minipage}
\begin{minipage}[t]{.25\linewidth}
\centering
\includegraphics[width=0.9\textwidth,height=3.5cm,keepaspectratio]{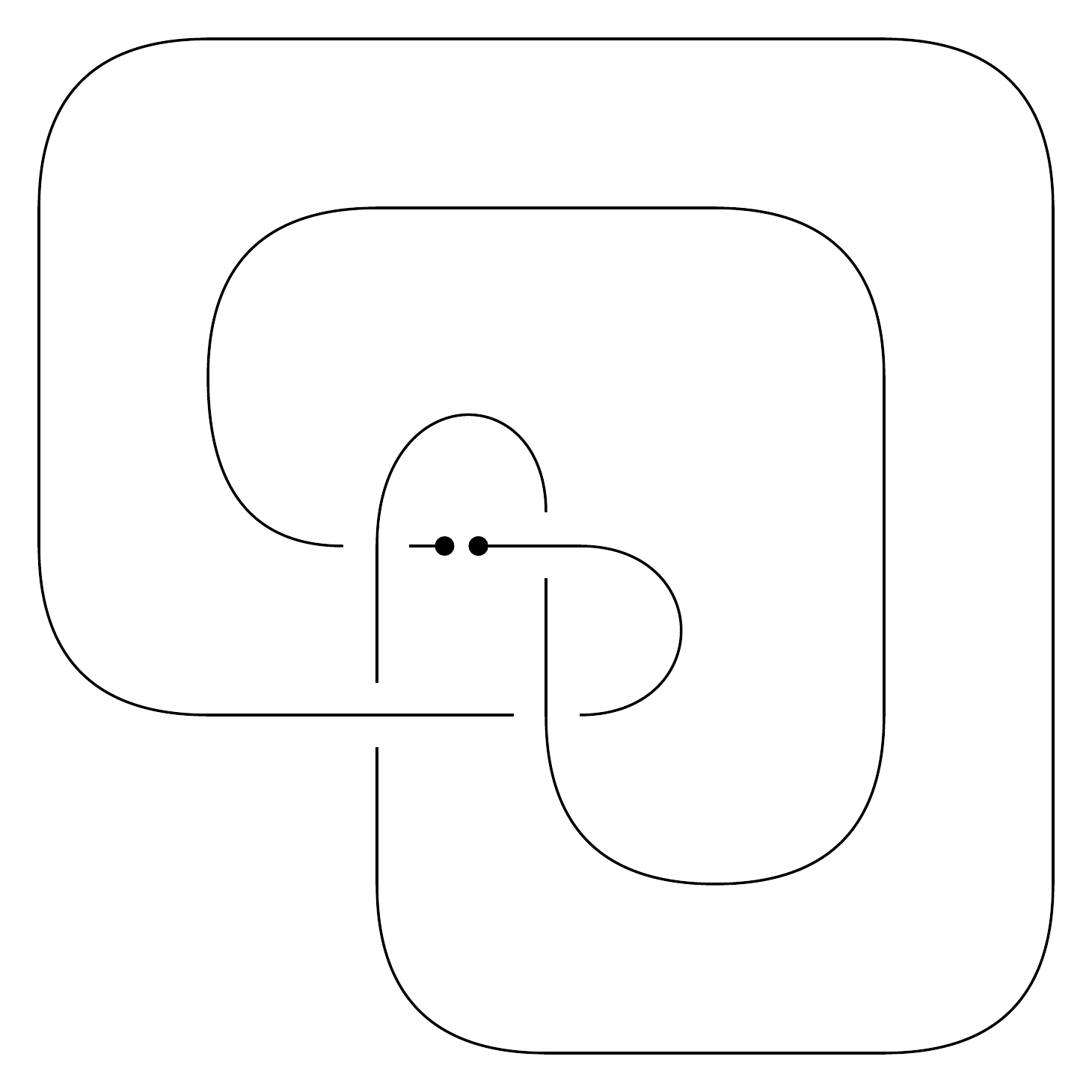}\\
\textcolor{black}{$4_{95}$}
\vspace{1cm}
\end{minipage}
\begin{minipage}[t]{.25\linewidth}
\centering
\includegraphics[width=0.9\textwidth,height=3.5cm,keepaspectratio]{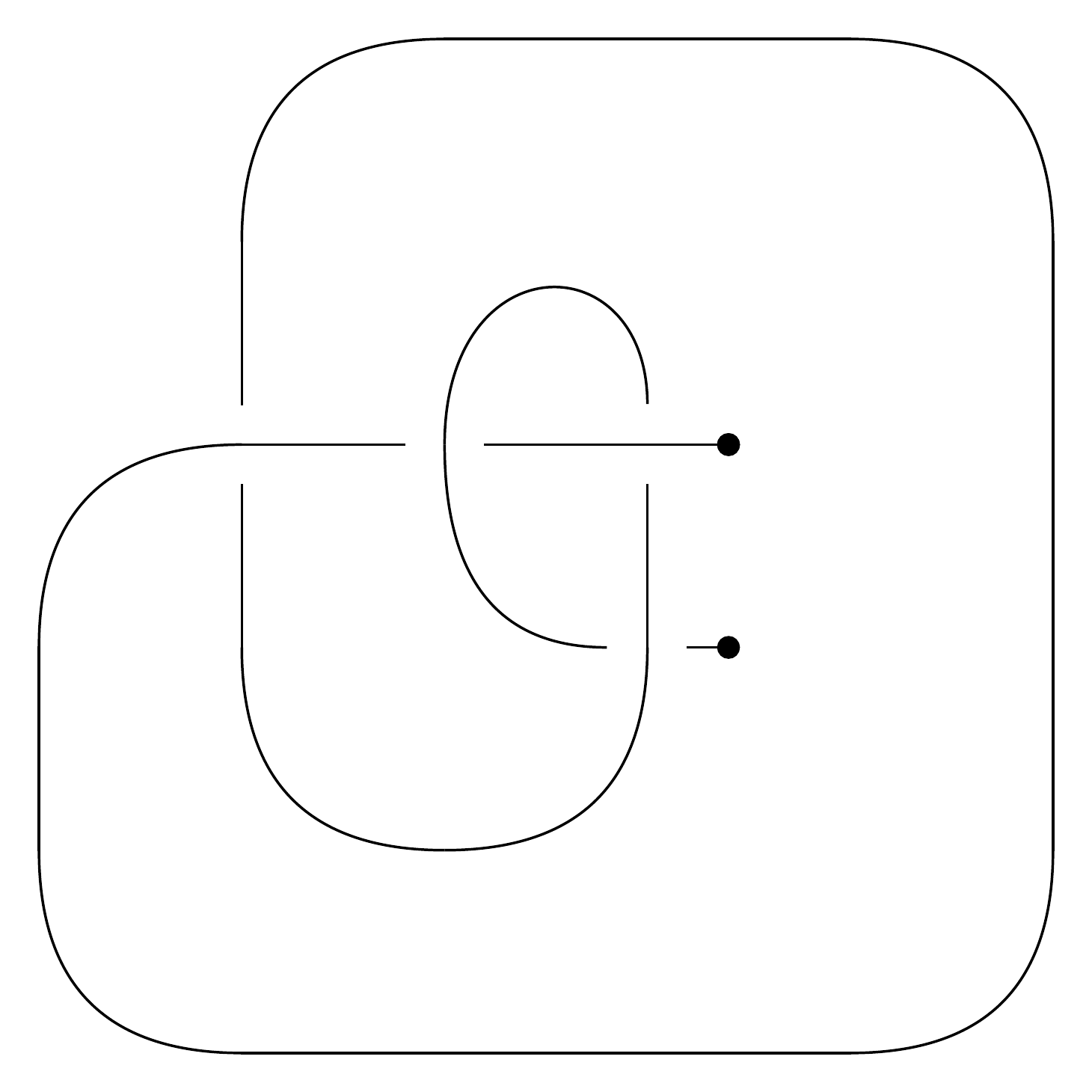}\\
\textcolor{black}{$4_{96}$}
\vspace{1cm}
\end{minipage}
\begin{minipage}[t]{.25\linewidth}
\centering
\includegraphics[width=0.9\textwidth,height=3.5cm,keepaspectratio]{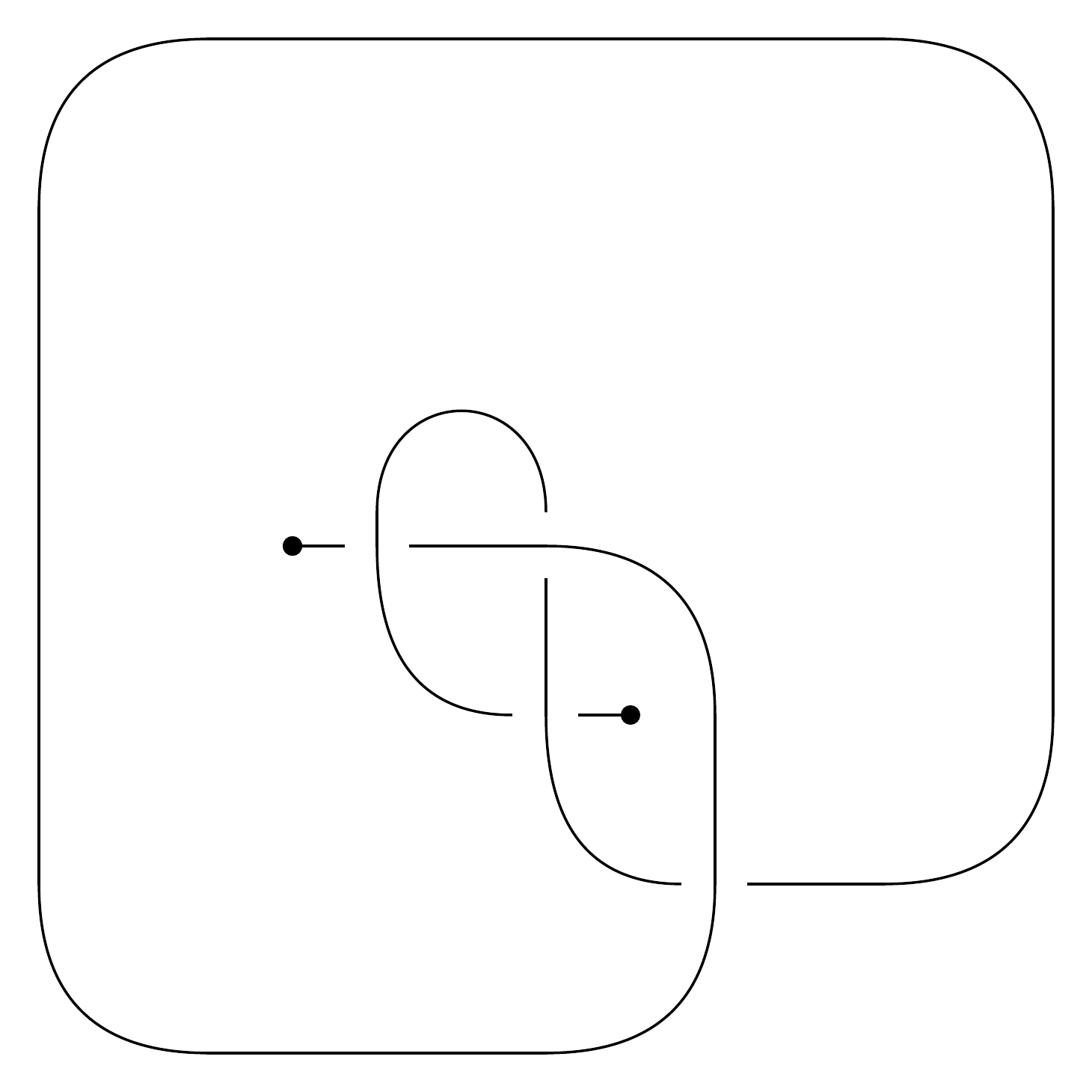}\\
\textcolor{black}{$4_{97}$}
\vspace{1cm}
\end{minipage}
\begin{minipage}[t]{.25\linewidth}
\centering
\includegraphics[width=0.9\textwidth,height=3.5cm,keepaspectratio]{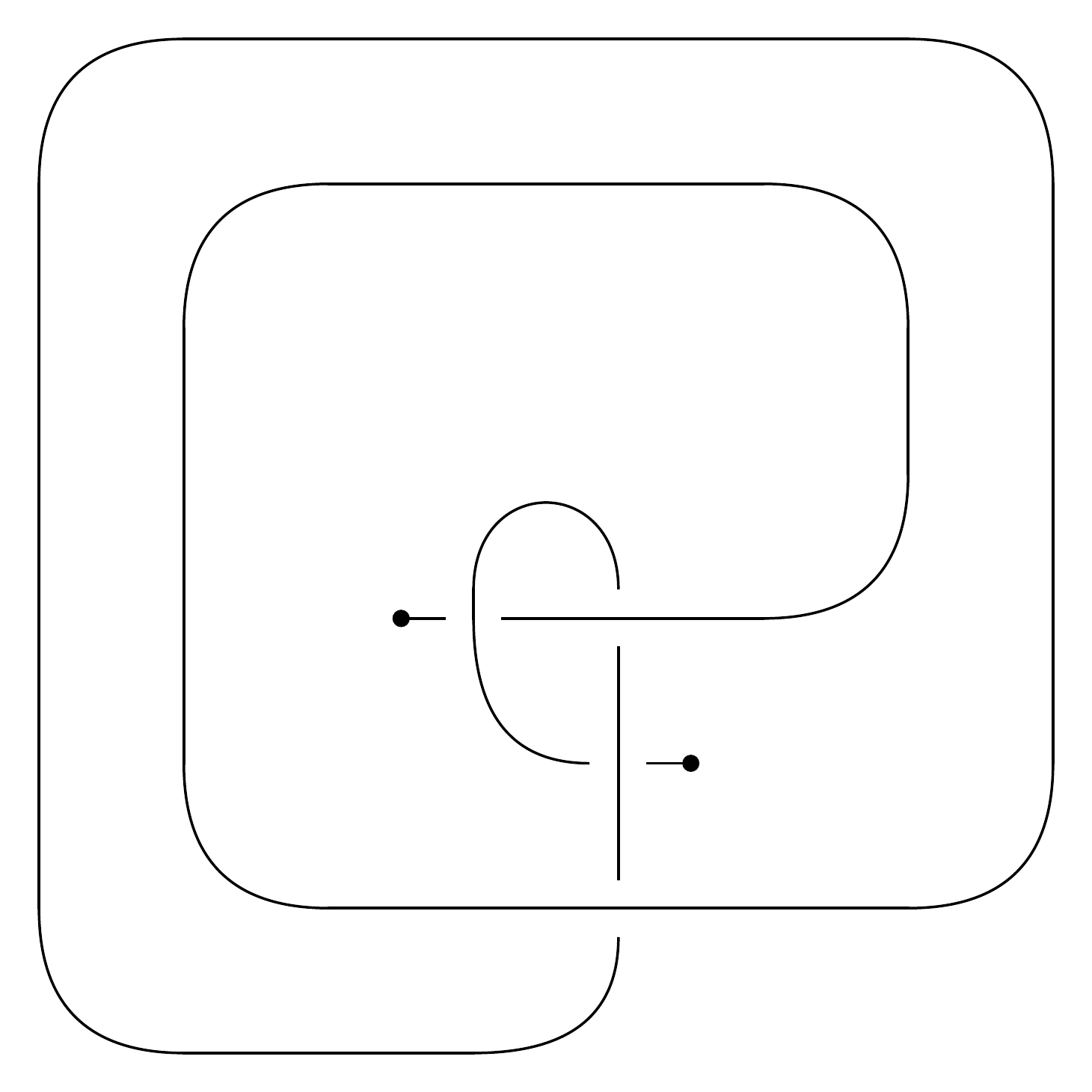}\\
\textcolor{black}{$4_{98}$}
\vspace{1cm}
\end{minipage}
\begin{minipage}[t]{.25\linewidth}
\centering
\includegraphics[width=0.9\textwidth,height=3.5cm,keepaspectratio]{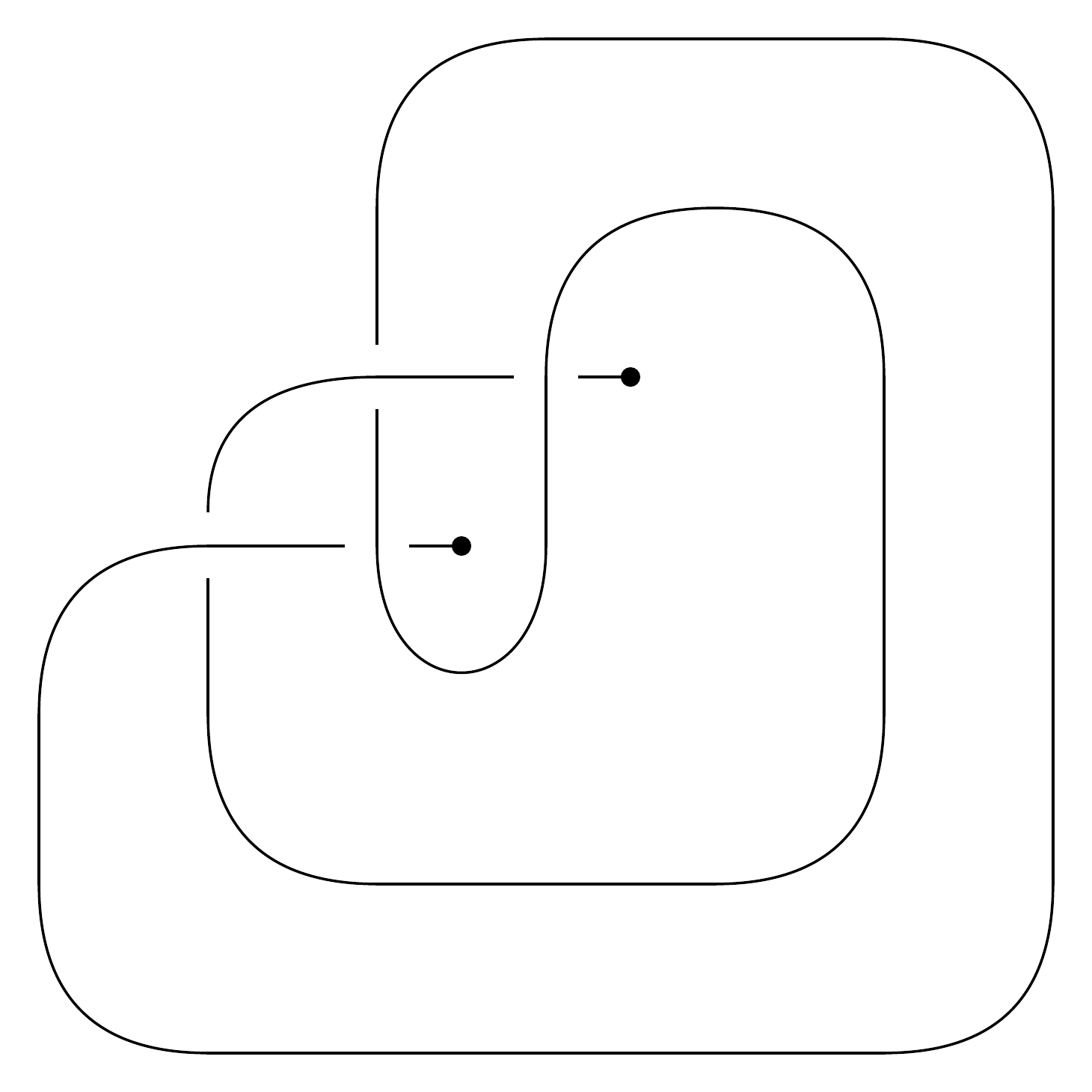}\\
\textcolor{black}{$4_{99}$}
\vspace{1cm}
\end{minipage}
\begin{minipage}[t]{.25\linewidth}
\centering
\includegraphics[width=0.9\textwidth,height=3.5cm,keepaspectratio]{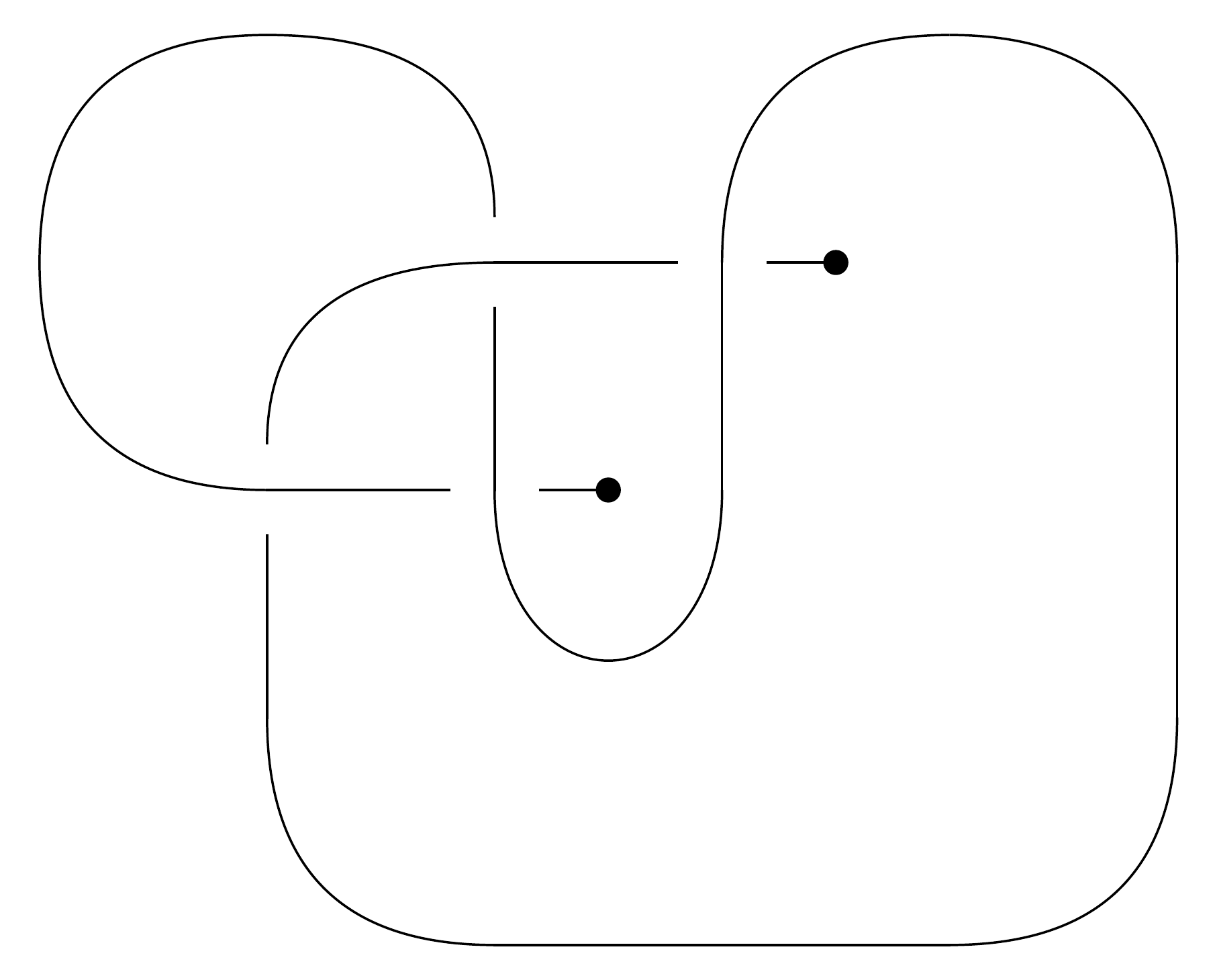}\\
\textcolor{black}{$4_{100}$}
\vspace{1cm}
\end{minipage}
\begin{minipage}[t]{.25\linewidth}
\centering
\includegraphics[width=0.9\textwidth,height=3.5cm,keepaspectratio]{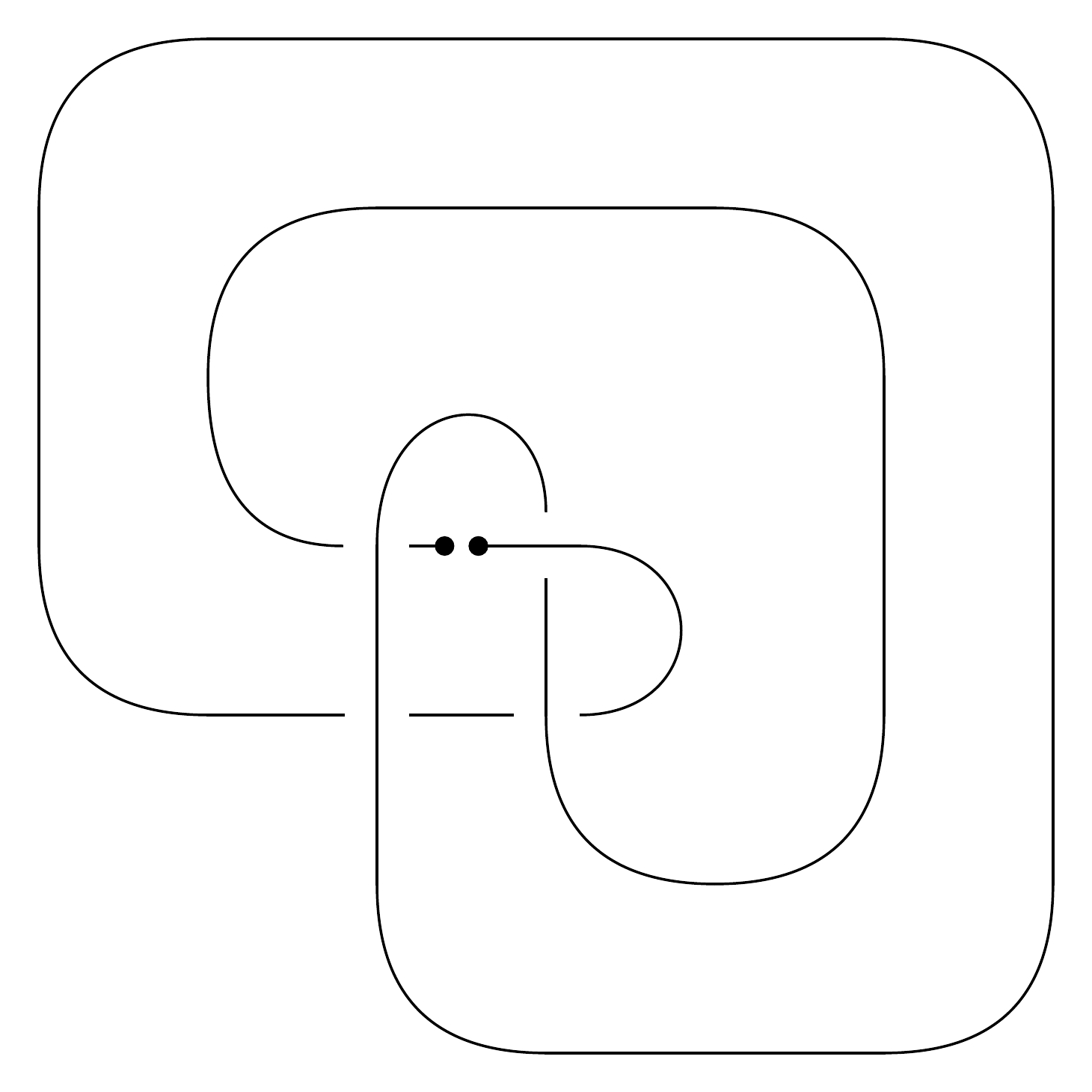}\\
\textcolor{black}{$4_{101}$}
\vspace{1cm}
\end{minipage}
\begin{minipage}[t]{.25\linewidth}
\centering
\includegraphics[width=0.9\textwidth,height=3.5cm,keepaspectratio]{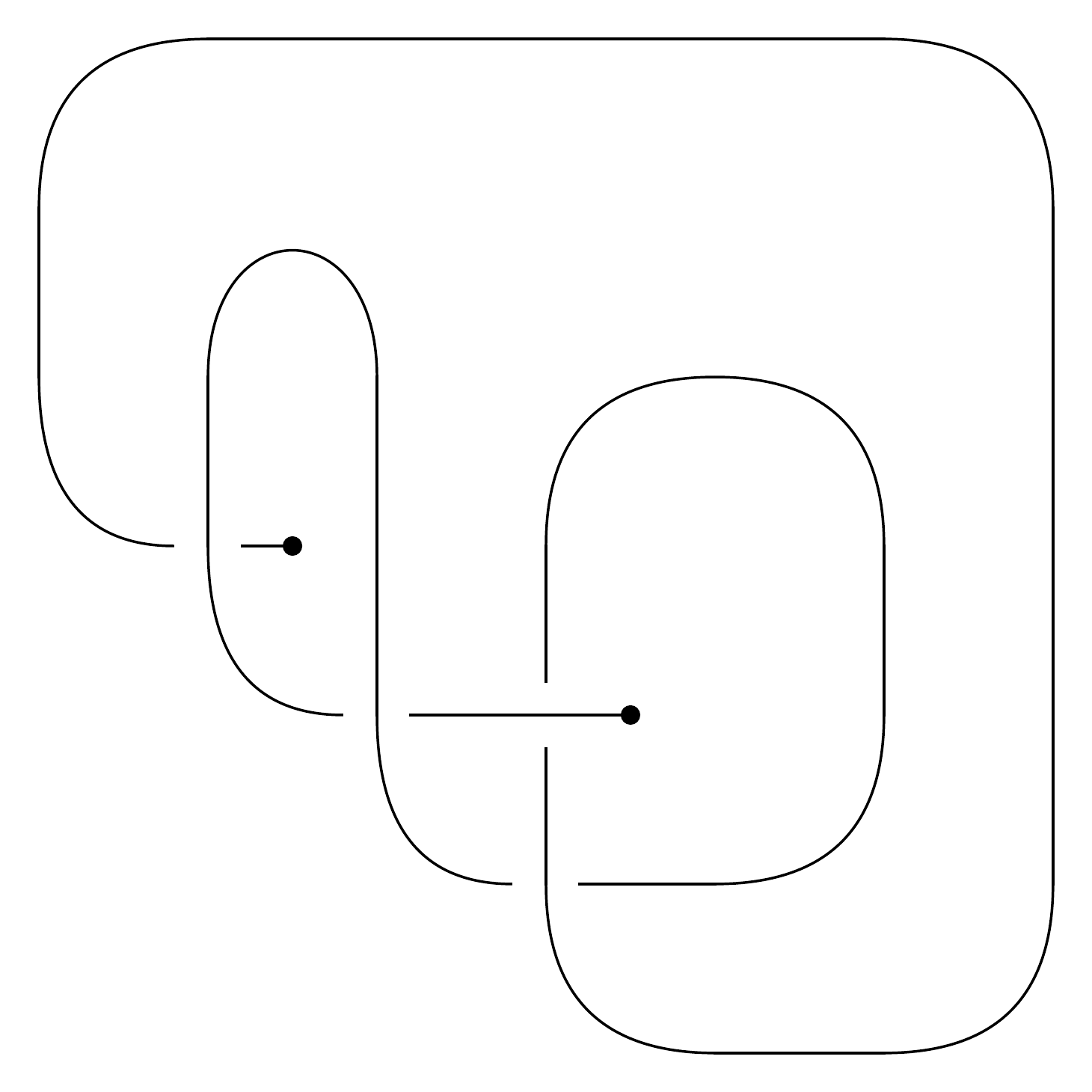}\\
\textcolor{black}{$4_{102}$}
\vspace{1cm}
\end{minipage}
\begin{minipage}[t]{.25\linewidth}
\centering
\includegraphics[width=0.9\textwidth,height=3.5cm,keepaspectratio]{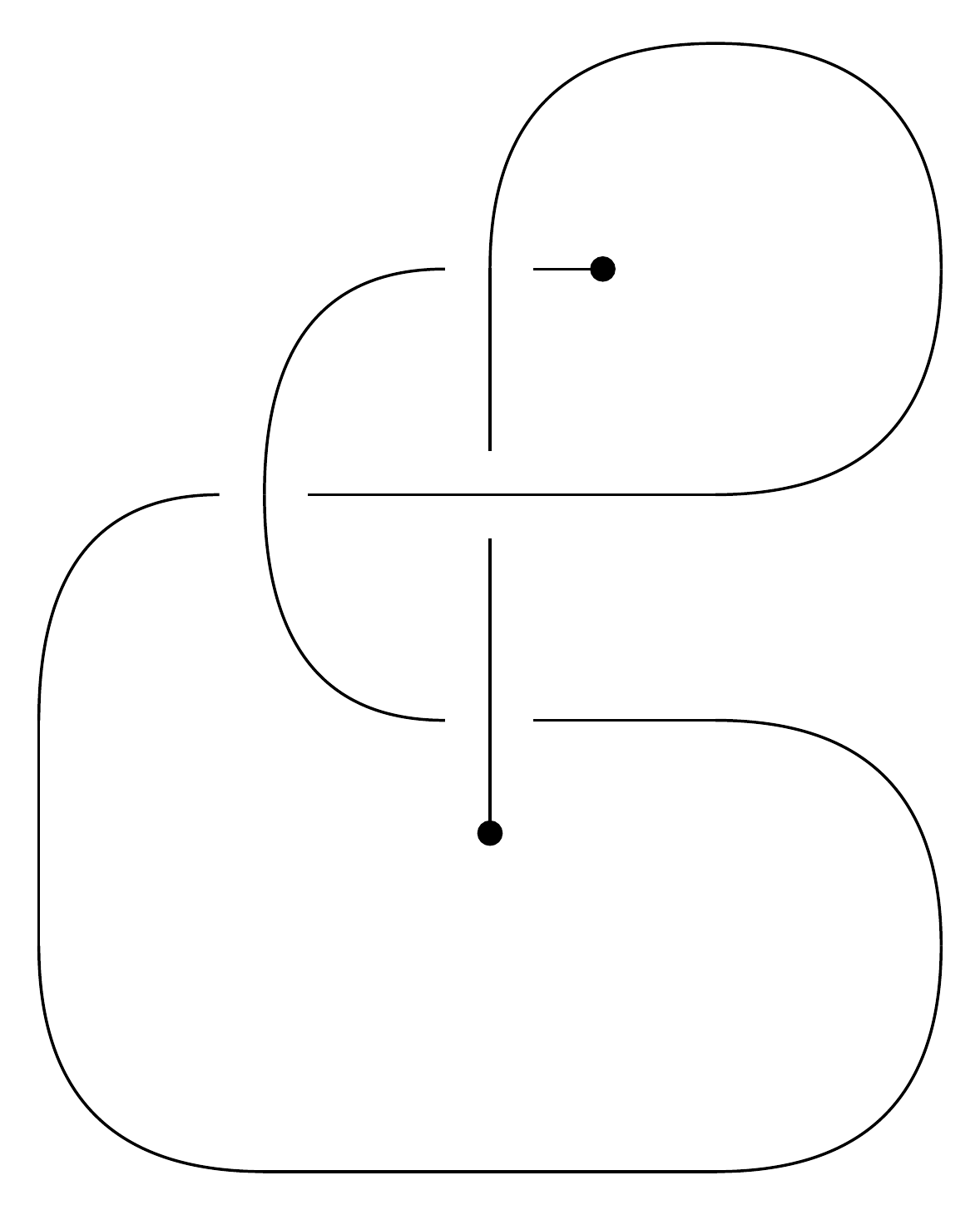}\\
\textcolor{black}{$4_{103}$}
\vspace{1cm}
\end{minipage}
\begin{minipage}[t]{.25\linewidth}
\centering
\includegraphics[width=0.9\textwidth,height=3.5cm,keepaspectratio]{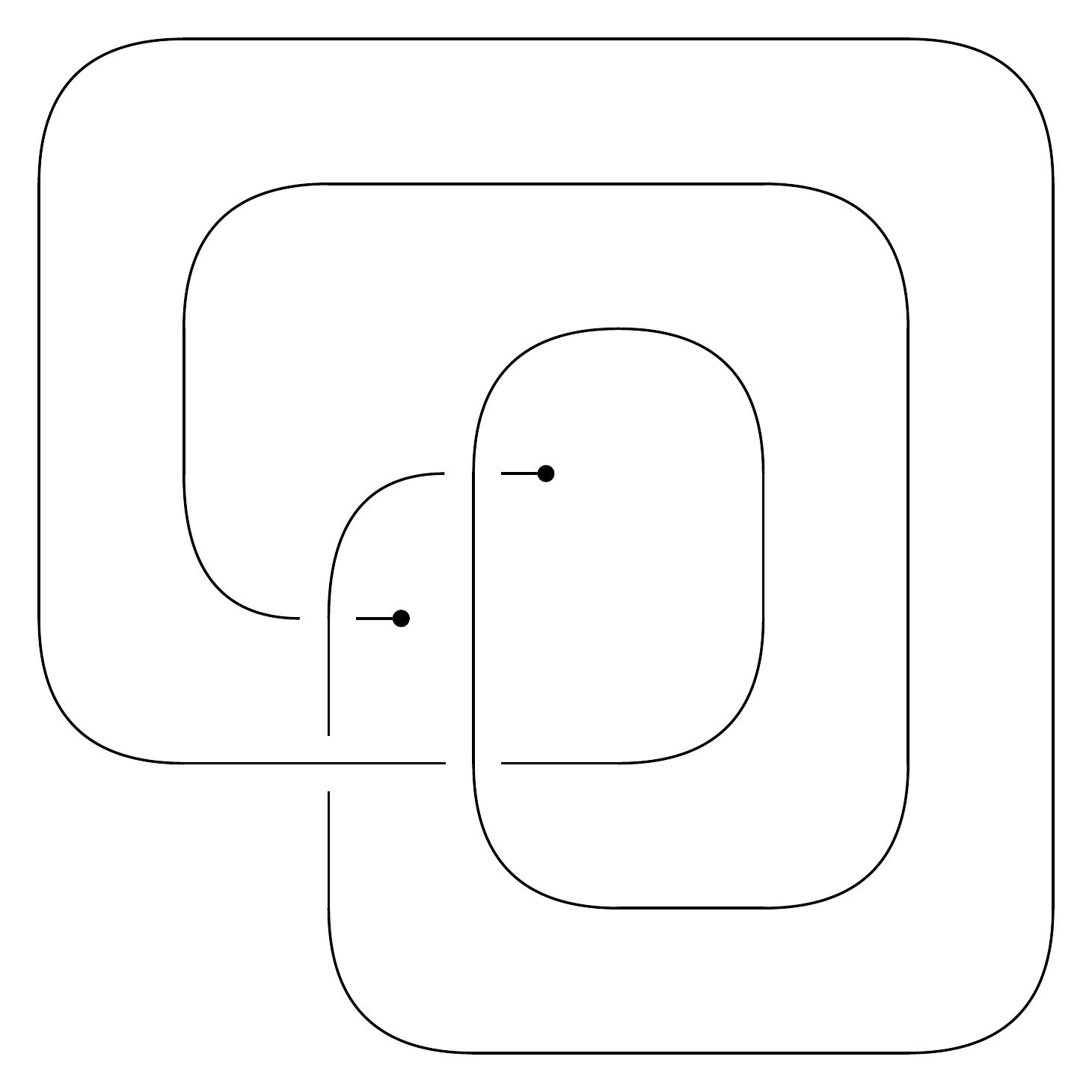}\\
\textcolor{black}{$4_{104}$}
\vspace{1cm}
\end{minipage}
\begin{minipage}[t]{.25\linewidth}
\centering
\includegraphics[width=0.9\textwidth,height=3.5cm,keepaspectratio]{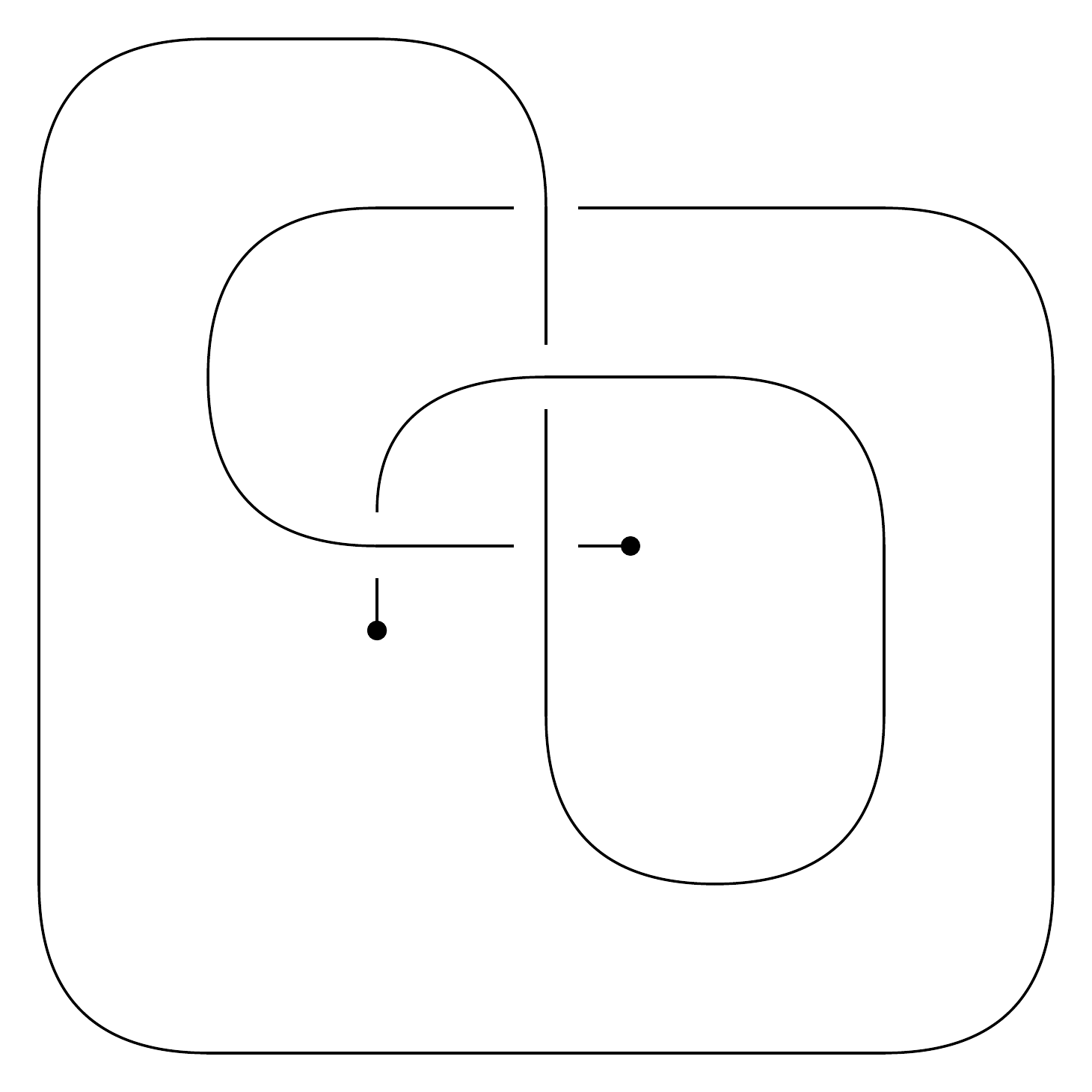}\\
\textcolor{black}{$4_{105}$}
\vspace{1cm}
\end{minipage}
\begin{minipage}[t]{.25\linewidth}
\centering
\includegraphics[width=0.9\textwidth,height=3.5cm,keepaspectratio]{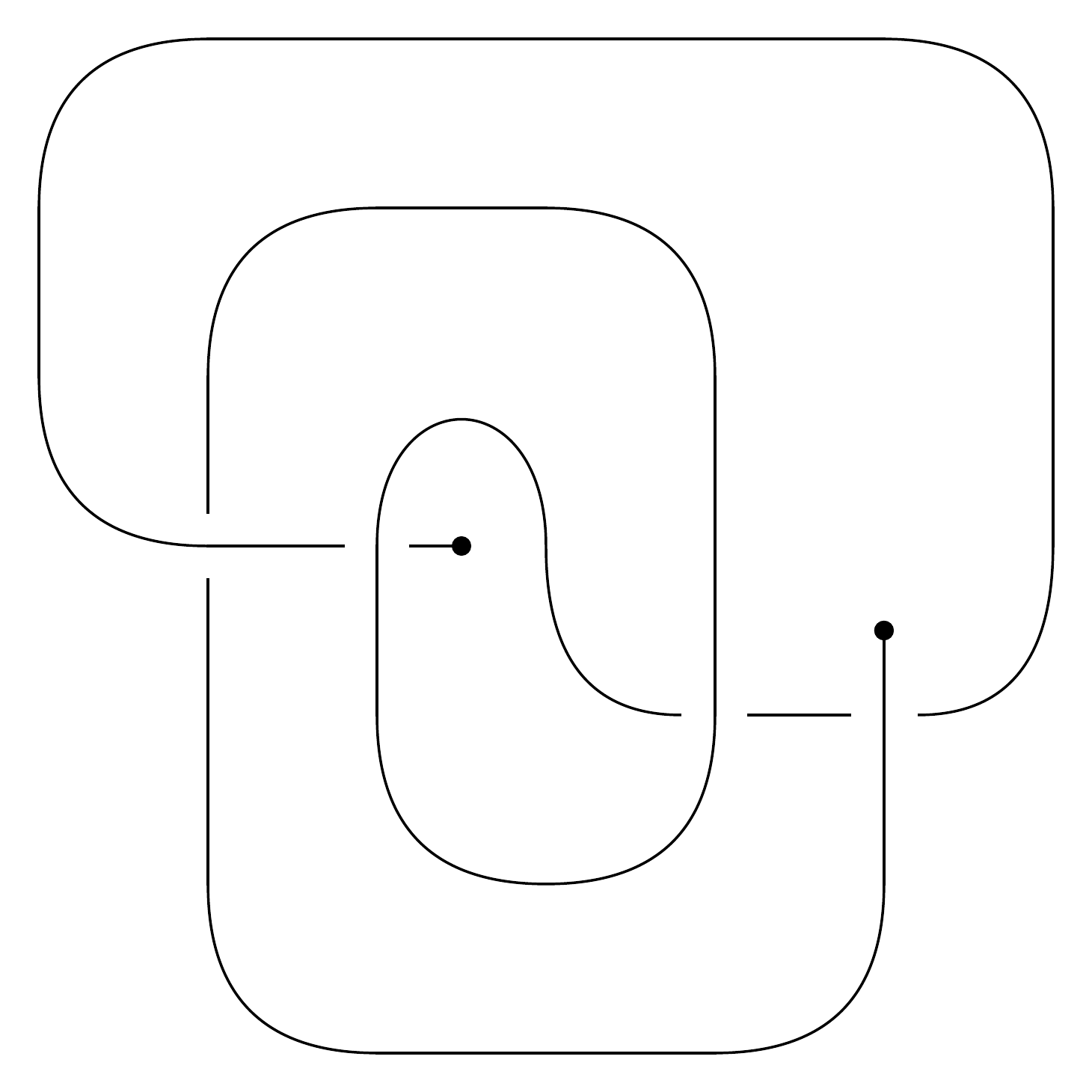}\\
\textcolor{black}{$4_{106}$}
\vspace{1cm}
\end{minipage}
\begin{minipage}[t]{.25\linewidth}
\centering
\includegraphics[width=0.9\textwidth,height=3.5cm,keepaspectratio]{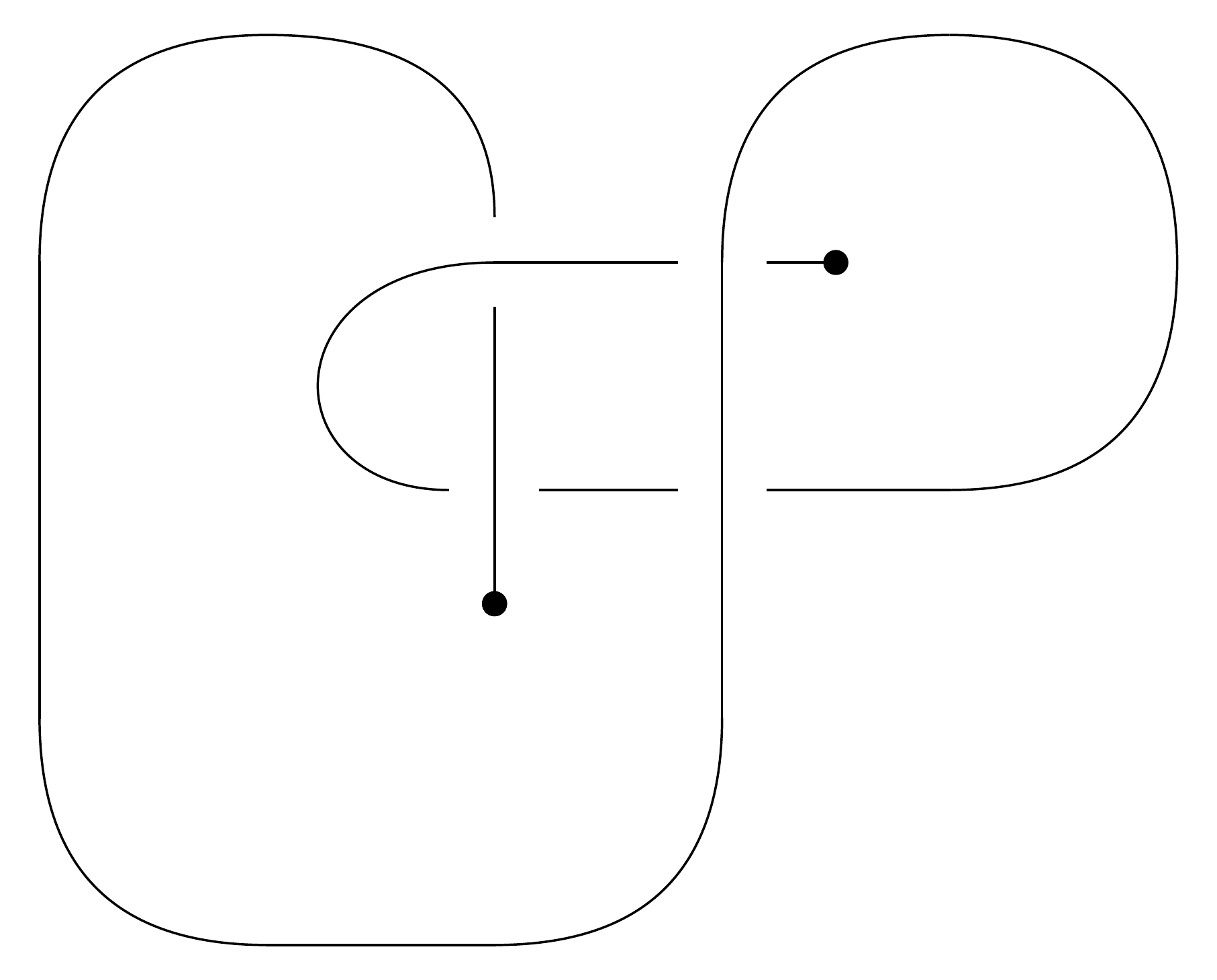}\\
\textcolor{black}{$4_{107}$}
\vspace{1cm}
\end{minipage}
\begin{minipage}[t]{.25\linewidth}
\centering
\includegraphics[width=0.9\textwidth,height=3.5cm,keepaspectratio]{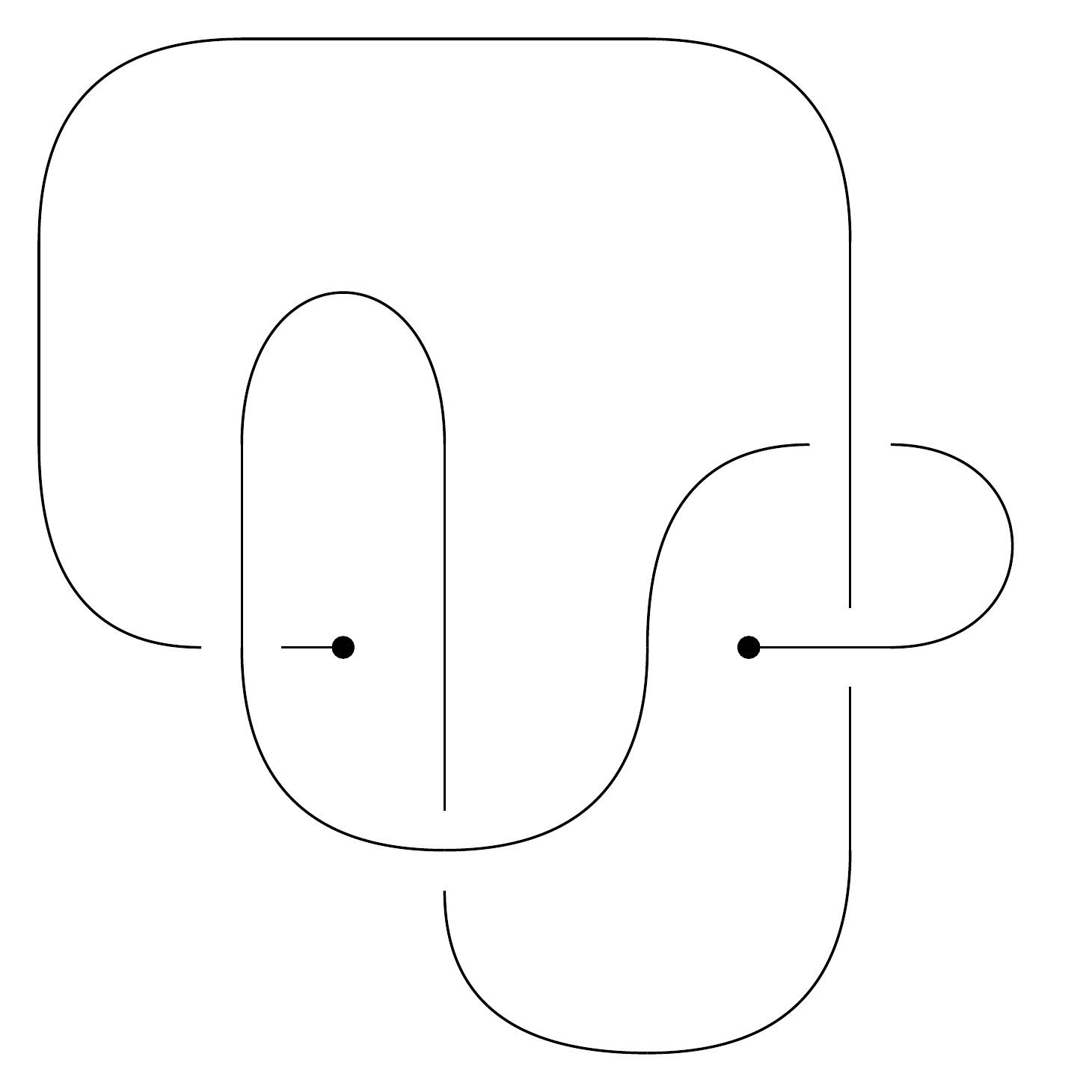}\\
\textcolor{black}{$4_{108}$}
\vspace{1cm}
\end{minipage}
\begin{minipage}[t]{.25\linewidth}
\centering
\includegraphics[width=0.9\textwidth,height=3.5cm,keepaspectratio]{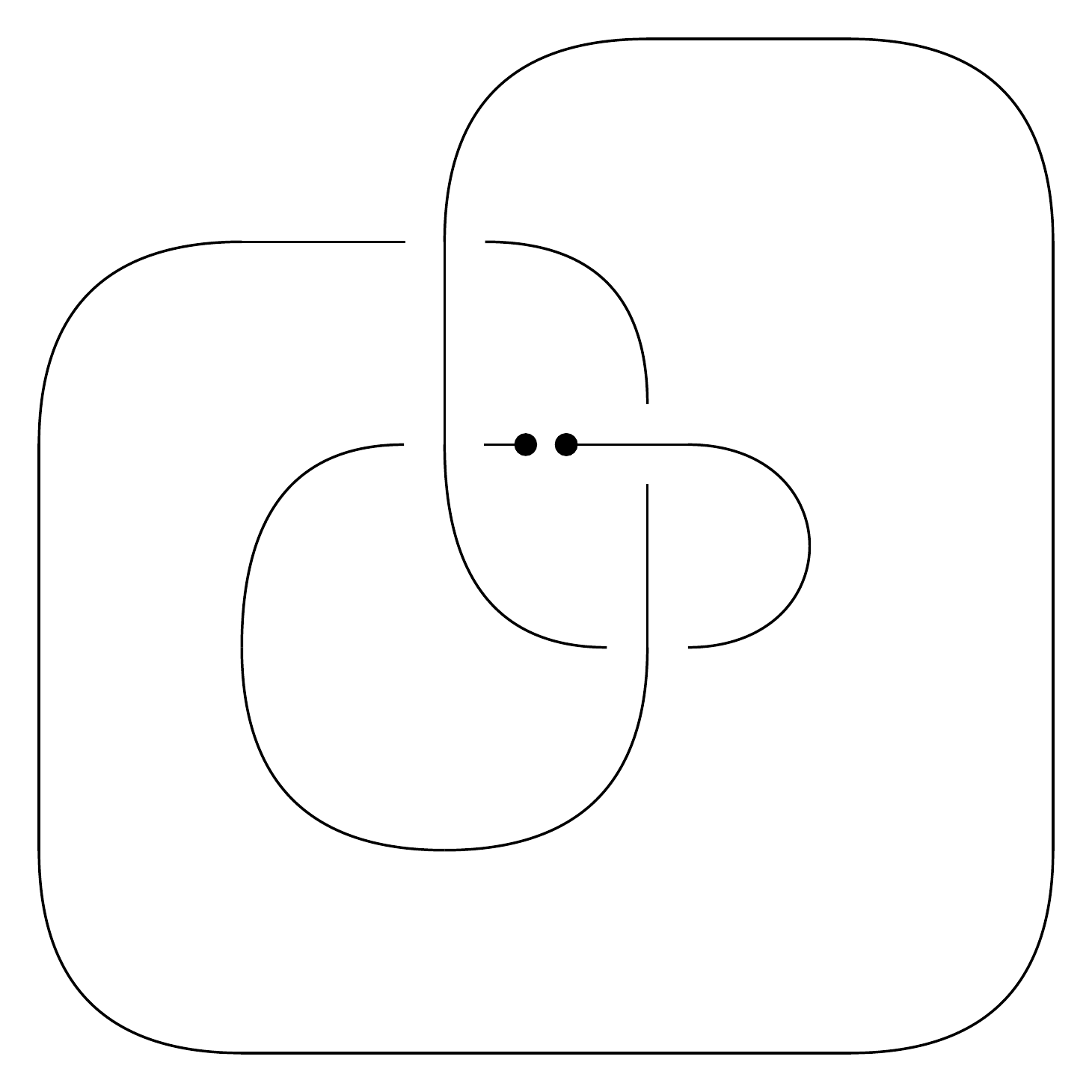}\\
\textcolor{black}{$4_{109}$}
\vspace{1cm}
\end{minipage}
\begin{minipage}[t]{.25\linewidth}
\centering
\includegraphics[width=0.9\textwidth,height=3.5cm,keepaspectratio]{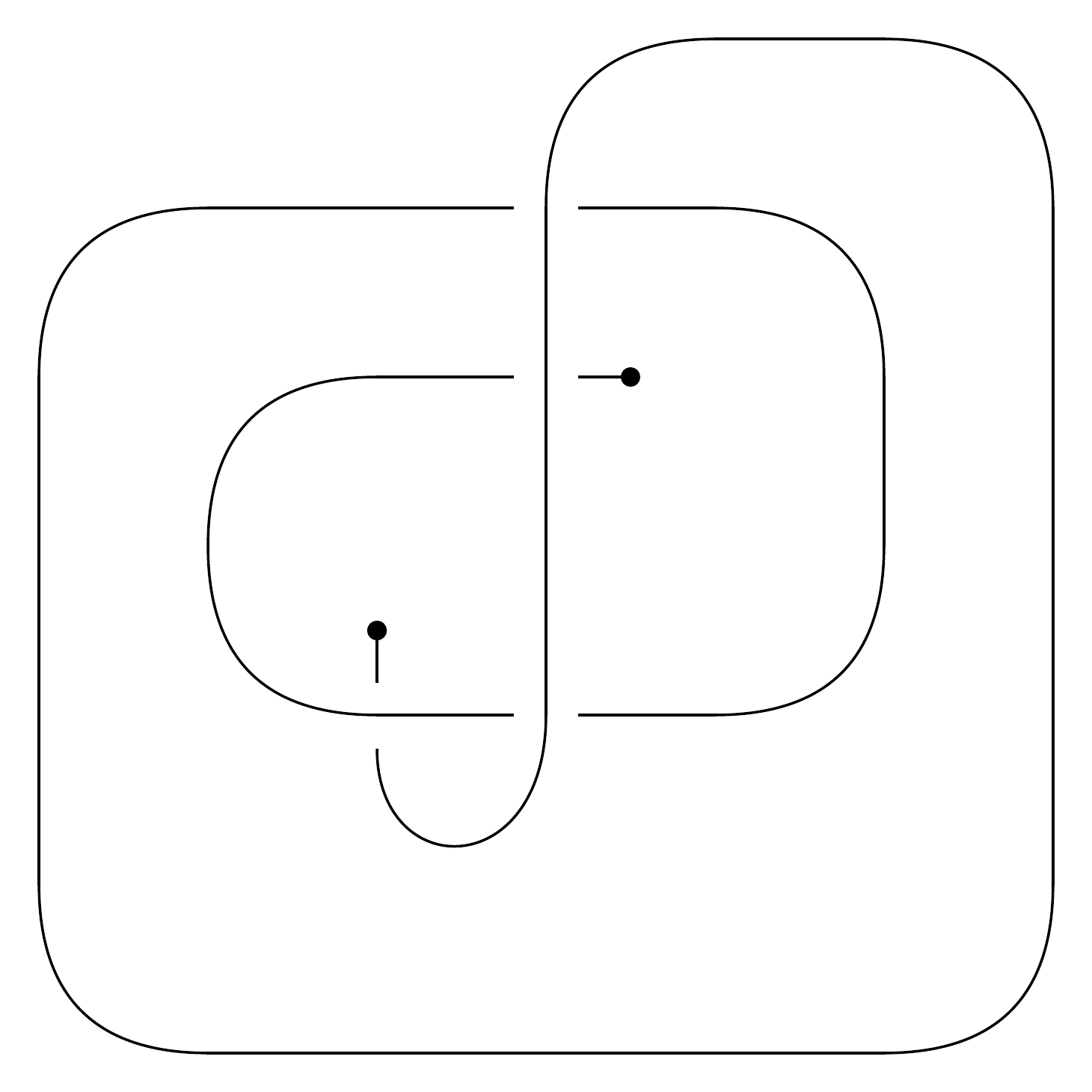}\\
\textcolor{black}{$4_{110}$}
\vspace{1cm}
\end{minipage}
\begin{minipage}[t]{.25\linewidth}
\centering
\includegraphics[width=0.9\textwidth,height=3.5cm,keepaspectratio]{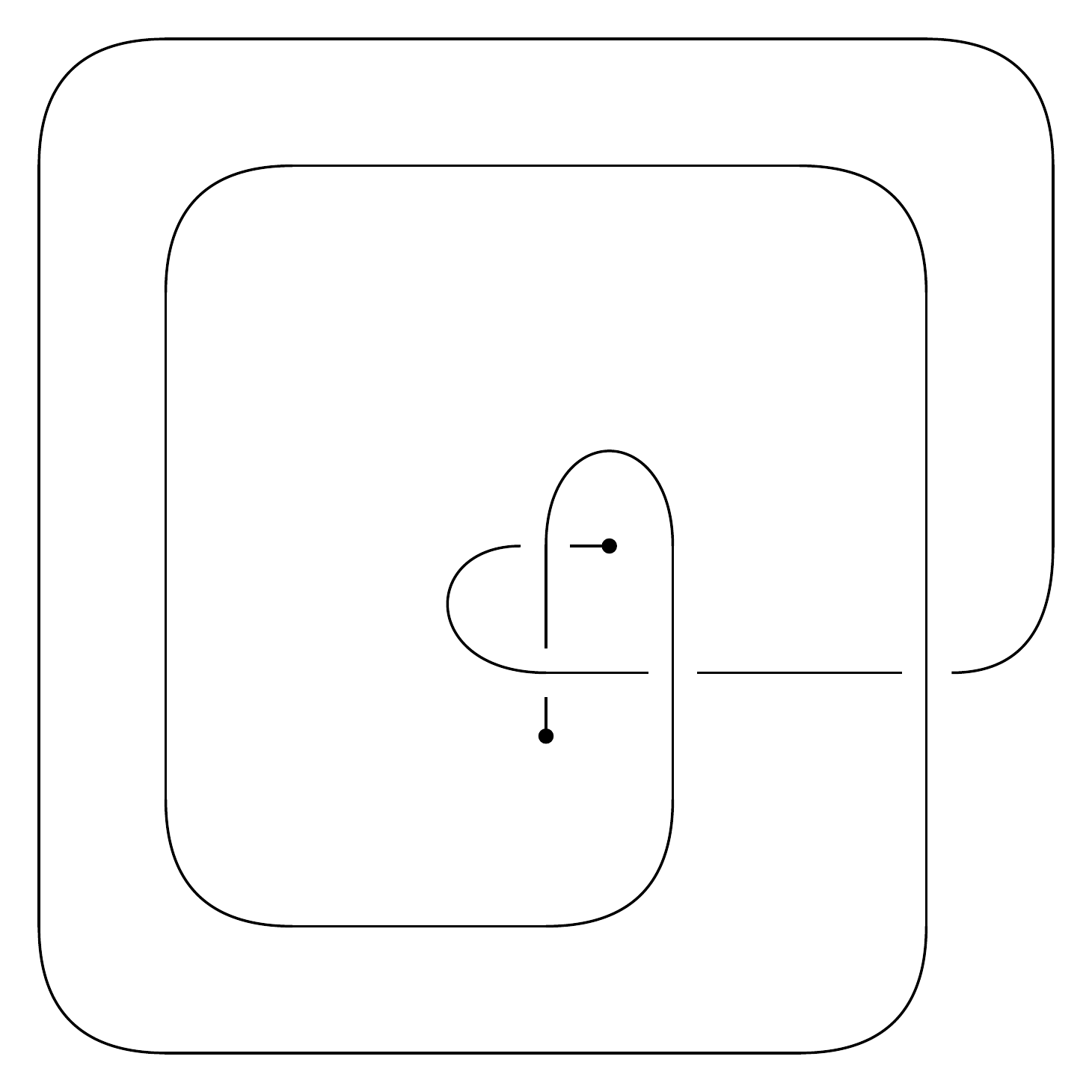}\\
\textcolor{black}{$4_{111}$}
\vspace{1cm}
\end{minipage}
\begin{minipage}[t]{.25\linewidth}
\centering
\includegraphics[width=0.9\textwidth,height=3.5cm,keepaspectratio]{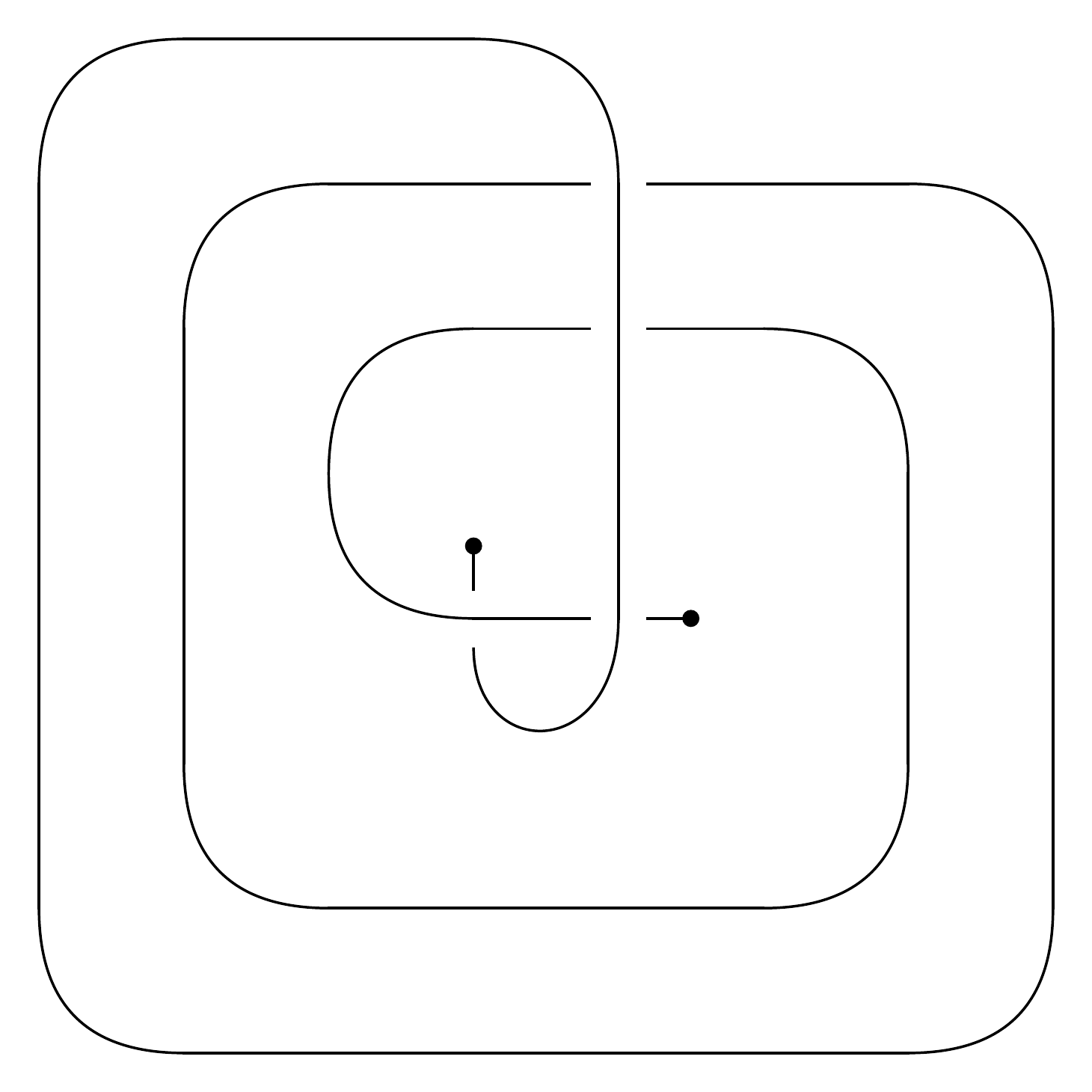}\\
\textcolor{black}{$4_{112}$}
\vspace{1cm}
\end{minipage}
\begin{minipage}[t]{.25\linewidth}
\centering
\includegraphics[width=0.9\textwidth,height=3.5cm,keepaspectratio]{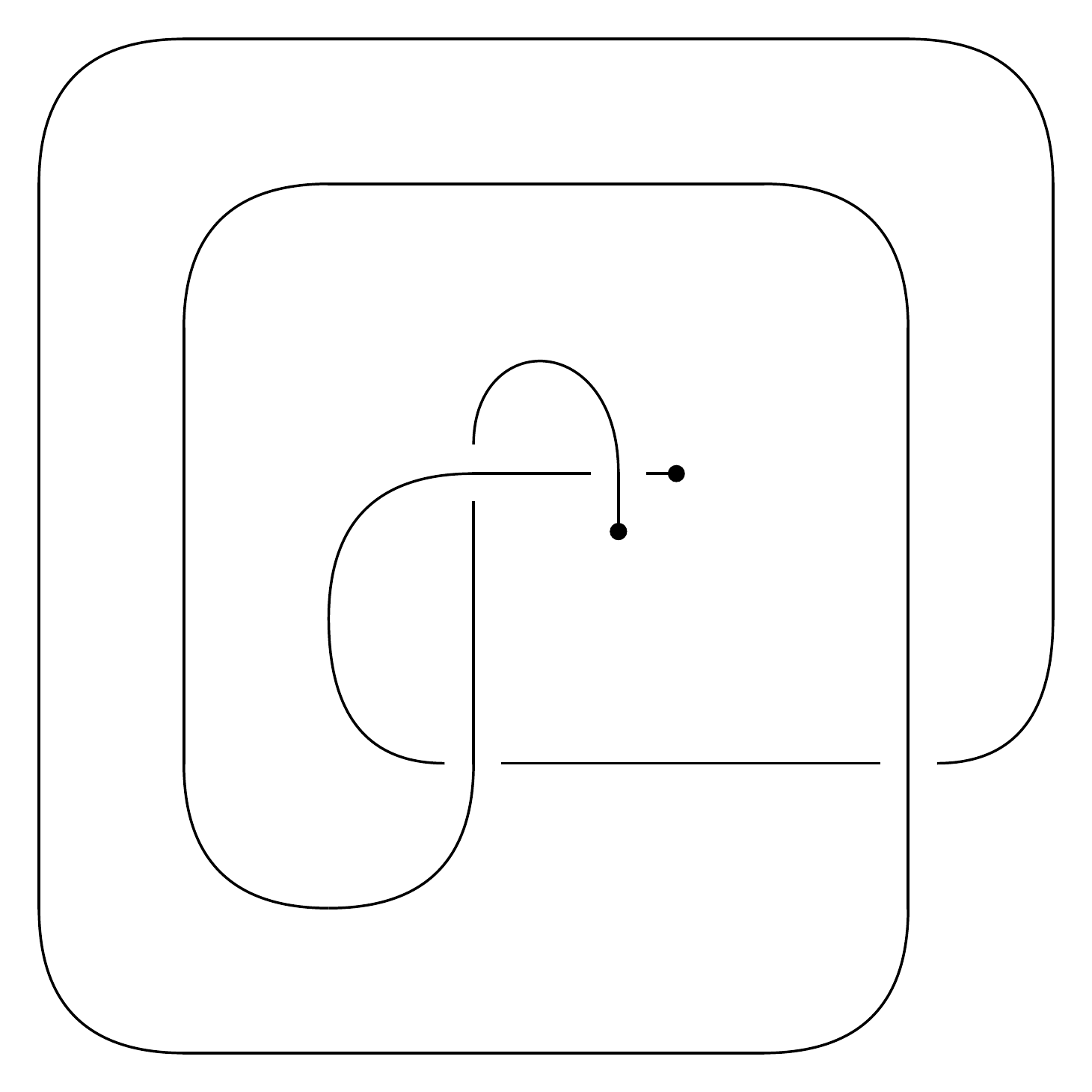}\\
\textcolor{black}{$4_{113}$}
\vspace{1cm}
\end{minipage}
\begin{minipage}[t]{.25\linewidth}
\centering
\includegraphics[width=0.9\textwidth,height=3.5cm,keepaspectratio]{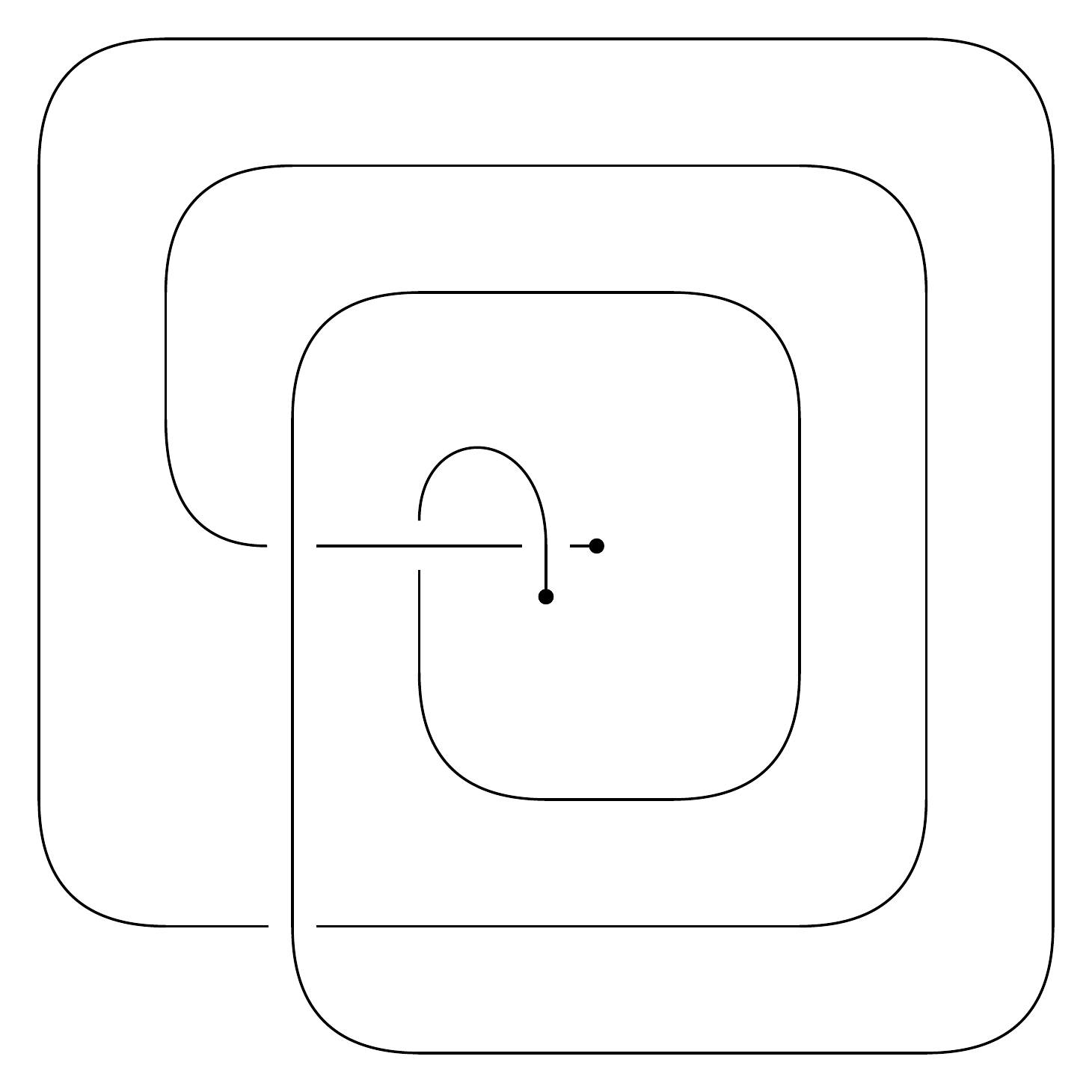}\\
\textcolor{black}{$4_{114}$}
\vspace{1cm}
\end{minipage}
\begin{minipage}[t]{.25\linewidth}
\centering
\includegraphics[width=0.9\textwidth,height=3.5cm,keepaspectratio]{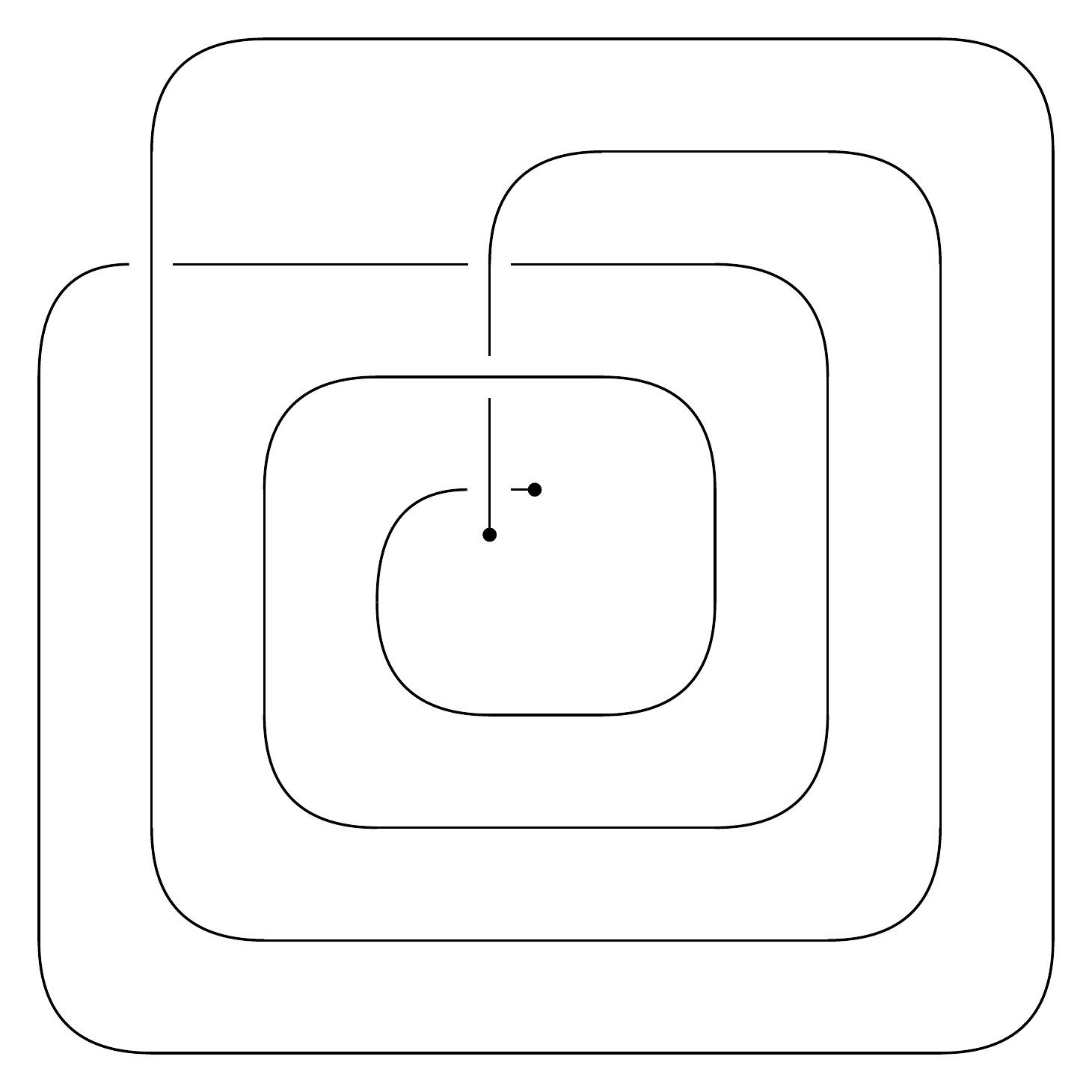}\\
\textcolor{black}{$4_{115}$}
\vspace{1cm}
\end{minipage}
\begin{minipage}[t]{.25\linewidth}
\centering
\includegraphics[width=0.9\textwidth,height=3.5cm,keepaspectratio]{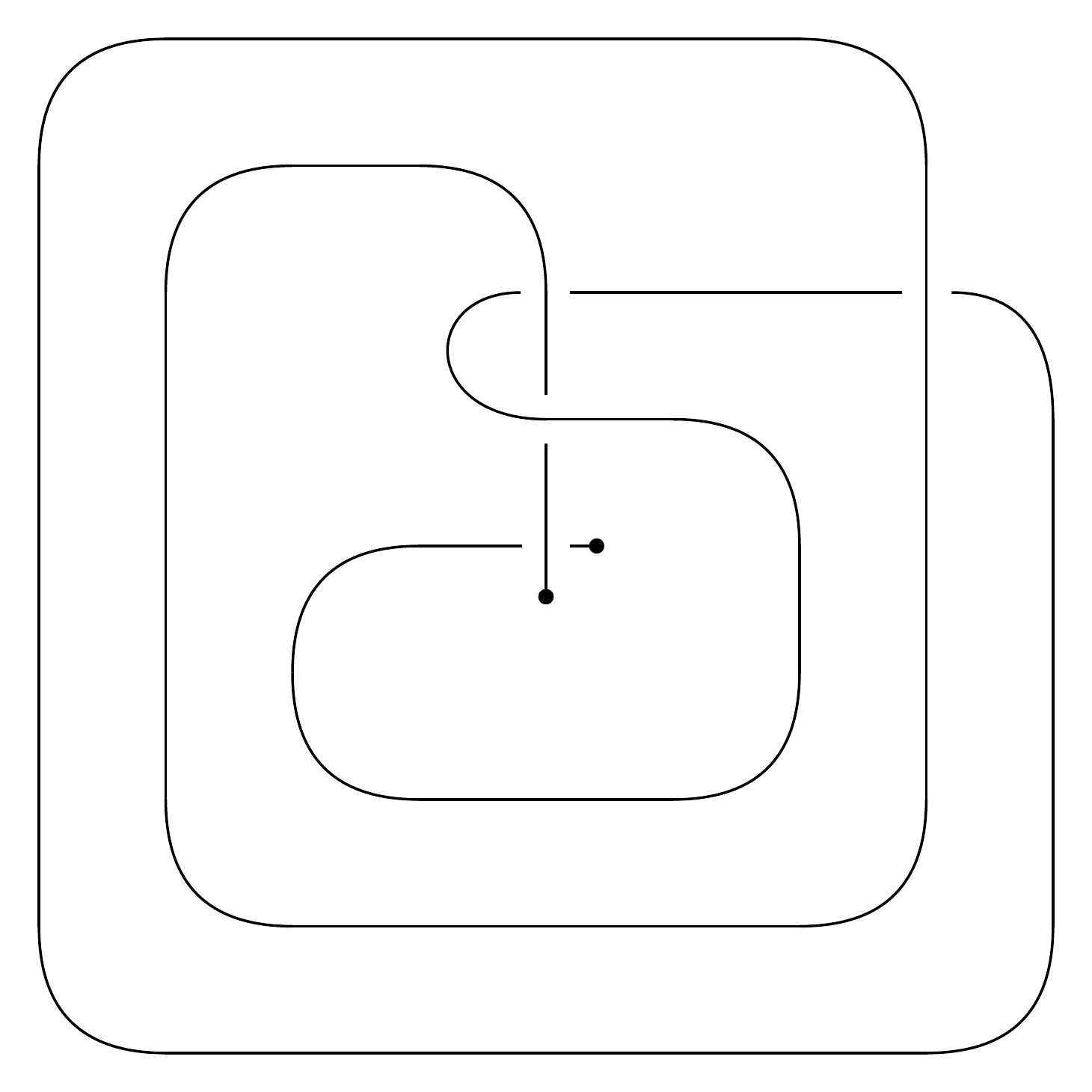}\\
\textcolor{black}{$4_{116}$}
\vspace{1cm}
\end{minipage}
\begin{minipage}[t]{.25\linewidth}
\centering
\includegraphics[width=0.9\textwidth,height=3.5cm,keepaspectratio]{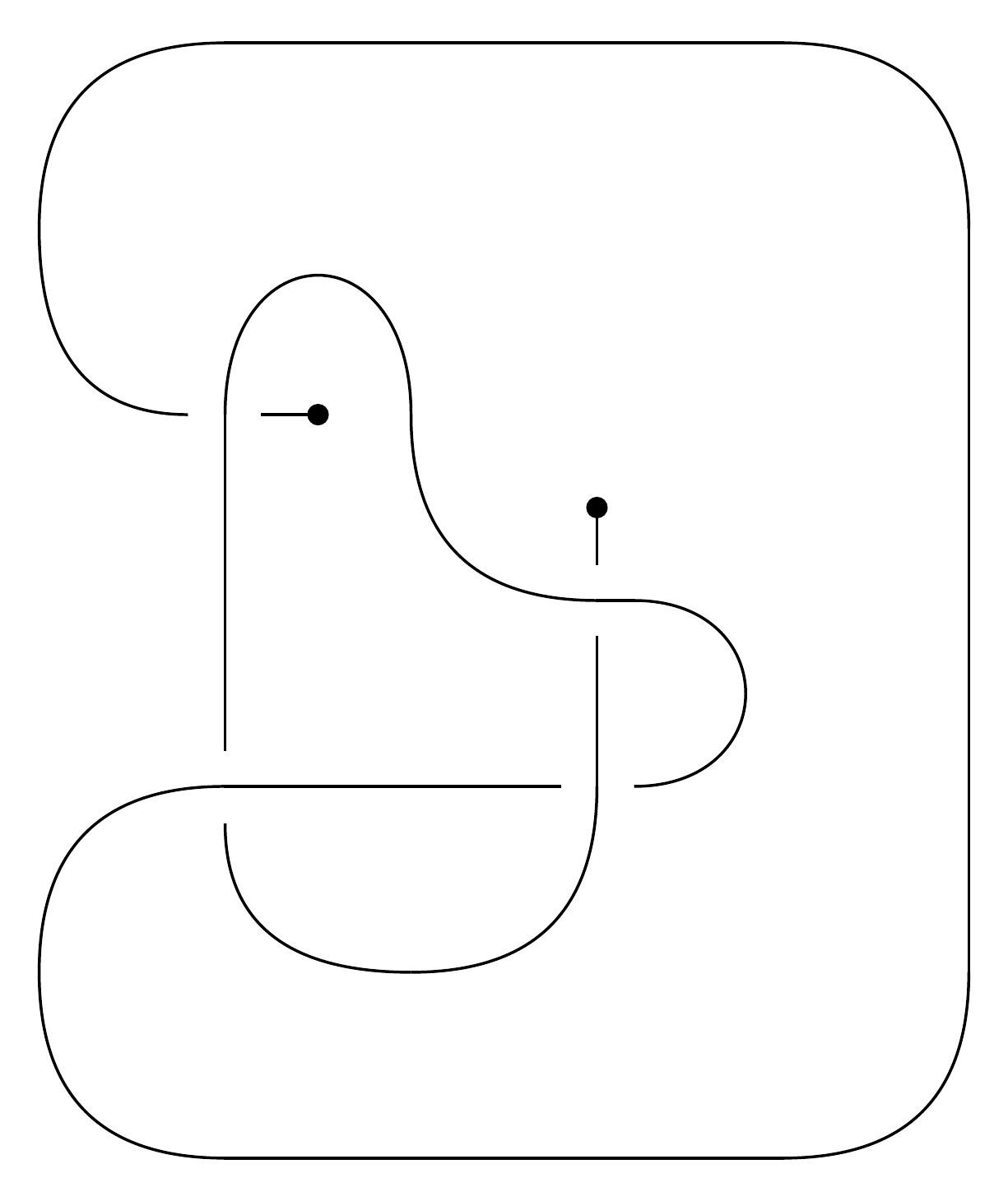}\\
\textcolor{black}{$4_{117}$}
\vspace{1cm}
\end{minipage}
\begin{minipage}[t]{.25\linewidth}
\centering
\includegraphics[width=0.9\textwidth,height=3.5cm,keepaspectratio]{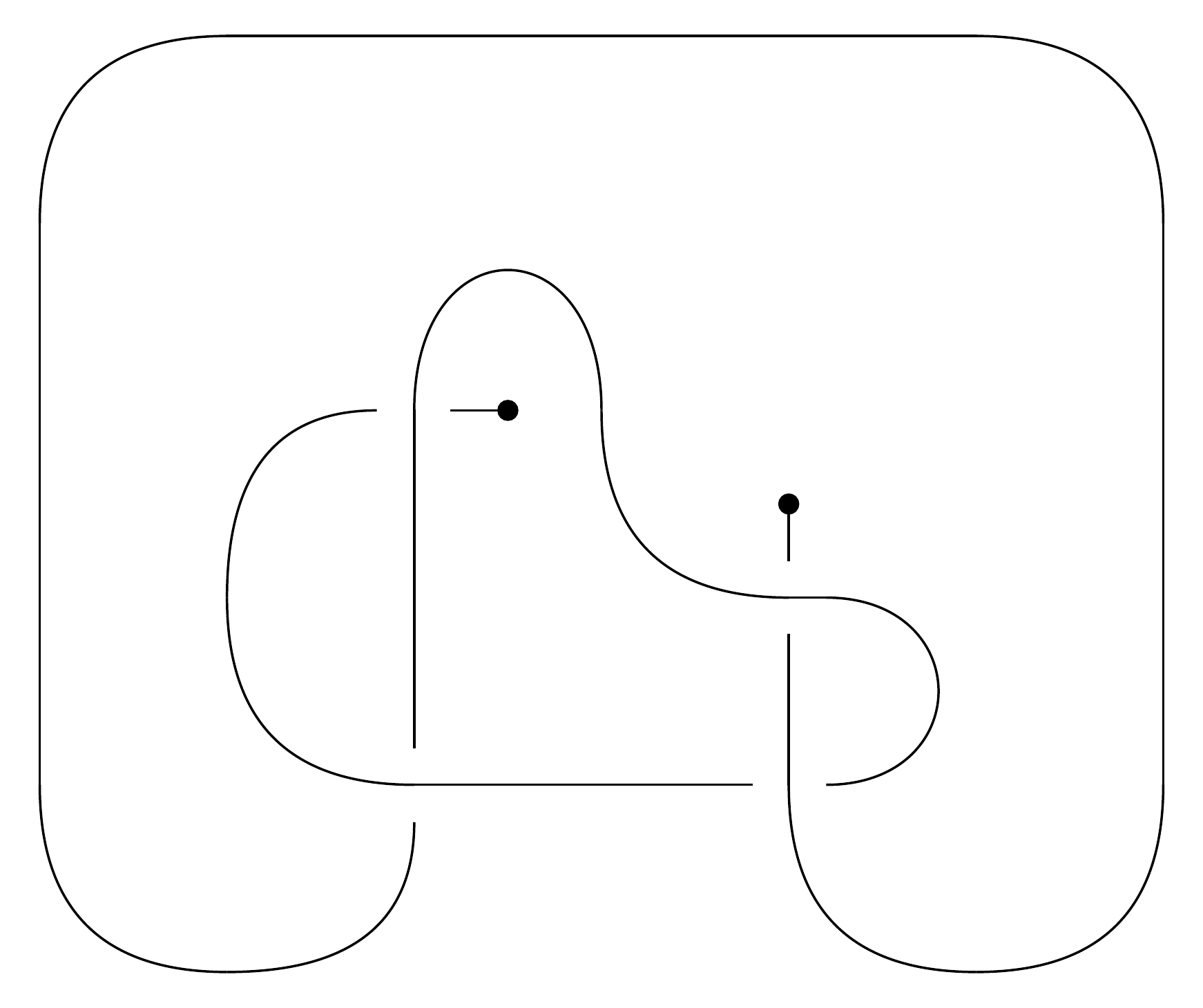}\\
\textcolor{black}{$4_{118}$}
\vspace{1cm}
\end{minipage}
\begin{minipage}[t]{.25\linewidth}
\centering
\includegraphics[width=0.9\textwidth,height=3.5cm,keepaspectratio]{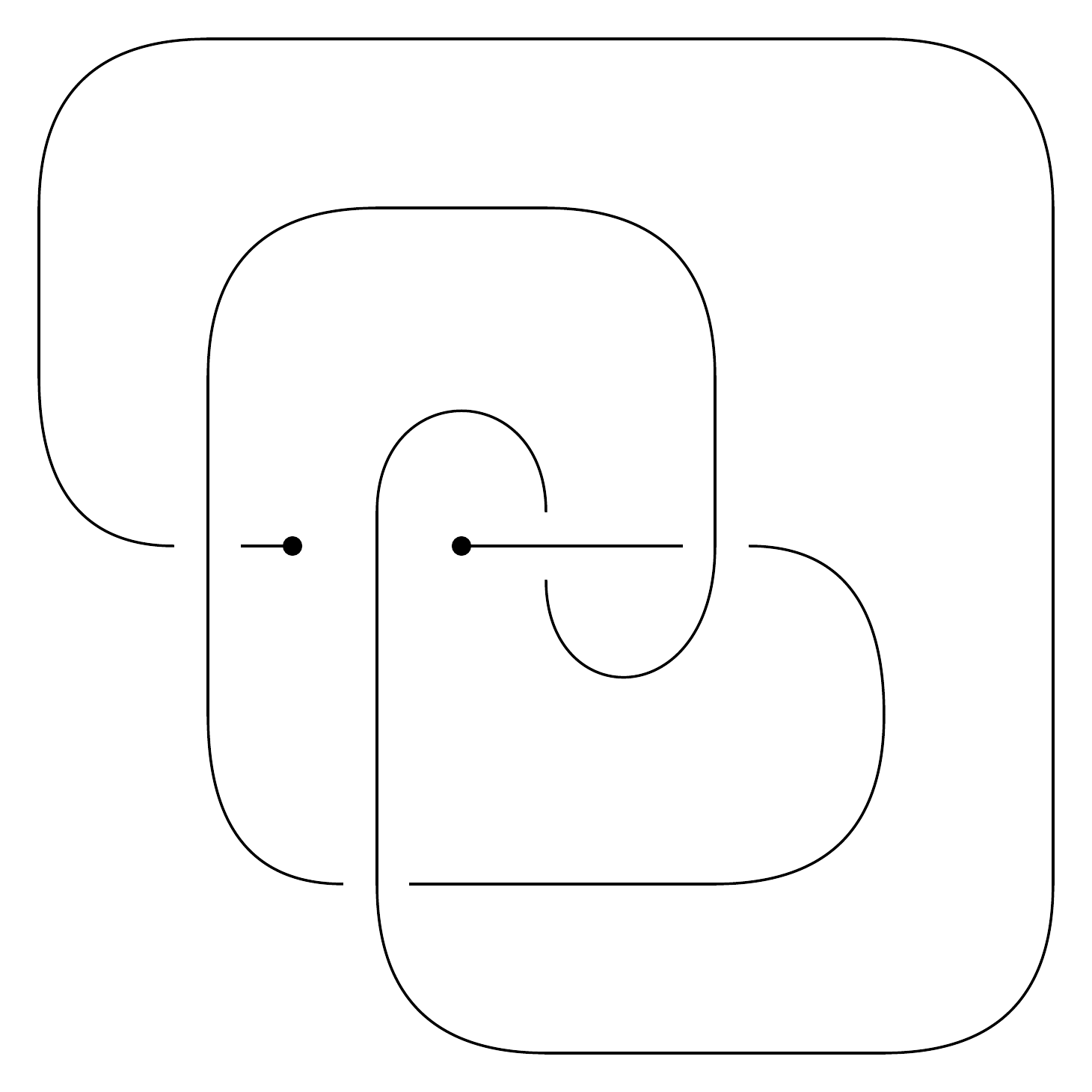}\\
\textcolor{black}{$4_{119}$}
\vspace{1cm}
\end{minipage}
\begin{minipage}[t]{.25\linewidth}
\centering
\includegraphics[width=0.9\textwidth,height=3.5cm,keepaspectratio]{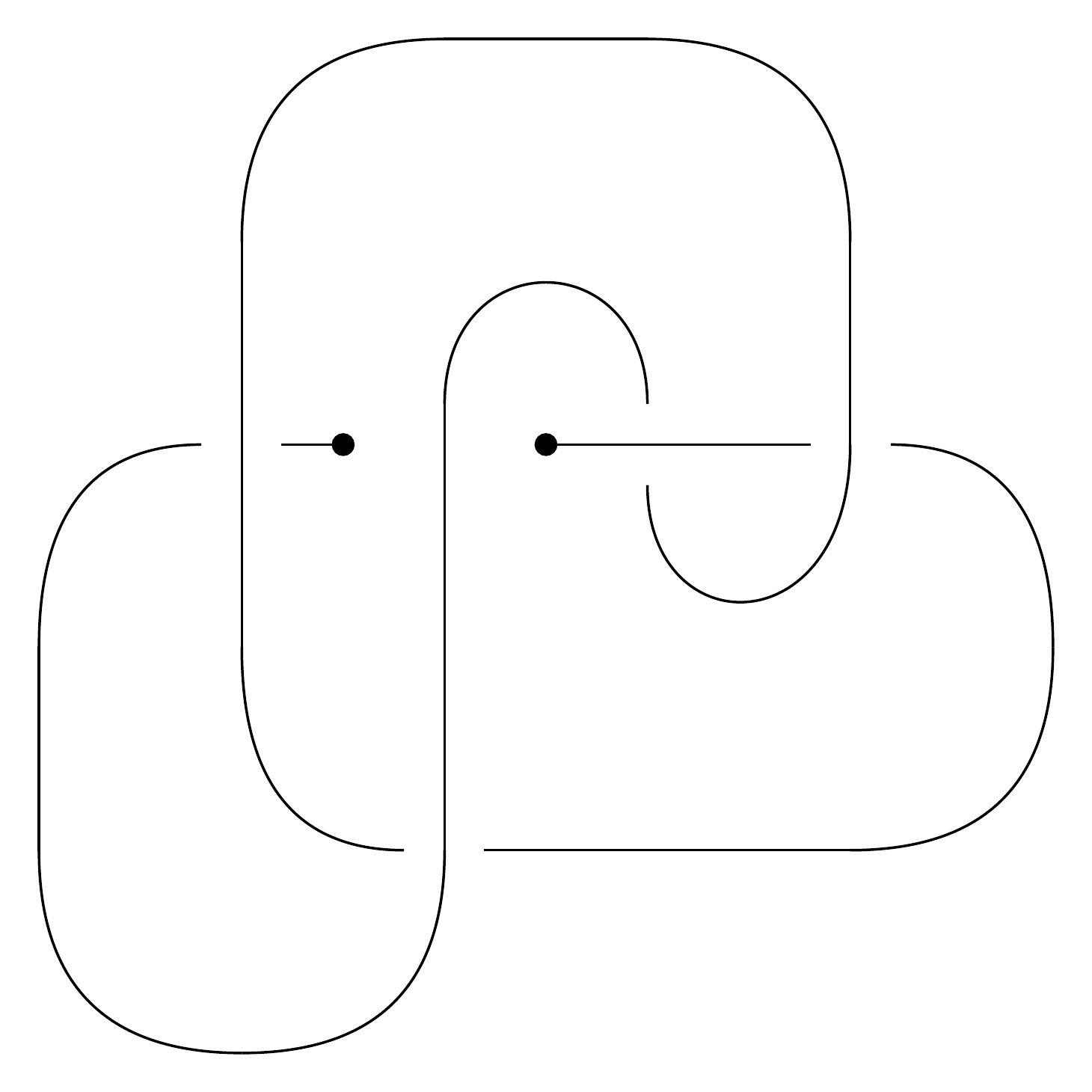}\\
\textcolor{black}{$4_{120}$}
\vspace{1cm}
\end{minipage}
\begin{minipage}[t]{.25\linewidth}
\centering
\includegraphics[width=0.9\textwidth,height=3.5cm,keepaspectratio]{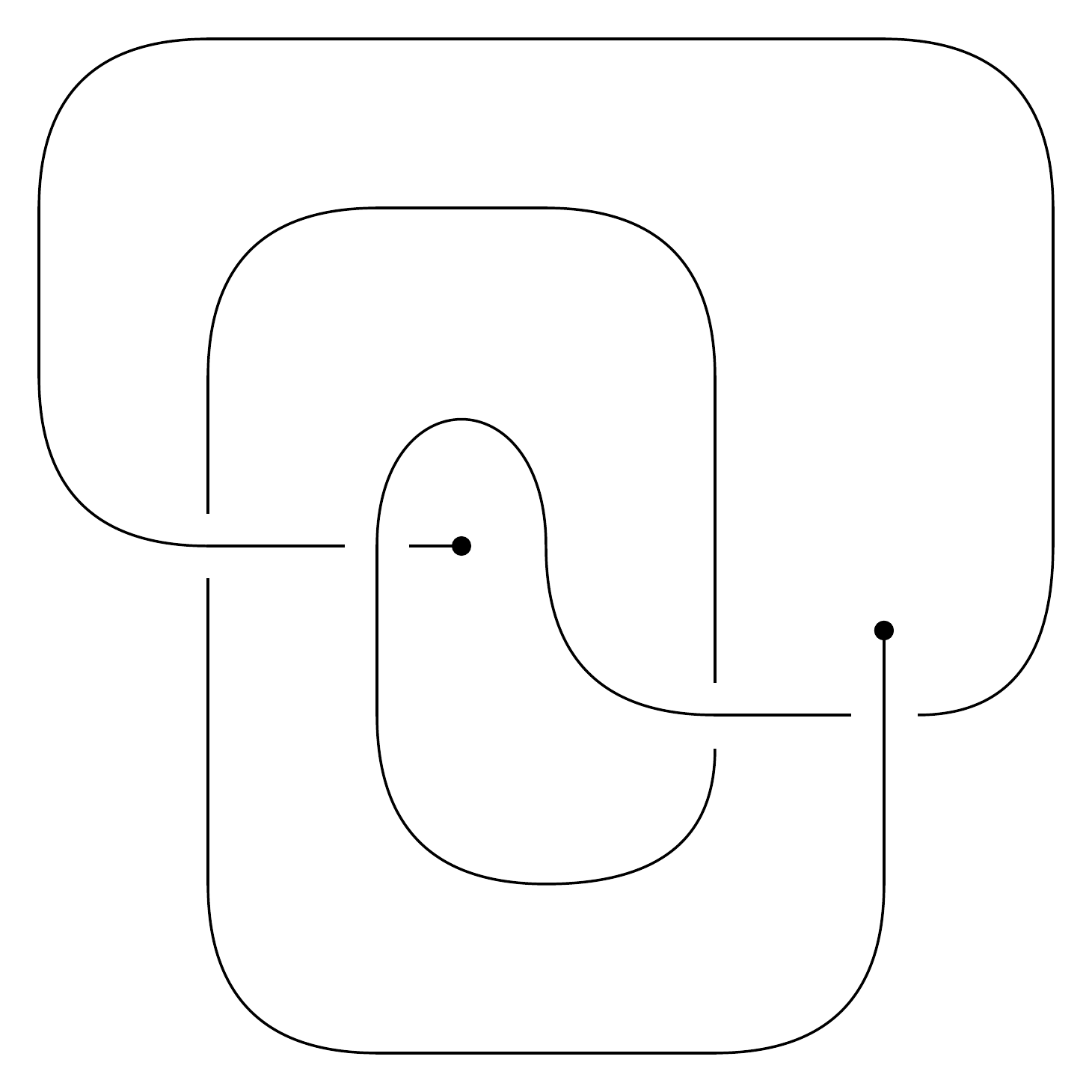}\\
\textcolor{black}{$4_{121}$}
\vspace{1cm}
\end{minipage}
\begin{minipage}[t]{.25\linewidth}
\centering
\includegraphics[width=0.9\textwidth,height=3.5cm,keepaspectratio]{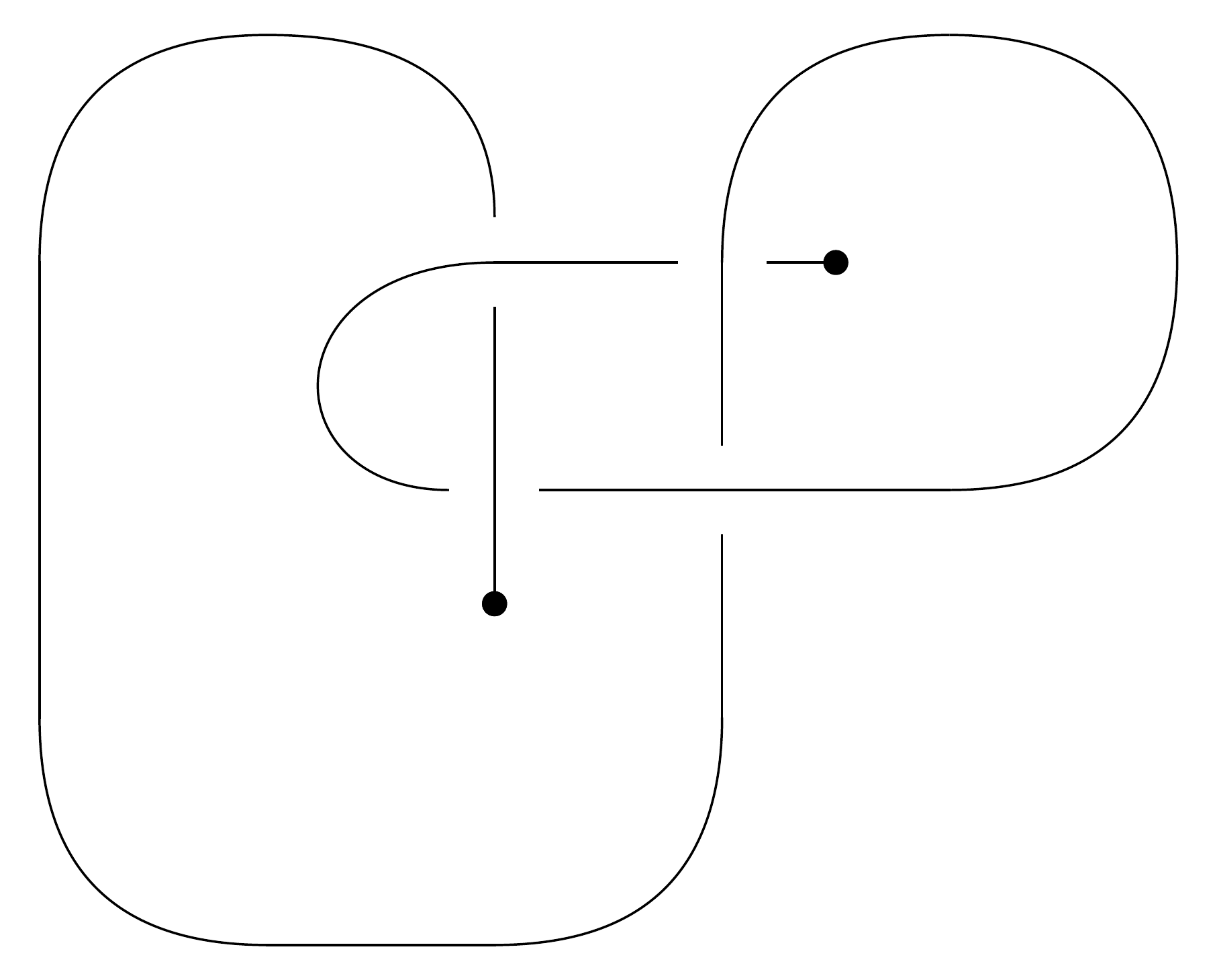}\\
\textcolor{black}{$4_{122}$}
\vspace{1cm}
\end{minipage}
\begin{minipage}[t]{.25\linewidth}
\centering
\includegraphics[width=0.9\textwidth,height=3.5cm,keepaspectratio]{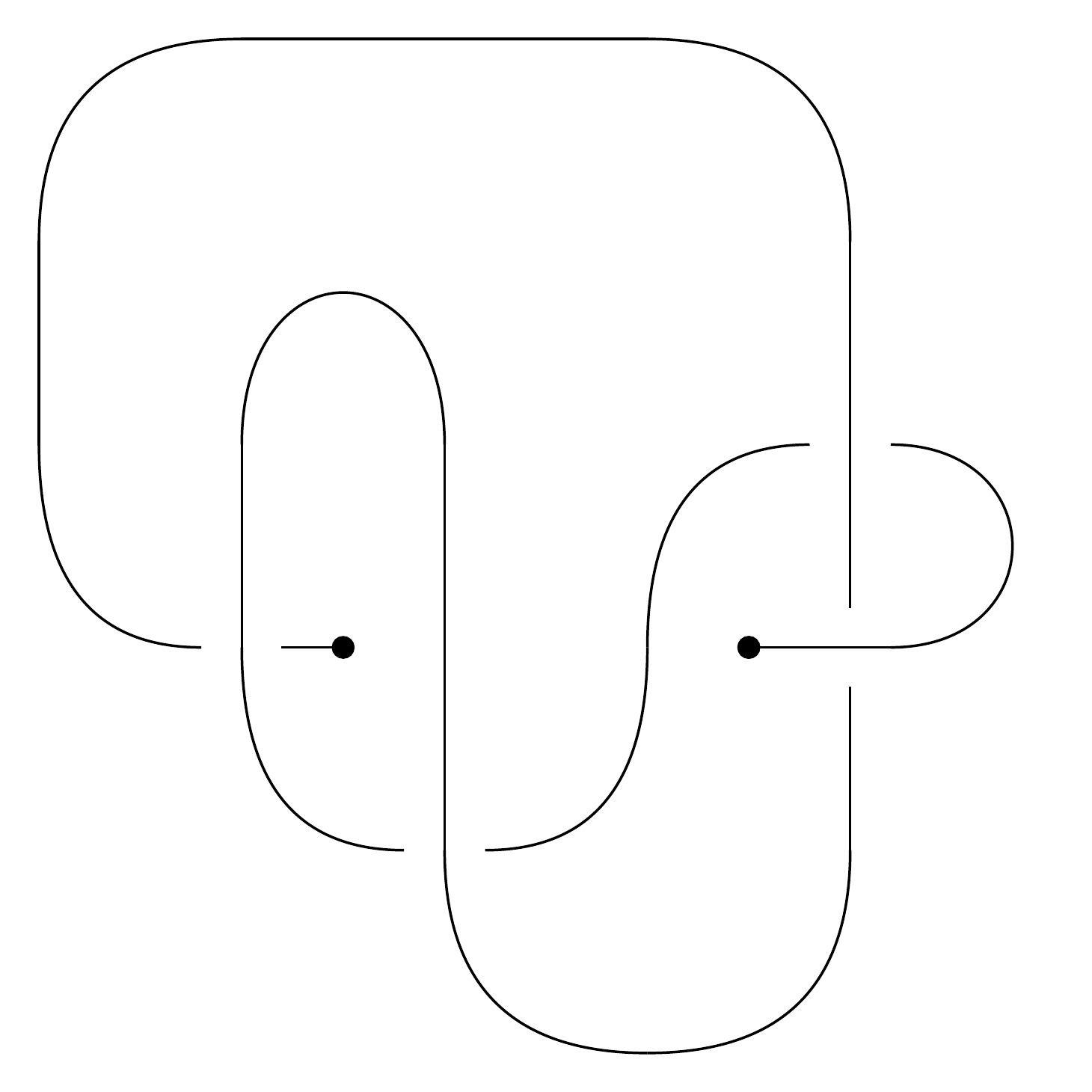}\\
\textcolor{black}{$4_{123}$}
\vspace{1cm}
\end{minipage}
\begin{minipage}[t]{.25\linewidth}
\centering
\includegraphics[width=0.9\textwidth,height=3.5cm,keepaspectratio]{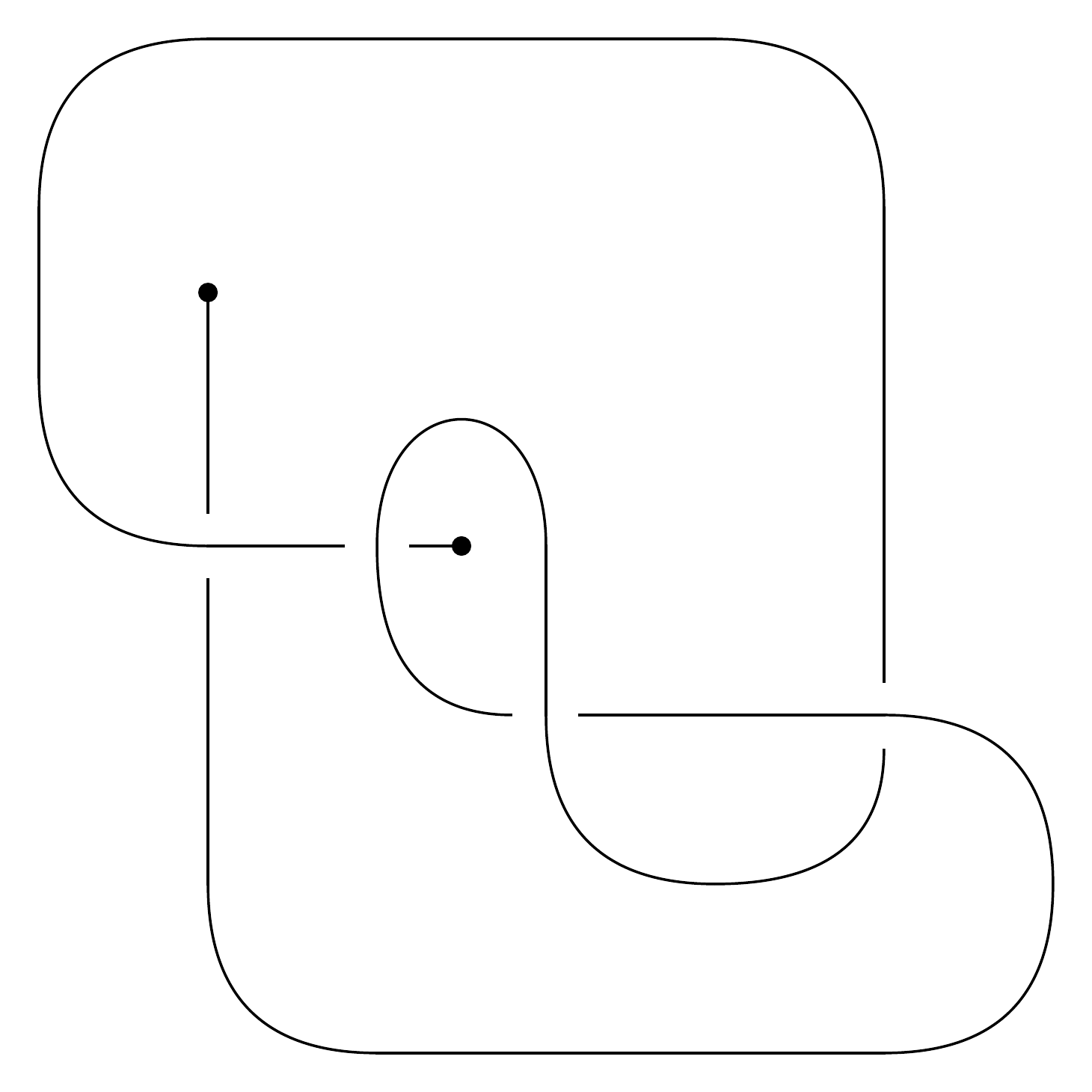}\\
\textcolor{black}{$4_{124}$}
\vspace{1cm}
\end{minipage}
\begin{minipage}[t]{.25\linewidth}
\centering
\includegraphics[width=0.9\textwidth,height=3.5cm,keepaspectratio]{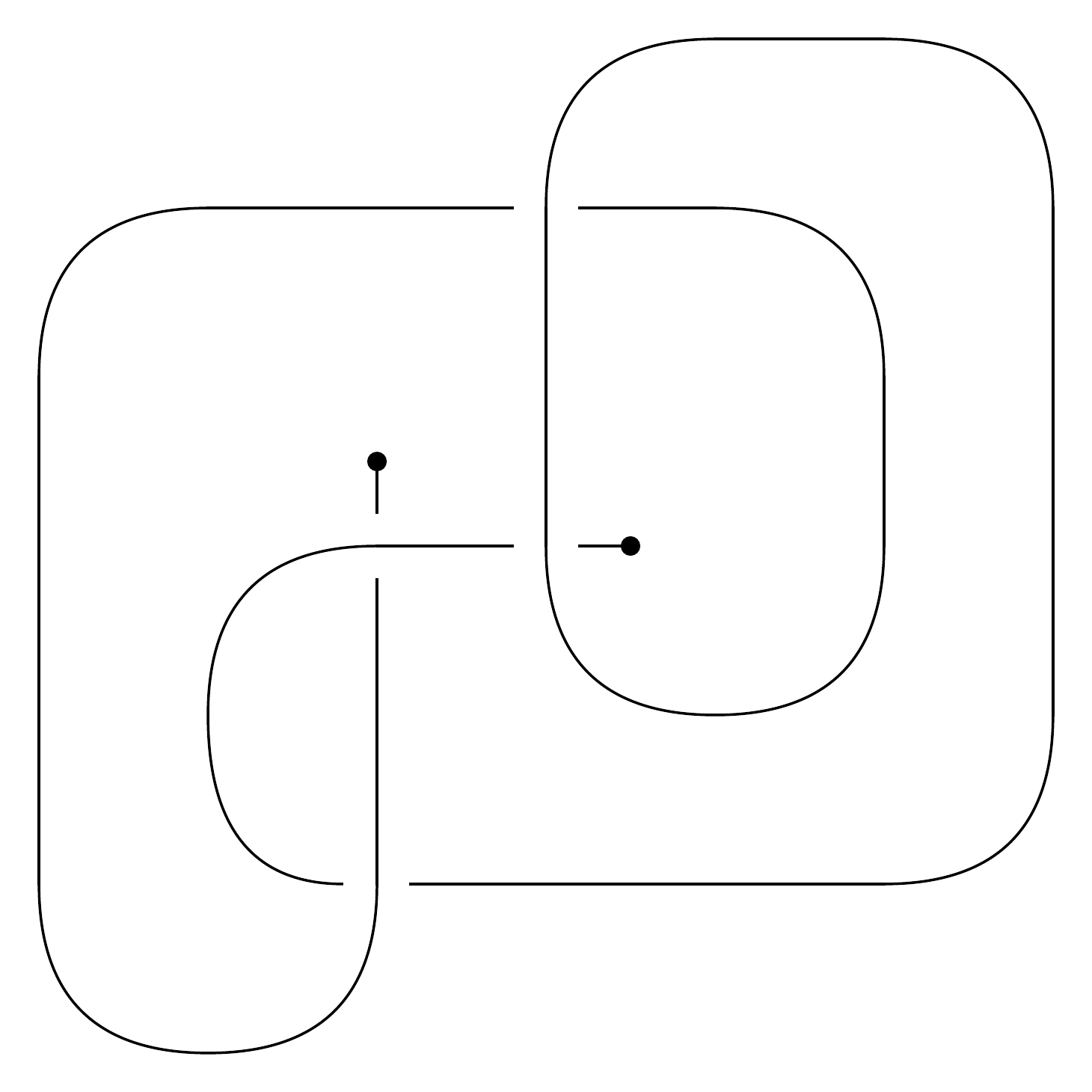}\\
\textcolor{black}{$4_{125}$}
\vspace{1cm}
\end{minipage}
\begin{minipage}[t]{.25\linewidth}
\centering
\includegraphics[width=0.9\textwidth,height=3.5cm,keepaspectratio]{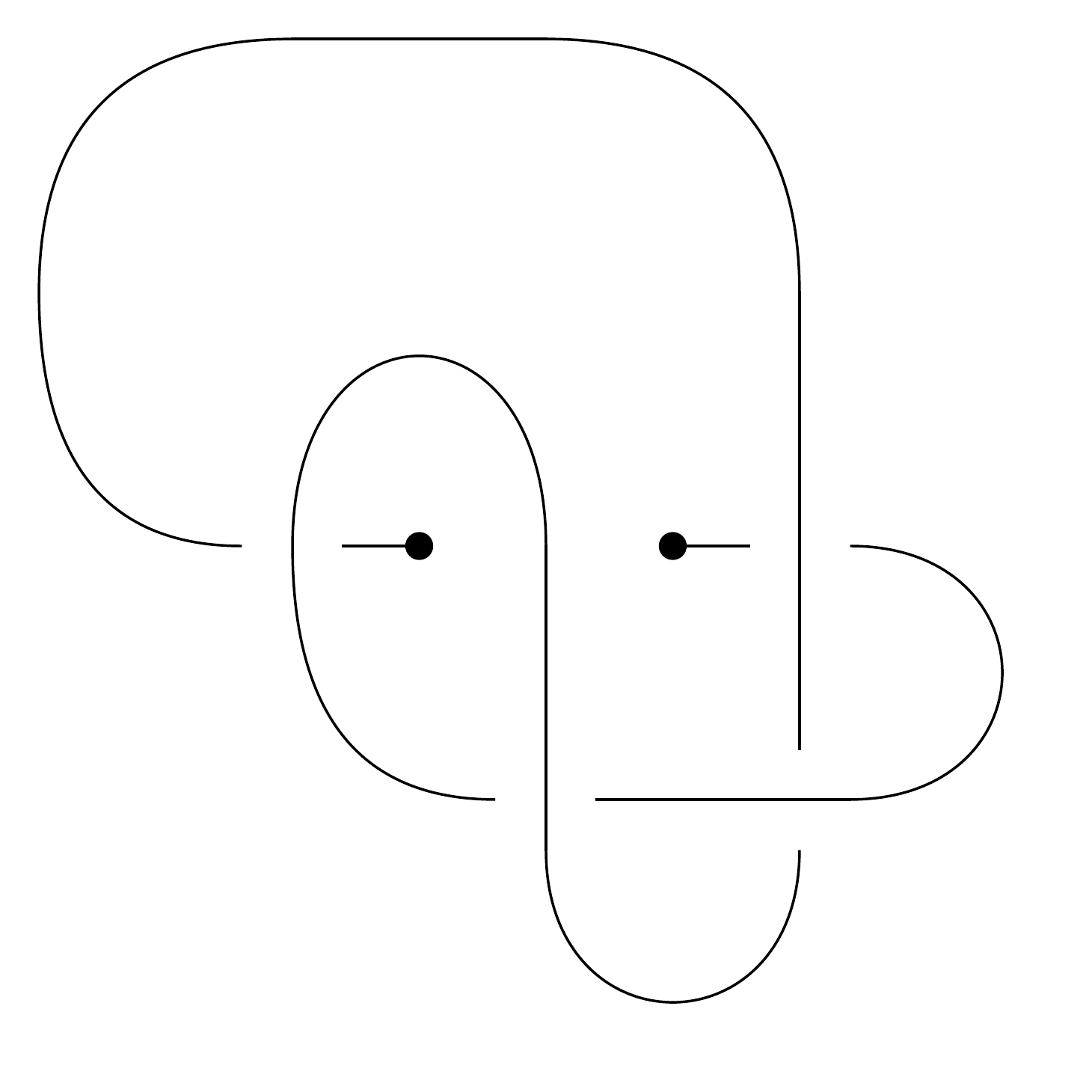}\\
\textcolor{black}{$4_{126}$}
\vspace{1cm}
\end{minipage}
\begin{minipage}[t]{.25\linewidth}
\centering
\includegraphics[width=0.9\textwidth,height=3.5cm,keepaspectratio]{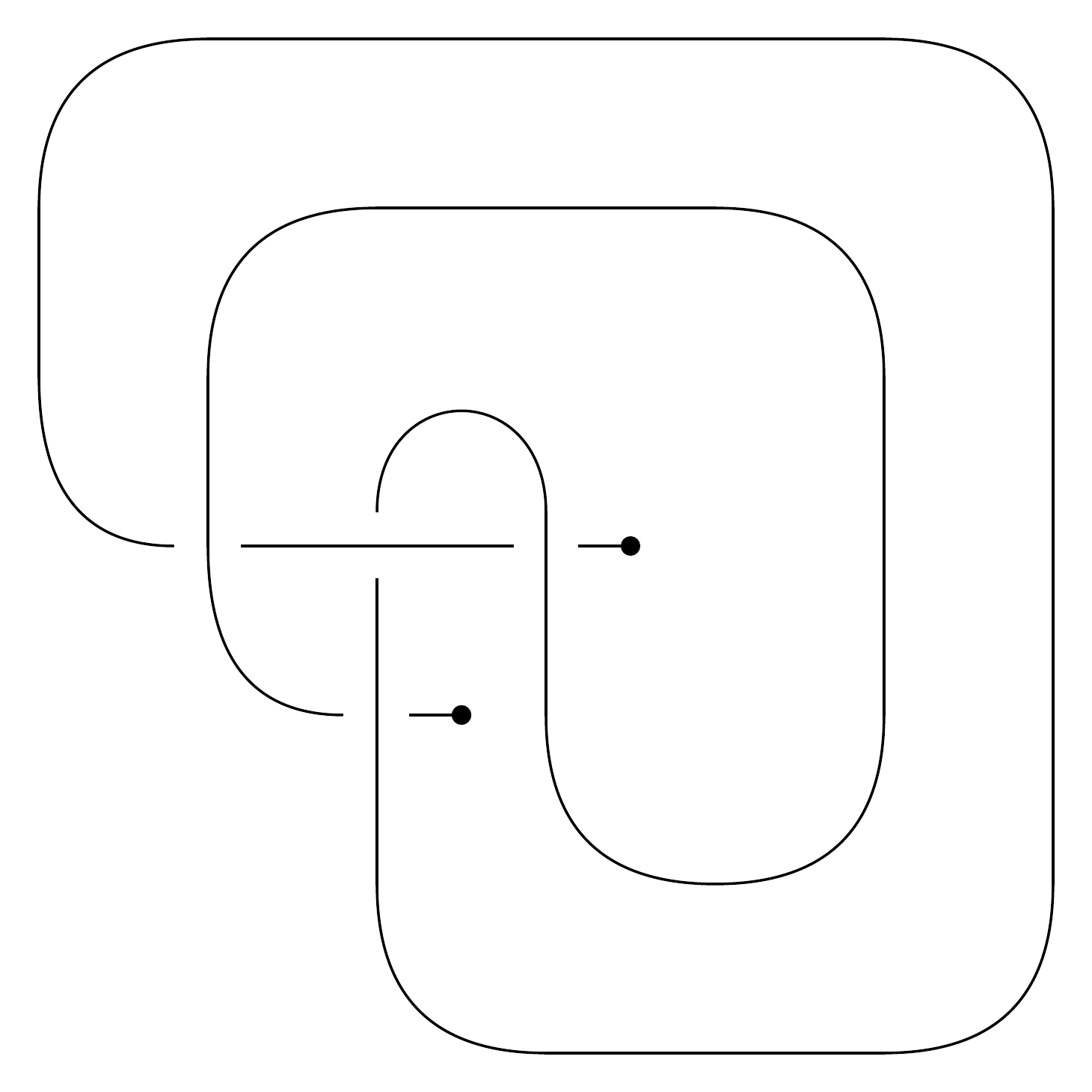}\\
\textcolor{black}{$4_{127}$}
\vspace{1cm}
\end{minipage}
\begin{minipage}[t]{.25\linewidth}
\centering
\includegraphics[width=0.9\textwidth,height=3.5cm,keepaspectratio]{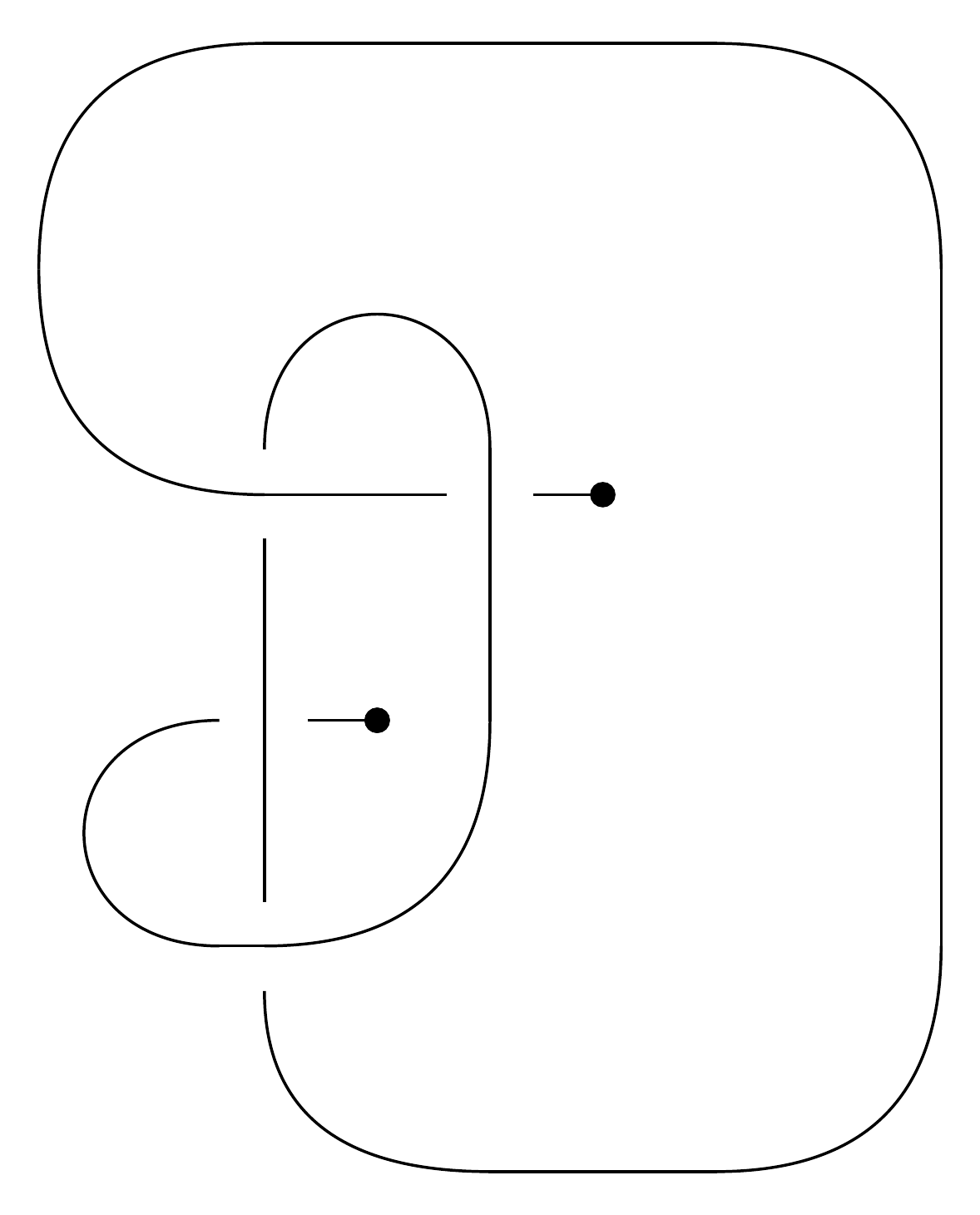}\\
\textcolor{black}{$4_{128}$}
\vspace{1cm}
\end{minipage}
\begin{minipage}[t]{.25\linewidth}
\centering
\includegraphics[width=0.9\textwidth,height=3.5cm,keepaspectratio]{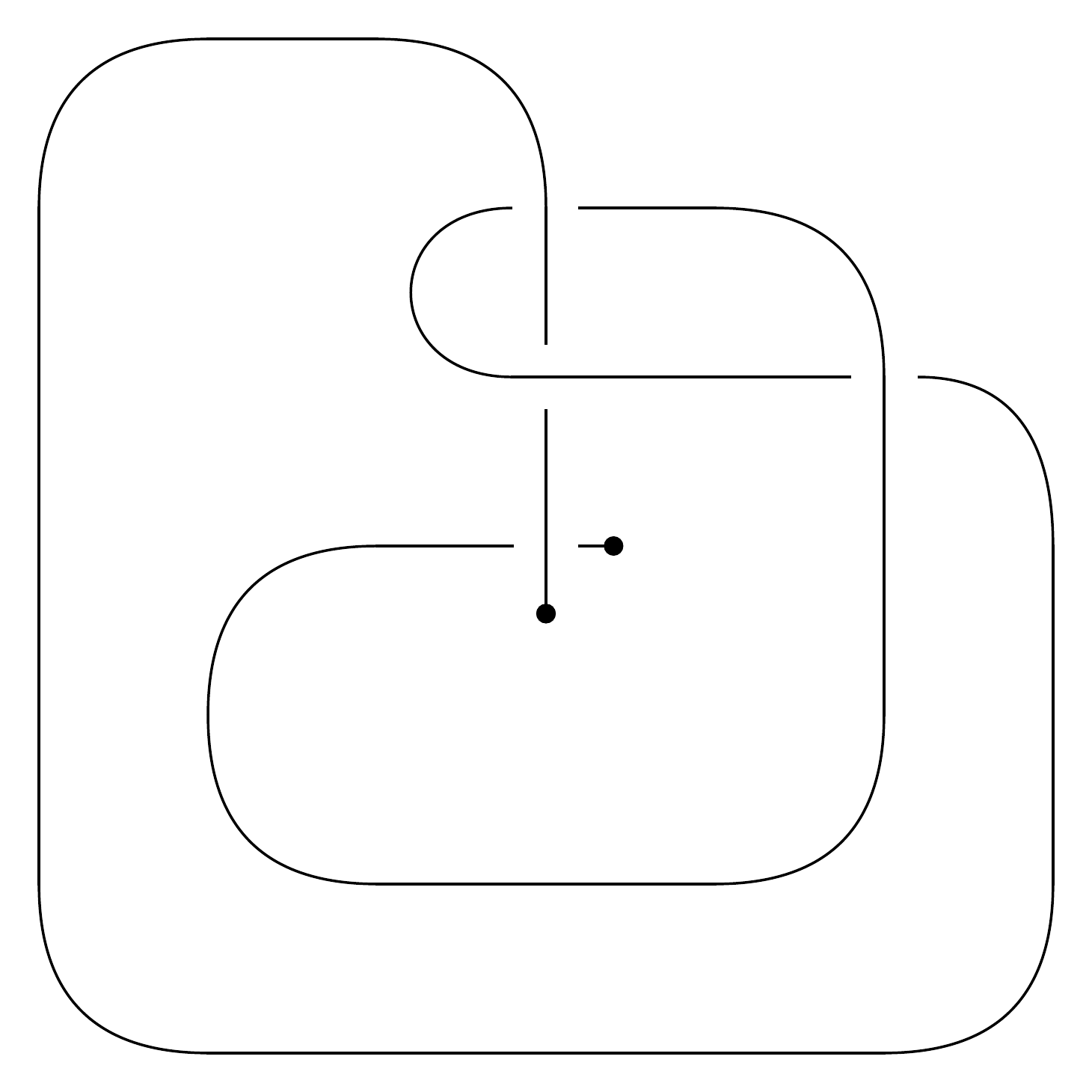}\\
\textcolor{black}{$4_{129}$}
\vspace{1cm}
\end{minipage}
\begin{minipage}[t]{.25\linewidth}
\centering
\includegraphics[width=0.9\textwidth,height=3.5cm,keepaspectratio]{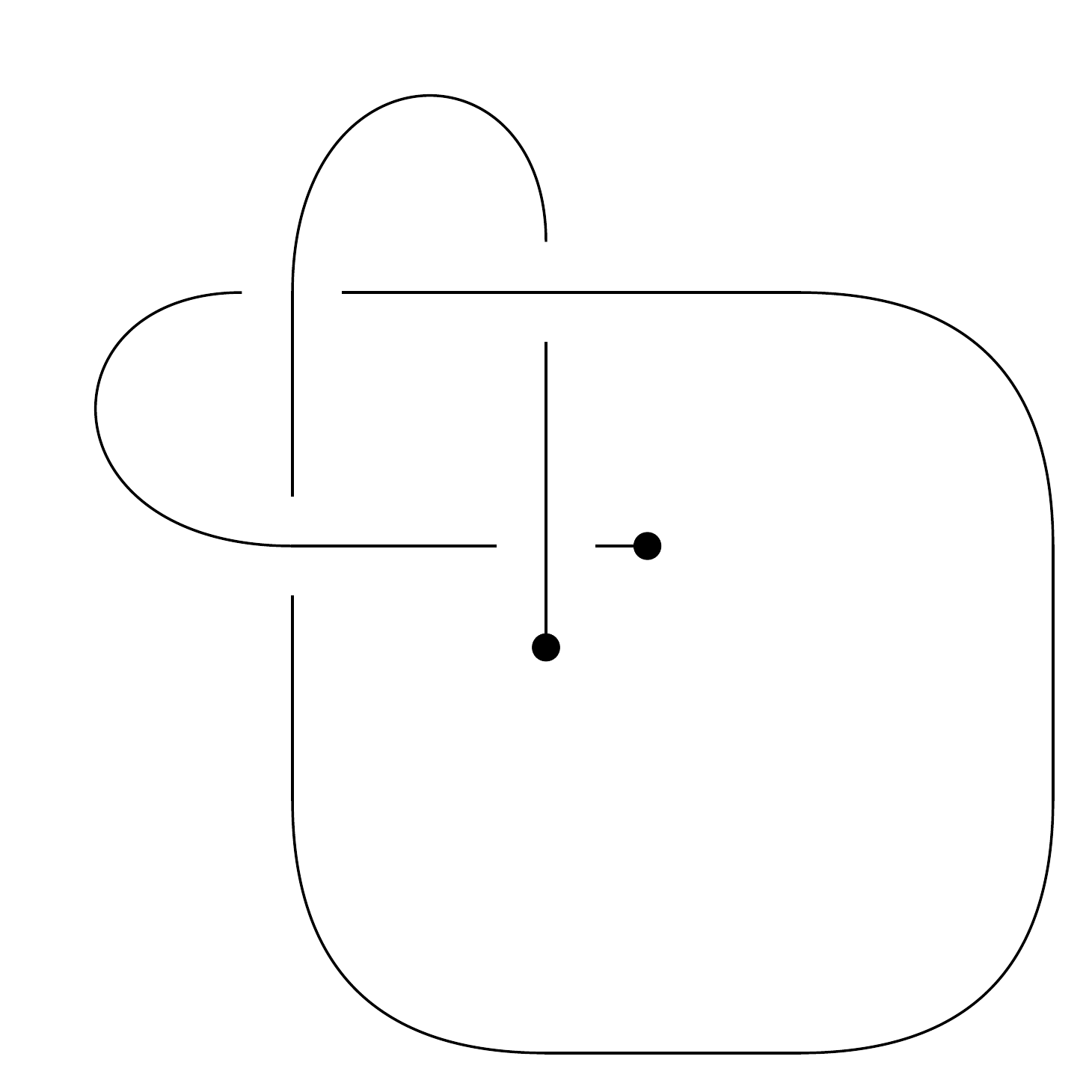}\\
\textcolor{black}{$4_{130}$}
\vspace{1cm}
\end{minipage}
\begin{minipage}[t]{.25\linewidth}
\centering
\includegraphics[width=0.9\textwidth,height=3.5cm,keepaspectratio]{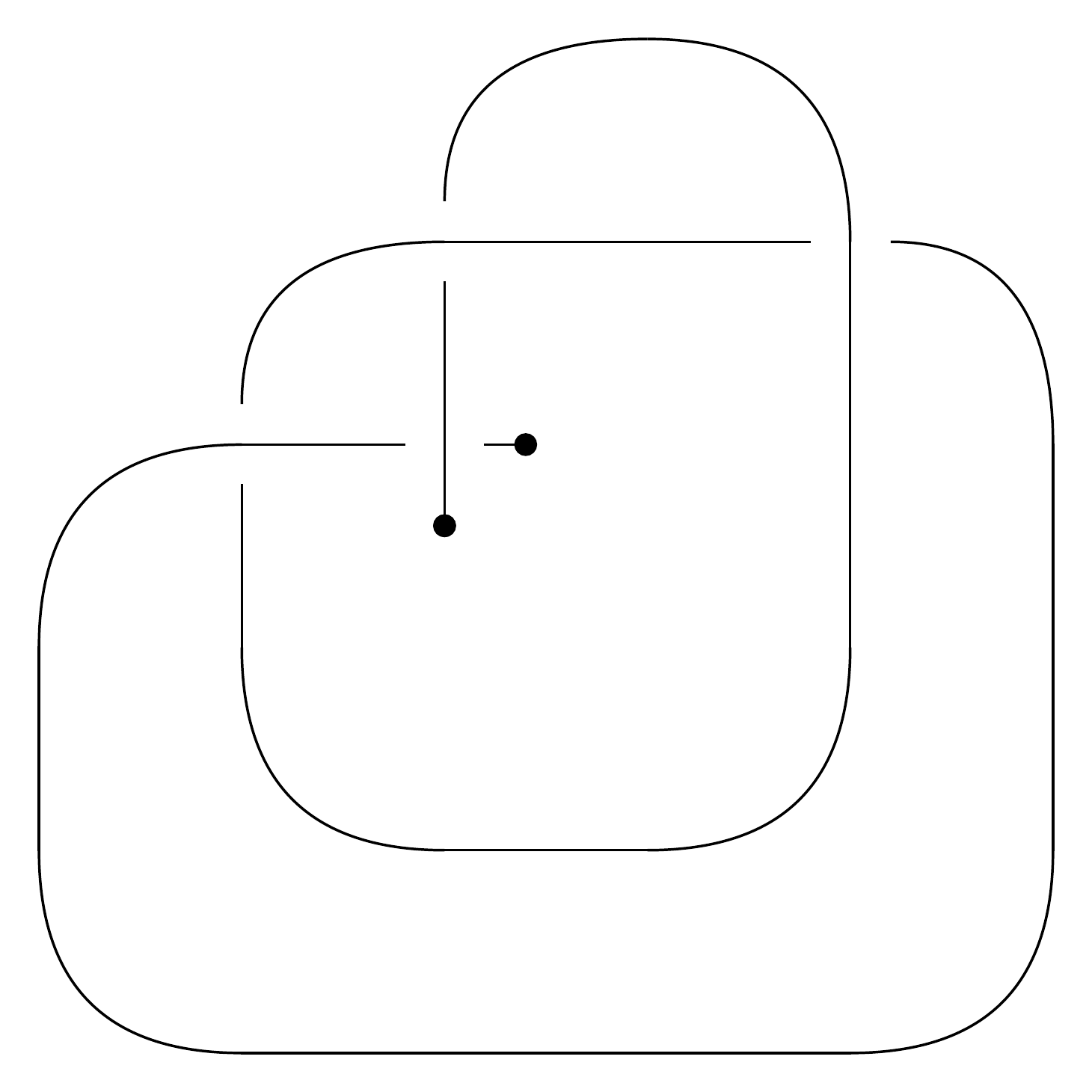}\\
\textcolor{black}{$4_{131}$}
\vspace{1cm}
\end{minipage}
\begin{minipage}[t]{.25\linewidth}
\centering
\includegraphics[width=0.9\textwidth,height=3.5cm,keepaspectratio]{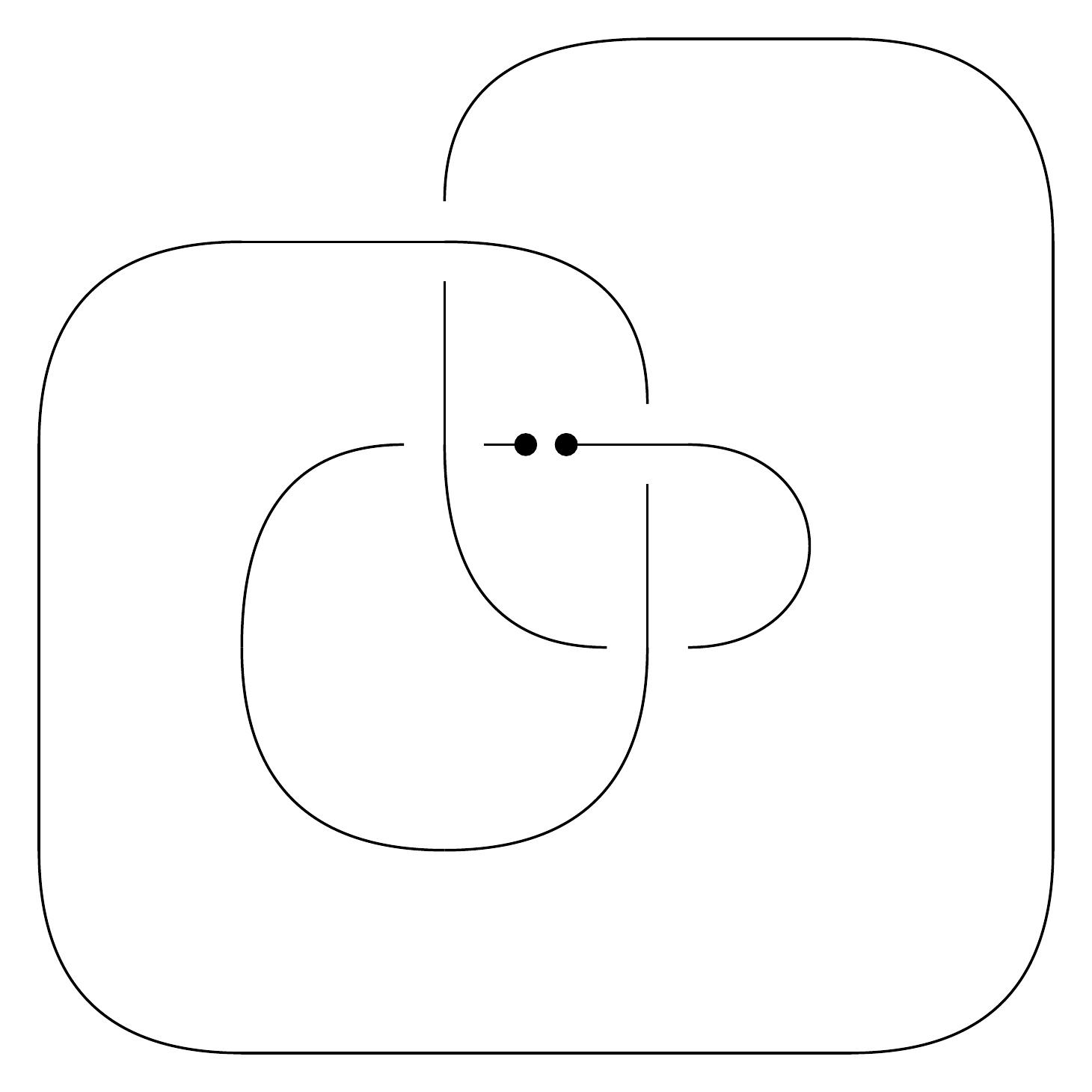}\\
\textcolor{black}{$4_{132}$}
\vspace{1cm}
\end{minipage}
\begin{minipage}[t]{.25\linewidth}
\centering
\includegraphics[width=0.9\textwidth,height=3.5cm,keepaspectratio]{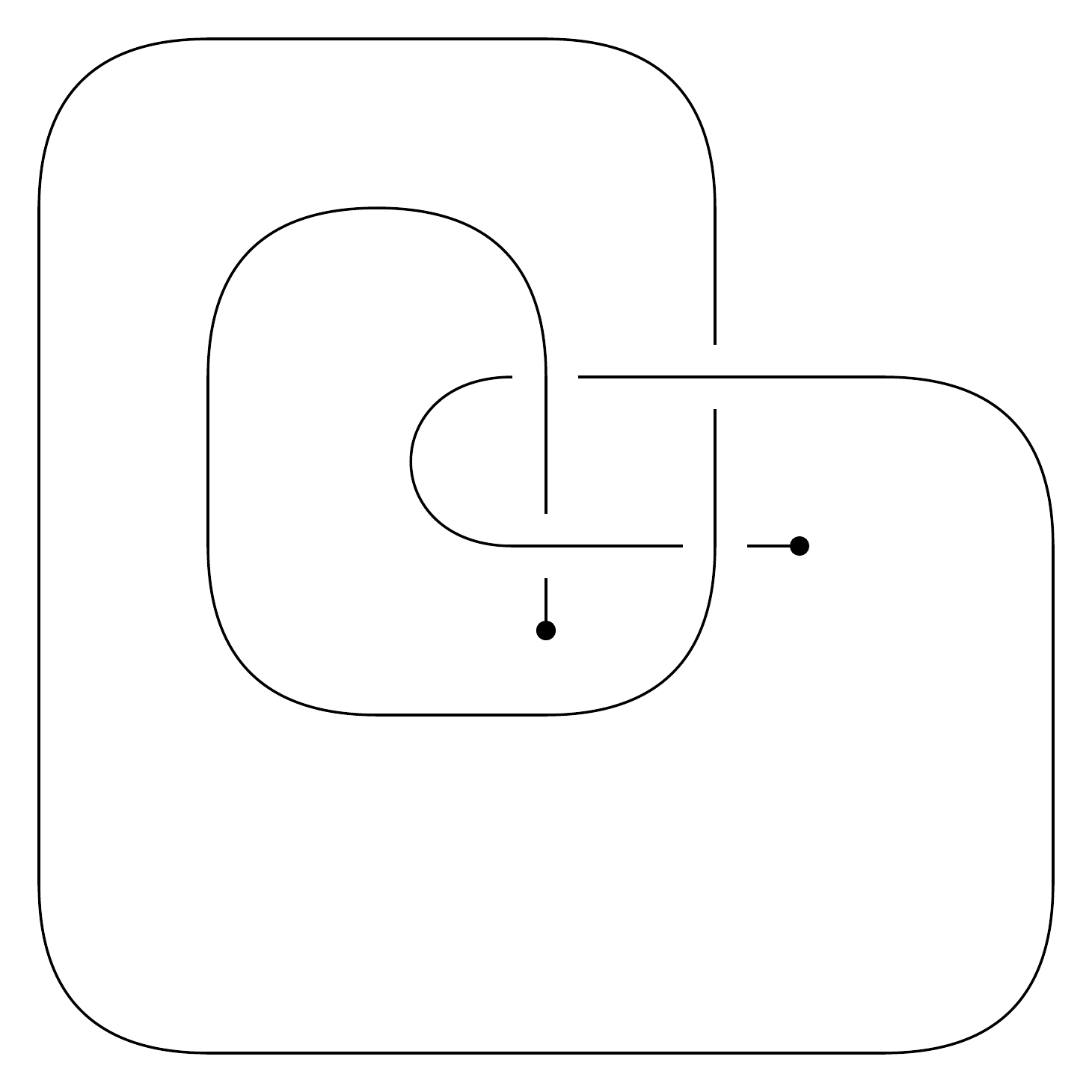}\\
\textcolor{black}{$4_{133}$}
\vspace{1cm}
\end{minipage}
\begin{minipage}[t]{.25\linewidth}
\centering
\includegraphics[width=0.9\textwidth,height=3.5cm,keepaspectratio]{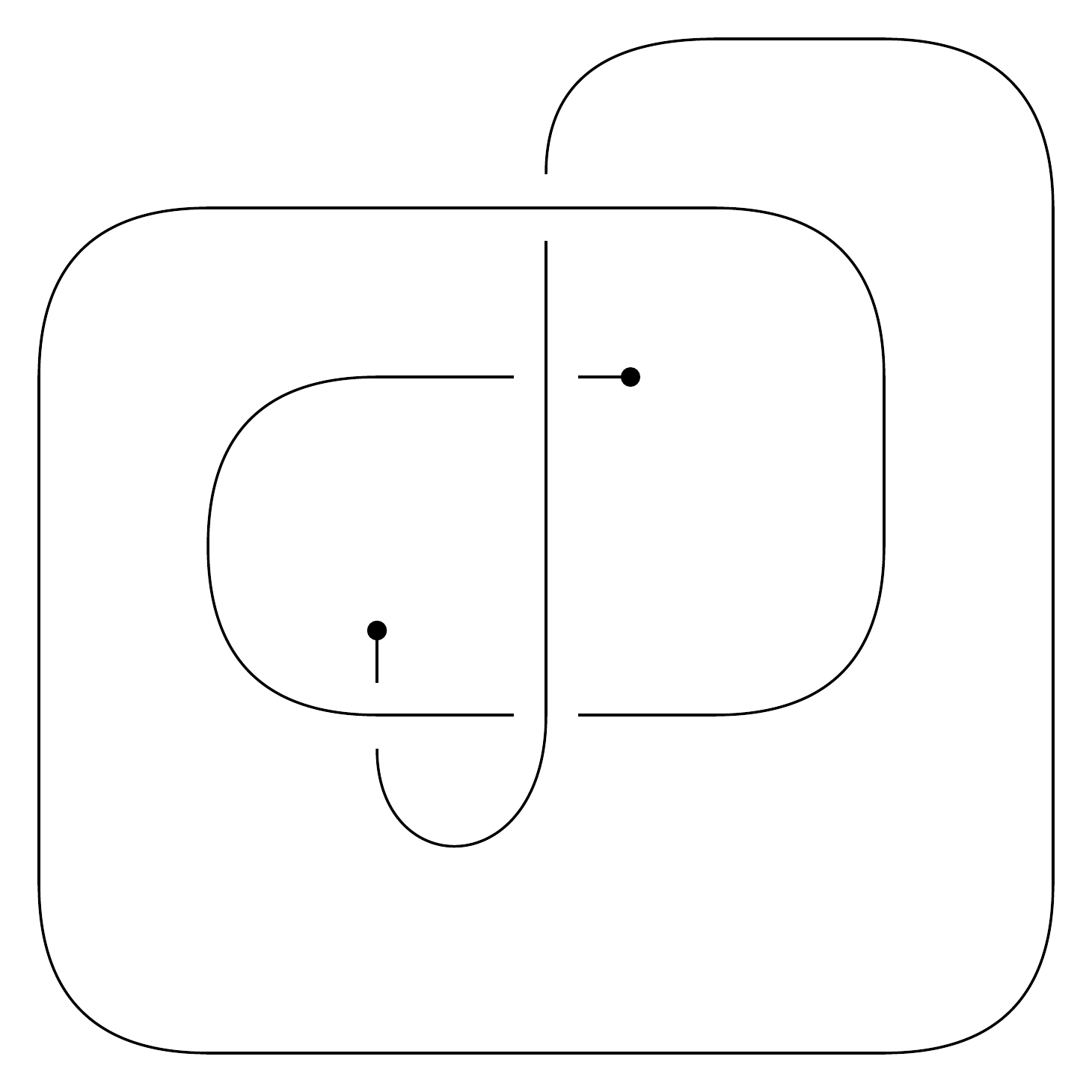}\\
\textcolor{black}{$4_{134}$}
\vspace{1cm}
\end{minipage}
\begin{minipage}[t]{.25\linewidth}
\centering
\includegraphics[width=0.9\textwidth,height=3.5cm,keepaspectratio]{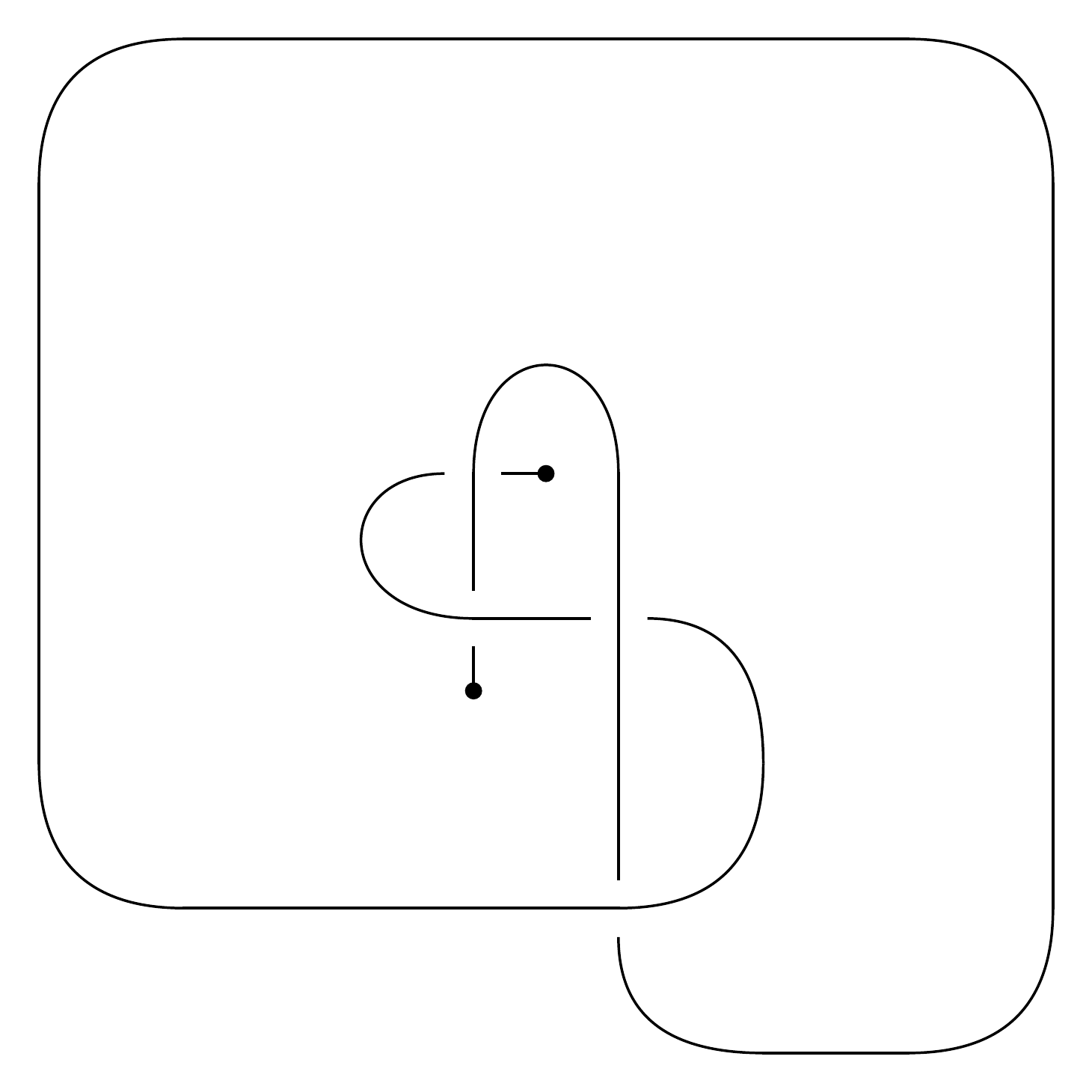}\\
\textcolor{black}{$4_{135}$}
\vspace{1cm}
\end{minipage}
\begin{minipage}[t]{.25\linewidth}
\centering
\includegraphics[width=0.9\textwidth,height=3.5cm,keepaspectratio]{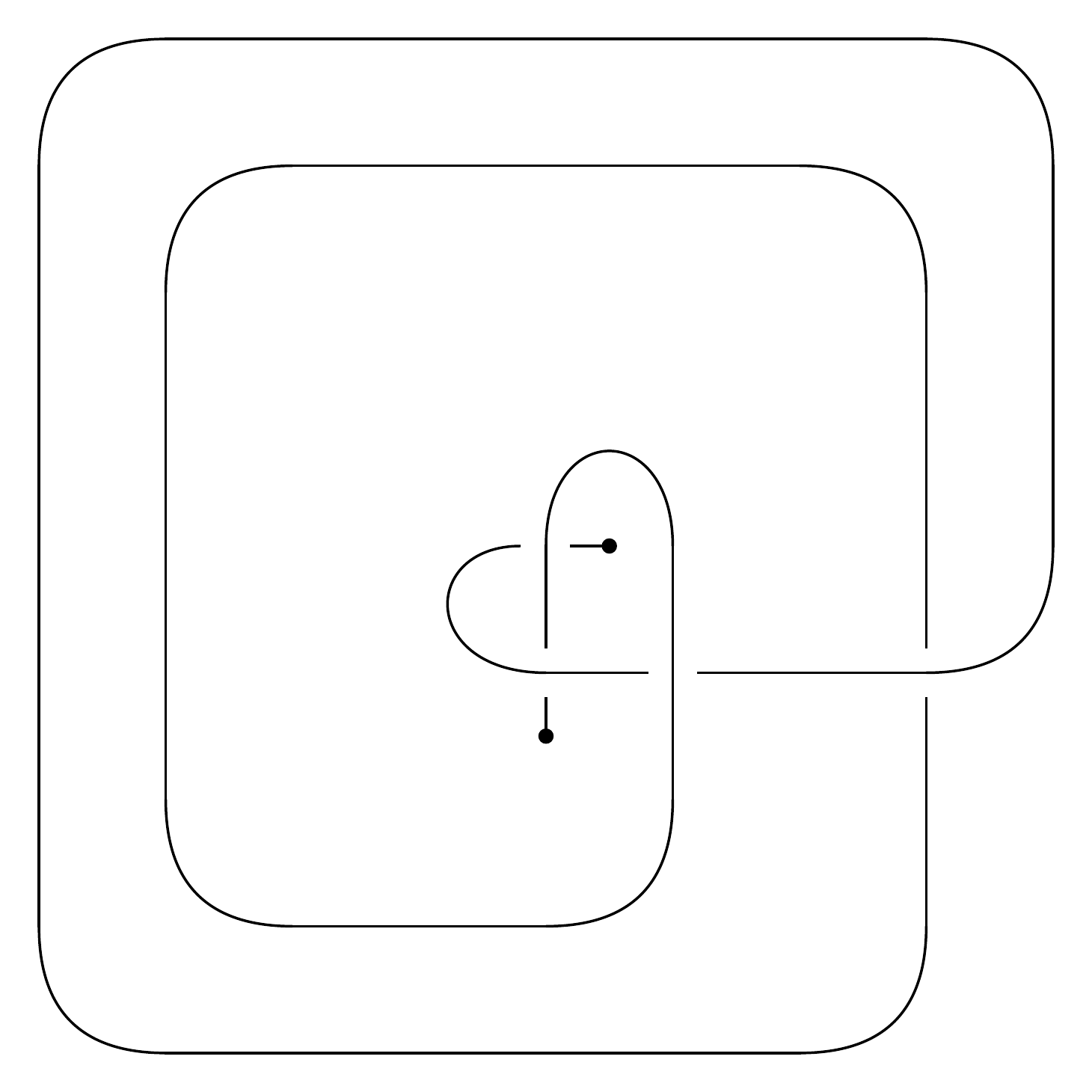}\\
\textcolor{black}{$4_{136}$}
\vspace{1cm}
\end{minipage}
\begin{minipage}[t]{.25\linewidth}
\centering
\includegraphics[width=0.9\textwidth,height=3.5cm,keepaspectratio]{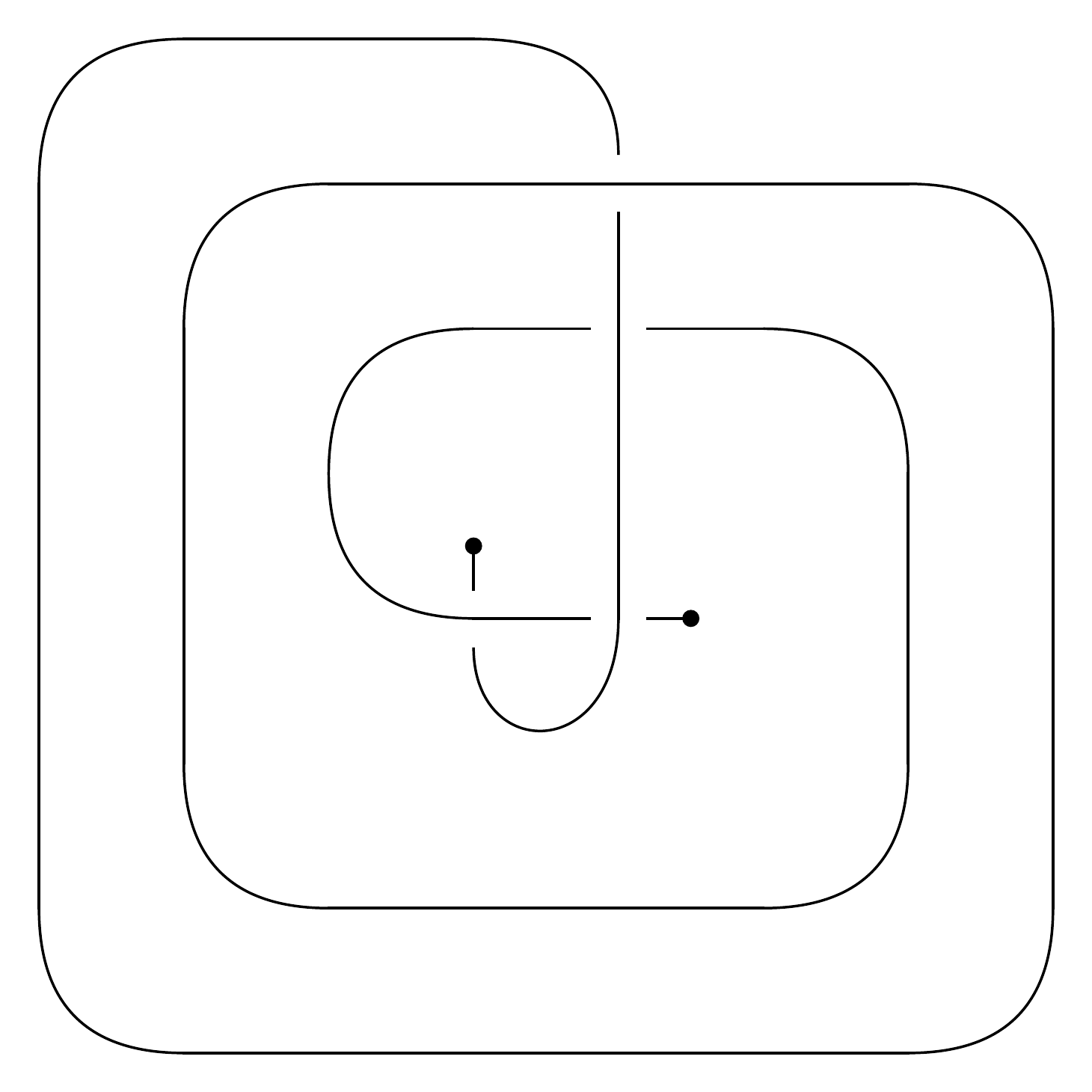}\\
\textcolor{black}{$4_{137}$}
\vspace{1cm}
\end{minipage}
\begin{minipage}[t]{.25\linewidth}
\centering
\includegraphics[width=0.9\textwidth,height=3.5cm,keepaspectratio]{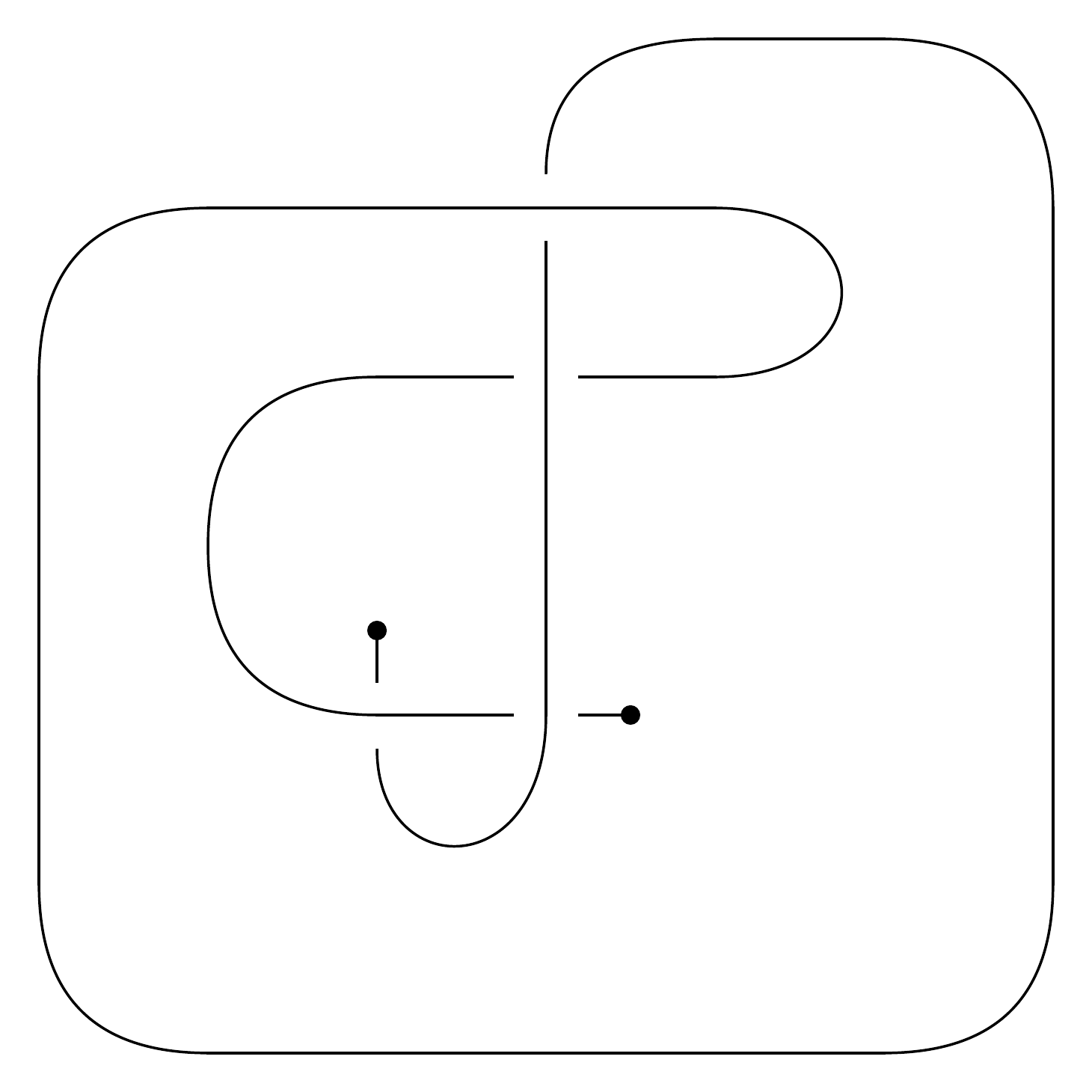}\\
\textcolor{black}{$4_{138}$}
\vspace{1cm}
\end{minipage}
\begin{minipage}[t]{.25\linewidth}
\centering
\includegraphics[width=0.9\textwidth,height=3.5cm,keepaspectratio]{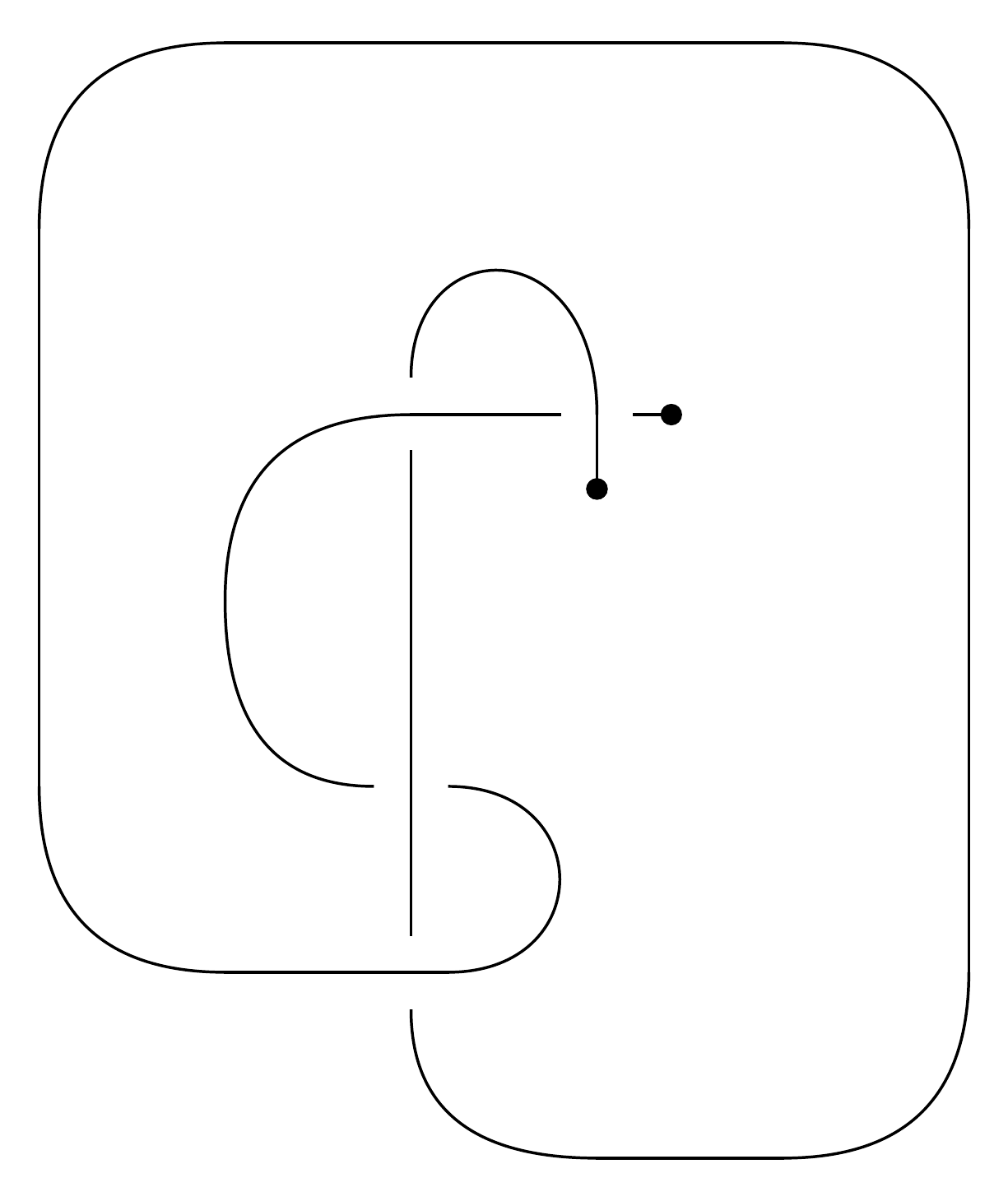}\\
\textcolor{black}{$4_{139}$}
\vspace{1cm}
\end{minipage}
\begin{minipage}[t]{.25\linewidth}
\centering
\includegraphics[width=0.9\textwidth,height=3.5cm,keepaspectratio]{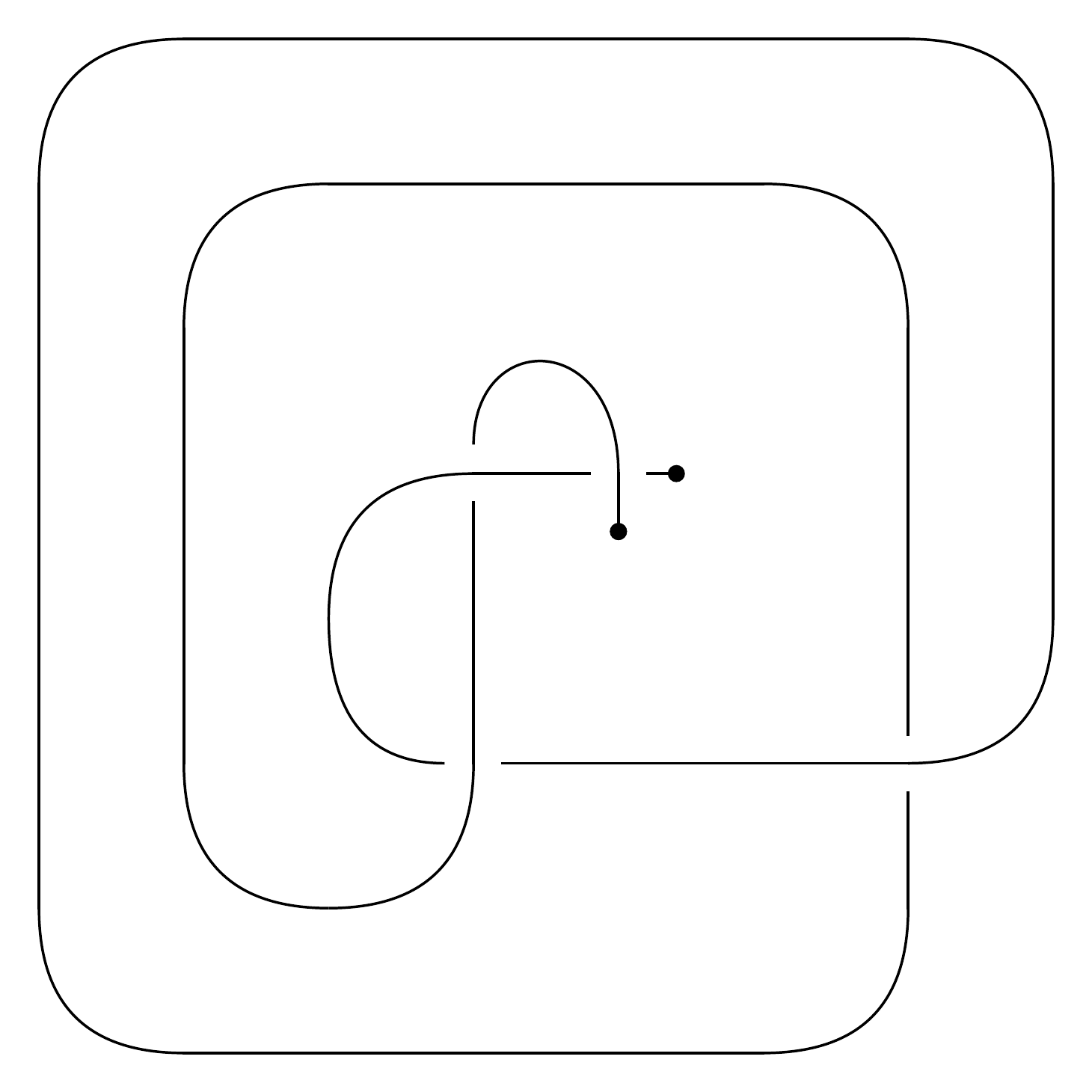}\\
\textcolor{black}{$4_{140}$}
\vspace{1cm}
\end{minipage}
\begin{minipage}[t]{.25\linewidth}
\centering
\includegraphics[width=0.9\textwidth,height=3.5cm,keepaspectratio]{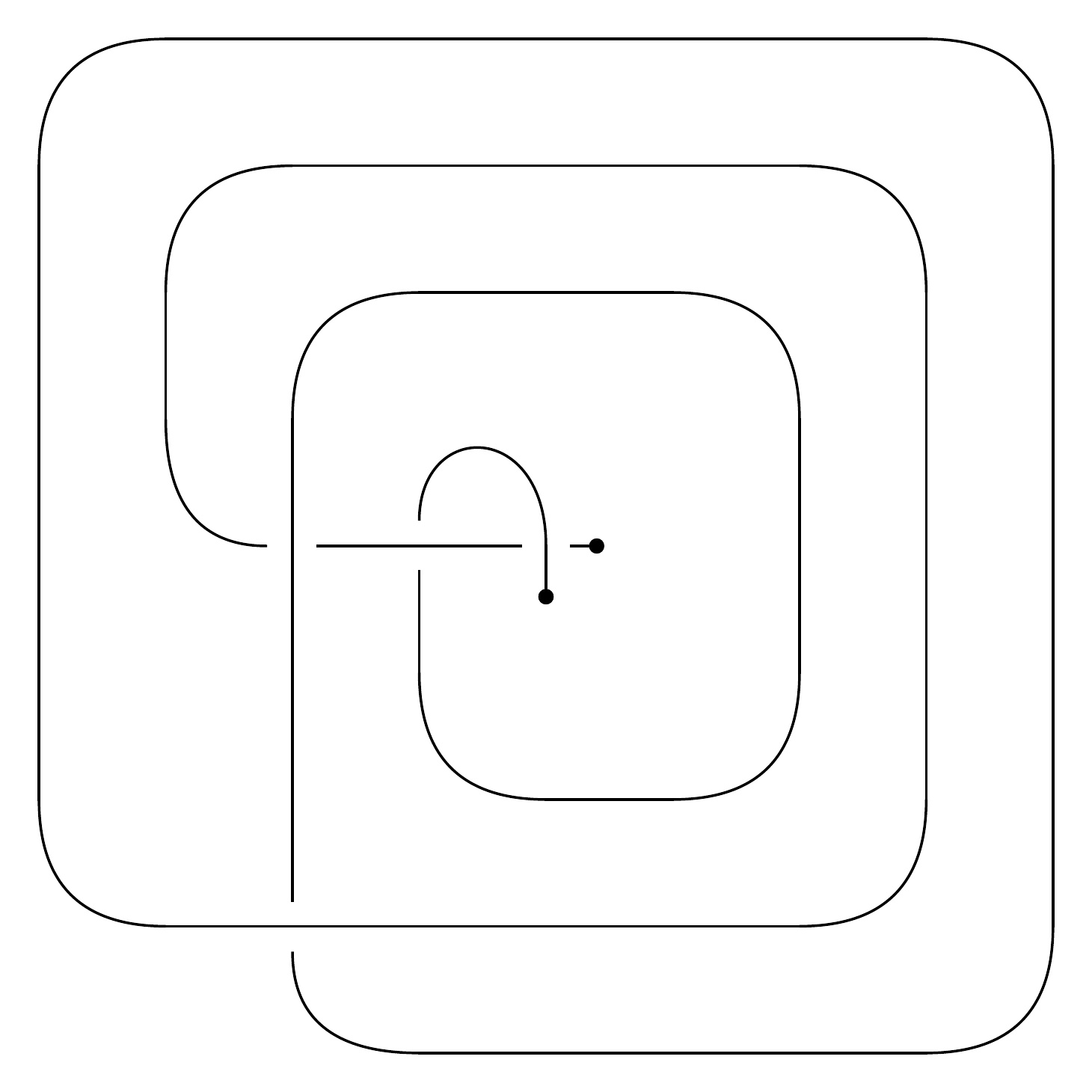}\\
\textcolor{black}{$4_{141}$}
\vspace{1cm}
\end{minipage}
\begin{minipage}[t]{.25\linewidth}
\centering
\includegraphics[width=0.9\textwidth,height=3.5cm,keepaspectratio]{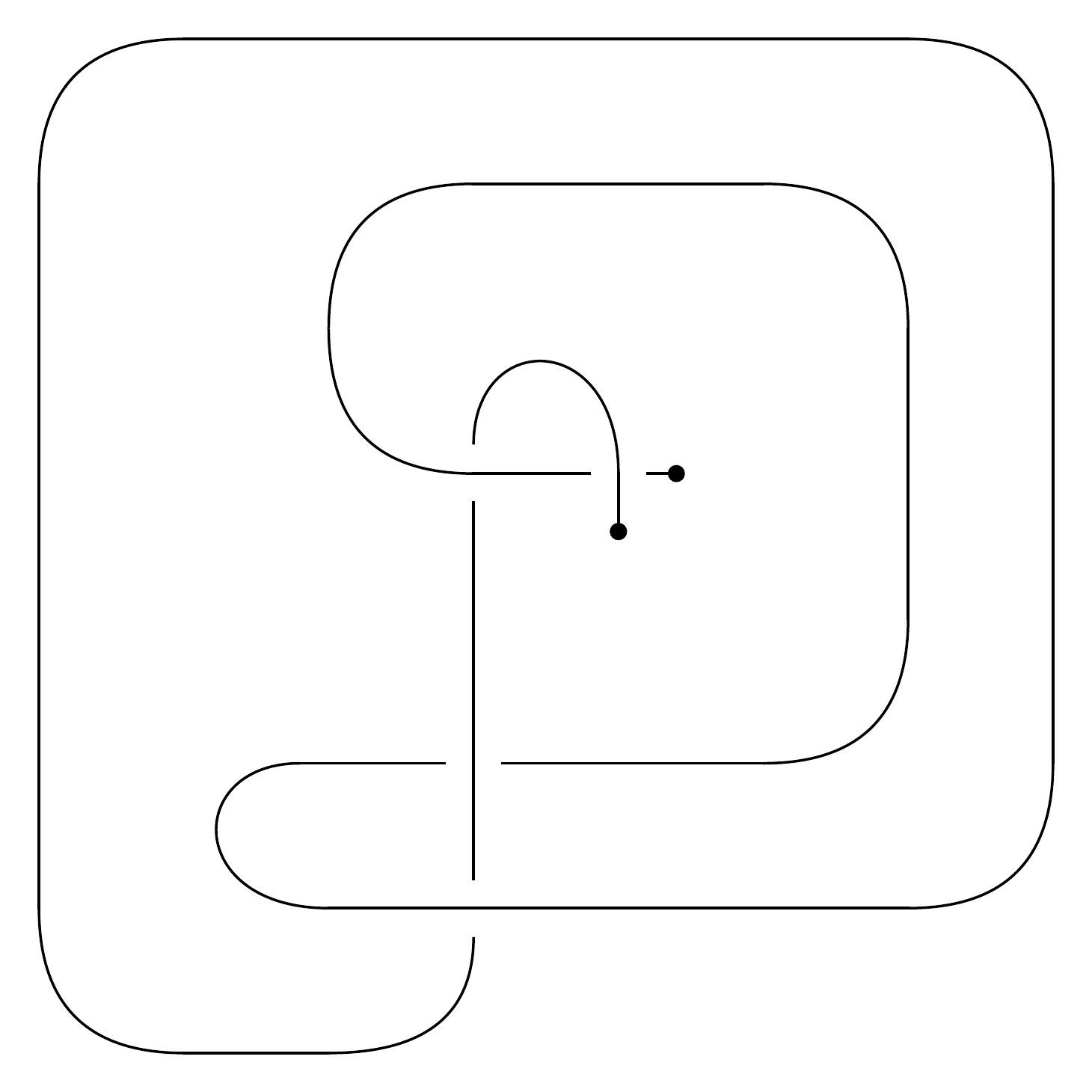}\\
\textcolor{black}{$4_{142}$}
\vspace{1cm}
\end{minipage}
\begin{minipage}[t]{.25\linewidth}
\centering
\includegraphics[width=0.9\textwidth,height=3.5cm,keepaspectratio]{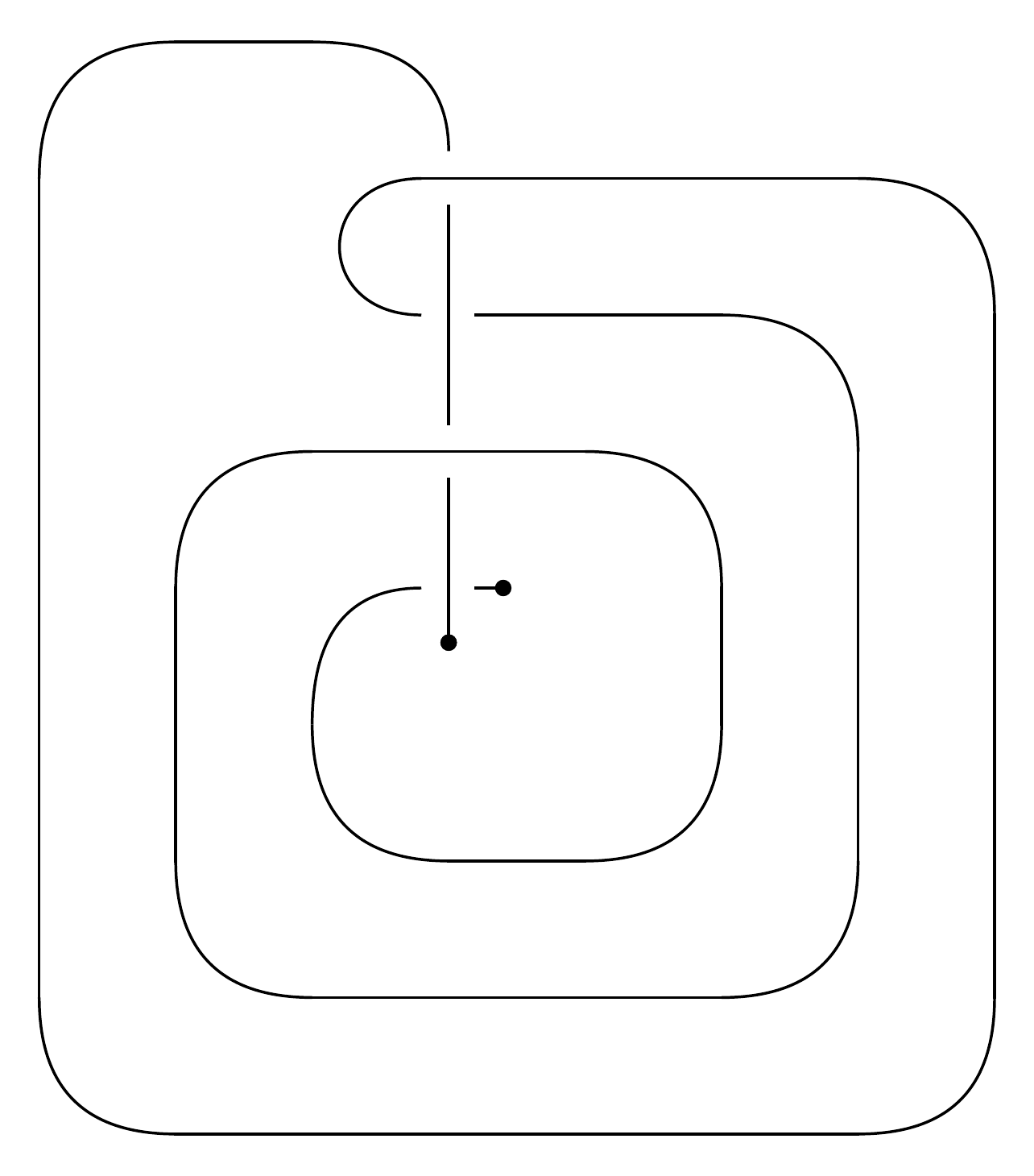}\\
\textcolor{black}{$4_{143}$}
\vspace{1cm}
\end{minipage}
\begin{minipage}[t]{.25\linewidth}
\centering
\includegraphics[width=0.9\textwidth,height=3.5cm,keepaspectratio]{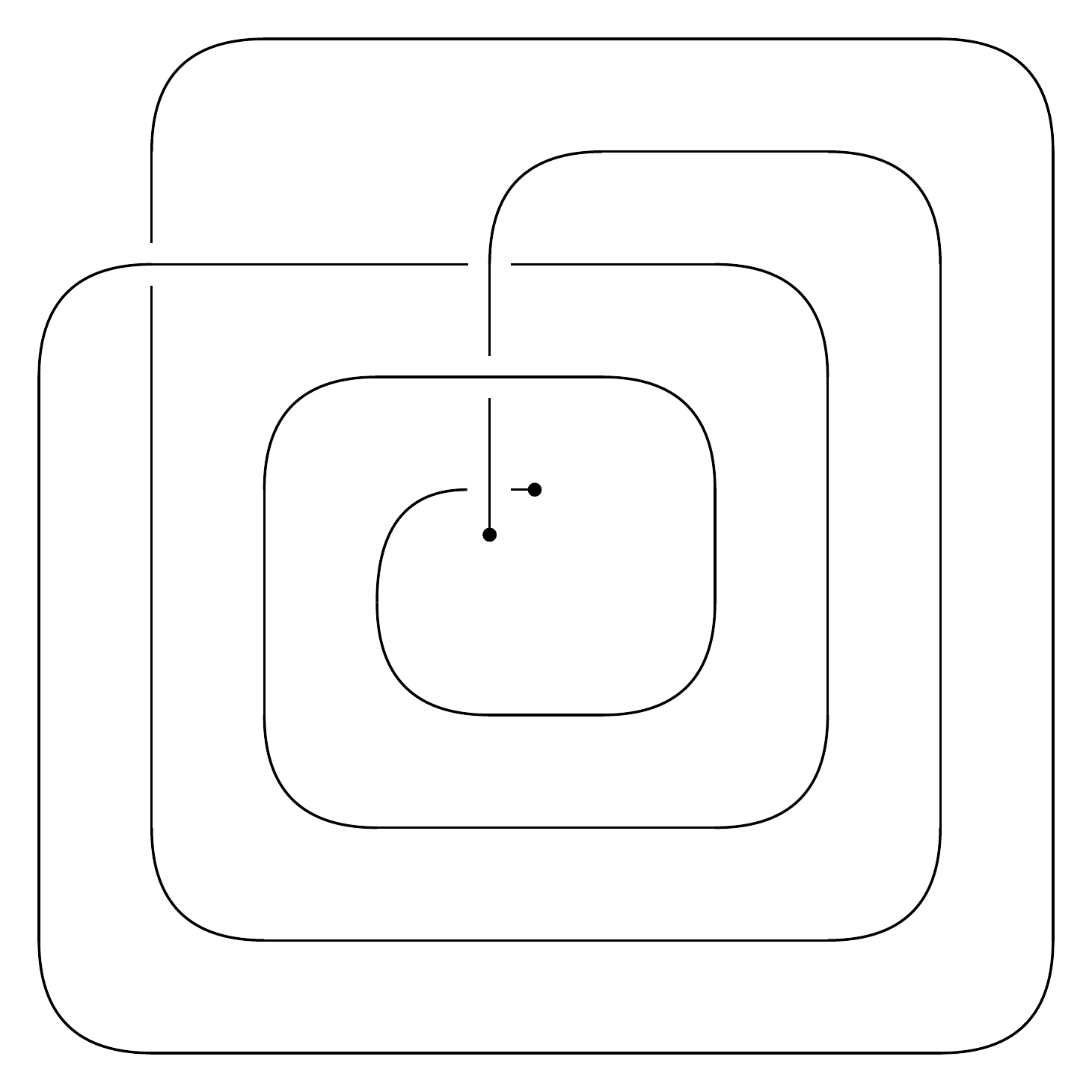}\\
\textcolor{black}{$4_{144}$}
\vspace{1cm}
\end{minipage}
\begin{minipage}[t]{.25\linewidth}
\centering
\includegraphics[width=0.9\textwidth,height=3.5cm,keepaspectratio]{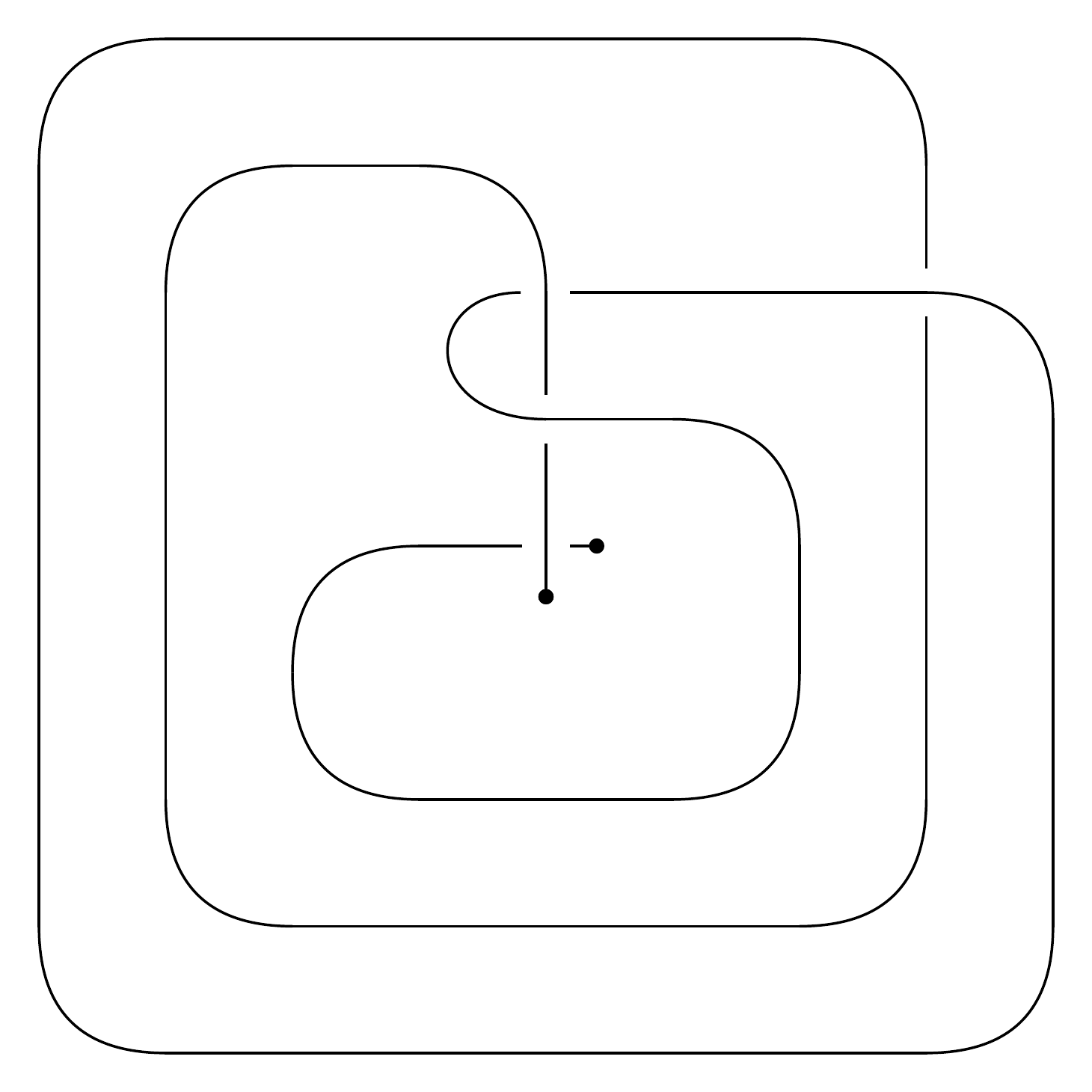}\\
\textcolor{black}{$4_{145}$}
\vspace{1cm}
\end{minipage}
\begin{minipage}[t]{.25\linewidth}
\centering
\includegraphics[width=0.9\textwidth,height=3.5cm,keepaspectratio]{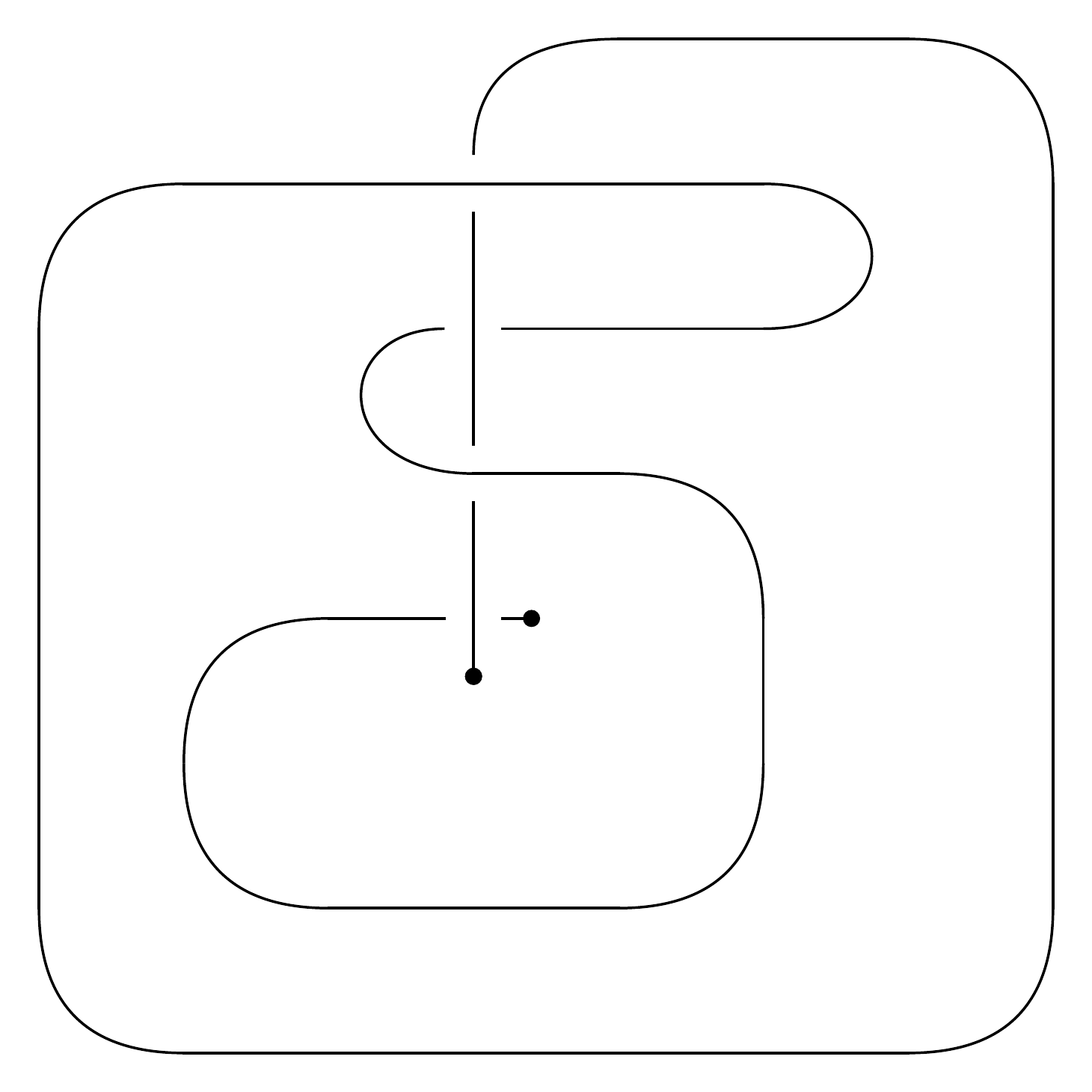}\\
\textcolor{black}{$4_{146}$}
\vspace{1cm}
\end{minipage}
\begin{minipage}[t]{.25\linewidth}
\centering
\includegraphics[width=0.9\textwidth,height=3.5cm,keepaspectratio]{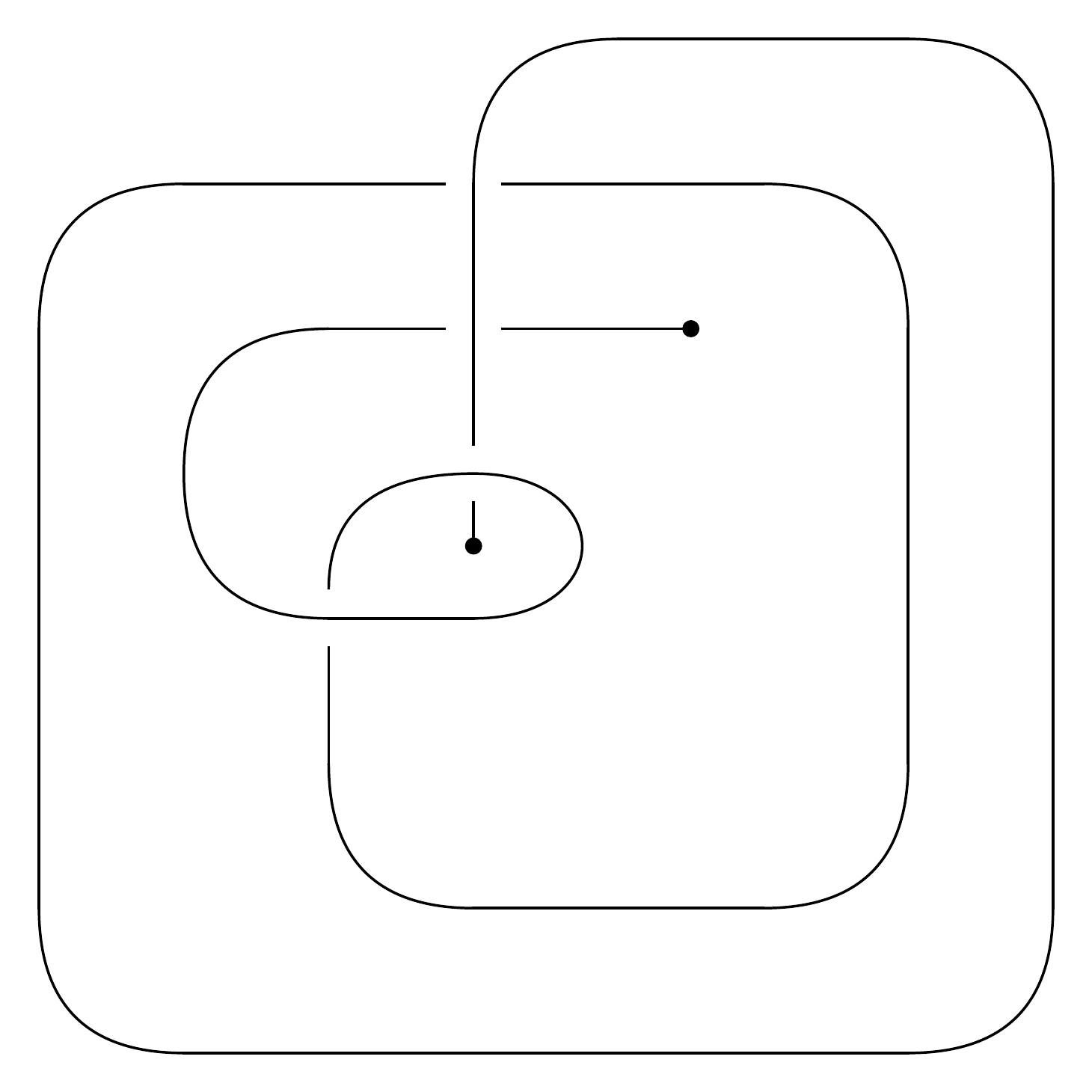}\\
\textcolor{black}{$4_{147}$}
\vspace{1cm}
\end{minipage}
\begin{minipage}[t]{.25\linewidth}
\centering
\includegraphics[width=0.9\textwidth,height=3.5cm,keepaspectratio]{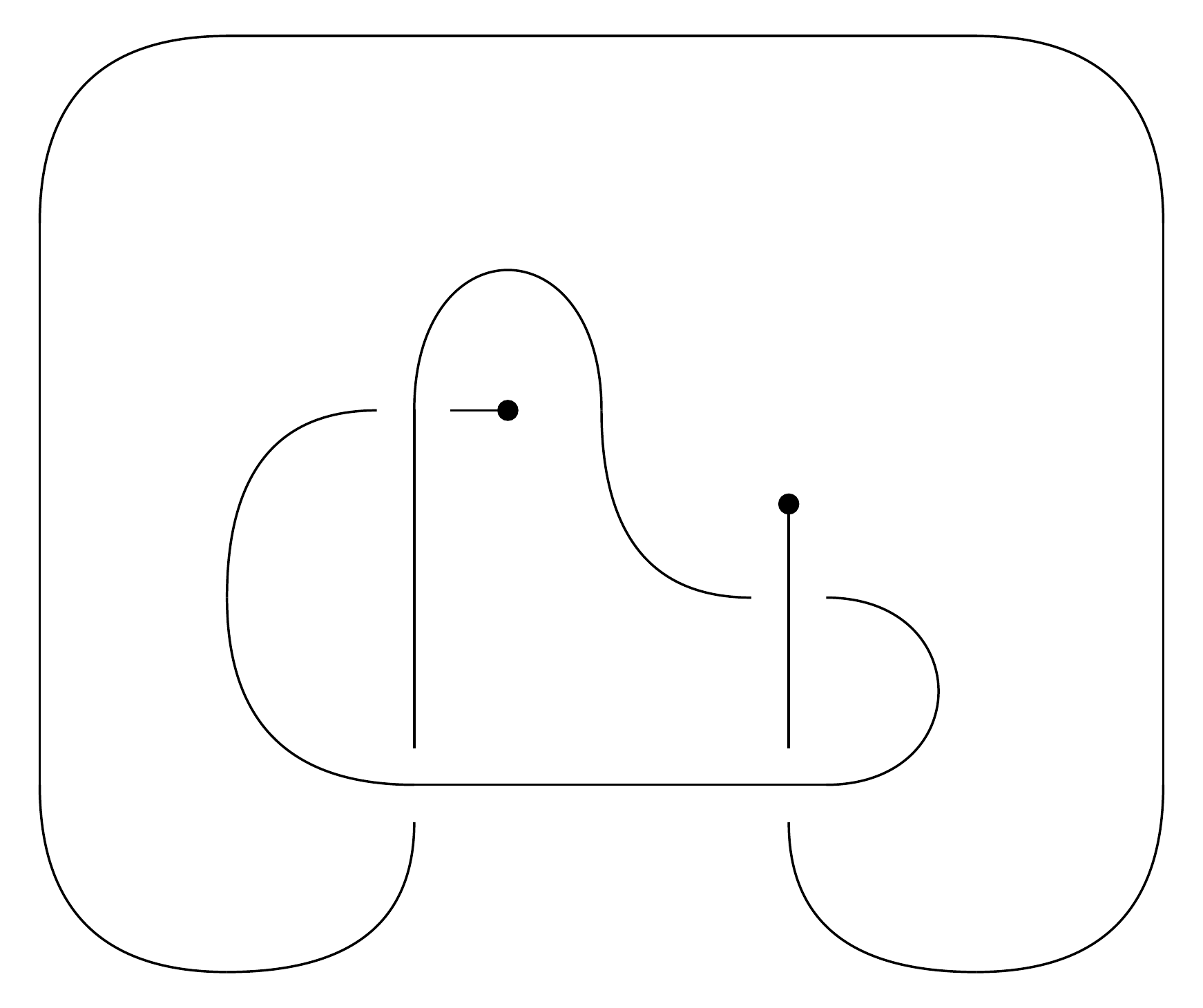}\\
\textcolor{black}{$4_{148}$}
\vspace{1cm}
\end{minipage}
\begin{minipage}[t]{.25\linewidth}
\centering
\includegraphics[width=0.9\textwidth,height=3.5cm,keepaspectratio]{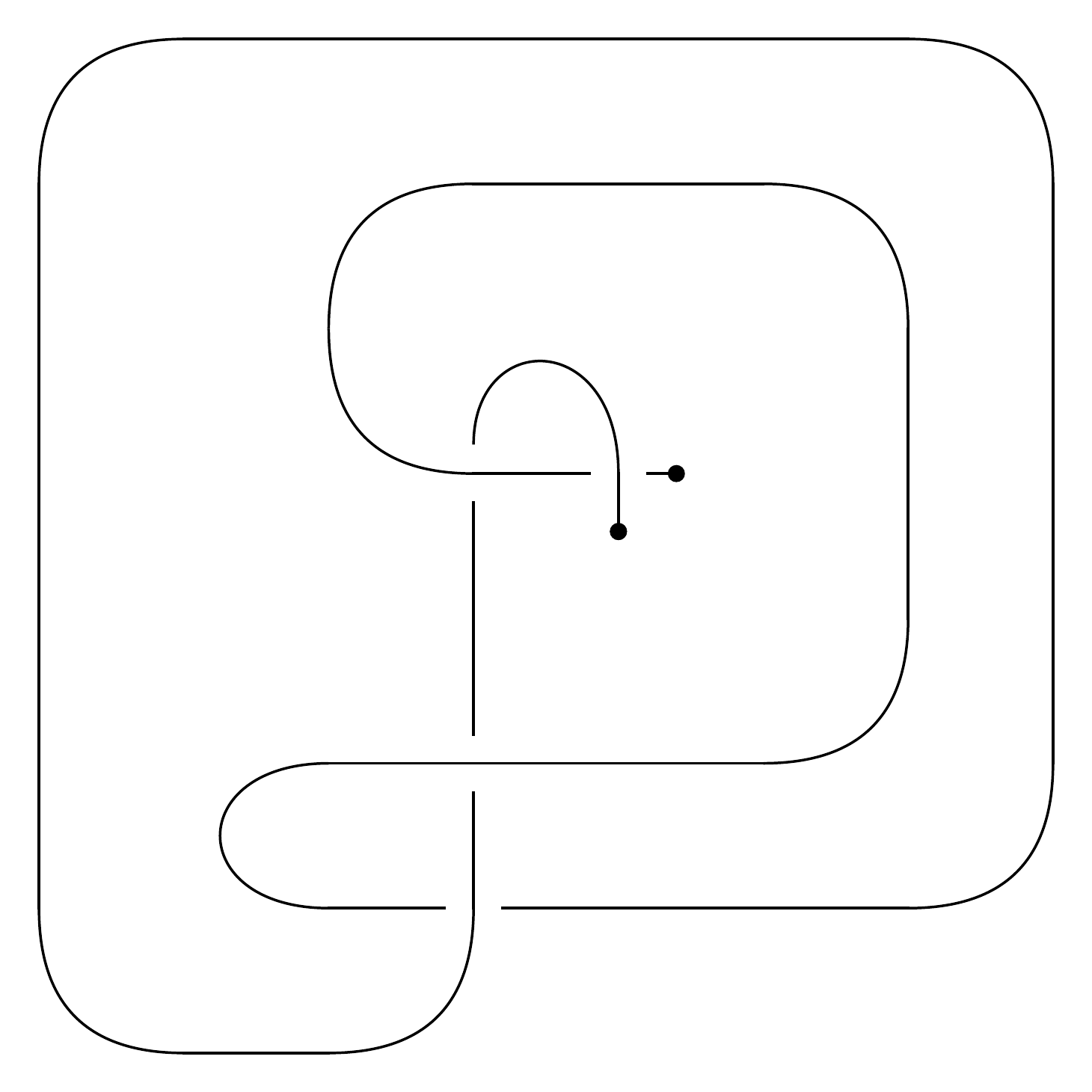}\\
\textcolor{black}{$4_{149}$}
\vspace{1cm}
\end{minipage}
\begin{minipage}[t]{.25\linewidth}
\centering
\includegraphics[width=0.9\textwidth,height=3.5cm,keepaspectratio]{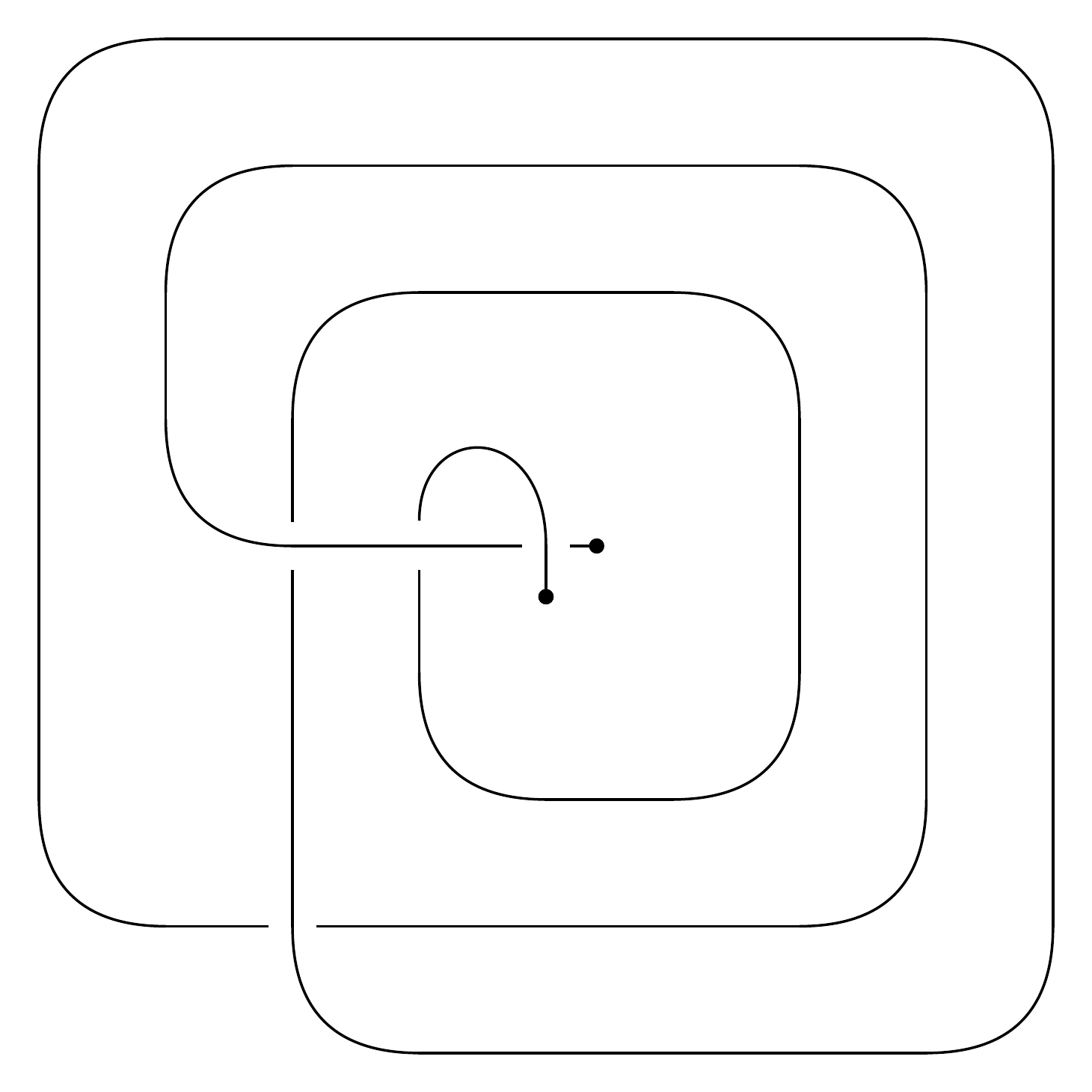}\\
\textcolor{black}{$4_{150}$}
\vspace{1cm}
\end{minipage}
\begin{minipage}[t]{.25\linewidth}
\centering
\includegraphics[width=0.9\textwidth,height=3.5cm,keepaspectratio]{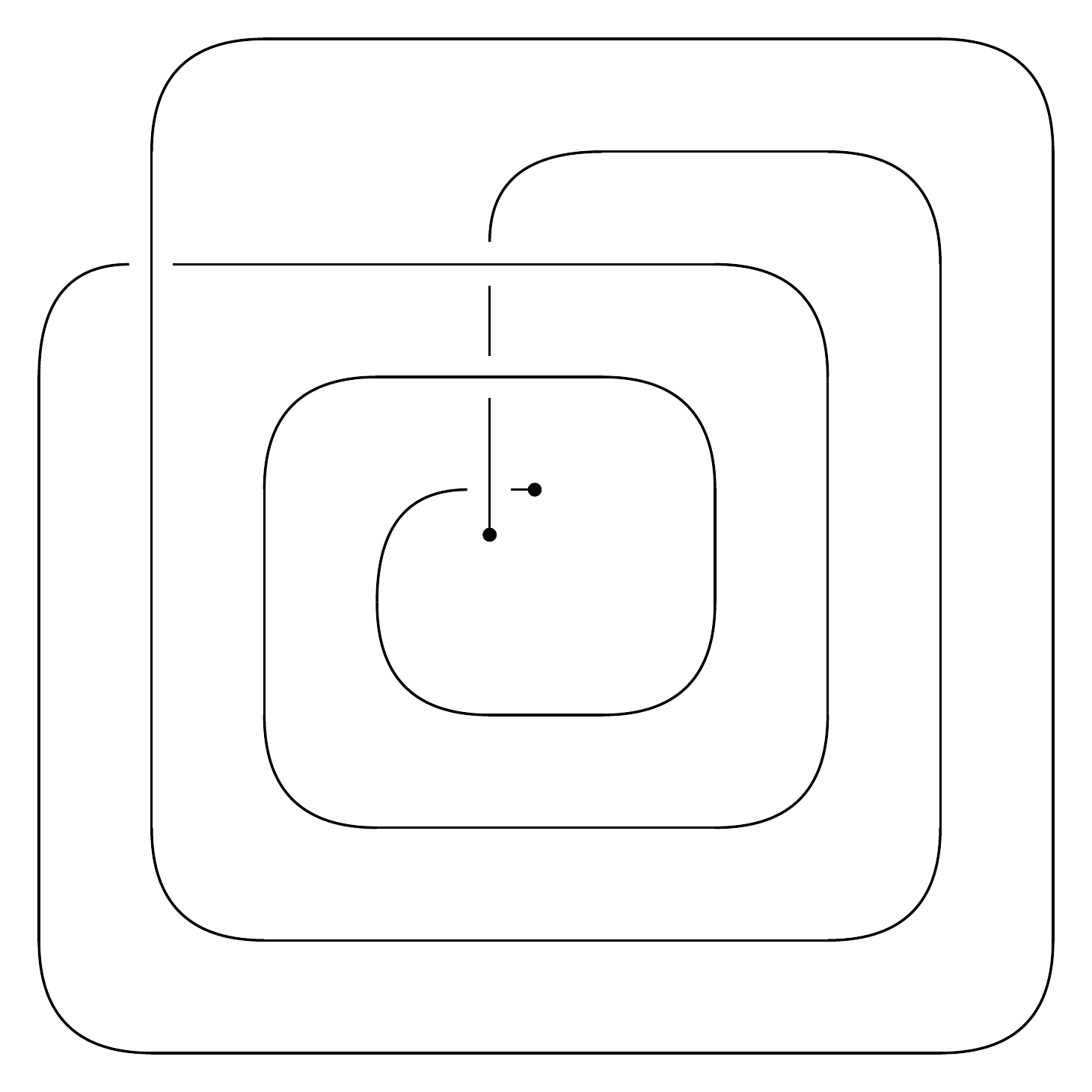}\\
\textcolor{black}{$4_{151}$}
\vspace{1cm}
\end{minipage}
\begin{minipage}[t]{.25\linewidth}
\centering
\includegraphics[width=0.9\textwidth,height=3.5cm,keepaspectratio]{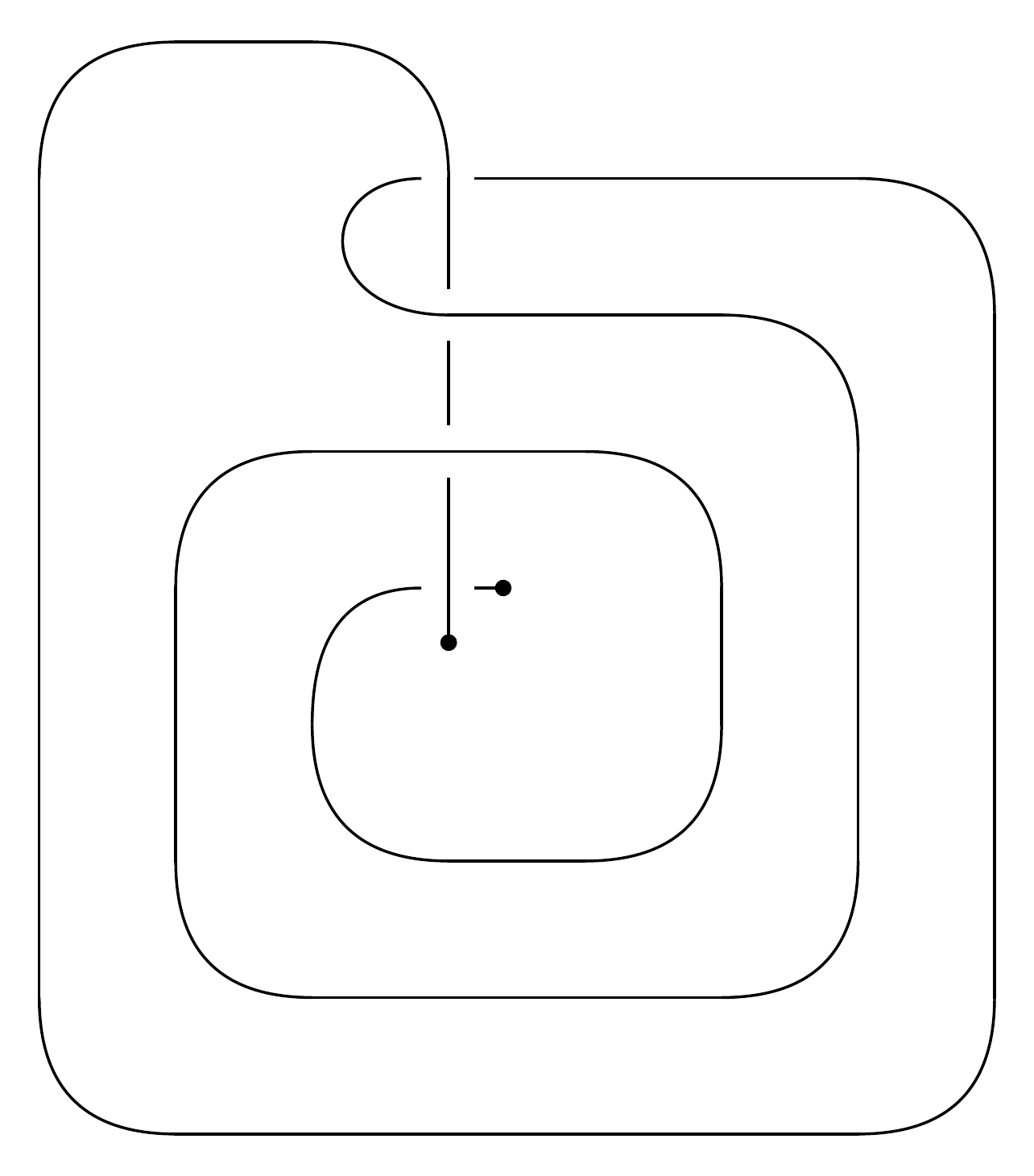}\\
\textcolor{black}{$4_{152}$}
\vspace{1cm}
\end{minipage}
\begin{minipage}[t]{.25\linewidth}
\centering
\includegraphics[width=0.9\textwidth,height=3.5cm,keepaspectratio]{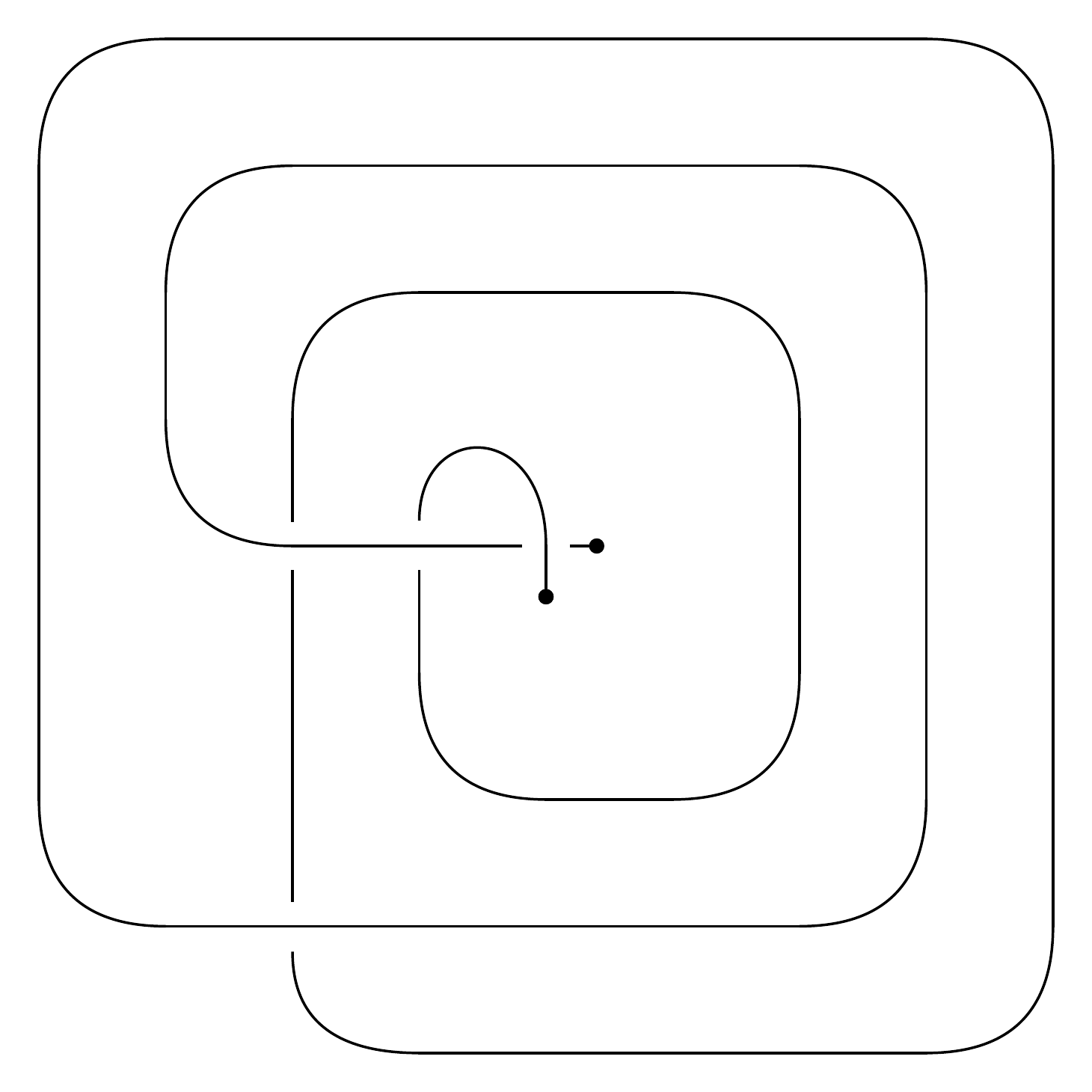}\\
\textcolor{black}{$4_{153}$}
\vspace{1cm}
\end{minipage}
\begin{minipage}[t]{.25\linewidth}
\centering
\includegraphics[width=0.9\textwidth,height=3.5cm,keepaspectratio]{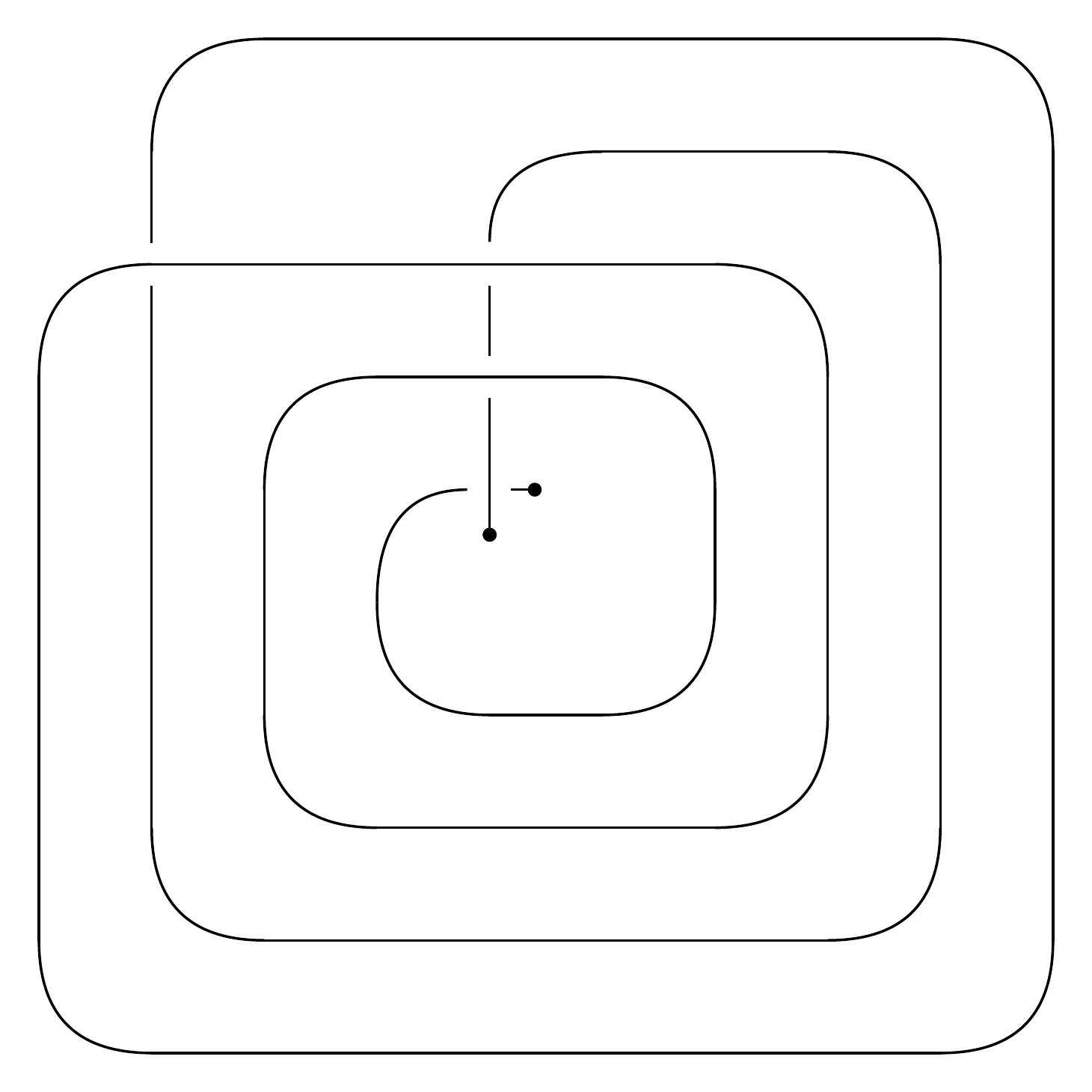}\\
\textcolor{black}{$4_{154}$}
\vspace{1cm}
\end{minipage}
\begin{minipage}[t]{.25\linewidth}
\centering
\includegraphics[width=0.9\textwidth,height=3.5cm,keepaspectratio]{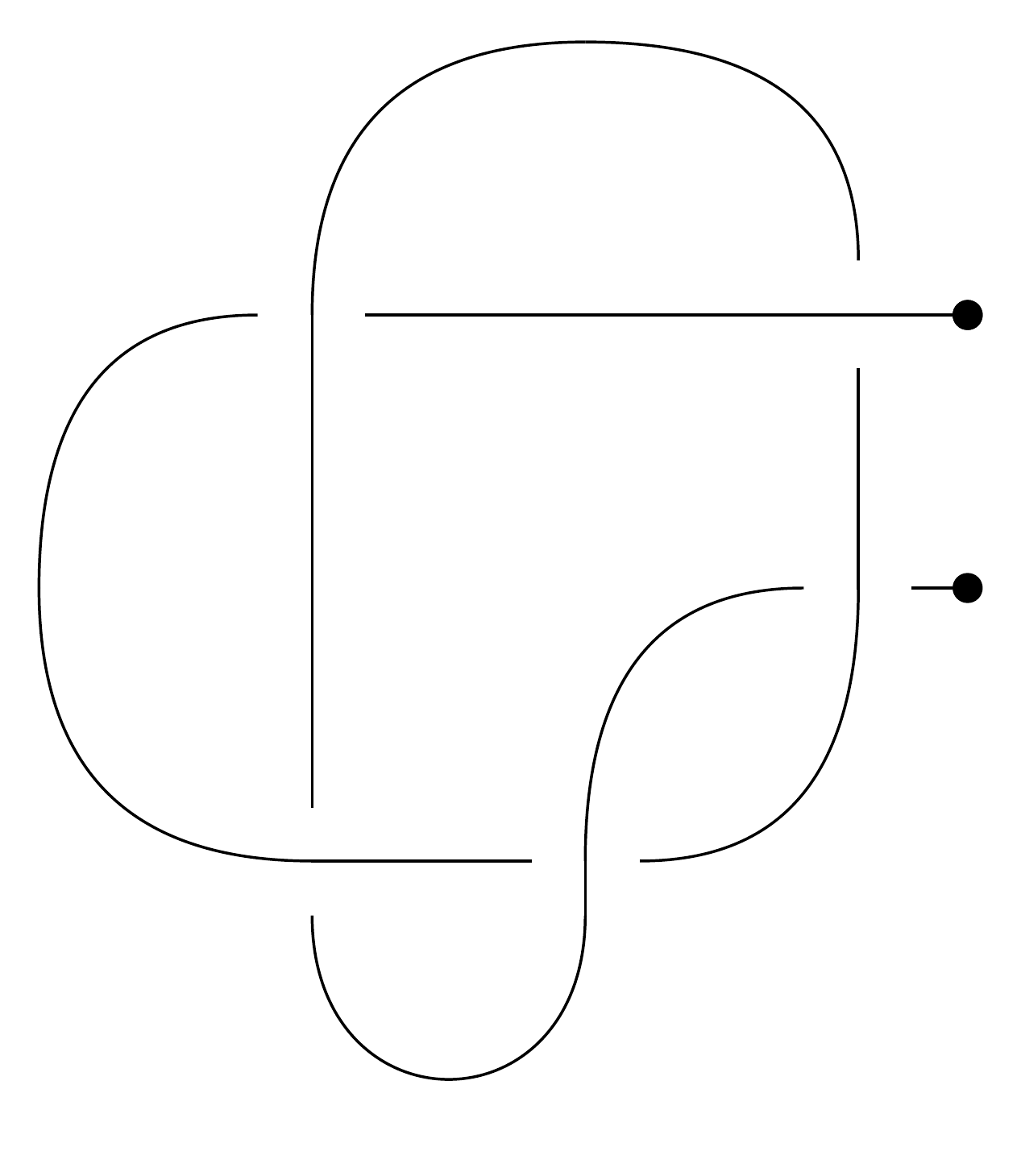}\\
\textcolor{red}{$5_{1}$}
\vspace{1cm}
\end{minipage}
\begin{minipage}[t]{.25\linewidth}
\centering
\includegraphics[width=0.9\textwidth,height=3.5cm,keepaspectratio]{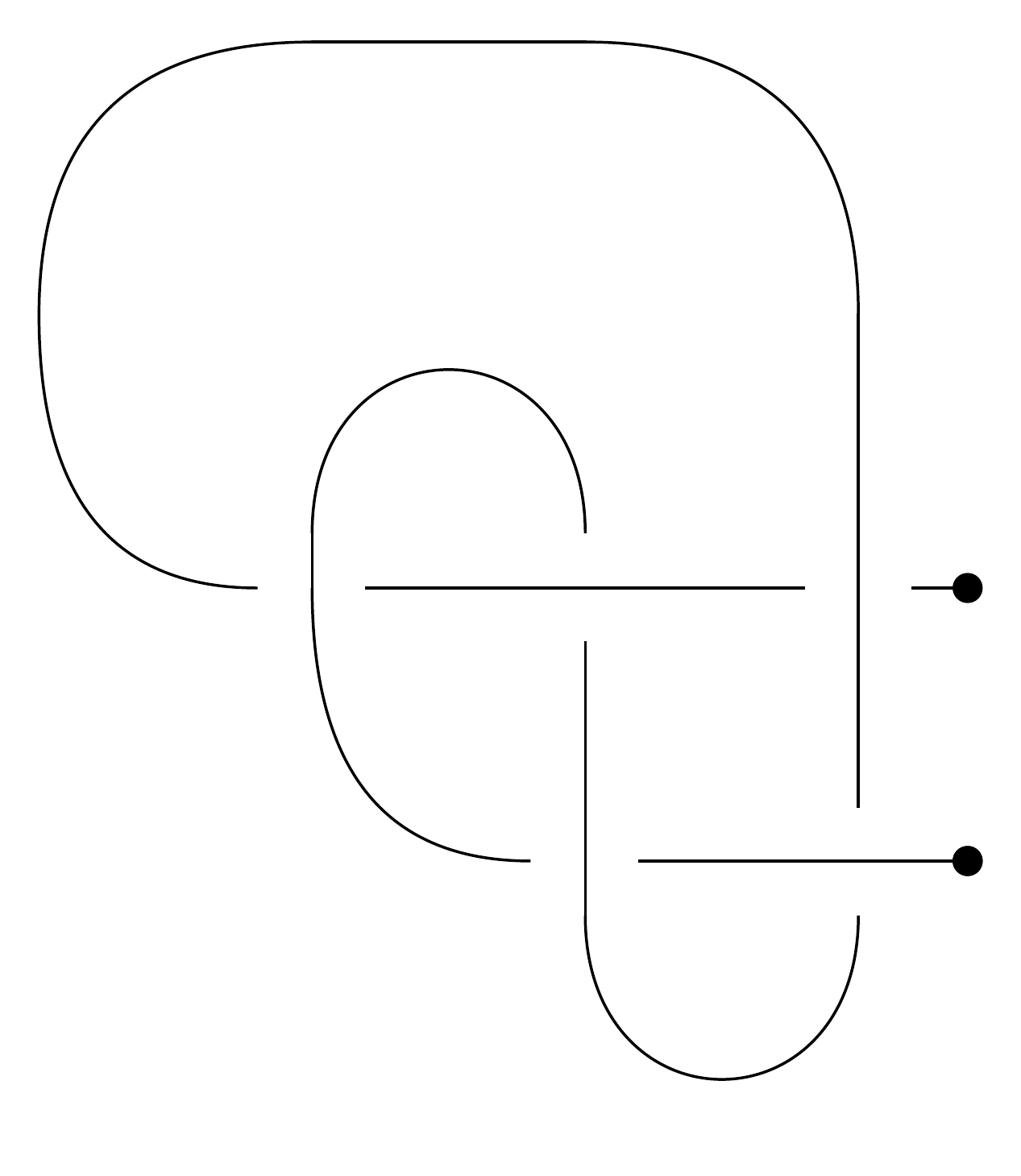}\\
\textcolor{red}{$5_{2}$}
\vspace{1cm}
\end{minipage}
\begin{minipage}[t]{.25\linewidth}
\centering
\includegraphics[width=0.9\textwidth,height=3.5cm,keepaspectratio]{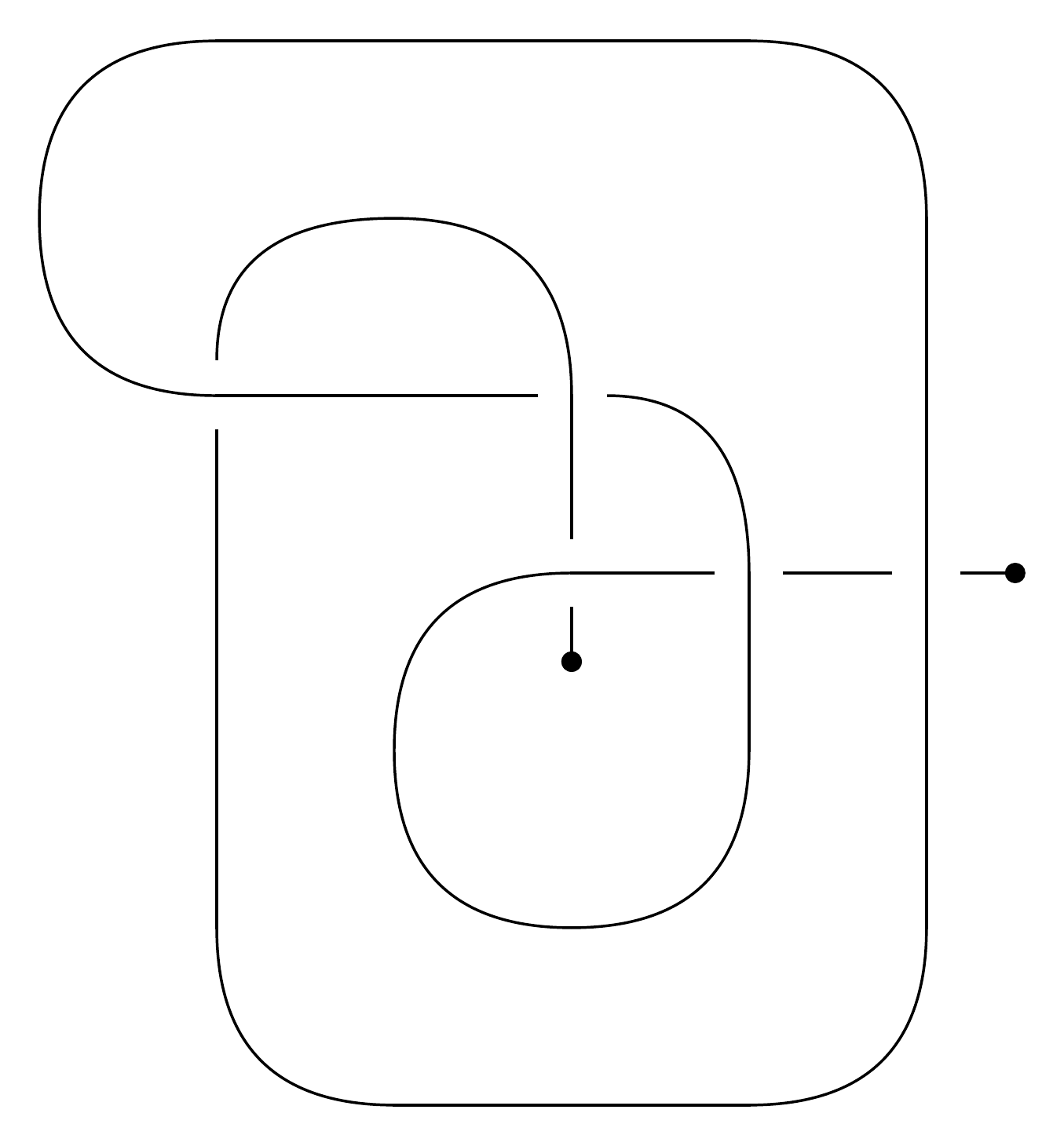}\\
\textcolor{blue}{$5_{3}$}
\vspace{1cm}
\end{minipage}
\begin{minipage}[t]{.25\linewidth}
\centering
\includegraphics[width=0.9\textwidth,height=3.5cm,keepaspectratio]{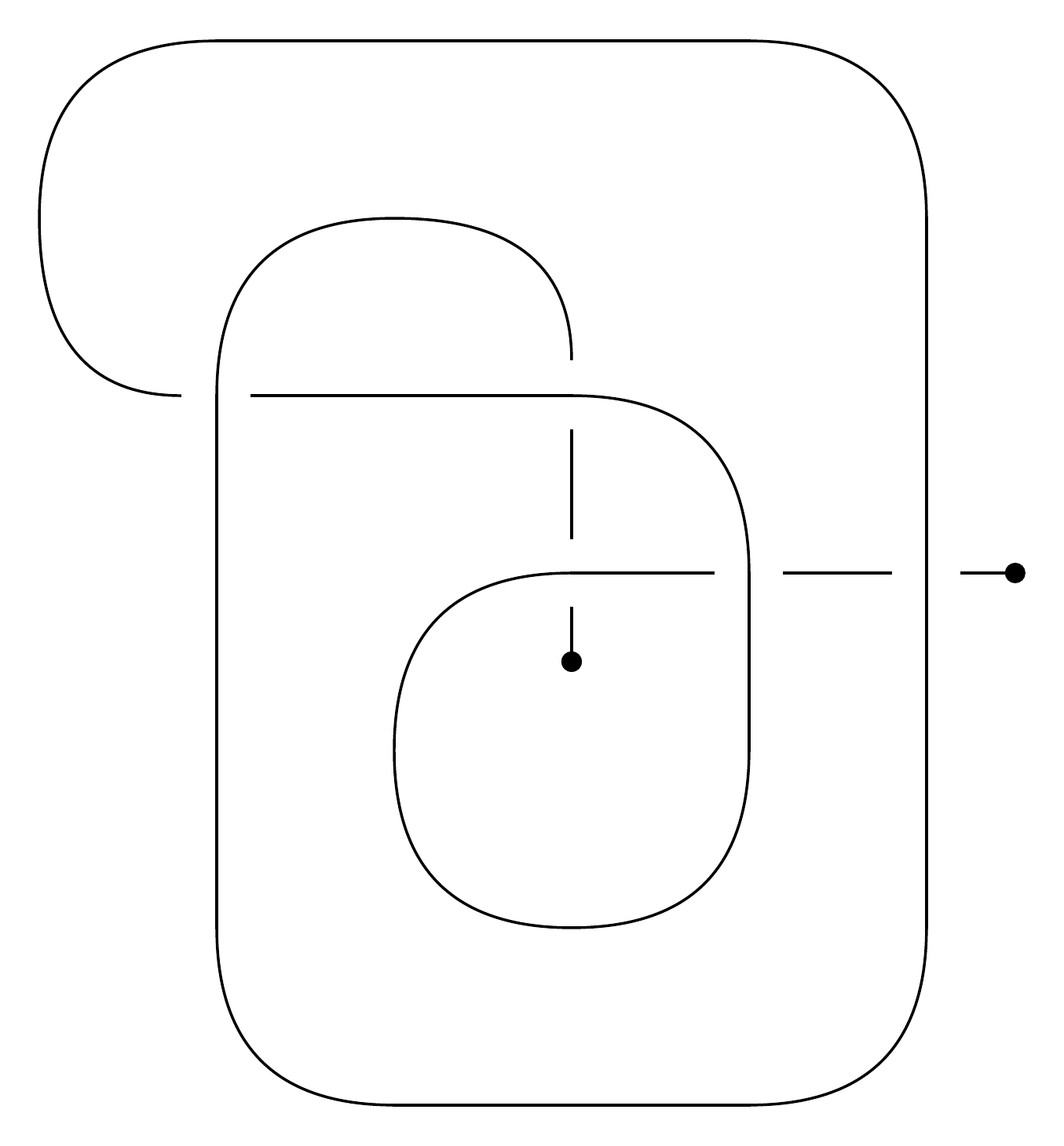}\\
\textcolor{blue}{$5_{4}$}
\vspace{1cm}
\end{minipage}
\begin{minipage}[t]{.25\linewidth}
\centering
\includegraphics[width=0.9\textwidth,height=3.5cm,keepaspectratio]{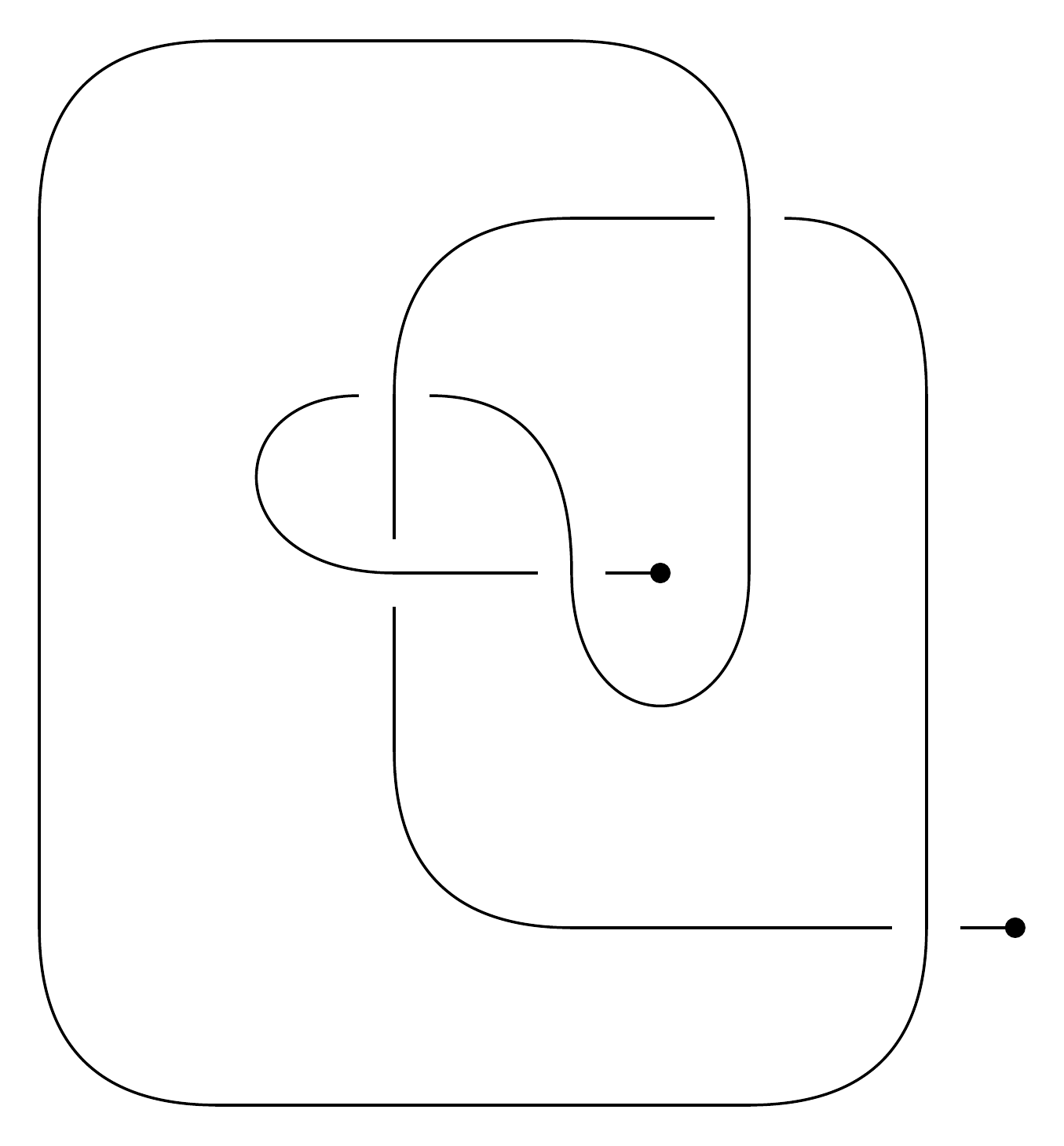}\\
\textcolor{blue}{$5_{5}$}
\vspace{1cm}
\end{minipage}
\begin{minipage}[t]{.25\linewidth}
\centering
\includegraphics[width=0.9\textwidth,height=3.5cm,keepaspectratio]{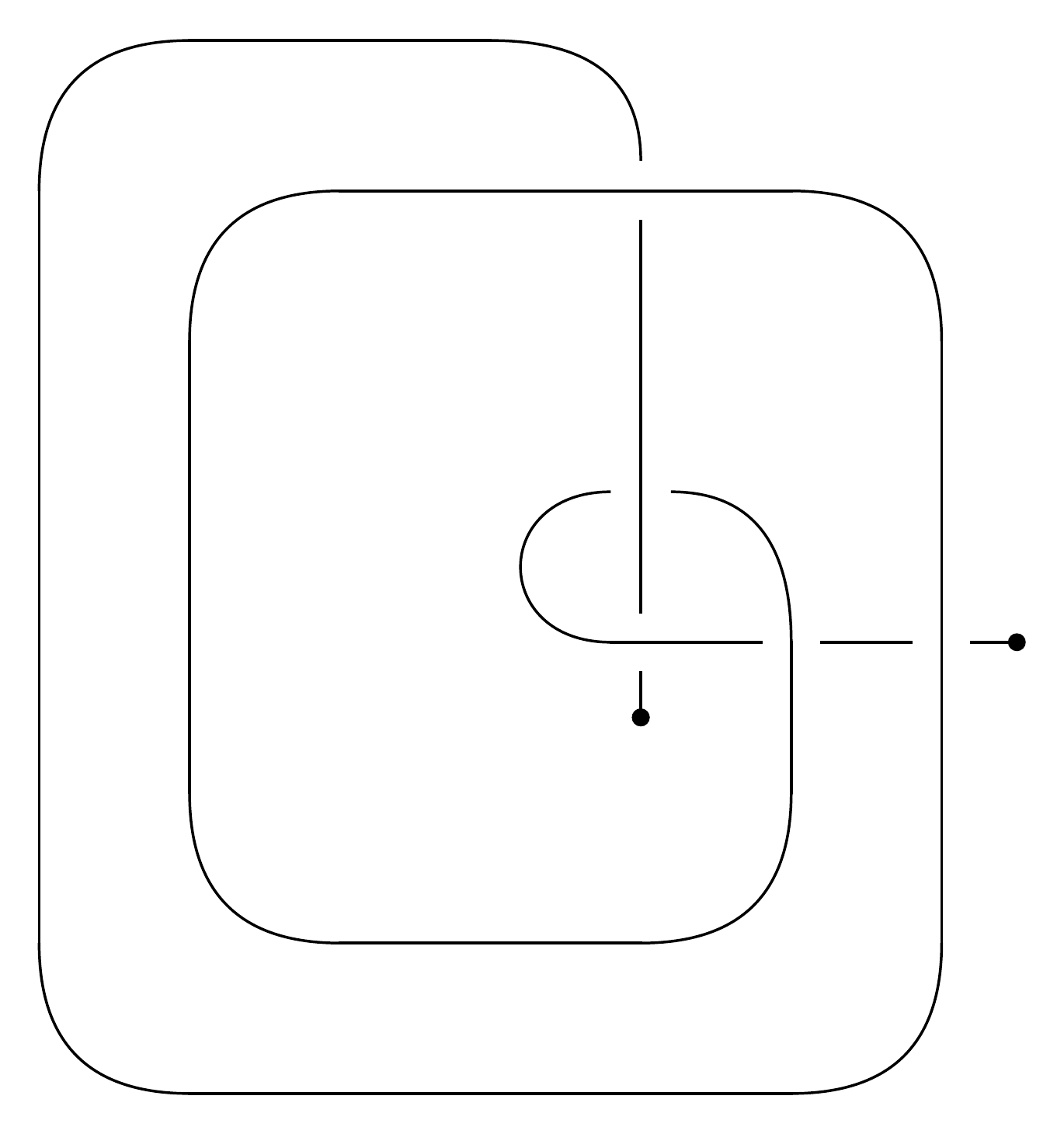}\\
\textcolor{blue}{$5_{6}$}
\vspace{1cm}
\end{minipage}
\begin{minipage}[t]{.25\linewidth}
\centering
\includegraphics[width=0.9\textwidth,height=3.5cm,keepaspectratio]{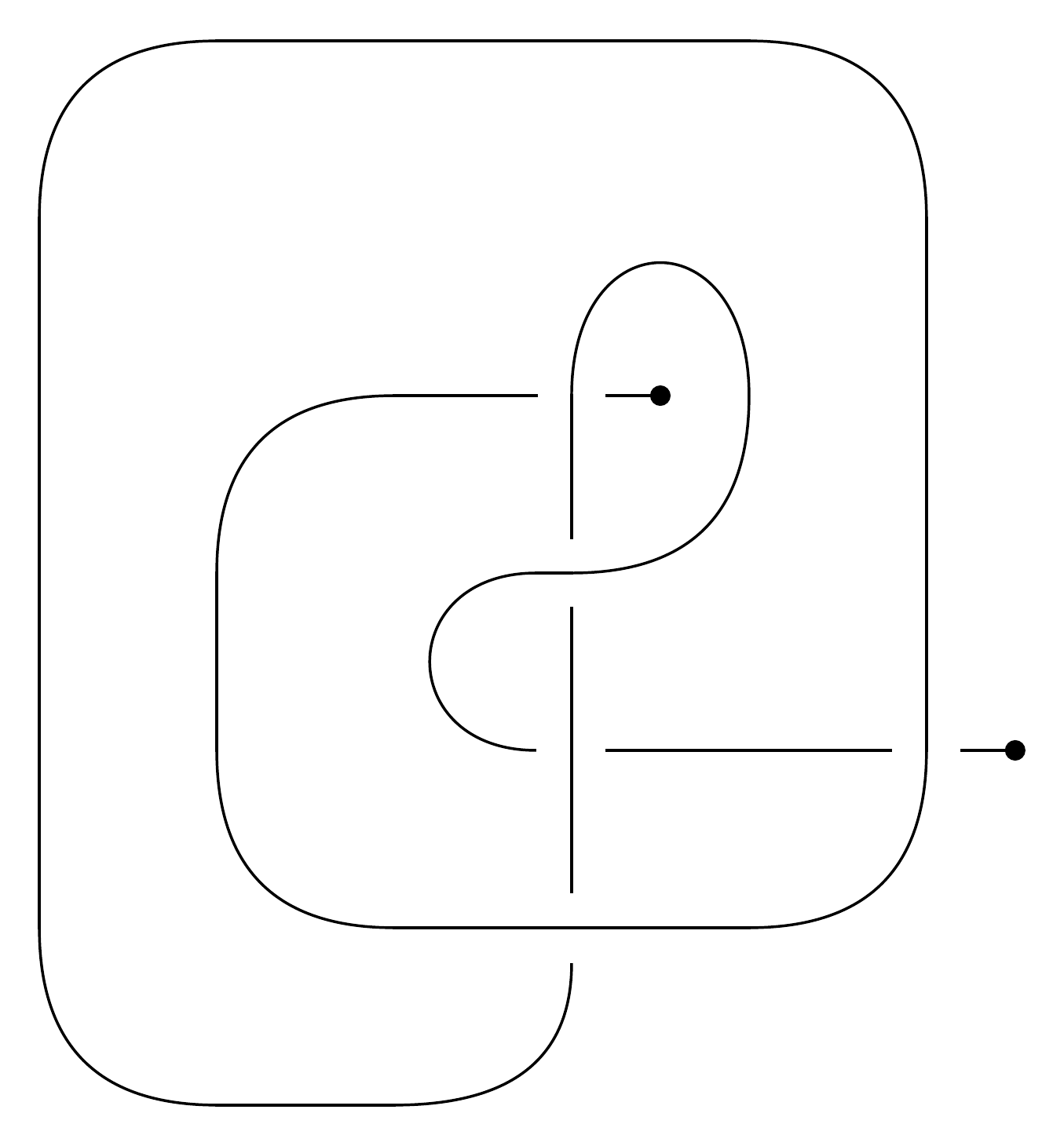}\\
\textcolor{blue}{$5_{7}$}
\vspace{1cm}
\end{minipage}
\begin{minipage}[t]{.25\linewidth}
\centering
\includegraphics[width=0.9\textwidth,height=3.5cm,keepaspectratio]{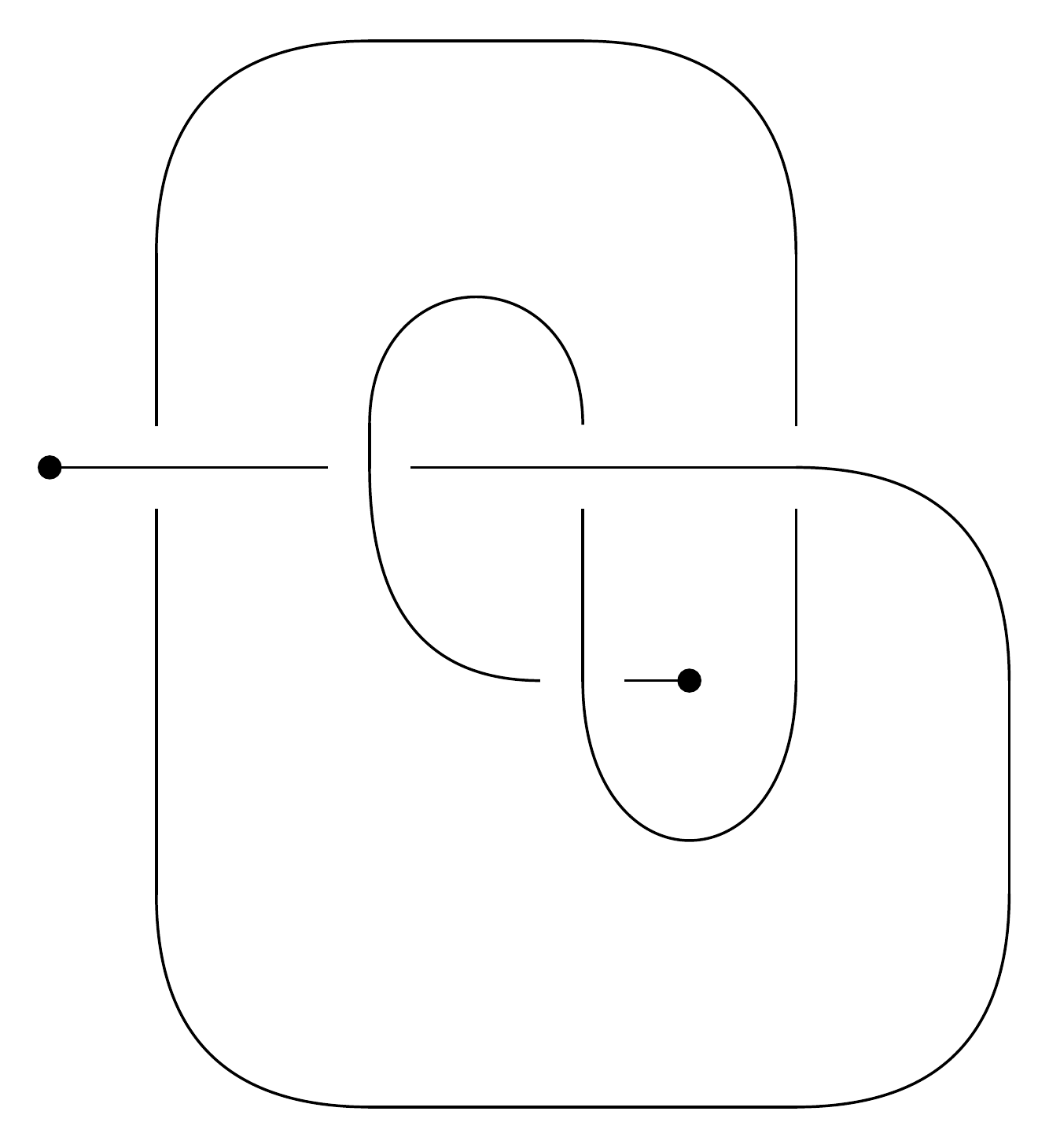}\\
\textcolor{blue}{$5_{8}$}
\vspace{1cm}
\end{minipage}
\begin{minipage}[t]{.25\linewidth}
\centering
\includegraphics[width=0.9\textwidth,height=3.5cm,keepaspectratio]{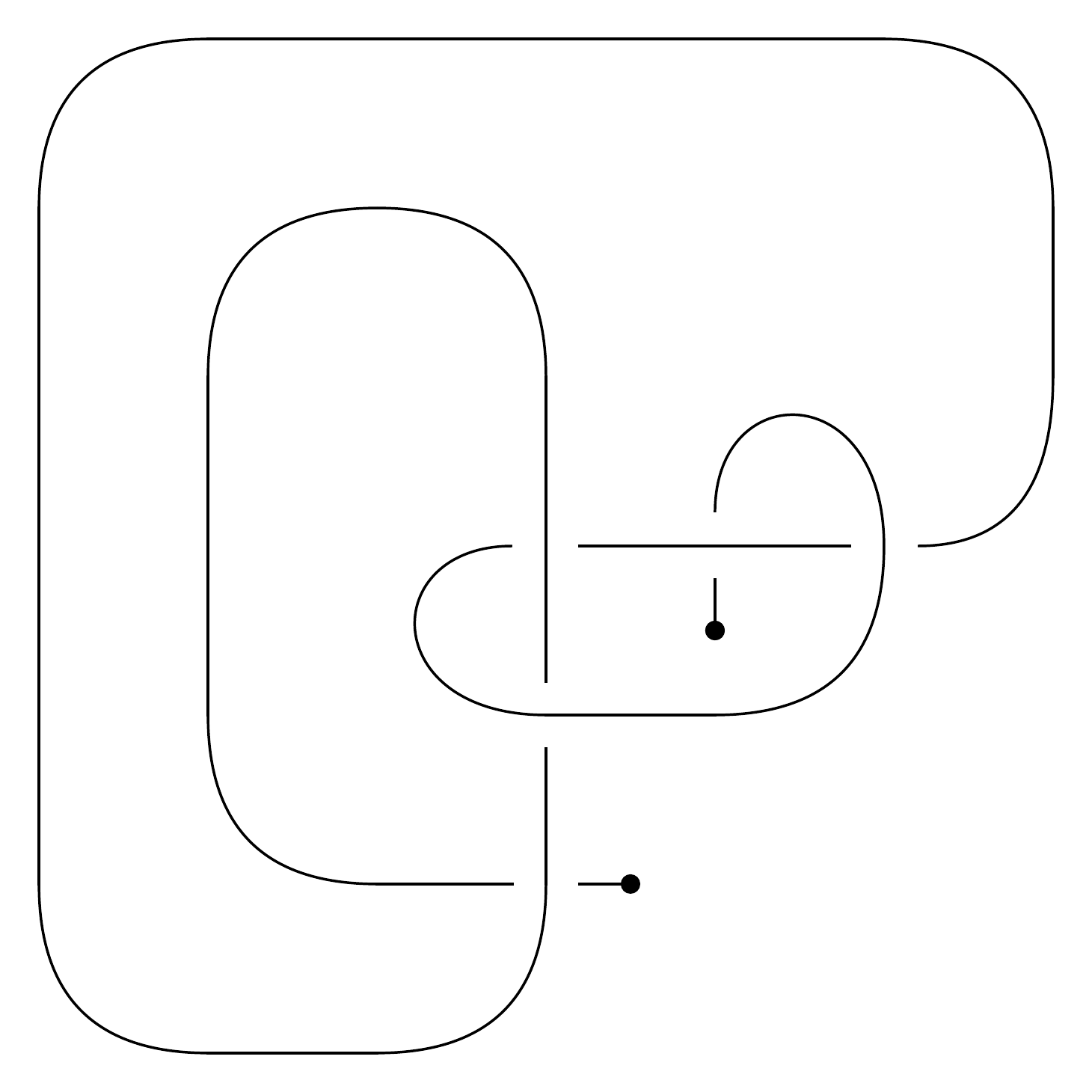}\\
\textcolor{blue}{$5_{9}$}
\vspace{1cm}
\end{minipage}
\begin{minipage}[t]{.25\linewidth}
\centering
\includegraphics[width=0.9\textwidth,height=3.5cm,keepaspectratio]{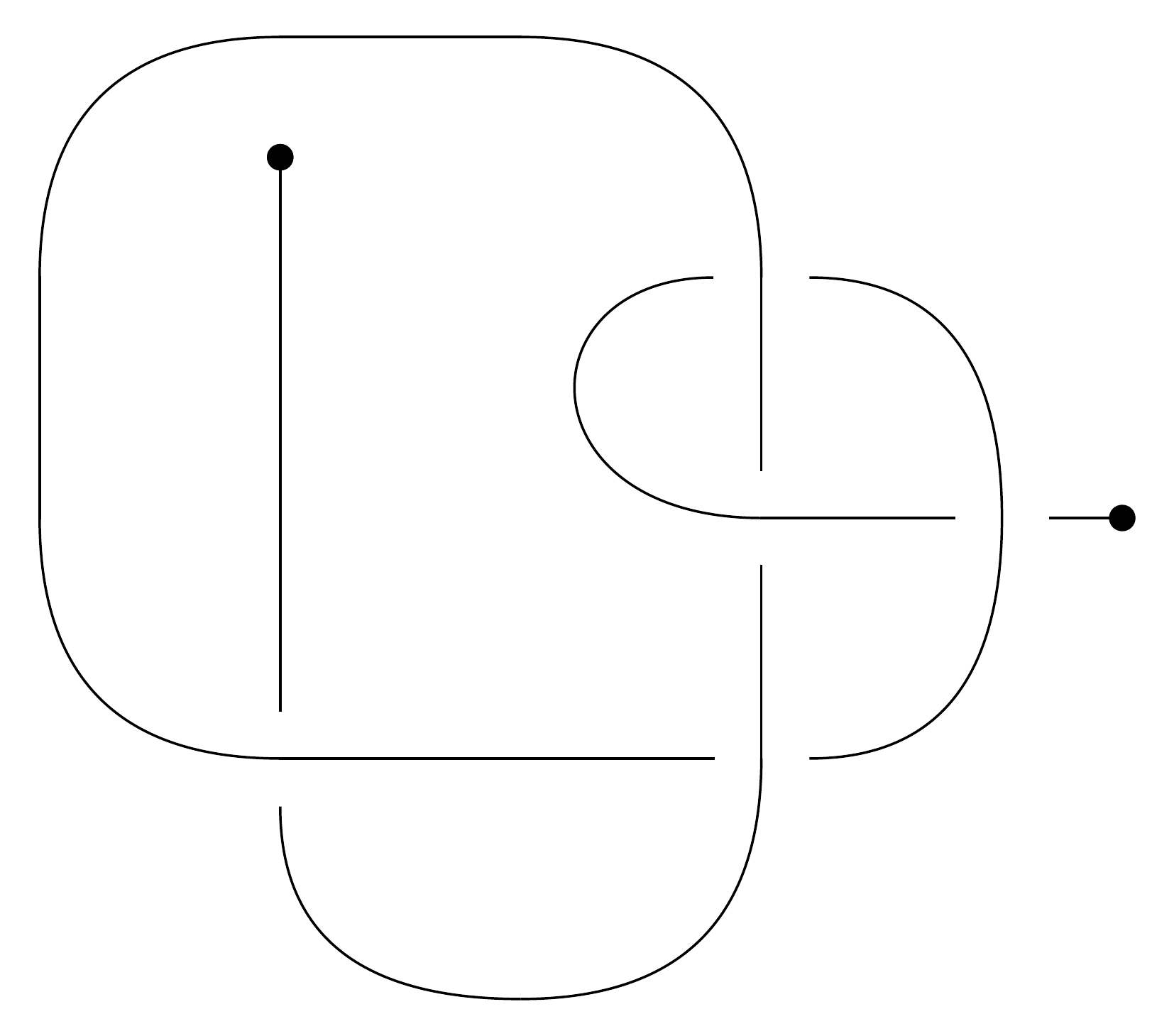}\\
\textcolor{blue}{$5_{10}$}
\vspace{1cm}
\end{minipage}
\begin{minipage}[t]{.25\linewidth}
\centering
\includegraphics[width=0.9\textwidth,height=3.5cm,keepaspectratio]{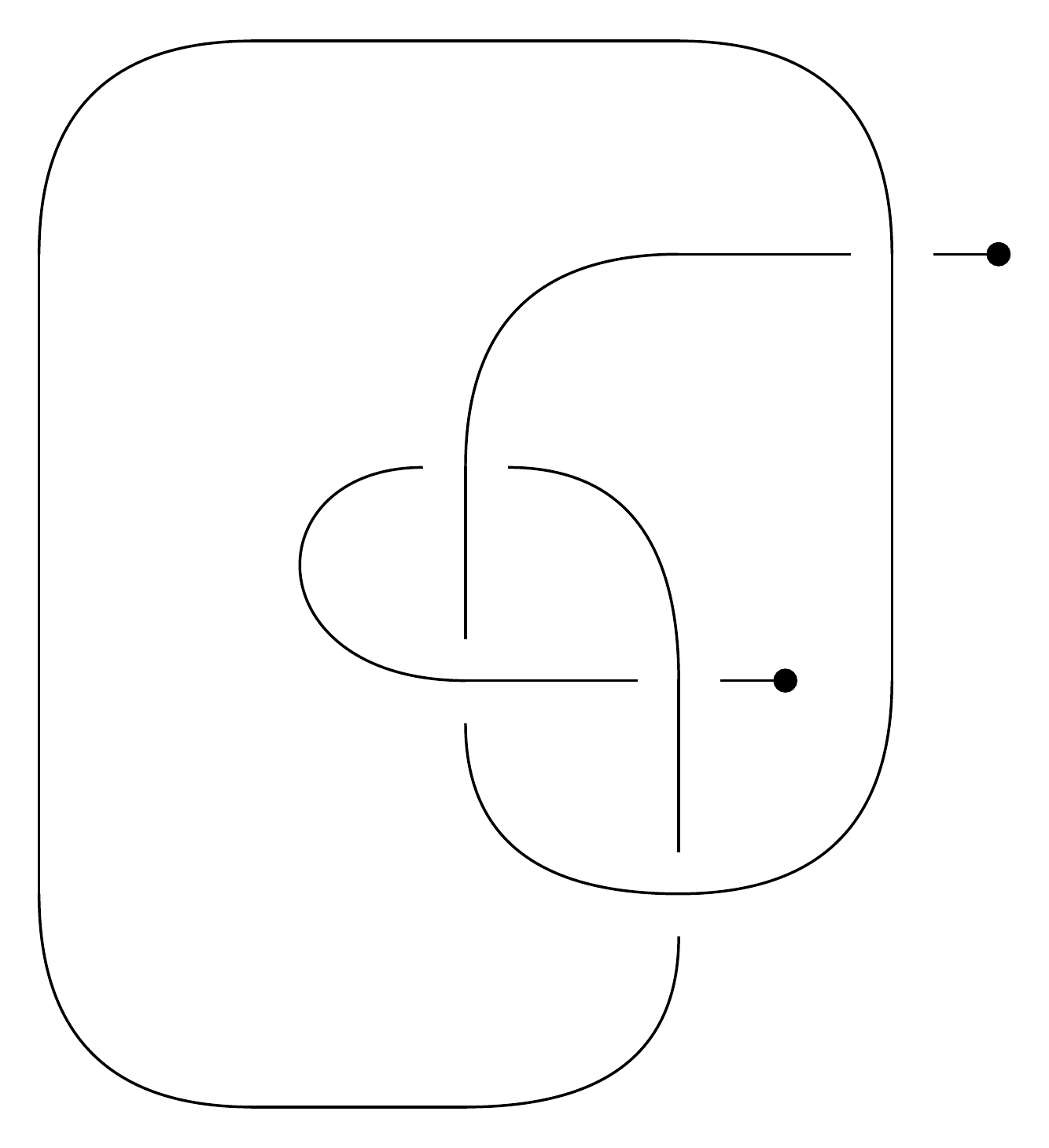}\\
\textcolor{blue}{$5_{11}$}
\vspace{1cm}
\end{minipage}
\begin{minipage}[t]{.25\linewidth}
\centering
\includegraphics[width=0.9\textwidth,height=3.5cm,keepaspectratio]{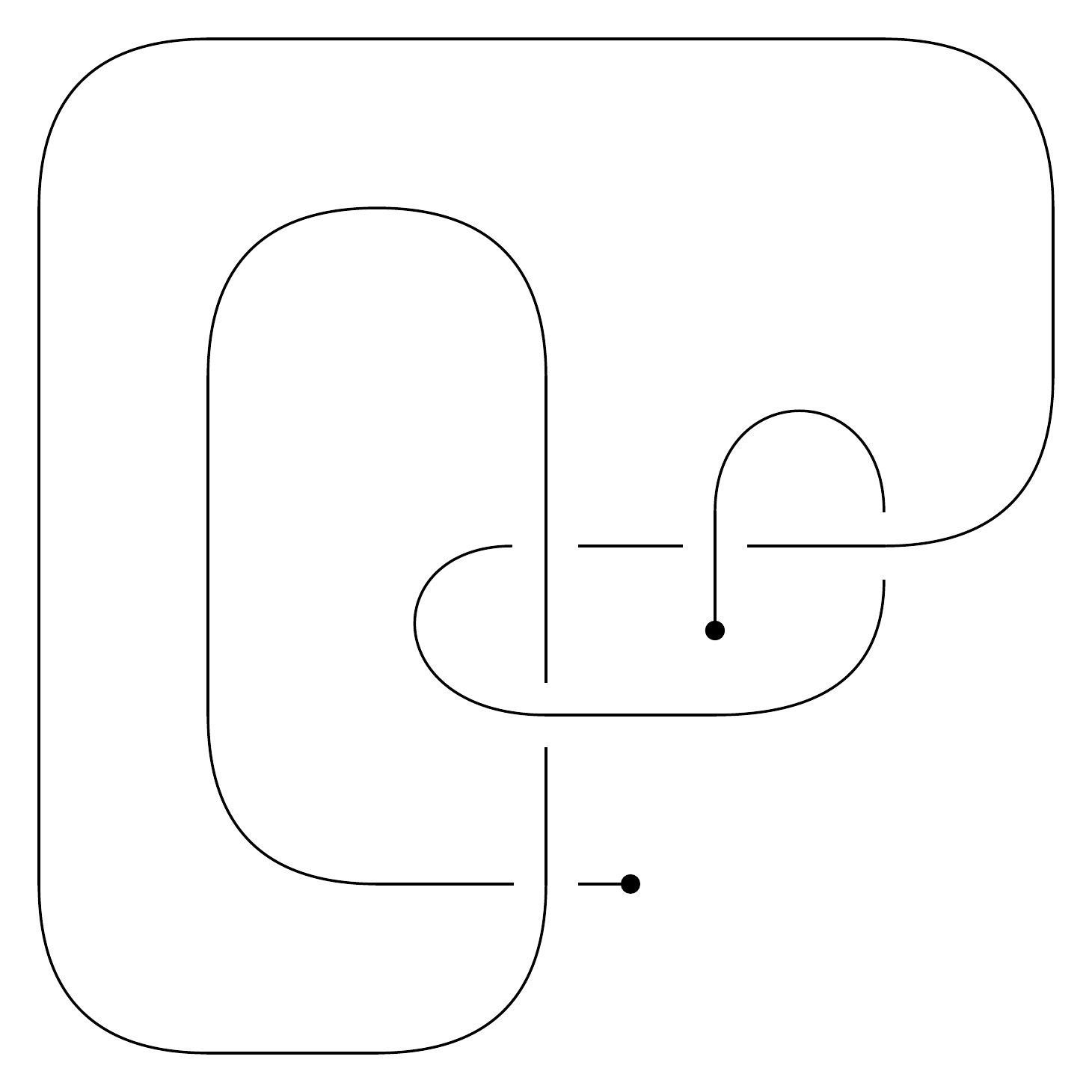}\\
\textcolor{blue}{$5_{12}$}
\vspace{1cm}
\end{minipage}
\begin{minipage}[t]{.25\linewidth}
\centering
\includegraphics[width=0.9\textwidth,height=3.5cm,keepaspectratio]{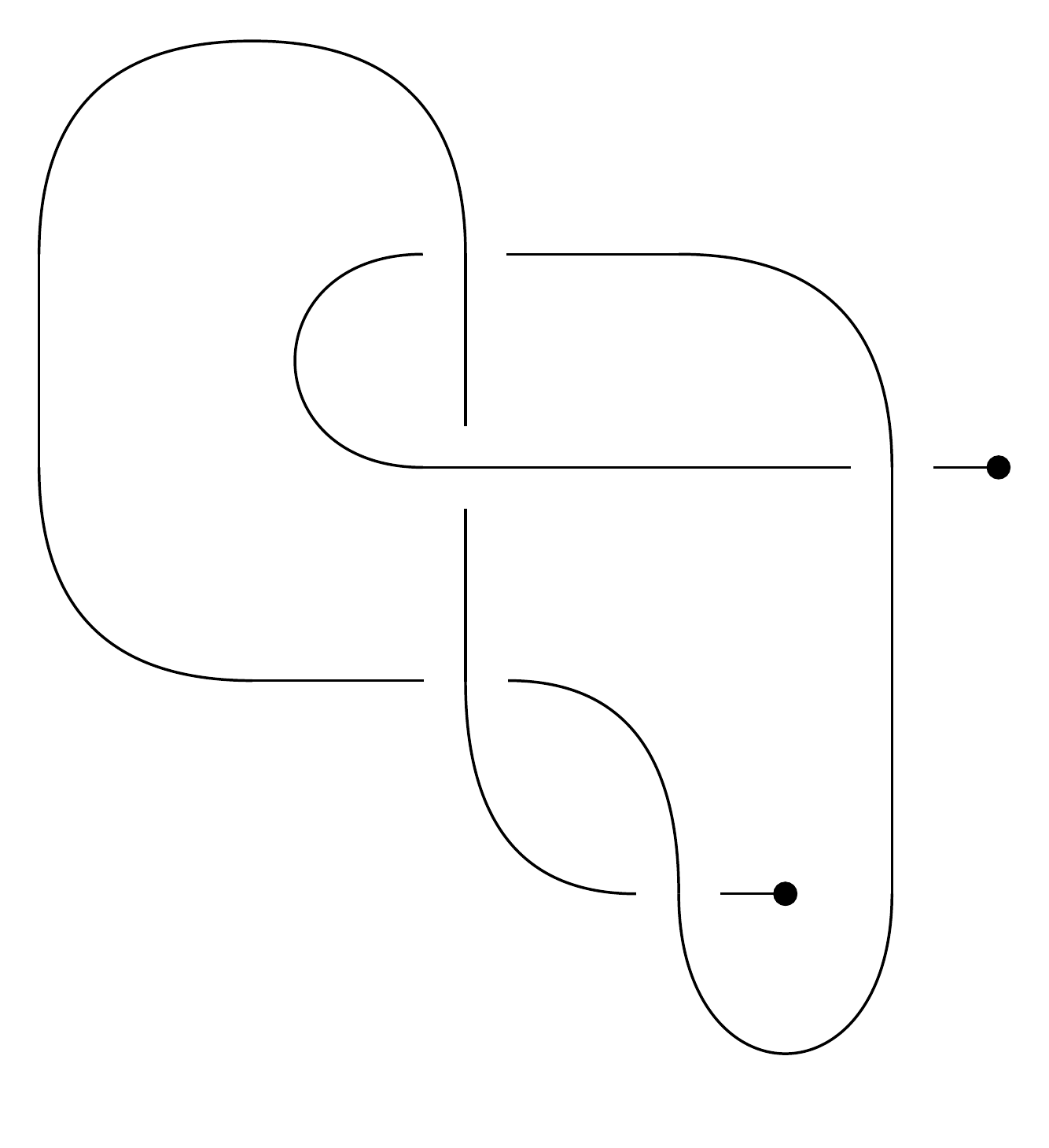}\\
\textcolor{blue}{$5_{13}$}
\vspace{1cm}
\end{minipage}
\begin{minipage}[t]{.25\linewidth}
\centering
\includegraphics[width=0.9\textwidth,height=3.5cm,keepaspectratio]{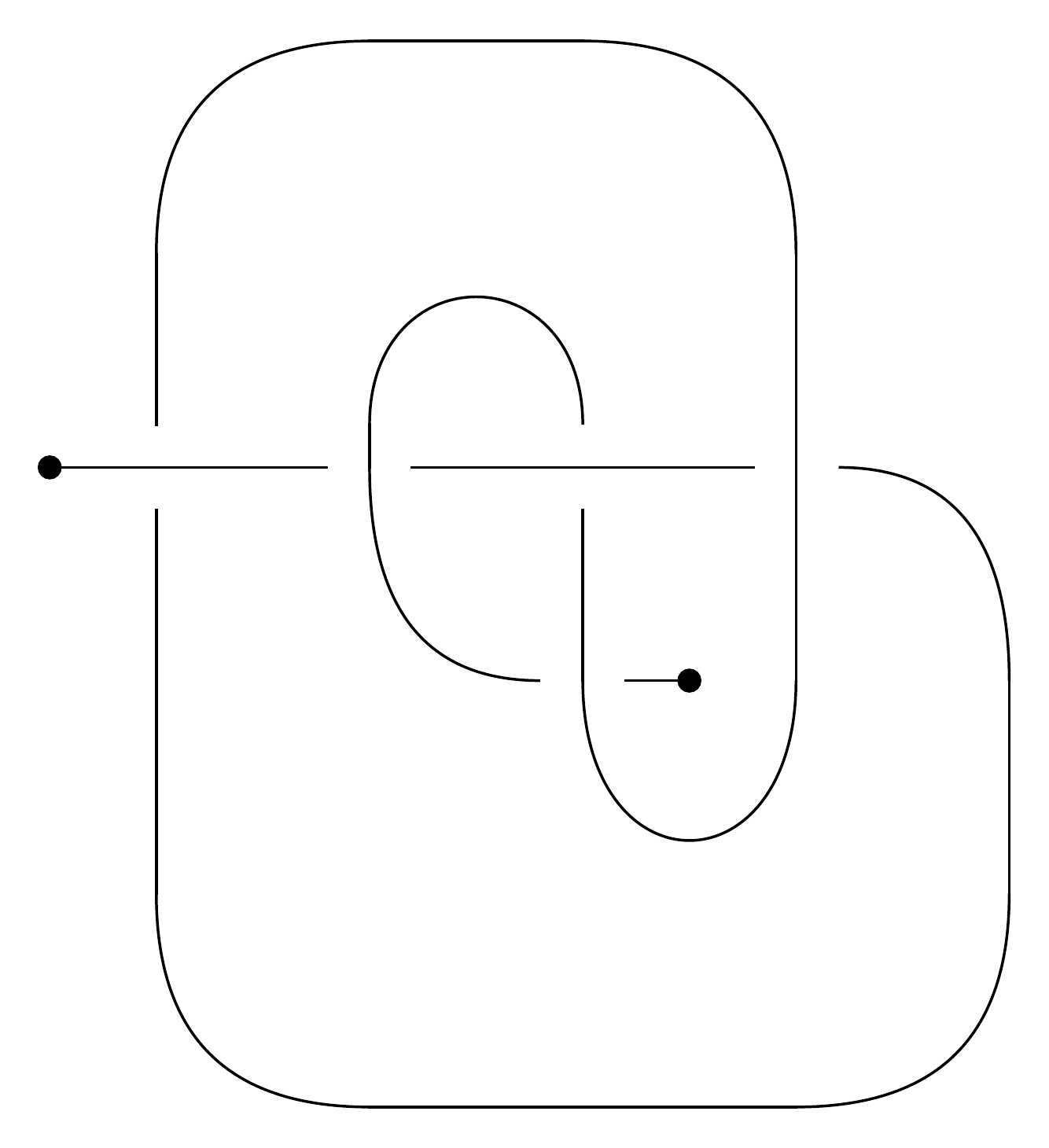}\\
\textcolor{blue}{$5_{14}$}
\vspace{1cm}
\end{minipage}
\begin{minipage}[t]{.25\linewidth}
\centering
\includegraphics[width=0.9\textwidth,height=3.5cm,keepaspectratio]{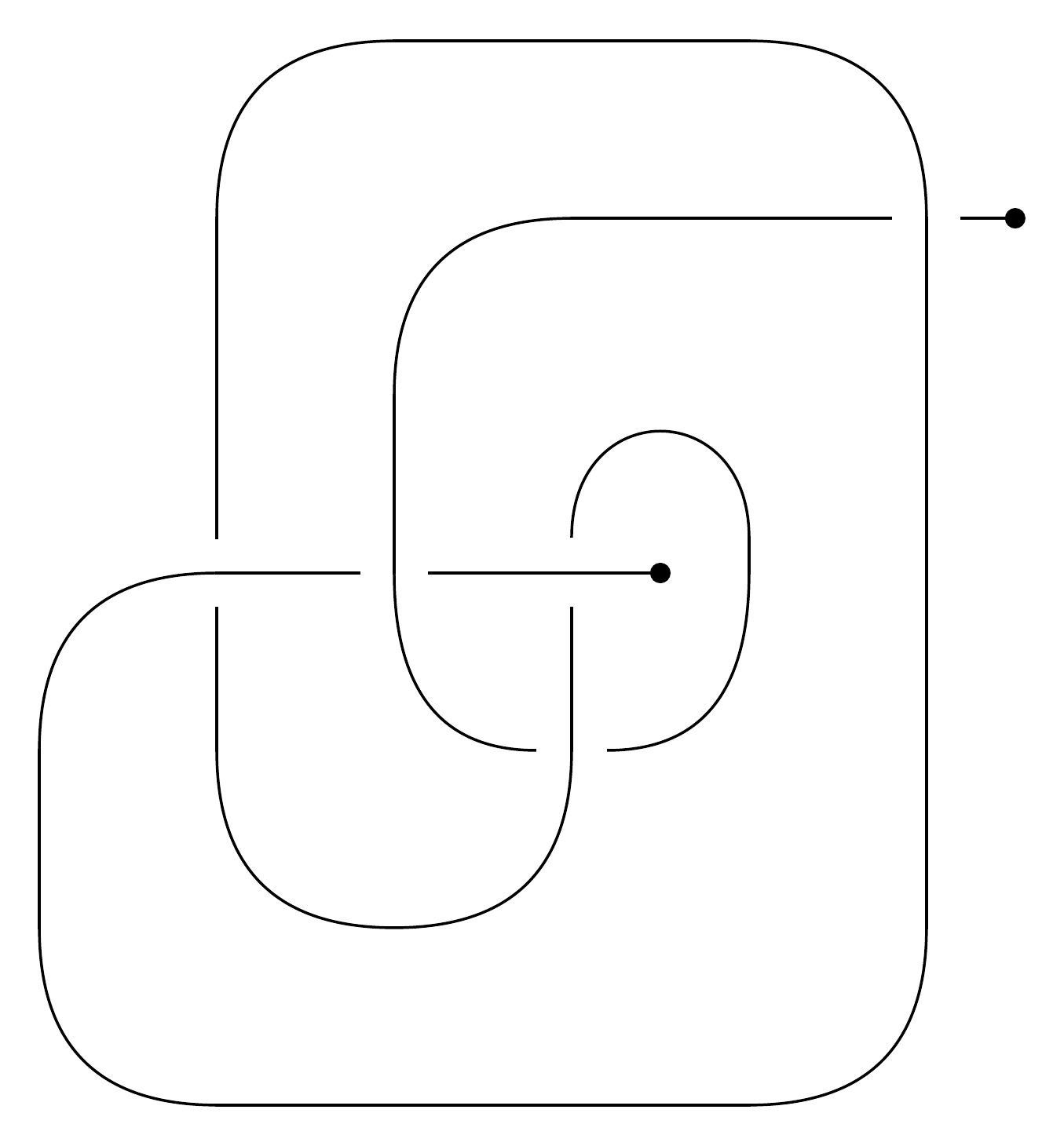}\\
\textcolor{blue}{$5_{15}$}
\vspace{1cm}
\end{minipage}
\begin{minipage}[t]{.25\linewidth}
\centering
\includegraphics[width=0.9\textwidth,height=3.5cm,keepaspectratio]{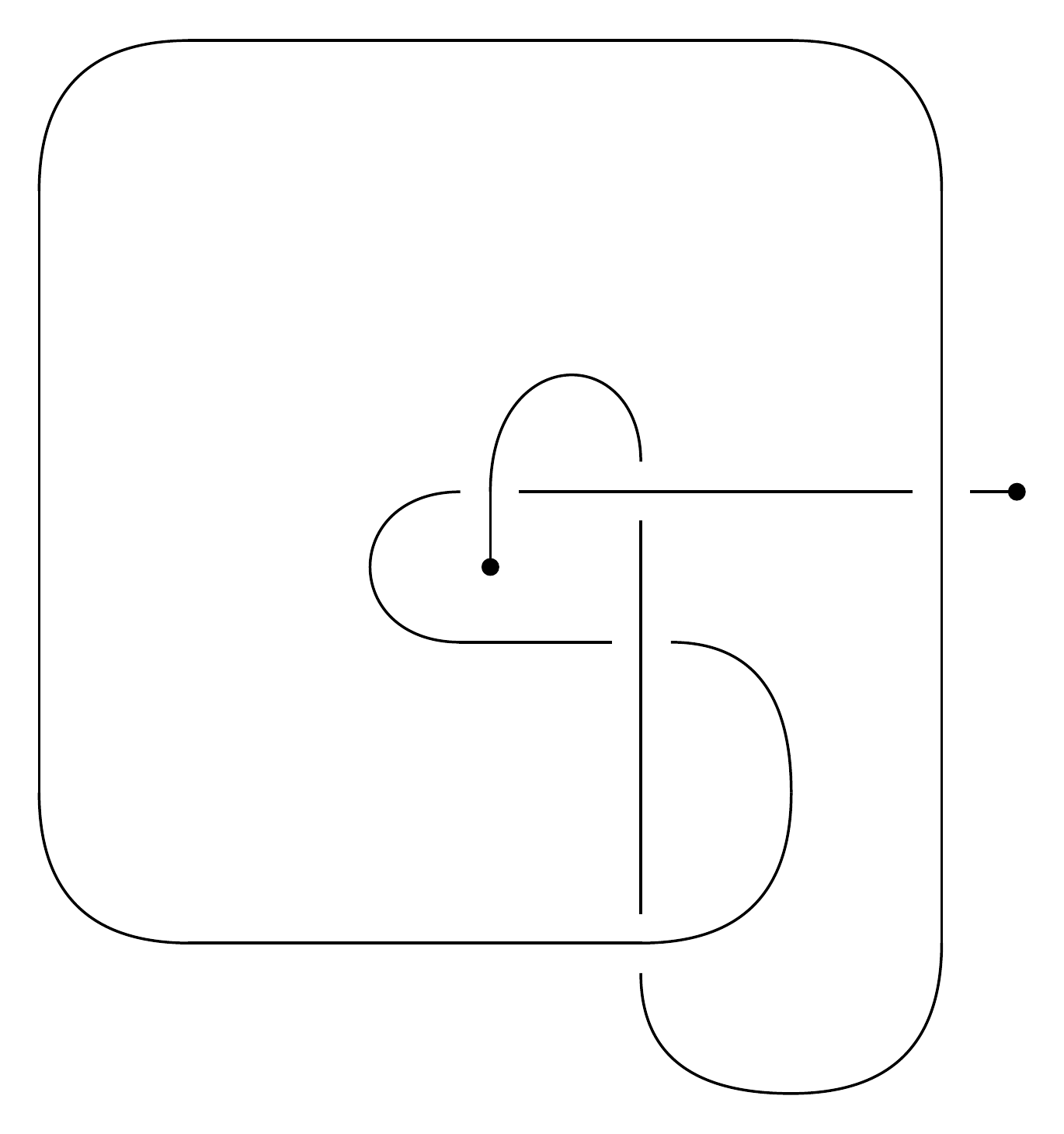}\\
\textcolor{blue}{$5_{16}$}
\vspace{1cm}
\end{minipage}
\begin{minipage}[t]{.25\linewidth}
\centering
\includegraphics[width=0.9\textwidth,height=3.5cm,keepaspectratio]{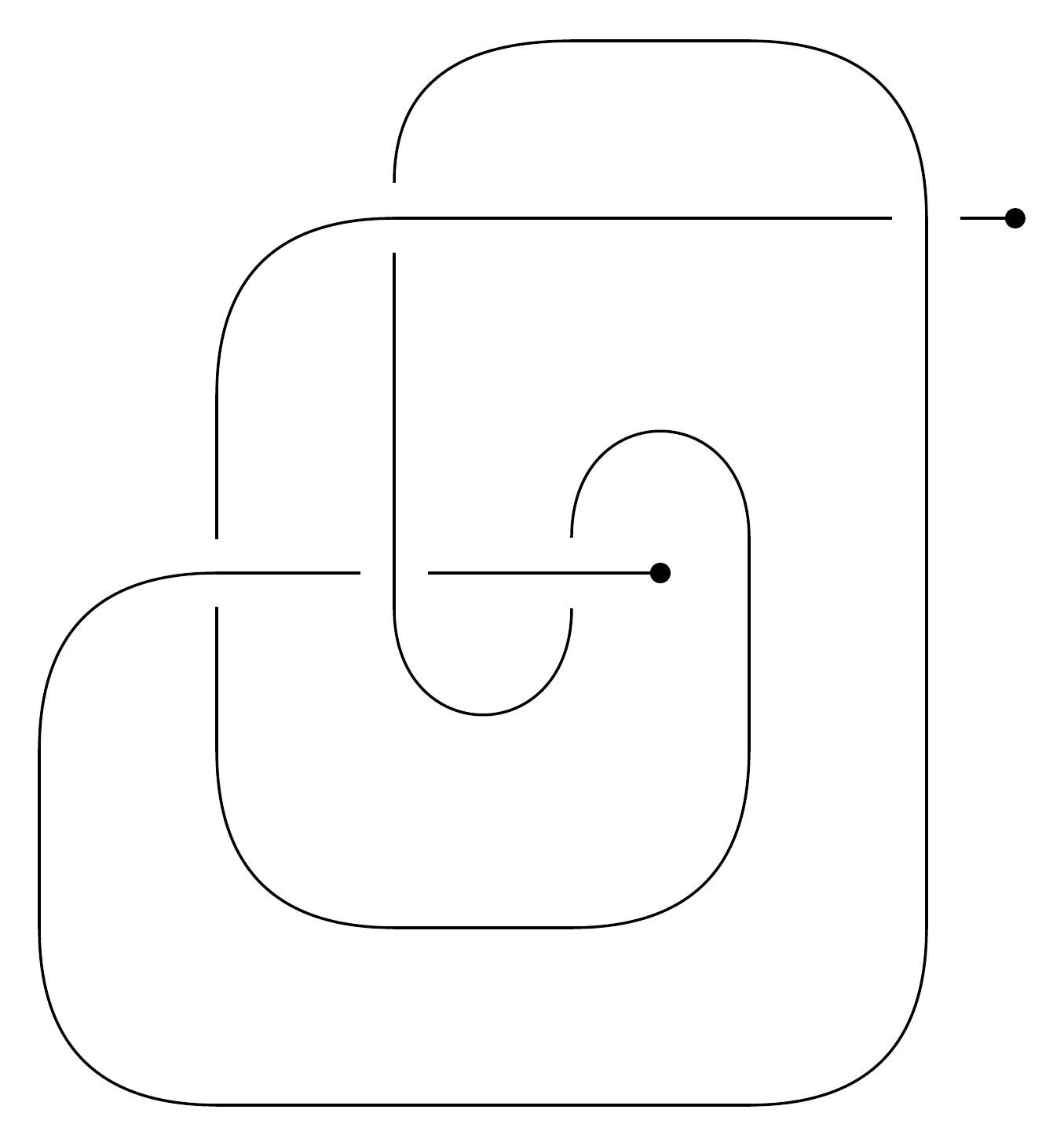}\\
\textcolor{blue}{$5_{17}$}
\vspace{1cm}
\end{minipage}
\begin{minipage}[t]{.25\linewidth}
\centering
\includegraphics[width=0.9\textwidth,height=3.5cm,keepaspectratio]{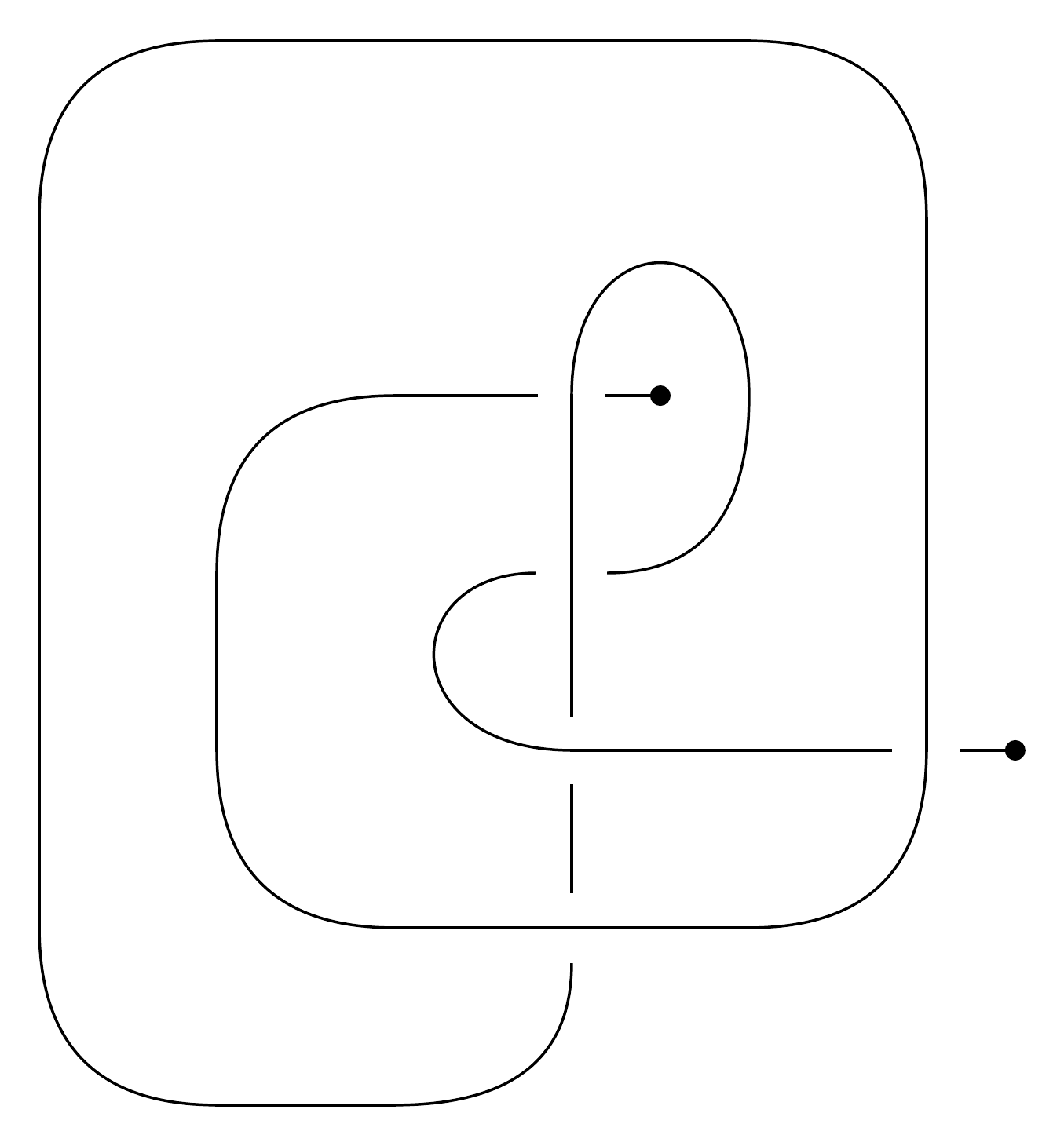}\\
\textcolor{blue}{$5_{18}$}
\vspace{1cm}
\end{minipage}
\begin{minipage}[t]{.25\linewidth}
\centering
\includegraphics[width=0.9\textwidth,height=3.5cm,keepaspectratio]{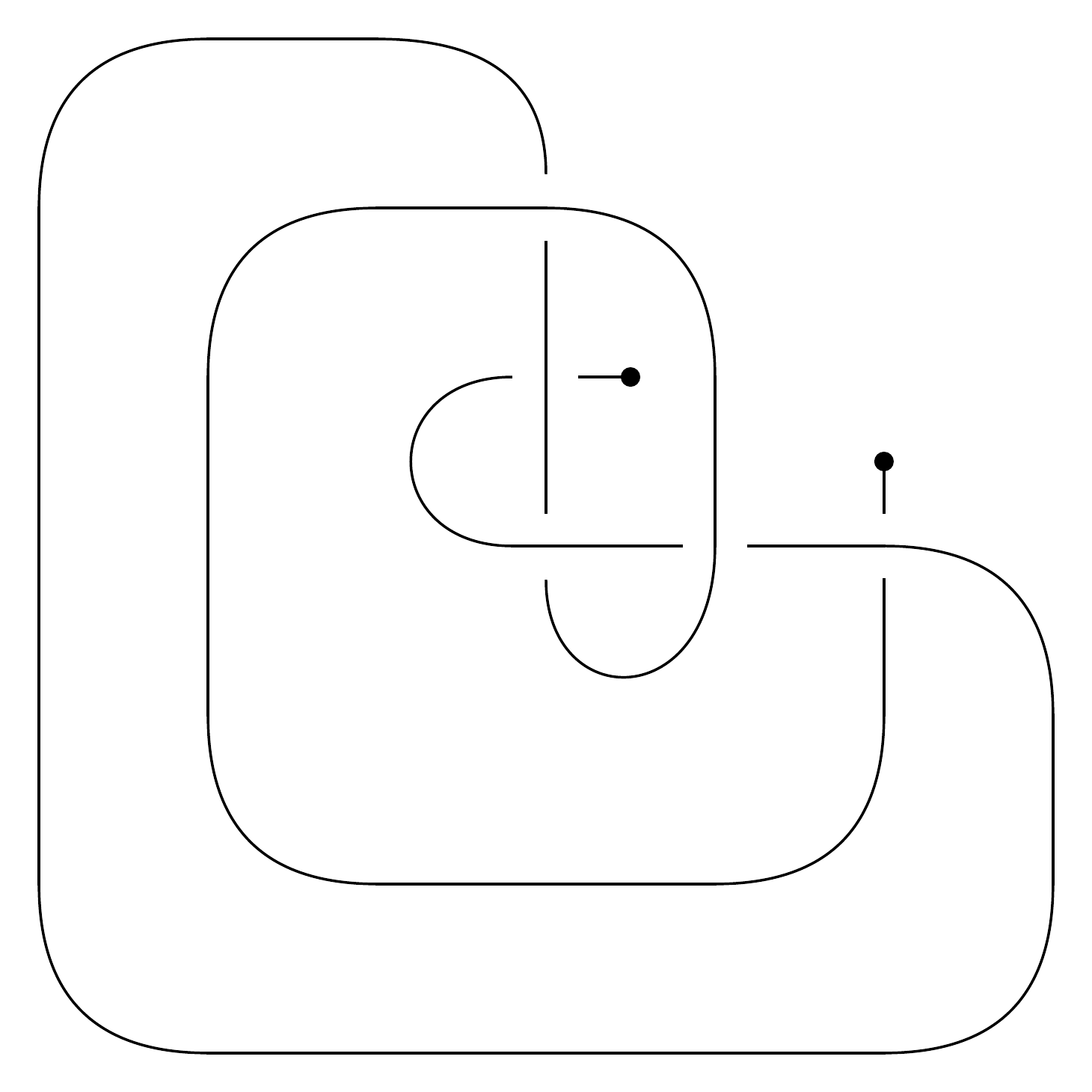}\\
\textcolor{blue}{$5_{19}$}
\vspace{1cm}
\end{minipage}
\begin{minipage}[t]{.25\linewidth}
\centering
\includegraphics[width=0.9\textwidth,height=3.5cm,keepaspectratio]{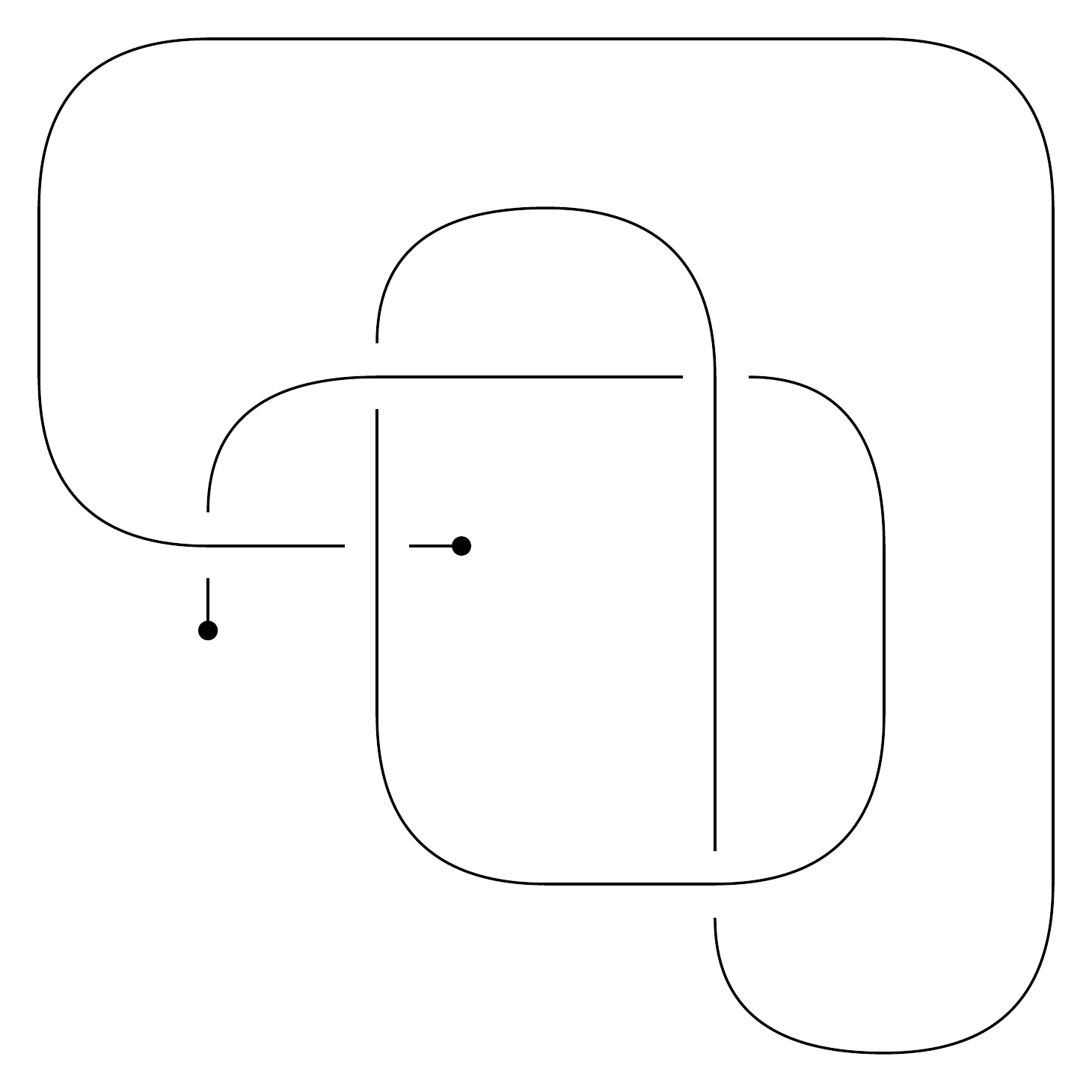}\\
\textcolor{blue}{$5_{20}$}
\vspace{1cm}
\end{minipage}
\begin{minipage}[t]{.25\linewidth}
\centering
\includegraphics[width=0.9\textwidth,height=3.5cm,keepaspectratio]{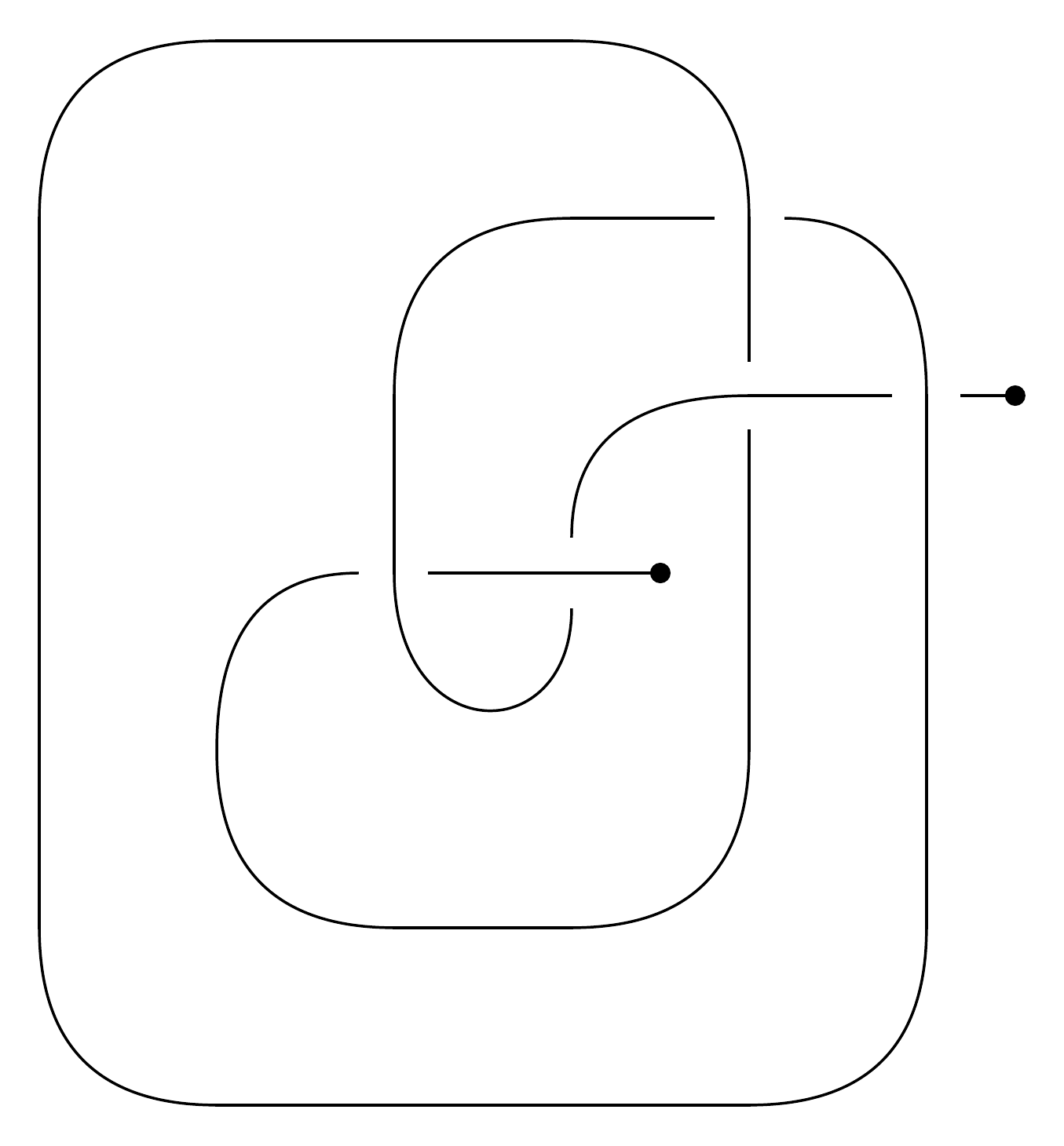}\\
\textcolor{blue}{$5_{21}$}
\vspace{1cm}
\end{minipage}
\begin{minipage}[t]{.25\linewidth}
\centering
\includegraphics[width=0.9\textwidth,height=3.5cm,keepaspectratio]{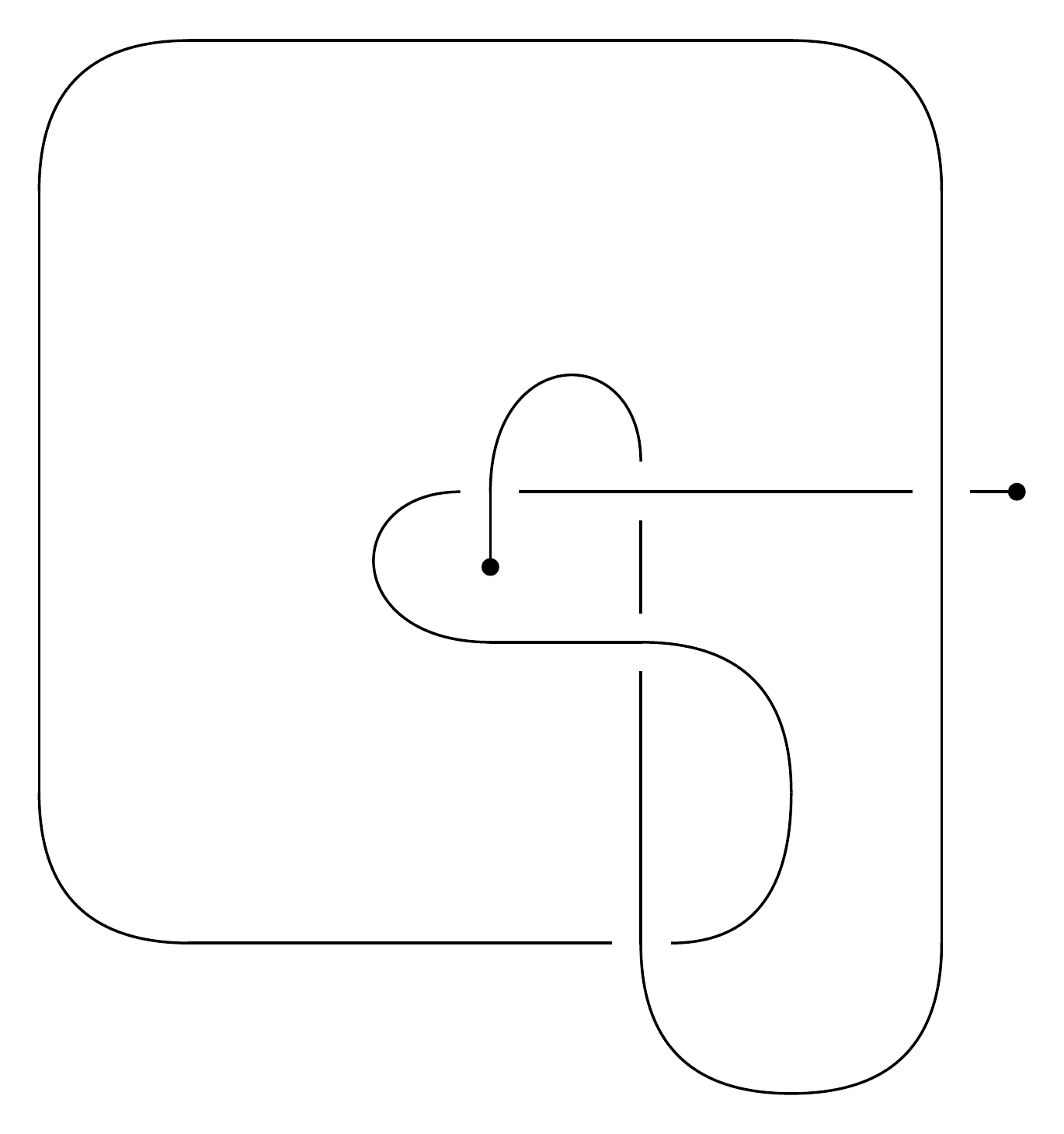}\\
\textcolor{blue}{$5_{22}$}
\vspace{1cm}
\end{minipage}
\begin{minipage}[t]{.25\linewidth}
\centering
\includegraphics[width=0.9\textwidth,height=3.5cm,keepaspectratio]{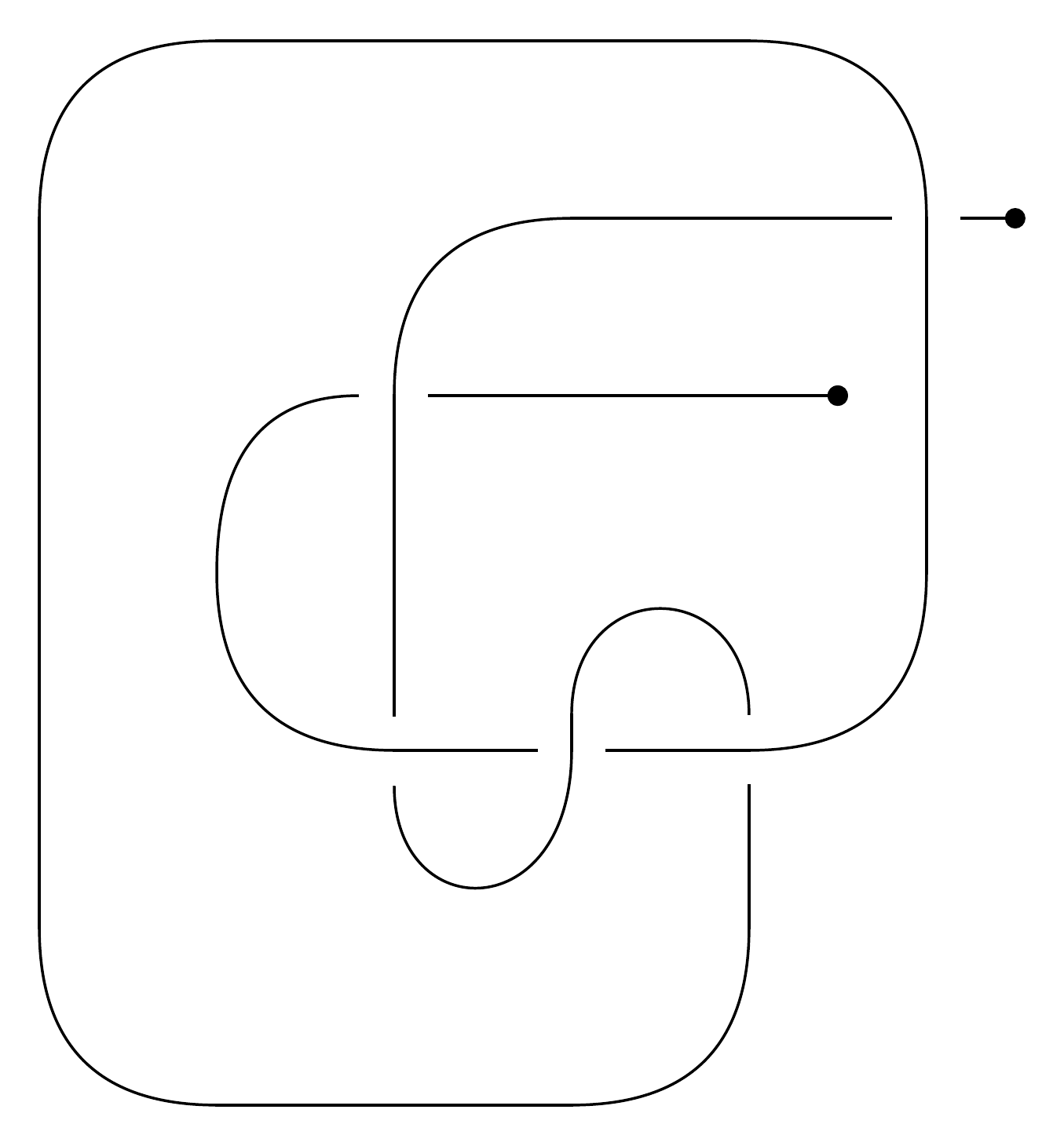}\\
\textcolor{blue}{$5_{23}$}
\vspace{1cm}
\end{minipage}
\begin{minipage}[t]{.25\linewidth}
\centering
\includegraphics[width=0.9\textwidth,height=3.5cm,keepaspectratio]{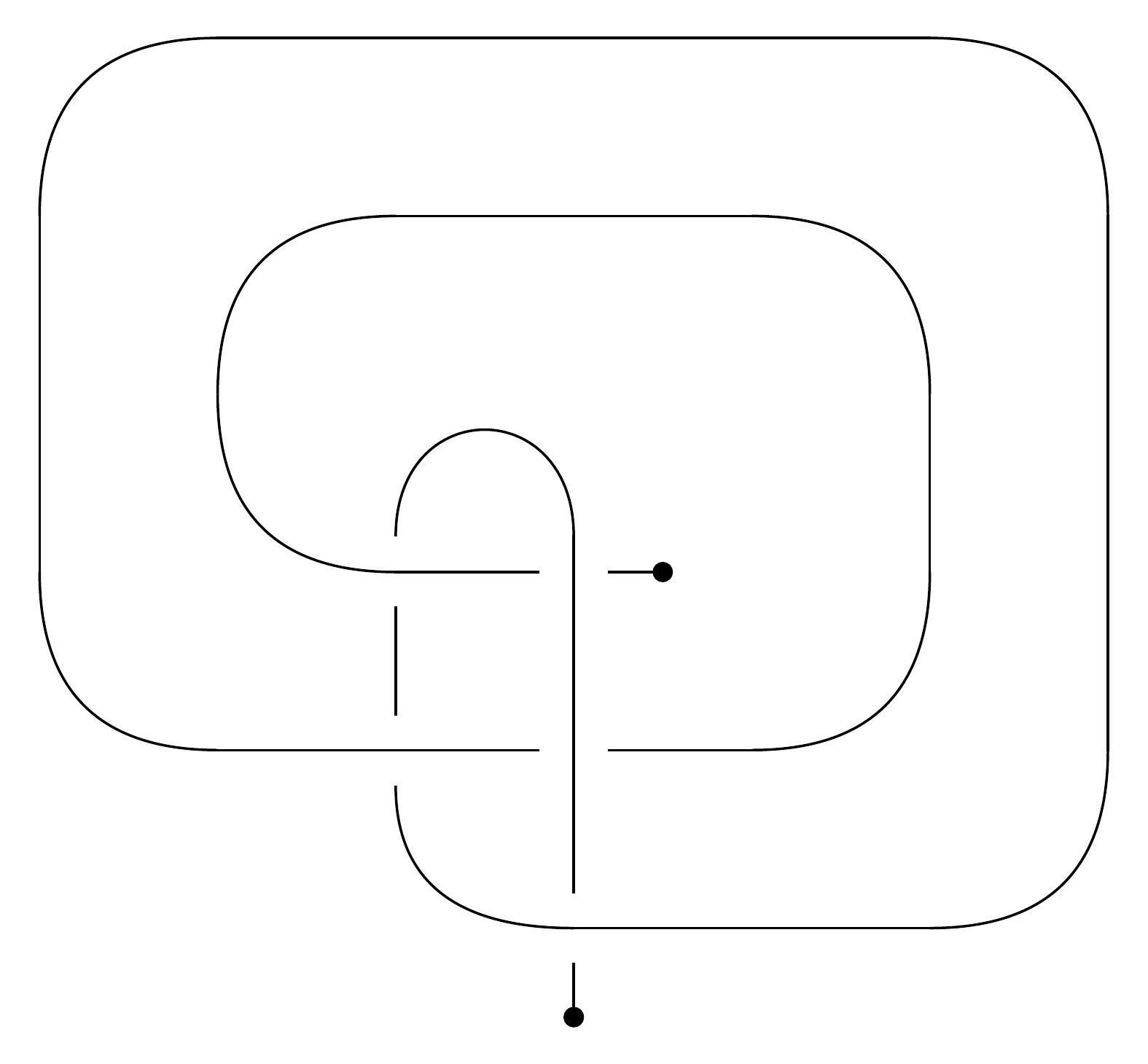}\\
\textcolor{blue}{$5_{24}$}
\vspace{1cm}
\end{minipage}
\begin{minipage}[t]{.25\linewidth}
\centering
\includegraphics[width=0.9\textwidth,height=3.5cm,keepaspectratio]{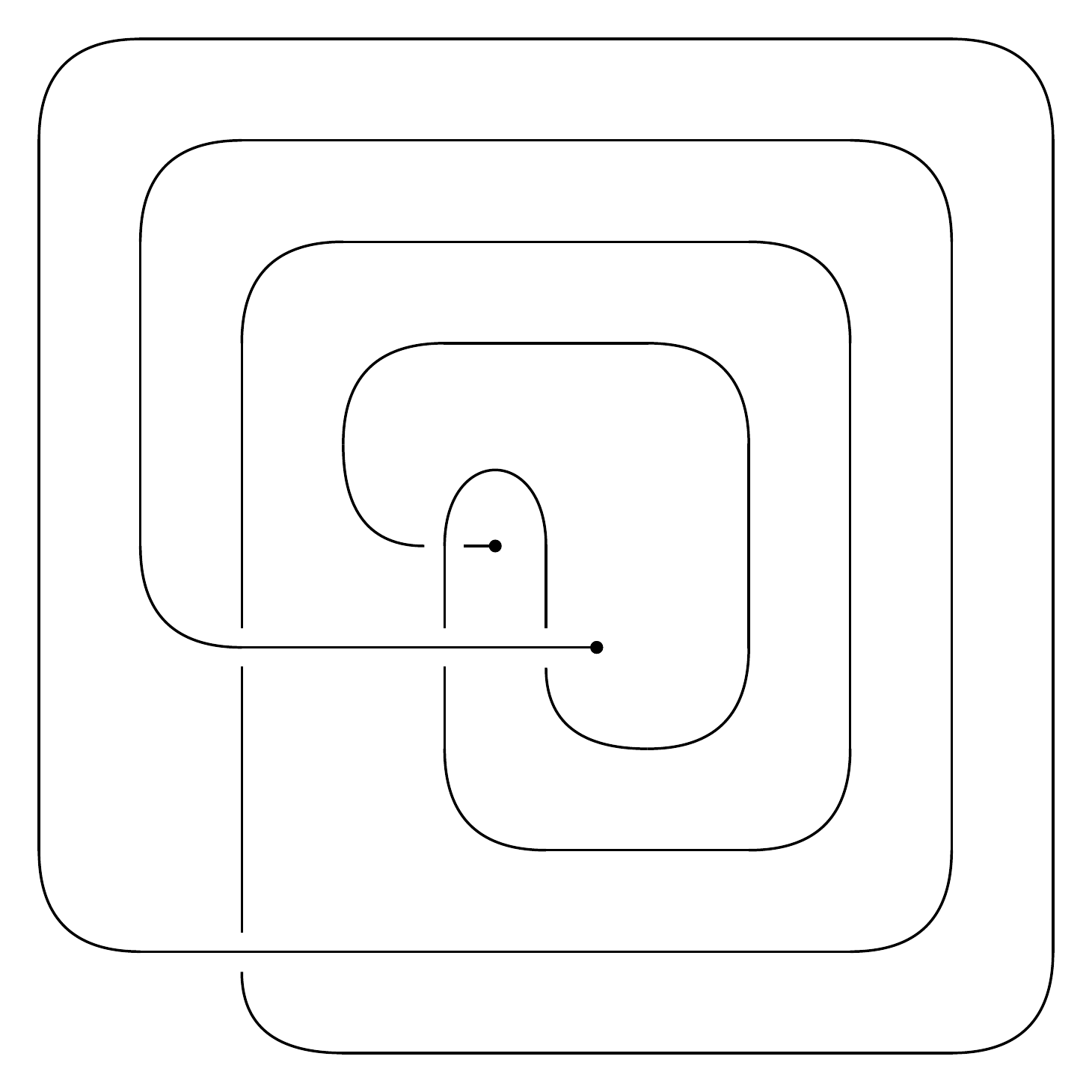}\\
\textcolor{black}{$5_{25}$}
\vspace{1cm}
\end{minipage}
\begin{minipage}[t]{.25\linewidth}
\centering
\includegraphics[width=0.9\textwidth,height=3.5cm,keepaspectratio]{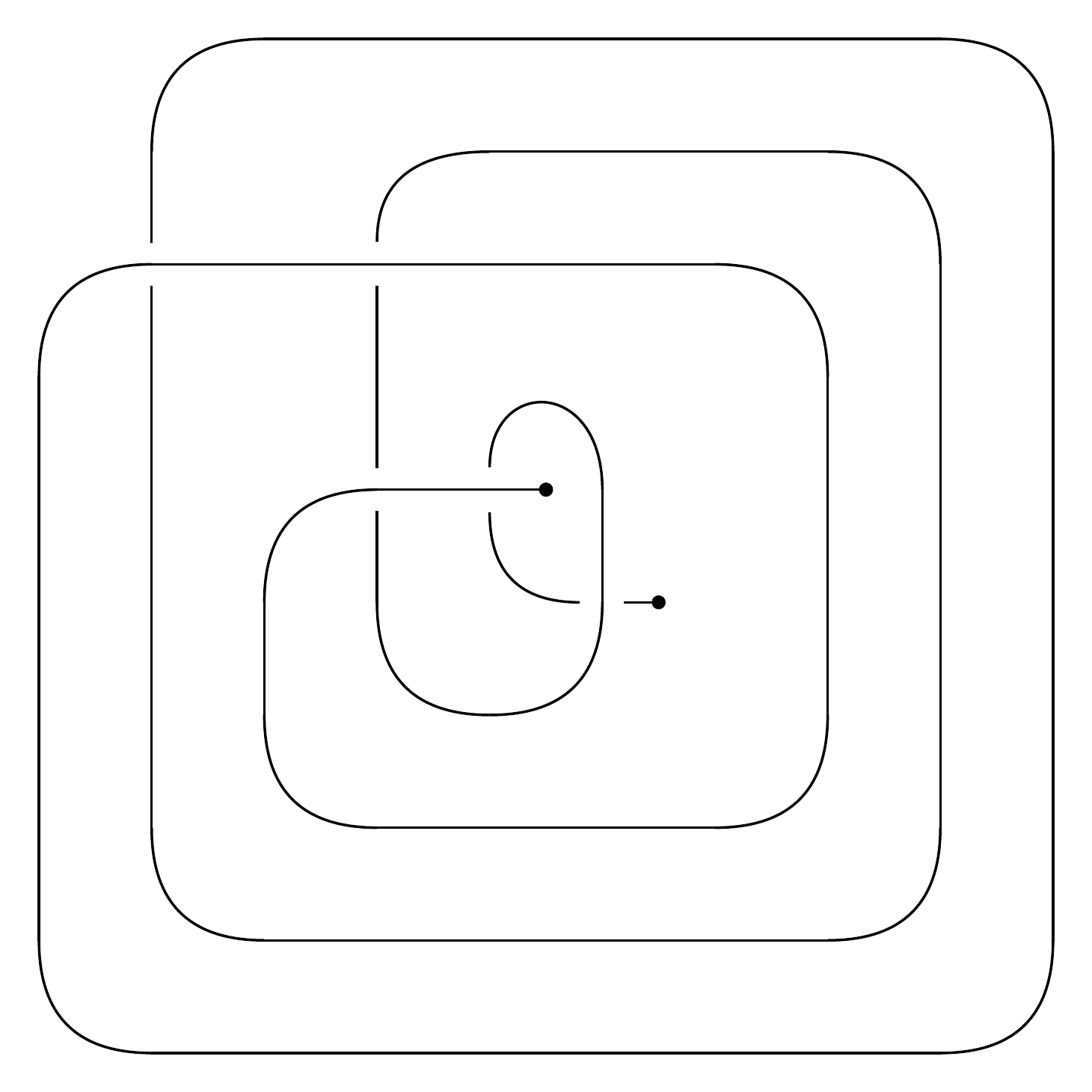}\\
\textcolor{black}{$5_{26}$}
\vspace{1cm}
\end{minipage}
\begin{minipage}[t]{.25\linewidth}
\centering
\includegraphics[width=0.9\textwidth,height=3.5cm,keepaspectratio]{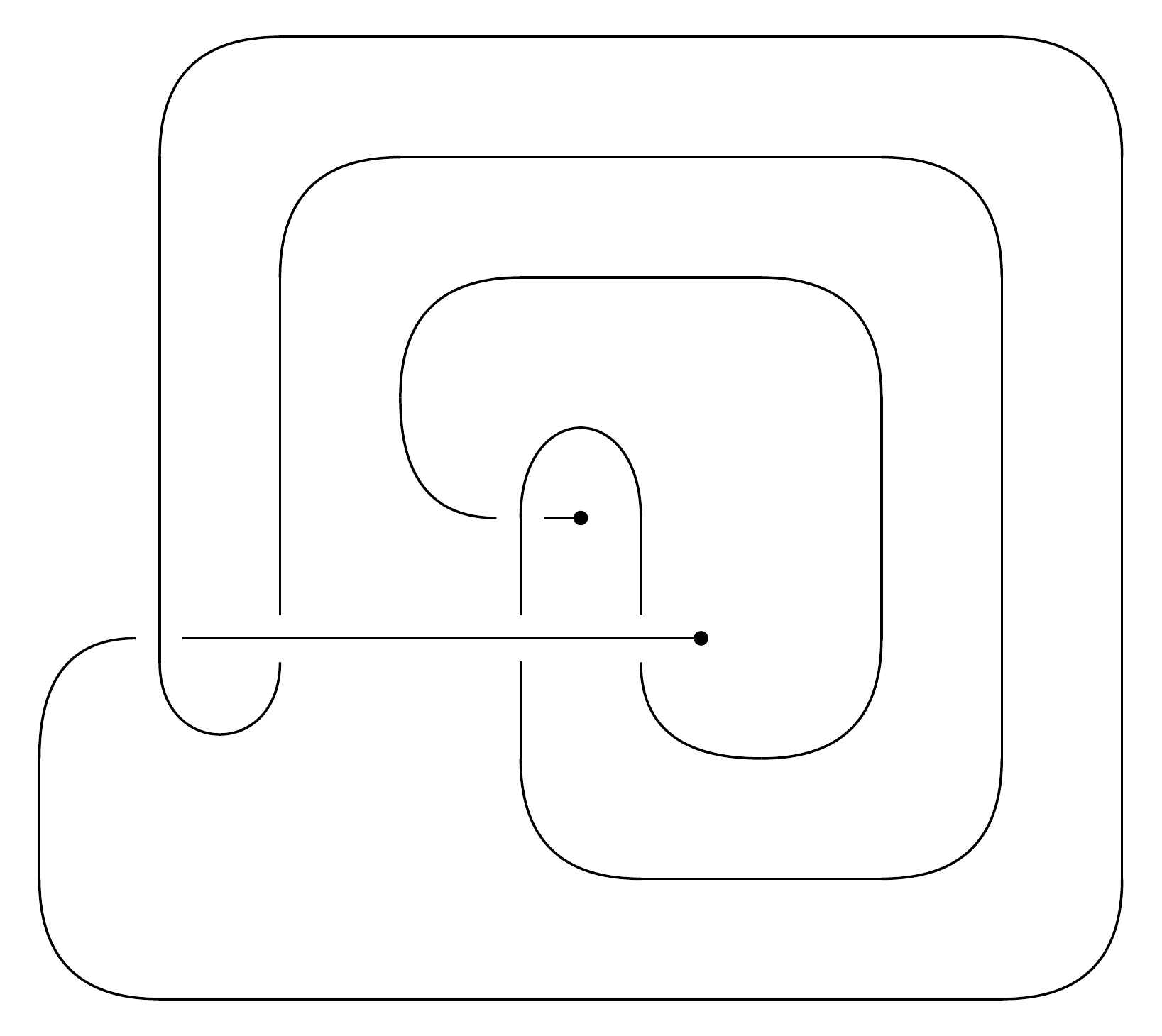}\\
\textcolor{black}{$5_{27}$}
\vspace{1cm}
\end{minipage}
\begin{minipage}[t]{.25\linewidth}
\centering
\includegraphics[width=0.9\textwidth,height=3.5cm,keepaspectratio]{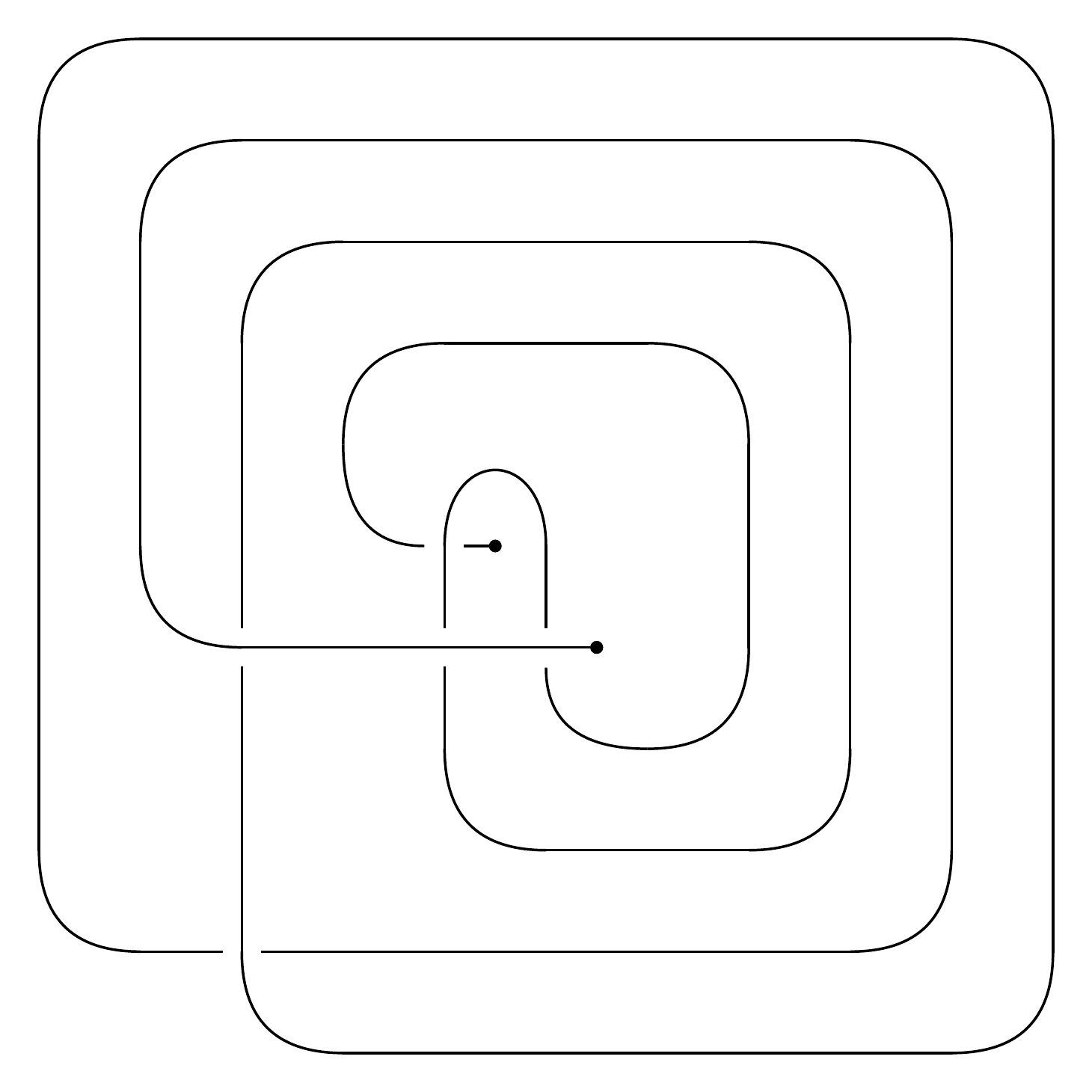}\\
\textcolor{black}{$5_{28}$}
\vspace{1cm}
\end{minipage}
\begin{minipage}[t]{.25\linewidth}
\centering
\includegraphics[width=0.9\textwidth,height=3.5cm,keepaspectratio]{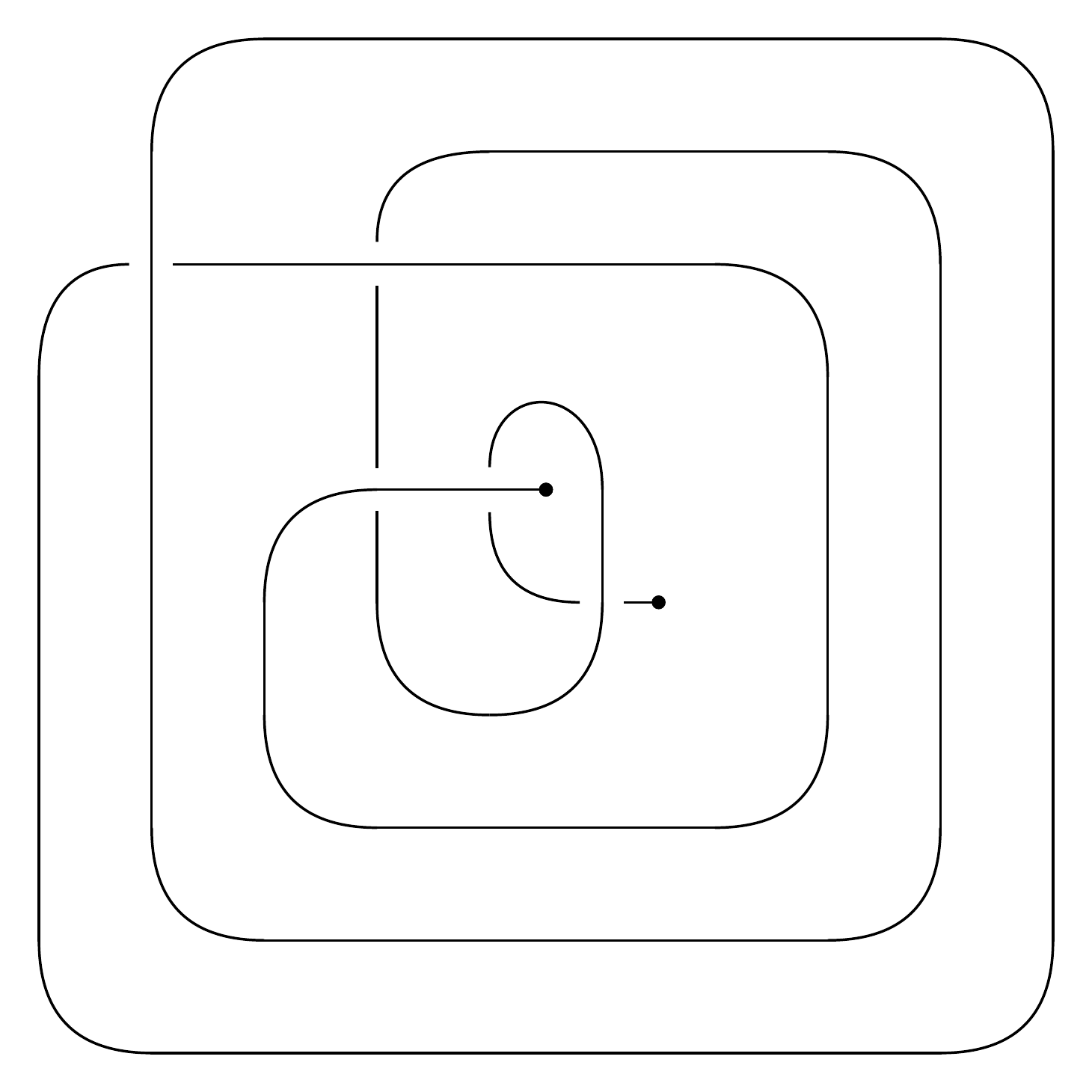}\\
\textcolor{black}{$5_{29}$}
\vspace{1cm}
\end{minipage}
\begin{minipage}[t]{.25\linewidth}
\centering
\includegraphics[width=0.9\textwidth,height=3.5cm,keepaspectratio]{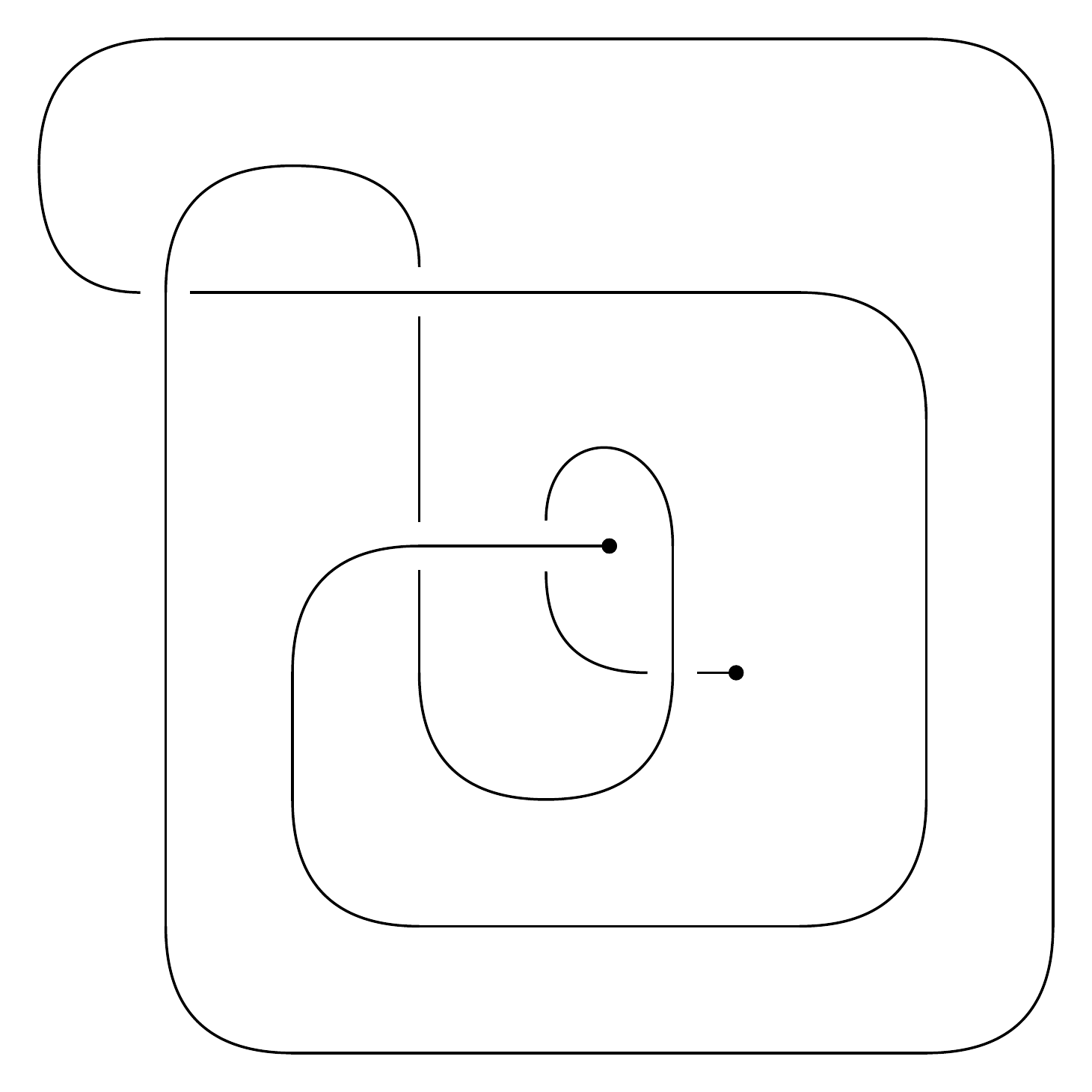}\\
\textcolor{black}{$5_{30}$}
\vspace{1cm}
\end{minipage}
\begin{minipage}[t]{.25\linewidth}
\centering
\includegraphics[width=0.9\textwidth,height=3.5cm,keepaspectratio]{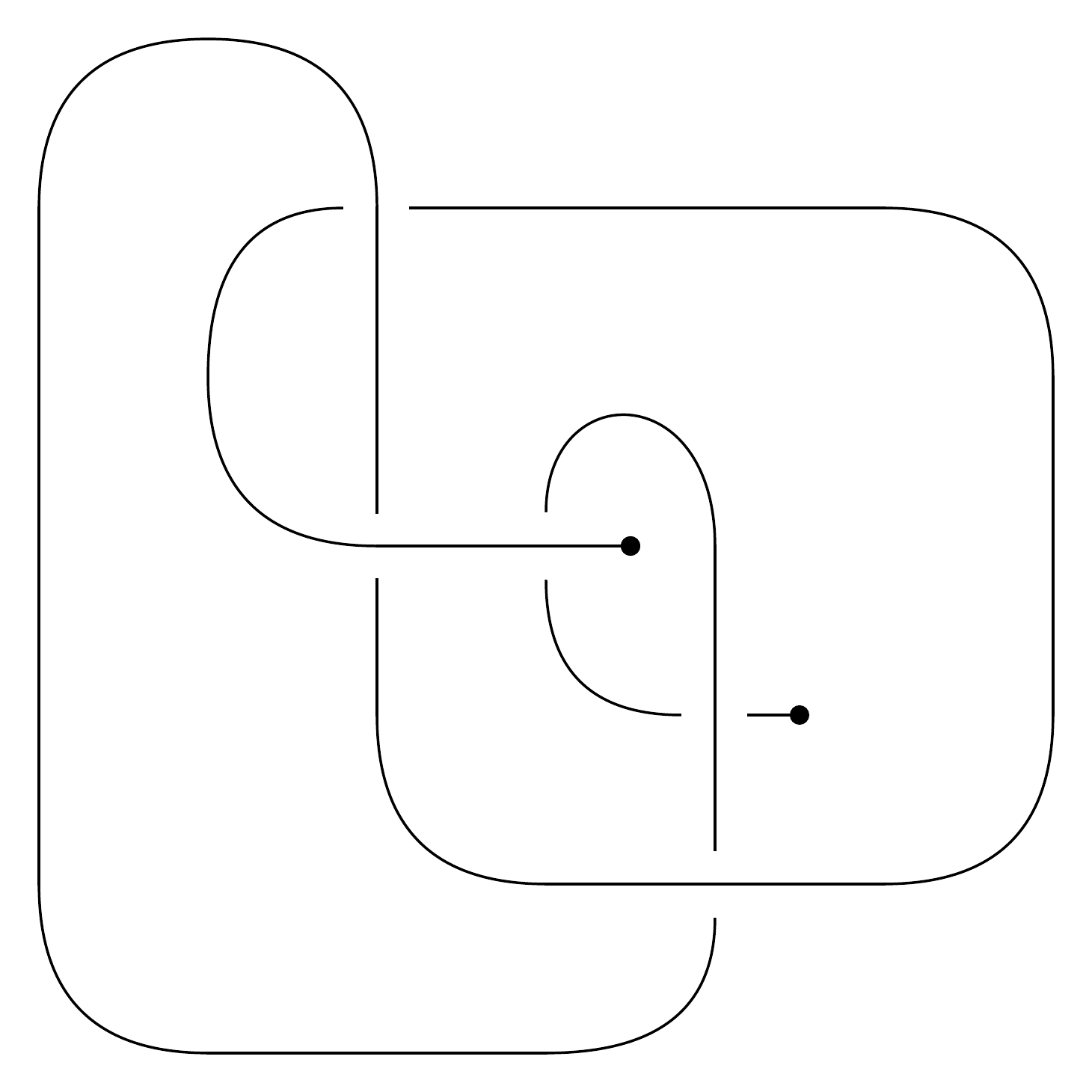}\\
\textcolor{black}{$5_{31}$}
\vspace{1cm}
\end{minipage}
\begin{minipage}[t]{.25\linewidth}
\centering
\includegraphics[width=0.9\textwidth,height=3.5cm,keepaspectratio]{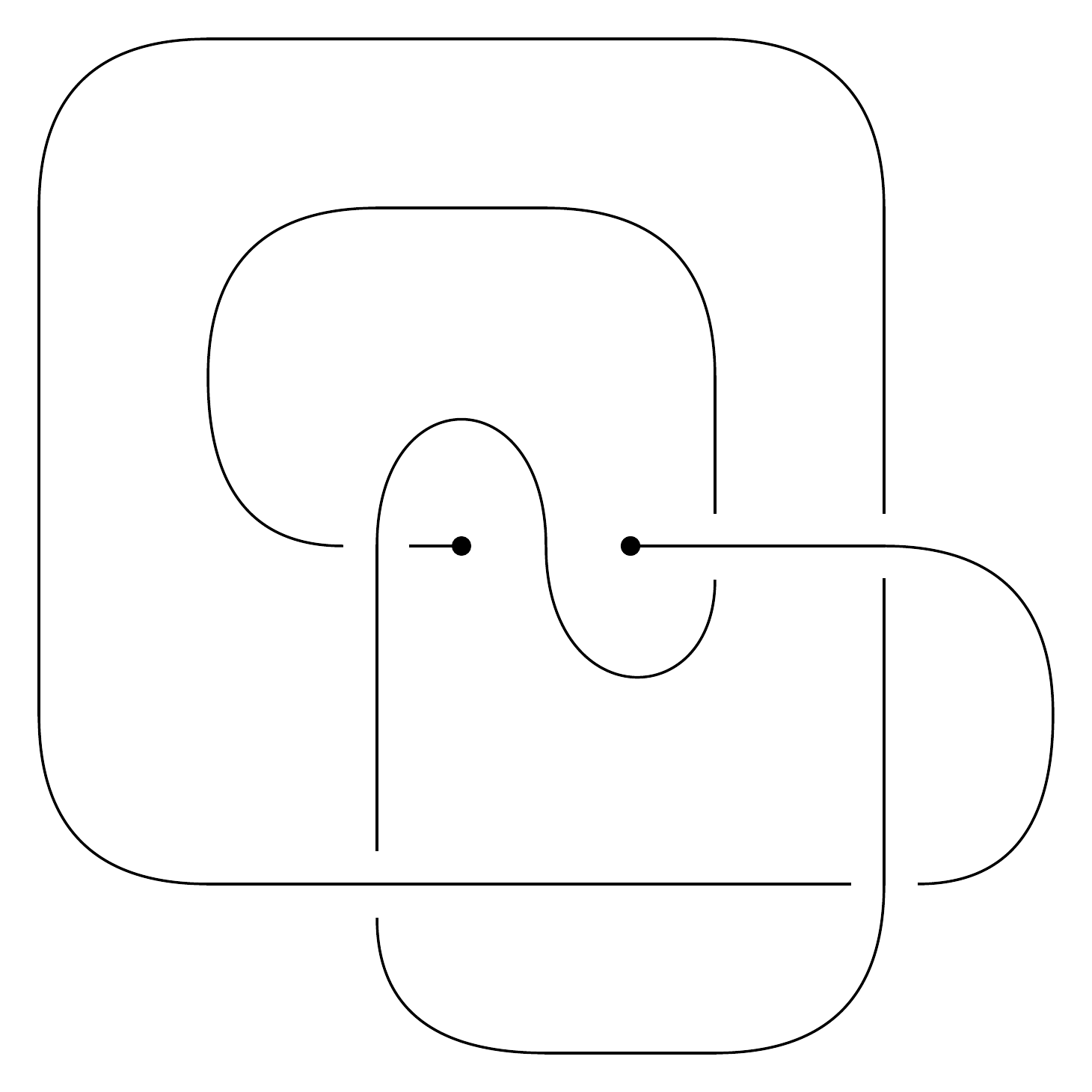}\\
\textcolor{black}{$5_{32}$}
\vspace{1cm}
\end{minipage}
\begin{minipage}[t]{.25\linewidth}
\centering
\includegraphics[width=0.9\textwidth,height=3.5cm,keepaspectratio]{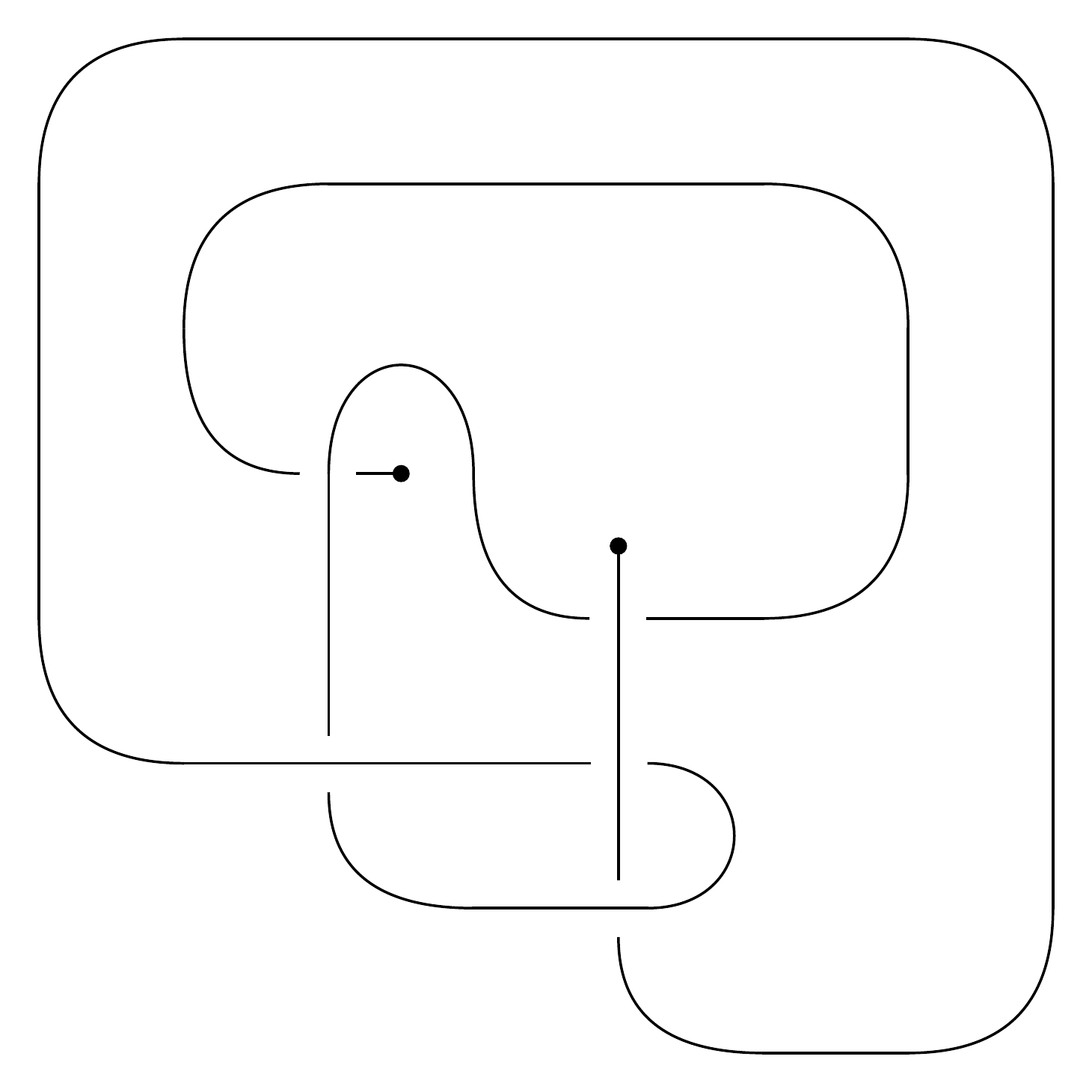}\\
\textcolor{black}{$5_{33}$}
\vspace{1cm}
\end{minipage}
\begin{minipage}[t]{.25\linewidth}
\centering
\includegraphics[width=0.9\textwidth,height=3.5cm,keepaspectratio]{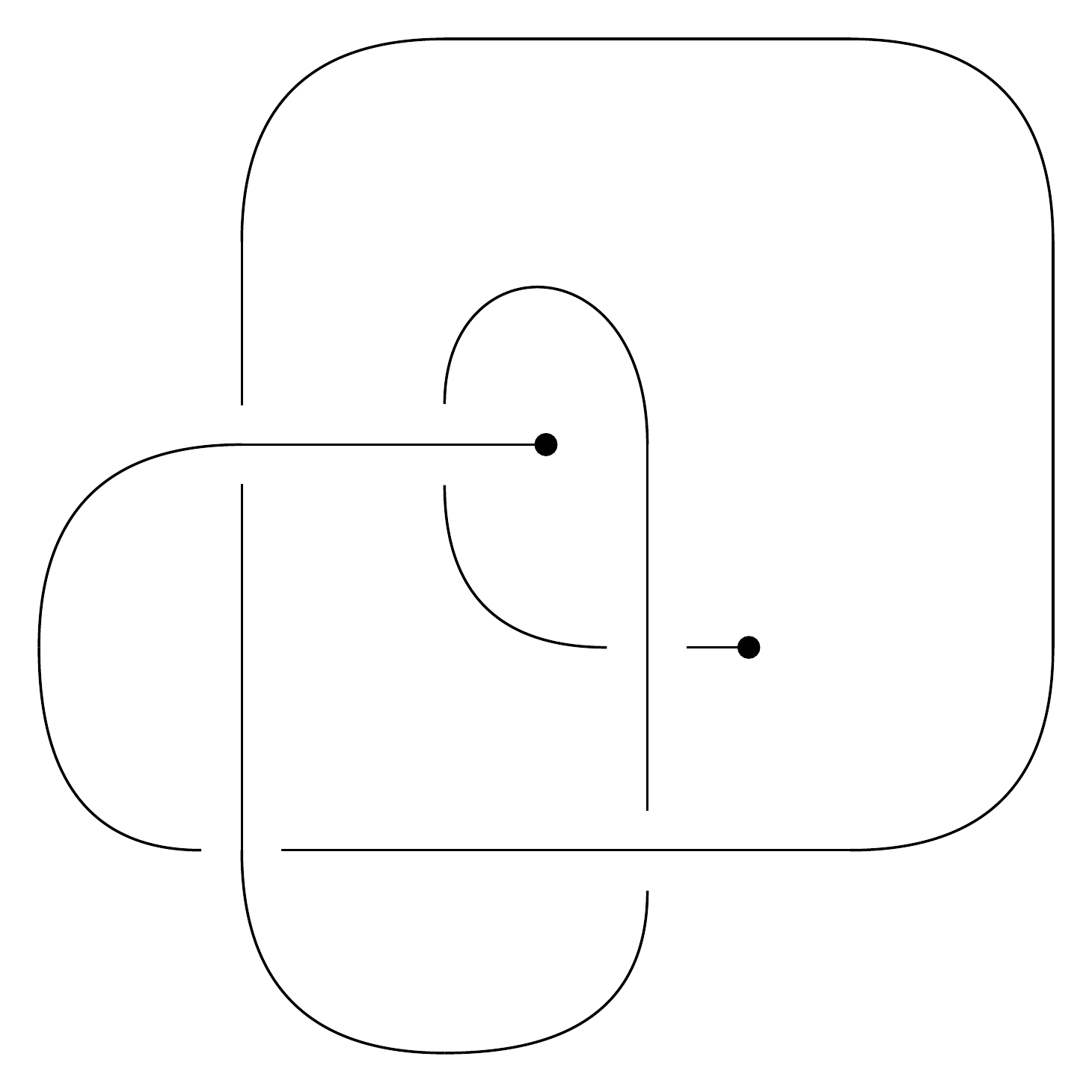}\\
\textcolor{black}{$5_{34}$}
\vspace{1cm}
\end{minipage}
\begin{minipage}[t]{.25\linewidth}
\centering
\includegraphics[width=0.9\textwidth,height=3.5cm,keepaspectratio]{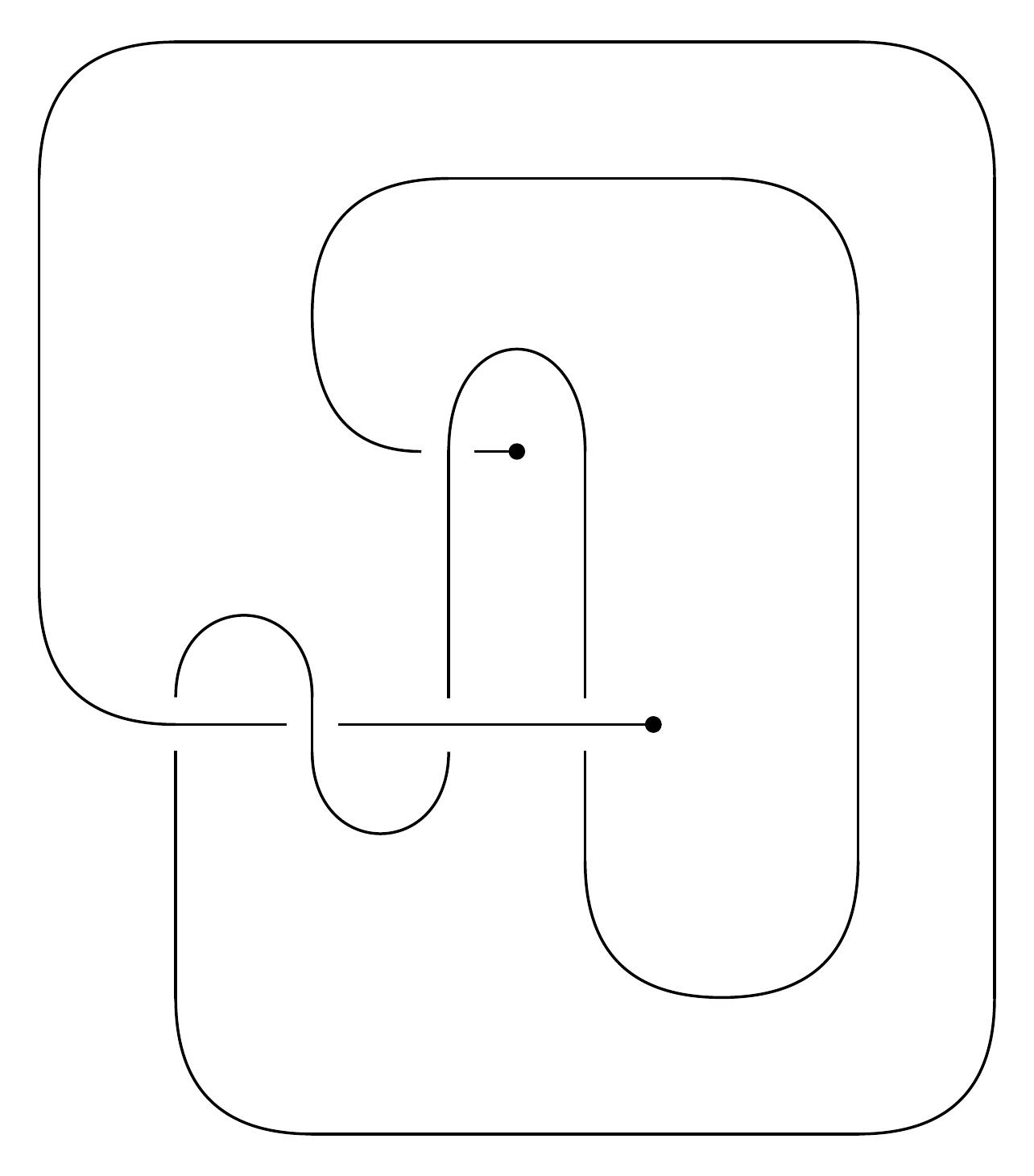}\\
\textcolor{black}{$5_{35}$}
\vspace{1cm}
\end{minipage}
\begin{minipage}[t]{.25\linewidth}
\centering
\includegraphics[width=0.9\textwidth,height=3.5cm,keepaspectratio]{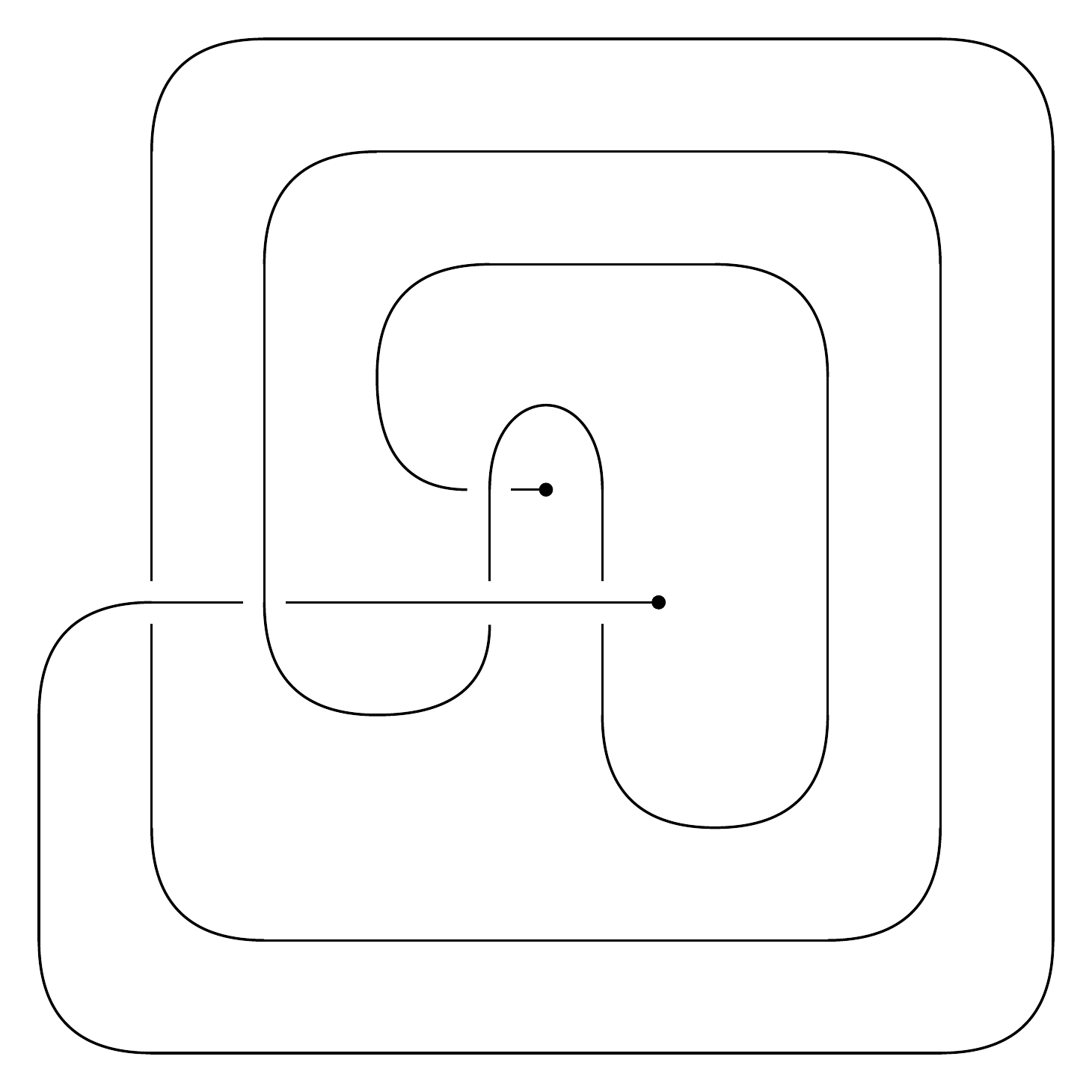}\\
\textcolor{black}{$5_{36}$}
\vspace{1cm}
\end{minipage}
\begin{minipage}[t]{.25\linewidth}
\centering
\includegraphics[width=0.9\textwidth,height=3.5cm,keepaspectratio]{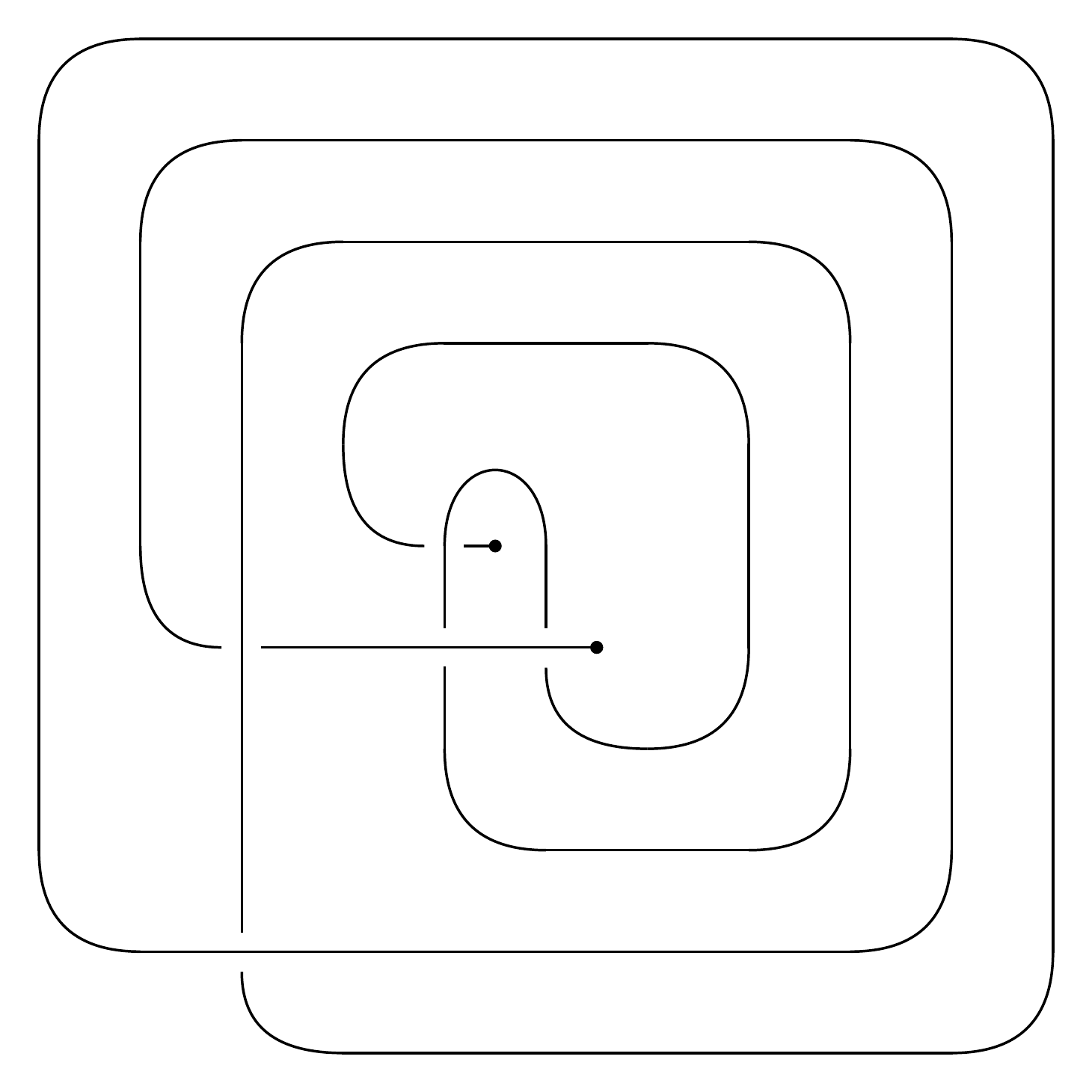}\\
\textcolor{black}{$5_{37}$}
\vspace{1cm}
\end{minipage}
\begin{minipage}[t]{.25\linewidth}
\centering
\includegraphics[width=0.9\textwidth,height=3.5cm,keepaspectratio]{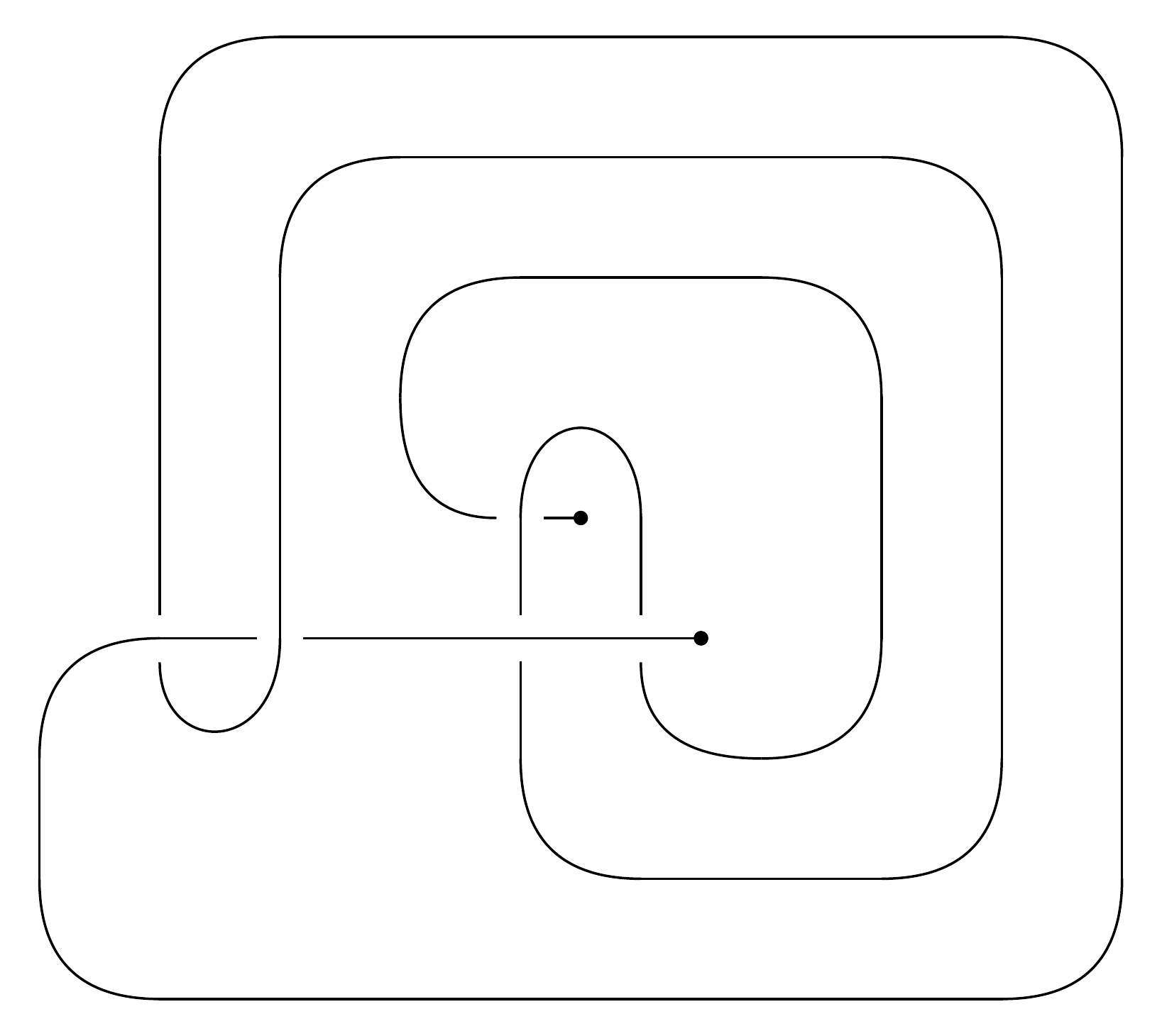}\\
\textcolor{black}{$5_{38}$}
\vspace{1cm}
\end{minipage}
\begin{minipage}[t]{.25\linewidth}
\centering
\includegraphics[width=0.9\textwidth,height=3.5cm,keepaspectratio]{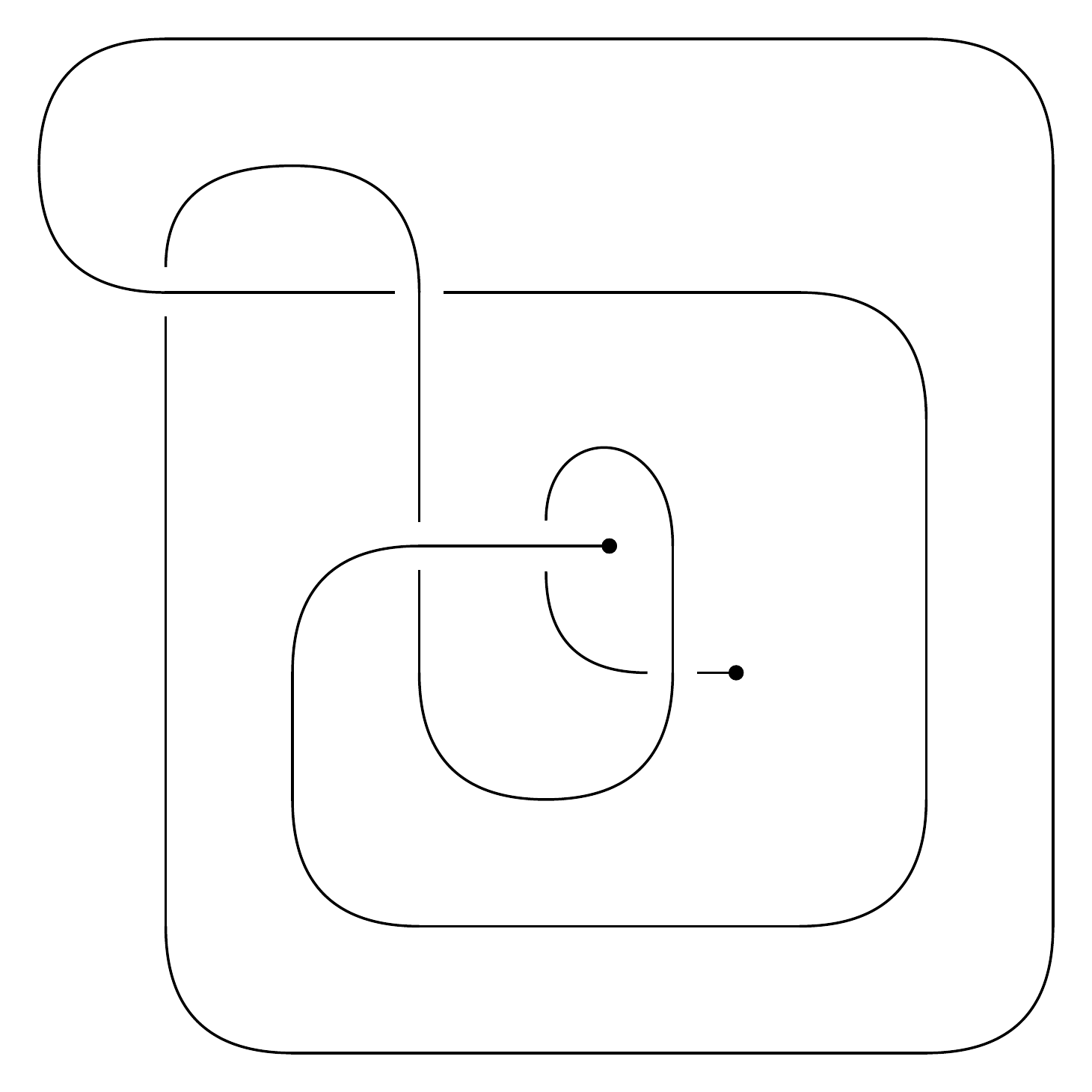}\\
\textcolor{black}{$5_{39}$}
\vspace{1cm}
\end{minipage}
\begin{minipage}[t]{.25\linewidth}
\centering
\includegraphics[width=0.9\textwidth,height=3.5cm,keepaspectratio]{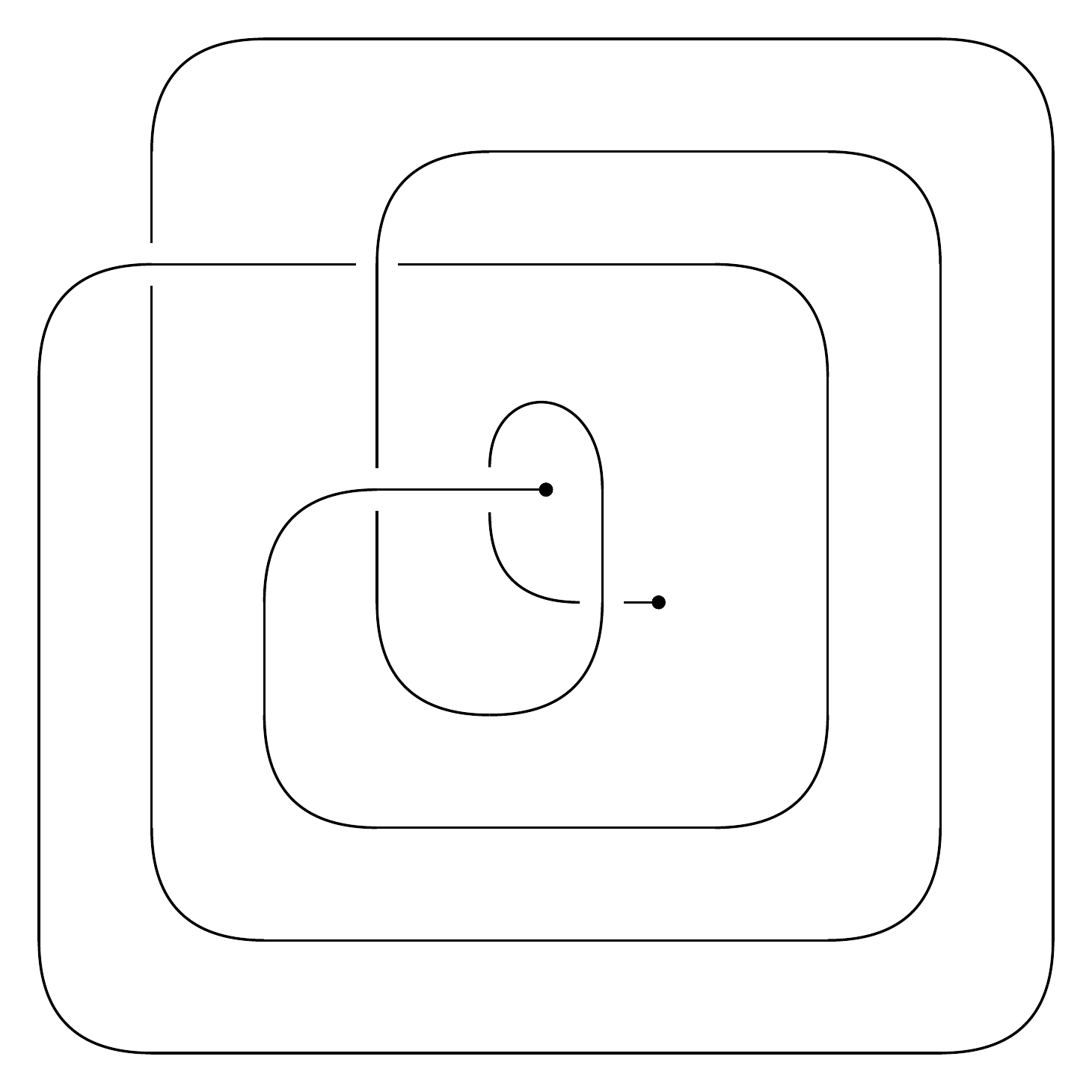}\\
\textcolor{black}{$5_{40}$}
\vspace{1cm}
\end{minipage}
\begin{minipage}[t]{.25\linewidth}
\centering
\includegraphics[width=0.9\textwidth,height=3.5cm,keepaspectratio]{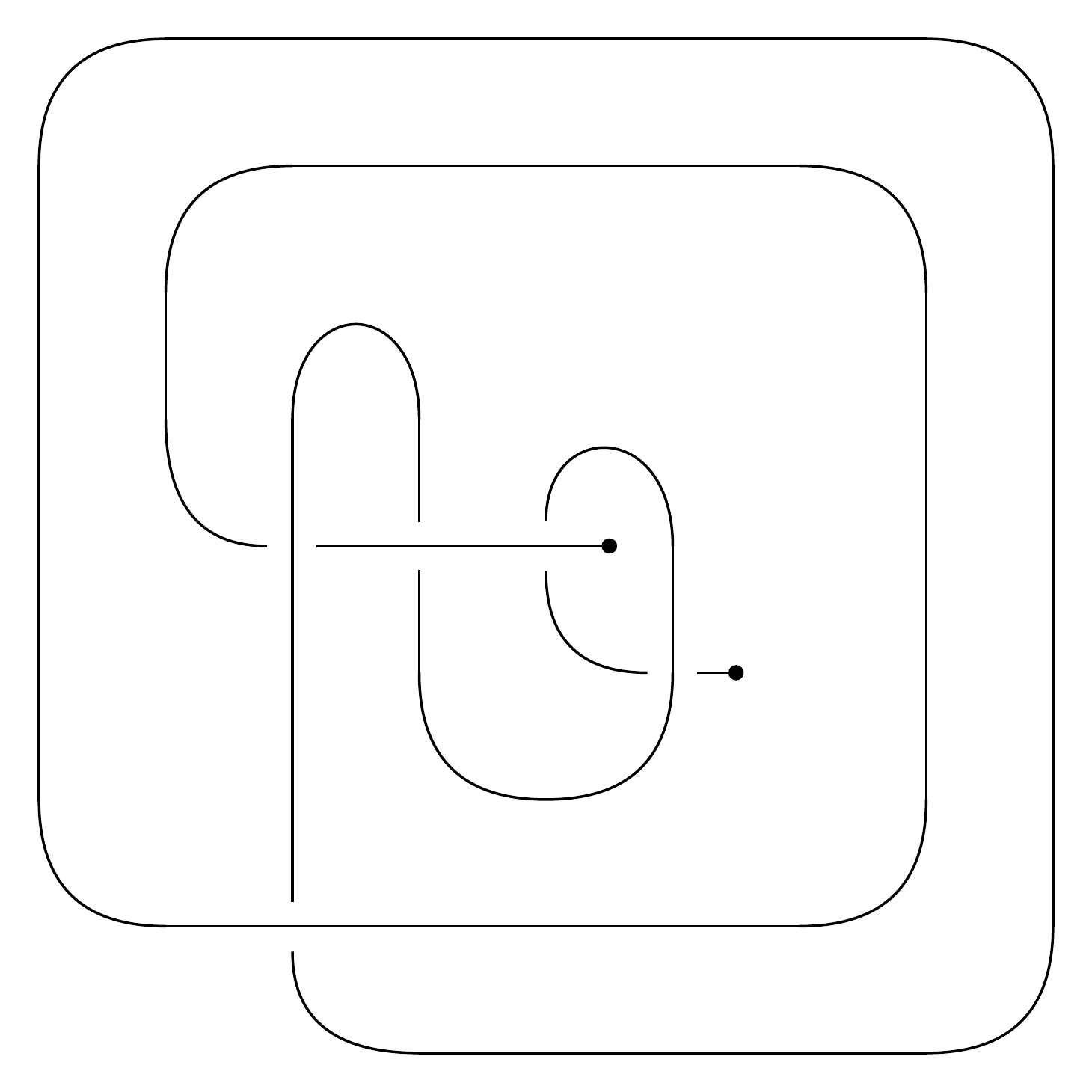}\\
\textcolor{black}{$5_{41}$}
\vspace{1cm}
\end{minipage}
\begin{minipage}[t]{.25\linewidth}
\centering
\includegraphics[width=0.9\textwidth,height=3.5cm,keepaspectratio]{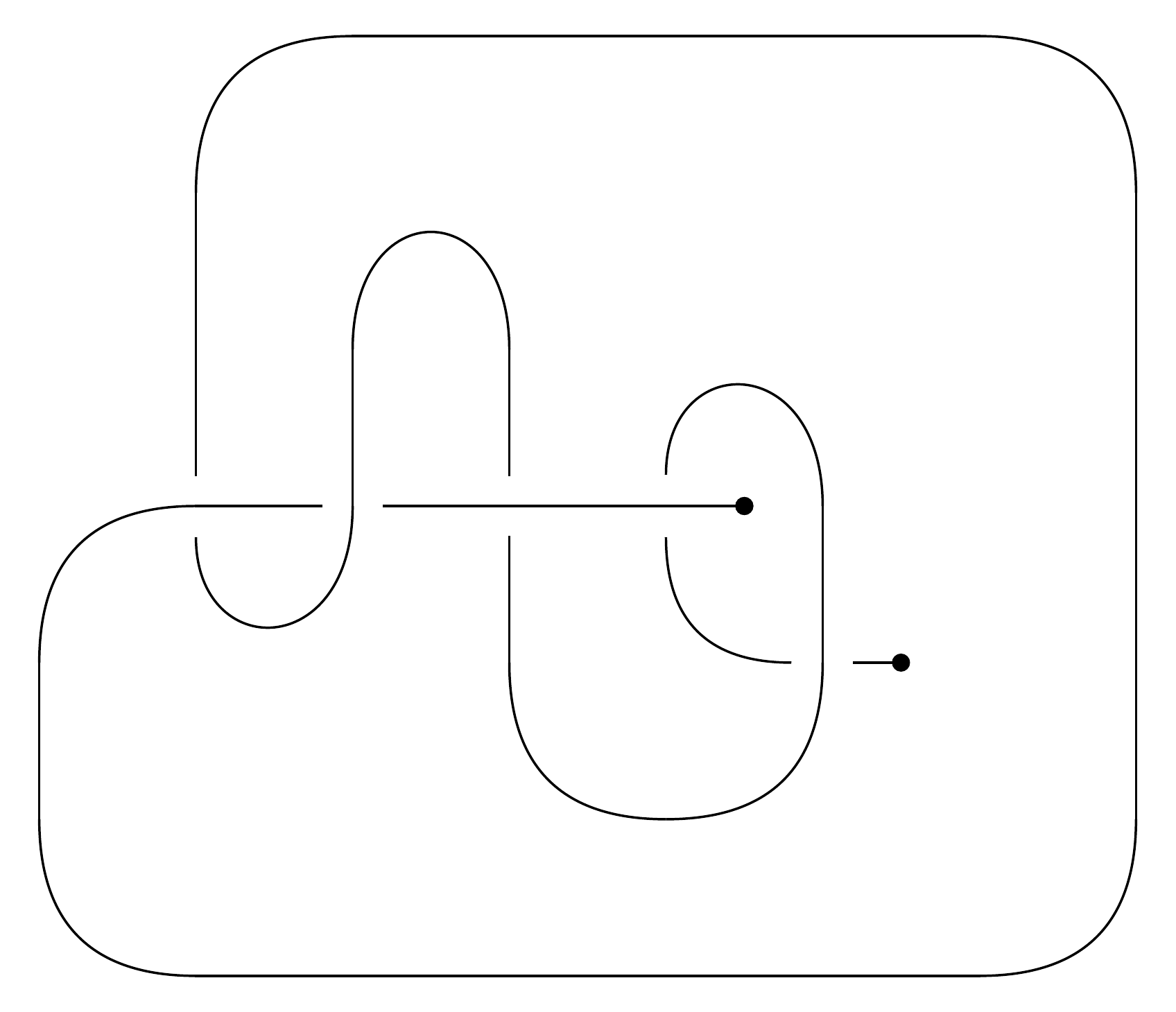}\\
\textcolor{black}{$5_{42}$}
\vspace{1cm}
\end{minipage}
\begin{minipage}[t]{.25\linewidth}
\centering
\includegraphics[width=0.9\textwidth,height=3.5cm,keepaspectratio]{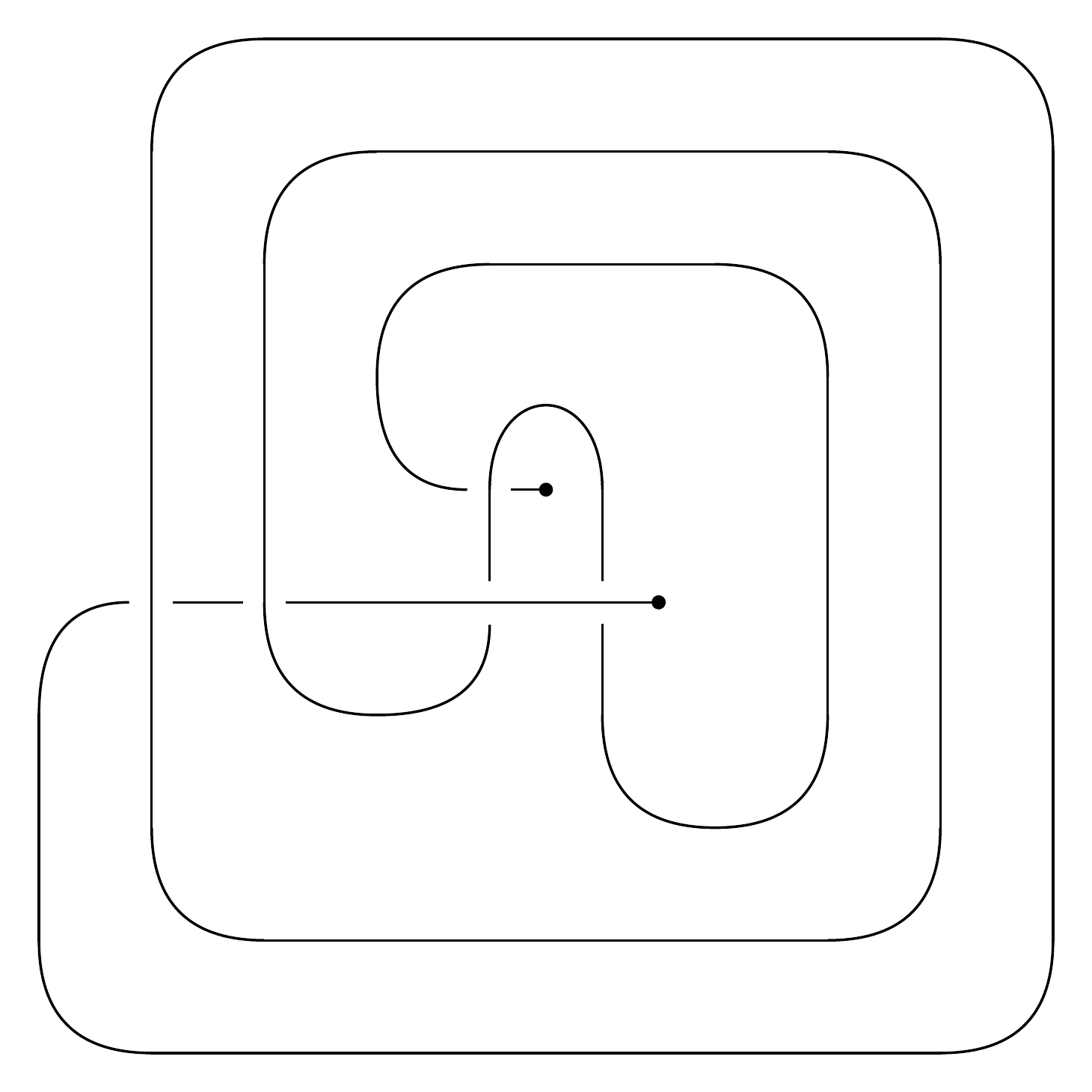}\\
\textcolor{black}{$5_{43}$}
\vspace{1cm}
\end{minipage}
\begin{minipage}[t]{.25\linewidth}
\centering
\includegraphics[width=0.9\textwidth,height=3.5cm,keepaspectratio]{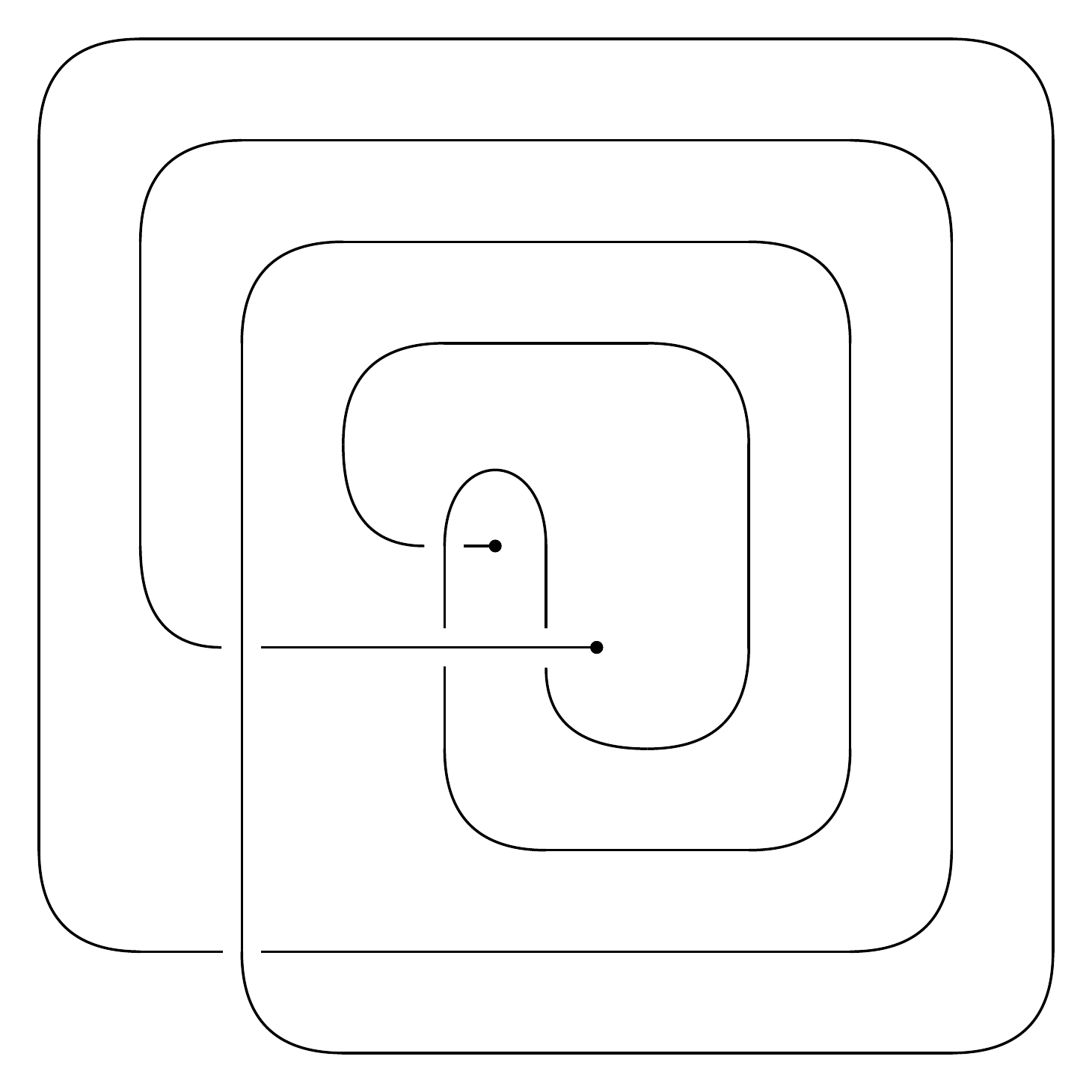}\\
\textcolor{black}{$5_{44}$}
\vspace{1cm}
\end{minipage}
\begin{minipage}[t]{.25\linewidth}
\centering
\includegraphics[width=0.9\textwidth,height=3.5cm,keepaspectratio]{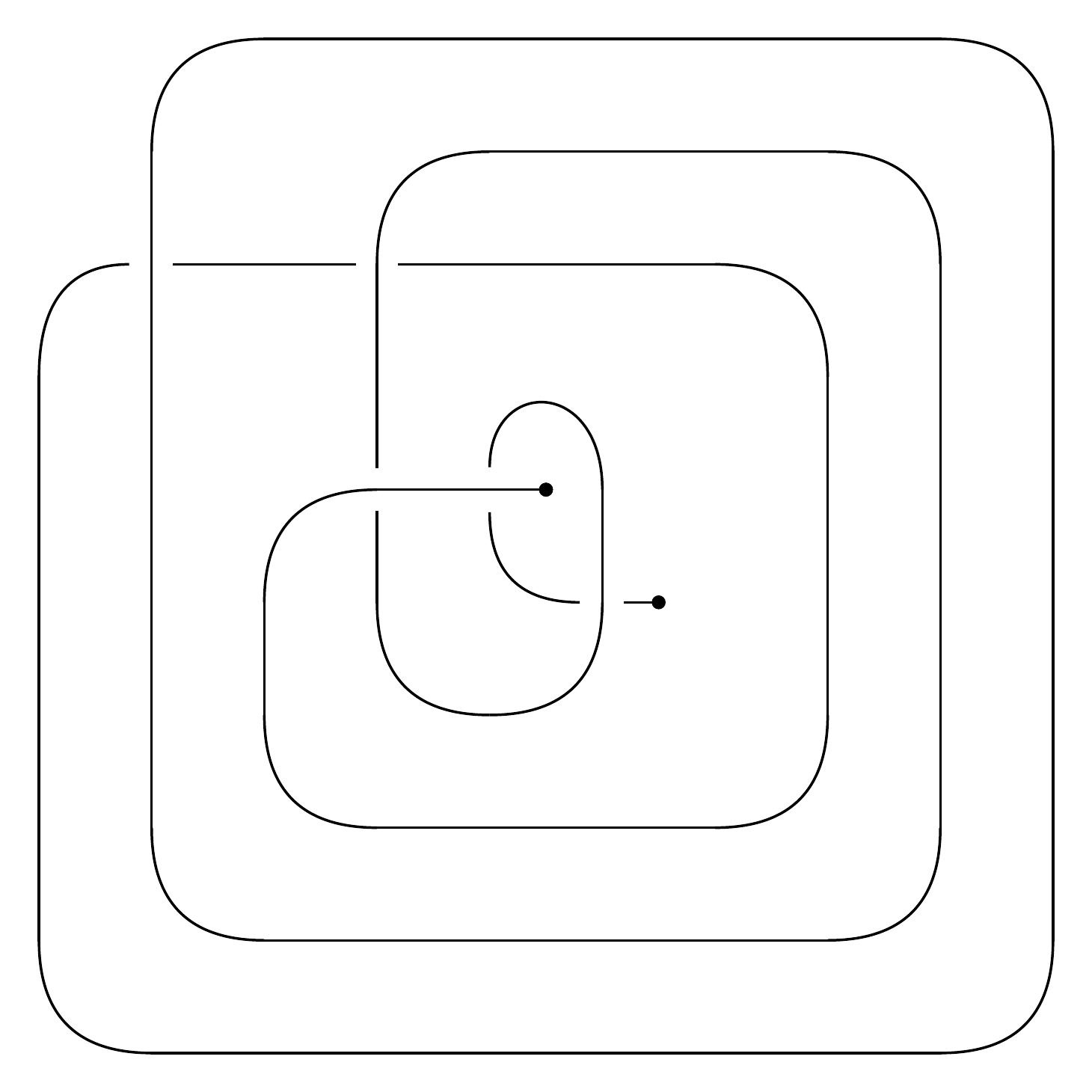}\\
\textcolor{black}{$5_{45}$}
\vspace{1cm}
\end{minipage}
\begin{minipage}[t]{.25\linewidth}
\centering
\includegraphics[width=0.9\textwidth,height=3.5cm,keepaspectratio]{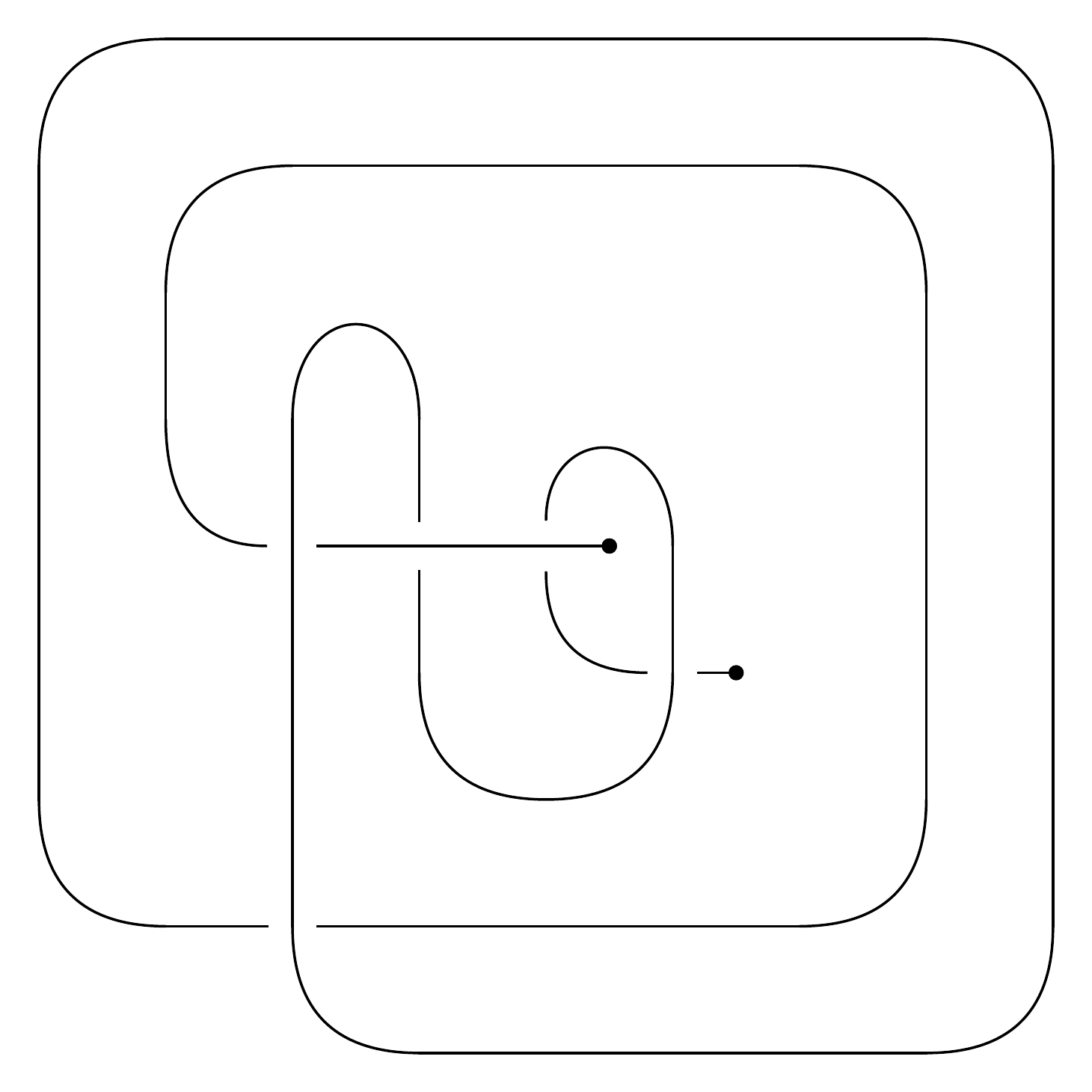}\\
\textcolor{black}{$5_{46}$}
\vspace{1cm}
\end{minipage}
\begin{minipage}[t]{.25\linewidth}
\centering
\includegraphics[width=0.9\textwidth,height=3.5cm,keepaspectratio]{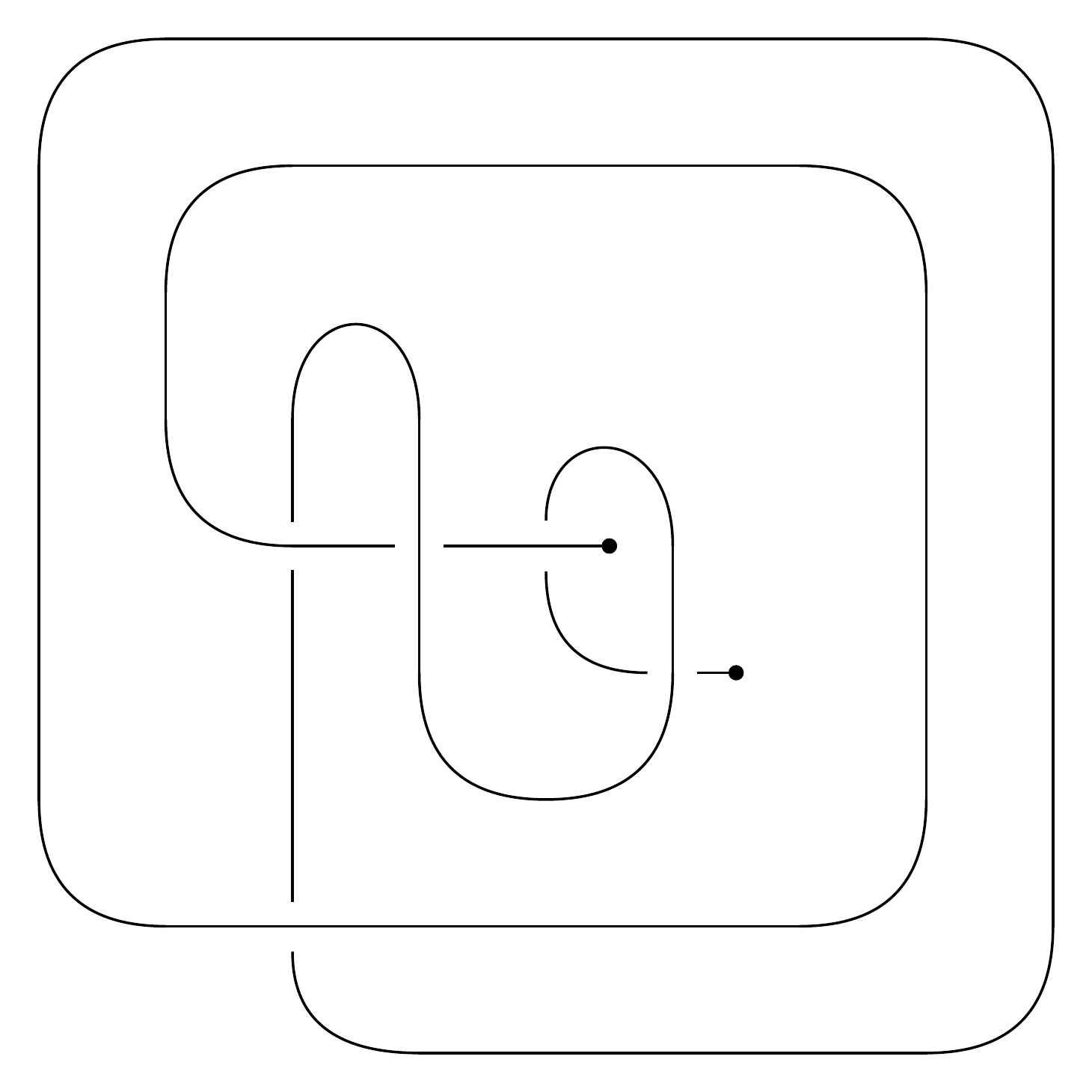}\\
\textcolor{black}{$5_{47}$}
\vspace{1cm}
\end{minipage}
\begin{minipage}[t]{.25\linewidth}
\centering
\includegraphics[width=0.9\textwidth,height=3.5cm,keepaspectratio]{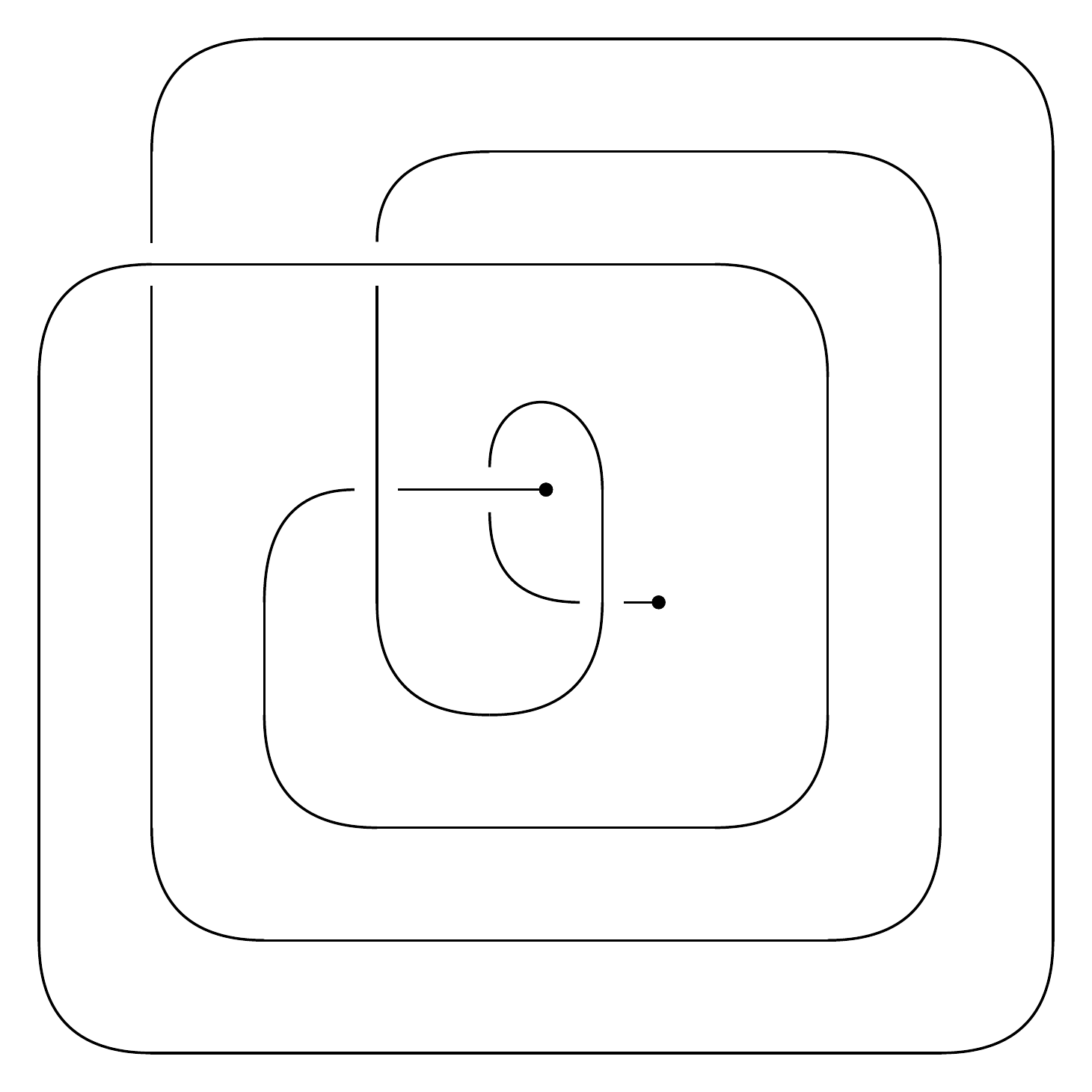}\\
\textcolor{black}{$5_{48}$}
\vspace{1cm}
\end{minipage}
\begin{minipage}[t]{.25\linewidth}
\centering
\includegraphics[width=0.9\textwidth,height=3.5cm,keepaspectratio]{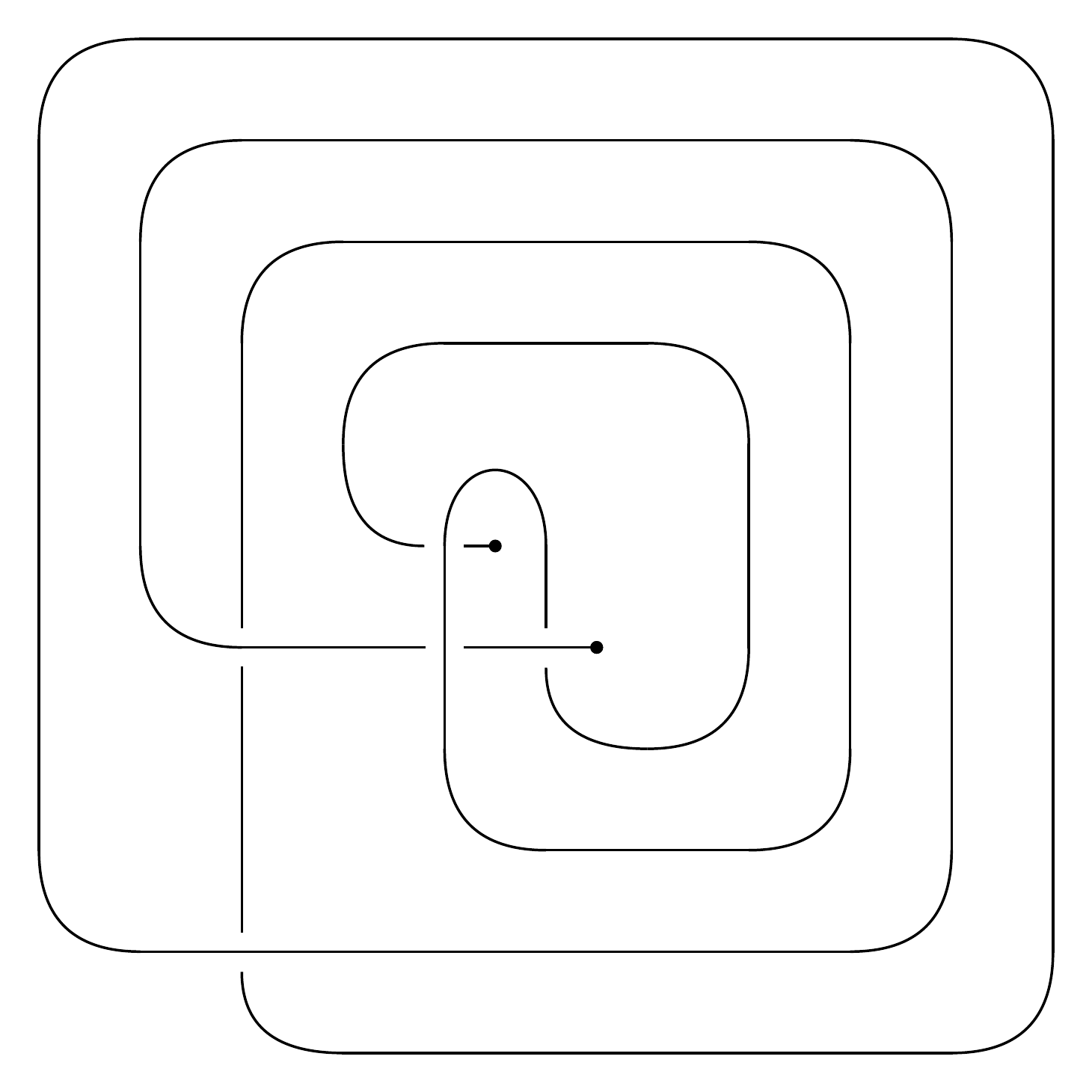}\\
\textcolor{black}{$5_{49}$}
\vspace{1cm}
\end{minipage}
\begin{minipage}[t]{.25\linewidth}
\centering
\includegraphics[width=0.9\textwidth,height=3.5cm,keepaspectratio]{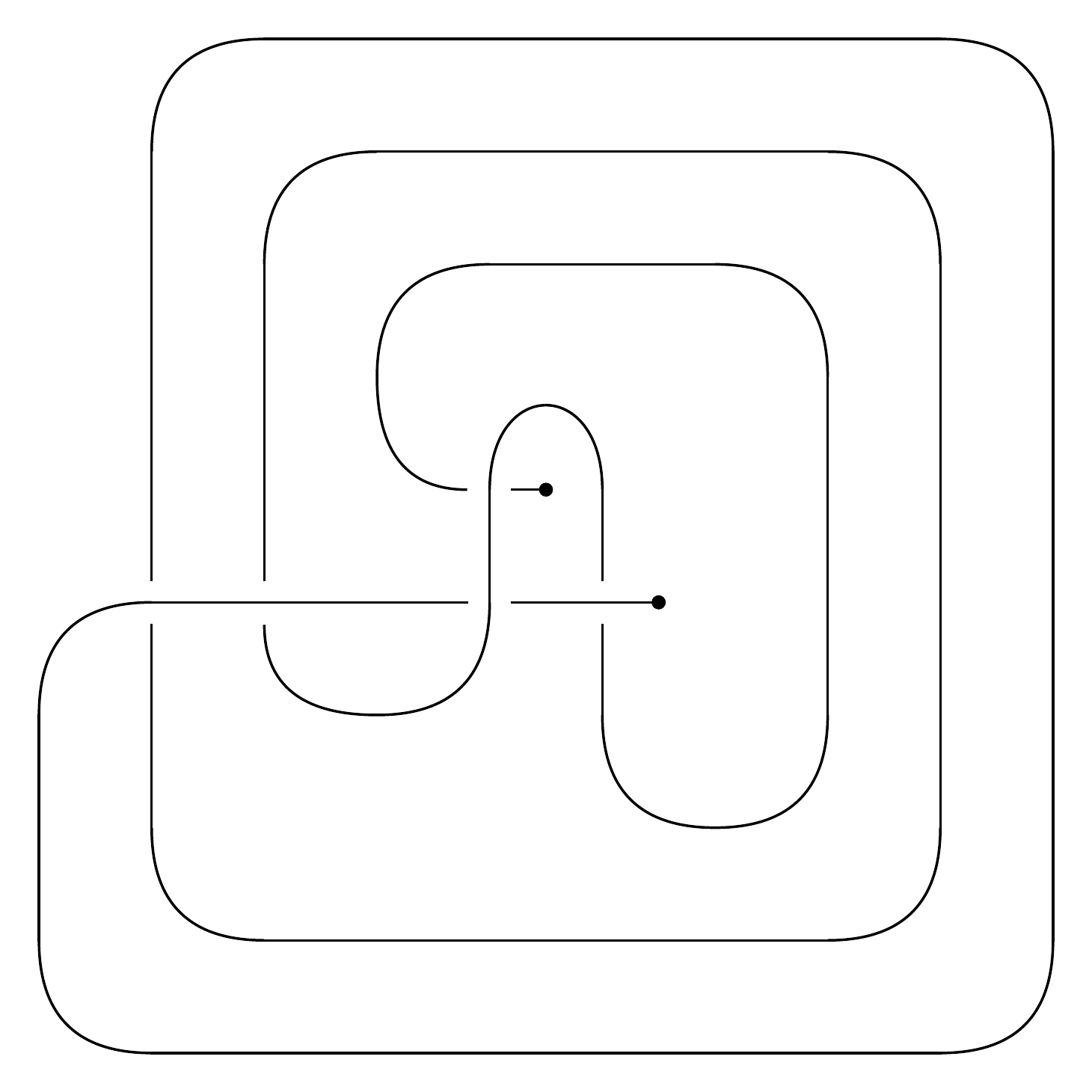}\\
\textcolor{black}{$5_{50}$}
\vspace{1cm}
\end{minipage}
\begin{minipage}[t]{.25\linewidth}
\centering
\includegraphics[width=0.9\textwidth,height=3.5cm,keepaspectratio]{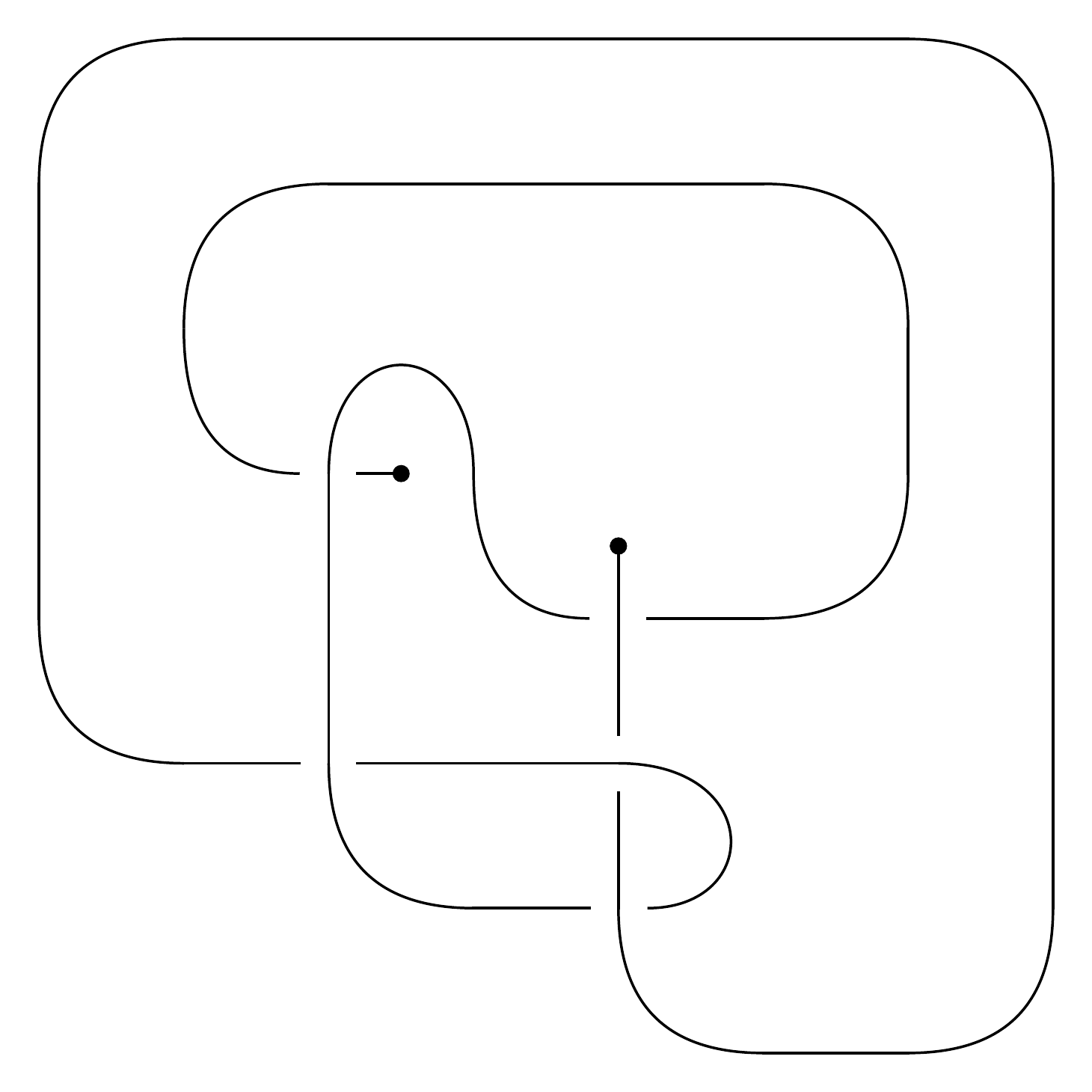}\\
\textcolor{black}{$5_{51}$}
\vspace{1cm}
\end{minipage}
\begin{minipage}[t]{.25\linewidth}
\centering
\includegraphics[width=0.9\textwidth,height=3.5cm,keepaspectratio]{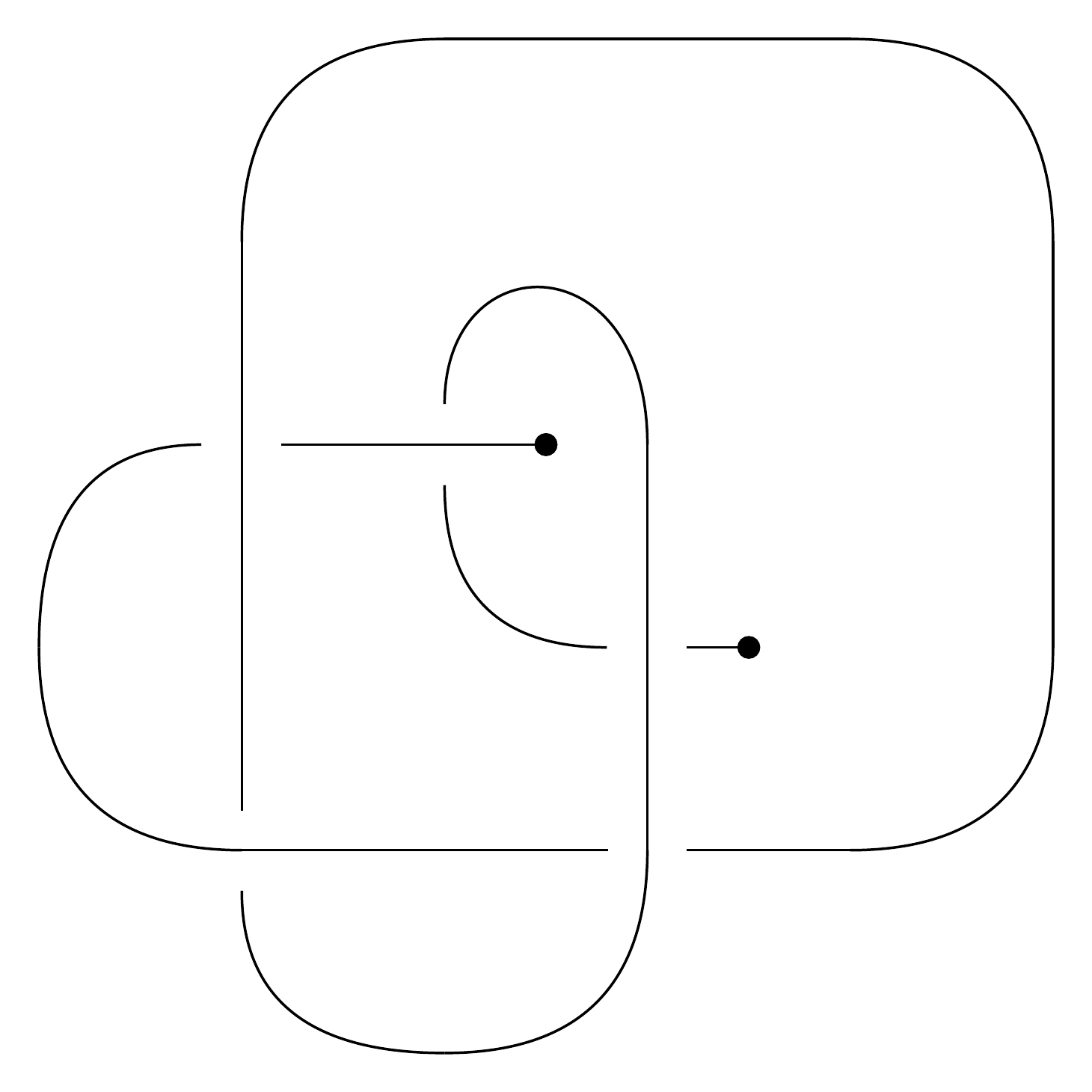}\\
\textcolor{black}{$5_{52}$}
\vspace{1cm}
\end{minipage}
\begin{minipage}[t]{.25\linewidth}
\centering
\includegraphics[width=0.9\textwidth,height=3.5cm,keepaspectratio]{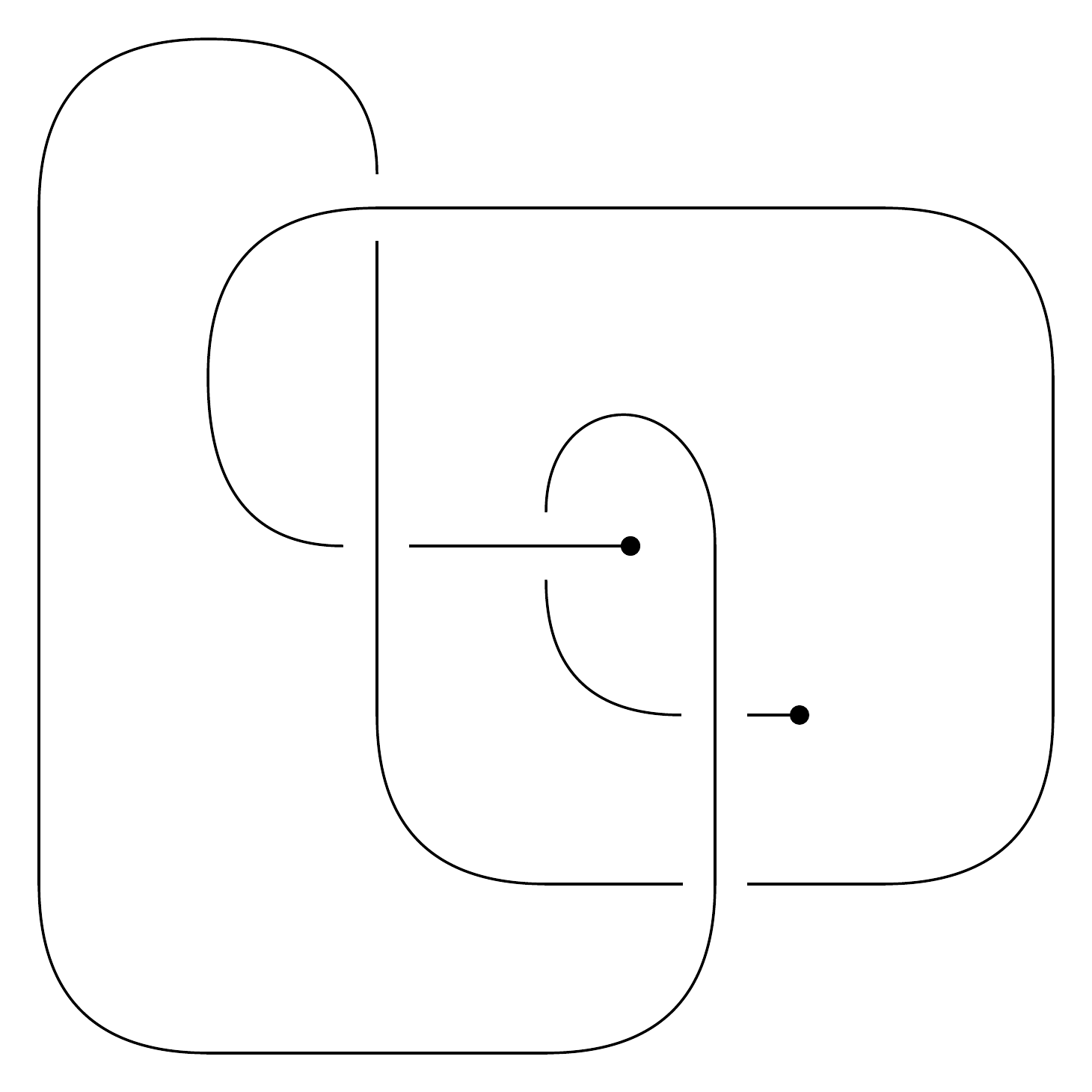}\\
\textcolor{black}{$5_{53}$}
\vspace{1cm}
\end{minipage}
\begin{minipage}[t]{.25\linewidth}
\centering
\includegraphics[width=0.9\textwidth,height=3.5cm,keepaspectratio]{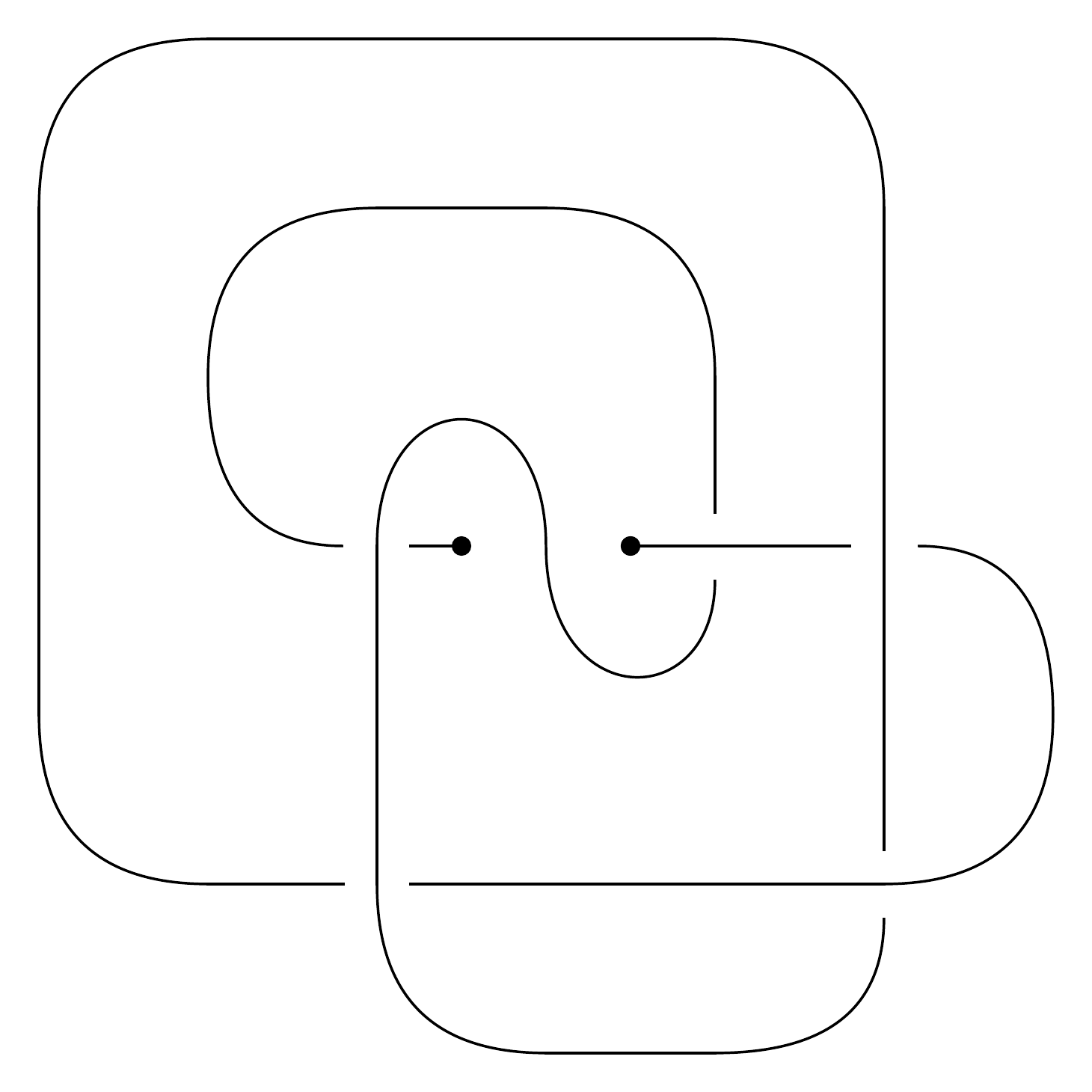}\\
\textcolor{black}{$5_{54}$}
\vspace{1cm}
\end{minipage}
\begin{minipage}[t]{.25\linewidth}
\centering
\includegraphics[width=0.9\textwidth,height=3.5cm,keepaspectratio]{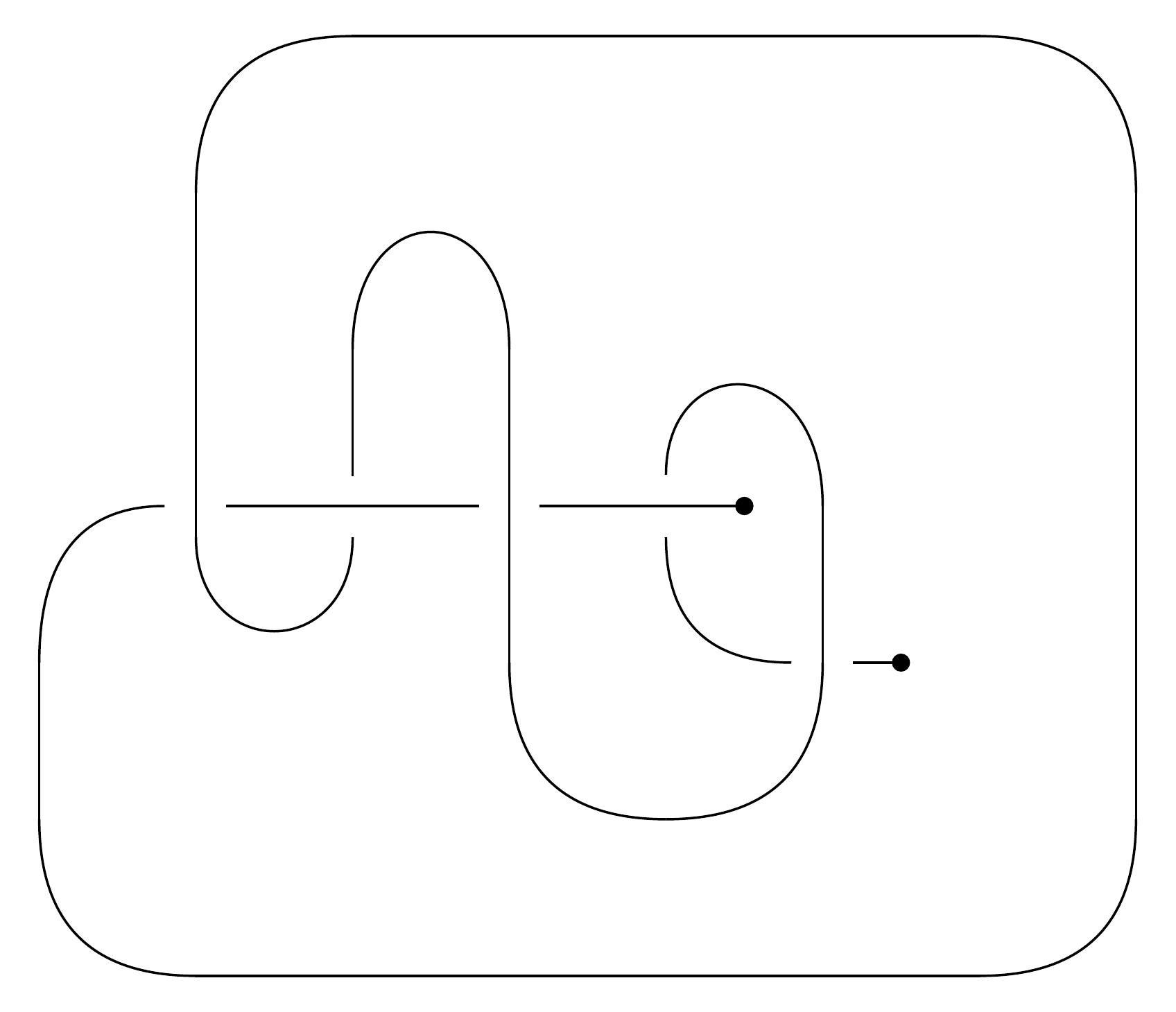}\\
\textcolor{black}{$5_{55}$}
\vspace{1cm}
\end{minipage}
\begin{minipage}[t]{.25\linewidth}
\centering
\includegraphics[width=0.9\textwidth,height=3.5cm,keepaspectratio]{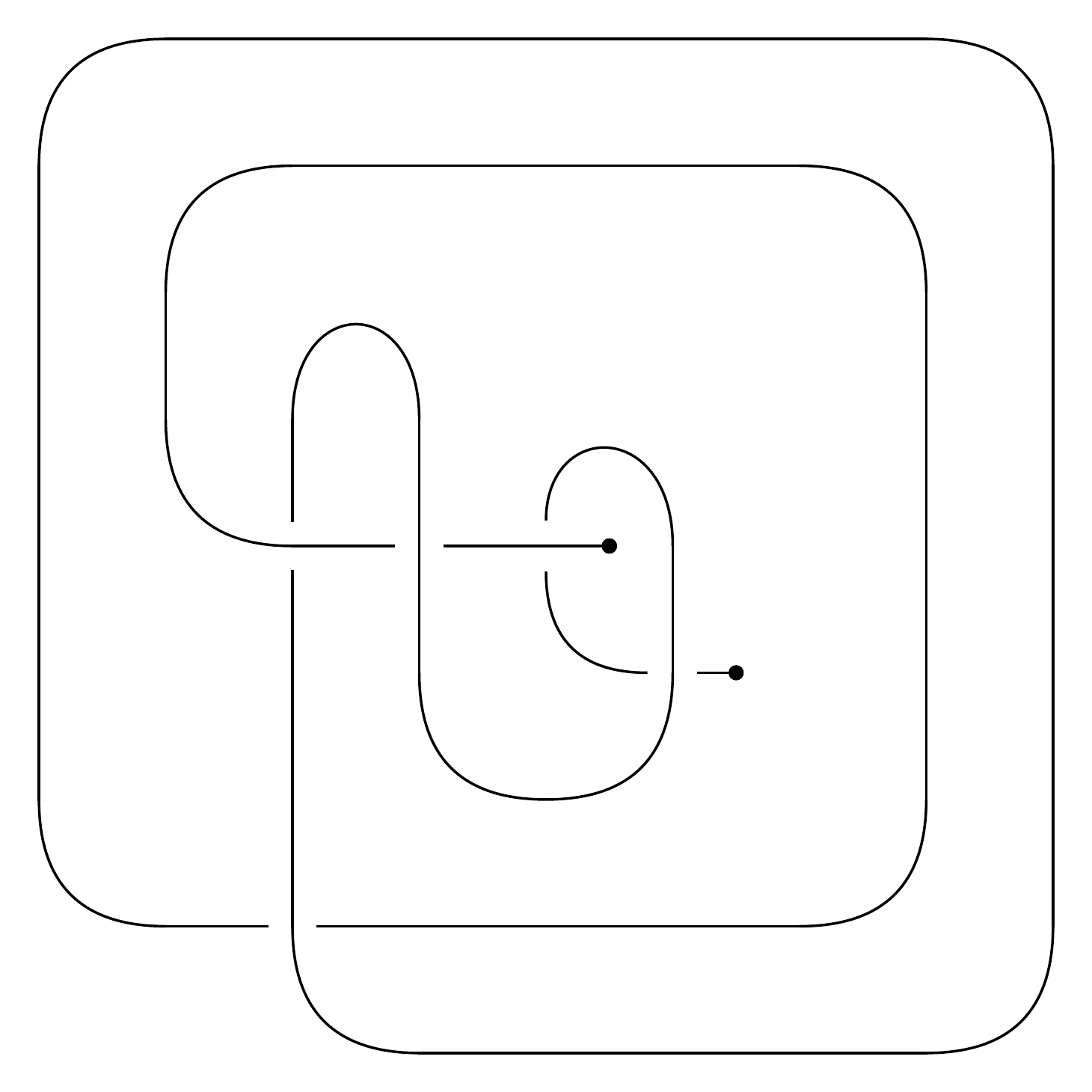}\\
\textcolor{black}{$5_{56}$}
\vspace{1cm}
\end{minipage}
\begin{minipage}[t]{.25\linewidth}
\centering
\includegraphics[width=0.9\textwidth,height=3.5cm,keepaspectratio]{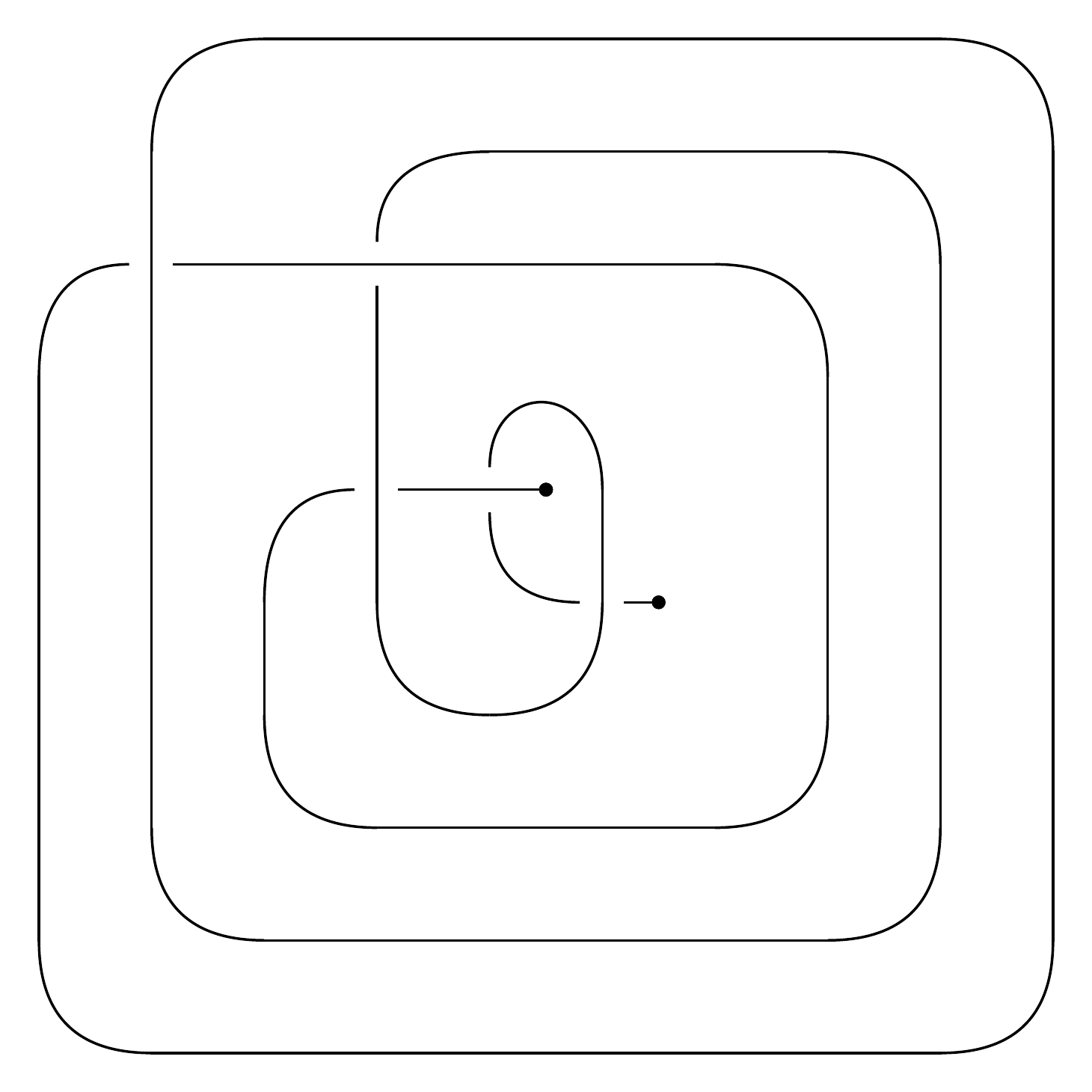}\\
\textcolor{black}{$5_{57}$}
\vspace{1cm}
\end{minipage}
\begin{minipage}[t]{.25\linewidth}
\centering
\includegraphics[width=0.9\textwidth,height=3.5cm,keepaspectratio]{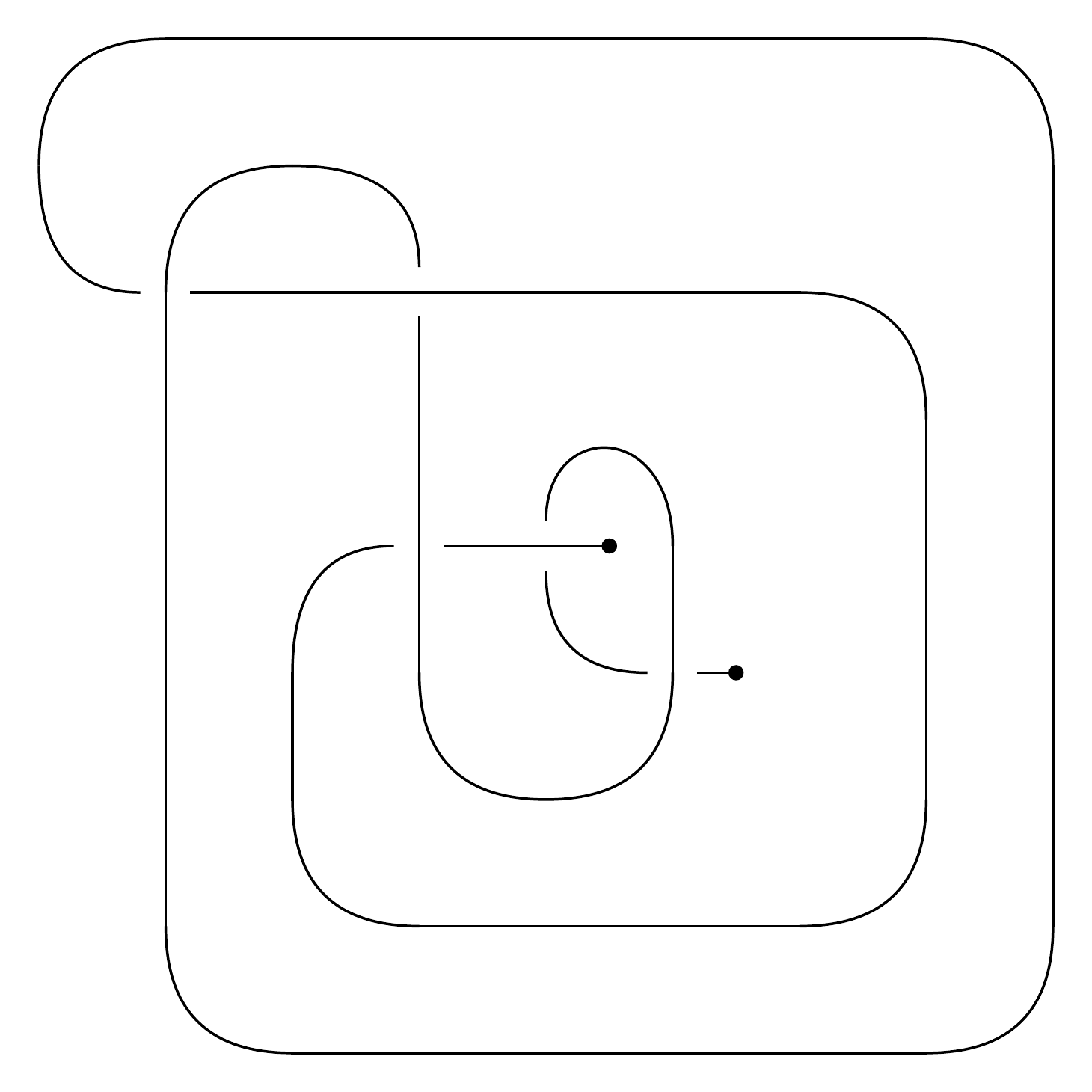}\\
\textcolor{black}{$5_{58}$}
\vspace{1cm}
\end{minipage}
\begin{minipage}[t]{.25\linewidth}
\centering
\includegraphics[width=0.9\textwidth,height=3.5cm,keepaspectratio]{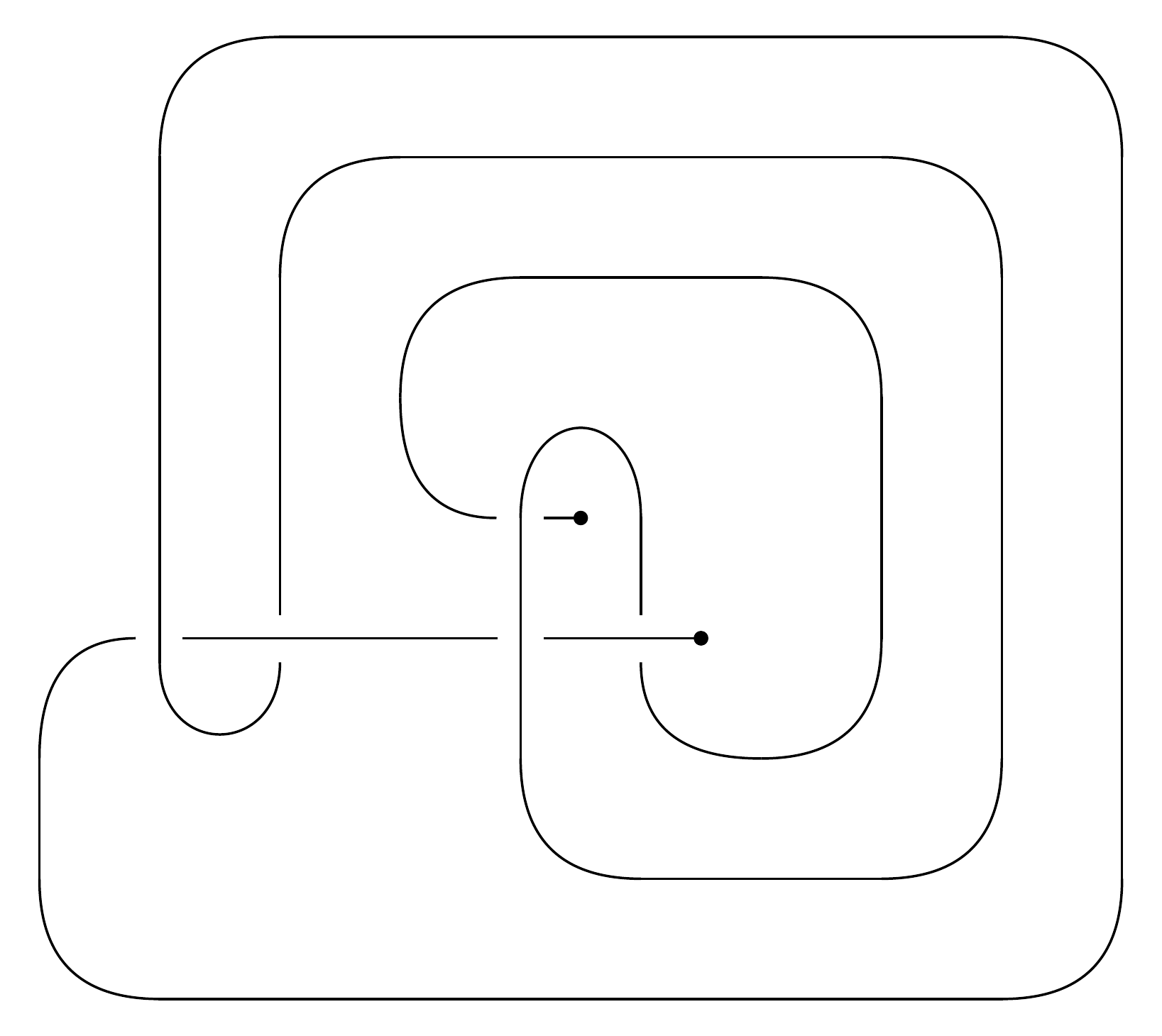}\\
\textcolor{black}{$5_{59}$}
\vspace{1cm}
\end{minipage}
\begin{minipage}[t]{.25\linewidth}
\centering
\includegraphics[width=0.9\textwidth,height=3.5cm,keepaspectratio]{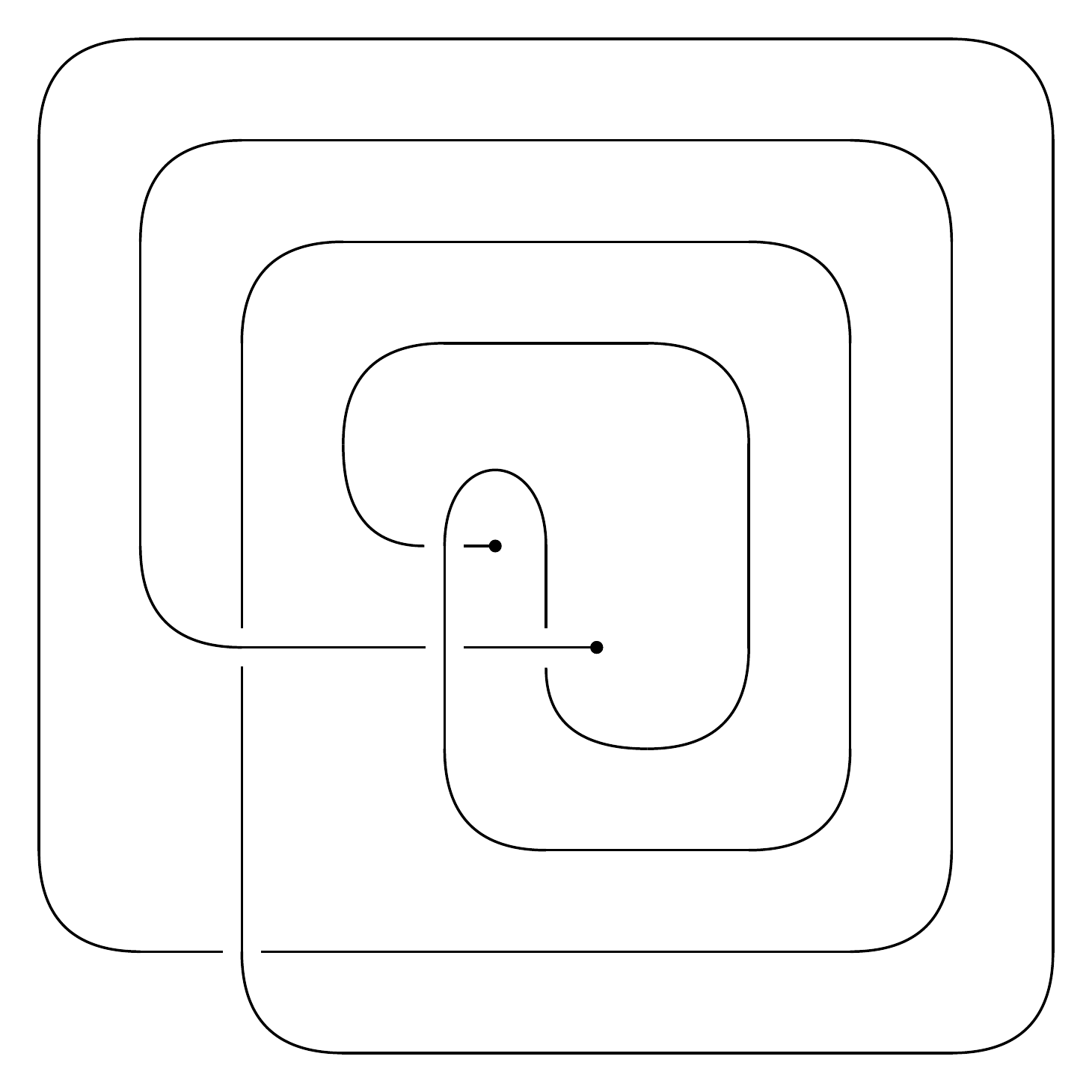}\\
\textcolor{black}{$5_{60}$}
\vspace{1cm}
\end{minipage}
\begin{minipage}[t]{.25\linewidth}
\centering
\includegraphics[width=0.9\textwidth,height=3.5cm,keepaspectratio]{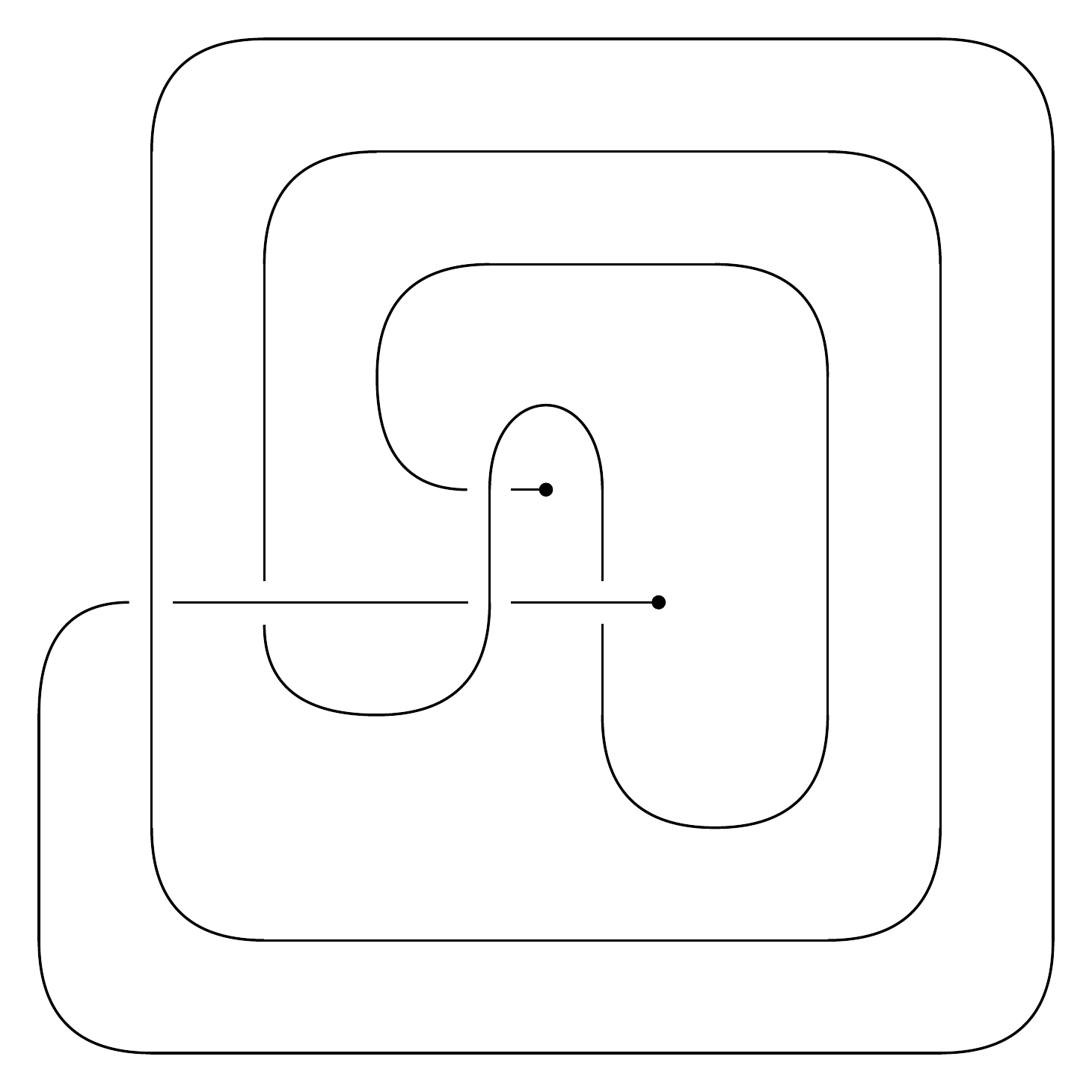}\\
\textcolor{black}{$5_{61}$}
\vspace{1cm}
\end{minipage}
\begin{minipage}[t]{.25\linewidth}
\centering
\includegraphics[width=0.9\textwidth,height=3.5cm,keepaspectratio]{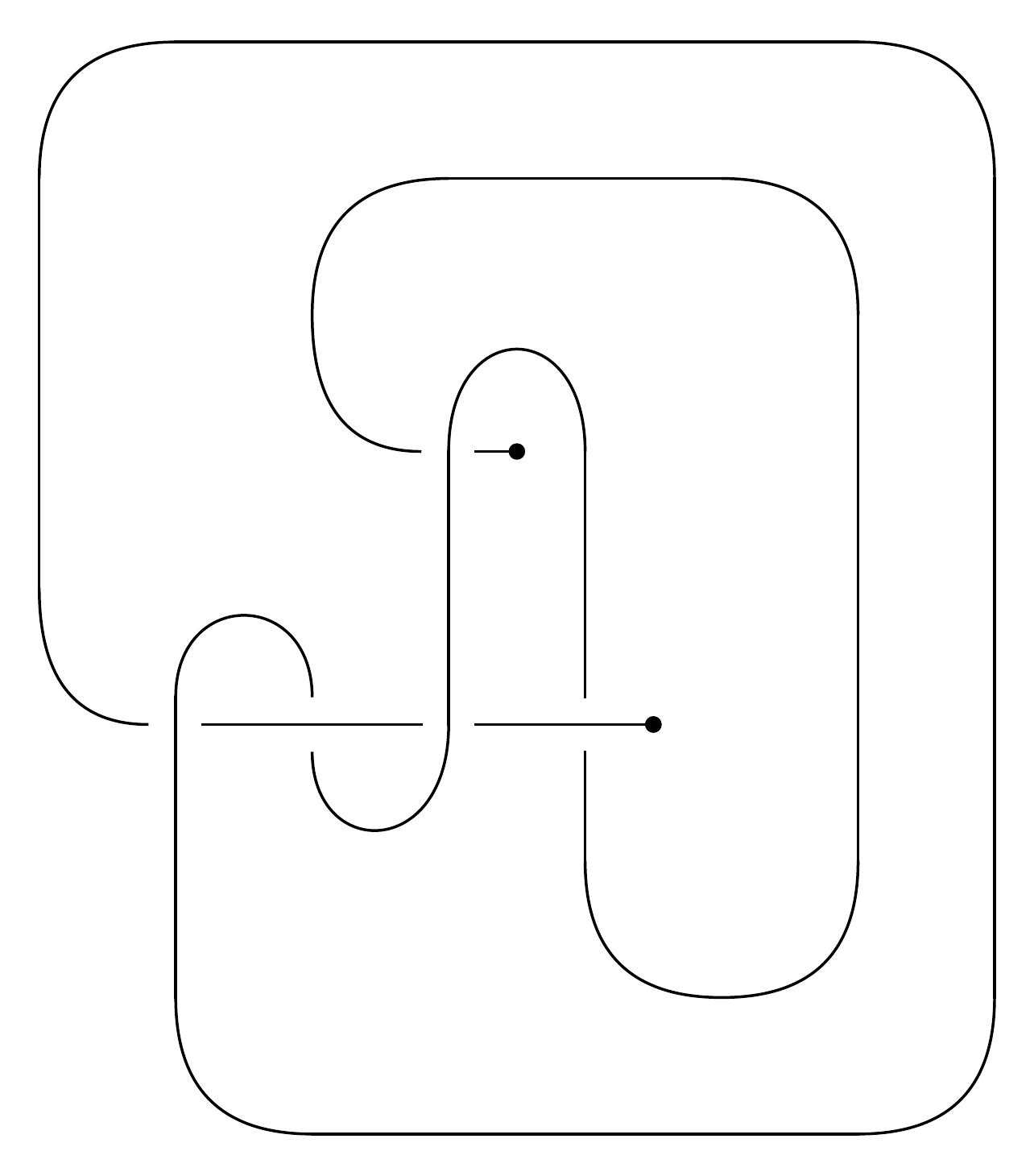}\\
\textcolor{black}{$5_{62}$}
\vspace{1cm}
\end{minipage}
\begin{minipage}[t]{.25\linewidth}
\centering
\includegraphics[width=0.9\textwidth,height=3.5cm,keepaspectratio]{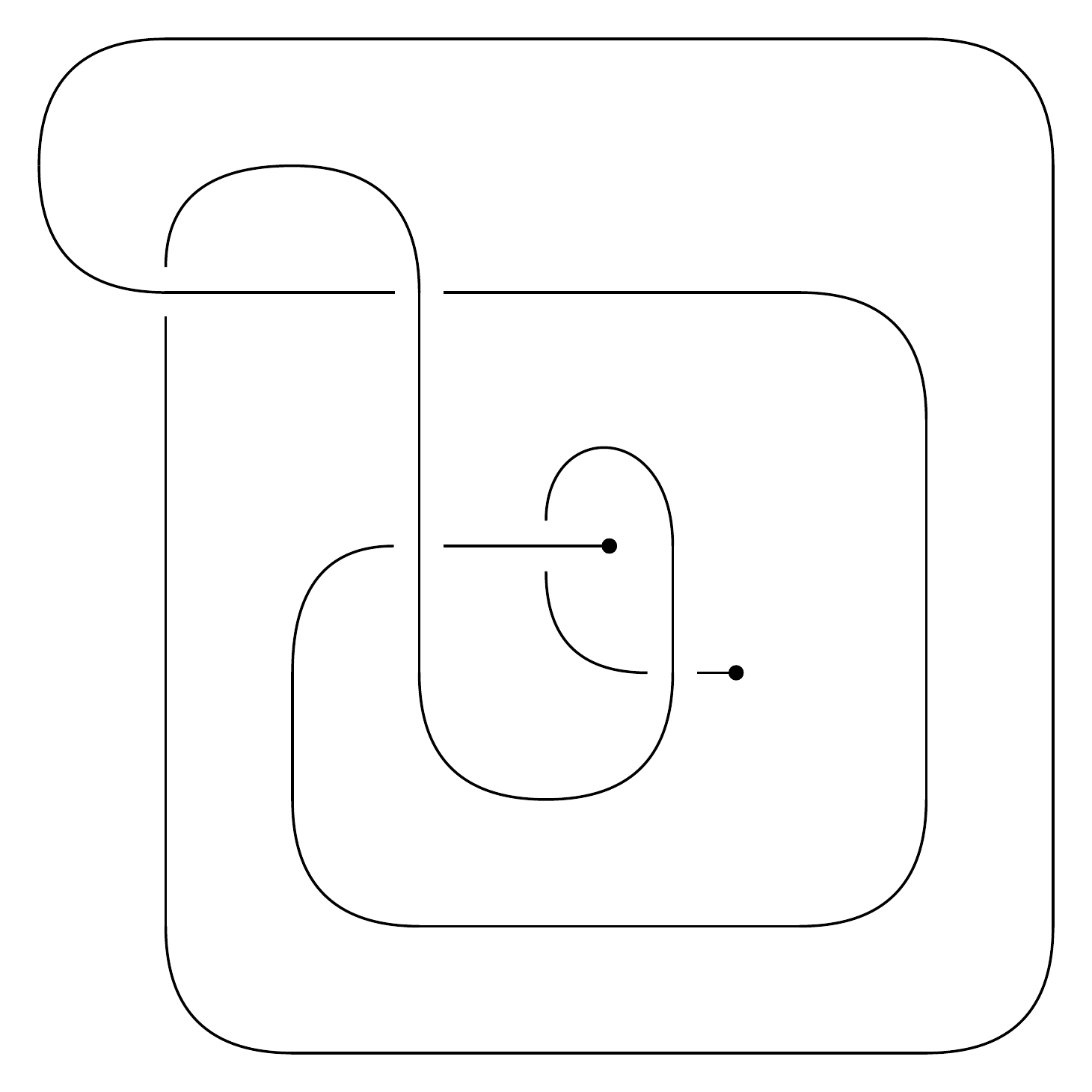}\\
\textcolor{black}{$5_{63}$}
\vspace{1cm}
\end{minipage}
\begin{minipage}[t]{.25\linewidth}
\centering
\includegraphics[width=0.9\textwidth,height=3.5cm,keepaspectratio]{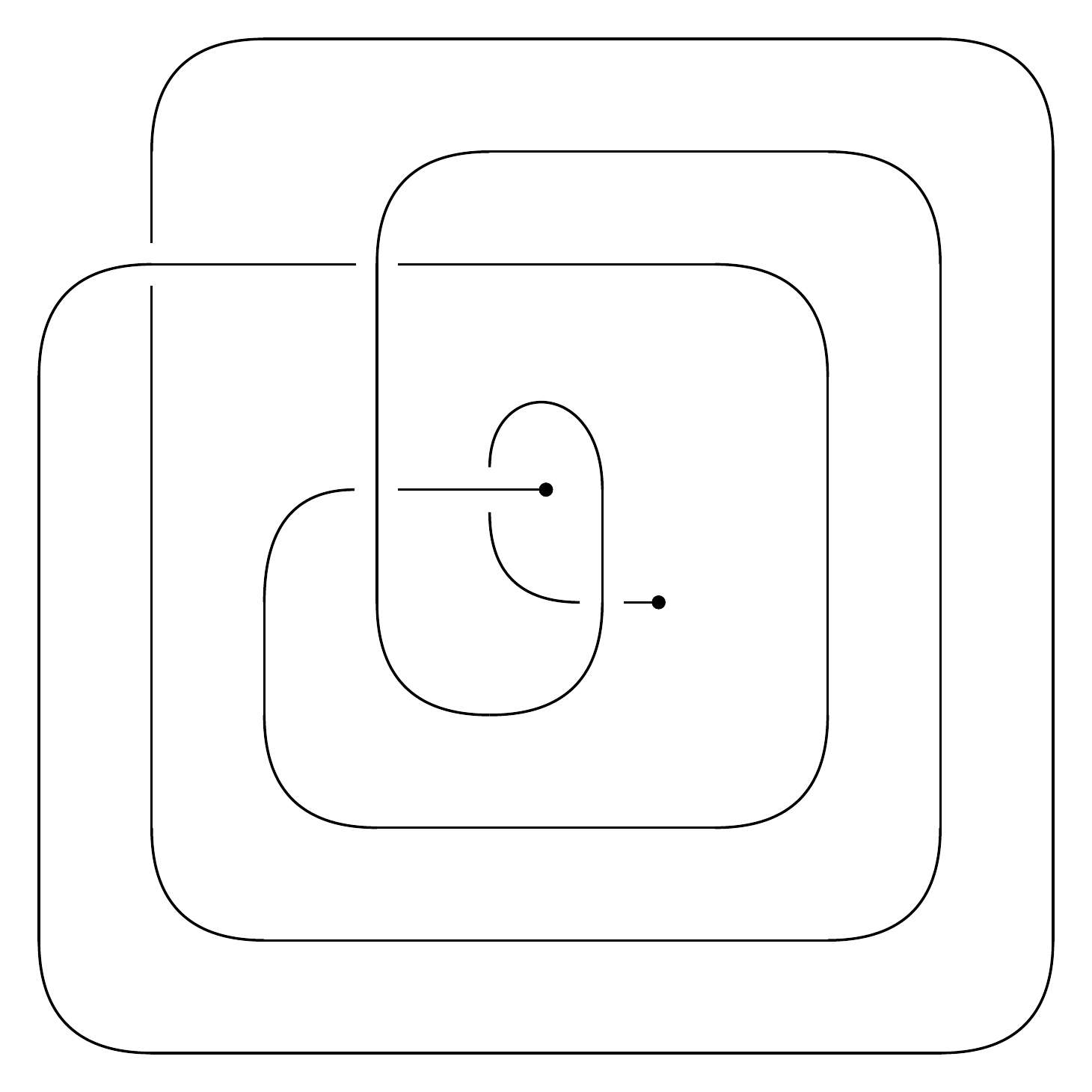}\\
\textcolor{black}{$5_{64}$}
\vspace{1cm}
\end{minipage}
\begin{minipage}[t]{.25\linewidth}
\centering
\includegraphics[width=0.9\textwidth,height=3.5cm,keepaspectratio]{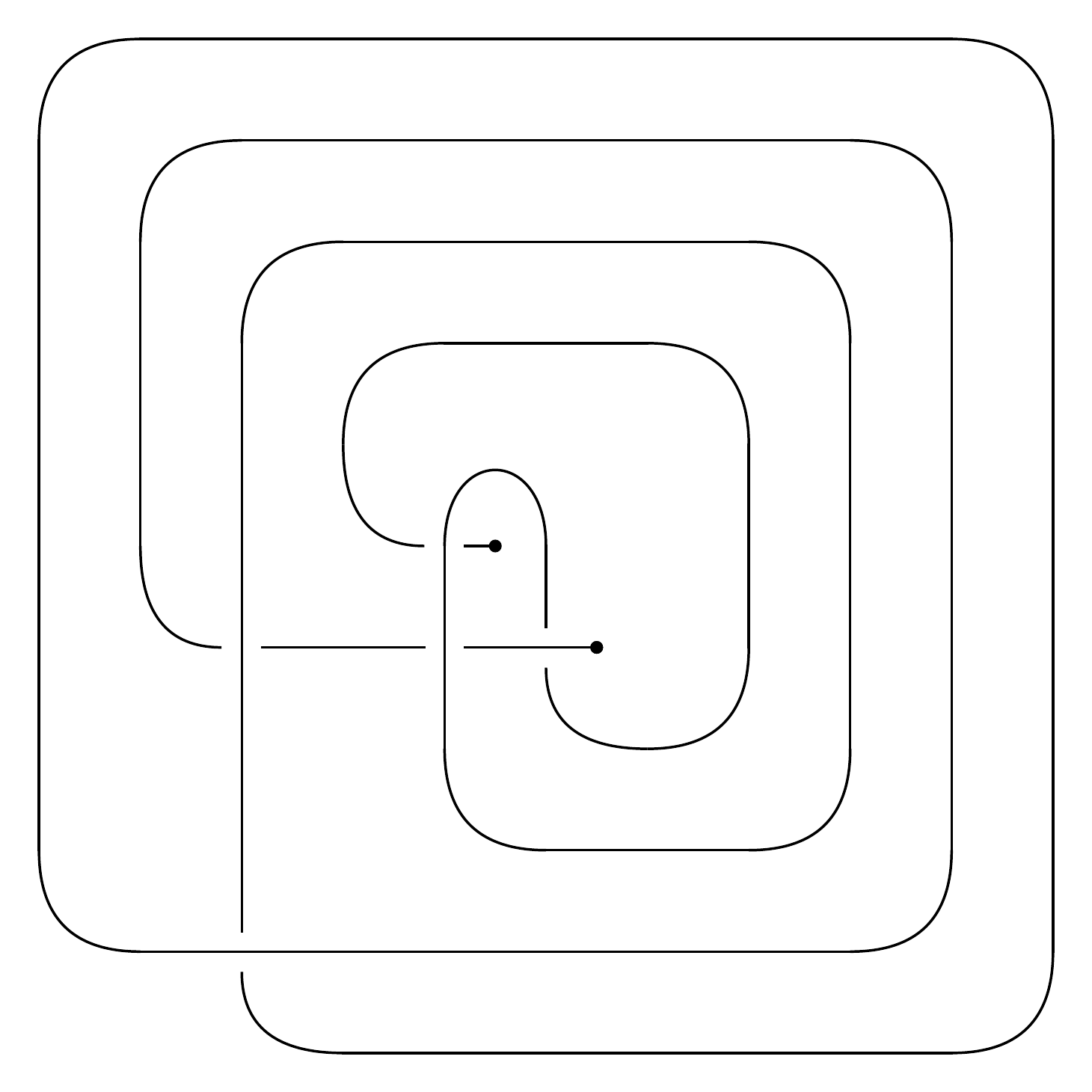}\\
\textcolor{black}{$5_{65}$}
\vspace{1cm}
\end{minipage}
\begin{minipage}[t]{.25\linewidth}
\centering
\includegraphics[width=0.9\textwidth,height=3.5cm,keepaspectratio]{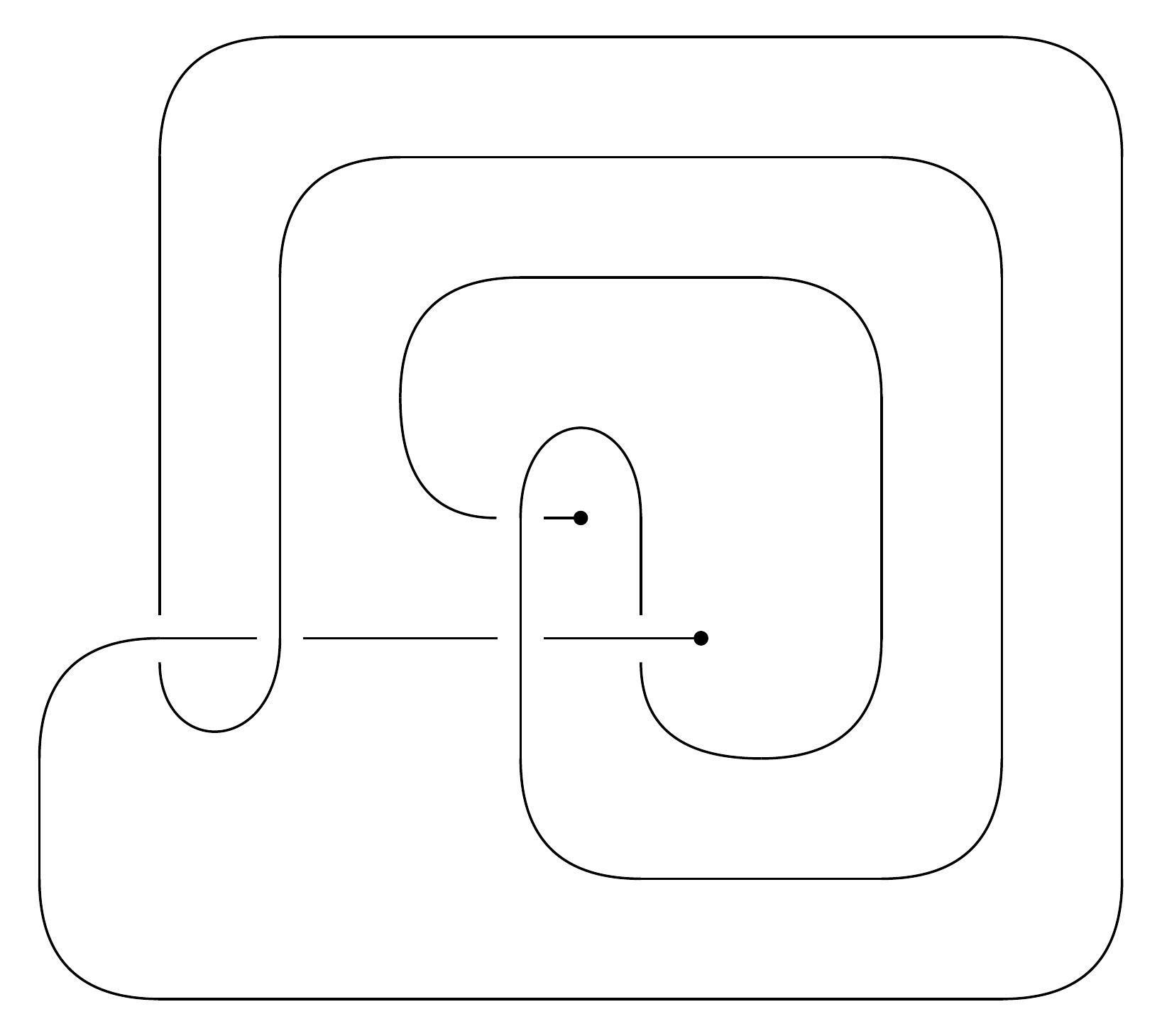}\\
\textcolor{black}{$5_{66}$}
\vspace{1cm}
\end{minipage}
\begin{minipage}[t]{.25\linewidth}
\centering
\includegraphics[width=0.9\textwidth,height=3.5cm,keepaspectratio]{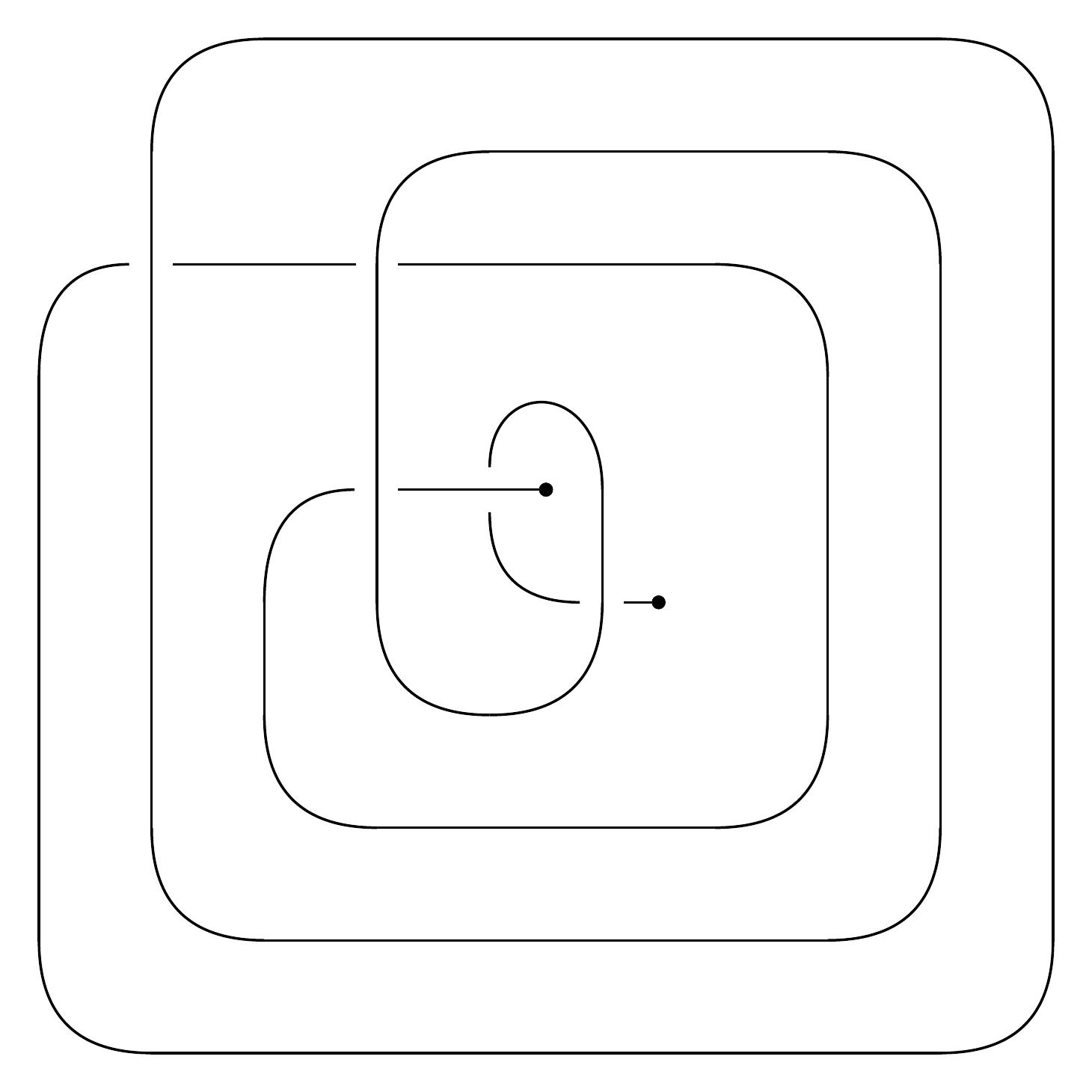}\\
\textcolor{black}{$5_{67}$}
\vspace{1cm}
\end{minipage}
\begin{minipage}[t]{.25\linewidth}
\centering
\includegraphics[width=0.9\textwidth,height=3.5cm,keepaspectratio]{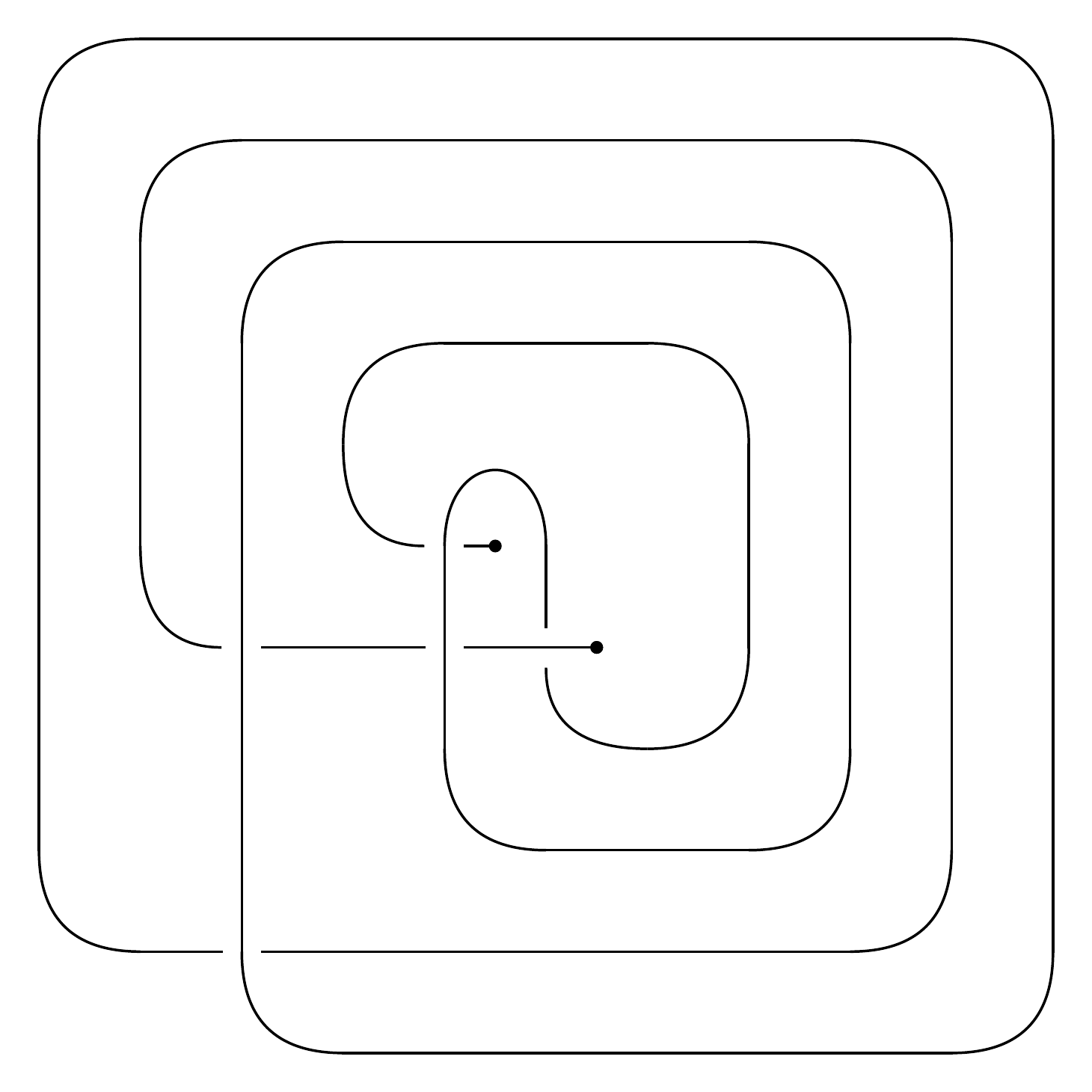}\\
\textcolor{black}{$5_{68}$}
\vspace{1cm}
\end{minipage}
\begin{minipage}[t]{.25\linewidth}
\centering
\includegraphics[width=0.9\textwidth,height=3.5cm,keepaspectratio]{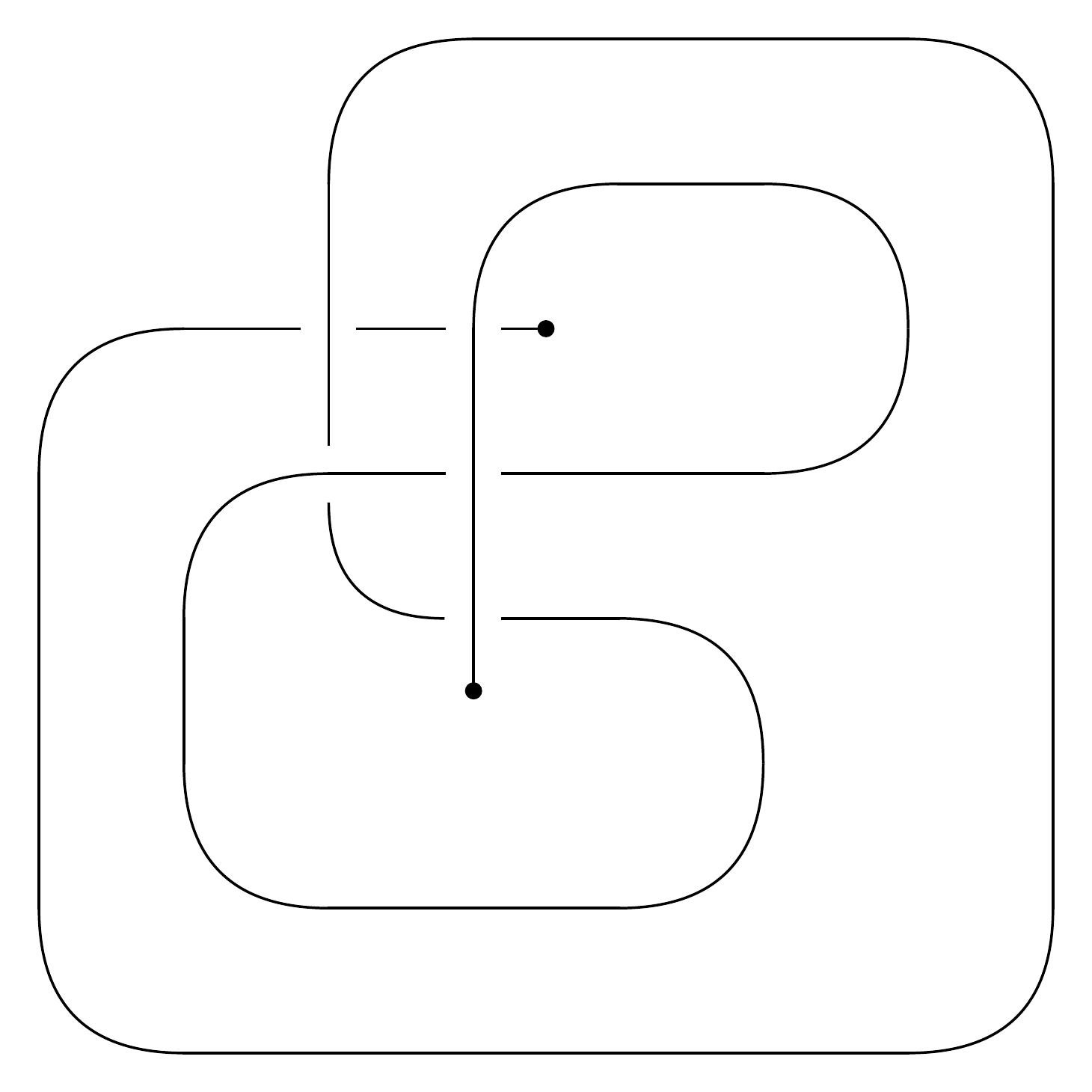}\\
\textcolor{black}{$5_{69}$}
\vspace{1cm}
\end{minipage}
\begin{minipage}[t]{.25\linewidth}
\centering
\includegraphics[width=0.9\textwidth,height=3.5cm,keepaspectratio]{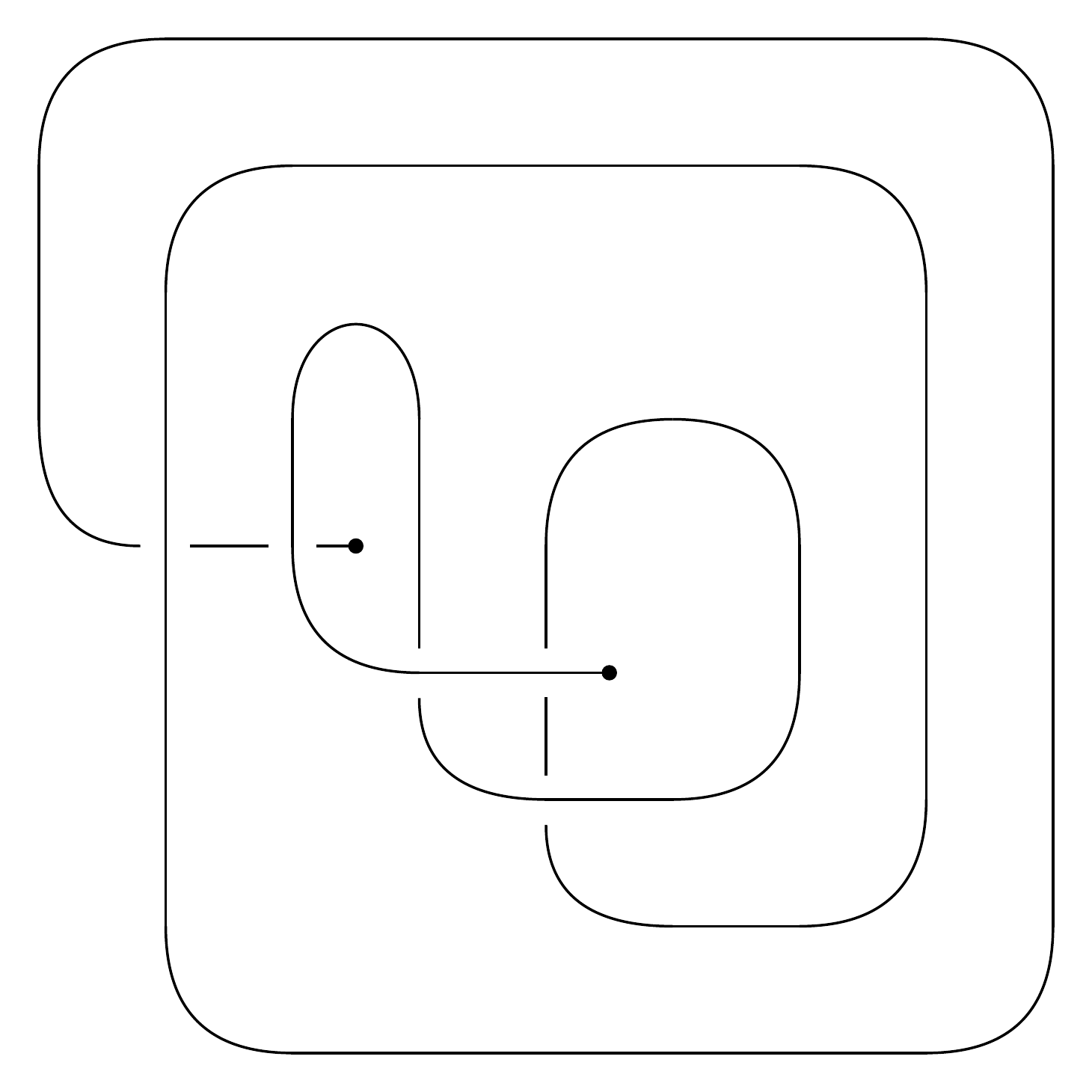}\\
\textcolor{black}{$5_{70}$}
\vspace{1cm}
\end{minipage}
\begin{minipage}[t]{.25\linewidth}
\centering
\includegraphics[width=0.9\textwidth,height=3.5cm,keepaspectratio]{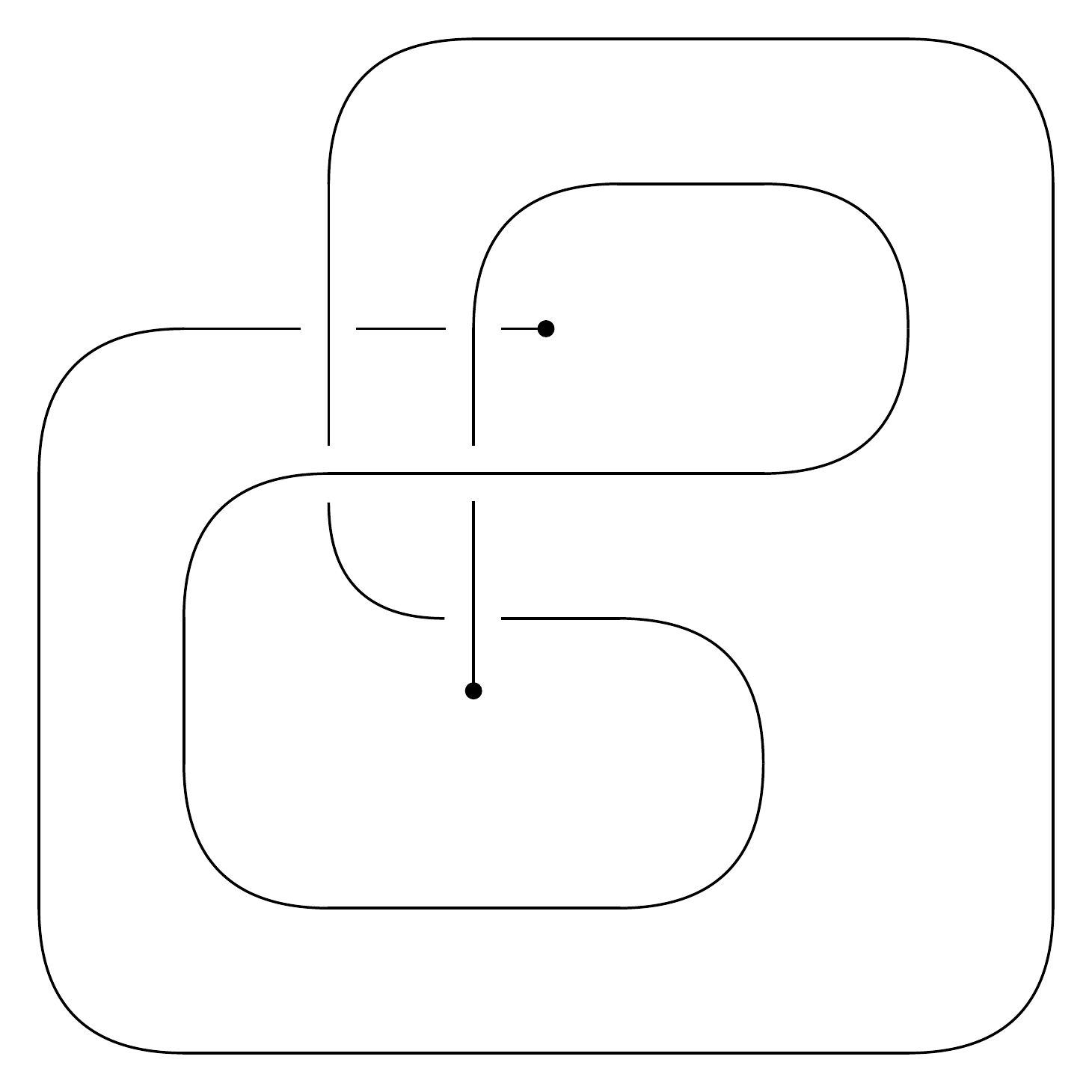}\\
\textcolor{black}{$5_{71}$}
\vspace{1cm}
\end{minipage}
\begin{minipage}[t]{.25\linewidth}
\centering
\includegraphics[width=0.9\textwidth,height=3.5cm,keepaspectratio]{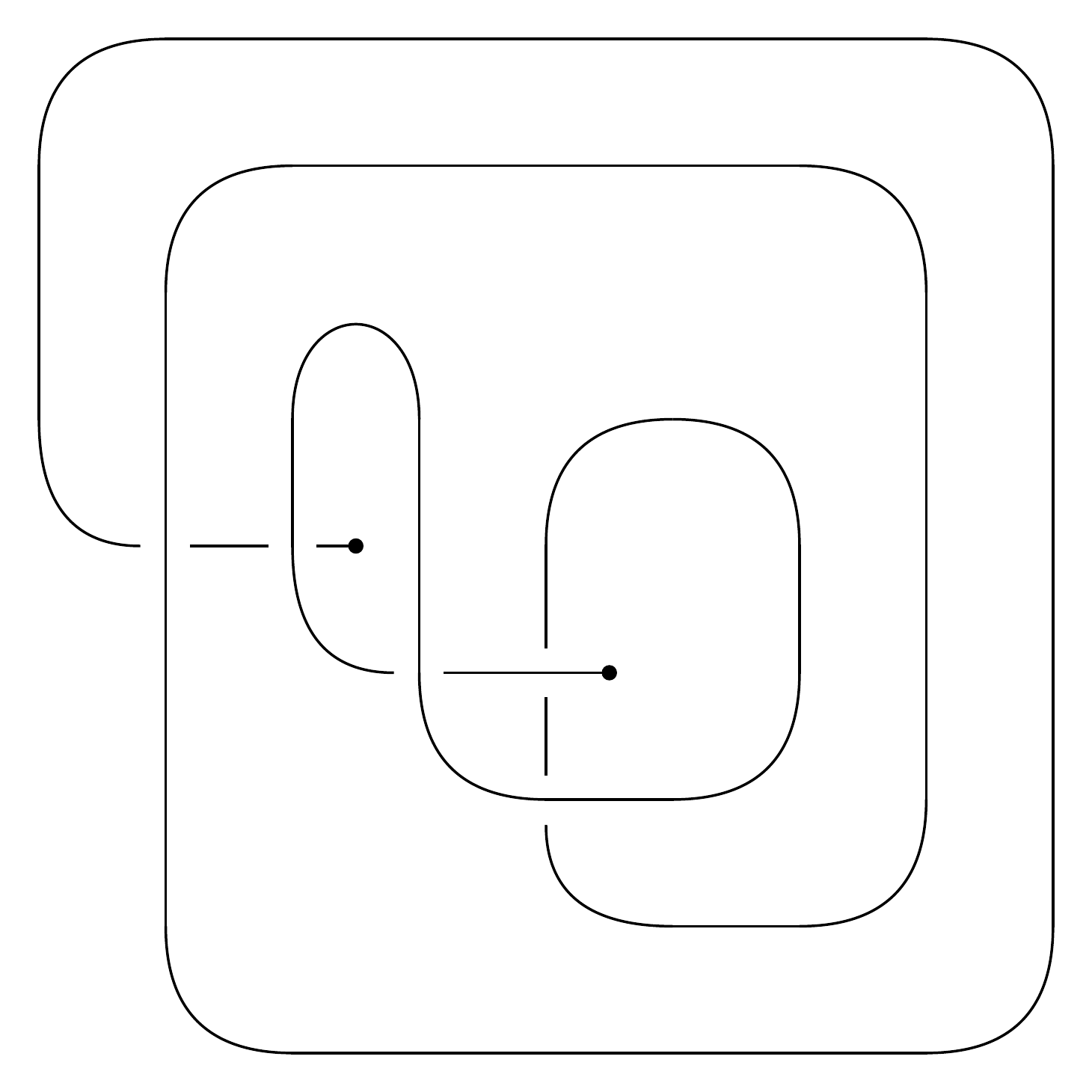}\\
\textcolor{black}{$5_{72}$}
\vspace{1cm}
\end{minipage}
\begin{minipage}[t]{.25\linewidth}
\centering
\includegraphics[width=0.9\textwidth,height=3.5cm,keepaspectratio]{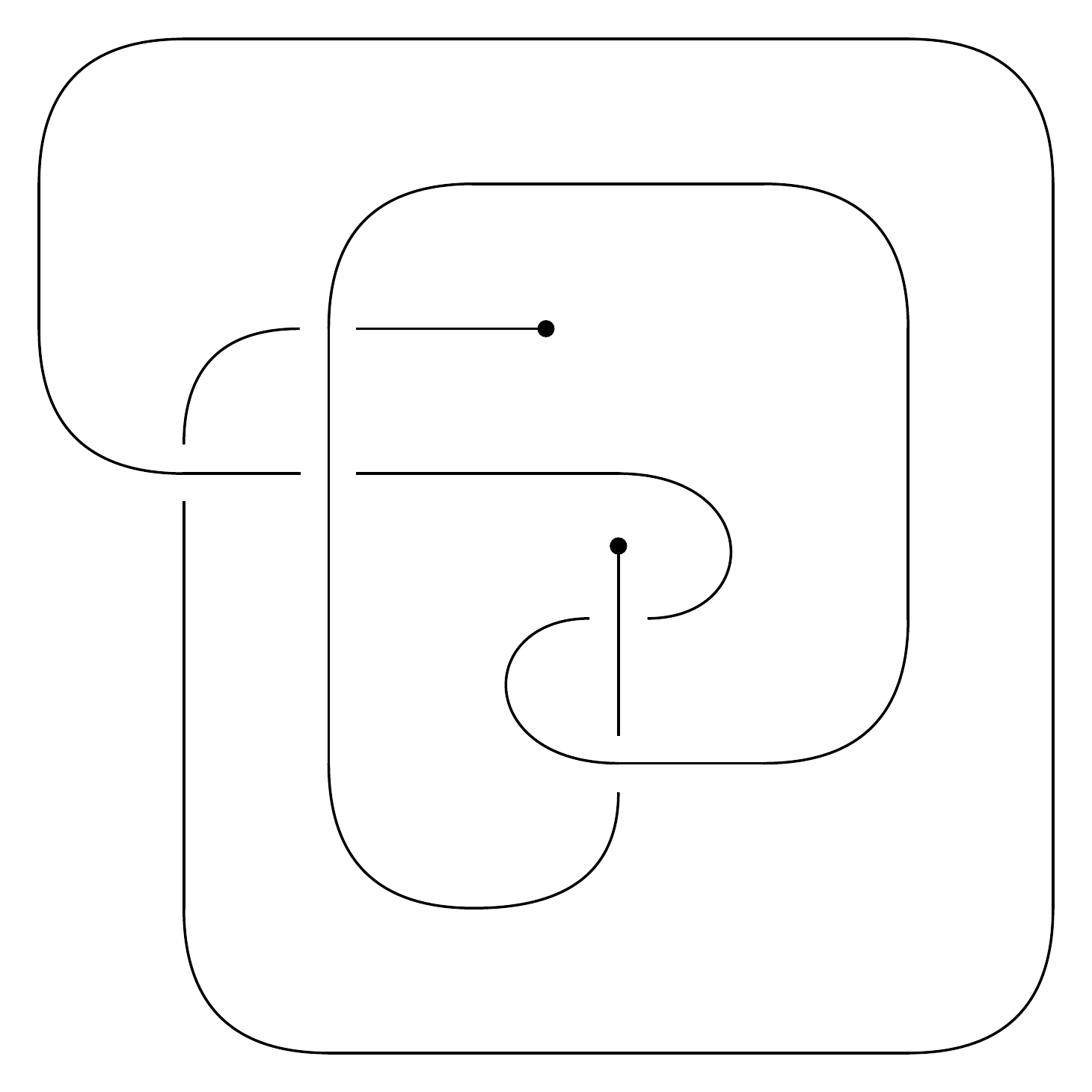}\\
\textcolor{black}{$5_{73}$}
\vspace{1cm}
\end{minipage}
\begin{minipage}[t]{.25\linewidth}
\centering
\includegraphics[width=0.9\textwidth,height=3.5cm,keepaspectratio]{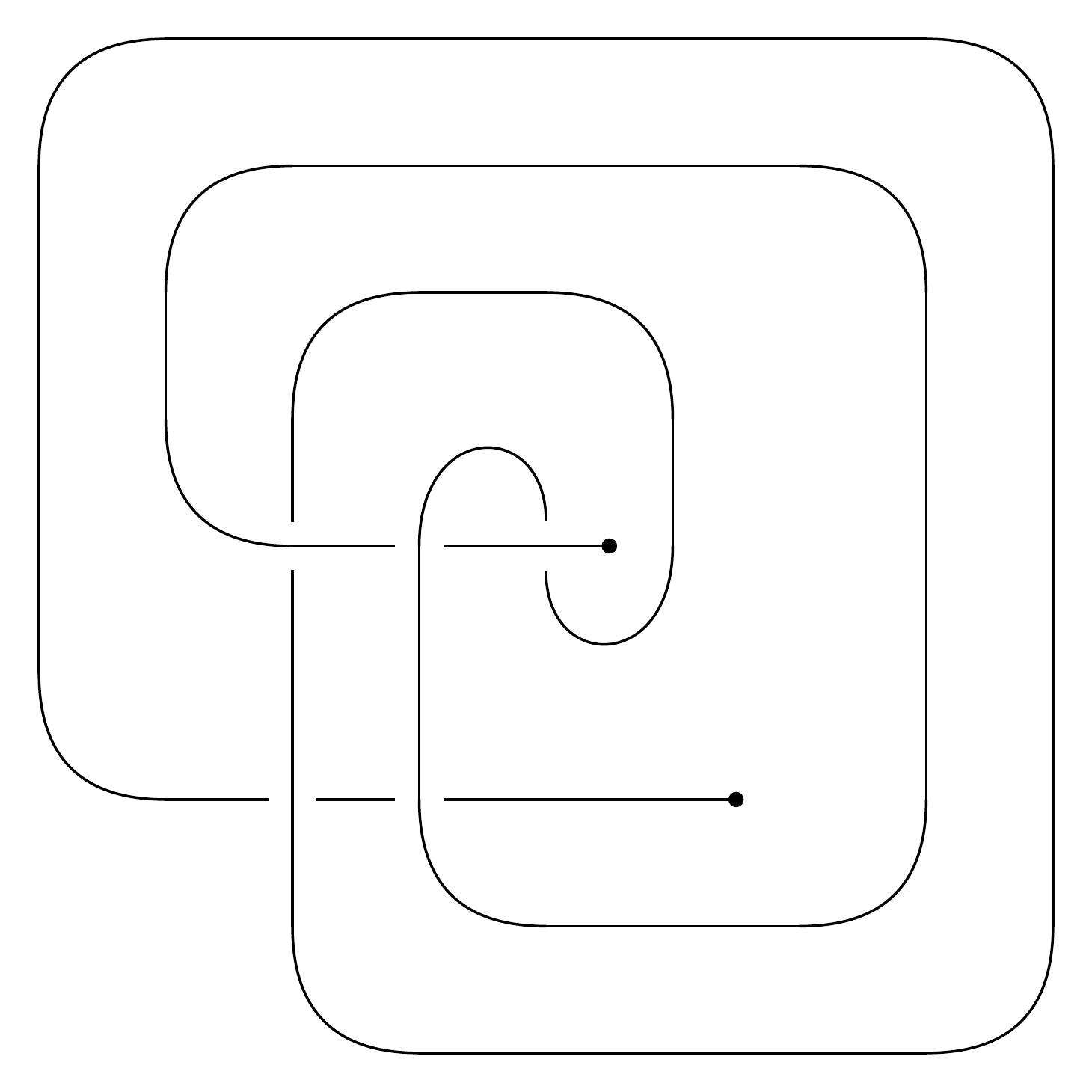}\\
\textcolor{black}{$5_{74}$}
\vspace{1cm}
\end{minipage}
\begin{minipage}[t]{.25\linewidth}
\centering
\includegraphics[width=0.9\textwidth,height=3.5cm,keepaspectratio]{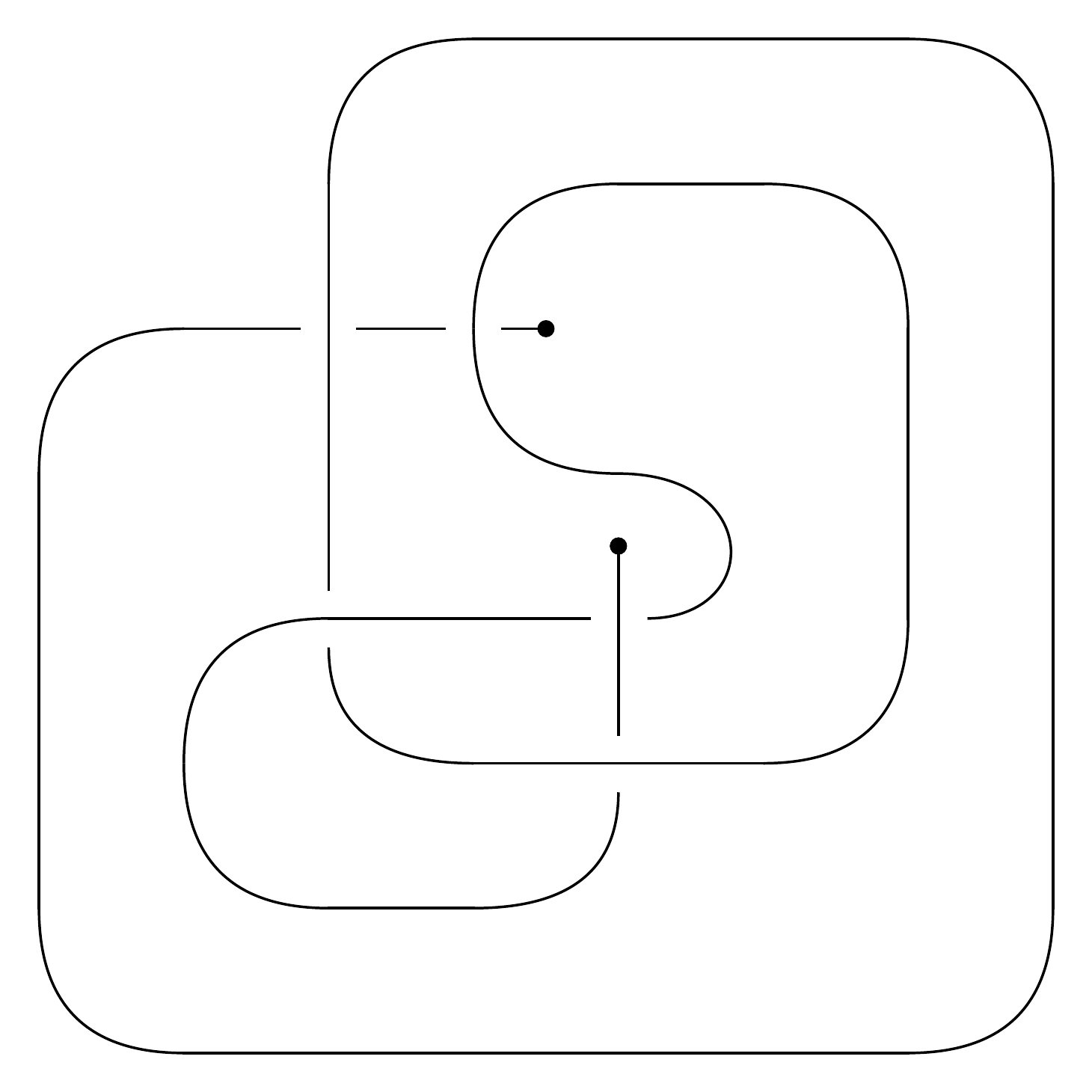}\\
\textcolor{black}{$5_{75}$}
\vspace{1cm}
\end{minipage}
\begin{minipage}[t]{.25\linewidth}
\centering
\includegraphics[width=0.9\textwidth,height=3.5cm,keepaspectratio]{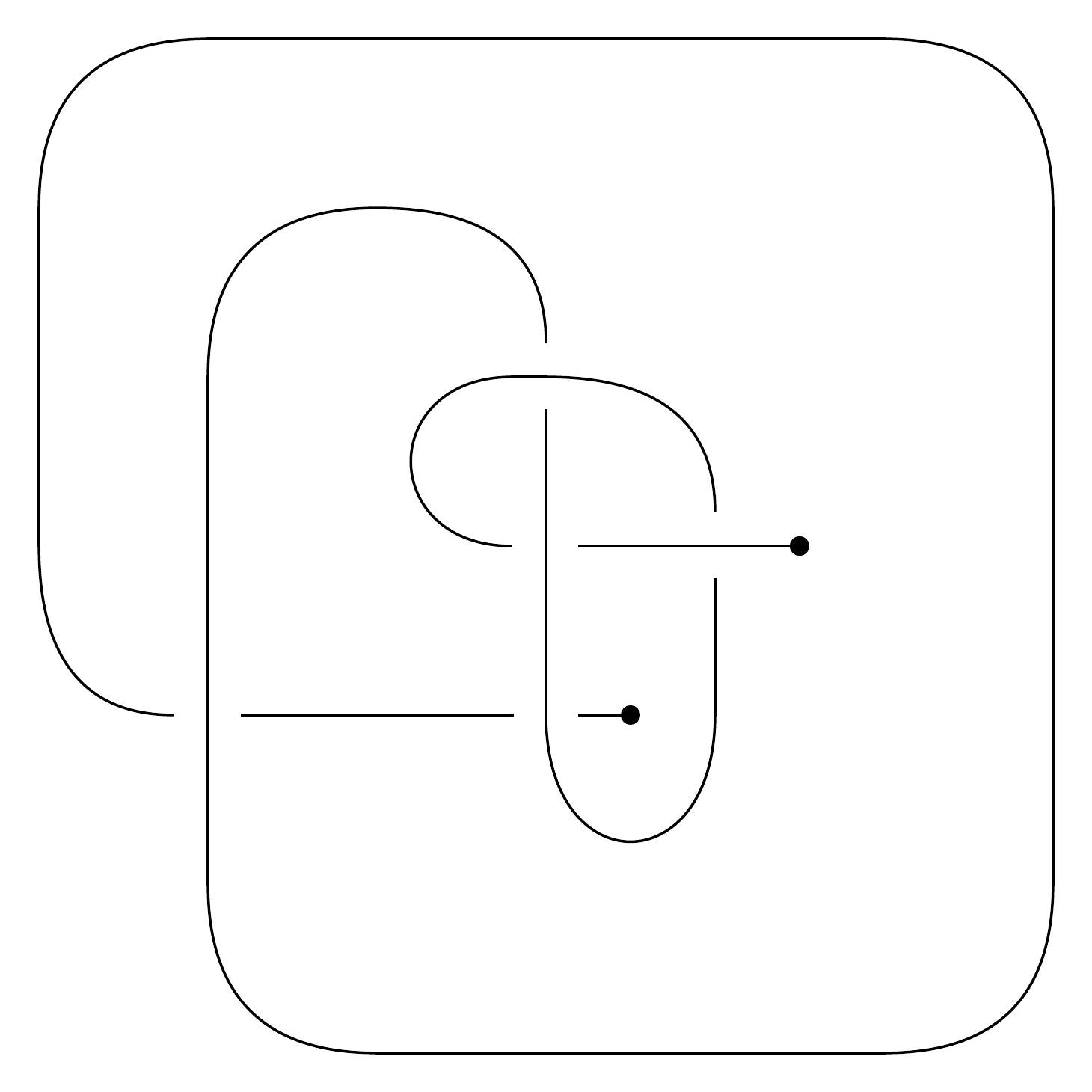}\\
\textcolor{black}{$5_{76}$}
\vspace{1cm}
\end{minipage}
\begin{minipage}[t]{.25\linewidth}
\centering
\includegraphics[width=0.9\textwidth,height=3.5cm,keepaspectratio]{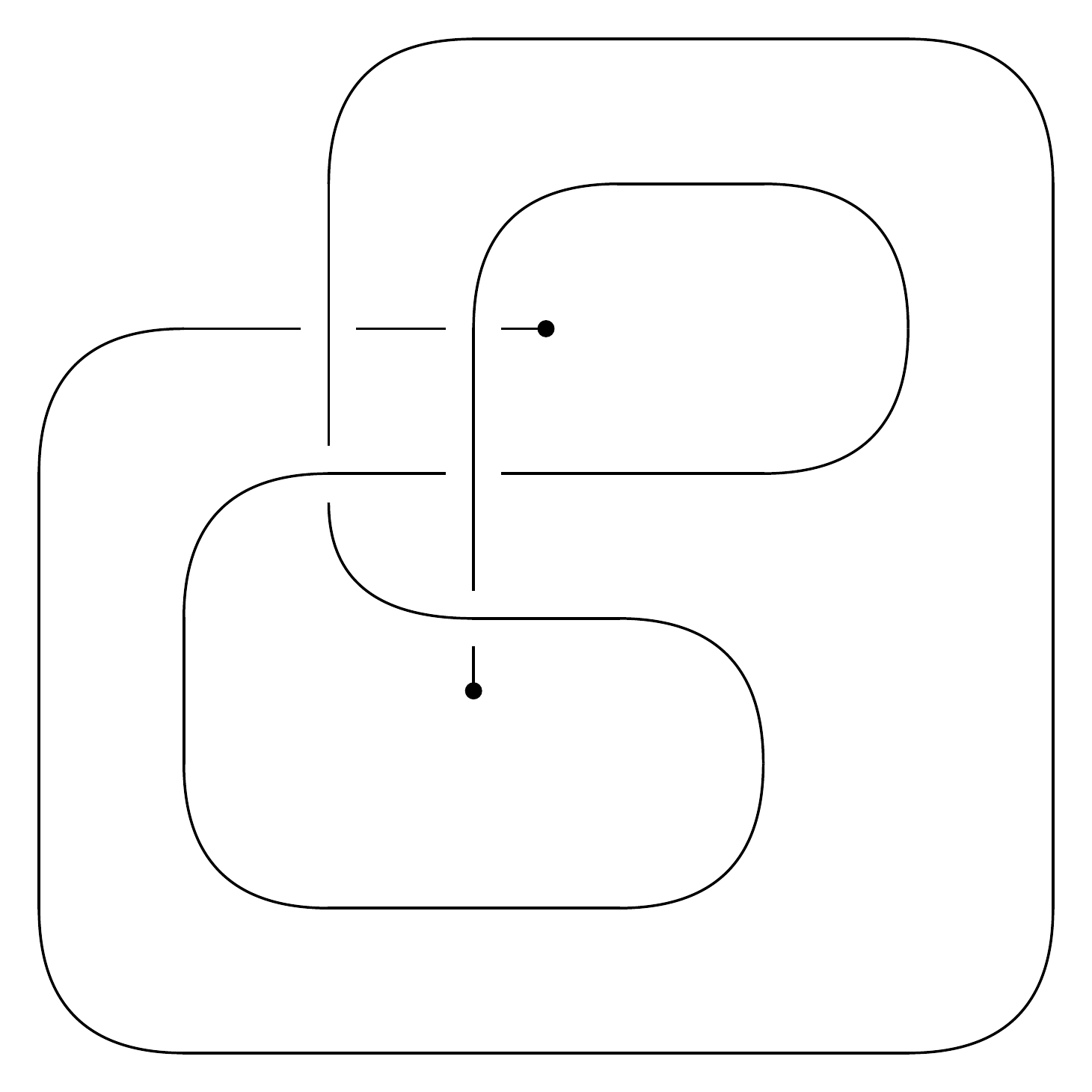}\\
\textcolor{black}{$5_{77}$}
\vspace{1cm}
\end{minipage}
\begin{minipage}[t]{.25\linewidth}
\centering
\includegraphics[width=0.9\textwidth,height=3.5cm,keepaspectratio]{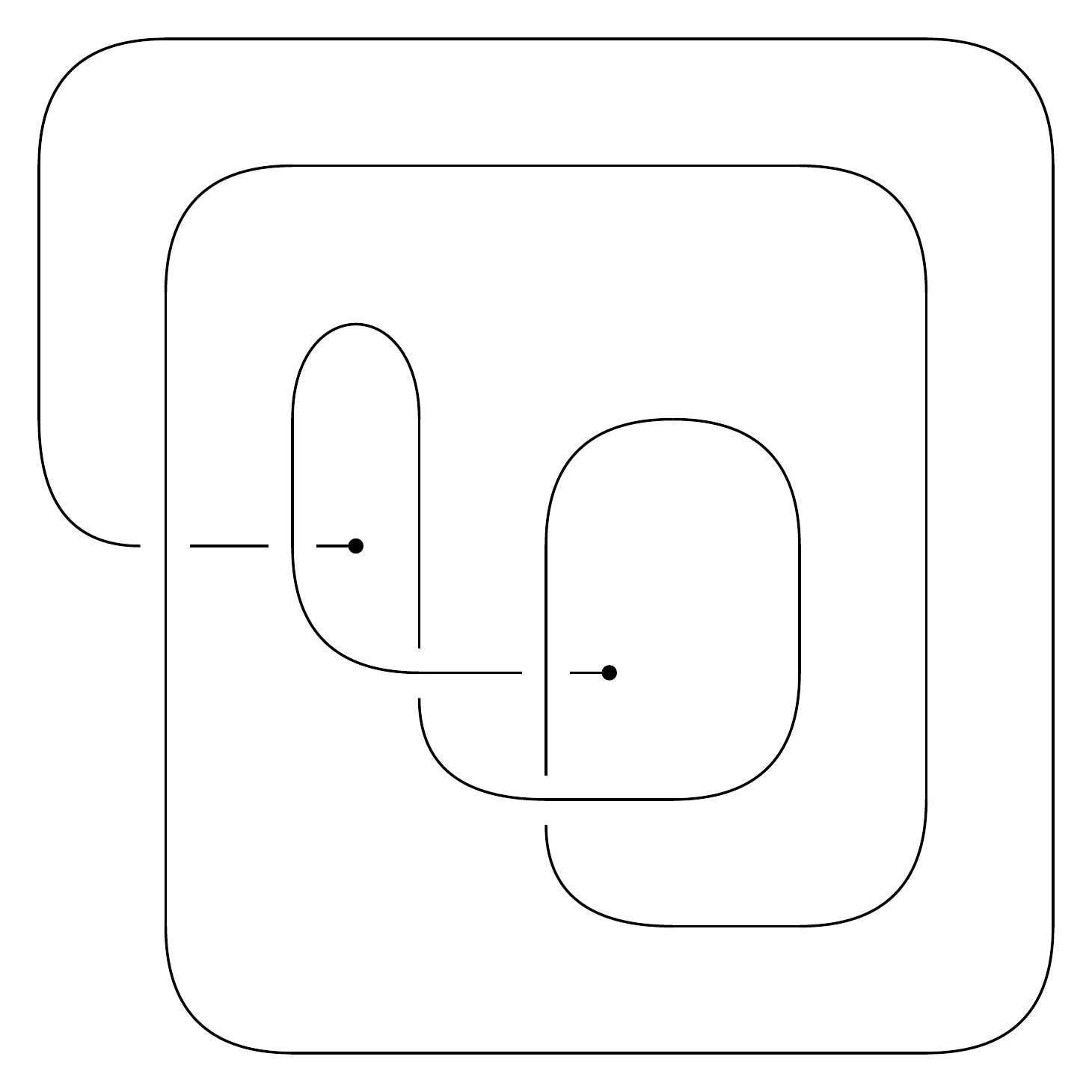}\\
\textcolor{black}{$5_{78}$}
\vspace{1cm}
\end{minipage}
\begin{minipage}[t]{.25\linewidth}
\centering
\includegraphics[width=0.9\textwidth,height=3.5cm,keepaspectratio]{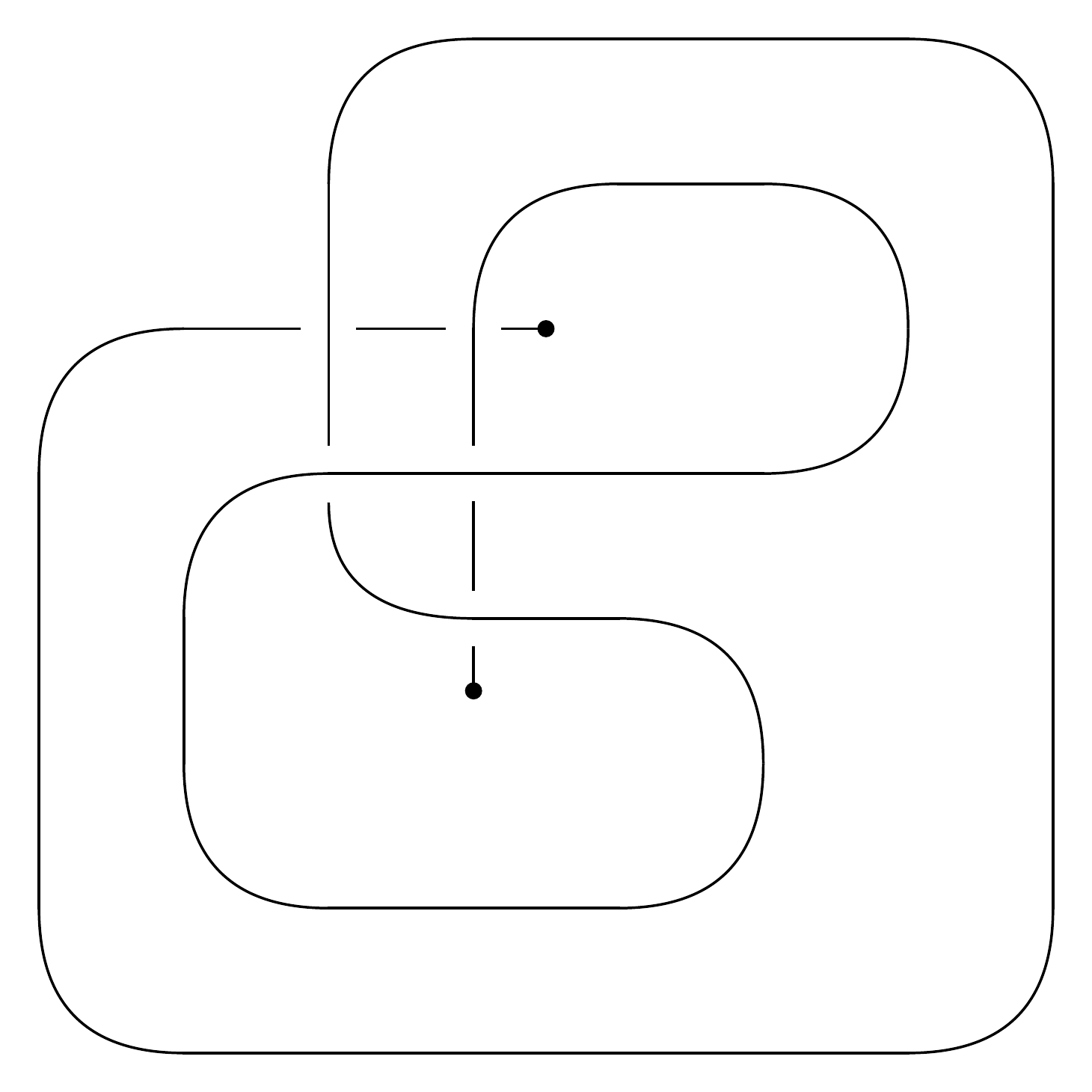}\\
\textcolor{black}{$5_{79}$}
\vspace{1cm}
\end{minipage}
\begin{minipage}[t]{.25\linewidth}
\centering
\includegraphics[width=0.9\textwidth,height=3.5cm,keepaspectratio]{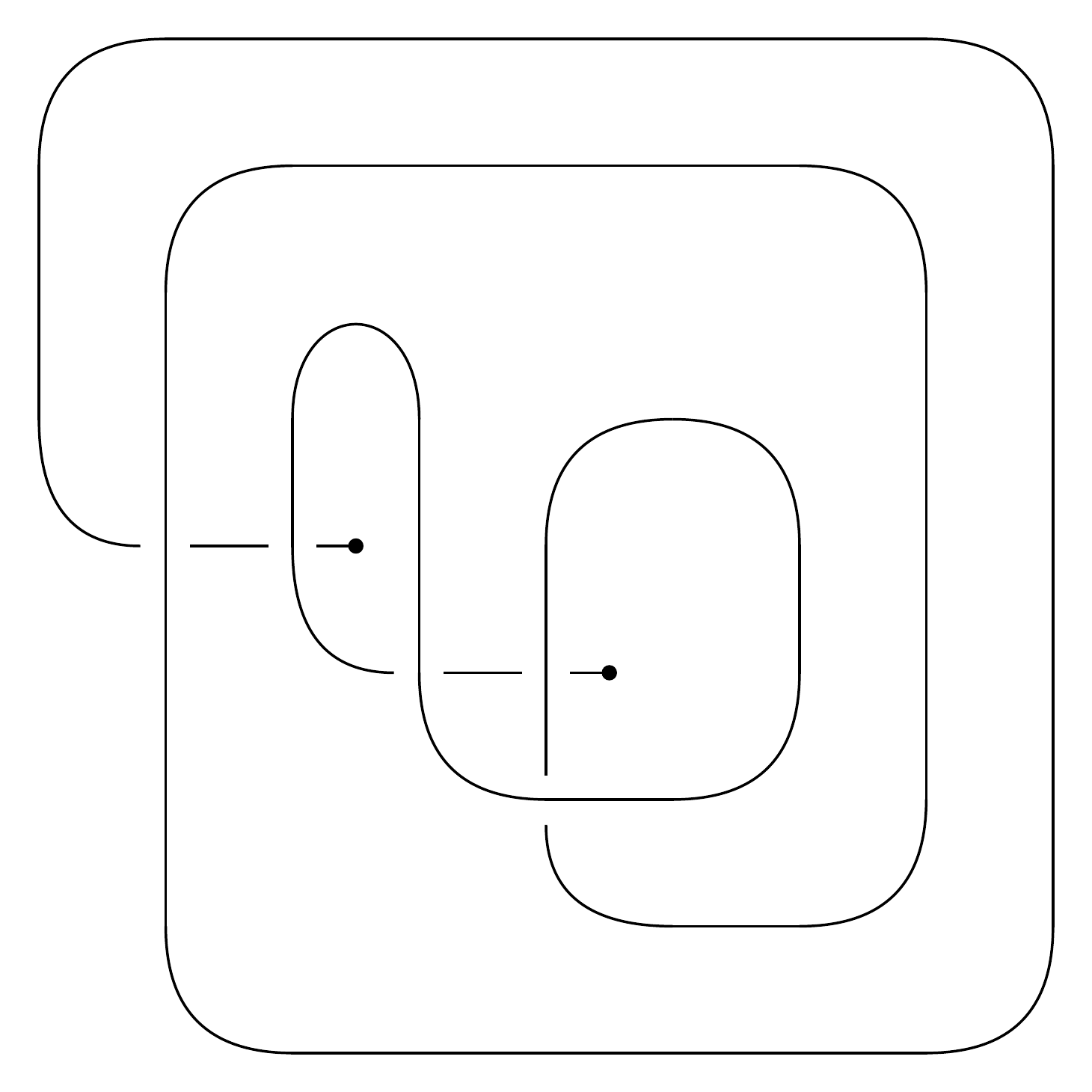}\\
\textcolor{black}{$5_{80}$}
\vspace{1cm}
\end{minipage}
\begin{minipage}[t]{.25\linewidth}
\centering
\includegraphics[width=0.9\textwidth,height=3.5cm,keepaspectratio]{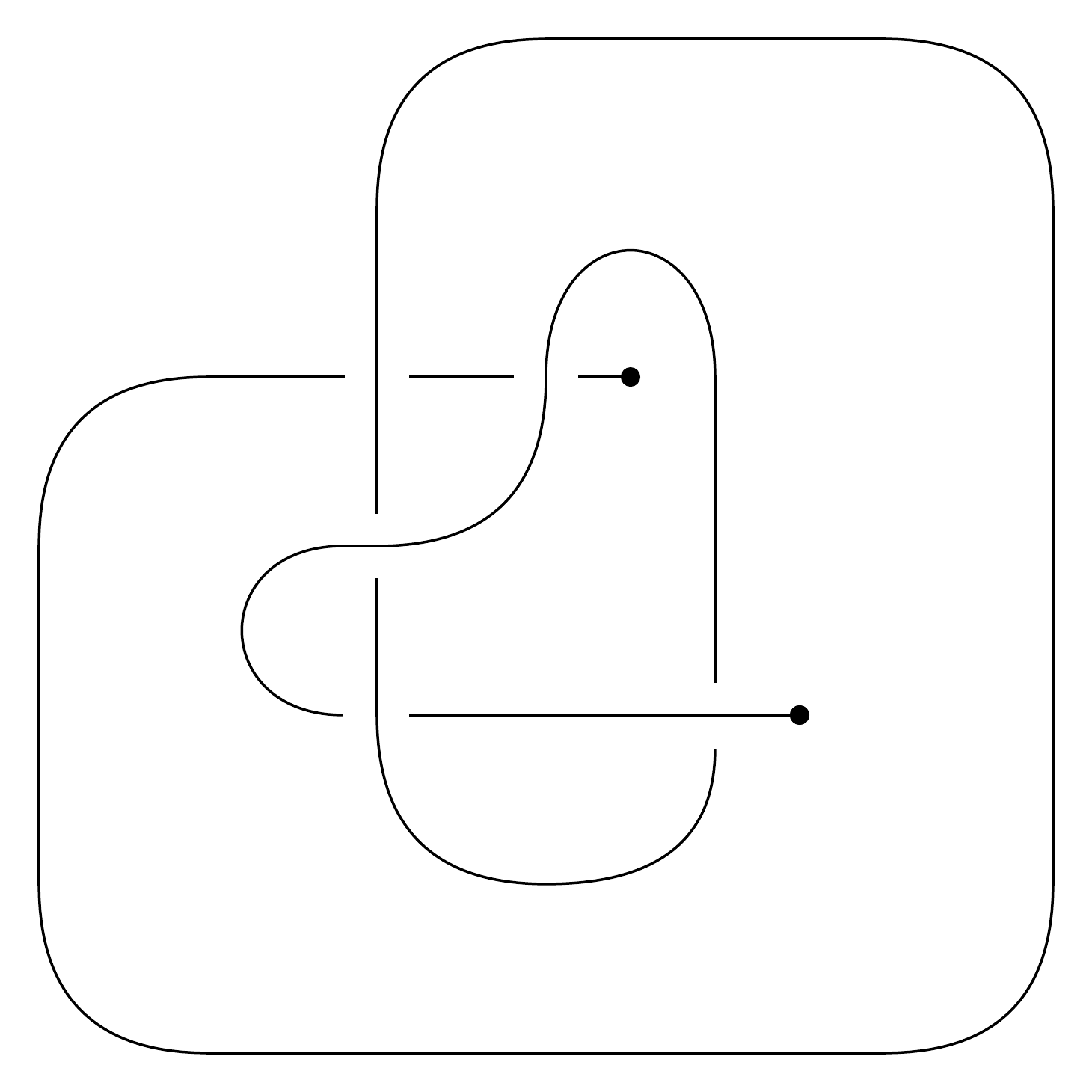}\\
\textcolor{black}{$5_{81}$}
\vspace{1cm}
\end{minipage}
\begin{minipage}[t]{.25\linewidth}
\centering
\includegraphics[width=0.9\textwidth,height=3.5cm,keepaspectratio]{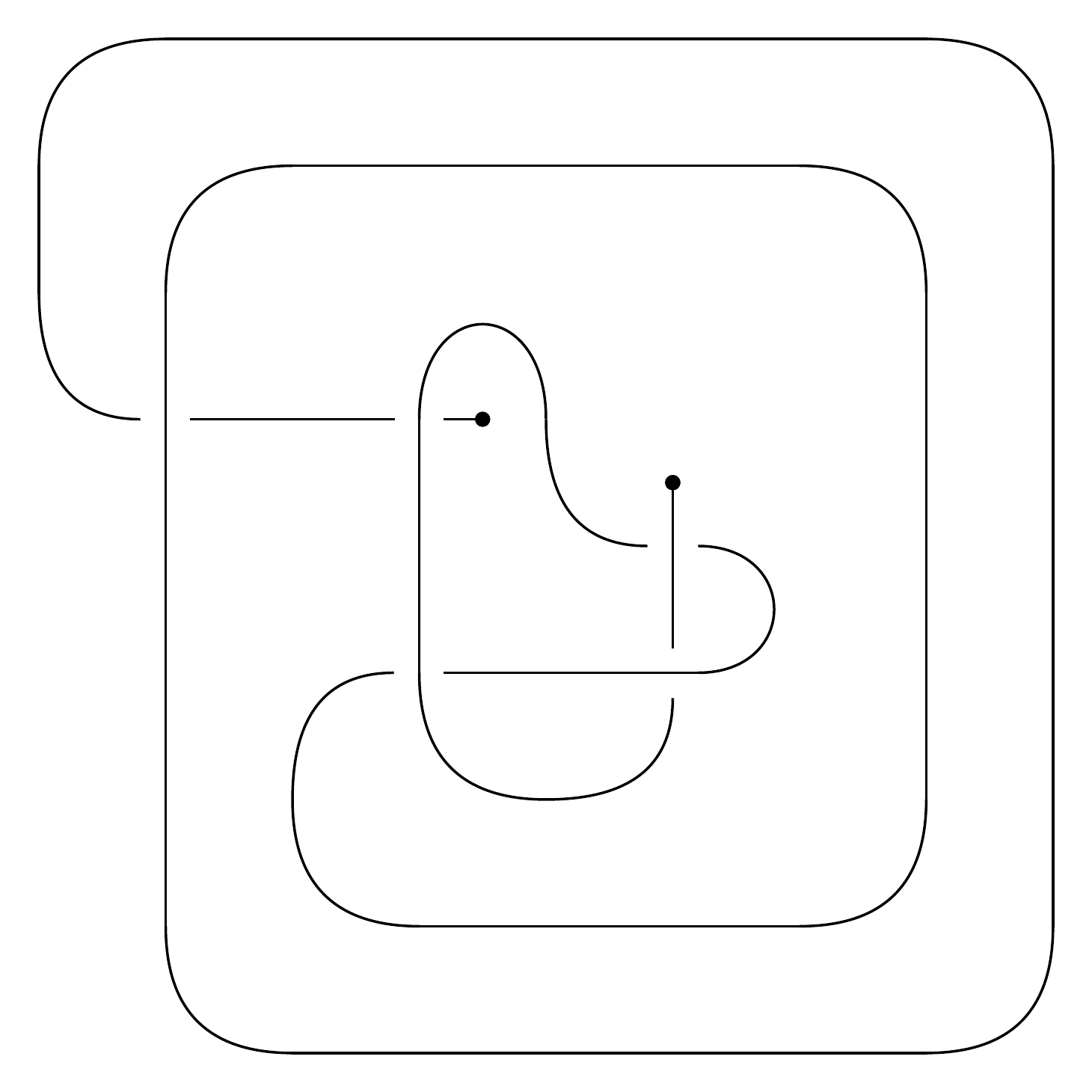}\\
\textcolor{black}{$5_{82}$}
\vspace{1cm}
\end{minipage}
\begin{minipage}[t]{.25\linewidth}
\centering
\includegraphics[width=0.9\textwidth,height=3.5cm,keepaspectratio]{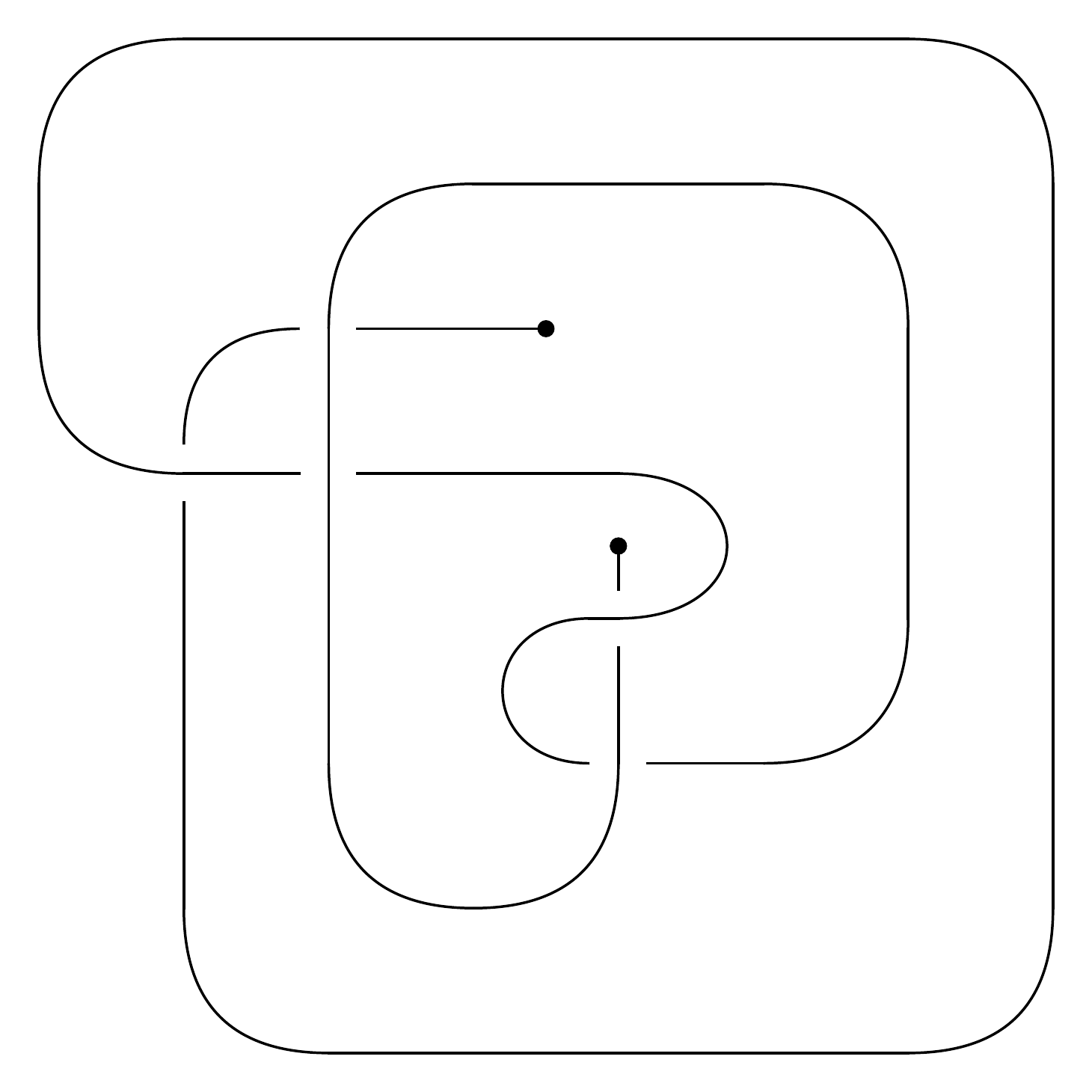}\\
\textcolor{black}{$5_{83}$}
\vspace{1cm}
\end{minipage}
\begin{minipage}[t]{.25\linewidth}
\centering
\includegraphics[width=0.9\textwidth,height=3.5cm,keepaspectratio]{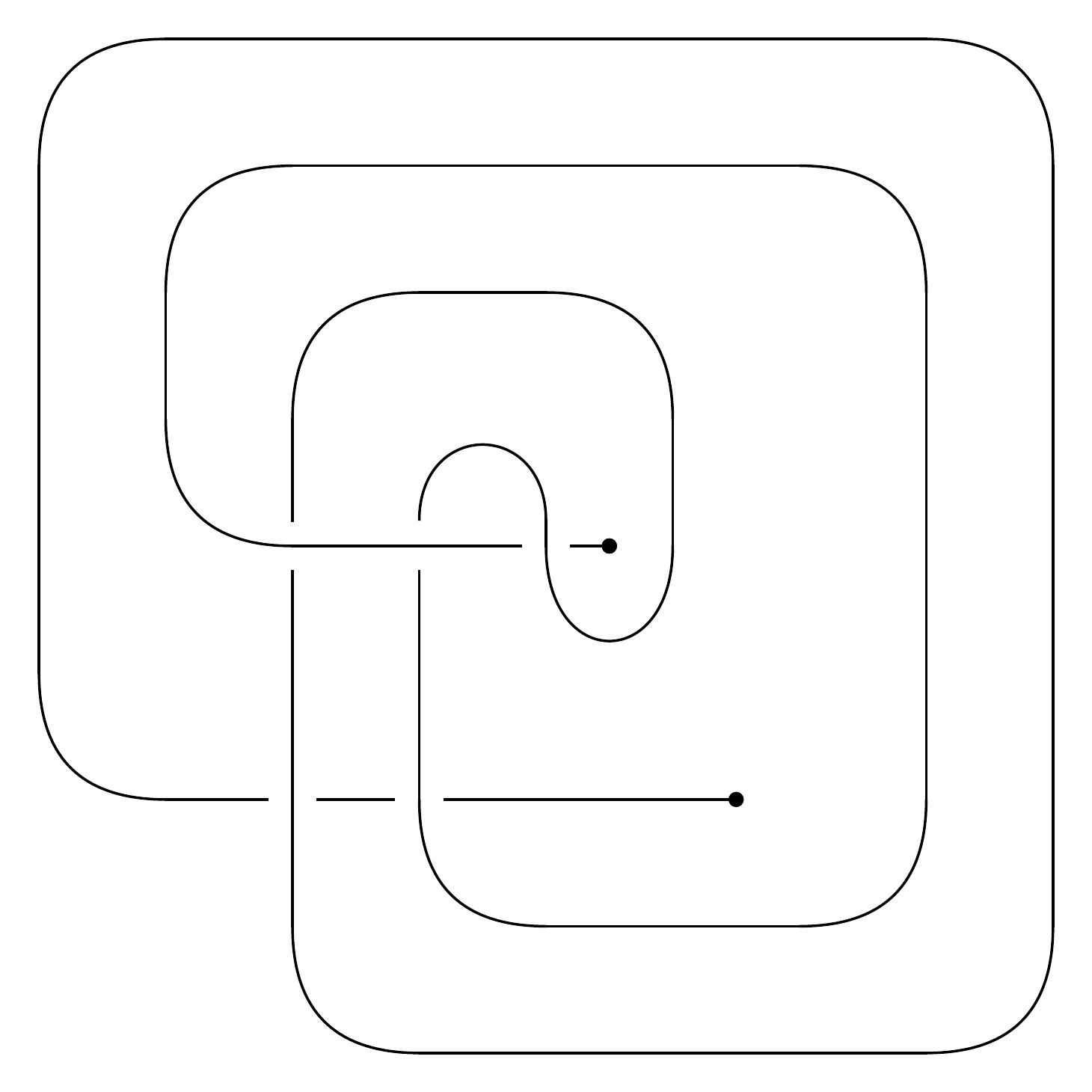}\\
\textcolor{black}{$5_{84}$}
\vspace{1cm}
\end{minipage}
\begin{minipage}[t]{.25\linewidth}
\centering
\includegraphics[width=0.9\textwidth,height=3.5cm,keepaspectratio]{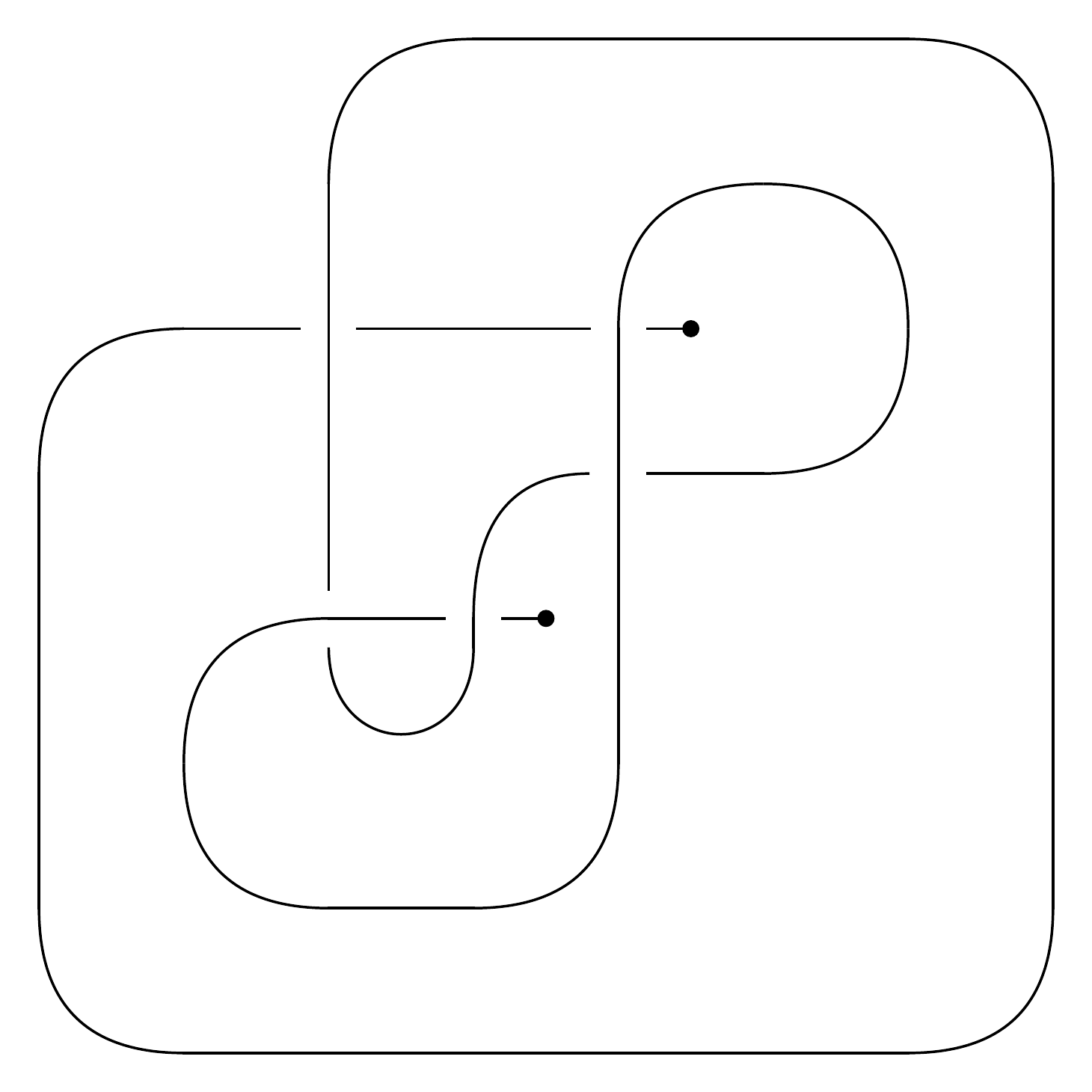}\\
\textcolor{black}{$5_{85}$}
\vspace{1cm}
\end{minipage}
\begin{minipage}[t]{.25\linewidth}
\centering
\includegraphics[width=0.9\textwidth,height=3.5cm,keepaspectratio]{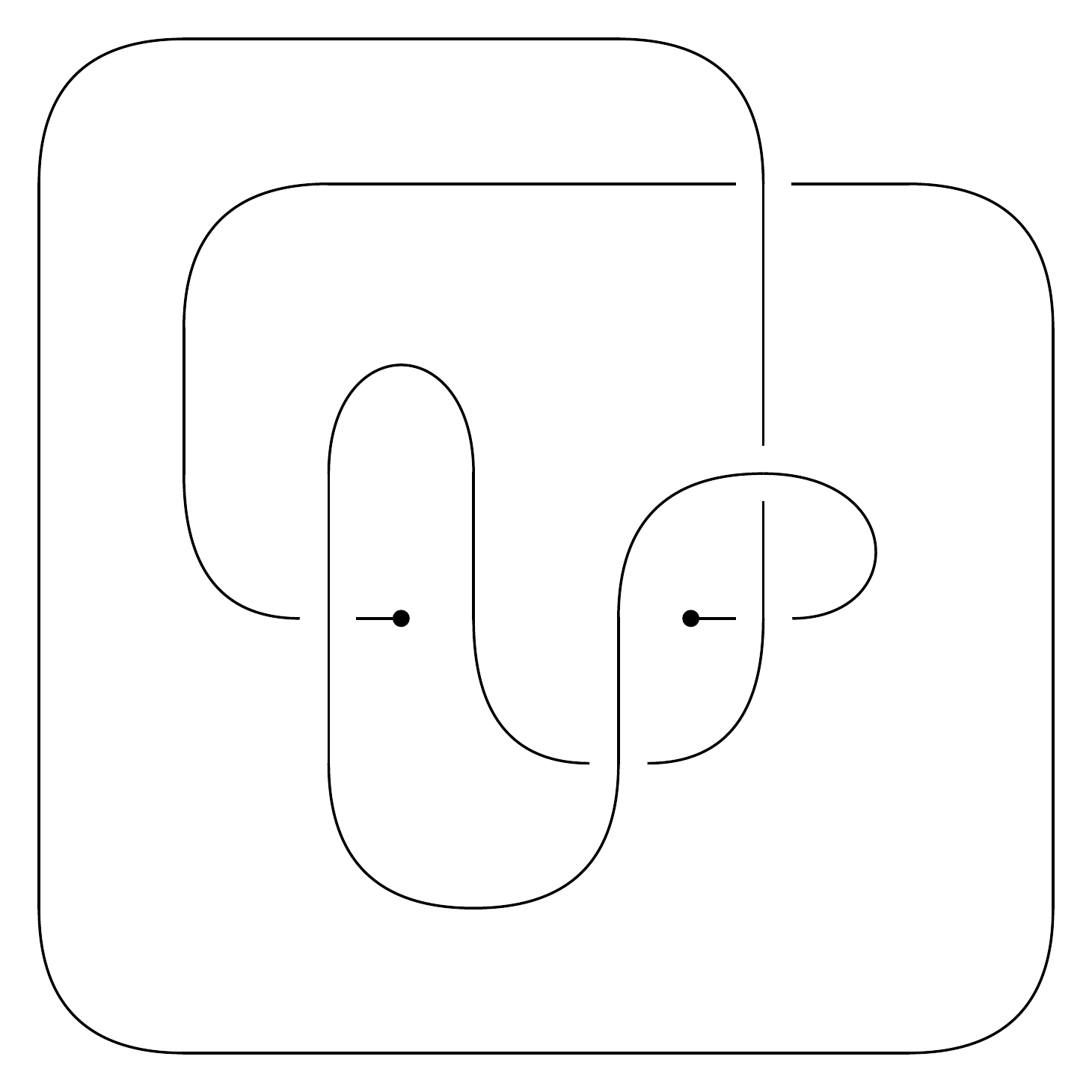}\\
\textcolor{black}{$5_{86}$}
\vspace{1cm}
\end{minipage}
\begin{minipage}[t]{.25\linewidth}
\centering
\includegraphics[width=0.9\textwidth,height=3.5cm,keepaspectratio]{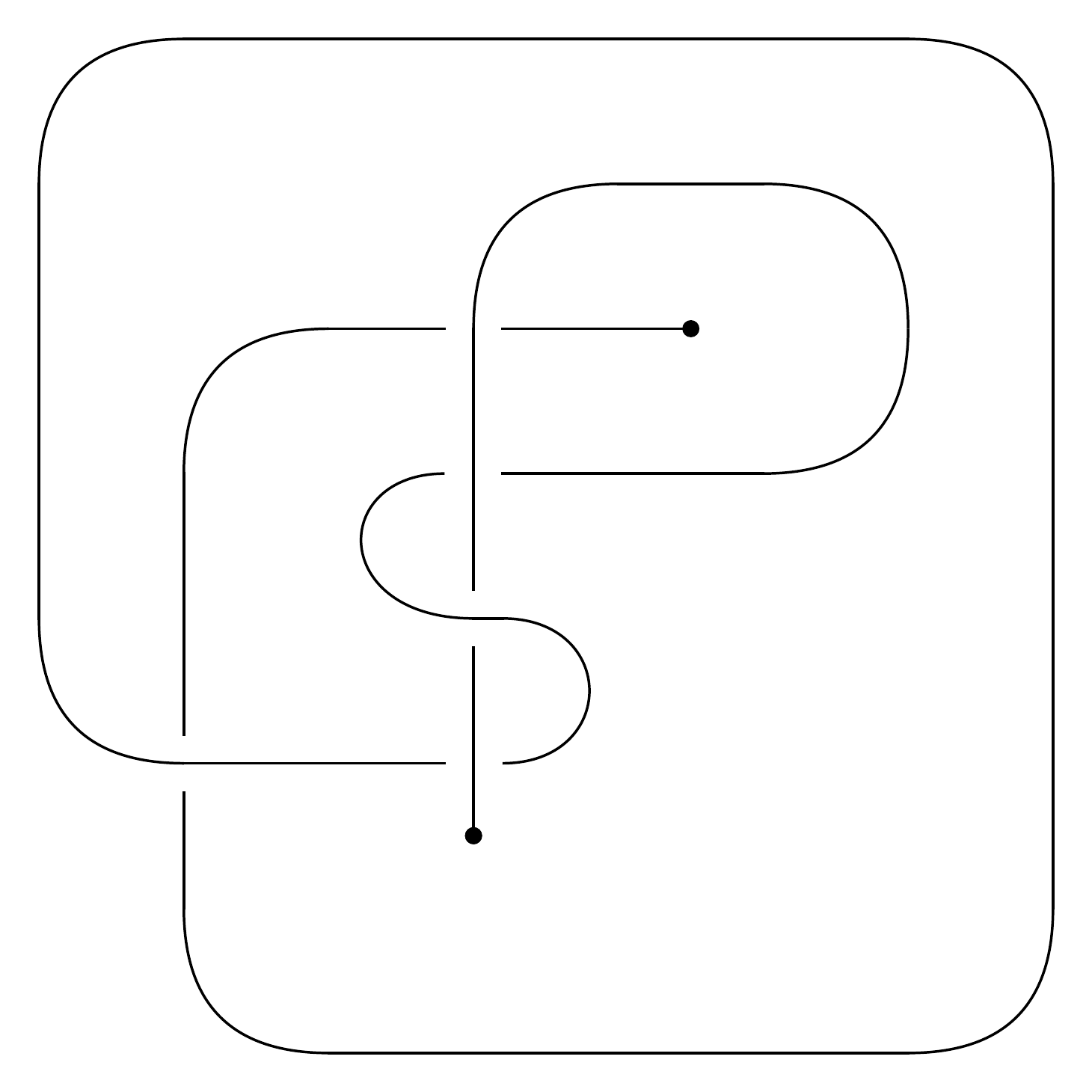}\\
\textcolor{black}{$5_{87}$}
\vspace{1cm}
\end{minipage}
\begin{minipage}[t]{.25\linewidth}
\centering
\includegraphics[width=0.9\textwidth,height=3.5cm,keepaspectratio]{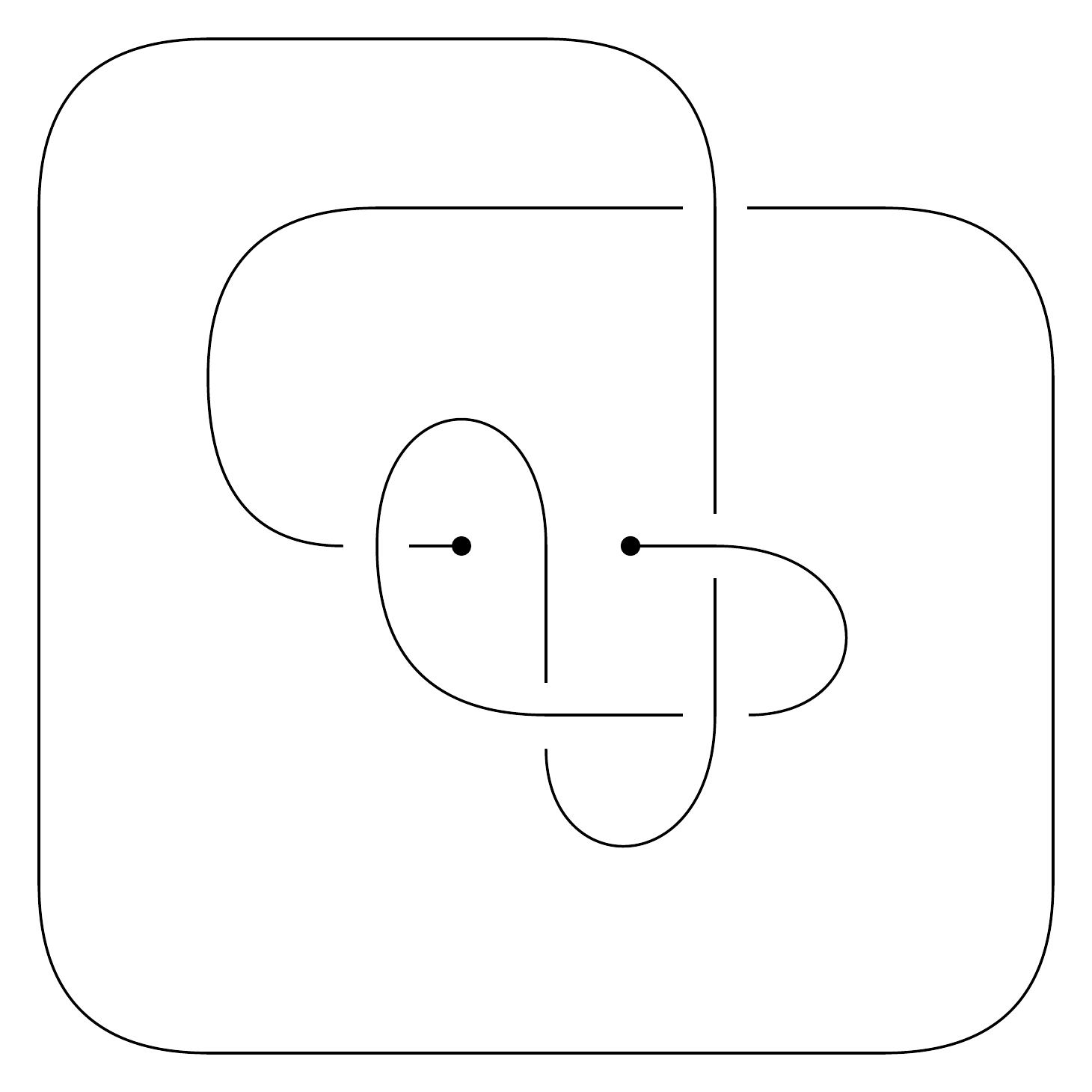}\\
\textcolor{black}{$5_{88}$}
\vspace{1cm}
\end{minipage}
\begin{minipage}[t]{.25\linewidth}
\centering
\includegraphics[width=0.9\textwidth,height=3.5cm,keepaspectratio]{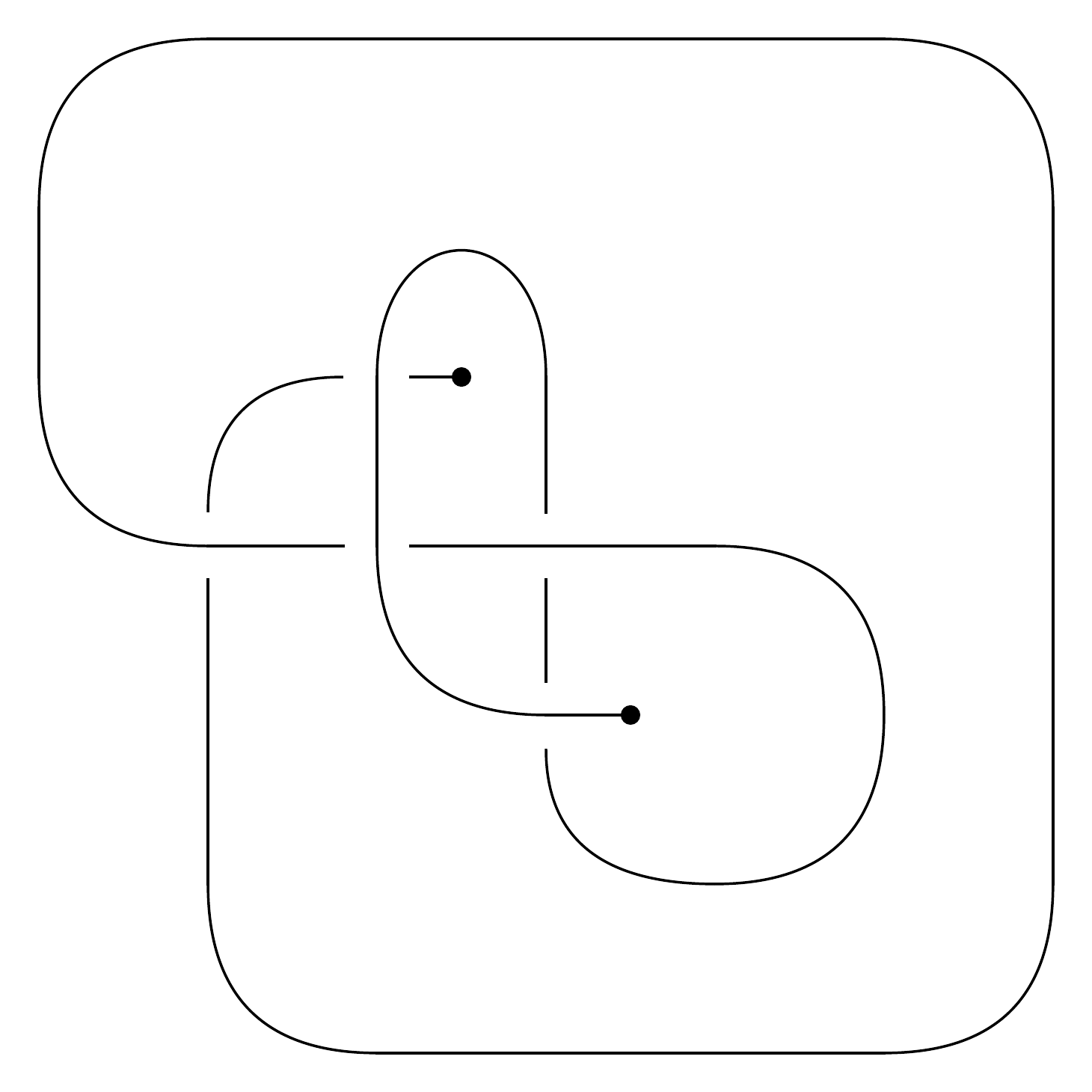}\\
\textcolor{black}{$5_{89}$}
\vspace{1cm}
\end{minipage}
\begin{minipage}[t]{.25\linewidth}
\centering
\includegraphics[width=0.9\textwidth,height=3.5cm,keepaspectratio]{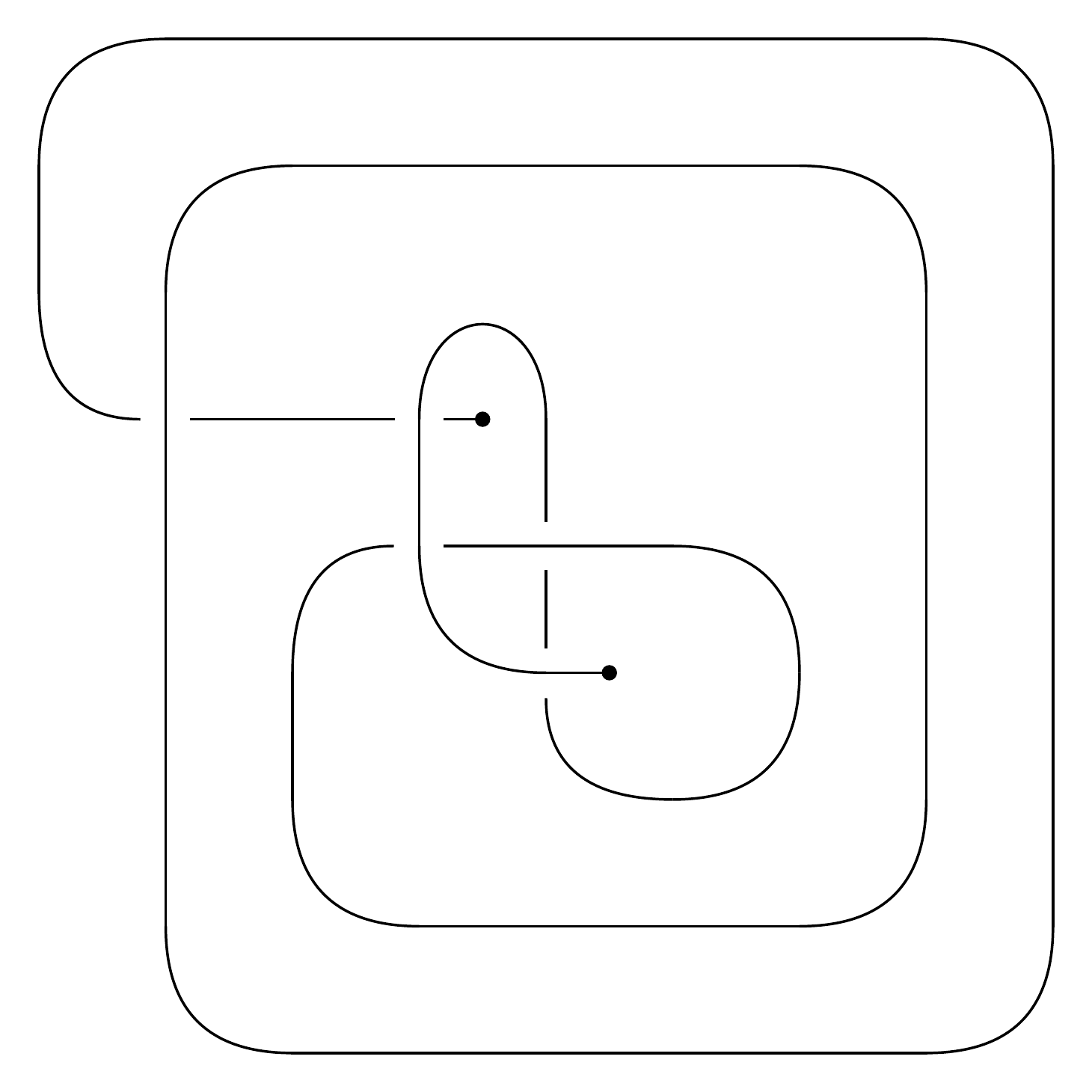}\\
\textcolor{black}{$5_{90}$}
\vspace{1cm}
\end{minipage}
\begin{minipage}[t]{.25\linewidth}
\centering
\includegraphics[width=0.9\textwidth,height=3.5cm,keepaspectratio]{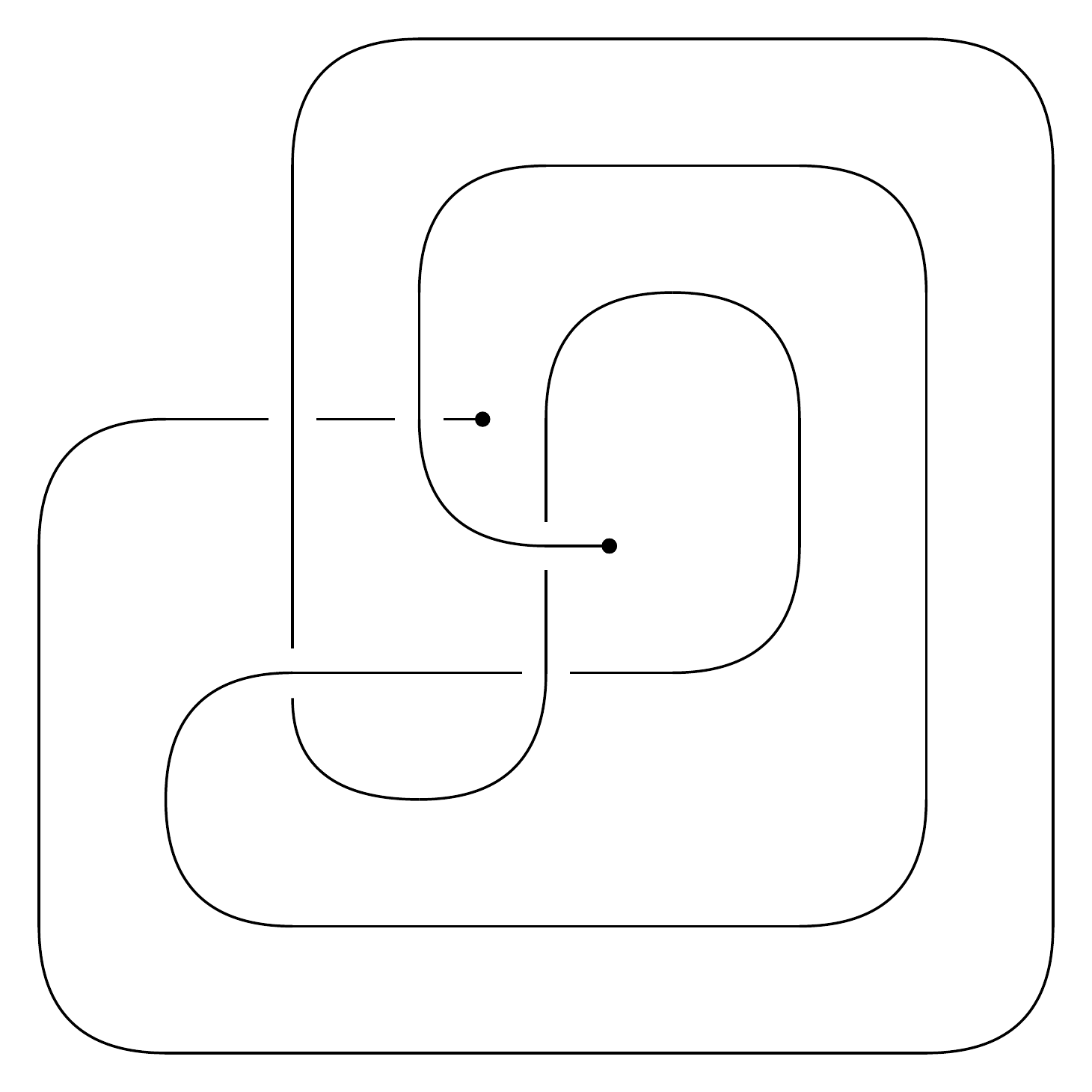}\\
\textcolor{black}{$5_{91}$}
\vspace{1cm}
\end{minipage}
\begin{minipage}[t]{.25\linewidth}
\centering
\includegraphics[width=0.9\textwidth,height=3.5cm,keepaspectratio]{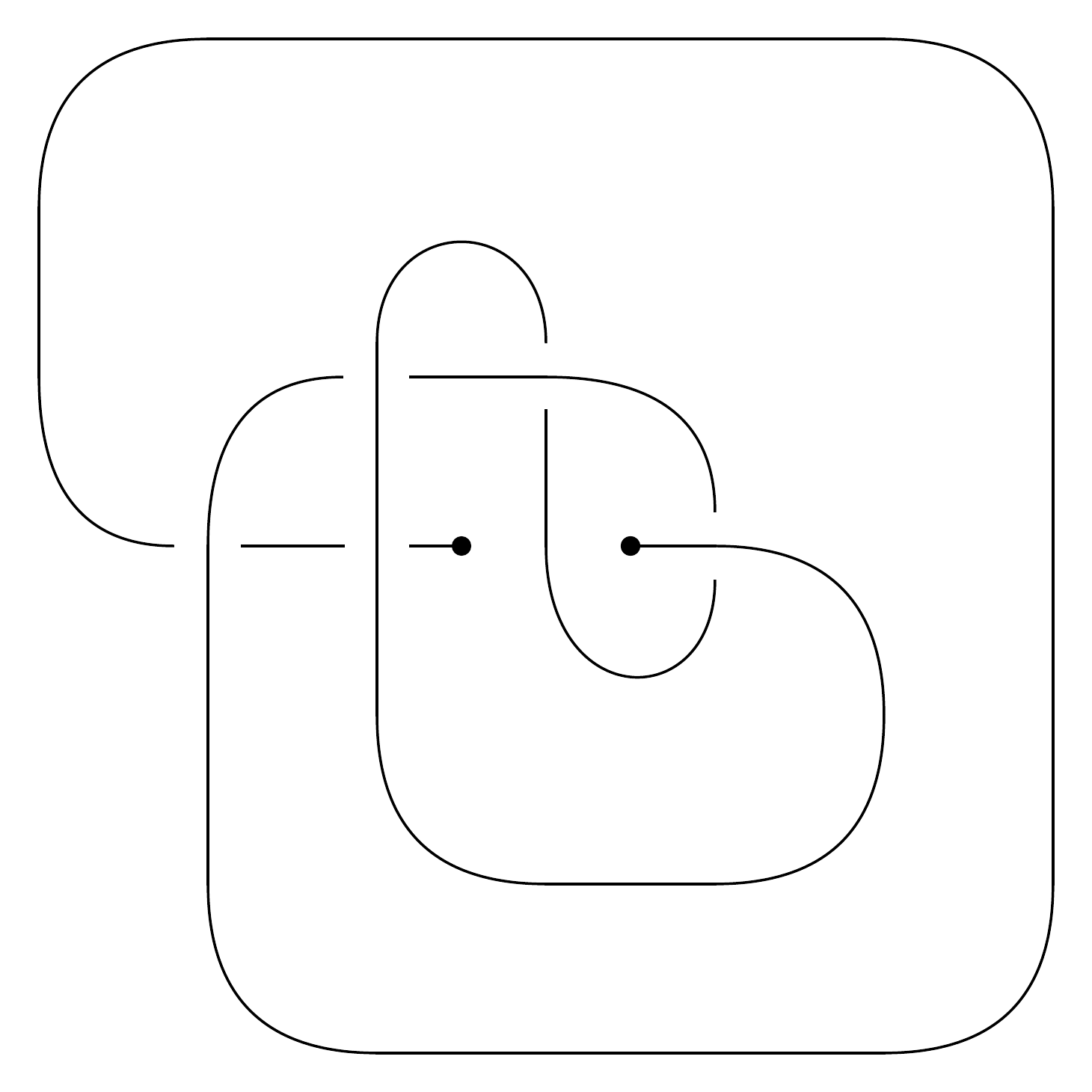}\\
\textcolor{black}{$5_{92}$}
\vspace{1cm}
\end{minipage}
\begin{minipage}[t]{.25\linewidth}
\centering
\includegraphics[width=0.9\textwidth,height=3.5cm,keepaspectratio]{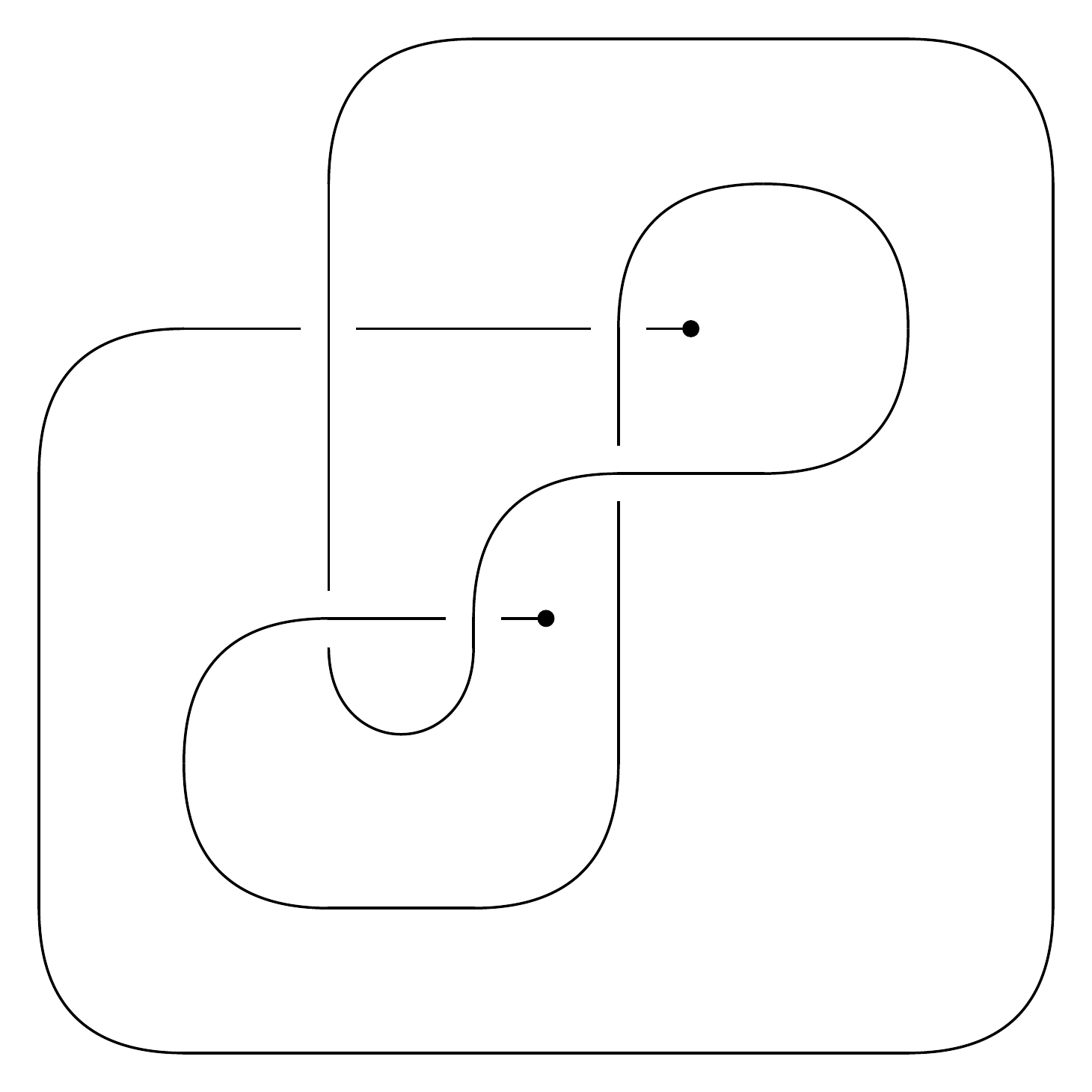}\\
\textcolor{black}{$5_{93}$}
\vspace{1cm}
\end{minipage}
\begin{minipage}[t]{.25\linewidth}
\centering
\includegraphics[width=0.9\textwidth,height=3.5cm,keepaspectratio]{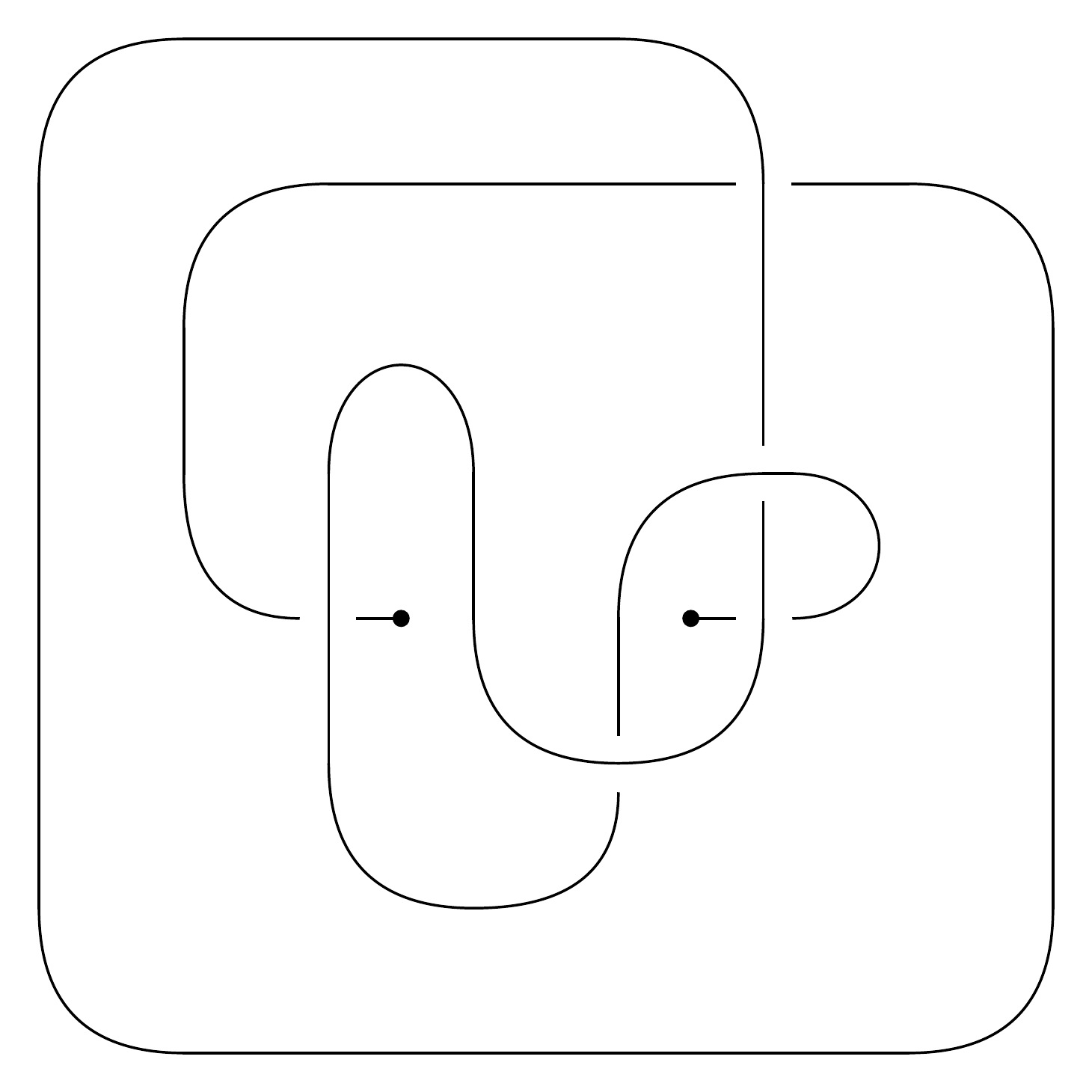}\\
\textcolor{black}{$5_{94}$}
\vspace{1cm}
\end{minipage}
\begin{minipage}[t]{.25\linewidth}
\centering
\includegraphics[width=0.9\textwidth,height=3.5cm,keepaspectratio]{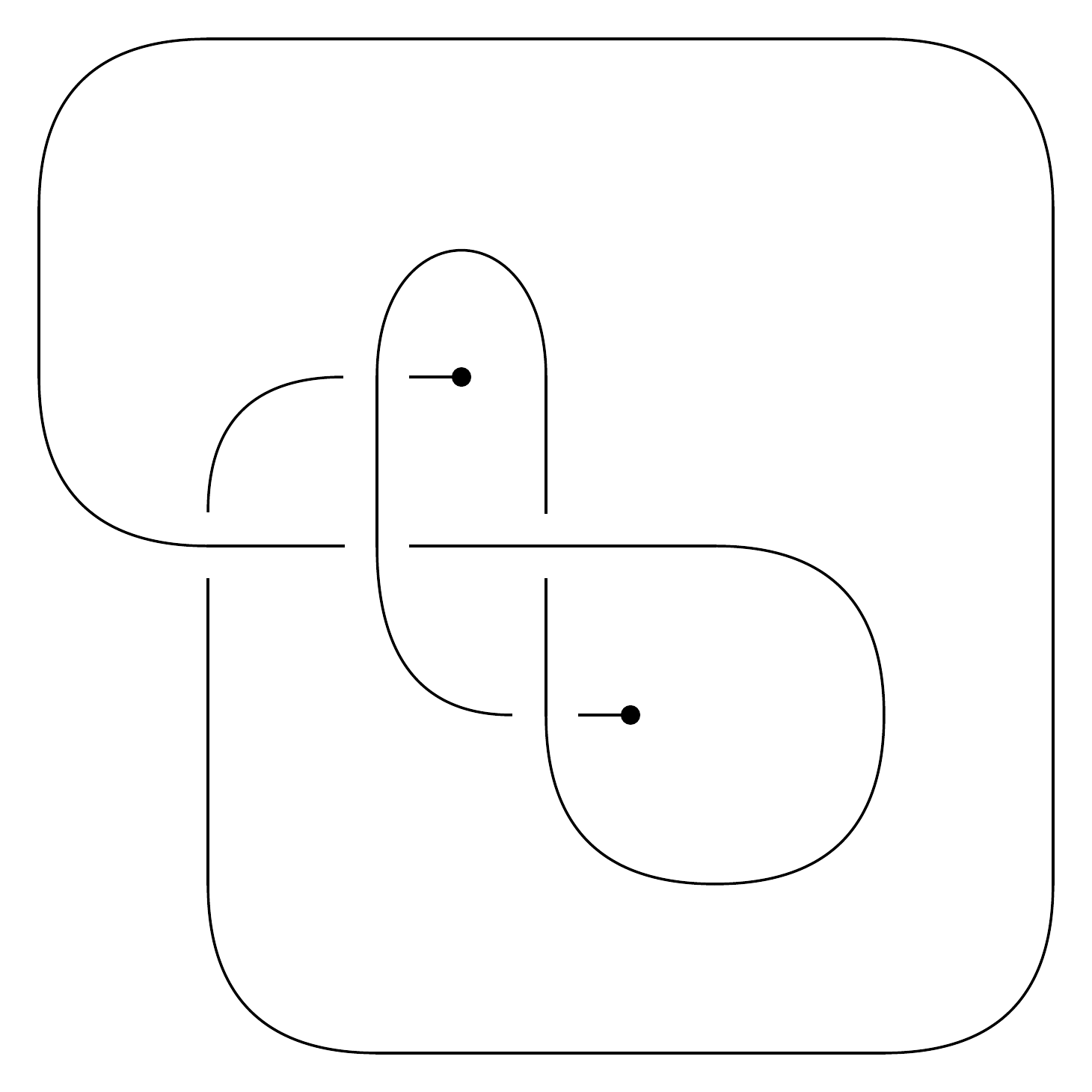}\\
\textcolor{black}{$5_{95}$}
\vspace{1cm}
\end{minipage}
\begin{minipage}[t]{.25\linewidth}
\centering
\includegraphics[width=0.9\textwidth,height=3.5cm,keepaspectratio]{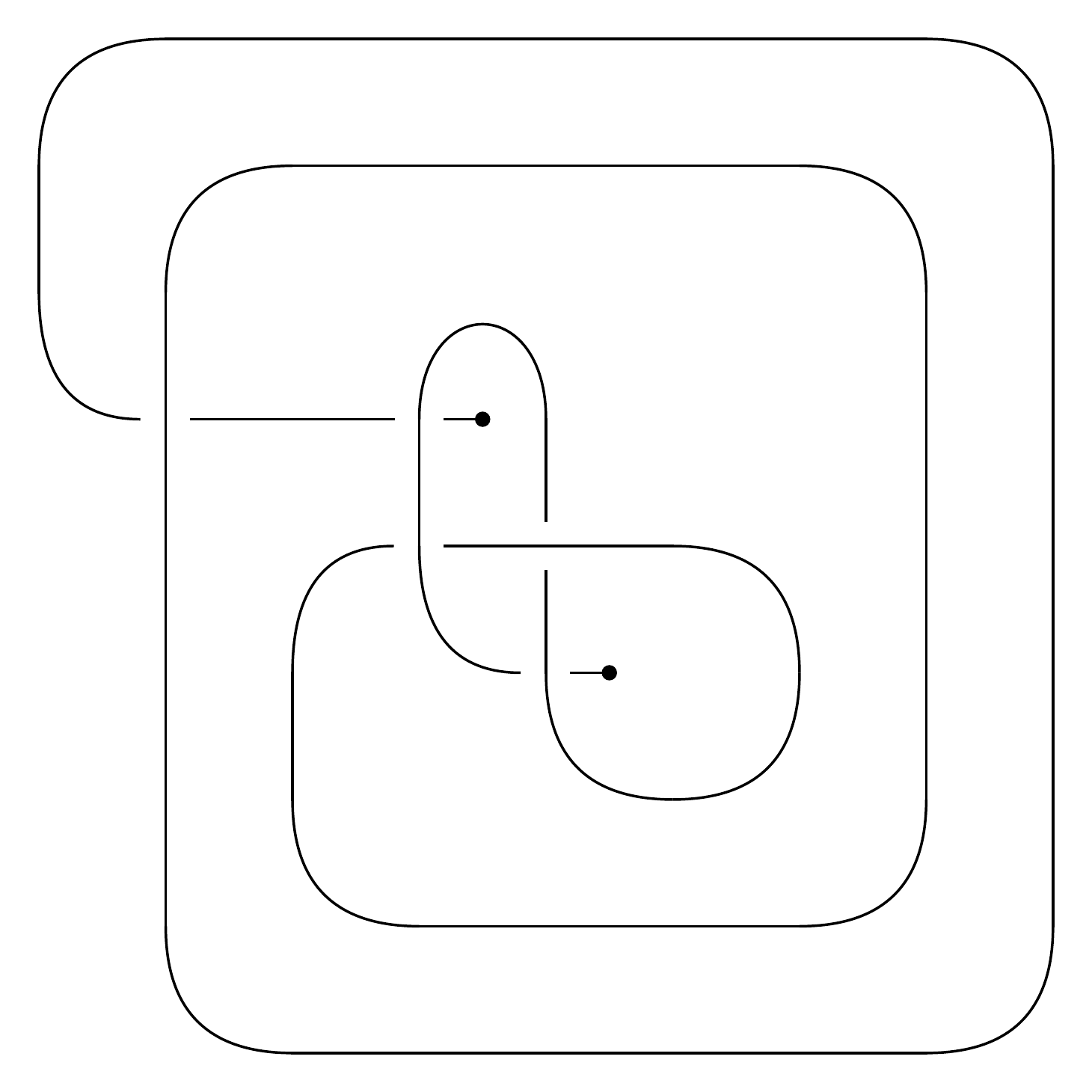}\\
\textcolor{black}{$5_{96}$}
\vspace{1cm}
\end{minipage}
\begin{minipage}[t]{.25\linewidth}
\centering
\includegraphics[width=0.9\textwidth,height=3.5cm,keepaspectratio]{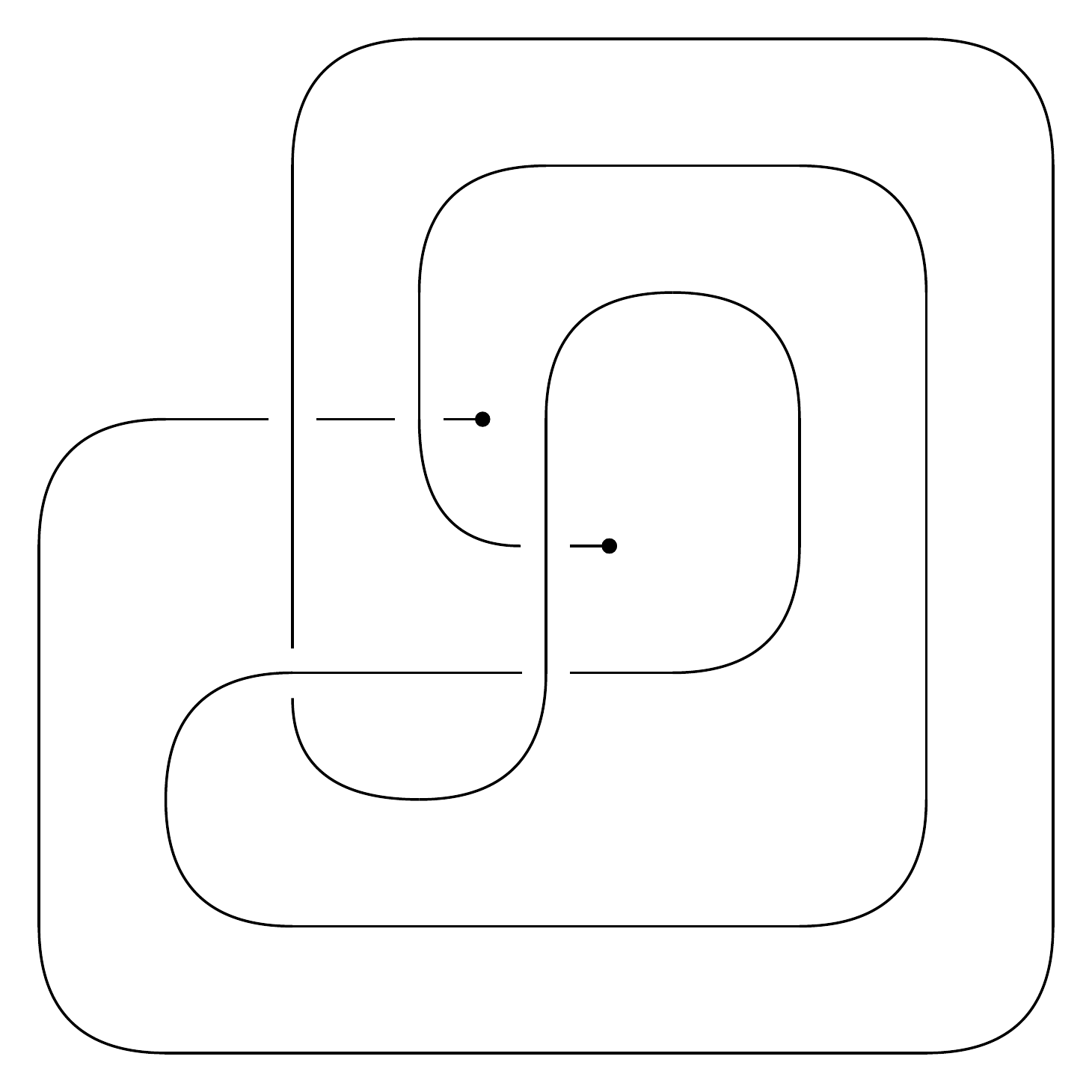}\\
\textcolor{black}{$5_{97}$}
\vspace{1cm}
\end{minipage}
\begin{minipage}[t]{.25\linewidth}
\centering
\includegraphics[width=0.9\textwidth,height=3.5cm,keepaspectratio]{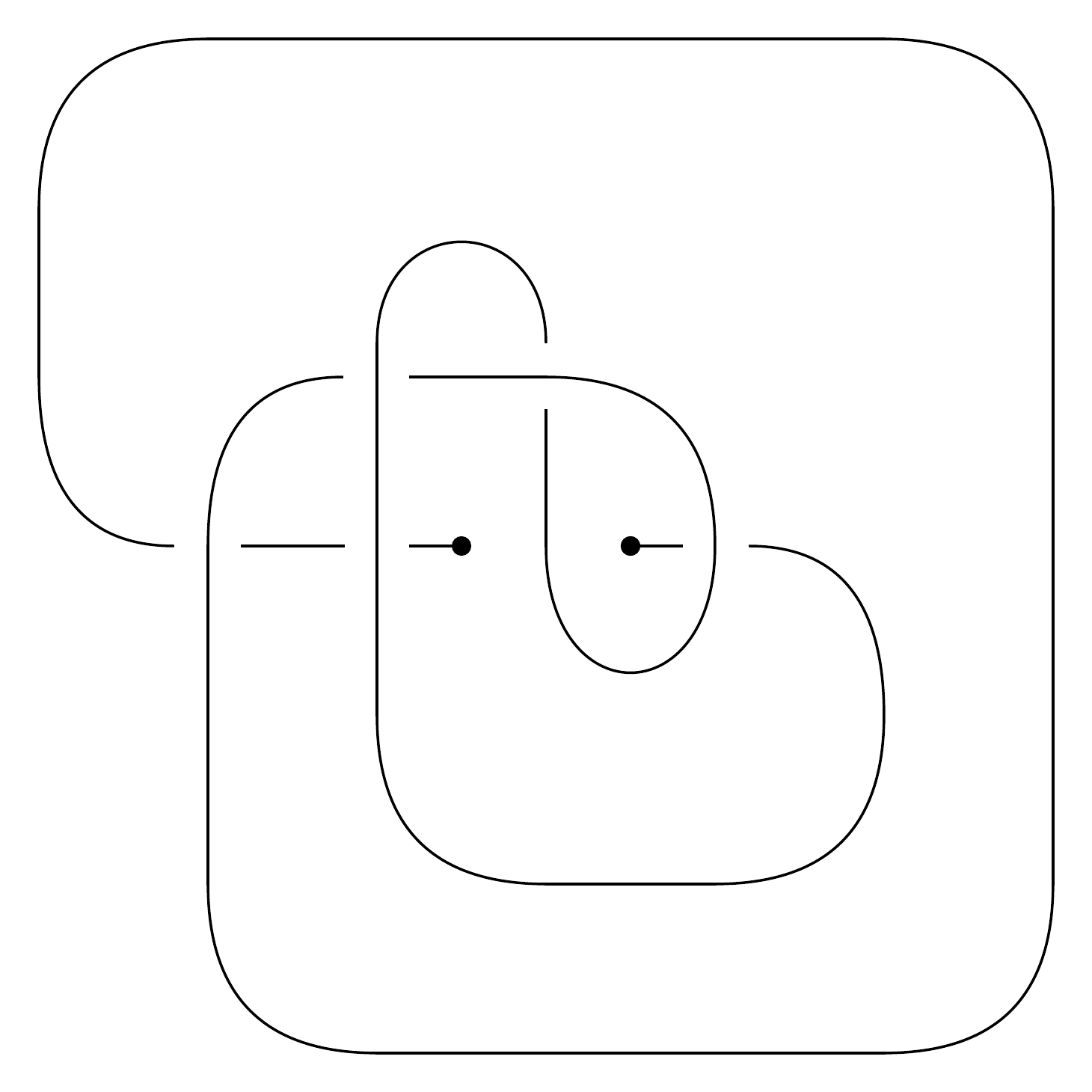}\\
\textcolor{black}{$5_{98}$}
\vspace{1cm}
\end{minipage}
\begin{minipage}[t]{.25\linewidth}
\centering
\includegraphics[width=0.9\textwidth,height=3.5cm,keepaspectratio]{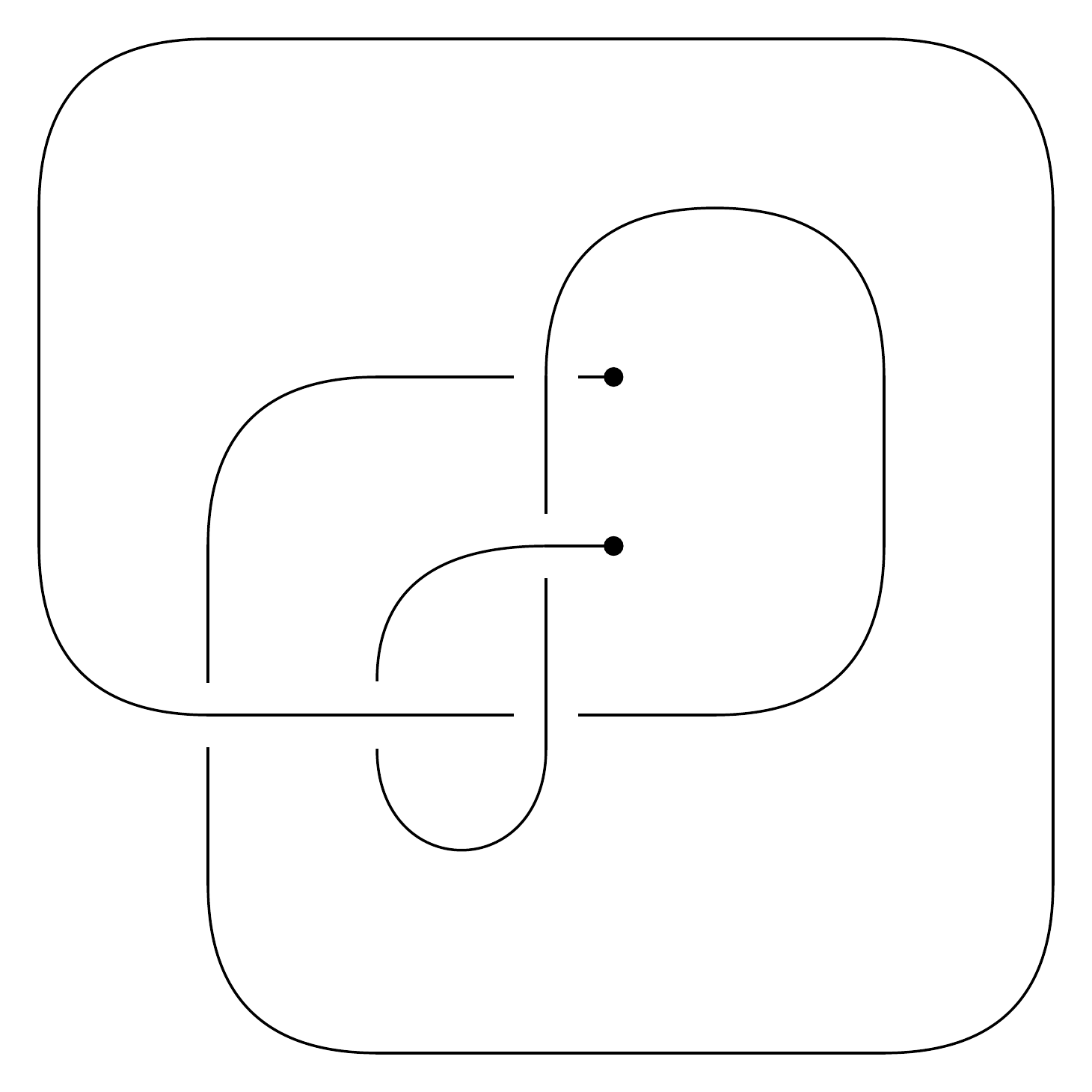}\\
\textcolor{black}{$5_{99}$}
\vspace{1cm}
\end{minipage}
\begin{minipage}[t]{.25\linewidth}
\centering
\includegraphics[width=0.9\textwidth,height=3.5cm,keepaspectratio]{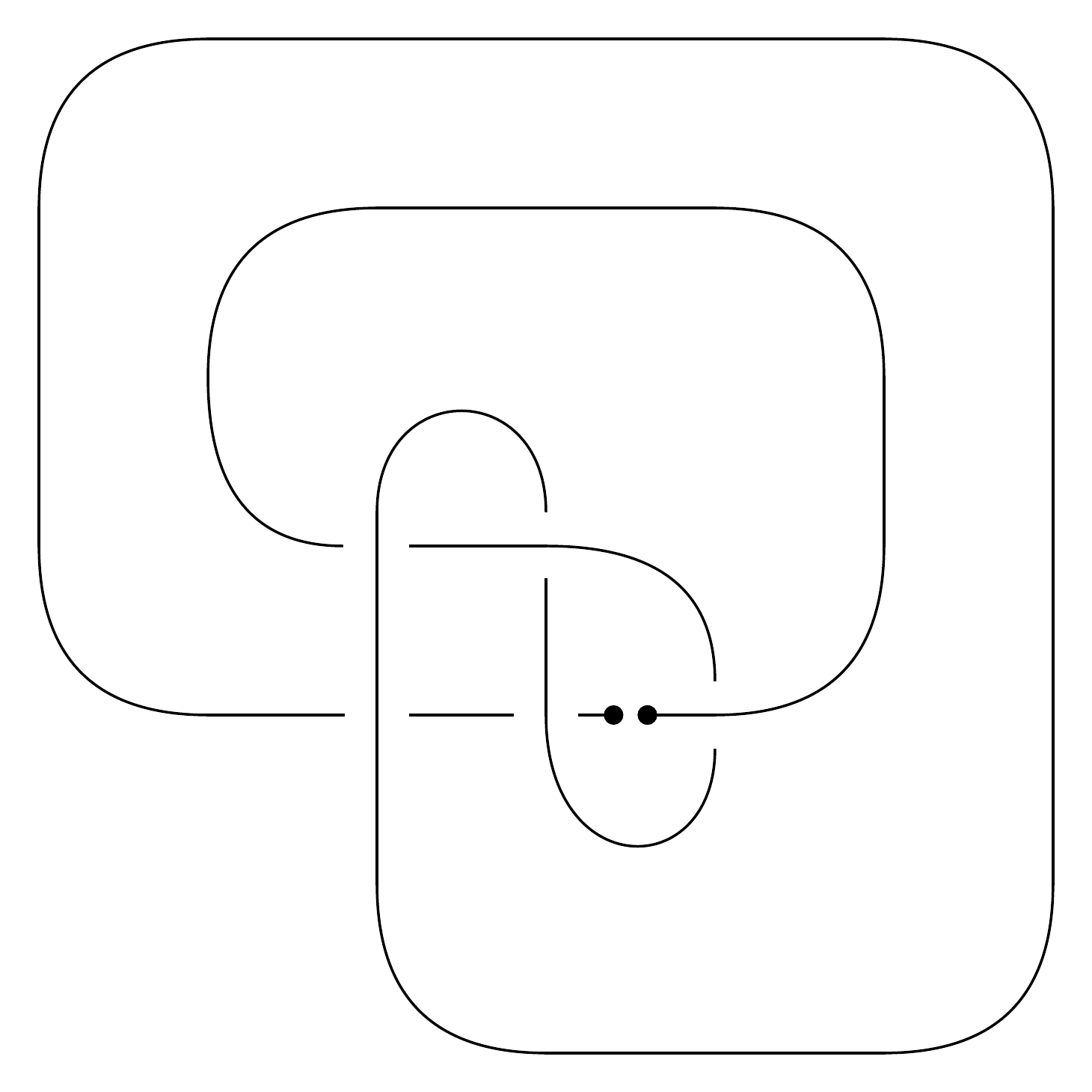}\\
\textcolor{black}{$5_{100}$}
\vspace{1cm}
\end{minipage}
\begin{minipage}[t]{.25\linewidth}
\centering
\includegraphics[width=0.9\textwidth,height=3.5cm,keepaspectratio]{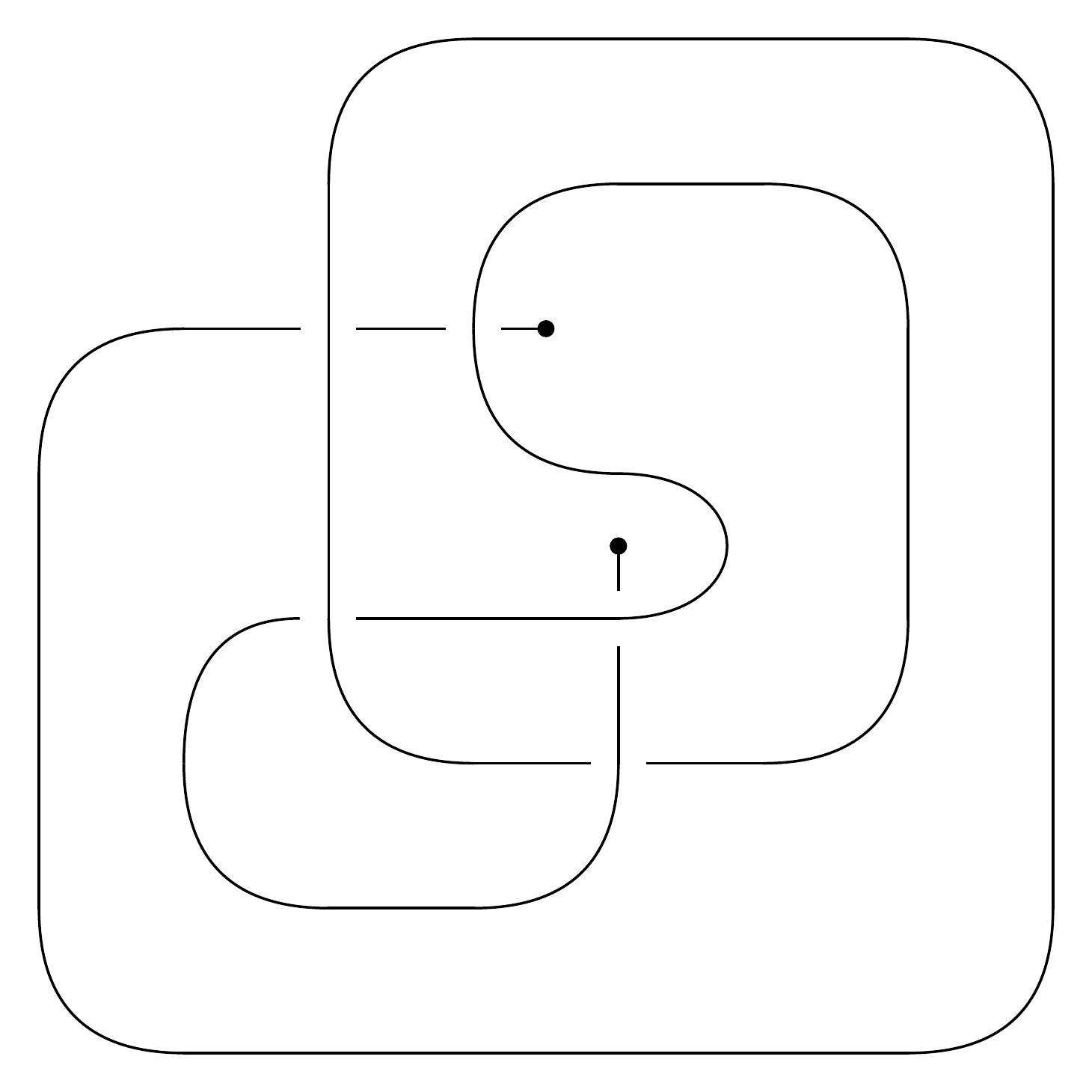}\\
\textcolor{black}{$5_{101}$}
\vspace{1cm}
\end{minipage}
\begin{minipage}[t]{.25\linewidth}
\centering
\includegraphics[width=0.9\textwidth,height=3.5cm,keepaspectratio]{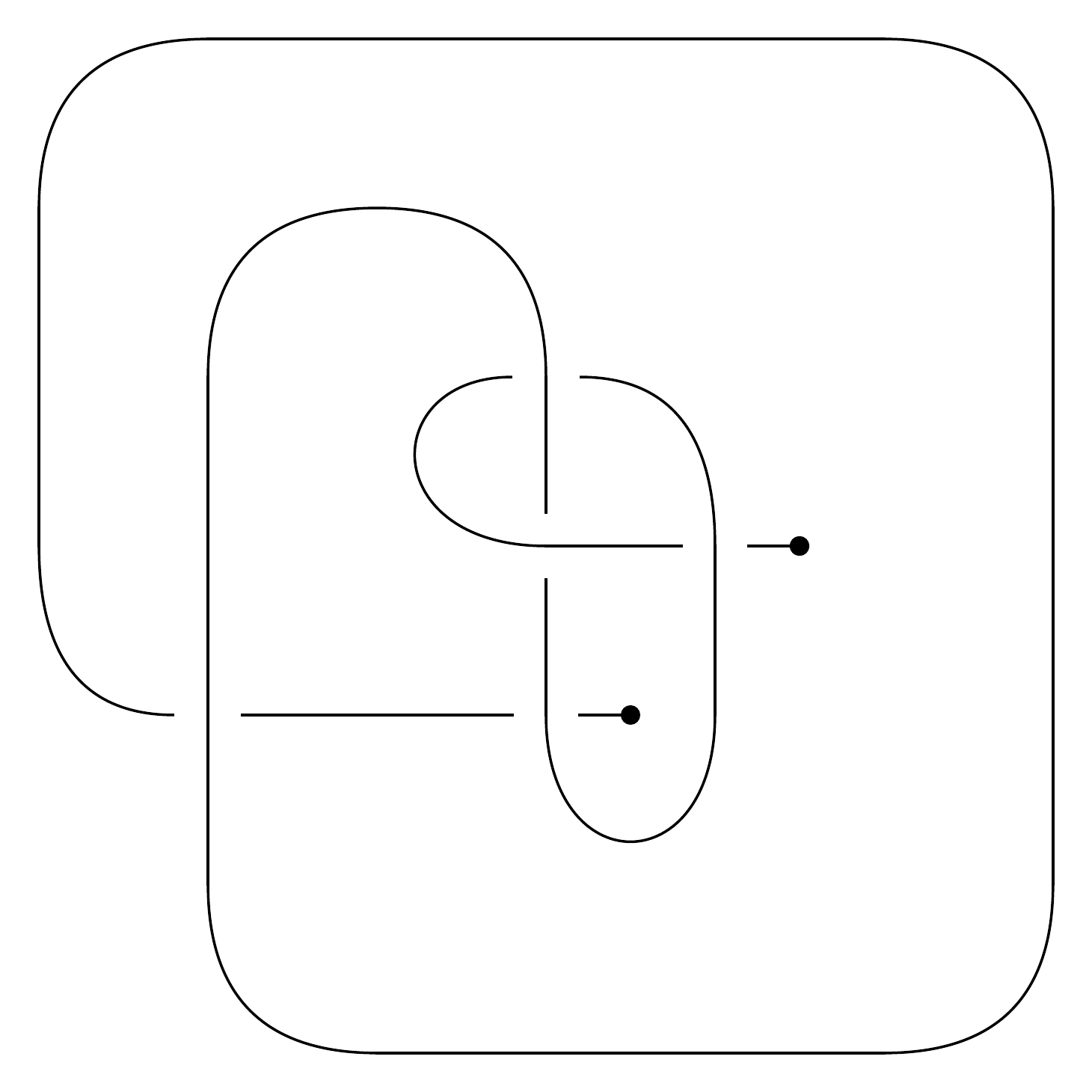}\\
\textcolor{black}{$5_{102}$}
\vspace{1cm}
\end{minipage}
\begin{minipage}[t]{.25\linewidth}
\centering
\includegraphics[width=0.9\textwidth,height=3.5cm,keepaspectratio]{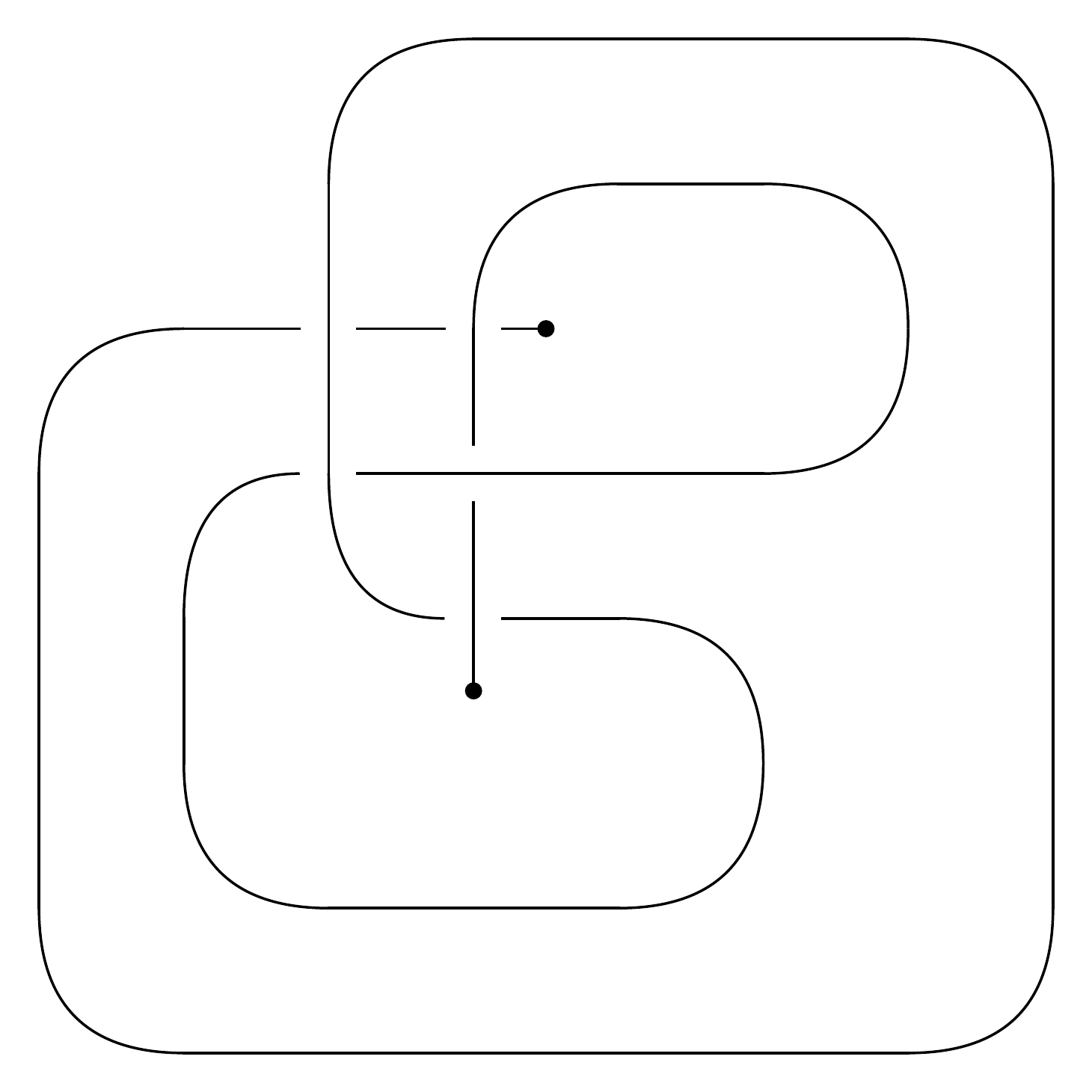}\\
\textcolor{black}{$5_{103}$}
\vspace{1cm}
\end{minipage}
\begin{minipage}[t]{.25\linewidth}
\centering
\includegraphics[width=0.9\textwidth,height=3.5cm,keepaspectratio]{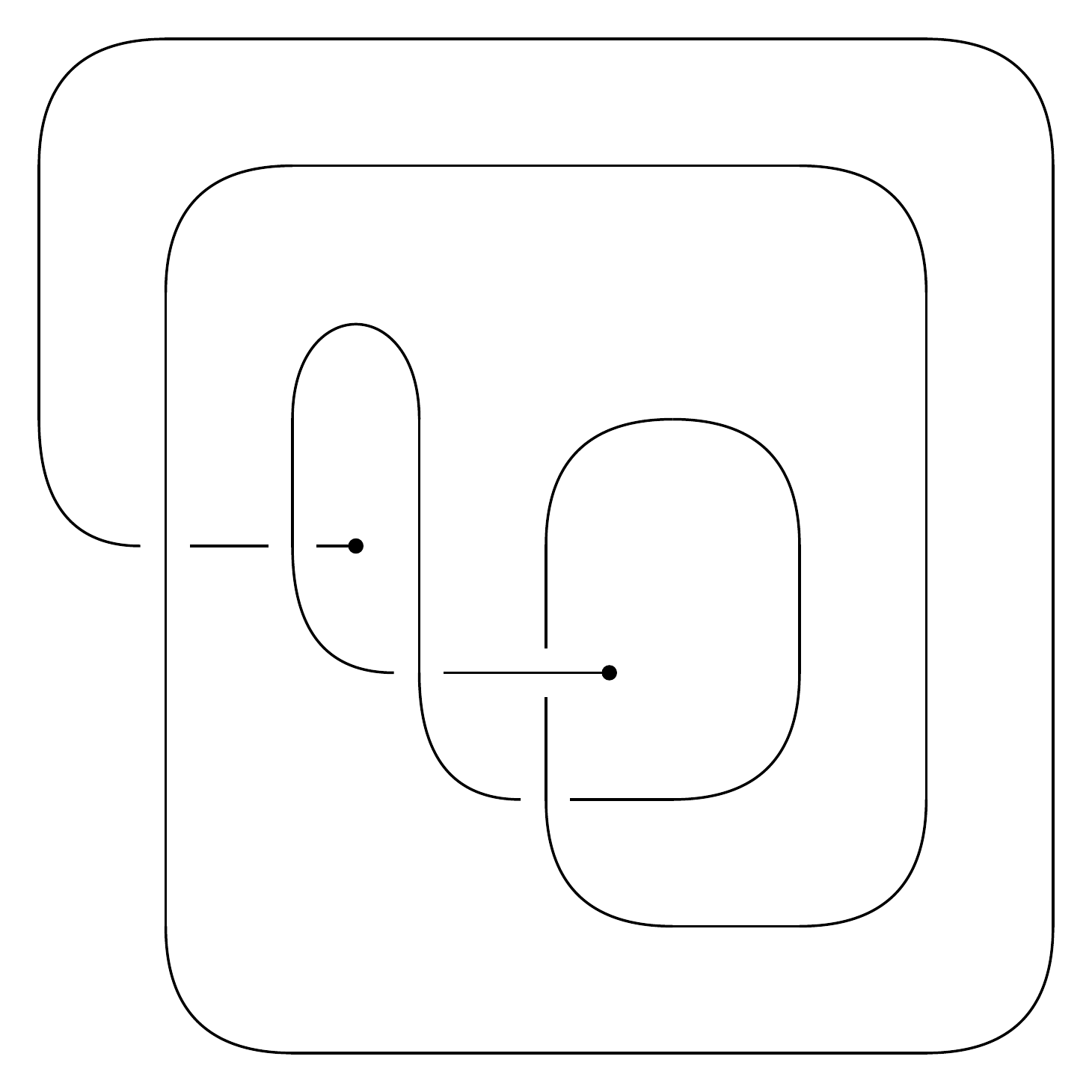}\\
\textcolor{black}{$5_{104}$}
\vspace{1cm}
\end{minipage}
\begin{minipage}[t]{.25\linewidth}
\centering
\includegraphics[width=0.9\textwidth,height=3.5cm,keepaspectratio]{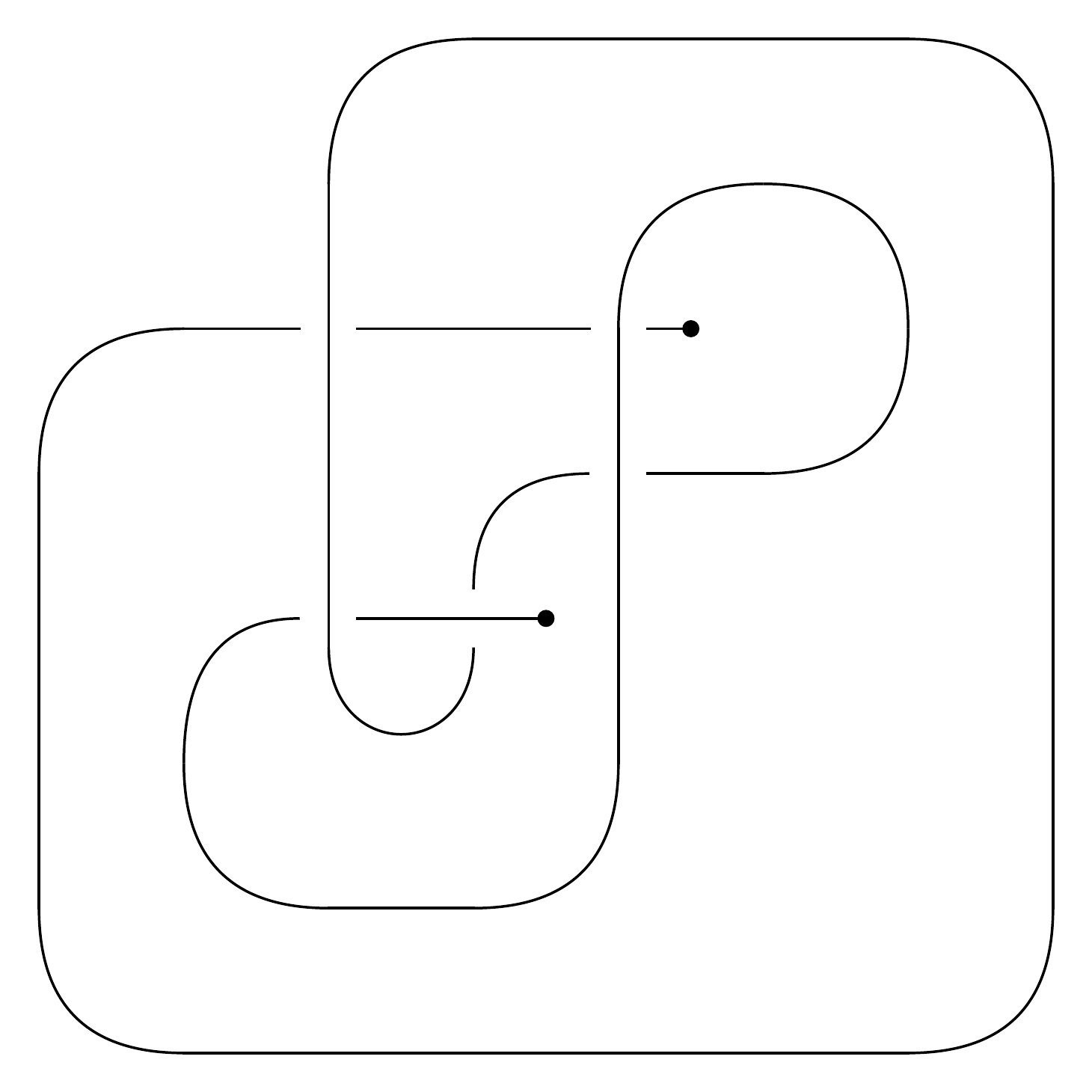}\\
\textcolor{black}{$5_{105}$}
\vspace{1cm}
\end{minipage}
\begin{minipage}[t]{.25\linewidth}
\centering
\includegraphics[width=0.9\textwidth,height=3.5cm,keepaspectratio]{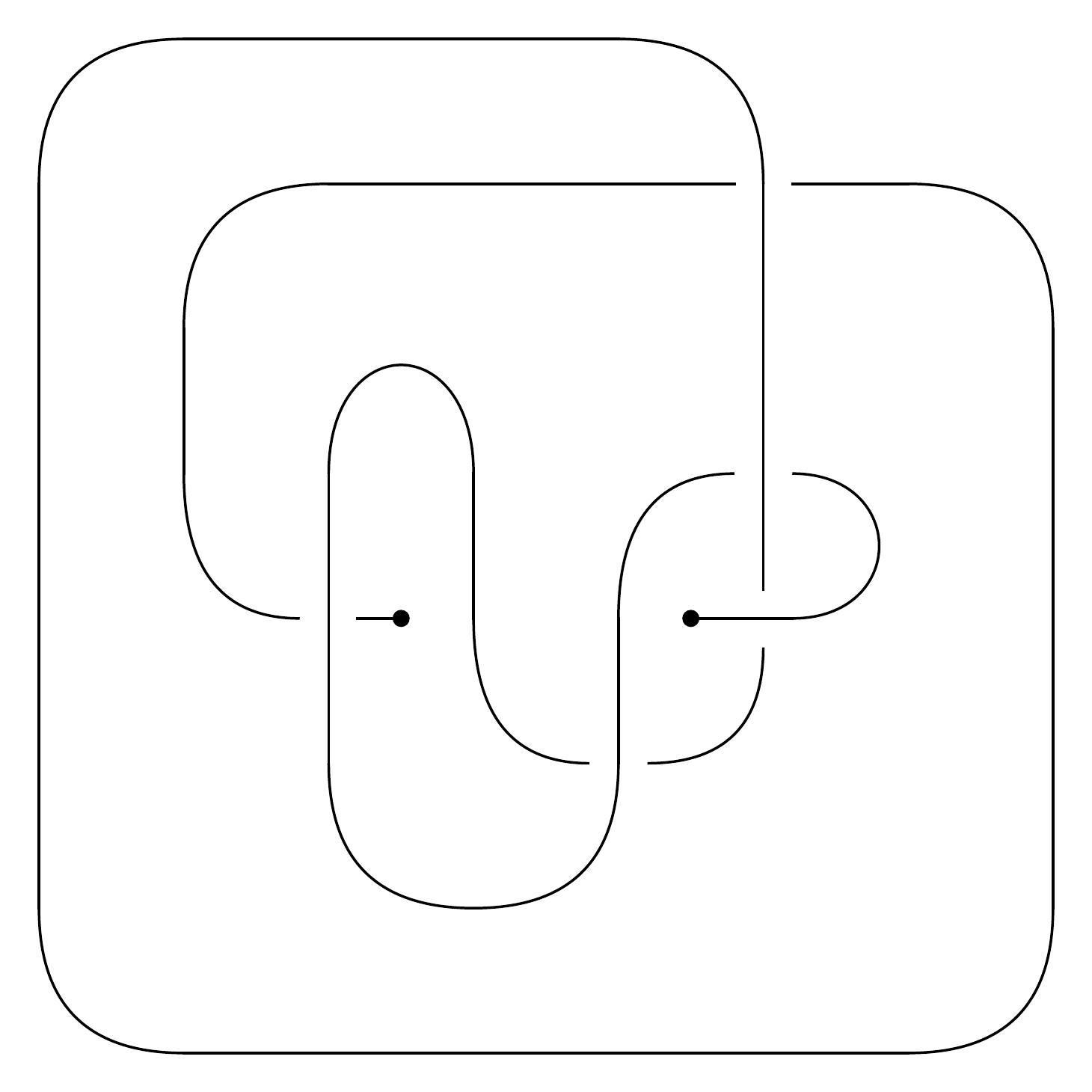}\\
\textcolor{black}{$5_{106}$}
\vspace{1cm}
\end{minipage}
\begin{minipage}[t]{.25\linewidth}
\centering
\includegraphics[width=0.9\textwidth,height=3.5cm,keepaspectratio]{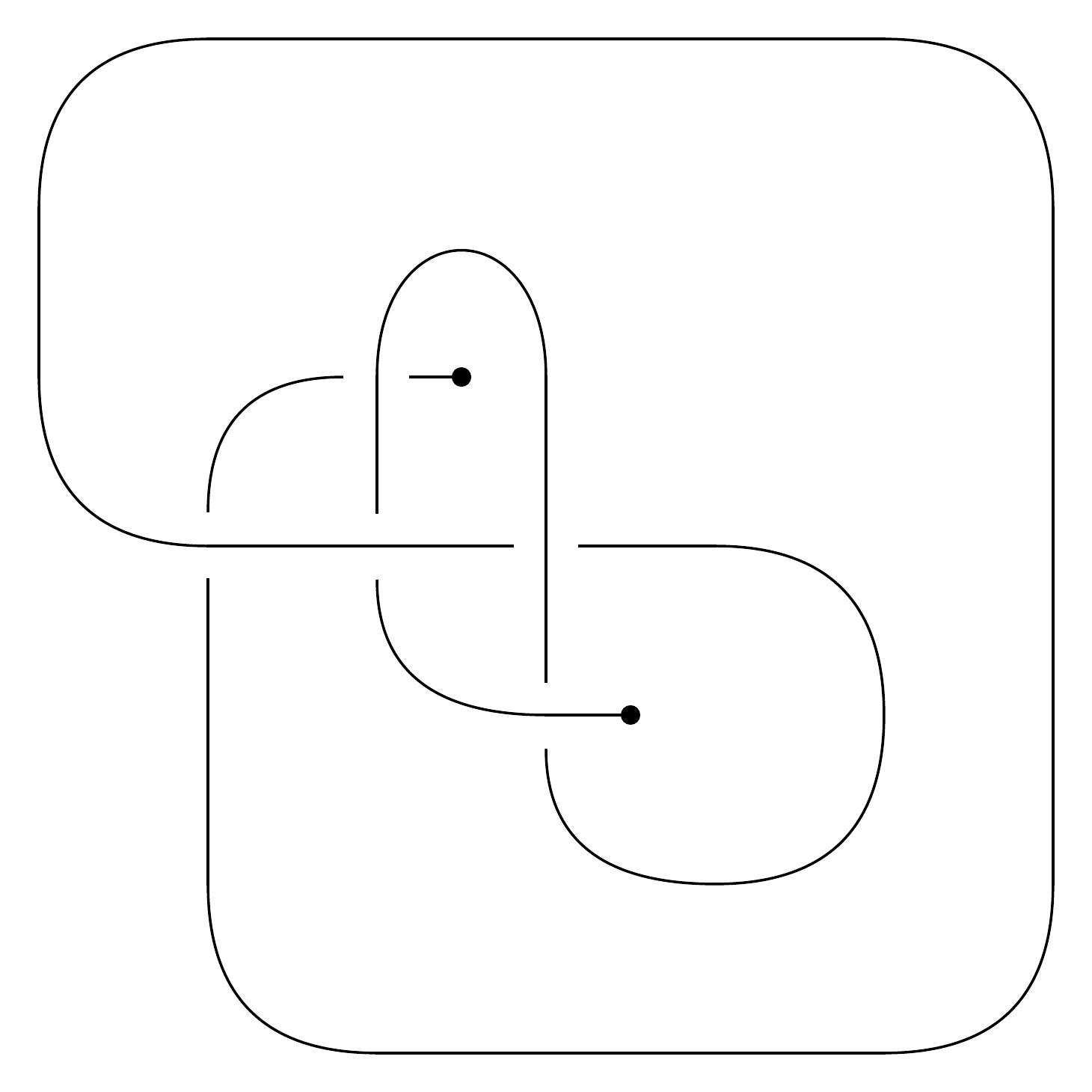}\\
\textcolor{black}{$5_{107}$}
\vspace{1cm}
\end{minipage}
\begin{minipage}[t]{.25\linewidth}
\centering
\includegraphics[width=0.9\textwidth,height=3.5cm,keepaspectratio]{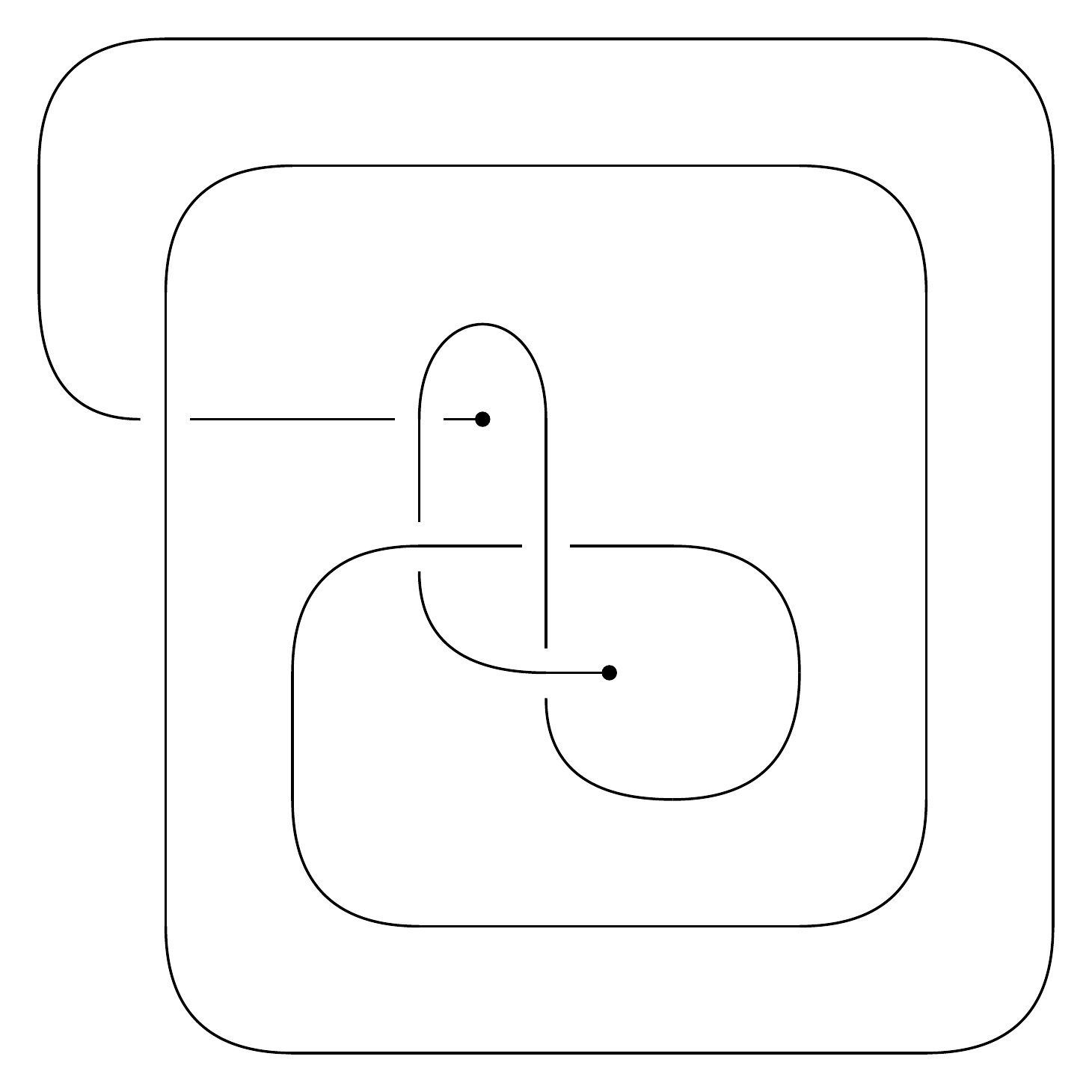}\\
\textcolor{black}{$5_{108}$}
\vspace{1cm}
\end{minipage}
\begin{minipage}[t]{.25\linewidth}
\centering
\includegraphics[width=0.9\textwidth,height=3.5cm,keepaspectratio]{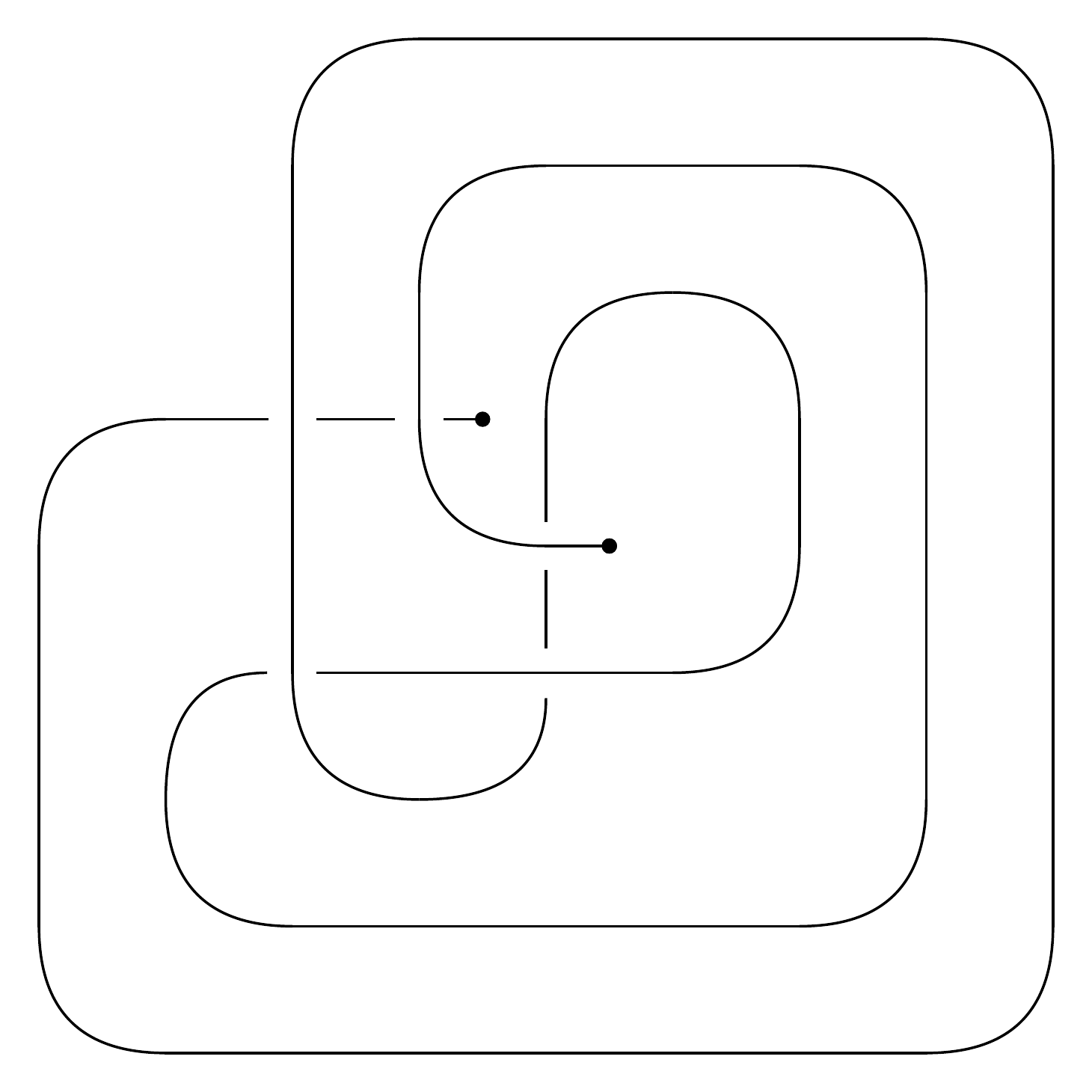}\\
\textcolor{black}{$5_{109}$}
\vspace{1cm}
\end{minipage}
\begin{minipage}[t]{.25\linewidth}
\centering
\includegraphics[width=0.9\textwidth,height=3.5cm,keepaspectratio]{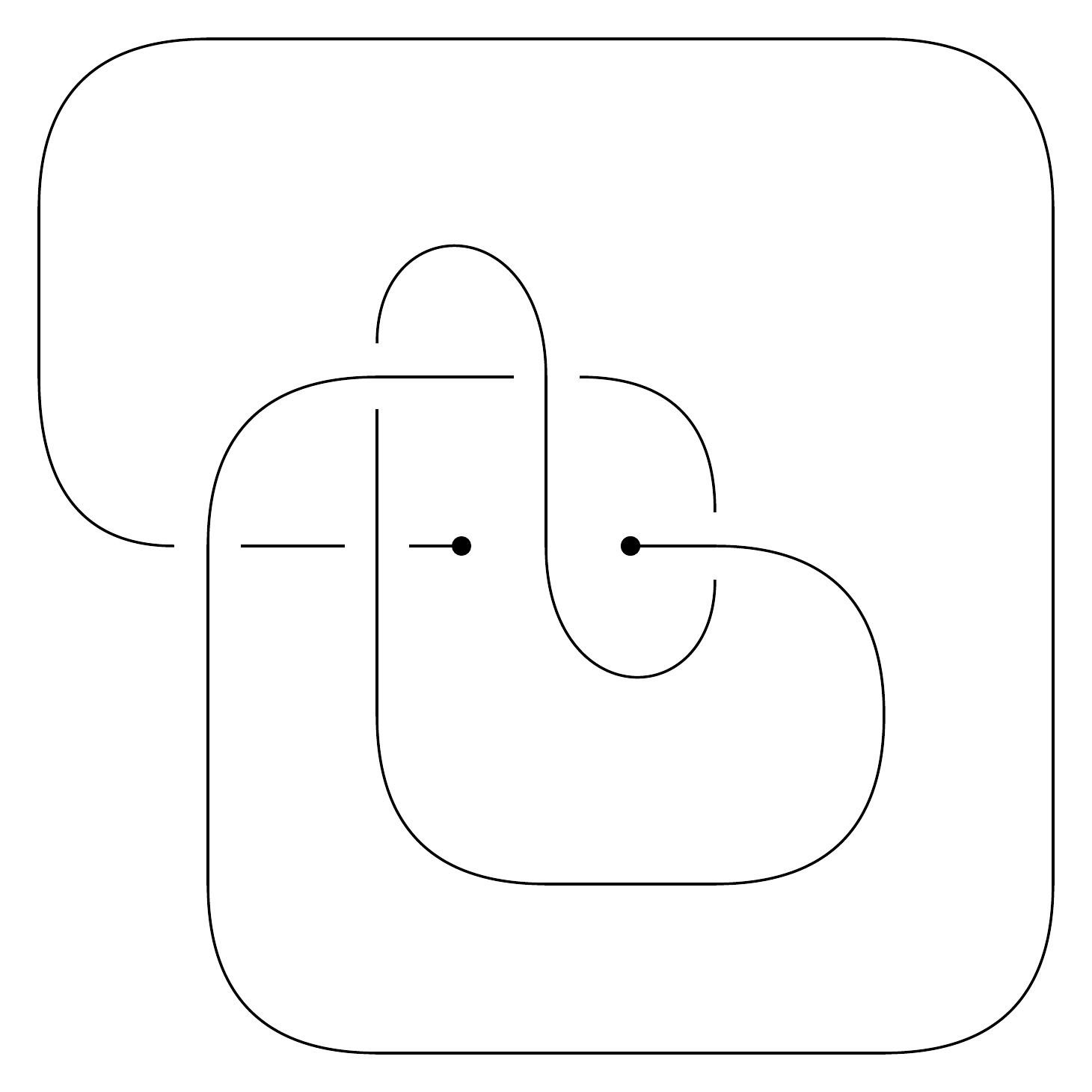}\\
\textcolor{black}{$5_{110}$}
\vspace{1cm}
\end{minipage}
\begin{minipage}[t]{.25\linewidth}
\centering
\includegraphics[width=0.9\textwidth,height=3.5cm,keepaspectratio]{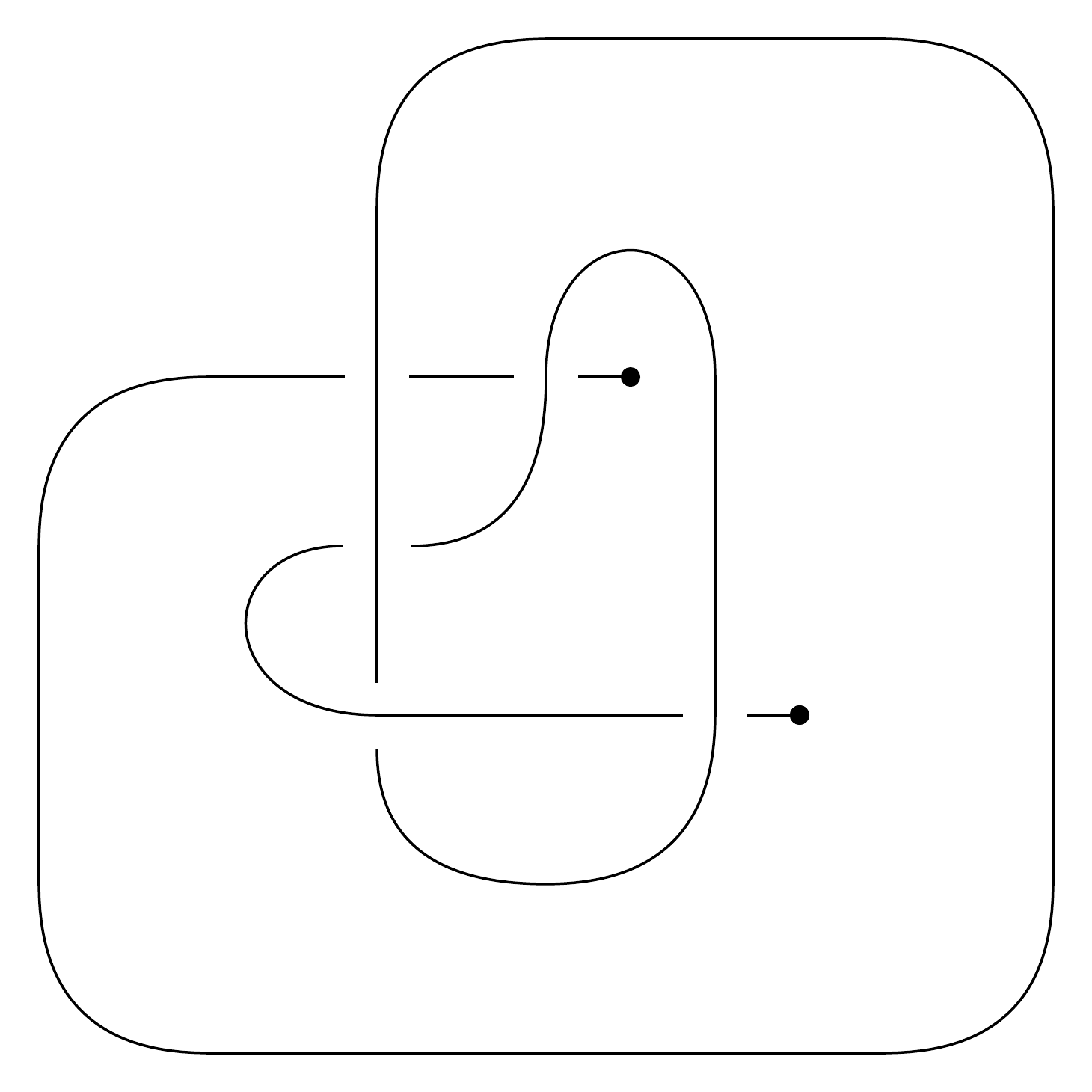}\\
\textcolor{black}{$5_{111}$}
\vspace{1cm}
\end{minipage}
\begin{minipage}[t]{.25\linewidth}
\centering
\includegraphics[width=0.9\textwidth,height=3.5cm,keepaspectratio]{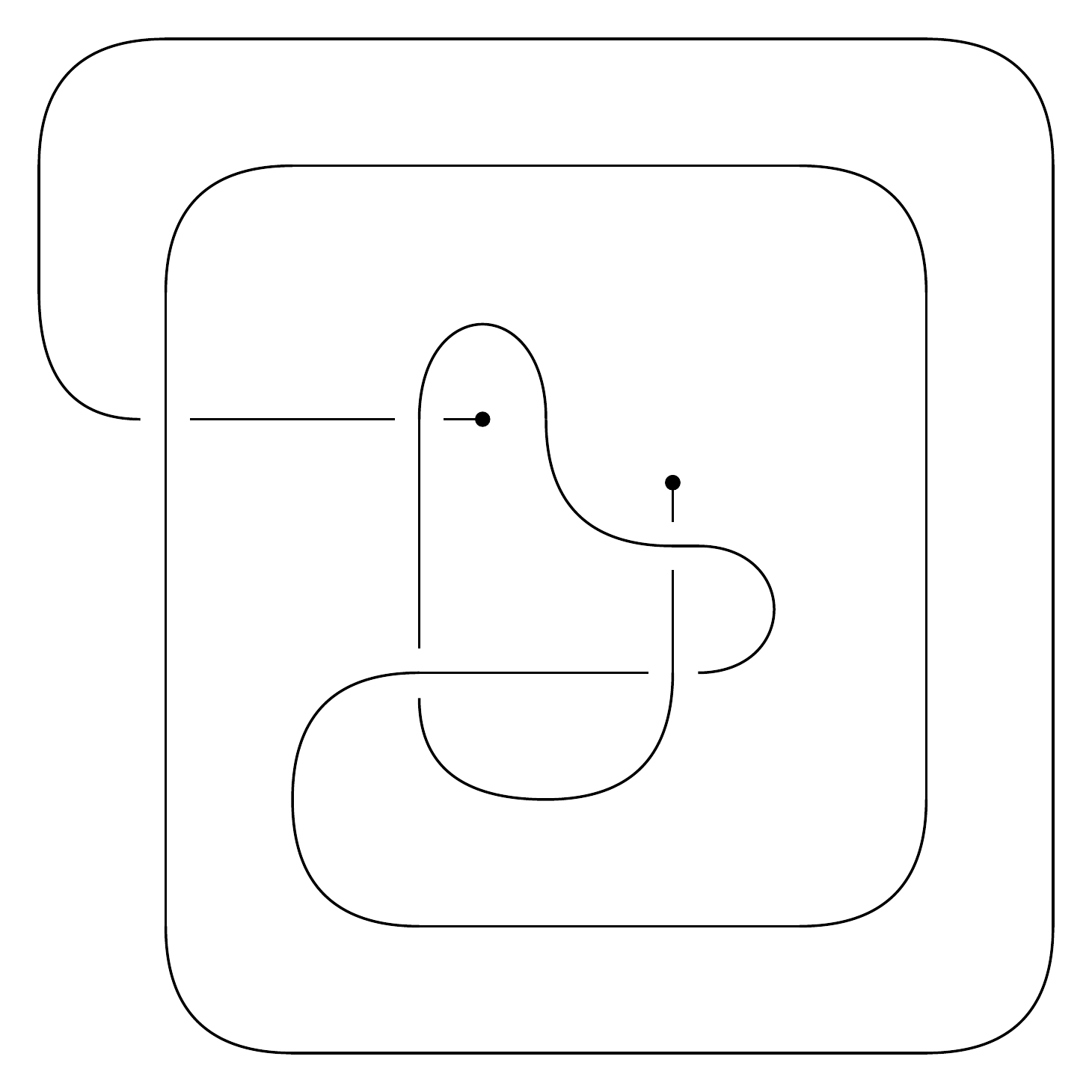}\\
\textcolor{black}{$5_{112}$}
\vspace{1cm}
\end{minipage}
\begin{minipage}[t]{.25\linewidth}
\centering
\includegraphics[width=0.9\textwidth,height=3.5cm,keepaspectratio]{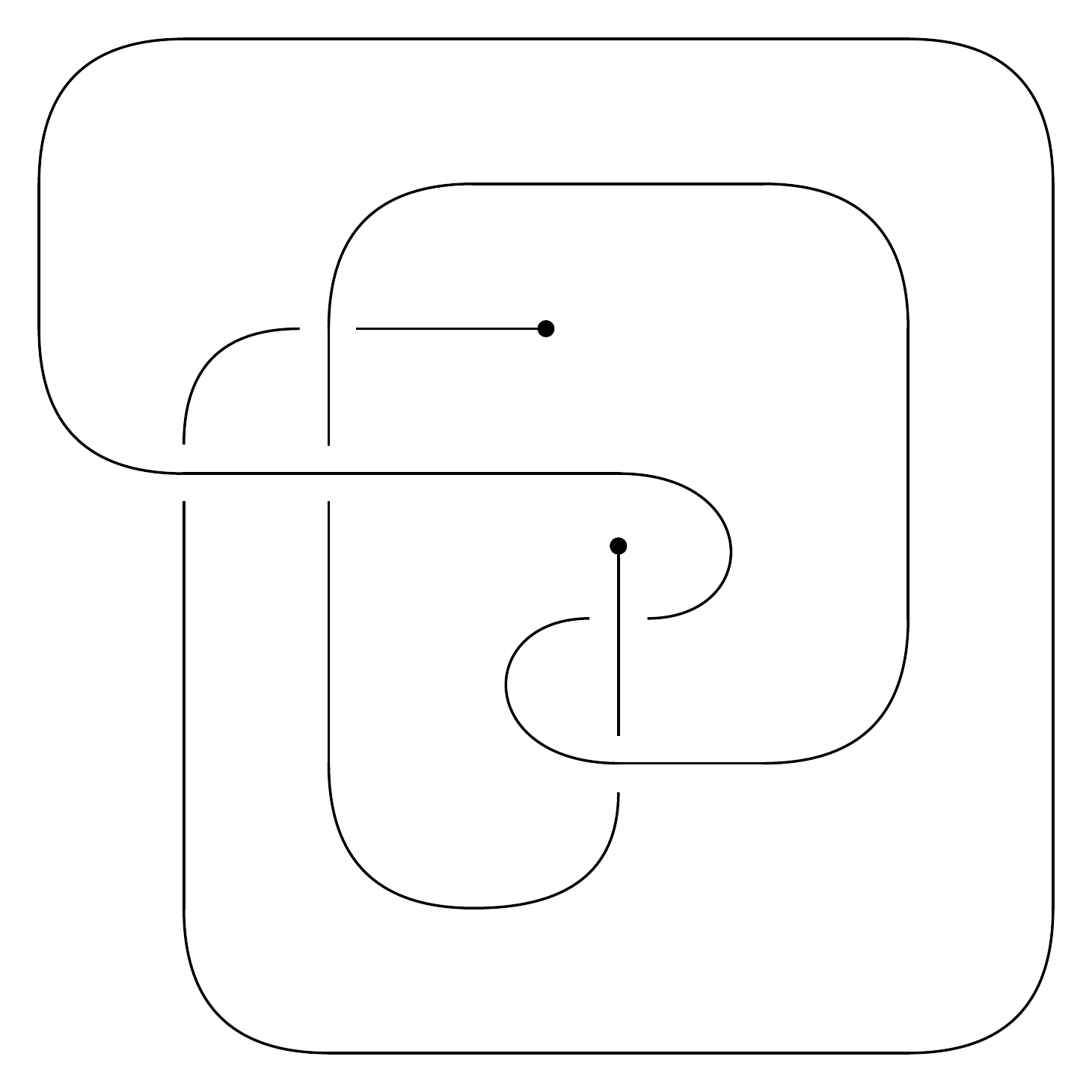}\\
\textcolor{black}{$5_{113}$}
\vspace{1cm}
\end{minipage}
\begin{minipage}[t]{.25\linewidth}
\centering
\includegraphics[width=0.9\textwidth,height=3.5cm,keepaspectratio]{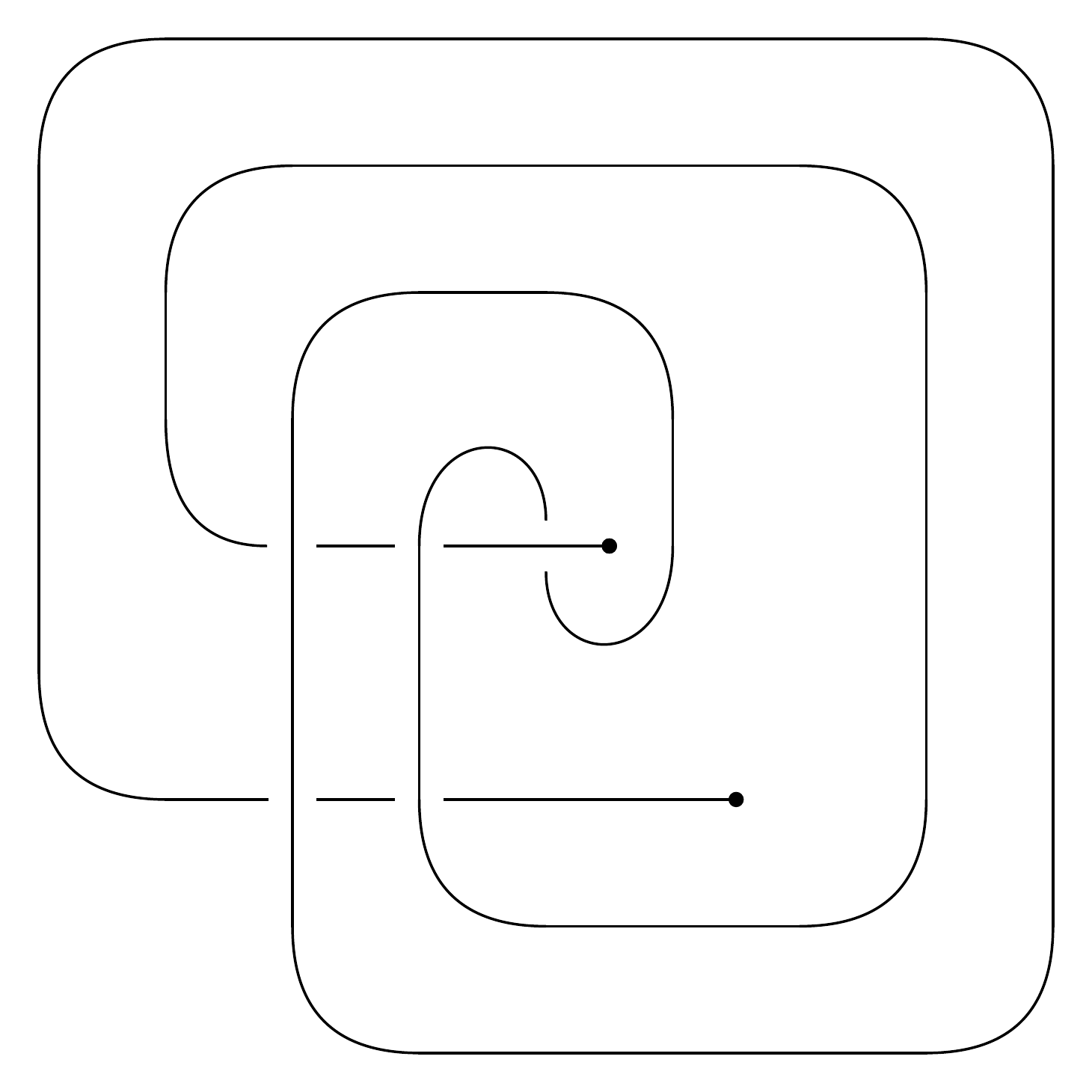}\\
\textcolor{black}{$5_{114}$}
\vspace{1cm}
\end{minipage}
\begin{minipage}[t]{.25\linewidth}
\centering
\includegraphics[width=0.9\textwidth,height=3.5cm,keepaspectratio]{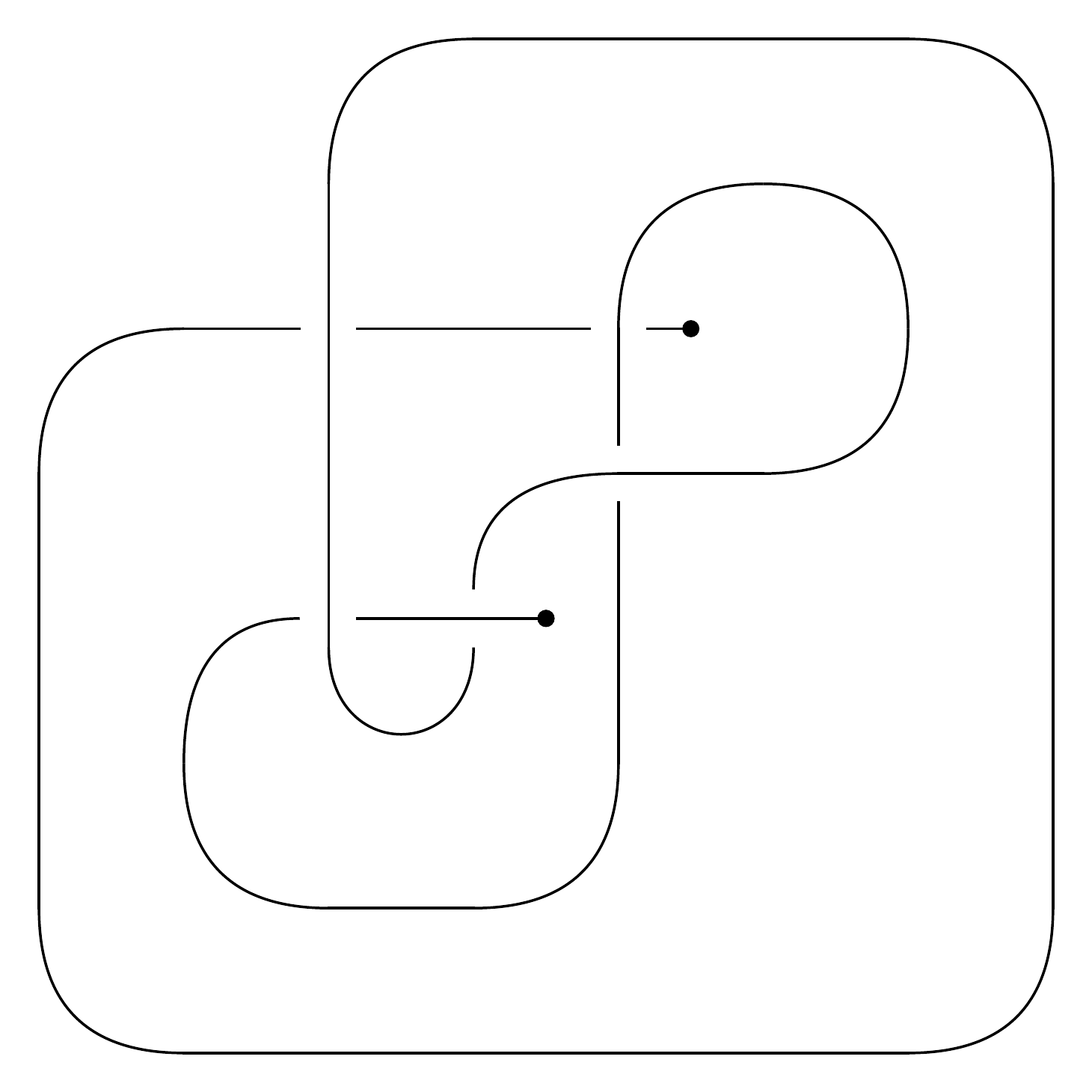}\\
\textcolor{black}{$5_{115}$}
\vspace{1cm}
\end{minipage}
\begin{minipage}[t]{.25\linewidth}
\centering
\includegraphics[width=0.9\textwidth,height=3.5cm,keepaspectratio]{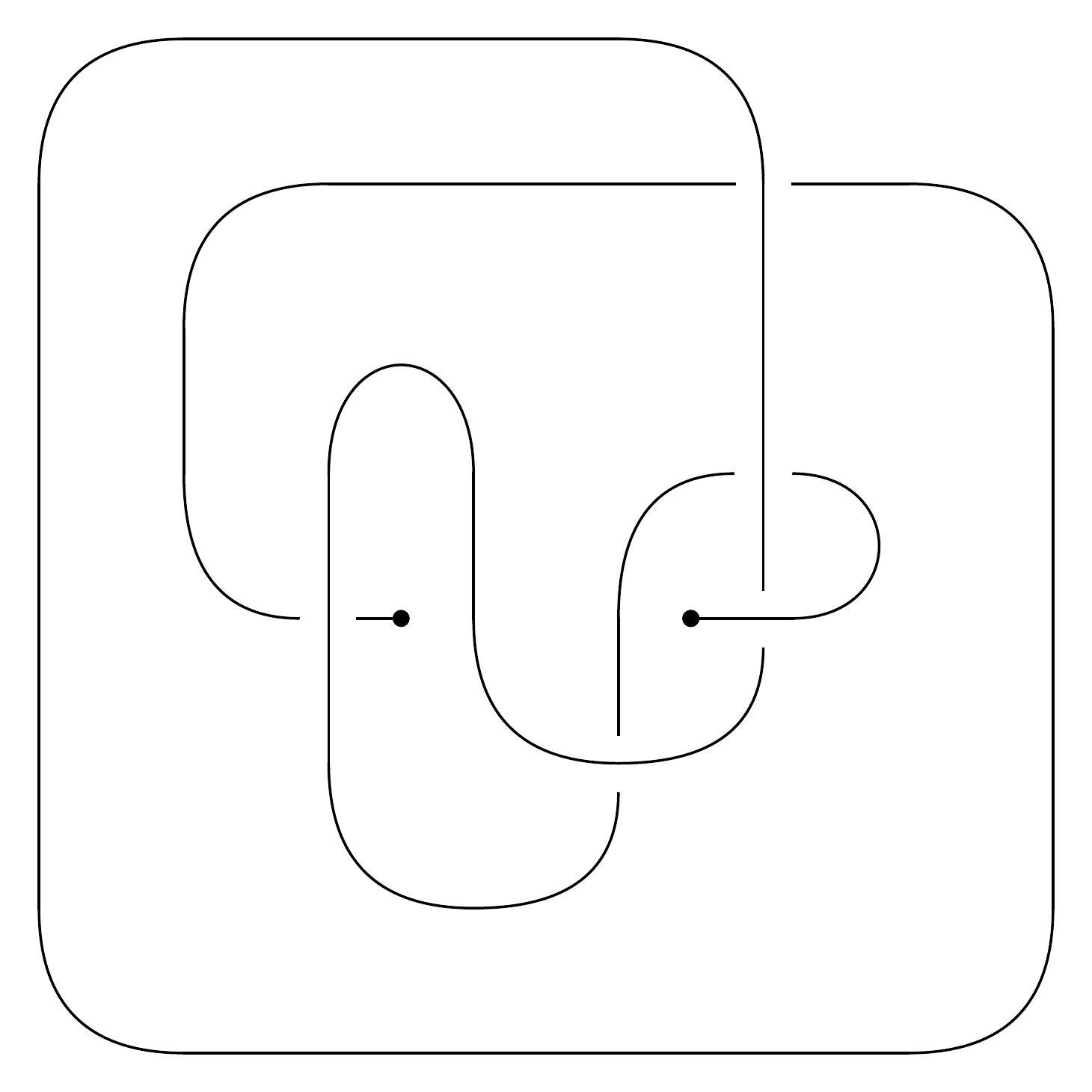}\\
\textcolor{black}{$5_{116}$}
\vspace{1cm}
\end{minipage}
\begin{minipage}[t]{.25\linewidth}
\centering
\includegraphics[width=0.9\textwidth,height=3.5cm,keepaspectratio]{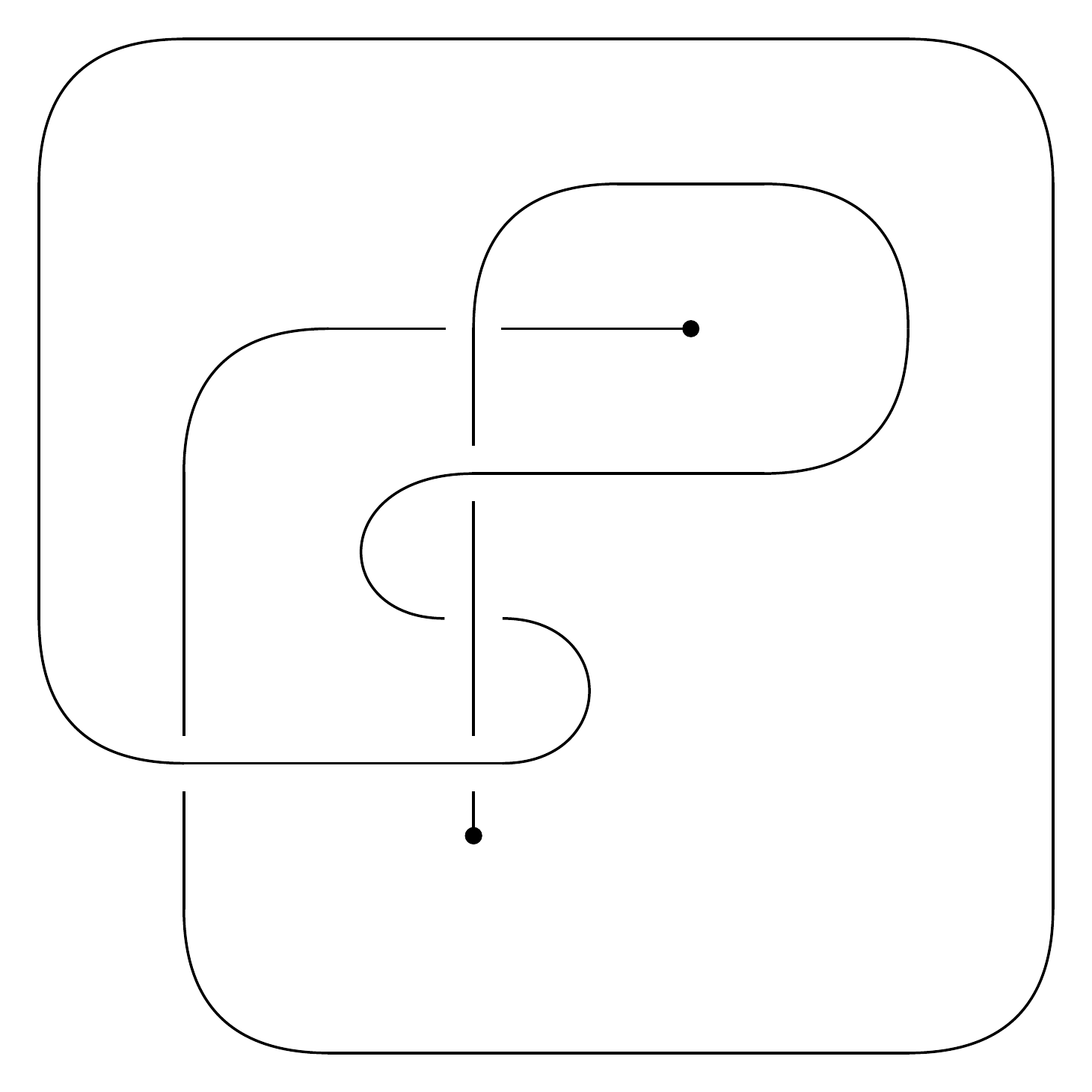}\\
\textcolor{black}{$5_{117}$}
\vspace{1cm}
\end{minipage}
\begin{minipage}[t]{.25\linewidth}
\centering
\includegraphics[width=0.9\textwidth,height=3.5cm,keepaspectratio]{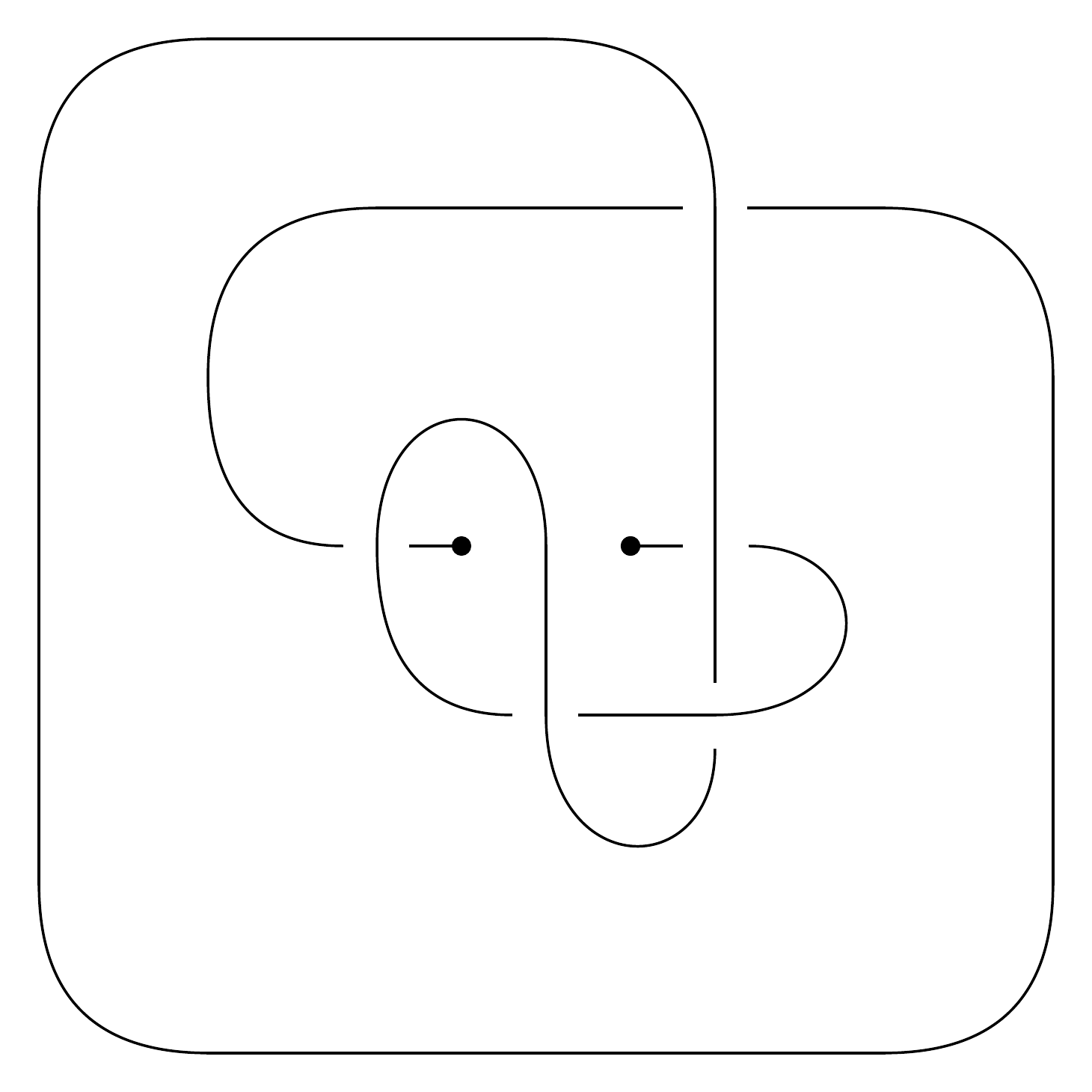}\\
\textcolor{black}{$5_{118}$}
\vspace{1cm}
\end{minipage}
\begin{minipage}[t]{.25\linewidth}
\centering
\includegraphics[width=0.9\textwidth,height=3.5cm,keepaspectratio]{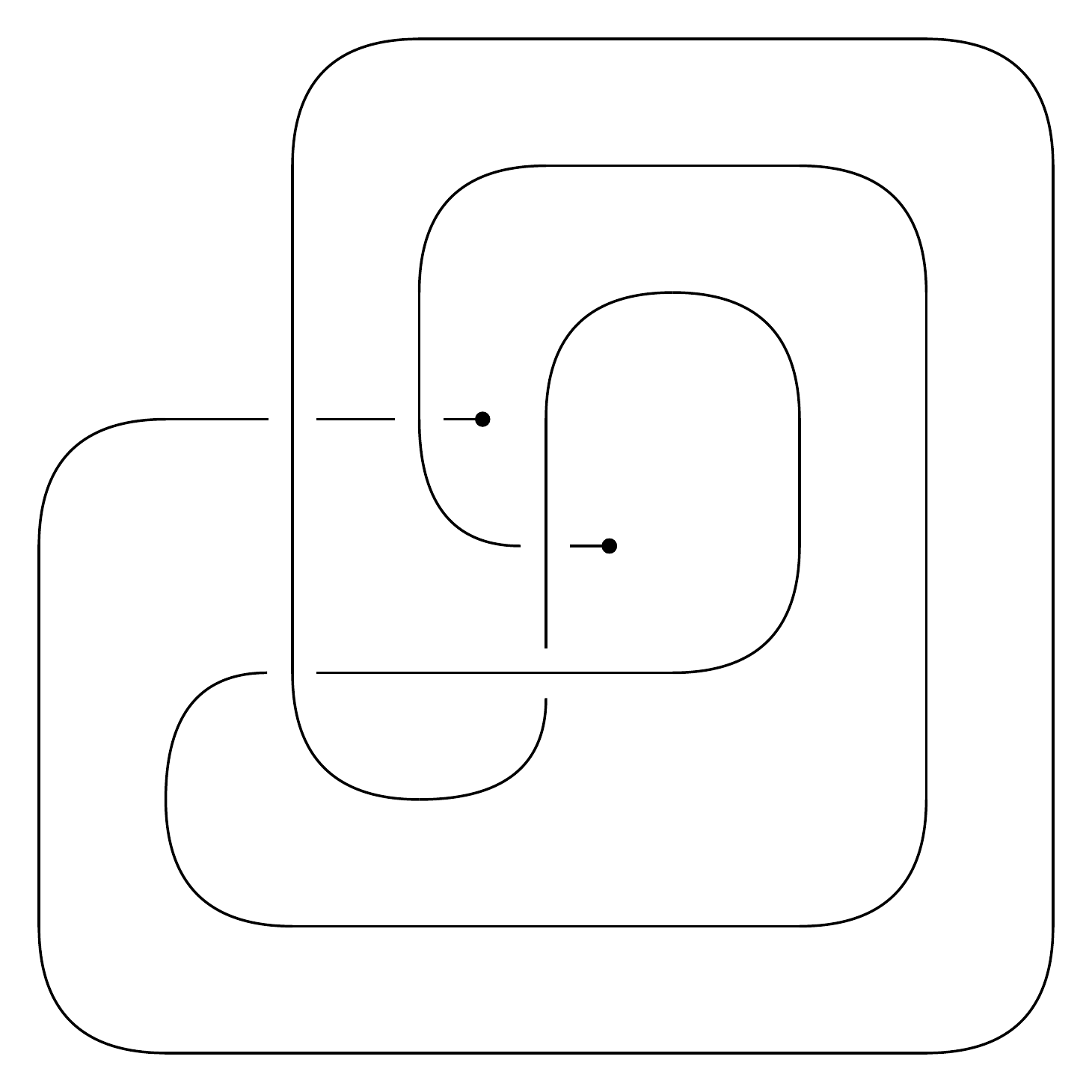}\\
\textcolor{black}{$5_{119}$}
\vspace{1cm}
\end{minipage}
\begin{minipage}[t]{.25\linewidth}
\centering
\includegraphics[width=0.9\textwidth,height=3.5cm,keepaspectratio]{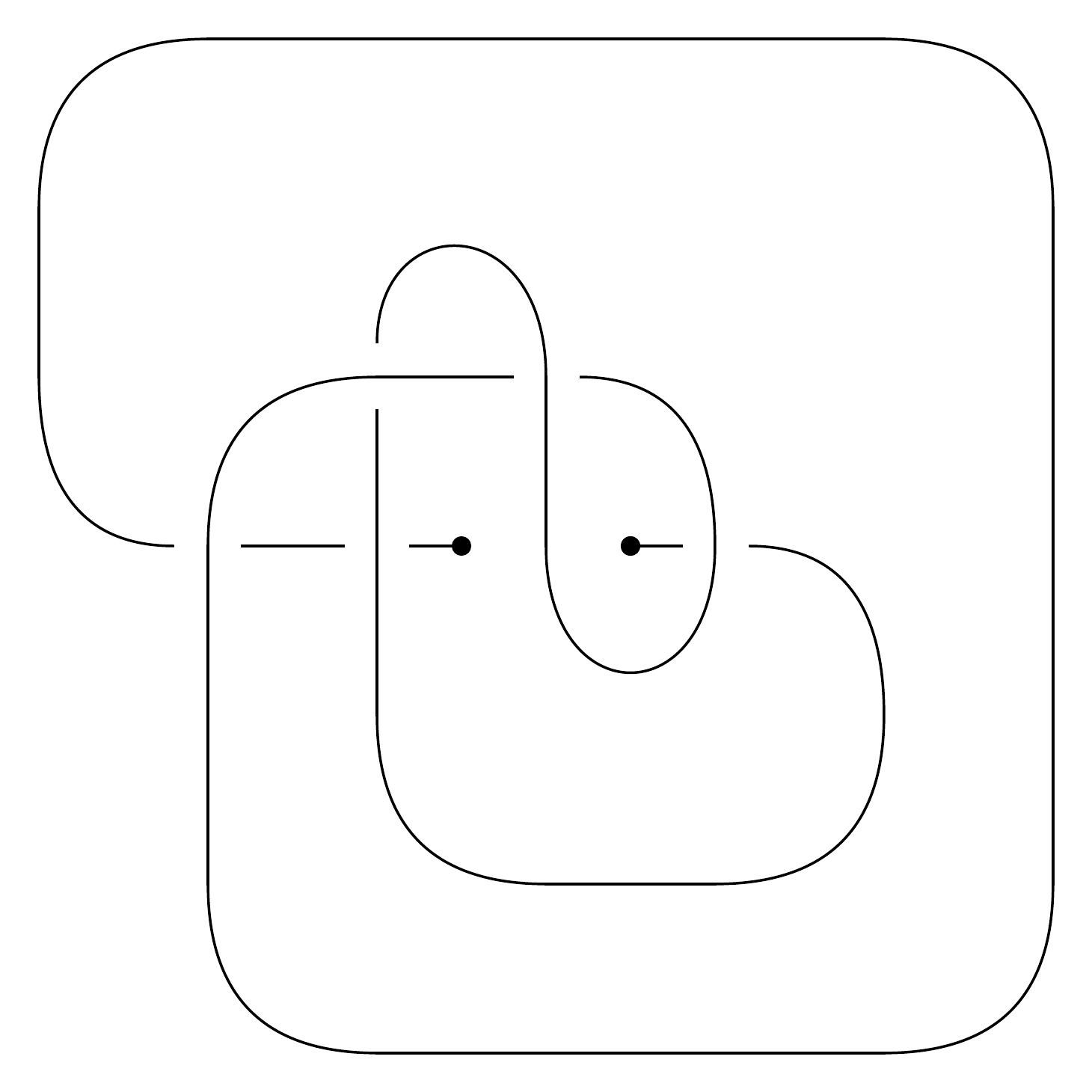}\\
\textcolor{black}{$5_{120}$}
\vspace{1cm}
\end{minipage}
\begin{minipage}[t]{.25\linewidth}
\centering
\includegraphics[width=0.9\textwidth,height=3.5cm,keepaspectratio]{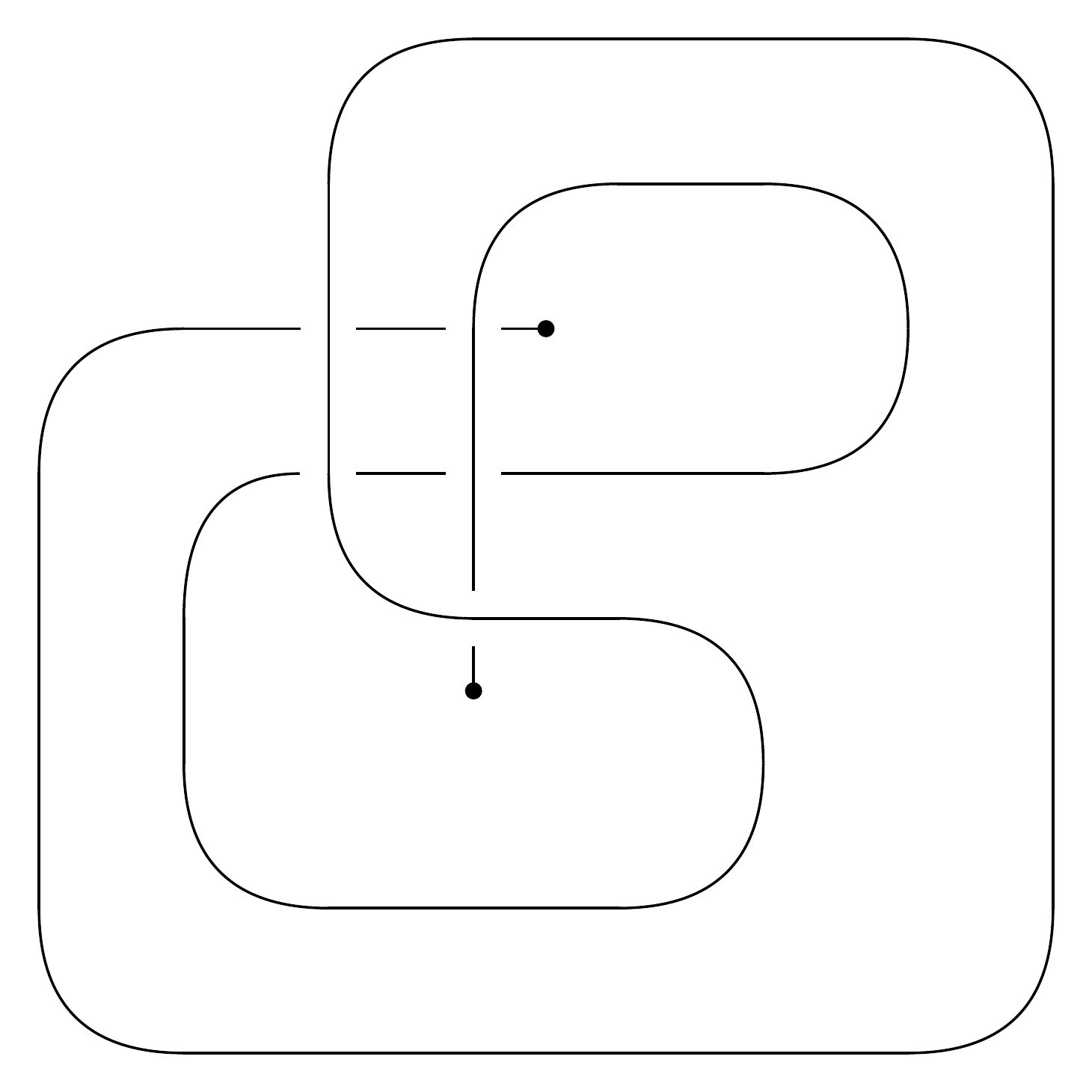}\\
\textcolor{black}{$5_{121}$}
\vspace{1cm}
\end{minipage}
\begin{minipage}[t]{.25\linewidth}
\centering
\includegraphics[width=0.9\textwidth,height=3.5cm,keepaspectratio]{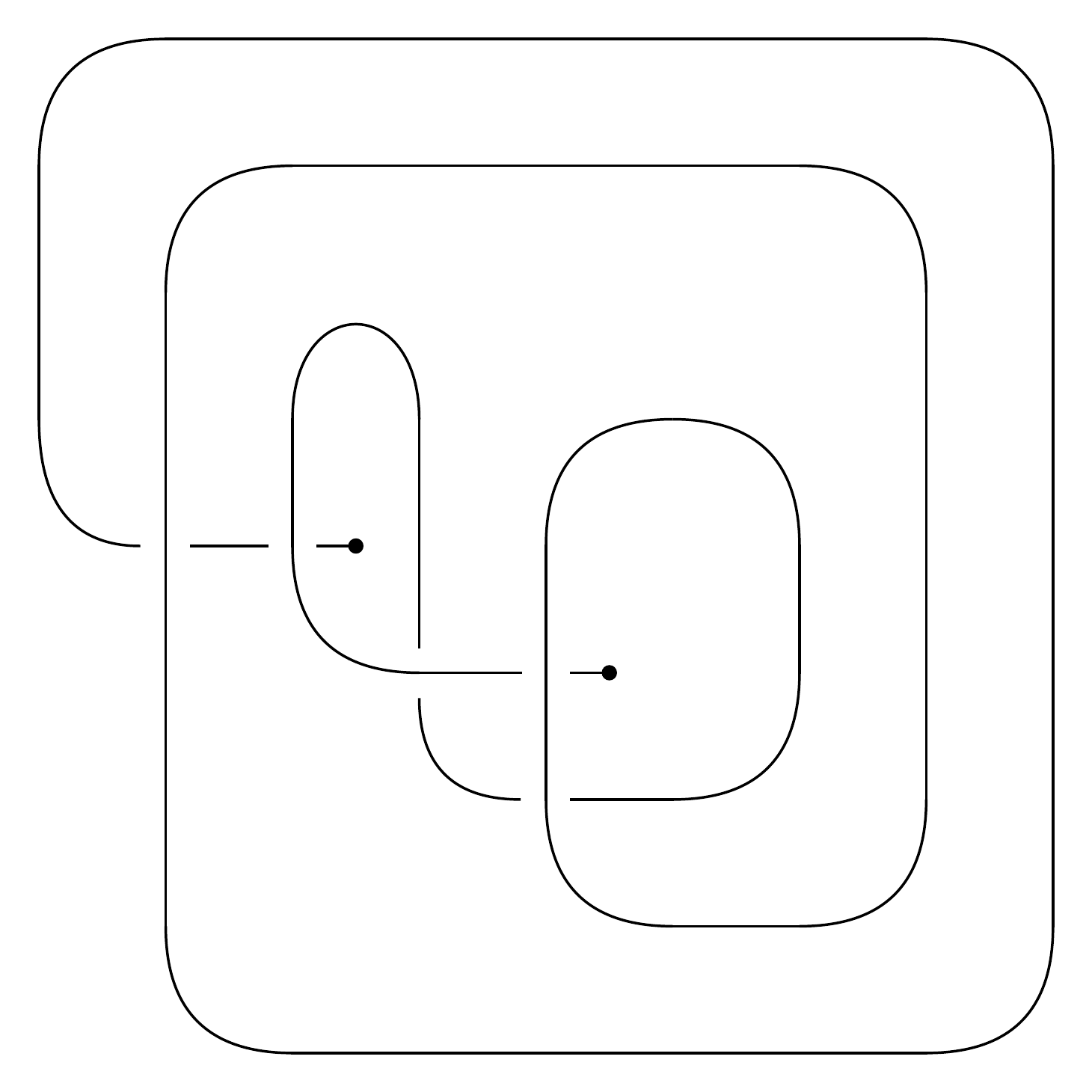}\\
\textcolor{black}{$5_{122}$}
\vspace{1cm}
\end{minipage}
\begin{minipage}[t]{.25\linewidth}
\centering
\includegraphics[width=0.9\textwidth,height=3.5cm,keepaspectratio]{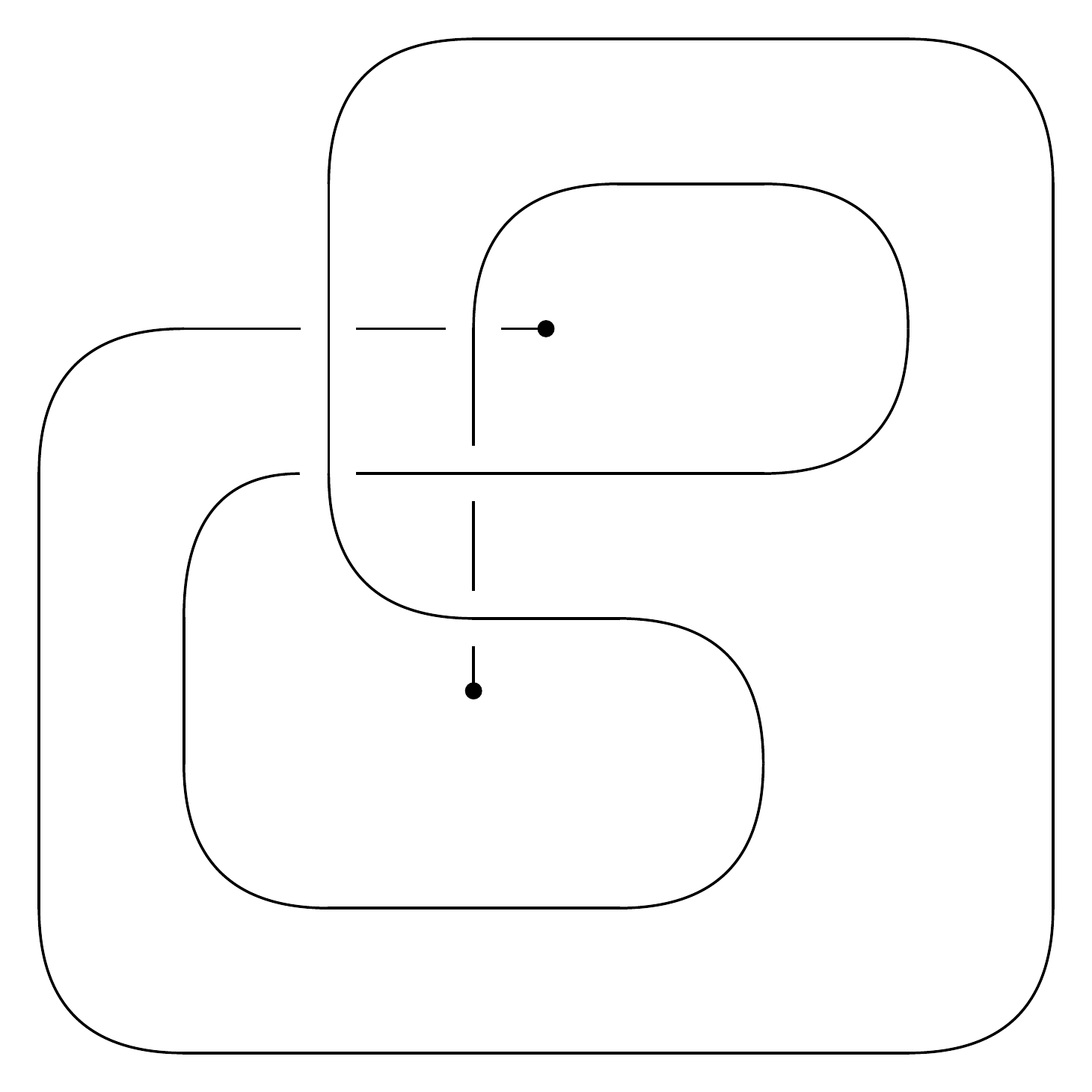}\\
\textcolor{black}{$5_{123}$}
\vspace{1cm}
\end{minipage}
\begin{minipage}[t]{.25\linewidth}
\centering
\includegraphics[width=0.9\textwidth,height=3.5cm,keepaspectratio]{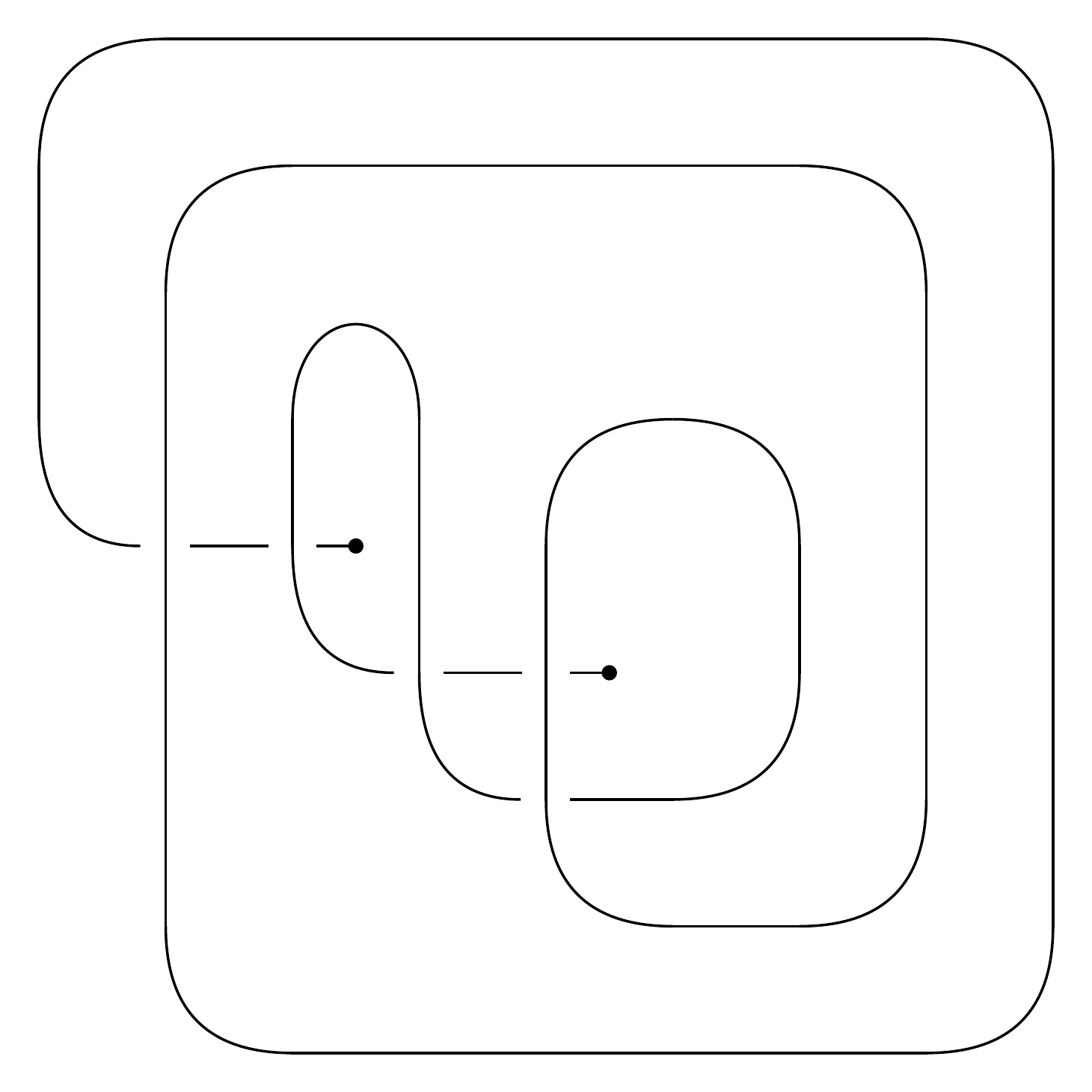}\\
\textcolor{black}{$5_{124}$}
\vspace{1cm}
\end{minipage}
\begin{minipage}[t]{.25\linewidth}
\centering
\includegraphics[width=0.9\textwidth,height=3.5cm,keepaspectratio]{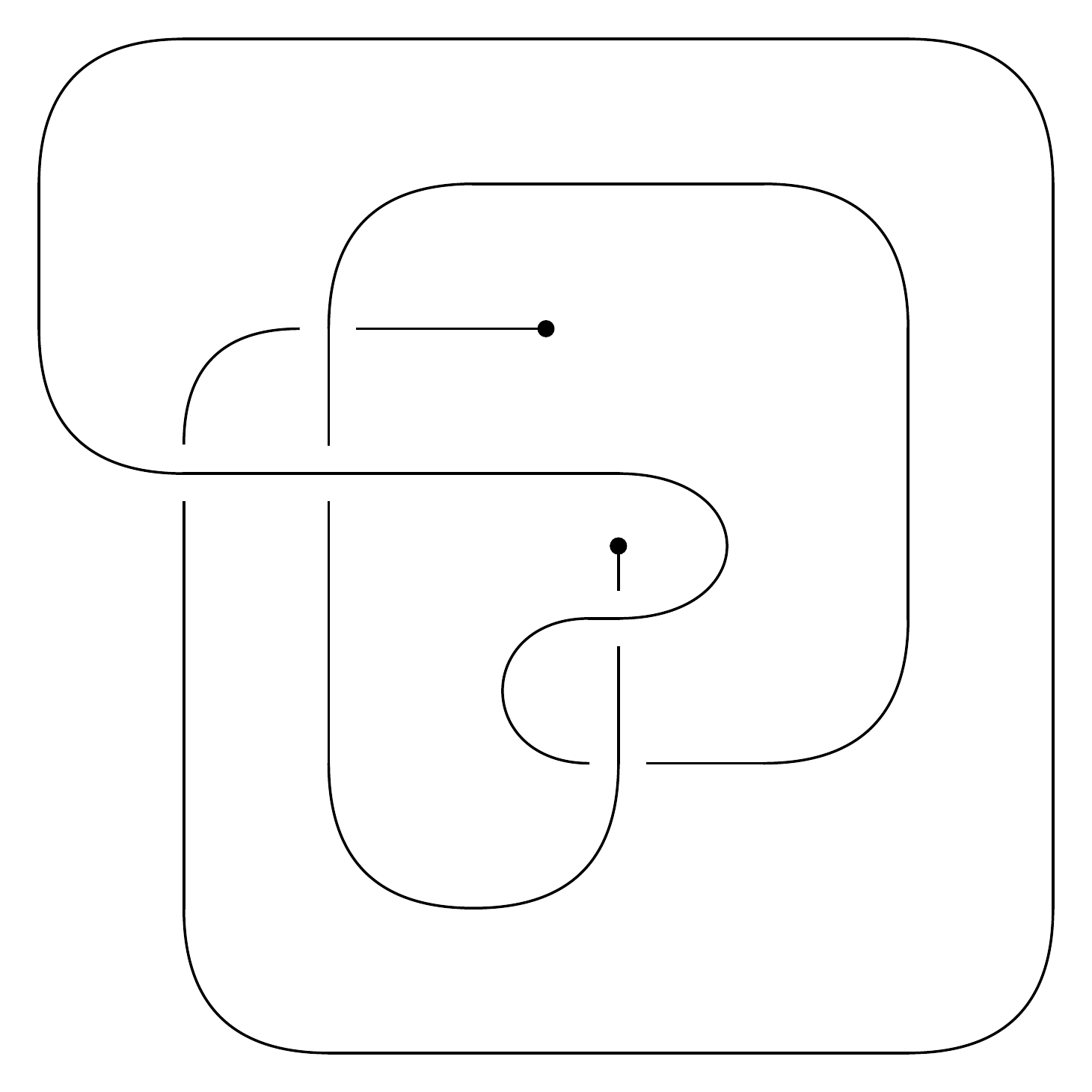}\\
\textcolor{black}{$5_{125}$}
\vspace{1cm}
\end{minipage}
\begin{minipage}[t]{.25\linewidth}
\centering
\includegraphics[width=0.9\textwidth,height=3.5cm,keepaspectratio]{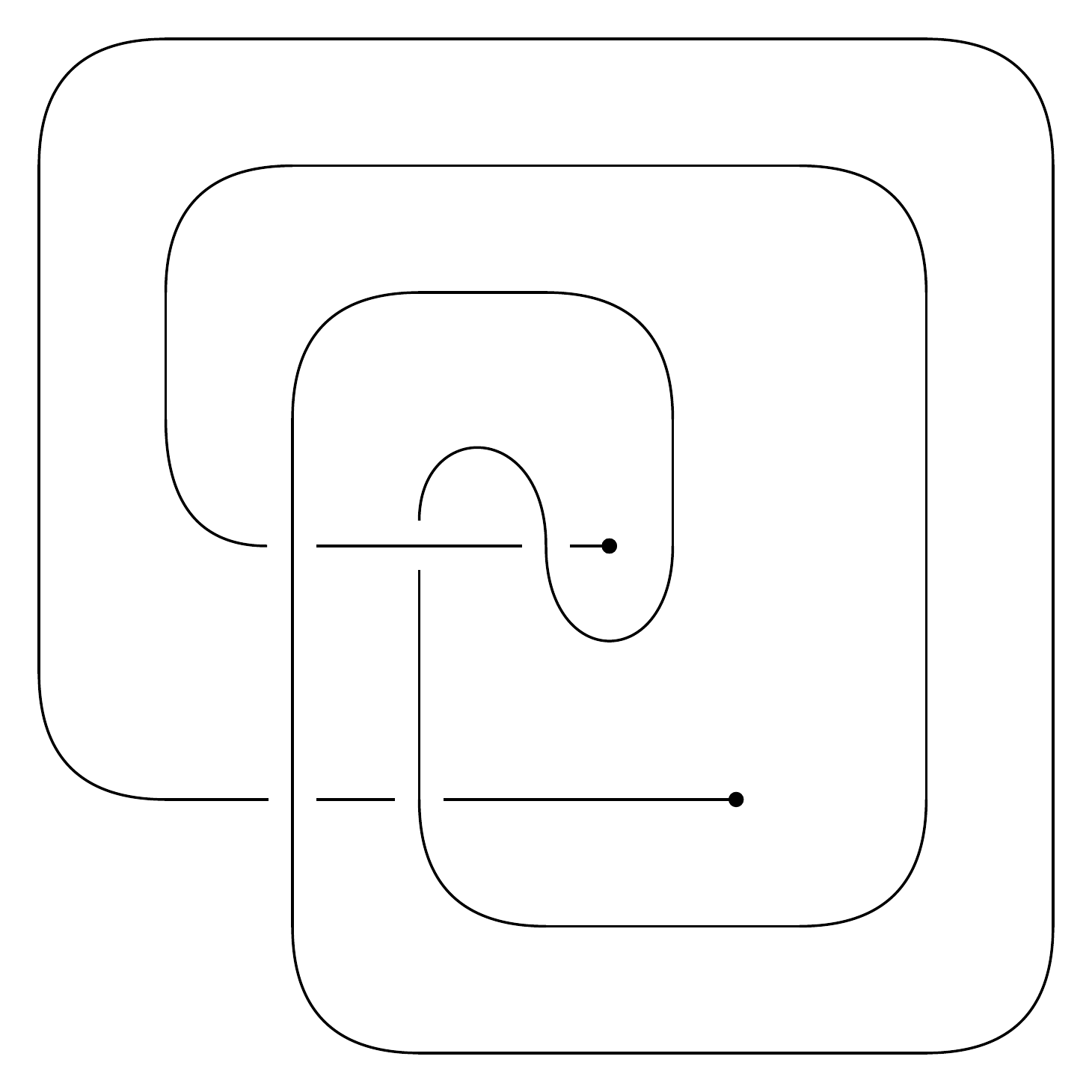}\\
\textcolor{black}{$5_{126}$}
\vspace{1cm}
\end{minipage}
\begin{minipage}[t]{.25\linewidth}
\centering
\includegraphics[width=0.9\textwidth,height=3.5cm,keepaspectratio]{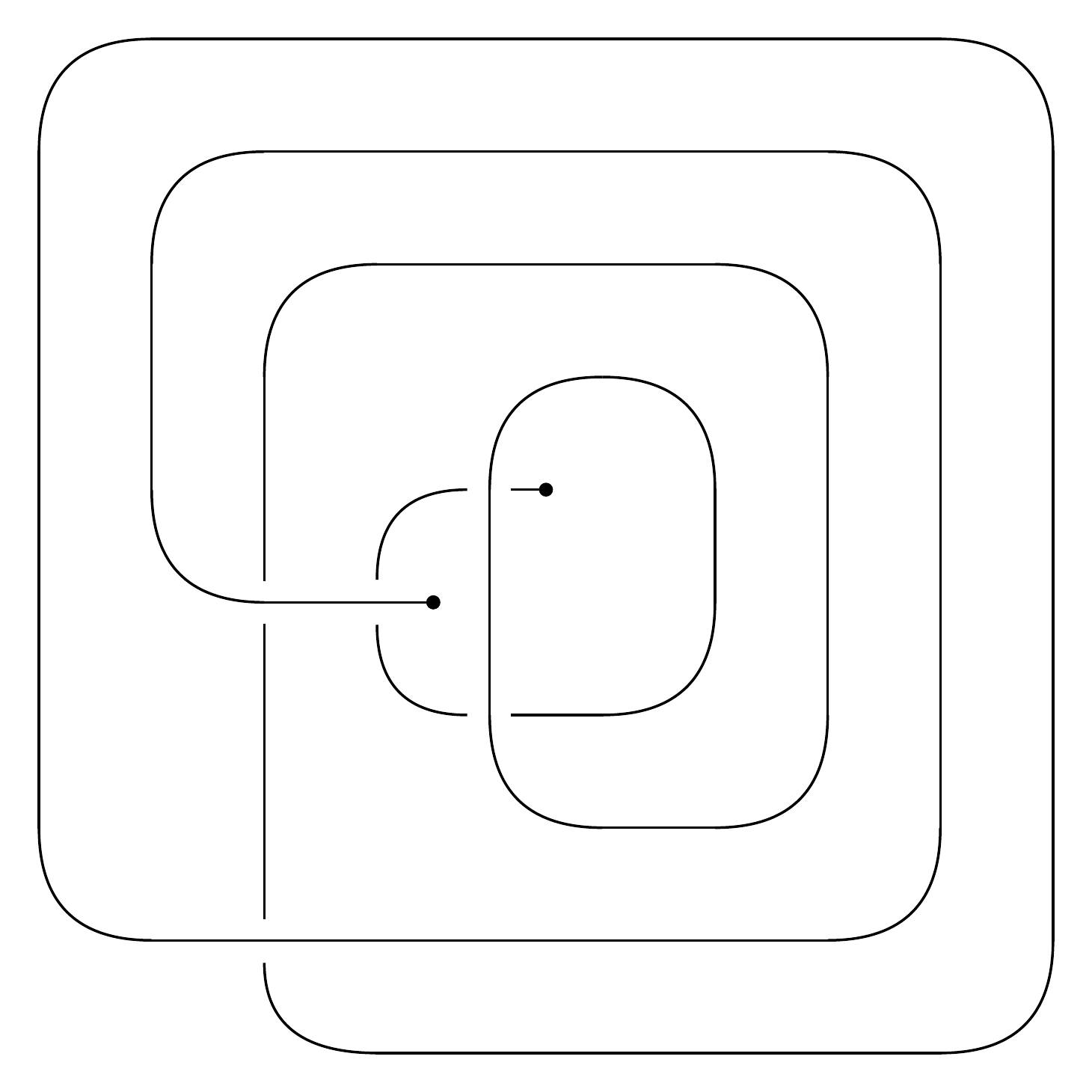}\\
\textcolor{black}{$5_{127}$}
\vspace{1cm}
\end{minipage}
\begin{minipage}[t]{.25\linewidth}
\centering
\includegraphics[width=0.9\textwidth,height=3.5cm,keepaspectratio]{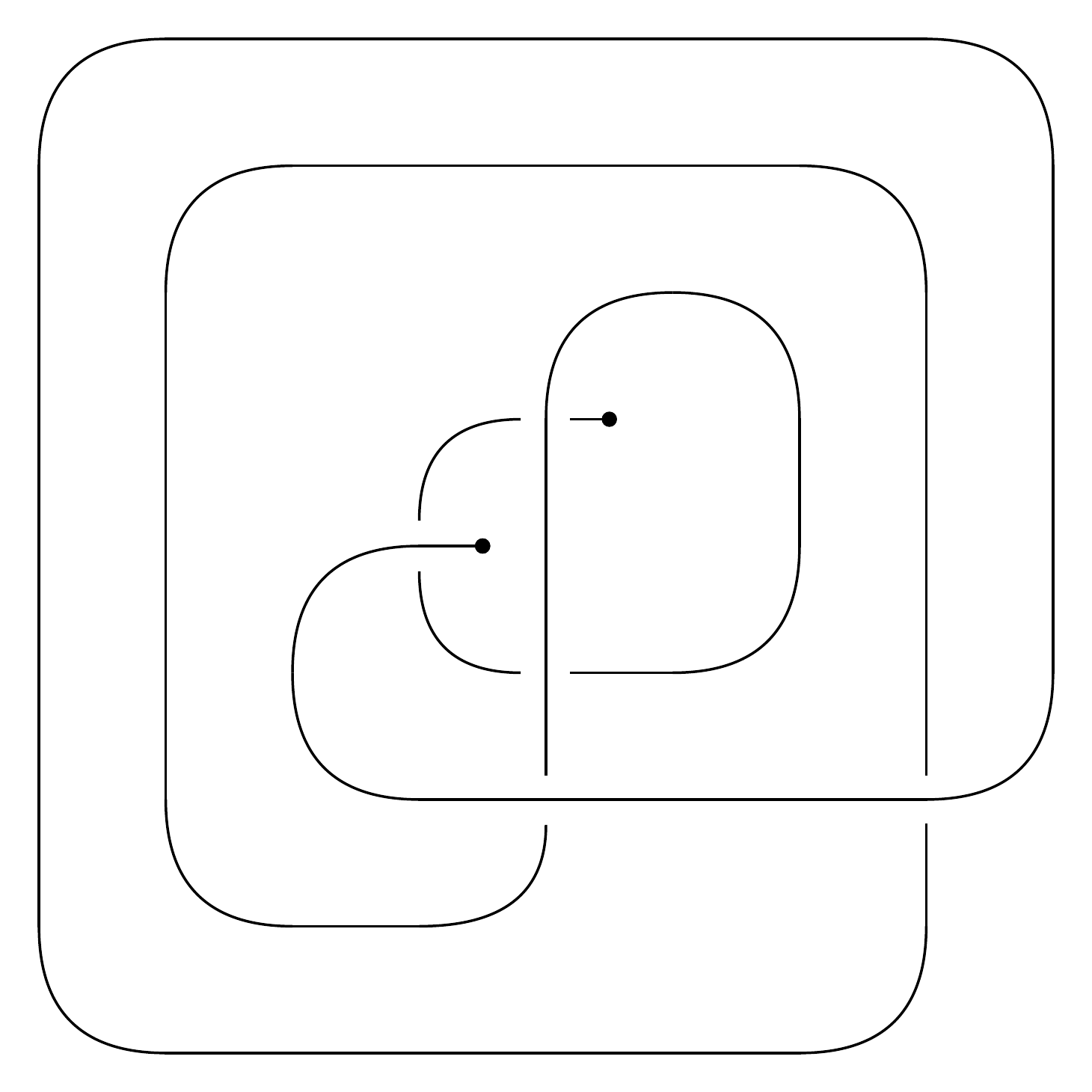}\\
\textcolor{black}{$5_{128}$}
\vspace{1cm}
\end{minipage}
\begin{minipage}[t]{.25\linewidth}
\centering
\includegraphics[width=0.9\textwidth,height=3.5cm,keepaspectratio]{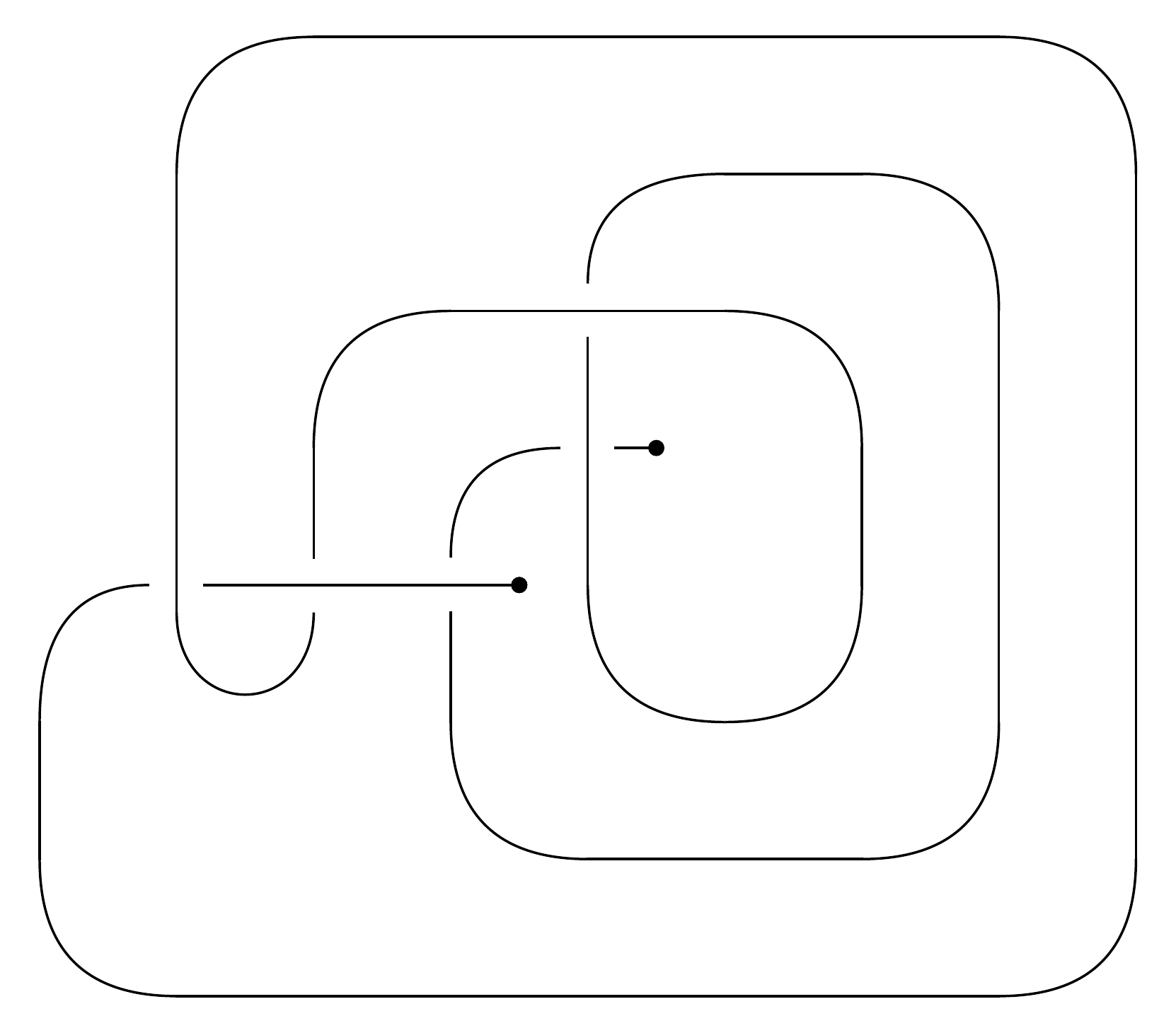}\\
\textcolor{black}{$5_{129}$}
\vspace{1cm}
\end{minipage}
\begin{minipage}[t]{.25\linewidth}
\centering
\includegraphics[width=0.9\textwidth,height=3.5cm,keepaspectratio]{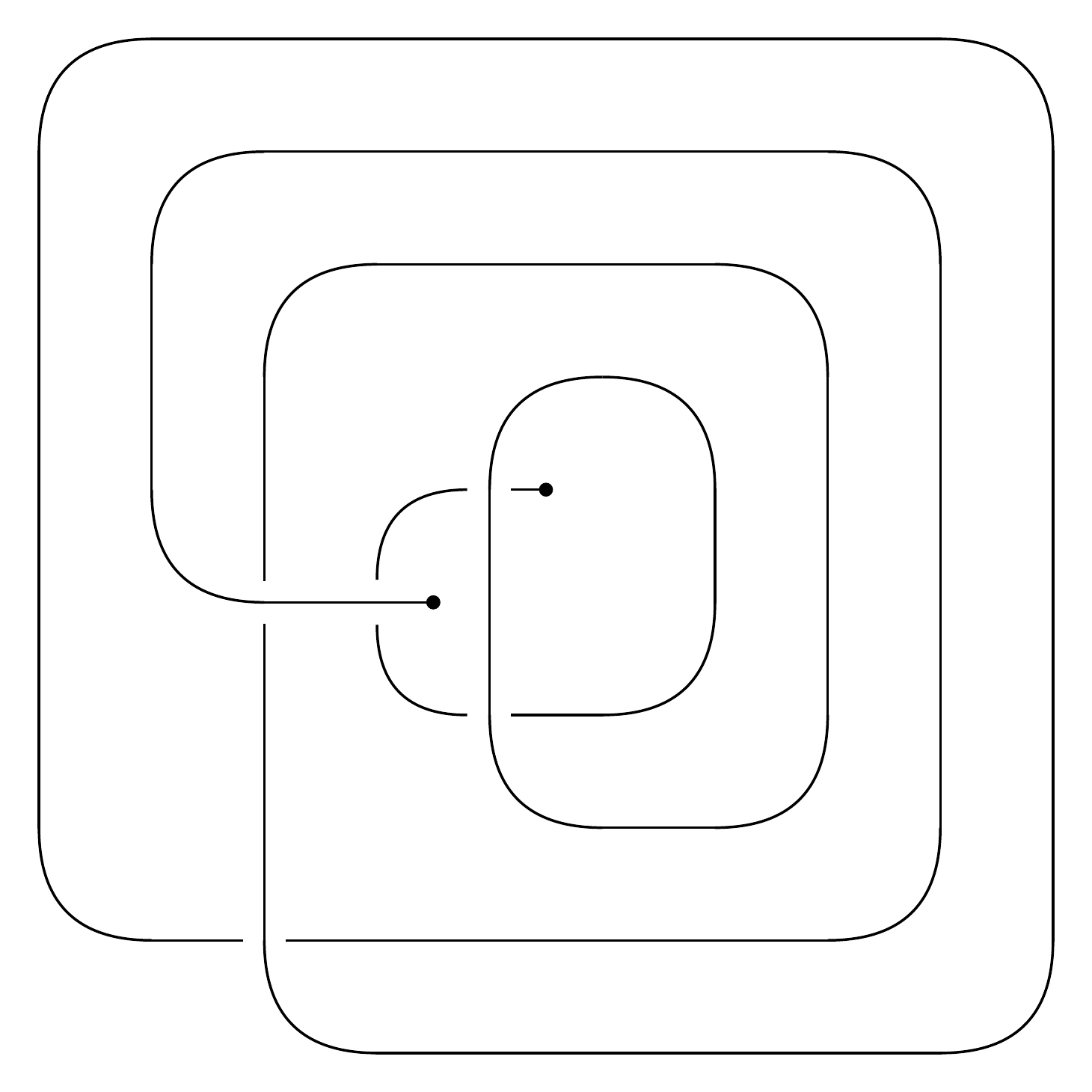}\\
\textcolor{black}{$5_{130}$}
\vspace{1cm}
\end{minipage}
\begin{minipage}[t]{.25\linewidth}
\centering
\includegraphics[width=0.9\textwidth,height=3.5cm,keepaspectratio]{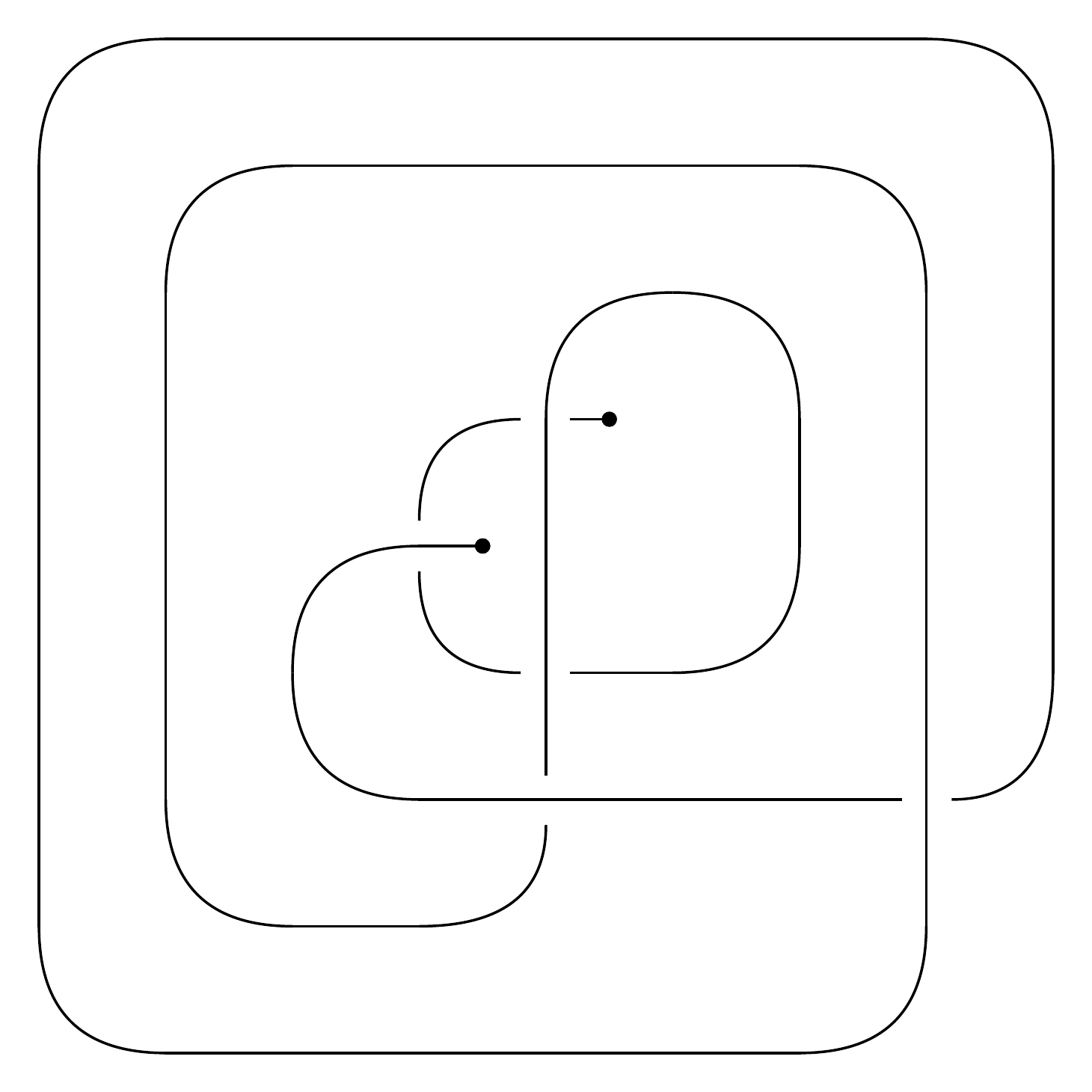}\\
\textcolor{black}{$5_{131}$}
\vspace{1cm}
\end{minipage}
\begin{minipage}[t]{.25\linewidth}
\centering
\includegraphics[width=0.9\textwidth,height=3.5cm,keepaspectratio]{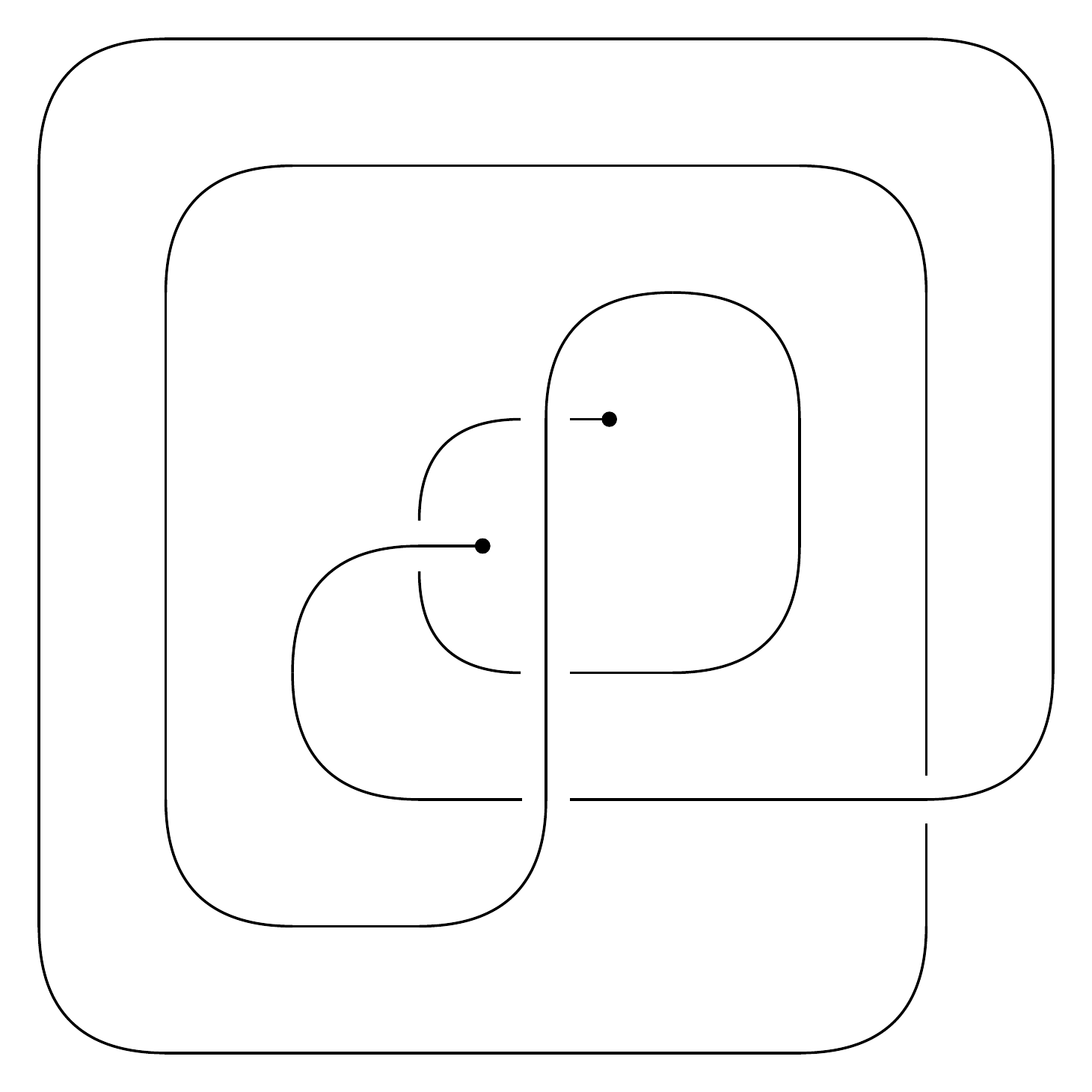}\\
\textcolor{black}{$5_{132}$}
\vspace{1cm}
\end{minipage}
\begin{minipage}[t]{.25\linewidth}
\centering
\includegraphics[width=0.9\textwidth,height=3.5cm,keepaspectratio]{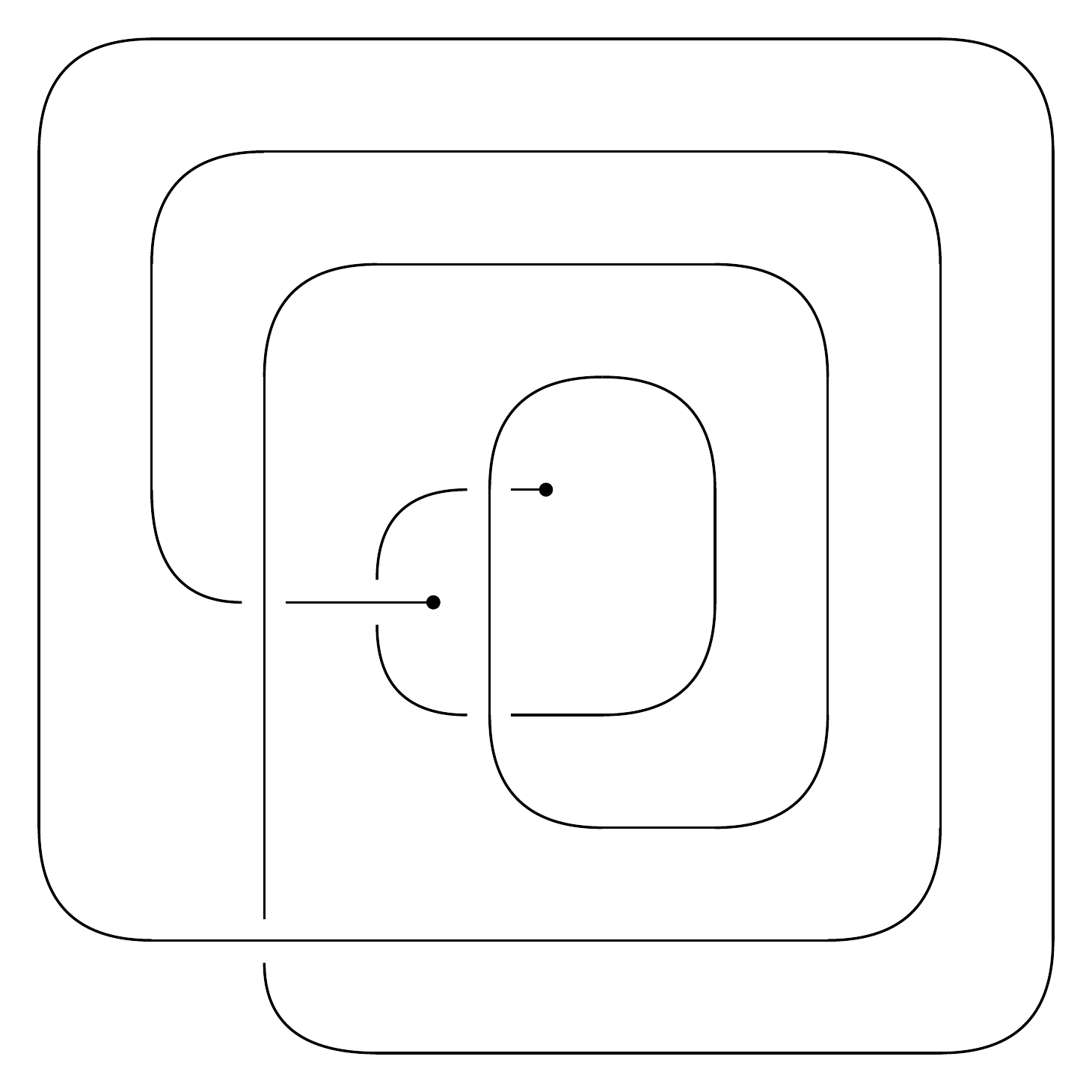}\\
\textcolor{black}{$5_{133}$}
\vspace{1cm}
\end{minipage}
\begin{minipage}[t]{.25\linewidth}
\centering
\includegraphics[width=0.9\textwidth,height=3.5cm,keepaspectratio]{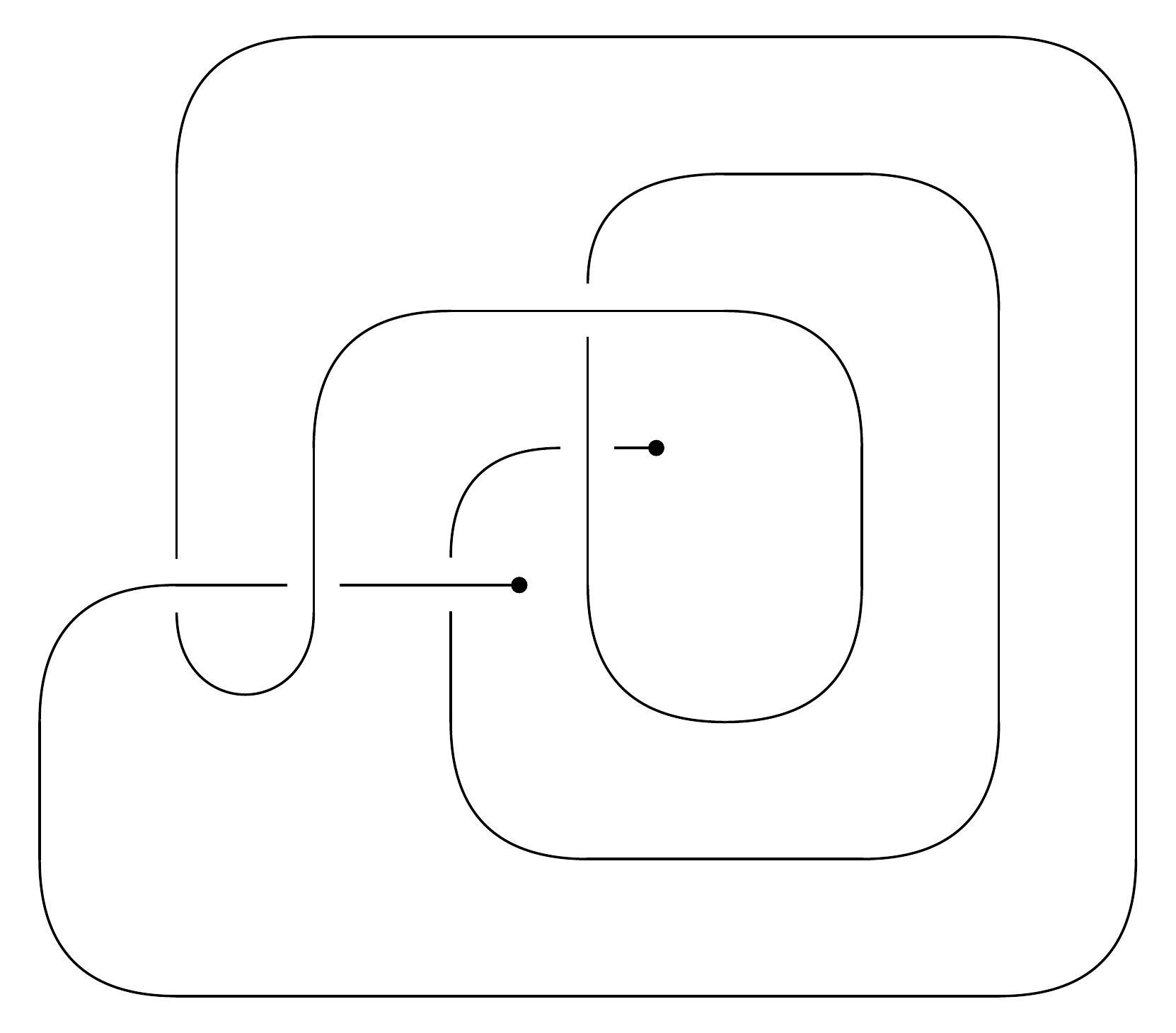}\\
\textcolor{black}{$5_{134}$}
\vspace{1cm}
\end{minipage}
\begin{minipage}[t]{.25\linewidth}
\centering
\includegraphics[width=0.9\textwidth,height=3.5cm,keepaspectratio]{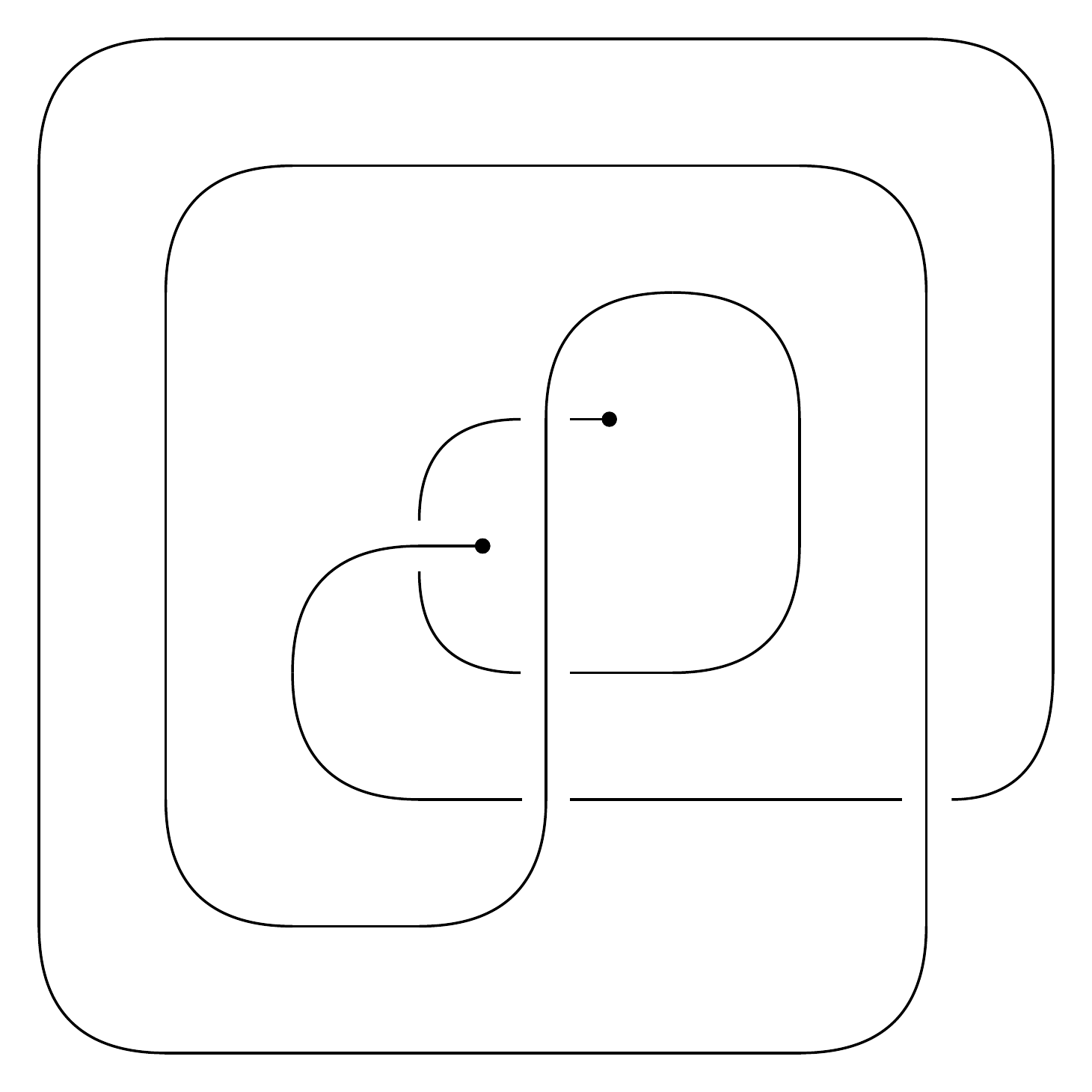}\\
\textcolor{black}{$5_{135}$}
\vspace{1cm}
\end{minipage}
\begin{minipage}[t]{.25\linewidth}
\centering
\includegraphics[width=0.9\textwidth,height=3.5cm,keepaspectratio]{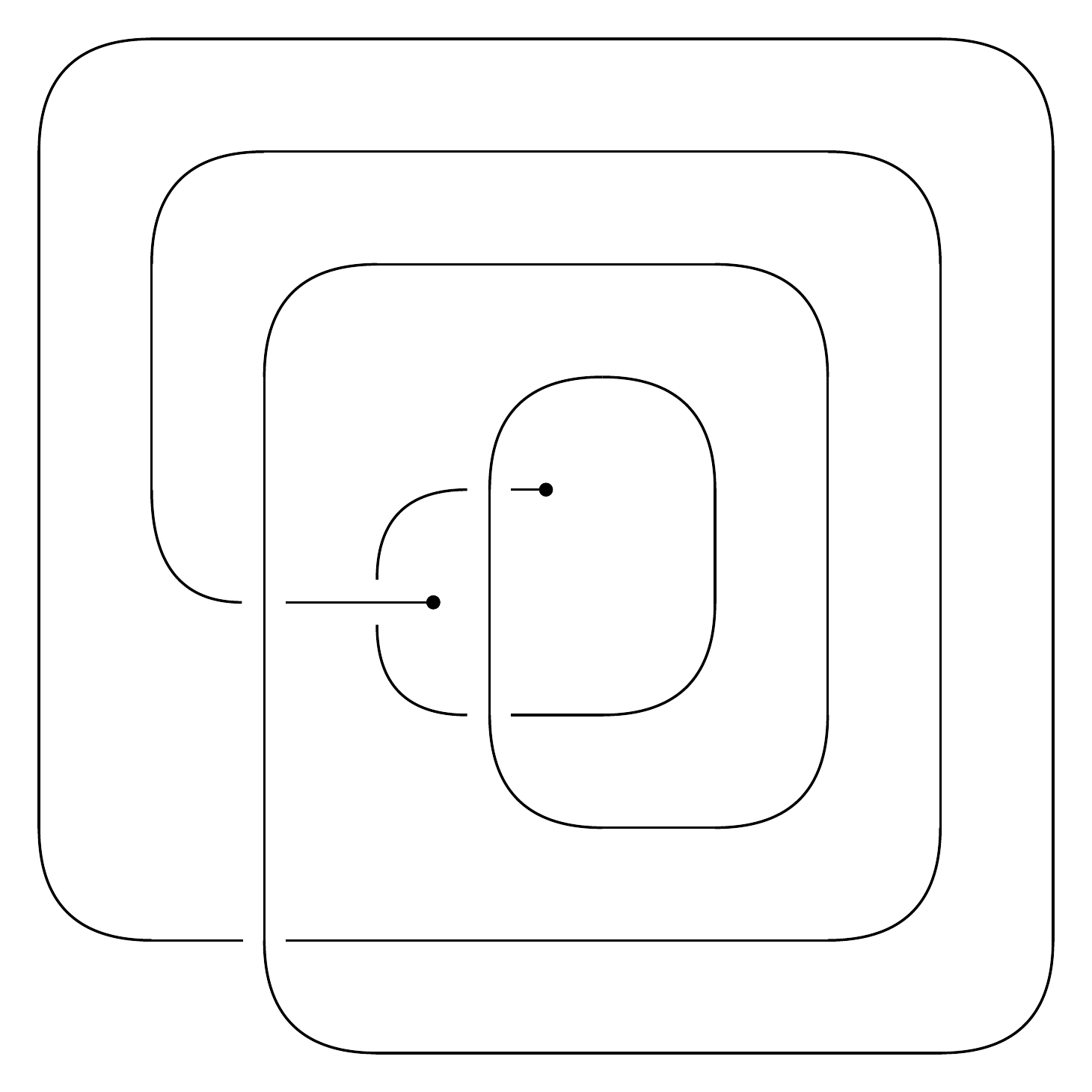}\\
\textcolor{black}{$5_{136}$}
\vspace{1cm}
\end{minipage}
\begin{minipage}[t]{.25\linewidth}
\centering
\includegraphics[width=0.9\textwidth,height=3.5cm,keepaspectratio]{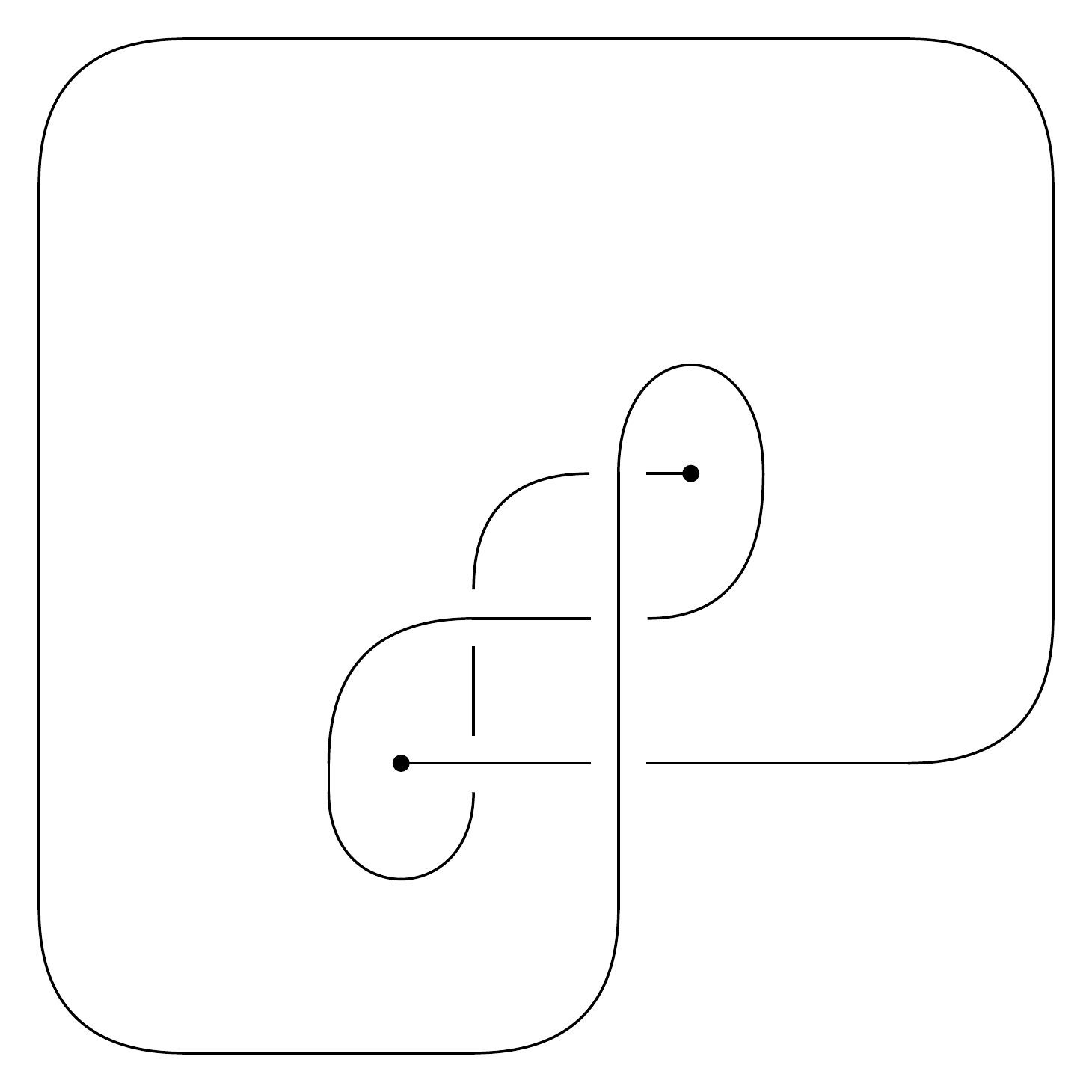}\\
\textcolor{black}{$5_{137}$}
\vspace{1cm}
\end{minipage}
\begin{minipage}[t]{.25\linewidth}
\centering
\includegraphics[width=0.9\textwidth,height=3.5cm,keepaspectratio]{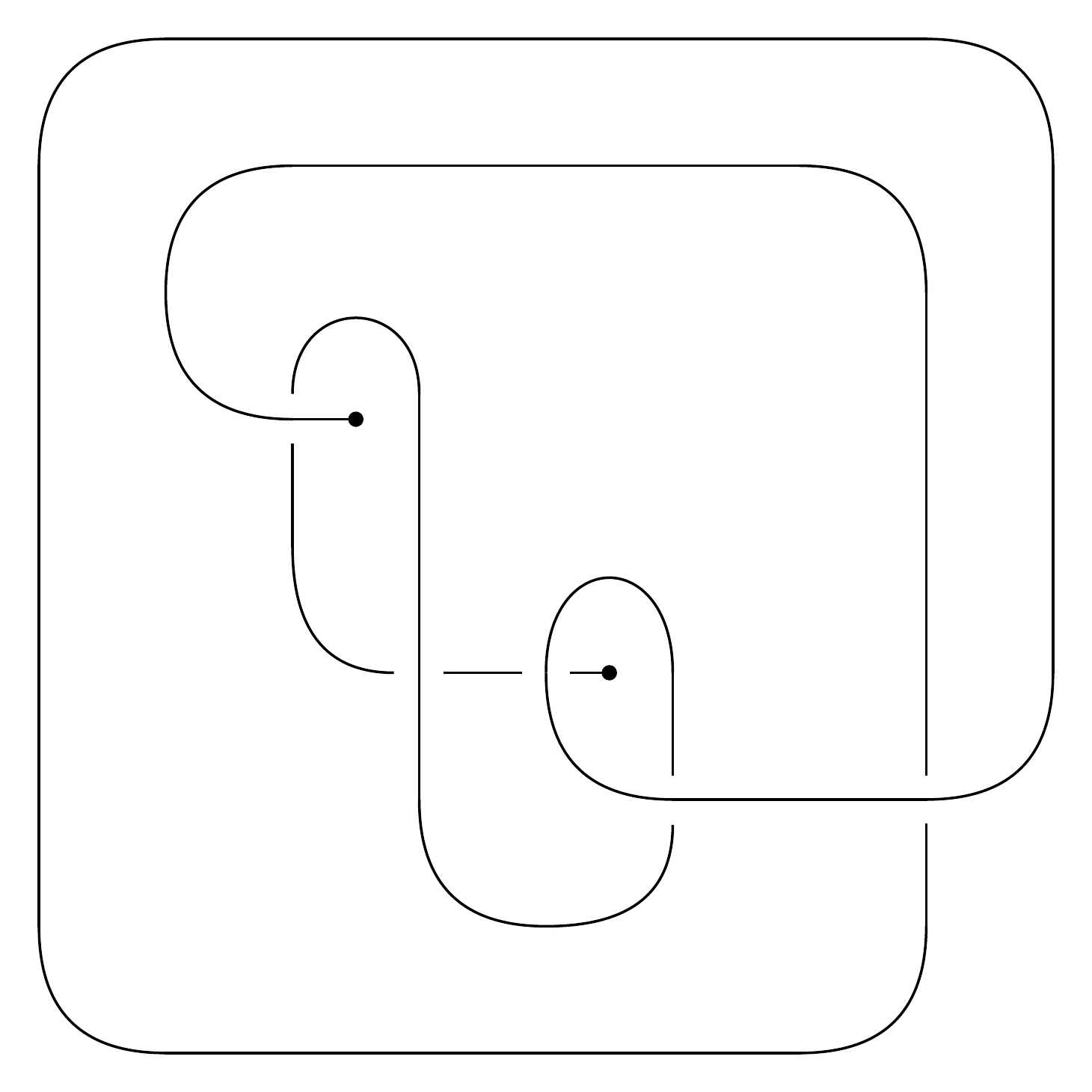}\\
\textcolor{black}{$5_{138}$}
\vspace{1cm}
\end{minipage}
\begin{minipage}[t]{.25\linewidth}
\centering
\includegraphics[width=0.9\textwidth,height=3.5cm,keepaspectratio]{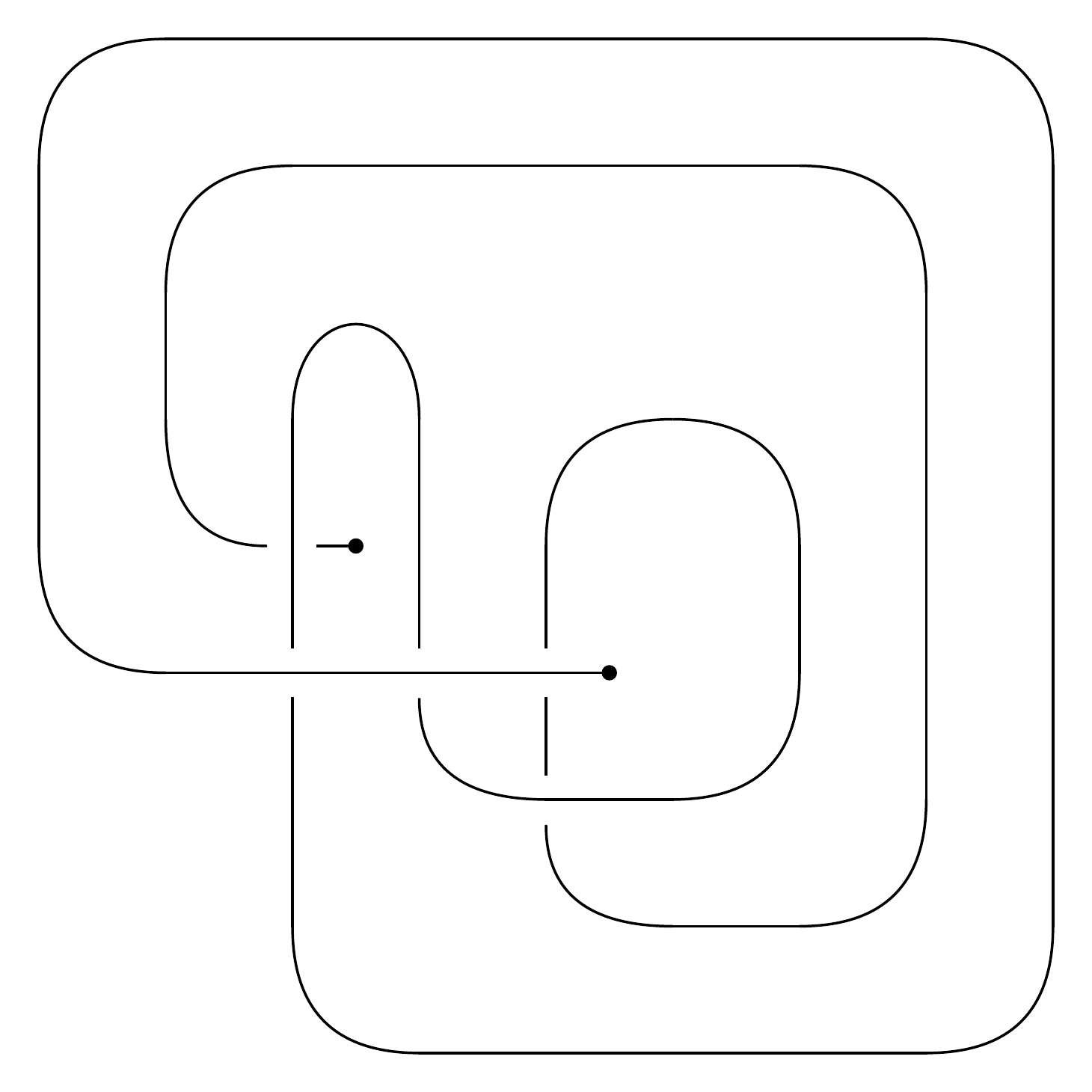}\\
\textcolor{black}{$5_{139}$}
\vspace{1cm}
\end{minipage}
\begin{minipage}[t]{.25\linewidth}
\centering
\includegraphics[width=0.9\textwidth,height=3.5cm,keepaspectratio]{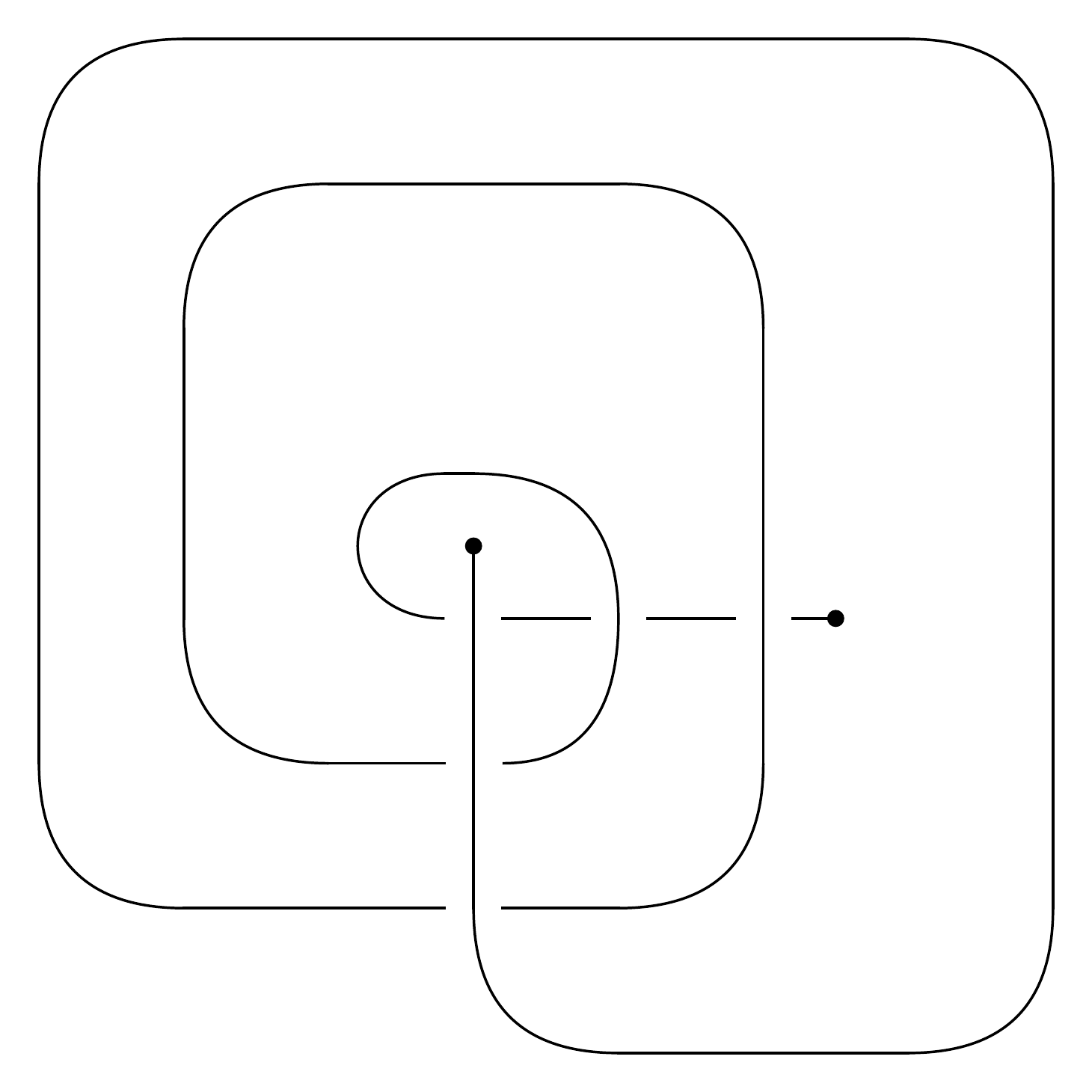}\\
\textcolor{black}{$5_{140}$}
\vspace{1cm}
\end{minipage}
\begin{minipage}[t]{.25\linewidth}
\centering
\includegraphics[width=0.9\textwidth,height=3.5cm,keepaspectratio]{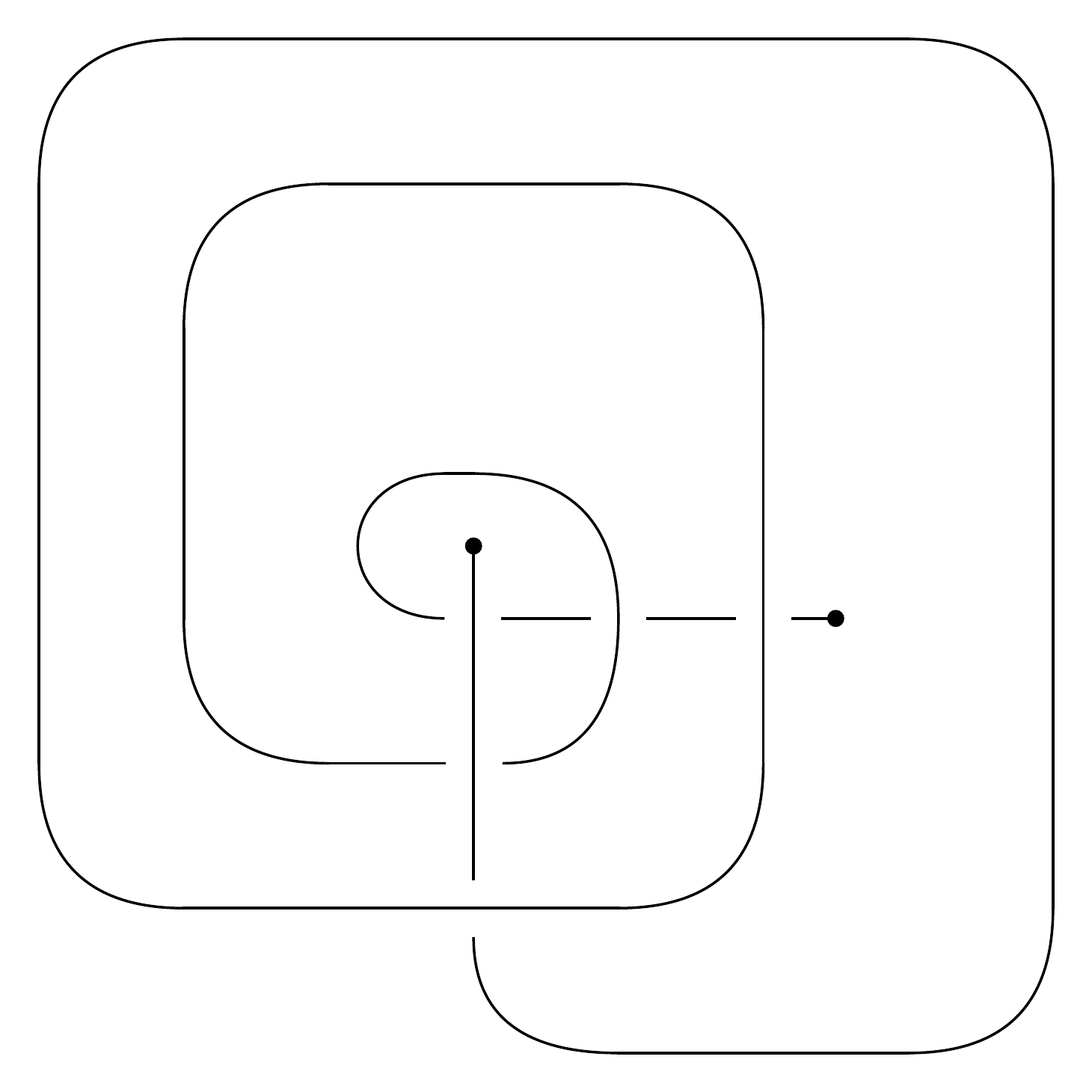}\\
\textcolor{black}{$5_{141}$}
\vspace{1cm}
\end{minipage}
\begin{minipage}[t]{.25\linewidth}
\centering
\includegraphics[width=0.9\textwidth,height=3.5cm,keepaspectratio]{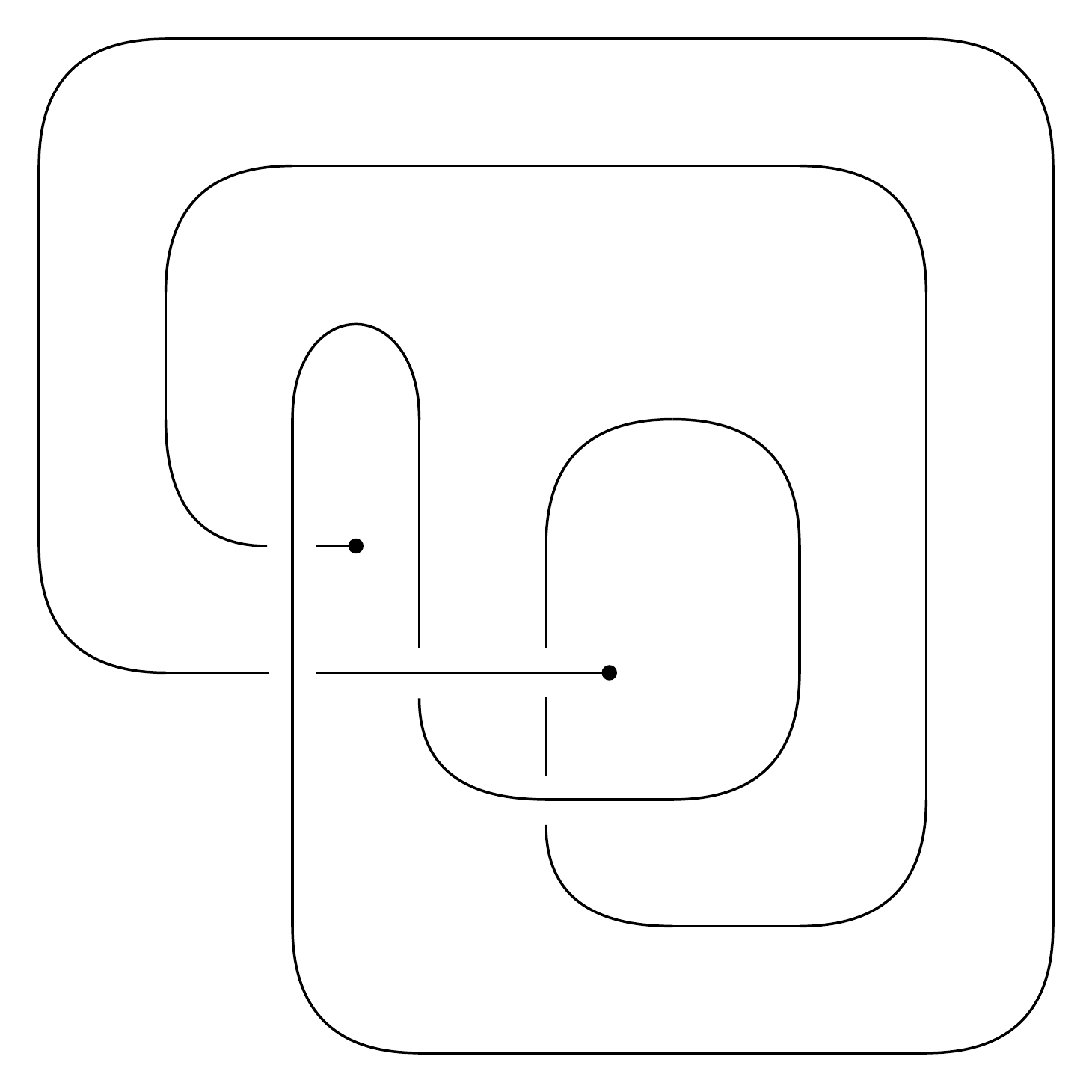}\\
\textcolor{black}{$5_{142}$}
\vspace{1cm}
\end{minipage}
\begin{minipage}[t]{.25\linewidth}
\centering
\includegraphics[width=0.9\textwidth,height=3.5cm,keepaspectratio]{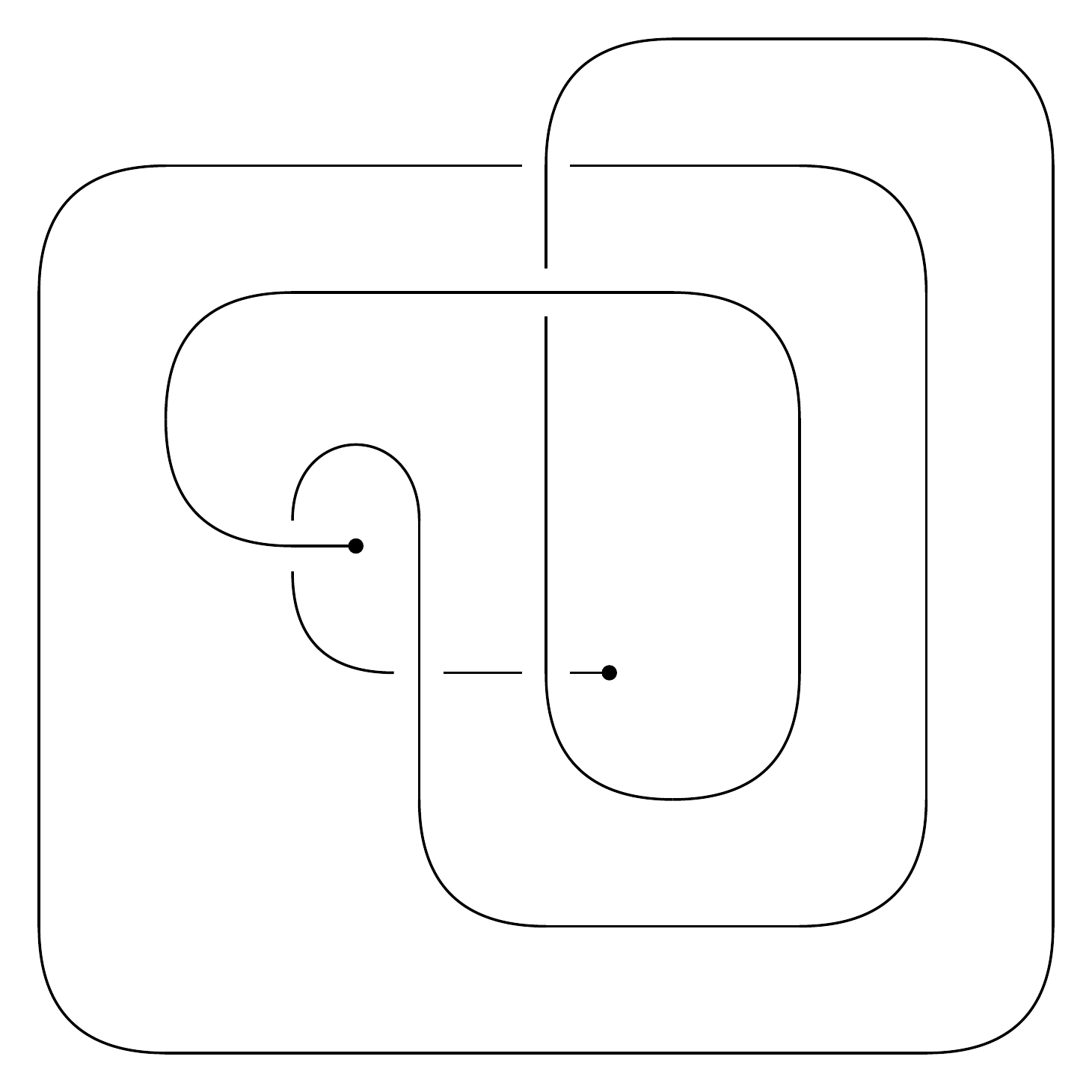}\\
\textcolor{black}{$5_{143}$}
\vspace{1cm}
\end{minipage}
\begin{minipage}[t]{.25\linewidth}
\centering
\includegraphics[width=0.9\textwidth,height=3.5cm,keepaspectratio]{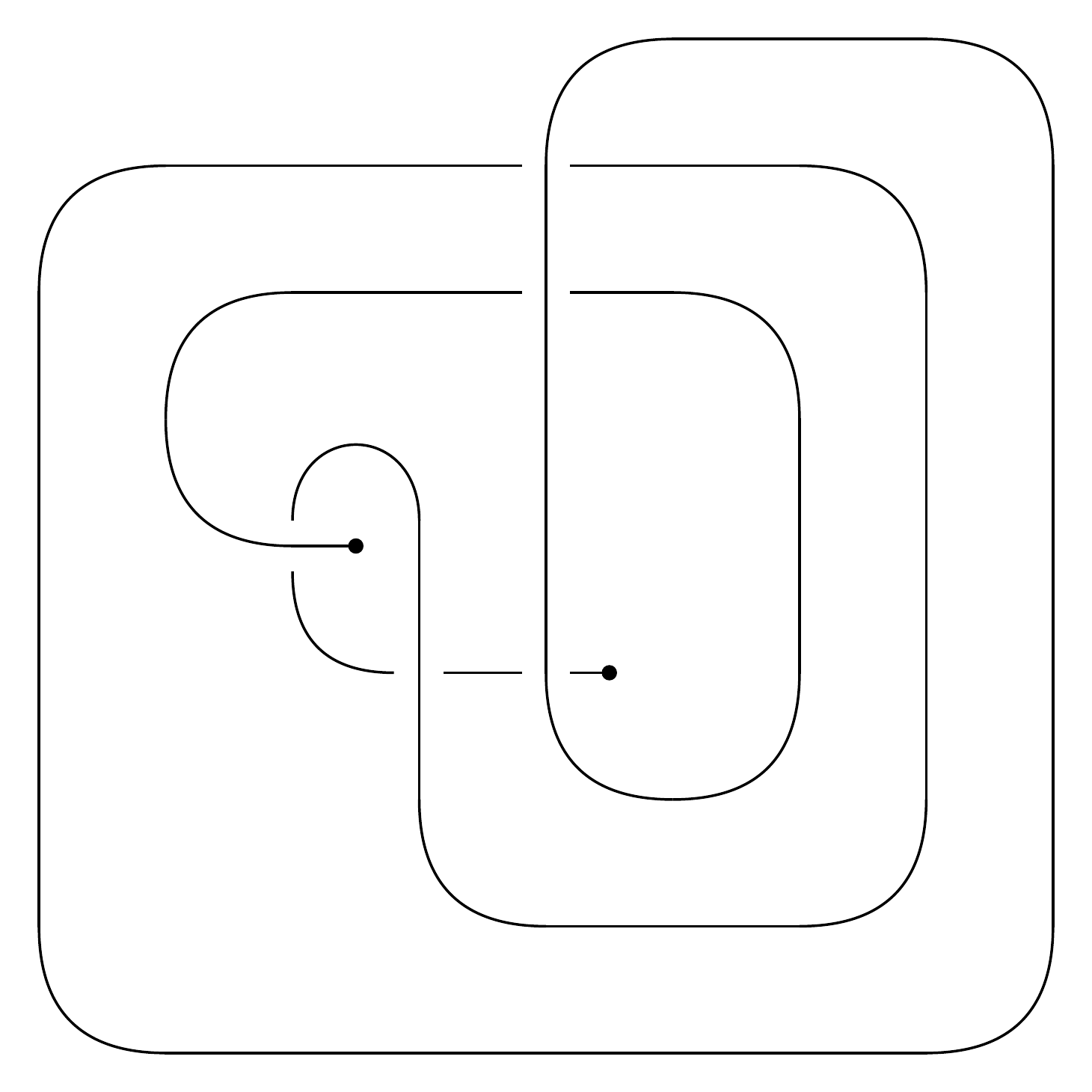}\\
\textcolor{black}{$5_{144}$}
\vspace{1cm}
\end{minipage}
\begin{minipage}[t]{.25\linewidth}
\centering
\includegraphics[width=0.9\textwidth,height=3.5cm,keepaspectratio]{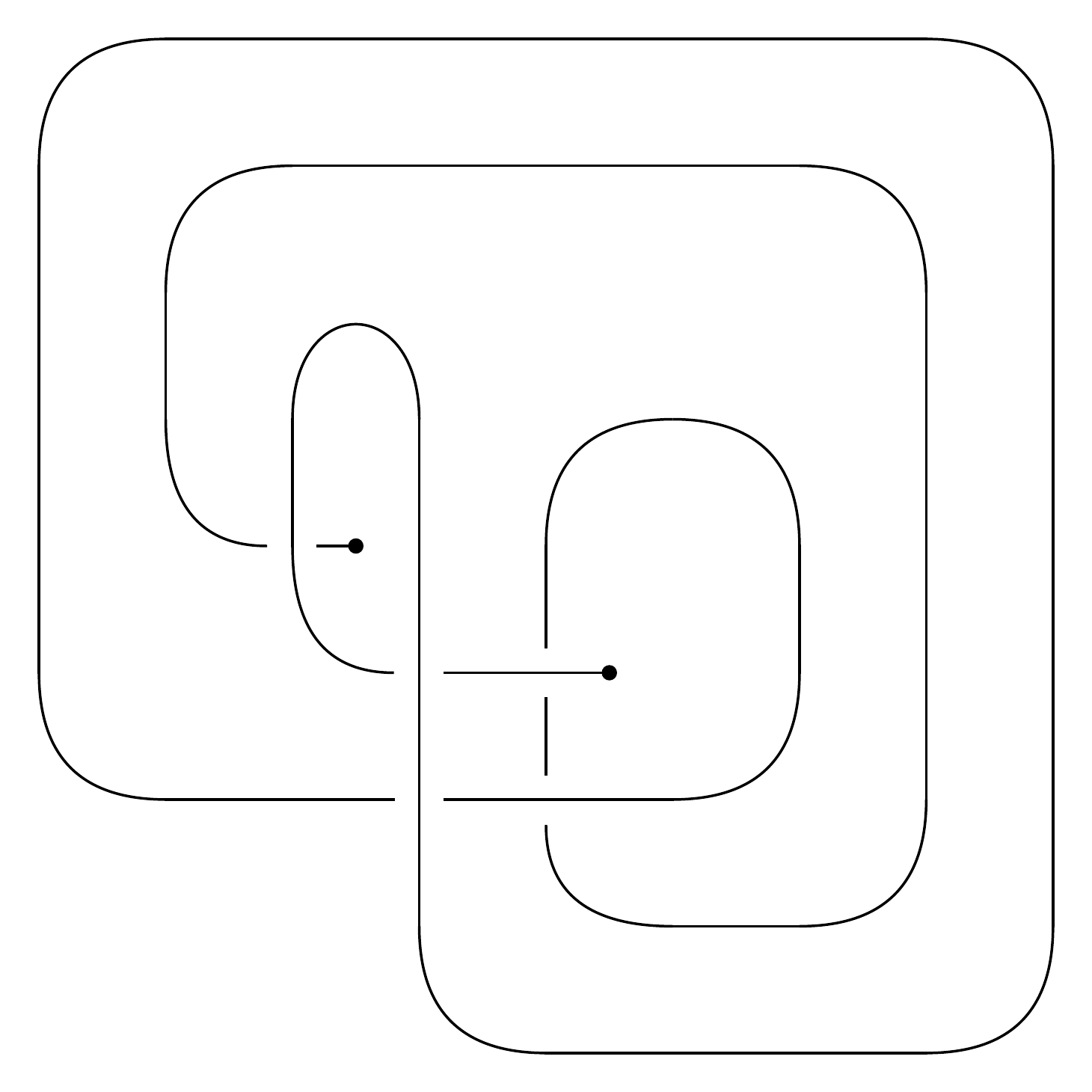}\\
\textcolor{black}{$5_{145}$}
\vspace{1cm}
\end{minipage}
\begin{minipage}[t]{.25\linewidth}
\centering
\includegraphics[width=0.9\textwidth,height=3.5cm,keepaspectratio]{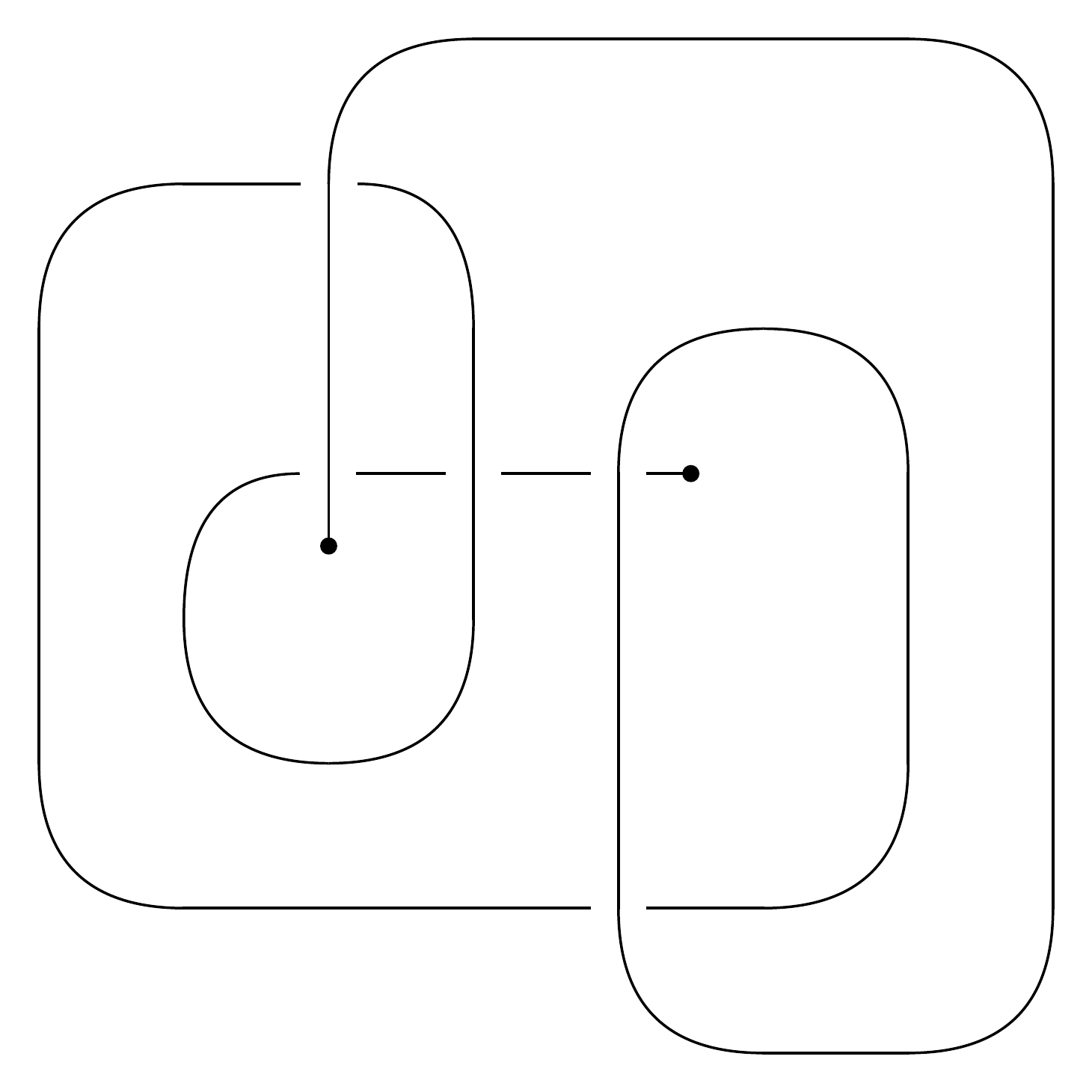}\\
\textcolor{black}{$5_{146}$}
\vspace{1cm}
\end{minipage}
\begin{minipage}[t]{.25\linewidth}
\centering
\includegraphics[width=0.9\textwidth,height=3.5cm,keepaspectratio]{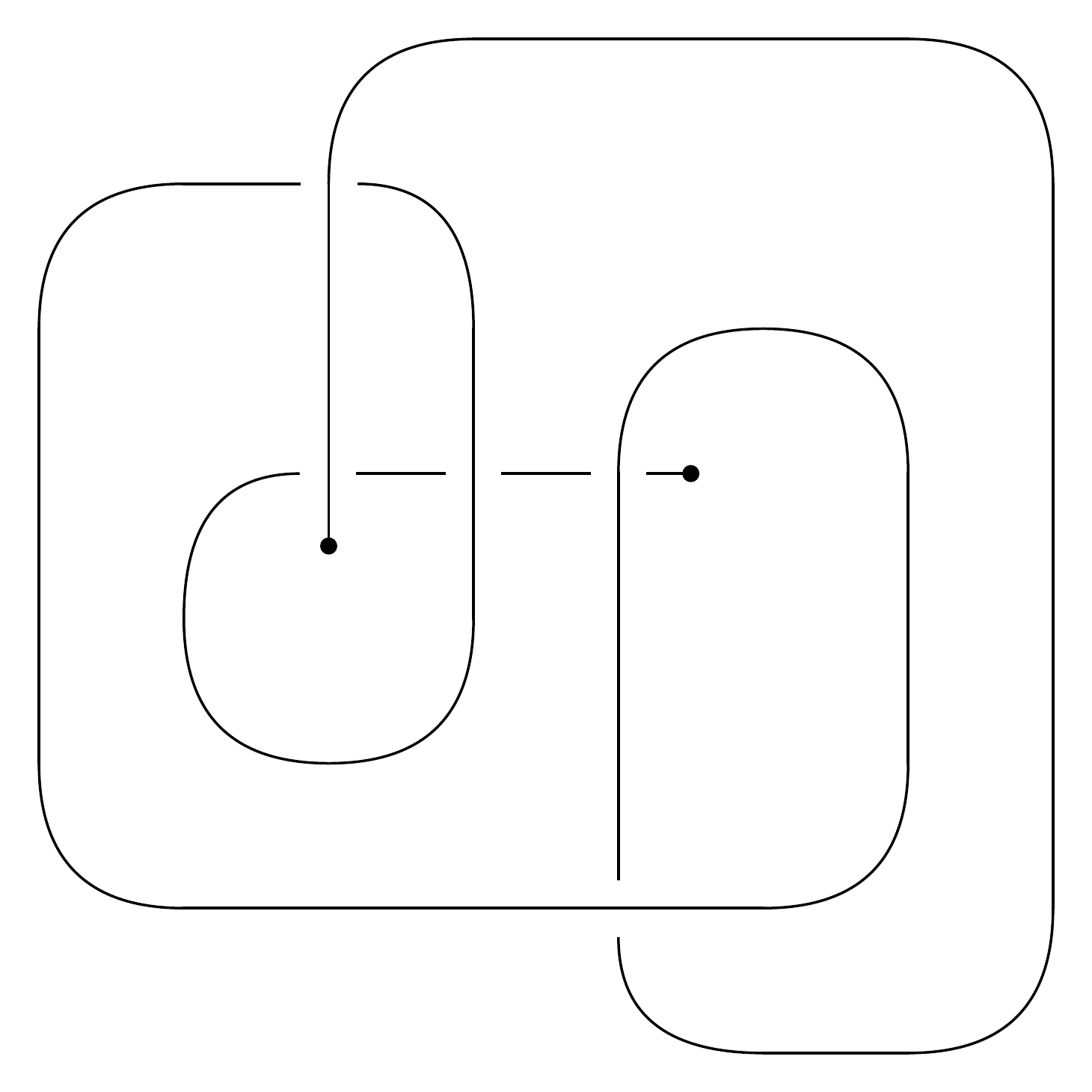}\\
\textcolor{black}{$5_{147}$}
\vspace{1cm}
\end{minipage}
\begin{minipage}[t]{.25\linewidth}
\centering
\includegraphics[width=0.9\textwidth,height=3.5cm,keepaspectratio]{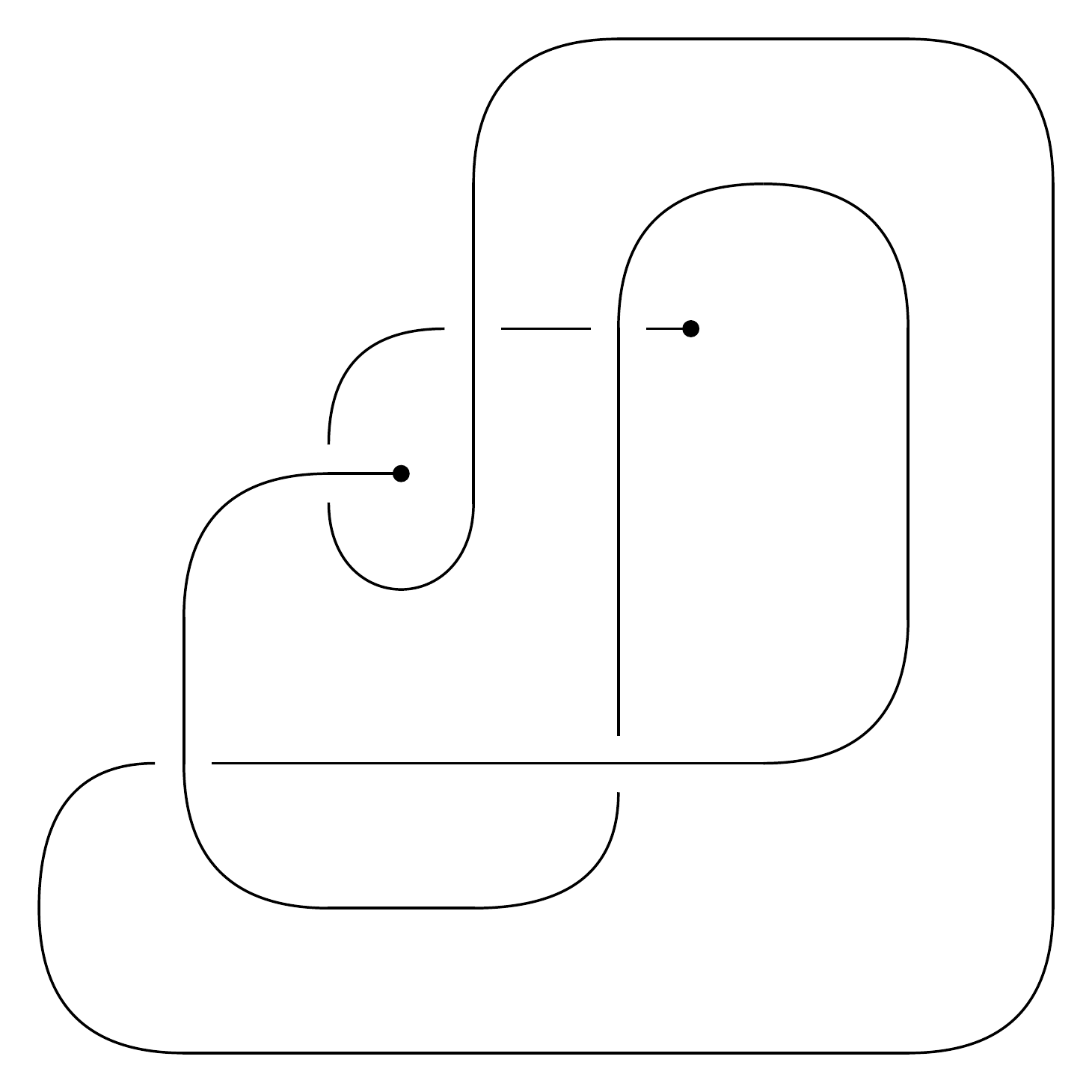}\\
\textcolor{black}{$5_{148}$}
\vspace{1cm}
\end{minipage}
\begin{minipage}[t]{.25\linewidth}
\centering
\includegraphics[width=0.9\textwidth,height=3.5cm,keepaspectratio]{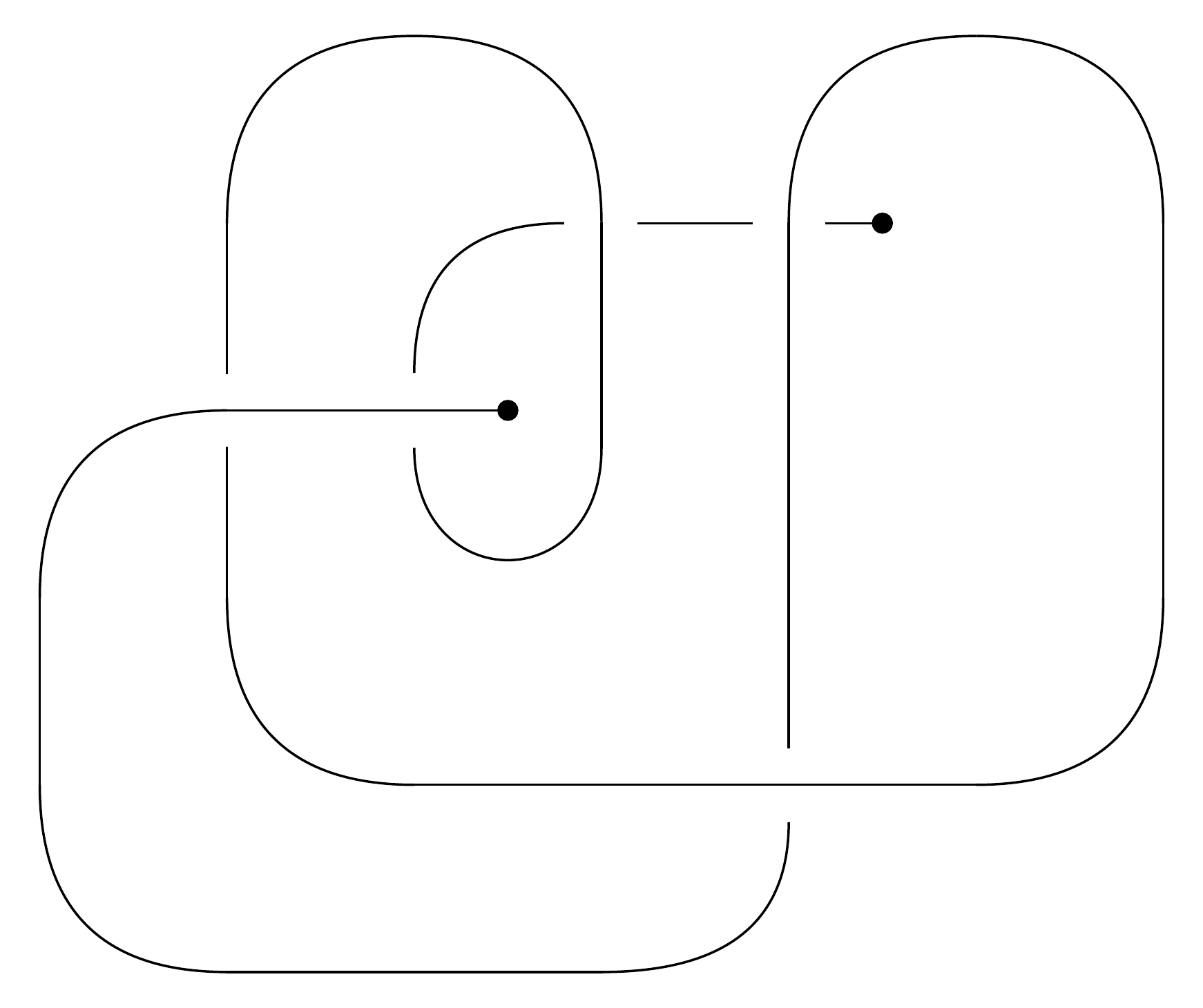}\\
\textcolor{black}{$5_{149}$}
\vspace{1cm}
\end{minipage}
\begin{minipage}[t]{.25\linewidth}
\centering
\includegraphics[width=0.9\textwidth,height=3.5cm,keepaspectratio]{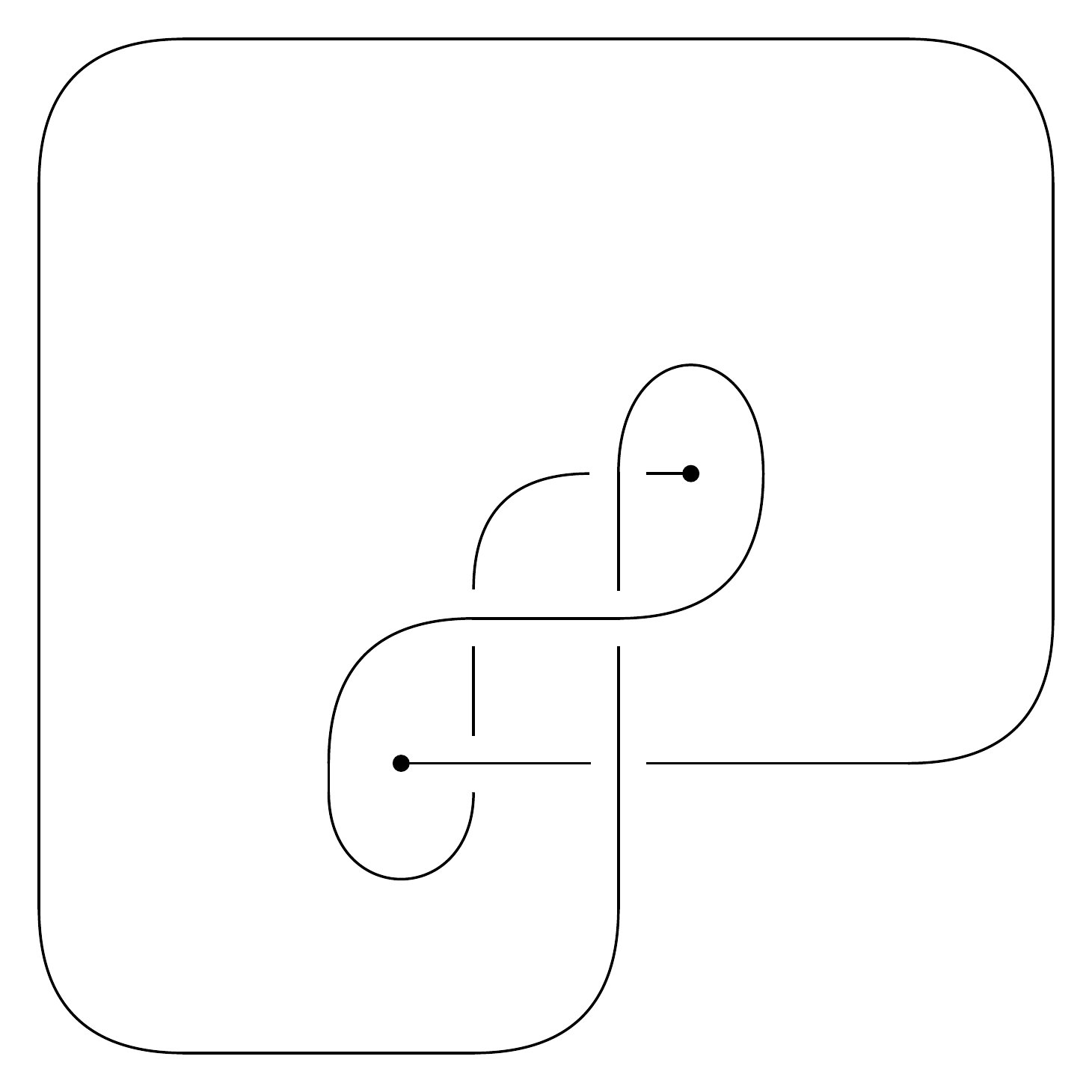}\\
\textcolor{black}{$5_{150}$}
\vspace{1cm}
\end{minipage}
\begin{minipage}[t]{.25\linewidth}
\centering
\includegraphics[width=0.9\textwidth,height=3.5cm,keepaspectratio]{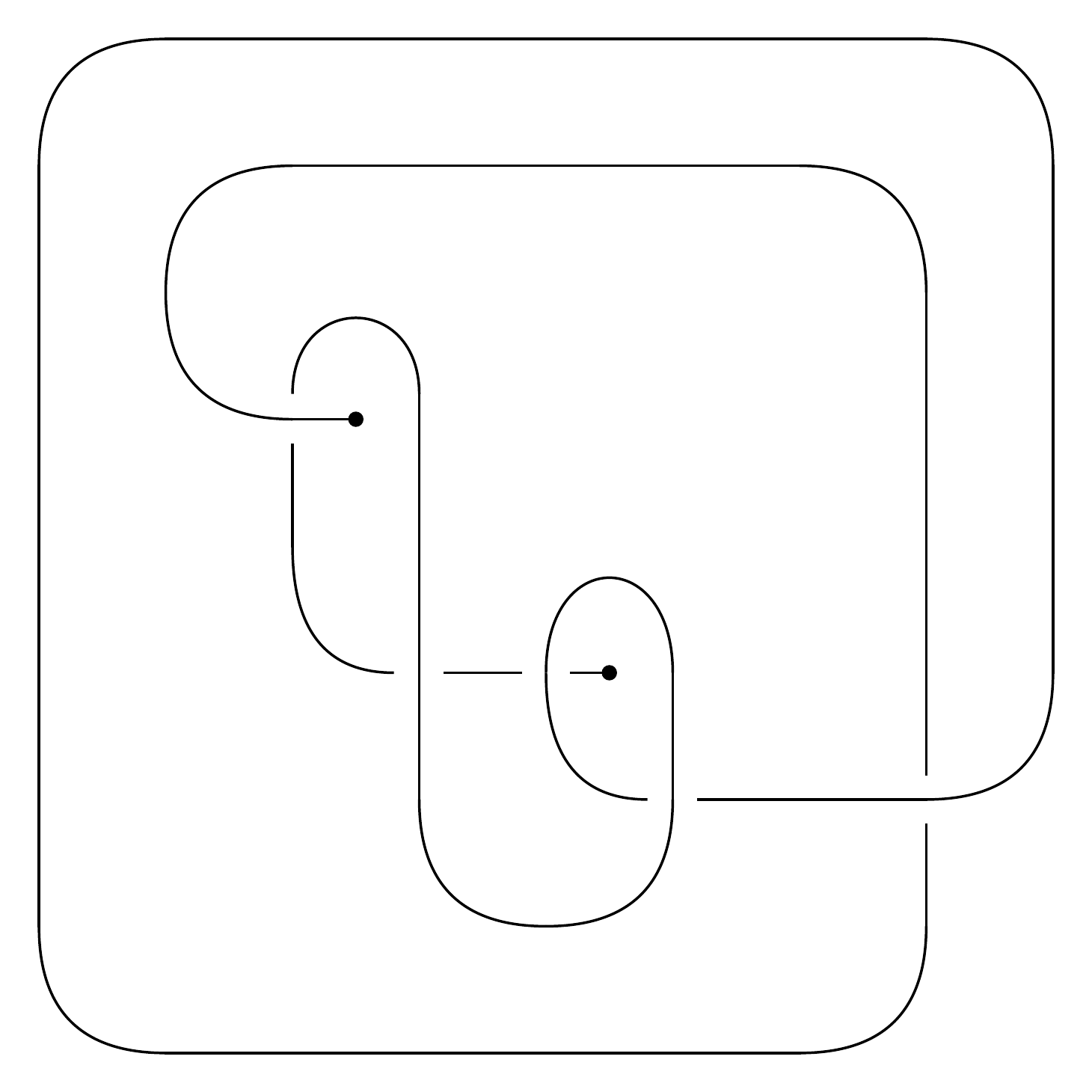}\\
\textcolor{black}{$5_{151}$}
\vspace{1cm}
\end{minipage}
\begin{minipage}[t]{.25\linewidth}
\centering
\includegraphics[width=0.9\textwidth,height=3.5cm,keepaspectratio]{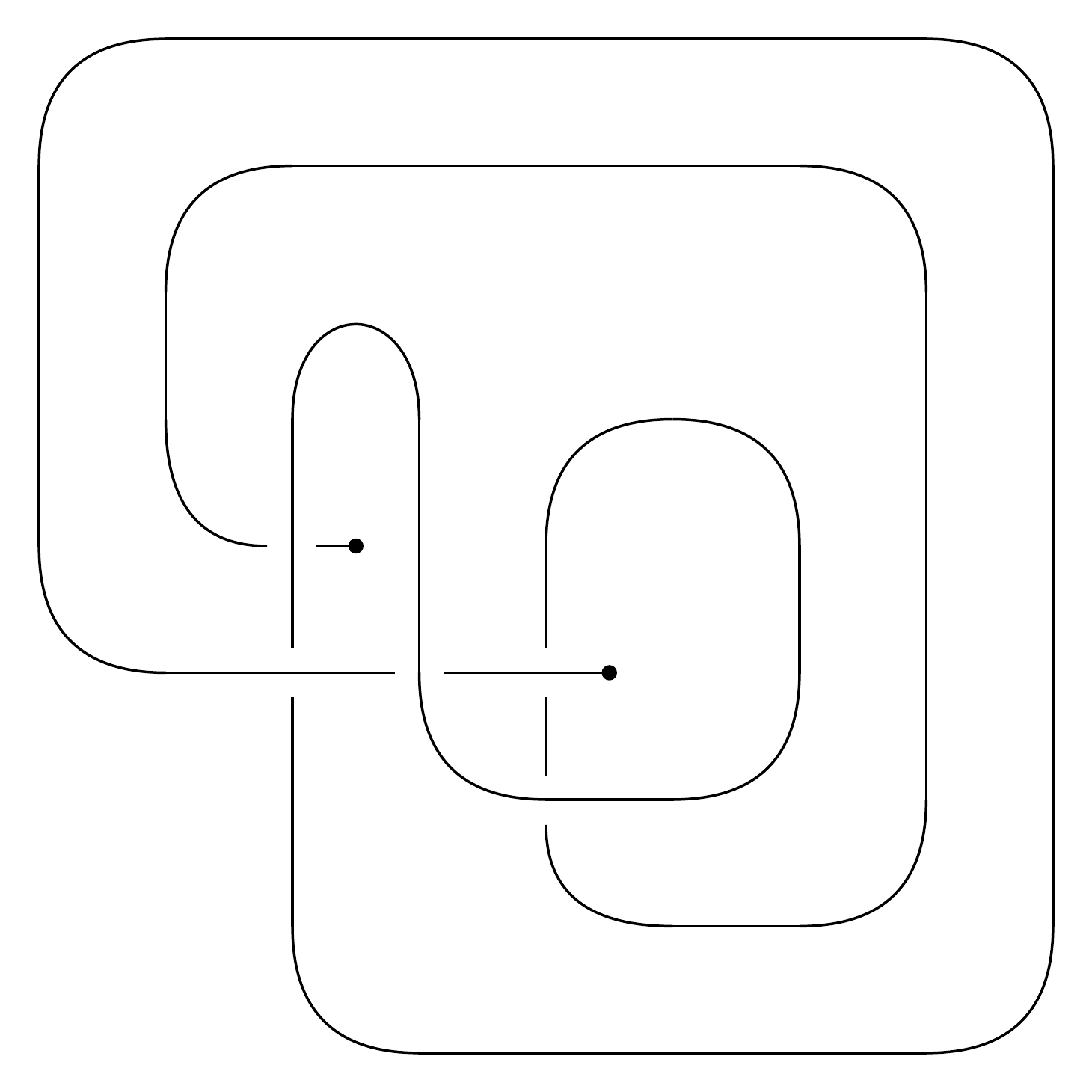}\\
\textcolor{black}{$5_{152}$}
\vspace{1cm}
\end{minipage}
\begin{minipage}[t]{.25\linewidth}
\centering
\includegraphics[width=0.9\textwidth,height=3.5cm,keepaspectratio]{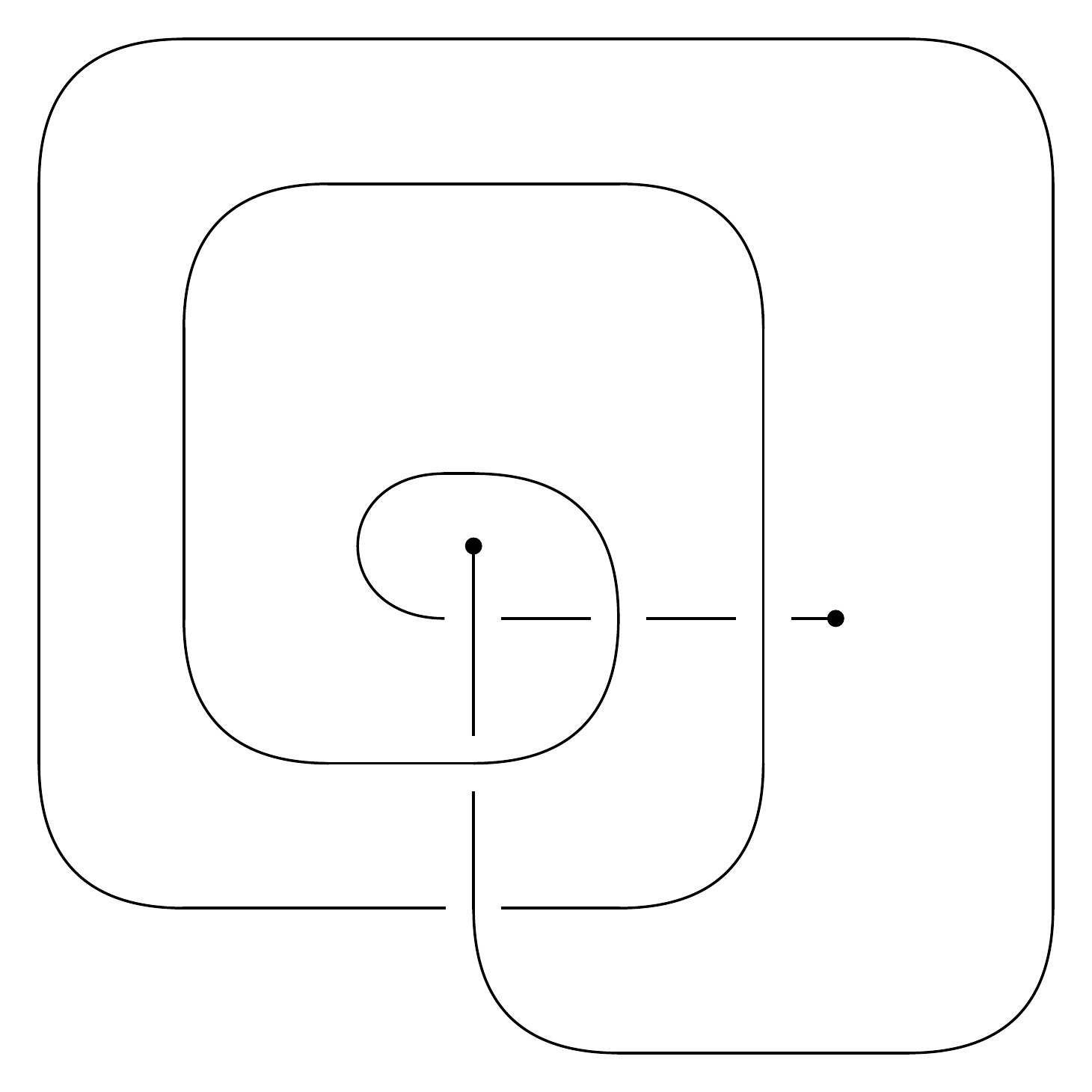}\\
\textcolor{black}{$5_{153}$}
\vspace{1cm}
\end{minipage}
\begin{minipage}[t]{.25\linewidth}
\centering
\includegraphics[width=0.9\textwidth,height=3.5cm,keepaspectratio]{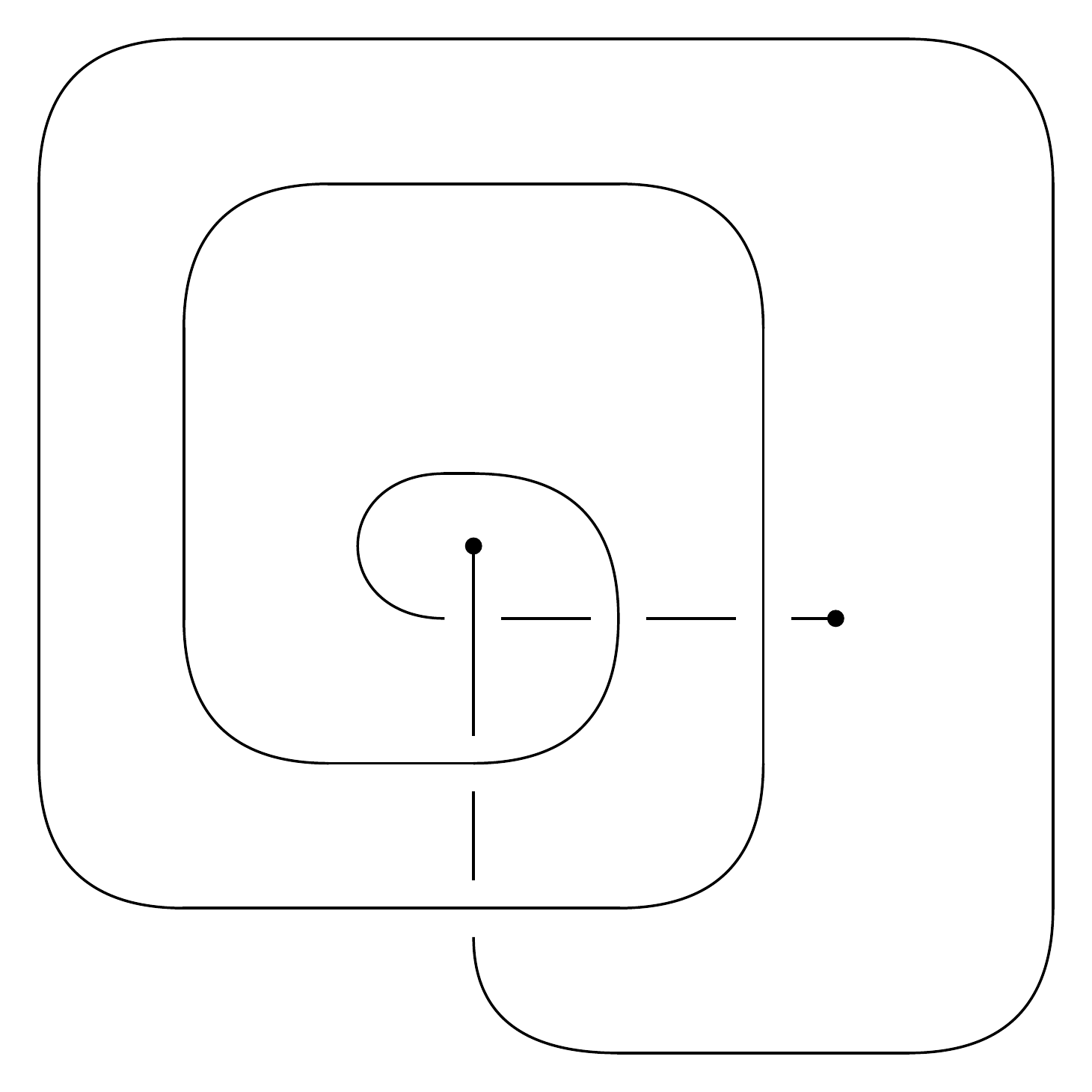}\\
\textcolor{black}{$5_{154}$}
\vspace{1cm}
\end{minipage}
\begin{minipage}[t]{.25\linewidth}
\centering
\includegraphics[width=0.9\textwidth,height=3.5cm,keepaspectratio]{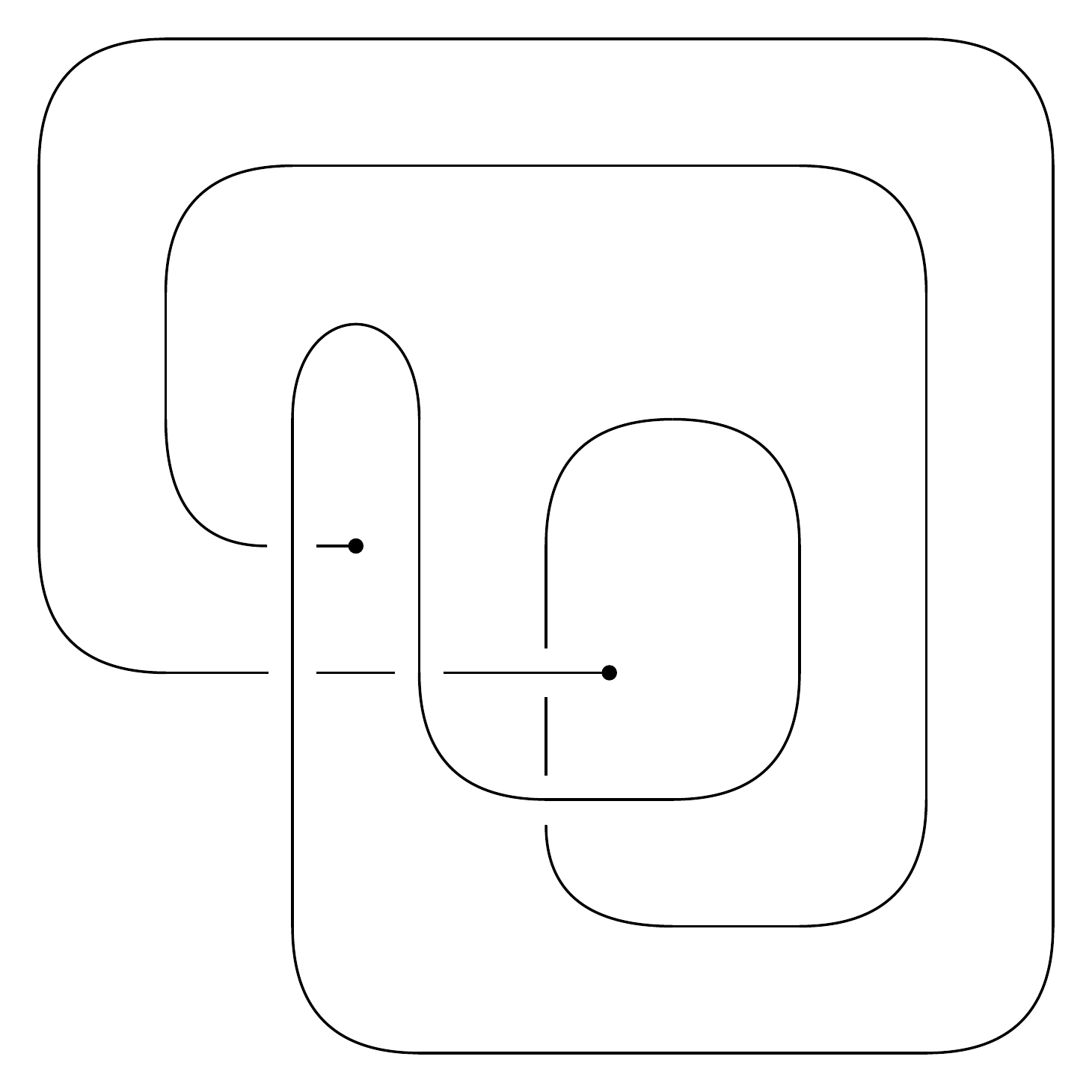}\\
\textcolor{black}{$5_{155}$}
\vspace{1cm}
\end{minipage}
\begin{minipage}[t]{.25\linewidth}
\centering
\includegraphics[width=0.9\textwidth,height=3.5cm,keepaspectratio]{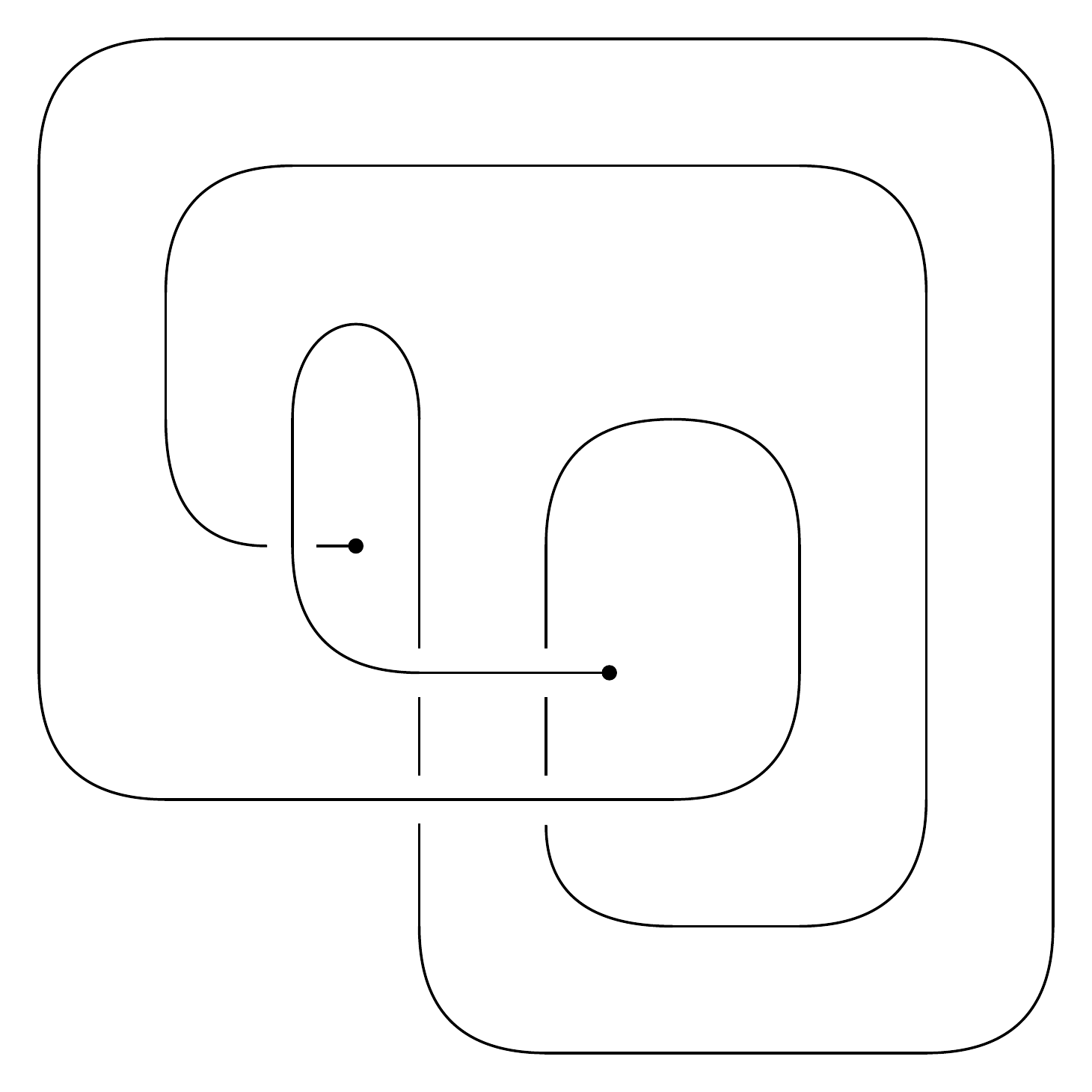}\\
\textcolor{black}{$5_{156}$}
\vspace{1cm}
\end{minipage}
\begin{minipage}[t]{.25\linewidth}
\centering
\includegraphics[width=0.9\textwidth,height=3.5cm,keepaspectratio]{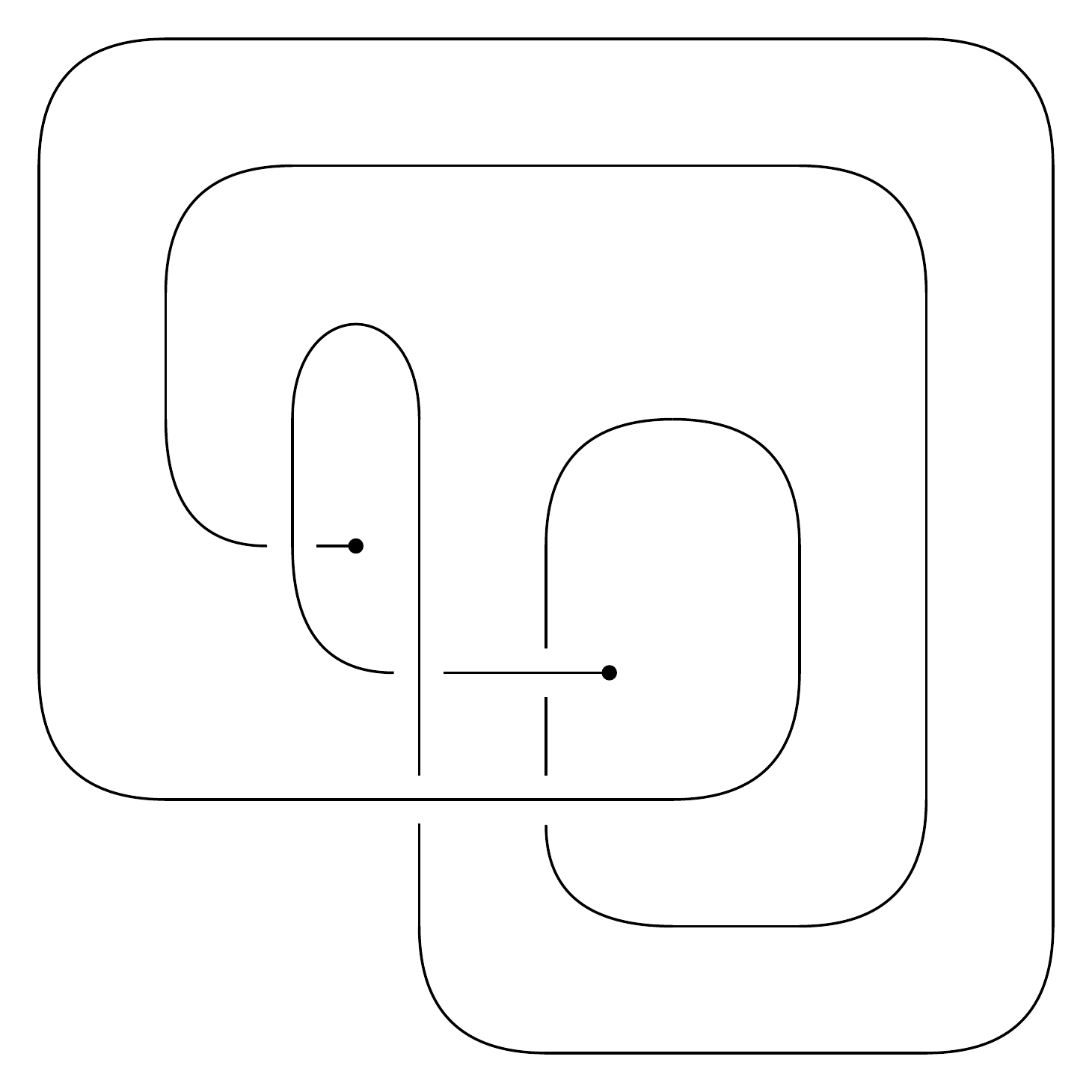}\\
\textcolor{black}{$5_{157}$}
\vspace{1cm}
\end{minipage}
\begin{minipage}[t]{.25\linewidth}
\centering
\includegraphics[width=0.9\textwidth,height=3.5cm,keepaspectratio]{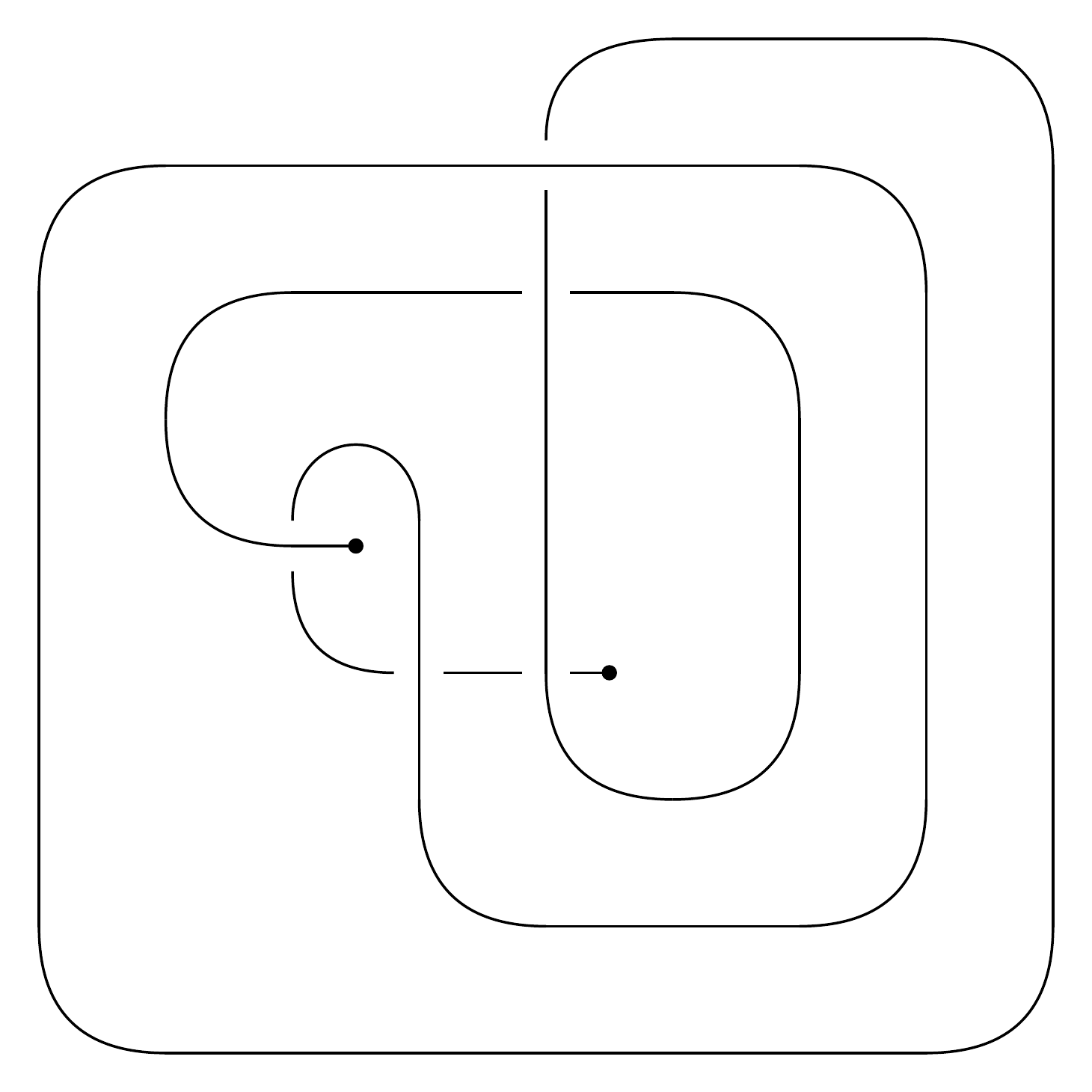}\\
\textcolor{black}{$5_{158}$}
\vspace{1cm}
\end{minipage}
\begin{minipage}[t]{.25\linewidth}
\centering
\includegraphics[width=0.9\textwidth,height=3.5cm,keepaspectratio]{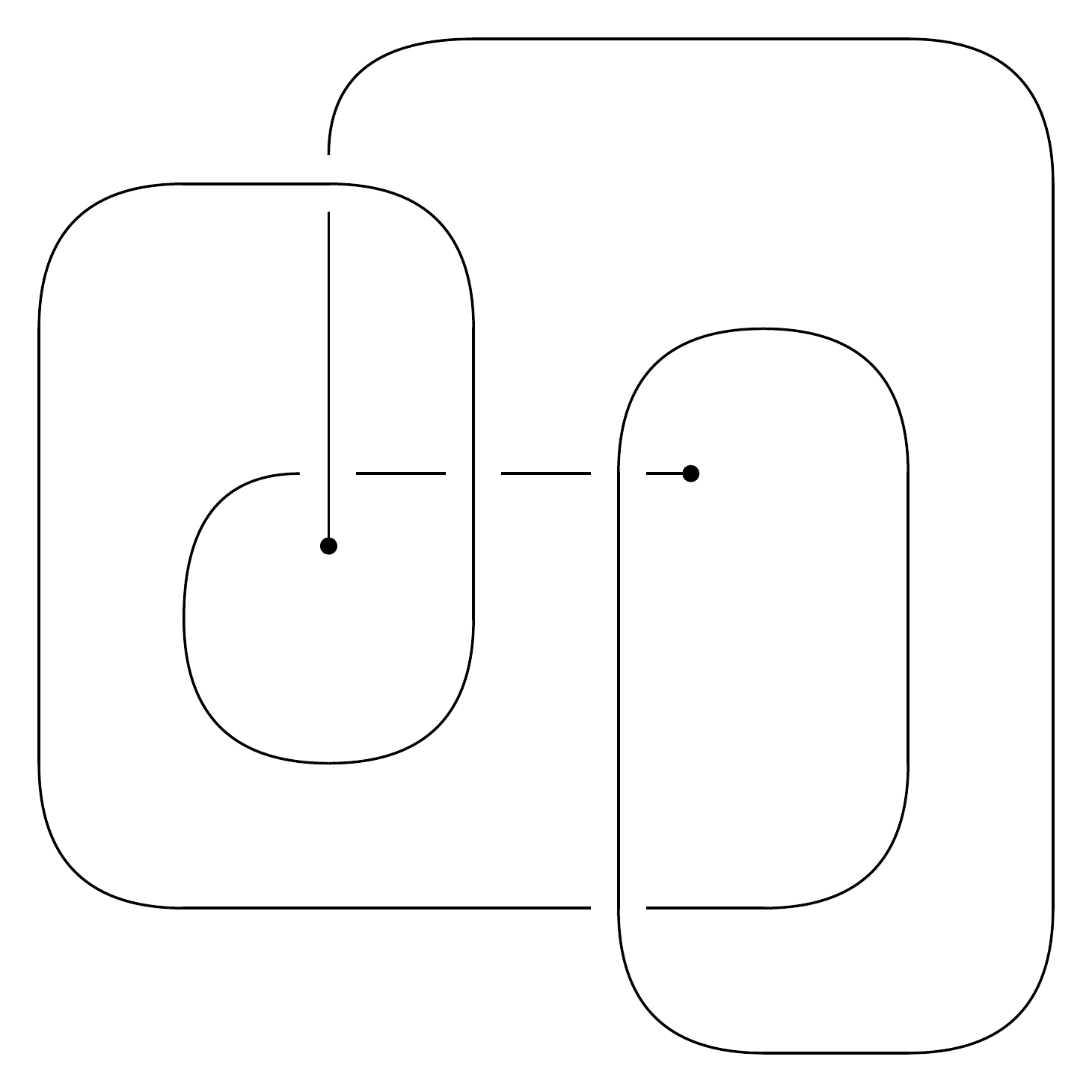}\\
\textcolor{black}{$5_{159}$}
\vspace{1cm}
\end{minipage}
\begin{minipage}[t]{.25\linewidth}
\centering
\includegraphics[width=0.9\textwidth,height=3.5cm,keepaspectratio]{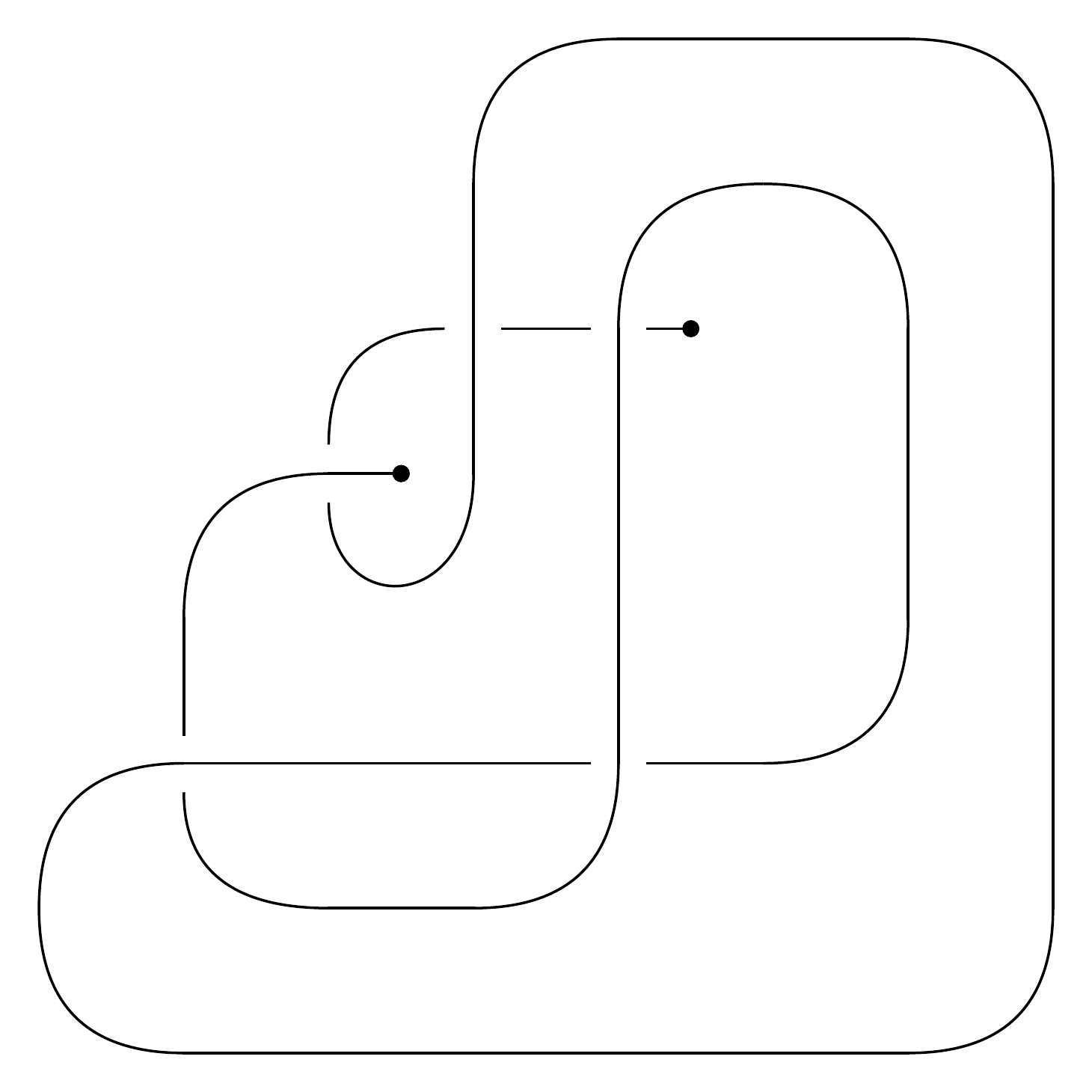}\\
\textcolor{black}{$5_{160}$}
\vspace{1cm}
\end{minipage}
\begin{minipage}[t]{.25\linewidth}
\centering
\includegraphics[width=0.9\textwidth,height=3.5cm,keepaspectratio]{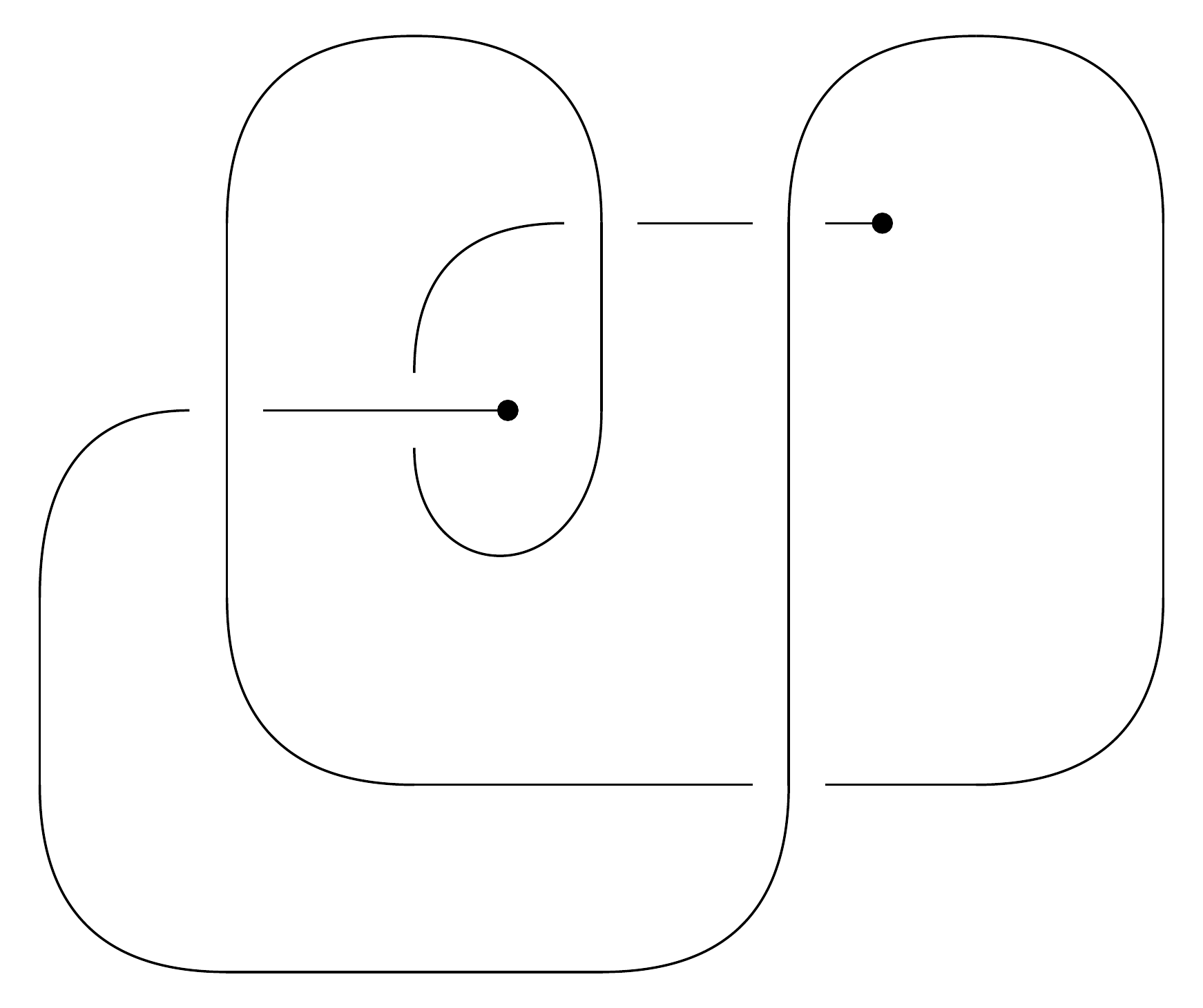}\\
\textcolor{black}{$5_{161}$}
\vspace{1cm}
\end{minipage}
\begin{minipage}[t]{.25\linewidth}
\centering
\includegraphics[width=0.9\textwidth,height=3.5cm,keepaspectratio]{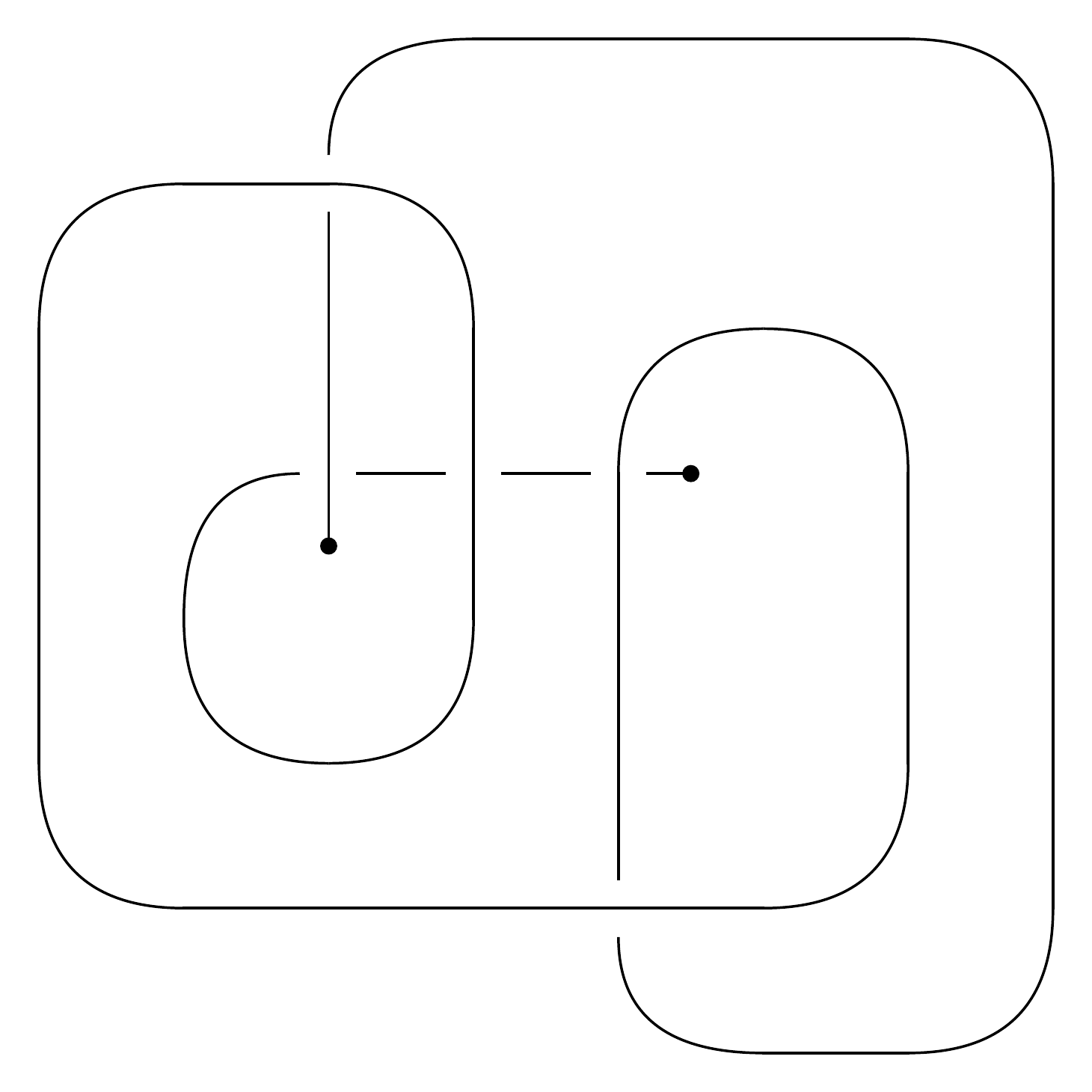}\\
\textcolor{black}{$5_{162}$}
\vspace{1cm}
\end{minipage}
\begin{minipage}[t]{.25\linewidth}
\centering
\includegraphics[width=0.9\textwidth,height=3.5cm,keepaspectratio]{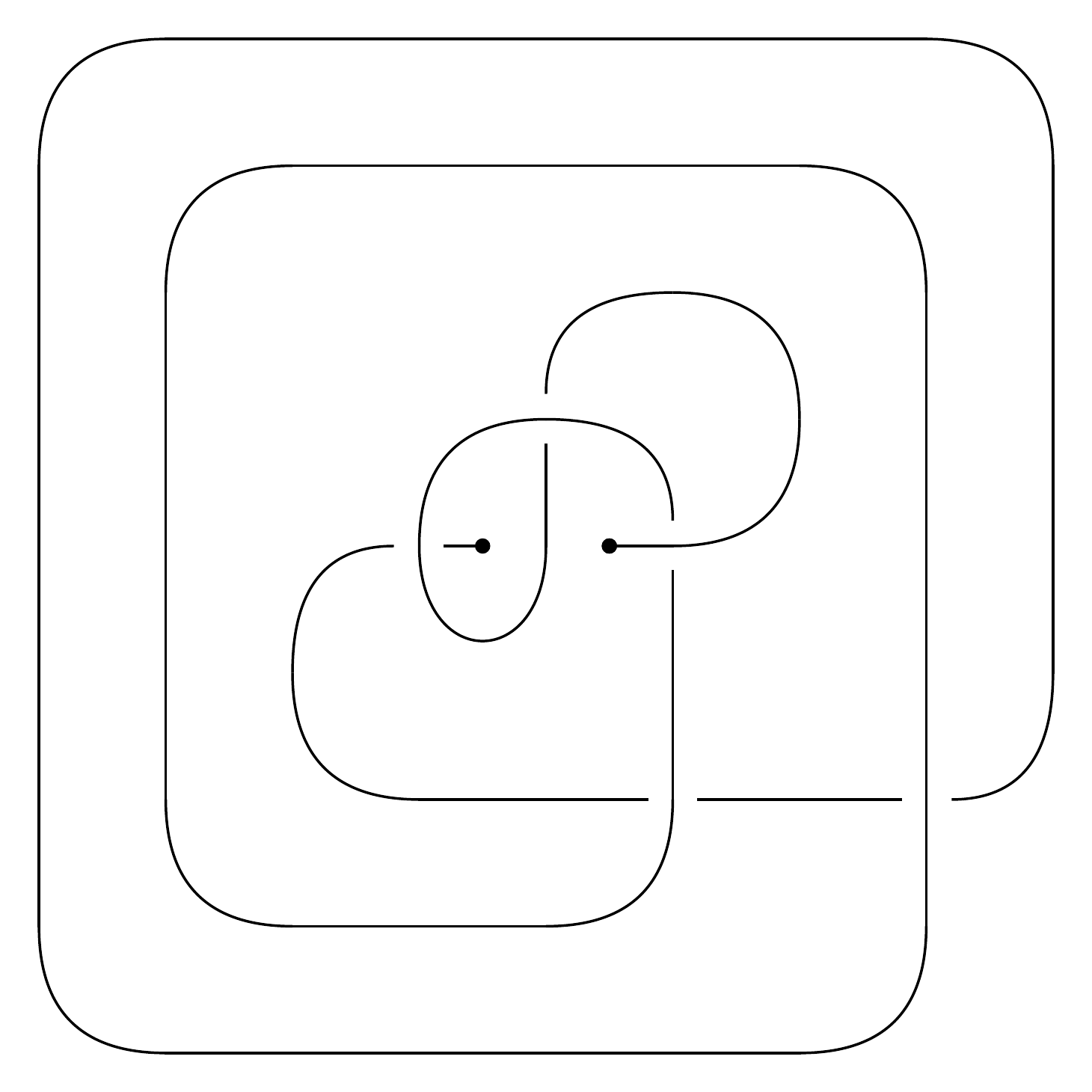}\\
\textcolor{black}{$5_{163}$}
\vspace{1cm}
\end{minipage}
\begin{minipage}[t]{.25\linewidth}
\centering
\includegraphics[width=0.9\textwidth,height=3.5cm,keepaspectratio]{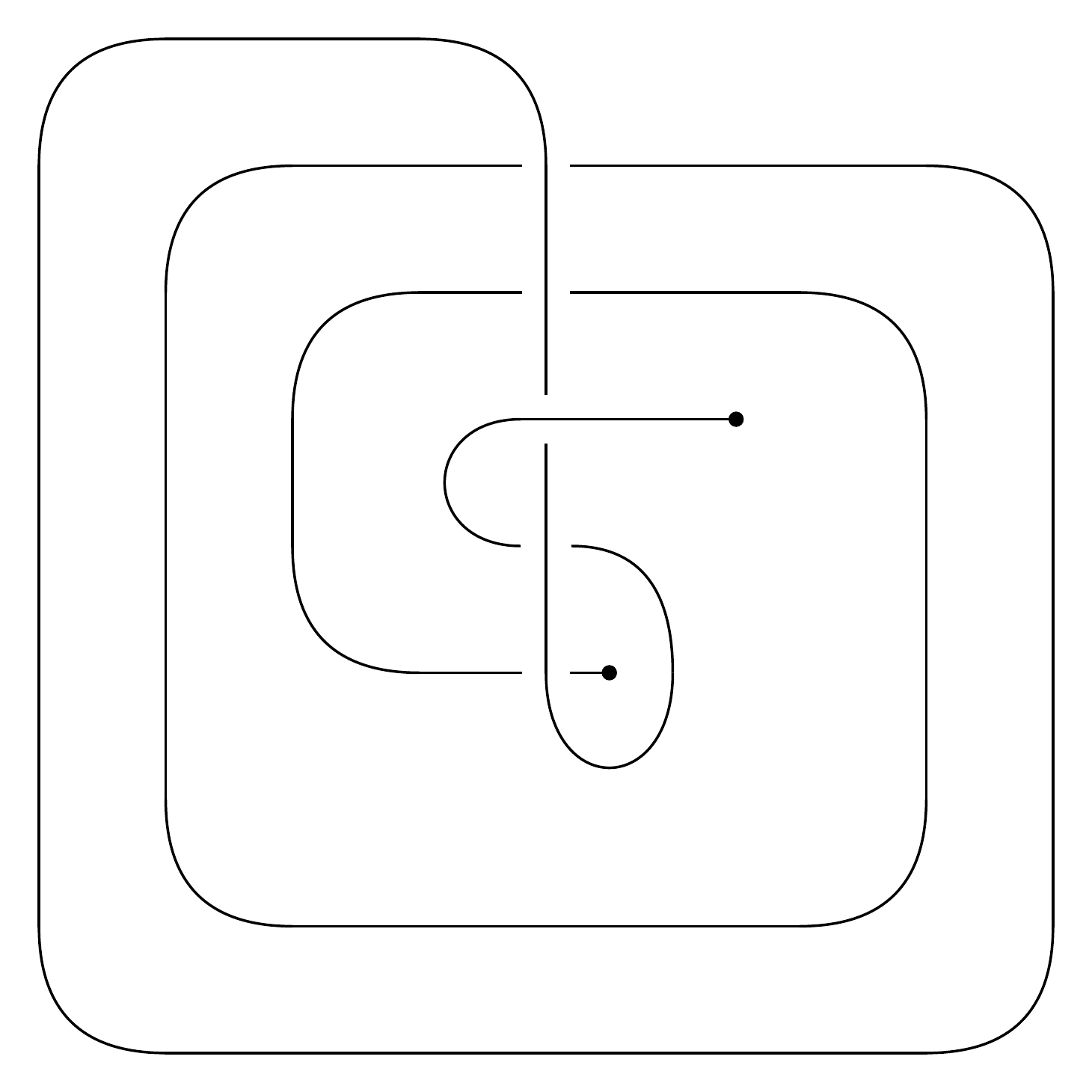}\\
\textcolor{black}{$5_{164}$}
\vspace{1cm}
\end{minipage}
\begin{minipage}[t]{.25\linewidth}
\centering
\includegraphics[width=0.9\textwidth,height=3.5cm,keepaspectratio]{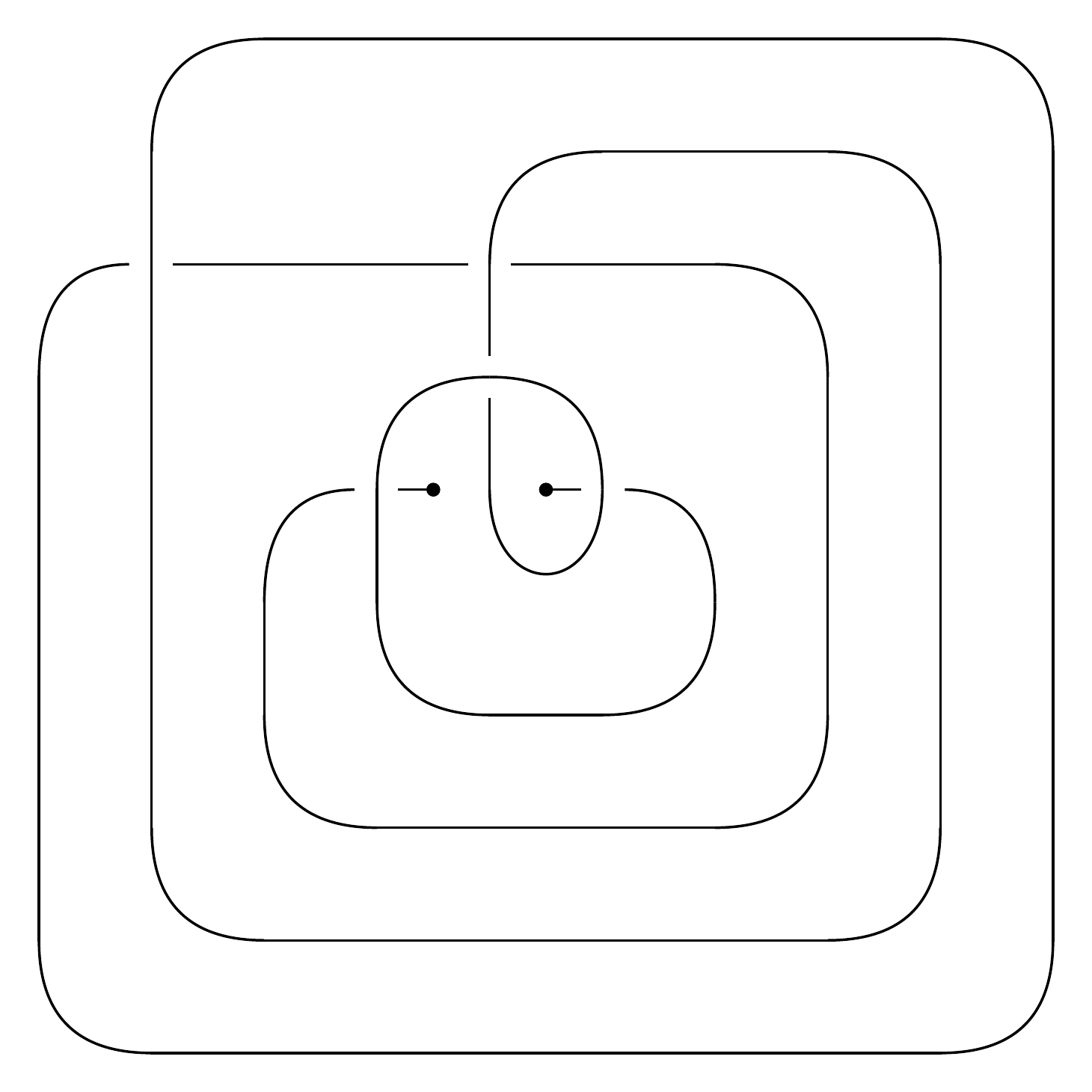}\\
\textcolor{black}{$5_{165}$}
\vspace{1cm}
\end{minipage}
\begin{minipage}[t]{.25\linewidth}
\centering
\includegraphics[width=0.9\textwidth,height=3.5cm,keepaspectratio]{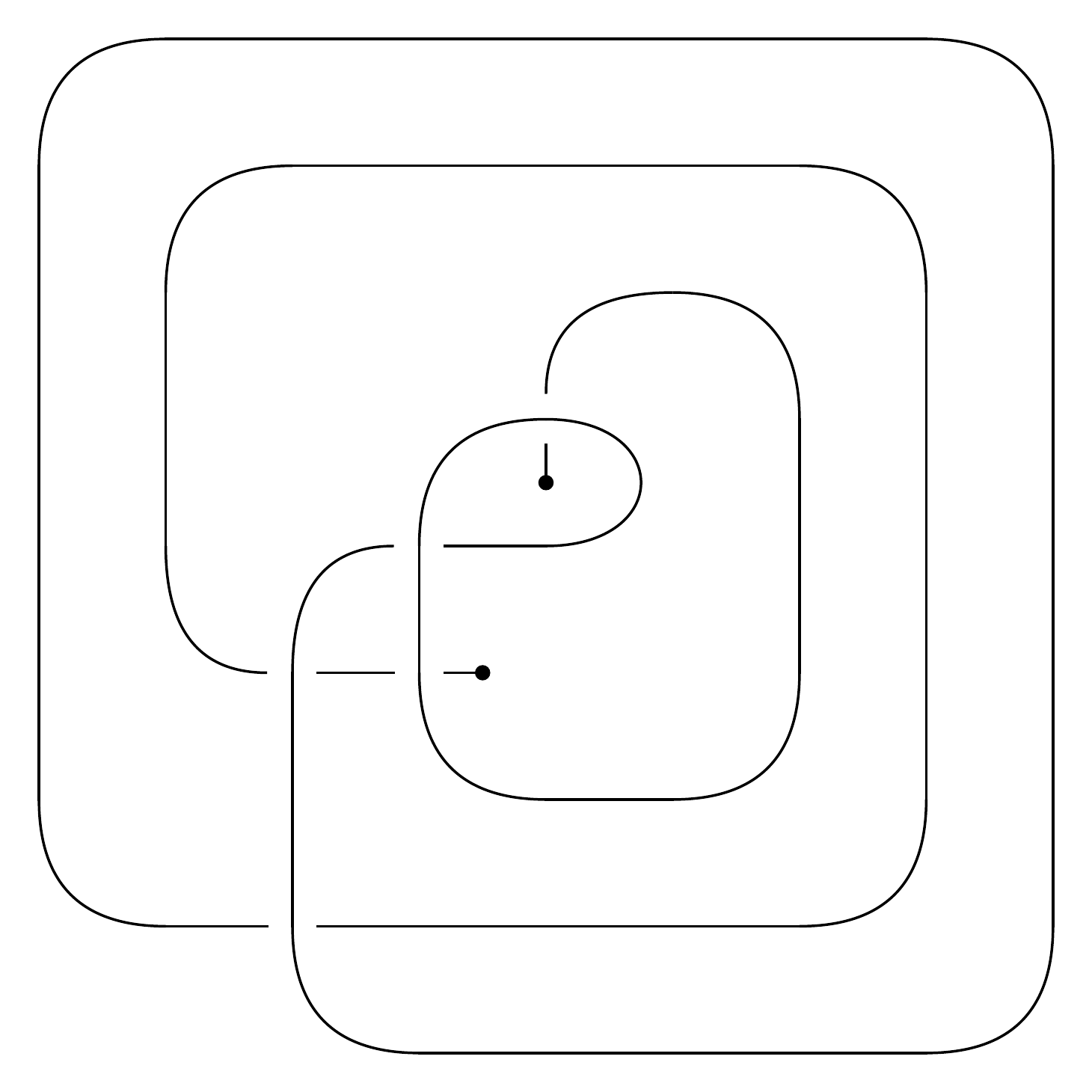}\\
\textcolor{black}{$5_{166}$}
\vspace{1cm}
\end{minipage}
\begin{minipage}[t]{.25\linewidth}
\centering
\includegraphics[width=0.9\textwidth,height=3.5cm,keepaspectratio]{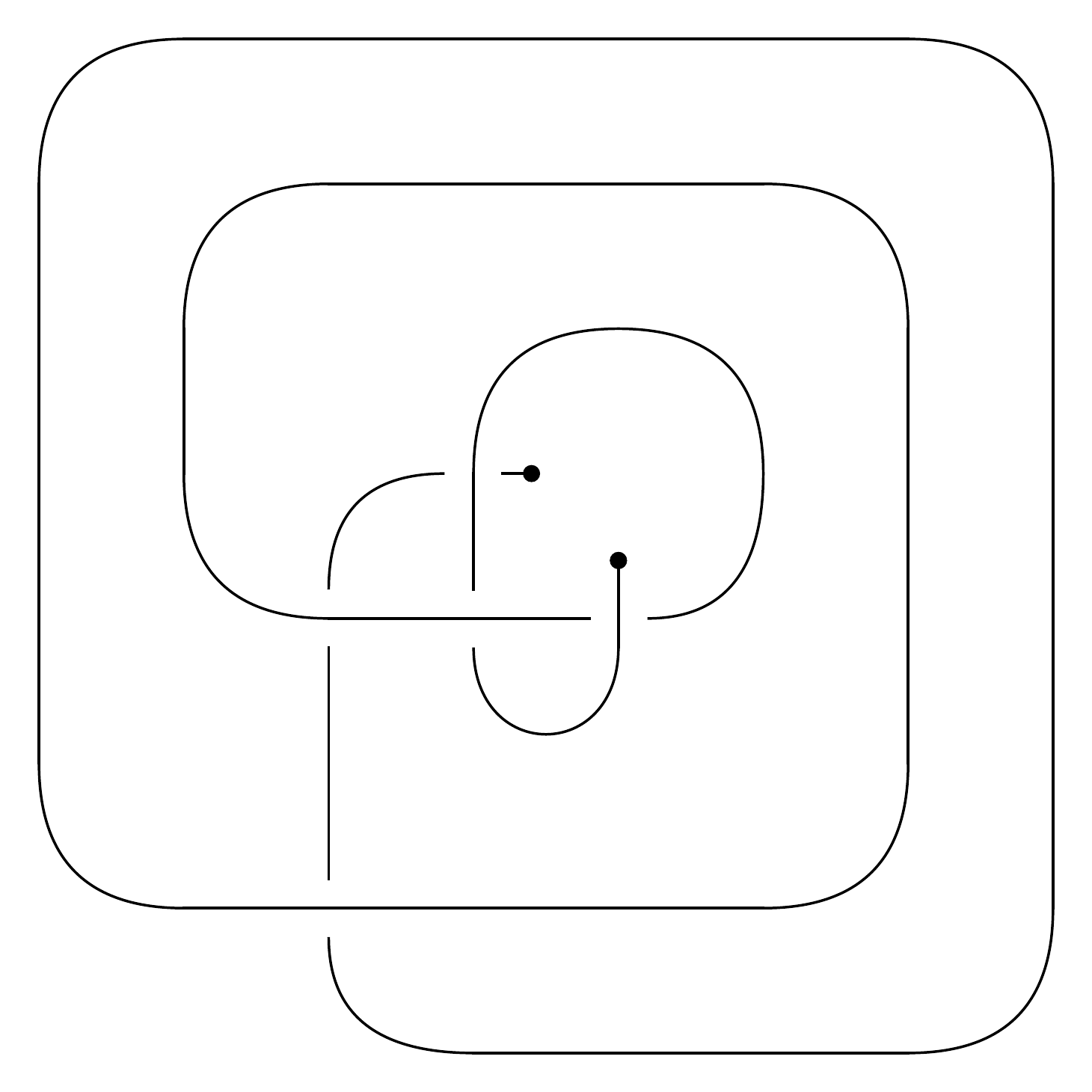}\\
\textcolor{black}{$5_{167}$}
\vspace{1cm}
\end{minipage}
\begin{minipage}[t]{.25\linewidth}
\centering
\includegraphics[width=0.9\textwidth,height=3.5cm,keepaspectratio]{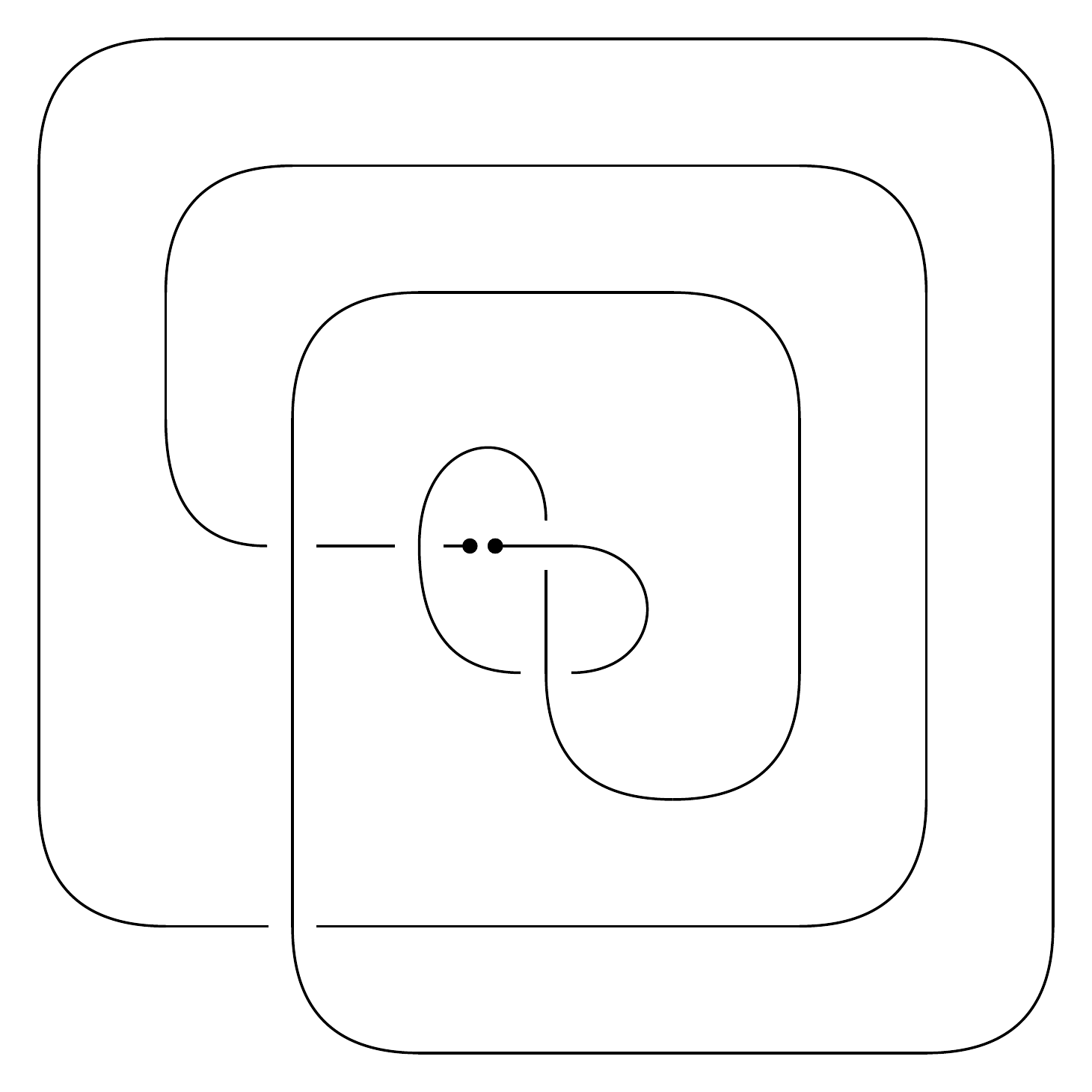}\\
\textcolor{black}{$5_{168}$}
\vspace{1cm}
\end{minipage}
\begin{minipage}[t]{.25\linewidth}
\centering
\includegraphics[width=0.9\textwidth,height=3.5cm,keepaspectratio]{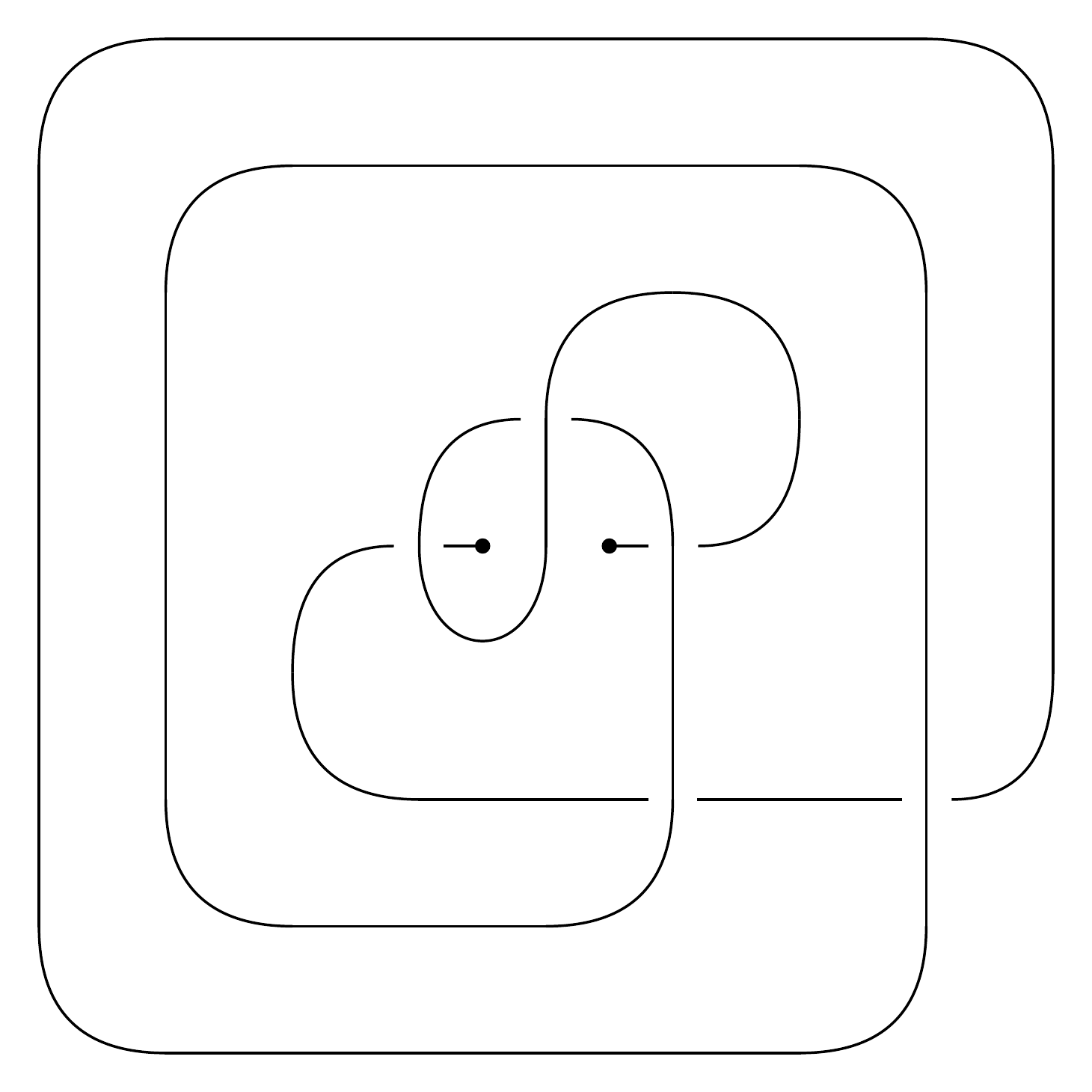}\\
\textcolor{black}{$5_{169}$}
\vspace{1cm}
\end{minipage}
\begin{minipage}[t]{.25\linewidth}
\centering
\includegraphics[width=0.9\textwidth,height=3.5cm,keepaspectratio]{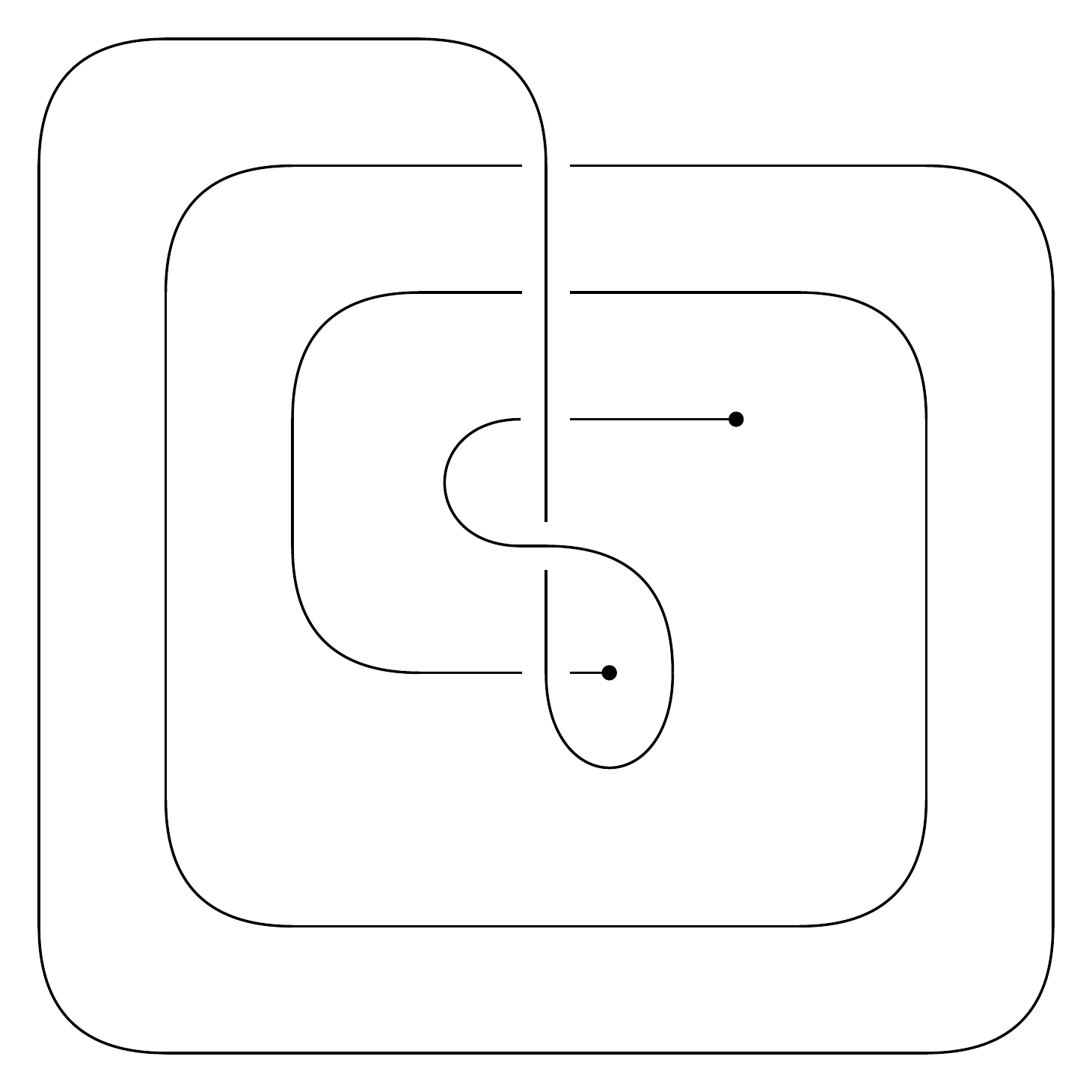}\\
\textcolor{black}{$5_{170}$}
\vspace{1cm}
\end{minipage}
\begin{minipage}[t]{.25\linewidth}
\centering
\includegraphics[width=0.9\textwidth,height=3.5cm,keepaspectratio]{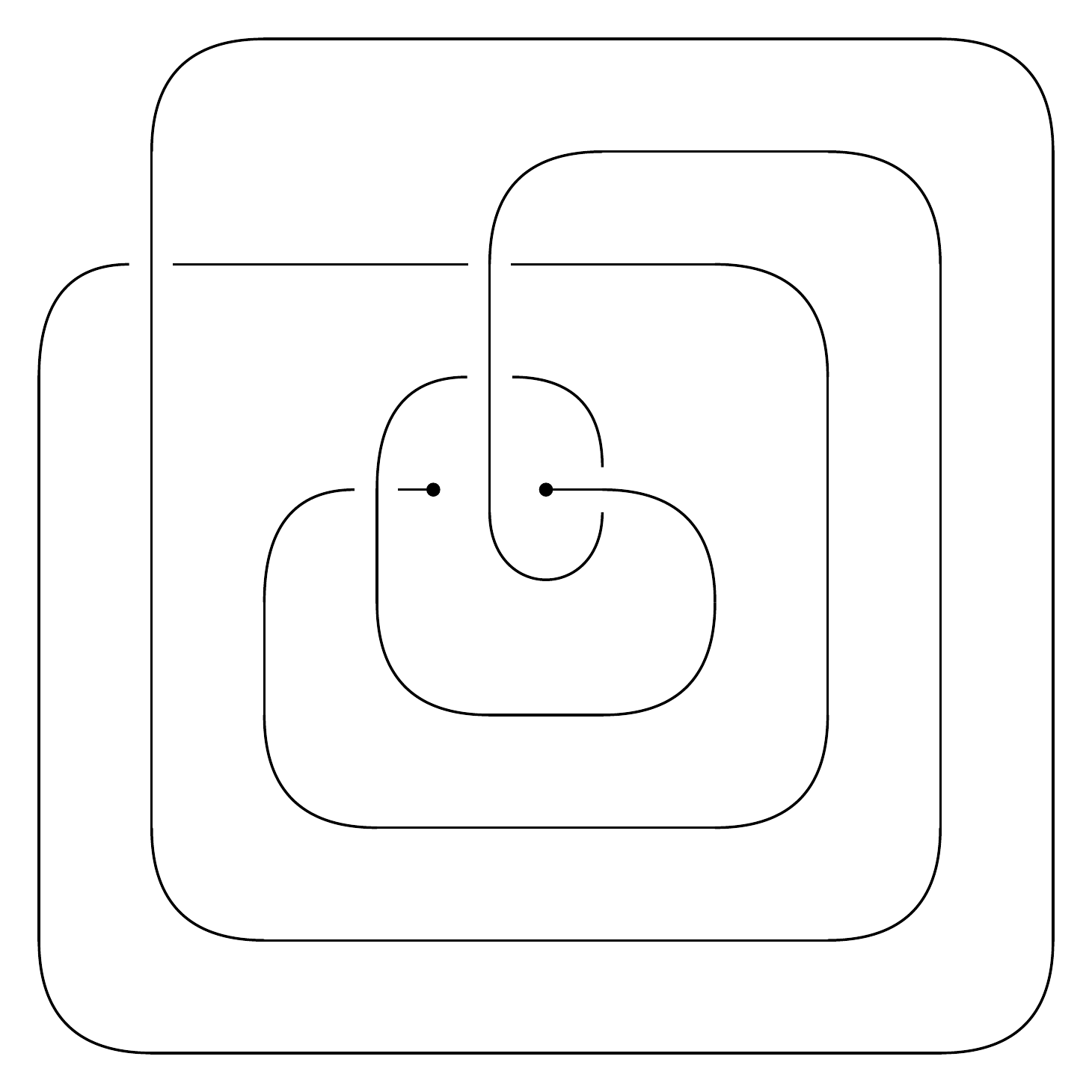}\\
\textcolor{black}{$5_{171}$}
\vspace{1cm}
\end{minipage}
\begin{minipage}[t]{.25\linewidth}
\centering
\includegraphics[width=0.9\textwidth,height=3.5cm,keepaspectratio]{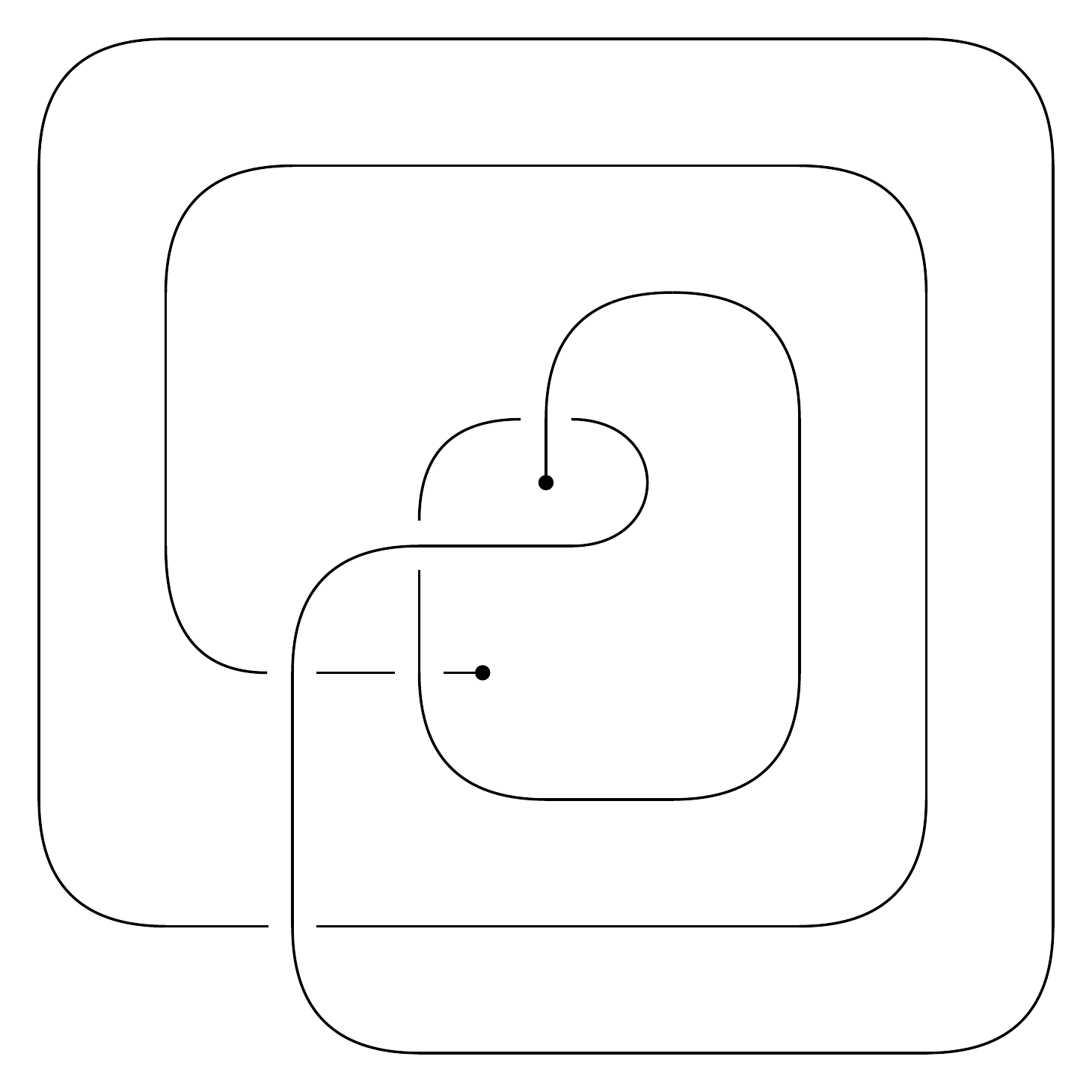}\\
\textcolor{black}{$5_{172}$}
\vspace{1cm}
\end{minipage}
\begin{minipage}[t]{.25\linewidth}
\centering
\includegraphics[width=0.9\textwidth,height=3.5cm,keepaspectratio]{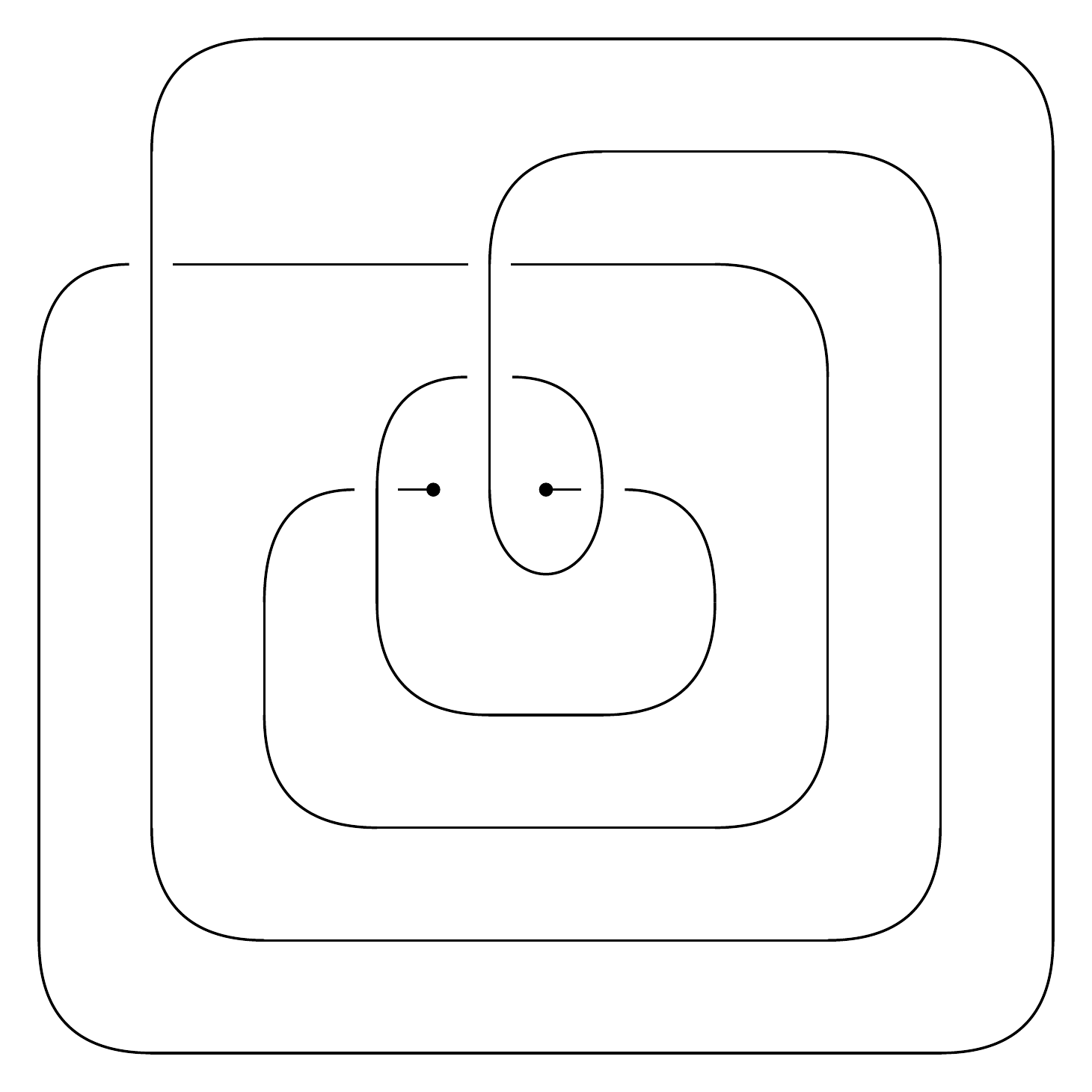}\\
\textcolor{black}{$5_{173}$}
\vspace{1cm}
\end{minipage}
\begin{minipage}[t]{.25\linewidth}
\centering
\includegraphics[width=0.9\textwidth,height=3.5cm,keepaspectratio]{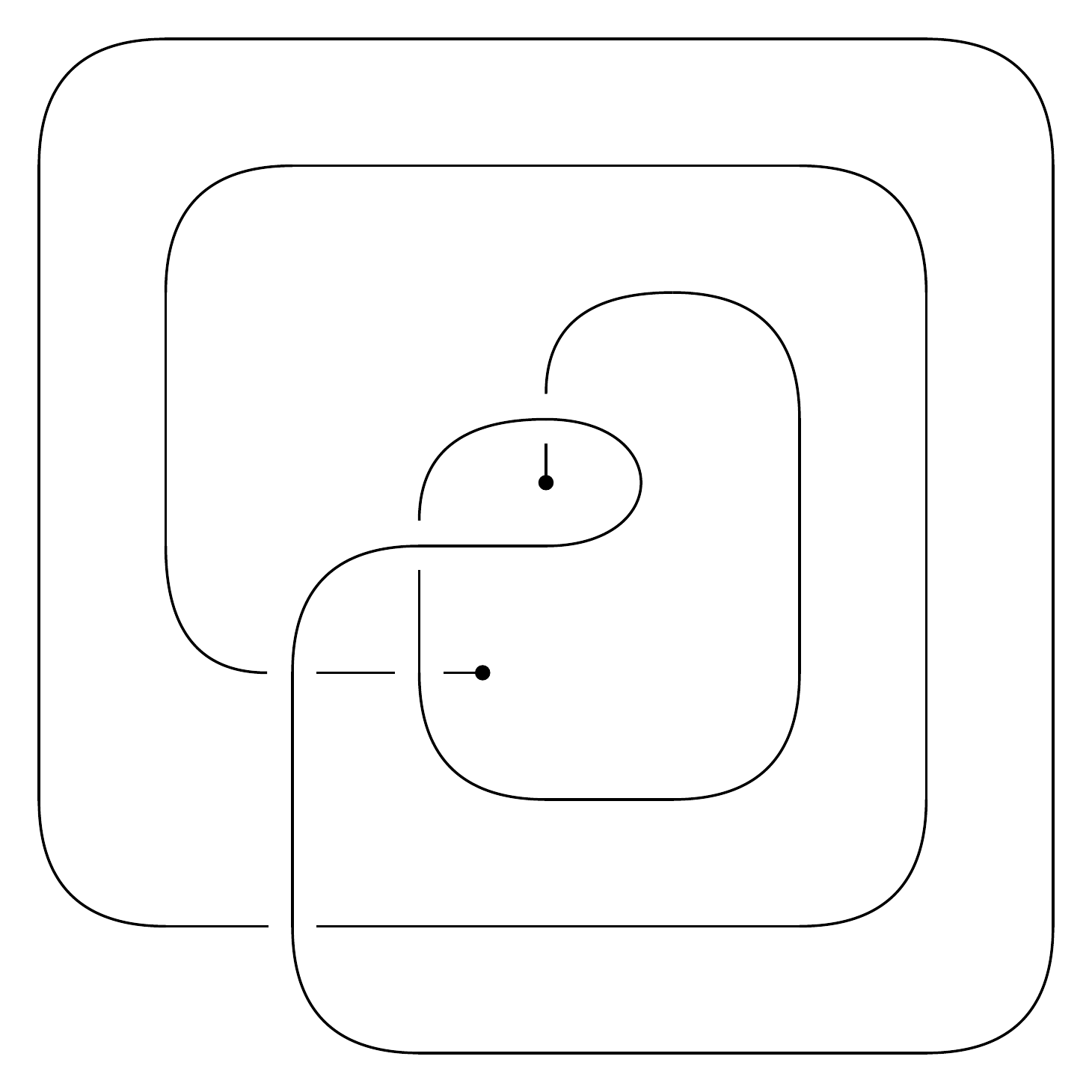}\\
\textcolor{black}{$5_{174}$}
\vspace{1cm}
\end{minipage}
\begin{minipage}[t]{.25\linewidth}
\centering
\includegraphics[width=0.9\textwidth,height=3.5cm,keepaspectratio]{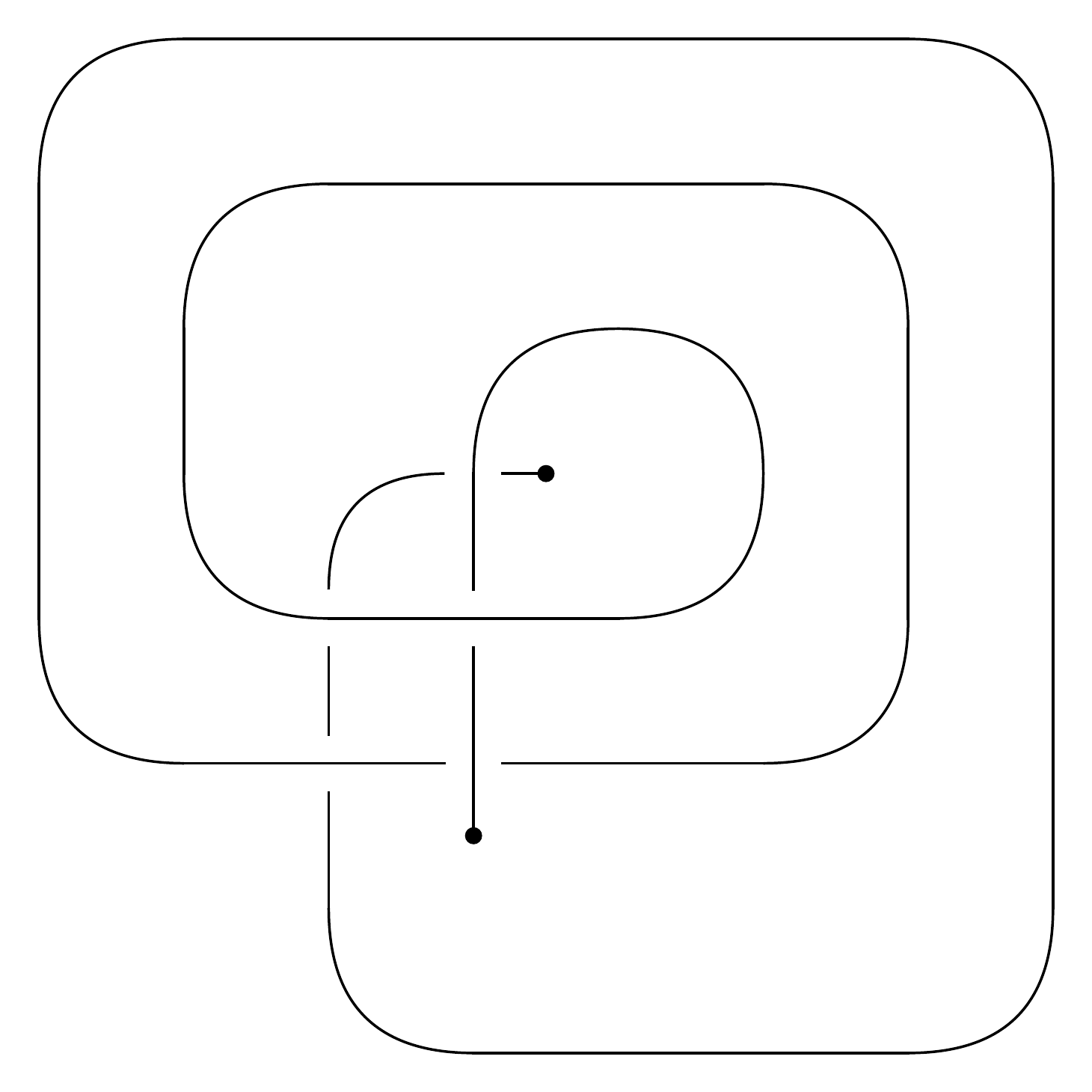}\\
\textcolor{black}{$5_{175}$}
\vspace{1cm}
\end{minipage}
\begin{minipage}[t]{.25\linewidth}
\centering
\includegraphics[width=0.9\textwidth,height=3.5cm,keepaspectratio]{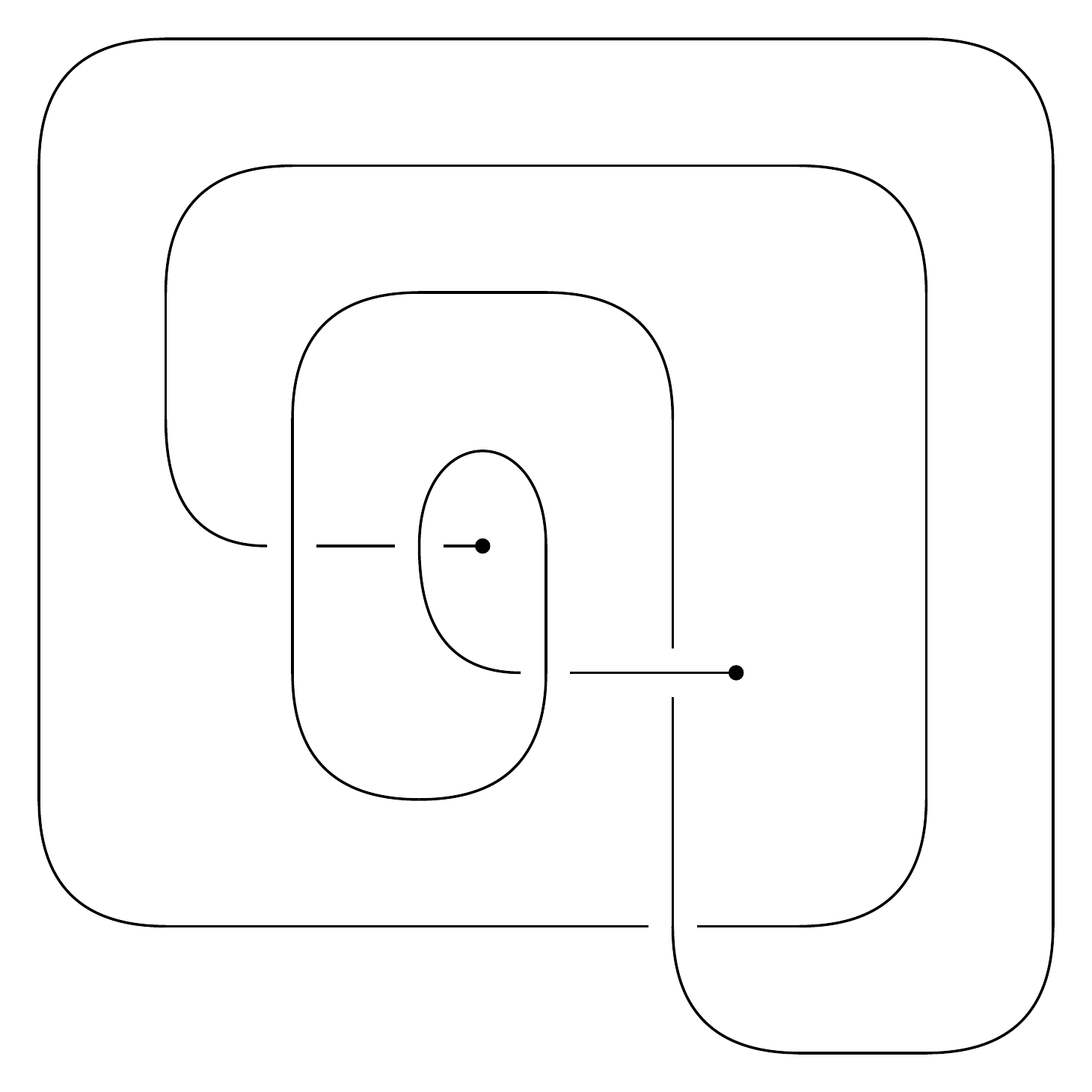}\\
\textcolor{black}{$5_{176}$}
\vspace{1cm}
\end{minipage}
\begin{minipage}[t]{.25\linewidth}
\centering
\includegraphics[width=0.9\textwidth,height=3.5cm,keepaspectratio]{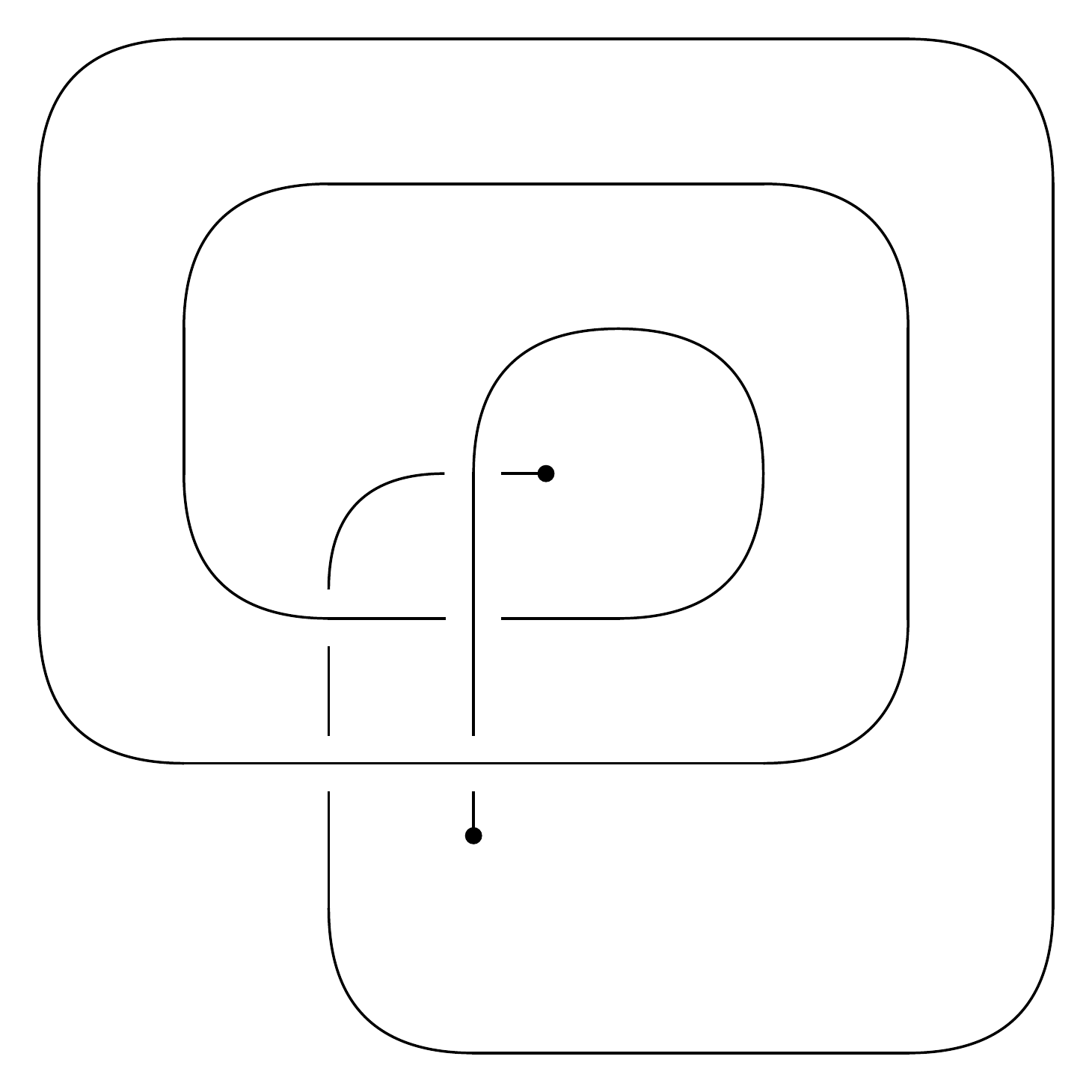}\\
\textcolor{black}{$5_{177}$}
\vspace{1cm}
\end{minipage}
\begin{minipage}[t]{.25\linewidth}
\centering
\includegraphics[width=0.9\textwidth,height=3.5cm,keepaspectratio]{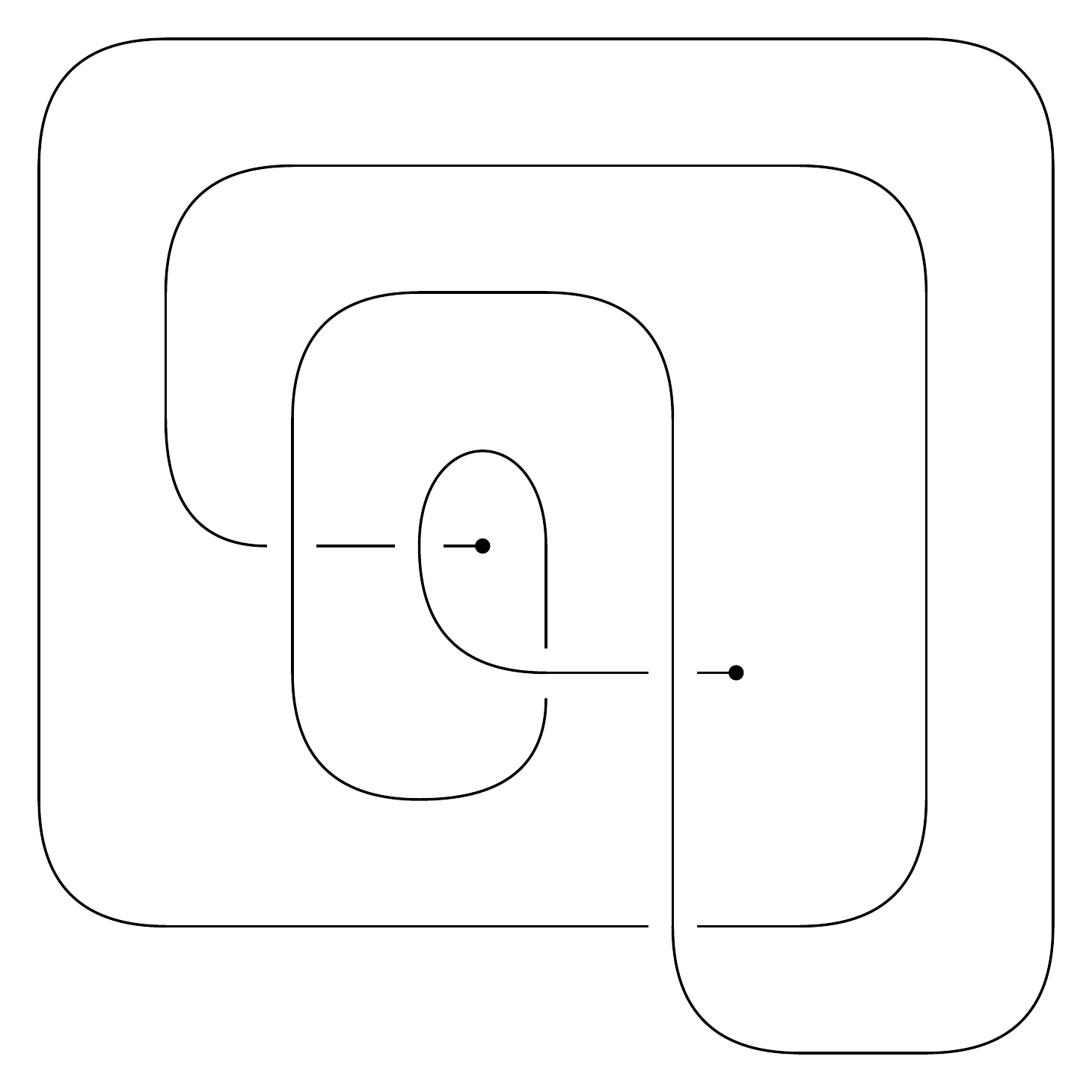}\\
\textcolor{black}{$5_{178}$}
\vspace{1cm}
\end{minipage}
\begin{minipage}[t]{.25\linewidth}
\centering
\includegraphics[width=0.9\textwidth,height=3.5cm,keepaspectratio]{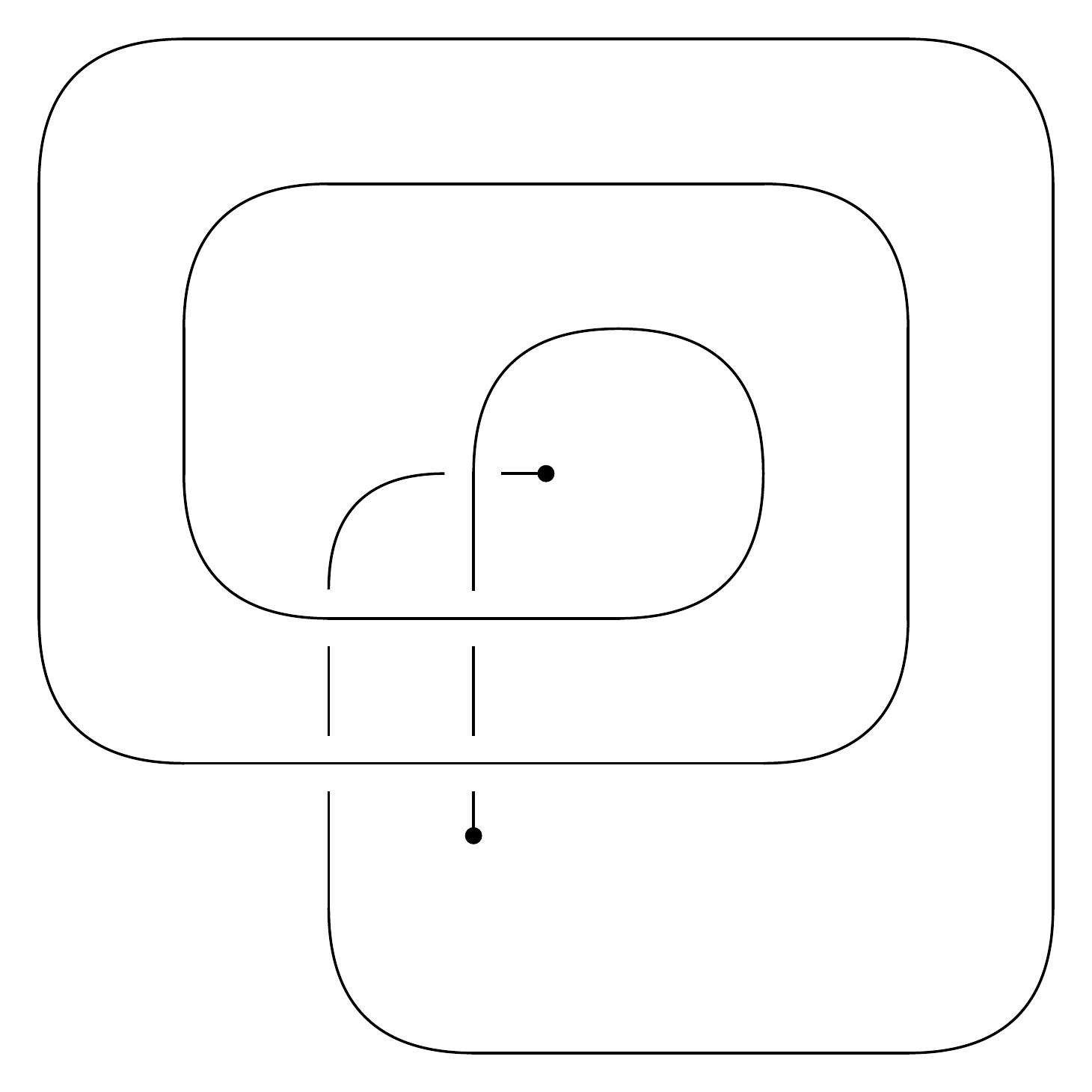}\\
\textcolor{black}{$5_{179}$}
\vspace{1cm}
\end{minipage}
\begin{minipage}[t]{.25\linewidth}
\centering
\includegraphics[width=0.9\textwidth,height=3.5cm,keepaspectratio]{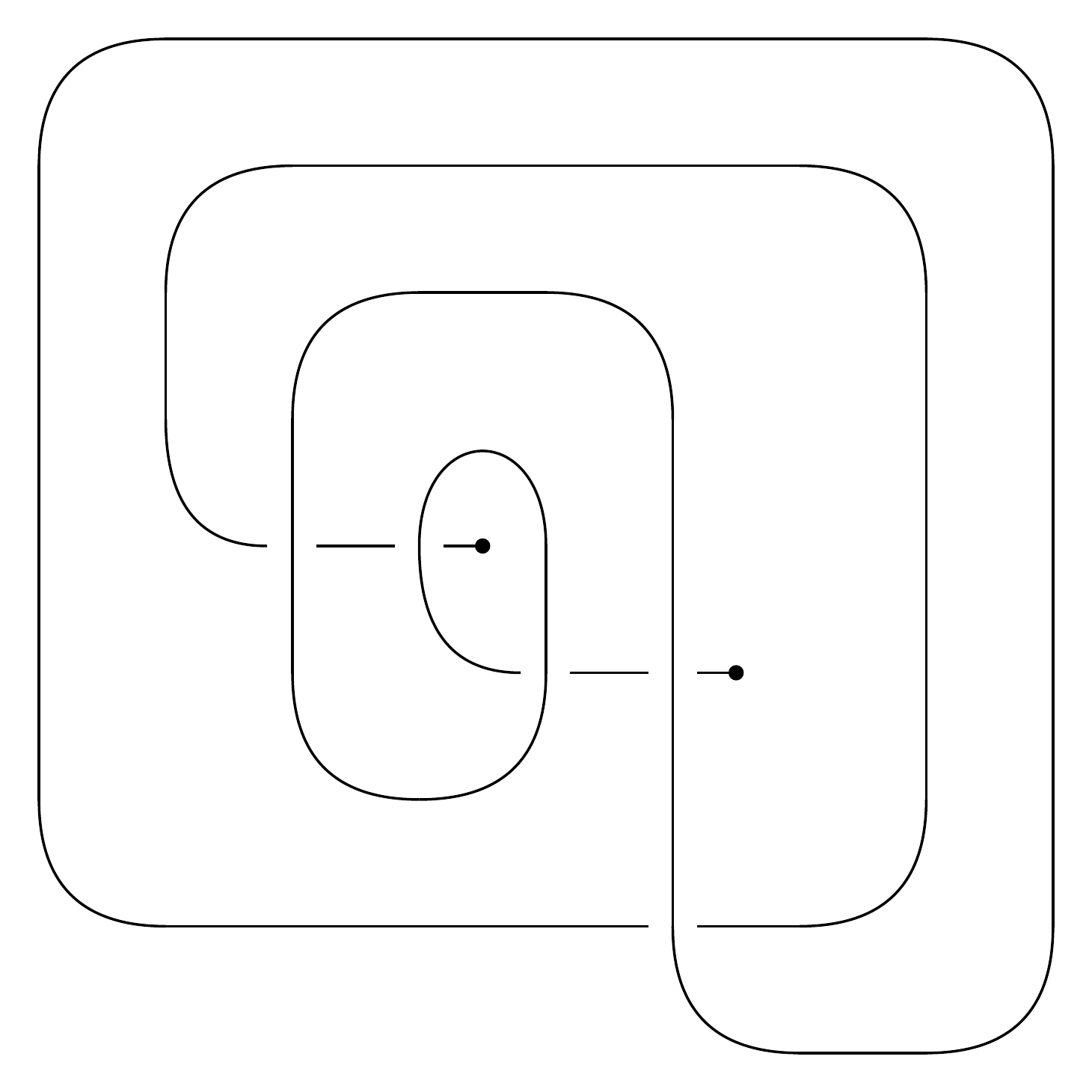}\\
\textcolor{black}{$5_{180}$}
\vspace{1cm}
\end{minipage}
\begin{minipage}[t]{.25\linewidth}
\centering
\includegraphics[width=0.9\textwidth,height=3.5cm,keepaspectratio]{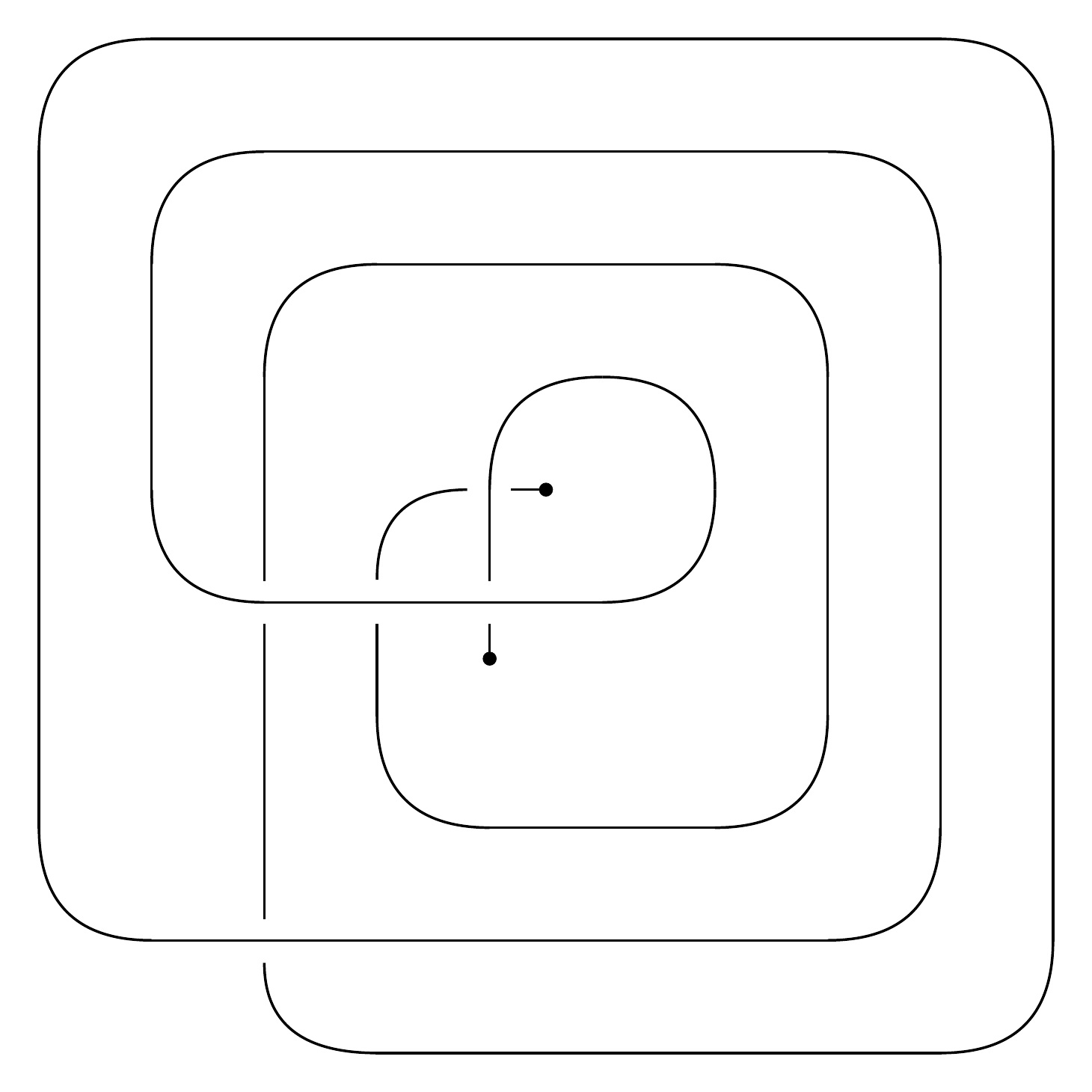}\\
\textcolor{black}{$5_{181}$}
\vspace{1cm}
\end{minipage}
\begin{minipage}[t]{.25\linewidth}
\centering
\includegraphics[width=0.9\textwidth,height=3.5cm,keepaspectratio]{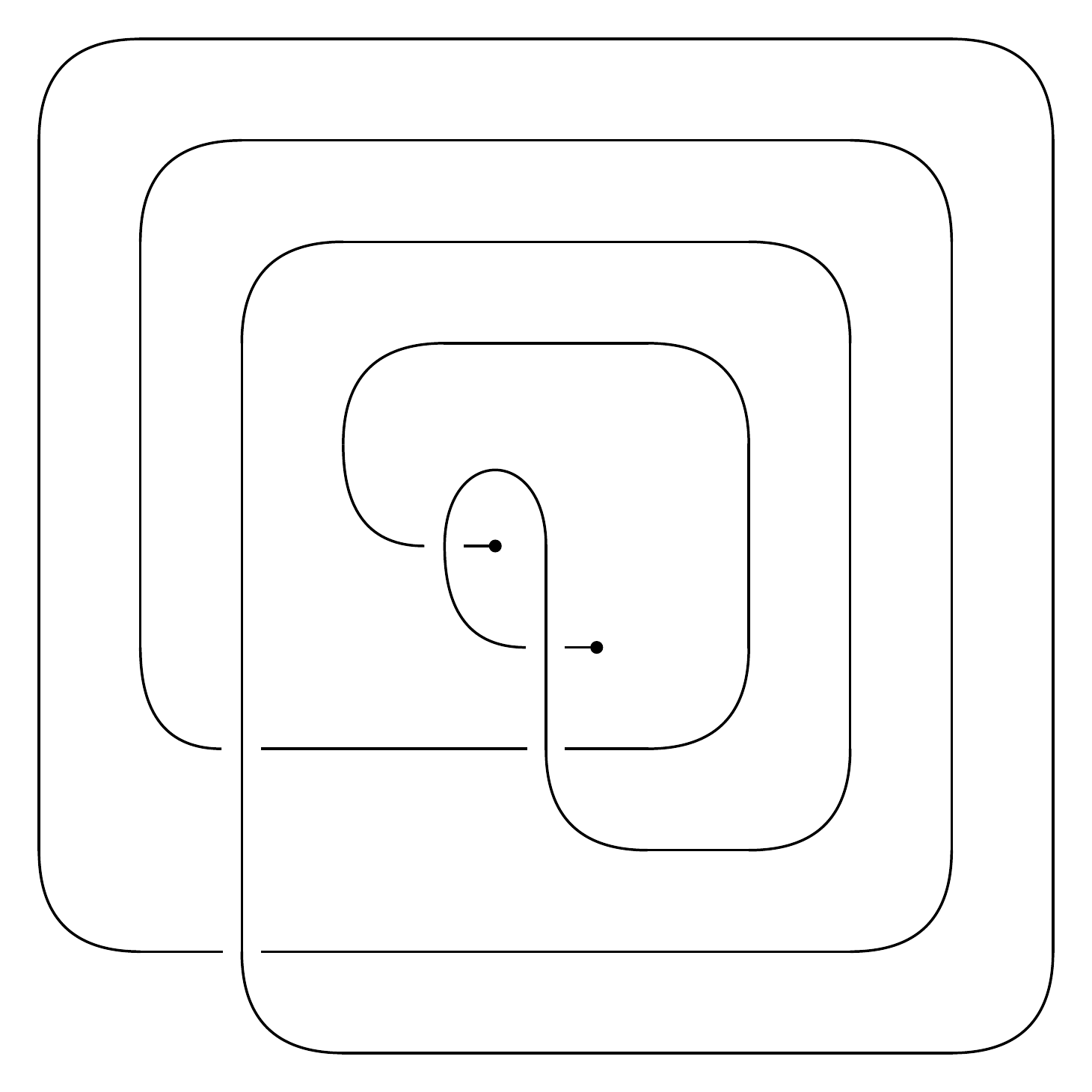}\\
\textcolor{black}{$5_{182}$}
\vspace{1cm}
\end{minipage}
\begin{minipage}[t]{.25\linewidth}
\centering
\includegraphics[width=0.9\textwidth,height=3.5cm,keepaspectratio]{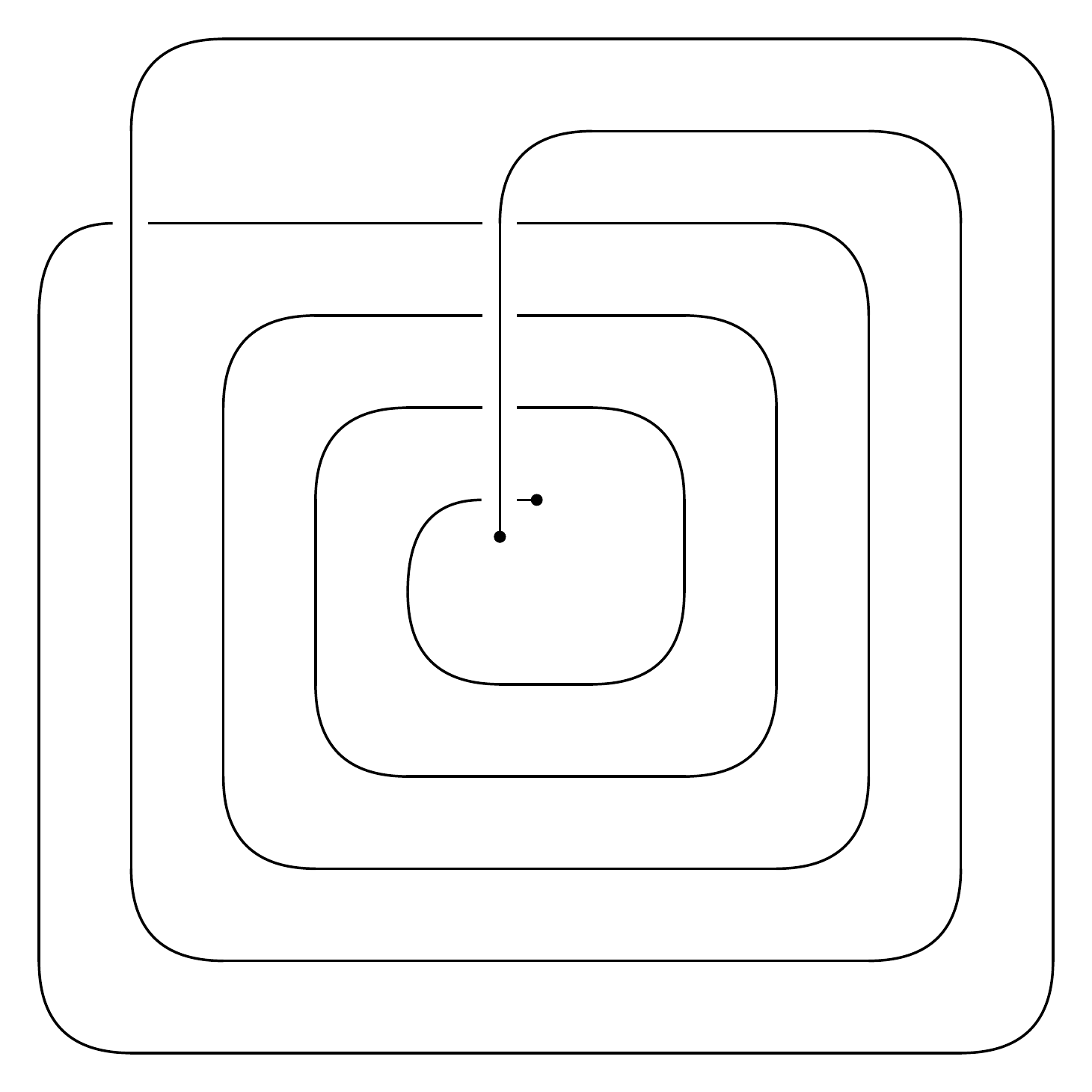}\\
\textcolor{black}{$5_{183}$}
\vspace{1cm}
\end{minipage}
\begin{minipage}[t]{.25\linewidth}
\centering
\includegraphics[width=0.9\textwidth,height=3.5cm,keepaspectratio]{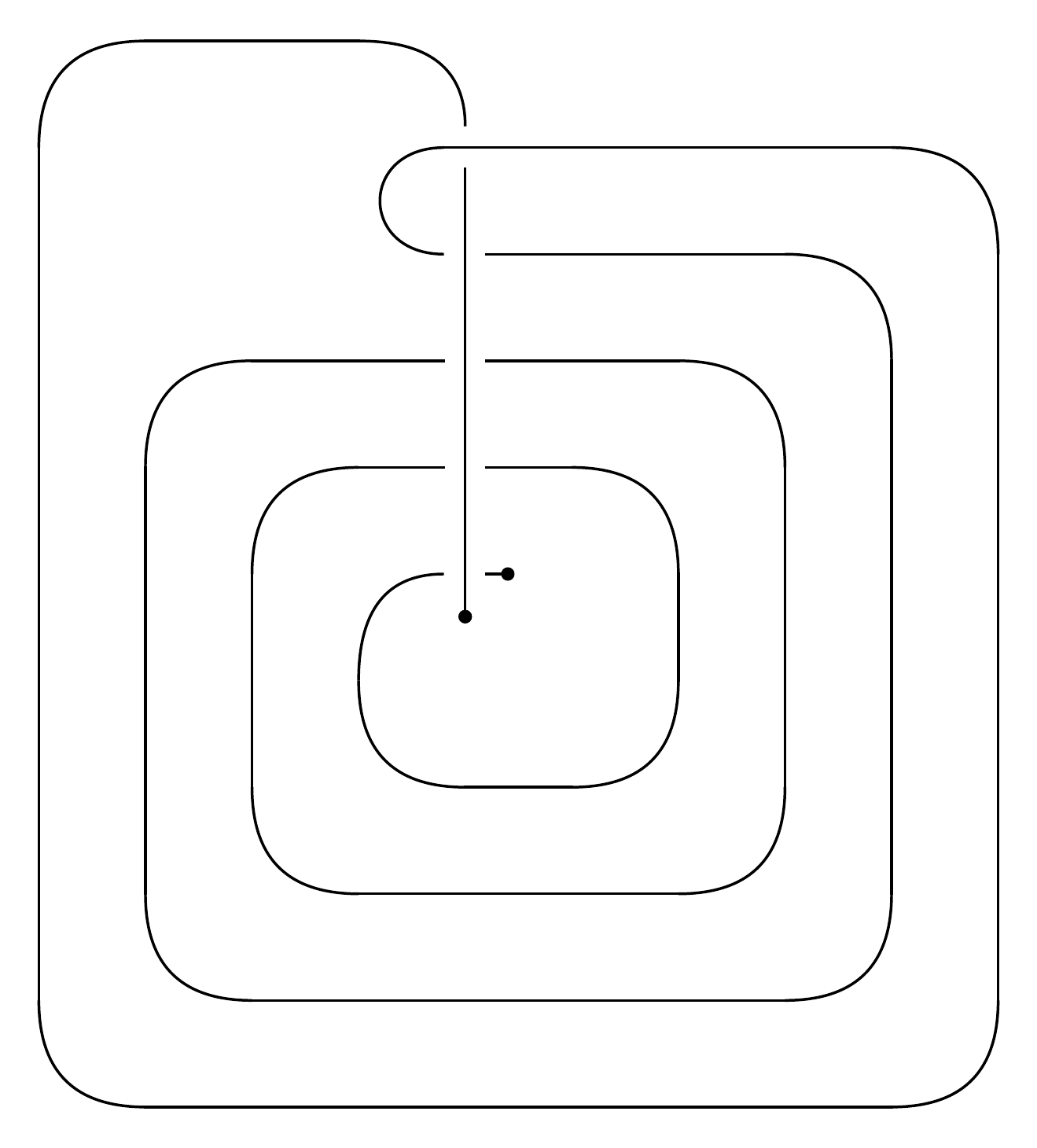}\\
\textcolor{black}{$5_{184}$}
\vspace{1cm}
\end{minipage}
\begin{minipage}[t]{.25\linewidth}
\centering
\includegraphics[width=0.9\textwidth,height=3.5cm,keepaspectratio]{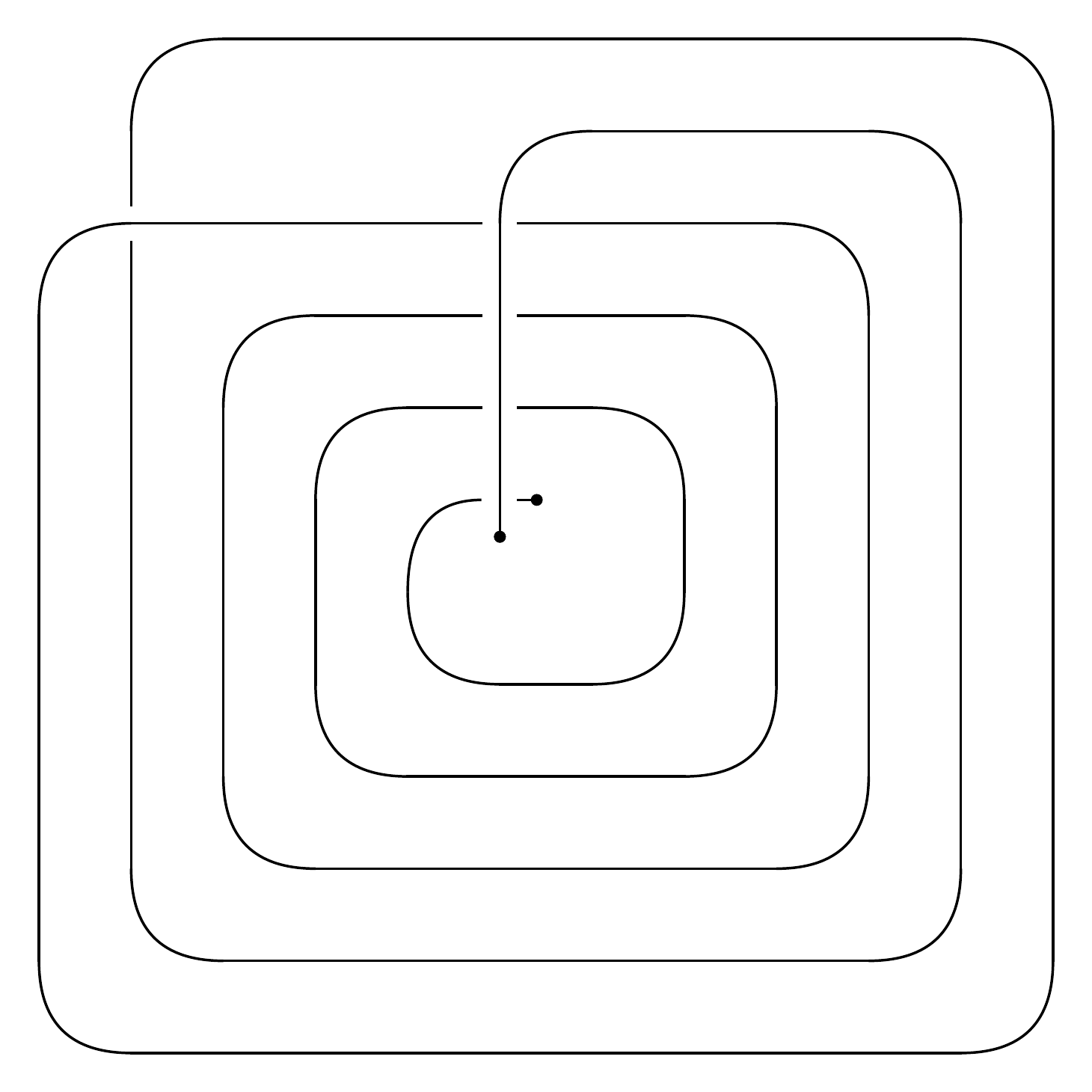}\\
\textcolor{black}{$5_{185}$}
\vspace{1cm}
\end{minipage}
\begin{minipage}[t]{.25\linewidth}
\centering
\includegraphics[width=0.9\textwidth,height=3.5cm,keepaspectratio]{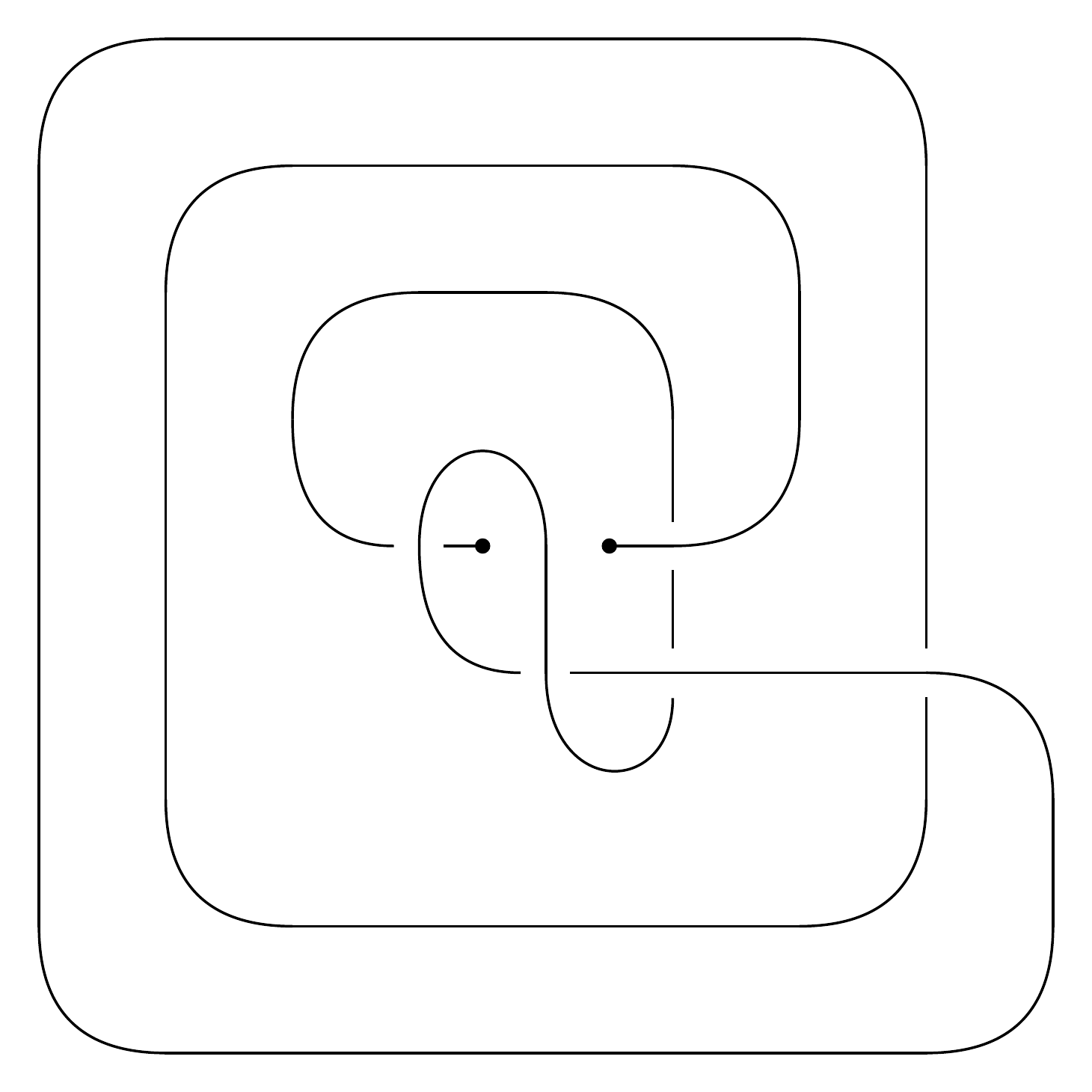}\\
\textcolor{black}{$5_{186}$}
\vspace{1cm}
\end{minipage}
\begin{minipage}[t]{.25\linewidth}
\centering
\includegraphics[width=0.9\textwidth,height=3.5cm,keepaspectratio]{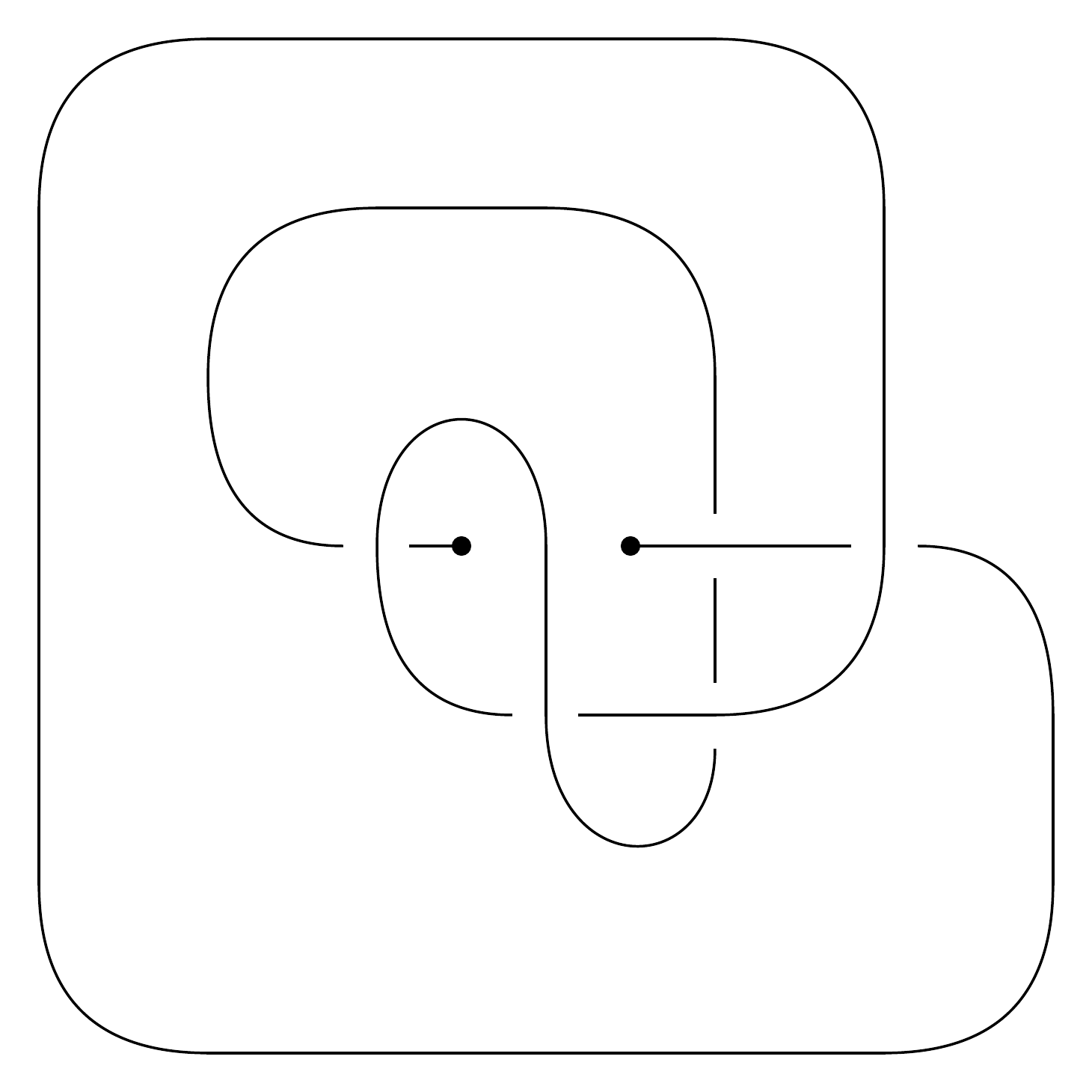}\\
\textcolor{black}{$5_{187}$}
\vspace{1cm}
\end{minipage}
\begin{minipage}[t]{.25\linewidth}
\centering
\includegraphics[width=0.9\textwidth,height=3.5cm,keepaspectratio]{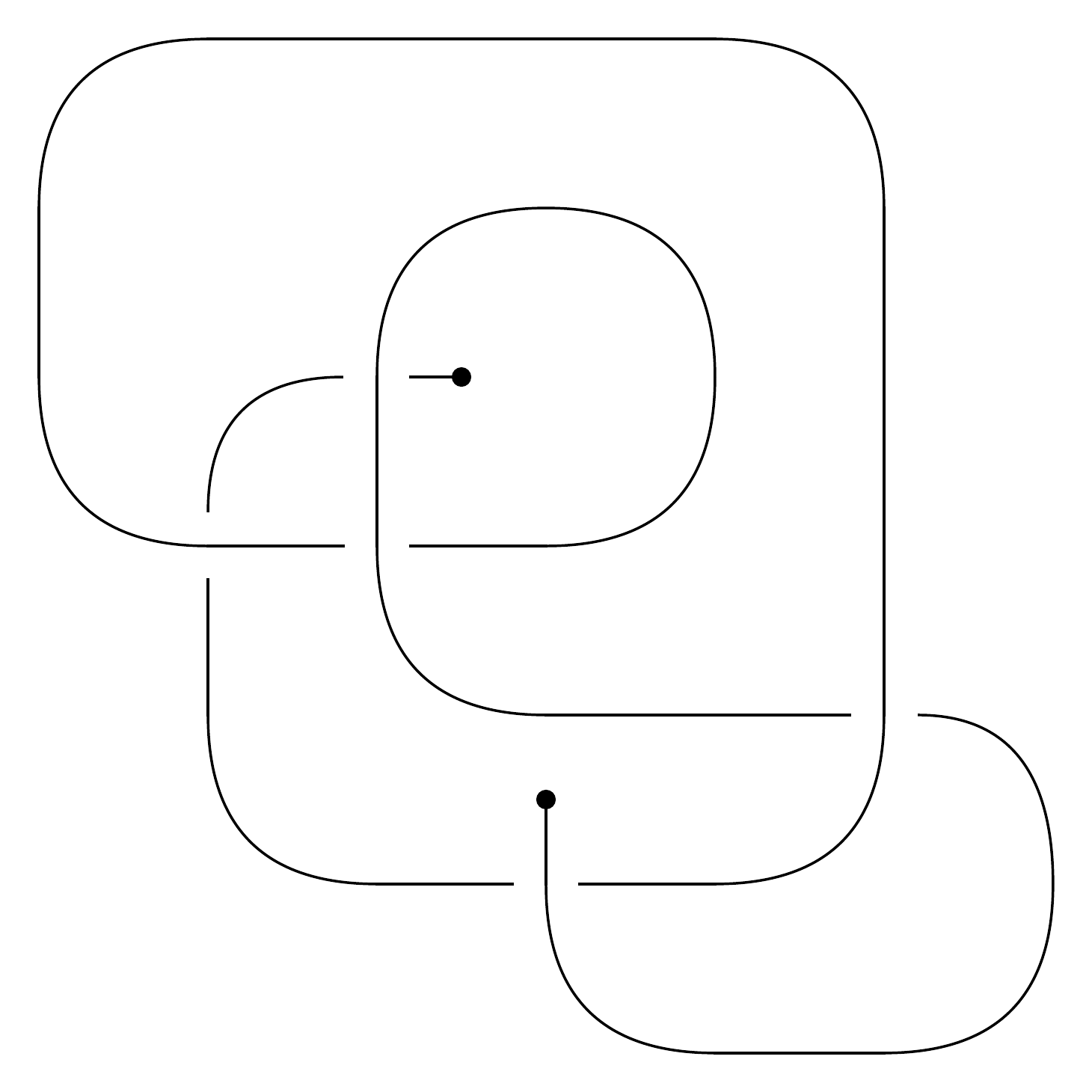}\\
\textcolor{black}{$5_{188}$}
\vspace{1cm}
\end{minipage}
\begin{minipage}[t]{.25\linewidth}
\centering
\includegraphics[width=0.9\textwidth,height=3.5cm,keepaspectratio]{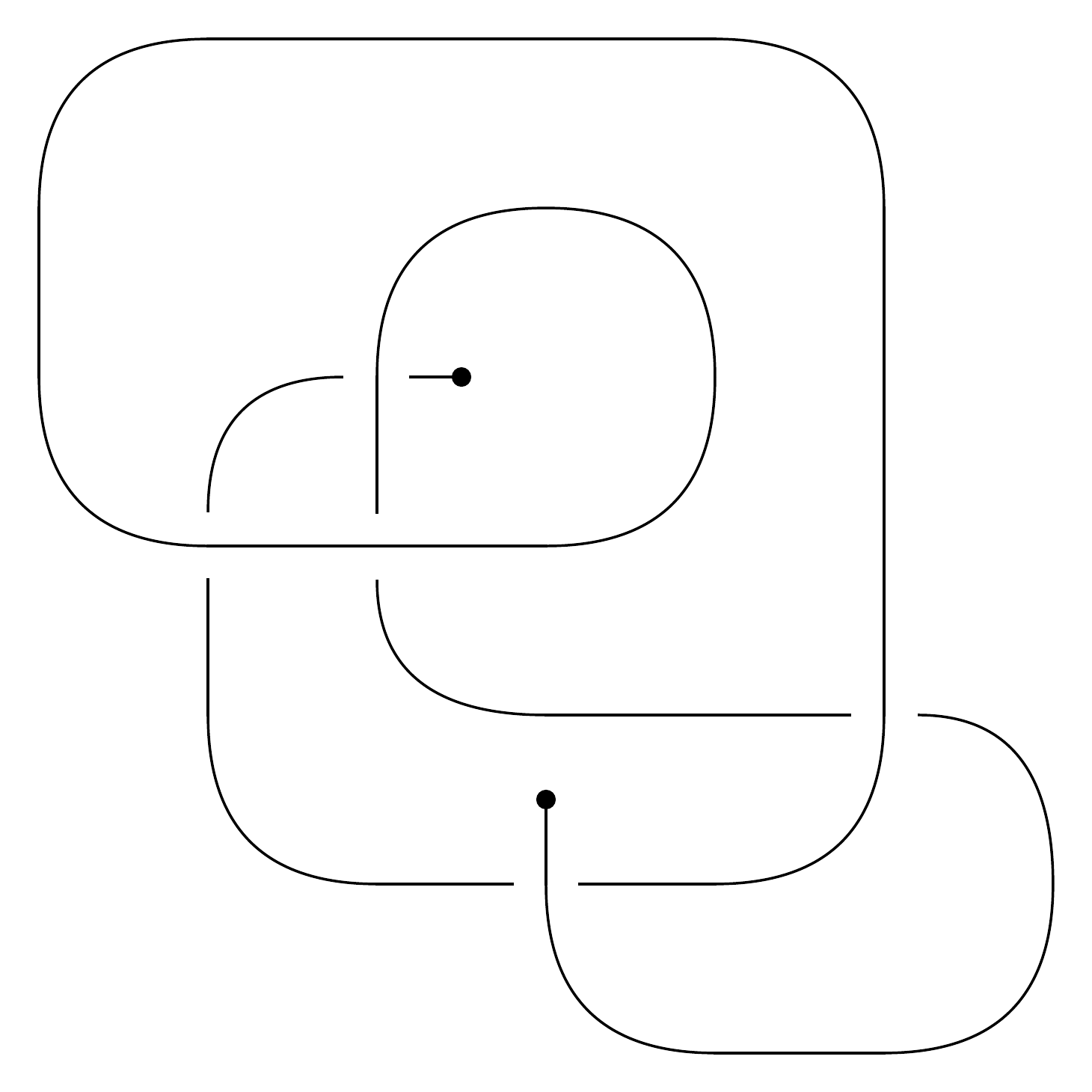}\\
\textcolor{black}{$5_{189}$}
\vspace{1cm}
\end{minipage}
\begin{minipage}[t]{.25\linewidth}
\centering
\includegraphics[width=0.9\textwidth,height=3.5cm,keepaspectratio]{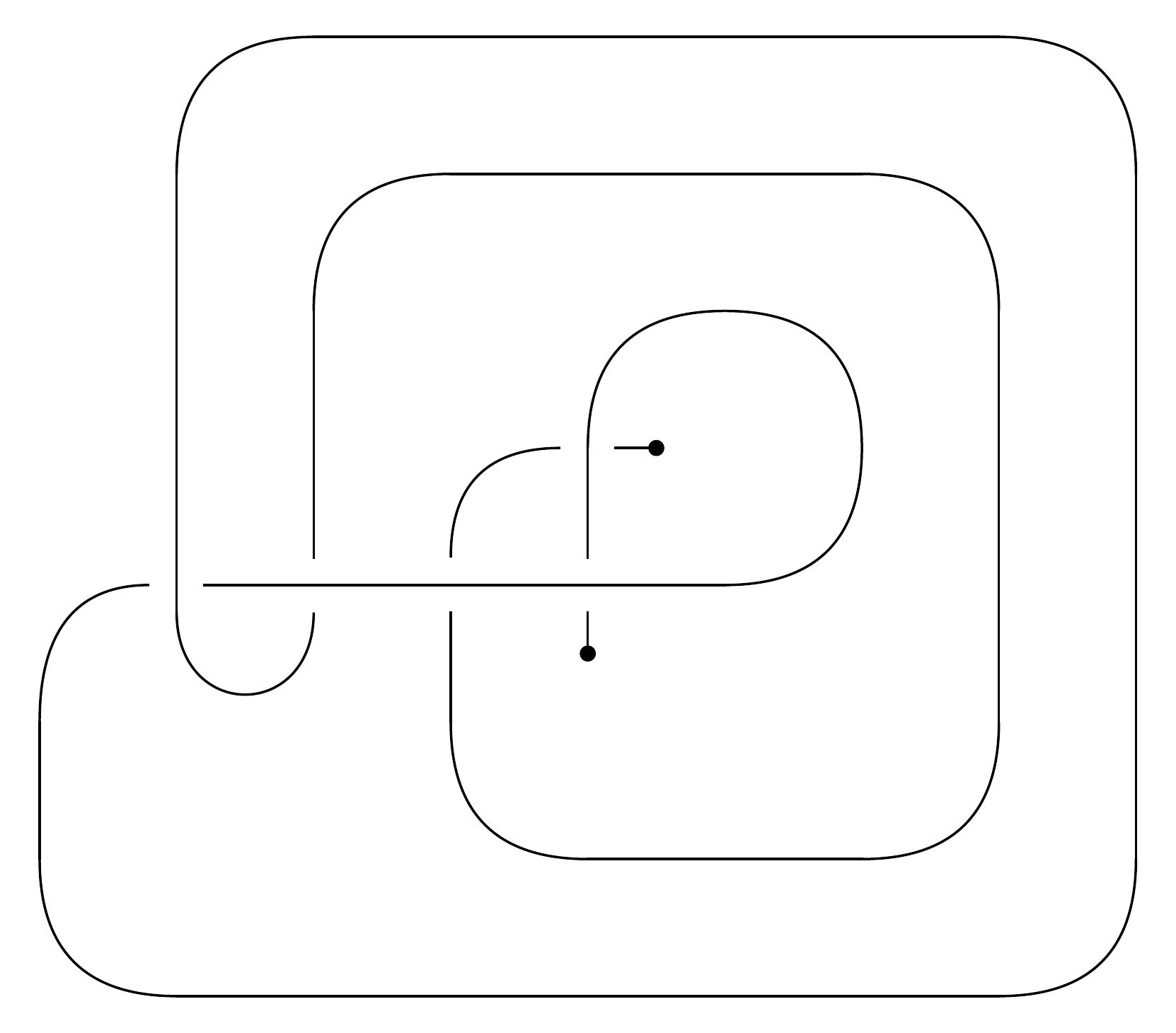}\\
\textcolor{black}{$5_{190}$}
\vspace{1cm}
\end{minipage}
\begin{minipage}[t]{.25\linewidth}
\centering
\includegraphics[width=0.9\textwidth,height=3.5cm,keepaspectratio]{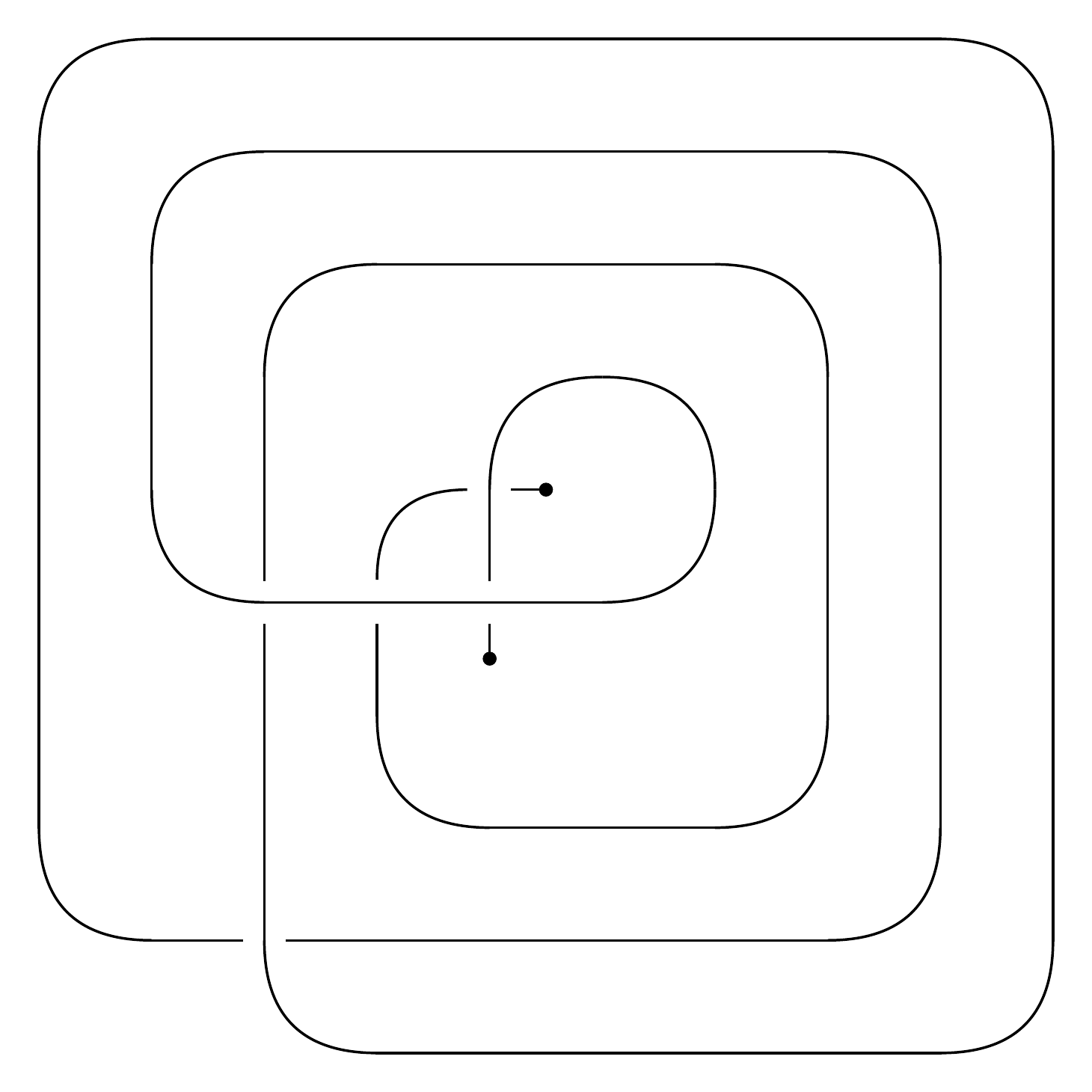}\\
\textcolor{black}{$5_{191}$}
\vspace{1cm}
\end{minipage}
\begin{minipage}[t]{.25\linewidth}
\centering
\includegraphics[width=0.9\textwidth,height=3.5cm,keepaspectratio]{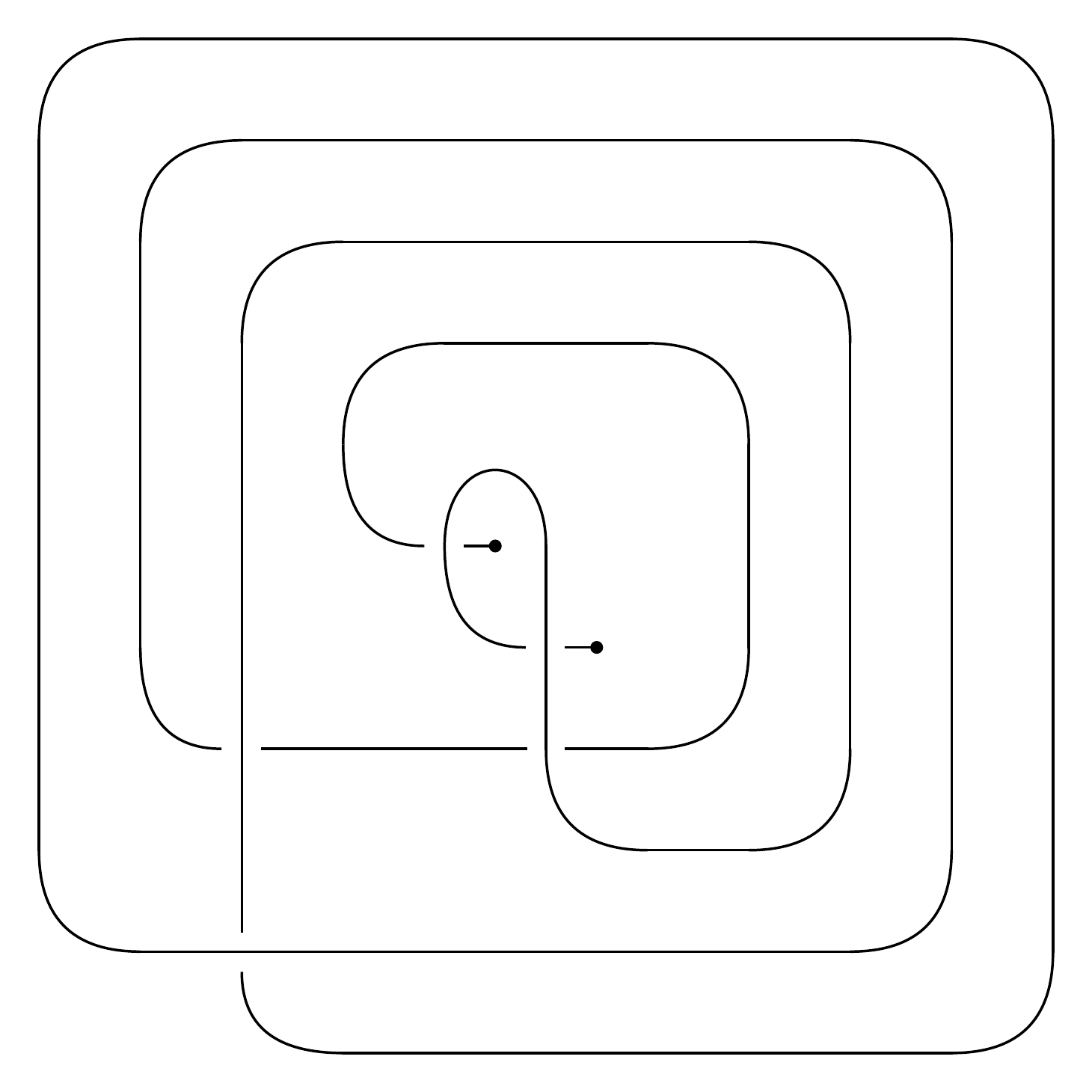}\\
\textcolor{black}{$5_{192}$}
\vspace{1cm}
\end{minipage}
\begin{minipage}[t]{.25\linewidth}
\centering
\includegraphics[width=0.9\textwidth,height=3.5cm,keepaspectratio]{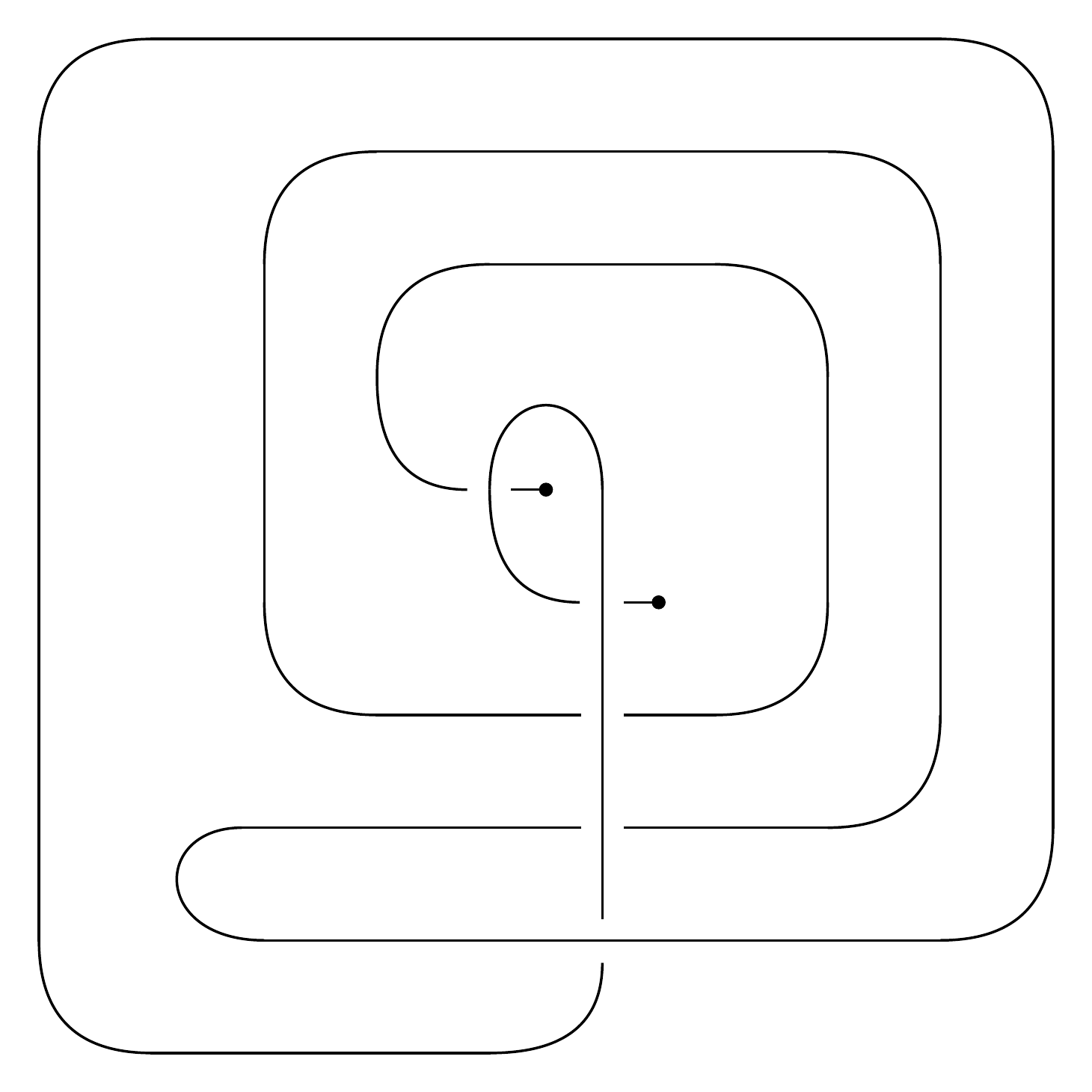}\\
\textcolor{black}{$5_{193}$}
\vspace{1cm}
\end{minipage}
\begin{minipage}[t]{.25\linewidth}
\centering
\includegraphics[width=0.9\textwidth,height=3.5cm,keepaspectratio]{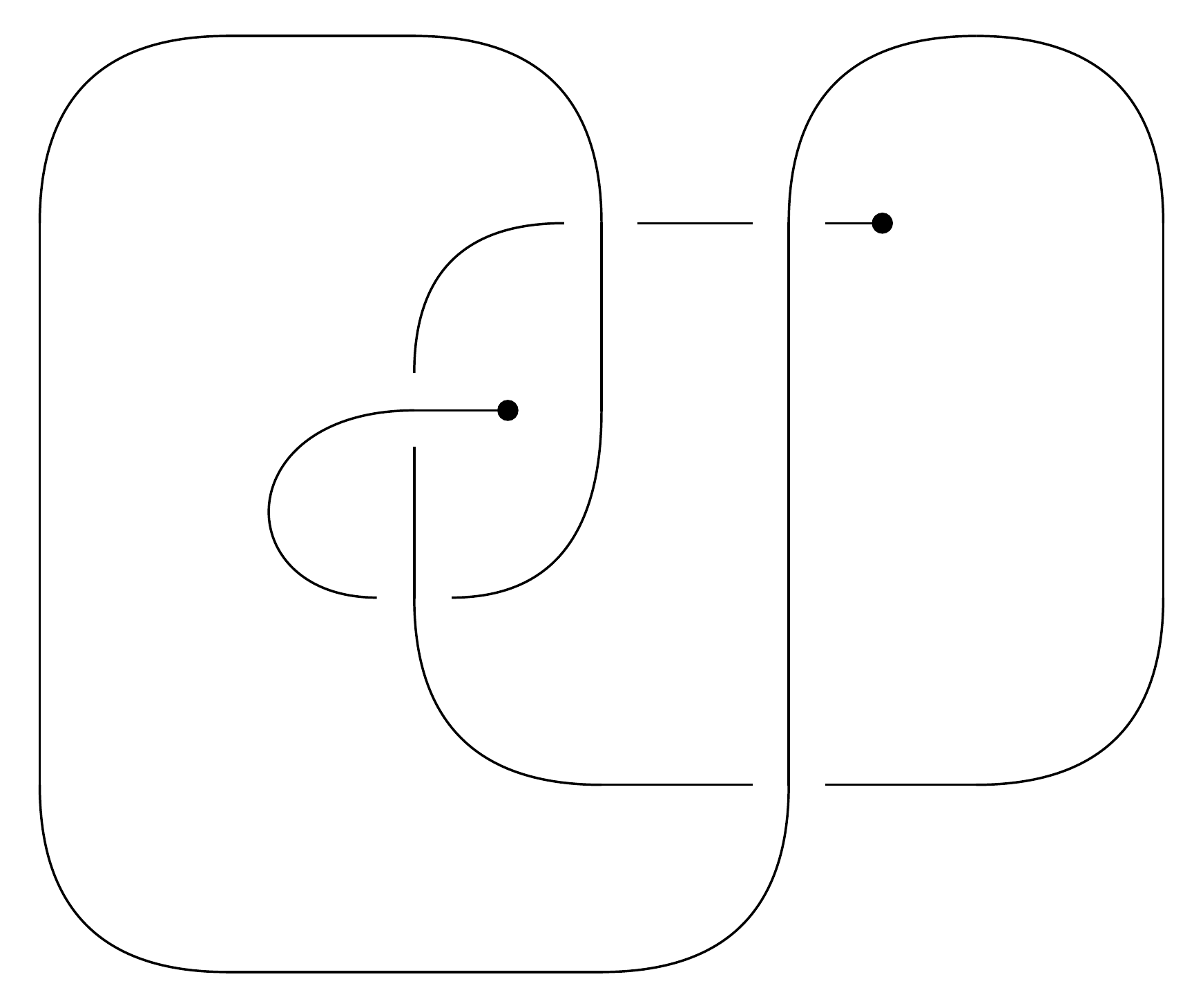}\\
\textcolor{black}{$5_{194}$}
\vspace{1cm}
\end{minipage}
\begin{minipage}[t]{.25\linewidth}
\centering
\includegraphics[width=0.9\textwidth,height=3.5cm,keepaspectratio]{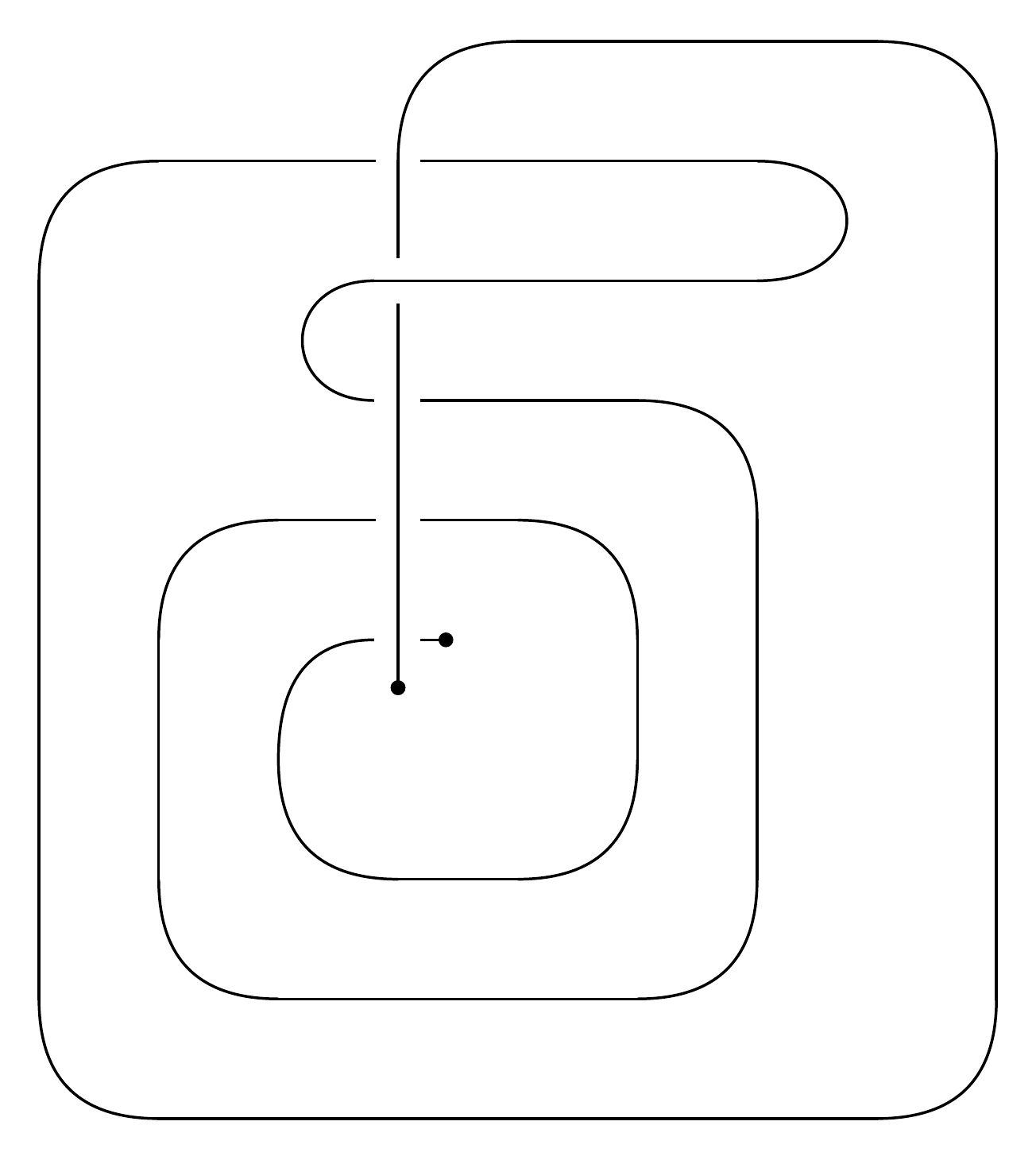}\\
\textcolor{black}{$5_{195}$}
\vspace{1cm}
\end{minipage}
\begin{minipage}[t]{.25\linewidth}
\centering
\includegraphics[width=0.9\textwidth,height=3.5cm,keepaspectratio]{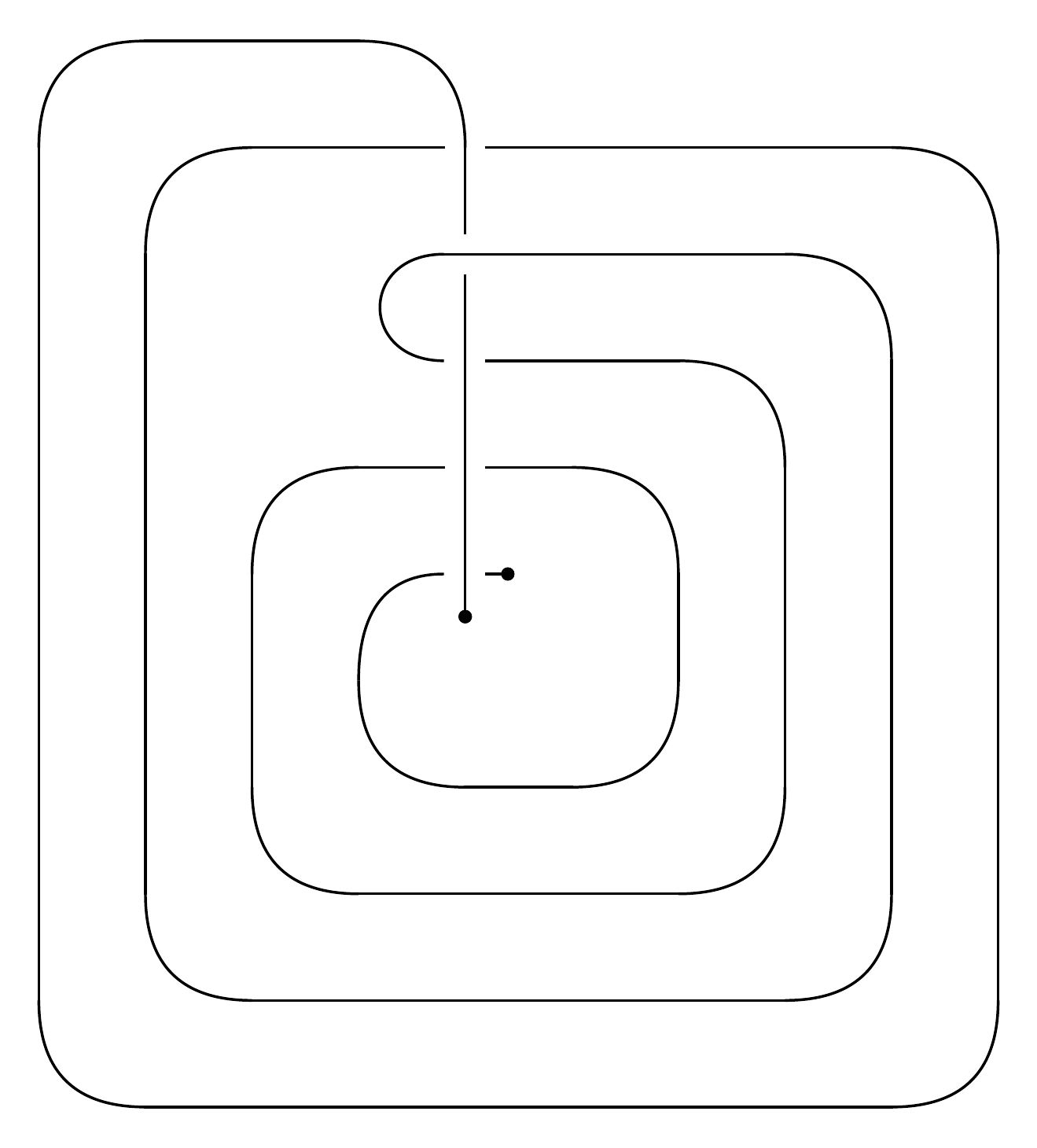}\\
\textcolor{black}{$5_{196}$}
\vspace{1cm}
\end{minipage}
\begin{minipage}[t]{.25\linewidth}
\centering
\includegraphics[width=0.9\textwidth,height=3.5cm,keepaspectratio]{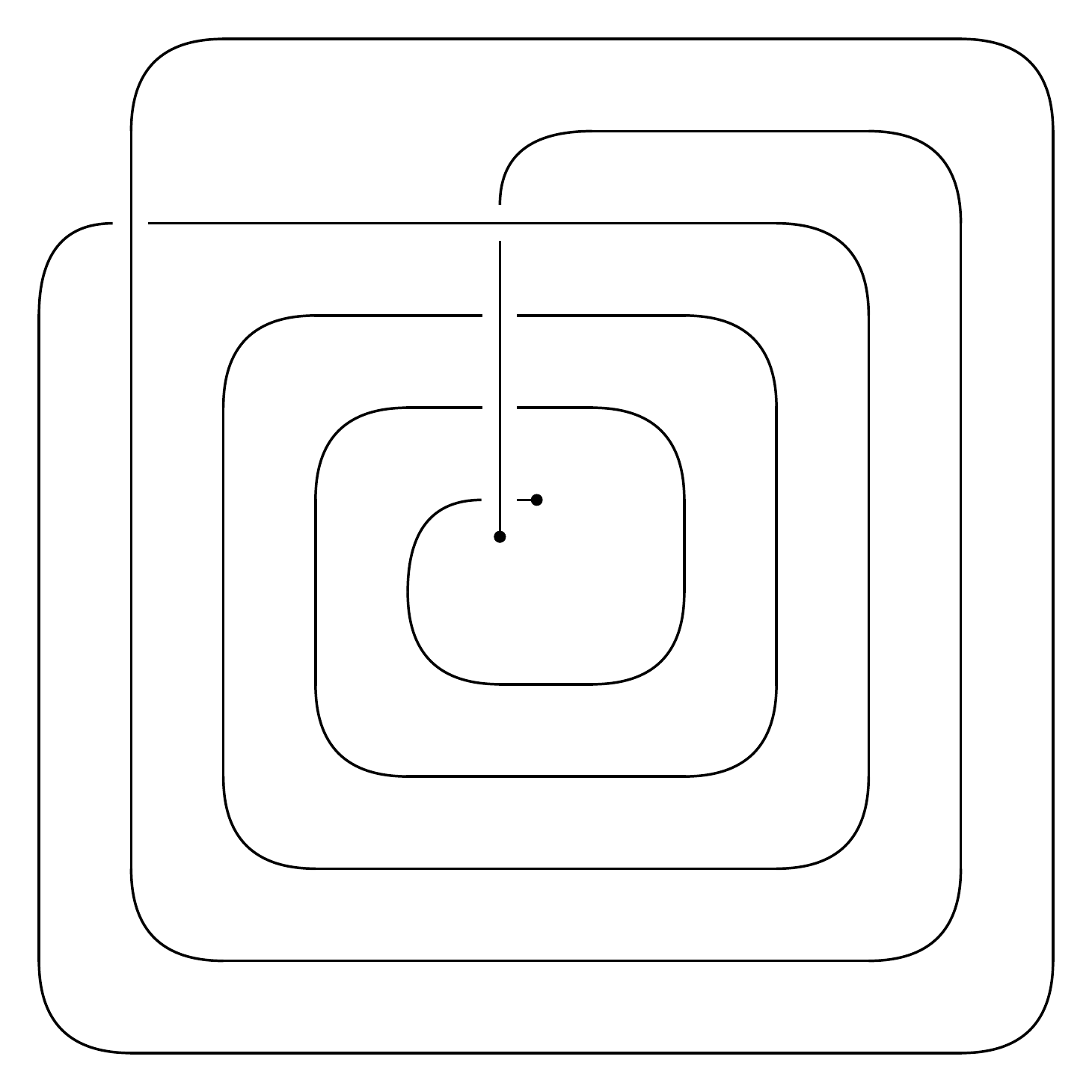}\\
\textcolor{black}{$5_{197}$}
\vspace{1cm}
\end{minipage}
\begin{minipage}[t]{.25\linewidth}
\centering
\includegraphics[width=0.9\textwidth,height=3.5cm,keepaspectratio]{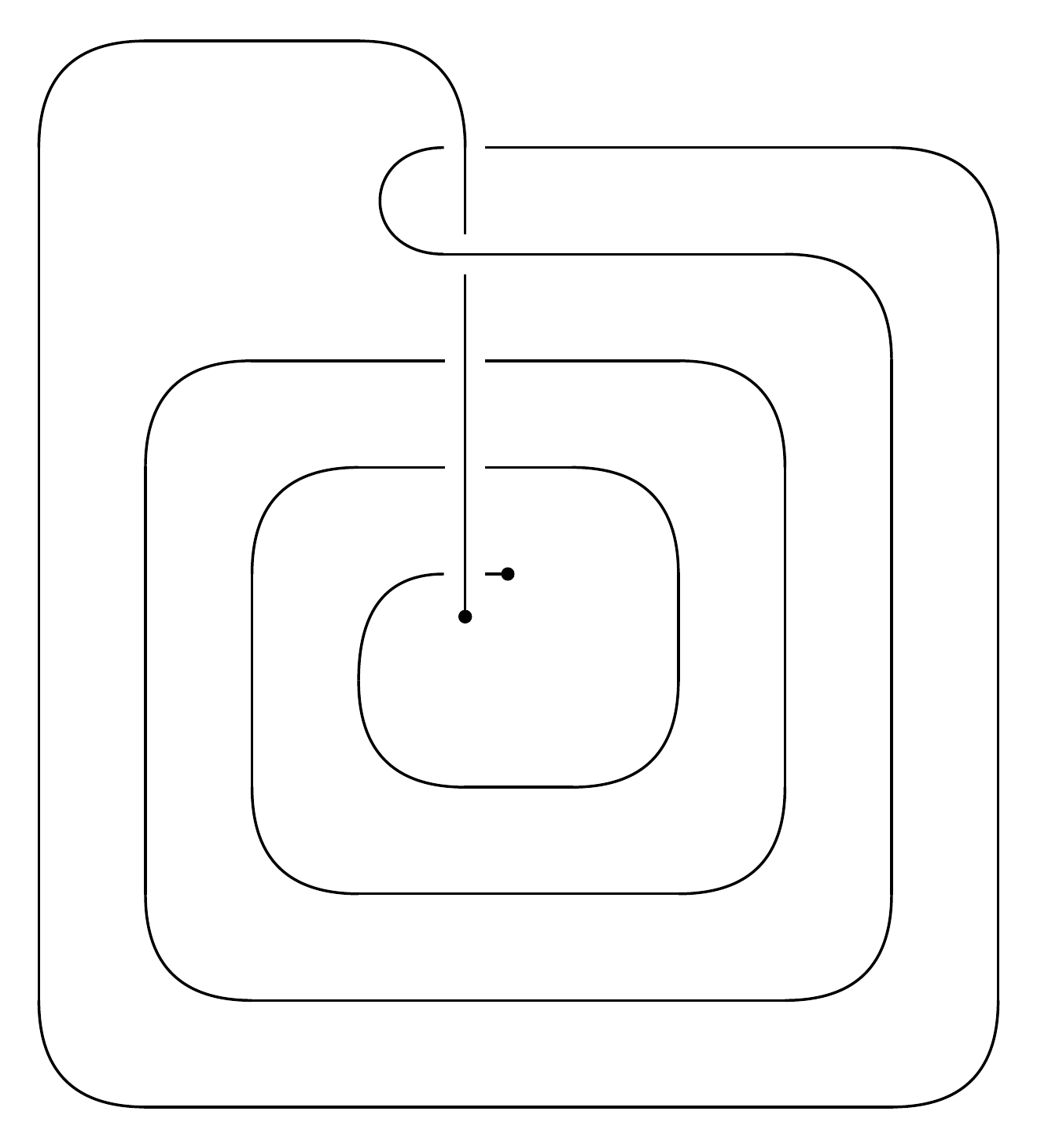}\\
\textcolor{black}{$5_{198}$}
\vspace{1cm}
\end{minipage}
\begin{minipage}[t]{.25\linewidth}
\centering
\includegraphics[width=0.9\textwidth,height=3.5cm,keepaspectratio]{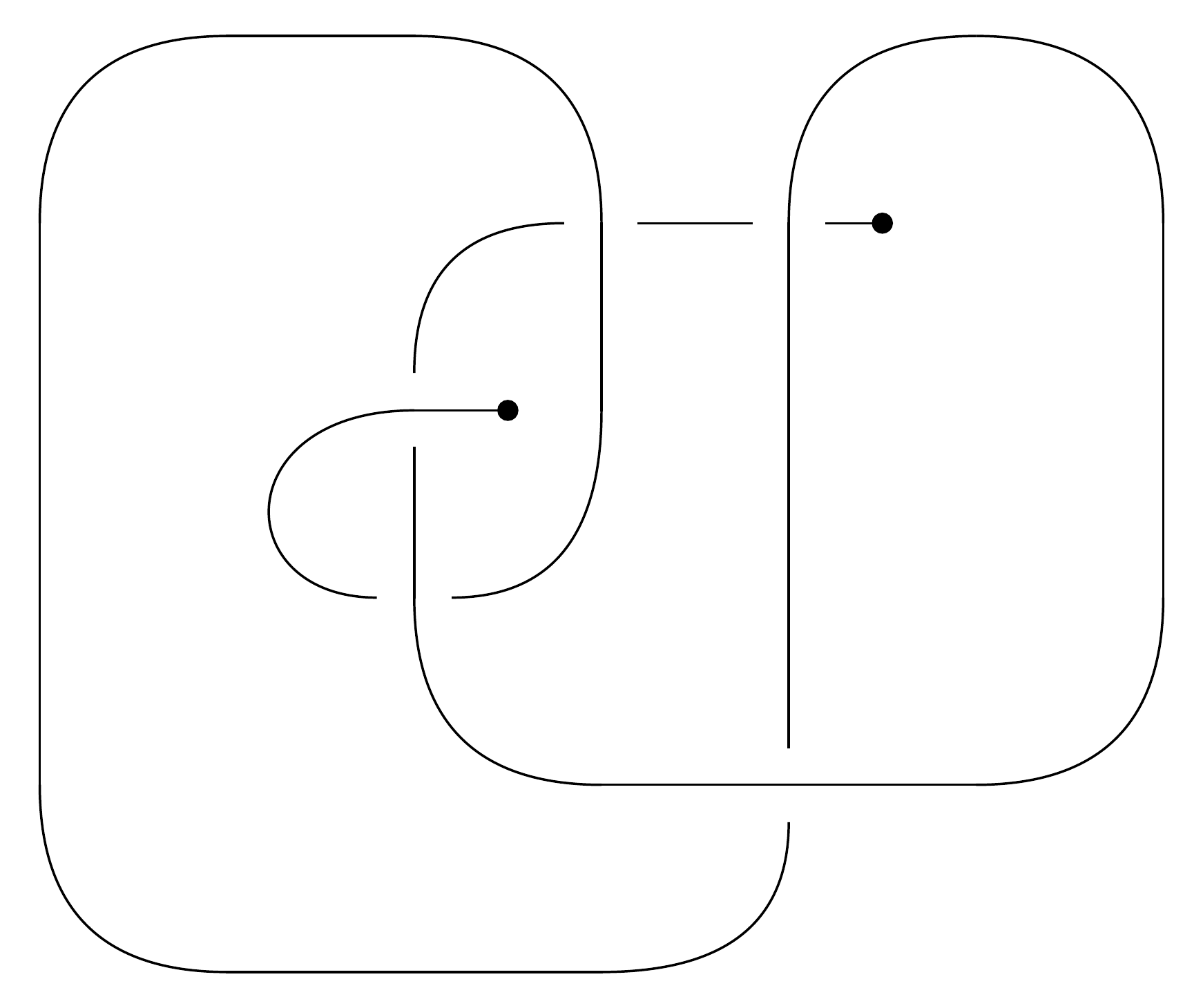}\\
\textcolor{black}{$5_{199}$}
\vspace{1cm}
\end{minipage}
\begin{minipage}[t]{.25\linewidth}
\centering
\includegraphics[width=0.9\textwidth,height=3.5cm,keepaspectratio]{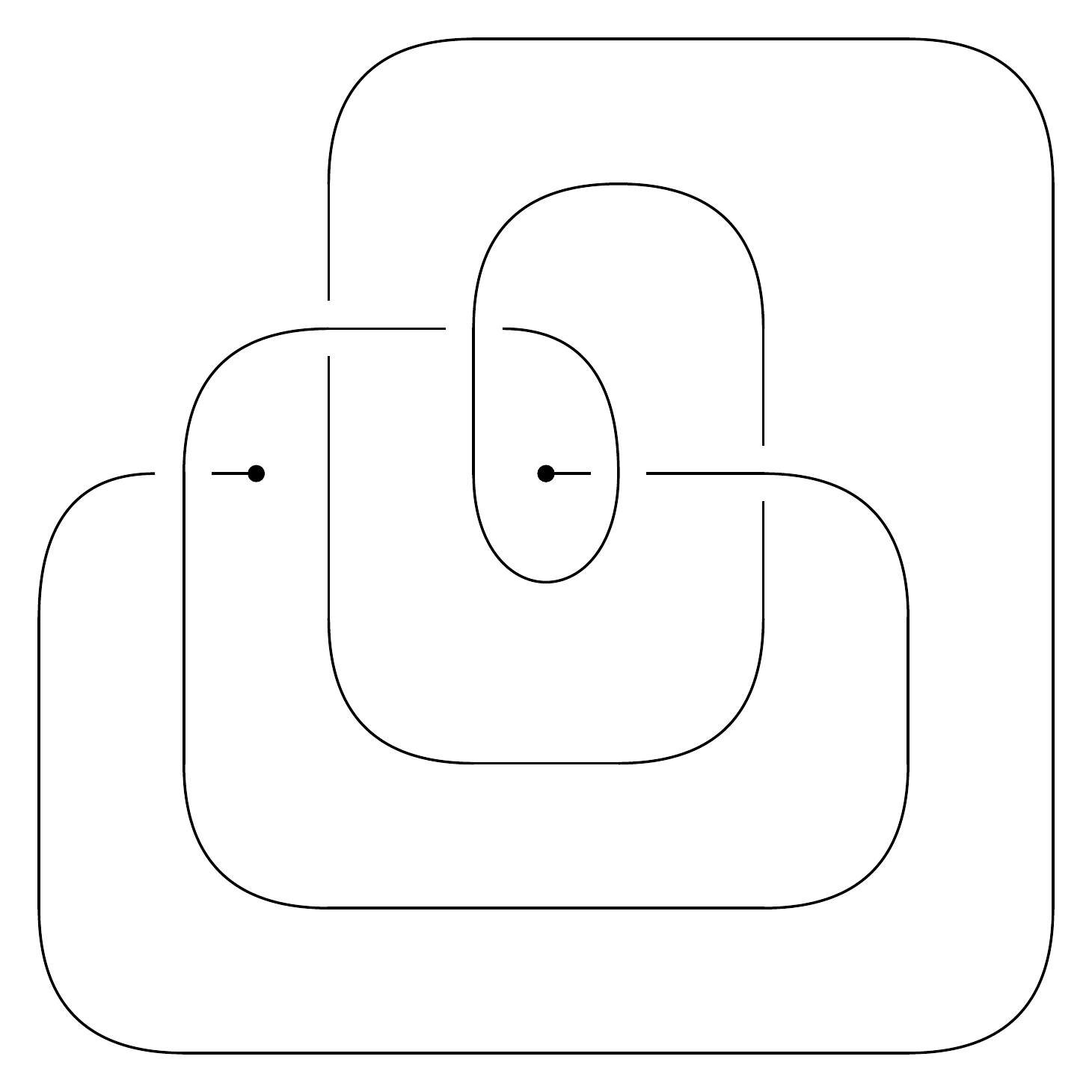}\\
\textcolor{black}{$5_{200}$}
\vspace{1cm}
\end{minipage}
\begin{minipage}[t]{.25\linewidth}
\centering
\includegraphics[width=0.9\textwidth,height=3.5cm,keepaspectratio]{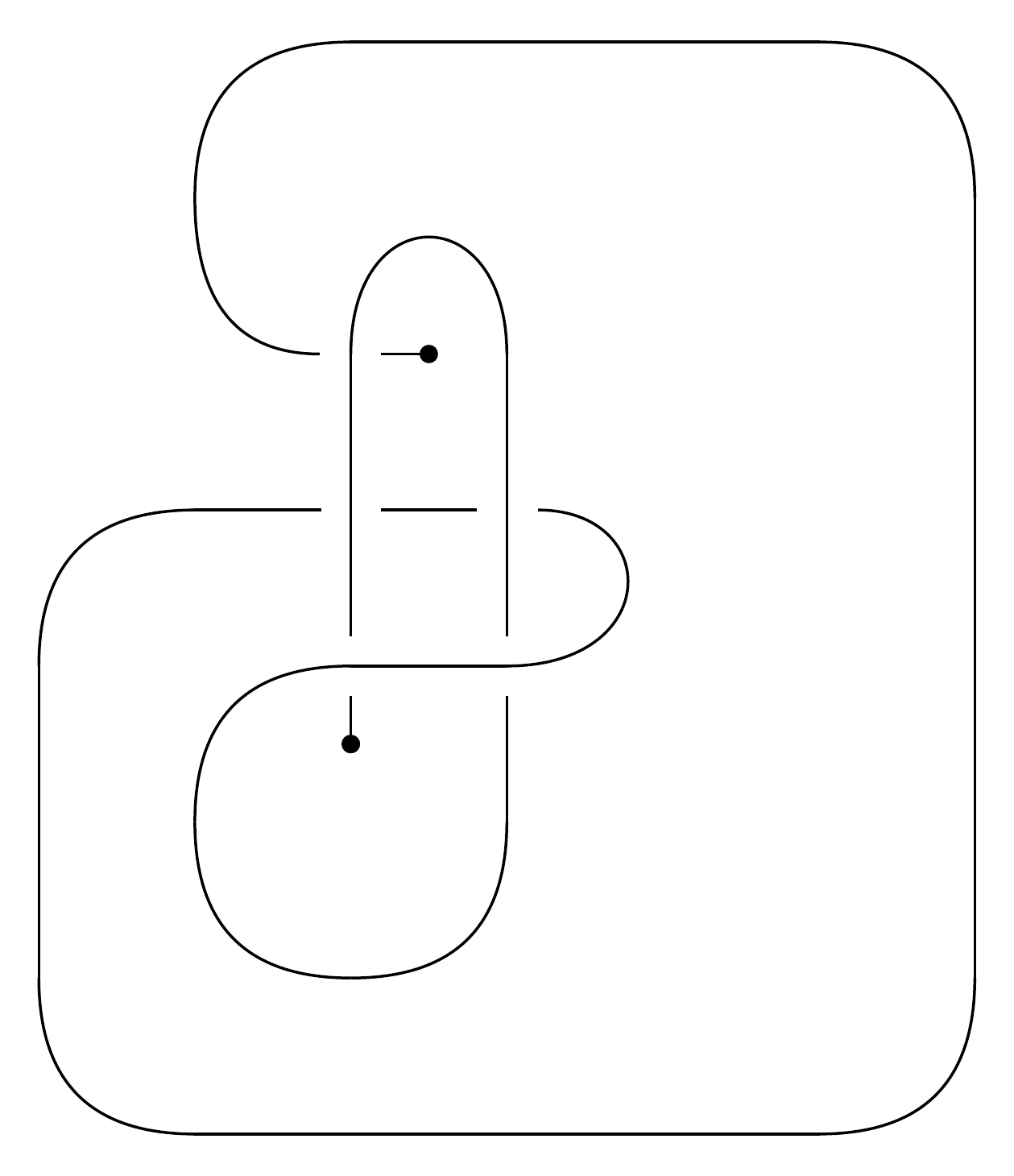}\\
\textcolor{black}{$5_{201}$}
\vspace{1cm}
\end{minipage}
\begin{minipage}[t]{.25\linewidth}
\centering
\includegraphics[width=0.9\textwidth,height=3.5cm,keepaspectratio]{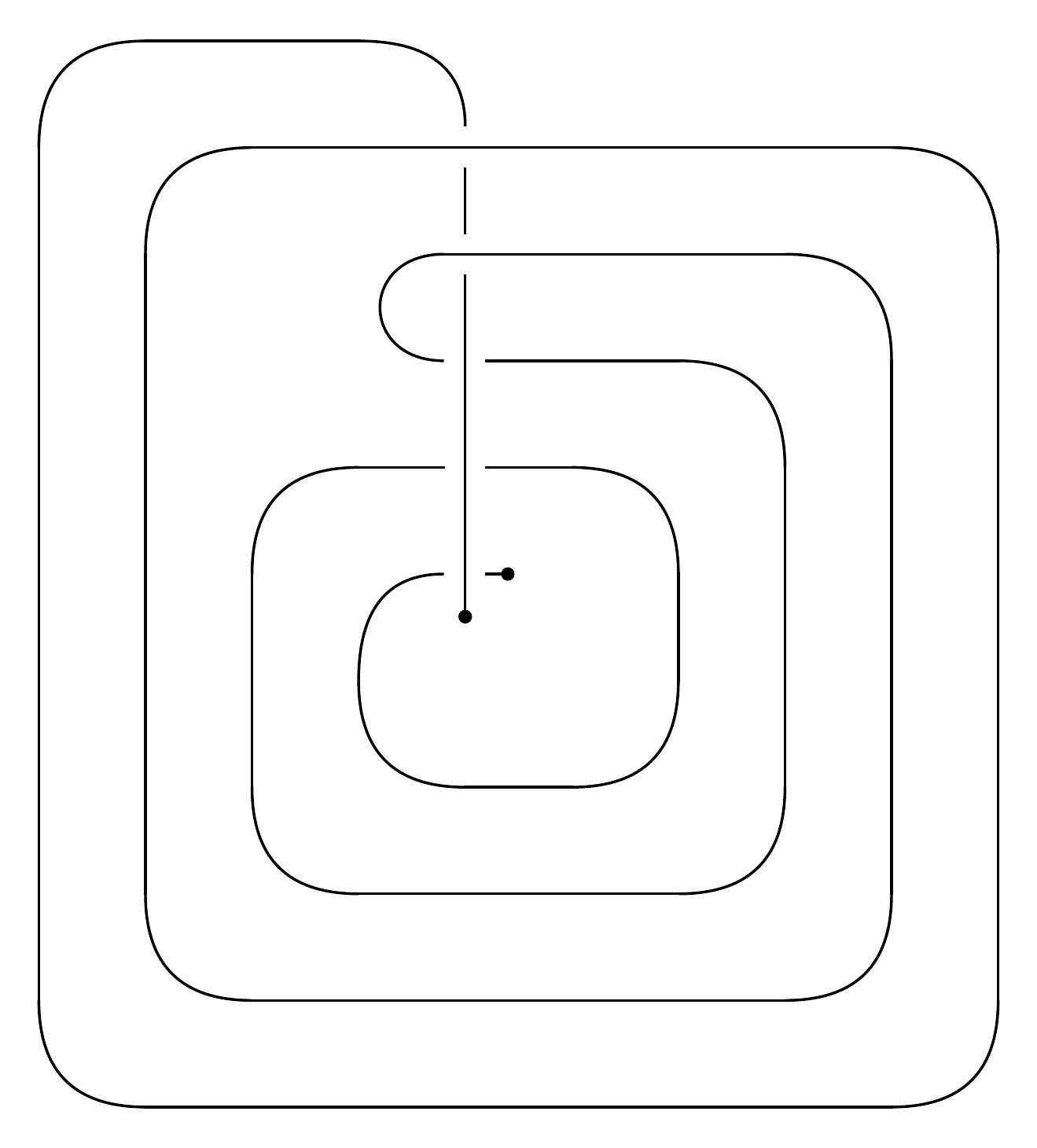}\\
\textcolor{black}{$5_{202}$}
\vspace{1cm}
\end{minipage}
\begin{minipage}[t]{.25\linewidth}
\centering
\includegraphics[width=0.9\textwidth,height=3.5cm,keepaspectratio]{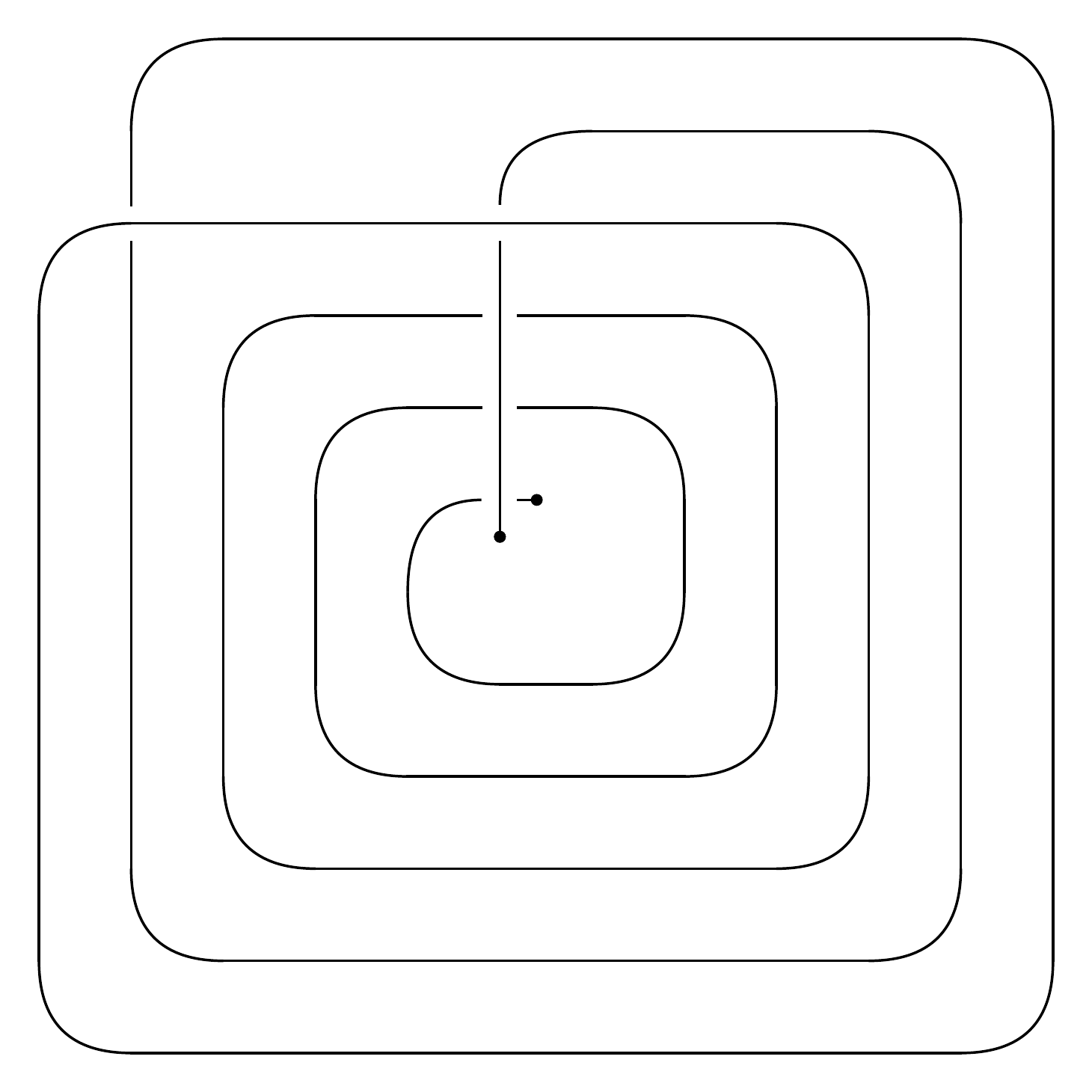}\\
\textcolor{black}{$5_{203}$}
\vspace{1cm}
\end{minipage}
\begin{minipage}[t]{.25\linewidth}
\centering
\includegraphics[width=0.9\textwidth,height=3.5cm,keepaspectratio]{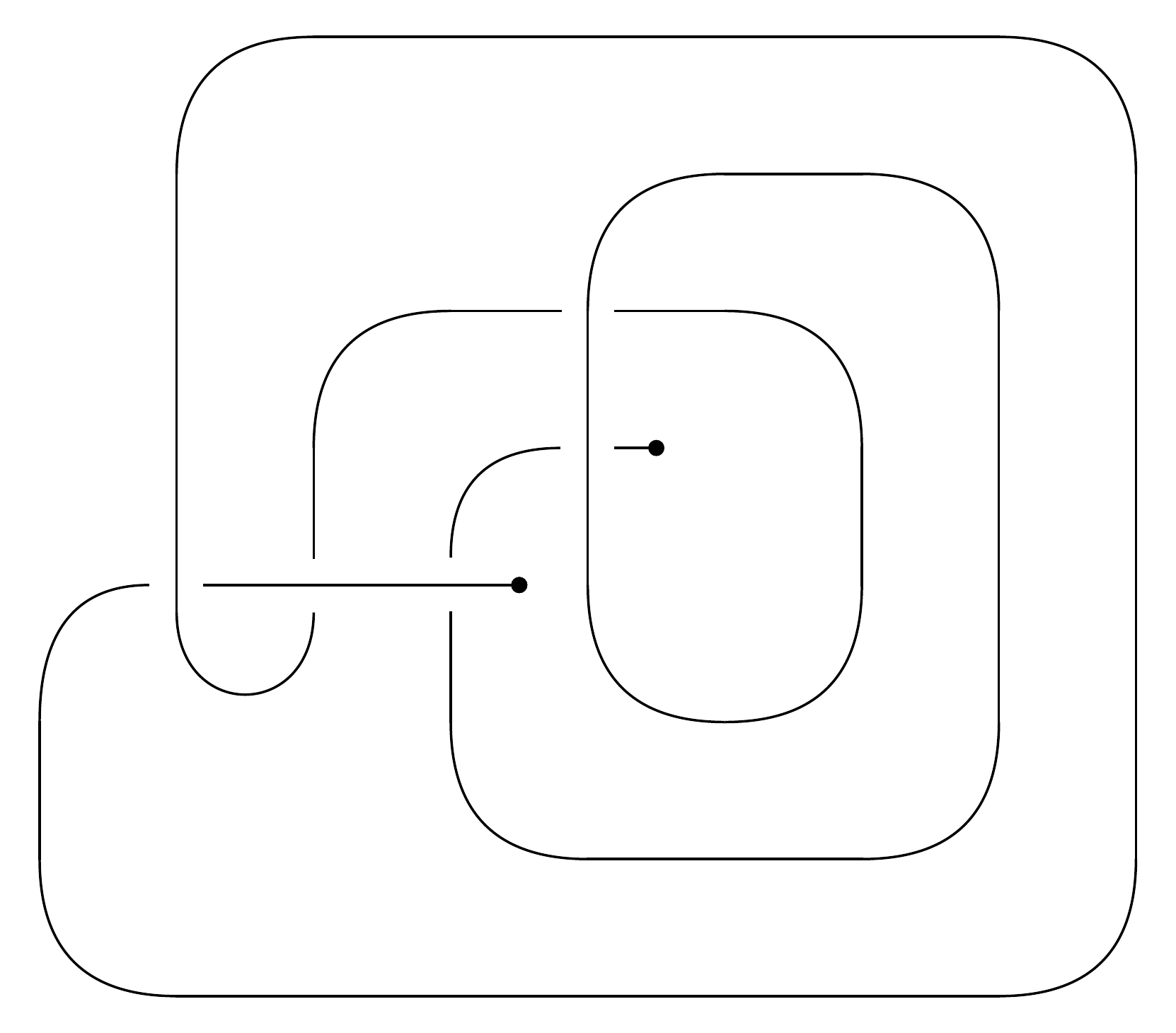}\\
\textcolor{black}{$5_{204}$}
\vspace{1cm}
\end{minipage}
\begin{minipage}[t]{.25\linewidth}
\centering
\includegraphics[width=0.9\textwidth,height=3.5cm,keepaspectratio]{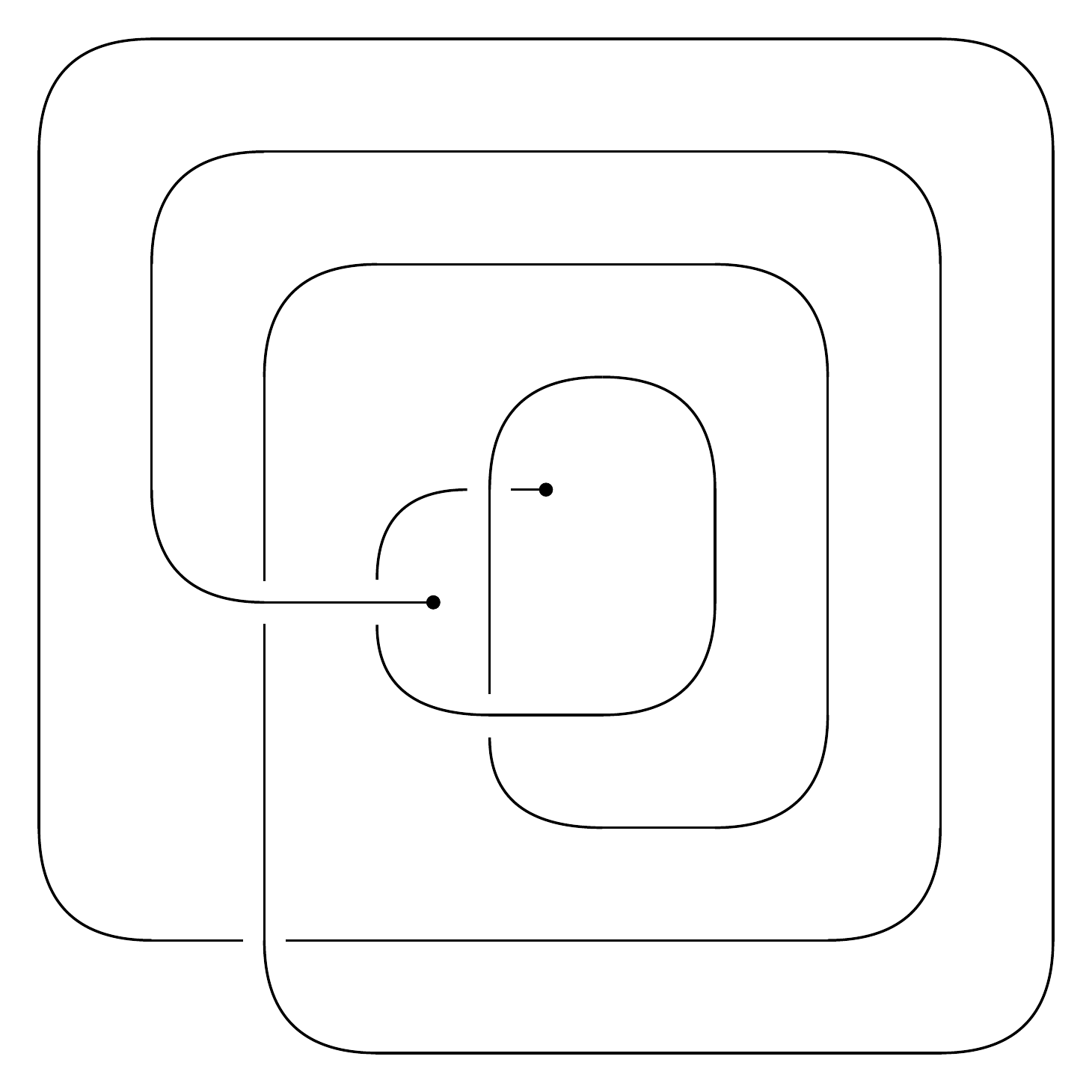}\\
\textcolor{black}{$5_{205}$}
\vspace{1cm}
\end{minipage}
\begin{minipage}[t]{.25\linewidth}
\centering
\includegraphics[width=0.9\textwidth,height=3.5cm,keepaspectratio]{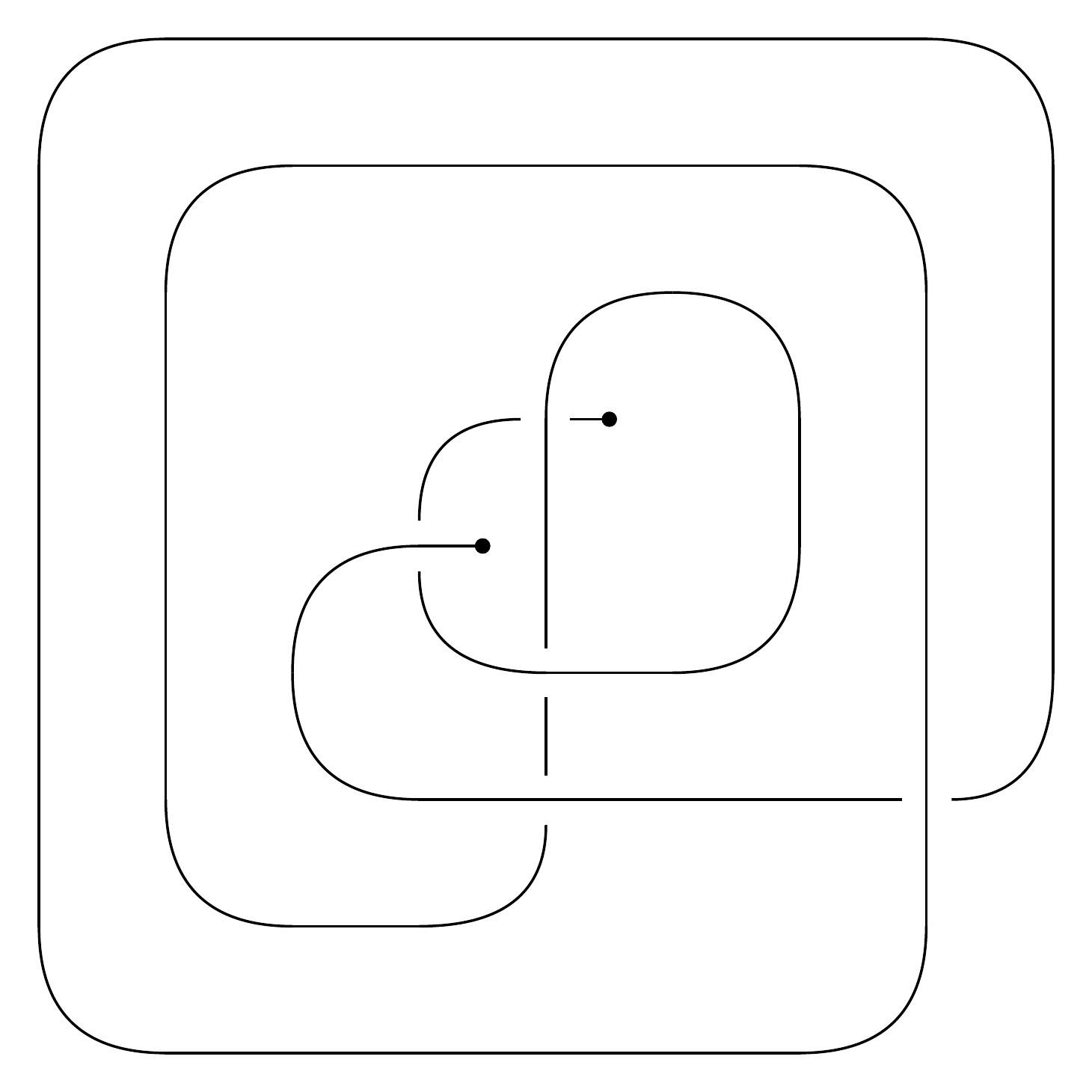}\\
\textcolor{black}{$5_{206}$}
\vspace{1cm}
\end{minipage}
\begin{minipage}[t]{.25\linewidth}
\centering
\includegraphics[width=0.9\textwidth,height=3.5cm,keepaspectratio]{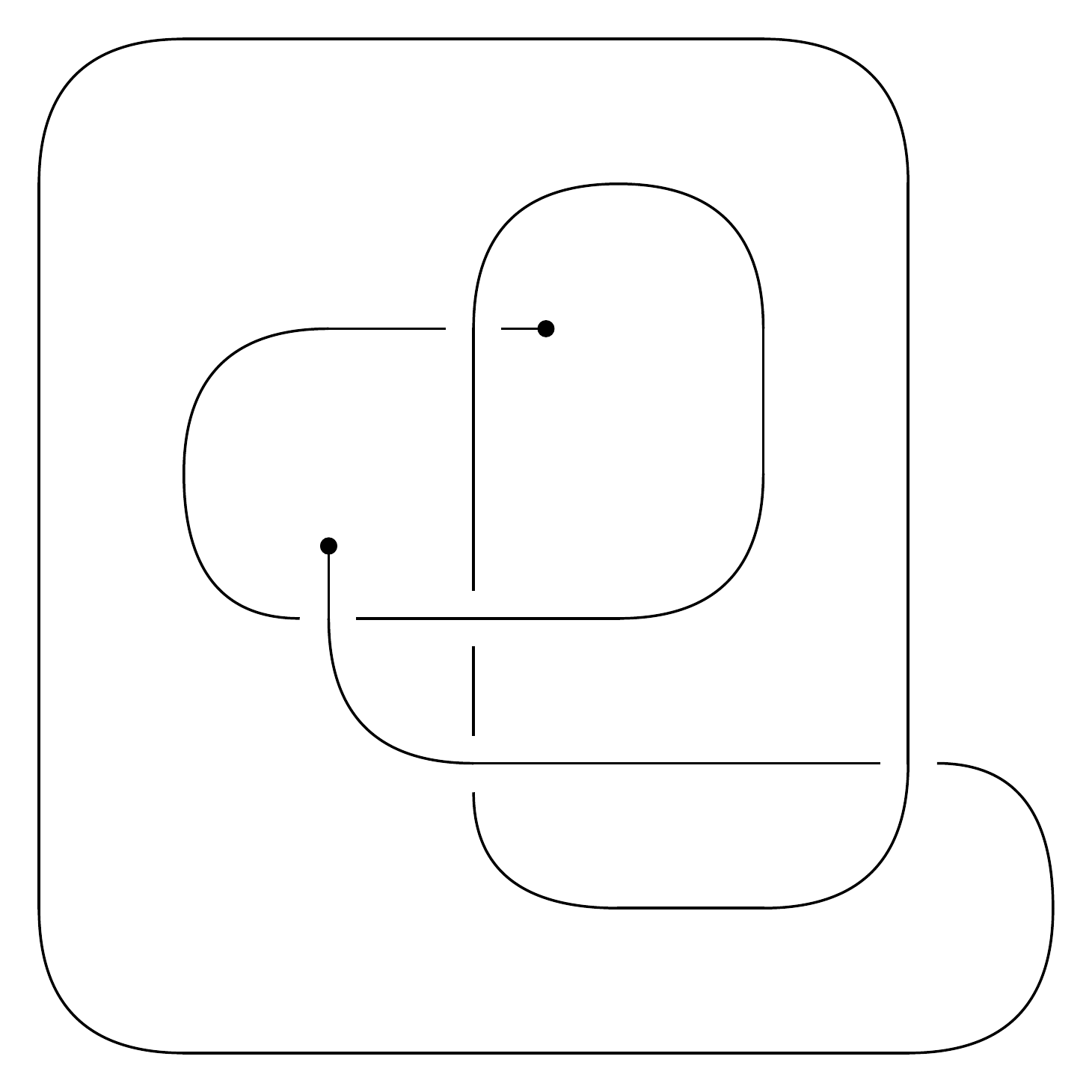}\\
\textcolor{black}{$5_{207}$}
\vspace{1cm}
\end{minipage}
\begin{minipage}[t]{.25\linewidth}
\centering
\includegraphics[width=0.9\textwidth,height=3.5cm,keepaspectratio]{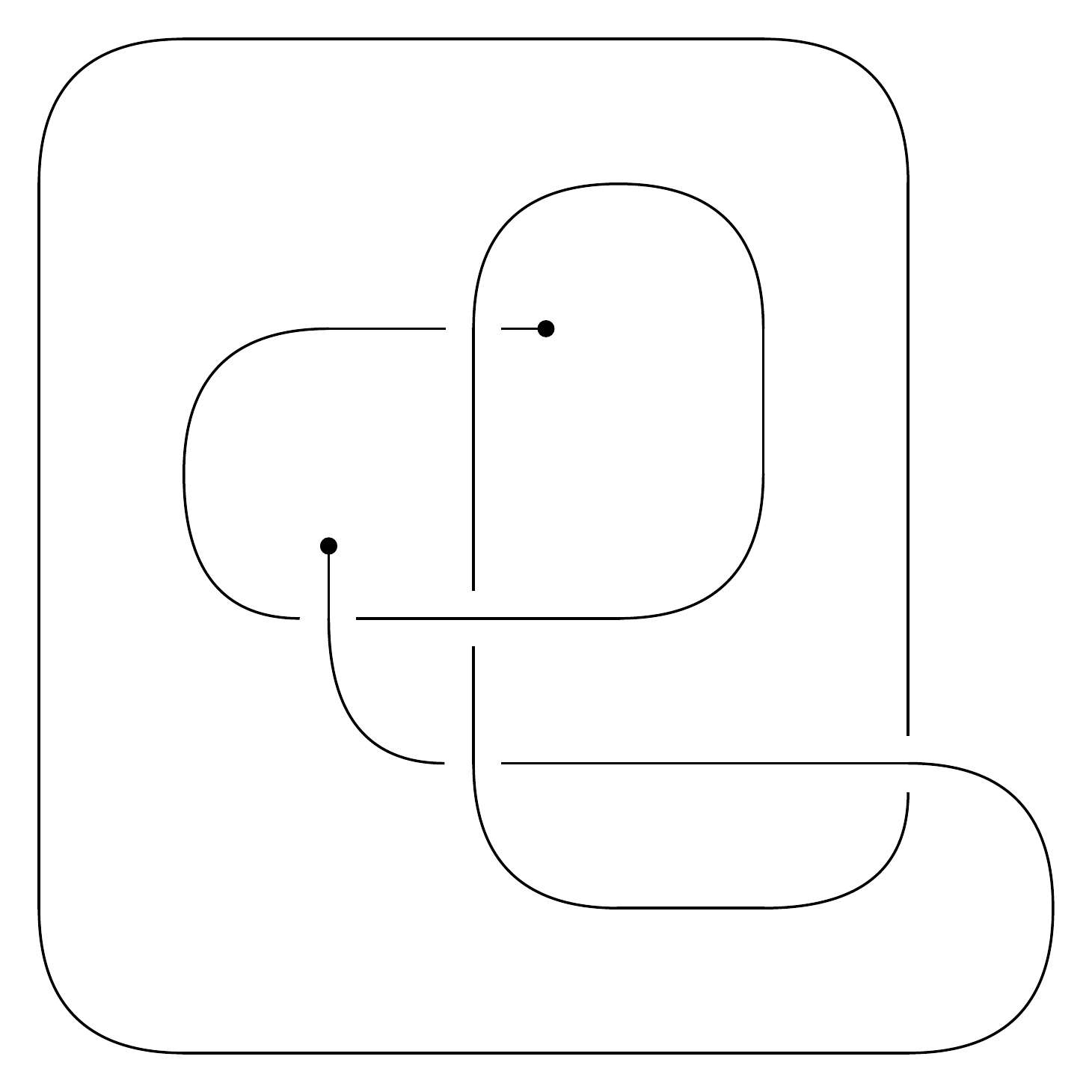}\\
\textcolor{black}{$5_{208}$}
\vspace{1cm}
\end{minipage}
\begin{minipage}[t]{.25\linewidth}
\centering
\includegraphics[width=0.9\textwidth,height=3.5cm,keepaspectratio]{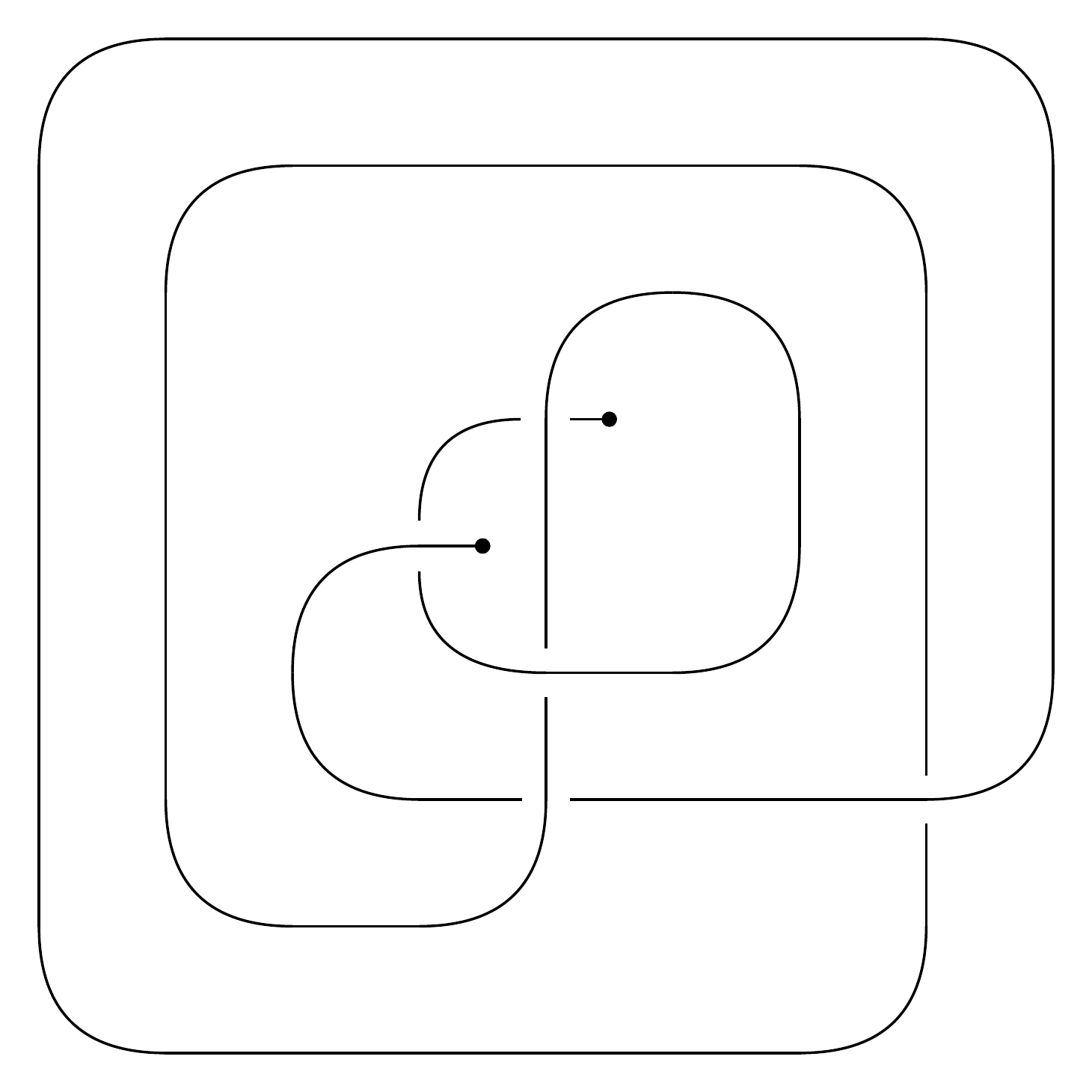}\\
\textcolor{black}{$5_{209}$}
\vspace{1cm}
\end{minipage}
\begin{minipage}[t]{.25\linewidth}
\centering
\includegraphics[width=0.9\textwidth,height=3.5cm,keepaspectratio]{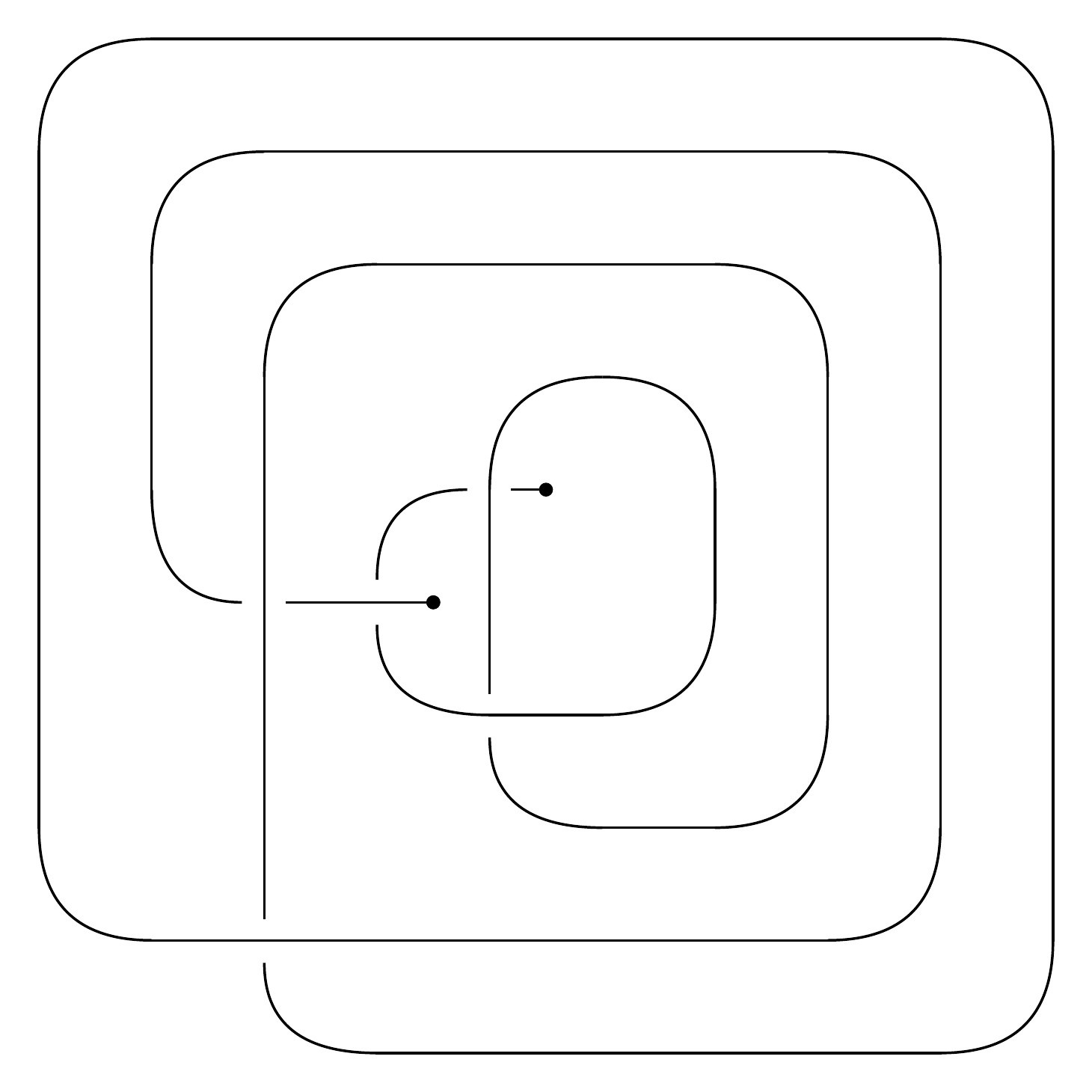}\\
\textcolor{black}{$5_{210}$}
\vspace{1cm}
\end{minipage}
\begin{minipage}[t]{.25\linewidth}
\centering
\includegraphics[width=0.9\textwidth,height=3.5cm,keepaspectratio]{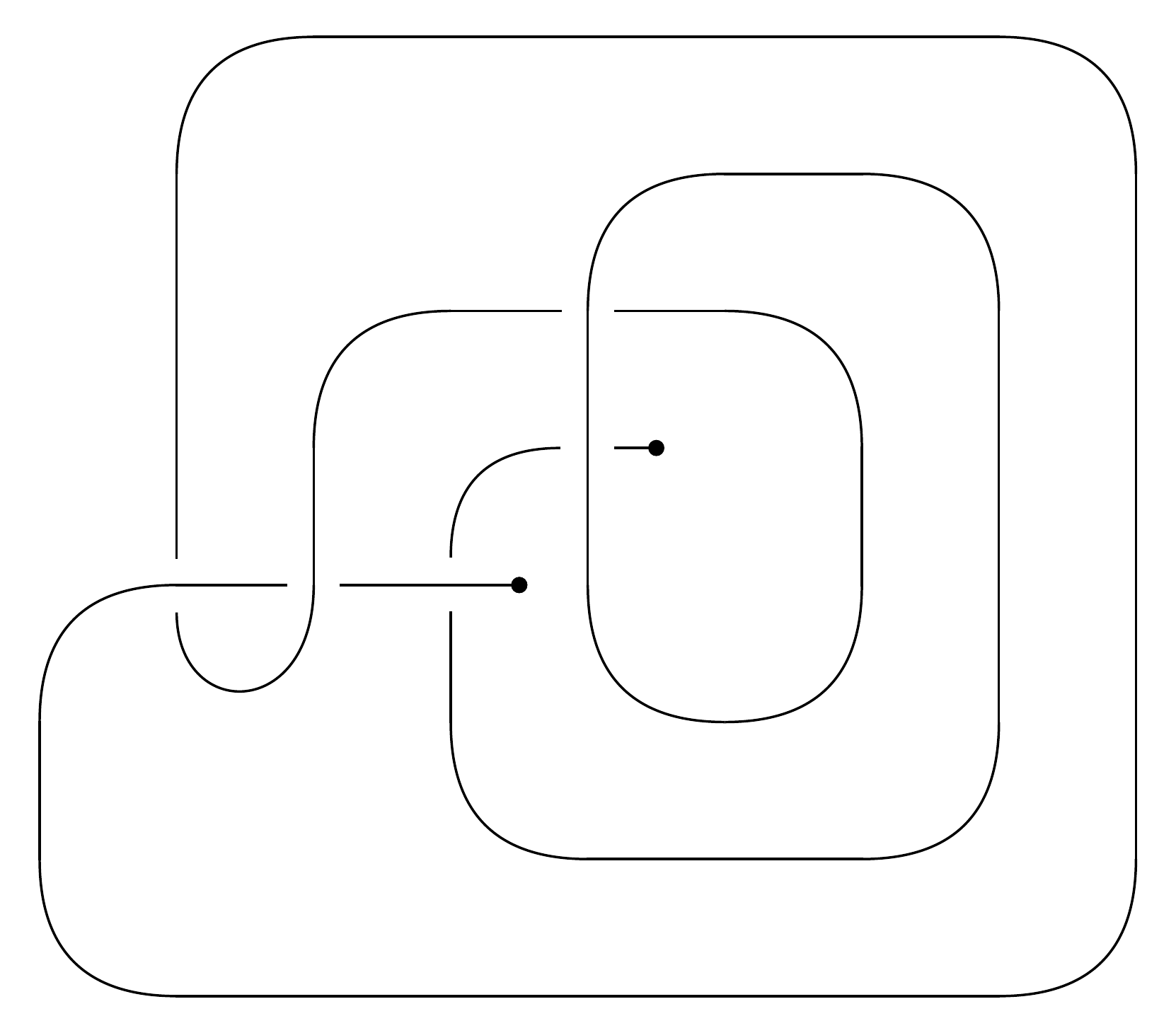}\\
\textcolor{black}{$5_{211}$}
\vspace{1cm}
\end{minipage}
\begin{minipage}[t]{.25\linewidth}
\centering
\includegraphics[width=0.9\textwidth,height=3.5cm,keepaspectratio]{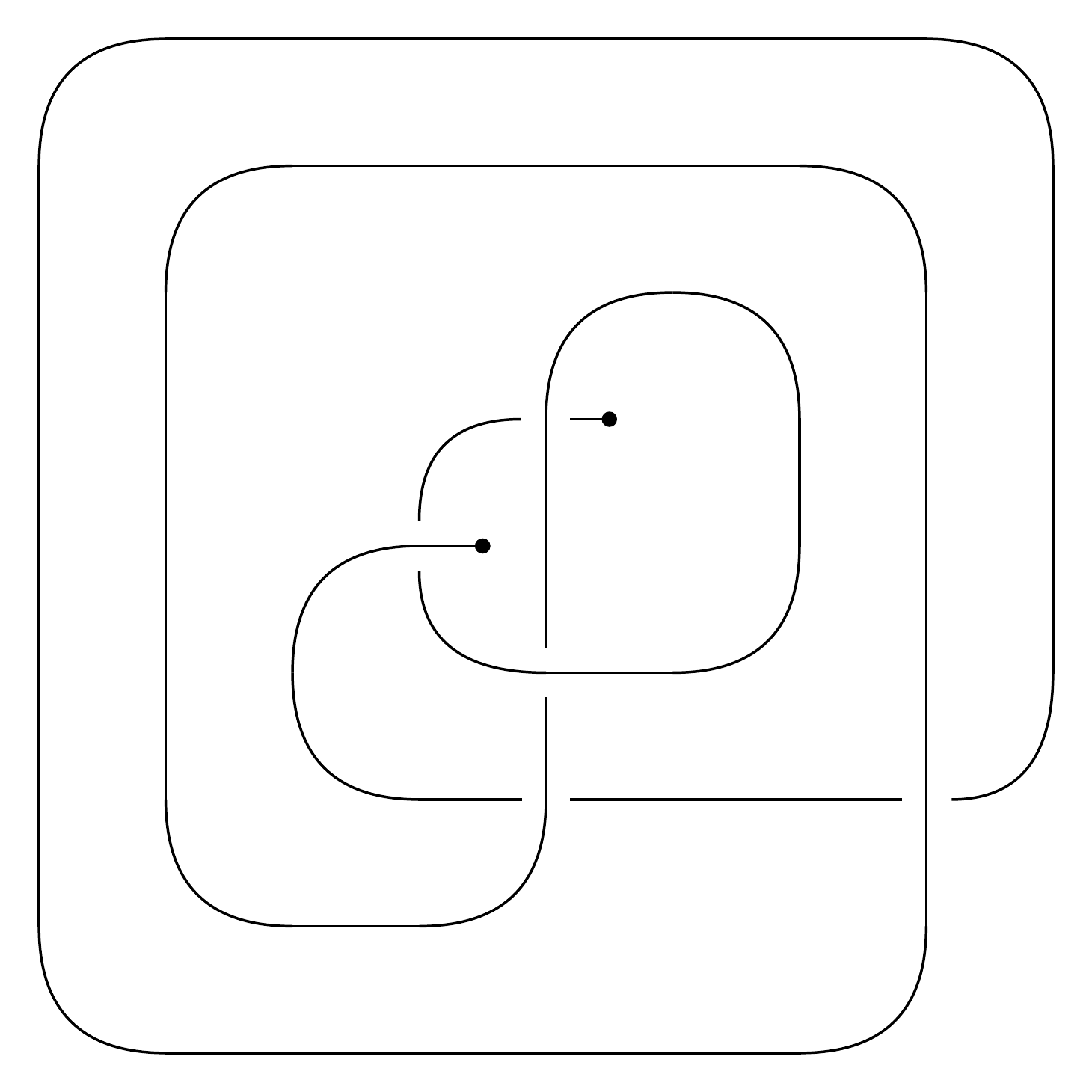}\\
\textcolor{black}{$5_{212}$}
\vspace{1cm}
\end{minipage}
\begin{minipage}[t]{.25\linewidth}
\centering
\includegraphics[width=0.9\textwidth,height=3.5cm,keepaspectratio]{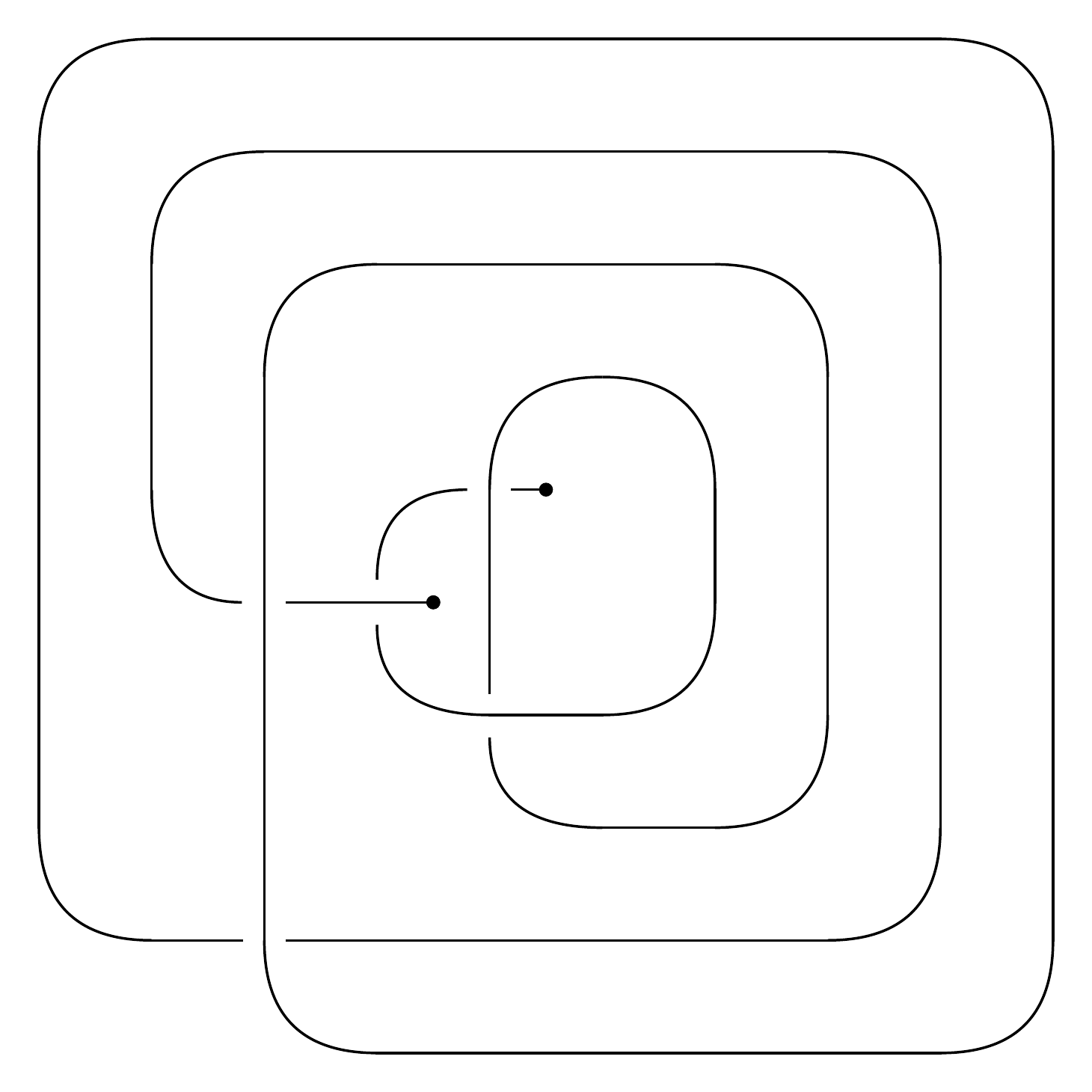}\\
\textcolor{black}{$5_{213}$}
\vspace{1cm}
\end{minipage}
\begin{minipage}[t]{.25\linewidth}
\centering
\includegraphics[width=0.9\textwidth,height=3.5cm,keepaspectratio]{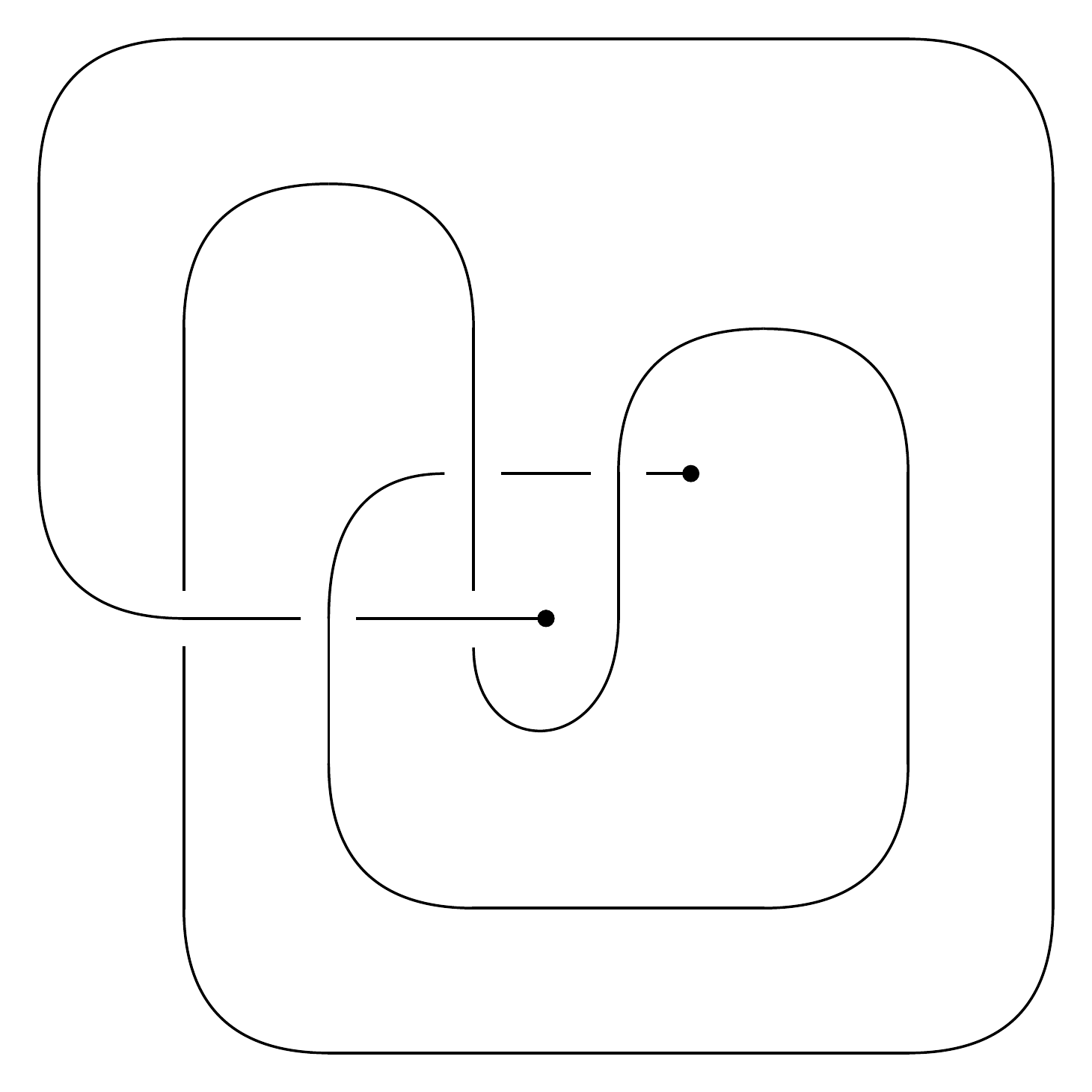}\\
\textcolor{black}{$5_{214}$}
\vspace{1cm}
\end{minipage}
\begin{minipage}[t]{.25\linewidth}
\centering
\includegraphics[width=0.9\textwidth,height=3.5cm,keepaspectratio]{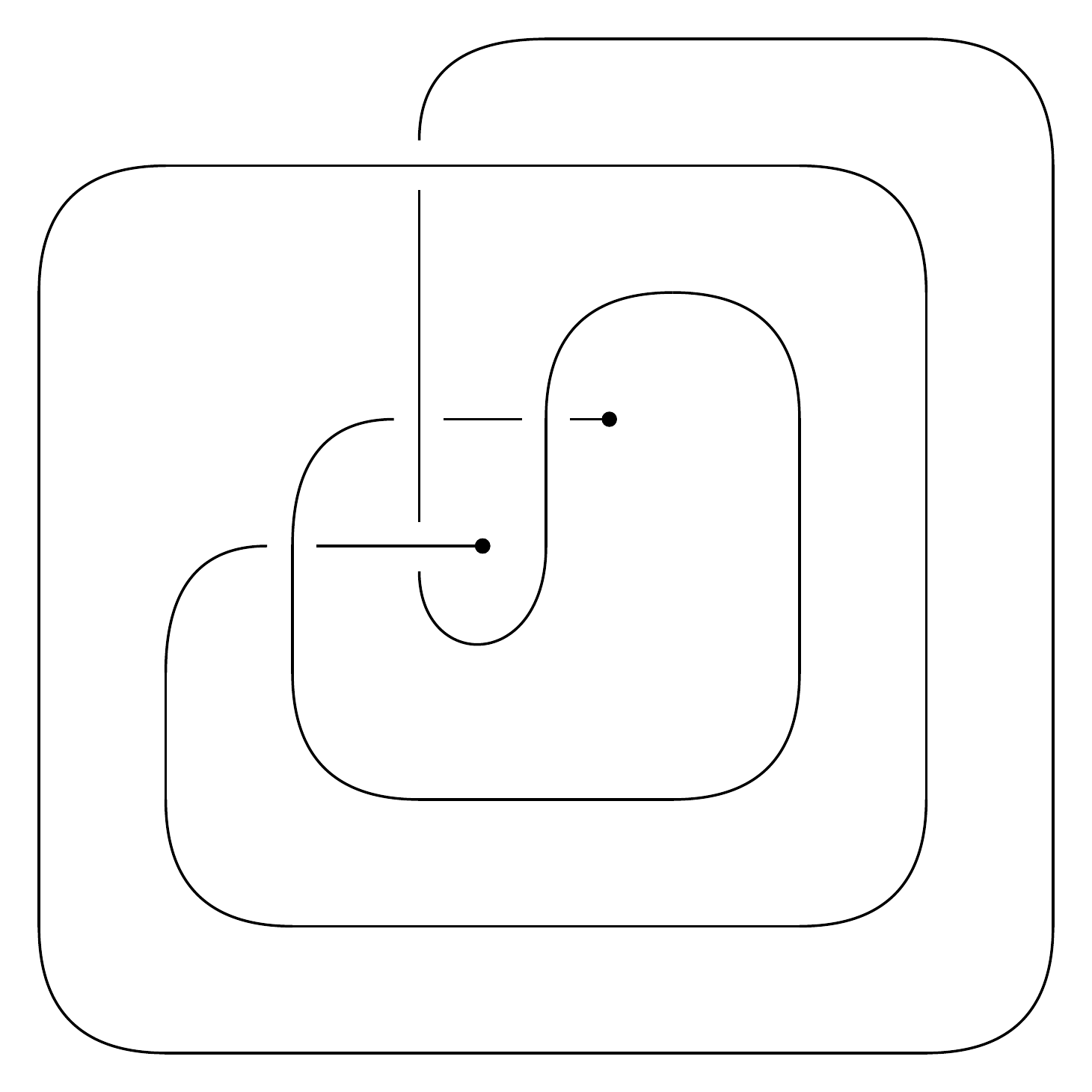}\\
\textcolor{black}{$5_{215}$}
\vspace{1cm}
\end{minipage}
\begin{minipage}[t]{.25\linewidth}
\centering
\includegraphics[width=0.9\textwidth,height=3.5cm,keepaspectratio]{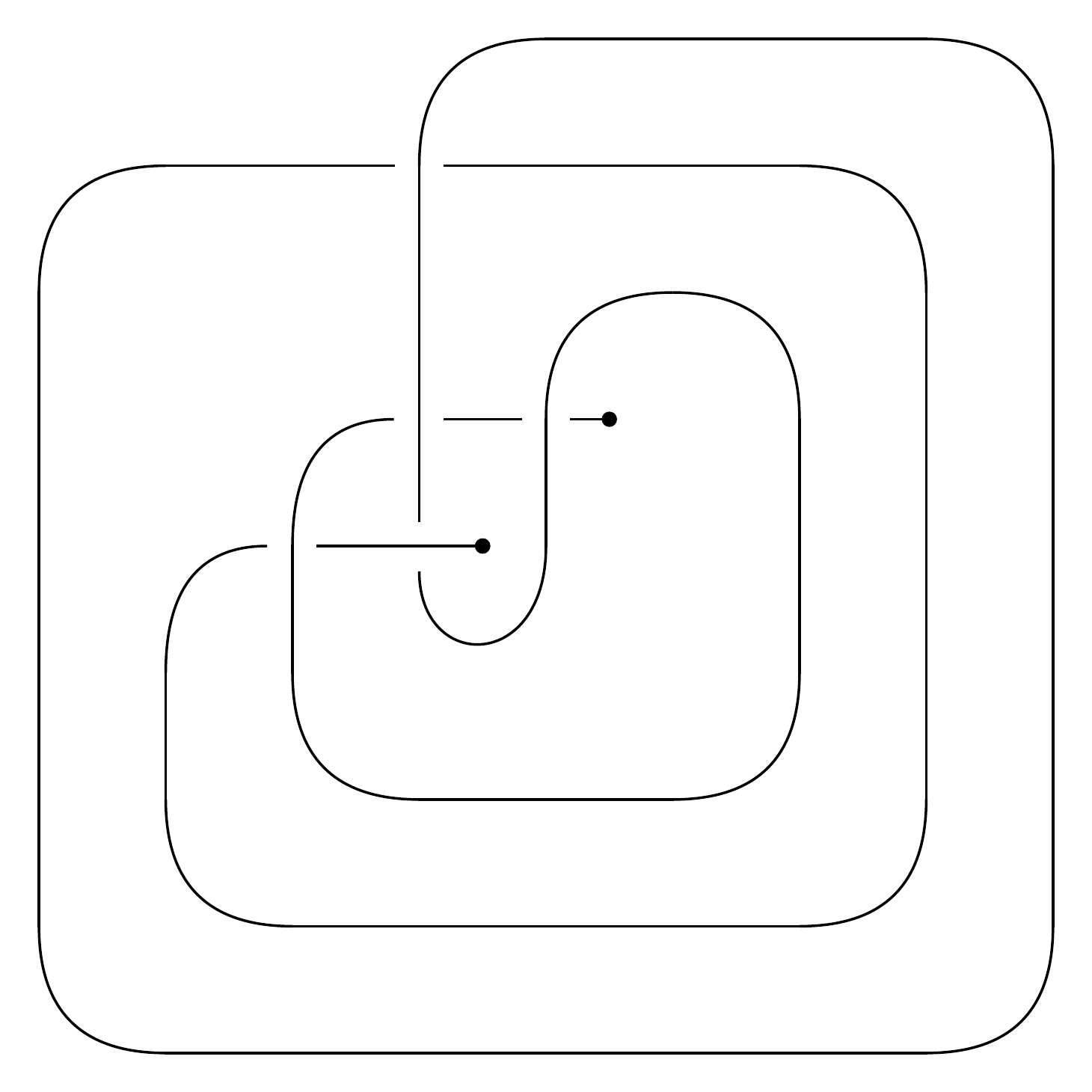}\\
\textcolor{black}{$5_{216}$}
\vspace{1cm}
\end{minipage}
\begin{minipage}[t]{.25\linewidth}
\centering
\includegraphics[width=0.9\textwidth,height=3.5cm,keepaspectratio]{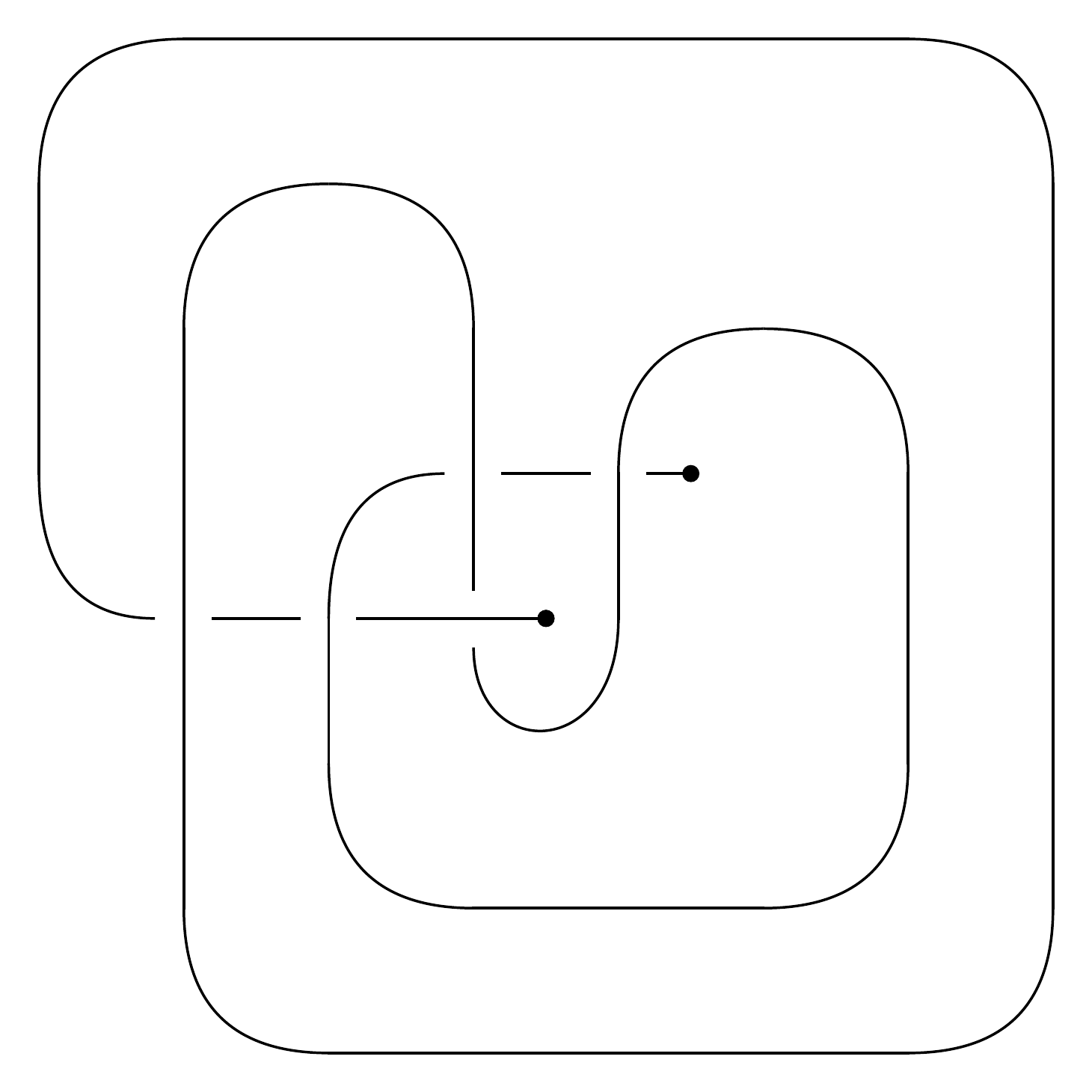}\\
\textcolor{black}{$5_{217}$}
\vspace{1cm}
\end{minipage}
\begin{minipage}[t]{.25\linewidth}
\centering
\includegraphics[width=0.9\textwidth,height=3.5cm,keepaspectratio]{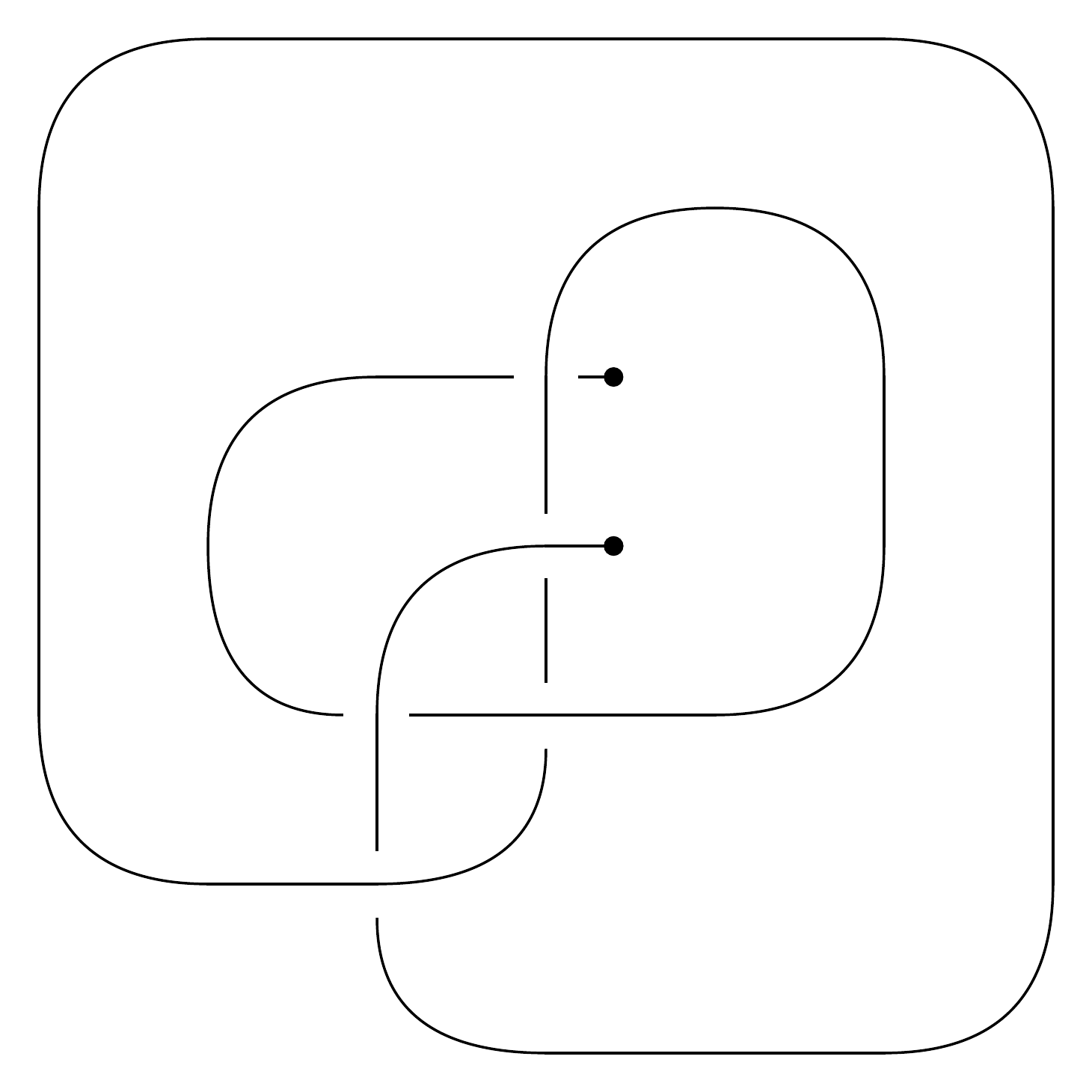}\\
\textcolor{black}{$5_{218}$}
\vspace{1cm}
\end{minipage}
\begin{minipage}[t]{.25\linewidth}
\centering
\includegraphics[width=0.9\textwidth,height=3.5cm,keepaspectratio]{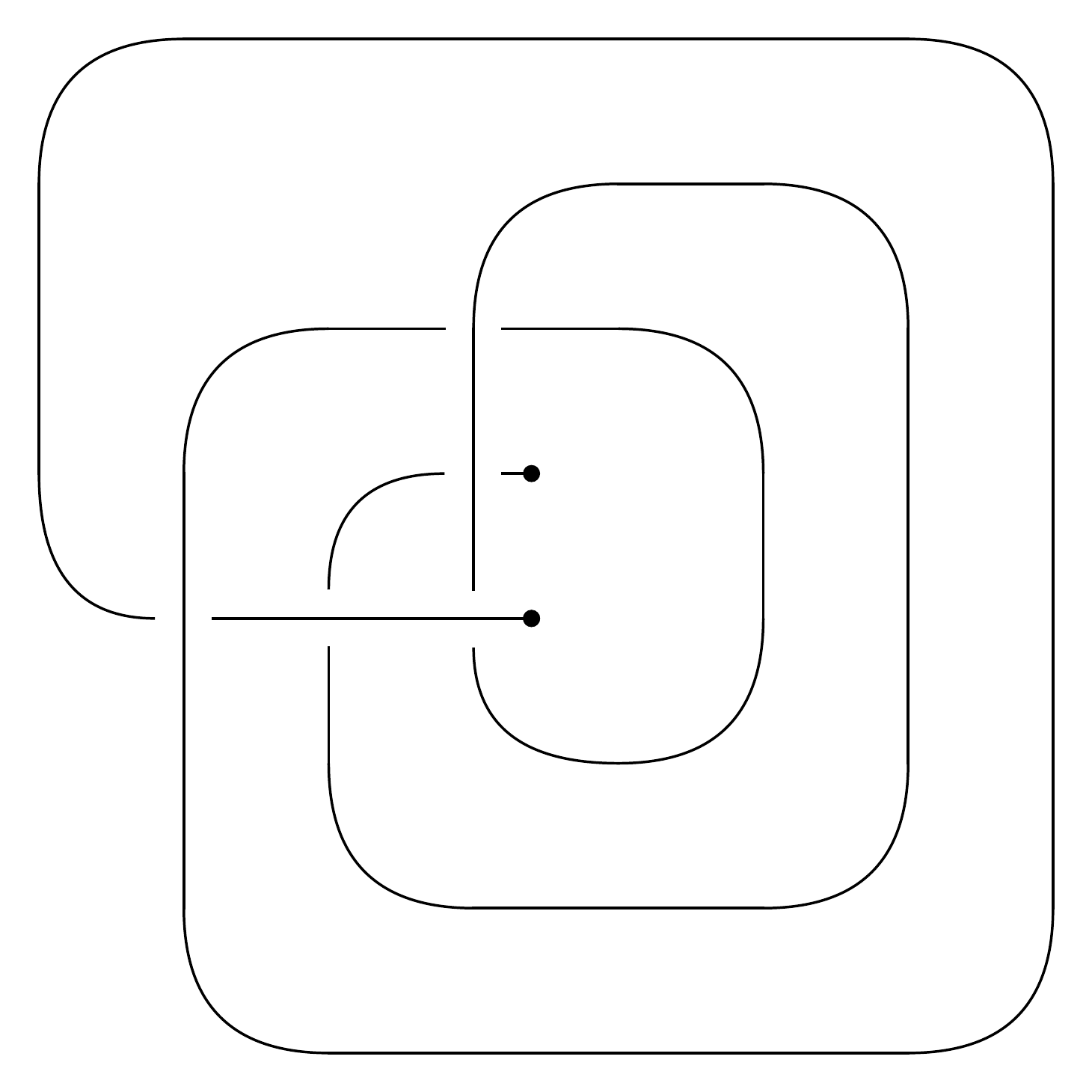}\\
\textcolor{black}{$5_{219}$}
\vspace{1cm}
\end{minipage}
\begin{minipage}[t]{.25\linewidth}
\centering
\includegraphics[width=0.9\textwidth,height=3.5cm,keepaspectratio]{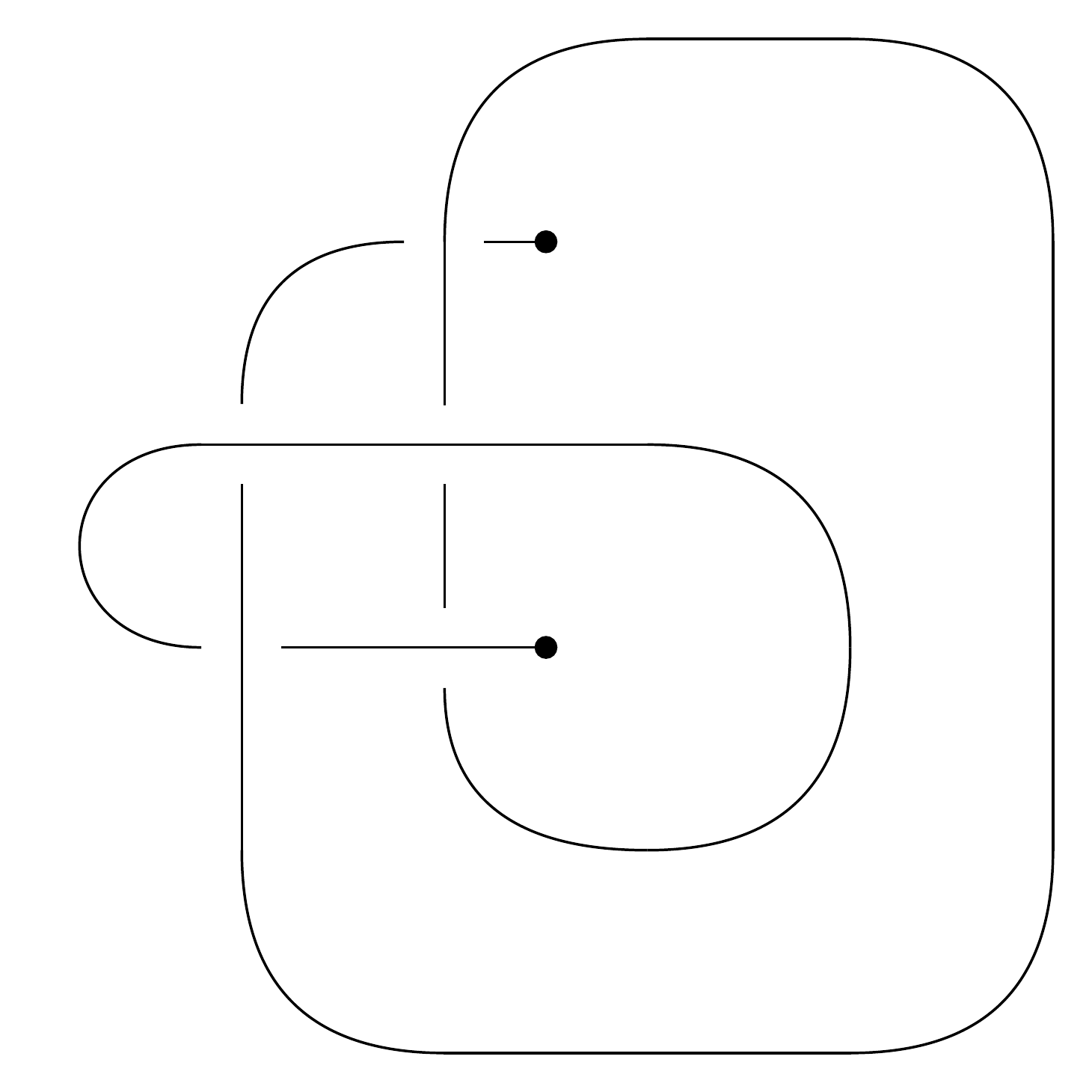}\\
\textcolor{black}{$5_{220}$}
\vspace{1cm}
\end{minipage}
\begin{minipage}[t]{.25\linewidth}
\centering
\includegraphics[width=0.9\textwidth,height=3.5cm,keepaspectratio]{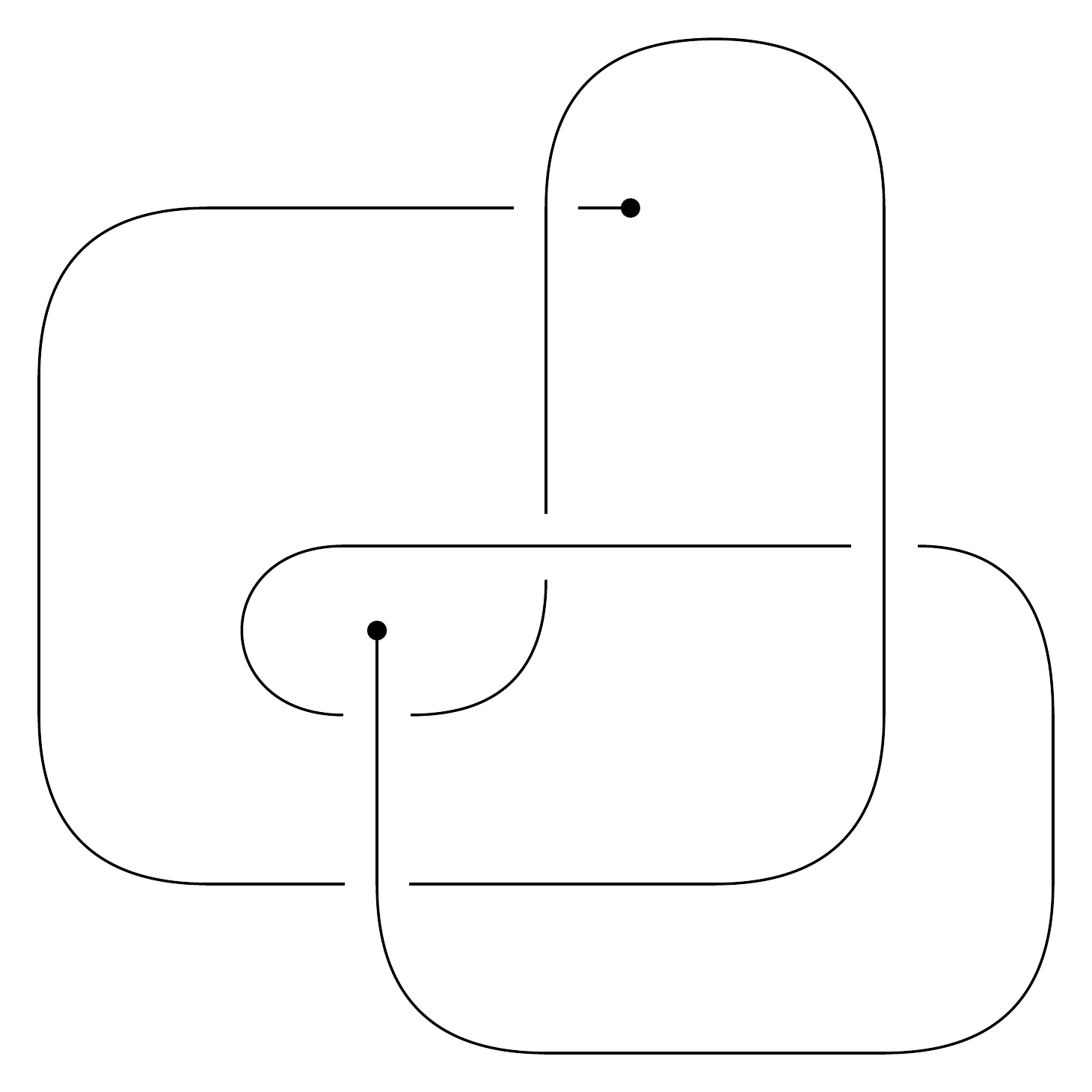}\\
\textcolor{black}{$5_{221}$}
\vspace{1cm}
\end{minipage}
\begin{minipage}[t]{.25\linewidth}
\centering
\includegraphics[width=0.9\textwidth,height=3.5cm,keepaspectratio]{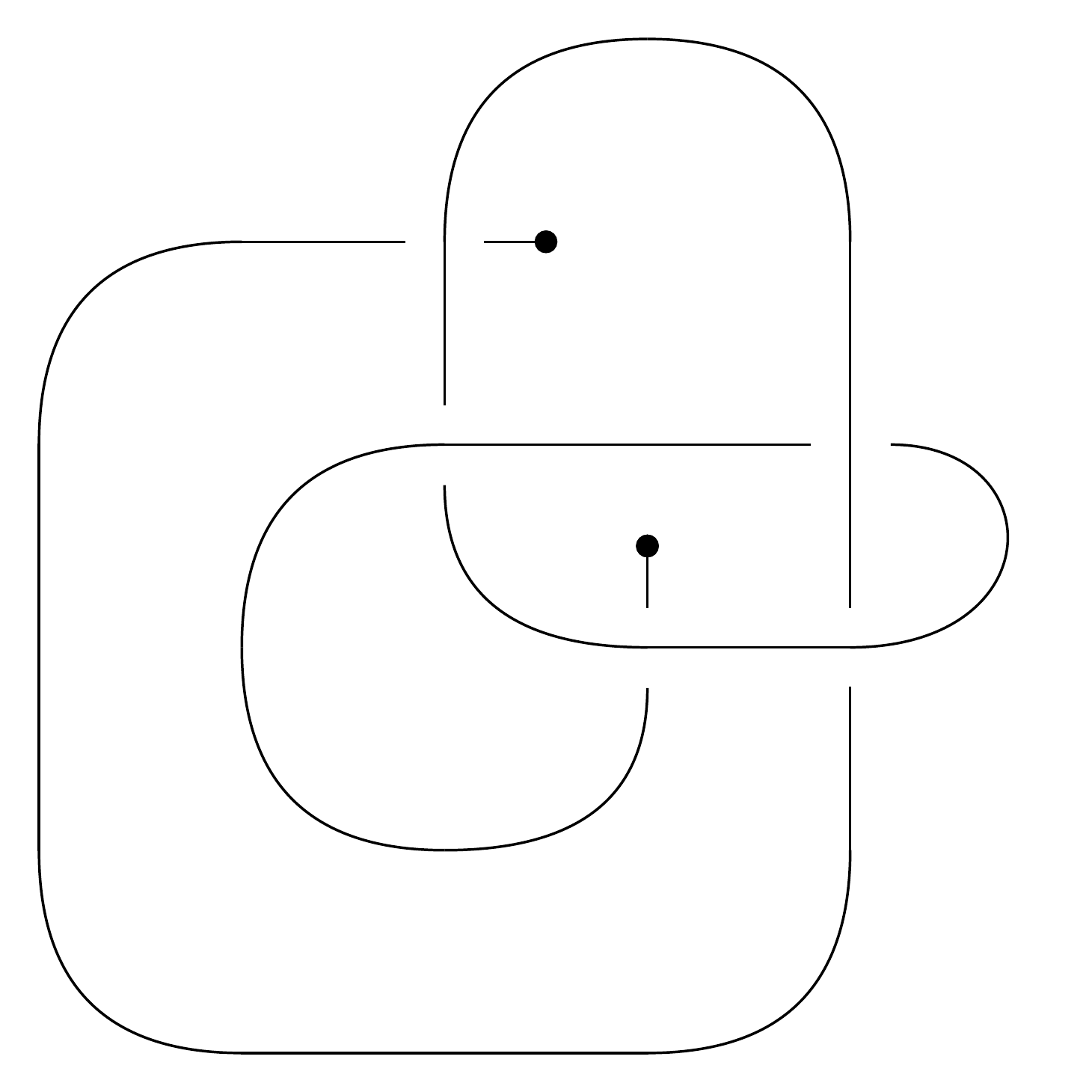}\\
\textcolor{black}{$5_{222}$}
\vspace{1cm}
\end{minipage}
\begin{minipage}[t]{.25\linewidth}
\centering
\includegraphics[width=0.9\textwidth,height=3.5cm,keepaspectratio]{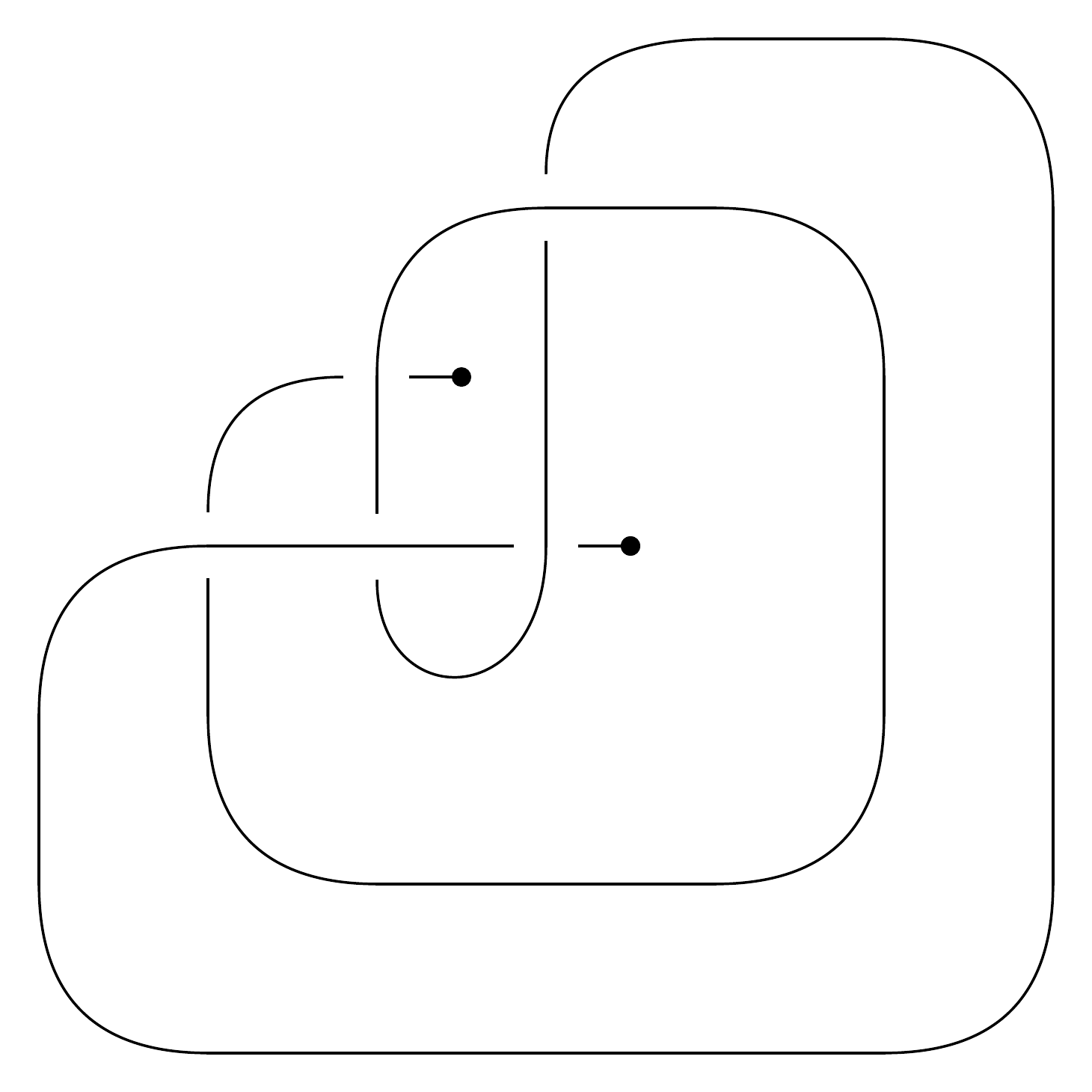}\\
\textcolor{black}{$5_{223}$}
\vspace{1cm}
\end{minipage}
\begin{minipage}[t]{.25\linewidth}
\centering
\includegraphics[width=0.9\textwidth,height=3.5cm,keepaspectratio]{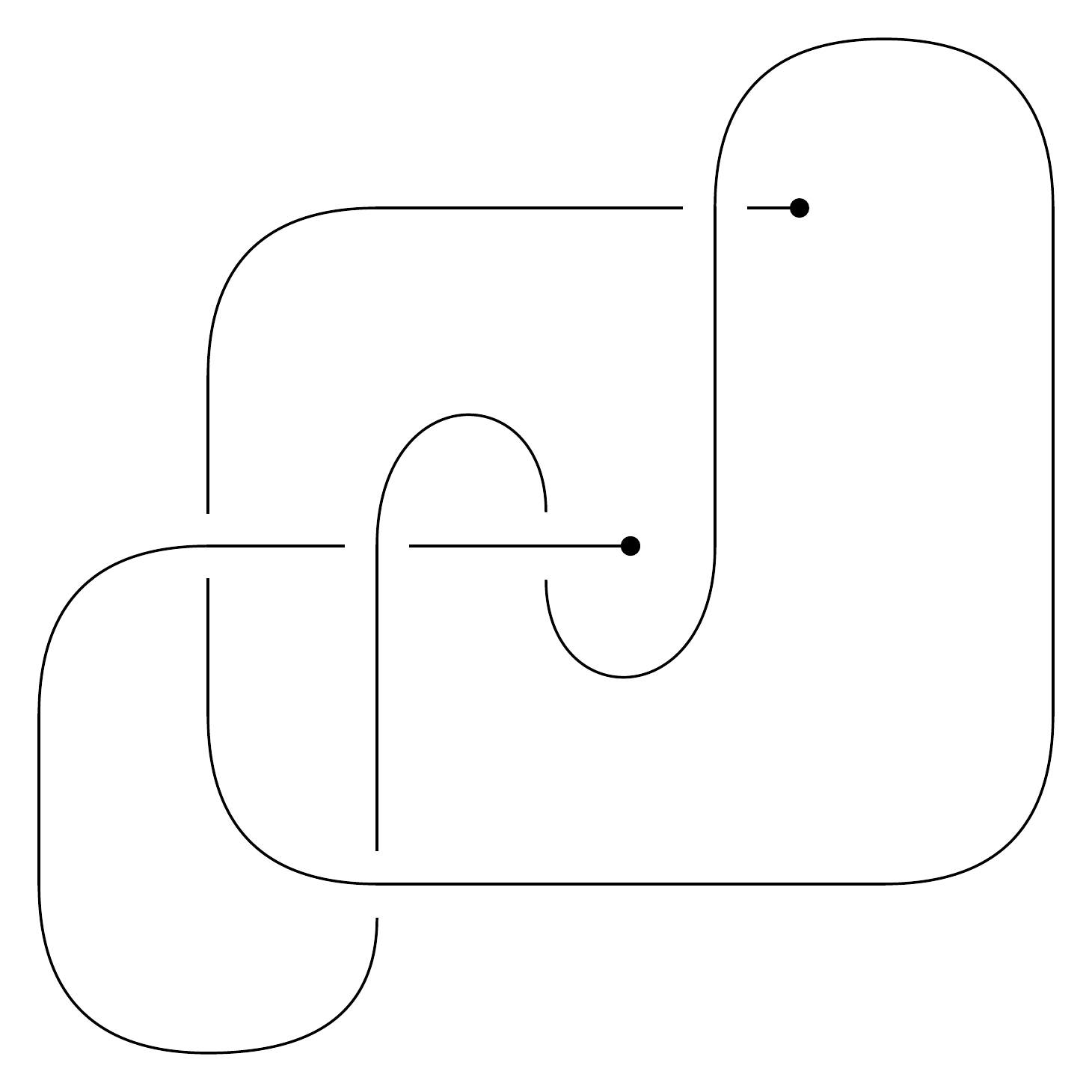}\\
\textcolor{black}{$5_{224}$}
\vspace{1cm}
\end{minipage}
\begin{minipage}[t]{.25\linewidth}
\centering
\includegraphics[width=0.9\textwidth,height=3.5cm,keepaspectratio]{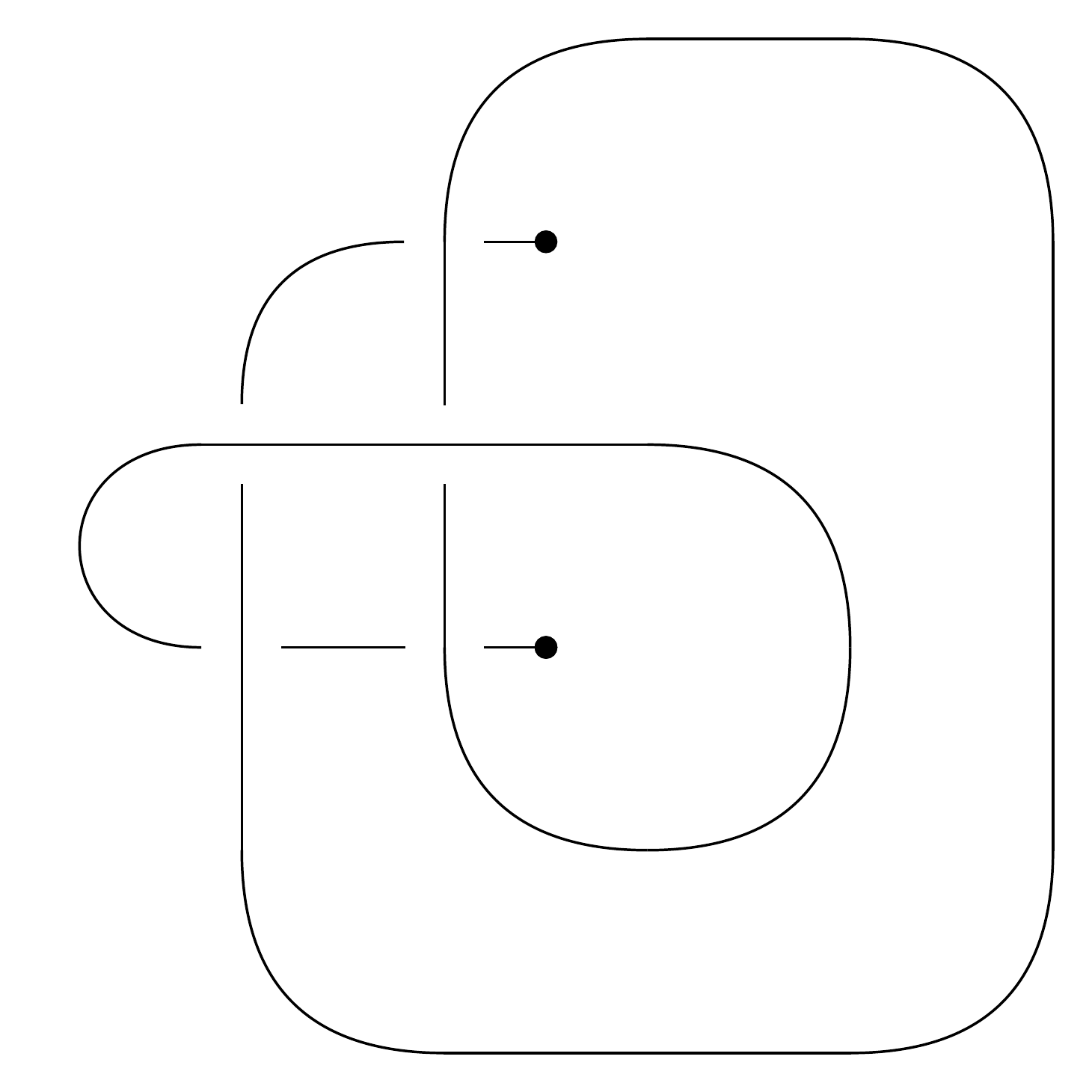}\\
\textcolor{black}{$5_{225}$}
\vspace{1cm}
\end{minipage}
\begin{minipage}[t]{.25\linewidth}
\centering
\includegraphics[width=0.9\textwidth,height=3.5cm,keepaspectratio]{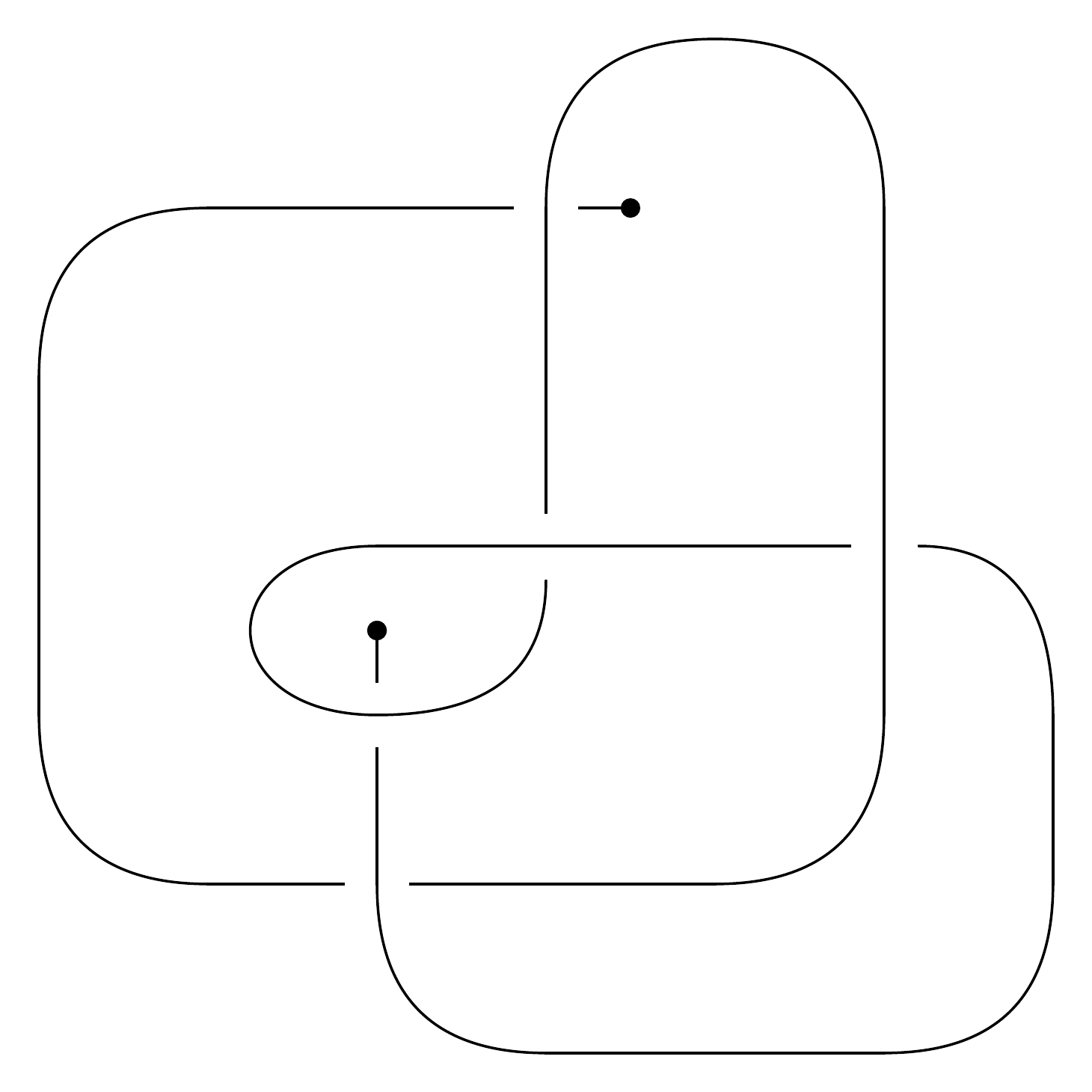}\\
\textcolor{black}{$5_{226}$}
\vspace{1cm}
\end{minipage}
\begin{minipage}[t]{.25\linewidth}
\centering
\includegraphics[width=0.9\textwidth,height=3.5cm,keepaspectratio]{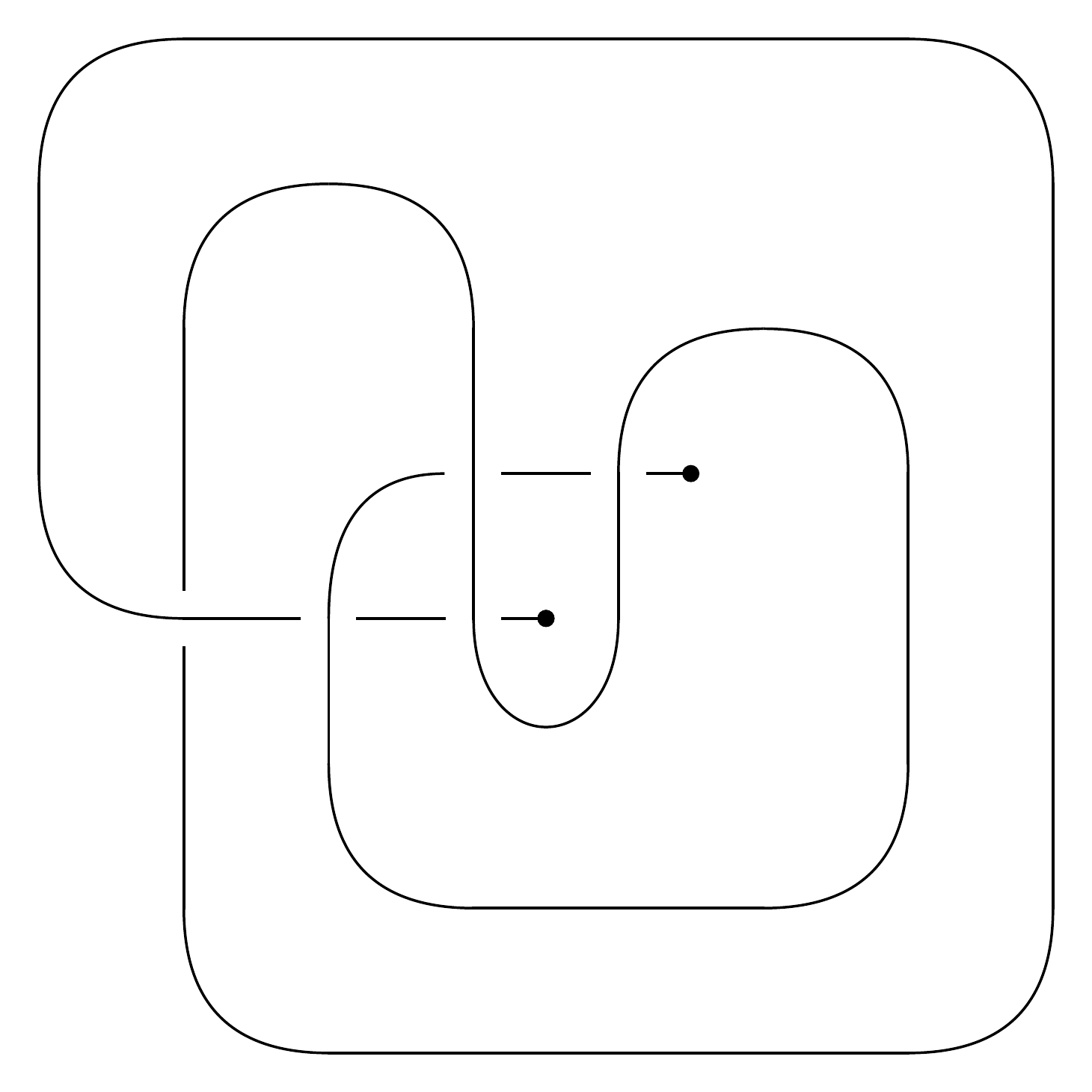}\\
\textcolor{black}{$5_{227}$}
\vspace{1cm}
\end{minipage}
\begin{minipage}[t]{.25\linewidth}
\centering
\includegraphics[width=0.9\textwidth,height=3.5cm,keepaspectratio]{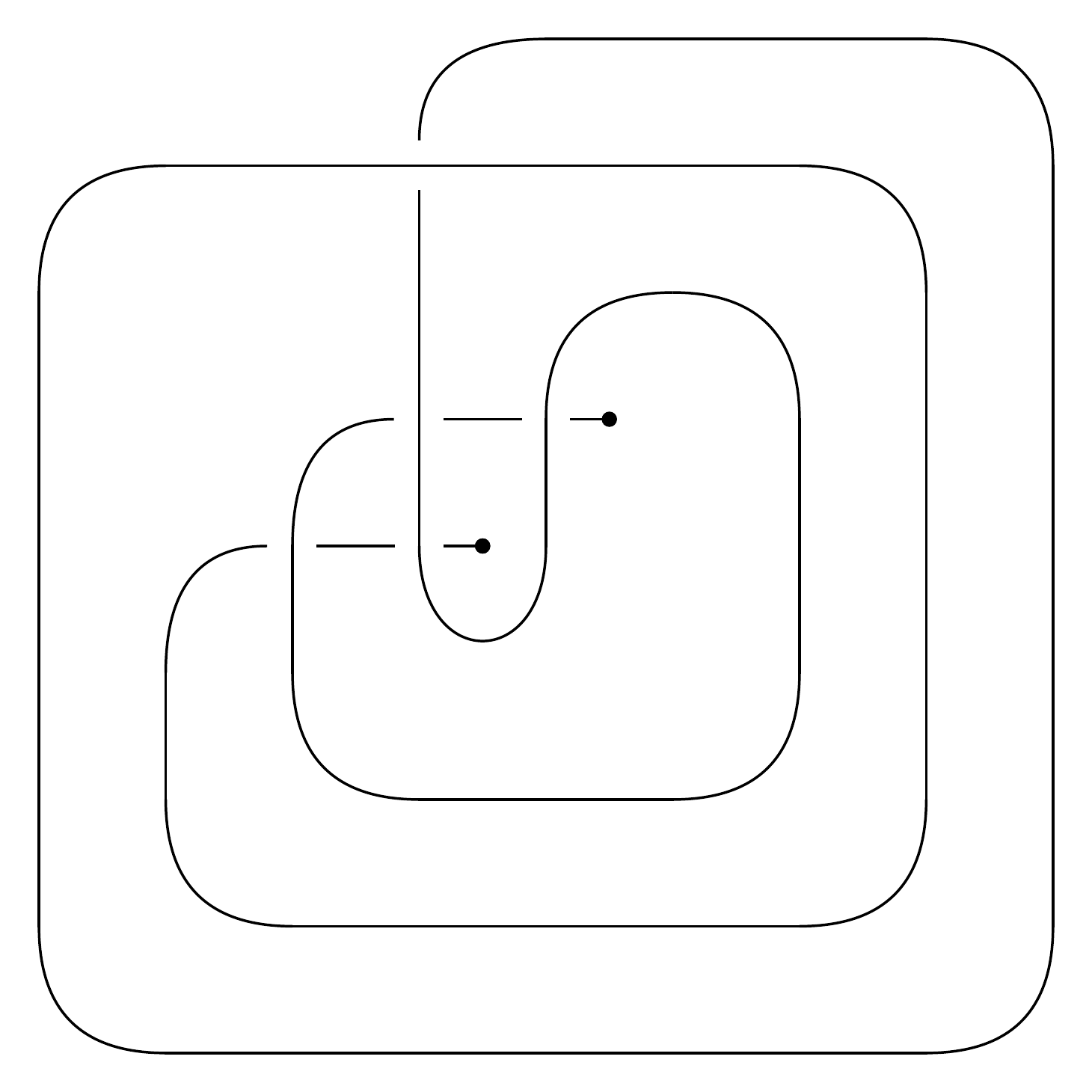}\\
\textcolor{black}{$5_{228}$}
\vspace{1cm}
\end{minipage}
\begin{minipage}[t]{.25\linewidth}
\centering
\includegraphics[width=0.9\textwidth,height=3.5cm,keepaspectratio]{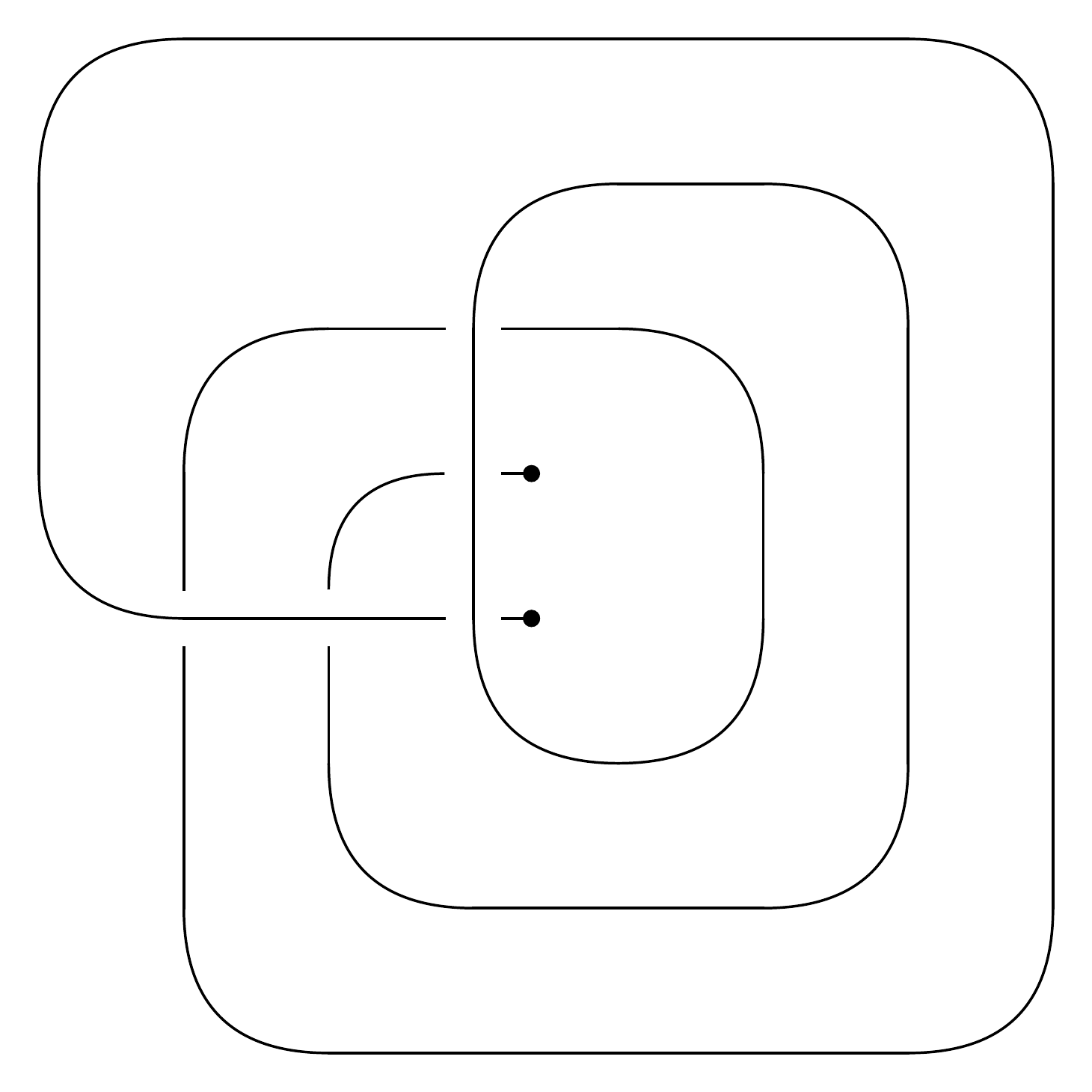}\\
\textcolor{black}{$5_{229}$}
\vspace{1cm}
\end{minipage}
\begin{minipage}[t]{.25\linewidth}
\centering
\includegraphics[width=0.9\textwidth,height=3.5cm,keepaspectratio]{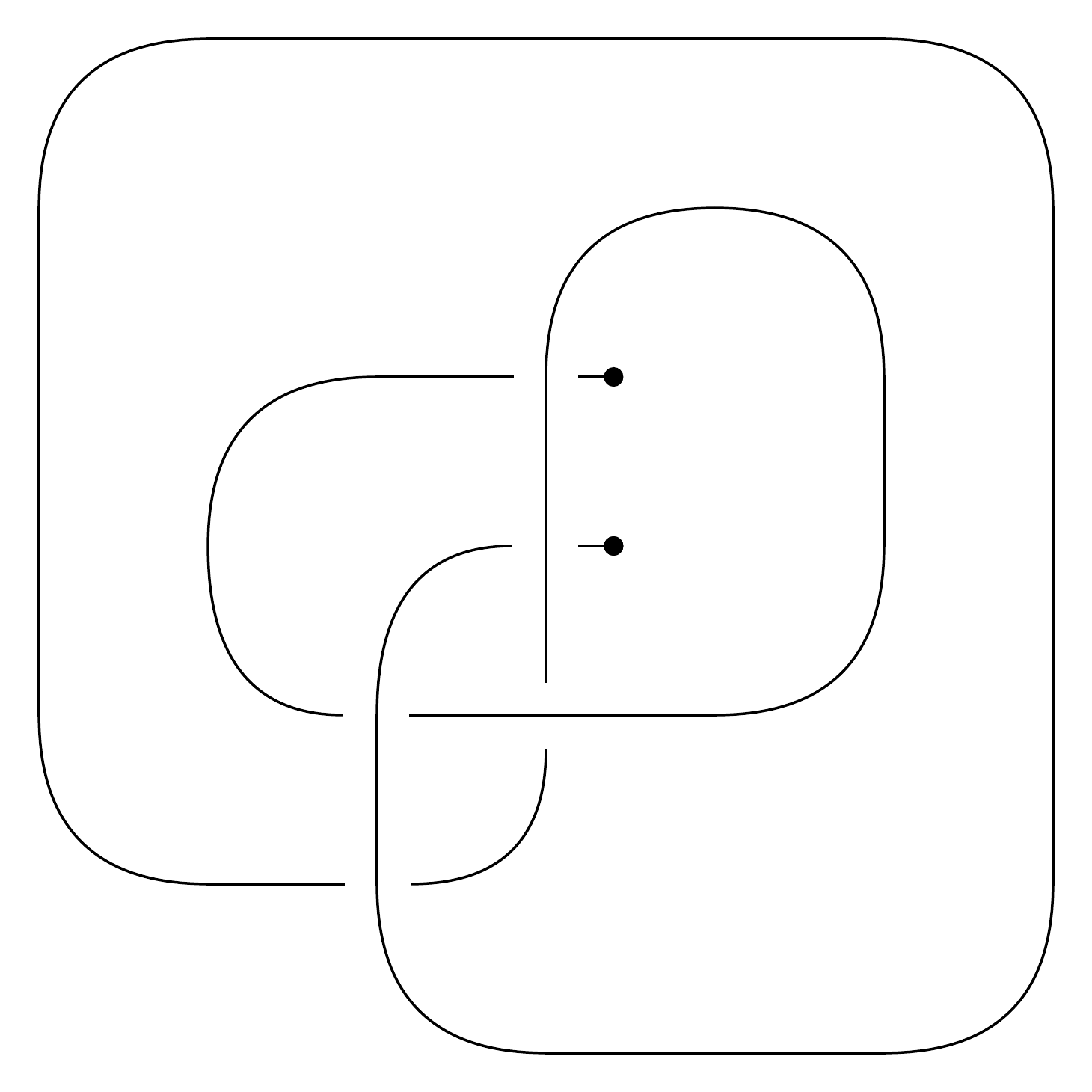}\\
\textcolor{black}{$5_{230}$}
\vspace{1cm}
\end{minipage}
\begin{minipage}[t]{.25\linewidth}
\centering
\includegraphics[width=0.9\textwidth,height=3.5cm,keepaspectratio]{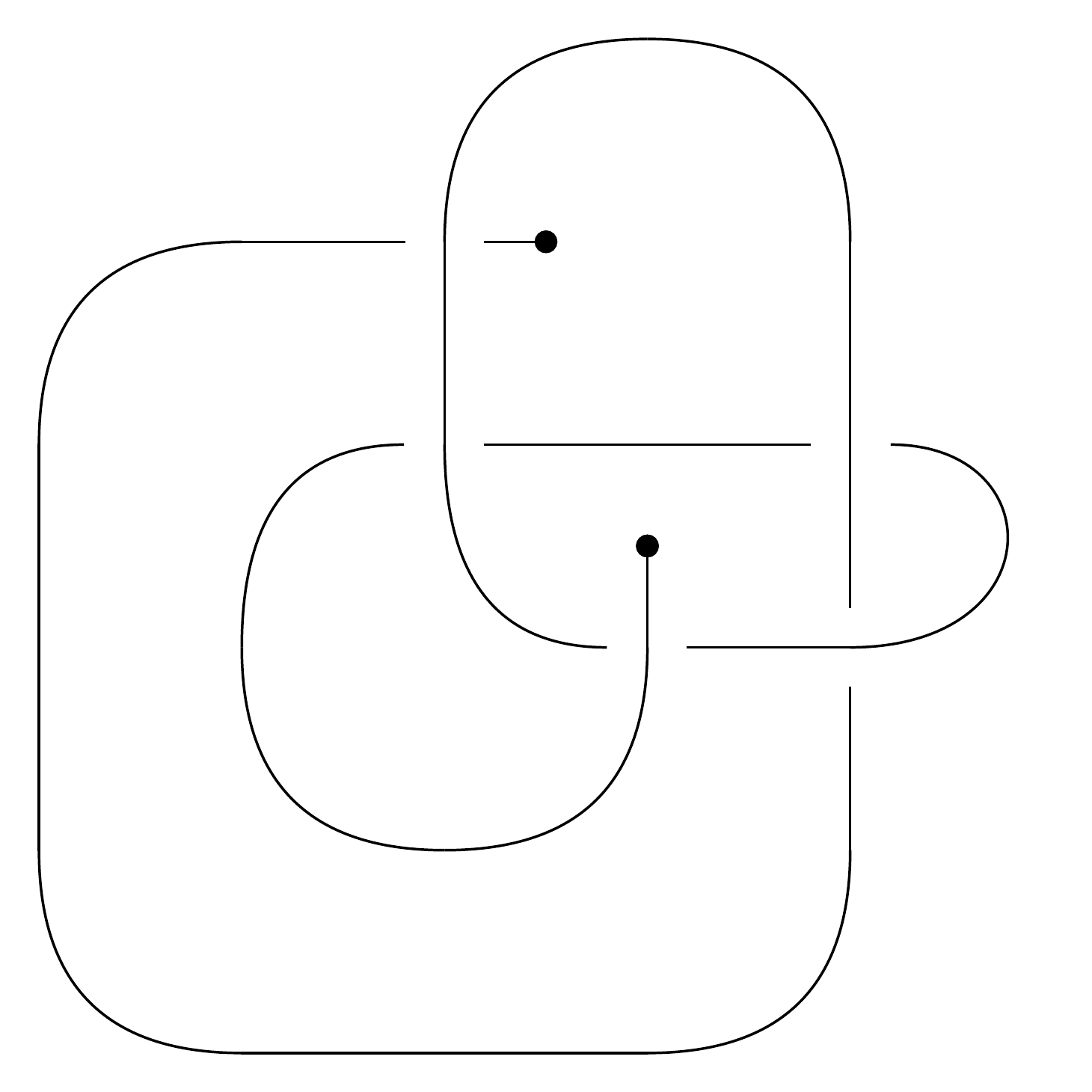}\\
\textcolor{black}{$5_{231}$}
\vspace{1cm}
\end{minipage}
\begin{minipage}[t]{.25\linewidth}
\centering
\includegraphics[width=0.9\textwidth,height=3.5cm,keepaspectratio]{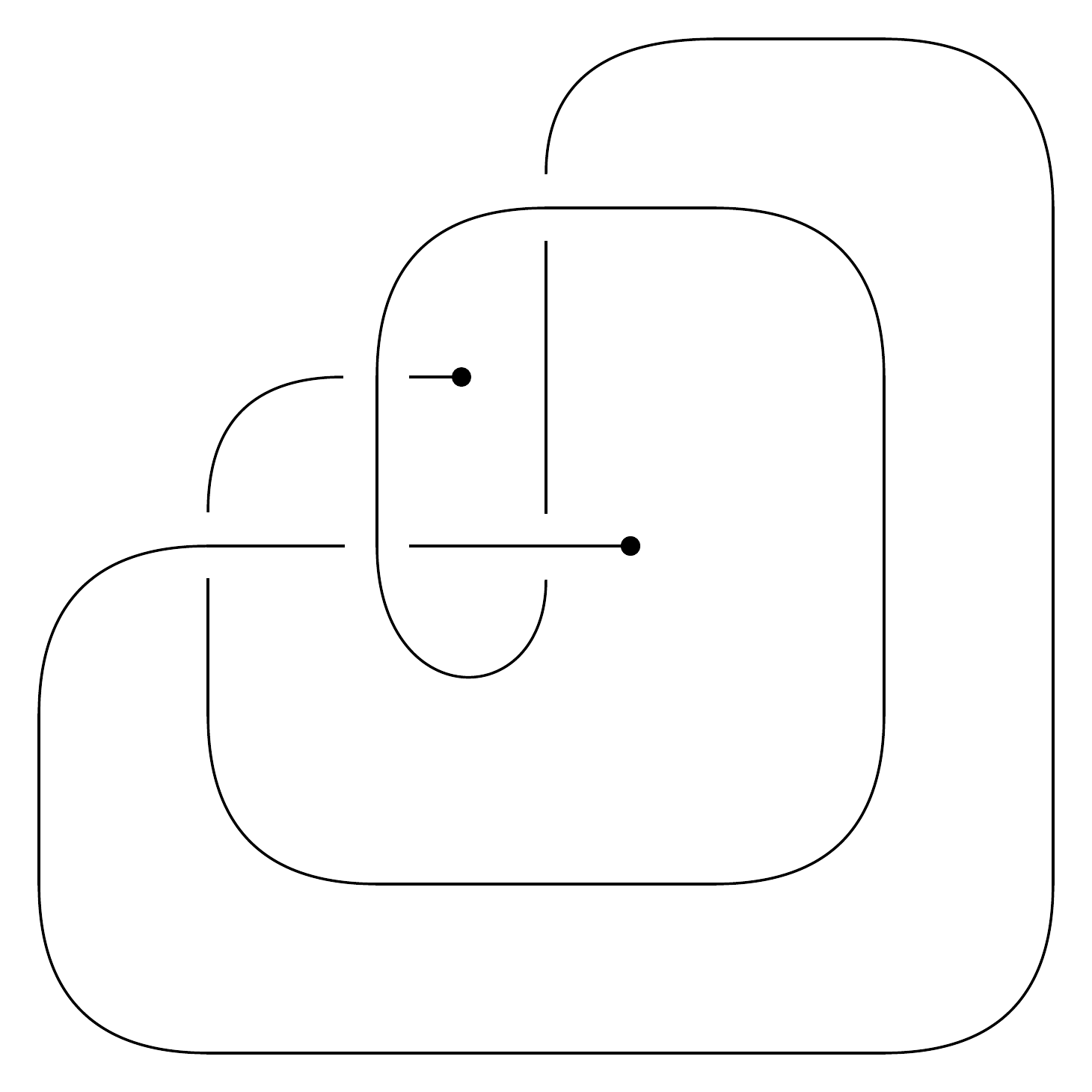}\\
\textcolor{black}{$5_{232}$}
\vspace{1cm}
\end{minipage}
\begin{minipage}[t]{.25\linewidth}
\centering
\includegraphics[width=0.9\textwidth,height=3.5cm,keepaspectratio]{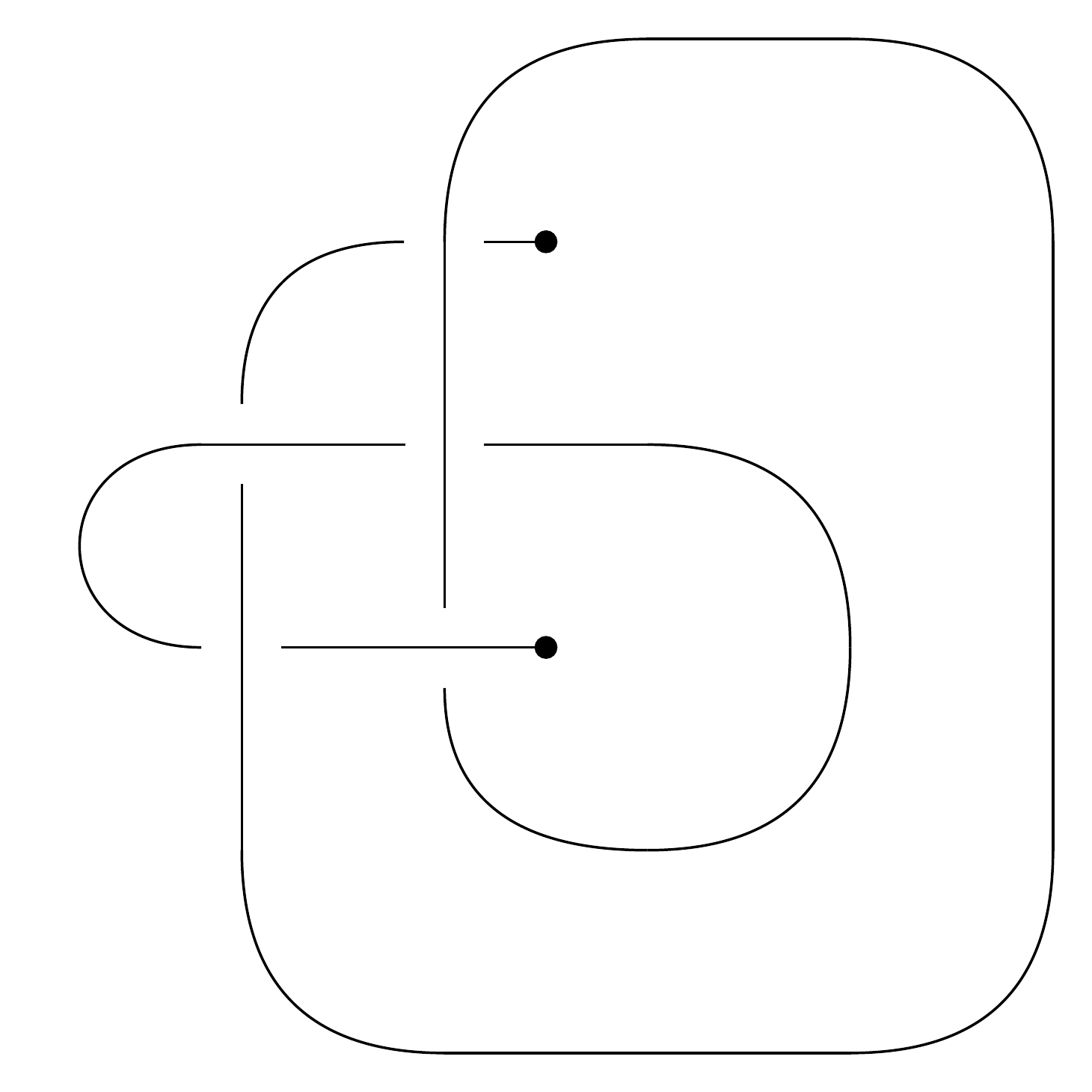}\\
\textcolor{black}{$5_{233}$}
\vspace{1cm}
\end{minipage}
\begin{minipage}[t]{.25\linewidth}
\centering
\includegraphics[width=0.9\textwidth,height=3.5cm,keepaspectratio]{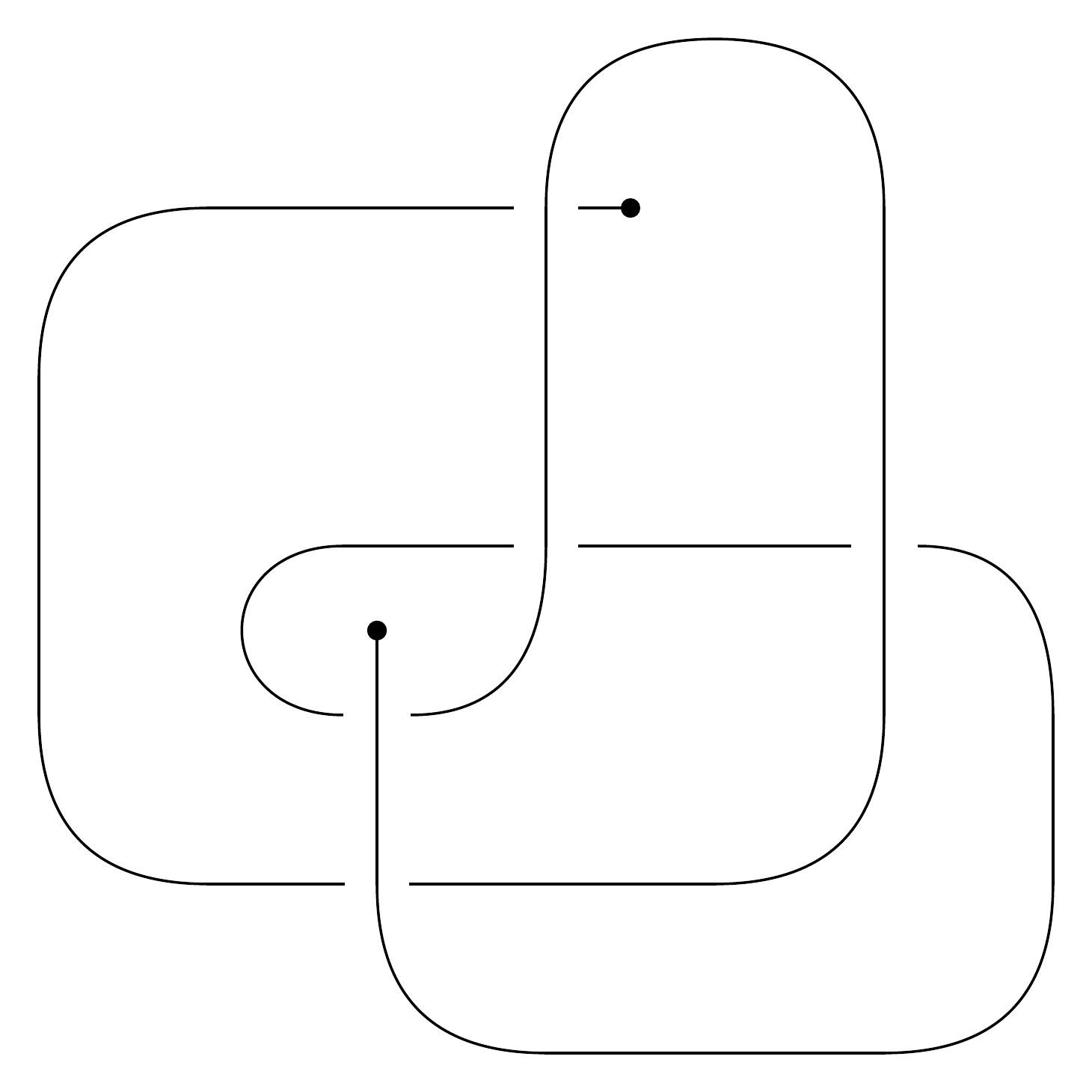}\\
\textcolor{black}{$5_{234}$}
\vspace{1cm}
\end{minipage}
\begin{minipage}[t]{.25\linewidth}
\centering
\includegraphics[width=0.9\textwidth,height=3.5cm,keepaspectratio]{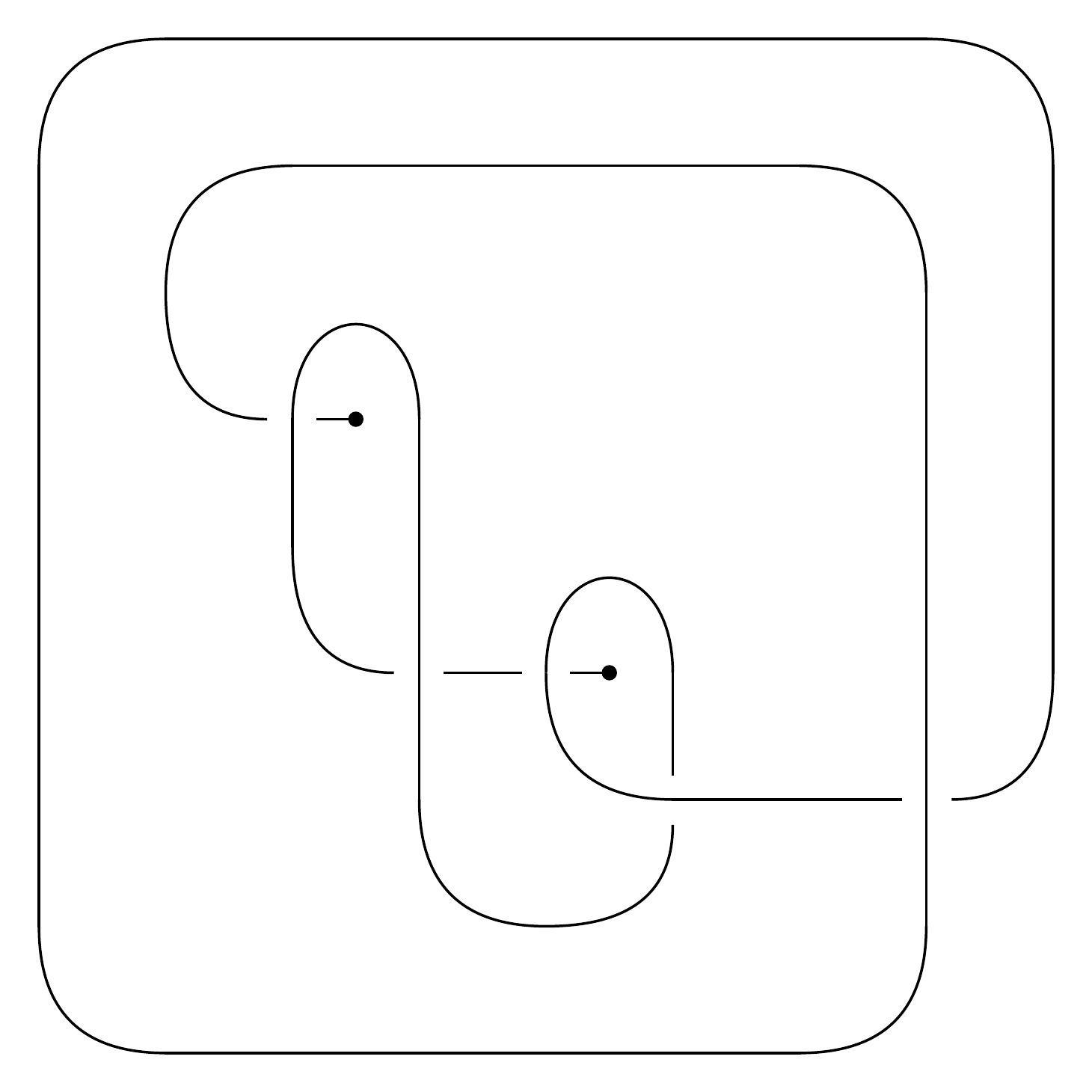}\\
\textcolor{black}{$5_{235}$}
\vspace{1cm}
\end{minipage}
\begin{minipage}[t]{.25\linewidth}
\centering
\includegraphics[width=0.9\textwidth,height=3.5cm,keepaspectratio]{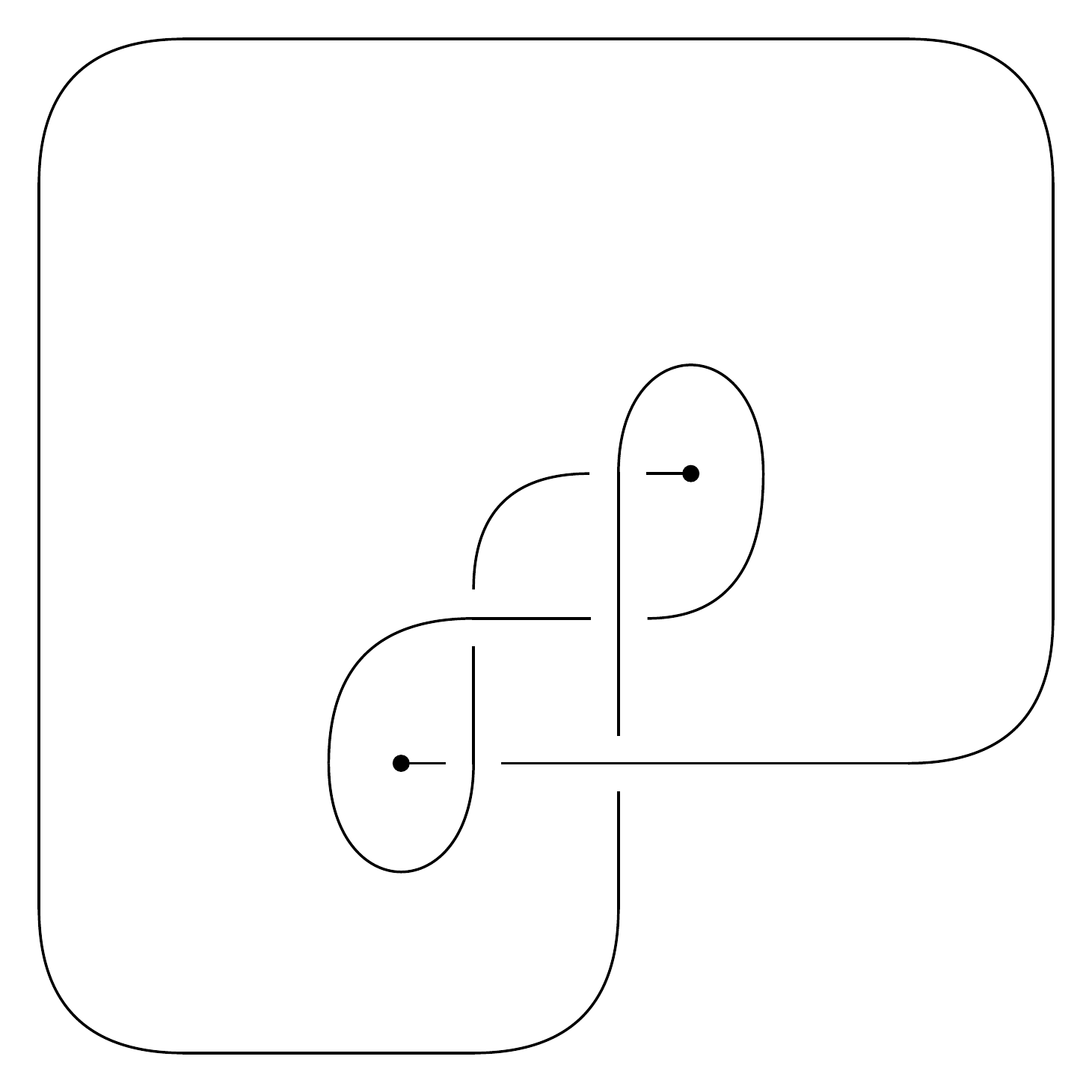}\\
\textcolor{black}{$5_{236}$}
\vspace{1cm}
\end{minipage}
\begin{minipage}[t]{.25\linewidth}
\centering
\includegraphics[width=0.9\textwidth,height=3.5cm,keepaspectratio]{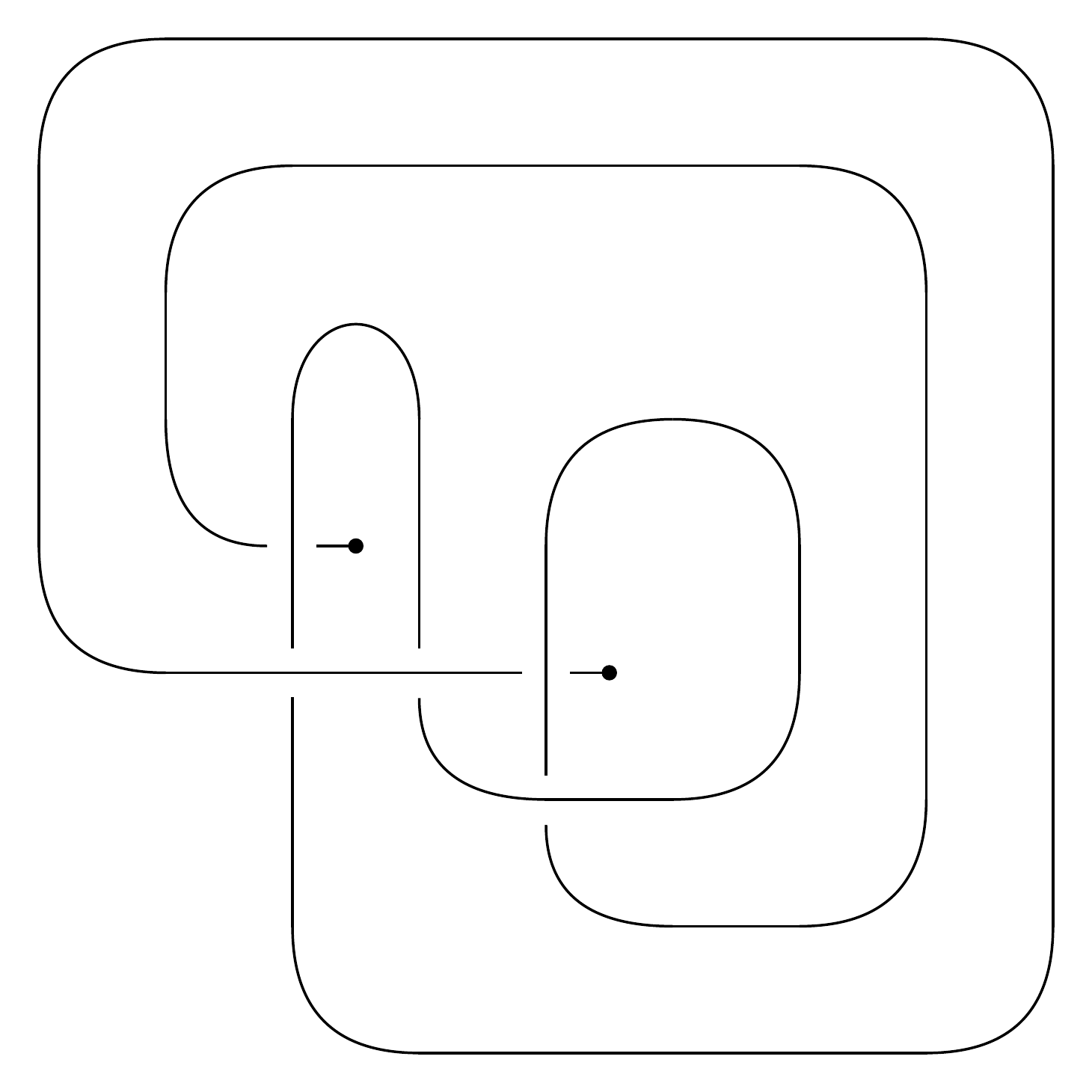}\\
\textcolor{black}{$5_{237}$}
\vspace{1cm}
\end{minipage}
\begin{minipage}[t]{.25\linewidth}
\centering
\includegraphics[width=0.9\textwidth,height=3.5cm,keepaspectratio]{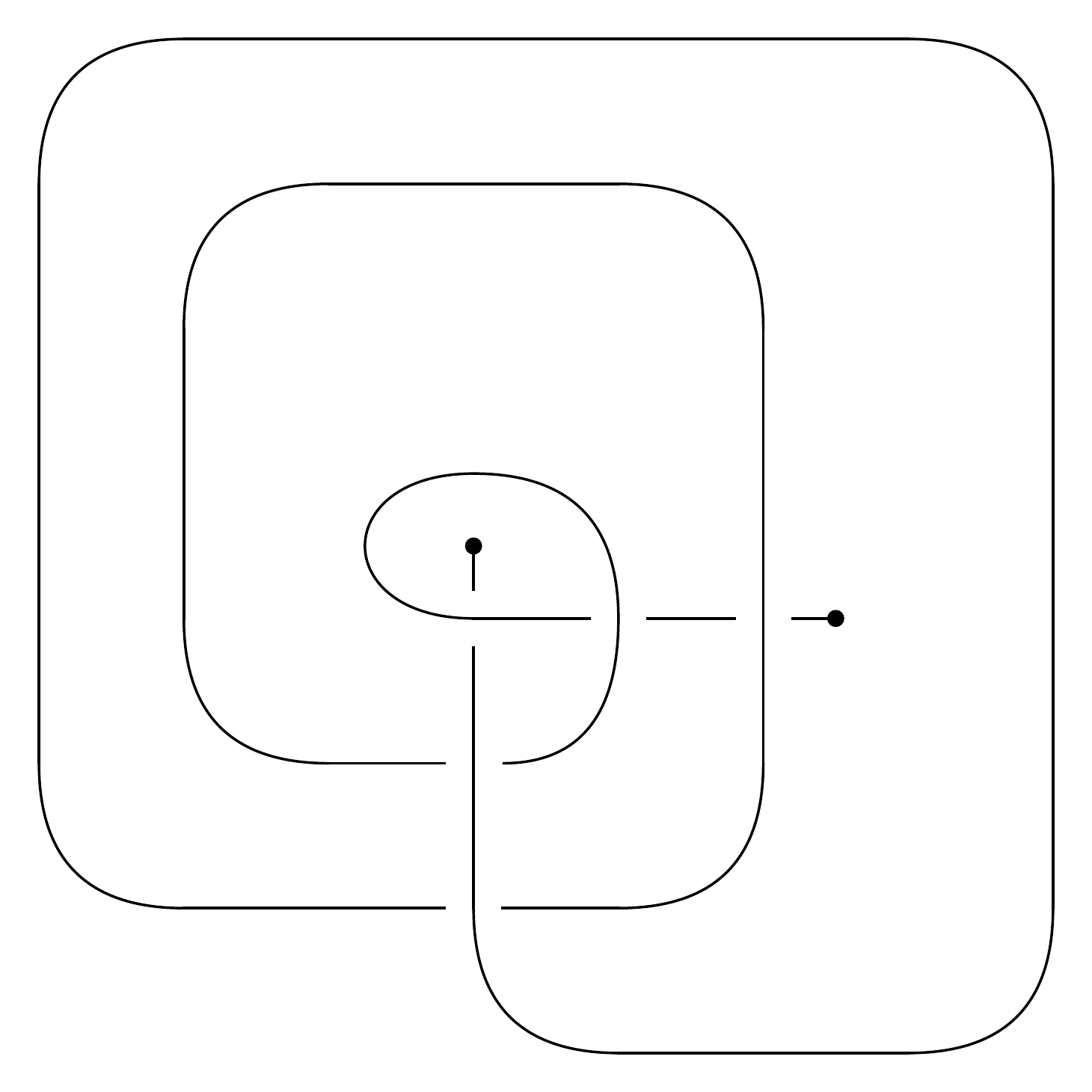}\\
\textcolor{black}{$5_{238}$}
\vspace{1cm}
\end{minipage}
\begin{minipage}[t]{.25\linewidth}
\centering
\includegraphics[width=0.9\textwidth,height=3.5cm,keepaspectratio]{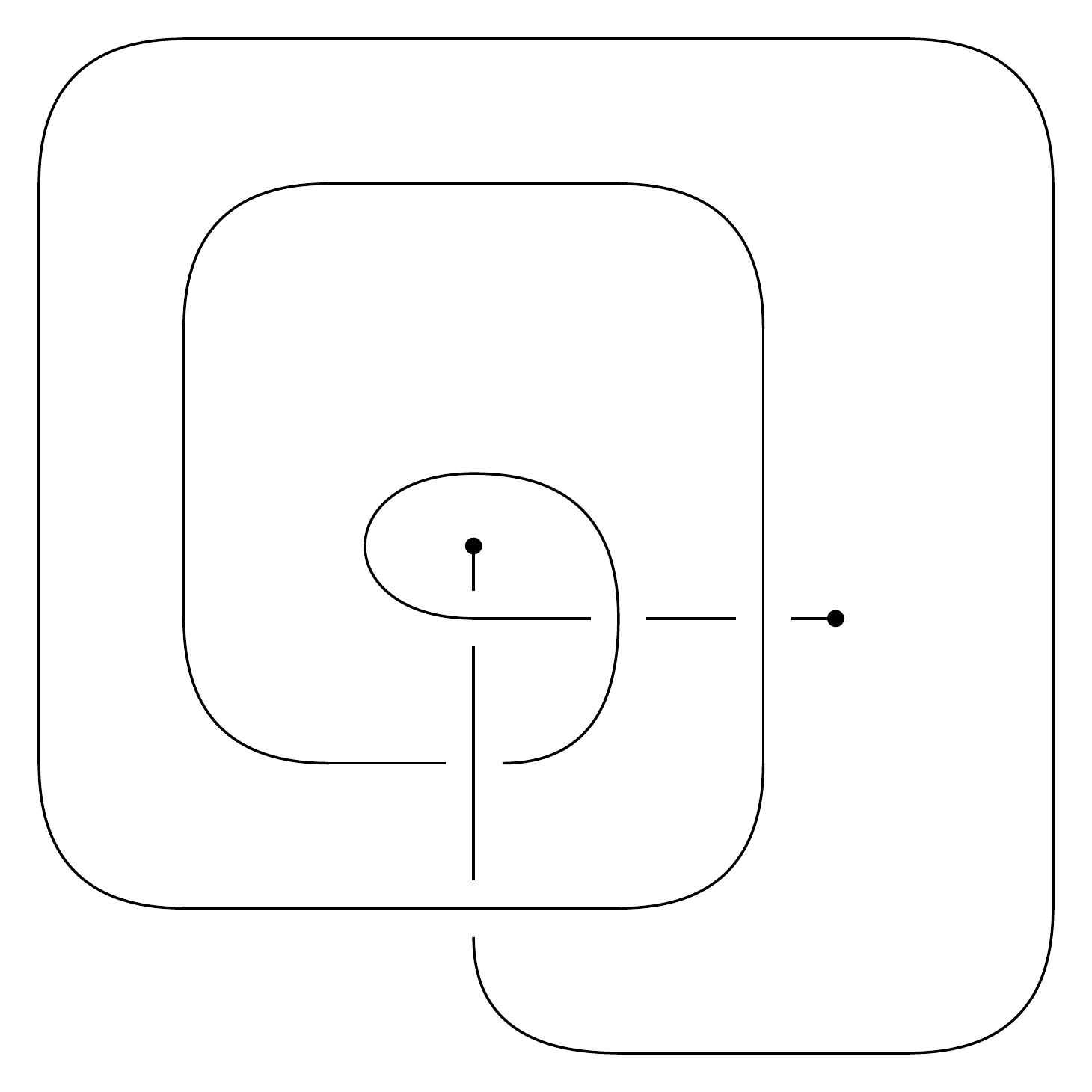}\\
\textcolor{black}{$5_{239}$}
\vspace{1cm}
\end{minipage}
\begin{minipage}[t]{.25\linewidth}
\centering
\includegraphics[width=0.9\textwidth,height=3.5cm,keepaspectratio]{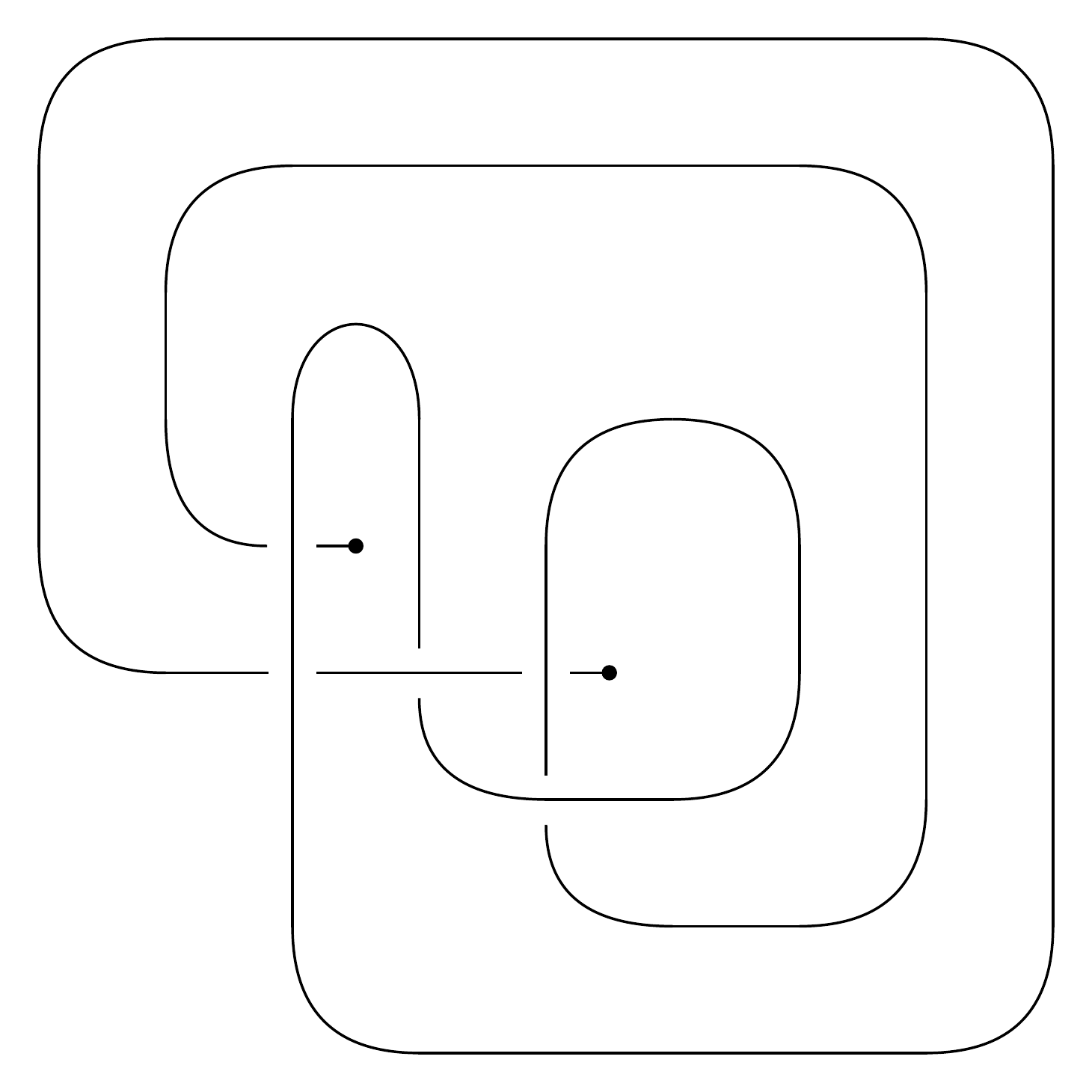}\\
\textcolor{black}{$5_{240}$}
\vspace{1cm}
\end{minipage}
\begin{minipage}[t]{.25\linewidth}
\centering
\includegraphics[width=0.9\textwidth,height=3.5cm,keepaspectratio]{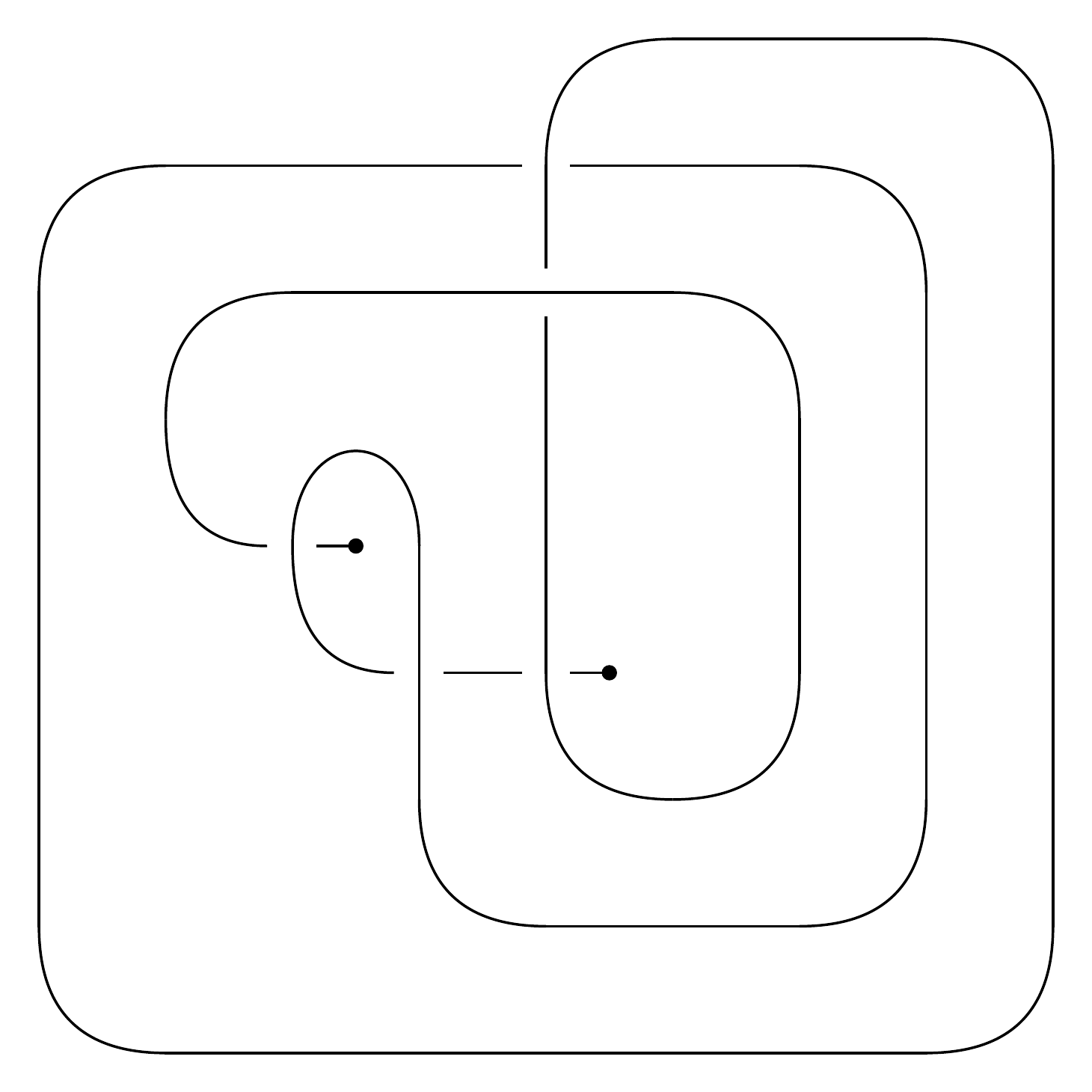}\\
\textcolor{black}{$5_{241}$}
\vspace{1cm}
\end{minipage}
\begin{minipage}[t]{.25\linewidth}
\centering
\includegraphics[width=0.9\textwidth,height=3.5cm,keepaspectratio]{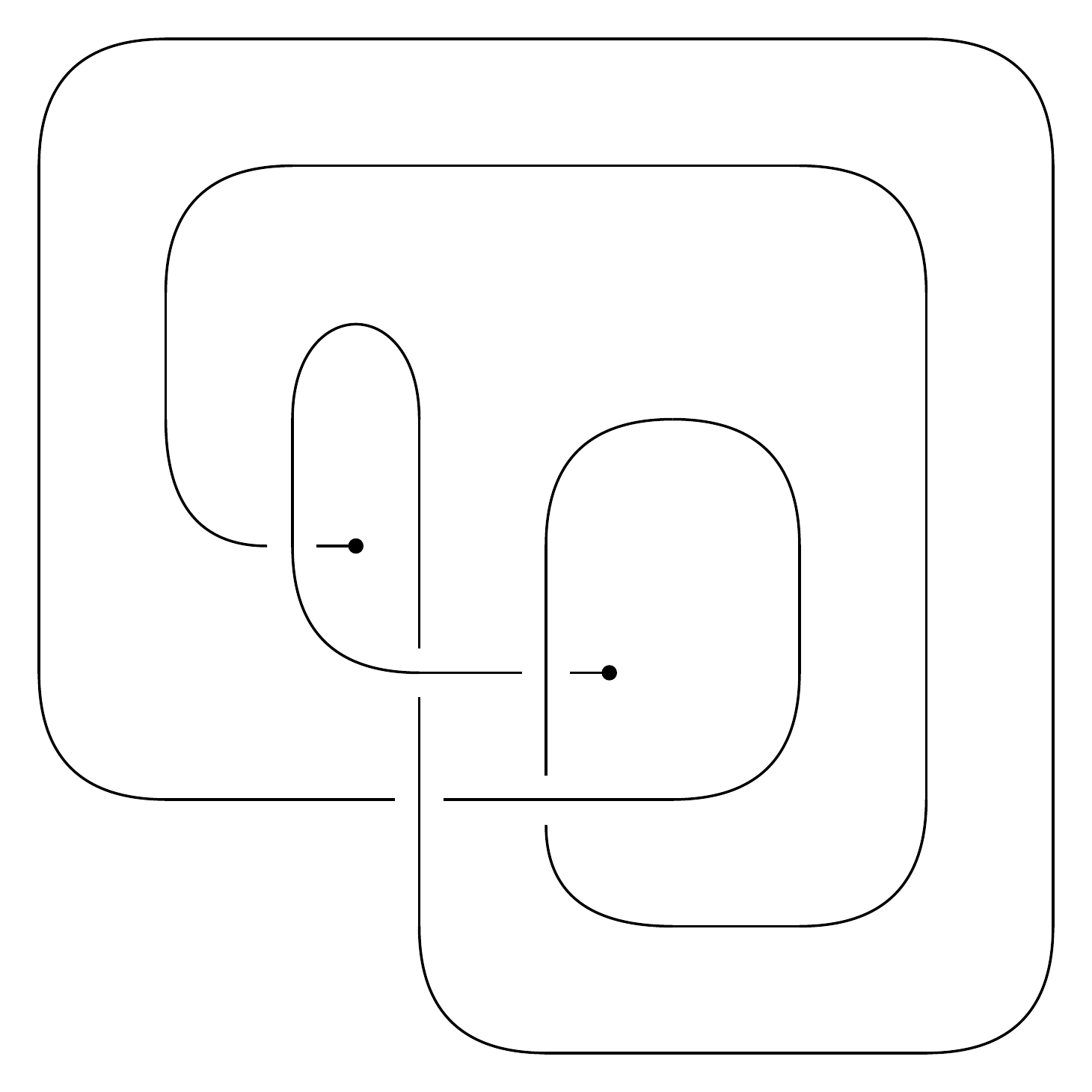}\\
\textcolor{black}{$5_{242}$}
\vspace{1cm}
\end{minipage}
\begin{minipage}[t]{.25\linewidth}
\centering
\includegraphics[width=0.9\textwidth,height=3.5cm,keepaspectratio]{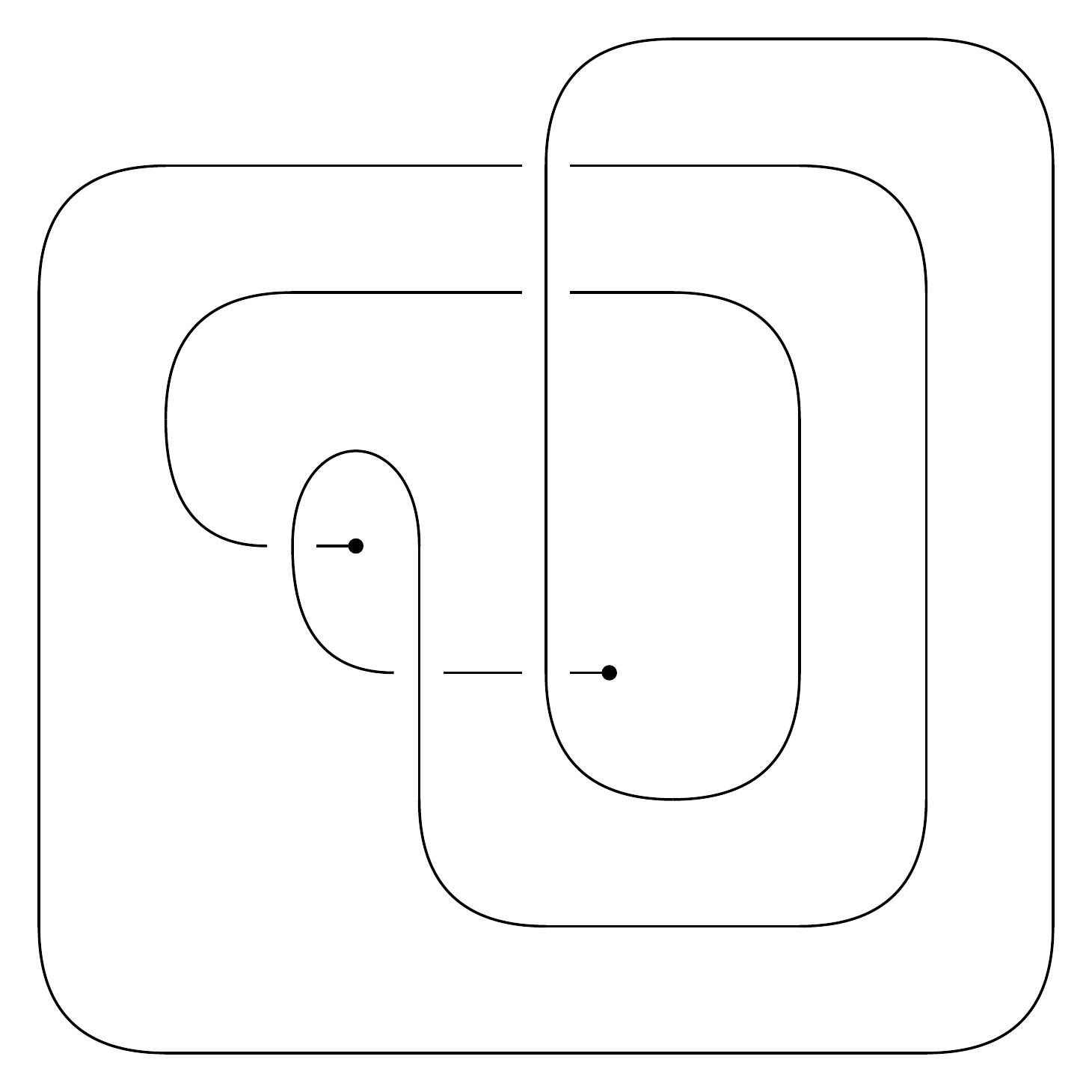}\\
\textcolor{black}{$5_{243}$}
\vspace{1cm}
\end{minipage}
\begin{minipage}[t]{.25\linewidth}
\centering
\includegraphics[width=0.9\textwidth,height=3.5cm,keepaspectratio]{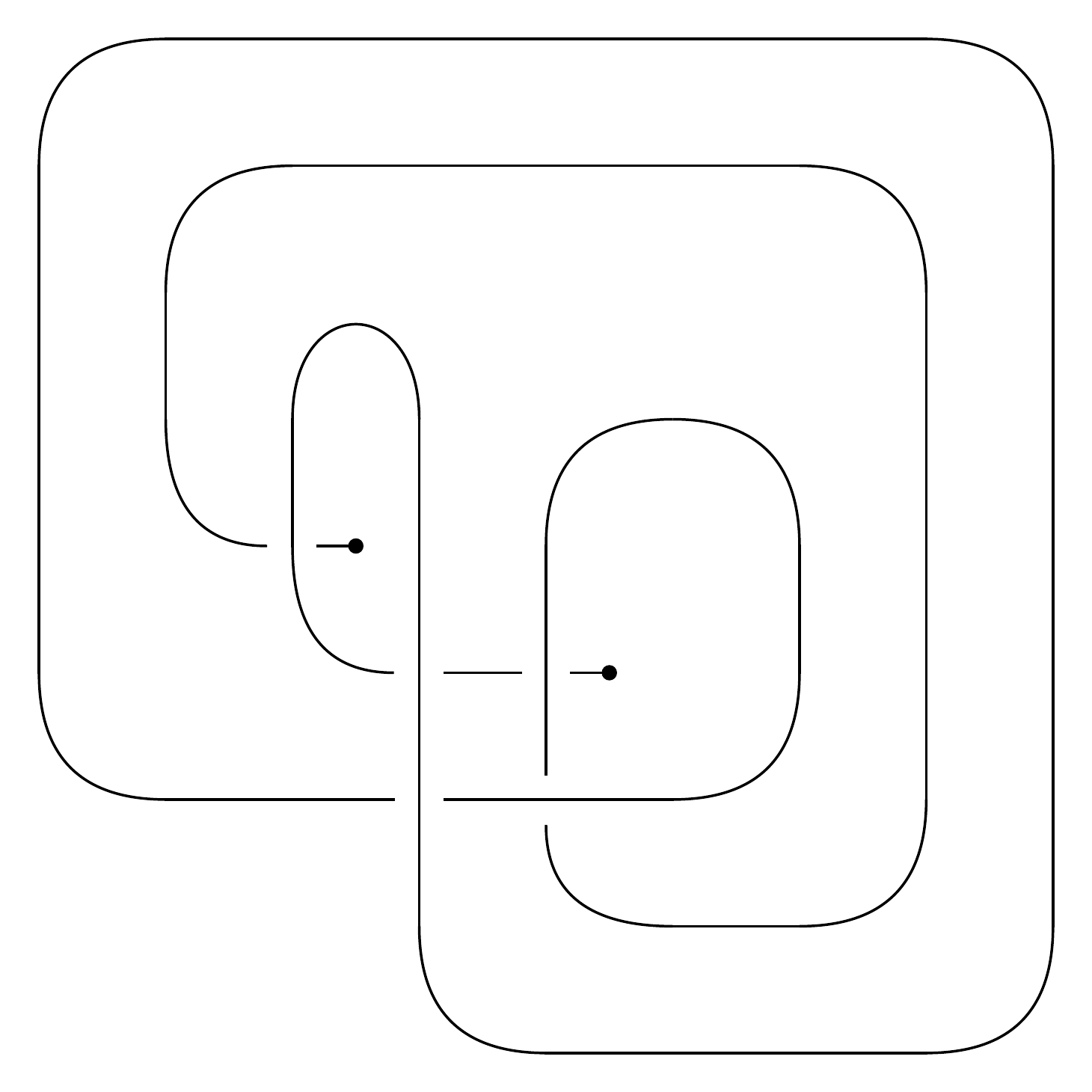}\\
\textcolor{black}{$5_{244}$}
\vspace{1cm}
\end{minipage}
\begin{minipage}[t]{.25\linewidth}
\centering
\includegraphics[width=0.9\textwidth,height=3.5cm,keepaspectratio]{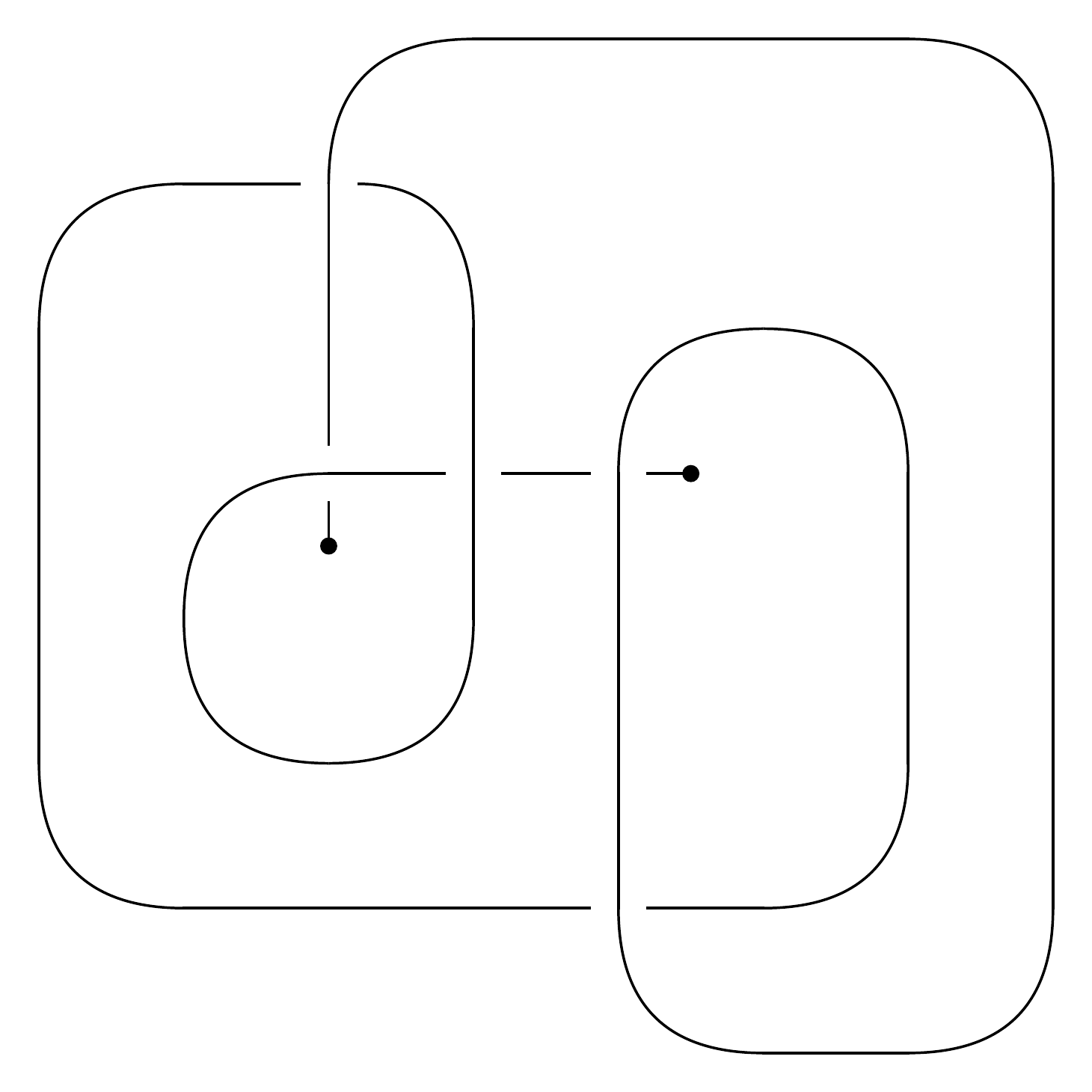}\\
\textcolor{black}{$5_{245}$}
\vspace{1cm}
\end{minipage}
\begin{minipage}[t]{.25\linewidth}
\centering
\includegraphics[width=0.9\textwidth,height=3.5cm,keepaspectratio]{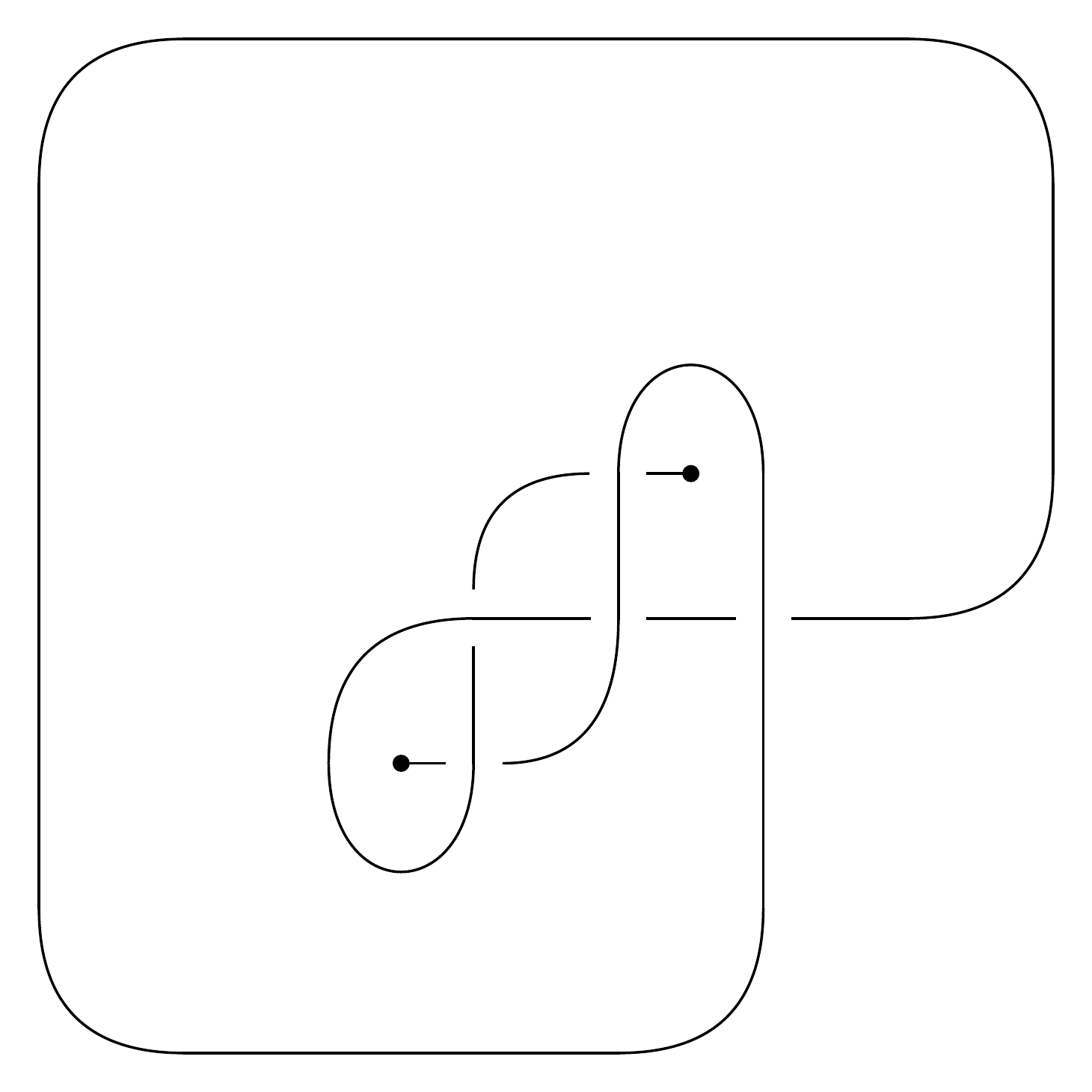}\\
\textcolor{black}{$5_{246}$}
\vspace{1cm}
\end{minipage}
\begin{minipage}[t]{.25\linewidth}
\centering
\includegraphics[width=0.9\textwidth,height=3.5cm,keepaspectratio]{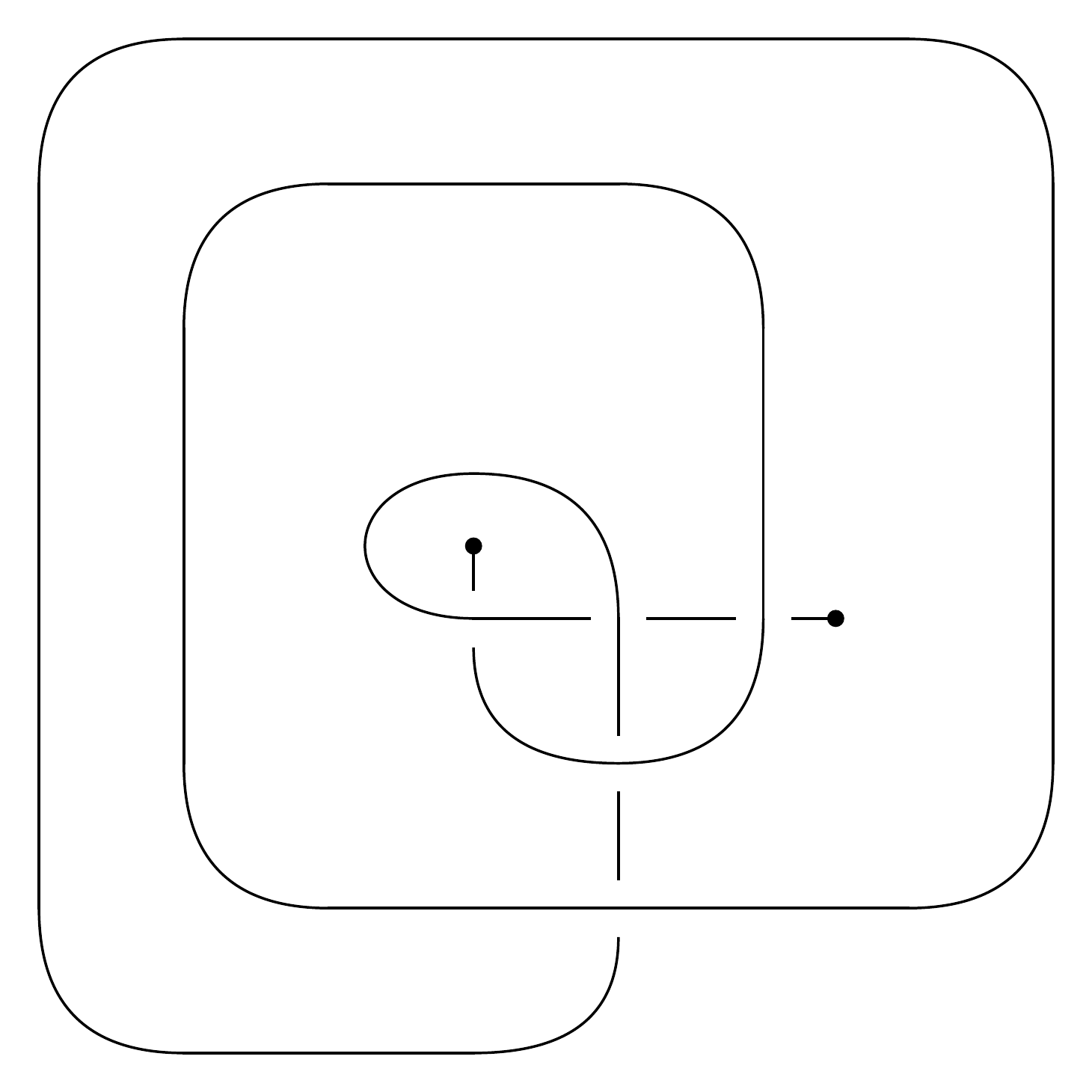}\\
\textcolor{black}{$5_{247}$}
\vspace{1cm}
\end{minipage}
\begin{minipage}[t]{.25\linewidth}
\centering
\includegraphics[width=0.9\textwidth,height=3.5cm,keepaspectratio]{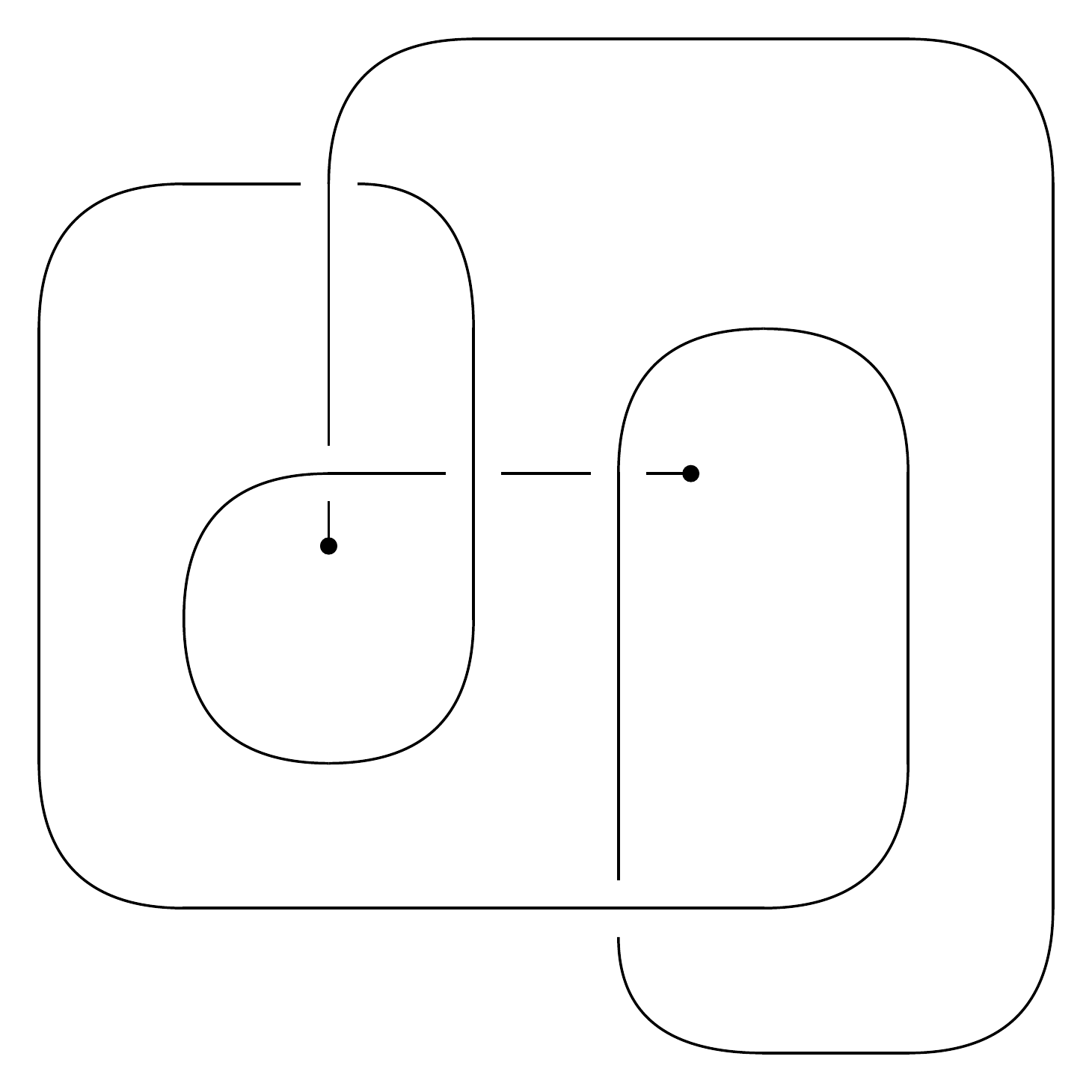}\\
\textcolor{black}{$5_{248}$}
\vspace{1cm}
\end{minipage}
\begin{minipage}[t]{.25\linewidth}
\centering
\includegraphics[width=0.9\textwidth,height=3.5cm,keepaspectratio]{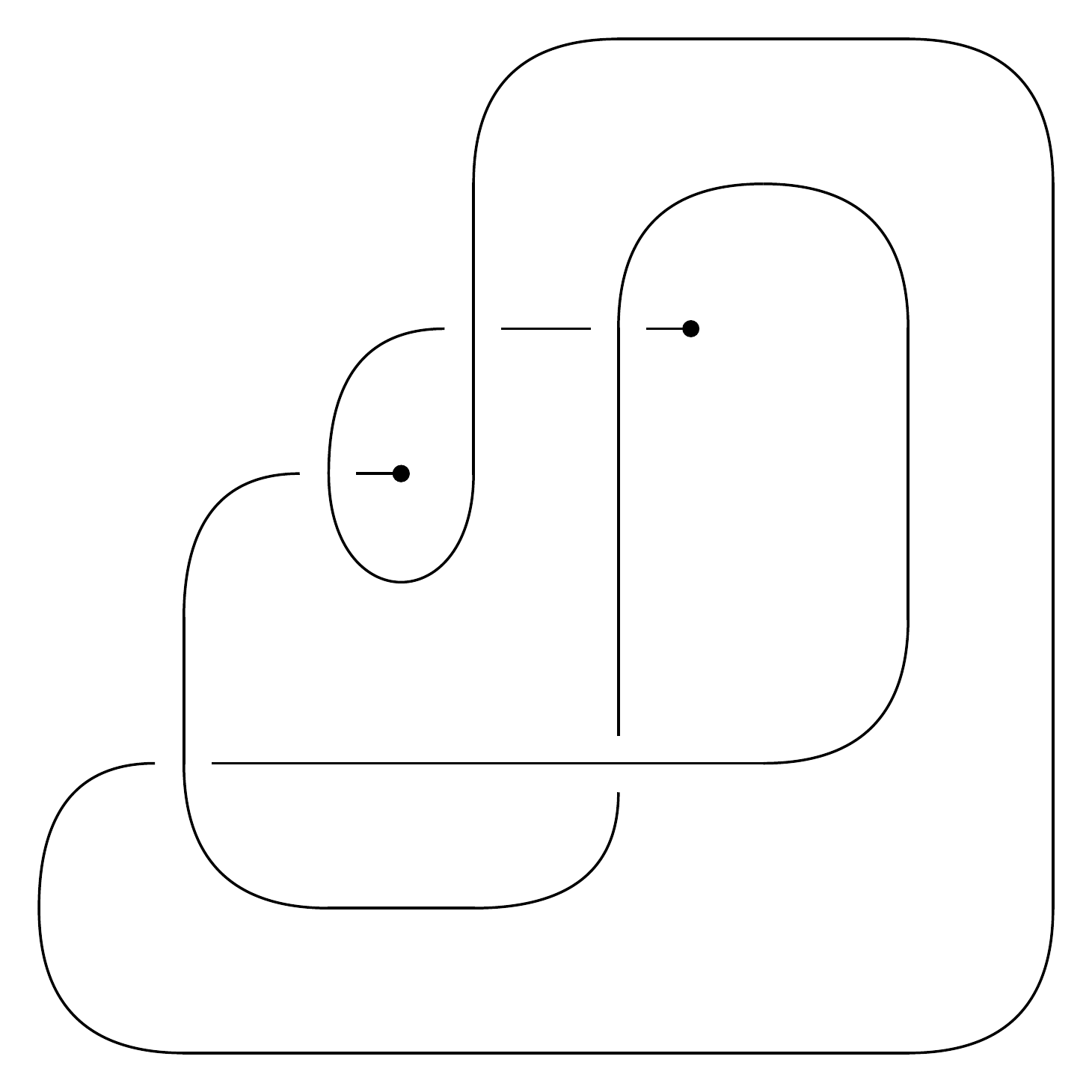}\\
\textcolor{black}{$5_{249}$}
\vspace{1cm}
\end{minipage}
\begin{minipage}[t]{.25\linewidth}
\centering
\includegraphics[width=0.9\textwidth,height=3.5cm,keepaspectratio]{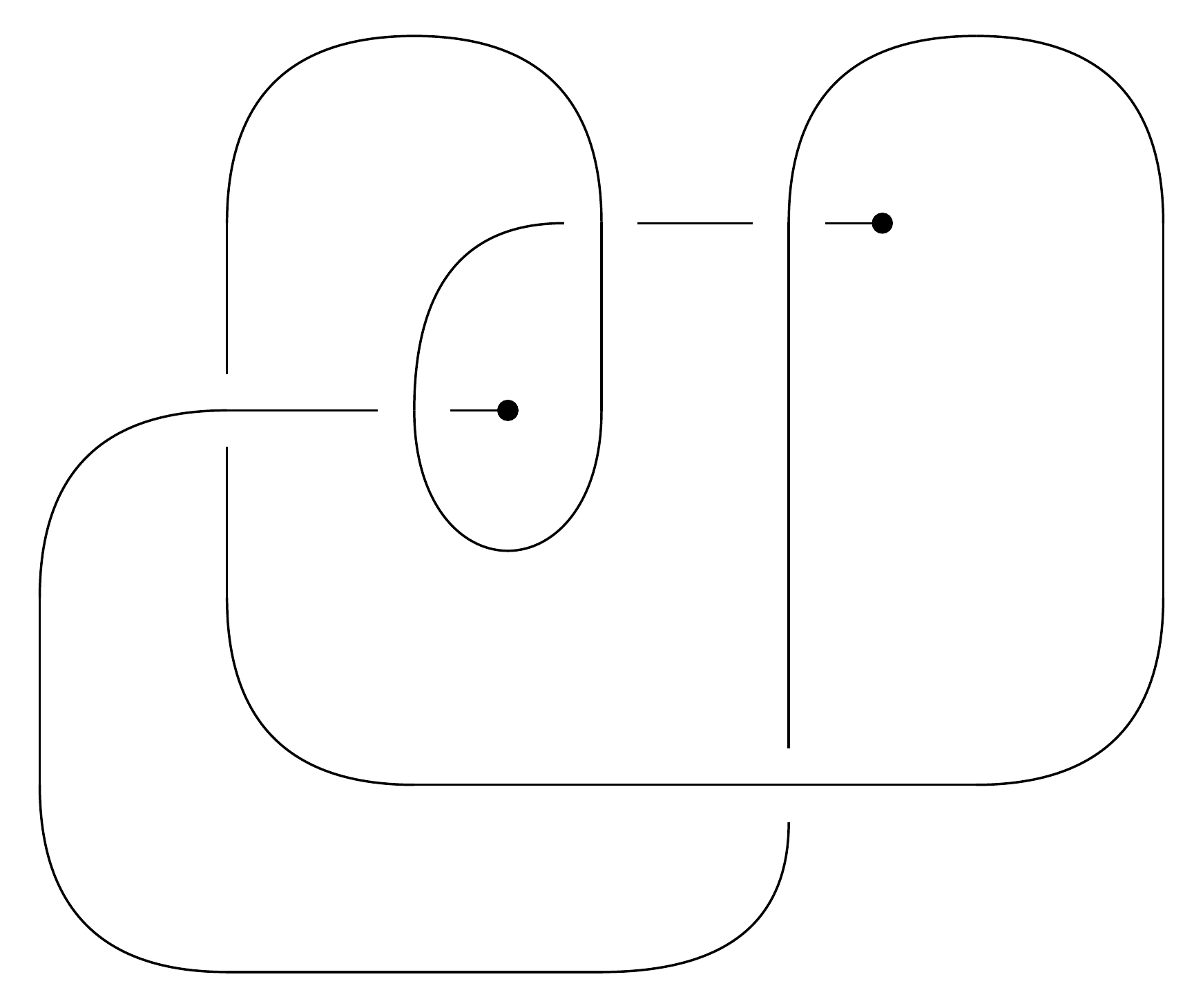}\\
\textcolor{black}{$5_{250}$}
\vspace{1cm}
\end{minipage}
\begin{minipage}[t]{.25\linewidth}
\centering
\includegraphics[width=0.9\textwidth,height=3.5cm,keepaspectratio]{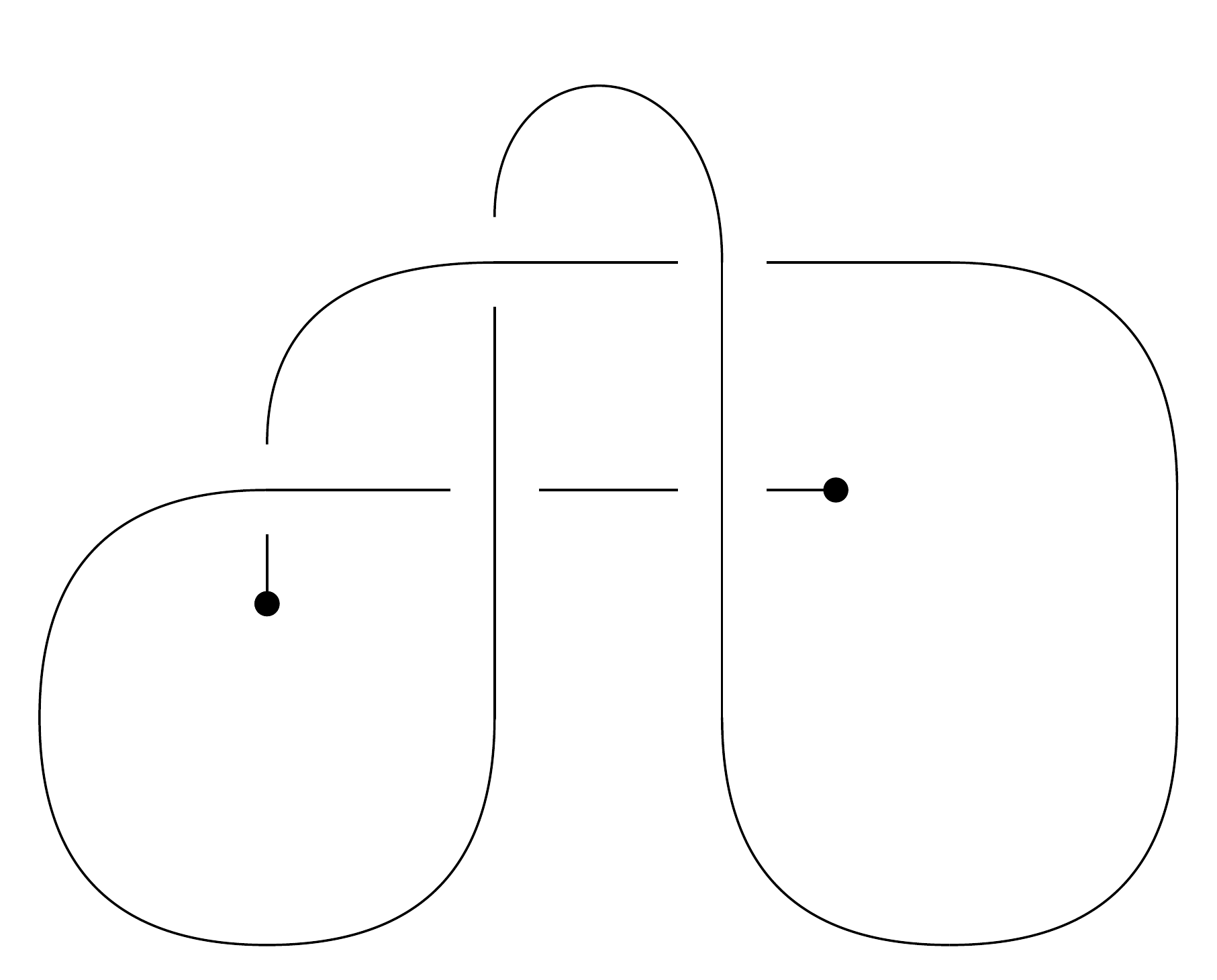}\\
\textcolor{black}{$5_{251}$}
\vspace{1cm}
\end{minipage}
\begin{minipage}[t]{.25\linewidth}
\centering
\includegraphics[width=0.9\textwidth,height=3.5cm,keepaspectratio]{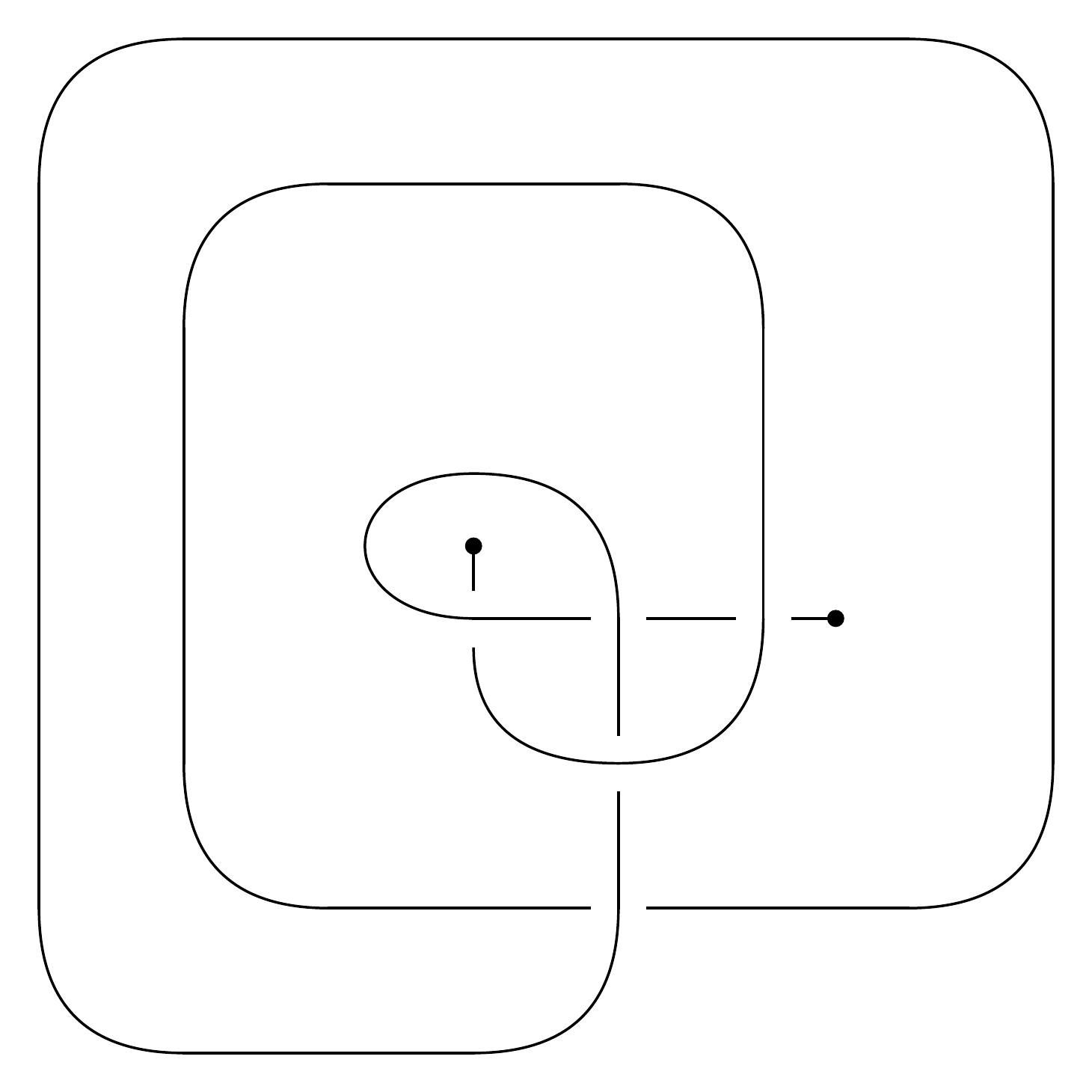}\\
\textcolor{black}{$5_{252}$}
\vspace{1cm}
\end{minipage}
\begin{minipage}[t]{.25\linewidth}
\centering
\includegraphics[width=0.9\textwidth,height=3.5cm,keepaspectratio]{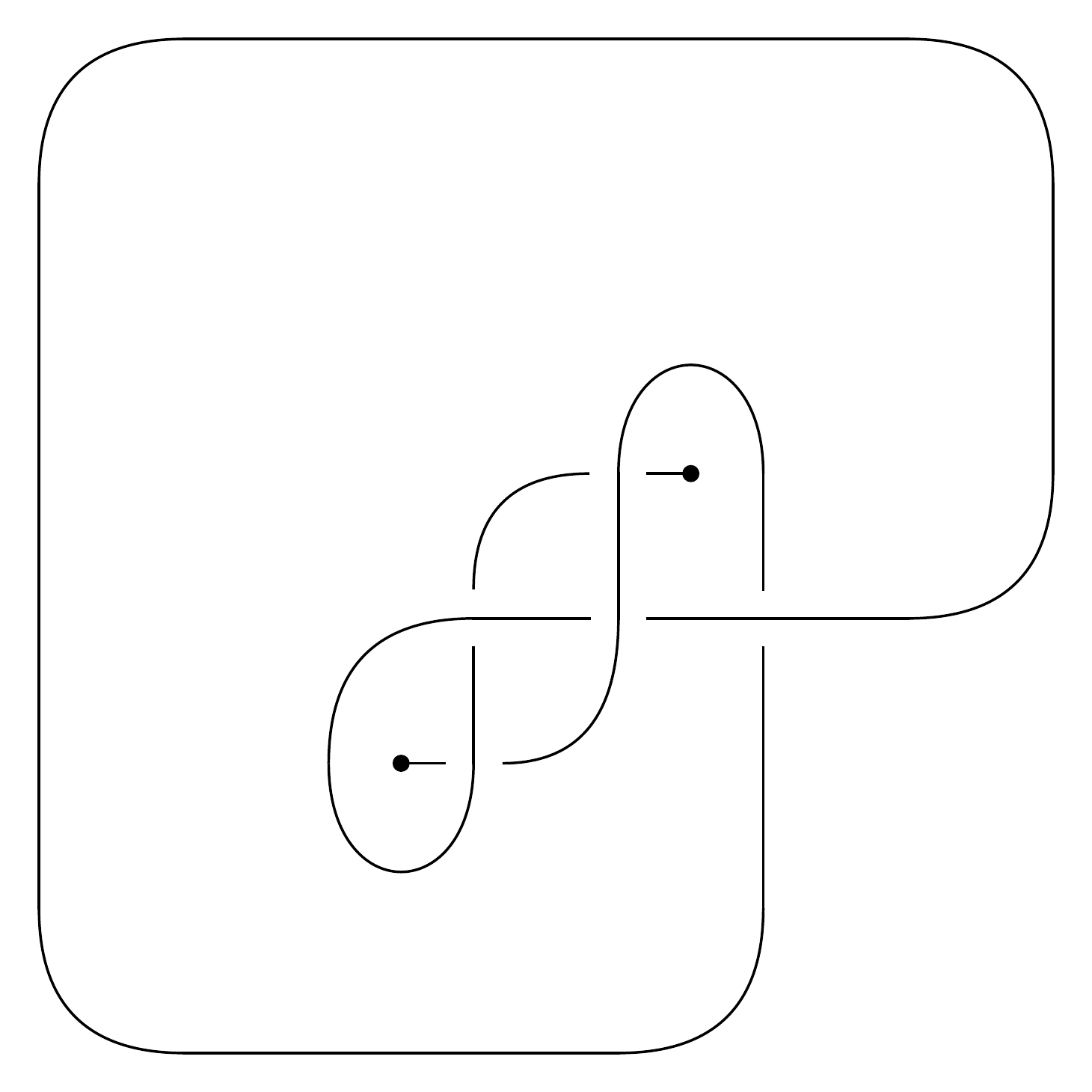}\\
\textcolor{black}{$5_{253}$}
\vspace{1cm}
\end{minipage}
\begin{minipage}[t]{.25\linewidth}
\centering
\includegraphics[width=0.9\textwidth,height=3.5cm,keepaspectratio]{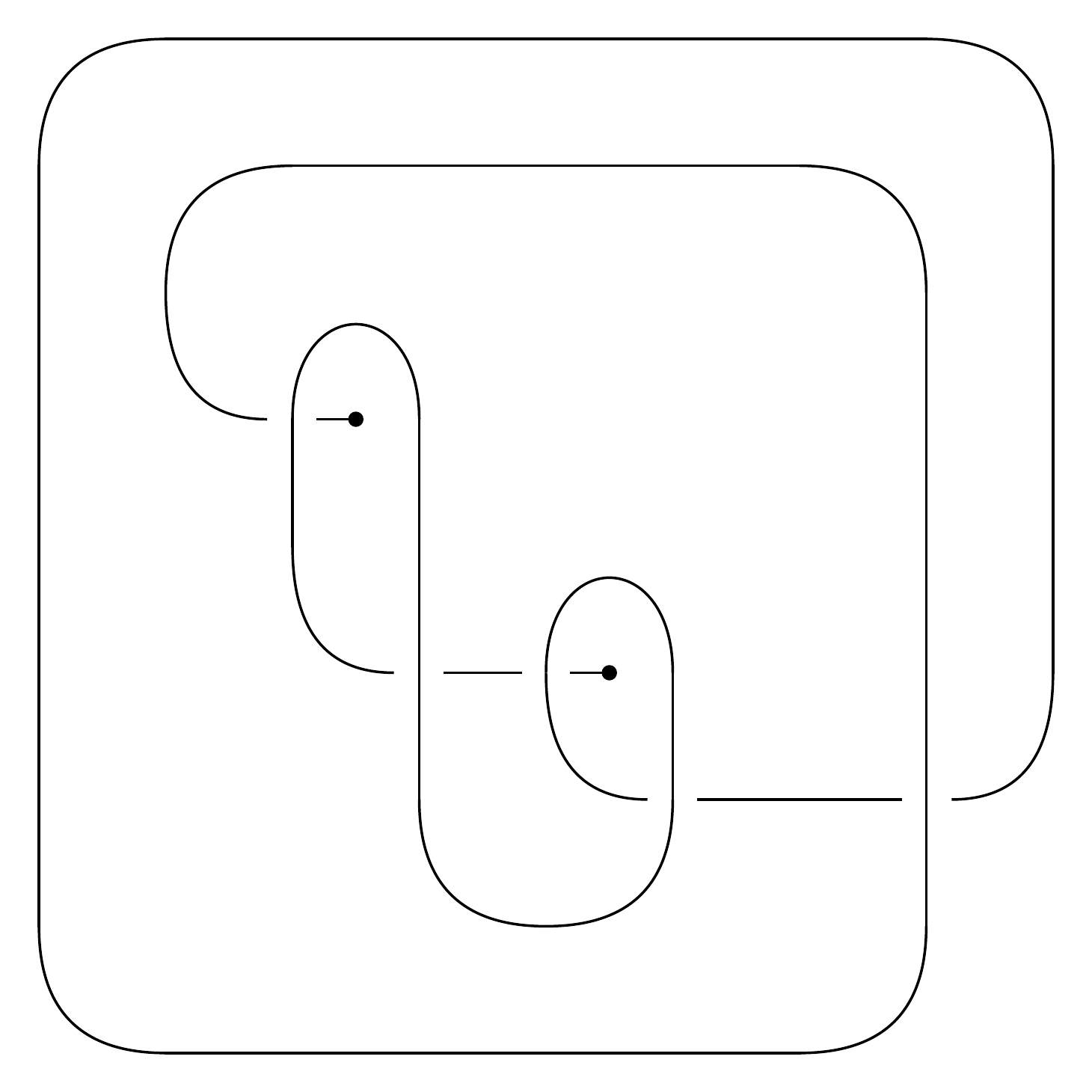}\\
\textcolor{black}{$5_{254}$}
\vspace{1cm}
\end{minipage}
\begin{minipage}[t]{.25\linewidth}
\centering
\includegraphics[width=0.9\textwidth,height=3.5cm,keepaspectratio]{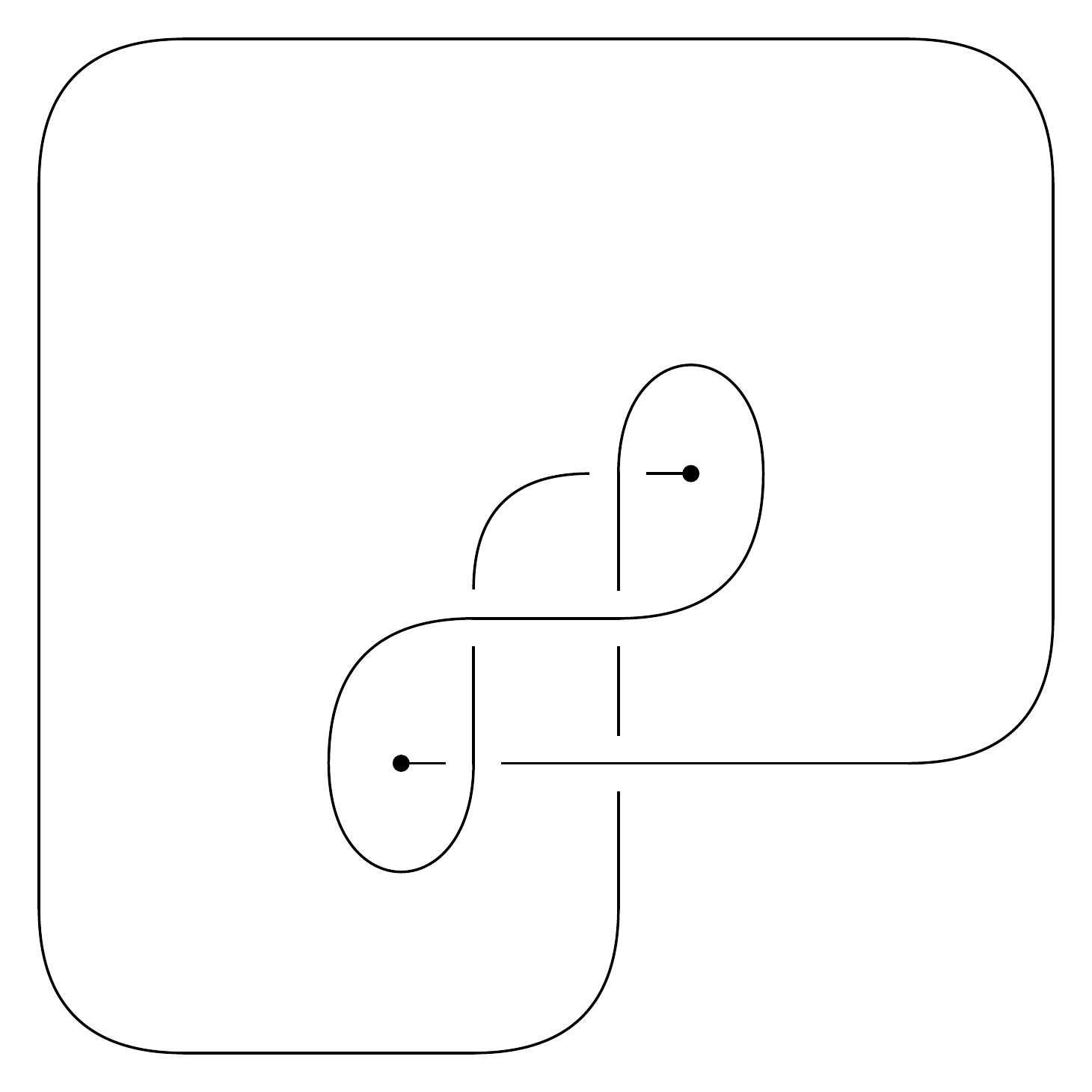}\\
\textcolor{black}{$5_{255}$}
\vspace{1cm}
\end{minipage}
\begin{minipage}[t]{.25\linewidth}
\centering
\includegraphics[width=0.9\textwidth,height=3.5cm,keepaspectratio]{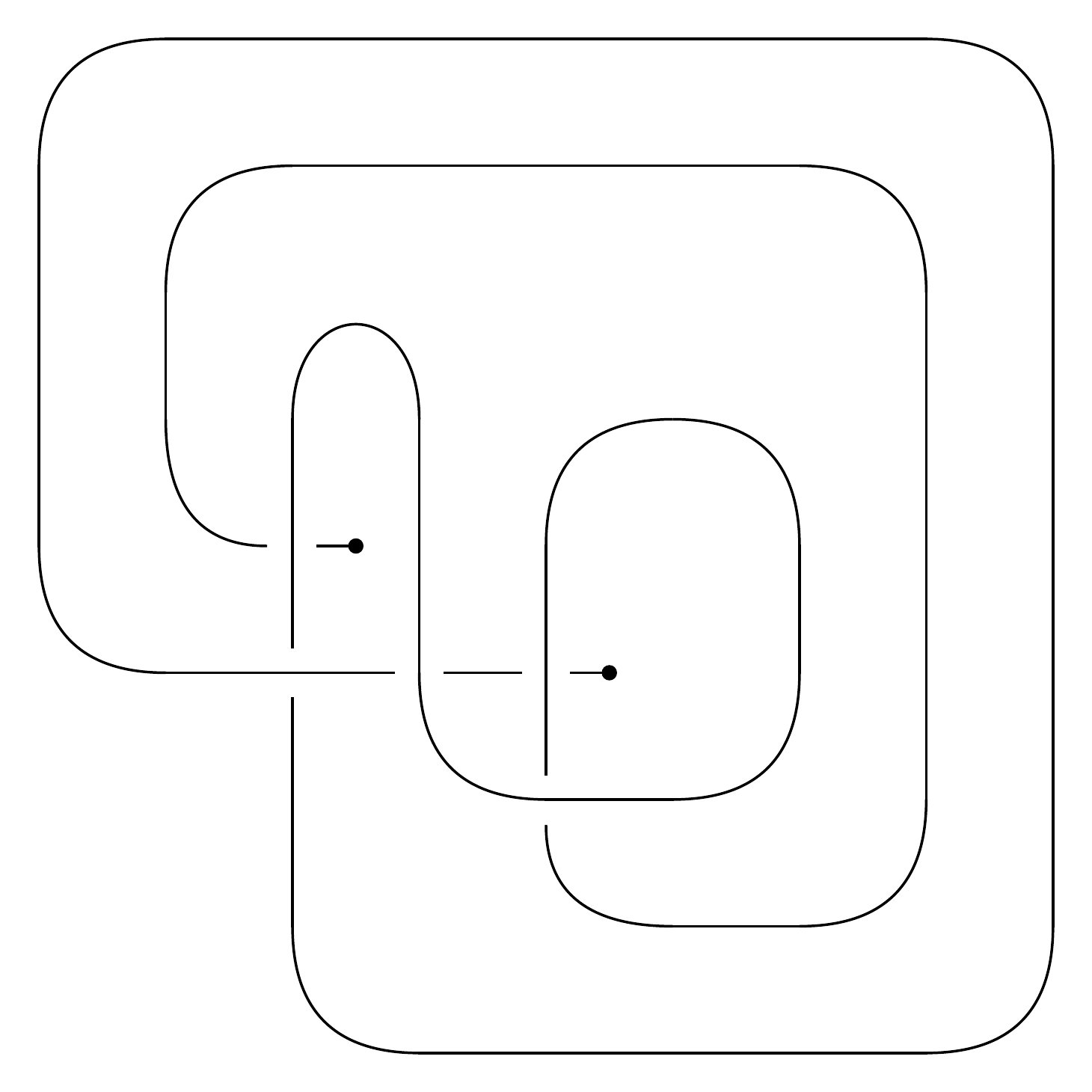}\\
\textcolor{black}{$5_{256}$}
\vspace{1cm}
\end{minipage}
\begin{minipage}[t]{.25\linewidth}
\centering
\includegraphics[width=0.9\textwidth,height=3.5cm,keepaspectratio]{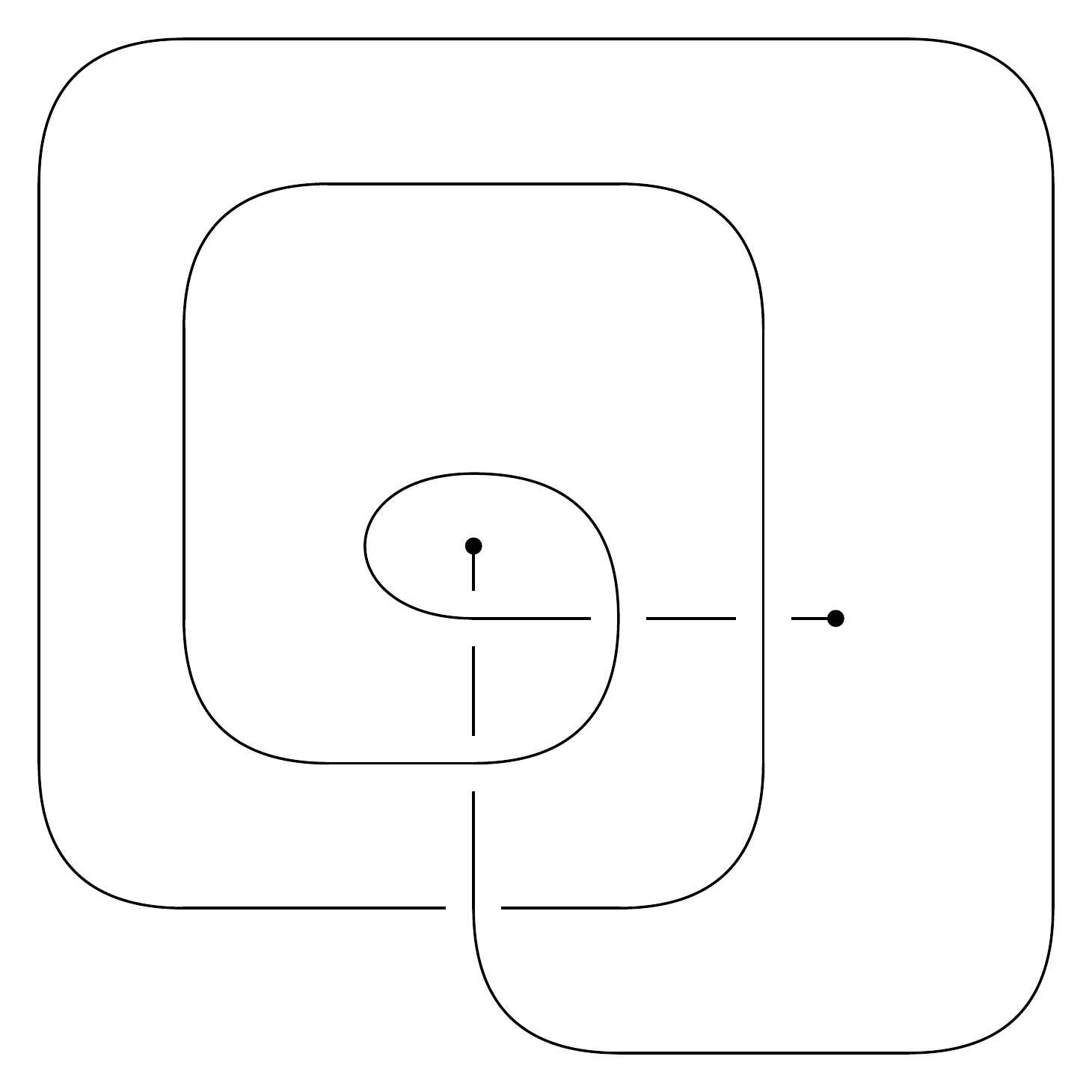}\\
\textcolor{black}{$5_{257}$}
\vspace{1cm}
\end{minipage}
\begin{minipage}[t]{.25\linewidth}
\centering
\includegraphics[width=0.9\textwidth,height=3.5cm,keepaspectratio]{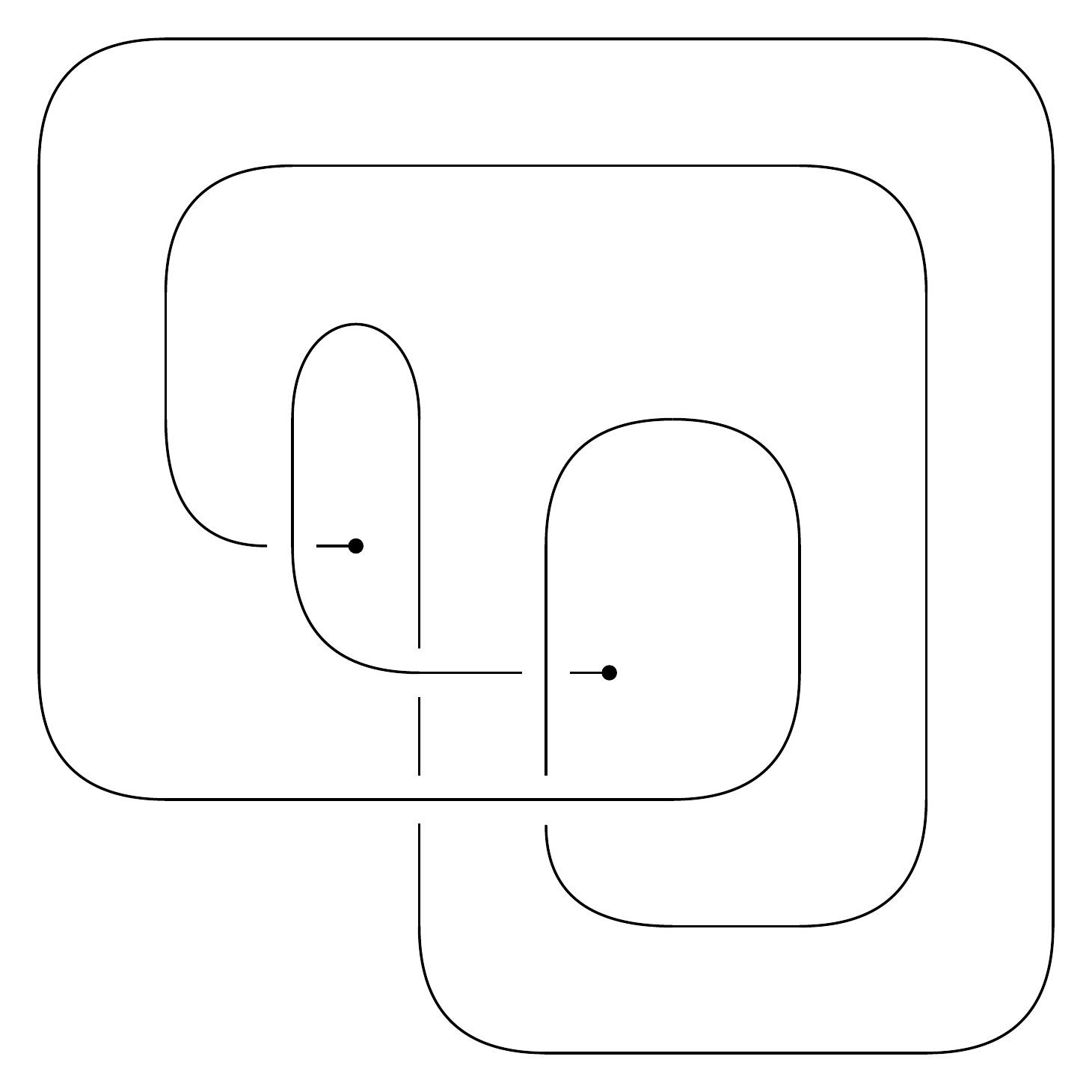}\\
\textcolor{black}{$5_{258}$}
\vspace{1cm}
\end{minipage}
\begin{minipage}[t]{.25\linewidth}
\centering
\includegraphics[width=0.9\textwidth,height=3.5cm,keepaspectratio]{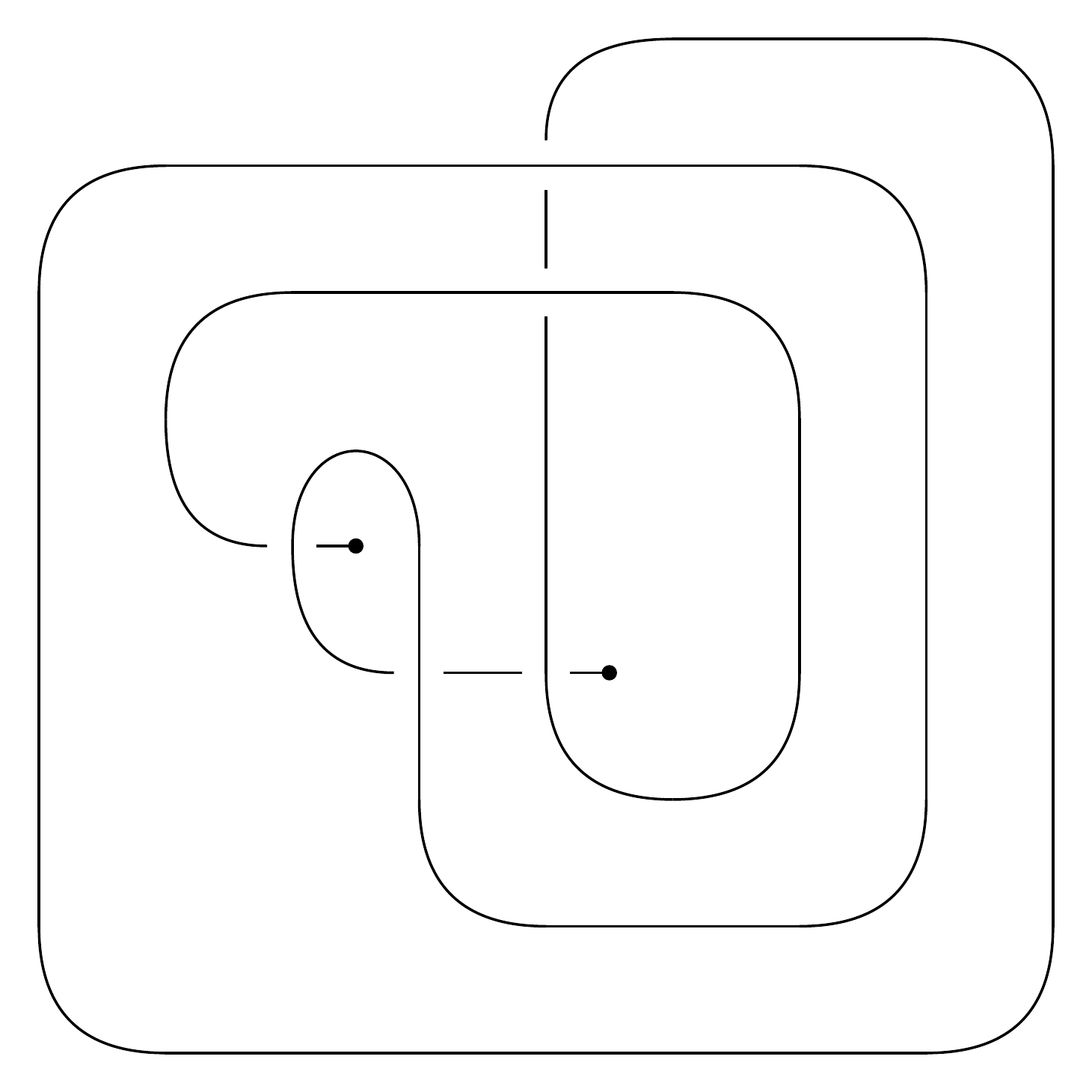}\\
\textcolor{black}{$5_{259}$}
\vspace{1cm}
\end{minipage}
\begin{minipage}[t]{.25\linewidth}
\centering
\includegraphics[width=0.9\textwidth,height=3.5cm,keepaspectratio]{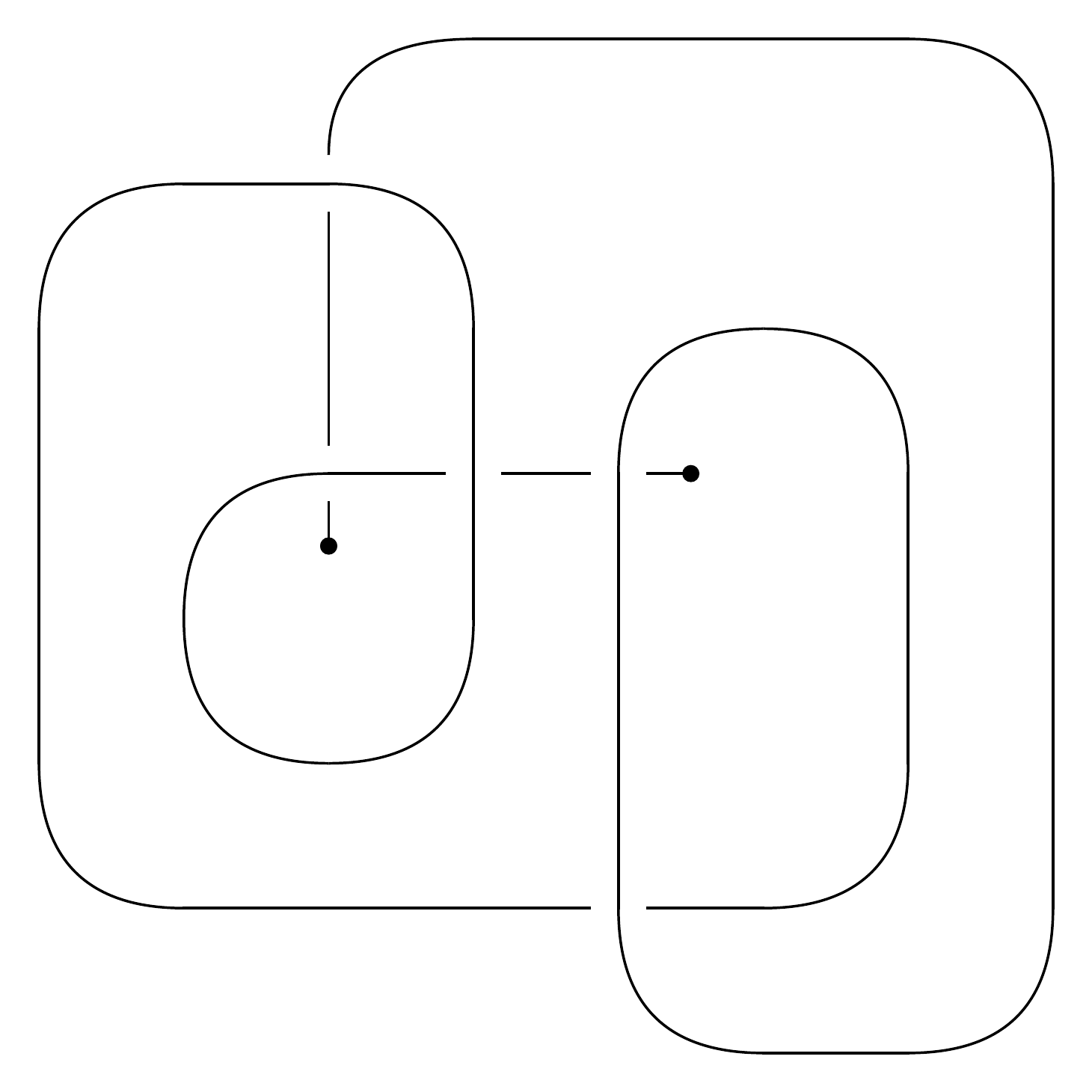}\\
\textcolor{black}{$5_{260}$}
\vspace{1cm}
\end{minipage}
\begin{minipage}[t]{.25\linewidth}
\centering
\includegraphics[width=0.9\textwidth,height=3.5cm,keepaspectratio]{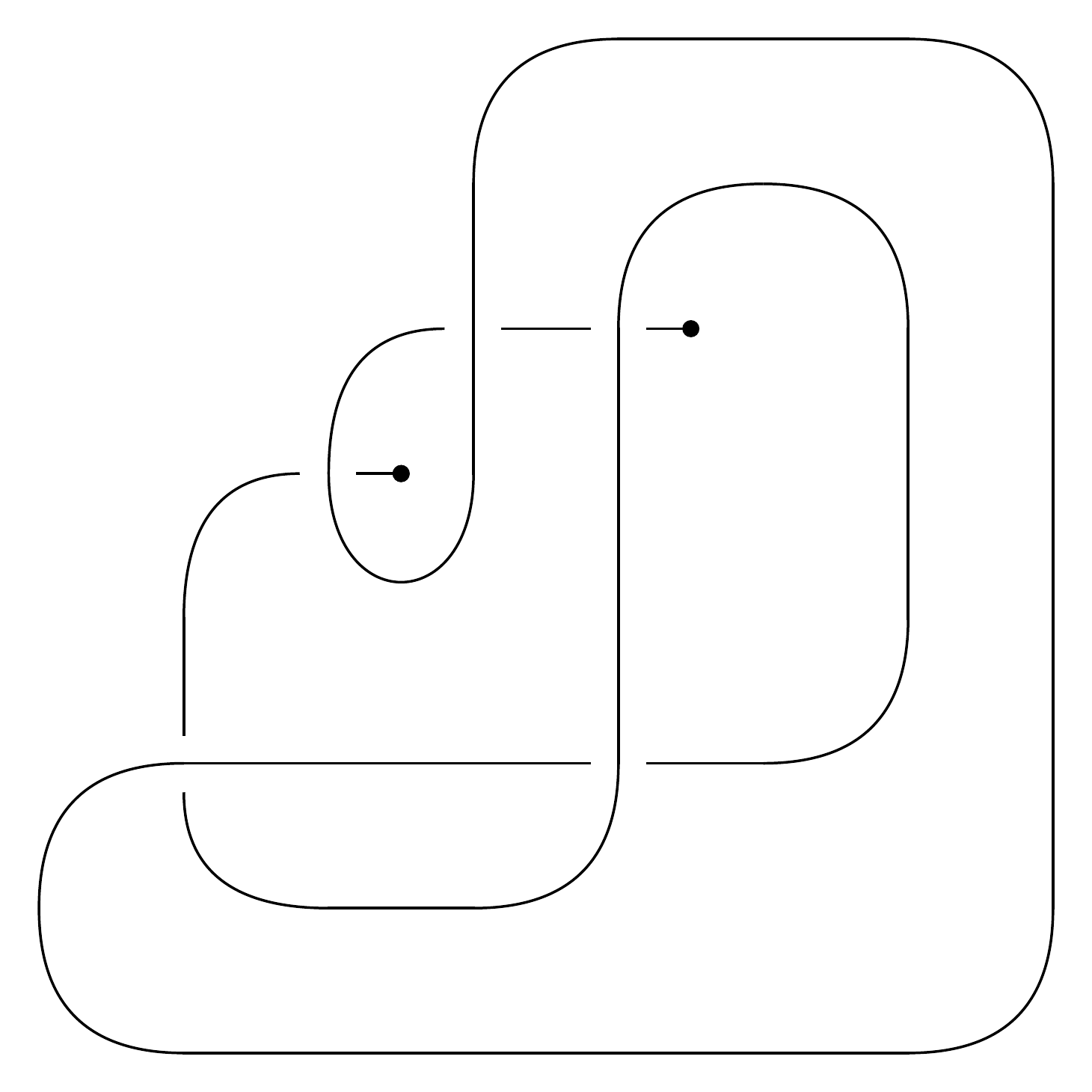}\\
\textcolor{black}{$5_{261}$}
\vspace{1cm}
\end{minipage}
\begin{minipage}[t]{.25\linewidth}
\centering
\includegraphics[width=0.9\textwidth,height=3.5cm,keepaspectratio]{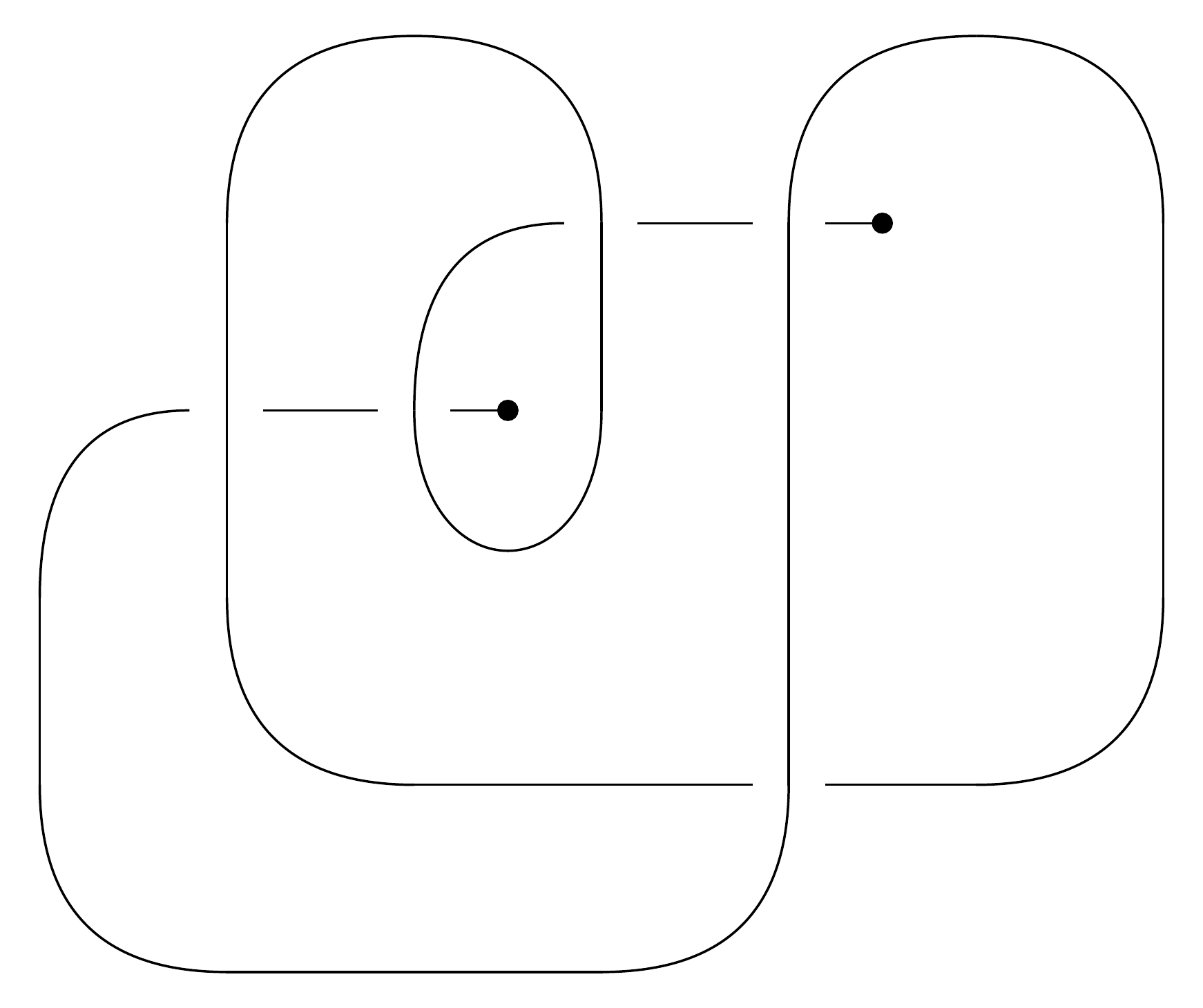}\\
\textcolor{black}{$5_{262}$}
\vspace{1cm}
\end{minipage}
\begin{minipage}[t]{.25\linewidth}
\centering
\includegraphics[width=0.9\textwidth,height=3.5cm,keepaspectratio]{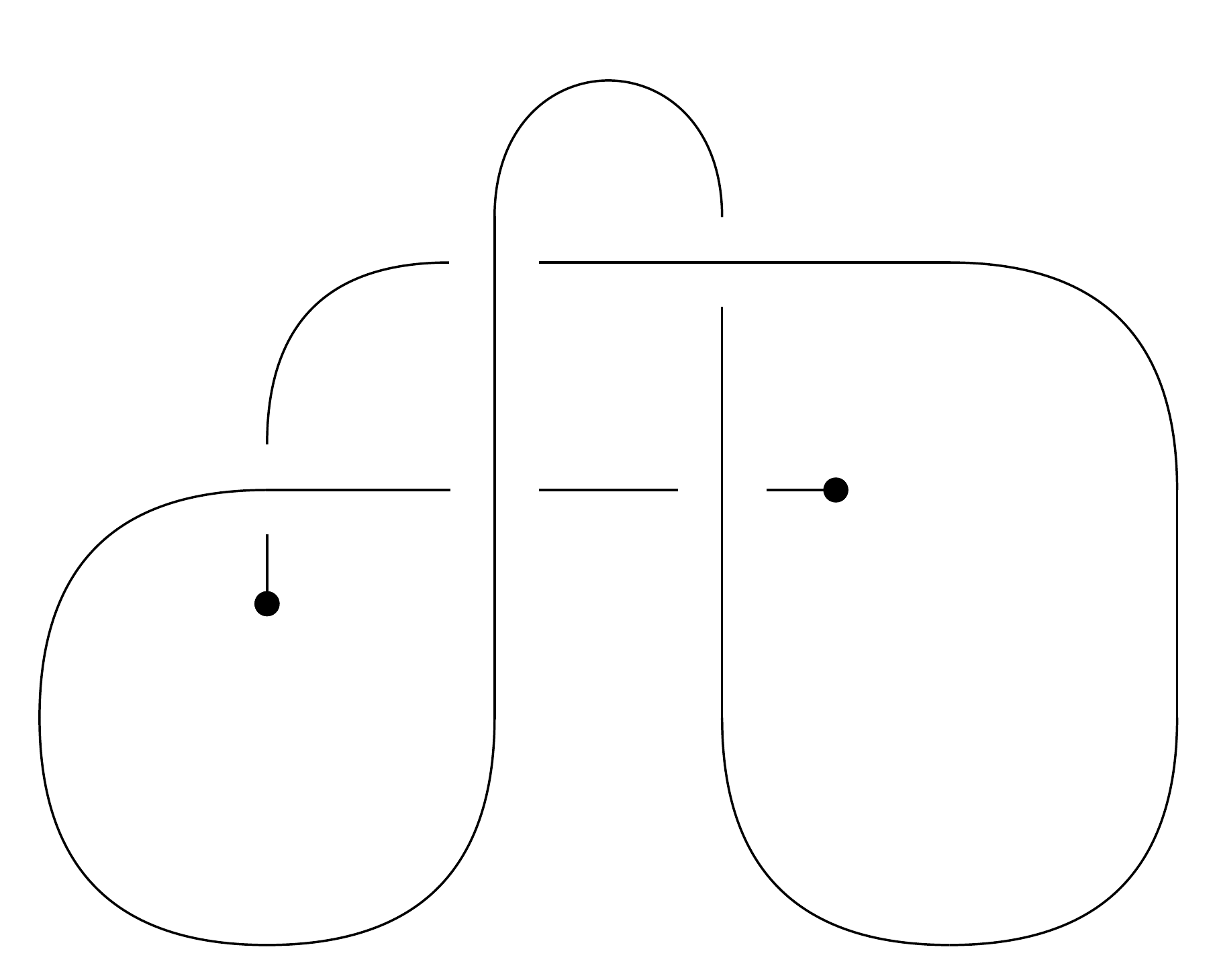}\\
\textcolor{black}{$5_{263}$}
\vspace{1cm}
\end{minipage}
\begin{minipage}[t]{.25\linewidth}
\centering
\includegraphics[width=0.9\textwidth,height=3.5cm,keepaspectratio]{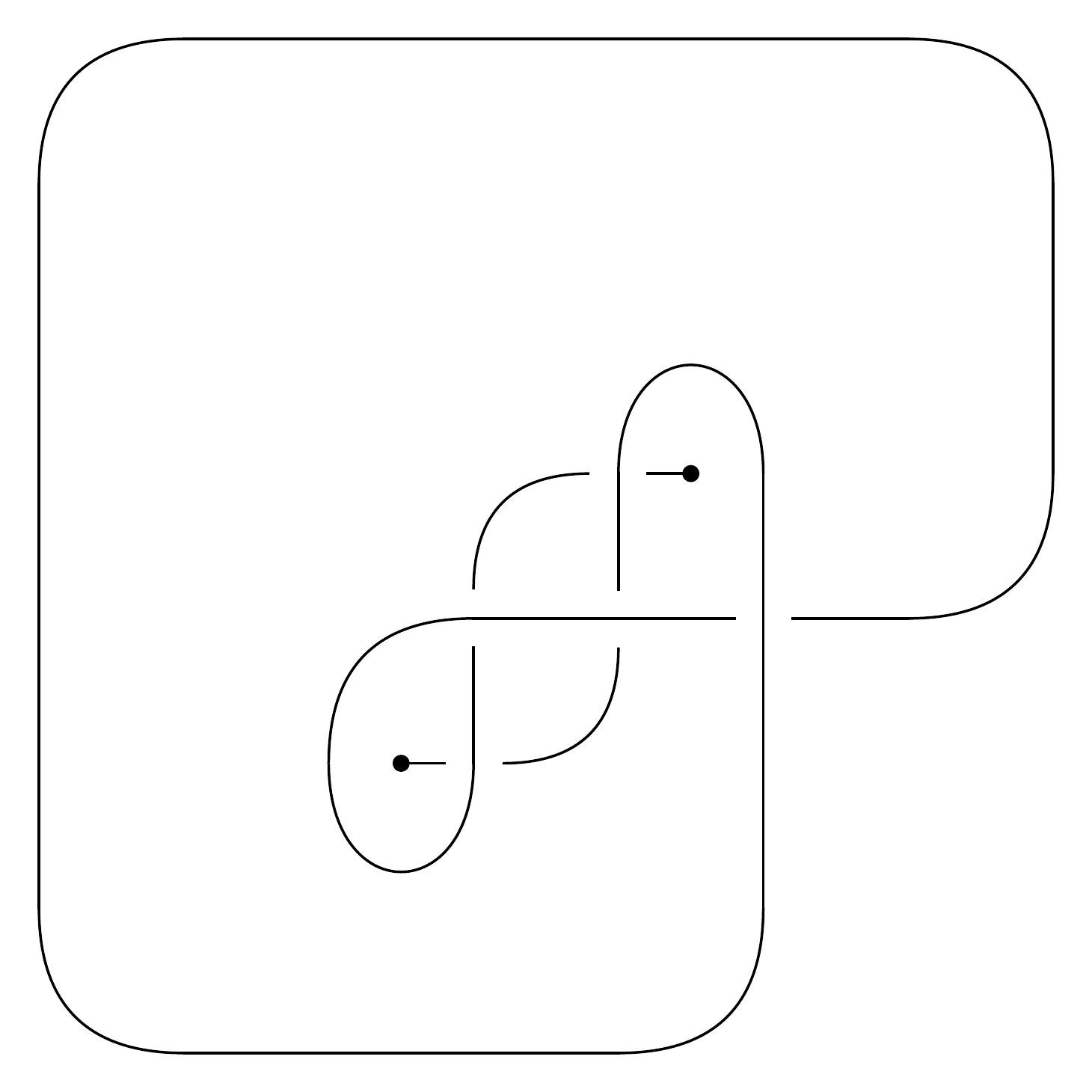}\\
\textcolor{black}{$5_{264}$}
\vspace{1cm}
\end{minipage}
\begin{minipage}[t]{.25\linewidth}
\centering
\includegraphics[width=0.9\textwidth,height=3.5cm,keepaspectratio]{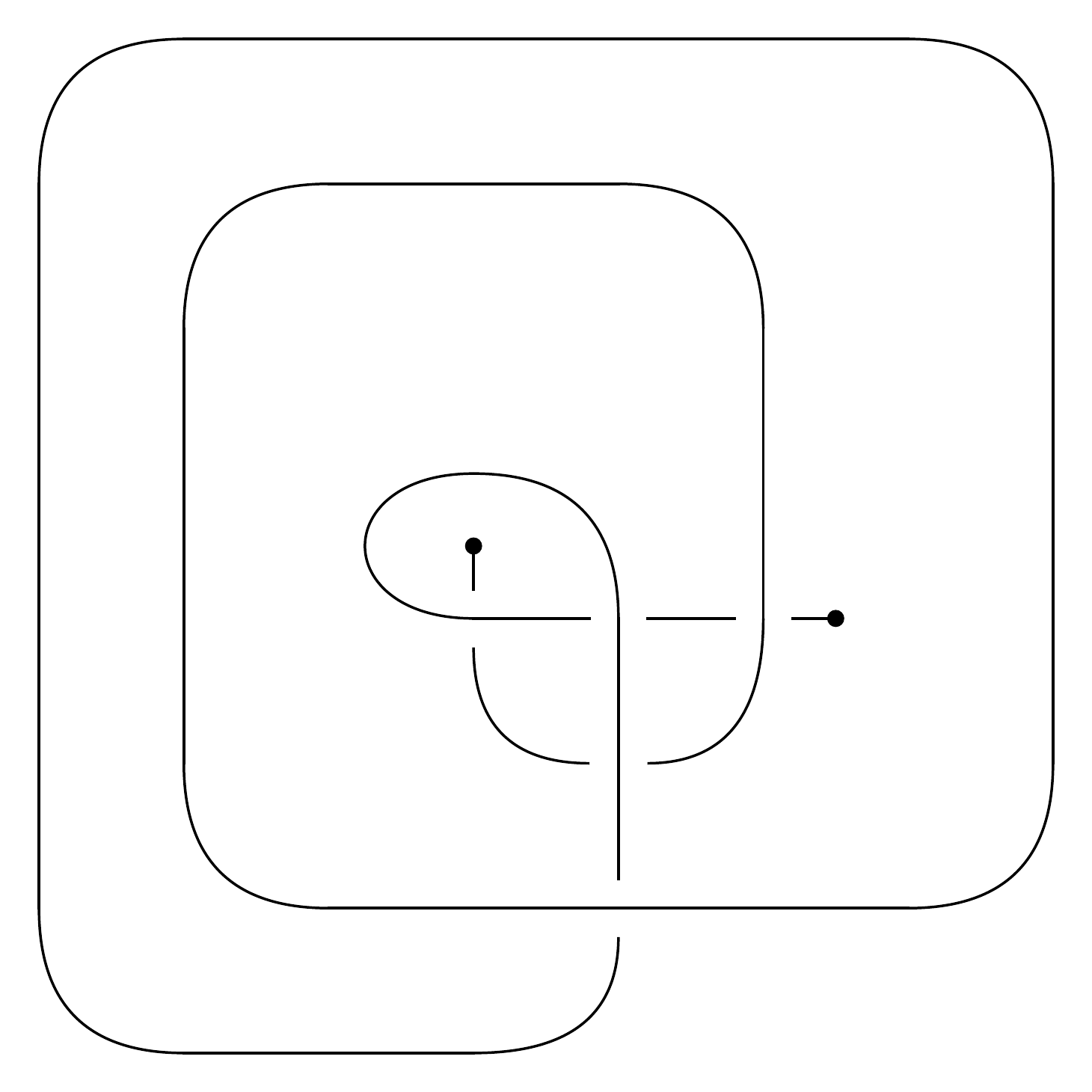}\\
\textcolor{black}{$5_{265}$}
\vspace{1cm}
\end{minipage}
\begin{minipage}[t]{.25\linewidth}
\centering
\includegraphics[width=0.9\textwidth,height=3.5cm,keepaspectratio]{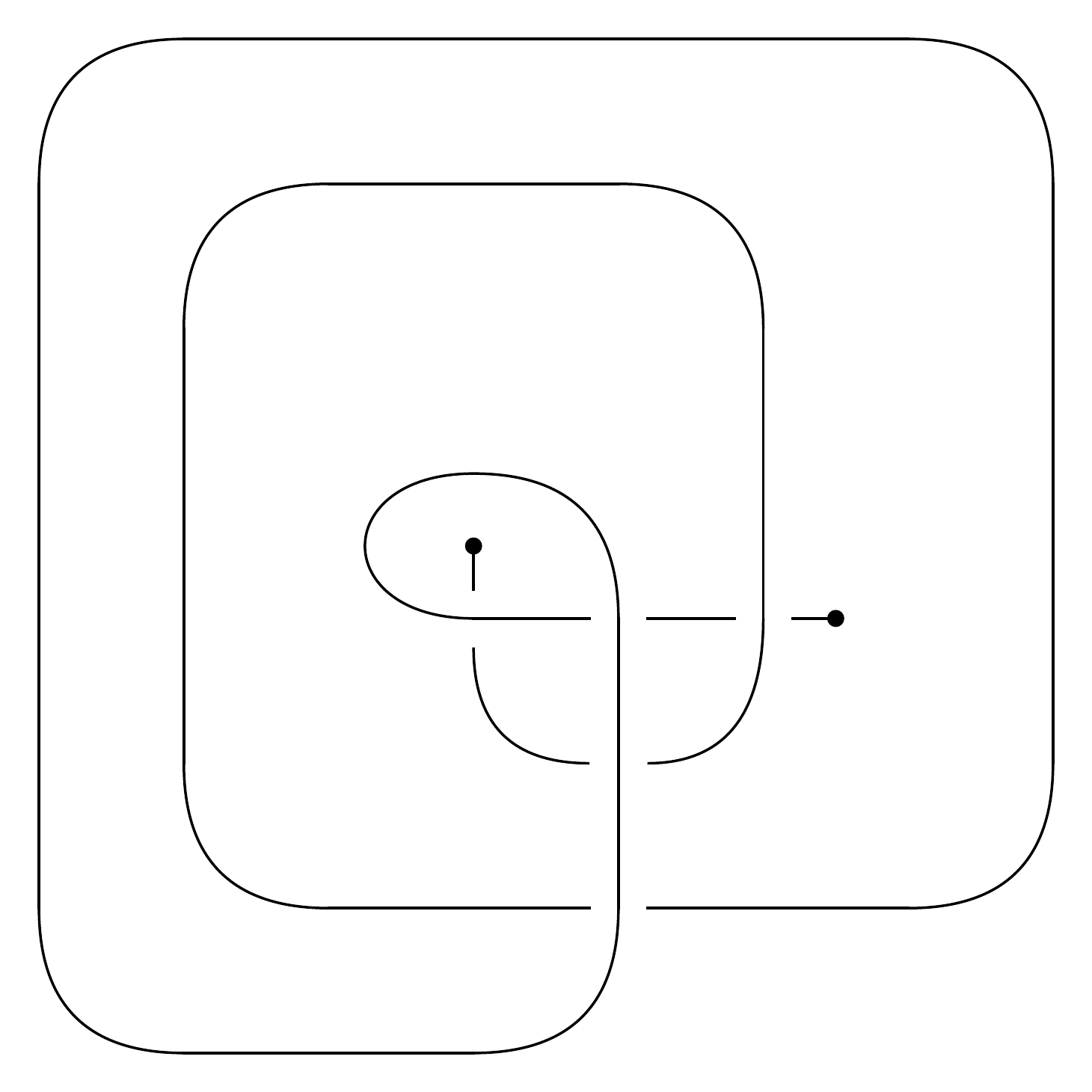}\\
\textcolor{black}{$5_{266}$}
\vspace{1cm}
\end{minipage}
\begin{minipage}[t]{.25\linewidth}
\centering
\includegraphics[width=0.9\textwidth,height=3.5cm,keepaspectratio]{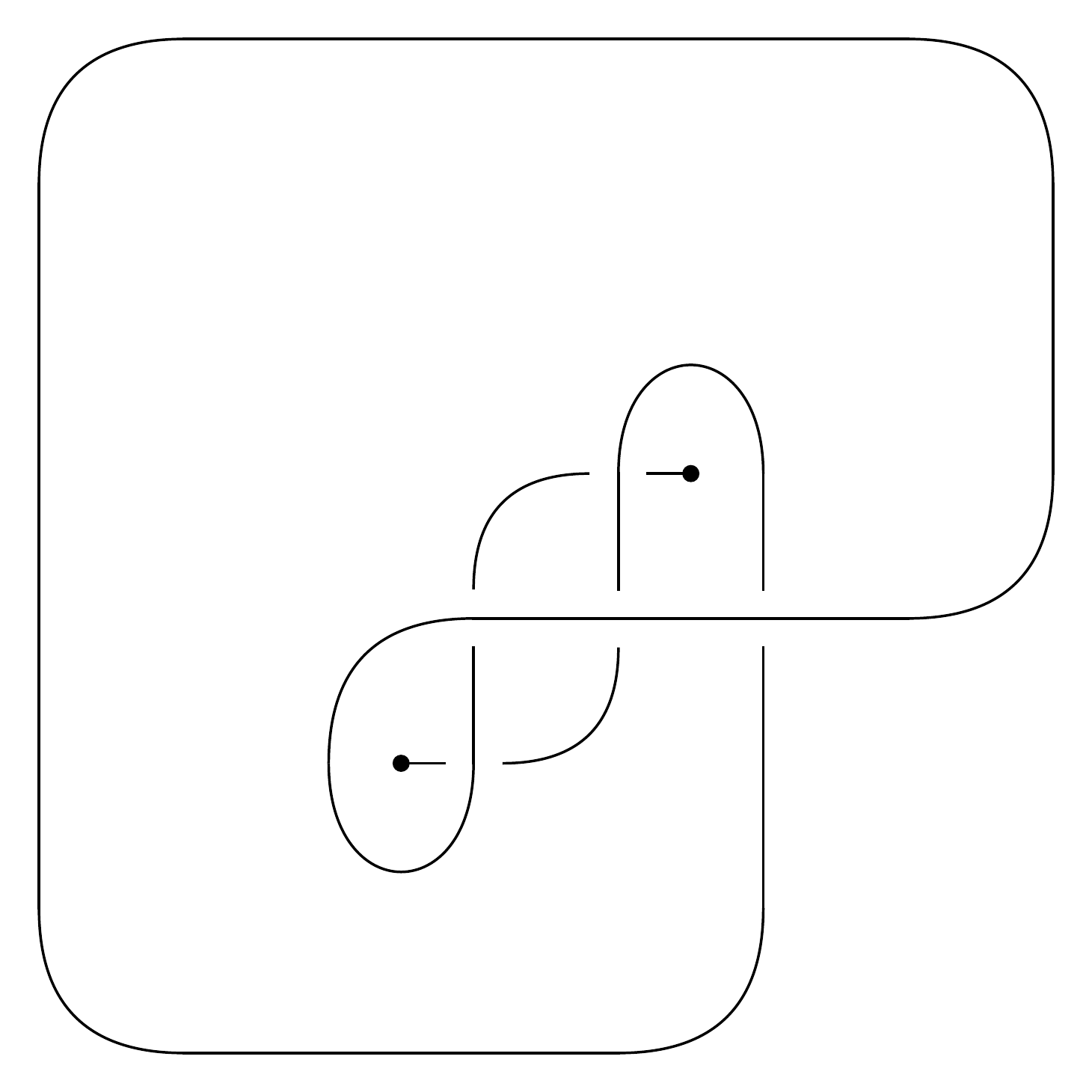}\\
\textcolor{black}{$5_{267}$}
\vspace{1cm}
\end{minipage}
\begin{minipage}[t]{.25\linewidth}
\centering
\includegraphics[width=0.9\textwidth,height=3.5cm,keepaspectratio]{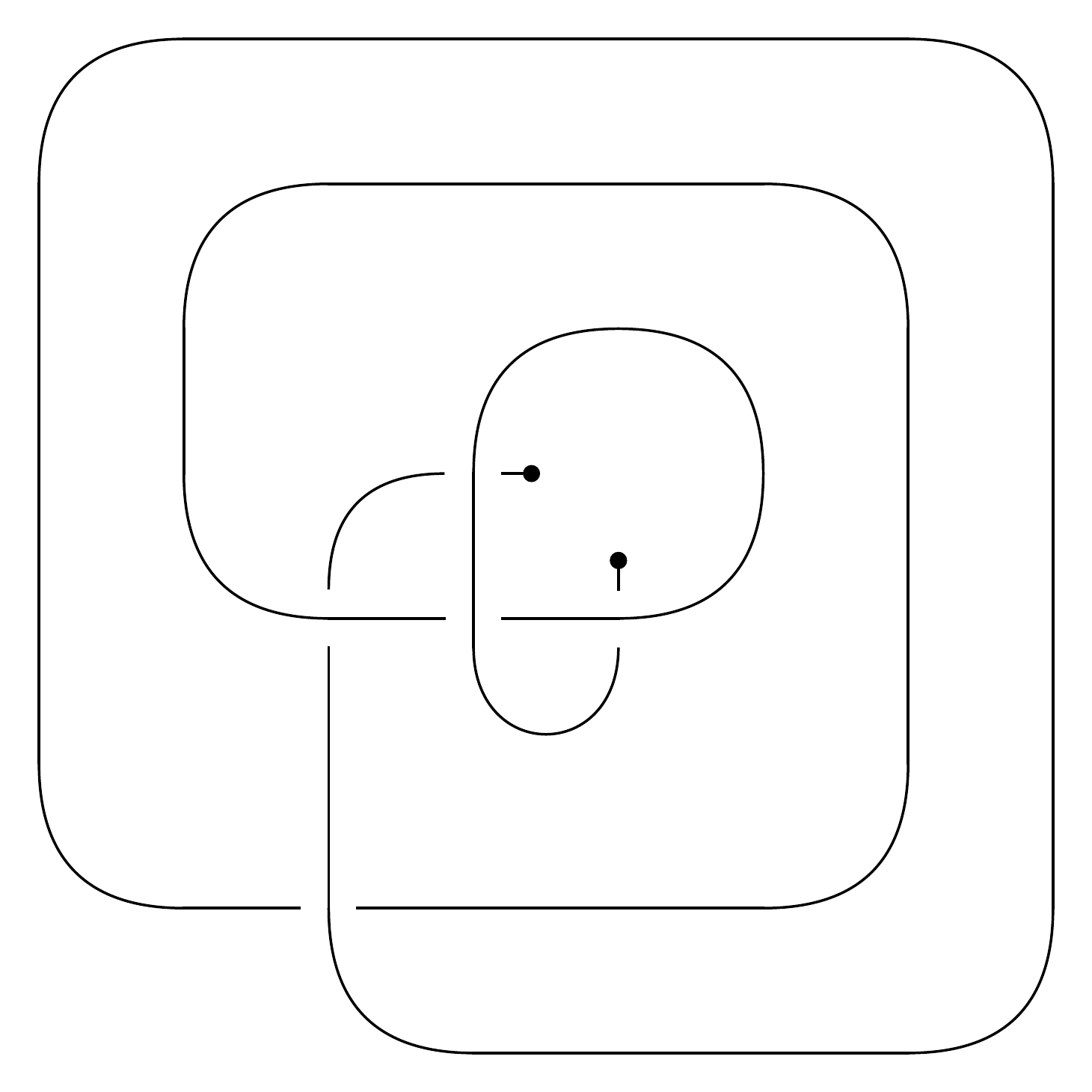}\\
\textcolor{black}{$5_{268}$}
\vspace{1cm}
\end{minipage}
\begin{minipage}[t]{.25\linewidth}
\centering
\includegraphics[width=0.9\textwidth,height=3.5cm,keepaspectratio]{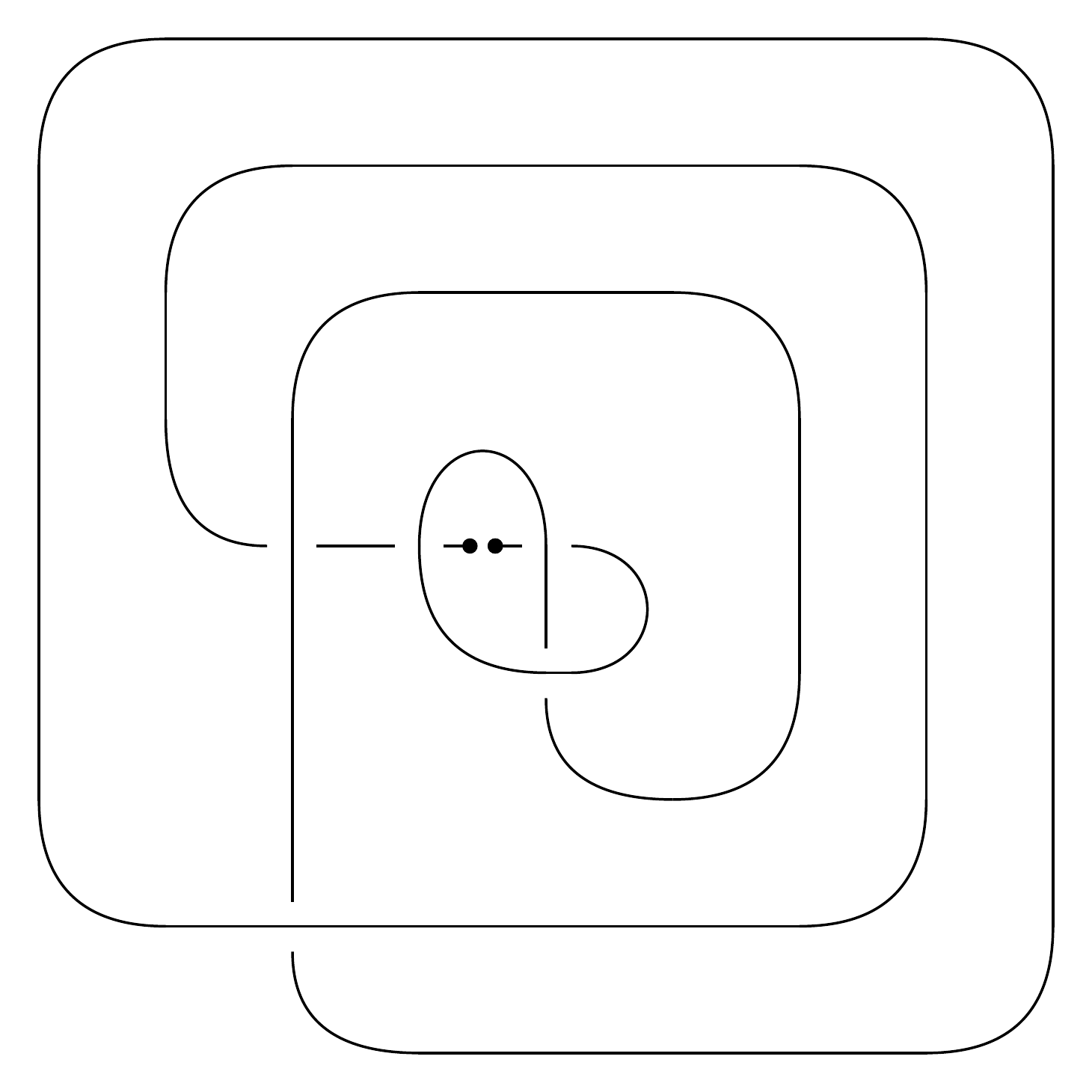}\\
\textcolor{black}{$5_{269}$}
\vspace{1cm}
\end{minipage}
\begin{minipage}[t]{.25\linewidth}
\centering
\includegraphics[width=0.9\textwidth,height=3.5cm,keepaspectratio]{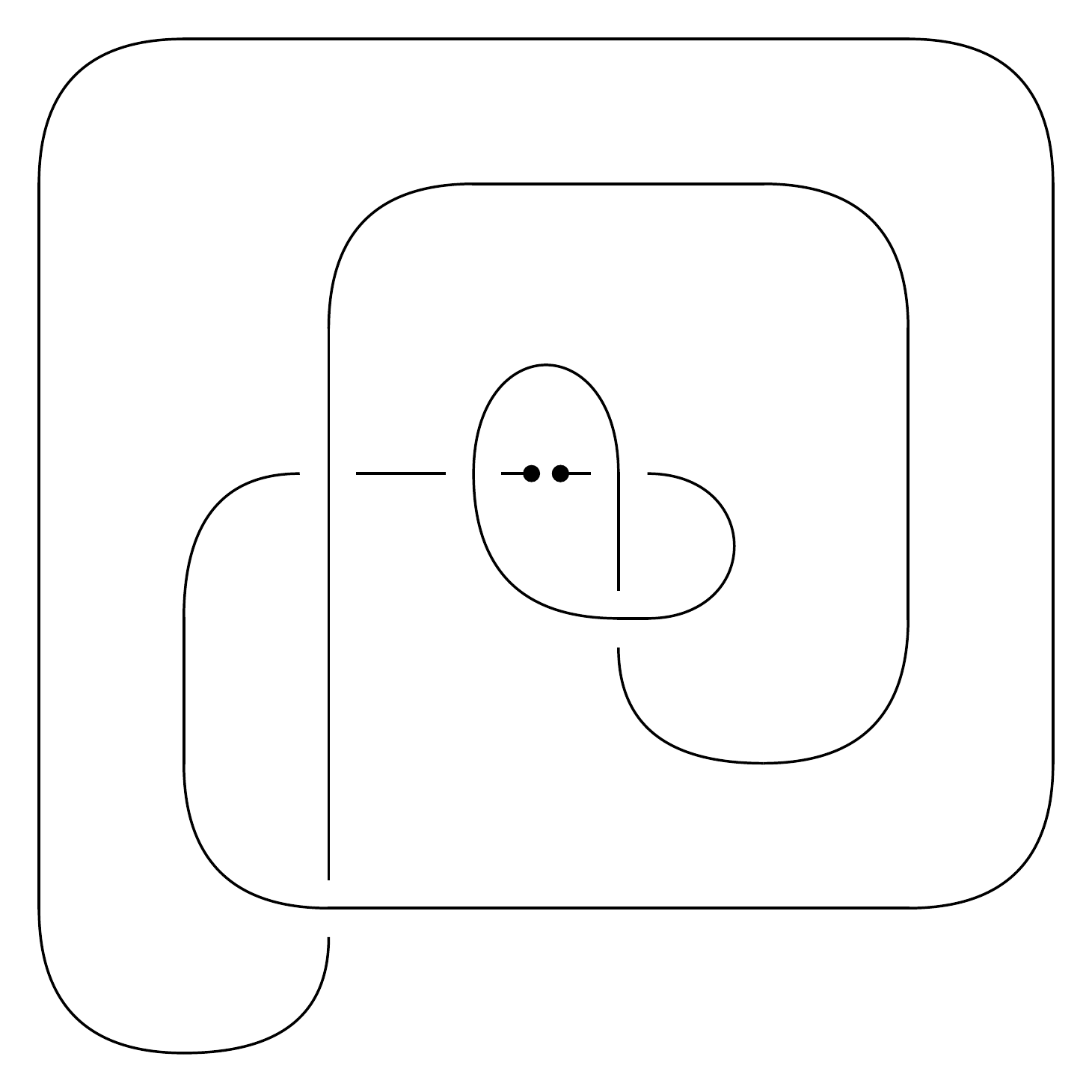}\\
\textcolor{black}{$5_{270}$}
\vspace{1cm}
\end{minipage}
\begin{minipage}[t]{.25\linewidth}
\centering
\includegraphics[width=0.9\textwidth,height=3.5cm,keepaspectratio]{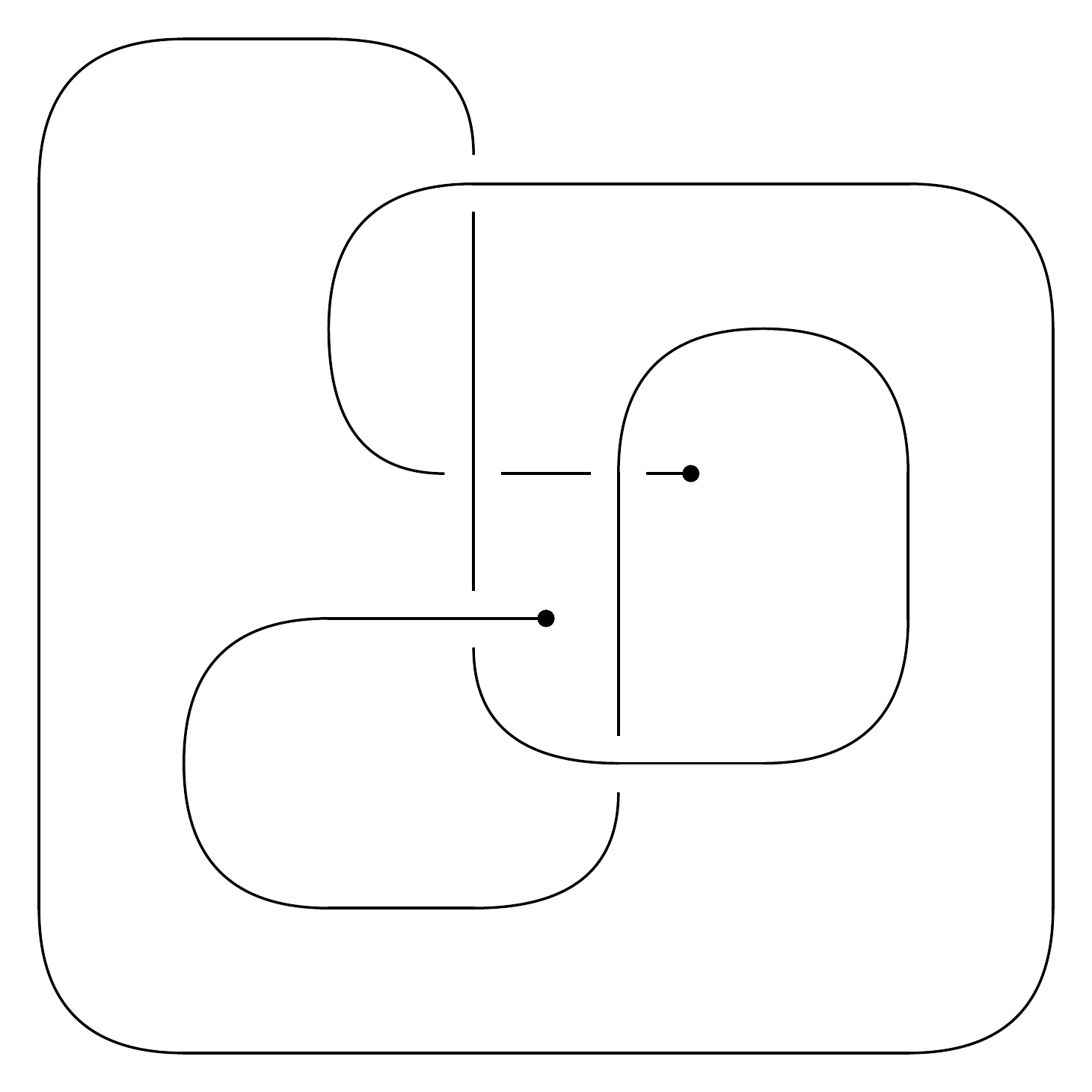}\\
\textcolor{black}{$5_{271}$}
\vspace{1cm}
\end{minipage}
\begin{minipage}[t]{.25\linewidth}
\centering
\includegraphics[width=0.9\textwidth,height=3.5cm,keepaspectratio]{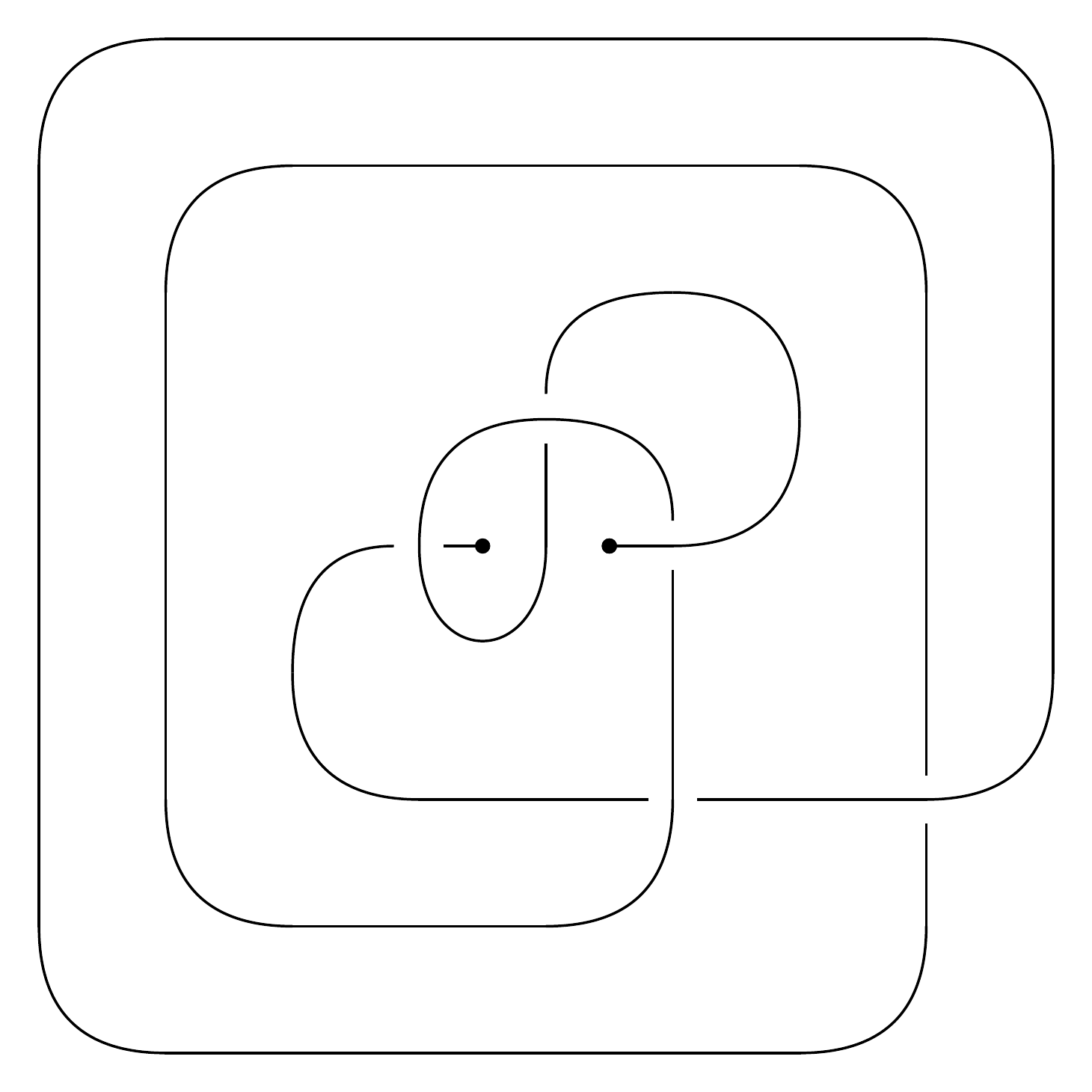}\\
\textcolor{black}{$5_{272}$}
\vspace{1cm}
\end{minipage}
\begin{minipage}[t]{.25\linewidth}
\centering
\includegraphics[width=0.9\textwidth,height=3.5cm,keepaspectratio]{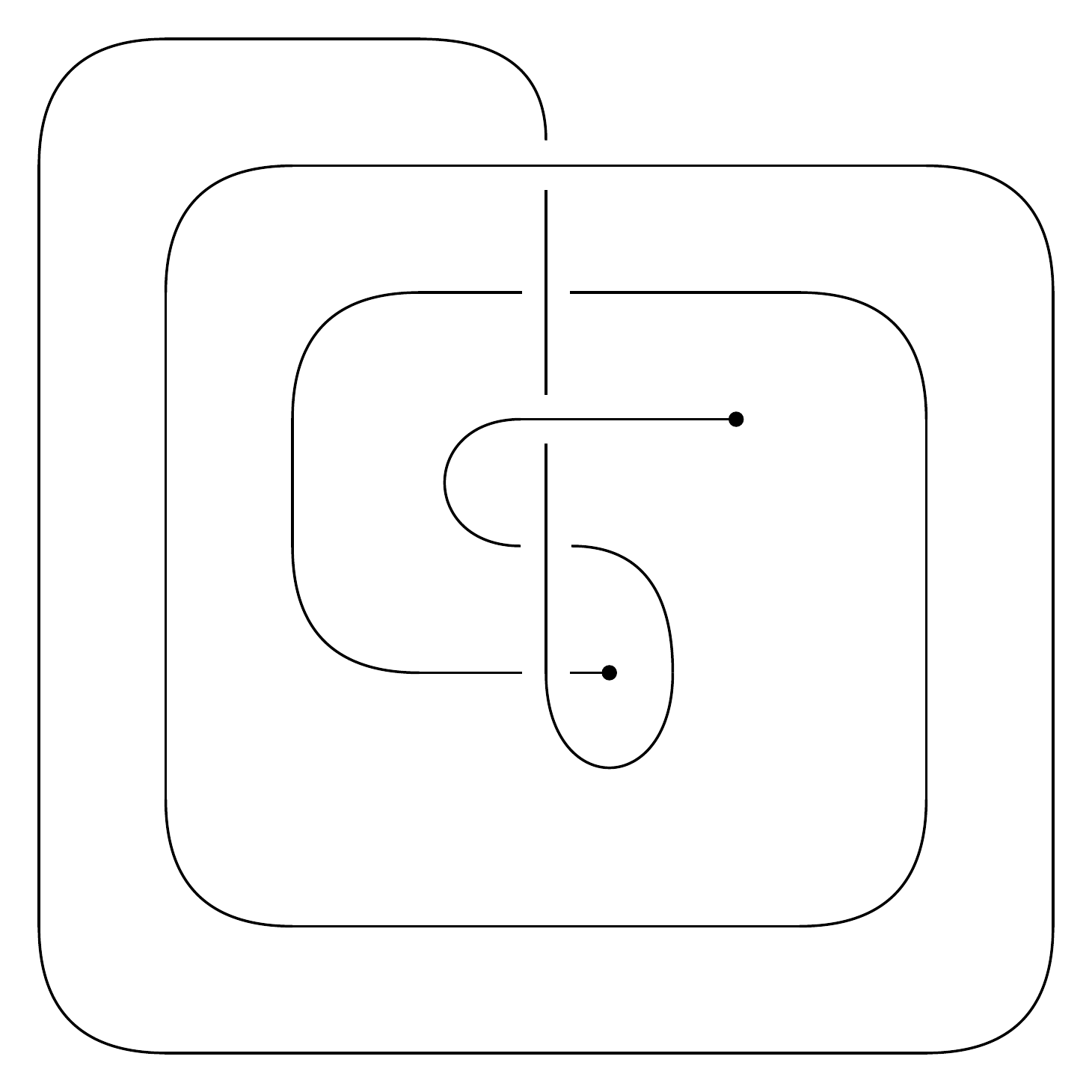}\\
\textcolor{black}{$5_{273}$}
\vspace{1cm}
\end{minipage}
\begin{minipage}[t]{.25\linewidth}
\centering
\includegraphics[width=0.9\textwidth,height=3.5cm,keepaspectratio]{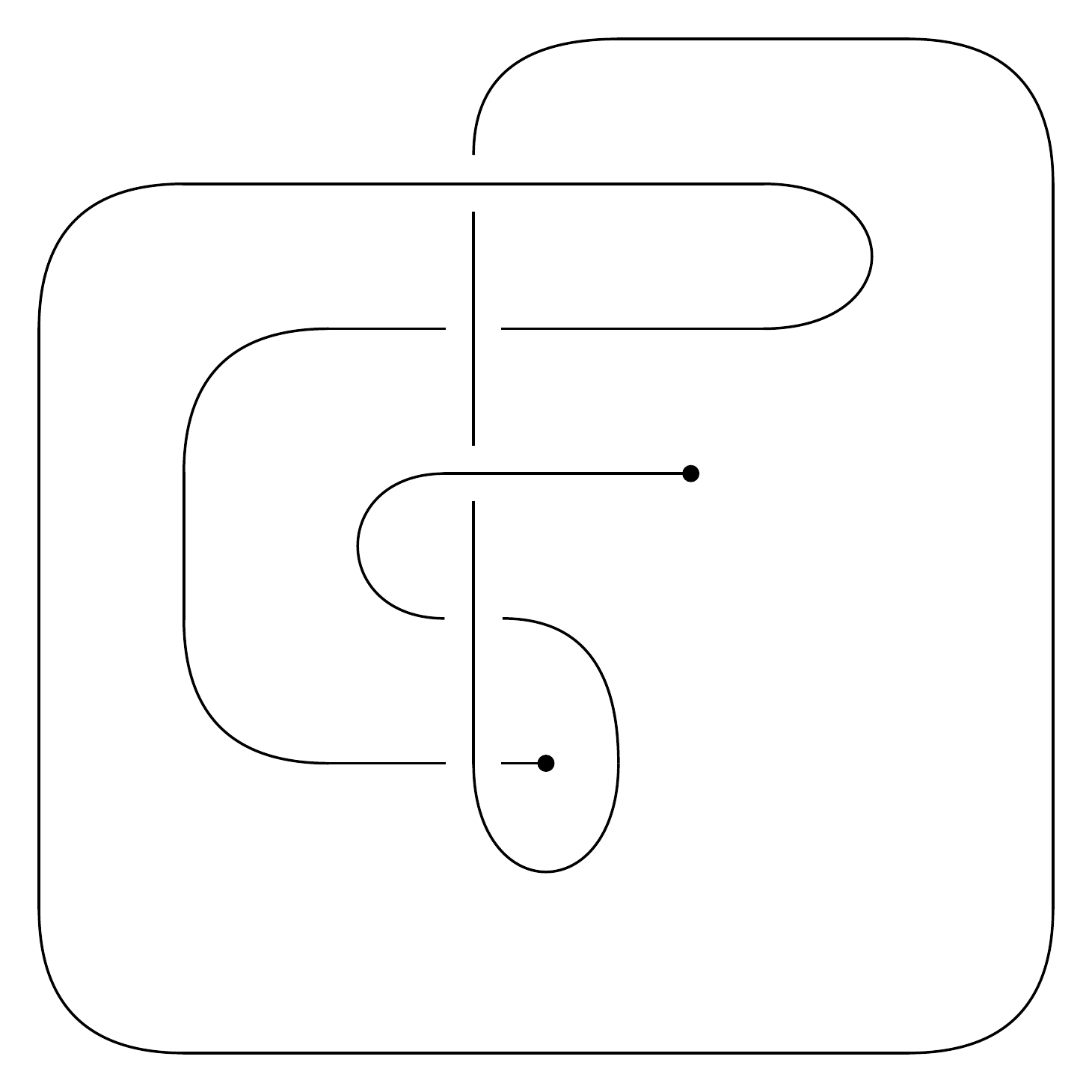}\\
\textcolor{black}{$5_{274}$}
\vspace{1cm}
\end{minipage}
\begin{minipage}[t]{.25\linewidth}
\centering
\includegraphics[width=0.9\textwidth,height=3.5cm,keepaspectratio]{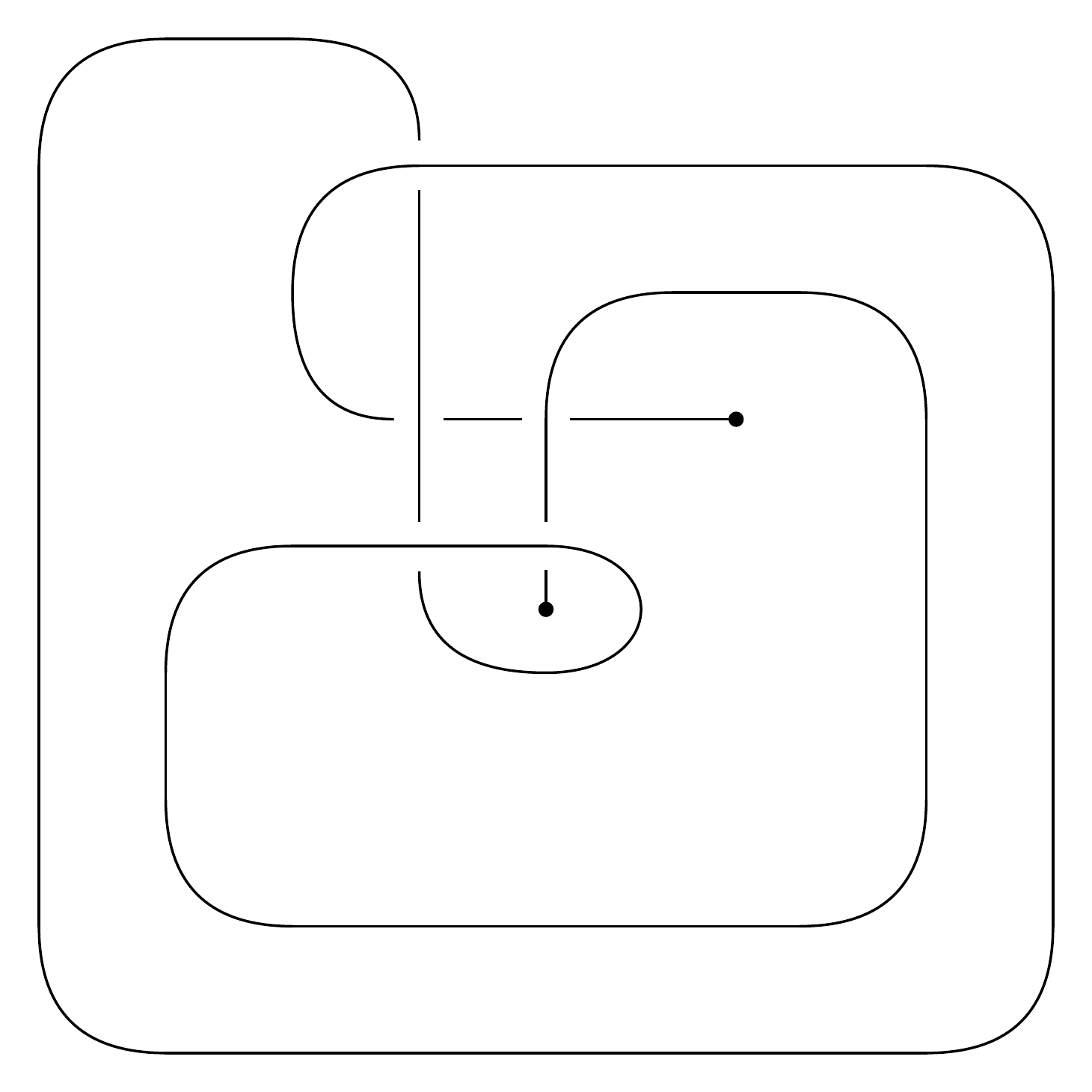}\\
\textcolor{black}{$5_{275}$}
\vspace{1cm}
\end{minipage}
\begin{minipage}[t]{.25\linewidth}
\centering
\includegraphics[width=0.9\textwidth,height=3.5cm,keepaspectratio]{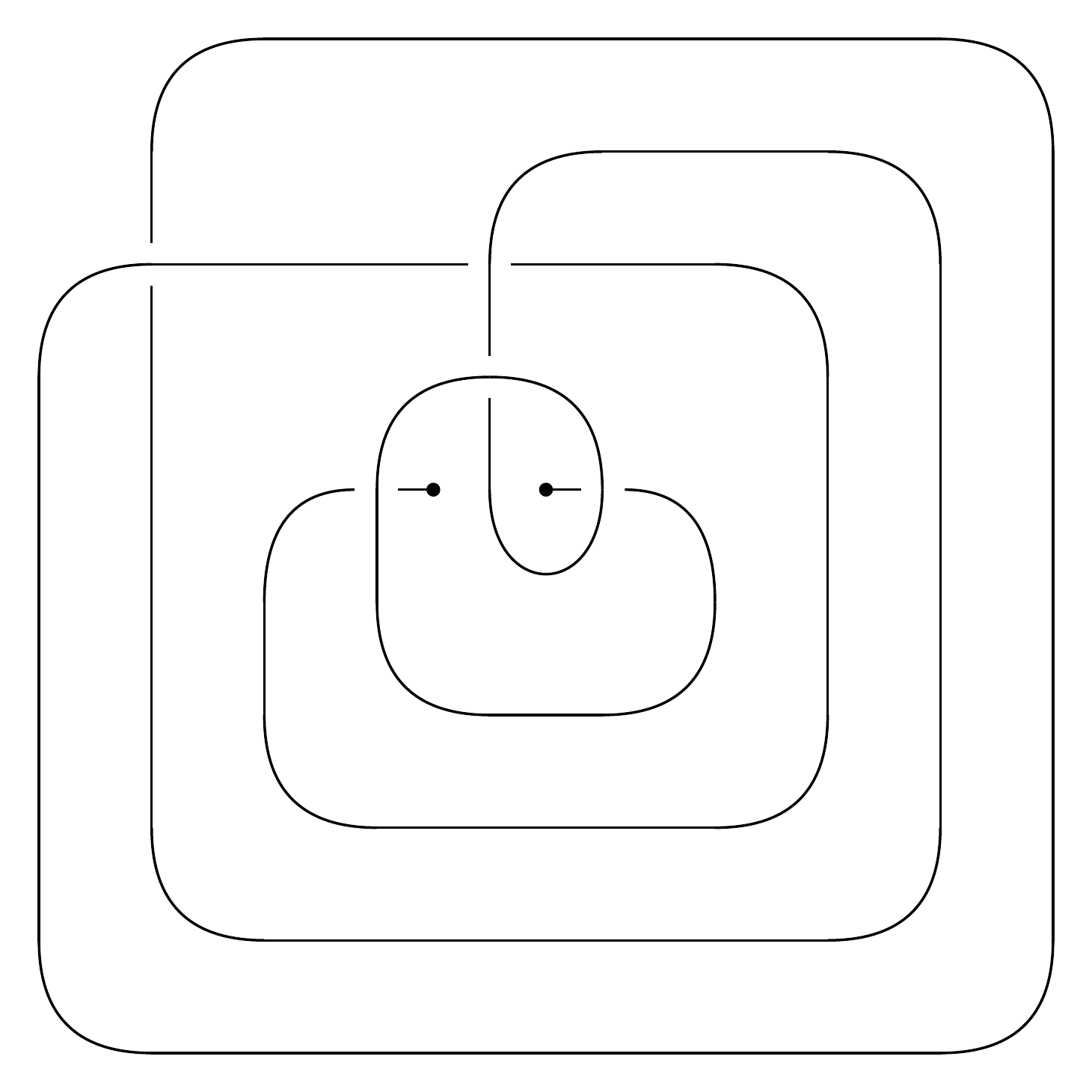}\\
\textcolor{black}{$5_{276}$}
\vspace{1cm}
\end{minipage}
\begin{minipage}[t]{.25\linewidth}
\centering
\includegraphics[width=0.9\textwidth,height=3.5cm,keepaspectratio]{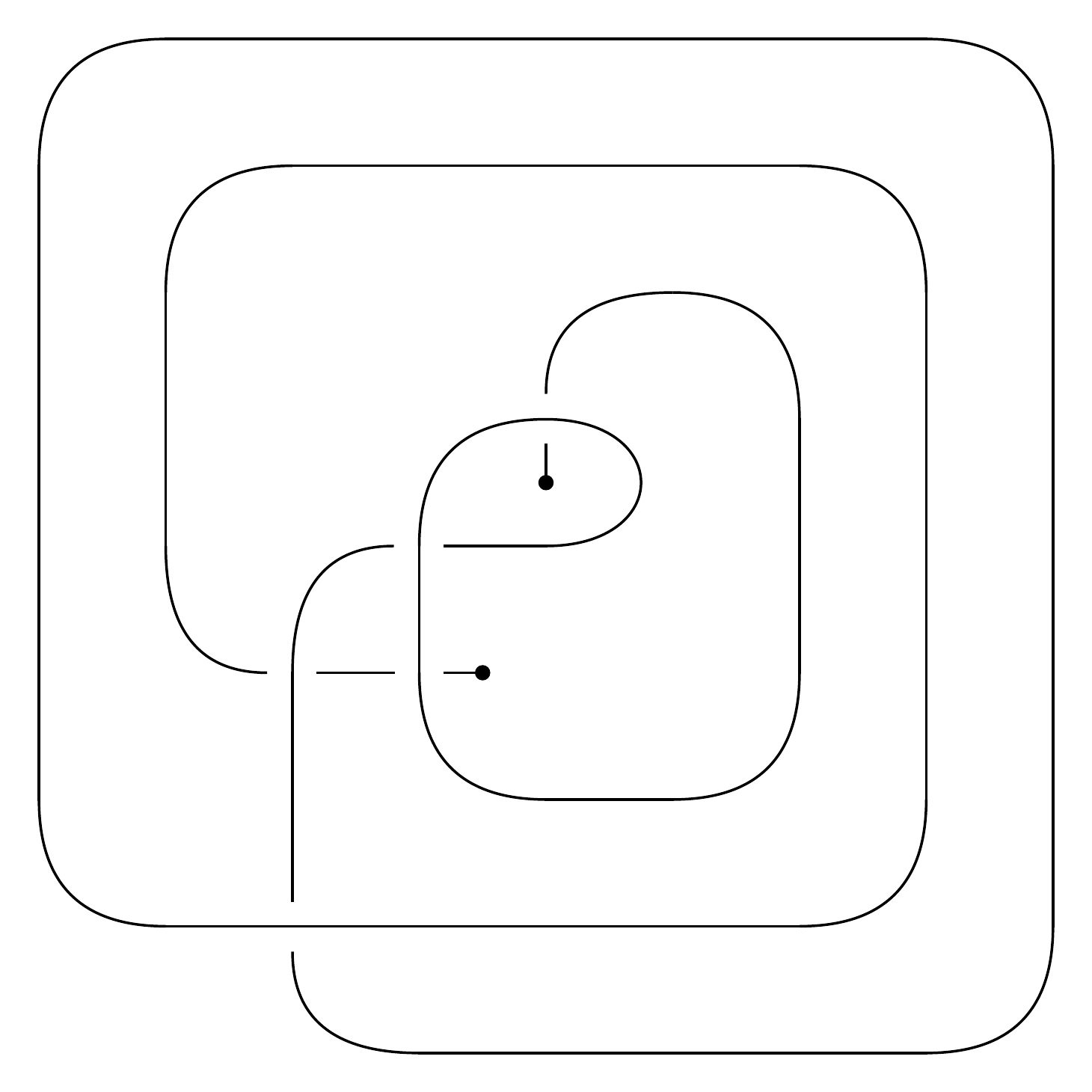}\\
\textcolor{black}{$5_{277}$}
\vspace{1cm}
\end{minipage}
\begin{minipage}[t]{.25\linewidth}
\centering
\includegraphics[width=0.9\textwidth,height=3.5cm,keepaspectratio]{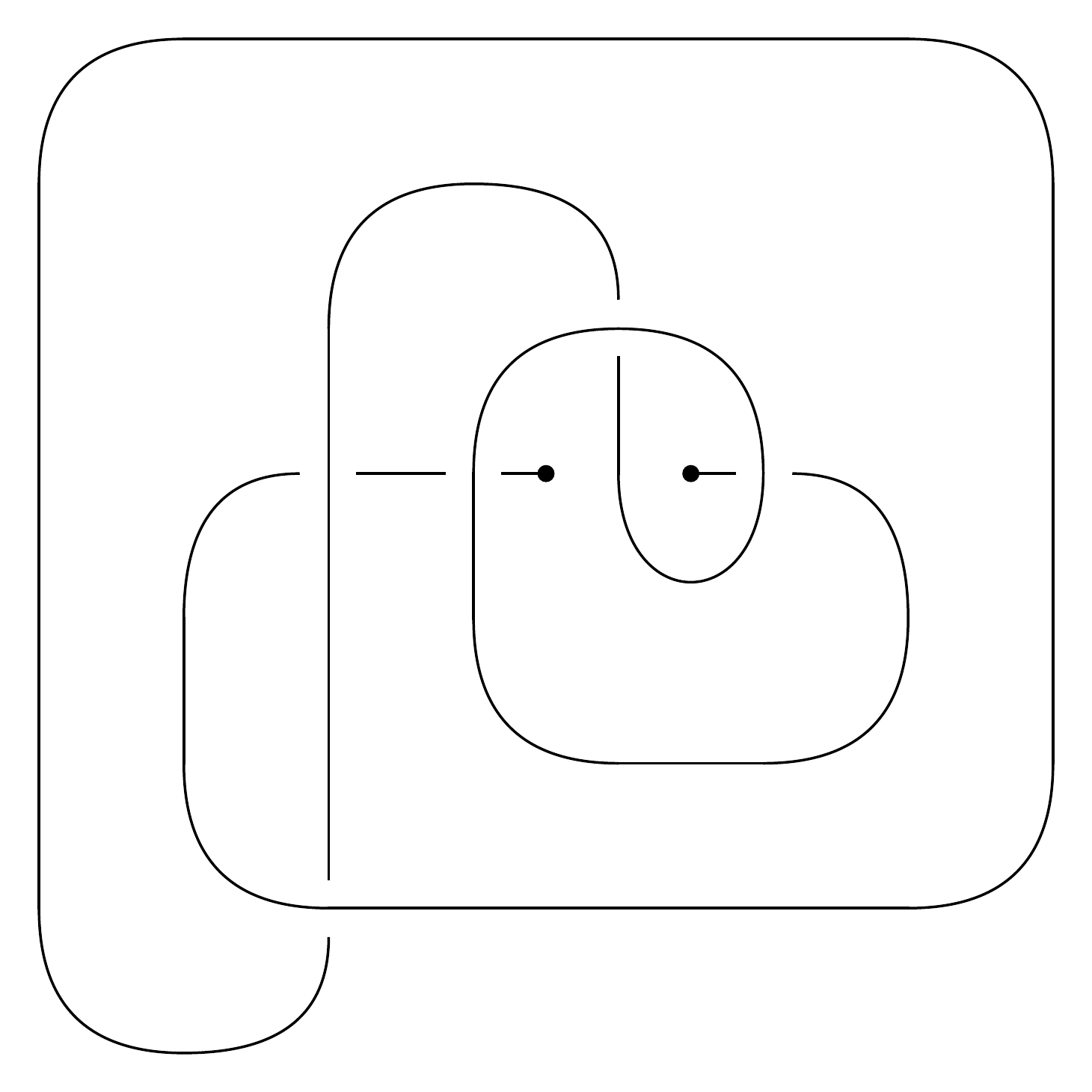}\\
\textcolor{black}{$5_{278}$}
\vspace{1cm}
\end{minipage}
\begin{minipage}[t]{.25\linewidth}
\centering
\includegraphics[width=0.9\textwidth,height=3.5cm,keepaspectratio]{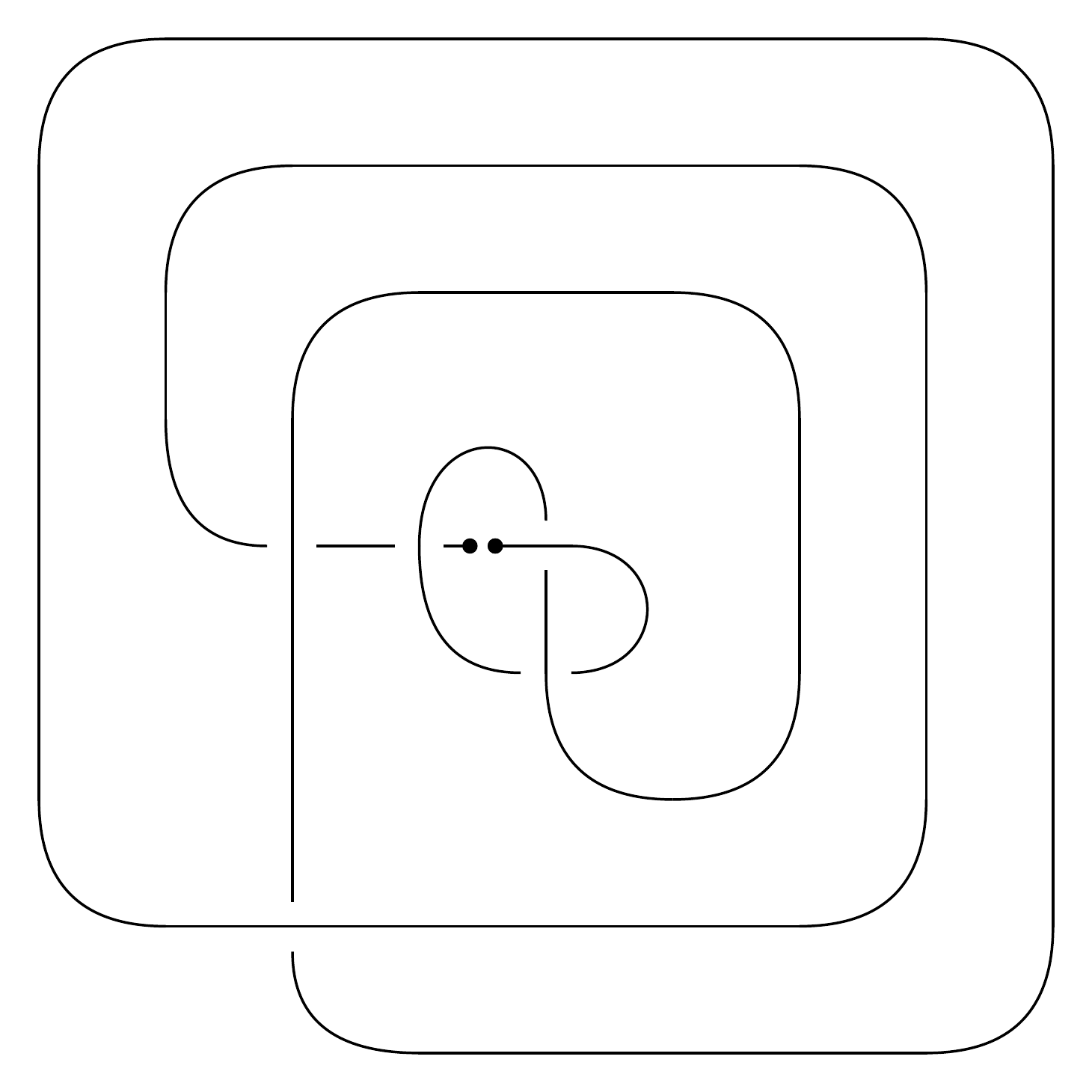}\\
\textcolor{black}{$5_{279}$}
\vspace{1cm}
\end{minipage}
\begin{minipage}[t]{.25\linewidth}
\centering
\includegraphics[width=0.9\textwidth,height=3.5cm,keepaspectratio]{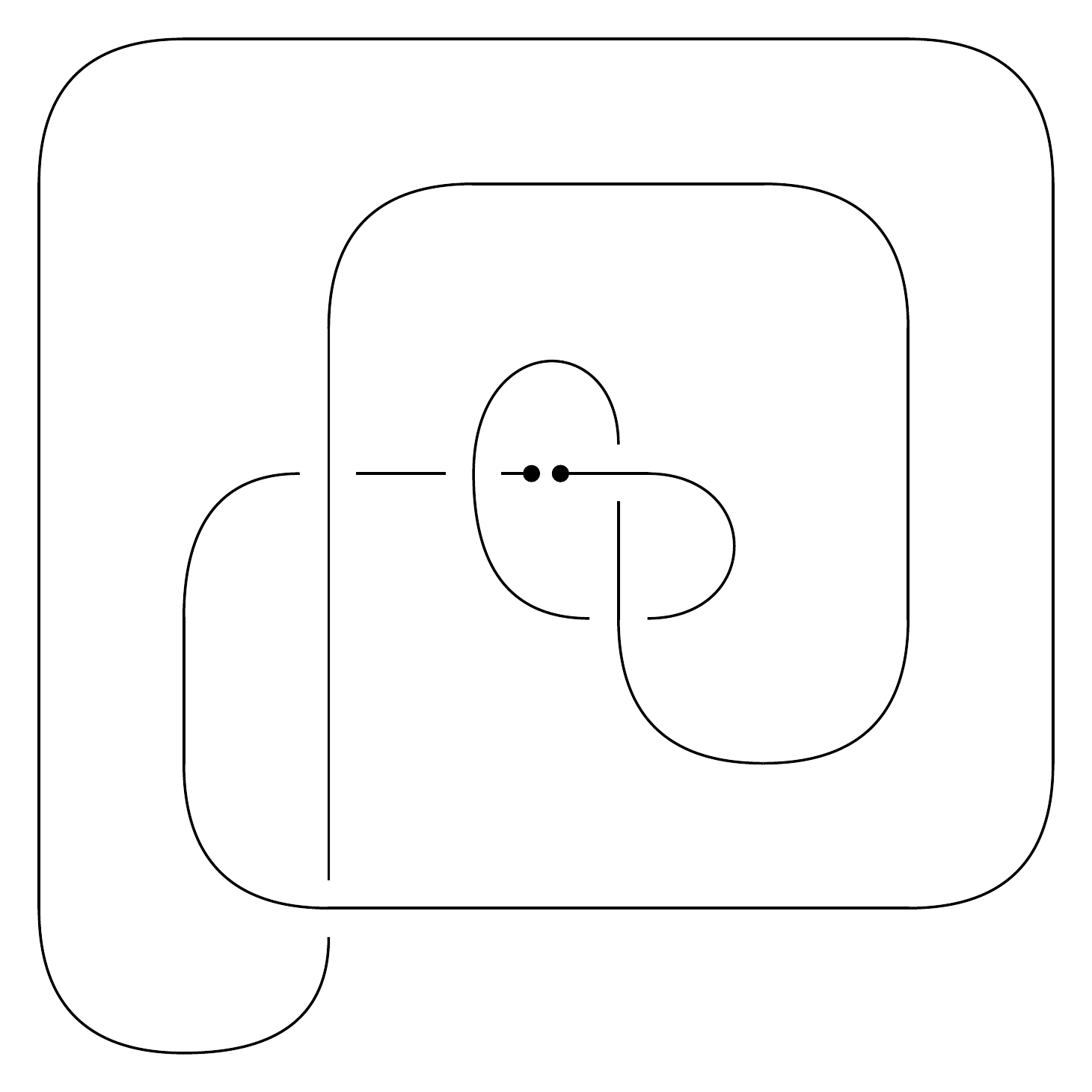}\\
\textcolor{black}{$5_{280}$}
\vspace{1cm}
\end{minipage}
\begin{minipage}[t]{.25\linewidth}
\centering
\includegraphics[width=0.9\textwidth,height=3.5cm,keepaspectratio]{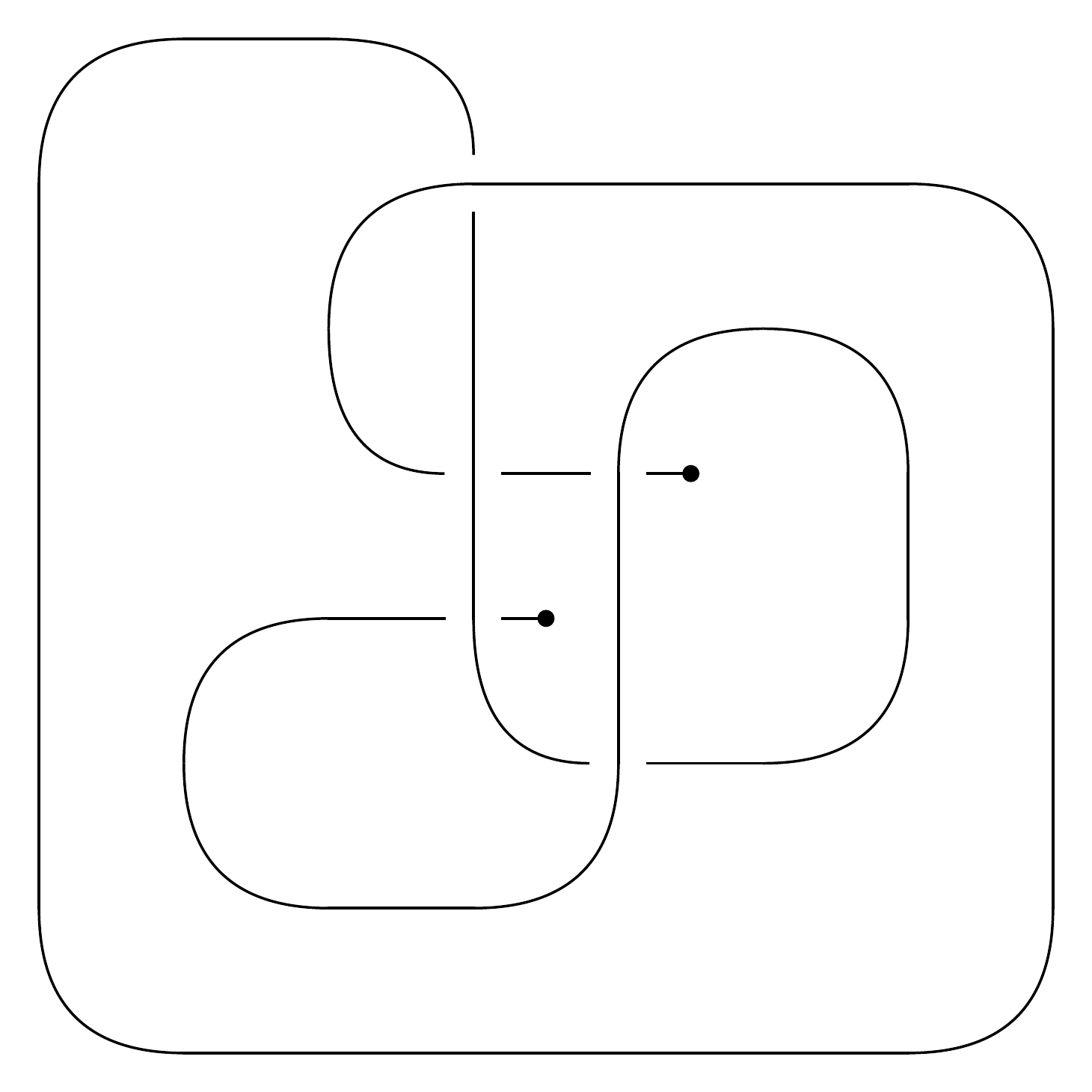}\\
\textcolor{black}{$5_{281}$}
\vspace{1cm}
\end{minipage}
\begin{minipage}[t]{.25\linewidth}
\centering
\includegraphics[width=0.9\textwidth,height=3.5cm,keepaspectratio]{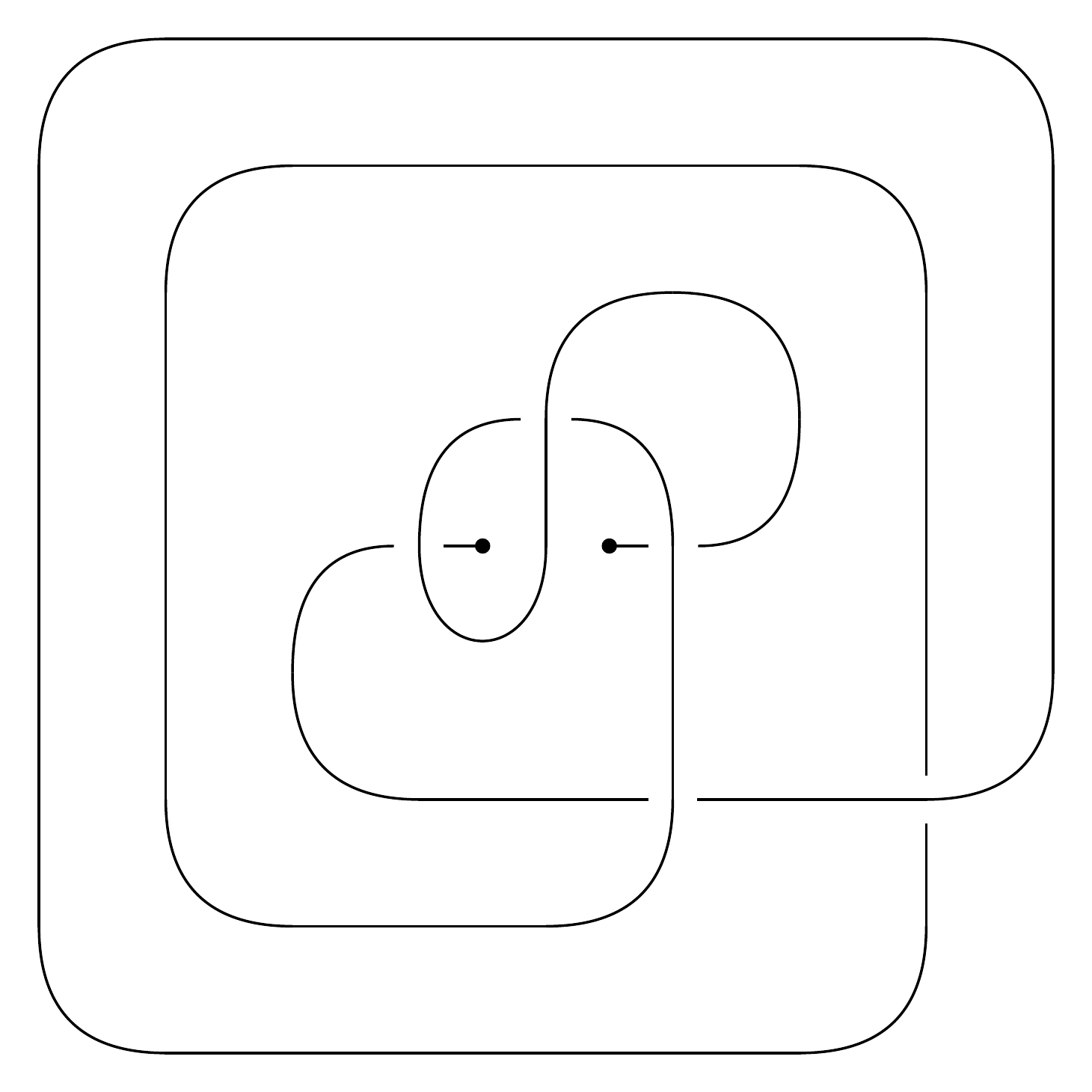}\\
\textcolor{black}{$5_{282}$}
\vspace{1cm}
\end{minipage}
\begin{minipage}[t]{.25\linewidth}
\centering
\includegraphics[width=0.9\textwidth,height=3.5cm,keepaspectratio]{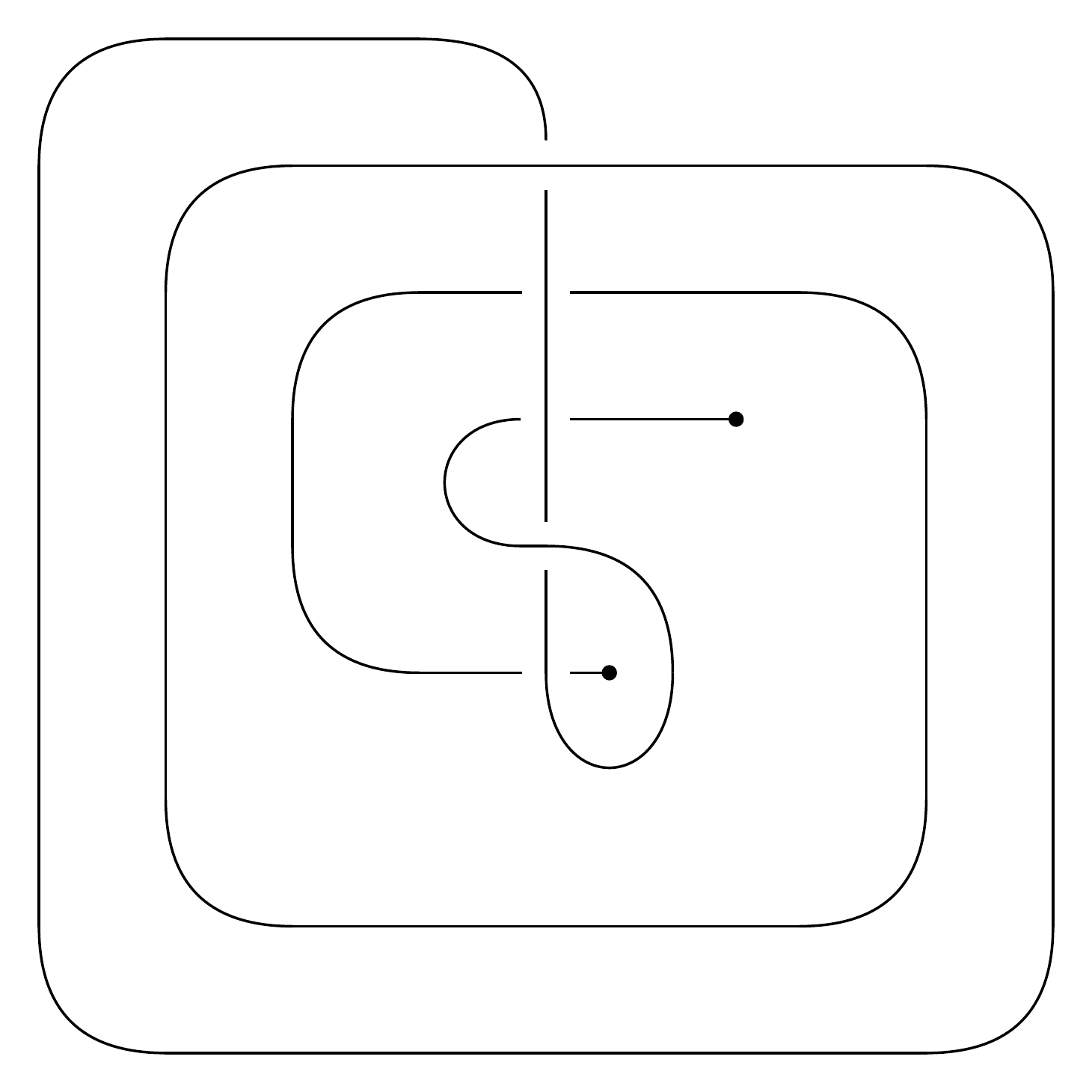}\\
\textcolor{black}{$5_{283}$}
\vspace{1cm}
\end{minipage}
\begin{minipage}[t]{.25\linewidth}
\centering
\includegraphics[width=0.9\textwidth,height=3.5cm,keepaspectratio]{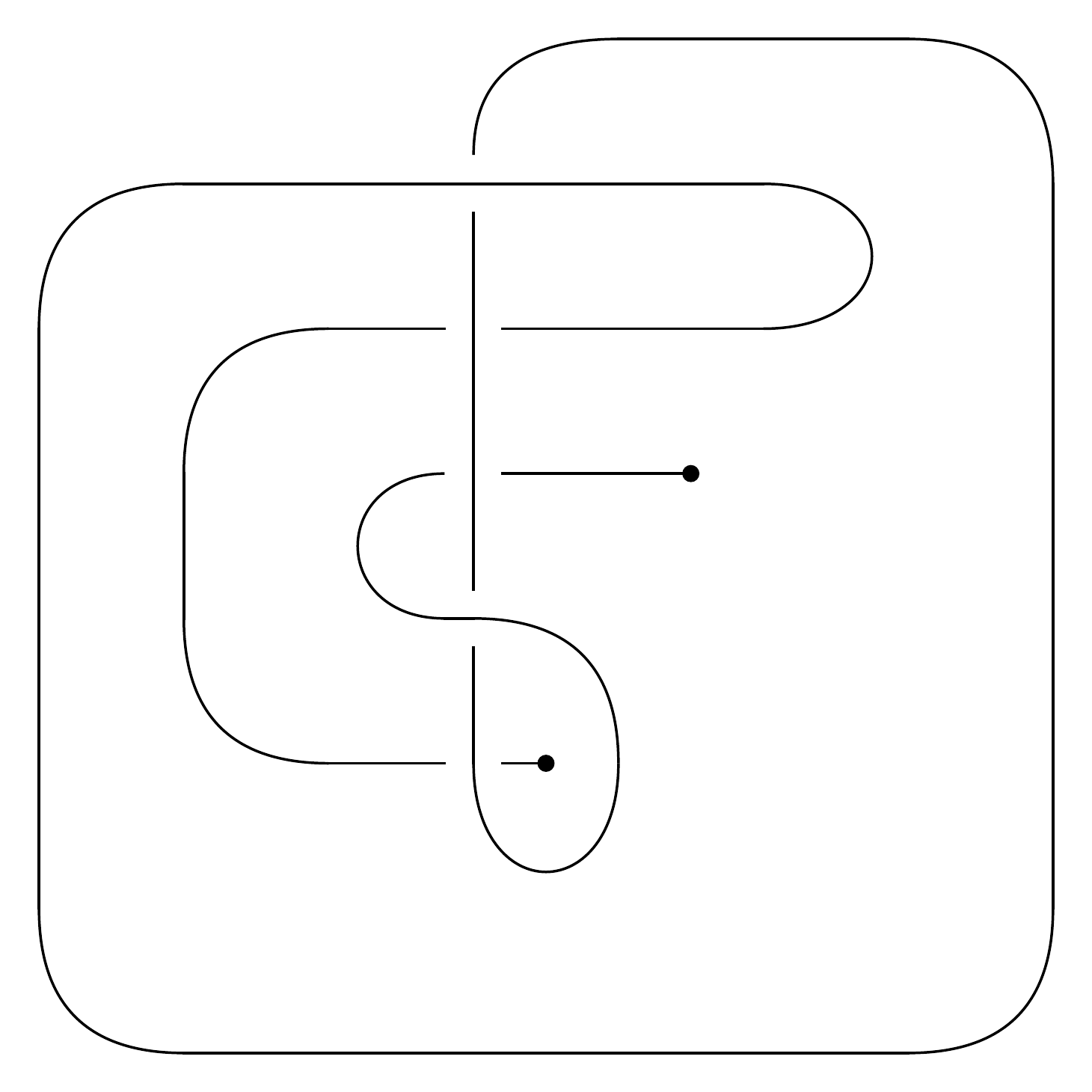}\\
\textcolor{black}{$5_{284}$}
\vspace{1cm}
\end{minipage}
\begin{minipage}[t]{.25\linewidth}
\centering
\includegraphics[width=0.9\textwidth,height=3.5cm,keepaspectratio]{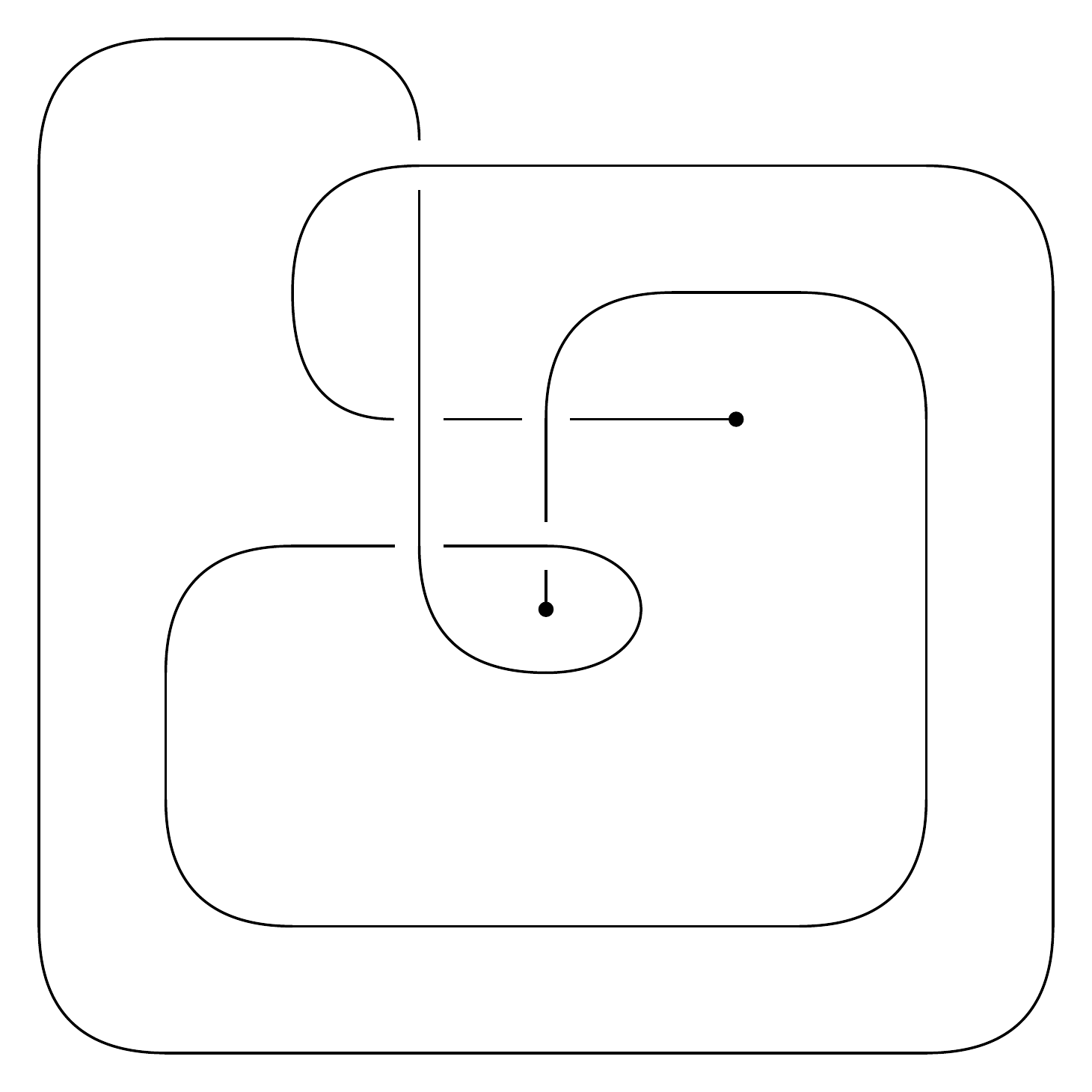}\\
\textcolor{black}{$5_{285}$}
\vspace{1cm}
\end{minipage}
\begin{minipage}[t]{.25\linewidth}
\centering
\includegraphics[width=0.9\textwidth,height=3.5cm,keepaspectratio]{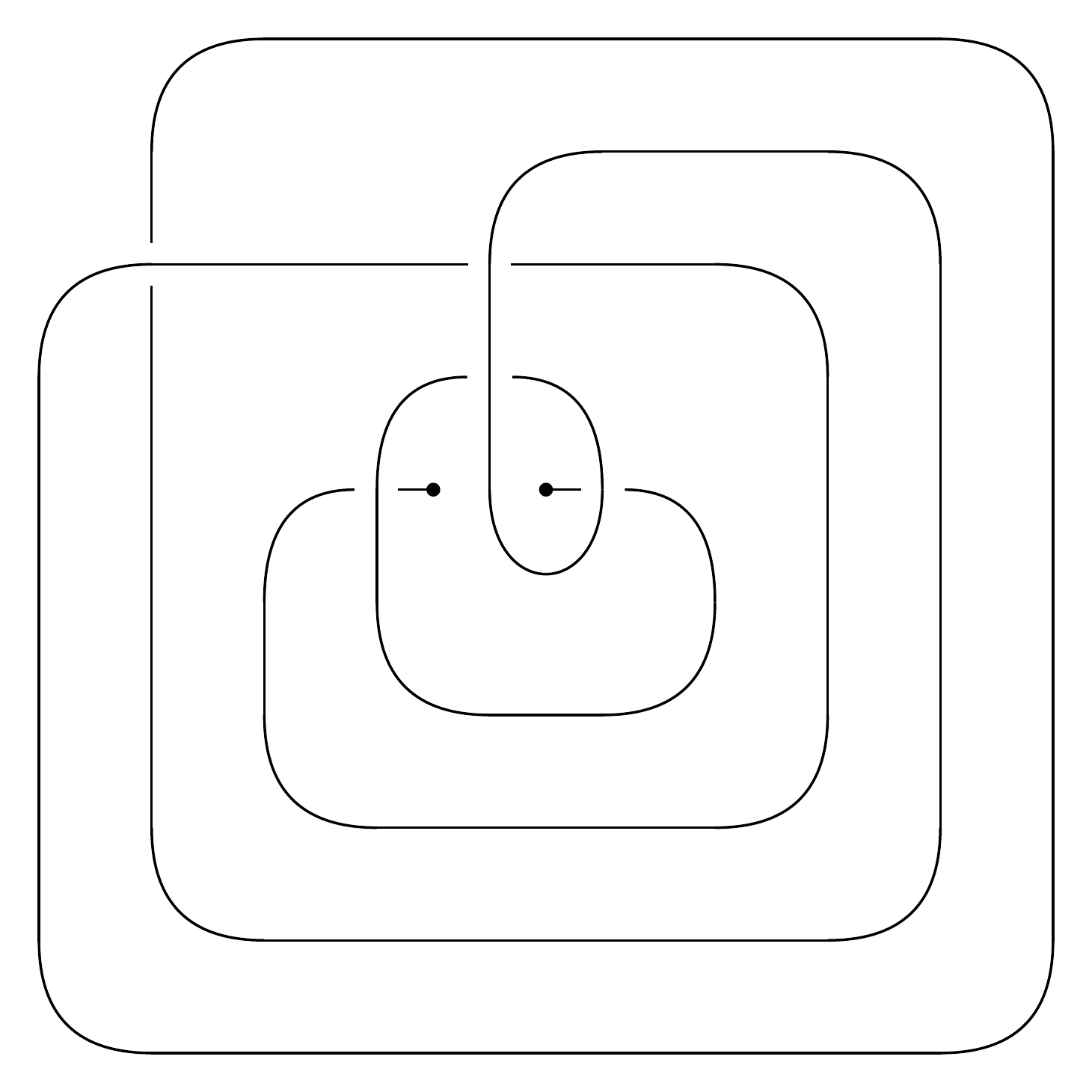}\\
\textcolor{black}{$5_{286}$}
\vspace{1cm}
\end{minipage}
\begin{minipage}[t]{.25\linewidth}
\centering
\includegraphics[width=0.9\textwidth,height=3.5cm,keepaspectratio]{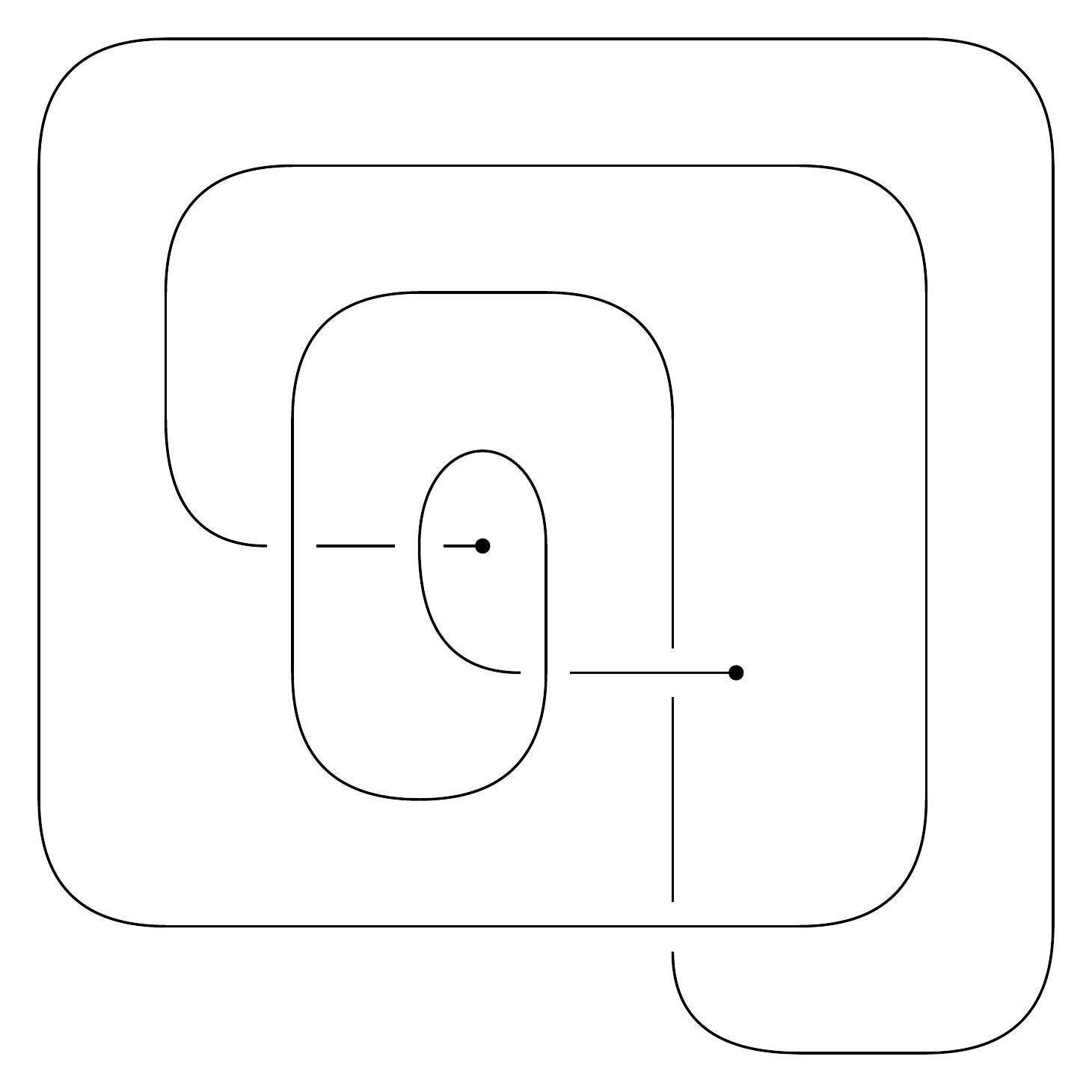}\\
\textcolor{black}{$5_{287}$}
\vspace{1cm}
\end{minipage}
\begin{minipage}[t]{.25\linewidth}
\centering
\includegraphics[width=0.9\textwidth,height=3.5cm,keepaspectratio]{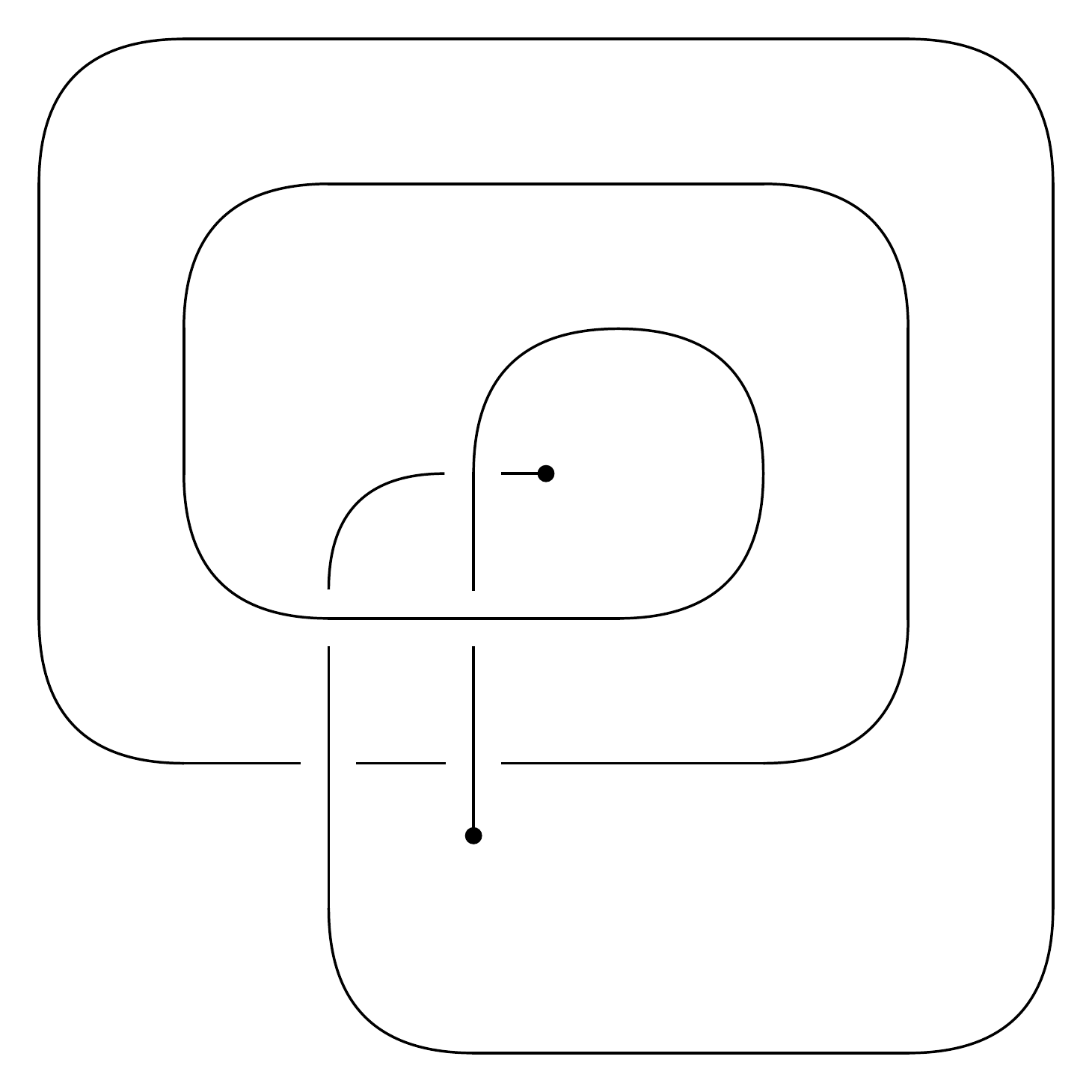}\\
\textcolor{black}{$5_{288}$}
\vspace{1cm}
\end{minipage}
\begin{minipage}[t]{.25\linewidth}
\centering
\includegraphics[width=0.9\textwidth,height=3.5cm,keepaspectratio]{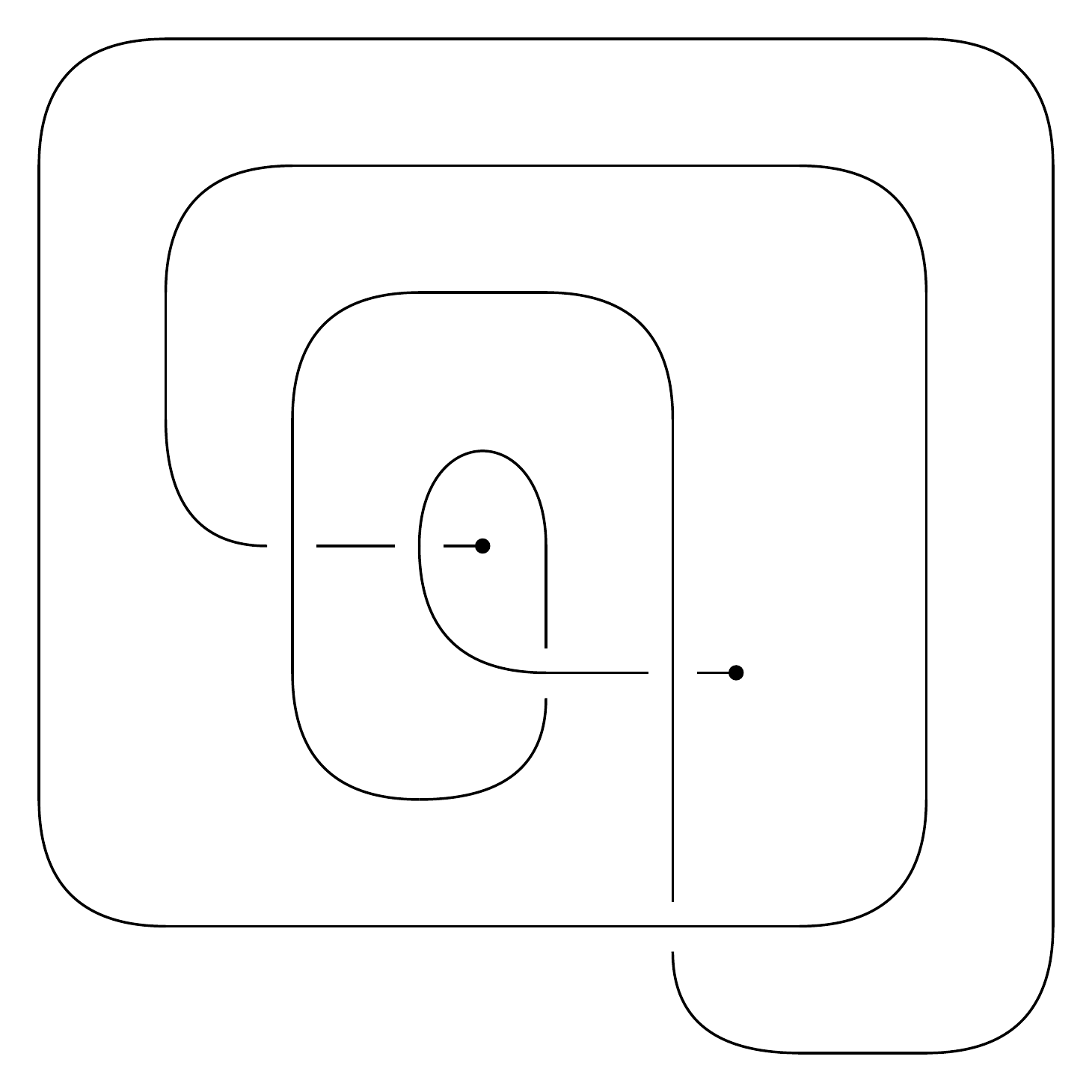}\\
\textcolor{black}{$5_{289}$}
\vspace{1cm}
\end{minipage}
\begin{minipage}[t]{.25\linewidth}
\centering
\includegraphics[width=0.9\textwidth,height=3.5cm,keepaspectratio]{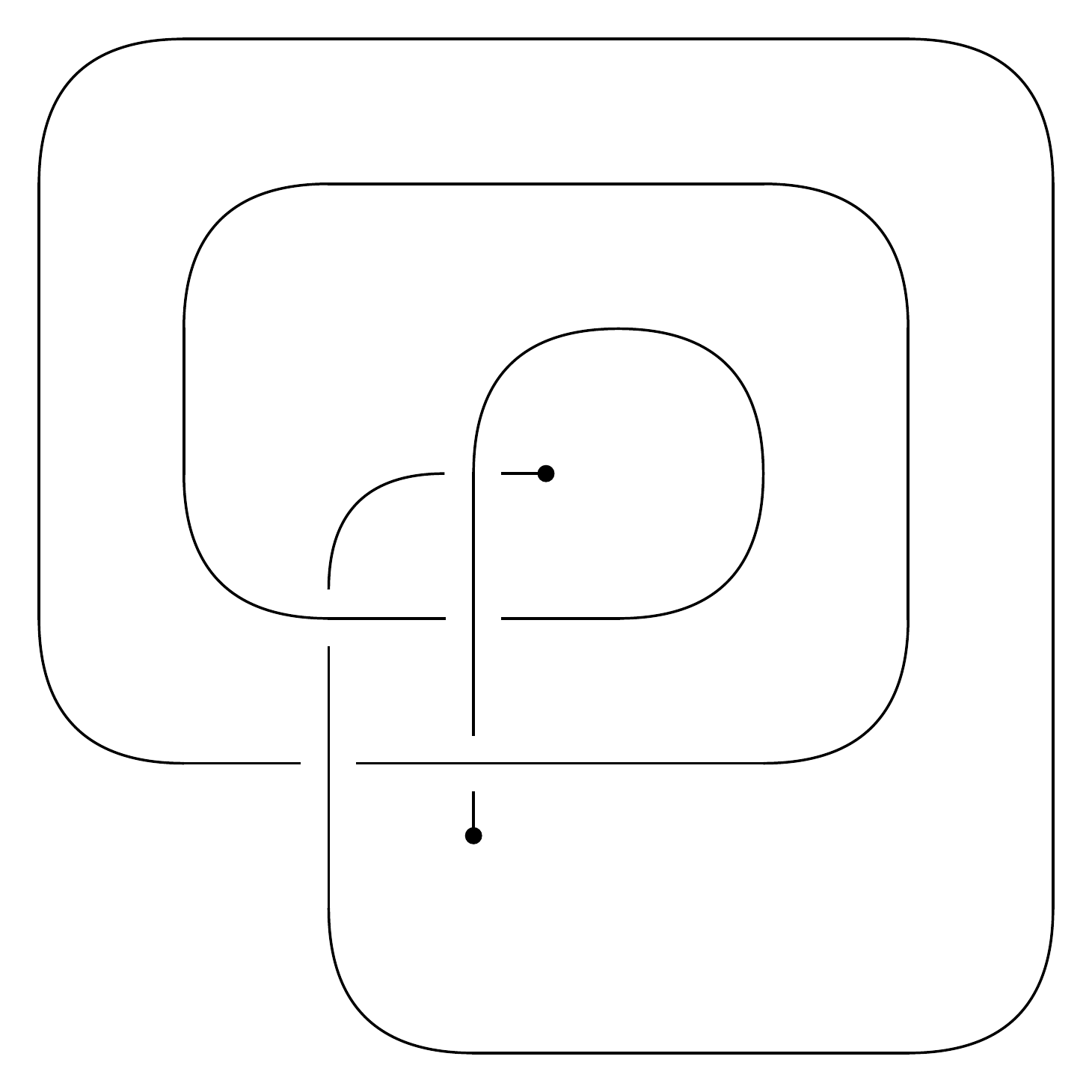}\\
\textcolor{black}{$5_{290}$}
\vspace{1cm}
\end{minipage}
\begin{minipage}[t]{.25\linewidth}
\centering
\includegraphics[width=0.9\textwidth,height=3.5cm,keepaspectratio]{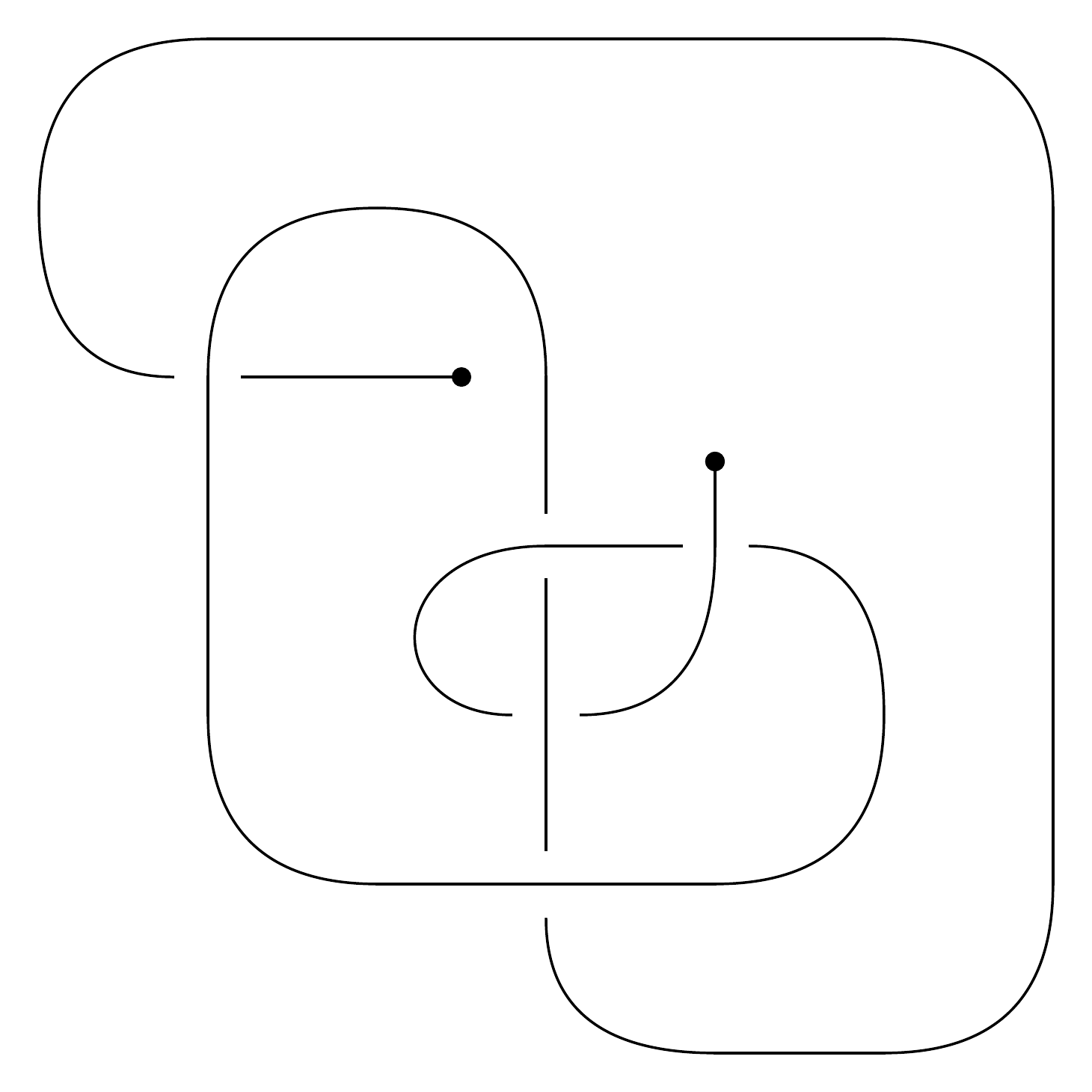}\\
\textcolor{black}{$5_{291}$}
\vspace{1cm}
\end{minipage}
\begin{minipage}[t]{.25\linewidth}
\centering
\includegraphics[width=0.9\textwidth,height=3.5cm,keepaspectratio]{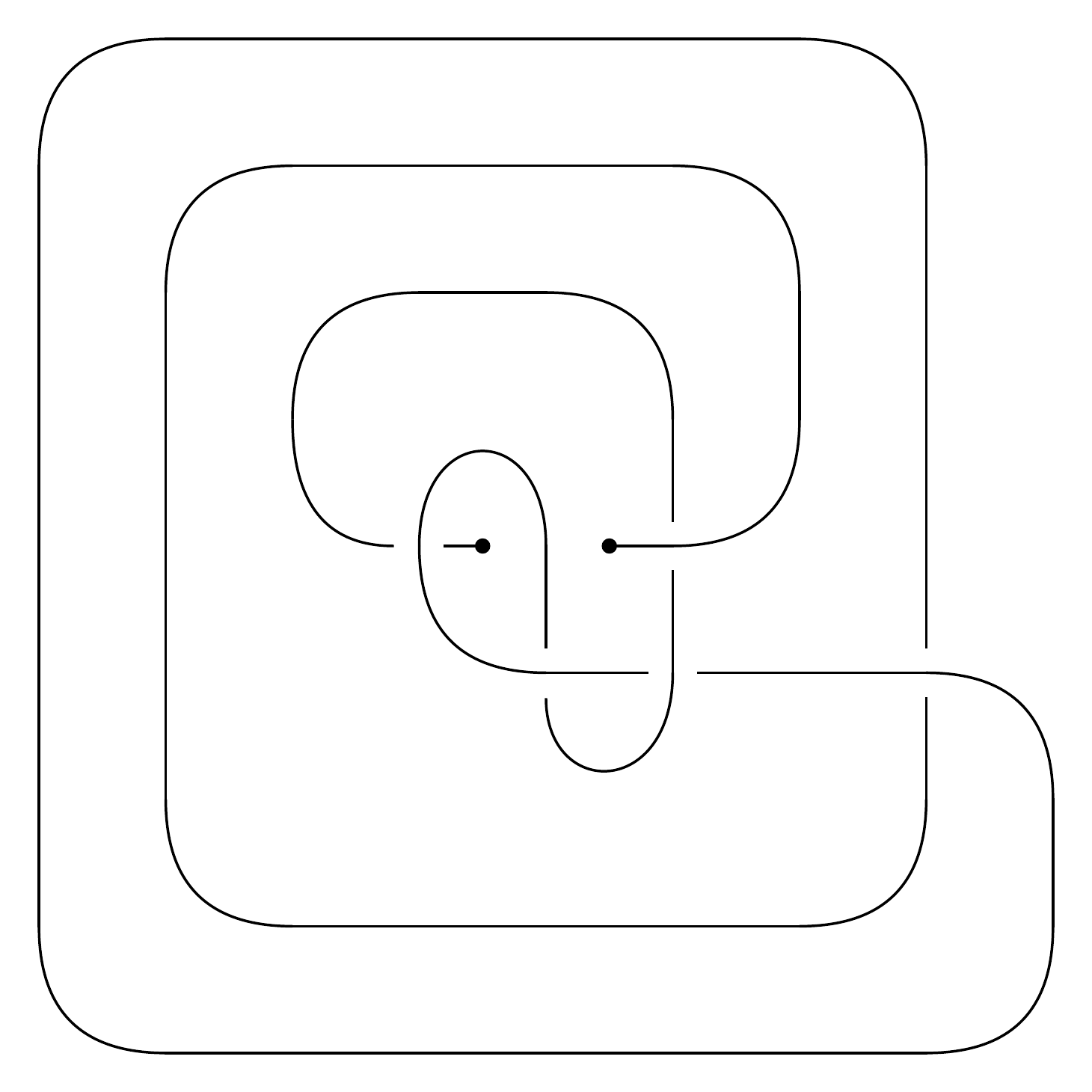}\\
\textcolor{black}{$5_{292}$}
\vspace{1cm}
\end{minipage}
\begin{minipage}[t]{.25\linewidth}
\centering
\includegraphics[width=0.9\textwidth,height=3.5cm,keepaspectratio]{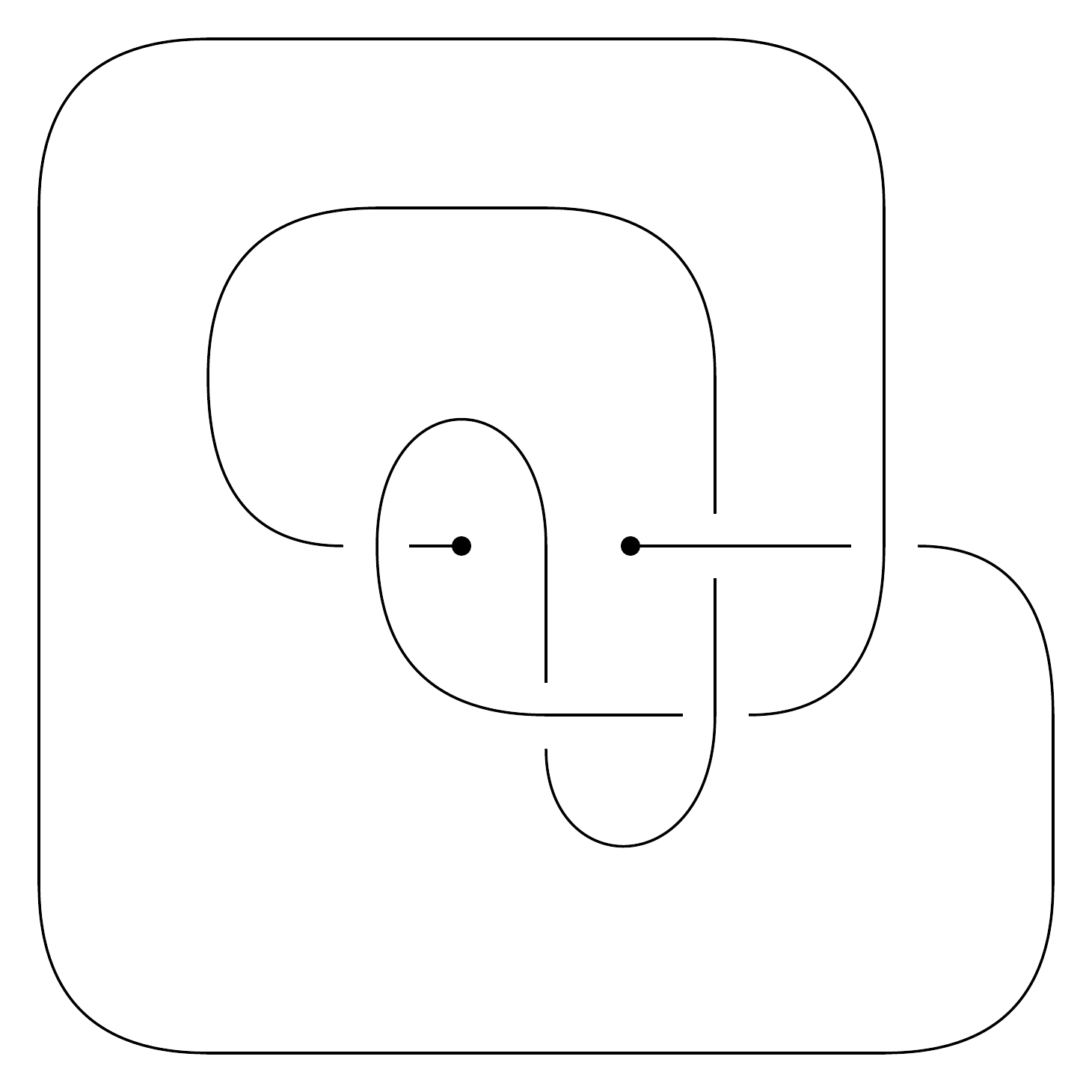}\\
\textcolor{black}{$5_{293}$}
\vspace{1cm}
\end{minipage}
\begin{minipage}[t]{.25\linewidth}
\centering
\includegraphics[width=0.9\textwidth,height=3.5cm,keepaspectratio]{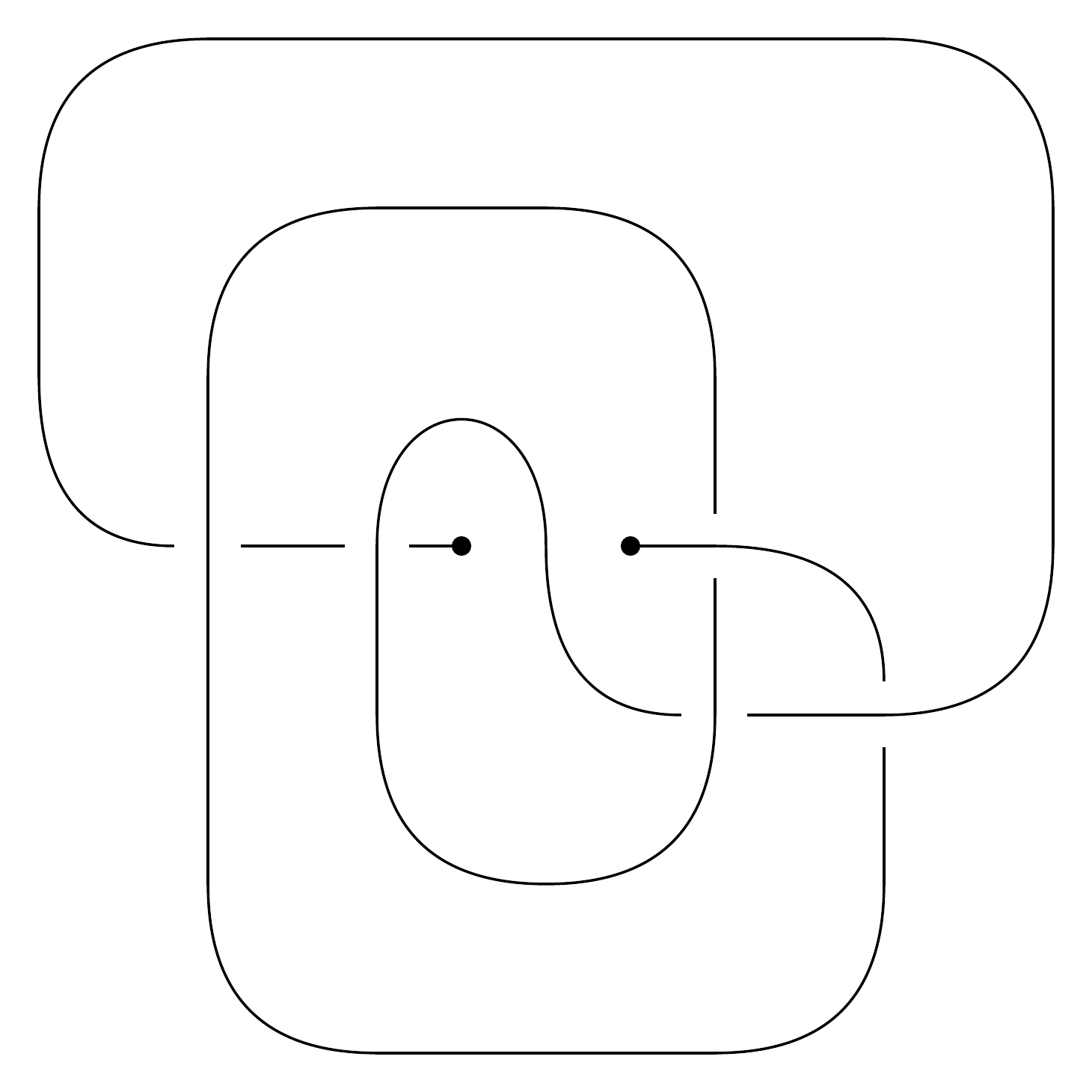}\\
\textcolor{black}{$5_{294}$}
\vspace{1cm}
\end{minipage}
\begin{minipage}[t]{.25\linewidth}
\centering
\includegraphics[width=0.9\textwidth,height=3.5cm,keepaspectratio]{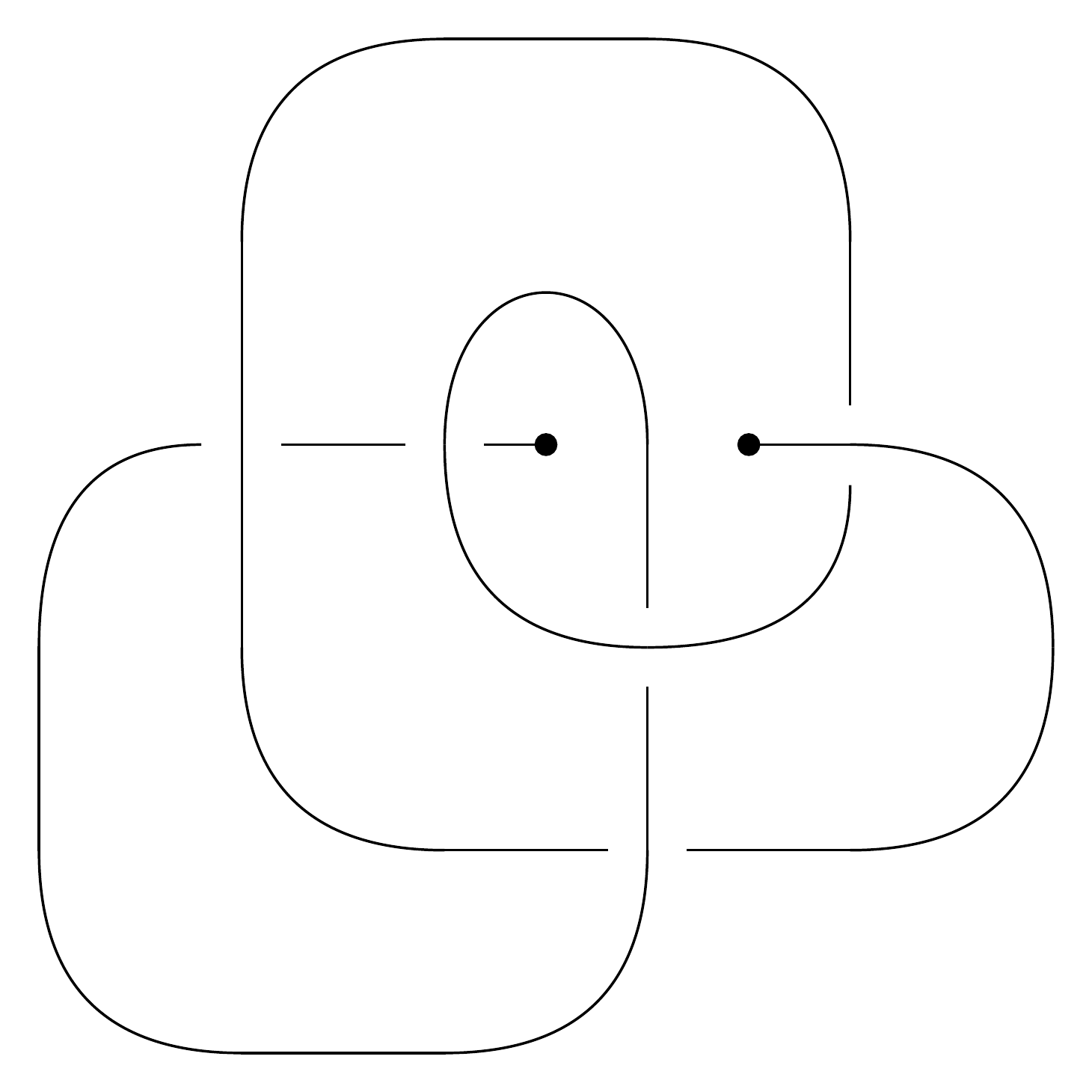}\\
\textcolor{black}{$5_{295}$}
\vspace{1cm}
\end{minipage}
\begin{minipage}[t]{.25\linewidth}
\centering
\includegraphics[width=0.9\textwidth,height=3.5cm,keepaspectratio]{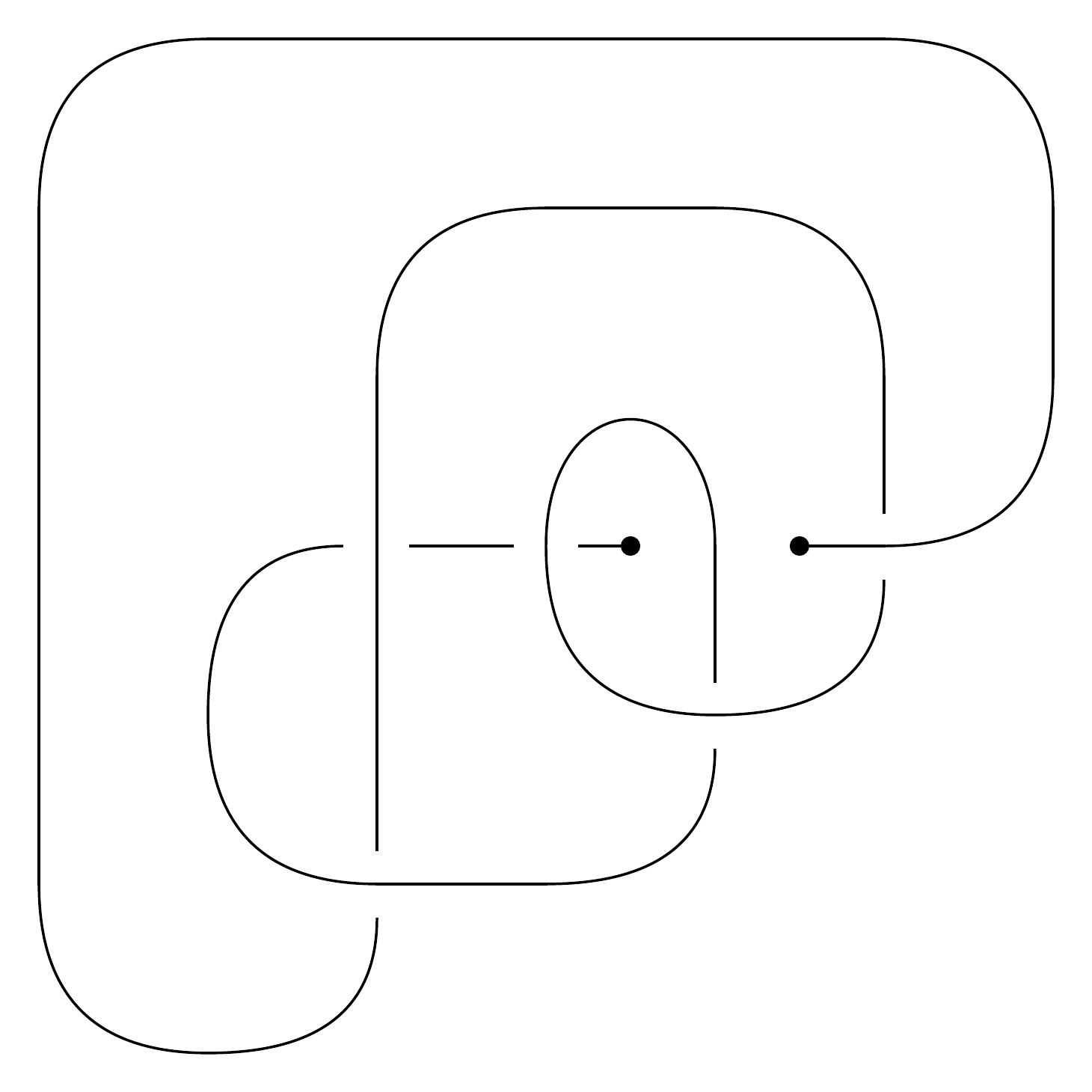}\\
\textcolor{black}{$5_{296}$}
\vspace{1cm}
\end{minipage}
\begin{minipage}[t]{.25\linewidth}
\centering
\includegraphics[width=0.9\textwidth,height=3.5cm,keepaspectratio]{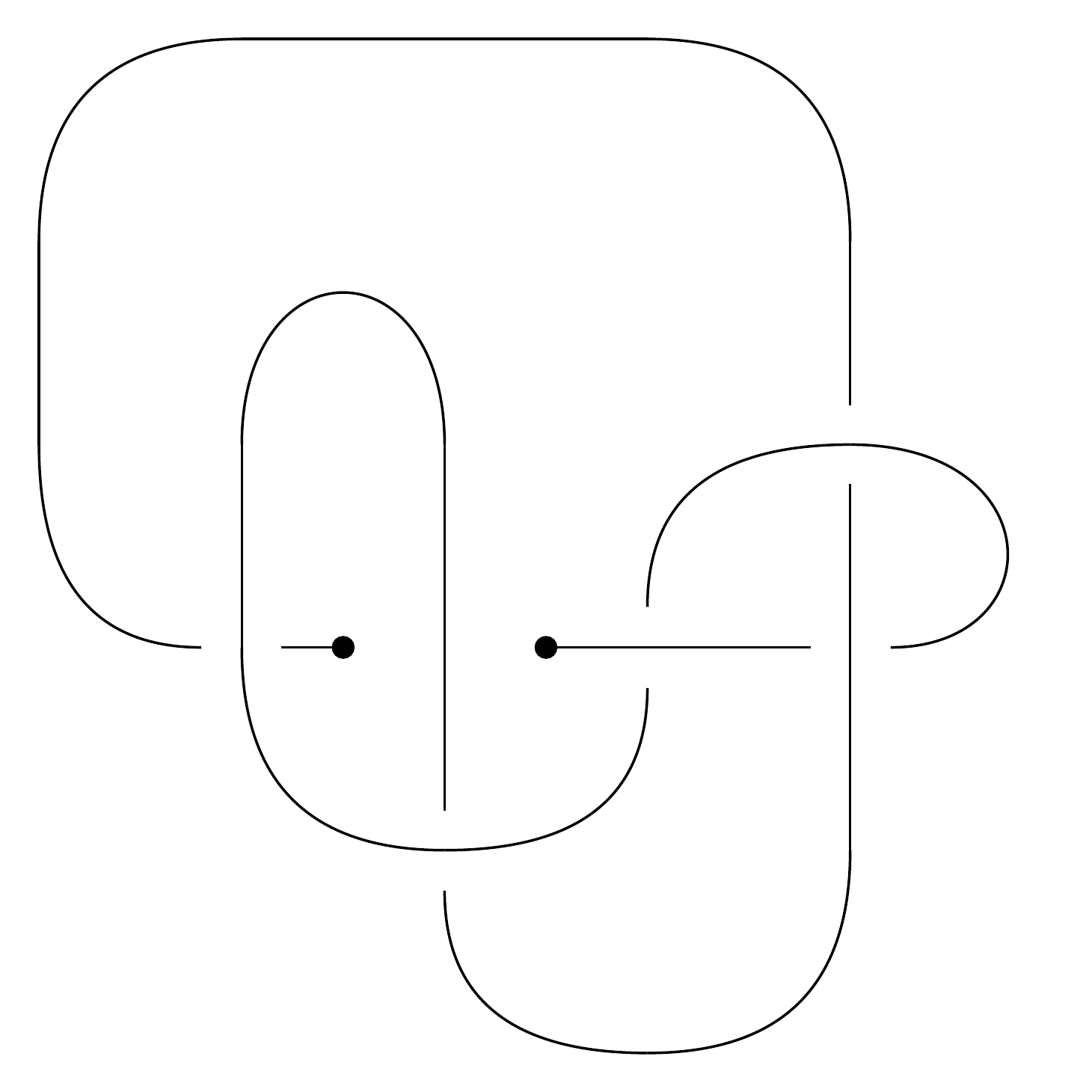}\\
\textcolor{black}{$5_{297}$}
\vspace{1cm}
\end{minipage}
\begin{minipage}[t]{.25\linewidth}
\centering
\includegraphics[width=0.9\textwidth,height=3.5cm,keepaspectratio]{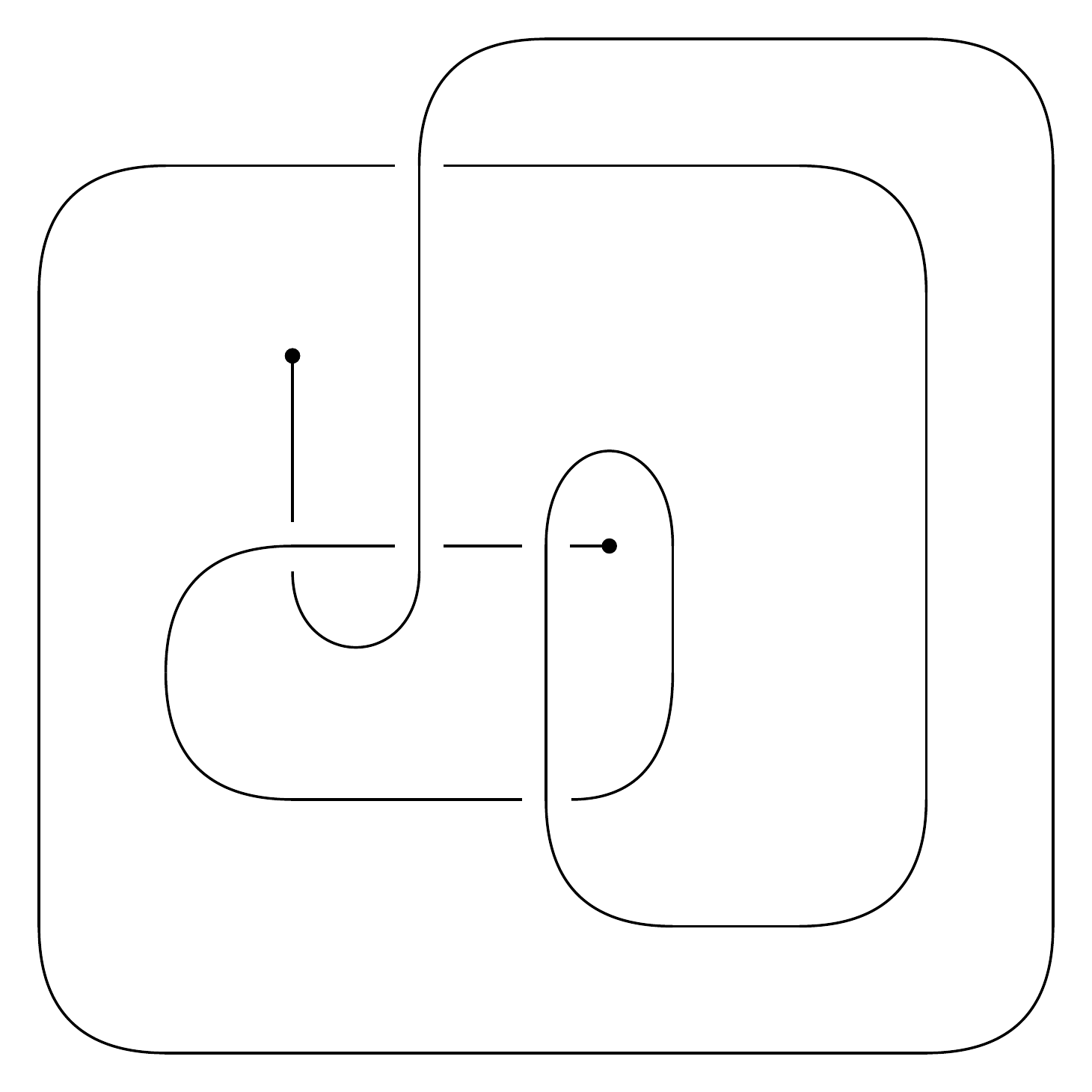}\\
\textcolor{black}{$5_{298}$}
\vspace{1cm}
\end{minipage}
\begin{minipage}[t]{.25\linewidth}
\centering
\includegraphics[width=0.9\textwidth,height=3.5cm,keepaspectratio]{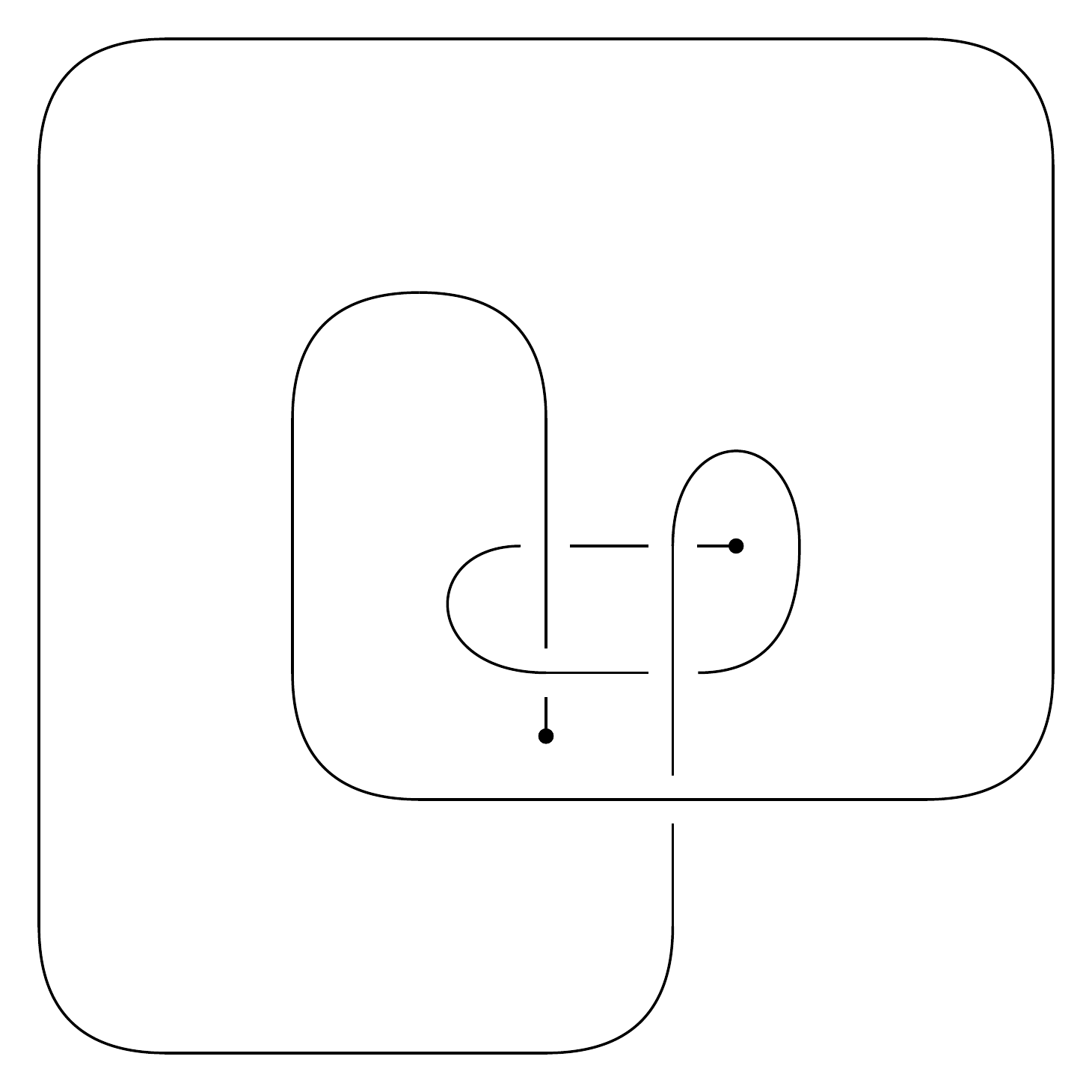}\\
\textcolor{black}{$5_{299}$}
\vspace{1cm}
\end{minipage}
\begin{minipage}[t]{.25\linewidth}
\centering
\includegraphics[width=0.9\textwidth,height=3.5cm,keepaspectratio]{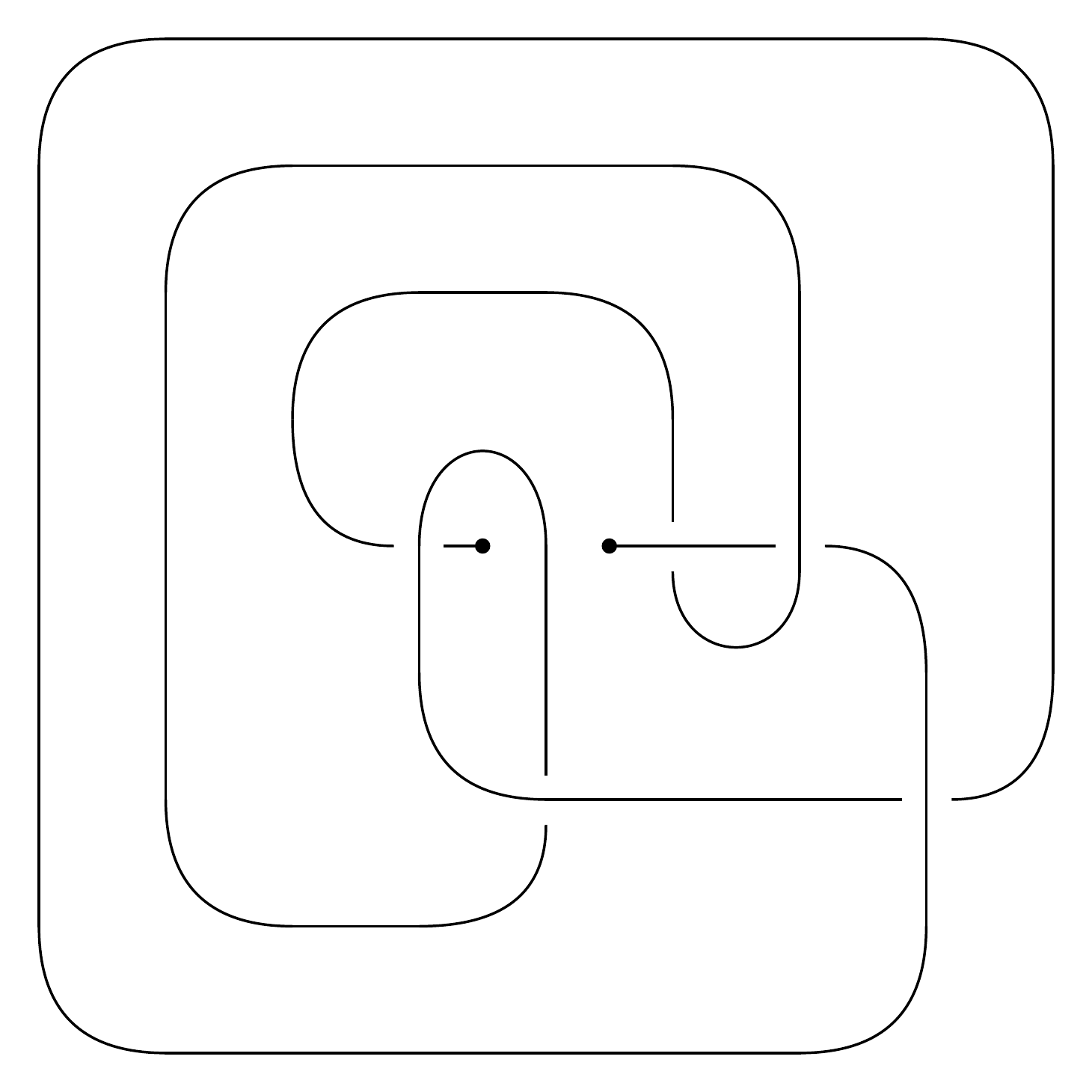}\\
\textcolor{black}{$5_{300}$}
\vspace{1cm}
\end{minipage}
\begin{minipage}[t]{.25\linewidth}
\centering
\includegraphics[width=0.9\textwidth,height=3.5cm,keepaspectratio]{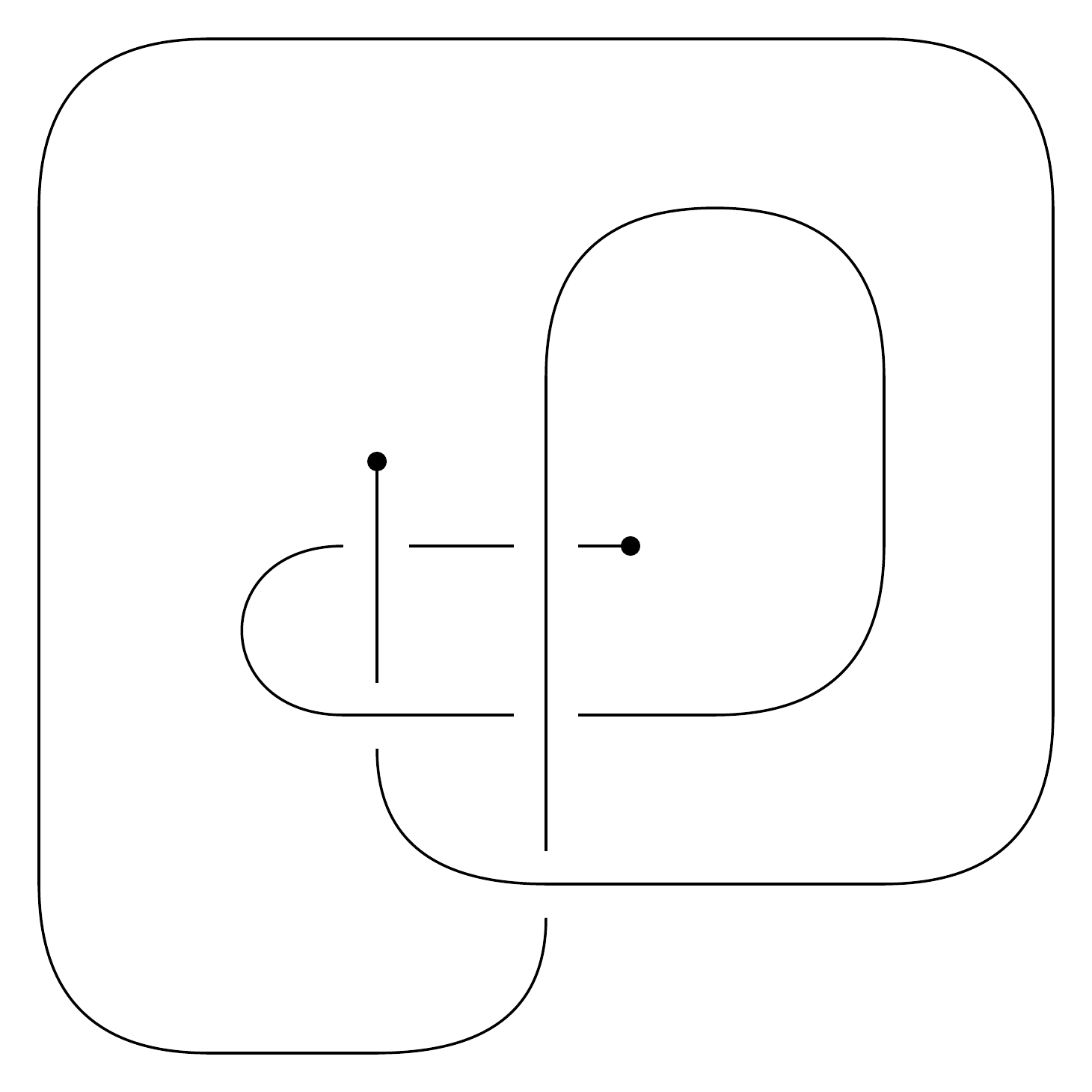}\\
\textcolor{black}{$5_{301}$}
\vspace{1cm}
\end{minipage}
\begin{minipage}[t]{.25\linewidth}
\centering
\includegraphics[width=0.9\textwidth,height=3.5cm,keepaspectratio]{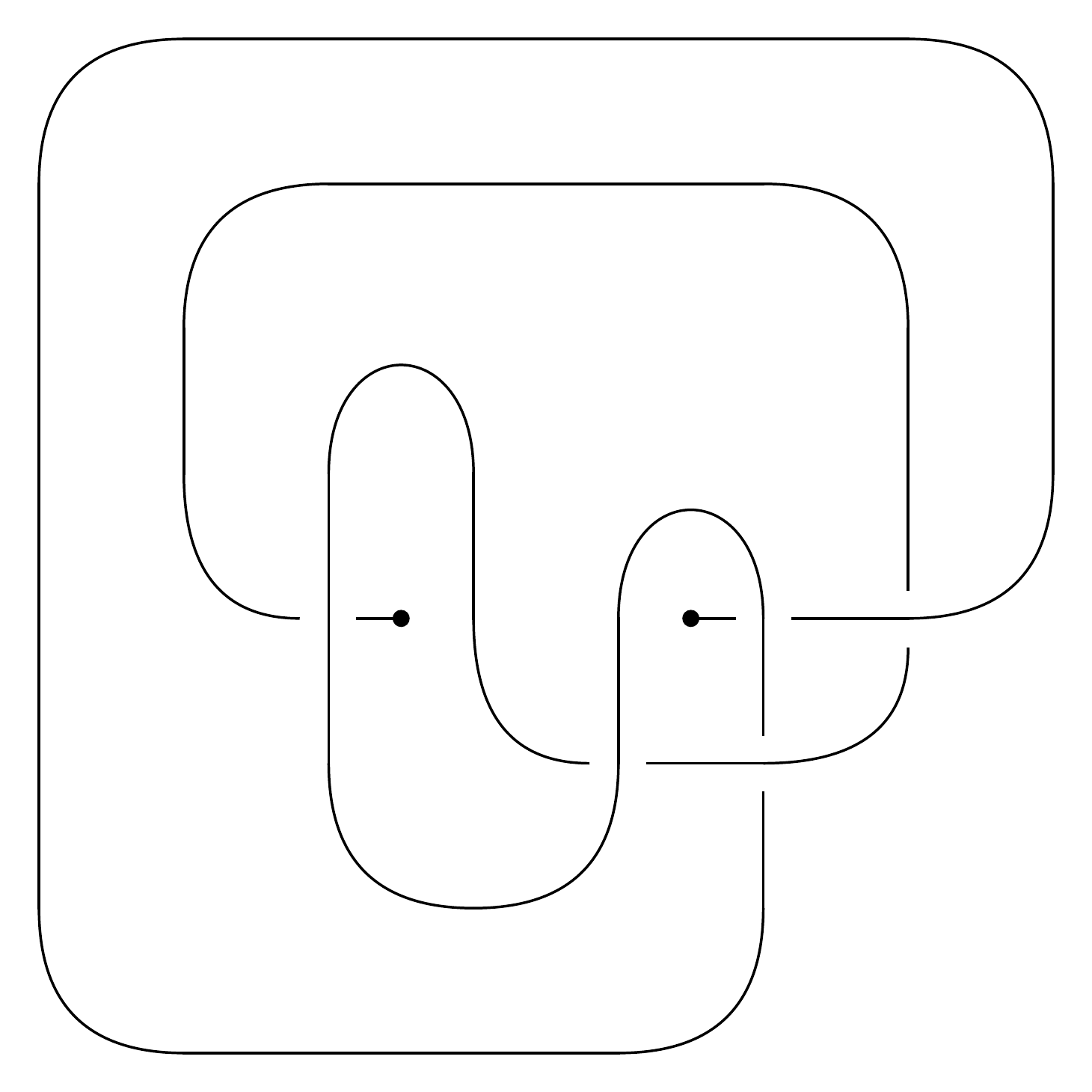}\\
\textcolor{black}{$5_{302}$}
\vspace{1cm}
\end{minipage}
\begin{minipage}[t]{.25\linewidth}
\centering
\includegraphics[width=0.9\textwidth,height=3.5cm,keepaspectratio]{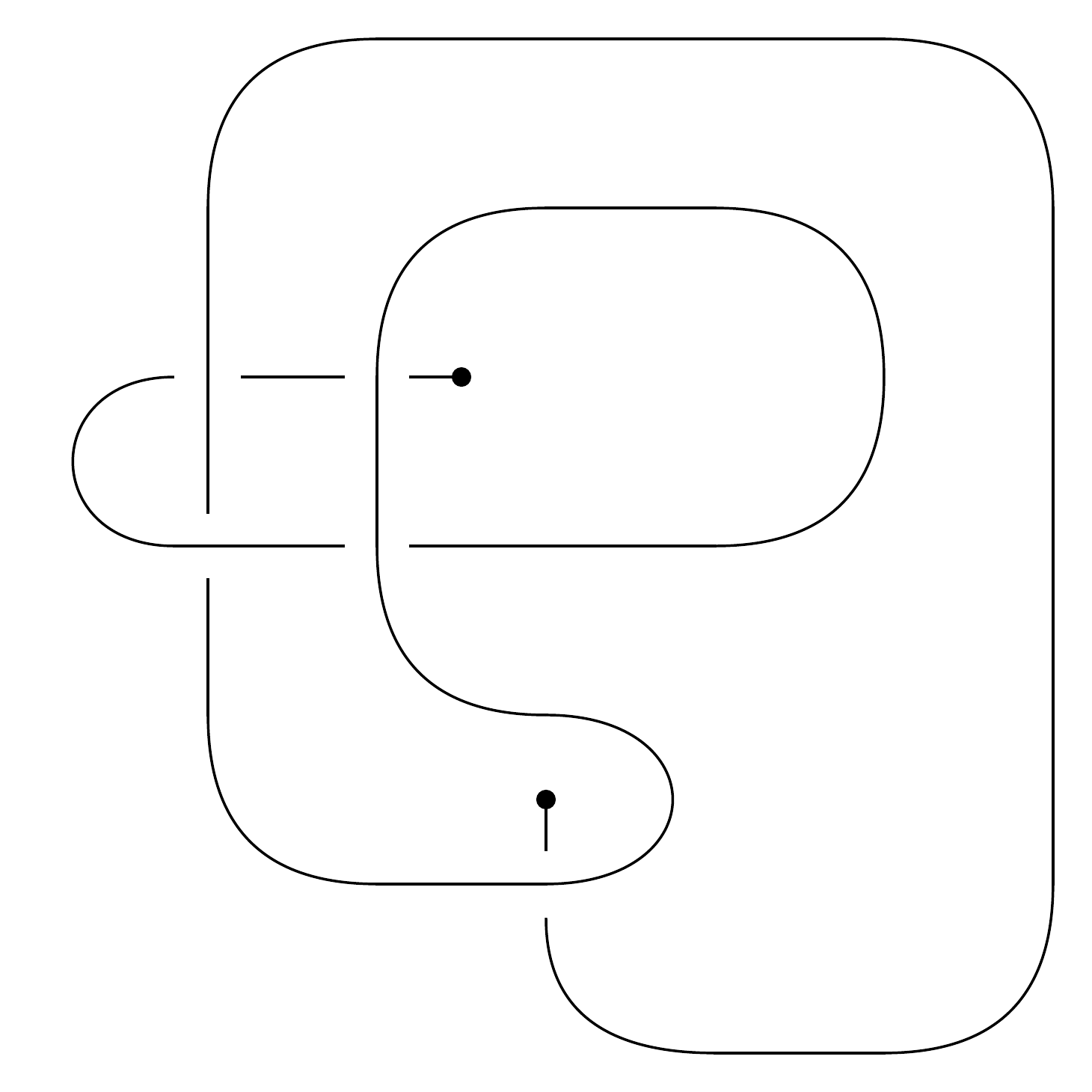}\\
\textcolor{black}{$5_{303}$}
\vspace{1cm}
\end{minipage}
\begin{minipage}[t]{.25\linewidth}
\centering
\includegraphics[width=0.9\textwidth,height=3.5cm,keepaspectratio]{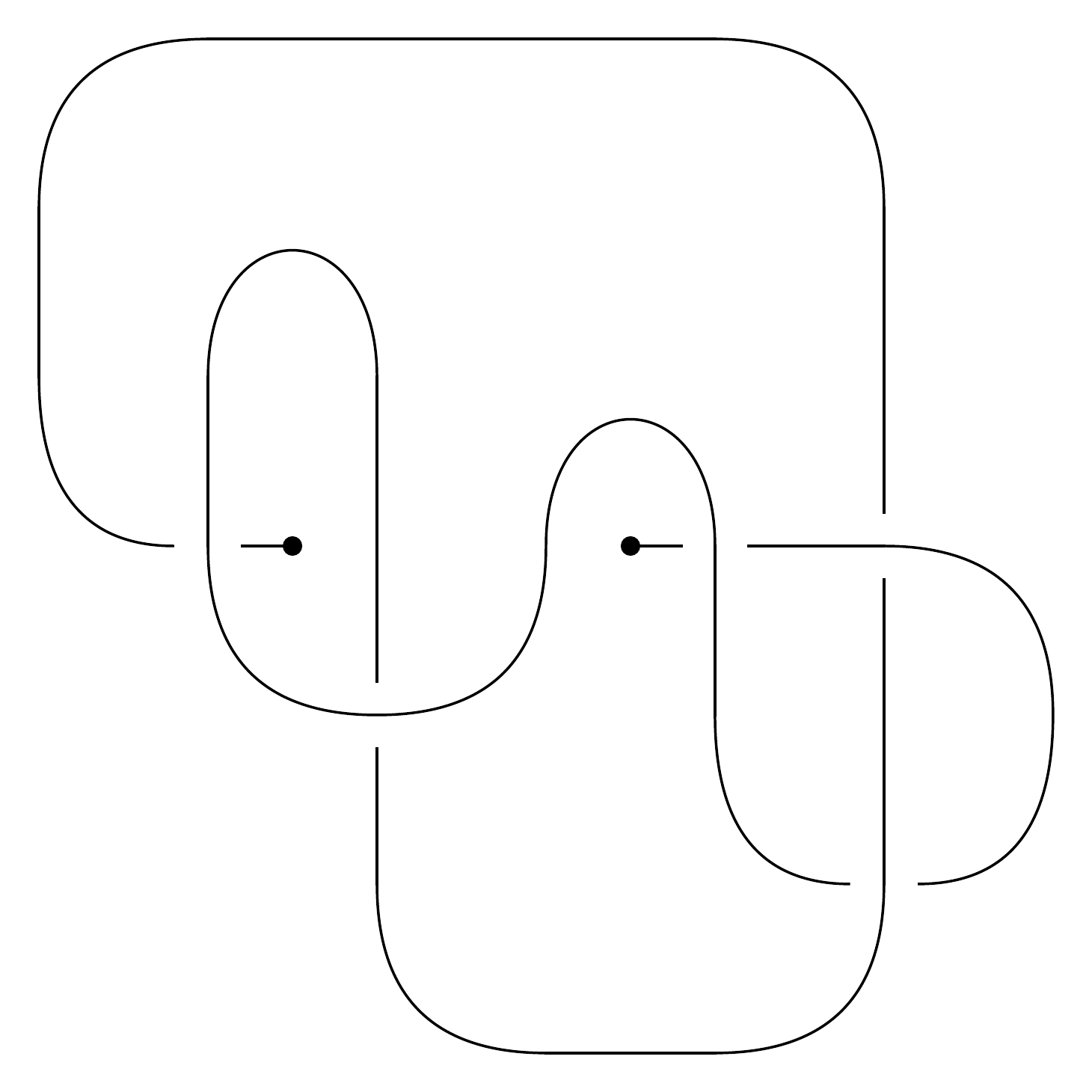}\\
\textcolor{black}{$5_{304}$}
\vspace{1cm}
\end{minipage}
\begin{minipage}[t]{.25\linewidth}
\centering
\includegraphics[width=0.9\textwidth,height=3.5cm,keepaspectratio]{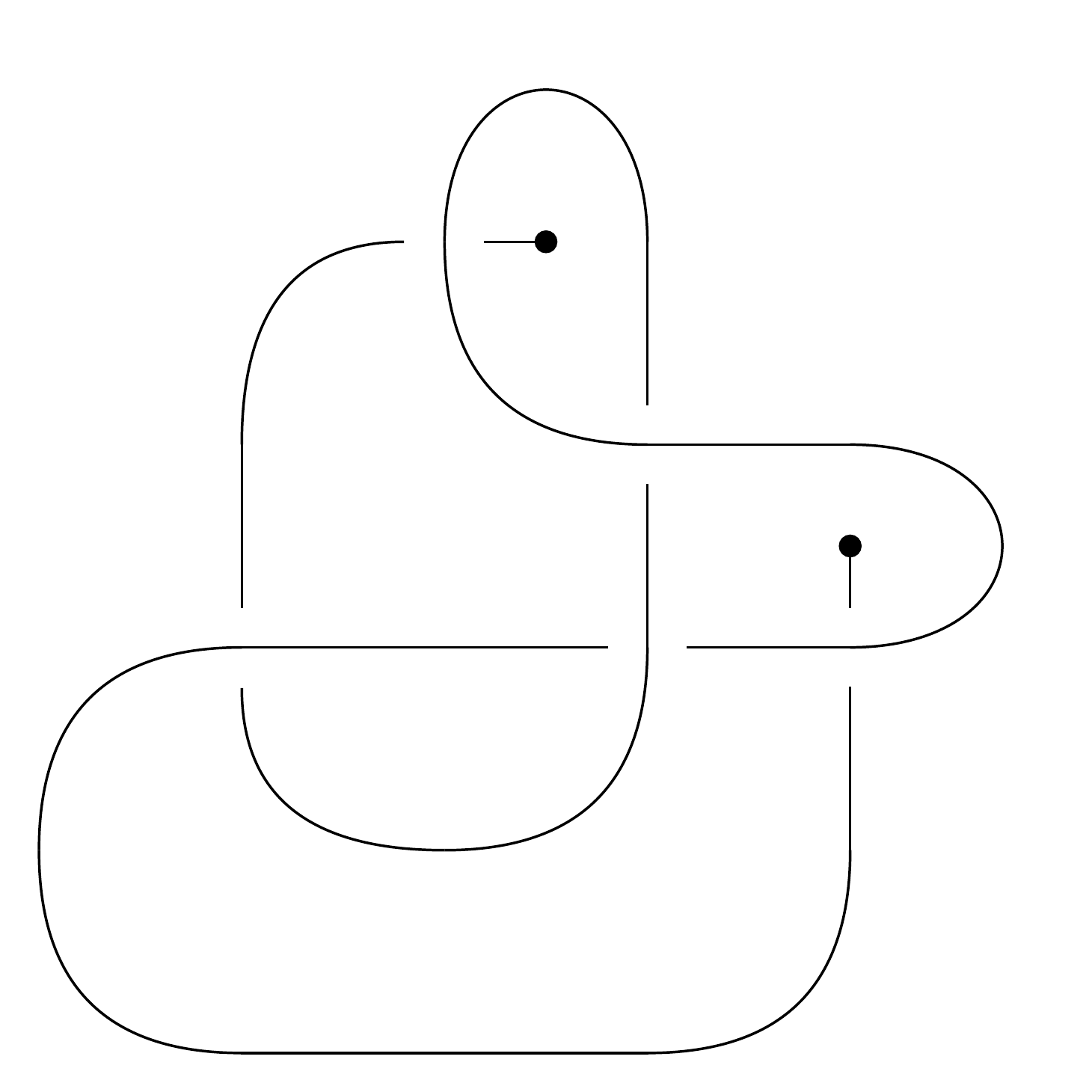}\\
\textcolor{black}{$5_{305}$}
\vspace{1cm}
\end{minipage}
\begin{minipage}[t]{.25\linewidth}
\centering
\includegraphics[width=0.9\textwidth,height=3.5cm,keepaspectratio]{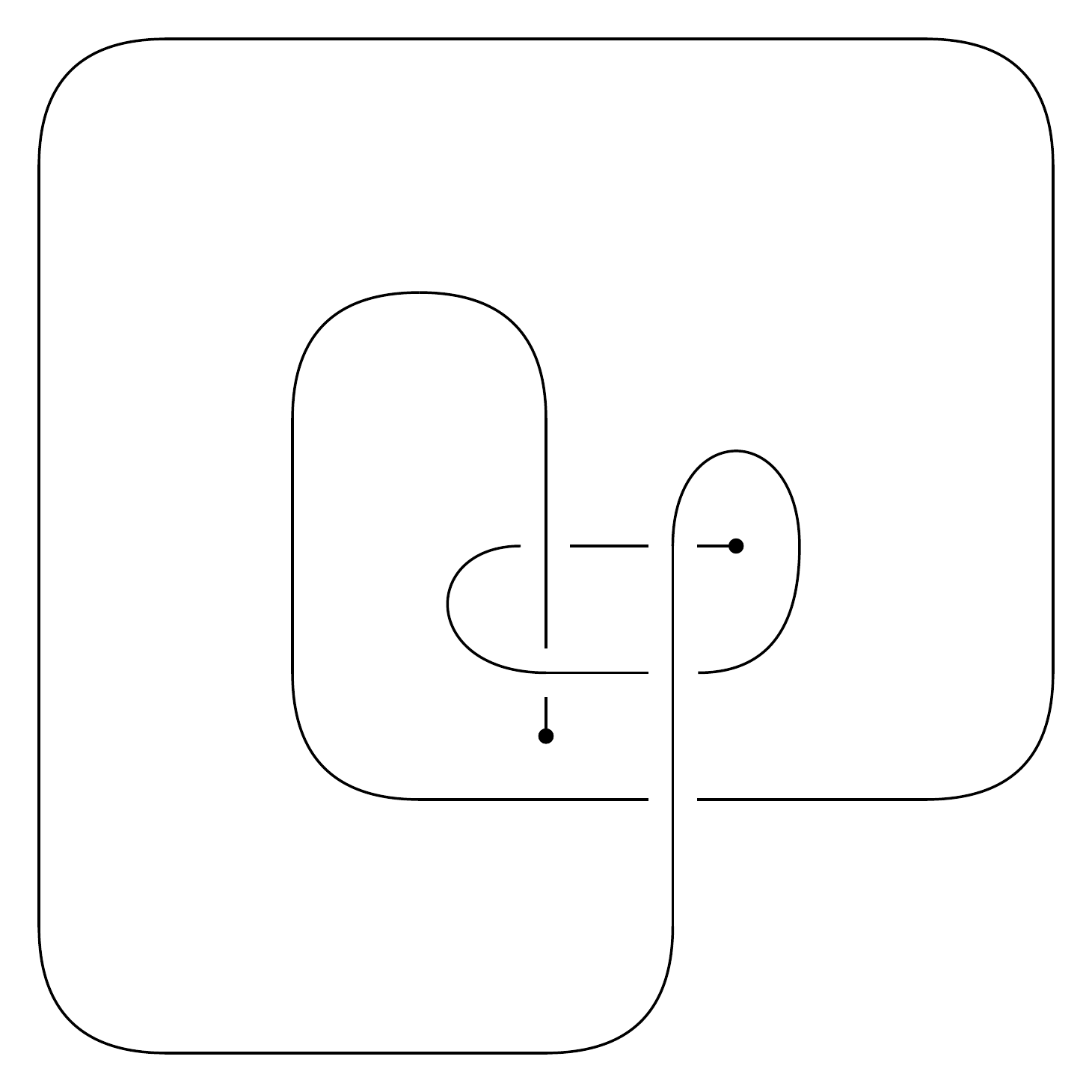}\\
\textcolor{black}{$5_{306}$}
\vspace{1cm}
\end{minipage}
\begin{minipage}[t]{.25\linewidth}
\centering
\includegraphics[width=0.9\textwidth,height=3.5cm,keepaspectratio]{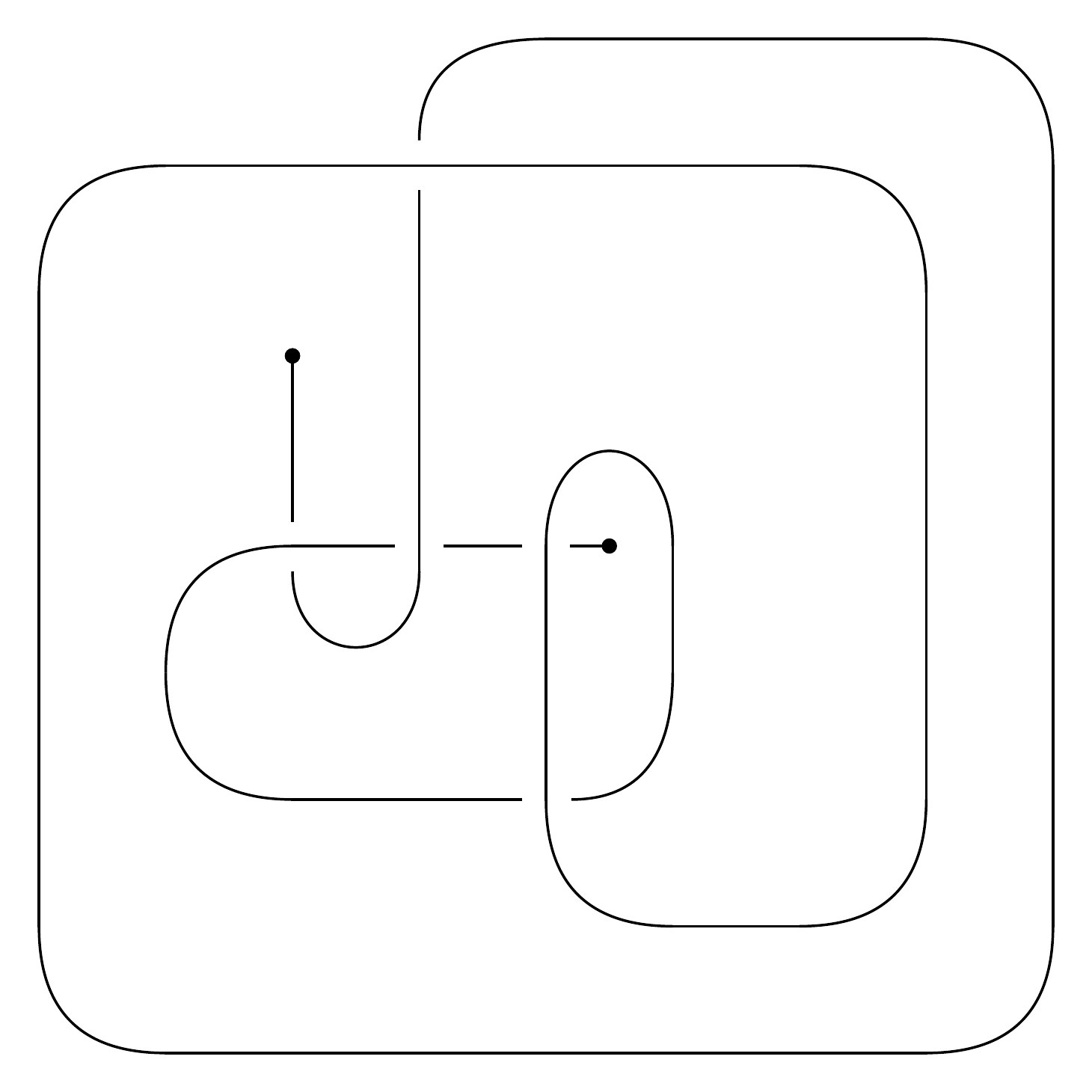}\\
\textcolor{black}{$5_{307}$}
\vspace{1cm}
\end{minipage}
\begin{minipage}[t]{.25\linewidth}
\centering
\includegraphics[width=0.9\textwidth,height=3.5cm,keepaspectratio]{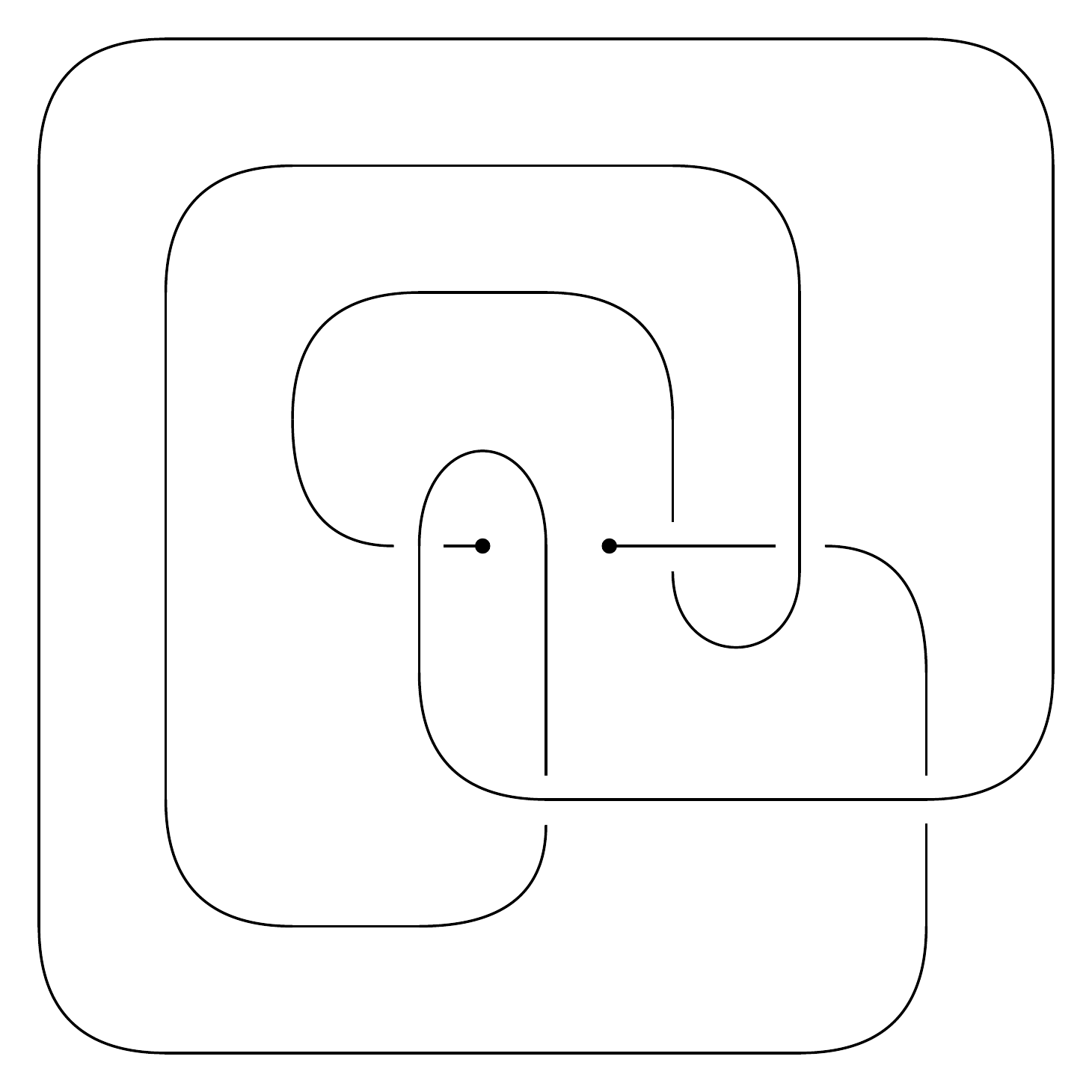}\\
\textcolor{black}{$5_{308}$}
\vspace{1cm}
\end{minipage}
\begin{minipage}[t]{.25\linewidth}
\centering
\includegraphics[width=0.9\textwidth,height=3.5cm,keepaspectratio]{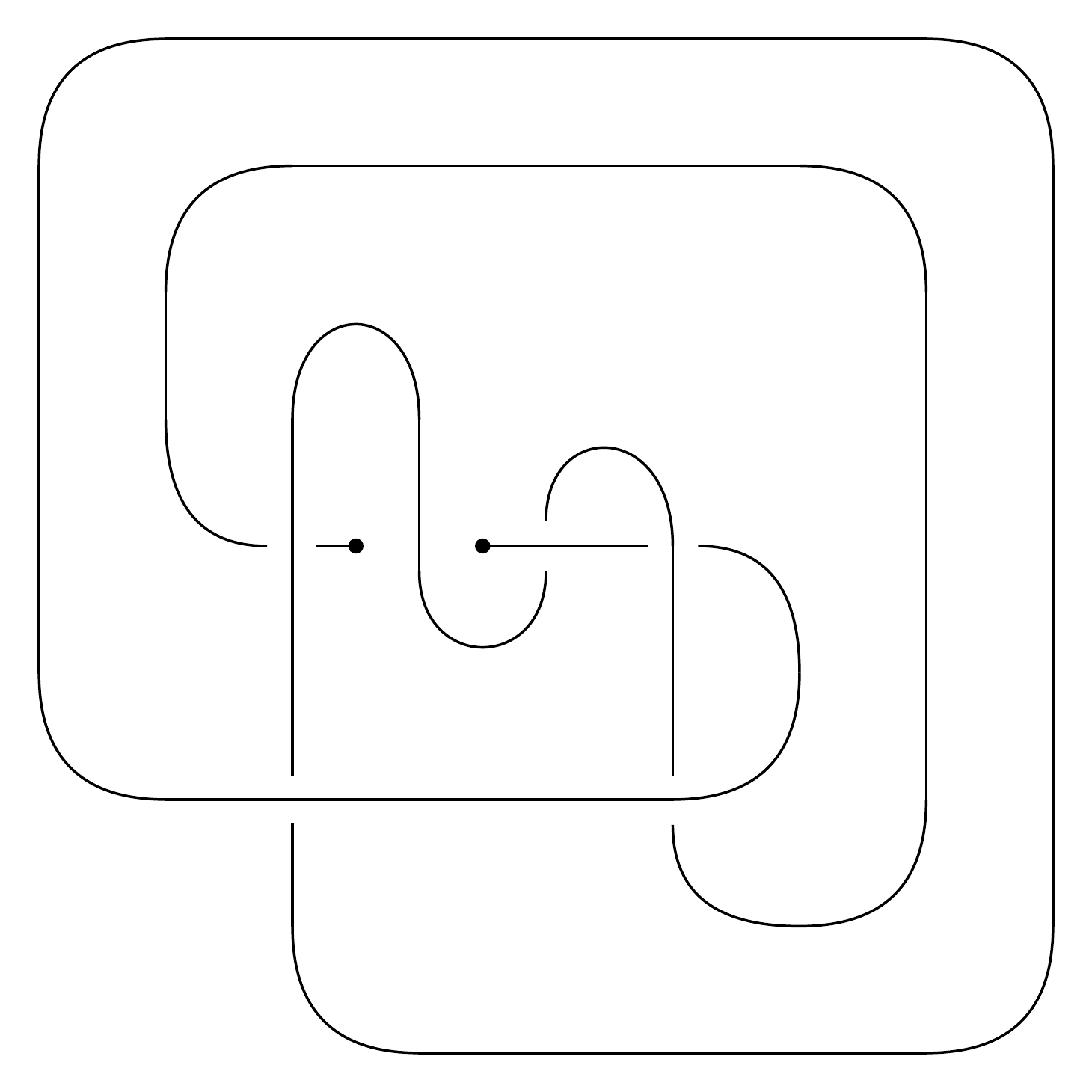}\\
\textcolor{black}{$5_{309}$}
\vspace{1cm}
\end{minipage}
\begin{minipage}[t]{.25\linewidth}
\centering
\includegraphics[width=0.9\textwidth,height=3.5cm,keepaspectratio]{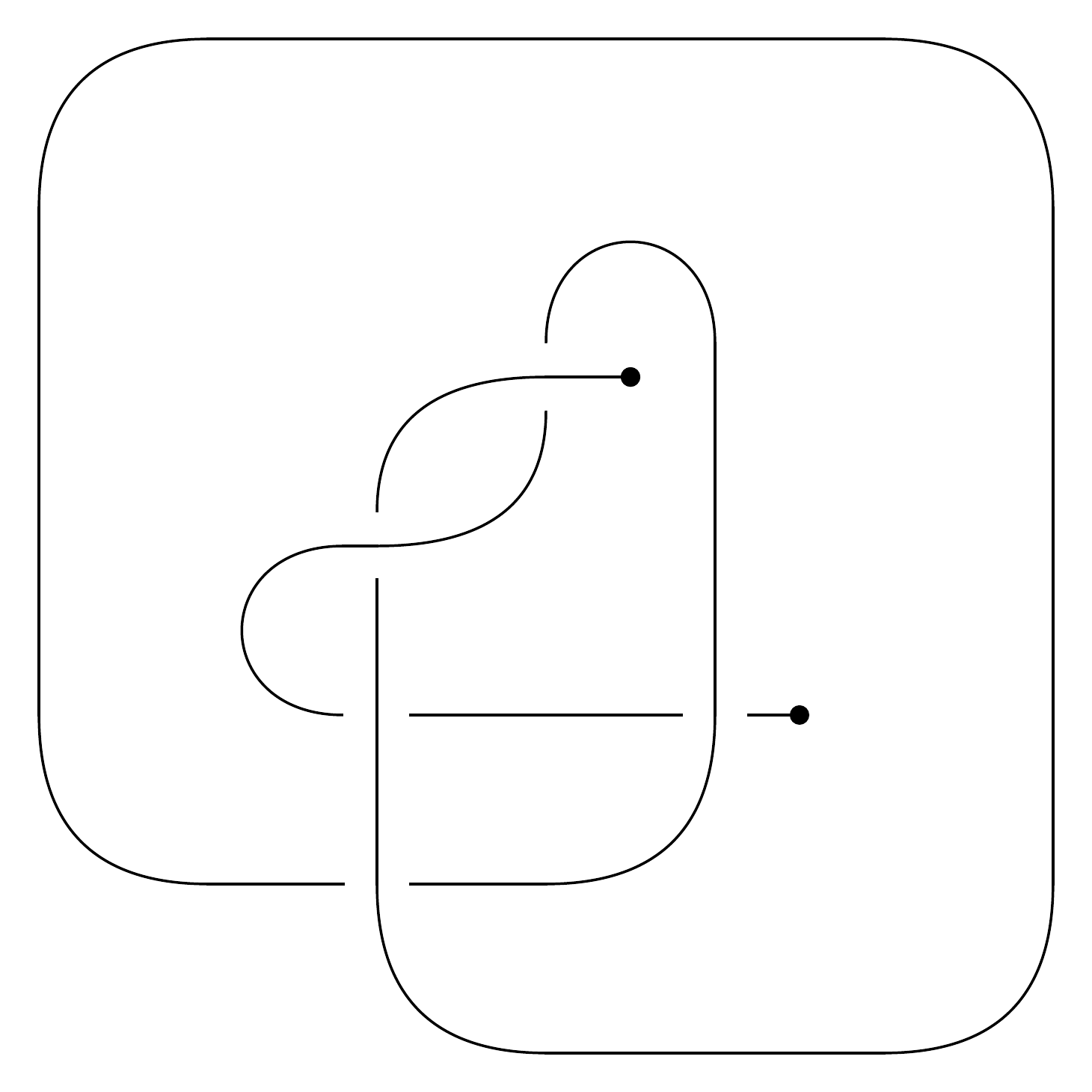}\\
\textcolor{black}{$5_{310}$}
\vspace{1cm}
\end{minipage}
\begin{minipage}[t]{.25\linewidth}
\centering
\includegraphics[width=0.9\textwidth,height=3.5cm,keepaspectratio]{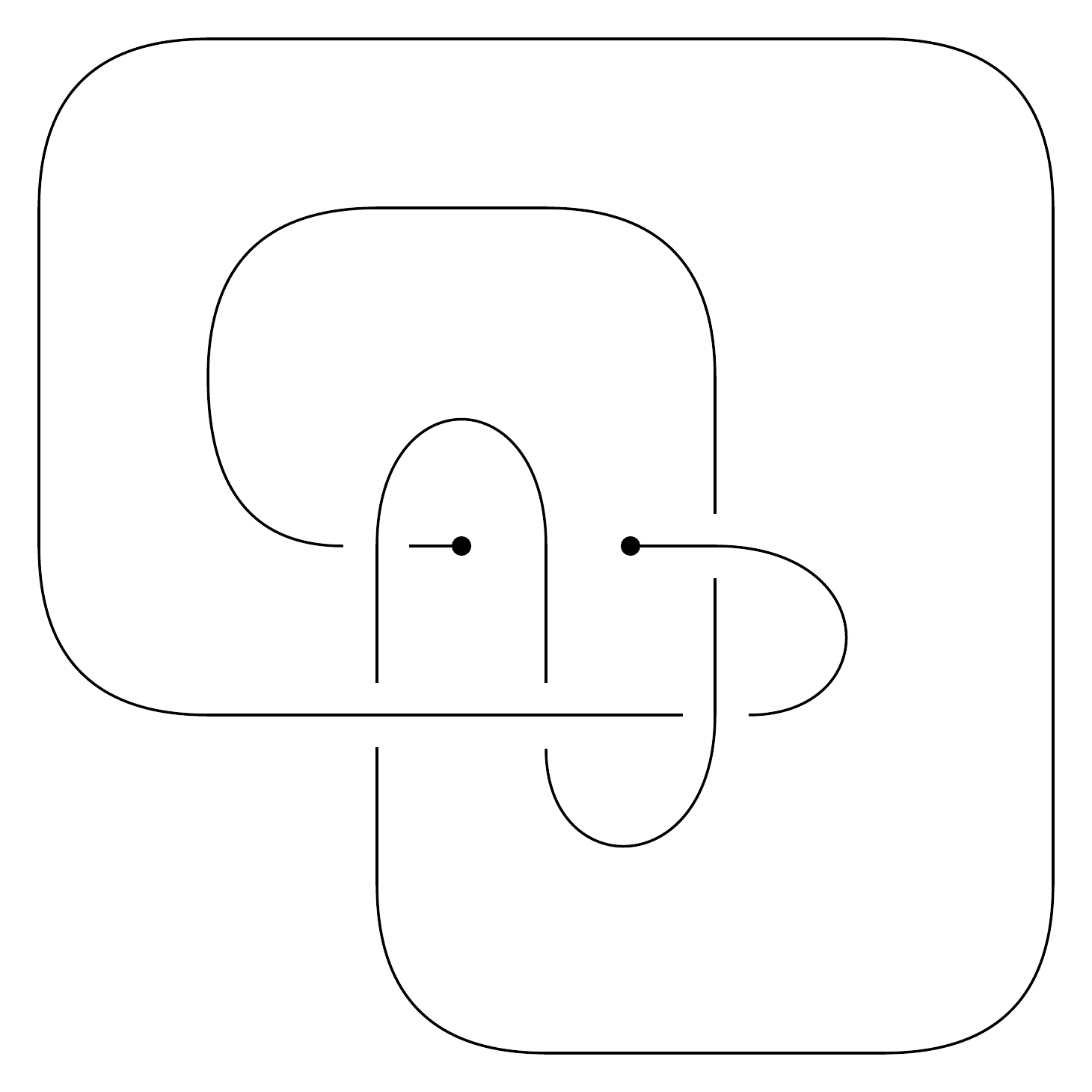}\\
\textcolor{black}{$5_{311}$}
\vspace{1cm}
\end{minipage}
\begin{minipage}[t]{.25\linewidth}
\centering
\includegraphics[width=0.9\textwidth,height=3.5cm,keepaspectratio]{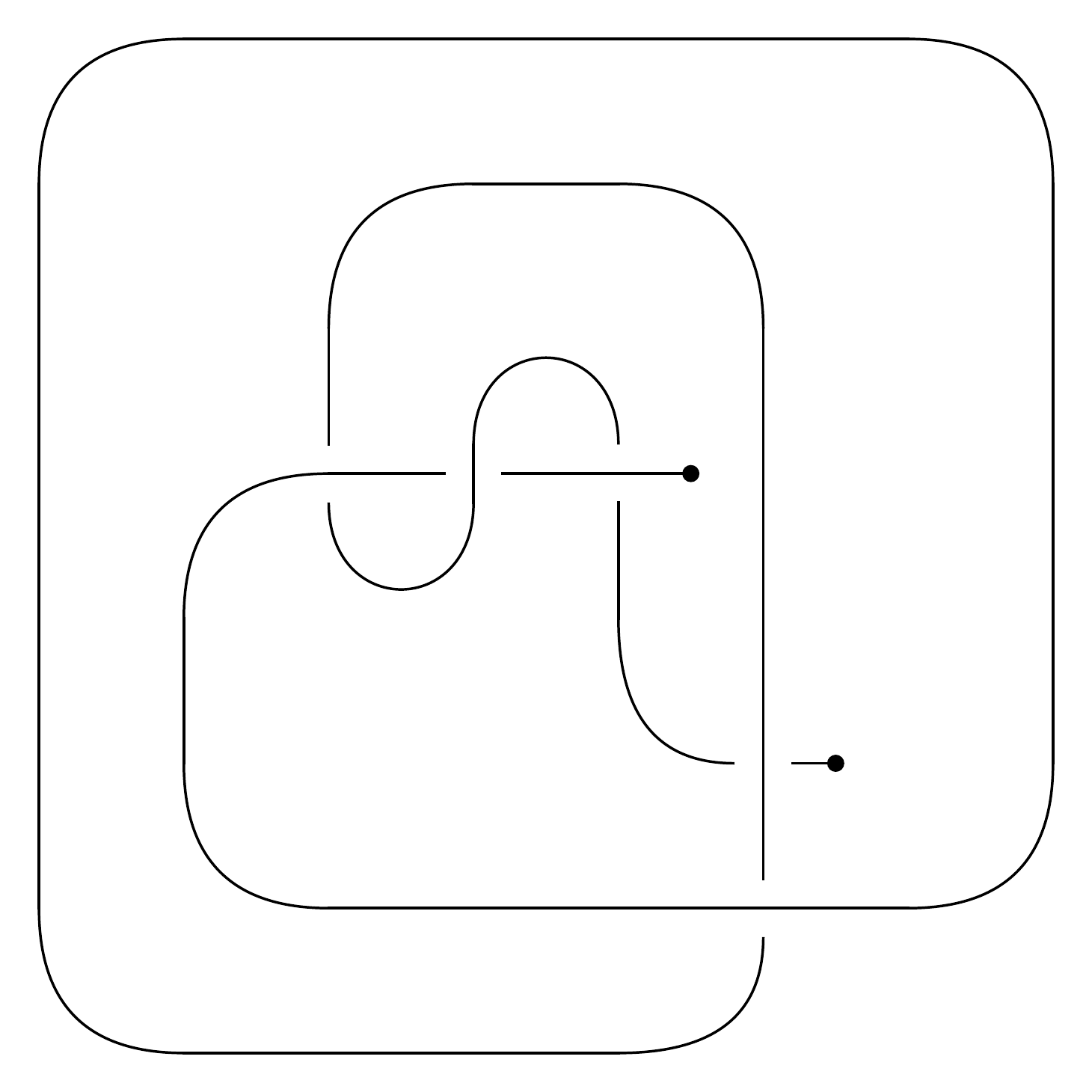}\\
\textcolor{black}{$5_{312}$}
\vspace{1cm}
\end{minipage}
\begin{minipage}[t]{.25\linewidth}
\centering
\includegraphics[width=0.9\textwidth,height=3.5cm,keepaspectratio]{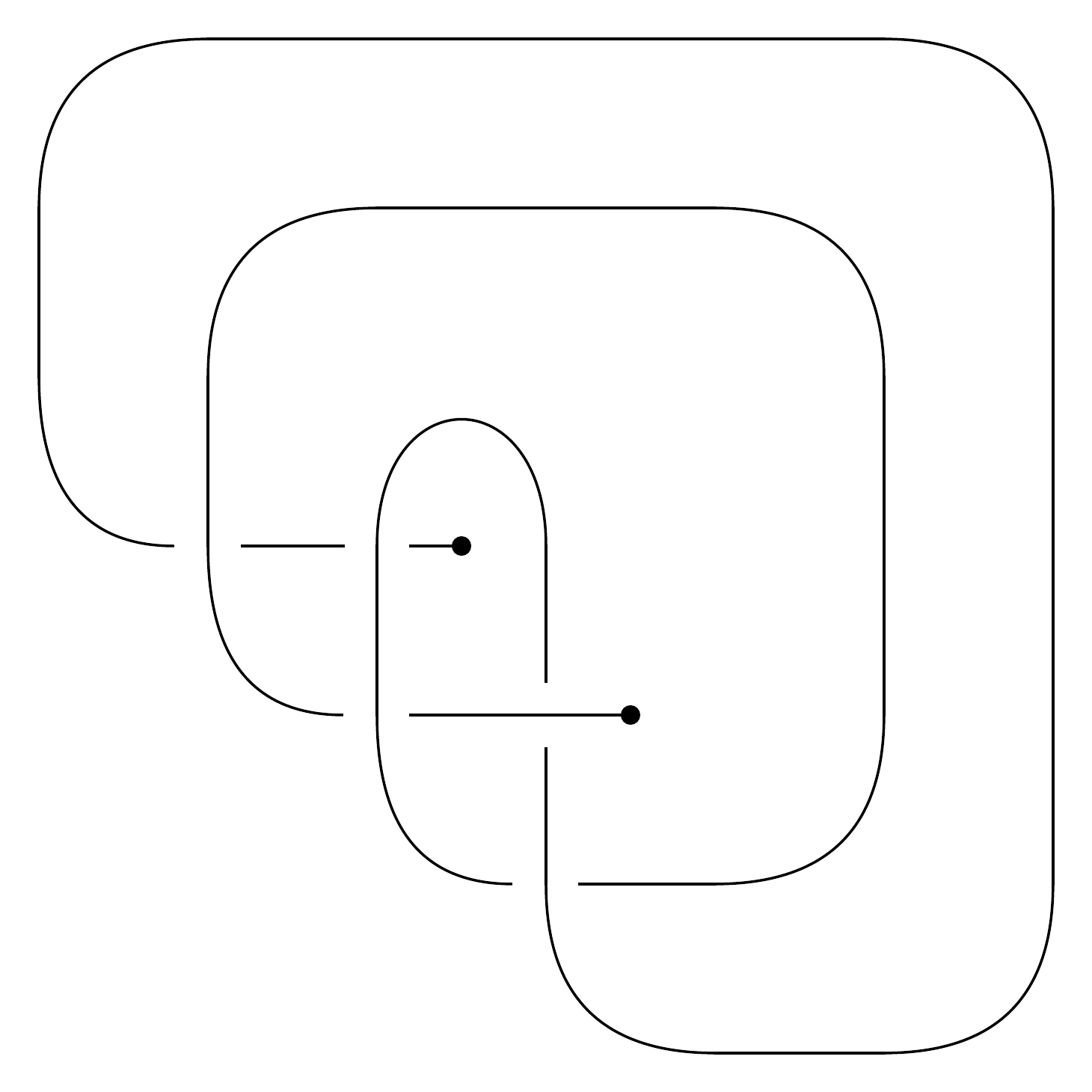}\\
\textcolor{black}{$5_{313}$}
\vspace{1cm}
\end{minipage}
\begin{minipage}[t]{.25\linewidth}
\centering
\includegraphics[width=0.9\textwidth,height=3.5cm,keepaspectratio]{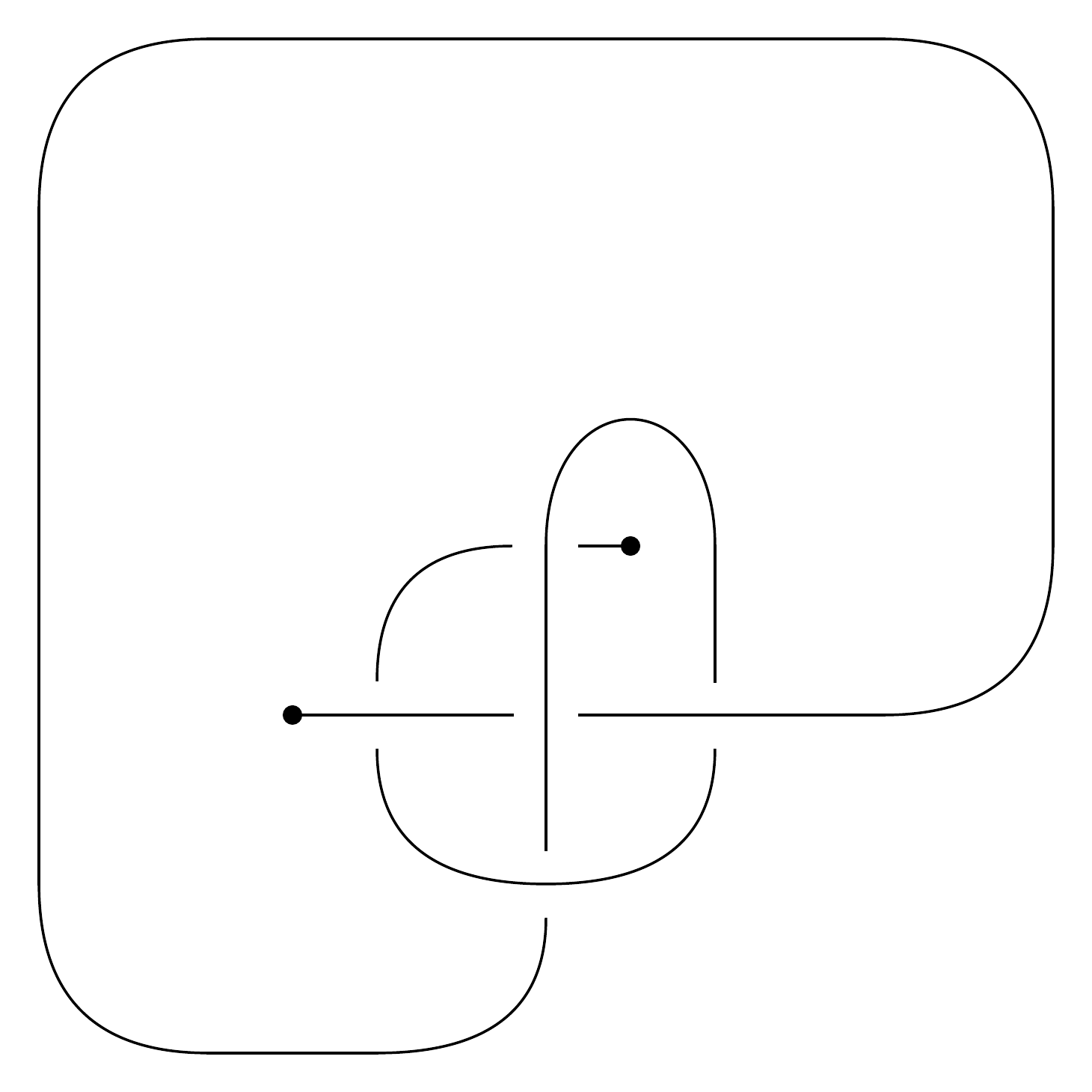}\\
\textcolor{black}{$5_{314}$}
\vspace{1cm}
\end{minipage}
\begin{minipage}[t]{.25\linewidth}
\centering
\includegraphics[width=0.9\textwidth,height=3.5cm,keepaspectratio]{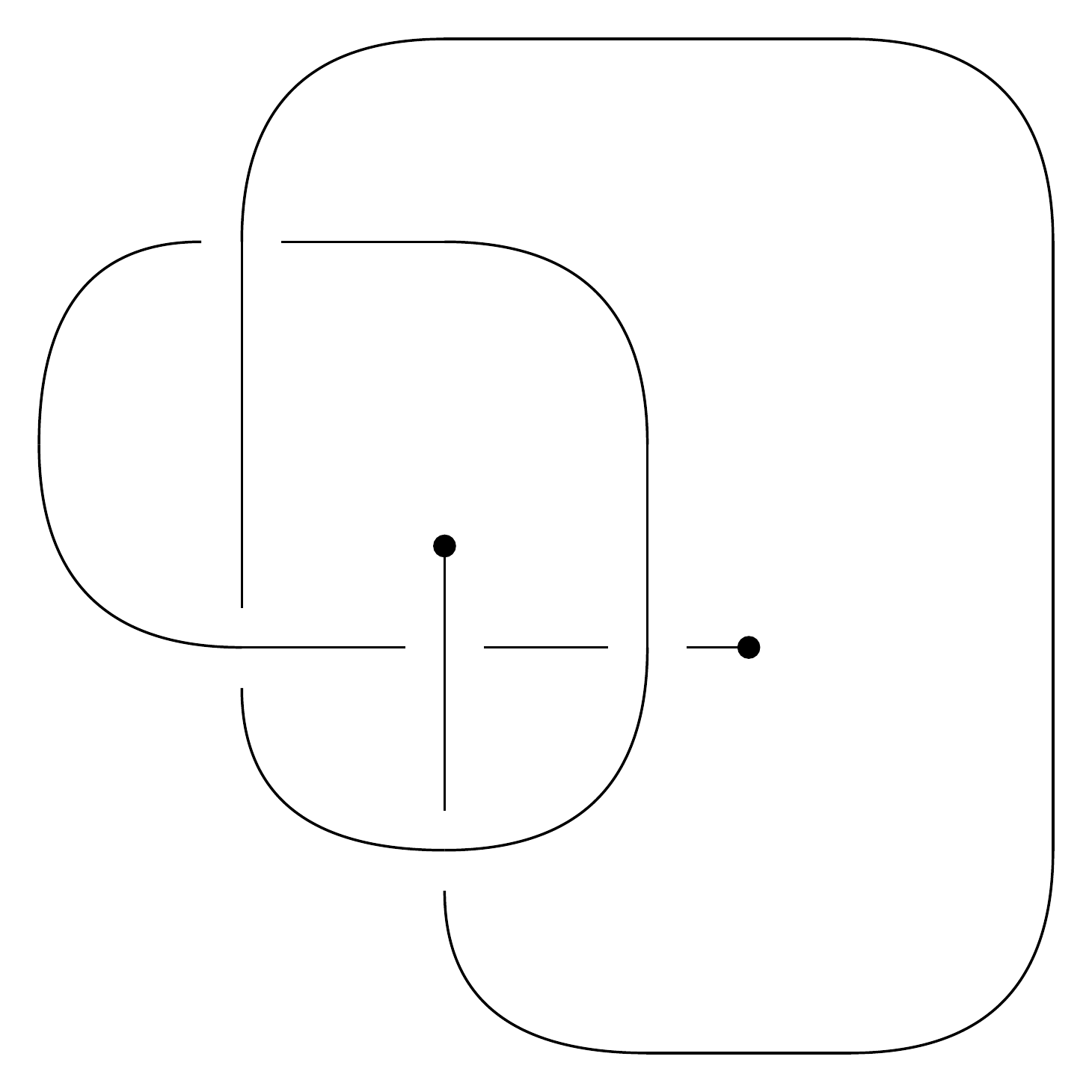}\\
\textcolor{black}{$5_{315}$}
\vspace{1cm}
\end{minipage}
\begin{minipage}[t]{.25\linewidth}
\centering
\includegraphics[width=0.9\textwidth,height=3.5cm,keepaspectratio]{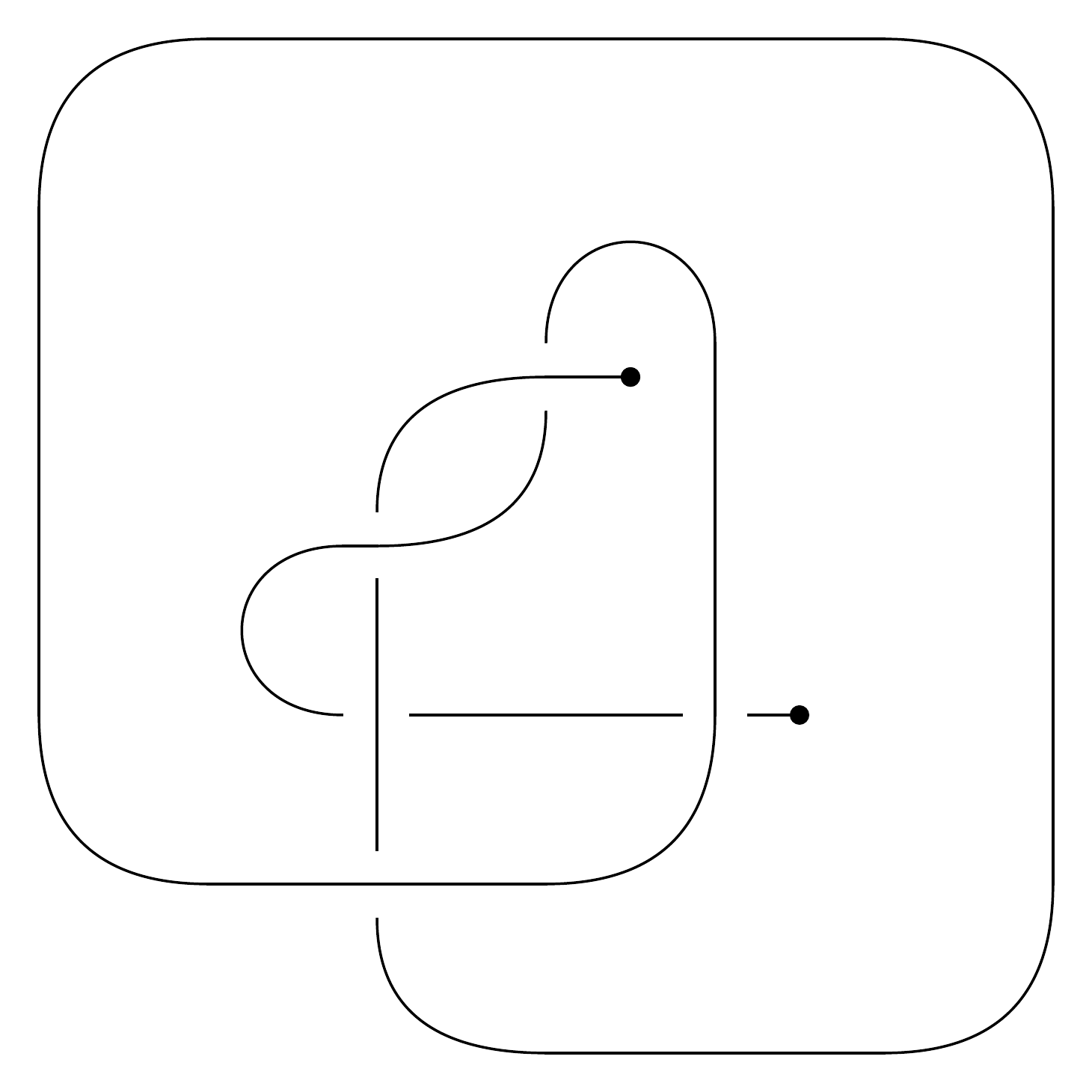}\\
\textcolor{black}{$5_{316}$}
\vspace{1cm}
\end{minipage}
\begin{minipage}[t]{.25\linewidth}
\centering
\includegraphics[width=0.9\textwidth,height=3.5cm,keepaspectratio]{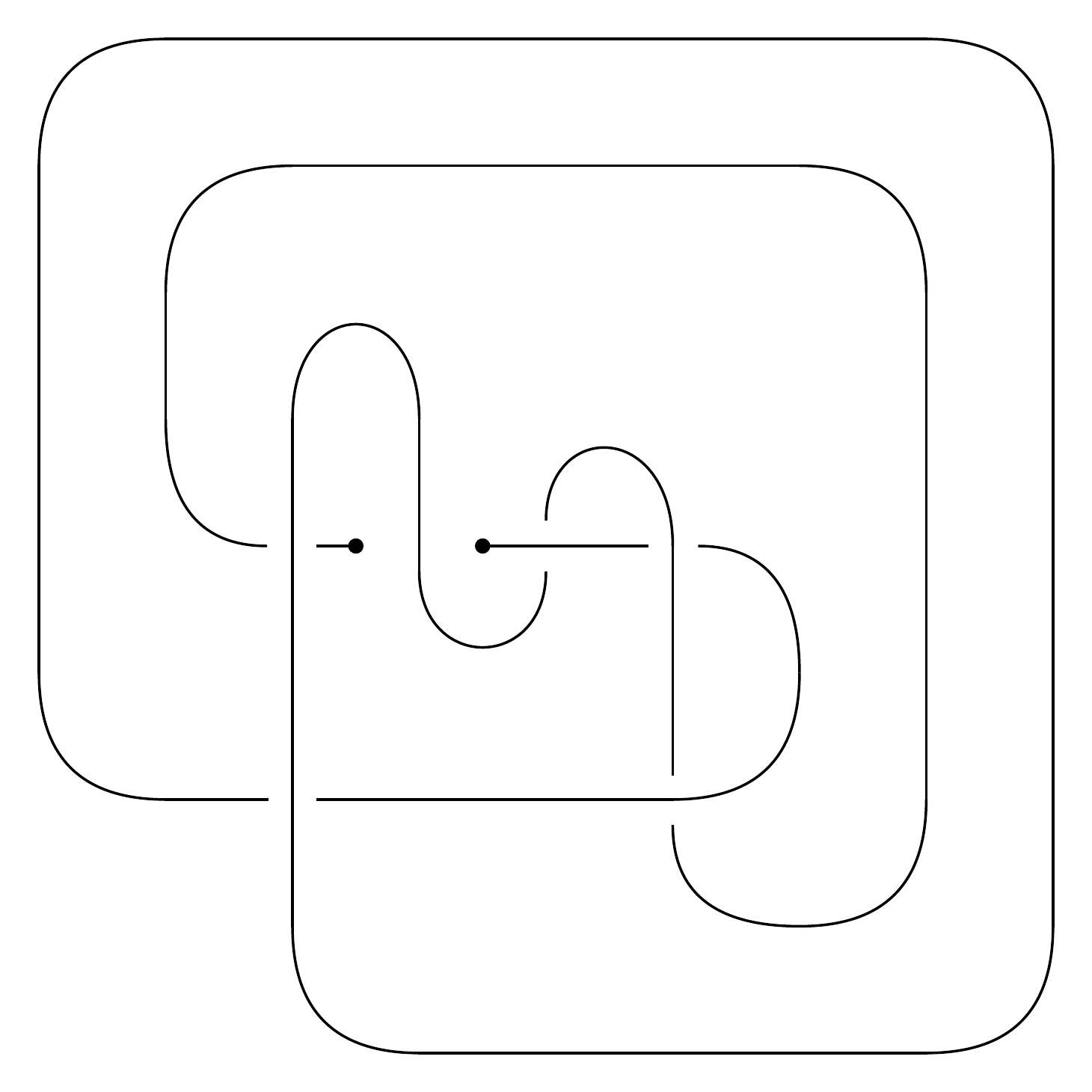}\\
\textcolor{black}{$5_{317}$}
\vspace{1cm}
\end{minipage}
\begin{minipage}[t]{.25\linewidth}
\centering
\includegraphics[width=0.9\textwidth,height=3.5cm,keepaspectratio]{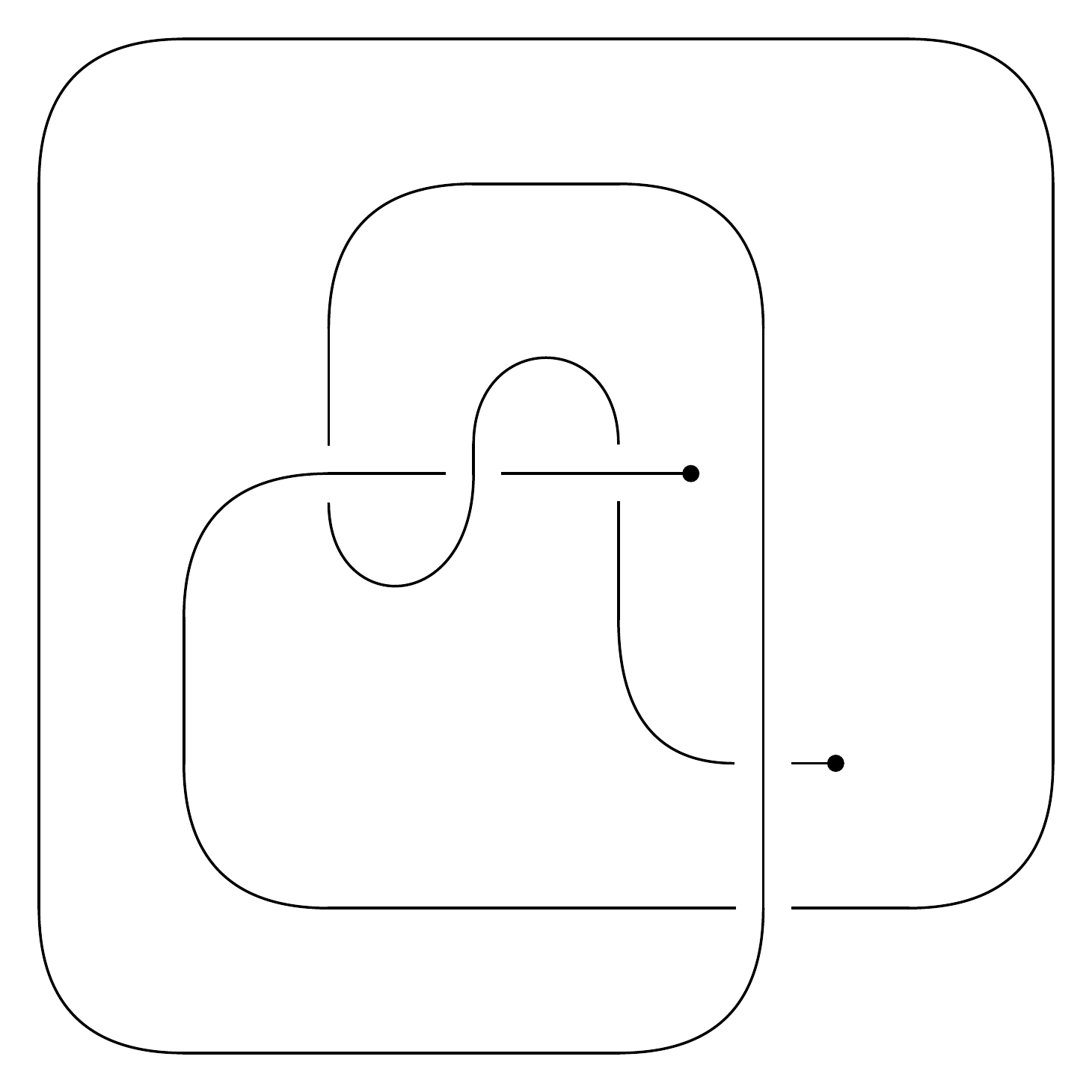}\\
\textcolor{black}{$5_{318}$}
\vspace{1cm}
\end{minipage}
\begin{minipage}[t]{.25\linewidth}
\centering
\includegraphics[width=0.9\textwidth,height=3.5cm,keepaspectratio]{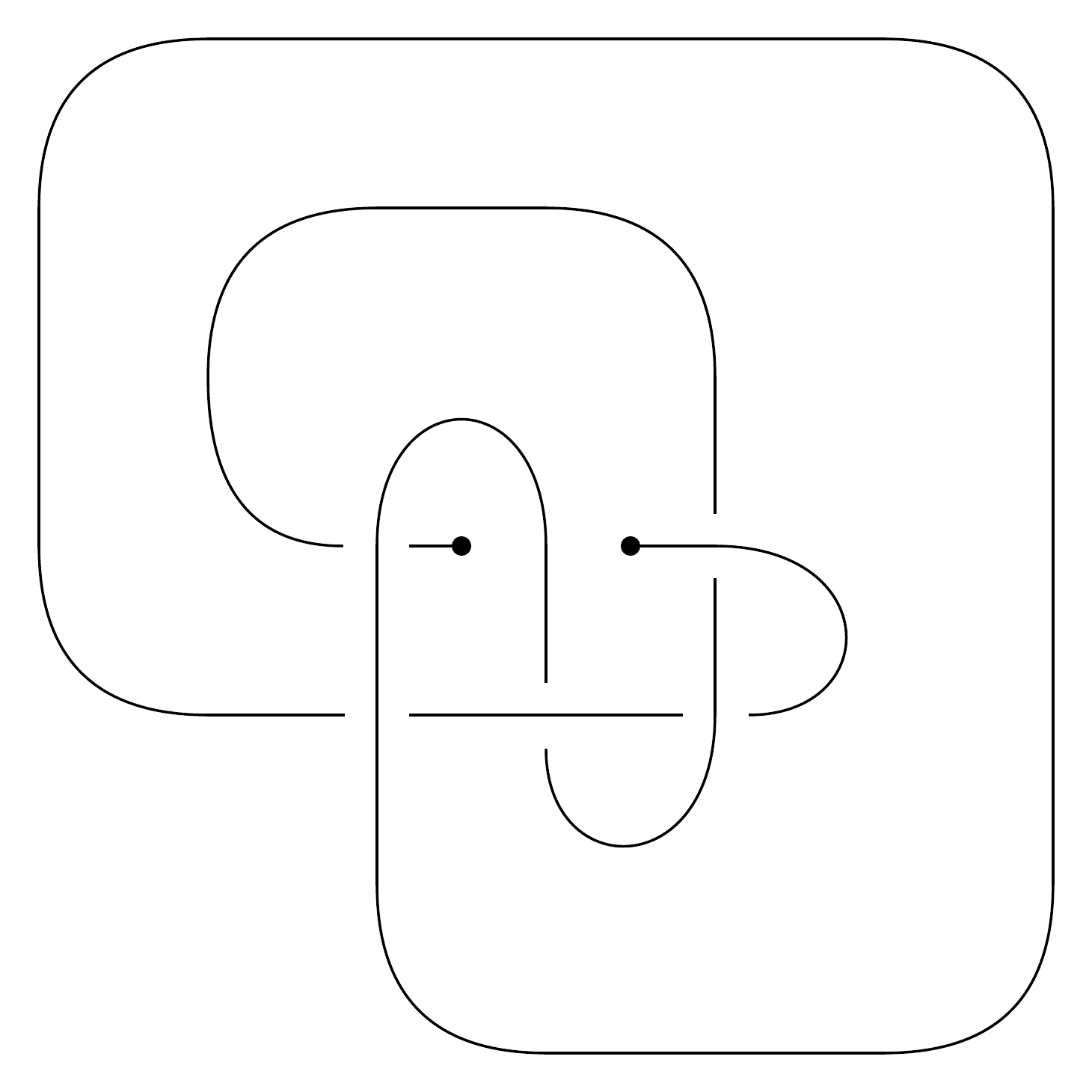}\\
\textcolor{black}{$5_{319}$}
\vspace{1cm}
\end{minipage}
\begin{minipage}[t]{.25\linewidth}
\centering
\includegraphics[width=0.9\textwidth,height=3.5cm,keepaspectratio]{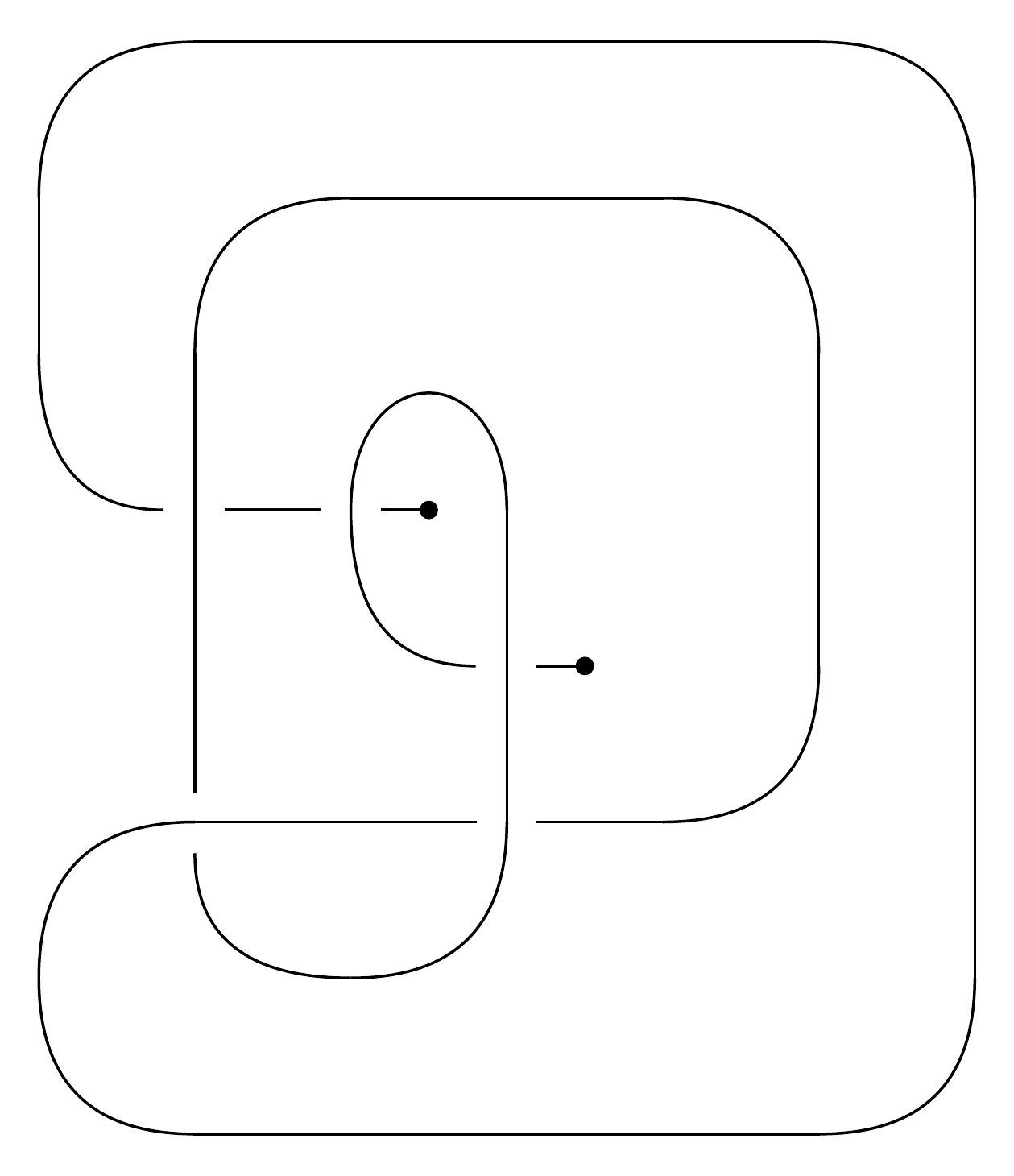}\\
\textcolor{black}{$5_{320}$}
\vspace{1cm}
\end{minipage}
\begin{minipage}[t]{.25\linewidth}
\centering
\includegraphics[width=0.9\textwidth,height=3.5cm,keepaspectratio]{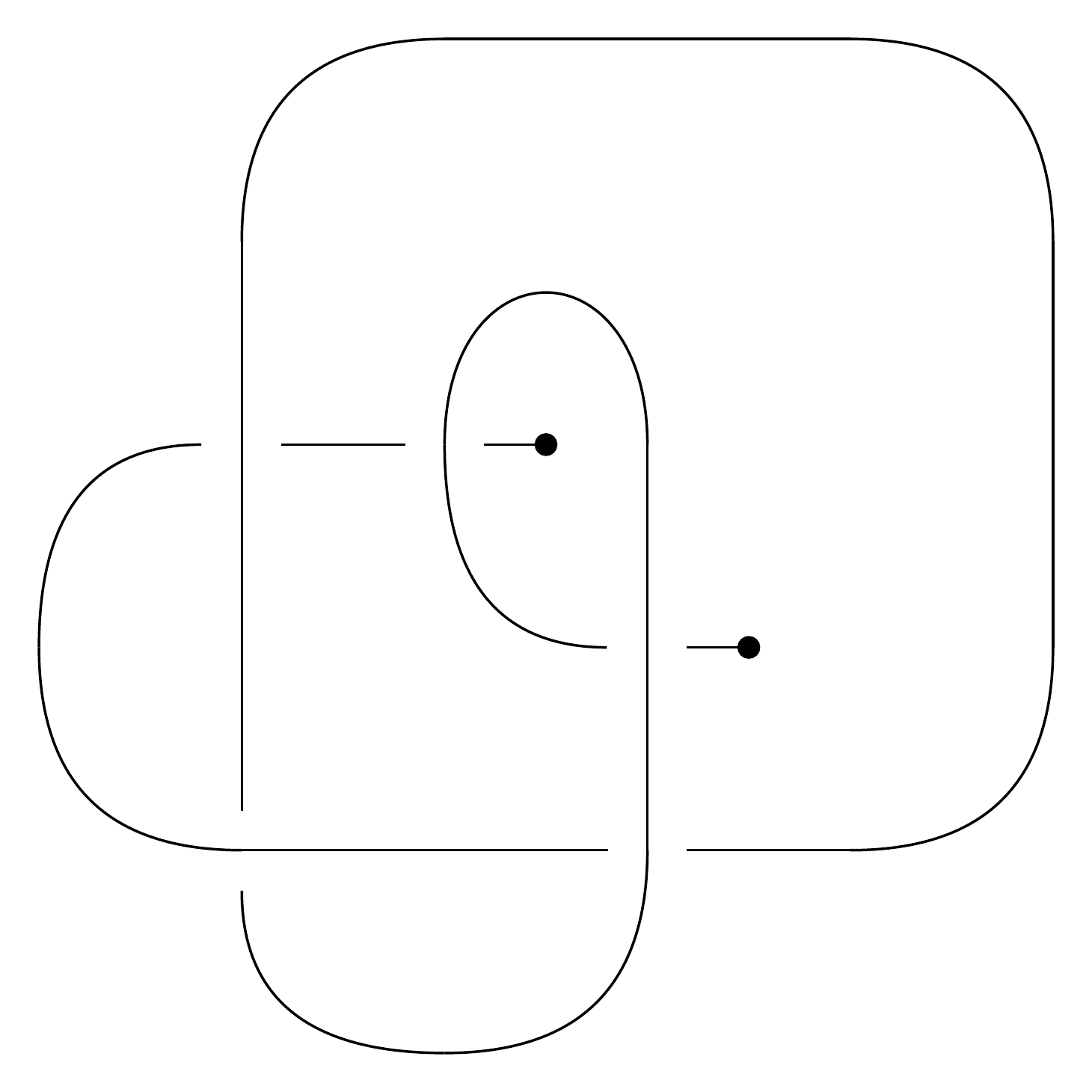}\\
\textcolor{black}{$5_{321}$}
\vspace{1cm}
\end{minipage}
\begin{minipage}[t]{.25\linewidth}
\centering
\includegraphics[width=0.9\textwidth,height=3.5cm,keepaspectratio]{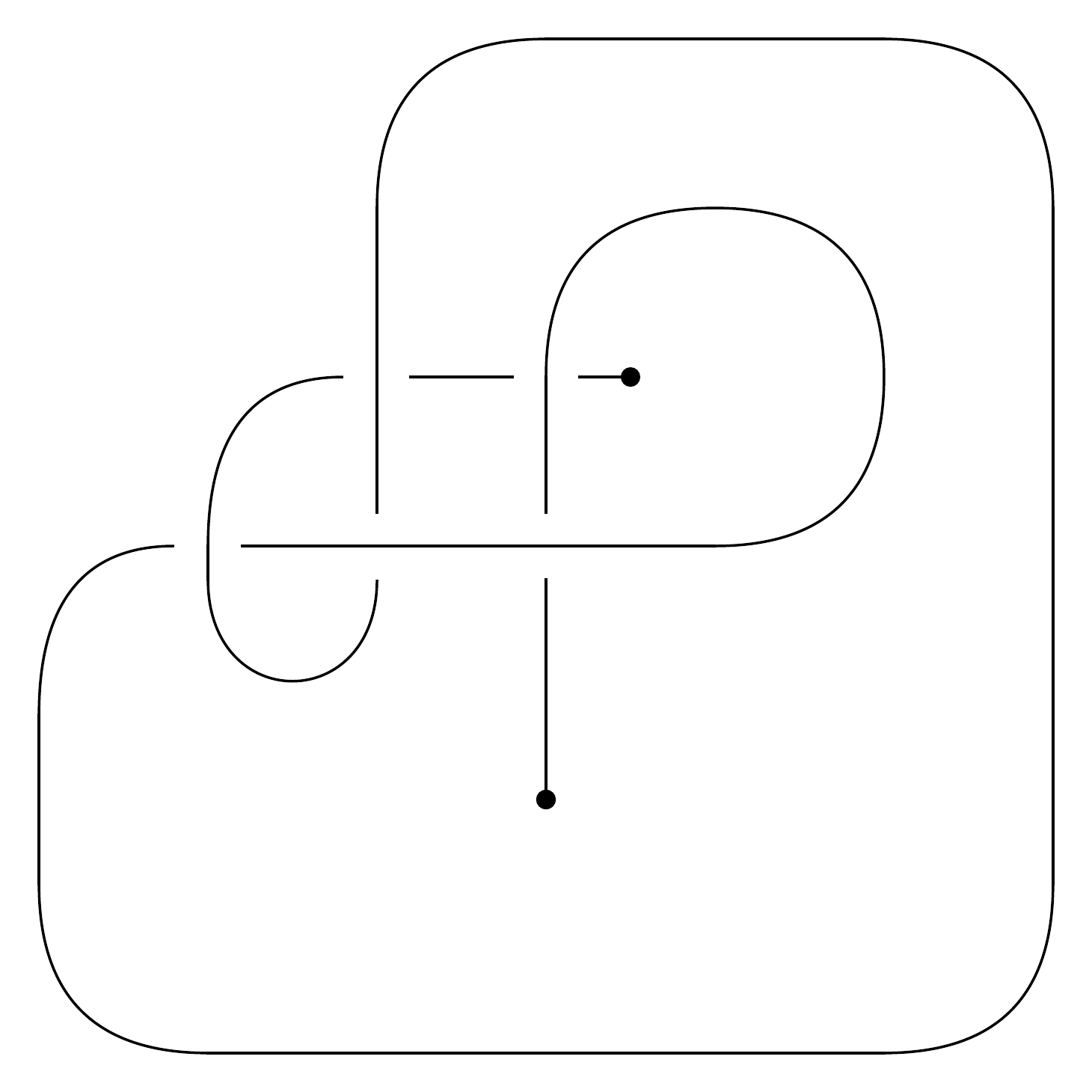}\\
\textcolor{black}{$5_{322}$}
\vspace{1cm}
\end{minipage}
\begin{minipage}[t]{.25\linewidth}
\centering
\includegraphics[width=0.9\textwidth,height=3.5cm,keepaspectratio]{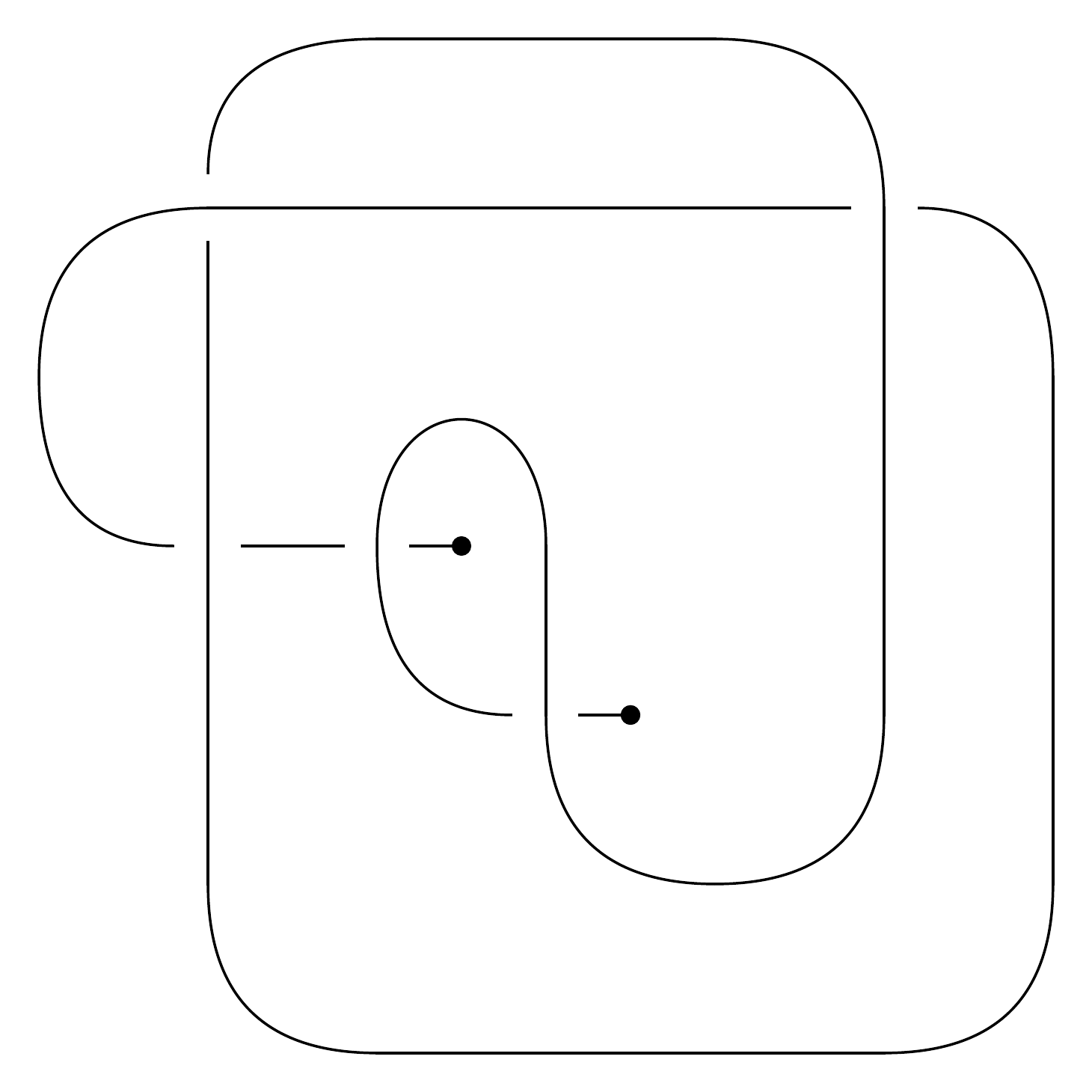}\\
\textcolor{black}{$5_{323}$}
\vspace{1cm}
\end{minipage}
\begin{minipage}[t]{.25\linewidth}
\centering
\includegraphics[width=0.9\textwidth,height=3.5cm,keepaspectratio]{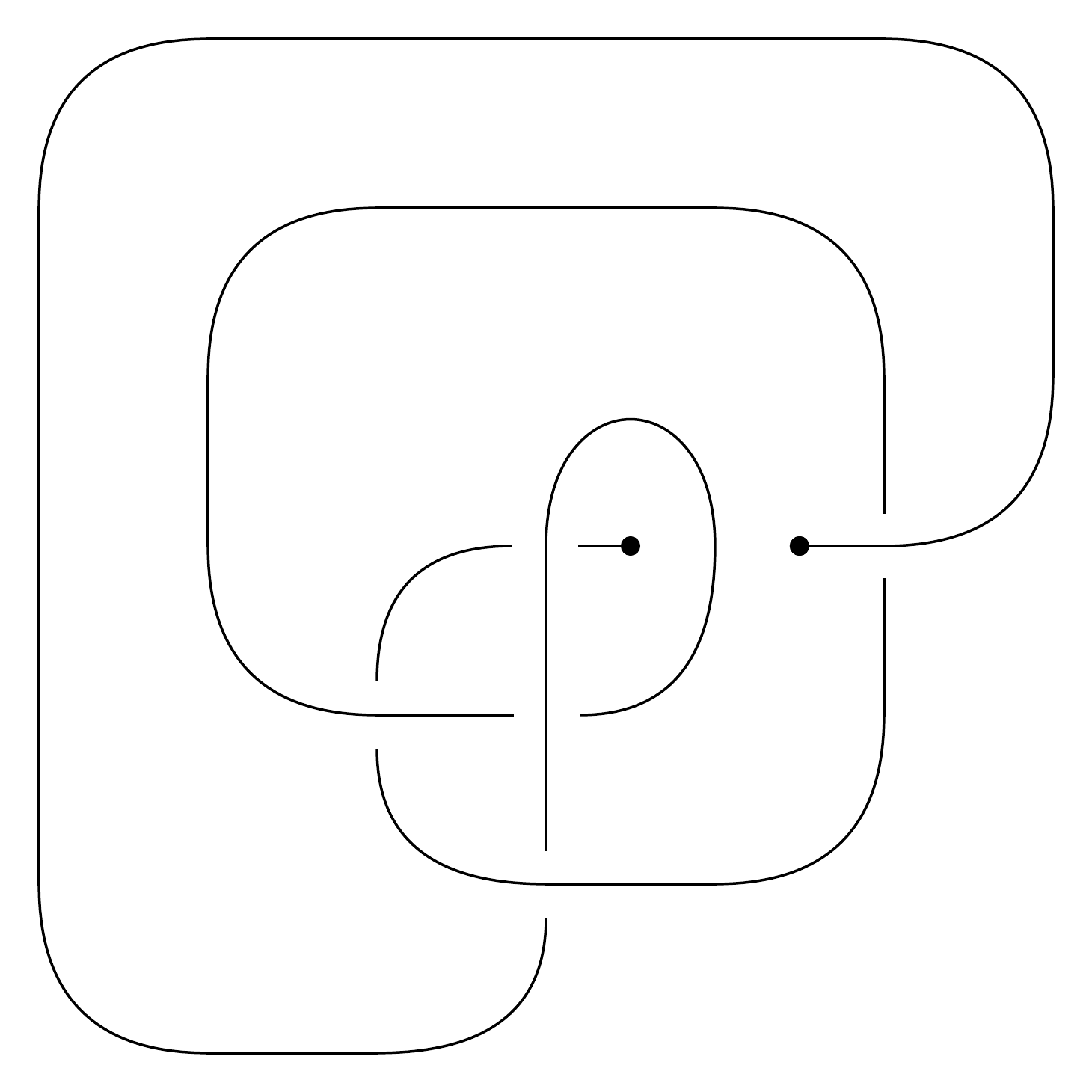}\\
\textcolor{black}{$5_{324}$}
\vspace{1cm}
\end{minipage}
\begin{minipage}[t]{.25\linewidth}
\centering
\includegraphics[width=0.9\textwidth,height=3.5cm,keepaspectratio]{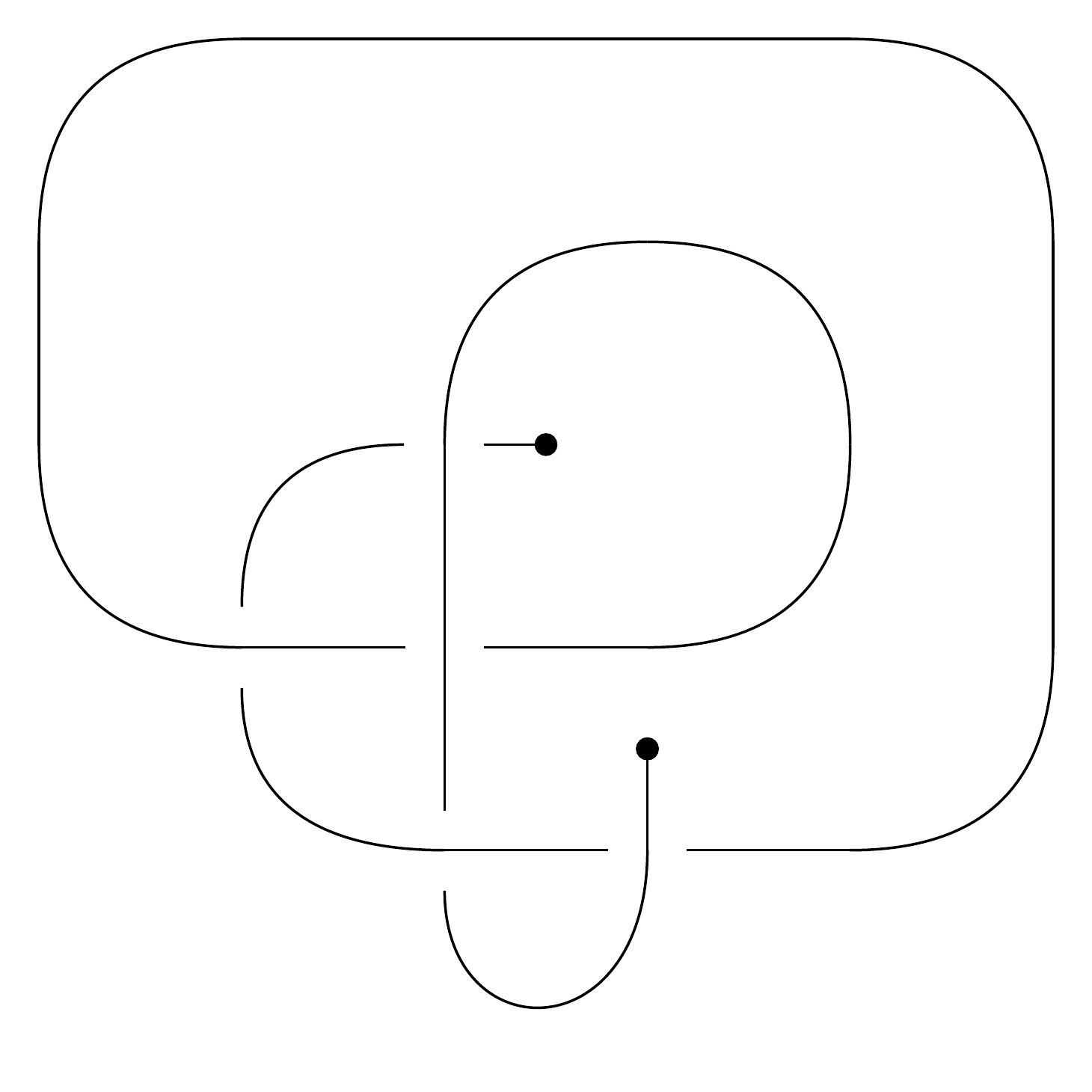}\\
\textcolor{black}{$5_{325}$}
\vspace{1cm}
\end{minipage}
\begin{minipage}[t]{.25\linewidth}
\centering
\includegraphics[width=0.9\textwidth,height=3.5cm,keepaspectratio]{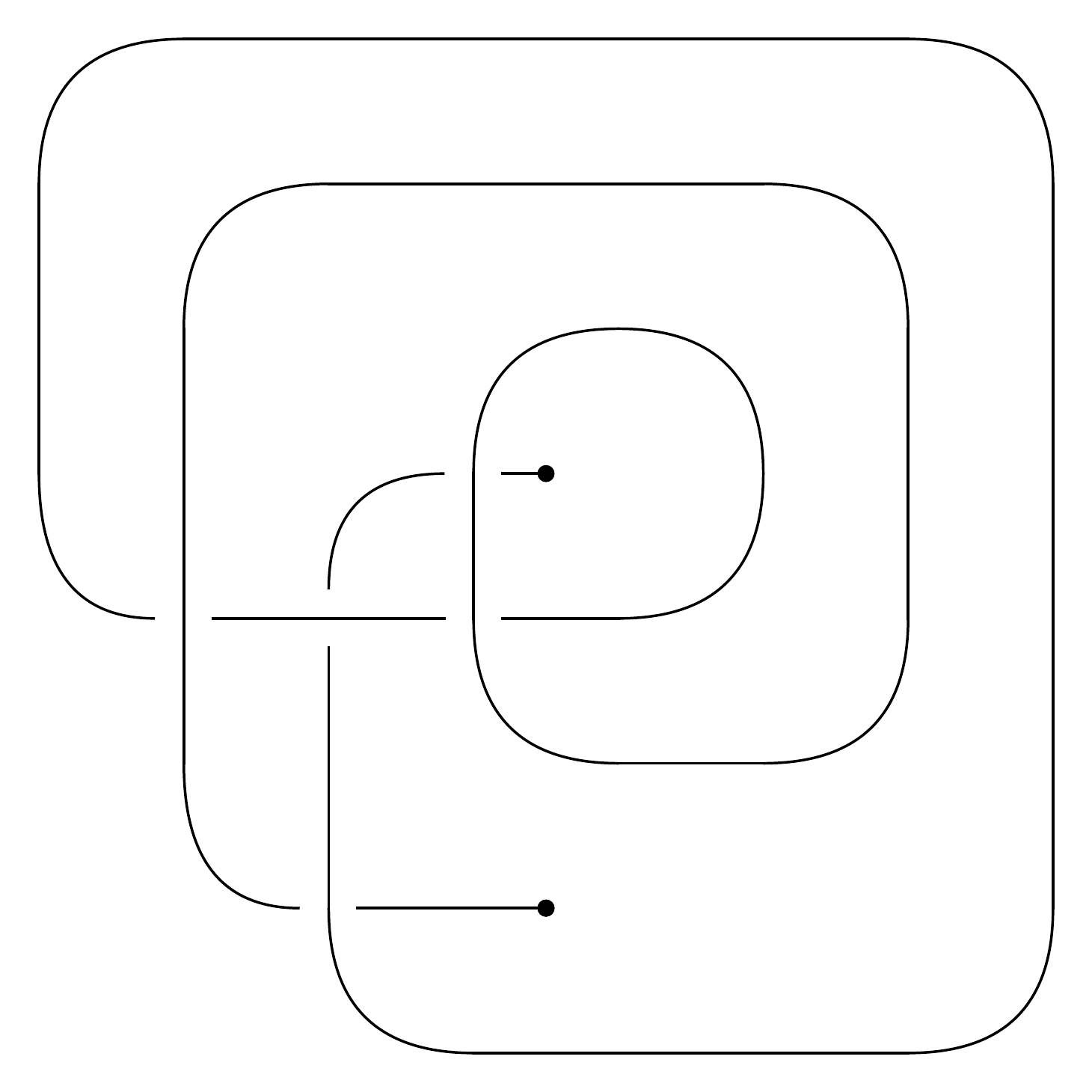}\\
\textcolor{black}{$5_{326}$}
\vspace{1cm}
\end{minipage}
\begin{minipage}[t]{.25\linewidth}
\centering
\includegraphics[width=0.9\textwidth,height=3.5cm,keepaspectratio]{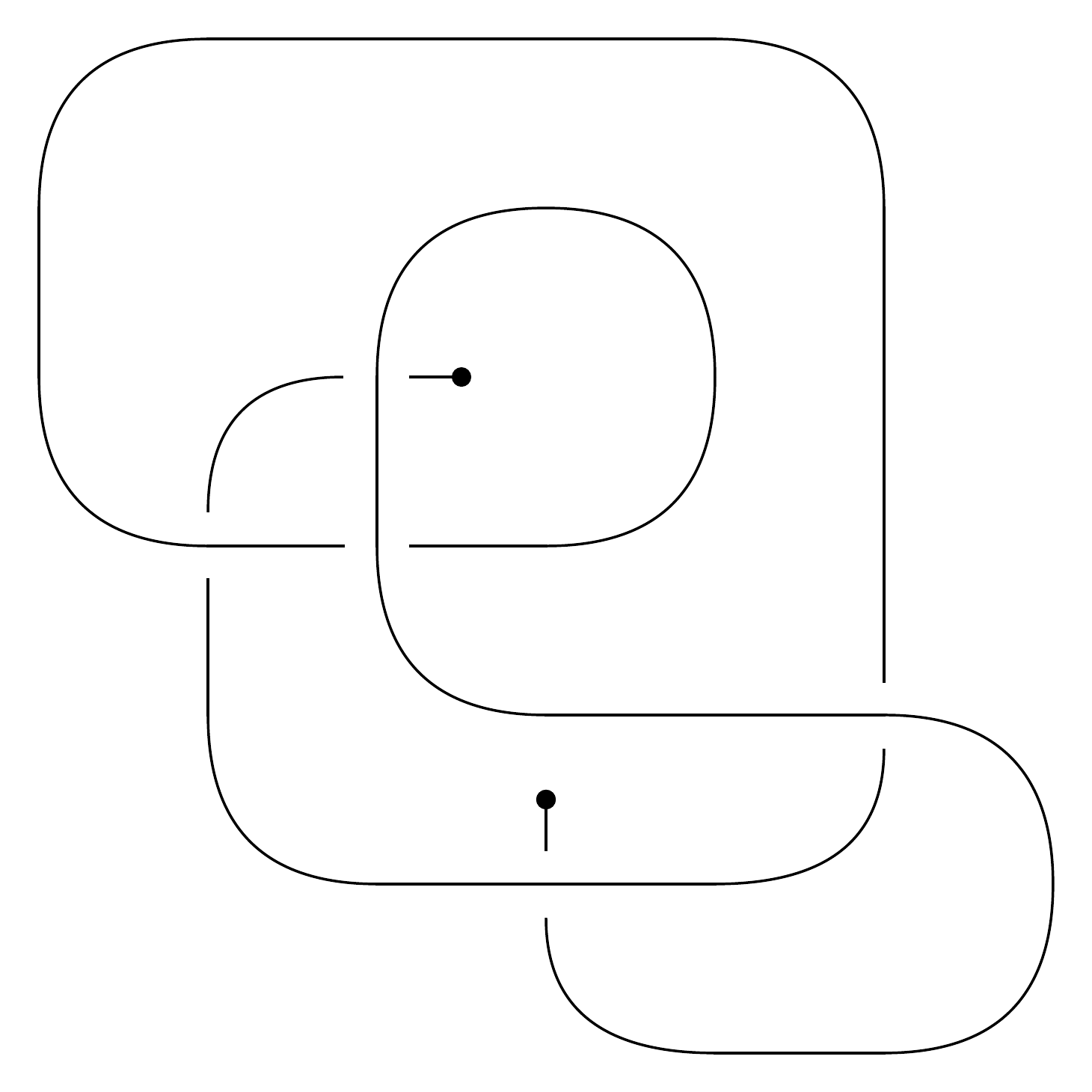}\\
\textcolor{black}{$5_{327}$}
\vspace{1cm}
\end{minipage}
\begin{minipage}[t]{.25\linewidth}
\centering
\includegraphics[width=0.9\textwidth,height=3.5cm,keepaspectratio]{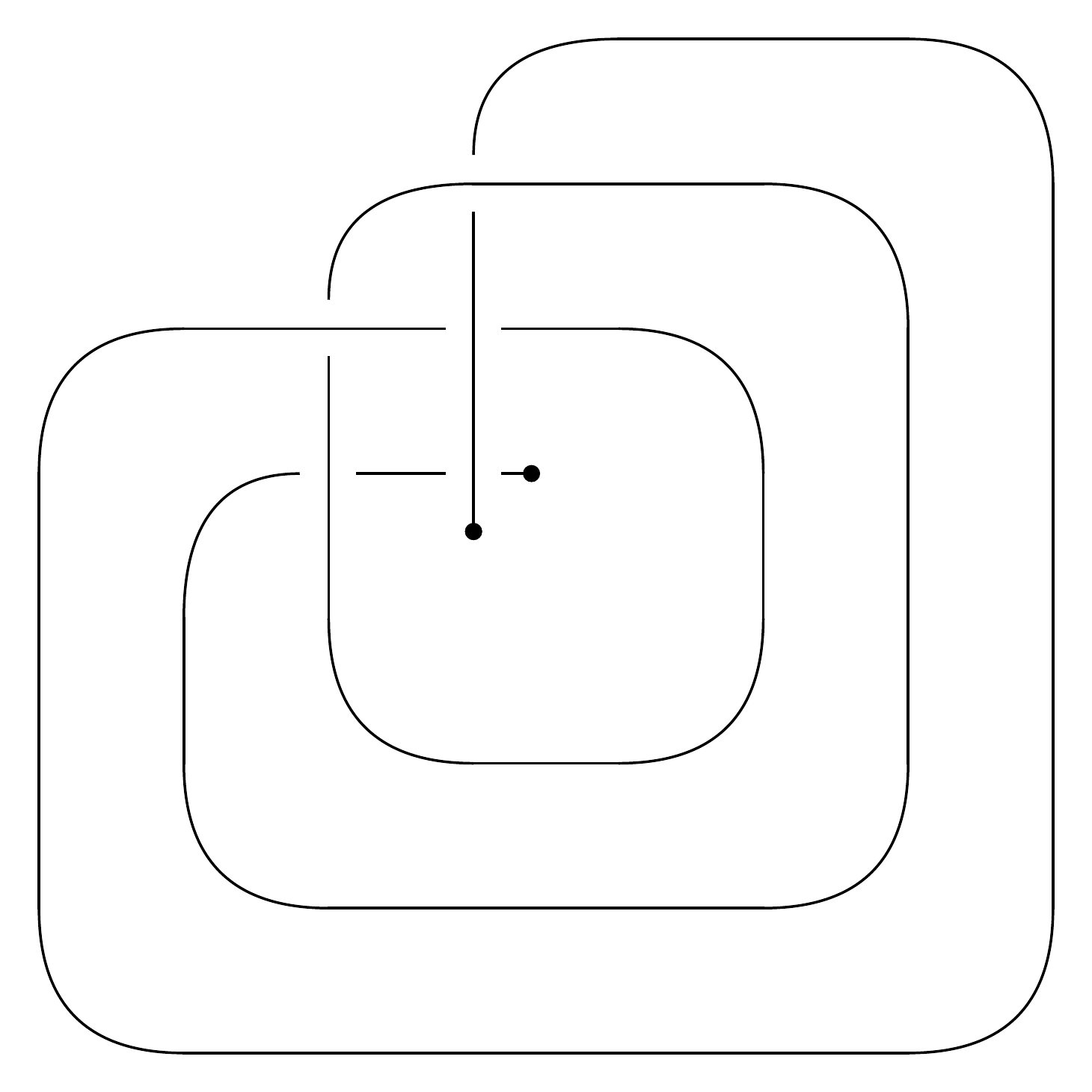}\\
\textcolor{black}{$5_{328}$}
\vspace{1cm}
\end{minipage}
\begin{minipage}[t]{.25\linewidth}
\centering
\includegraphics[width=0.9\textwidth,height=3.5cm,keepaspectratio]{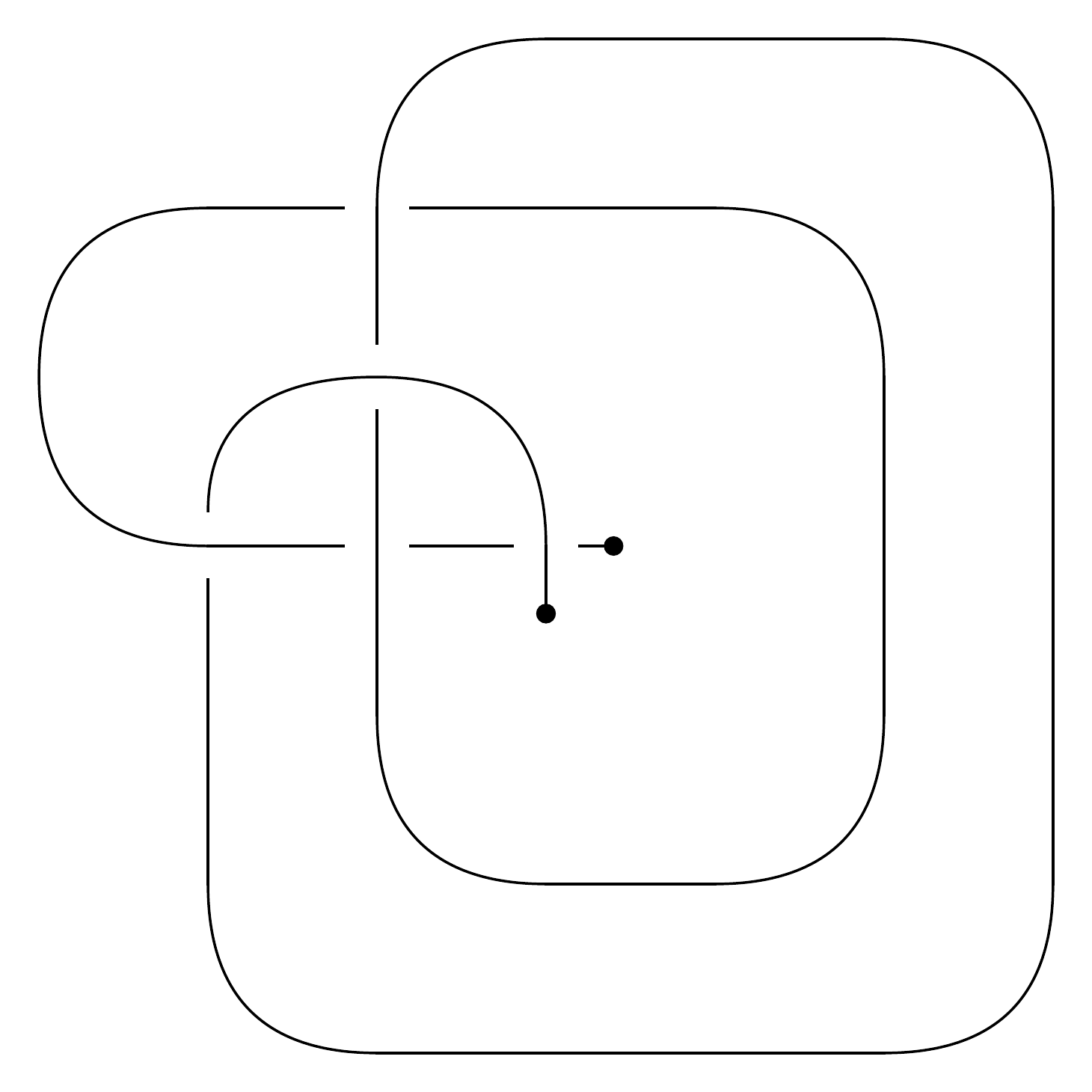}\\
\textcolor{black}{$5_{329}$}
\vspace{1cm}
\end{minipage}
\begin{minipage}[t]{.25\linewidth}
\centering
\includegraphics[width=0.9\textwidth,height=3.5cm,keepaspectratio]{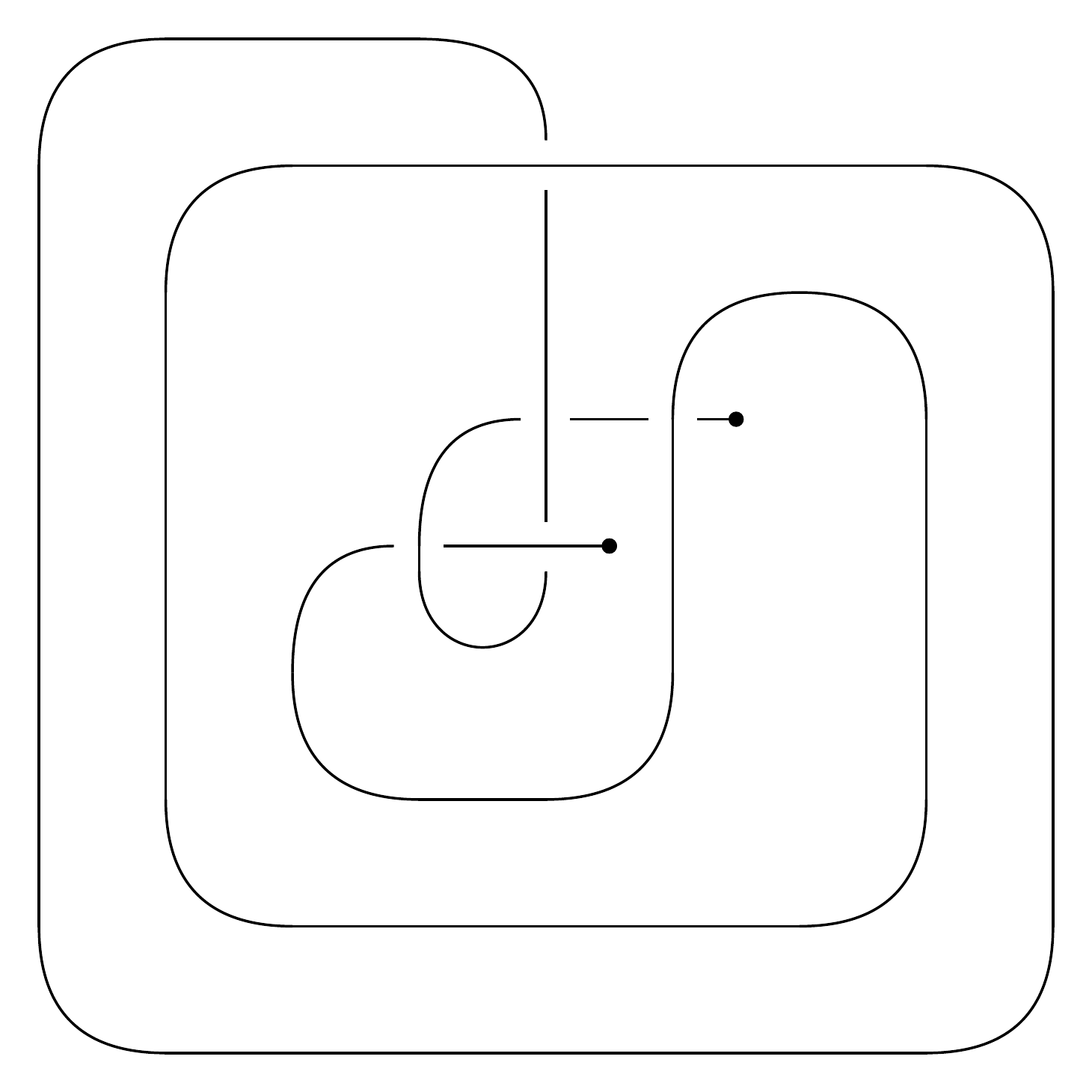}\\
\textcolor{black}{$5_{330}$}
\vspace{1cm}
\end{minipage}
\begin{minipage}[t]{.25\linewidth}
\centering
\includegraphics[width=0.9\textwidth,height=3.5cm,keepaspectratio]{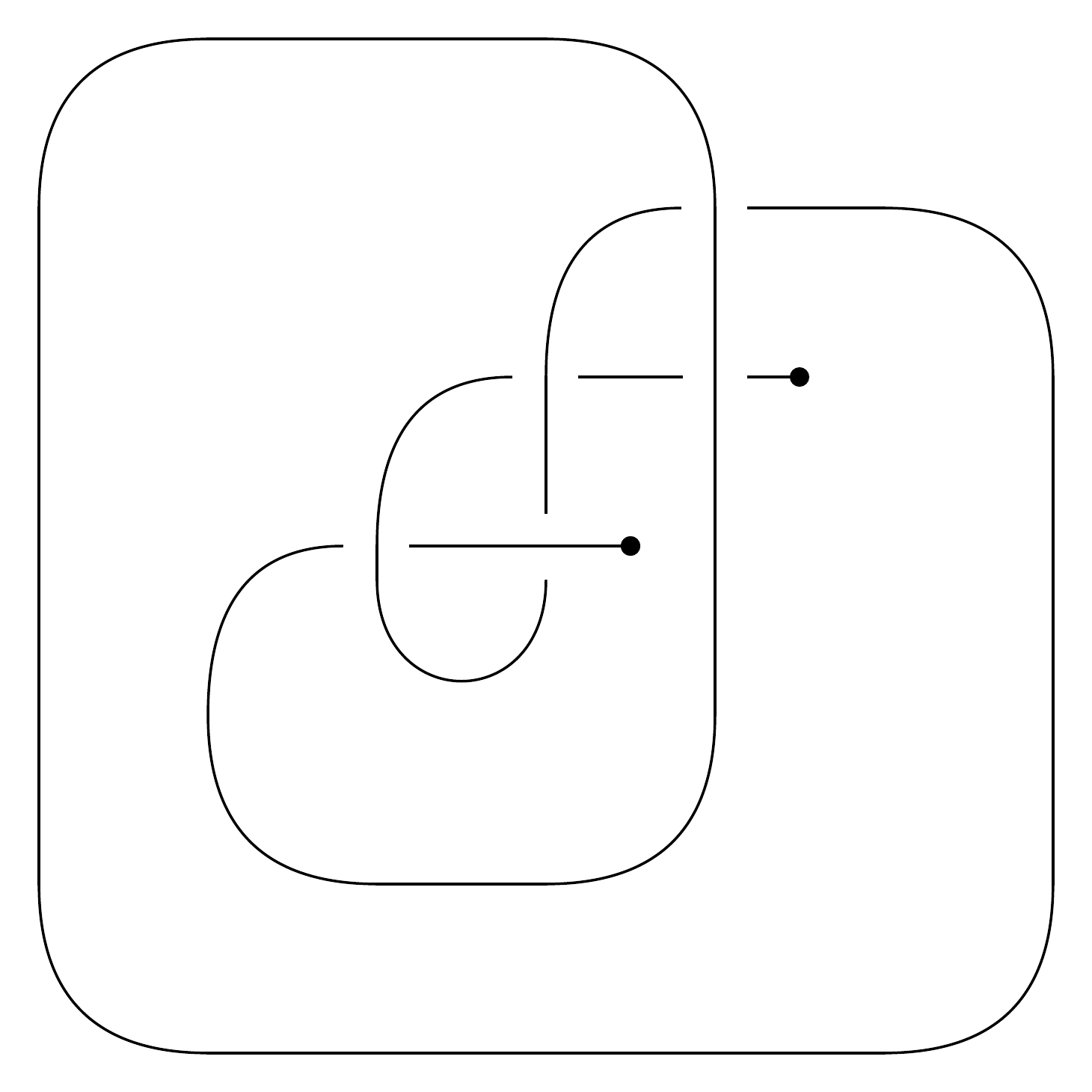}\\
\textcolor{black}{$5_{331}$}
\vspace{1cm}
\end{minipage}
\begin{minipage}[t]{.25\linewidth}
\centering
\includegraphics[width=0.9\textwidth,height=3.5cm,keepaspectratio]{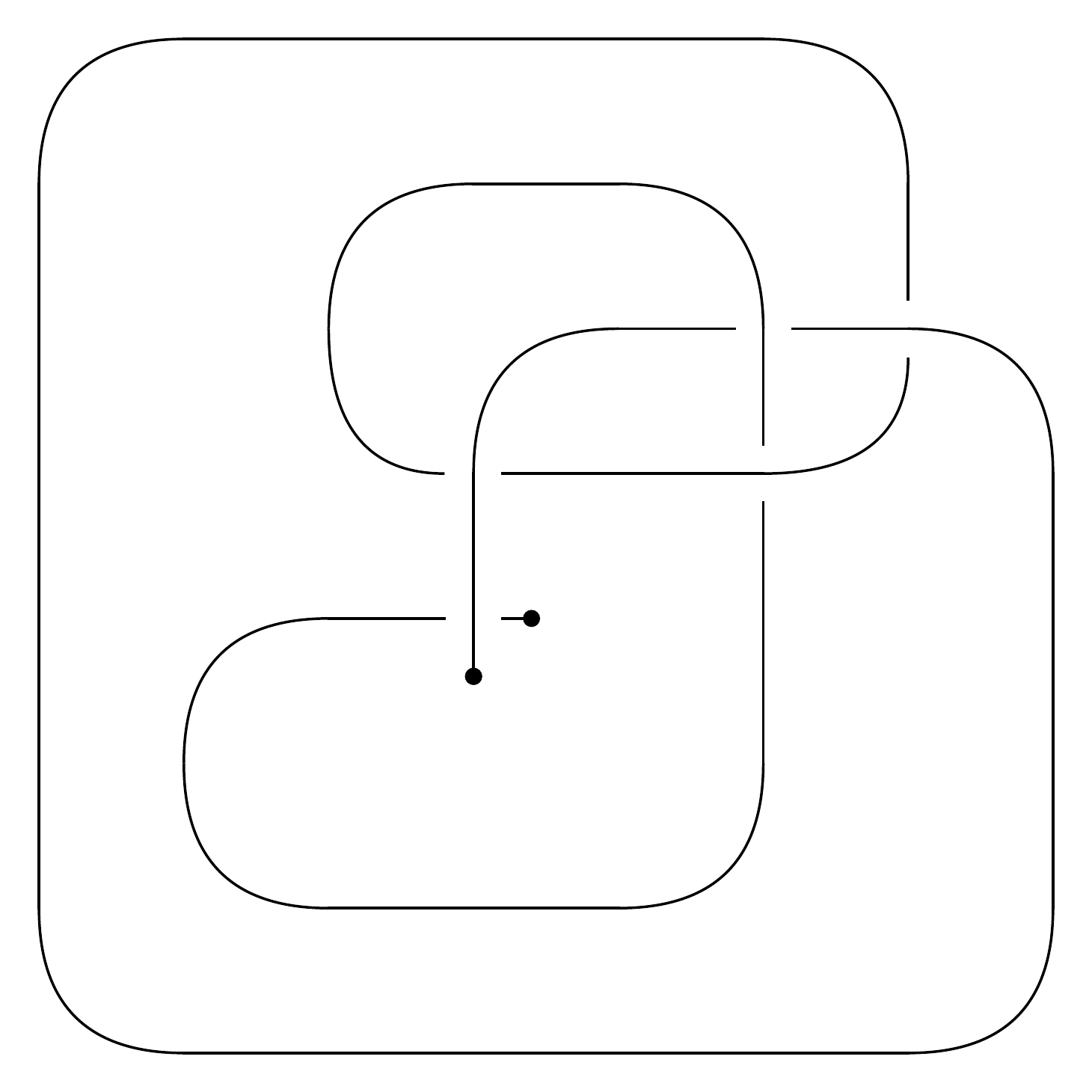}\\
\textcolor{black}{$5_{332}$}
\vspace{1cm}
\end{minipage}
\begin{minipage}[t]{.25\linewidth}
\centering
\includegraphics[width=0.9\textwidth,height=3.5cm,keepaspectratio]{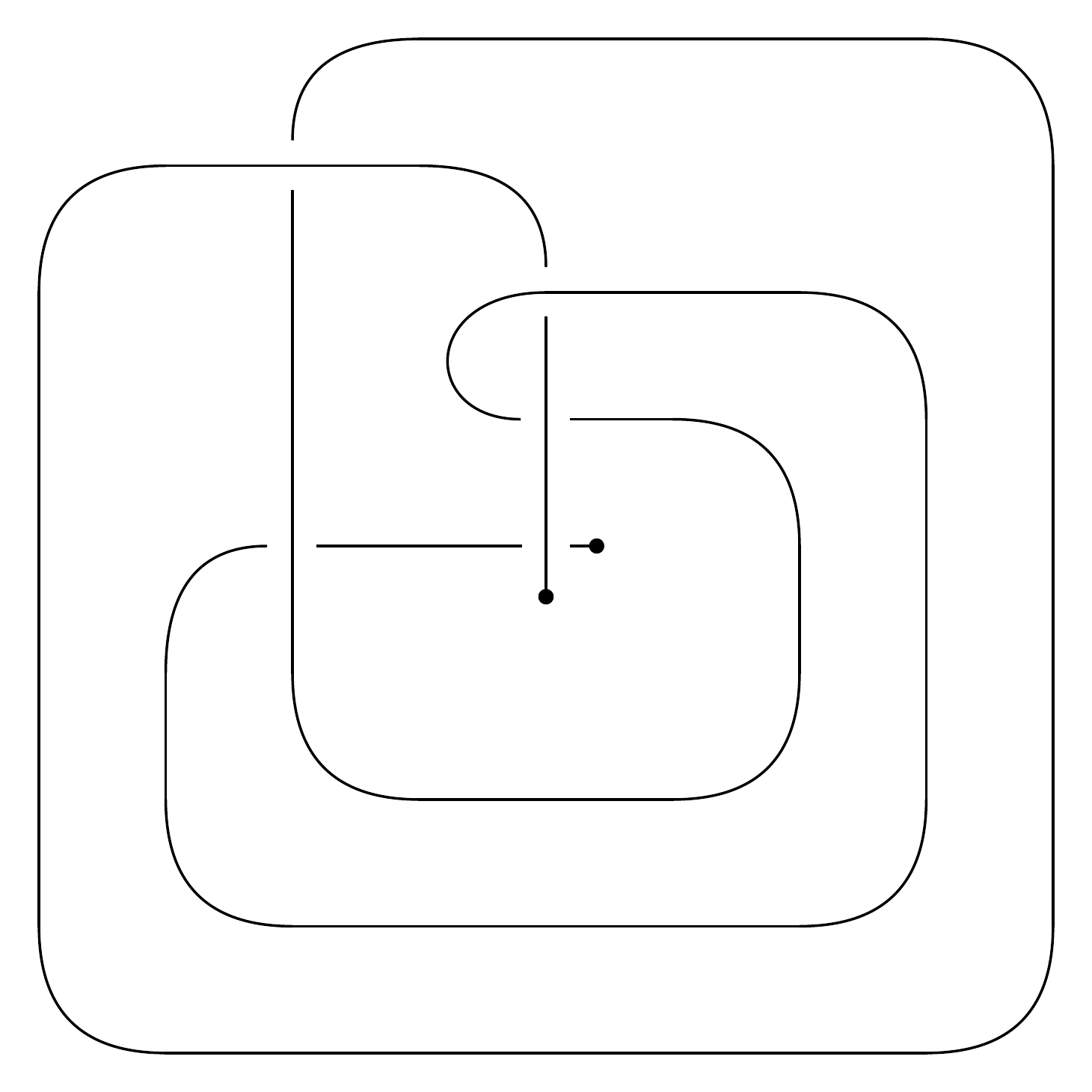}\\
\textcolor{black}{$5_{333}$}
\vspace{1cm}
\end{minipage}
\begin{minipage}[t]{.25\linewidth}
\centering
\includegraphics[width=0.9\textwidth,height=3.5cm,keepaspectratio]{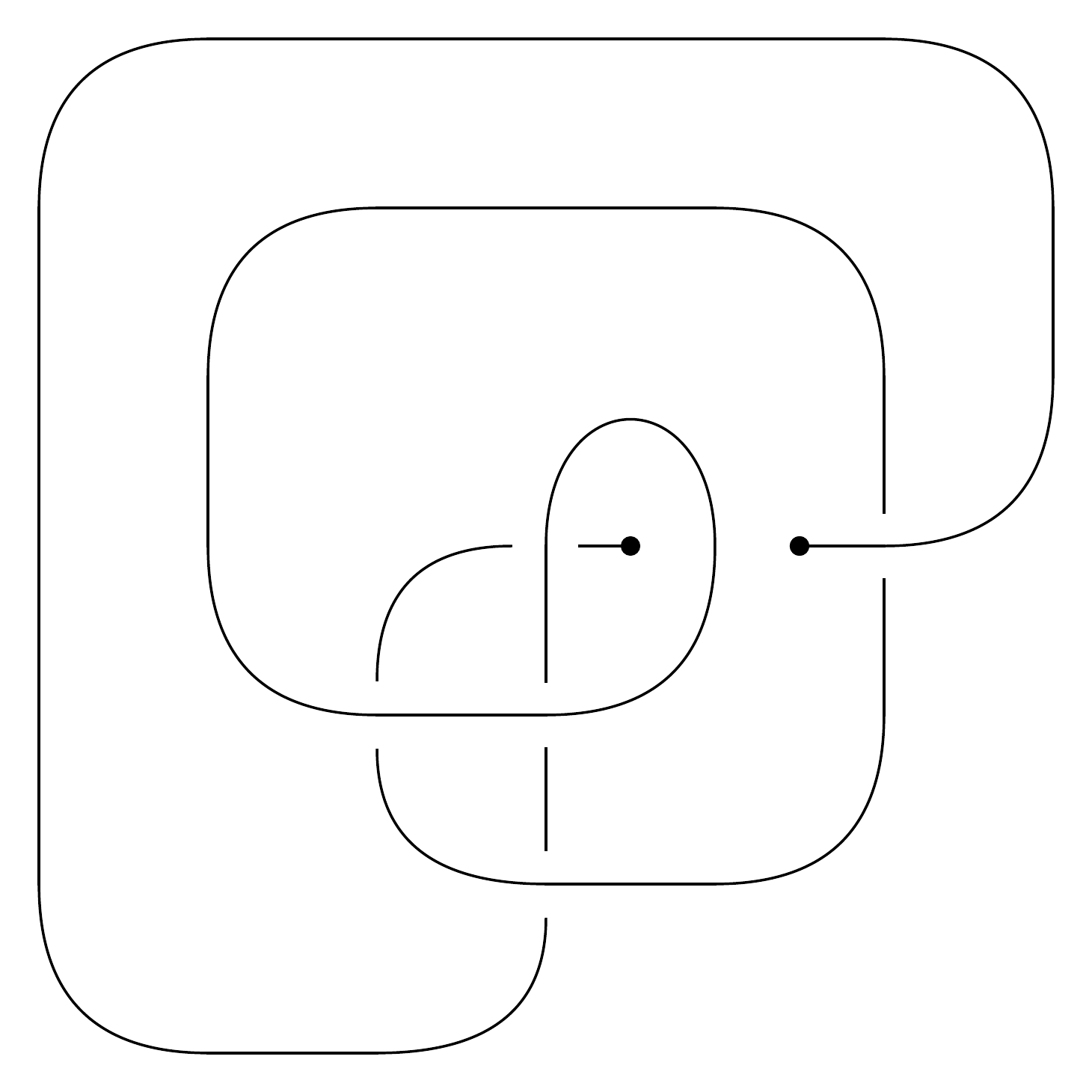}\\
\textcolor{black}{$5_{334}$}
\vspace{1cm}
\end{minipage}
\begin{minipage}[t]{.25\linewidth}
\centering
\includegraphics[width=0.9\textwidth,height=3.5cm,keepaspectratio]{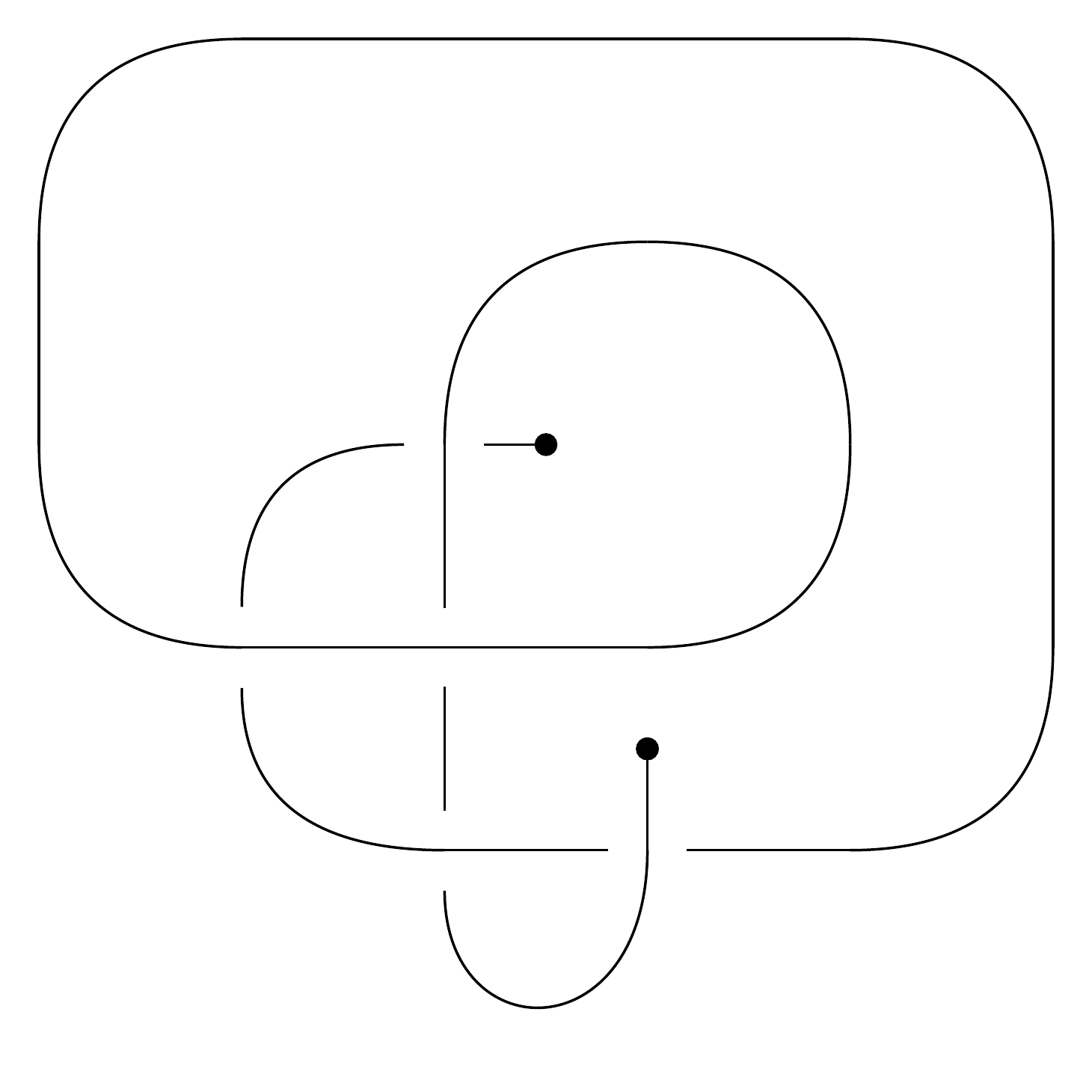}\\
\textcolor{black}{$5_{335}$}
\vspace{1cm}
\end{minipage}
\begin{minipage}[t]{.25\linewidth}
\centering
\includegraphics[width=0.9\textwidth,height=3.5cm,keepaspectratio]{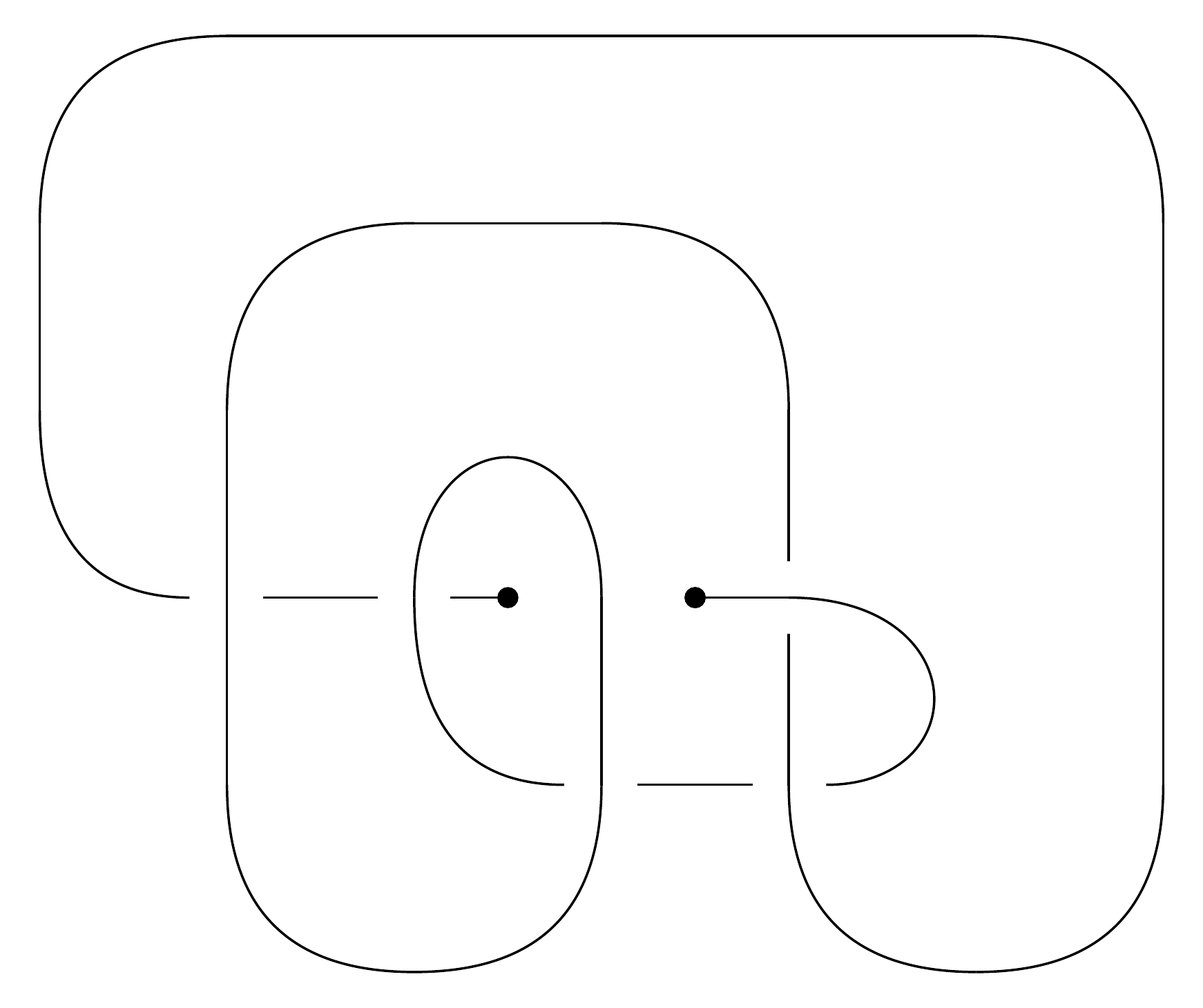}\\
\textcolor{black}{$5_{336}$}
\vspace{1cm}
\end{minipage}
\begin{minipage}[t]{.25\linewidth}
\centering
\includegraphics[width=0.9\textwidth,height=3.5cm,keepaspectratio]{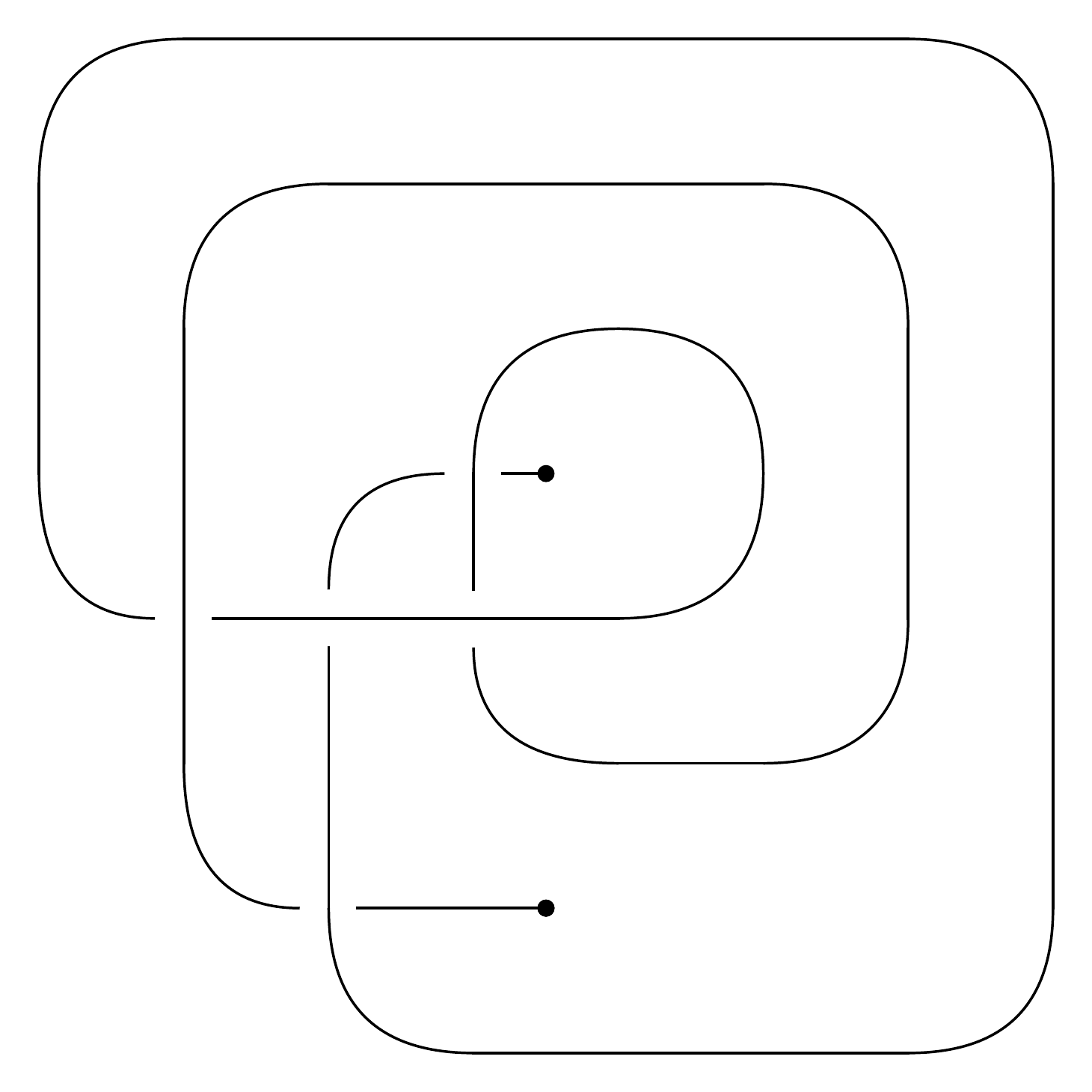}\\
\textcolor{black}{$5_{337}$}
\vspace{1cm}
\end{minipage}
\begin{minipage}[t]{.25\linewidth}
\centering
\includegraphics[width=0.9\textwidth,height=3.5cm,keepaspectratio]{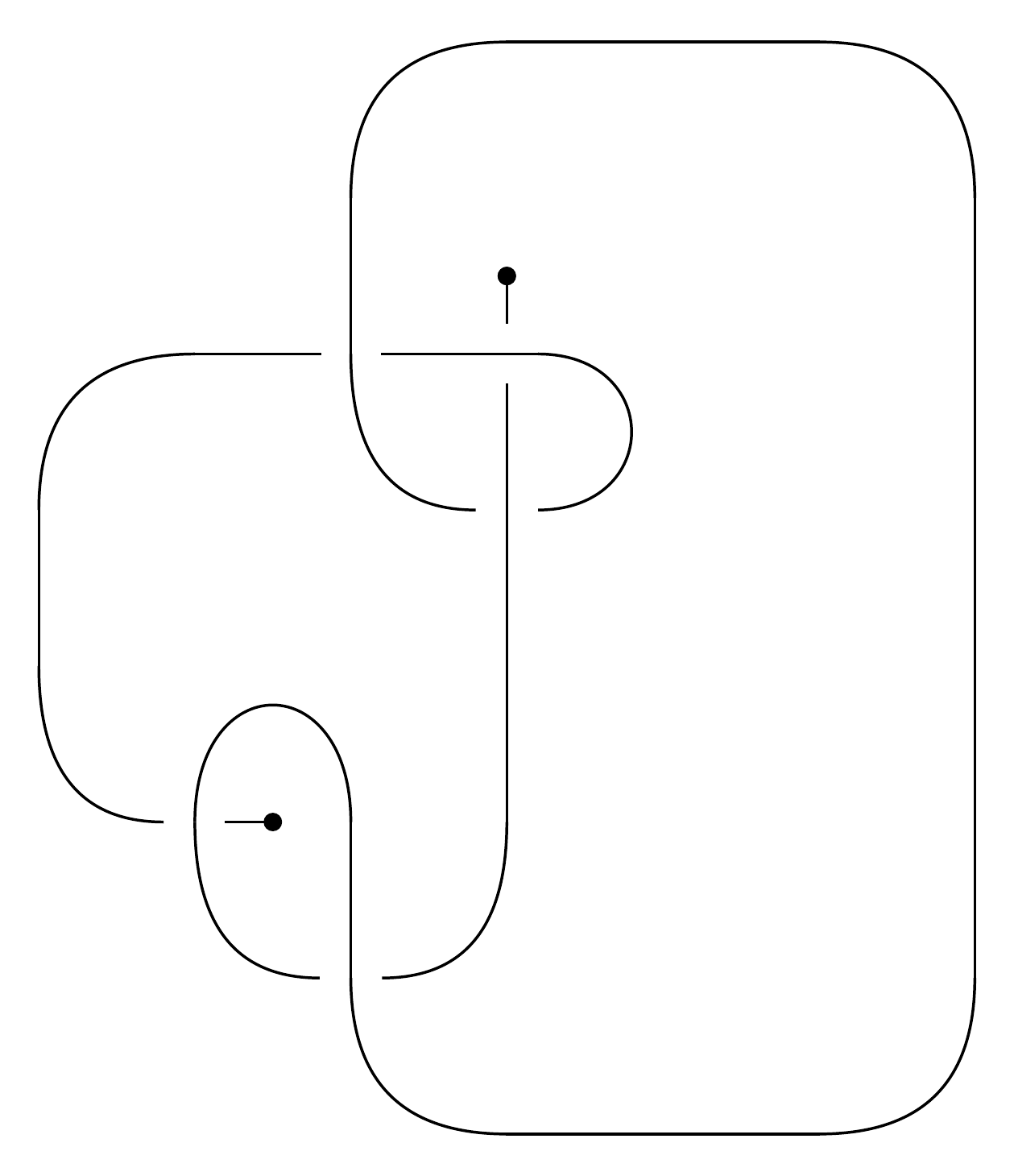}\\
\textcolor{black}{$5_{338}$}
\vspace{1cm}
\end{minipage}
\begin{minipage}[t]{.25\linewidth}
\centering
\includegraphics[width=0.9\textwidth,height=3.5cm,keepaspectratio]{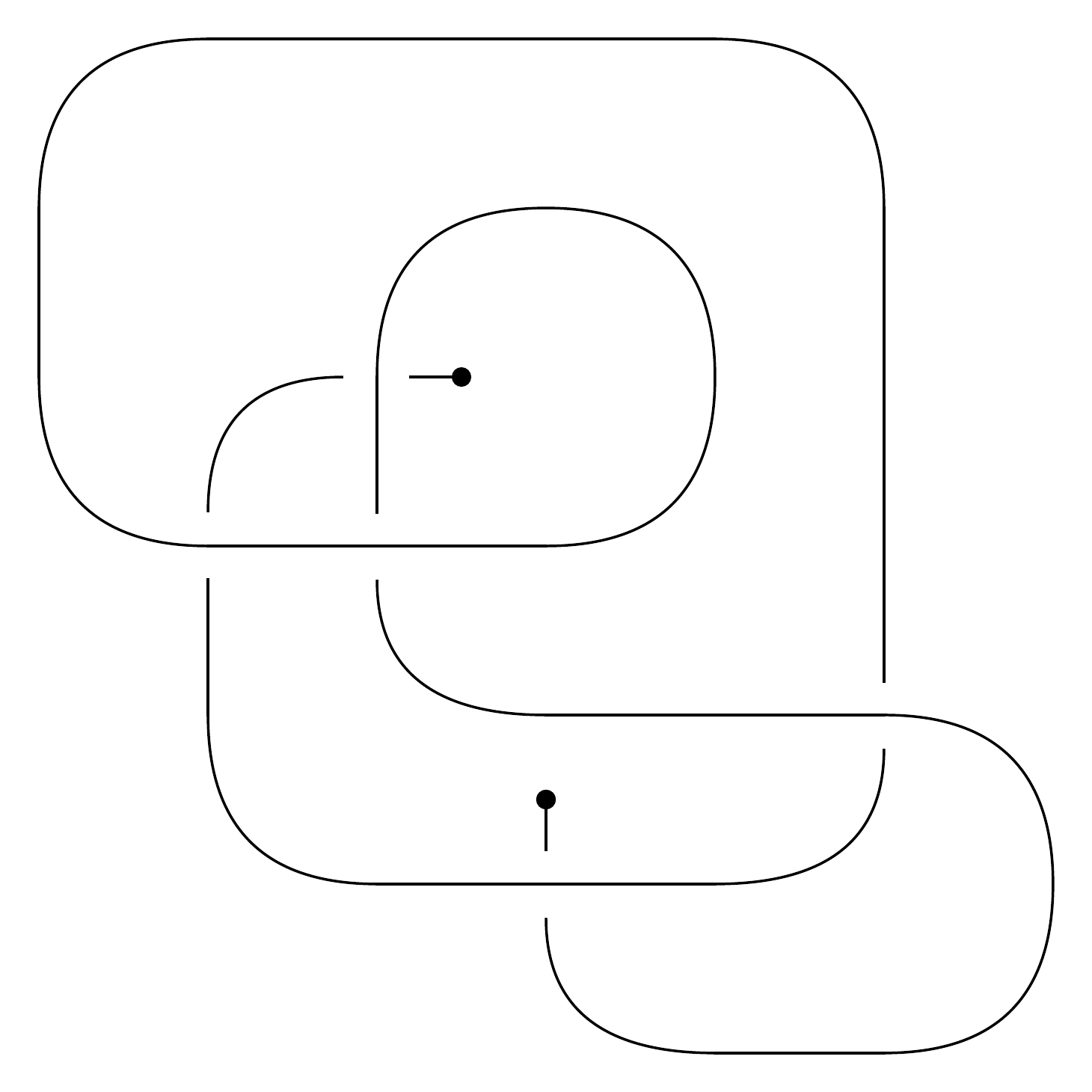}\\
\textcolor{black}{$5_{339}$}
\vspace{1cm}
\end{minipage}
\begin{minipage}[t]{.25\linewidth}
\centering
\includegraphics[width=0.9\textwidth,height=3.5cm,keepaspectratio]{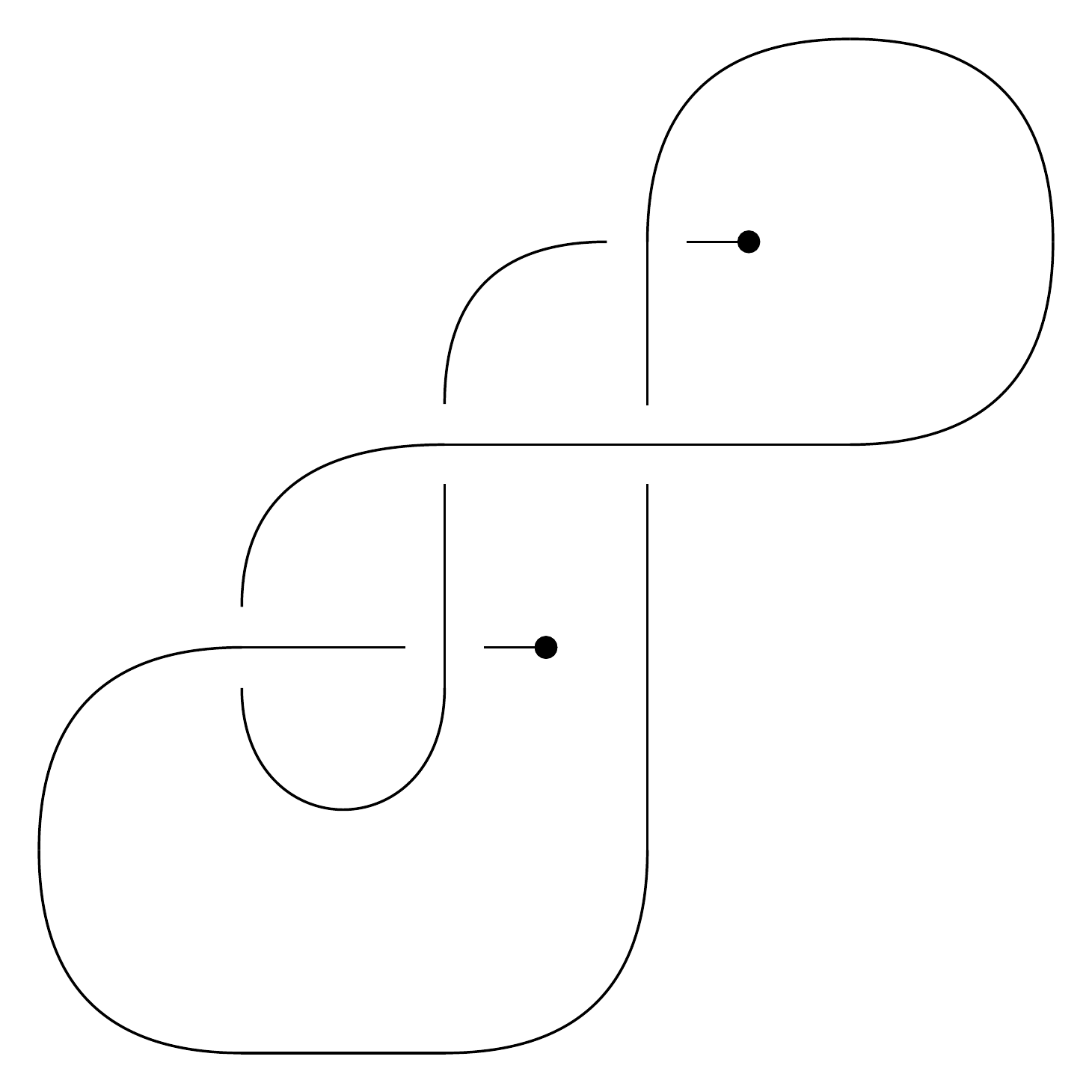}\\
\textcolor{black}{$5_{340}$}
\vspace{1cm}
\end{minipage}
\begin{minipage}[t]{.25\linewidth}
\centering
\includegraphics[width=0.9\textwidth,height=3.5cm,keepaspectratio]{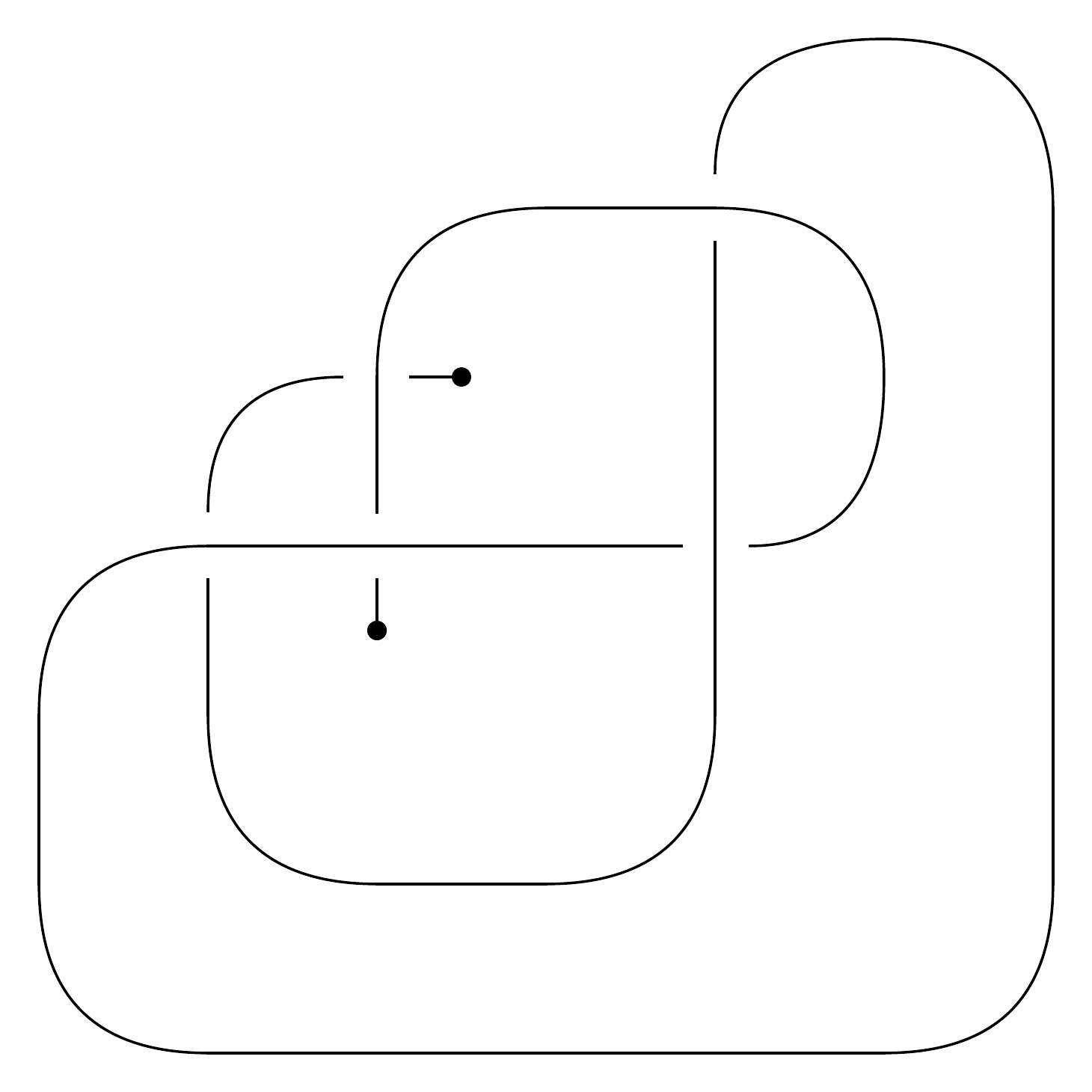}\\
\textcolor{black}{$5_{341}$}
\vspace{1cm}
\end{minipage}
\begin{minipage}[t]{.25\linewidth}
\centering
\includegraphics[width=0.9\textwidth,height=3.5cm,keepaspectratio]{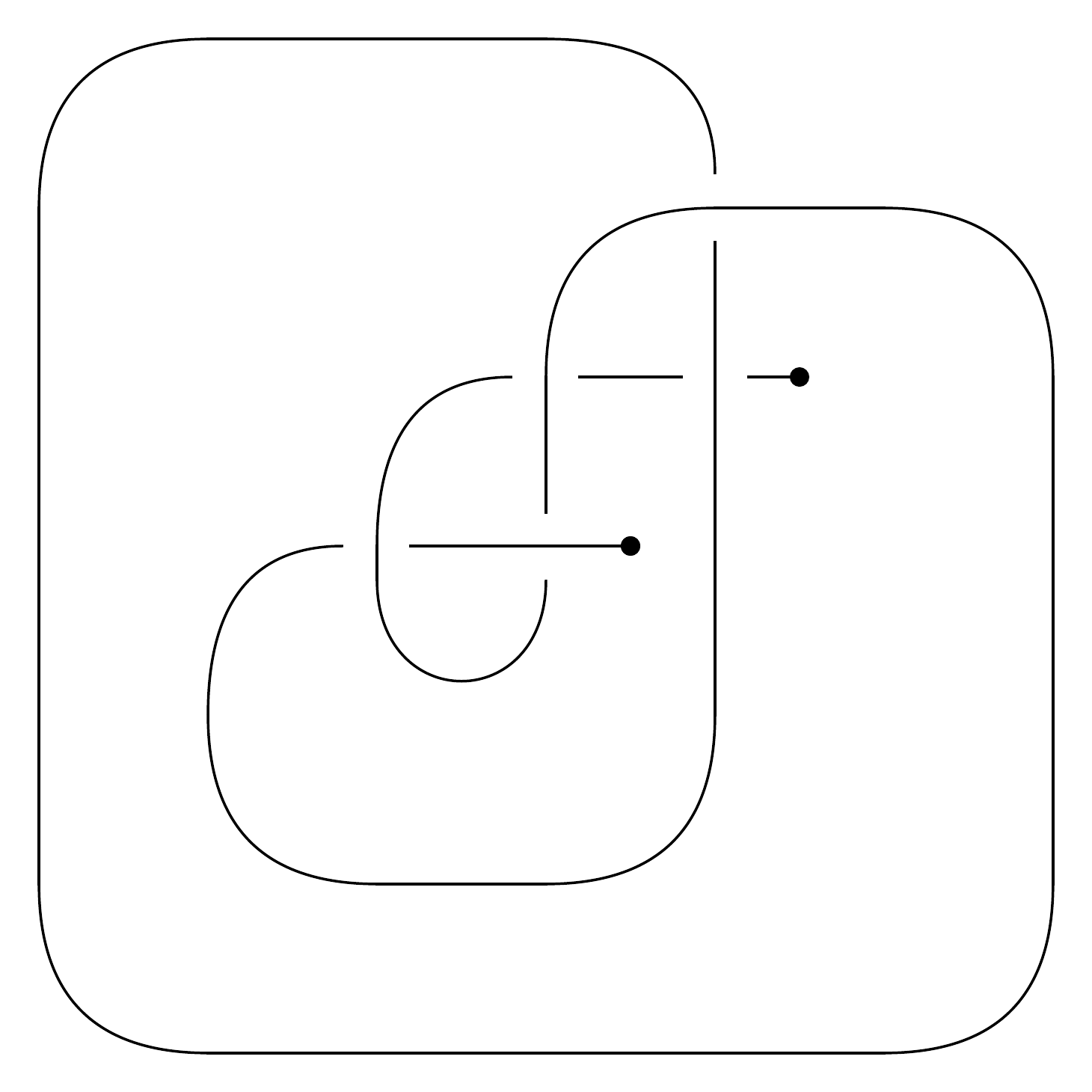}\\
\textcolor{black}{$5_{342}$}
\vspace{1cm}
\end{minipage}
\begin{minipage}[t]{.25\linewidth}
\centering
\includegraphics[width=0.9\textwidth,height=3.5cm,keepaspectratio]{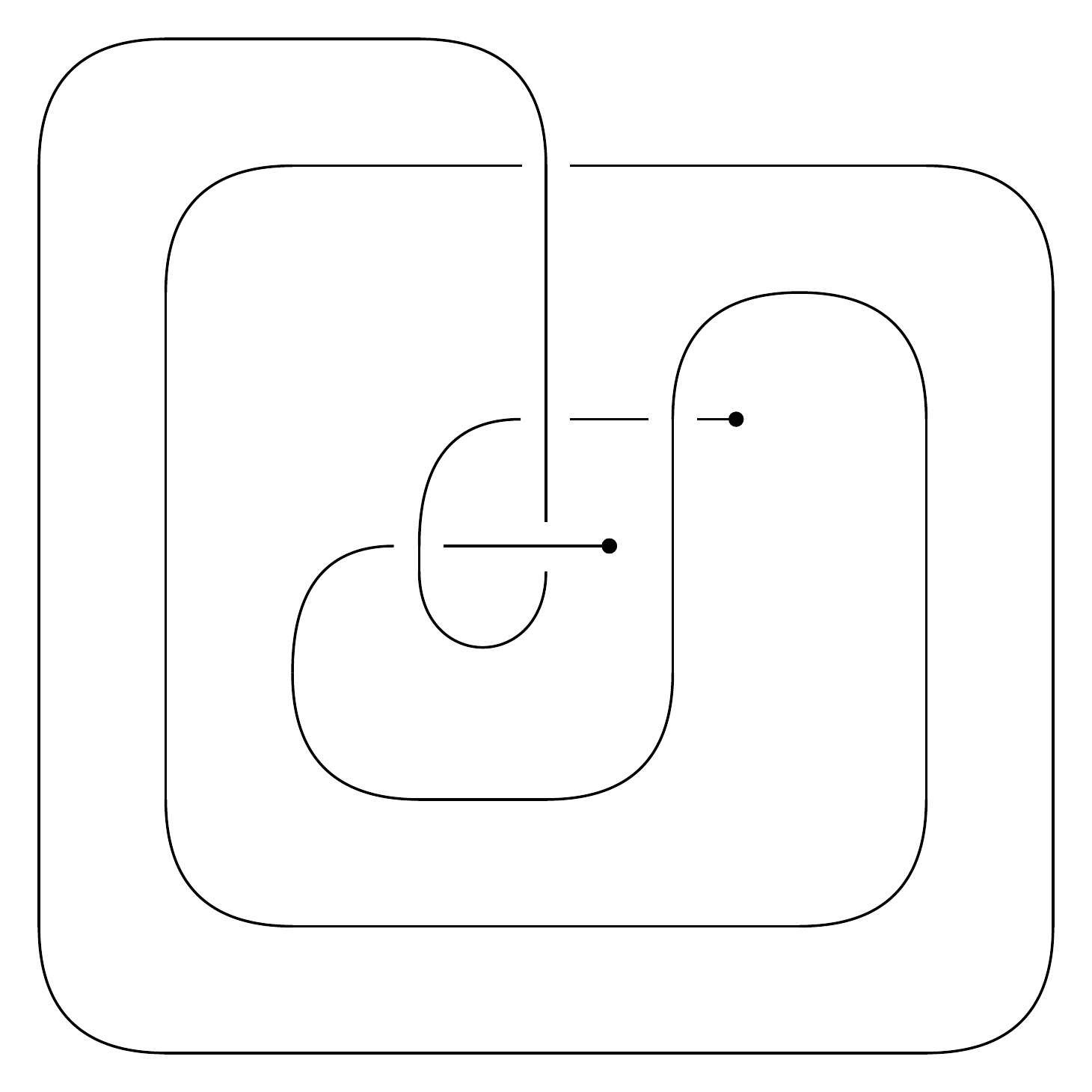}\\
\textcolor{black}{$5_{343}$}
\vspace{1cm}
\end{minipage}
\begin{minipage}[t]{.25\linewidth}
\centering
\includegraphics[width=0.9\textwidth,height=3.5cm,keepaspectratio]{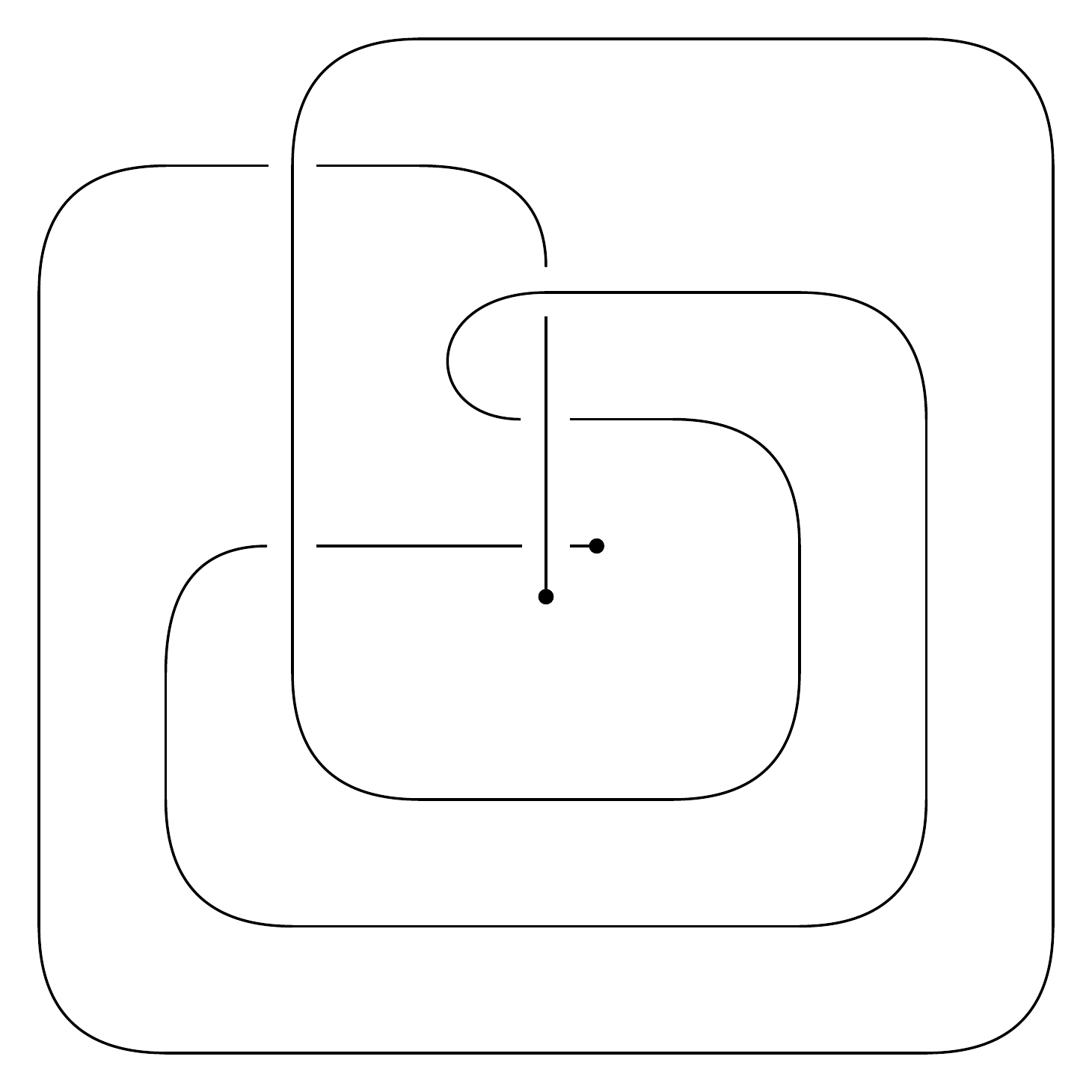}\\
\textcolor{black}{$5_{344}$}
\vspace{1cm}
\end{minipage}
\begin{minipage}[t]{.25\linewidth}
\centering
\includegraphics[width=0.9\textwidth,height=3.5cm,keepaspectratio]{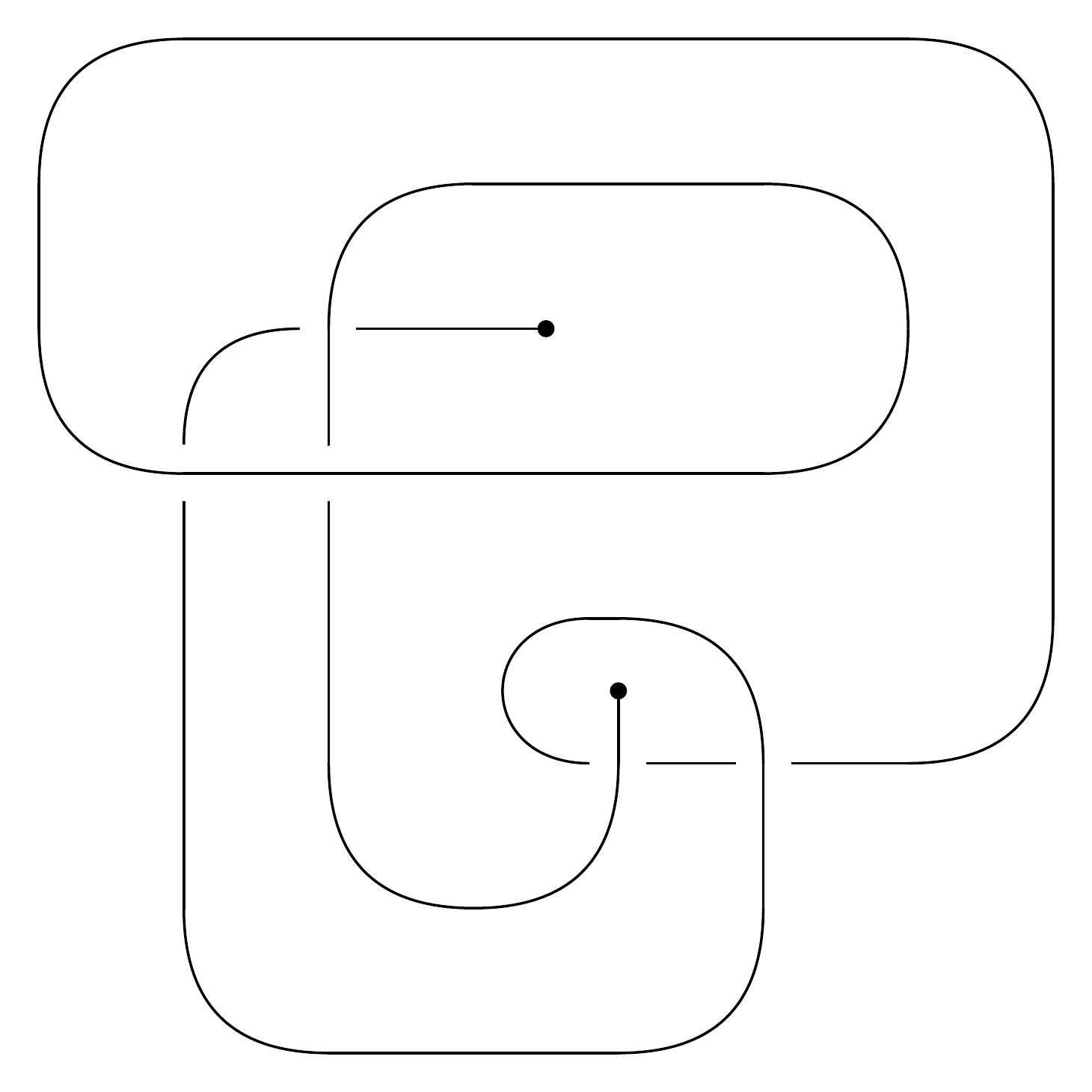}\\
\textcolor{black}{$5_{345}$}
\vspace{1cm}
\end{minipage}
\begin{minipage}[t]{.25\linewidth}
\centering
\includegraphics[width=0.9\textwidth,height=3.5cm,keepaspectratio]{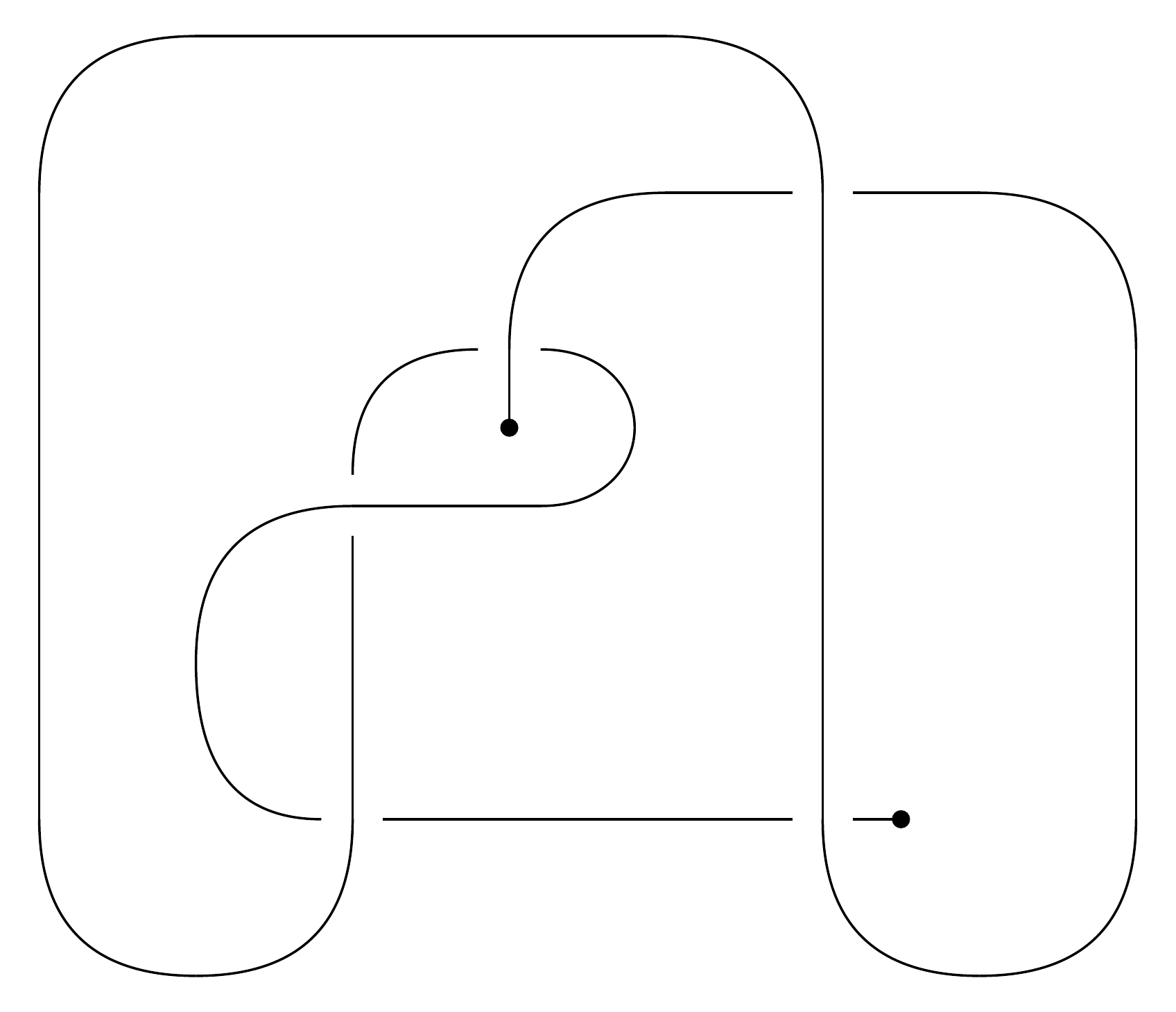}\\
\textcolor{black}{$5_{346}$}
\vspace{1cm}
\end{minipage}
\begin{minipage}[t]{.25\linewidth}
\centering
\includegraphics[width=0.9\textwidth,height=3.5cm,keepaspectratio]{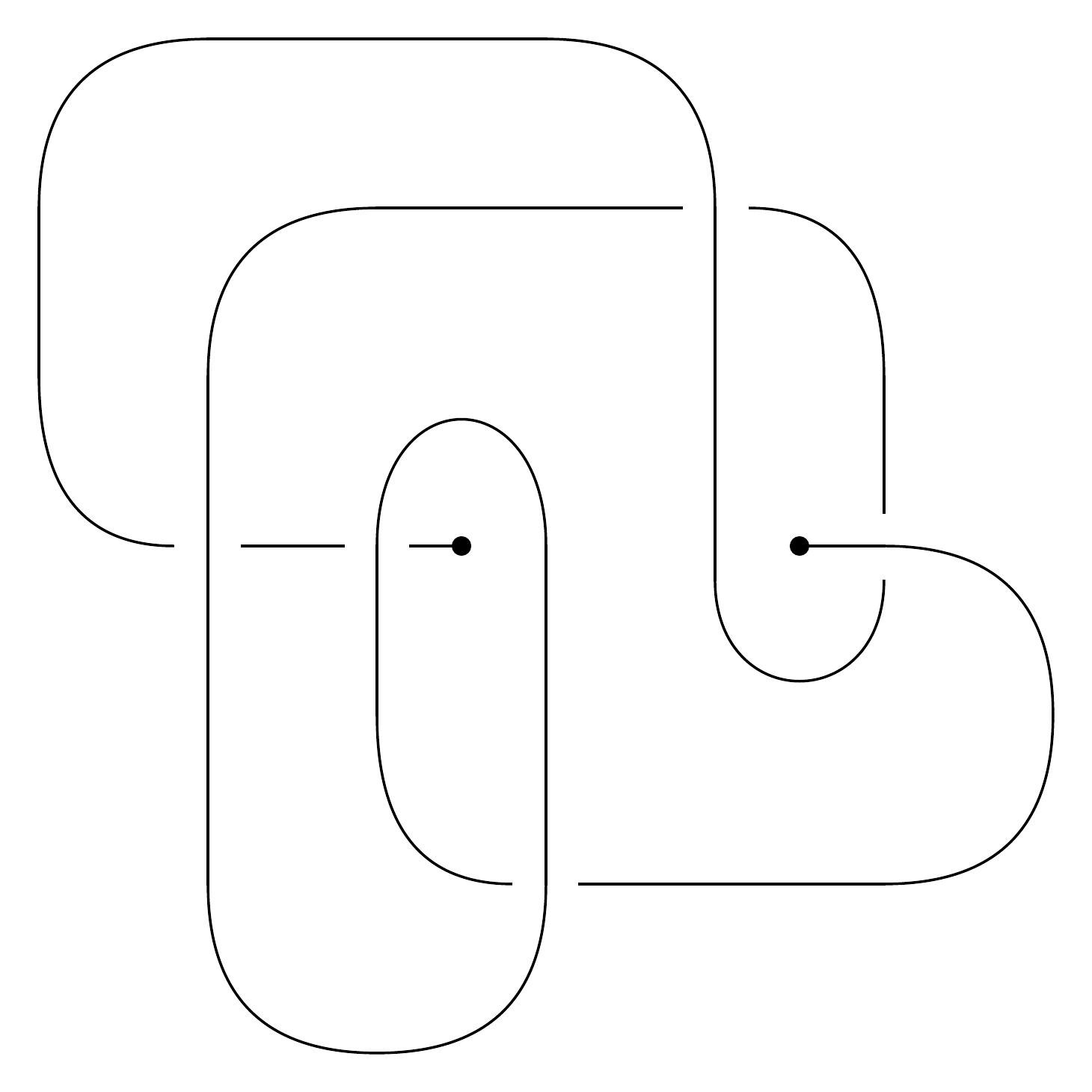}\\
\textcolor{black}{$5_{347}$}
\vspace{1cm}
\end{minipage}
\begin{minipage}[t]{.25\linewidth}
\centering
\includegraphics[width=0.9\textwidth,height=3.5cm,keepaspectratio]{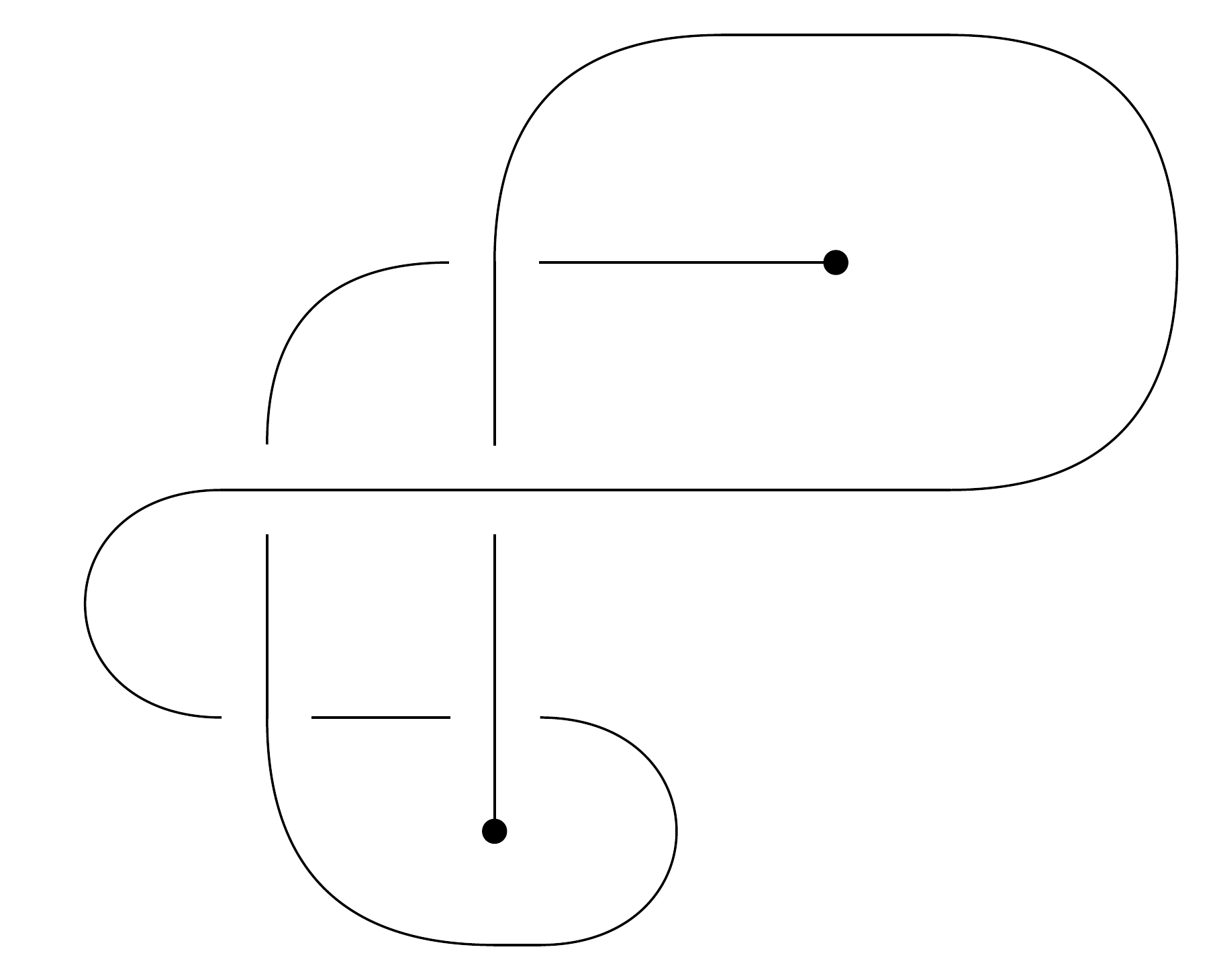}\\
\textcolor{black}{$5_{348}$}
\vspace{1cm}
\end{minipage}
\begin{minipage}[t]{.25\linewidth}
\centering
\includegraphics[width=0.9\textwidth,height=3.5cm,keepaspectratio]{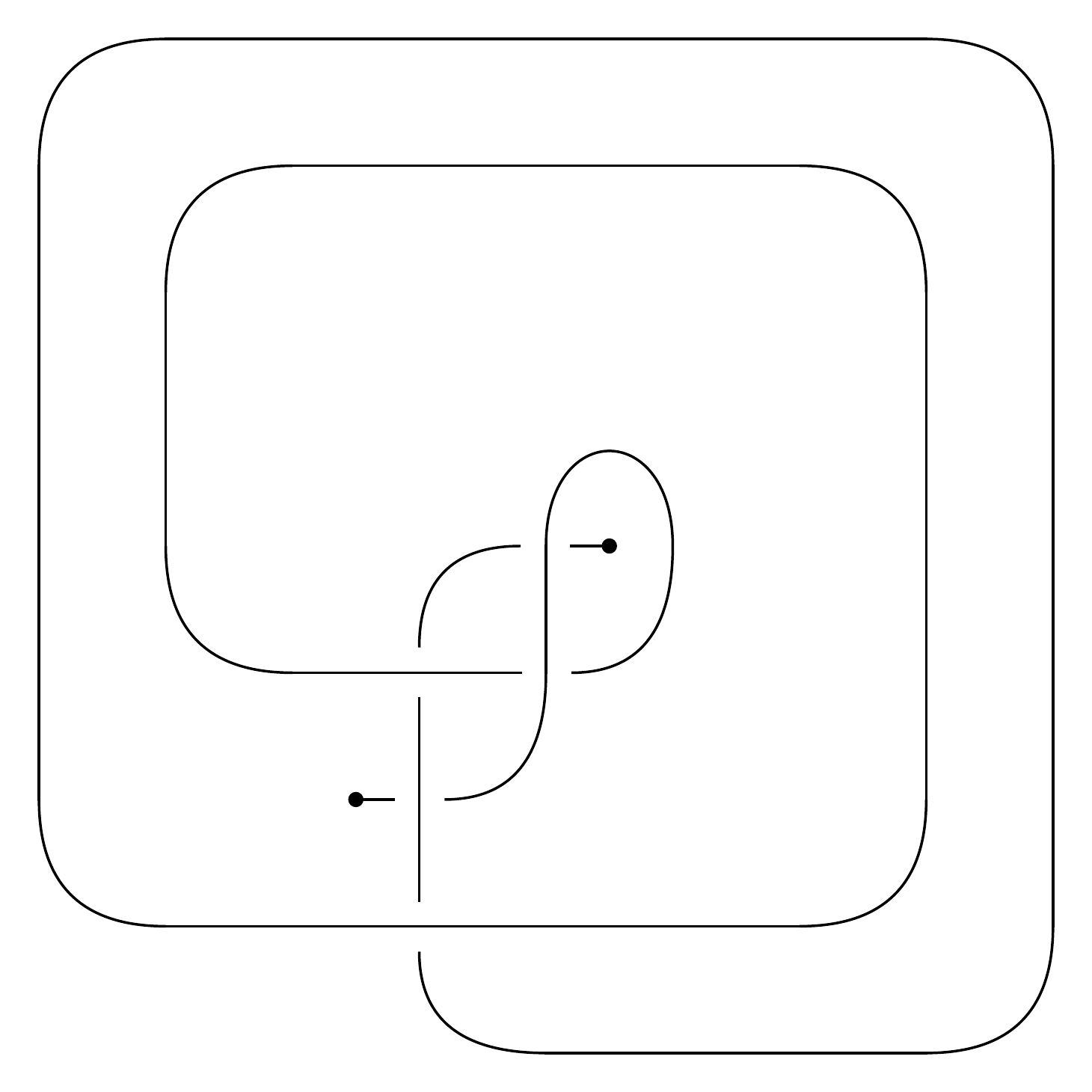}\\
\textcolor{black}{$5_{349}$}
\vspace{1cm}
\end{minipage}
\begin{minipage}[t]{.25\linewidth}
\centering
\includegraphics[width=0.9\textwidth,height=3.5cm,keepaspectratio]{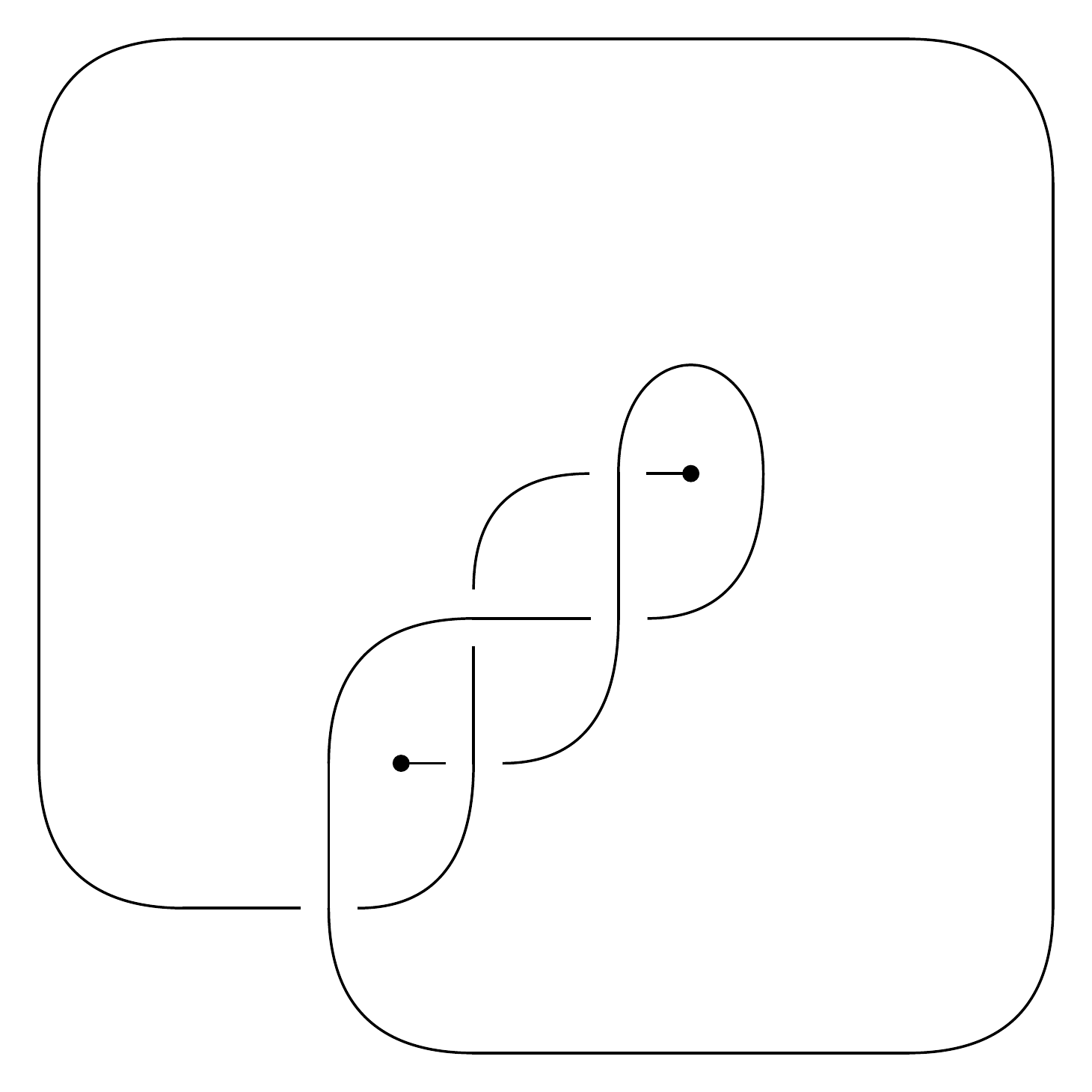}\\
\textcolor{black}{$5_{350}$}
\vspace{1cm}
\end{minipage}
\begin{minipage}[t]{.25\linewidth}
\centering
\includegraphics[width=0.9\textwidth,height=3.5cm,keepaspectratio]{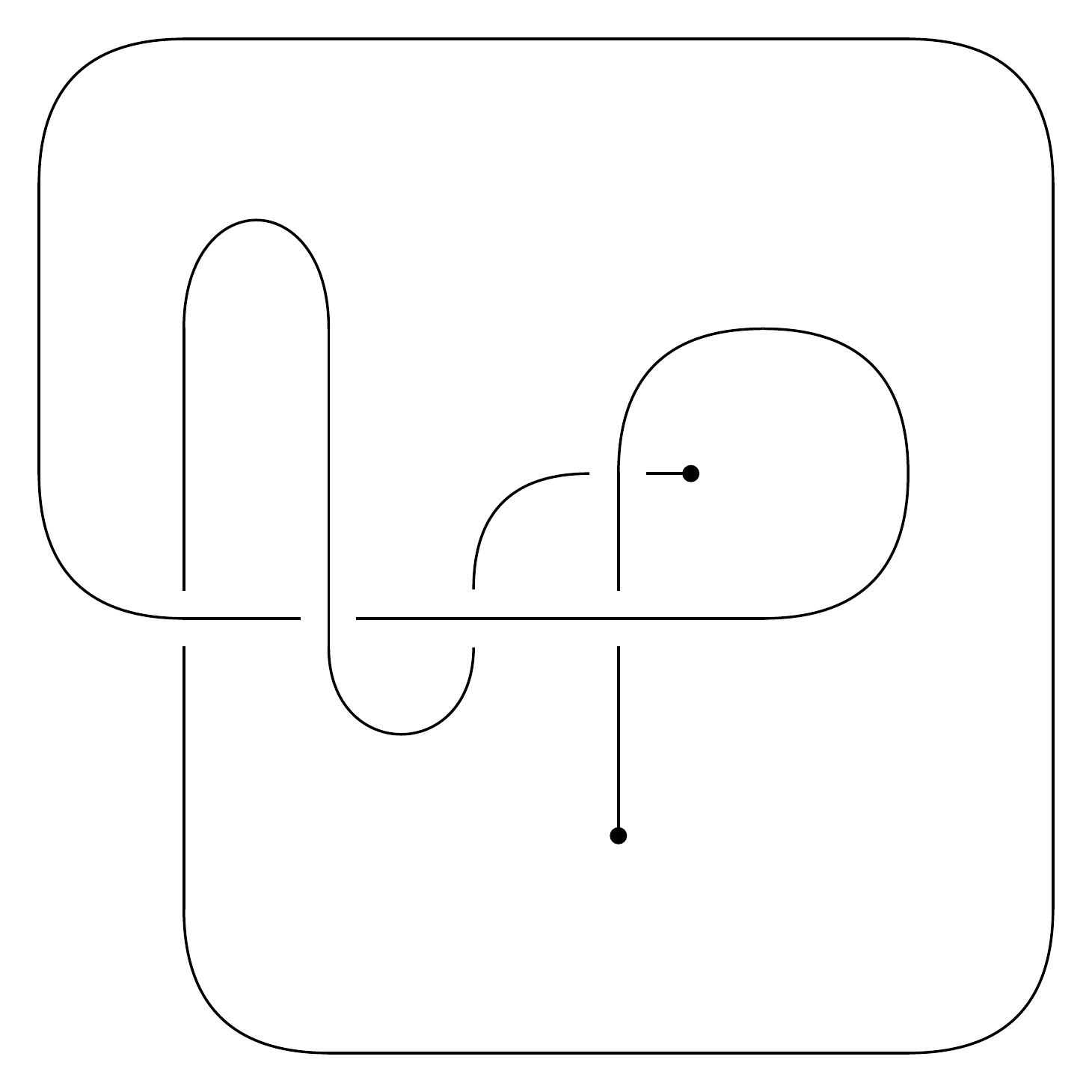}\\
\textcolor{black}{$5_{351}$}
\vspace{1cm}
\end{minipage}
\begin{minipage}[t]{.25\linewidth}
\centering
\includegraphics[width=0.9\textwidth,height=3.5cm,keepaspectratio]{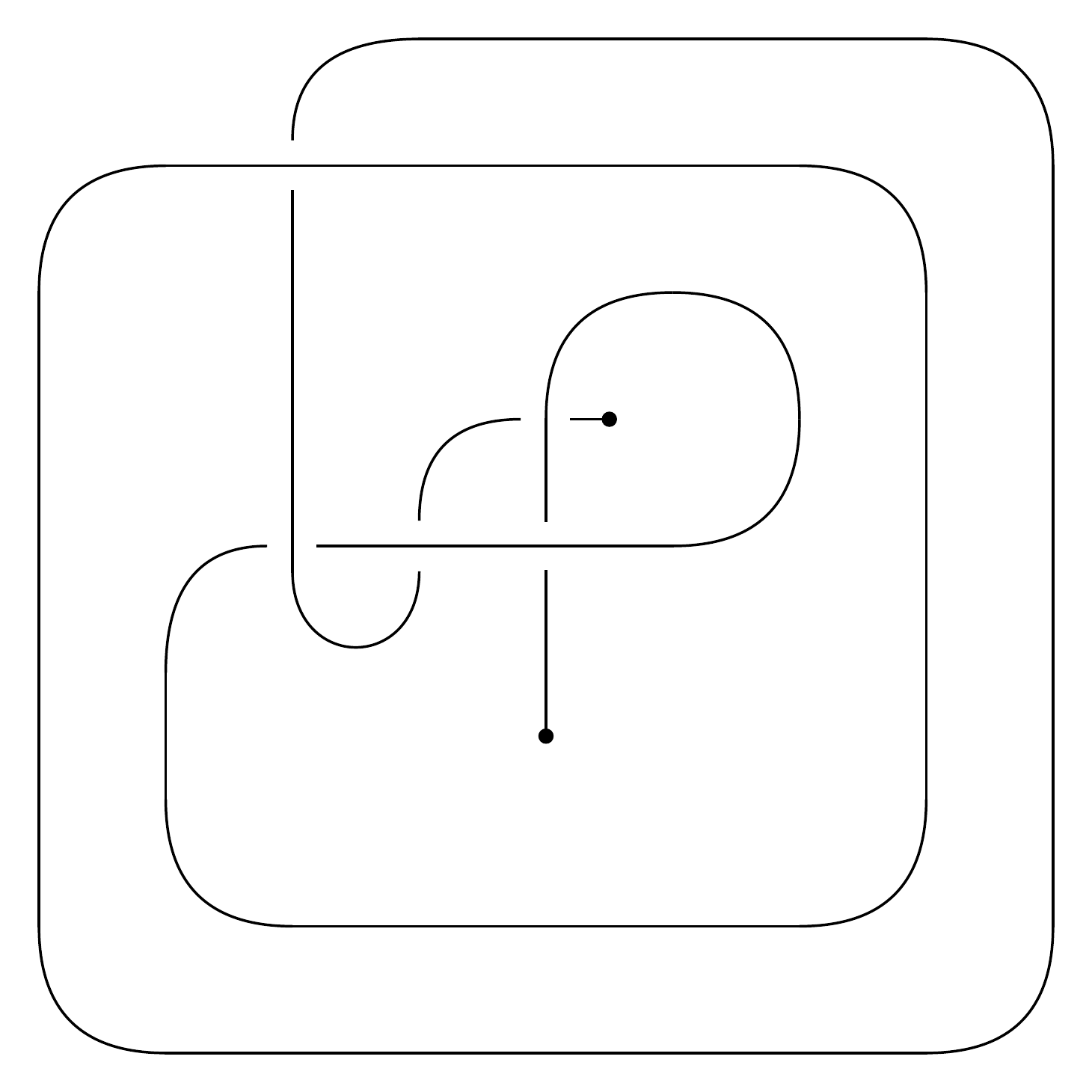}\\
\textcolor{black}{$5_{352}$}
\vspace{1cm}
\end{minipage}
\begin{minipage}[t]{.25\linewidth}
\centering
\includegraphics[width=0.9\textwidth,height=3.5cm,keepaspectratio]{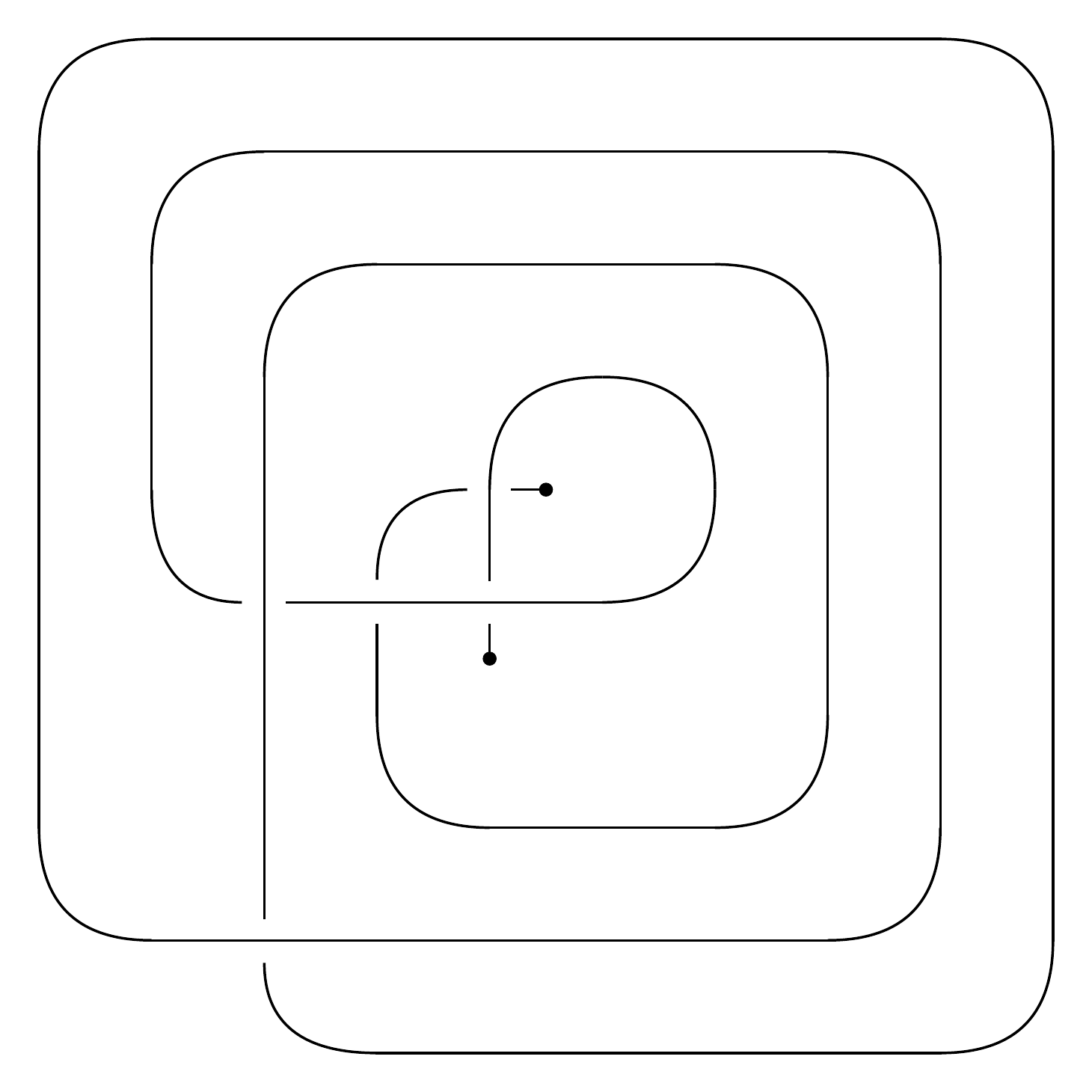}\\
\textcolor{black}{$5_{353}$}
\vspace{1cm}
\end{minipage}
\begin{minipage}[t]{.25\linewidth}
\centering
\includegraphics[width=0.9\textwidth,height=3.5cm,keepaspectratio]{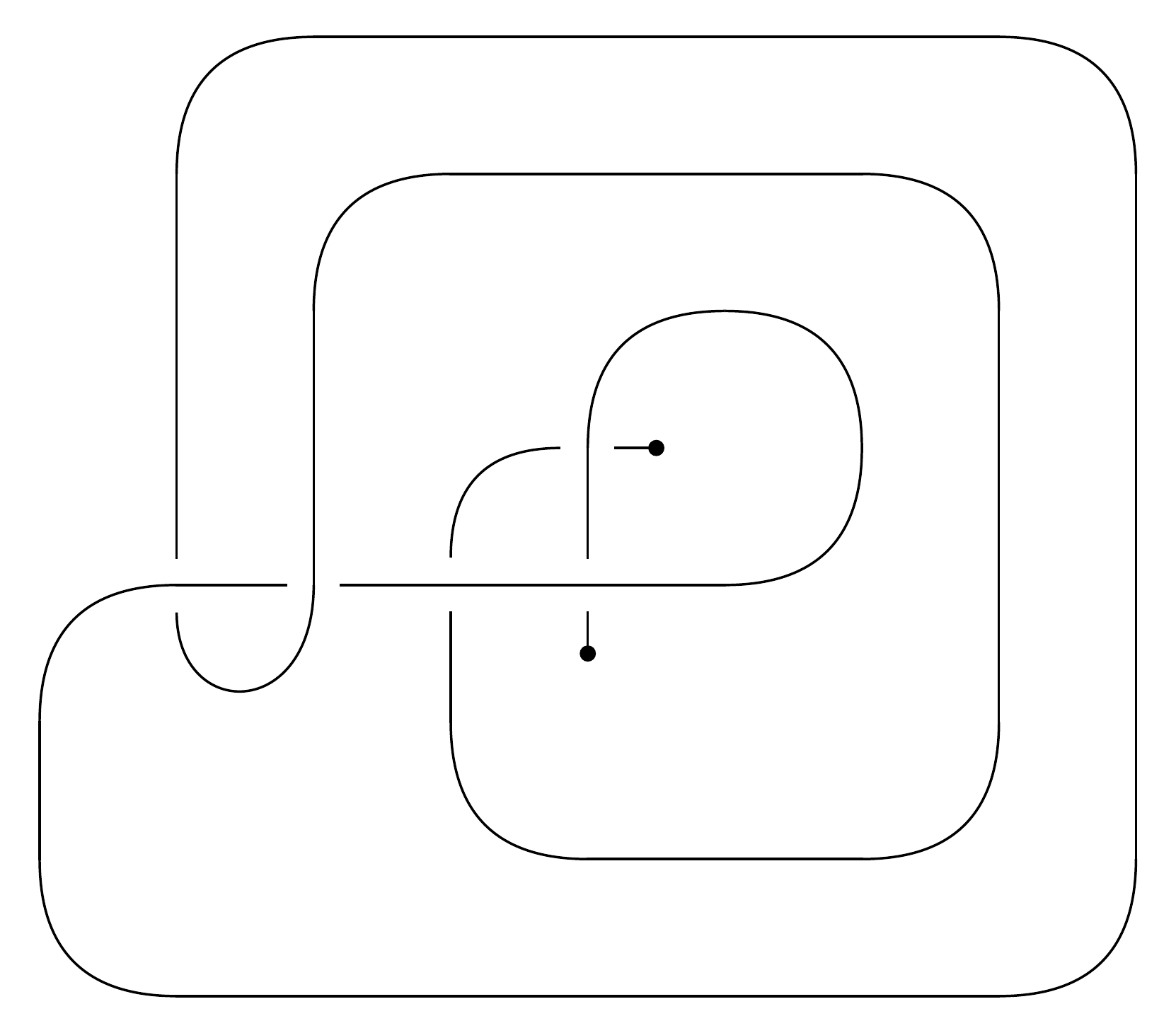}\\
\textcolor{black}{$5_{354}$}
\vspace{1cm}
\end{minipage}
\begin{minipage}[t]{.25\linewidth}
\centering
\includegraphics[width=0.9\textwidth,height=3.5cm,keepaspectratio]{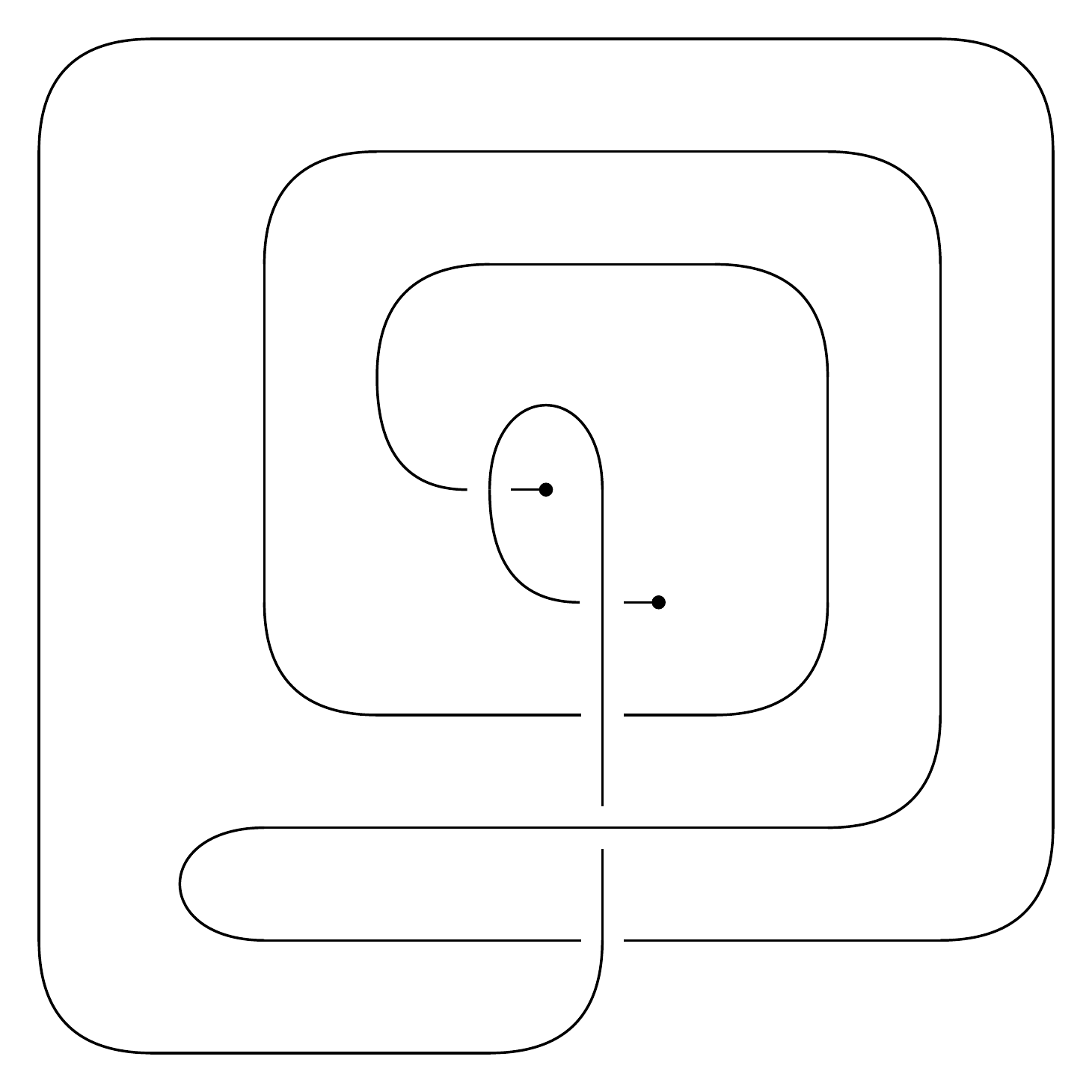}\\
\textcolor{black}{$5_{355}$}
\vspace{1cm}
\end{minipage}
\begin{minipage}[t]{.25\linewidth}
\centering
\includegraphics[width=0.9\textwidth,height=3.5cm,keepaspectratio]{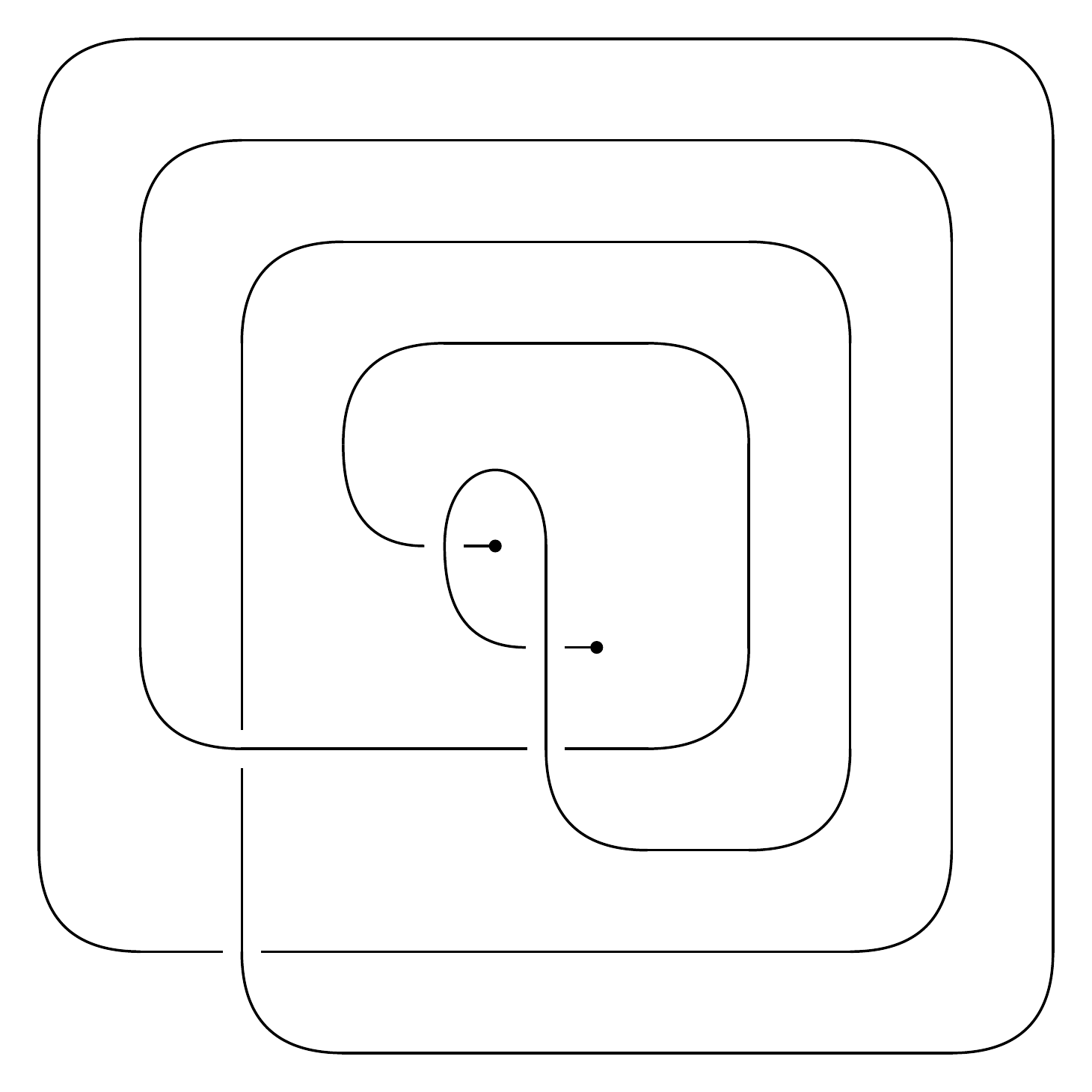}\\
\textcolor{black}{$5_{356}$}
\vspace{1cm}
\end{minipage}
\begin{minipage}[t]{.25\linewidth}
\centering
\includegraphics[width=0.9\textwidth,height=3.5cm,keepaspectratio]{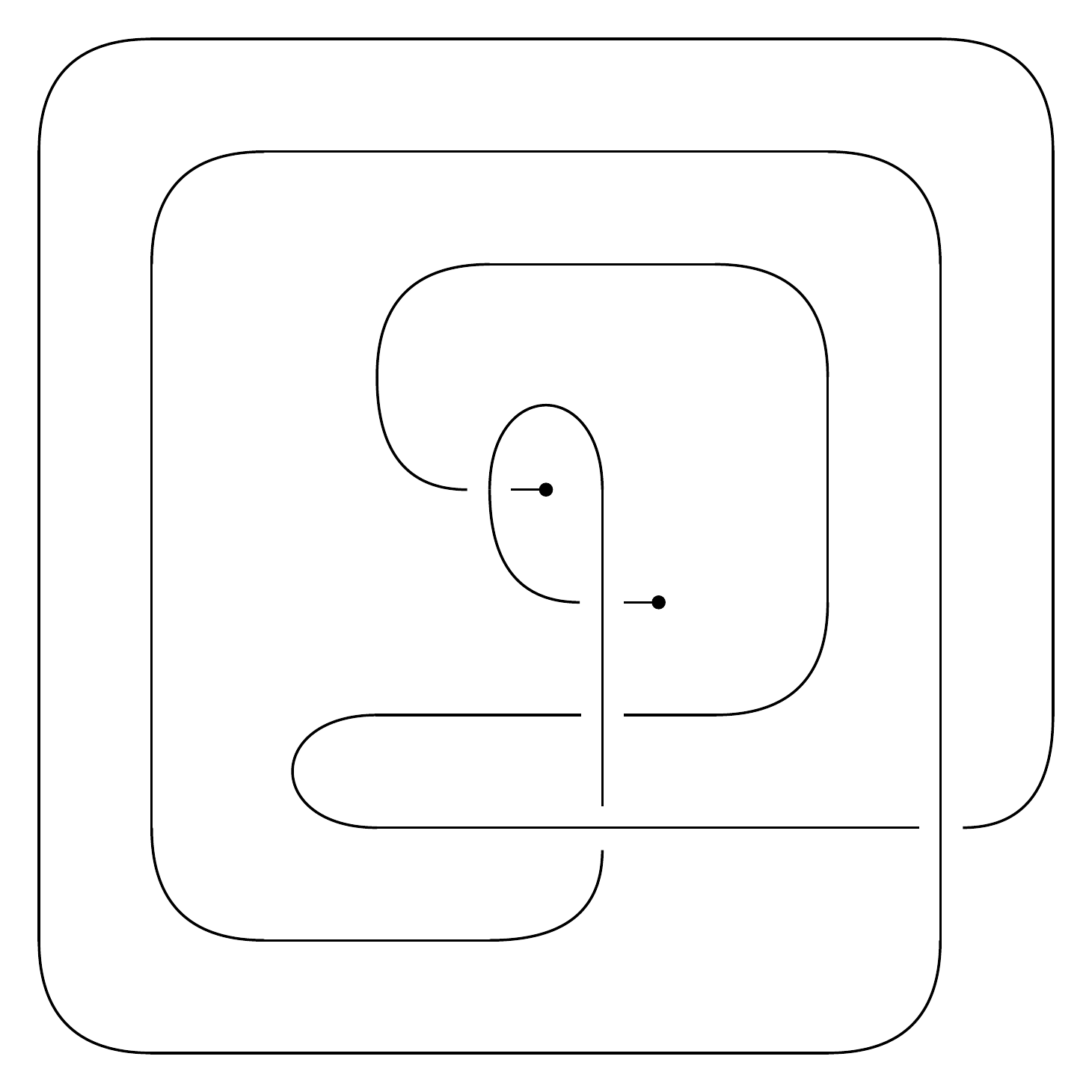}\\
\textcolor{black}{$5_{357}$}
\vspace{1cm}
\end{minipage}
\begin{minipage}[t]{.25\linewidth}
\centering
\includegraphics[width=0.9\textwidth,height=3.5cm,keepaspectratio]{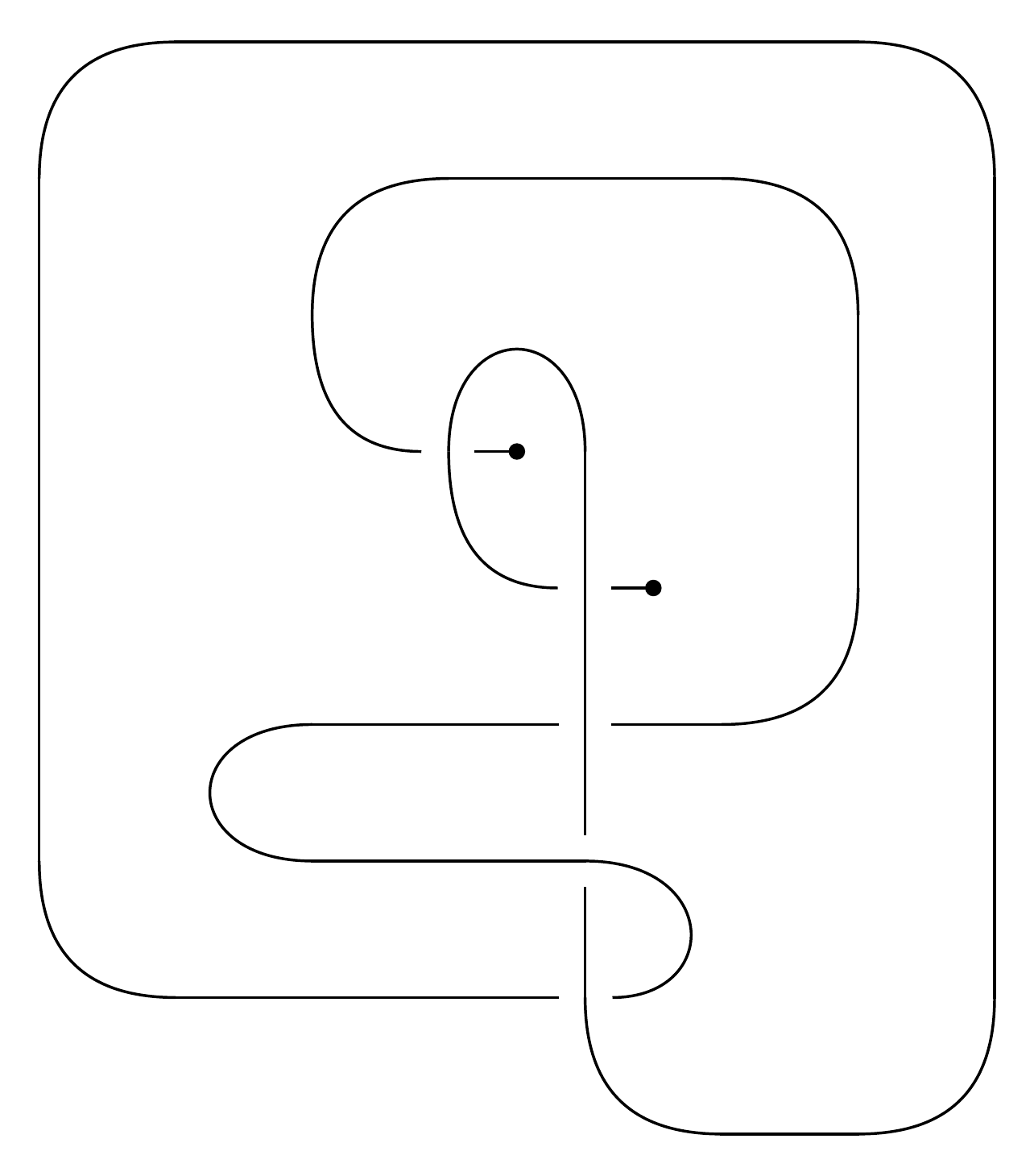}\\
\textcolor{black}{$5_{358}$}
\vspace{1cm}
\end{minipage}
\begin{minipage}[t]{.25\linewidth}
\centering
\includegraphics[width=0.9\textwidth,height=3.5cm,keepaspectratio]{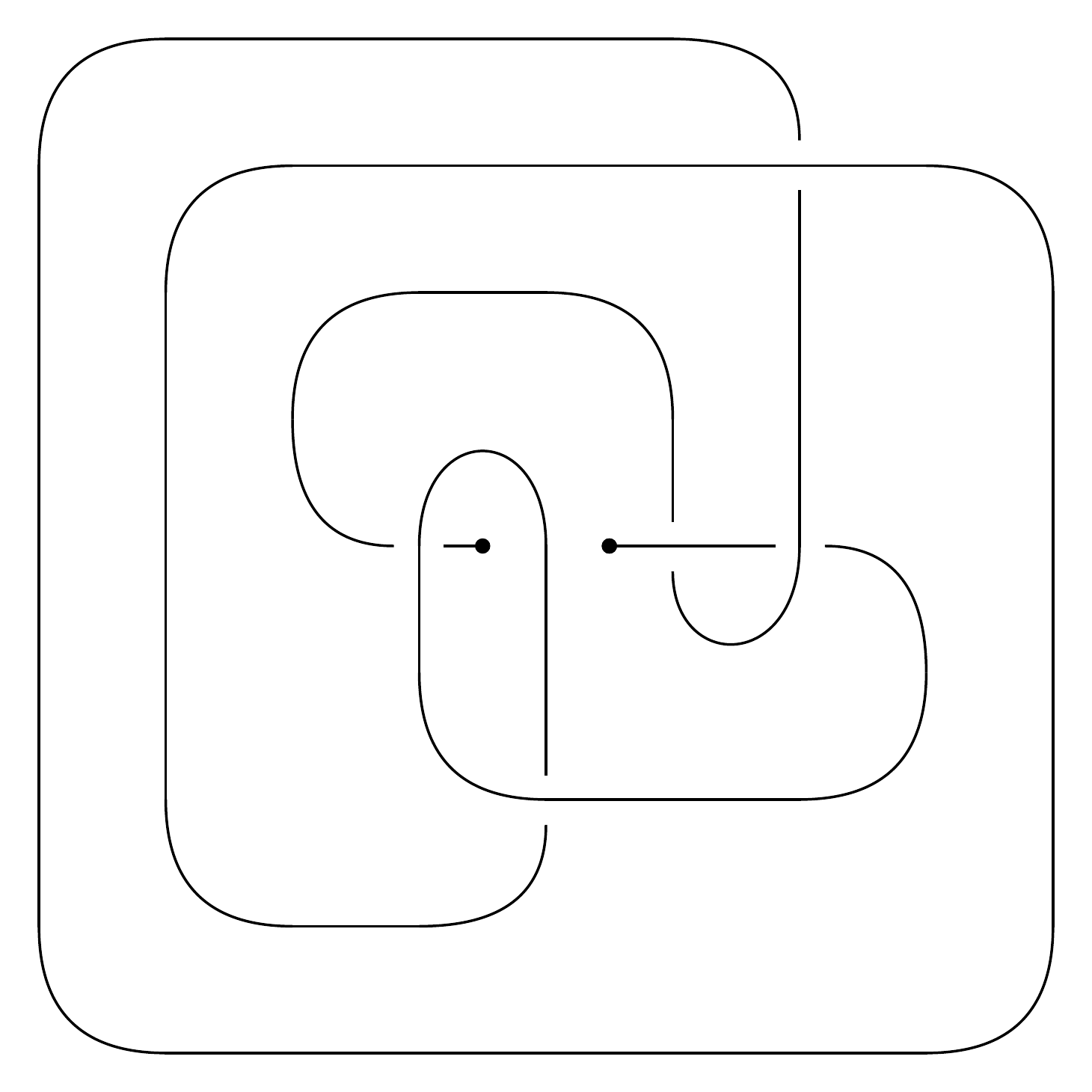}\\
\textcolor{black}{$5_{359}$}
\vspace{1cm}
\end{minipage}
\begin{minipage}[t]{.25\linewidth}
\centering
\includegraphics[width=0.9\textwidth,height=3.5cm,keepaspectratio]{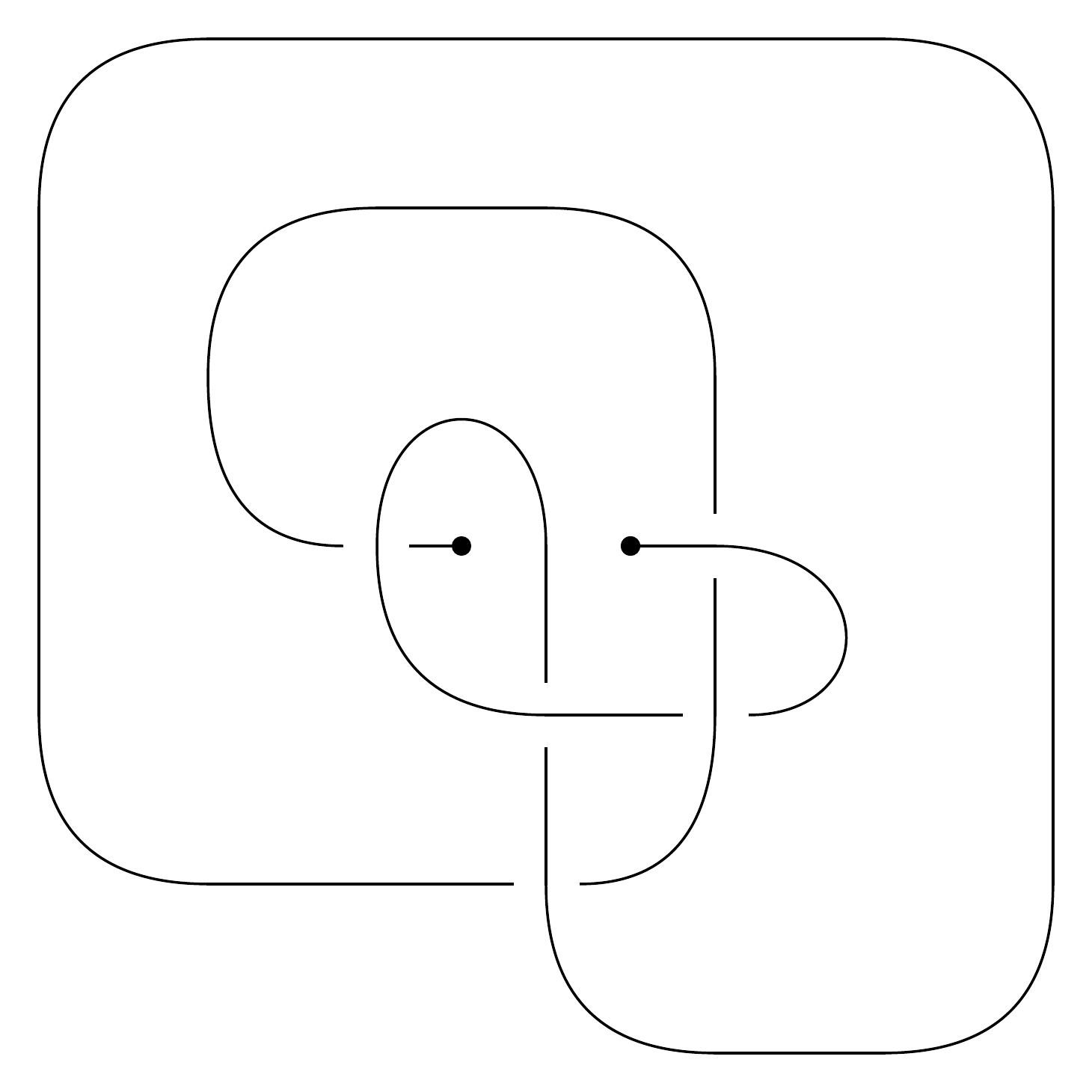}\\
\textcolor{black}{$5_{360}$}
\vspace{1cm}
\end{minipage}
\begin{minipage}[t]{.25\linewidth}
\centering
\includegraphics[width=0.9\textwidth,height=3.5cm,keepaspectratio]{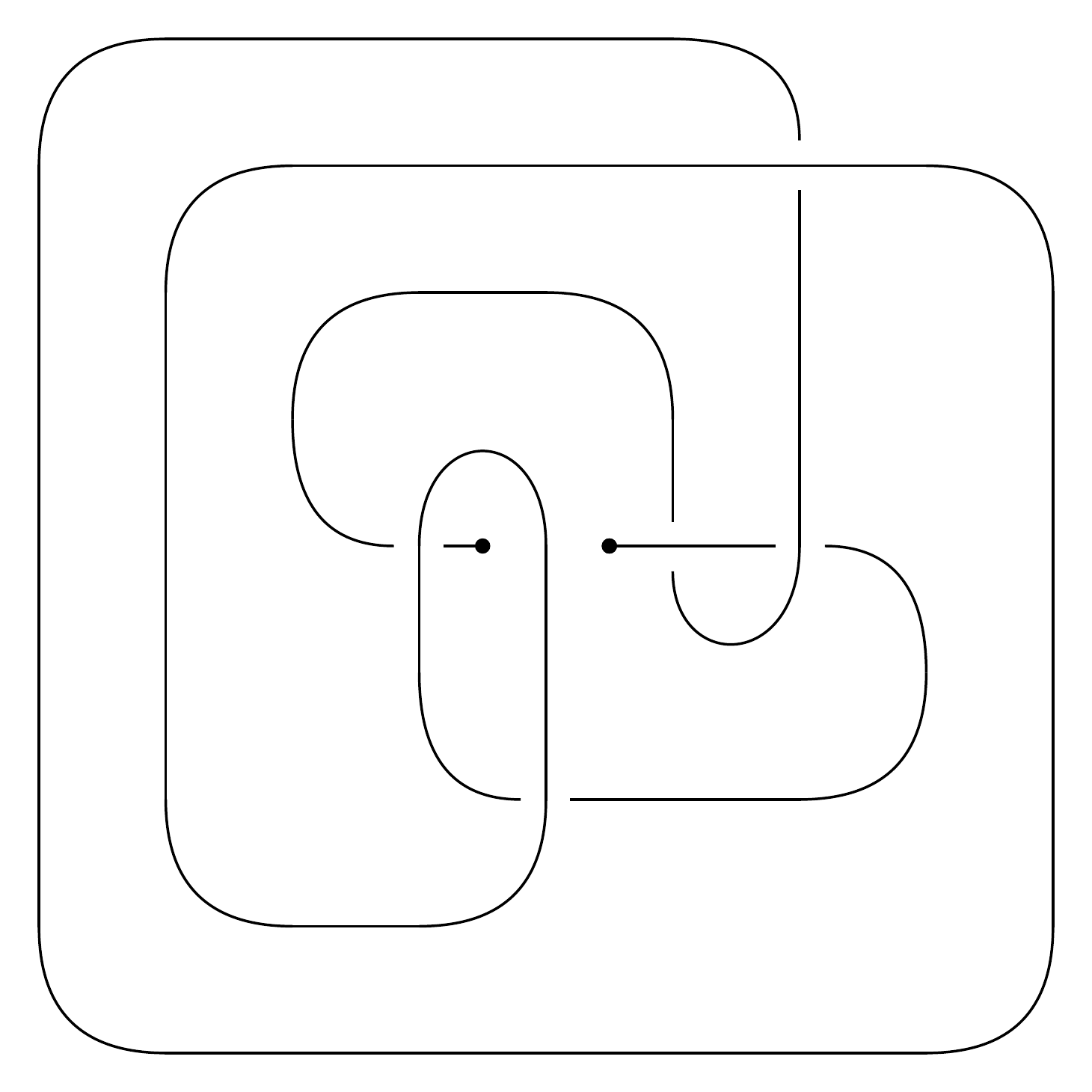}\\
\textcolor{black}{$5_{361}$}
\vspace{1cm}
\end{minipage}
\begin{minipage}[t]{.25\linewidth}
\centering
\includegraphics[width=0.9\textwidth,height=3.5cm,keepaspectratio]{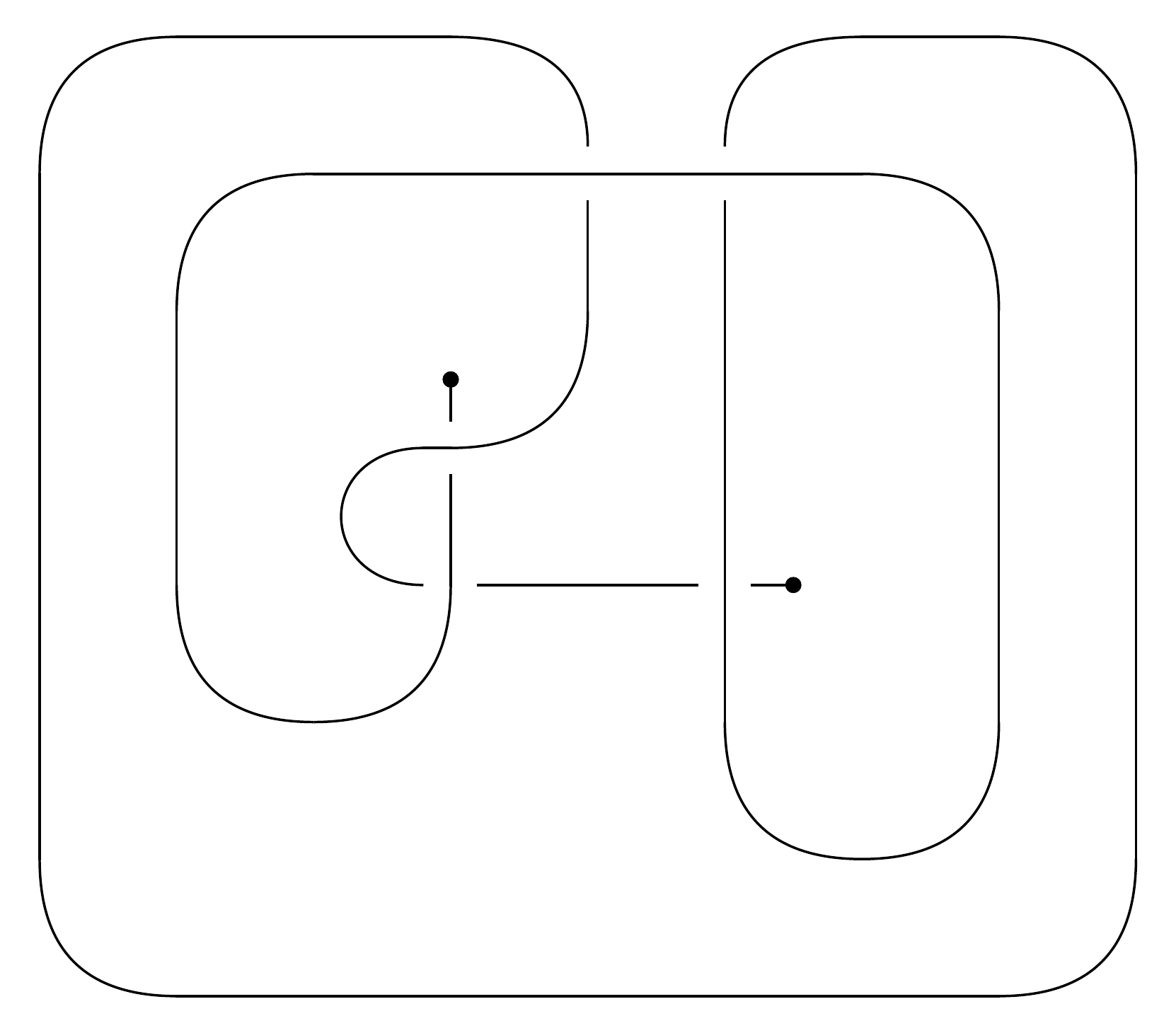}\\
\textcolor{black}{$5_{362}$}
\vspace{1cm}
\end{minipage}
\begin{minipage}[t]{.25\linewidth}
\centering
\includegraphics[width=0.9\textwidth,height=3.5cm,keepaspectratio]{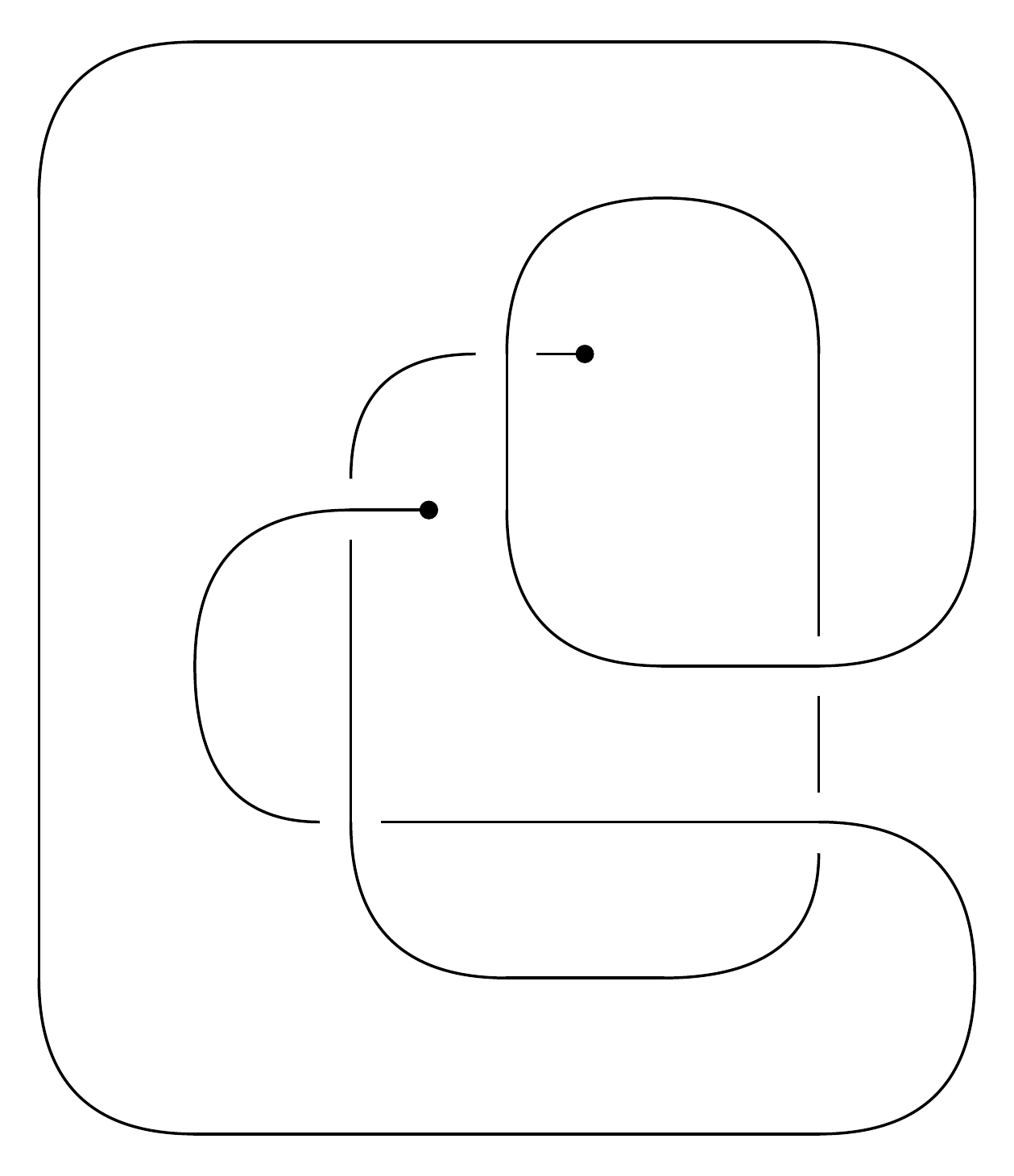}\\
\textcolor{black}{$5_{363}$}
\vspace{1cm}
\end{minipage}
\begin{minipage}[t]{.25\linewidth}
\centering
\includegraphics[width=0.9\textwidth,height=3.5cm,keepaspectratio]{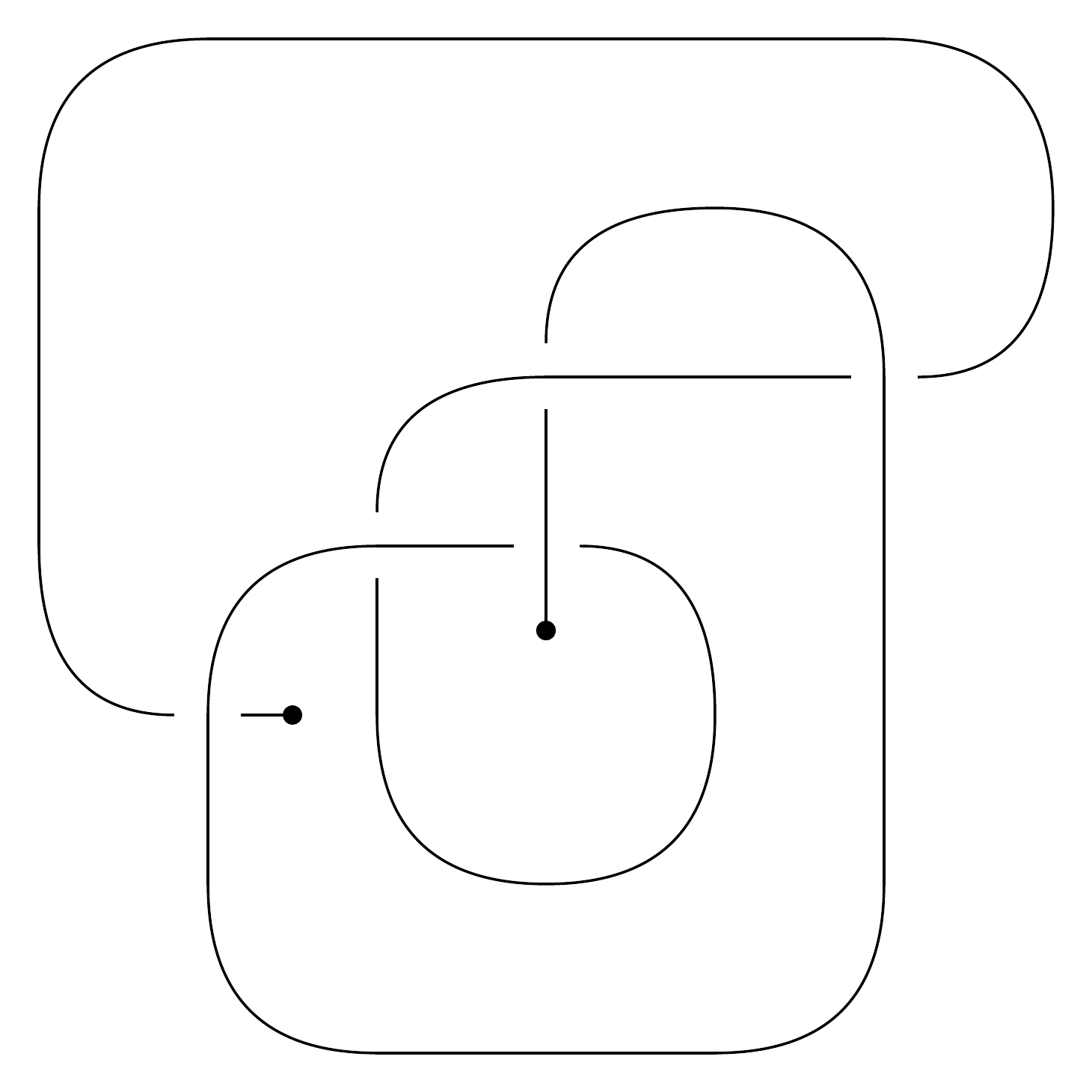}\\
\textcolor{black}{$5_{364}$}
\vspace{1cm}
\end{minipage}
\begin{minipage}[t]{.25\linewidth}
\centering
\includegraphics[width=0.9\textwidth,height=3.5cm,keepaspectratio]{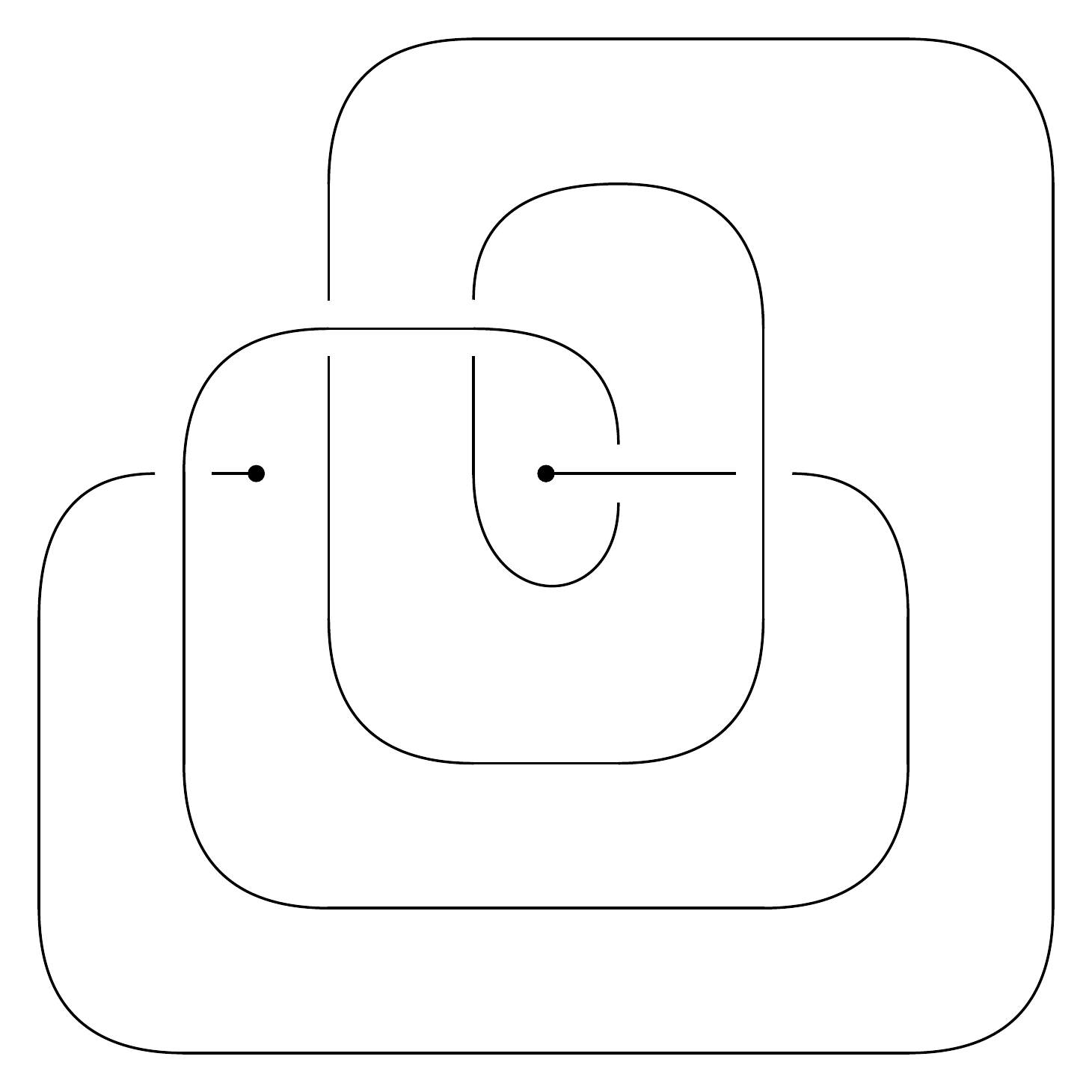}\\
\textcolor{black}{$5_{365}$}
\vspace{1cm}
\end{minipage}
\begin{minipage}[t]{.25\linewidth}
\centering
\includegraphics[width=0.9\textwidth,height=3.5cm,keepaspectratio]{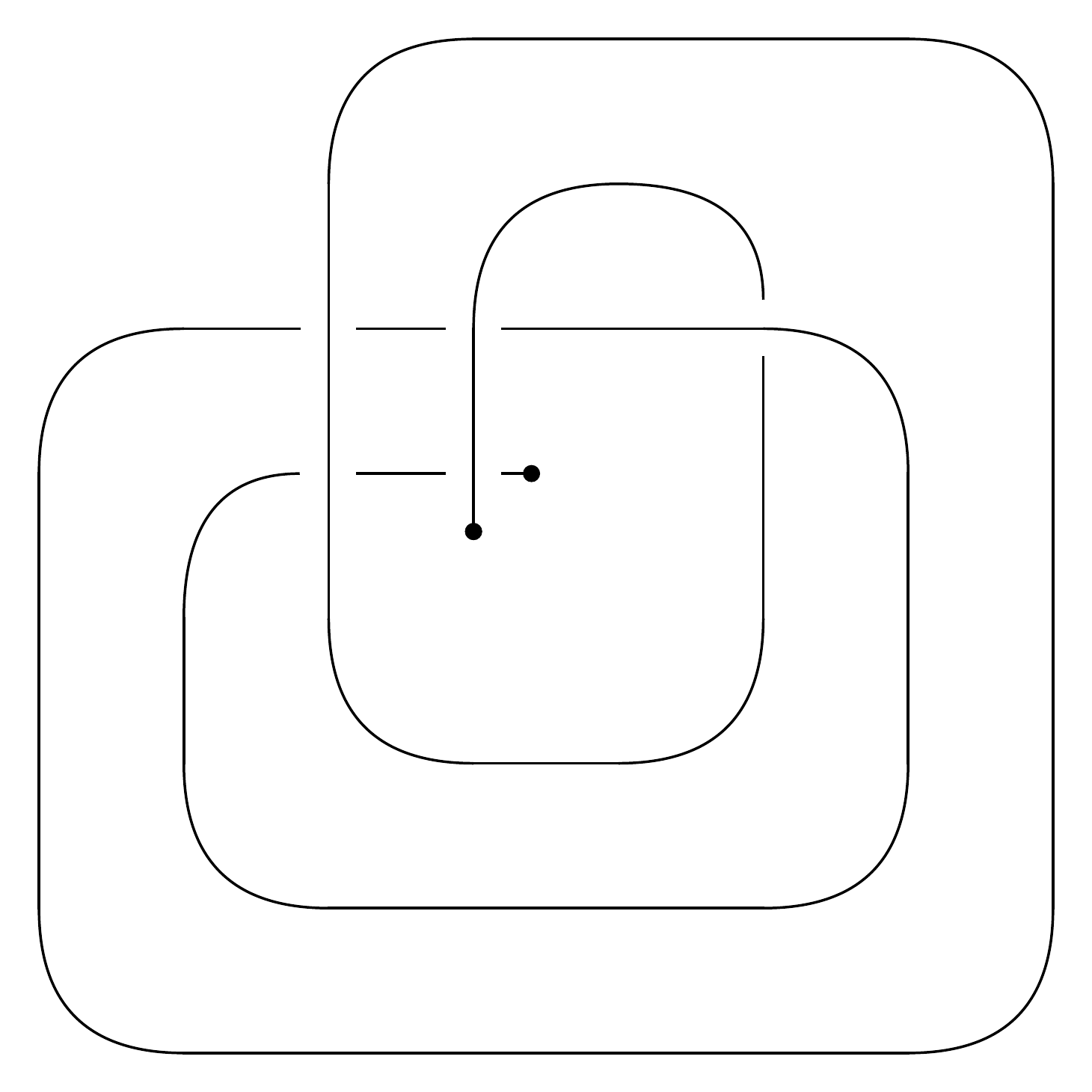}\\
\textcolor{black}{$5_{366}$}
\vspace{1cm}
\end{minipage}
\begin{minipage}[t]{.25\linewidth}
\centering
\includegraphics[width=0.9\textwidth,height=3.5cm,keepaspectratio]{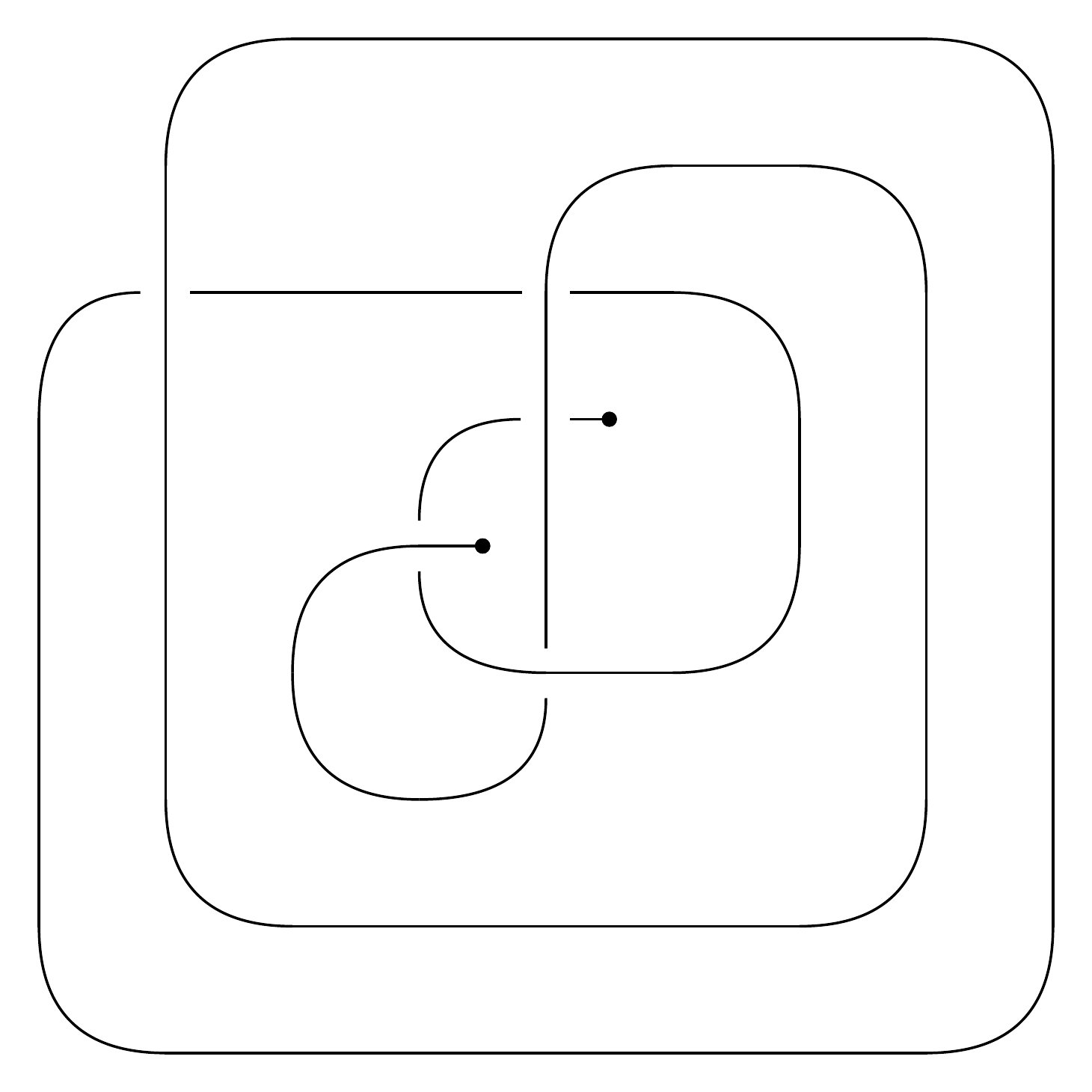}\\
\textcolor{black}{$5_{367}$}
\vspace{1cm}
\end{minipage}
\begin{minipage}[t]{.25\linewidth}
\centering
\includegraphics[width=0.9\textwidth,height=3.5cm,keepaspectratio]{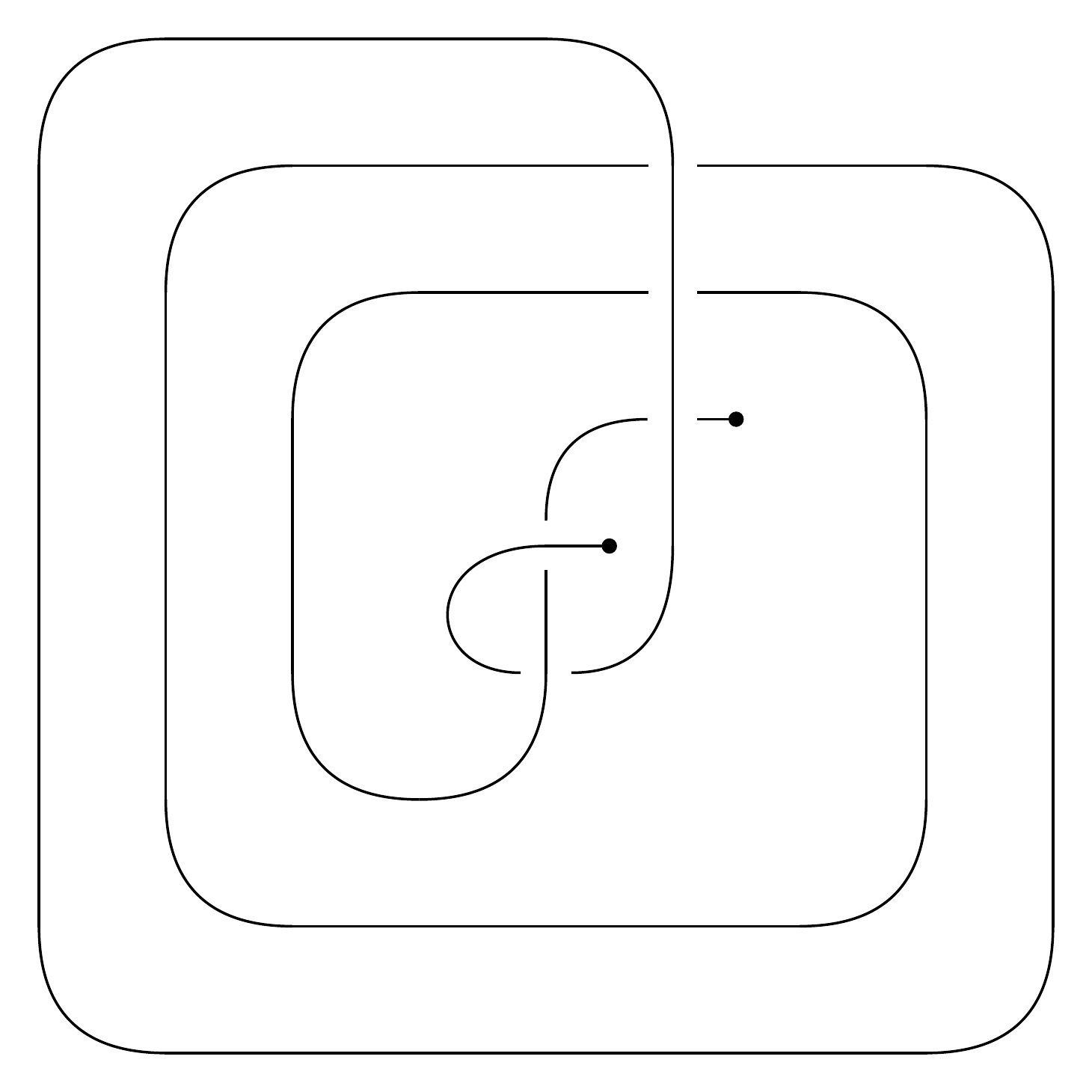}\\
\textcolor{black}{$5_{368}$}
\vspace{1cm}
\end{minipage}
\begin{minipage}[t]{.25\linewidth}
\centering
\includegraphics[width=0.9\textwidth,height=3.5cm,keepaspectratio]{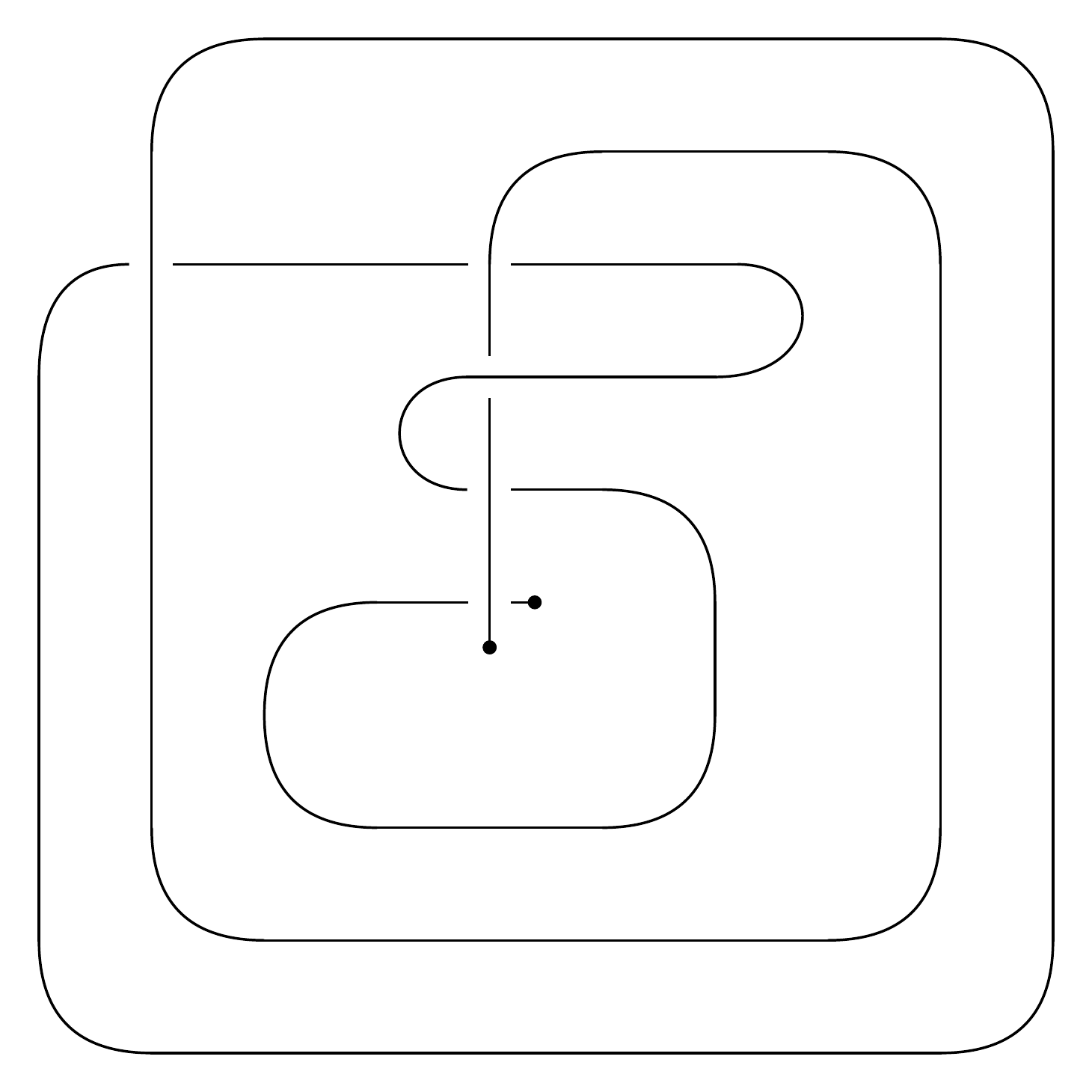}\\
\textcolor{black}{$5_{369}$}
\vspace{1cm}
\end{minipage}
\begin{minipage}[t]{.25\linewidth}
\centering
\includegraphics[width=0.9\textwidth,height=3.5cm,keepaspectratio]{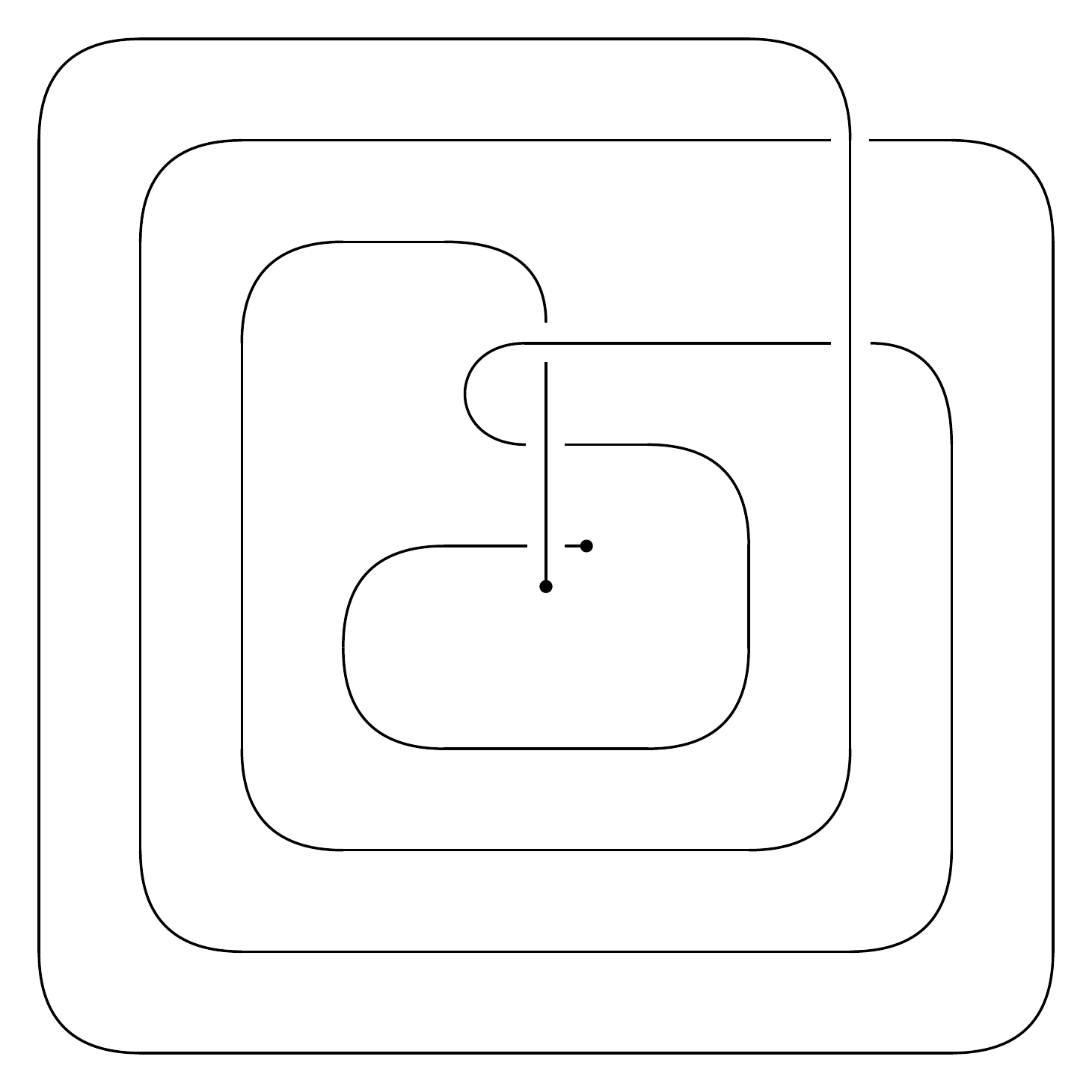}\\
\textcolor{black}{$5_{370}$}
\vspace{1cm}
\end{minipage}
\begin{minipage}[t]{.25\linewidth}
\centering
\includegraphics[width=0.9\textwidth,height=3.5cm,keepaspectratio]{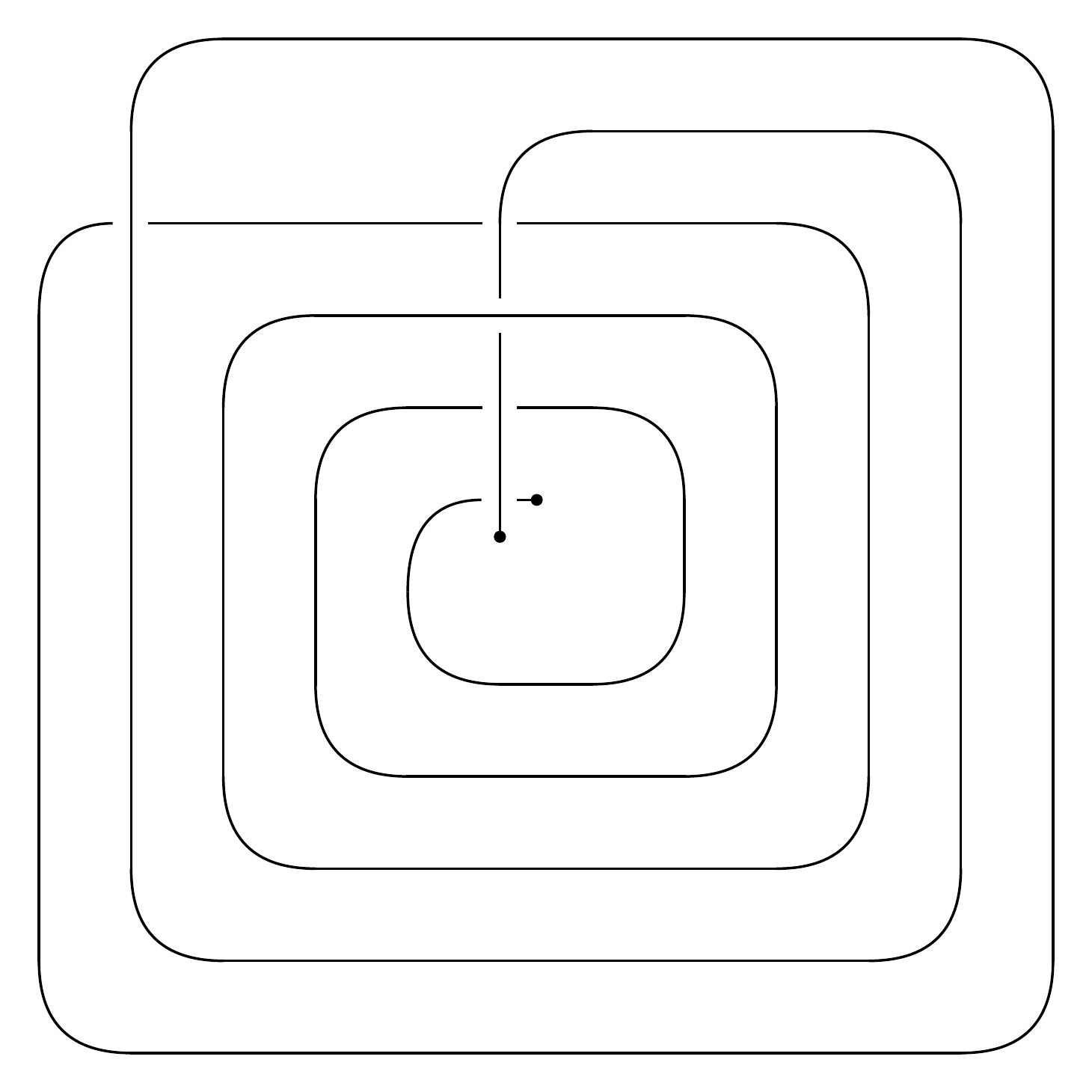}\\
\textcolor{black}{$5_{371}$}
\vspace{1cm}
\end{minipage}
\begin{minipage}[t]{.25\linewidth}
\centering
\includegraphics[width=0.9\textwidth,height=3.5cm,keepaspectratio]{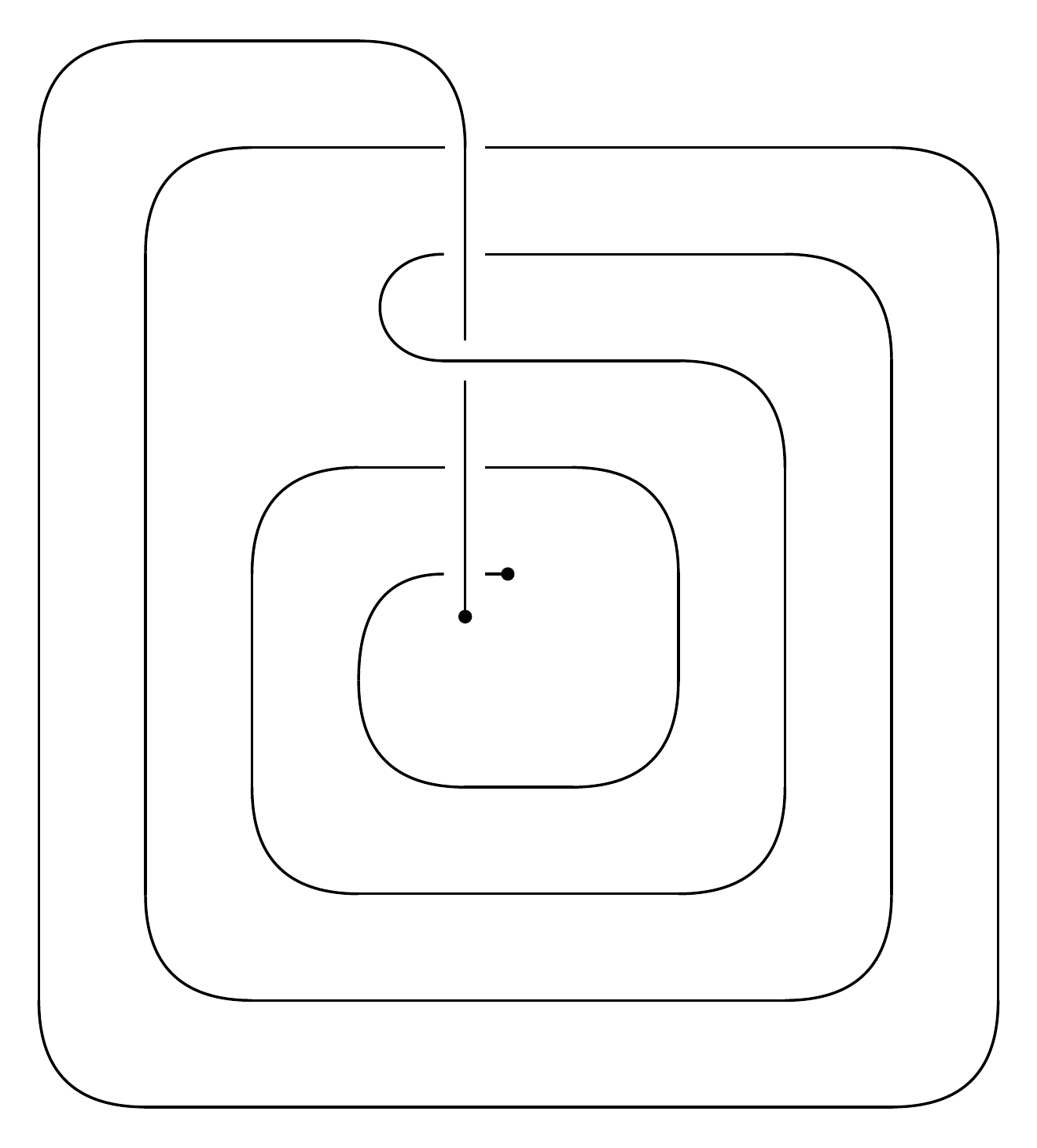}\\
\textcolor{black}{$5_{372}$}
\vspace{1cm}
\end{minipage}
\begin{minipage}[t]{.25\linewidth}
\centering
\includegraphics[width=0.9\textwidth,height=3.5cm,keepaspectratio]{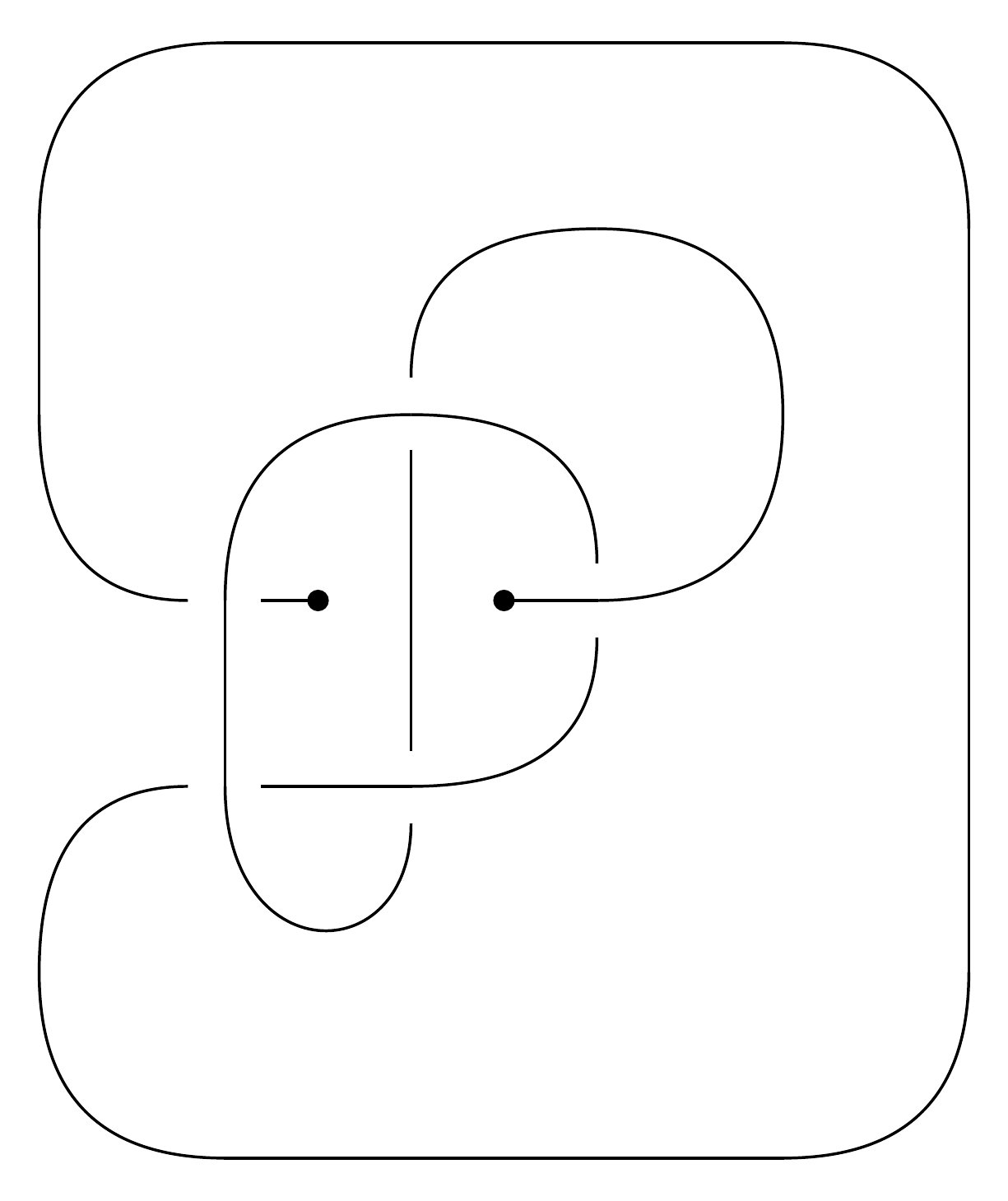}\\
\textcolor{black}{$5_{373}$}
\vspace{1cm}
\end{minipage}
\begin{minipage}[t]{.25\linewidth}
\centering
\includegraphics[width=0.9\textwidth,height=3.5cm,keepaspectratio]{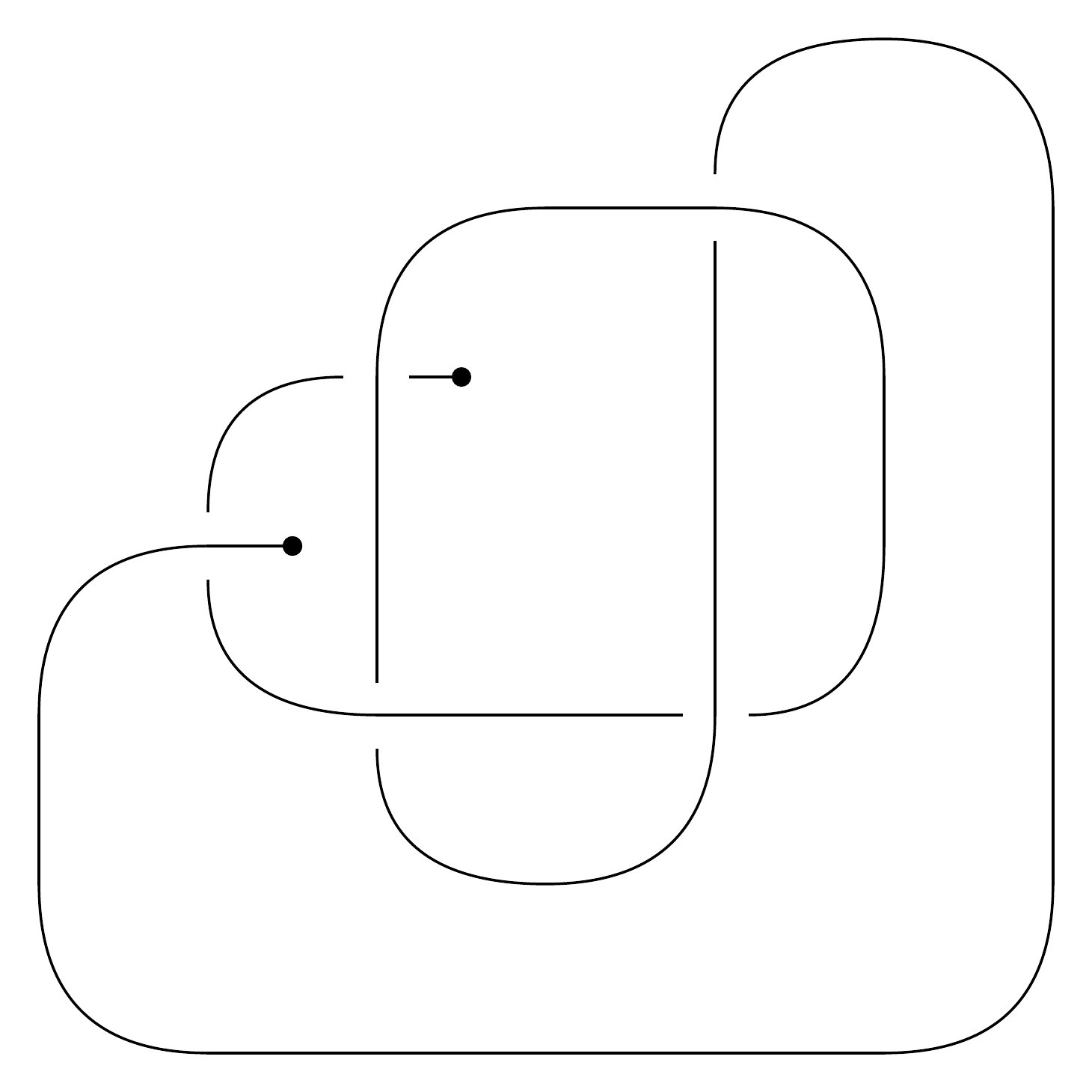}\\
\textcolor{black}{$5_{374}$}
\vspace{1cm}
\end{minipage}
\begin{minipage}[t]{.25\linewidth}
\centering
\includegraphics[width=0.9\textwidth,height=3.5cm,keepaspectratio]{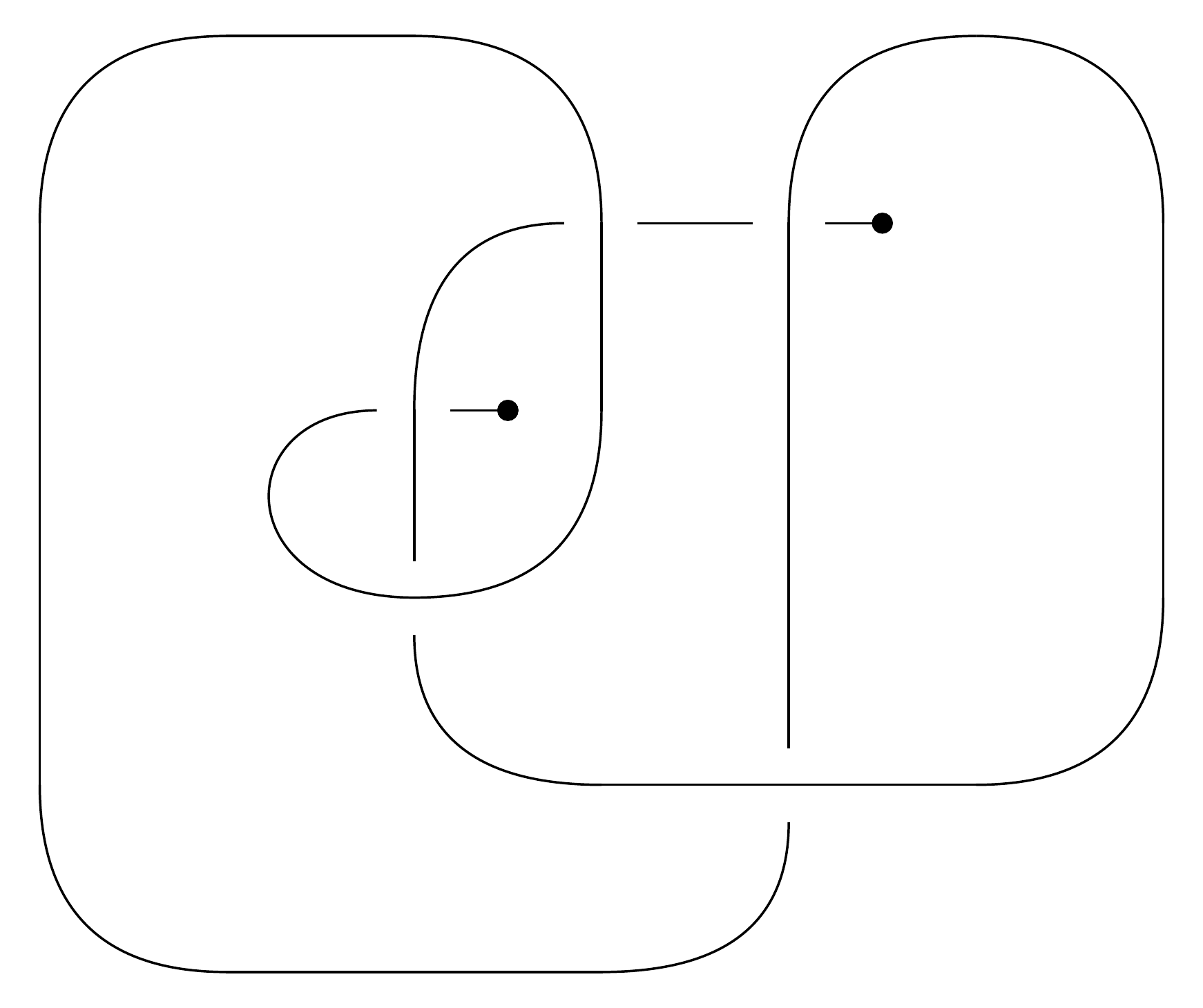}\\
\textcolor{black}{$5_{375}$}
\vspace{1cm}
\end{minipage}
\begin{minipage}[t]{.25\linewidth}
\centering
\includegraphics[width=0.9\textwidth,height=3.5cm,keepaspectratio]{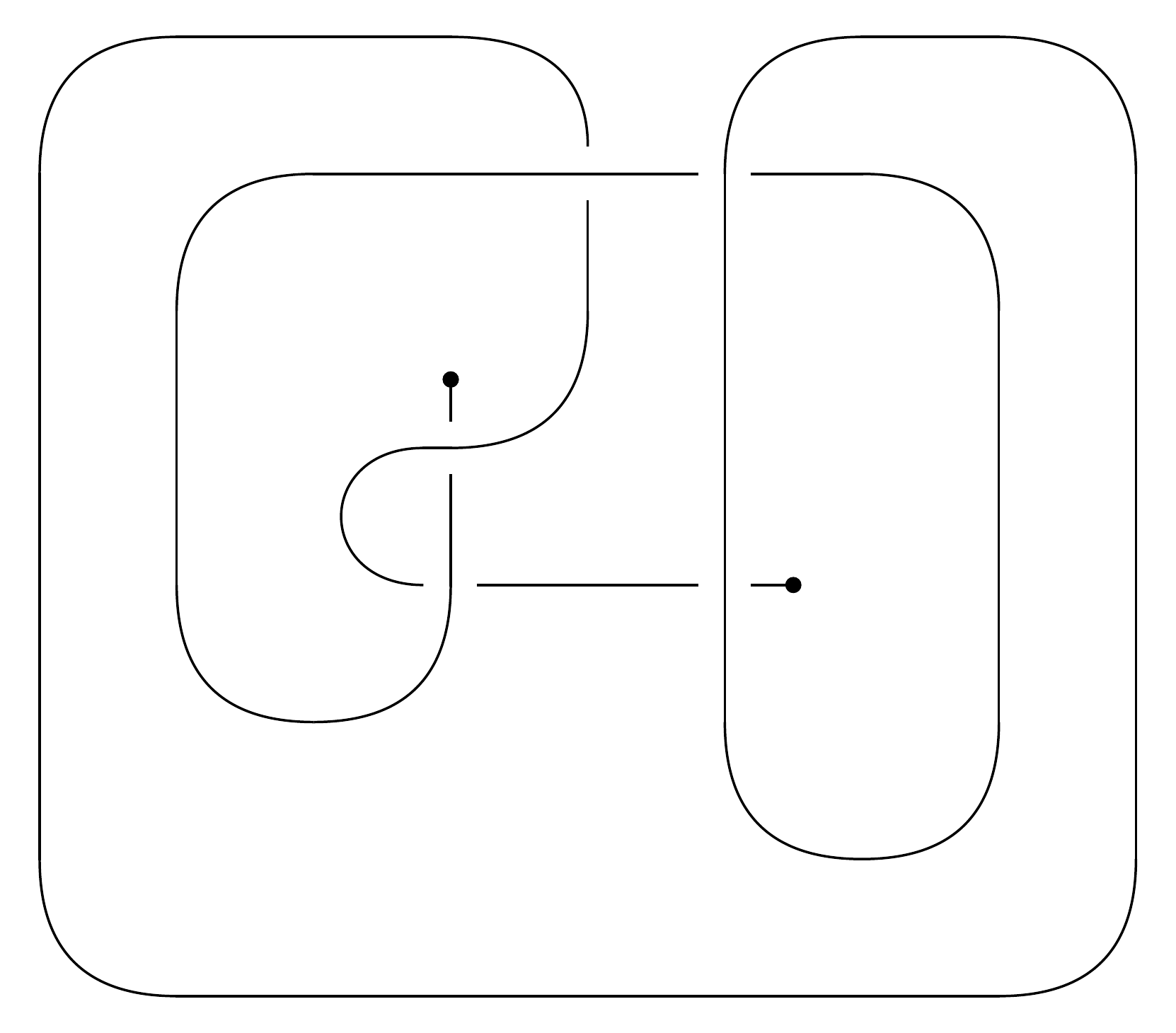}\\
\textcolor{black}{$5_{376}$}
\vspace{1cm}
\end{minipage}
\begin{minipage}[t]{.25\linewidth}
\centering
\includegraphics[width=0.9\textwidth,height=3.5cm,keepaspectratio]{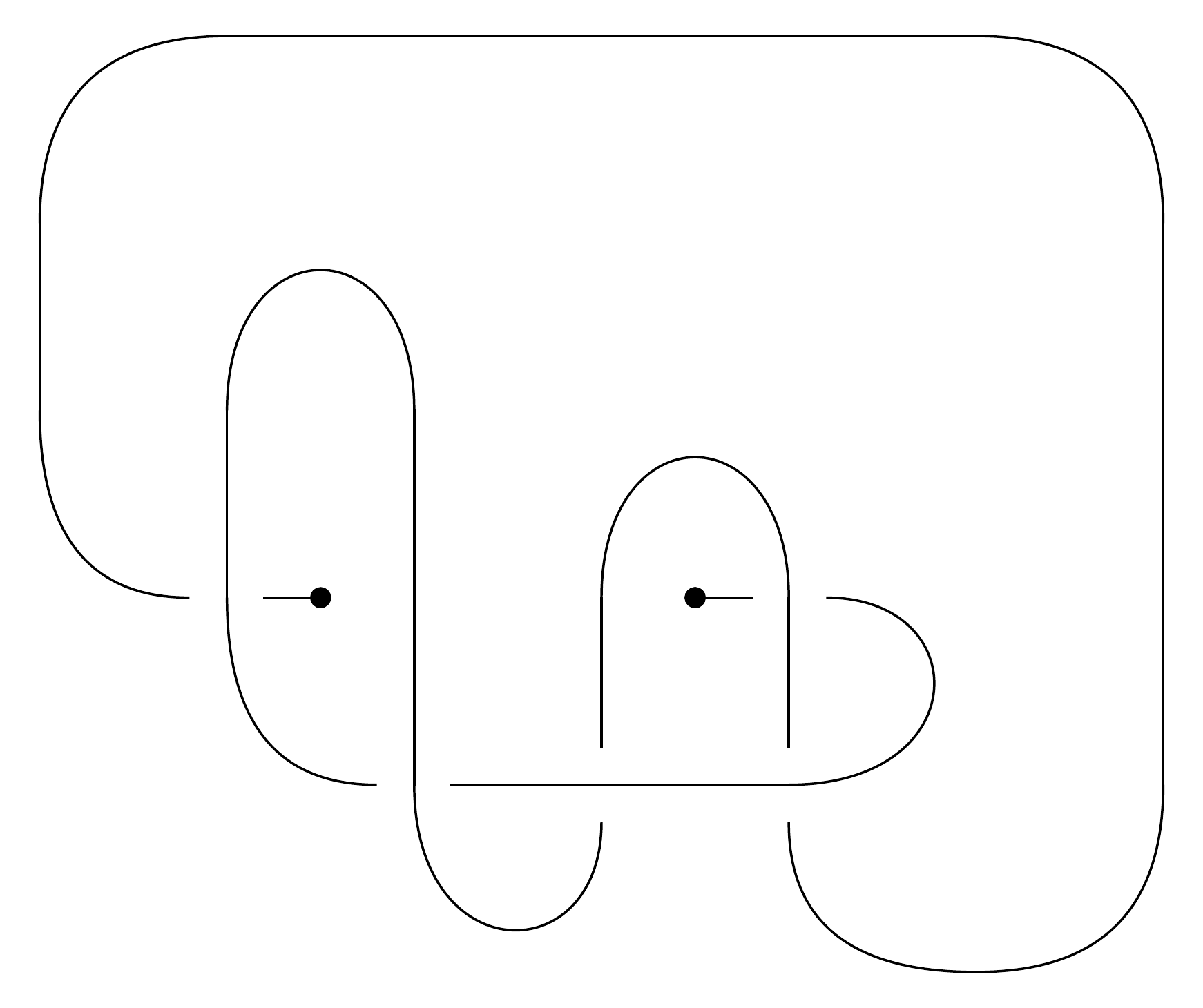}\\
\textcolor{black}{$5_{377}$}
\vspace{1cm}
\end{minipage}
\begin{minipage}[t]{.25\linewidth}
\centering
\includegraphics[width=0.9\textwidth,height=3.5cm,keepaspectratio]{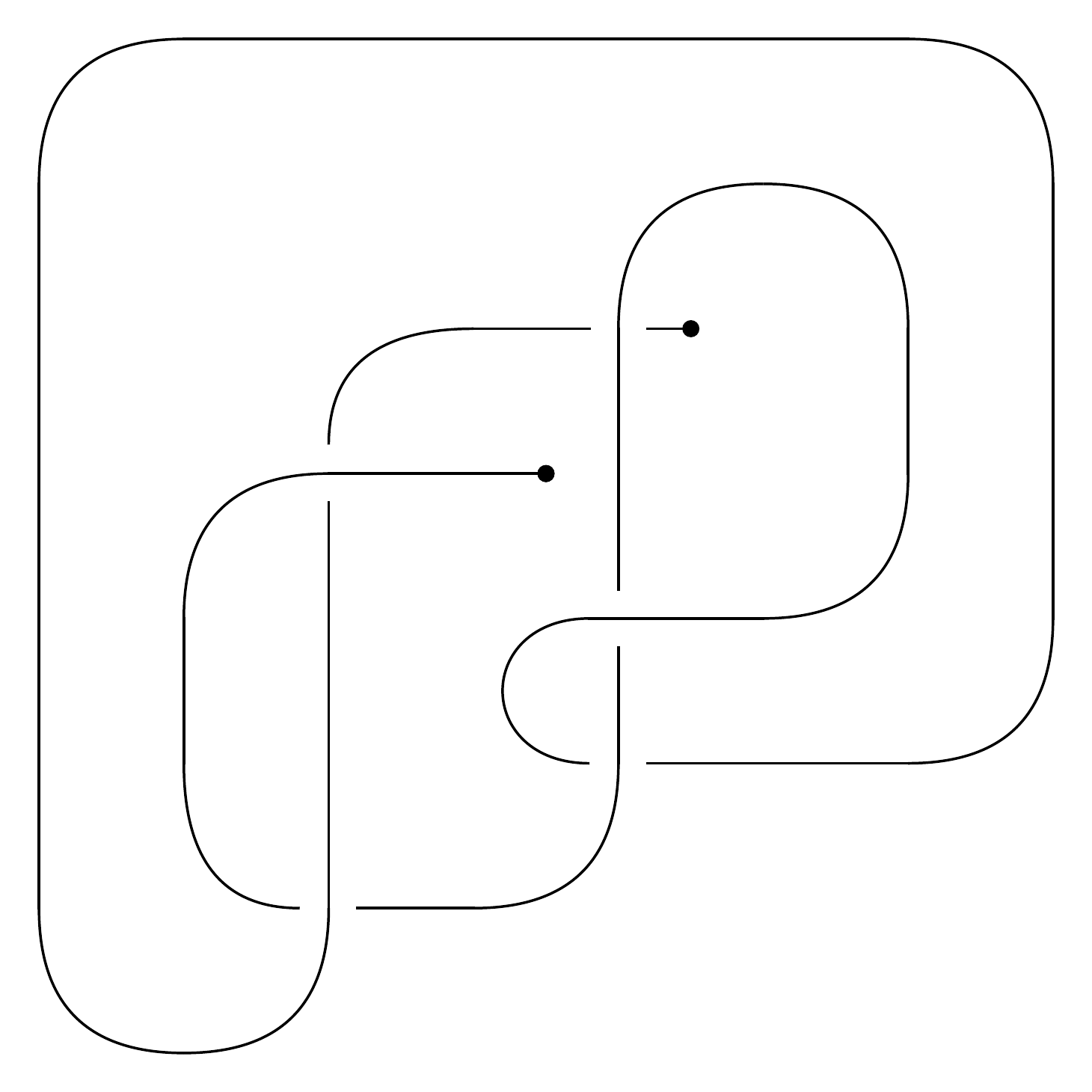}\\
\textcolor{black}{$5_{378}$}
\vspace{1cm}
\end{minipage}
\begin{minipage}[t]{.25\linewidth}
\centering
\includegraphics[width=0.9\textwidth,height=3.5cm,keepaspectratio]{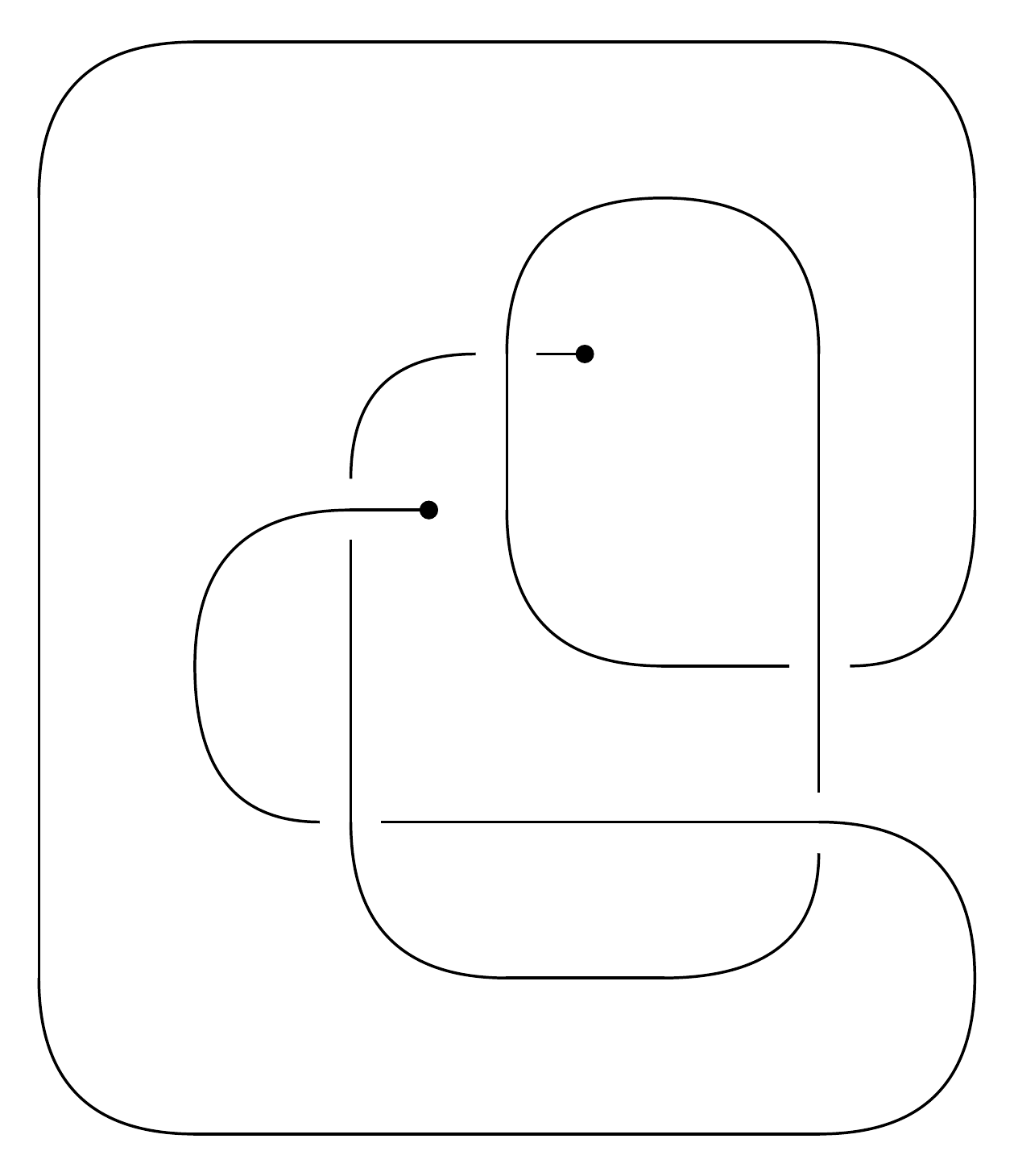}\\
\textcolor{black}{$5_{379}$}
\vspace{1cm}
\end{minipage}
\begin{minipage}[t]{.25\linewidth}
\centering
\includegraphics[width=0.9\textwidth,height=3.5cm,keepaspectratio]{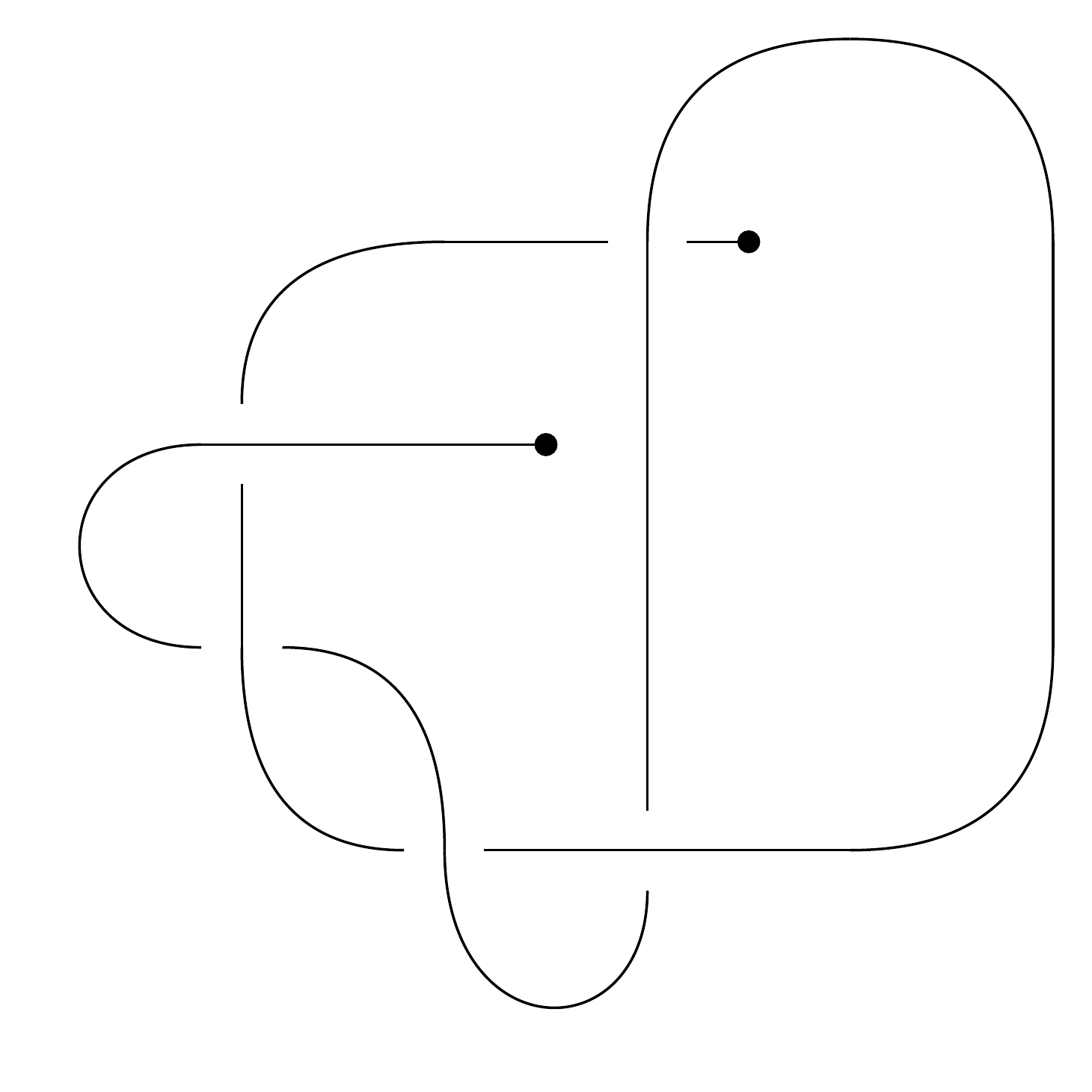}\\
\textcolor{black}{$5_{380}$}
\vspace{1cm}
\end{minipage}
\begin{minipage}[t]{.25\linewidth}
\centering
\includegraphics[width=0.9\textwidth,height=3.5cm,keepaspectratio]{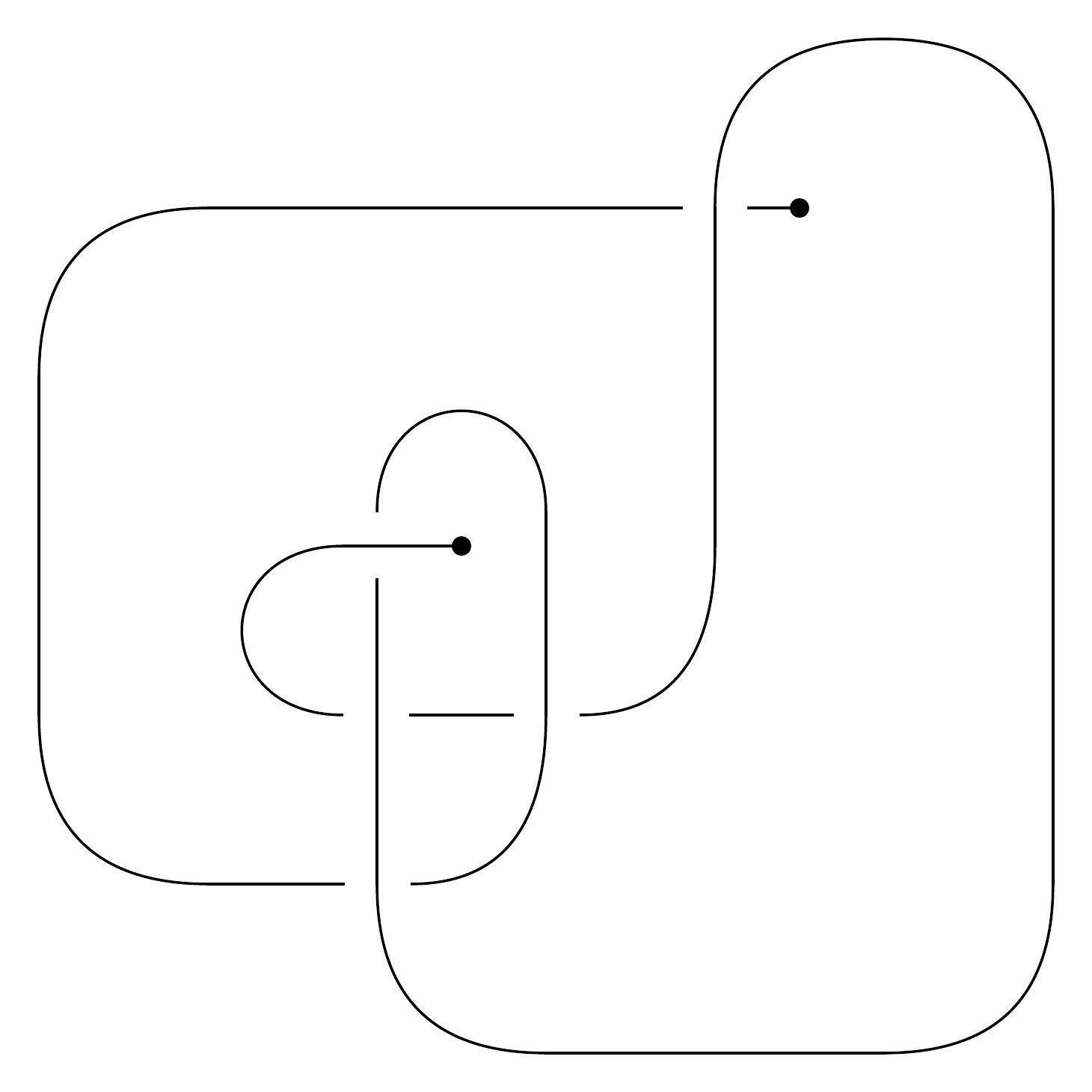}\\
\textcolor{black}{$5_{381}$}
\vspace{1cm}
\end{minipage}
\begin{minipage}[t]{.25\linewidth}
\centering
\includegraphics[width=0.9\textwidth,height=3.5cm,keepaspectratio]{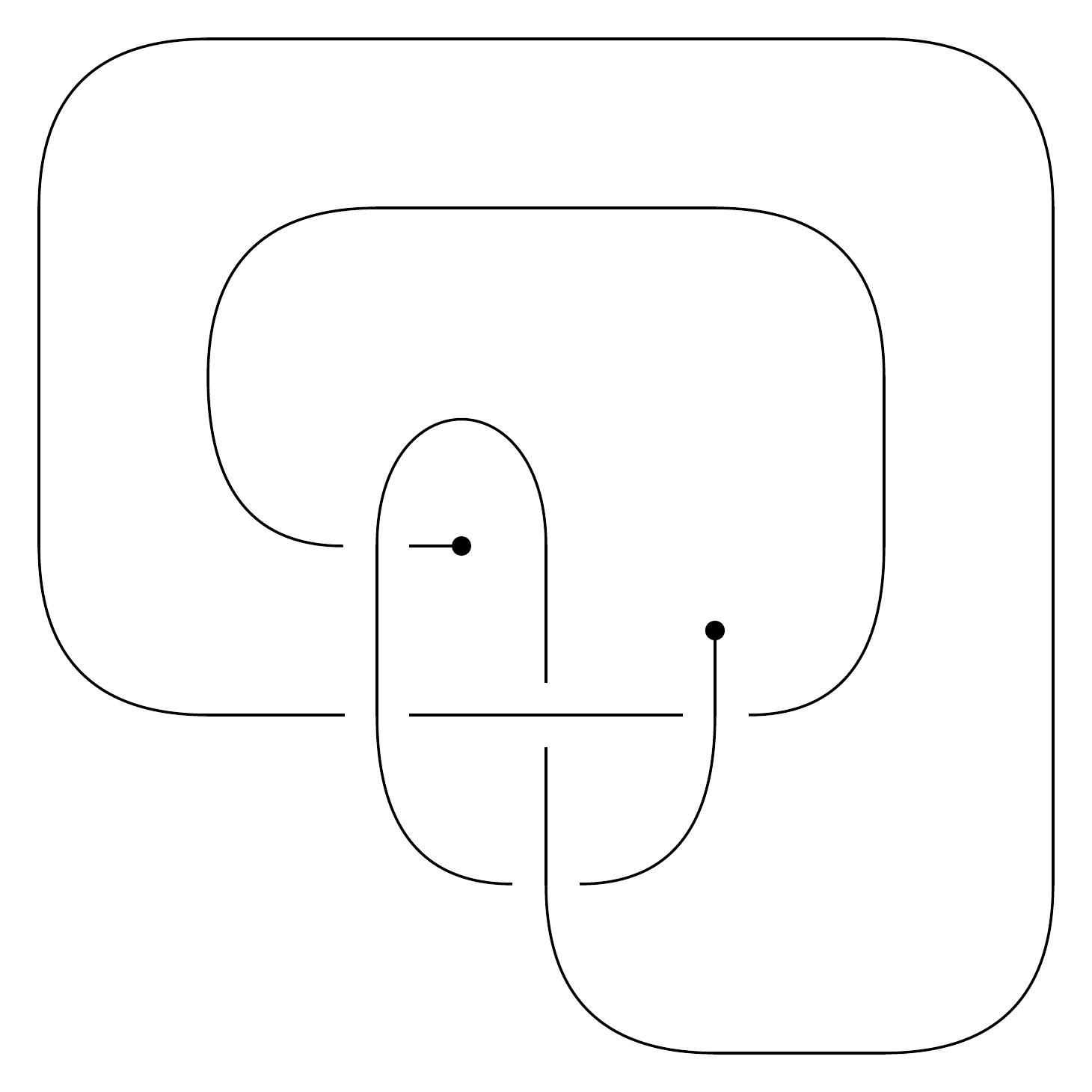}\\
\textcolor{black}{$5_{382}$}
\vspace{1cm}
\end{minipage}
\begin{minipage}[t]{.25\linewidth}
\centering
\includegraphics[width=0.9\textwidth,height=3.5cm,keepaspectratio]{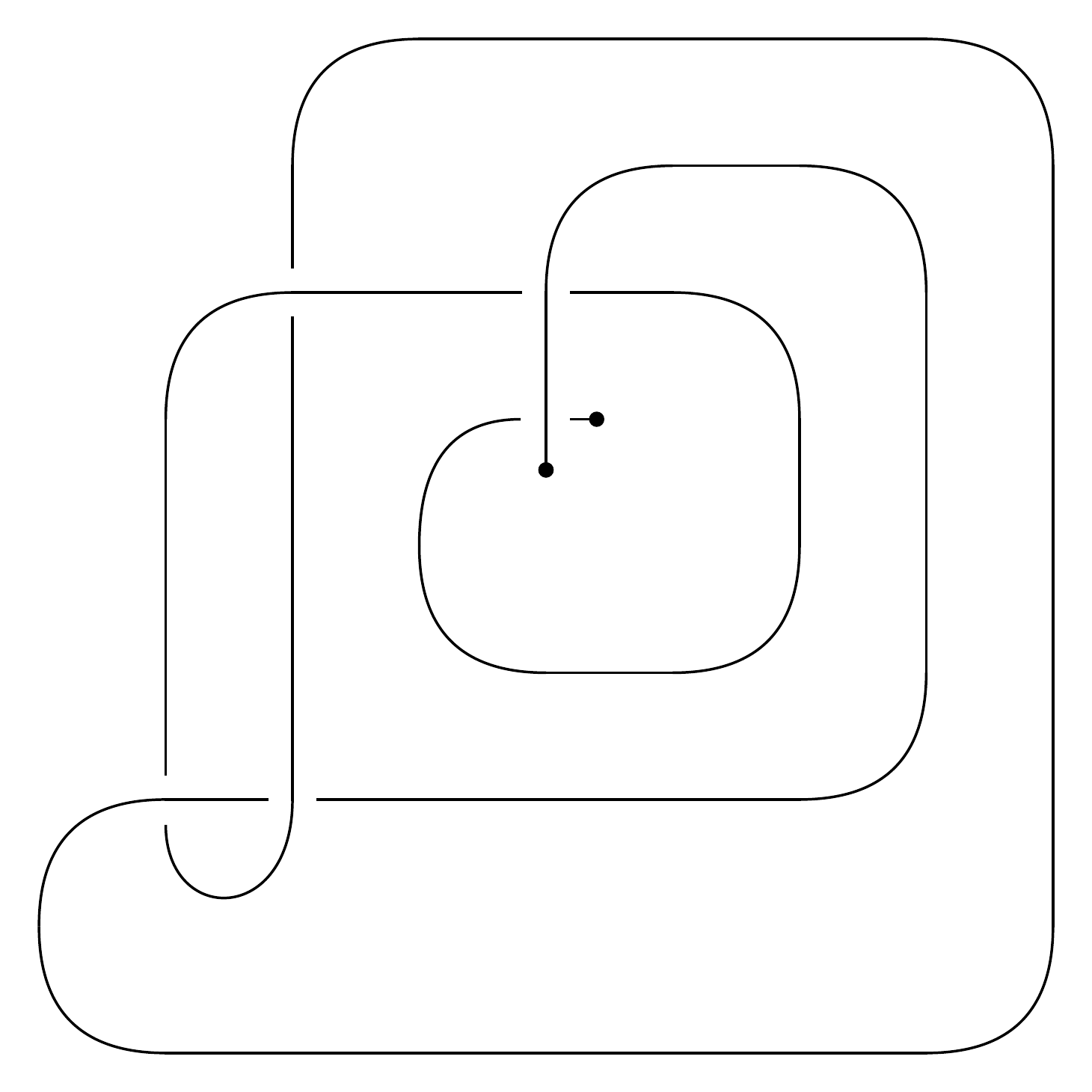}\\
\textcolor{black}{$5_{383}$}
\vspace{1cm}
\end{minipage}
\begin{minipage}[t]{.25\linewidth}
\centering
\includegraphics[width=0.9\textwidth,height=3.5cm,keepaspectratio]{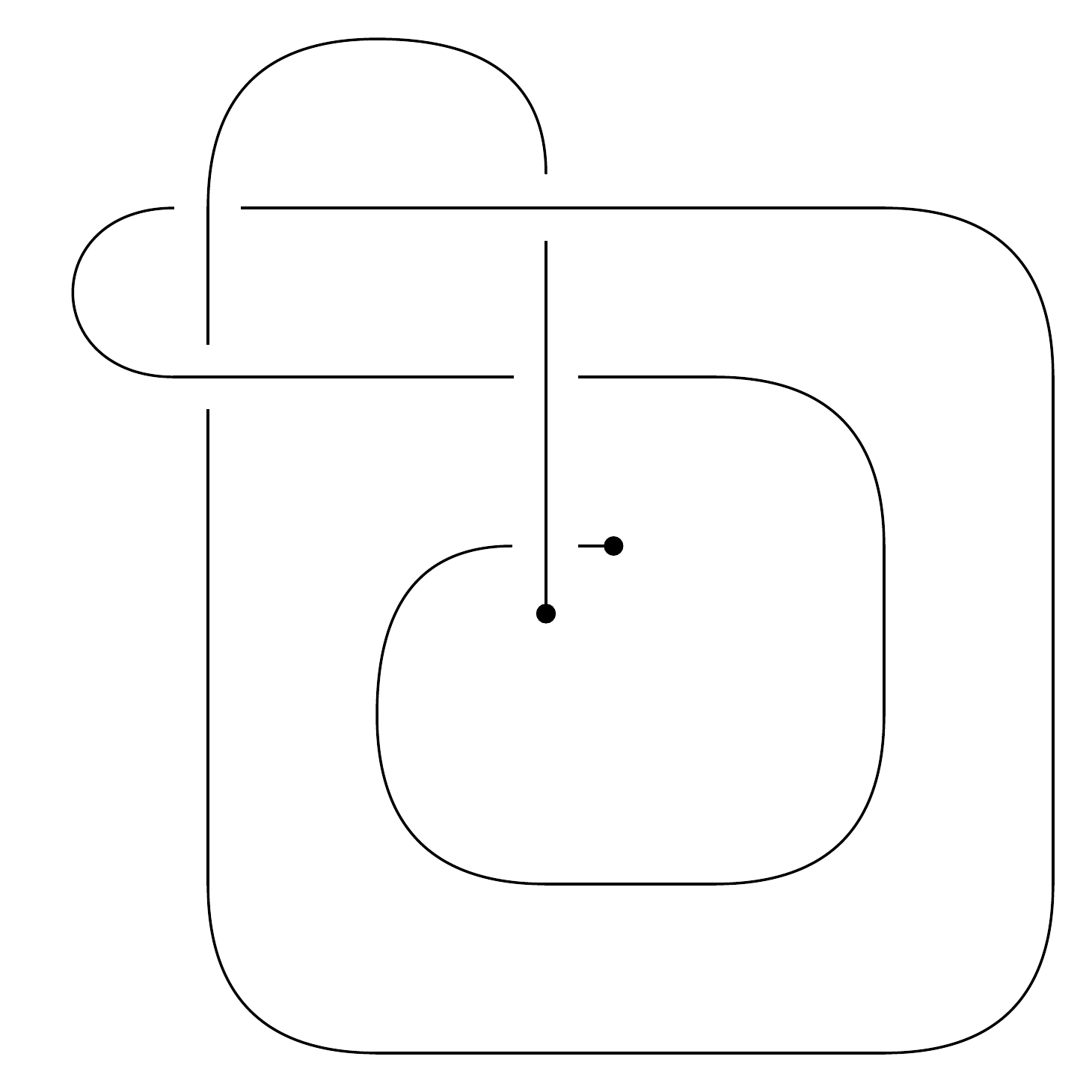}\\
\textcolor{black}{$5_{384}$}
\vspace{1cm}
\end{minipage}
\begin{minipage}[t]{.25\linewidth}
\centering
\includegraphics[width=0.9\textwidth,height=3.5cm,keepaspectratio]{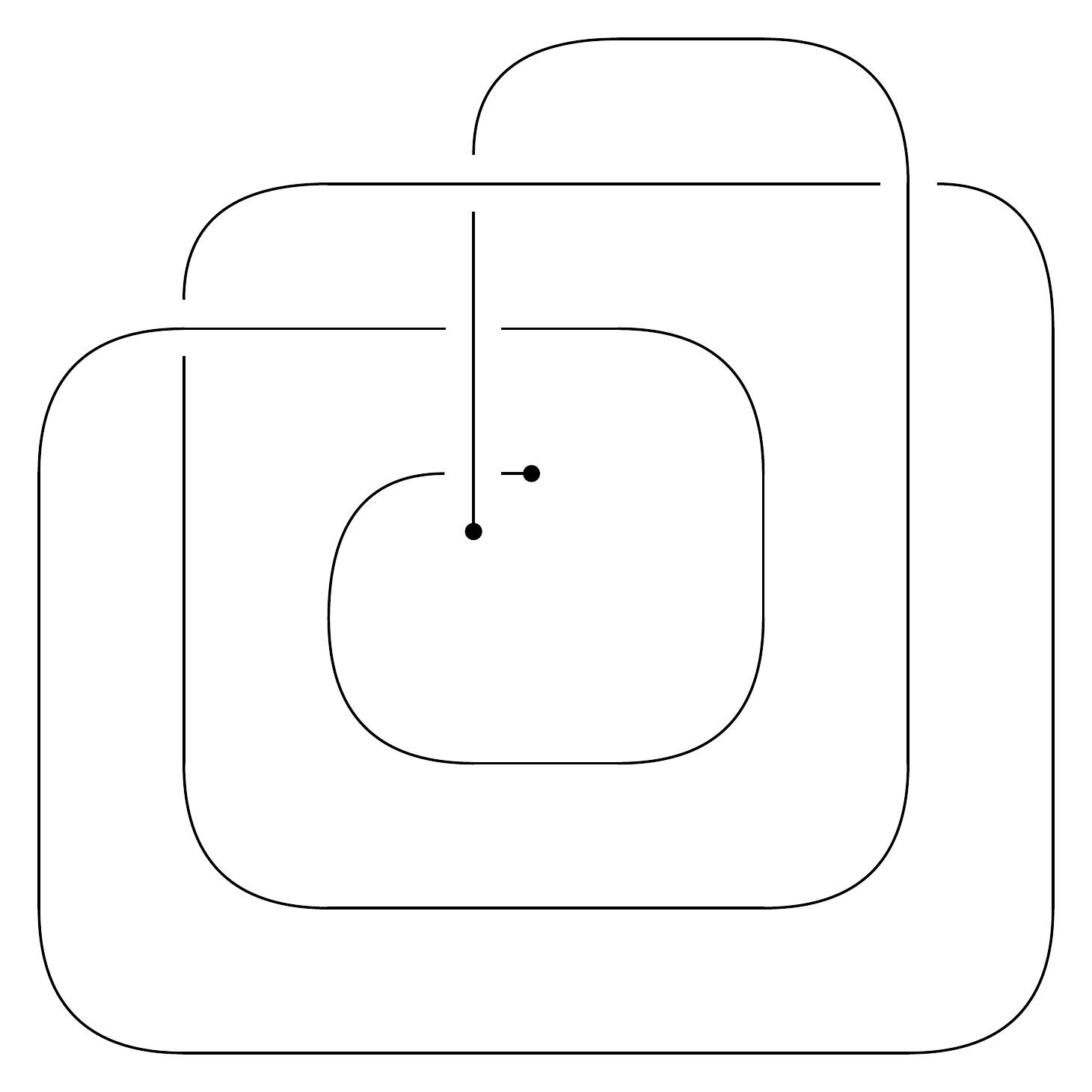}\\
\textcolor{black}{$5_{385}$}
\vspace{1cm}
\end{minipage}
\begin{minipage}[t]{.25\linewidth}
\centering
\includegraphics[width=0.9\textwidth,height=3.5cm,keepaspectratio]{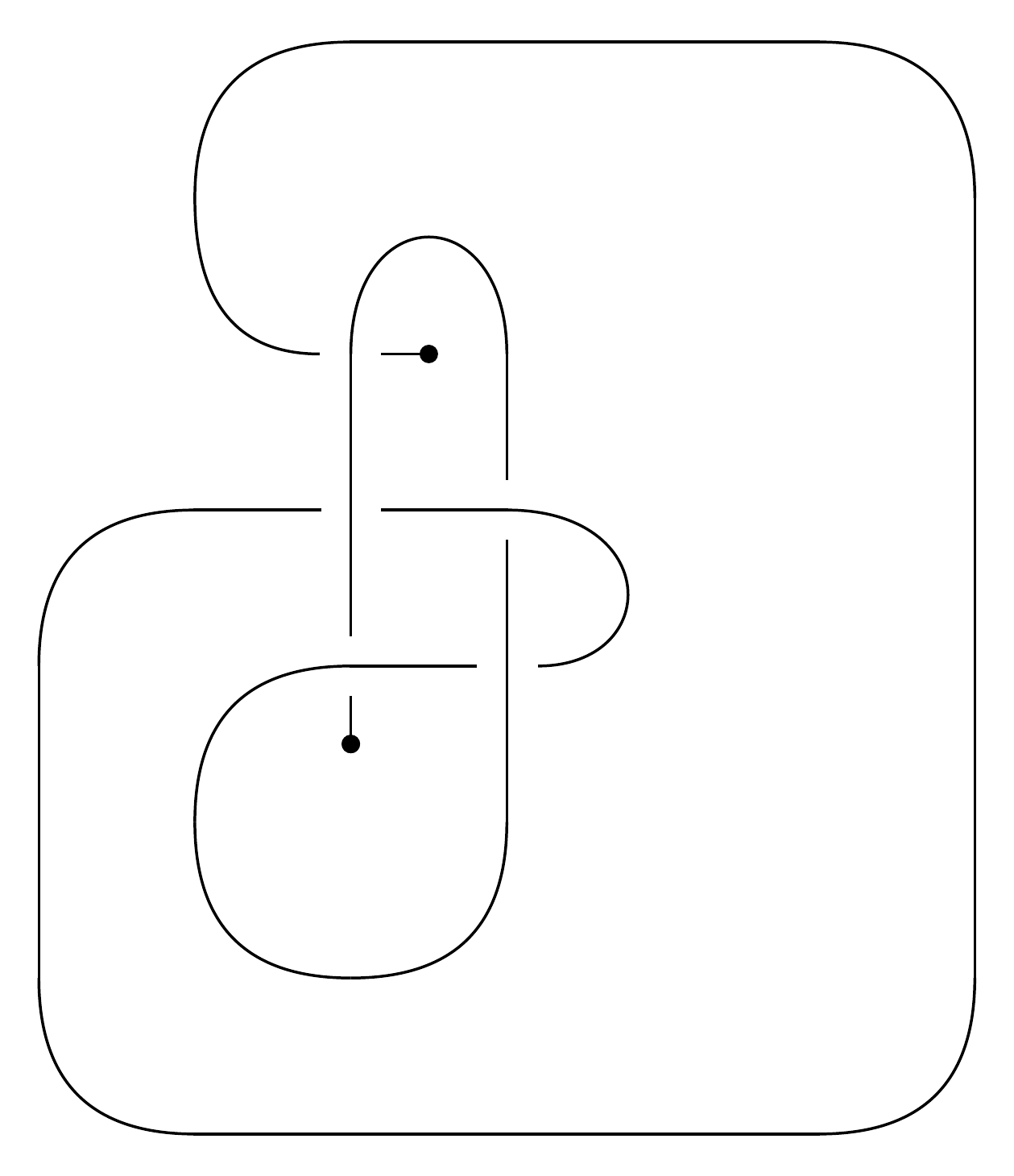}\\
\textcolor{black}{$5_{386}$}
\vspace{1cm}
\end{minipage}
\begin{minipage}[t]{.25\linewidth}
\centering
\includegraphics[width=0.9\textwidth,height=3.5cm,keepaspectratio]{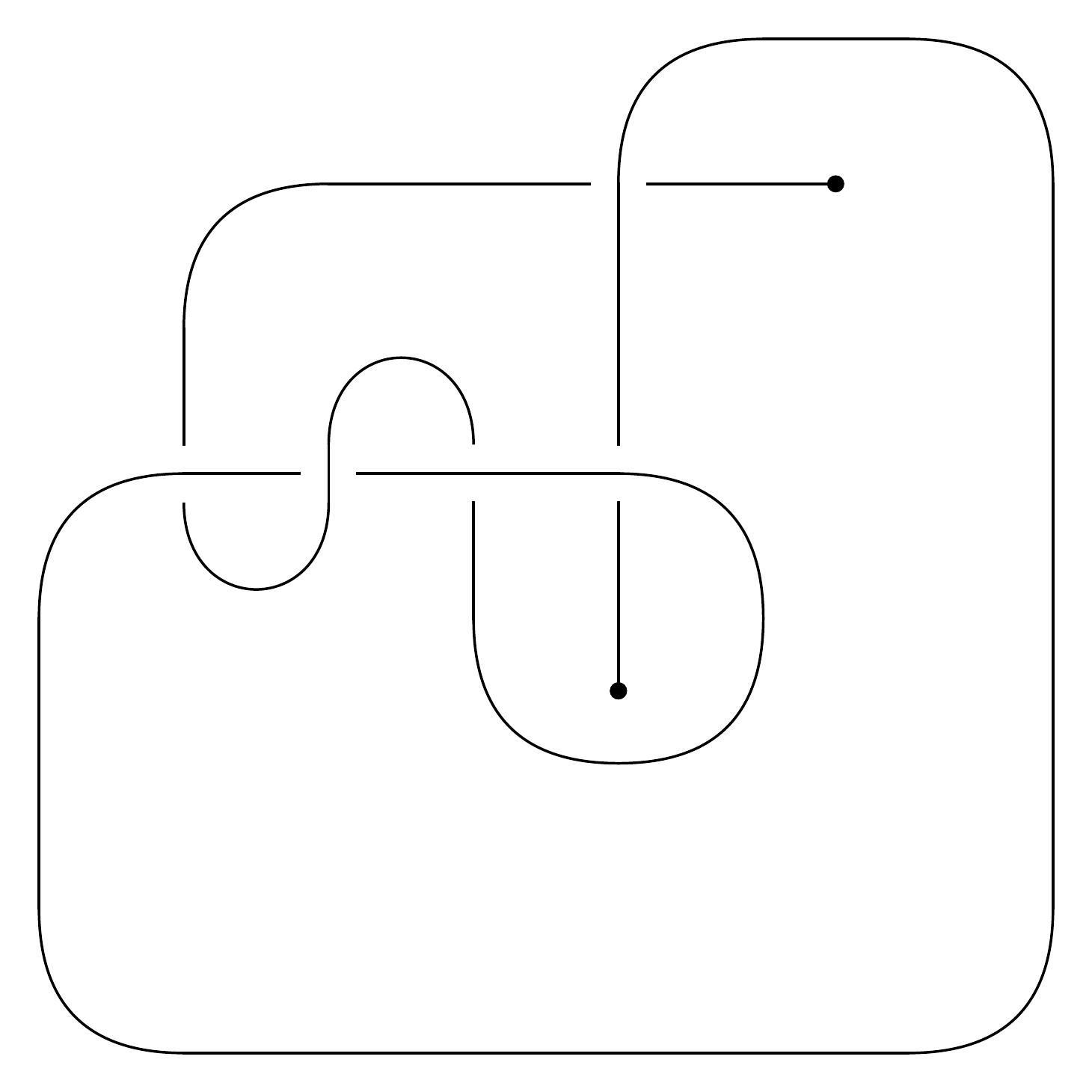}\\
\textcolor{black}{$5_{387}$}
\vspace{1cm}
\end{minipage}
\begin{minipage}[t]{.25\linewidth}
\centering
\includegraphics[width=0.9\textwidth,height=3.5cm,keepaspectratio]{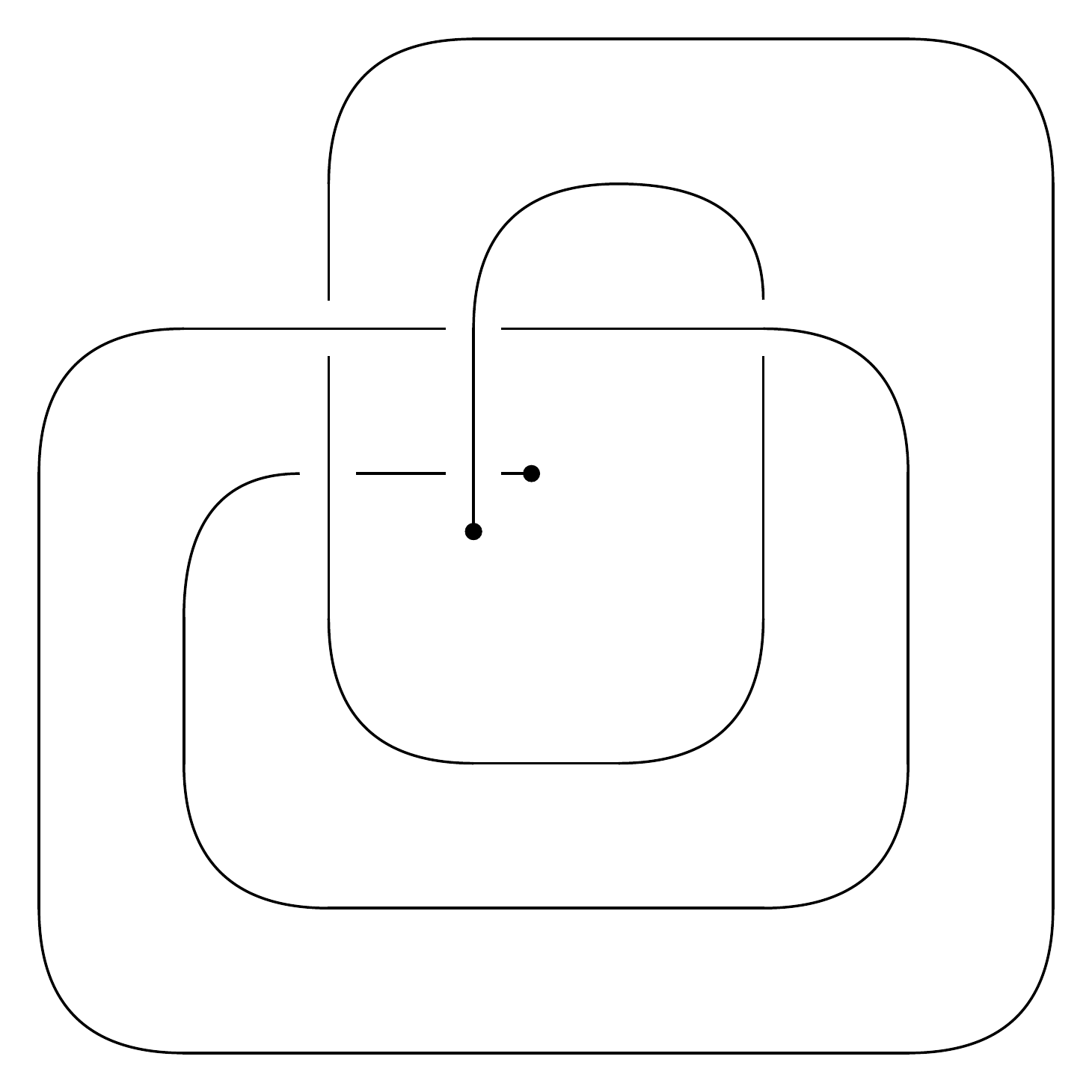}\\
\textcolor{black}{$5_{388}$}
\vspace{1cm}
\end{minipage}
\begin{minipage}[t]{.25\linewidth}
\centering
\includegraphics[width=0.9\textwidth,height=3.5cm,keepaspectratio]{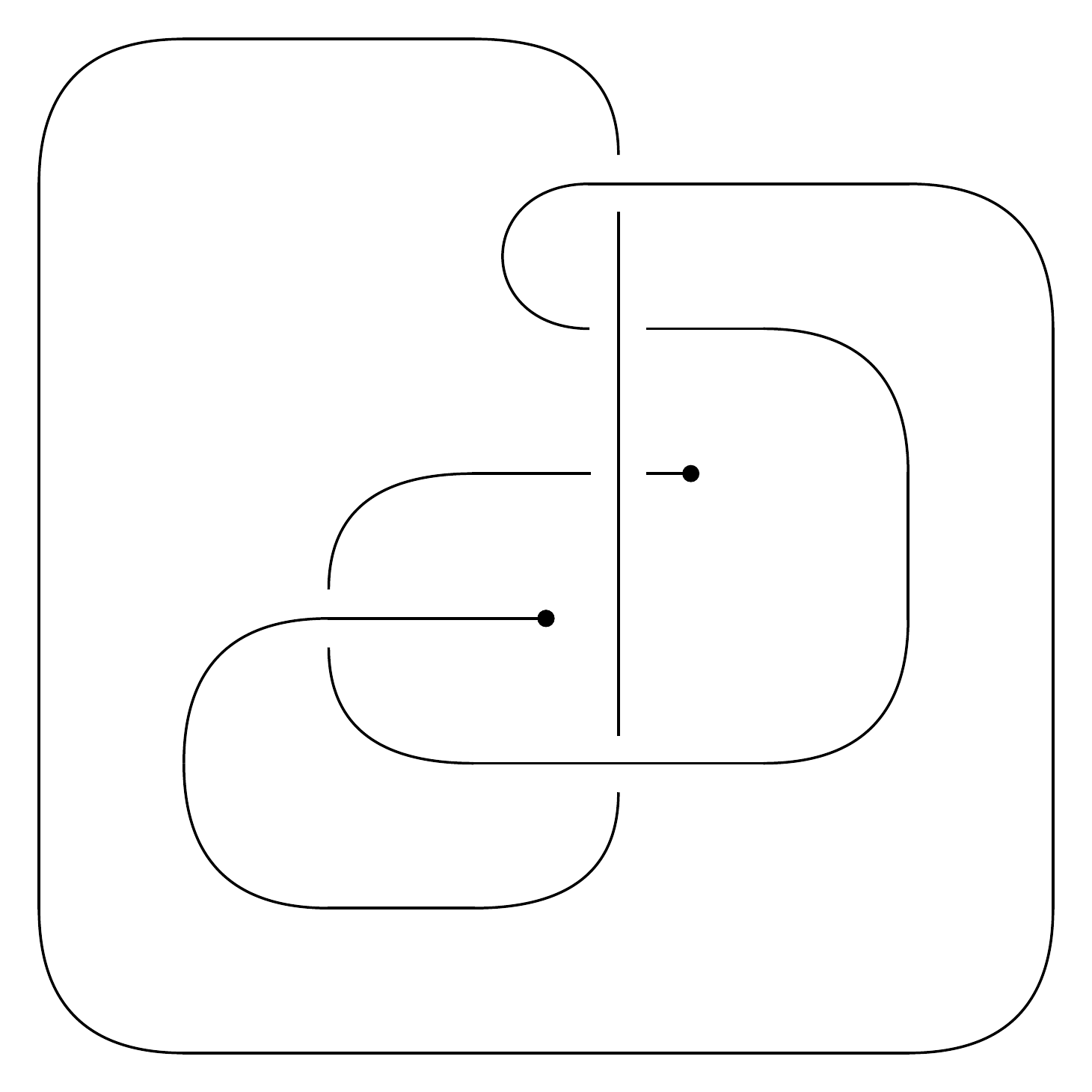}\\
\textcolor{black}{$5_{389}$}
\vspace{1cm}
\end{minipage}
\begin{minipage}[t]{.25\linewidth}
\centering
\includegraphics[width=0.9\textwidth,height=3.5cm,keepaspectratio]{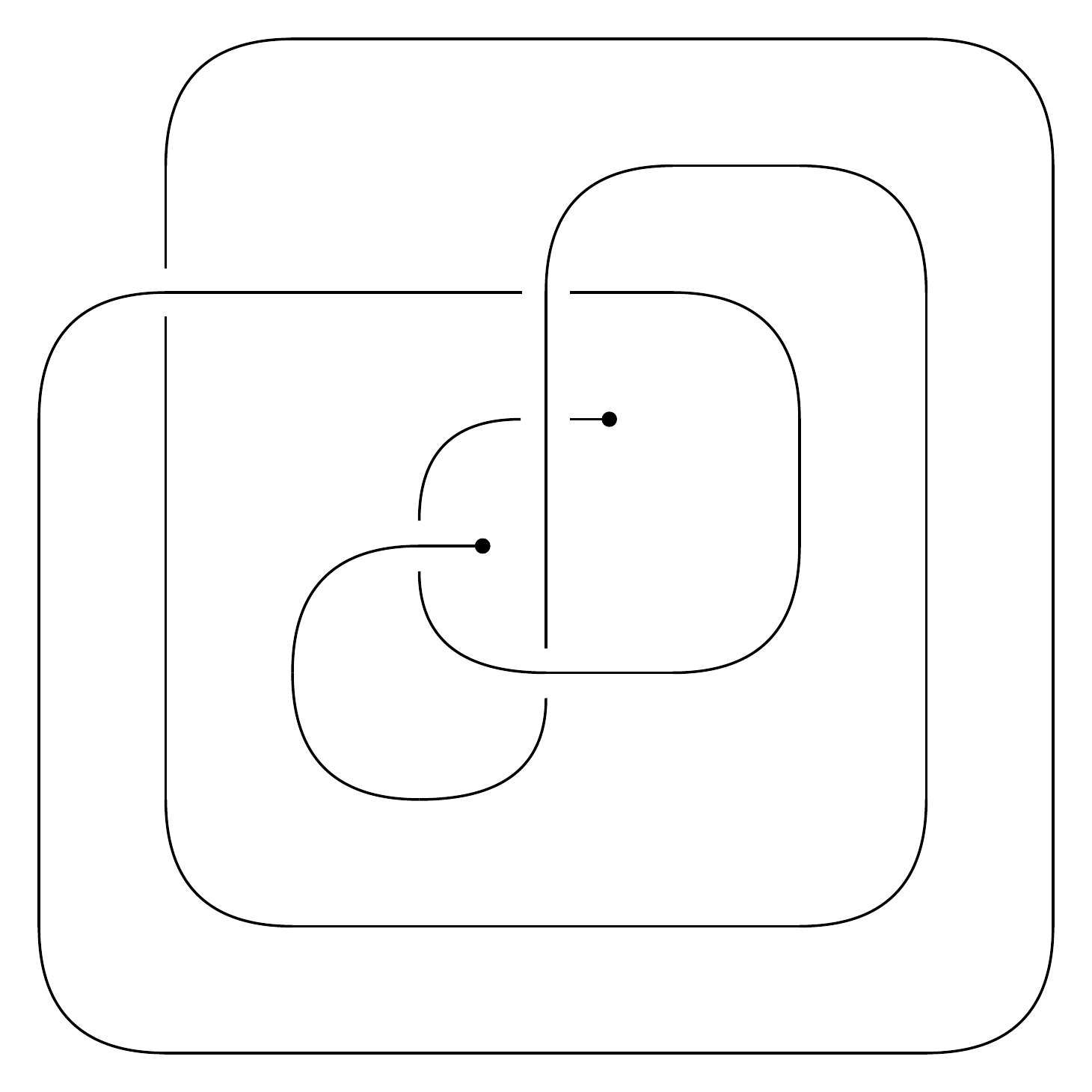}\\
\textcolor{black}{$5_{390}$}
\vspace{1cm}
\end{minipage}
\begin{minipage}[t]{.25\linewidth}
\centering
\includegraphics[width=0.9\textwidth,height=3.5cm,keepaspectratio]{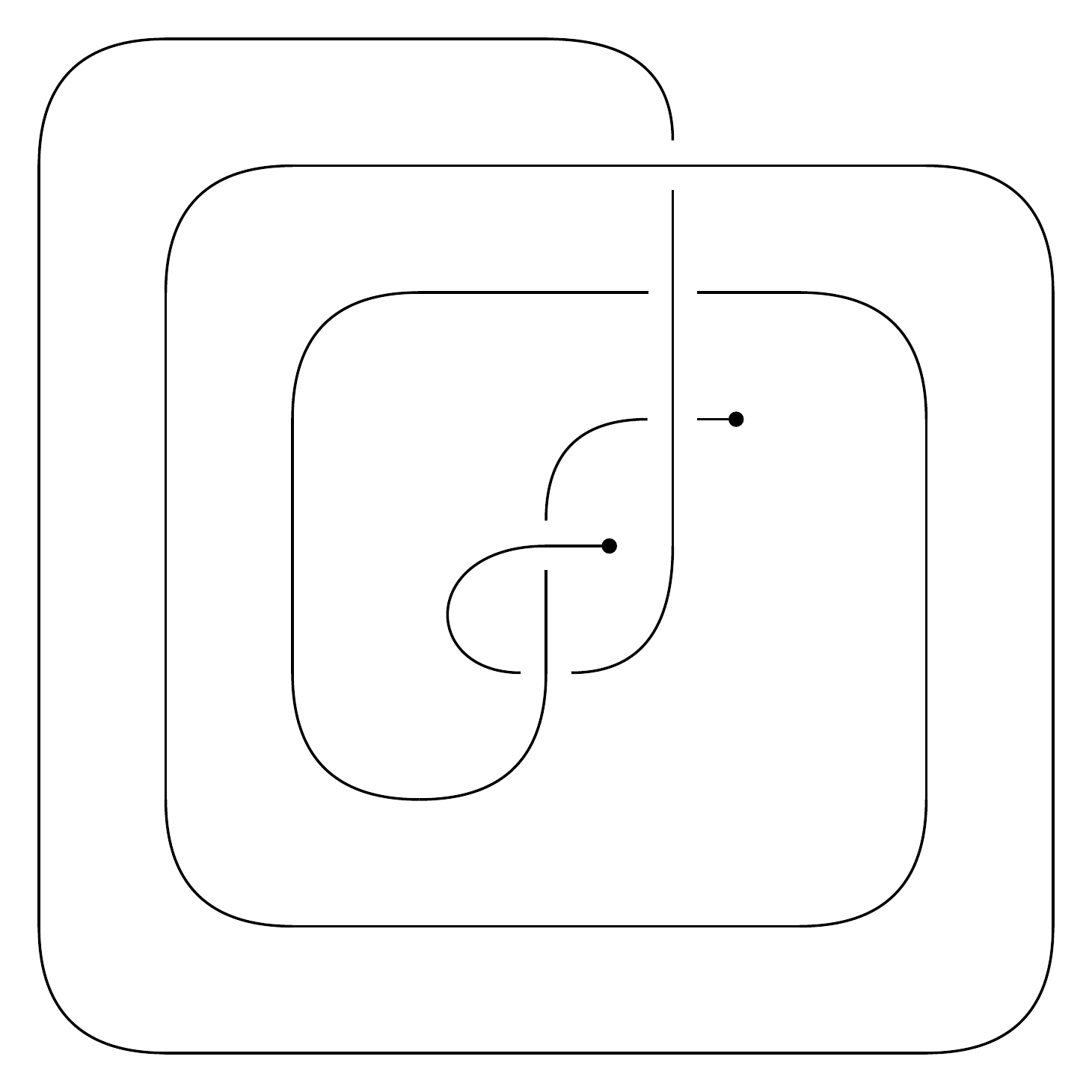}\\
\textcolor{black}{$5_{391}$}
\vspace{1cm}
\end{minipage}
\begin{minipage}[t]{.25\linewidth}
\centering
\includegraphics[width=0.9\textwidth,height=3.5cm,keepaspectratio]{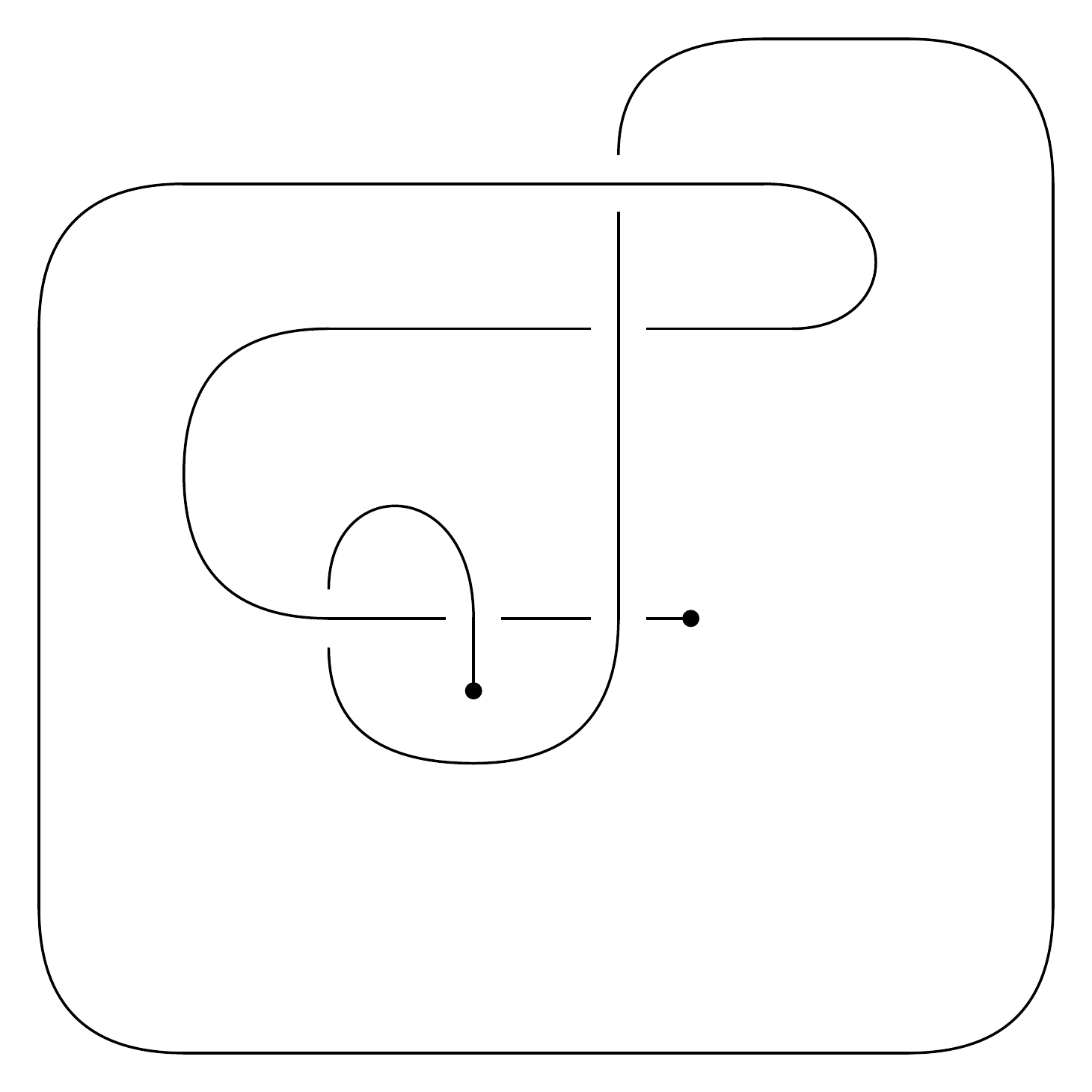}\\
\textcolor{black}{$5_{392}$}
\vspace{1cm}
\end{minipage}
\begin{minipage}[t]{.25\linewidth}
\centering
\includegraphics[width=0.9\textwidth,height=3.5cm,keepaspectratio]{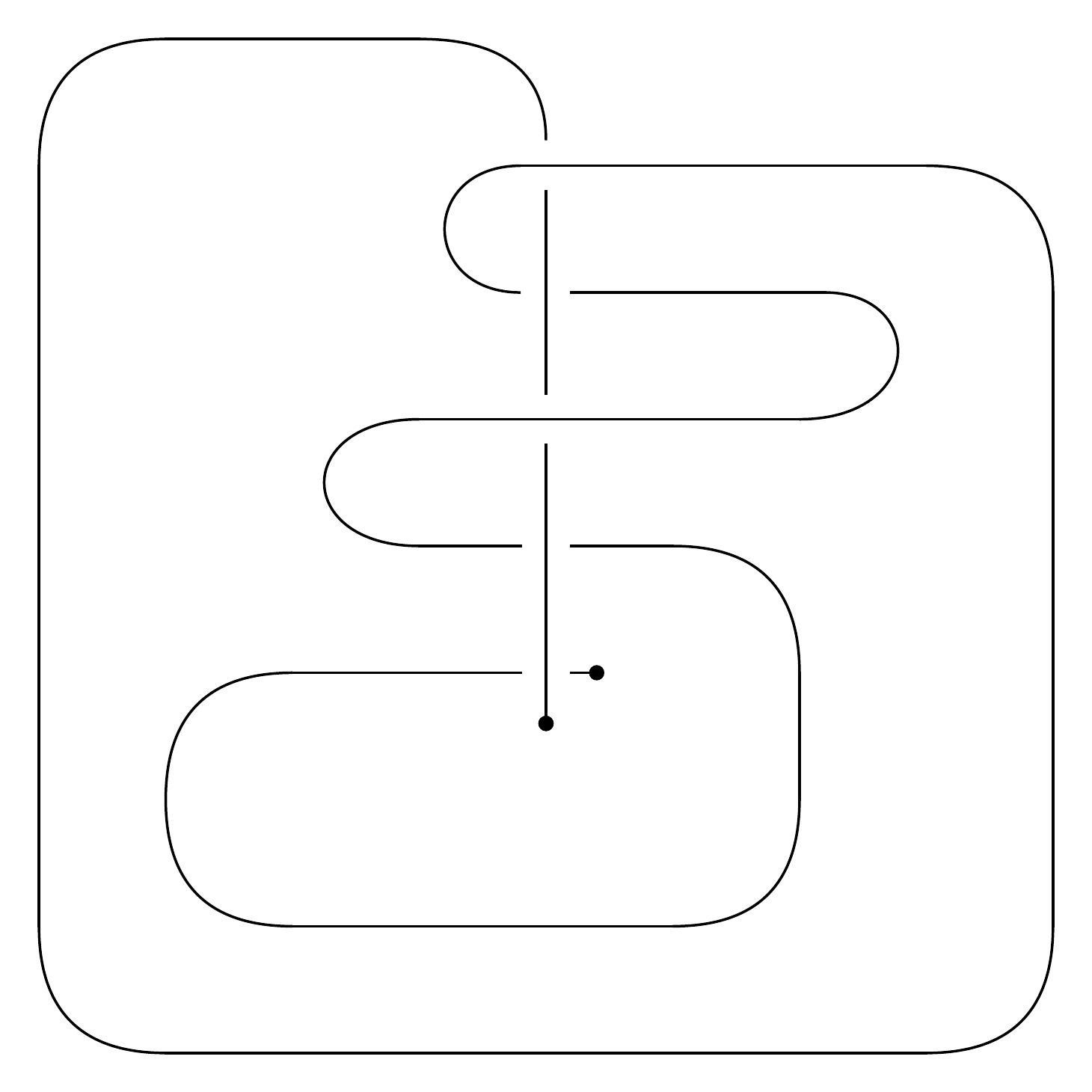}\\
\textcolor{black}{$5_{393}$}
\vspace{1cm}
\end{minipage}
\begin{minipage}[t]{.25\linewidth}
\centering
\includegraphics[width=0.9\textwidth,height=3.5cm,keepaspectratio]{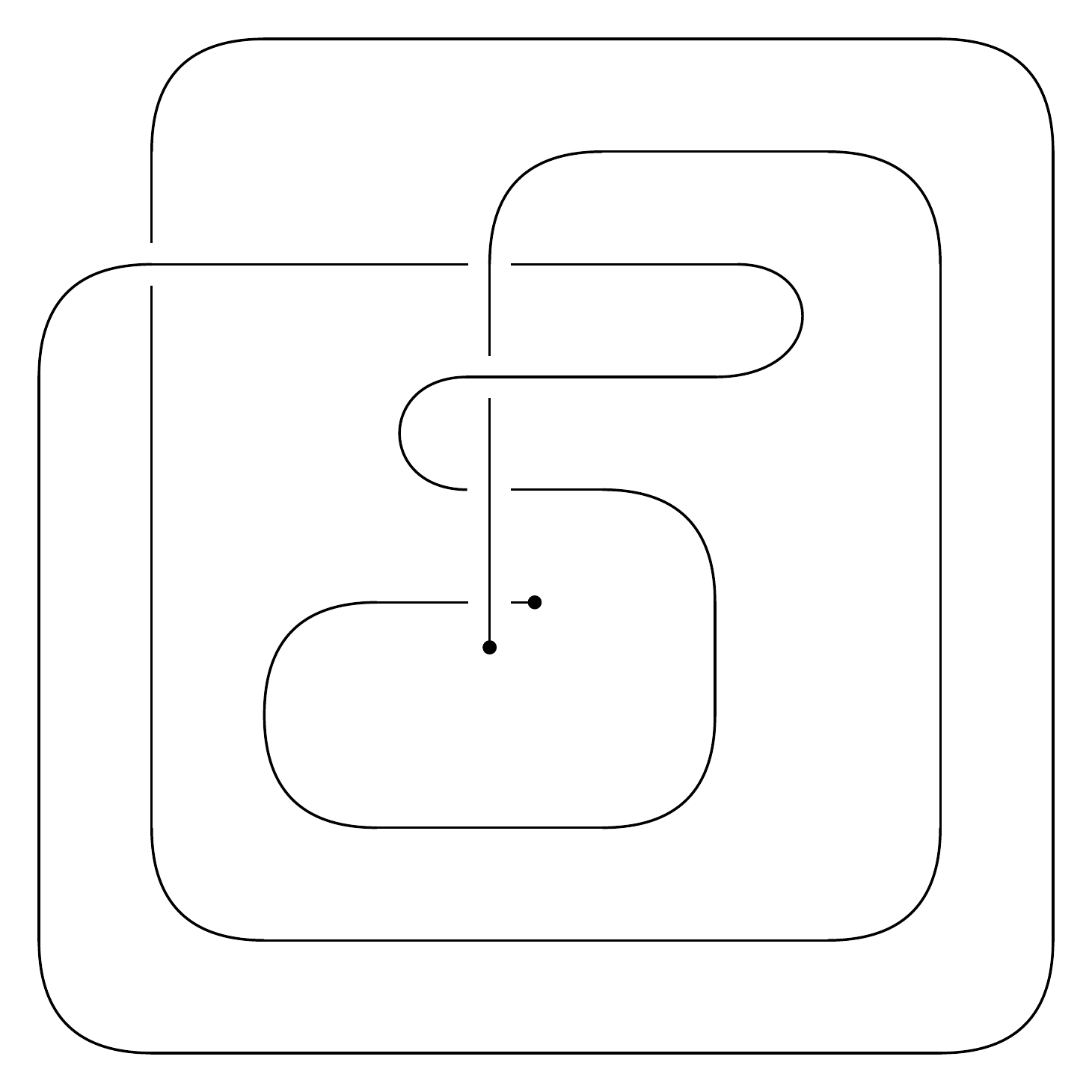}\\
\textcolor{black}{$5_{394}$}
\vspace{1cm}
\end{minipage}
\begin{minipage}[t]{.25\linewidth}
\centering
\includegraphics[width=0.9\textwidth,height=3.5cm,keepaspectratio]{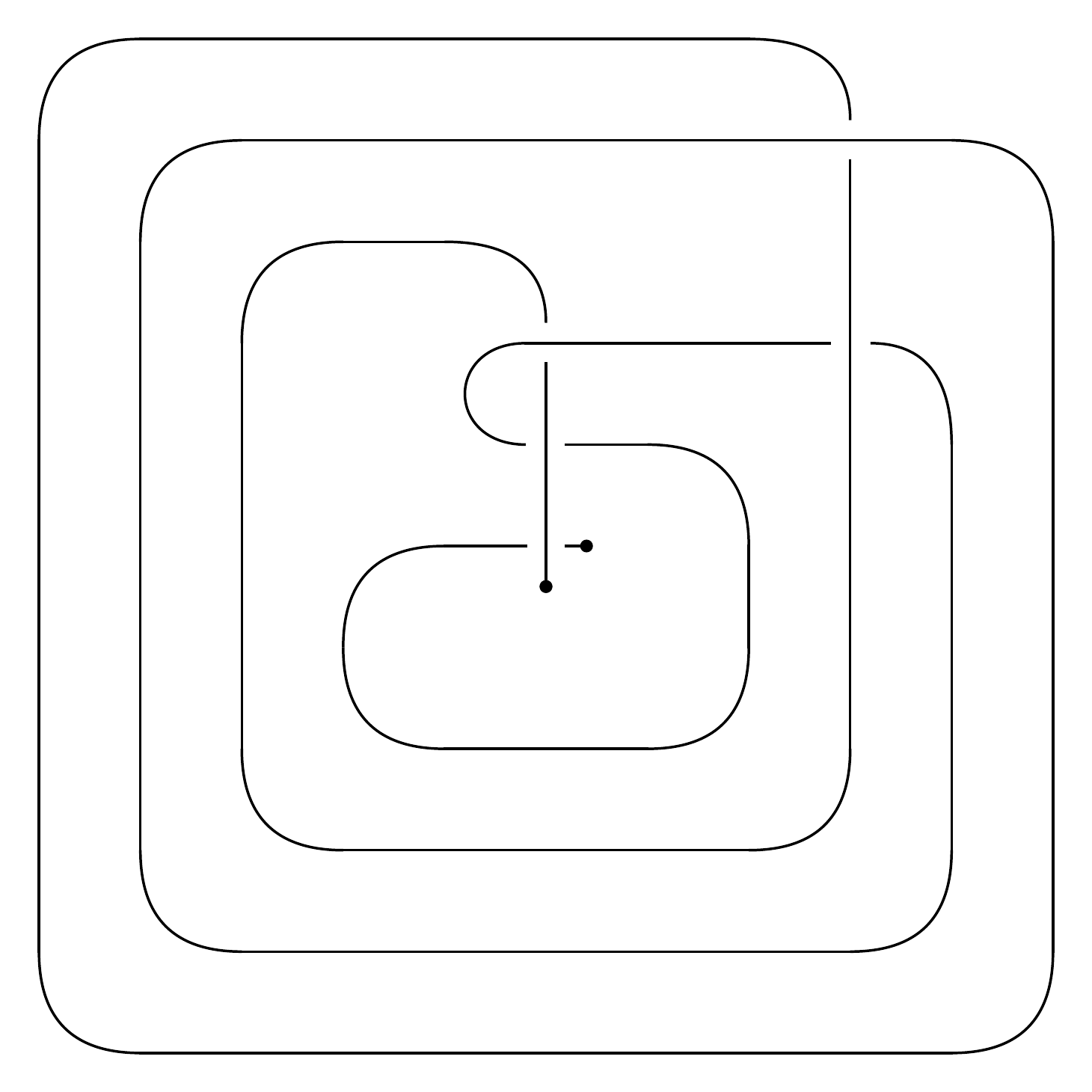}\\
\textcolor{black}{$5_{395}$}
\vspace{1cm}
\end{minipage}
\begin{minipage}[t]{.25\linewidth}
\centering
\includegraphics[width=0.9\textwidth,height=3.5cm,keepaspectratio]{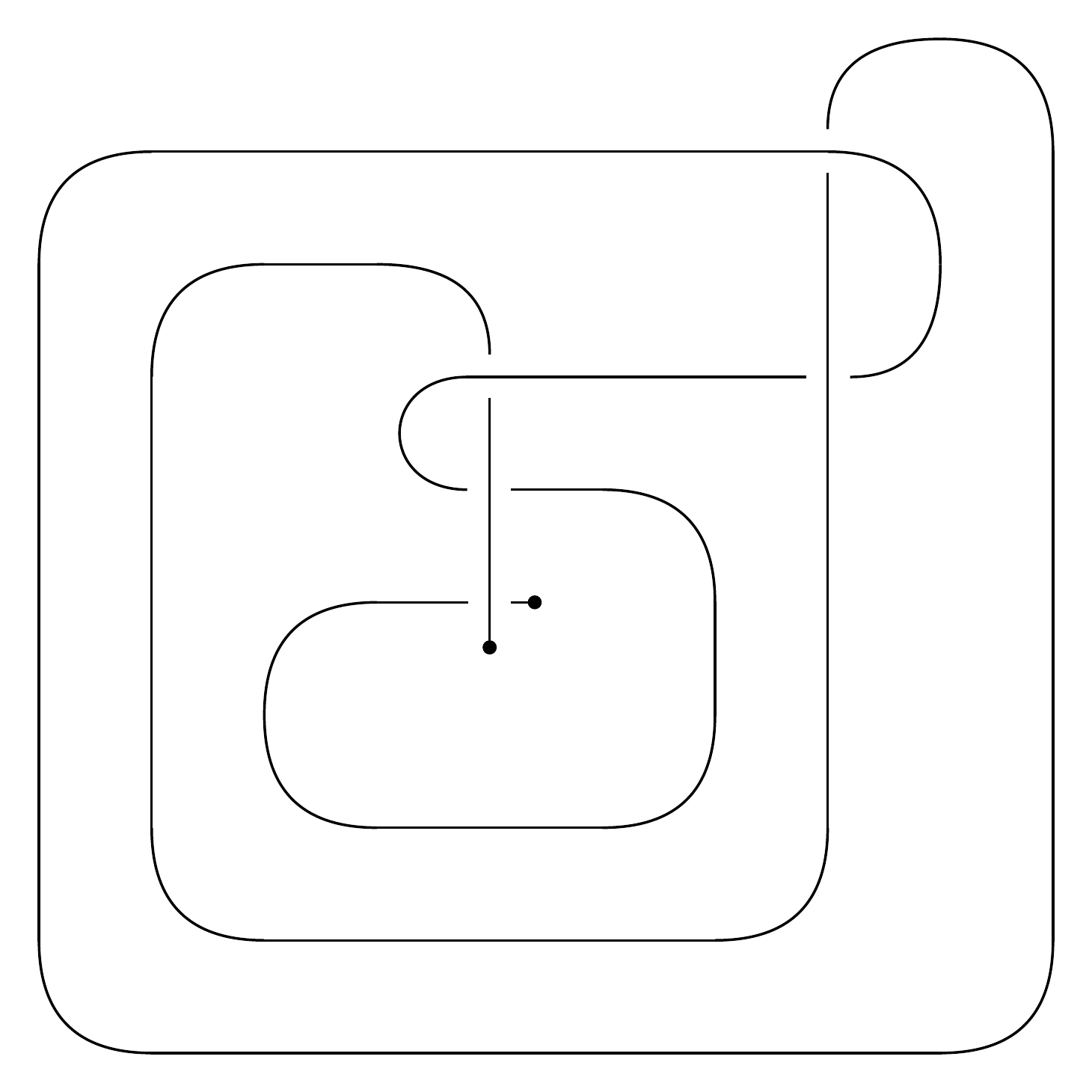}\\
\textcolor{black}{$5_{396}$}
\vspace{1cm}
\end{minipage}
\begin{minipage}[t]{.25\linewidth}
\centering
\includegraphics[width=0.9\textwidth,height=3.5cm,keepaspectratio]{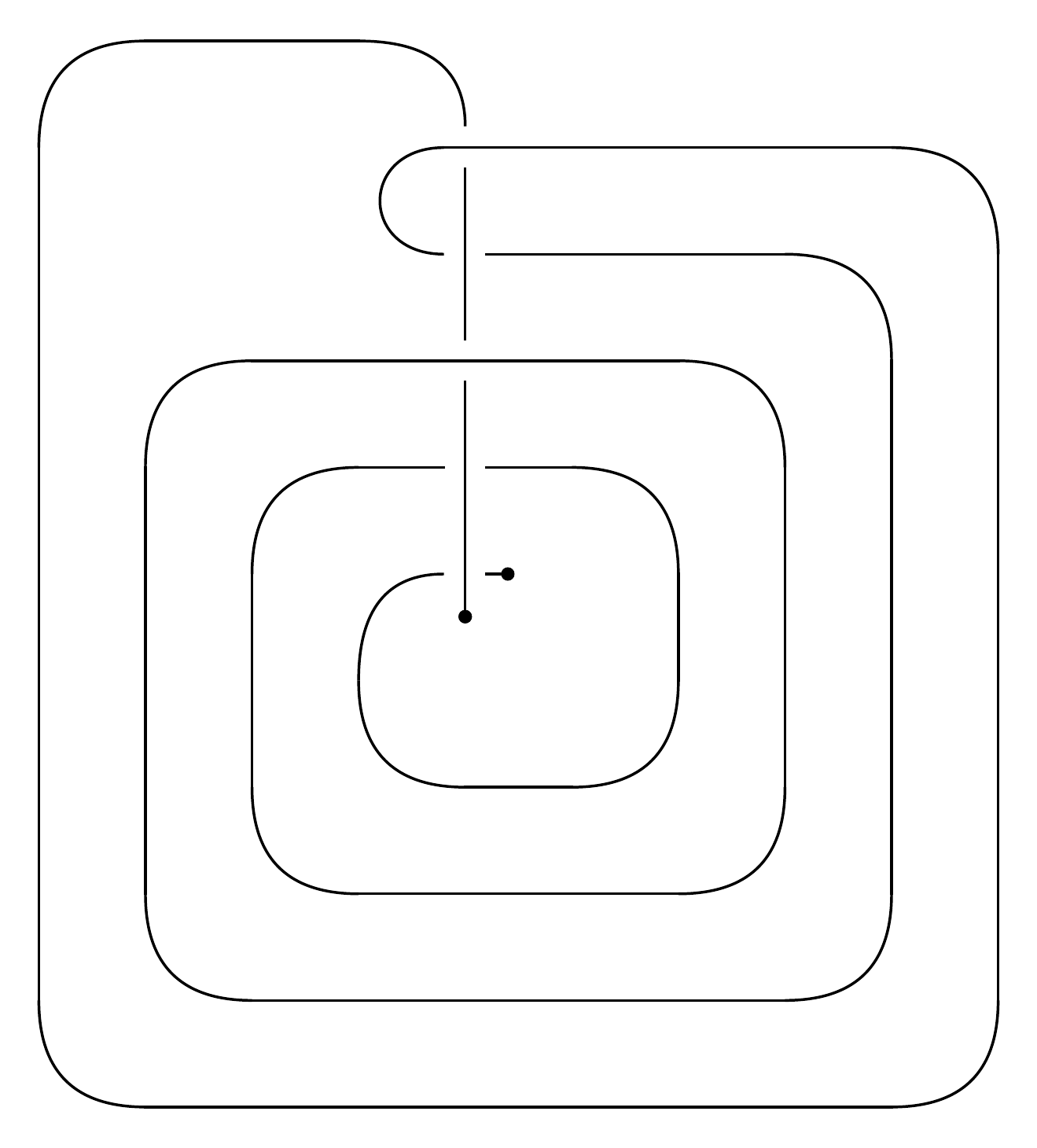}\\
\textcolor{black}{$5_{397}$}
\vspace{1cm}
\end{minipage}
\begin{minipage}[t]{.25\linewidth}
\centering
\includegraphics[width=0.9\textwidth,height=3.5cm,keepaspectratio]{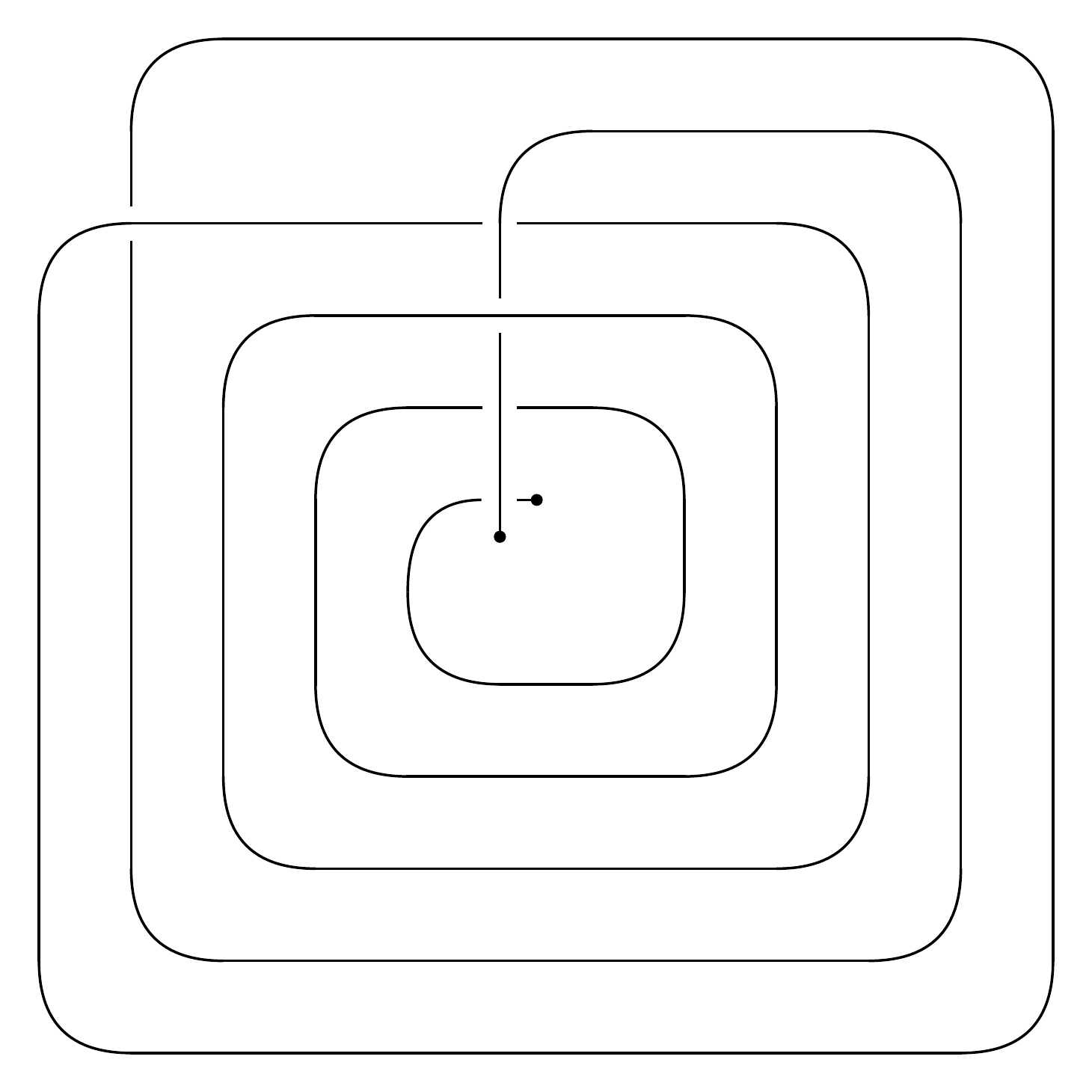}\\
\textcolor{black}{$5_{398}$}
\vspace{1cm}
\end{minipage}
\begin{minipage}[t]{.25\linewidth}
\centering
\includegraphics[width=0.9\textwidth,height=3.5cm,keepaspectratio]{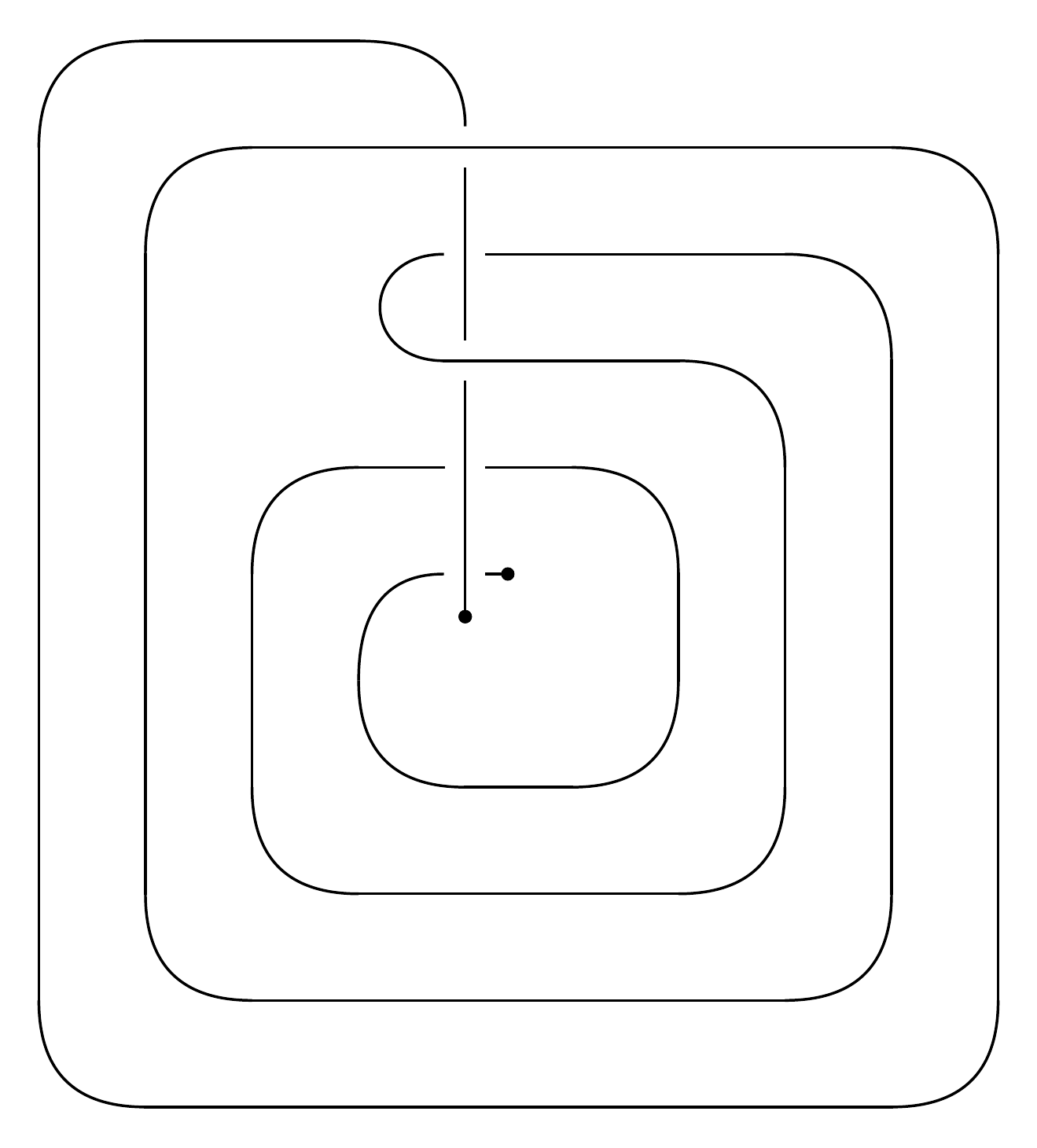}\\
\textcolor{black}{$5_{399}$}
\vspace{1cm}
\end{minipage}
\begin{minipage}[t]{.25\linewidth}
\centering
\includegraphics[width=0.9\textwidth,height=3.5cm,keepaspectratio]{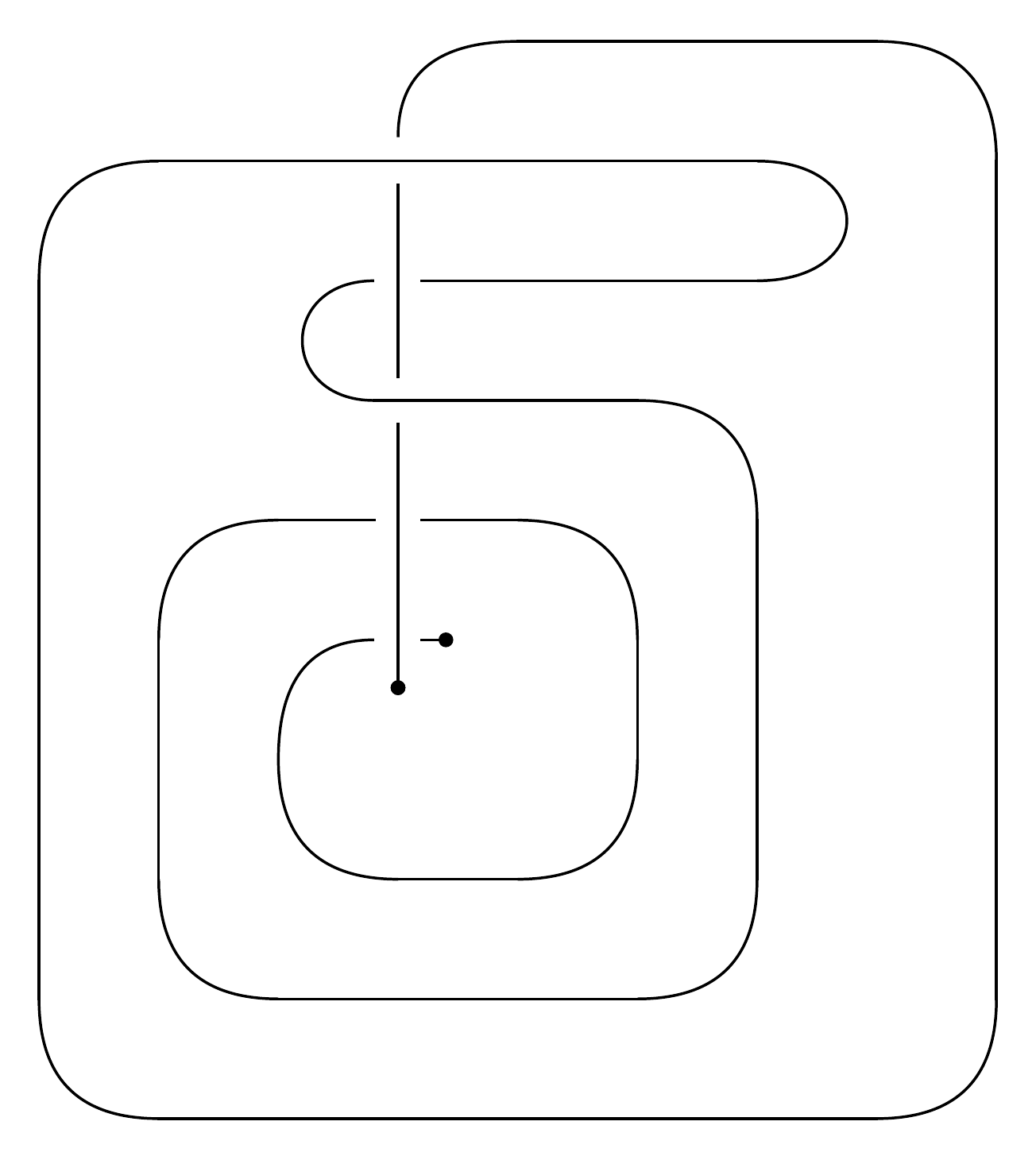}\\
\textcolor{black}{$5_{400}$}
\vspace{1cm}
\end{minipage}
\begin{minipage}[t]{.25\linewidth}
\centering
\includegraphics[width=0.9\textwidth,height=3.5cm,keepaspectratio]{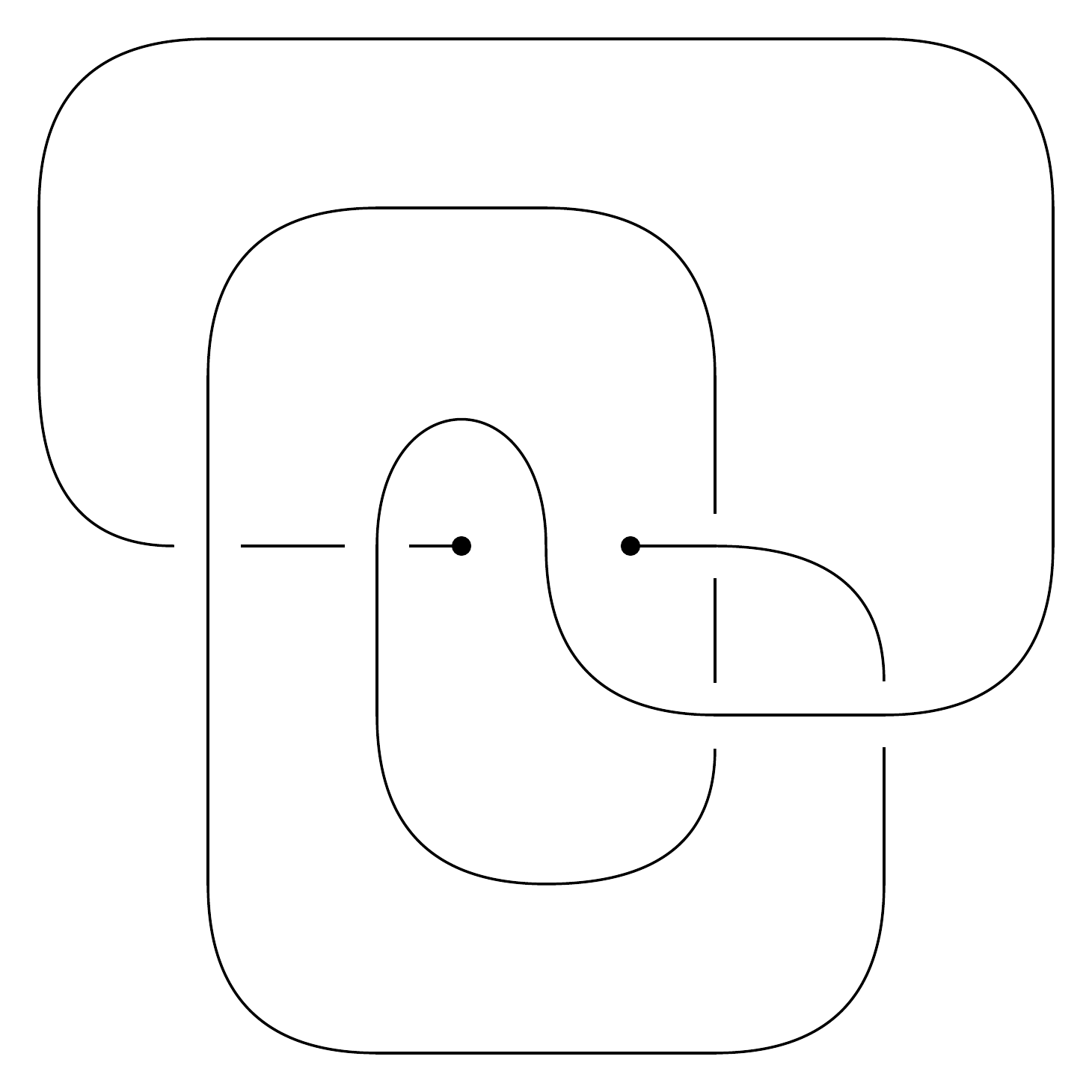}\\
\textcolor{black}{$5_{401}$}
\vspace{1cm}
\end{minipage}
\begin{minipage}[t]{.25\linewidth}
\centering
\includegraphics[width=0.9\textwidth,height=3.5cm,keepaspectratio]{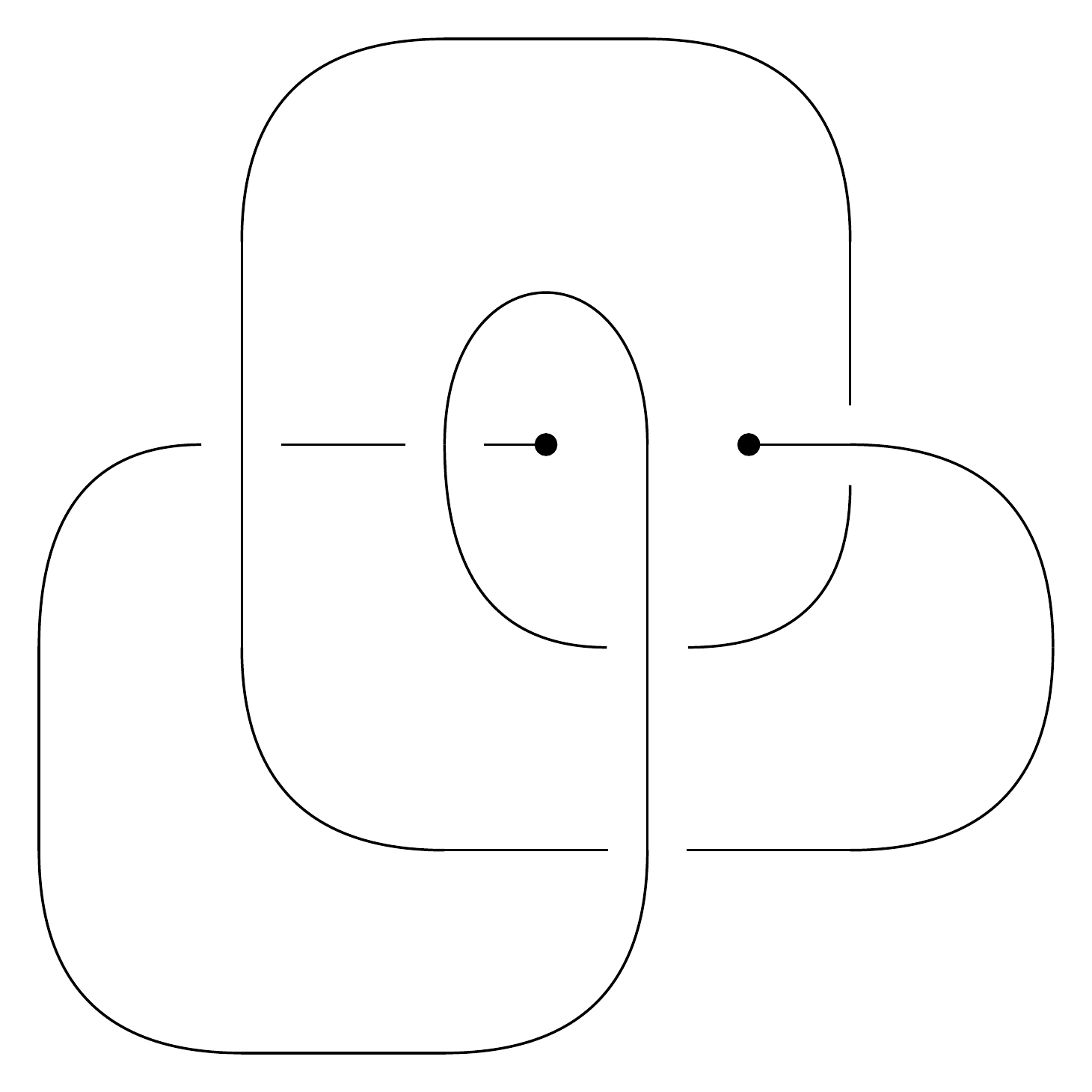}\\
\textcolor{black}{$5_{402}$}
\vspace{1cm}
\end{minipage}
\begin{minipage}[t]{.25\linewidth}
\centering
\includegraphics[width=0.9\textwidth,height=3.5cm,keepaspectratio]{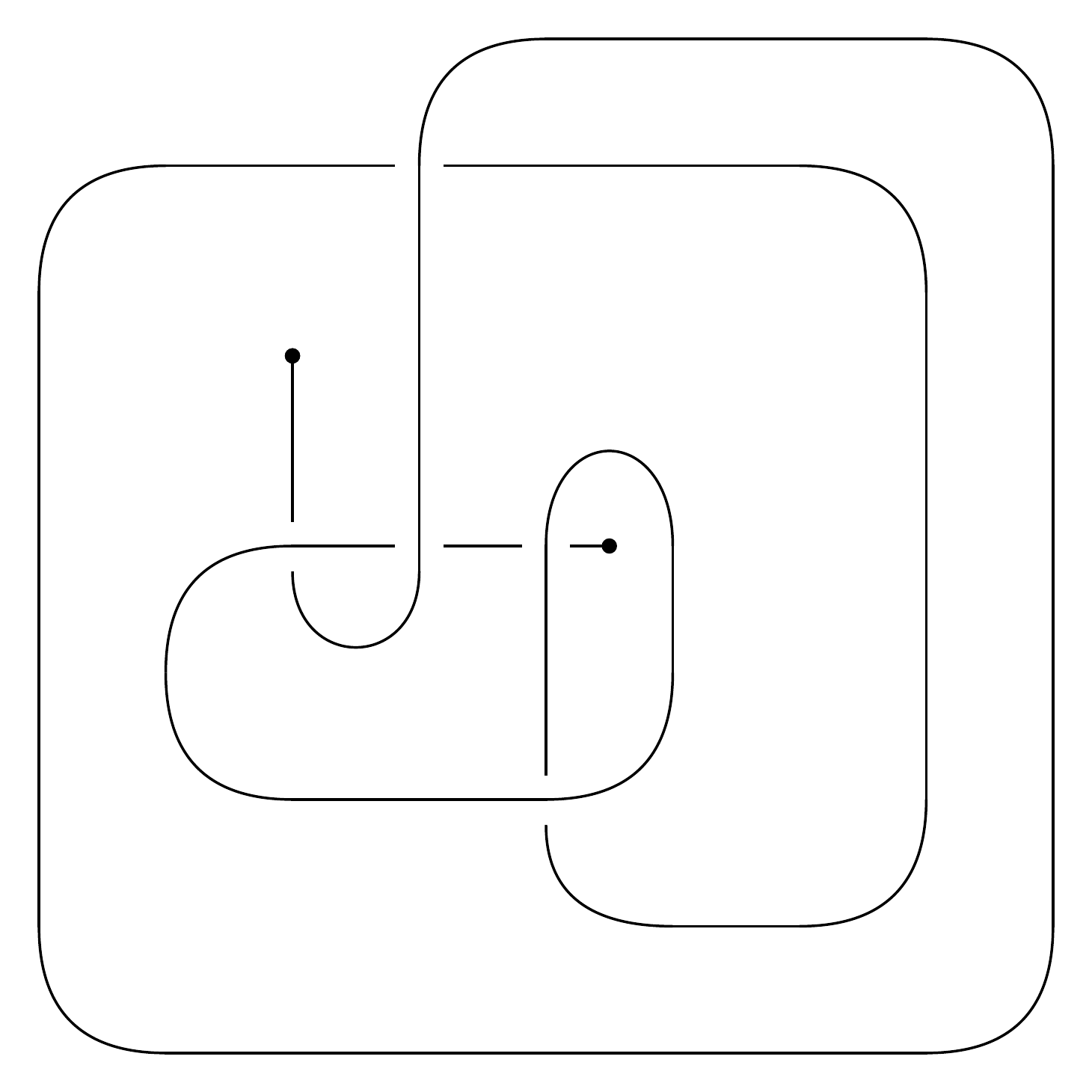}\\
\textcolor{black}{$5_{403}$}
\vspace{1cm}
\end{minipage}
\begin{minipage}[t]{.25\linewidth}
\centering
\includegraphics[width=0.9\textwidth,height=3.5cm,keepaspectratio]{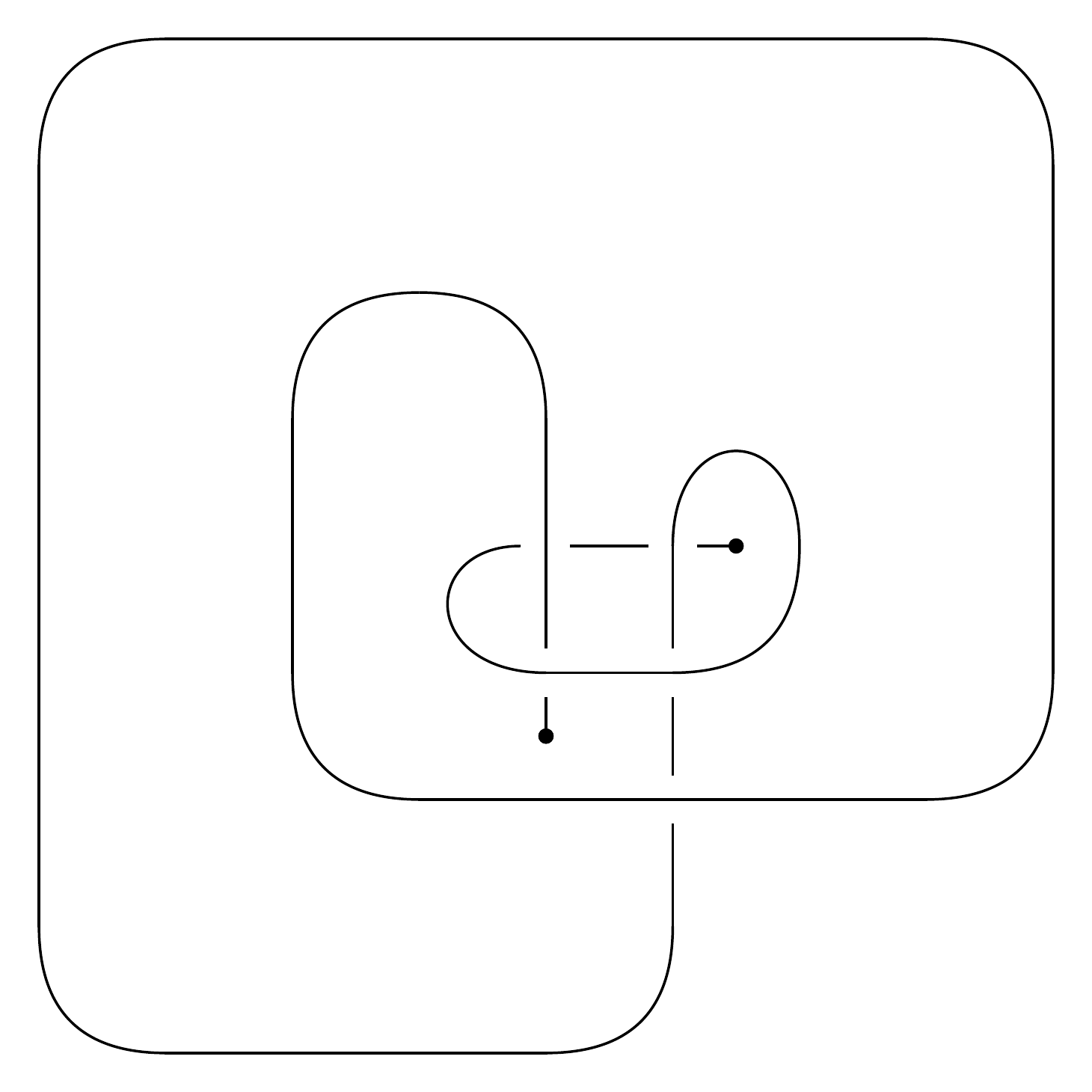}\\
\textcolor{black}{$5_{404}$}
\vspace{1cm}
\end{minipage}
\begin{minipage}[t]{.25\linewidth}
\centering
\includegraphics[width=0.9\textwidth,height=3.5cm,keepaspectratio]{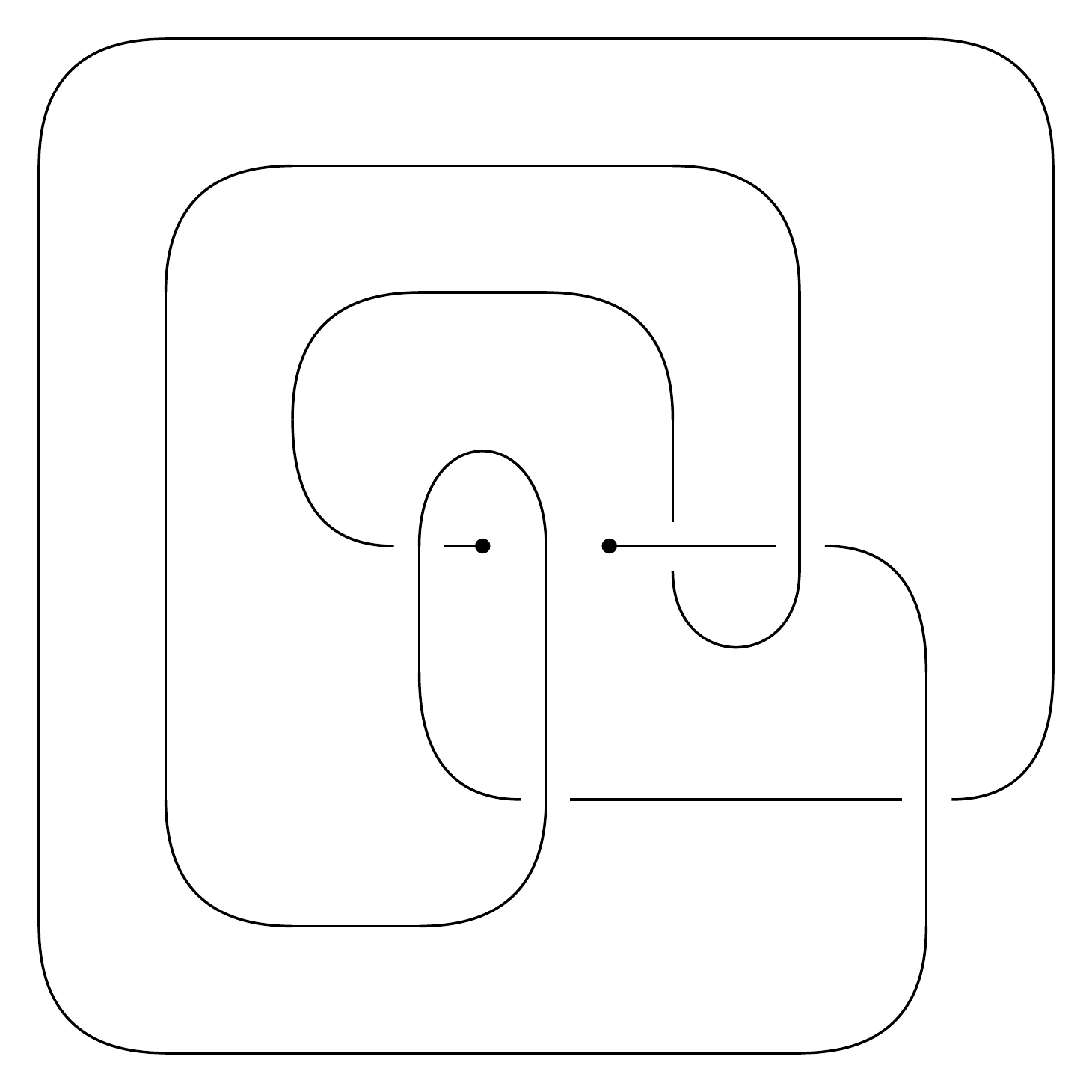}\\
\textcolor{black}{$5_{405}$}
\vspace{1cm}
\end{minipage}
\begin{minipage}[t]{.25\linewidth}
\centering
\includegraphics[width=0.9\textwidth,height=3.5cm,keepaspectratio]{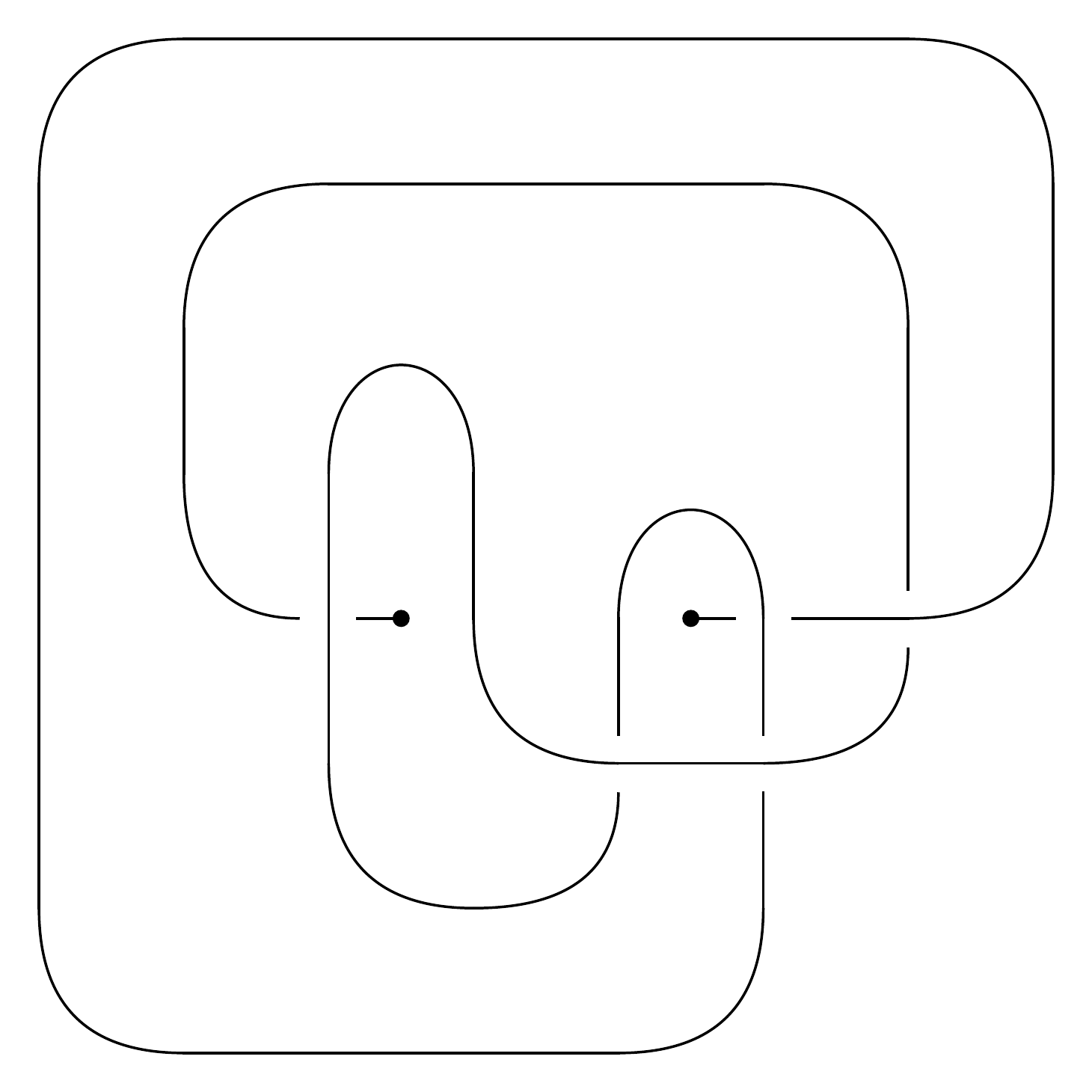}\\
\textcolor{black}{$5_{406}$}
\vspace{1cm}
\end{minipage}
\begin{minipage}[t]{.25\linewidth}
\centering
\includegraphics[width=0.9\textwidth,height=3.5cm,keepaspectratio]{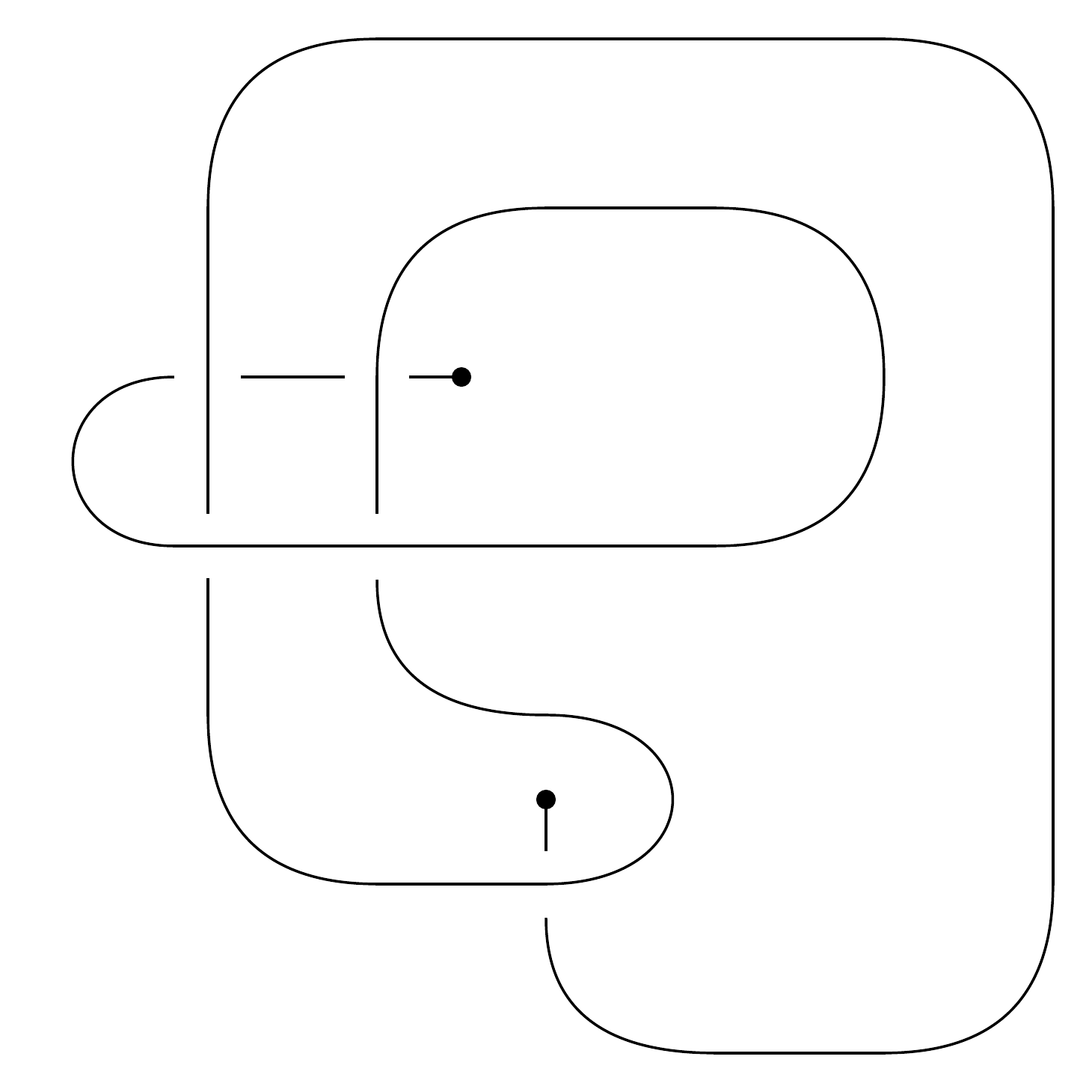}\\
\textcolor{black}{$5_{407}$}
\vspace{1cm}
\end{minipage}
\begin{minipage}[t]{.25\linewidth}
\centering
\includegraphics[width=0.9\textwidth,height=3.5cm,keepaspectratio]{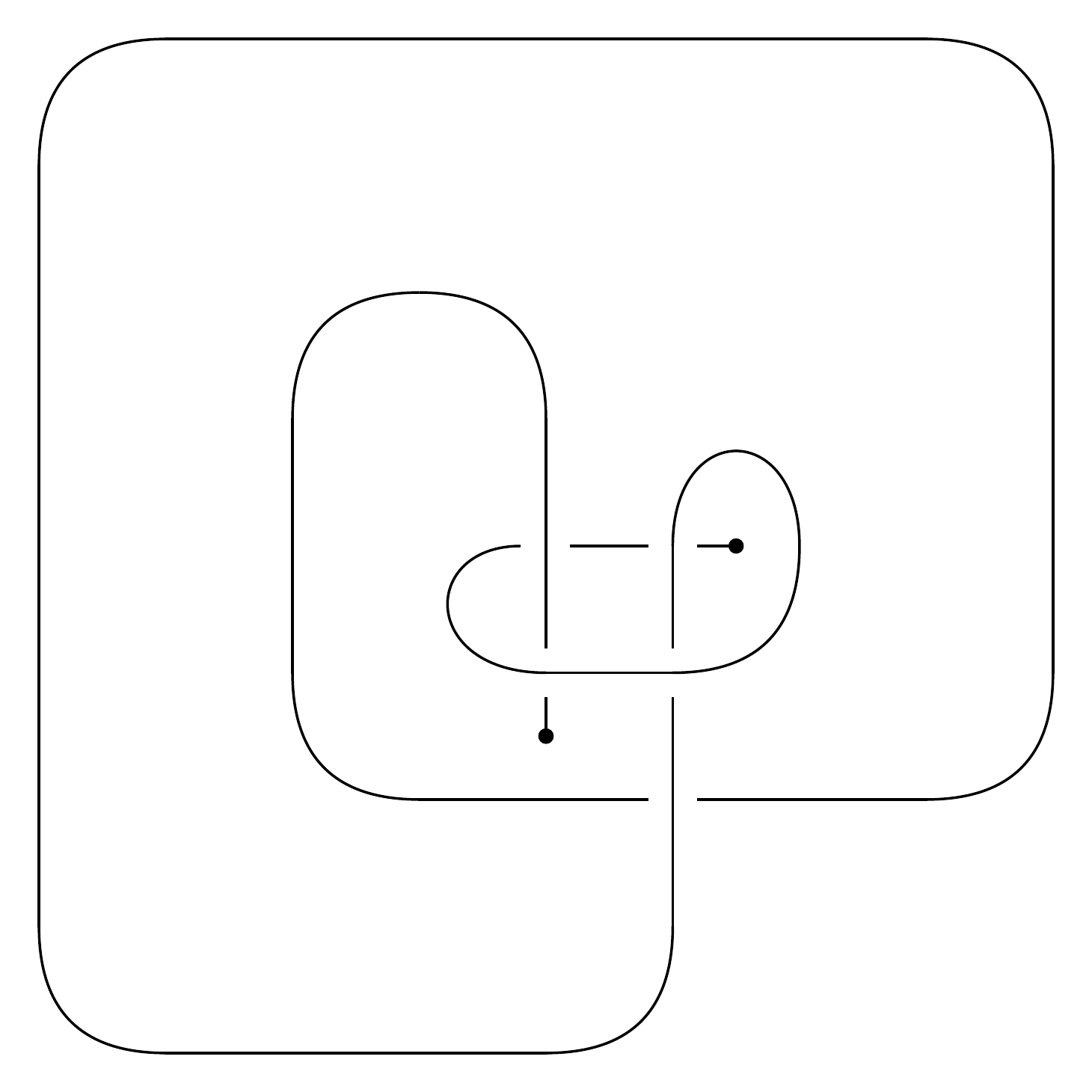}\\
\textcolor{black}{$5_{408}$}
\vspace{1cm}
\end{minipage}
\begin{minipage}[t]{.25\linewidth}
\centering
\includegraphics[width=0.9\textwidth,height=3.5cm,keepaspectratio]{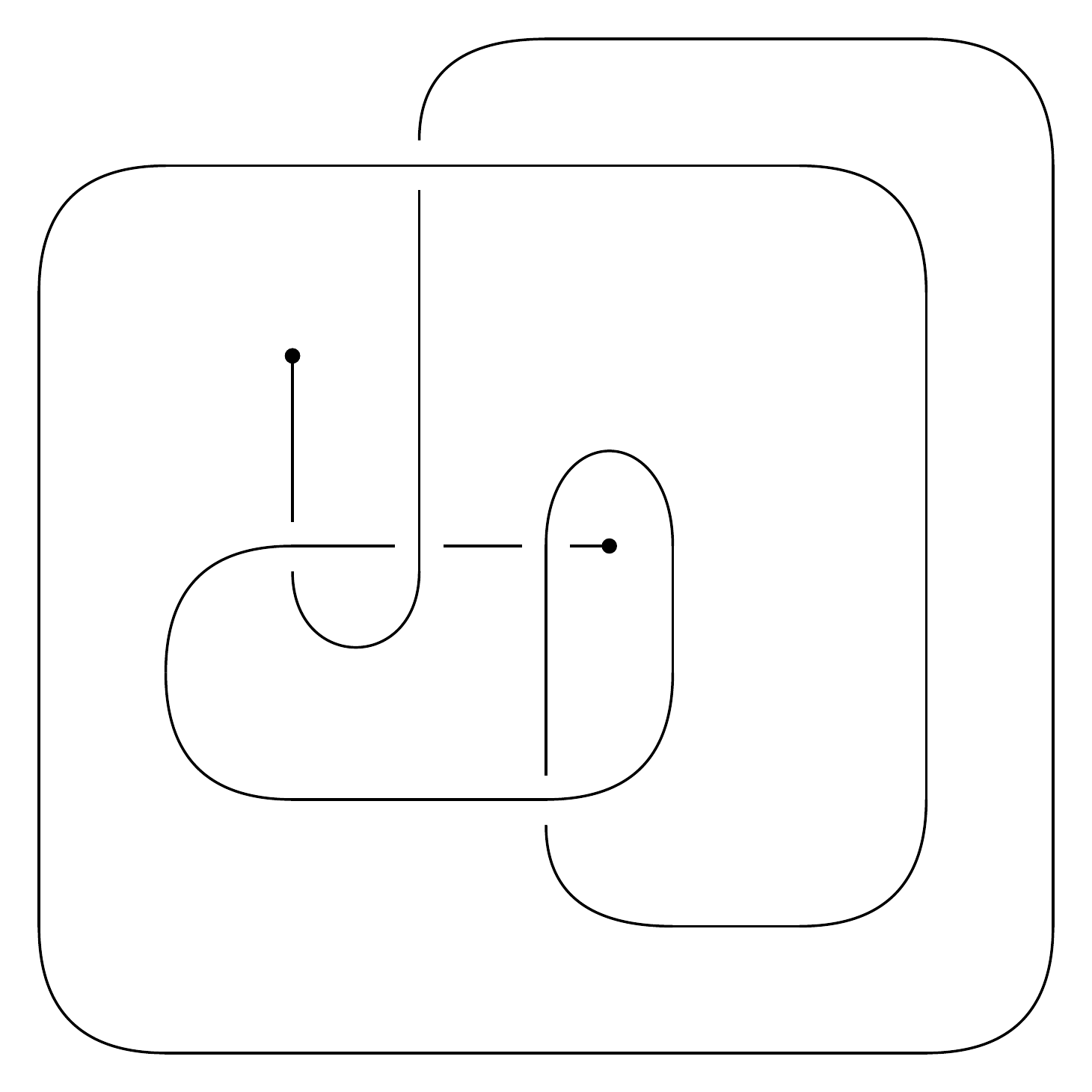}\\
\textcolor{black}{$5_{409}$}
\vspace{1cm}
\end{minipage}
\begin{minipage}[t]{.25\linewidth}
\centering
\includegraphics[width=0.9\textwidth,height=3.5cm,keepaspectratio]{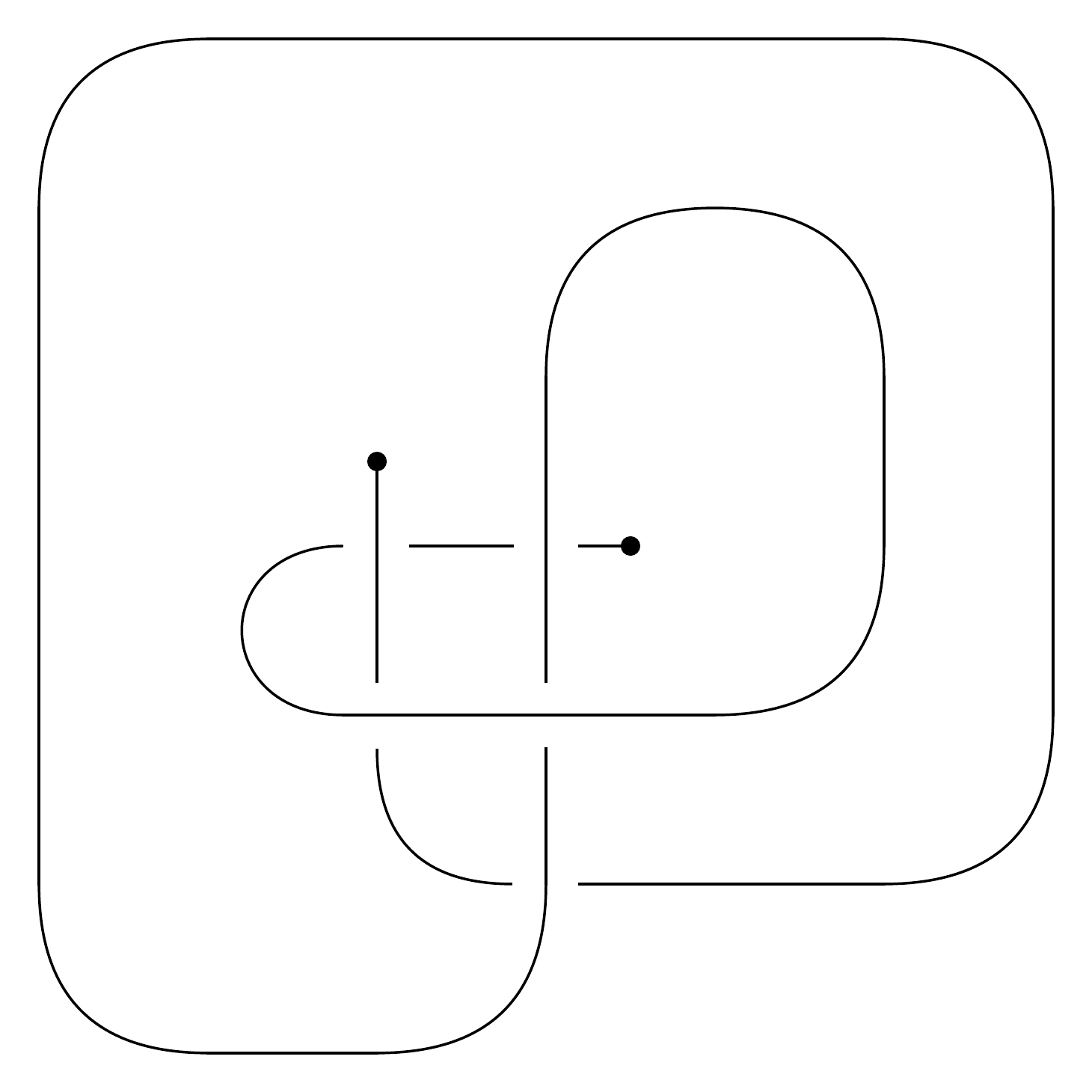}\\
\textcolor{black}{$5_{410}$}
\vspace{1cm}
\end{minipage}
\begin{minipage}[t]{.25\linewidth}
\centering
\includegraphics[width=0.9\textwidth,height=3.5cm,keepaspectratio]{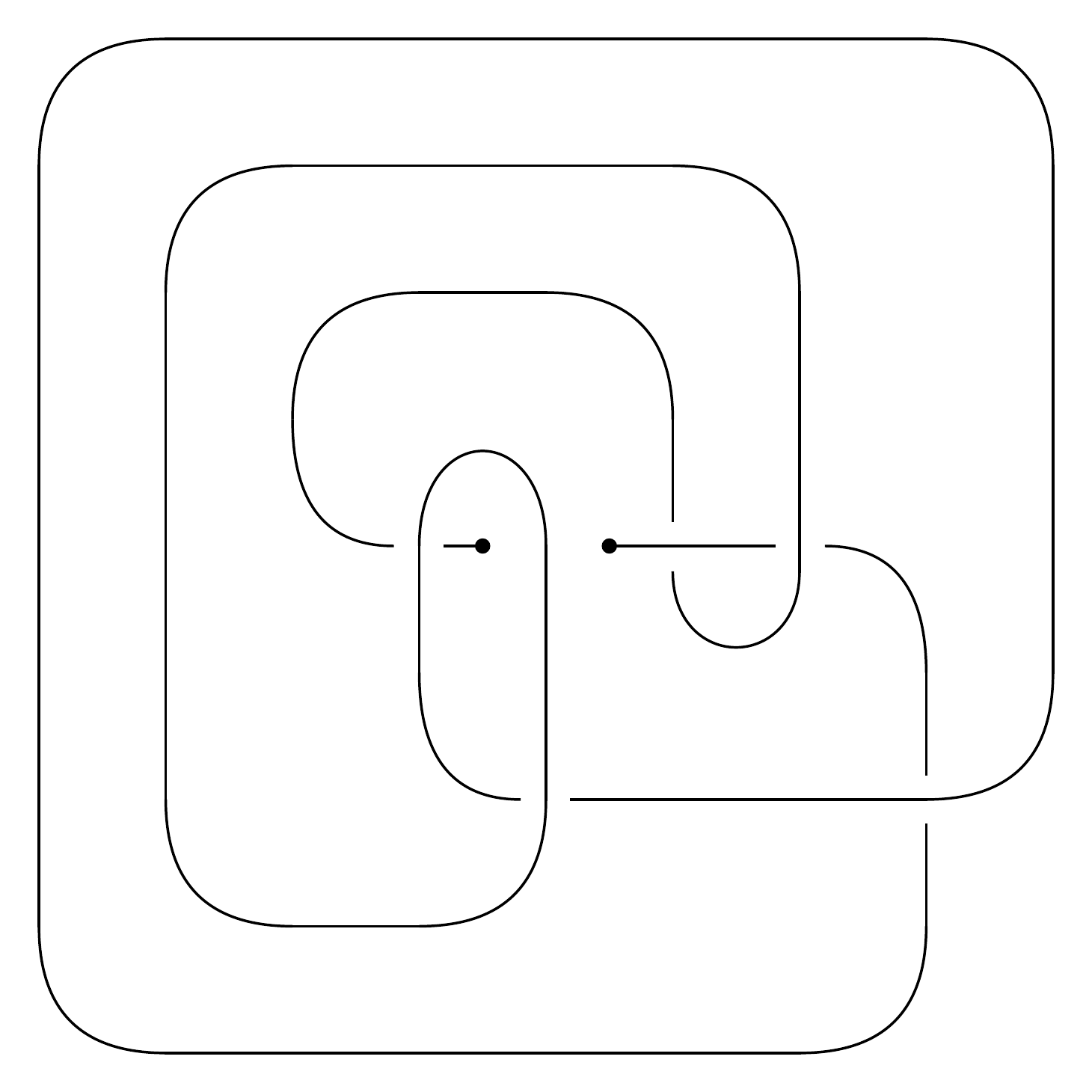}\\
\textcolor{black}{$5_{411}$}
\vspace{1cm}
\end{minipage}
\begin{minipage}[t]{.25\linewidth}
\centering
\includegraphics[width=0.9\textwidth,height=3.5cm,keepaspectratio]{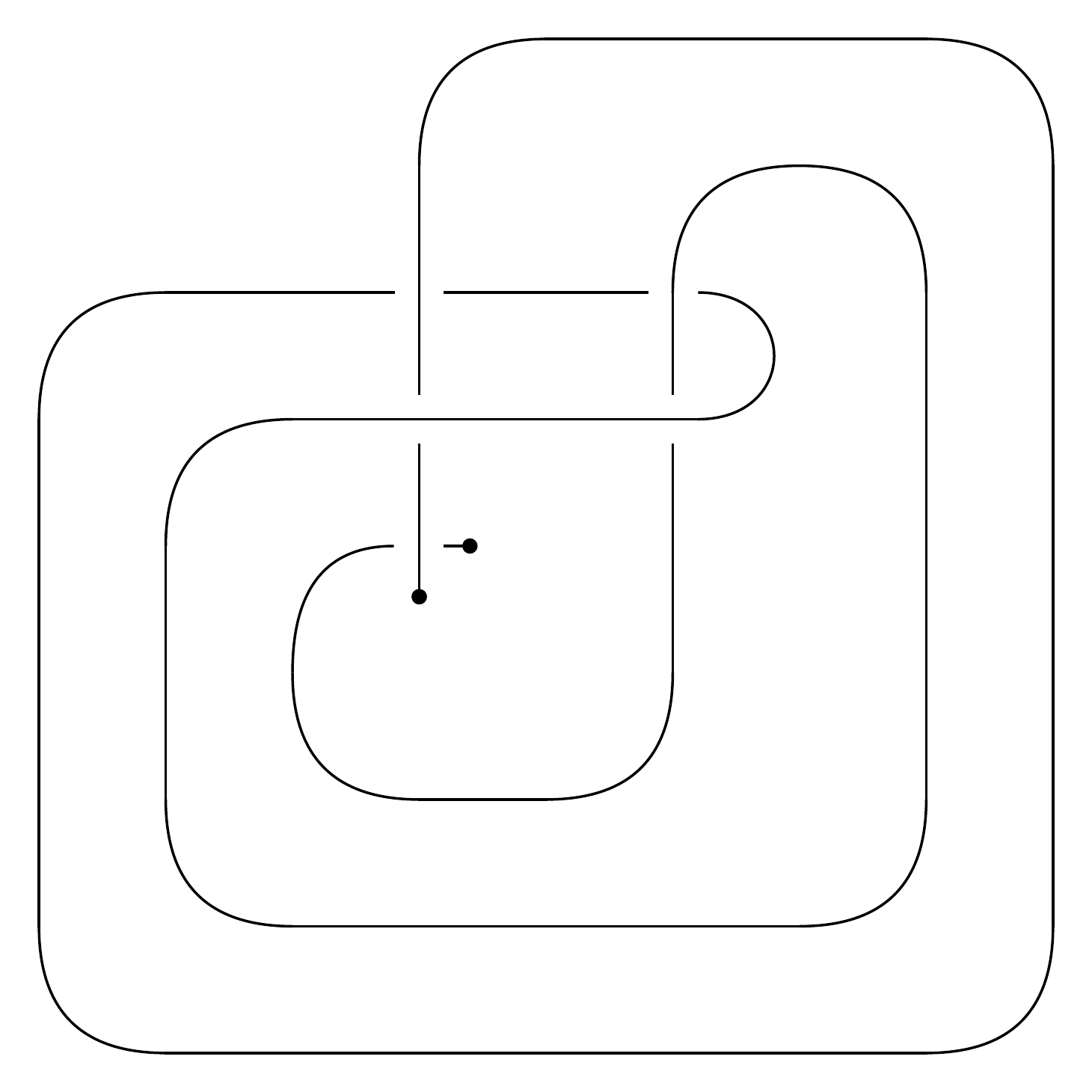}\\
\textcolor{black}{$5_{412}$}
\vspace{1cm}
\end{minipage}
\begin{minipage}[t]{.25\linewidth}
\centering
\includegraphics[width=0.9\textwidth,height=3.5cm,keepaspectratio]{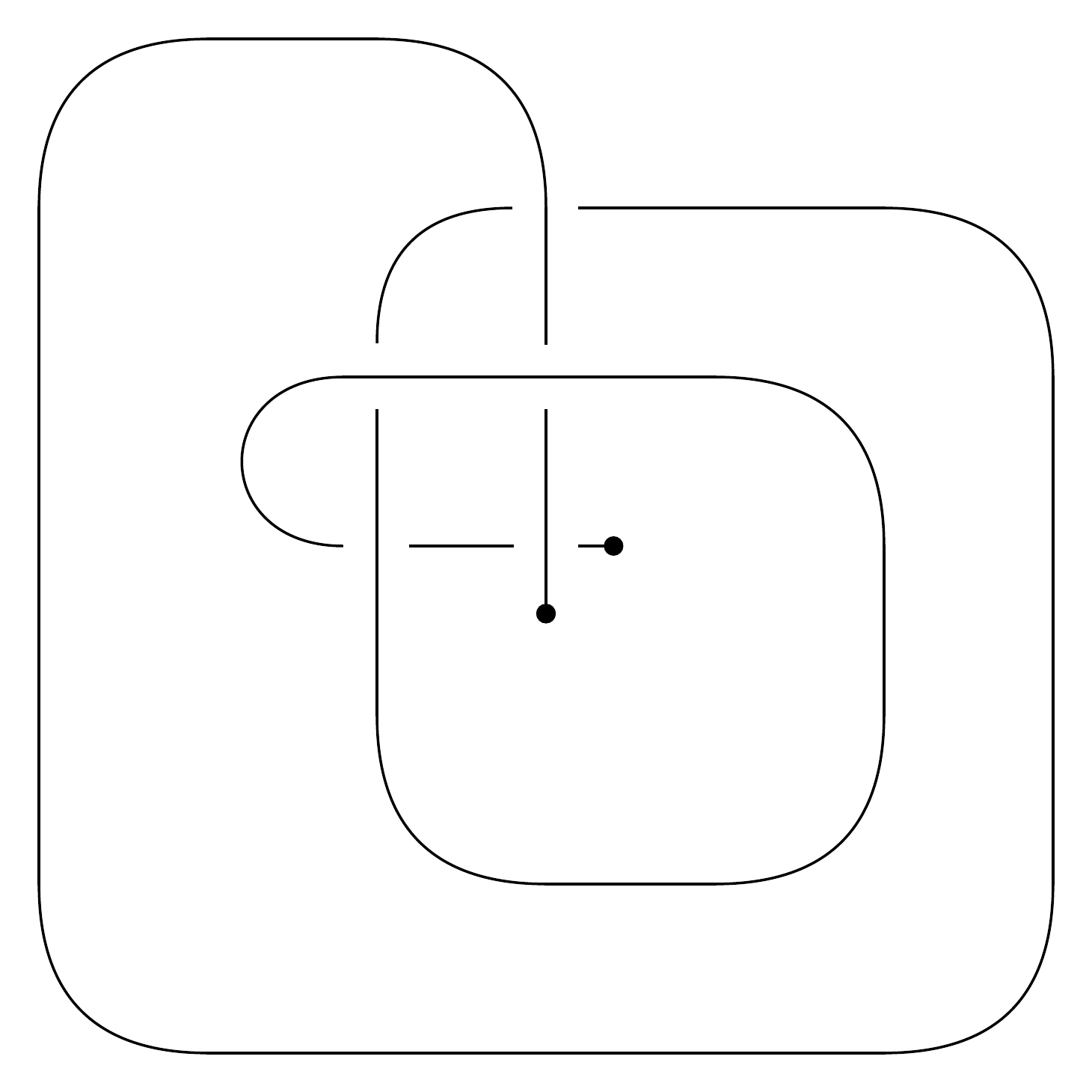}\\
\textcolor{black}{$5_{413}$}
\vspace{1cm}
\end{minipage}
\begin{minipage}[t]{.25\linewidth}
\centering
\includegraphics[width=0.9\textwidth,height=3.5cm,keepaspectratio]{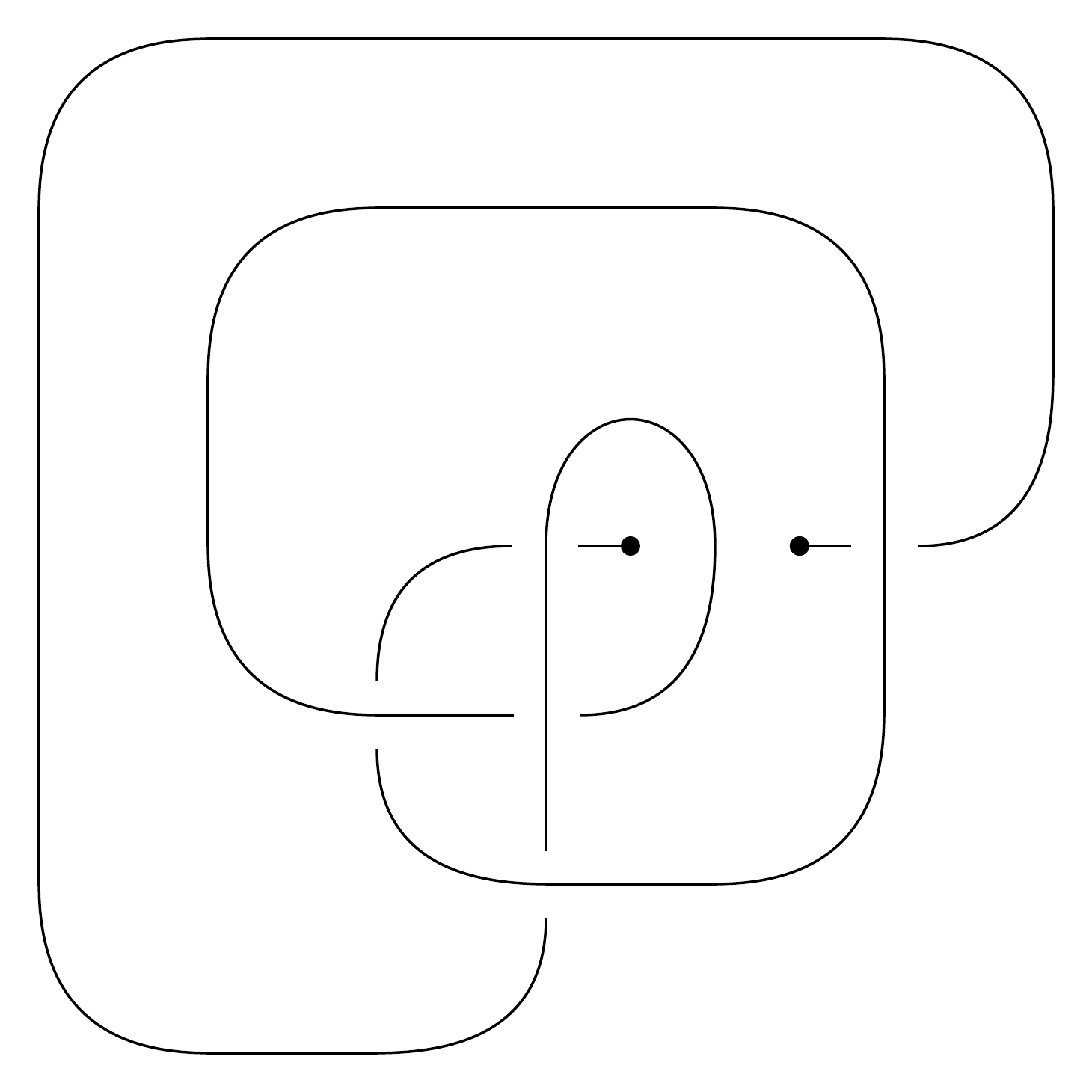}\\
\textcolor{black}{$5_{414}$}
\vspace{1cm}
\end{minipage}
\begin{minipage}[t]{.25\linewidth}
\centering
\includegraphics[width=0.9\textwidth,height=3.5cm,keepaspectratio]{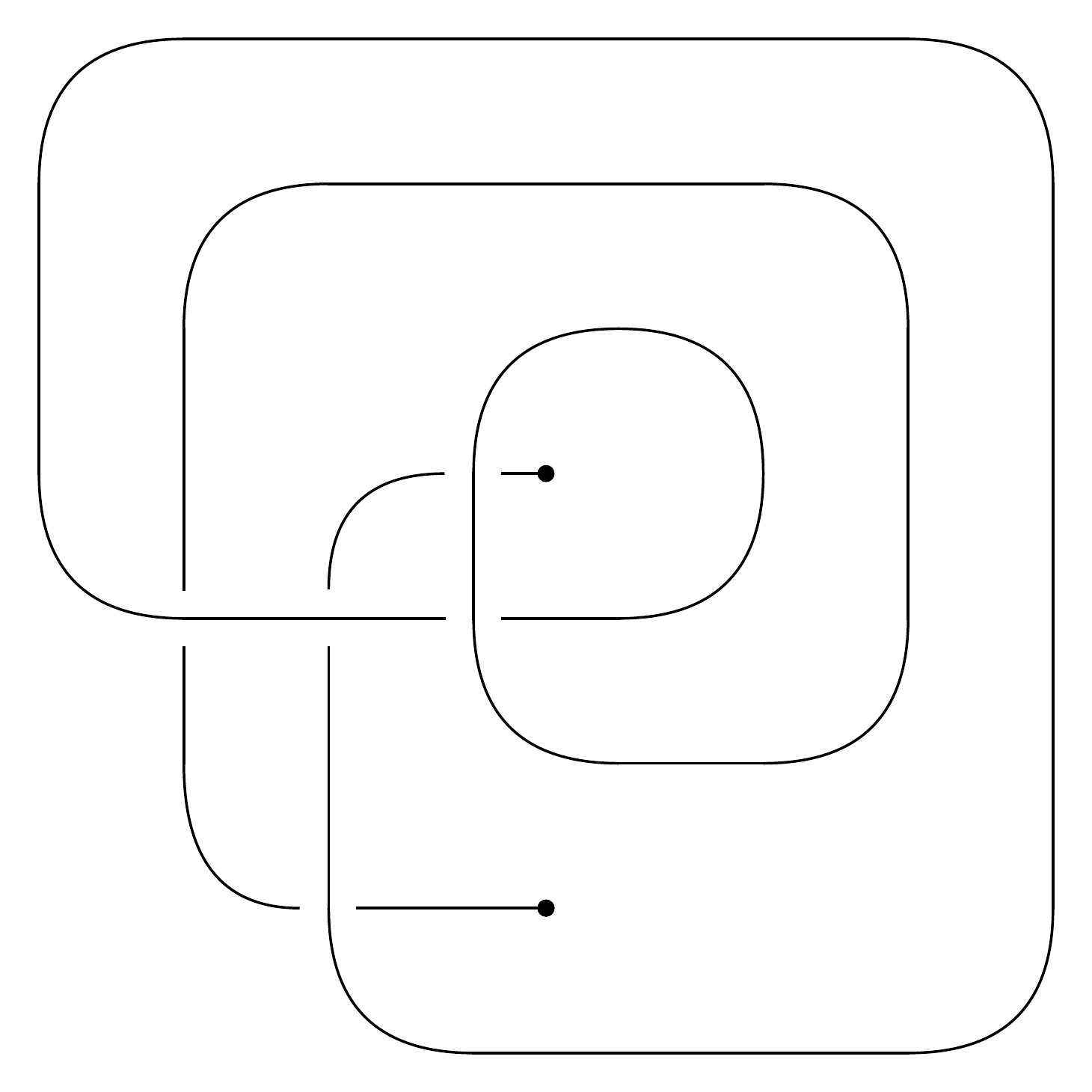}\\
\textcolor{black}{$5_{415}$}
\vspace{1cm}
\end{minipage}
\begin{minipage}[t]{.25\linewidth}
\centering
\includegraphics[width=0.9\textwidth,height=3.5cm,keepaspectratio]{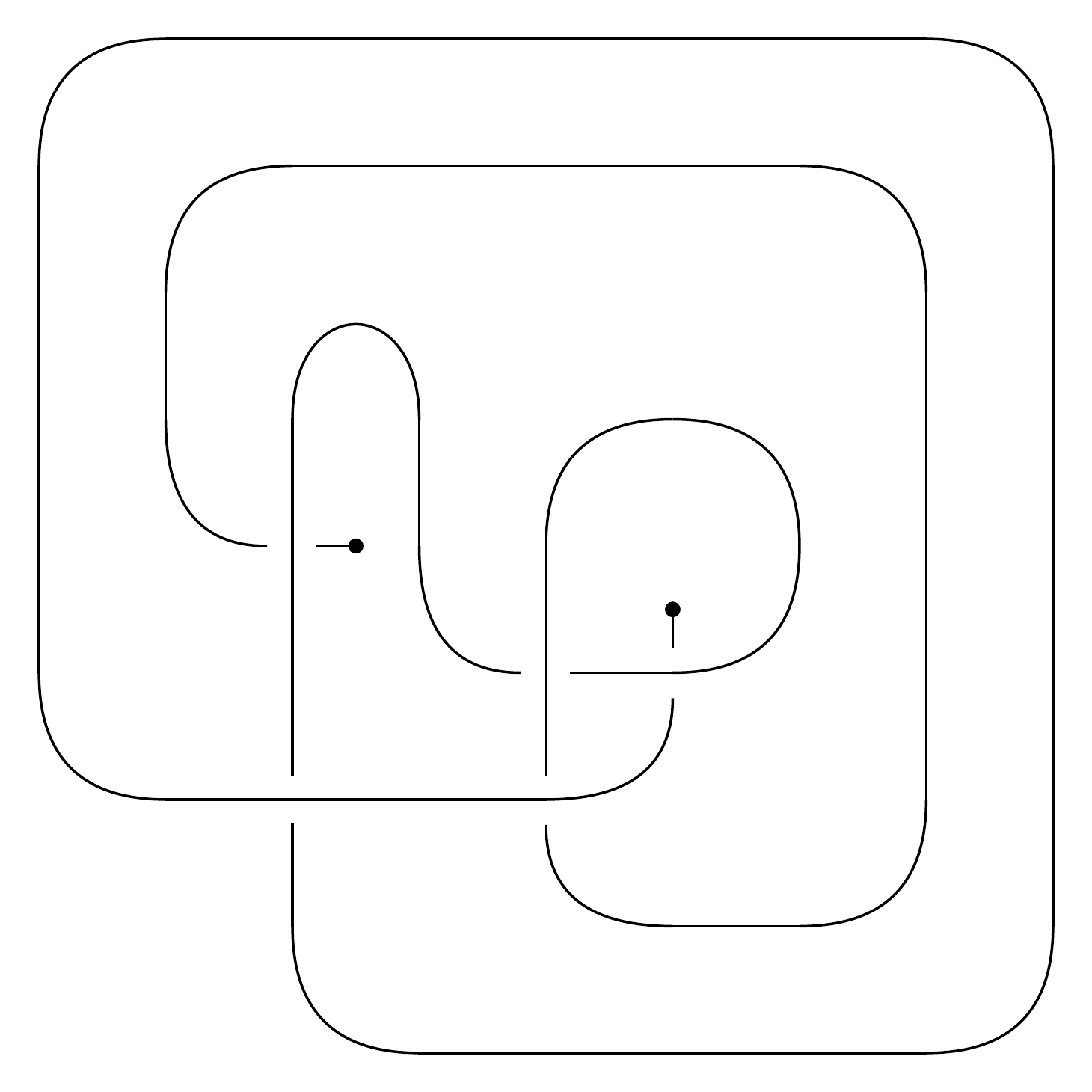}\\
\textcolor{black}{$5_{416}$}
\vspace{1cm}
\end{minipage}
\begin{minipage}[t]{.25\linewidth}
\centering
\includegraphics[width=0.9\textwidth,height=3.5cm,keepaspectratio]{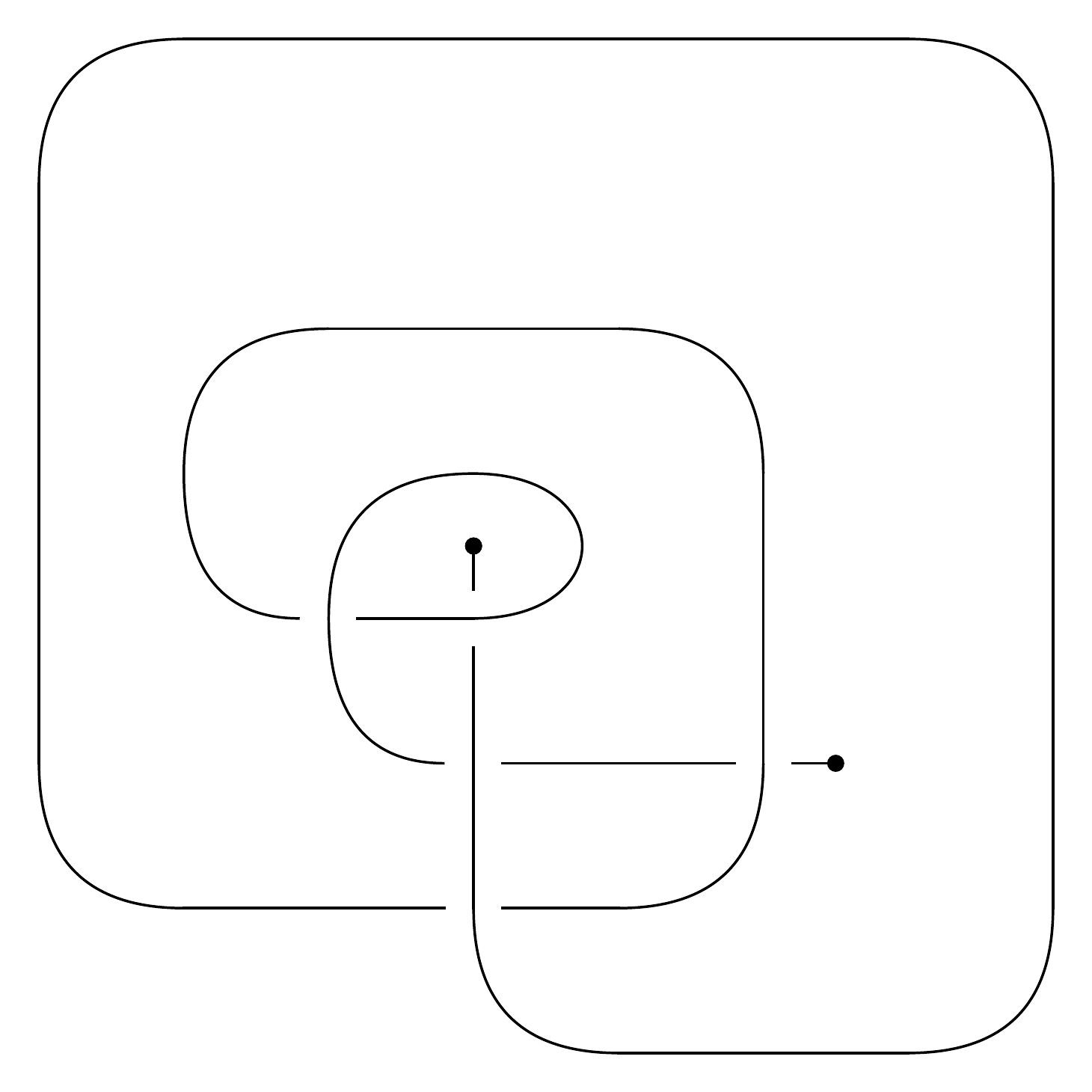}\\
\textcolor{black}{$5_{417}$}
\vspace{1cm}
\end{minipage}
\begin{minipage}[t]{.25\linewidth}
\centering
\includegraphics[width=0.9\textwidth,height=3.5cm,keepaspectratio]{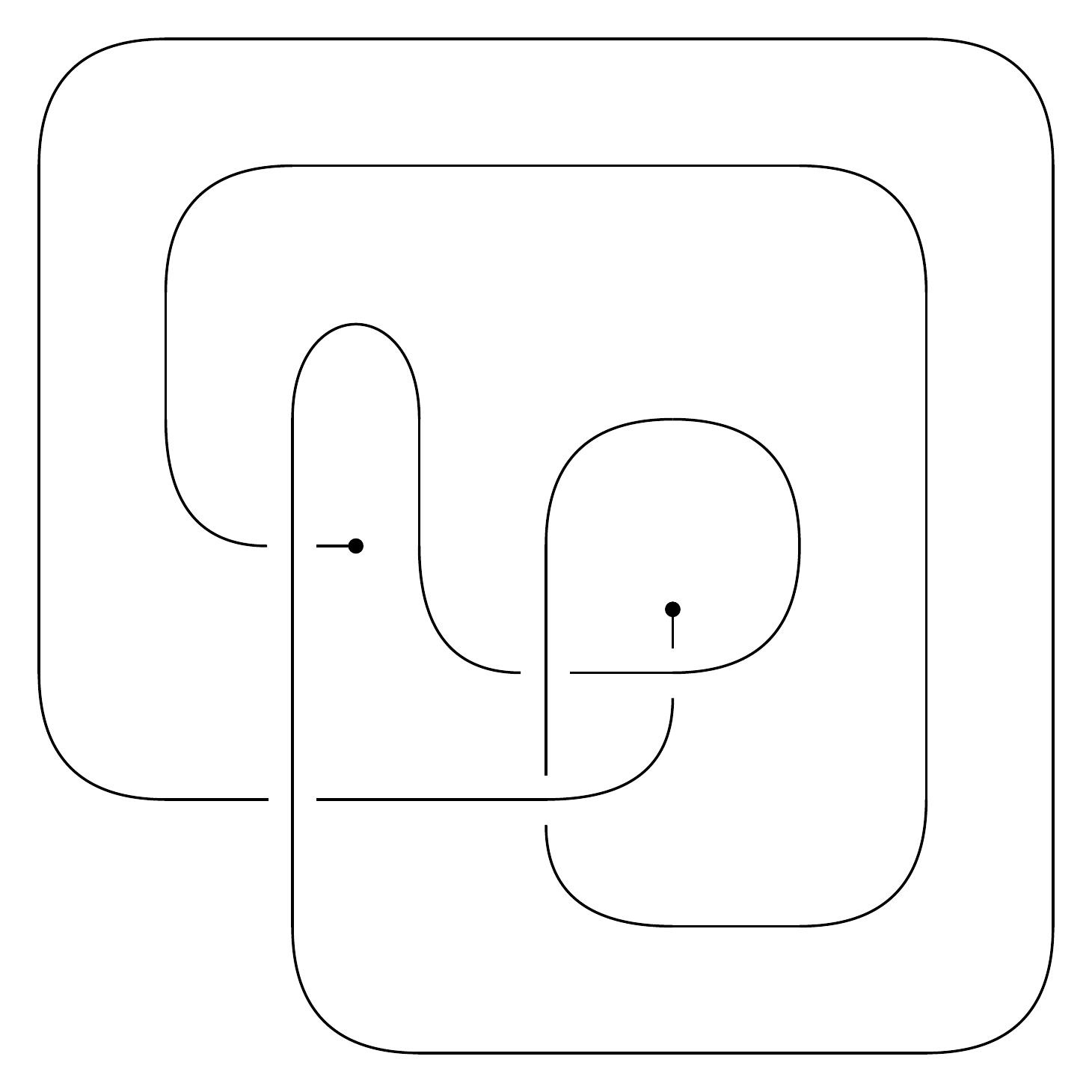}\\
\textcolor{black}{$5_{418}$}
\vspace{1cm}
\end{minipage}
\begin{minipage}[t]{.25\linewidth}
\centering
\includegraphics[width=0.9\textwidth,height=3.5cm,keepaspectratio]{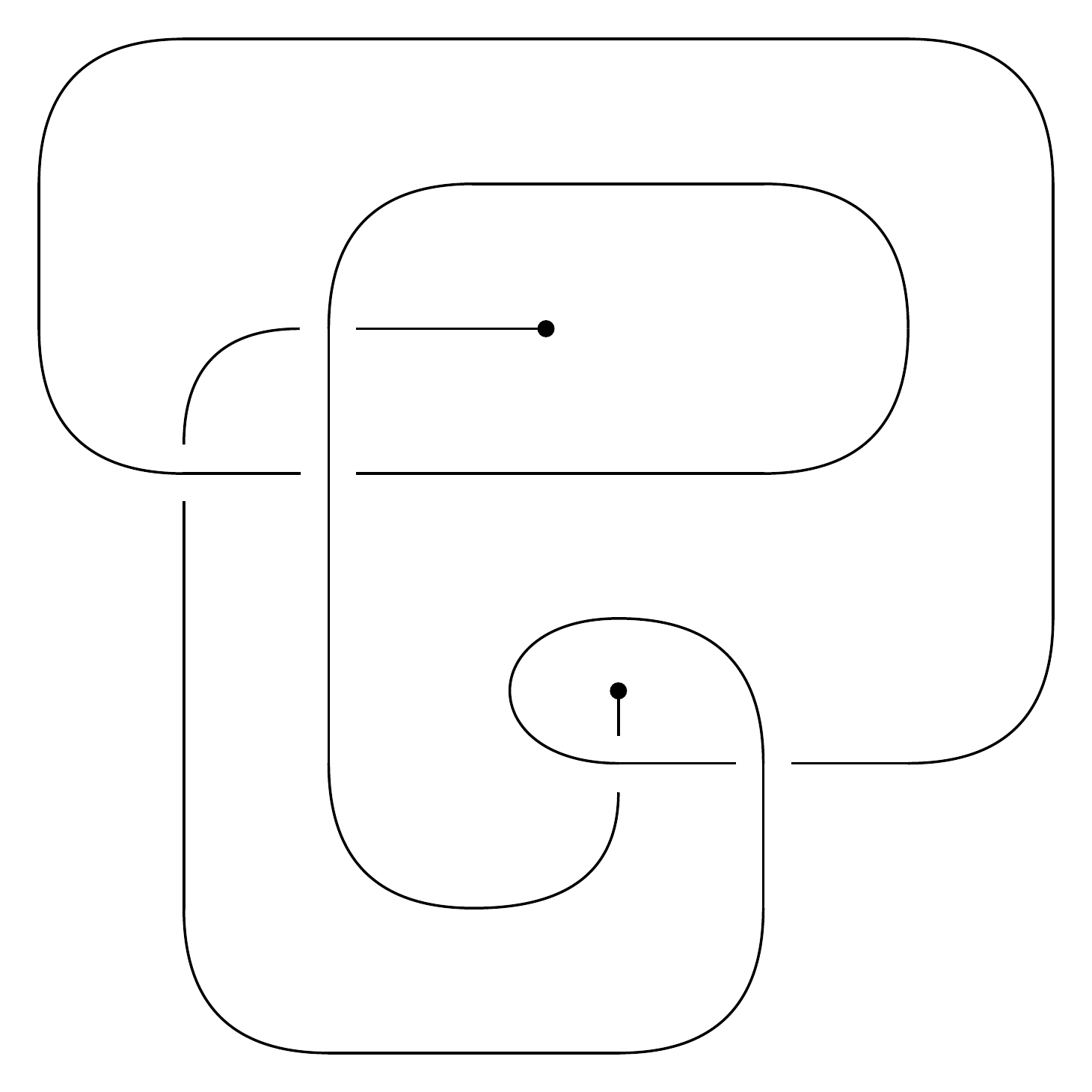}\\
\textcolor{black}{$5_{419}$}
\vspace{1cm}
\end{minipage}
\begin{minipage}[t]{.25\linewidth}
\centering
\includegraphics[width=0.9\textwidth,height=3.5cm,keepaspectratio]{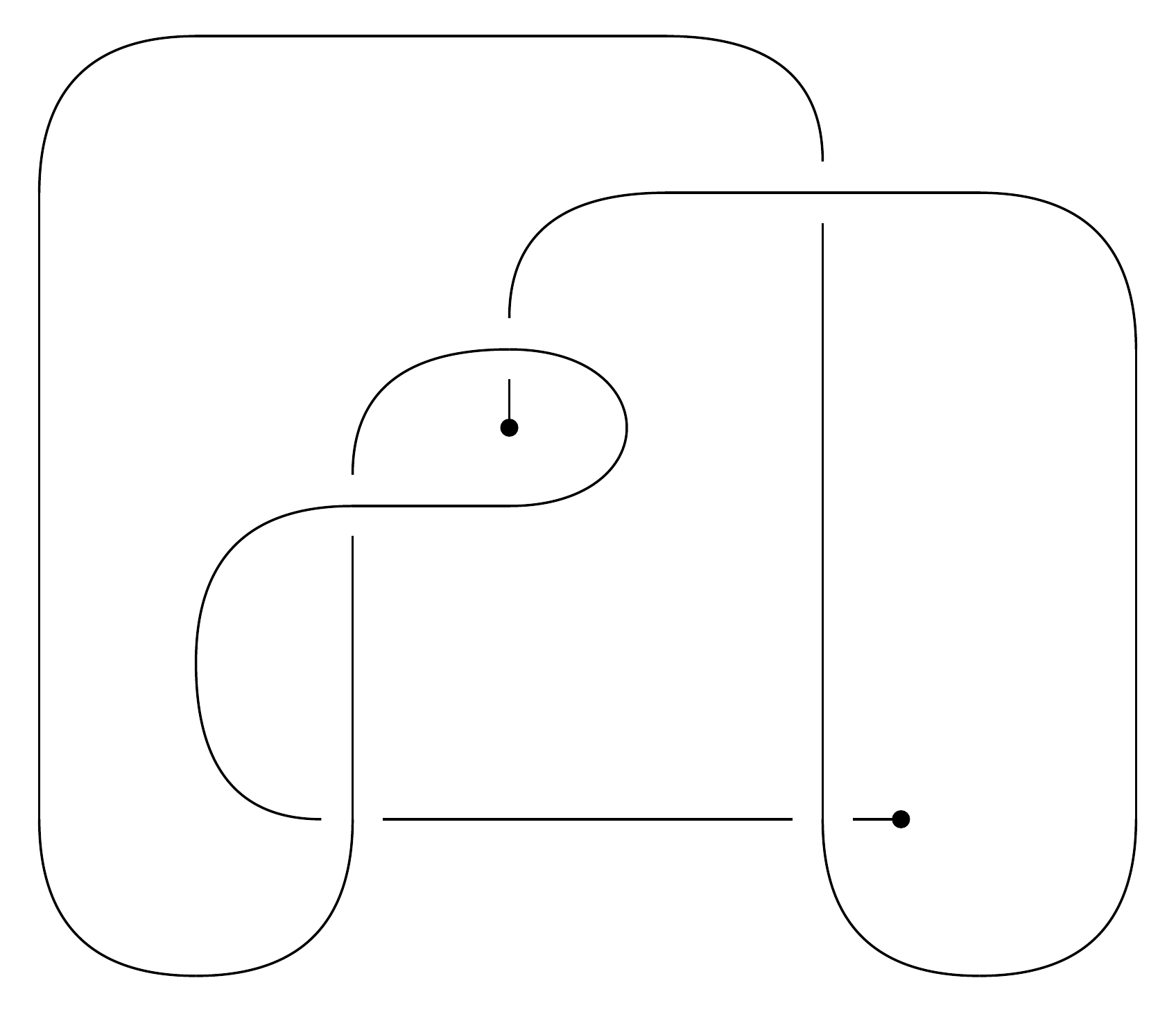}\\
\textcolor{black}{$5_{420}$}
\vspace{1cm}
\end{minipage}
\begin{minipage}[t]{.25\linewidth}
\centering
\includegraphics[width=0.9\textwidth,height=3.5cm,keepaspectratio]{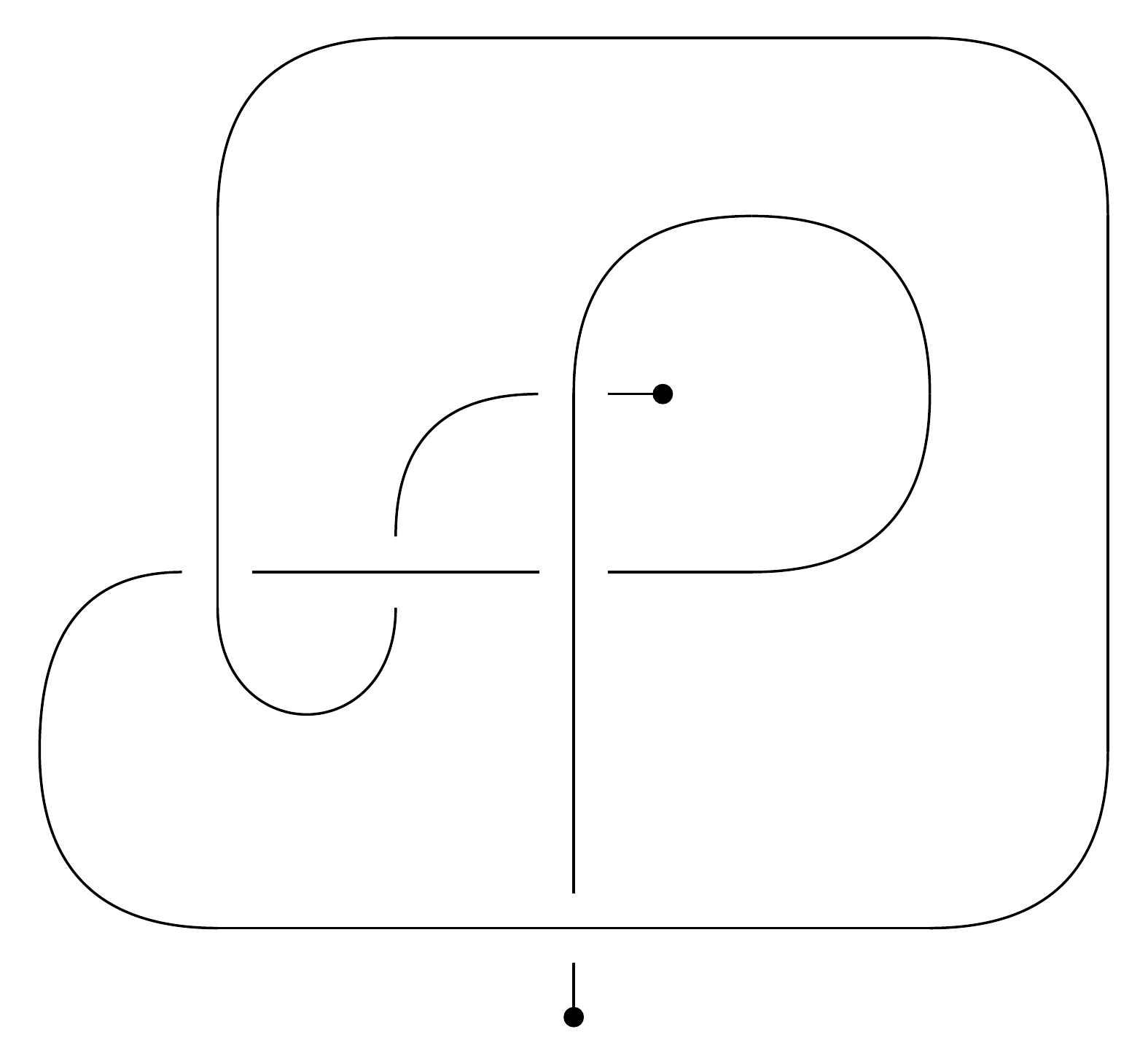}\\
\textcolor{black}{$5_{421}$}
\vspace{1cm}
\end{minipage}
\begin{minipage}[t]{.25\linewidth}
\centering
\includegraphics[width=0.9\textwidth,height=3.5cm,keepaspectratio]{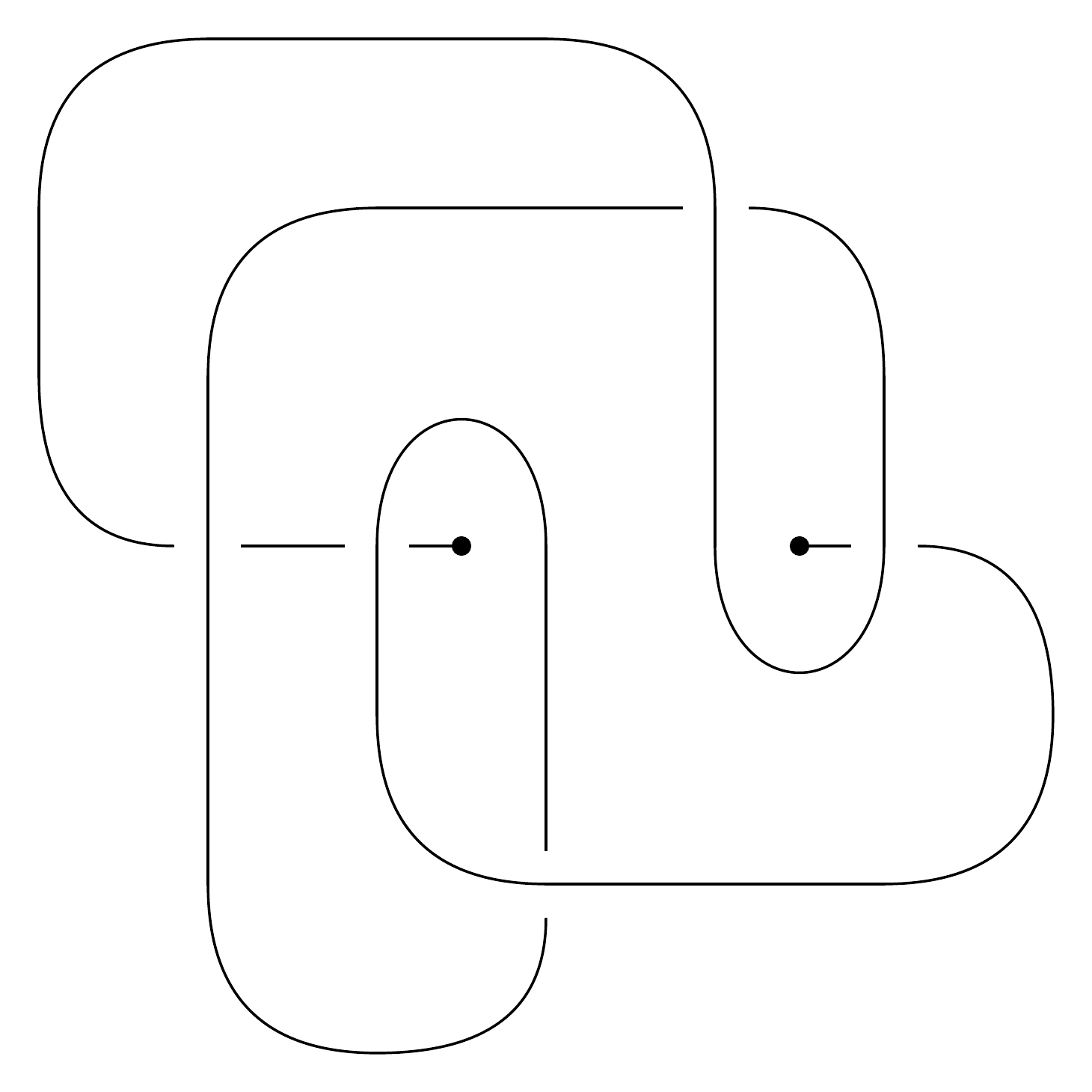}\\
\textcolor{black}{$5_{422}$}
\vspace{1cm}
\end{minipage}
\begin{minipage}[t]{.25\linewidth}
\centering
\includegraphics[width=0.9\textwidth,height=3.5cm,keepaspectratio]{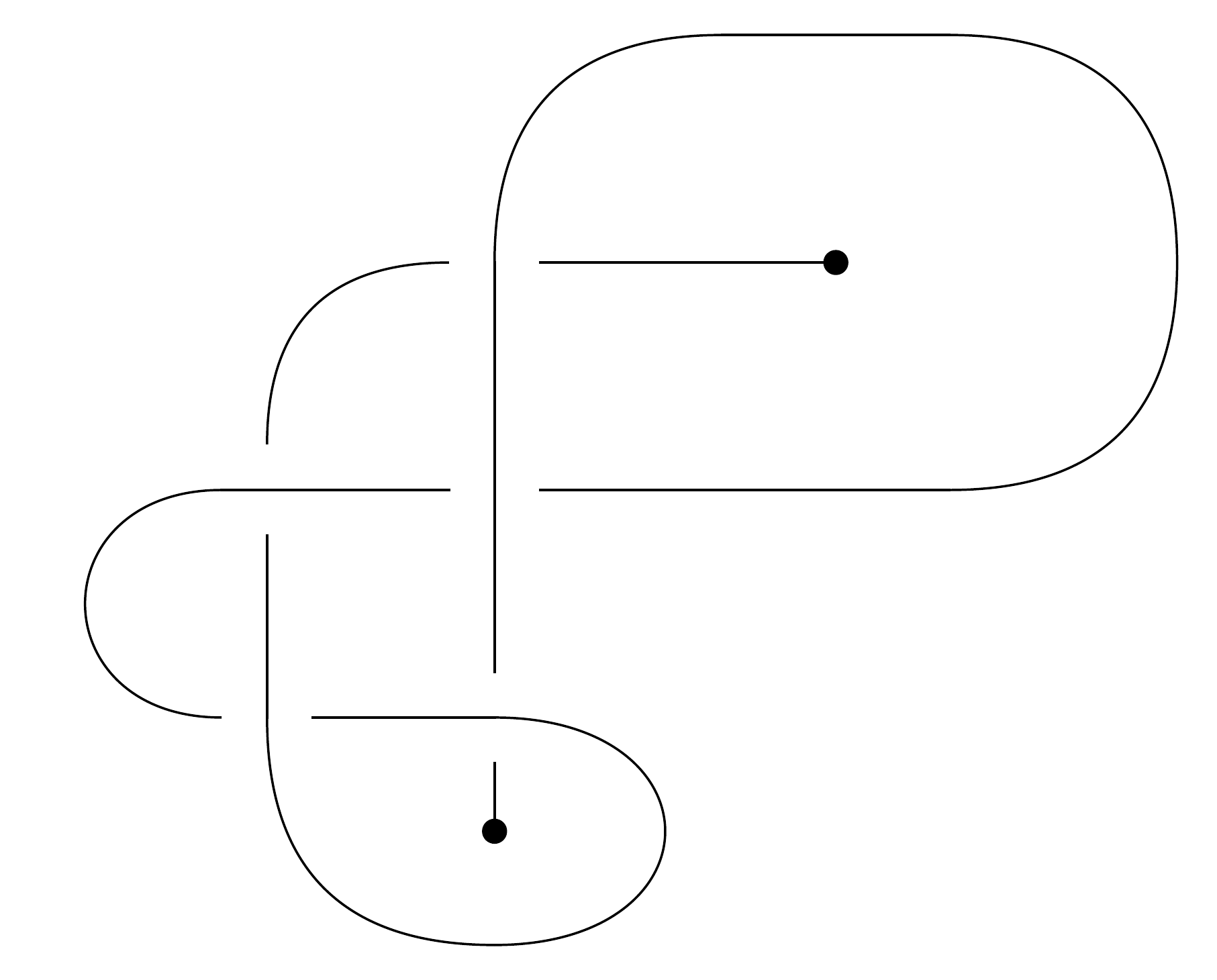}\\
\textcolor{black}{$5_{423}$}
\vspace{1cm}
\end{minipage}
\begin{minipage}[t]{.25\linewidth}
\centering
\includegraphics[width=0.9\textwidth,height=3.5cm,keepaspectratio]{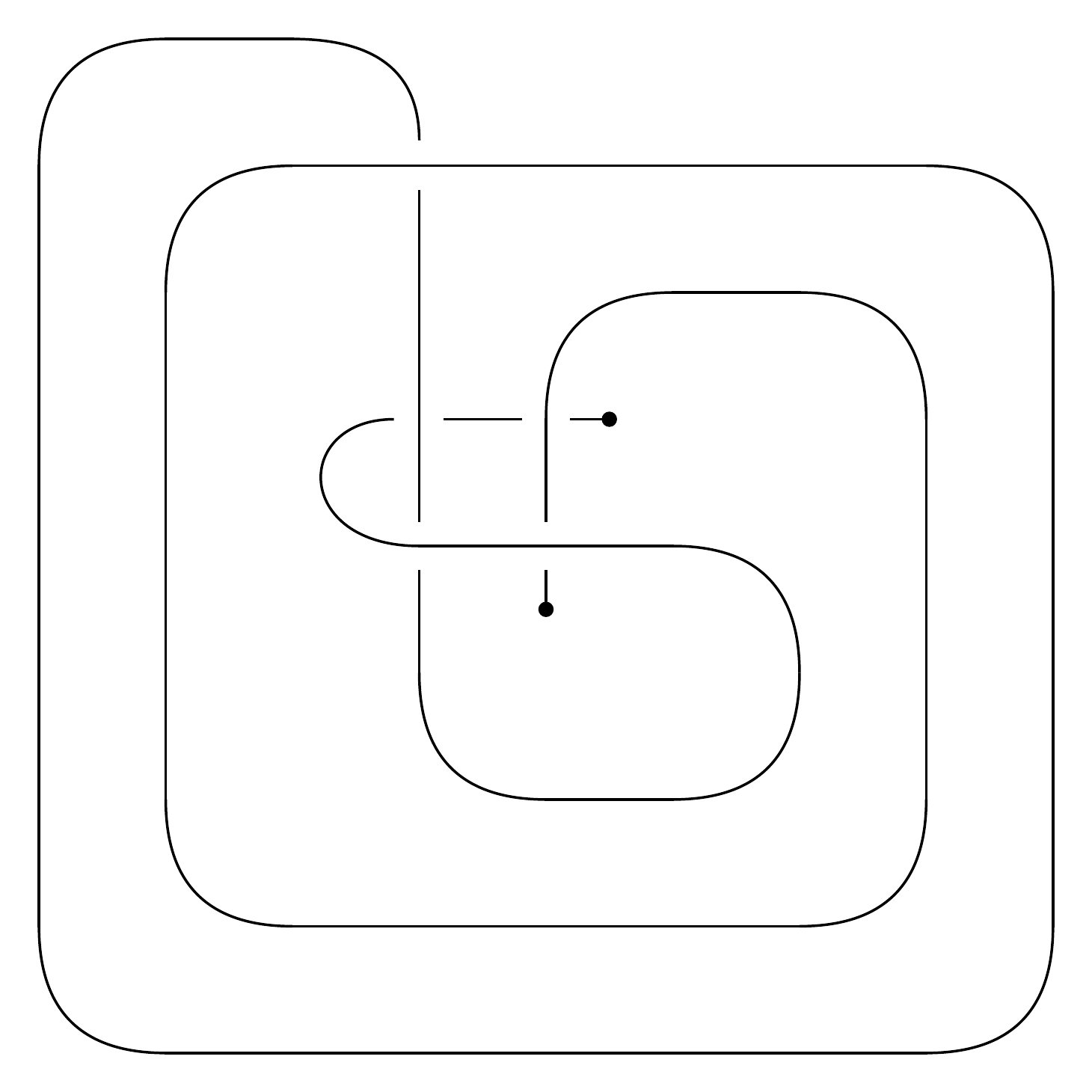}\\
\textcolor{black}{$5_{424}$}
\vspace{1cm}
\end{minipage}
\begin{minipage}[t]{.25\linewidth}
\centering
\includegraphics[width=0.9\textwidth,height=3.5cm,keepaspectratio]{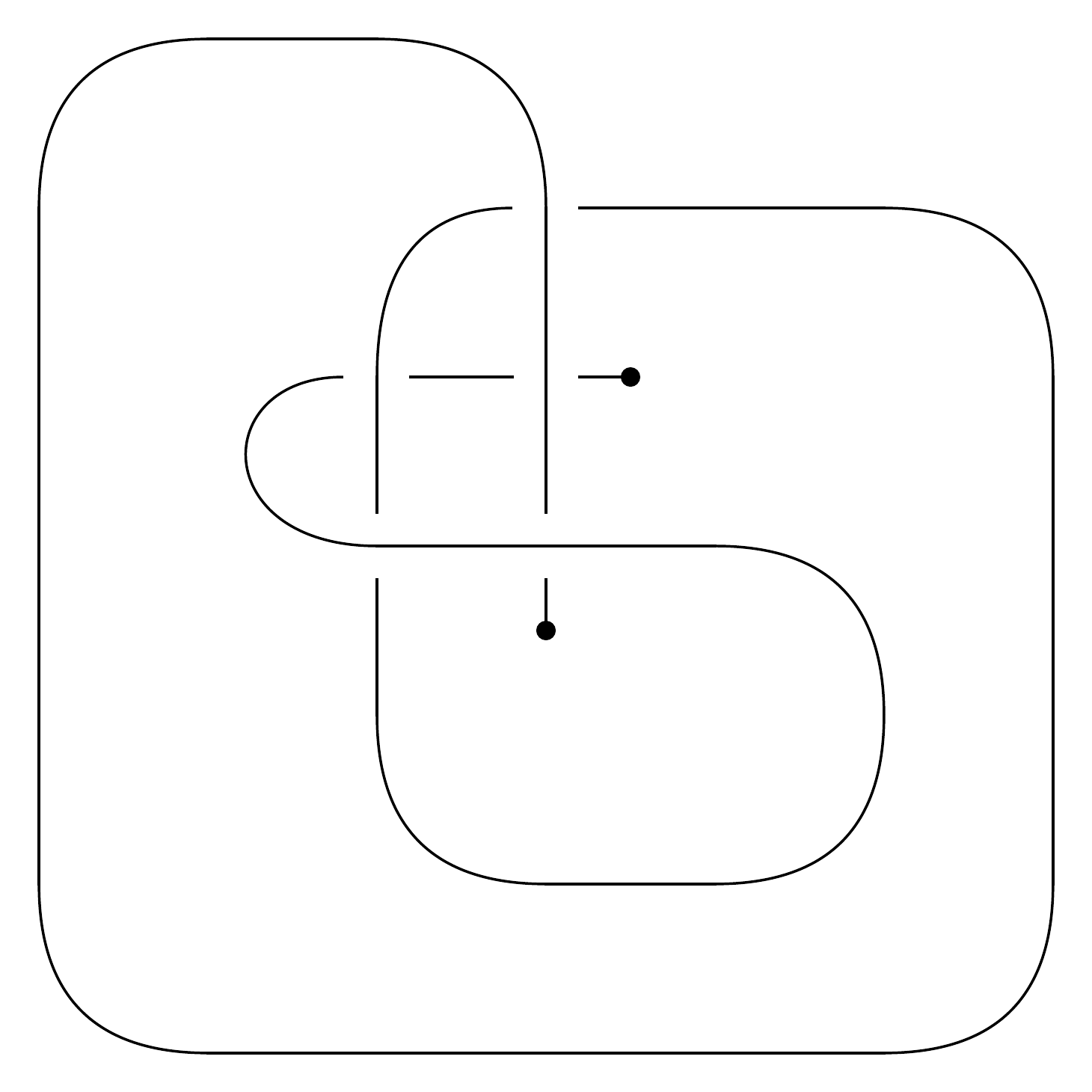}\\
\textcolor{black}{$5_{425}$}
\vspace{1cm}
\end{minipage}
\begin{minipage}[t]{.25\linewidth}
\centering
\includegraphics[width=0.9\textwidth,height=3.5cm,keepaspectratio]{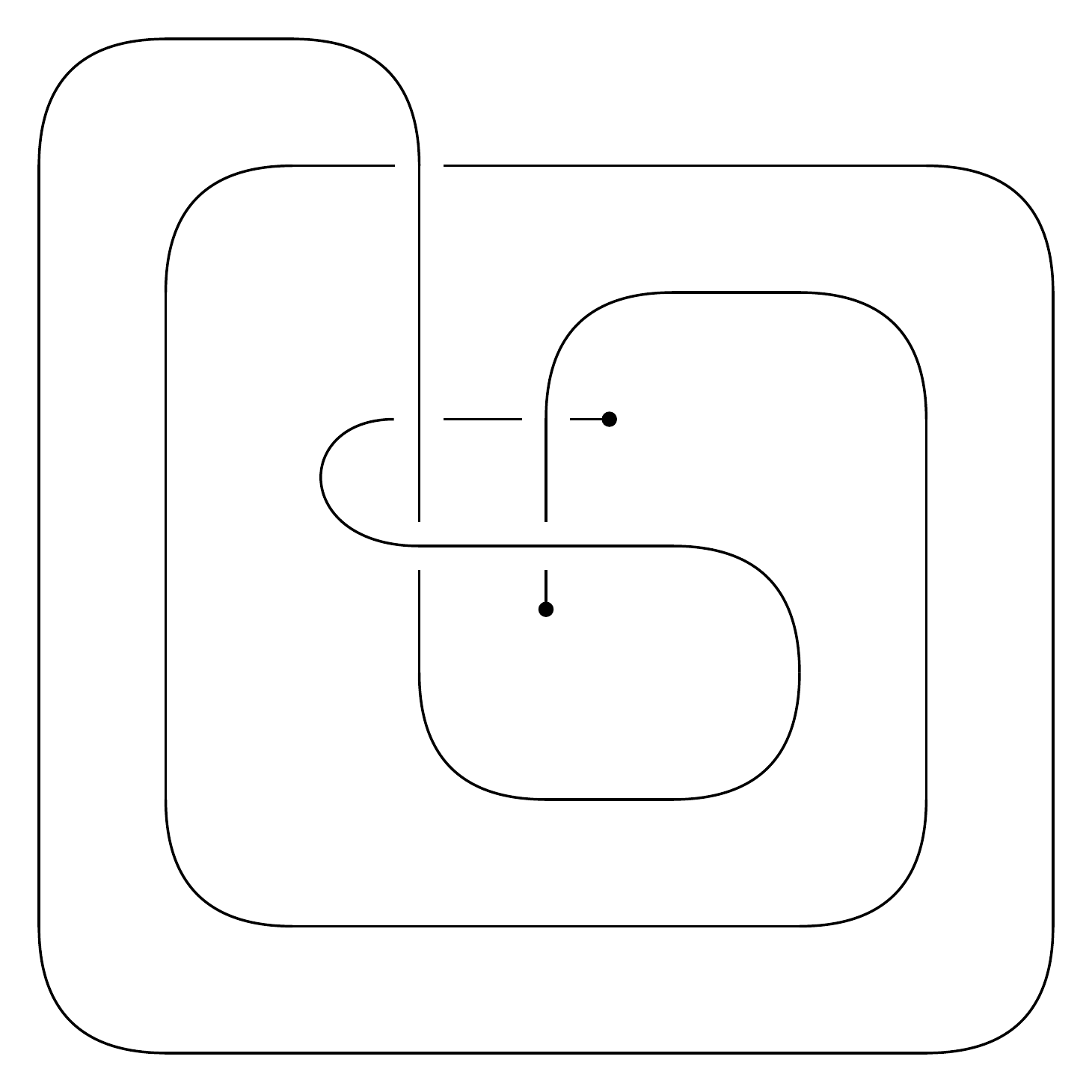}\\
\textcolor{black}{$5_{426}$}
\vspace{1cm}
\end{minipage}
\begin{minipage}[t]{.25\linewidth}
\centering
\includegraphics[width=0.9\textwidth,height=3.5cm,keepaspectratio]{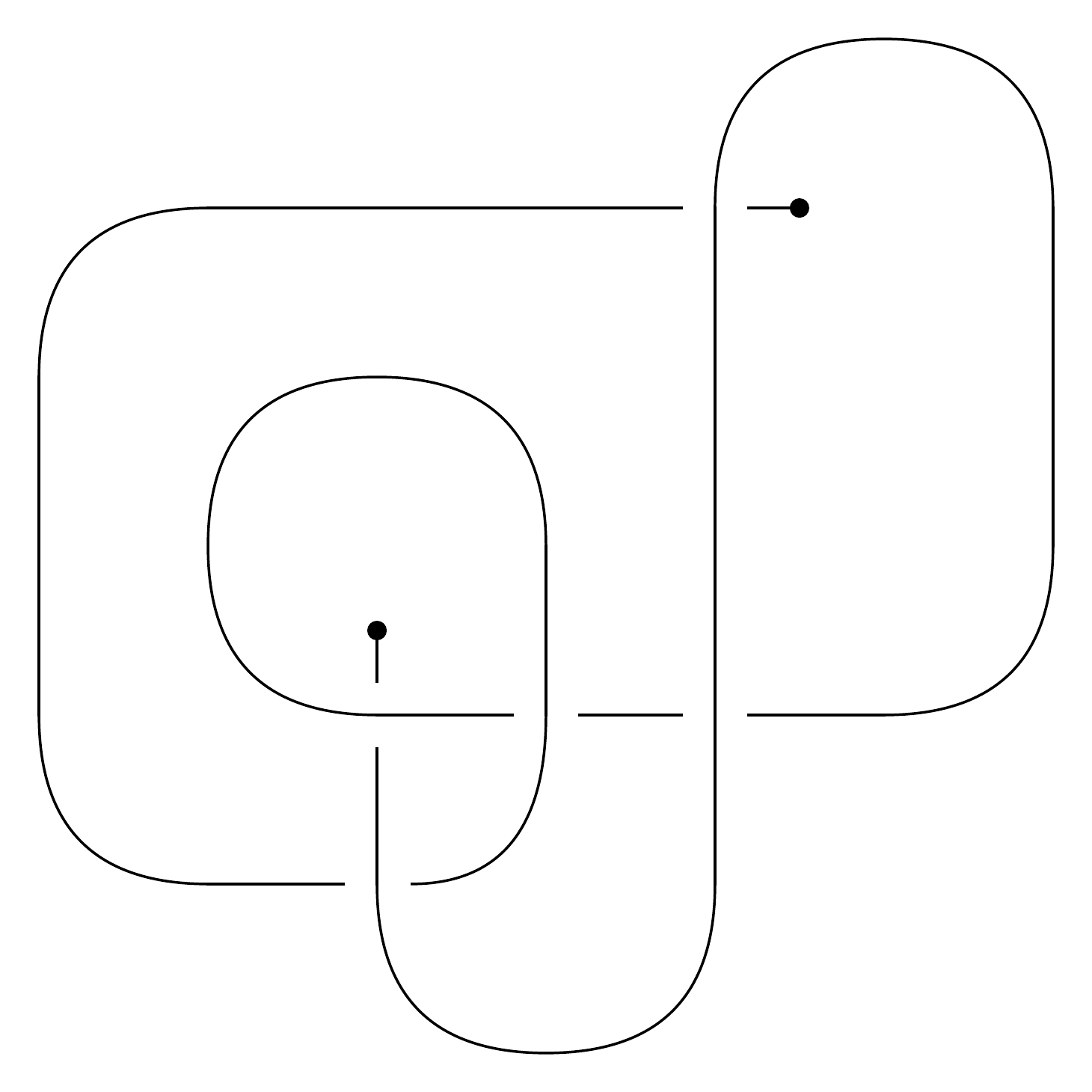}\\
\textcolor{black}{$5_{427}$}
\vspace{1cm}
\end{minipage}
\begin{minipage}[t]{.25\linewidth}
\centering
\includegraphics[width=0.9\textwidth,height=3.5cm,keepaspectratio]{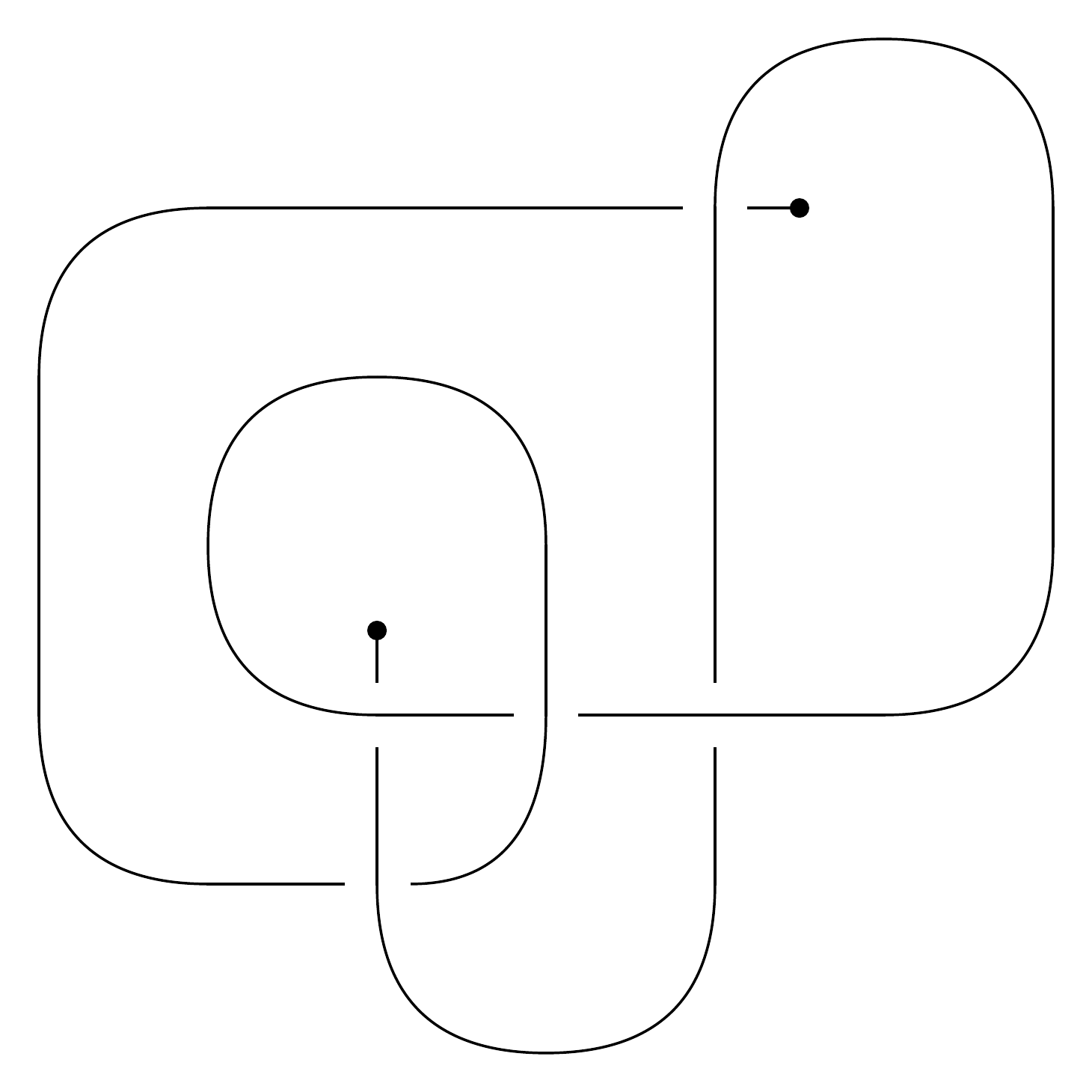}\\
\textcolor{black}{$5_{428}$}
\vspace{1cm}
\end{minipage}
\begin{minipage}[t]{.25\linewidth}
\centering
\includegraphics[width=0.9\textwidth,height=3.5cm,keepaspectratio]{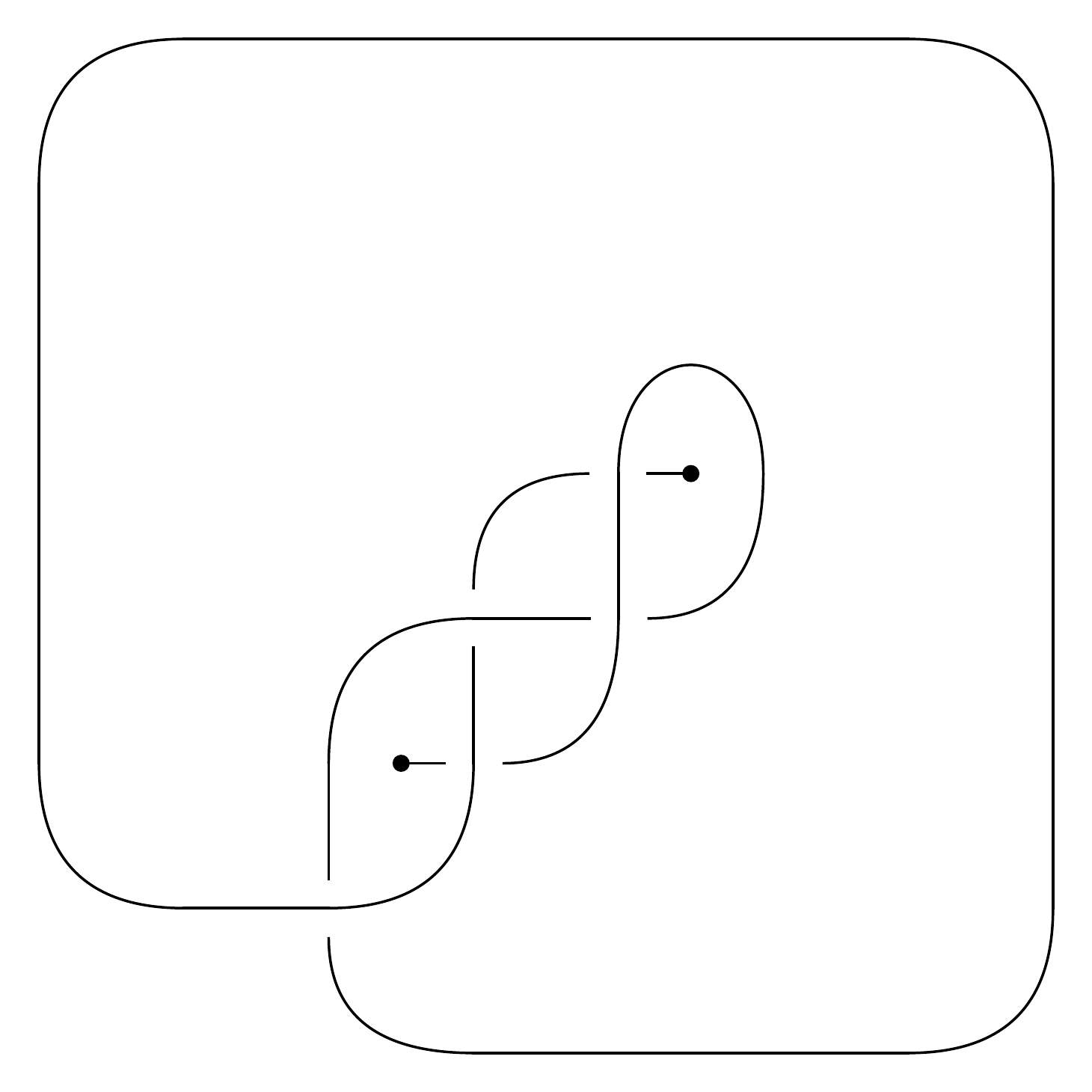}\\
\textcolor{black}{$5_{429}$}
\vspace{1cm}
\end{minipage}
\begin{minipage}[t]{.25\linewidth}
\centering
\includegraphics[width=0.9\textwidth,height=3.5cm,keepaspectratio]{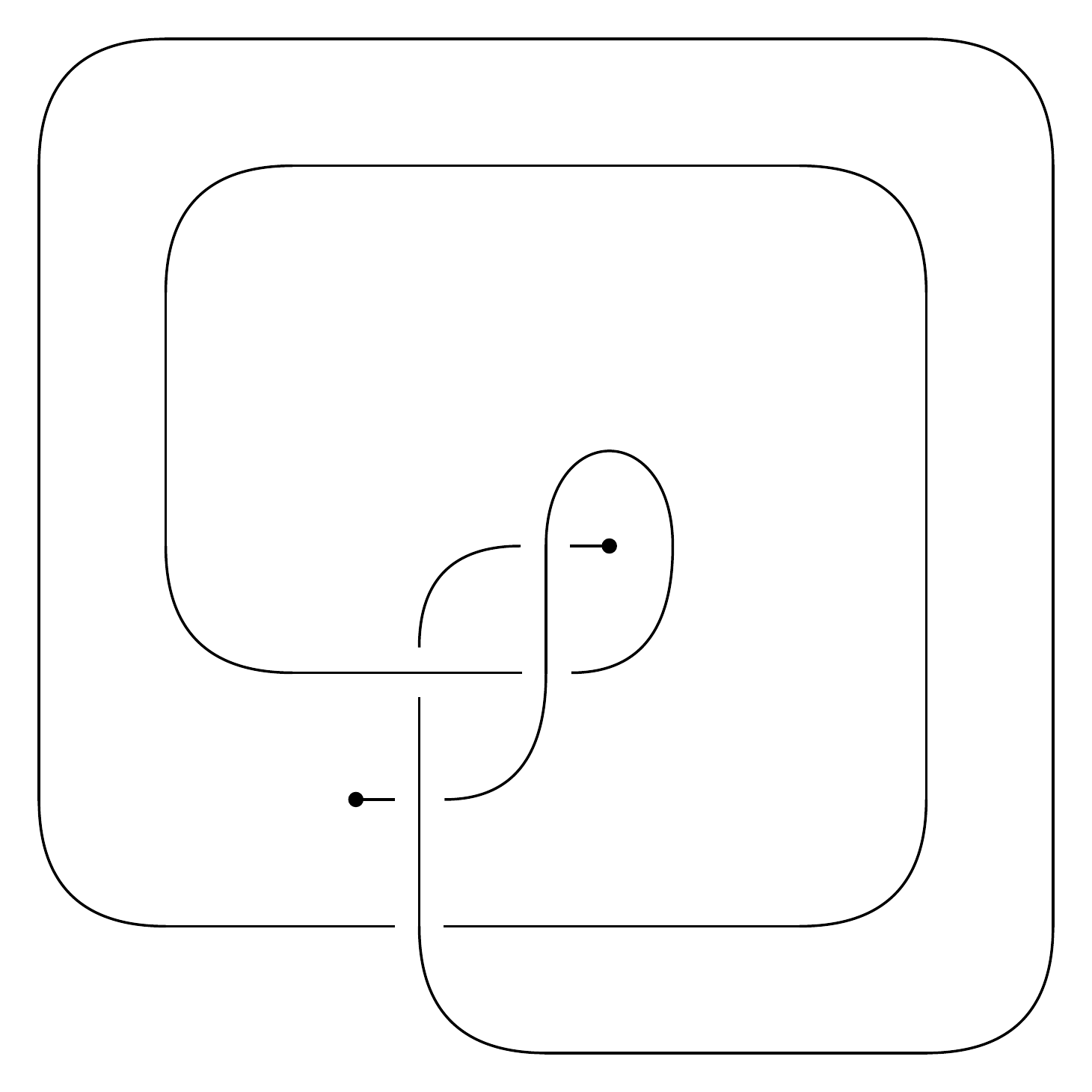}\\
\textcolor{black}{$5_{430}$}
\vspace{1cm}
\end{minipage}
\begin{minipage}[t]{.25\linewidth}
\centering
\includegraphics[width=0.9\textwidth,height=3.5cm,keepaspectratio]{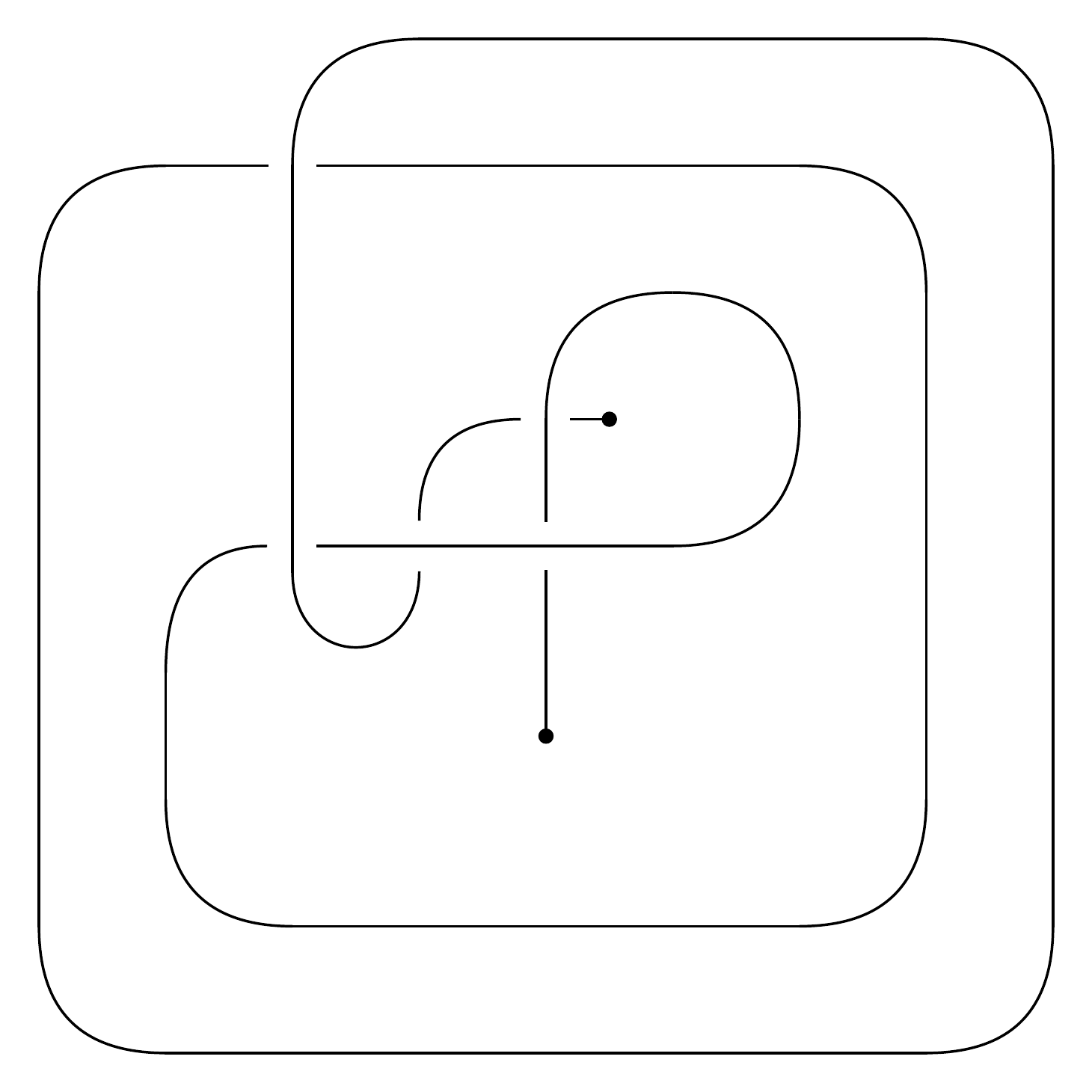}\\
\textcolor{black}{$5_{431}$}
\vspace{1cm}
\end{minipage}
\begin{minipage}[t]{.25\linewidth}
\centering
\includegraphics[width=0.9\textwidth,height=3.5cm,keepaspectratio]{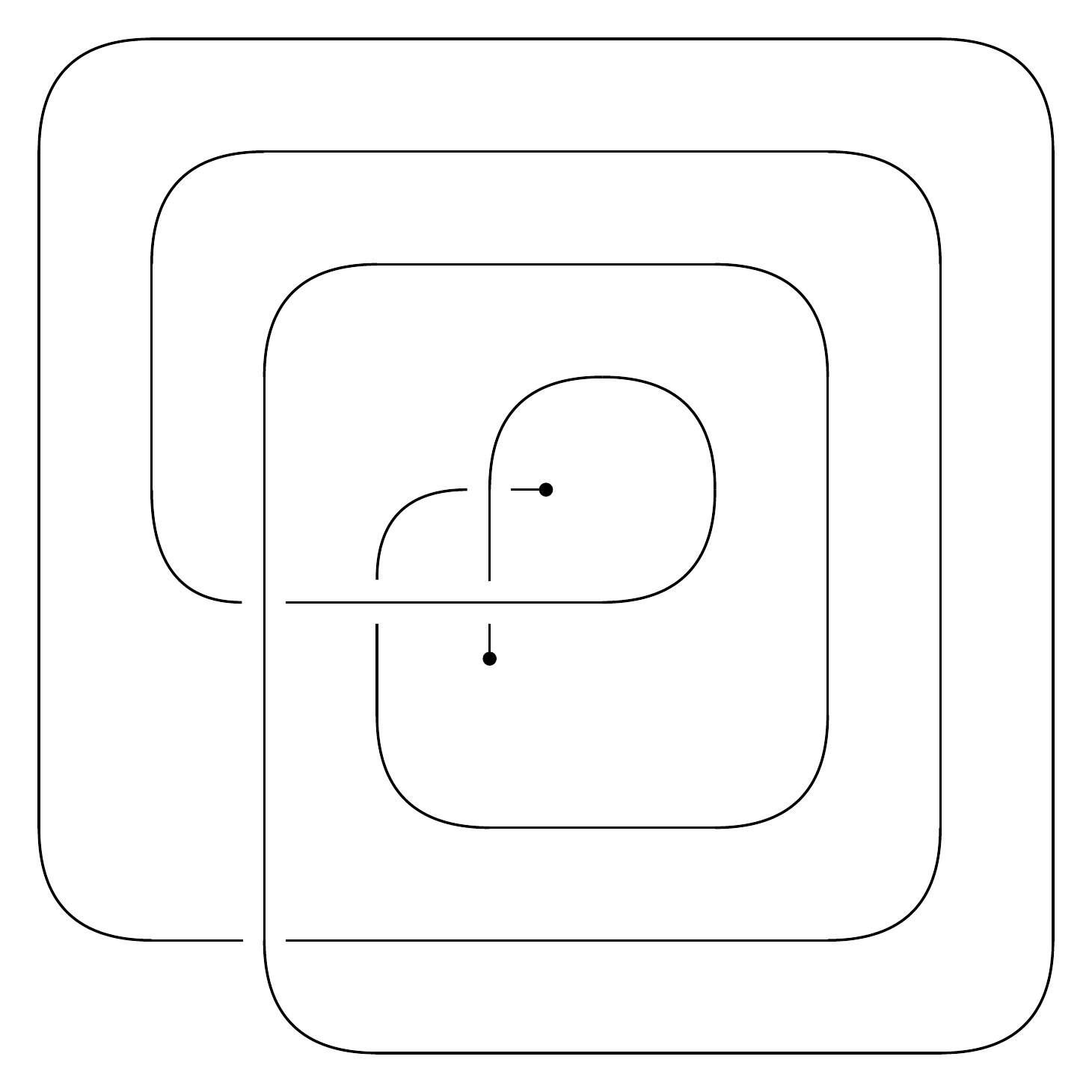}\\
\textcolor{black}{$5_{432}$}
\vspace{1cm}
\end{minipage}
\begin{minipage}[t]{.25\linewidth}
\centering
\includegraphics[width=0.9\textwidth,height=3.5cm,keepaspectratio]{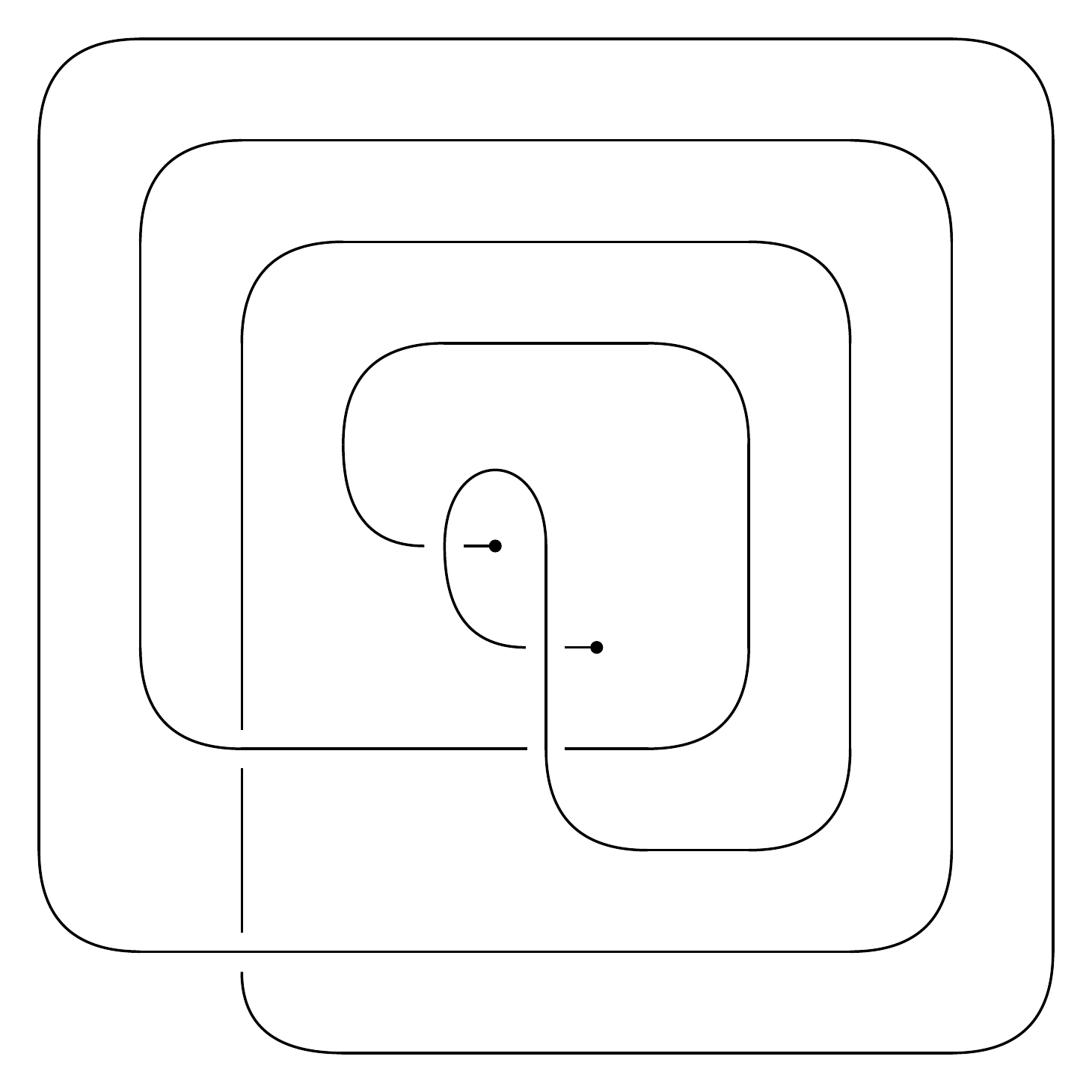}\\
\textcolor{black}{$5_{433}$}
\vspace{1cm}
\end{minipage}
\begin{minipage}[t]{.25\linewidth}
\centering
\includegraphics[width=0.9\textwidth,height=3.5cm,keepaspectratio]{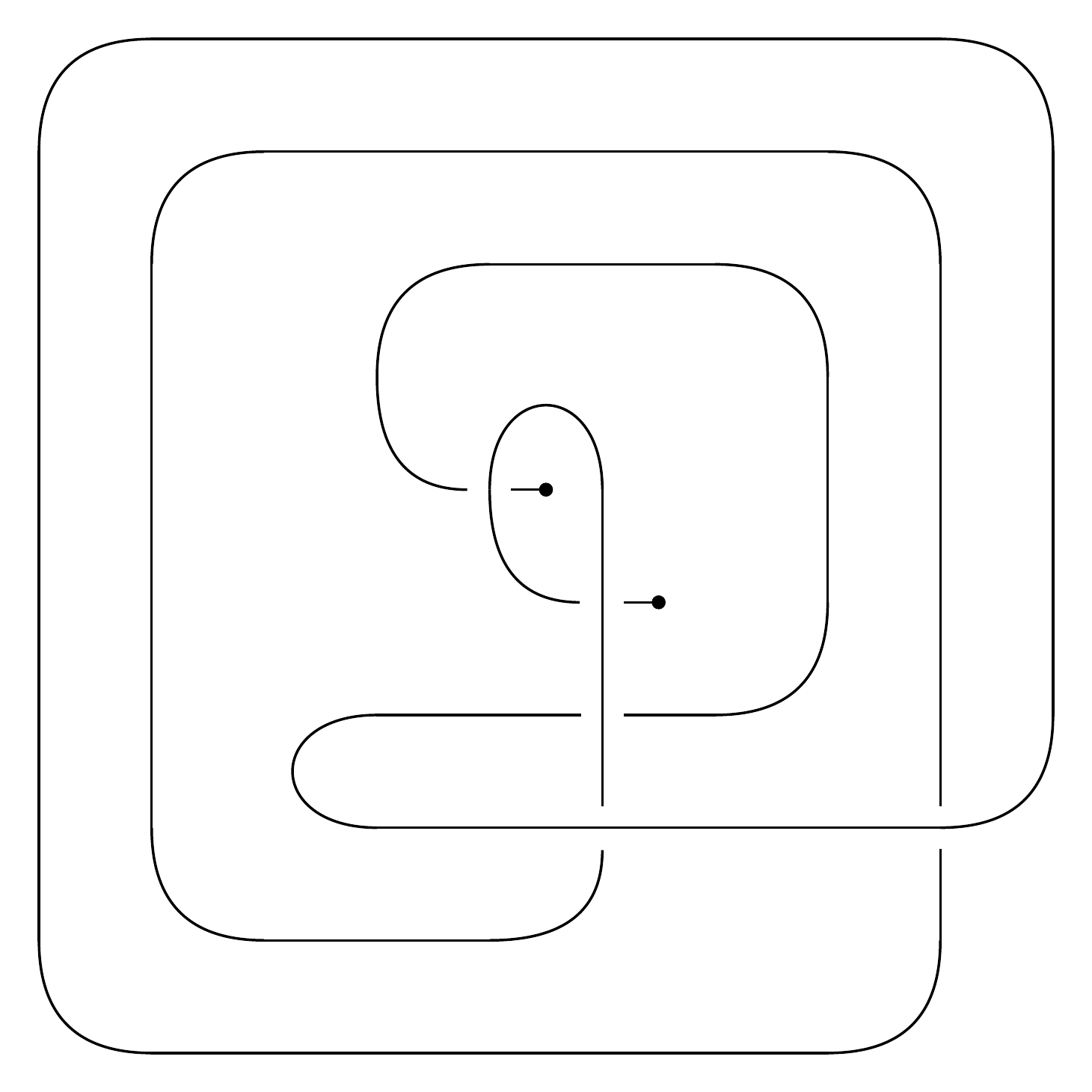}\\
\textcolor{black}{$5_{434}$}
\vspace{1cm}
\end{minipage}
\begin{minipage}[t]{.25\linewidth}
\centering
\includegraphics[width=0.9\textwidth,height=3.5cm,keepaspectratio]{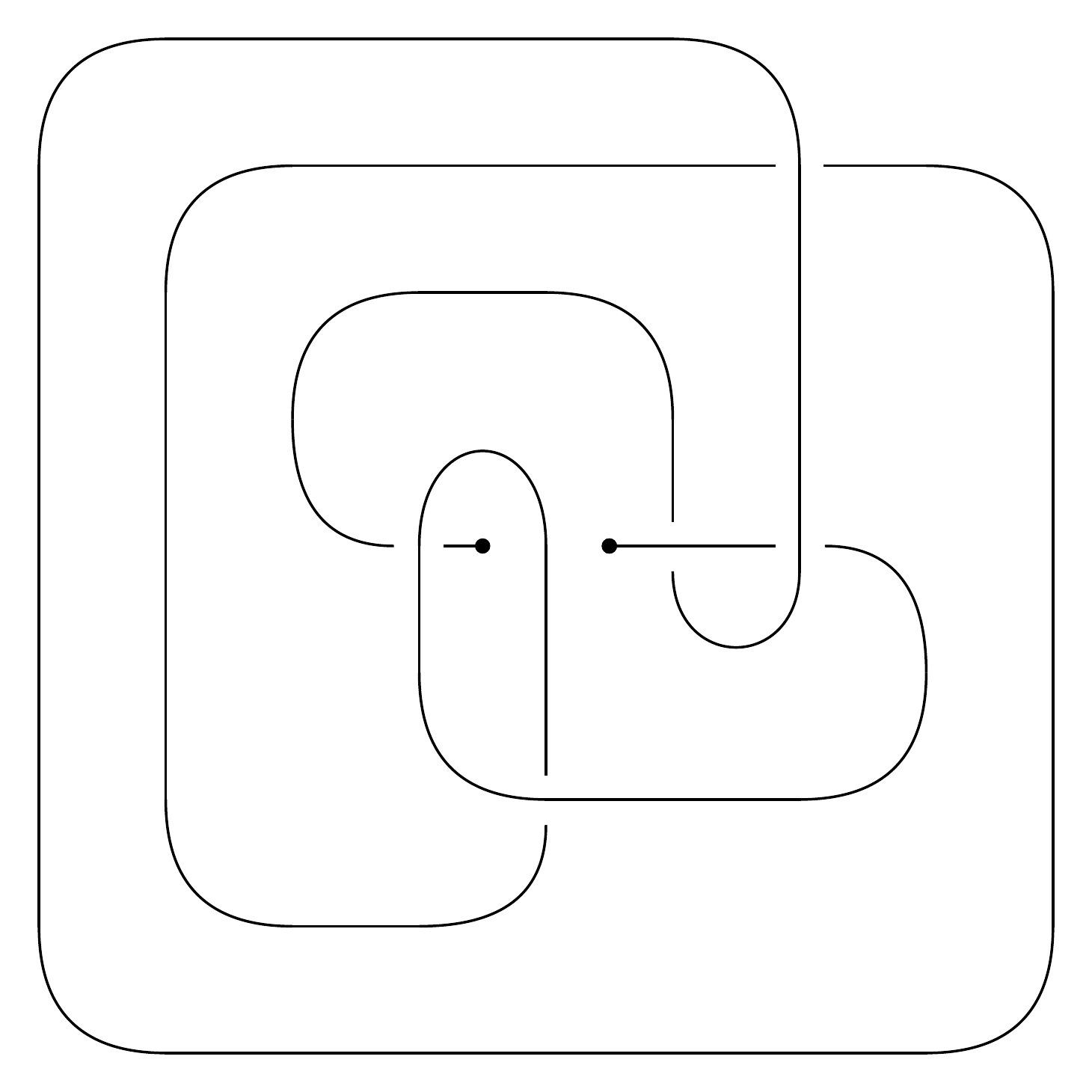}\\
\textcolor{black}{$5_{435}$}
\vspace{1cm}
\end{minipage}
\begin{minipage}[t]{.25\linewidth}
\centering
\includegraphics[width=0.9\textwidth,height=3.5cm,keepaspectratio]{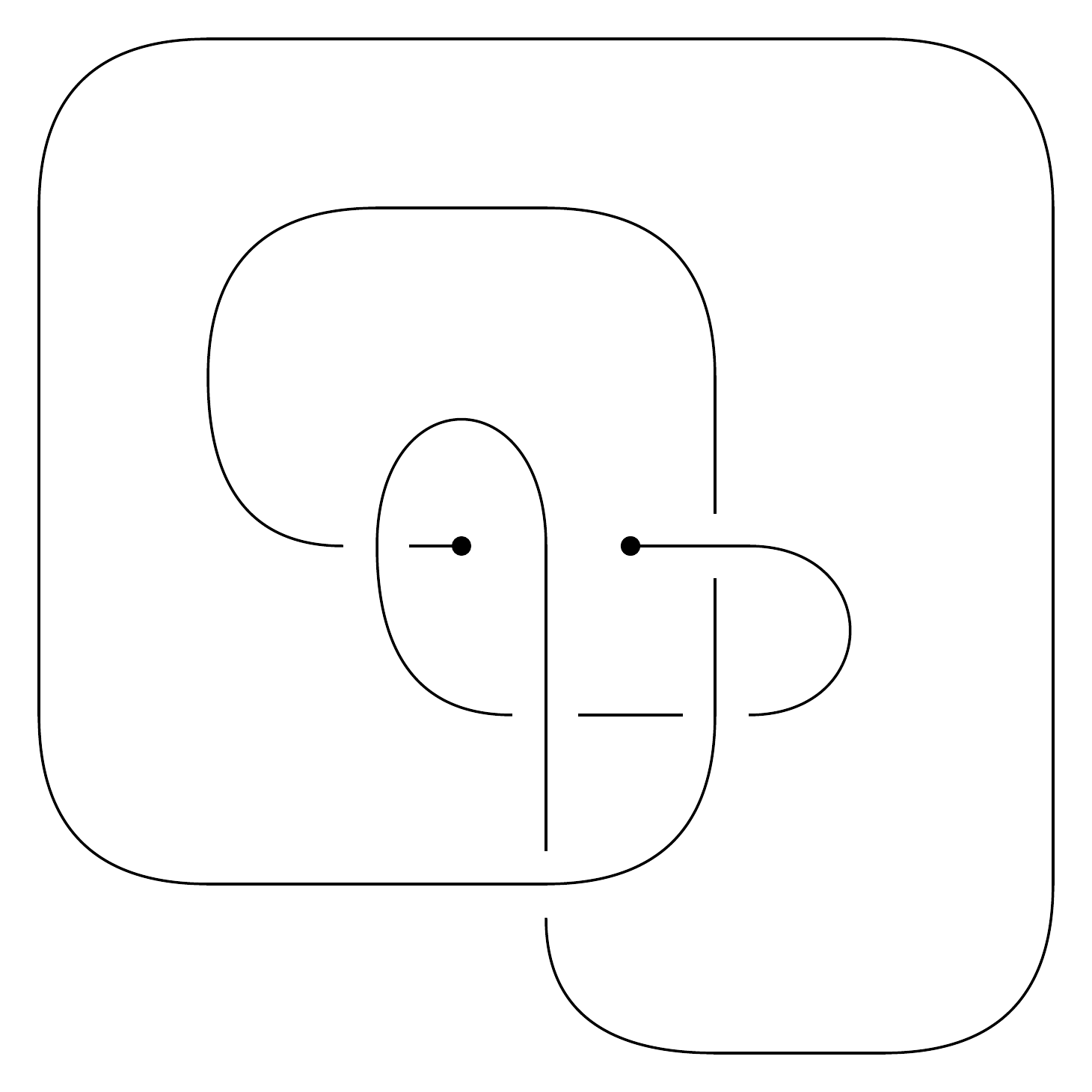}\\
\textcolor{black}{$5_{436}$}
\vspace{1cm}
\end{minipage}
\begin{minipage}[t]{.25\linewidth}
\centering
\includegraphics[width=0.9\textwidth,height=3.5cm,keepaspectratio]{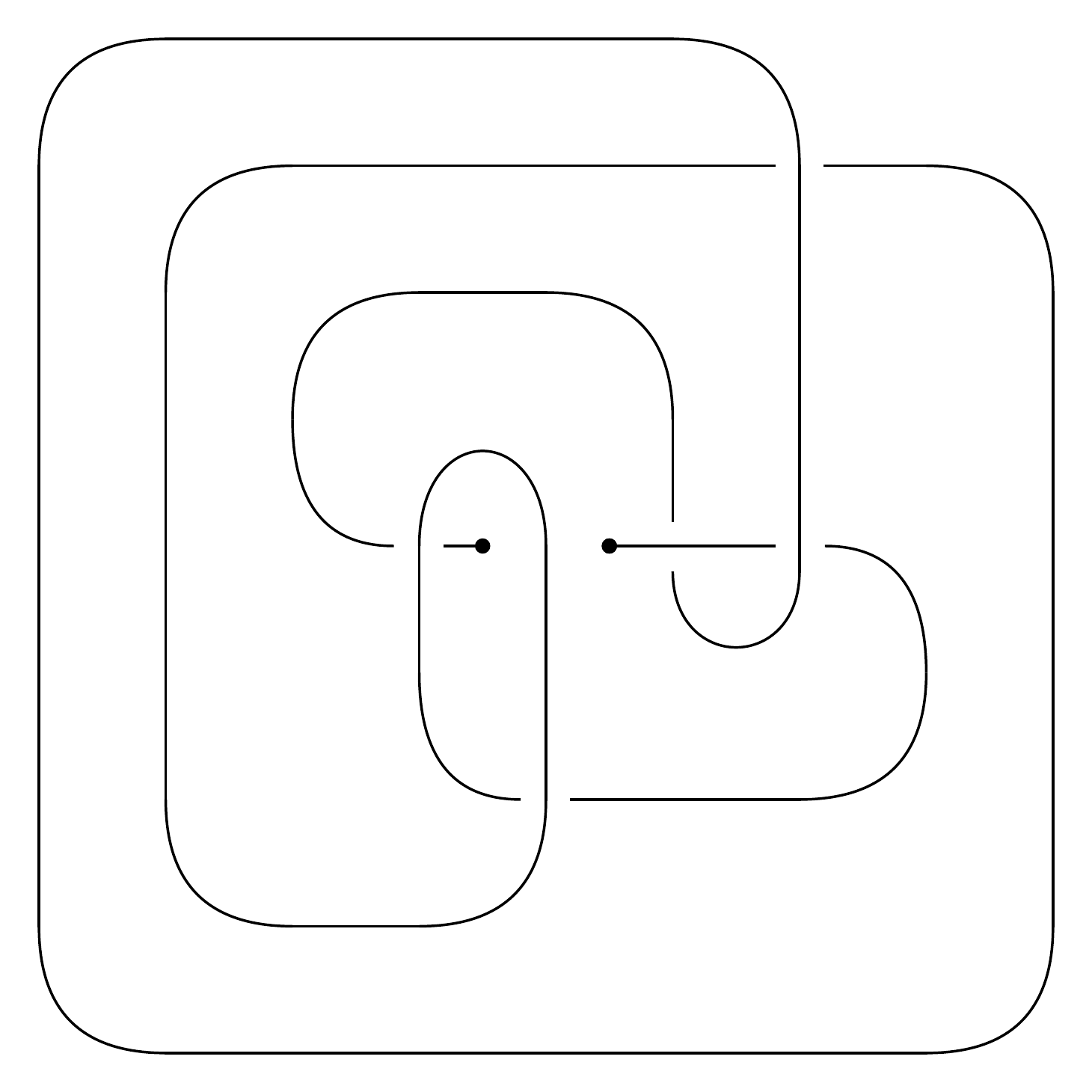}\\
\textcolor{black}{$5_{437}$}
\vspace{1cm}
\end{minipage}
\begin{minipage}[t]{.25\linewidth}
\centering
\includegraphics[width=0.9\textwidth,height=3.5cm,keepaspectratio]{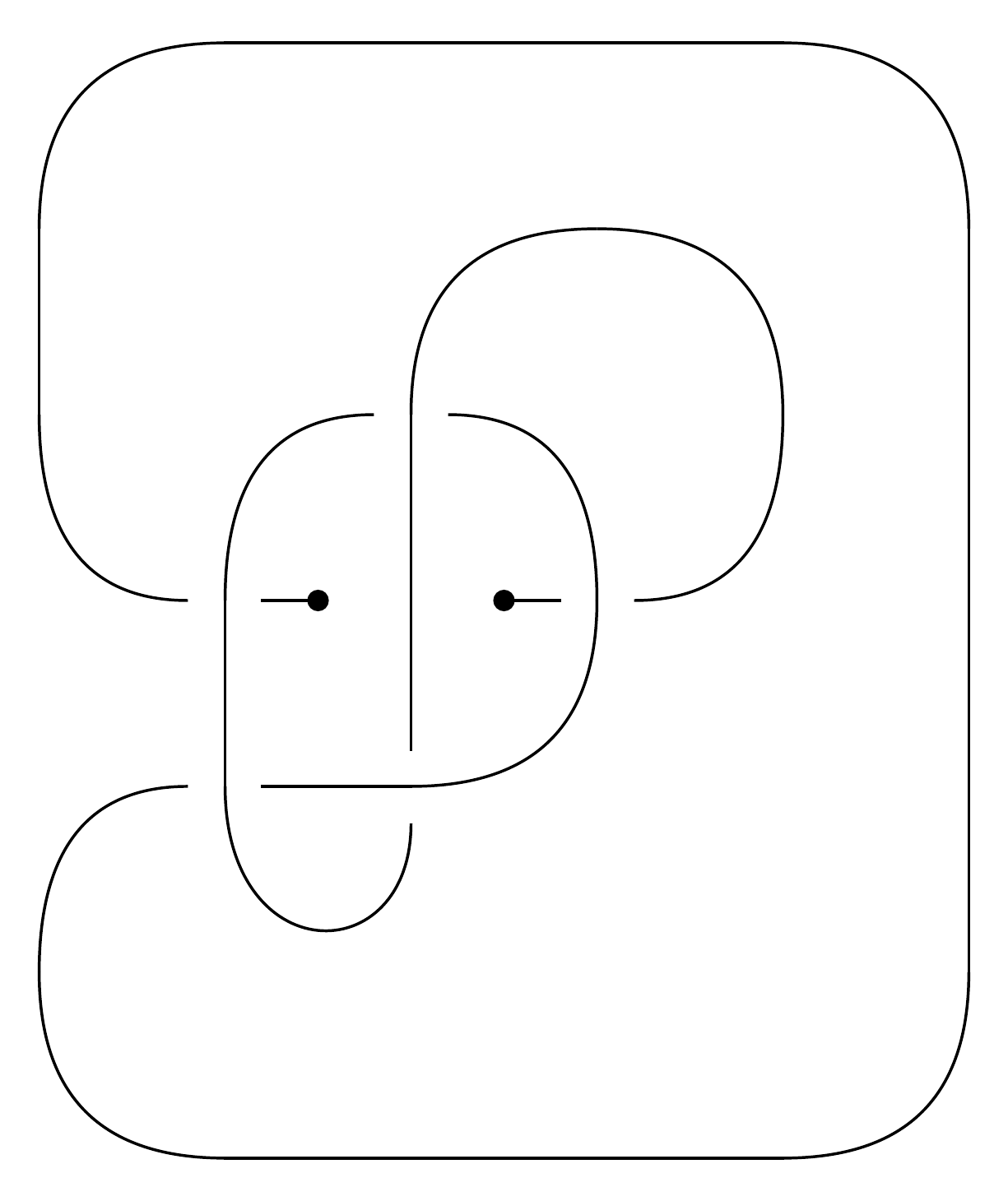}\\
\textcolor{black}{$5_{438}$}
\vspace{1cm}
\end{minipage}
\begin{minipage}[t]{.25\linewidth}
\centering
\includegraphics[width=0.9\textwidth,height=3.5cm,keepaspectratio]{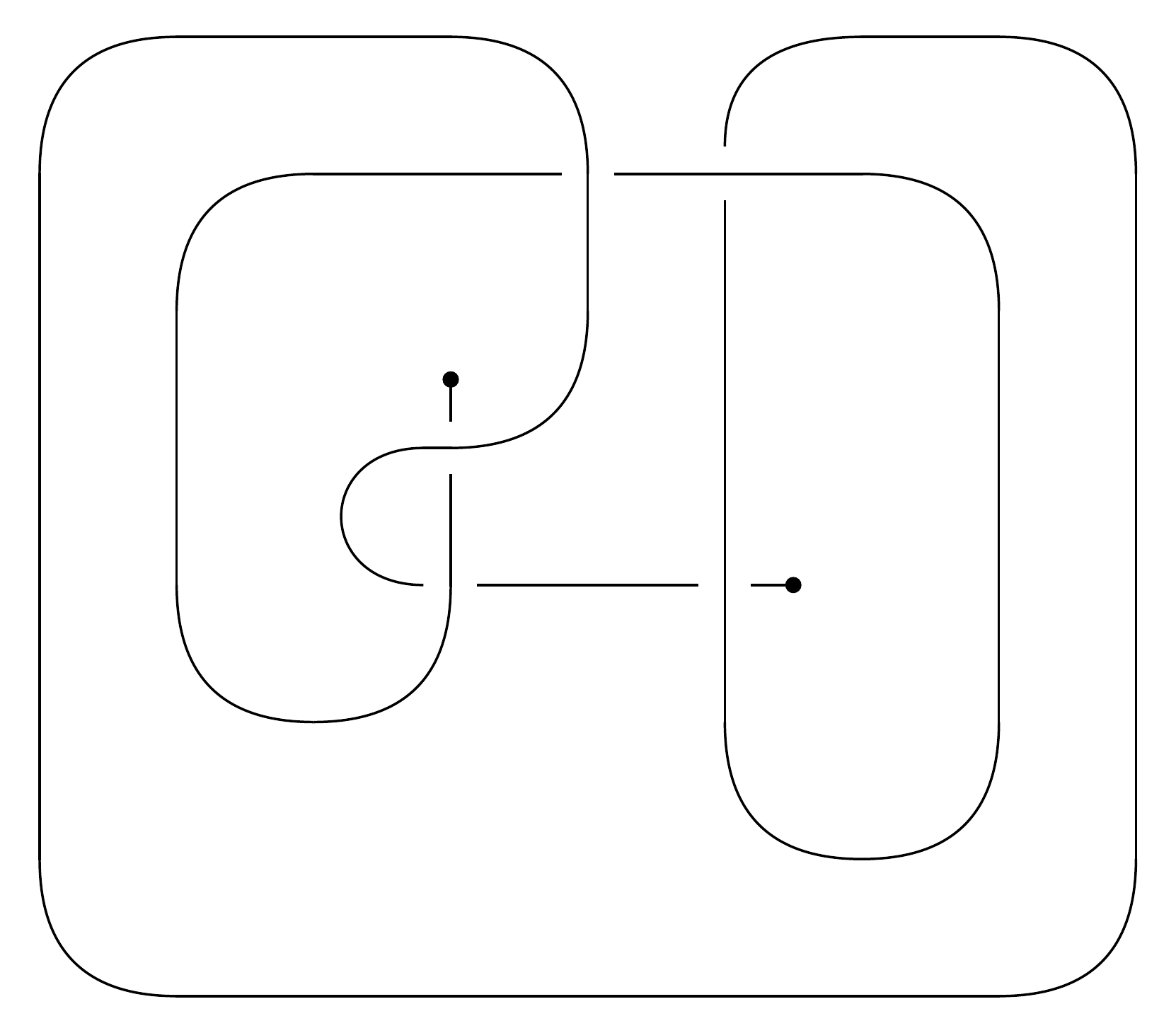}\\
\textcolor{black}{$5_{439}$}
\vspace{1cm}
\end{minipage}
\begin{minipage}[t]{.25\linewidth}
\centering
\includegraphics[width=0.9\textwidth,height=3.5cm,keepaspectratio]{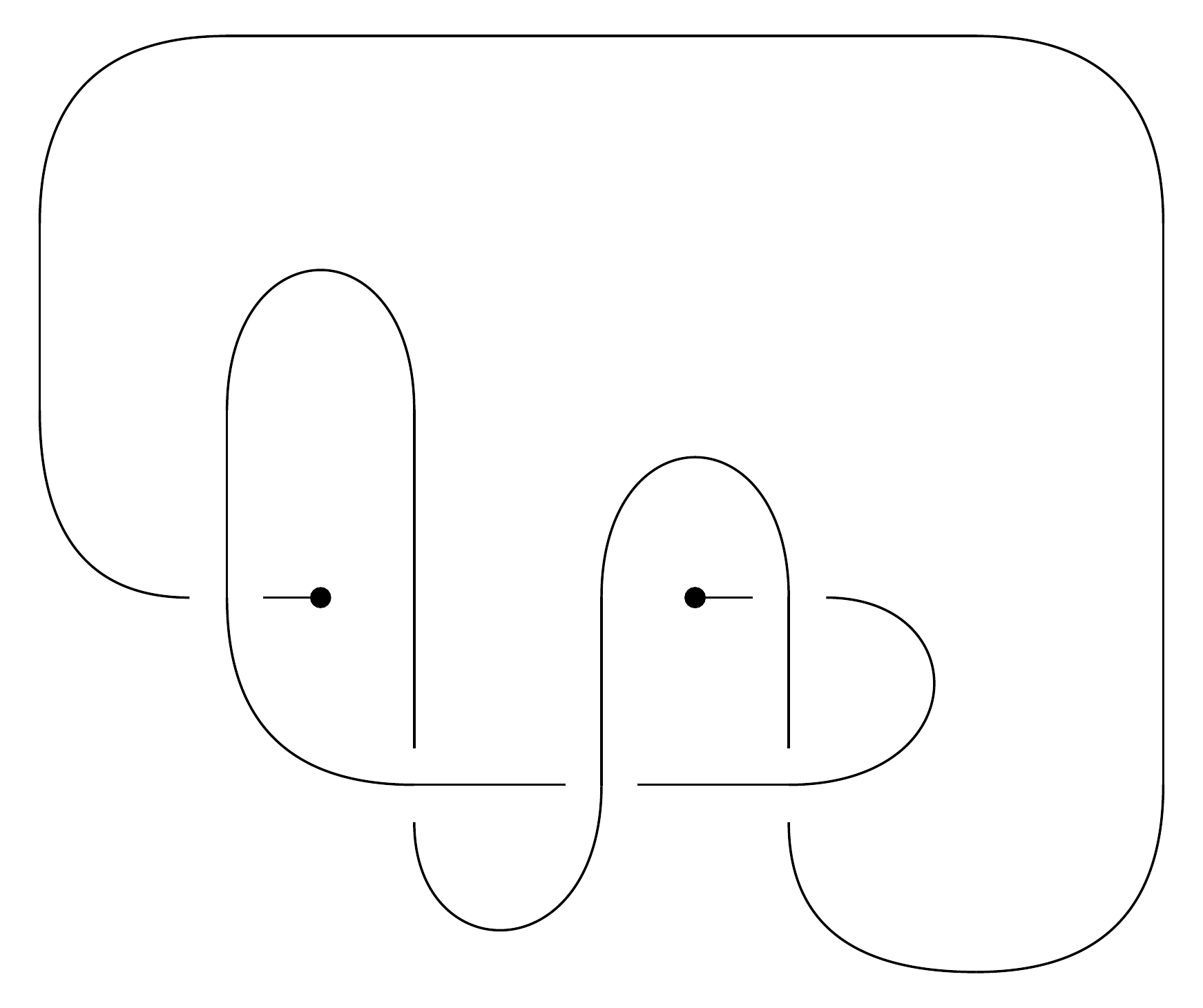}\\
\textcolor{black}{$5_{440}$}
\vspace{1cm}
\end{minipage}
\begin{minipage}[t]{.25\linewidth}
\centering
\includegraphics[width=0.9\textwidth,height=3.5cm,keepaspectratio]{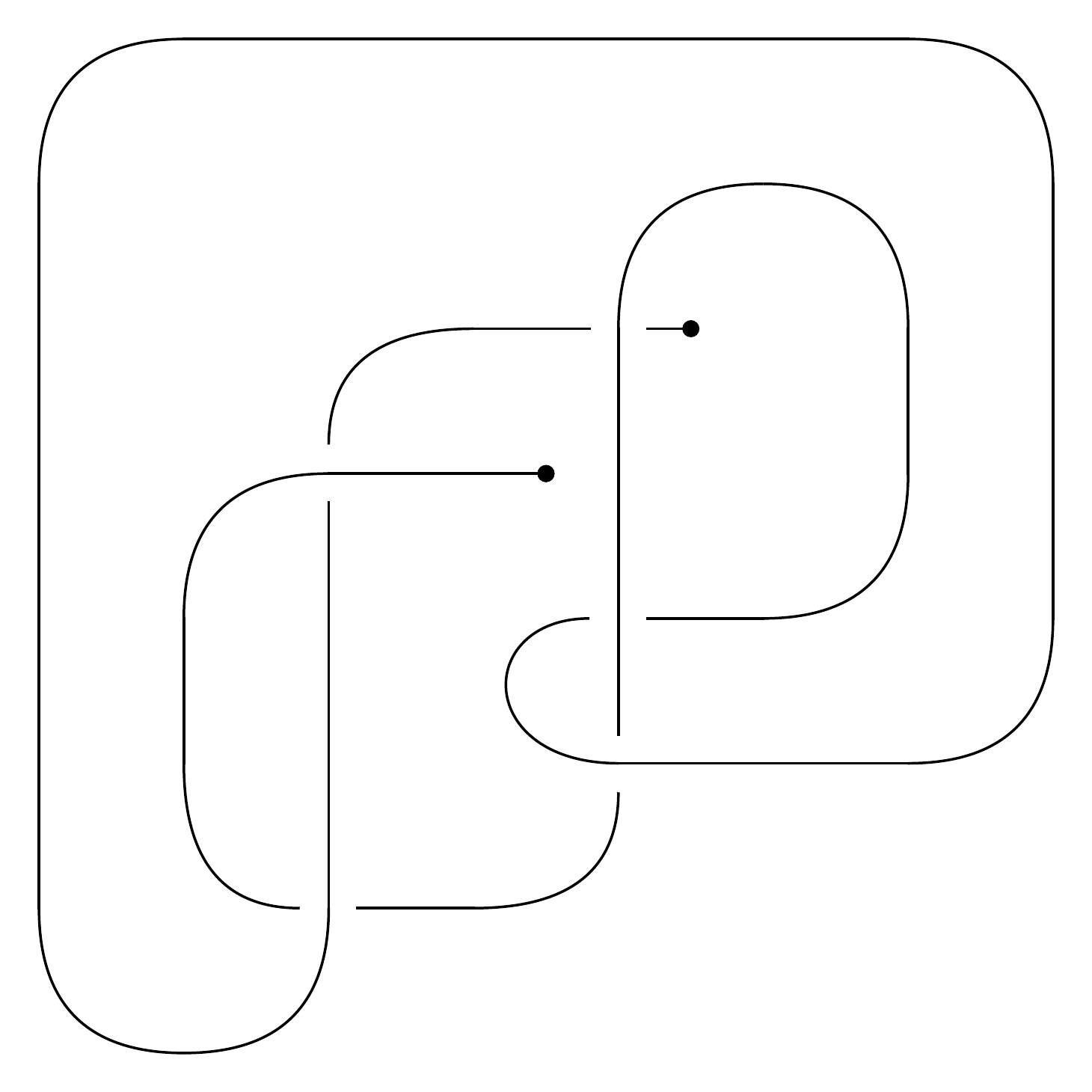}\\
\textcolor{black}{$5_{441}$}
\vspace{1cm}
\end{minipage}
\begin{minipage}[t]{.25\linewidth}
\centering
\includegraphics[width=0.9\textwidth,height=3.5cm,keepaspectratio]{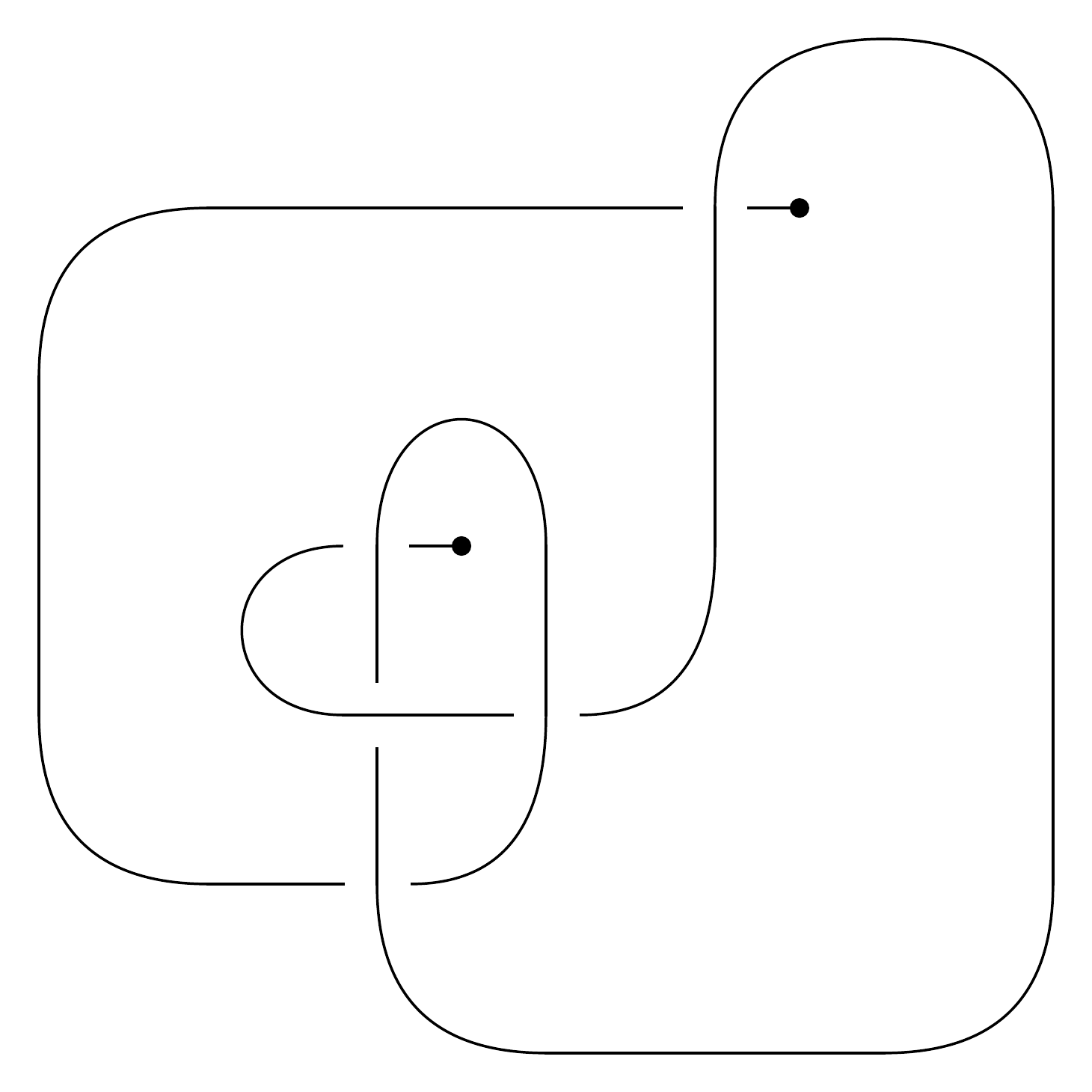}\\
\textcolor{black}{$5_{442}$}
\vspace{1cm}
\end{minipage}
\begin{minipage}[t]{.25\linewidth}
\centering
\includegraphics[width=0.9\textwidth,height=3.5cm,keepaspectratio]{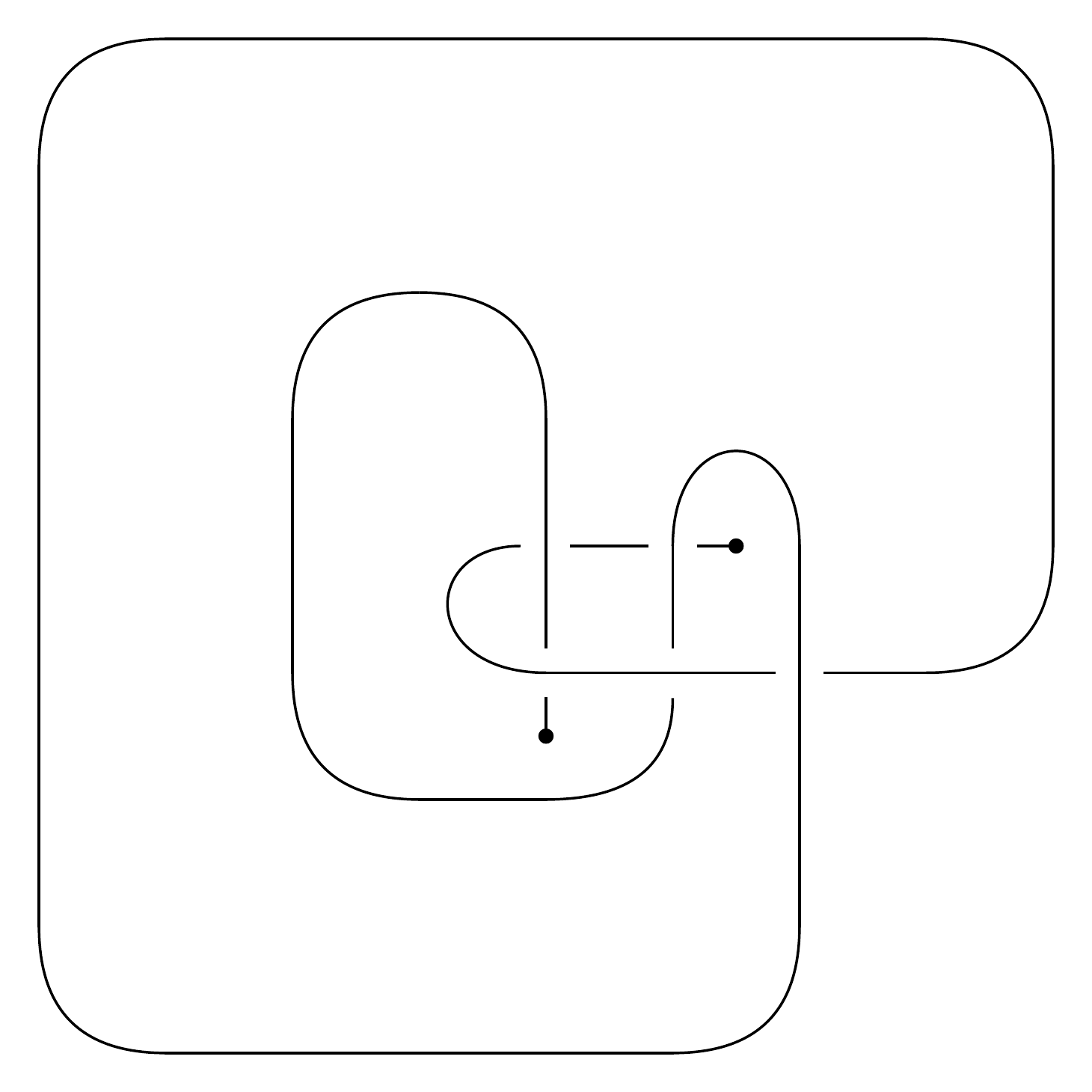}\\
\textcolor{black}{$5_{443}$}
\vspace{1cm}
\end{minipage}
\begin{minipage}[t]{.25\linewidth}
\centering
\includegraphics[width=0.9\textwidth,height=3.5cm,keepaspectratio]{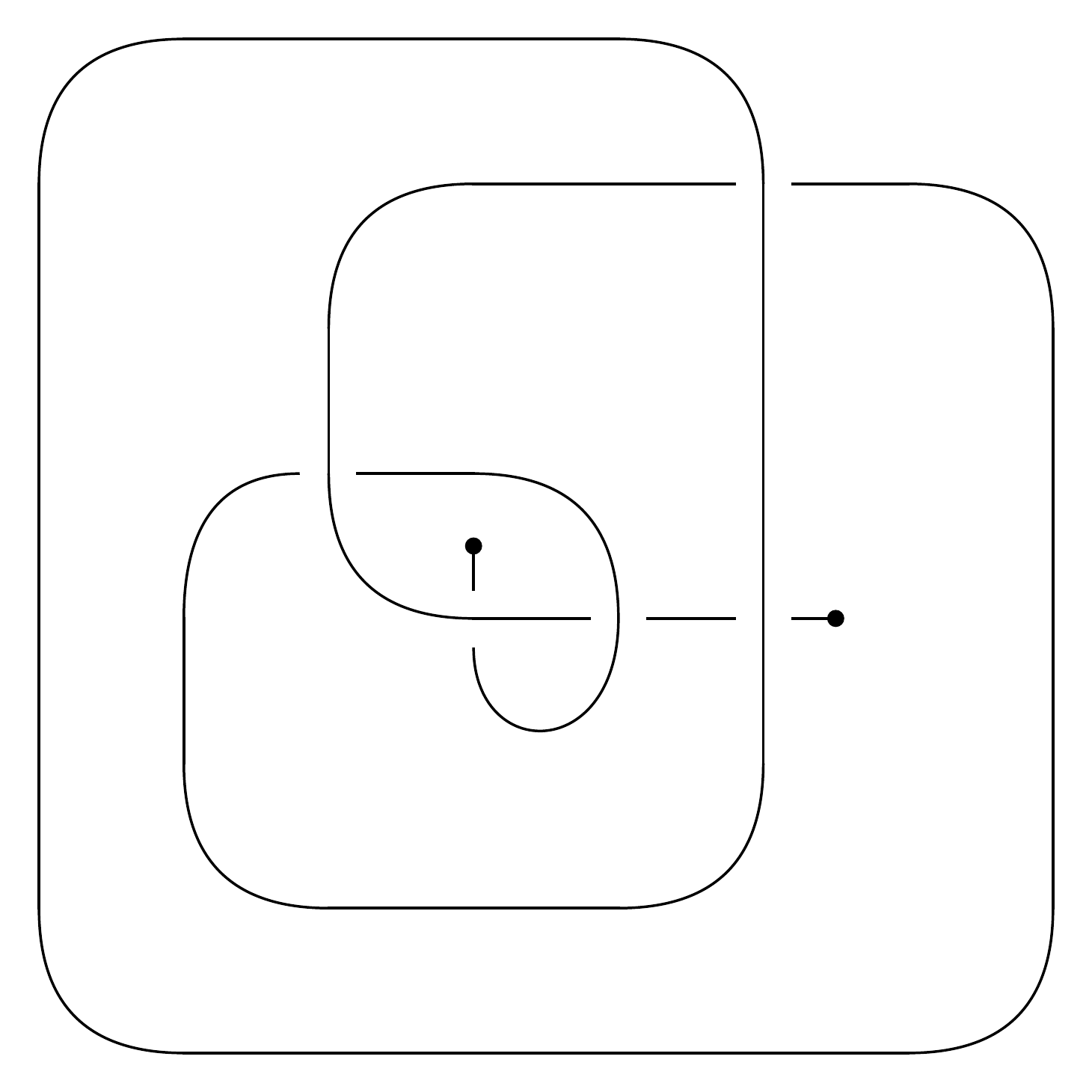}\\
\textcolor{black}{$5_{444}$}
\vspace{1cm}
\end{minipage}
\begin{minipage}[t]{.25\linewidth}
\centering
\includegraphics[width=0.9\textwidth,height=3.5cm,keepaspectratio]{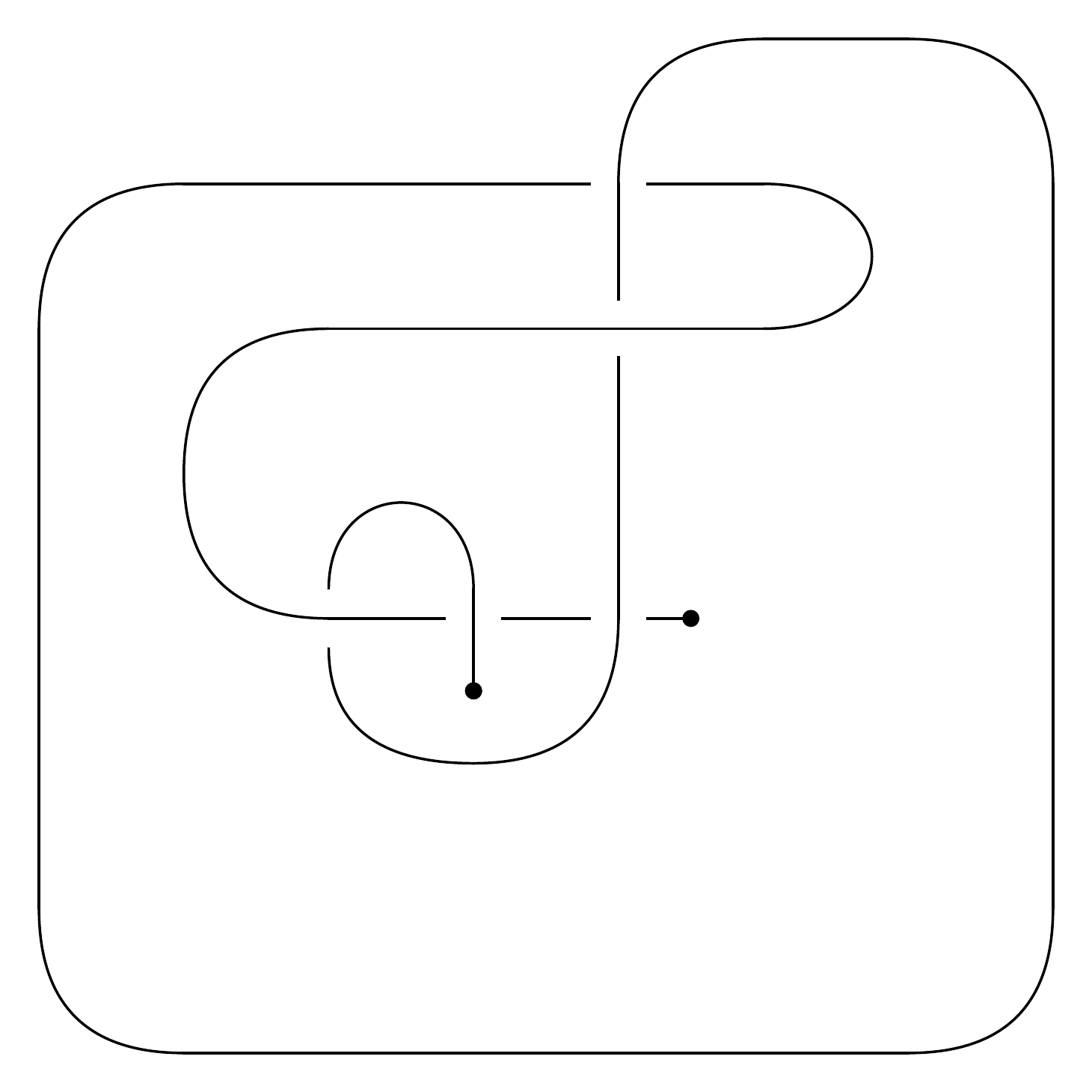}\\
\textcolor{black}{$5_{445}$}
\vspace{1cm}
\end{minipage}
\begin{minipage}[t]{.25\linewidth}
\centering
\includegraphics[width=0.9\textwidth,height=3.5cm,keepaspectratio]{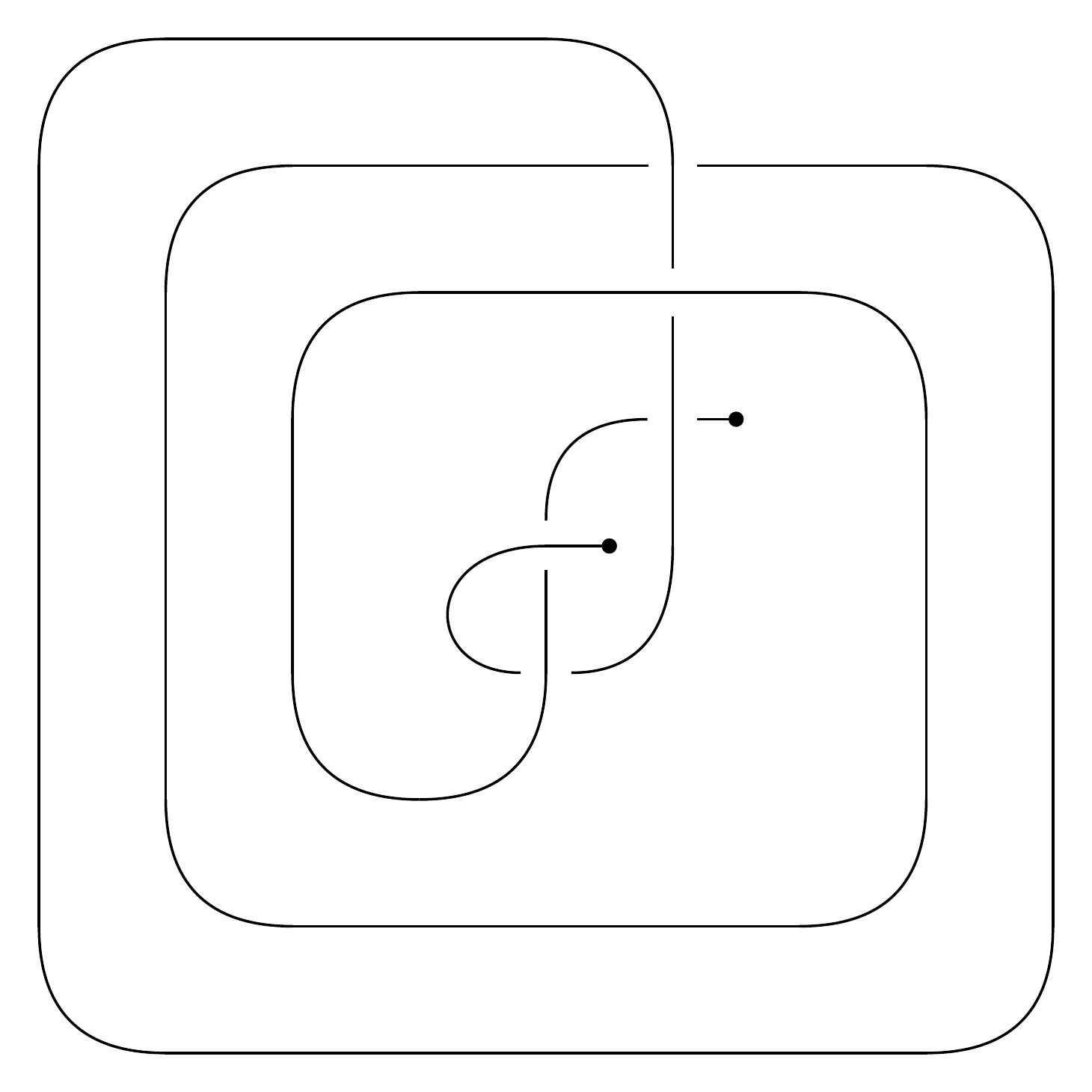}\\
\textcolor{black}{$5_{446}$}
\vspace{1cm}
\end{minipage}
\begin{minipage}[t]{.25\linewidth}
\centering
\includegraphics[width=0.9\textwidth,height=3.5cm,keepaspectratio]{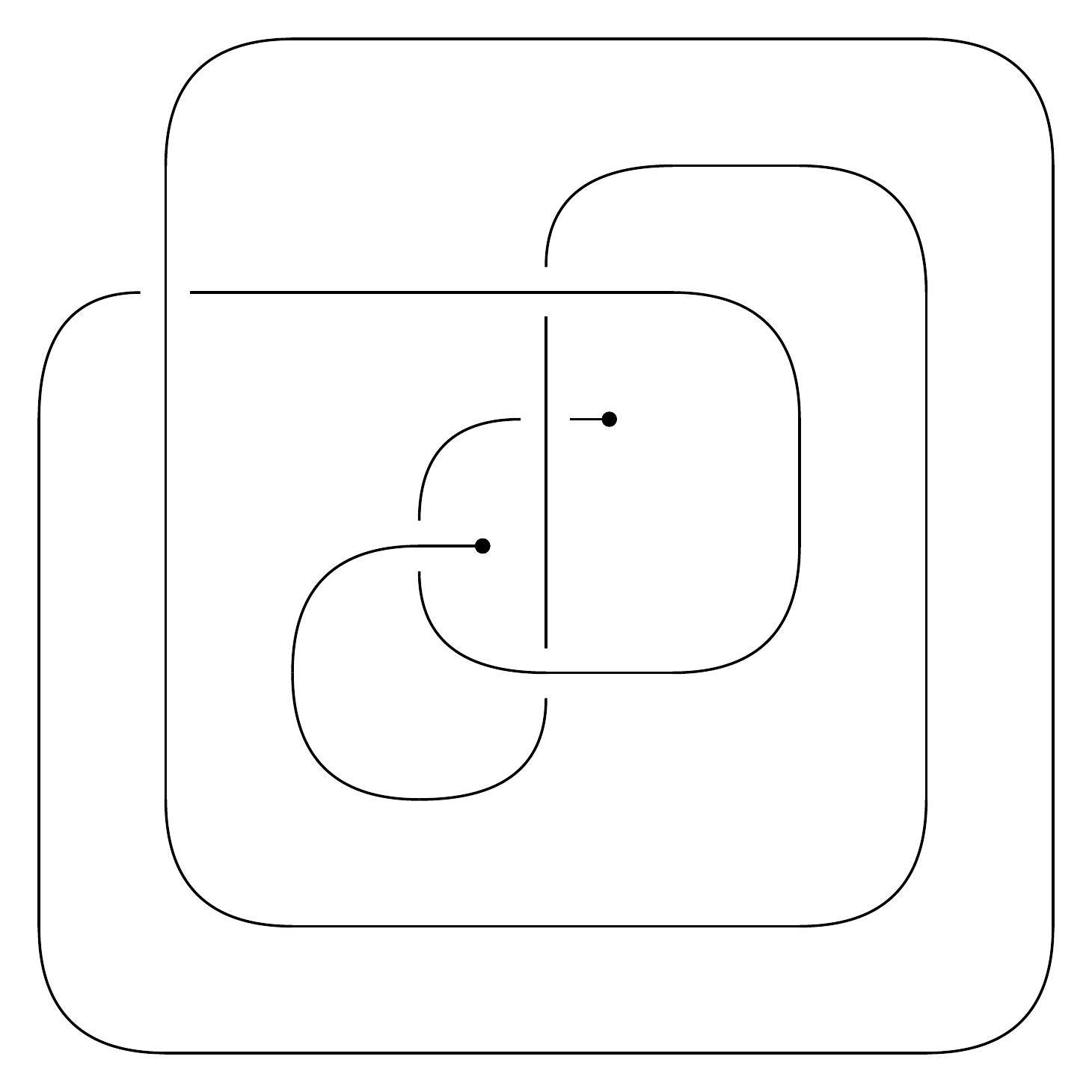}\\
\textcolor{black}{$5_{447}$}
\vspace{1cm}
\end{minipage}
\begin{minipage}[t]{.25\linewidth}
\centering
\includegraphics[width=0.9\textwidth,height=3.5cm,keepaspectratio]{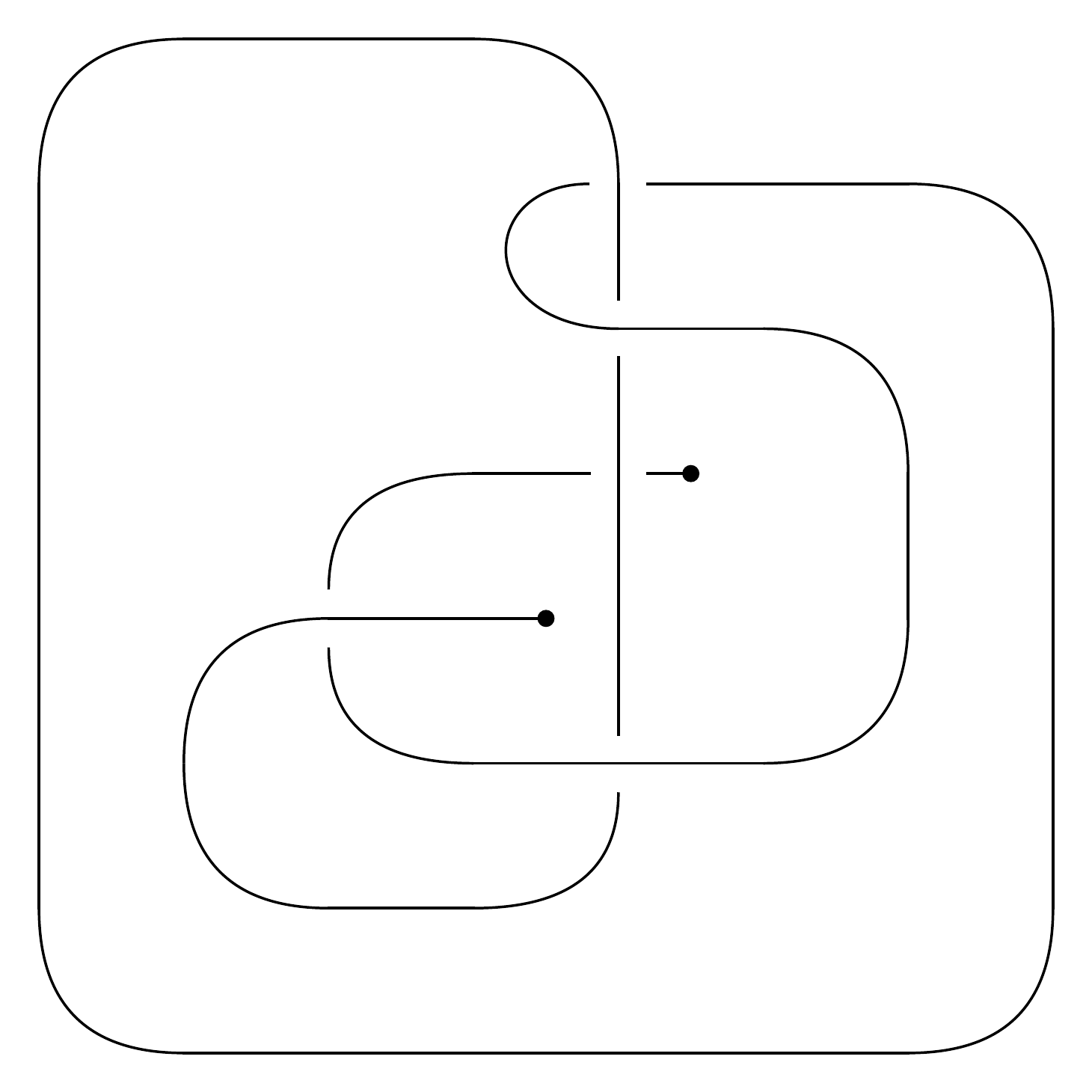}\\
\textcolor{black}{$5_{448}$}
\vspace{1cm}
\end{minipage}
\begin{minipage}[t]{.25\linewidth}
\centering
\includegraphics[width=0.9\textwidth,height=3.5cm,keepaspectratio]{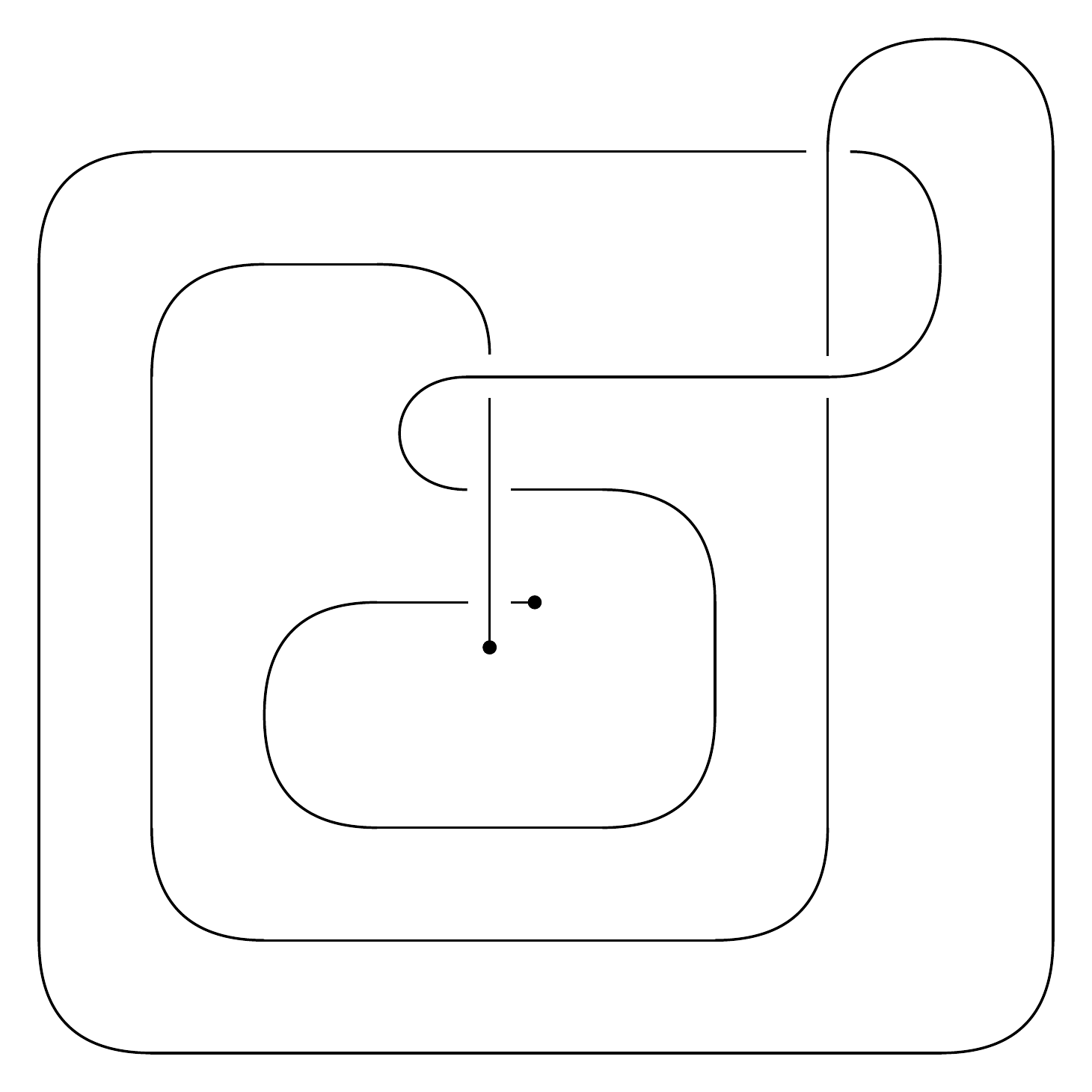}\\
\textcolor{black}{$5_{449}$}
\vspace{1cm}
\end{minipage}
\begin{minipage}[t]{.25\linewidth}
\centering
\includegraphics[width=0.9\textwidth,height=3.5cm,keepaspectratio]{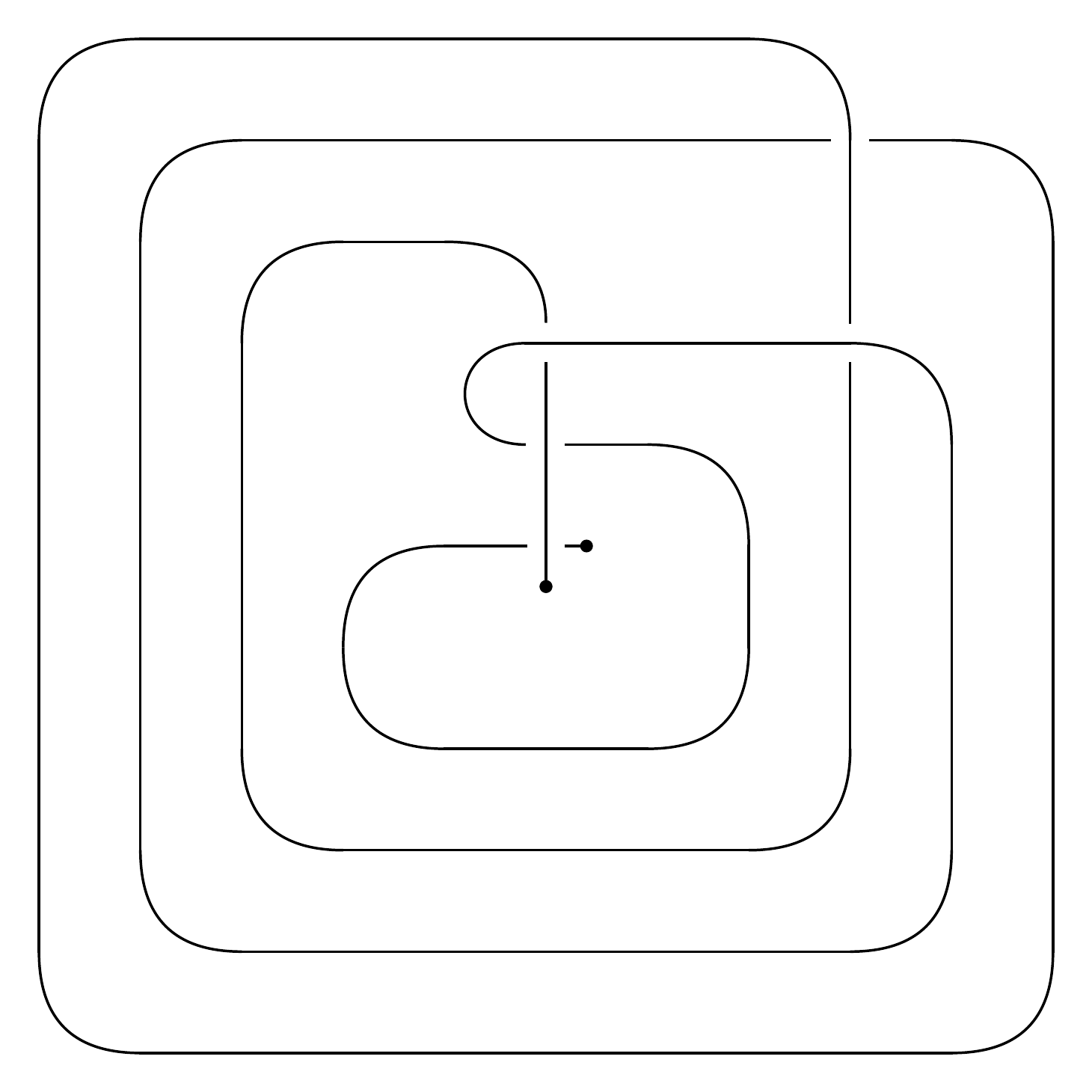}\\
\textcolor{black}{$5_{450}$}
\vspace{1cm}
\end{minipage}
\begin{minipage}[t]{.25\linewidth}
\centering
\includegraphics[width=0.9\textwidth,height=3.5cm,keepaspectratio]{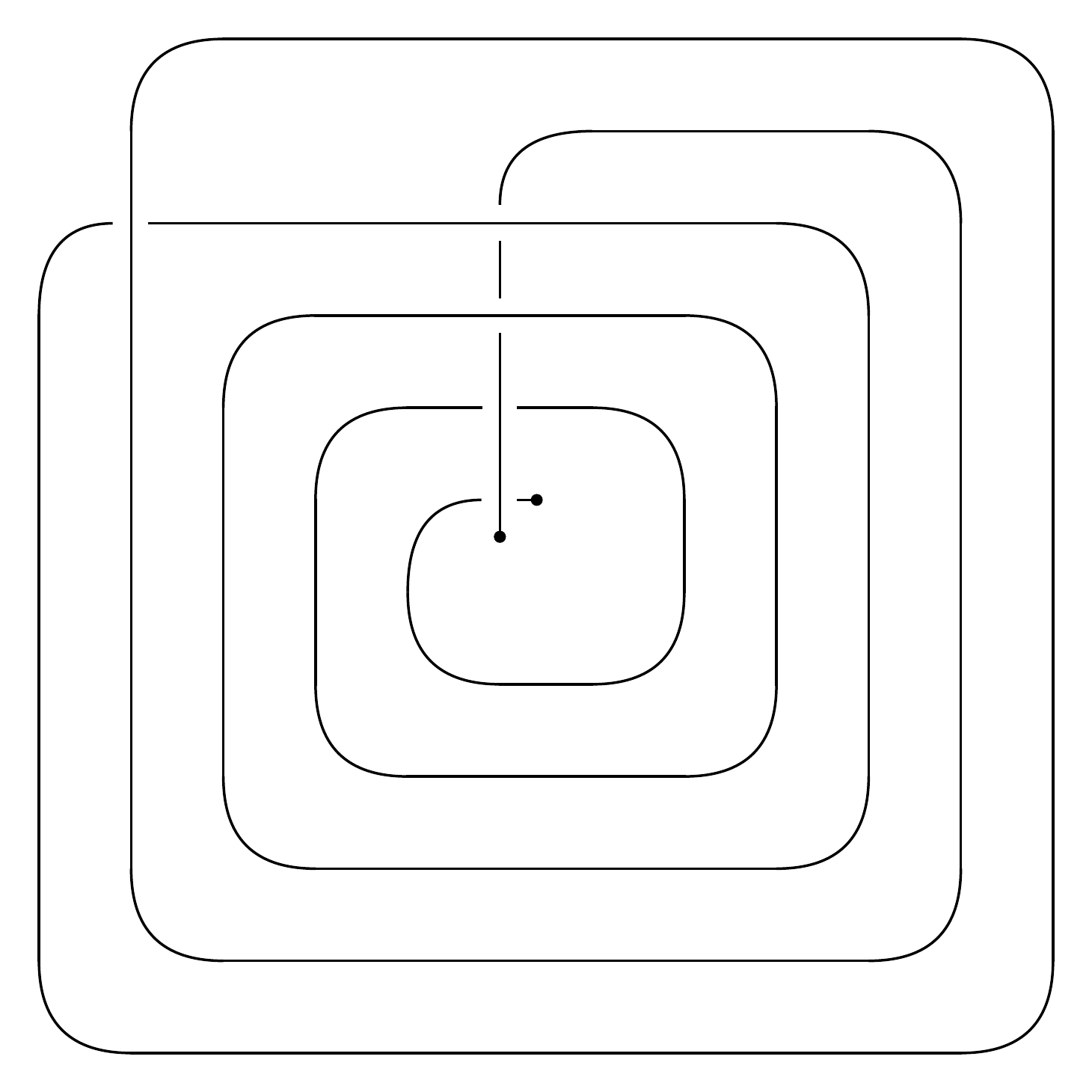}\\
\textcolor{black}{$5_{451}$}
\vspace{1cm}
\end{minipage}
\begin{minipage}[t]{.25\linewidth}
\centering
\includegraphics[width=0.9\textwidth,height=3.5cm,keepaspectratio]{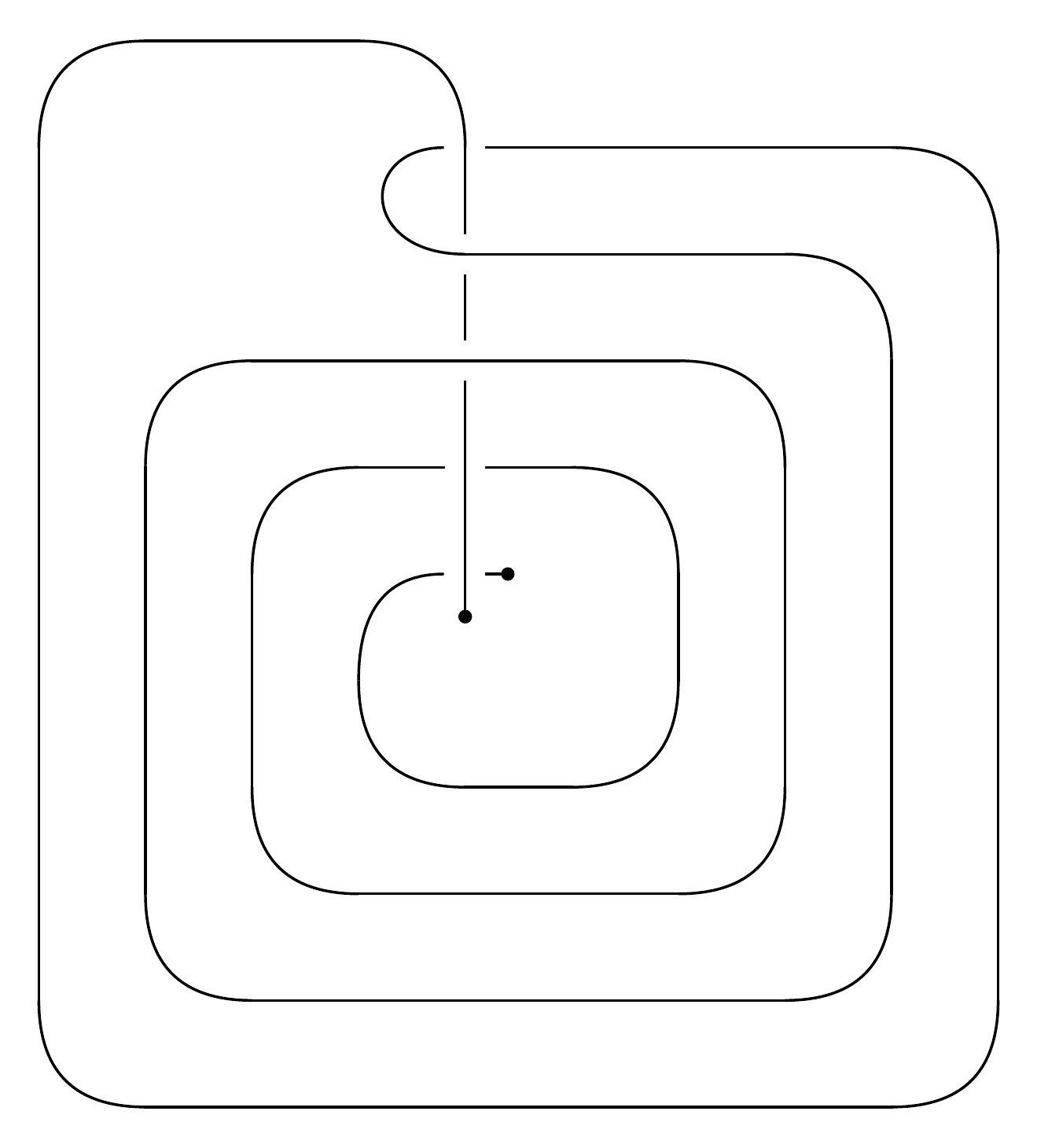}\\
\textcolor{black}{$5_{452}$}
\vspace{1cm}
\end{minipage}
\begin{minipage}[t]{.25\linewidth}
\centering
\includegraphics[width=0.9\textwidth,height=3.5cm,keepaspectratio]{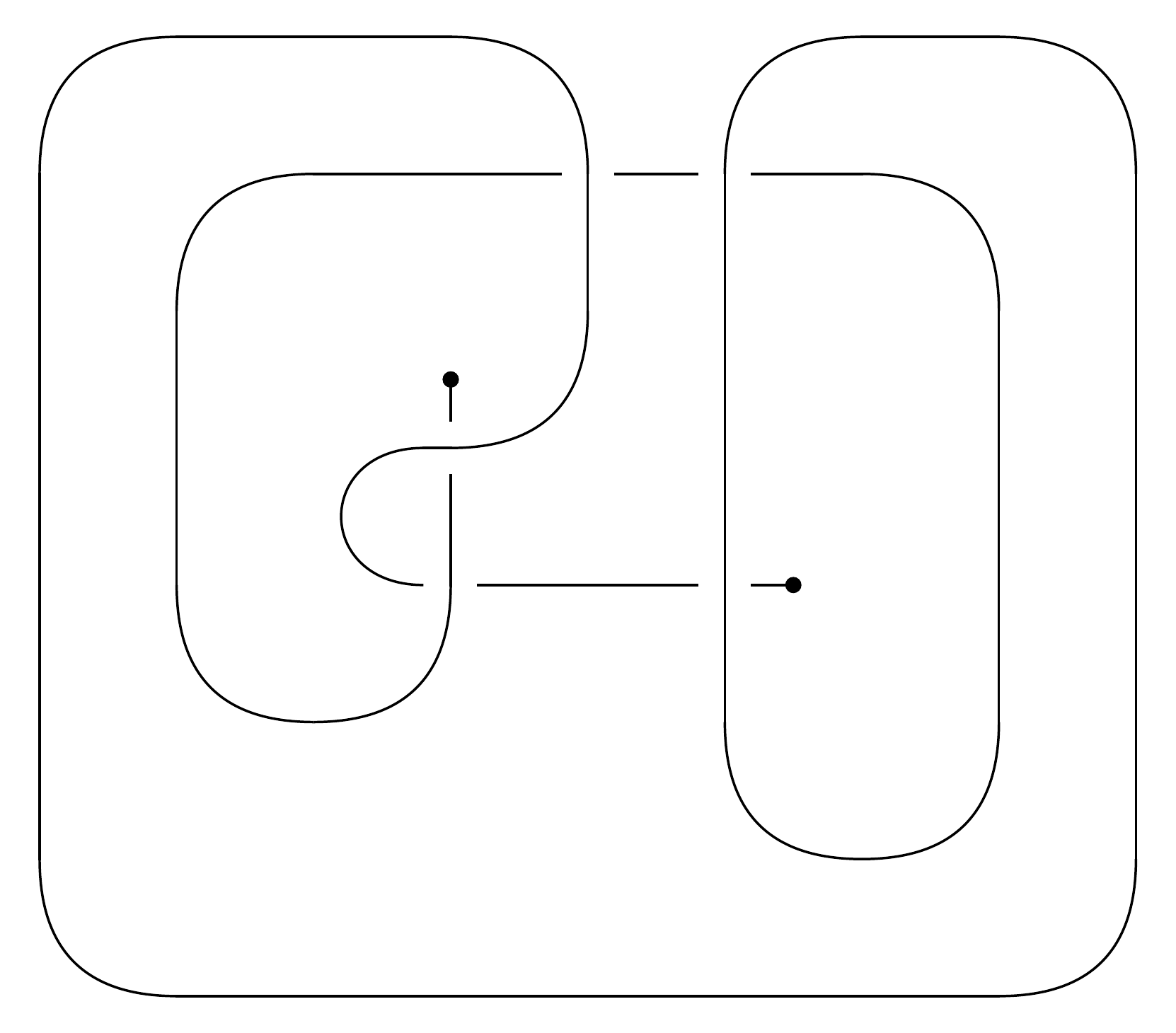}\\
\textcolor{black}{$5_{453}$}
\vspace{1cm}
\end{minipage}
\begin{minipage}[t]{.25\linewidth}
\centering
\includegraphics[width=0.9\textwidth,height=3.5cm,keepaspectratio]{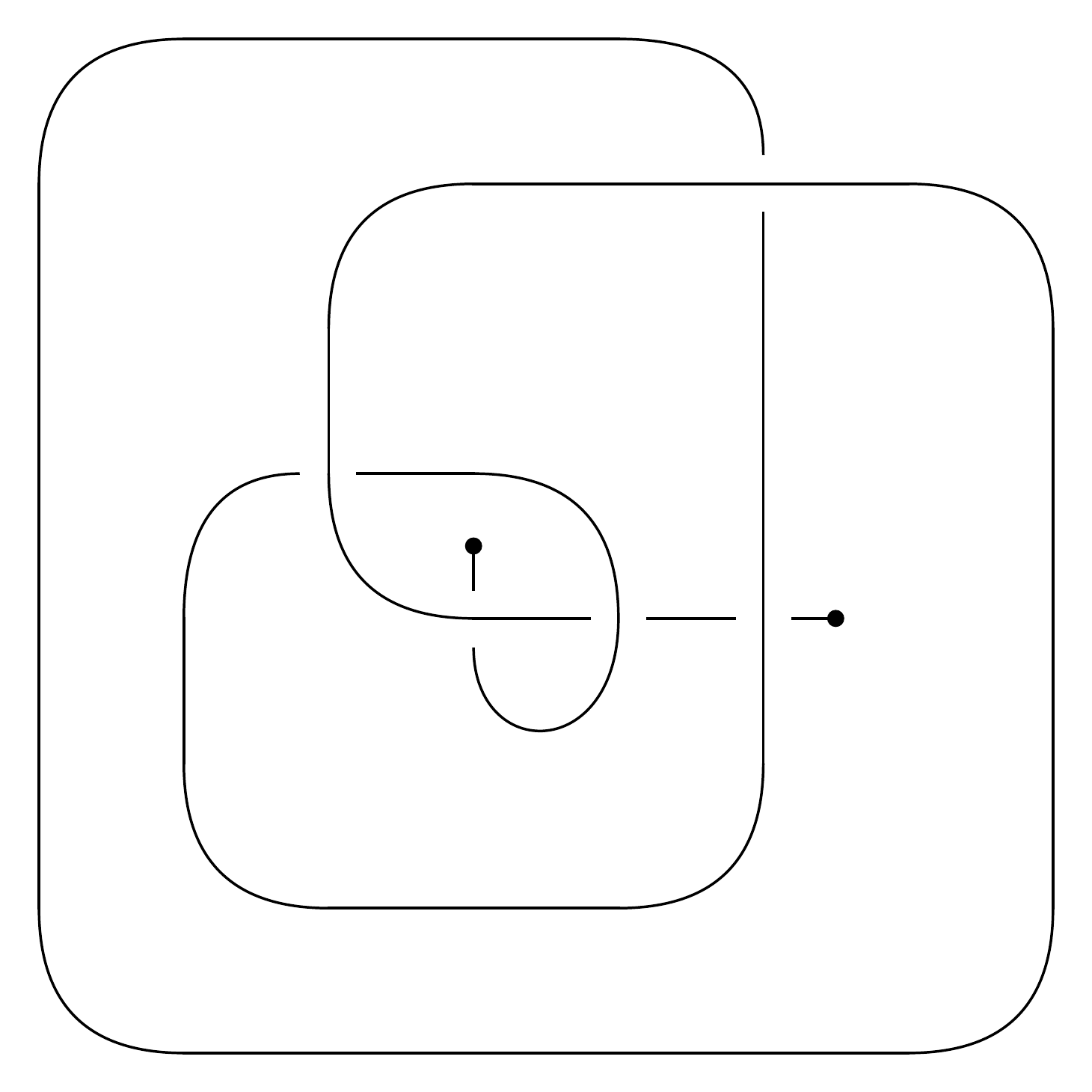}\\
\textcolor{black}{$5_{454}$}
\vspace{1cm}
\end{minipage}
\begin{minipage}[t]{.25\linewidth}
\centering
\includegraphics[width=0.9\textwidth,height=3.5cm,keepaspectratio]{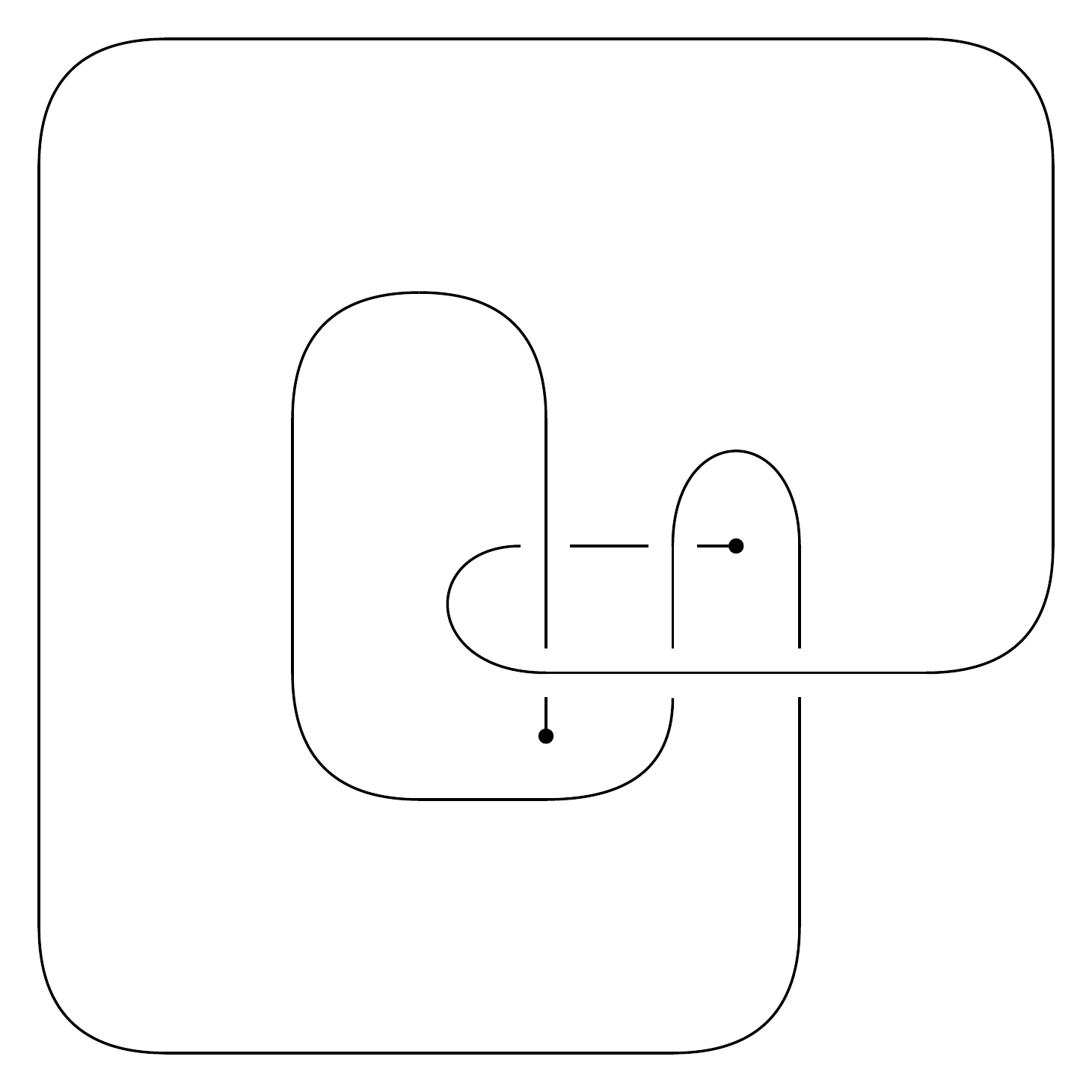}\\
\textcolor{black}{$5_{455}$}
\vspace{1cm}
\end{minipage}
\begin{minipage}[t]{.25\linewidth}
\centering
\includegraphics[width=0.9\textwidth,height=3.5cm,keepaspectratio]{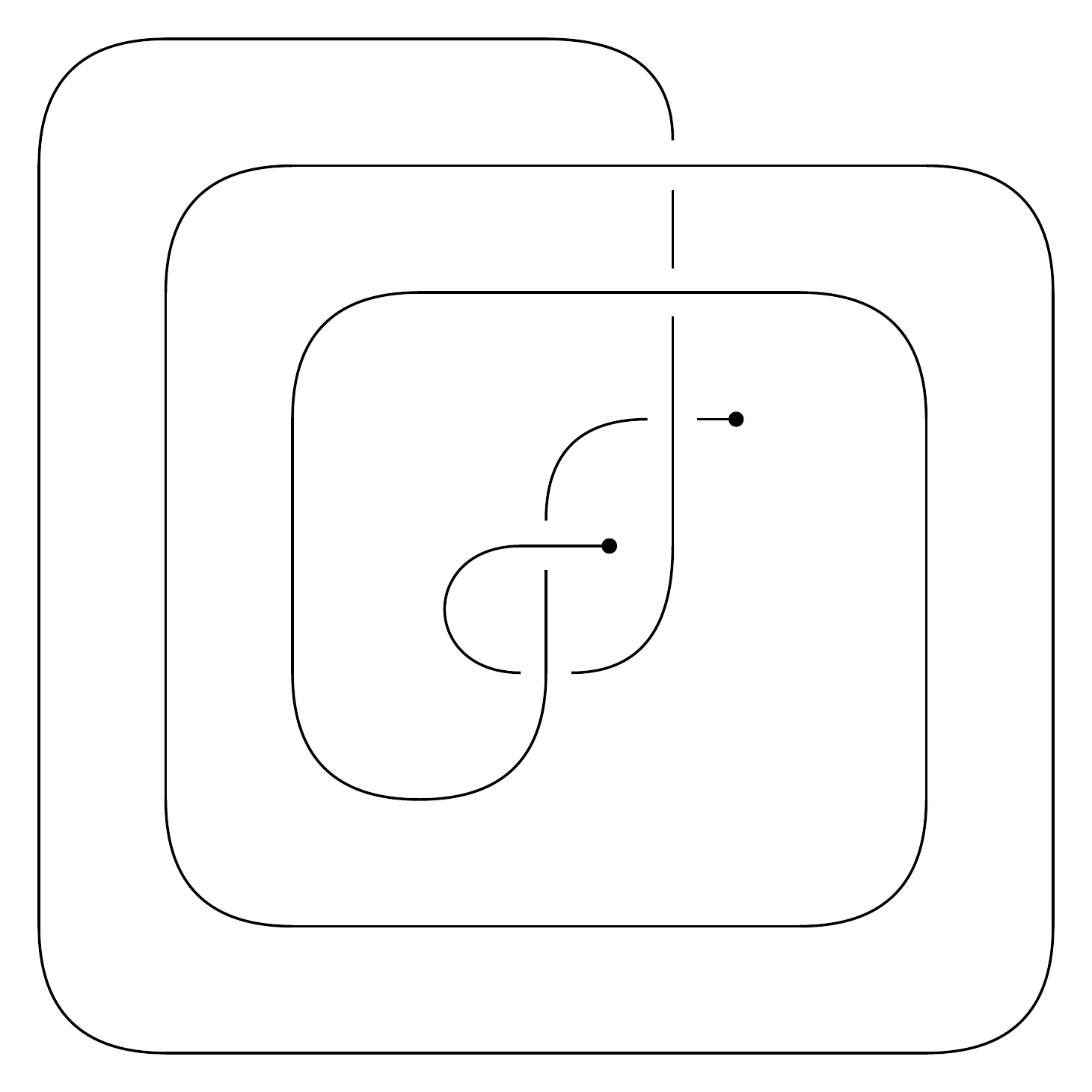}\\
\textcolor{black}{$5_{456}$}
\vspace{1cm}
\end{minipage}
\begin{minipage}[t]{.25\linewidth}
\centering
\includegraphics[width=0.9\textwidth,height=3.5cm,keepaspectratio]{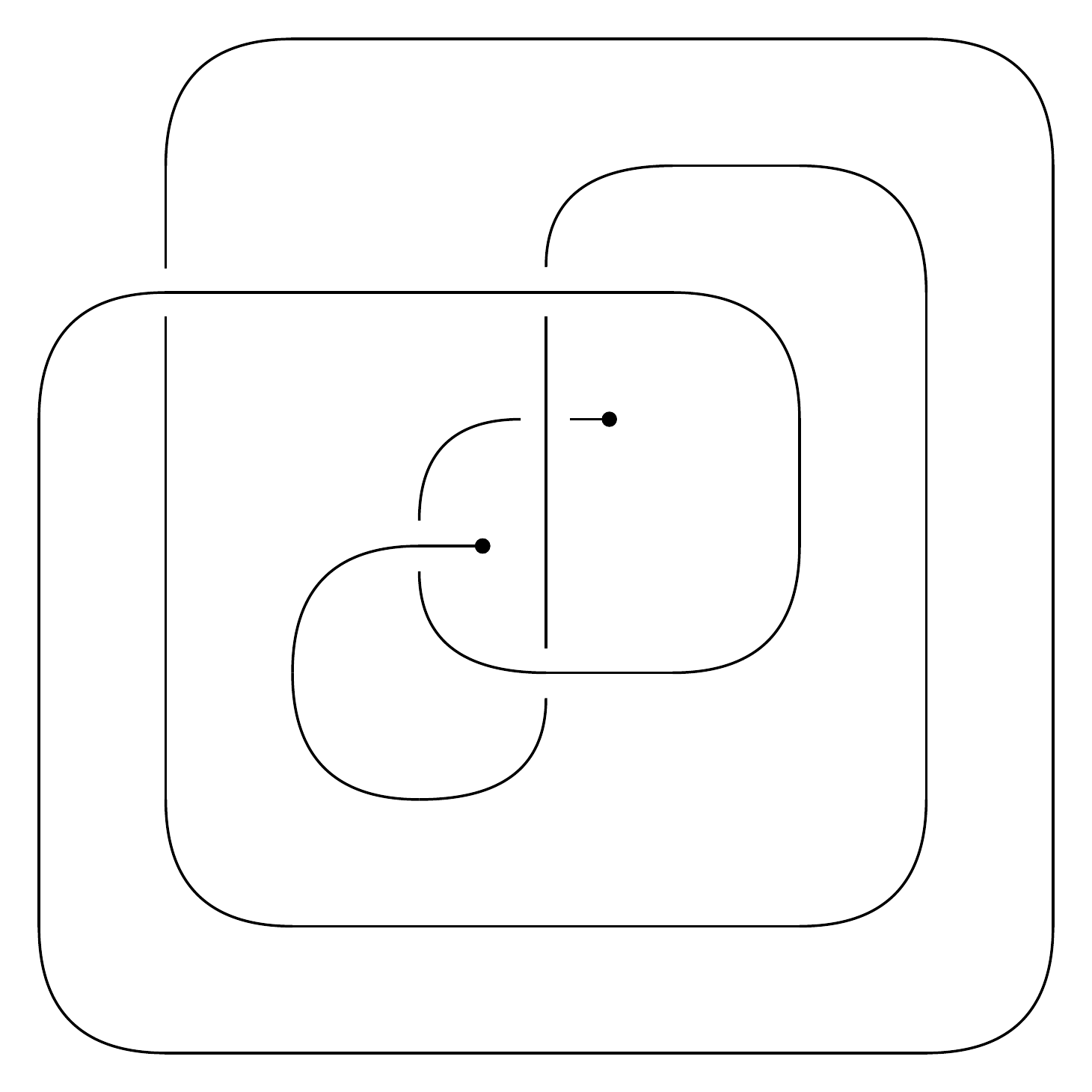}\\
\textcolor{black}{$5_{457}$}
\vspace{1cm}
\end{minipage}
\begin{minipage}[t]{.25\linewidth}
\centering
\includegraphics[width=0.9\textwidth,height=3.5cm,keepaspectratio]{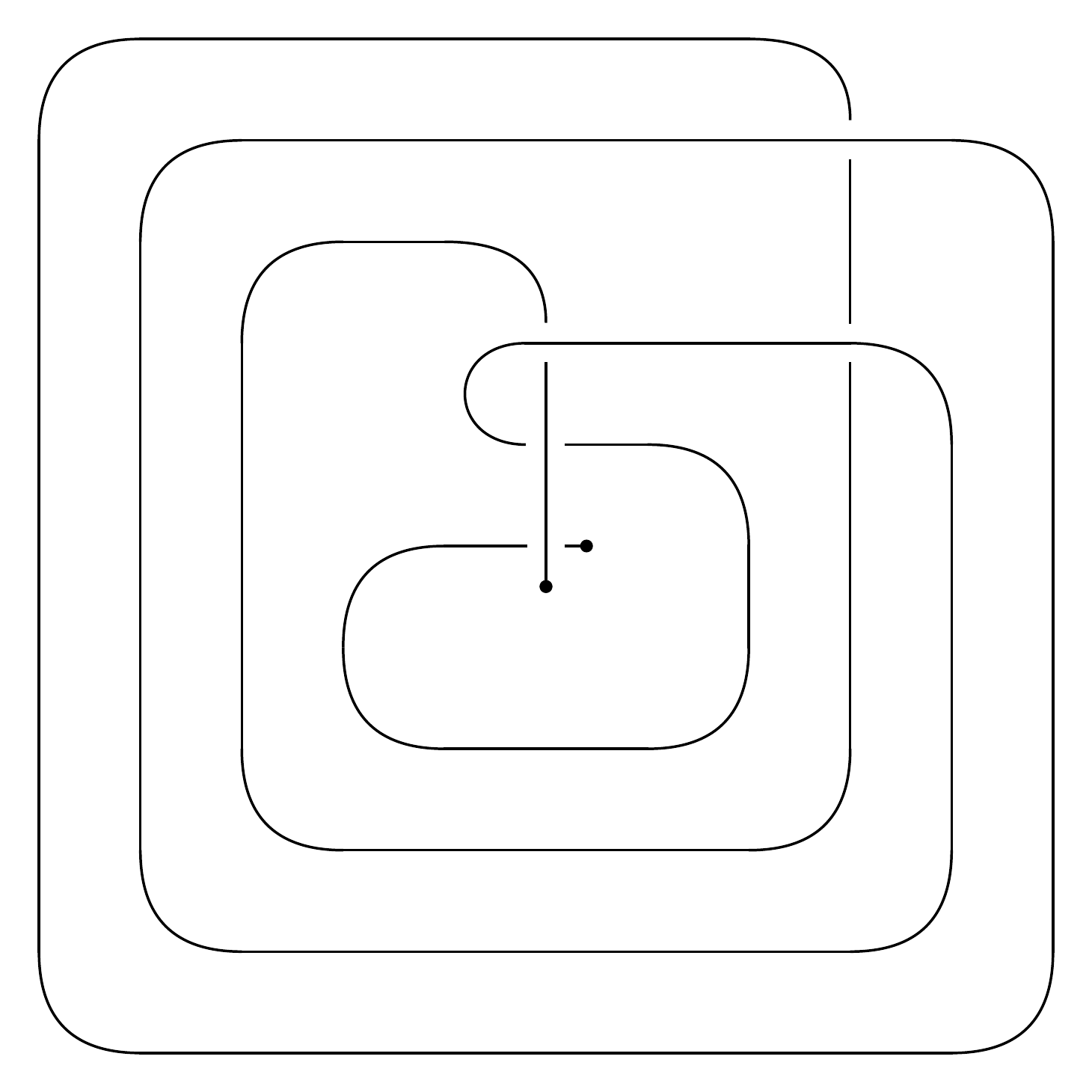}\\
\textcolor{black}{$5_{458}$}
\vspace{1cm}
\end{minipage}
\begin{minipage}[t]{.25\linewidth}
\centering
\includegraphics[width=0.9\textwidth,height=3.5cm,keepaspectratio]{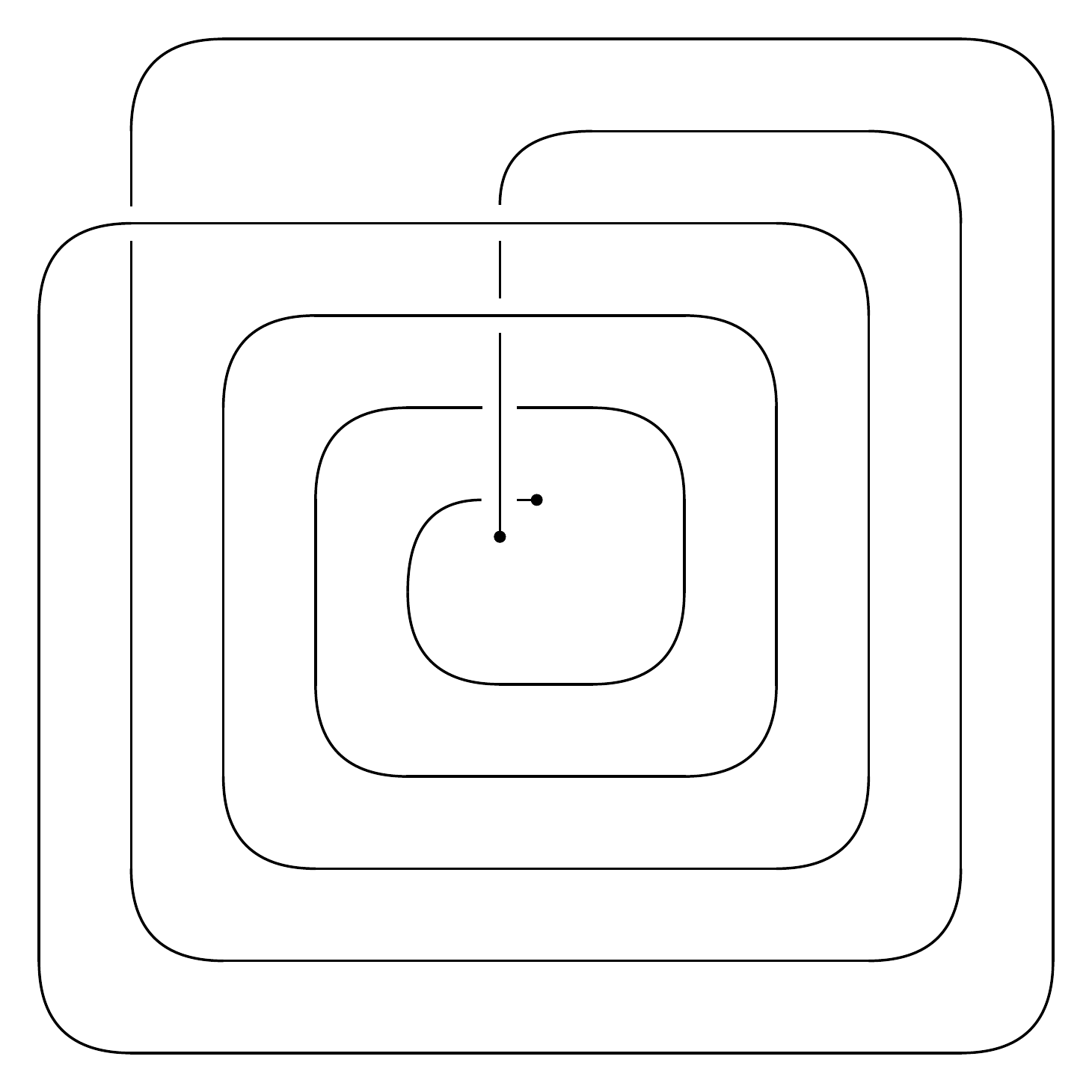}\\
\textcolor{black}{$5_{459}$}
\vspace{1cm}
\end{minipage}
\begin{minipage}[t]{.25\linewidth}
\centering
\includegraphics[width=0.9\textwidth,height=3.5cm,keepaspectratio]{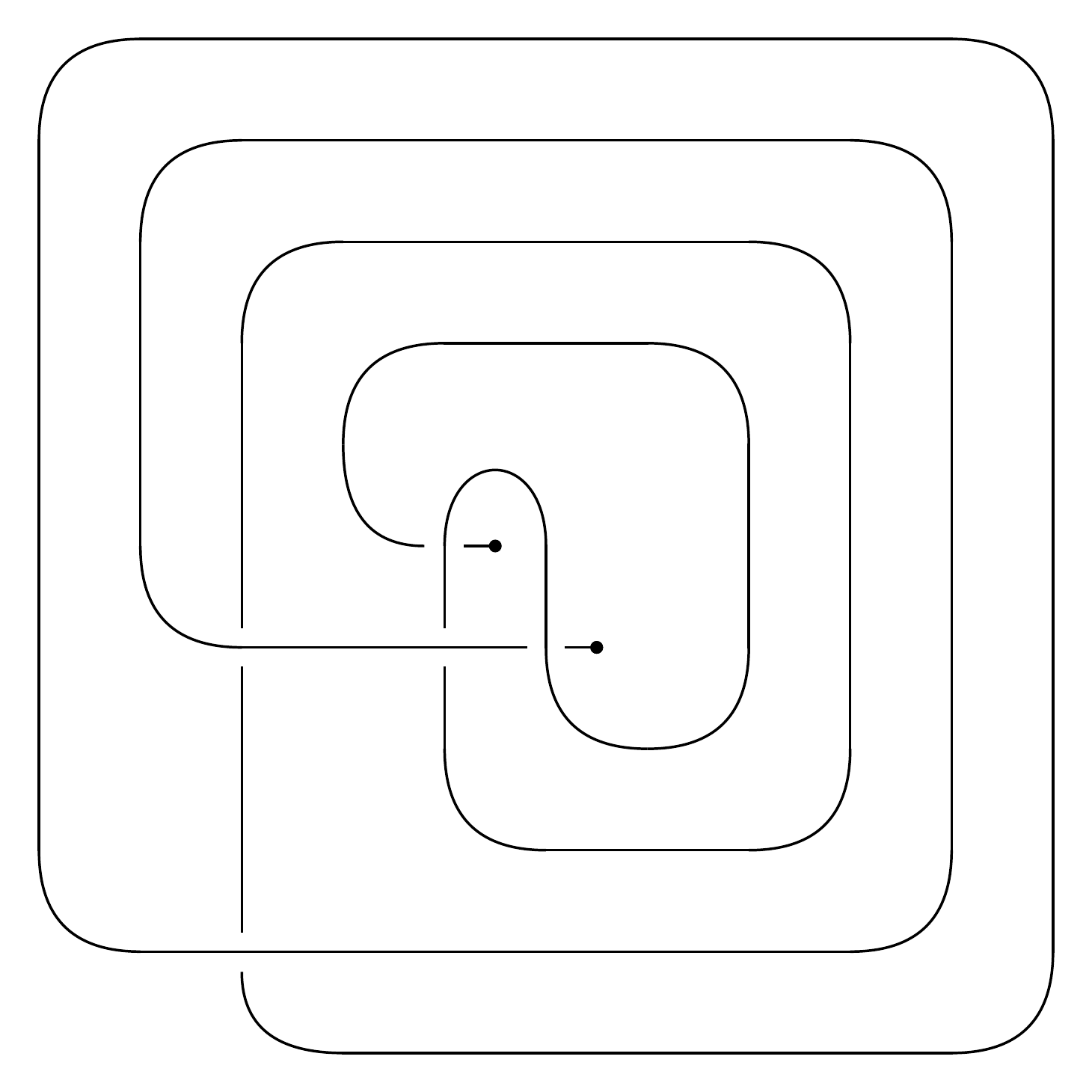}\\
\textcolor{black}{$5_{460}$}
\vspace{1cm}
\end{minipage}
\begin{minipage}[t]{.25\linewidth}
\centering
\includegraphics[width=0.9\textwidth,height=3.5cm,keepaspectratio]{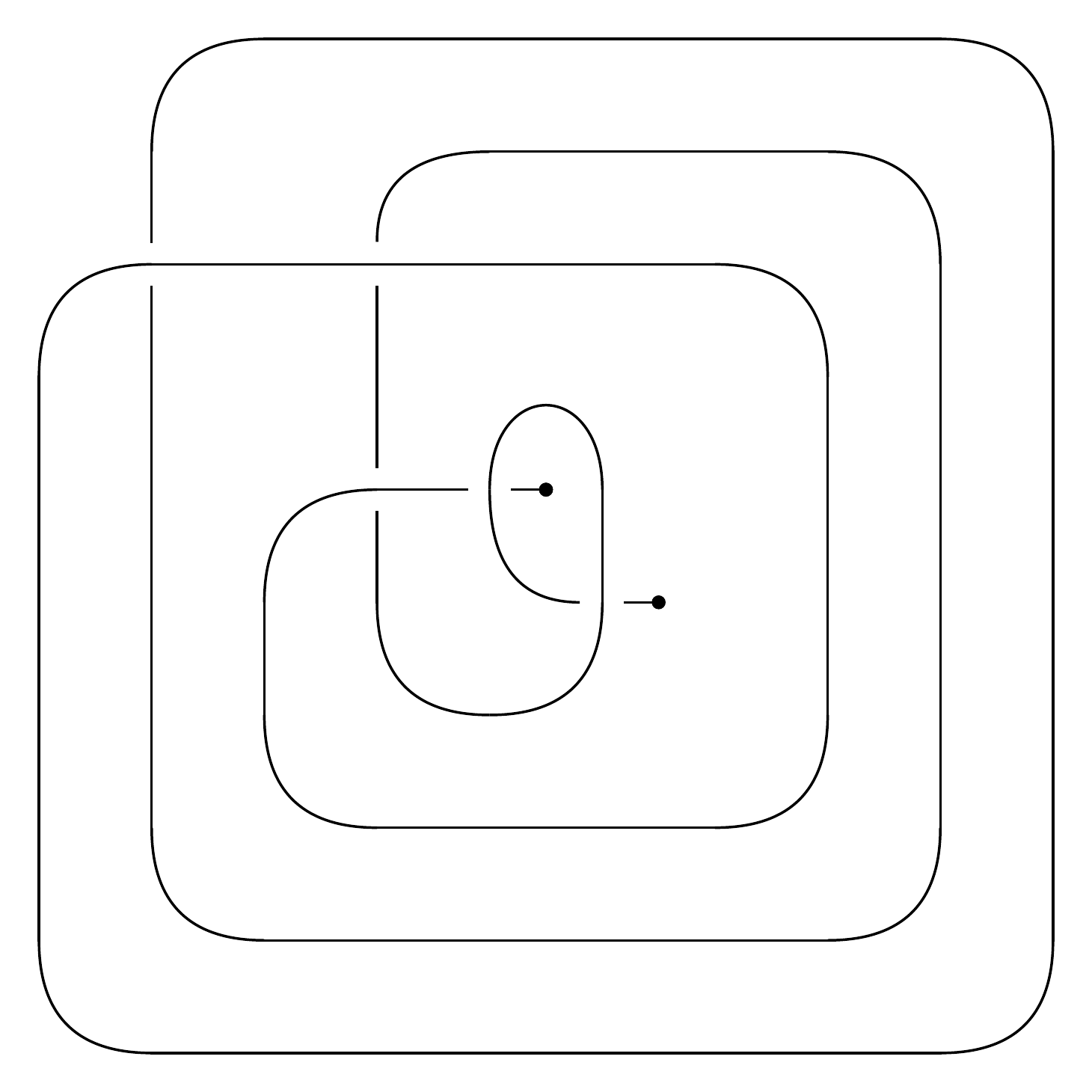}\\
\textcolor{black}{$5_{461}$}
\vspace{1cm}
\end{minipage}
\begin{minipage}[t]{.25\linewidth}
\centering
\includegraphics[width=0.9\textwidth,height=3.5cm,keepaspectratio]{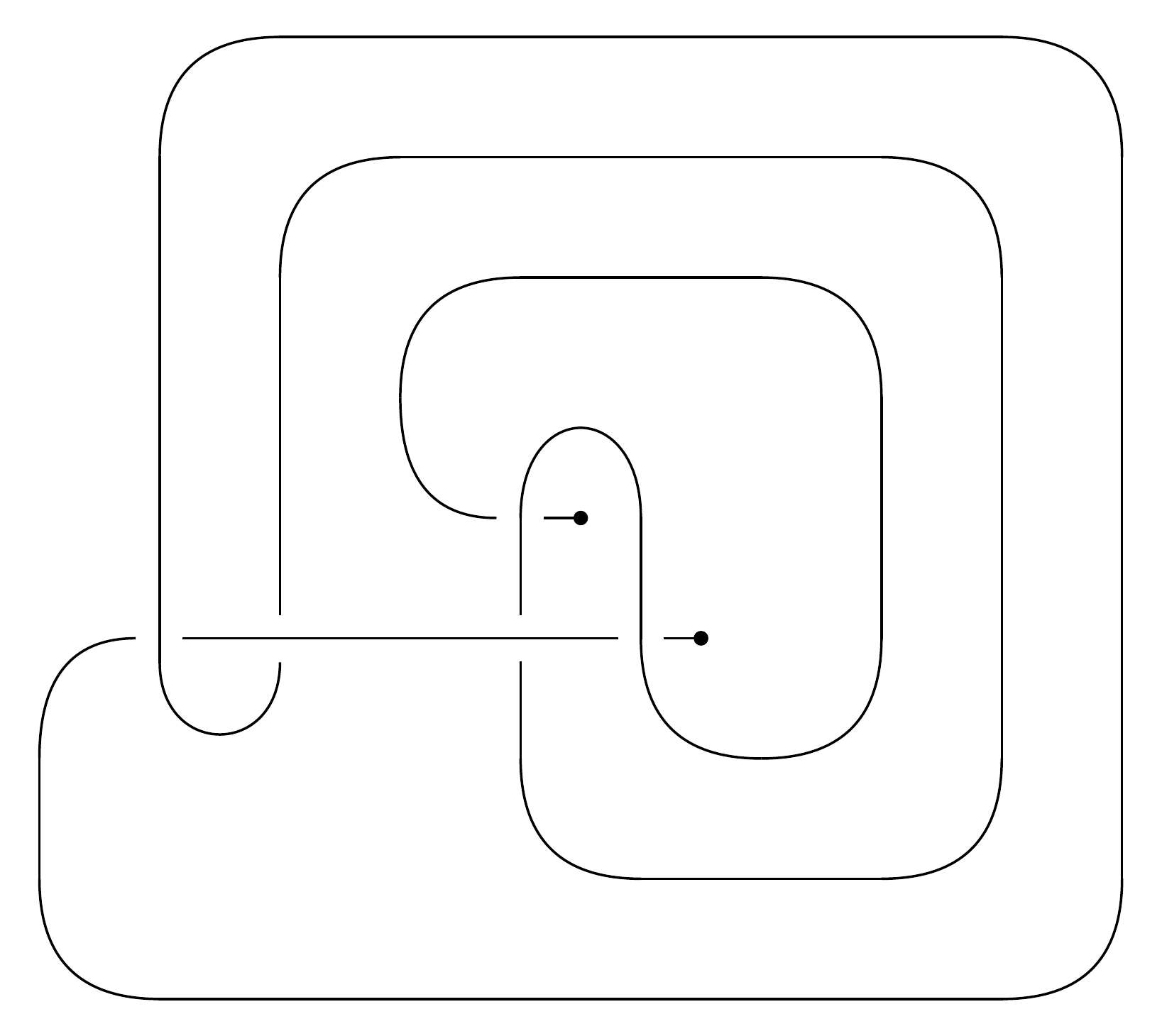}\\
\textcolor{black}{$5_{462}$}
\vspace{1cm}
\end{minipage}
\begin{minipage}[t]{.25\linewidth}
\centering
\includegraphics[width=0.9\textwidth,height=3.5cm,keepaspectratio]{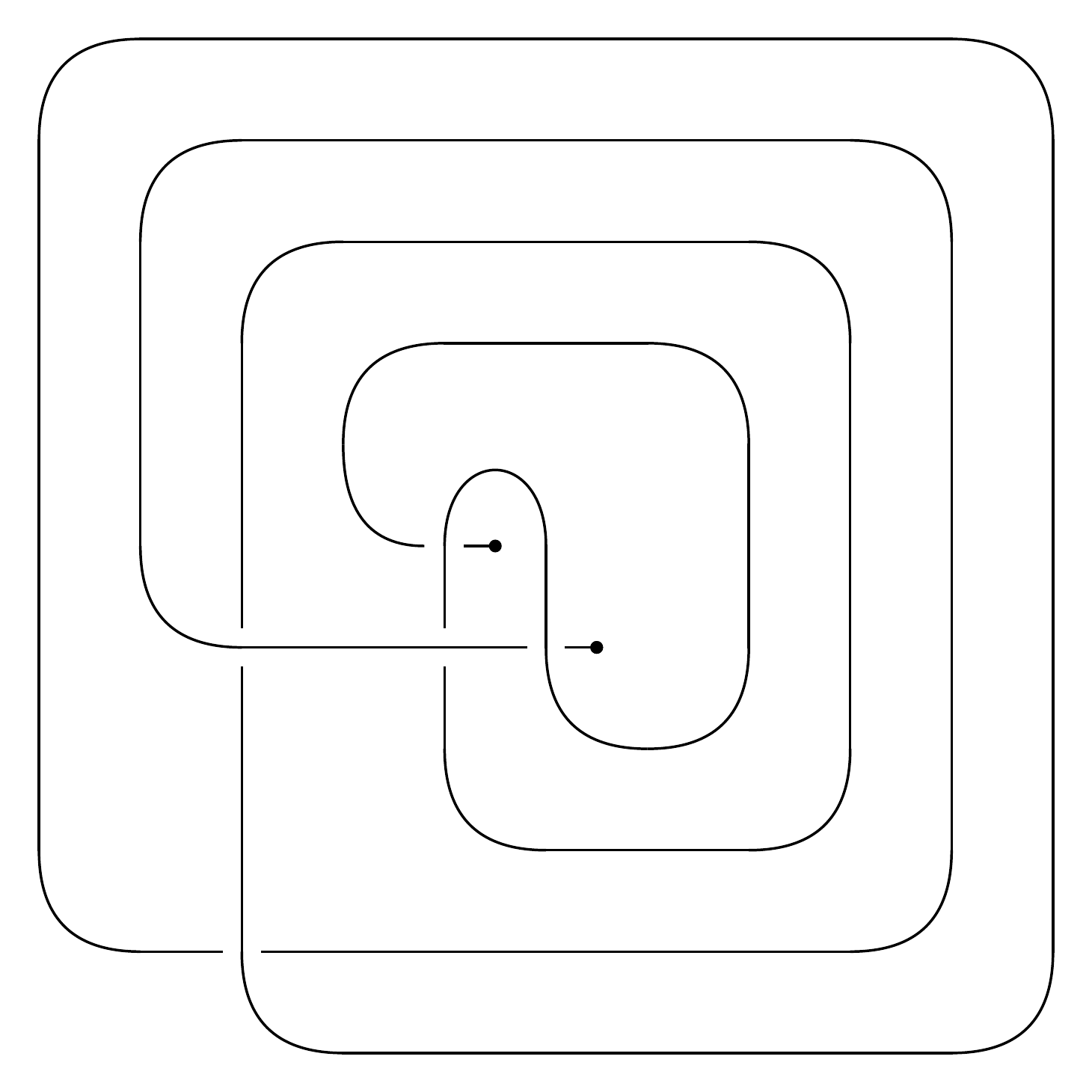}\\
\textcolor{black}{$5_{463}$}
\vspace{1cm}
\end{minipage}
\begin{minipage}[t]{.25\linewidth}
\centering
\includegraphics[width=0.9\textwidth,height=3.5cm,keepaspectratio]{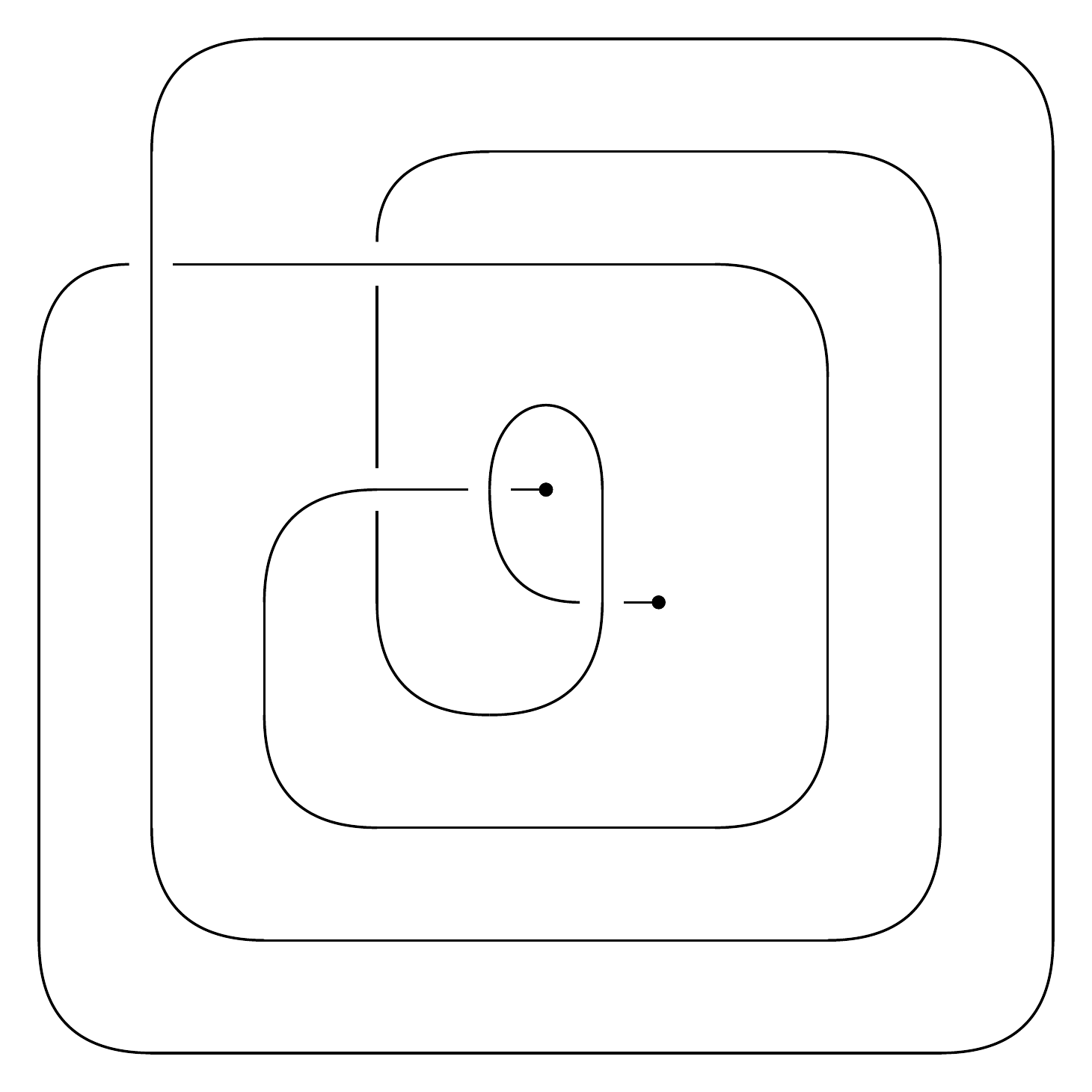}\\
\textcolor{black}{$5_{464}$}
\vspace{1cm}
\end{minipage}
\begin{minipage}[t]{.25\linewidth}
\centering
\includegraphics[width=0.9\textwidth,height=3.5cm,keepaspectratio]{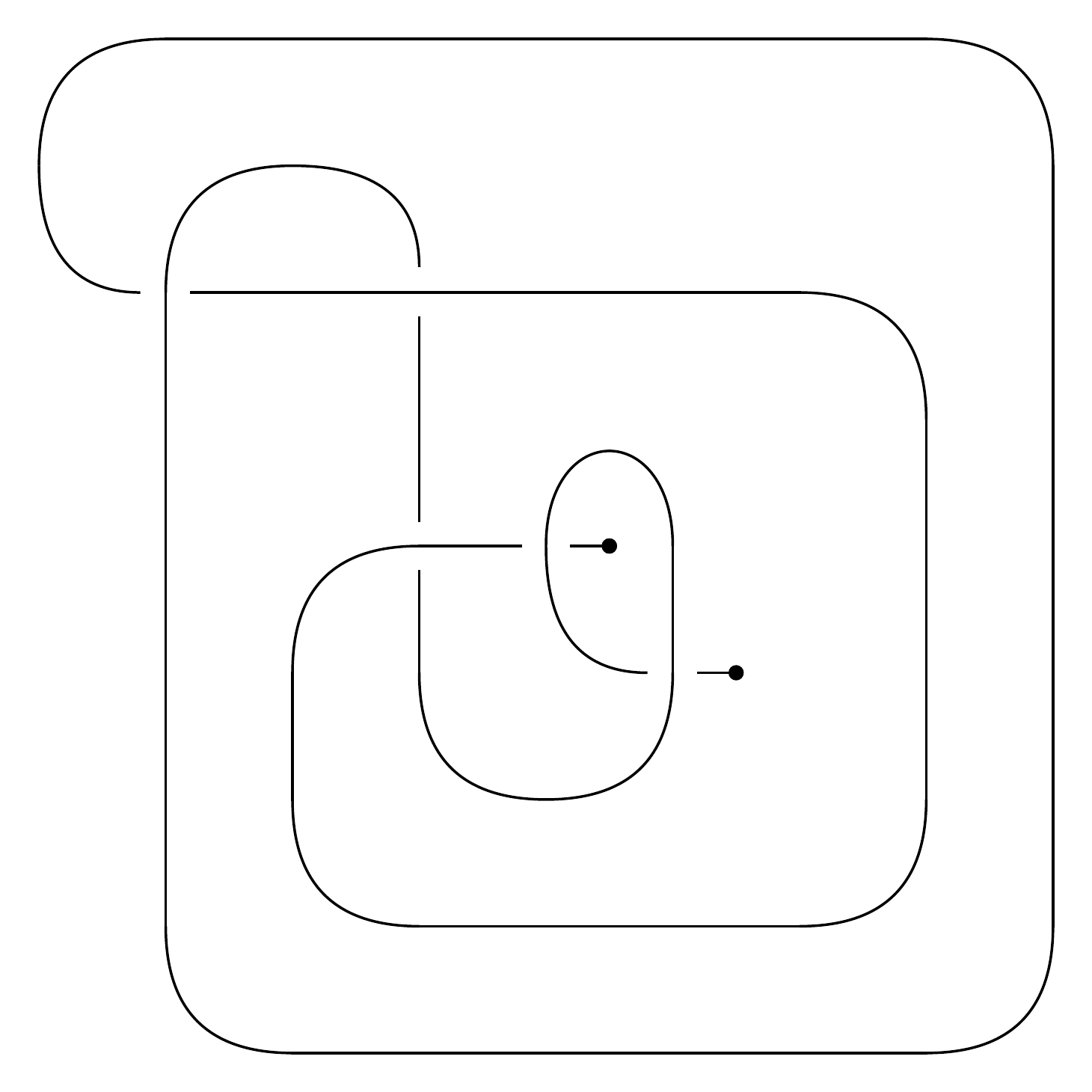}\\
\textcolor{black}{$5_{465}$}
\vspace{1cm}
\end{minipage}
\begin{minipage}[t]{.25\linewidth}
\centering
\includegraphics[width=0.9\textwidth,height=3.5cm,keepaspectratio]{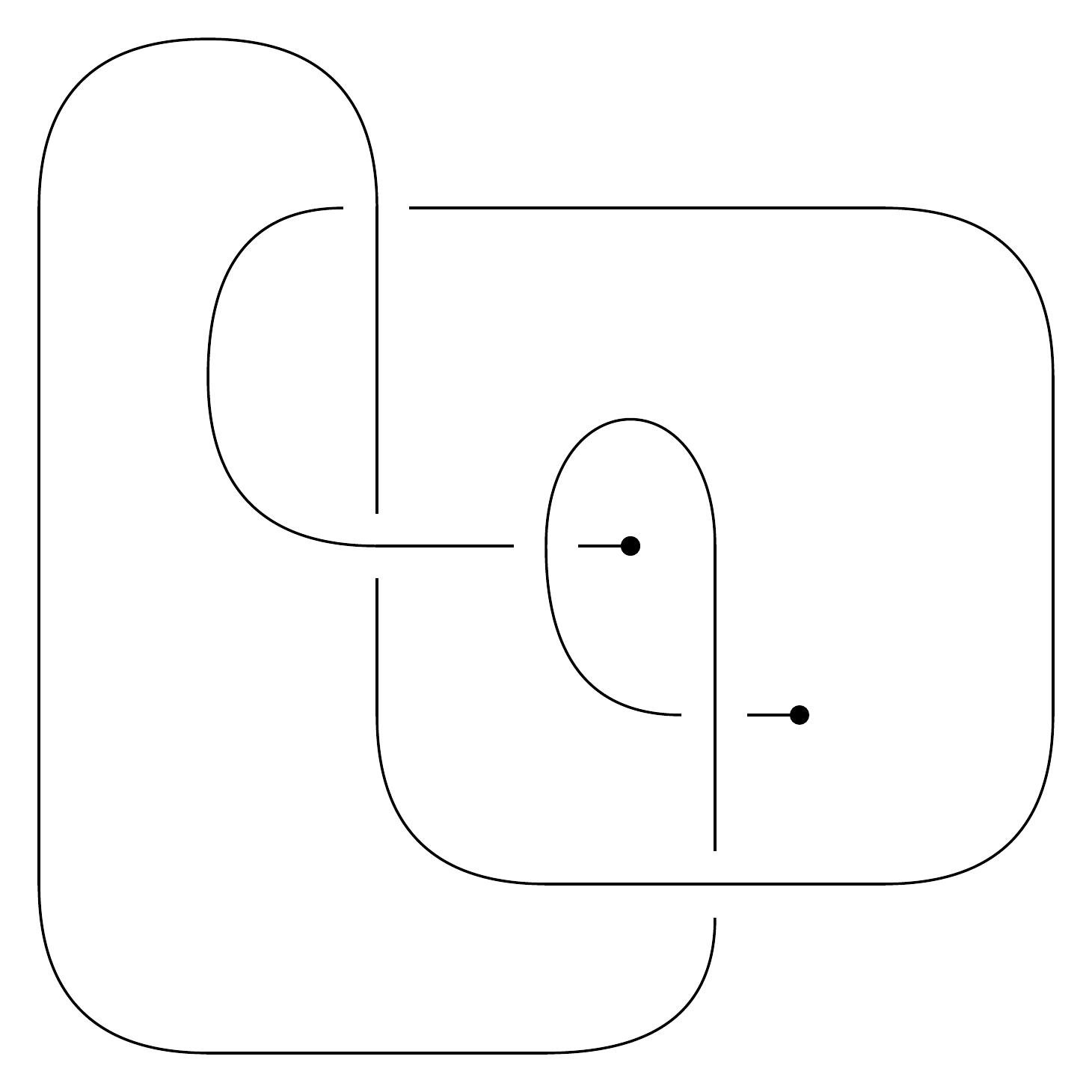}\\
\textcolor{black}{$5_{466}$}
\vspace{1cm}
\end{minipage}
\begin{minipage}[t]{.25\linewidth}
\centering
\includegraphics[width=0.9\textwidth,height=3.5cm,keepaspectratio]{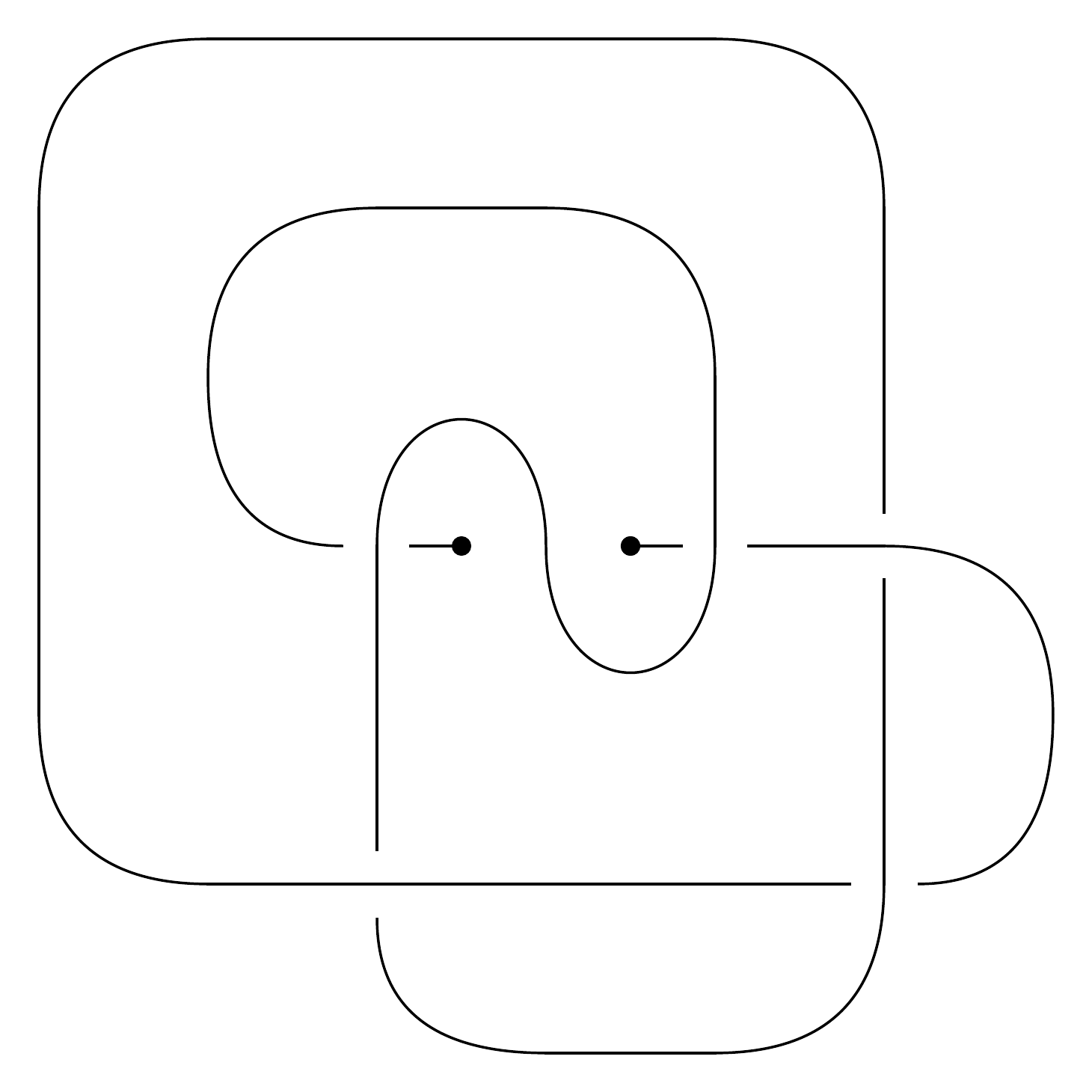}\\
\textcolor{black}{$5_{467}$}
\vspace{1cm}
\end{minipage}
\begin{minipage}[t]{.25\linewidth}
\centering
\includegraphics[width=0.9\textwidth,height=3.5cm,keepaspectratio]{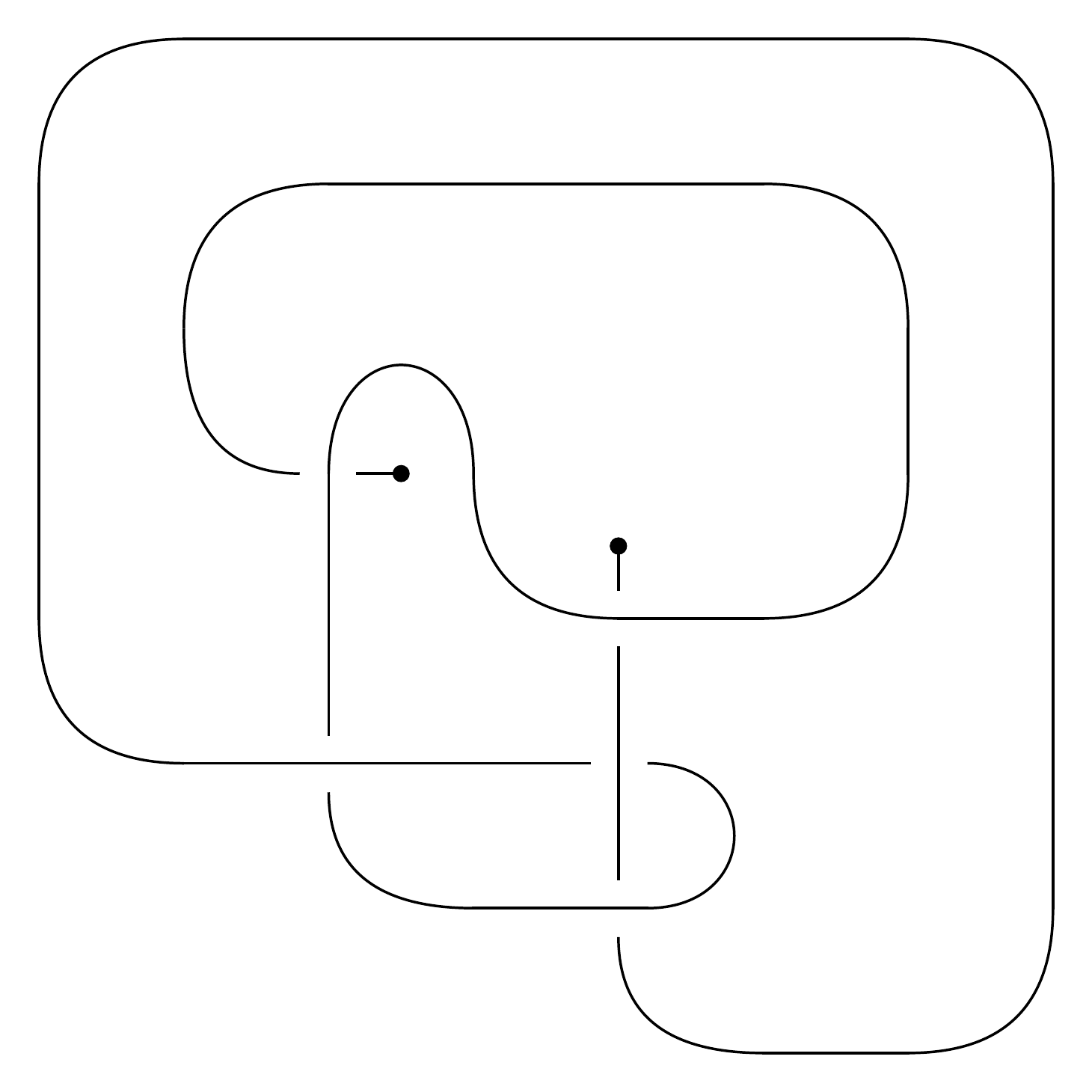}\\
\textcolor{black}{$5_{468}$}
\vspace{1cm}
\end{minipage}
\begin{minipage}[t]{.25\linewidth}
\centering
\includegraphics[width=0.9\textwidth,height=3.5cm,keepaspectratio]{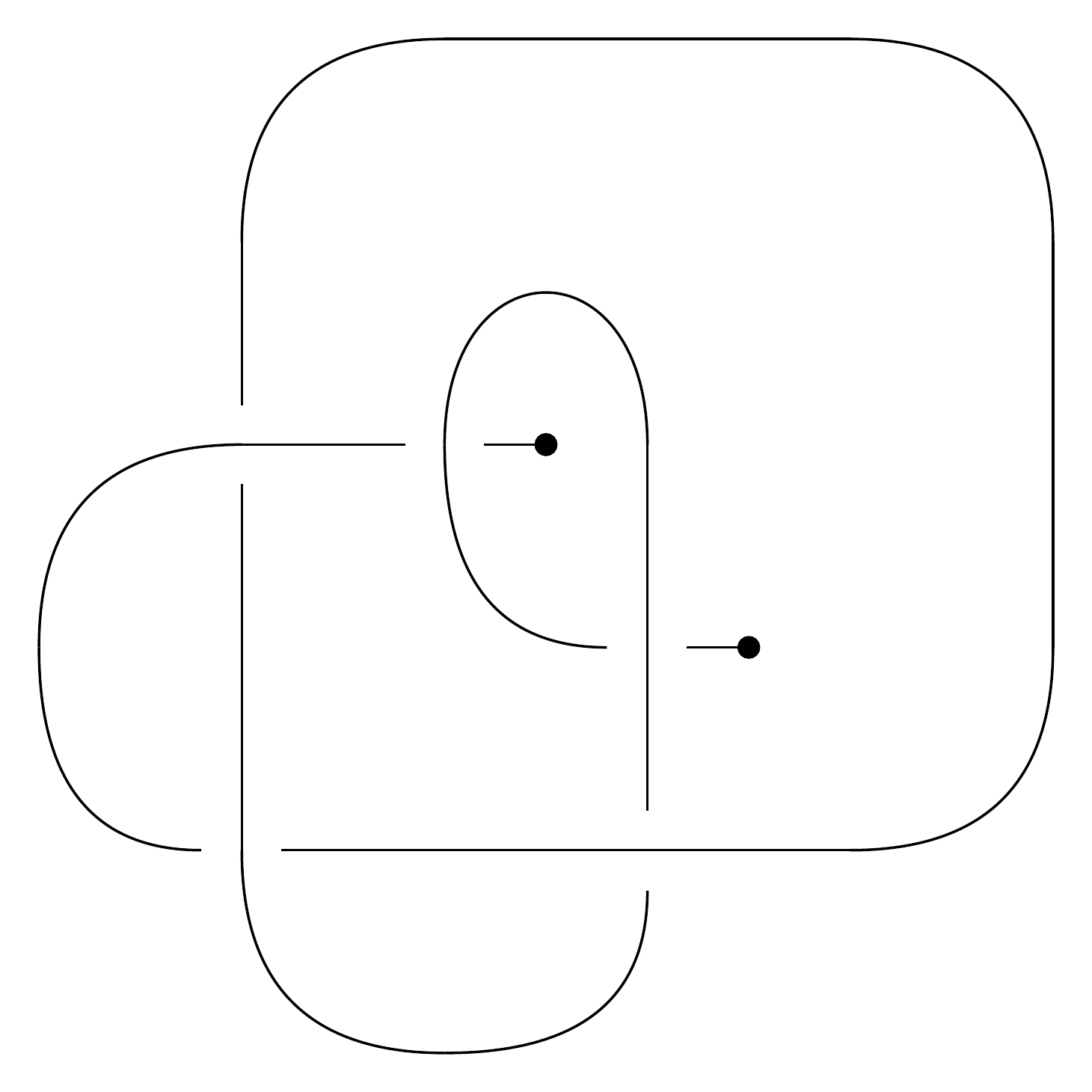}\\
\textcolor{black}{$5_{469}$}
\vspace{1cm}
\end{minipage}
\begin{minipage}[t]{.25\linewidth}
\centering
\includegraphics[width=0.9\textwidth,height=3.5cm,keepaspectratio]{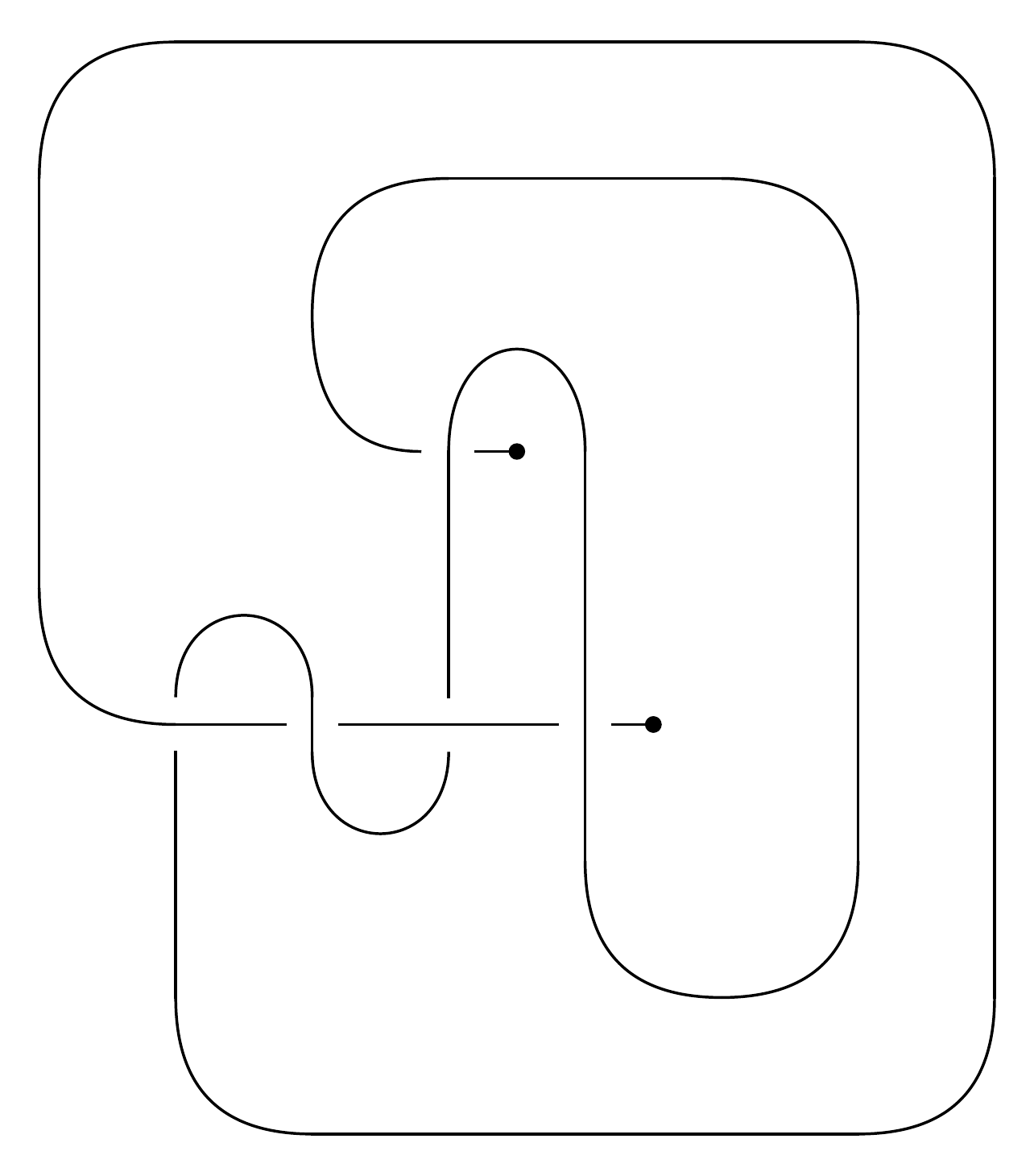}\\
\textcolor{black}{$5_{470}$}
\vspace{1cm}
\end{minipage}
\begin{minipage}[t]{.25\linewidth}
\centering
\includegraphics[width=0.9\textwidth,height=3.5cm,keepaspectratio]{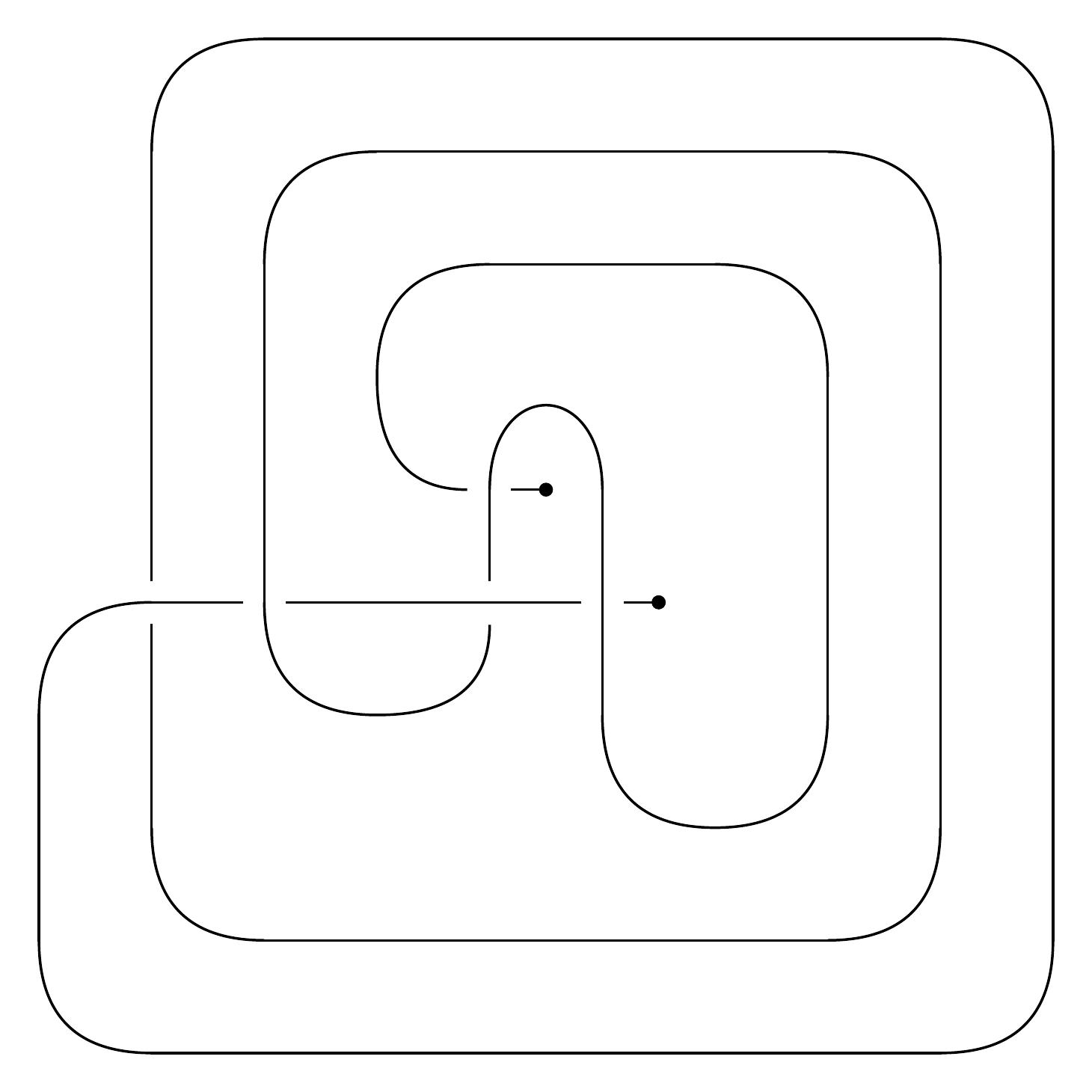}\\
\textcolor{black}{$5_{471}$}
\vspace{1cm}
\end{minipage}
\begin{minipage}[t]{.25\linewidth}
\centering
\includegraphics[width=0.9\textwidth,height=3.5cm,keepaspectratio]{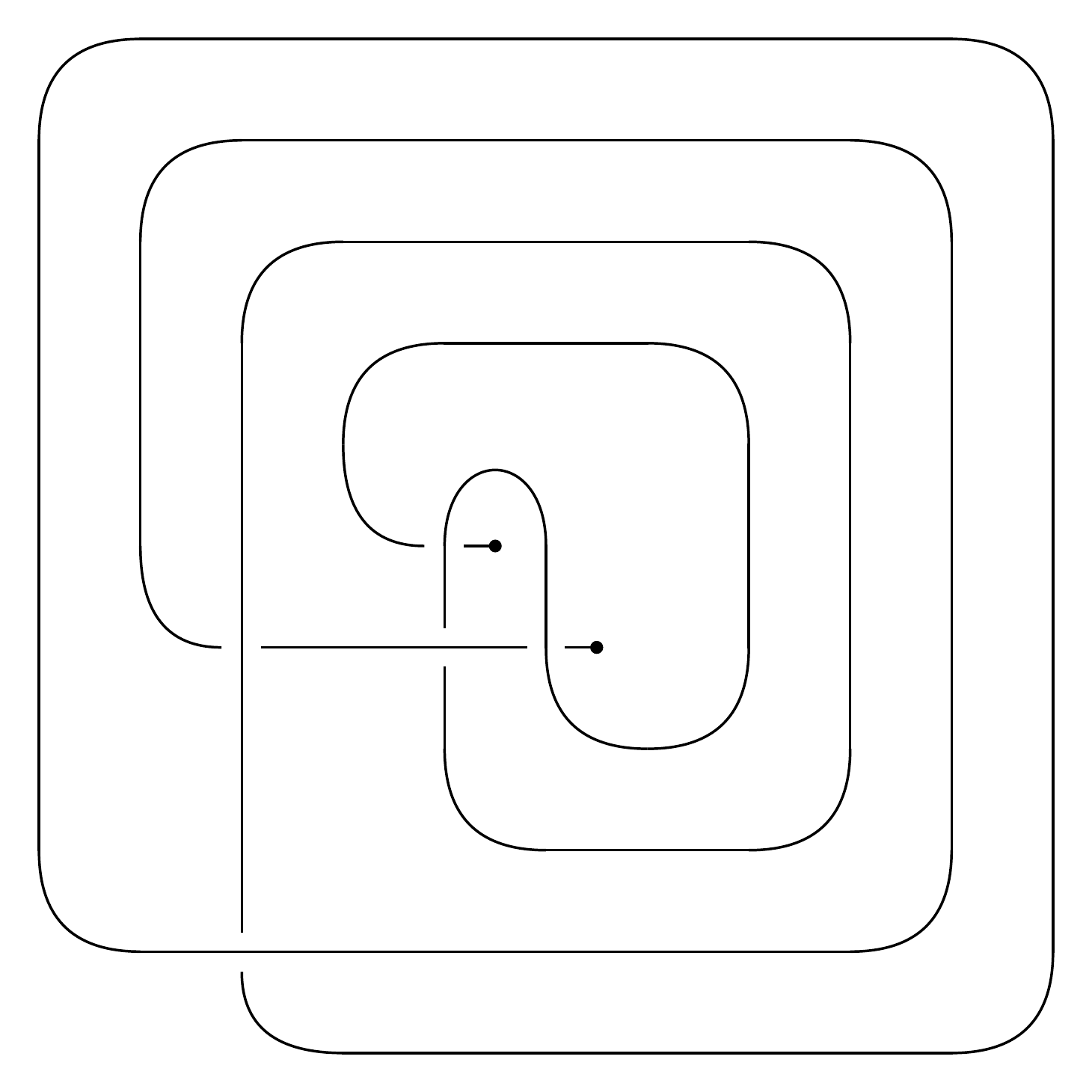}\\
\textcolor{black}{$5_{472}$}
\vspace{1cm}
\end{minipage}
\begin{minipage}[t]{.25\linewidth}
\centering
\includegraphics[width=0.9\textwidth,height=3.5cm,keepaspectratio]{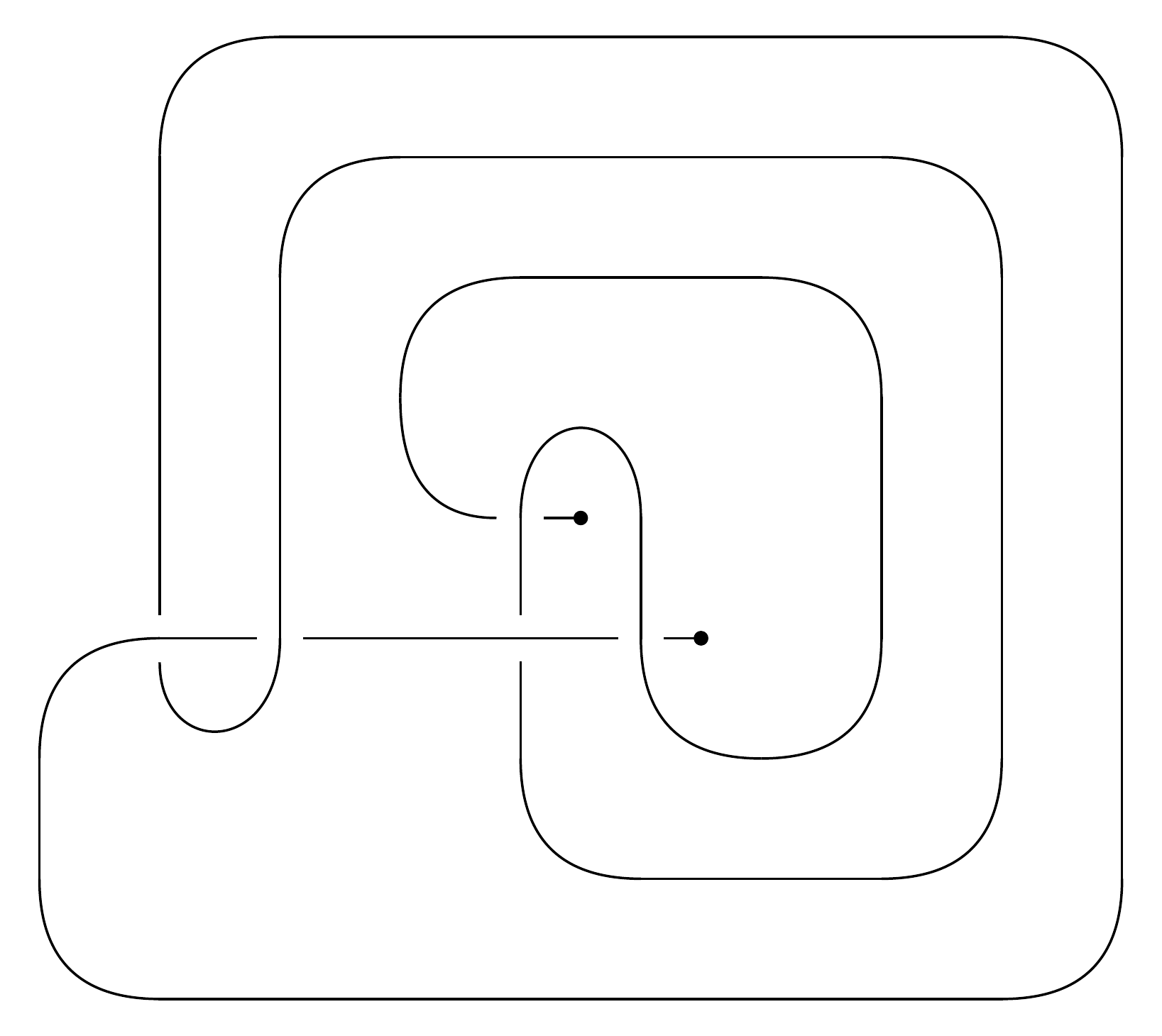}\\
\textcolor{black}{$5_{473}$}
\vspace{1cm}
\end{minipage}
\begin{minipage}[t]{.25\linewidth}
\centering
\includegraphics[width=0.9\textwidth,height=3.5cm,keepaspectratio]{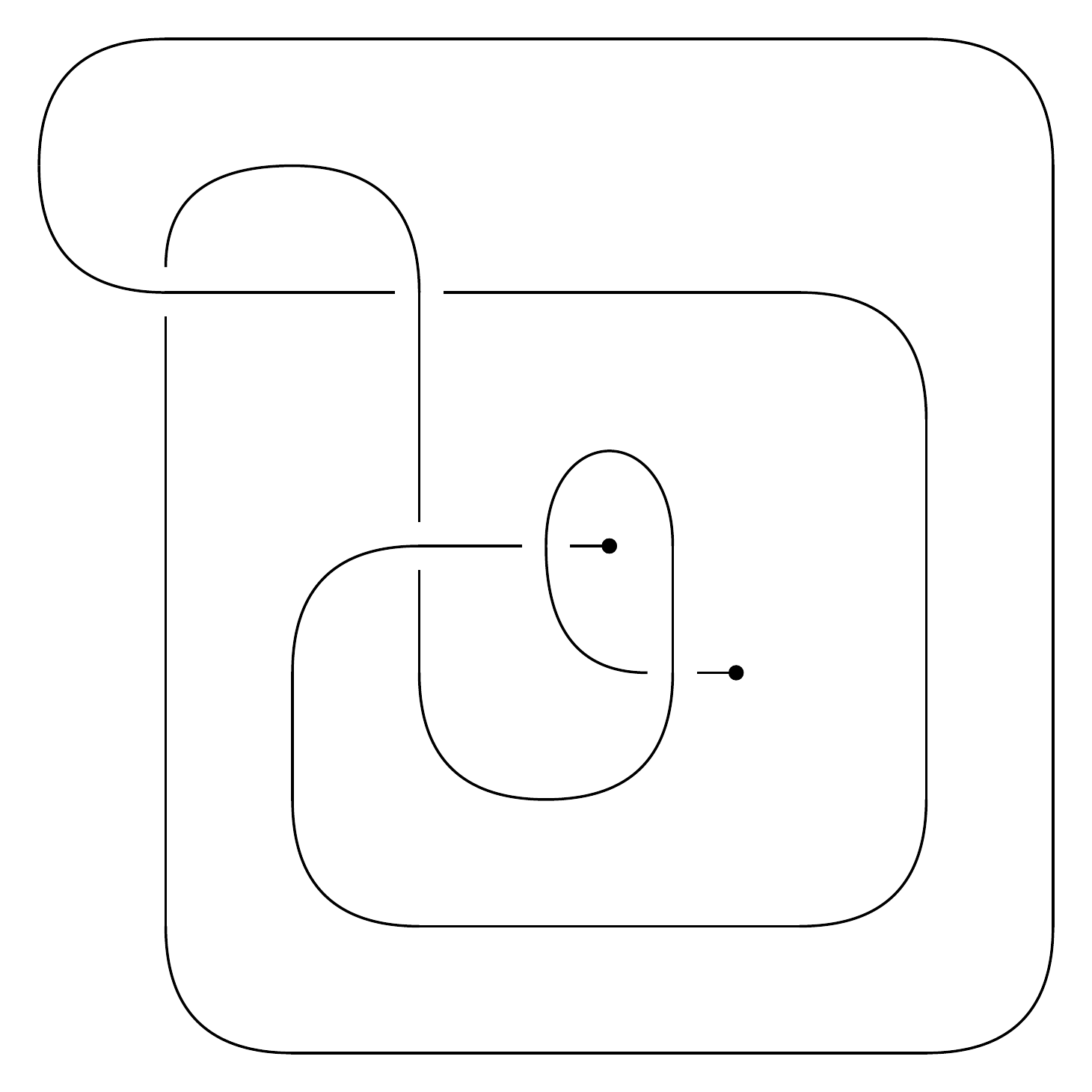}\\
\textcolor{black}{$5_{474}$}
\vspace{1cm}
\end{minipage}
\begin{minipage}[t]{.25\linewidth}
\centering
\includegraphics[width=0.9\textwidth,height=3.5cm,keepaspectratio]{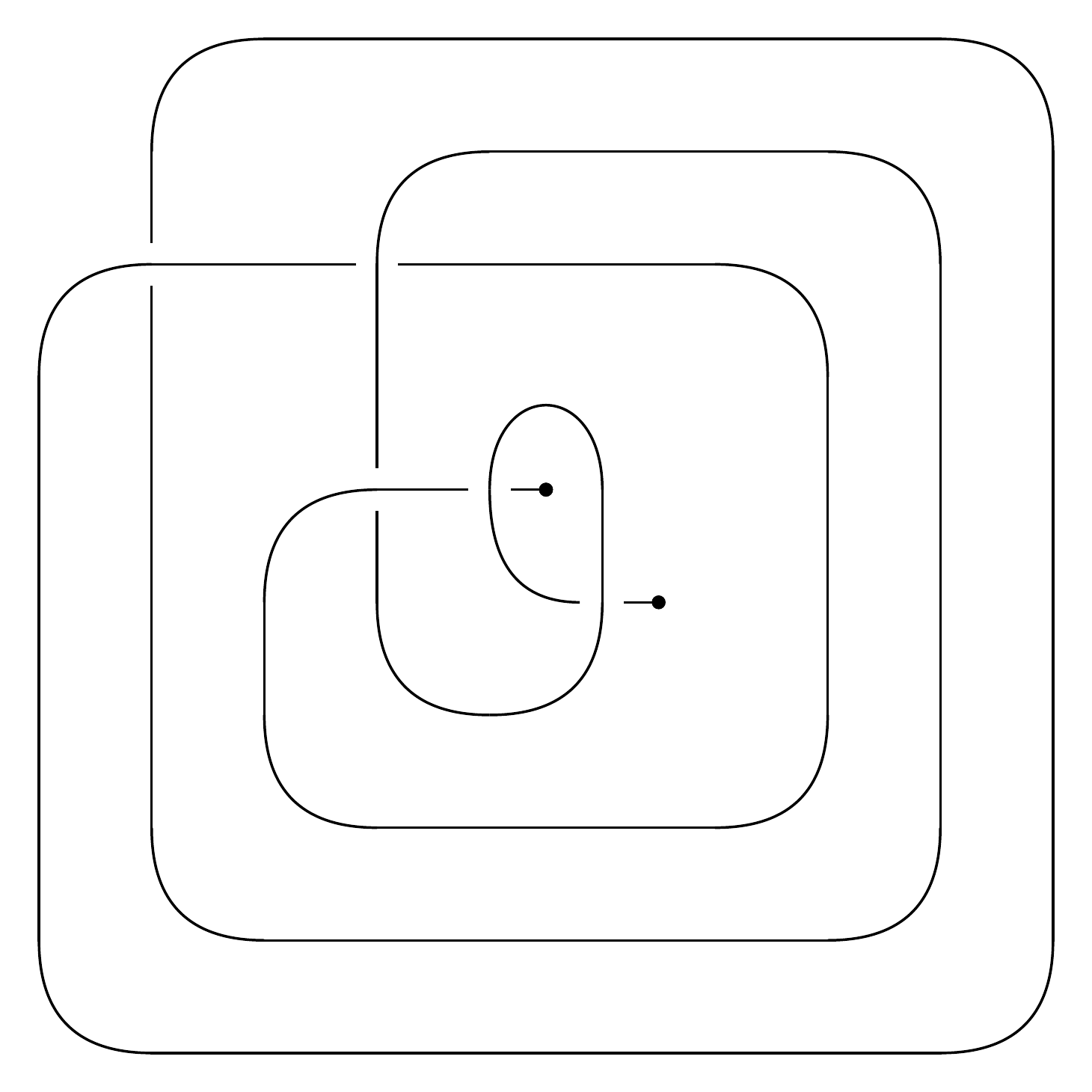}\\
\textcolor{black}{$5_{475}$}
\vspace{1cm}
\end{minipage}
\begin{minipage}[t]{.25\linewidth}
\centering
\includegraphics[width=0.9\textwidth,height=3.5cm,keepaspectratio]{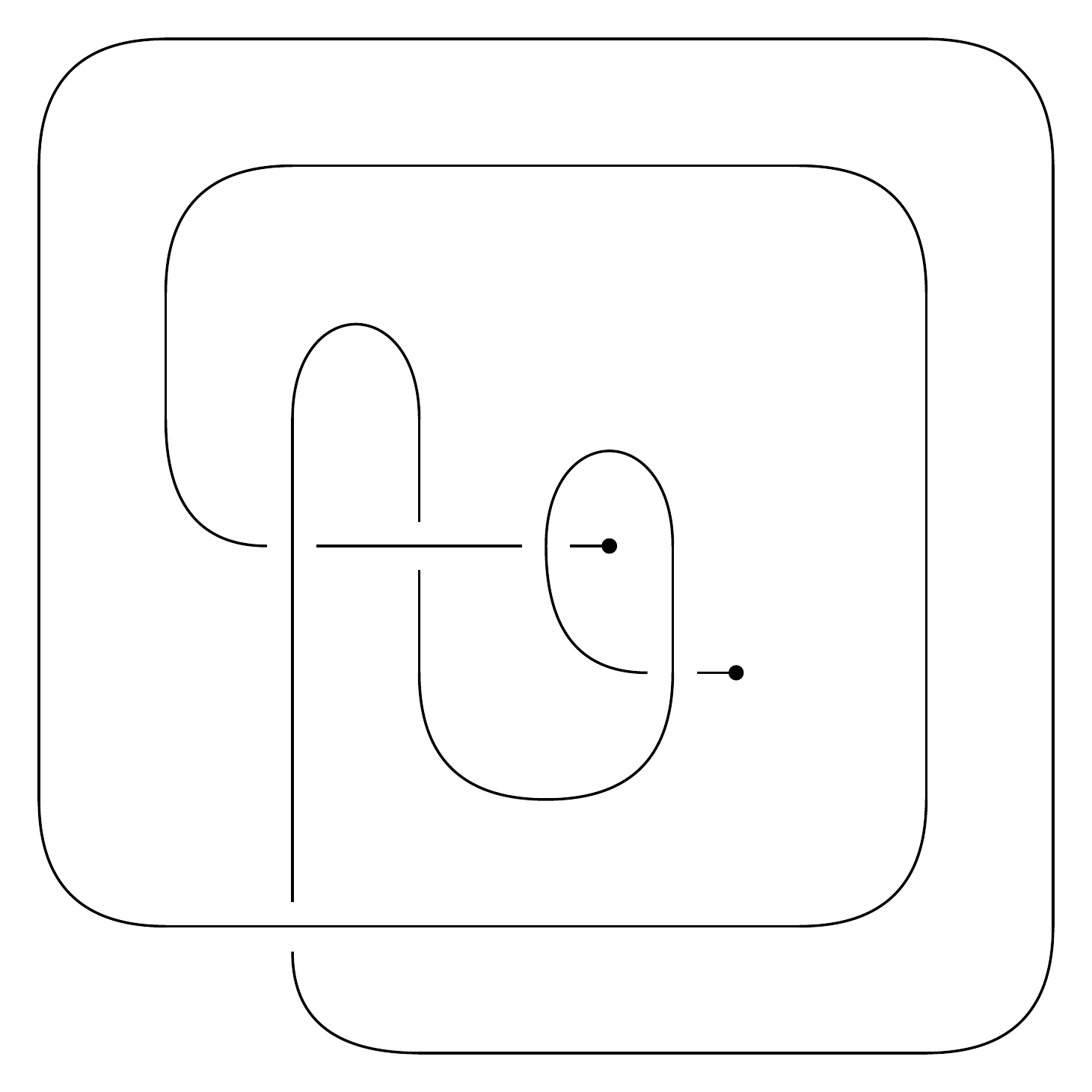}\\
\textcolor{black}{$5_{476}$}
\vspace{1cm}
\end{minipage}
\begin{minipage}[t]{.25\linewidth}
\centering
\includegraphics[width=0.9\textwidth,height=3.5cm,keepaspectratio]{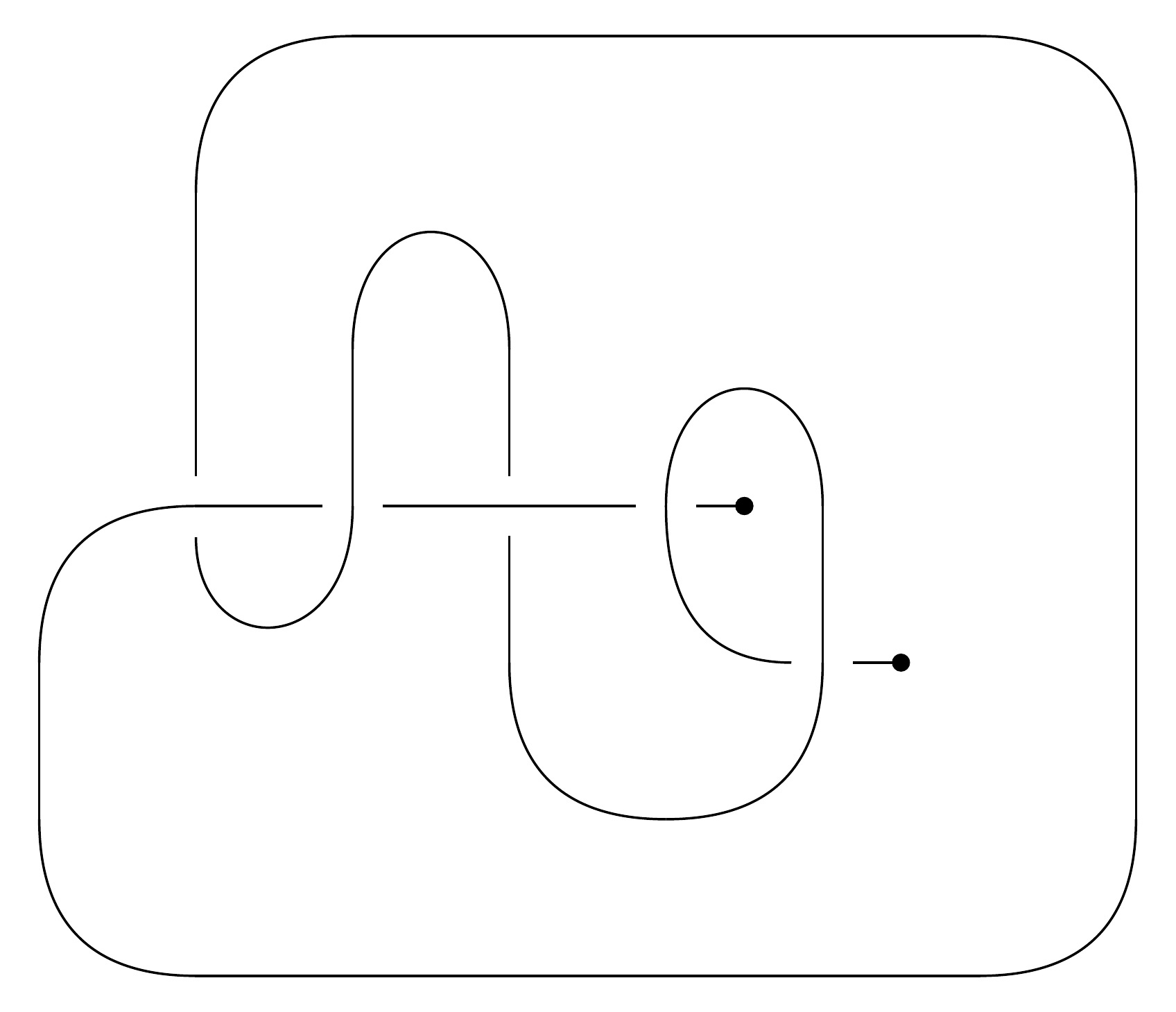}\\
\textcolor{black}{$5_{477}$}
\vspace{1cm}
\end{minipage}
\begin{minipage}[t]{.25\linewidth}
\centering
\includegraphics[width=0.9\textwidth,height=3.5cm,keepaspectratio]{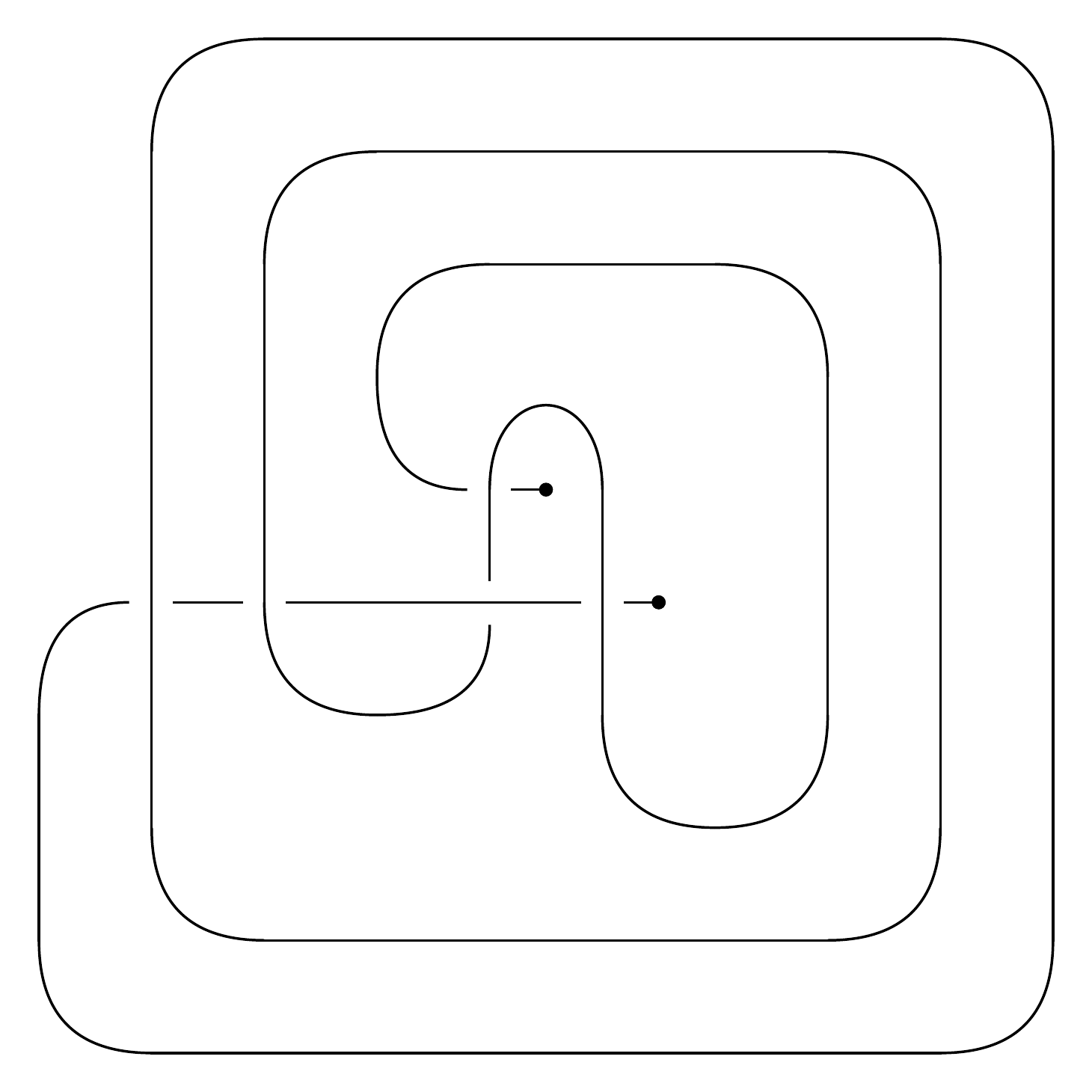}\\
\textcolor{black}{$5_{478}$}
\vspace{1cm}
\end{minipage}
\begin{minipage}[t]{.25\linewidth}
\centering
\includegraphics[width=0.9\textwidth,height=3.5cm,keepaspectratio]{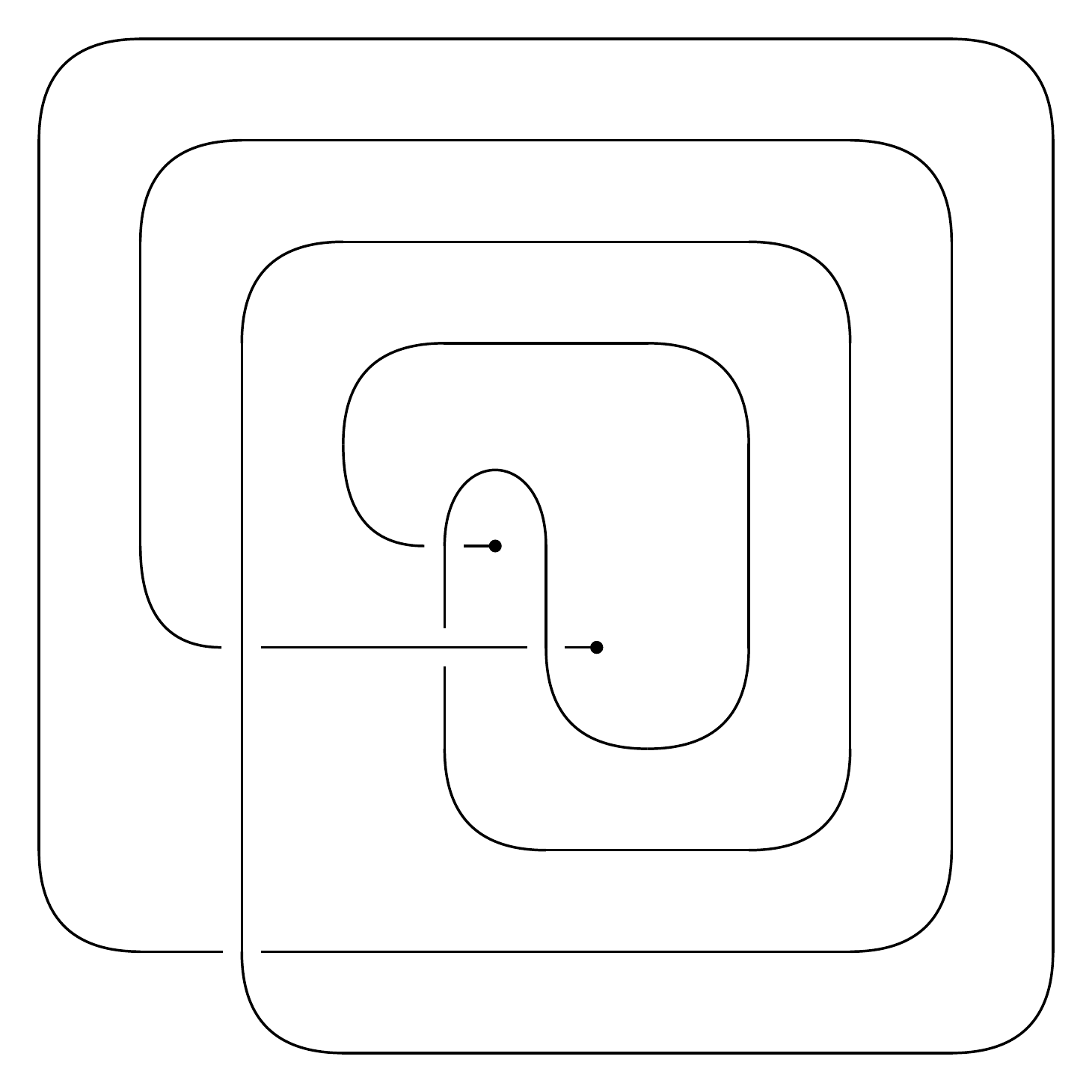}\\
\textcolor{black}{$5_{479}$}
\vspace{1cm}
\end{minipage}
\begin{minipage}[t]{.25\linewidth}
\centering
\includegraphics[width=0.9\textwidth,height=3.5cm,keepaspectratio]{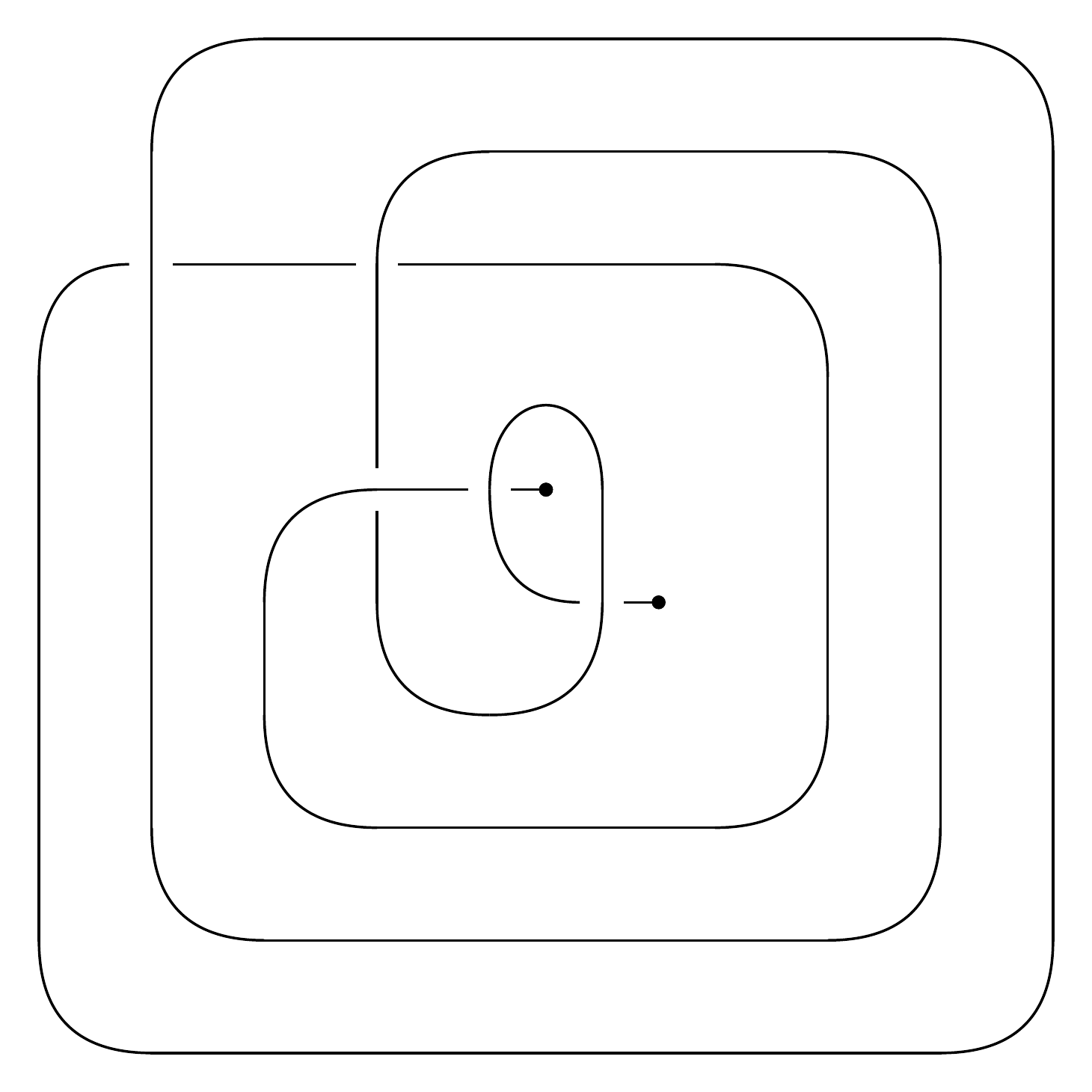}\\
\textcolor{black}{$5_{480}$}
\vspace{1cm}
\end{minipage}
\begin{minipage}[t]{.25\linewidth}
\centering
\includegraphics[width=0.9\textwidth,height=3.5cm,keepaspectratio]{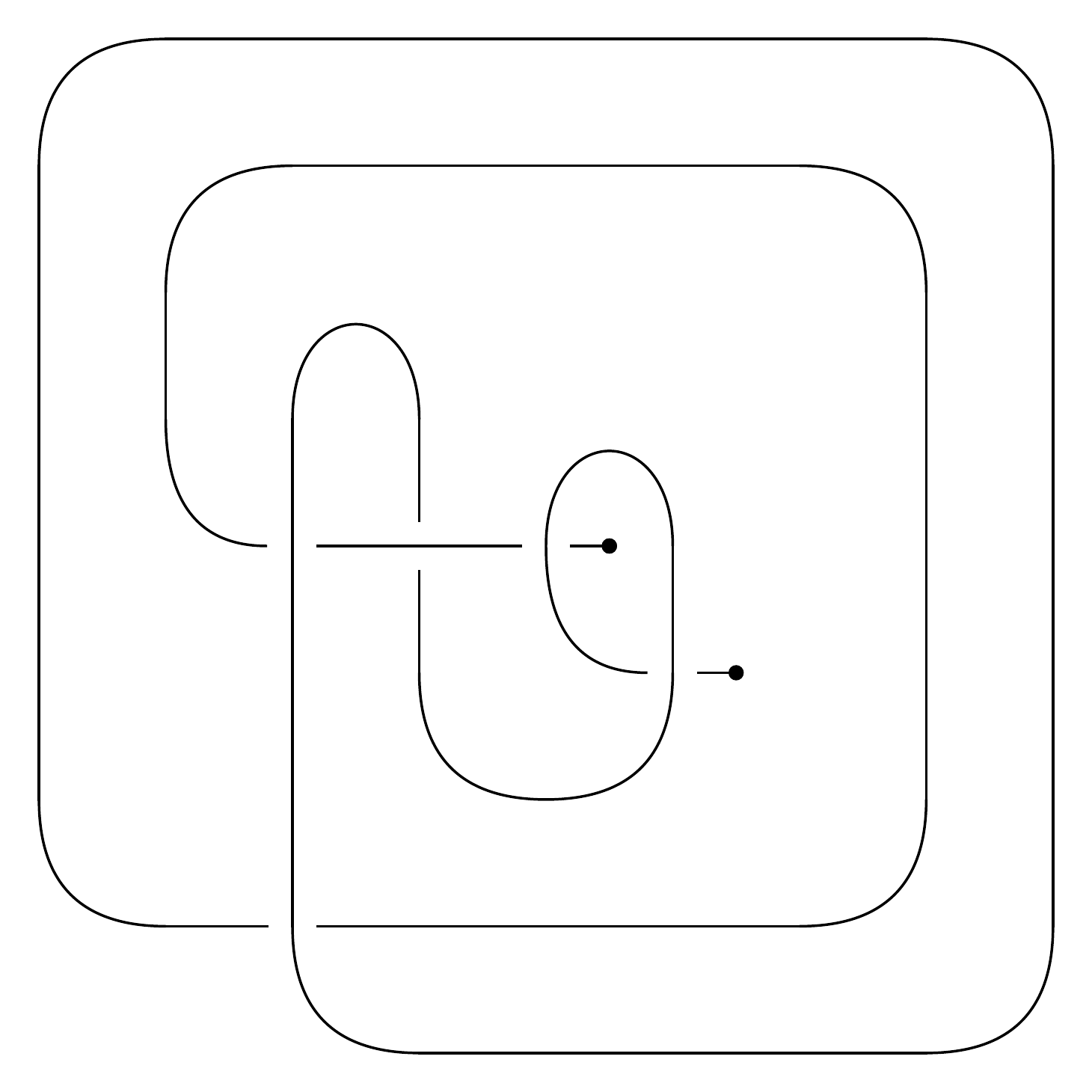}\\
\textcolor{black}{$5_{481}$}
\vspace{1cm}
\end{minipage}
\begin{minipage}[t]{.25\linewidth}
\centering
\includegraphics[width=0.9\textwidth,height=3.5cm,keepaspectratio]{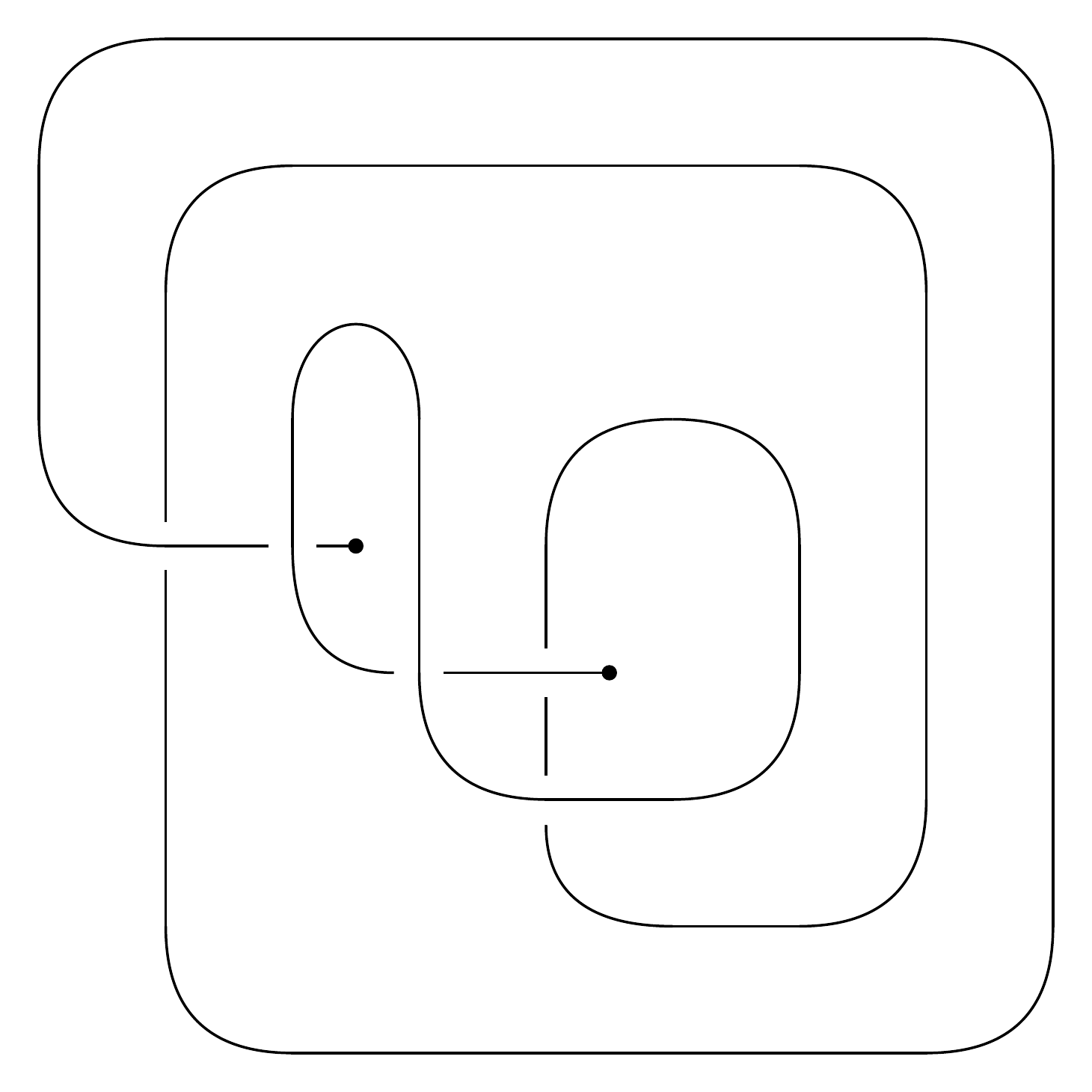}\\
\textcolor{black}{$5_{482}$}
\vspace{1cm}
\end{minipage}
\begin{minipage}[t]{.25\linewidth}
\centering
\includegraphics[width=0.9\textwidth,height=3.5cm,keepaspectratio]{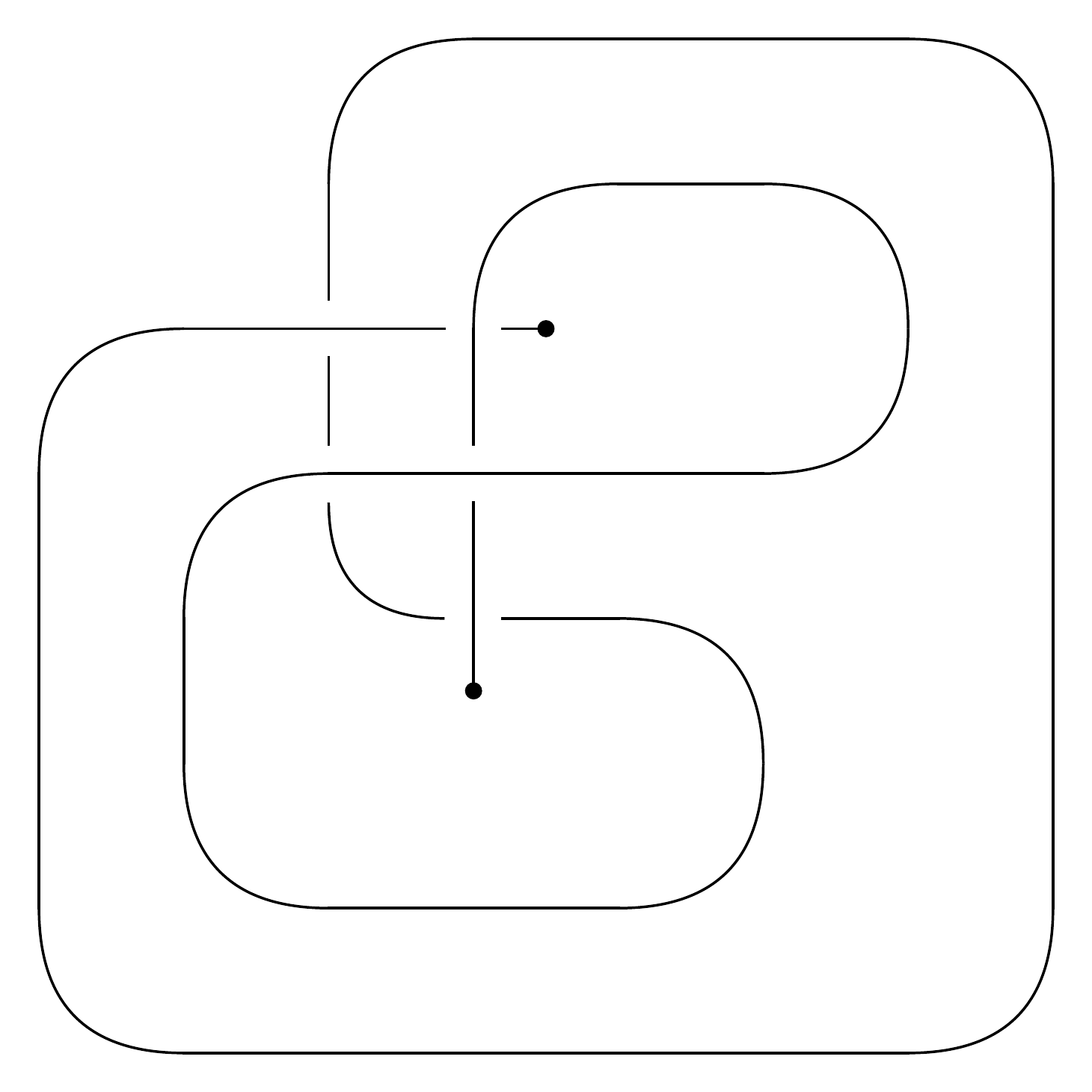}\\
\textcolor{black}{$5_{483}$}
\vspace{1cm}
\end{minipage}
\begin{minipage}[t]{.25\linewidth}
\centering
\includegraphics[width=0.9\textwidth,height=3.5cm,keepaspectratio]{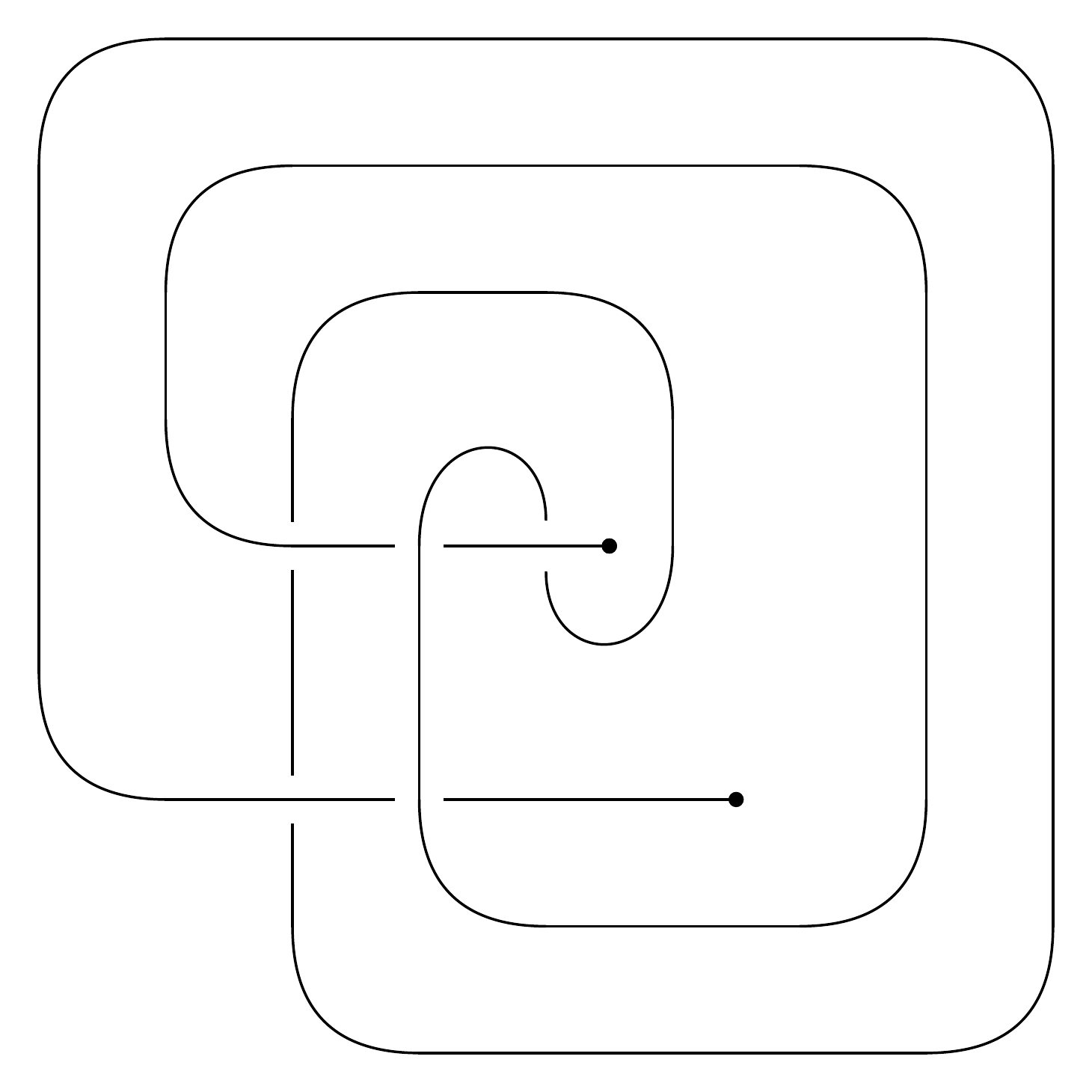}\\
\textcolor{black}{$5_{484}$}
\vspace{1cm}
\end{minipage}
\begin{minipage}[t]{.25\linewidth}
\centering
\includegraphics[width=0.9\textwidth,height=3.5cm,keepaspectratio]{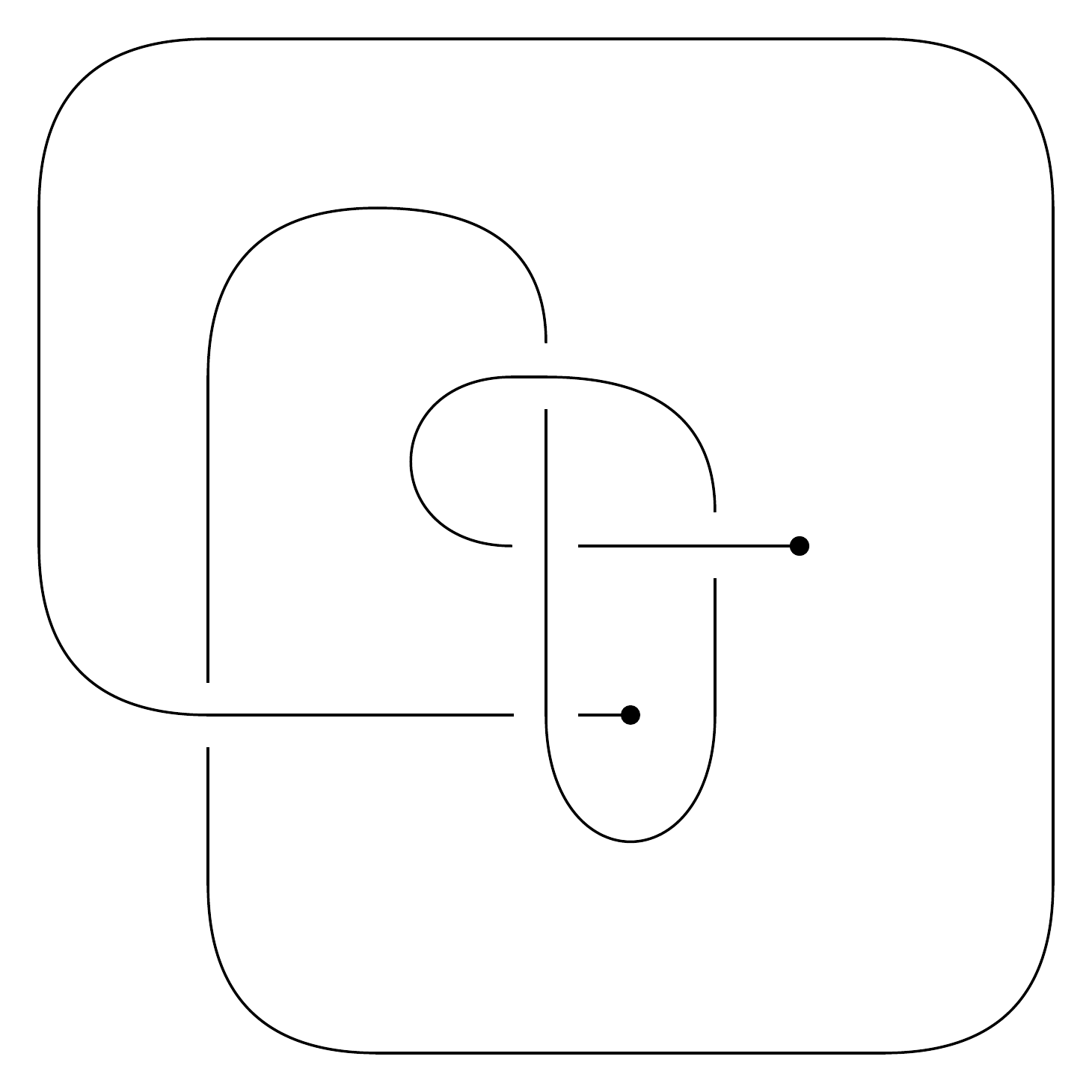}\\
\textcolor{black}{$5_{485}$}
\vspace{1cm}
\end{minipage}
\begin{minipage}[t]{.25\linewidth}
\centering
\includegraphics[width=0.9\textwidth,height=3.5cm,keepaspectratio]{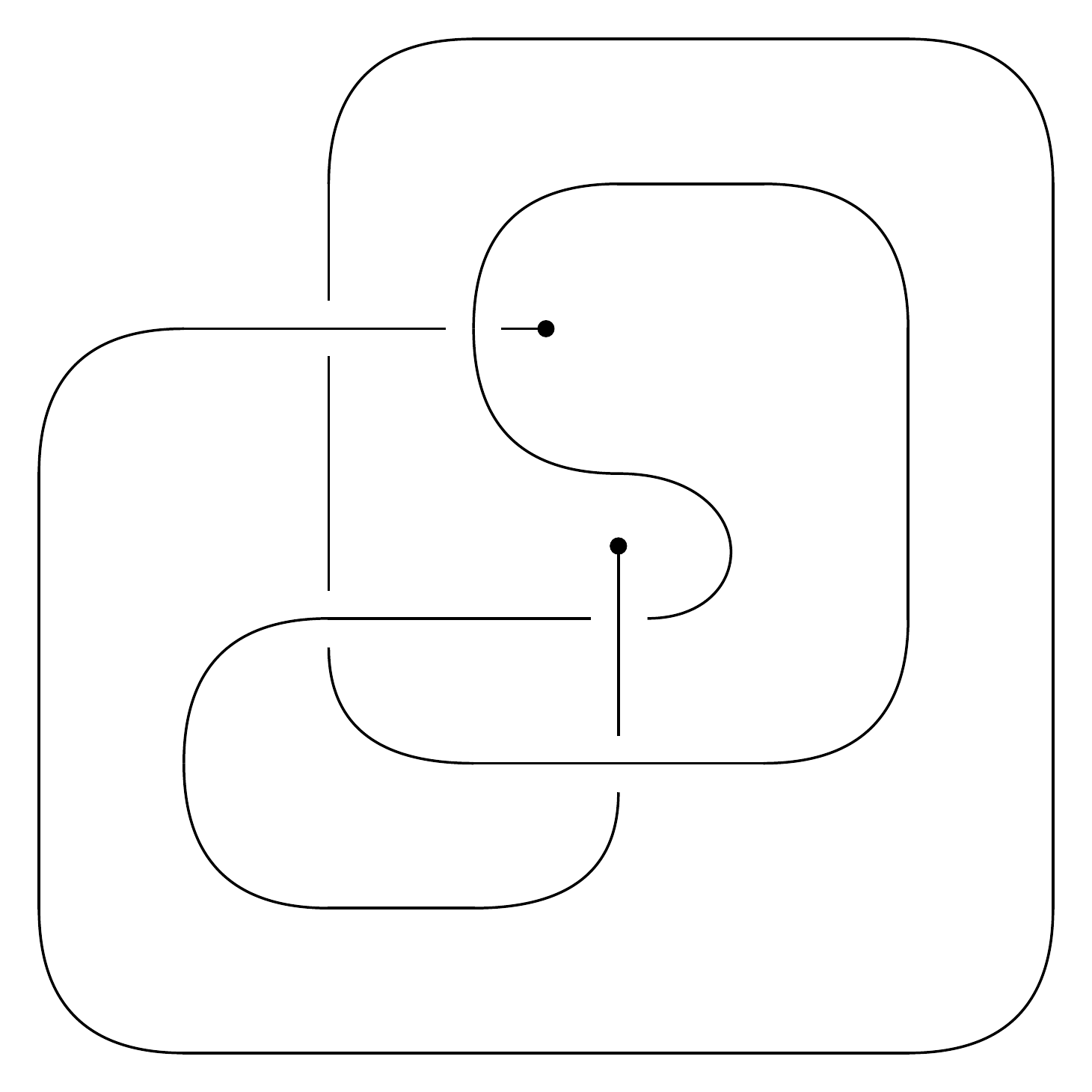}\\
\textcolor{black}{$5_{486}$}
\vspace{1cm}
\end{minipage}
\begin{minipage}[t]{.25\linewidth}
\centering
\includegraphics[width=0.9\textwidth,height=3.5cm,keepaspectratio]{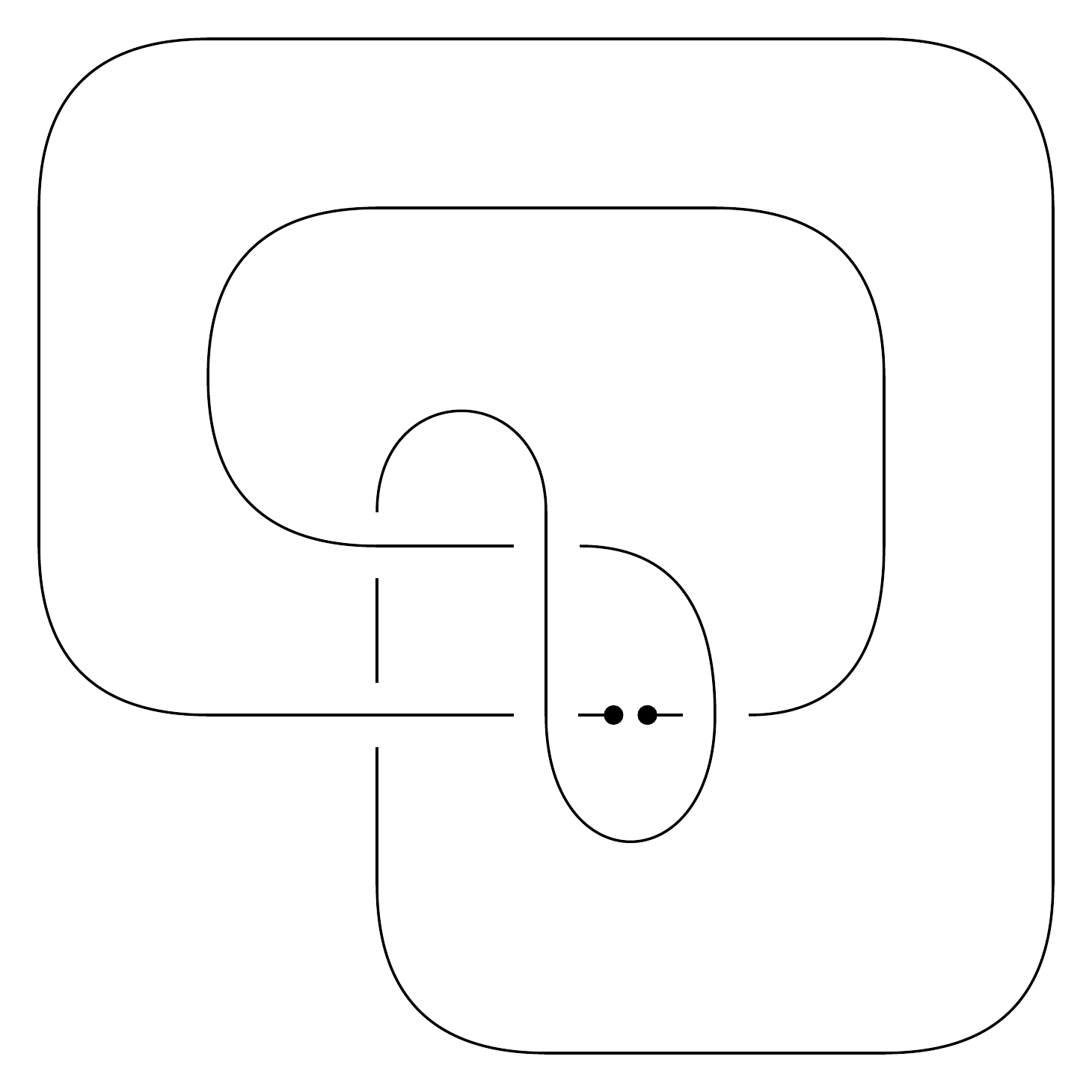}\\
\textcolor{black}{$5_{487}$}
\vspace{1cm}
\end{minipage}
\begin{minipage}[t]{.25\linewidth}
\centering
\includegraphics[width=0.9\textwidth,height=3.5cm,keepaspectratio]{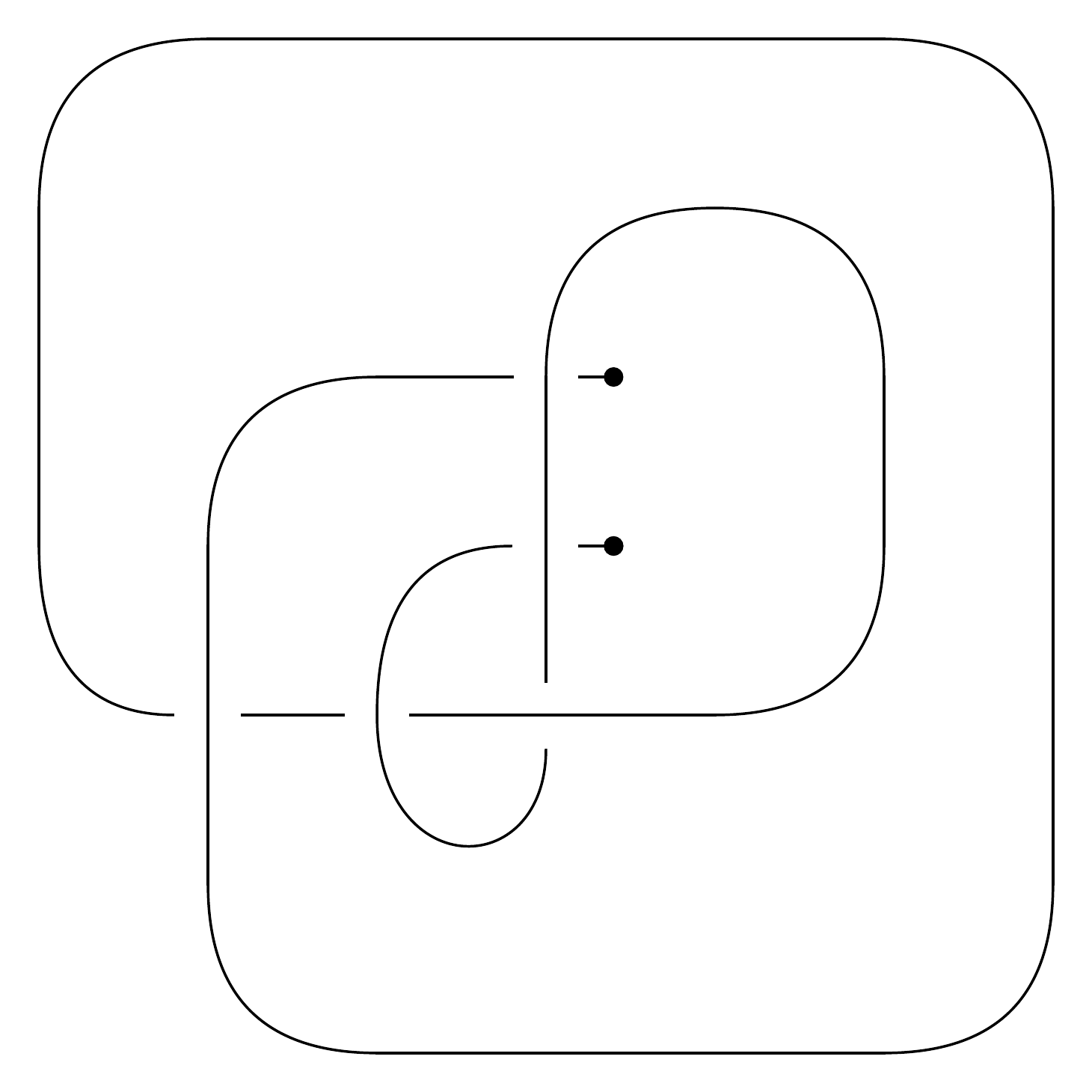}\\
\textcolor{black}{$5_{488}$}
\vspace{1cm}
\end{minipage}
\begin{minipage}[t]{.25\linewidth}
\centering
\includegraphics[width=0.9\textwidth,height=3.5cm,keepaspectratio]{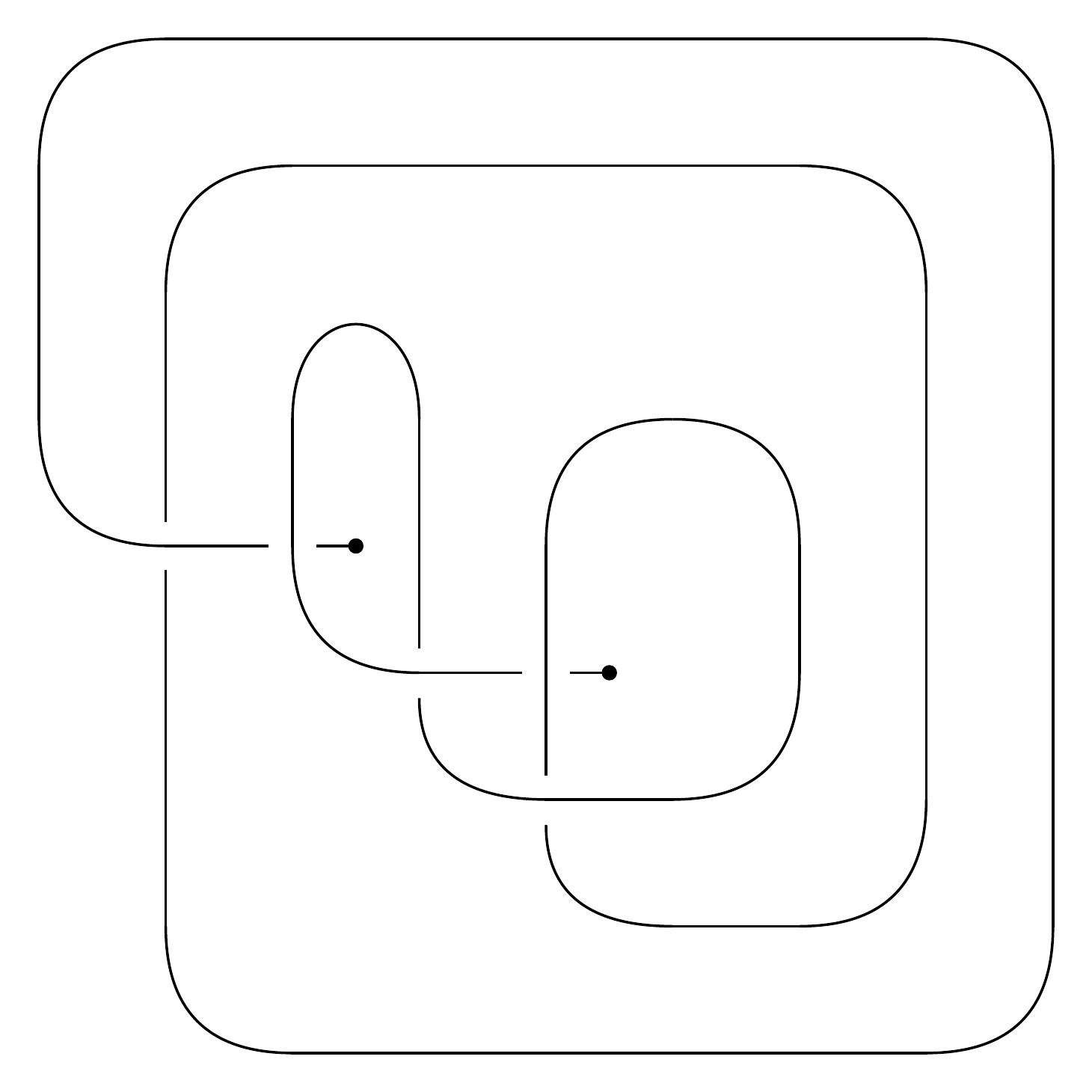}\\
\textcolor{black}{$5_{489}$}
\vspace{1cm}
\end{minipage}
\begin{minipage}[t]{.25\linewidth}
\centering
\includegraphics[width=0.9\textwidth,height=3.5cm,keepaspectratio]{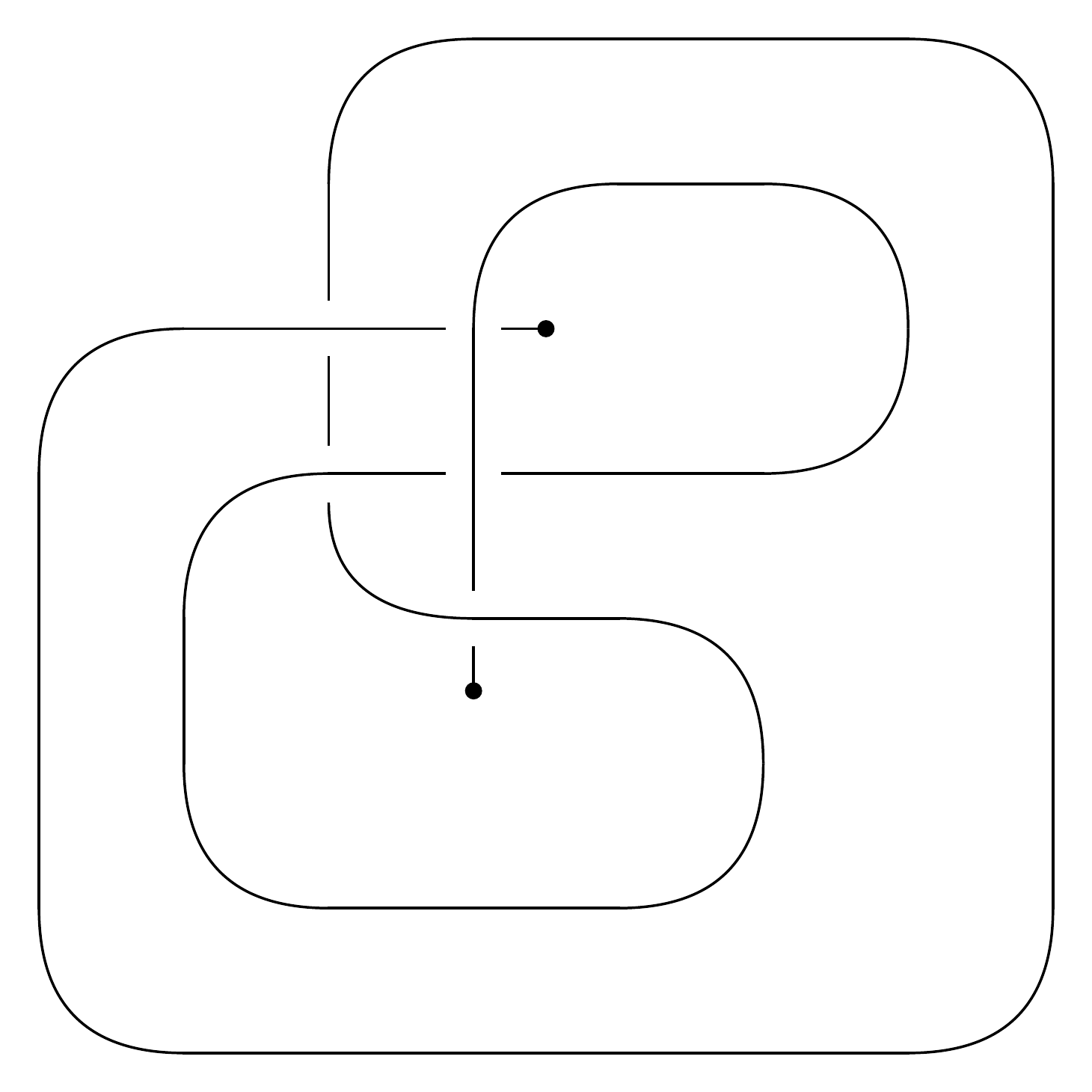}\\
\textcolor{black}{$5_{490}$}
\vspace{1cm}
\end{minipage}
\begin{minipage}[t]{.25\linewidth}
\centering
\includegraphics[width=0.9\textwidth,height=3.5cm,keepaspectratio]{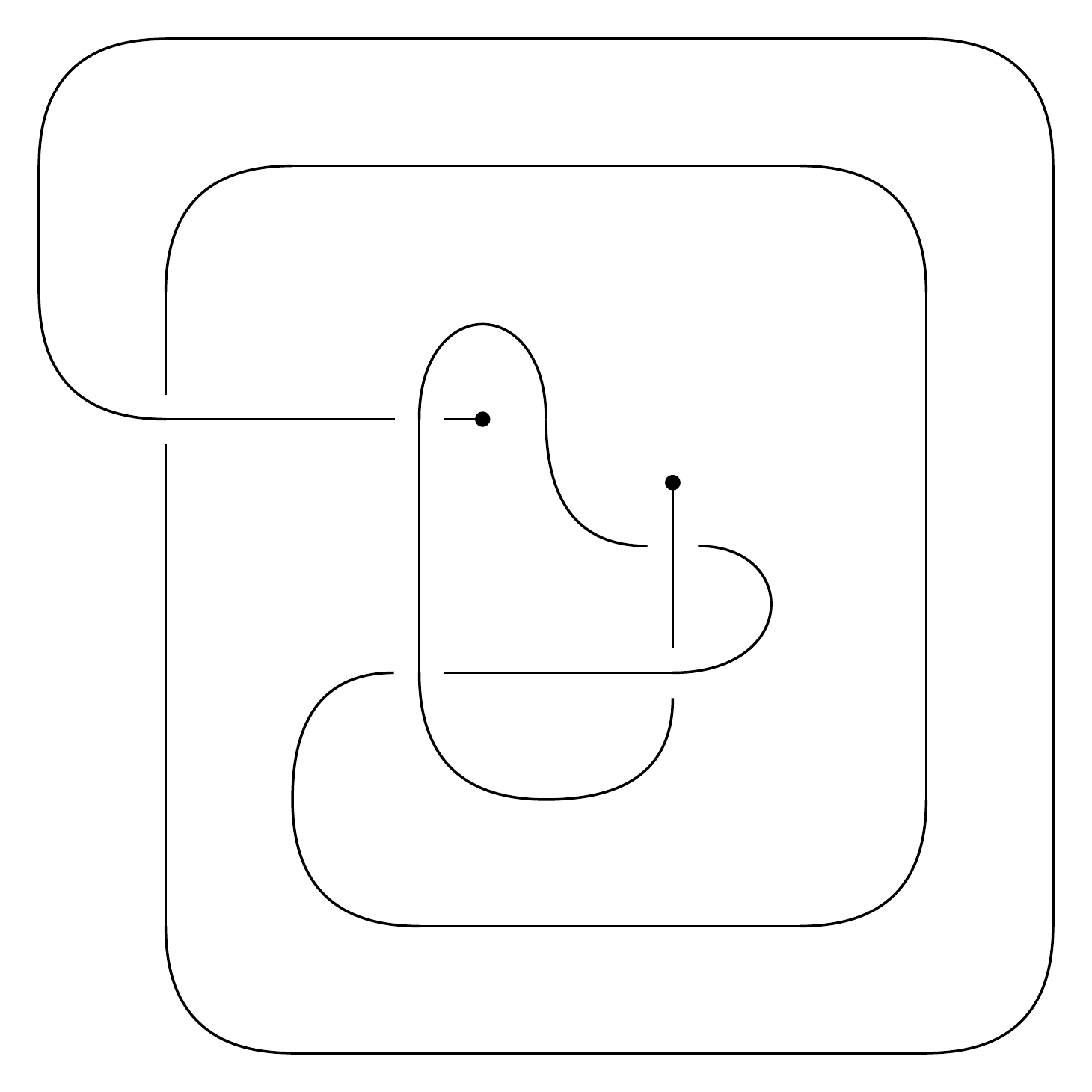}\\
\textcolor{black}{$5_{491}$}
\vspace{1cm}
\end{minipage}
\begin{minipage}[t]{.25\linewidth}
\centering
\includegraphics[width=0.9\textwidth,height=3.5cm,keepaspectratio]{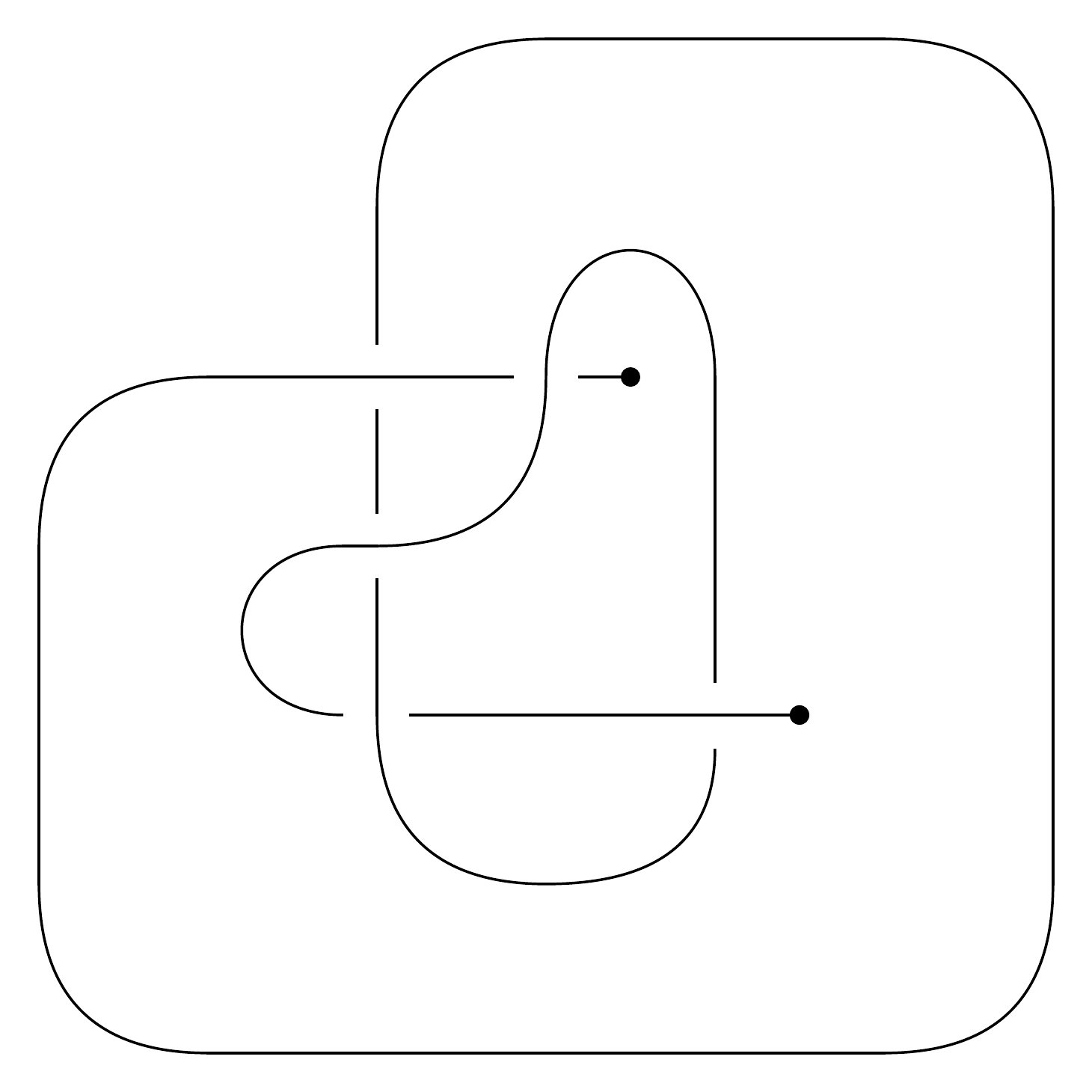}\\
\textcolor{black}{$5_{492}$}
\vspace{1cm}
\end{minipage}
\begin{minipage}[t]{.25\linewidth}
\centering
\includegraphics[width=0.9\textwidth,height=3.5cm,keepaspectratio]{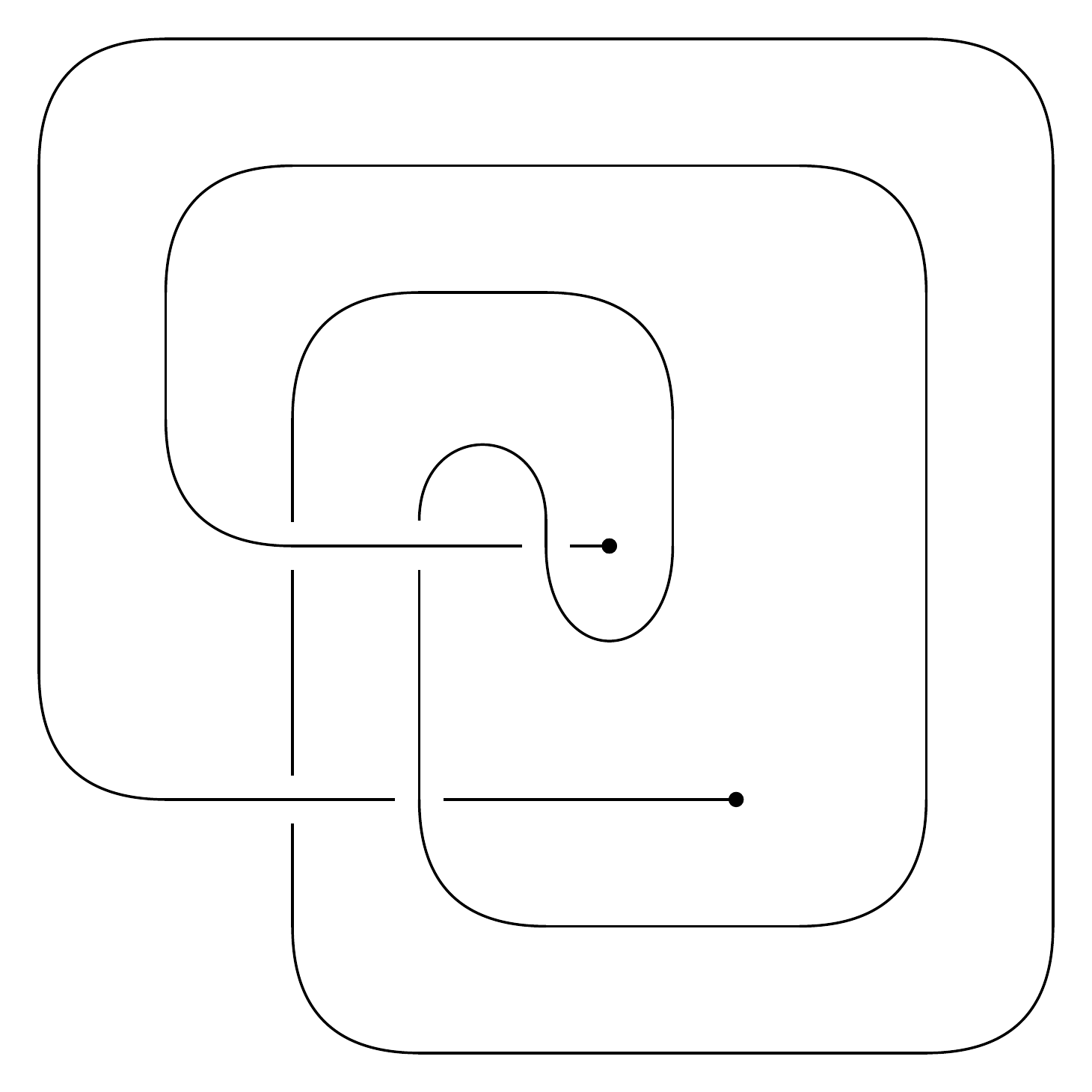}\\
\textcolor{black}{$5_{493}$}
\vspace{1cm}
\end{minipage}
\begin{minipage}[t]{.25\linewidth}
\centering
\includegraphics[width=0.9\textwidth,height=3.5cm,keepaspectratio]{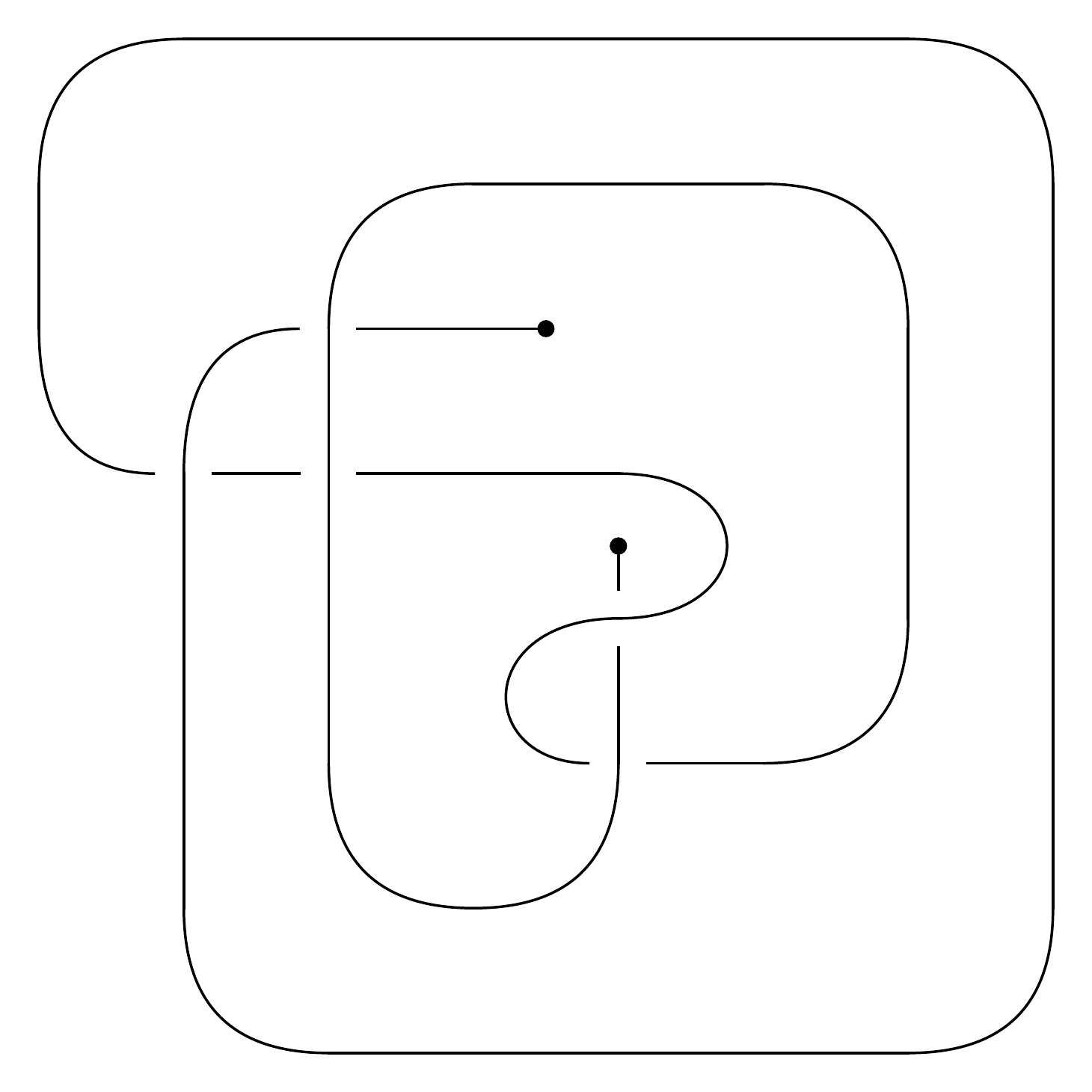}\\
\textcolor{black}{$5_{494}$}
\vspace{1cm}
\end{minipage}
\begin{minipage}[t]{.25\linewidth}
\centering
\includegraphics[width=0.9\textwidth,height=3.5cm,keepaspectratio]{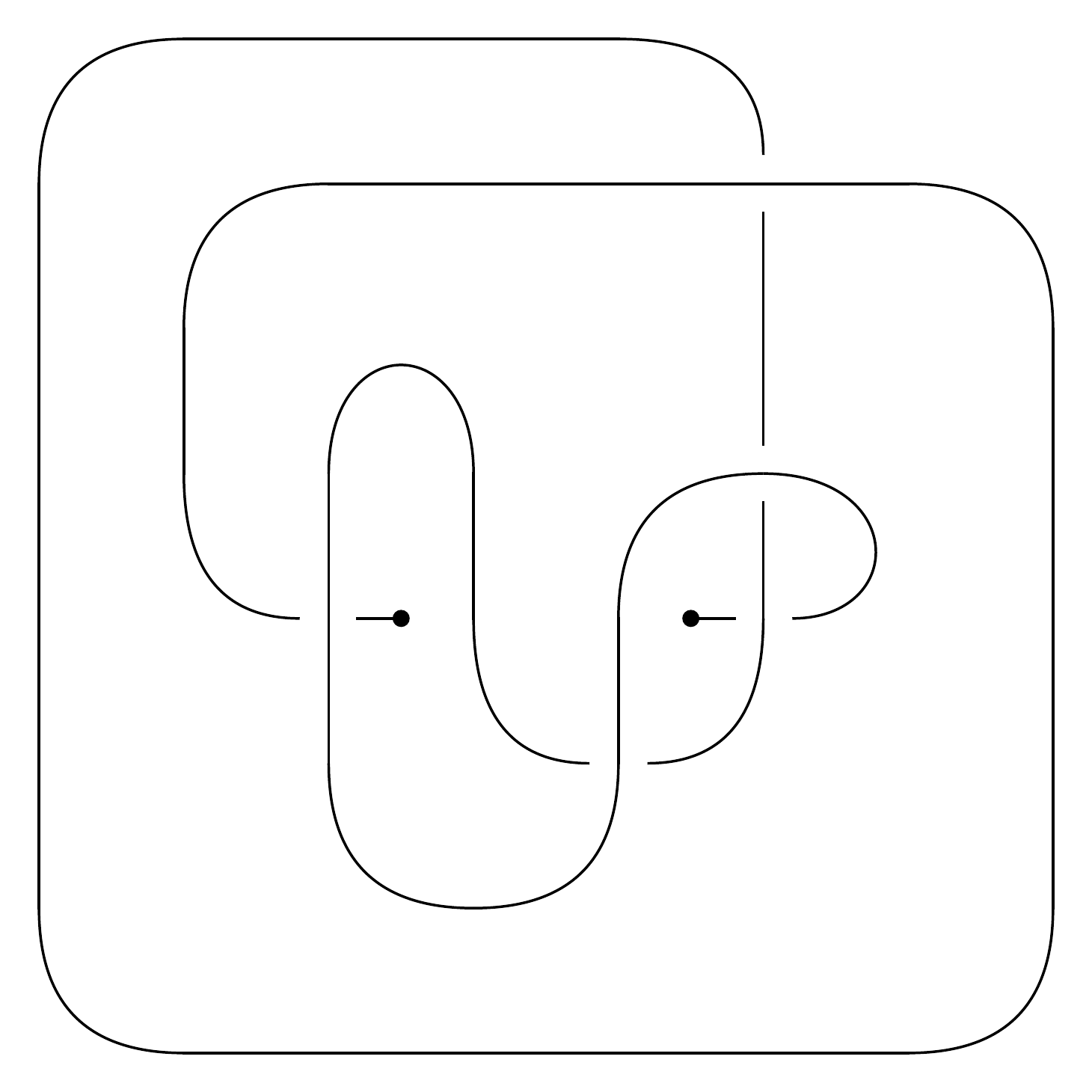}\\
\textcolor{black}{$5_{495}$}
\vspace{1cm}
\end{minipage}
\begin{minipage}[t]{.25\linewidth}
\centering
\includegraphics[width=0.9\textwidth,height=3.5cm,keepaspectratio]{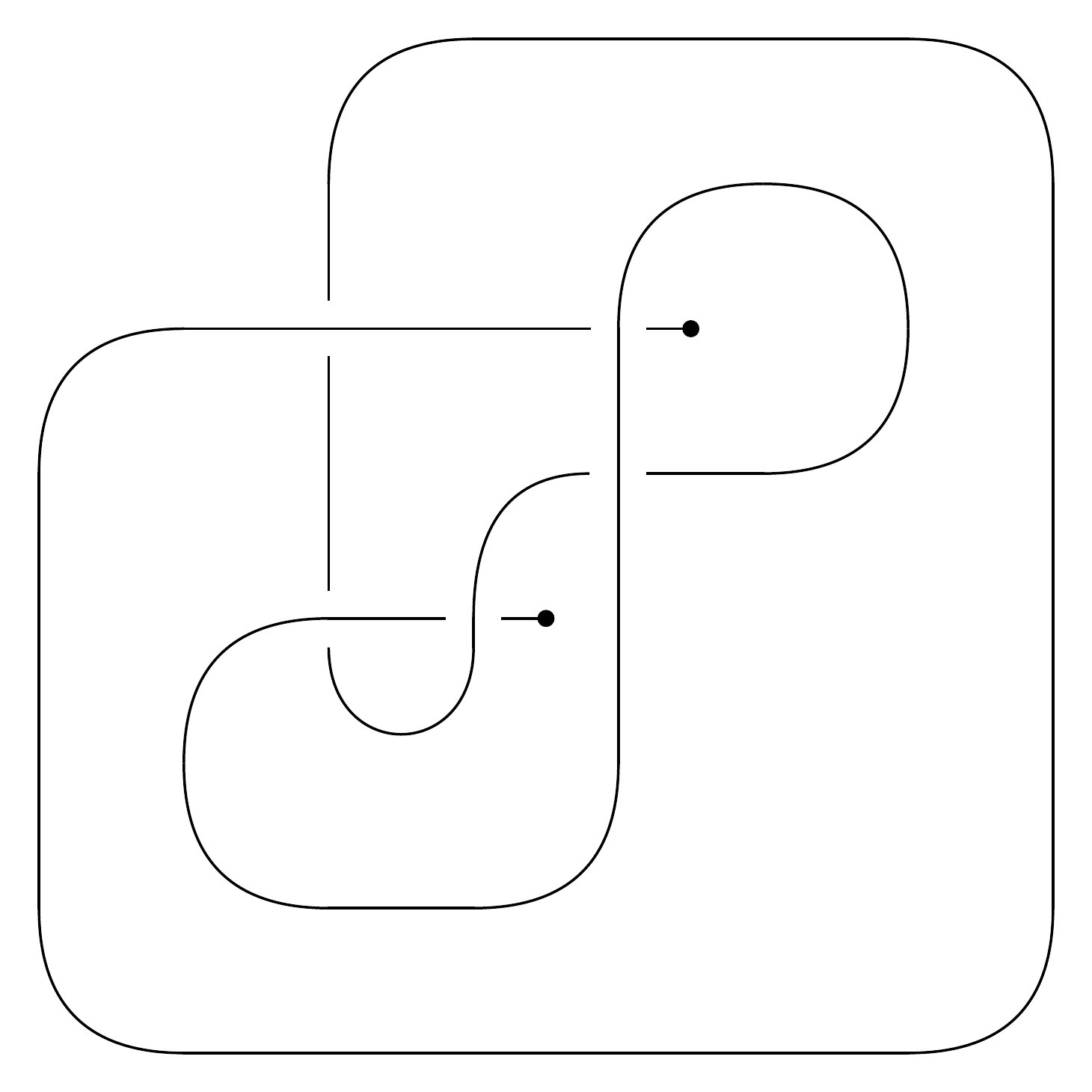}\\
\textcolor{black}{$5_{496}$}
\vspace{1cm}
\end{minipage}
\begin{minipage}[t]{.25\linewidth}
\centering
\includegraphics[width=0.9\textwidth,height=3.5cm,keepaspectratio]{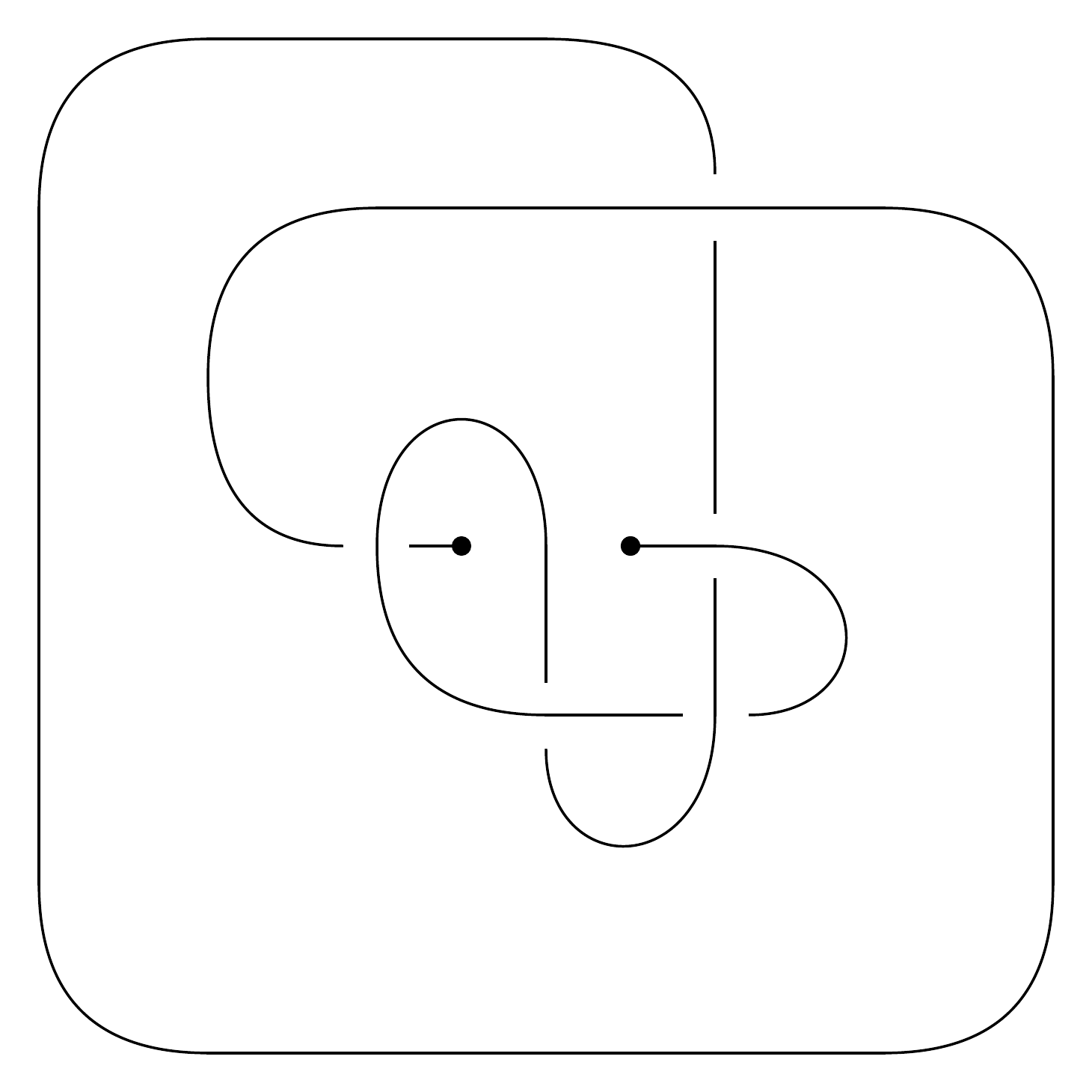}\\
\textcolor{black}{$5_{497}$}
\vspace{1cm}
\end{minipage}
\begin{minipage}[t]{.25\linewidth}
\centering
\includegraphics[width=0.9\textwidth,height=3.5cm,keepaspectratio]{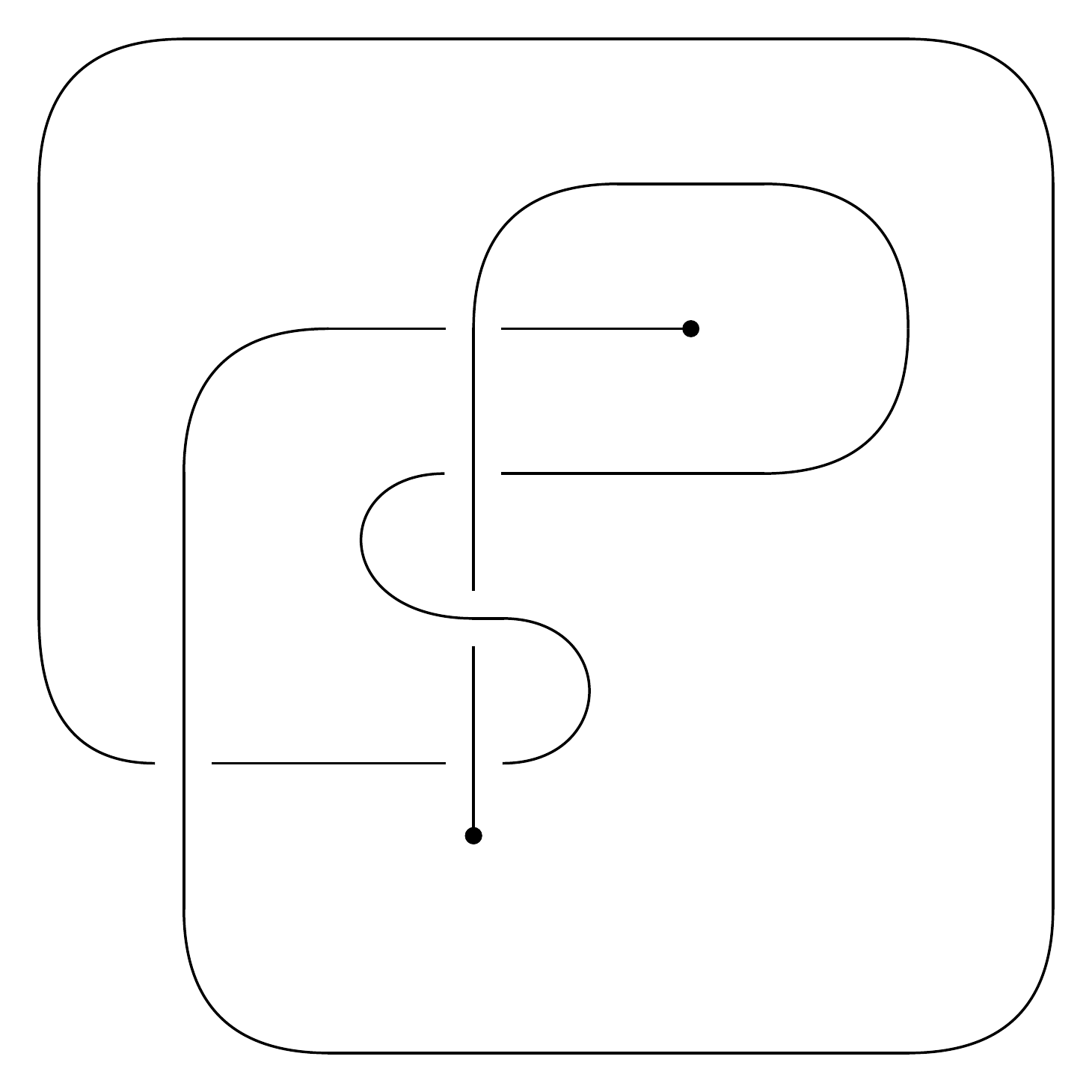}\\
\textcolor{black}{$5_{498}$}
\vspace{1cm}
\end{minipage}
\begin{minipage}[t]{.25\linewidth}
\centering
\includegraphics[width=0.9\textwidth,height=3.5cm,keepaspectratio]{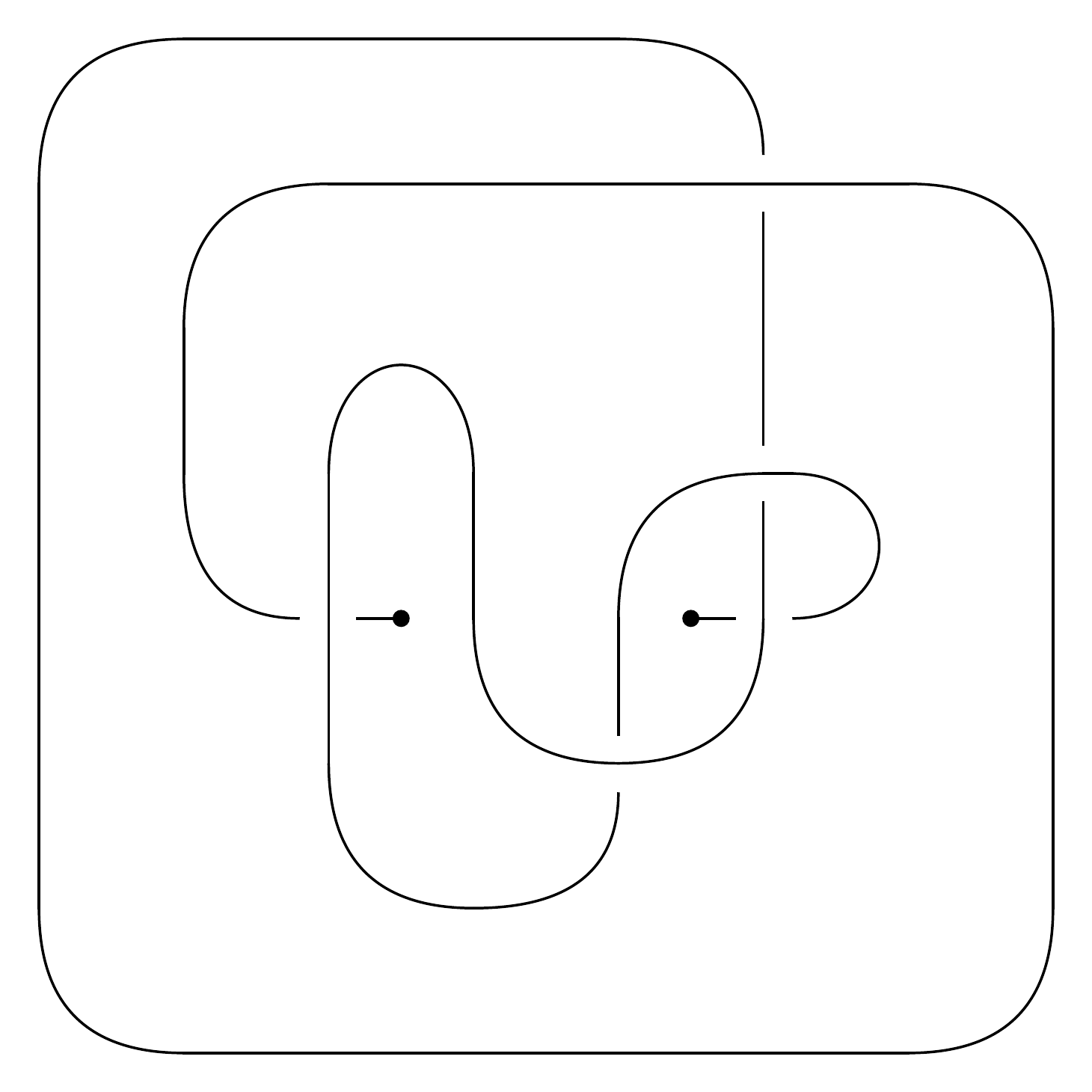}\\
\textcolor{black}{$5_{499}$}
\vspace{1cm}
\end{minipage}
\begin{minipage}[t]{.25\linewidth}
\centering
\includegraphics[width=0.9\textwidth,height=3.5cm,keepaspectratio]{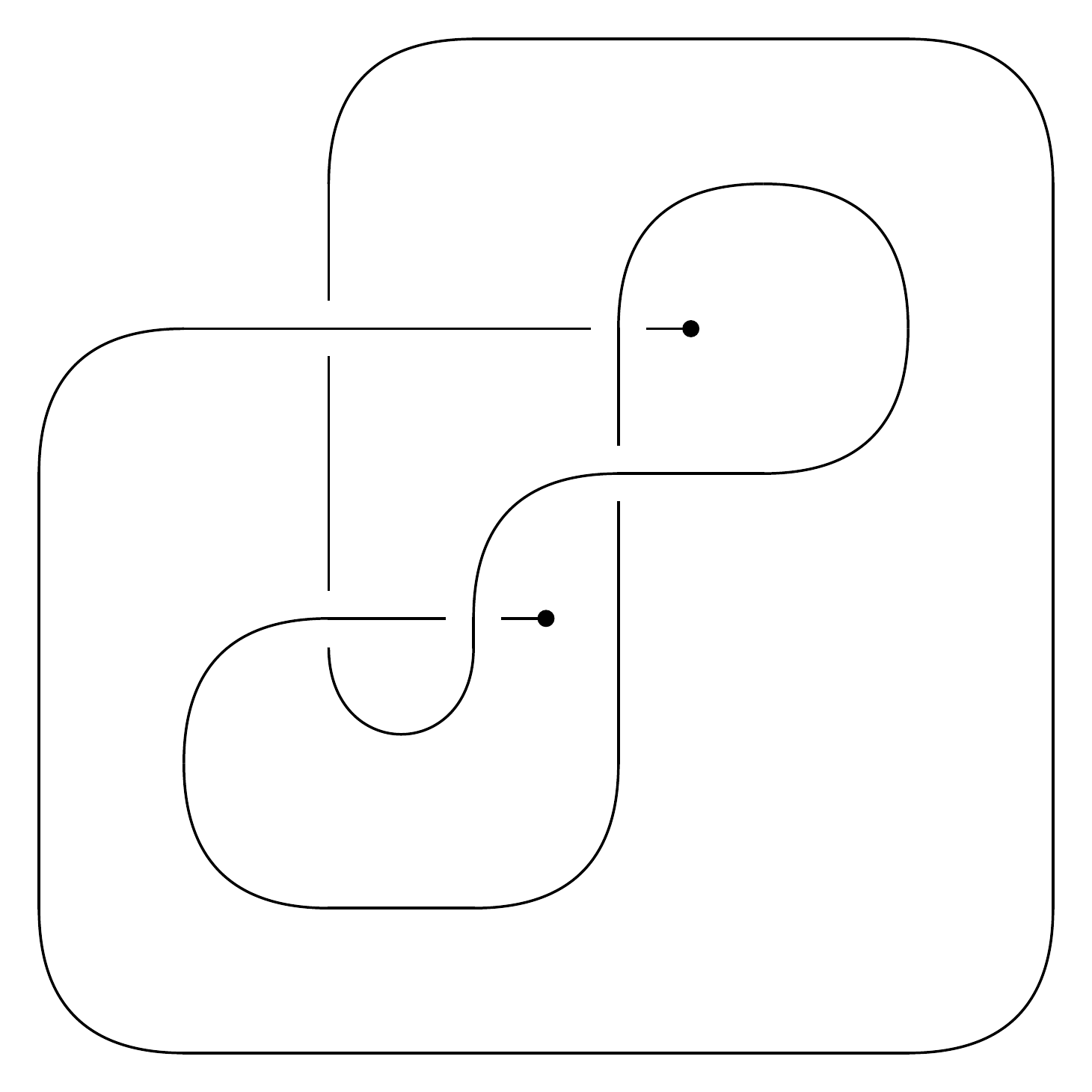}\\
\textcolor{black}{$5_{500}$}
\vspace{1cm}
\end{minipage}
\begin{minipage}[t]{.25\linewidth}
\centering
\includegraphics[width=0.9\textwidth,height=3.5cm,keepaspectratio]{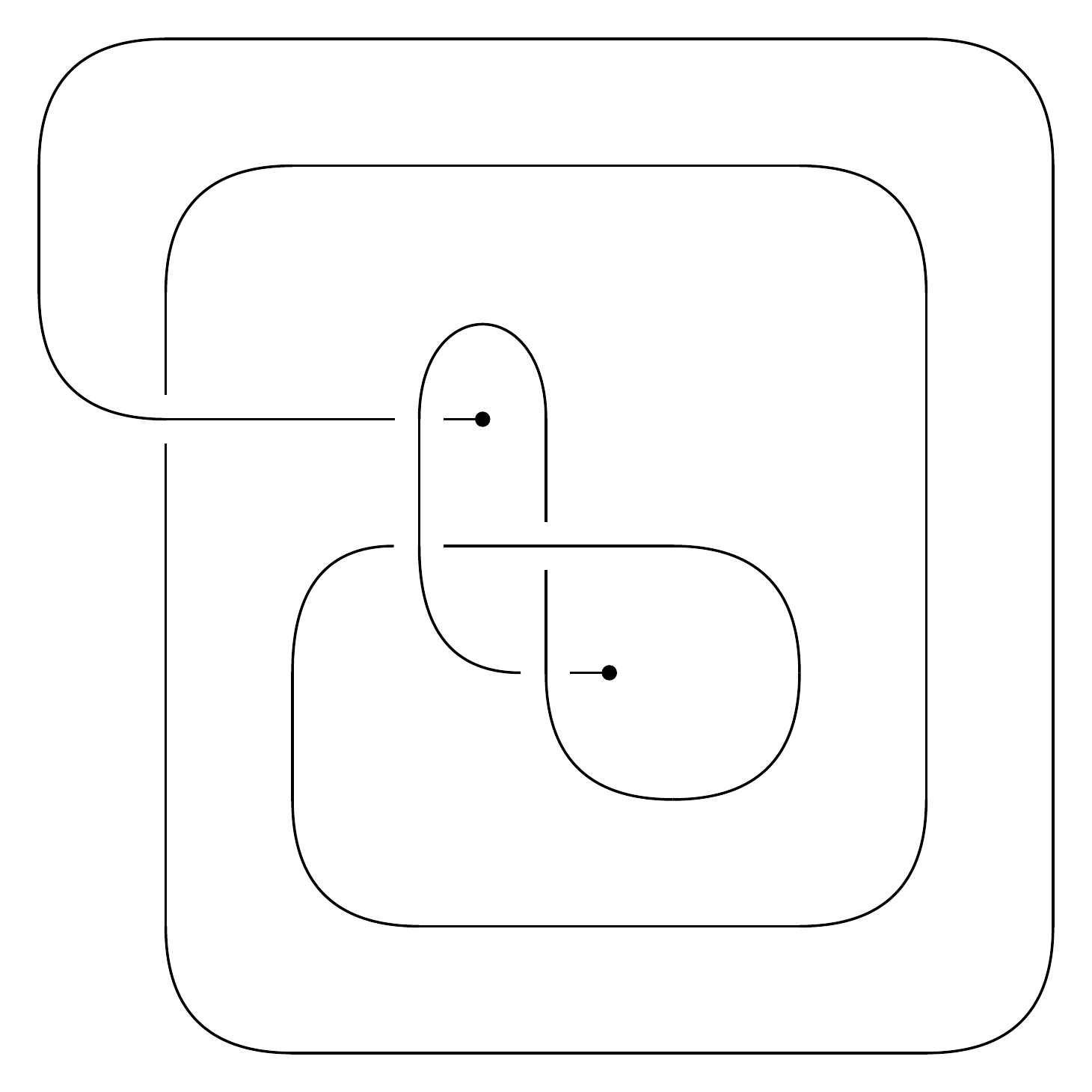}\\
\textcolor{black}{$5_{501}$}
\vspace{1cm}
\end{minipage}
\begin{minipage}[t]{.25\linewidth}
\centering
\includegraphics[width=0.9\textwidth,height=3.5cm,keepaspectratio]{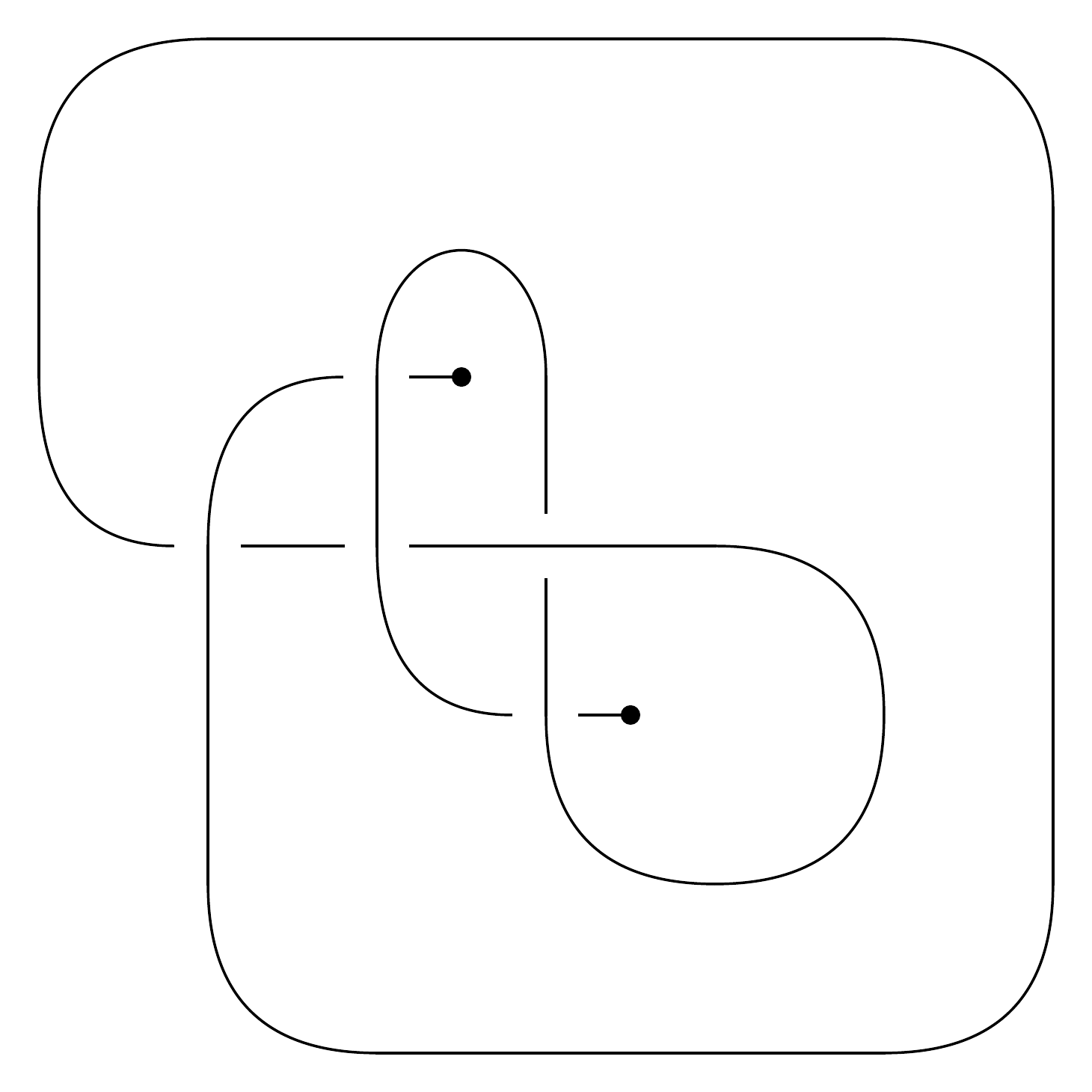}\\
\textcolor{black}{$5_{502}$}
\vspace{1cm}
\end{minipage}
\begin{minipage}[t]{.25\linewidth}
\centering
\includegraphics[width=0.9\textwidth,height=3.5cm,keepaspectratio]{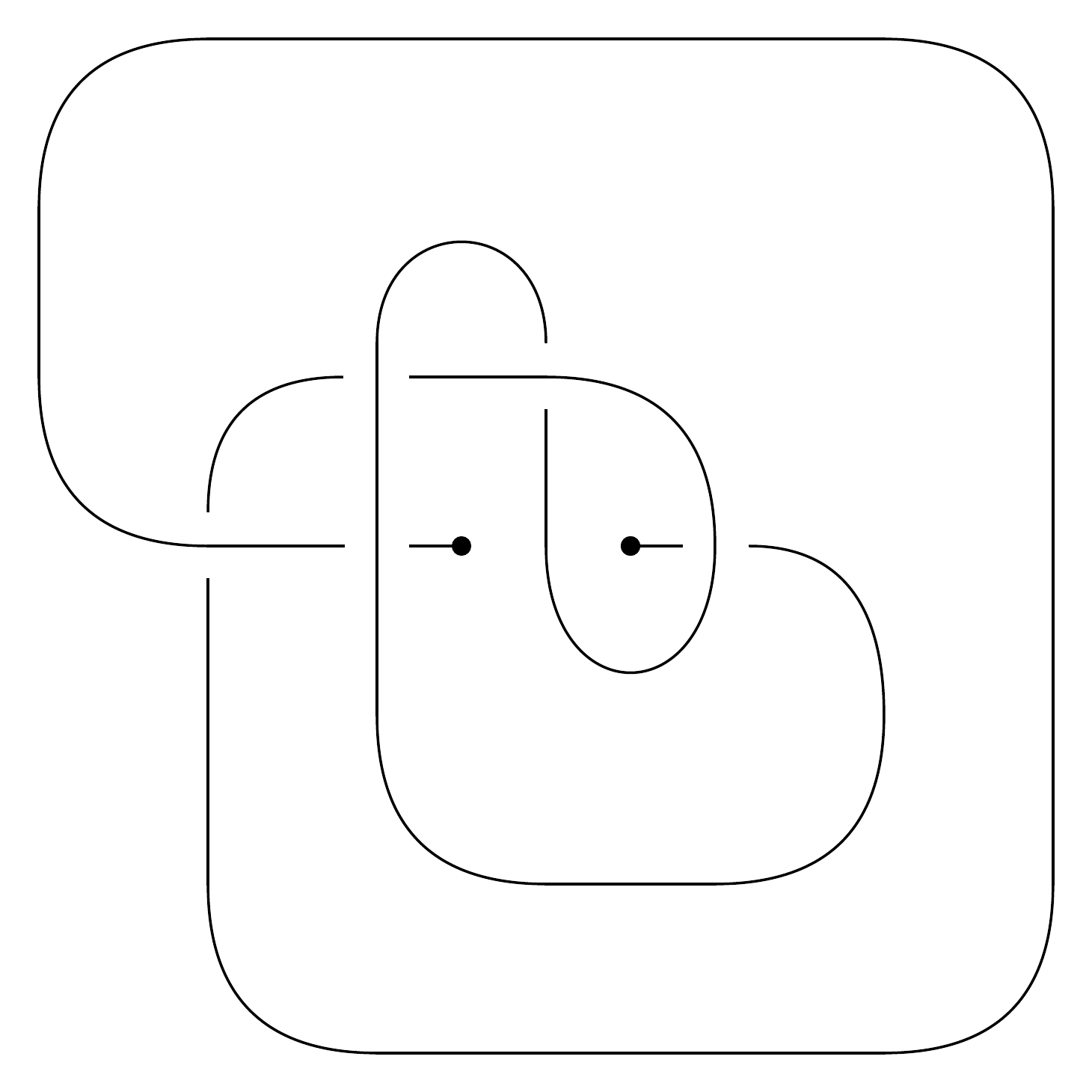}\\
\textcolor{black}{$5_{503}$}
\vspace{1cm}
\end{minipage}
\begin{minipage}[t]{.25\linewidth}
\centering
\includegraphics[width=0.9\textwidth,height=3.5cm,keepaspectratio]{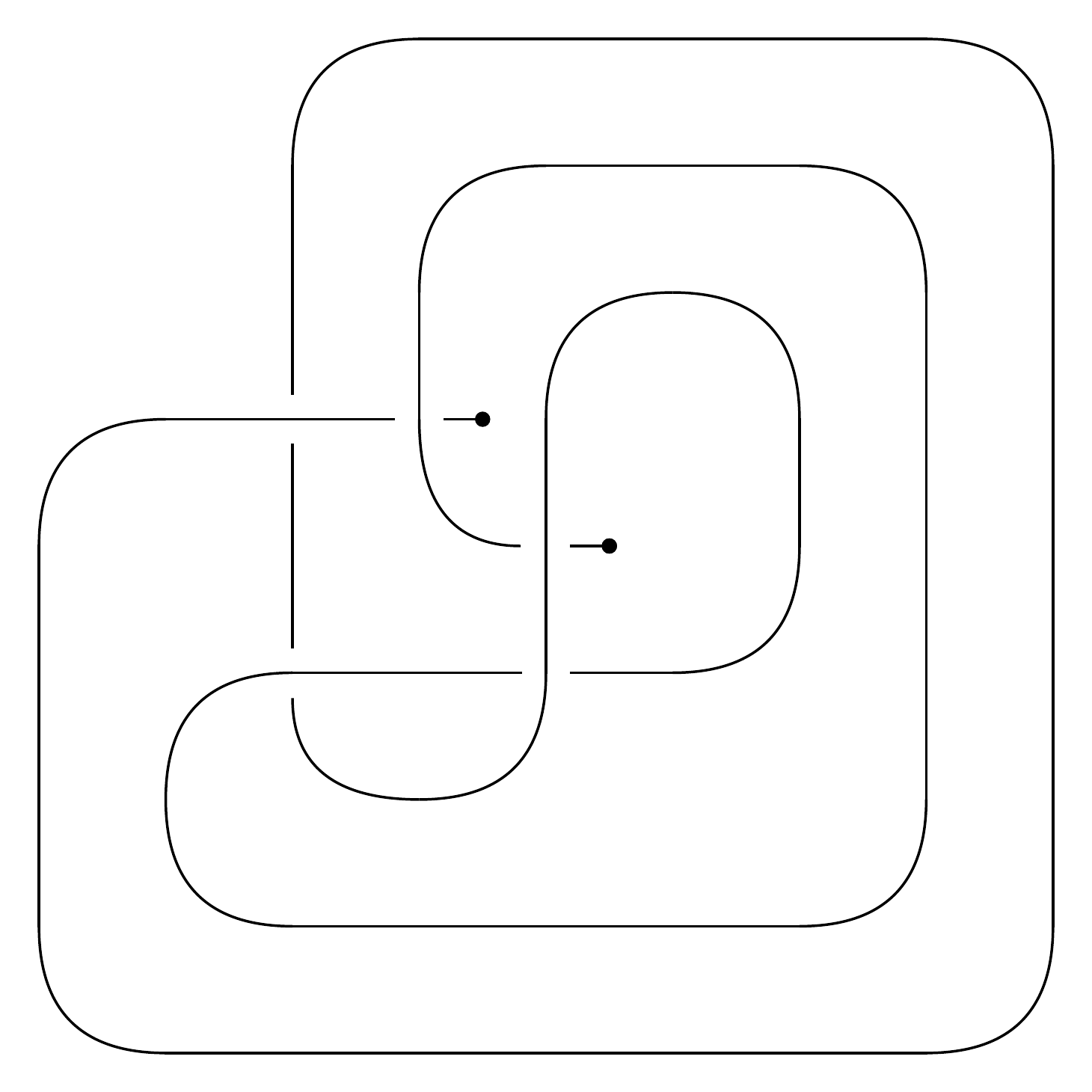}\\
\textcolor{black}{$5_{504}$}
\vspace{1cm}
\end{minipage}
\begin{minipage}[t]{.25\linewidth}
\centering
\includegraphics[width=0.9\textwidth,height=3.5cm,keepaspectratio]{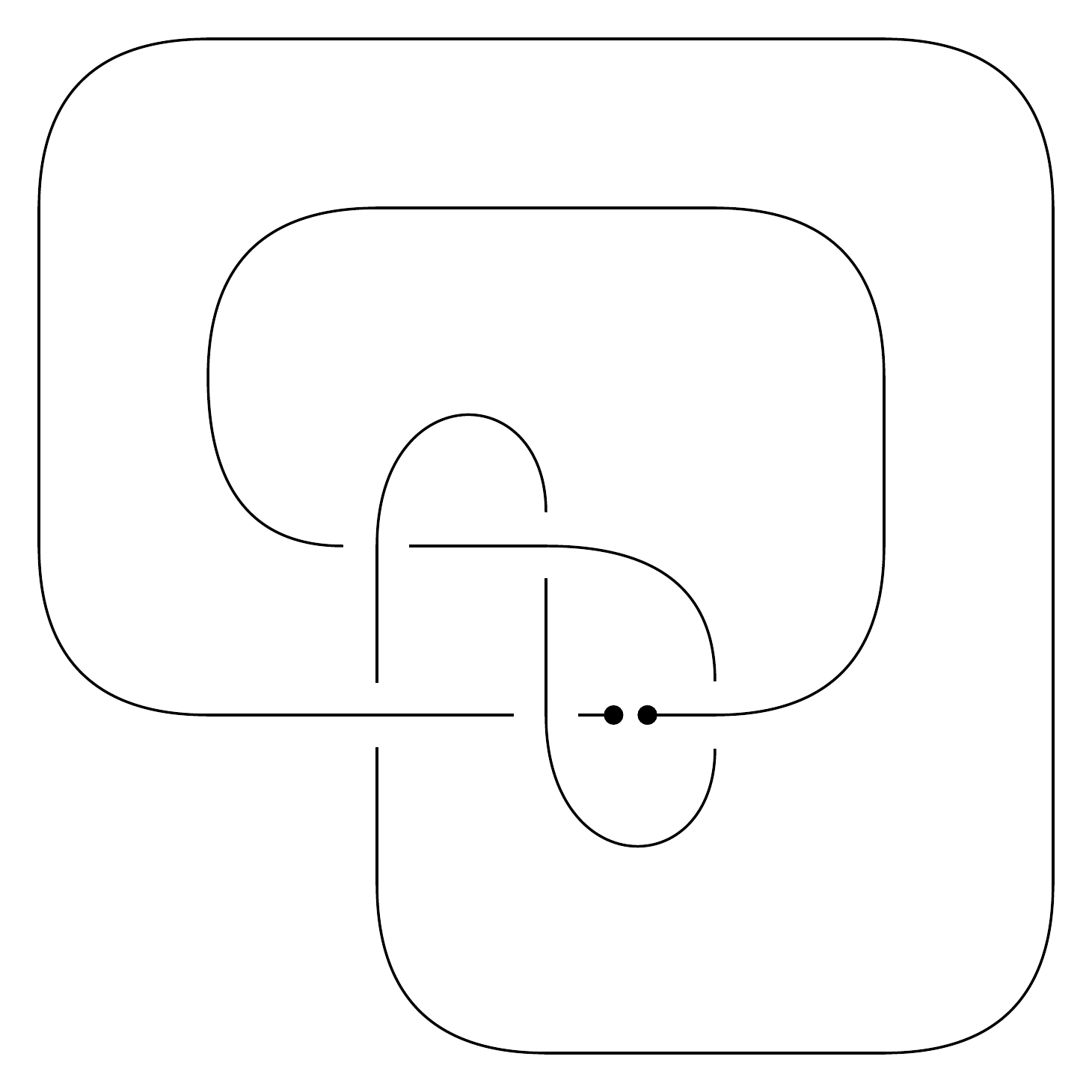}\\
\textcolor{black}{$5_{505}$}
\vspace{1cm}
\end{minipage}
\begin{minipage}[t]{.25\linewidth}
\centering
\includegraphics[width=0.9\textwidth,height=3.5cm,keepaspectratio]{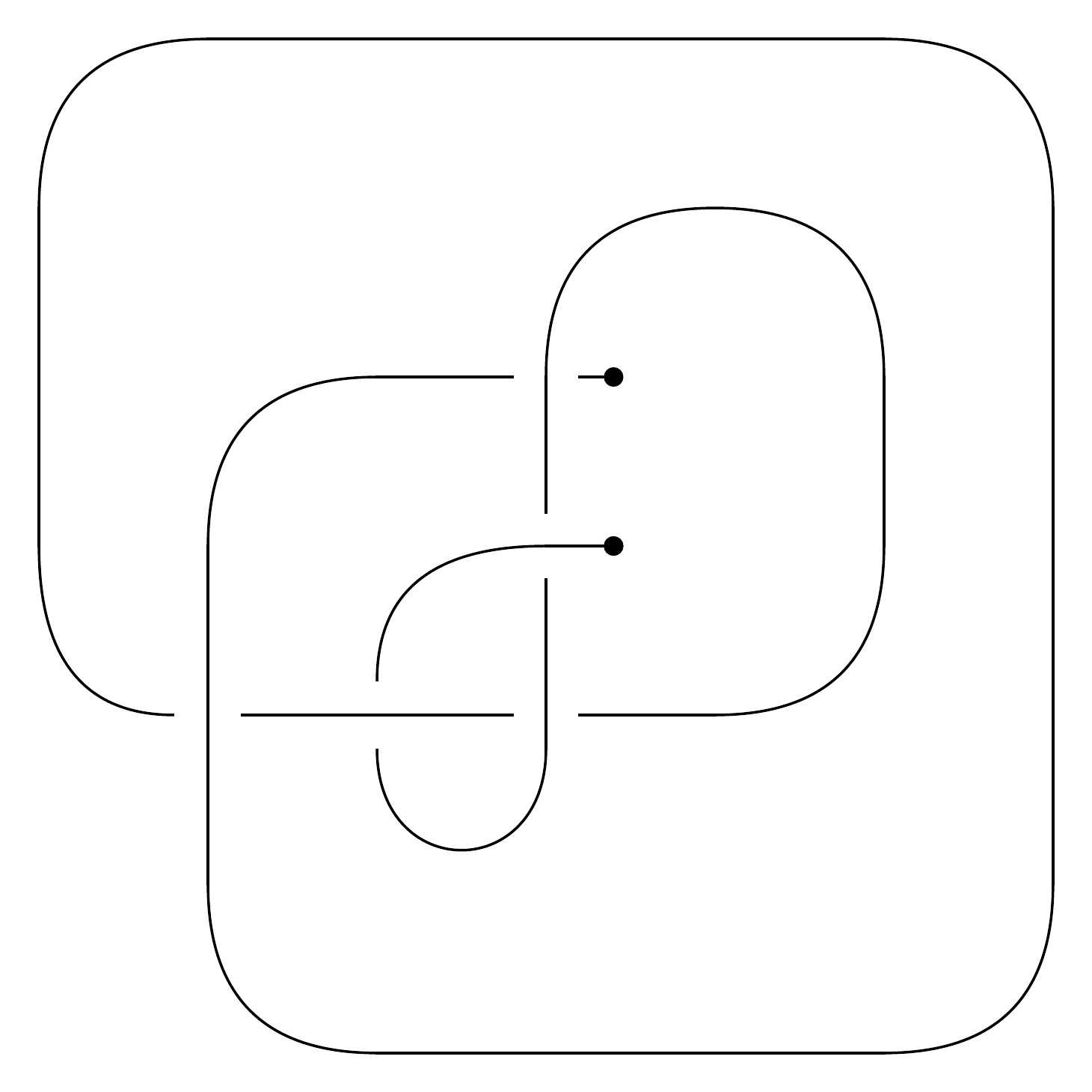}\\
\textcolor{black}{$5_{506}$}
\vspace{1cm}
\end{minipage}
\begin{minipage}[t]{.25\linewidth}
\centering
\includegraphics[width=0.9\textwidth,height=3.5cm,keepaspectratio]{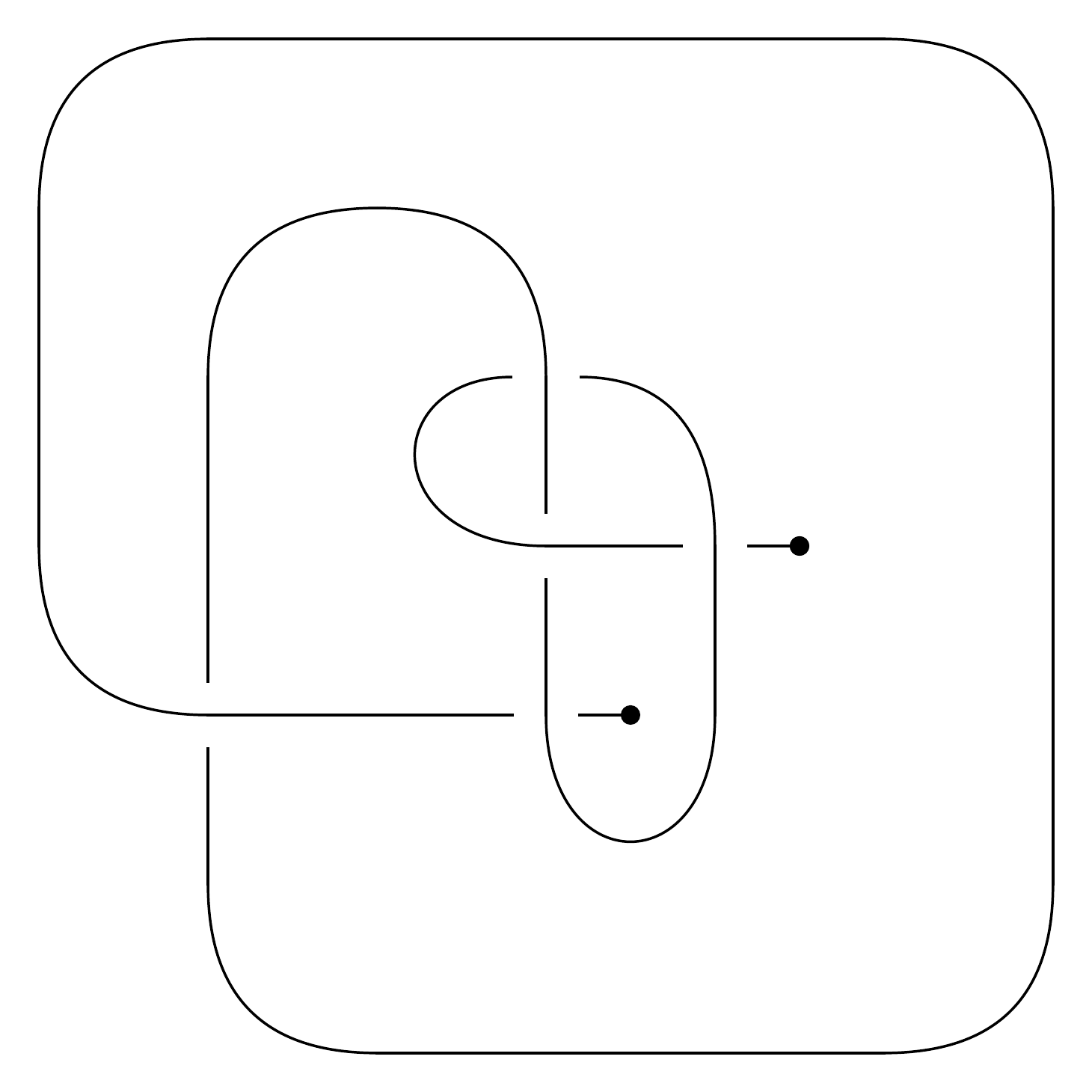}\\
\textcolor{black}{$5_{507}$}
\vspace{1cm}
\end{minipage}
\begin{minipage}[t]{.25\linewidth}
\centering
\includegraphics[width=0.9\textwidth,height=3.5cm,keepaspectratio]{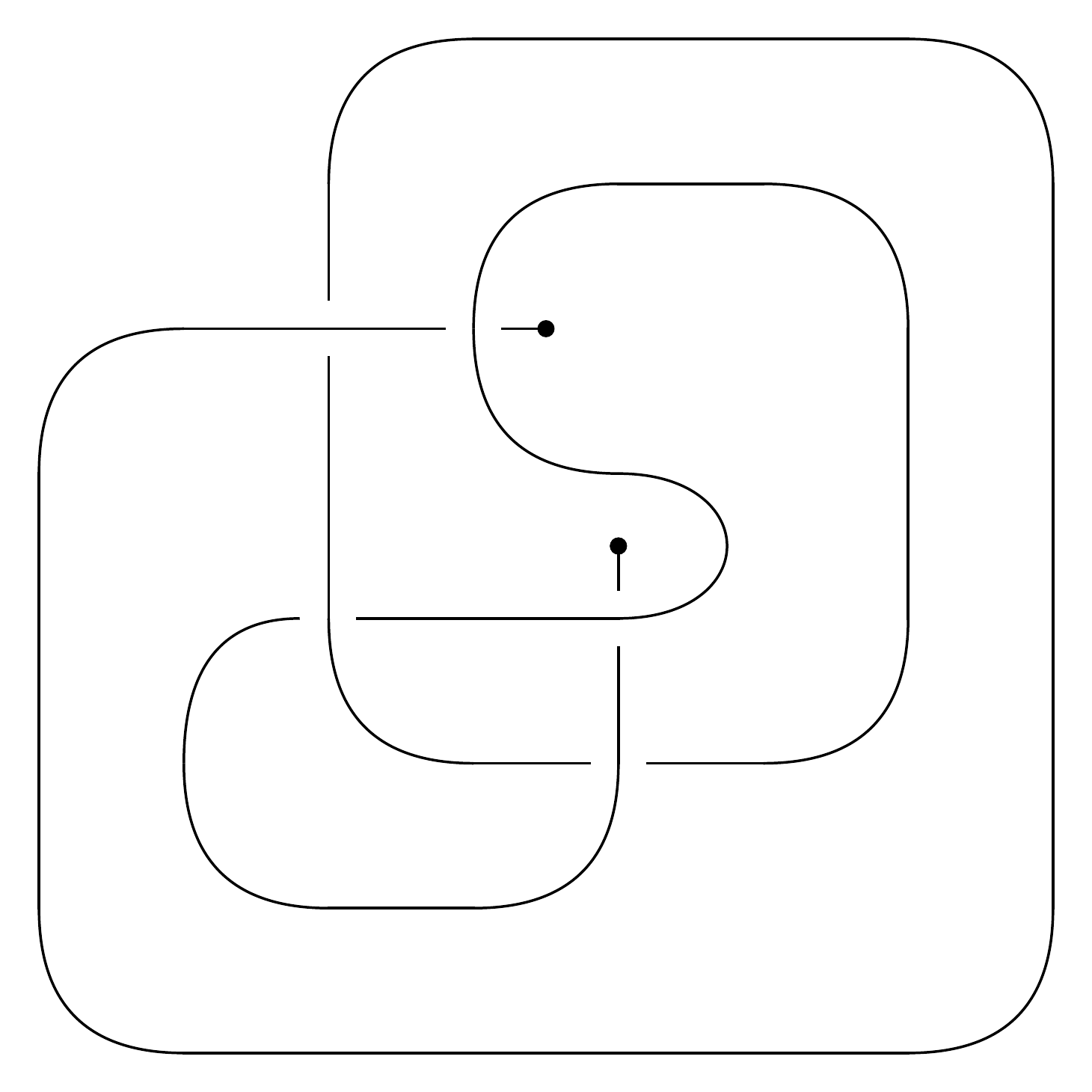}\\
\textcolor{black}{$5_{508}$}
\vspace{1cm}
\end{minipage}
\begin{minipage}[t]{.25\linewidth}
\centering
\includegraphics[width=0.9\textwidth,height=3.5cm,keepaspectratio]{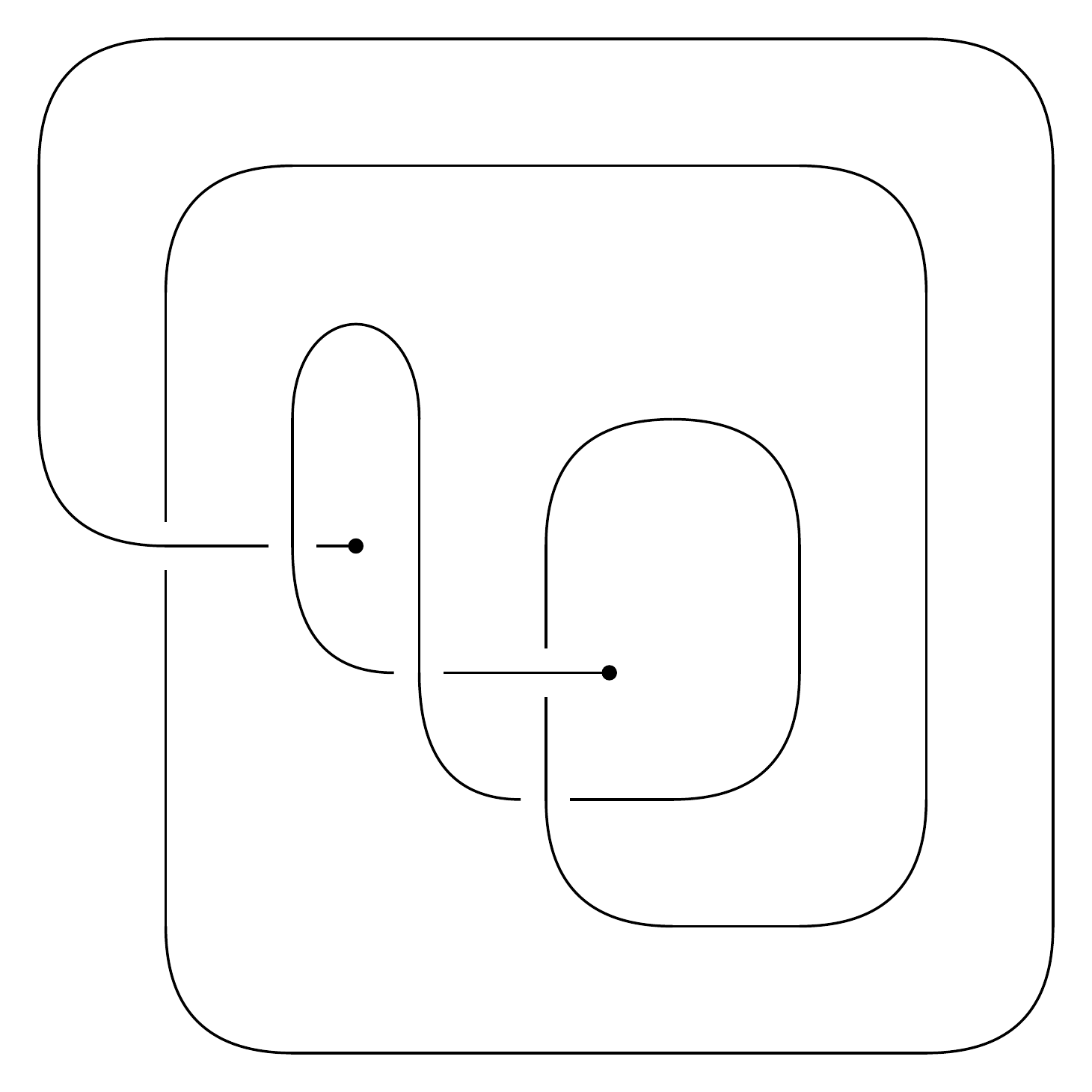}\\
\textcolor{black}{$5_{509}$}
\vspace{1cm}
\end{minipage}
\begin{minipage}[t]{.25\linewidth}
\centering
\includegraphics[width=0.9\textwidth,height=3.5cm,keepaspectratio]{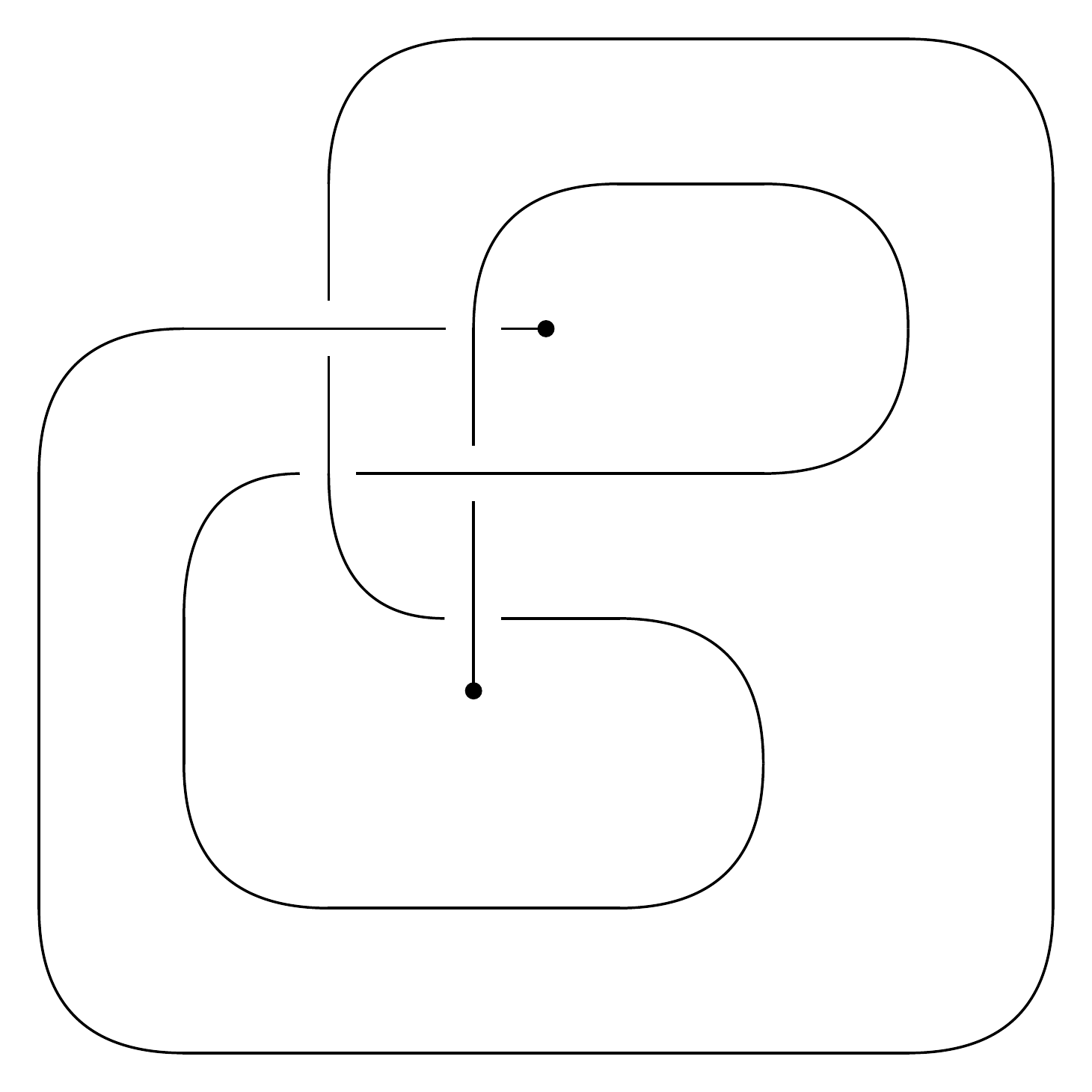}\\
\textcolor{black}{$5_{510}$}
\vspace{1cm}
\end{minipage}
\begin{minipage}[t]{.25\linewidth}
\centering
\includegraphics[width=0.9\textwidth,height=3.5cm,keepaspectratio]{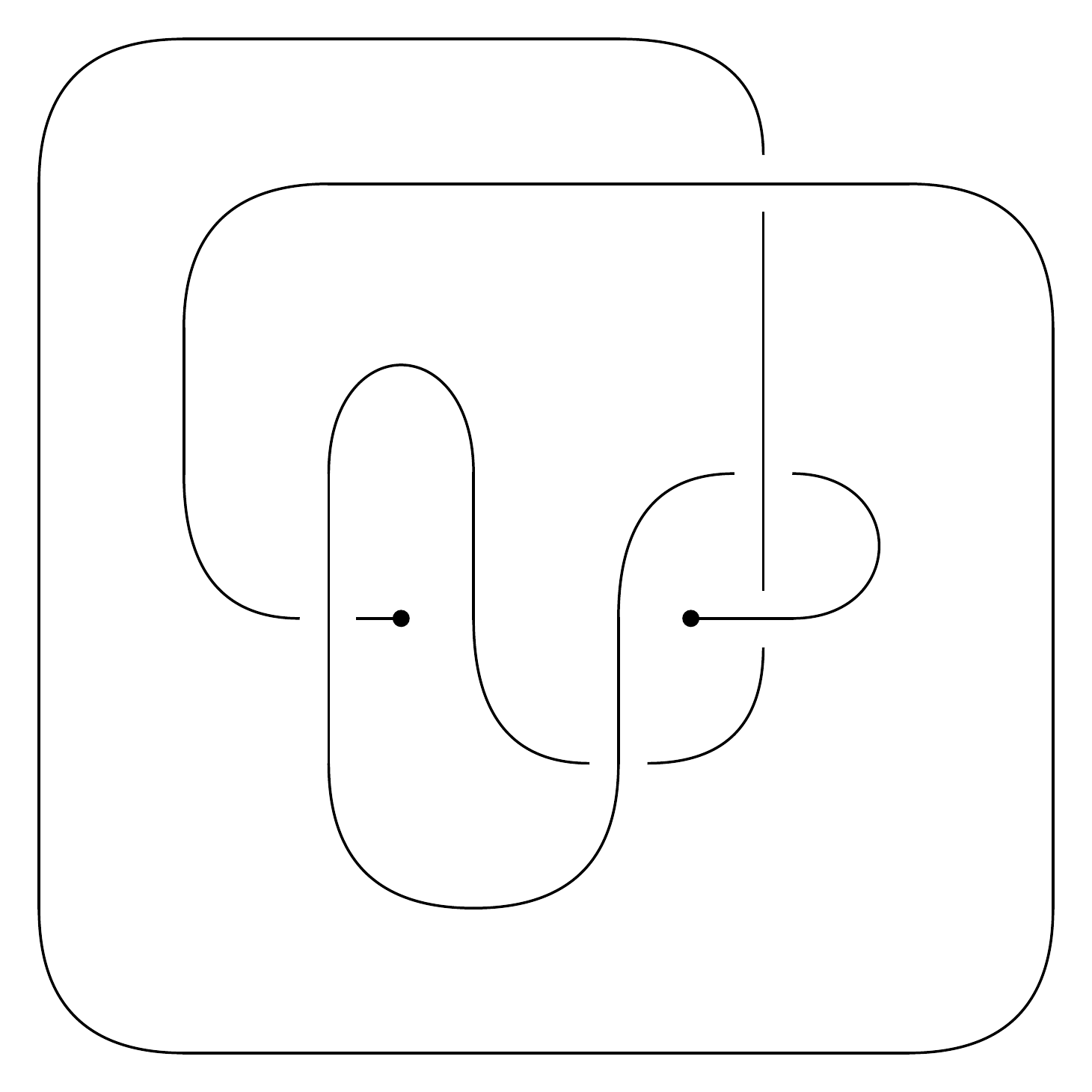}\\
\textcolor{black}{$5_{511}$}
\vspace{1cm}
\end{minipage}
\begin{minipage}[t]{.25\linewidth}
\centering
\includegraphics[width=0.9\textwidth,height=3.5cm,keepaspectratio]{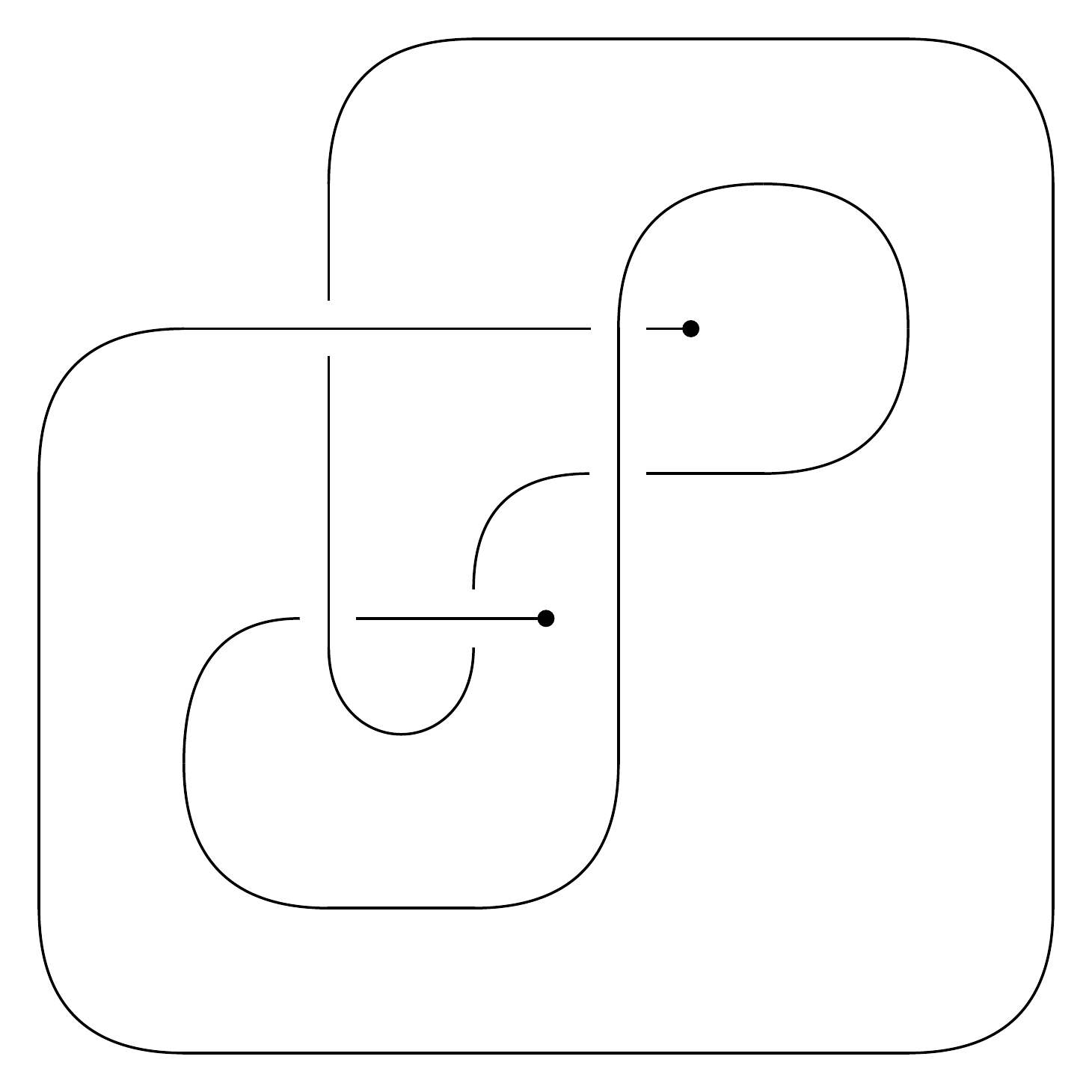}\\
\textcolor{black}{$5_{512}$}
\vspace{1cm}
\end{minipage}
\begin{minipage}[t]{.25\linewidth}
\centering
\includegraphics[width=0.9\textwidth,height=3.5cm,keepaspectratio]{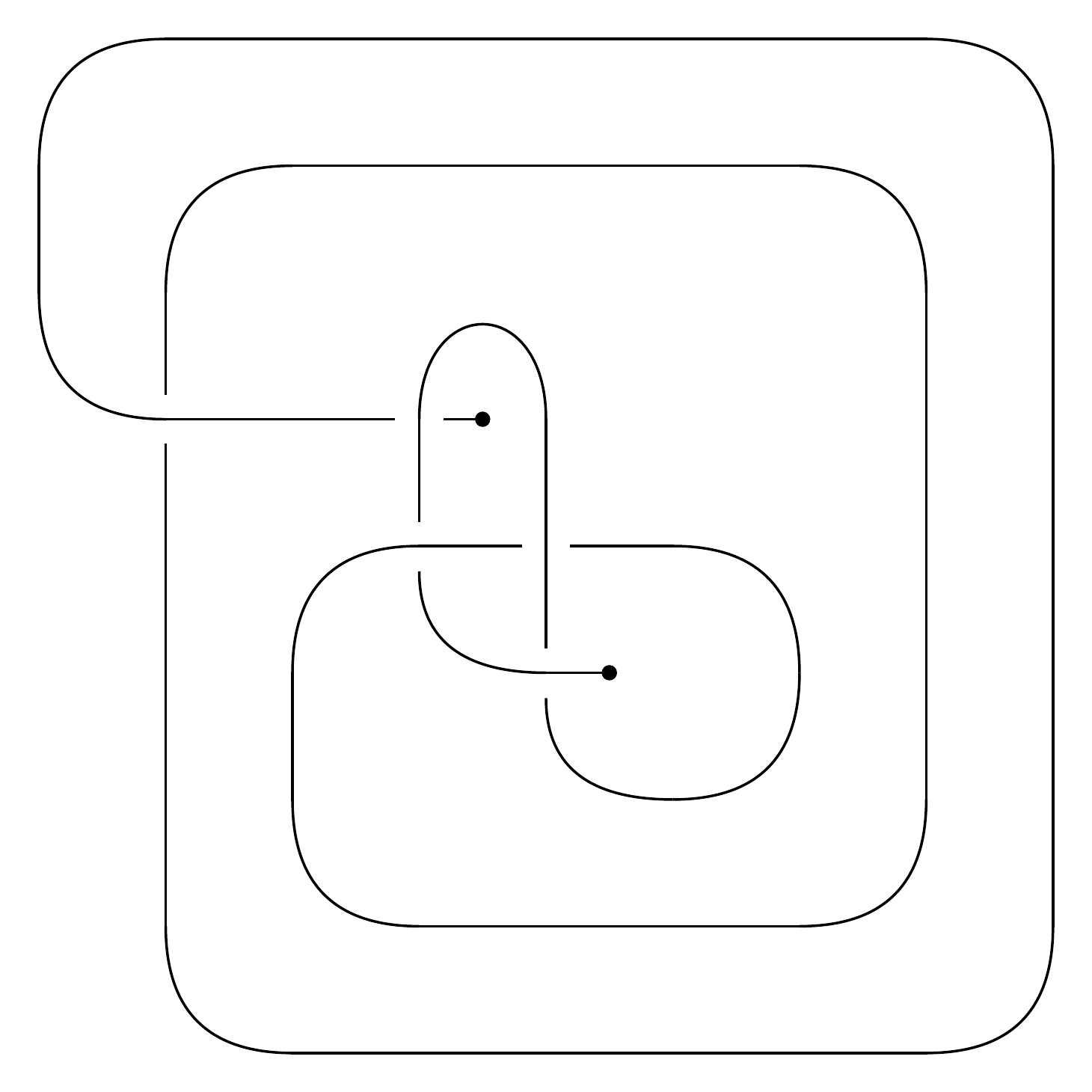}\\
\textcolor{black}{$5_{513}$}
\vspace{1cm}
\end{minipage}
\begin{minipage}[t]{.25\linewidth}
\centering
\includegraphics[width=0.9\textwidth,height=3.5cm,keepaspectratio]{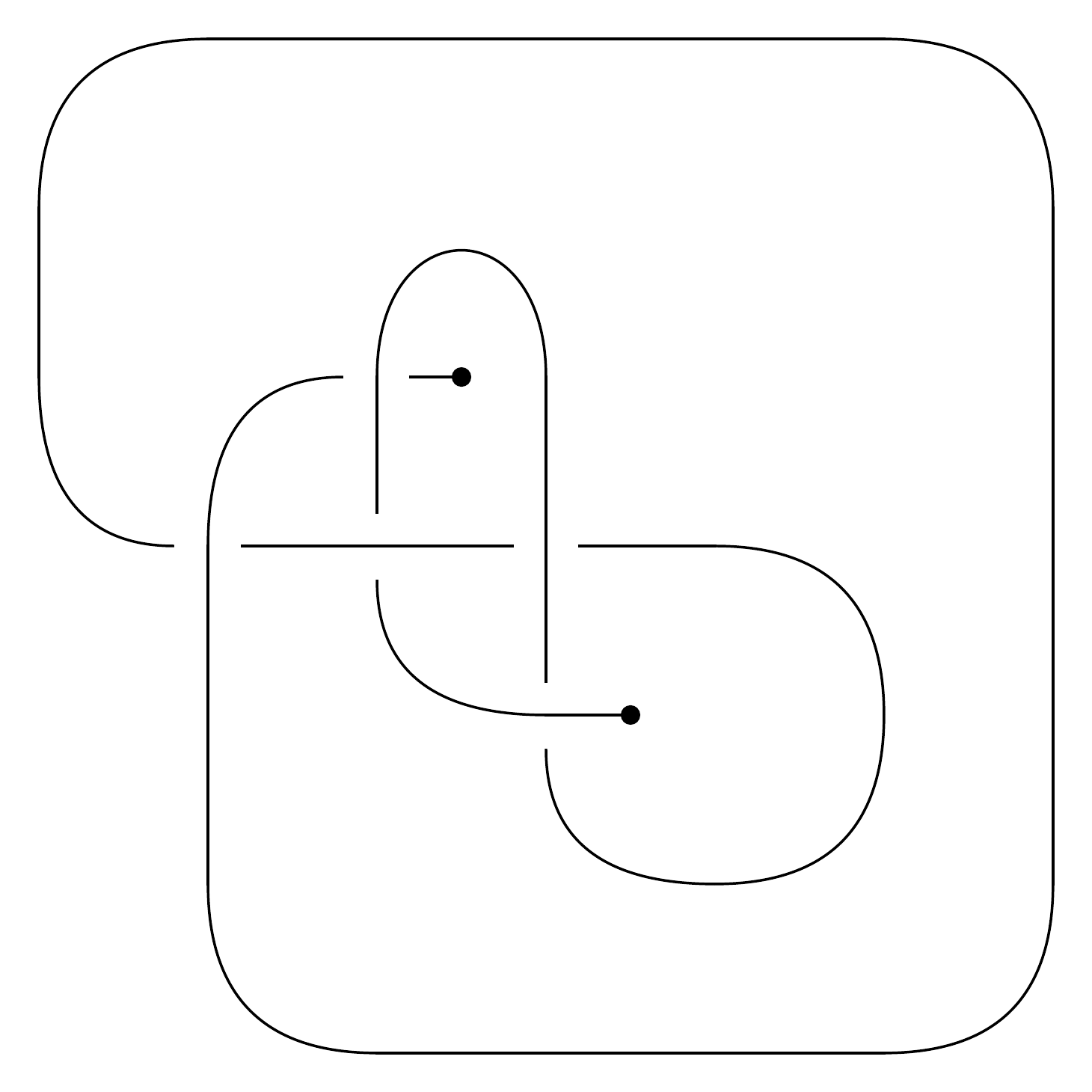}\\
\textcolor{black}{$5_{514}$}
\vspace{1cm}
\end{minipage}
\begin{minipage}[t]{.25\linewidth}
\centering
\includegraphics[width=0.9\textwidth,height=3.5cm,keepaspectratio]{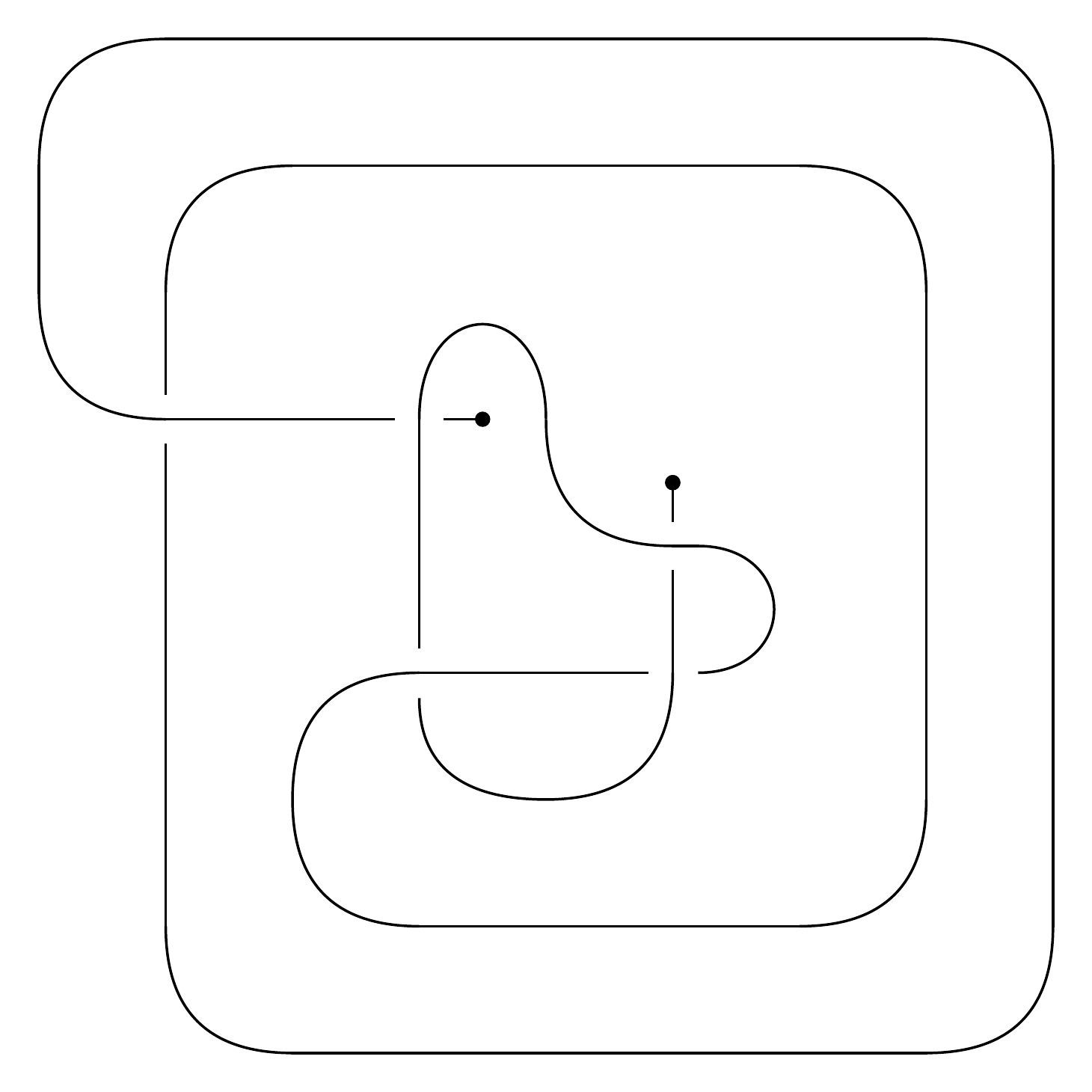}\\
\textcolor{black}{$5_{515}$}
\vspace{1cm}
\end{minipage}
\begin{minipage}[t]{.25\linewidth}
\centering
\includegraphics[width=0.9\textwidth,height=3.5cm,keepaspectratio]{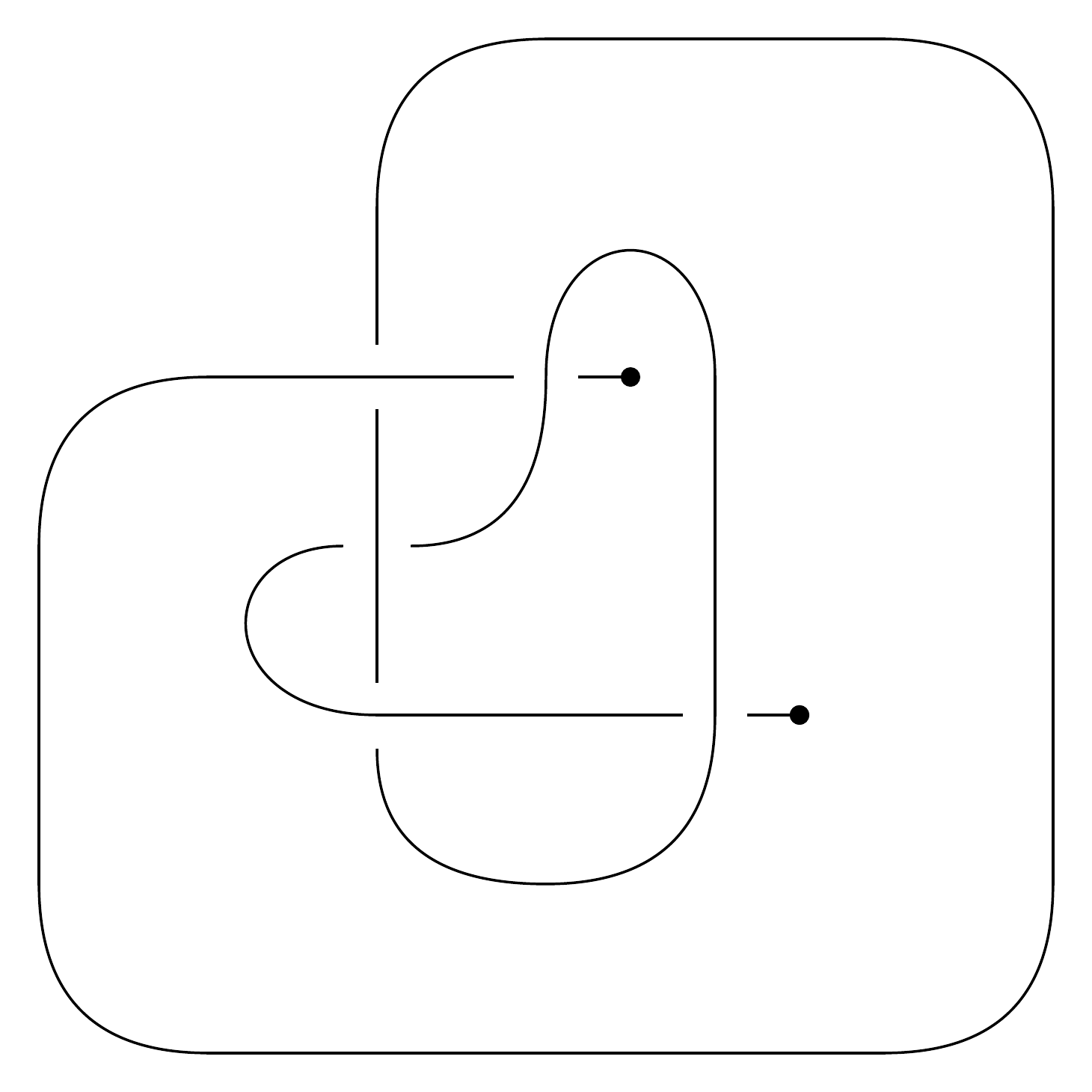}\\
\textcolor{black}{$5_{516}$}
\vspace{1cm}
\end{minipage}
\begin{minipage}[t]{.25\linewidth}
\centering
\includegraphics[width=0.9\textwidth,height=3.5cm,keepaspectratio]{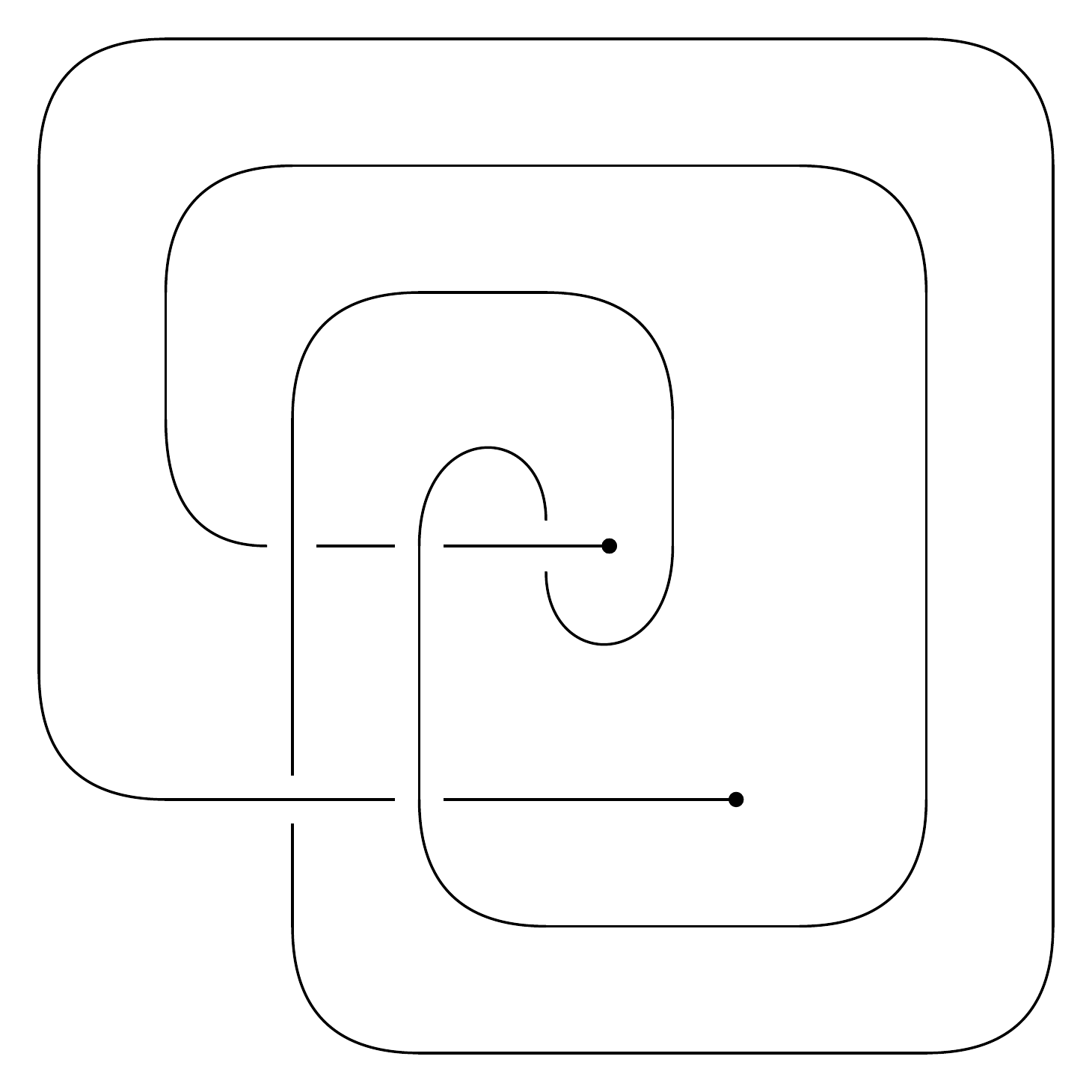}\\
\textcolor{black}{$5_{517}$}
\vspace{1cm}
\end{minipage}
\begin{minipage}[t]{.25\linewidth}
\centering
\includegraphics[width=0.9\textwidth,height=3.5cm,keepaspectratio]{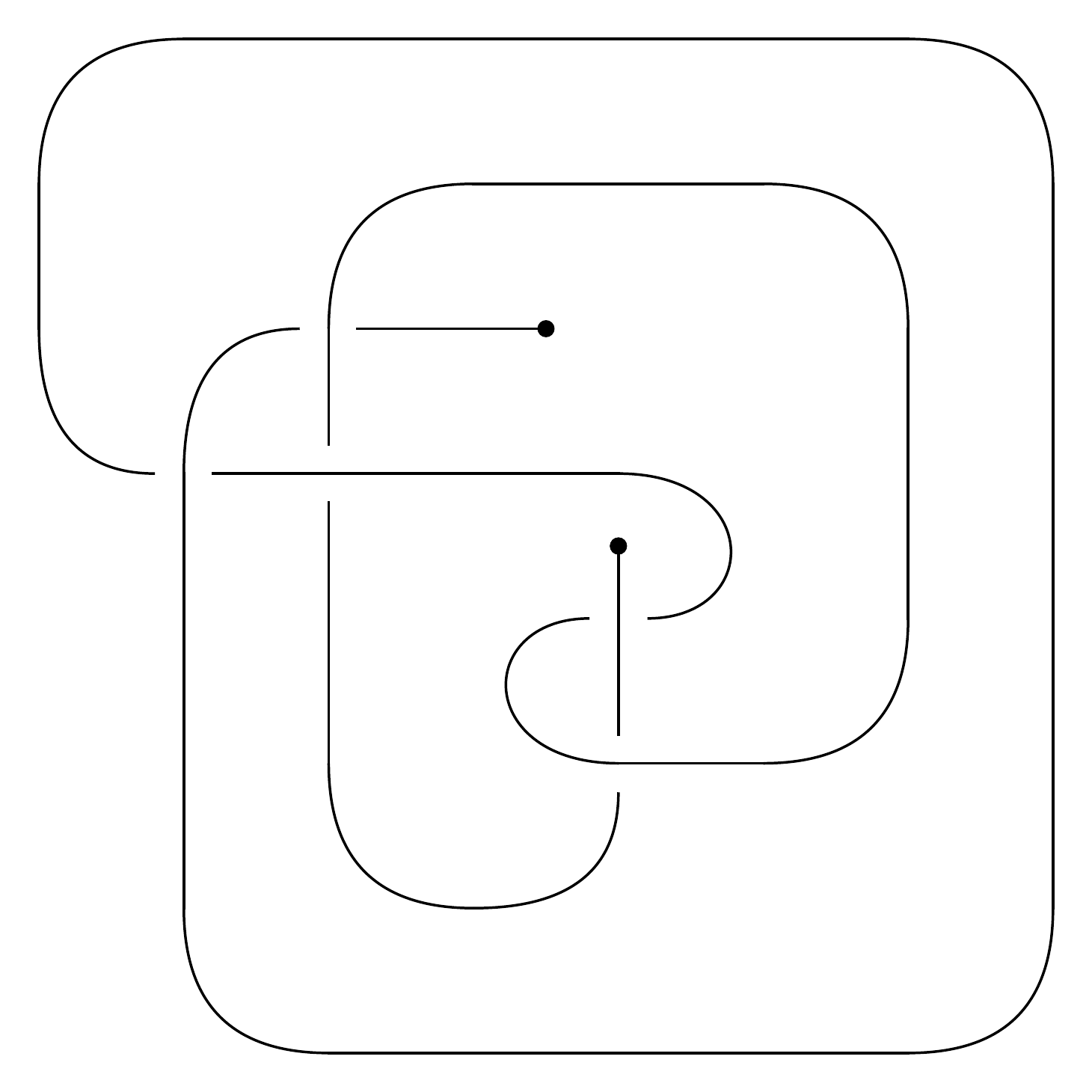}\\
\textcolor{black}{$5_{518}$}
\vspace{1cm}
\end{minipage}
\begin{minipage}[t]{.25\linewidth}
\centering
\includegraphics[width=0.9\textwidth,height=3.5cm,keepaspectratio]{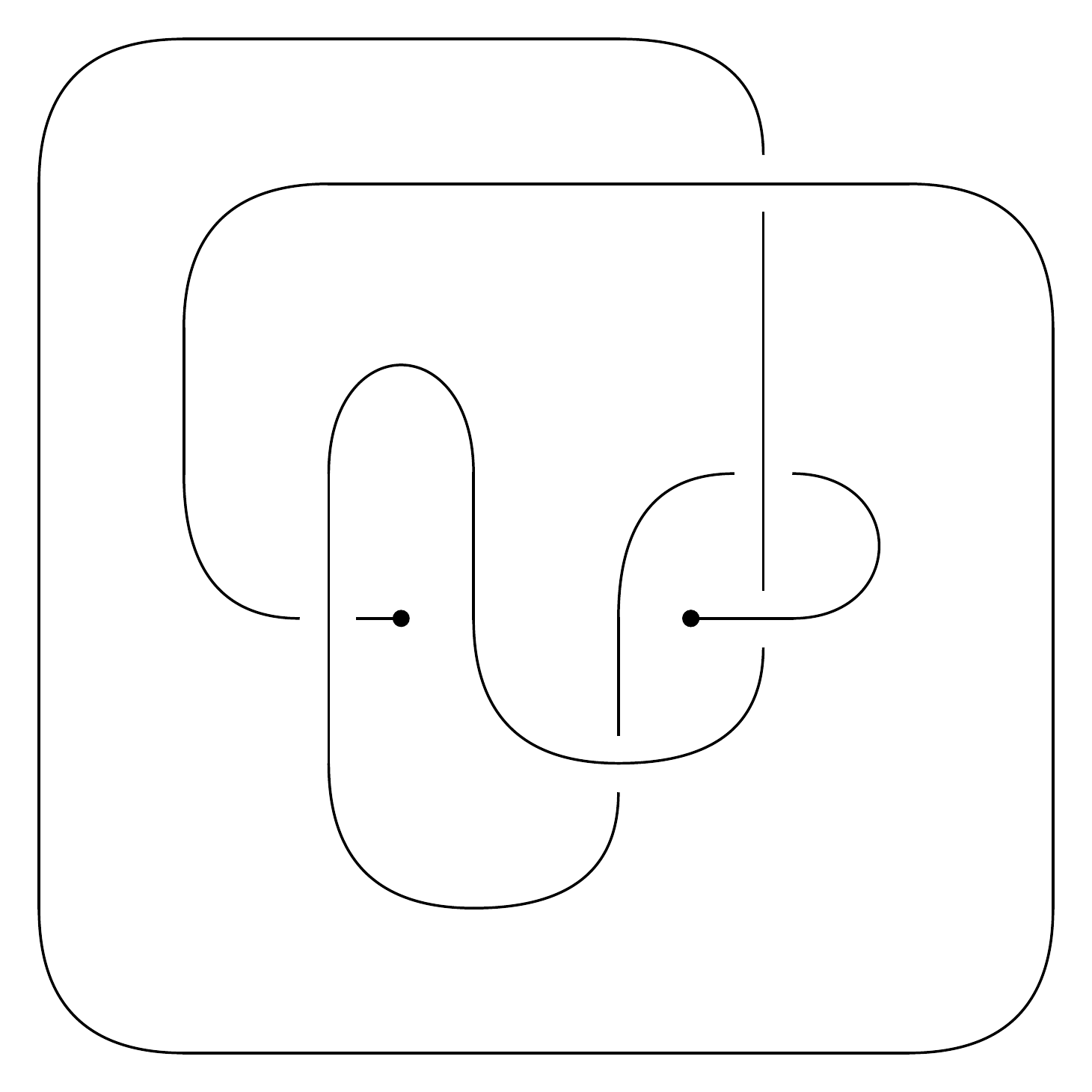}\\
\textcolor{black}{$5_{519}$}
\vspace{1cm}
\end{minipage}
\begin{minipage}[t]{.25\linewidth}
\centering
\includegraphics[width=0.9\textwidth,height=3.5cm,keepaspectratio]{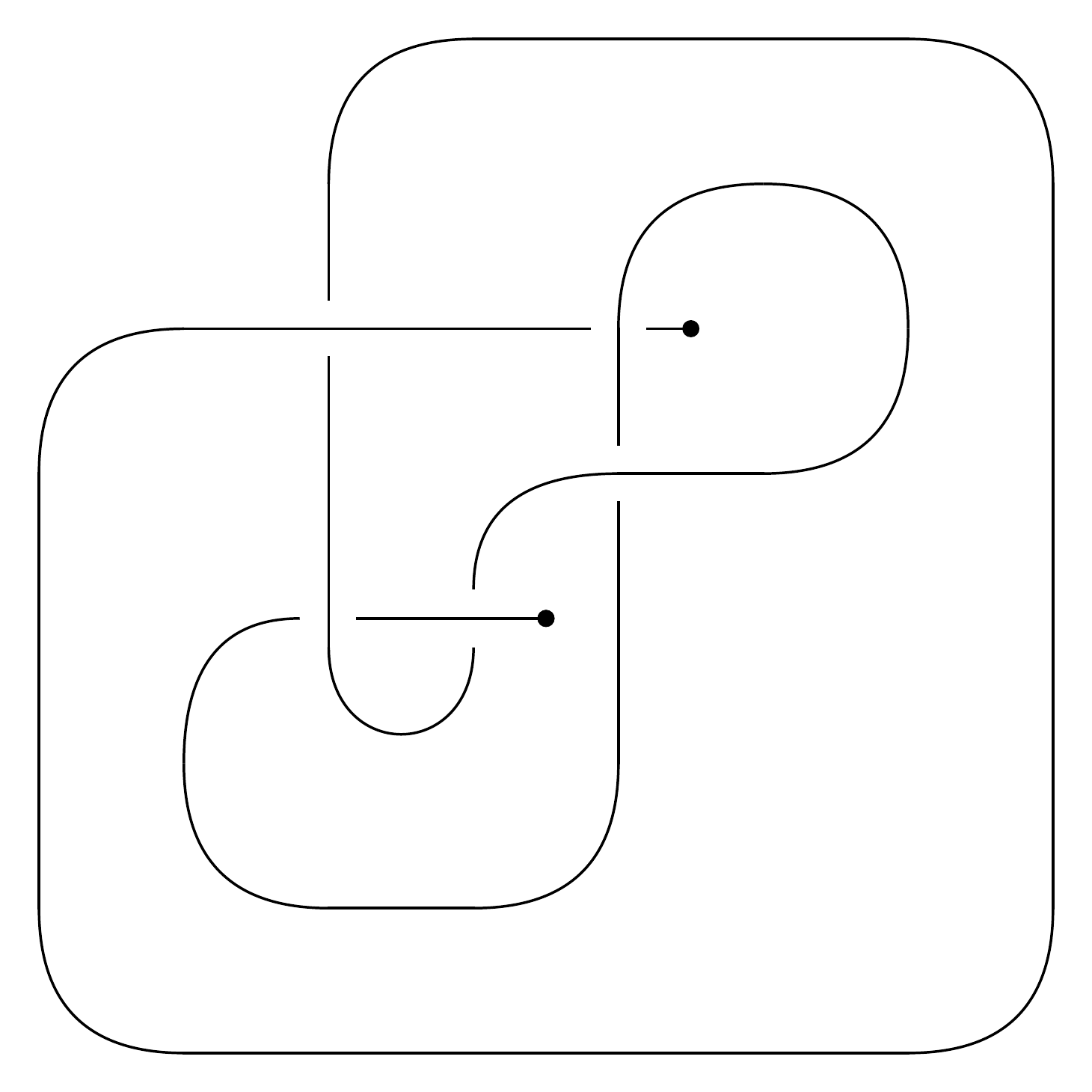}\\
\textcolor{black}{$5_{520}$}
\vspace{1cm}
\end{minipage}
\begin{minipage}[t]{.25\linewidth}
\centering
\includegraphics[width=0.9\textwidth,height=3.5cm,keepaspectratio]{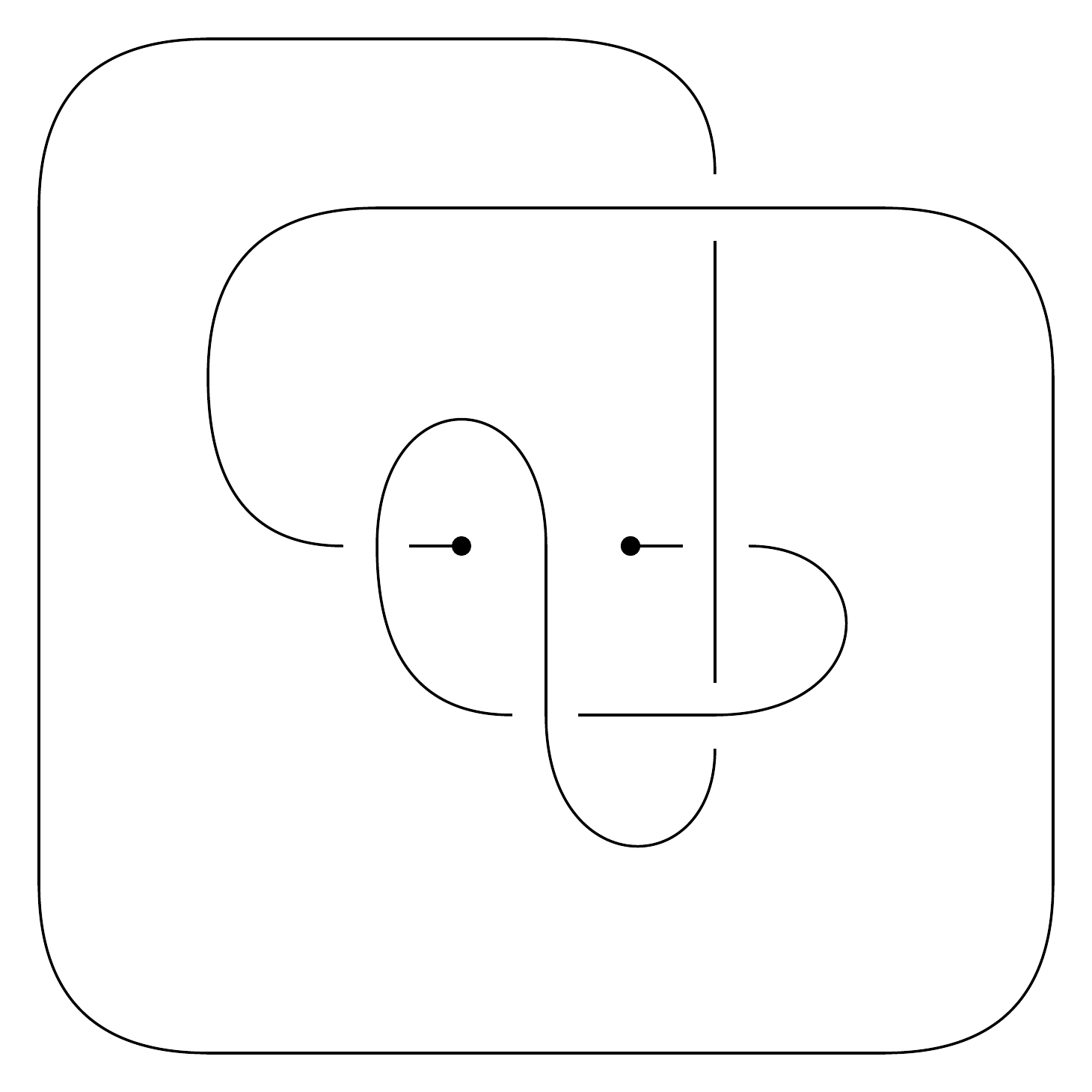}\\
\textcolor{black}{$5_{521}$}
\vspace{1cm}
\end{minipage}
\begin{minipage}[t]{.25\linewidth}
\centering
\includegraphics[width=0.9\textwidth,height=3.5cm,keepaspectratio]{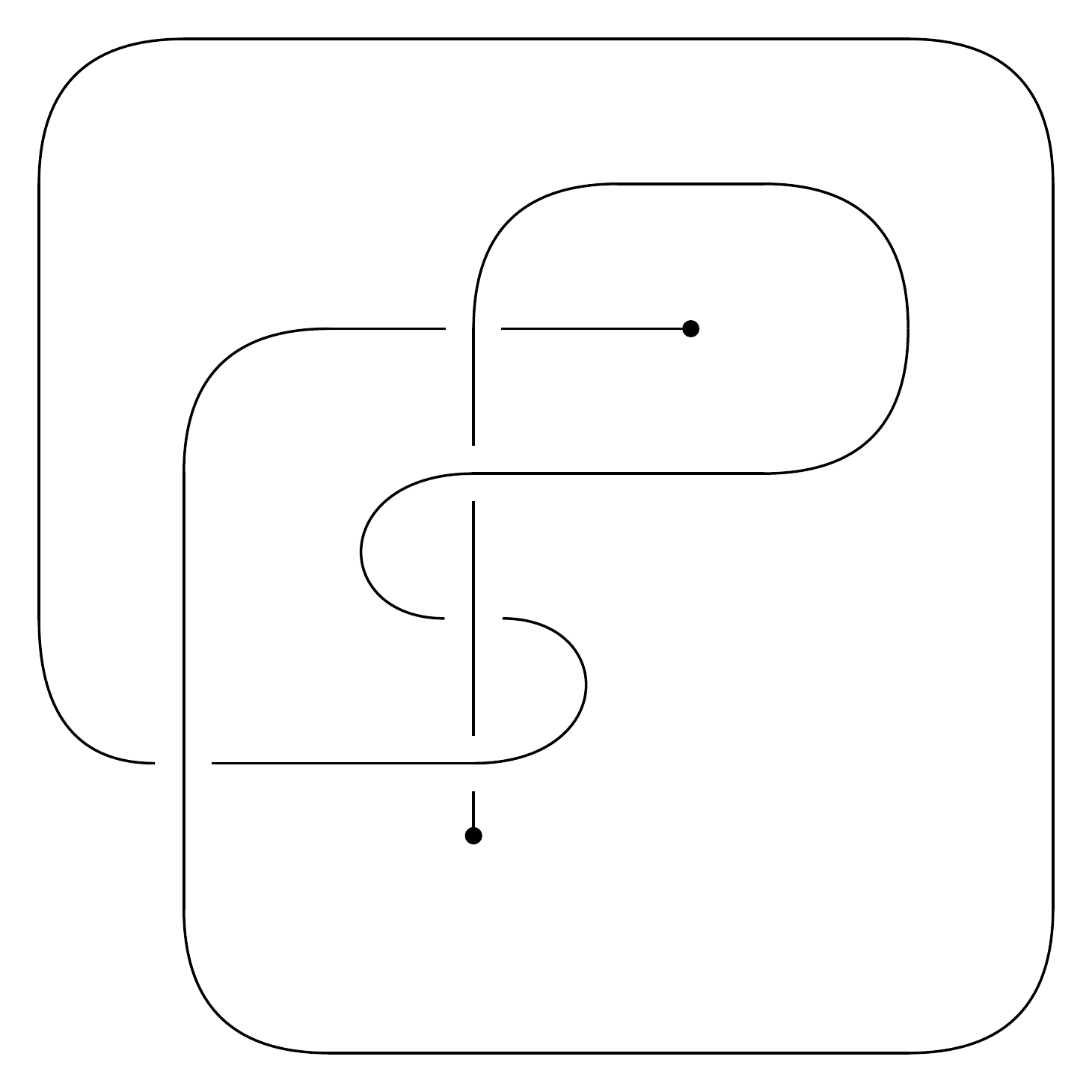}\\
\textcolor{black}{$5_{522}$}
\vspace{1cm}
\end{minipage}
\begin{minipage}[t]{.25\linewidth}
\centering
\includegraphics[width=0.9\textwidth,height=3.5cm,keepaspectratio]{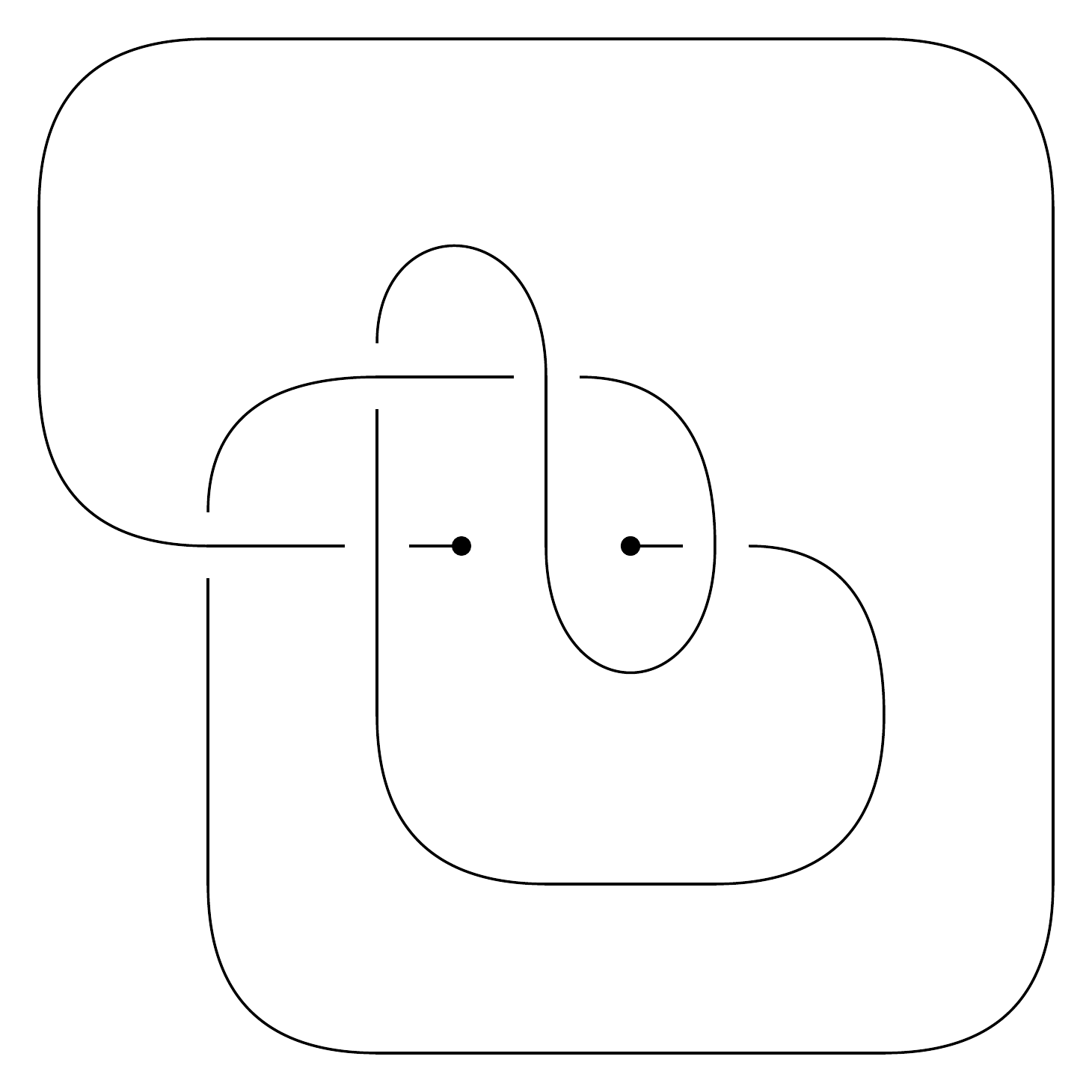}\\
\textcolor{black}{$5_{523}$}
\vspace{1cm}
\end{minipage}
\begin{minipage}[t]{.25\linewidth}
\centering
\includegraphics[width=0.9\textwidth,height=3.5cm,keepaspectratio]{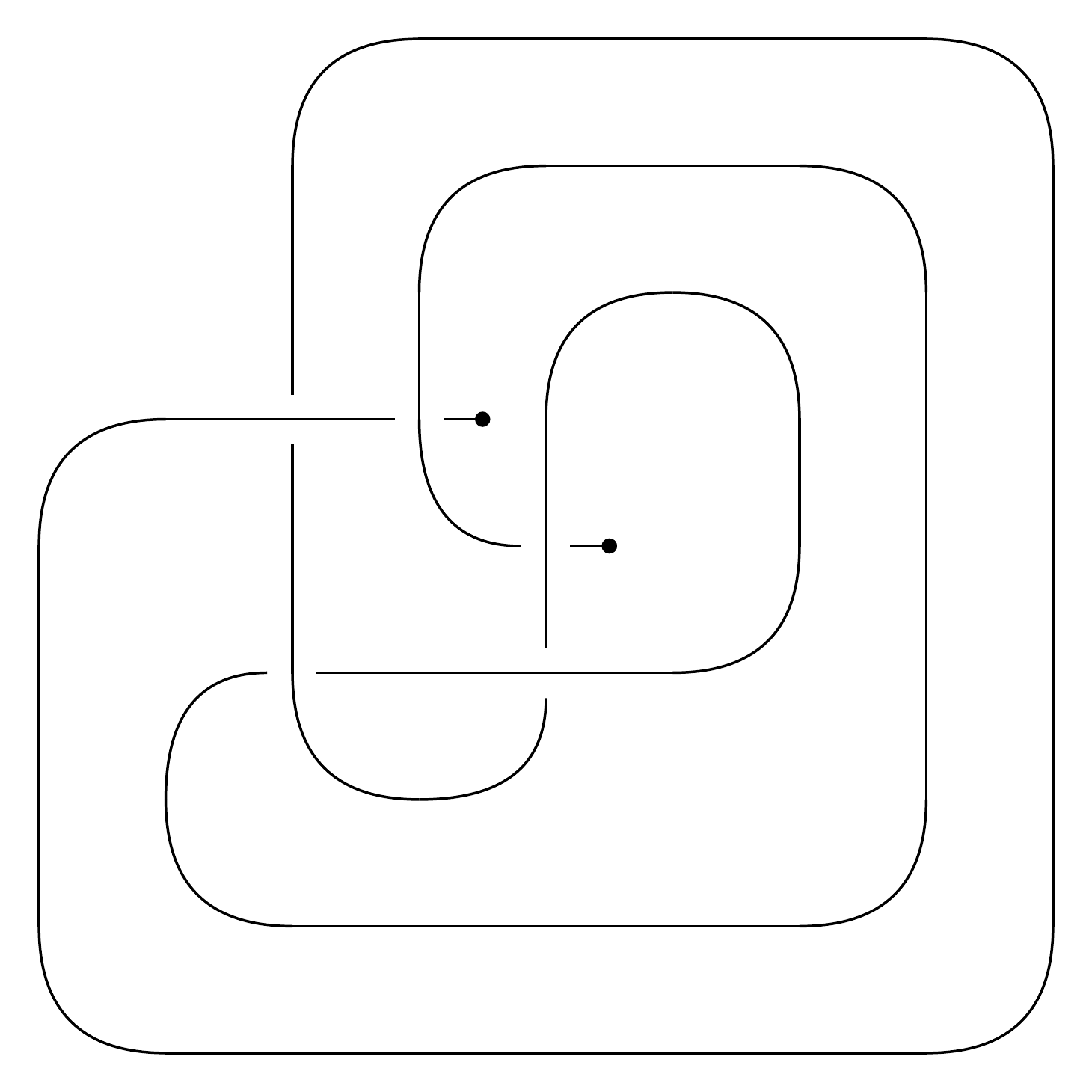}\\
\textcolor{black}{$5_{524}$}
\vspace{1cm}
\end{minipage}
\begin{minipage}[t]{.25\linewidth}
\centering
\includegraphics[width=0.9\textwidth,height=3.5cm,keepaspectratio]{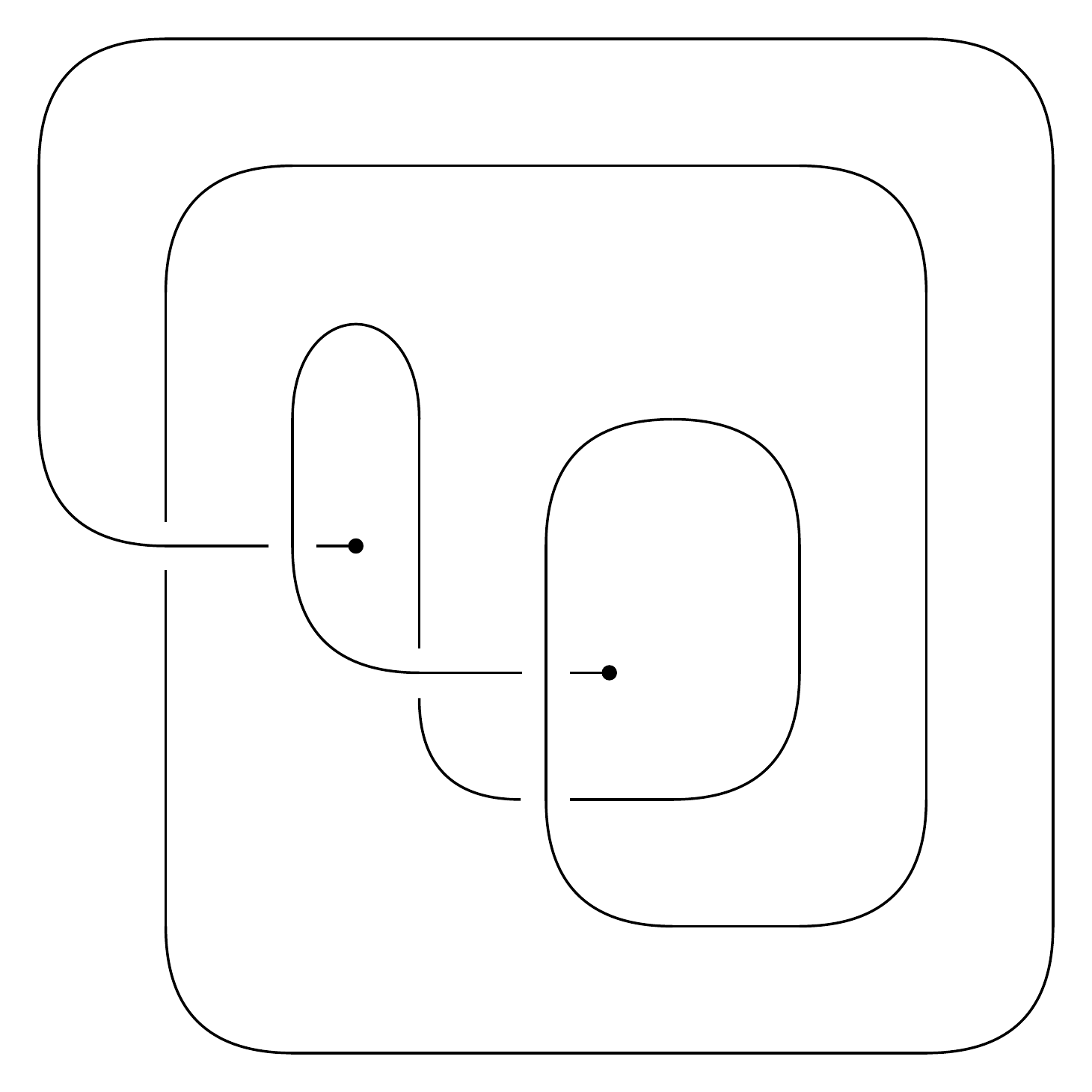}\\
\textcolor{black}{$5_{525}$}
\vspace{1cm}
\end{minipage}
\begin{minipage}[t]{.25\linewidth}
\centering
\includegraphics[width=0.9\textwidth,height=3.5cm,keepaspectratio]{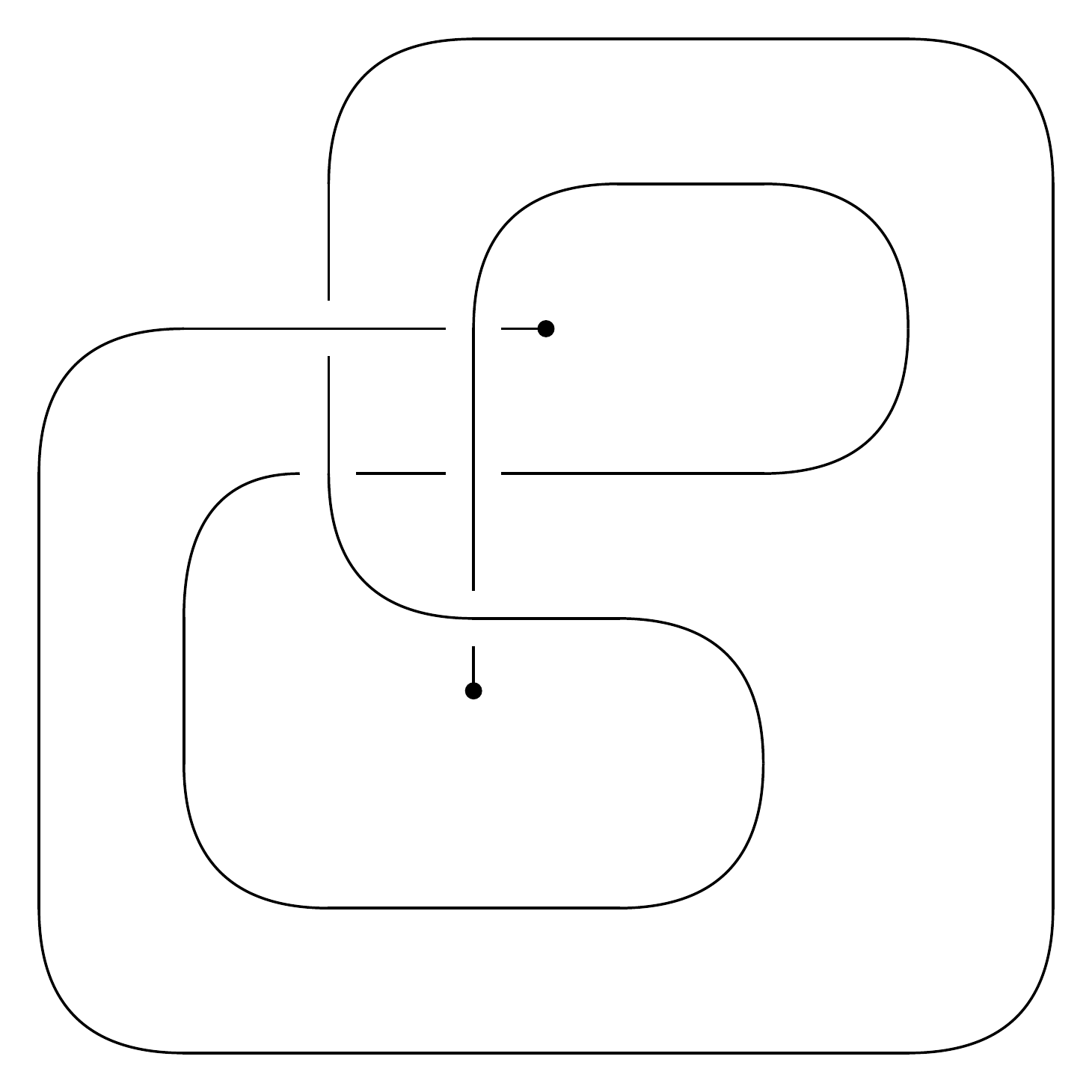}\\
\textcolor{black}{$5_{526}$}
\vspace{1cm}
\end{minipage}
\begin{minipage}[t]{.25\linewidth}
\centering
\includegraphics[width=0.9\textwidth,height=3.5cm,keepaspectratio]{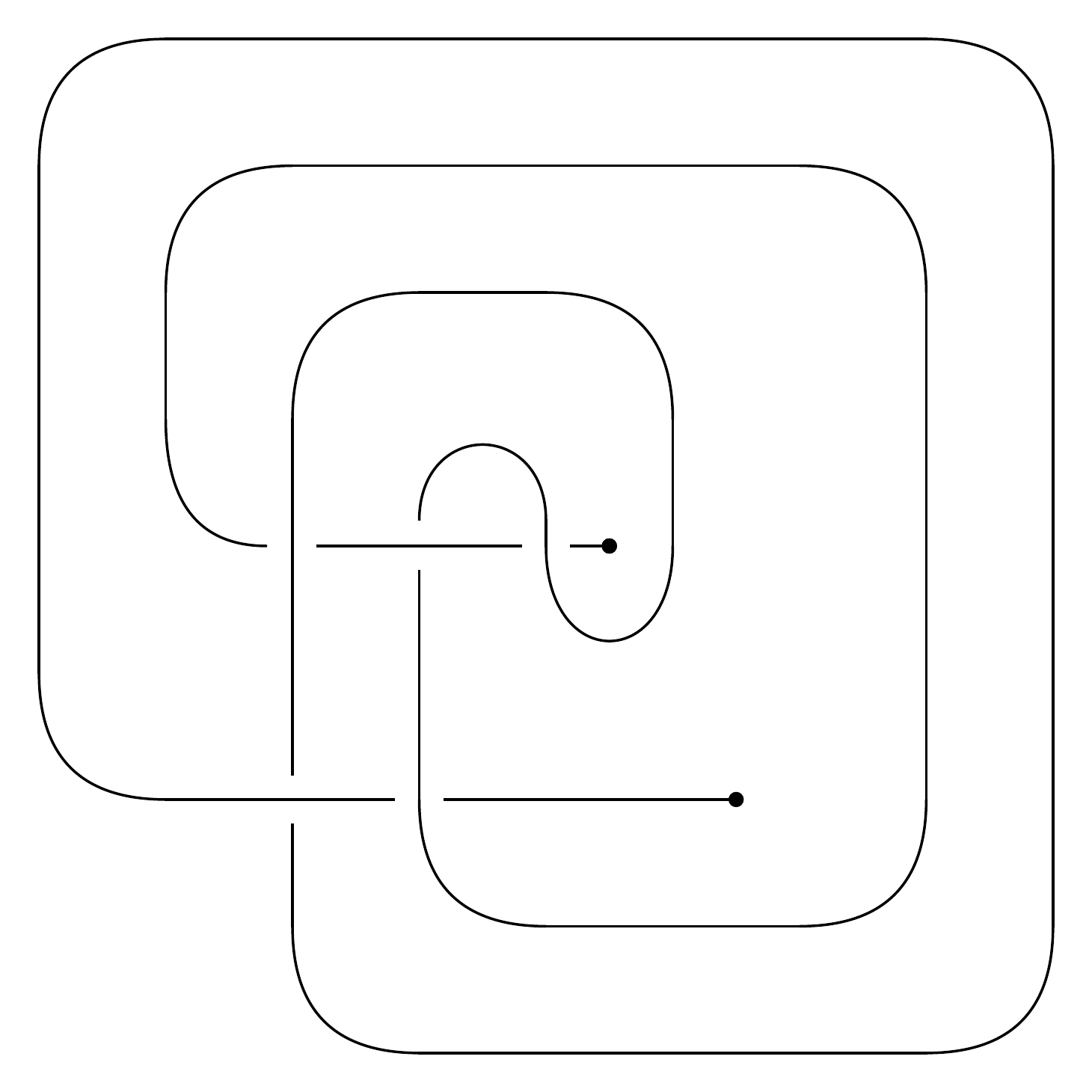}\\
\textcolor{black}{$5_{527}$}
\vspace{1cm}
\end{minipage}
\begin{minipage}[t]{.25\linewidth}
\centering
\includegraphics[width=0.9\textwidth,height=3.5cm,keepaspectratio]{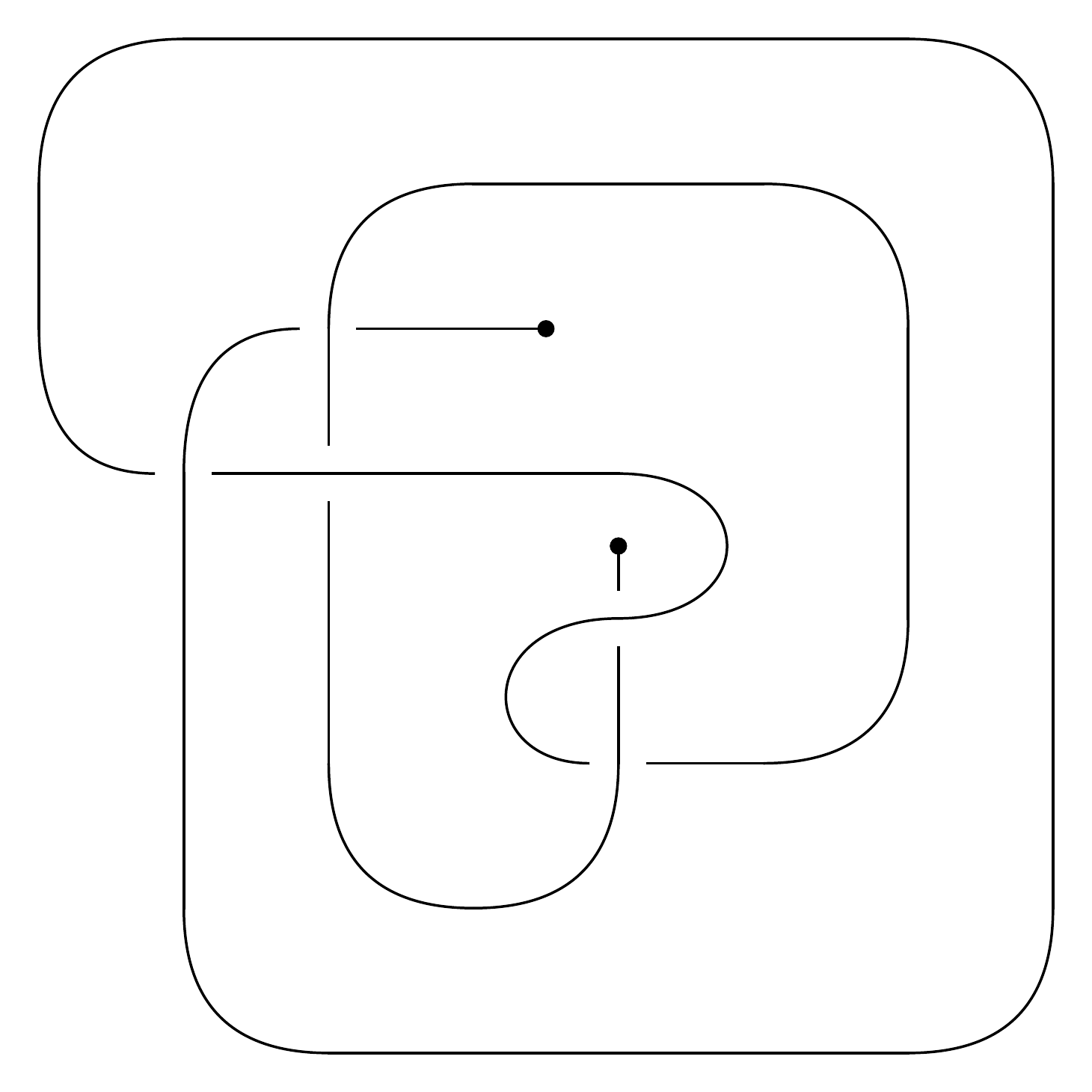}\\
\textcolor{black}{$5_{528}$}
\vspace{1cm}
\end{minipage}
\begin{minipage}[t]{.25\linewidth}
\centering
\includegraphics[width=0.9\textwidth,height=3.5cm,keepaspectratio]{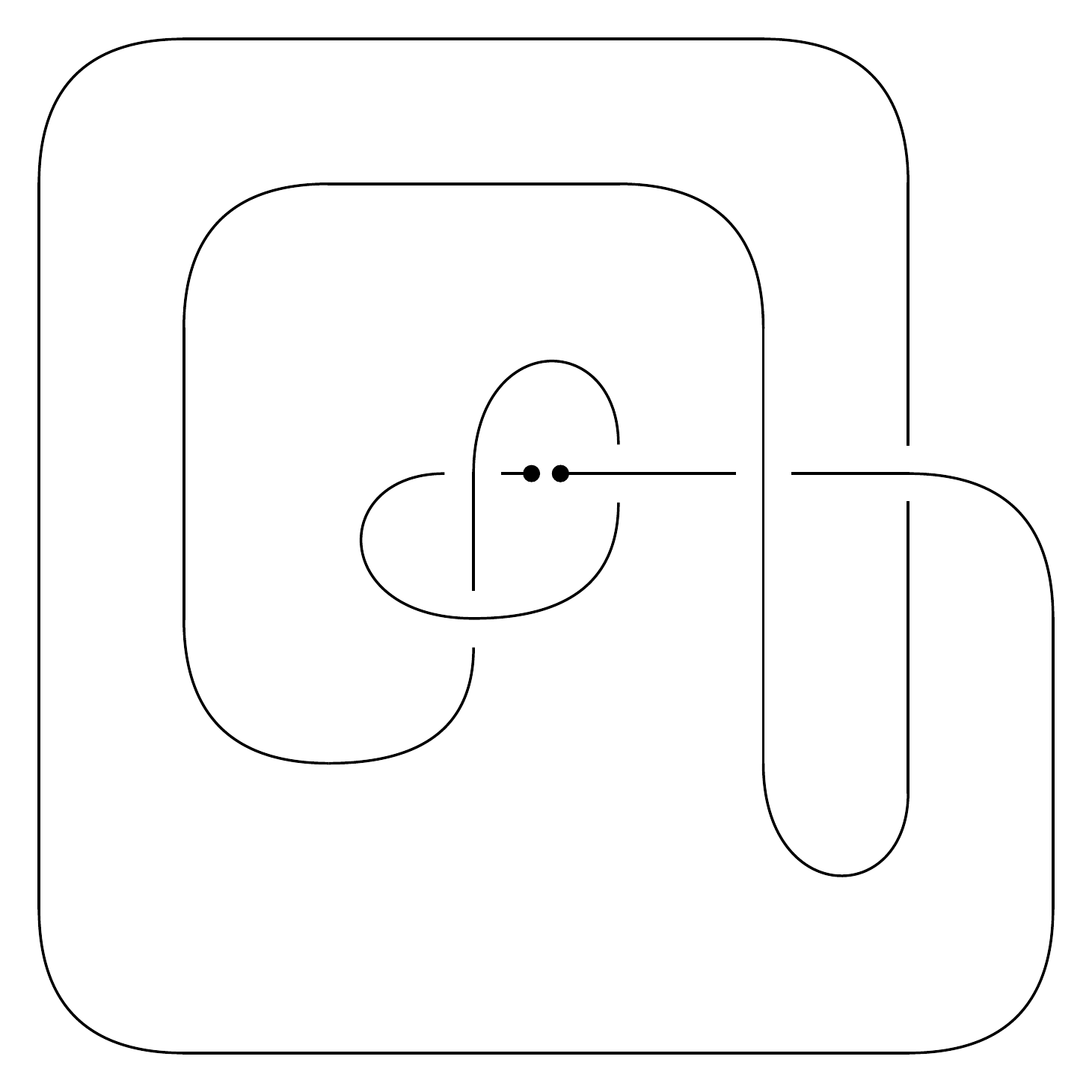}\\
\textcolor{black}{$5_{529}$}
\vspace{1cm}
\end{minipage}
\begin{minipage}[t]{.25\linewidth}
\centering
\includegraphics[width=0.9\textwidth,height=3.5cm,keepaspectratio]{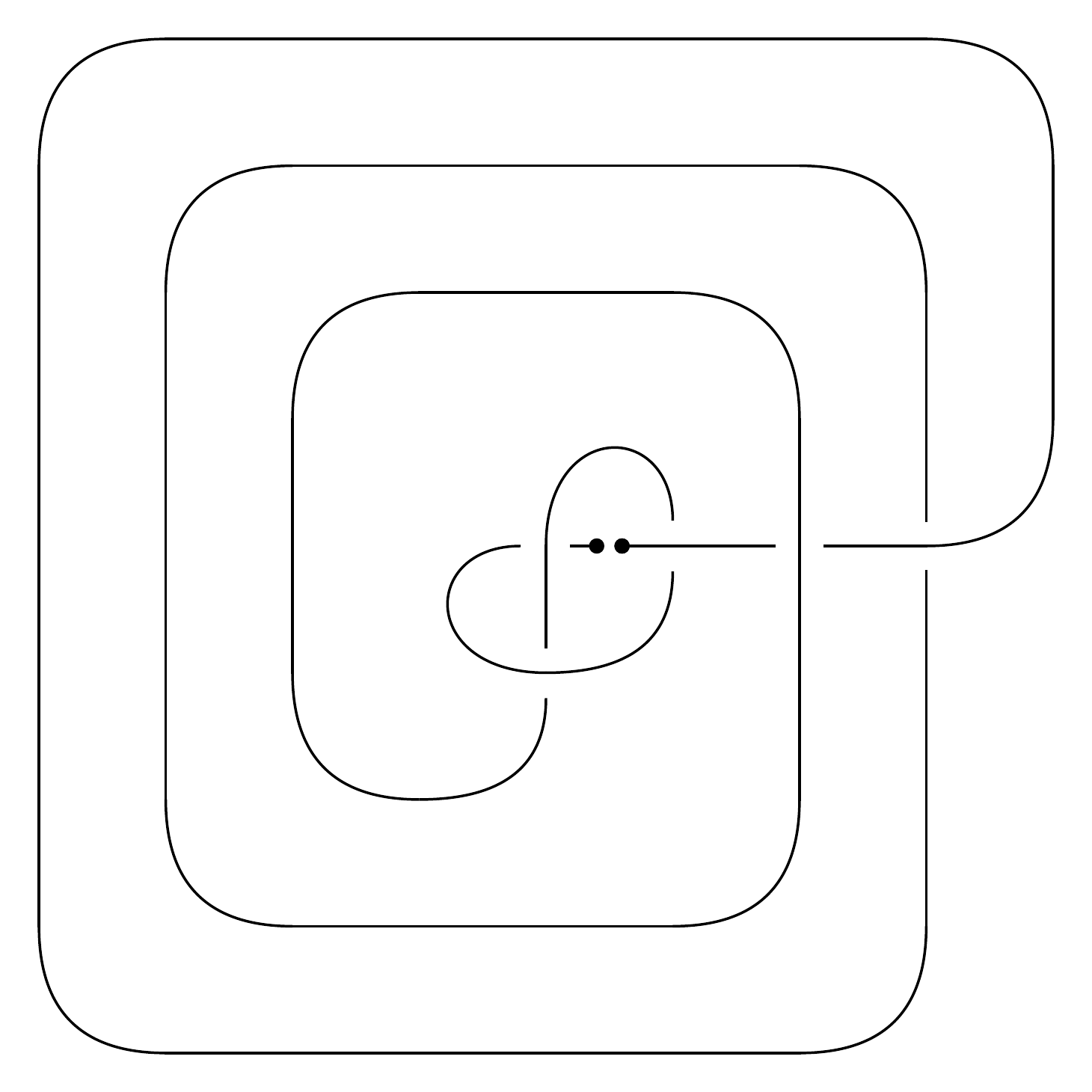}\\
\textcolor{black}{$5_{530}$}
\vspace{1cm}
\end{minipage}
\begin{minipage}[t]{.25\linewidth}
\centering
\includegraphics[width=0.9\textwidth,height=3.5cm,keepaspectratio]{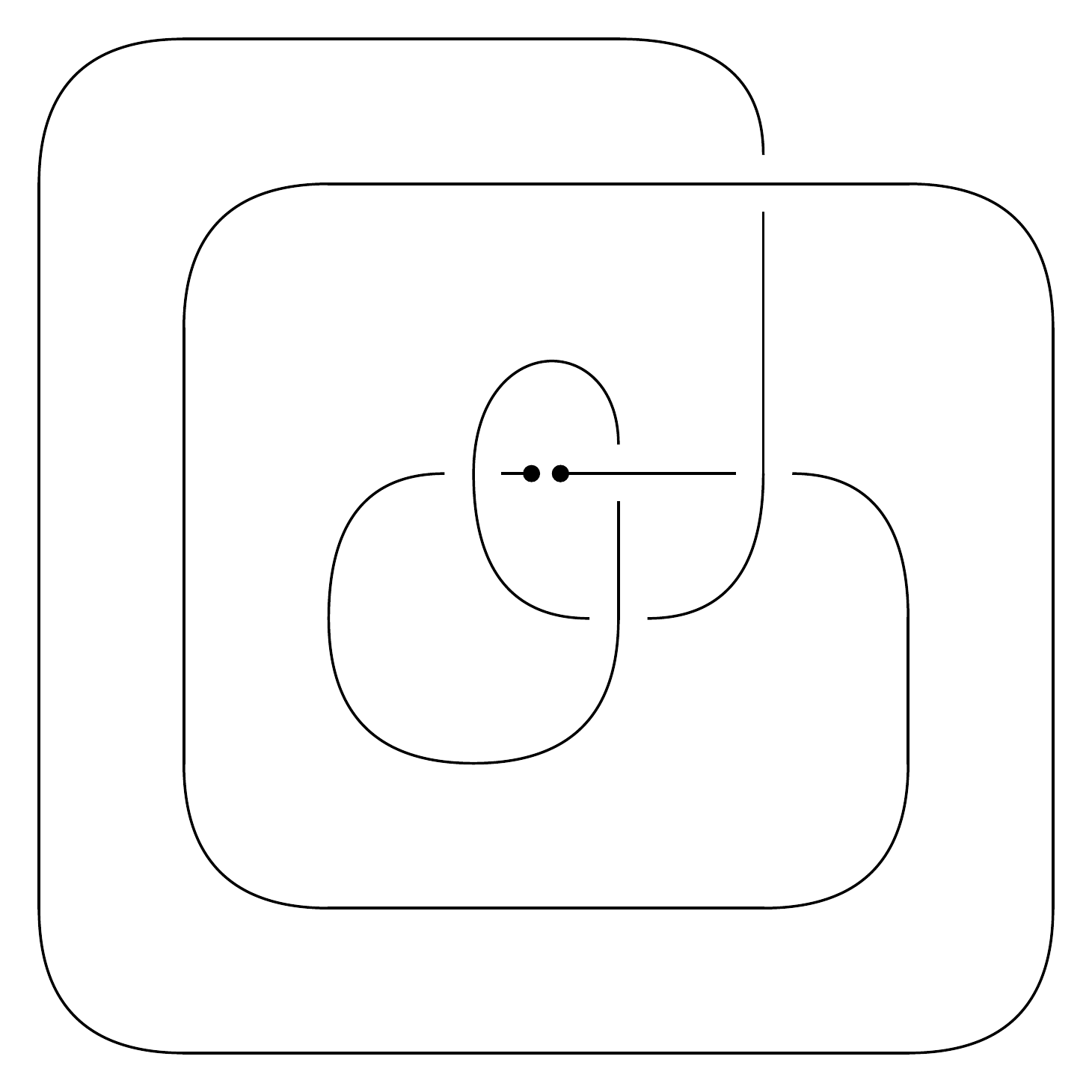}\\
\textcolor{black}{$5_{531}$}
\vspace{1cm}
\end{minipage}
\begin{minipage}[t]{.25\linewidth}
\centering
\includegraphics[width=0.9\textwidth,height=3.5cm,keepaspectratio]{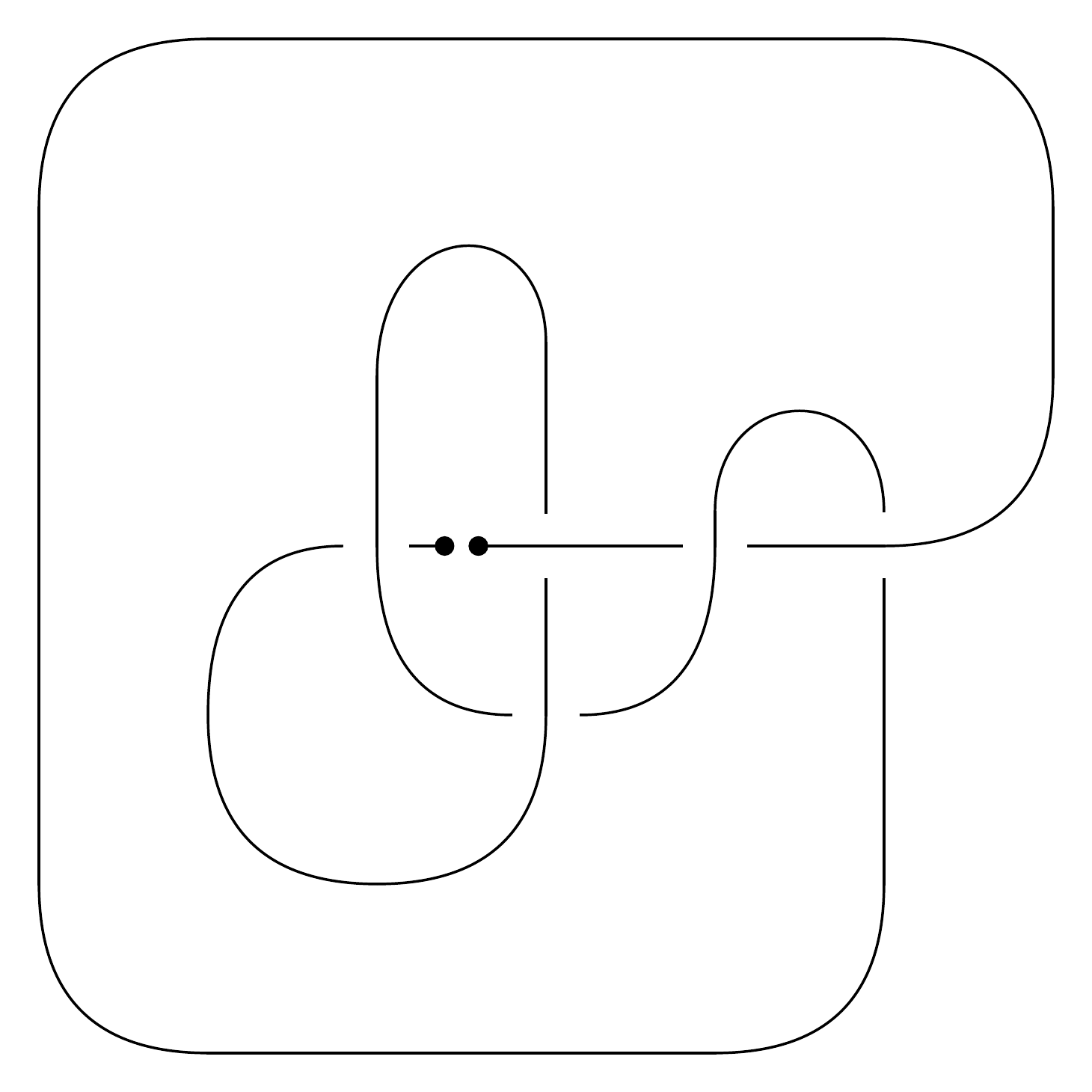}\\
\textcolor{black}{$5_{532}$}
\vspace{1cm}
\end{minipage}
\begin{minipage}[t]{.25\linewidth}
\centering
\includegraphics[width=0.9\textwidth,height=3.5cm,keepaspectratio]{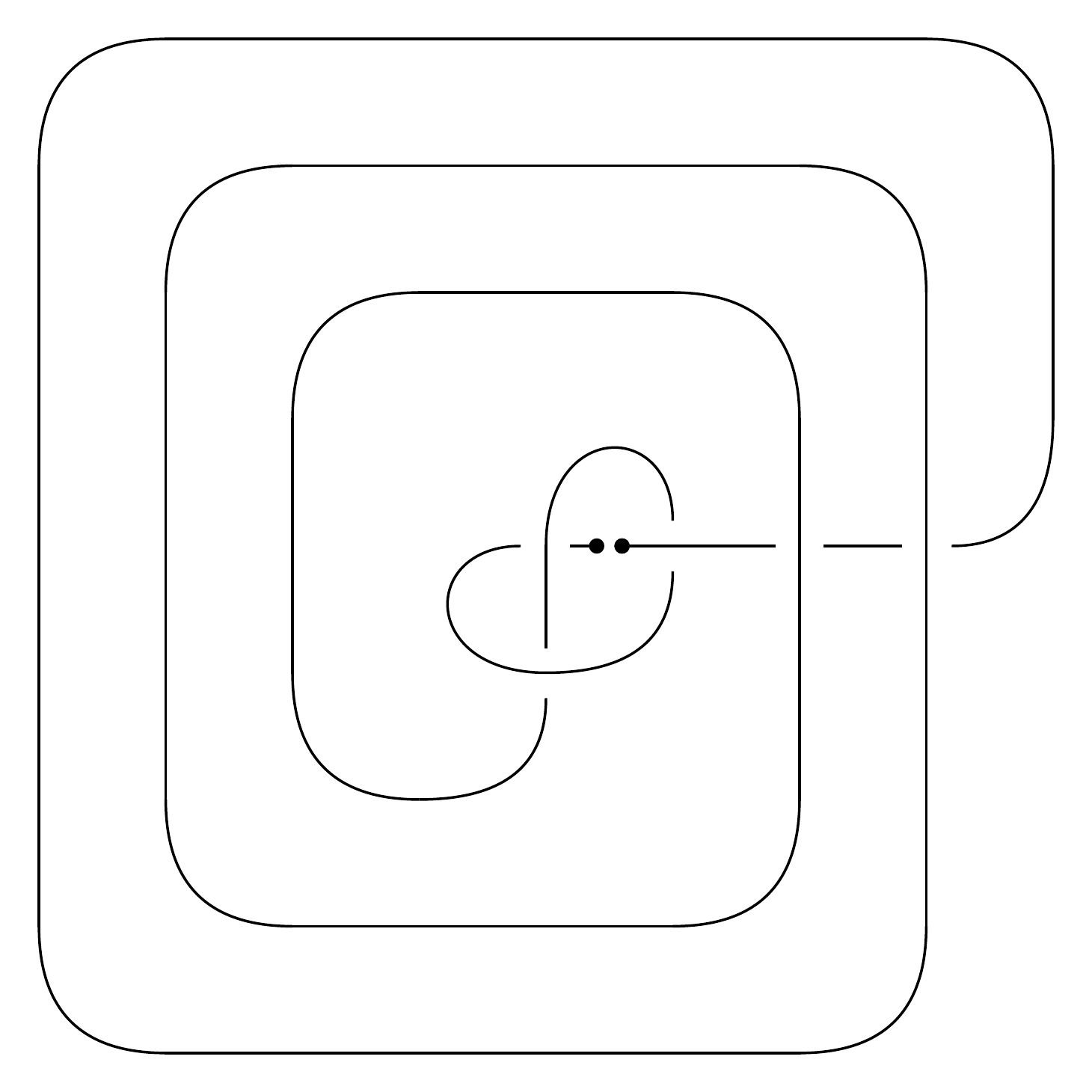}\\
\textcolor{black}{$5_{533}$}
\vspace{1cm}
\end{minipage}
\begin{minipage}[t]{.25\linewidth}
\centering
\includegraphics[width=0.9\textwidth,height=3.5cm,keepaspectratio]{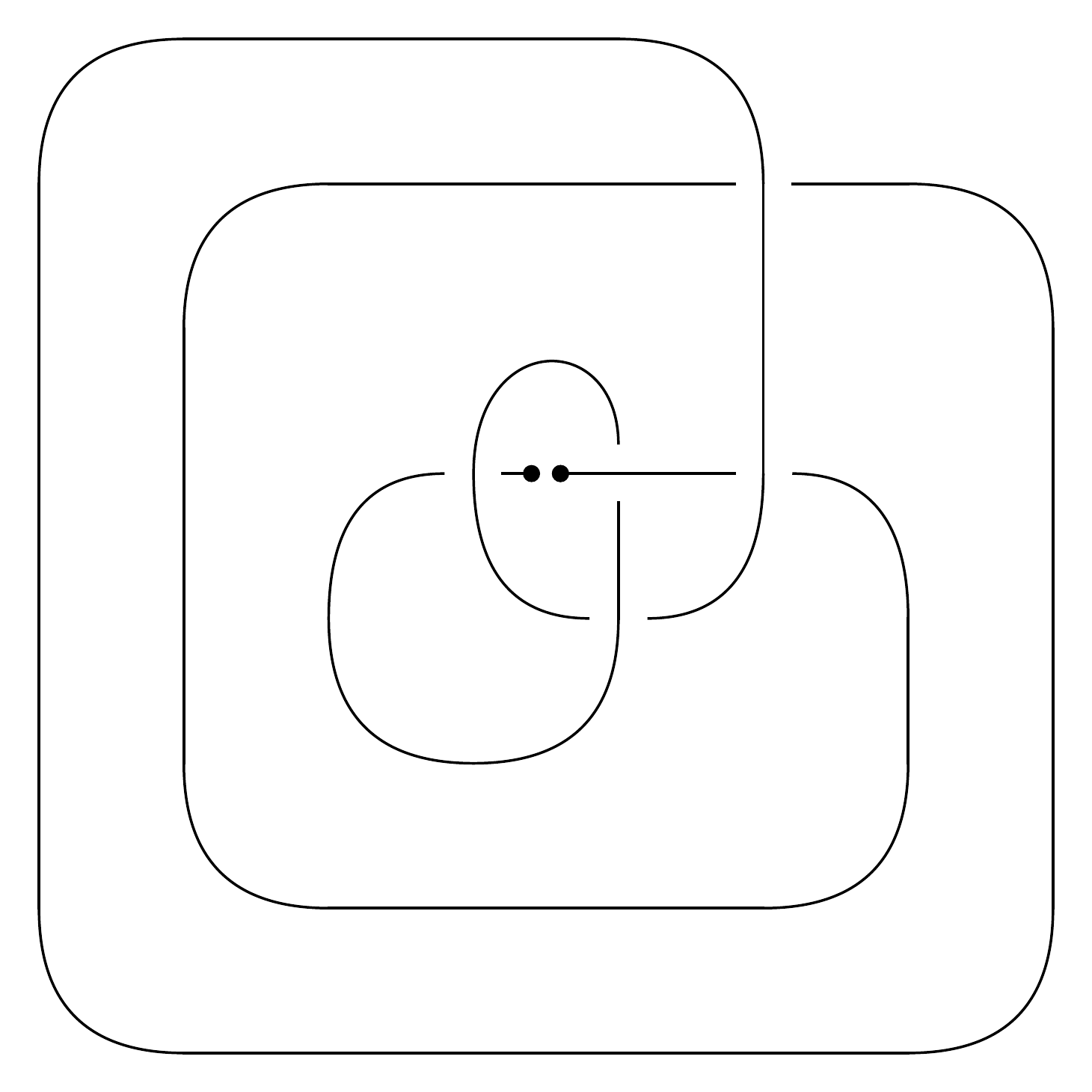}\\
\textcolor{black}{$5_{534}$}
\vspace{1cm}
\end{minipage}
\begin{minipage}[t]{.25\linewidth}
\centering
\includegraphics[width=0.9\textwidth,height=3.5cm,keepaspectratio]{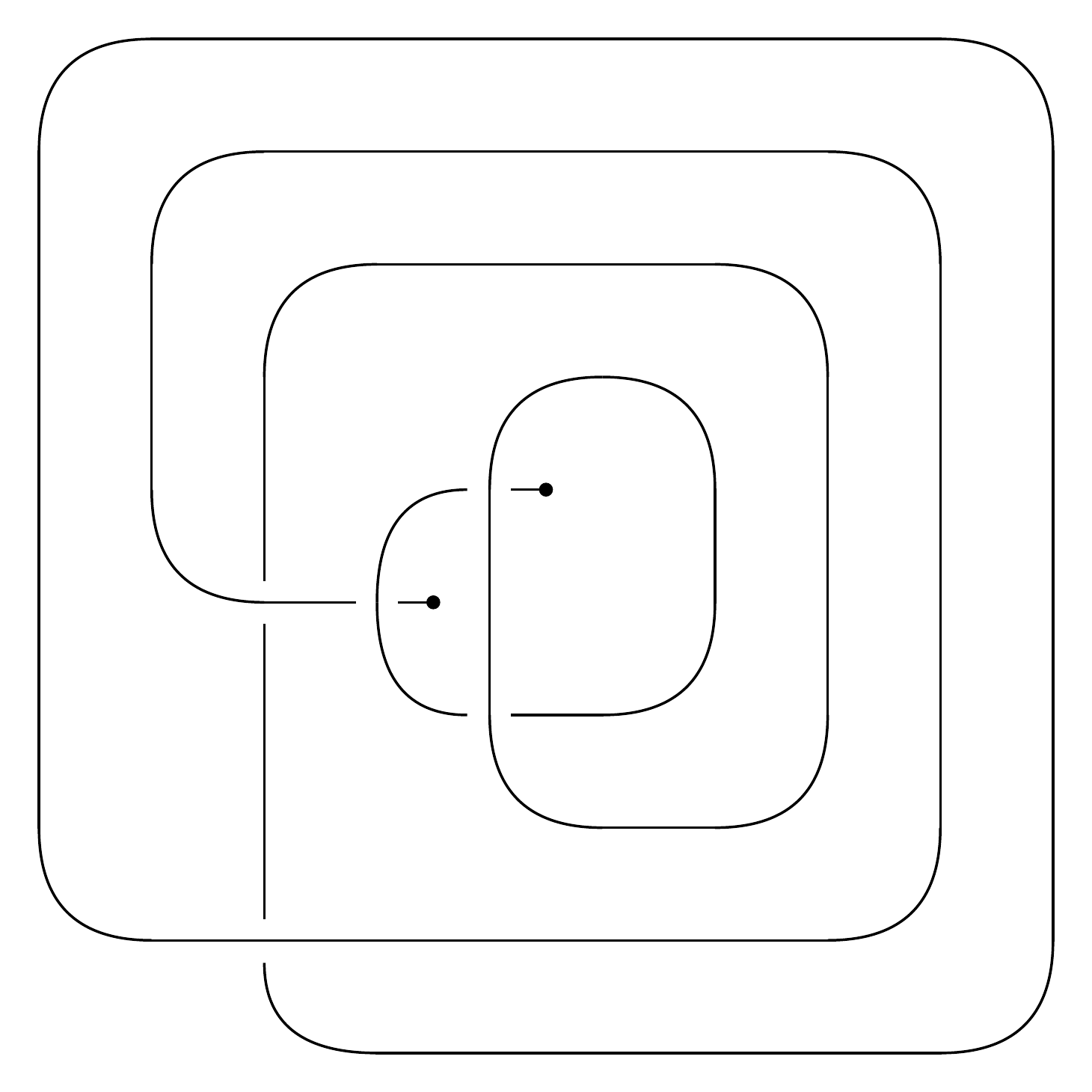}\\
\textcolor{black}{$5_{535}$}
\vspace{1cm}
\end{minipage}
\begin{minipage}[t]{.25\linewidth}
\centering
\includegraphics[width=0.9\textwidth,height=3.5cm,keepaspectratio]{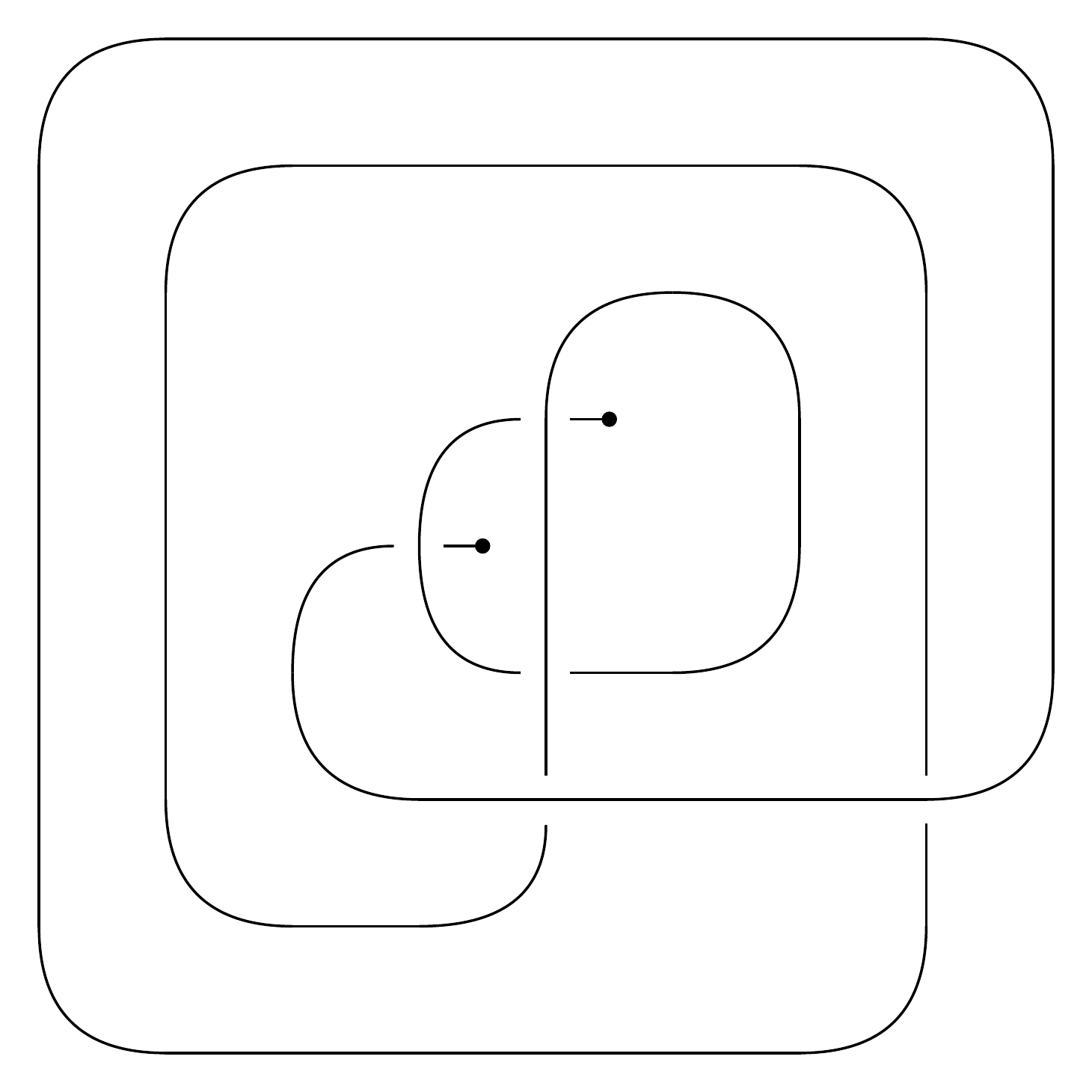}\\
\textcolor{black}{$5_{536}$}
\vspace{1cm}
\end{minipage}
\begin{minipage}[t]{.25\linewidth}
\centering
\includegraphics[width=0.9\textwidth,height=3.5cm,keepaspectratio]{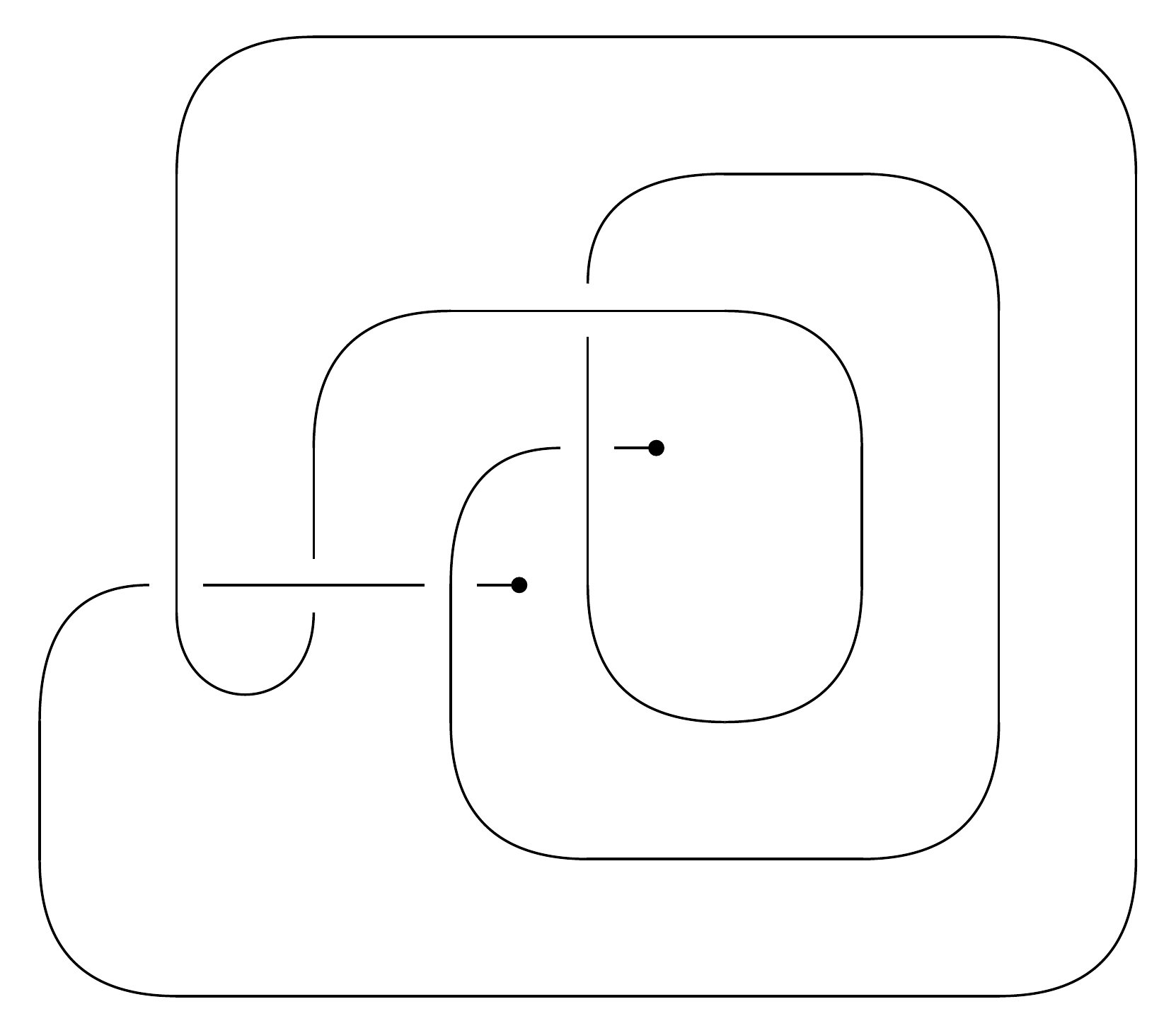}\\
\textcolor{black}{$5_{537}$}
\vspace{1cm}
\end{minipage}
\begin{minipage}[t]{.25\linewidth}
\centering
\includegraphics[width=0.9\textwidth,height=3.5cm,keepaspectratio]{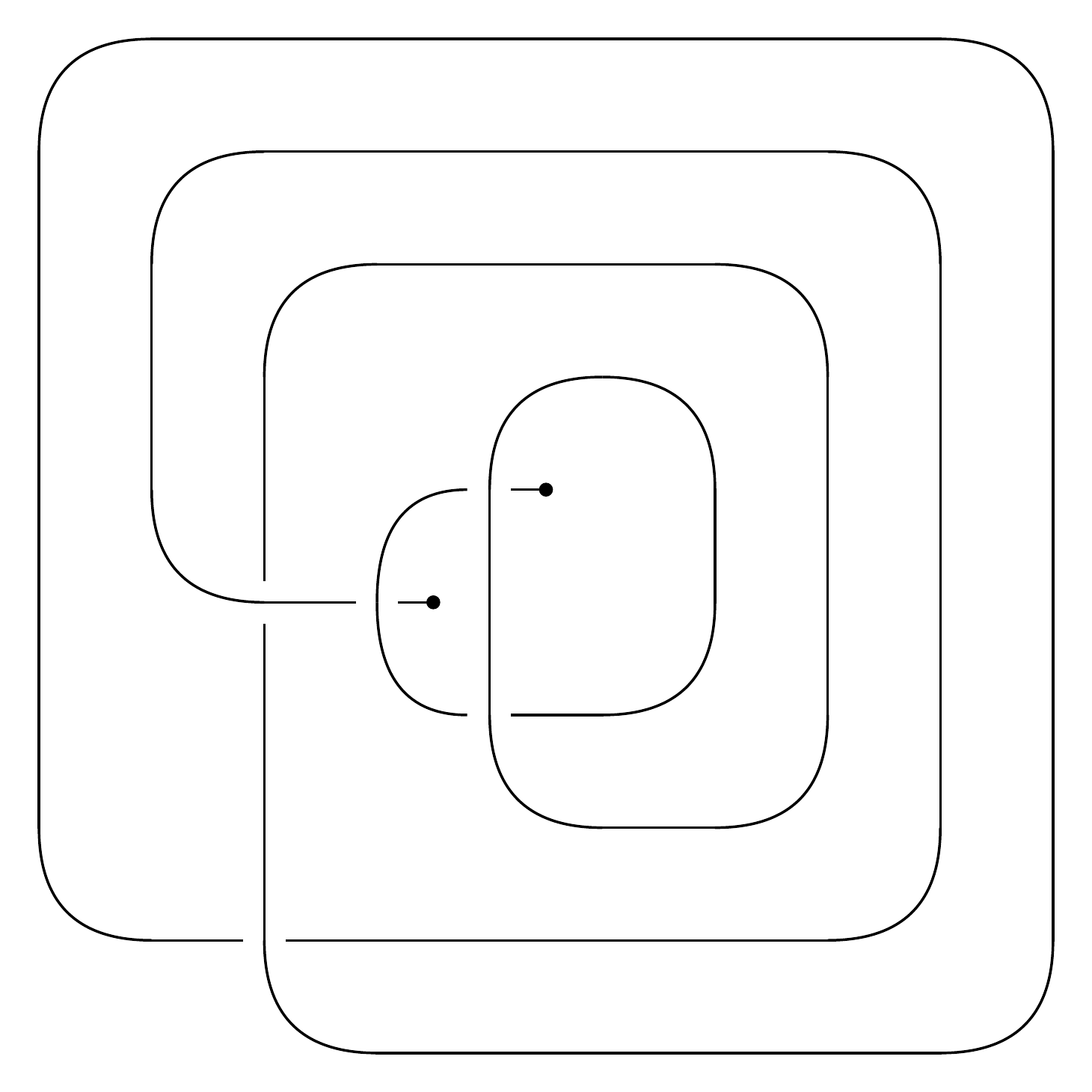}\\
\textcolor{black}{$5_{538}$}
\vspace{1cm}
\end{minipage}
\begin{minipage}[t]{.25\linewidth}
\centering
\includegraphics[width=0.9\textwidth,height=3.5cm,keepaspectratio]{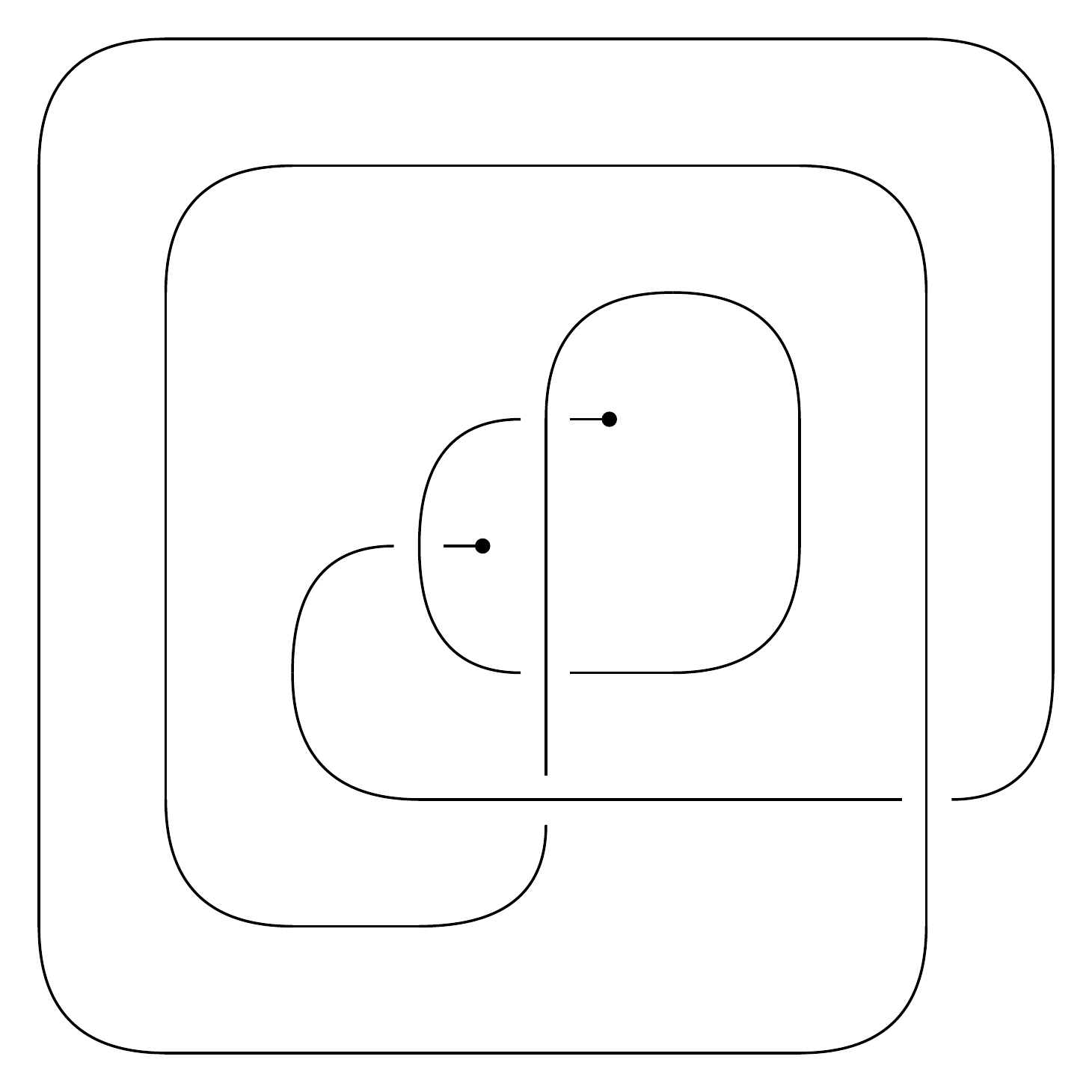}\\
\textcolor{black}{$5_{539}$}
\vspace{1cm}
\end{minipage}
\begin{minipage}[t]{.25\linewidth}
\centering
\includegraphics[width=0.9\textwidth,height=3.5cm,keepaspectratio]{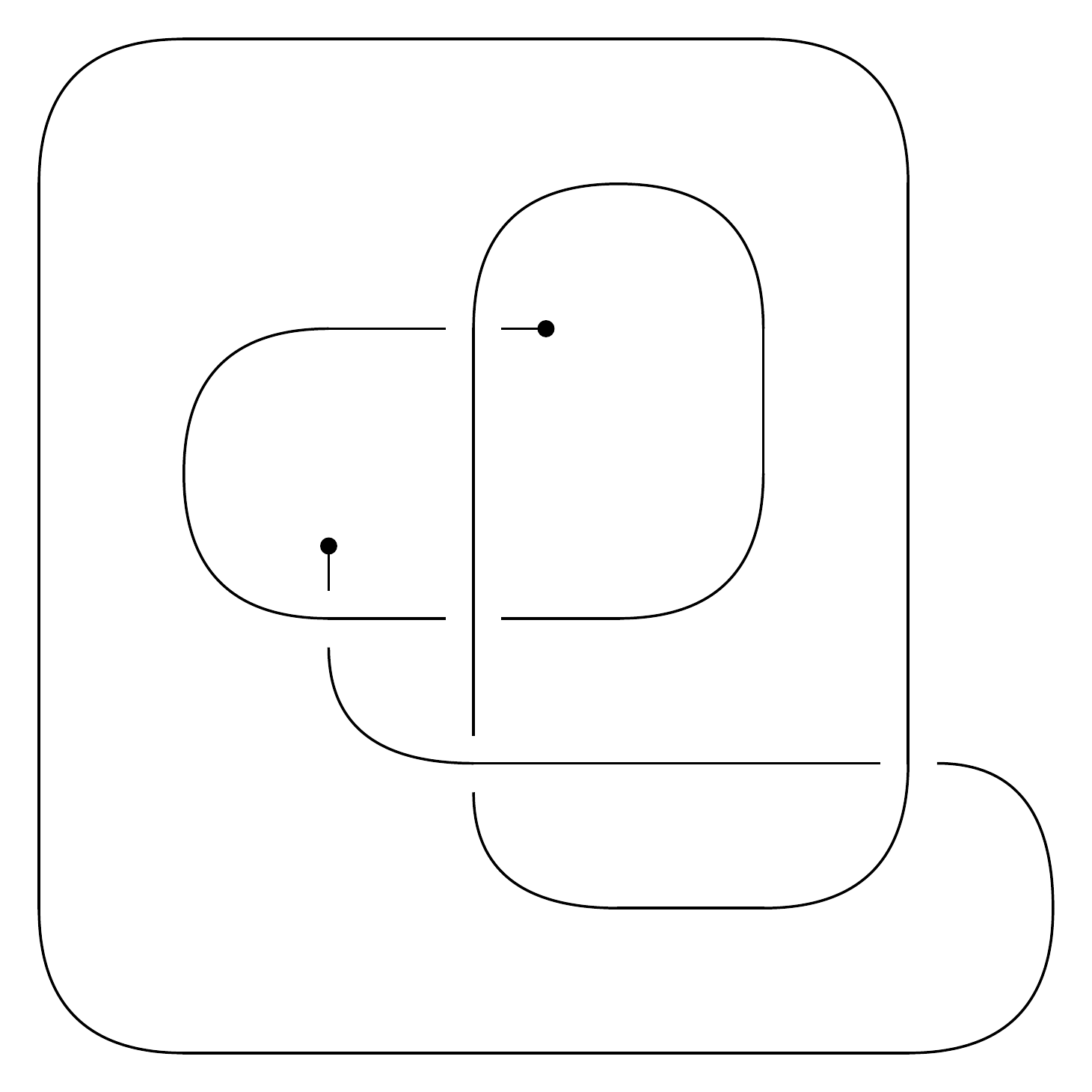}\\
\textcolor{black}{$5_{540}$}
\vspace{1cm}
\end{minipage}
\begin{minipage}[t]{.25\linewidth}
\centering
\includegraphics[width=0.9\textwidth,height=3.5cm,keepaspectratio]{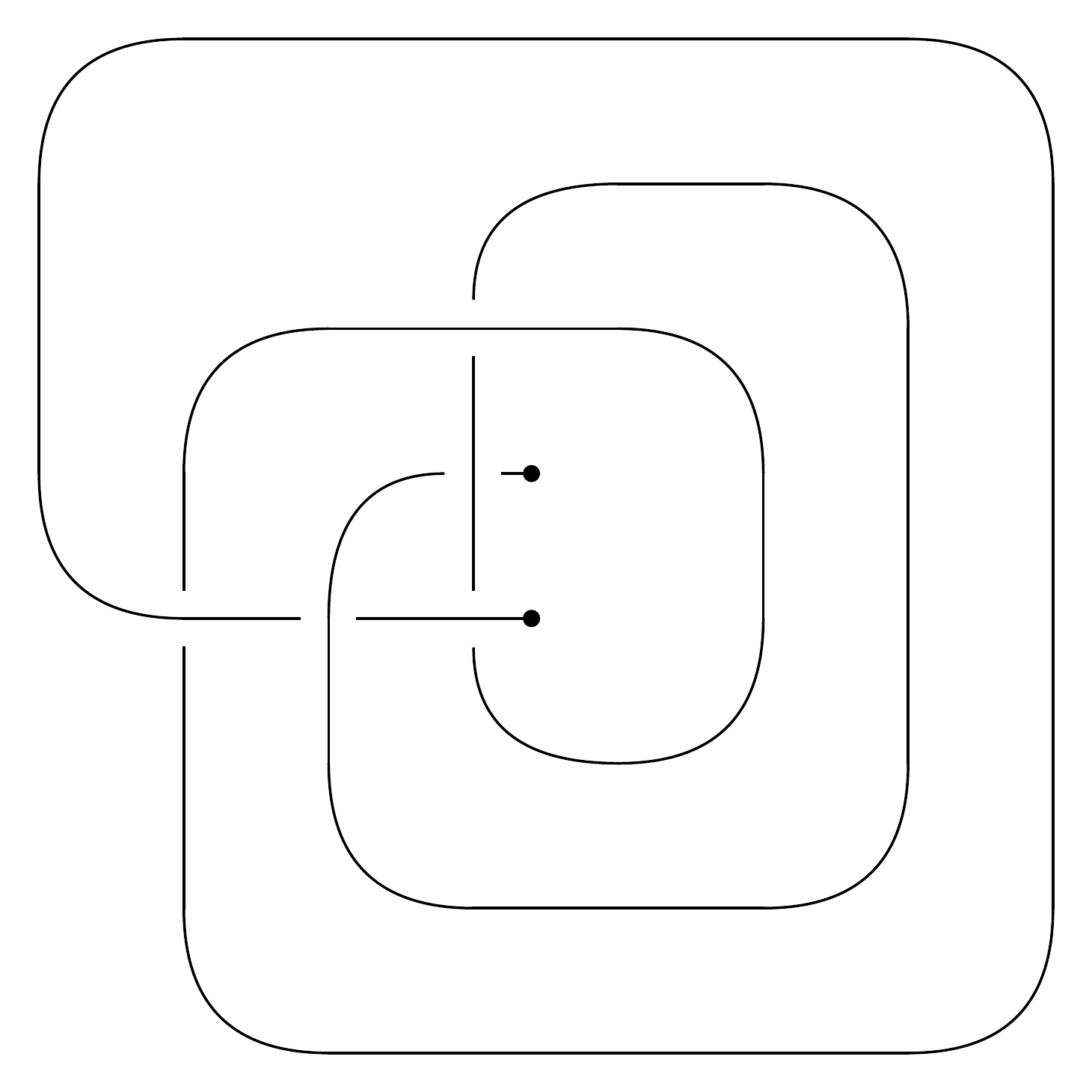}\\
\textcolor{black}{$5_{541}$}
\vspace{1cm}
\end{minipage}
\begin{minipage}[t]{.25\linewidth}
\centering
\includegraphics[width=0.9\textwidth,height=3.5cm,keepaspectratio]{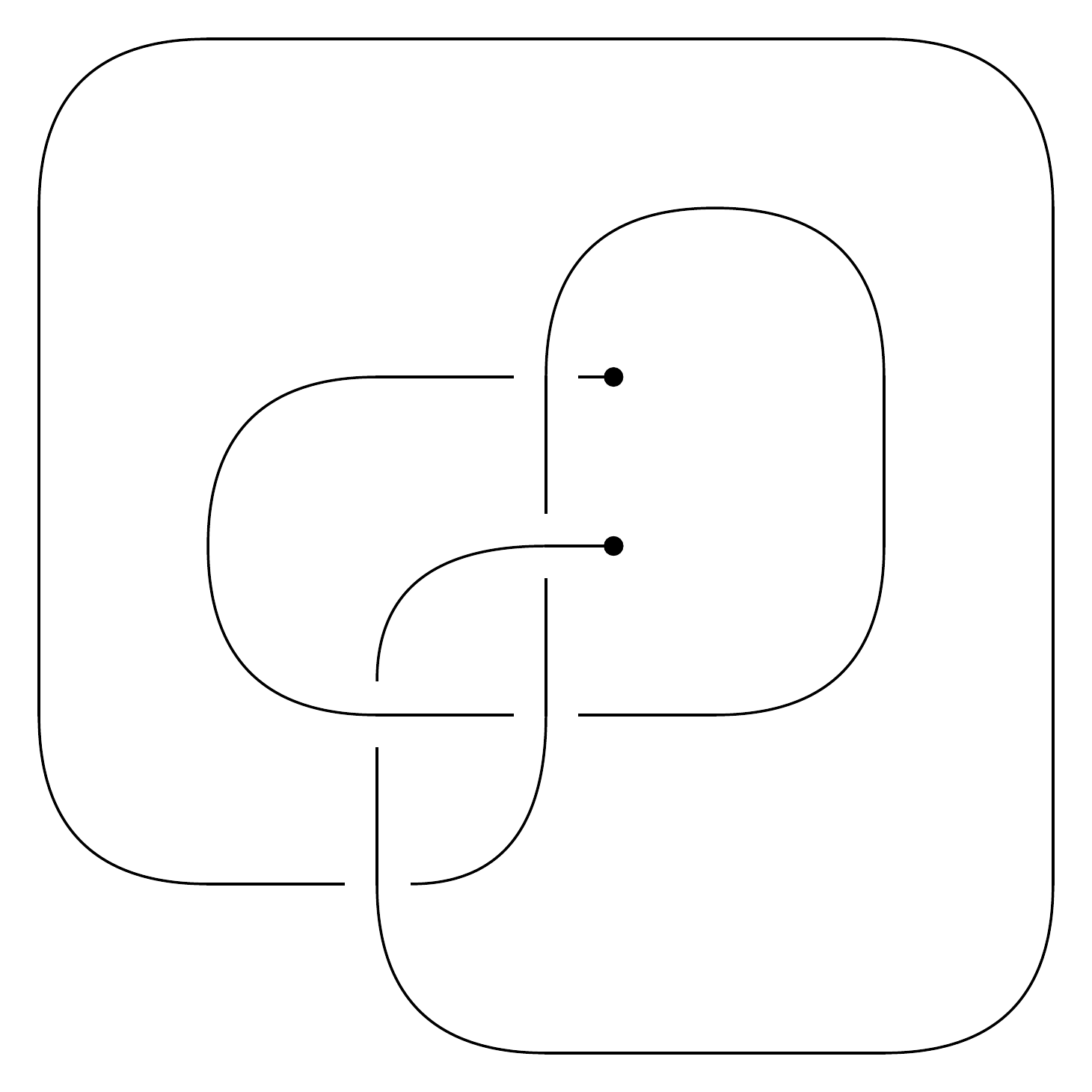}\\
\textcolor{black}{$5_{542}$}
\vspace{1cm}
\end{minipage}
\begin{minipage}[t]{.25\linewidth}
\centering
\includegraphics[width=0.9\textwidth,height=3.5cm,keepaspectratio]{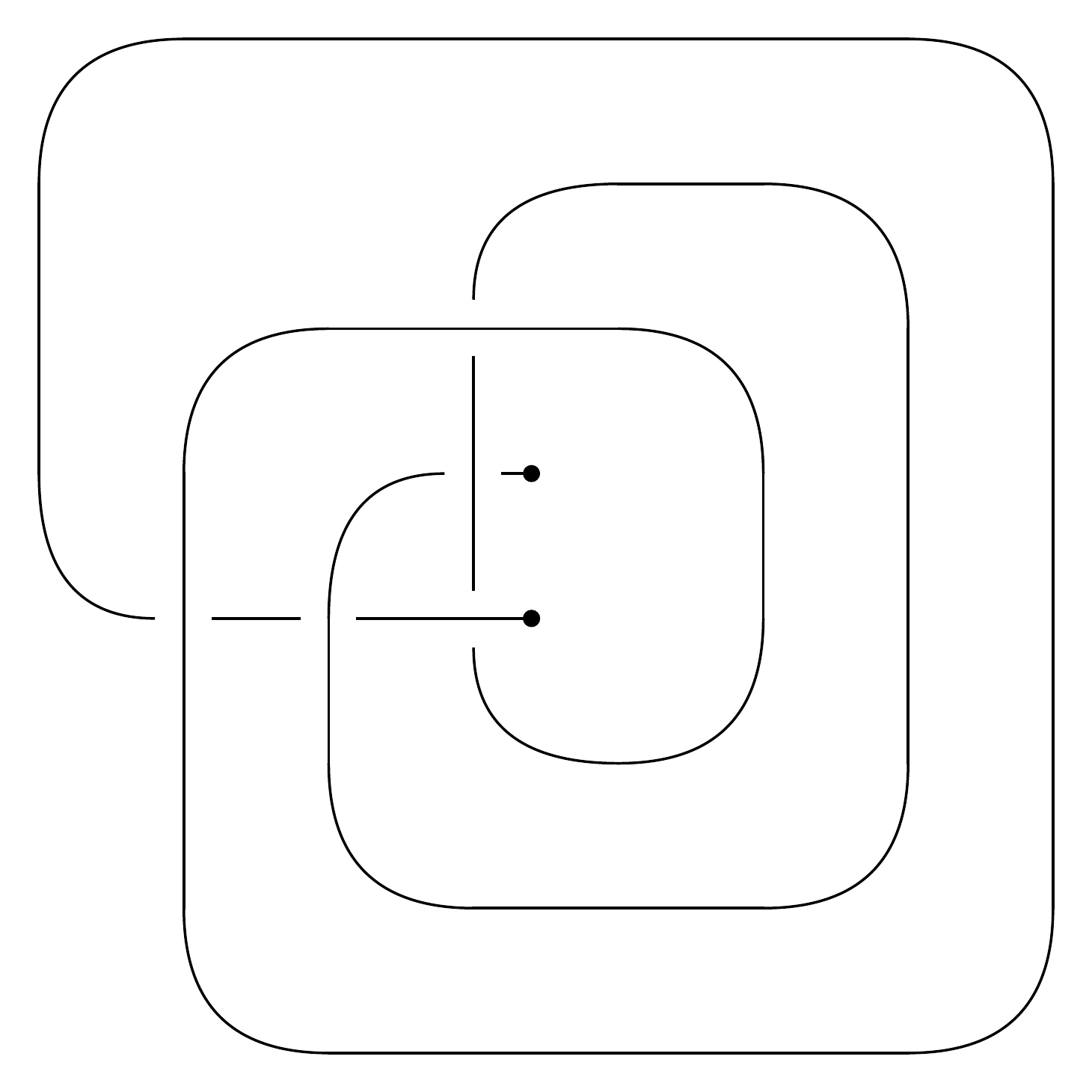}\\
\textcolor{black}{$5_{543}$}
\vspace{1cm}
\end{minipage}
\begin{minipage}[t]{.25\linewidth}
\centering
\includegraphics[width=0.9\textwidth,height=3.5cm,keepaspectratio]{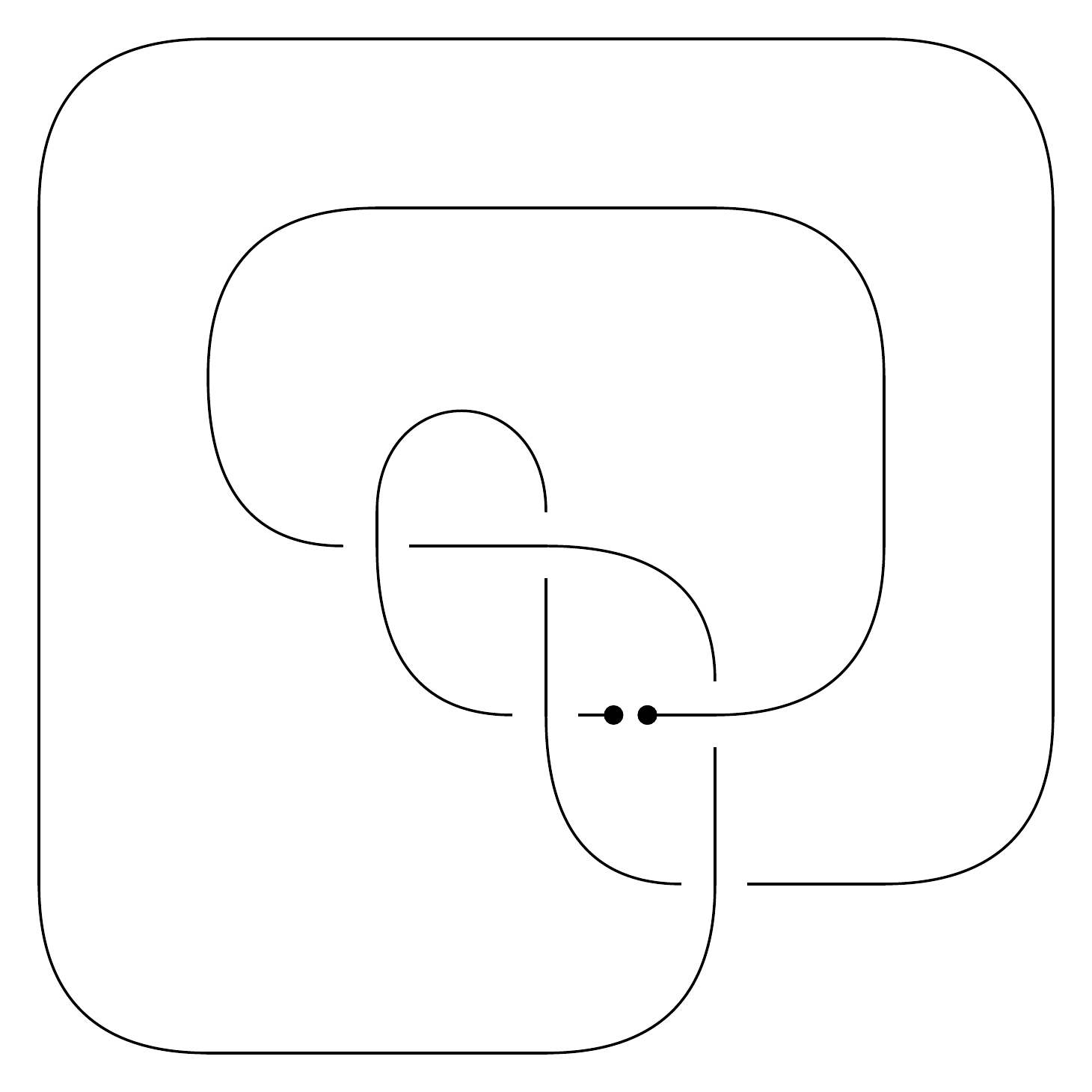}\\
\textcolor{black}{$5_{544}$}
\vspace{1cm}
\end{minipage}
\begin{minipage}[t]{.25\linewidth}
\centering
\includegraphics[width=0.9\textwidth,height=3.5cm,keepaspectratio]{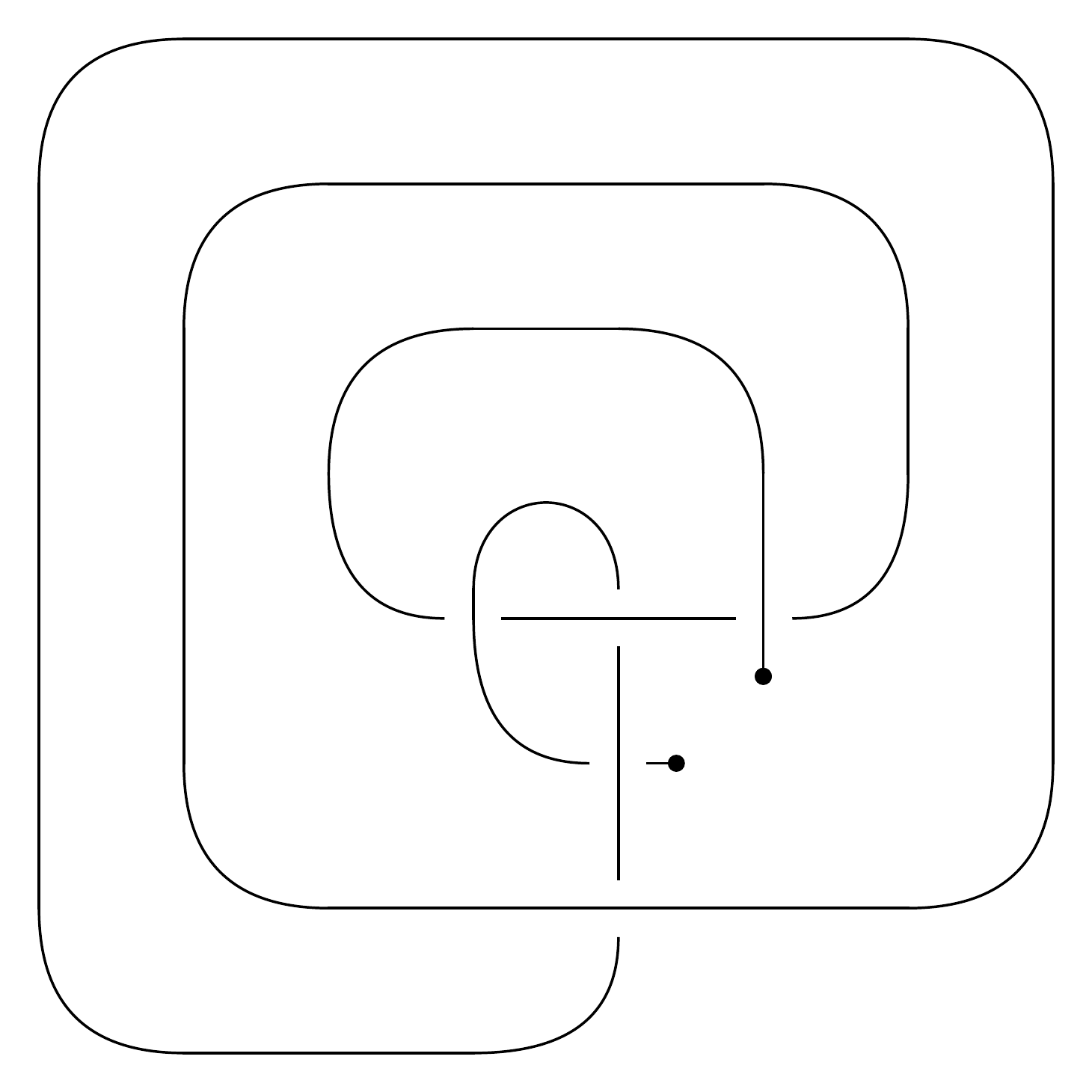}\\
\textcolor{black}{$5_{545}$}
\vspace{1cm}
\end{minipage}
\begin{minipage}[t]{.25\linewidth}
\centering
\includegraphics[width=0.9\textwidth,height=3.5cm,keepaspectratio]{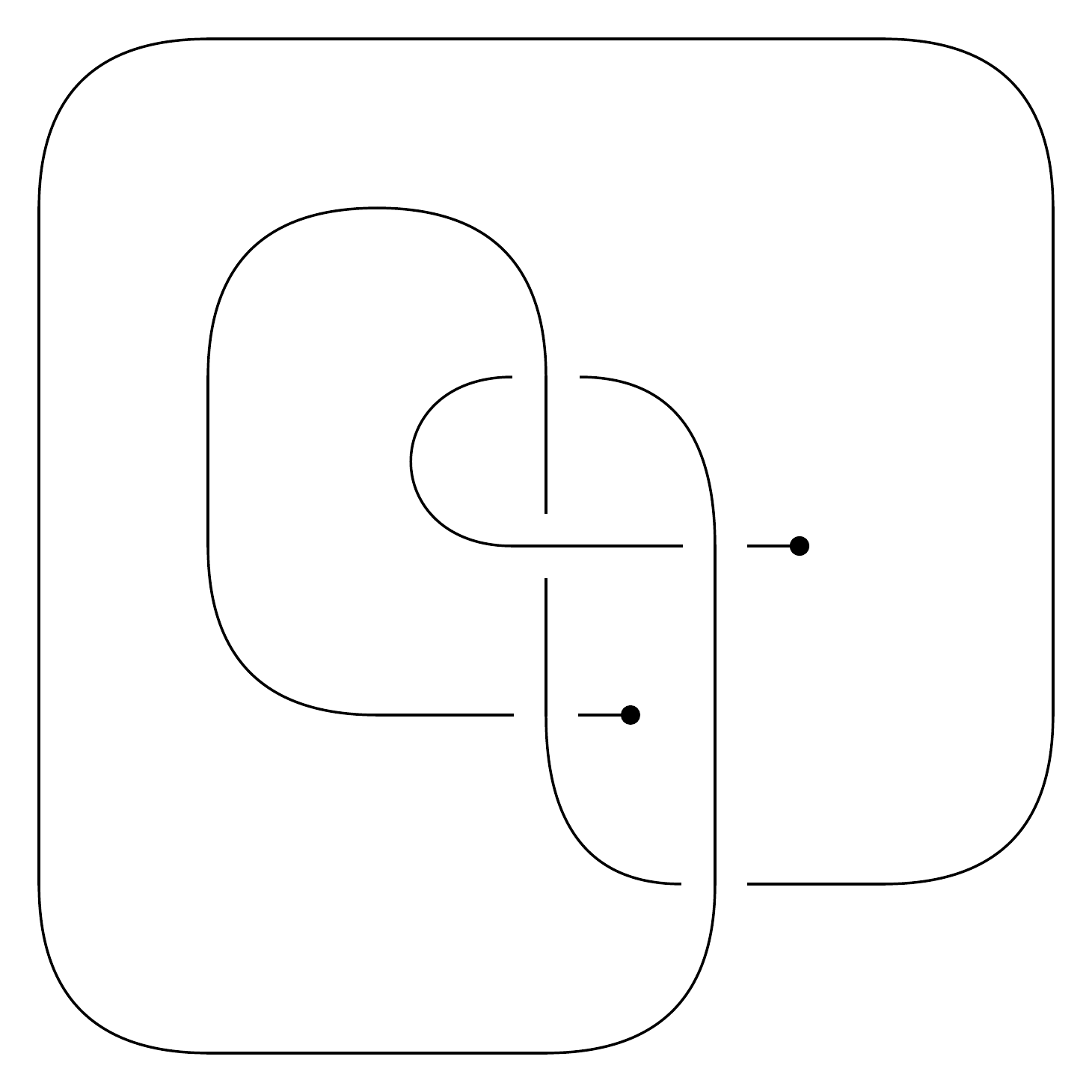}\\
\textcolor{black}{$5_{546}$}
\vspace{1cm}
\end{minipage}
\begin{minipage}[t]{.25\linewidth}
\centering
\includegraphics[width=0.9\textwidth,height=3.5cm,keepaspectratio]{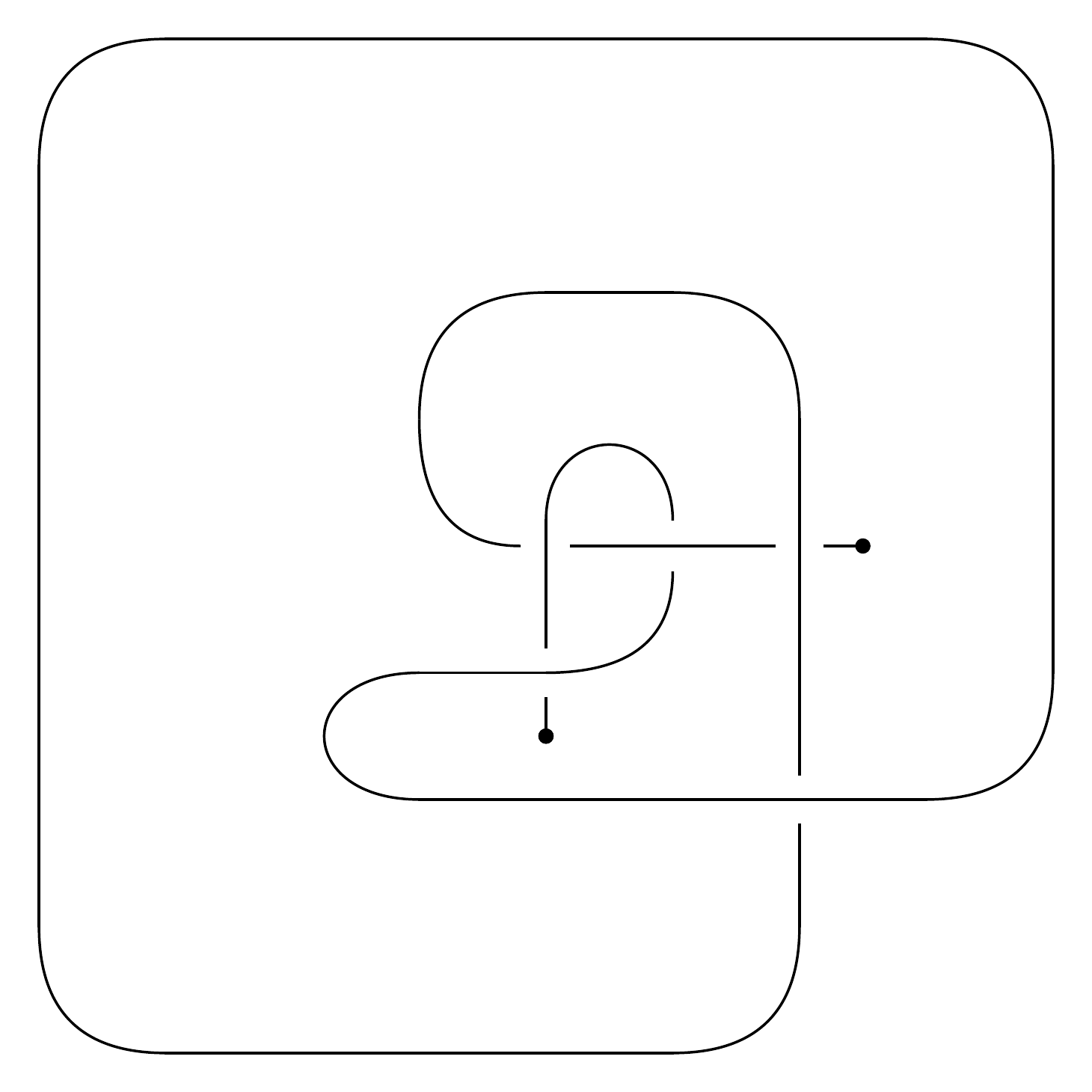}\\
\textcolor{black}{$5_{547}$}
\vspace{1cm}
\end{minipage}
\begin{minipage}[t]{.25\linewidth}
\centering
\includegraphics[width=0.9\textwidth,height=3.5cm,keepaspectratio]{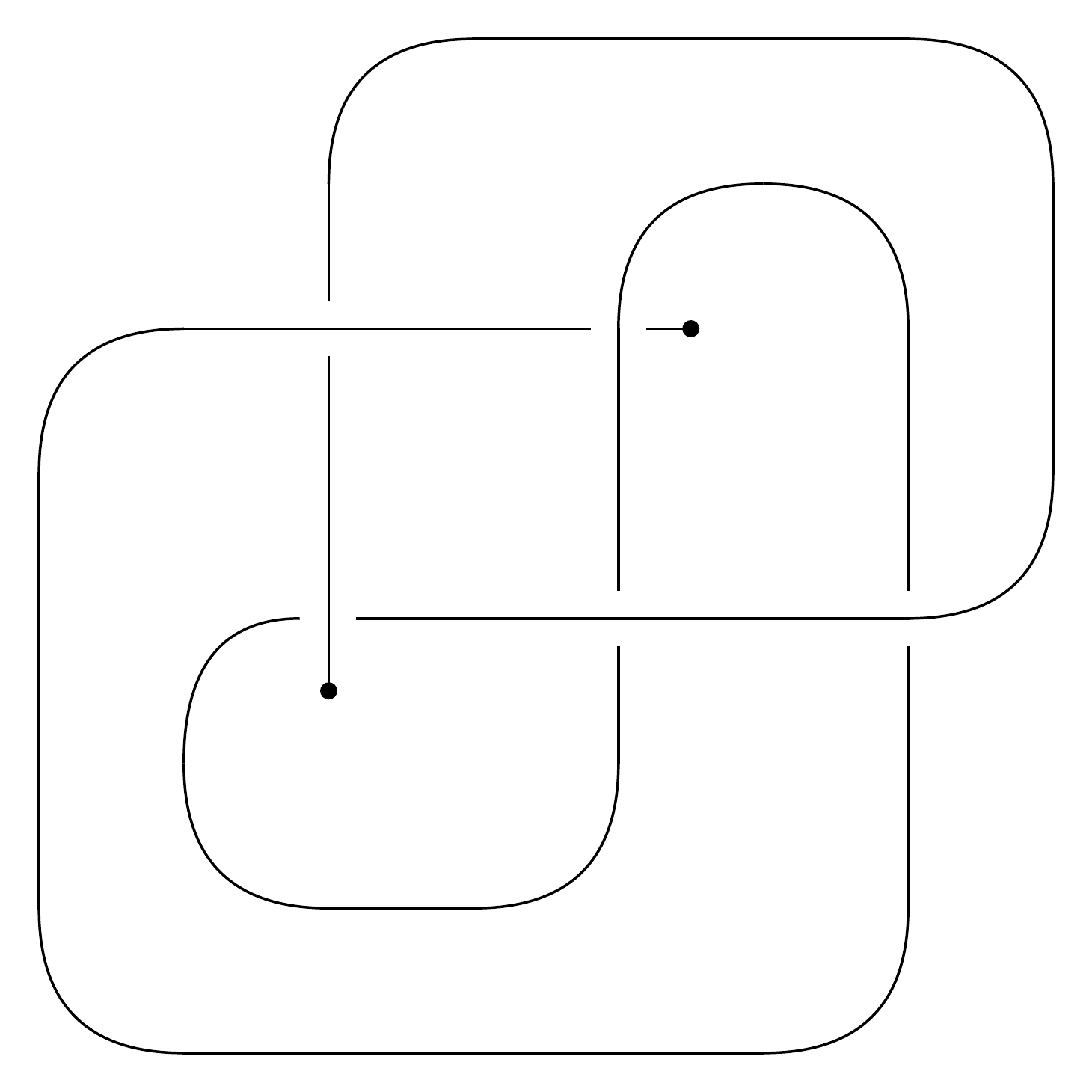}\\
\textcolor{black}{$5_{548}$}
\vspace{1cm}
\end{minipage}
\begin{minipage}[t]{.25\linewidth}
\centering
\includegraphics[width=0.9\textwidth,height=3.5cm,keepaspectratio]{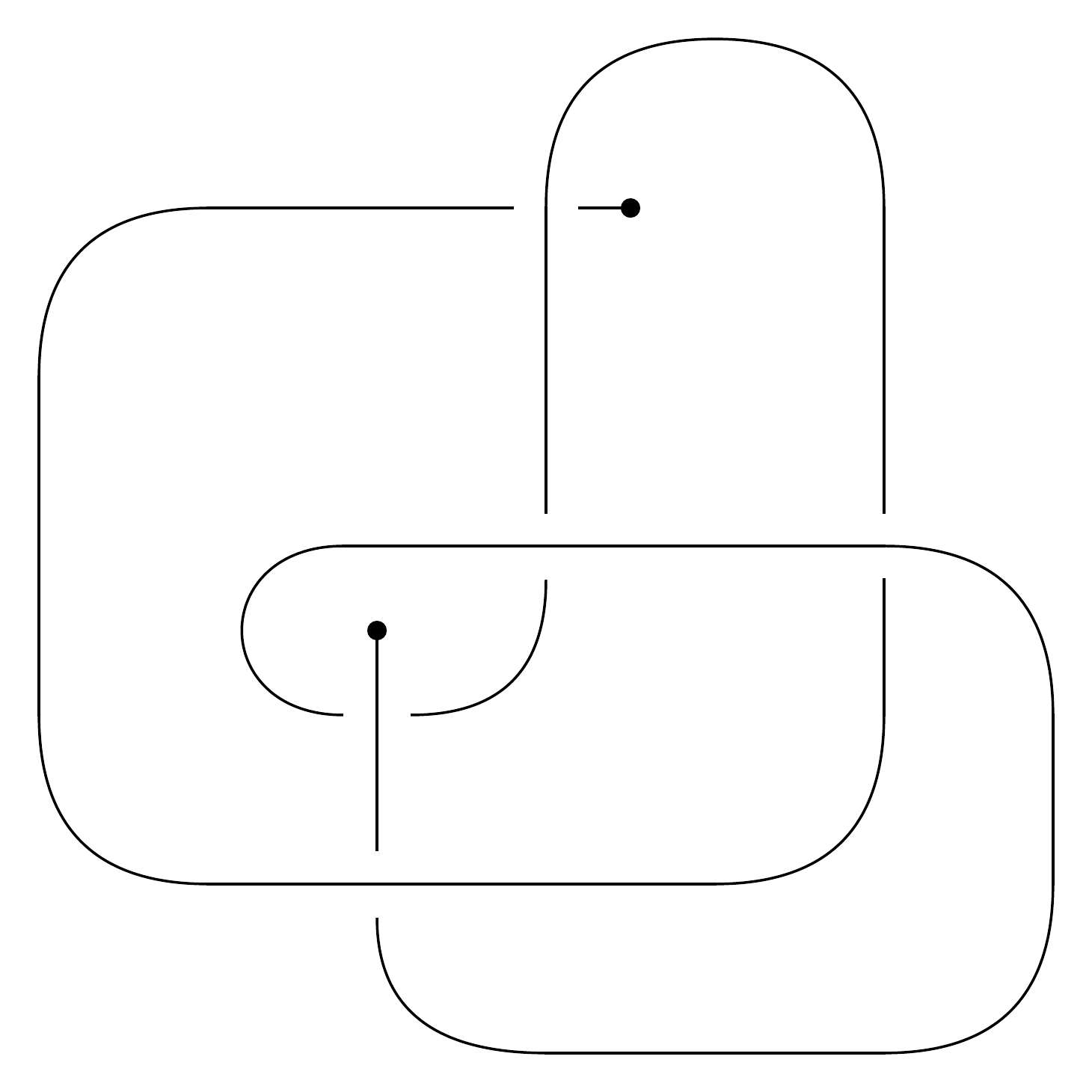}\\
\textcolor{black}{$5_{549}$}
\vspace{1cm}
\end{minipage}
\begin{minipage}[t]{.25\linewidth}
\centering
\includegraphics[width=0.9\textwidth,height=3.5cm,keepaspectratio]{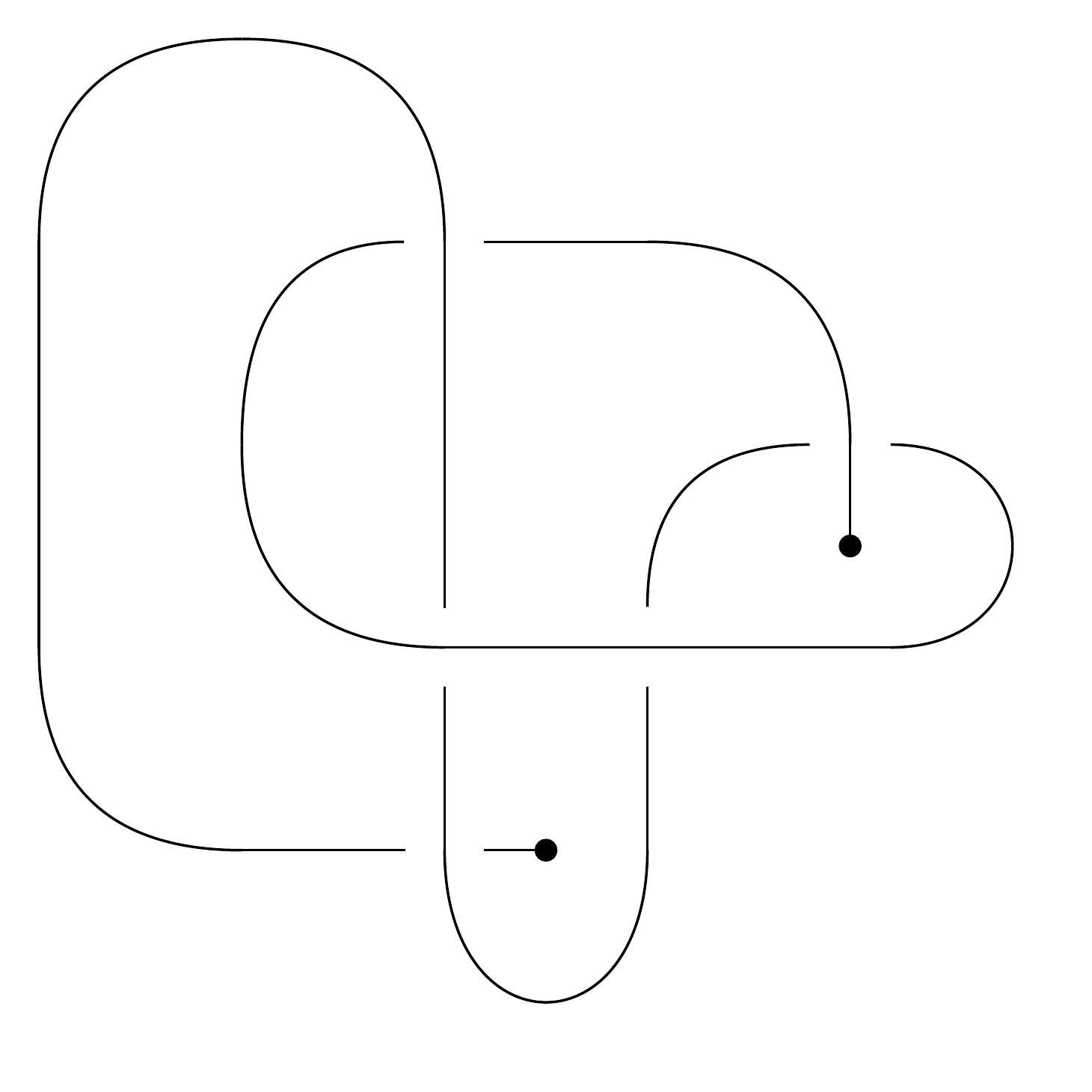}\\
\textcolor{black}{$5_{550}$}
\vspace{1cm}
\end{minipage}
\begin{minipage}[t]{.25\linewidth}
\centering
\includegraphics[width=0.9\textwidth,height=3.5cm,keepaspectratio]{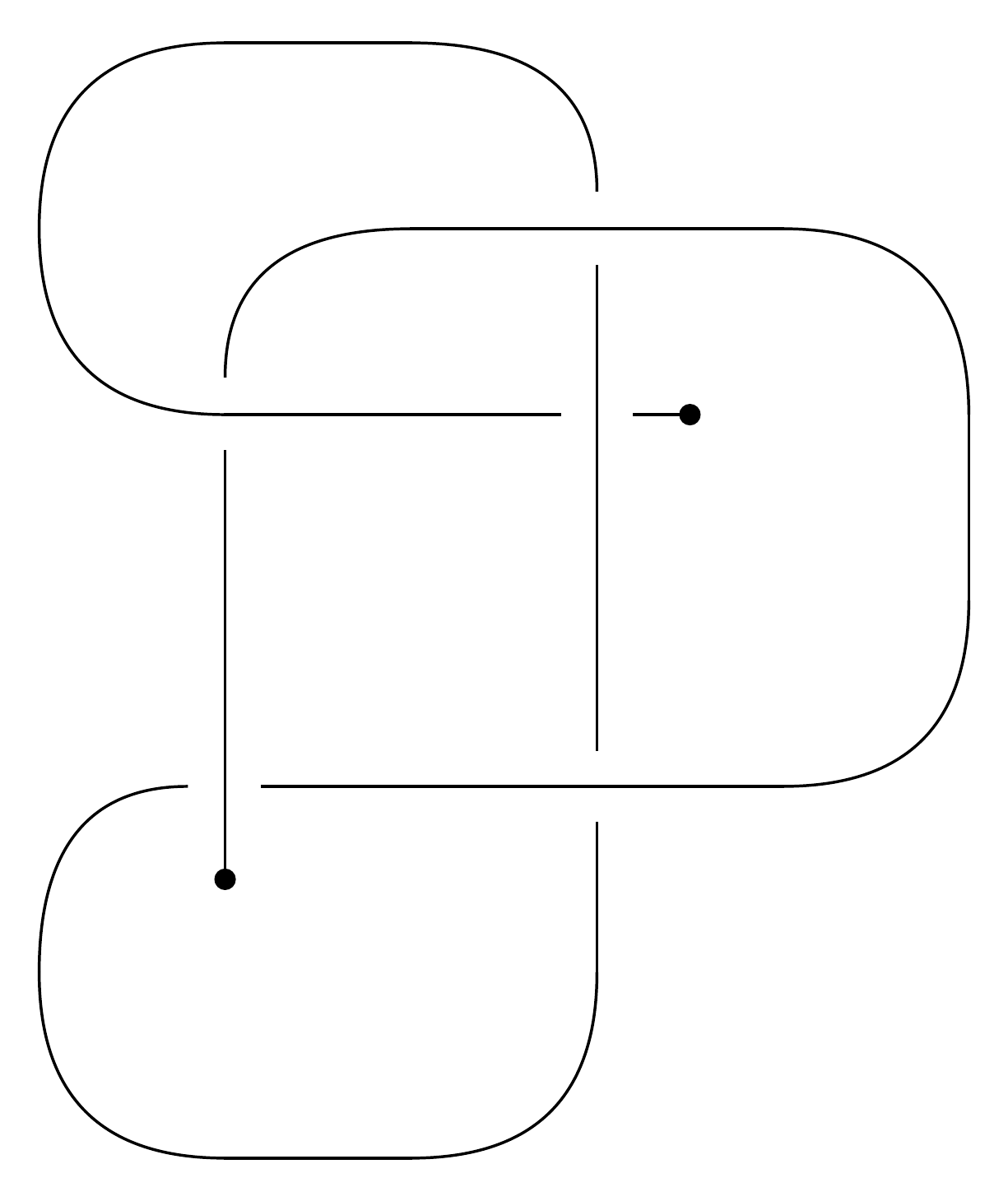}\\
\textcolor{black}{$5_{551}$}
\vspace{1cm}
\end{minipage}
\begin{minipage}[t]{.25\linewidth}
\centering
\includegraphics[width=0.9\textwidth,height=3.5cm,keepaspectratio]{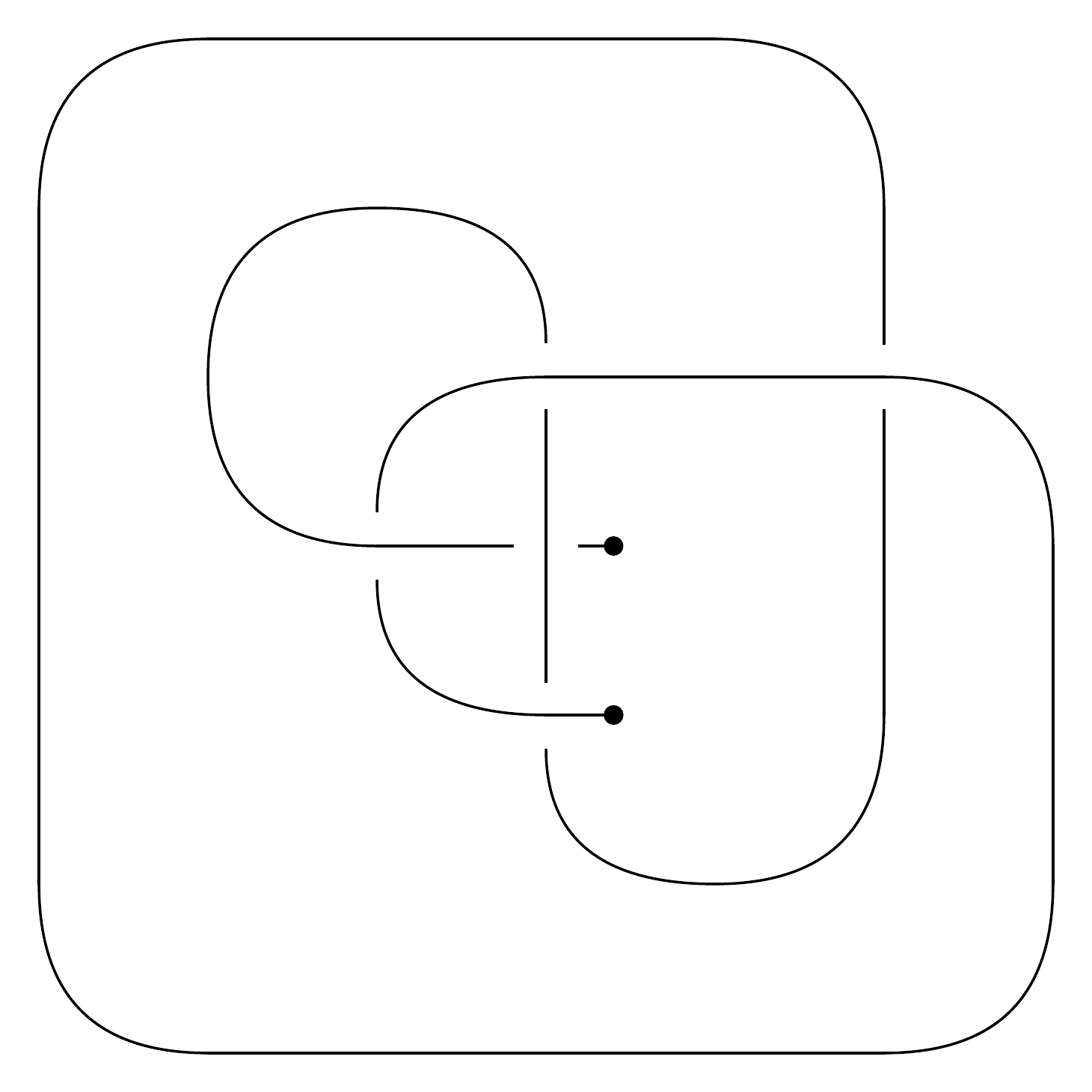}\\
\textcolor{black}{$5_{552}$}
\vspace{1cm}
\end{minipage}
\begin{minipage}[t]{.25\linewidth}
\centering
\includegraphics[width=0.9\textwidth,height=3.5cm,keepaspectratio]{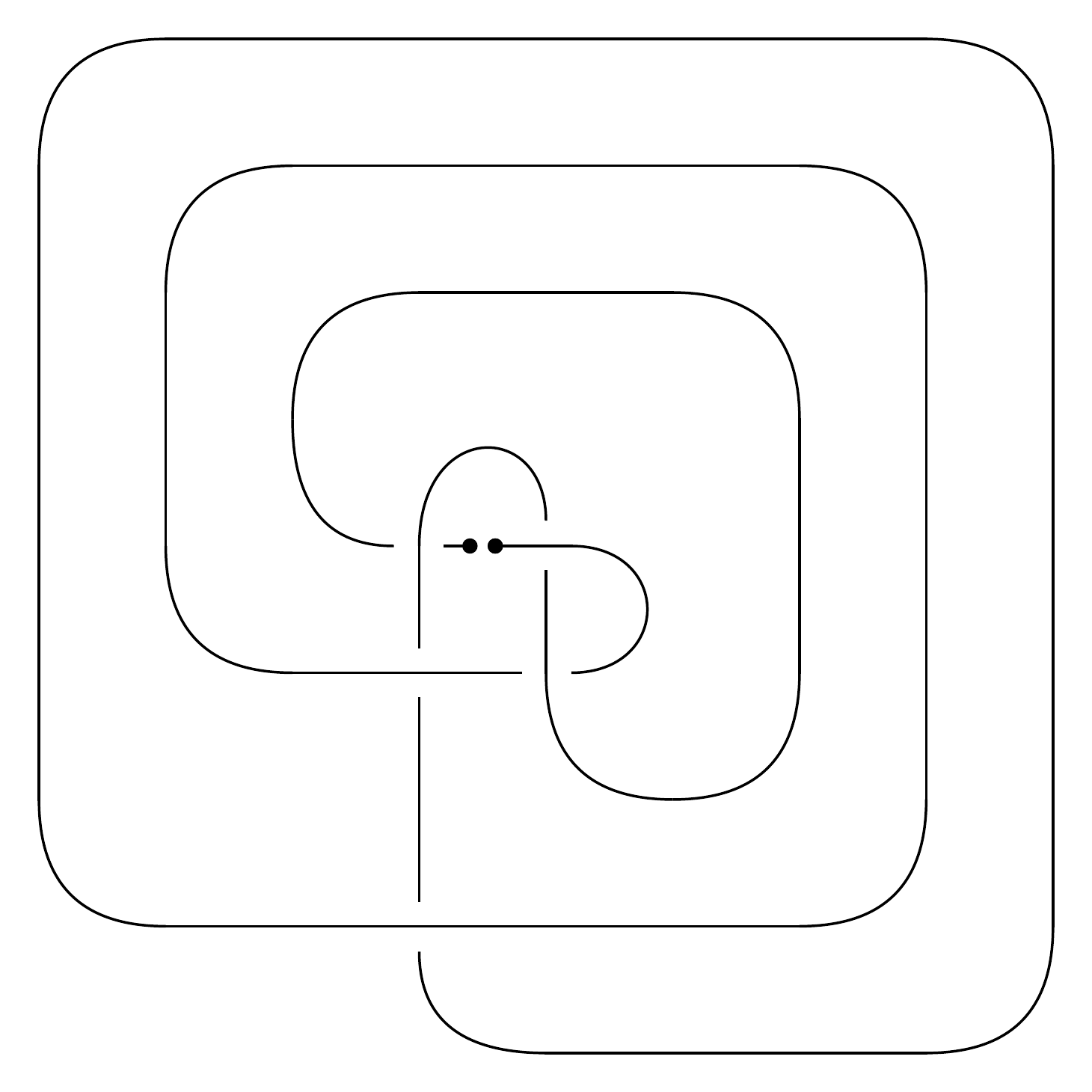}\\
\textcolor{black}{$5_{553}$}
\vspace{1cm}
\end{minipage}
\begin{minipage}[t]{.25\linewidth}
\centering
\includegraphics[width=0.9\textwidth,height=3.5cm,keepaspectratio]{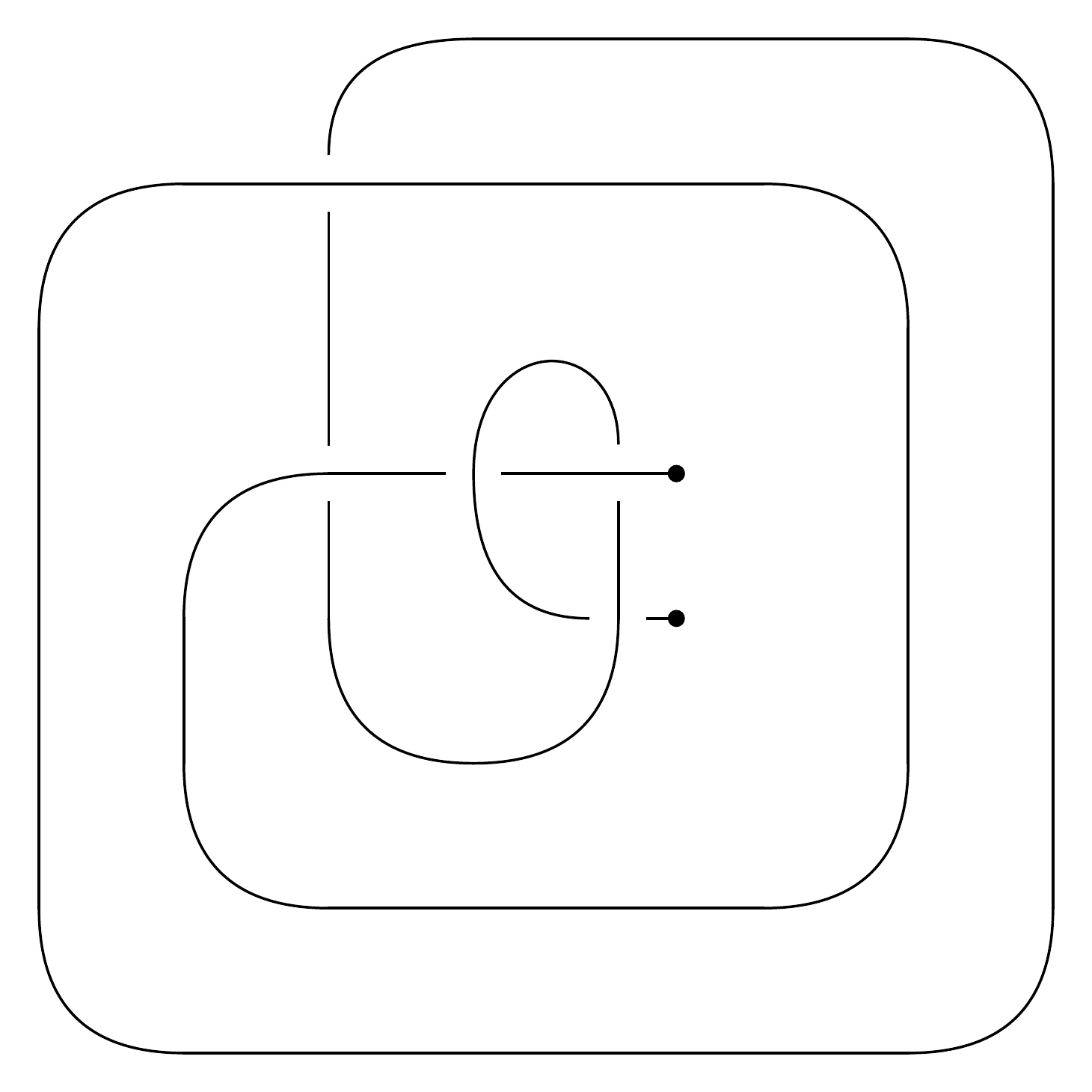}\\
\textcolor{black}{$5_{554}$}
\vspace{1cm}
\end{minipage}
\begin{minipage}[t]{.25\linewidth}
\centering
\includegraphics[width=0.9\textwidth,height=3.5cm,keepaspectratio]{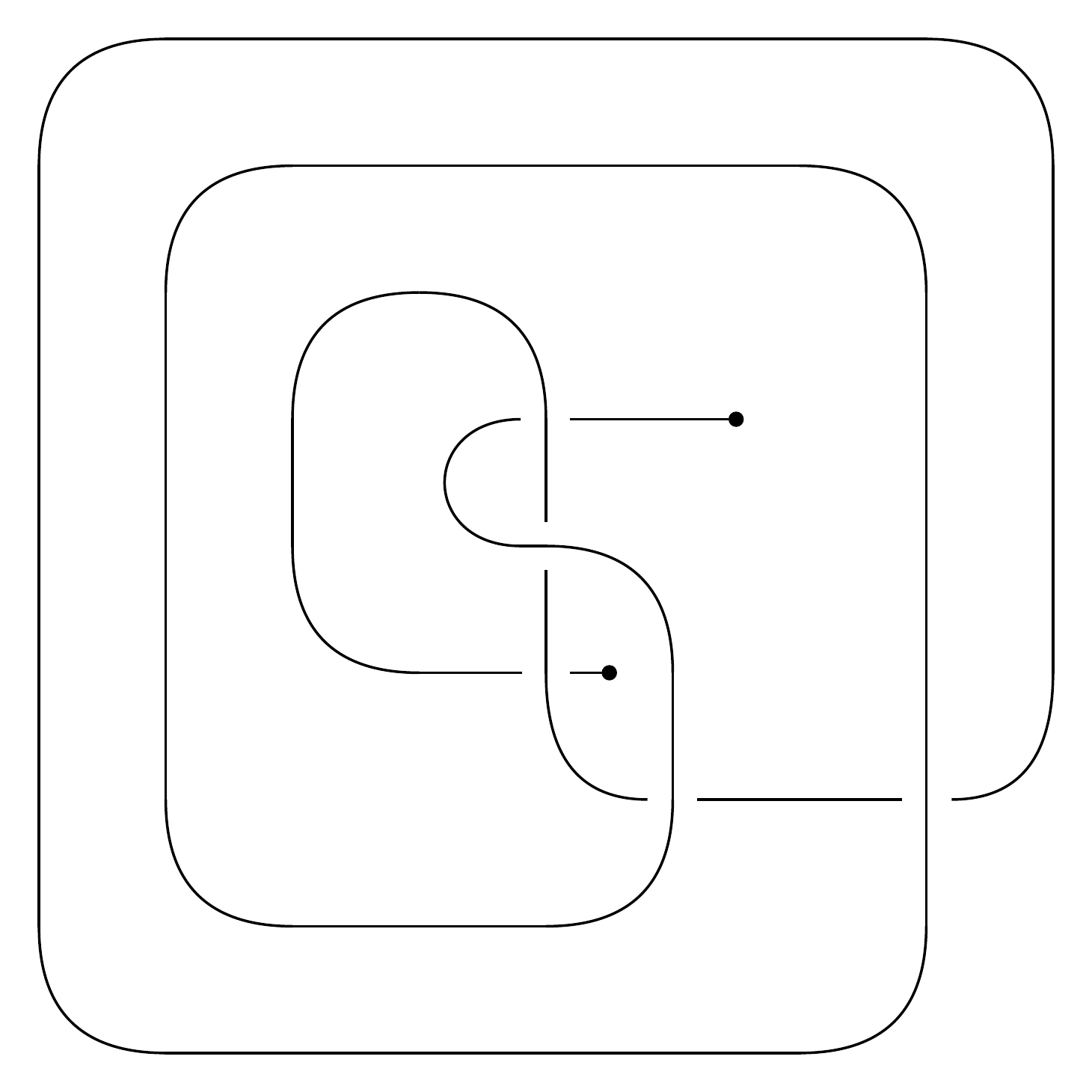}\\
\textcolor{black}{$5_{555}$}
\vspace{1cm}
\end{minipage}
\begin{minipage}[t]{.25\linewidth}
\centering
\includegraphics[width=0.9\textwidth,height=3.5cm,keepaspectratio]{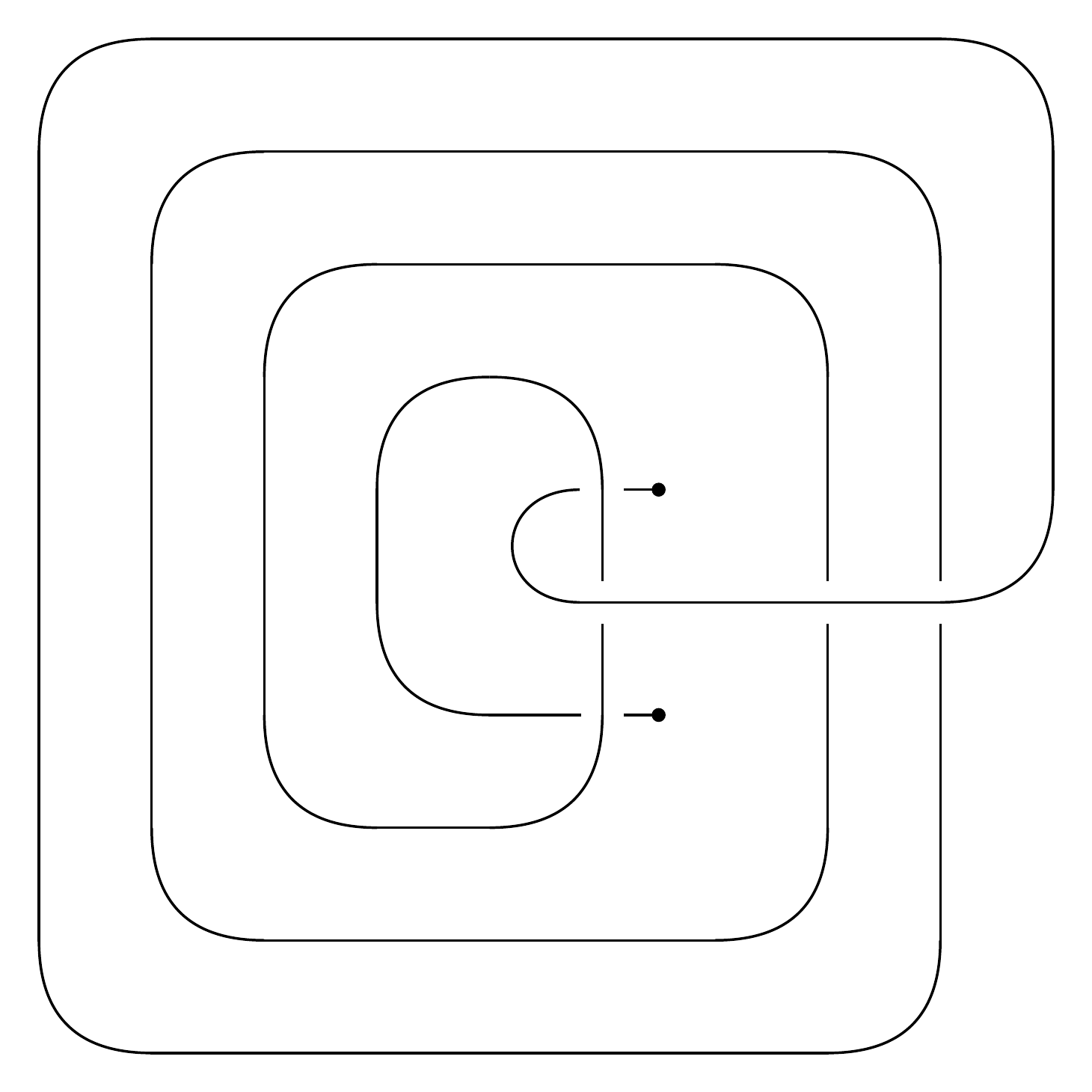}\\
\textcolor{black}{$5_{556}$}
\vspace{1cm}
\end{minipage}
\begin{minipage}[t]{.25\linewidth}
\centering
\includegraphics[width=0.9\textwidth,height=3.5cm,keepaspectratio]{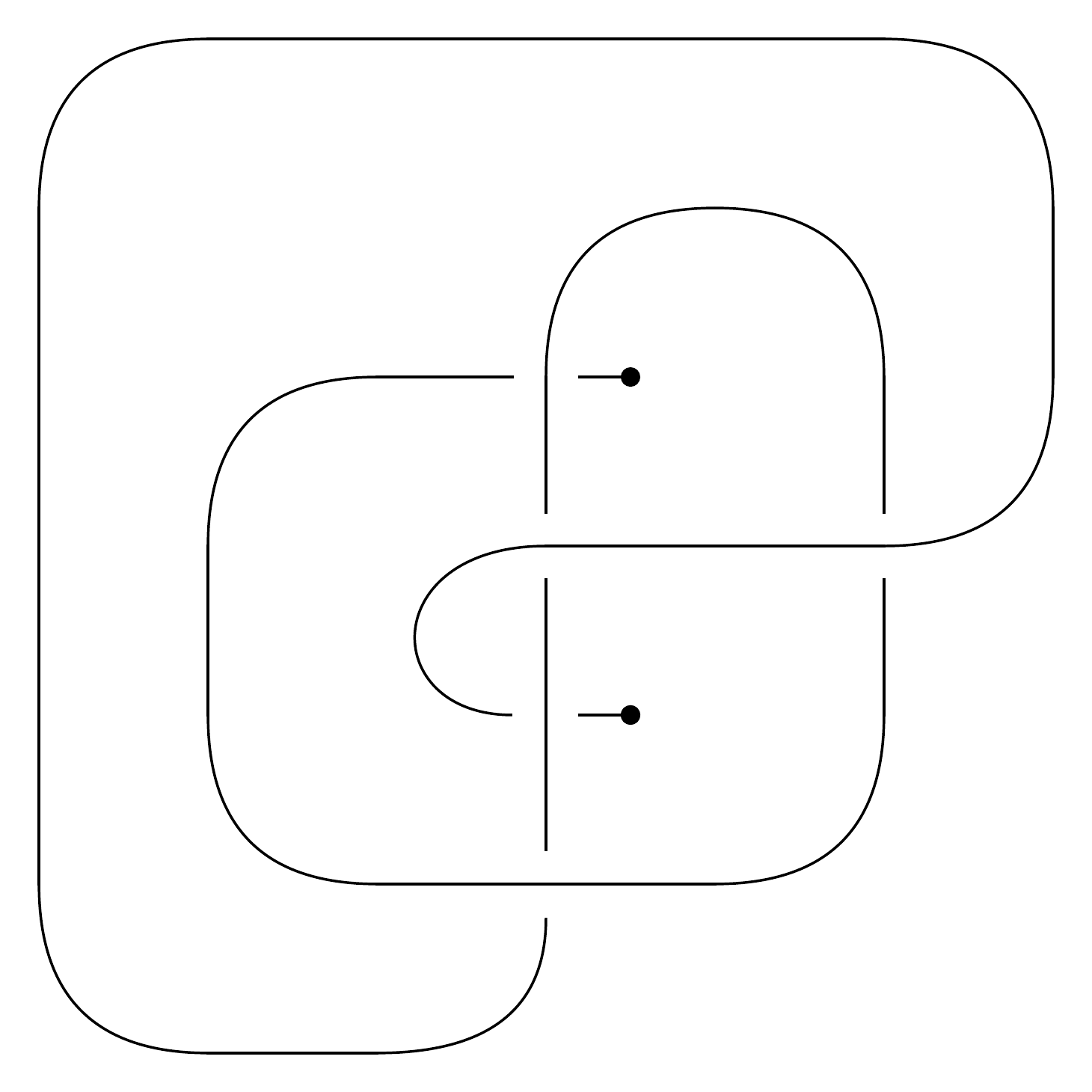}\\
\textcolor{black}{$5_{557}$}
\vspace{1cm}
\end{minipage}
\begin{minipage}[t]{.25\linewidth}
\centering
\includegraphics[width=0.9\textwidth,height=3.5cm,keepaspectratio]{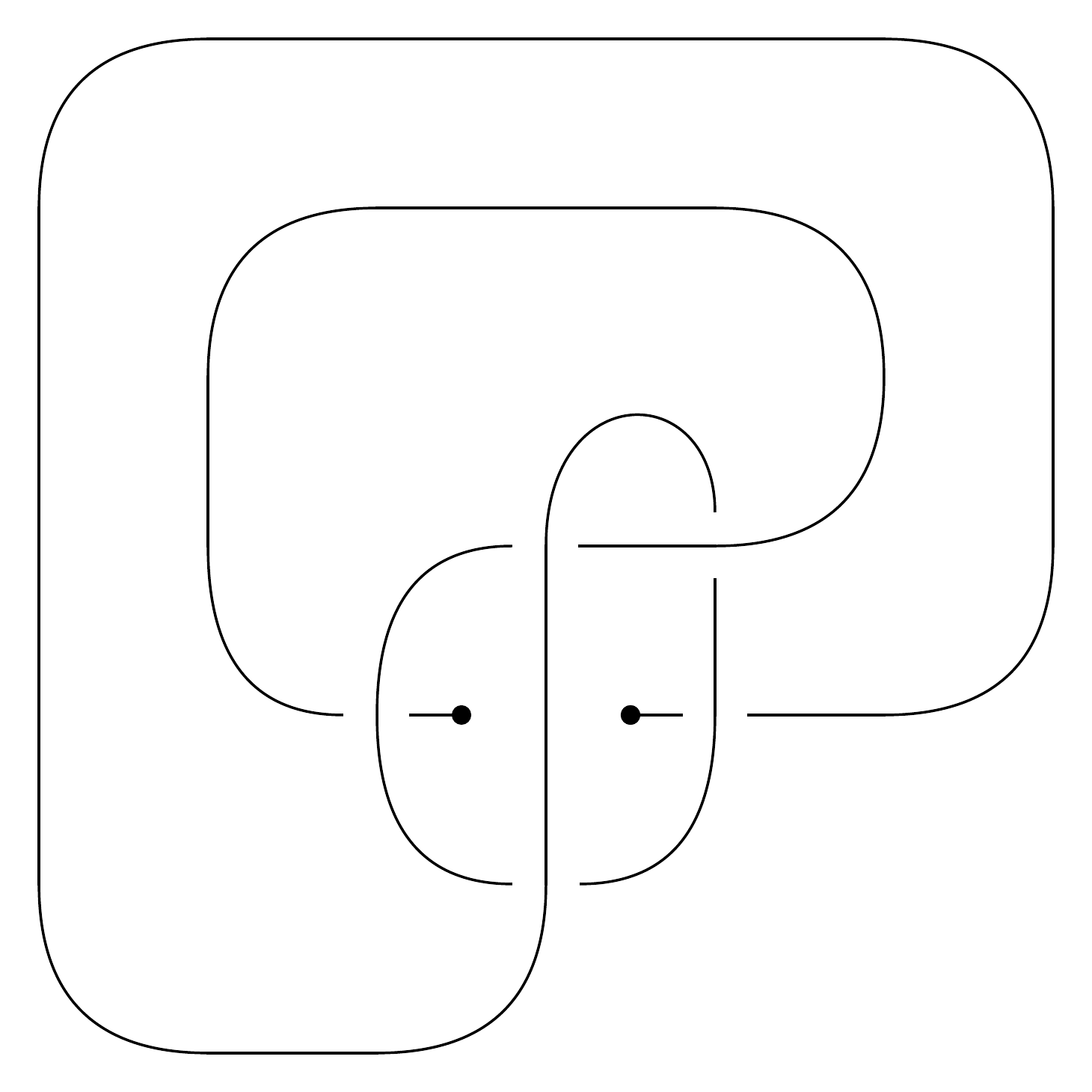}\\
\textcolor{black}{$5_{558}$}
\vspace{1cm}
\end{minipage}
\begin{minipage}[t]{.25\linewidth}
\centering
\includegraphics[width=0.9\textwidth,height=3.5cm,keepaspectratio]{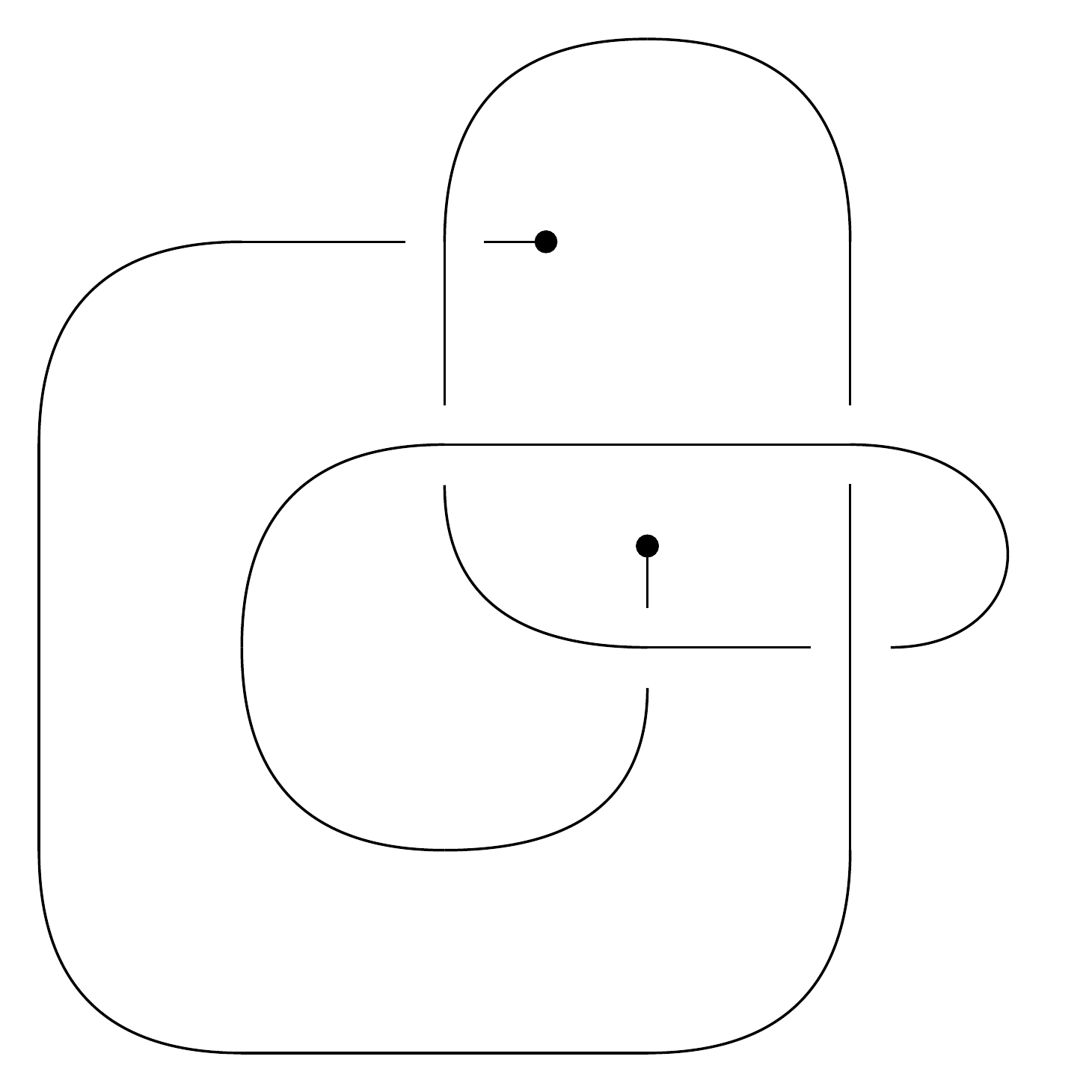}\\
\textcolor{black}{$5_{559}$}
\vspace{1cm}
\end{minipage}
\begin{minipage}[t]{.25\linewidth}
\centering
\includegraphics[width=0.9\textwidth,height=3.5cm,keepaspectratio]{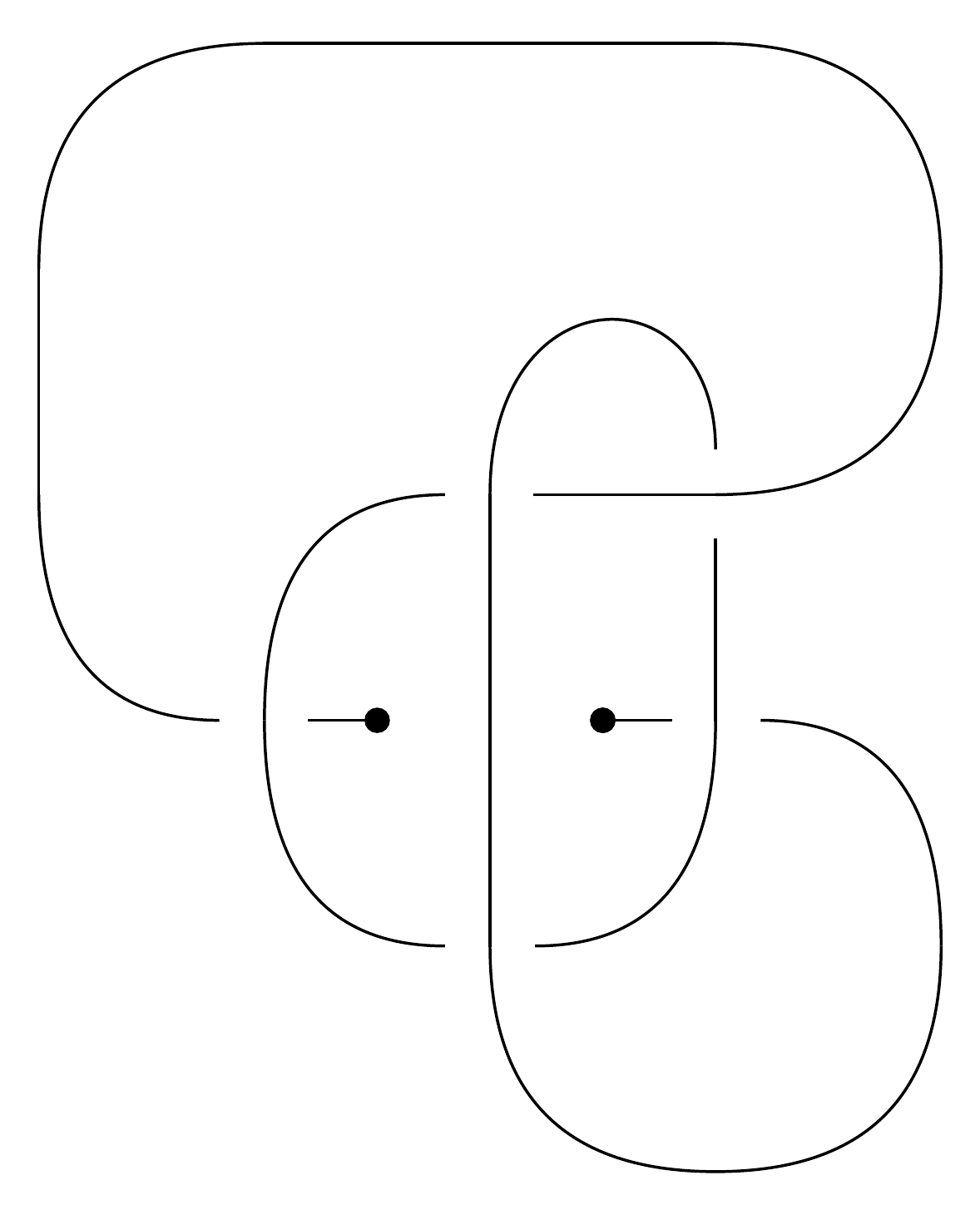}\\
\textcolor{black}{$5_{560}$}
\vspace{1cm}
\end{minipage}
\begin{minipage}[t]{.25\linewidth}
\centering
\includegraphics[width=0.9\textwidth,height=3.5cm,keepaspectratio]{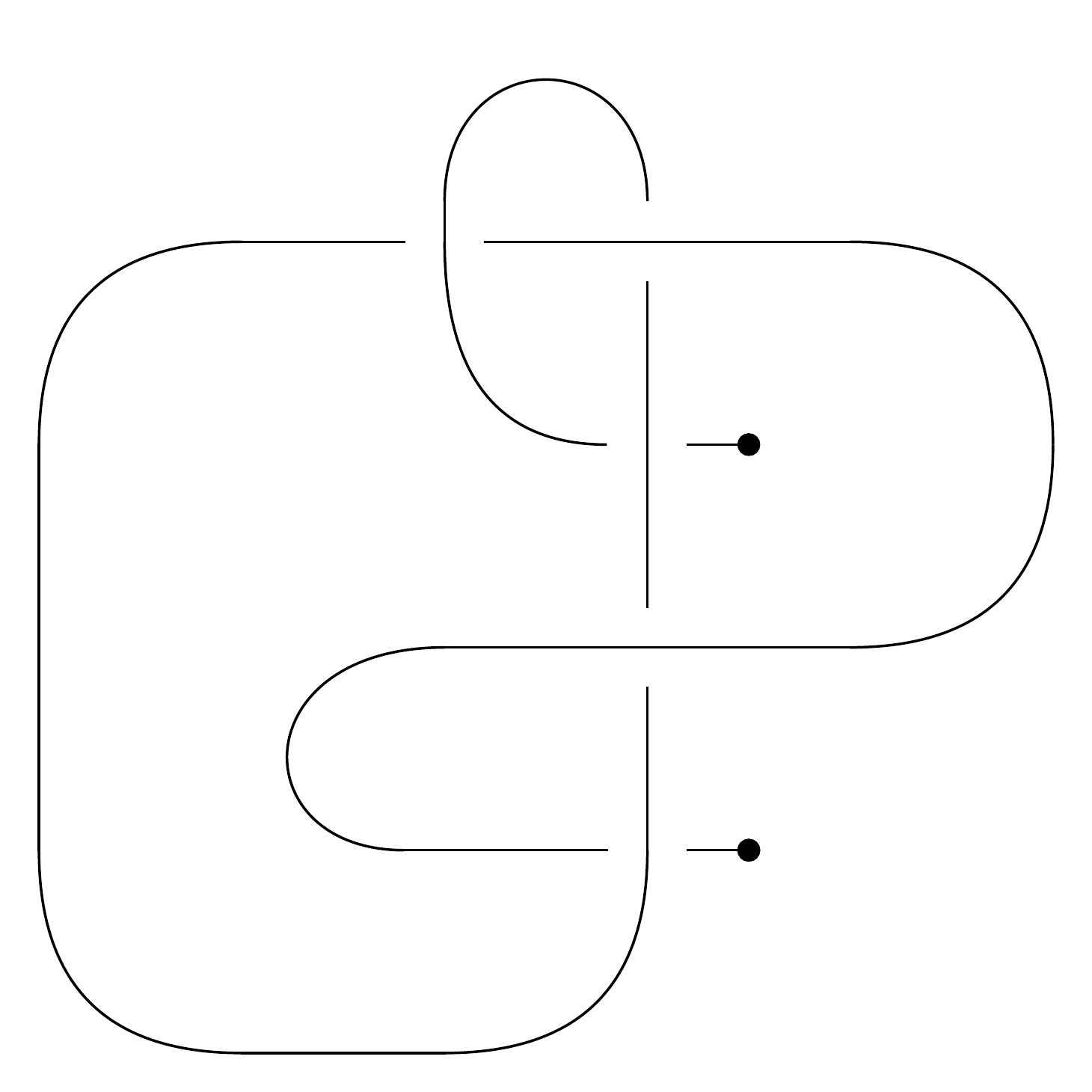}\\
\textcolor{black}{$5_{561}$}
\vspace{1cm}
\end{minipage}
\begin{minipage}[t]{.25\linewidth}
\centering
\includegraphics[width=0.9\textwidth,height=3.5cm,keepaspectratio]{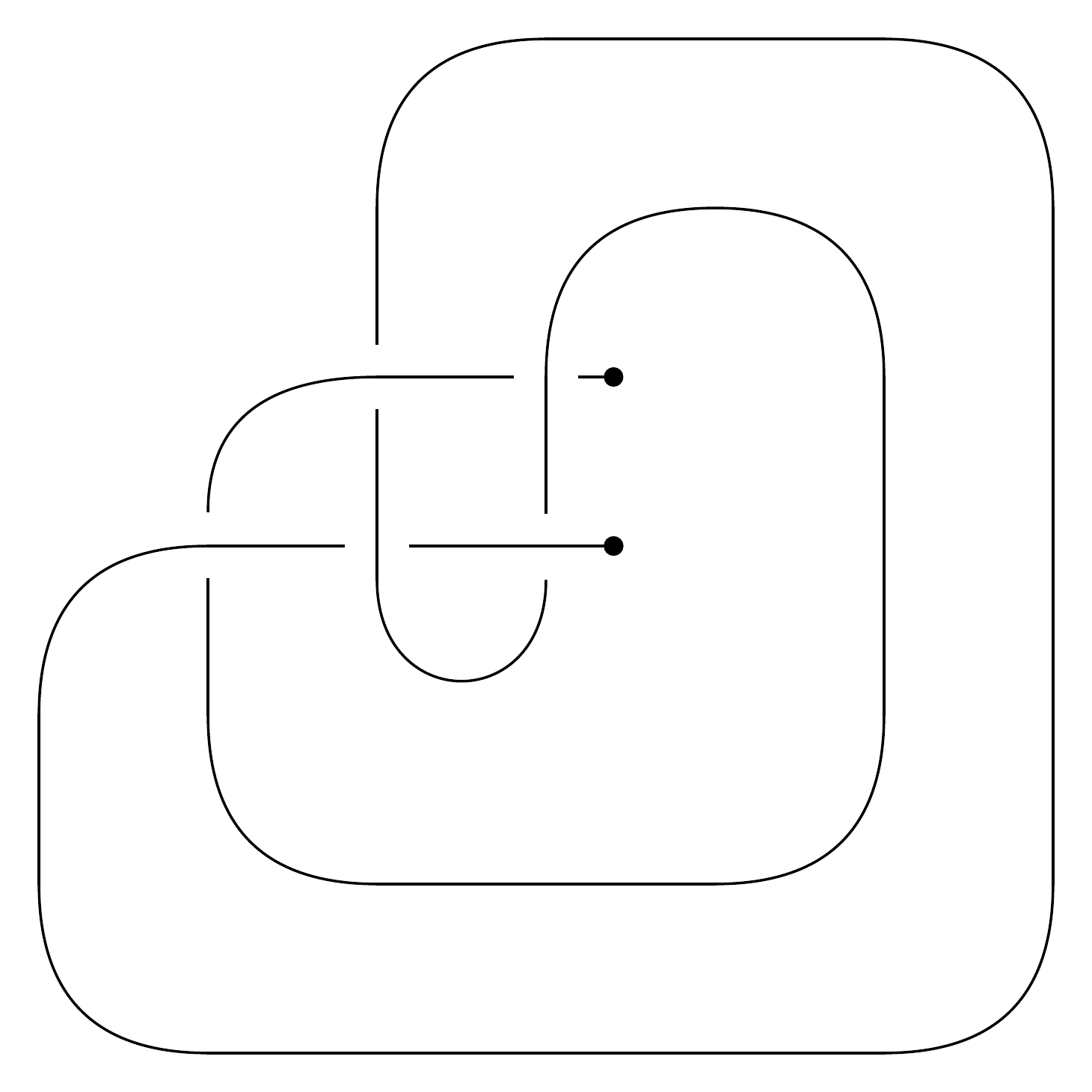}\\
\textcolor{black}{$5_{562}$}
\vspace{1cm}
\end{minipage}
\begin{minipage}[t]{.25\linewidth}
\centering
\includegraphics[width=0.9\textwidth,height=3.5cm,keepaspectratio]{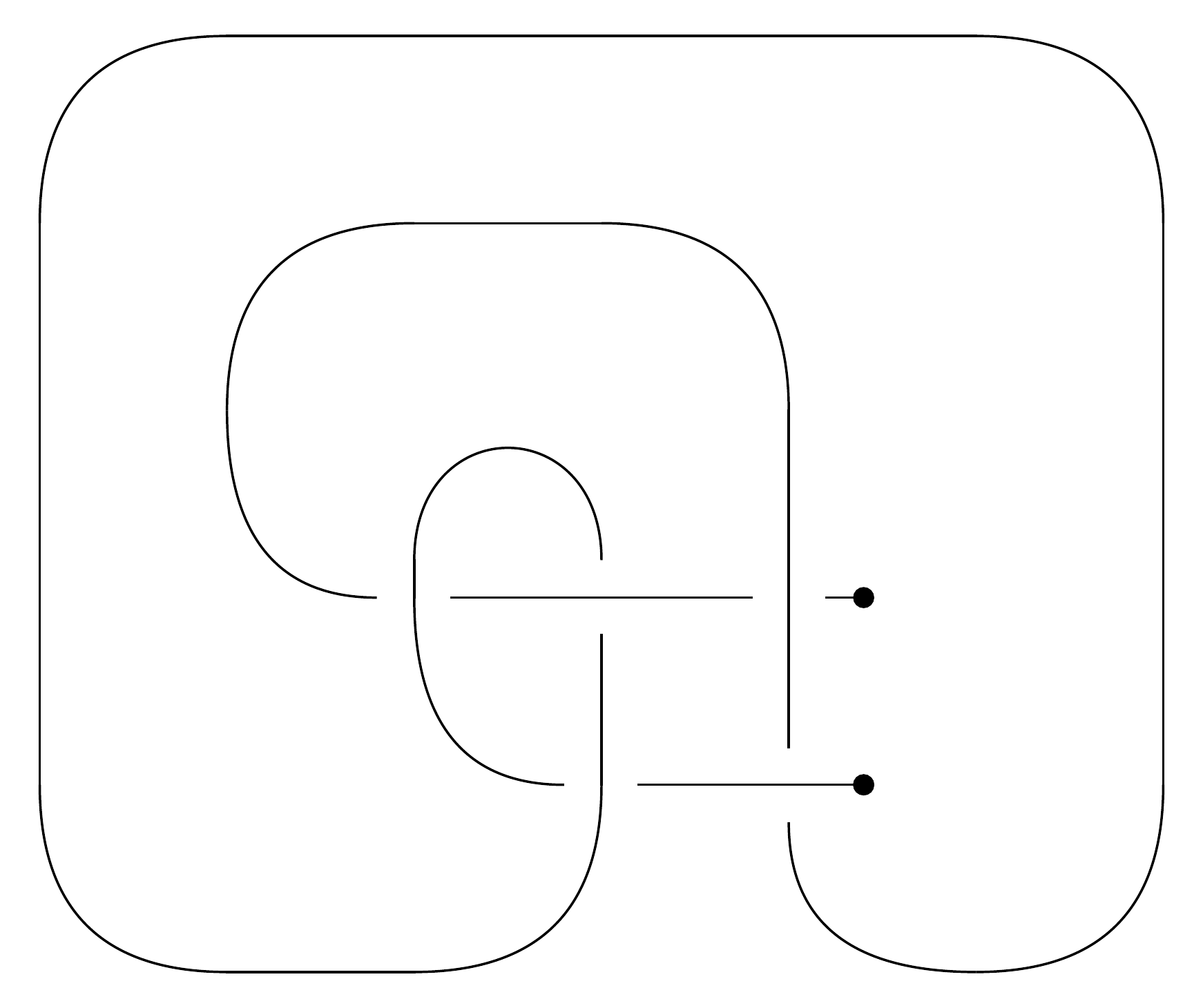}\\
\textcolor{black}{$5_{563}$}
\vspace{1cm}
\end{minipage}
\begin{minipage}[t]{.25\linewidth}
\centering
\includegraphics[width=0.9\textwidth,height=3.5cm,keepaspectratio]{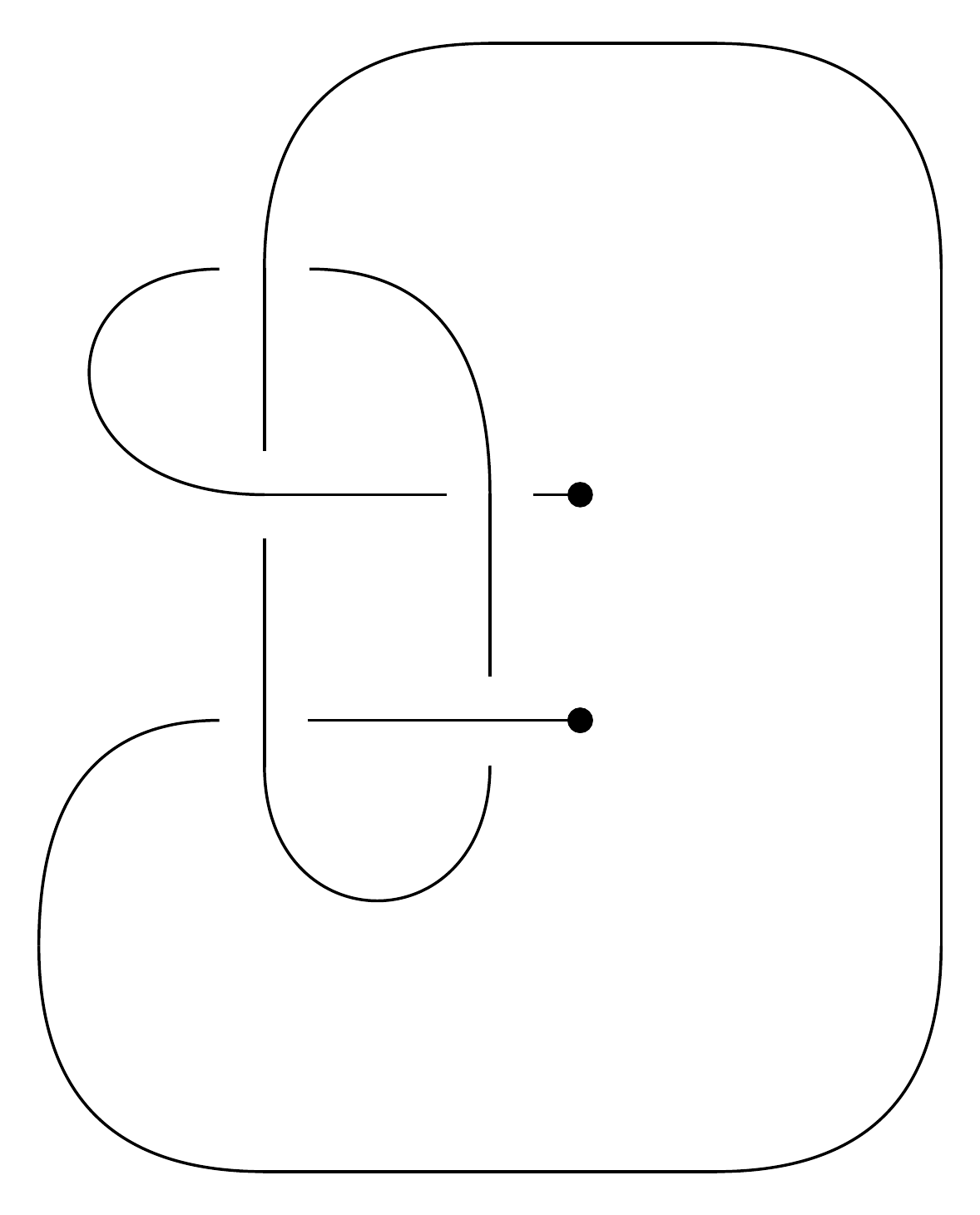}\\
\textcolor{black}{$5_{564}$}
\vspace{1cm}
\end{minipage}
\begin{minipage}[t]{.25\linewidth}
\centering
\includegraphics[width=0.9\textwidth,height=3.5cm,keepaspectratio]{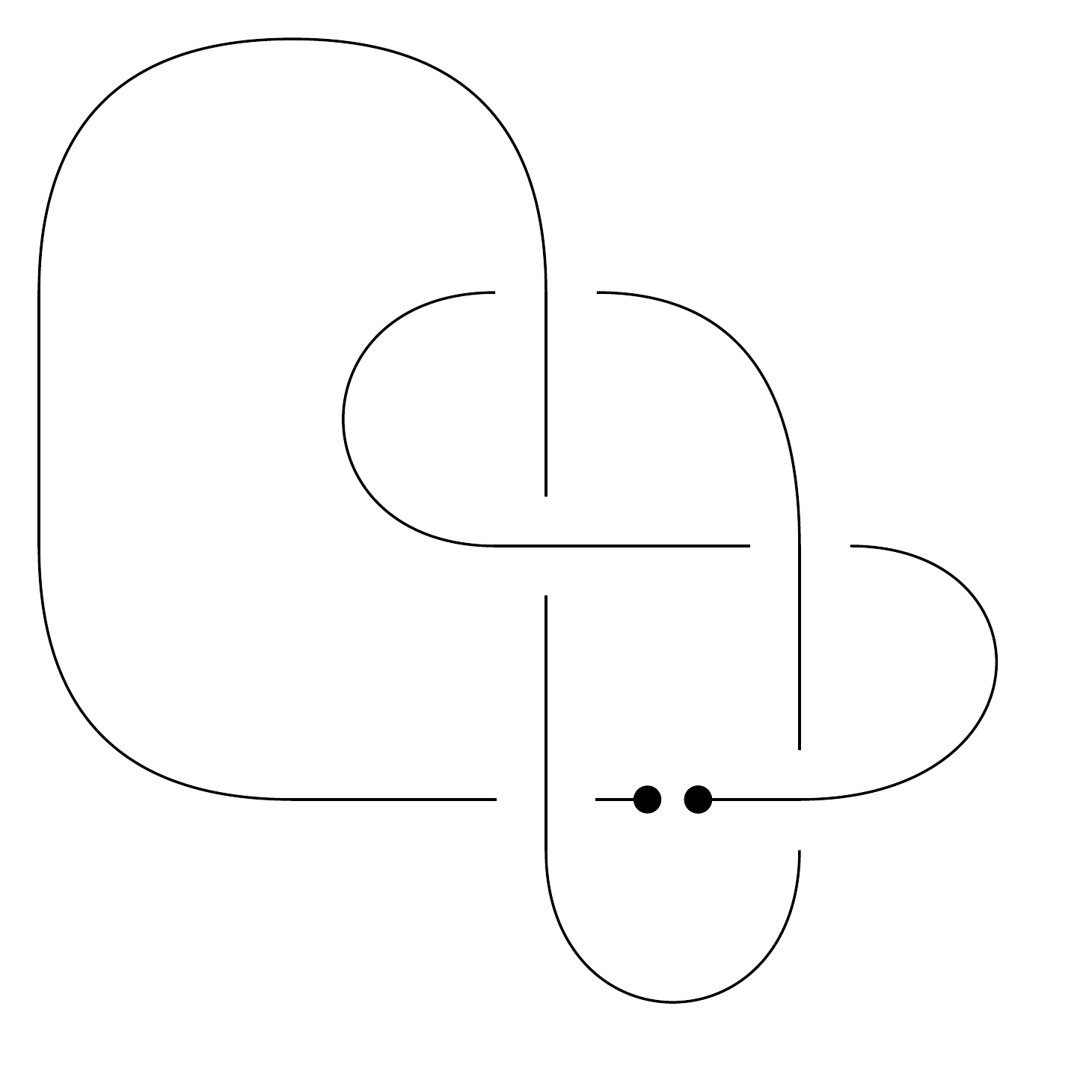}\\
\textcolor{black}{$5_{565}$}
\vspace{1cm}
\end{minipage}
\begin{minipage}[t]{.25\linewidth}
\centering
\includegraphics[width=0.9\textwidth,height=3.5cm,keepaspectratio]{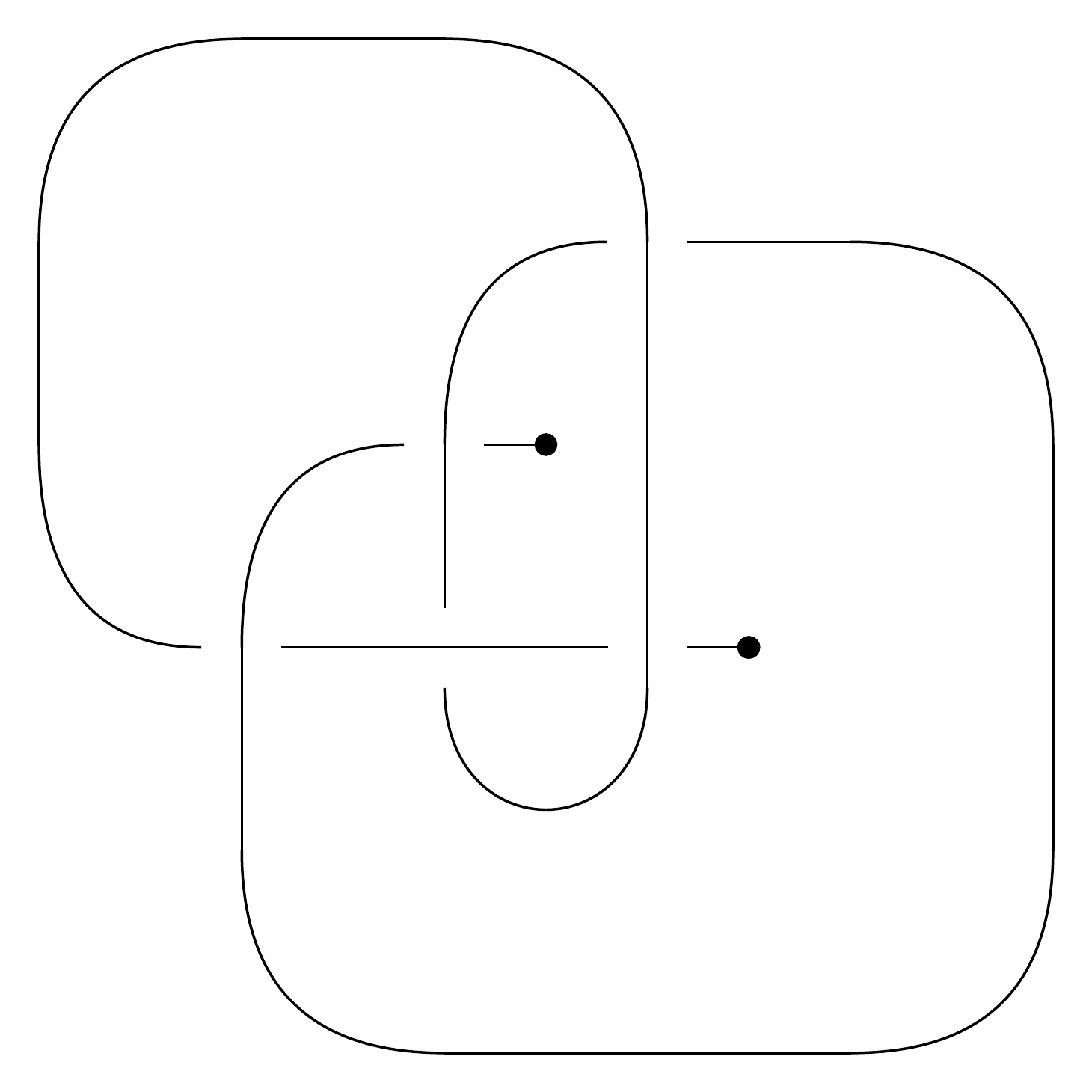}\\
\textcolor{black}{$5_{566}$}
\vspace{1cm}
\end{minipage}
\begin{minipage}[t]{.25\linewidth}
\centering
\includegraphics[width=0.9\textwidth,height=3.5cm,keepaspectratio]{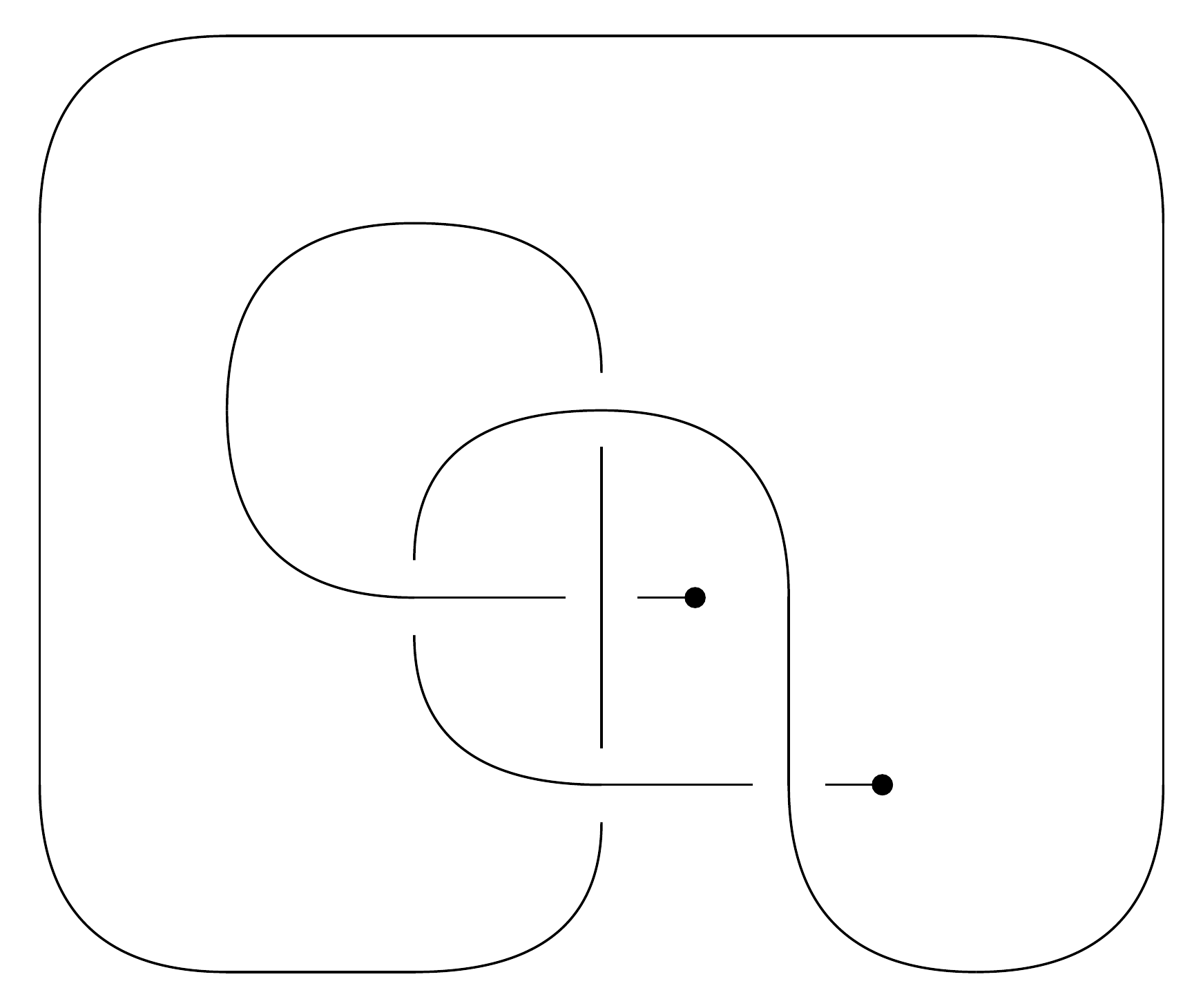}\\
\textcolor{black}{$5_{567}$}
\vspace{1cm}
\end{minipage}
\begin{minipage}[t]{.25\linewidth}
\centering
\includegraphics[width=0.9\textwidth,height=3.5cm,keepaspectratio]{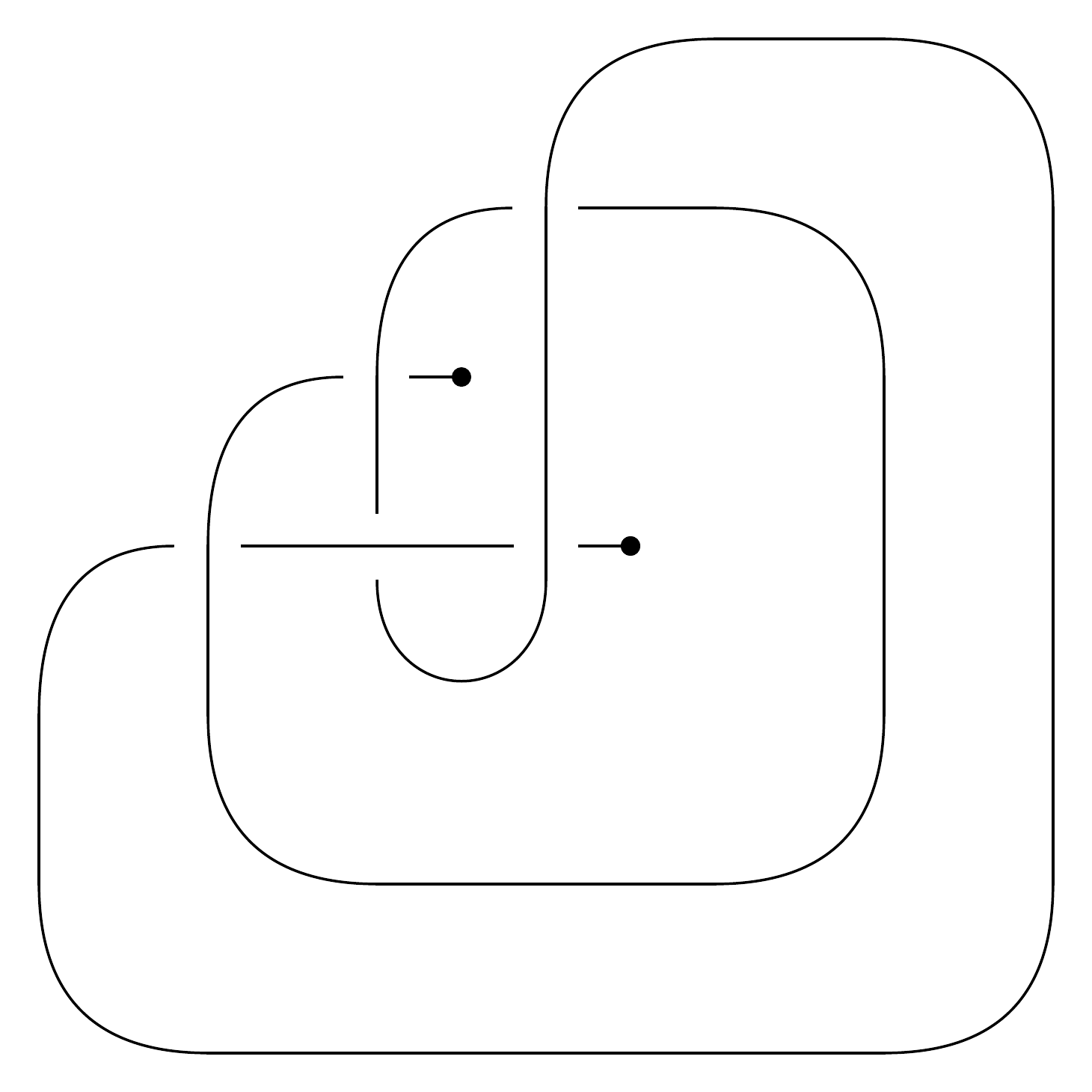}\\
\textcolor{black}{$5_{568}$}
\vspace{1cm}
\end{minipage}
\begin{minipage}[t]{.25\linewidth}
\centering
\includegraphics[width=0.9\textwidth,height=3.5cm,keepaspectratio]{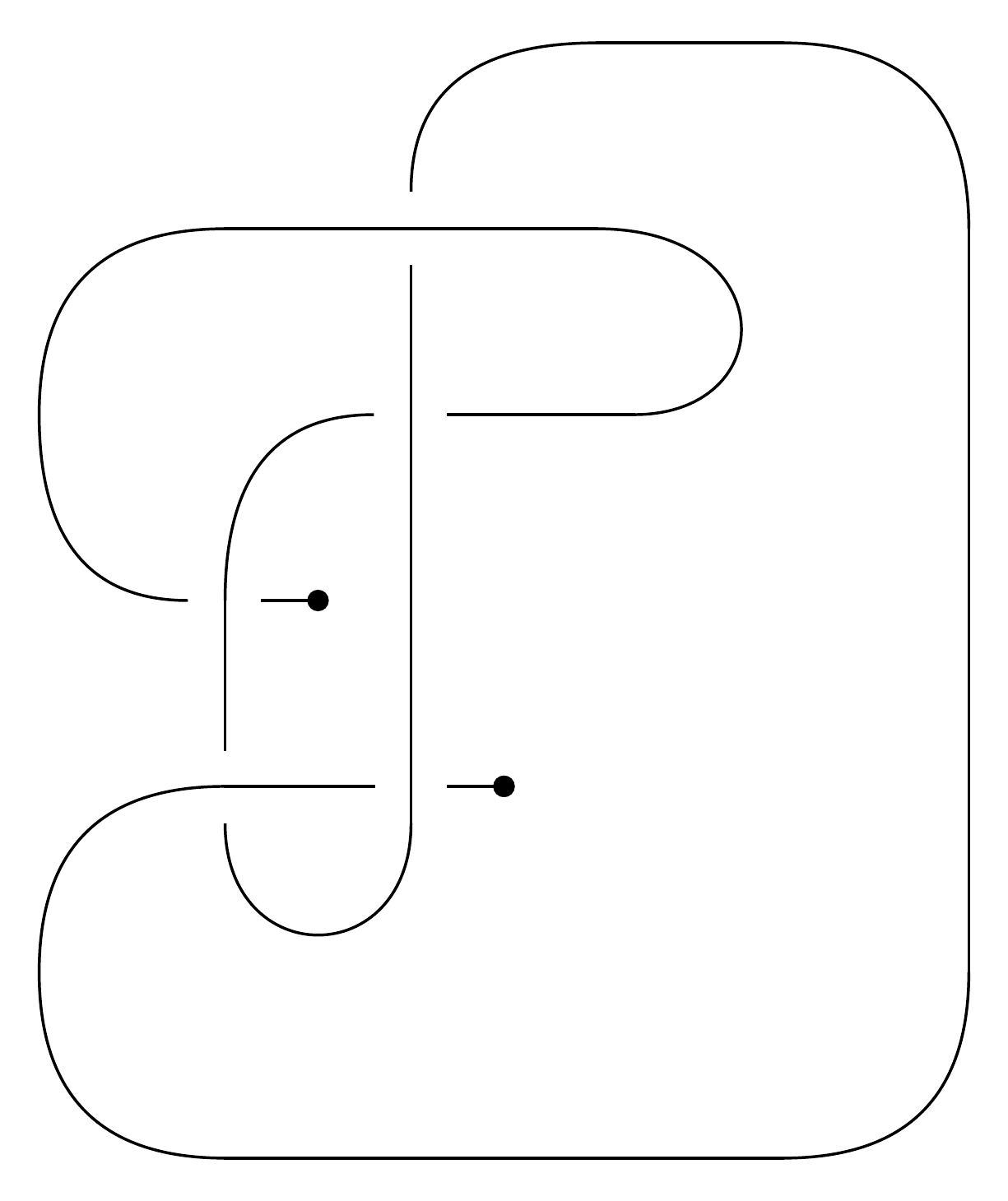}\\
\textcolor{black}{$5_{569}$}
\vspace{1cm}
\end{minipage}
\begin{minipage}[t]{.25\linewidth}
\centering
\includegraphics[width=0.9\textwidth,height=3.5cm,keepaspectratio]{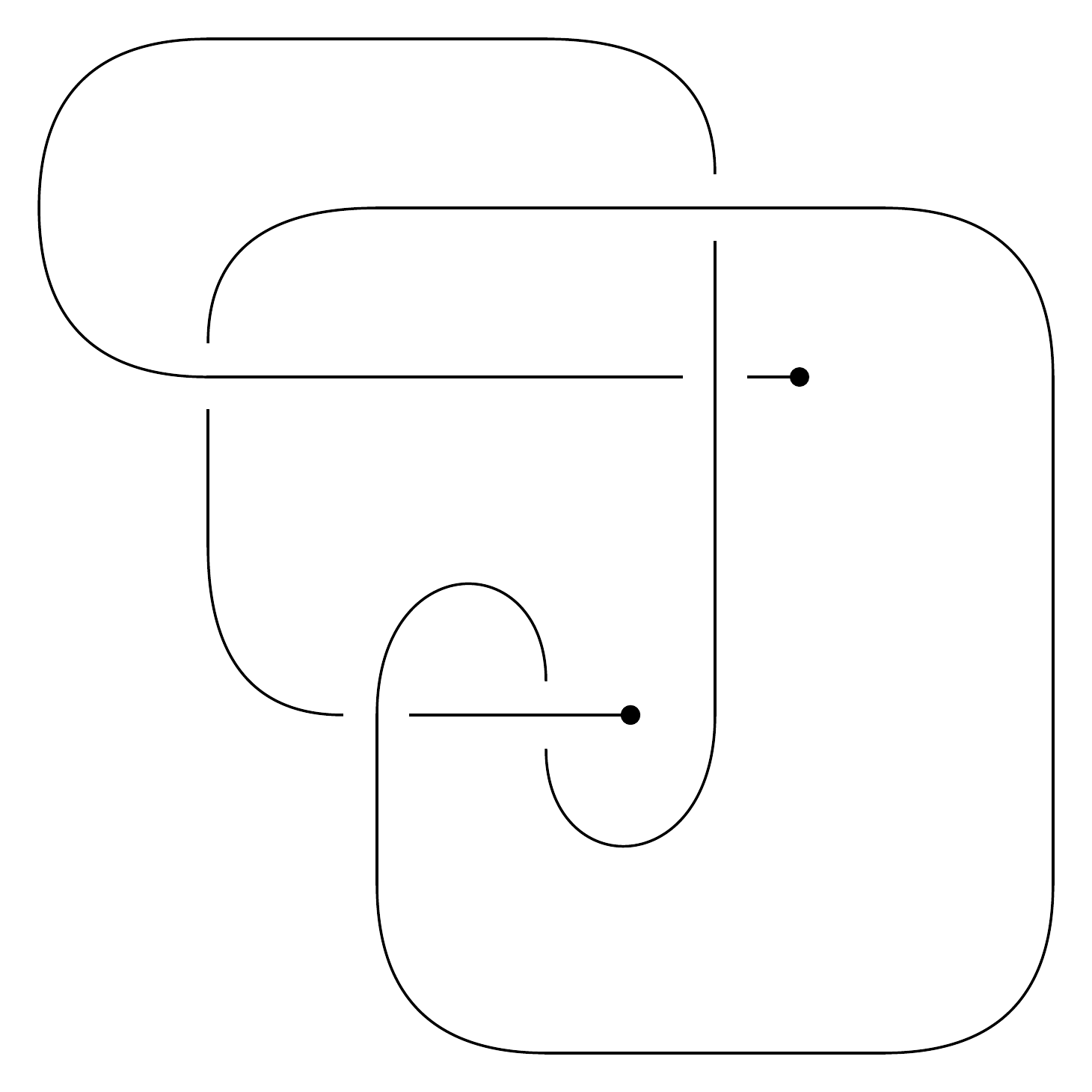}\\
\textcolor{black}{$5_{570}$}
\vspace{1cm}
\end{minipage}
\begin{minipage}[t]{.25\linewidth}
\centering
\includegraphics[width=0.9\textwidth,height=3.5cm,keepaspectratio]{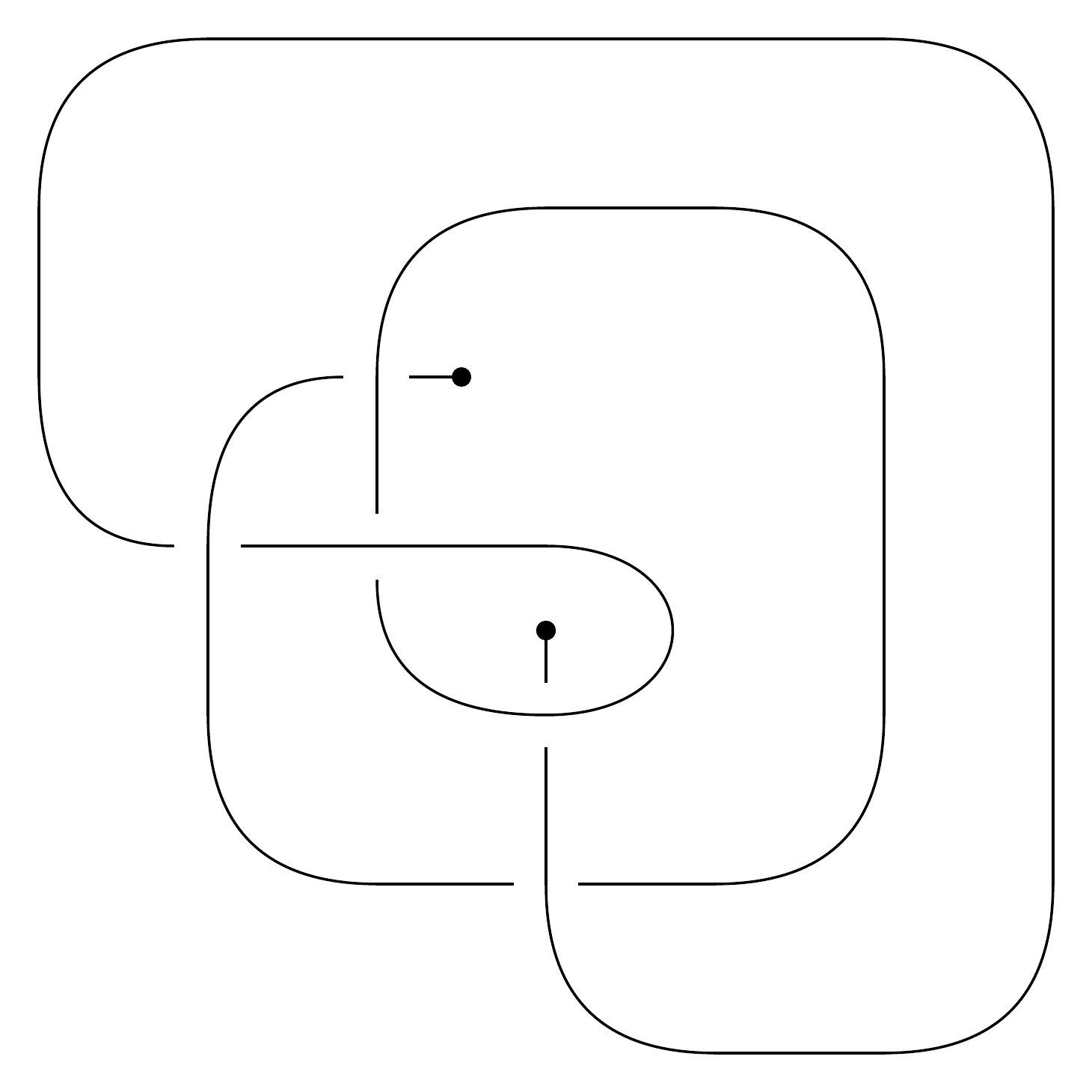}\\
\textcolor{black}{$5_{571}$}
\vspace{1cm}
\end{minipage}
\begin{minipage}[t]{.25\linewidth}
\centering
\includegraphics[width=0.9\textwidth,height=3.5cm,keepaspectratio]{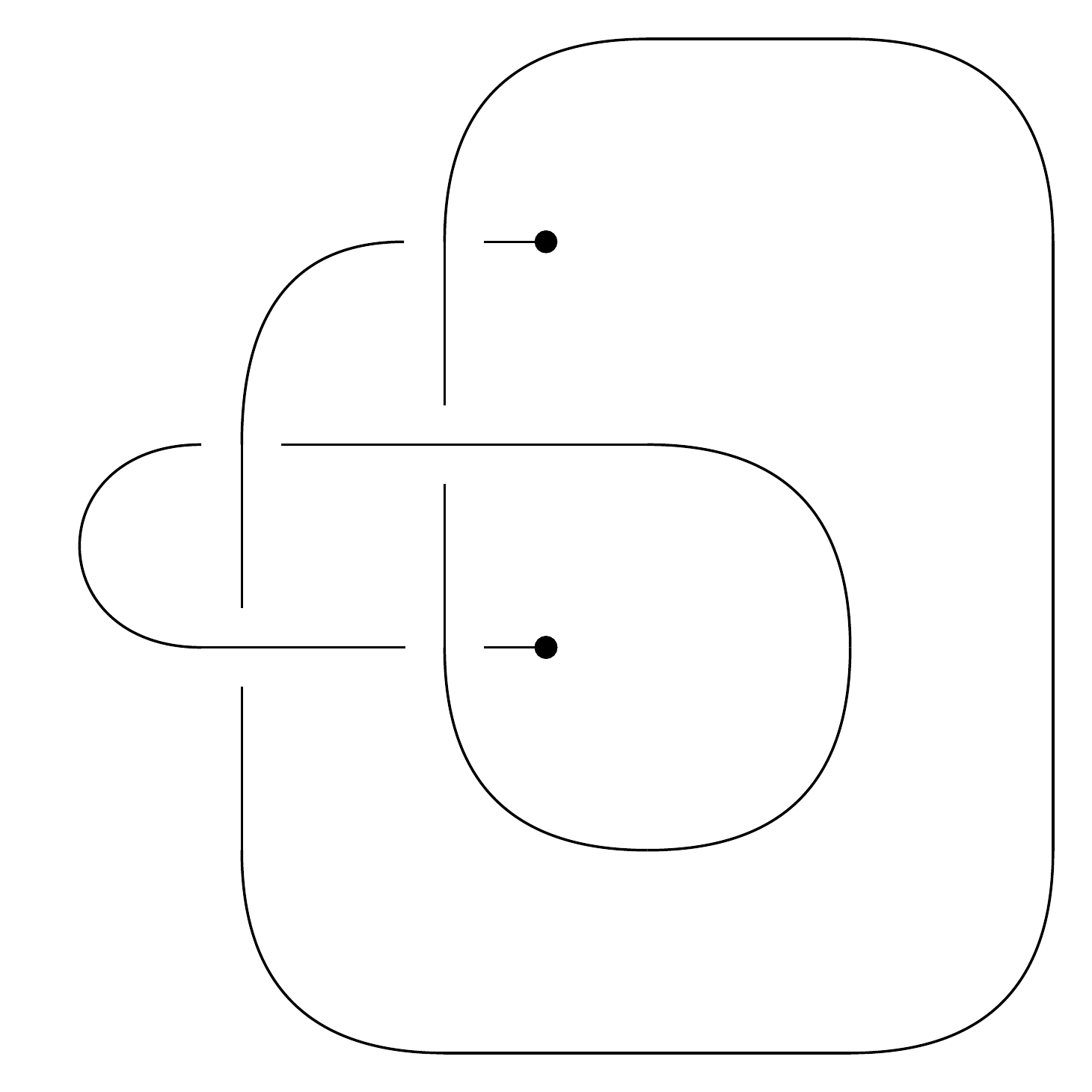}\\
\textcolor{black}{$5_{572}$}
\vspace{1cm}
\end{minipage}
\begin{minipage}[t]{.25\linewidth}
\centering
\includegraphics[width=0.9\textwidth,height=3.5cm,keepaspectratio]{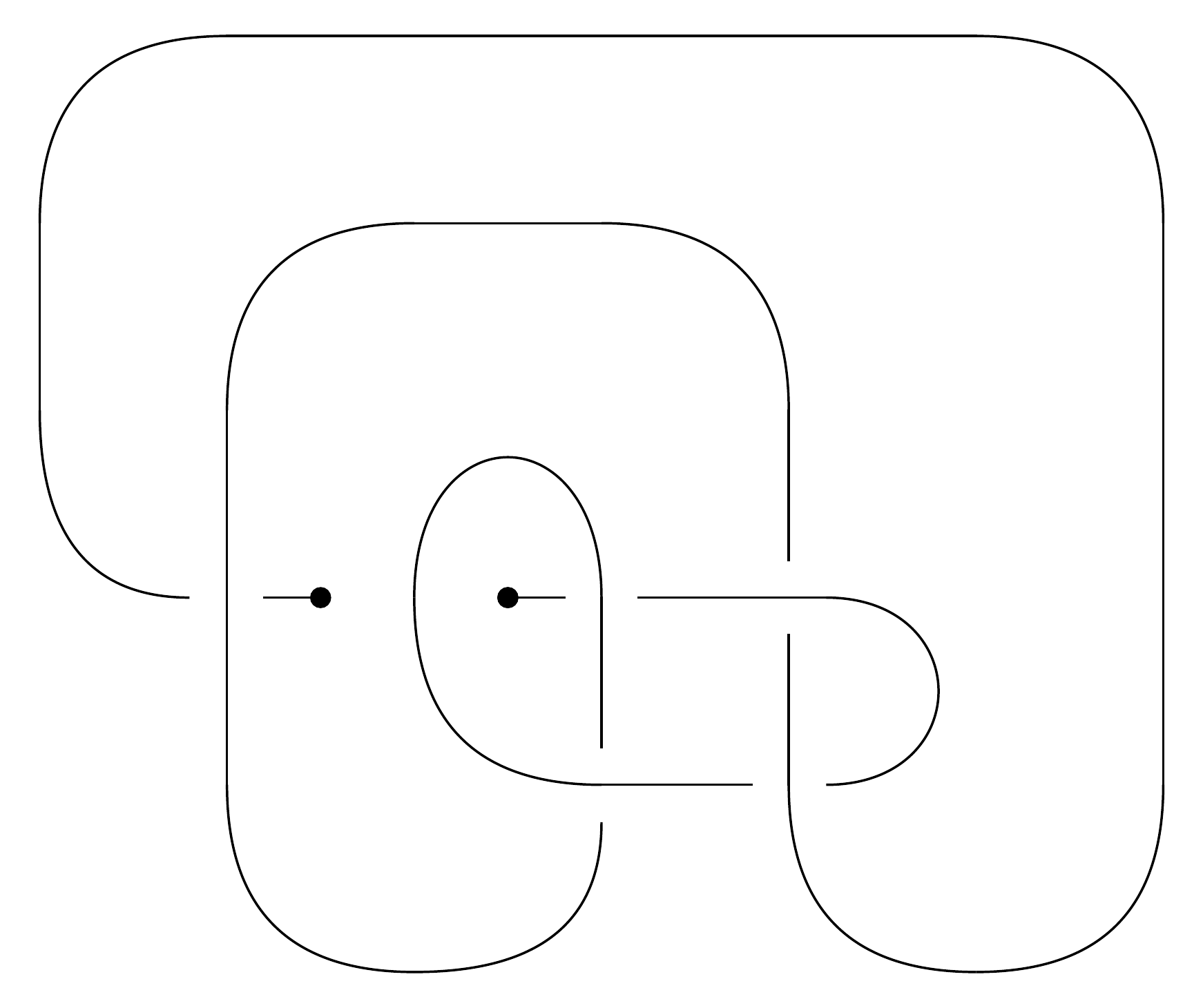}\\
\textcolor{black}{$5_{573}$}
\vspace{1cm}
\end{minipage}
\begin{minipage}[t]{.25\linewidth}
\centering
\includegraphics[width=0.9\textwidth,height=3.5cm,keepaspectratio]{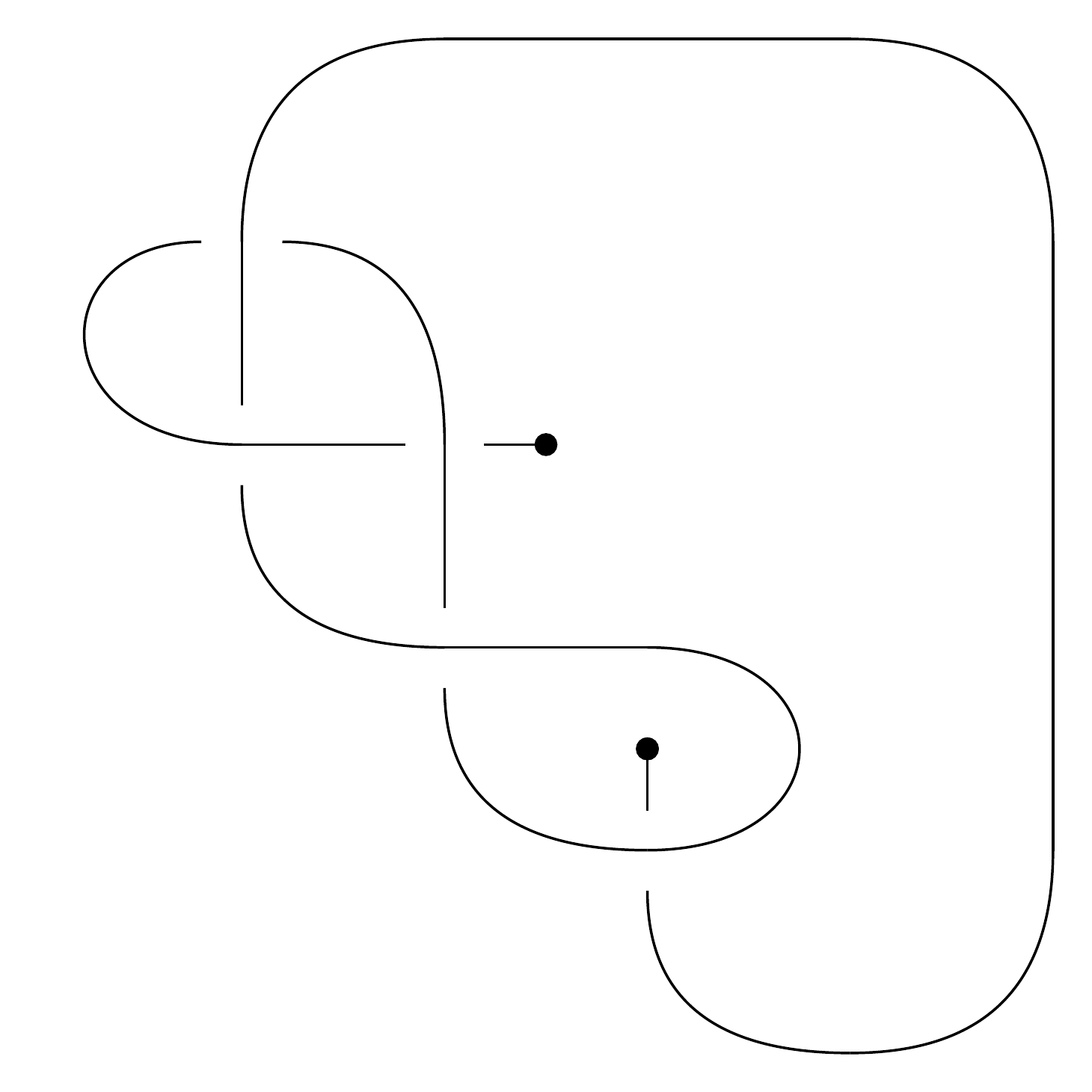}\\
\textcolor{black}{$5_{574}$}
\vspace{1cm}
\end{minipage}
\begin{minipage}[t]{.25\linewidth}
\centering
\includegraphics[width=0.9\textwidth,height=3.5cm,keepaspectratio]{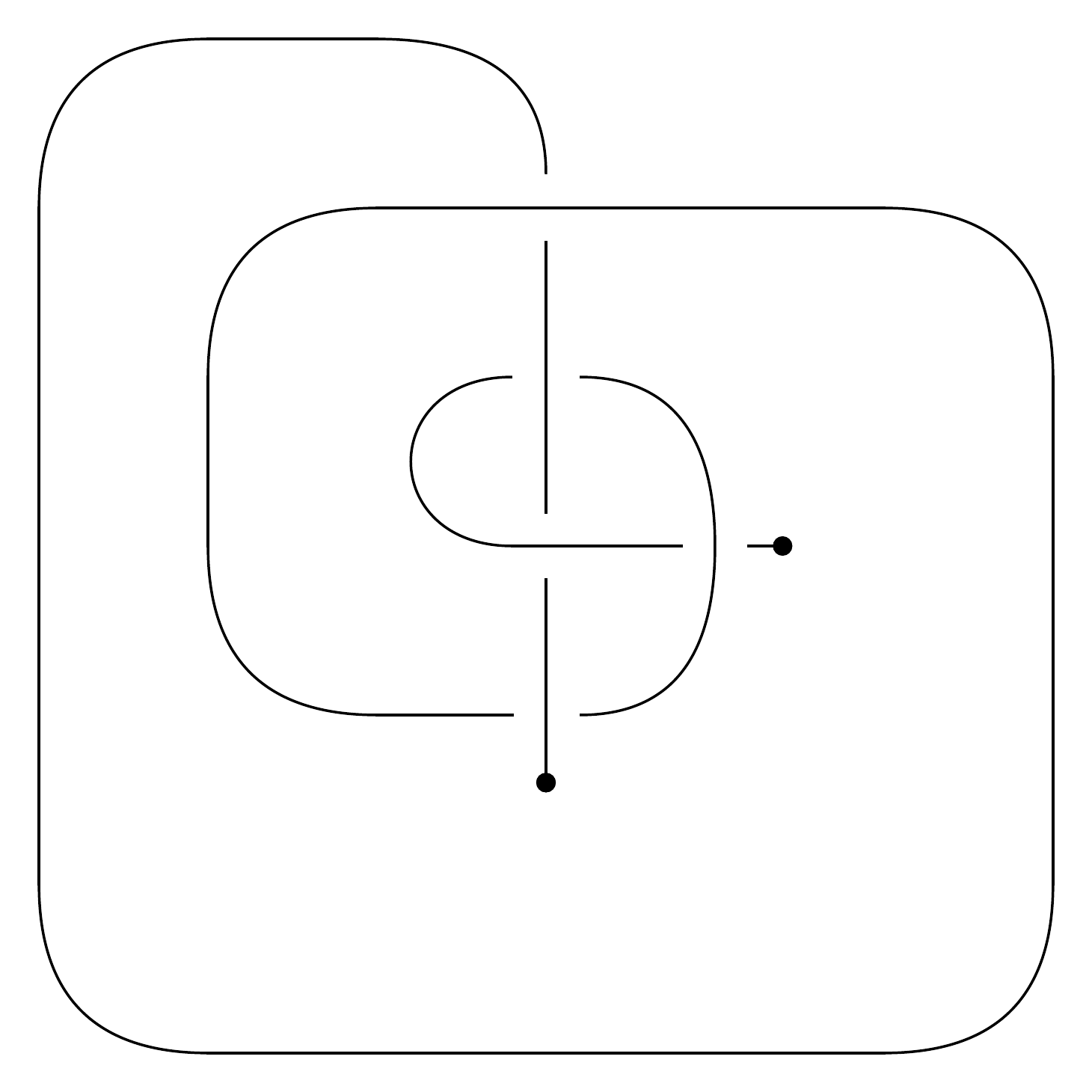}\\
\textcolor{black}{$5_{575}$}
\vspace{1cm}
\end{minipage}
\begin{minipage}[t]{.25\linewidth}
\centering
\includegraphics[width=0.9\textwidth,height=3.5cm,keepaspectratio]{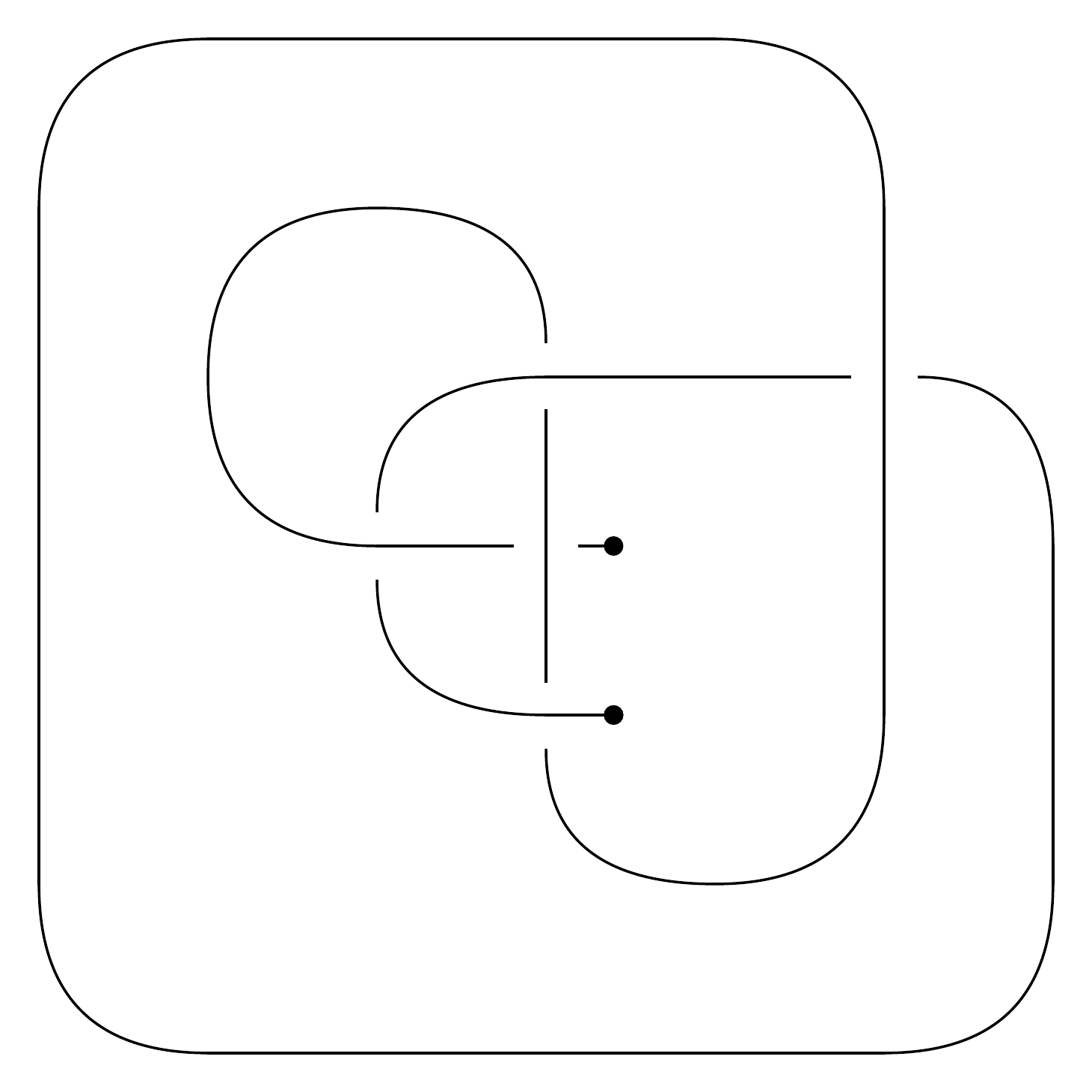}\\
\textcolor{black}{$5_{576}$}
\vspace{1cm}
\end{minipage}
\begin{minipage}[t]{.25\linewidth}
\centering
\includegraphics[width=0.9\textwidth,height=3.5cm,keepaspectratio]{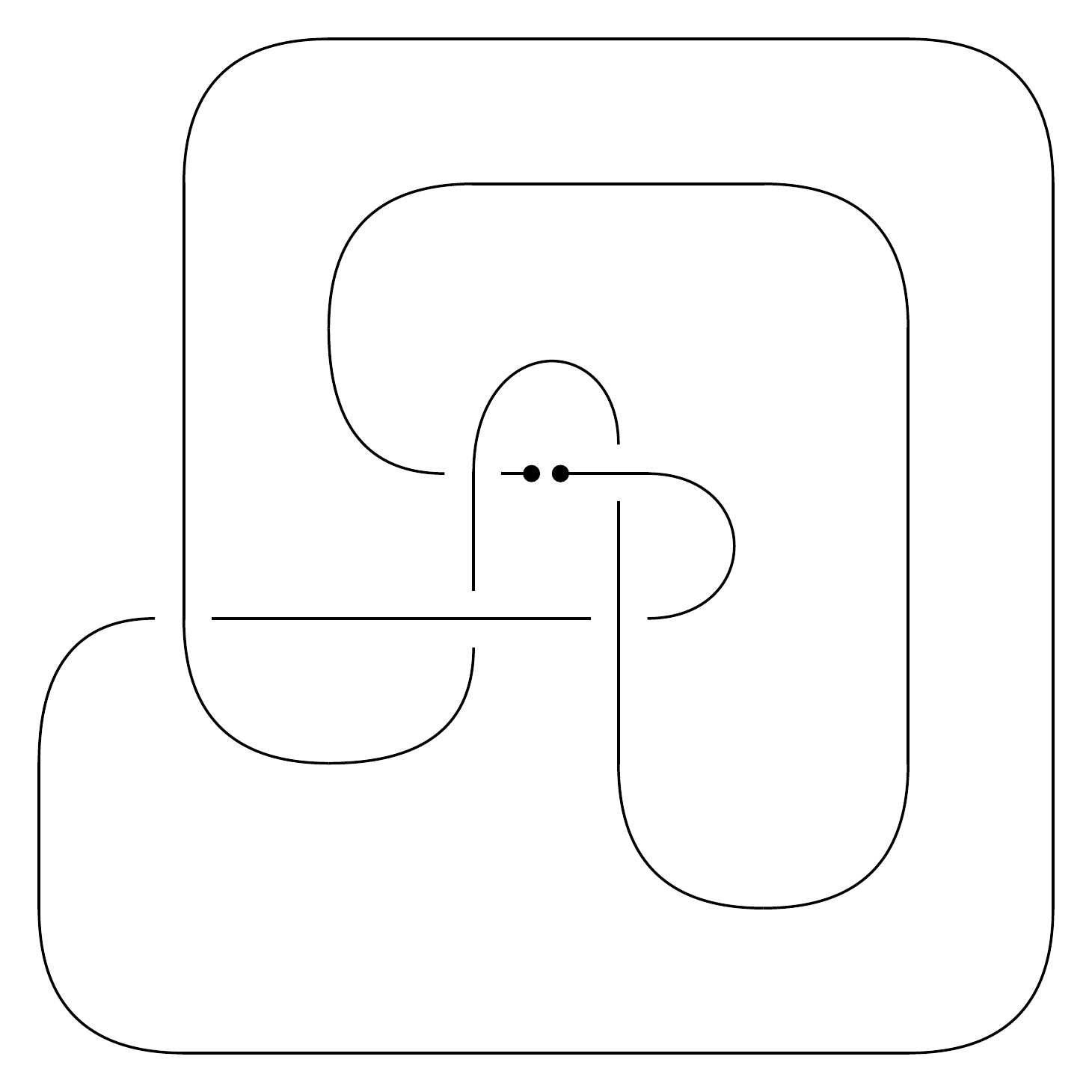}\\
\textcolor{black}{$5_{577}$}
\vspace{1cm}
\end{minipage}
\begin{minipage}[t]{.25\linewidth}
\centering
\includegraphics[width=0.9\textwidth,height=3.5cm,keepaspectratio]{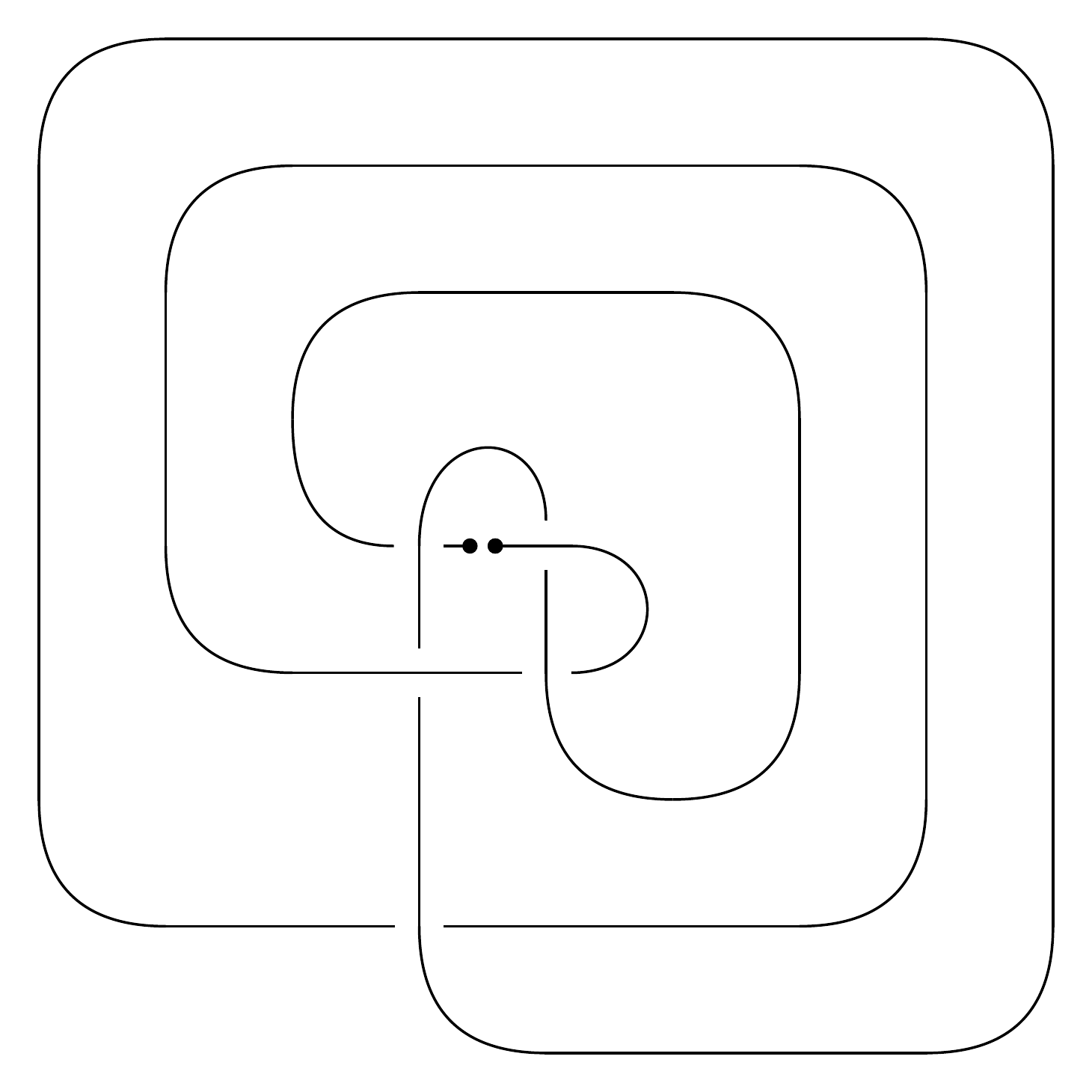}\\
\textcolor{black}{$5_{578}$}
\vspace{1cm}
\end{minipage}
\begin{minipage}[t]{.25\linewidth}
\centering
\includegraphics[width=0.9\textwidth,height=3.5cm,keepaspectratio]{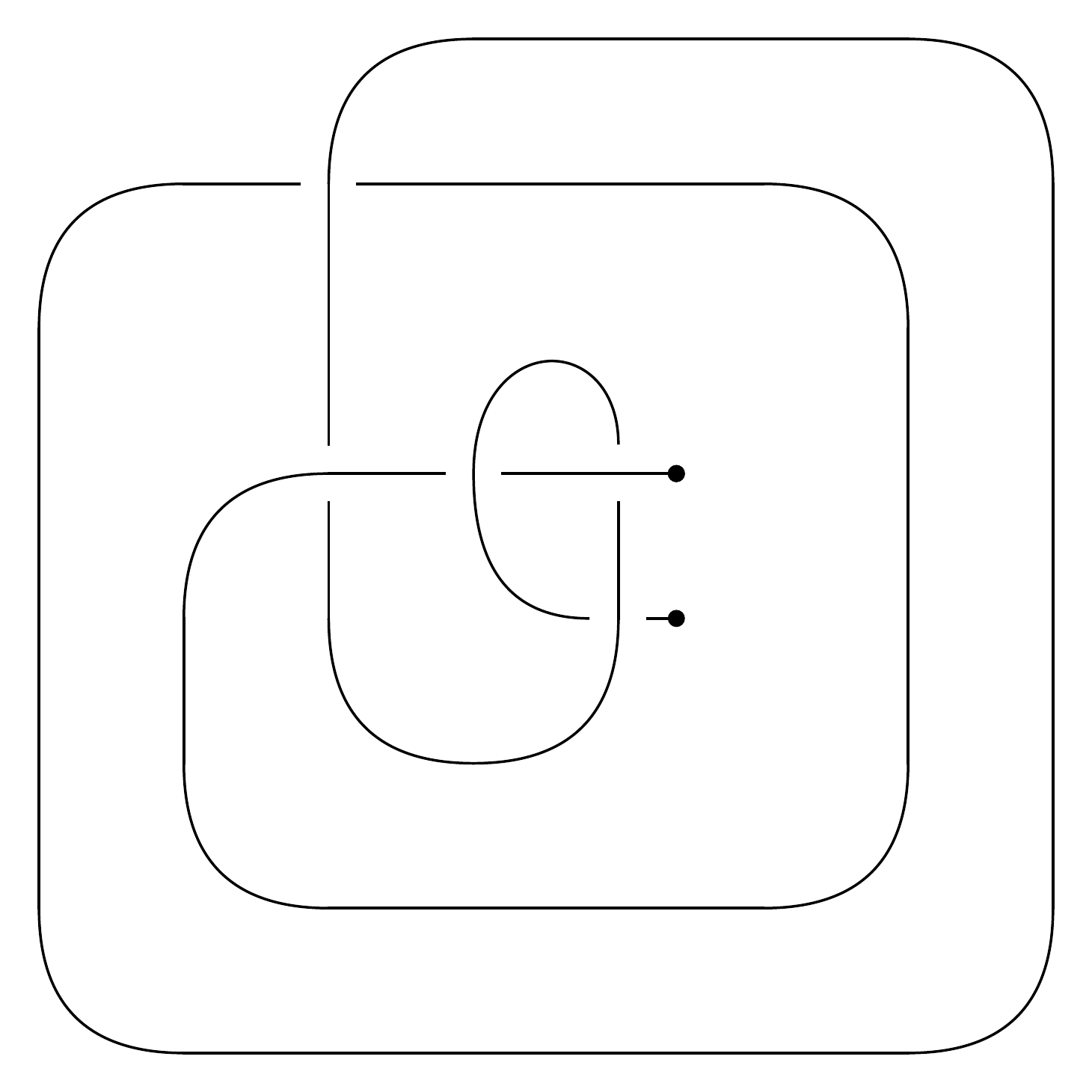}\\
\textcolor{black}{$5_{579}$}
\vspace{1cm}
\end{minipage}
\begin{minipage}[t]{.25\linewidth}
\centering
\includegraphics[width=0.9\textwidth,height=3.5cm,keepaspectratio]{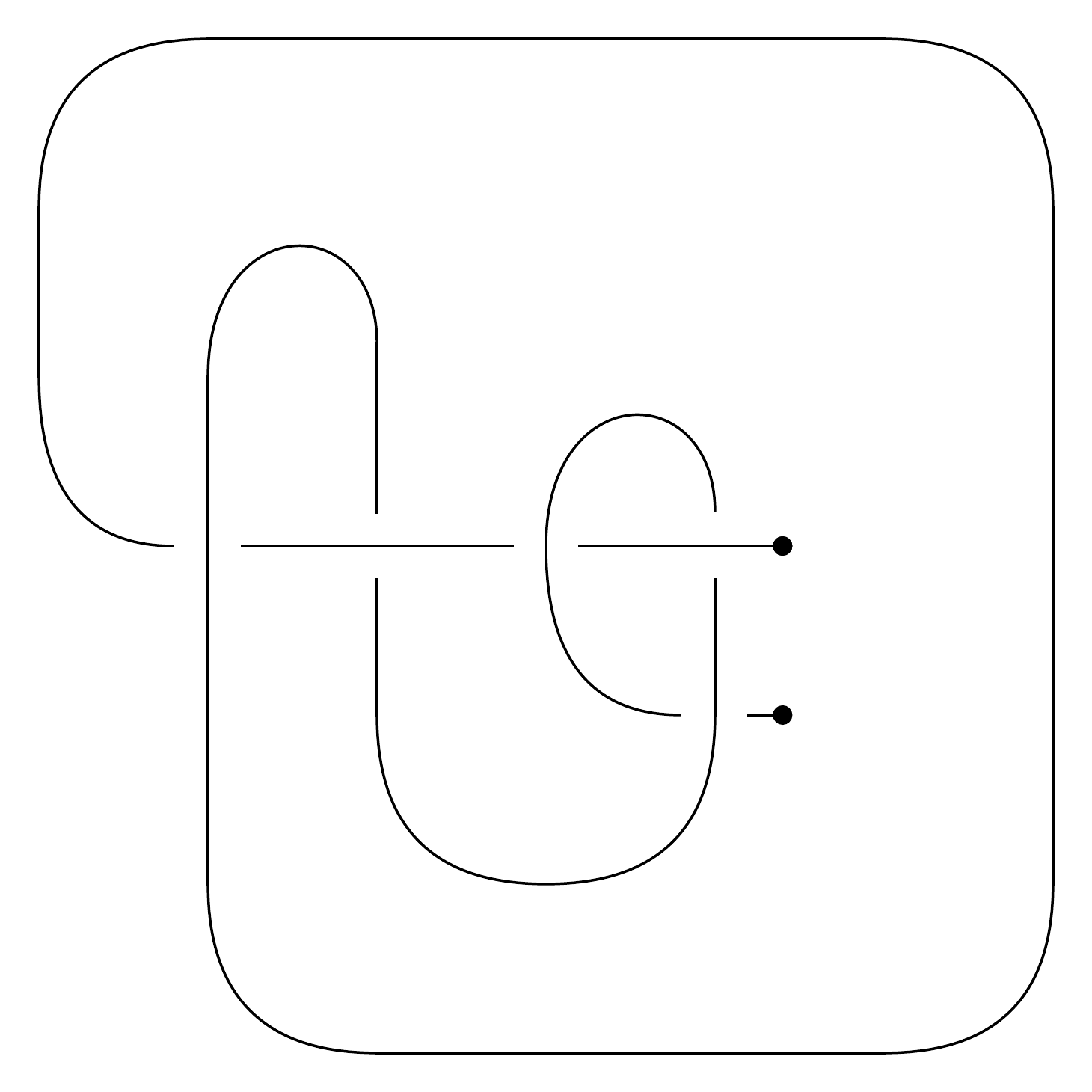}\\
\textcolor{black}{$5_{580}$}
\vspace{1cm}
\end{minipage}
\begin{minipage}[t]{.25\linewidth}
\centering
\includegraphics[width=0.9\textwidth,height=3.5cm,keepaspectratio]{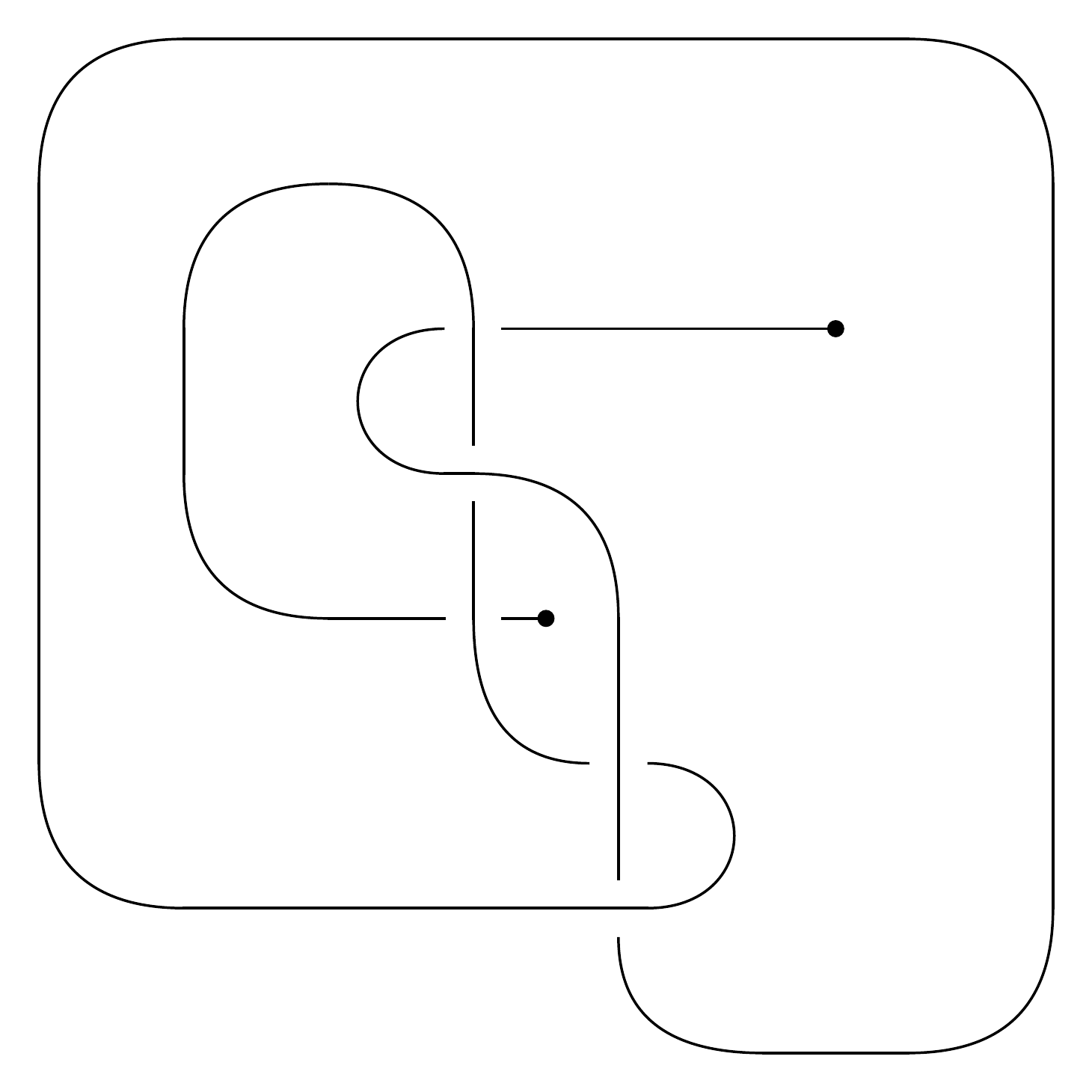}\\
\textcolor{black}{$5_{581}$}
\vspace{1cm}
\end{minipage}
\begin{minipage}[t]{.25\linewidth}
\centering
\includegraphics[width=0.9\textwidth,height=3.5cm,keepaspectratio]{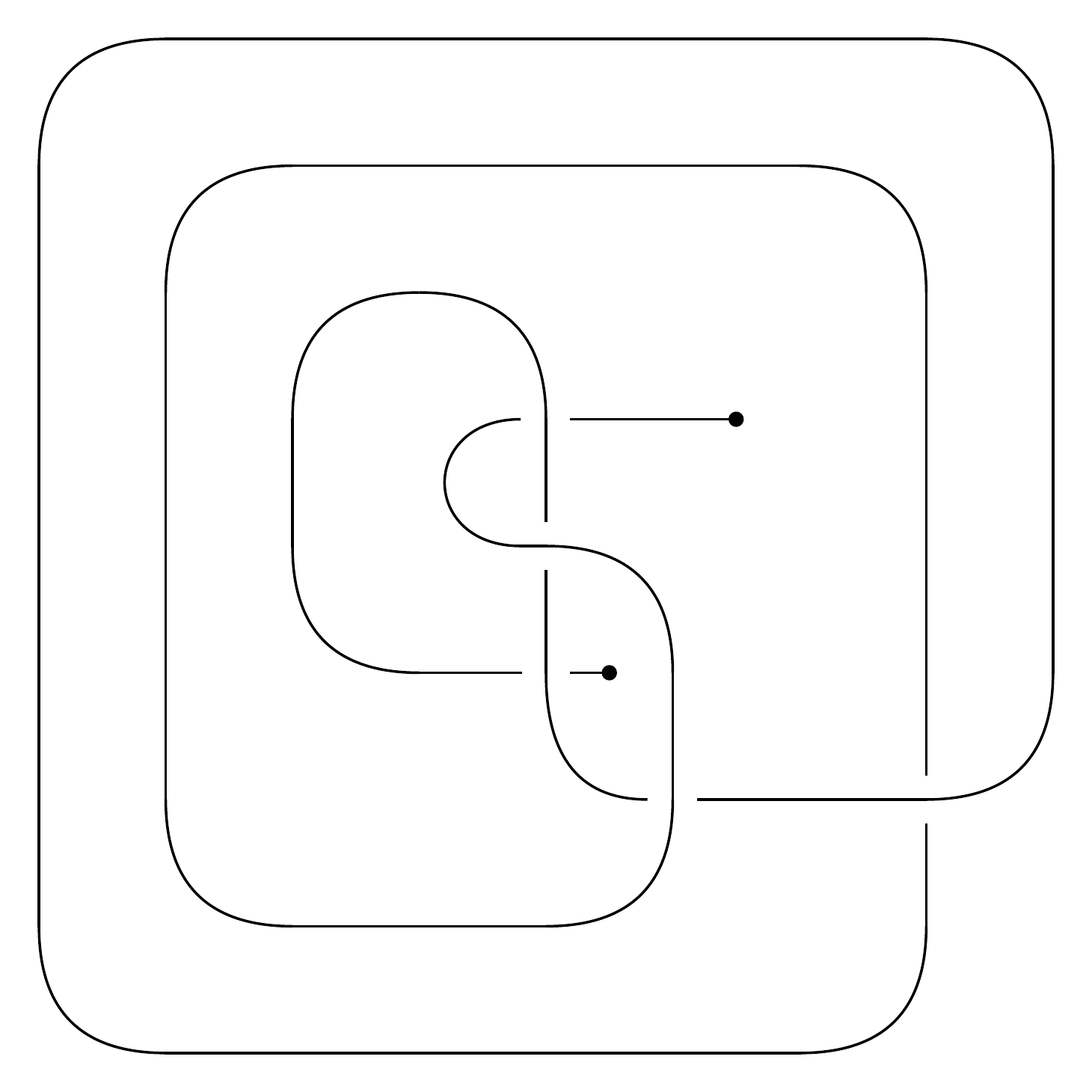}\\
\textcolor{black}{$5_{582}$}
\vspace{1cm}
\end{minipage}
\begin{minipage}[t]{.25\linewidth}
\centering
\includegraphics[width=0.9\textwidth,height=3.5cm,keepaspectratio]{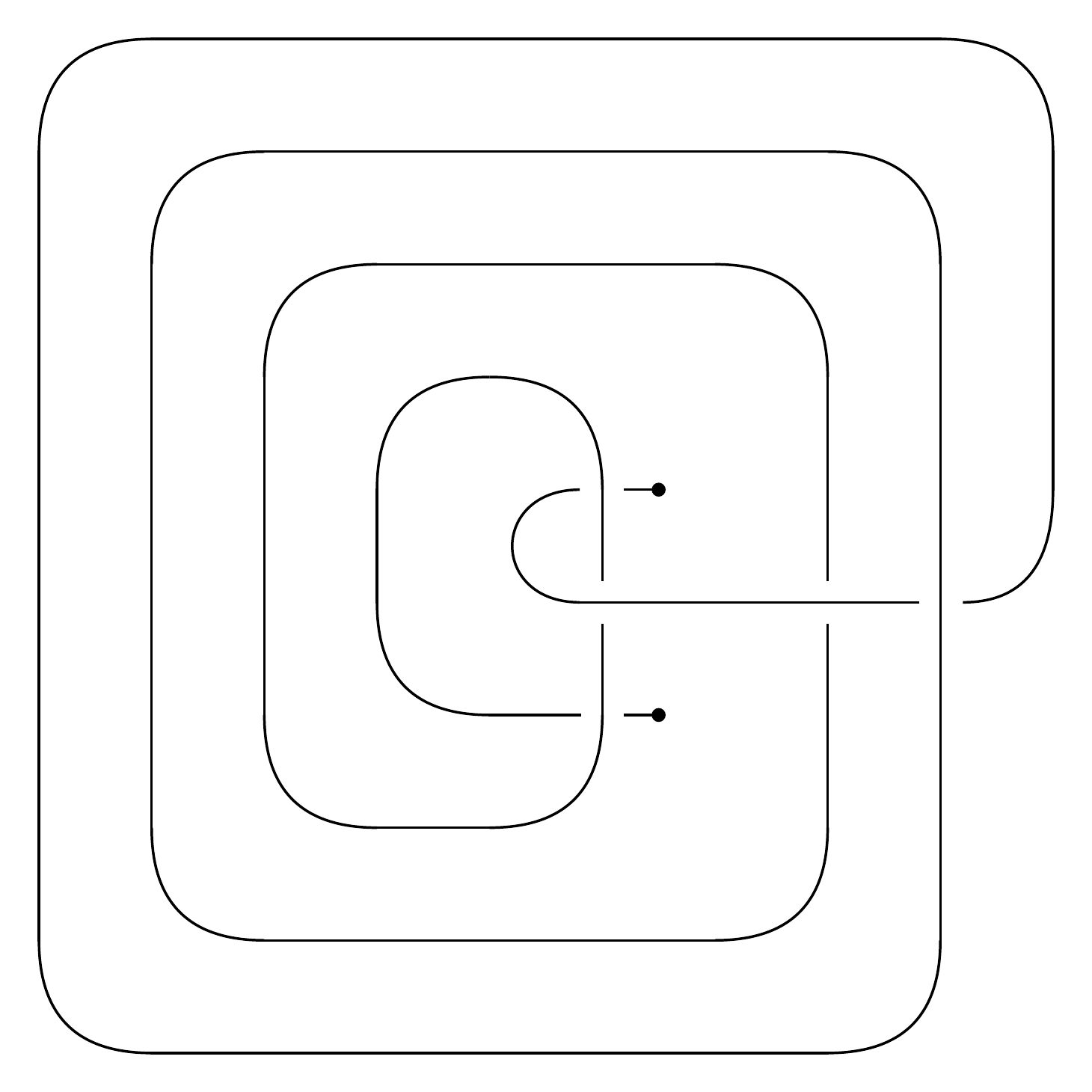}\\
\textcolor{black}{$5_{583}$}
\vspace{1cm}
\end{minipage}
\begin{minipage}[t]{.25\linewidth}
\centering
\includegraphics[width=0.9\textwidth,height=3.5cm,keepaspectratio]{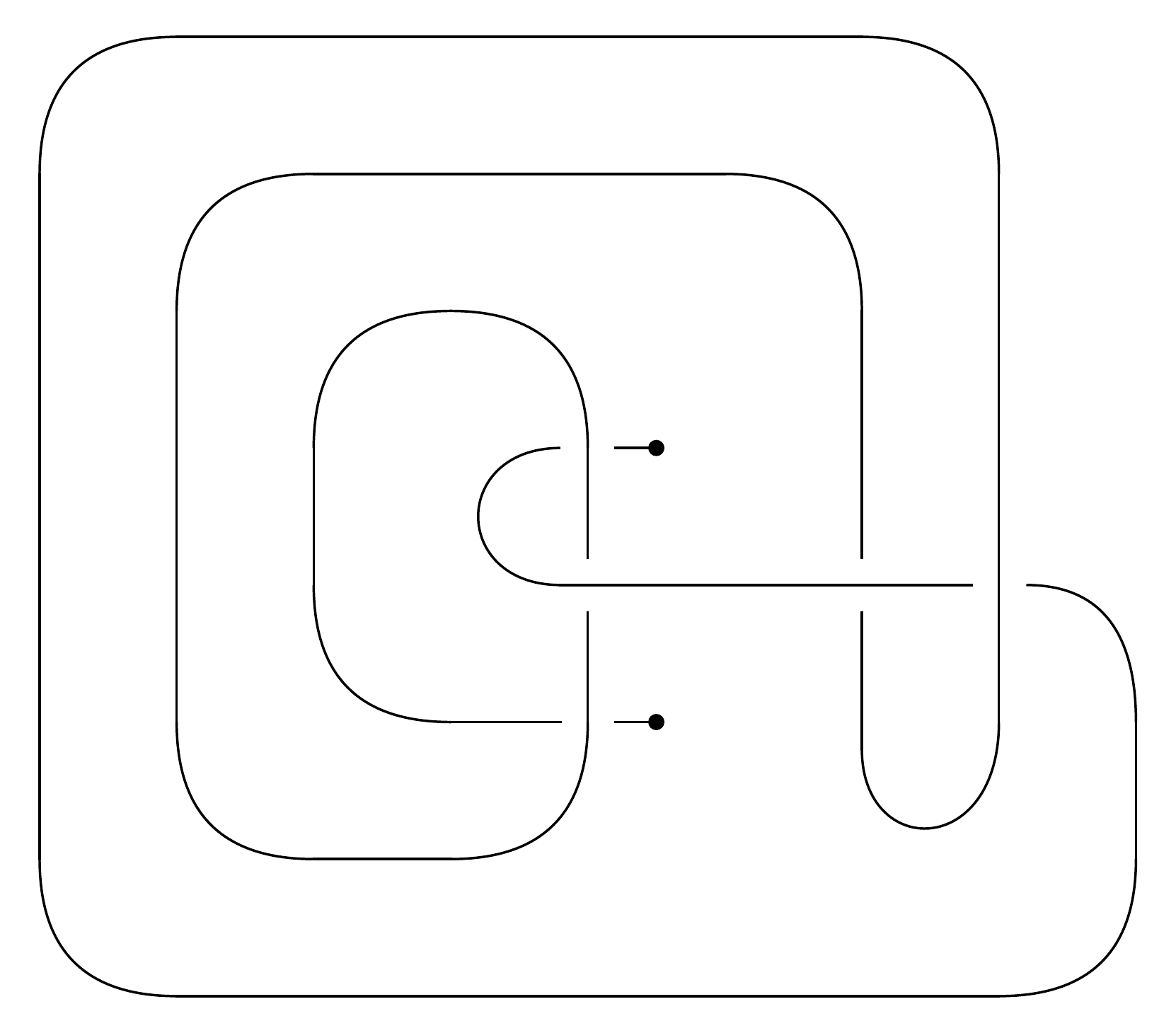}\\
\textcolor{black}{$5_{584}$}
\vspace{1cm}
\end{minipage}
\begin{minipage}[t]{.25\linewidth}
\centering
\includegraphics[width=0.9\textwidth,height=3.5cm,keepaspectratio]{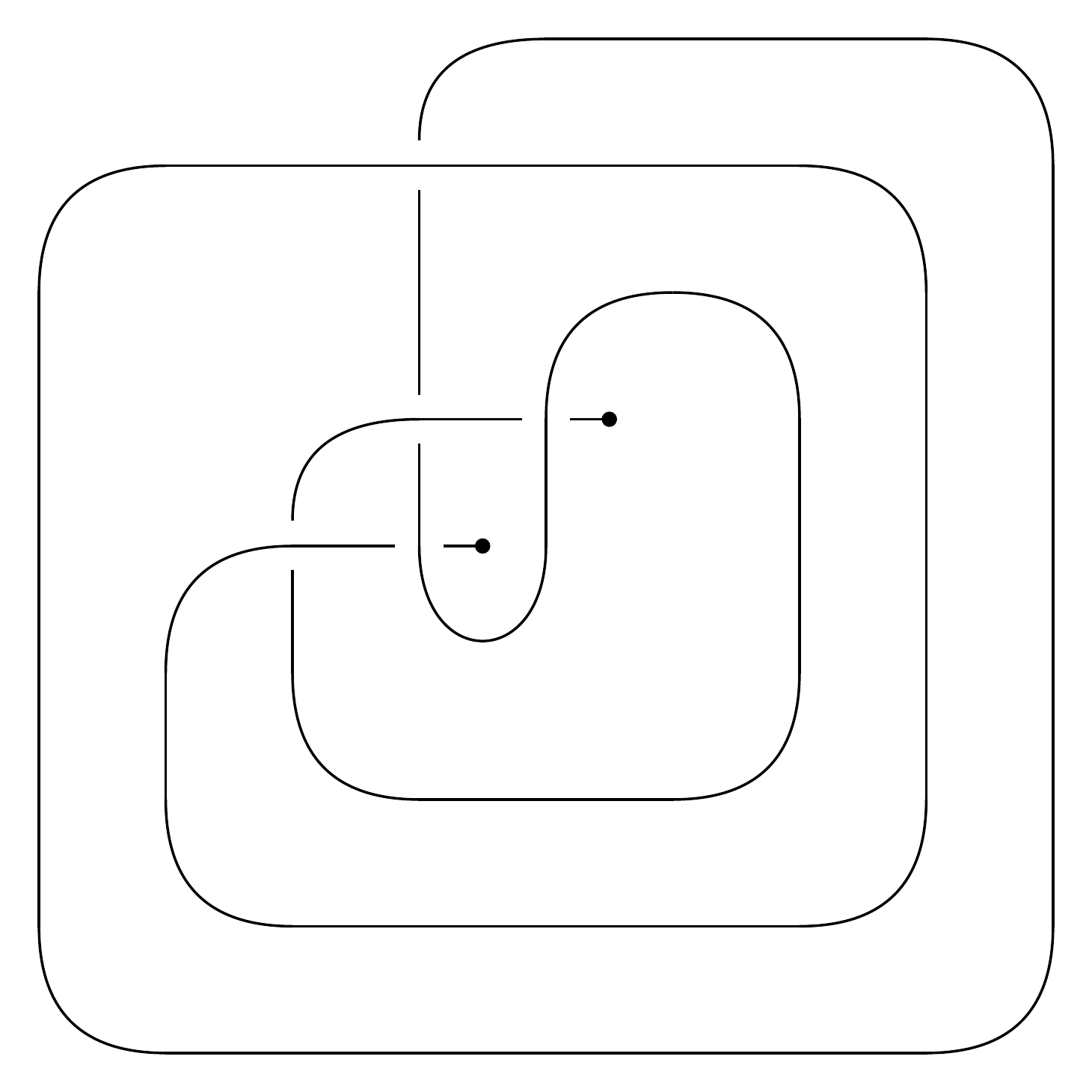}\\
\textcolor{black}{$5_{585}$}
\vspace{1cm}
\end{minipage}
\begin{minipage}[t]{.25\linewidth}
\centering
\includegraphics[width=0.9\textwidth,height=3.5cm,keepaspectratio]{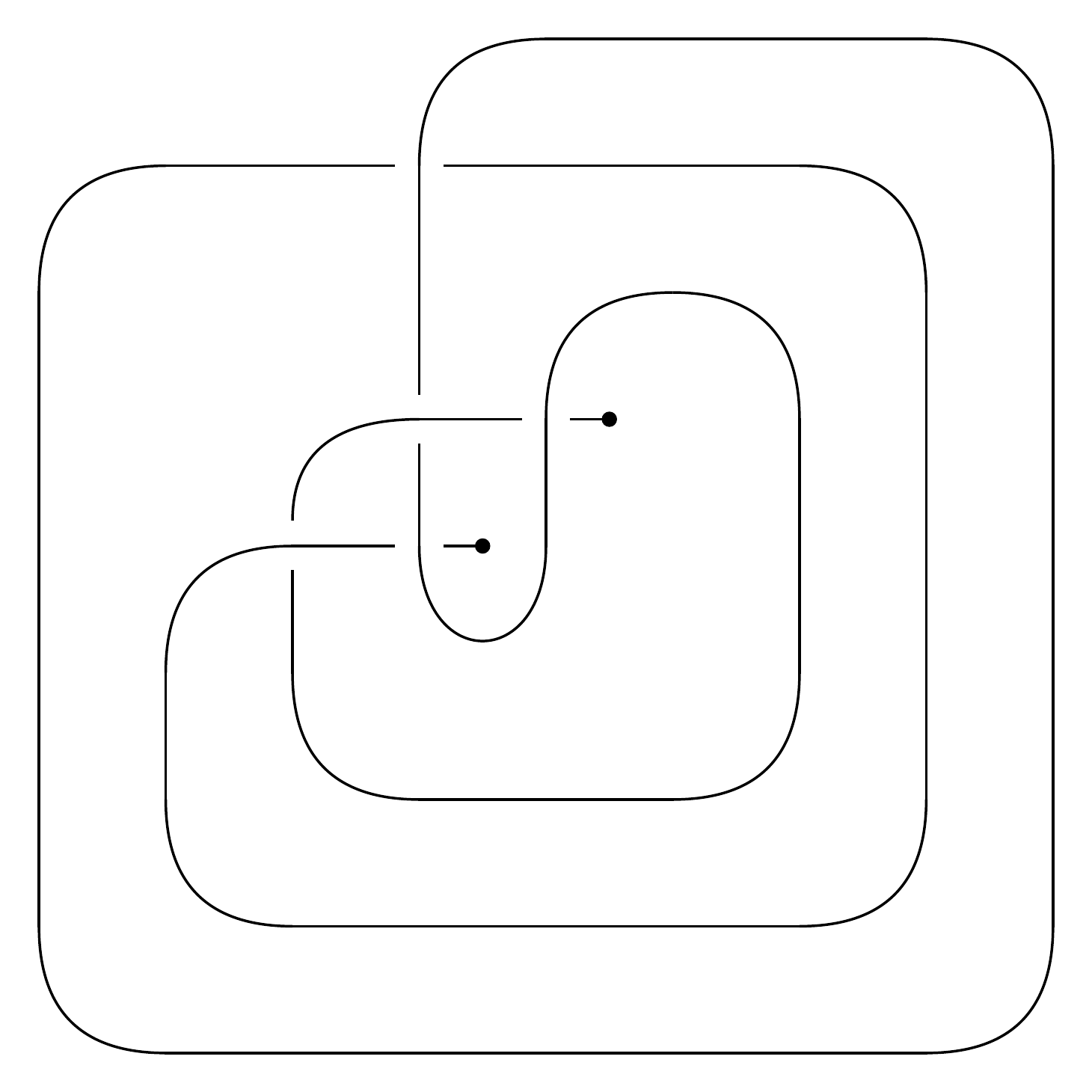}\\
\textcolor{black}{$5_{586}$}
\vspace{1cm}
\end{minipage}
\begin{minipage}[t]{.25\linewidth}
\centering
\includegraphics[width=0.9\textwidth,height=3.5cm,keepaspectratio]{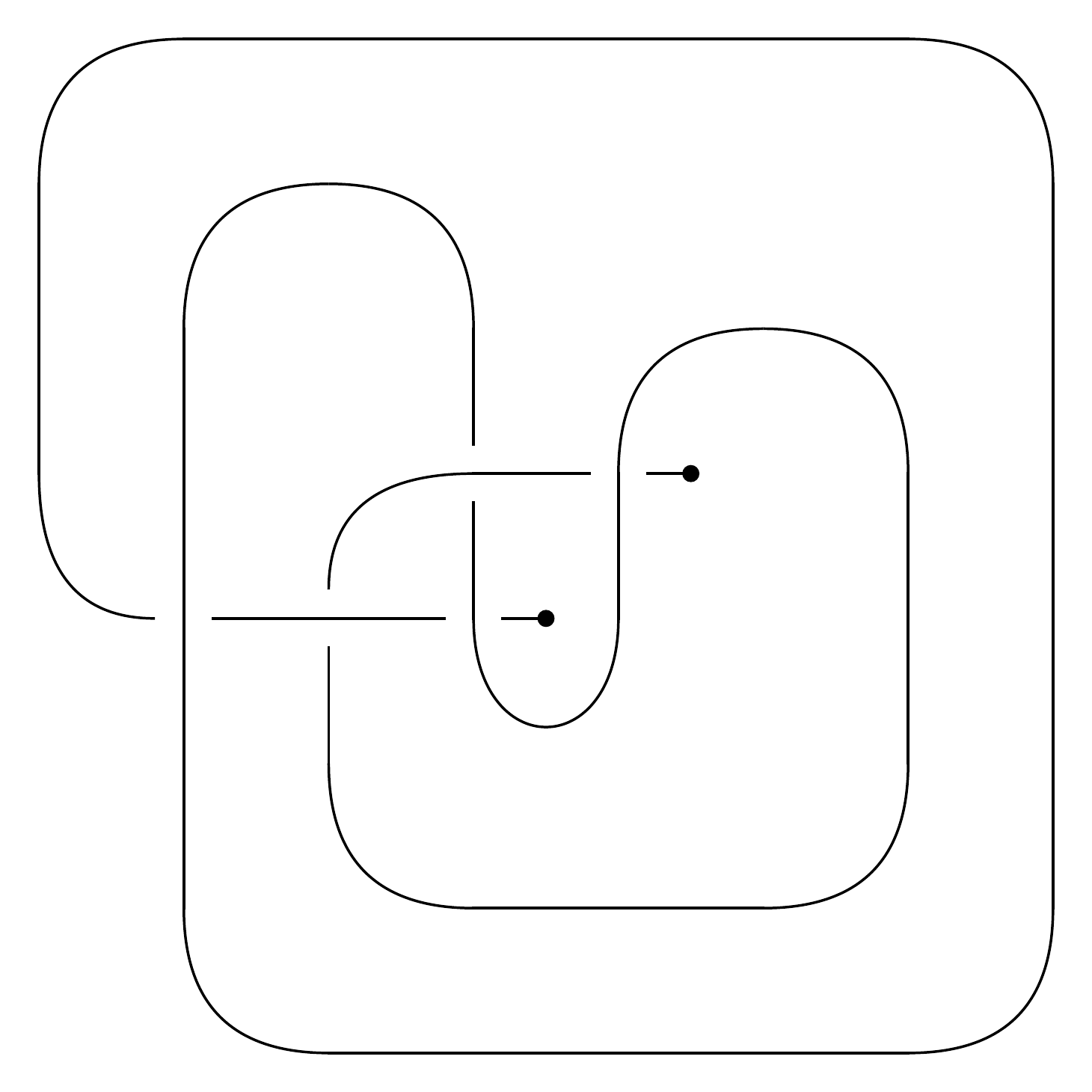}\\
\textcolor{black}{$5_{587}$}
\vspace{1cm}
\end{minipage}
\begin{minipage}[t]{.25\linewidth}
\centering
\includegraphics[width=0.9\textwidth,height=3.5cm,keepaspectratio]{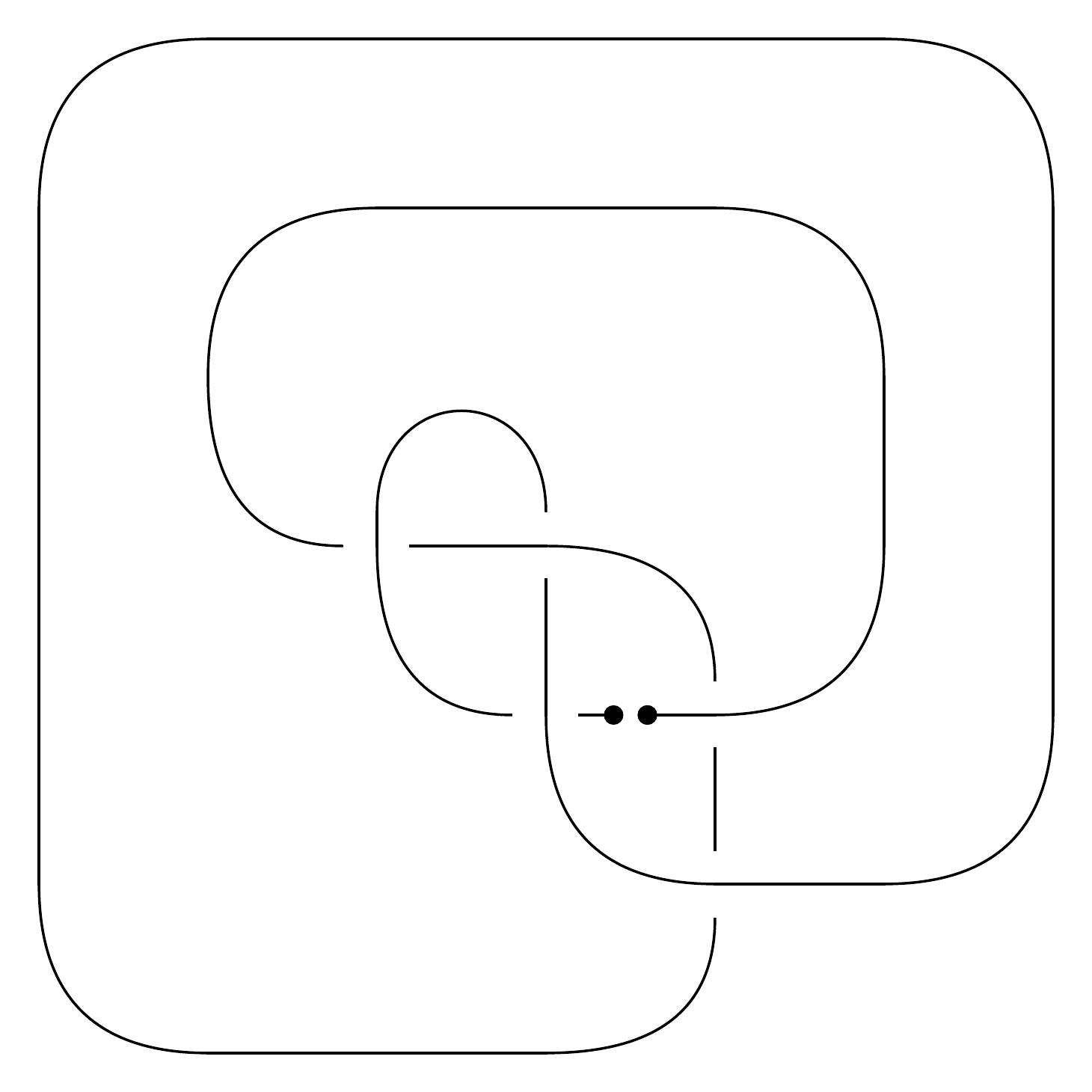}\\
\textcolor{black}{$5_{588}$}
\vspace{1cm}
\end{minipage}
\begin{minipage}[t]{.25\linewidth}
\centering
\includegraphics[width=0.9\textwidth,height=3.5cm,keepaspectratio]{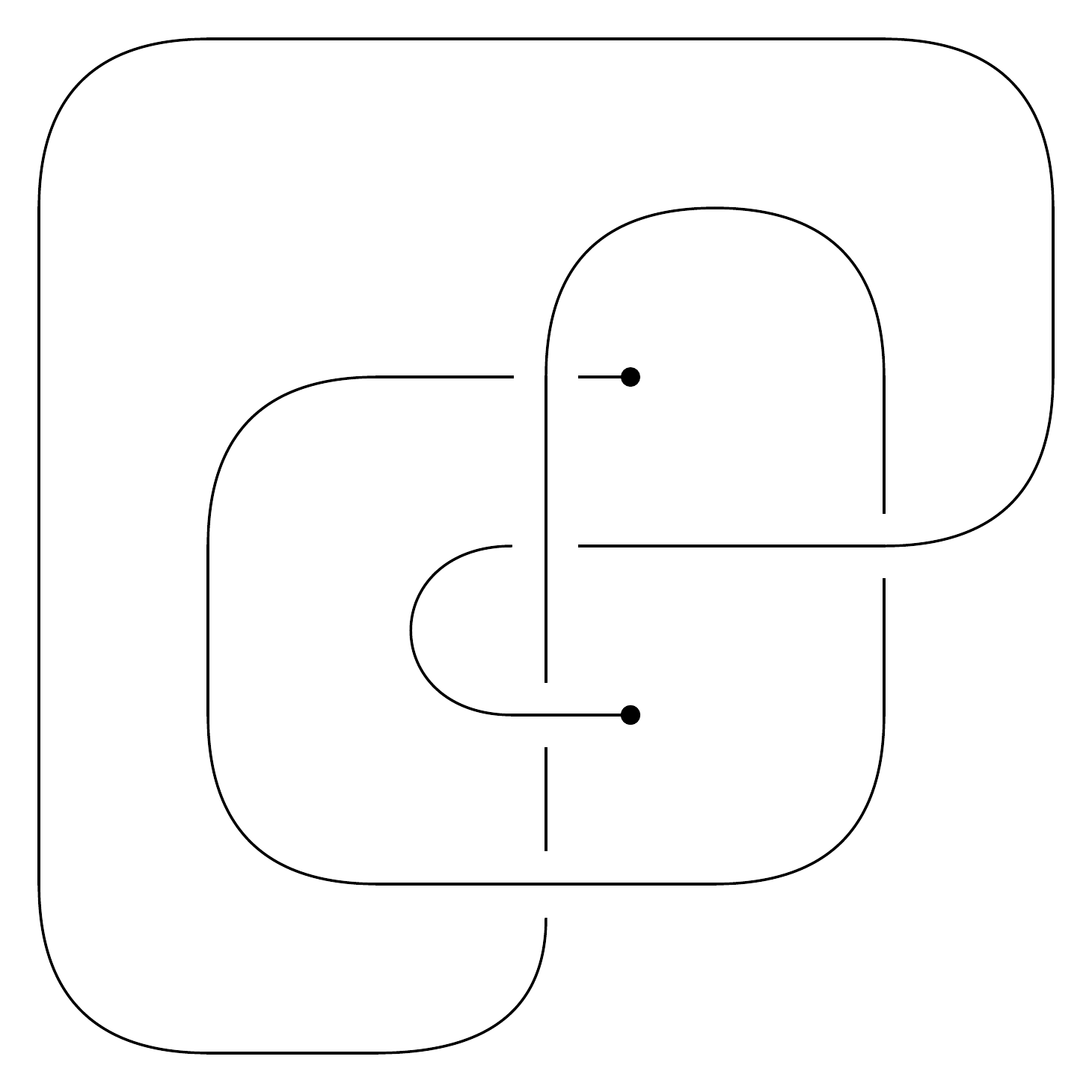}\\
\textcolor{black}{$5_{589}$}
\vspace{1cm}
\end{minipage}
\begin{minipage}[t]{.25\linewidth}
\centering
\includegraphics[width=0.9\textwidth,height=3.5cm,keepaspectratio]{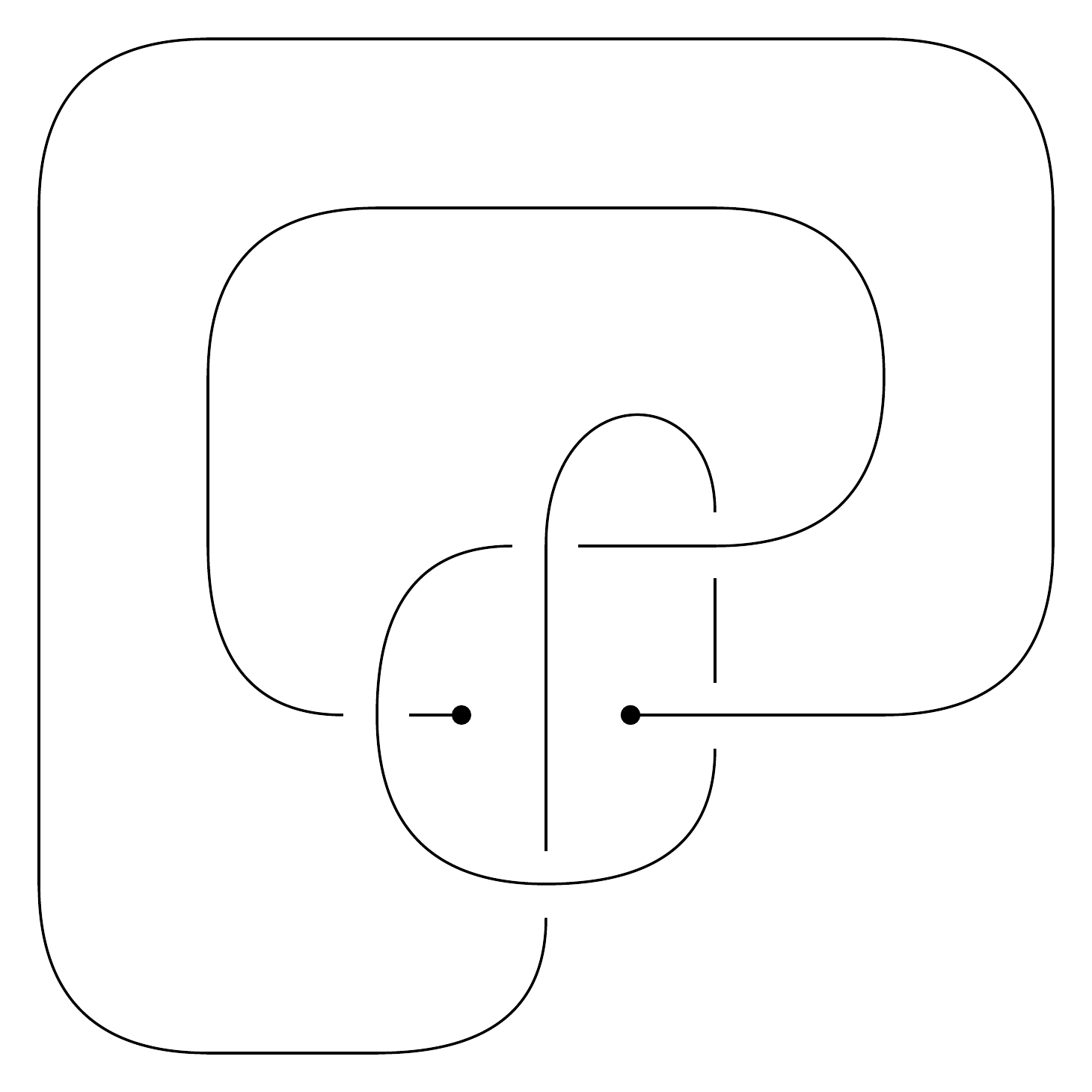}\\
\textcolor{black}{$5_{590}$}
\vspace{1cm}
\end{minipage}
\begin{minipage}[t]{.25\linewidth}
\centering
\includegraphics[width=0.9\textwidth,height=3.5cm,keepaspectratio]{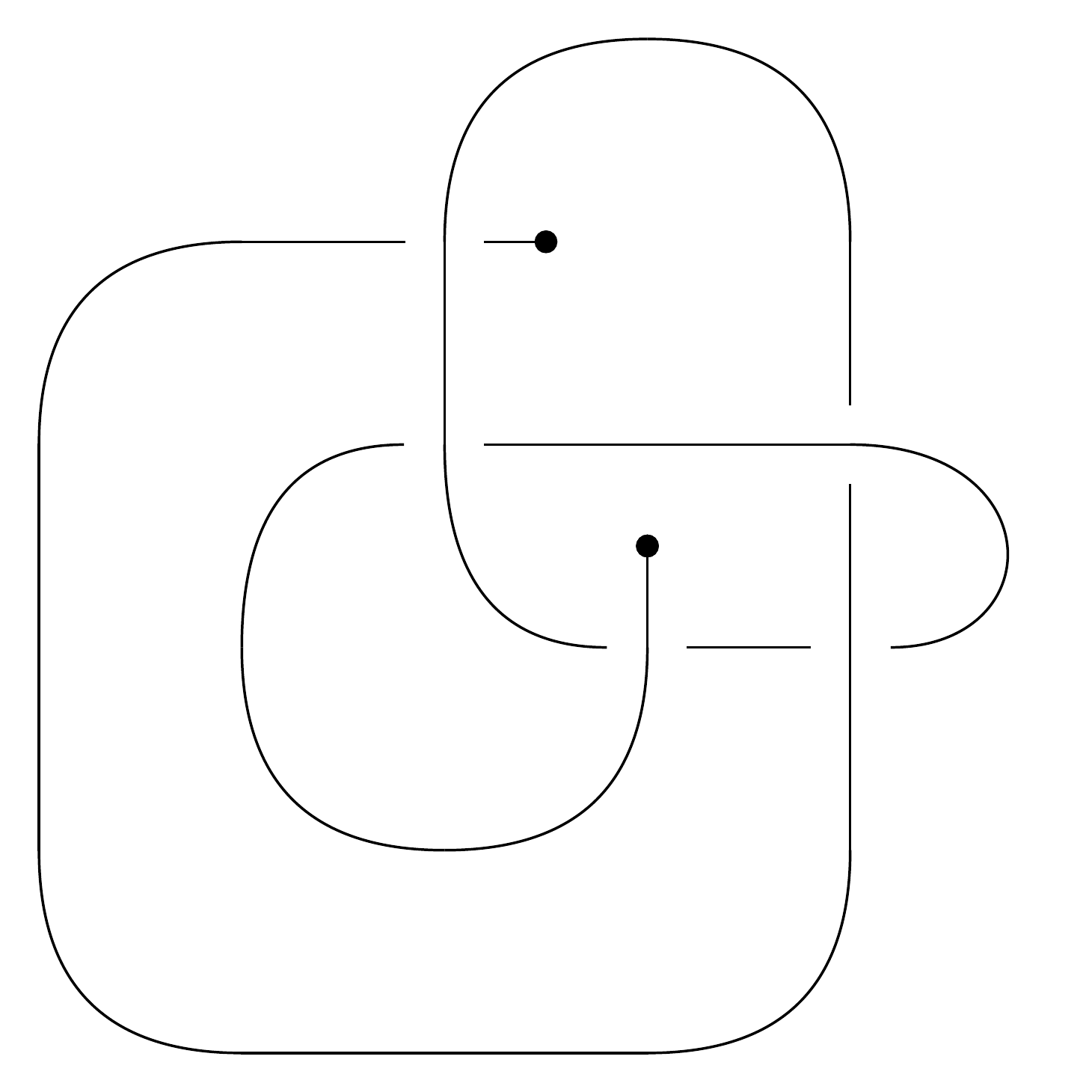}\\
\textcolor{black}{$5_{591}$}
\vspace{1cm}
\end{minipage}
\begin{minipage}[t]{.25\linewidth}
\centering
\includegraphics[width=0.9\textwidth,height=3.5cm,keepaspectratio]{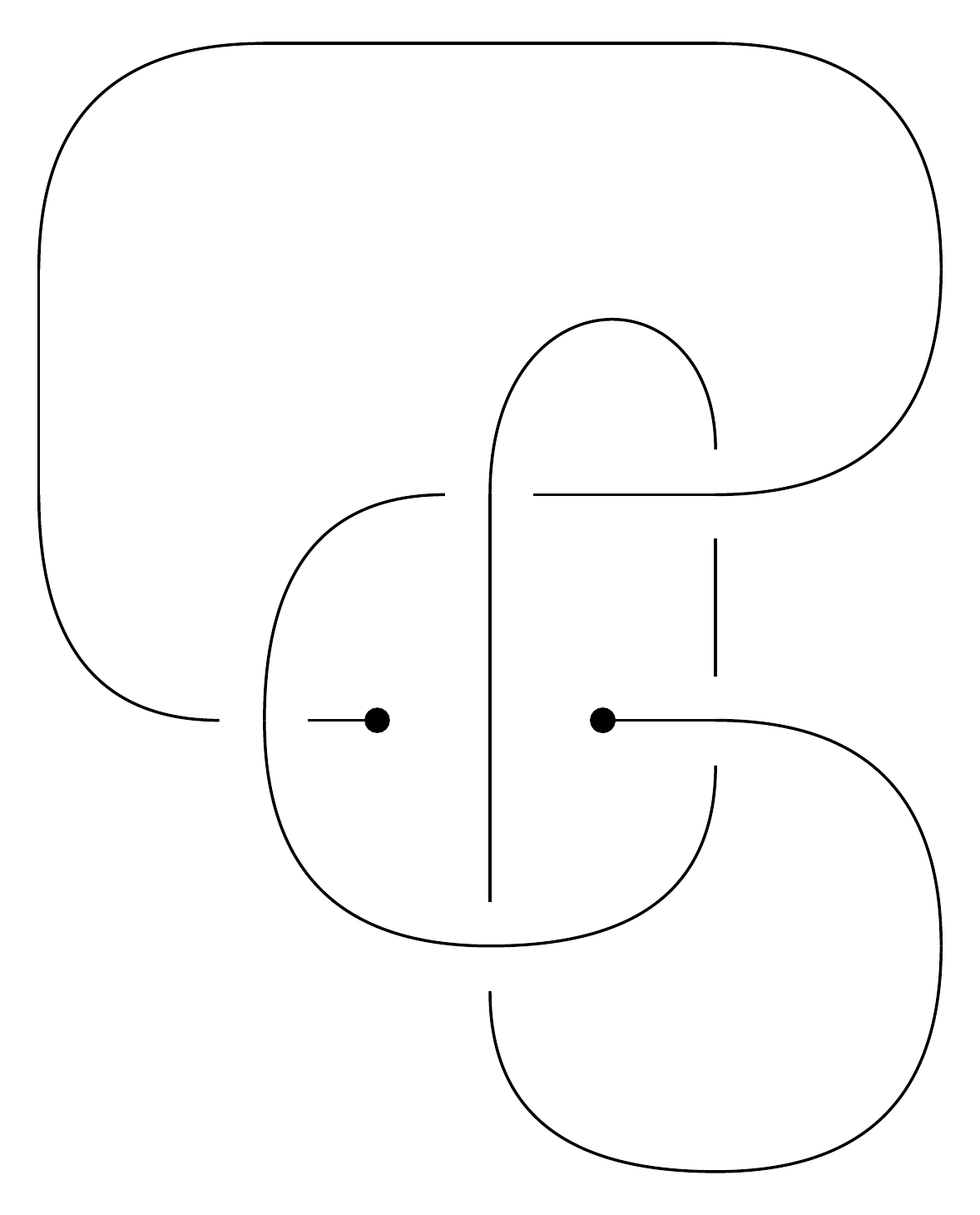}\\
\textcolor{black}{$5_{592}$}
\vspace{1cm}
\end{minipage}
\begin{minipage}[t]{.25\linewidth}
\centering
\includegraphics[width=0.9\textwidth,height=3.5cm,keepaspectratio]{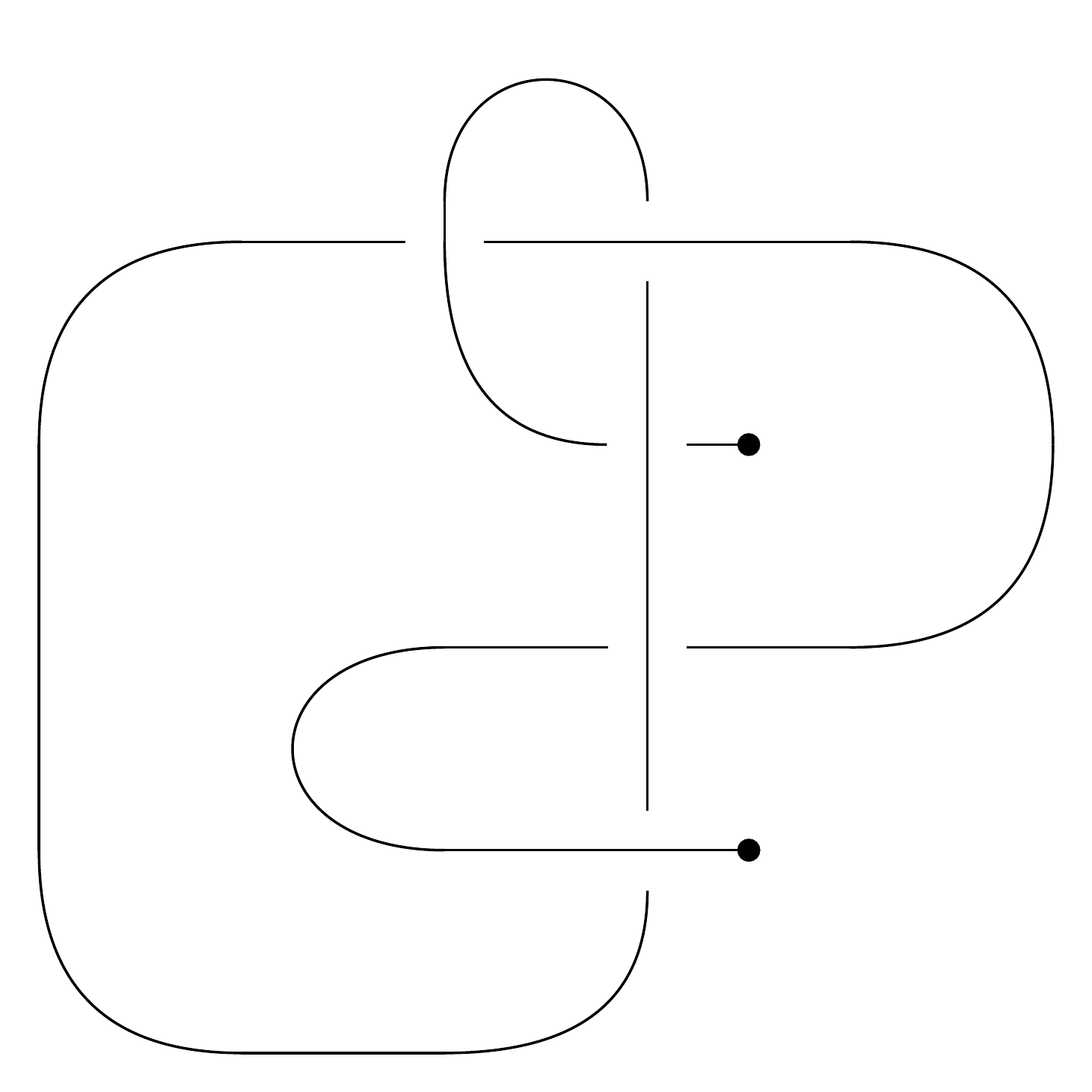}\\
\textcolor{black}{$5_{593}$}
\vspace{1cm}
\end{minipage}
\begin{minipage}[t]{.25\linewidth}
\centering
\includegraphics[width=0.9\textwidth,height=3.5cm,keepaspectratio]{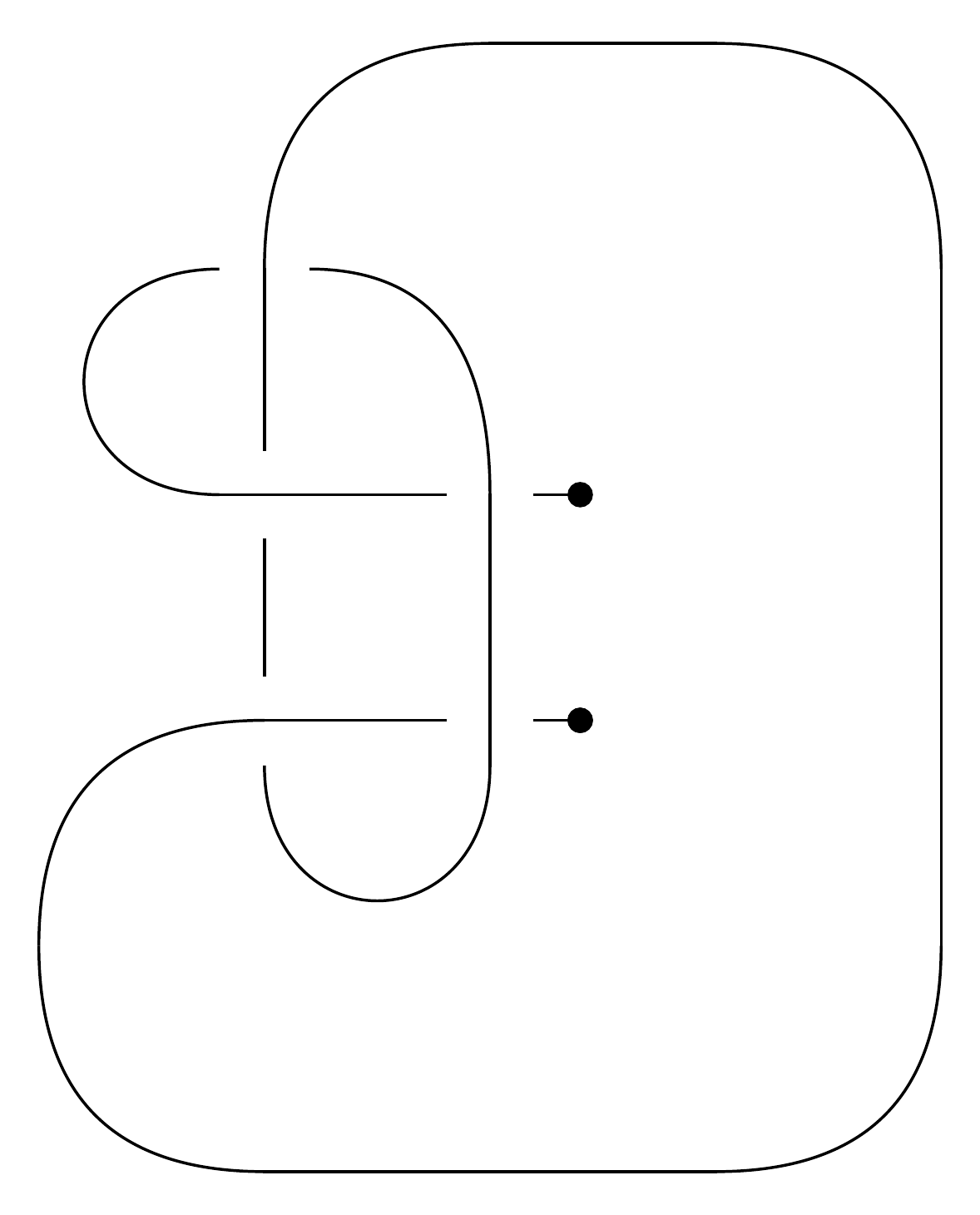}\\
\textcolor{black}{$5_{594}$}
\vspace{1cm}
\end{minipage}
\begin{minipage}[t]{.25\linewidth}
\centering
\includegraphics[width=0.9\textwidth,height=3.5cm,keepaspectratio]{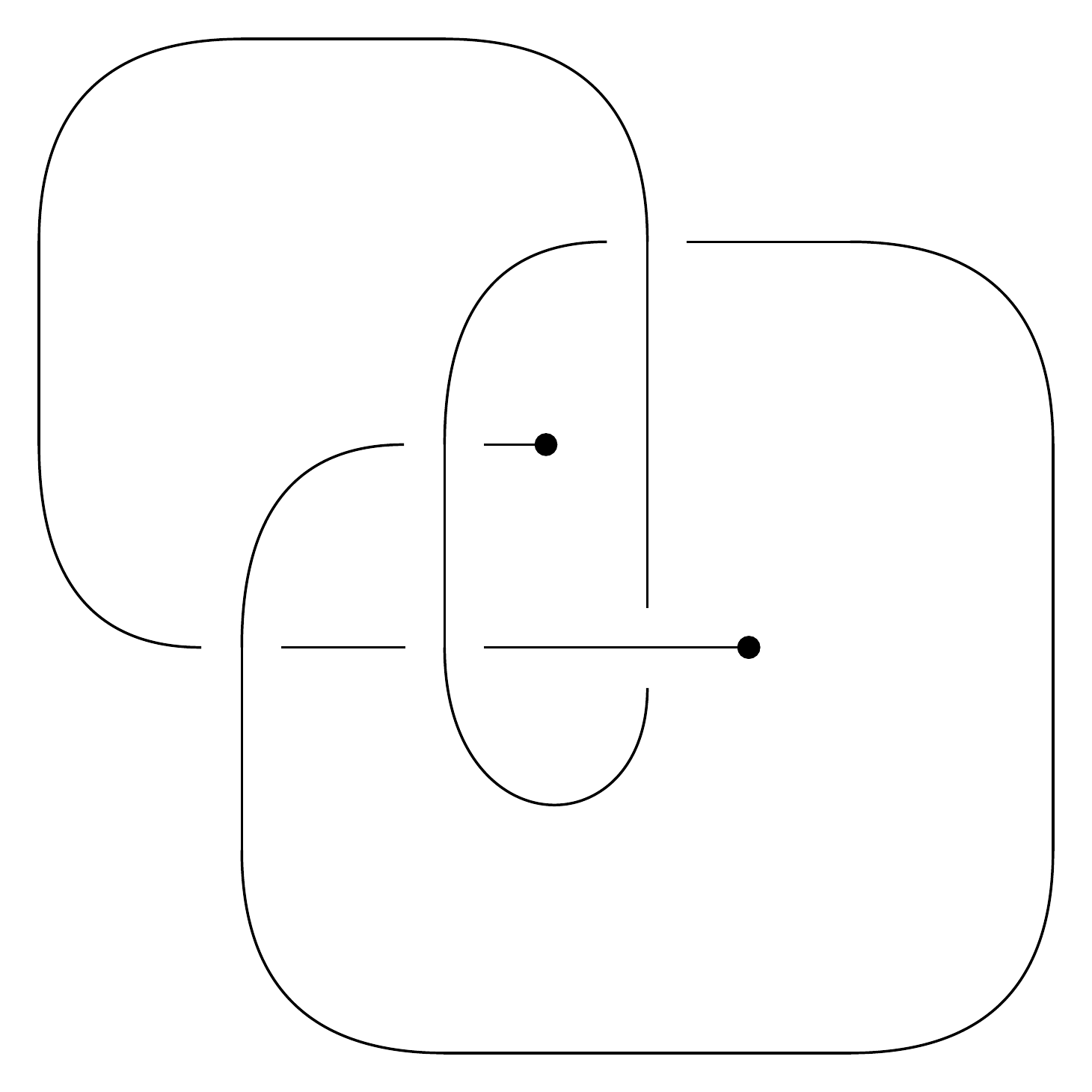}\\
\textcolor{black}{$5_{595}$}
\vspace{1cm}
\end{minipage}
\begin{minipage}[t]{.25\linewidth}
\centering
\includegraphics[width=0.9\textwidth,height=3.5cm,keepaspectratio]{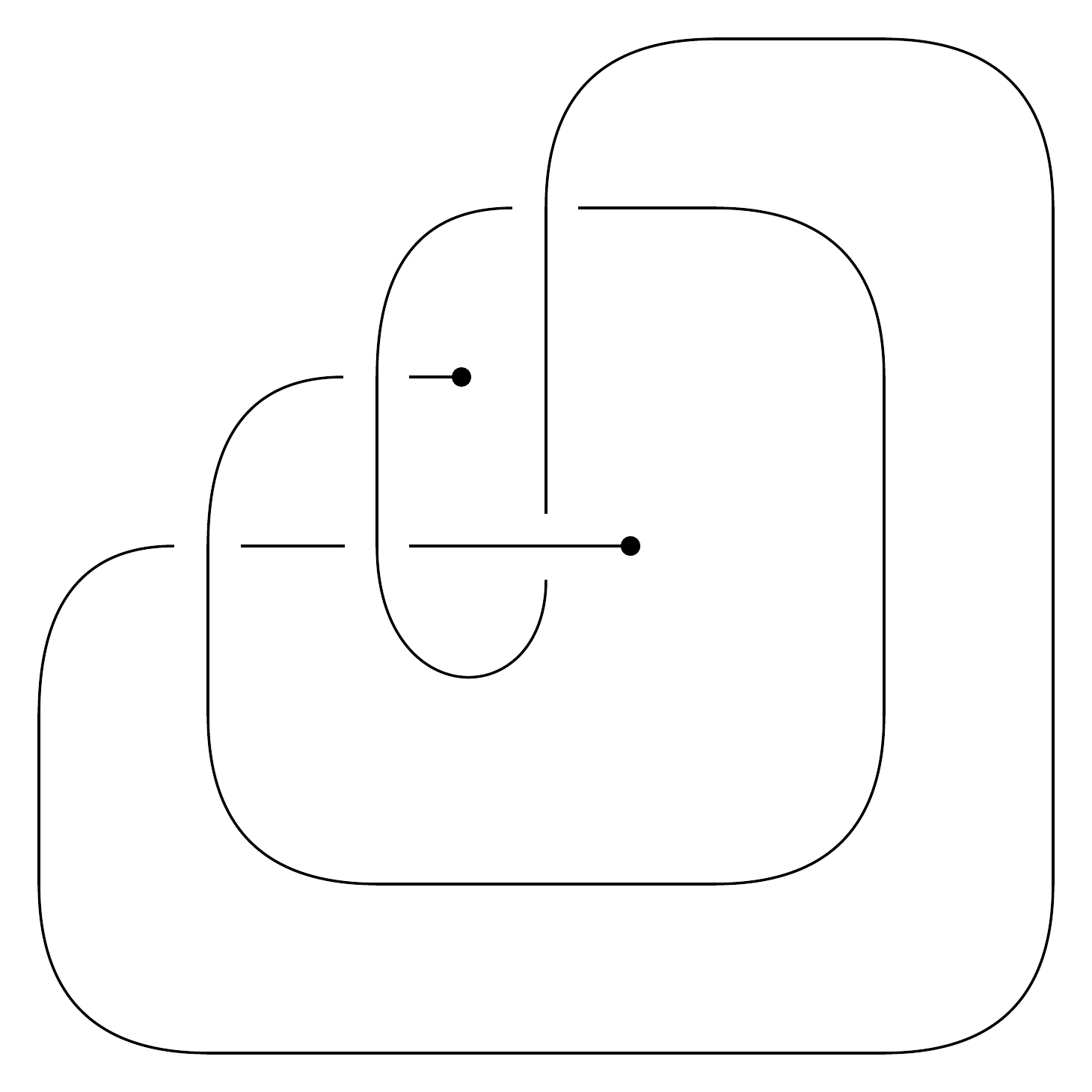}\\
\textcolor{black}{$5_{596}$}
\vspace{1cm}
\end{minipage}
\begin{minipage}[t]{.25\linewidth}
\centering
\includegraphics[width=0.9\textwidth,height=3.5cm,keepaspectratio]{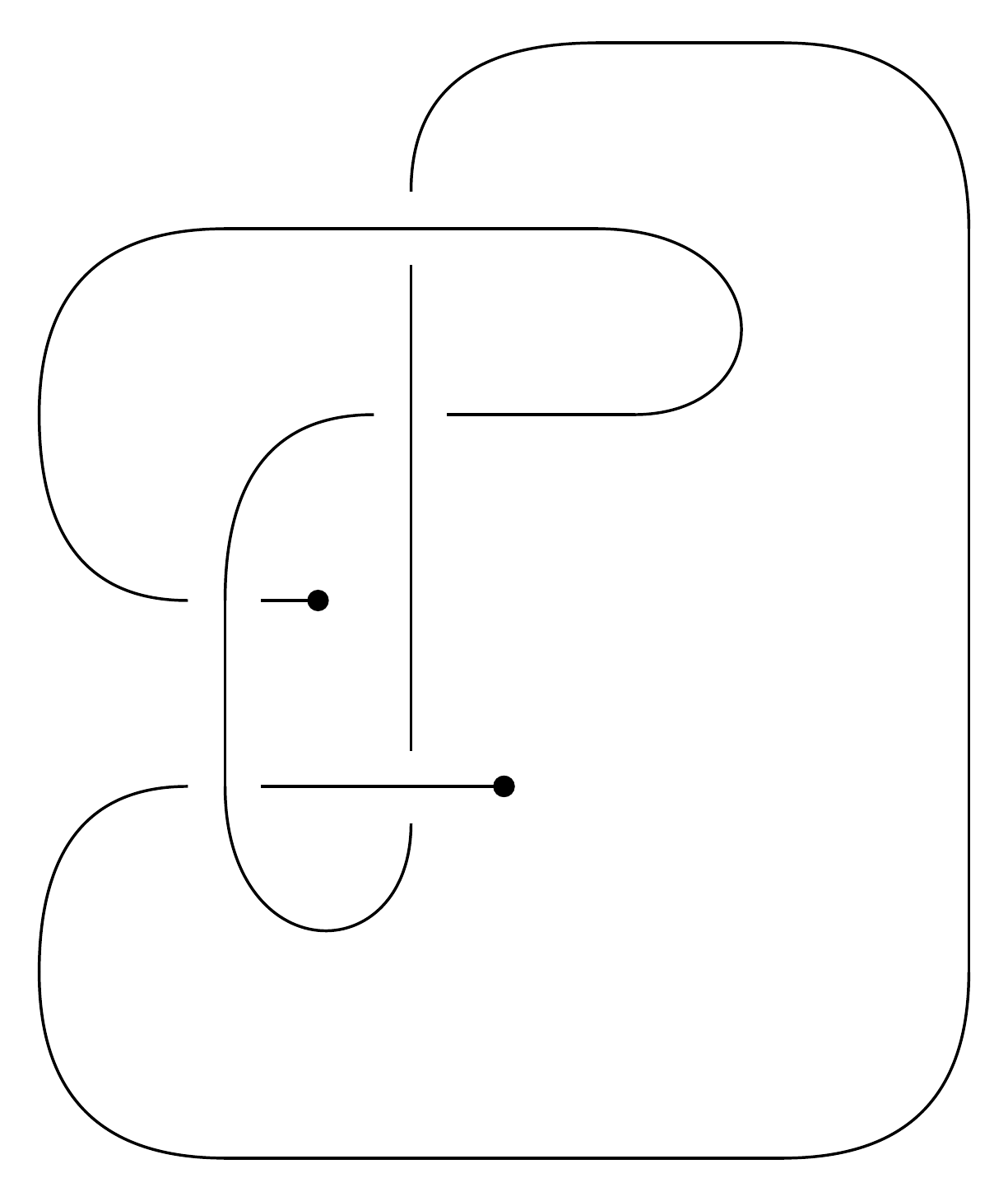}\\
\textcolor{black}{$5_{597}$}
\vspace{1cm}
\end{minipage}
\begin{minipage}[t]{.25\linewidth}
\centering
\includegraphics[width=0.9\textwidth,height=3.5cm,keepaspectratio]{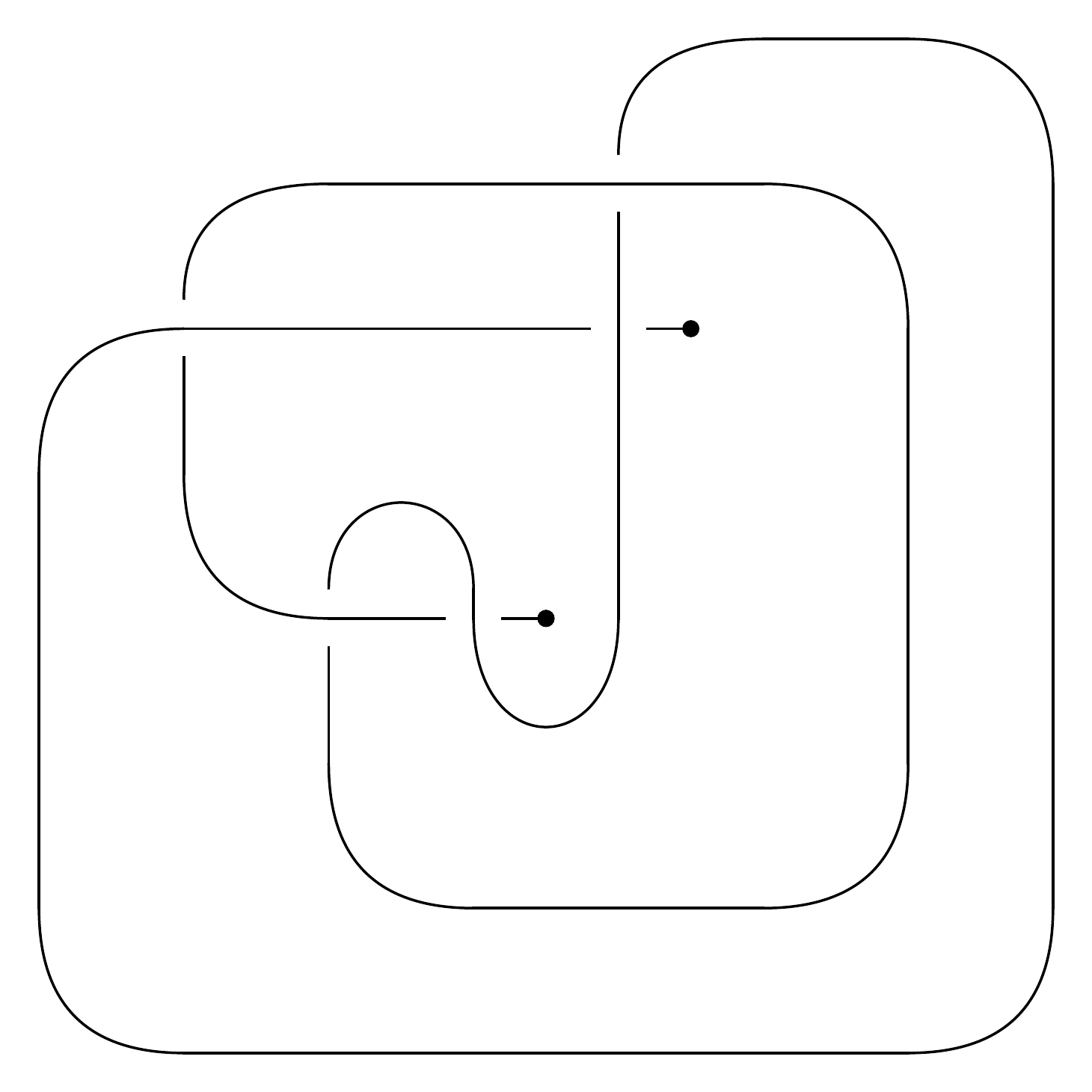}\\
\textcolor{black}{$5_{598}$}
\vspace{1cm}
\end{minipage}
\begin{minipage}[t]{.25\linewidth}
\centering
\includegraphics[width=0.9\textwidth,height=3.5cm,keepaspectratio]{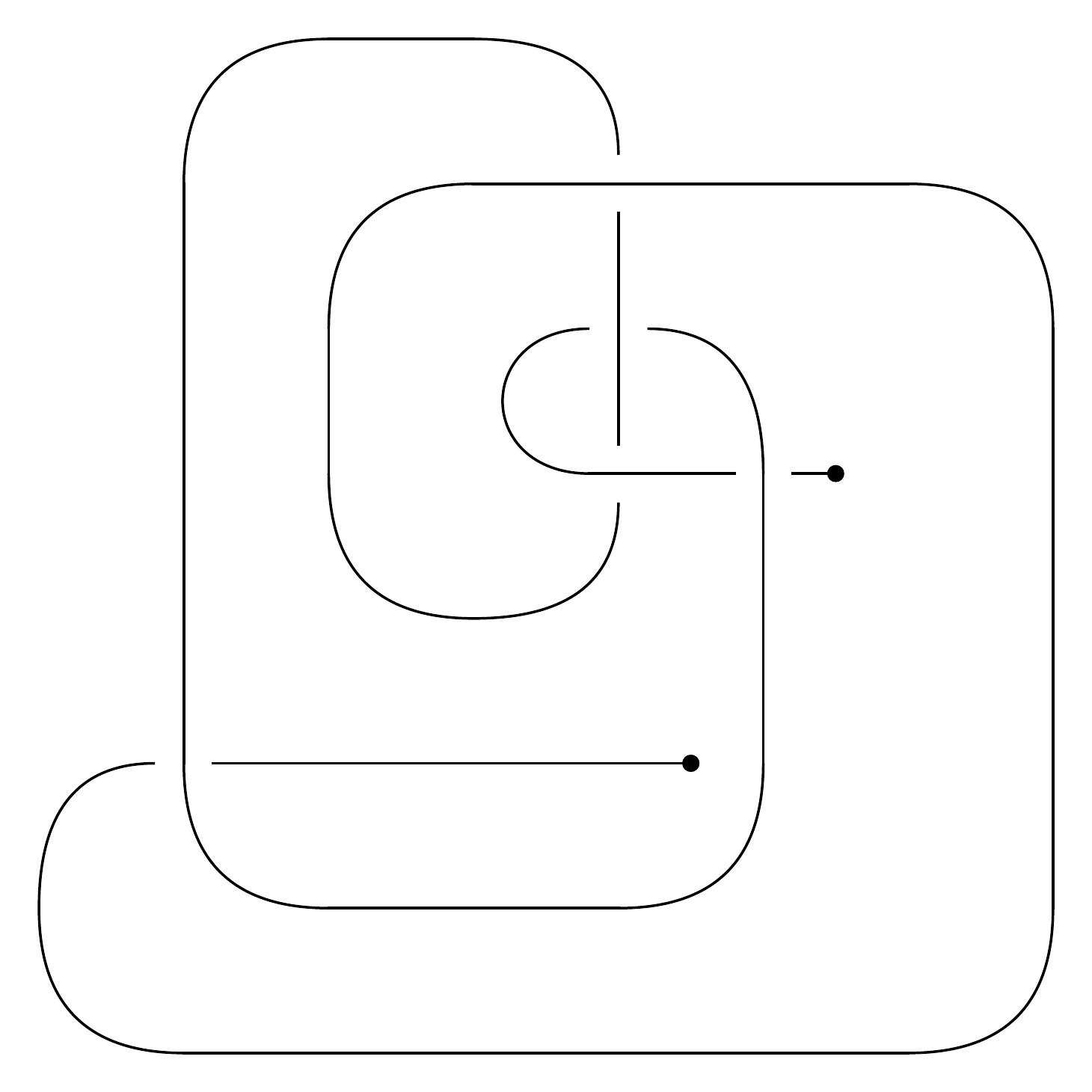}\\
\textcolor{black}{$5_{599}$}
\vspace{1cm}
\end{minipage}
\begin{minipage}[t]{.25\linewidth}
\centering
\includegraphics[width=0.9\textwidth,height=3.5cm,keepaspectratio]{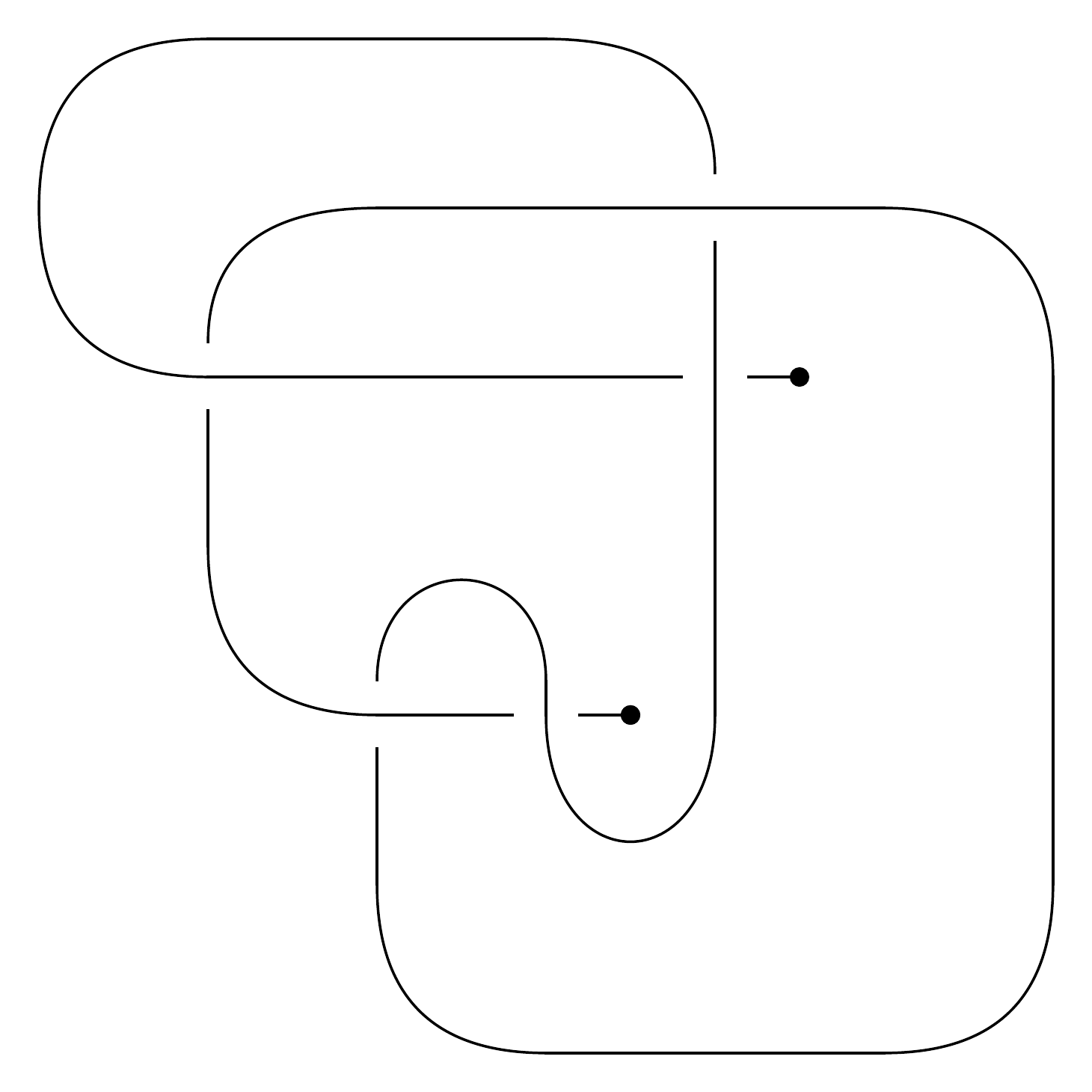}\\
\textcolor{black}{$5_{600}$}
\vspace{1cm}
\end{minipage}
\begin{minipage}[t]{.25\linewidth}
\centering
\includegraphics[width=0.9\textwidth,height=3.5cm,keepaspectratio]{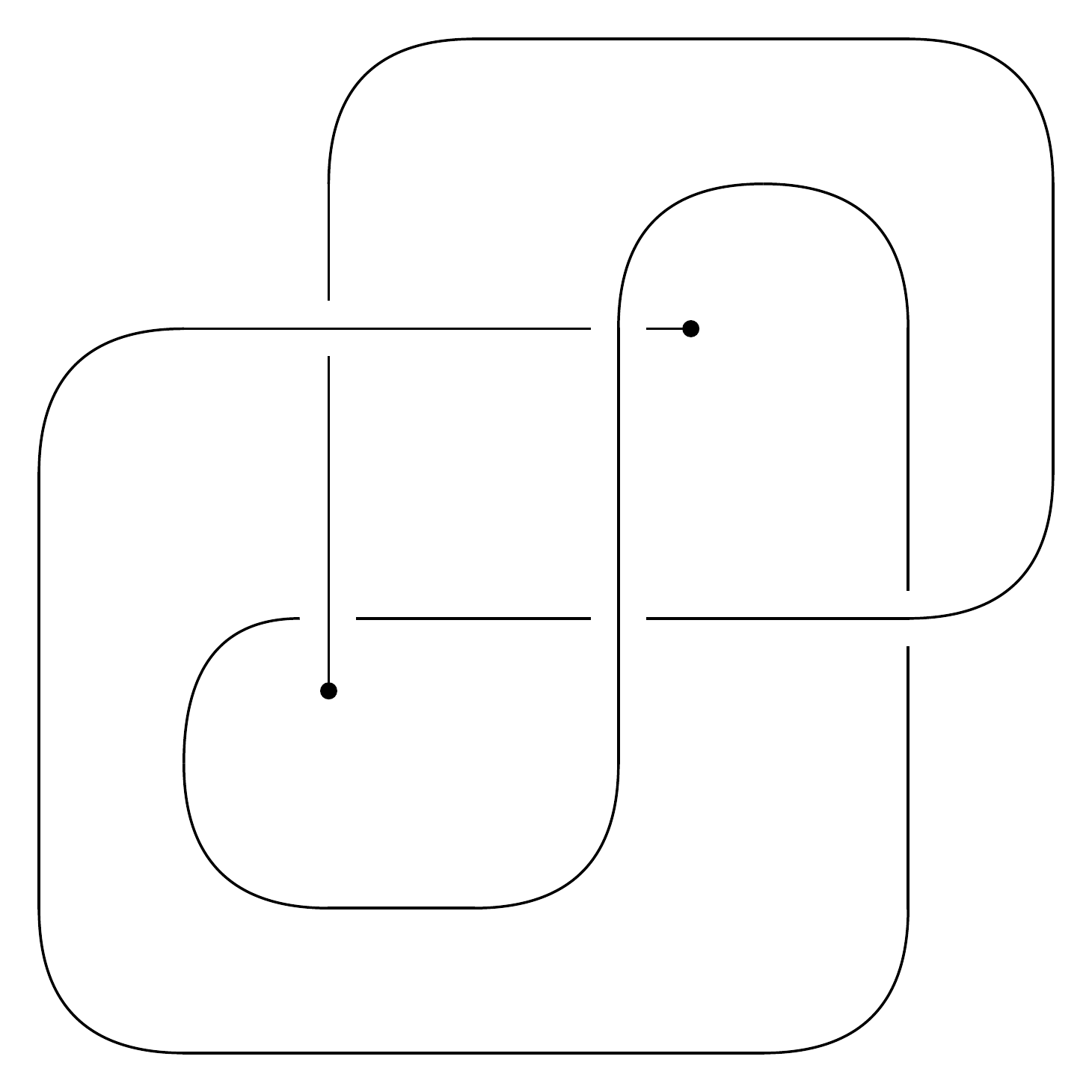}\\
\textcolor{black}{$5_{601}$}
\vspace{1cm}
\end{minipage}
\begin{minipage}[t]{.25\linewidth}
\centering
\includegraphics[width=0.9\textwidth,height=3.5cm,keepaspectratio]{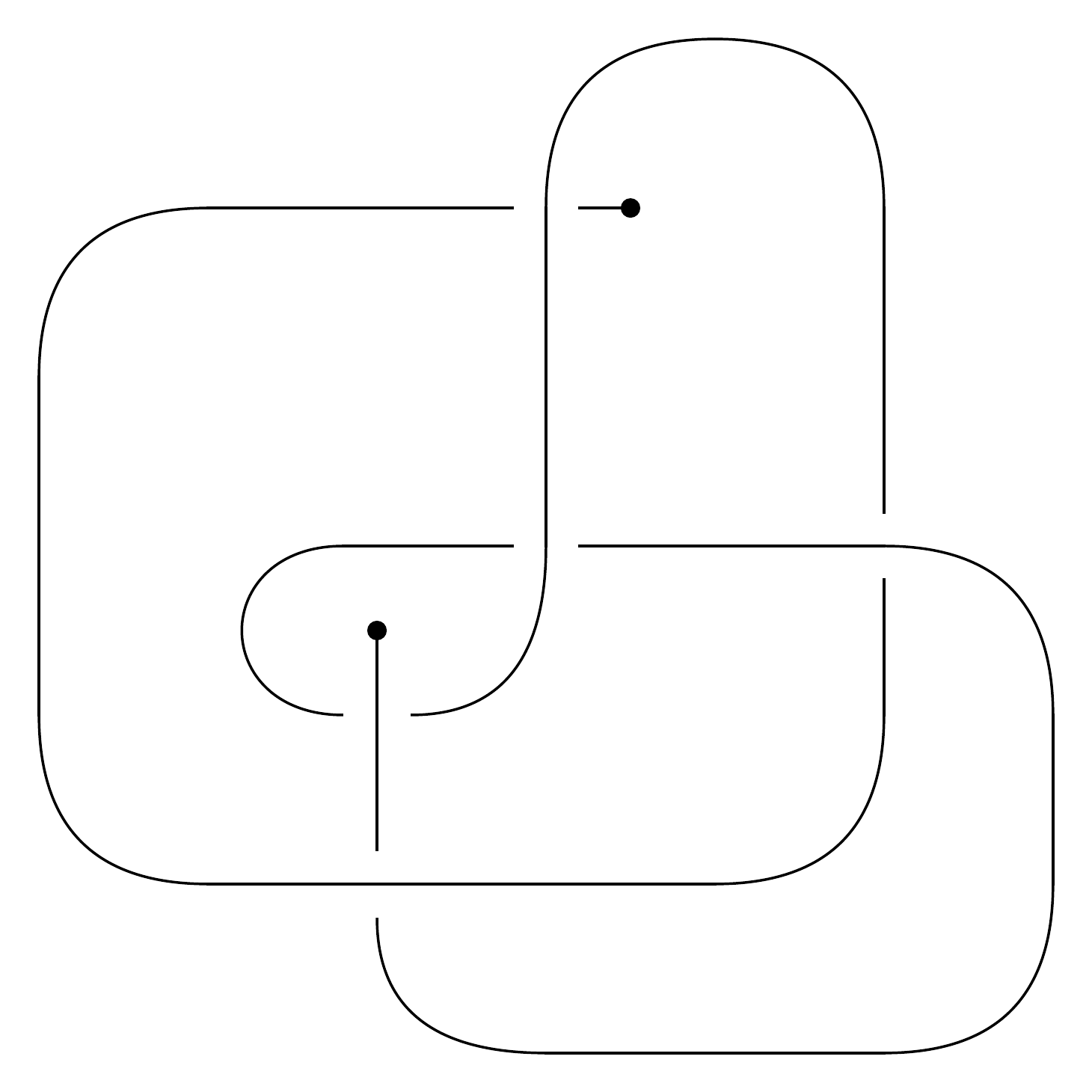}\\
\textcolor{black}{$5_{602}$}
\vspace{1cm}
\end{minipage}
\begin{minipage}[t]{.25\linewidth}
\centering
\includegraphics[width=0.9\textwidth,height=3.5cm,keepaspectratio]{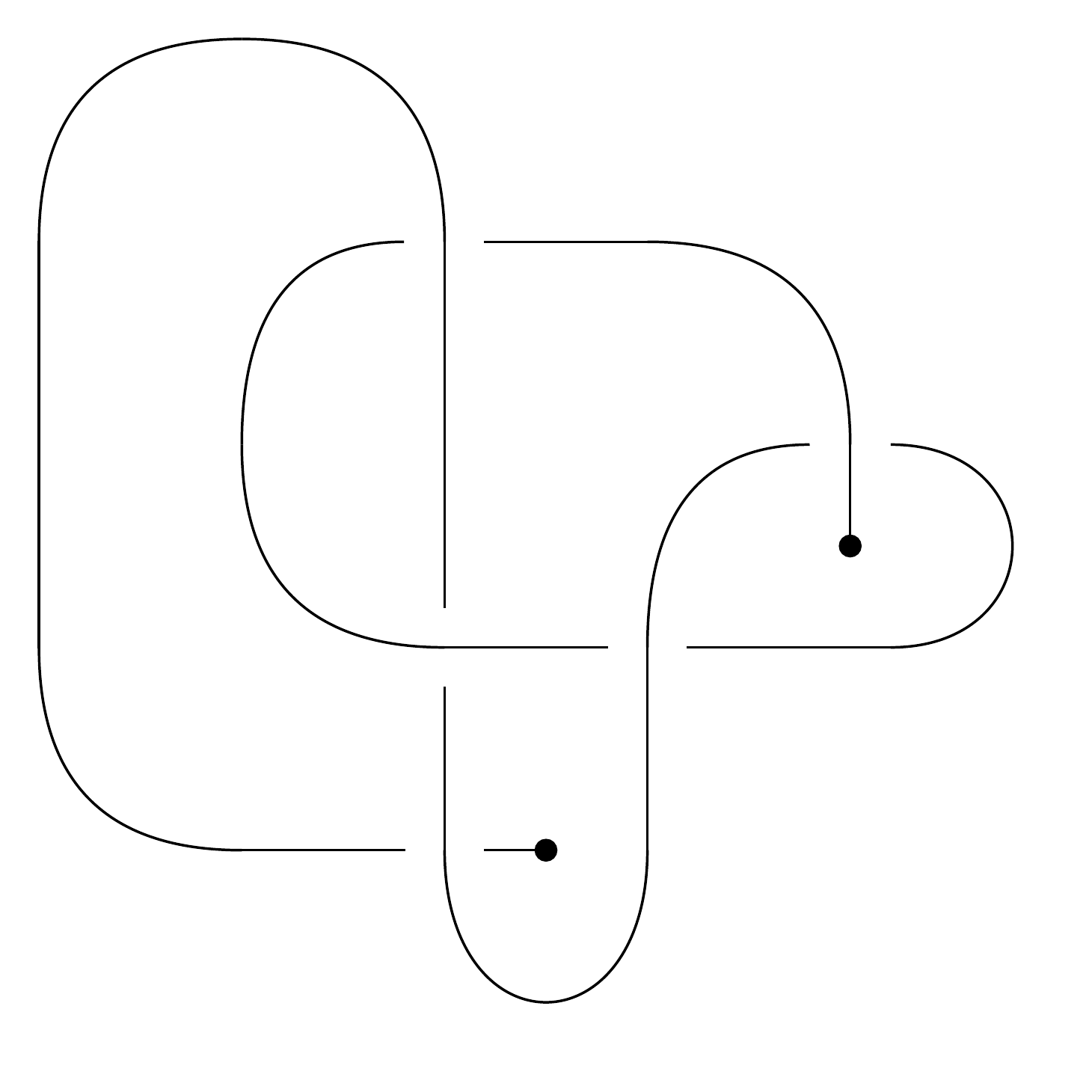}\\
\textcolor{black}{$5_{603}$}
\vspace{1cm}
\end{minipage}
\begin{minipage}[t]{.25\linewidth}
\centering
\includegraphics[width=0.9\textwidth,height=3.5cm,keepaspectratio]{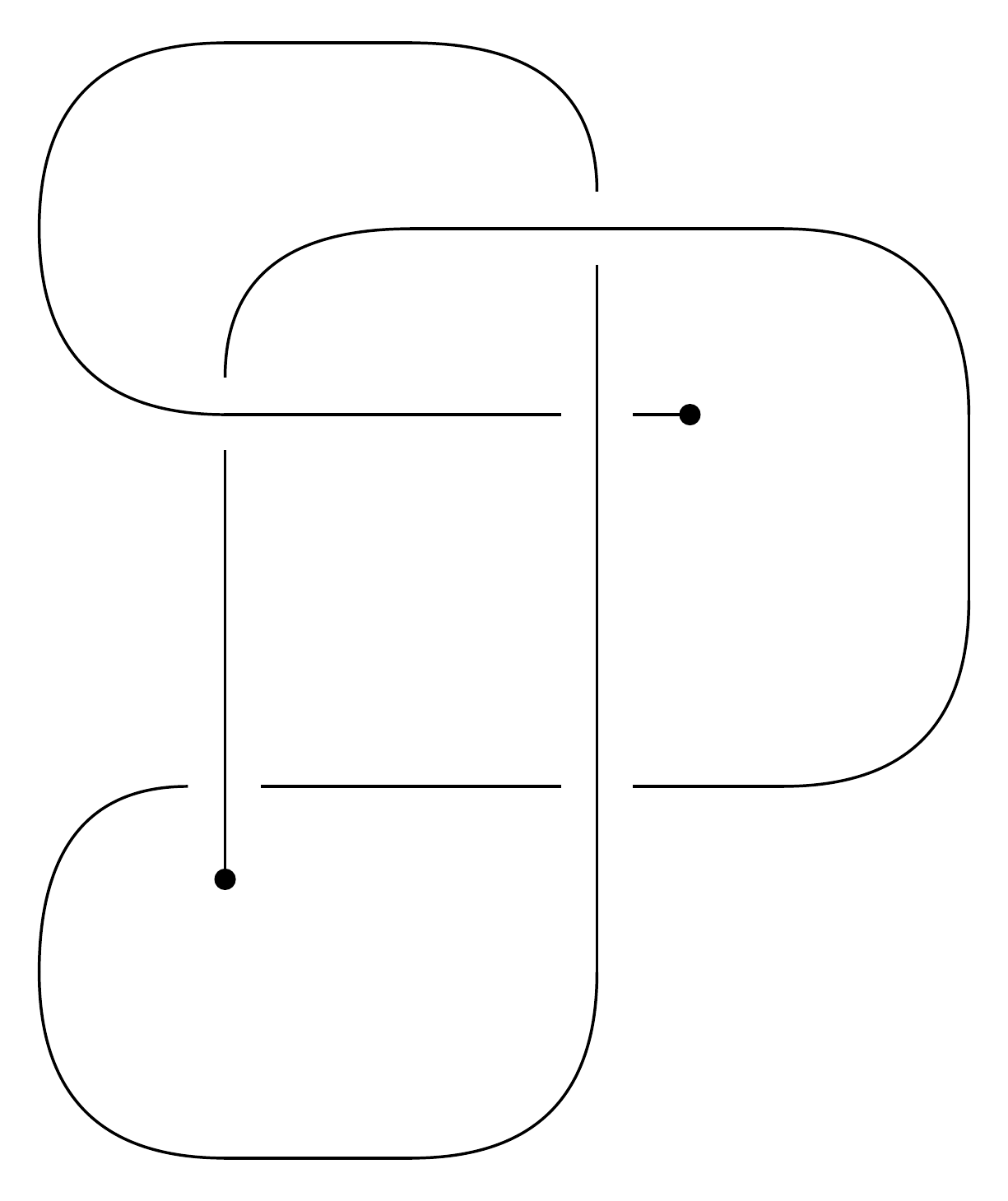}\\
\textcolor{black}{$5_{604}$}
\vspace{1cm}
\end{minipage}
\begin{minipage}[t]{.25\linewidth}
\centering
\includegraphics[width=0.9\textwidth,height=3.5cm,keepaspectratio]{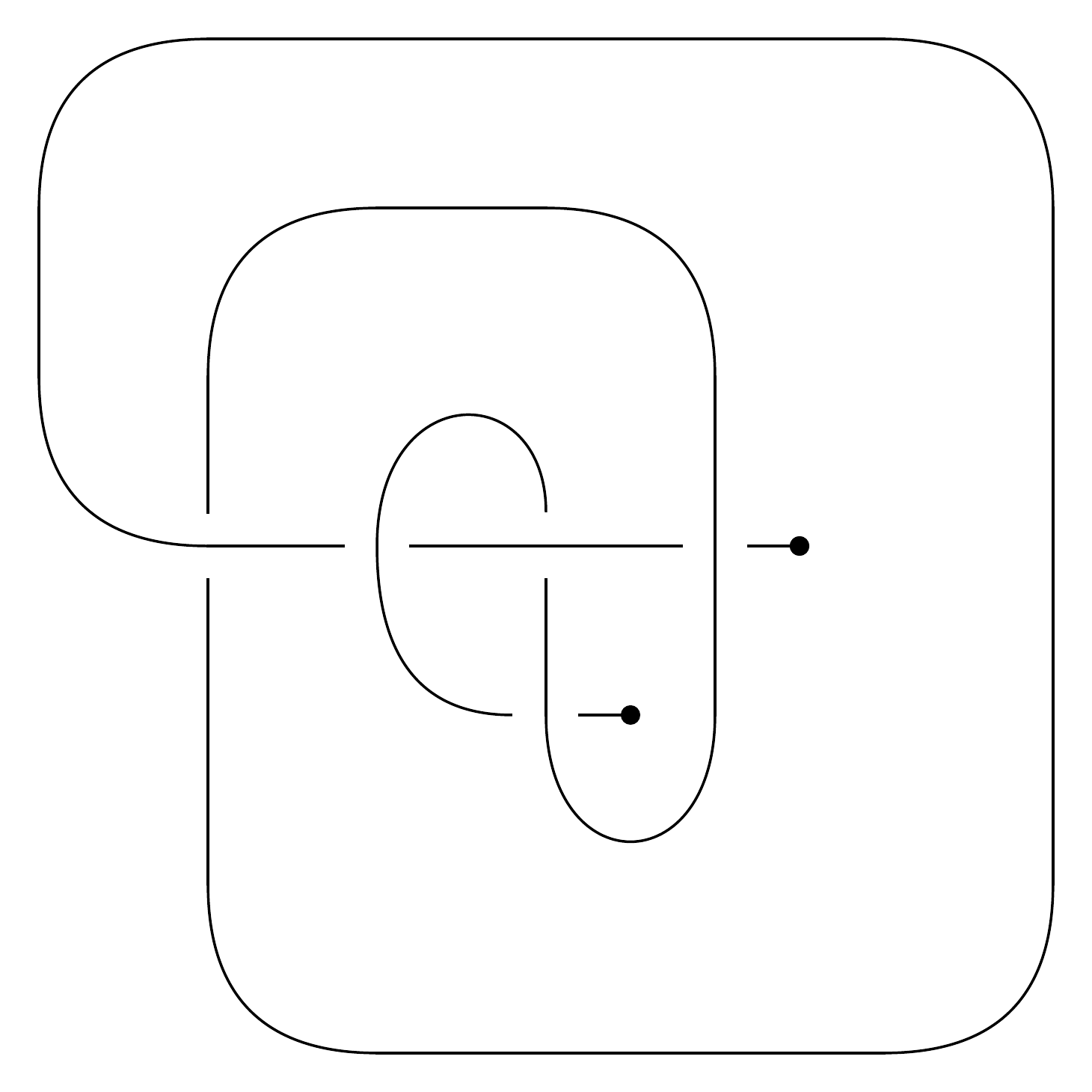}\\
\textcolor{black}{$5_{605}$}
\vspace{1cm}
\end{minipage}
\begin{minipage}[t]{.25\linewidth}
\centering
\includegraphics[width=0.9\textwidth,height=3.5cm,keepaspectratio]{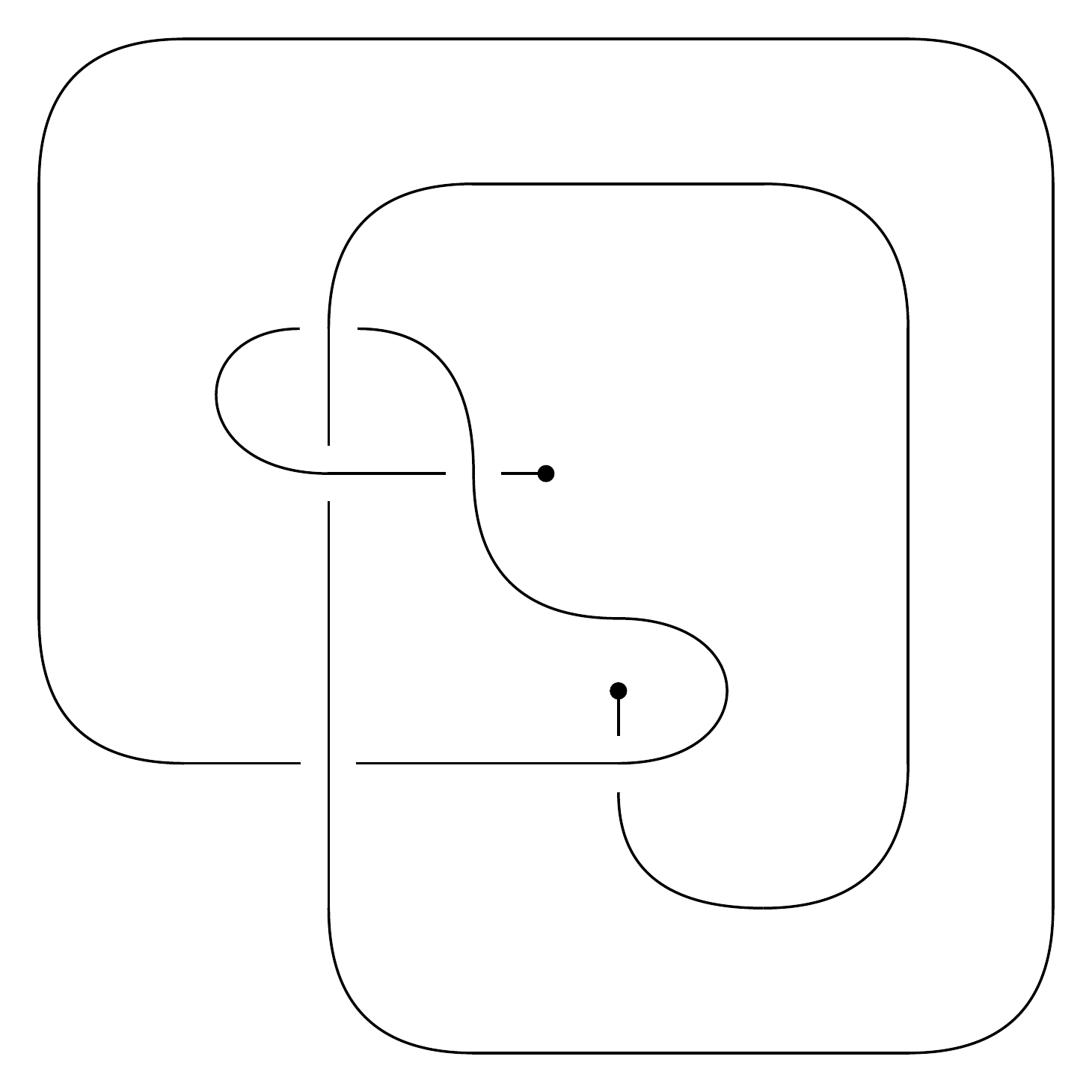}\\
\textcolor{black}{$5_{606}$}
\vspace{1cm}
\end{minipage}
\begin{minipage}[t]{.25\linewidth}
\centering
\includegraphics[width=0.9\textwidth,height=3.5cm,keepaspectratio]{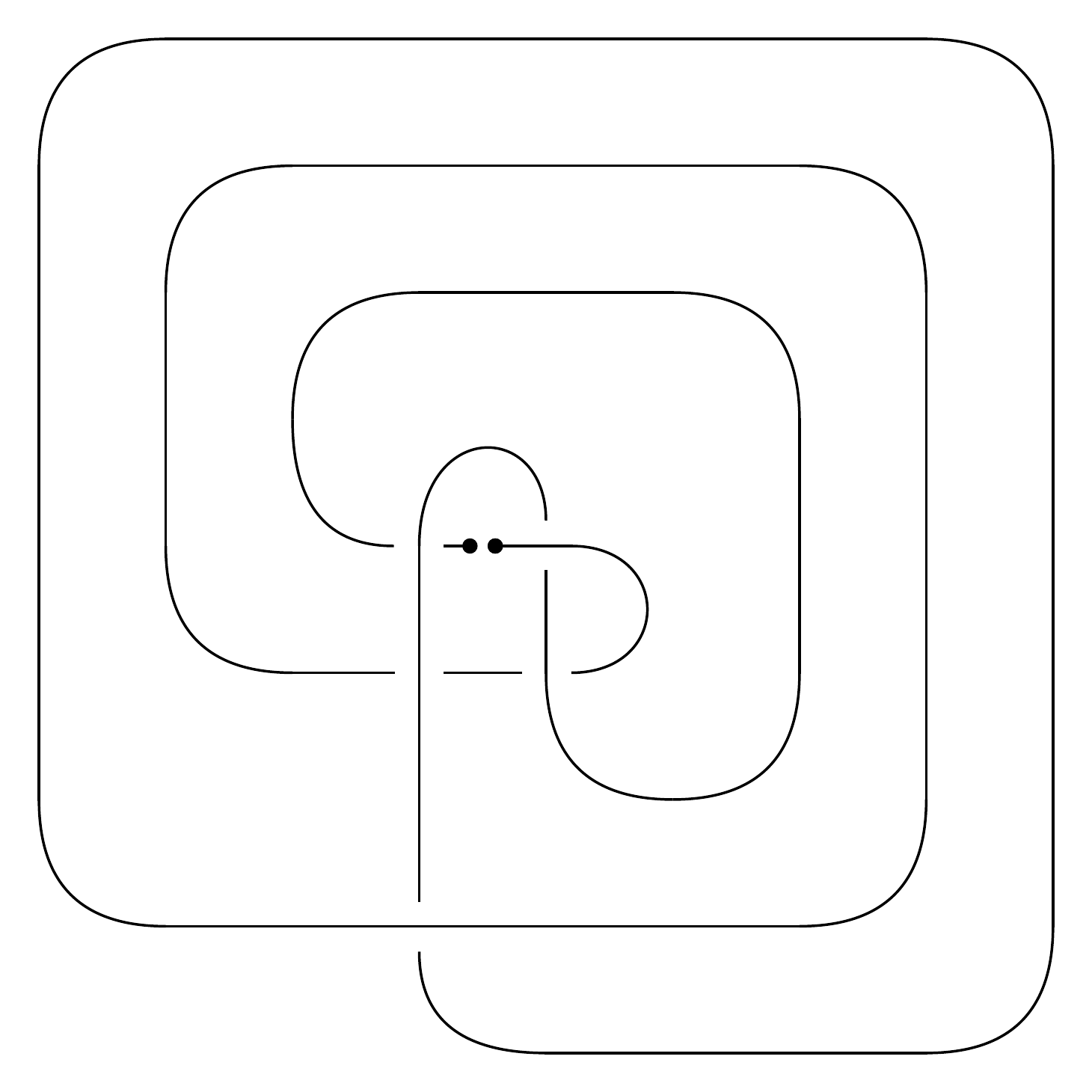}\\
\textcolor{black}{$5_{607}$}
\vspace{1cm}
\end{minipage}
\begin{minipage}[t]{.25\linewidth}
\centering
\includegraphics[width=0.9\textwidth,height=3.5cm,keepaspectratio]{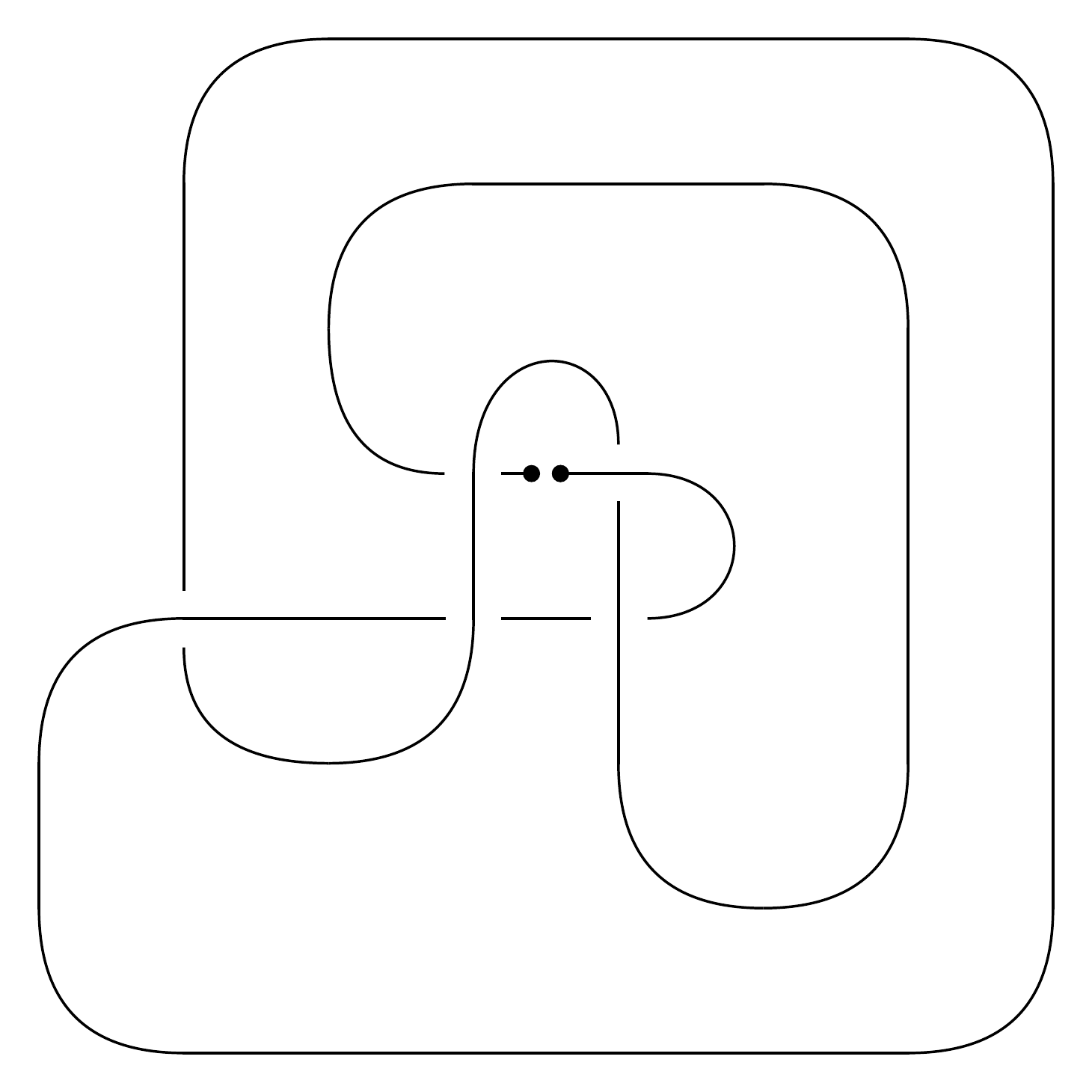}\\
\textcolor{black}{$5_{608}$}
\vspace{1cm}
\end{minipage}
\begin{minipage}[t]{.25\linewidth}
\centering
\includegraphics[width=0.9\textwidth,height=3.5cm,keepaspectratio]{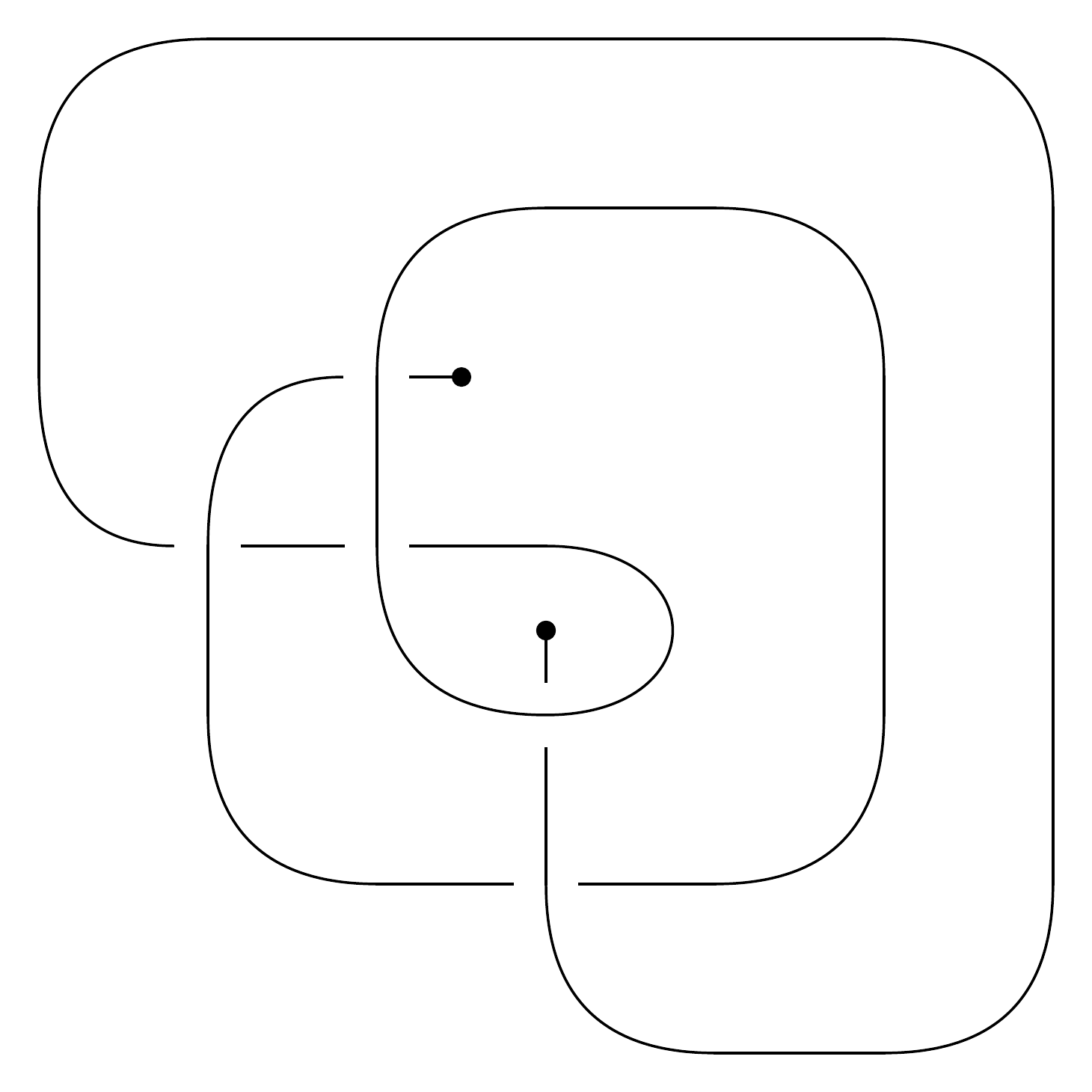}\\
\textcolor{black}{$5_{609}$}
\vspace{1cm}
\end{minipage}
\begin{minipage}[t]{.25\linewidth}
\centering
\includegraphics[width=0.9\textwidth,height=3.5cm,keepaspectratio]{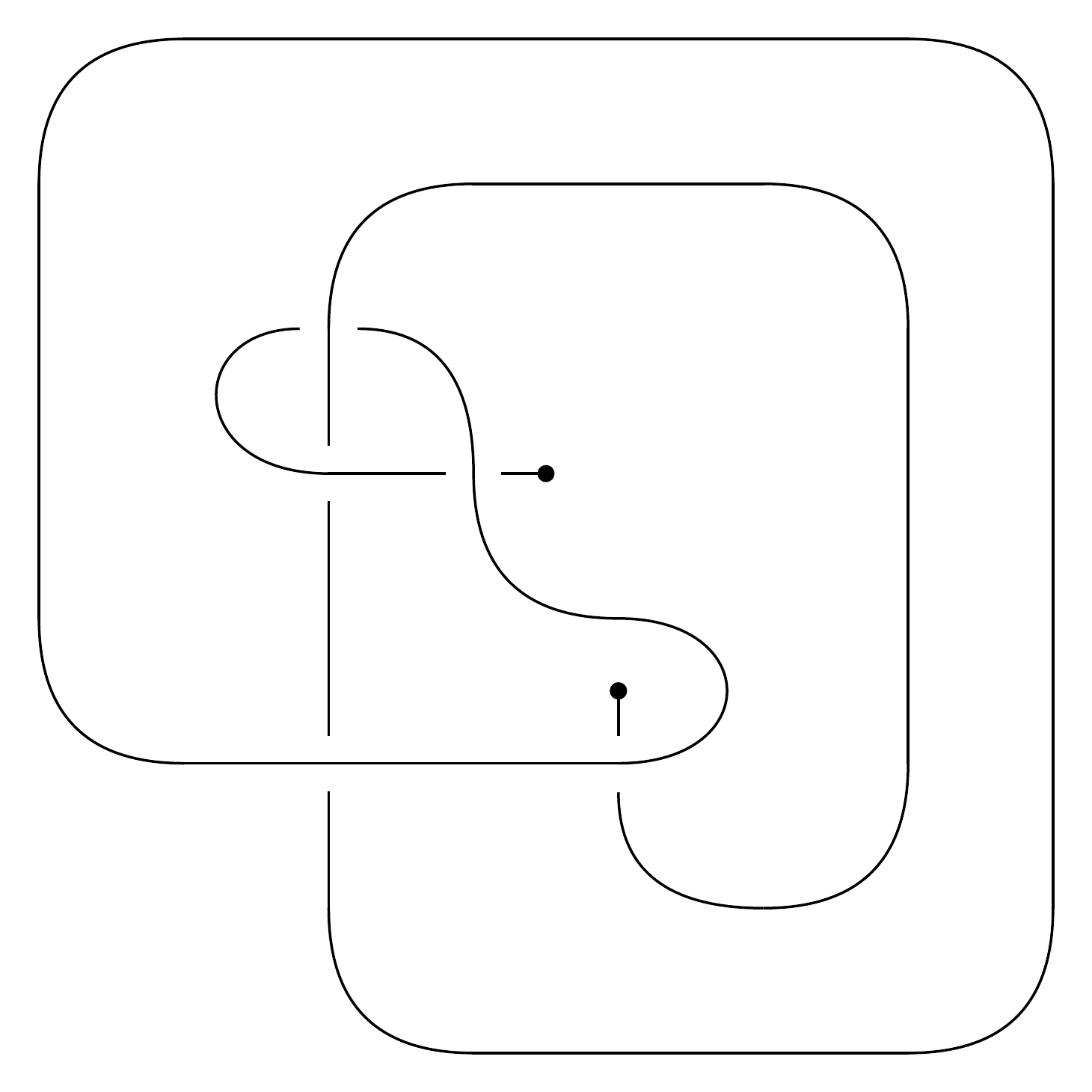}\\
\textcolor{black}{$5_{610}$}
\vspace{1cm}
\end{minipage}
\begin{minipage}[t]{.25\linewidth}
\centering
\includegraphics[width=0.9\textwidth,height=3.5cm,keepaspectratio]{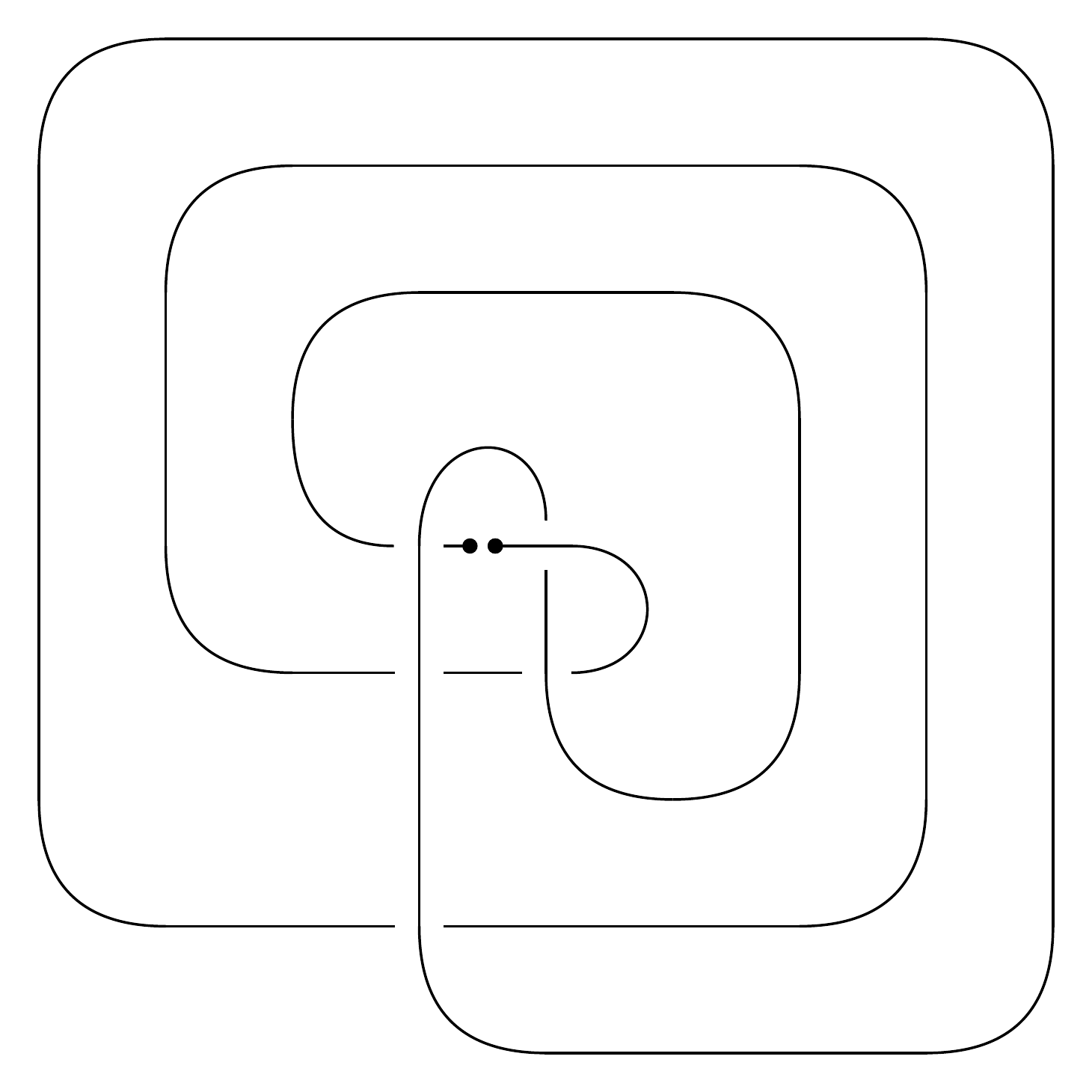}\\
\textcolor{black}{$5_{611}$}
\vspace{1cm}
\end{minipage}
\begin{minipage}[t]{.25\linewidth}
\centering
\includegraphics[width=0.9\textwidth,height=3.5cm,keepaspectratio]{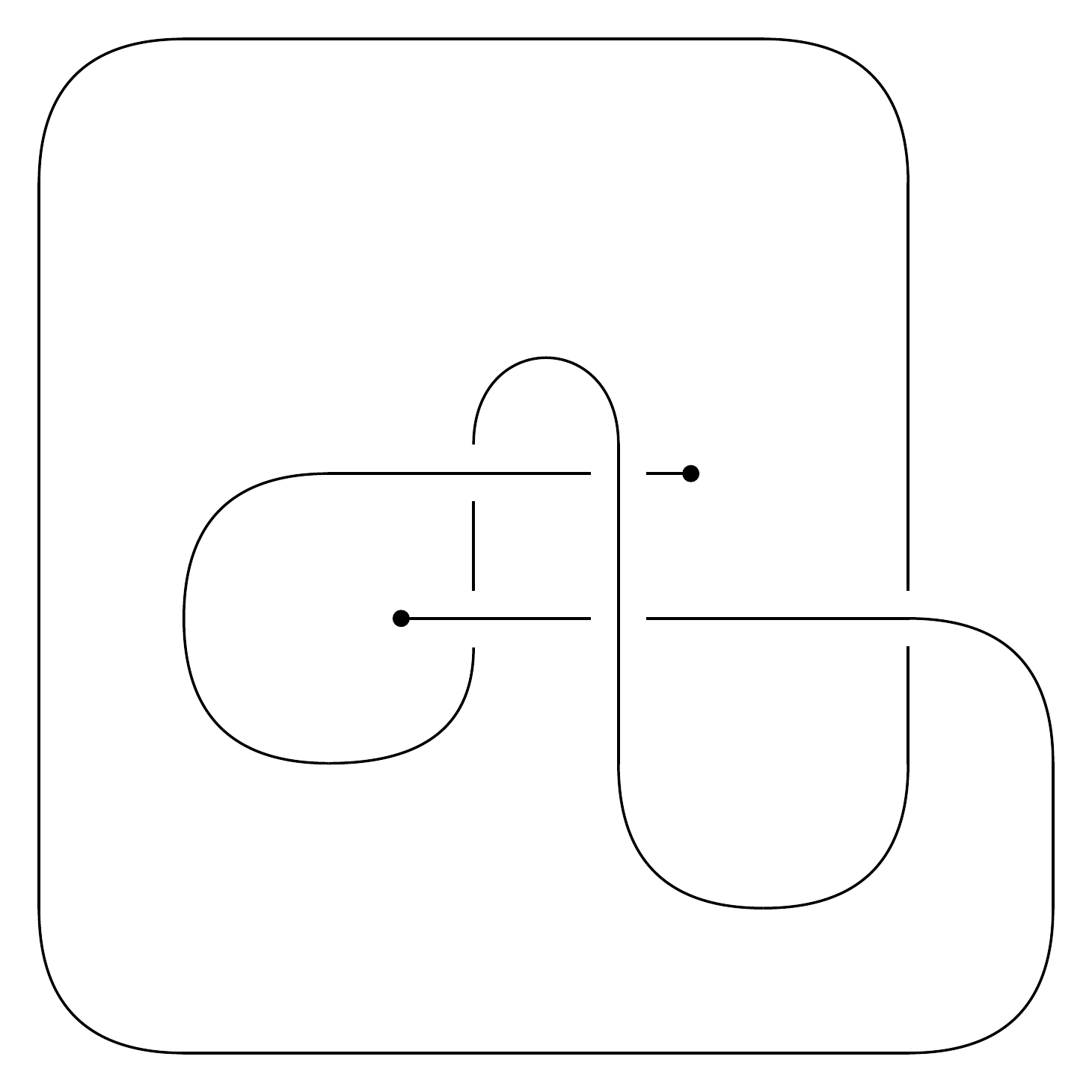}\\
\textcolor{black}{$5_{612}$}
\vspace{1cm}
\end{minipage}
\begin{minipage}[t]{.25\linewidth}
\centering
\includegraphics[width=0.9\textwidth,height=3.5cm,keepaspectratio]{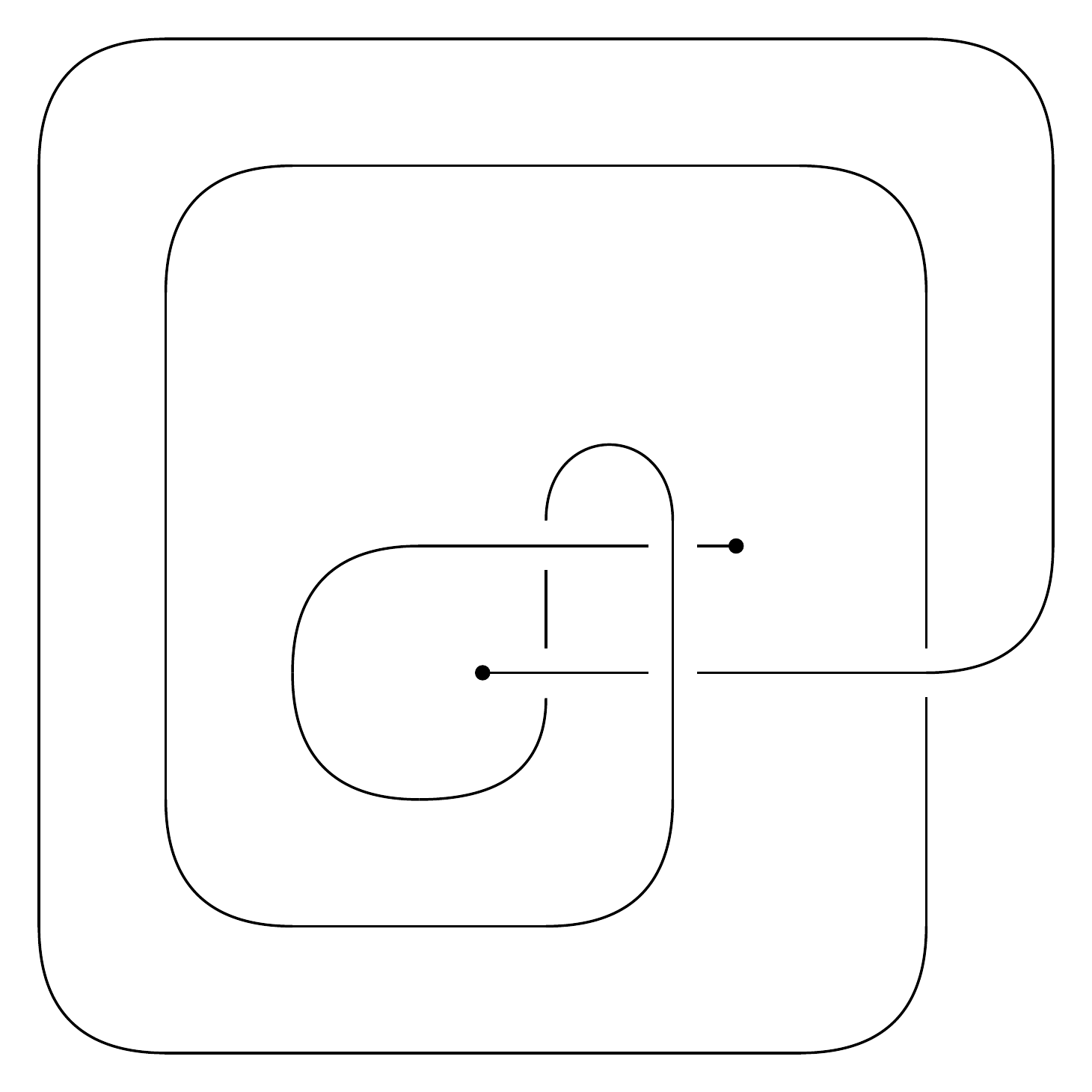}\\
\textcolor{black}{$5_{613}$}
\vspace{1cm}
\end{minipage}
\begin{minipage}[t]{.25\linewidth}
\centering
\includegraphics[width=0.9\textwidth,height=3.5cm,keepaspectratio]{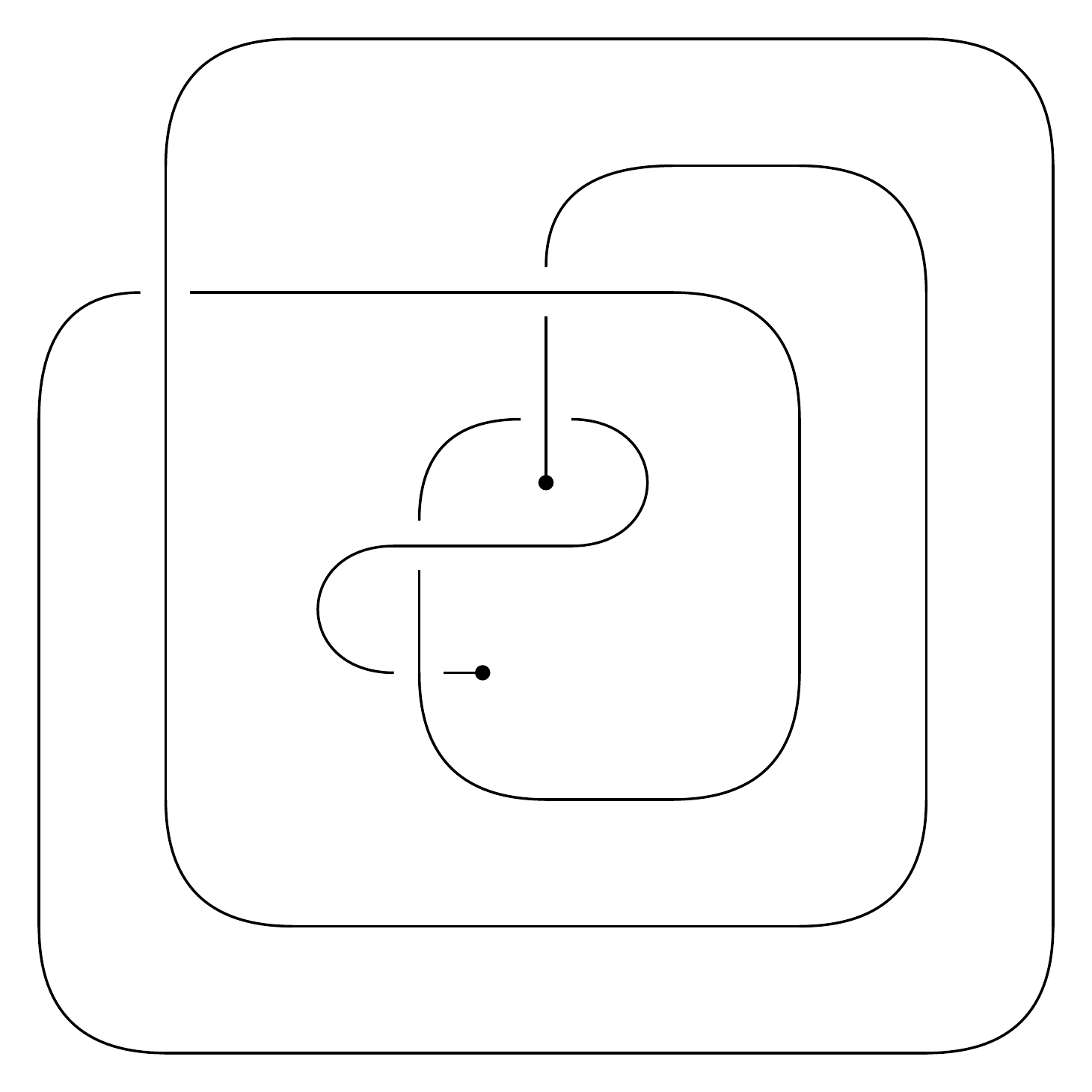}\\
\textcolor{black}{$5_{614}$}
\vspace{1cm}
\end{minipage}
\begin{minipage}[t]{.25\linewidth}
\centering
\includegraphics[width=0.9\textwidth,height=3.5cm,keepaspectratio]{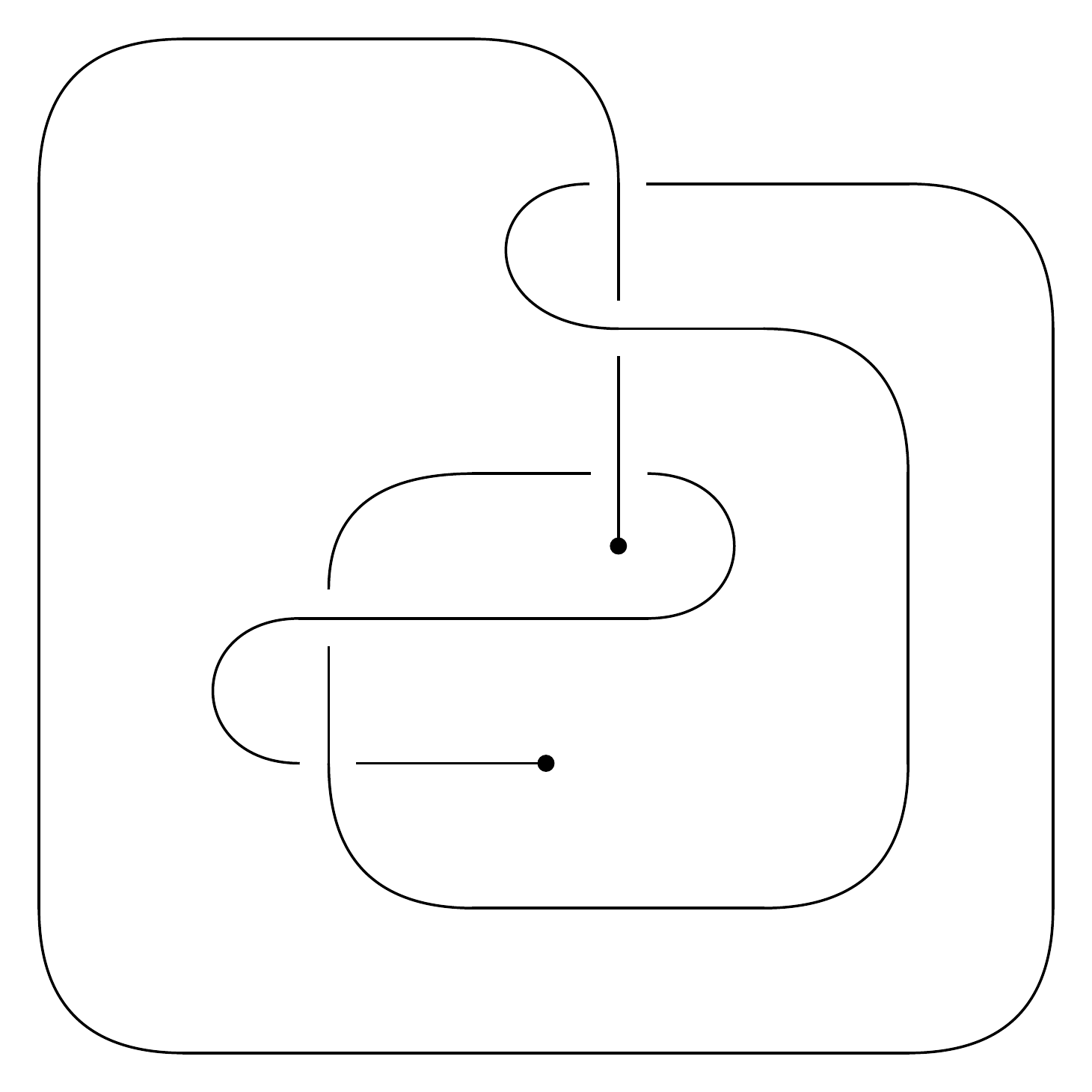}\\
\textcolor{black}{$5_{615}$}
\vspace{1cm}
\end{minipage}
\begin{minipage}[t]{.25\linewidth}
\centering
\includegraphics[width=0.9\textwidth,height=3.5cm,keepaspectratio]{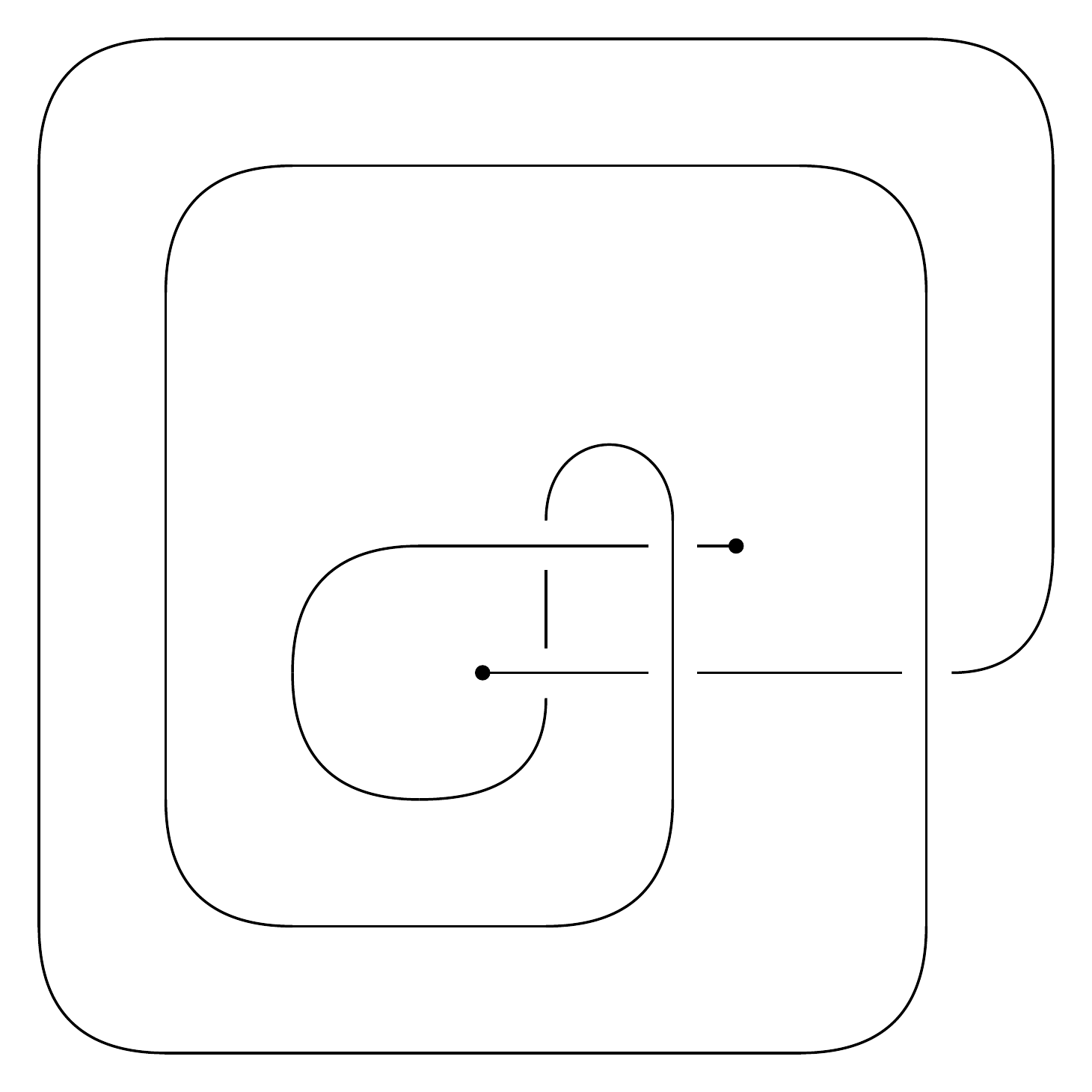}\\
\textcolor{black}{$5_{616}$}
\vspace{1cm}
\end{minipage}
\begin{minipage}[t]{.25\linewidth}
\centering
\includegraphics[width=0.9\textwidth,height=3.5cm,keepaspectratio]{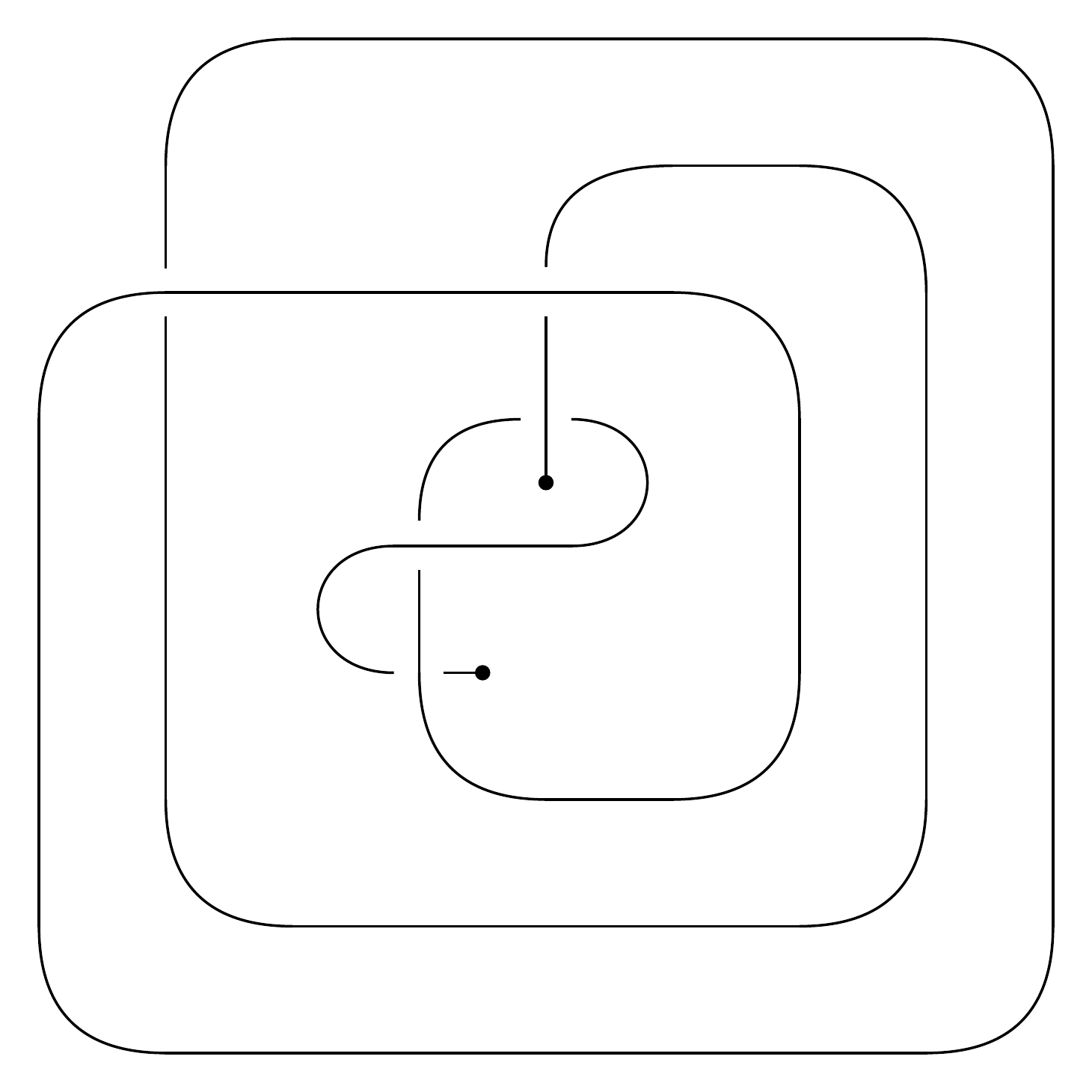}\\
\textcolor{black}{$5_{617}$}
\vspace{1cm}
\end{minipage}
\begin{minipage}[t]{.25\linewidth}
\centering
\includegraphics[width=0.9\textwidth,height=3.5cm,keepaspectratio]{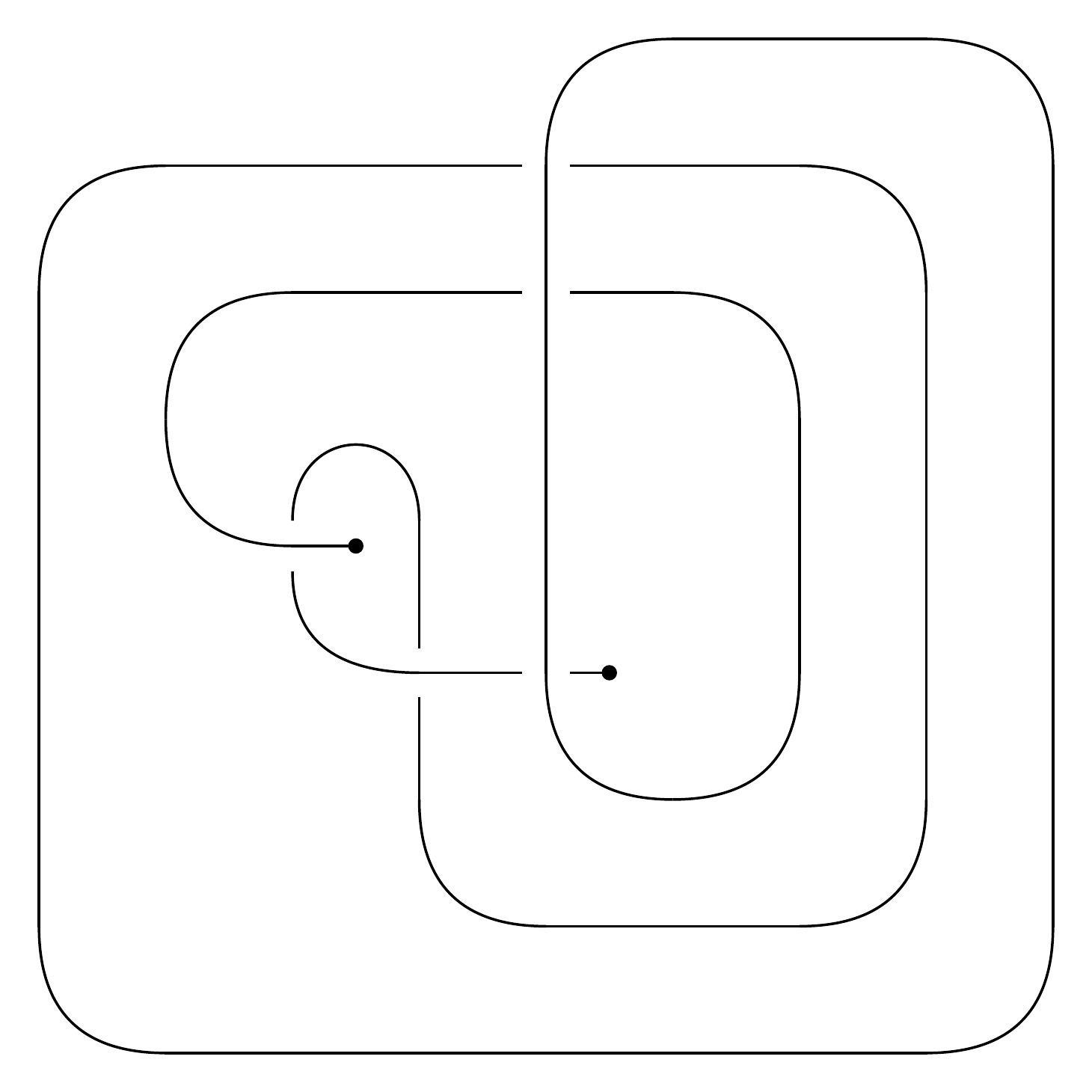}\\
\textcolor{black}{$5_{618}$}
\vspace{1cm}
\end{minipage}
\begin{minipage}[t]{.25\linewidth}
\centering
\includegraphics[width=0.9\textwidth,height=3.5cm,keepaspectratio]{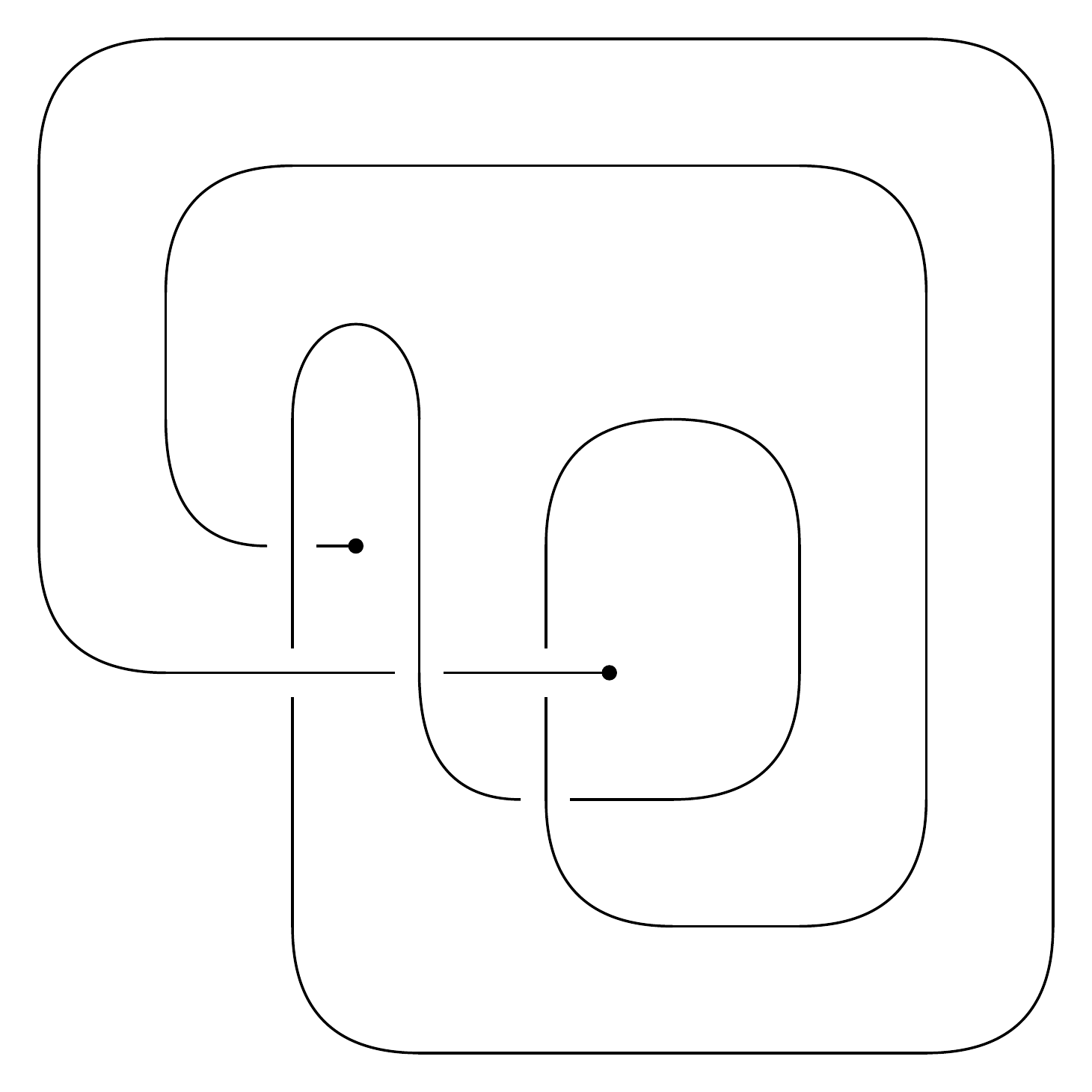}\\
\textcolor{black}{$5_{619}$}
\vspace{1cm}
\end{minipage}
\begin{minipage}[t]{.25\linewidth}
\centering
\includegraphics[width=0.9\textwidth,height=3.5cm,keepaspectratio]{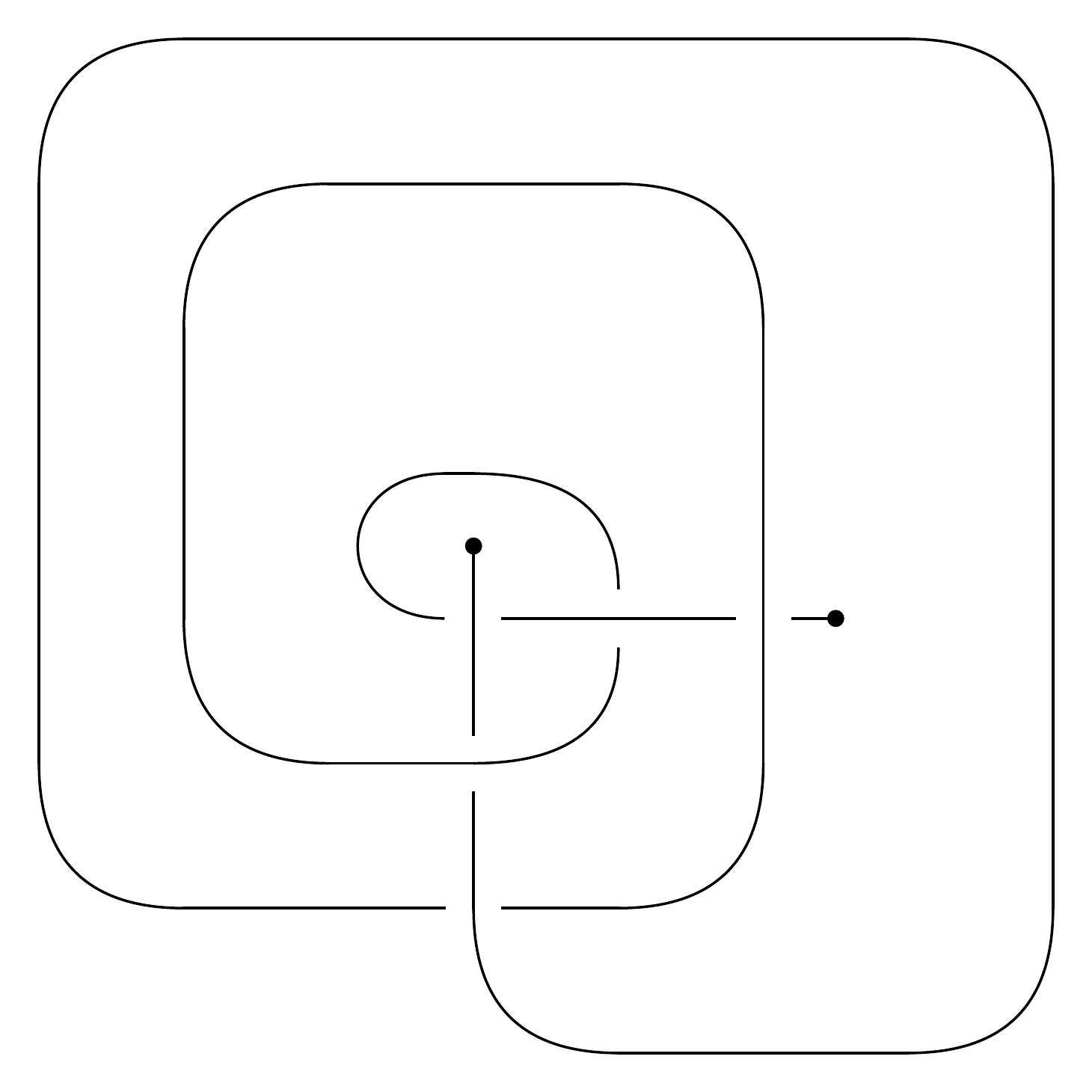}\\
\textcolor{black}{$5_{620}$}
\vspace{1cm}
\end{minipage}
\begin{minipage}[t]{.25\linewidth}
\centering
\includegraphics[width=0.9\textwidth,height=3.5cm,keepaspectratio]{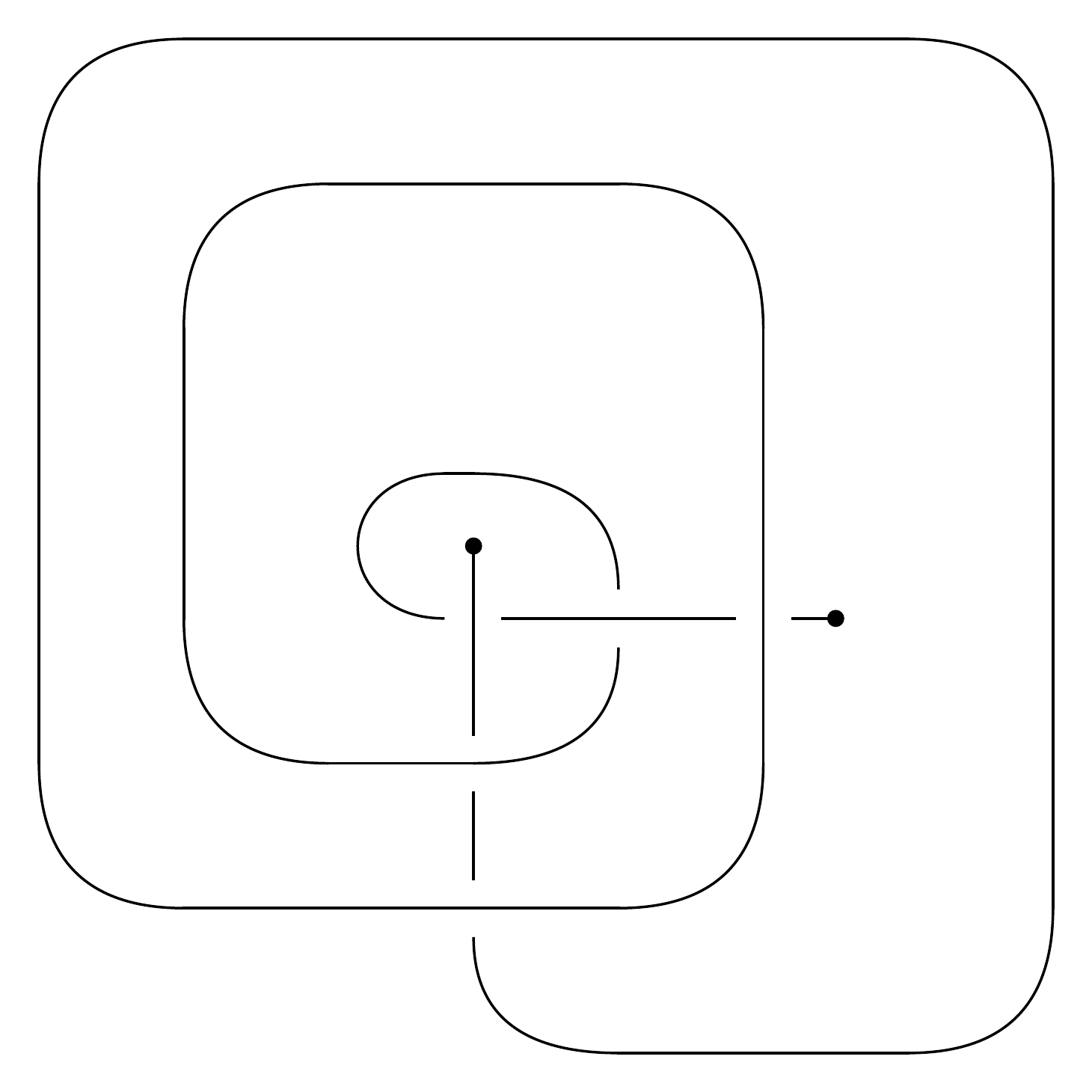}\\
\textcolor{black}{$5_{621}$}
\vspace{1cm}
\end{minipage}
\begin{minipage}[t]{.25\linewidth}
\centering
\includegraphics[width=0.9\textwidth,height=3.5cm,keepaspectratio]{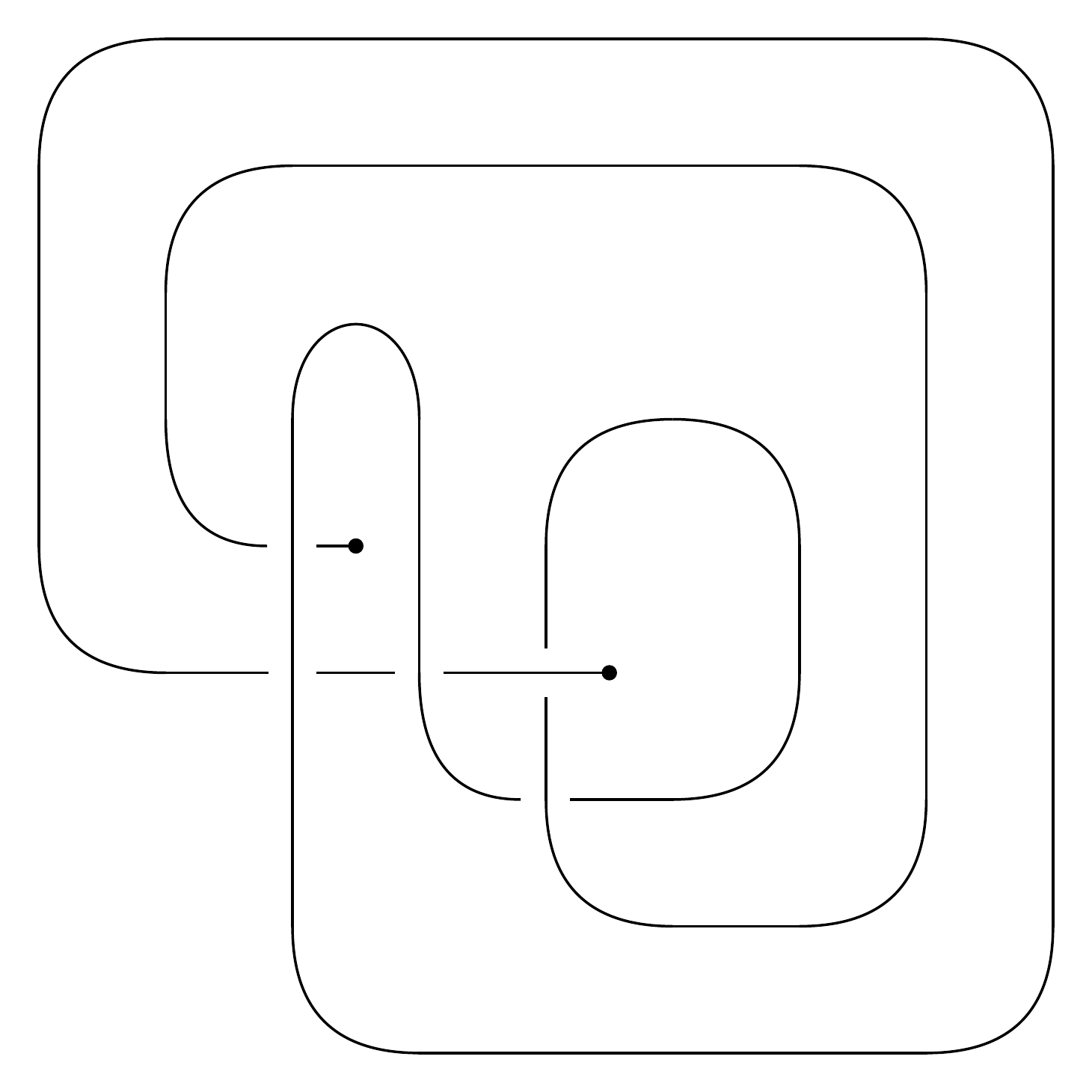}\\
\textcolor{black}{$5_{622}$}
\vspace{1cm}
\end{minipage}
\begin{minipage}[t]{.25\linewidth}
\centering
\includegraphics[width=0.9\textwidth,height=3.5cm,keepaspectratio]{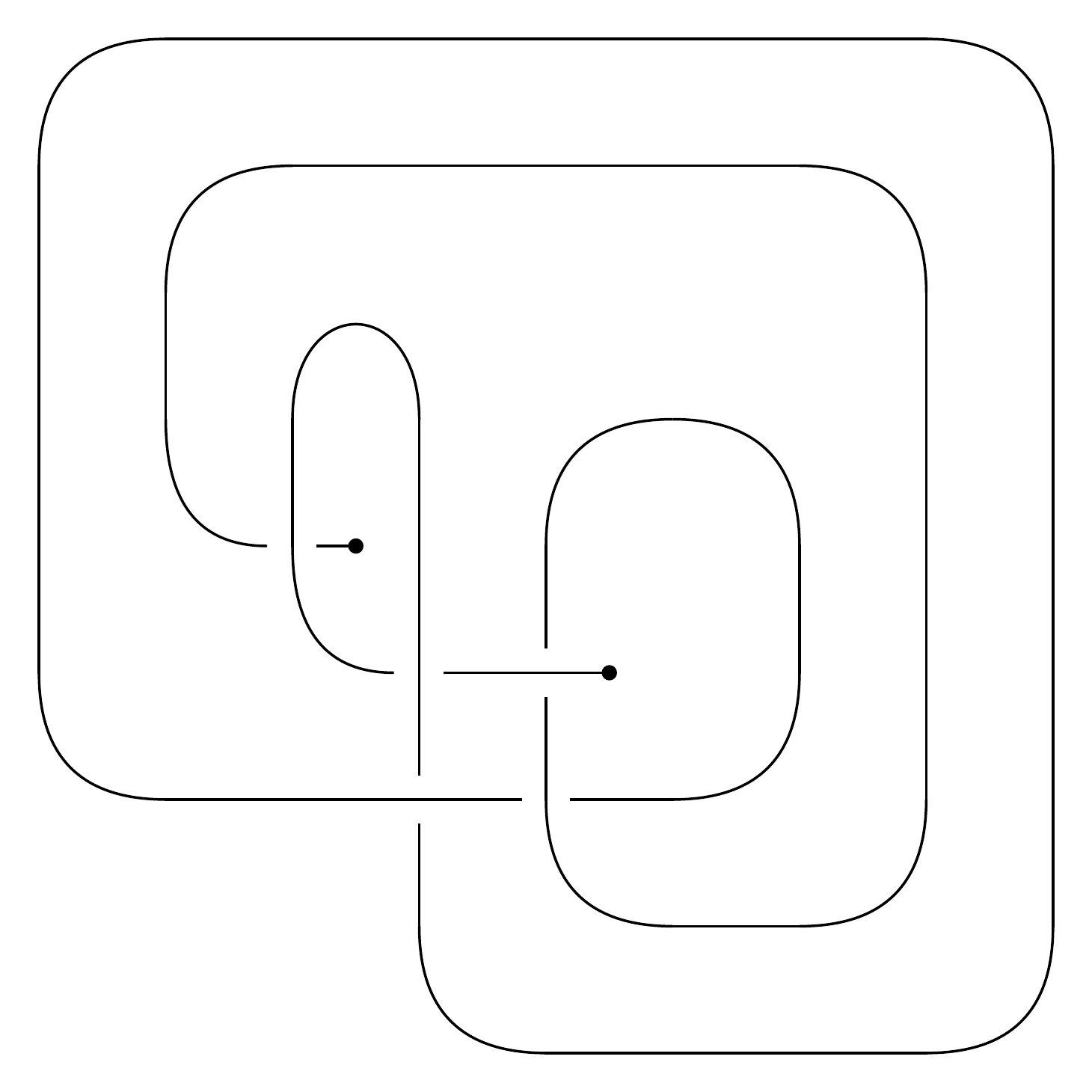}\\
\textcolor{black}{$5_{623}$}
\vspace{1cm}
\end{minipage}
\begin{minipage}[t]{.25\linewidth}
\centering
\includegraphics[width=0.9\textwidth,height=3.5cm,keepaspectratio]{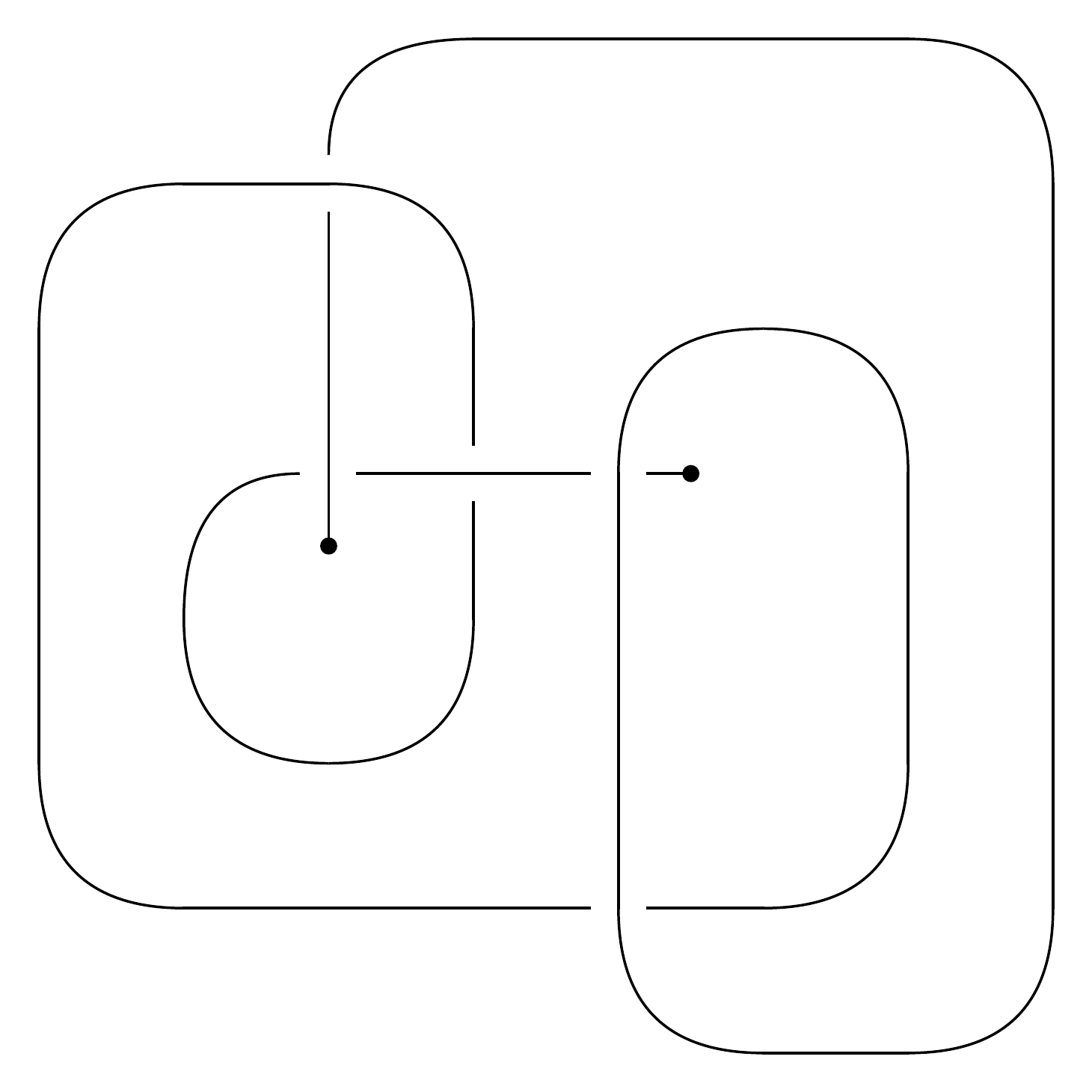}\\
\textcolor{black}{$5_{624}$}
\vspace{1cm}
\end{minipage}
\begin{minipage}[t]{.25\linewidth}
\centering
\includegraphics[width=0.9\textwidth,height=3.5cm,keepaspectratio]{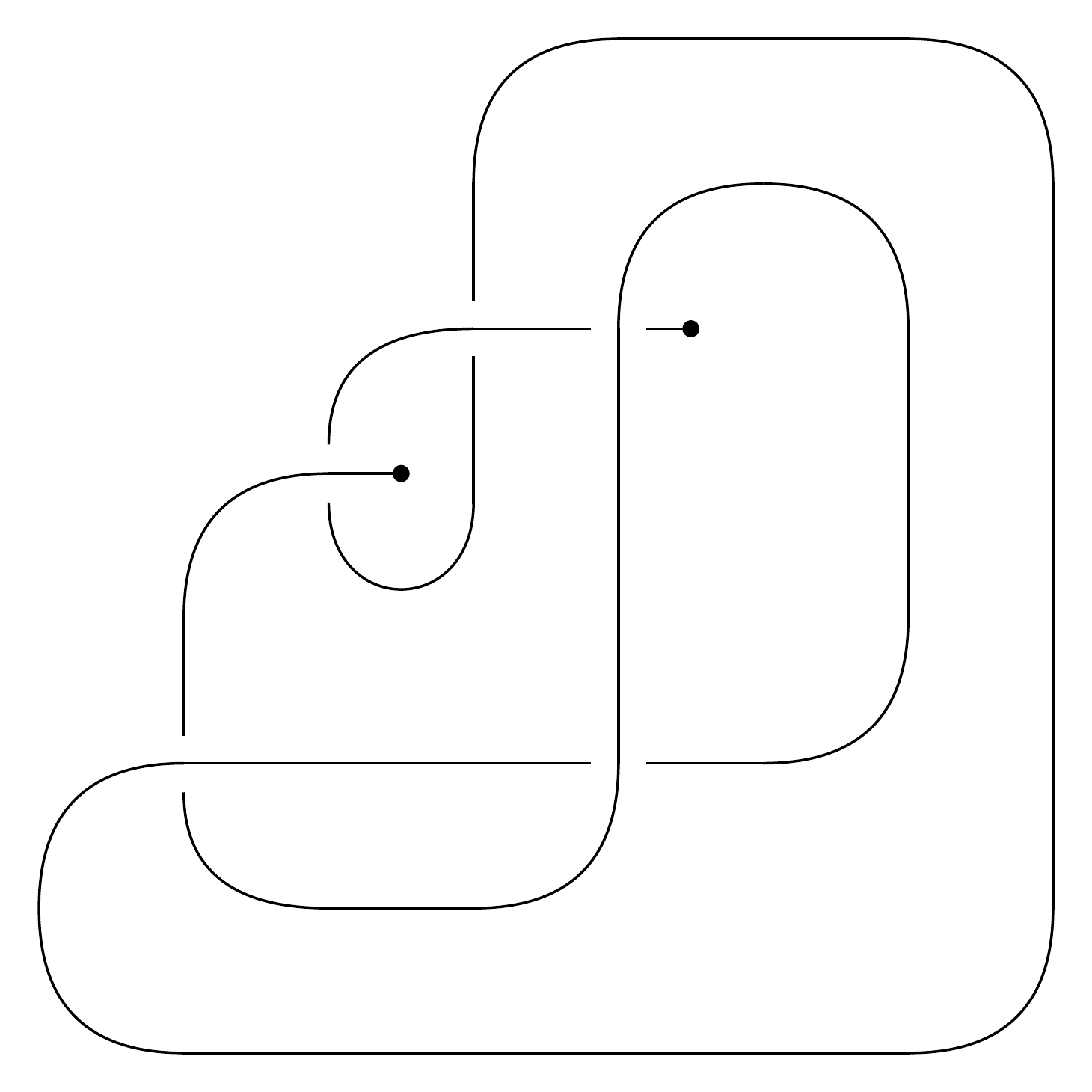}\\
\textcolor{black}{$5_{625}$}
\vspace{1cm}
\end{minipage}
\begin{minipage}[t]{.25\linewidth}
\centering
\includegraphics[width=0.9\textwidth,height=3.5cm,keepaspectratio]{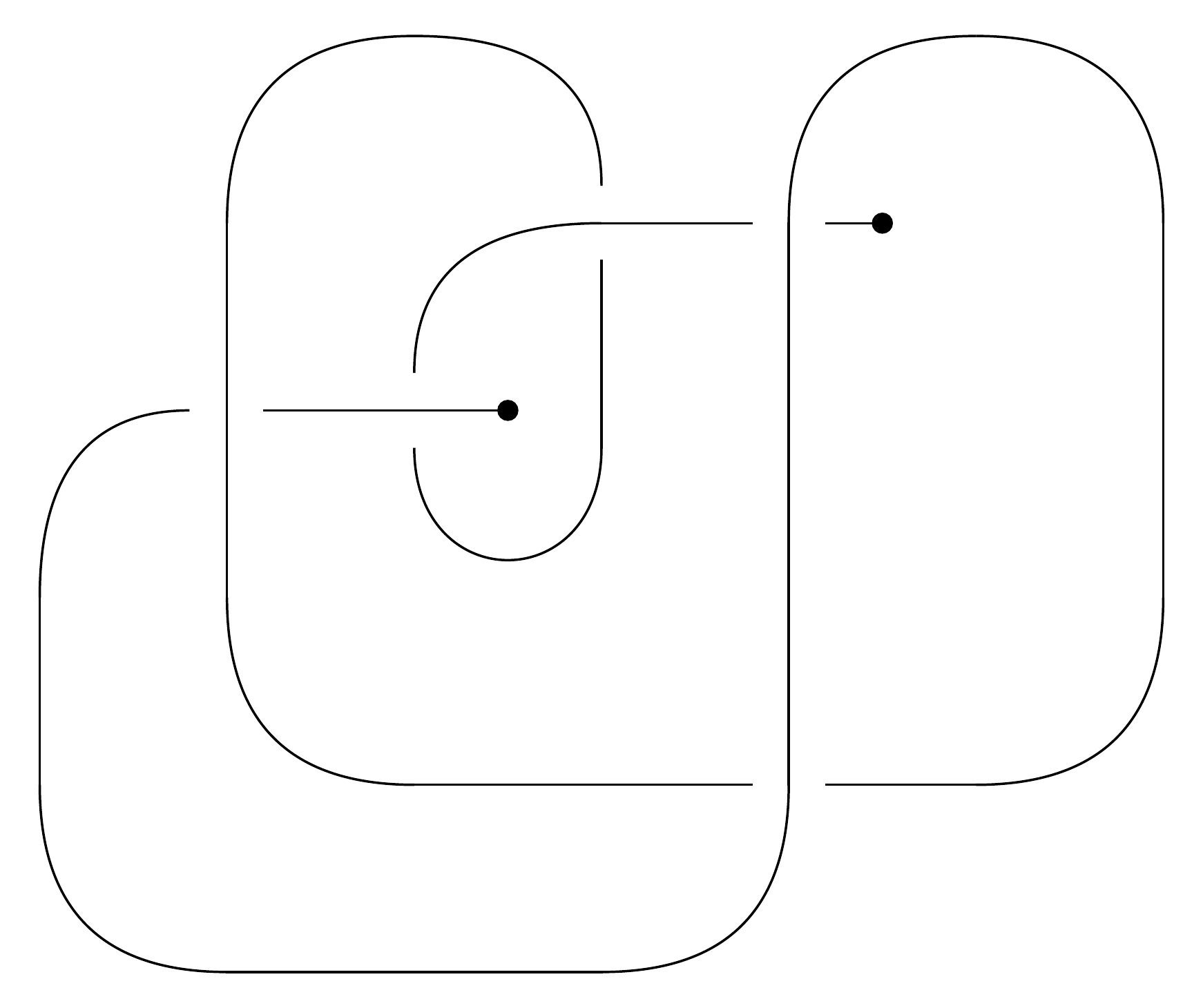}\\
\textcolor{black}{$5_{626}$}
\vspace{1cm}
\end{minipage}
\begin{minipage}[t]{.25\linewidth}
\centering
\includegraphics[width=0.9\textwidth,height=3.5cm,keepaspectratio]{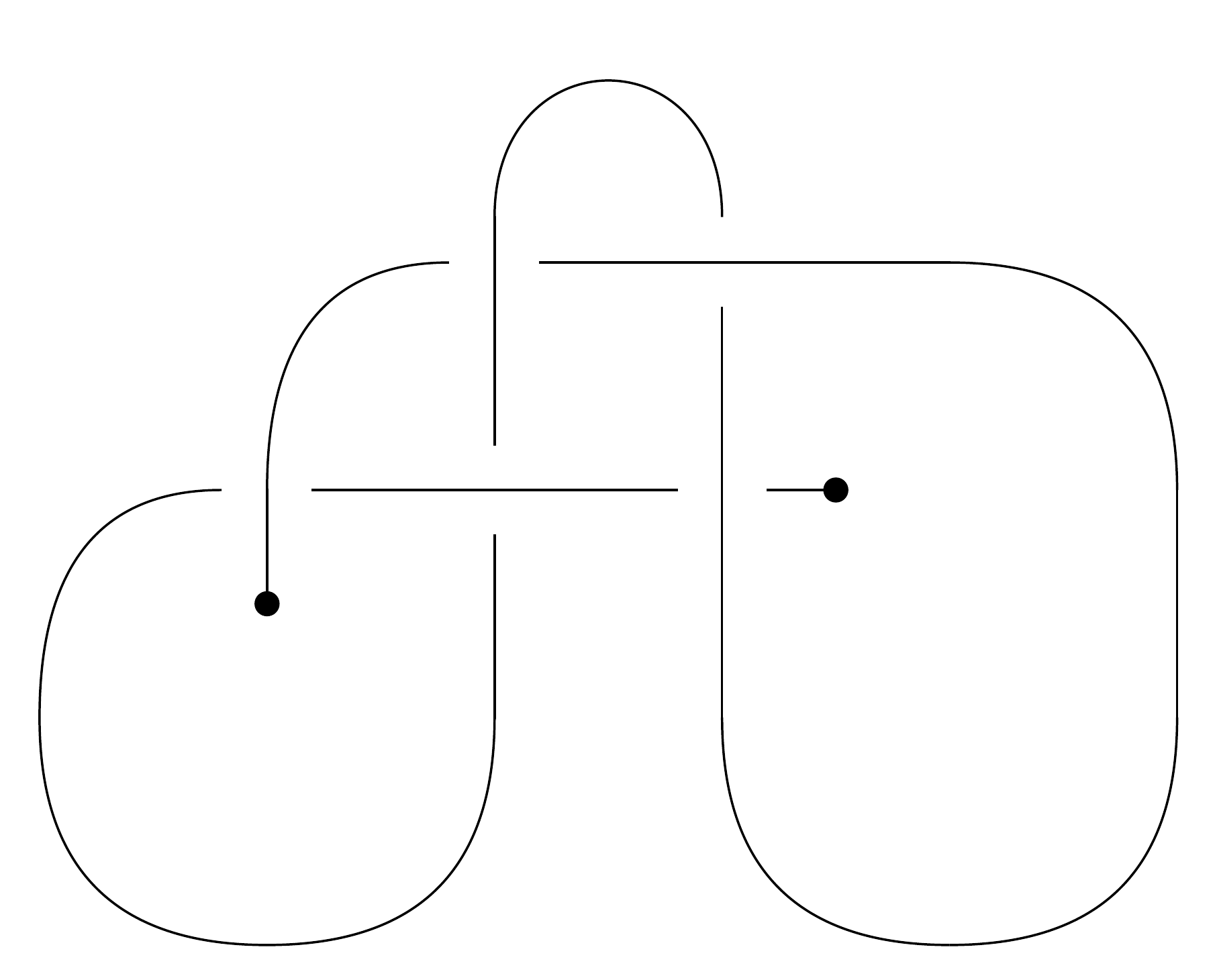}\\
\textcolor{black}{$5_{627}$}
\vspace{1cm}
\end{minipage}
\begin{minipage}[t]{.25\linewidth}
\centering
\includegraphics[width=0.9\textwidth,height=3.5cm,keepaspectratio]{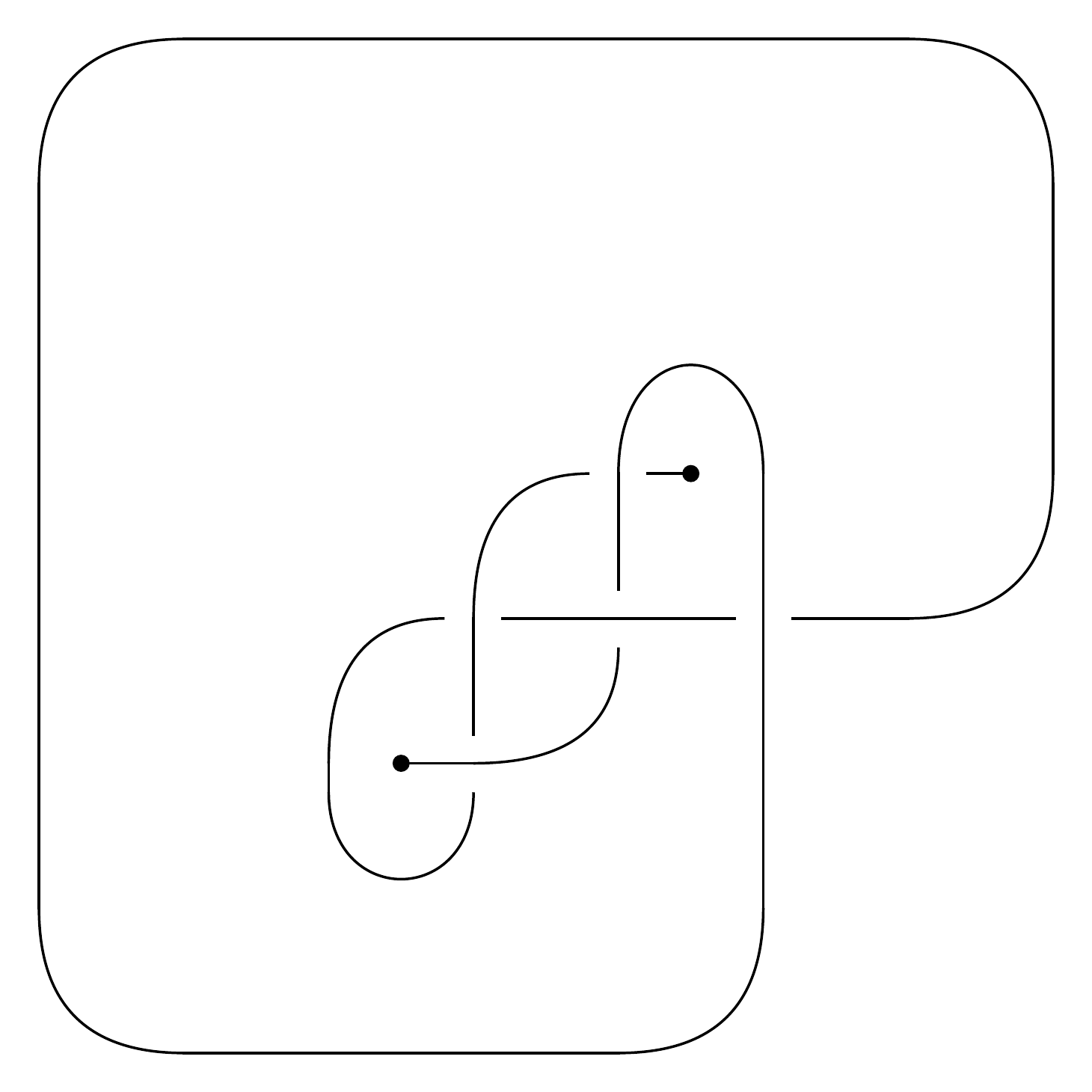}\\
\textcolor{black}{$5_{628}$}
\vspace{1cm}
\end{minipage}
\begin{minipage}[t]{.25\linewidth}
\centering
\includegraphics[width=0.9\textwidth,height=3.5cm,keepaspectratio]{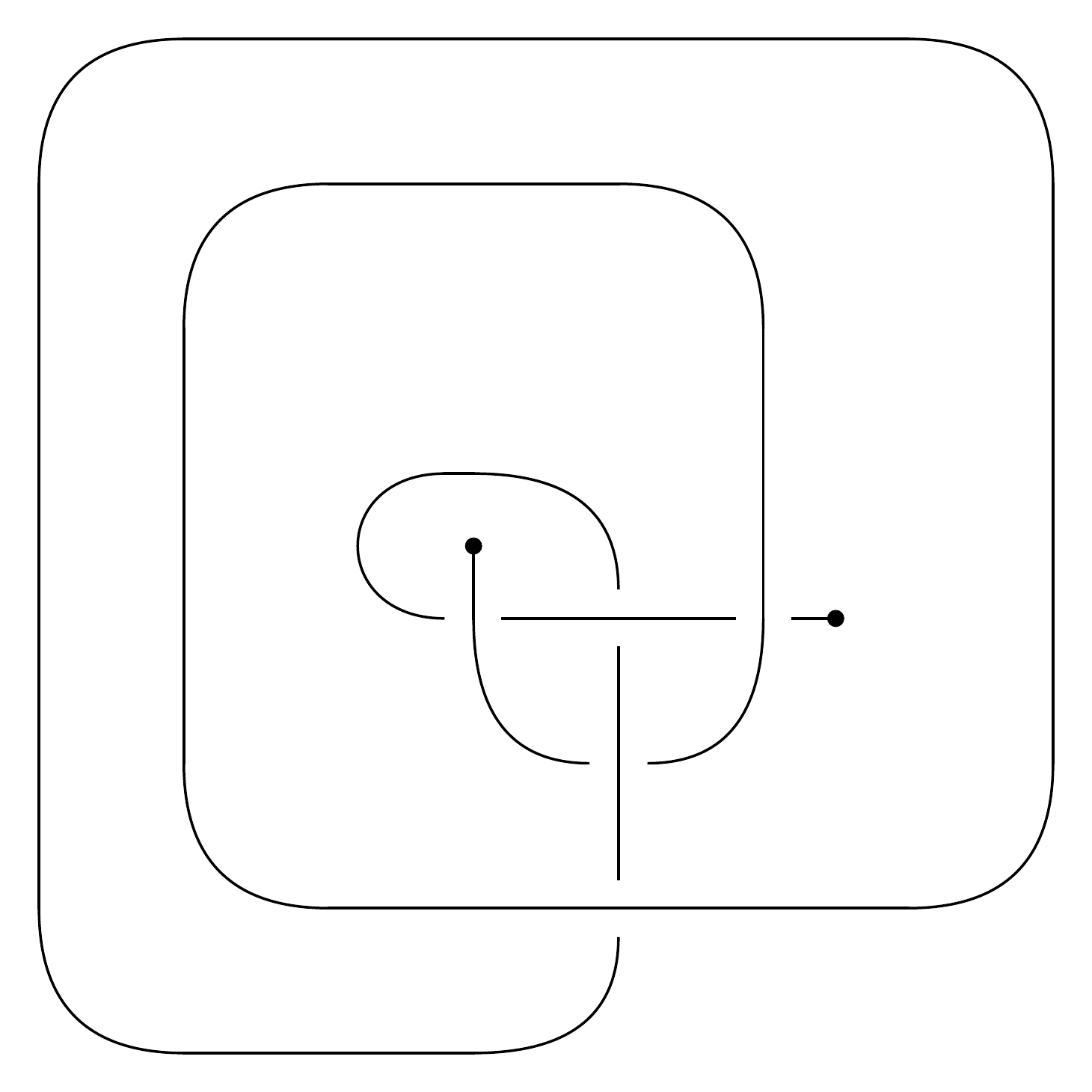}\\
\textcolor{black}{$5_{629}$}
\vspace{1cm}
\end{minipage}
\begin{minipage}[t]{.25\linewidth}
\centering
\includegraphics[width=0.9\textwidth,height=3.5cm,keepaspectratio]{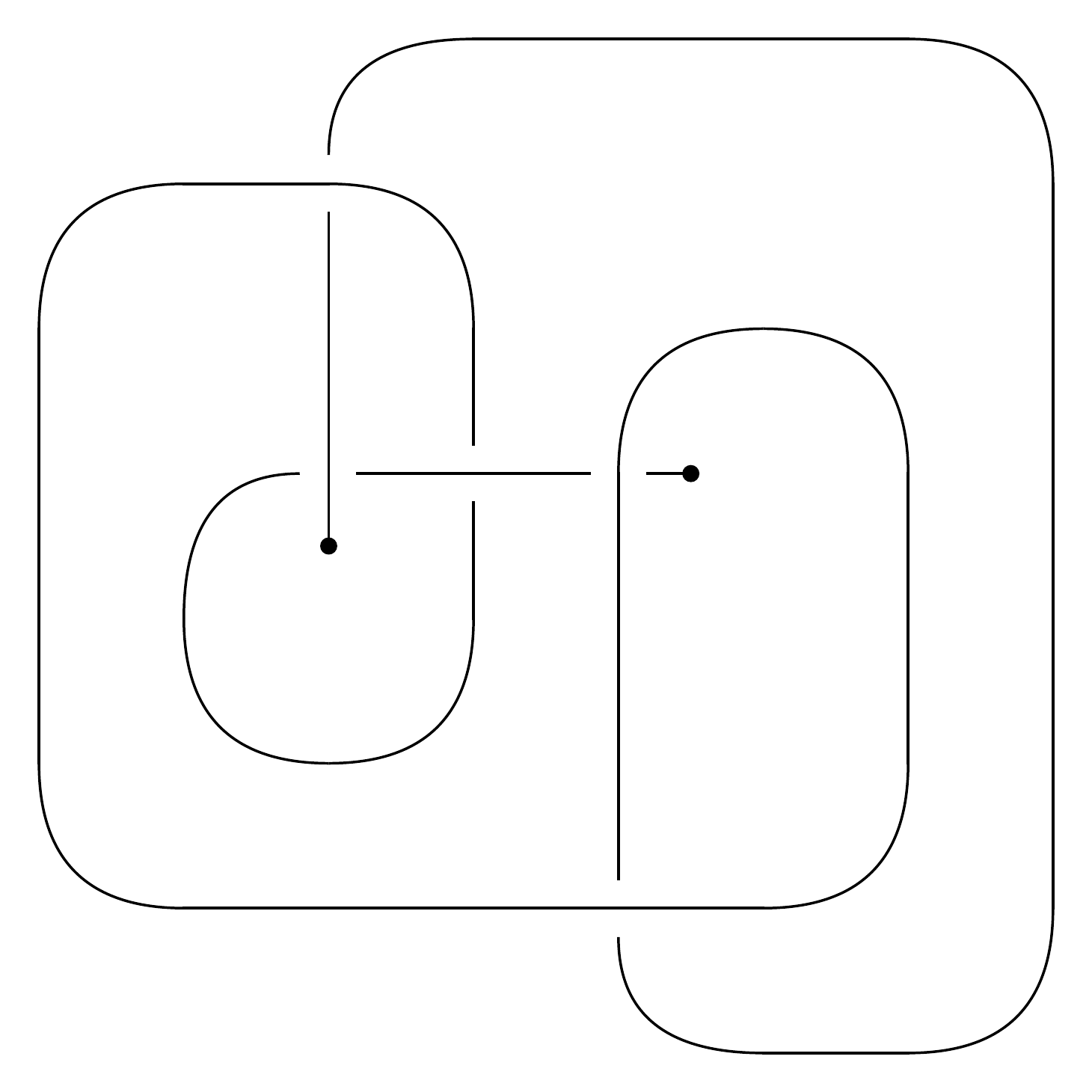}\\
\textcolor{black}{$5_{630}$}
\vspace{1cm}
\end{minipage}
\begin{minipage}[t]{.25\linewidth}
\centering
\includegraphics[width=0.9\textwidth,height=3.5cm,keepaspectratio]{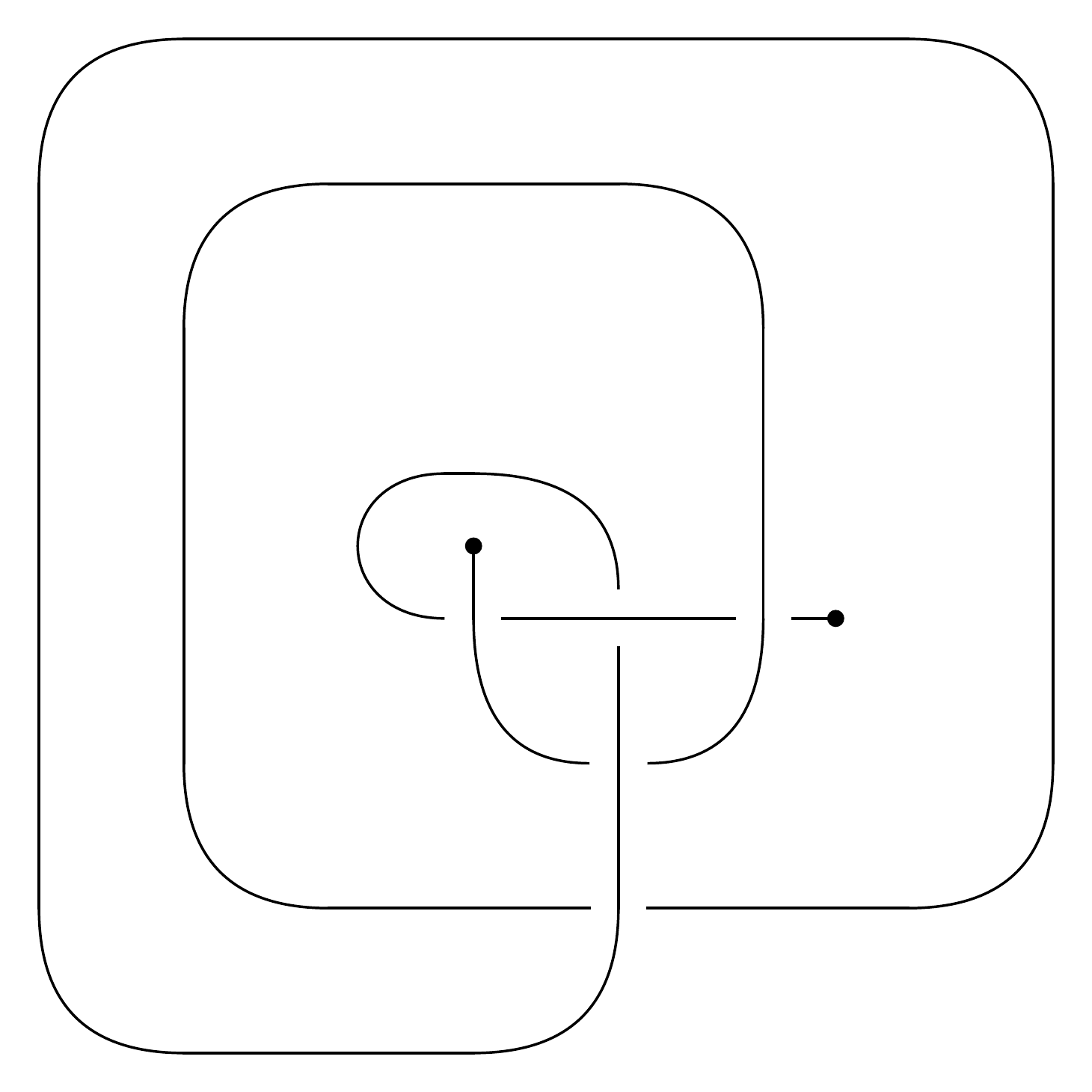}\\
\textcolor{black}{$5_{631}$}
\vspace{1cm}
\end{minipage}
\begin{minipage}[t]{.25\linewidth}
\centering
\includegraphics[width=0.9\textwidth,height=3.5cm,keepaspectratio]{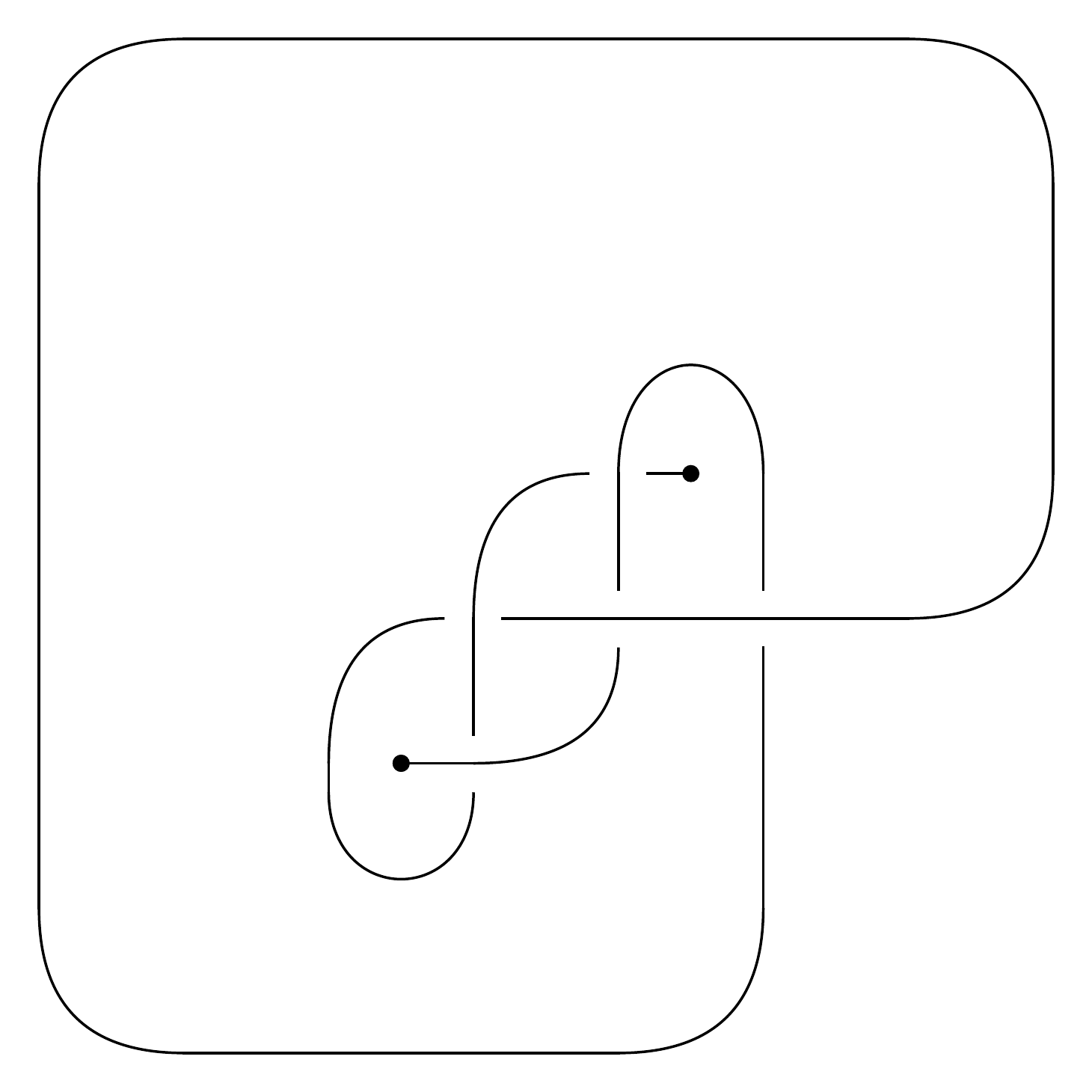}\\
\textcolor{black}{$5_{632}$}
\vspace{1cm}
\end{minipage}
\begin{minipage}[t]{.25\linewidth}
\centering
\includegraphics[width=0.9\textwidth,height=3.5cm,keepaspectratio]{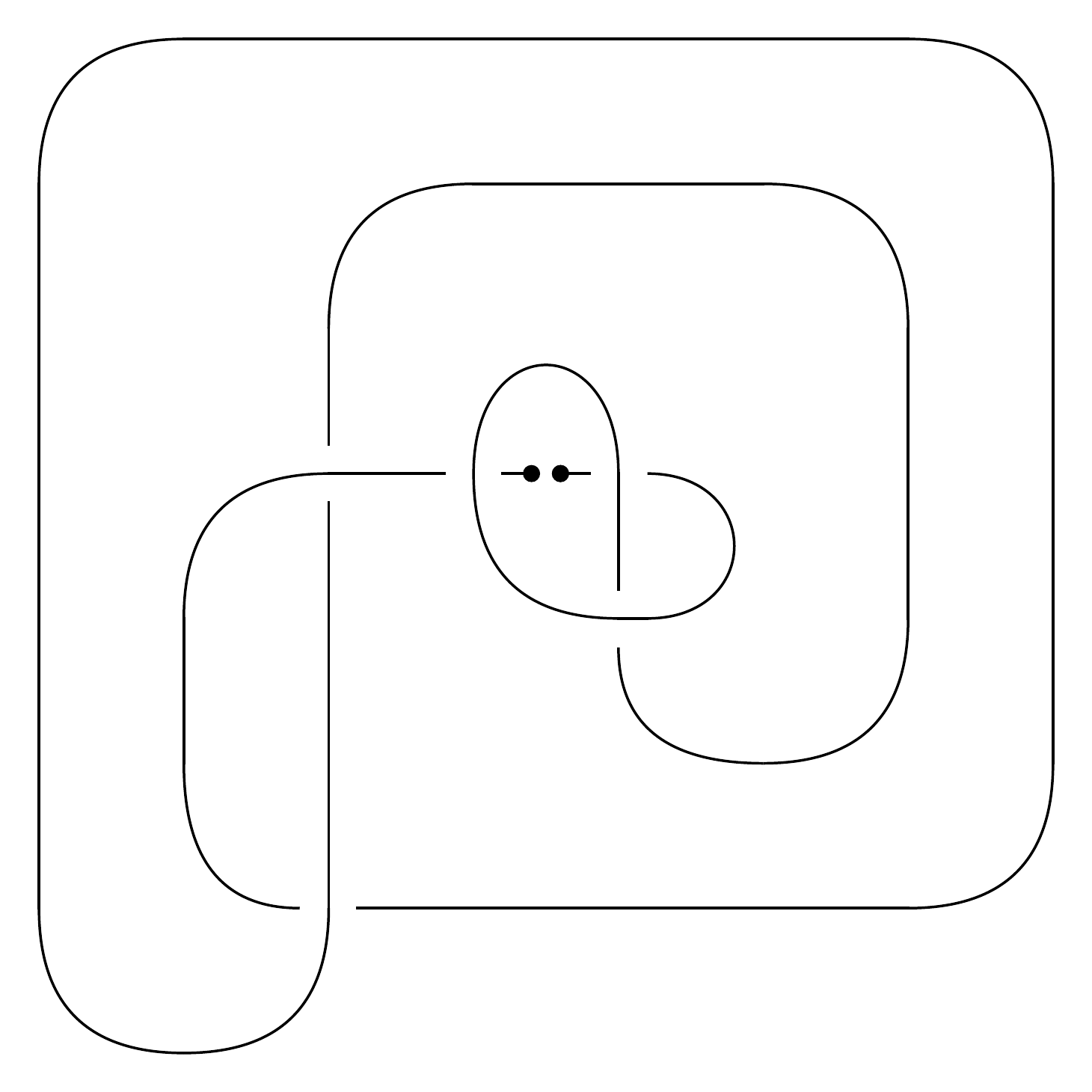}\\
\textcolor{black}{$5_{633}$}
\vspace{1cm}
\end{minipage}
\begin{minipage}[t]{.25\linewidth}
\centering
\includegraphics[width=0.9\textwidth,height=3.5cm,keepaspectratio]{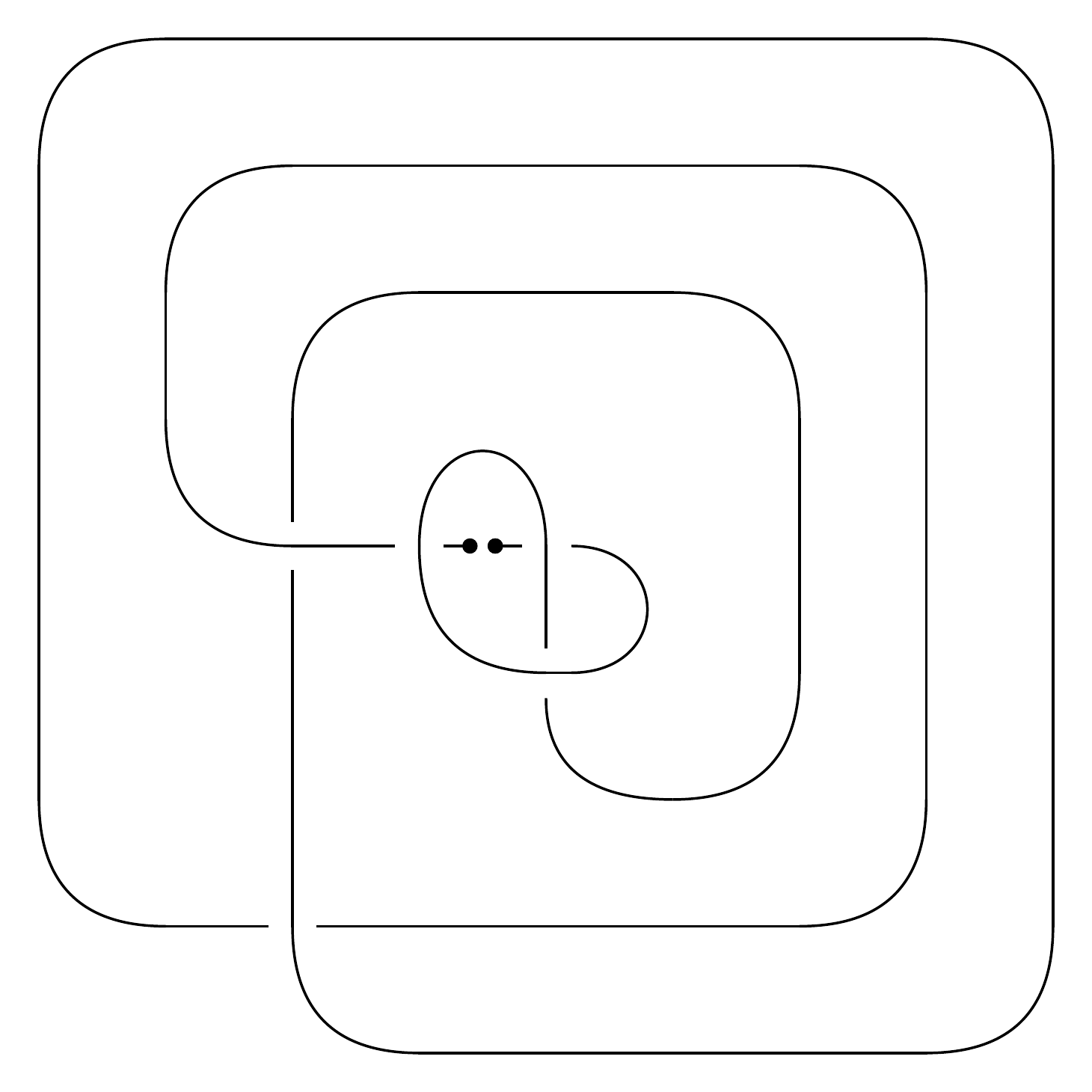}\\
\textcolor{black}{$5_{634}$}
\vspace{1cm}
\end{minipage}
\begin{minipage}[t]{.25\linewidth}
\centering
\includegraphics[width=0.9\textwidth,height=3.5cm,keepaspectratio]{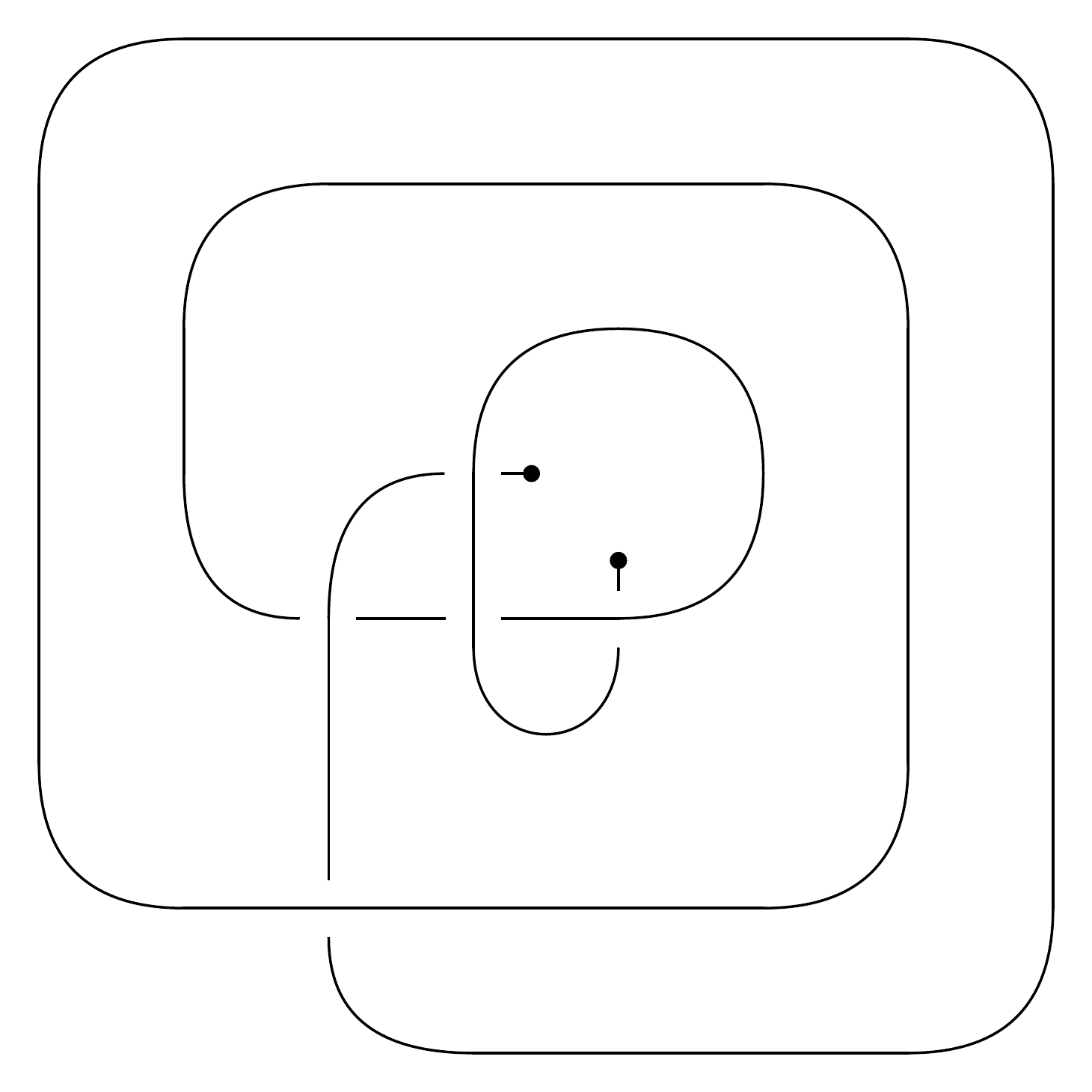}\\
\textcolor{black}{$5_{635}$}
\vspace{1cm}
\end{minipage}
\begin{minipage}[t]{.25\linewidth}
\centering
\includegraphics[width=0.9\textwidth,height=3.5cm,keepaspectratio]{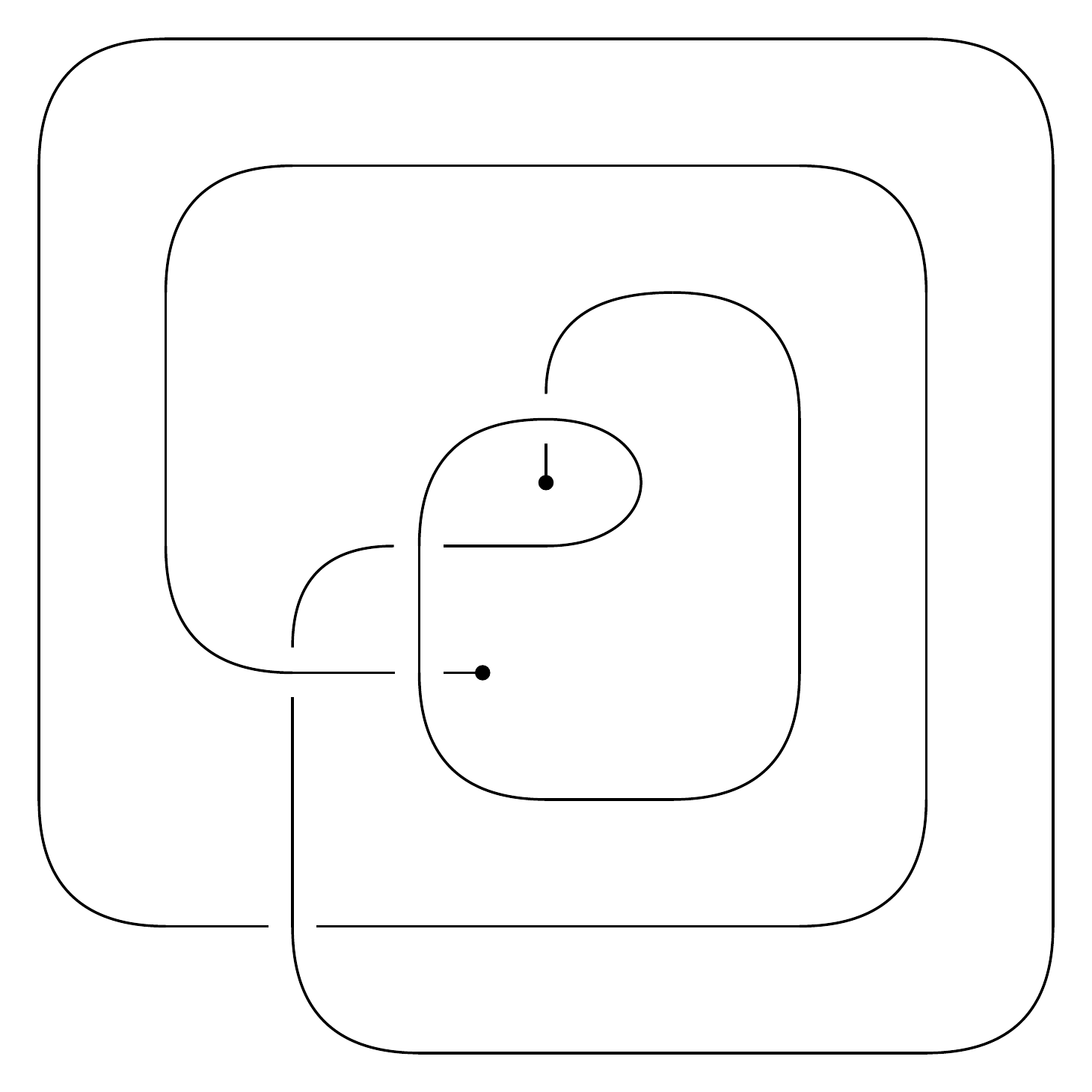}\\
\textcolor{black}{$5_{636}$}
\vspace{1cm}
\end{minipage}
\begin{minipage}[t]{.25\linewidth}
\centering
\includegraphics[width=0.9\textwidth,height=3.5cm,keepaspectratio]{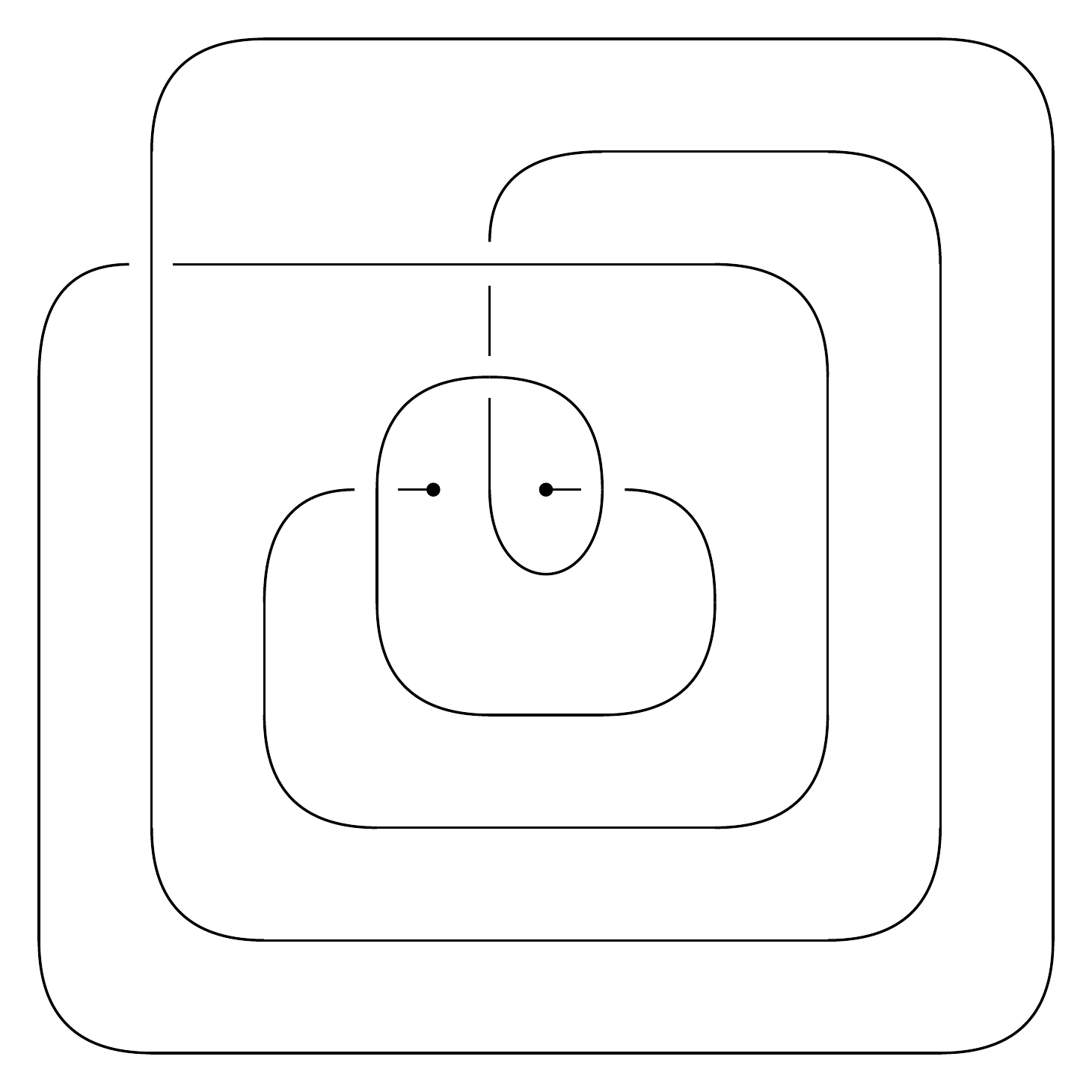}\\
\textcolor{black}{$5_{637}$}
\vspace{1cm}
\end{minipage}
\begin{minipage}[t]{.25\linewidth}
\centering
\includegraphics[width=0.9\textwidth,height=3.5cm,keepaspectratio]{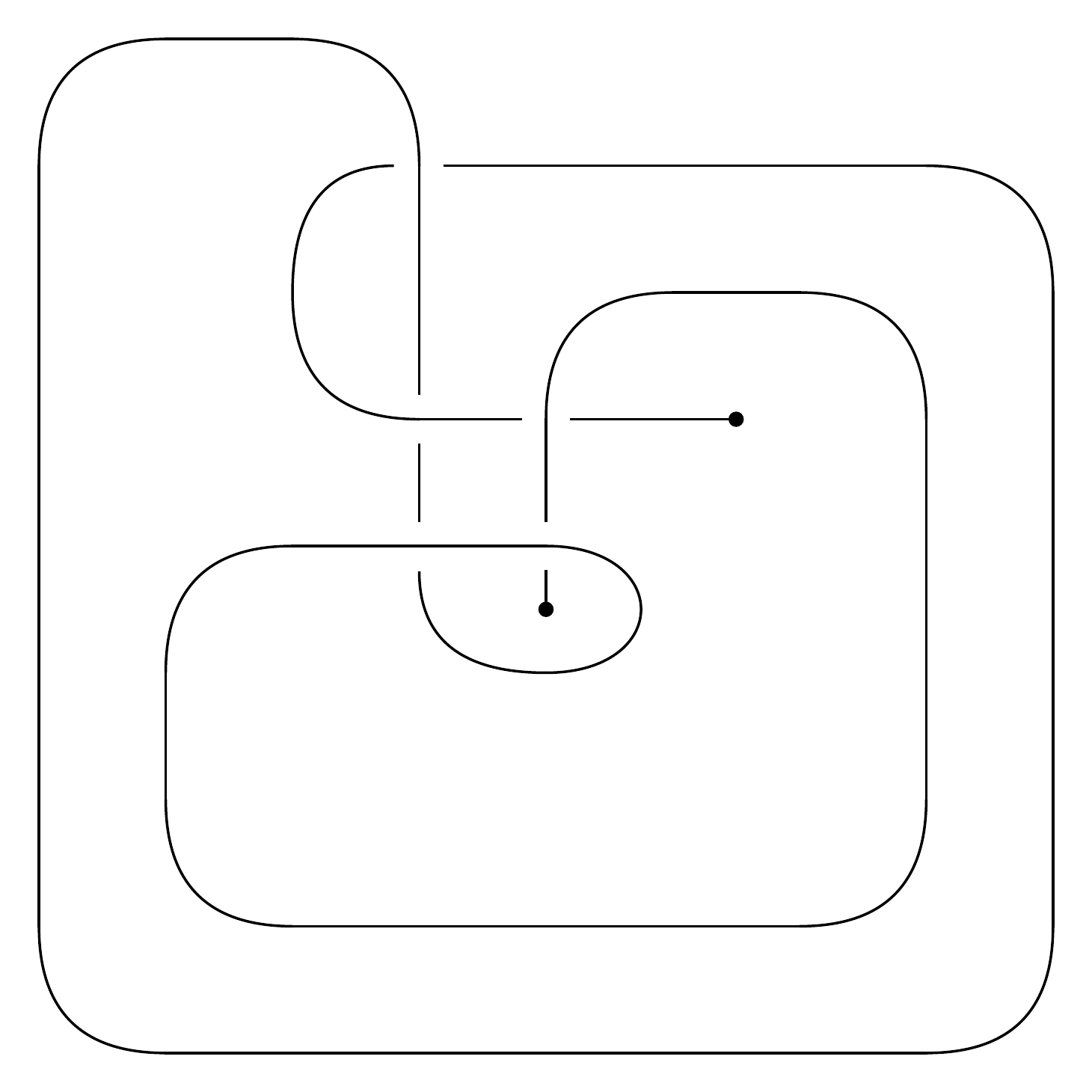}\\
\textcolor{black}{$5_{638}$}
\vspace{1cm}
\end{minipage}
\begin{minipage}[t]{.25\linewidth}
\centering
\includegraphics[width=0.9\textwidth,height=3.5cm,keepaspectratio]{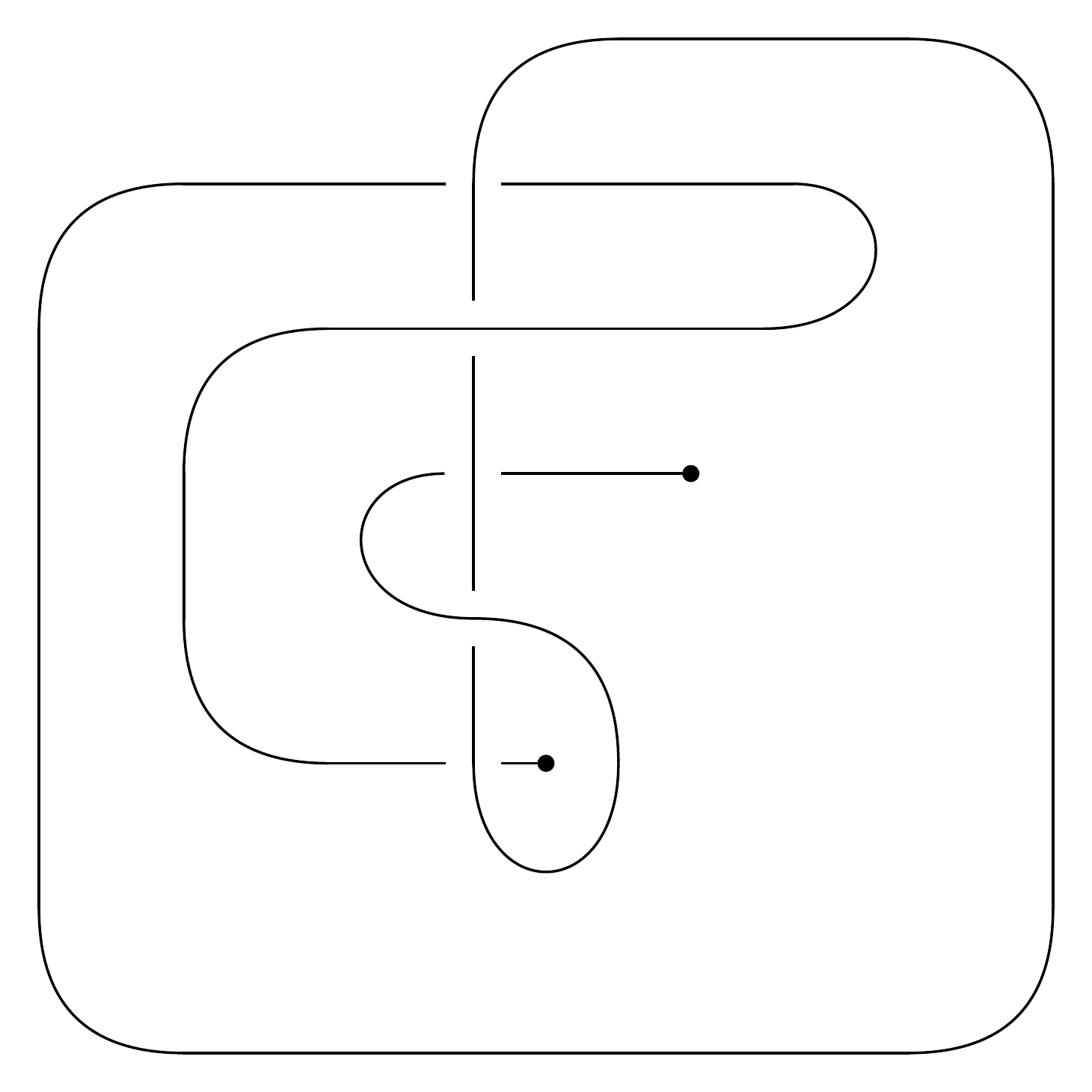}\\
\textcolor{black}{$5_{639}$}
\vspace{1cm}
\end{minipage}
\begin{minipage}[t]{.25\linewidth}
\centering
\includegraphics[width=0.9\textwidth,height=3.5cm,keepaspectratio]{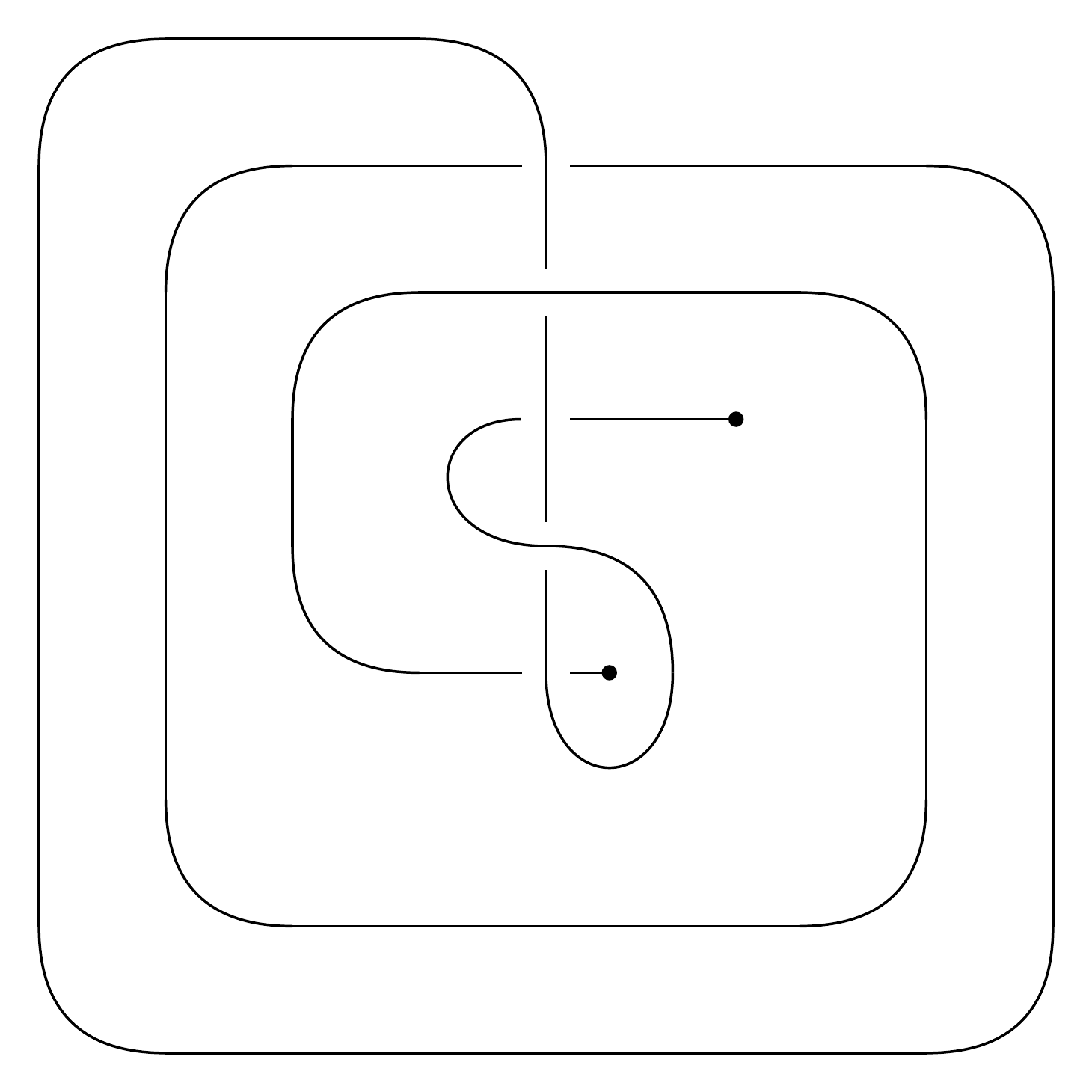}\\
\textcolor{black}{$5_{640}$}
\vspace{1cm}
\end{minipage}
\begin{minipage}[t]{.25\linewidth}
\centering
\includegraphics[width=0.9\textwidth,height=3.5cm,keepaspectratio]{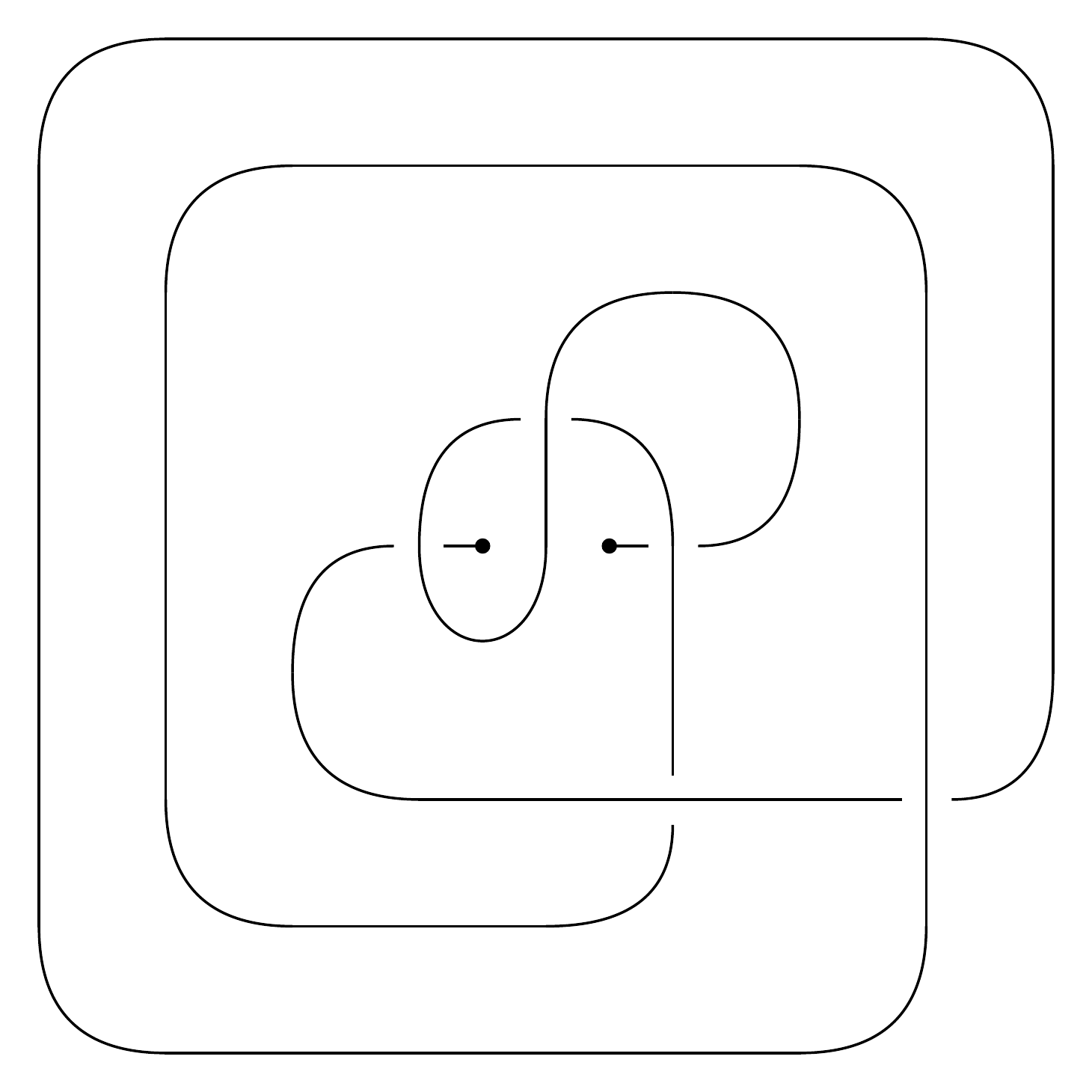}\\
\textcolor{black}{$5_{641}$}
\vspace{1cm}
\end{minipage}
\begin{minipage}[t]{.25\linewidth}
\centering
\includegraphics[width=0.9\textwidth,height=3.5cm,keepaspectratio]{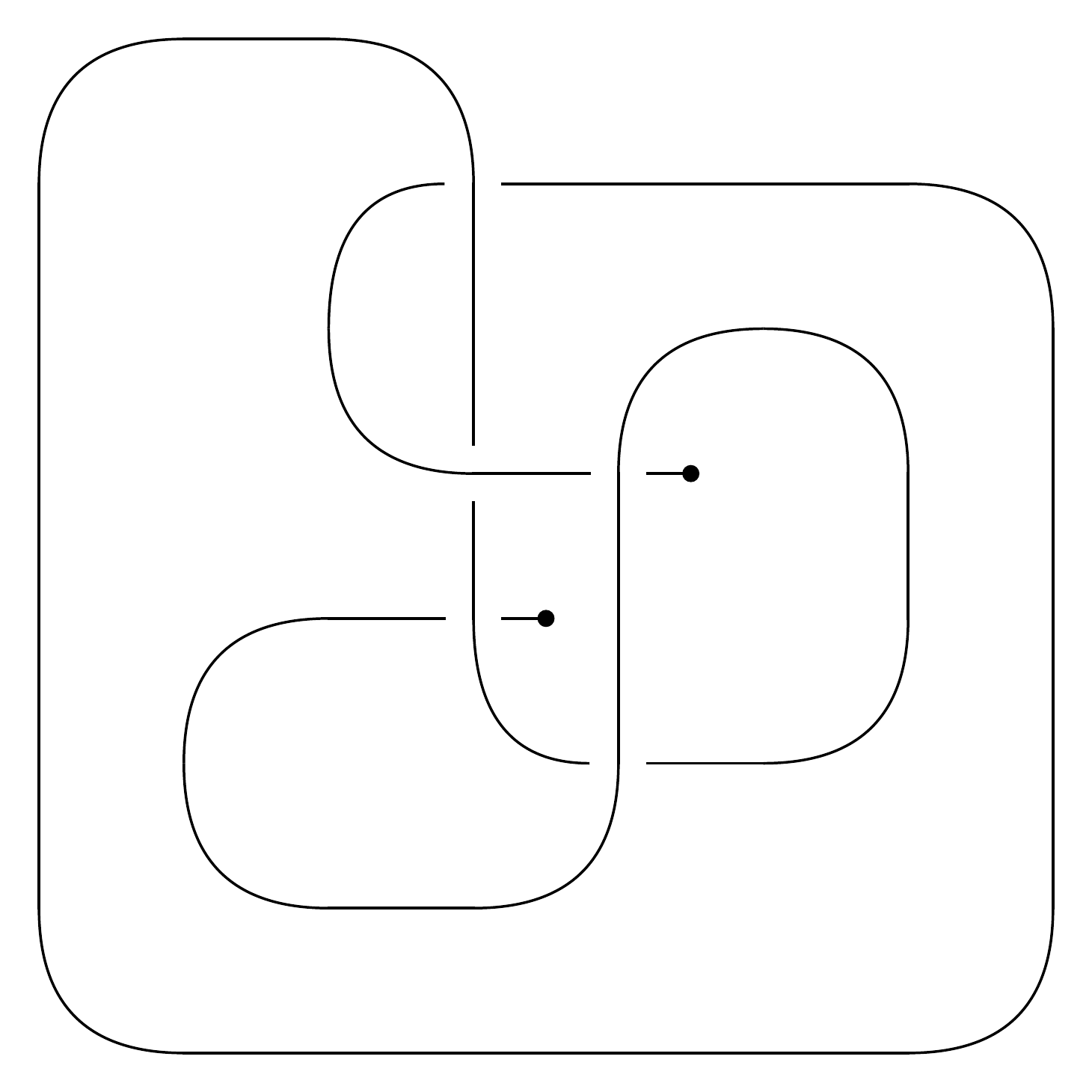}\\
\textcolor{black}{$5_{642}$}
\vspace{1cm}
\end{minipage}
\begin{minipage}[t]{.25\linewidth}
\centering
\includegraphics[width=0.9\textwidth,height=3.5cm,keepaspectratio]{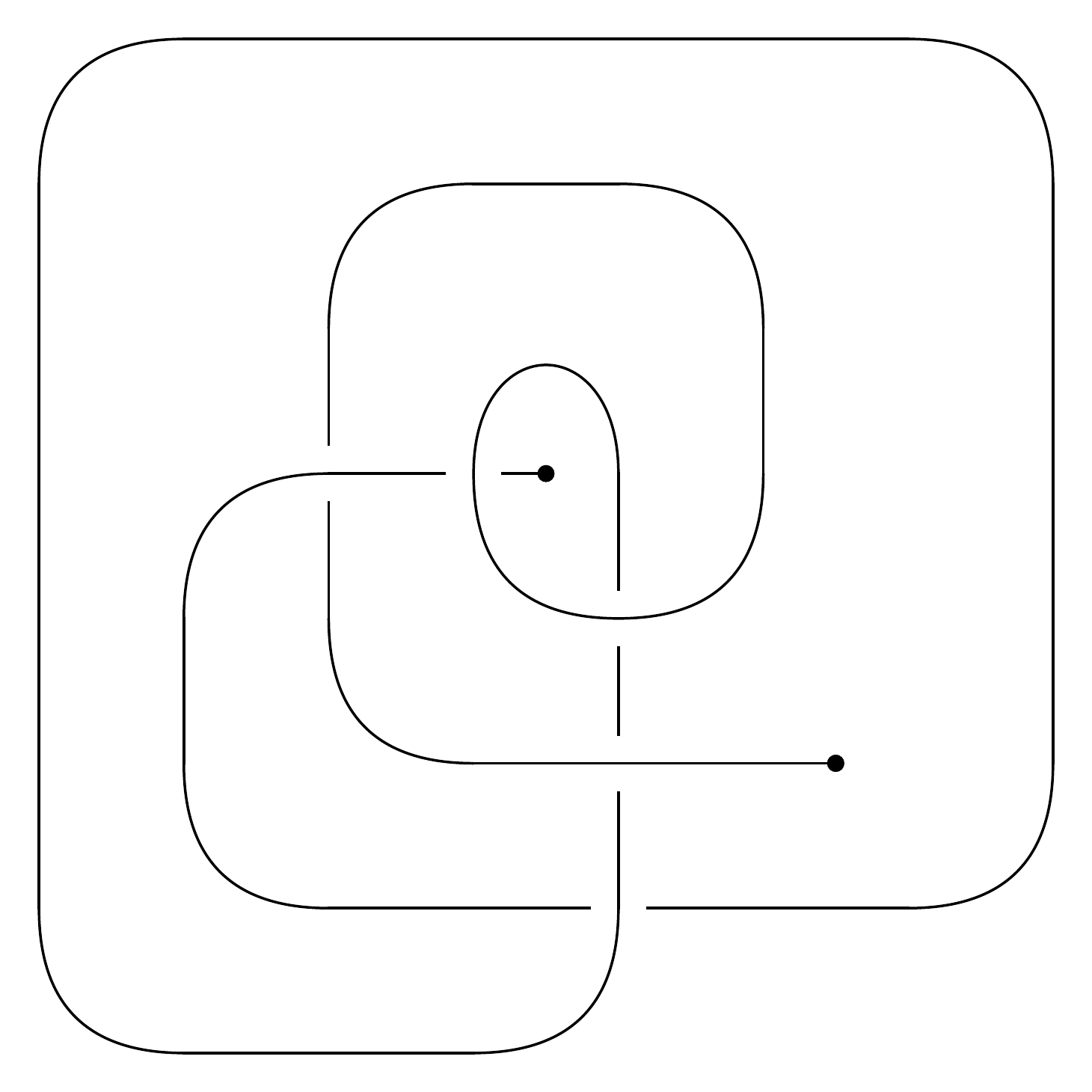}\\
\textcolor{black}{$5_{643}$}
\vspace{1cm}
\end{minipage}
\begin{minipage}[t]{.25\linewidth}
\centering
\includegraphics[width=0.9\textwidth,height=3.5cm,keepaspectratio]{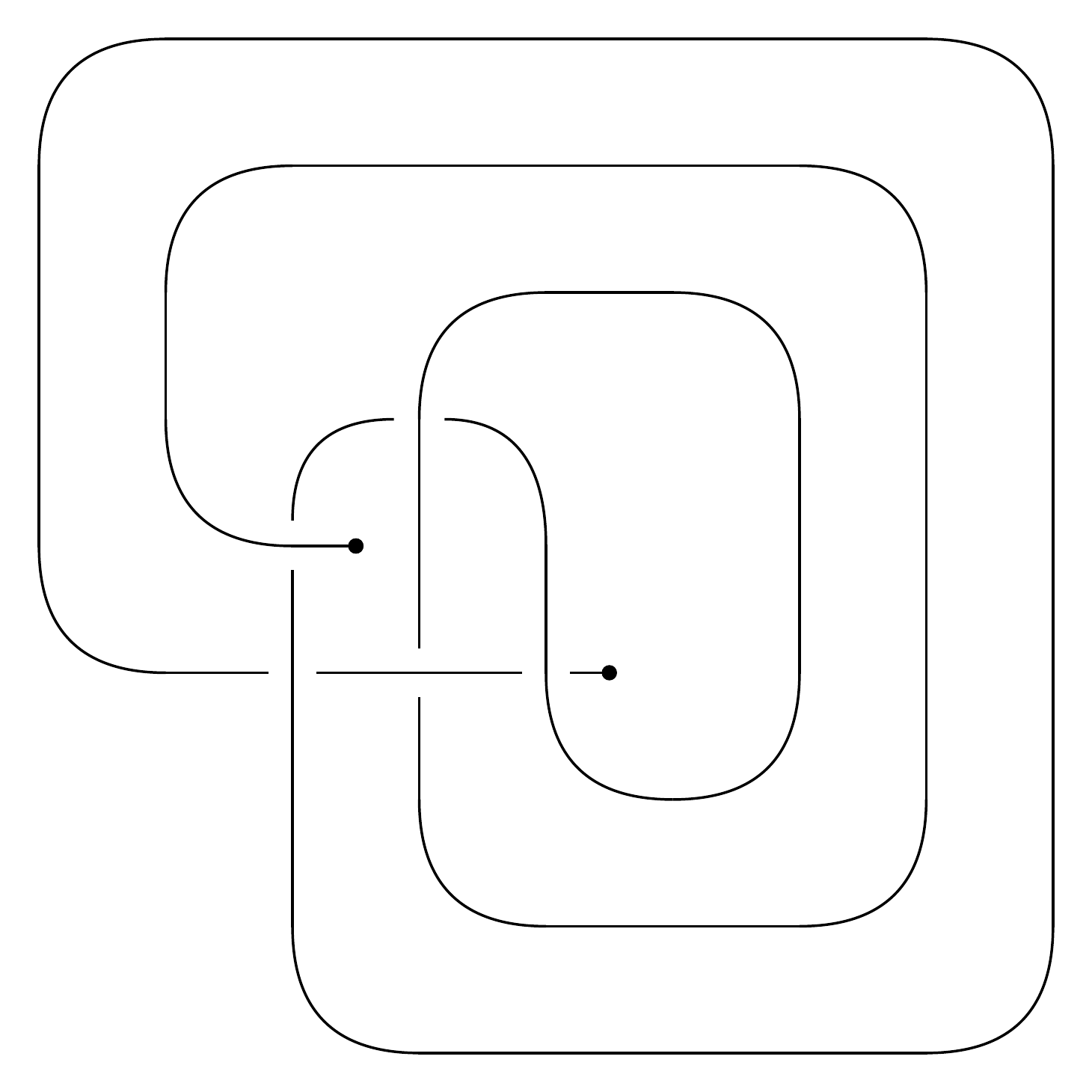}\\
\textcolor{black}{$5_{644}$}
\vspace{1cm}
\end{minipage}
\begin{minipage}[t]{.25\linewidth}
\centering
\includegraphics[width=0.9\textwidth,height=3.5cm,keepaspectratio]{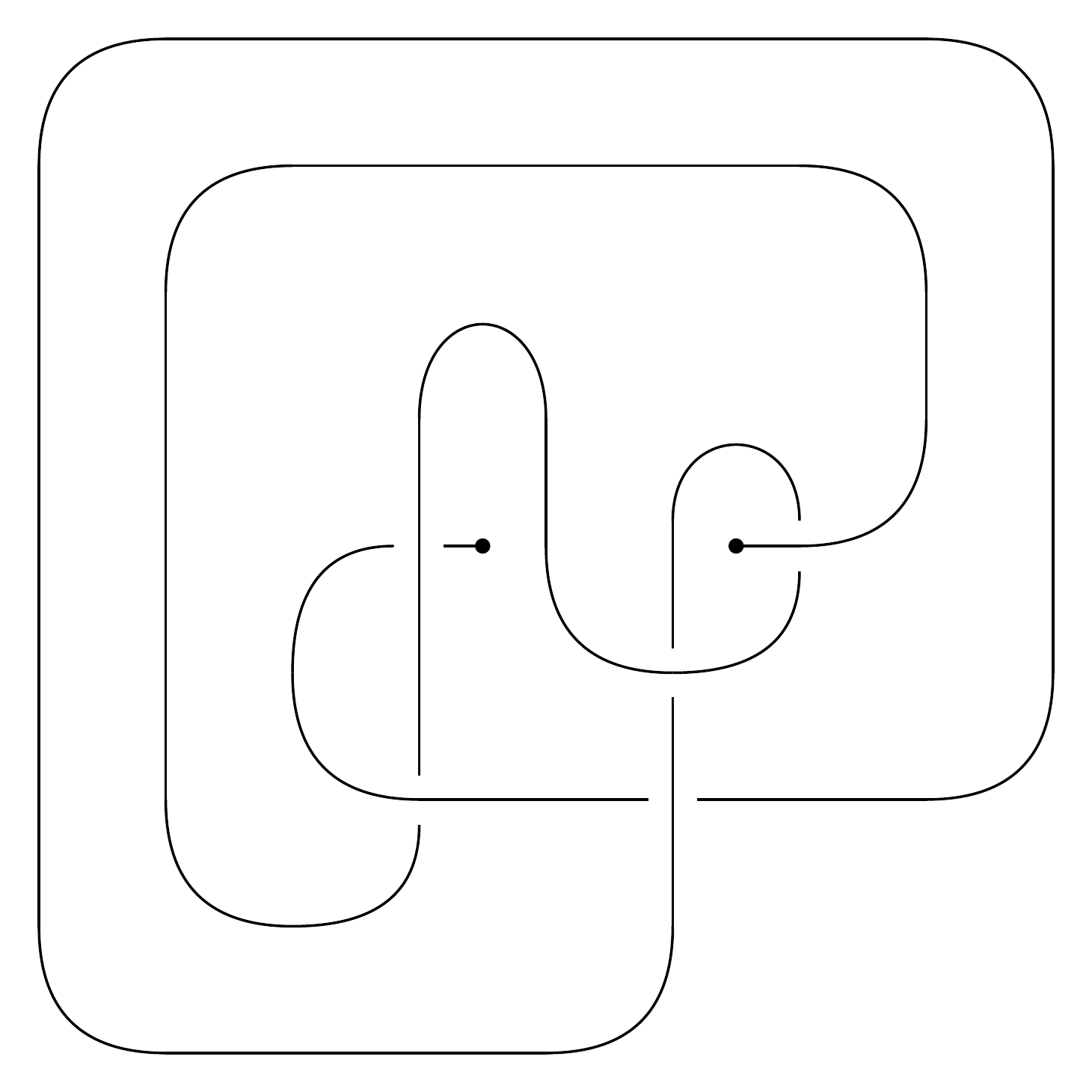}\\
\textcolor{black}{$5_{645}$}
\vspace{1cm}
\end{minipage}
\begin{minipage}[t]{.25\linewidth}
\centering
\includegraphics[width=0.9\textwidth,height=3.5cm,keepaspectratio]{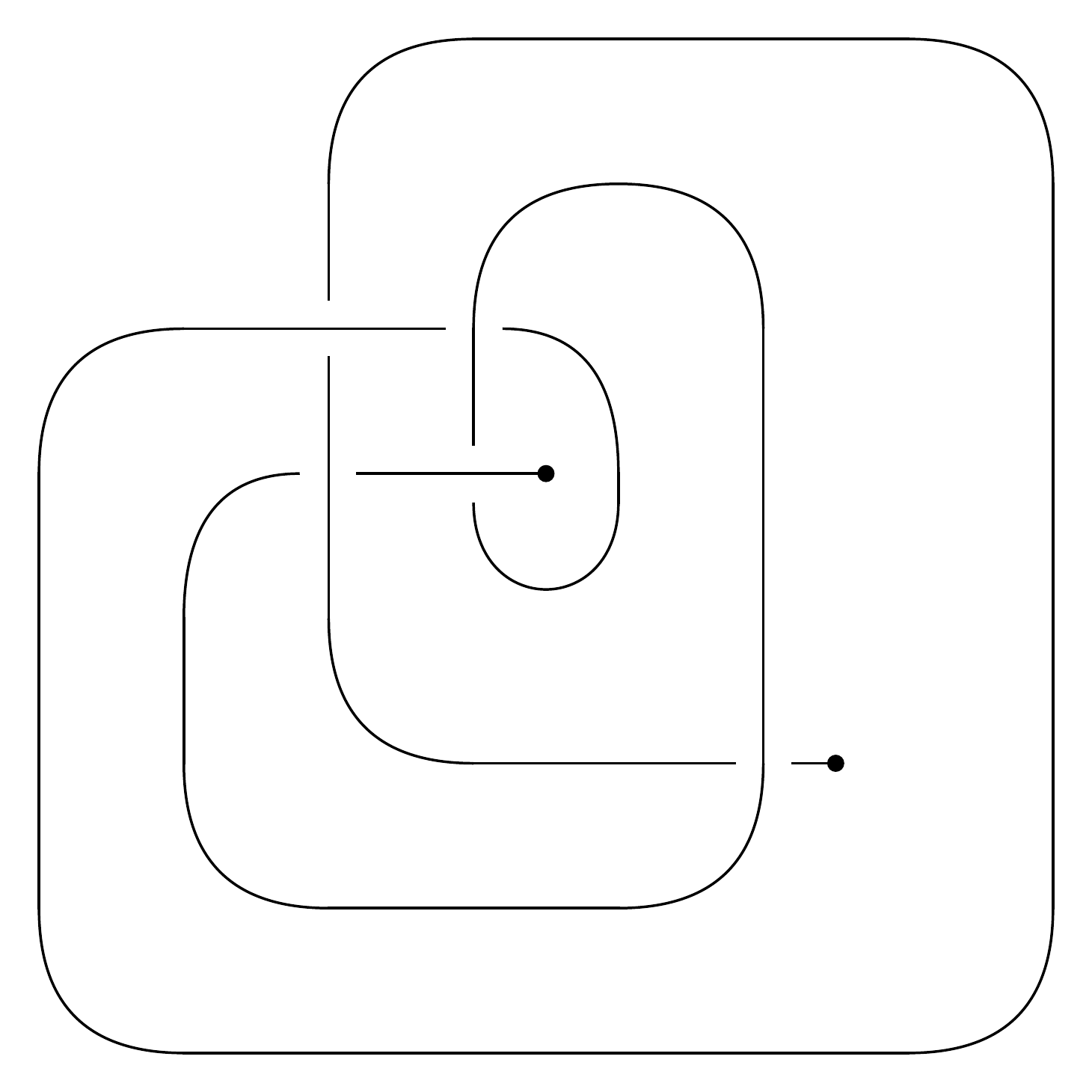}\\
\textcolor{black}{$5_{646}$}
\vspace{1cm}
\end{minipage}
\begin{minipage}[t]{.25\linewidth}
\centering
\includegraphics[width=0.9\textwidth,height=3.5cm,keepaspectratio]{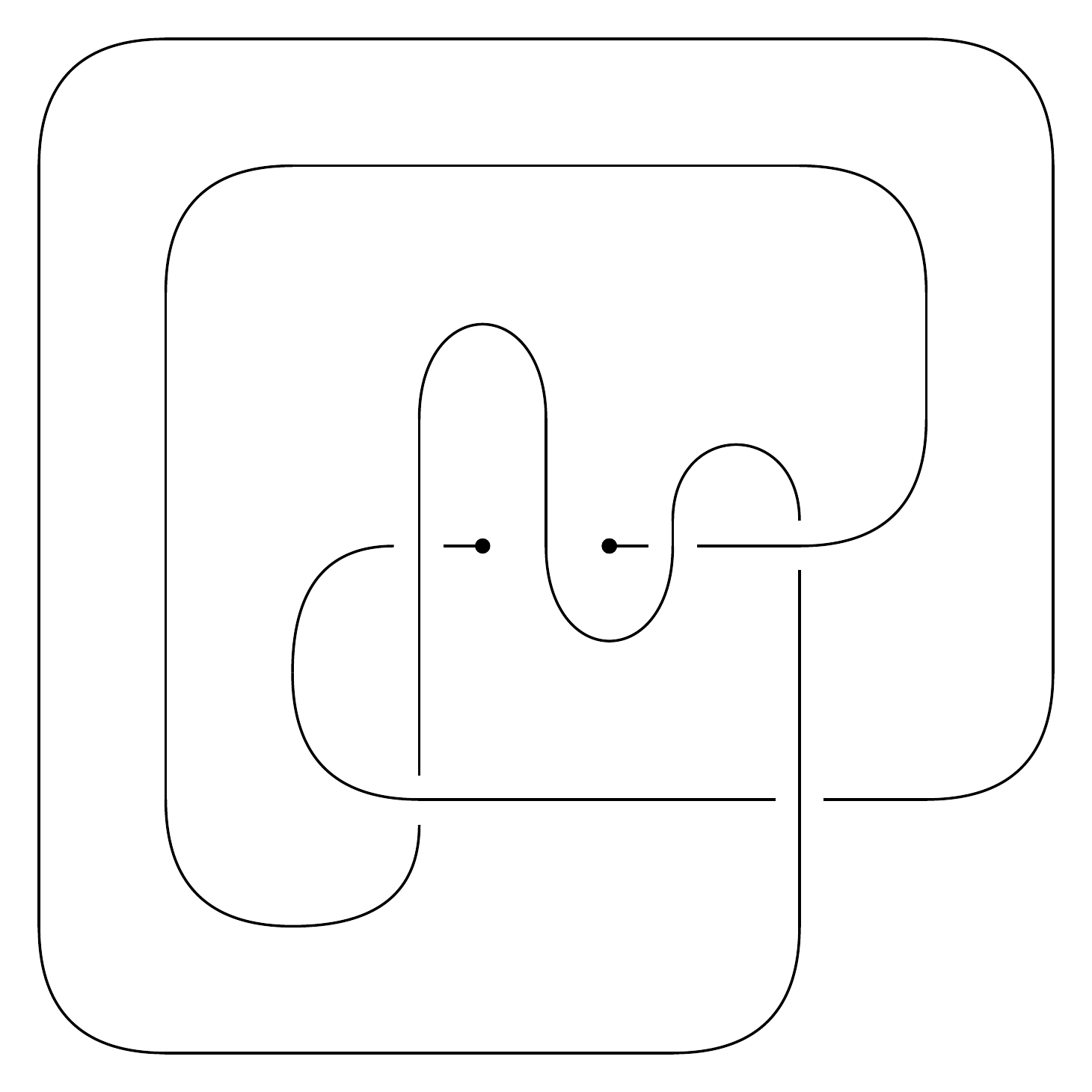}\\
\textcolor{black}{$5_{647}$}
\vspace{1cm}
\end{minipage}
\begin{minipage}[t]{.25\linewidth}
\centering
\includegraphics[width=0.9\textwidth,height=3.5cm,keepaspectratio]{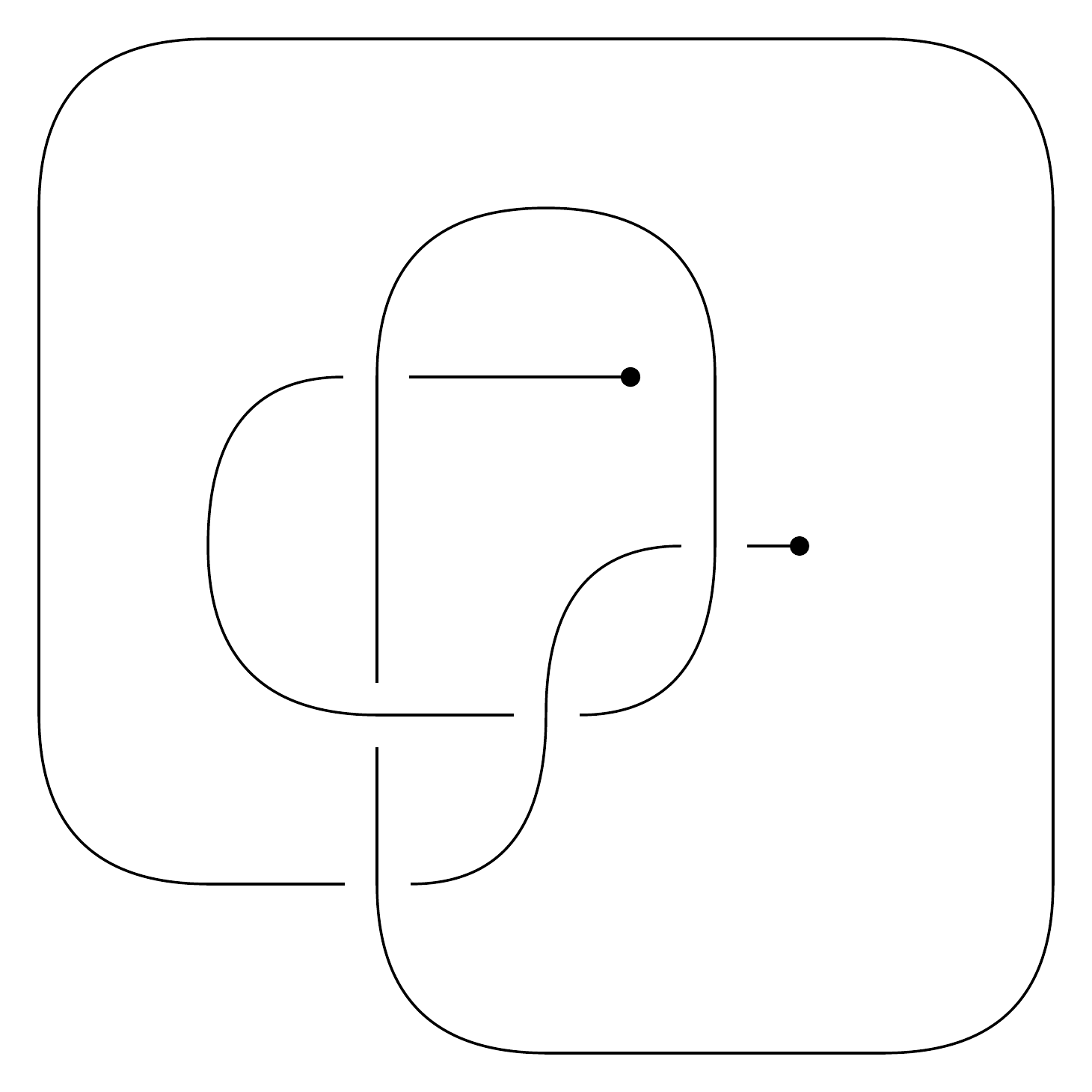}\\
\textcolor{black}{$5_{648}$}
\vspace{1cm}
\end{minipage}
\begin{minipage}[t]{.25\linewidth}
\centering
\includegraphics[width=0.9\textwidth,height=3.5cm,keepaspectratio]{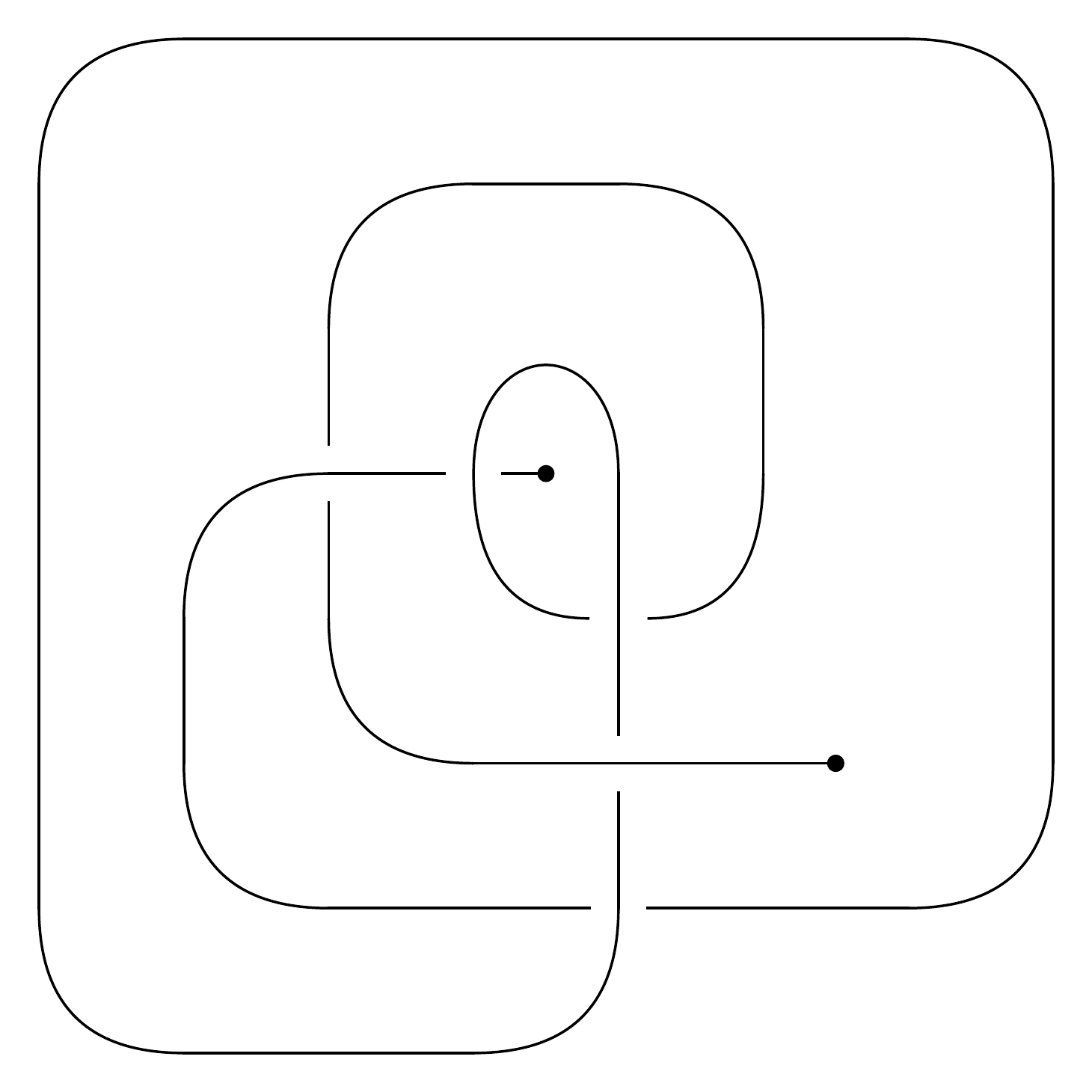}\\
\textcolor{black}{$5_{649}$}
\vspace{1cm}
\end{minipage}
\begin{minipage}[t]{.25\linewidth}
\centering
\includegraphics[width=0.9\textwidth,height=3.5cm,keepaspectratio]{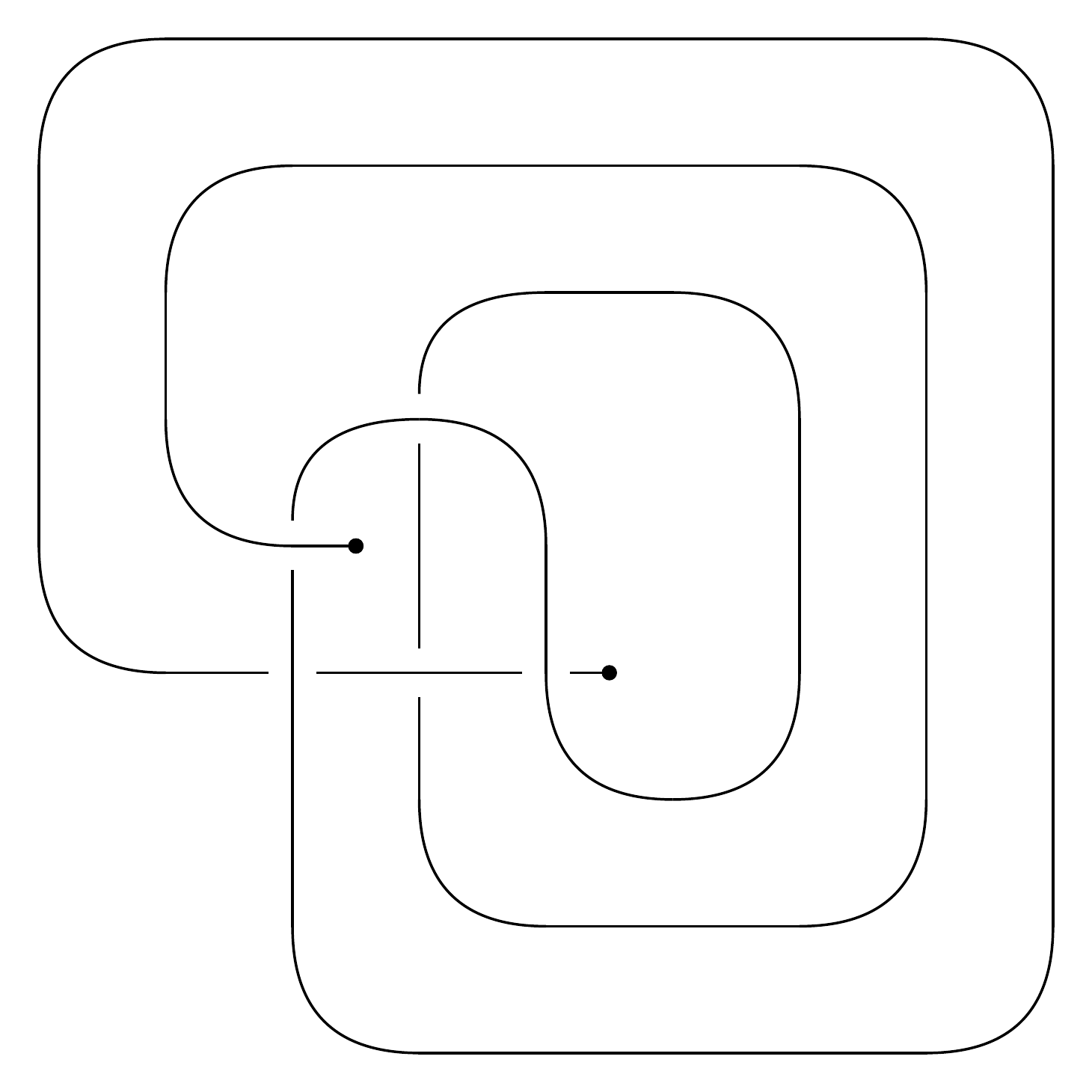}\\
\textcolor{black}{$5_{650}$}
\vspace{1cm}
\end{minipage}
\begin{minipage}[t]{.25\linewidth}
\centering
\includegraphics[width=0.9\textwidth,height=3.5cm,keepaspectratio]{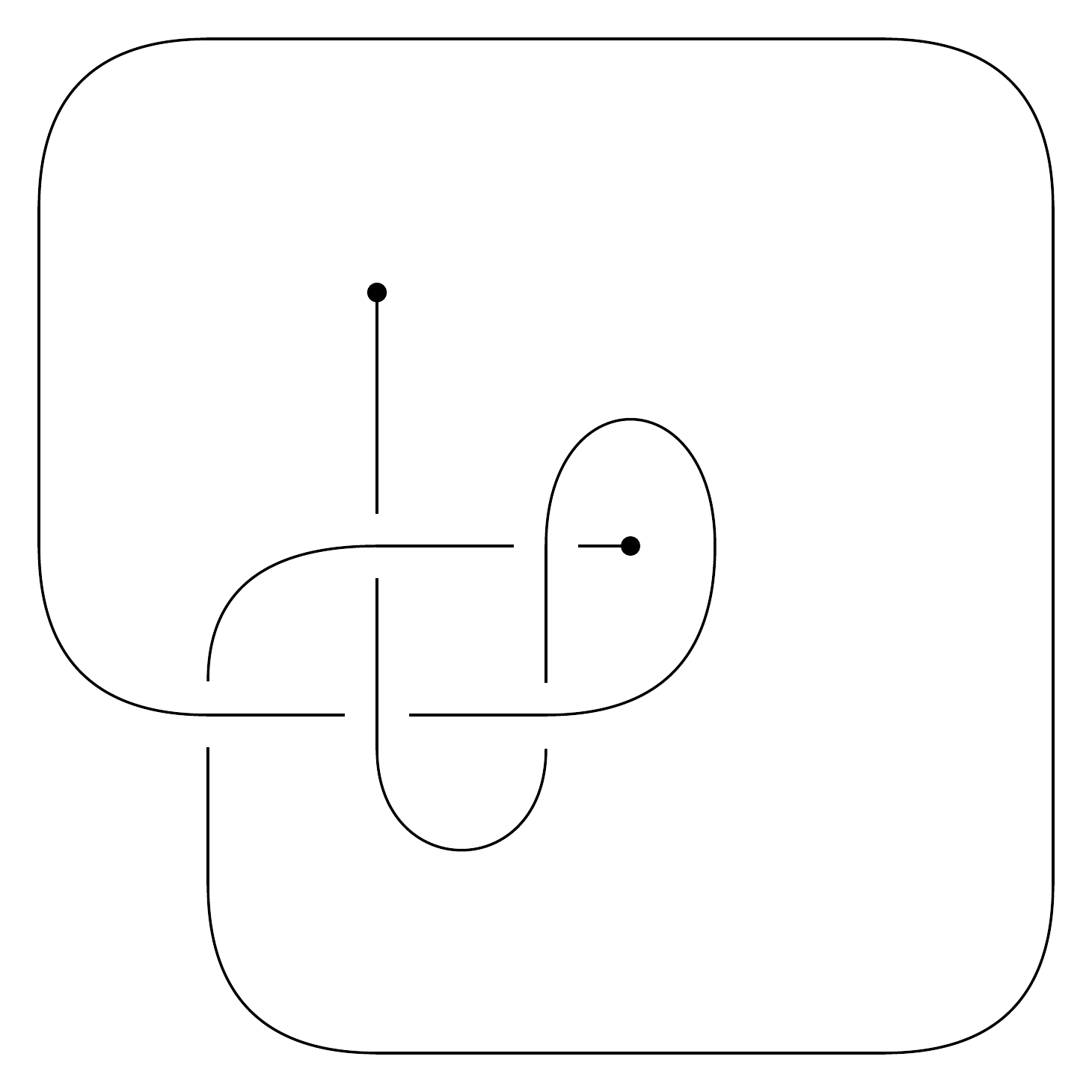}\\
\textcolor{black}{$5_{651}$}
\vspace{1cm}
\end{minipage}
\begin{minipage}[t]{.25\linewidth}
\centering
\includegraphics[width=0.9\textwidth,height=3.5cm,keepaspectratio]{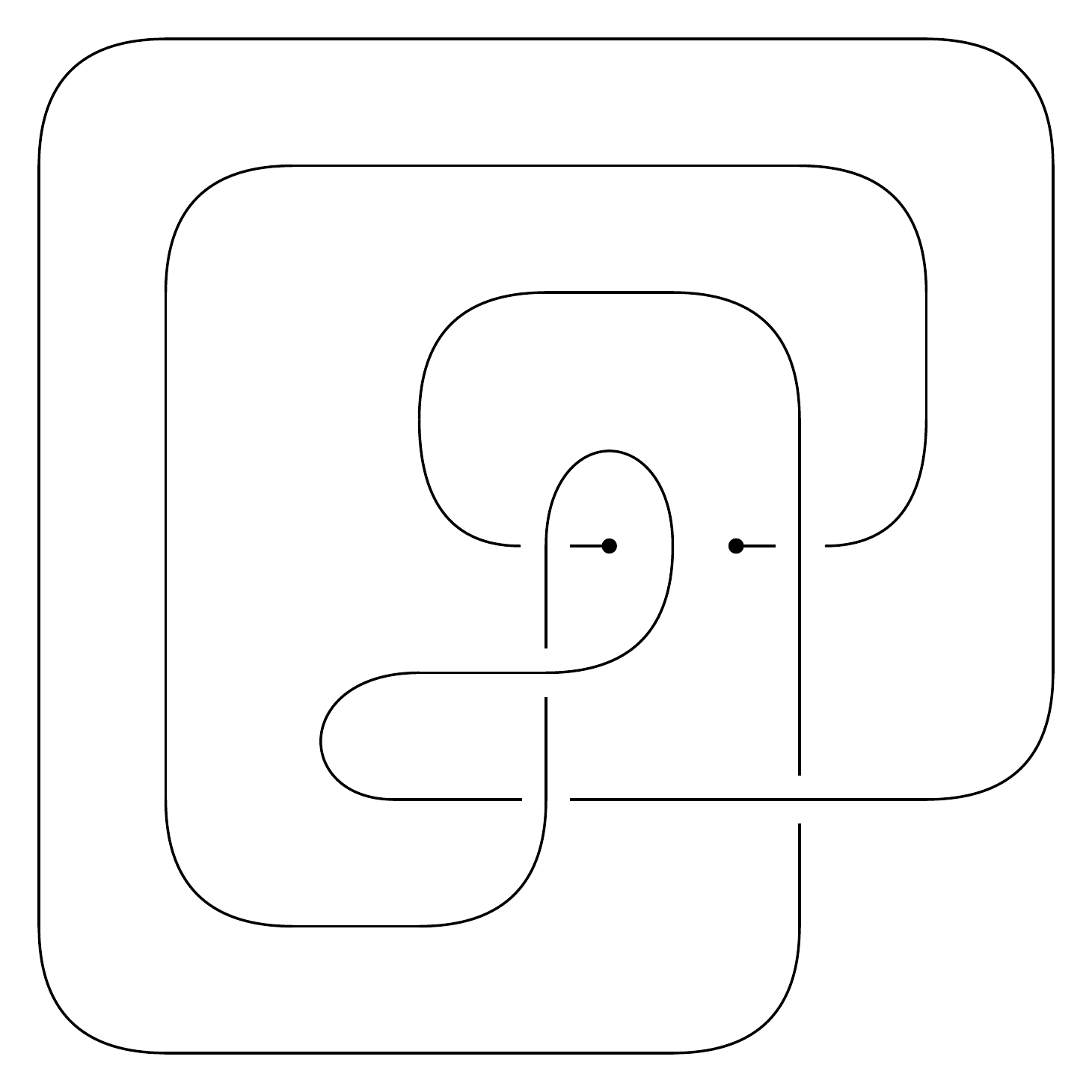}\\
\textcolor{black}{$5_{652}$}
\vspace{1cm}
\end{minipage}
\begin{minipage}[t]{.25\linewidth}
\centering
\includegraphics[width=0.9\textwidth,height=3.5cm,keepaspectratio]{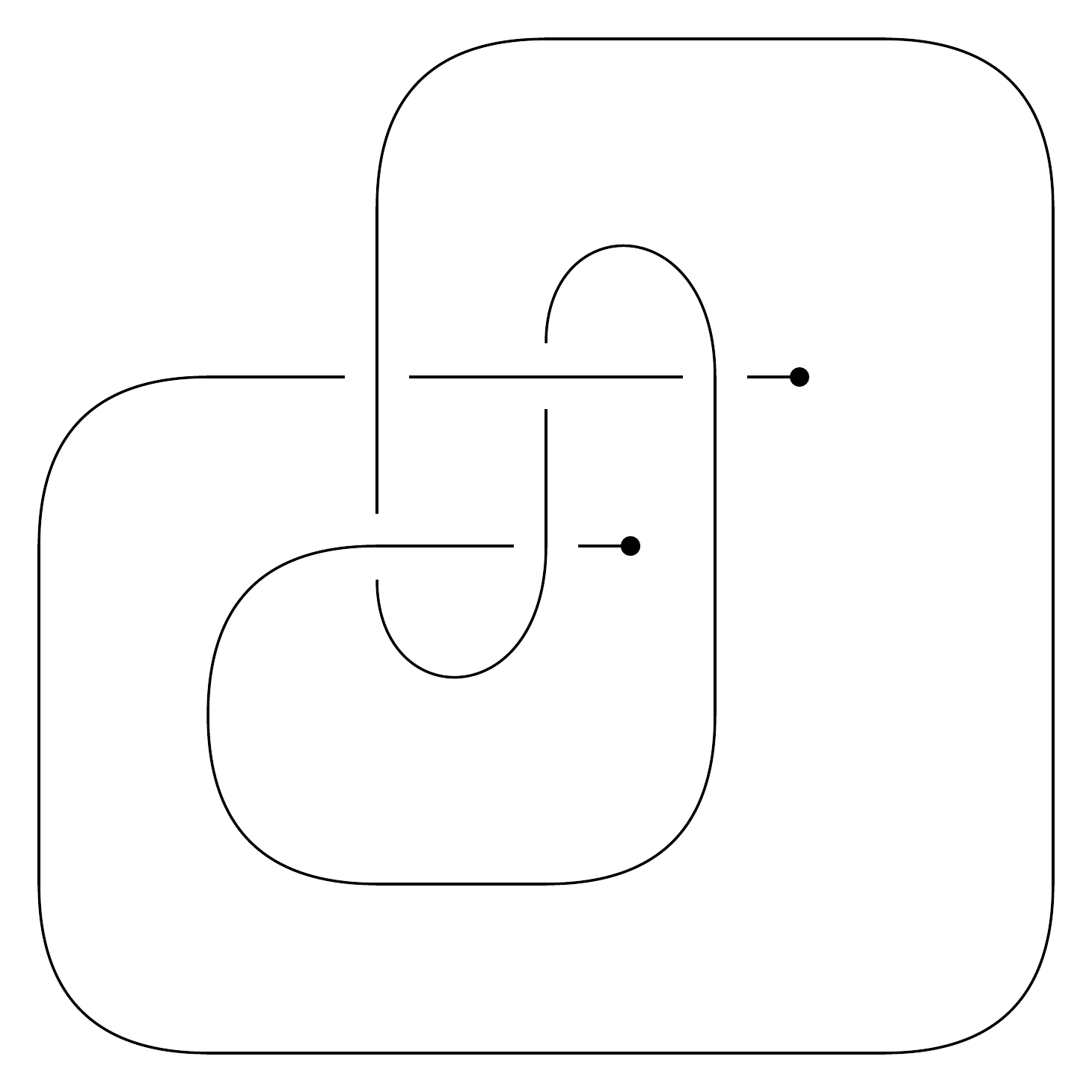}\\
\textcolor{black}{$5_{653}$}
\vspace{1cm}
\end{minipage}
\begin{minipage}[t]{.25\linewidth}
\centering
\includegraphics[width=0.9\textwidth,height=3.5cm,keepaspectratio]{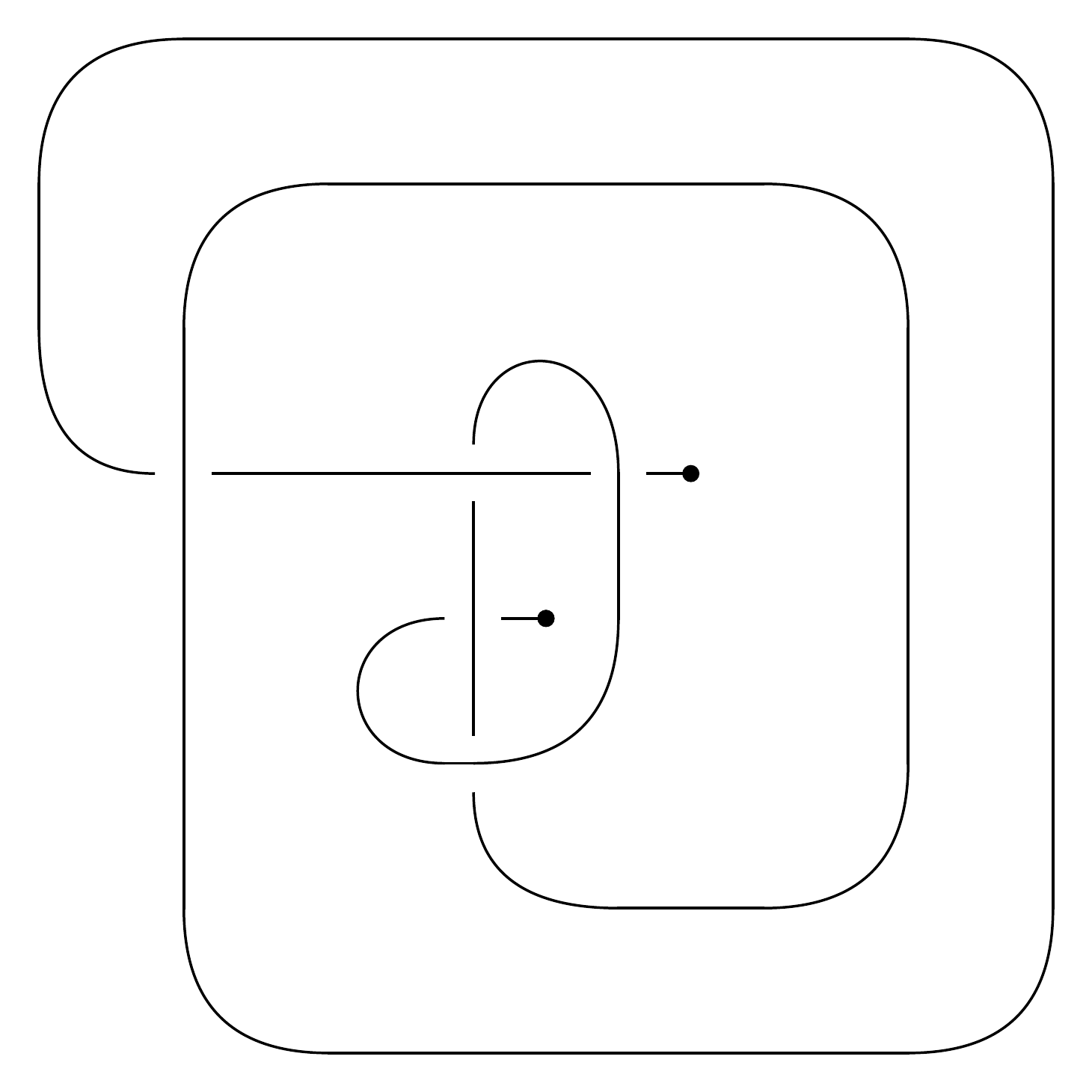}\\
\textcolor{black}{$5_{654}$}
\vspace{1cm}
\end{minipage}
\begin{minipage}[t]{.25\linewidth}
\centering
\includegraphics[width=0.9\textwidth,height=3.5cm,keepaspectratio]{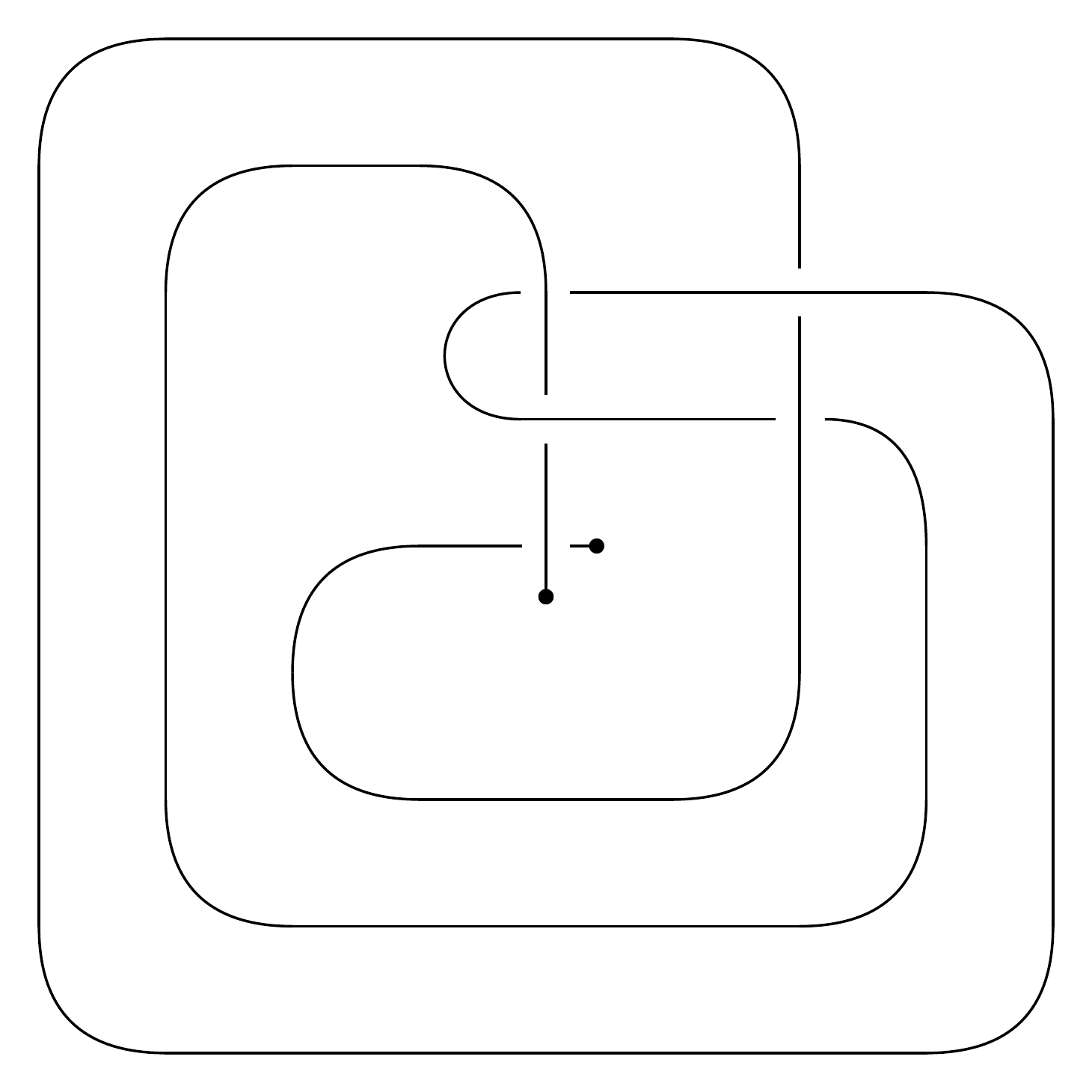}\\
\textcolor{black}{$5_{655}$}
\vspace{1cm}
\end{minipage}
\begin{minipage}[t]{.25\linewidth}
\centering
\includegraphics[width=0.9\textwidth,height=3.5cm,keepaspectratio]{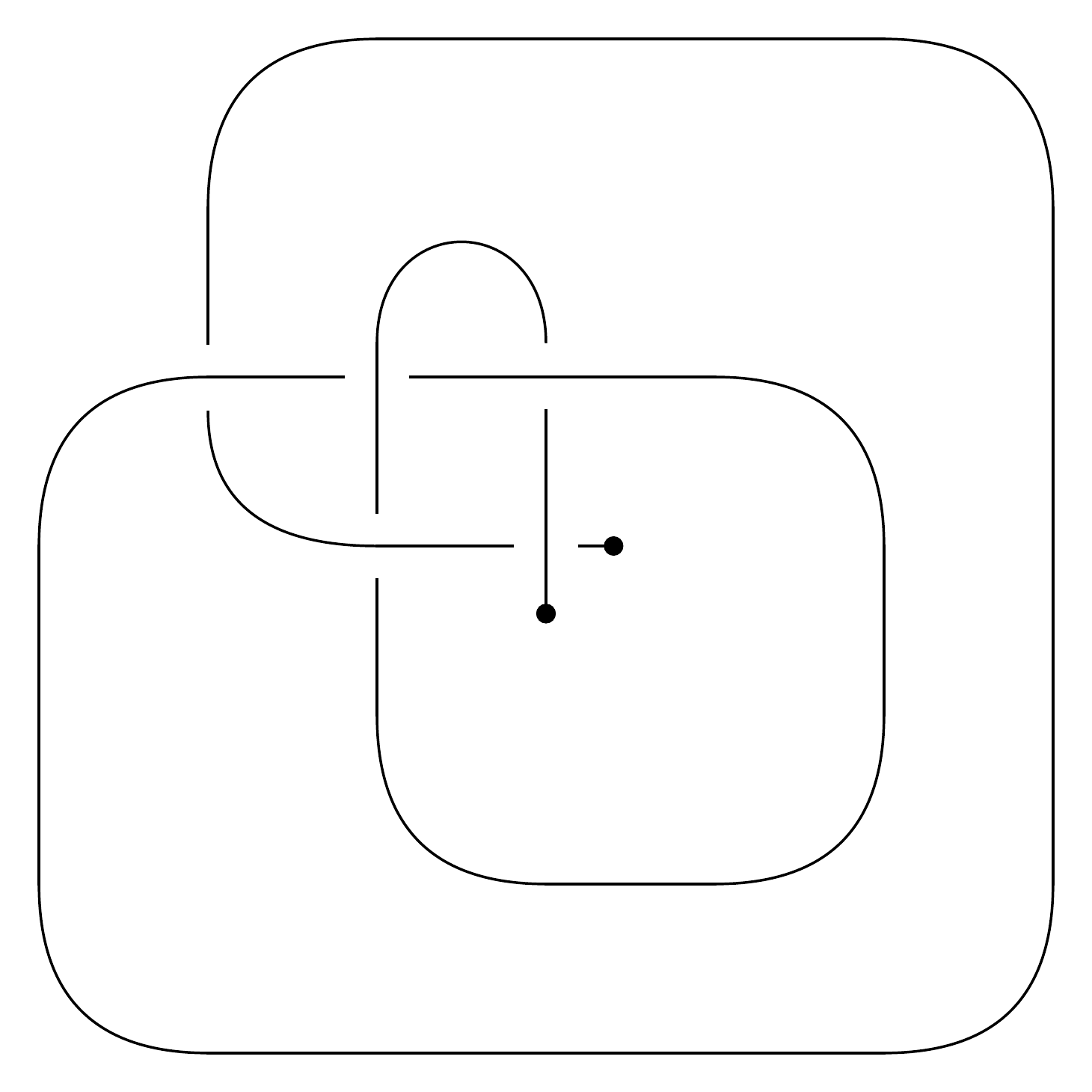}\\
\textcolor{black}{$5_{656}$}
\vspace{1cm}
\end{minipage}
\begin{minipage}[t]{.25\linewidth}
\centering
\includegraphics[width=0.9\textwidth,height=3.5cm,keepaspectratio]{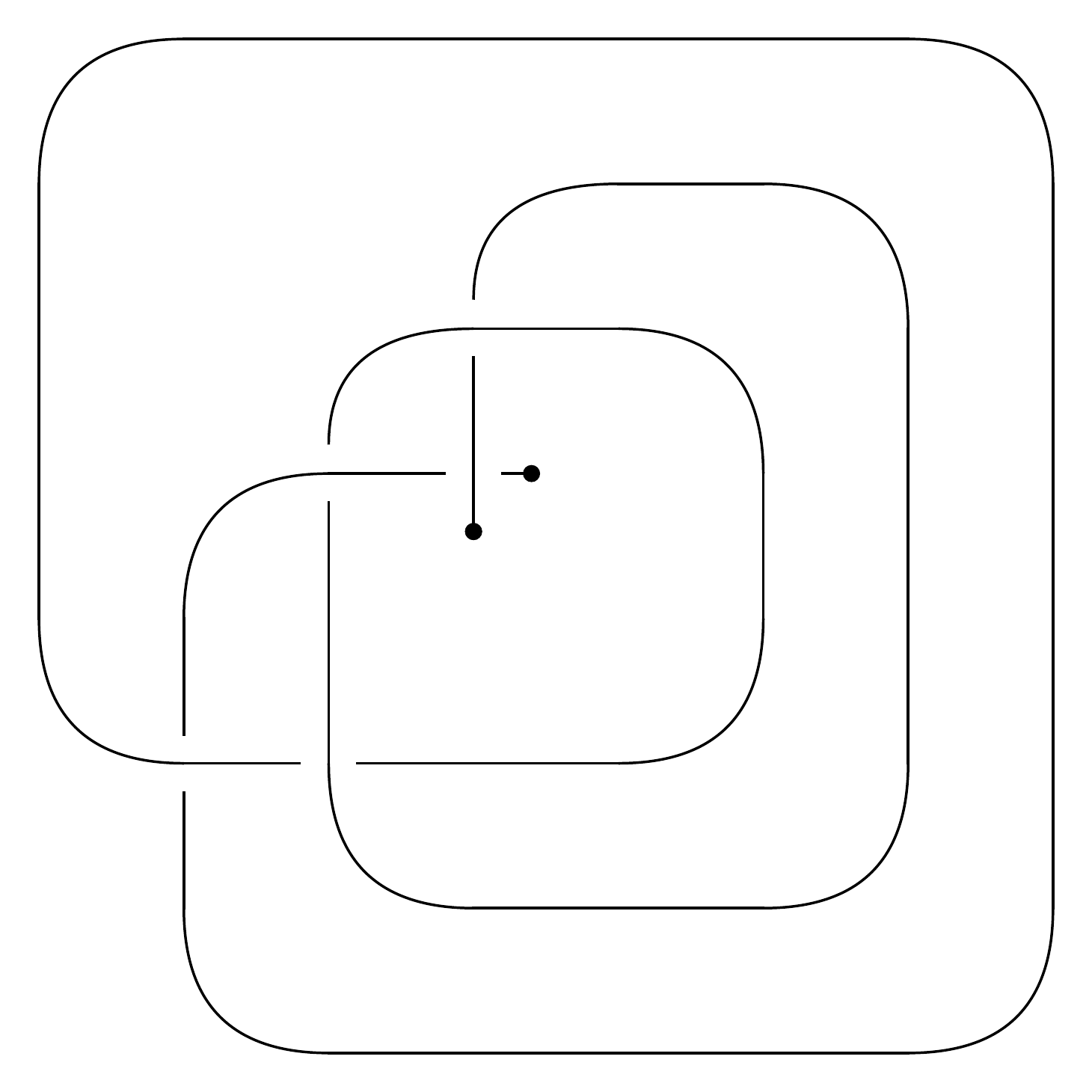}\\
\textcolor{black}{$5_{657}$}
\vspace{1cm}
\end{minipage}
\begin{minipage}[t]{.25\linewidth}
\centering
\includegraphics[width=0.9\textwidth,height=3.5cm,keepaspectratio]{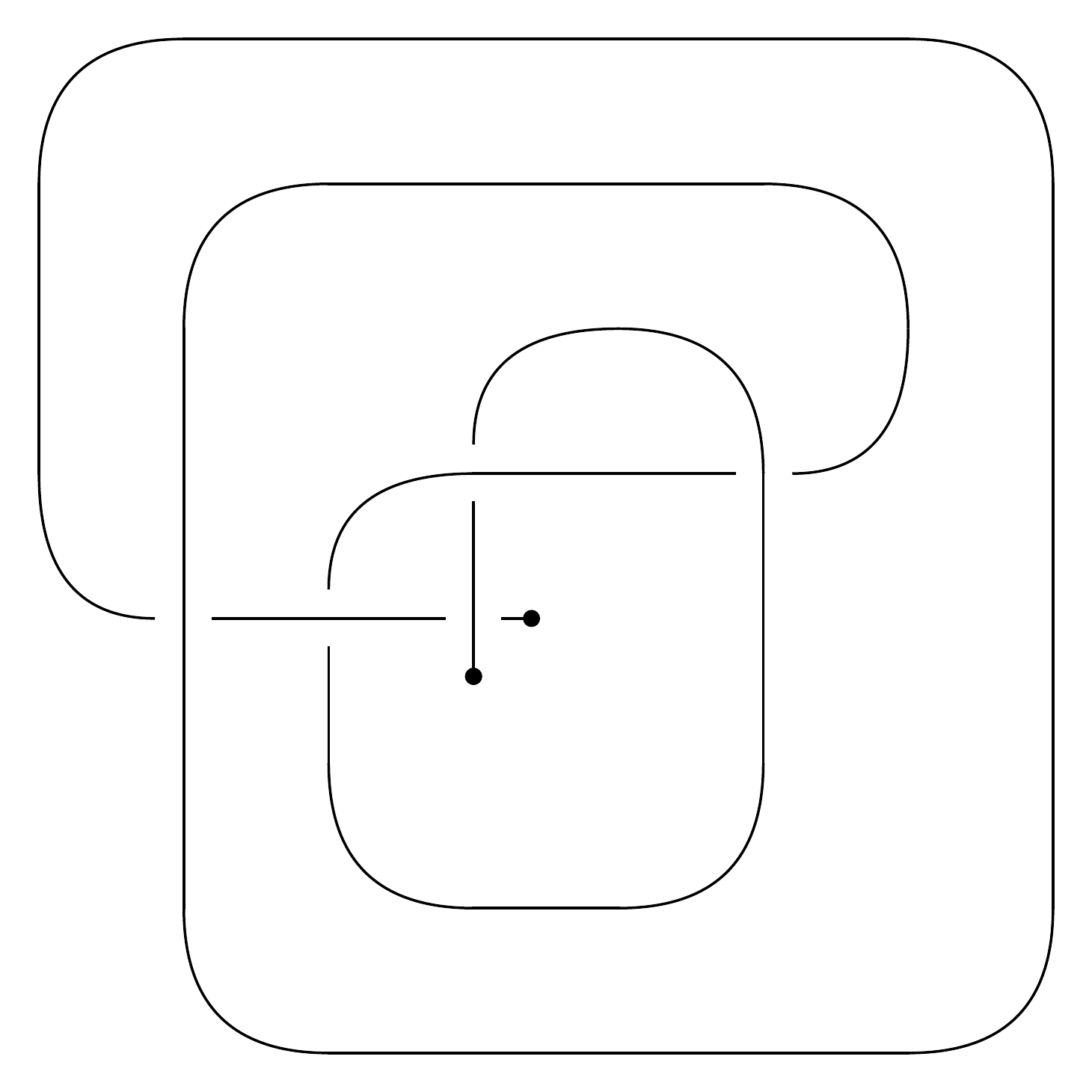}\\
\textcolor{black}{$5_{658}$}
\vspace{1cm}
\end{minipage}
\begin{minipage}[t]{.25\linewidth}
\centering
\includegraphics[width=0.9\textwidth,height=3.5cm,keepaspectratio]{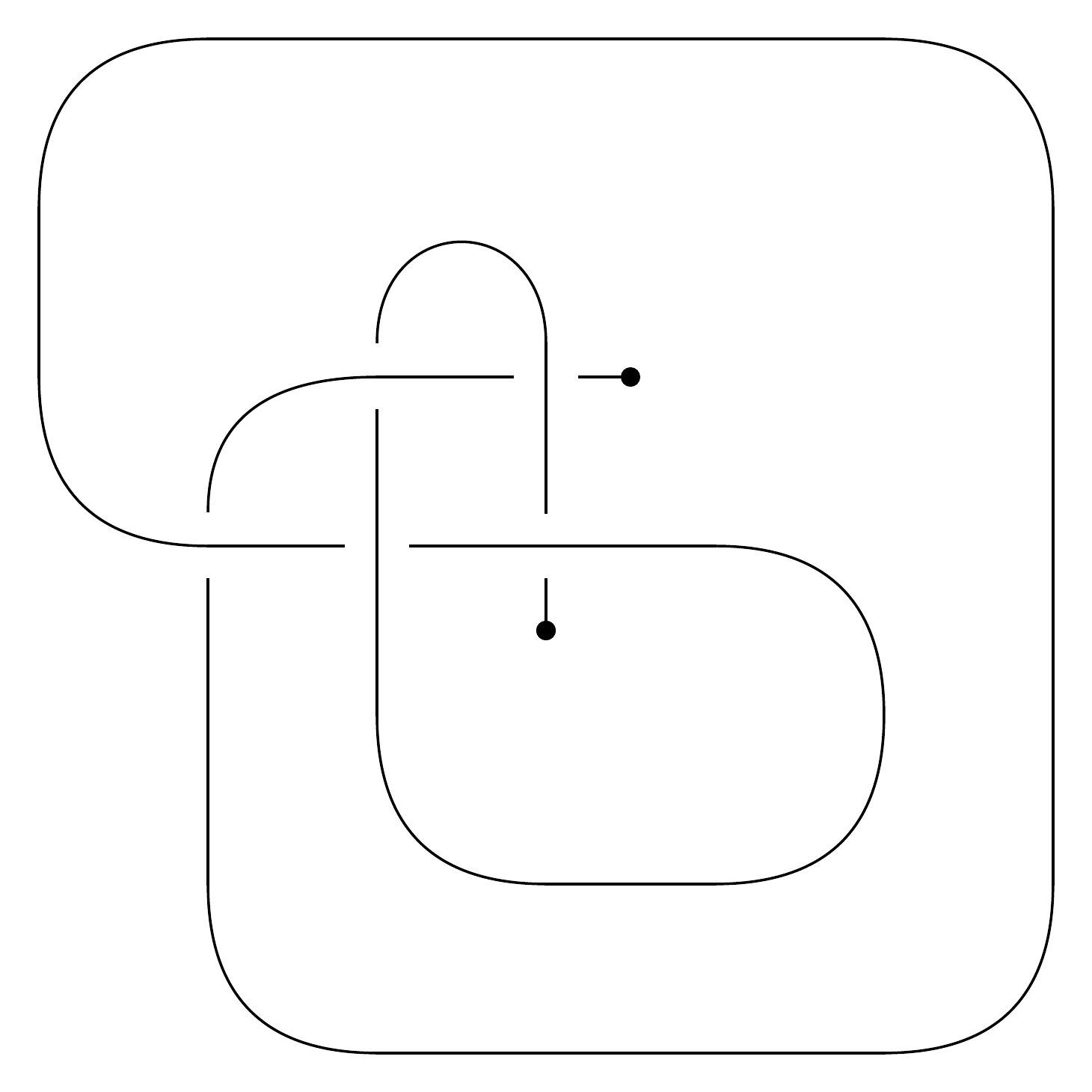}\\
\textcolor{black}{$5_{659}$}
\vspace{1cm}
\end{minipage}
\begin{minipage}[t]{.25\linewidth}
\centering
\includegraphics[width=0.9\textwidth,height=3.5cm,keepaspectratio]{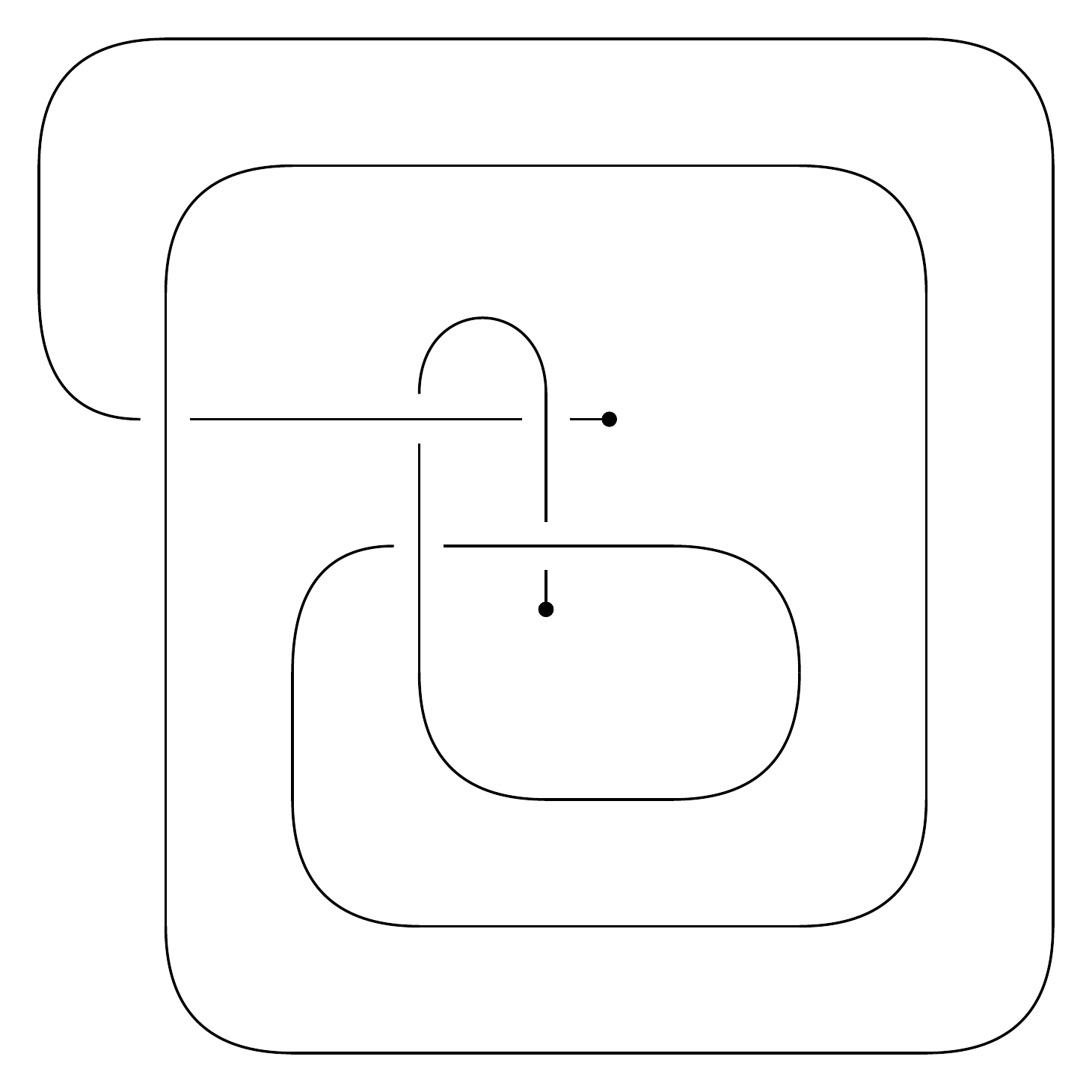}\\
\textcolor{black}{$5_{660}$}
\vspace{1cm}
\end{minipage}
\begin{minipage}[t]{.25\linewidth}
\centering
\includegraphics[width=0.9\textwidth,height=3.5cm,keepaspectratio]{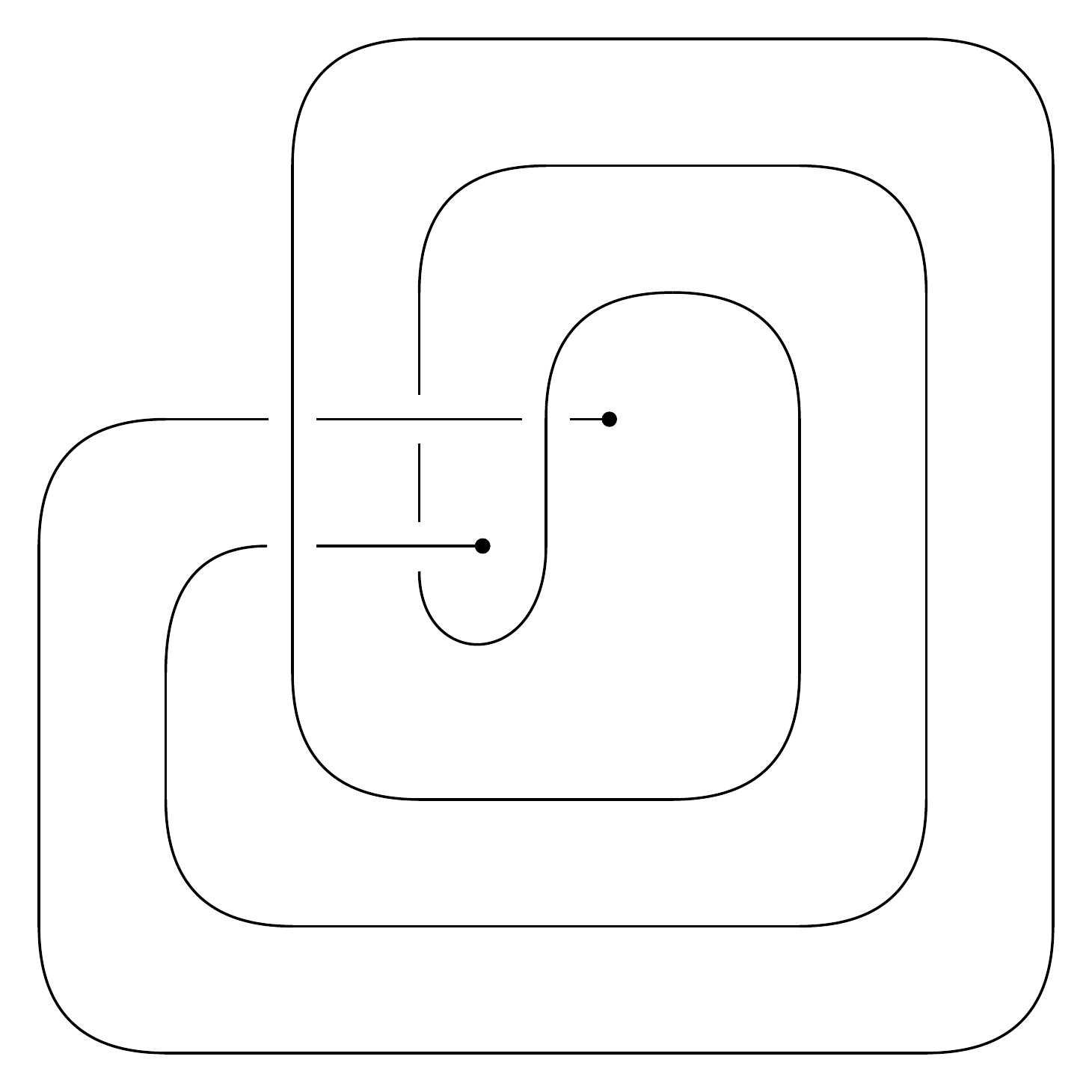}\\
\textcolor{black}{$5_{661}$}
\vspace{1cm}
\end{minipage}
\begin{minipage}[t]{.25\linewidth}
\centering
\includegraphics[width=0.9\textwidth,height=3.5cm,keepaspectratio]{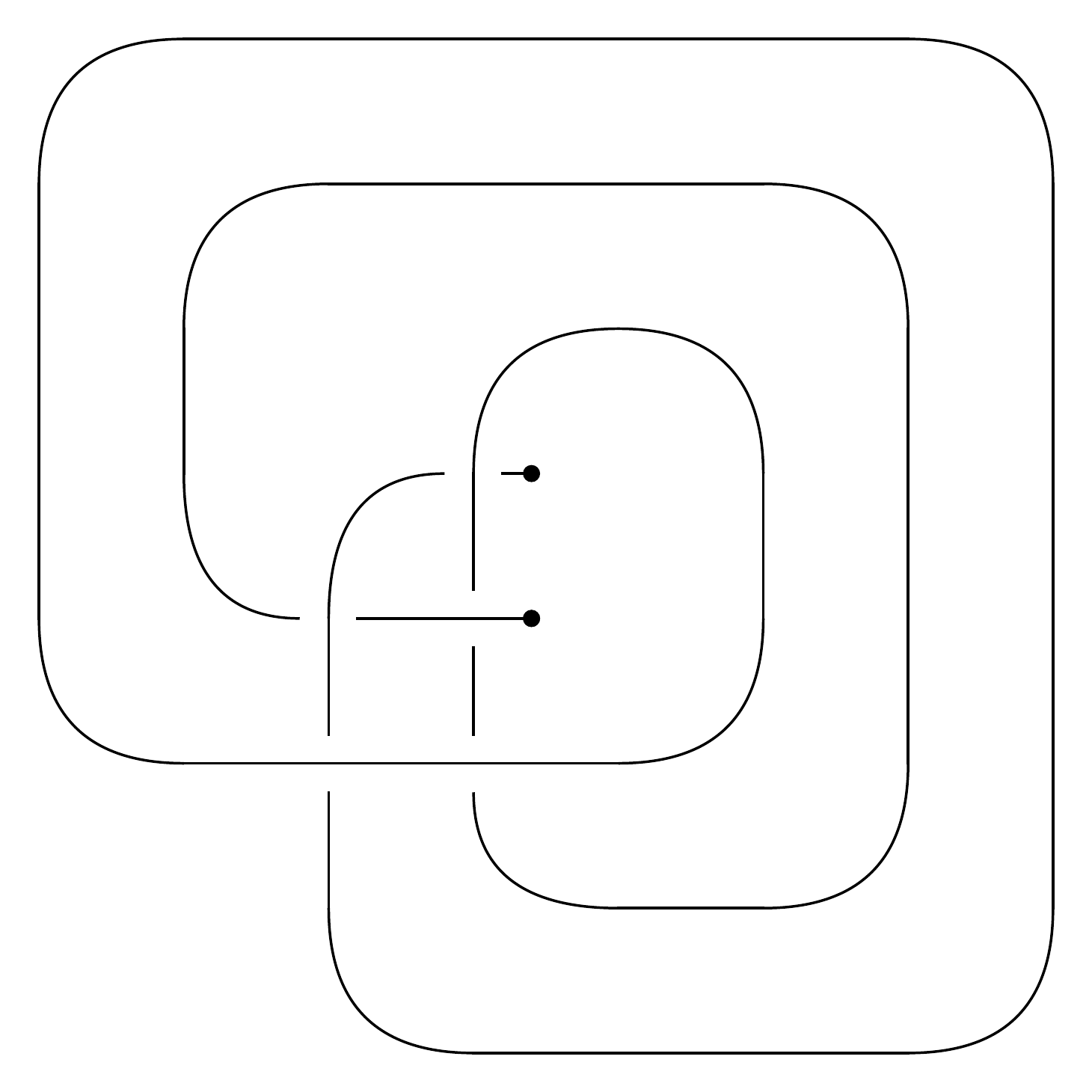}\\
\textcolor{black}{$5_{662}$}
\vspace{1cm}
\end{minipage}
\begin{minipage}[t]{.25\linewidth}
\centering
\includegraphics[width=0.9\textwidth,height=3.5cm,keepaspectratio]{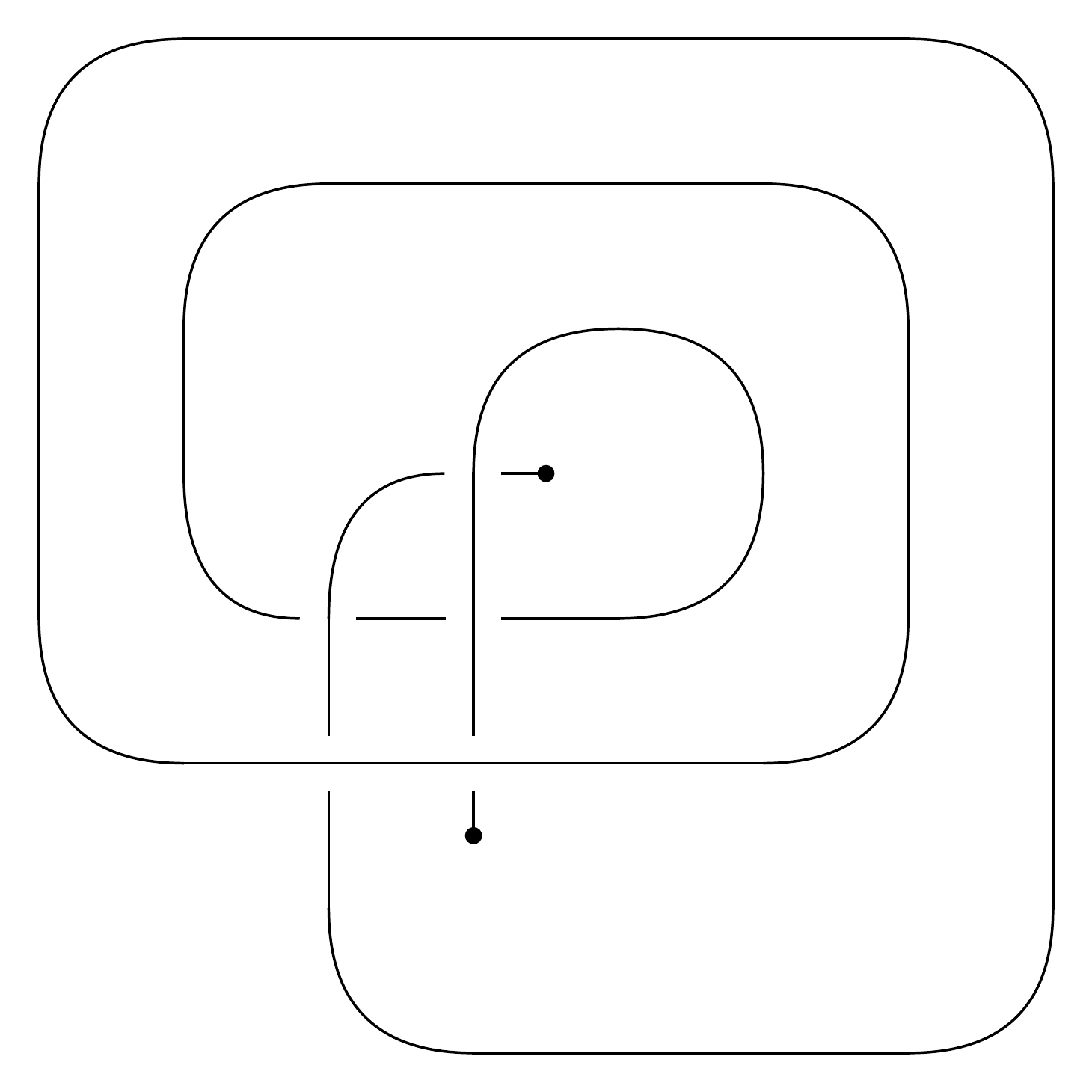}\\
\textcolor{black}{$5_{663}$}
\vspace{1cm}
\end{minipage}
\begin{minipage}[t]{.25\linewidth}
\centering
\includegraphics[width=0.9\textwidth,height=3.5cm,keepaspectratio]{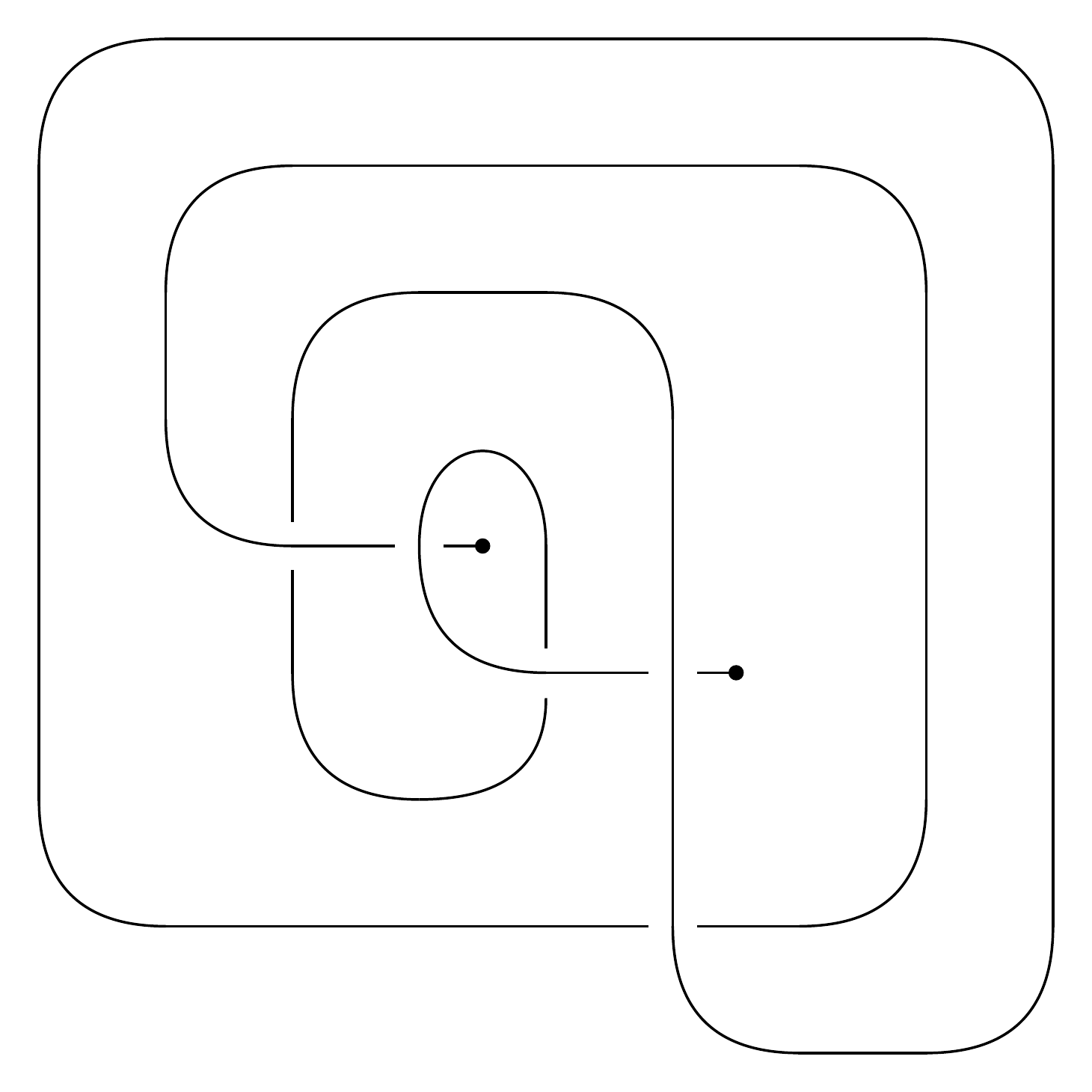}\\
\textcolor{black}{$5_{664}$}
\vspace{1cm}
\end{minipage}
\begin{minipage}[t]{.25\linewidth}
\centering
\includegraphics[width=0.9\textwidth,height=3.5cm,keepaspectratio]{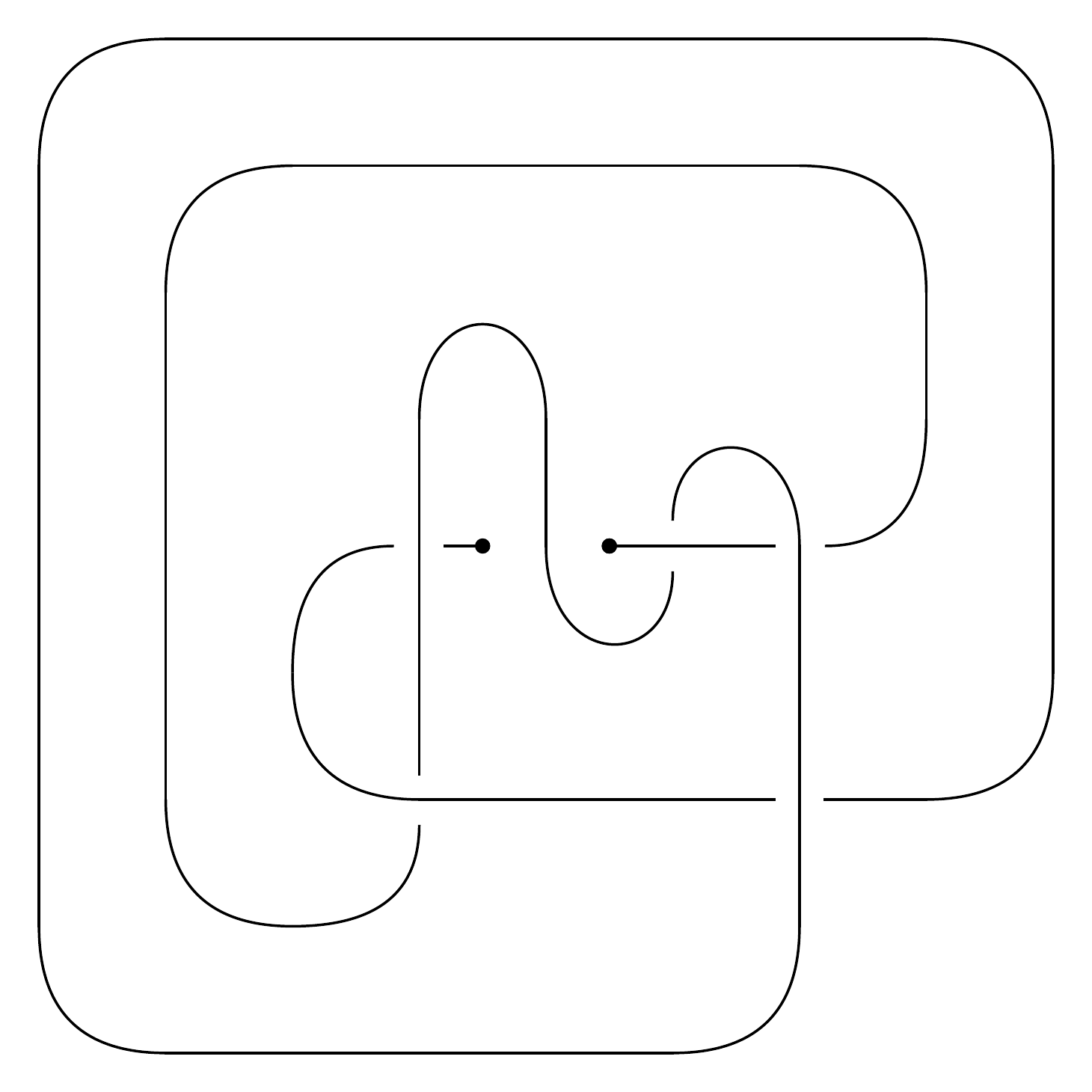}\\
\textcolor{black}{$5_{665}$}
\vspace{1cm}
\end{minipage}
\begin{minipage}[t]{.25\linewidth}
\centering
\includegraphics[width=0.9\textwidth,height=3.5cm,keepaspectratio]{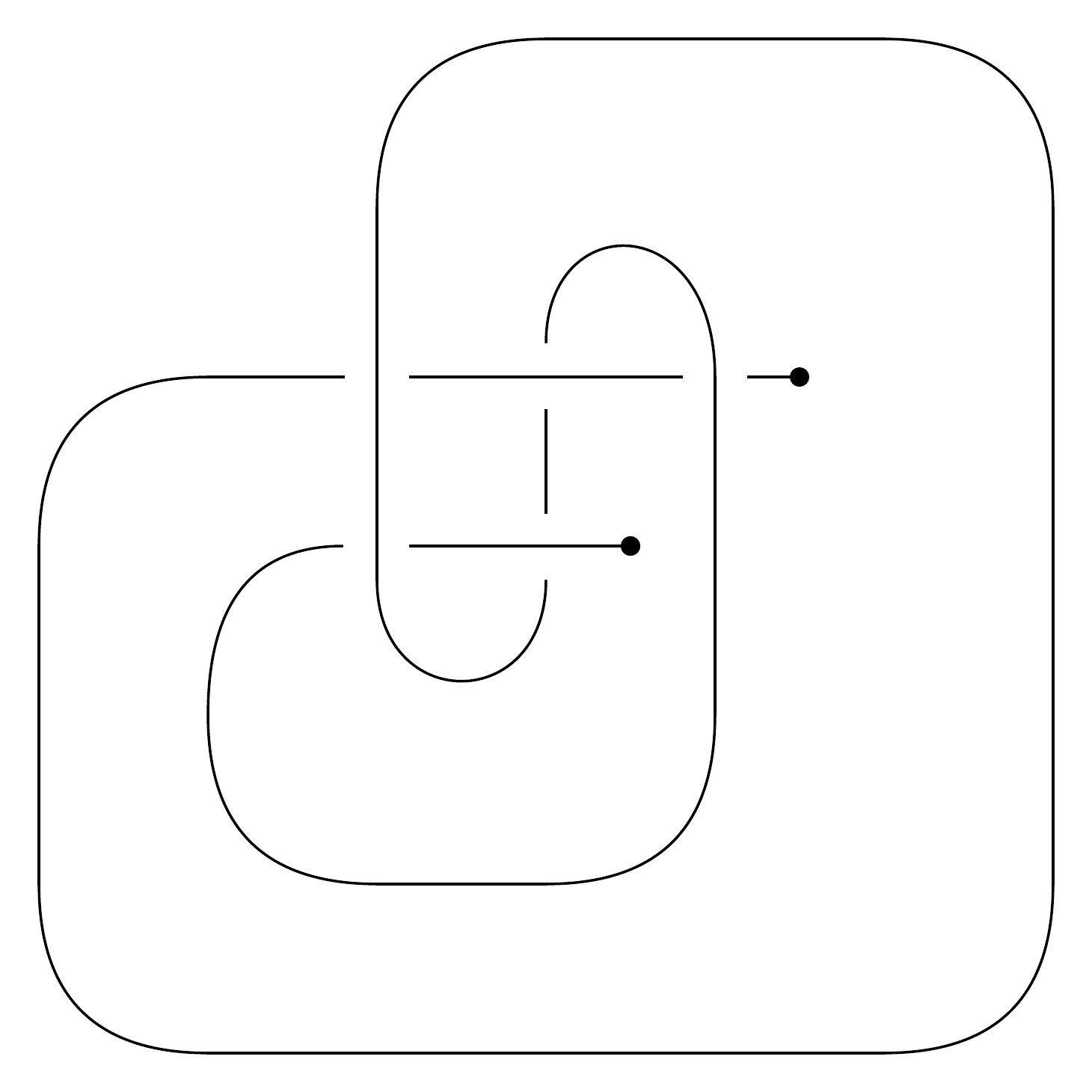}\\
\textcolor{black}{$5_{666}$}
\vspace{1cm}
\end{minipage}
\begin{minipage}[t]{.25\linewidth}
\centering
\includegraphics[width=0.9\textwidth,height=3.5cm,keepaspectratio]{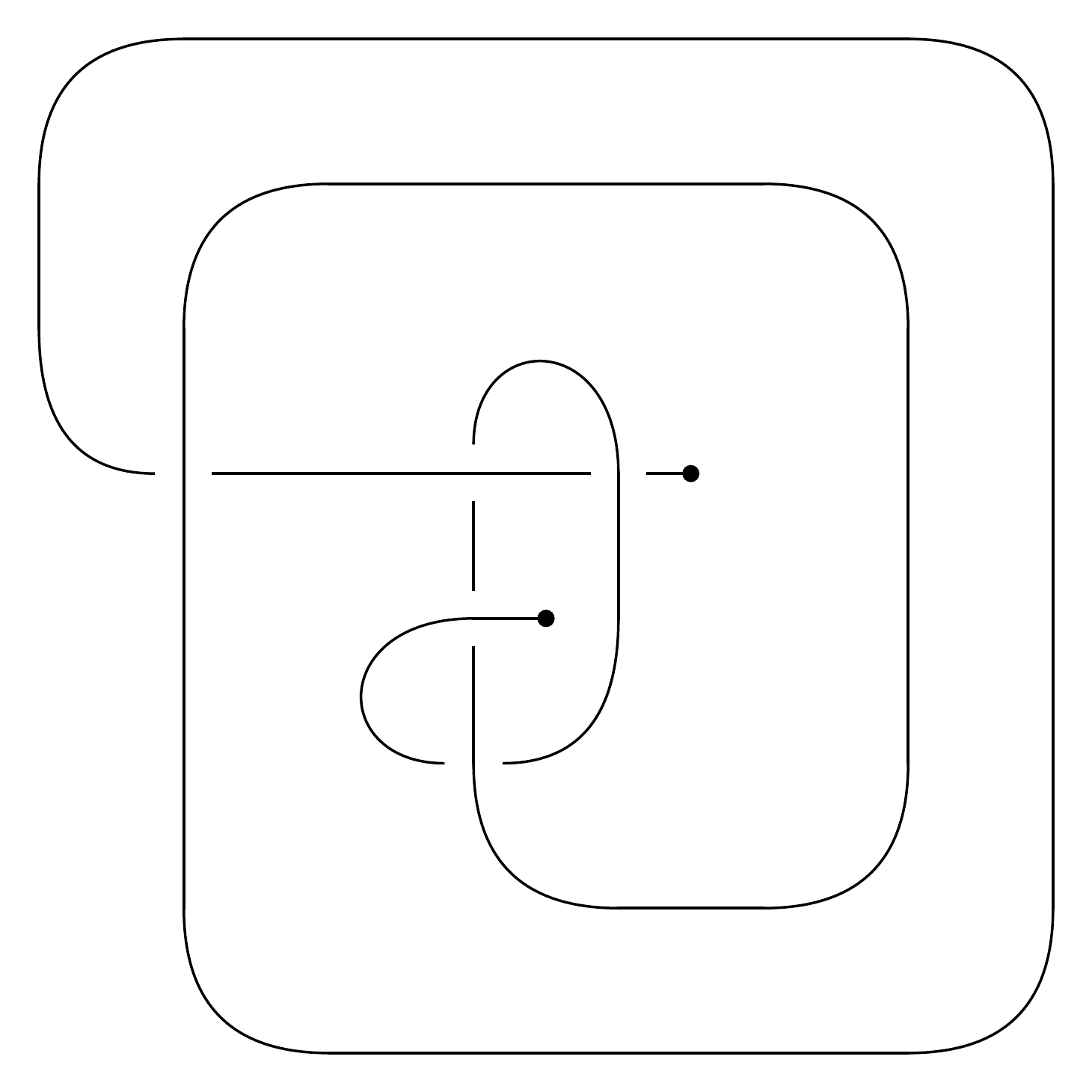}\\
\textcolor{black}{$5_{667}$}
\vspace{1cm}
\end{minipage}
\begin{minipage}[t]{.25\linewidth}
\centering
\includegraphics[width=0.9\textwidth,height=3.5cm,keepaspectratio]{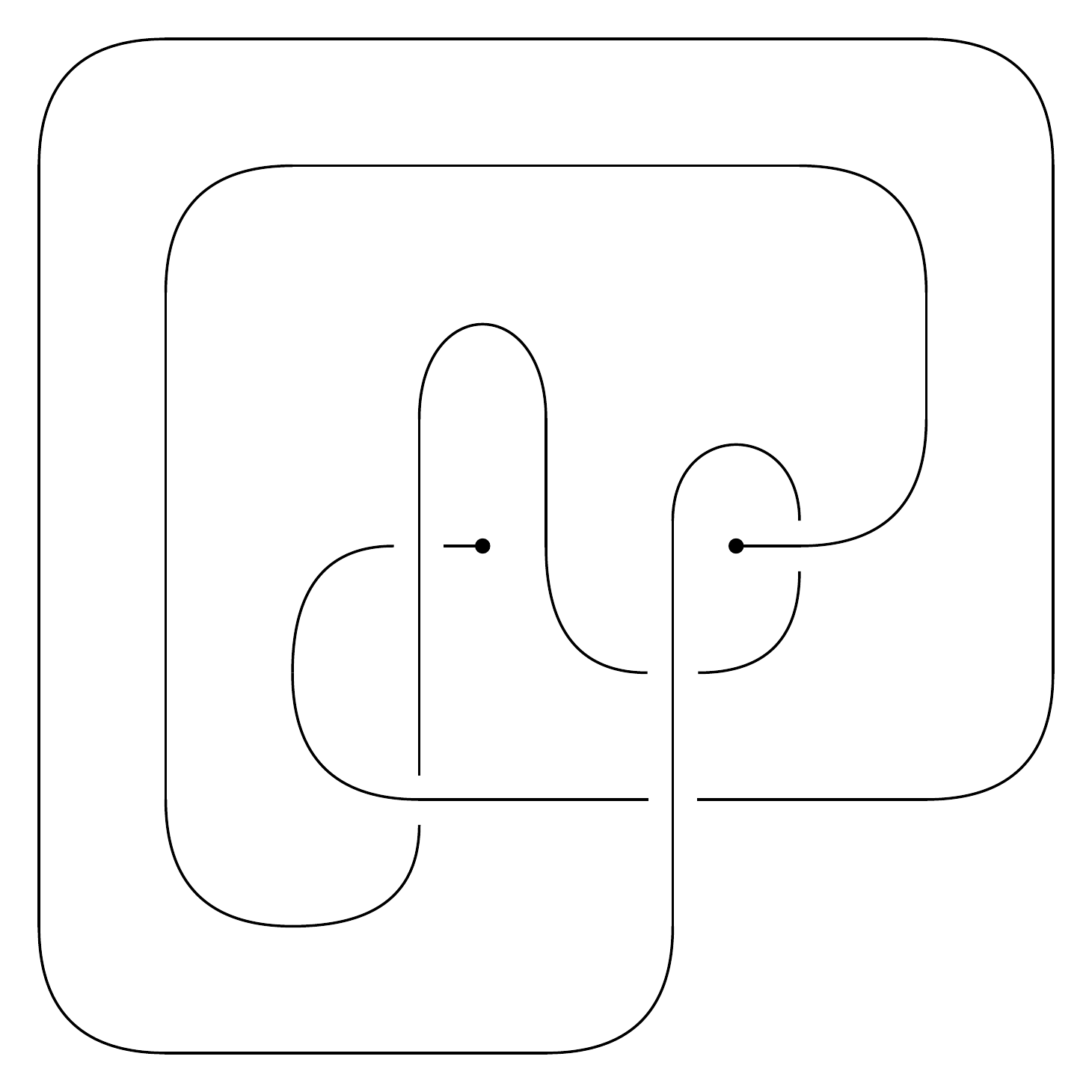}\\
\textcolor{black}{$5_{668}$}
\vspace{1cm}
\end{minipage}
\begin{minipage}[t]{.25\linewidth}
\centering
\includegraphics[width=0.9\textwidth,height=3.5cm,keepaspectratio]{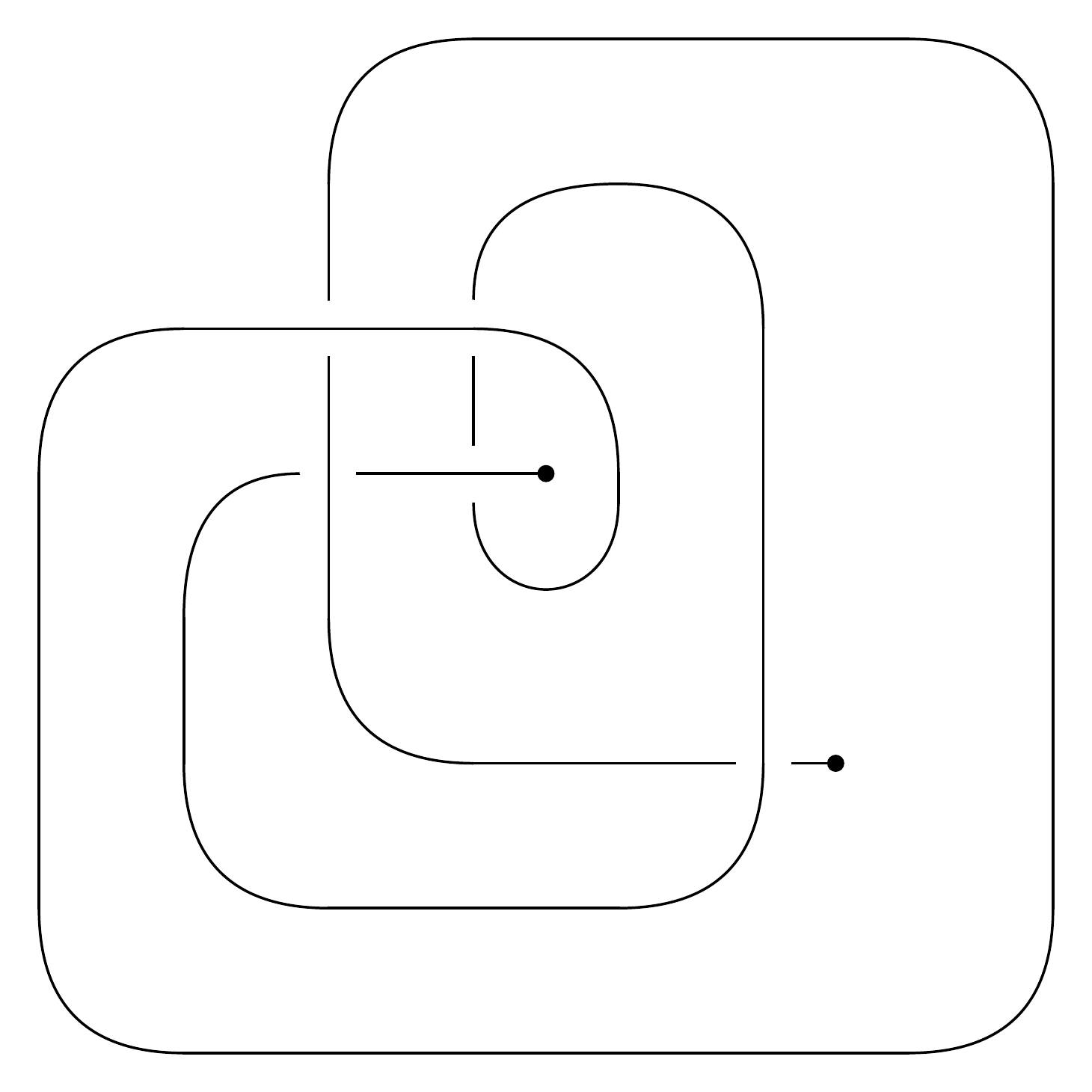}\\
\textcolor{black}{$5_{669}$}
\vspace{1cm}
\end{minipage}
\begin{minipage}[t]{.25\linewidth}
\centering
\includegraphics[width=0.9\textwidth,height=3.5cm,keepaspectratio]{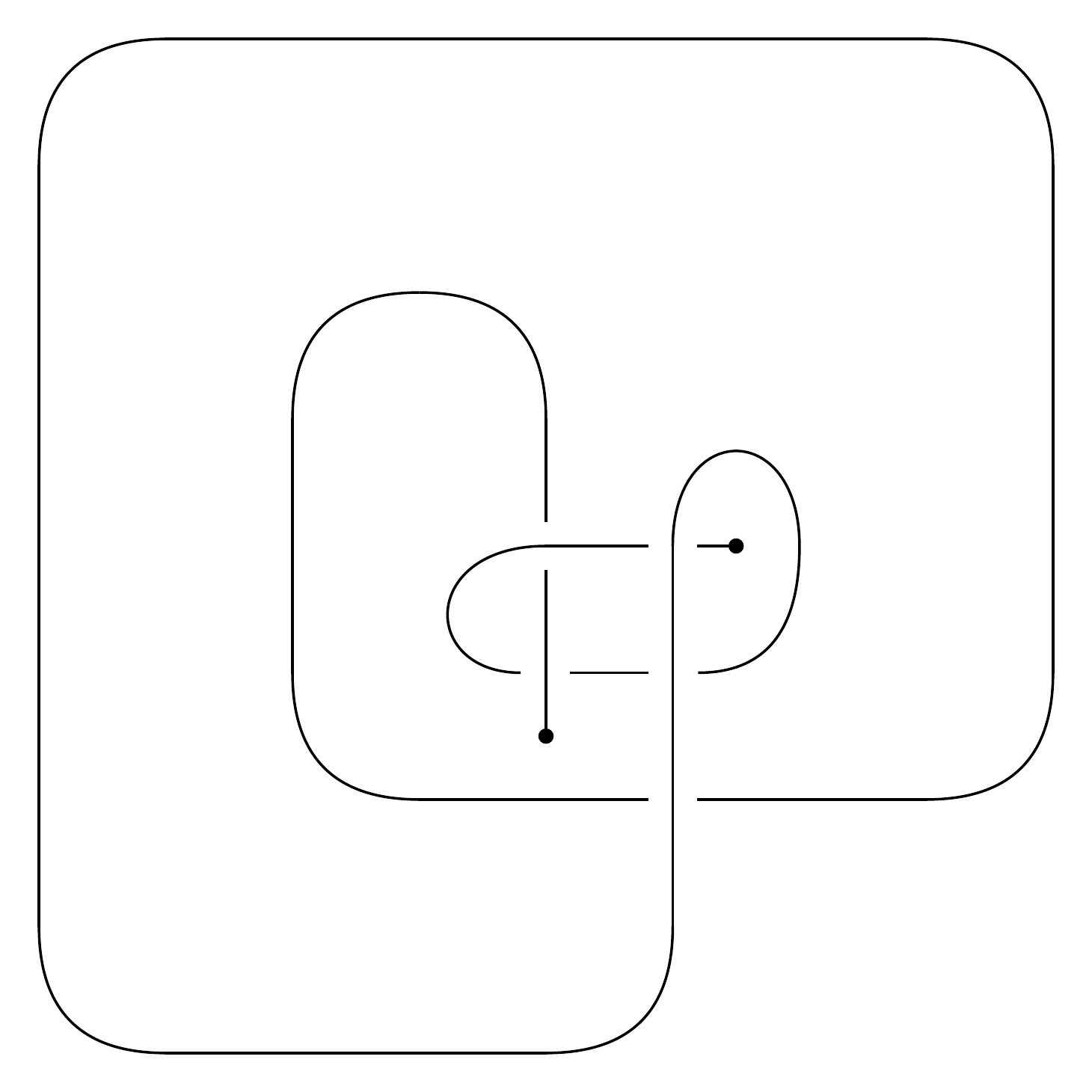}\\
\textcolor{black}{$5_{670}$}
\vspace{1cm}
\end{minipage}
\begin{minipage}[t]{.25\linewidth}
\centering
\includegraphics[width=0.9\textwidth,height=3.5cm,keepaspectratio]{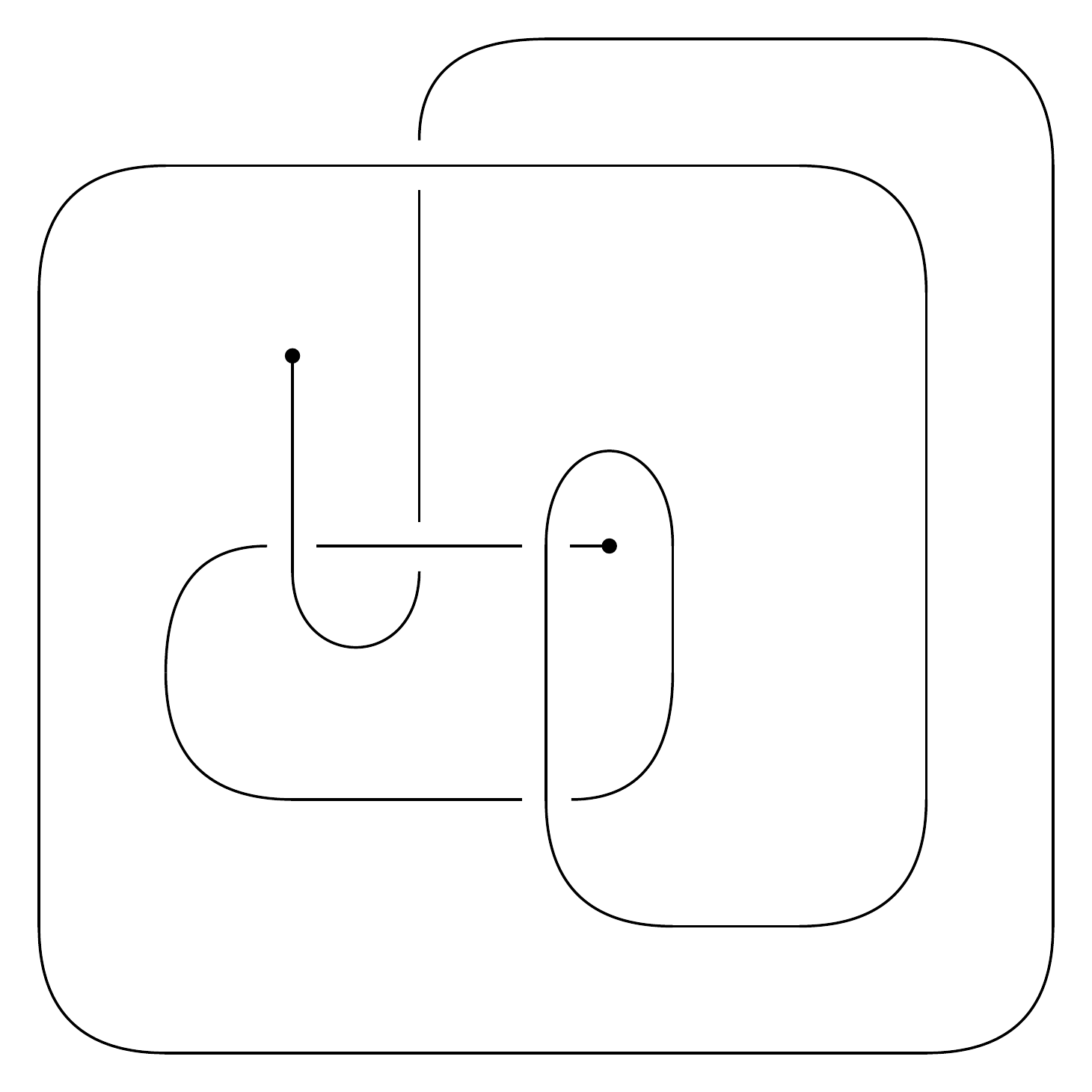}\\
\textcolor{black}{$5_{671}$}
\vspace{1cm}
\end{minipage}
\begin{minipage}[t]{.25\linewidth}
\centering
\includegraphics[width=0.9\textwidth,height=3.5cm,keepaspectratio]{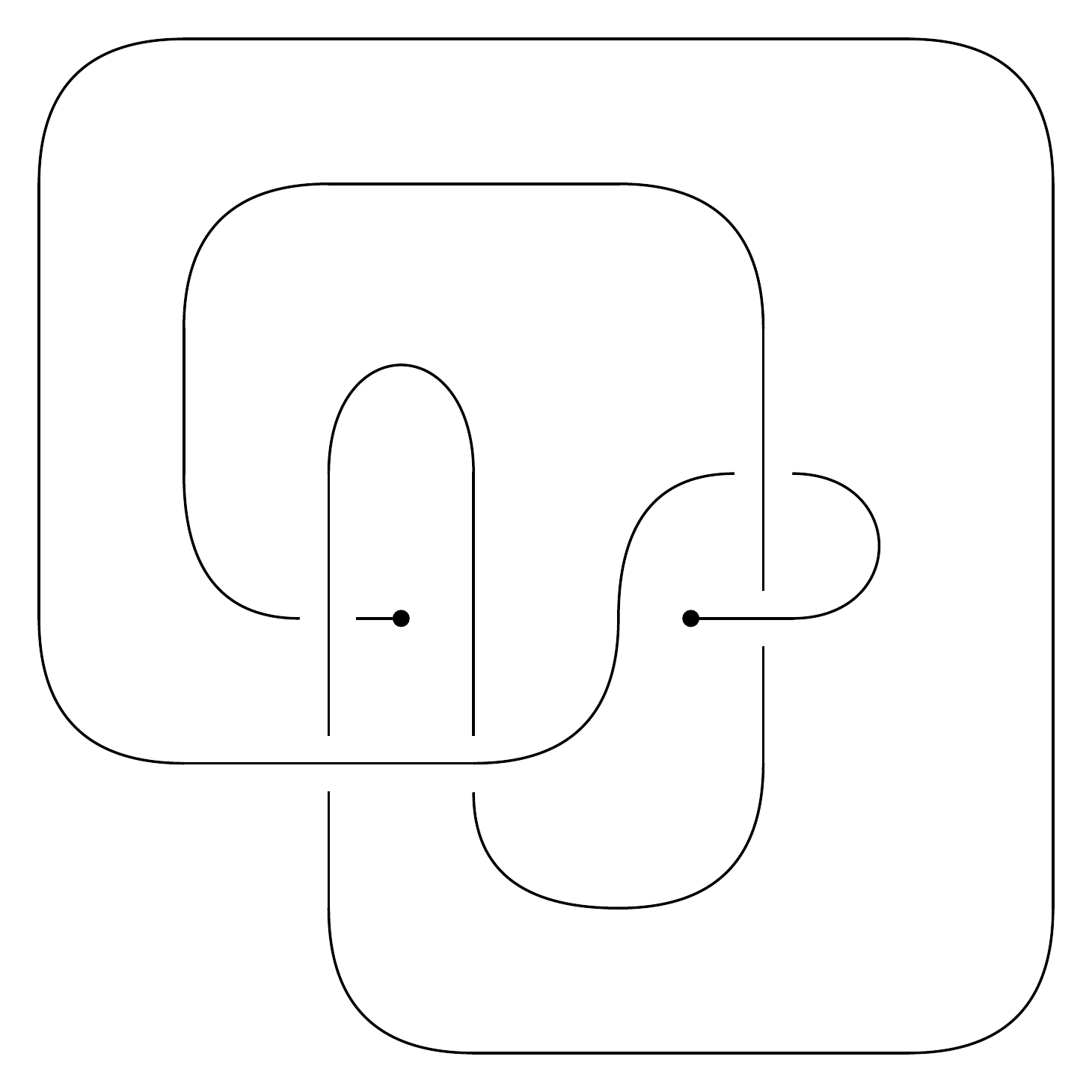}\\
\textcolor{black}{$5_{672}$}
\vspace{1cm}
\end{minipage}
\begin{minipage}[t]{.25\linewidth}
\centering
\includegraphics[width=0.9\textwidth,height=3.5cm,keepaspectratio]{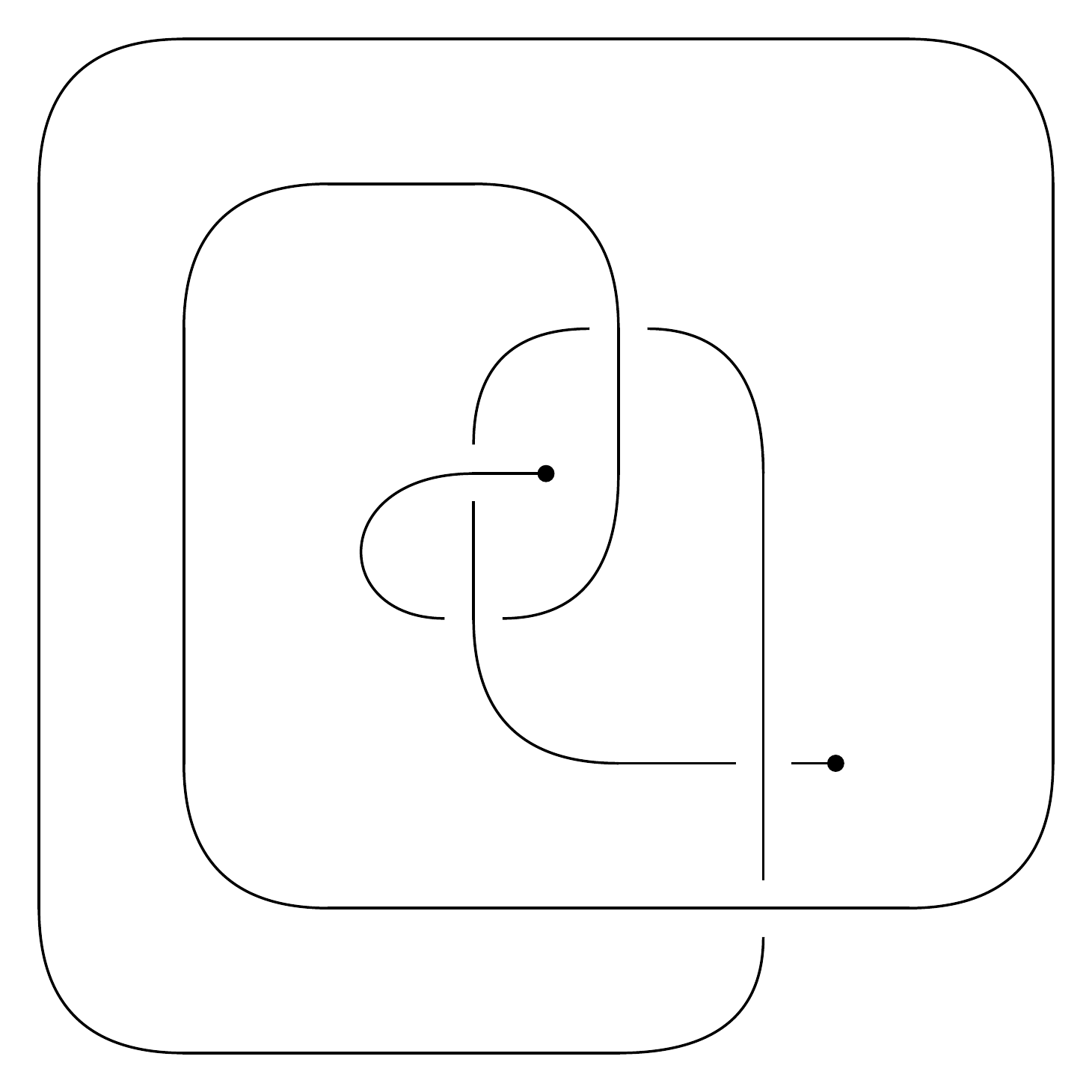}\\
\textcolor{black}{$5_{673}$}
\vspace{1cm}
\end{minipage}
\begin{minipage}[t]{.25\linewidth}
\centering
\includegraphics[width=0.9\textwidth,height=3.5cm,keepaspectratio]{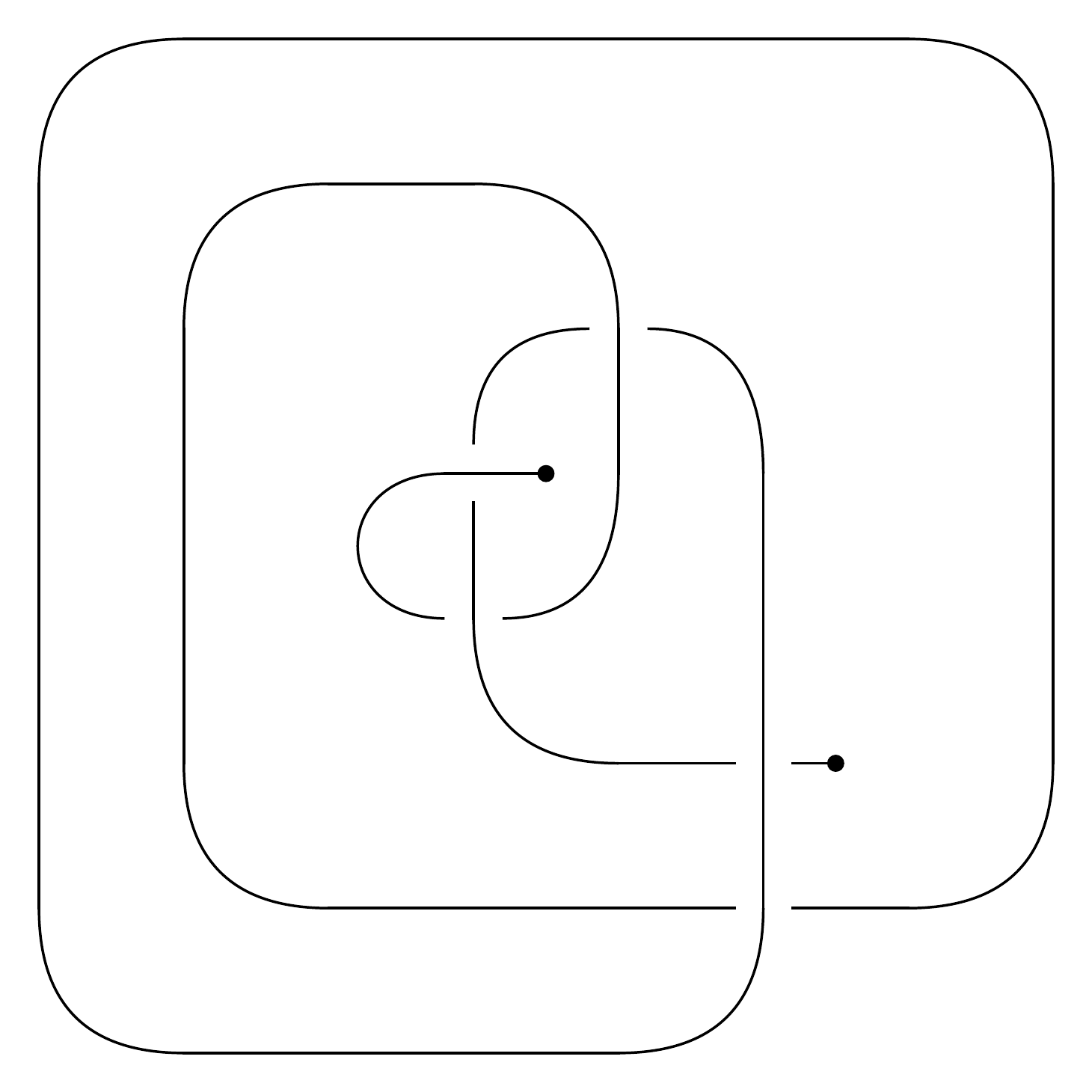}\\
\textcolor{black}{$5_{674}$}
\vspace{1cm}
\end{minipage}
\begin{minipage}[t]{.25\linewidth}
\centering
\includegraphics[width=0.9\textwidth,height=3.5cm,keepaspectratio]{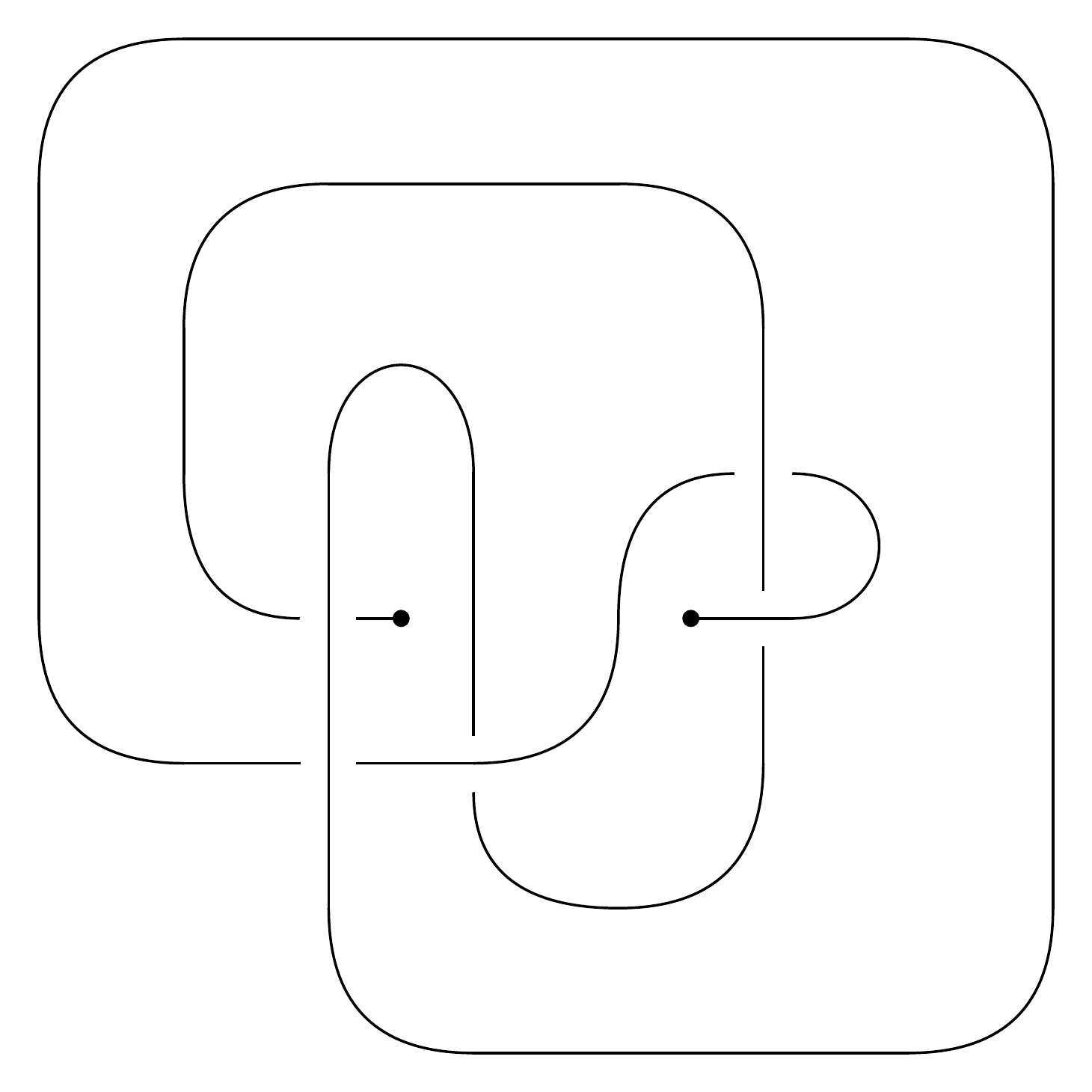}\\
\textcolor{black}{$5_{675}$}
\vspace{1cm}
\end{minipage}
\begin{minipage}[t]{.25\linewidth}
\centering
\includegraphics[width=0.9\textwidth,height=3.5cm,keepaspectratio]{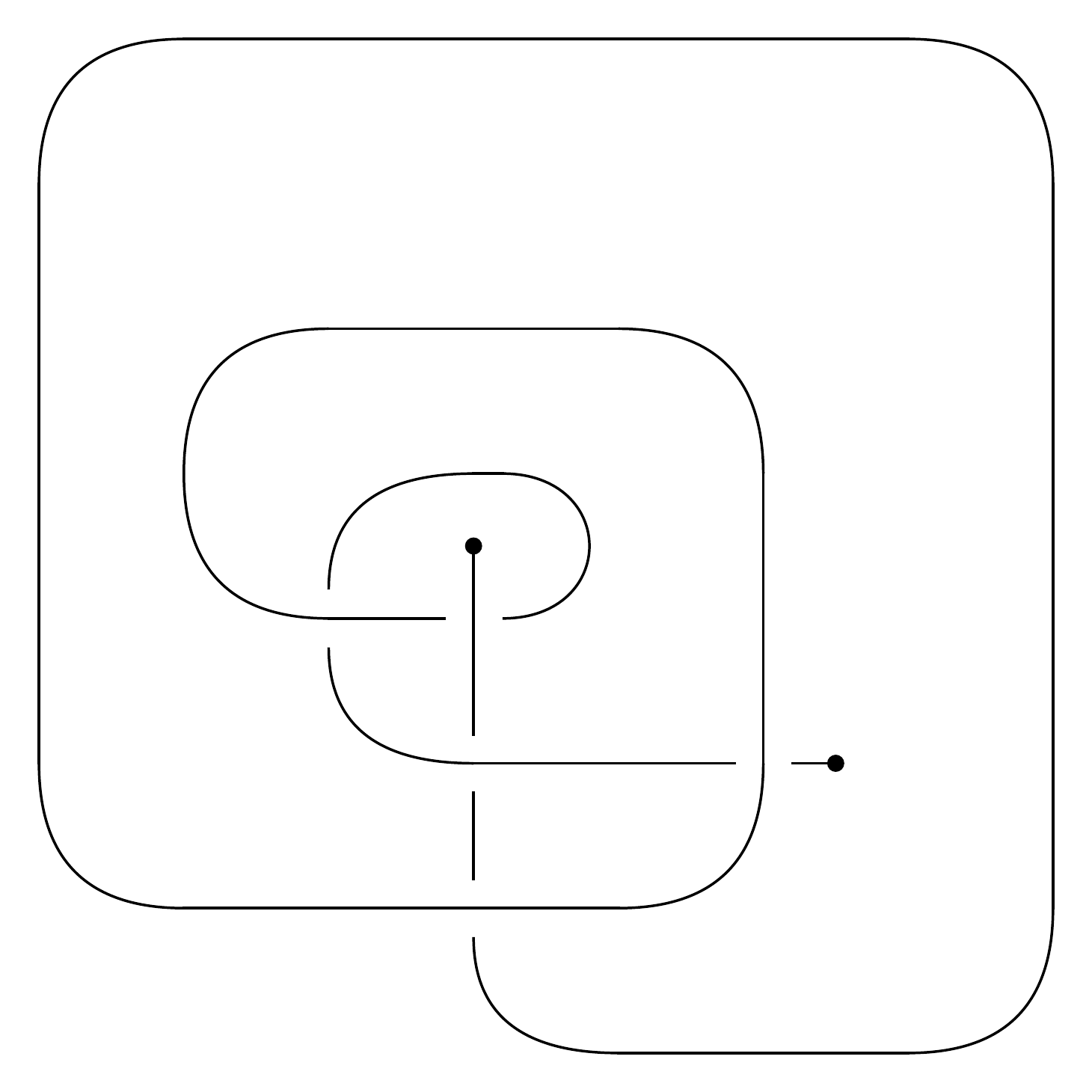}\\
\textcolor{black}{$5_{676}$}
\vspace{1cm}
\end{minipage}
\begin{minipage}[t]{.25\linewidth}
\centering
\includegraphics[width=0.9\textwidth,height=3.5cm,keepaspectratio]{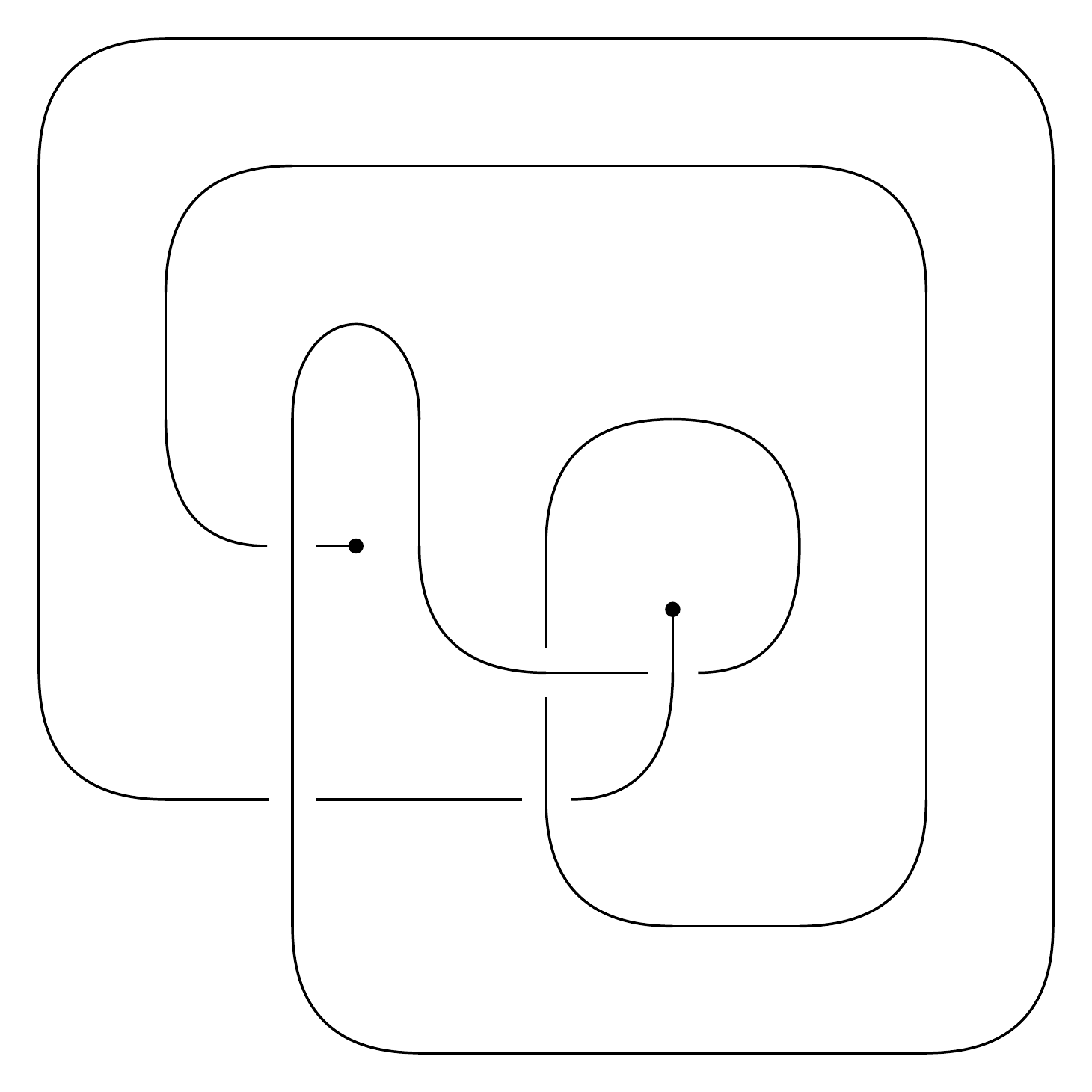}\\
\textcolor{black}{$5_{677}$}
\vspace{1cm}
\end{minipage}
\begin{minipage}[t]{.25\linewidth}
\centering
\includegraphics[width=0.9\textwidth,height=3.5cm,keepaspectratio]{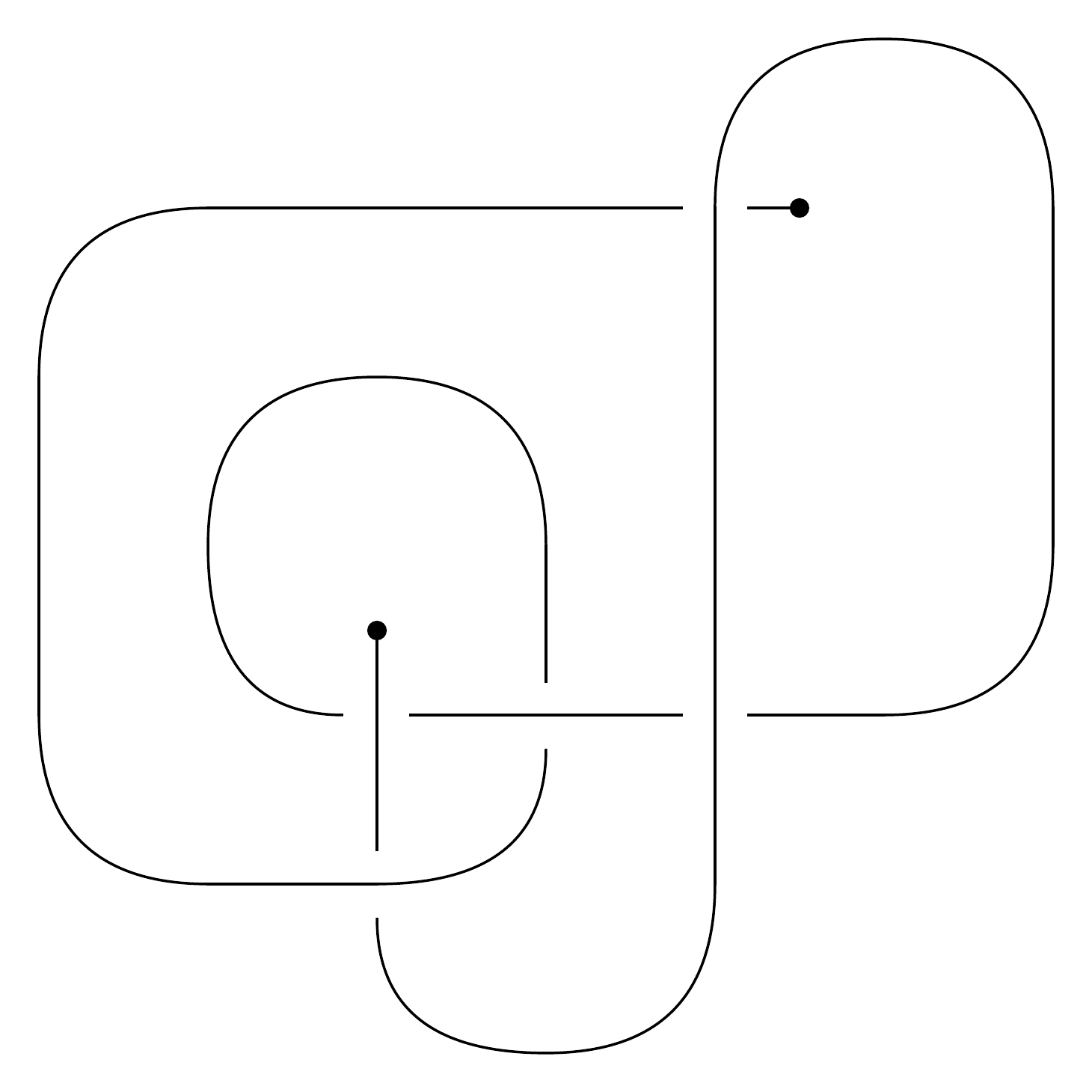}\\
\textcolor{black}{$5_{678}$}
\vspace{1cm}
\end{minipage}
\begin{minipage}[t]{.25\linewidth}
\centering
\includegraphics[width=0.9\textwidth,height=3.5cm,keepaspectratio]{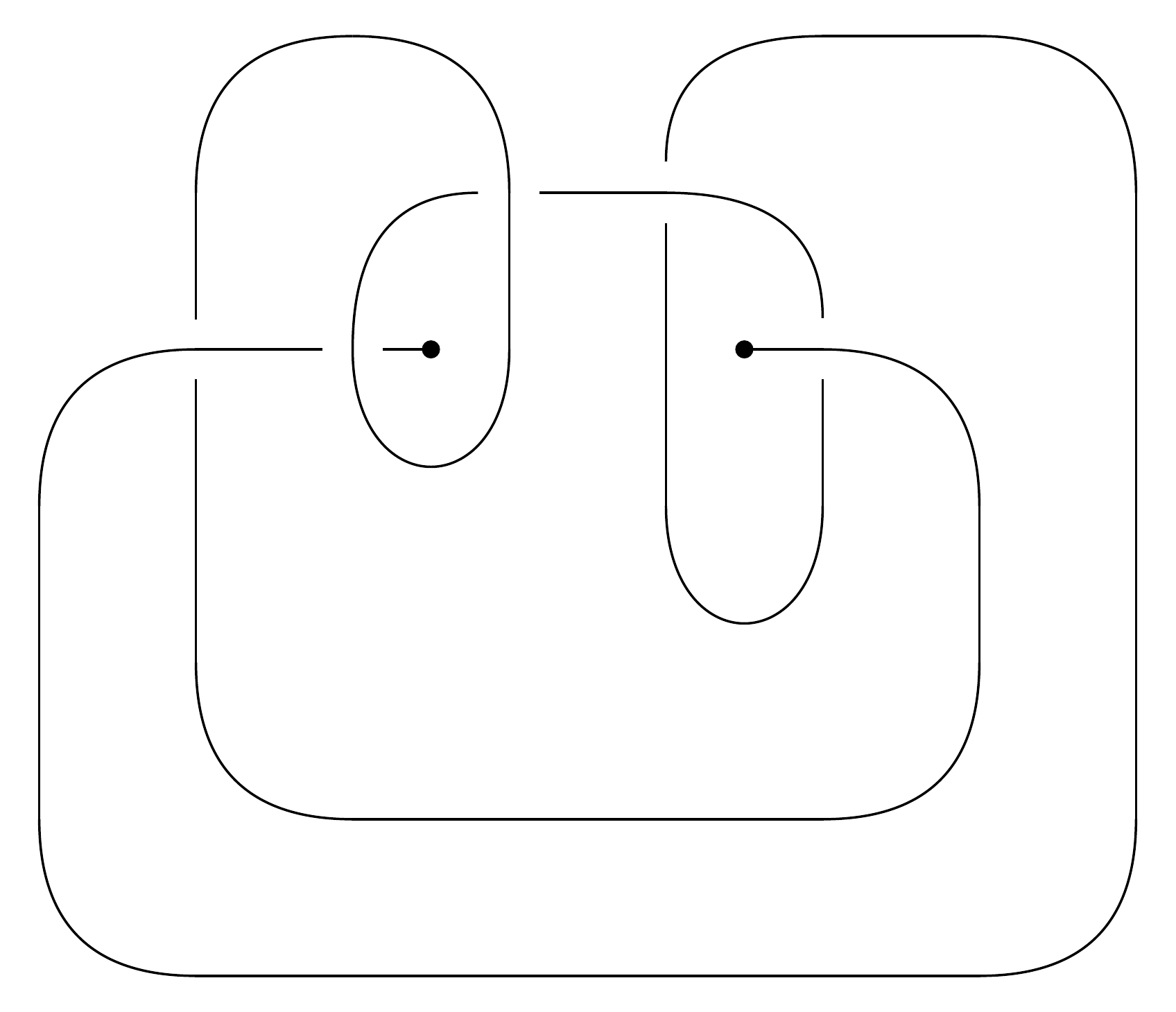}\\
\textcolor{black}{$5_{679}$}
\vspace{1cm}
\end{minipage}
\begin{minipage}[t]{.25\linewidth}
\centering
\includegraphics[width=0.9\textwidth,height=3.5cm,keepaspectratio]{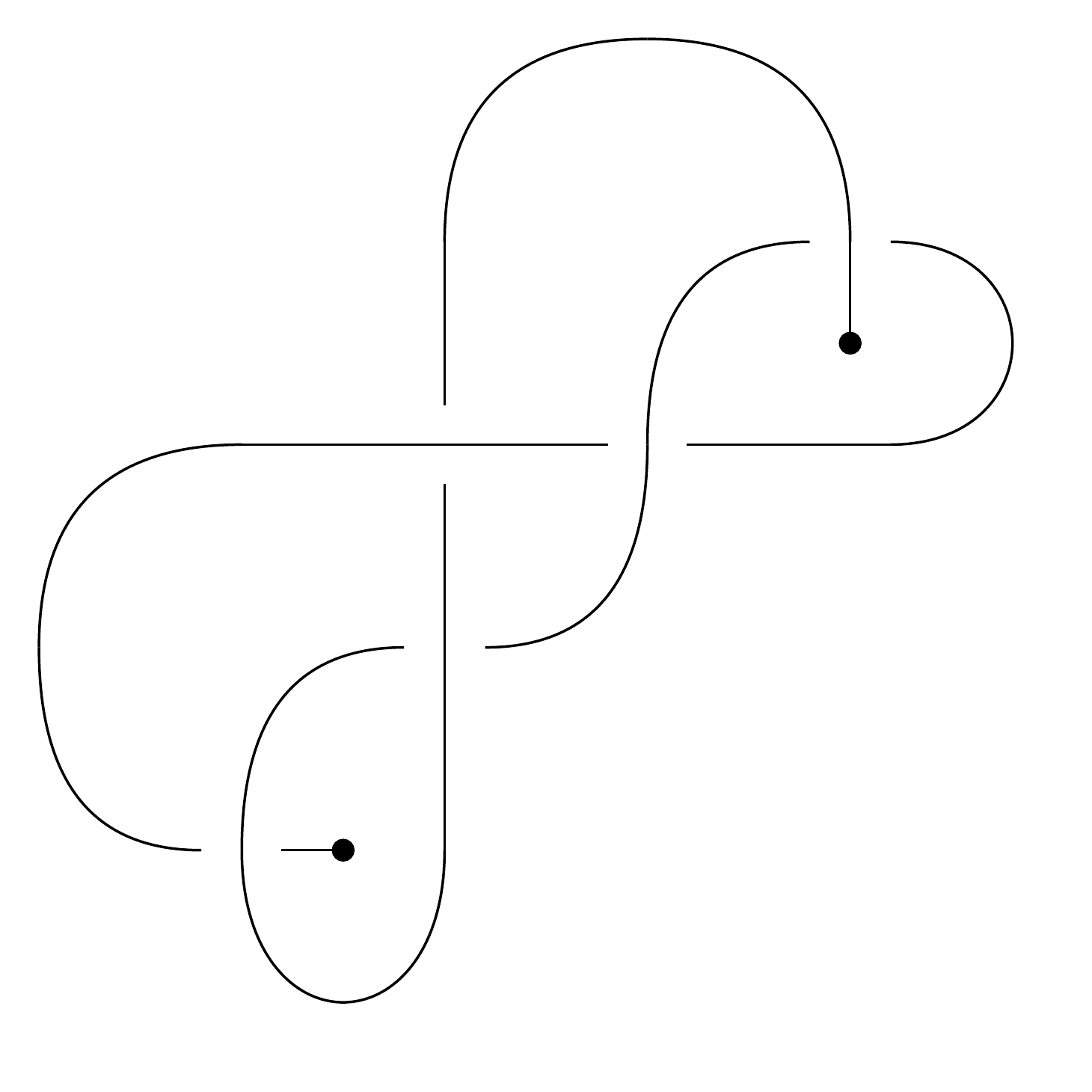}\\
\textcolor{black}{$5_{680}$}
\vspace{1cm}
\end{minipage}
\begin{minipage}[t]{.25\linewidth}
\centering
\includegraphics[width=0.9\textwidth,height=3.5cm,keepaspectratio]{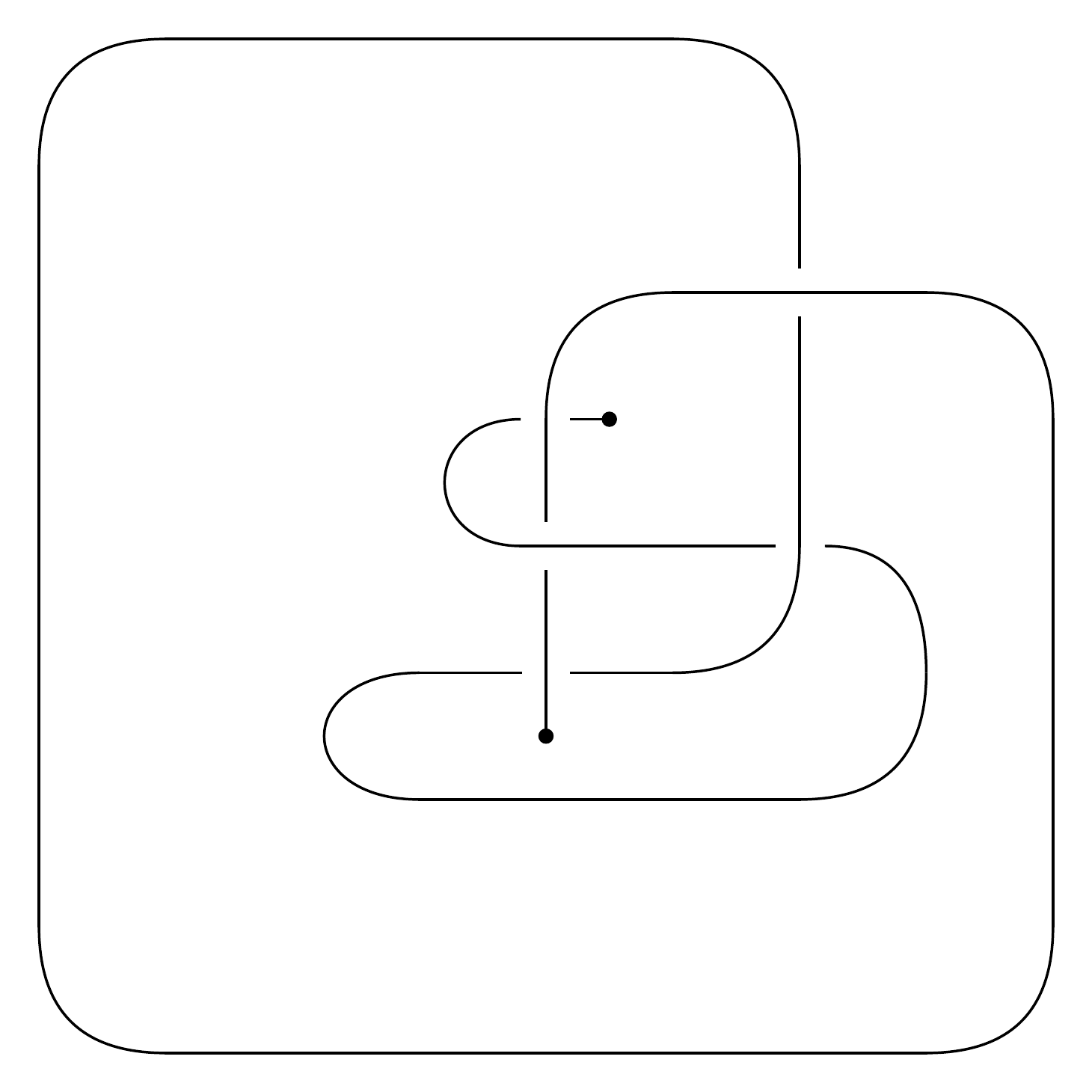}\\
\textcolor{black}{$5_{681}$}
\vspace{1cm}
\end{minipage}
\begin{minipage}[t]{.25\linewidth}
\centering
\includegraphics[width=0.9\textwidth,height=3.5cm,keepaspectratio]{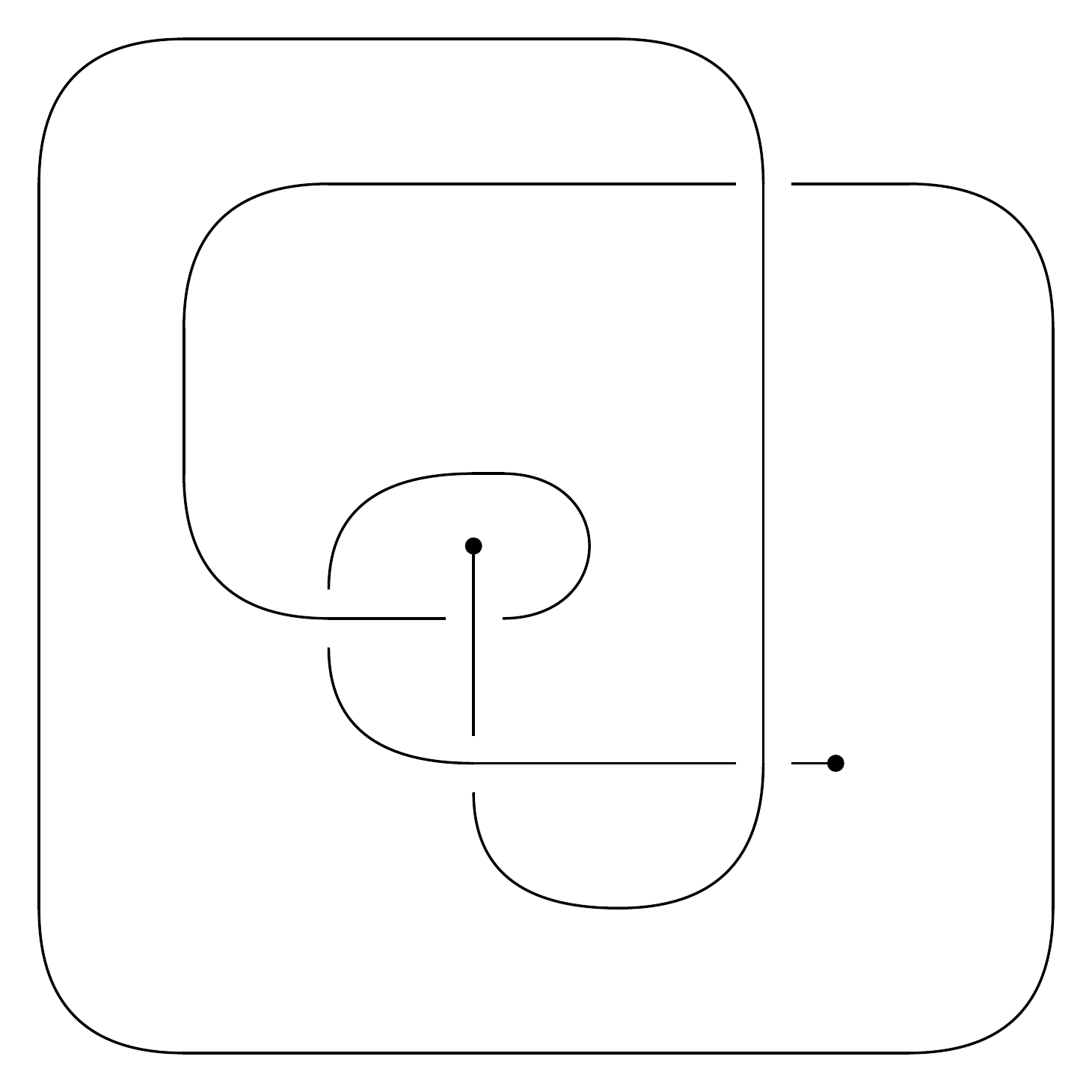}\\
\textcolor{black}{$5_{682}$}
\vspace{1cm}
\end{minipage}
\begin{minipage}[t]{.25\linewidth}
\centering
\includegraphics[width=0.9\textwidth,height=3.5cm,keepaspectratio]{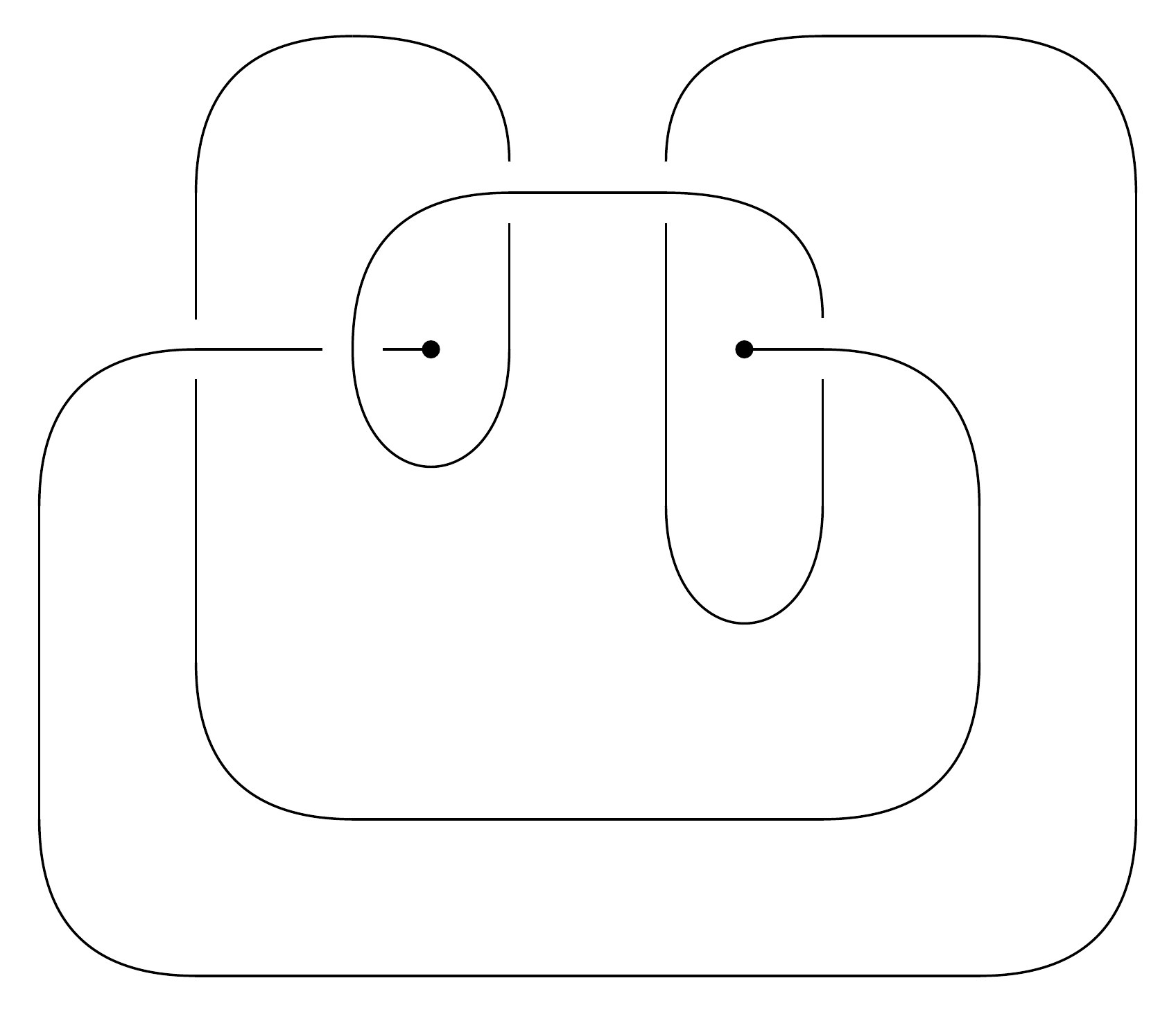}\\
\textcolor{black}{$5_{683}$}
\vspace{1cm}
\end{minipage}
\begin{minipage}[t]{.25\linewidth}
\centering
\includegraphics[width=0.9\textwidth,height=3.5cm,keepaspectratio]{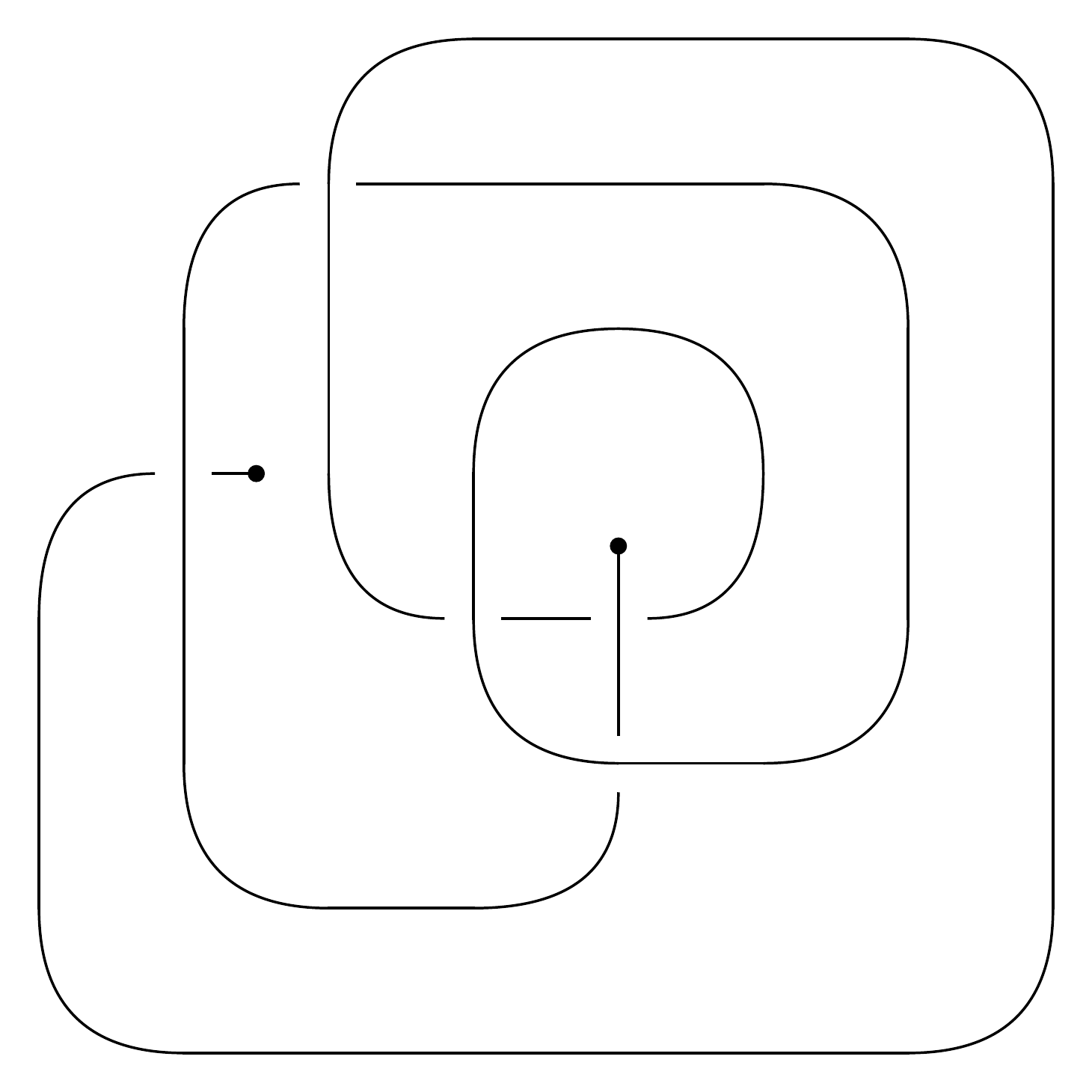}\\
\textcolor{black}{$5_{684}$}
\vspace{1cm}
\end{minipage}
\begin{minipage}[t]{.25\linewidth}
\centering
\includegraphics[width=0.9\textwidth,height=3.5cm,keepaspectratio]{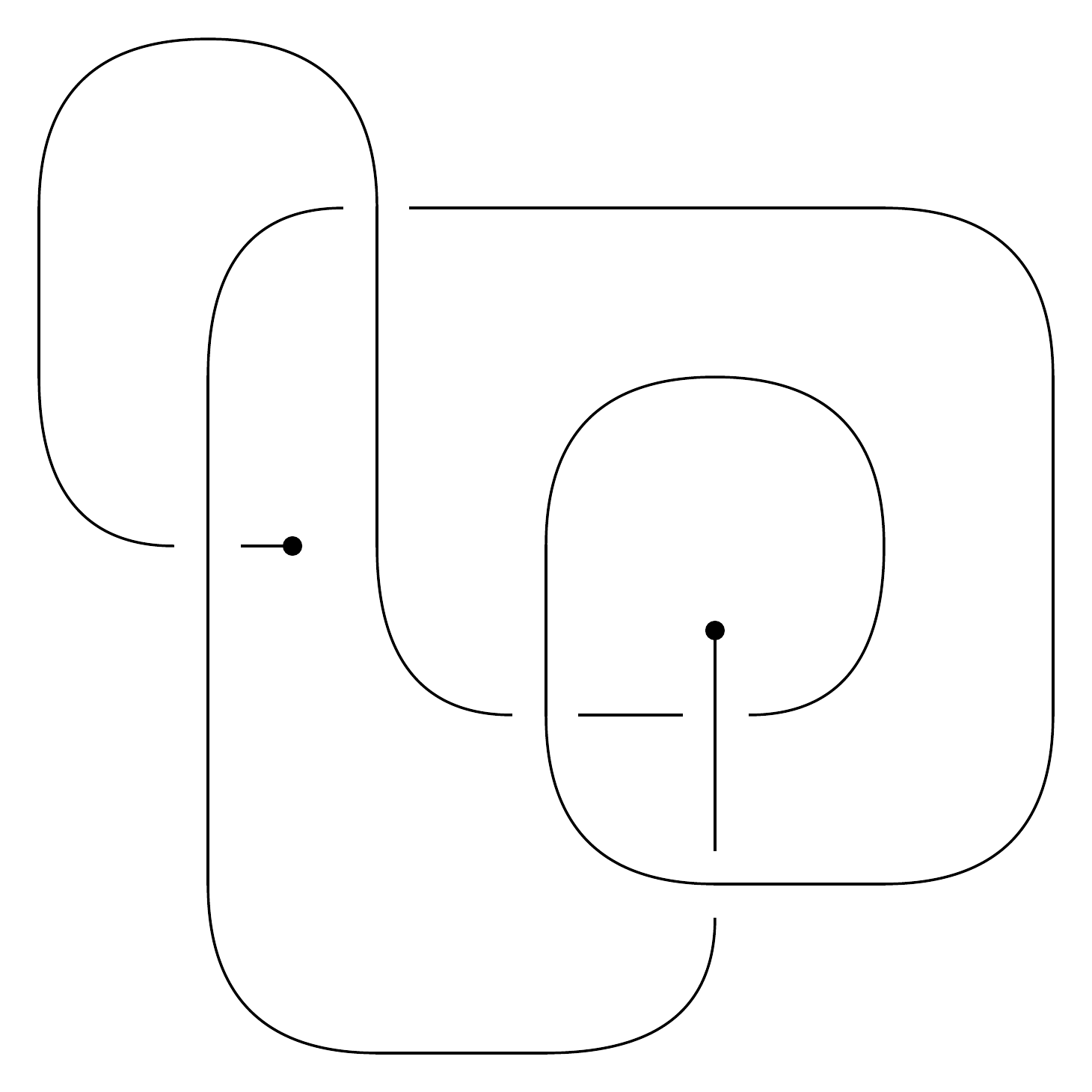}\\
\textcolor{black}{$5_{685}$}
\vspace{1cm}
\end{minipage}
\begin{minipage}[t]{.25\linewidth}
\centering
\includegraphics[width=0.9\textwidth,height=3.5cm,keepaspectratio]{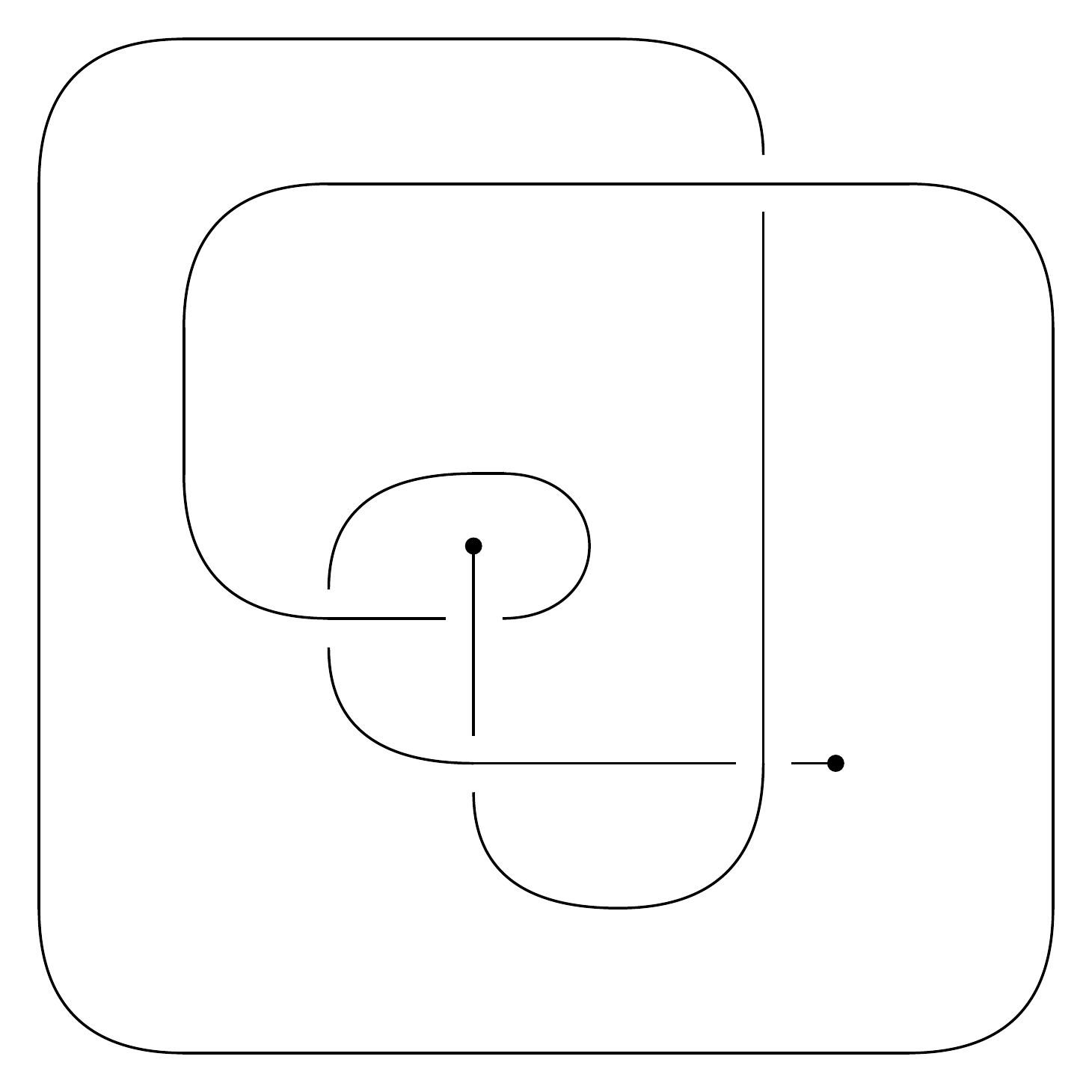}\\
\textcolor{black}{$5_{686}$}
\vspace{1cm}
\end{minipage}
\begin{minipage}[t]{.25\linewidth}
\centering
\includegraphics[width=0.9\textwidth,height=3.5cm,keepaspectratio]{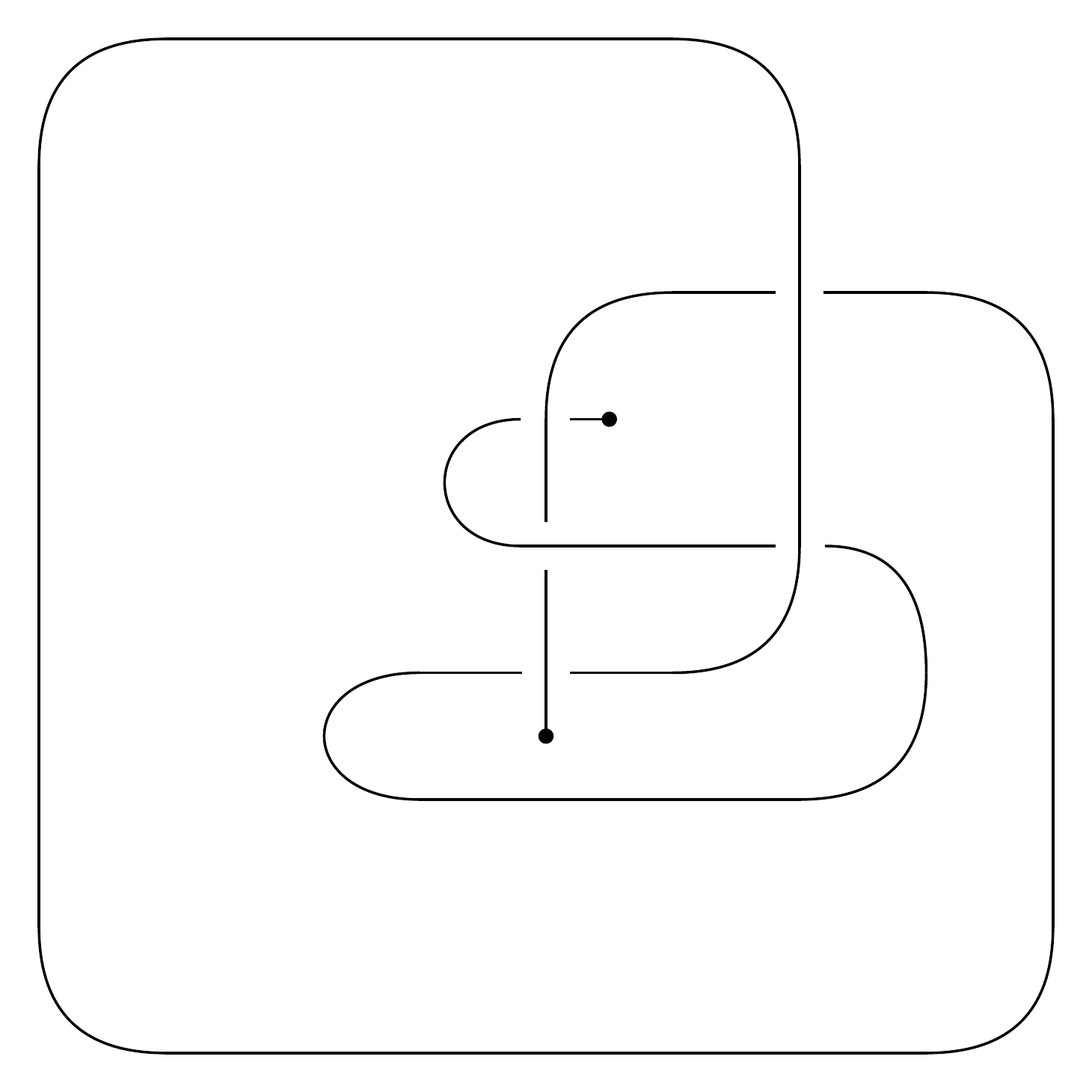}\\
\textcolor{black}{$5_{687}$}
\vspace{1cm}
\end{minipage}
\begin{minipage}[t]{.25\linewidth}
\centering
\includegraphics[width=0.9\textwidth,height=3.5cm,keepaspectratio]{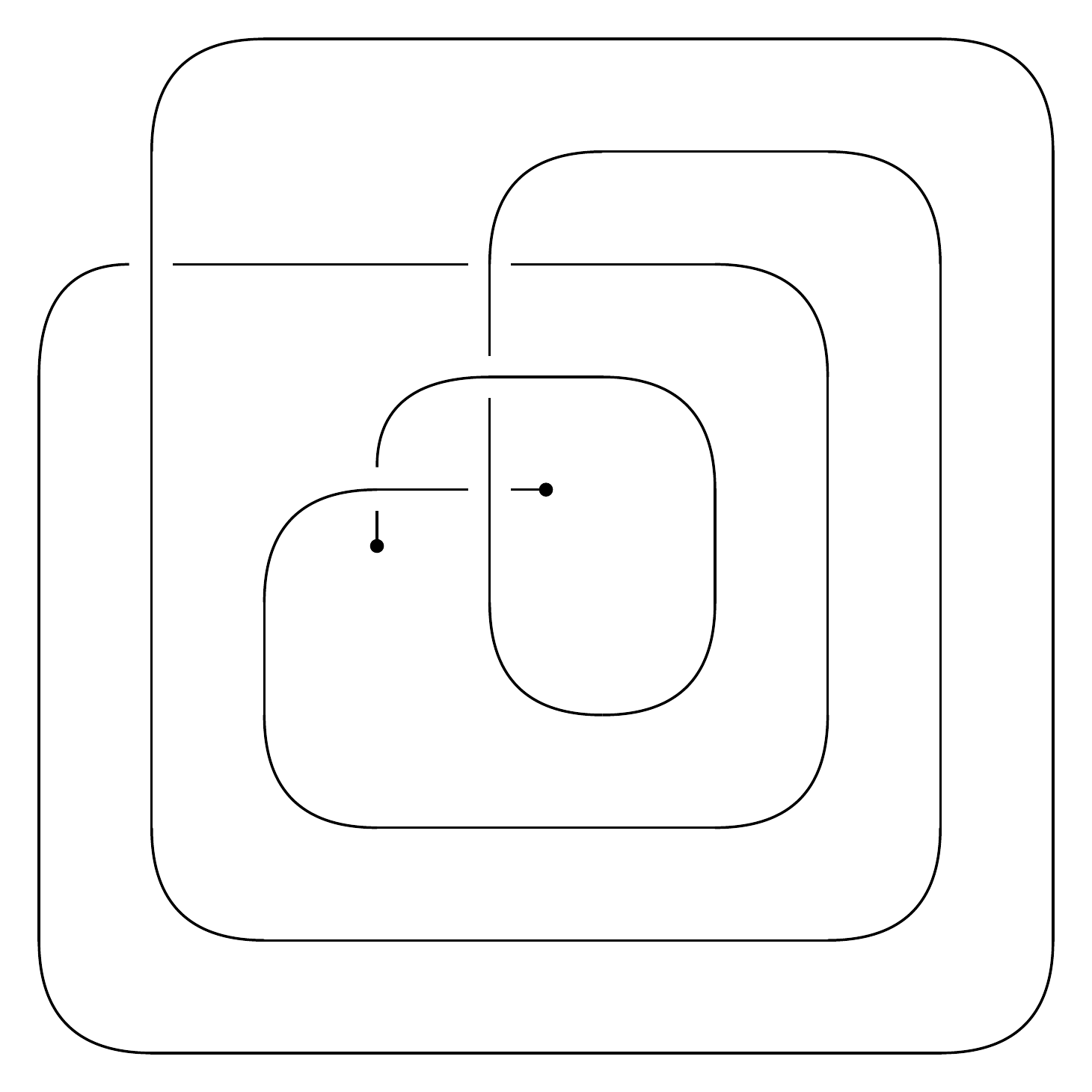}\\
\textcolor{black}{$5_{688}$}
\vspace{1cm}
\end{minipage}
\begin{minipage}[t]{.25\linewidth}
\centering
\includegraphics[width=0.9\textwidth,height=3.5cm,keepaspectratio]{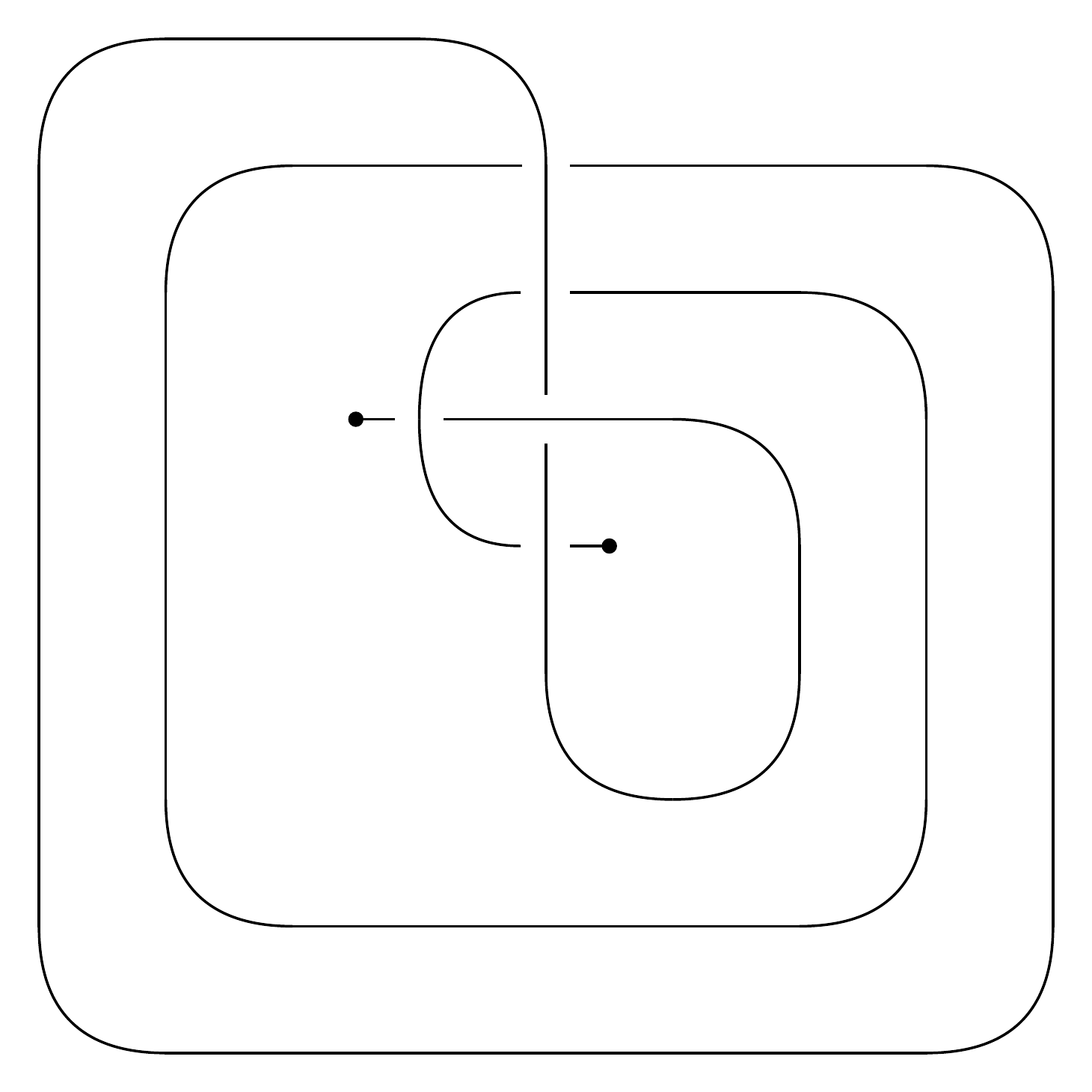}\\
\textcolor{black}{$5_{689}$}
\vspace{1cm}
\end{minipage}
\begin{minipage}[t]{.25\linewidth}
\centering
\includegraphics[width=0.9\textwidth,height=3.5cm,keepaspectratio]{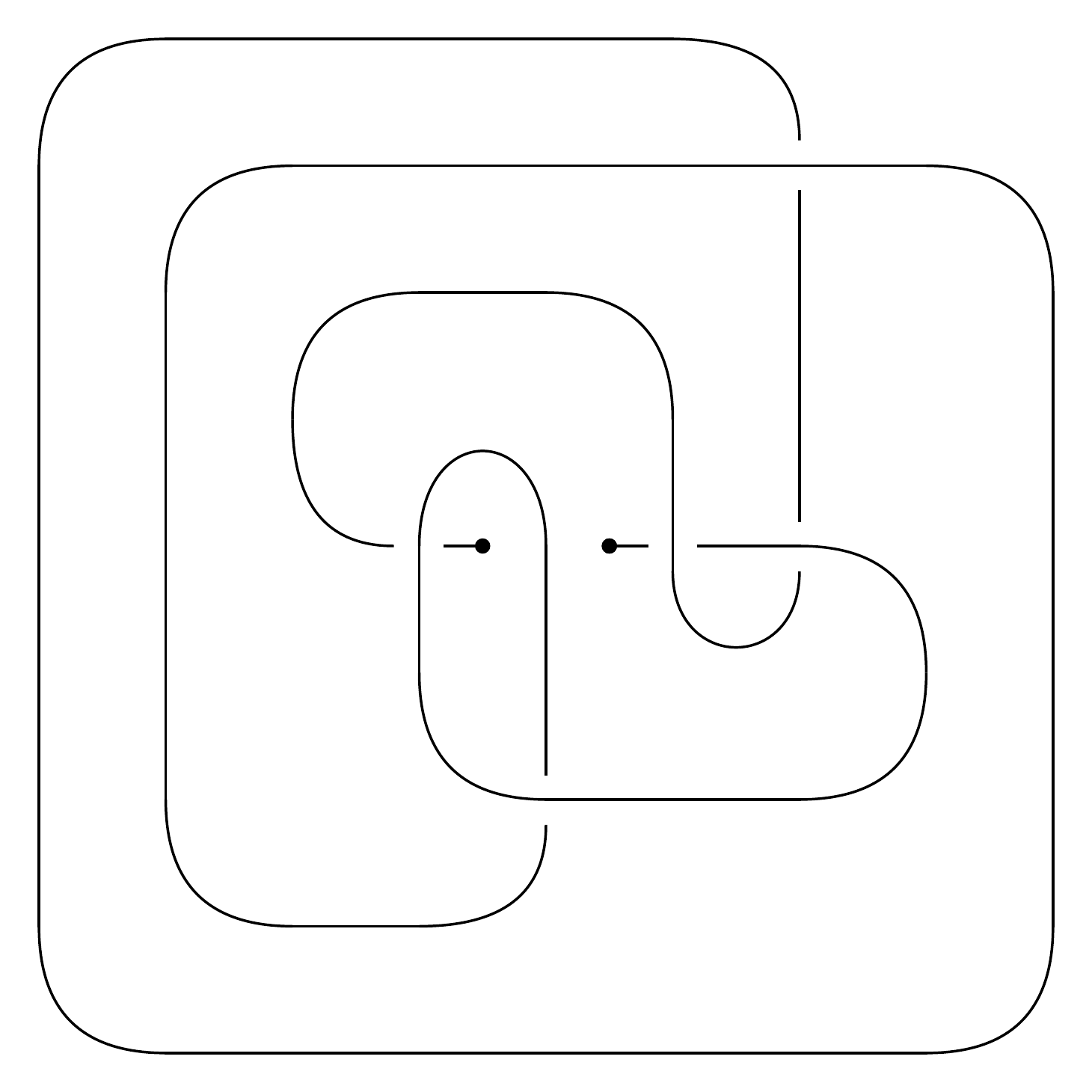}\\
\textcolor{black}{$5_{690}$}
\vspace{1cm}
\end{minipage}
\begin{minipage}[t]{.25\linewidth}
\centering
\includegraphics[width=0.9\textwidth,height=3.5cm,keepaspectratio]{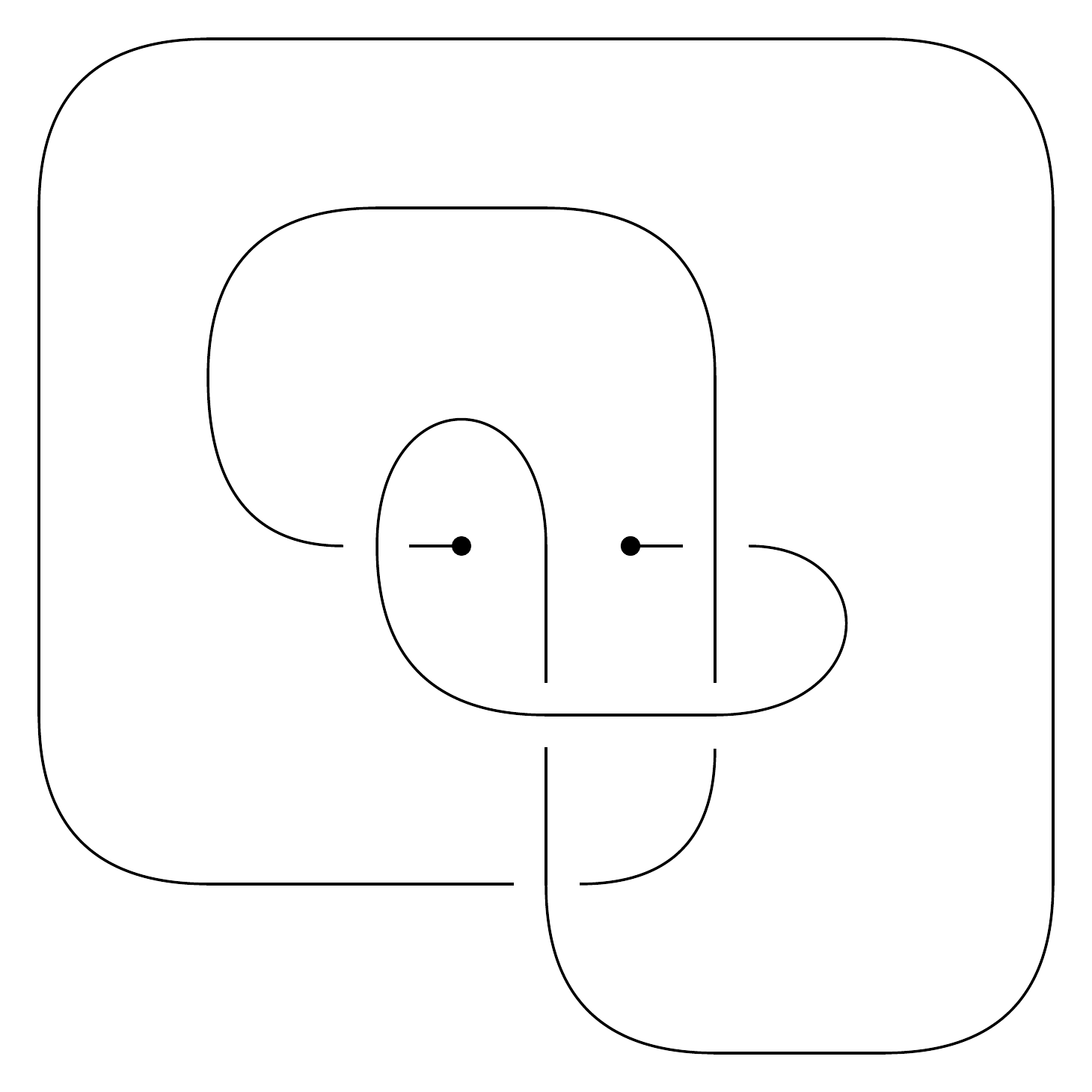}\\
\textcolor{black}{$5_{691}$}
\vspace{1cm}
\end{minipage}
\begin{minipage}[t]{.25\linewidth}
\centering
\includegraphics[width=0.9\textwidth,height=3.5cm,keepaspectratio]{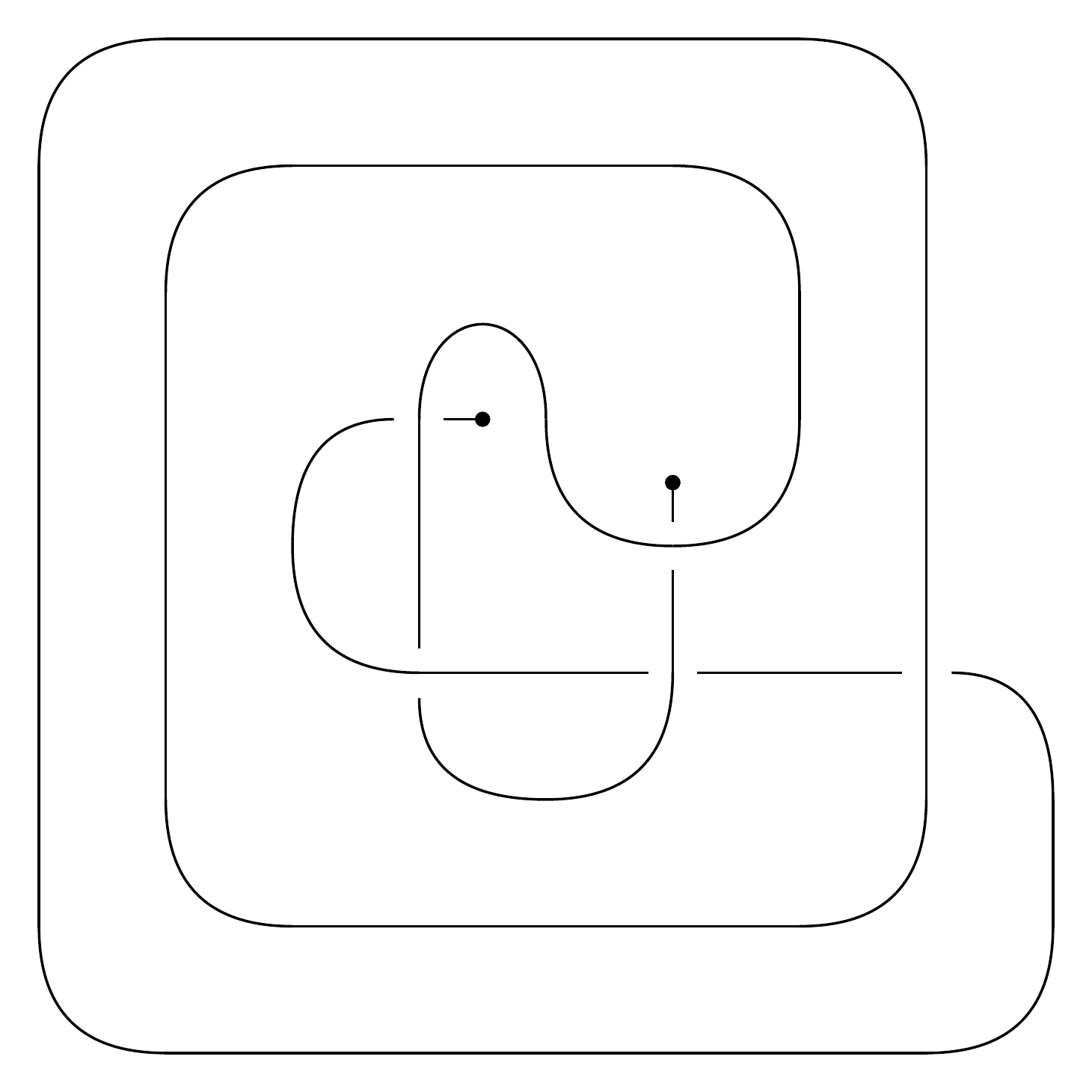}\\
\textcolor{black}{$5_{692}$}
\vspace{1cm}
\end{minipage}
\begin{minipage}[t]{.25\linewidth}
\centering
\includegraphics[width=0.9\textwidth,height=3.5cm,keepaspectratio]{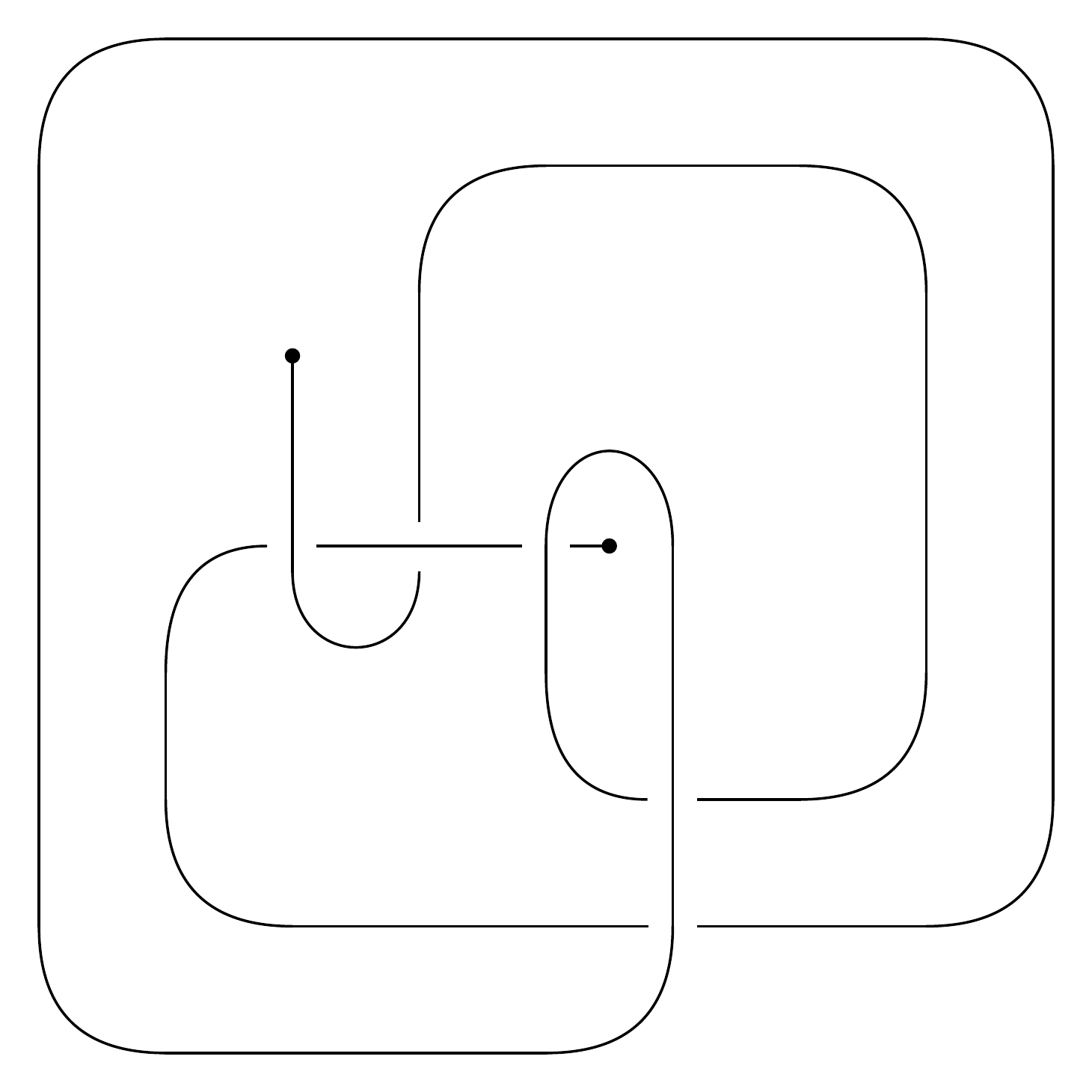}\\
\textcolor{black}{$5_{693}$}
\vspace{1cm}
\end{minipage}
\begin{minipage}[t]{.25\linewidth}
\centering
\includegraphics[width=0.9\textwidth,height=3.5cm,keepaspectratio]{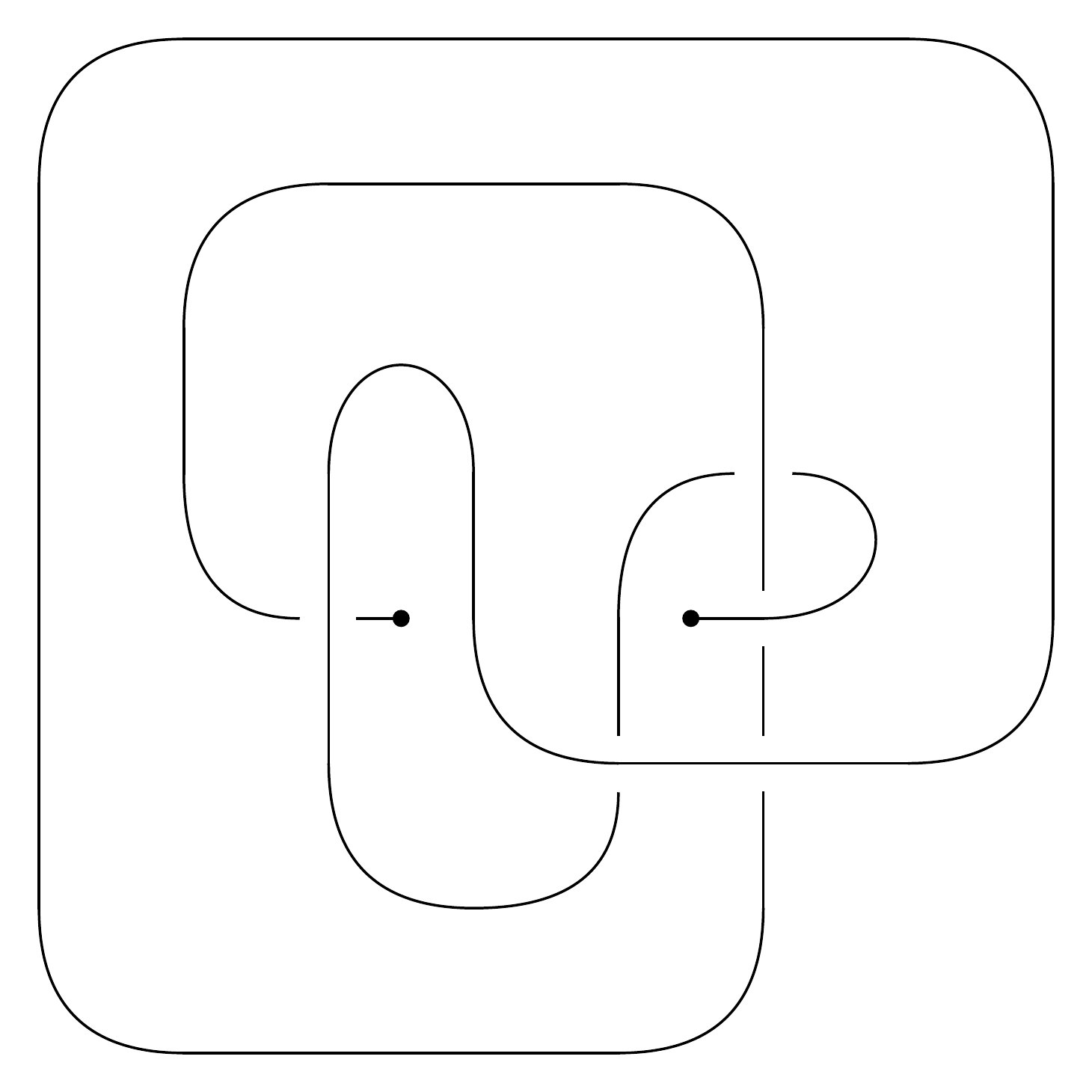}\\
\textcolor{black}{$5_{694}$}
\vspace{1cm}
\end{minipage}
\begin{minipage}[t]{.25\linewidth}
\centering
\includegraphics[width=0.9\textwidth,height=3.5cm,keepaspectratio]{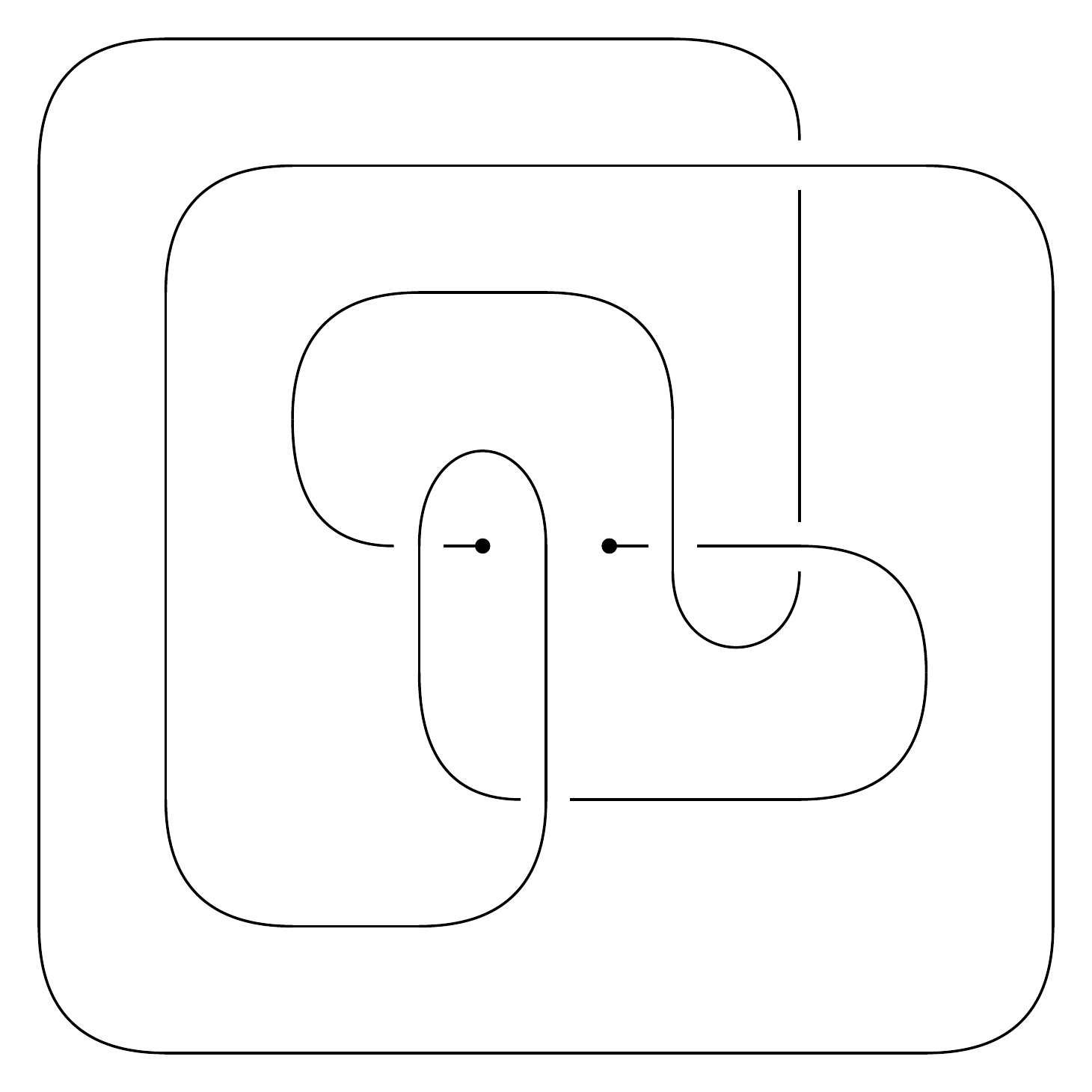}\\
\textcolor{black}{$5_{695}$}
\vspace{1cm}
\end{minipage}
\begin{minipage}[t]{.25\linewidth}
\centering
\includegraphics[width=0.9\textwidth,height=3.5cm,keepaspectratio]{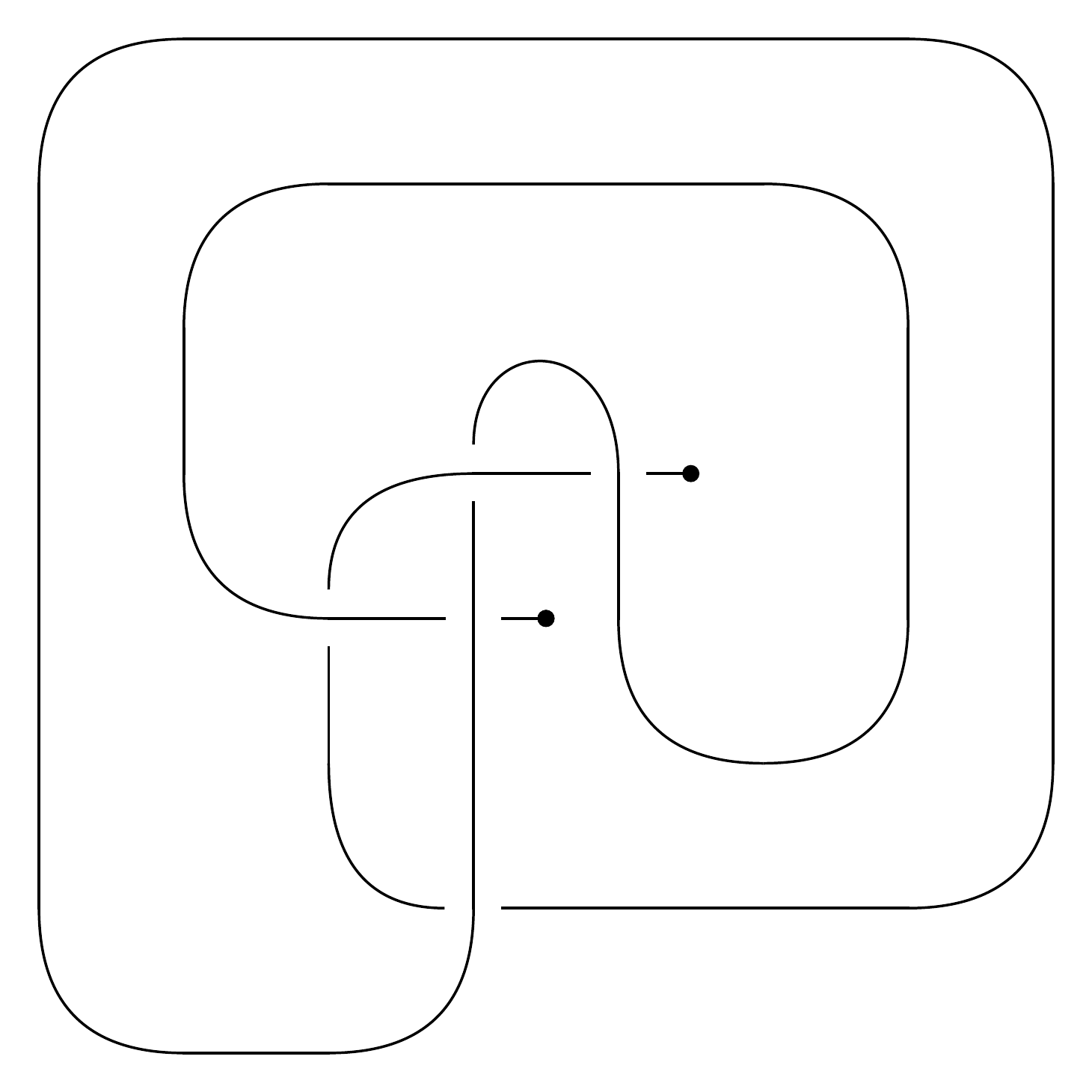}\\
\textcolor{black}{$5_{696}$}
\vspace{1cm}
\end{minipage}
\begin{minipage}[t]{.25\linewidth}
\centering
\includegraphics[width=0.9\textwidth,height=3.5cm,keepaspectratio]{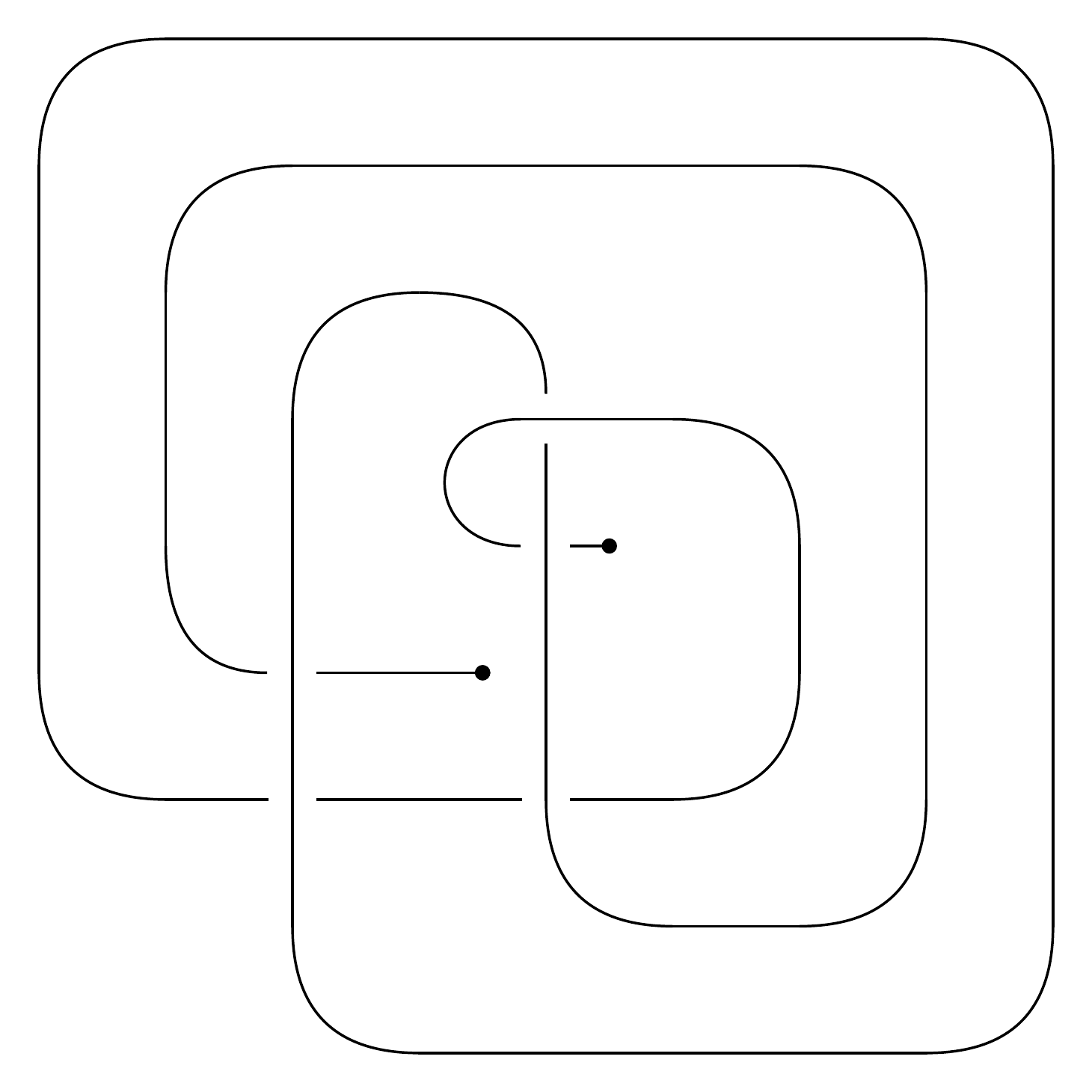}\\
\textcolor{black}{$5_{697}$}
\vspace{1cm}
\end{minipage}
\begin{minipage}[t]{.25\linewidth}
\centering
\includegraphics[width=0.9\textwidth,height=3.5cm,keepaspectratio]{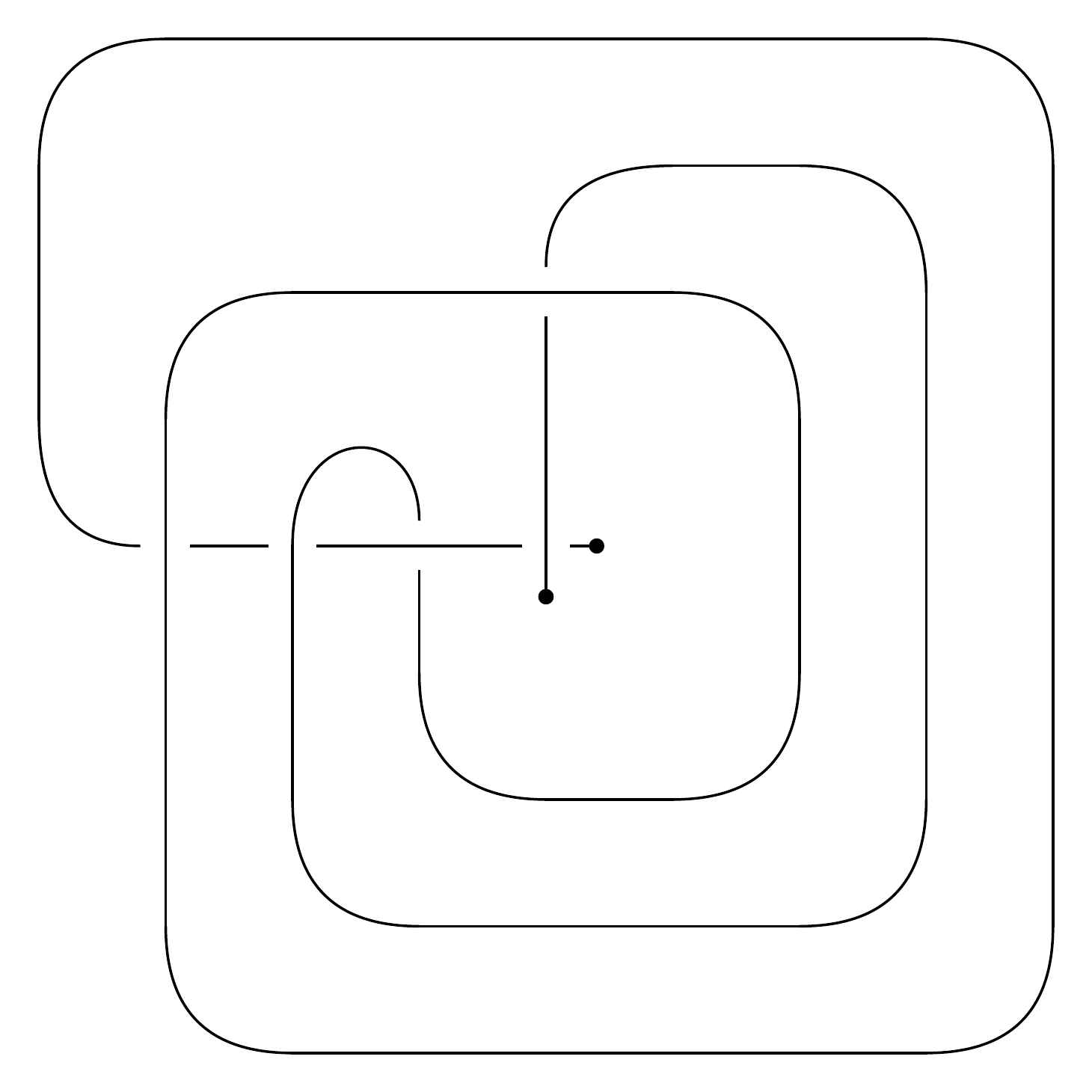}\\
\textcolor{black}{$5_{698}$}
\vspace{1cm}
\end{minipage}
\begin{minipage}[t]{.25\linewidth}
\centering
\includegraphics[width=0.9\textwidth,height=3.5cm,keepaspectratio]{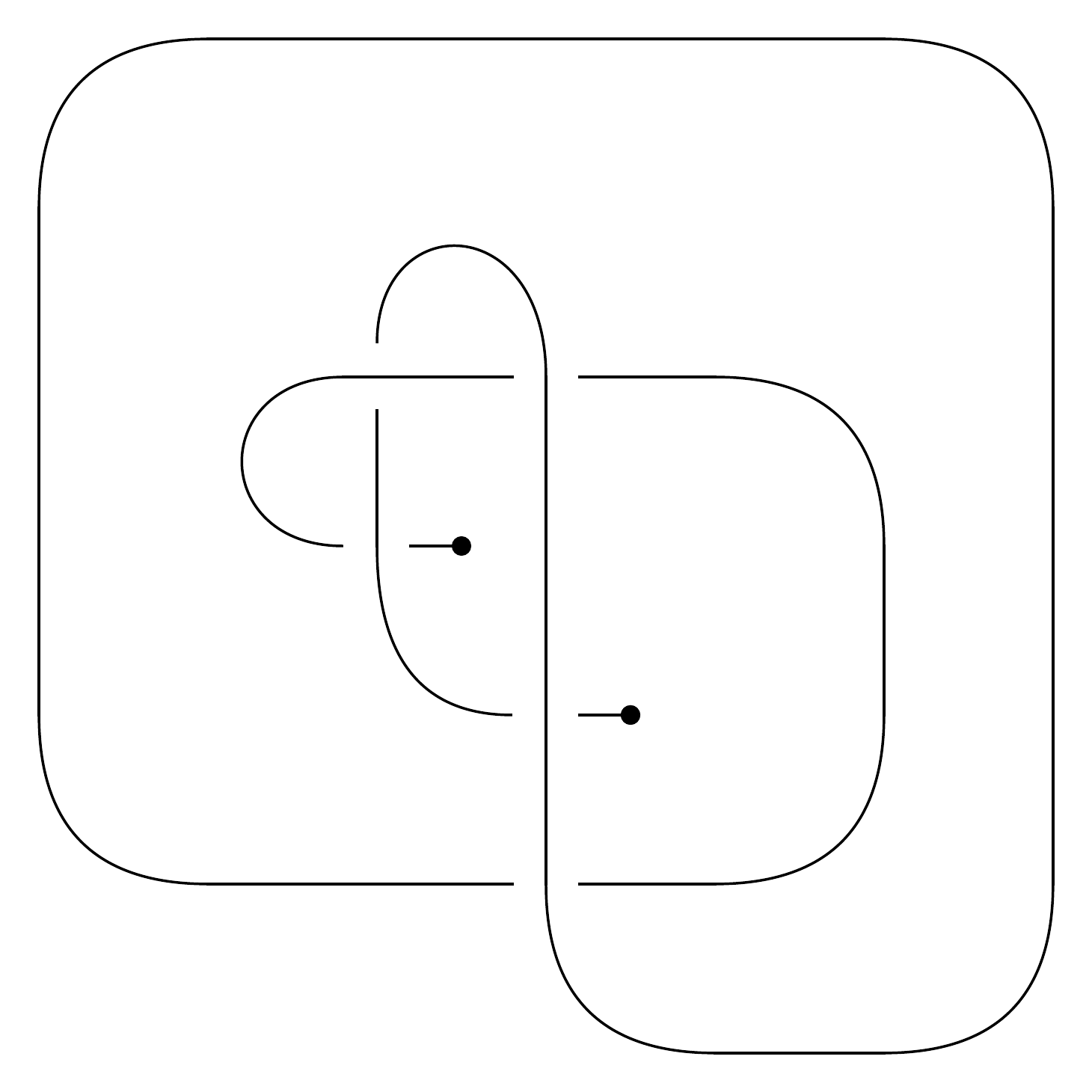}\\
\textcolor{black}{$5_{699}$}
\vspace{1cm}
\end{minipage}
\begin{minipage}[t]{.25\linewidth}
\centering
\includegraphics[width=0.9\textwidth,height=3.5cm,keepaspectratio]{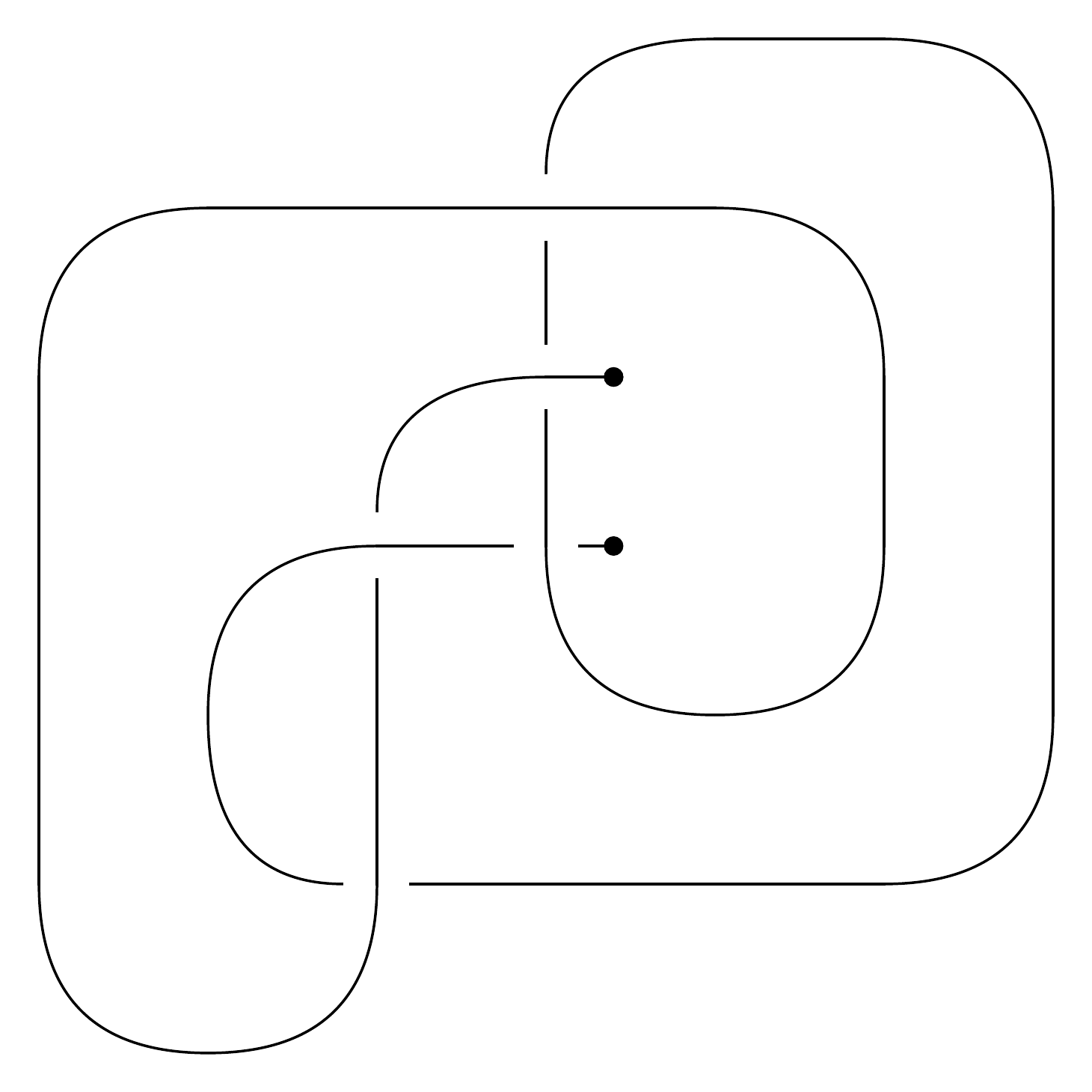}\\
\textcolor{black}{$5_{700}$}
\vspace{1cm}
\end{minipage}
\begin{minipage}[t]{.25\linewidth}
\centering
\includegraphics[width=0.9\textwidth,height=3.5cm,keepaspectratio]{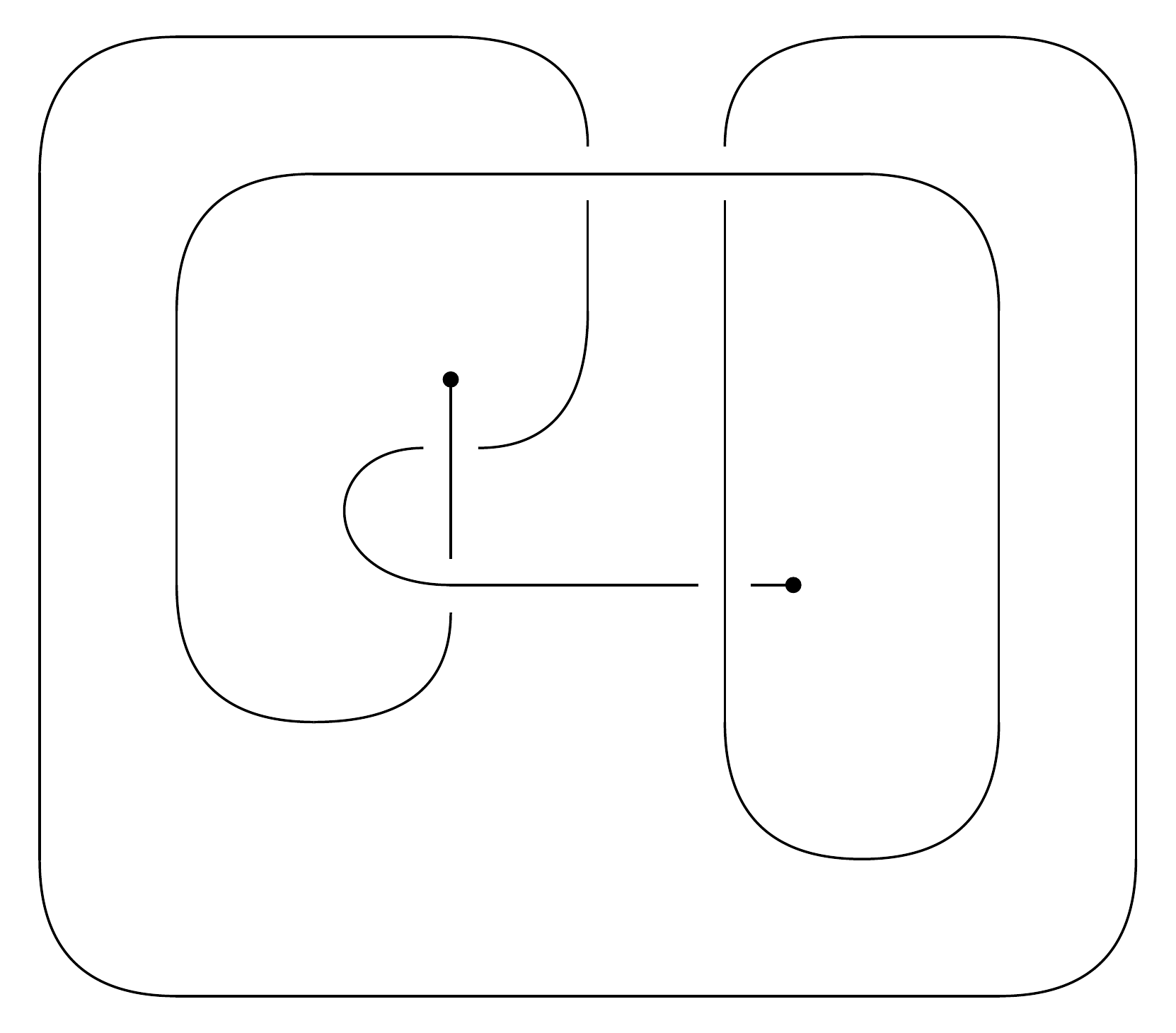}\\
\textcolor{black}{$5_{701}$}
\vspace{1cm}
\end{minipage}
\begin{minipage}[t]{.25\linewidth}
\centering
\includegraphics[width=0.9\textwidth,height=3.5cm,keepaspectratio]{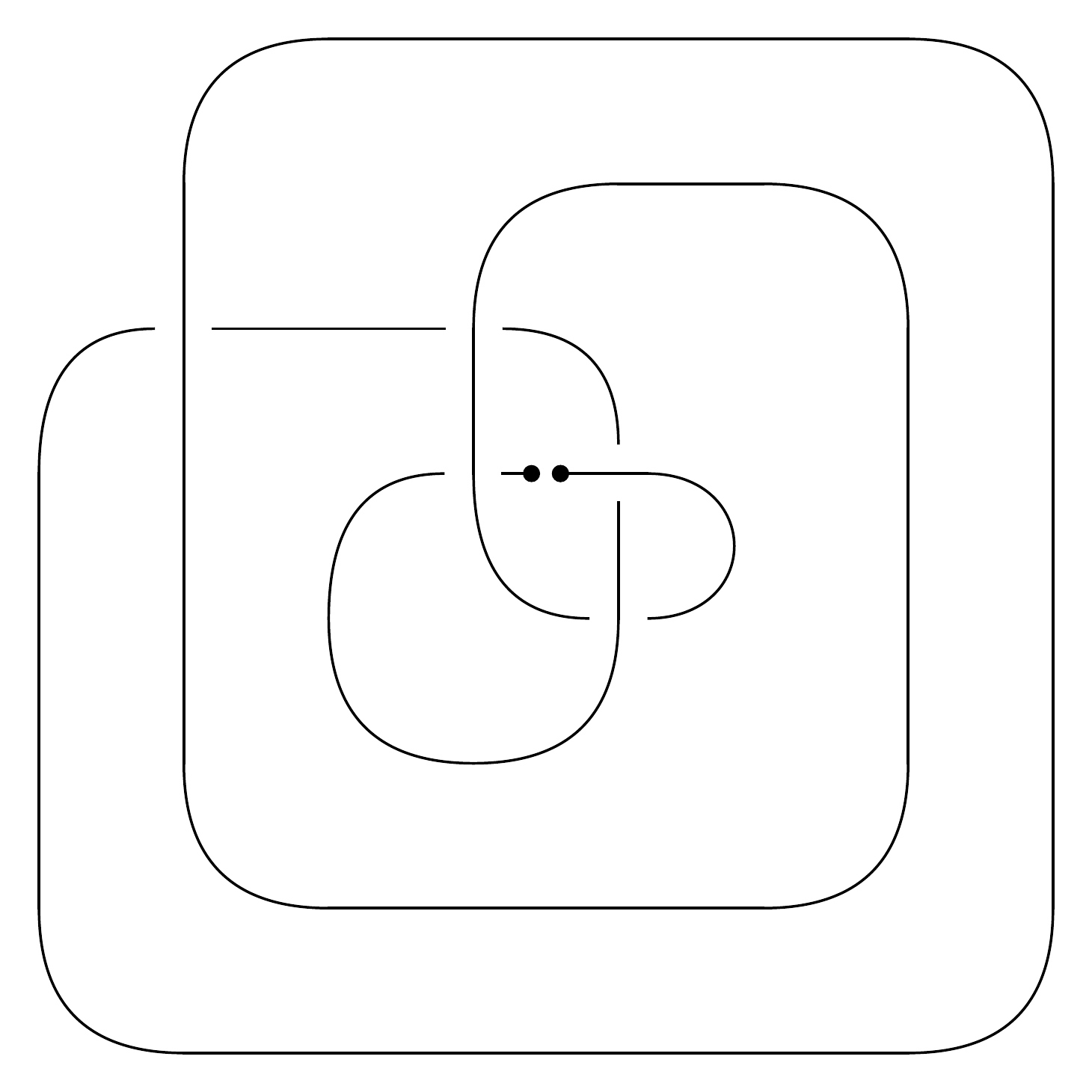}\\
\textcolor{black}{$5_{702}$}
\vspace{1cm}
\end{minipage}
\begin{minipage}[t]{.25\linewidth}
\centering
\includegraphics[width=0.9\textwidth,height=3.5cm,keepaspectratio]{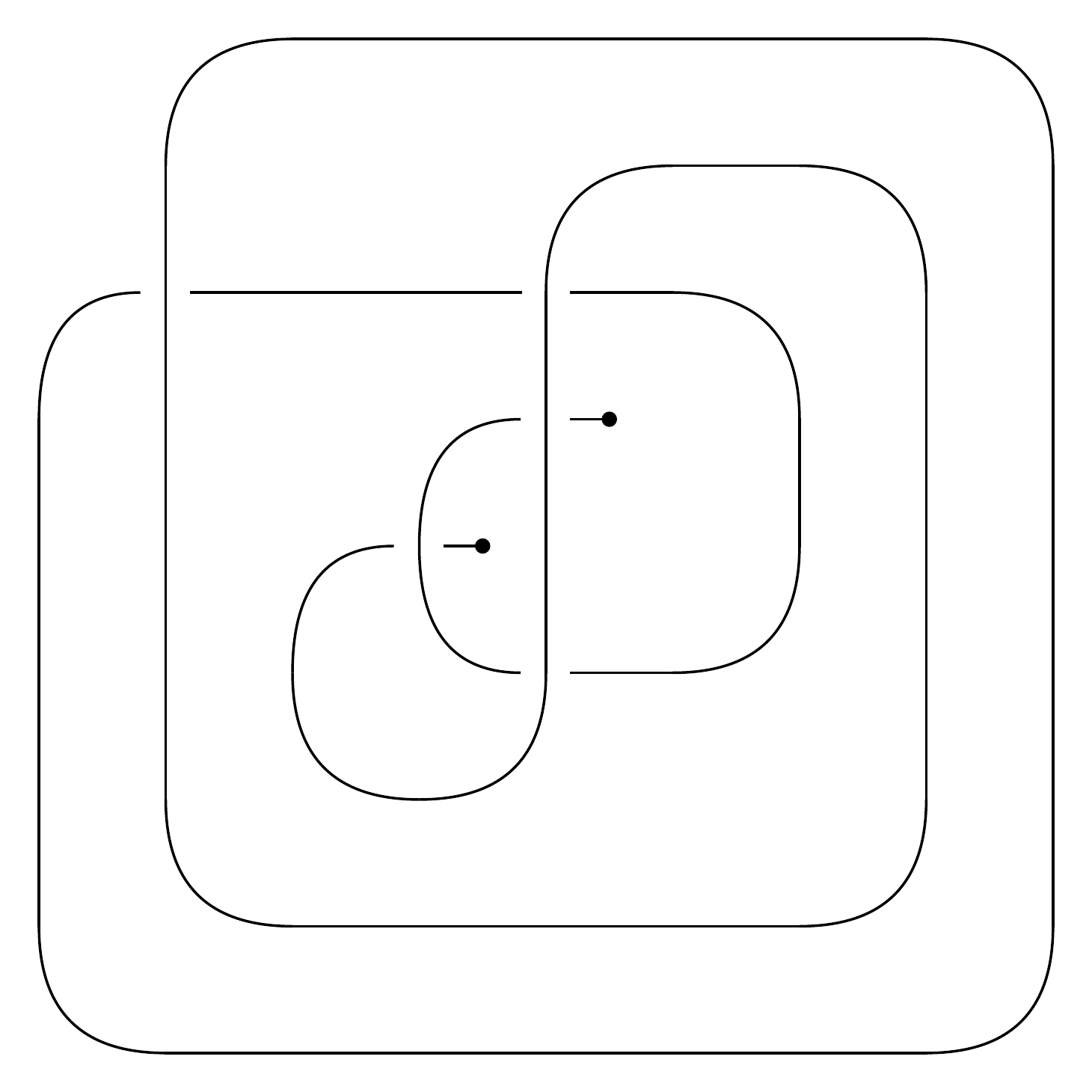}\\
\textcolor{black}{$5_{703}$}
\vspace{1cm}
\end{minipage}
\begin{minipage}[t]{.25\linewidth}
\centering
\includegraphics[width=0.9\textwidth,height=3.5cm,keepaspectratio]{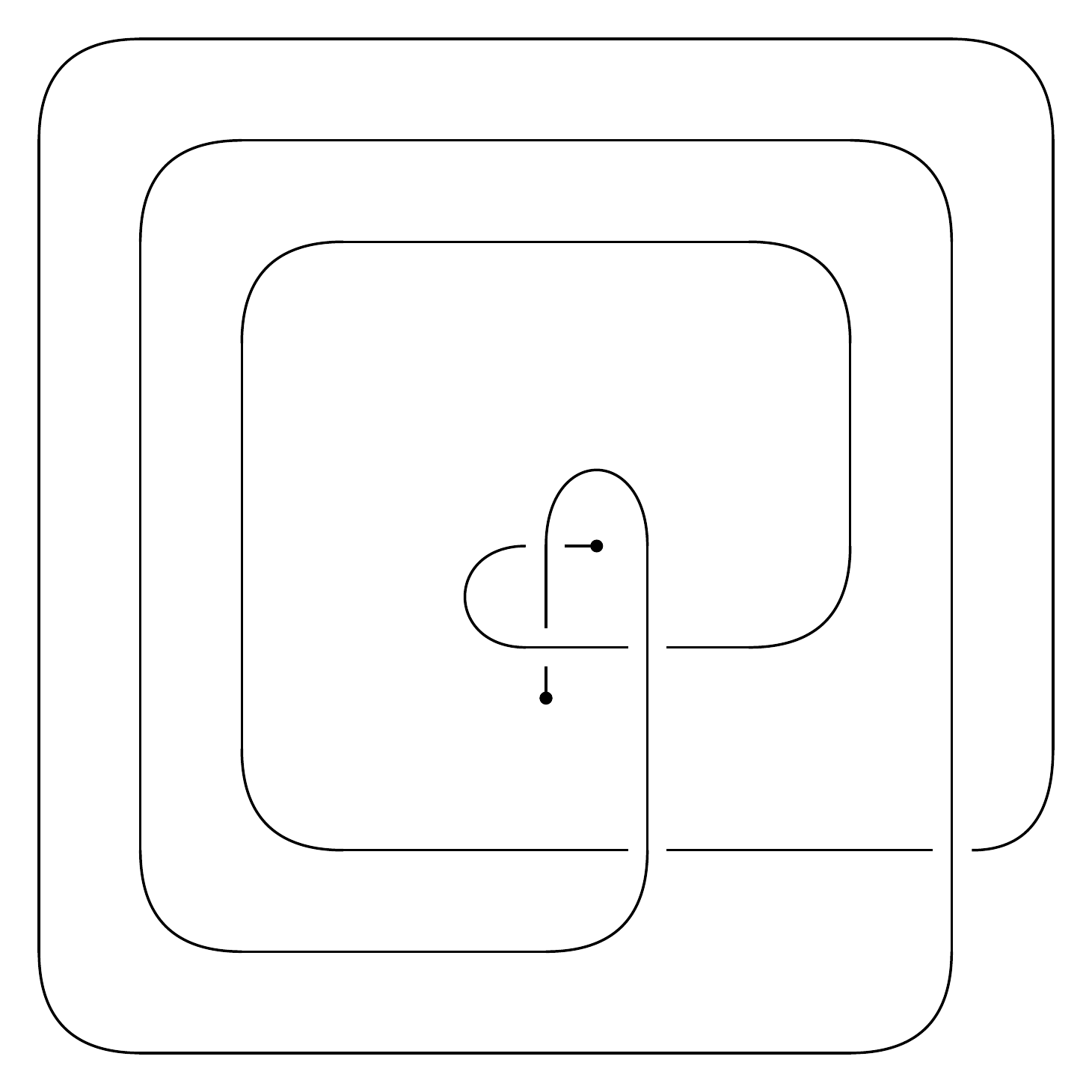}\\
\textcolor{black}{$5_{704}$}
\vspace{1cm}
\end{minipage}
\begin{minipage}[t]{.25\linewidth}
\centering
\includegraphics[width=0.9\textwidth,height=3.5cm,keepaspectratio]{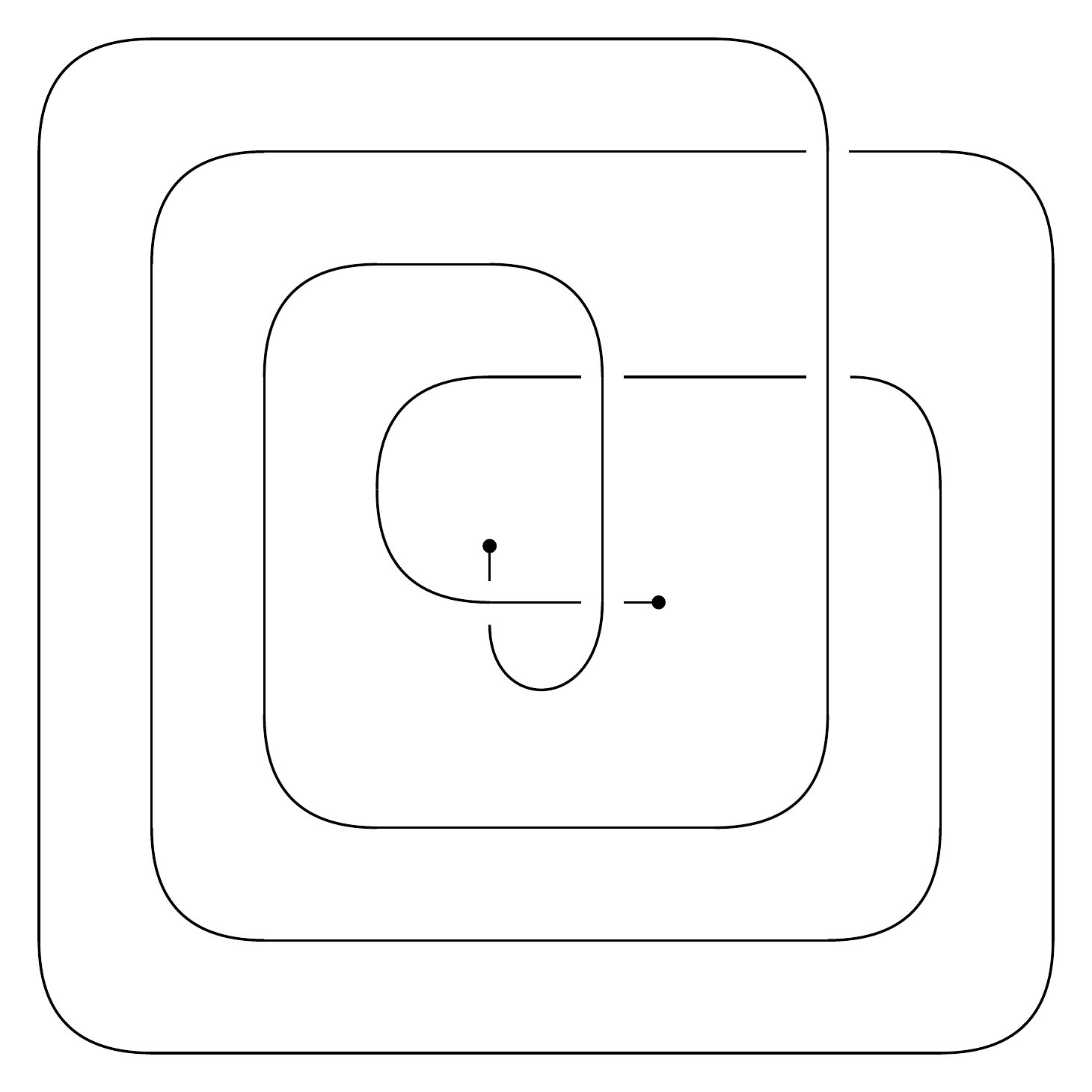}\\
\textcolor{black}{$5_{705}$}
\vspace{1cm}
\end{minipage}
\begin{minipage}[t]{.25\linewidth}
\centering
\includegraphics[width=0.9\textwidth,height=3.5cm,keepaspectratio]{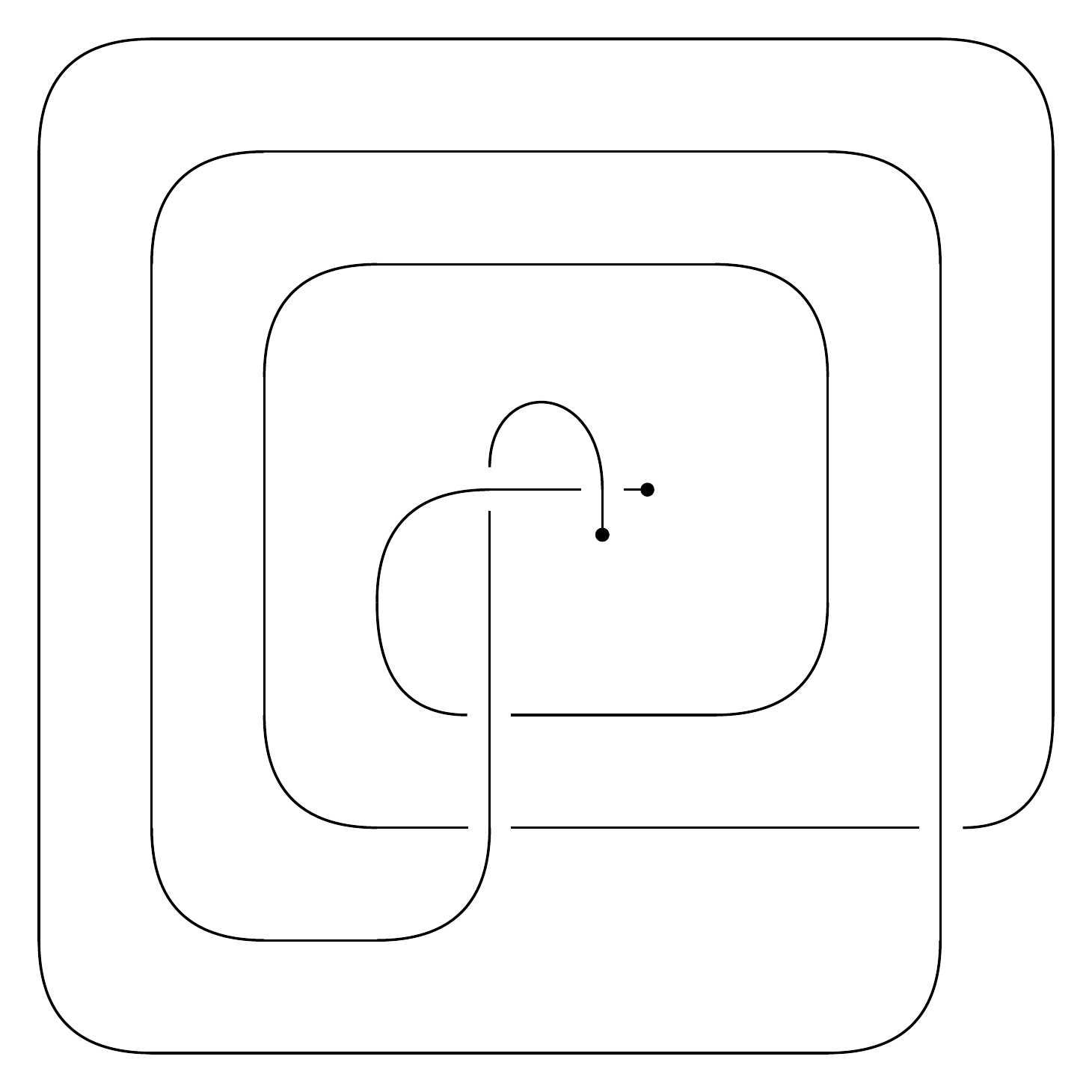}\\
\textcolor{black}{$5_{706}$}
\vspace{1cm}
\end{minipage}
\begin{minipage}[t]{.25\linewidth}
\centering
\includegraphics[width=0.9\textwidth,height=3.5cm,keepaspectratio]{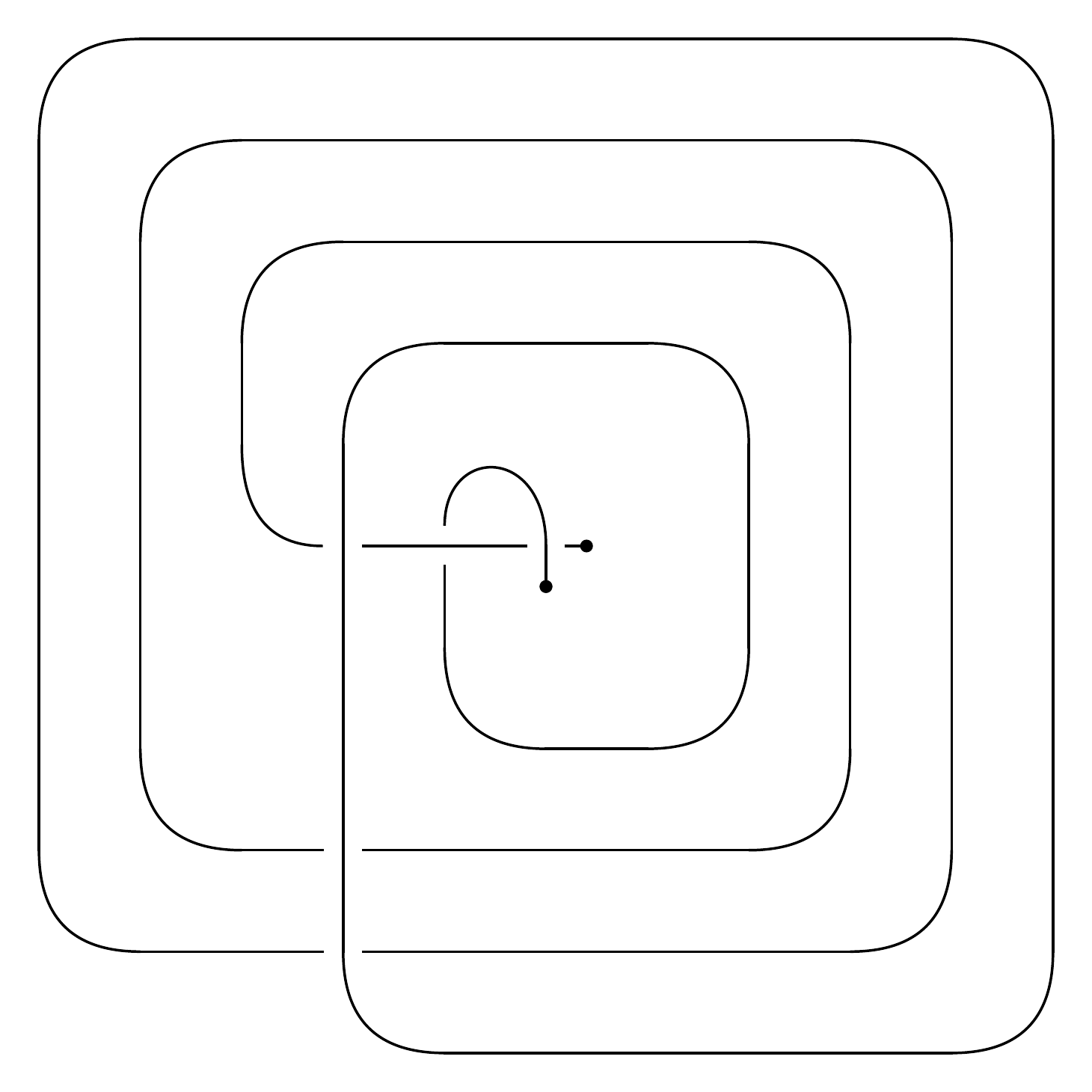}\\
\textcolor{black}{$5_{707}$}
\vspace{1cm}
\end{minipage}
\begin{minipage}[t]{.25\linewidth}
\centering
\includegraphics[width=0.9\textwidth,height=3.5cm,keepaspectratio]{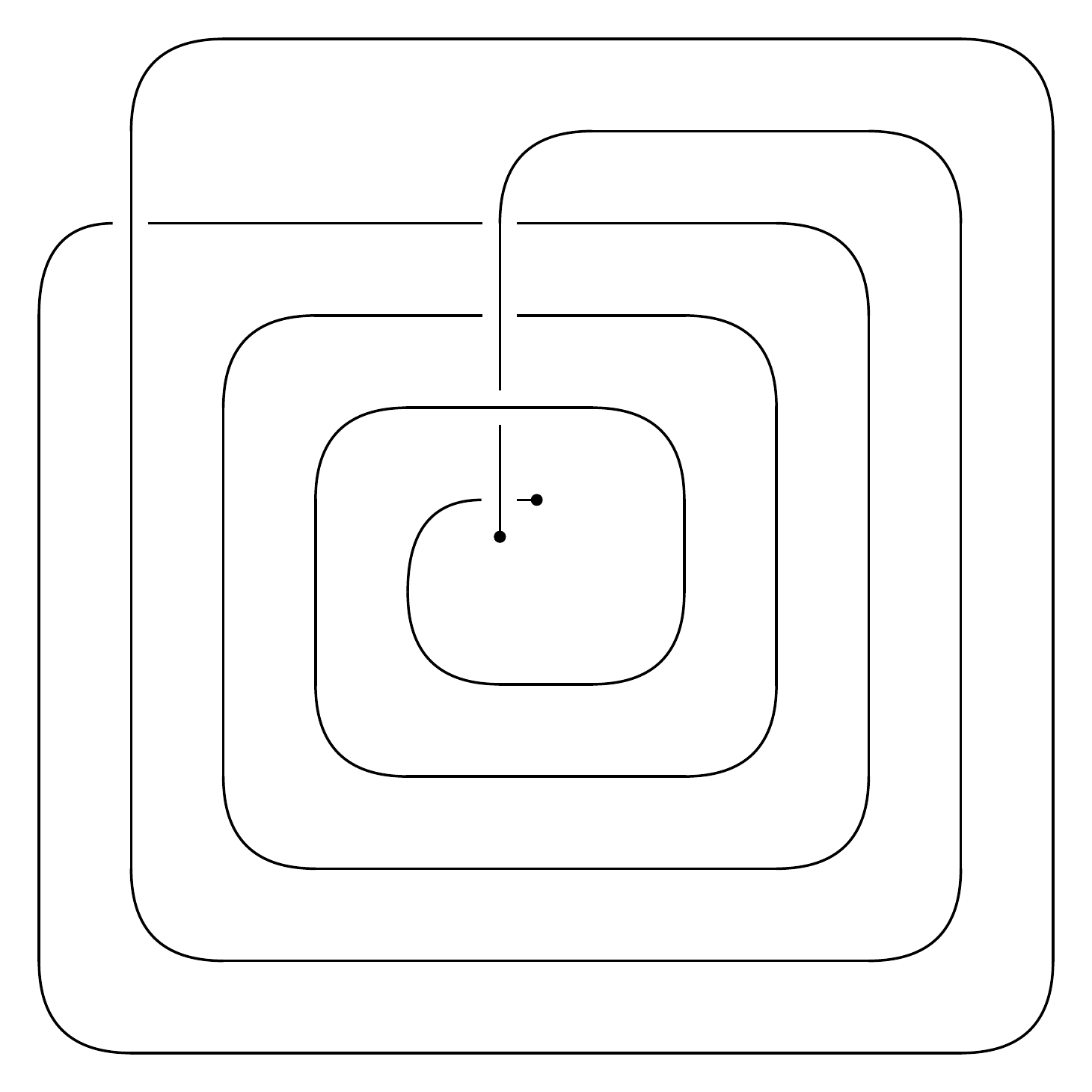}\\
\textcolor{black}{$5_{708}$}
\vspace{1cm}
\end{minipage}
\begin{minipage}[t]{.25\linewidth}
\centering
\includegraphics[width=0.9\textwidth,height=3.5cm,keepaspectratio]{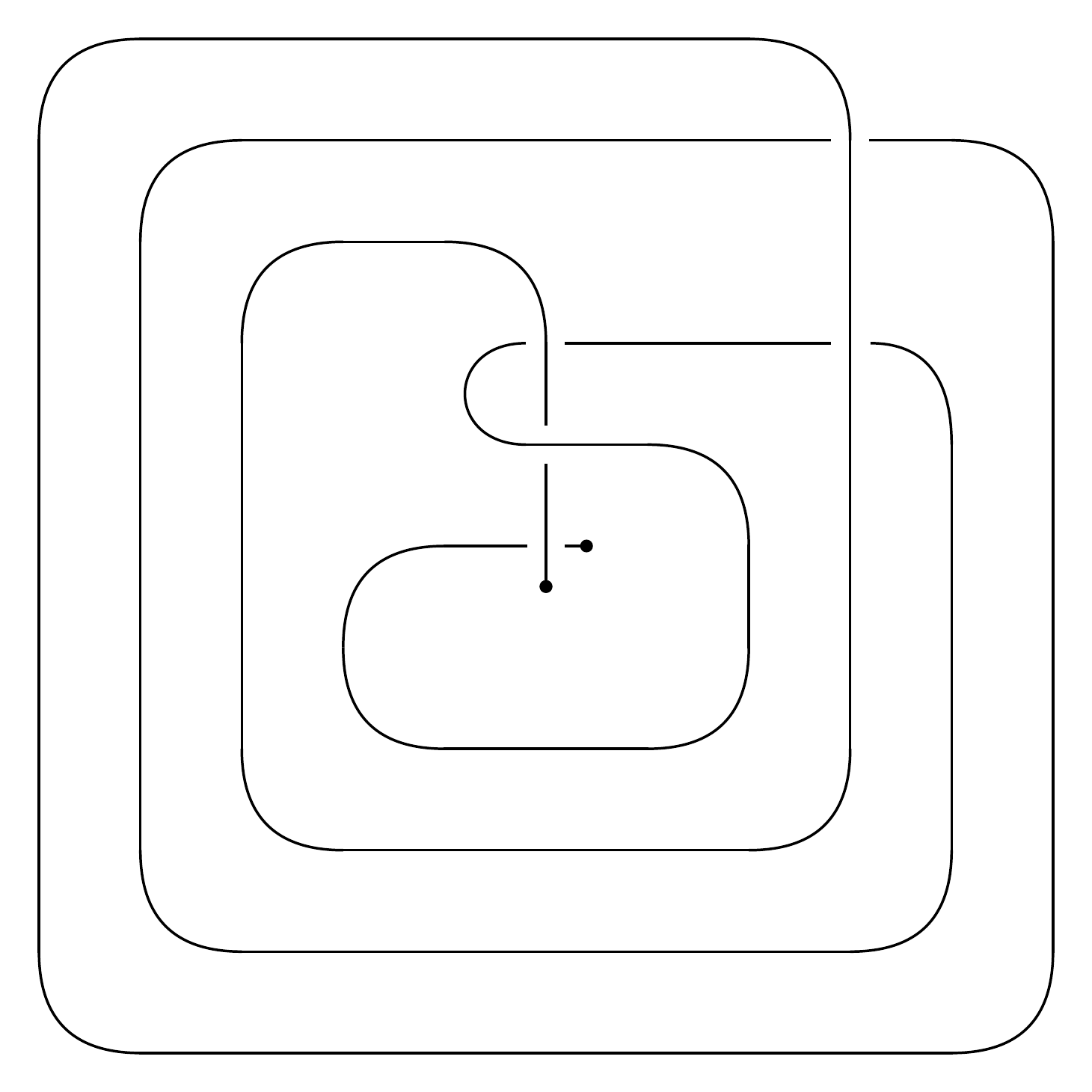}\\
\textcolor{black}{$5_{709}$}
\vspace{1cm}
\end{minipage}
\begin{minipage}[t]{.25\linewidth}
\centering
\includegraphics[width=0.9\textwidth,height=3.5cm,keepaspectratio]{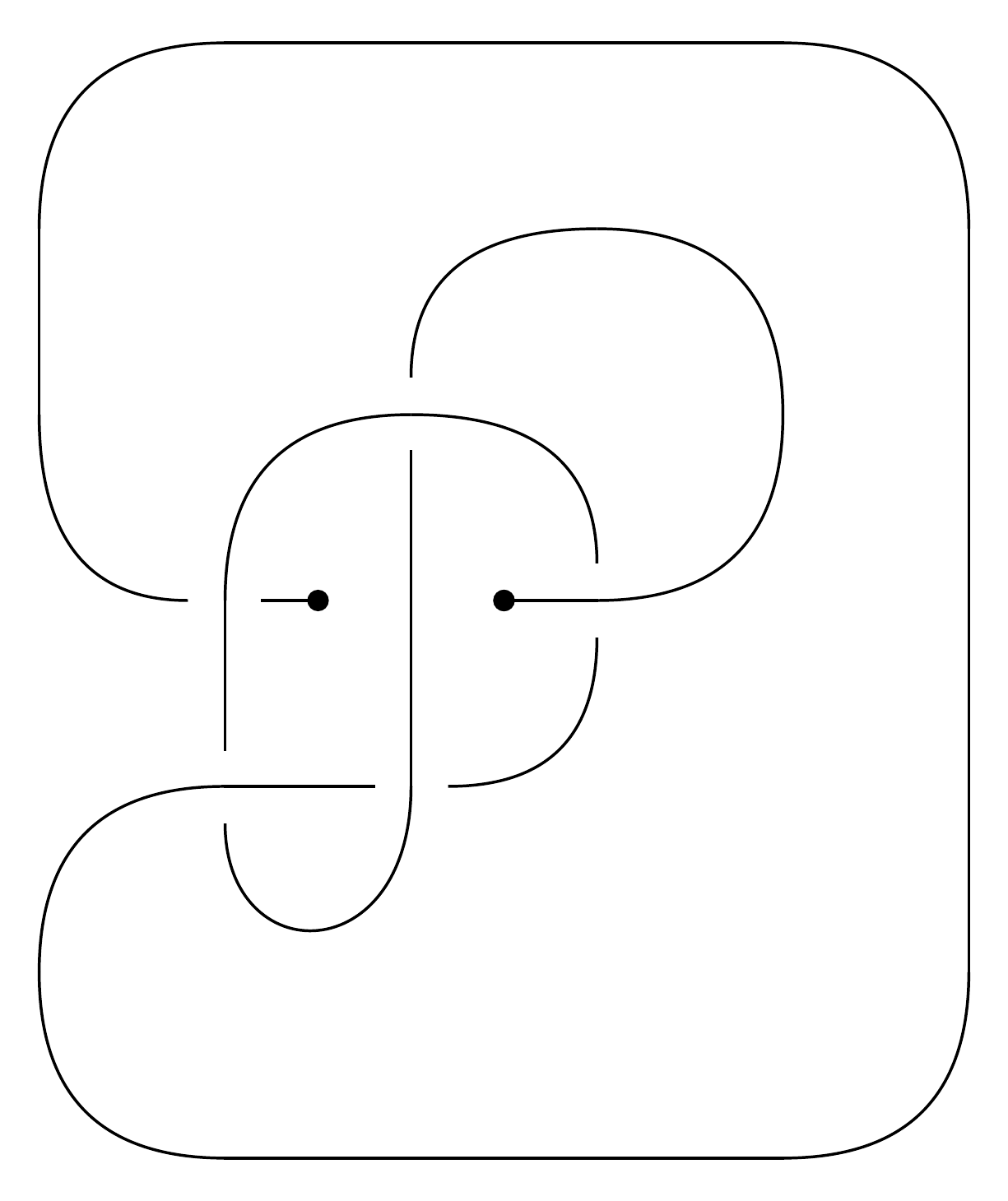}\\
\textcolor{black}{$5_{710}$}
\vspace{1cm}
\end{minipage}
\begin{minipage}[t]{.25\linewidth}
\centering
\includegraphics[width=0.9\textwidth,height=3.5cm,keepaspectratio]{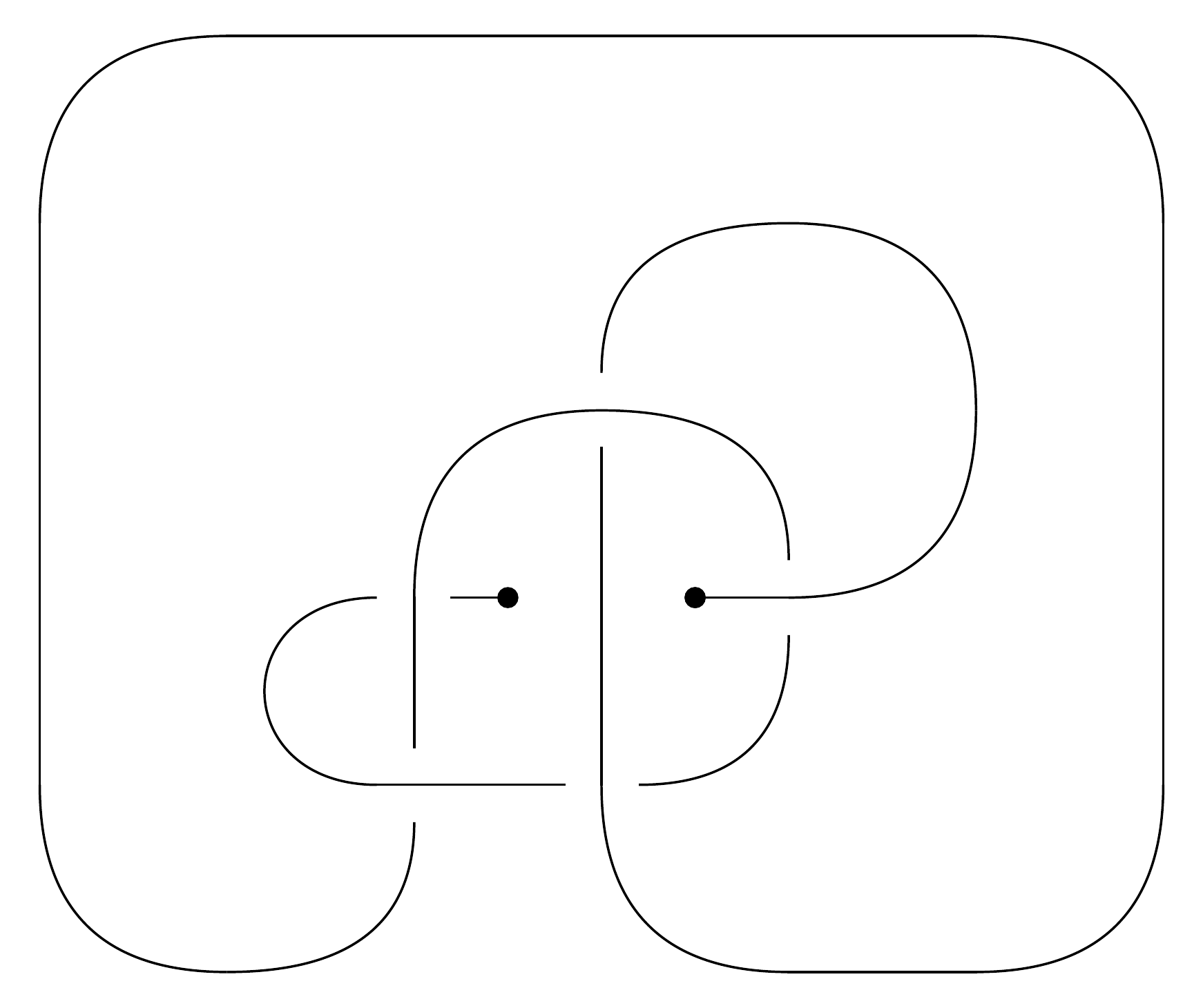}\\
\textcolor{black}{$5_{711}$}
\vspace{1cm}
\end{minipage}
\begin{minipage}[t]{.25\linewidth}
\centering
\includegraphics[width=0.9\textwidth,height=3.5cm,keepaspectratio]{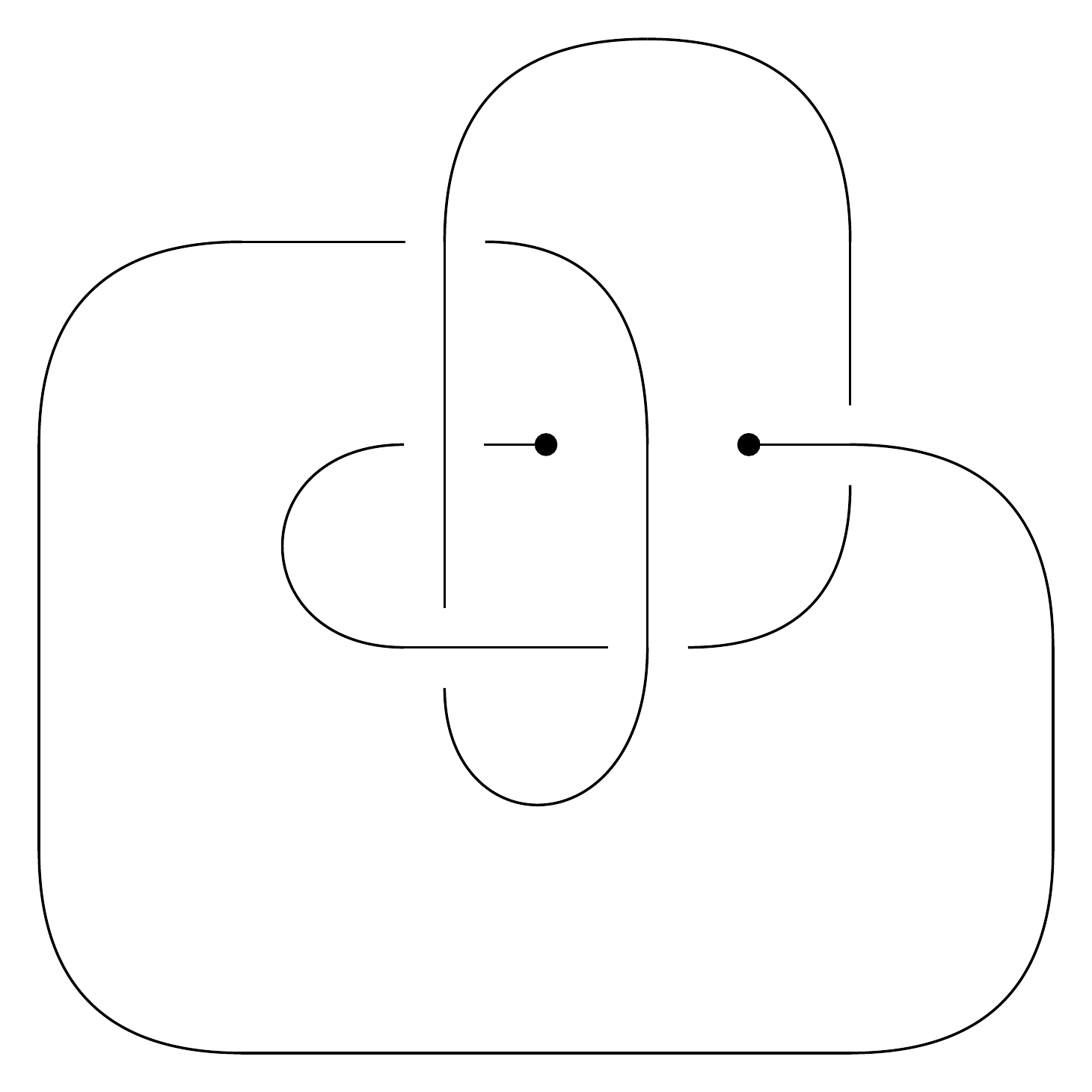}\\
\textcolor{black}{$5_{712}$}
\vspace{1cm}
\end{minipage}
\begin{minipage}[t]{.25\linewidth}
\centering
\includegraphics[width=0.9\textwidth,height=3.5cm,keepaspectratio]{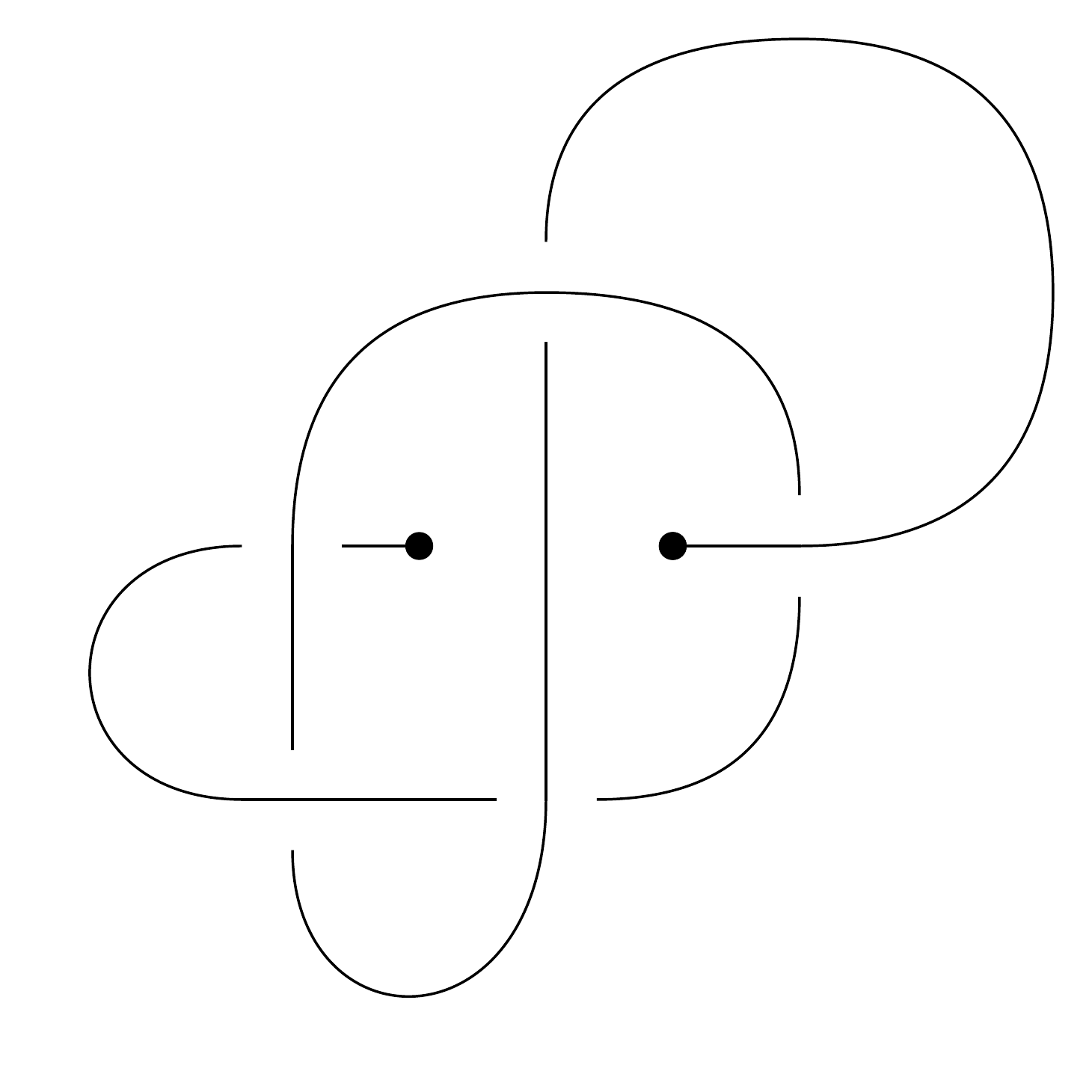}\\
\textcolor{black}{$5_{713}$}
\vspace{1cm}
\end{minipage}
\begin{minipage}[t]{.25\linewidth}
\centering
\includegraphics[width=0.9\textwidth,height=3.5cm,keepaspectratio]{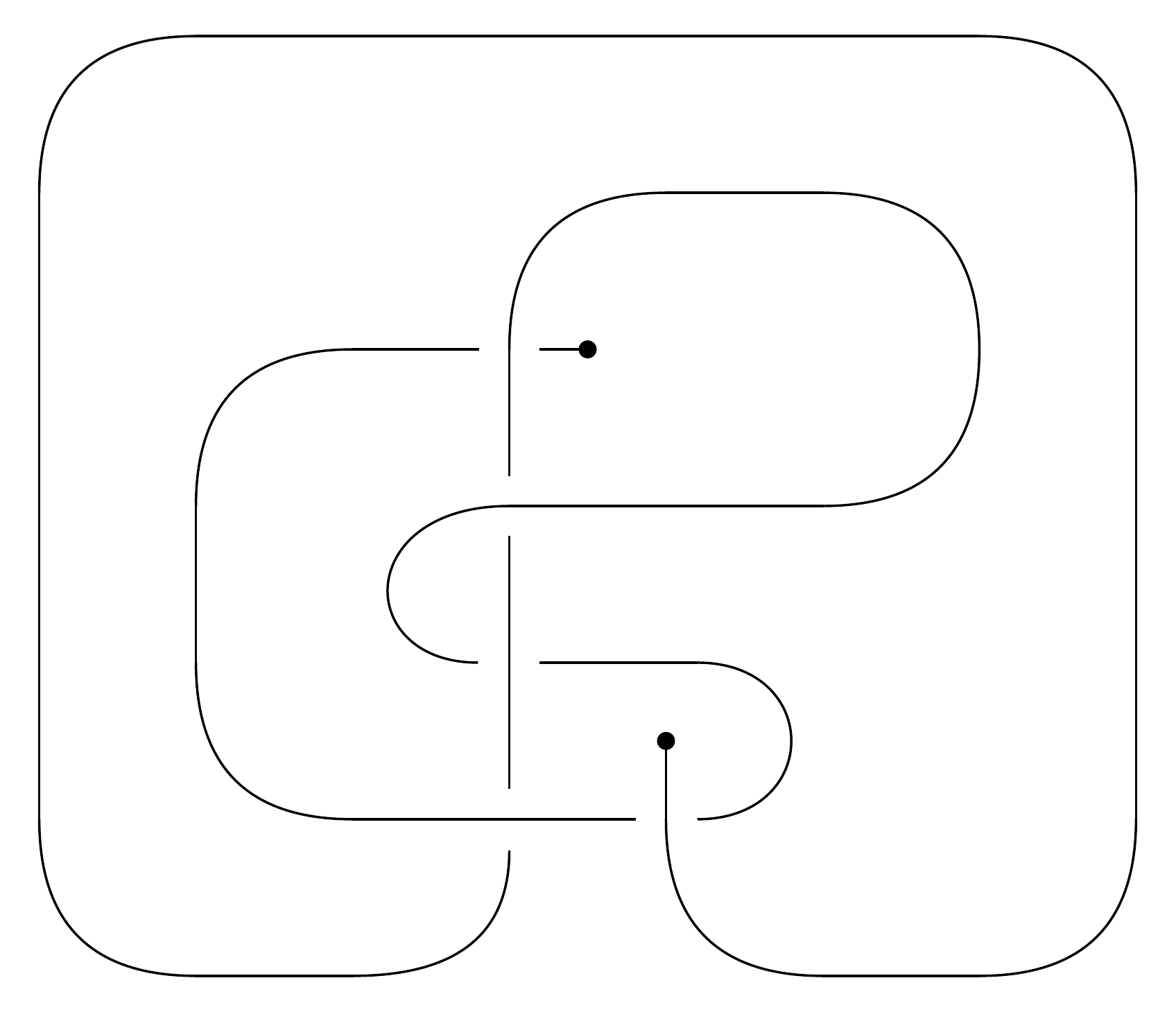}\\
\textcolor{black}{$5_{714}$}
\vspace{1cm}
\end{minipage}
\begin{minipage}[t]{.25\linewidth}
\centering
\includegraphics[width=0.9\textwidth,height=3.5cm,keepaspectratio]{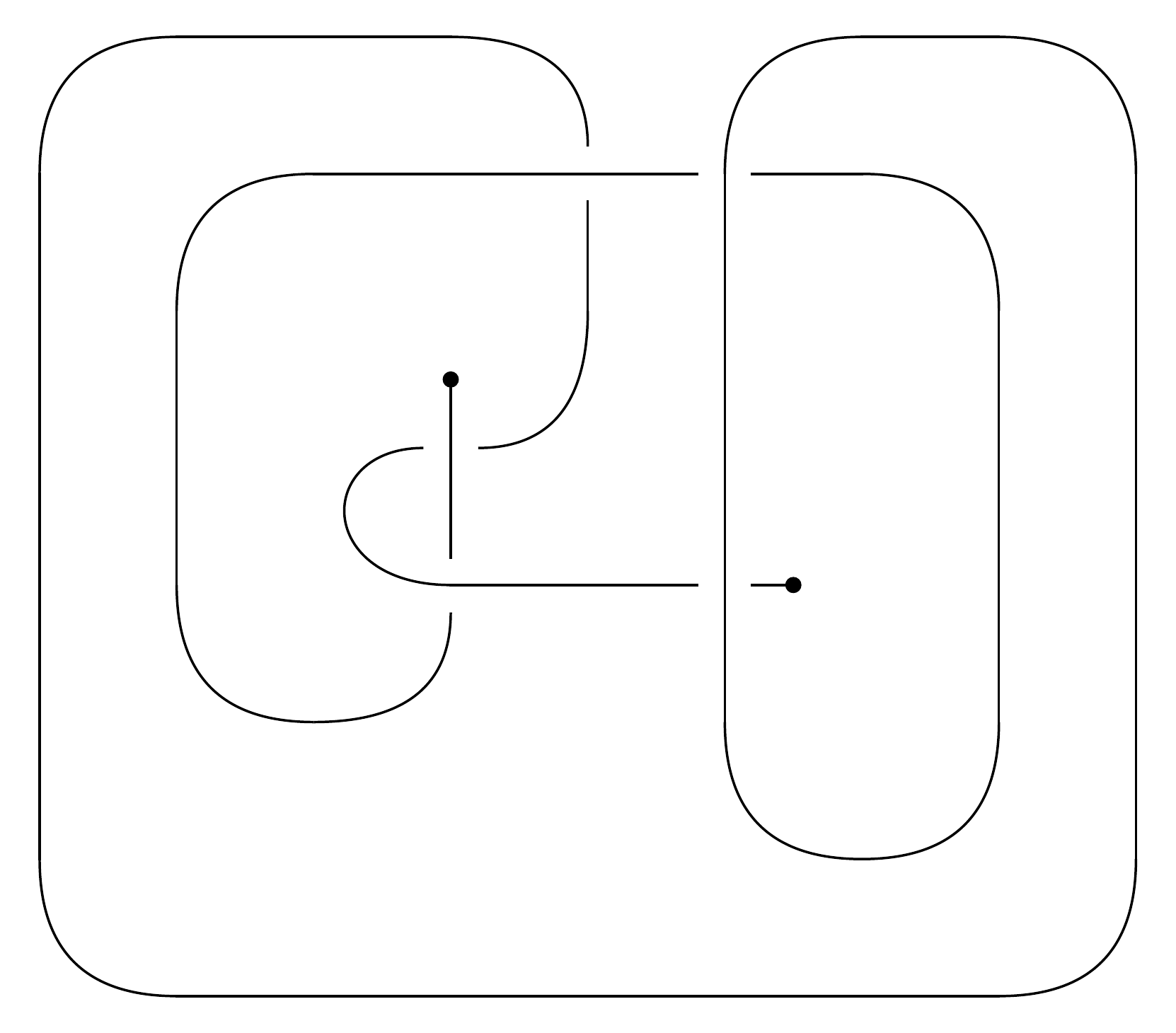}\\
\textcolor{black}{$5_{715}$}
\vspace{1cm}
\end{minipage}
\begin{minipage}[t]{.25\linewidth}
\centering
\includegraphics[width=0.9\textwidth,height=3.5cm,keepaspectratio]{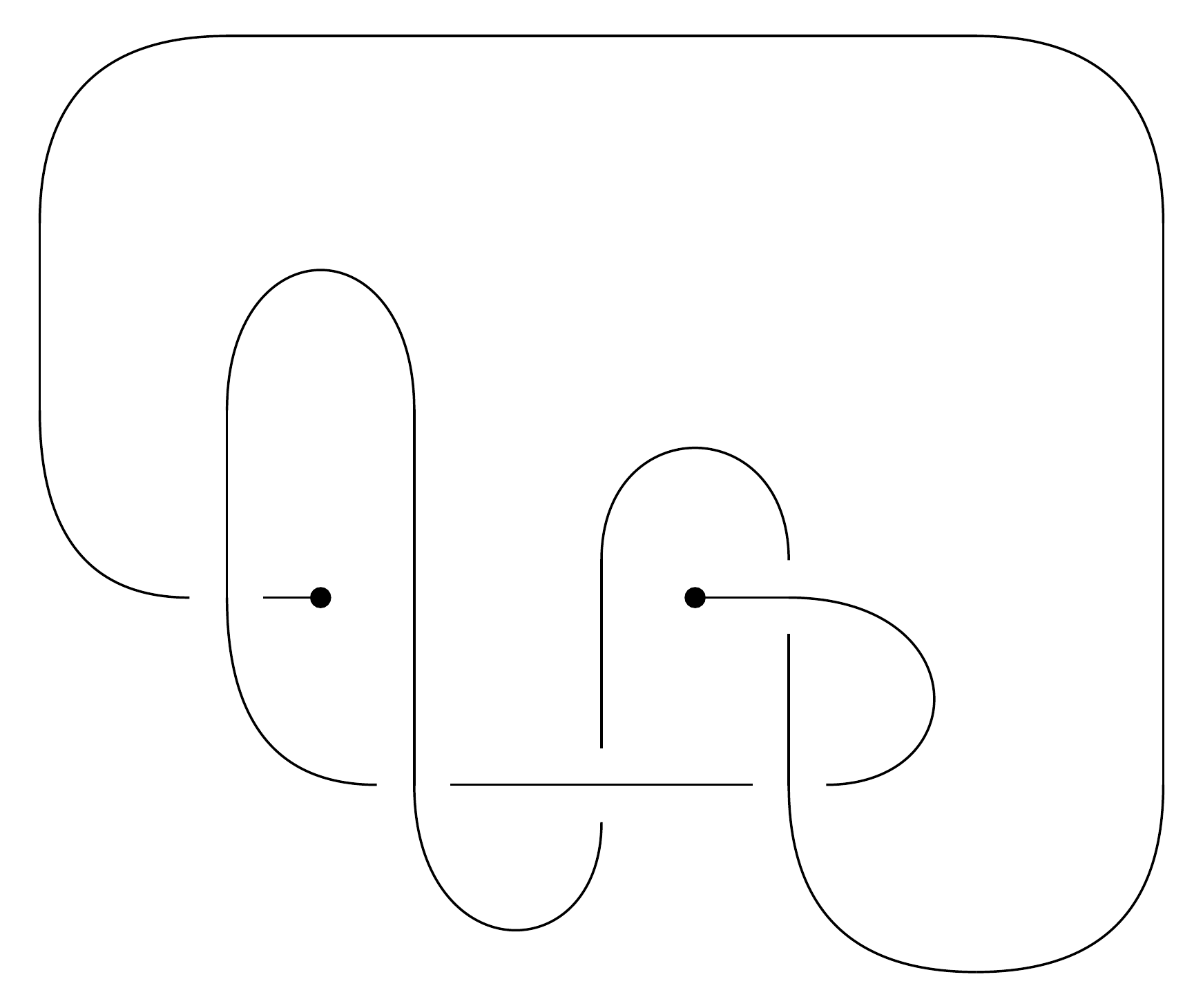}\\
\textcolor{black}{$5_{716}$}
\vspace{1cm}
\end{minipage}
\begin{minipage}[t]{.25\linewidth}
\centering
\includegraphics[width=0.9\textwidth,height=3.5cm,keepaspectratio]{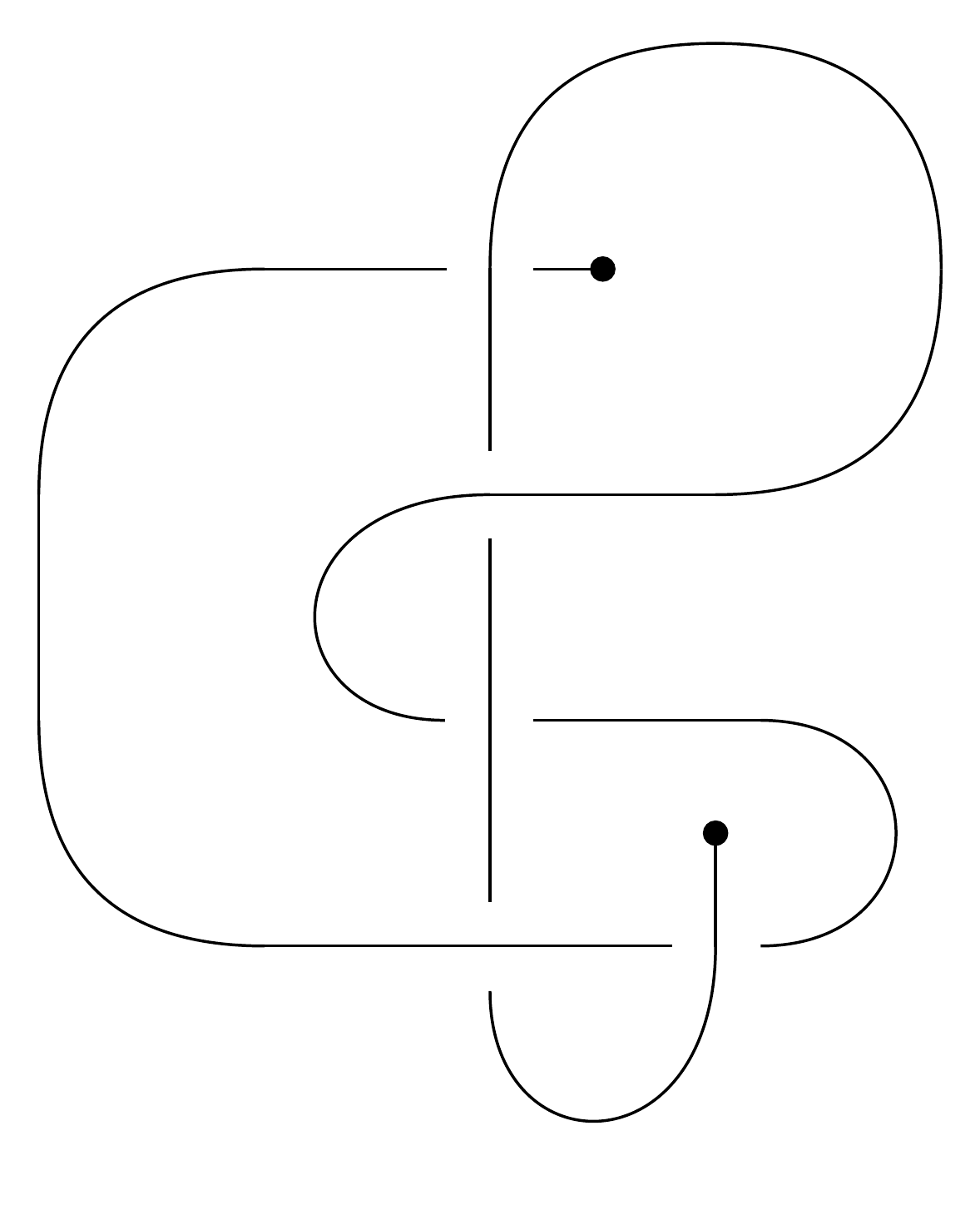}\\
\textcolor{black}{$5_{717}$}
\vspace{1cm}
\end{minipage}
\begin{minipage}[t]{.25\linewidth}
\centering
\includegraphics[width=0.9\textwidth,height=3.5cm,keepaspectratio]{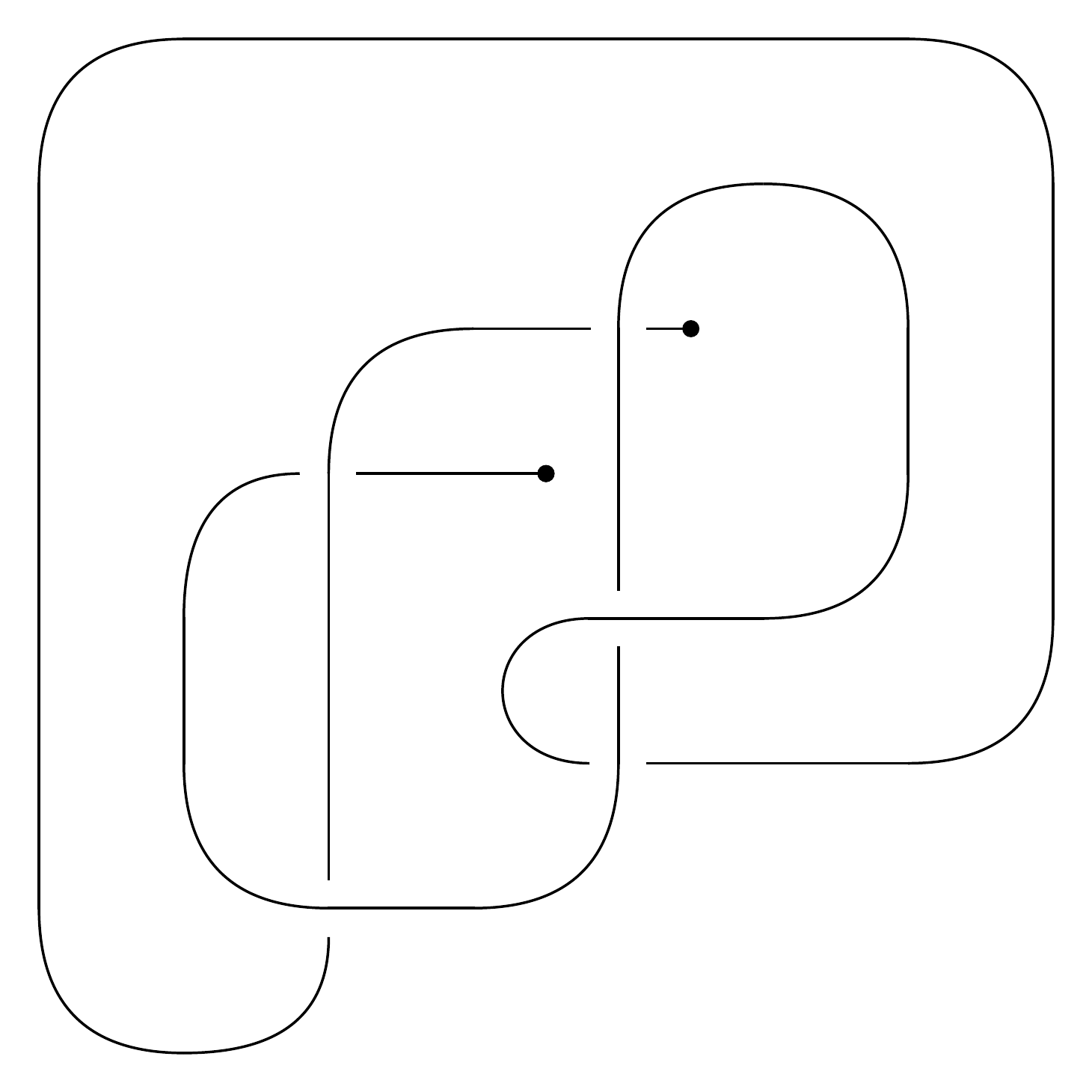}\\
\textcolor{black}{$5_{718}$}
\vspace{1cm}
\end{minipage}
\begin{minipage}[t]{.25\linewidth}
\centering
\includegraphics[width=0.9\textwidth,height=3.5cm,keepaspectratio]{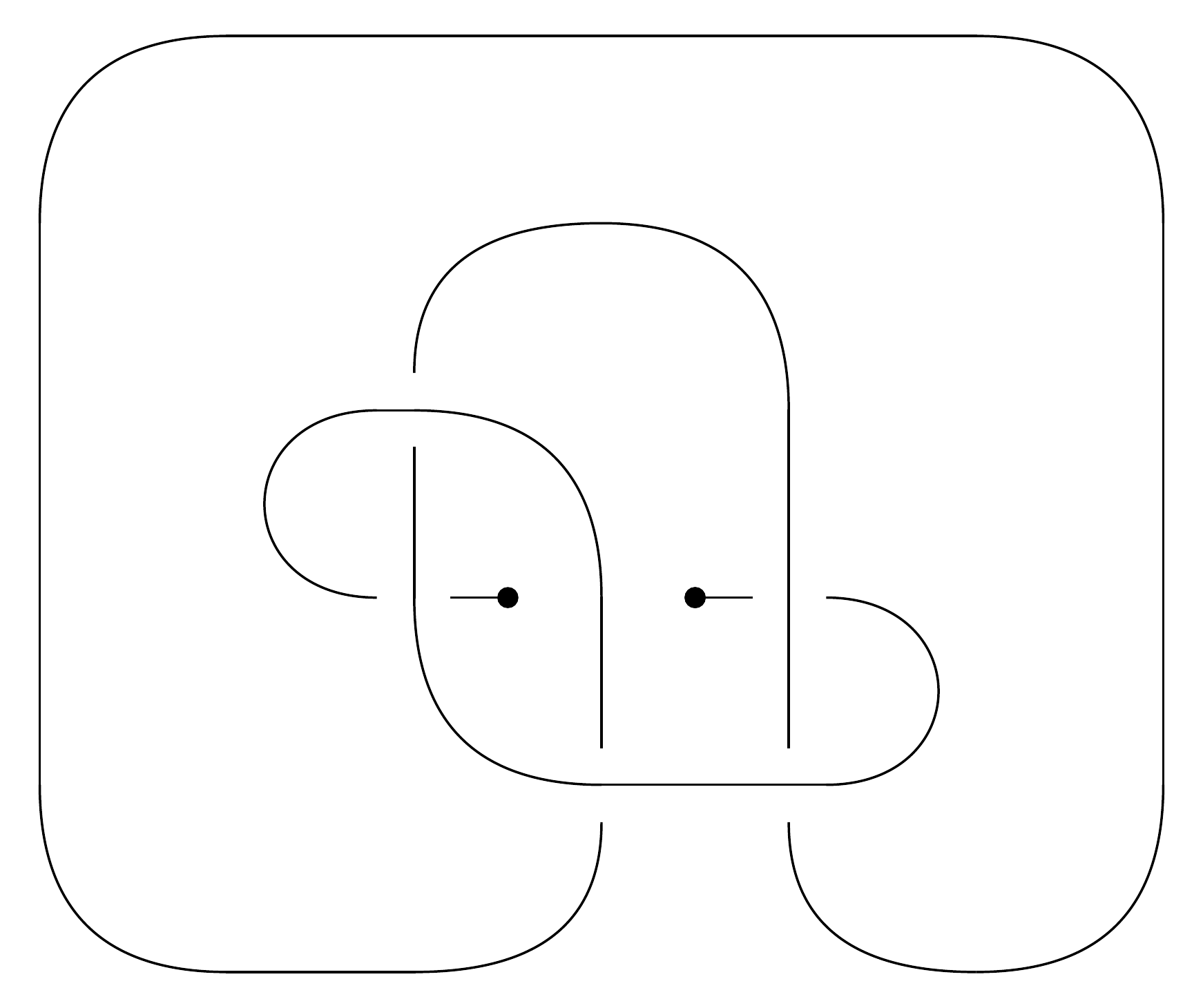}\\
\textcolor{black}{$5_{719}$}
\vspace{1cm}
\end{minipage}
\begin{minipage}[t]{.25\linewidth}
\centering
\includegraphics[width=0.9\textwidth,height=3.5cm,keepaspectratio]{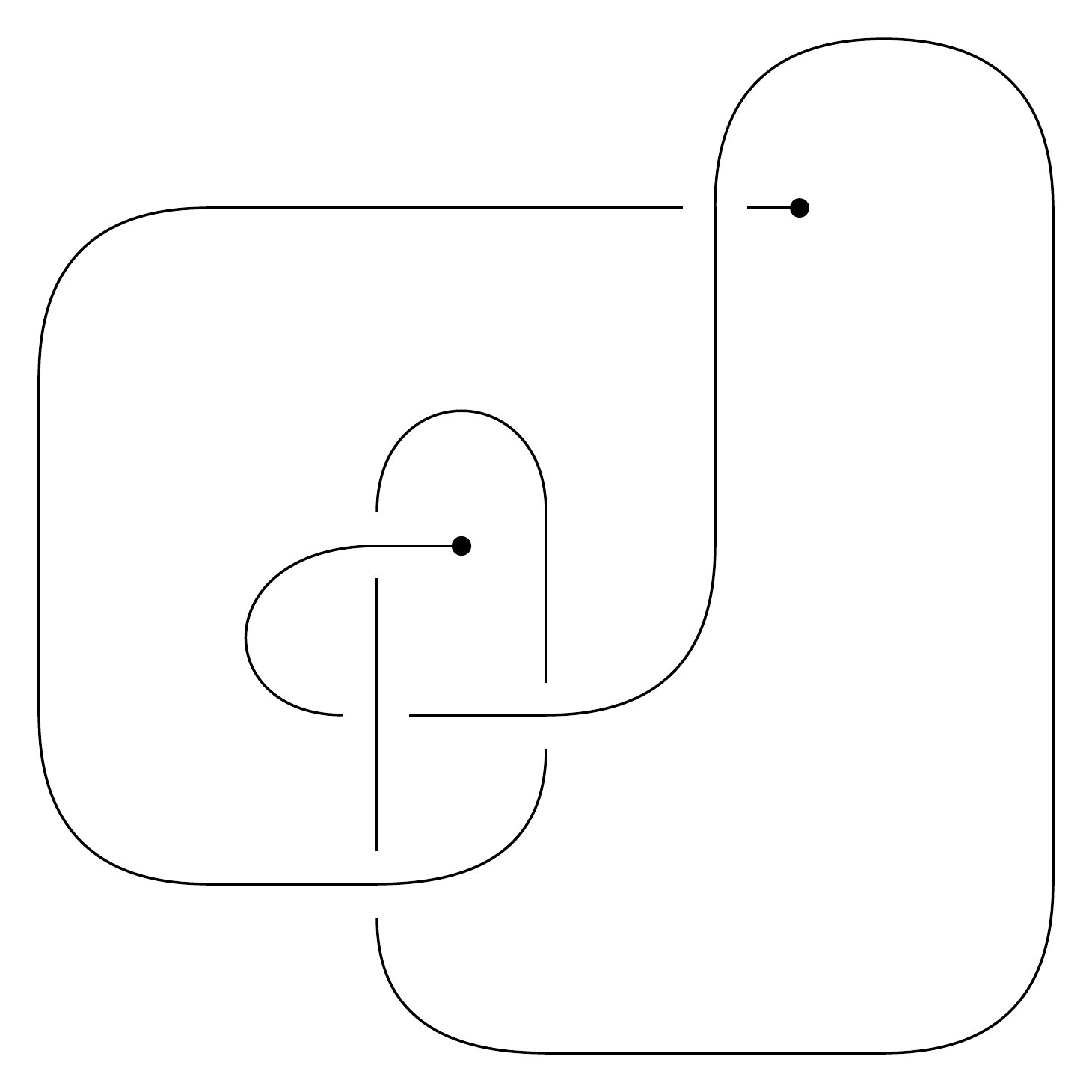}\\
\textcolor{black}{$5_{720}$}
\vspace{1cm}
\end{minipage}
\begin{minipage}[t]{.25\linewidth}
\centering
\includegraphics[width=0.9\textwidth,height=3.5cm,keepaspectratio]{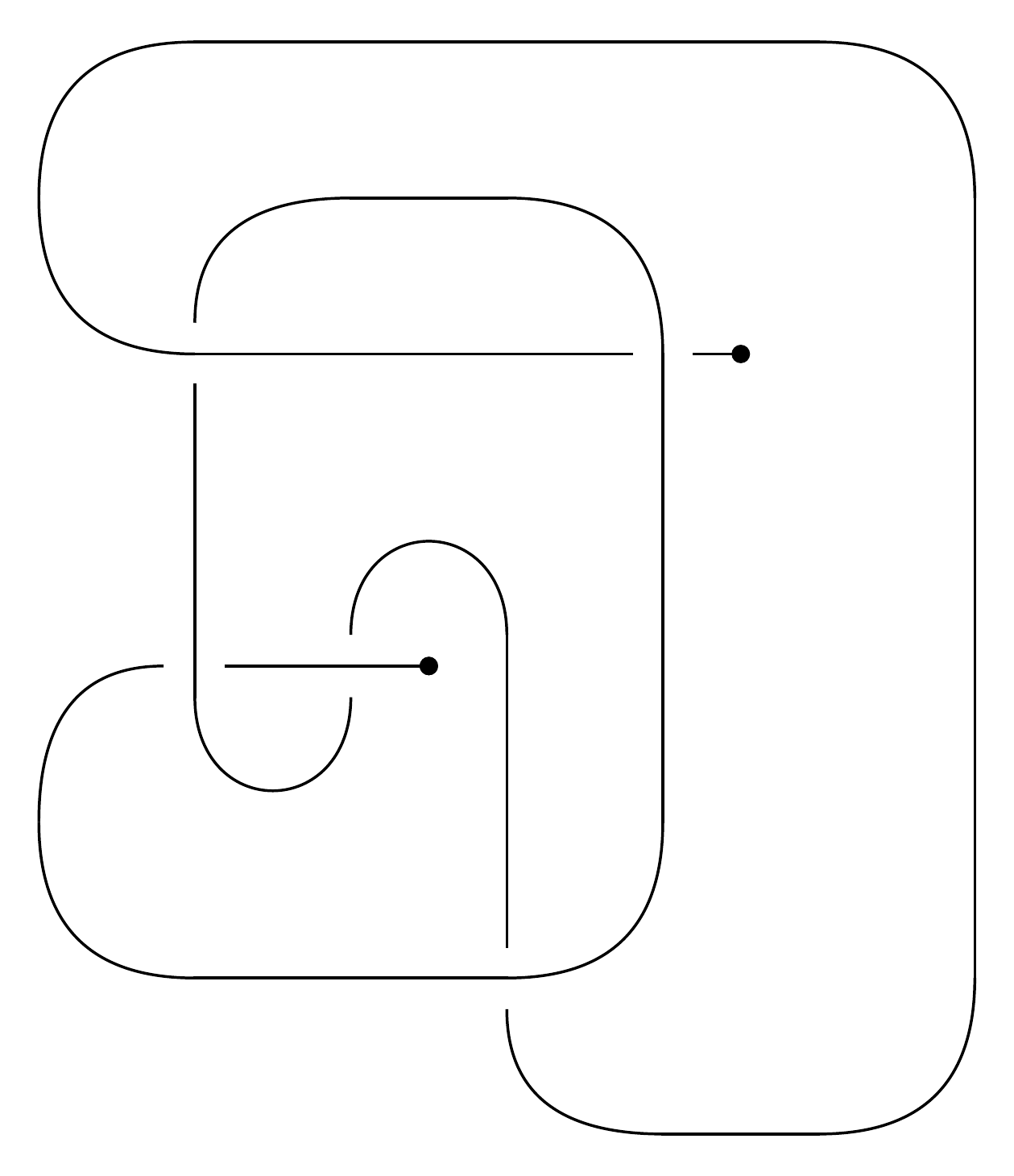}\\
\textcolor{black}{$5_{721}$}
\vspace{1cm}
\end{minipage}
\begin{minipage}[t]{.25\linewidth}
\centering
\includegraphics[width=0.9\textwidth,height=3.5cm,keepaspectratio]{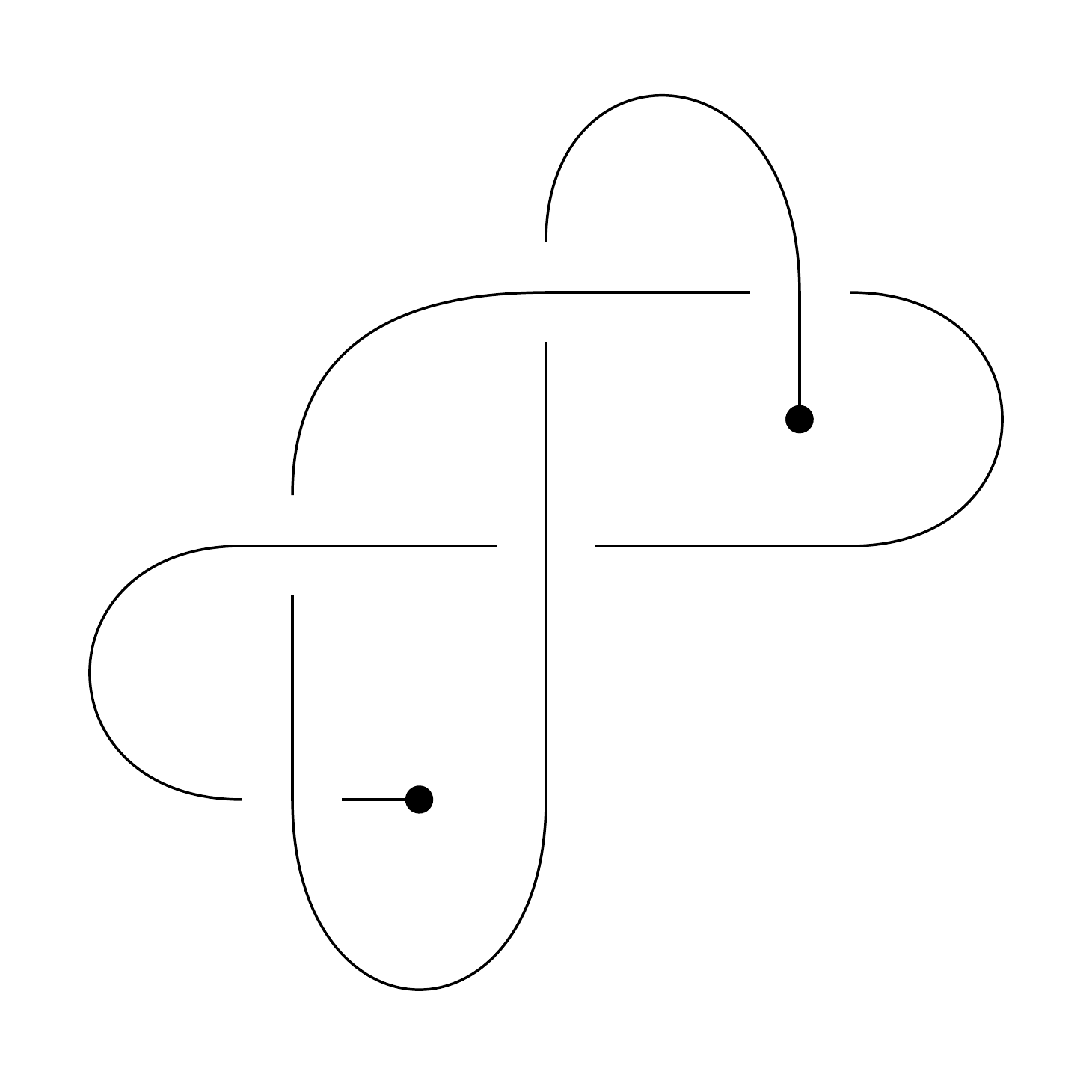}\\
\textcolor{black}{$5_{722}$}
\vspace{1cm}
\end{minipage}
\begin{minipage}[t]{.25\linewidth}
\centering
\includegraphics[width=0.9\textwidth,height=3.5cm,keepaspectratio]{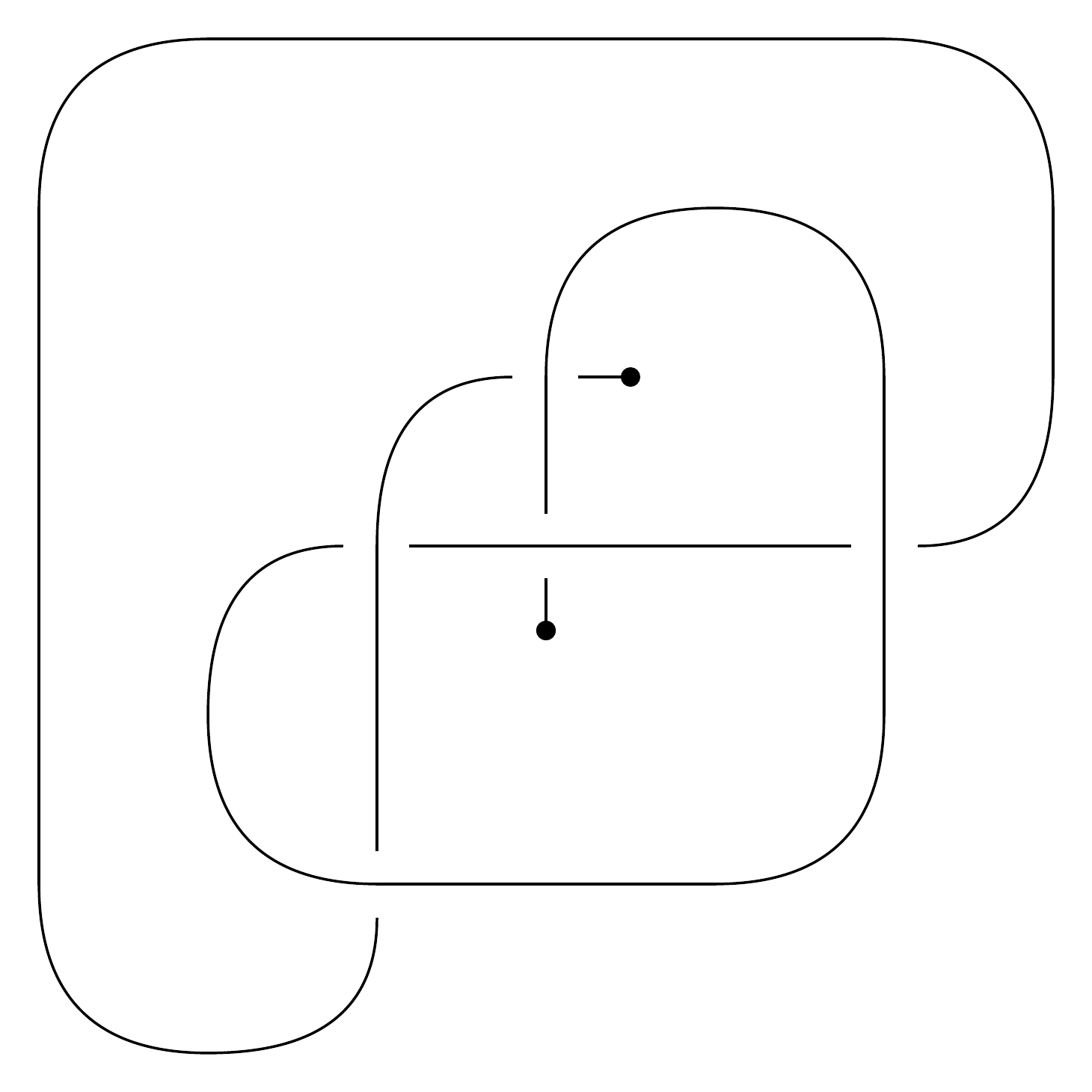}\\
\textcolor{black}{$5_{723}$}
\vspace{1cm}
\end{minipage}
\begin{minipage}[t]{.25\linewidth}
\centering
\includegraphics[width=0.9\textwidth,height=3.5cm,keepaspectratio]{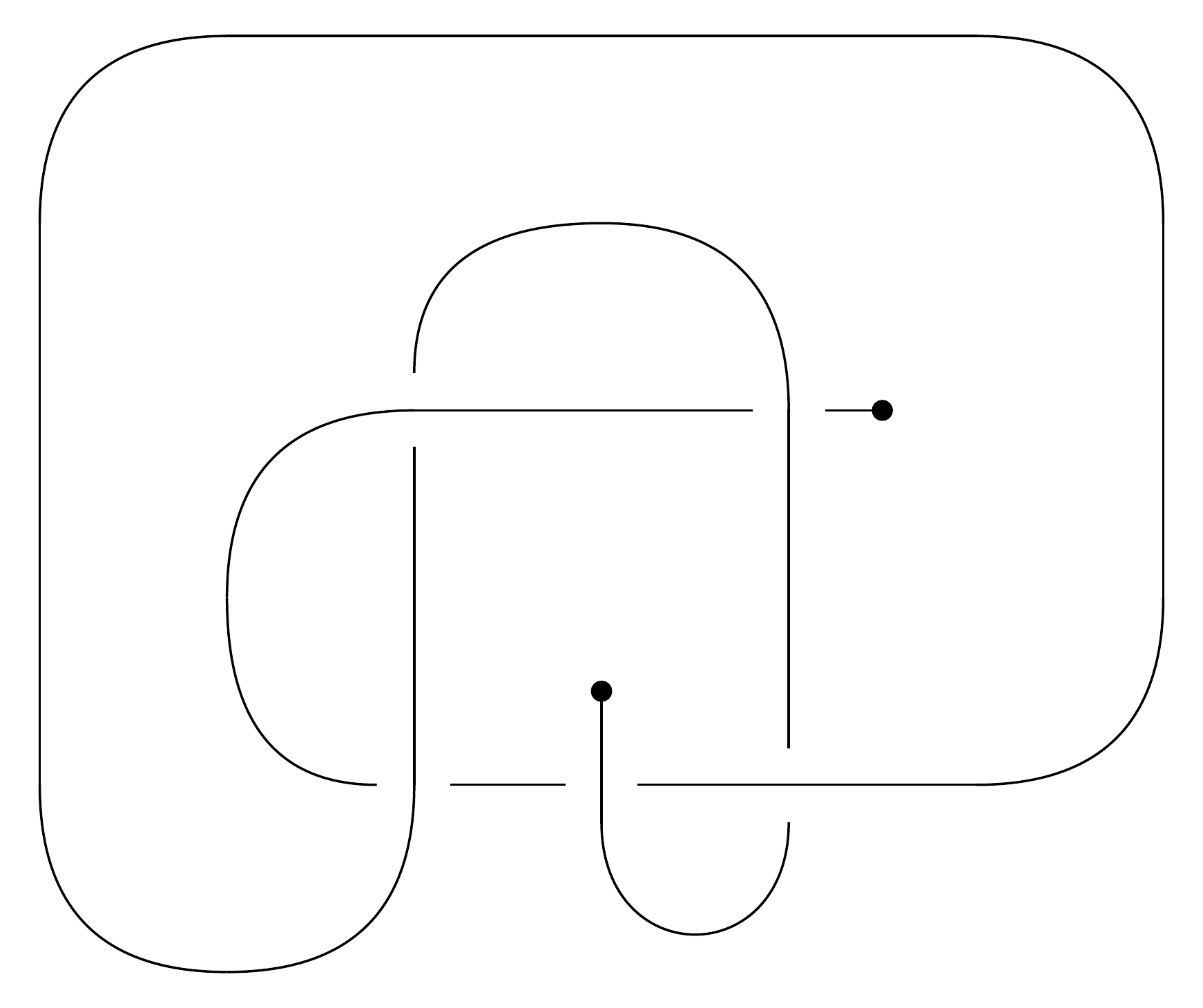}\\
\textcolor{black}{$5_{724}$}
\vspace{1cm}
\end{minipage}
\begin{minipage}[t]{.25\linewidth}
\centering
\includegraphics[width=0.9\textwidth,height=3.5cm,keepaspectratio]{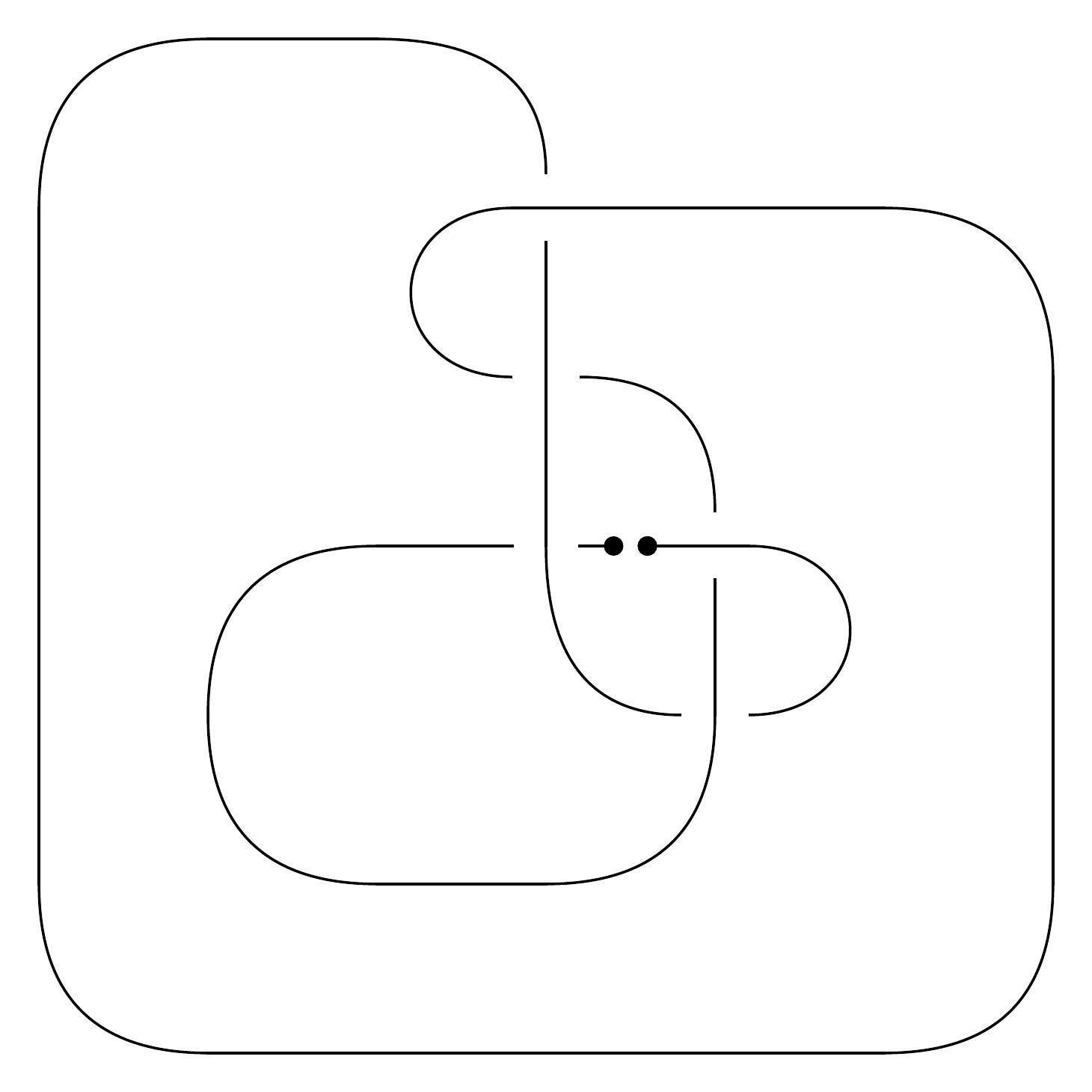}\\
\textcolor{black}{$5_{725}$}
\vspace{1cm}
\end{minipage}
\begin{minipage}[t]{.25\linewidth}
\centering
\includegraphics[width=0.9\textwidth,height=3.5cm,keepaspectratio]{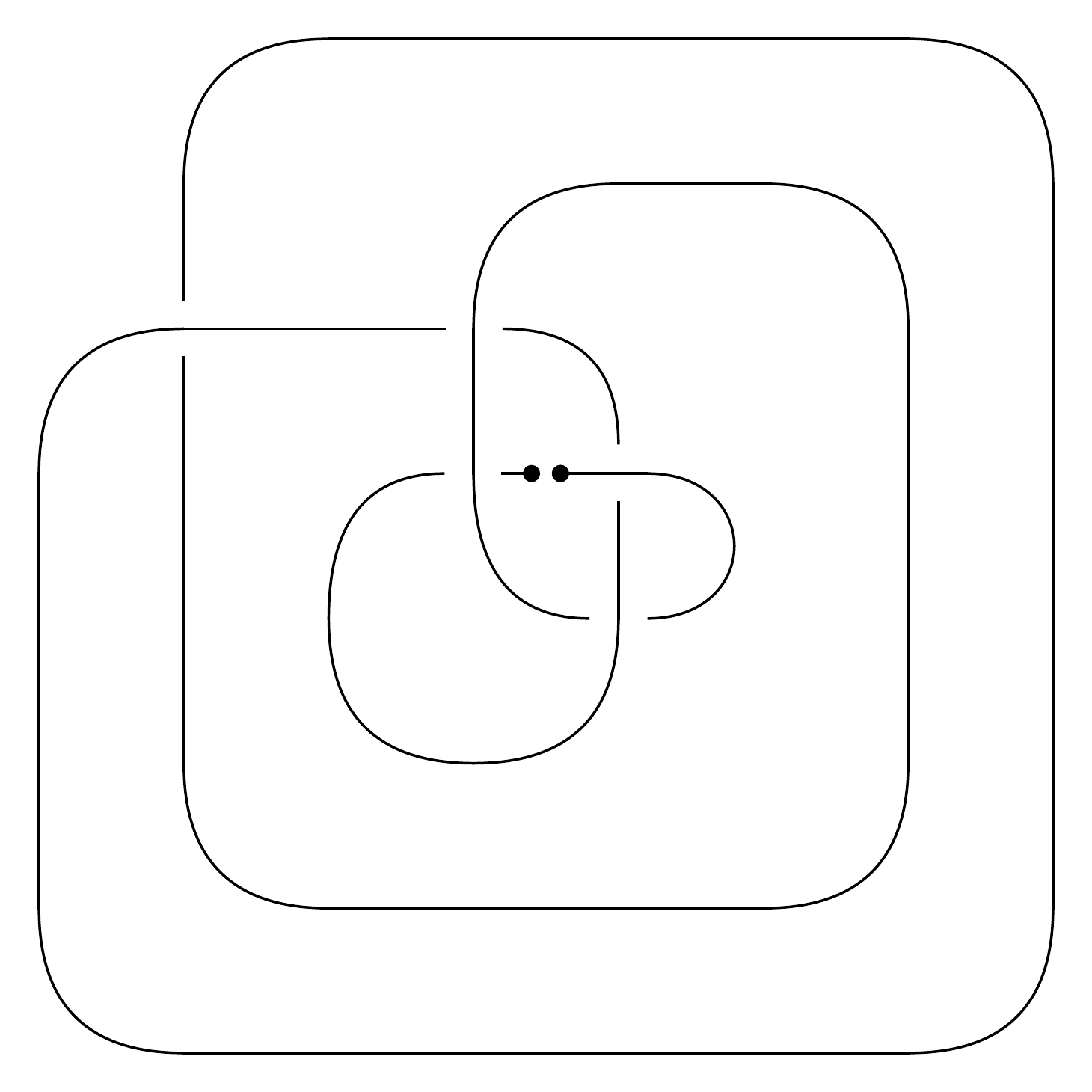}\\
\textcolor{black}{$5_{726}$}
\vspace{1cm}
\end{minipage}
\begin{minipage}[t]{.25\linewidth}
\centering
\includegraphics[width=0.9\textwidth,height=3.5cm,keepaspectratio]{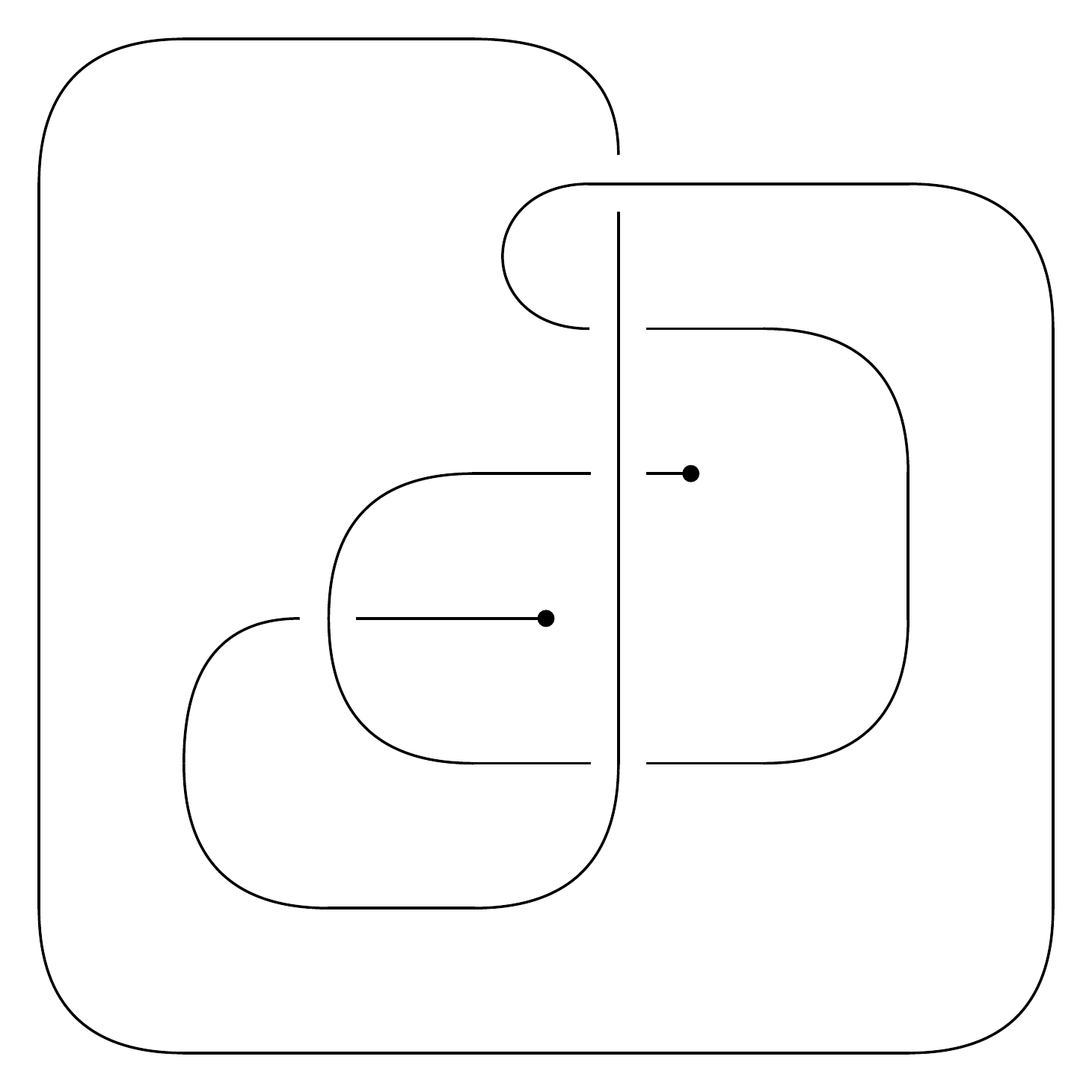}\\
\textcolor{black}{$5_{727}$}
\vspace{1cm}
\end{minipage}
\begin{minipage}[t]{.25\linewidth}
\centering
\includegraphics[width=0.9\textwidth,height=3.5cm,keepaspectratio]{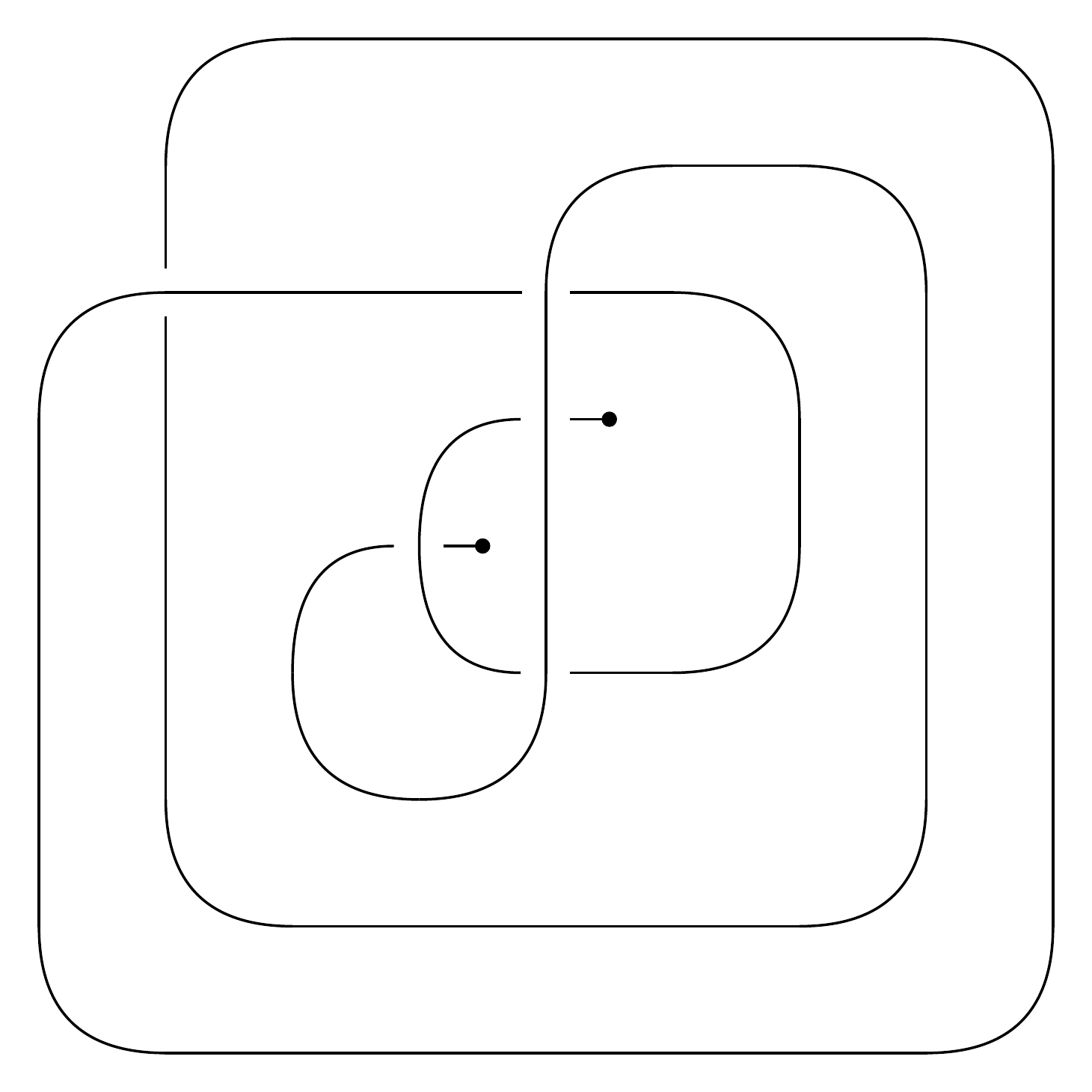}\\
\textcolor{black}{$5_{728}$}
\vspace{1cm}
\end{minipage}
\begin{minipage}[t]{.25\linewidth}
\centering
\includegraphics[width=0.9\textwidth,height=3.5cm,keepaspectratio]{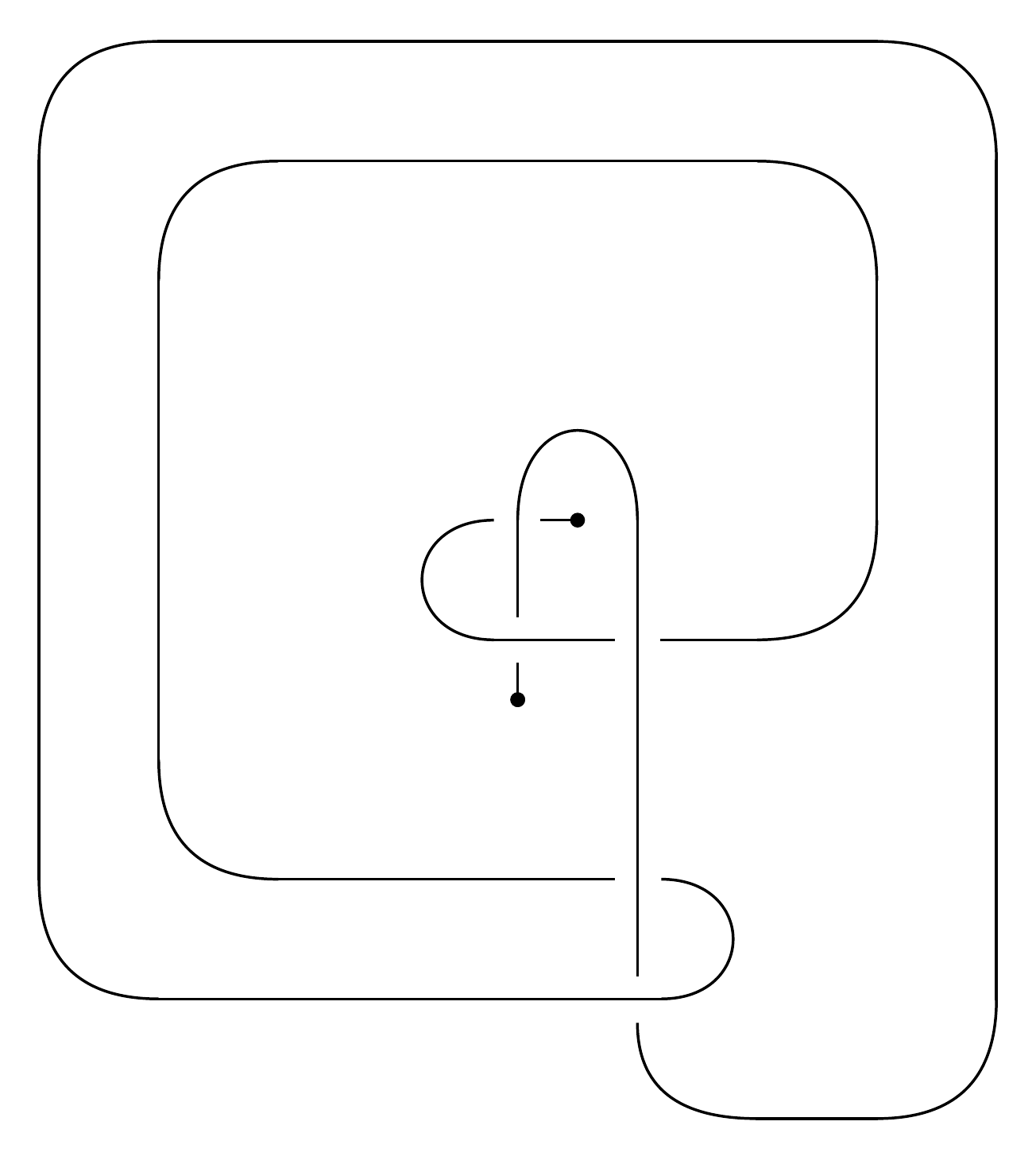}\\
\textcolor{black}{$5_{729}$}
\vspace{1cm}
\end{minipage}
\begin{minipage}[t]{.25\linewidth}
\centering
\includegraphics[width=0.9\textwidth,height=3.5cm,keepaspectratio]{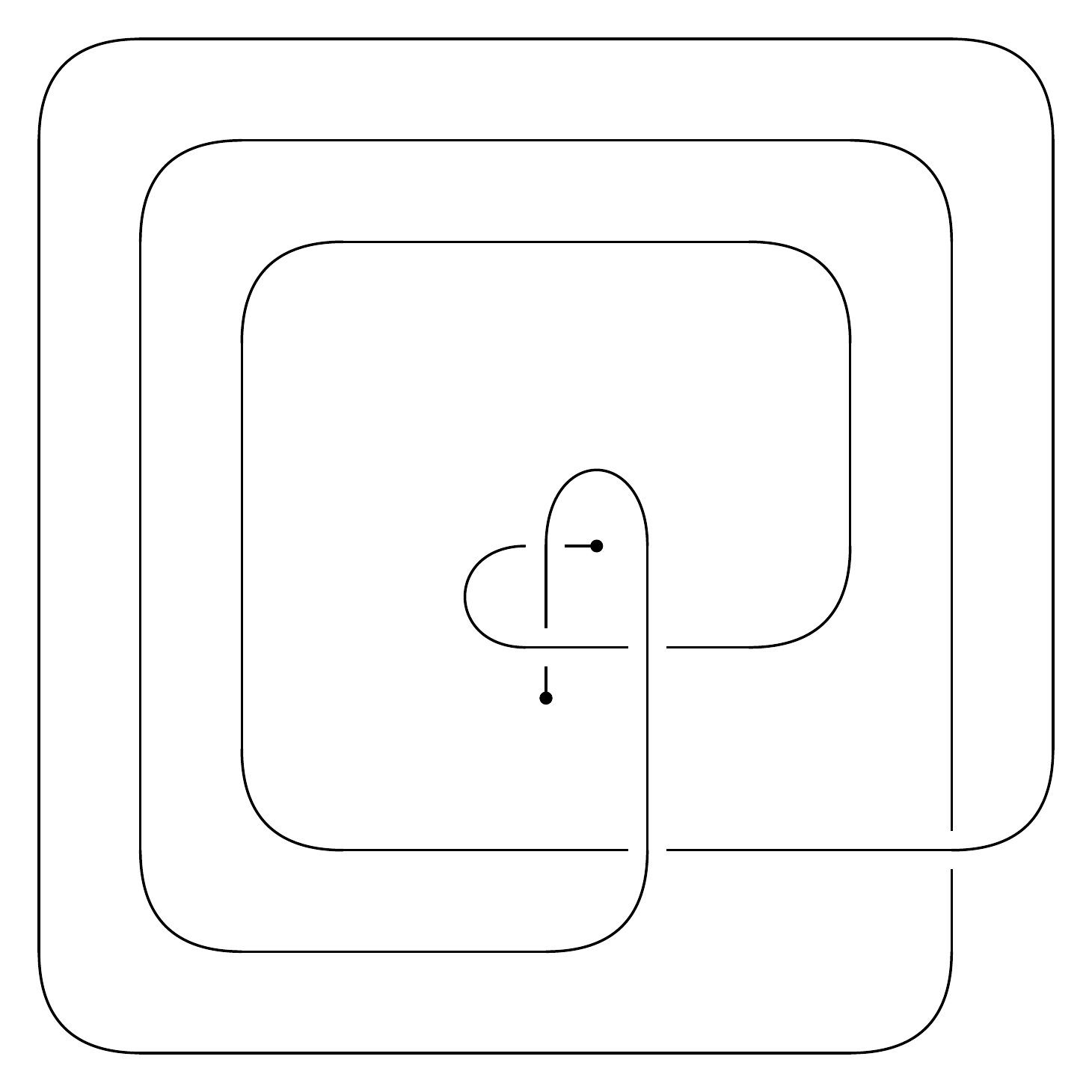}\\
\textcolor{black}{$5_{730}$}
\vspace{1cm}
\end{minipage}
\begin{minipage}[t]{.25\linewidth}
\centering
\includegraphics[width=0.9\textwidth,height=3.5cm,keepaspectratio]{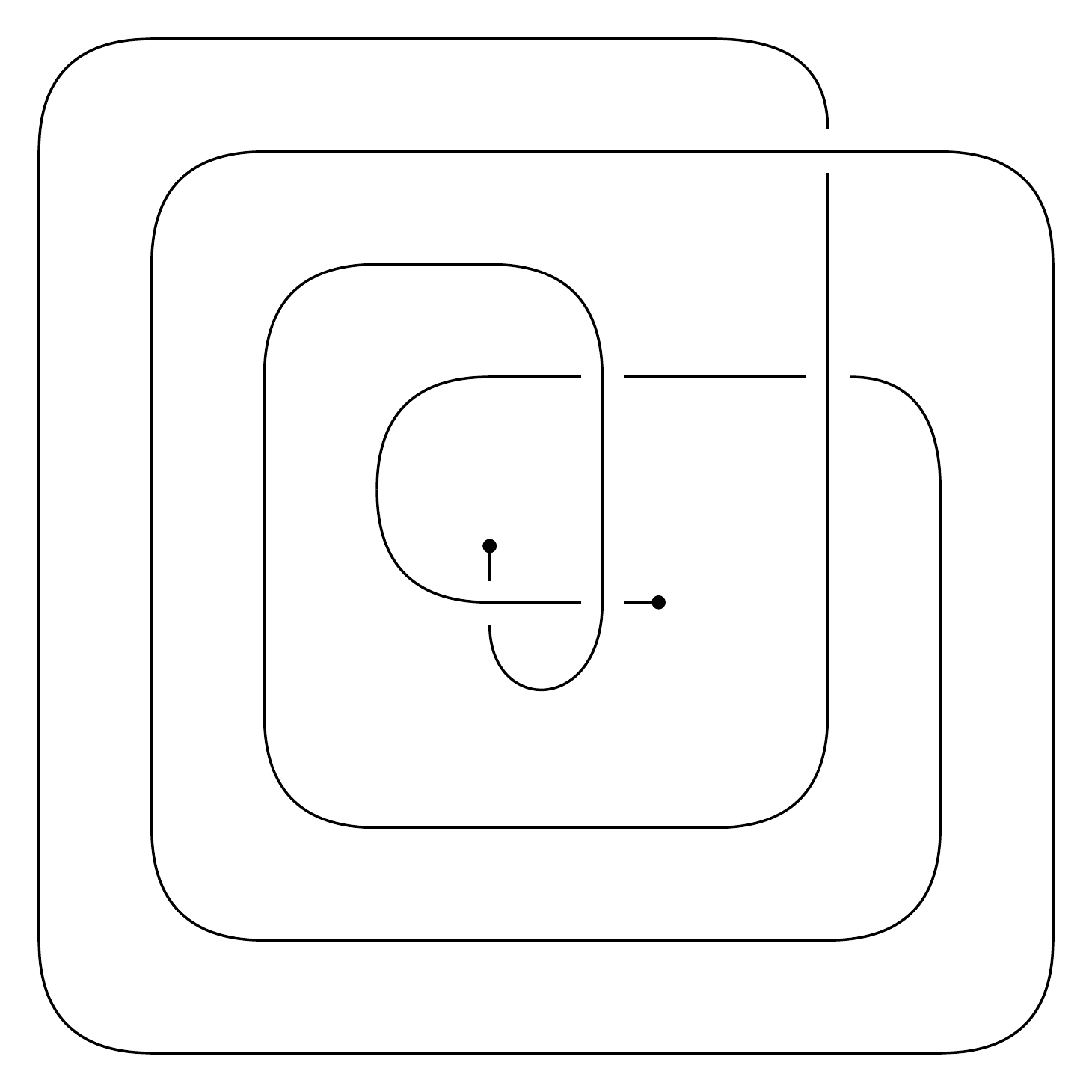}\\
\textcolor{black}{$5_{731}$}
\vspace{1cm}
\end{minipage}
\begin{minipage}[t]{.25\linewidth}
\centering
\includegraphics[width=0.9\textwidth,height=3.5cm,keepaspectratio]{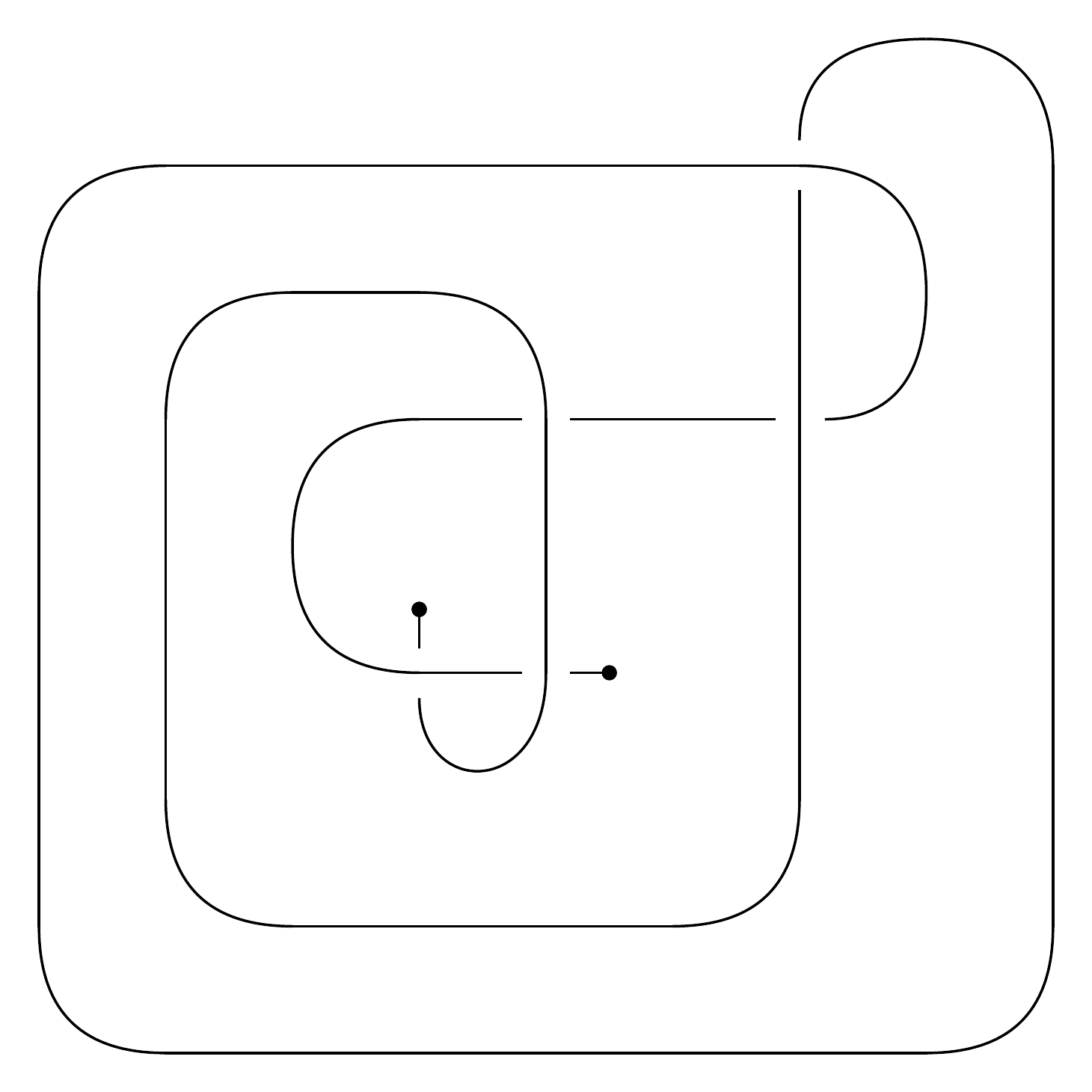}\\
\textcolor{black}{$5_{732}$}
\vspace{1cm}
\end{minipage}
\begin{minipage}[t]{.25\linewidth}
\centering
\includegraphics[width=0.9\textwidth,height=3.5cm,keepaspectratio]{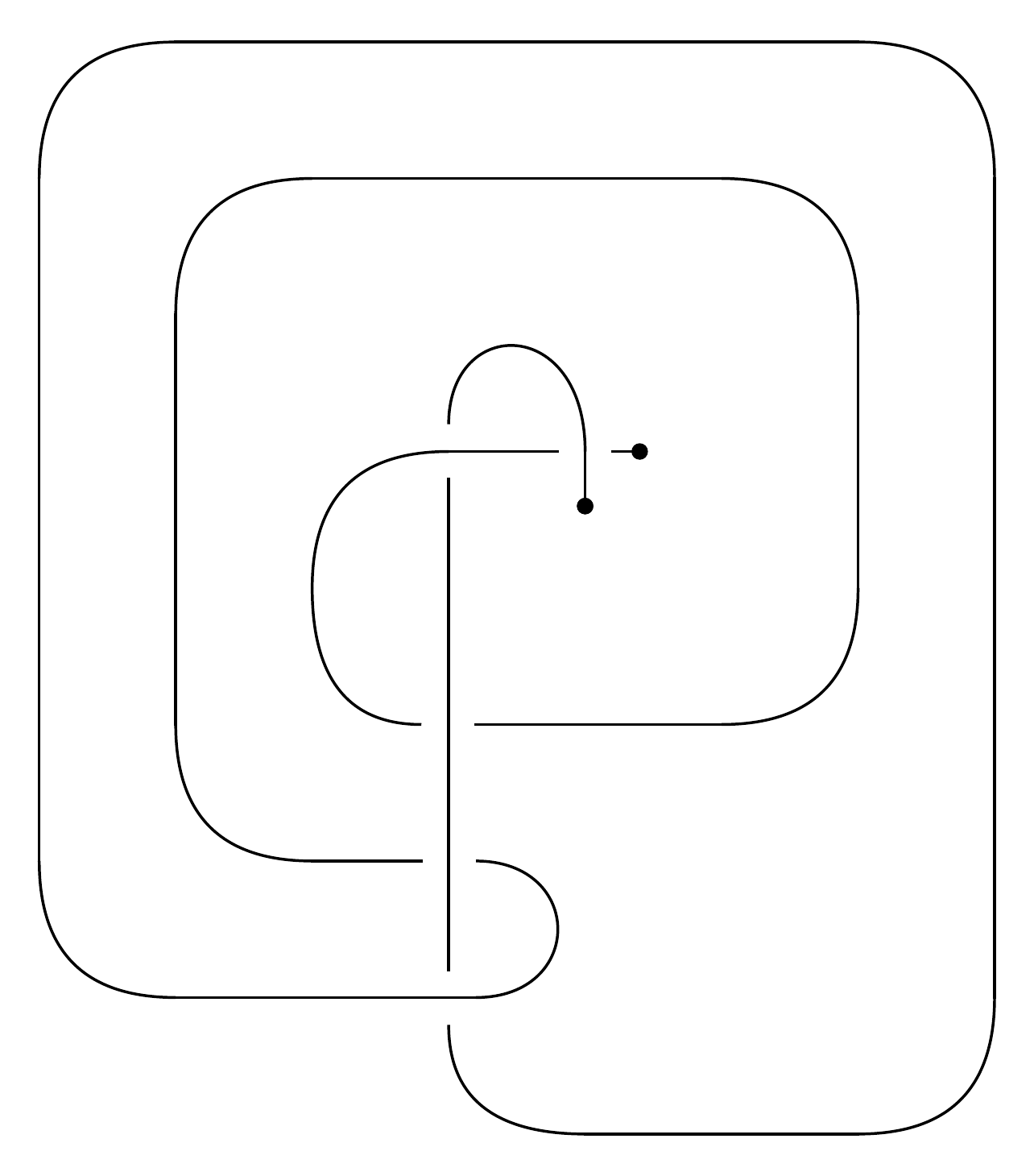}\\
\textcolor{black}{$5_{733}$}
\vspace{1cm}
\end{minipage}
\begin{minipage}[t]{.25\linewidth}
\centering
\includegraphics[width=0.9\textwidth,height=3.5cm,keepaspectratio]{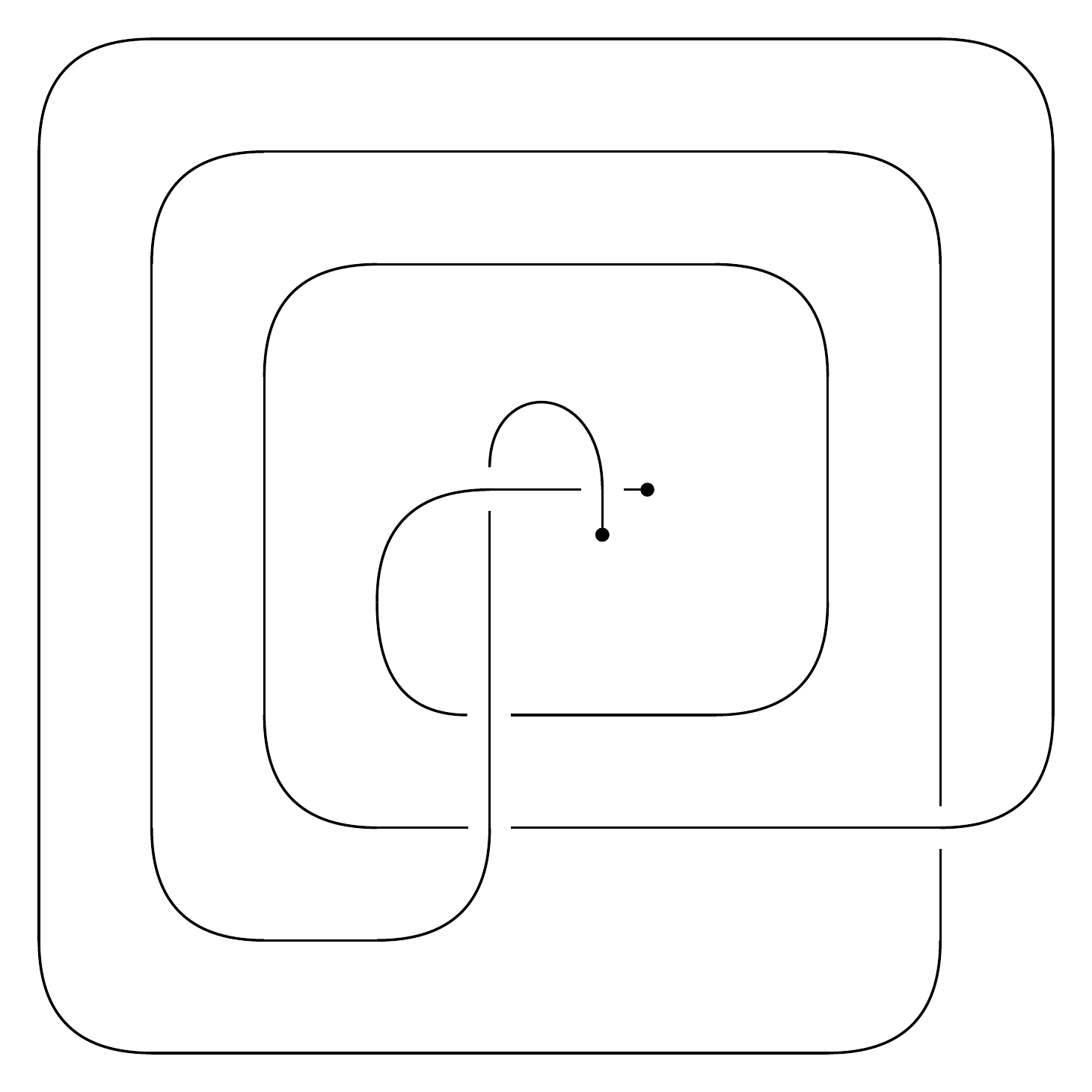}\\
\textcolor{black}{$5_{734}$}
\vspace{1cm}
\end{minipage}
\begin{minipage}[t]{.25\linewidth}
\centering
\includegraphics[width=0.9\textwidth,height=3.5cm,keepaspectratio]{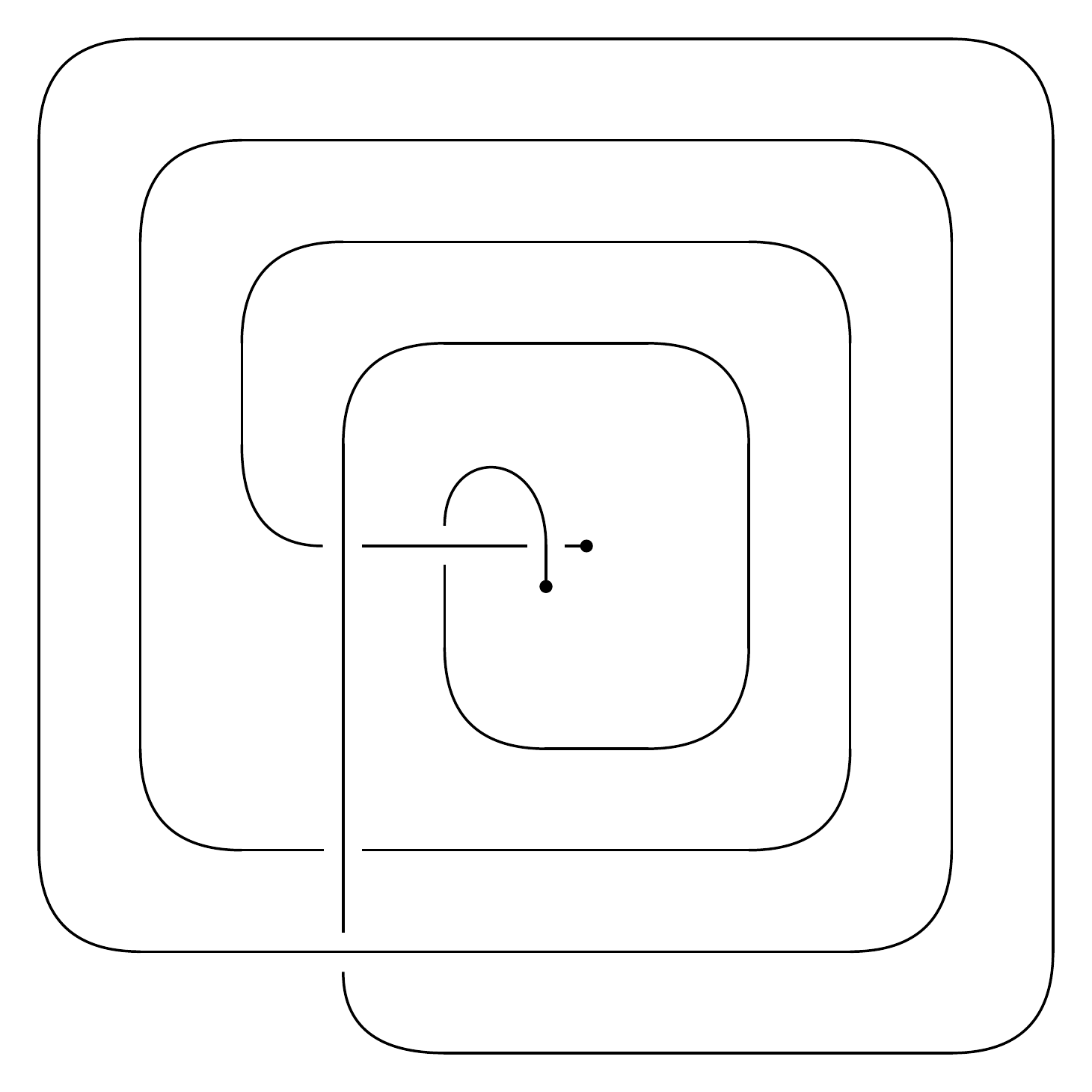}\\
\textcolor{black}{$5_{735}$}
\vspace{1cm}
\end{minipage}
\begin{minipage}[t]{.25\linewidth}
\centering
\includegraphics[width=0.9\textwidth,height=3.5cm,keepaspectratio]{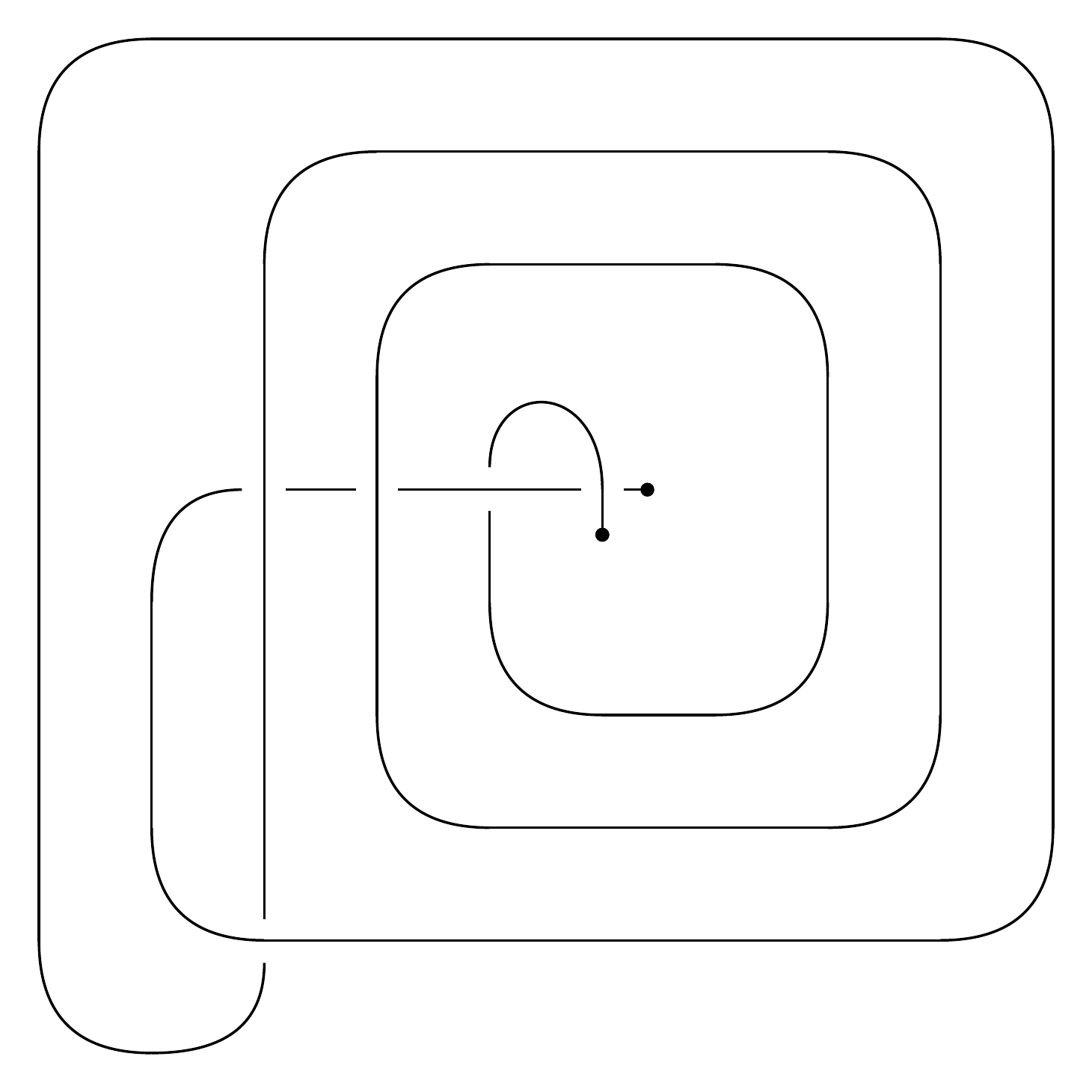}\\
\textcolor{black}{$5_{736}$}
\vspace{1cm}
\end{minipage}
\begin{minipage}[t]{.25\linewidth}
\centering
\includegraphics[width=0.9\textwidth,height=3.5cm,keepaspectratio]{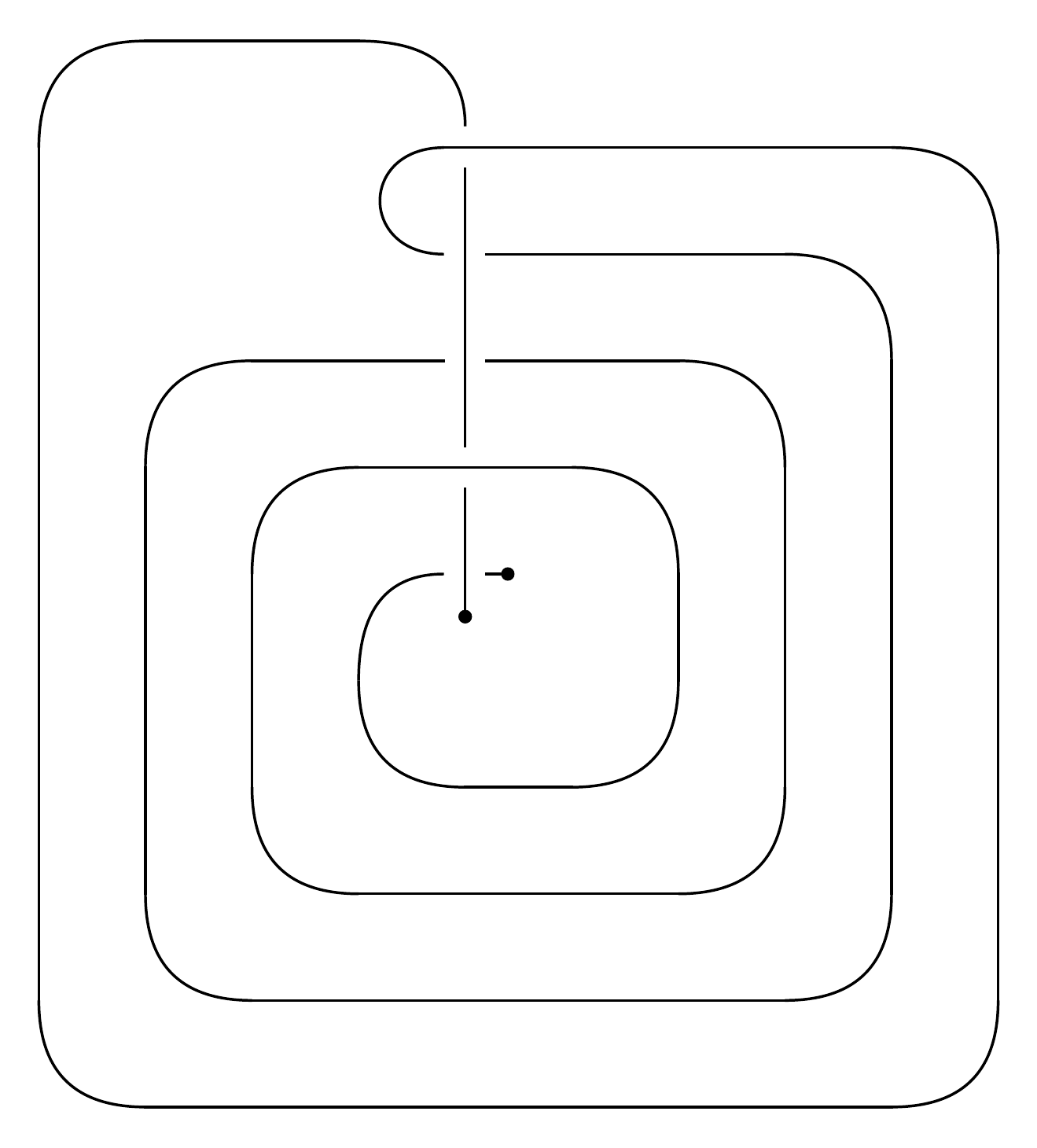}\\
\textcolor{black}{$5_{737}$}
\vspace{1cm}
\end{minipage}
\begin{minipage}[t]{.25\linewidth}
\centering
\includegraphics[width=0.9\textwidth,height=3.5cm,keepaspectratio]{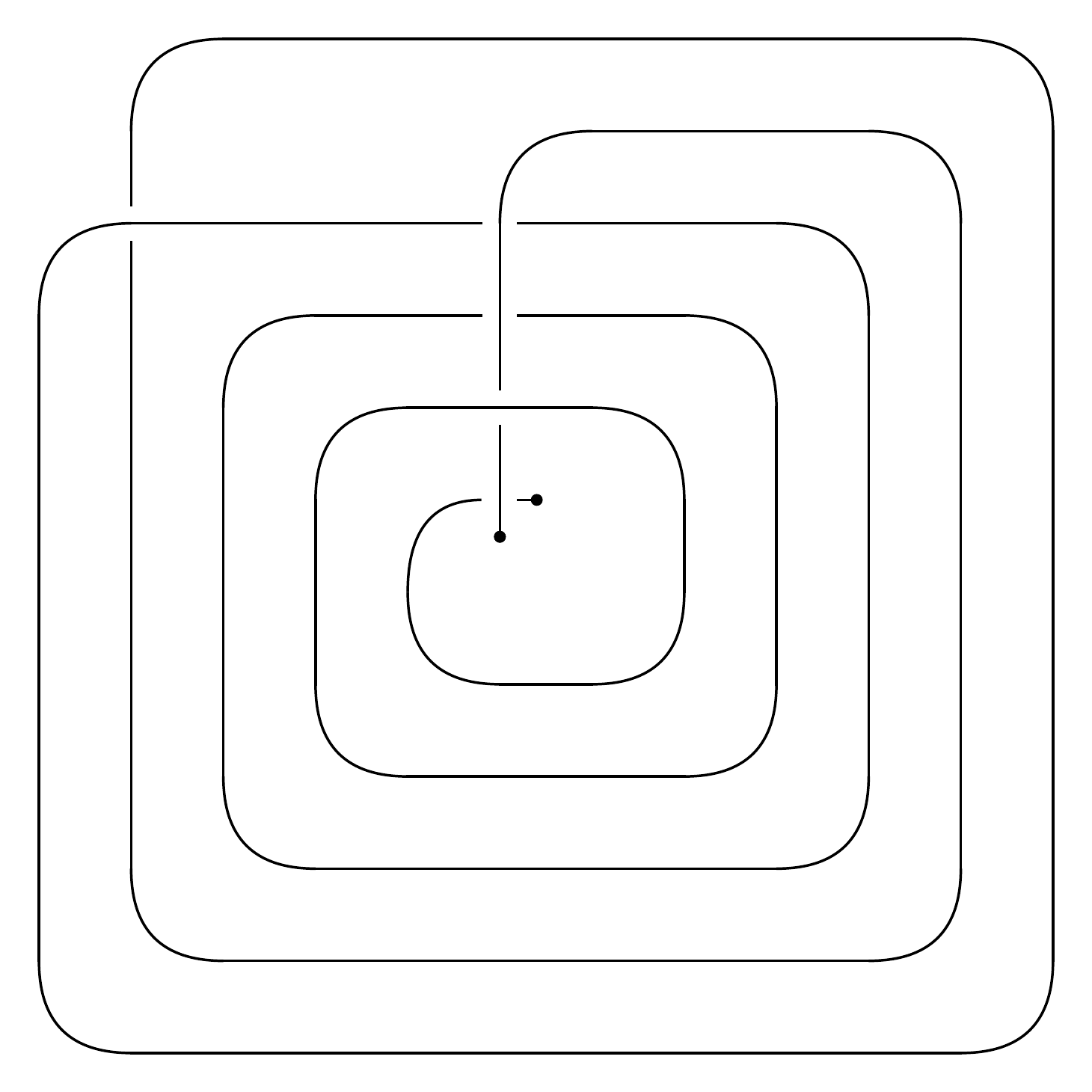}\\
\textcolor{black}{$5_{738}$}
\vspace{1cm}
\end{minipage}
\begin{minipage}[t]{.25\linewidth}
\centering
\includegraphics[width=0.9\textwidth,height=3.5cm,keepaspectratio]{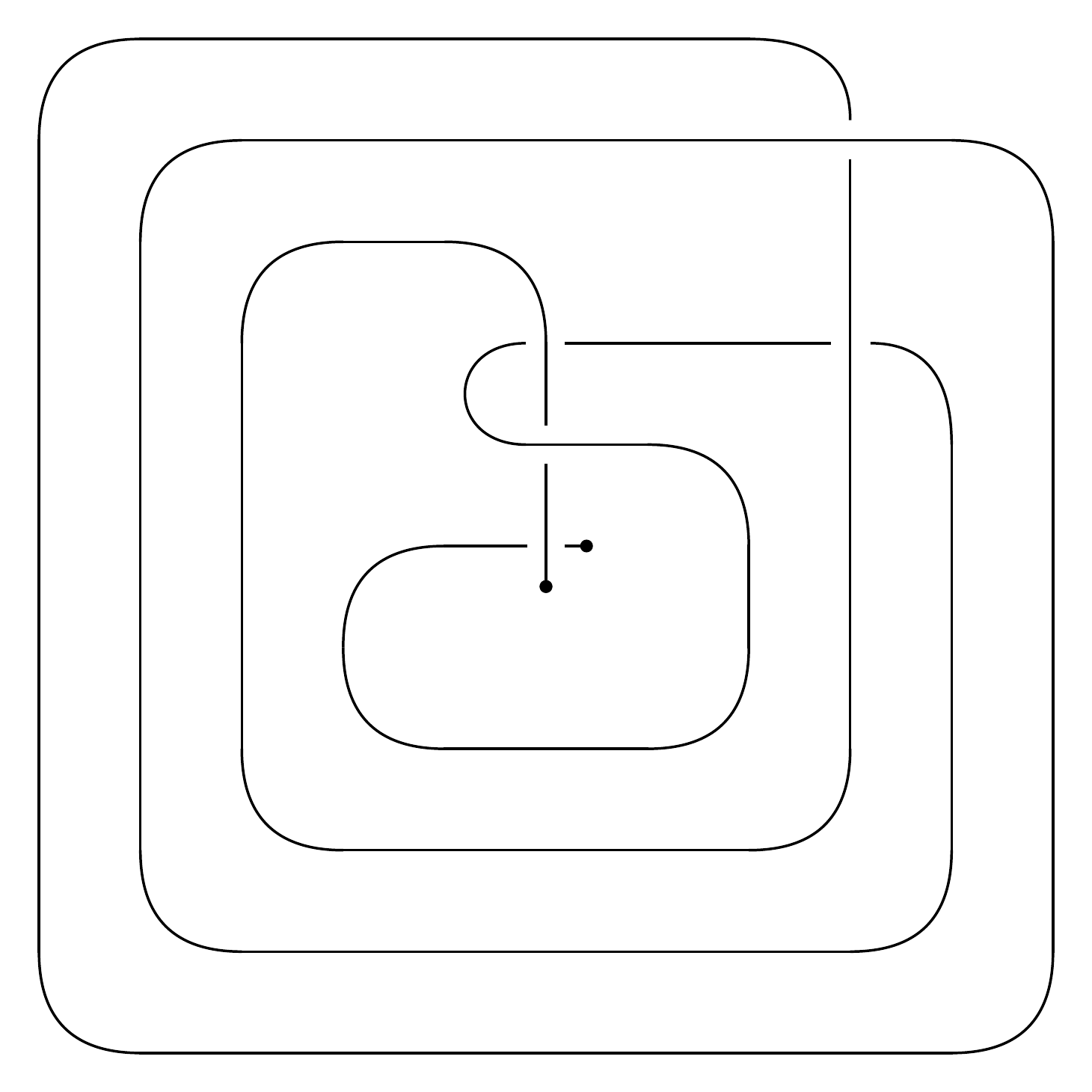}\\
\textcolor{black}{$5_{739}$}
\vspace{1cm}
\end{minipage}
\begin{minipage}[t]{.25\linewidth}
\centering
\includegraphics[width=0.9\textwidth,height=3.5cm,keepaspectratio]{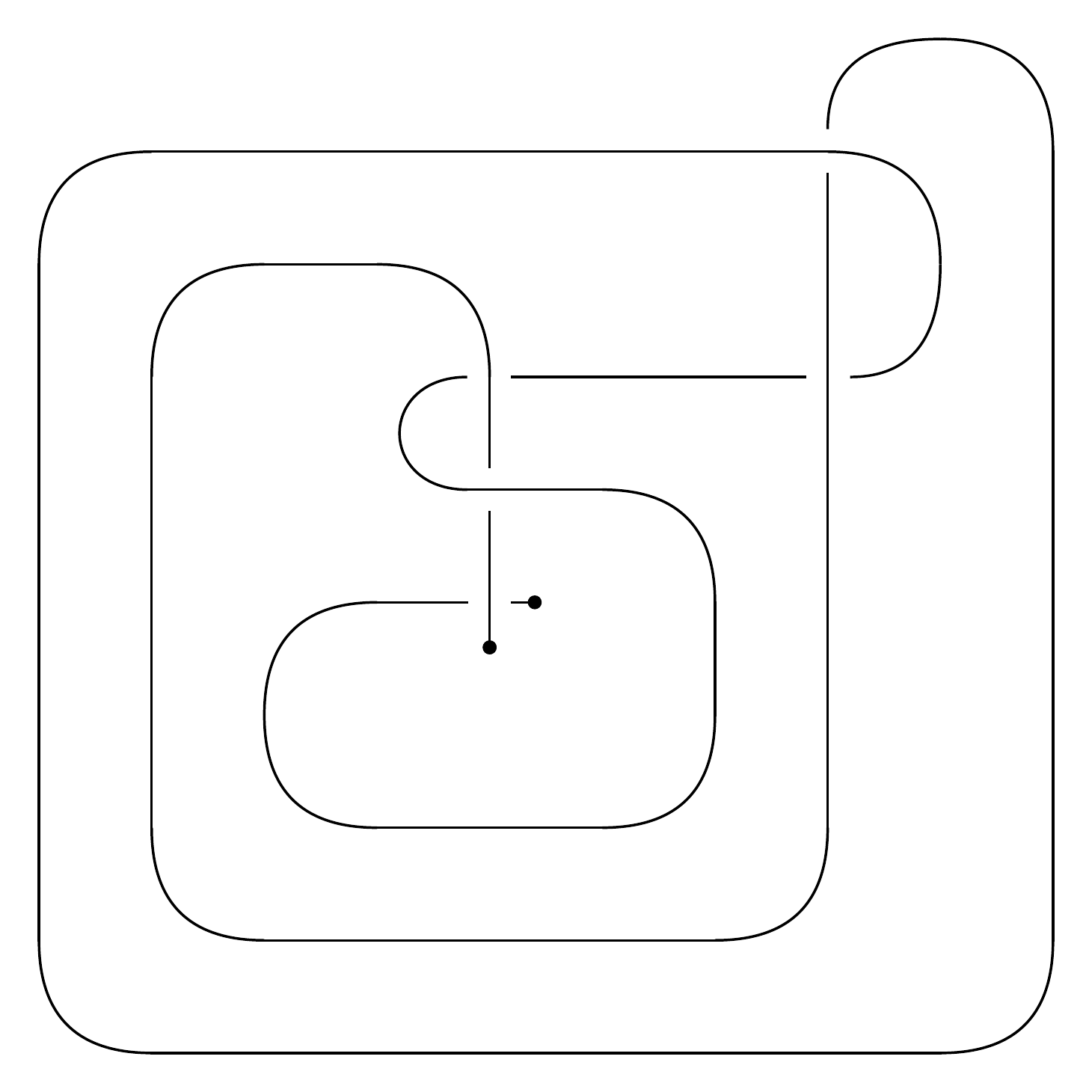}\\
\textcolor{black}{$5_{740}$}
\vspace{1cm}
\end{minipage}
\begin{minipage}[t]{.25\linewidth}
\centering
\includegraphics[width=0.9\textwidth,height=3.5cm,keepaspectratio]{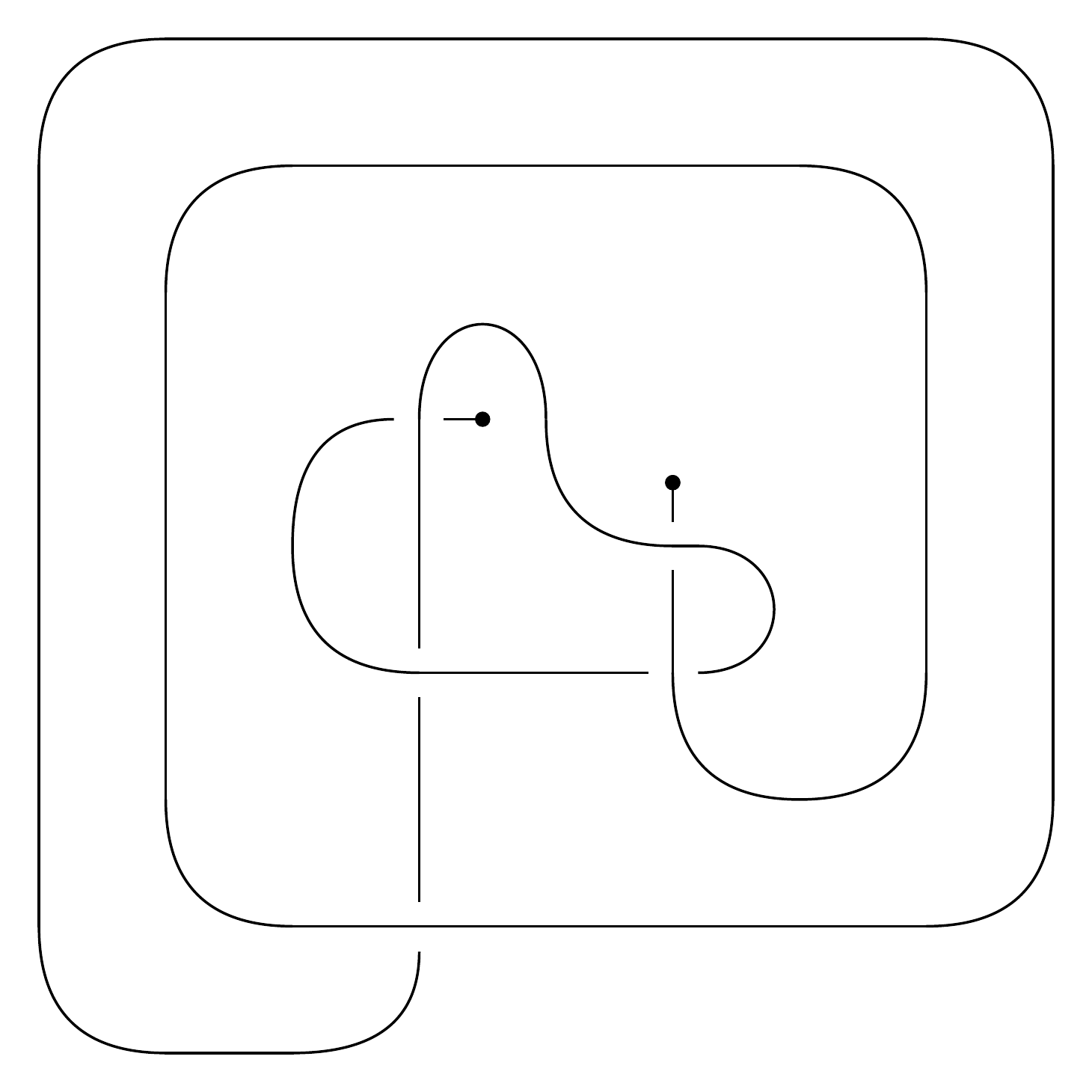}\\
\textcolor{black}{$5_{741}$}
\vspace{1cm}
\end{minipage}
\begin{minipage}[t]{.25\linewidth}
\centering
\includegraphics[width=0.9\textwidth,height=3.5cm,keepaspectratio]{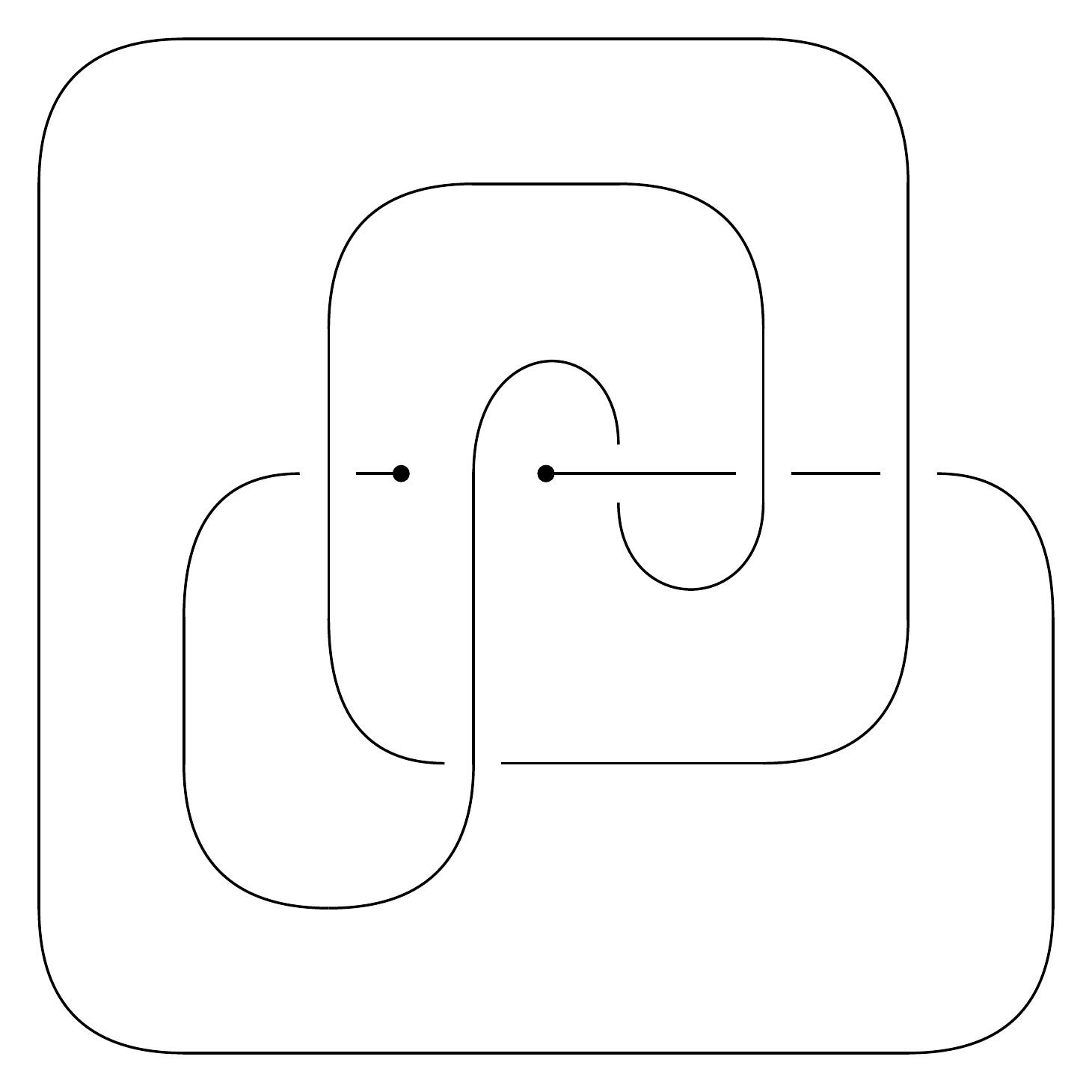}\\
\textcolor{black}{$5_{742}$}
\vspace{1cm}
\end{minipage}
\begin{minipage}[t]{.25\linewidth}
\centering
\includegraphics[width=0.9\textwidth,height=3.5cm,keepaspectratio]{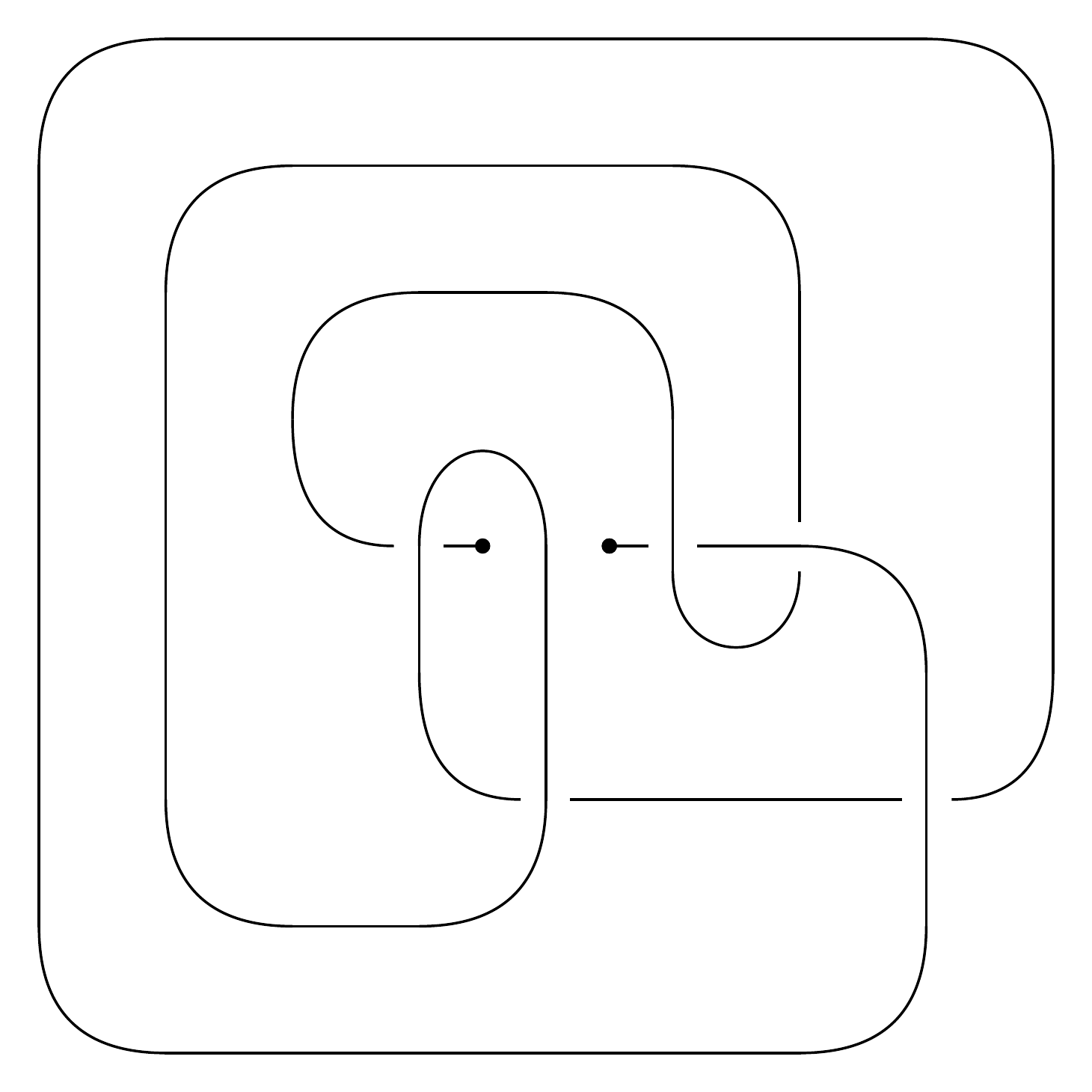}\\
\textcolor{black}{$5_{743}$}
\vspace{1cm}
\end{minipage}
\begin{minipage}[t]{.25\linewidth}
\centering
\includegraphics[width=0.9\textwidth,height=3.5cm,keepaspectratio]{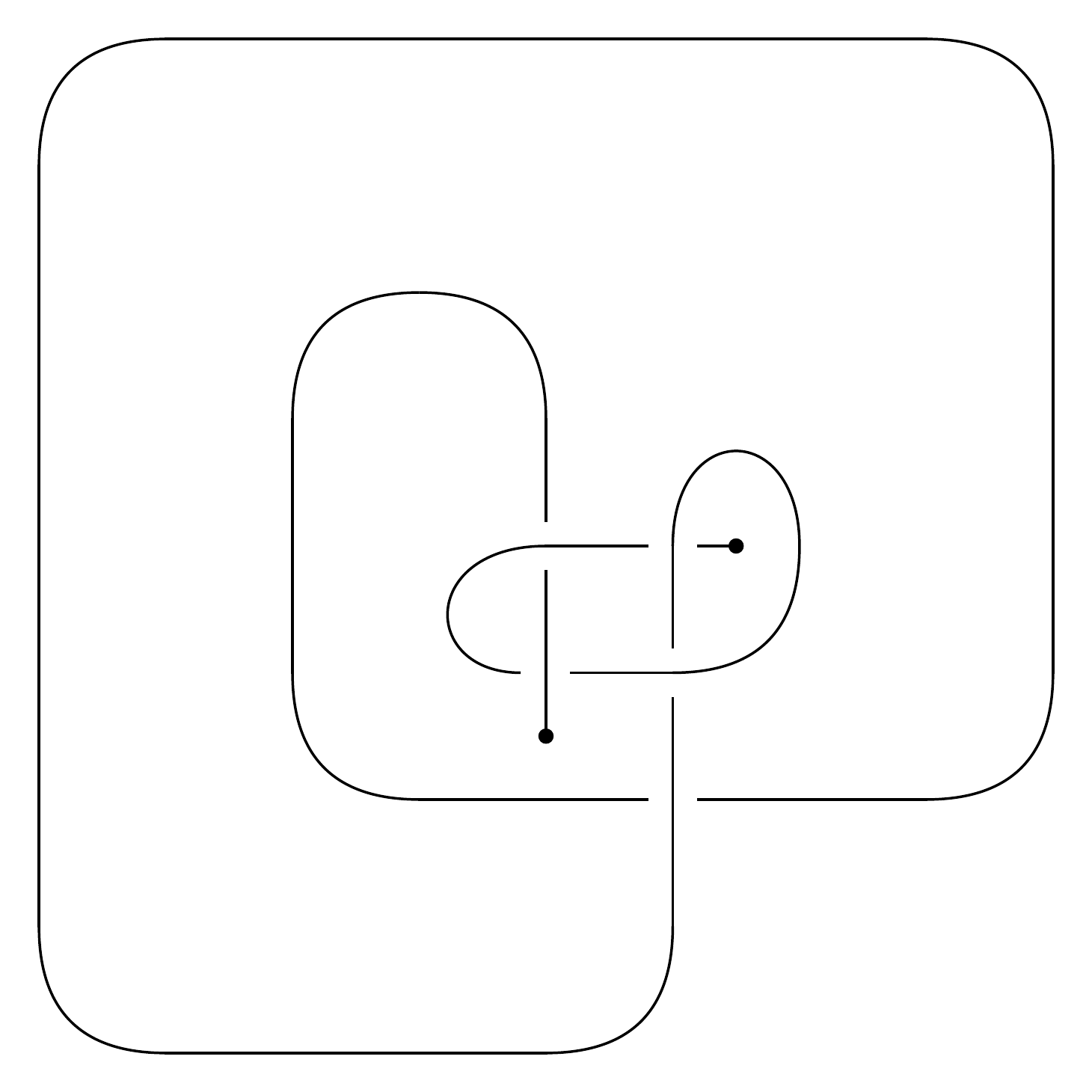}\\
\textcolor{black}{$5_{744}$}
\vspace{1cm}
\end{minipage}
\begin{minipage}[t]{.25\linewidth}
\centering
\includegraphics[width=0.9\textwidth,height=3.5cm,keepaspectratio]{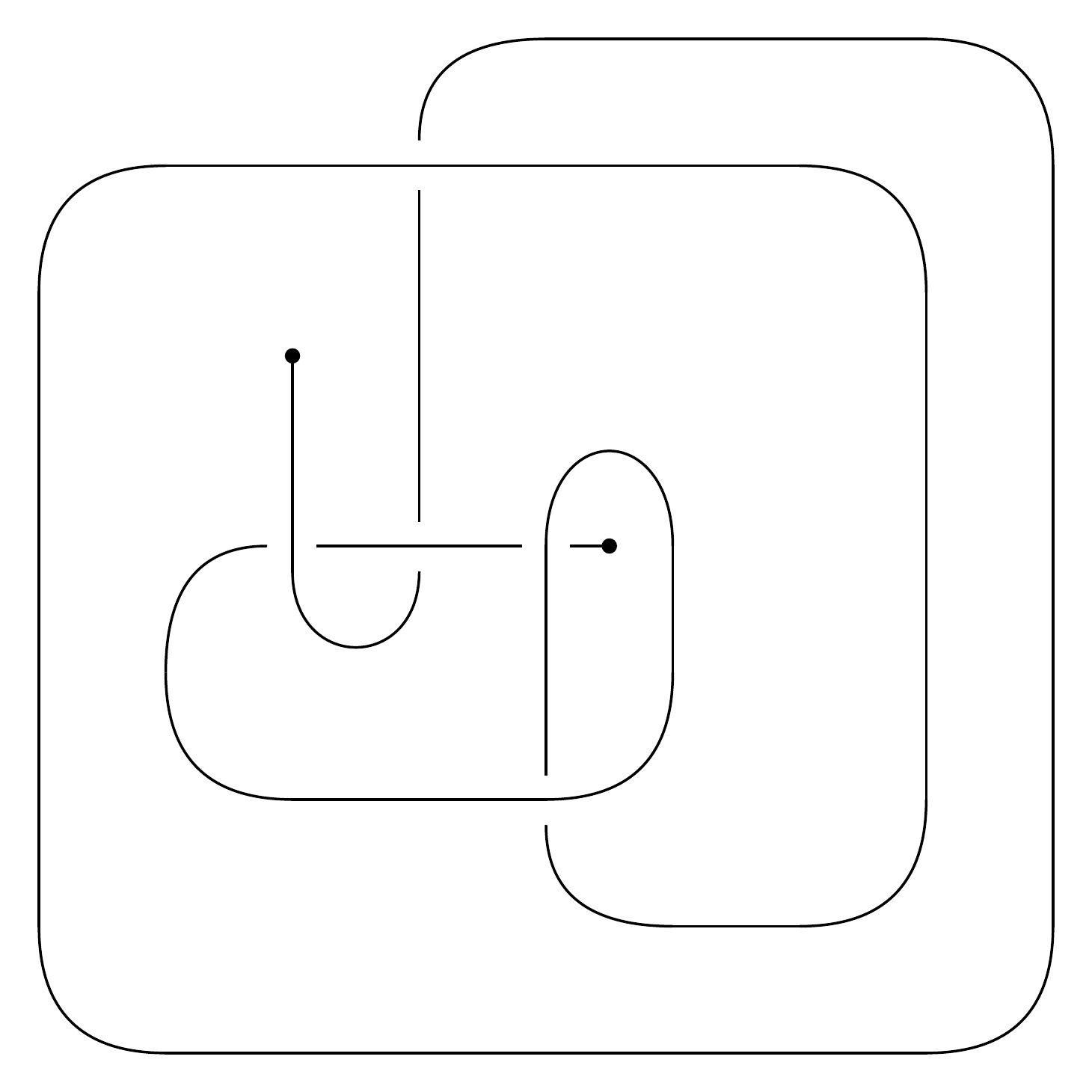}\\
\textcolor{black}{$5_{745}$}
\vspace{1cm}
\end{minipage}
\begin{minipage}[t]{.25\linewidth}
\centering
\includegraphics[width=0.9\textwidth,height=3.5cm,keepaspectratio]{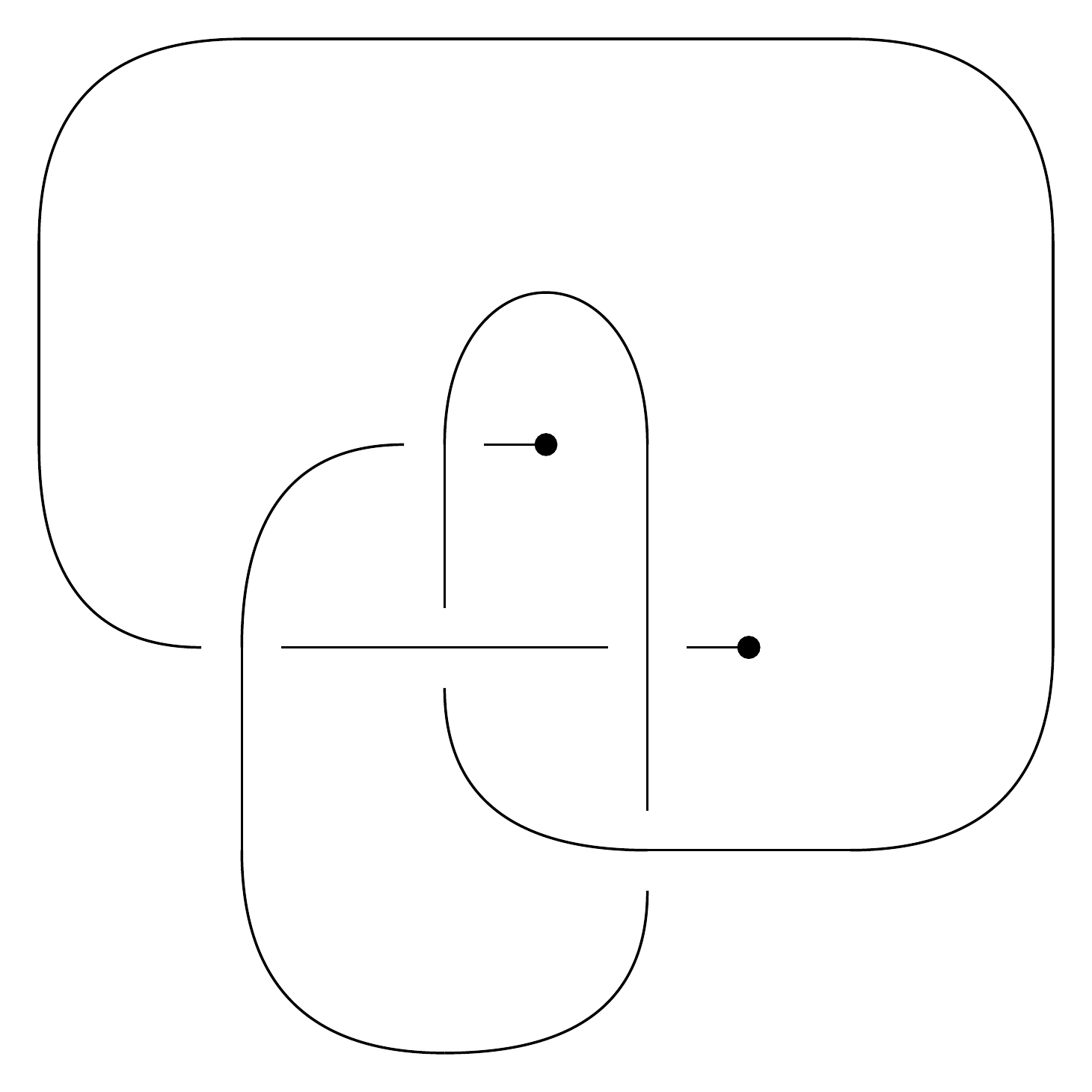}\\
\textcolor{black}{$5_{746}$}
\vspace{1cm}
\end{minipage}
\begin{minipage}[t]{.25\linewidth}
\centering
\includegraphics[width=0.9\textwidth,height=3.5cm,keepaspectratio]{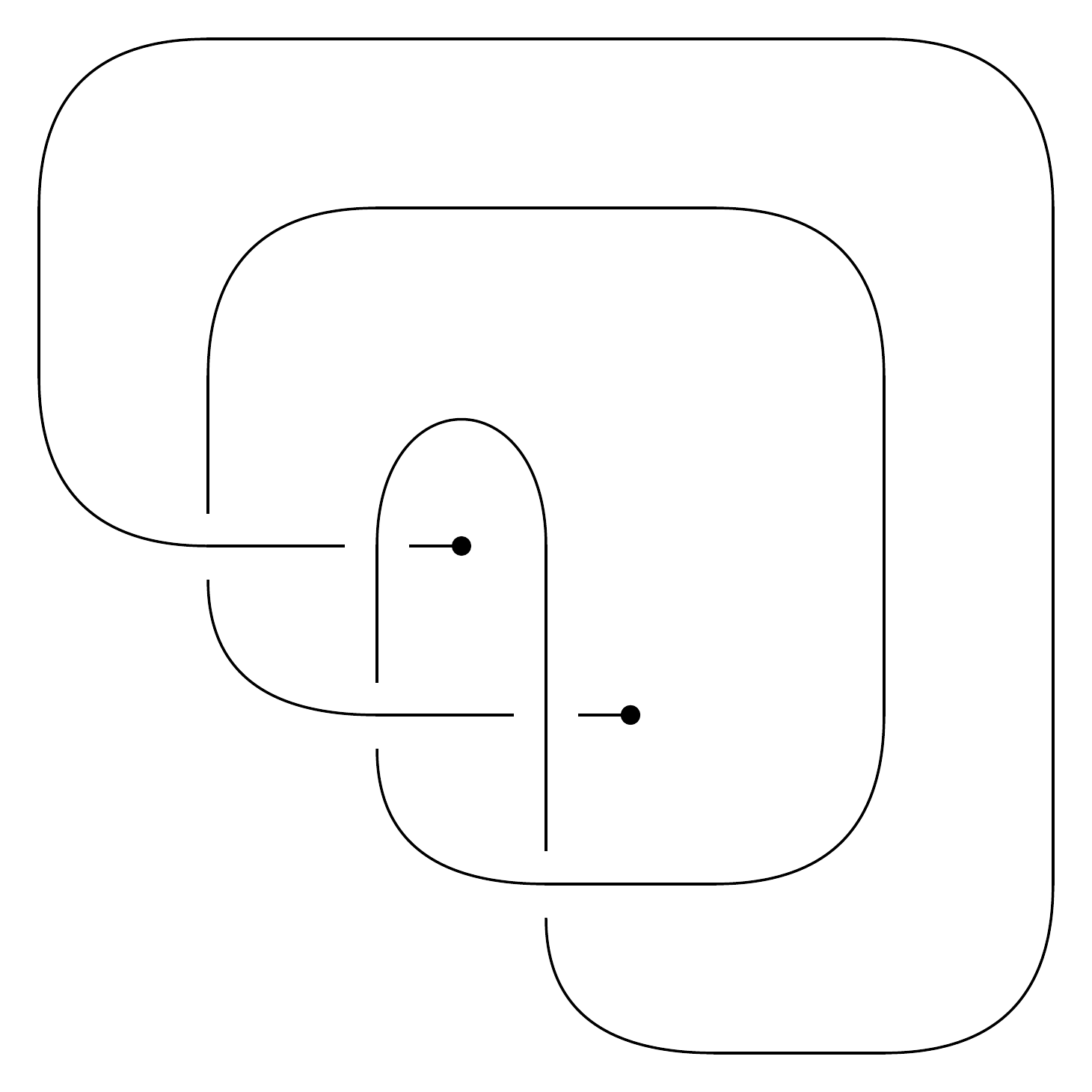}\\
\textcolor{black}{$5_{747}$}
\vspace{1cm}
\end{minipage}
\begin{minipage}[t]{.25\linewidth}
\centering
\includegraphics[width=0.9\textwidth,height=3.5cm,keepaspectratio]{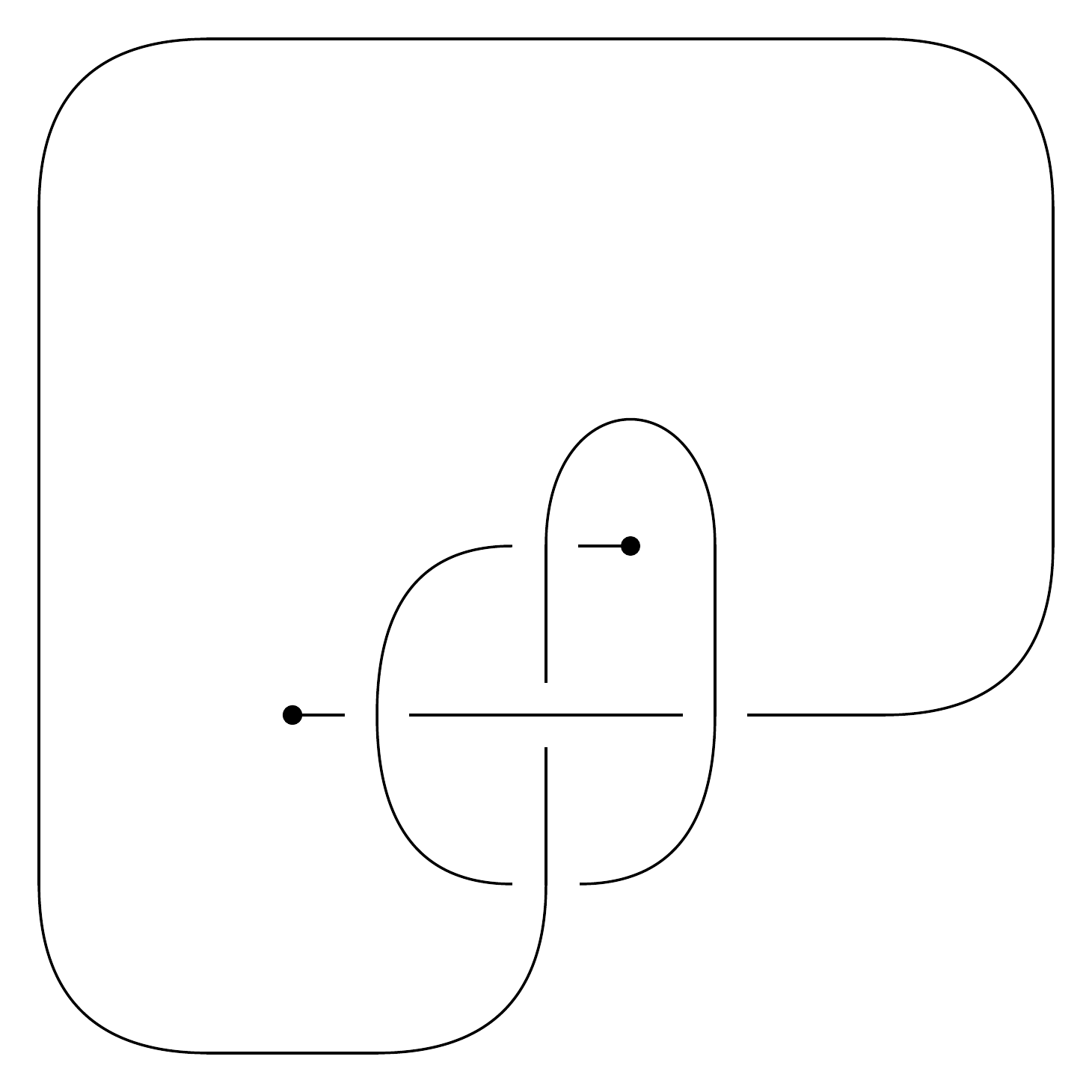}\\
\textcolor{black}{$5_{748}$}
\vspace{1cm}
\end{minipage}
\begin{minipage}[t]{.25\linewidth}
\centering
\includegraphics[width=0.9\textwidth,height=3.5cm,keepaspectratio]{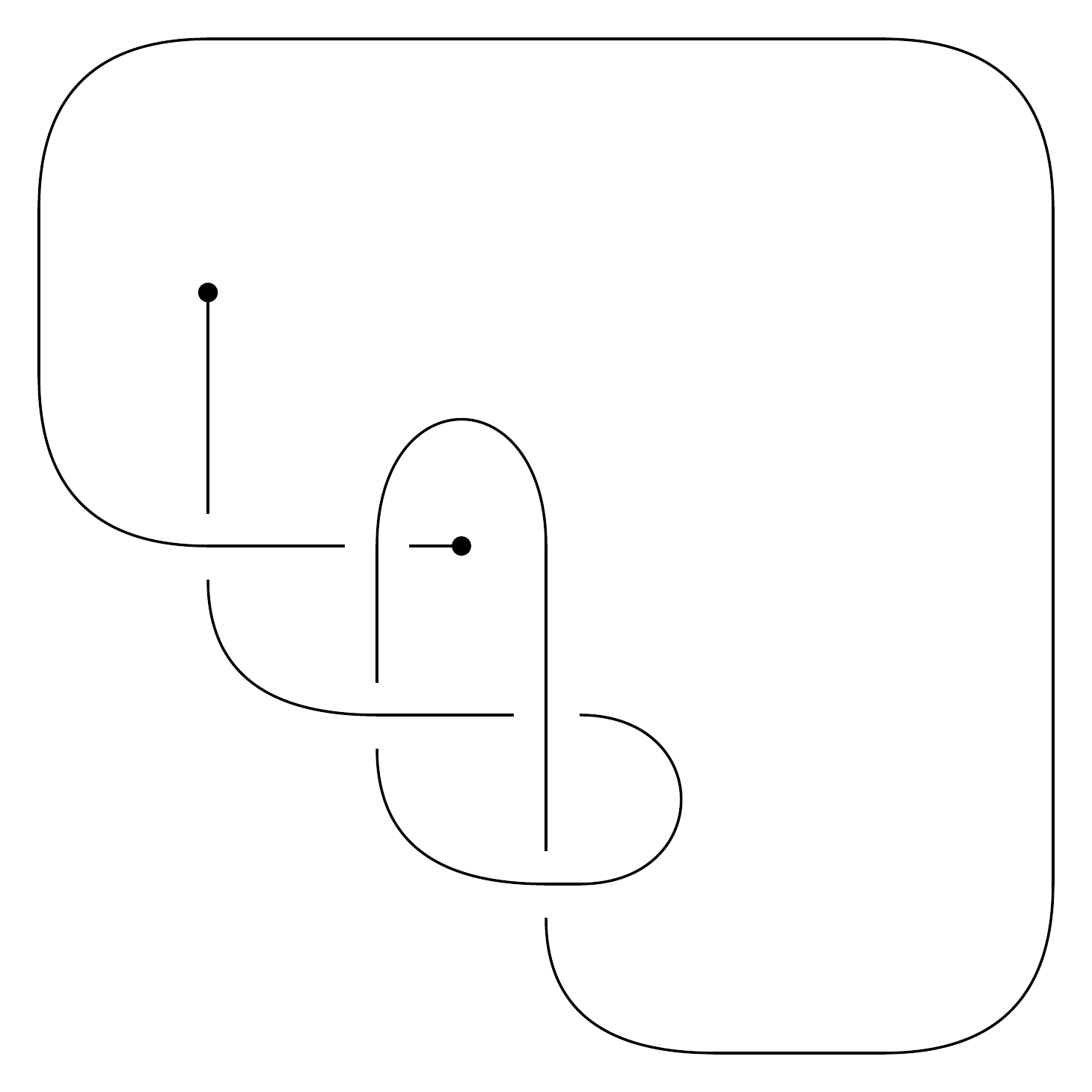}\\
\textcolor{black}{$5_{749}$}
\vspace{1cm}
\end{minipage}
\begin{minipage}[t]{.25\linewidth}
\centering
\includegraphics[width=0.9\textwidth,height=3.5cm,keepaspectratio]{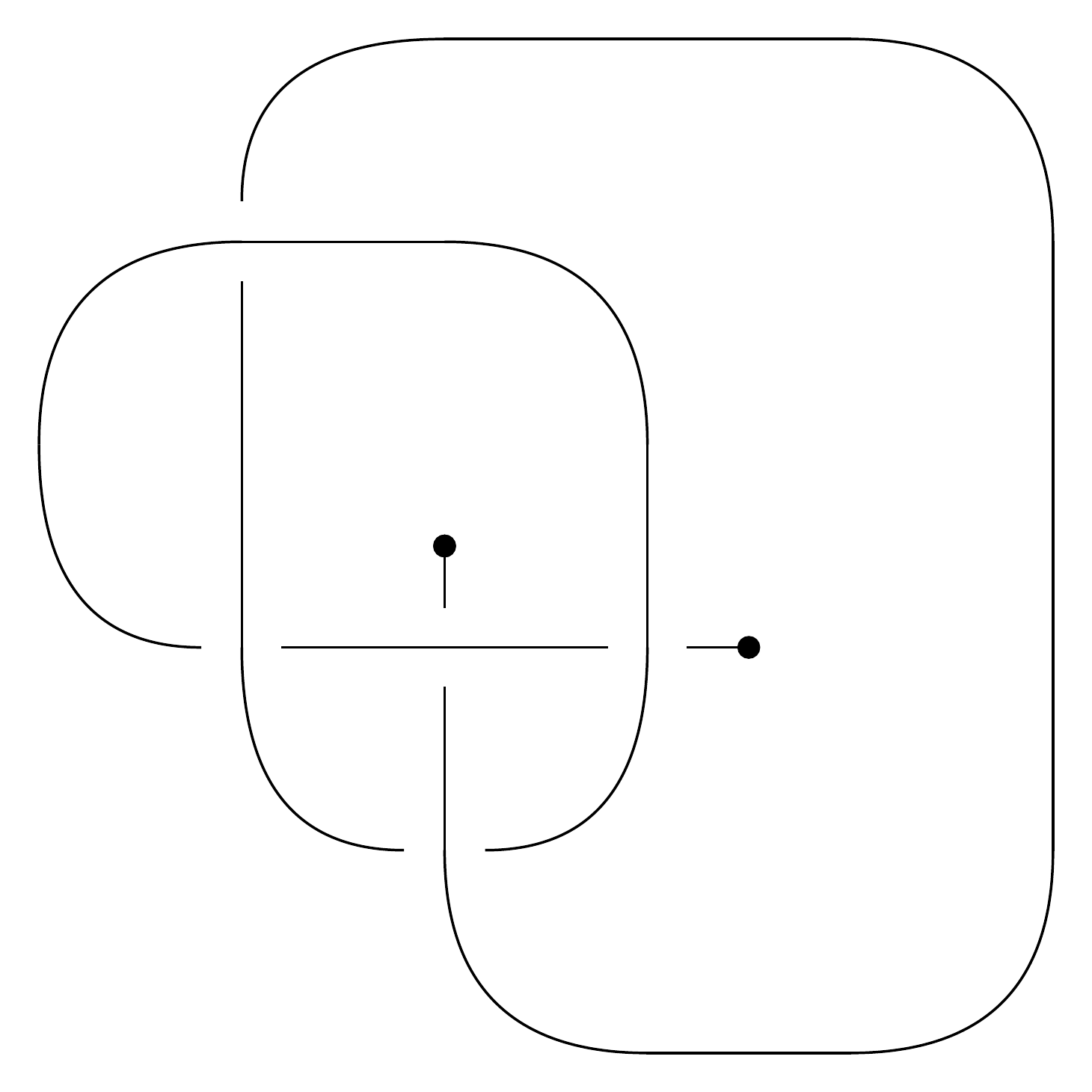}\\
\textcolor{black}{$5_{750}$}
\vspace{1cm}
\end{minipage}
\begin{minipage}[t]{.25\linewidth}
\centering
\includegraphics[width=0.9\textwidth,height=3.5cm,keepaspectratio]{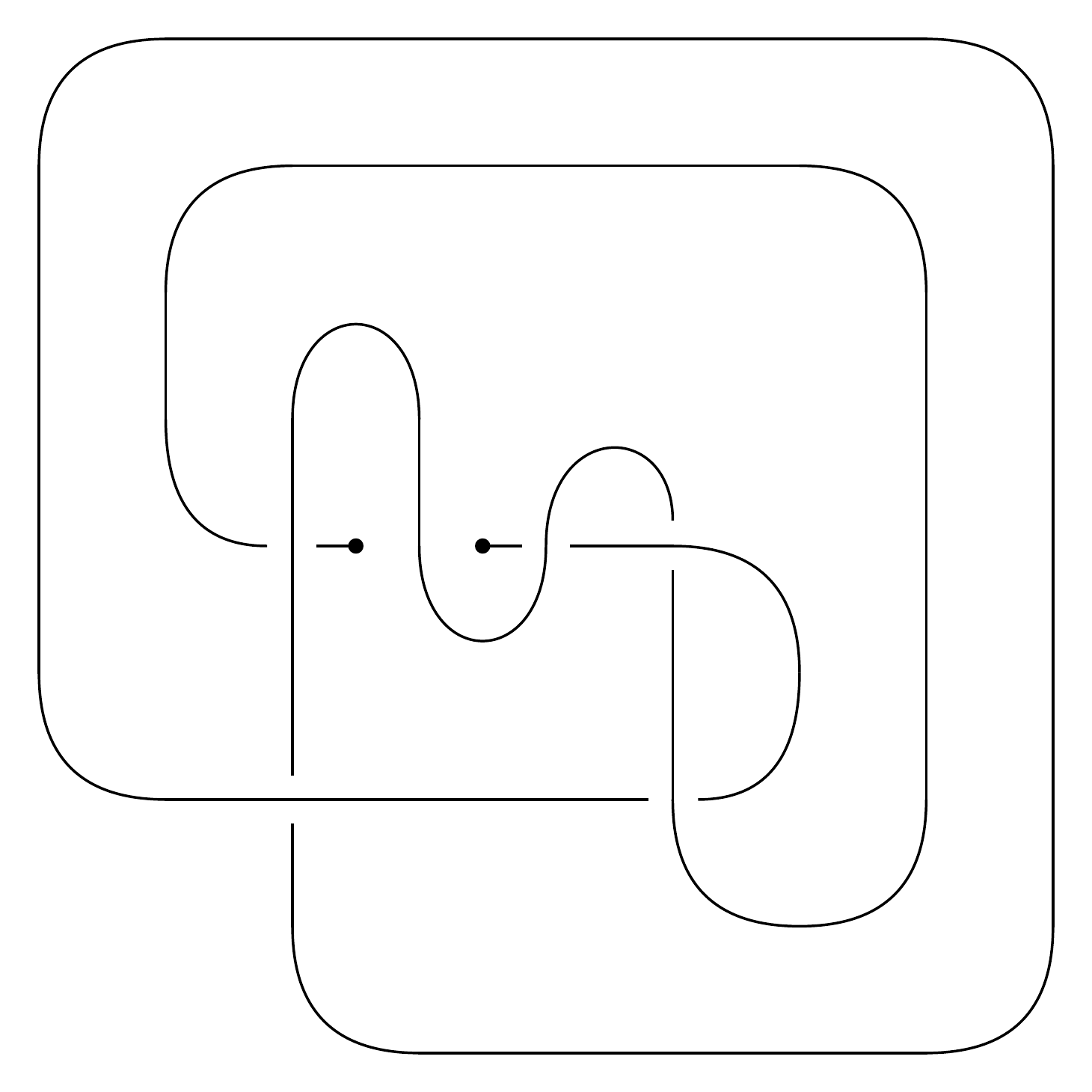}\\
\textcolor{black}{$5_{751}$}
\vspace{1cm}
\end{minipage}
\begin{minipage}[t]{.25\linewidth}
\centering
\includegraphics[width=0.9\textwidth,height=3.5cm,keepaspectratio]{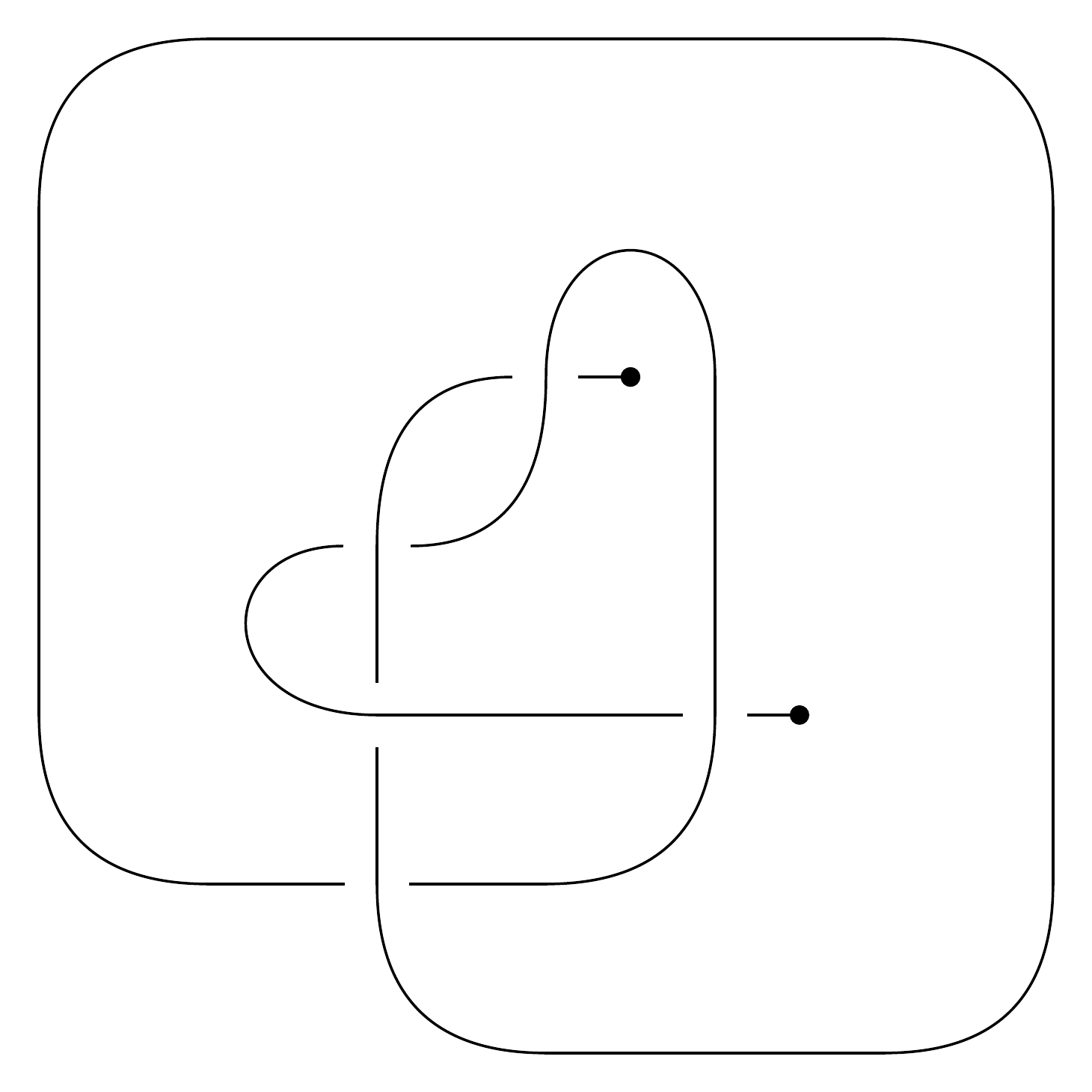}\\
\textcolor{black}{$5_{752}$}
\vspace{1cm}
\end{minipage}
\begin{minipage}[t]{.25\linewidth}
\centering
\includegraphics[width=0.9\textwidth,height=3.5cm,keepaspectratio]{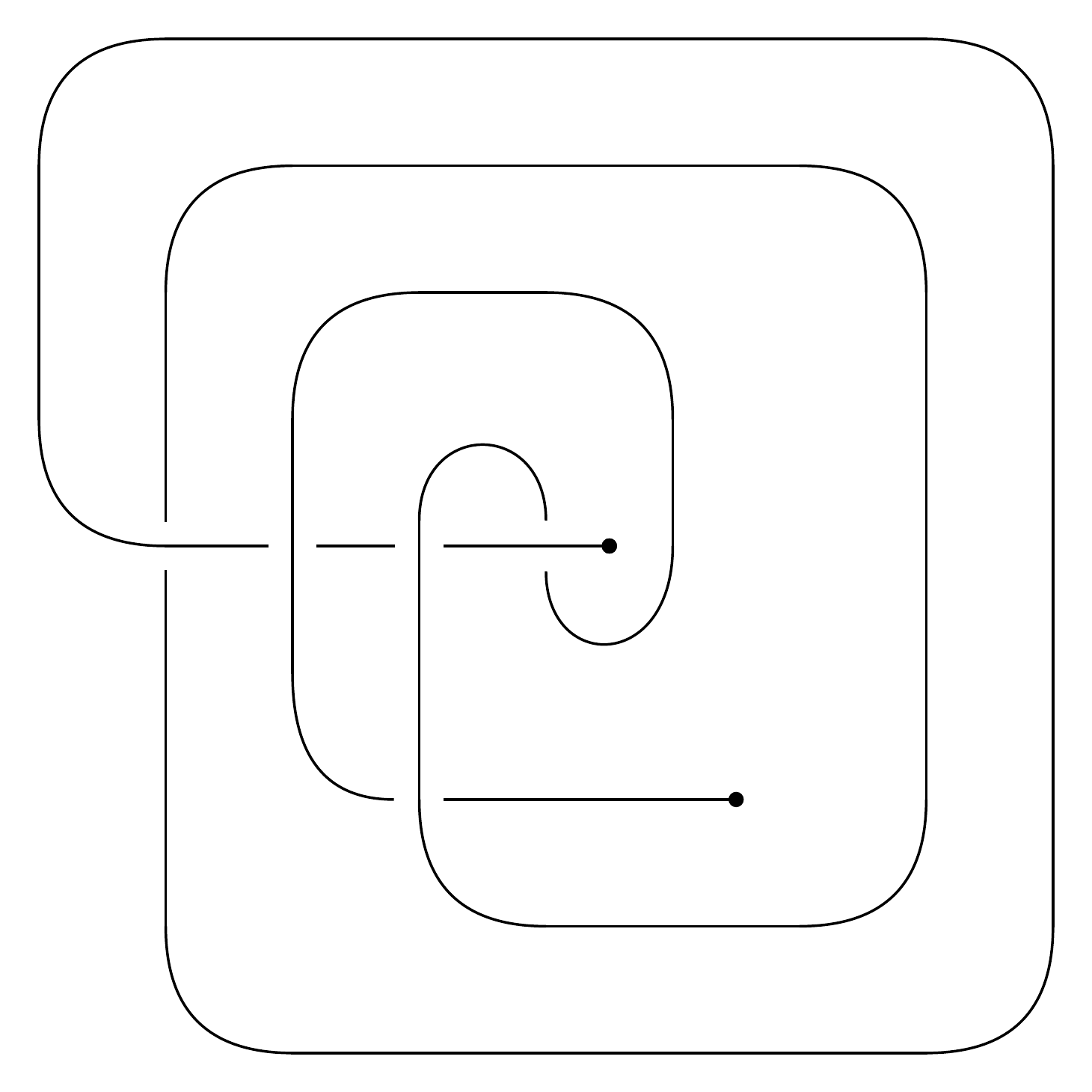}\\
\textcolor{black}{$5_{753}$}
\vspace{1cm}
\end{minipage}
\begin{minipage}[t]{.25\linewidth}
\centering
\includegraphics[width=0.9\textwidth,height=3.5cm,keepaspectratio]{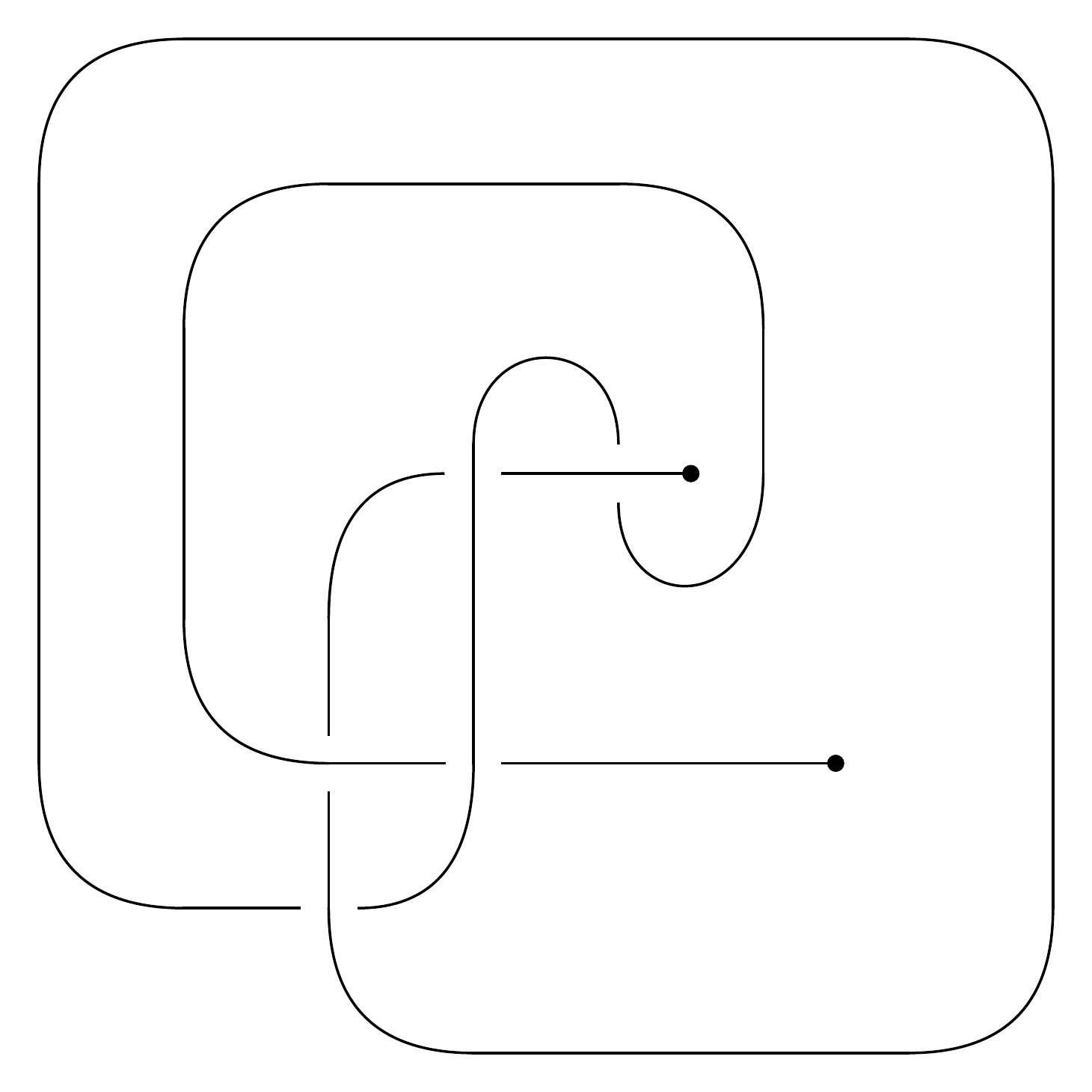}\\
\textcolor{black}{$5_{754}$}
\vspace{1cm}
\end{minipage}
\begin{minipage}[t]{.25\linewidth}
\centering
\includegraphics[width=0.9\textwidth,height=3.5cm,keepaspectratio]{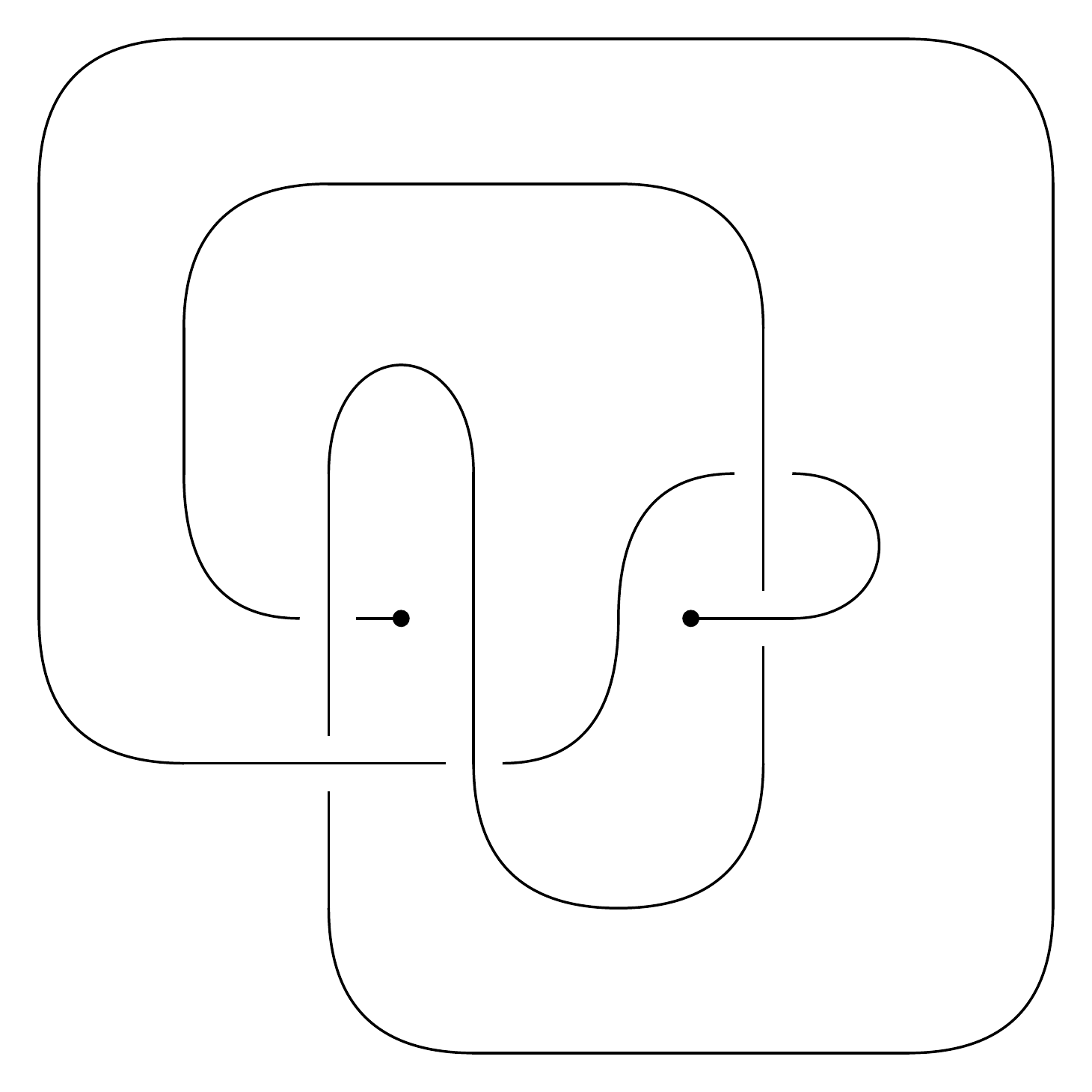}\\
\textcolor{black}{$5_{755}$}
\vspace{1cm}
\end{minipage}
\begin{minipage}[t]{.25\linewidth}
\centering
\includegraphics[width=0.9\textwidth,height=3.5cm,keepaspectratio]{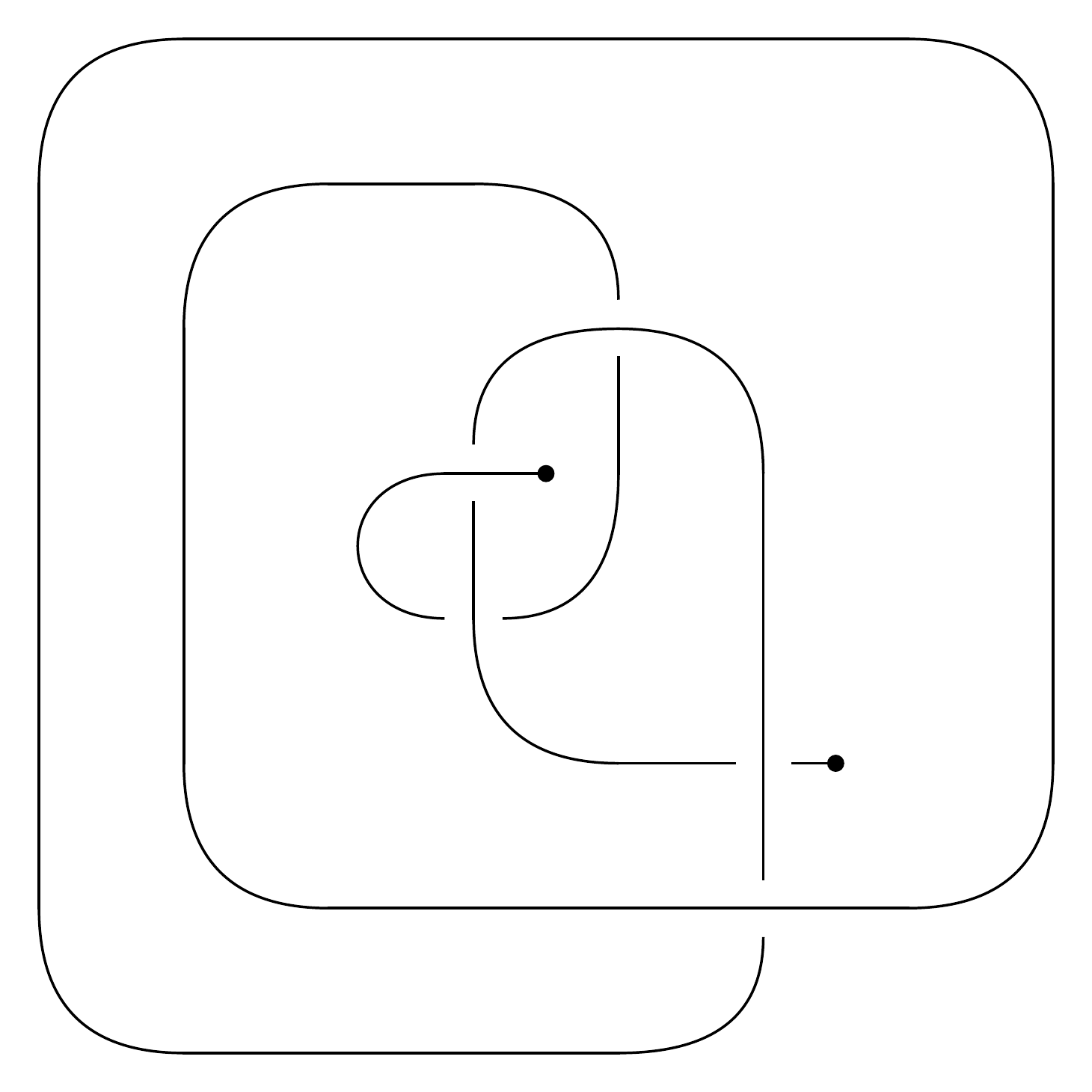}\\
\textcolor{black}{$5_{756}$}
\vspace{1cm}
\end{minipage}
\begin{minipage}[t]{.25\linewidth}
\centering
\includegraphics[width=0.9\textwidth,height=3.5cm,keepaspectratio]{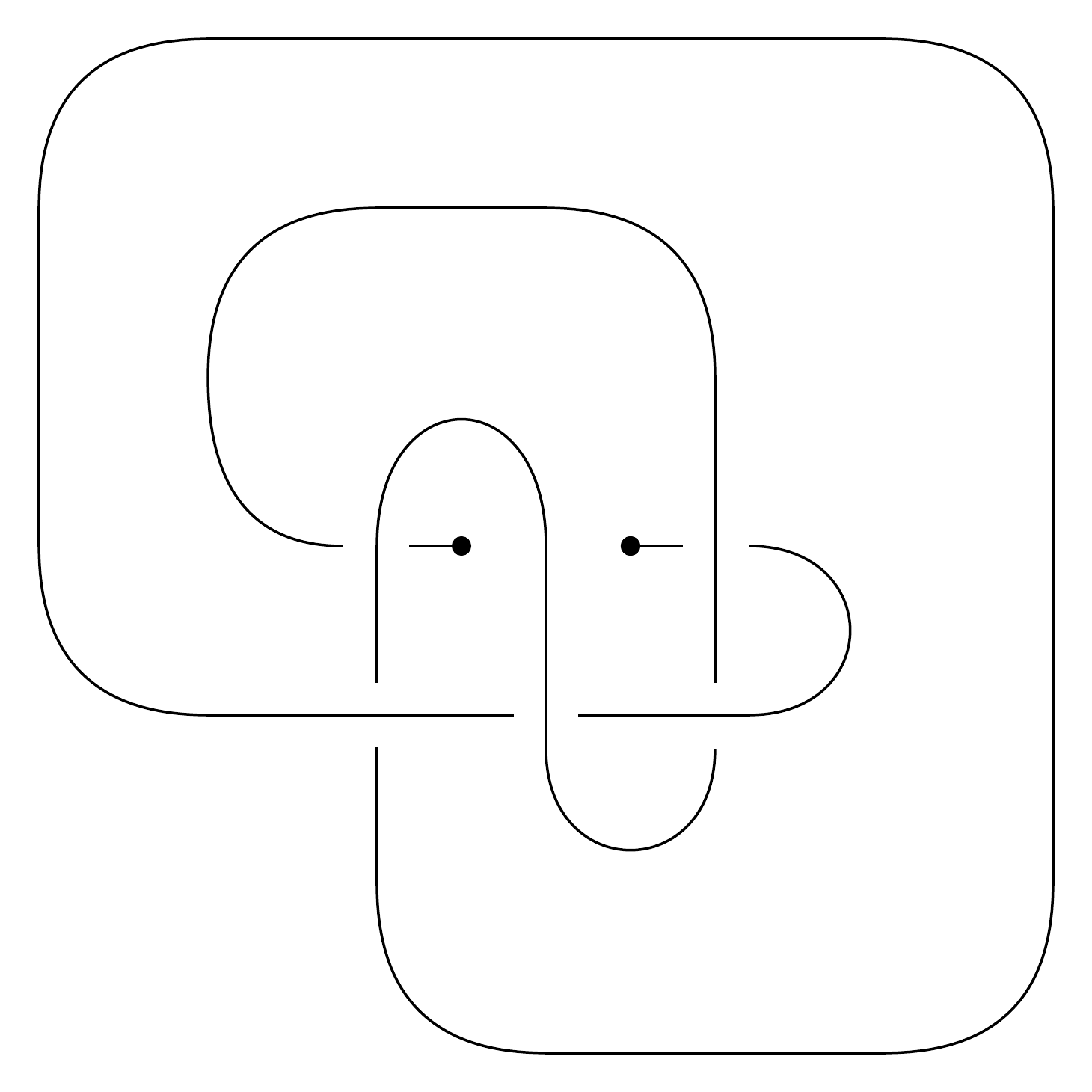}\\
\textcolor{black}{$5_{757}$}
\vspace{1cm}
\end{minipage}
\begin{minipage}[t]{.25\linewidth}
\centering
\includegraphics[width=0.9\textwidth,height=3.5cm,keepaspectratio]{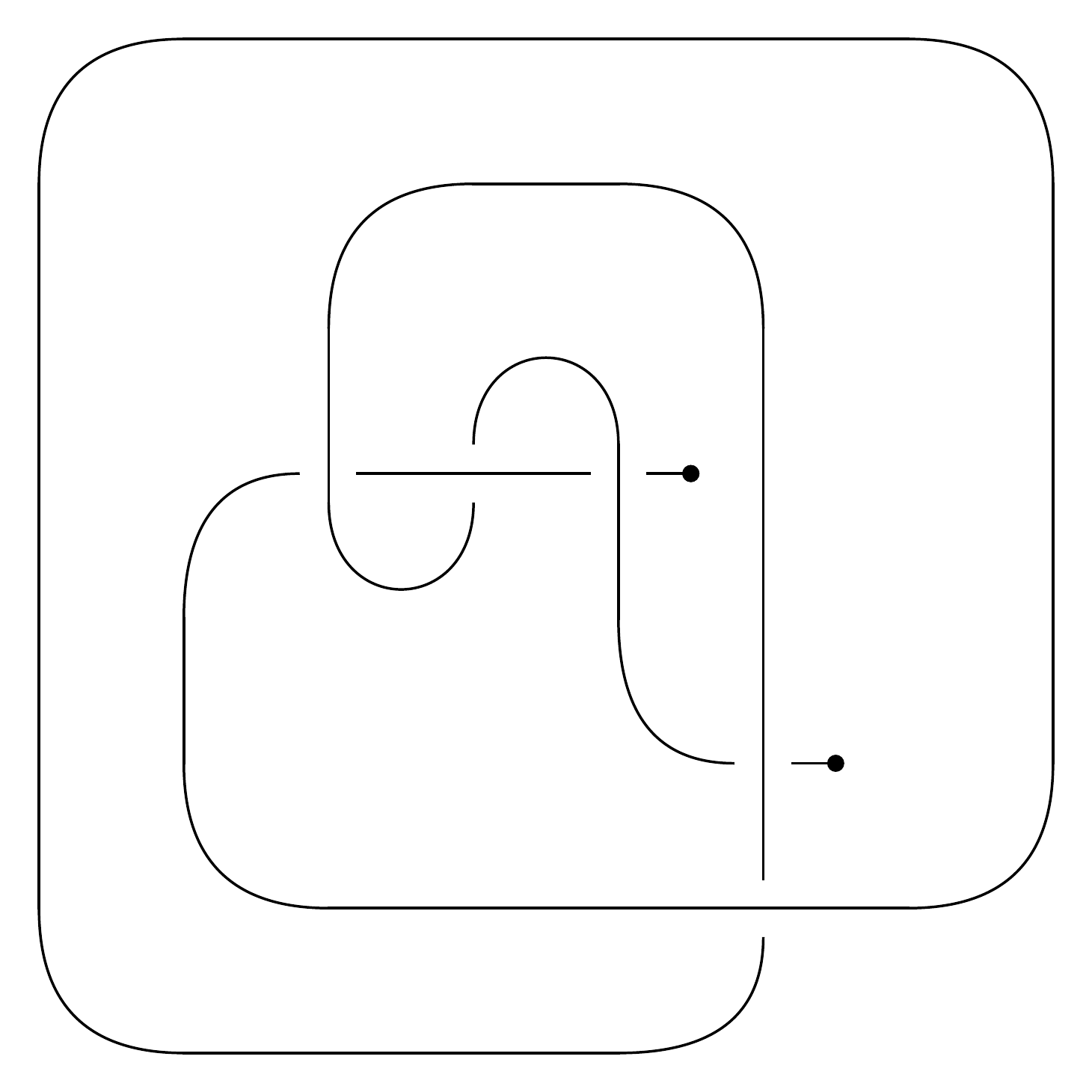}\\
\textcolor{black}{$5_{758}$}
\vspace{1cm}
\end{minipage}
\begin{minipage}[t]{.25\linewidth}
\centering
\includegraphics[width=0.9\textwidth,height=3.5cm,keepaspectratio]{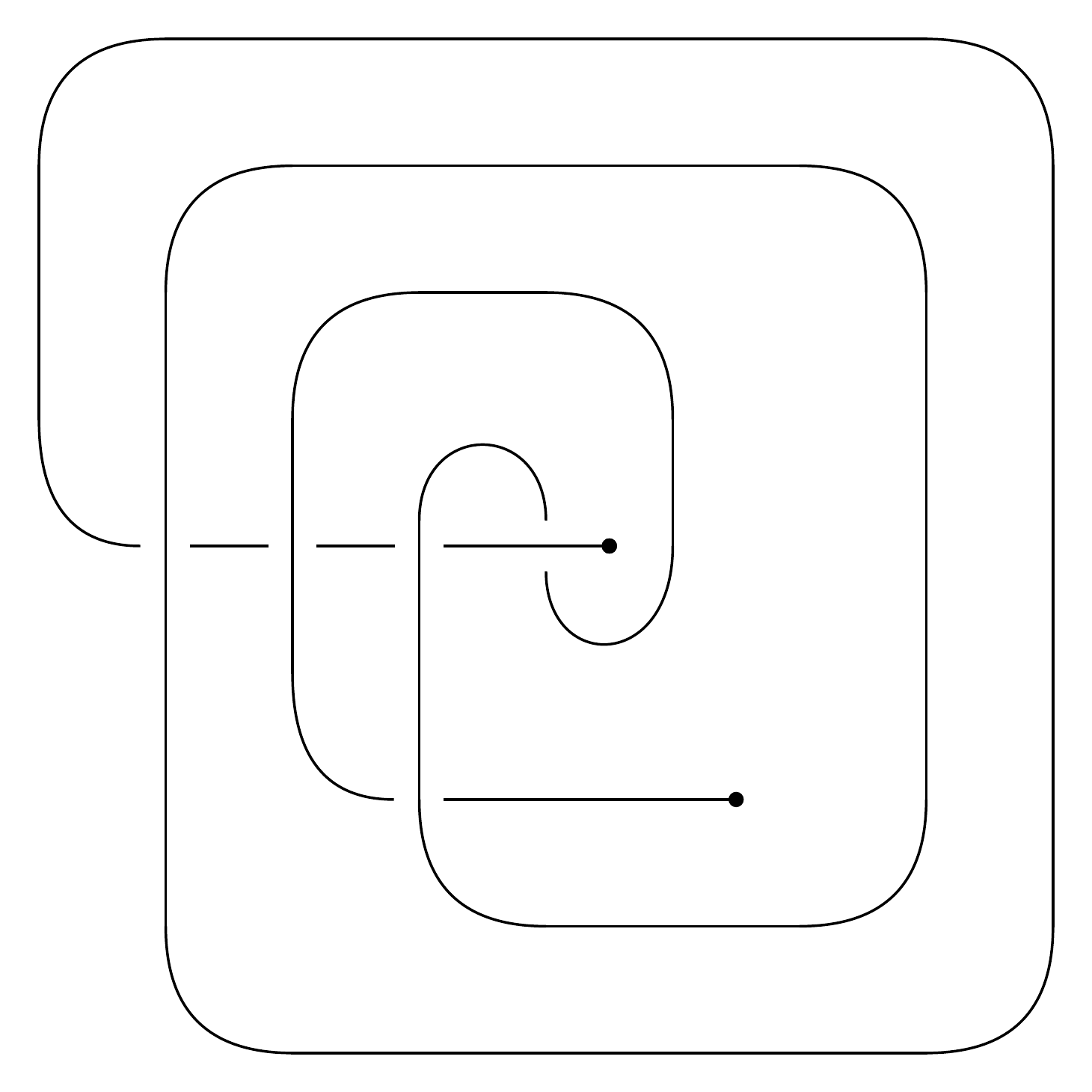}\\
\textcolor{black}{$5_{759}$}
\vspace{1cm}
\end{minipage}
\begin{minipage}[t]{.25\linewidth}
\centering
\includegraphics[width=0.9\textwidth,height=3.5cm,keepaspectratio]{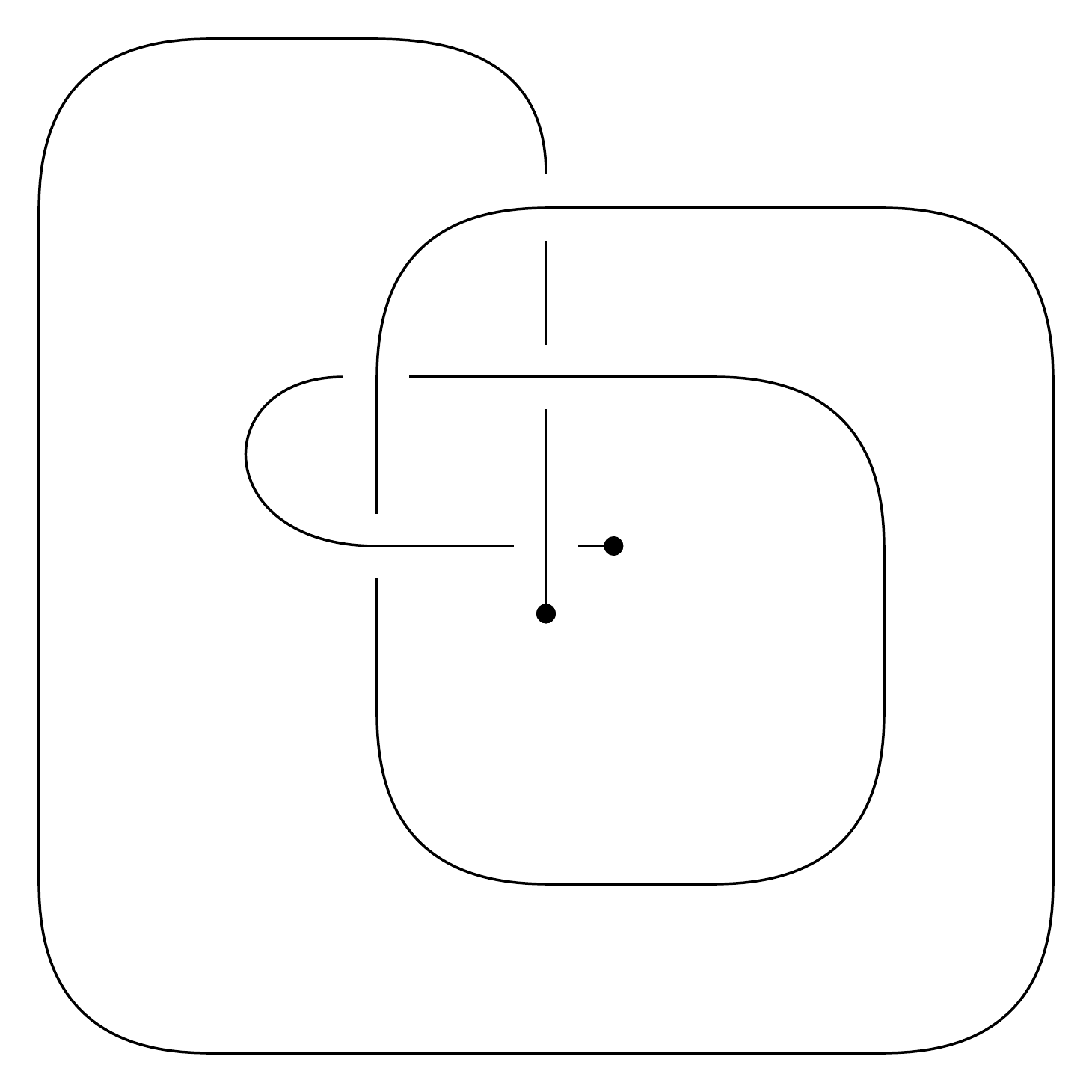}\\
\textcolor{black}{$5_{760}$}
\vspace{1cm}
\end{minipage}
\begin{minipage}[t]{.25\linewidth}
\centering
\includegraphics[width=0.9\textwidth,height=3.5cm,keepaspectratio]{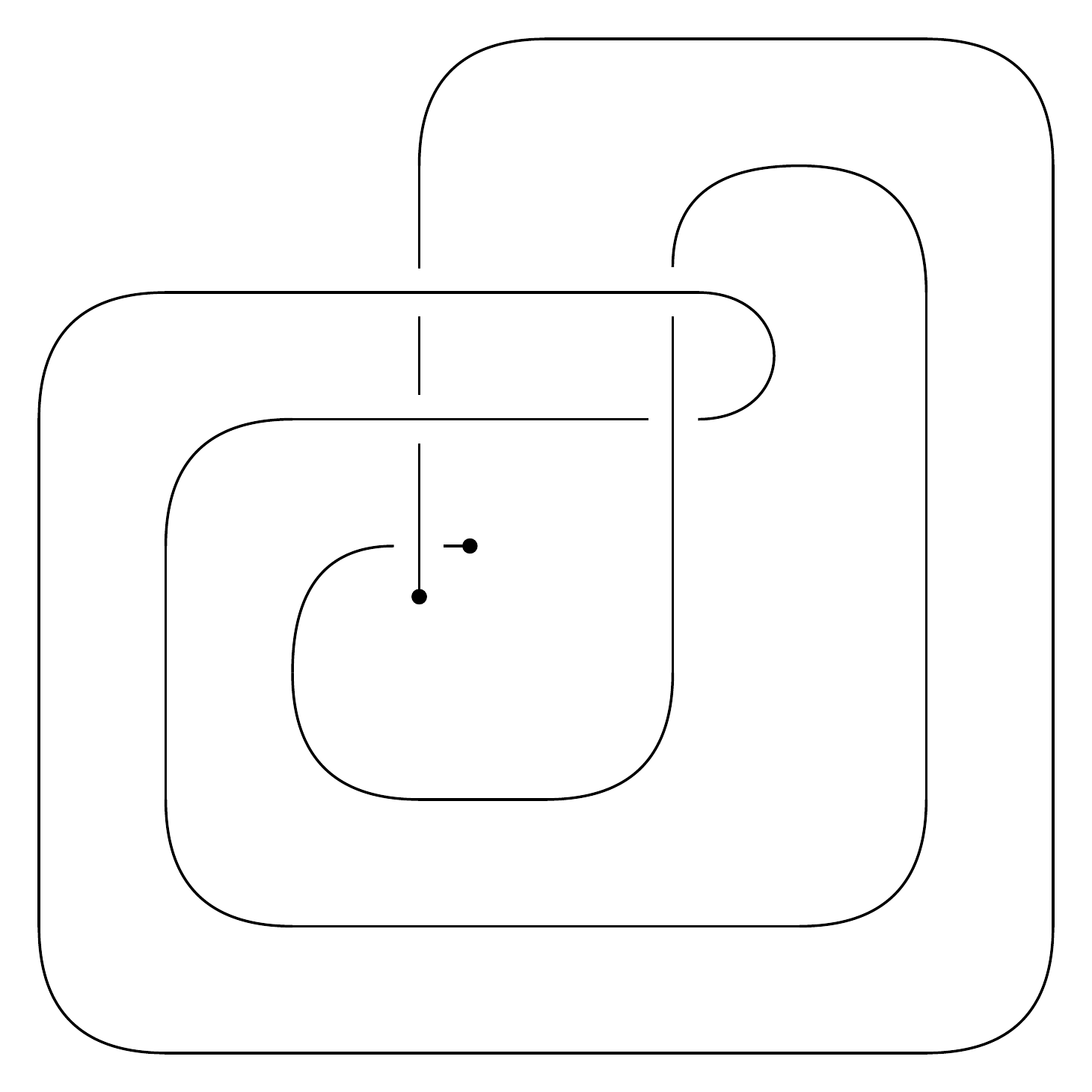}\\
\textcolor{black}{$5_{761}$}
\vspace{1cm}
\end{minipage}
\begin{minipage}[t]{.25\linewidth}
\centering
\includegraphics[width=0.9\textwidth,height=3.5cm,keepaspectratio]{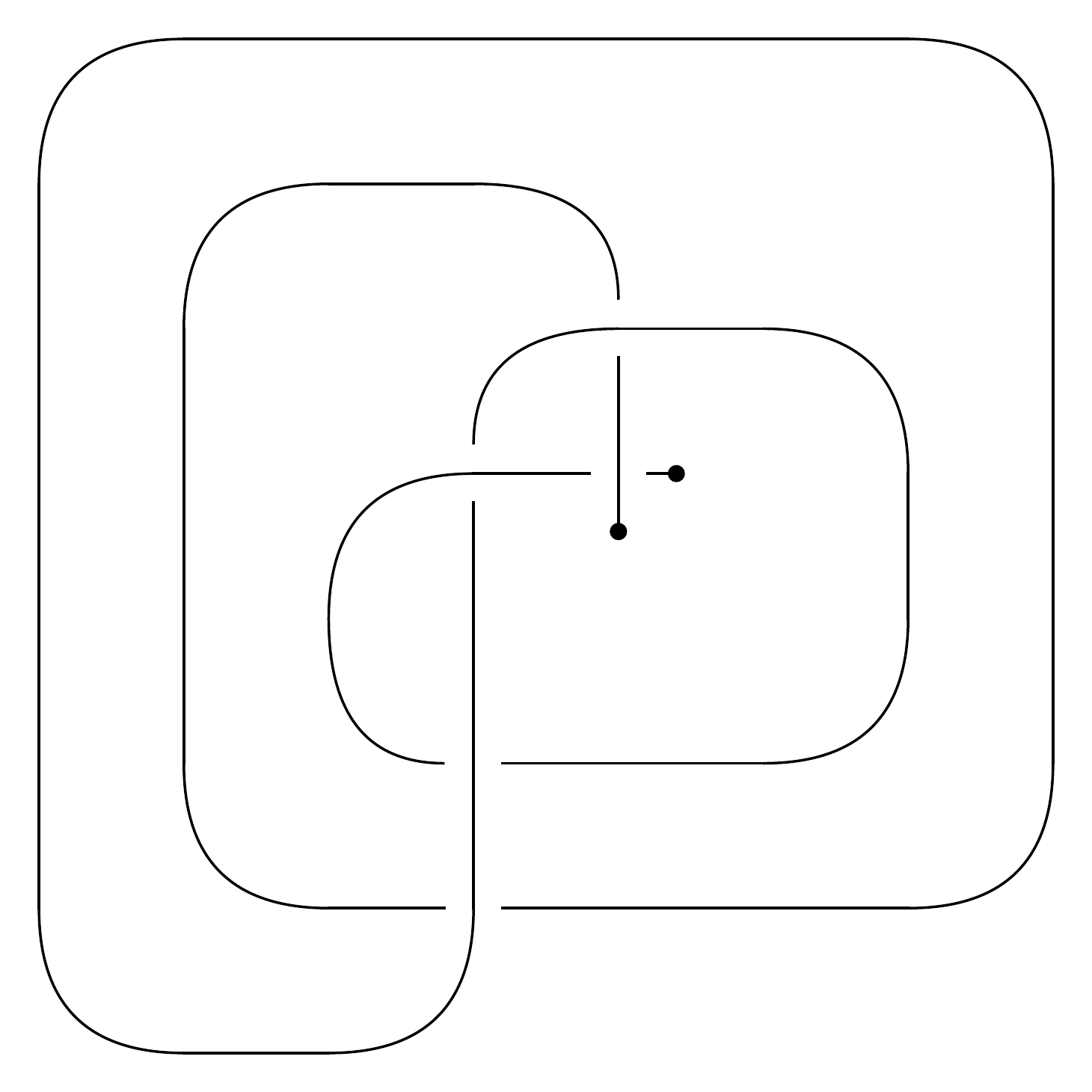}\\
\textcolor{black}{$5_{762}$}
\vspace{1cm}
\end{minipage}
\begin{minipage}[t]{.25\linewidth}
\centering
\includegraphics[width=0.9\textwidth,height=3.5cm,keepaspectratio]{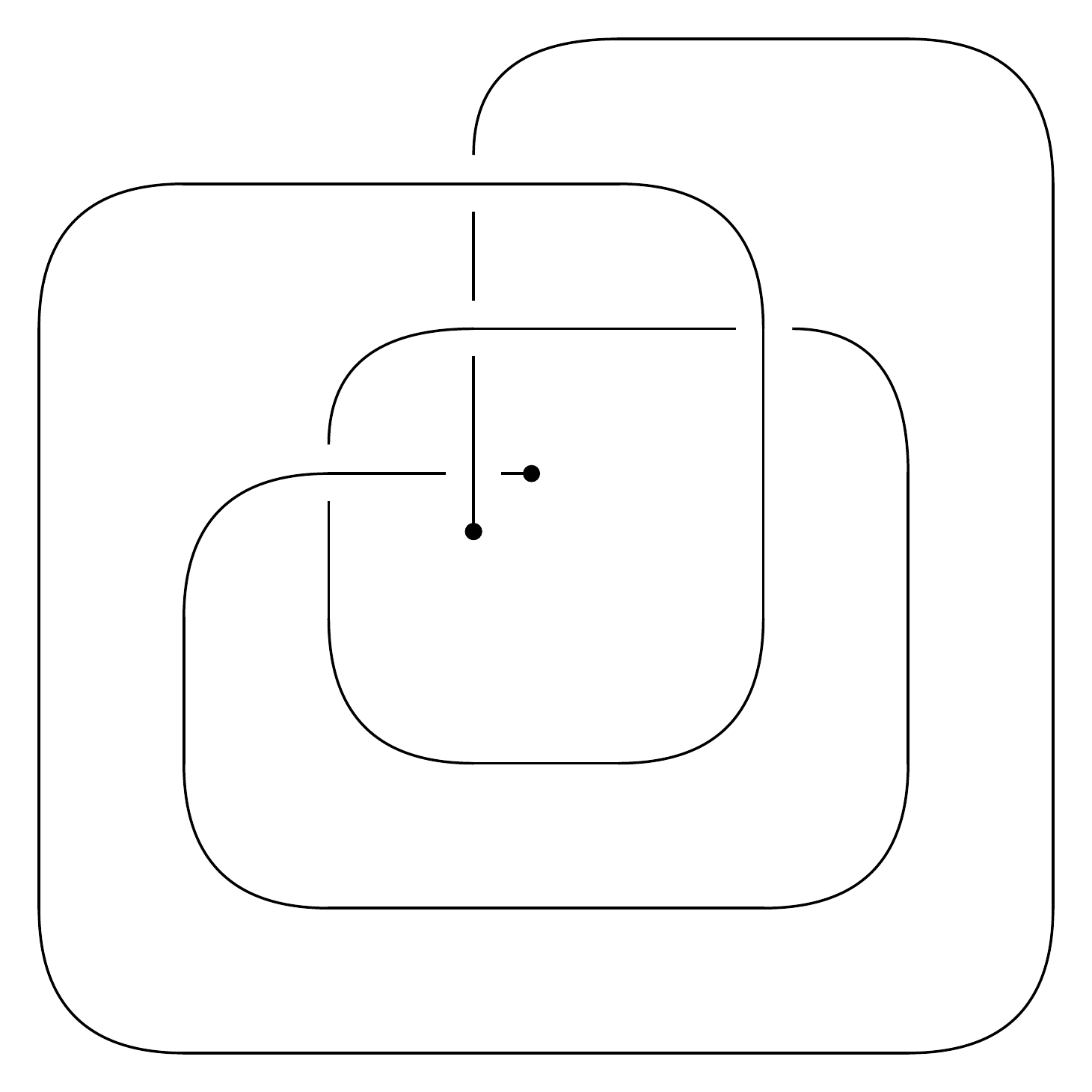}\\
\textcolor{black}{$5_{763}$}
\vspace{1cm}
\end{minipage}
\begin{minipage}[t]{.25\linewidth}
\centering
\includegraphics[width=0.9\textwidth,height=3.5cm,keepaspectratio]{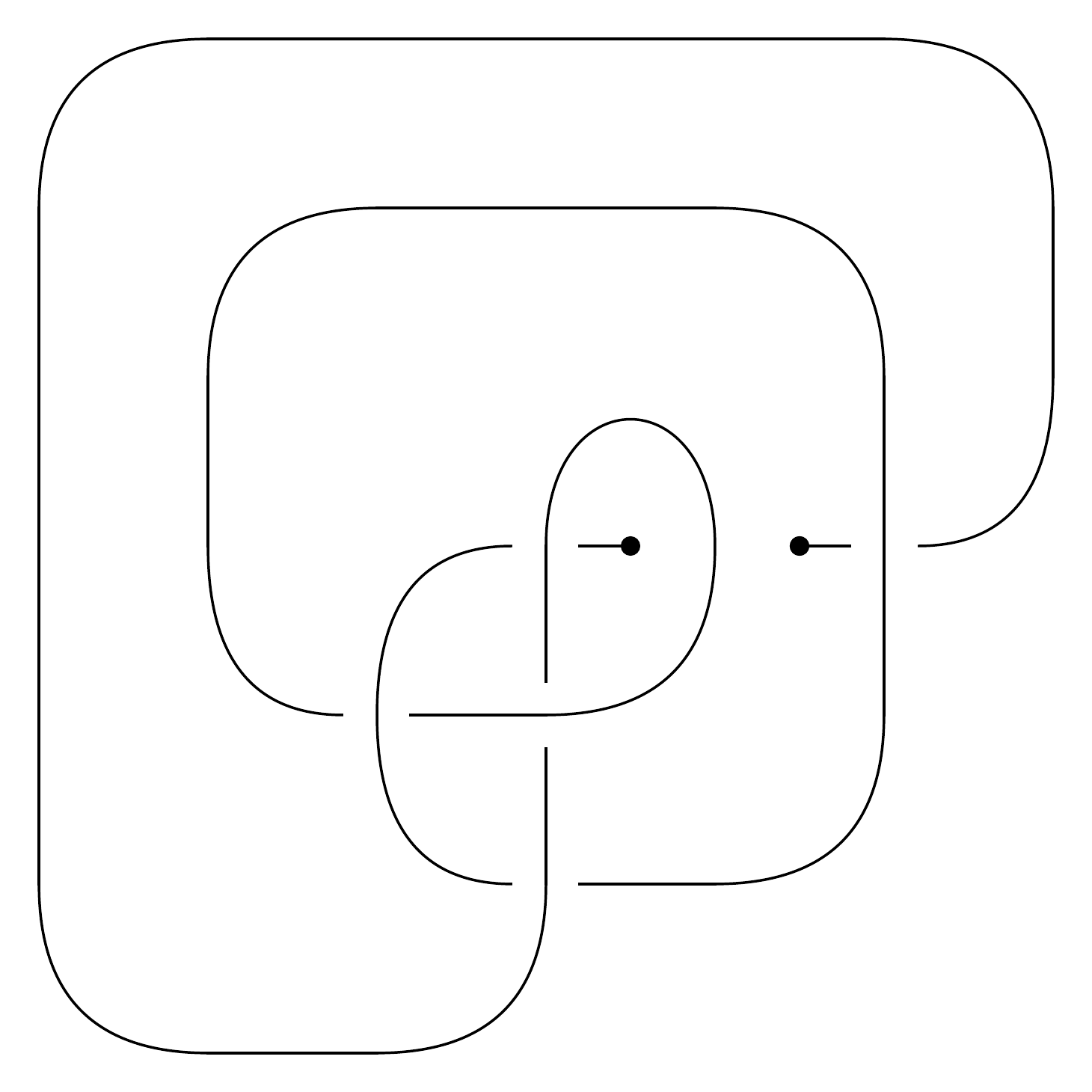}\\
\textcolor{black}{$5_{764}$}
\vspace{1cm}
\end{minipage}
\begin{minipage}[t]{.25\linewidth}
\centering
\includegraphics[width=0.9\textwidth,height=3.5cm,keepaspectratio]{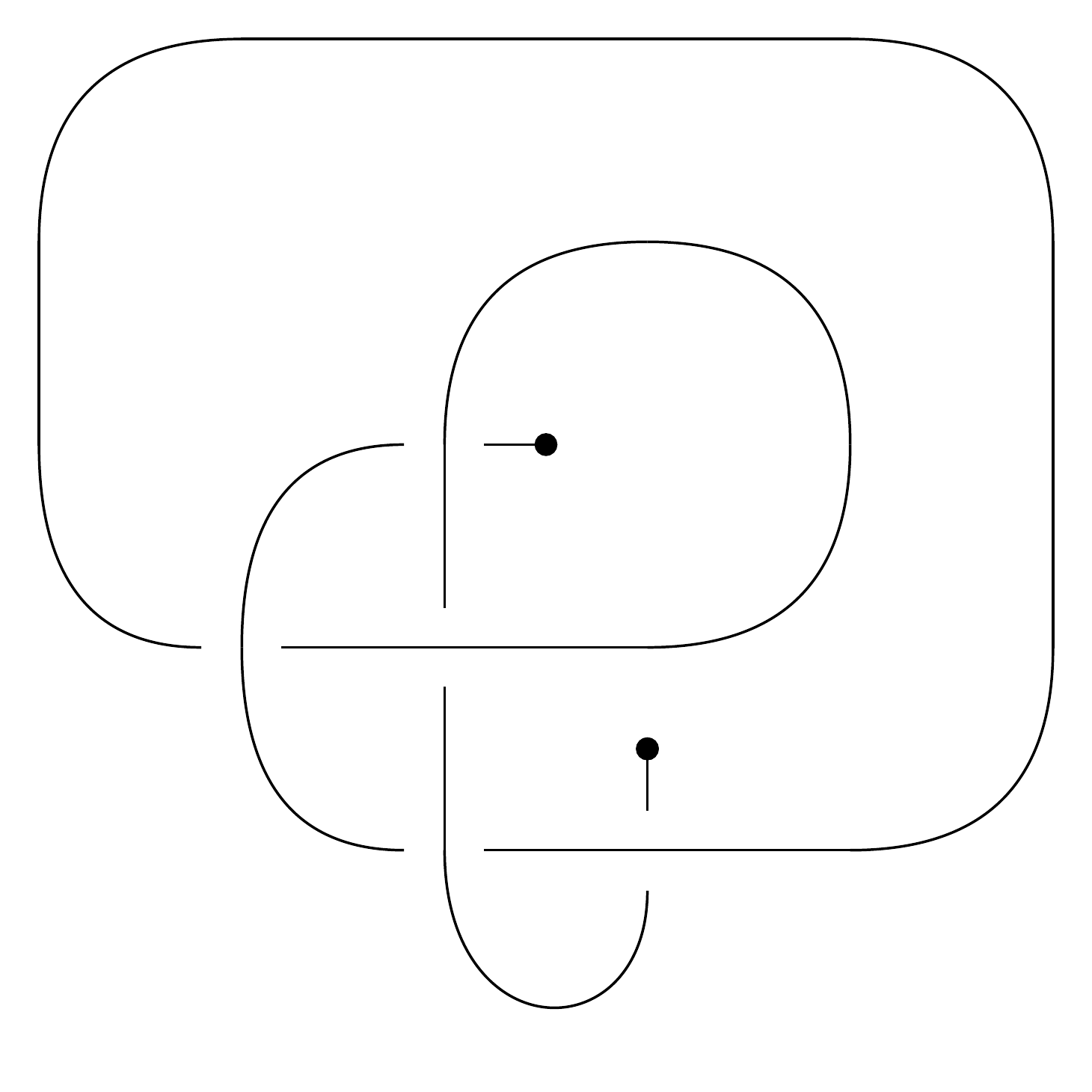}\\
\textcolor{black}{$5_{765}$}
\vspace{1cm}
\end{minipage}
\begin{minipage}[t]{.25\linewidth}
\centering
\includegraphics[width=0.9\textwidth,height=3.5cm,keepaspectratio]{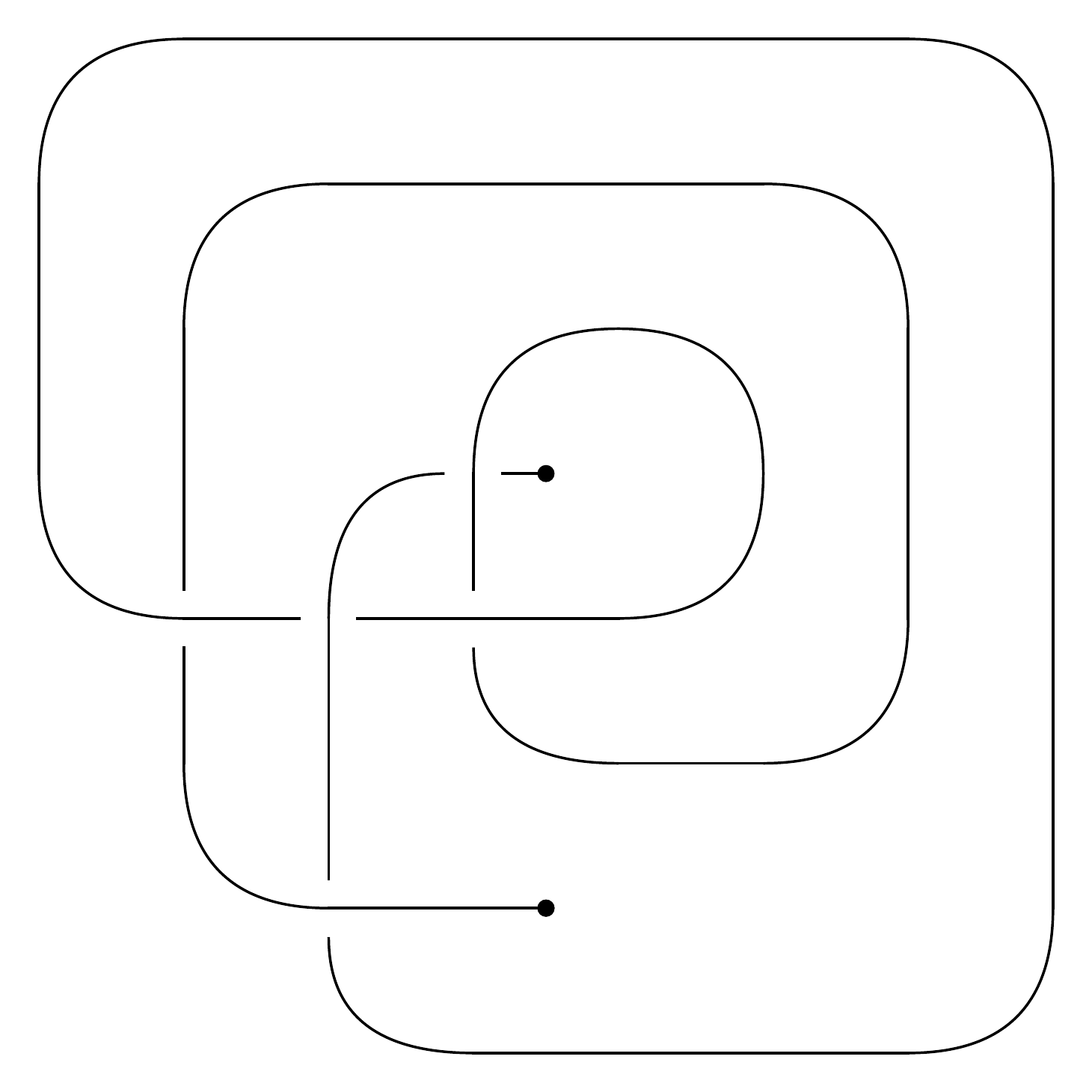}\\
\textcolor{black}{$5_{766}$}
\vspace{1cm}
\end{minipage}
\begin{minipage}[t]{.25\linewidth}
\centering
\includegraphics[width=0.9\textwidth,height=3.5cm,keepaspectratio]{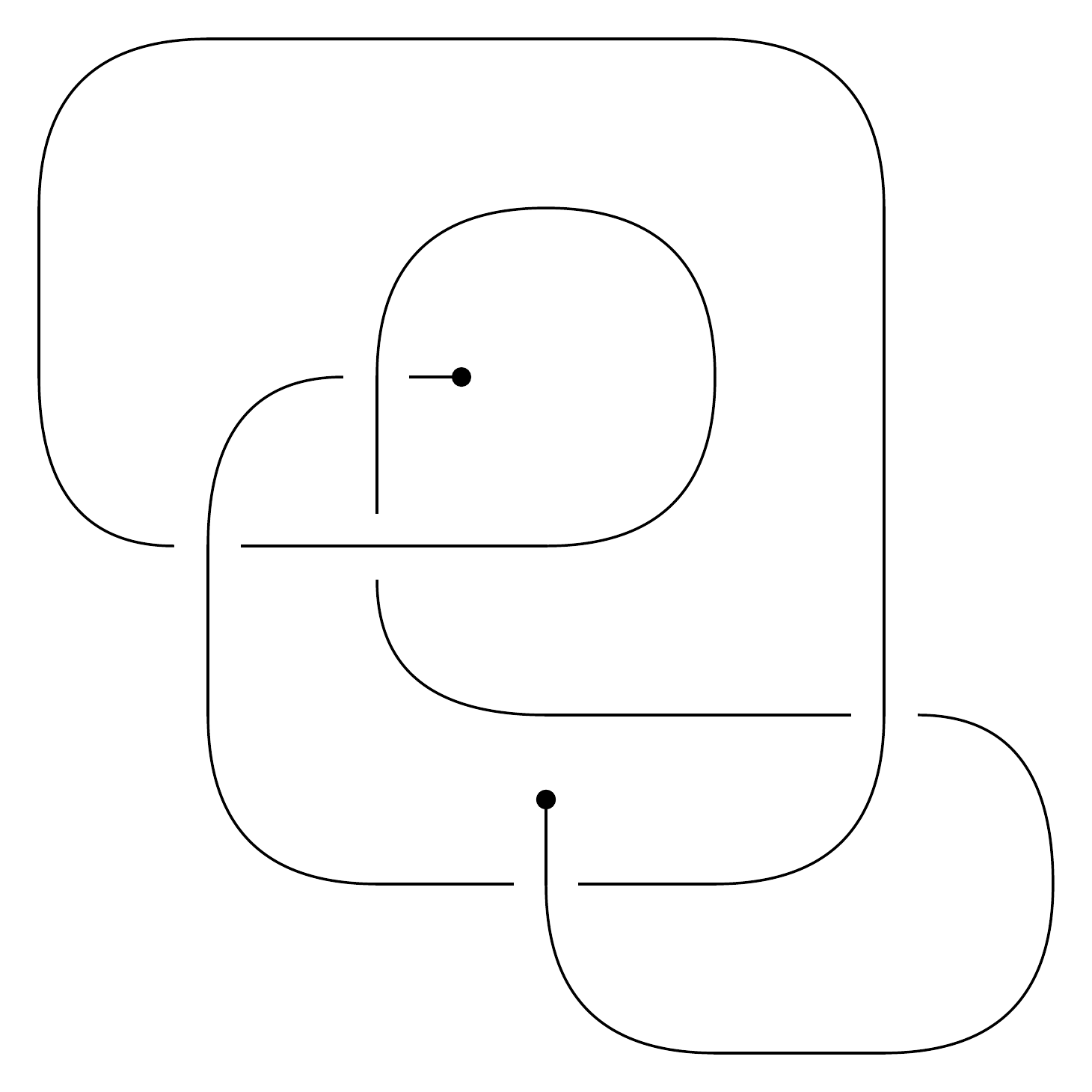}\\
\textcolor{black}{$5_{767}$}
\vspace{1cm}
\end{minipage}
\begin{minipage}[t]{.25\linewidth}
\centering
\includegraphics[width=0.9\textwidth,height=3.5cm,keepaspectratio]{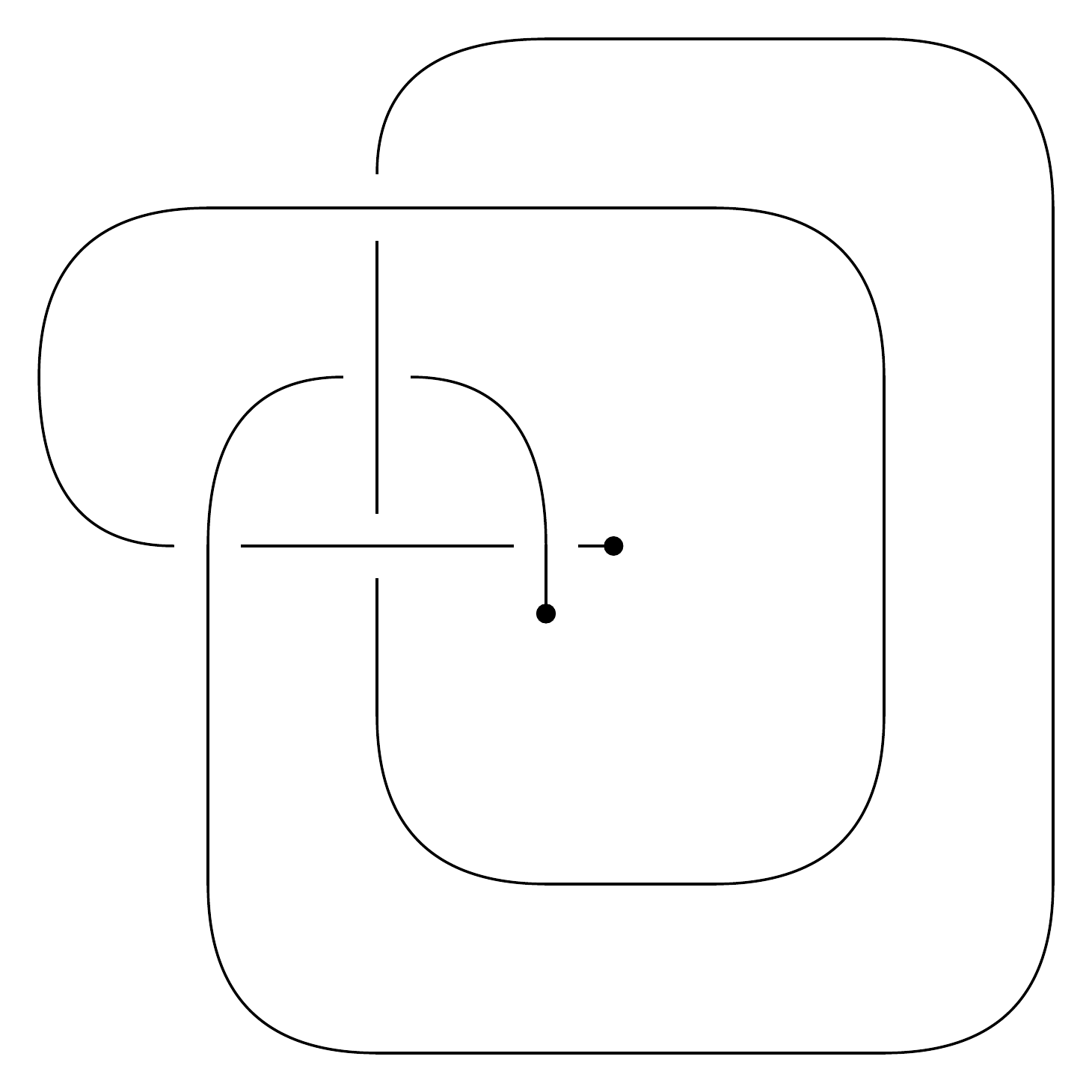}\\
\textcolor{black}{$5_{768}$}
\vspace{1cm}
\end{minipage}
\begin{minipage}[t]{.25\linewidth}
\centering
\includegraphics[width=0.9\textwidth,height=3.5cm,keepaspectratio]{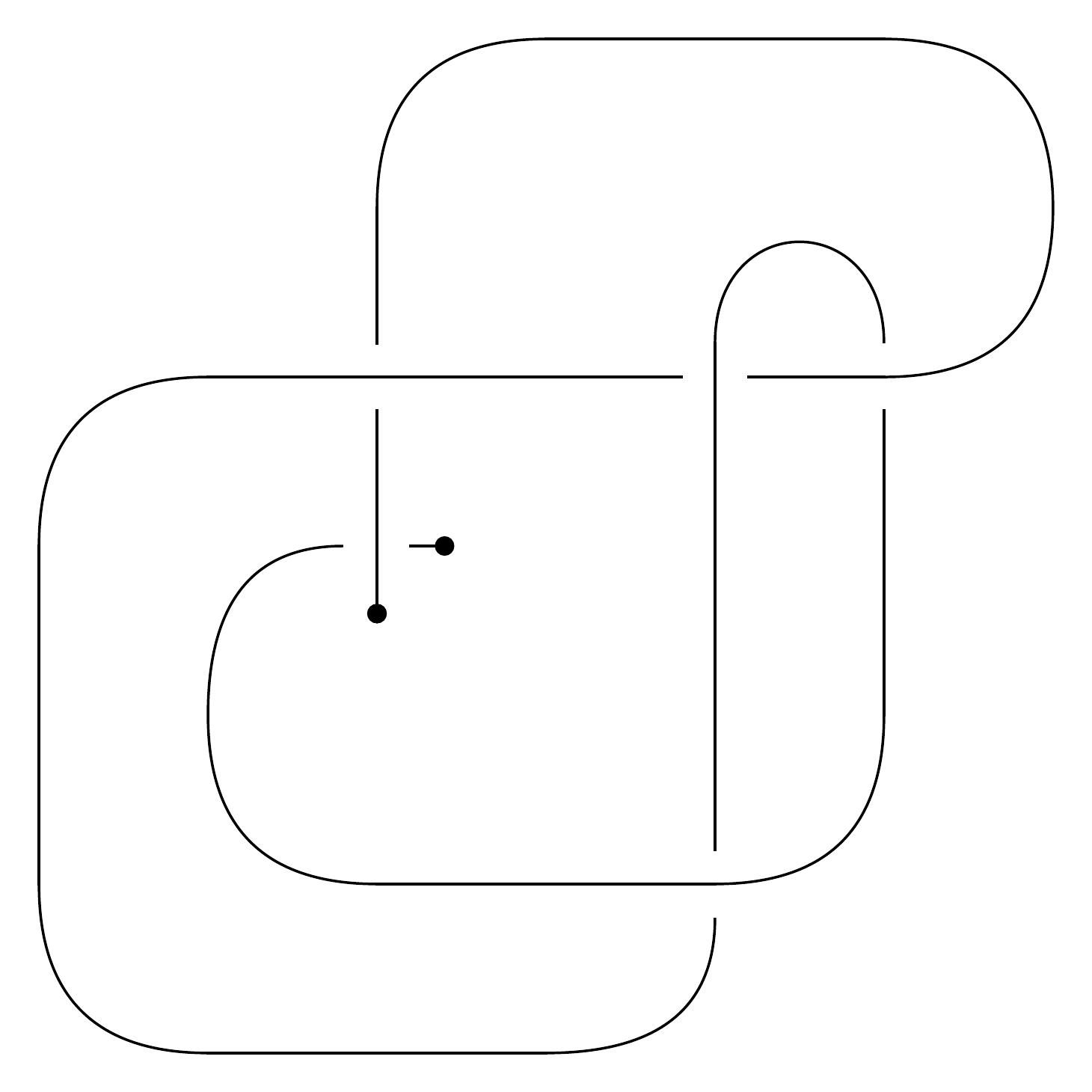}\\
\textcolor{black}{$5_{769}$}
\vspace{1cm}
\end{minipage}
\begin{minipage}[t]{.25\linewidth}
\centering
\includegraphics[width=0.9\textwidth,height=3.5cm,keepaspectratio]{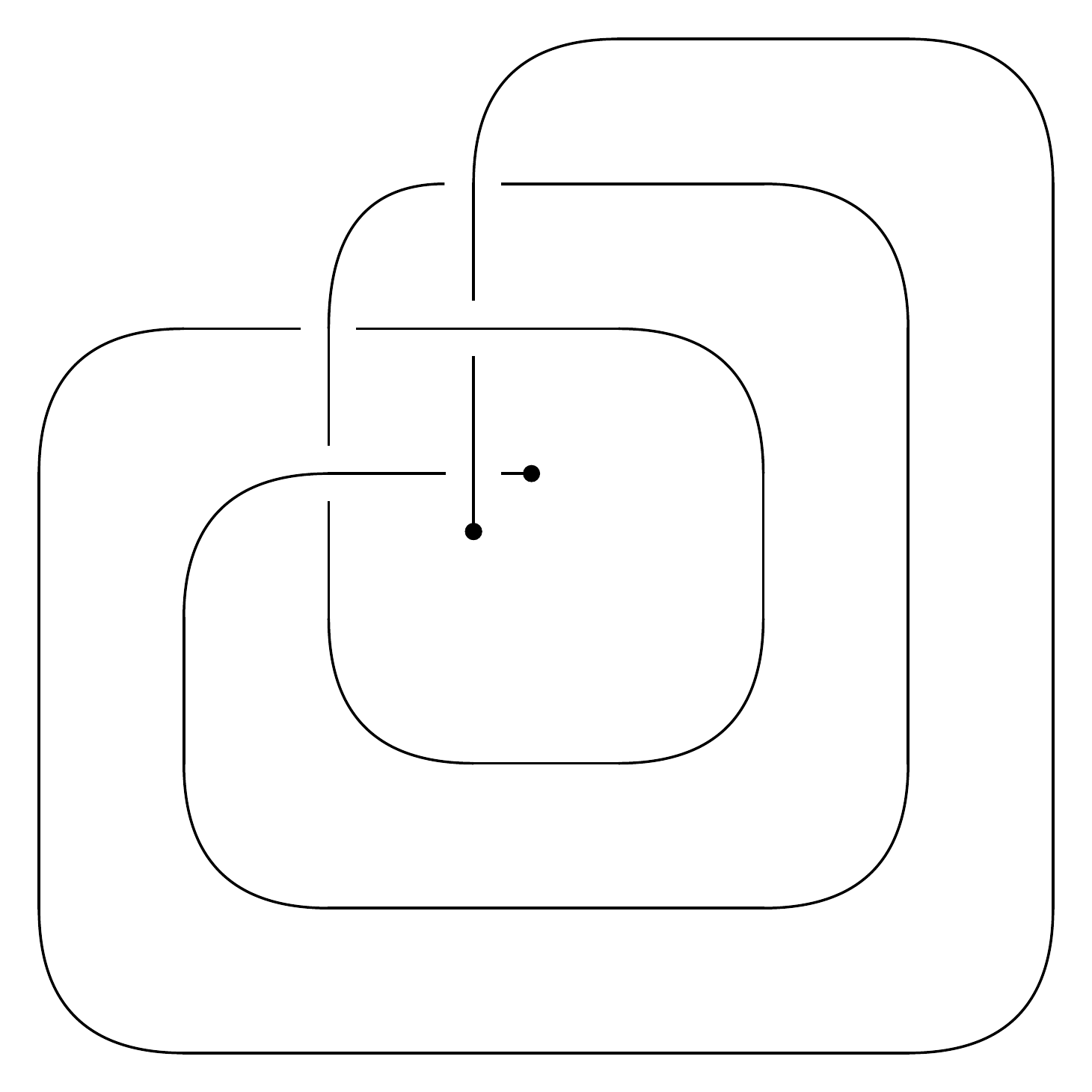}\\
\textcolor{black}{$5_{770}$}
\vspace{1cm}
\end{minipage}
\begin{minipage}[t]{.25\linewidth}
\centering
\includegraphics[width=0.9\textwidth,height=3.5cm,keepaspectratio]{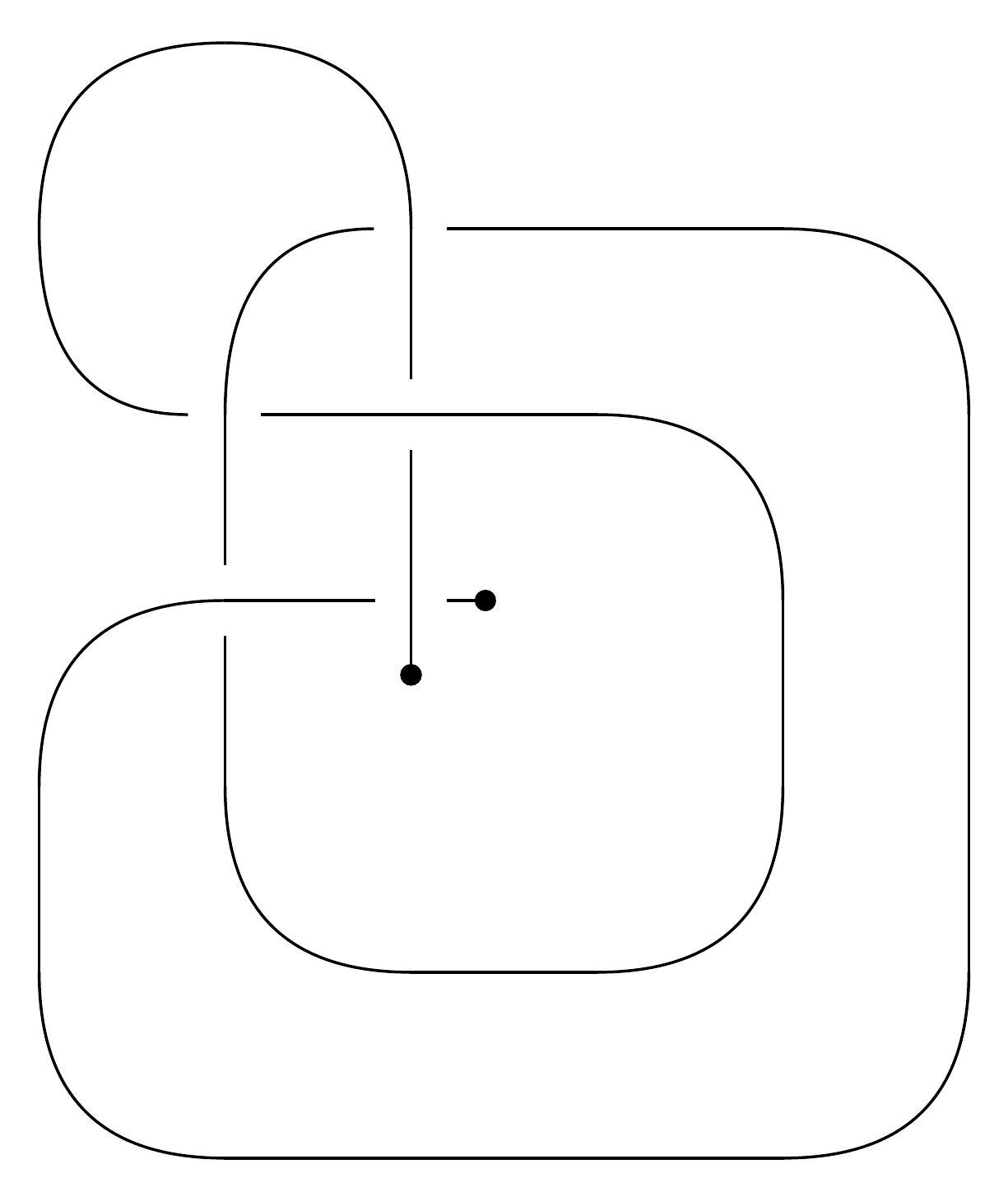}\\
\textcolor{black}{$5_{771}$}
\vspace{1cm}
\end{minipage}
\begin{minipage}[t]{.25\linewidth}
\centering
\includegraphics[width=0.9\textwidth,height=3.5cm,keepaspectratio]{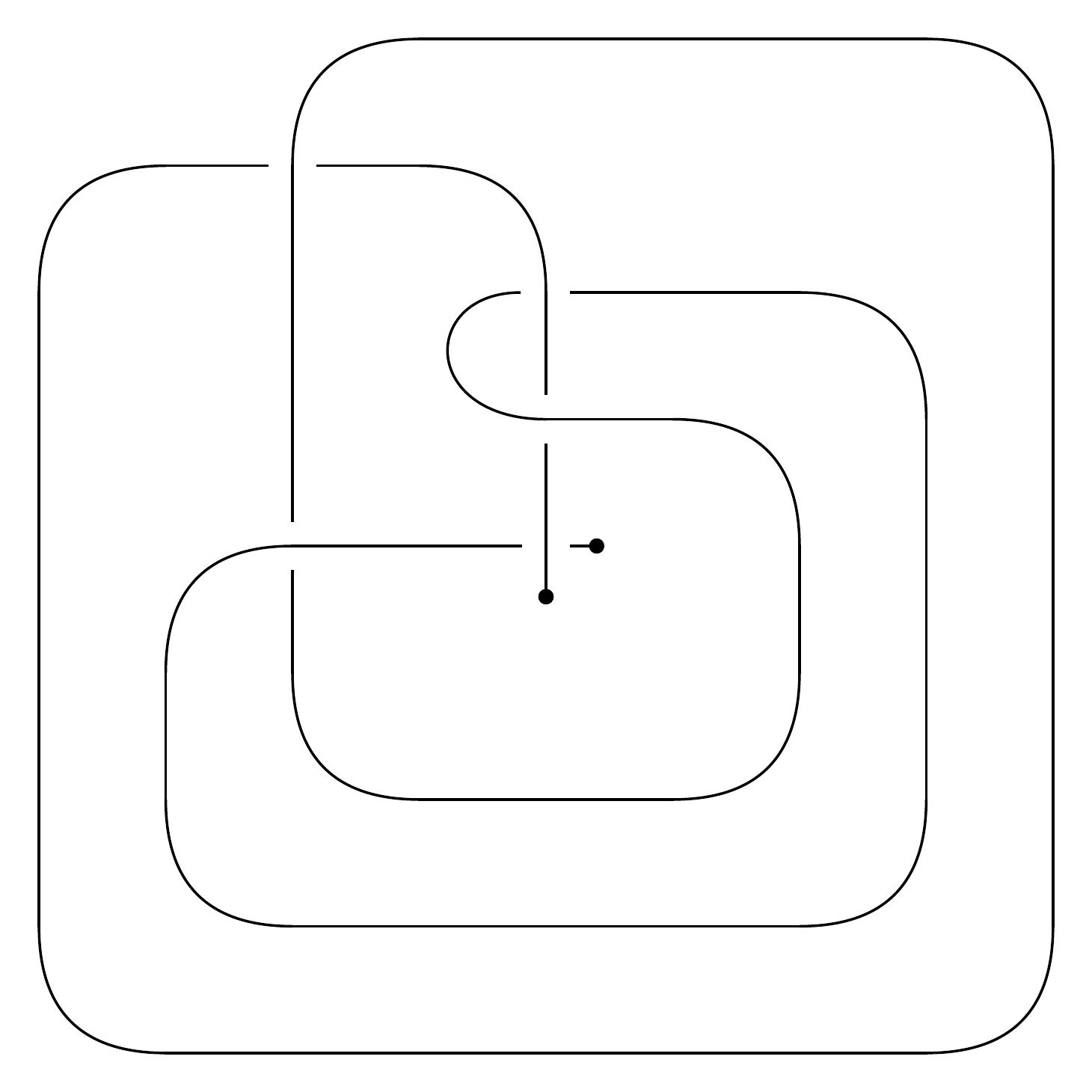}\\
\textcolor{black}{$5_{772}$}
\vspace{1cm}
\end{minipage}
\begin{minipage}[t]{.25\linewidth}
\centering
\includegraphics[width=0.9\textwidth,height=3.5cm,keepaspectratio]{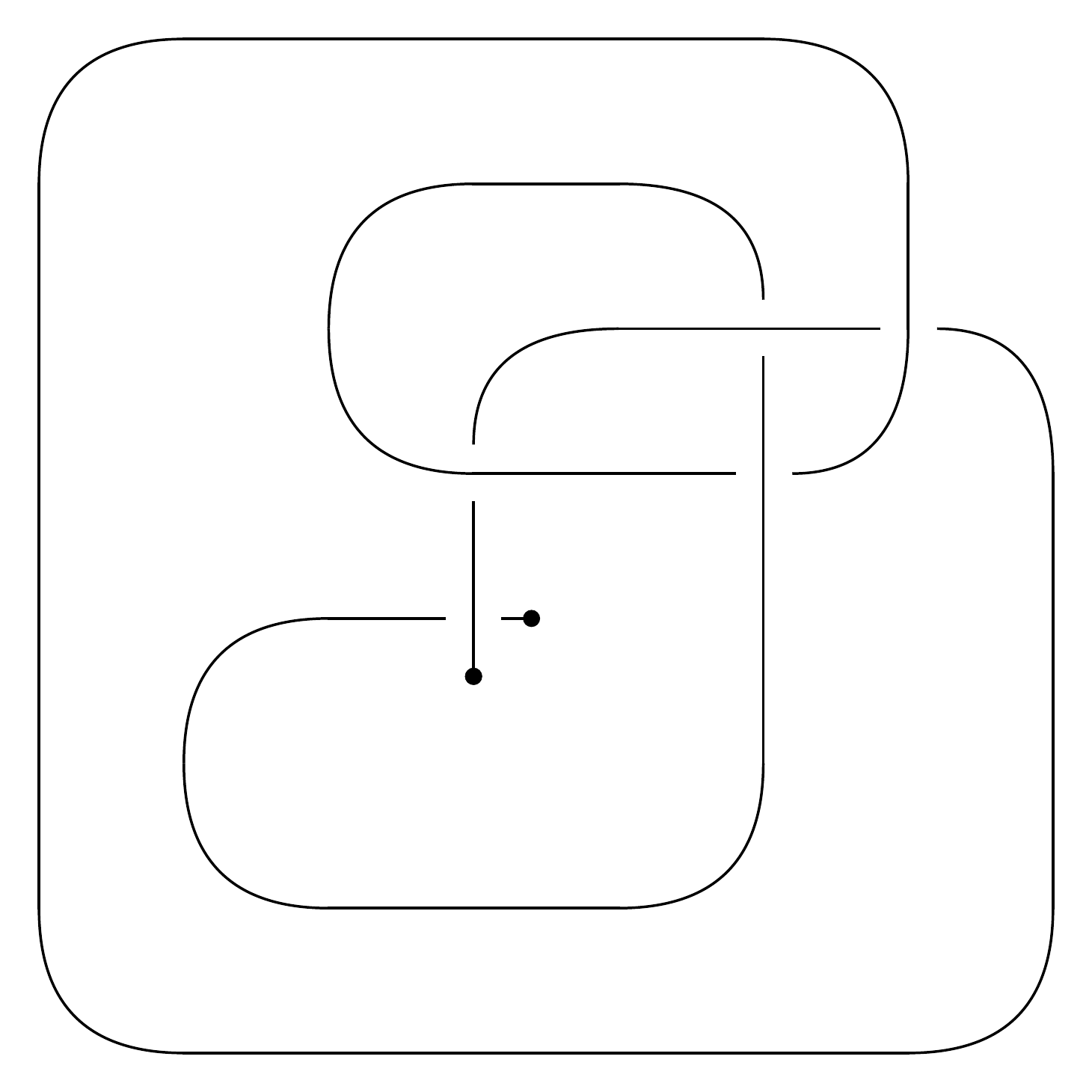}\\
\textcolor{black}{$5_{773}$}
\vspace{1cm}
\end{minipage}
\begin{minipage}[t]{.25\linewidth}
\centering
\includegraphics[width=0.9\textwidth,height=3.5cm,keepaspectratio]{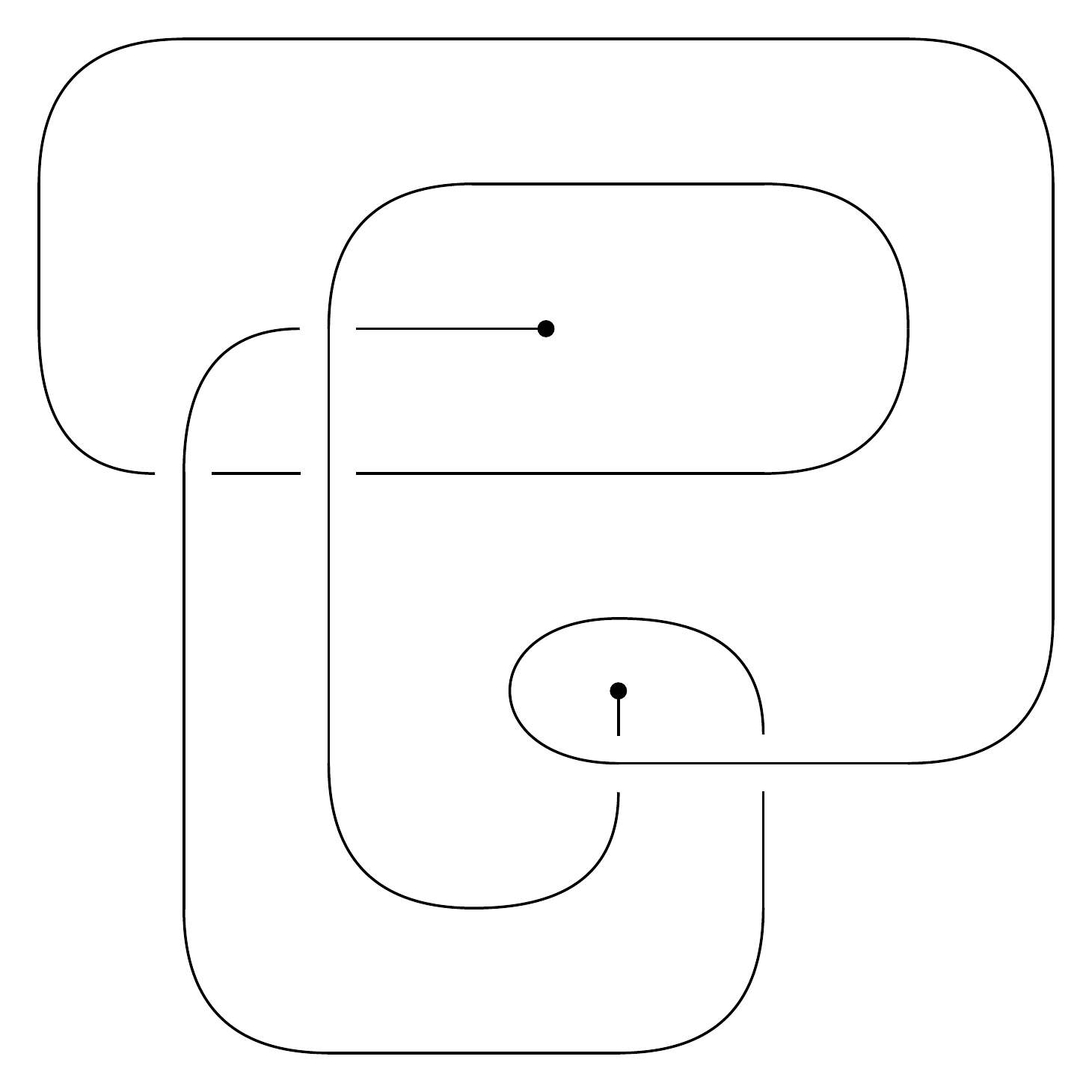}\\
\textcolor{black}{$5_{774}$}
\vspace{1cm}
\end{minipage}
\begin{minipage}[t]{.25\linewidth}
\centering
\includegraphics[width=0.9\textwidth,height=3.5cm,keepaspectratio]{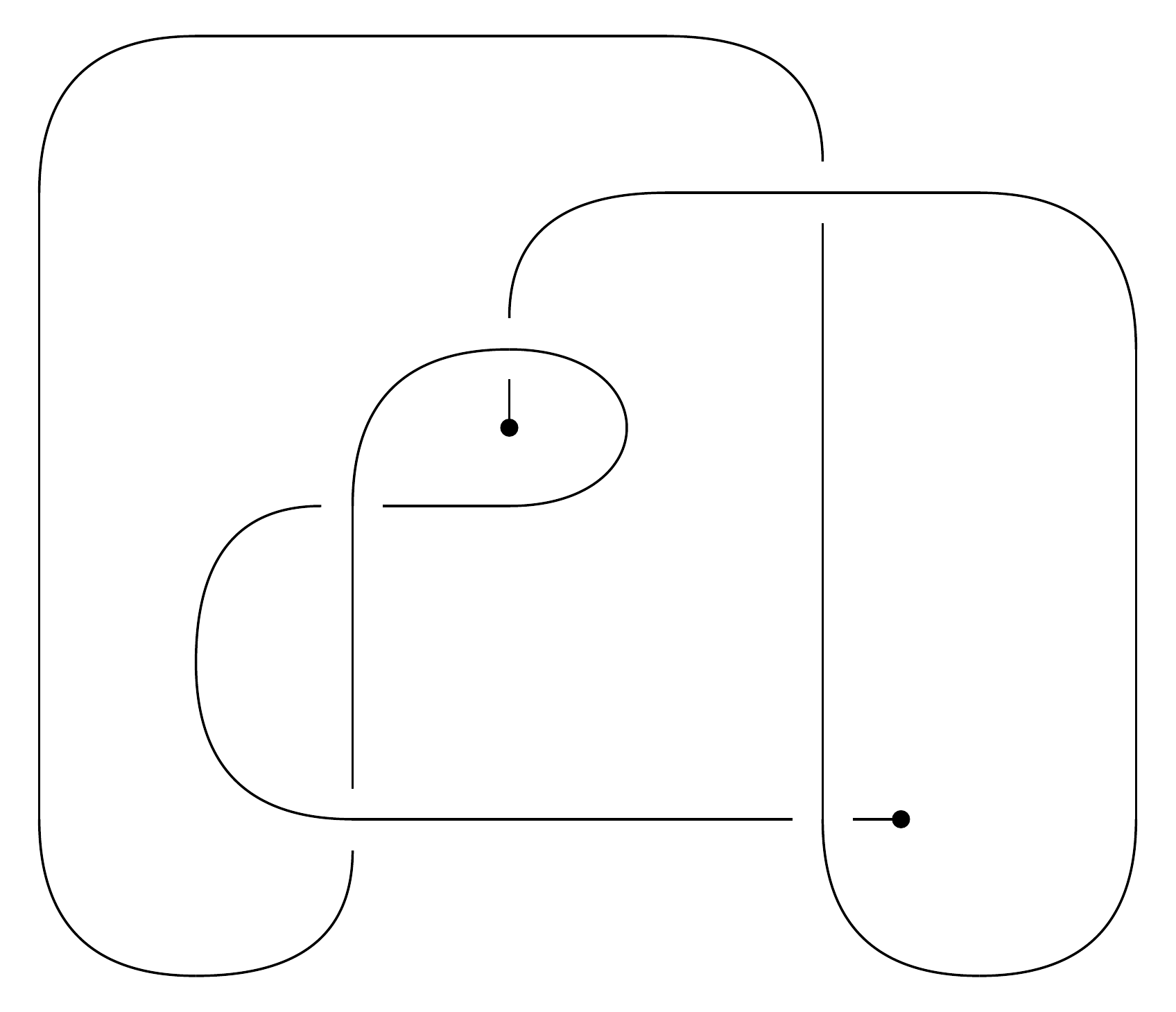}\\
\textcolor{black}{$5_{775}$}
\vspace{1cm}
\end{minipage}
\begin{minipage}[t]{.25\linewidth}
\centering
\includegraphics[width=0.9\textwidth,height=3.5cm,keepaspectratio]{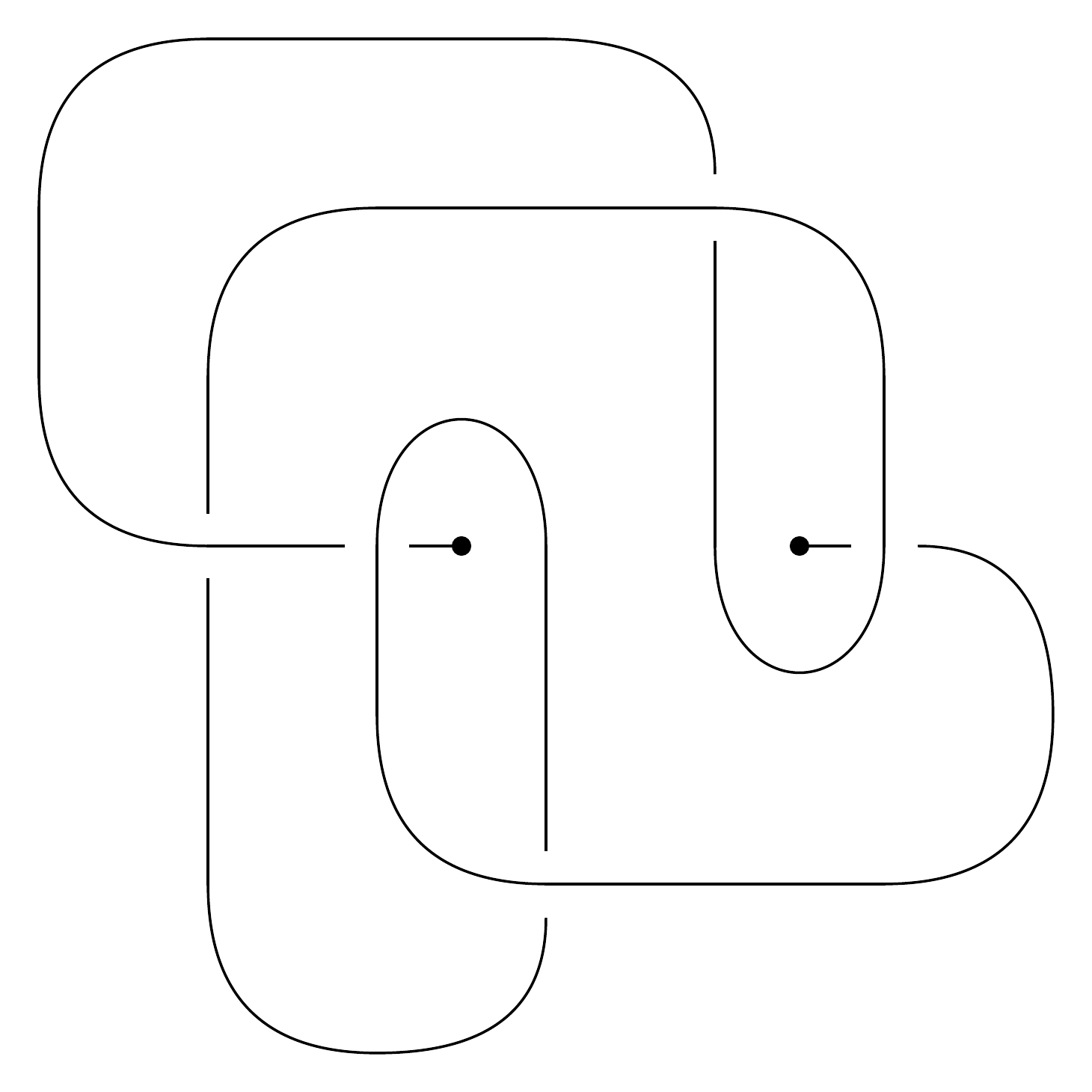}\\
\textcolor{black}{$5_{776}$}
\vspace{1cm}
\end{minipage}
\begin{minipage}[t]{.25\linewidth}
\centering
\includegraphics[width=0.9\textwidth,height=3.5cm,keepaspectratio]{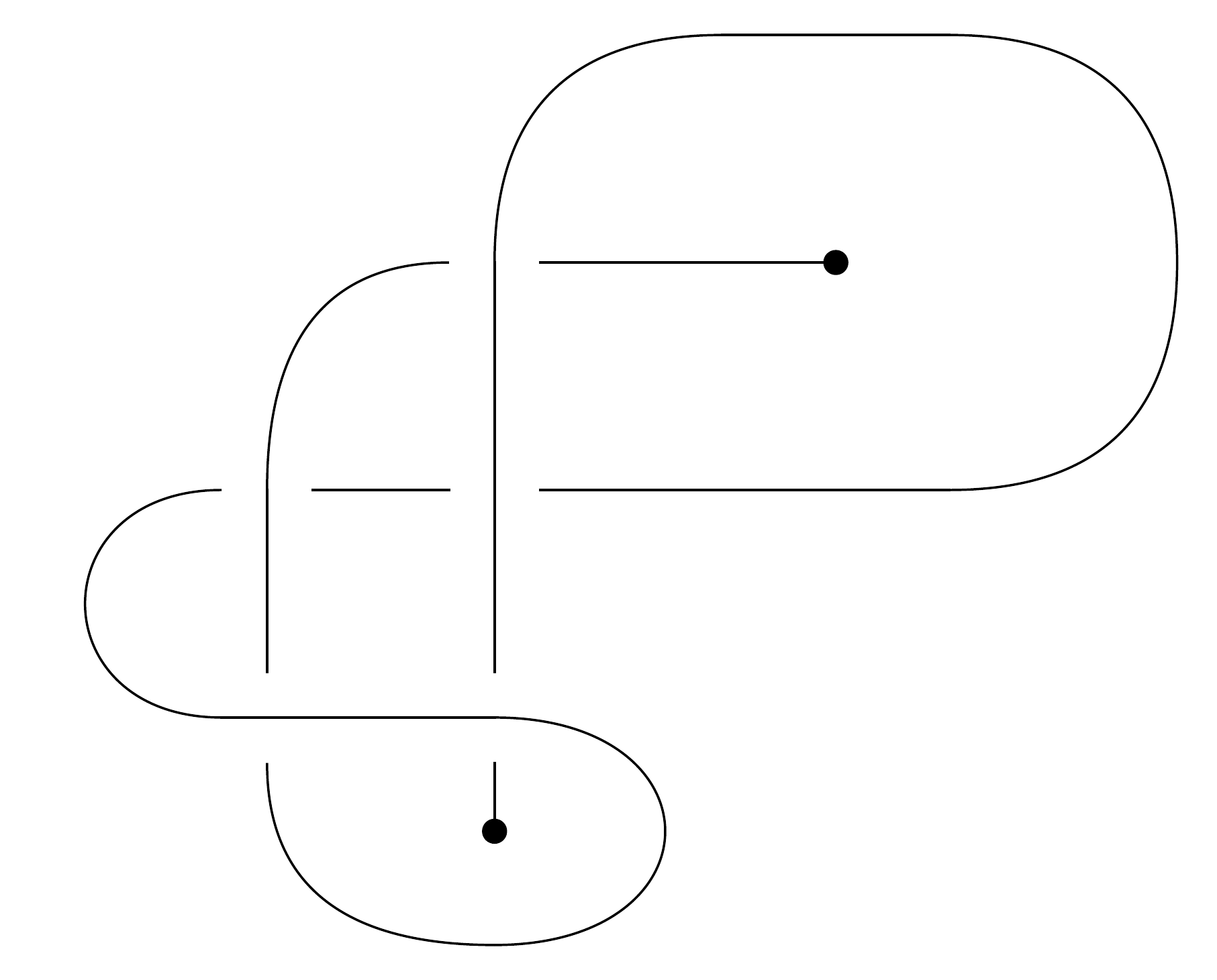}\\
\textcolor{black}{$5_{777}$}
\vspace{1cm}
\end{minipage}
\begin{minipage}[t]{.25\linewidth}
\centering
\includegraphics[width=0.9\textwidth,height=3.5cm,keepaspectratio]{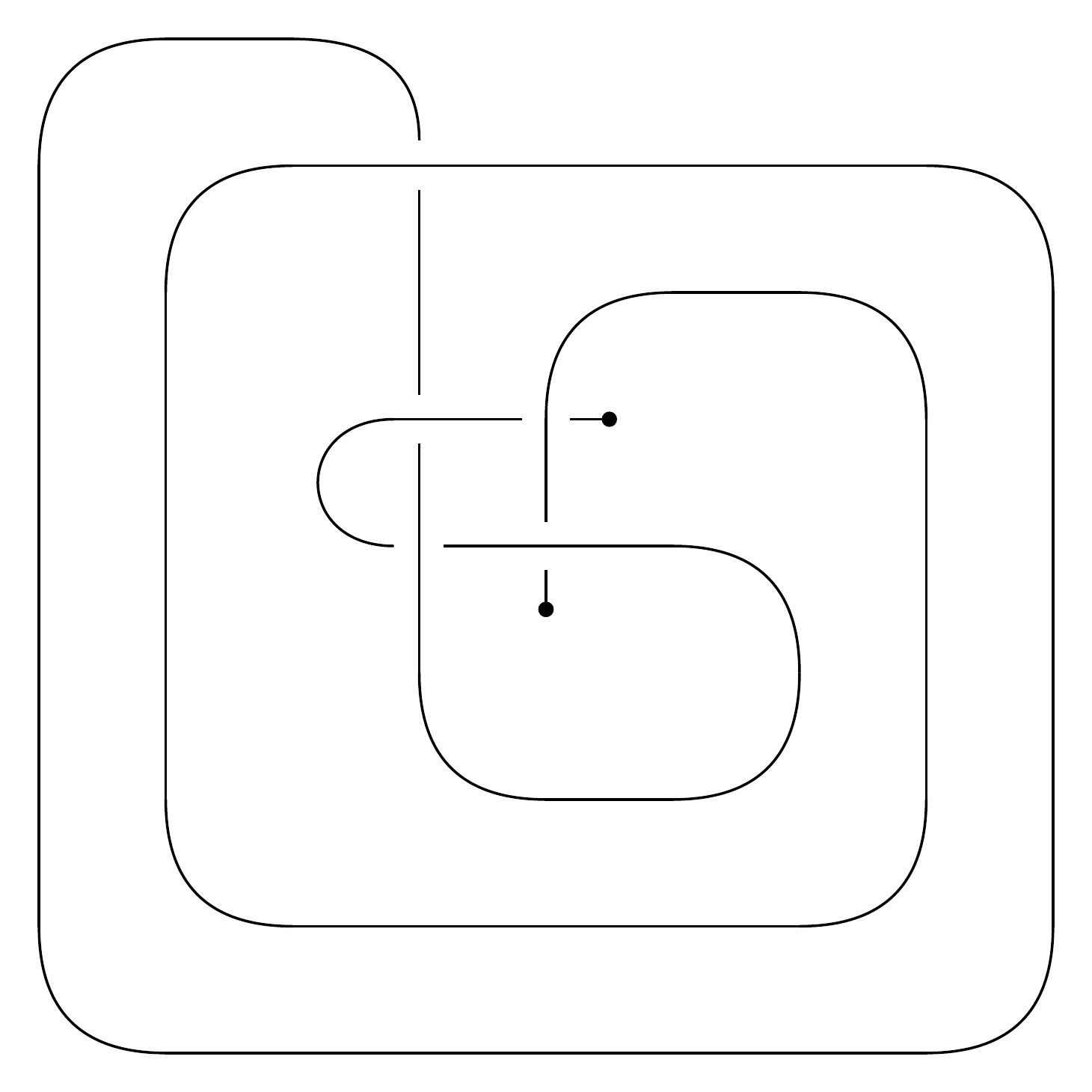}\\
\textcolor{black}{$5_{778}$}
\vspace{1cm}
\end{minipage}
\begin{minipage}[t]{.25\linewidth}
\centering
\includegraphics[width=0.9\textwidth,height=3.5cm,keepaspectratio]{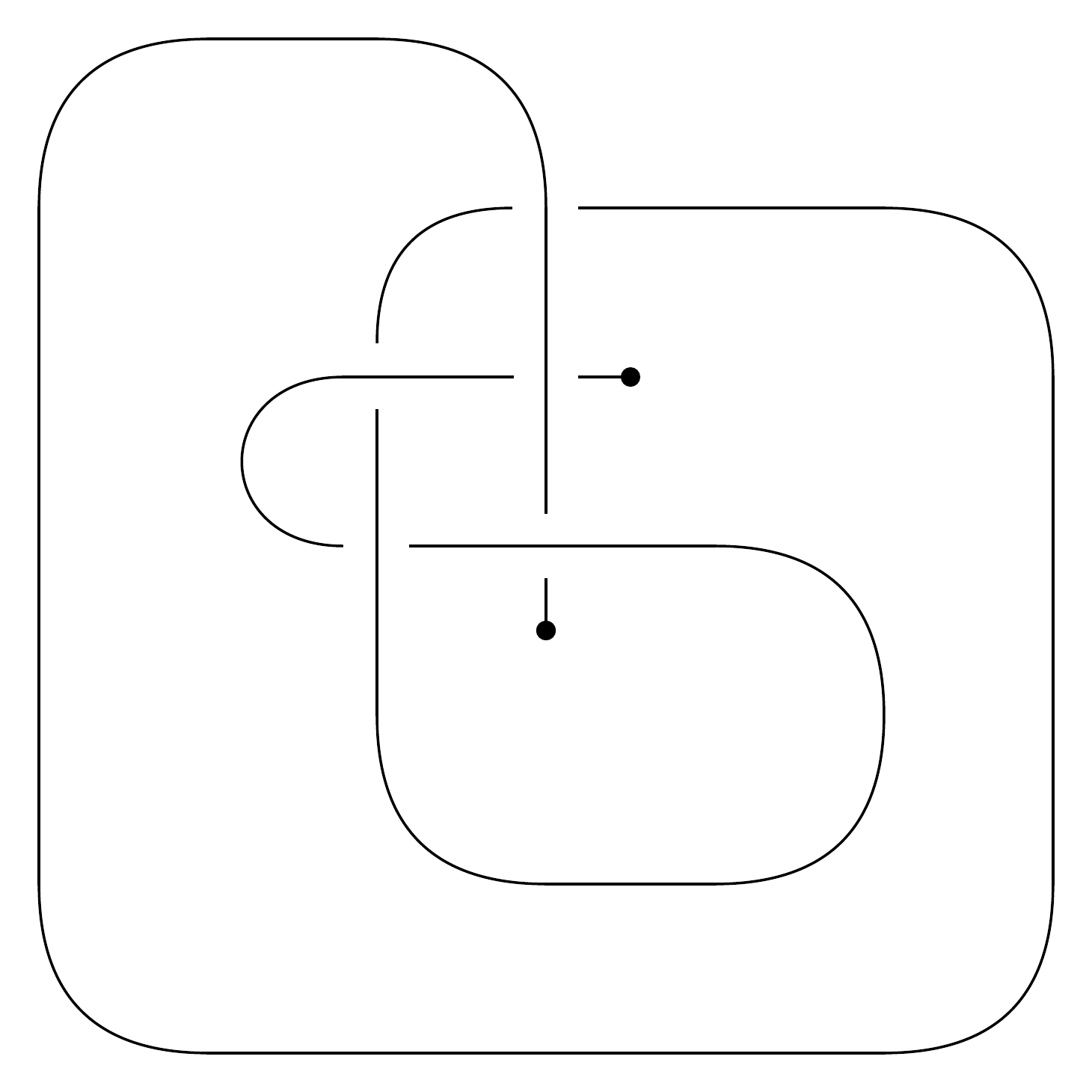}\\
\textcolor{black}{$5_{779}$}
\vspace{1cm}
\end{minipage}
\begin{minipage}[t]{.25\linewidth}
\centering
\includegraphics[width=0.9\textwidth,height=3.5cm,keepaspectratio]{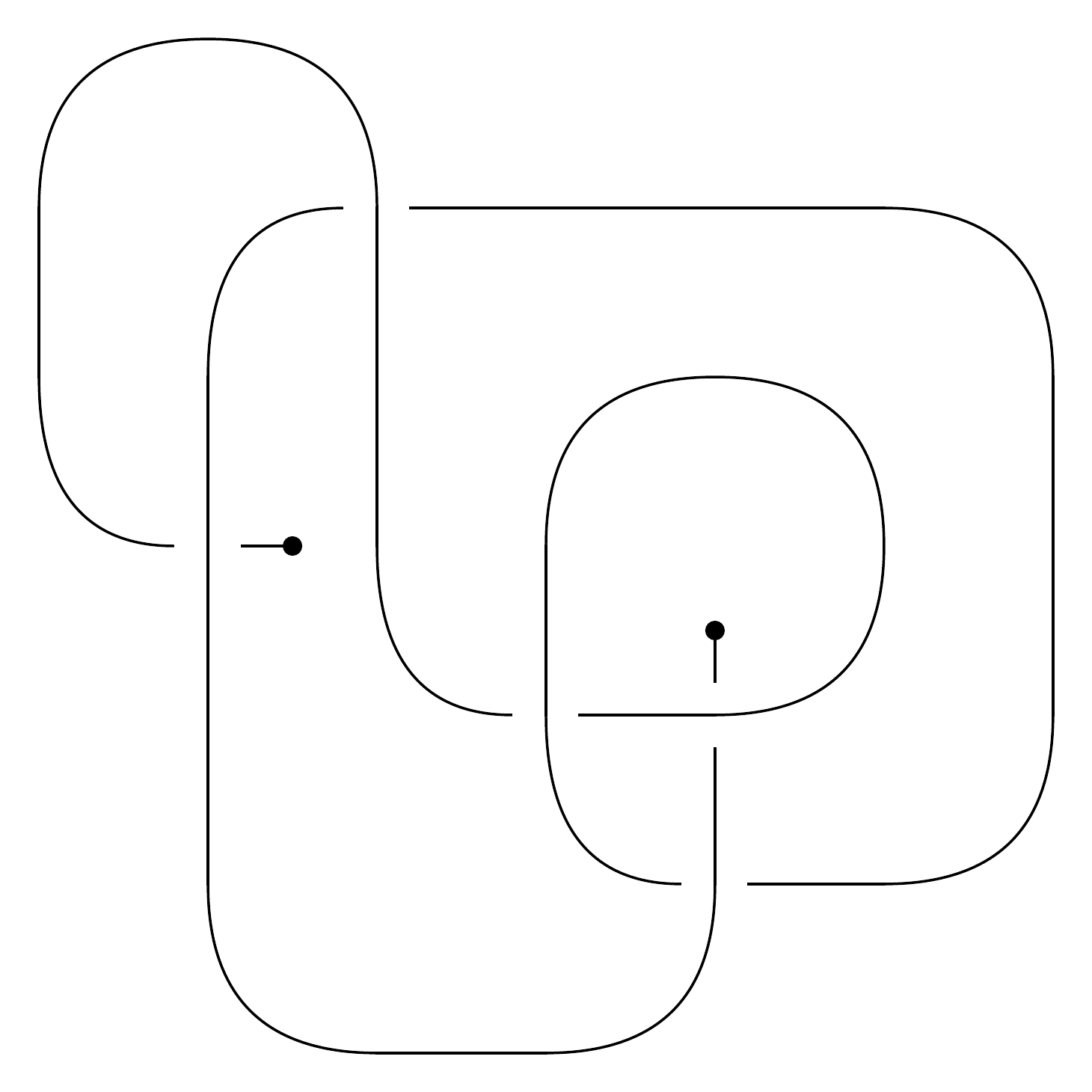}\\
\textcolor{black}{$5_{780}$}
\vspace{1cm}
\end{minipage}
\begin{minipage}[t]{.25\linewidth}
\centering
\includegraphics[width=0.9\textwidth,height=3.5cm,keepaspectratio]{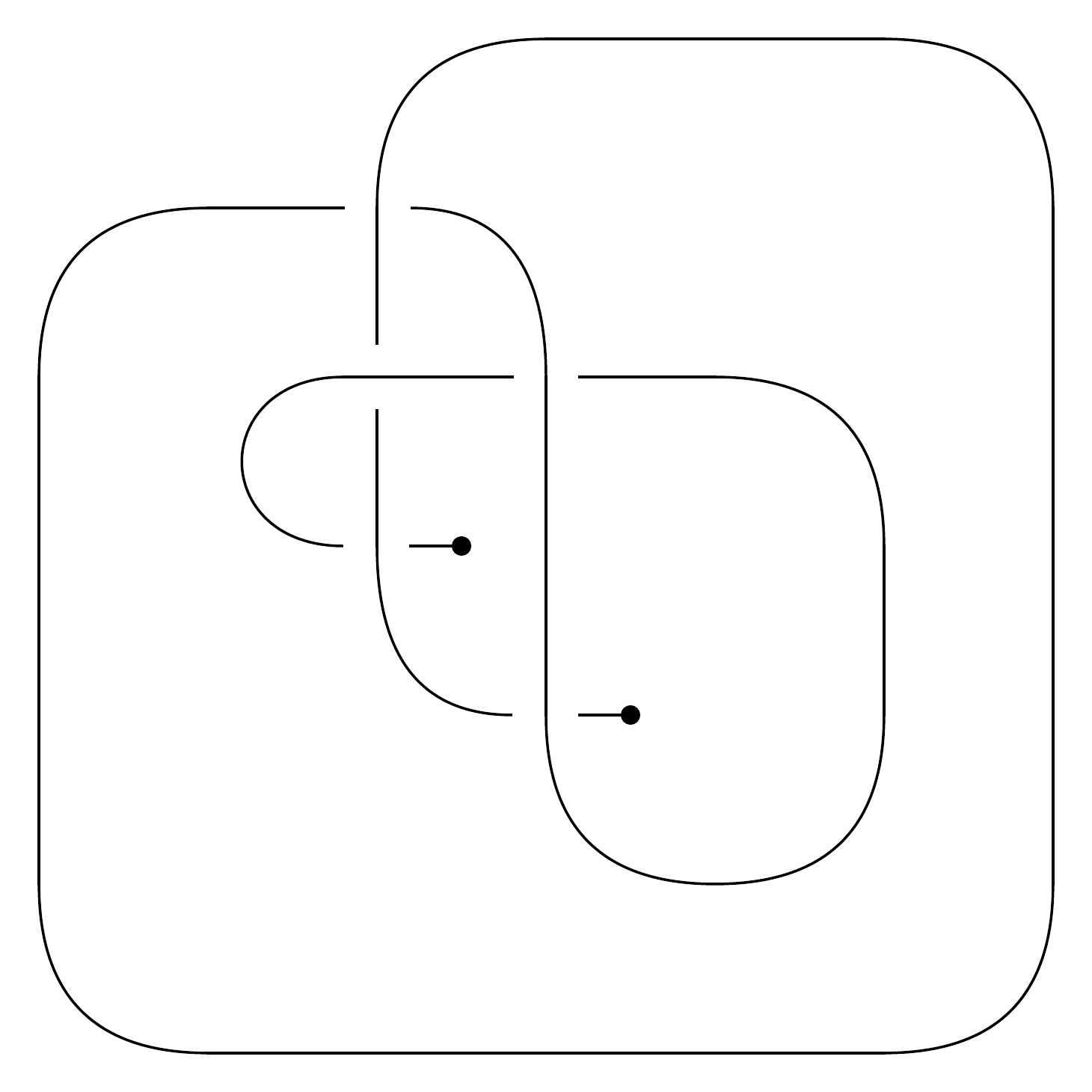}\\
\textcolor{black}{$5_{781}$}
\vspace{1cm}
\end{minipage}
\begin{minipage}[t]{.25\linewidth}
\centering
\includegraphics[width=0.9\textwidth,height=3.5cm,keepaspectratio]{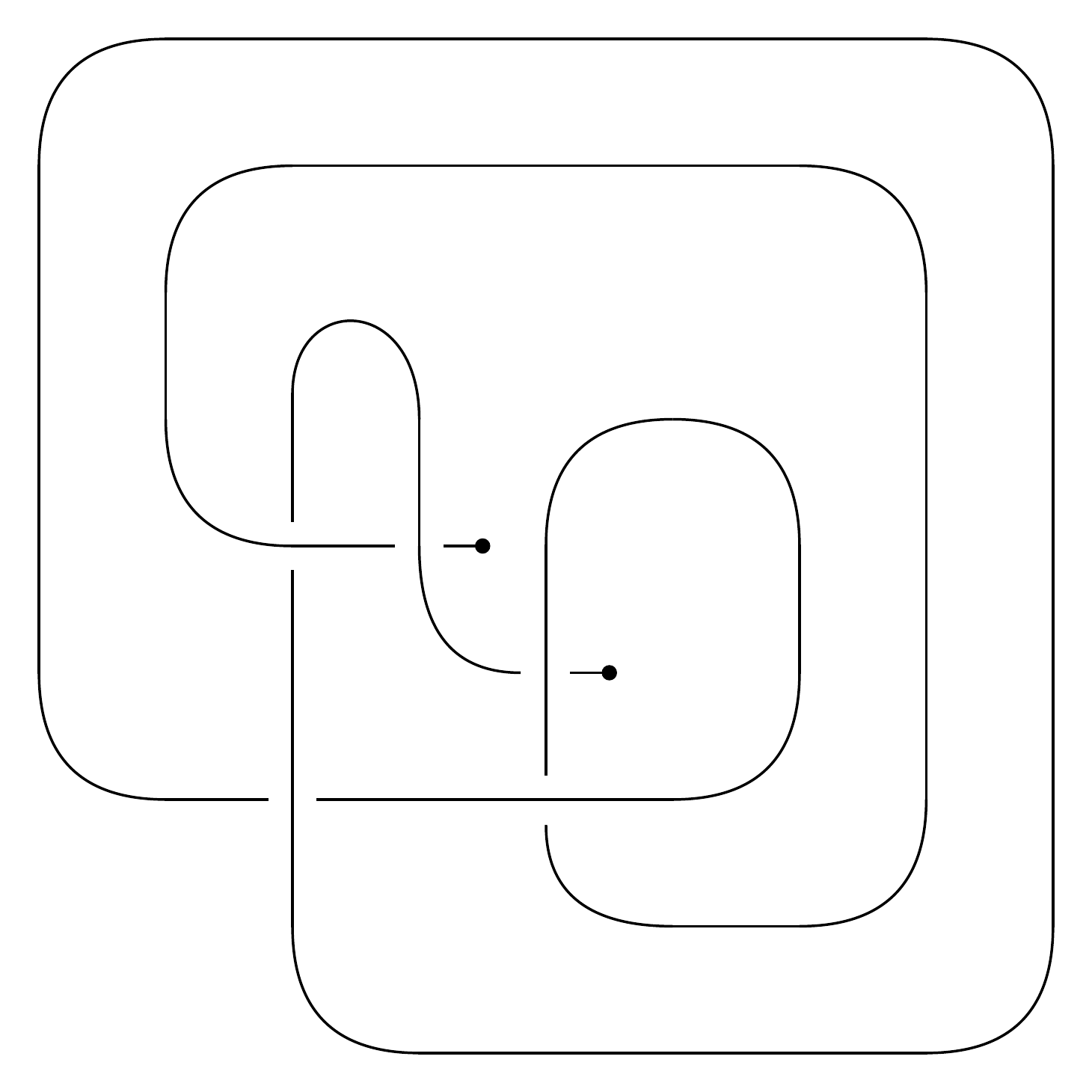}\\
\textcolor{black}{$5_{782}$}
\vspace{1cm}
\end{minipage}
\begin{minipage}[t]{.25\linewidth}
\centering
\includegraphics[width=0.9\textwidth,height=3.5cm,keepaspectratio]{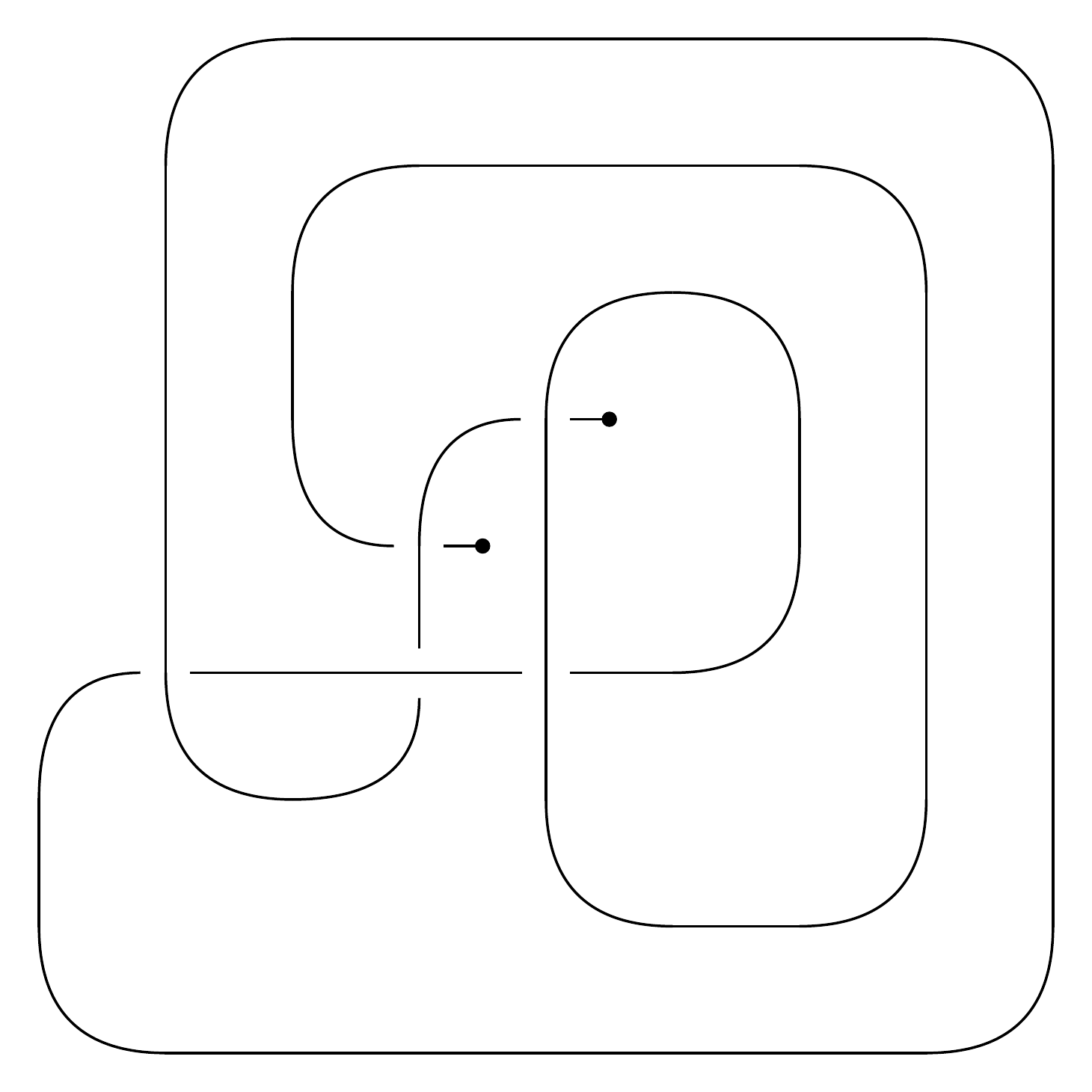}\\
\textcolor{black}{$5_{783}$}
\vspace{1cm}
\end{minipage}
\begin{minipage}[t]{.25\linewidth}
\centering
\includegraphics[width=0.9\textwidth,height=3.5cm,keepaspectratio]{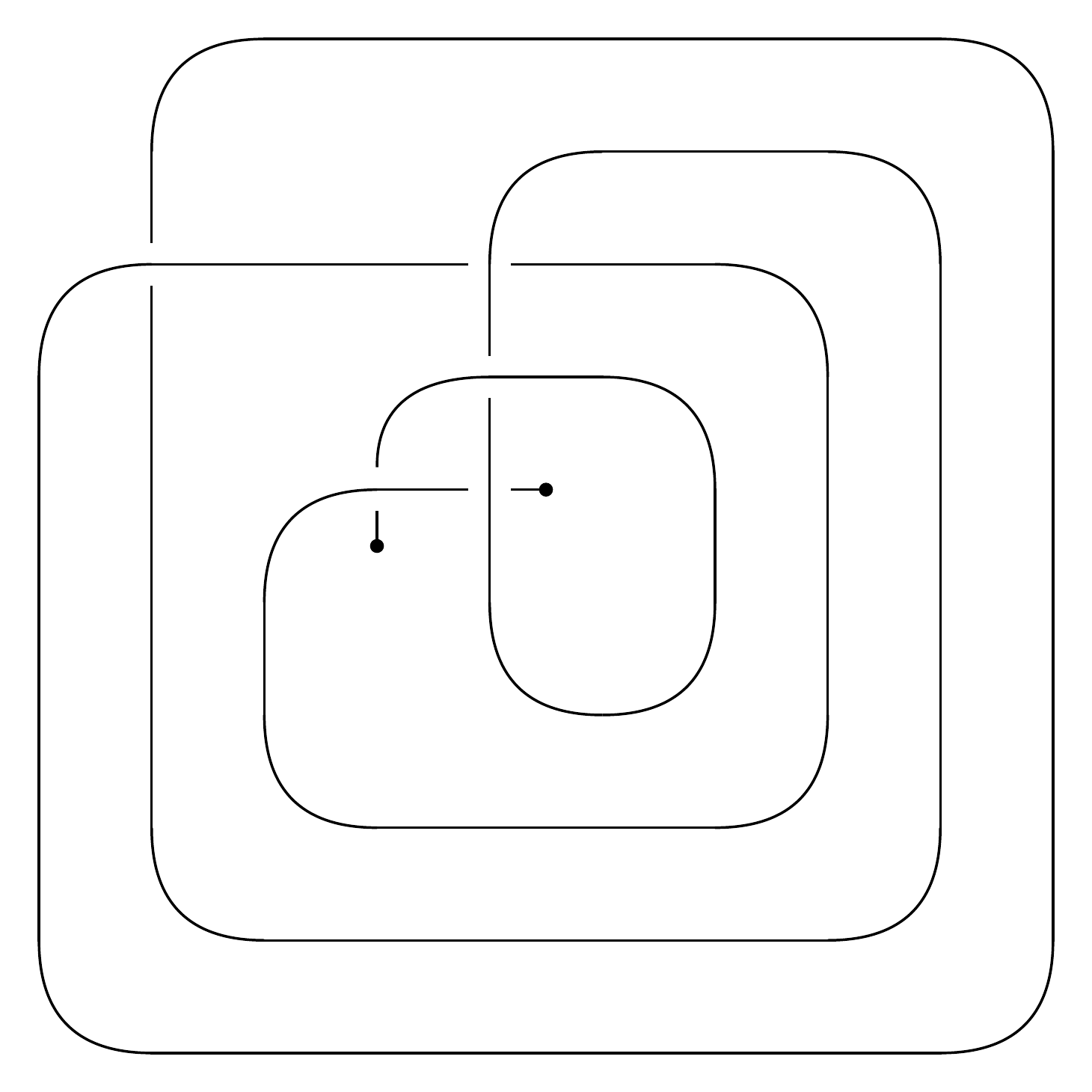}\\
\textcolor{black}{$5_{784}$}
\vspace{1cm}
\end{minipage}
\begin{minipage}[t]{.25\linewidth}
\centering
\includegraphics[width=0.9\textwidth,height=3.5cm,keepaspectratio]{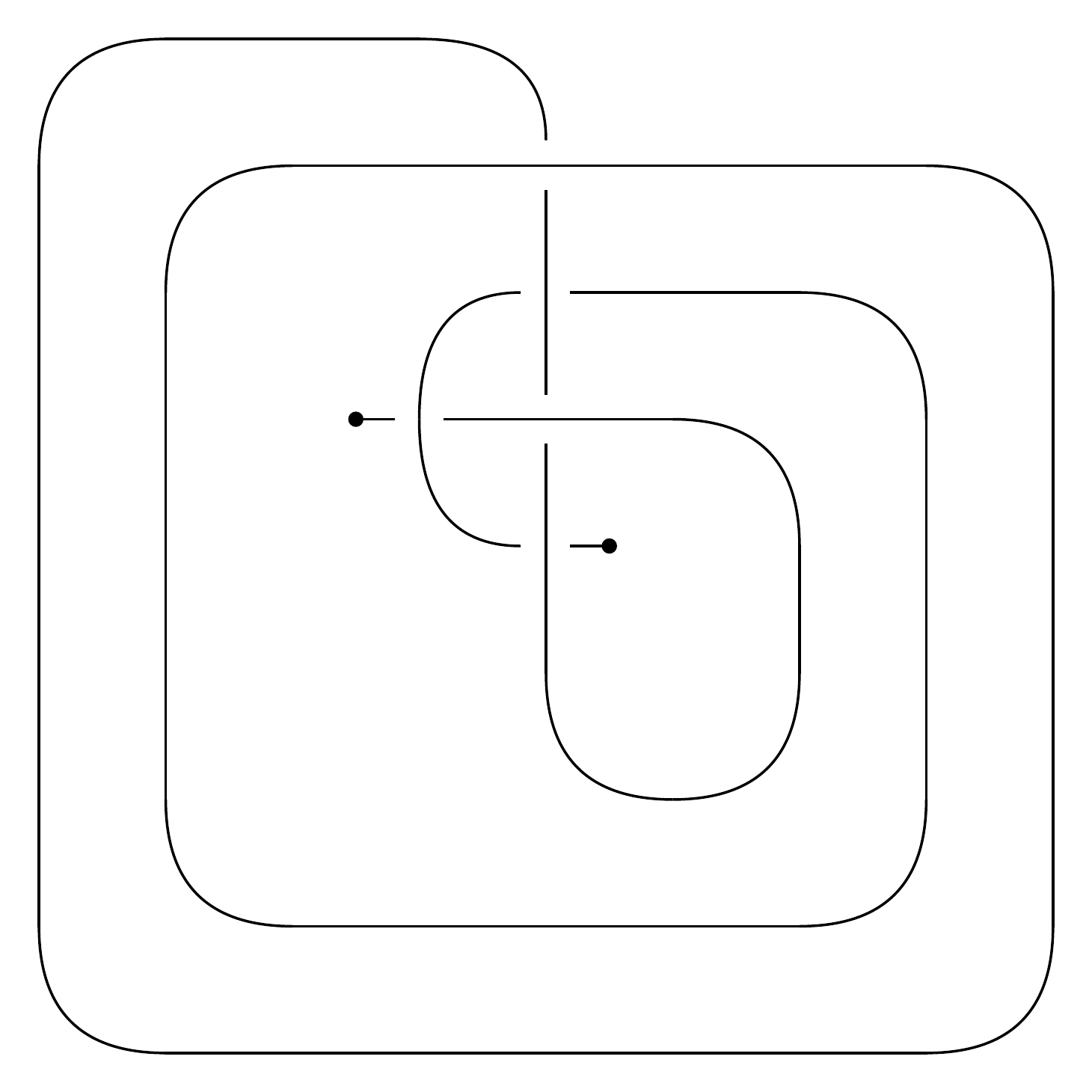}\\
\textcolor{black}{$5_{785}$}
\vspace{1cm}
\end{minipage}
\begin{minipage}[t]{.25\linewidth}
\centering
\includegraphics[width=0.9\textwidth,height=3.5cm,keepaspectratio]{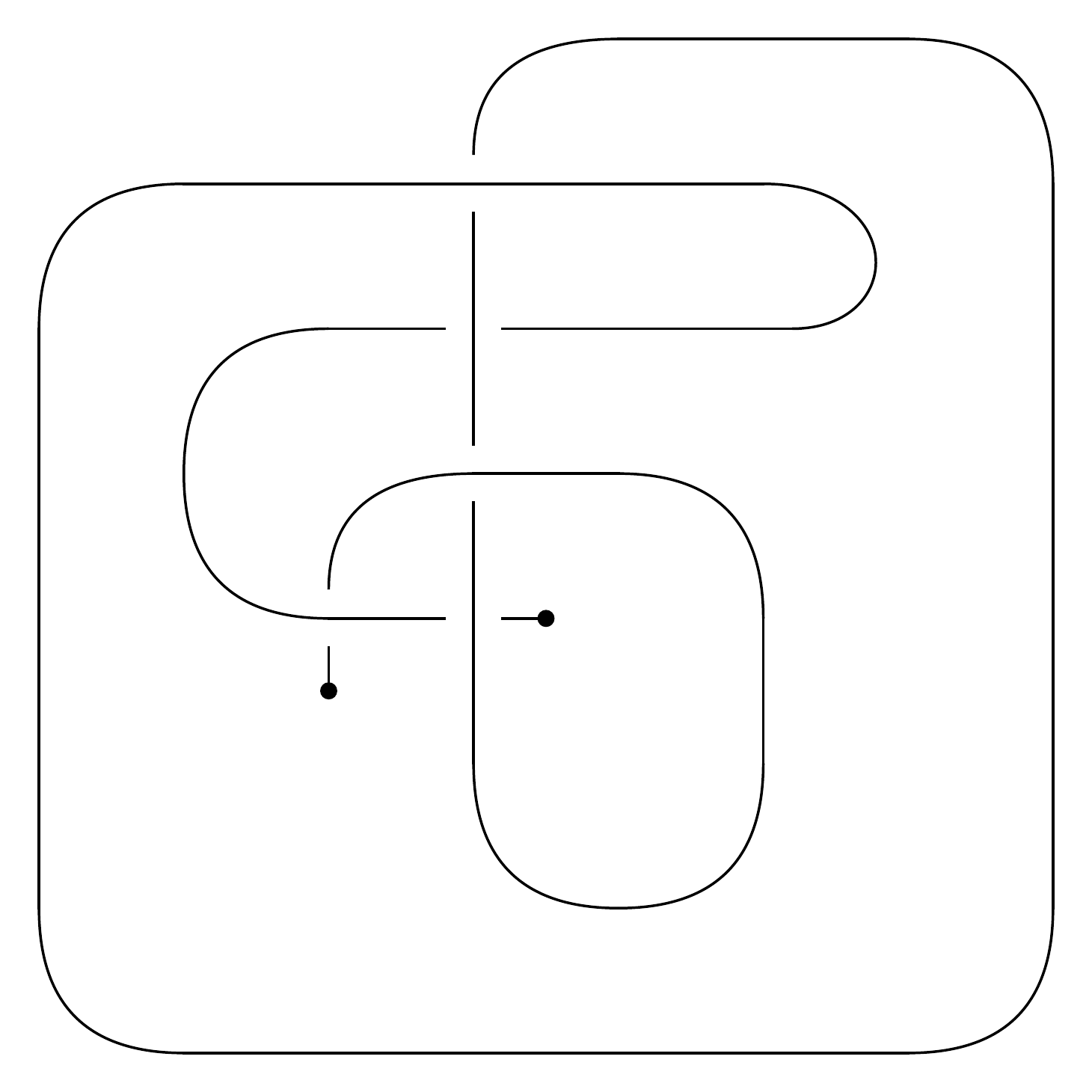}\\
\textcolor{black}{$5_{786}$}
\vspace{1cm}
\end{minipage}
\begin{minipage}[t]{.25\linewidth}
\centering
\includegraphics[width=0.9\textwidth,height=3.5cm,keepaspectratio]{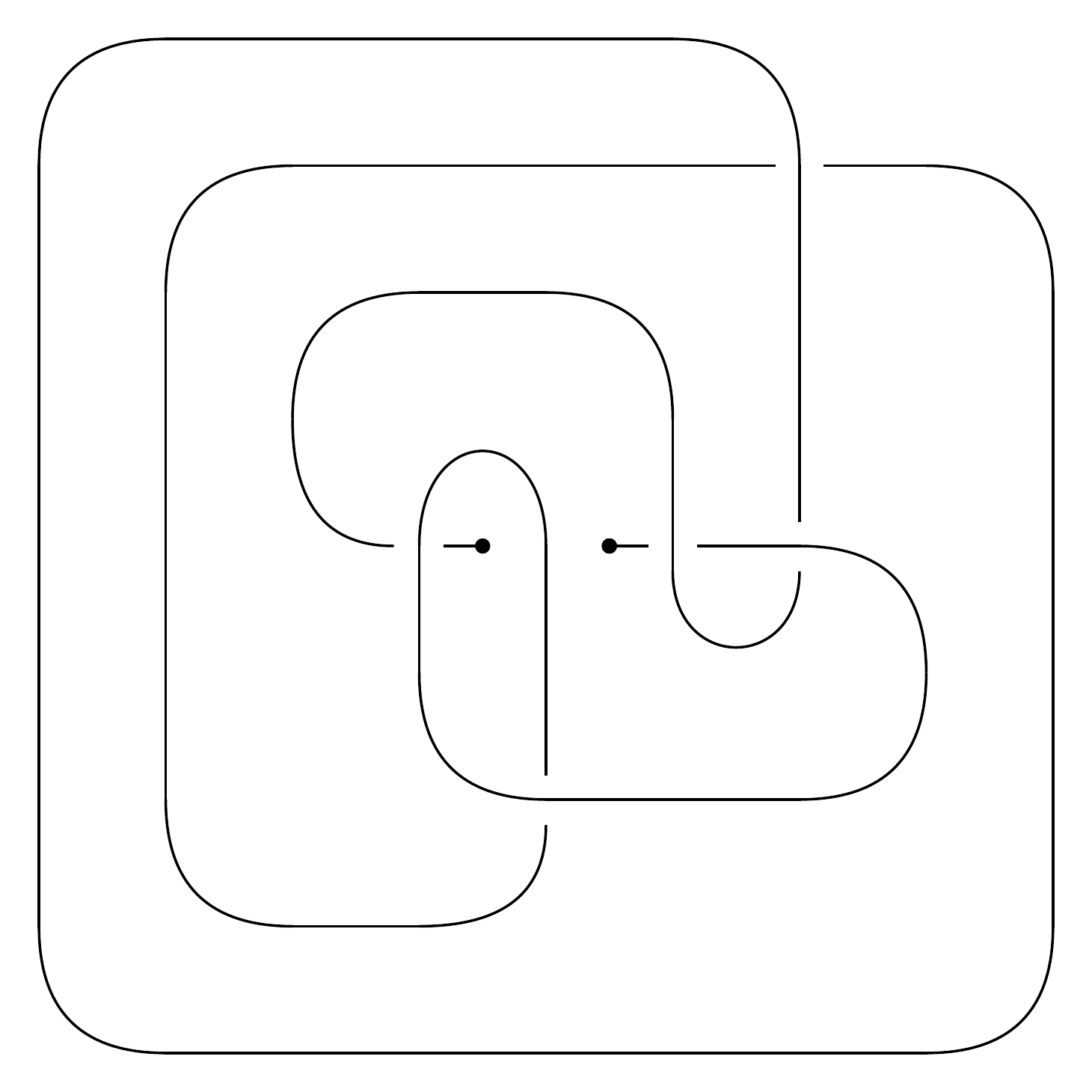}\\
\textcolor{black}{$5_{787}$}
\vspace{1cm}
\end{minipage}
\begin{minipage}[t]{.25\linewidth}
\centering
\includegraphics[width=0.9\textwidth,height=3.5cm,keepaspectratio]{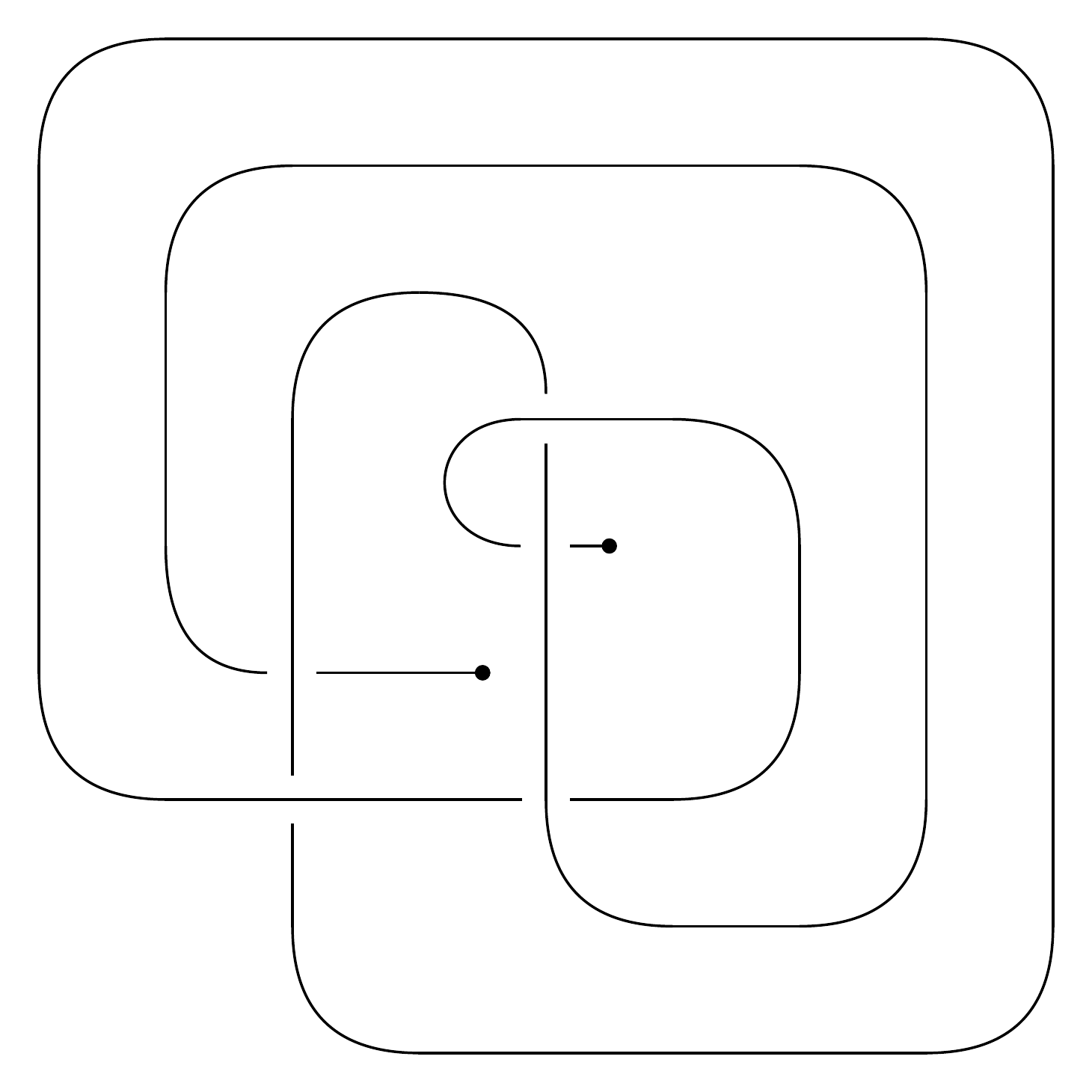}\\
\textcolor{black}{$5_{788}$}
\vspace{1cm}
\end{minipage}
\begin{minipage}[t]{.25\linewidth}
\centering
\includegraphics[width=0.9\textwidth,height=3.5cm,keepaspectratio]{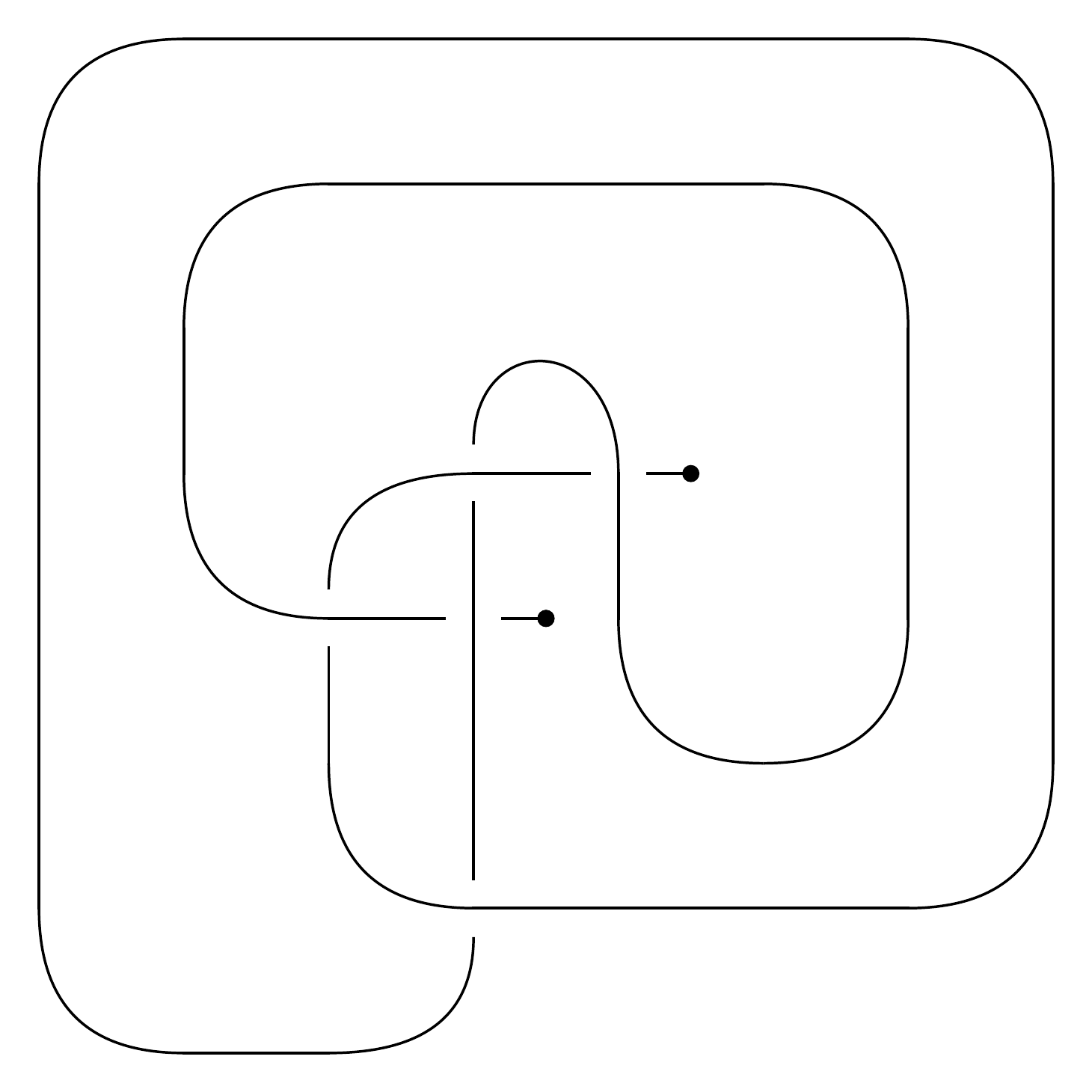}\\
\textcolor{black}{$5_{789}$}
\vspace{1cm}
\end{minipage}
\begin{minipage}[t]{.25\linewidth}
\centering
\includegraphics[width=0.9\textwidth,height=3.5cm,keepaspectratio]{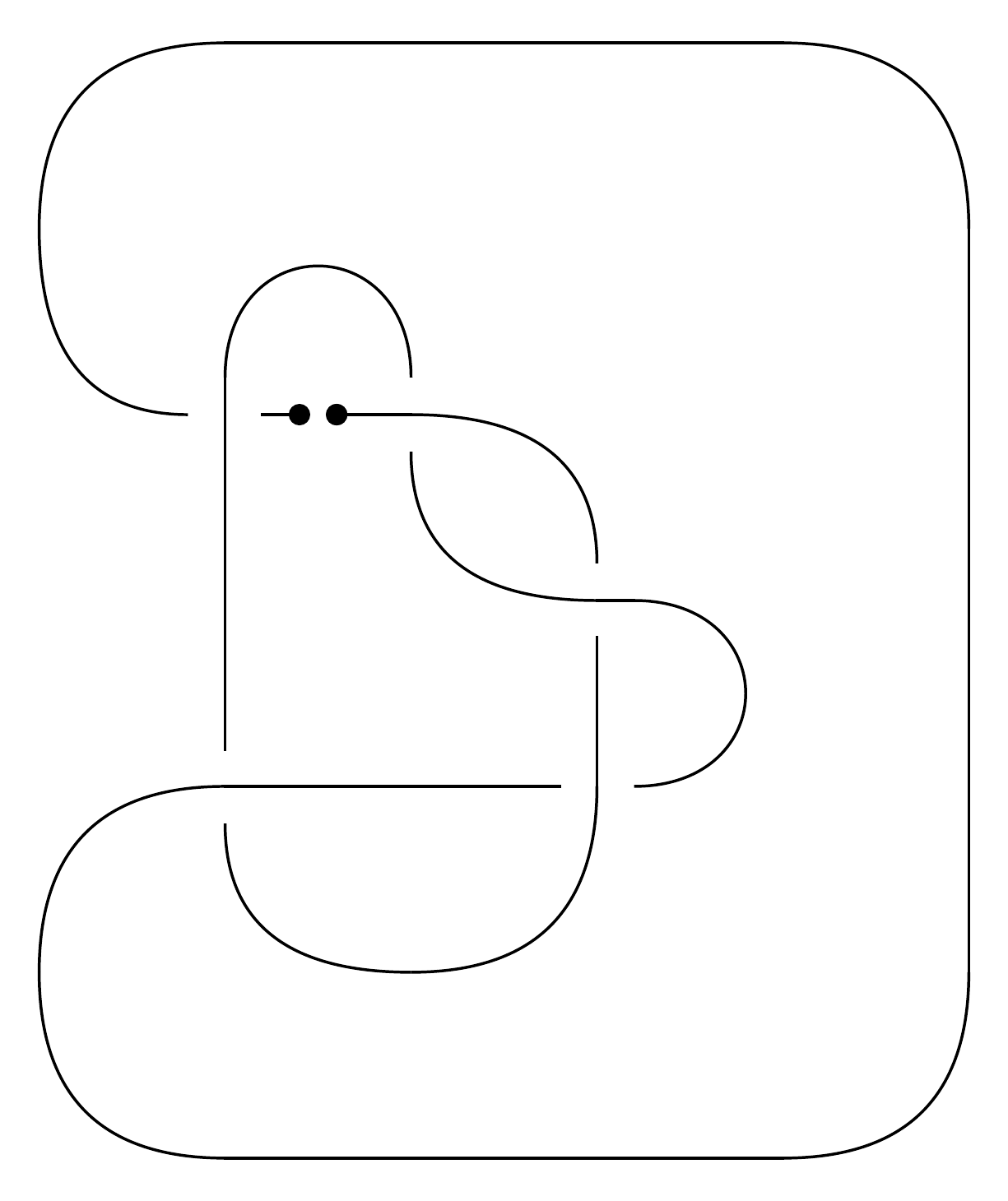}\\
\textcolor{black}{$5_{790}$}
\vspace{1cm}
\end{minipage}
\begin{minipage}[t]{.25\linewidth}
\centering
\includegraphics[width=0.9\textwidth,height=3.5cm,keepaspectratio]{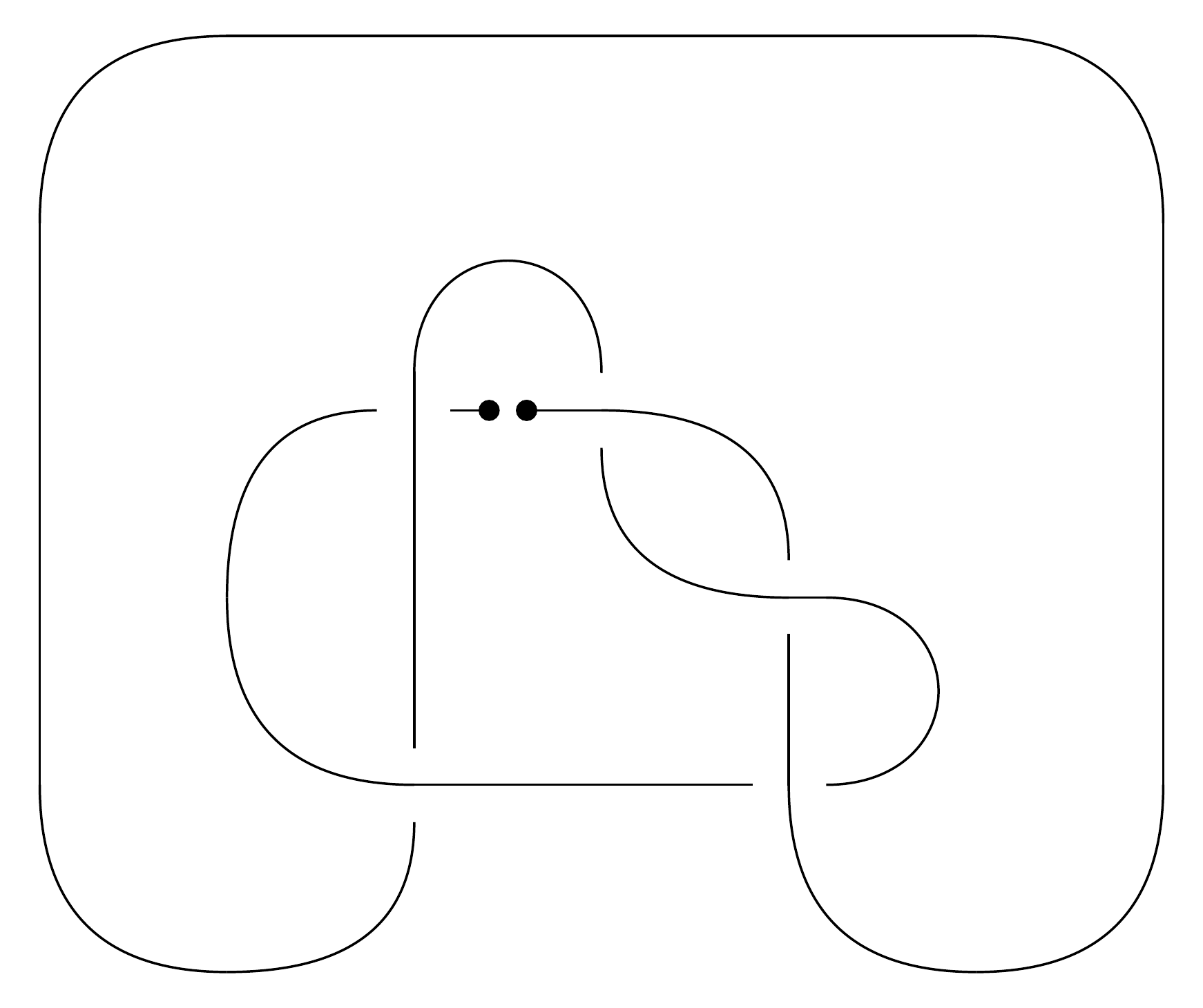}\\
\textcolor{black}{$5_{791}$}
\vspace{1cm}
\end{minipage}
\begin{minipage}[t]{.25\linewidth}
\centering
\includegraphics[width=0.9\textwidth,height=3.5cm,keepaspectratio]{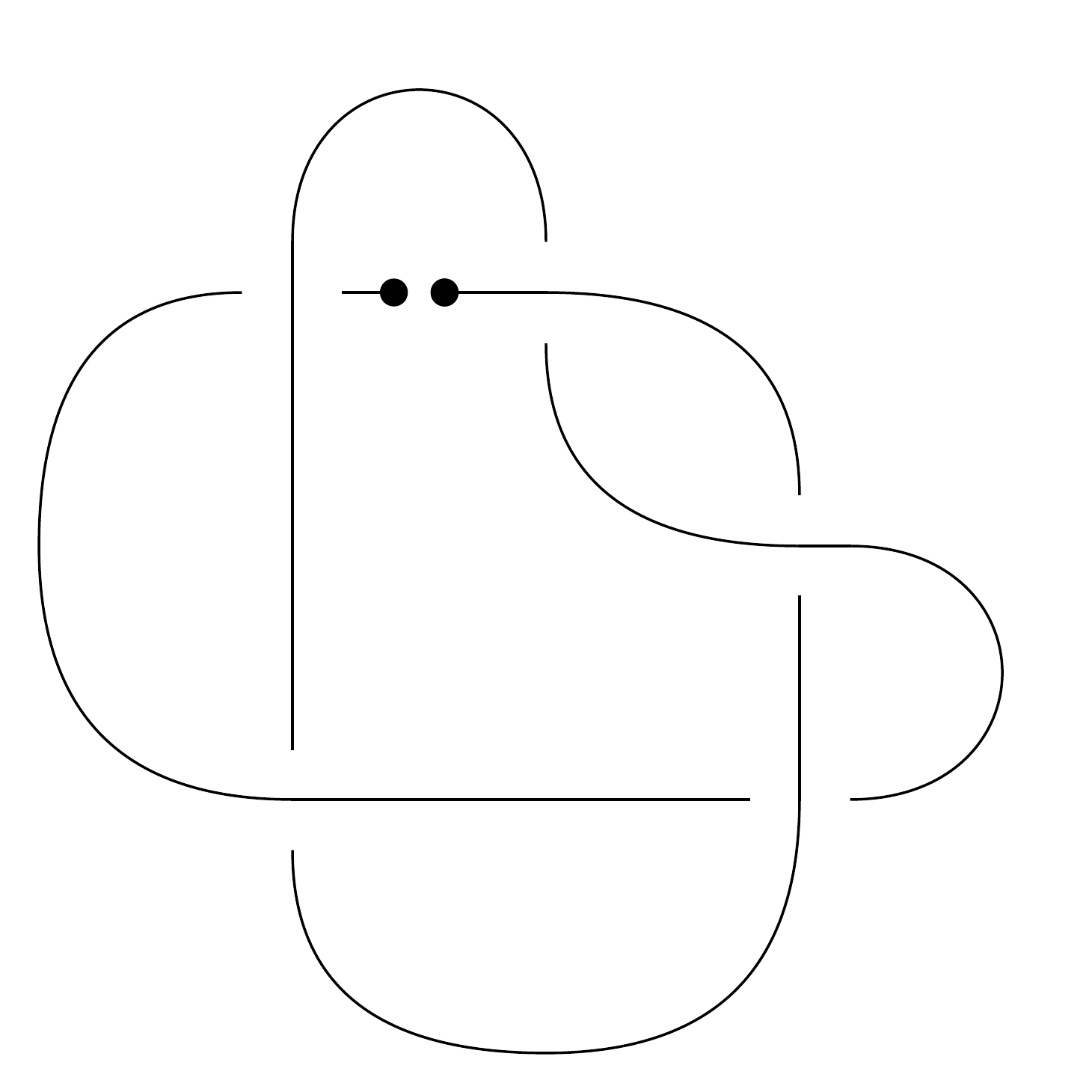}\\
\textcolor{black}{$5_{792}$}
\vspace{1cm}
\end{minipage}
\begin{minipage}[t]{.25\linewidth}
\centering
\includegraphics[width=0.9\textwidth,height=3.5cm,keepaspectratio]{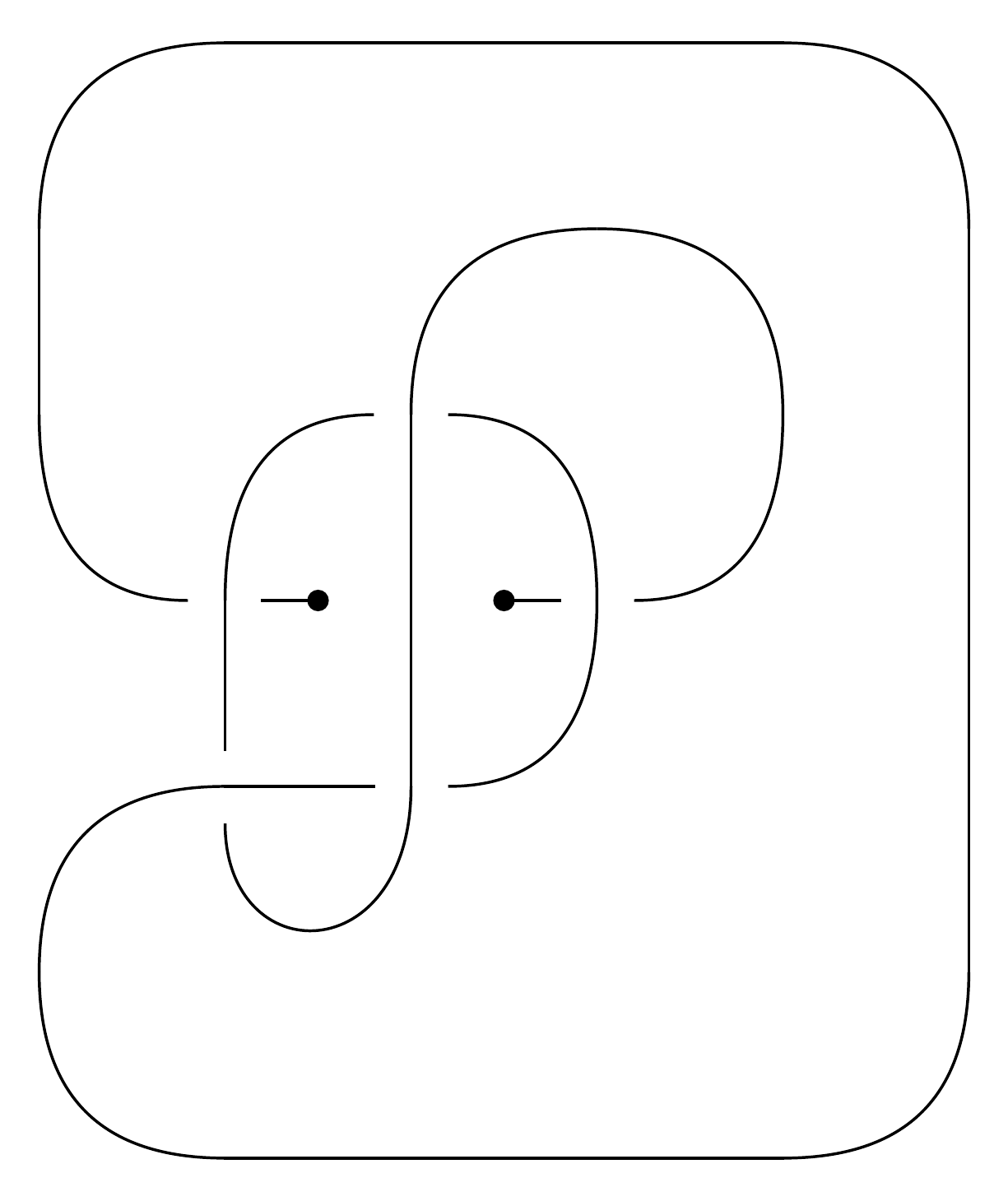}\\
\textcolor{black}{$5_{793}$}
\vspace{1cm}
\end{minipage}
\begin{minipage}[t]{.25\linewidth}
\centering
\includegraphics[width=0.9\textwidth,height=3.5cm,keepaspectratio]{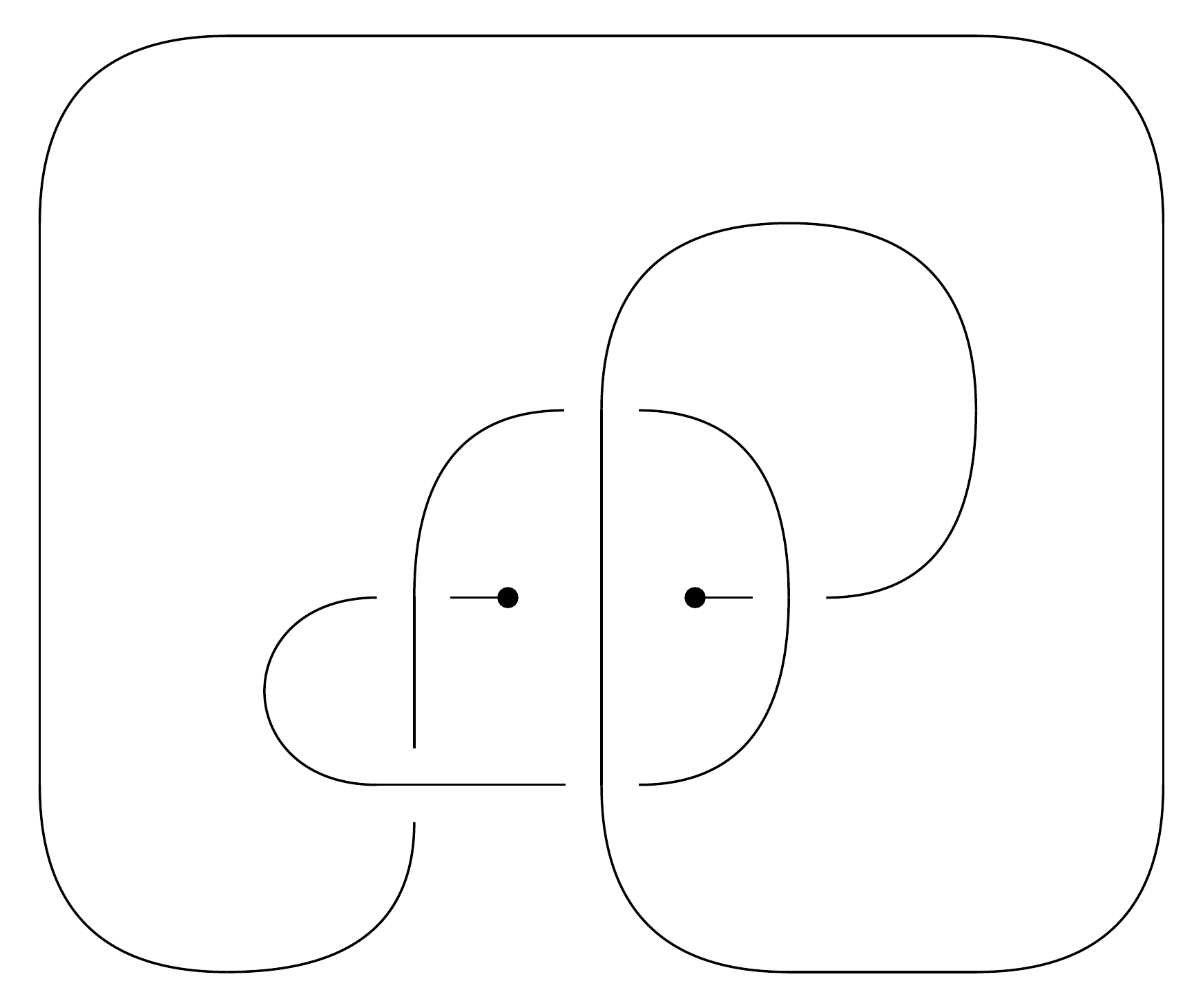}\\
\textcolor{black}{$5_{794}$}
\vspace{1cm}
\end{minipage}
\begin{minipage}[t]{.25\linewidth}
\centering
\includegraphics[width=0.9\textwidth,height=3.5cm,keepaspectratio]{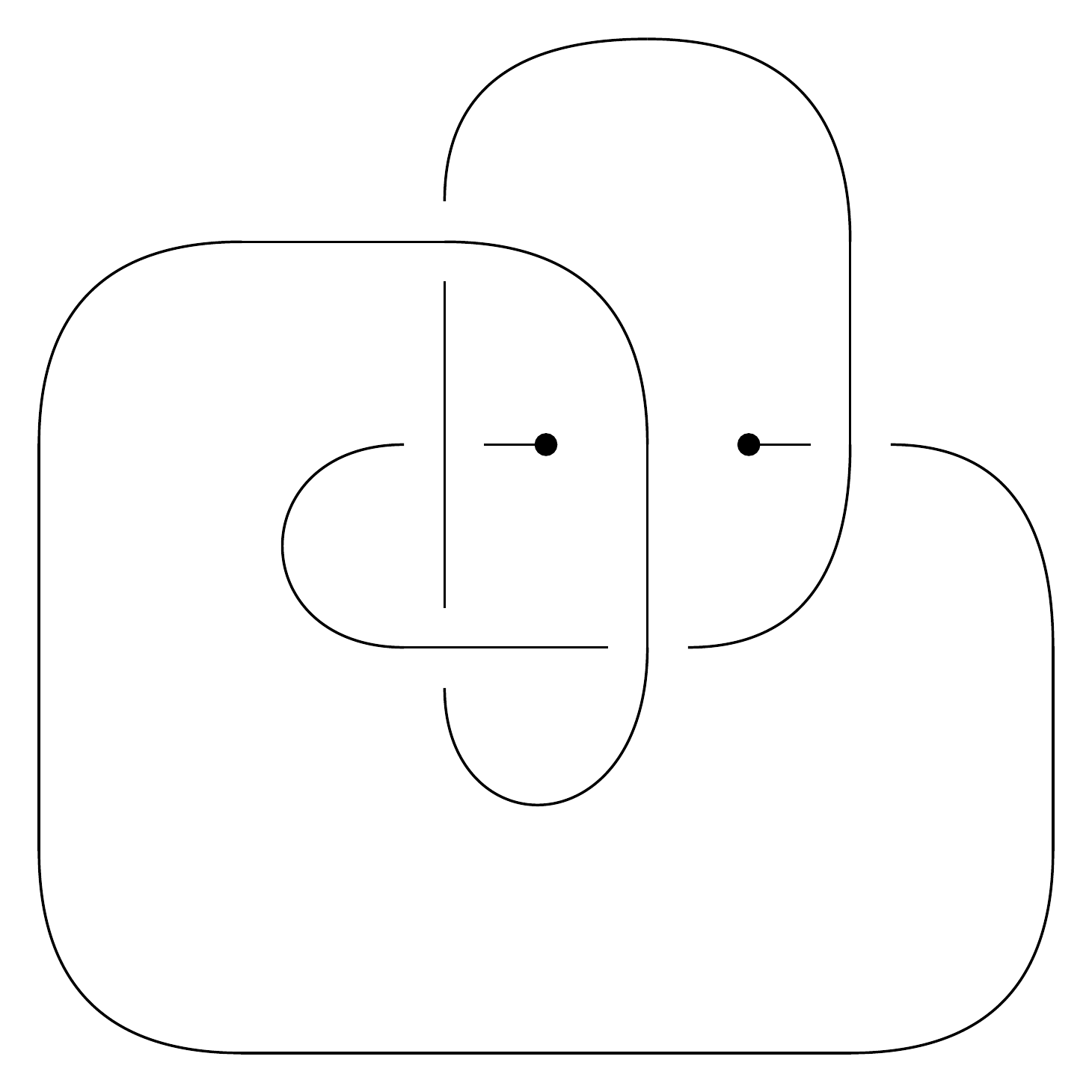}\\
\textcolor{black}{$5_{795}$}
\vspace{1cm}
\end{minipage}
\begin{minipage}[t]{.25\linewidth}
\centering
\includegraphics[width=0.9\textwidth,height=3.5cm,keepaspectratio]{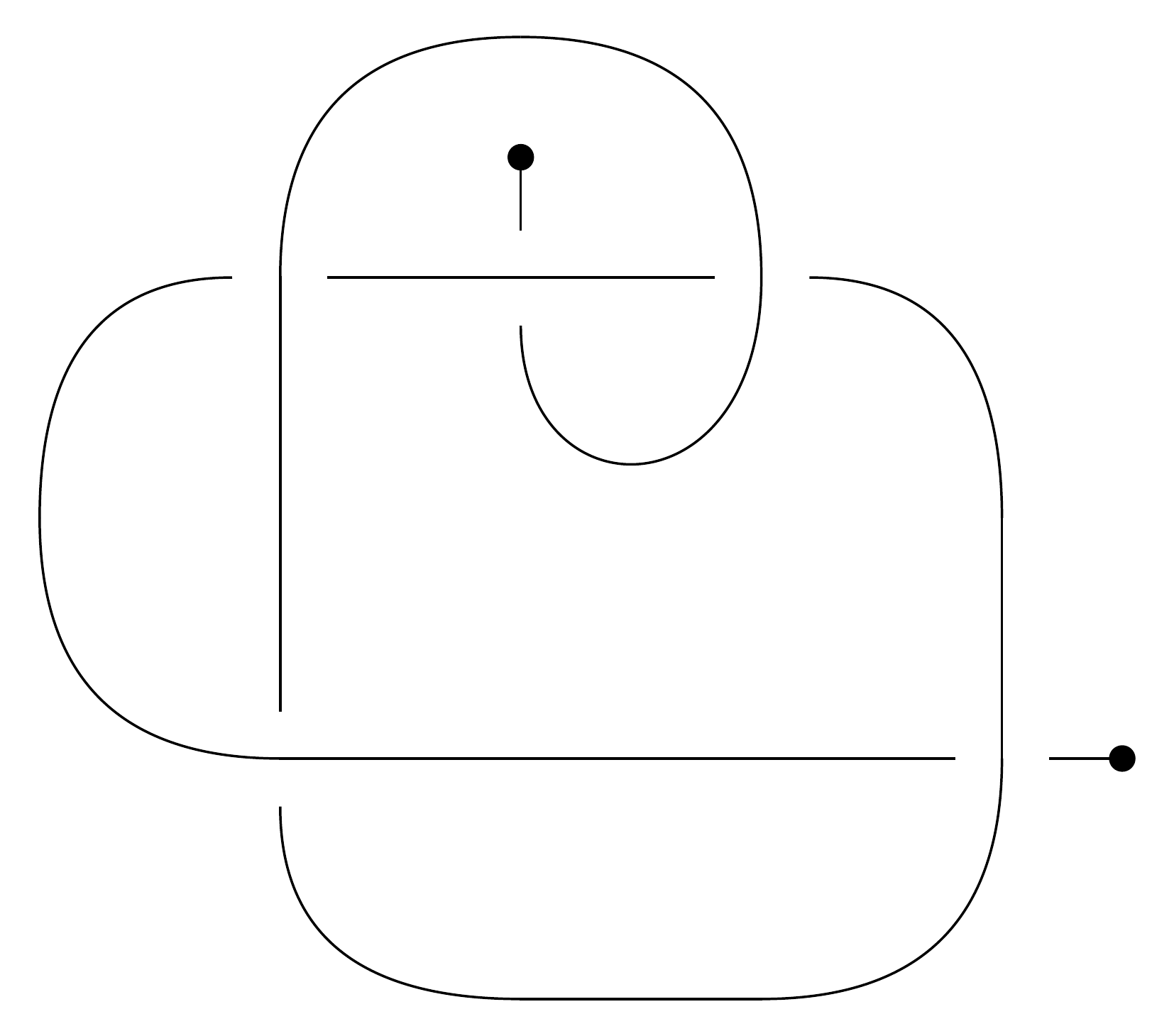}\\
\textcolor{black}{$5_{796}$}
\vspace{1cm}
\end{minipage}
\begin{minipage}[t]{.25\linewidth}
\centering
\includegraphics[width=0.9\textwidth,height=3.5cm,keepaspectratio]{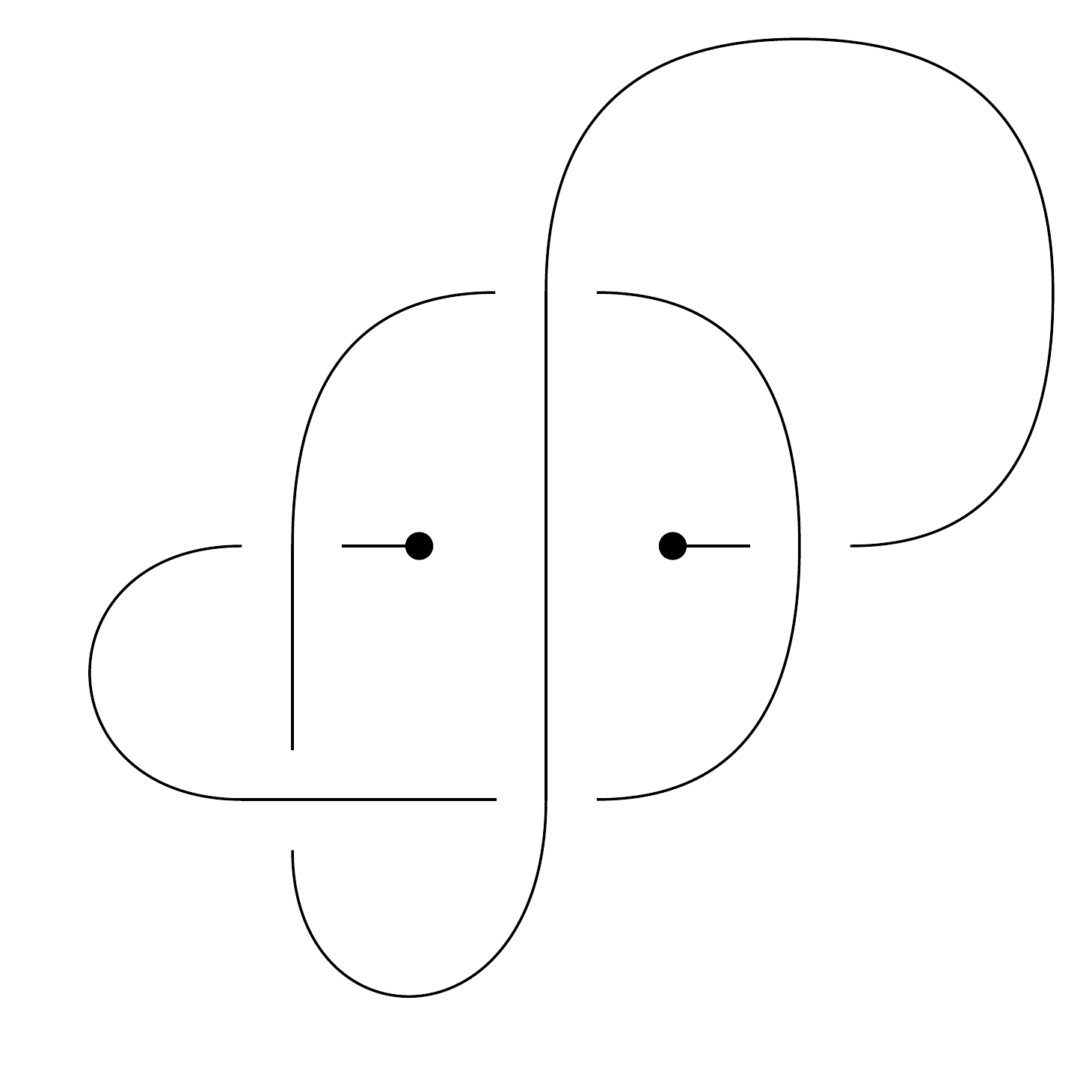}\\
\textcolor{black}{$5_{797}$}
\vspace{1cm}
\end{minipage}
\begin{minipage}[t]{.25\linewidth}
\centering
\includegraphics[width=0.9\textwidth,height=3.5cm,keepaspectratio]{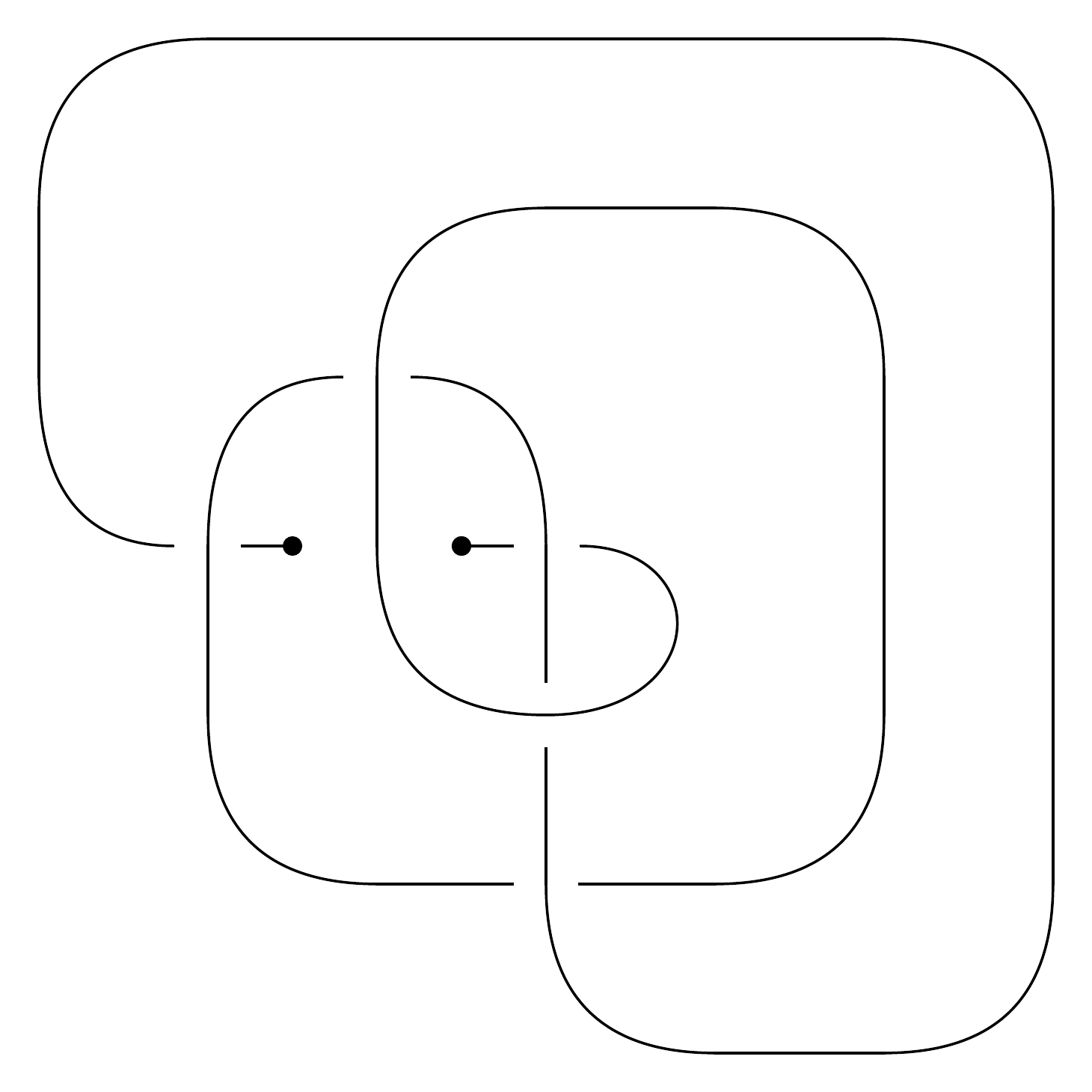}\\
\textcolor{black}{$5_{798}$}
\vspace{1cm}
\end{minipage}
\begin{minipage}[t]{.25\linewidth}
\centering
\includegraphics[width=0.9\textwidth,height=3.5cm,keepaspectratio]{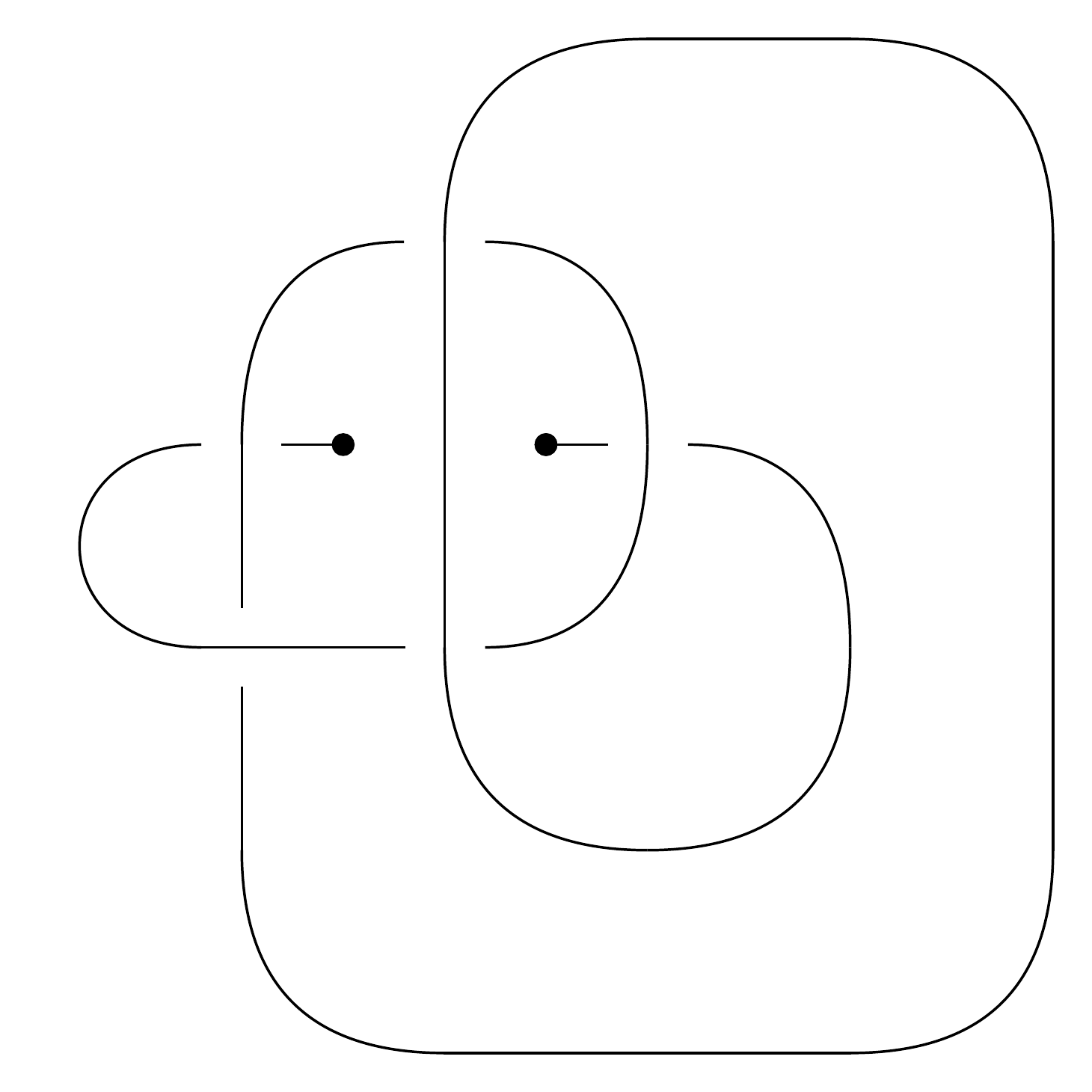}\\
\textcolor{black}{$5_{799}$}
\vspace{1cm}
\end{minipage}
\begin{minipage}[t]{.25\linewidth}
\centering
\includegraphics[width=0.9\textwidth,height=3.5cm,keepaspectratio]{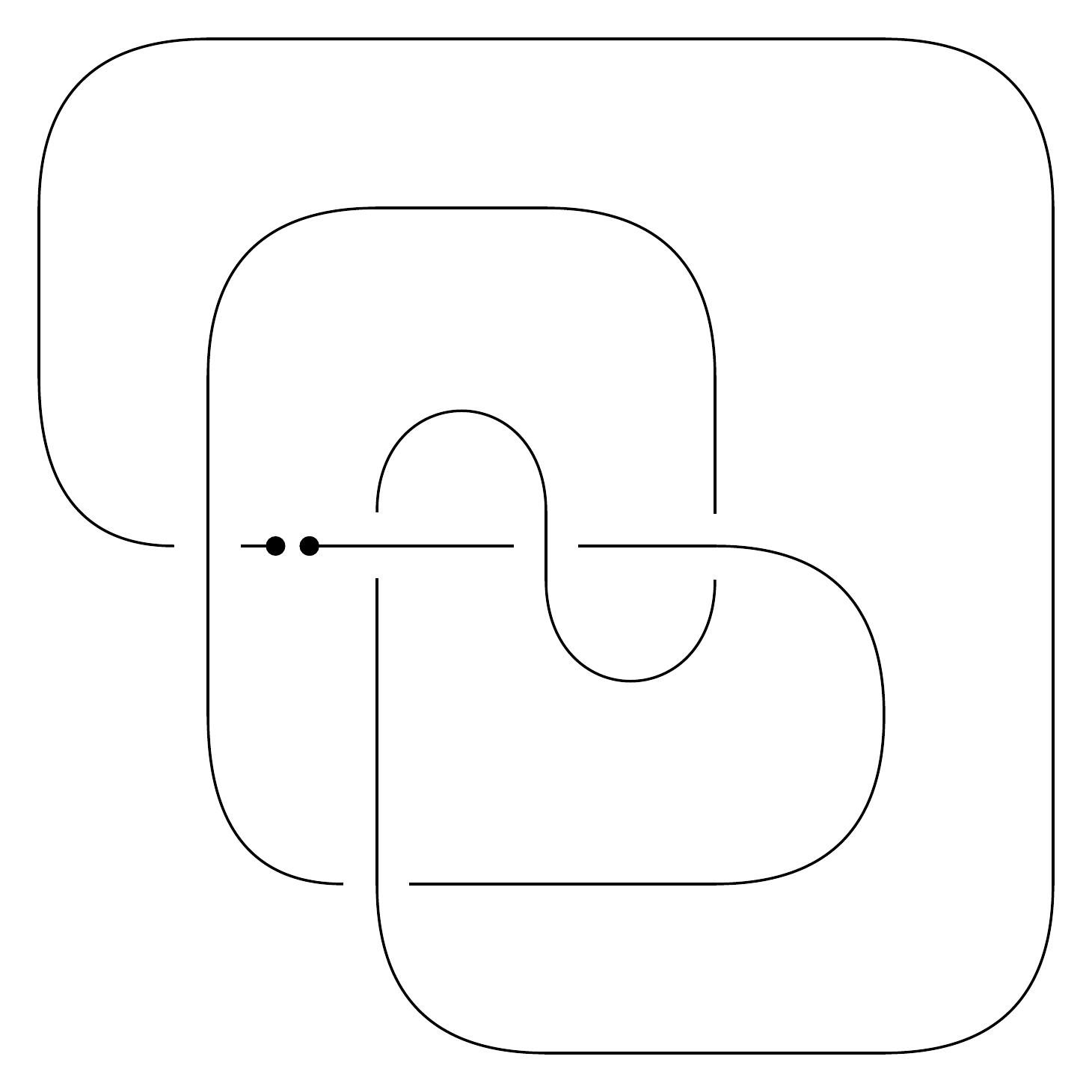}\\
\textcolor{black}{$5_{800}$}
\vspace{1cm}
\end{minipage}
\begin{minipage}[t]{.25\linewidth}
\centering
\includegraphics[width=0.9\textwidth,height=3.5cm,keepaspectratio]{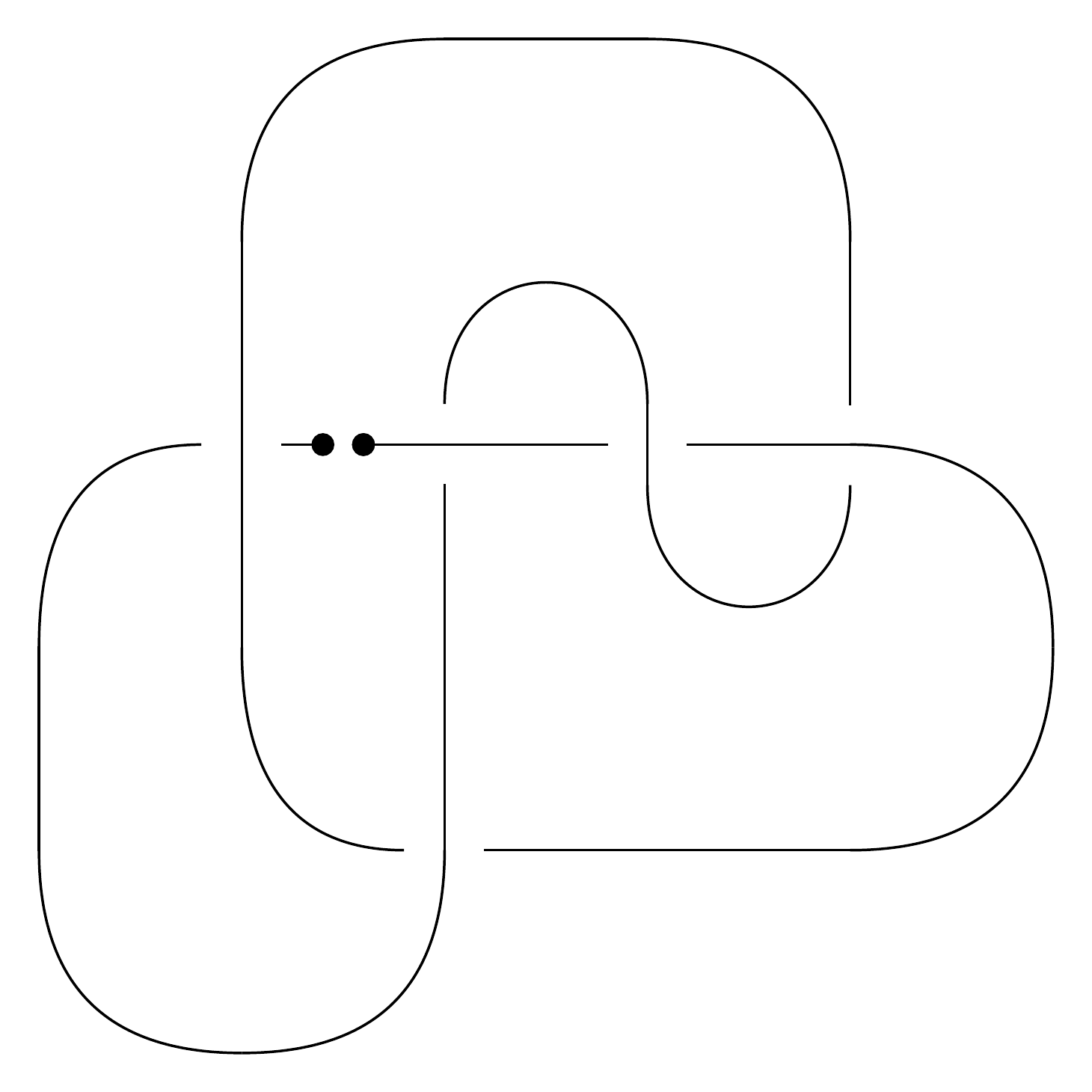}\\
\textcolor{black}{$5_{801}$}
\vspace{1cm}
\end{minipage}
\begin{minipage}[t]{.25\linewidth}
\centering
\includegraphics[width=0.9\textwidth,height=3.5cm,keepaspectratio]{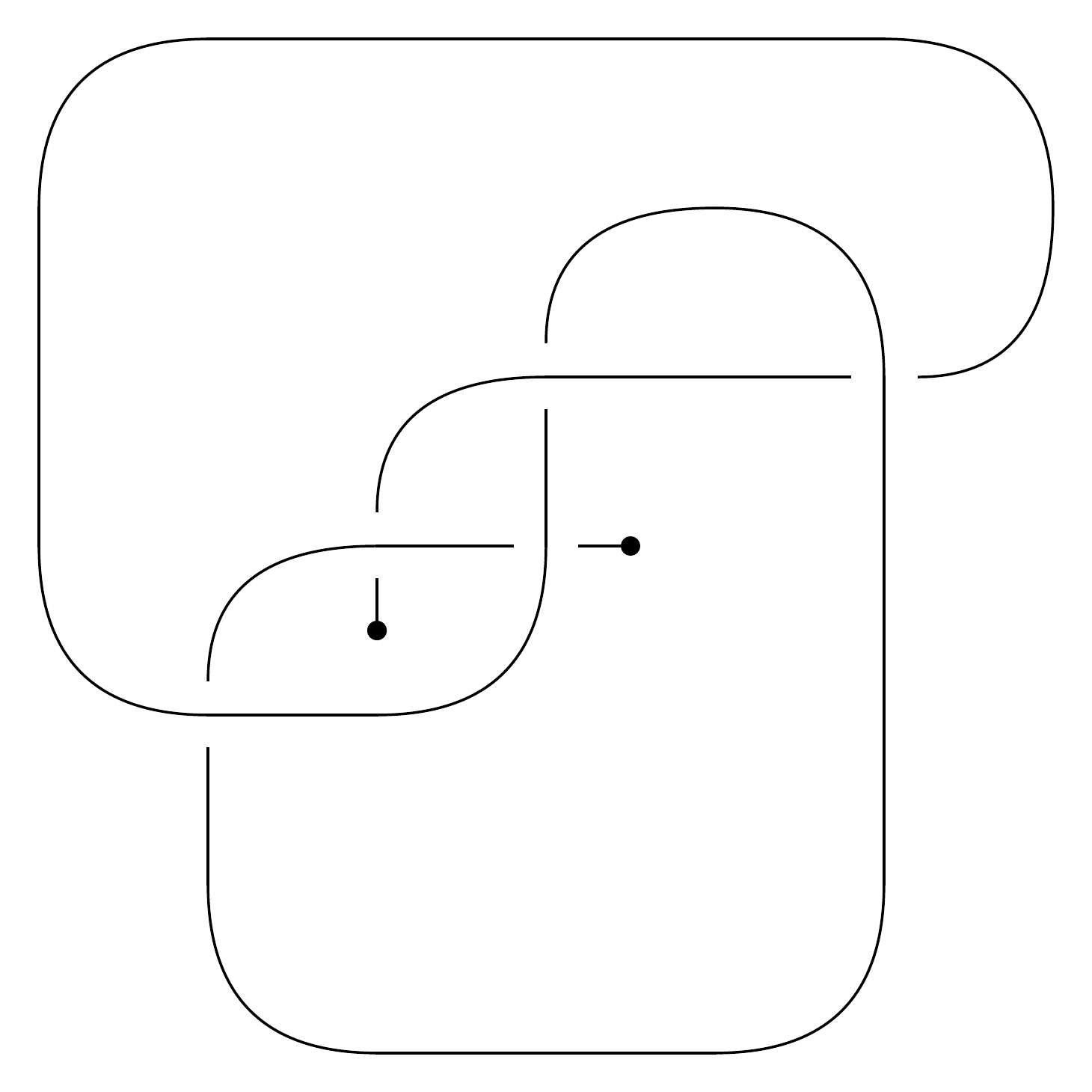}\\
\textcolor{black}{$5_{802}$}
\vspace{1cm}
\end{minipage}
\begin{minipage}[t]{.25\linewidth}
\centering
\includegraphics[width=0.9\textwidth,height=3.5cm,keepaspectratio]{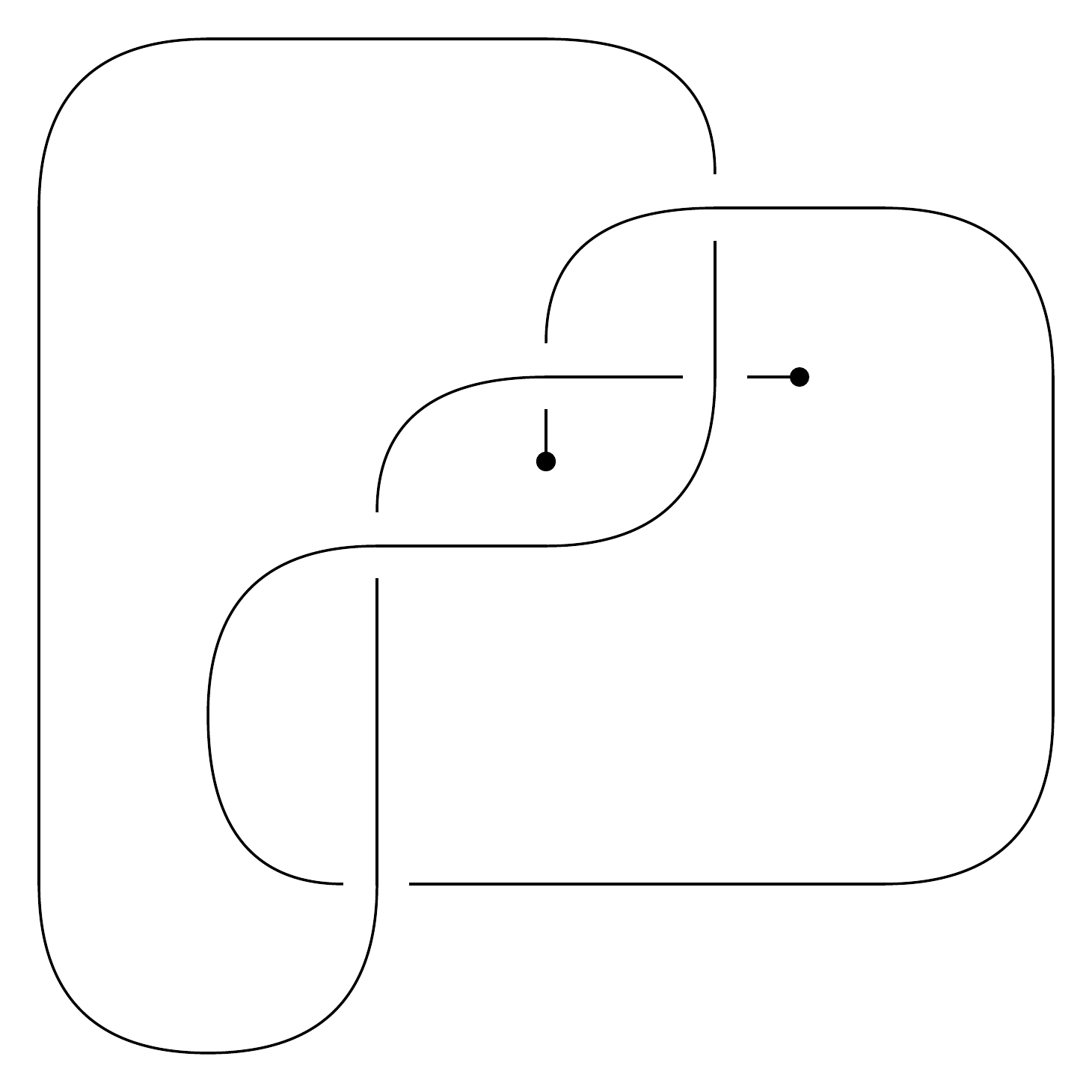}\\
\textcolor{black}{$5_{803}$}
\vspace{1cm}
\end{minipage}
\begin{minipage}[t]{.25\linewidth}
\centering
\includegraphics[width=0.9\textwidth,height=3.5cm,keepaspectratio]{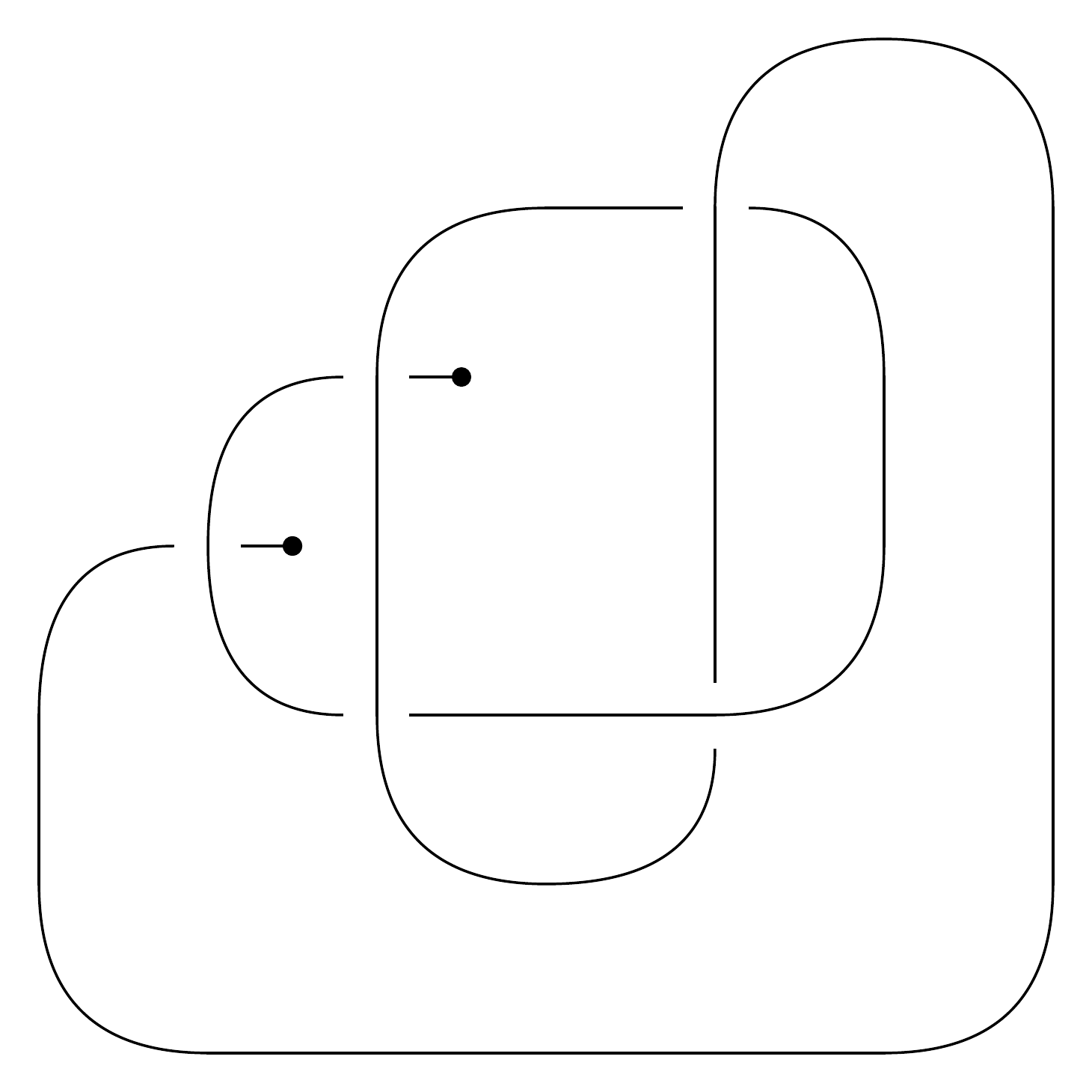}\\
\textcolor{black}{$5_{804}$}
\vspace{1cm}
\end{minipage}
\begin{minipage}[t]{.25\linewidth}
\centering
\includegraphics[width=0.9\textwidth,height=3.5cm,keepaspectratio]{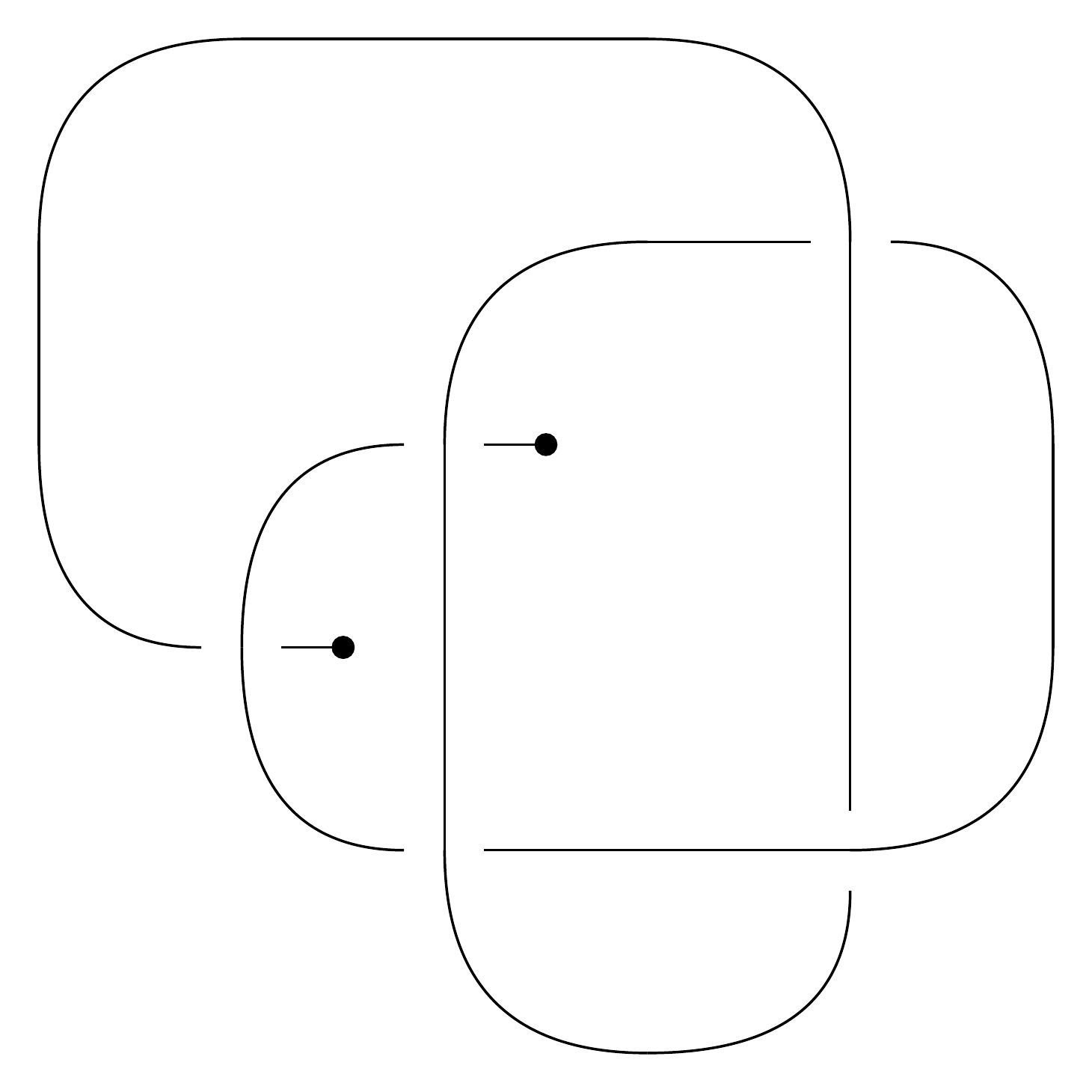}\\
\textcolor{black}{$5_{805}$}
\vspace{1cm}
\end{minipage}
\begin{minipage}[t]{.25\linewidth}
\centering
\includegraphics[width=0.9\textwidth,height=3.5cm,keepaspectratio]{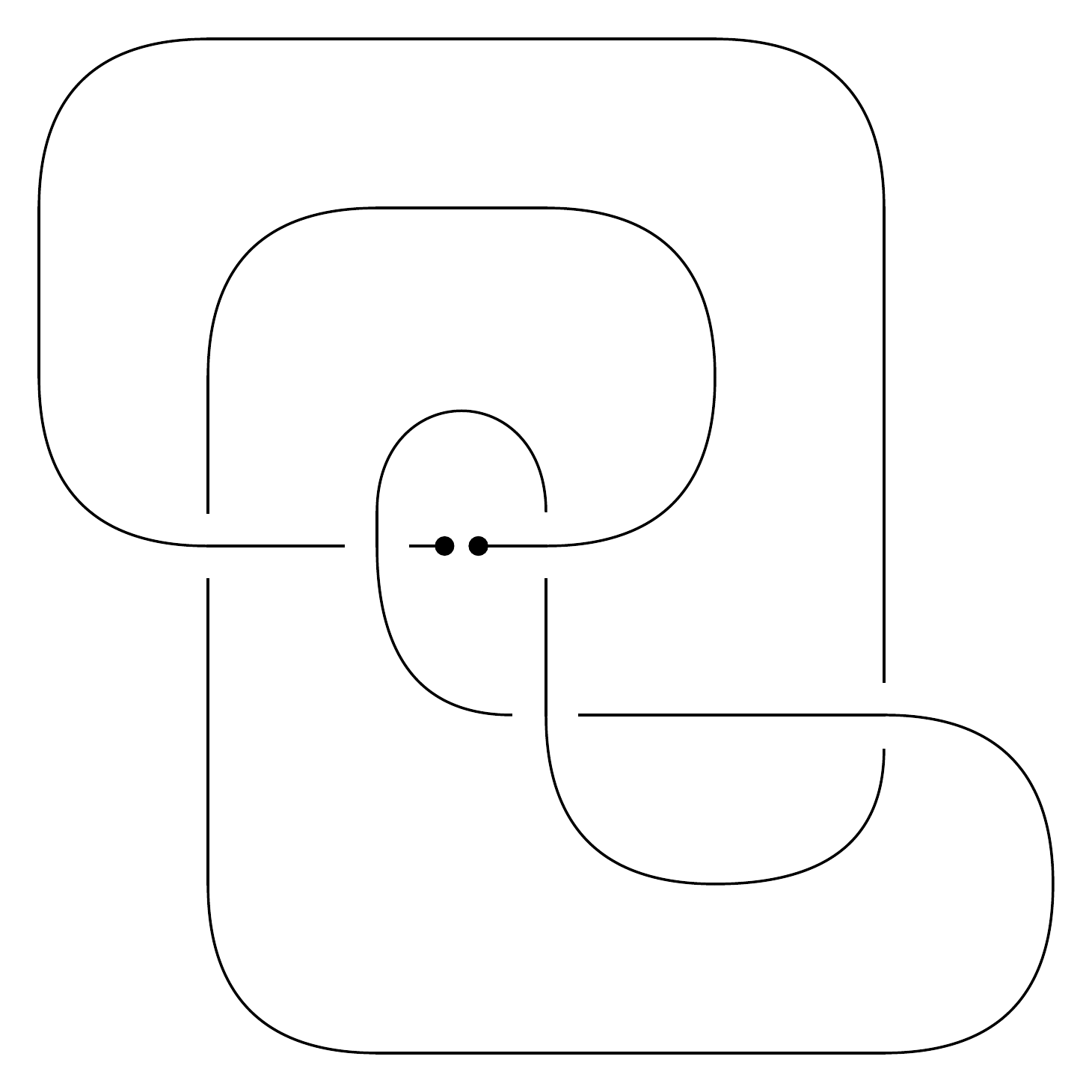}\\
\textcolor{black}{$5_{806}$}
\vspace{1cm}
\end{minipage}
\begin{minipage}[t]{.25\linewidth}
\centering
\includegraphics[width=0.9\textwidth,height=3.5cm,keepaspectratio]{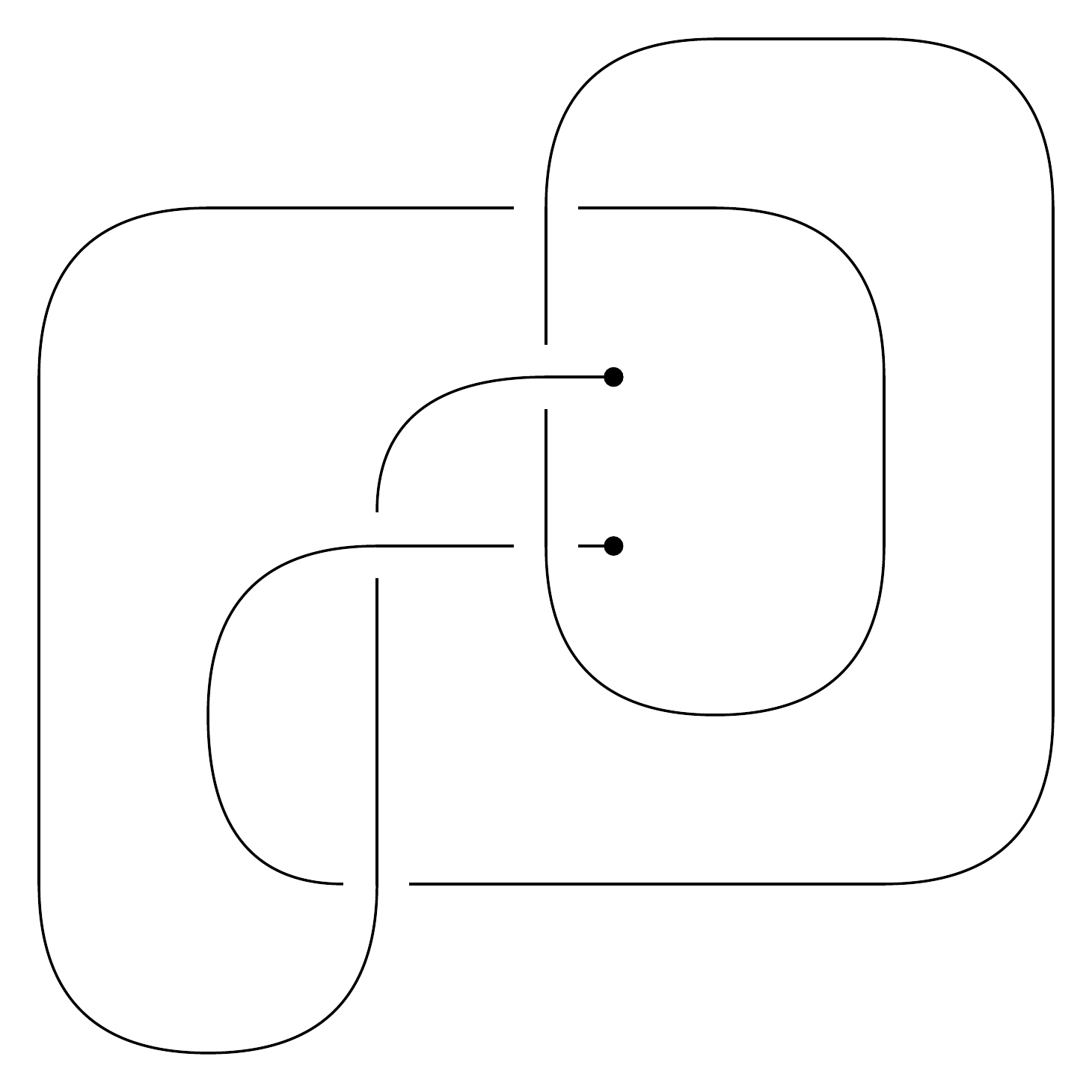}\\
\textcolor{black}{$5_{807}$}
\vspace{1cm}
\end{minipage}
\begin{minipage}[t]{.25\linewidth}
\centering
\includegraphics[width=0.9\textwidth,height=3.5cm,keepaspectratio]{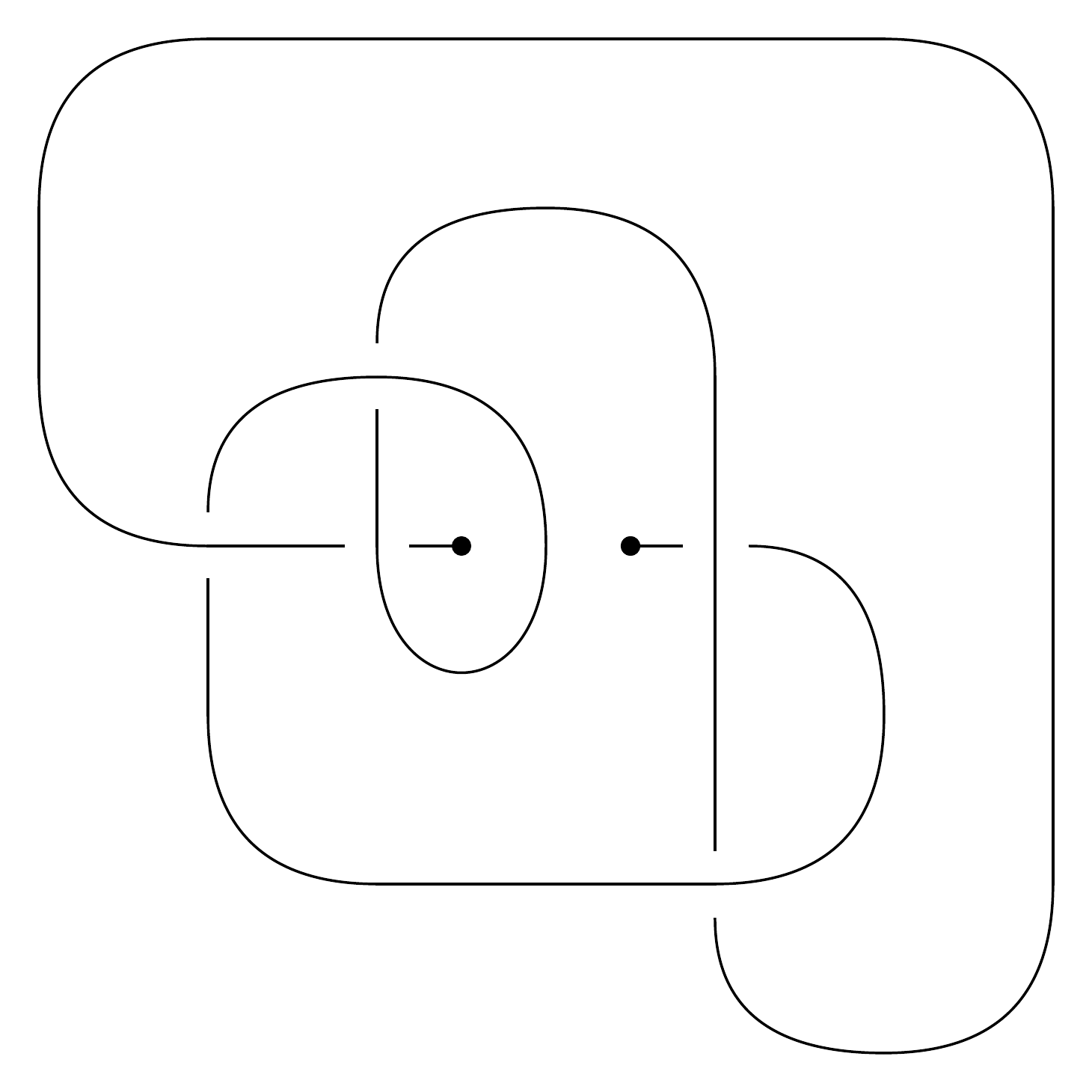}\\
\textcolor{black}{$5_{808}$}
\vspace{1cm}
\end{minipage}
\begin{minipage}[t]{.25\linewidth}
\centering
\includegraphics[width=0.9\textwidth,height=3.5cm,keepaspectratio]{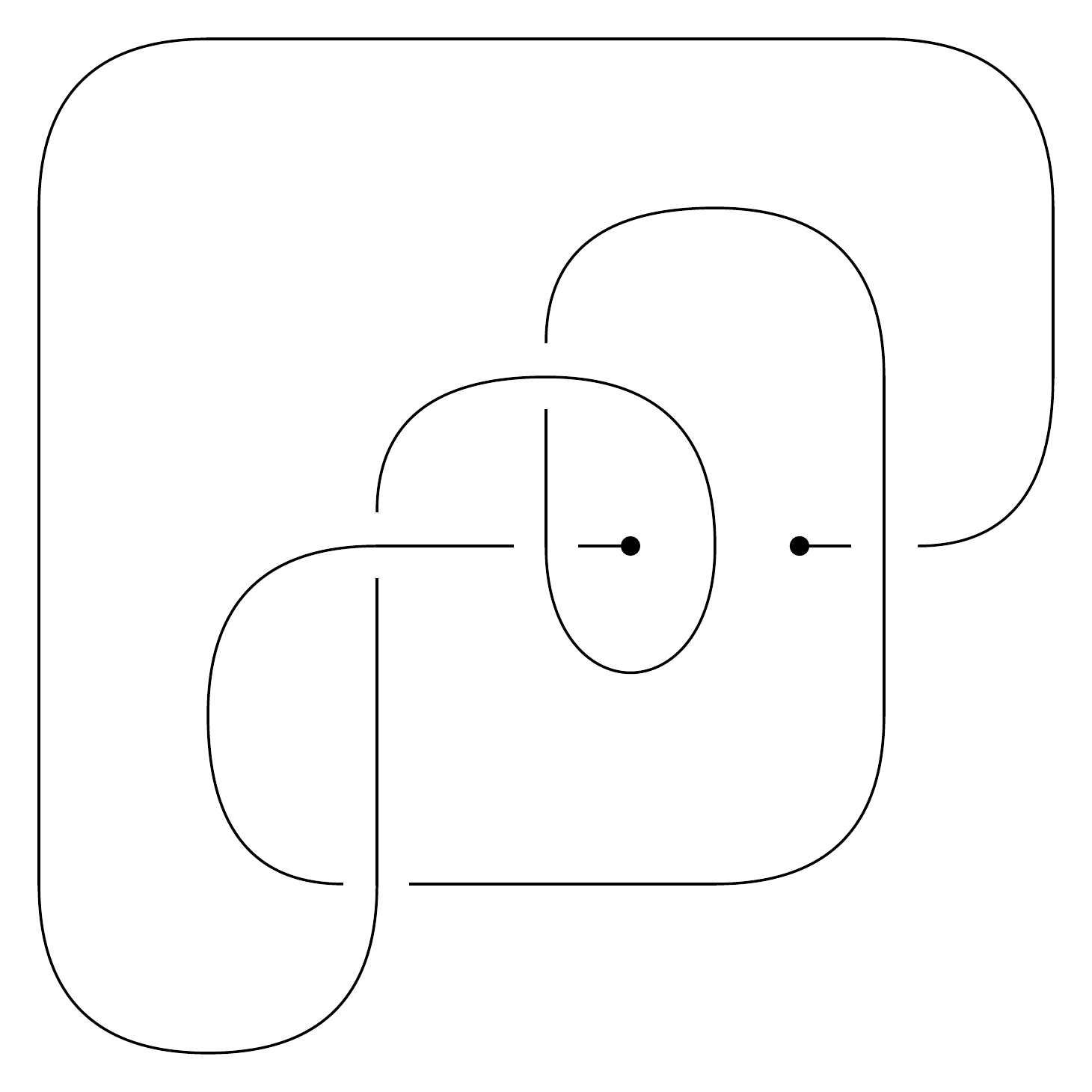}\\
\textcolor{black}{$5_{809}$}
\vspace{1cm}
\end{minipage}
\begin{minipage}[t]{.25\linewidth}
\centering
\includegraphics[width=0.9\textwidth,height=3.5cm,keepaspectratio]{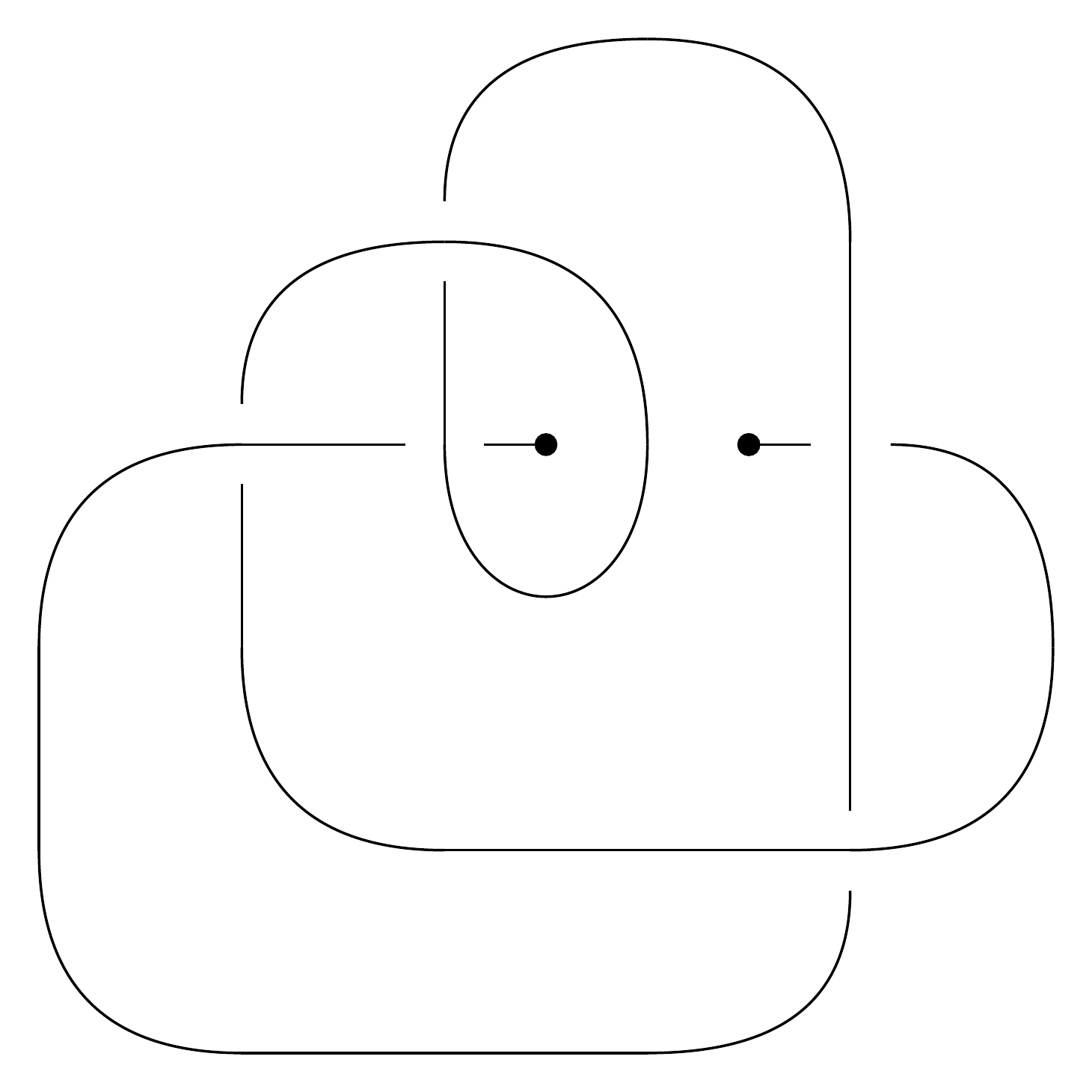}\\
\textcolor{black}{$5_{810}$}
\vspace{1cm}
\end{minipage}
\begin{minipage}[t]{.25\linewidth}
\centering
\includegraphics[width=0.9\textwidth,height=3.5cm,keepaspectratio]{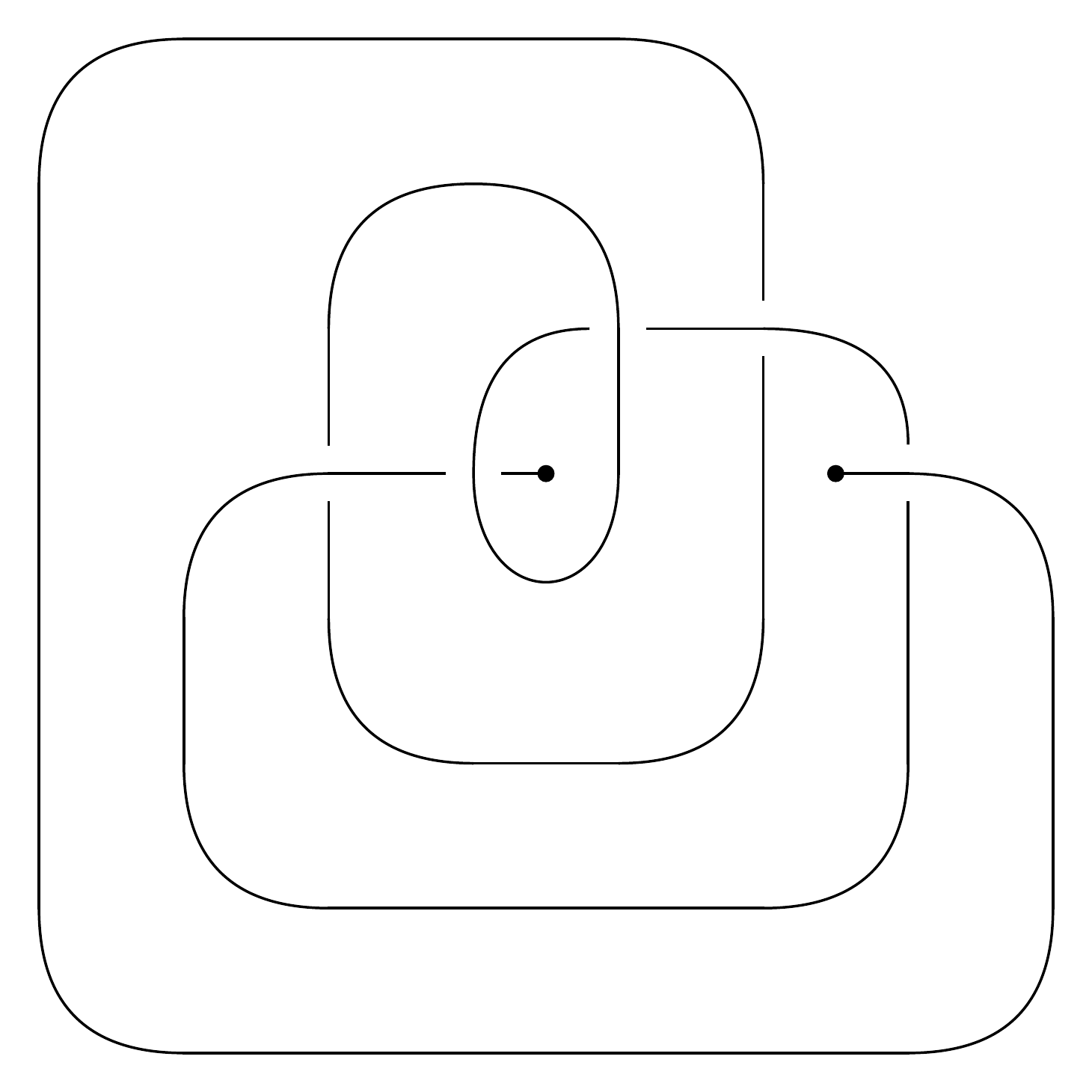}\\
\textcolor{black}{$5_{811}$}
\vspace{1cm}
\end{minipage}
\begin{minipage}[t]{.25\linewidth}
\centering
\includegraphics[width=0.9\textwidth,height=3.5cm,keepaspectratio]{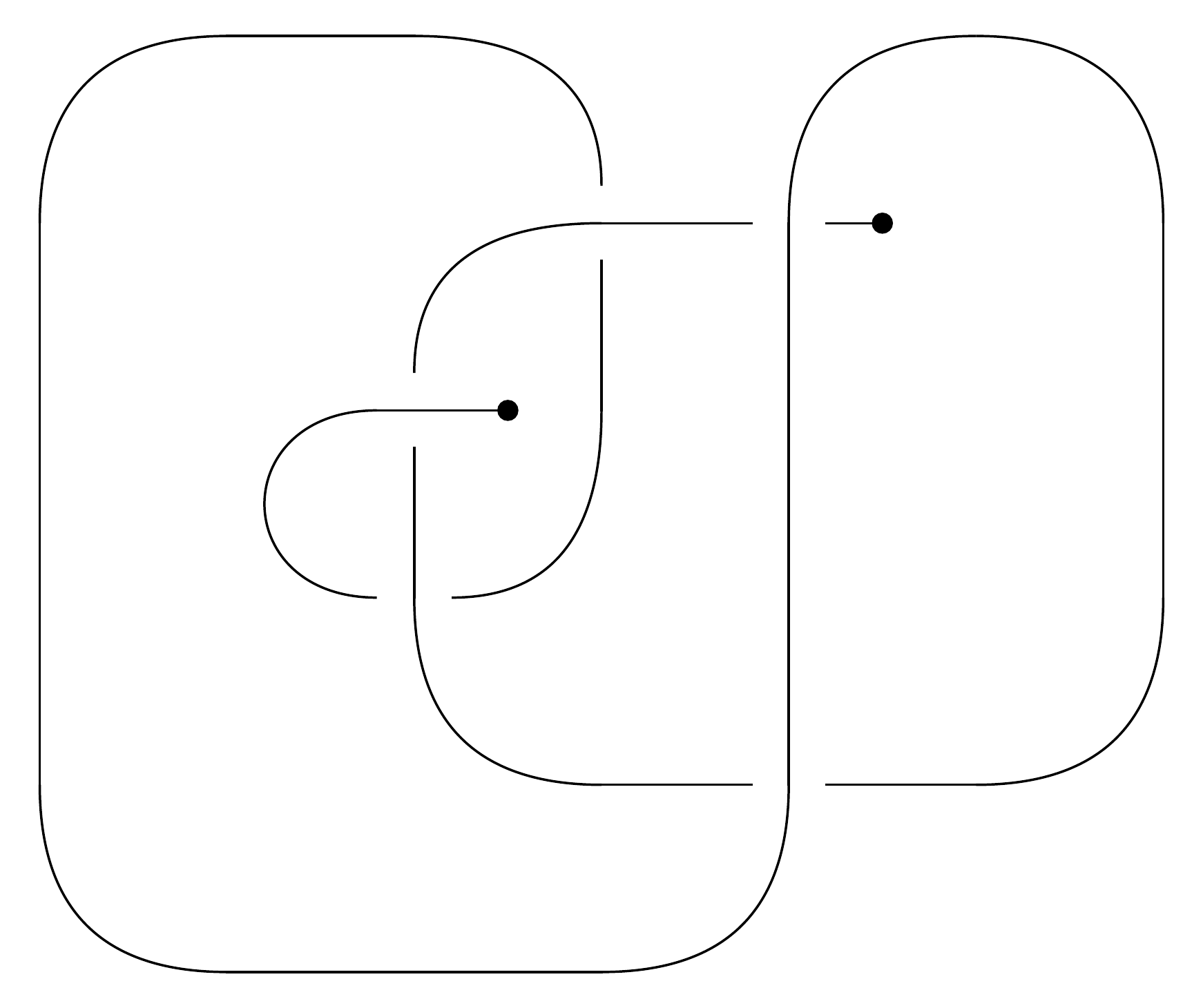}\\
\textcolor{black}{$5_{812}$}
\vspace{1cm}
\end{minipage}
\begin{minipage}[t]{.25\linewidth}
\centering
\includegraphics[width=0.9\textwidth,height=3.5cm,keepaspectratio]{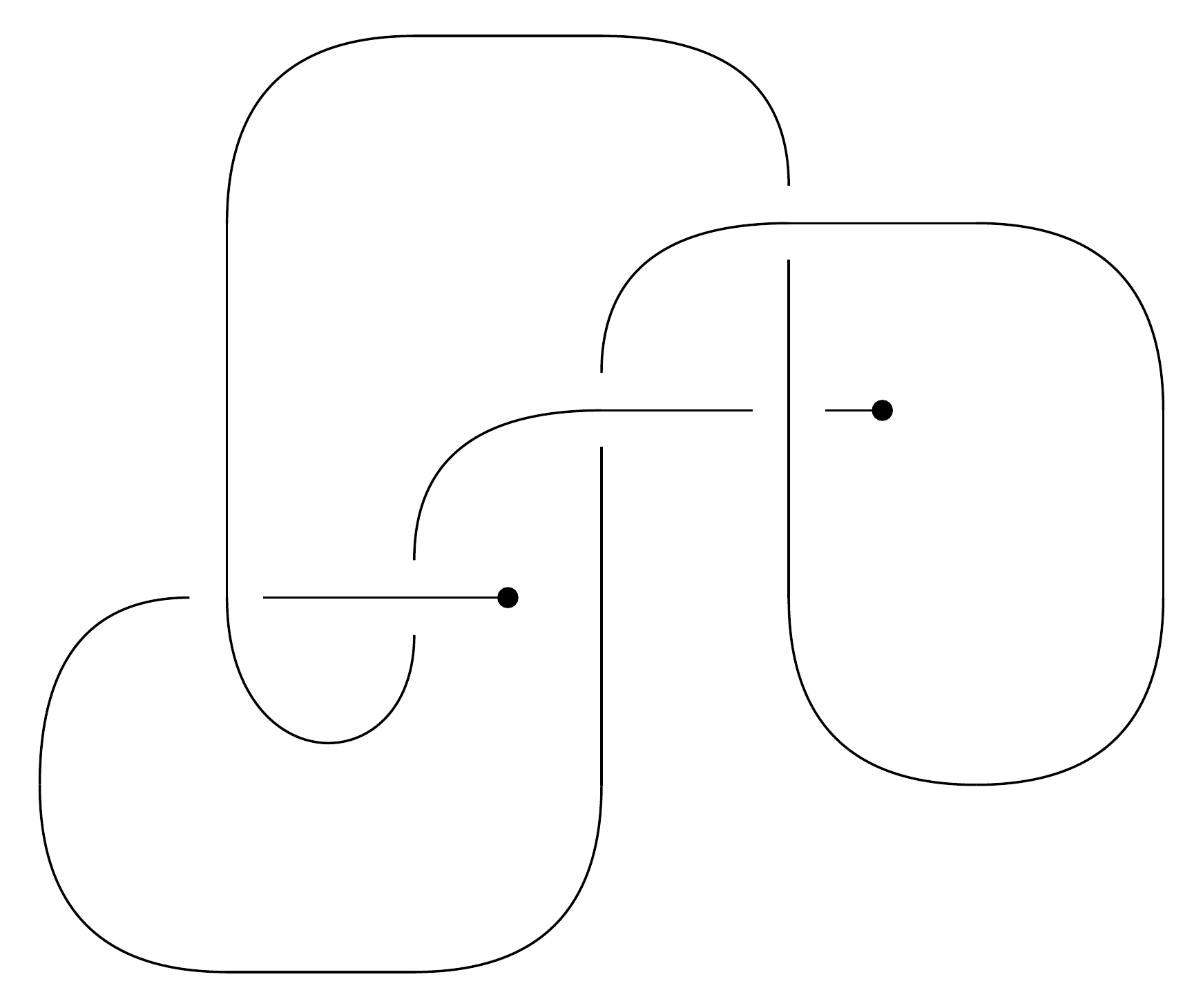}\\
\textcolor{black}{$5_{813}$}
\vspace{1cm}
\end{minipage}
\begin{minipage}[t]{.25\linewidth}
\centering
\includegraphics[width=0.9\textwidth,height=3.5cm,keepaspectratio]{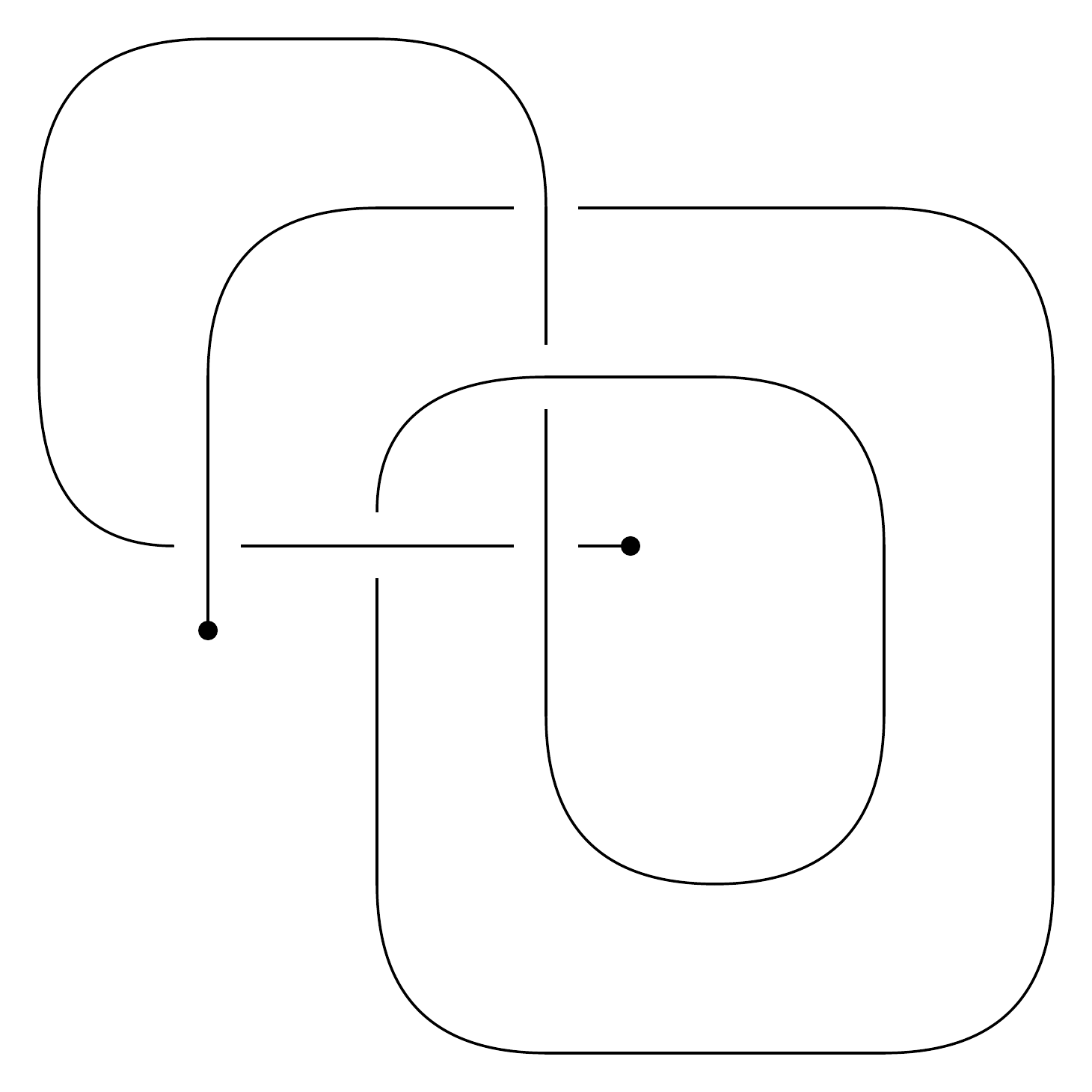}\\
\textcolor{black}{$5_{814}$}
\vspace{1cm}
\end{minipage}
\begin{minipage}[t]{.25\linewidth}
\centering
\includegraphics[width=0.9\textwidth,height=3.5cm,keepaspectratio]{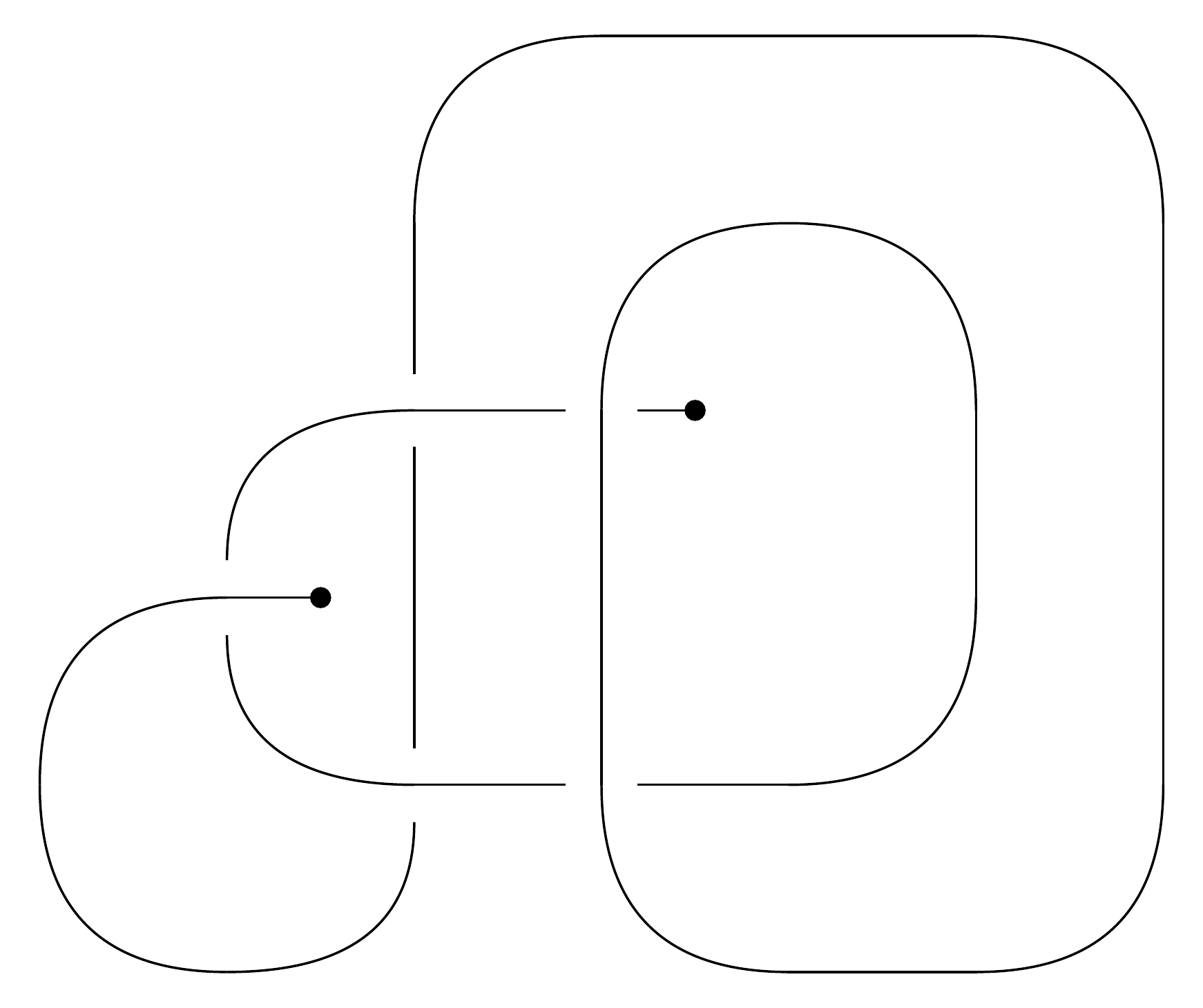}\\
\textcolor{black}{$5_{815}$}
\vspace{1cm}
\end{minipage}
\begin{minipage}[t]{.25\linewidth}
\centering
\includegraphics[width=0.9\textwidth,height=3.5cm,keepaspectratio]{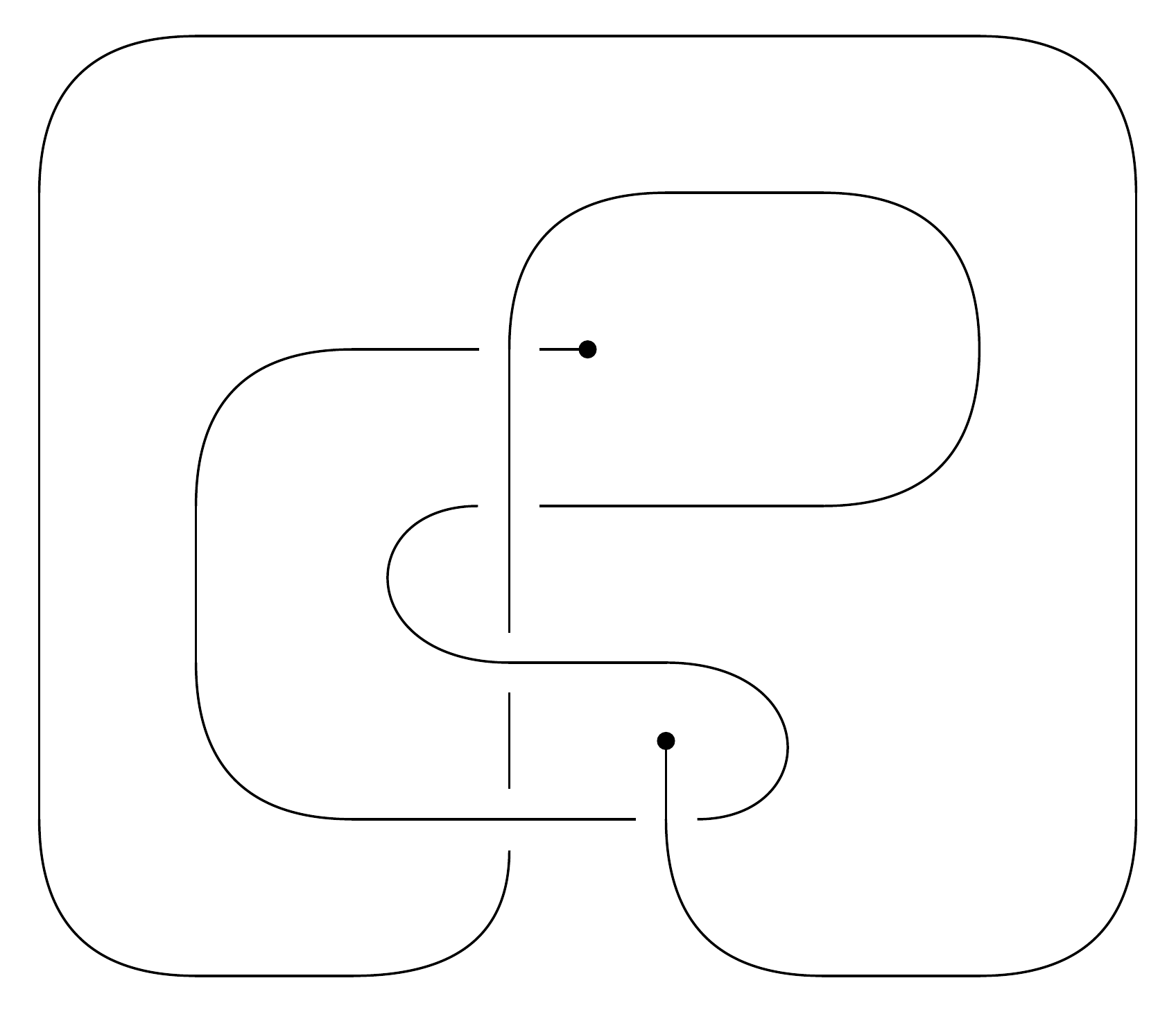}\\
\textcolor{black}{$5_{816}$}
\vspace{1cm}
\end{minipage}
\begin{minipage}[t]{.25\linewidth}
\centering
\includegraphics[width=0.9\textwidth,height=3.5cm,keepaspectratio]{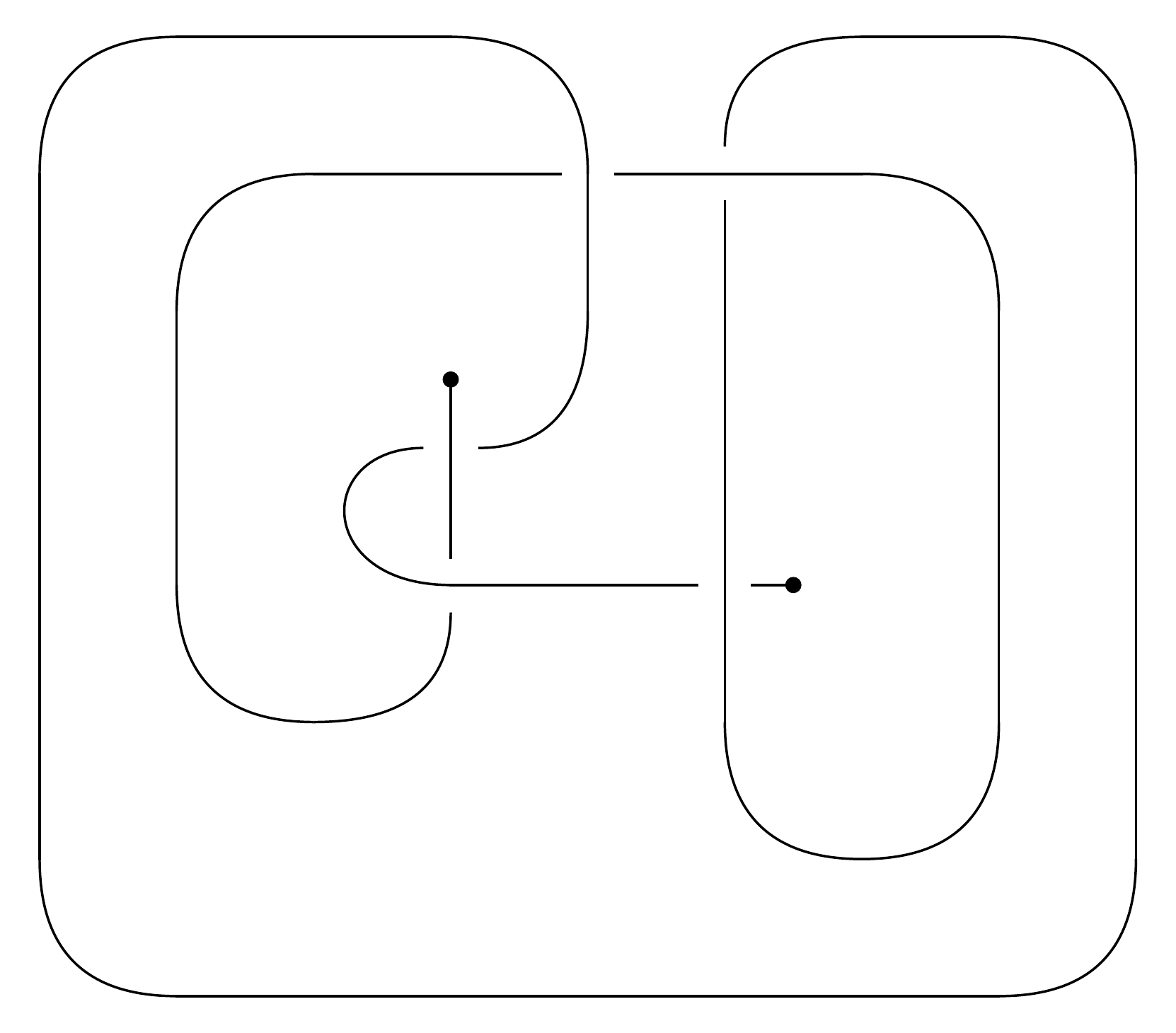}\\
\textcolor{black}{$5_{817}$}
\vspace{1cm}
\end{minipage}
\begin{minipage}[t]{.25\linewidth}
\centering
\includegraphics[width=0.9\textwidth,height=3.5cm,keepaspectratio]{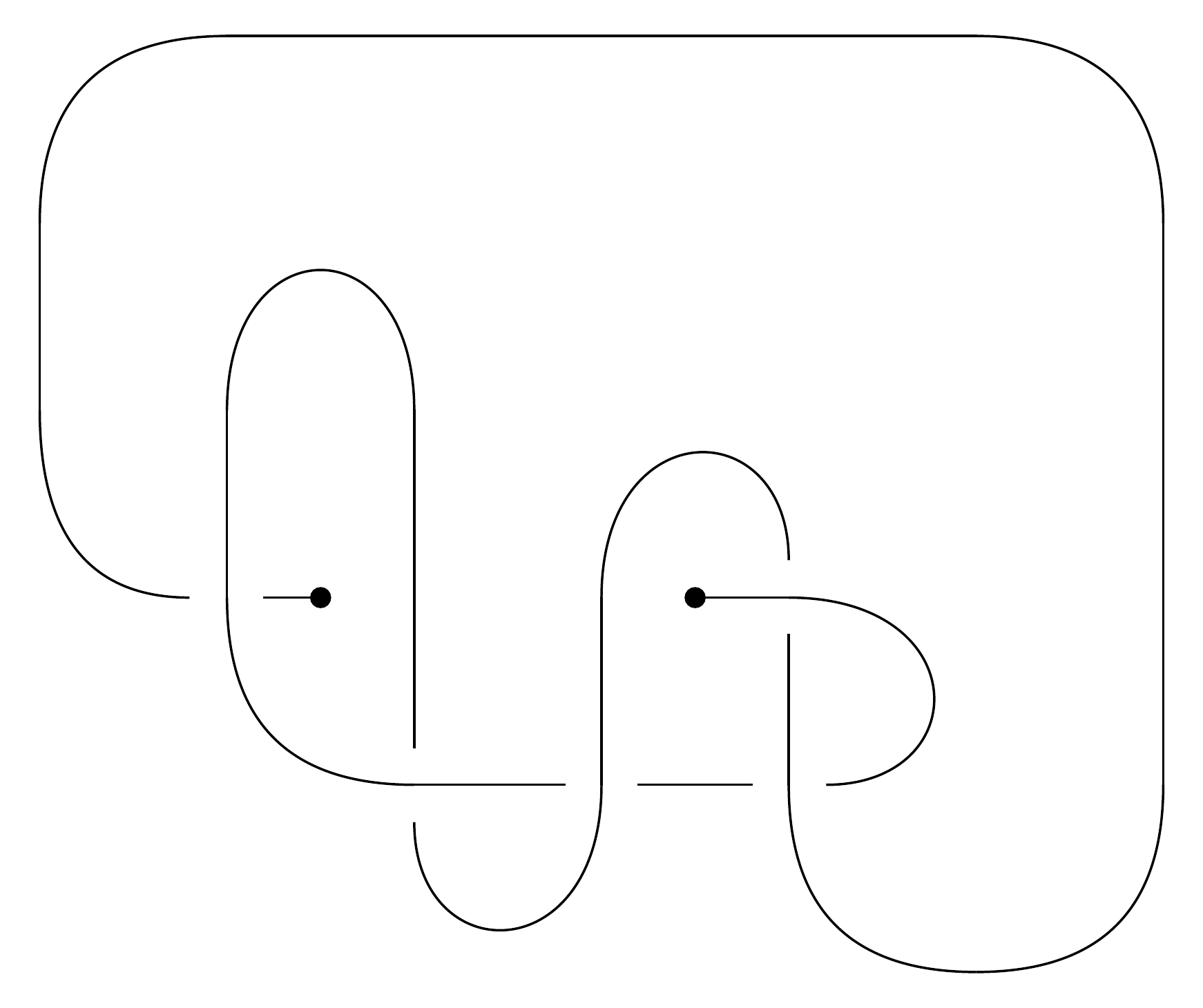}\\
\textcolor{black}{$5_{818}$}
\vspace{1cm}
\end{minipage}
\begin{minipage}[t]{.25\linewidth}
\centering
\includegraphics[width=0.9\textwidth,height=3.5cm,keepaspectratio]{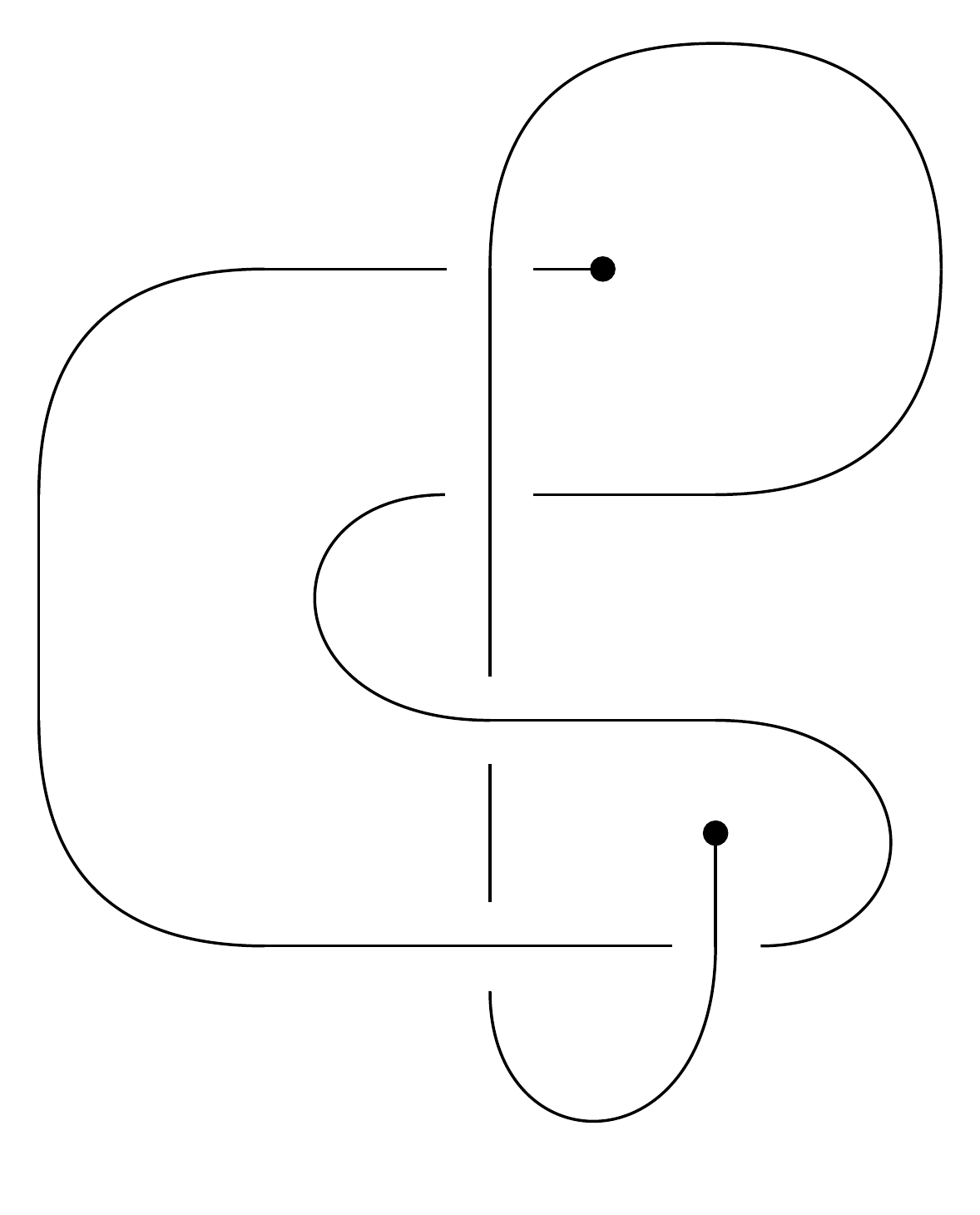}\\
\textcolor{black}{$5_{819}$}
\vspace{1cm}
\end{minipage}
\begin{minipage}[t]{.25\linewidth}
\centering
\includegraphics[width=0.9\textwidth,height=3.5cm,keepaspectratio]{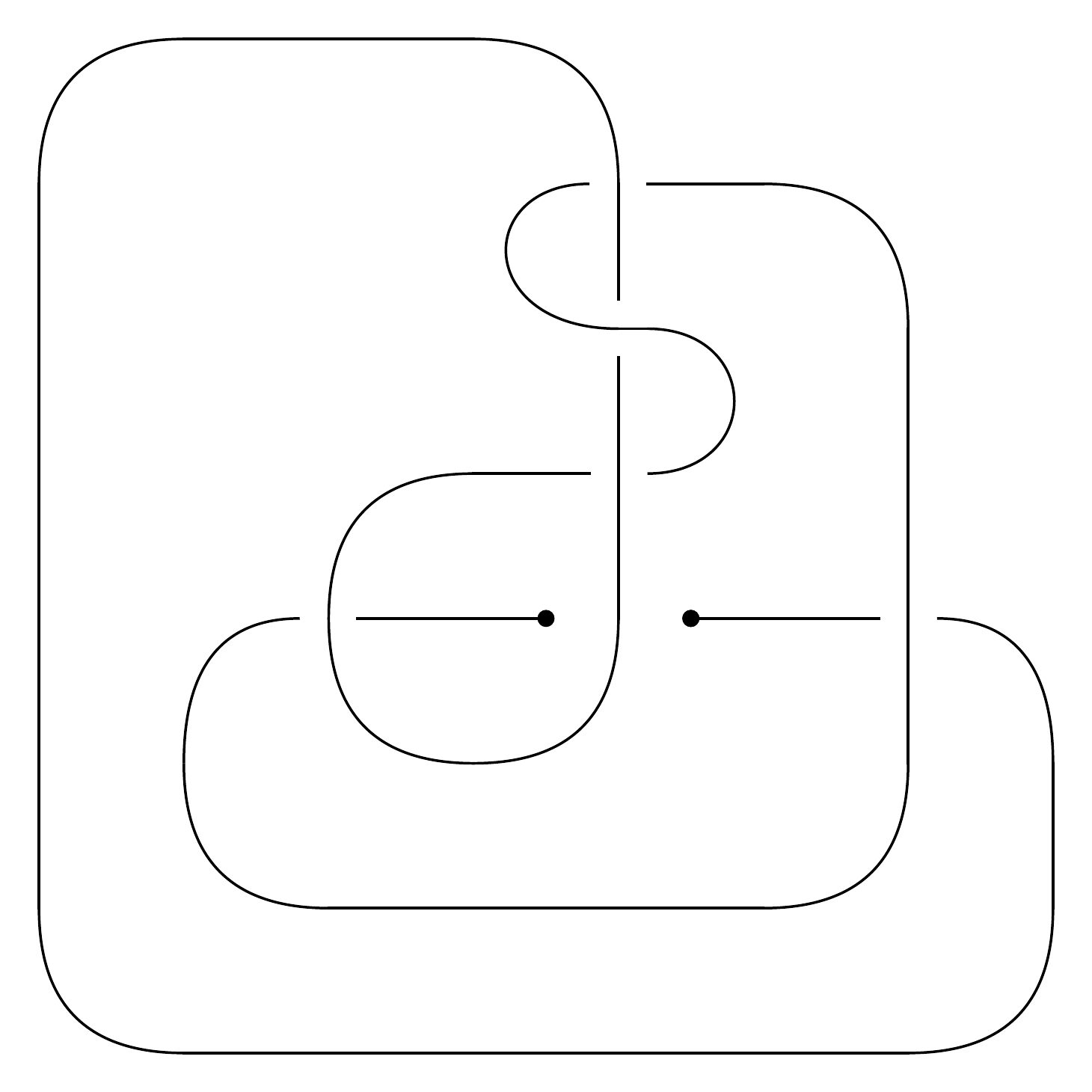}\\
\textcolor{black}{$5_{820}$}
\vspace{1cm}
\end{minipage}
\begin{minipage}[t]{.25\linewidth}
\centering
\includegraphics[width=0.9\textwidth,height=3.5cm,keepaspectratio]{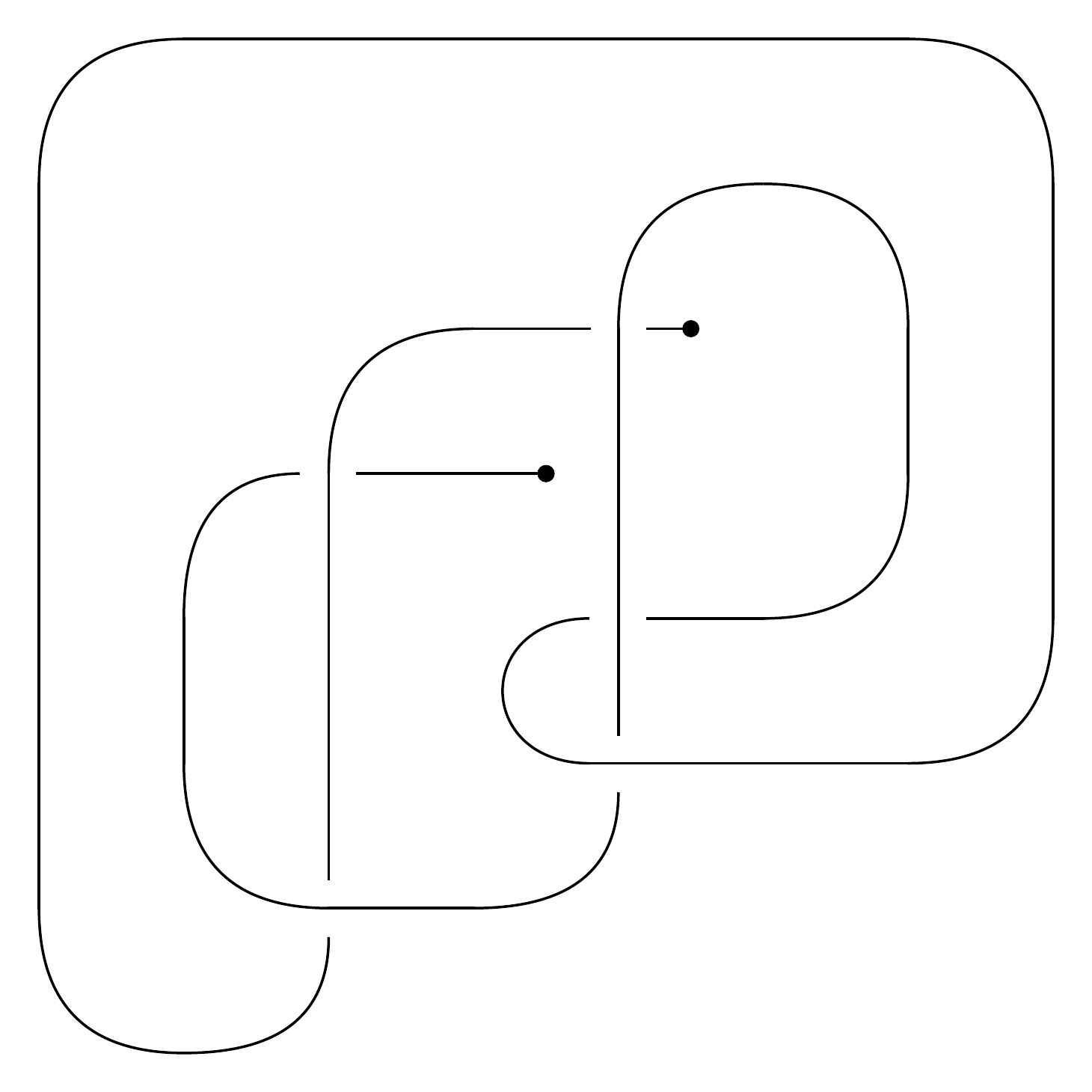}\\
\textcolor{black}{$5_{821}$}
\vspace{1cm}
\end{minipage}
\begin{minipage}[t]{.25\linewidth}
\centering
\includegraphics[width=0.9\textwidth,height=3.5cm,keepaspectratio]{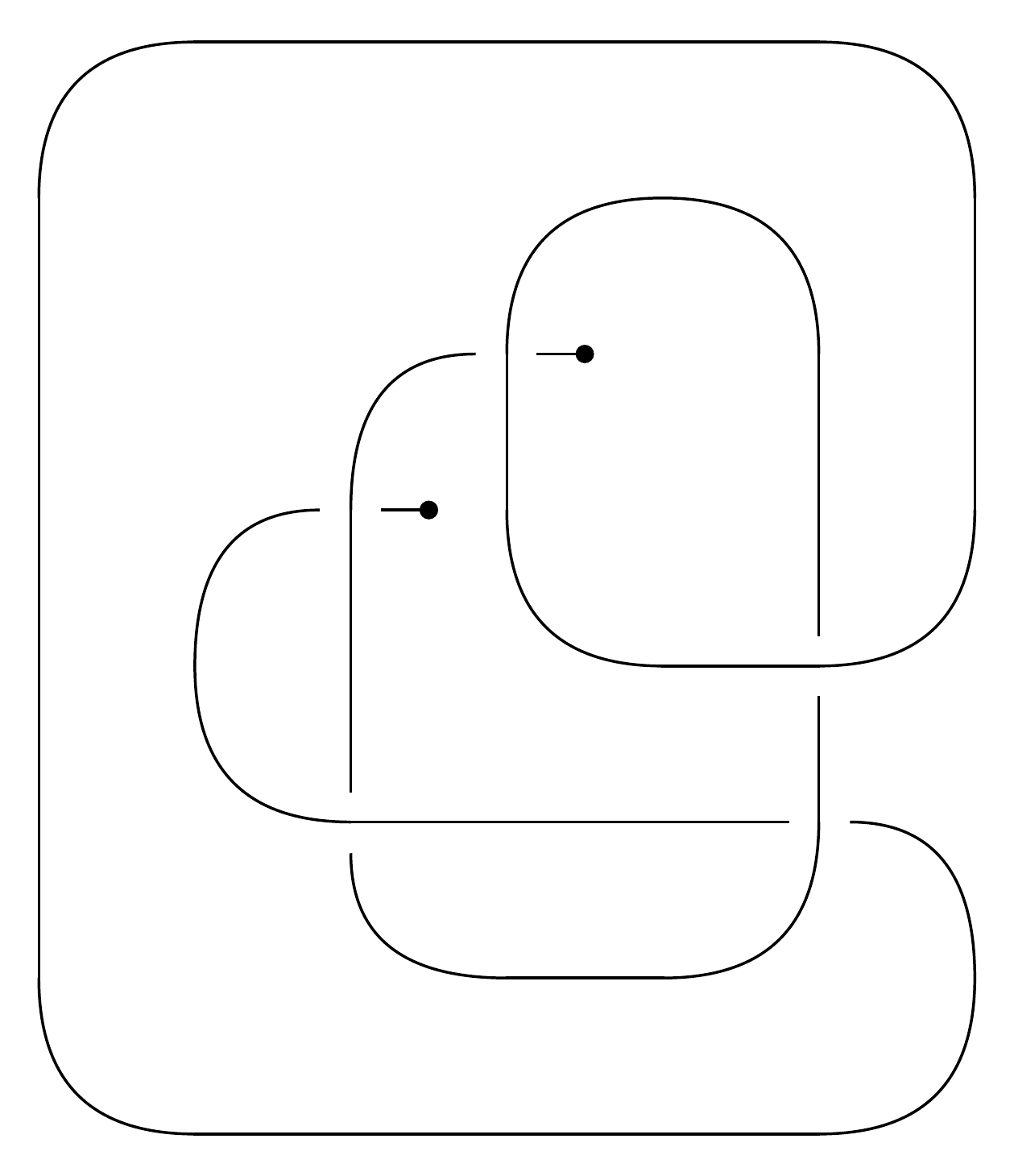}\\
\textcolor{black}{$5_{822}$}
\vspace{1cm}
\end{minipage}
\begin{minipage}[t]{.25\linewidth}
\centering
\includegraphics[width=0.9\textwidth,height=3.5cm,keepaspectratio]{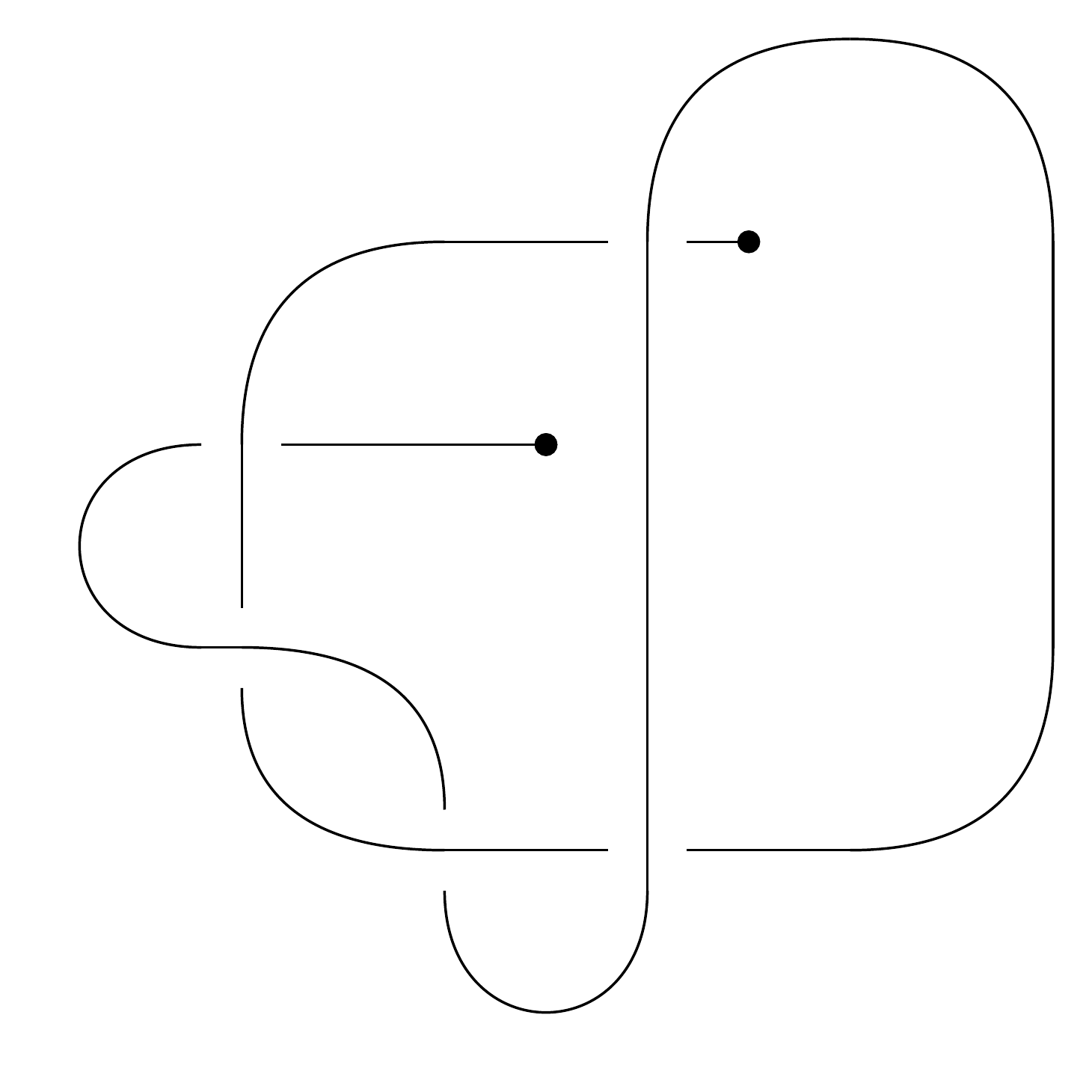}\\
\textcolor{black}{$5_{823}$}
\vspace{1cm}
\end{minipage}
\begin{minipage}[t]{.25\linewidth}
\centering
\includegraphics[width=0.9\textwidth,height=3.5cm,keepaspectratio]{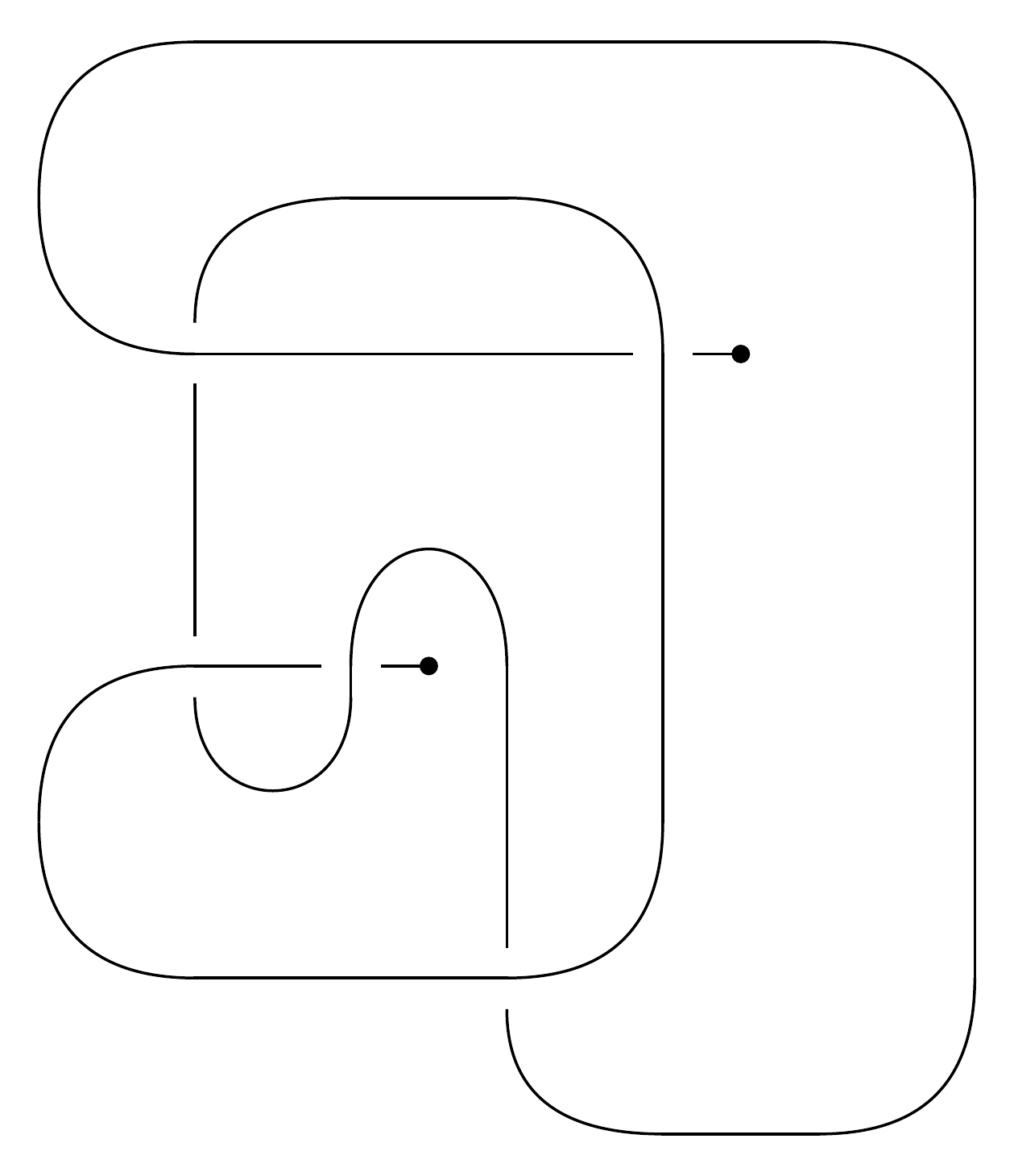}\\
\textcolor{black}{$5_{824}$}
\vspace{1cm}
\end{minipage}
\begin{minipage}[t]{.25\linewidth}
\centering
\includegraphics[width=0.9\textwidth,height=3.5cm,keepaspectratio]{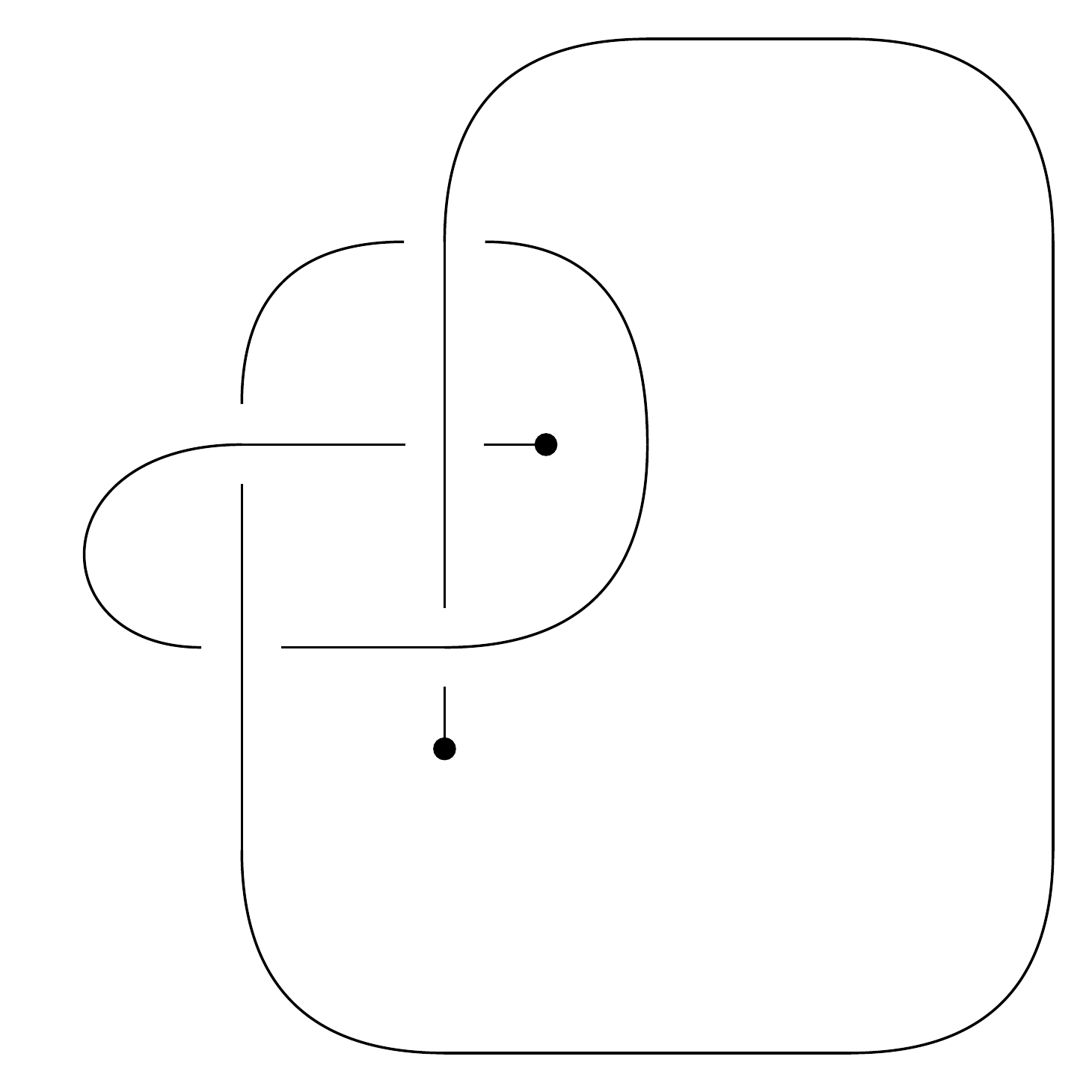}\\
\textcolor{black}{$5_{825}$}
\vspace{1cm}
\end{minipage}
\begin{minipage}[t]{.25\linewidth}
\centering
\includegraphics[width=0.9\textwidth,height=3.5cm,keepaspectratio]{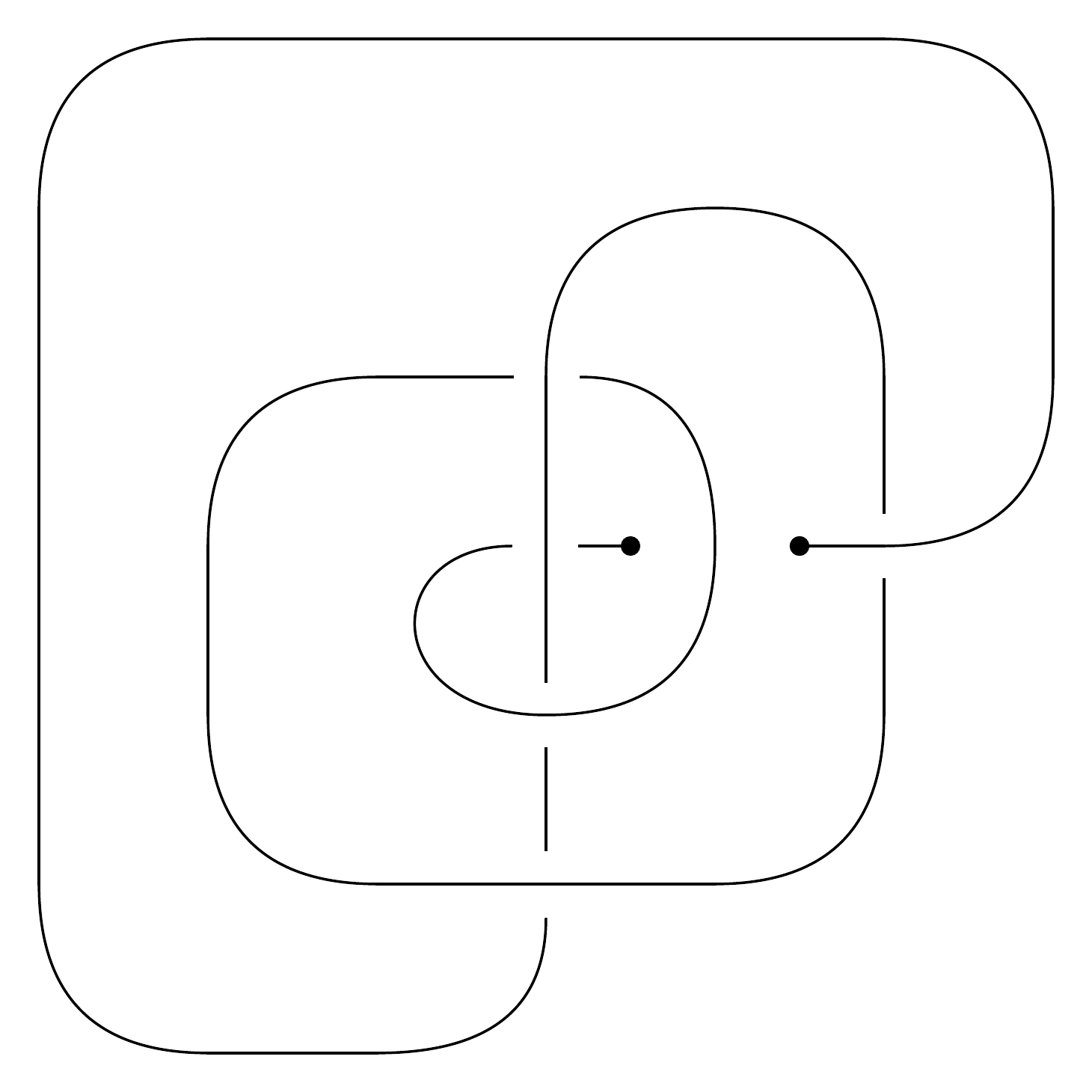}\\
\textcolor{black}{$5_{826}$}
\vspace{1cm}
\end{minipage}
\begin{minipage}[t]{.25\linewidth}
\centering
\includegraphics[width=0.9\textwidth,height=3.5cm,keepaspectratio]{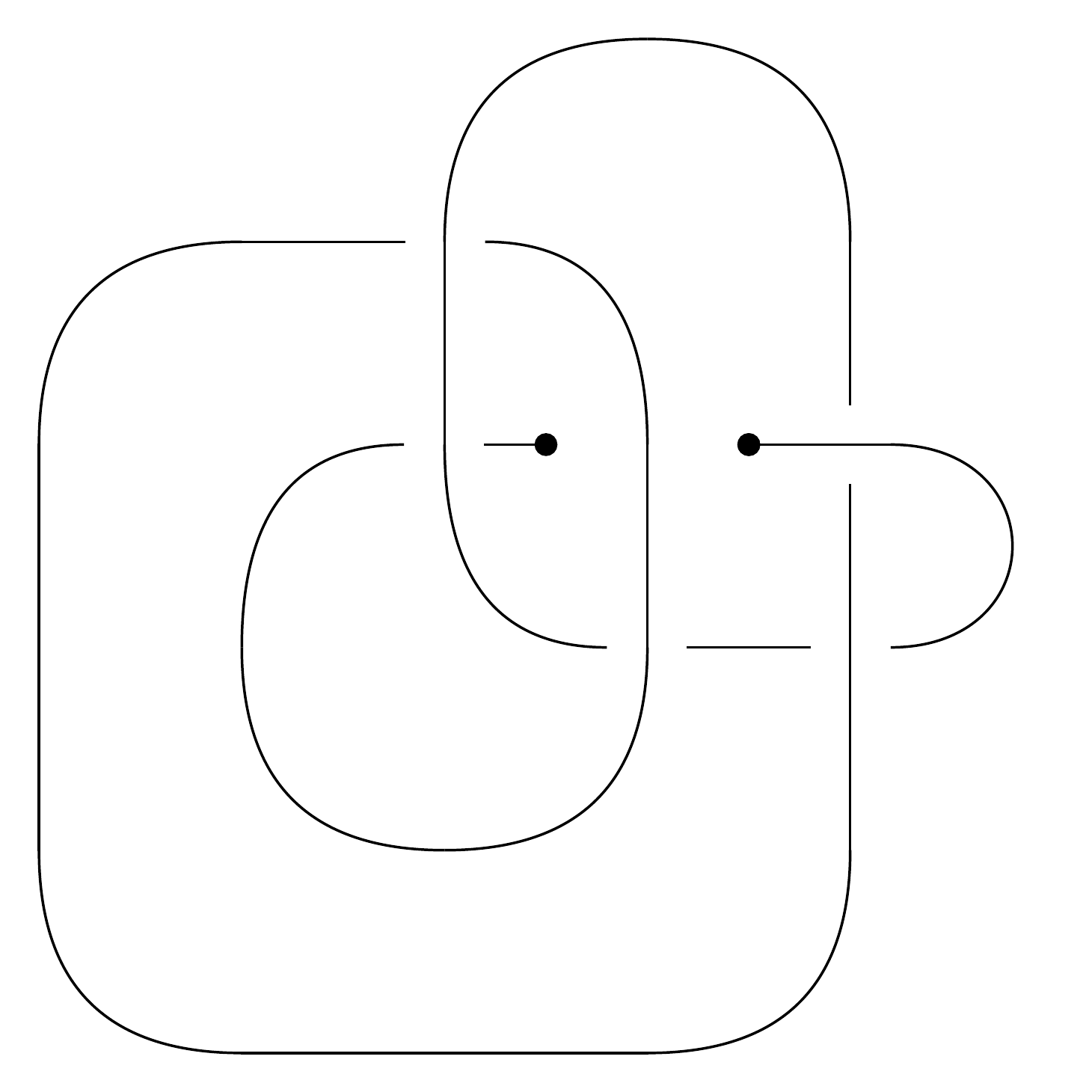}\\
\textcolor{black}{$5_{827}$}
\vspace{1cm}
\end{minipage}
\begin{minipage}[t]{.25\linewidth}
\centering
\includegraphics[width=0.9\textwidth,height=3.5cm,keepaspectratio]{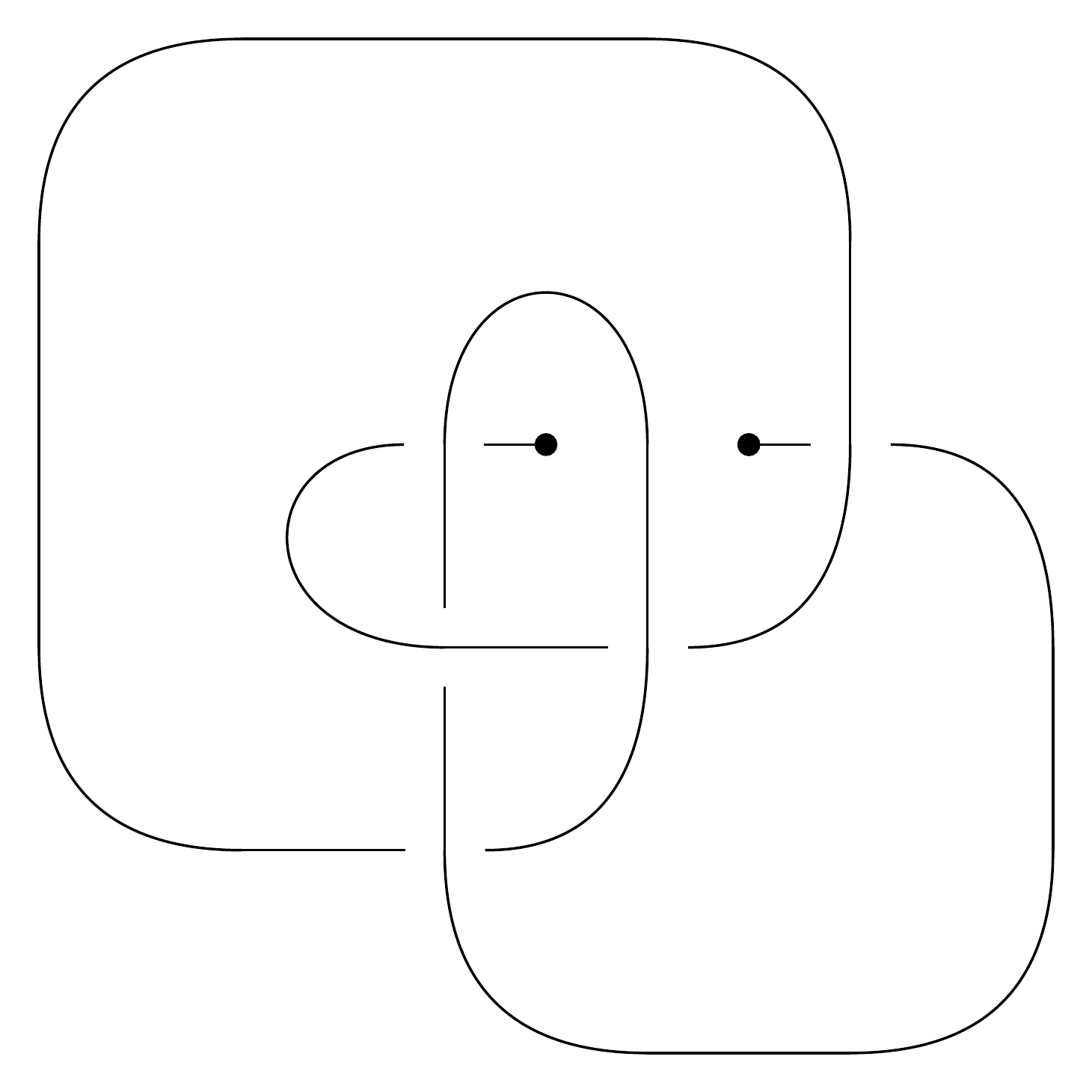}\\
\textcolor{black}{$5_{828}$}
\vspace{1cm}
\end{minipage}
\begin{minipage}[t]{.25\linewidth}
\centering
\includegraphics[width=0.9\textwidth,height=3.5cm,keepaspectratio]{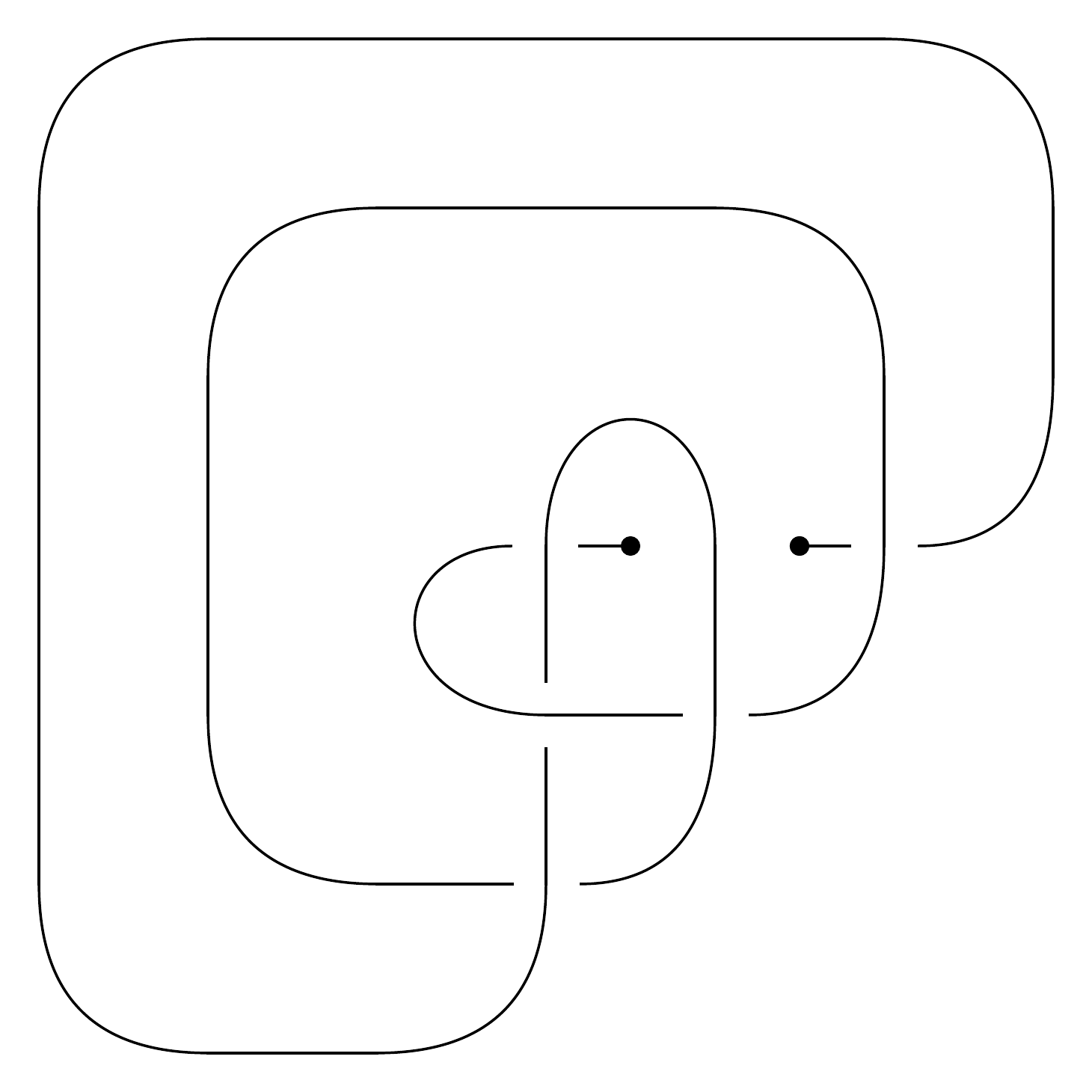}\\
\textcolor{black}{$5_{829}$}
\vspace{1cm}
\end{minipage}
\begin{minipage}[t]{.25\linewidth}
\centering
\includegraphics[width=0.9\textwidth,height=3.5cm,keepaspectratio]{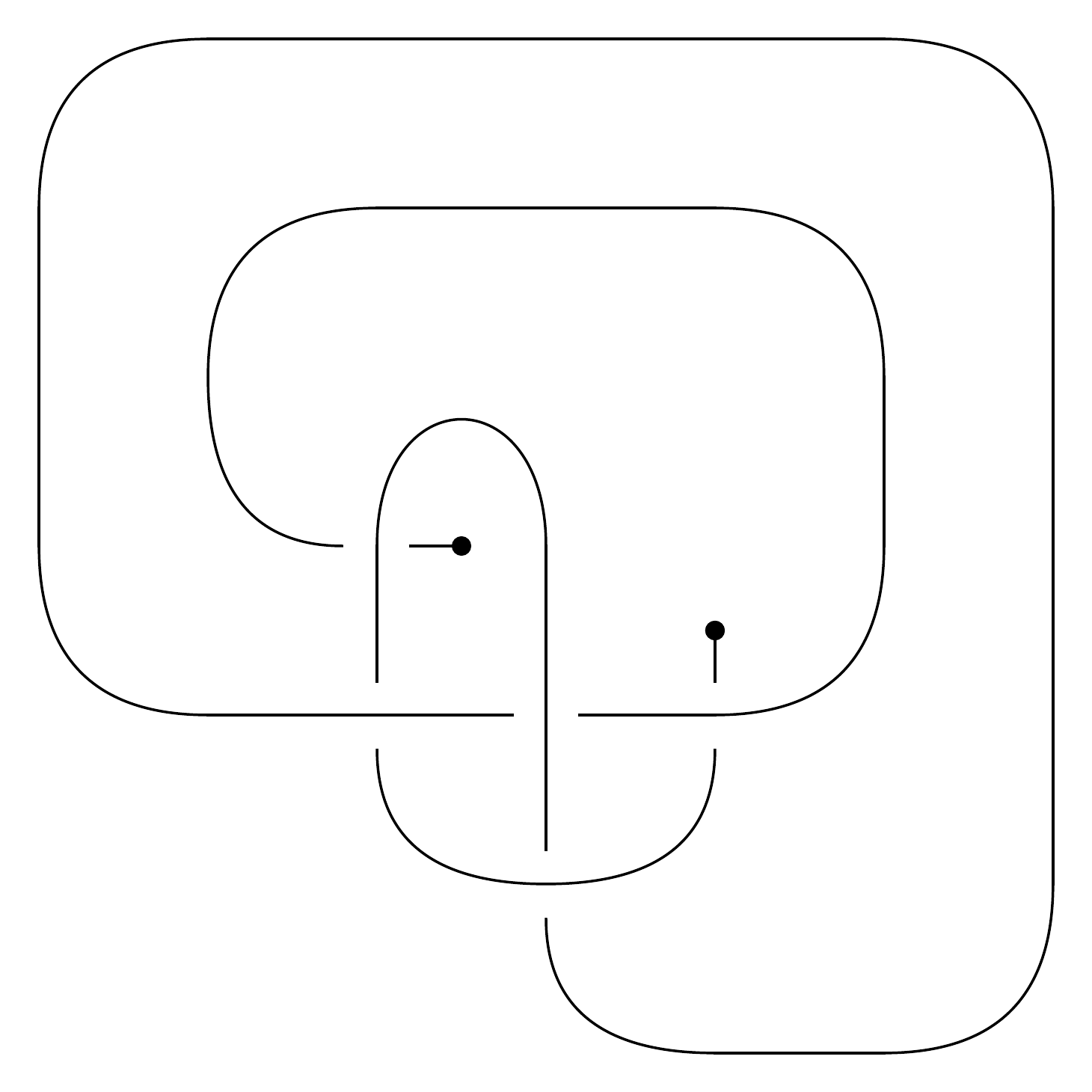}\\
\textcolor{black}{$5_{830}$}
\vspace{1cm}
\end{minipage}
\begin{minipage}[t]{.25\linewidth}
\centering
\includegraphics[width=0.9\textwidth,height=3.5cm,keepaspectratio]{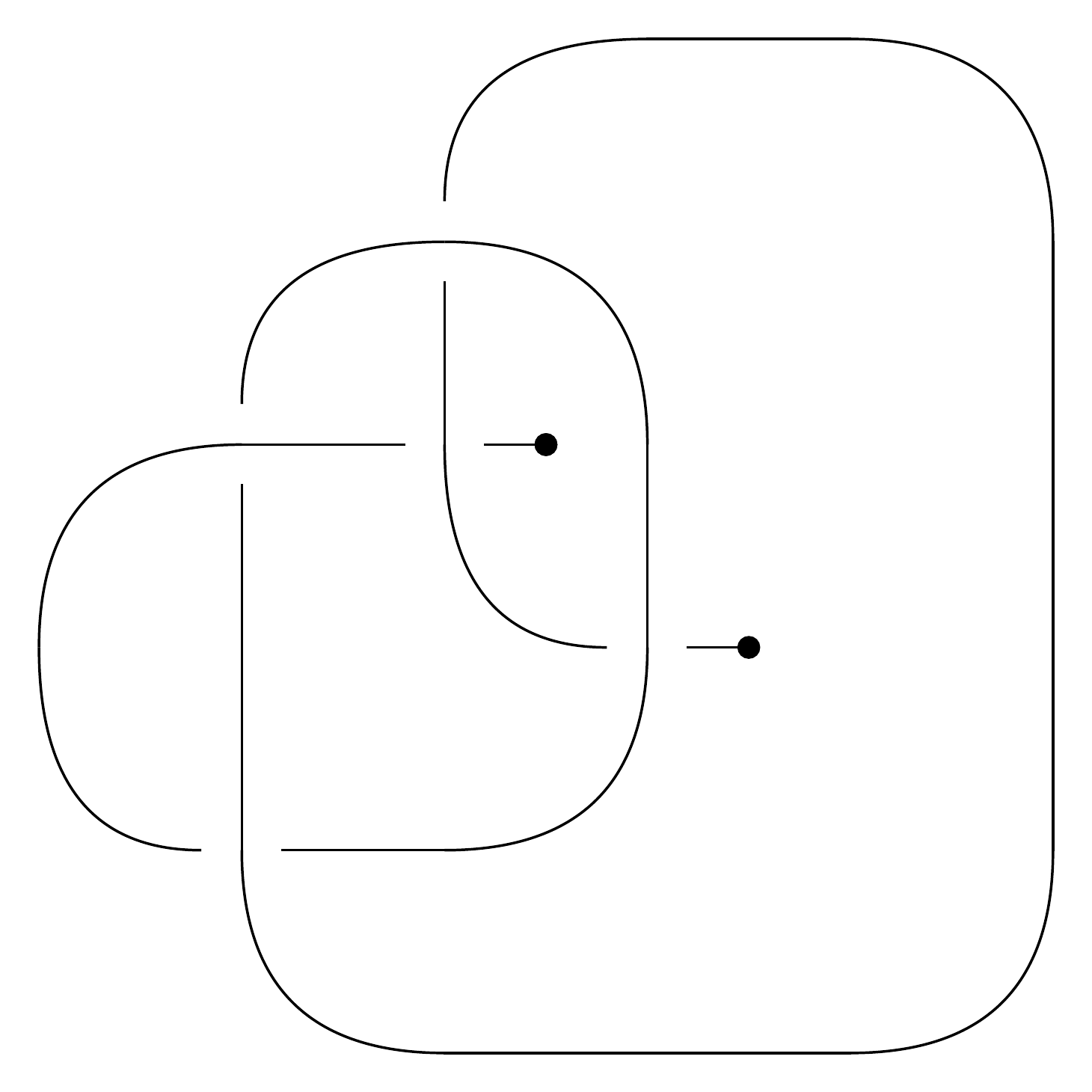}\\
\textcolor{black}{$5_{831}$}
\vspace{1cm}
\end{minipage}
\begin{minipage}[t]{.25\linewidth}
\centering
\includegraphics[width=0.9\textwidth,height=3.5cm,keepaspectratio]{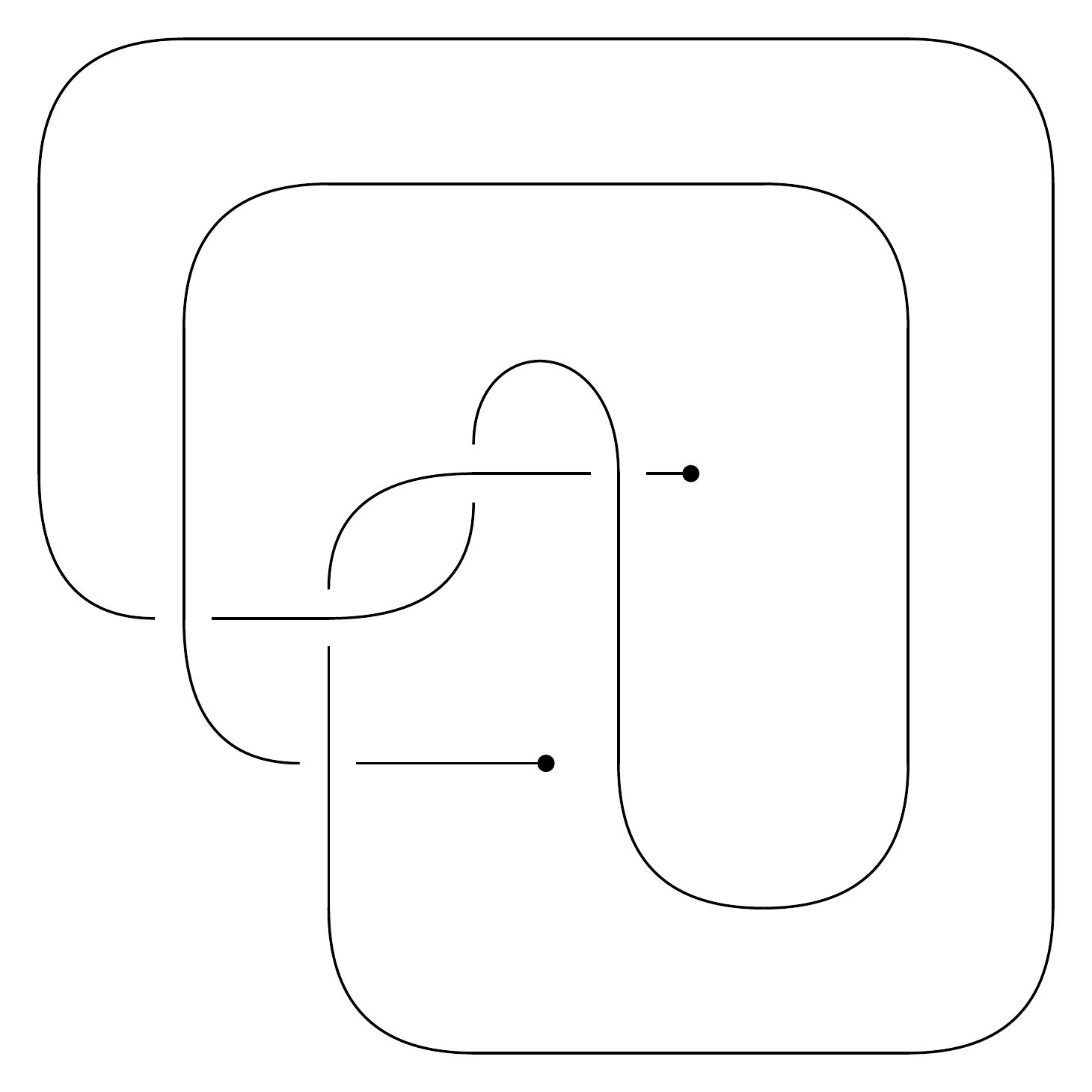}\\
\textcolor{black}{$5_{832}$}
\vspace{1cm}
\end{minipage}
\begin{minipage}[t]{.25\linewidth}
\centering
\includegraphics[width=0.9\textwidth,height=3.5cm,keepaspectratio]{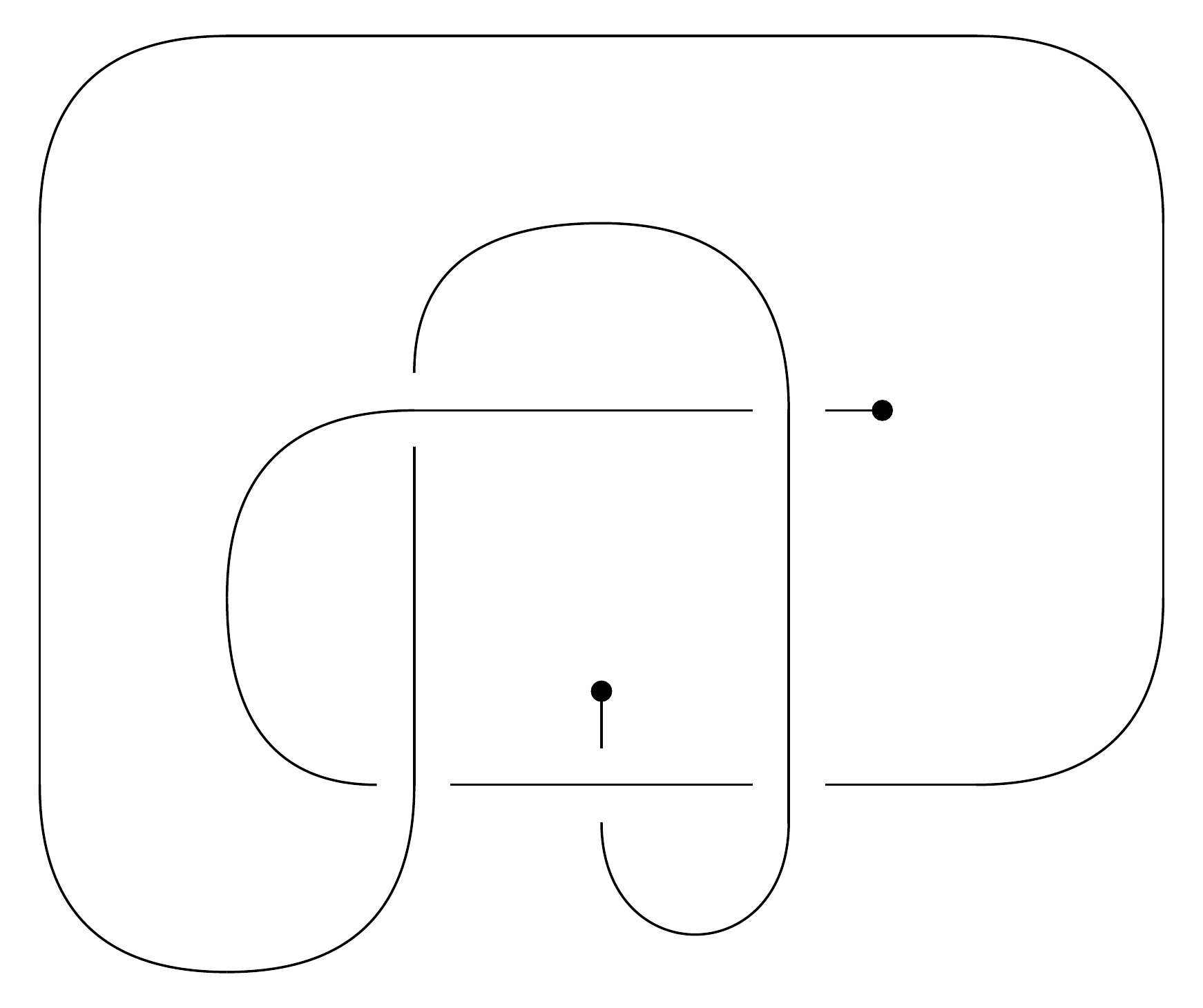}\\
\textcolor{black}{$5_{833}$}
\vspace{1cm}
\end{minipage}
\begin{minipage}[t]{.25\linewidth}
\centering
\includegraphics[width=0.9\textwidth,height=3.5cm,keepaspectratio]{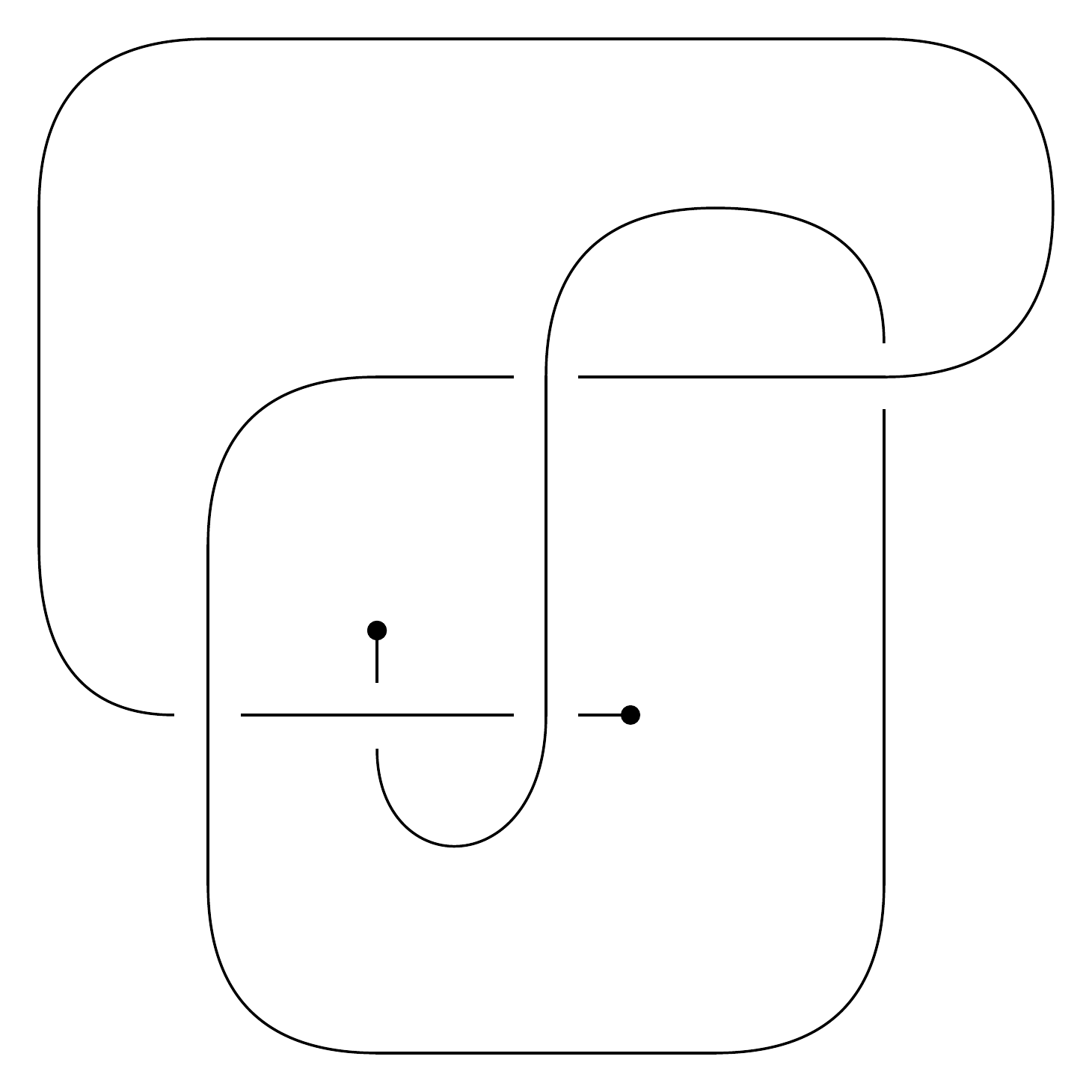}\\
\textcolor{black}{$5_{834}$}
\vspace{1cm}
\end{minipage}
\begin{minipage}[t]{.25\linewidth}
\centering
\includegraphics[width=0.9\textwidth,height=3.5cm,keepaspectratio]{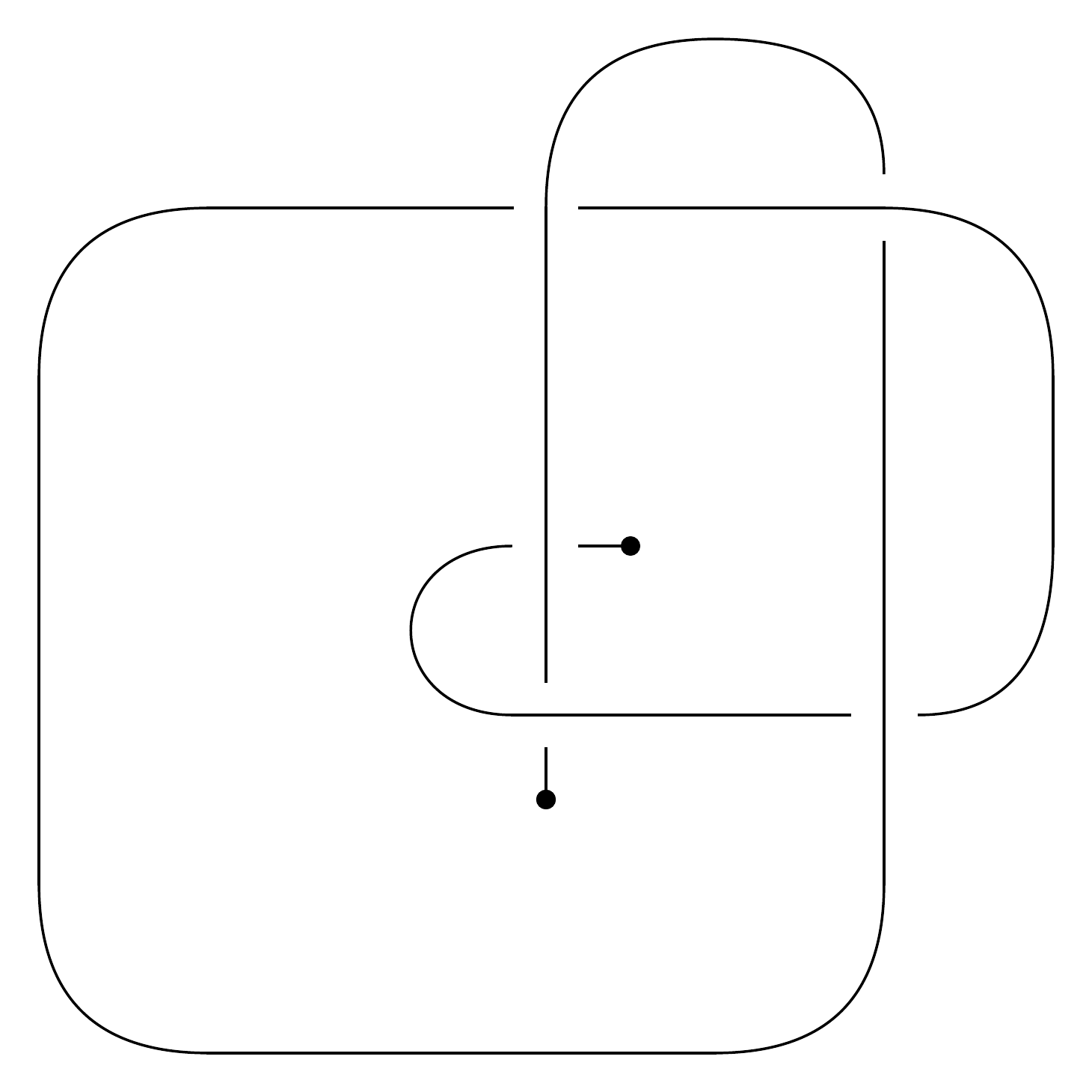}\\
\textcolor{black}{$5_{835}$}
\vspace{1cm}
\end{minipage}
\begin{minipage}[t]{.25\linewidth}
\centering
\includegraphics[width=0.9\textwidth,height=3.5cm,keepaspectratio]{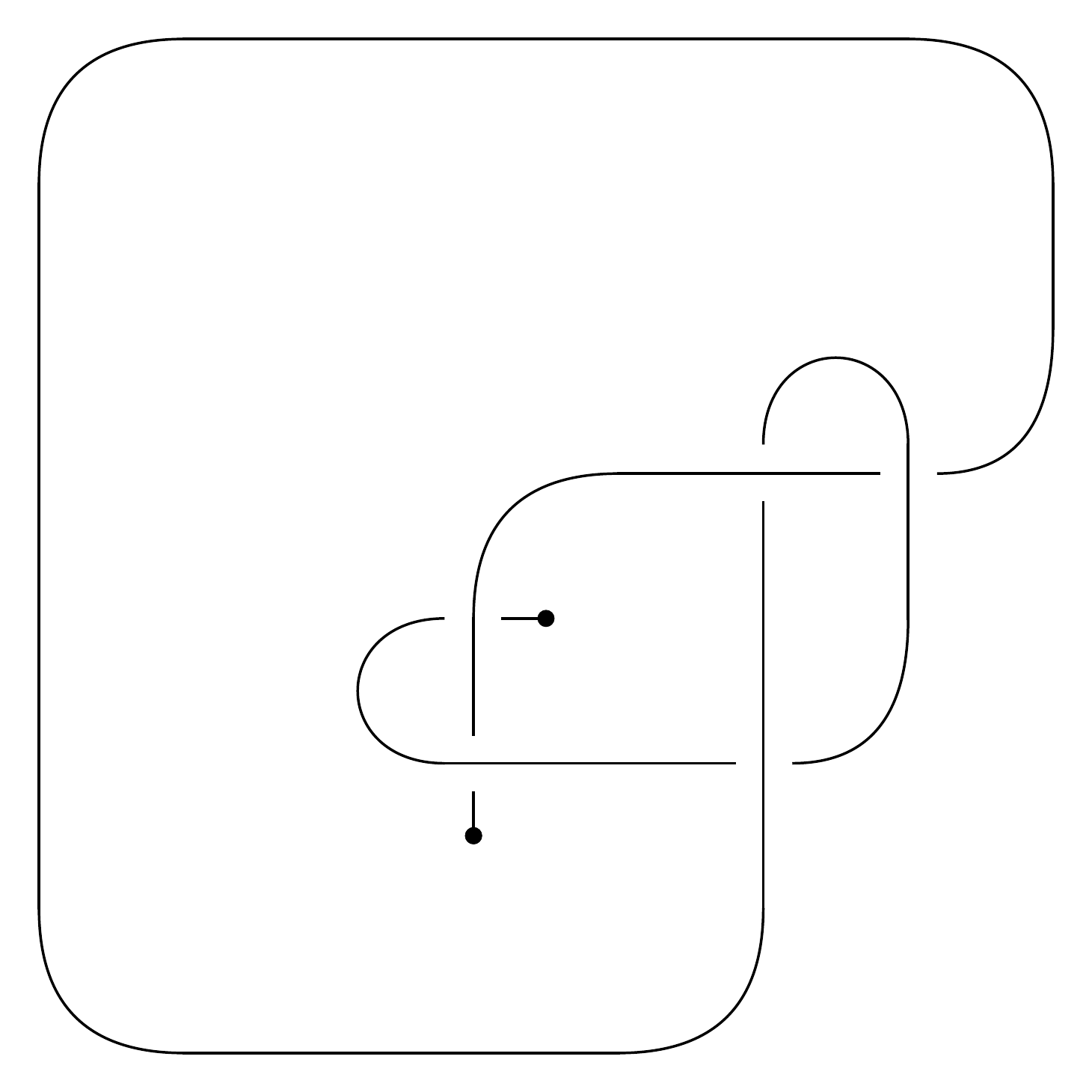}\\
\textcolor{black}{$5_{836}$}
\vspace{1cm}
\end{minipage}
\begin{minipage}[t]{.25\linewidth}
\centering
\includegraphics[width=0.9\textwidth,height=3.5cm,keepaspectratio]{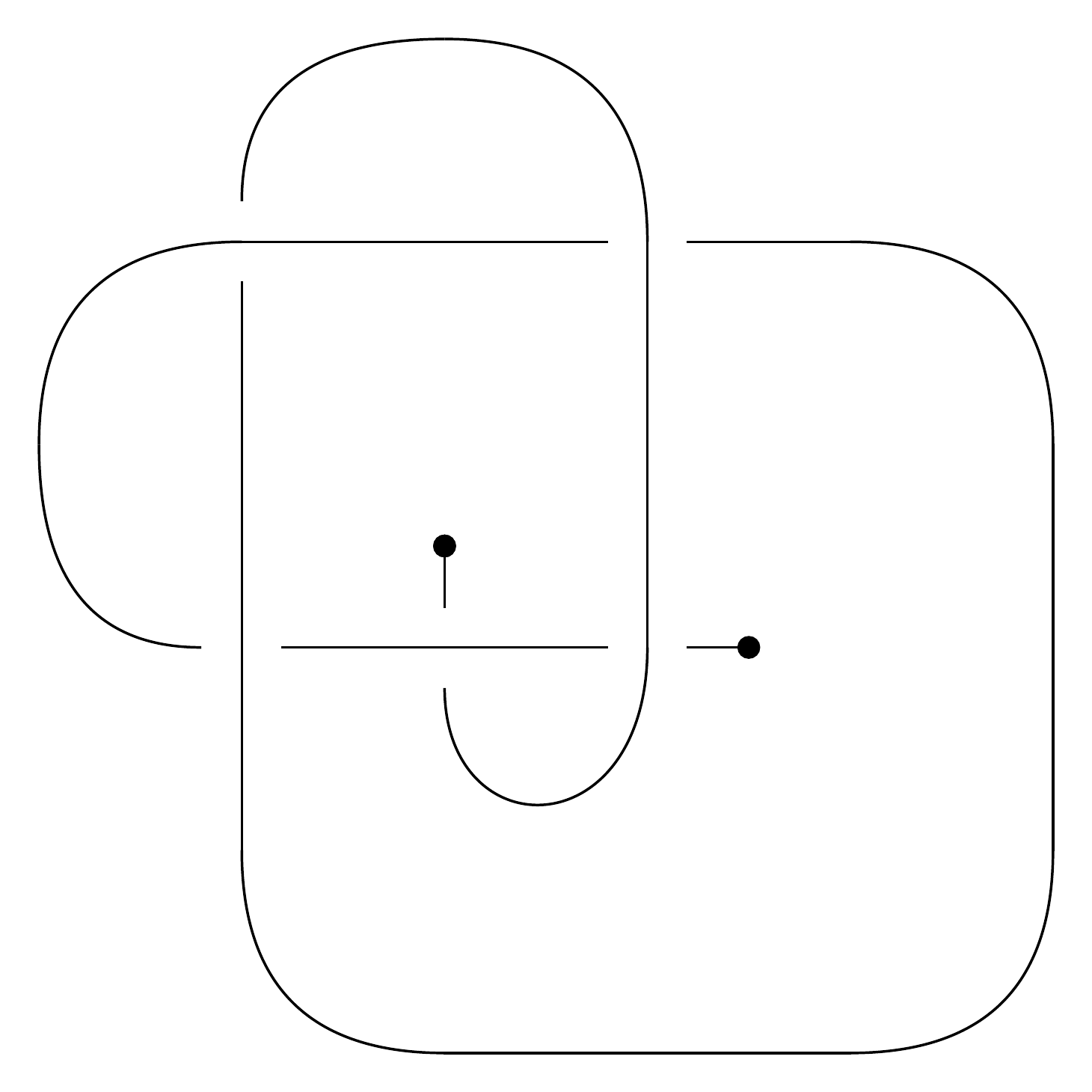}\\
\textcolor{black}{$5_{837}$}
\vspace{1cm}
\end{minipage}
\begin{minipage}[t]{.25\linewidth}
\centering
\includegraphics[width=0.9\textwidth,height=3.5cm,keepaspectratio]{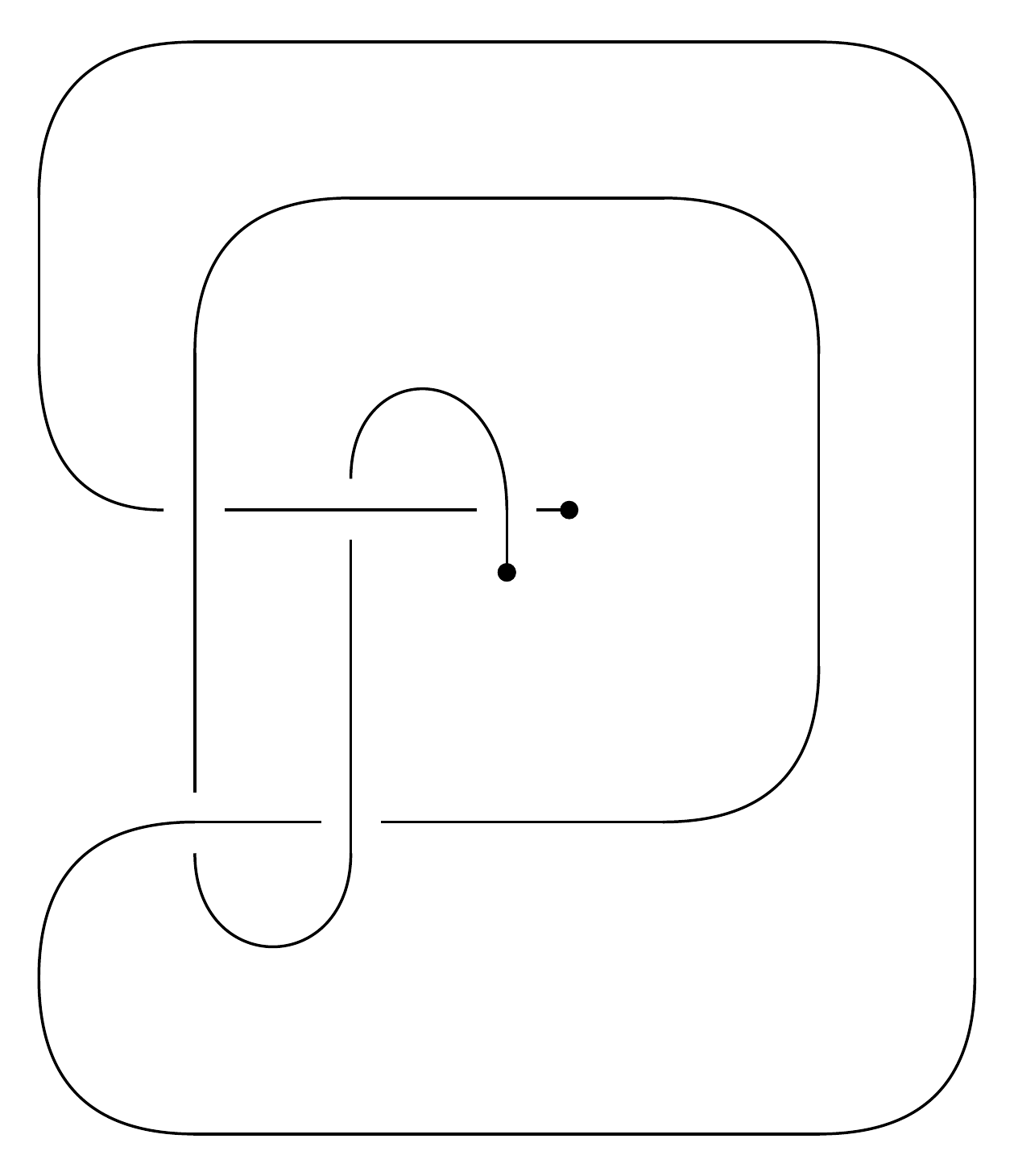}\\
\textcolor{black}{$5_{838}$}
\vspace{1cm}
\end{minipage}
\begin{minipage}[t]{.25\linewidth}
\centering
\includegraphics[width=0.9\textwidth,height=3.5cm,keepaspectratio]{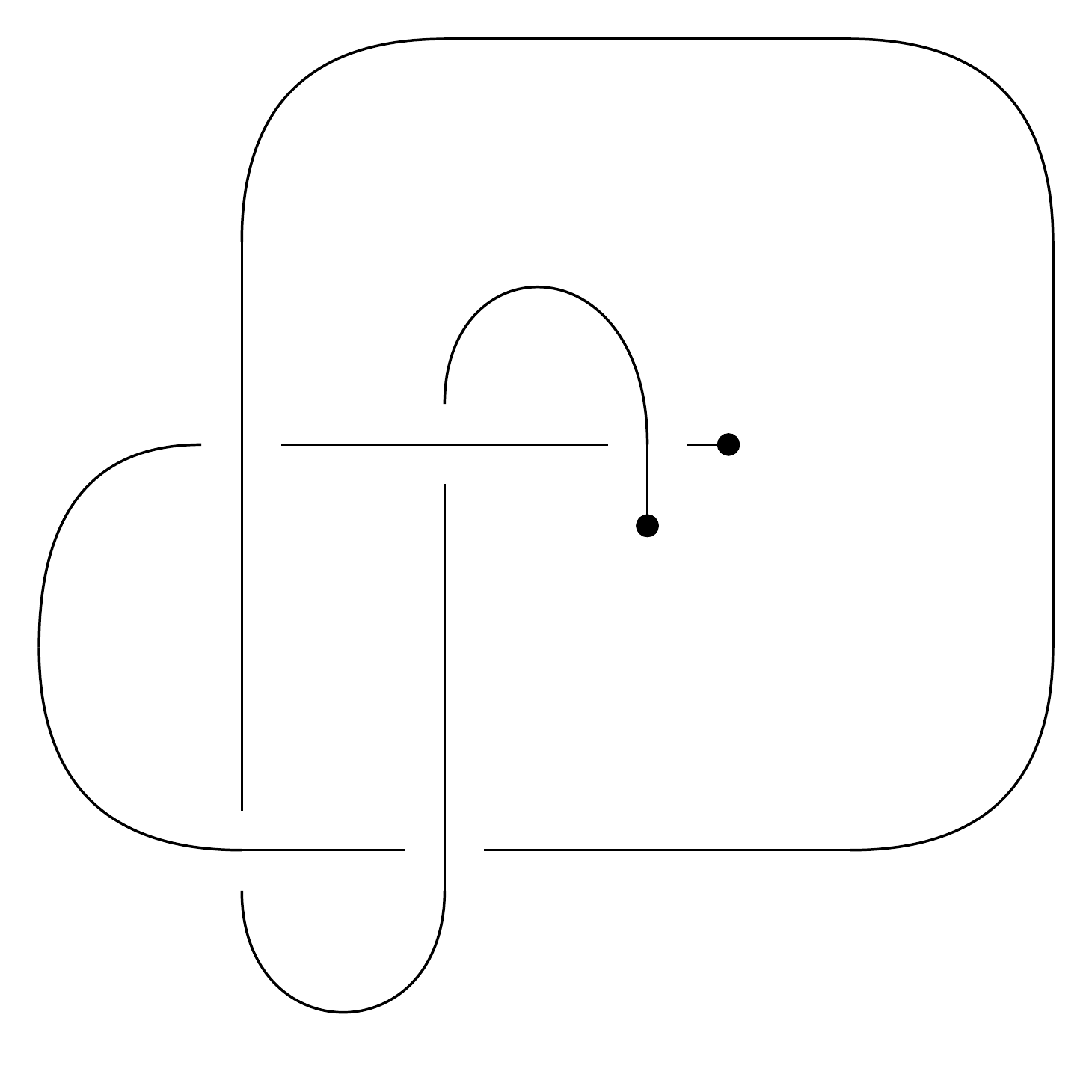}\\
\textcolor{black}{$5_{839}$}
\vspace{1cm}
\end{minipage}
\begin{minipage}[t]{.25\linewidth}
\centering
\includegraphics[width=0.9\textwidth,height=3.5cm,keepaspectratio]{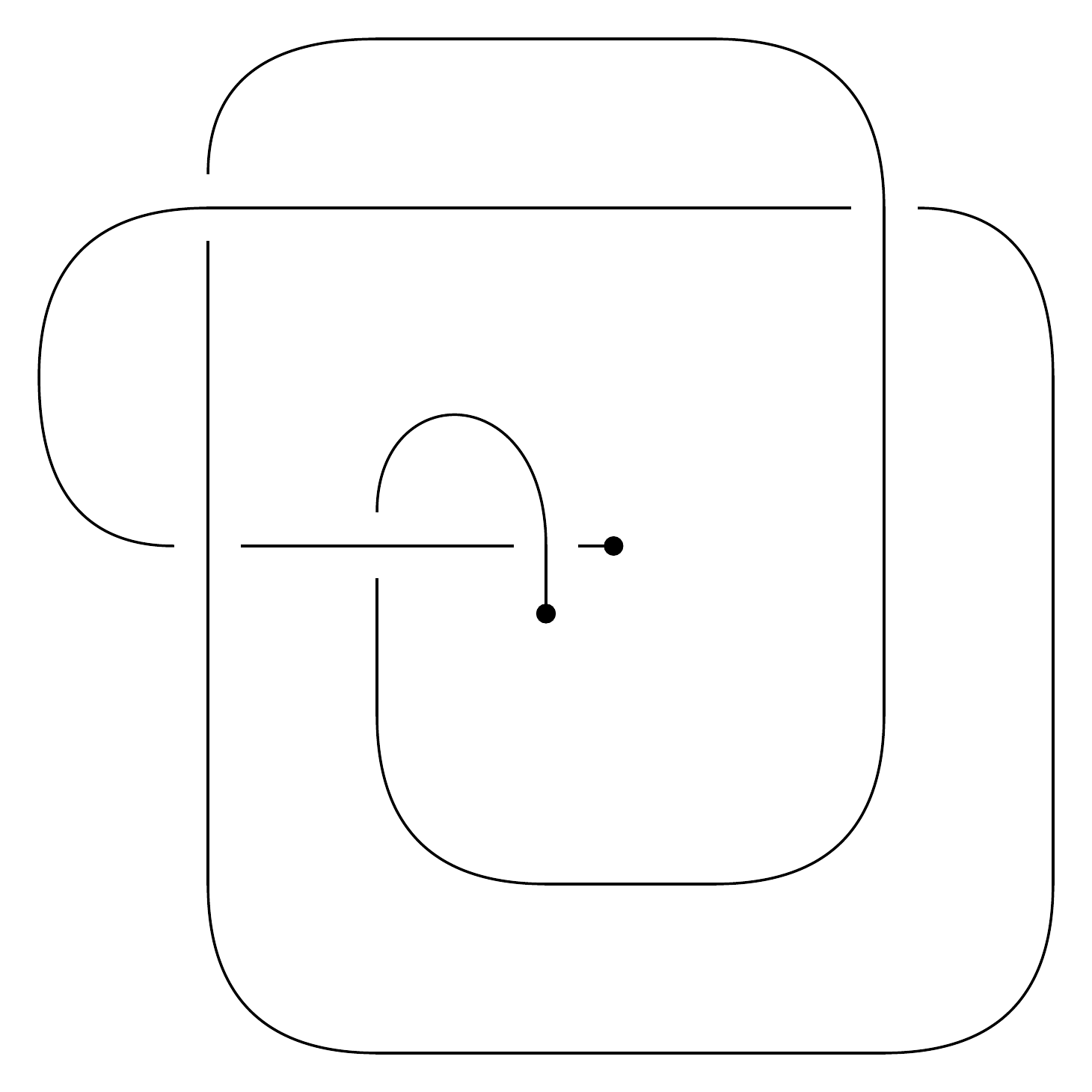}\\
\textcolor{black}{$5_{840}$}
\vspace{1cm}
\end{minipage}
\begin{minipage}[t]{.25\linewidth}
\centering
\includegraphics[width=0.9\textwidth,height=3.5cm,keepaspectratio]{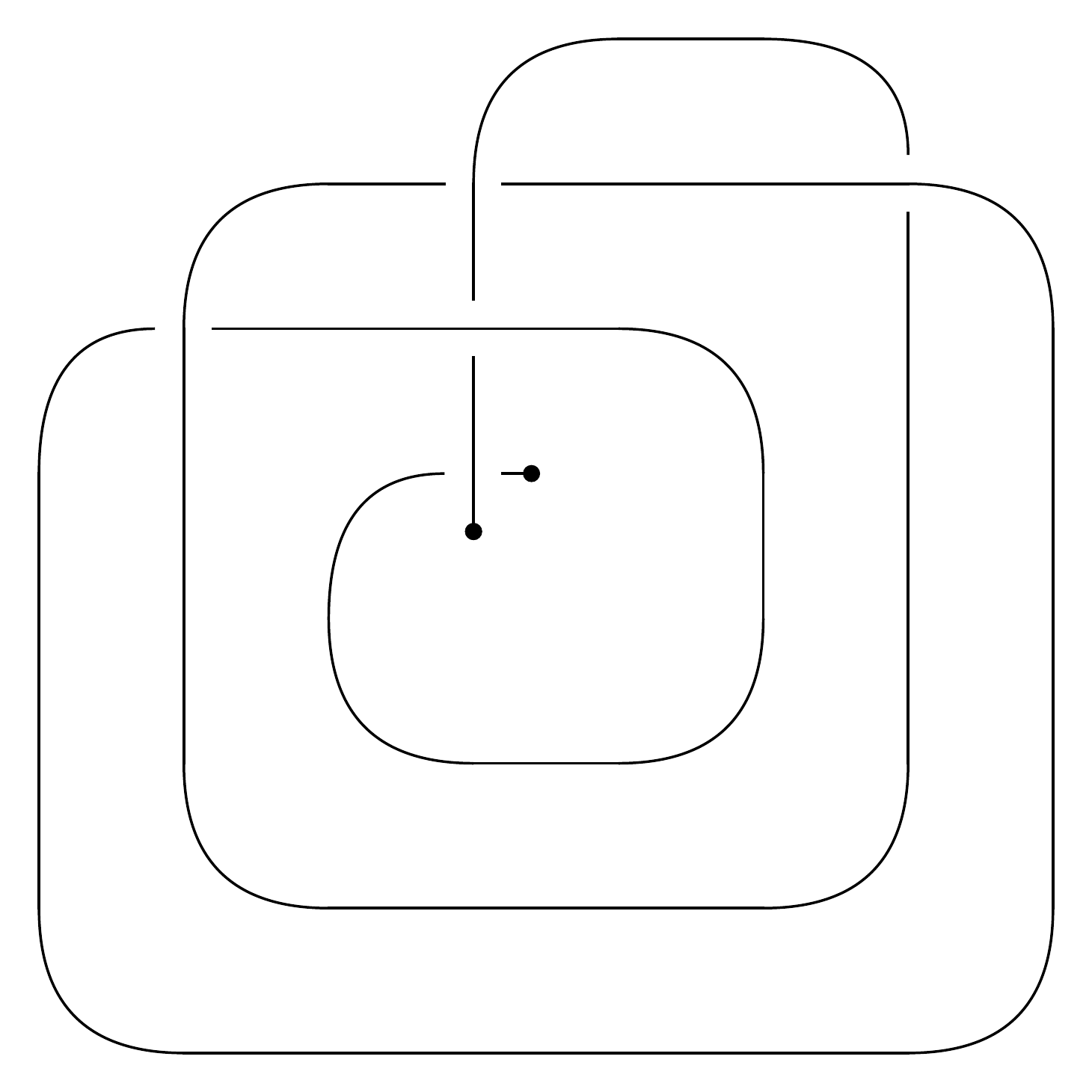}\\
\textcolor{black}{$5_{841}$}
\vspace{1cm}
\end{minipage}
\begin{minipage}[t]{.25\linewidth}
\centering
\includegraphics[width=0.9\textwidth,height=3.5cm,keepaspectratio]{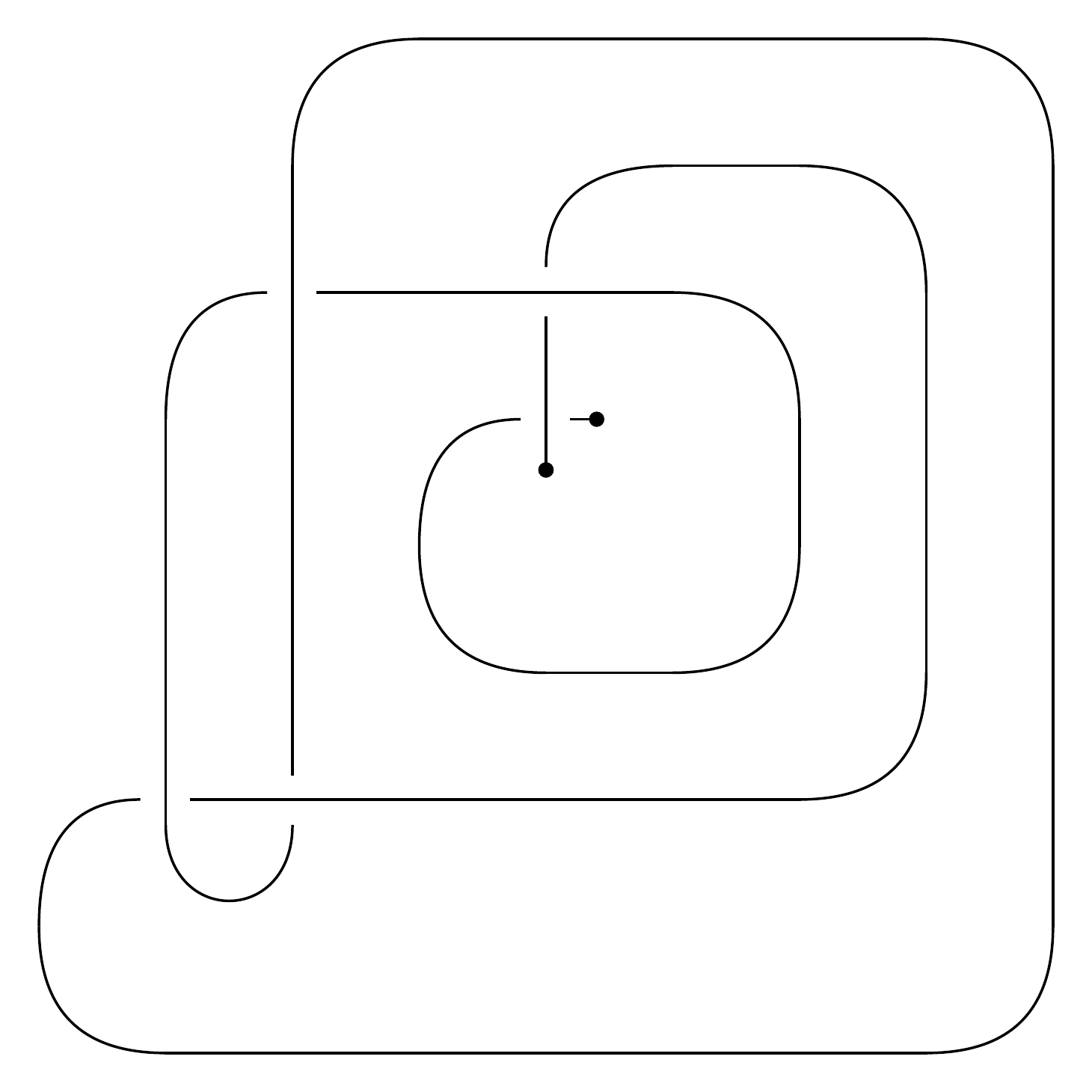}\\
\textcolor{black}{$5_{842}$}
\vspace{1cm}
\end{minipage}
\begin{minipage}[t]{.25\linewidth}
\centering
\includegraphics[width=0.9\textwidth,height=3.5cm,keepaspectratio]{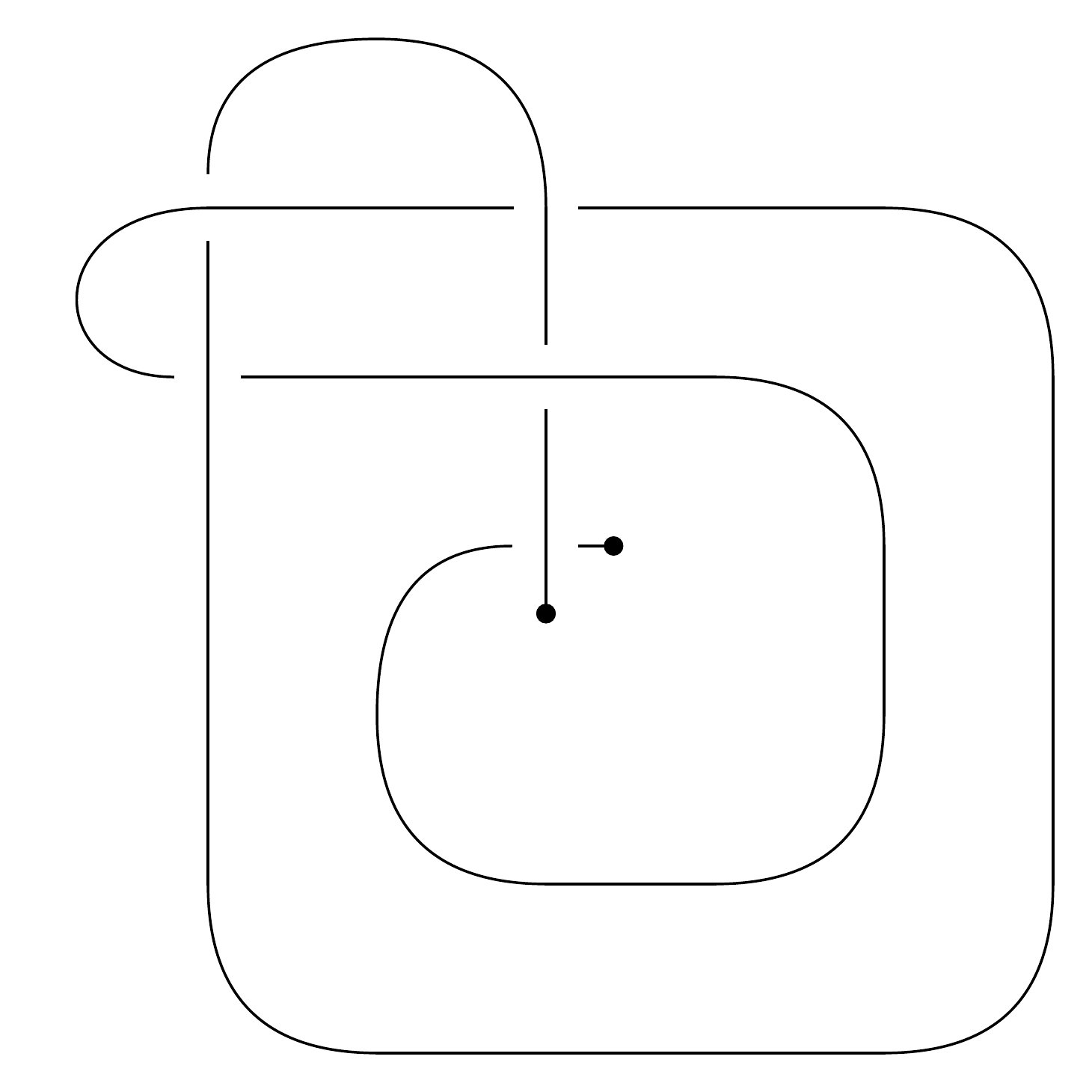}\\
\textcolor{black}{$5_{843}$}
\vspace{1cm}
\end{minipage}
\begin{minipage}[t]{.25\linewidth}
\centering
\includegraphics[width=0.9\textwidth,height=3.5cm,keepaspectratio]{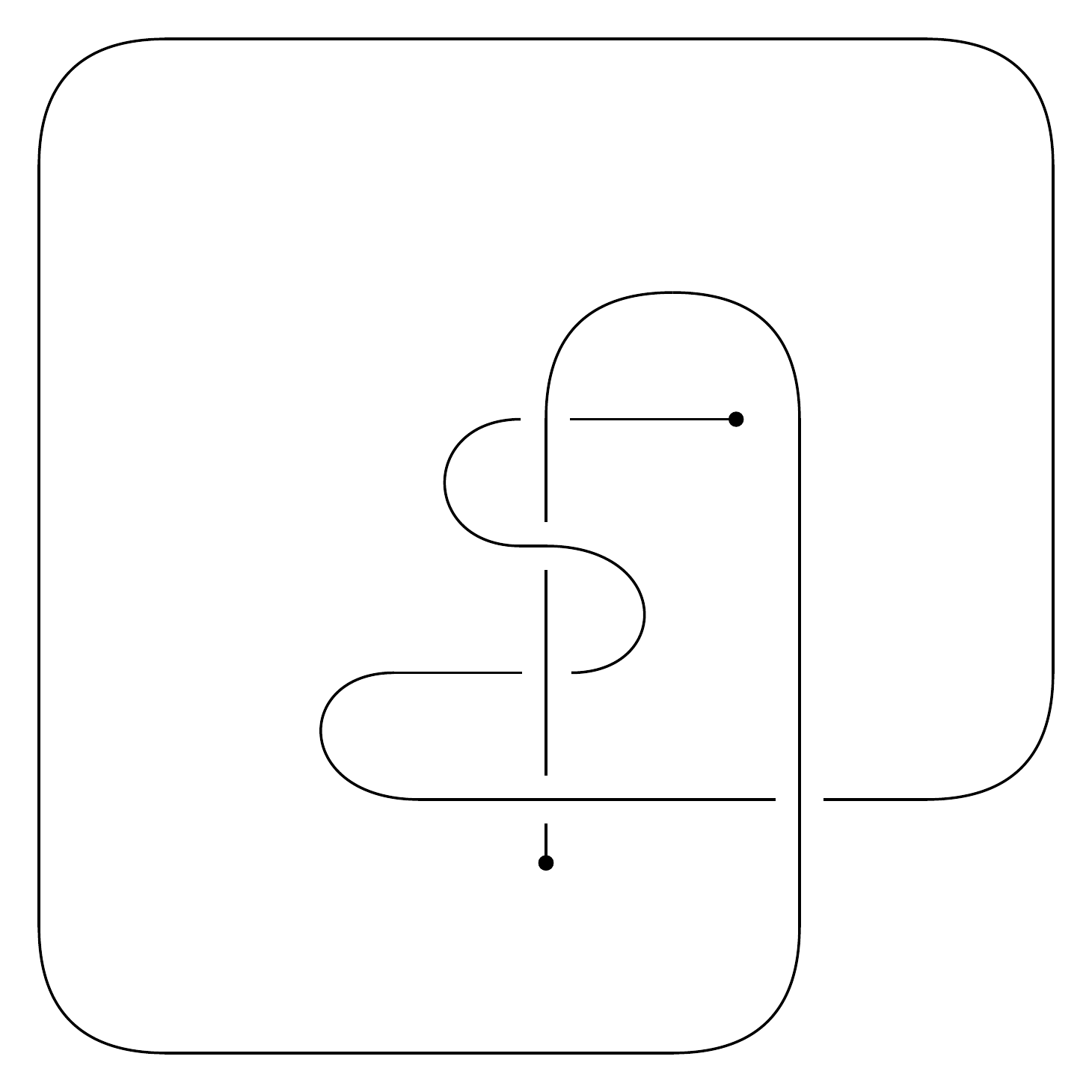}\\
\textcolor{black}{$5_{844}$}
\vspace{1cm}
\end{minipage}
\begin{minipage}[t]{.25\linewidth}
\centering
\includegraphics[width=0.9\textwidth,height=3.5cm,keepaspectratio]{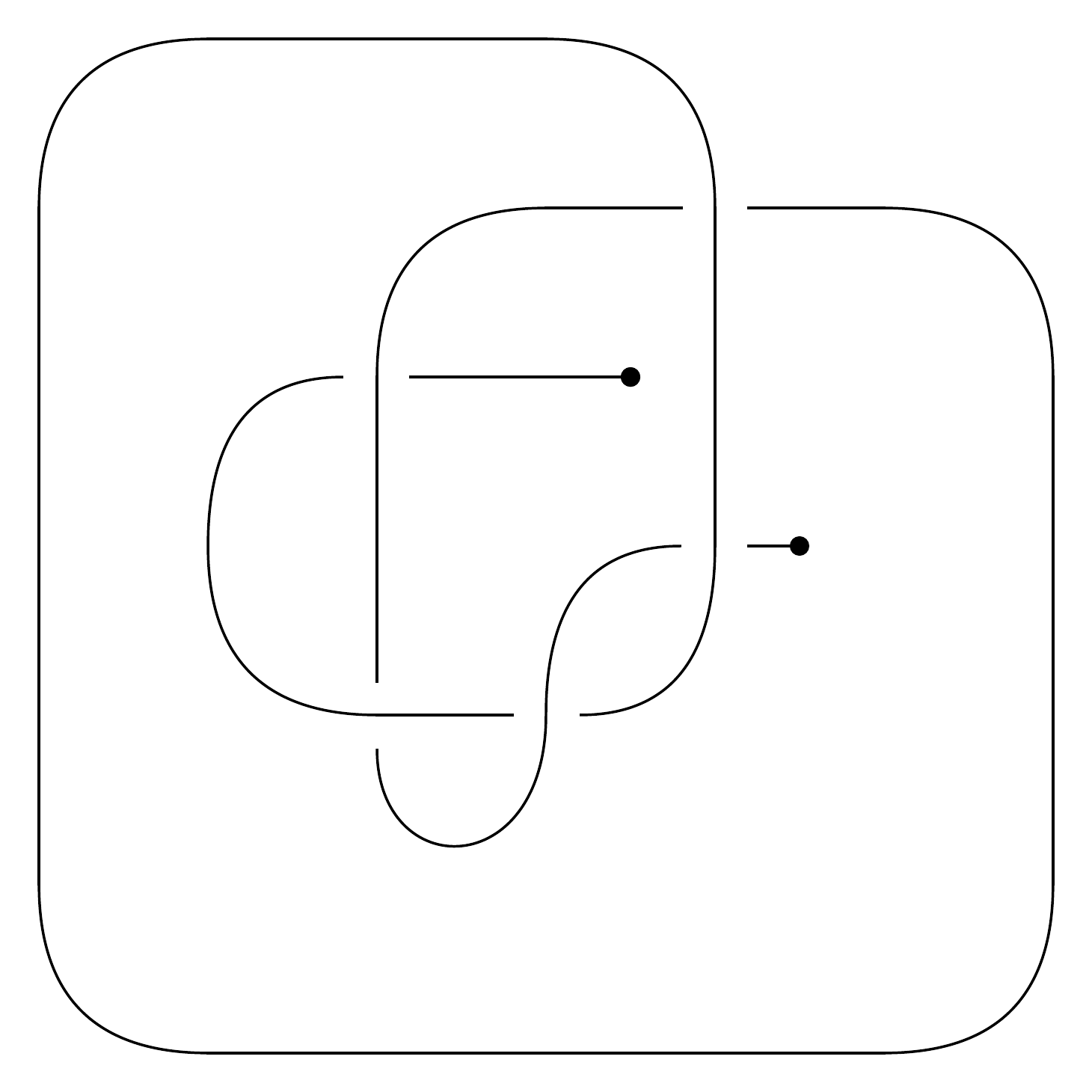}\\
\textcolor{black}{$5_{845}$}
\vspace{1cm}
\end{minipage}
\begin{minipage}[t]{.25\linewidth}
\centering
\includegraphics[width=0.9\textwidth,height=3.5cm,keepaspectratio]{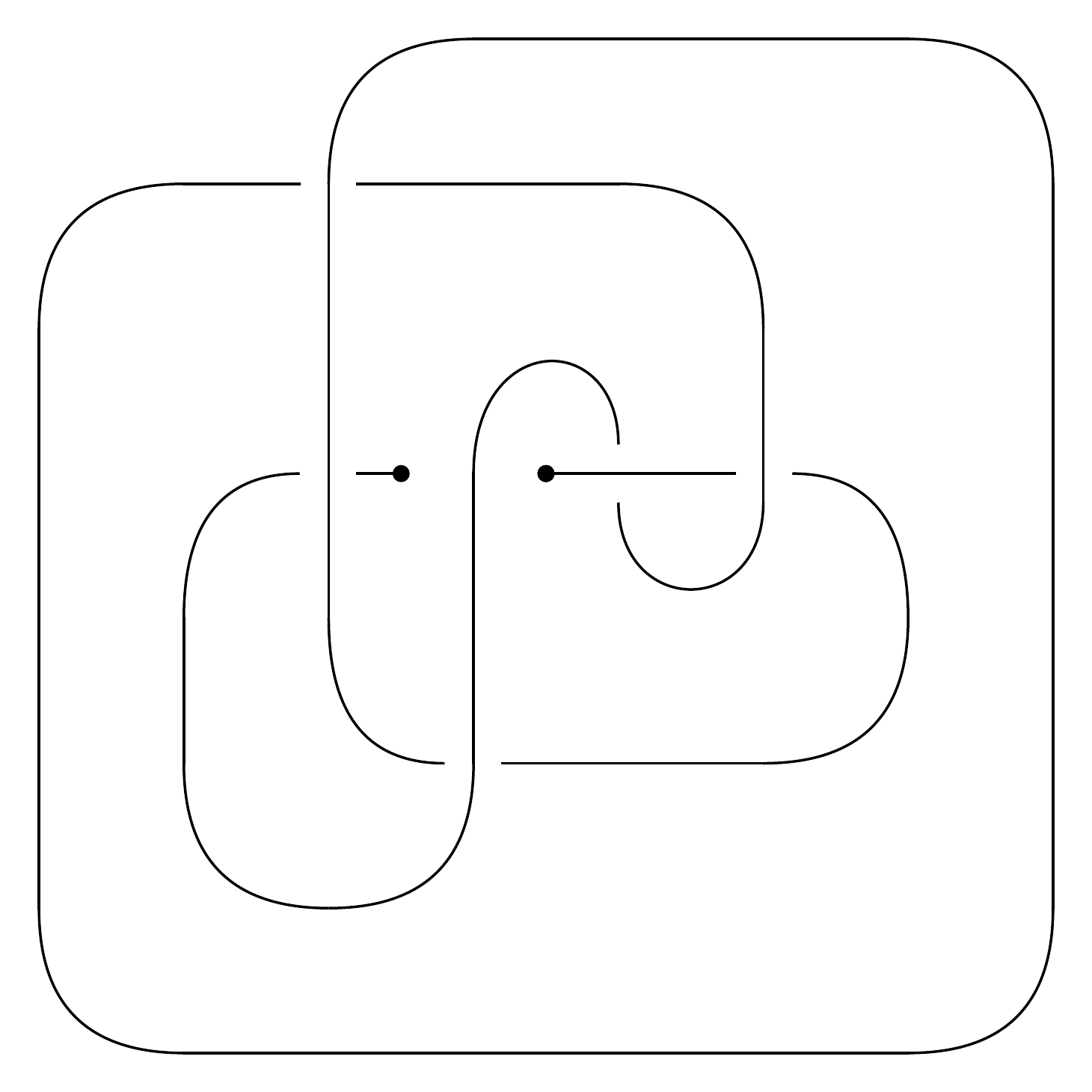}\\
\textcolor{black}{$5_{846}$}
\vspace{1cm}
\end{minipage}
\begin{minipage}[t]{.25\linewidth}
\centering
\includegraphics[width=0.9\textwidth,height=3.5cm,keepaspectratio]{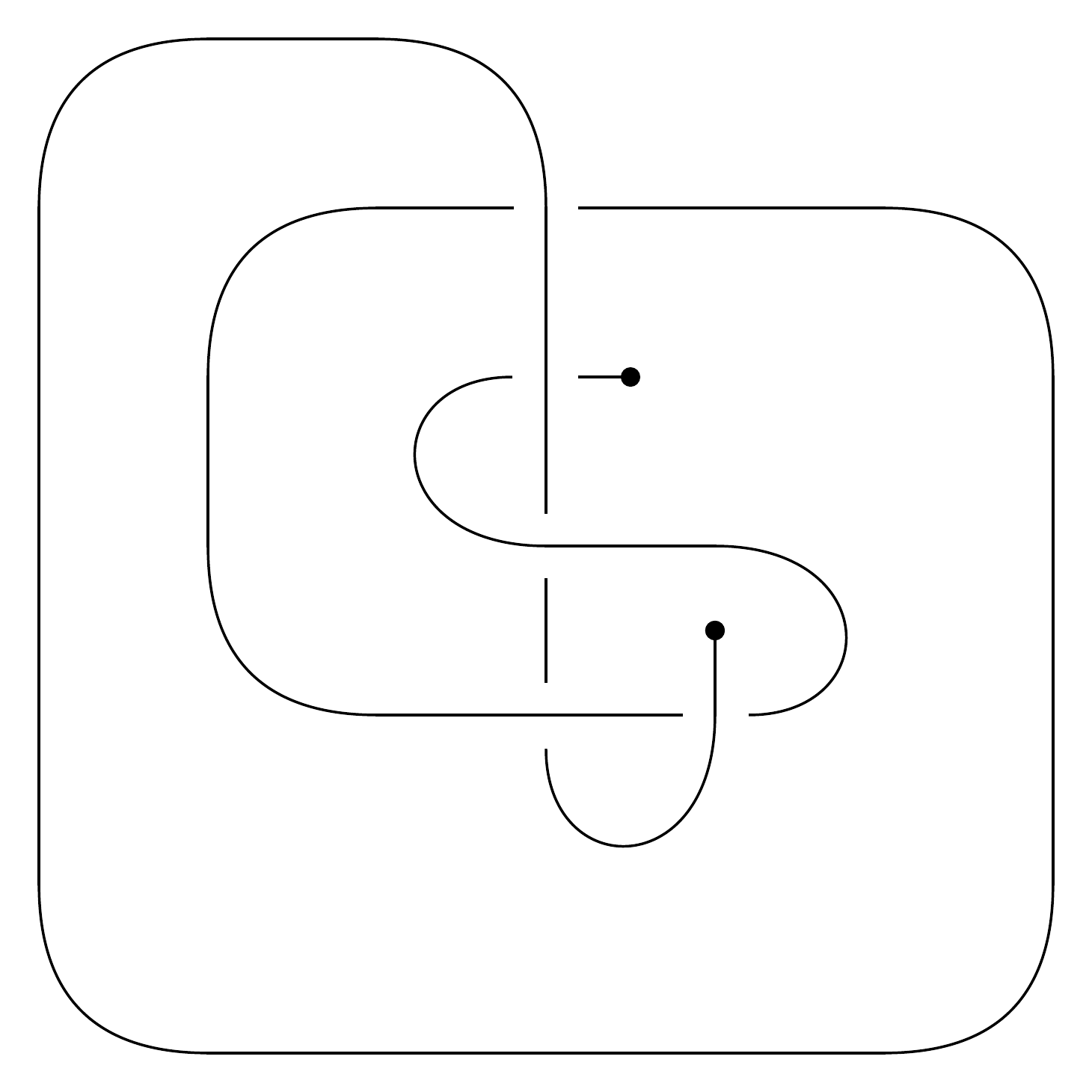}\\
\textcolor{black}{$5_{847}$}
\vspace{1cm}
\end{minipage}
\begin{minipage}[t]{.25\linewidth}
\centering
\includegraphics[width=0.9\textwidth,height=3.5cm,keepaspectratio]{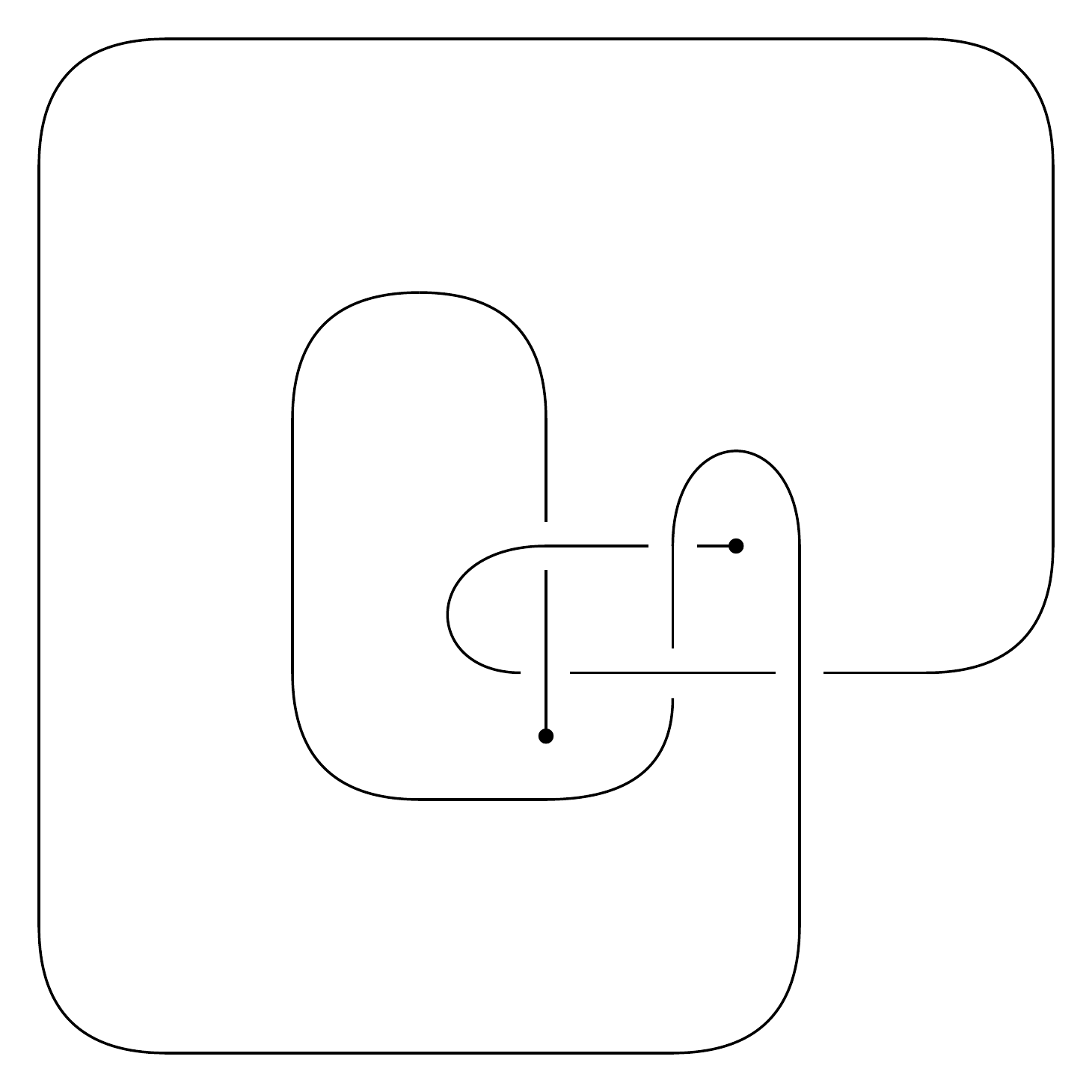}\\
\textcolor{black}{$5_{848}$}
\vspace{1cm}
\end{minipage}
\begin{minipage}[t]{.25\linewidth}
\centering
\includegraphics[width=0.9\textwidth,height=3.5cm,keepaspectratio]{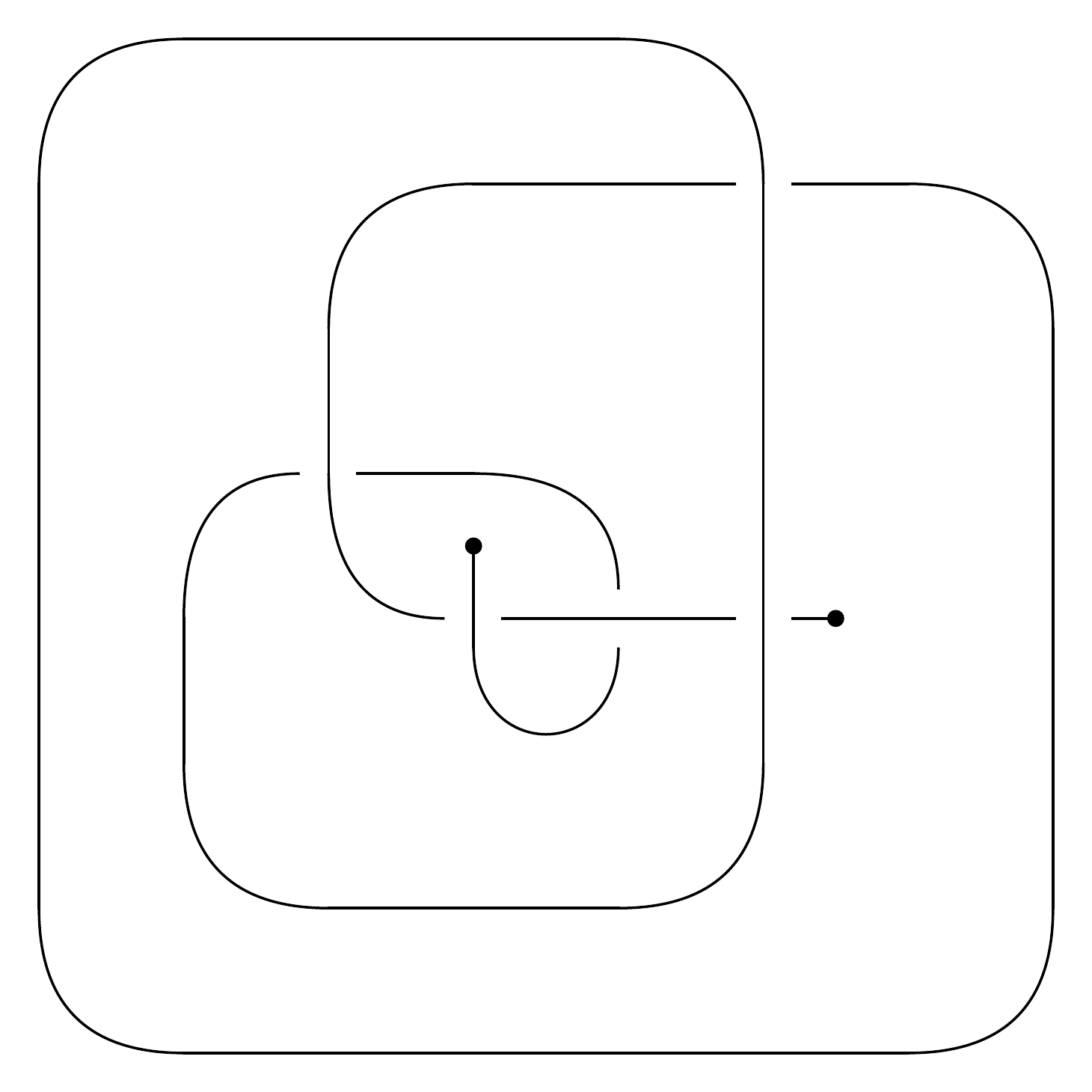}\\
\textcolor{black}{$5_{849}$}
\vspace{1cm}
\end{minipage}
\begin{minipage}[t]{.25\linewidth}
\centering
\includegraphics[width=0.9\textwidth,height=3.5cm,keepaspectratio]{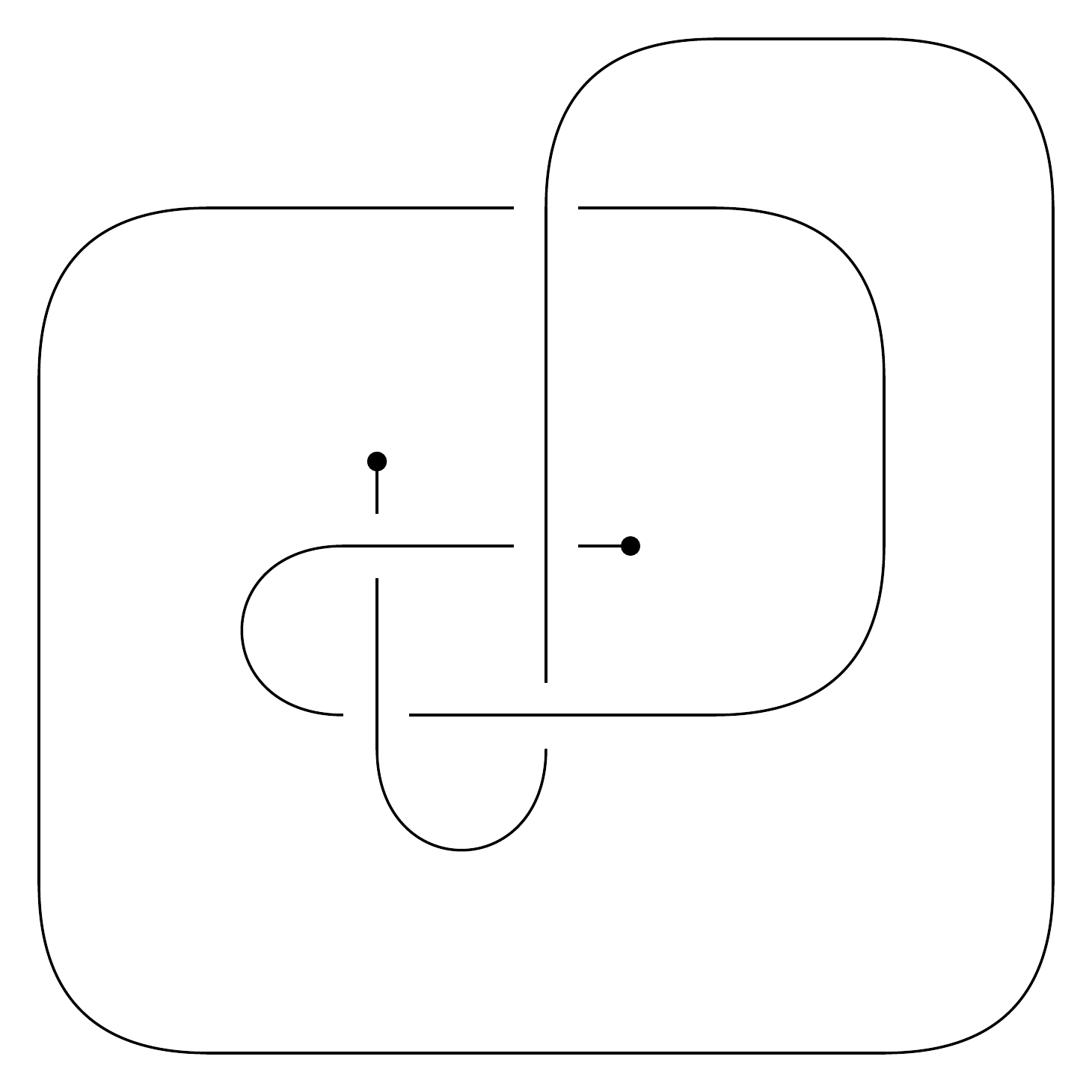}\\
\textcolor{black}{$5_{850}$}
\vspace{1cm}
\end{minipage}
\begin{minipage}[t]{.25\linewidth}
\centering
\includegraphics[width=0.9\textwidth,height=3.5cm,keepaspectratio]{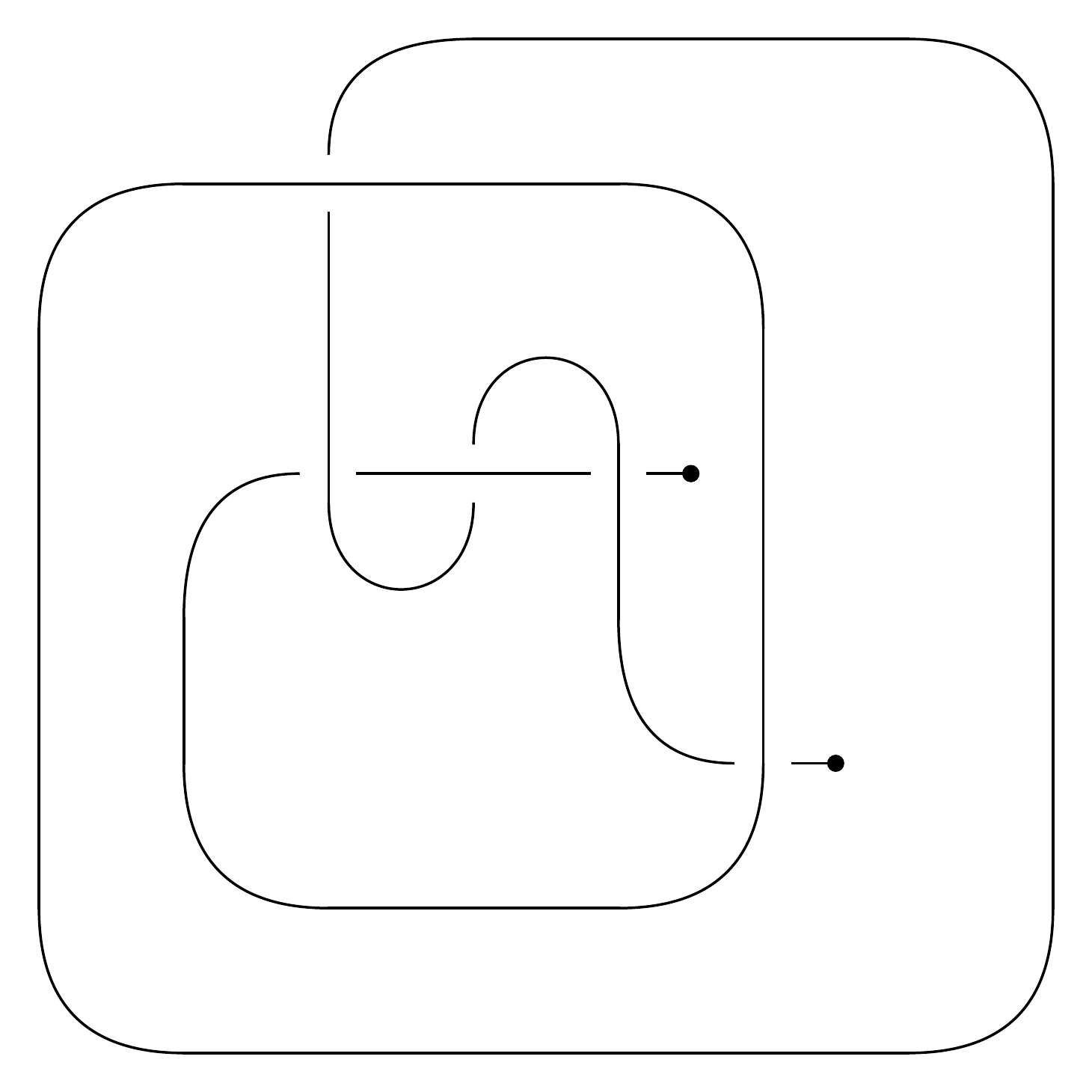}\\
\textcolor{black}{$5_{851}$}
\vspace{1cm}
\end{minipage}
\begin{minipage}[t]{.25\linewidth}
\centering
\includegraphics[width=0.9\textwidth,height=3.5cm,keepaspectratio]{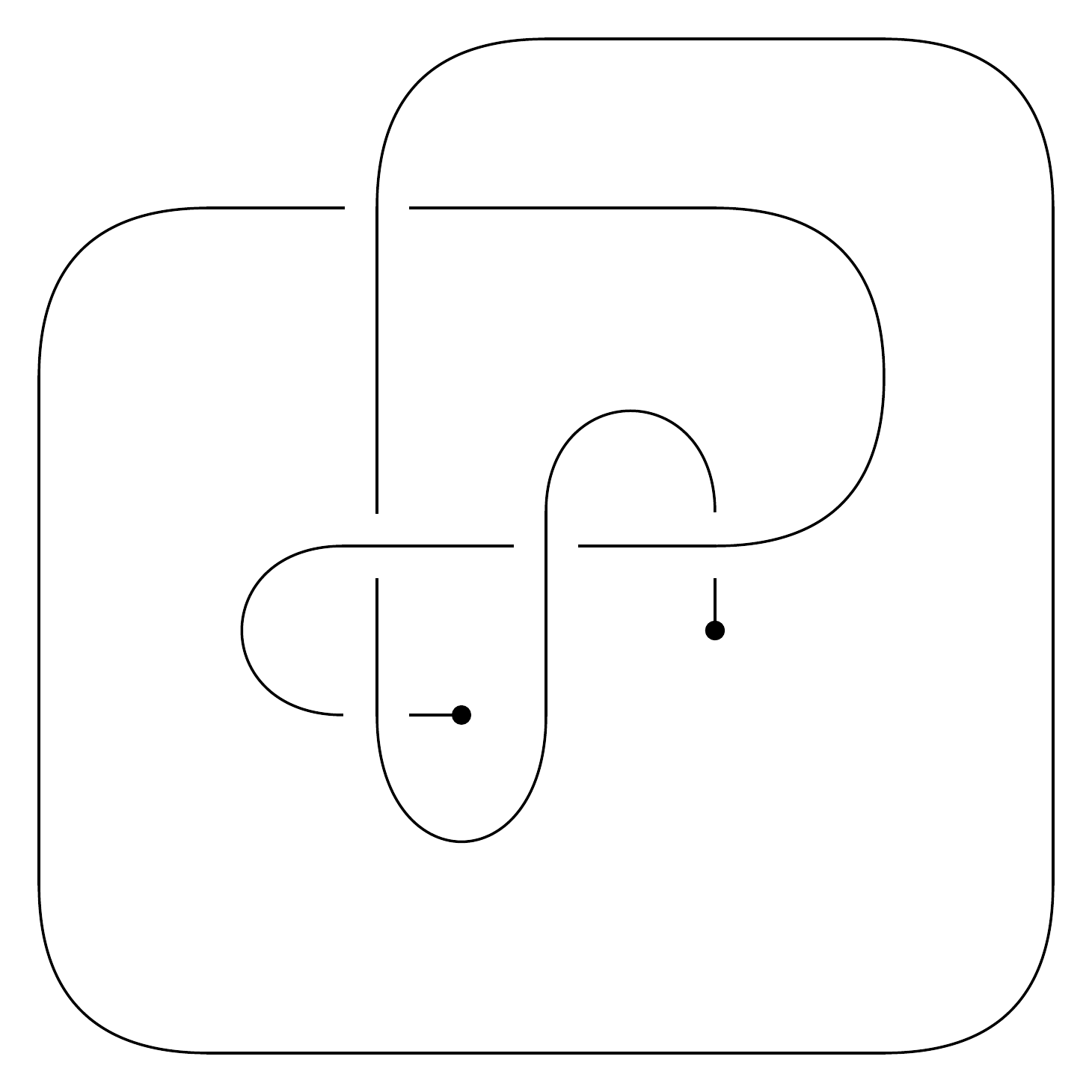}\\
\textcolor{black}{$5_{852}$}
\vspace{1cm}
\end{minipage}
\begin{minipage}[t]{.25\linewidth}
\centering
\includegraphics[width=0.9\textwidth,height=3.5cm,keepaspectratio]{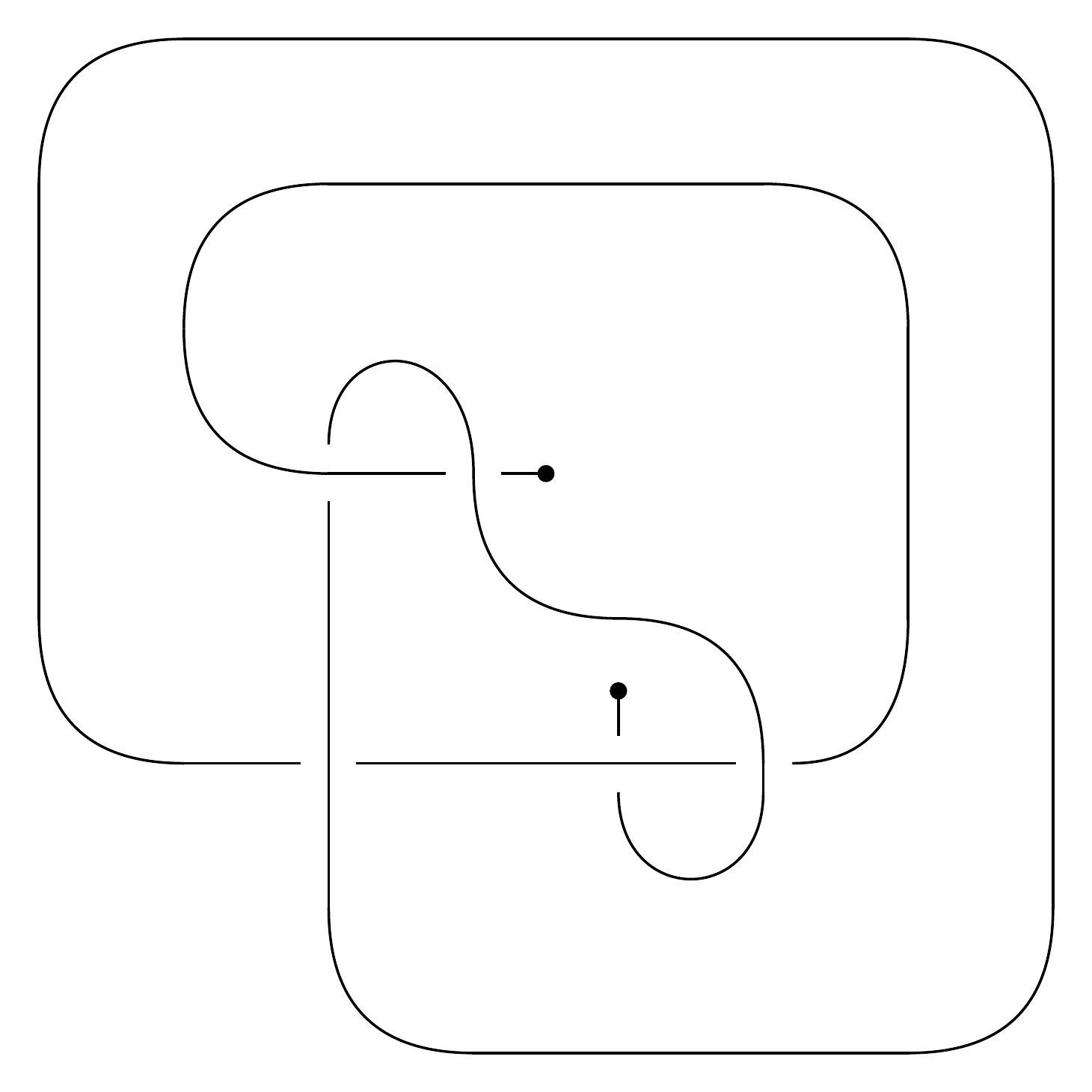}\\
\textcolor{black}{$5_{853}$}
\vspace{1cm}
\end{minipage}
\begin{minipage}[t]{.25\linewidth}
\centering
\includegraphics[width=0.9\textwidth,height=3.5cm,keepaspectratio]{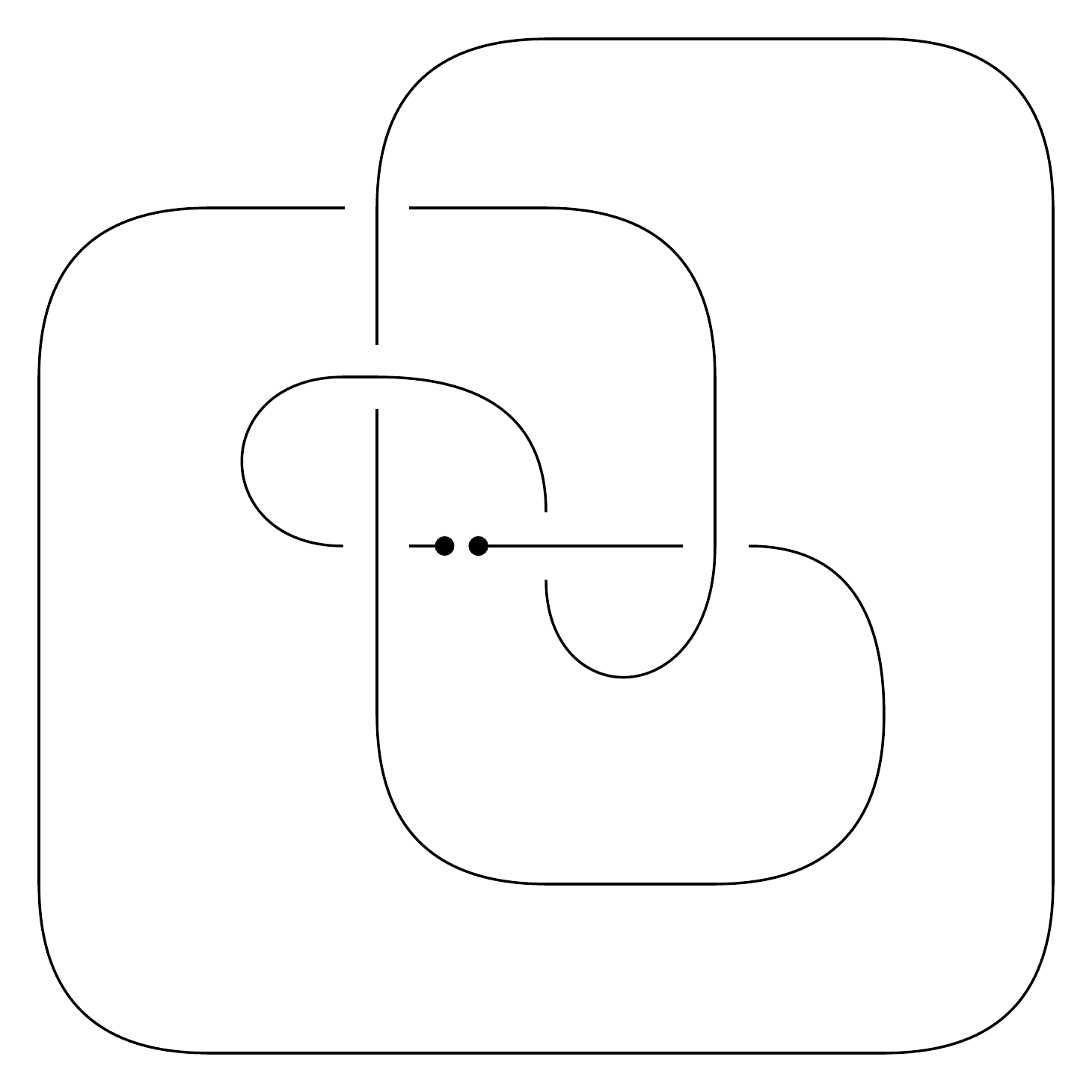}\\
\textcolor{black}{$5_{854}$}
\vspace{1cm}
\end{minipage}
\begin{minipage}[t]{.25\linewidth}
\centering
\includegraphics[width=0.9\textwidth,height=3.5cm,keepaspectratio]{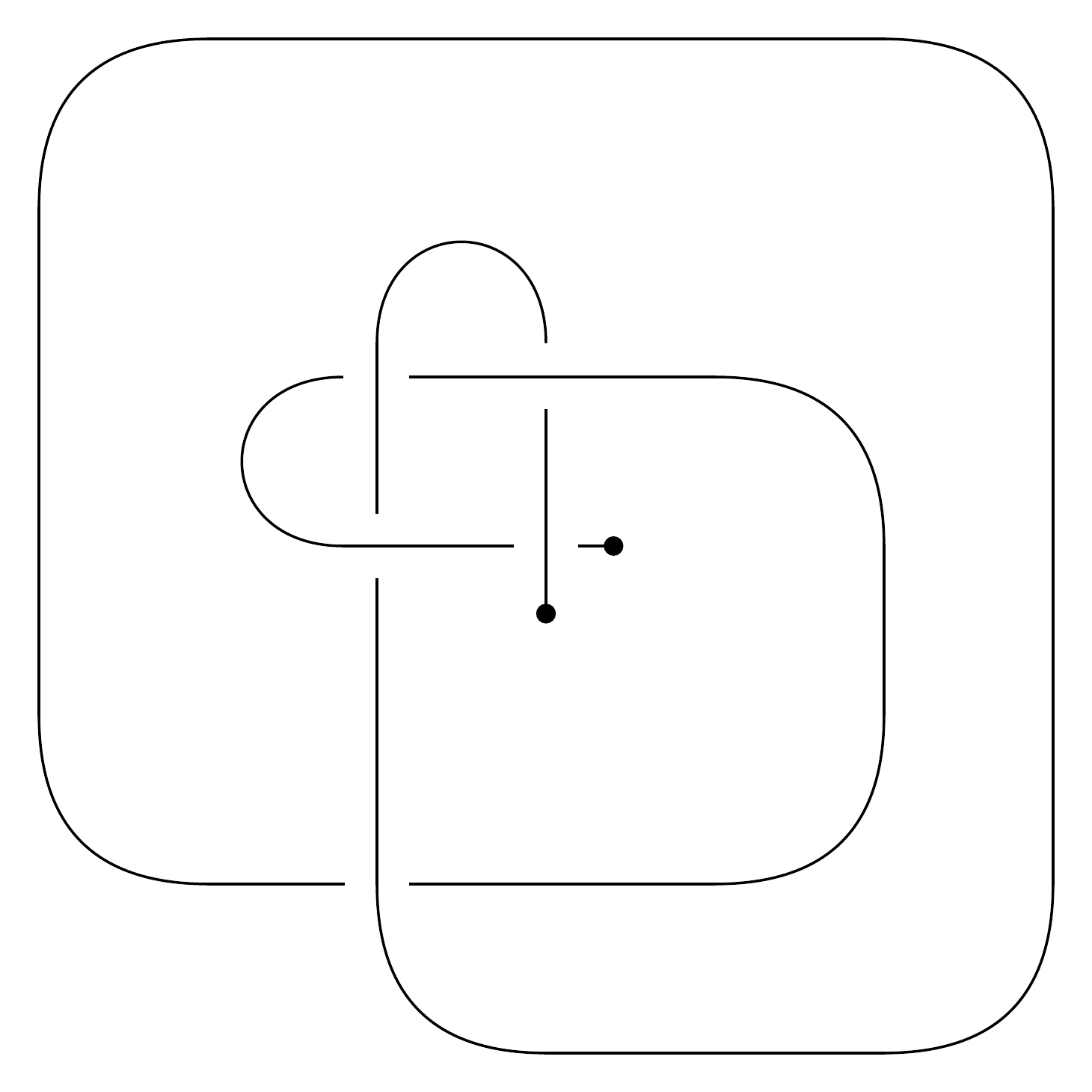}\\
\textcolor{black}{$5_{855}$}
\vspace{1cm}
\end{minipage}
\begin{minipage}[t]{.25\linewidth}
\centering
\includegraphics[width=0.9\textwidth,height=3.5cm,keepaspectratio]{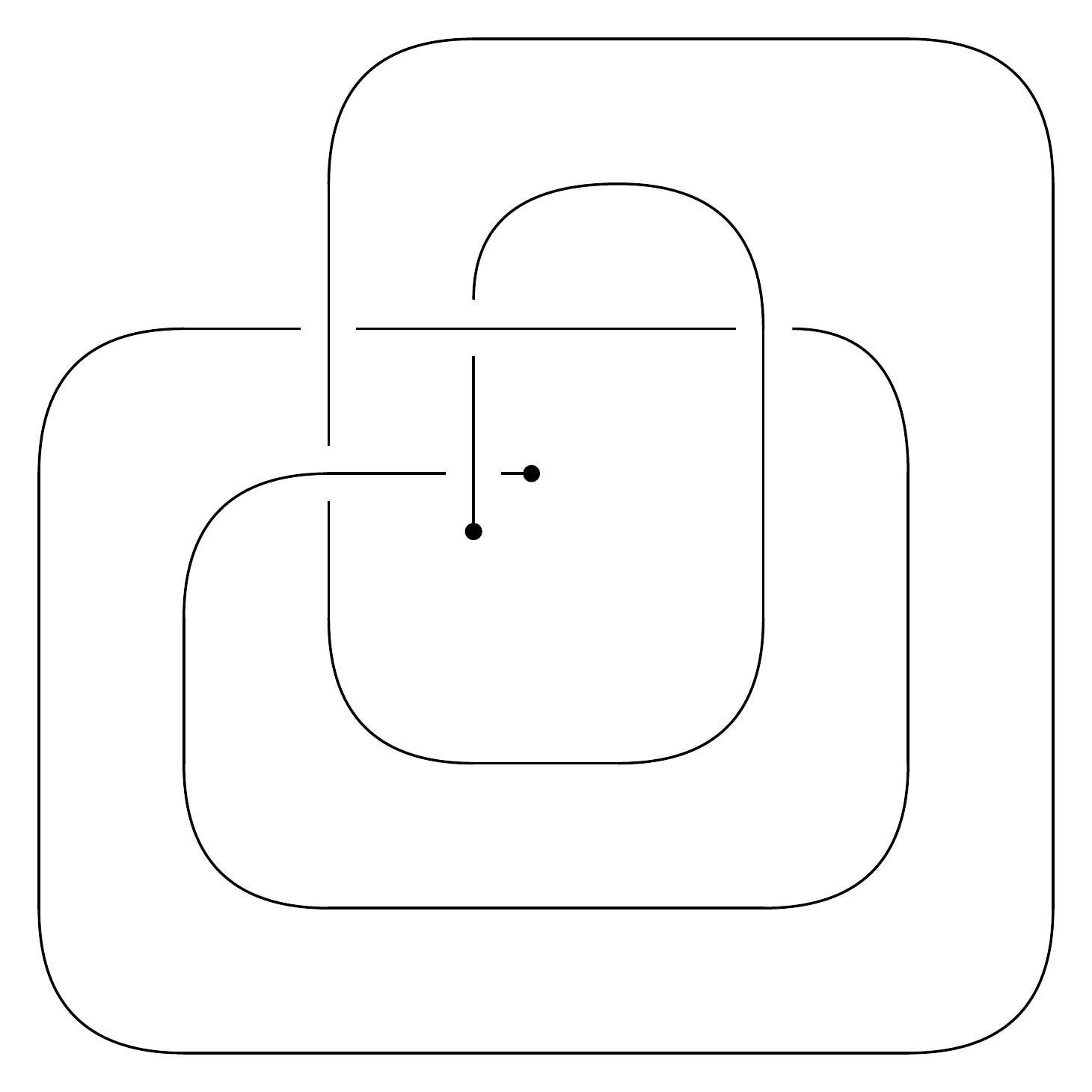}\\
\textcolor{black}{$5_{856}$}
\vspace{1cm}
\end{minipage}
\begin{minipage}[t]{.25\linewidth}
\centering
\includegraphics[width=0.9\textwidth,height=3.5cm,keepaspectratio]{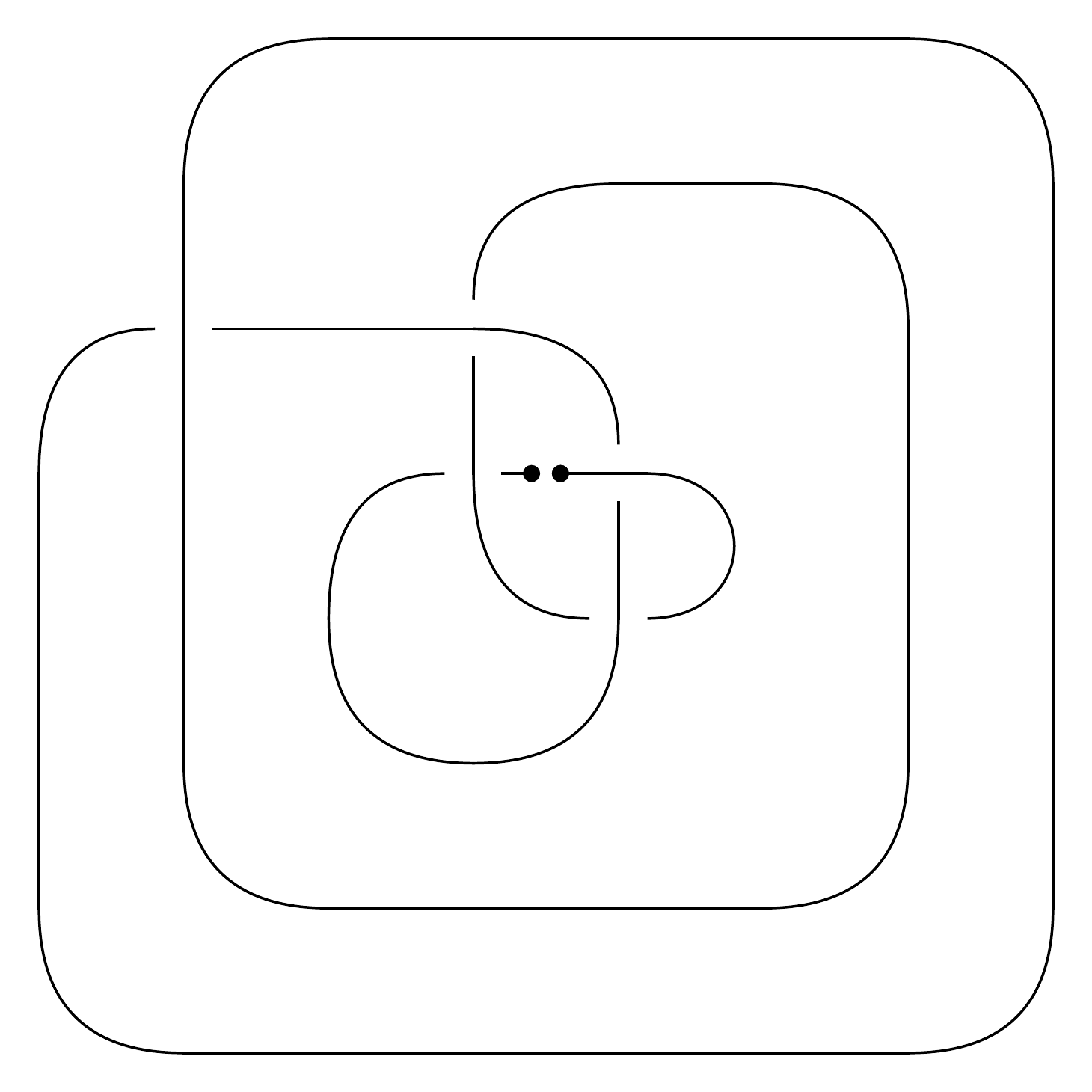}\\
\textcolor{black}{$5_{857}$}
\vspace{1cm}
\end{minipage}
\begin{minipage}[t]{.25\linewidth}
\centering
\includegraphics[width=0.9\textwidth,height=3.5cm,keepaspectratio]{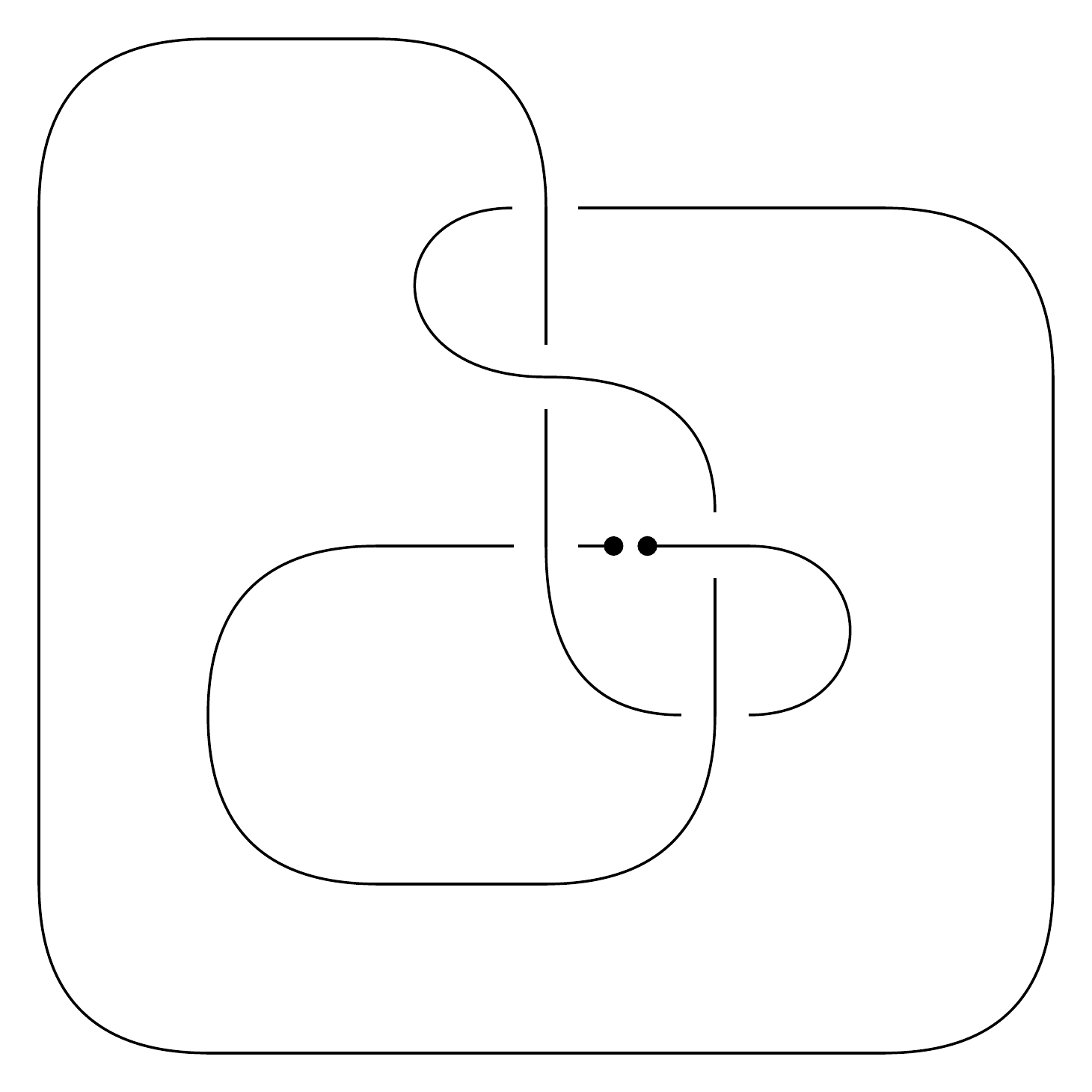}\\
\textcolor{black}{$5_{858}$}
\vspace{1cm}
\end{minipage}
\begin{minipage}[t]{.25\linewidth}
\centering
\includegraphics[width=0.9\textwidth,height=3.5cm,keepaspectratio]{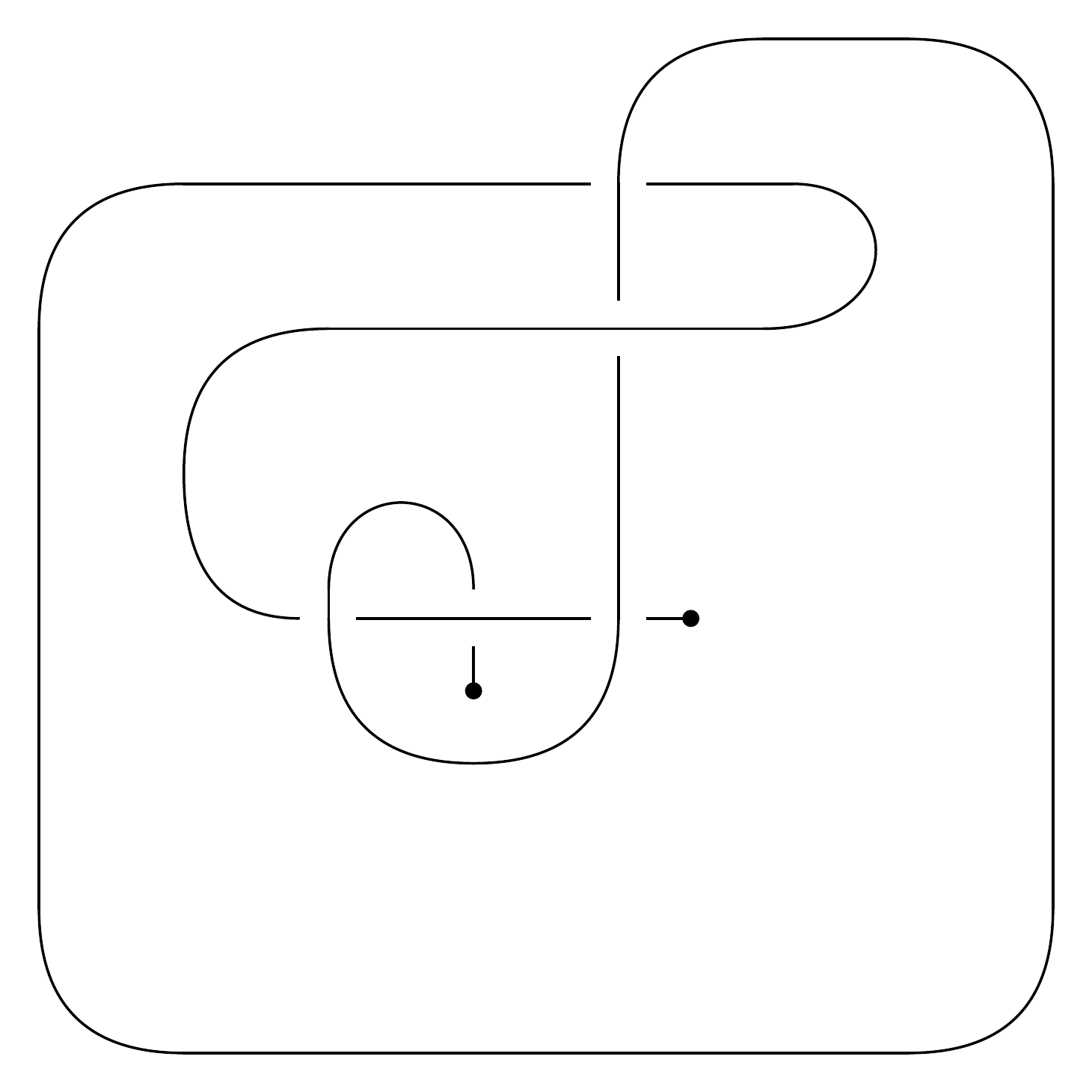}\\
\textcolor{black}{$5_{859}$}
\vspace{1cm}
\end{minipage}
\begin{minipage}[t]{.25\linewidth}
\centering
\includegraphics[width=0.9\textwidth,height=3.5cm,keepaspectratio]{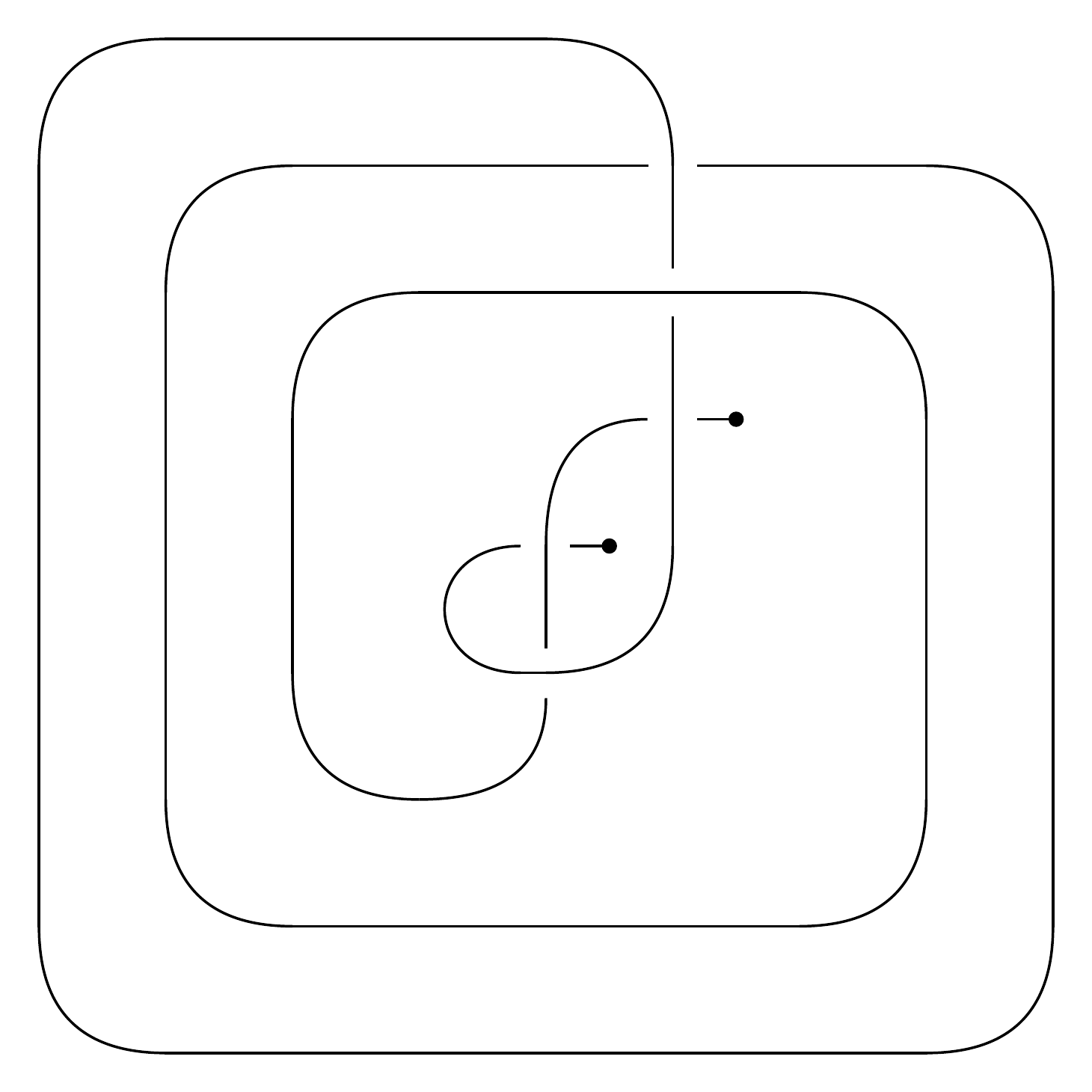}\\
\textcolor{black}{$5_{860}$}
\vspace{1cm}
\end{minipage}
\begin{minipage}[t]{.25\linewidth}
\centering
\includegraphics[width=0.9\textwidth,height=3.5cm,keepaspectratio]{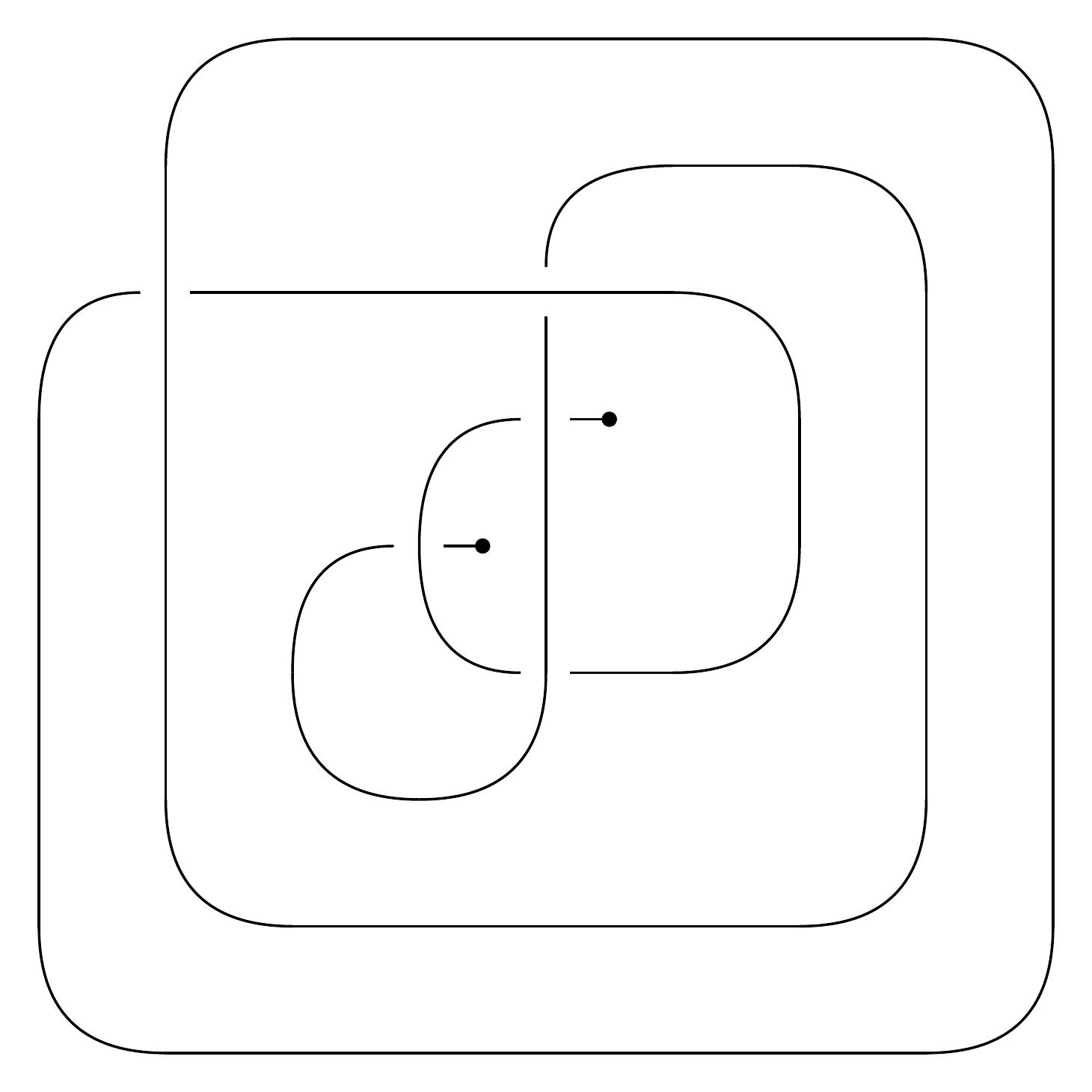}\\
\textcolor{black}{$5_{861}$}
\vspace{1cm}
\end{minipage}
\begin{minipage}[t]{.25\linewidth}
\centering
\includegraphics[width=0.9\textwidth,height=3.5cm,keepaspectratio]{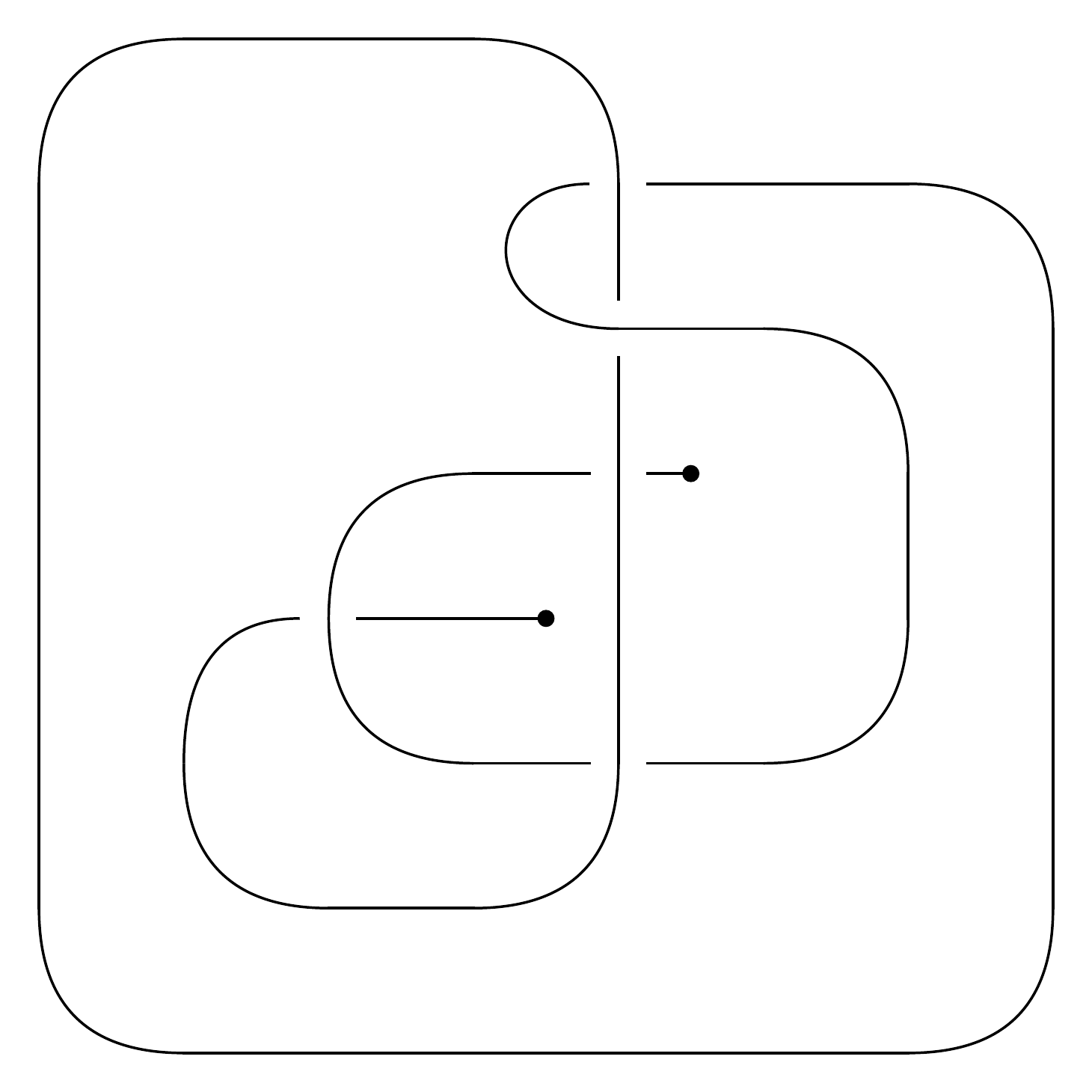}\\
\textcolor{black}{$5_{862}$}
\vspace{1cm}
\end{minipage}
\begin{minipage}[t]{.25\linewidth}
\centering
\includegraphics[width=0.9\textwidth,height=3.5cm,keepaspectratio]{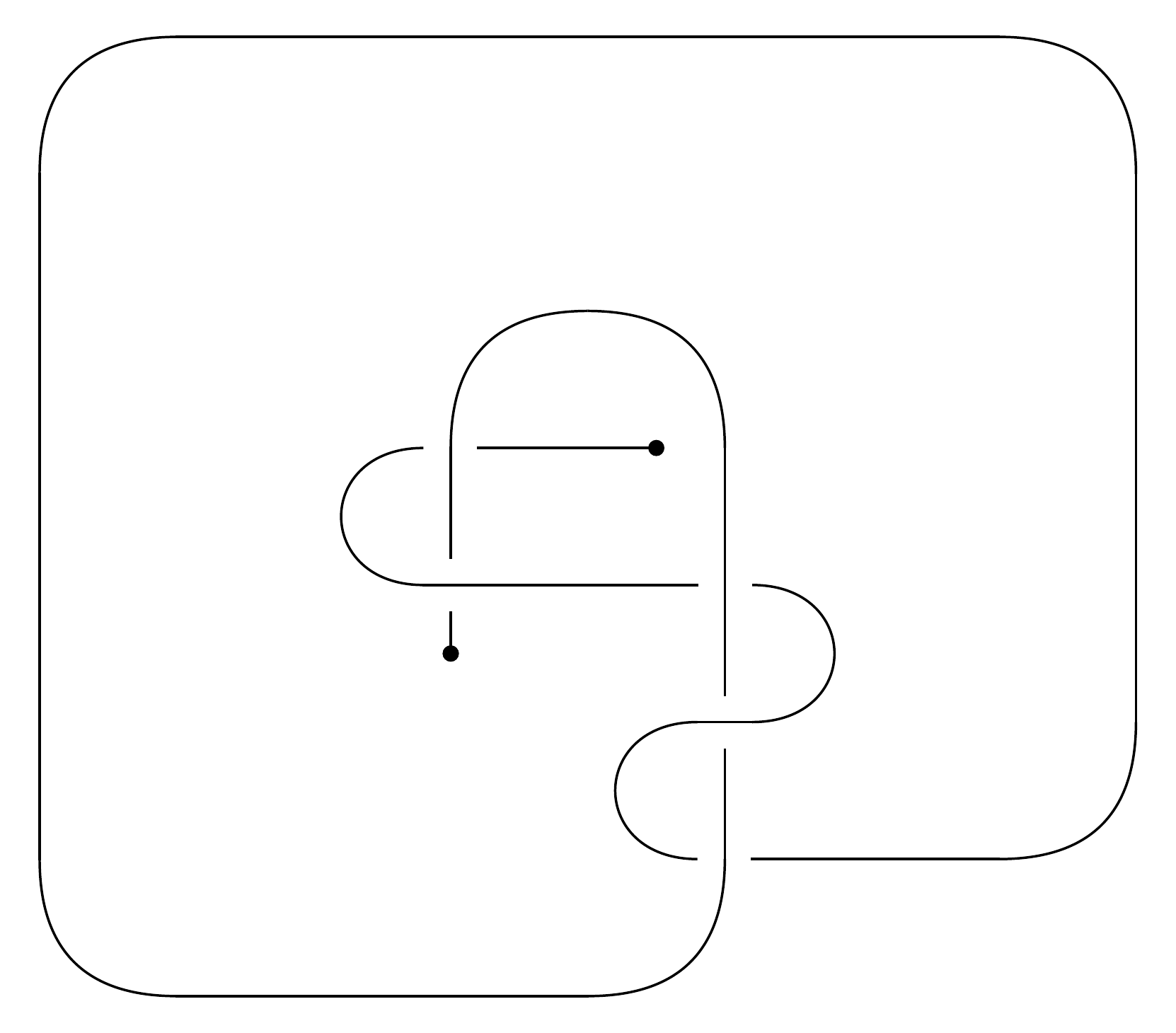}\\
\textcolor{black}{$5_{863}$}
\vspace{1cm}
\end{minipage}
\begin{minipage}[t]{.25\linewidth}
\centering
\includegraphics[width=0.9\textwidth,height=3.5cm,keepaspectratio]{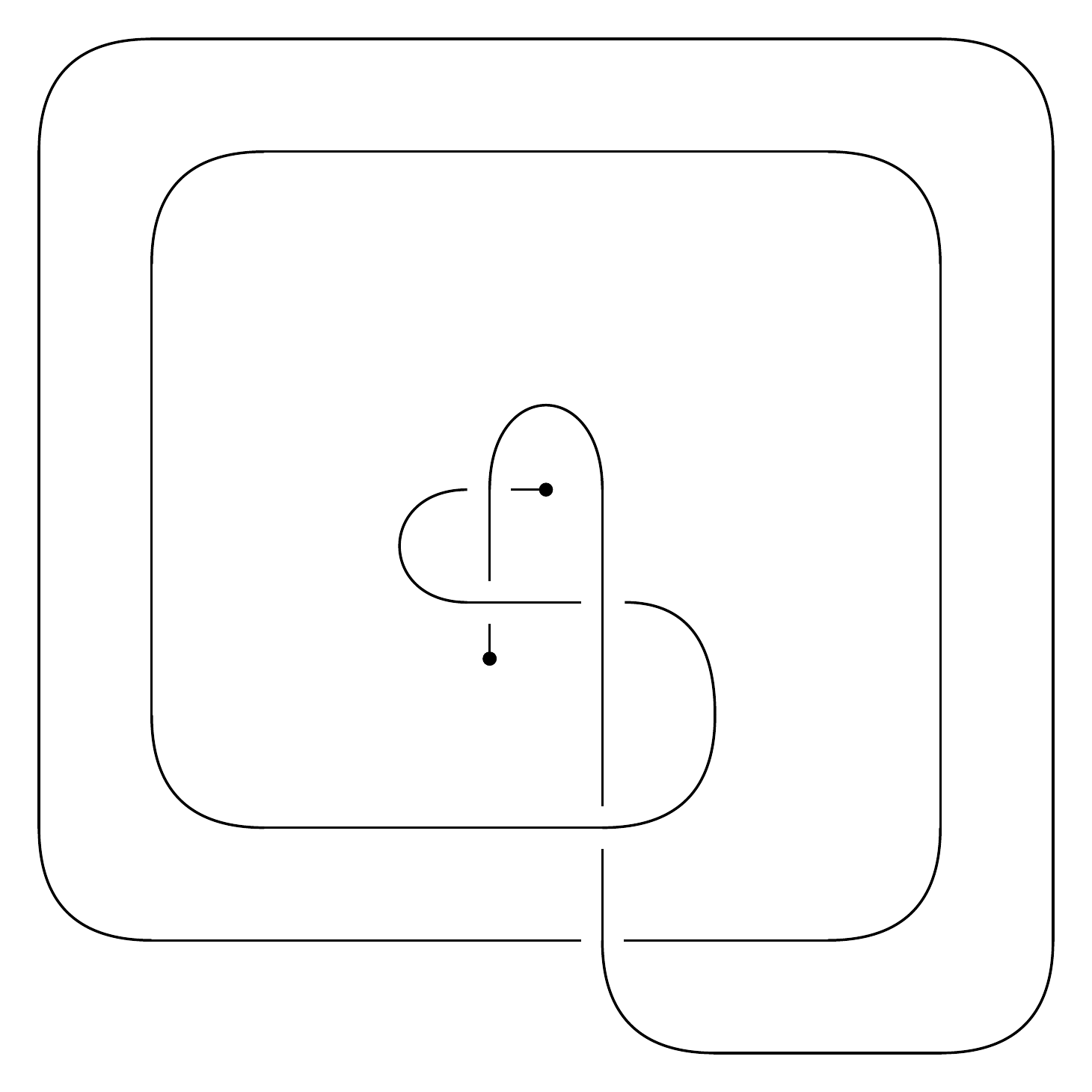}\\
\textcolor{black}{$5_{864}$}
\vspace{1cm}
\end{minipage}
\begin{minipage}[t]{.25\linewidth}
\centering
\includegraphics[width=0.9\textwidth,height=3.5cm,keepaspectratio]{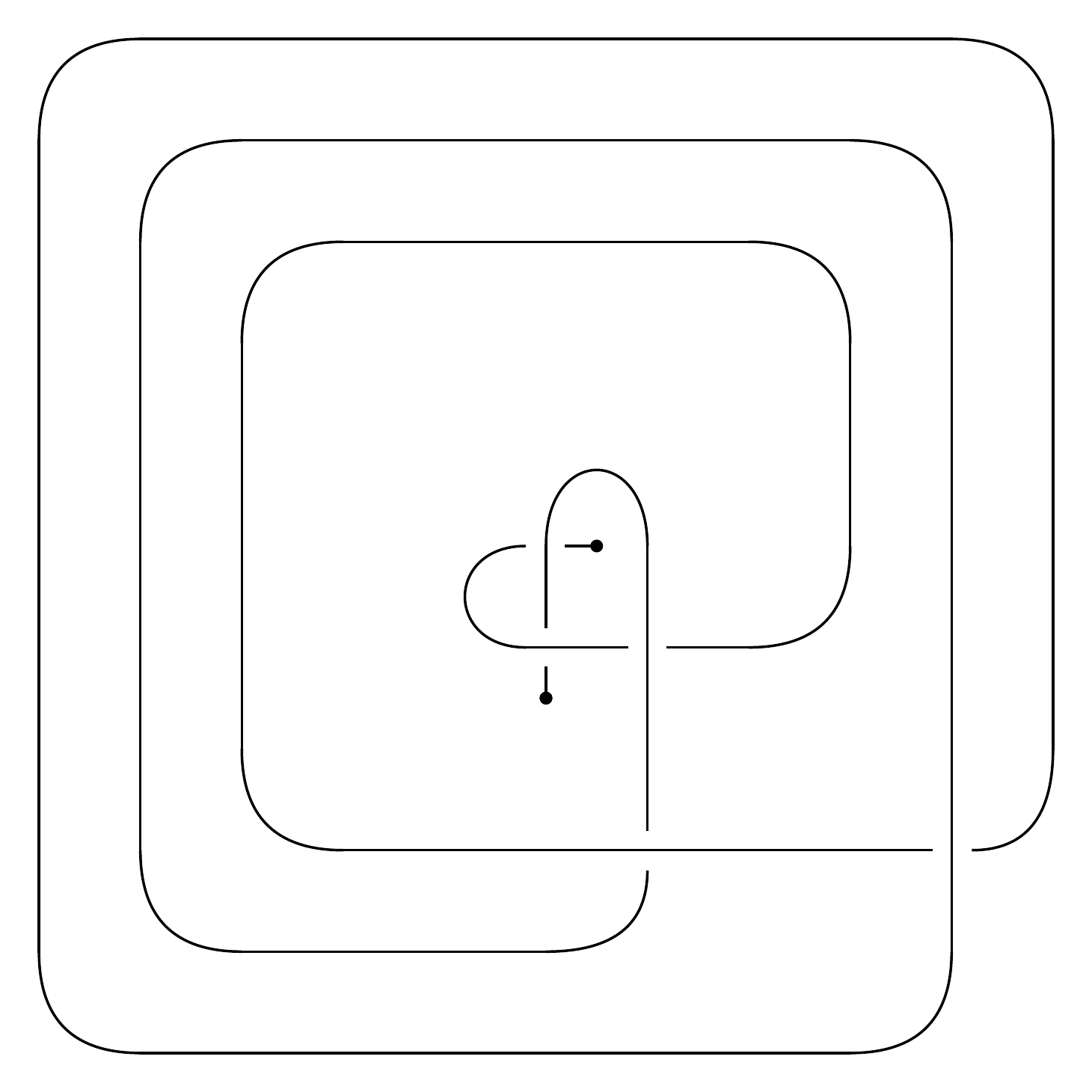}\\
\textcolor{black}{$5_{865}$}
\vspace{1cm}
\end{minipage}
\begin{minipage}[t]{.25\linewidth}
\centering
\includegraphics[width=0.9\textwidth,height=3.5cm,keepaspectratio]{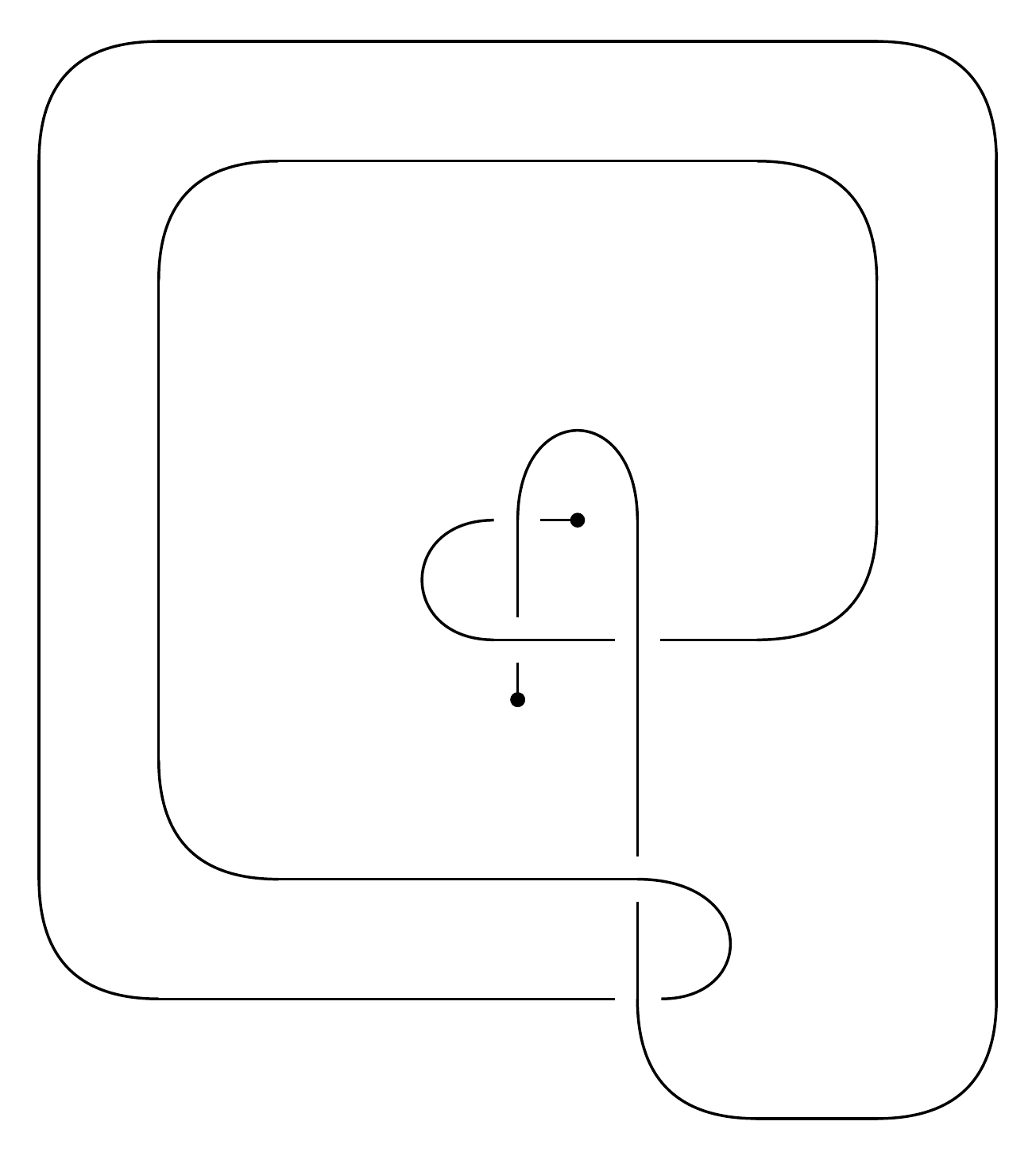}\\
\textcolor{black}{$5_{866}$}
\vspace{1cm}
\end{minipage}
\begin{minipage}[t]{.25\linewidth}
\centering
\includegraphics[width=0.9\textwidth,height=3.5cm,keepaspectratio]{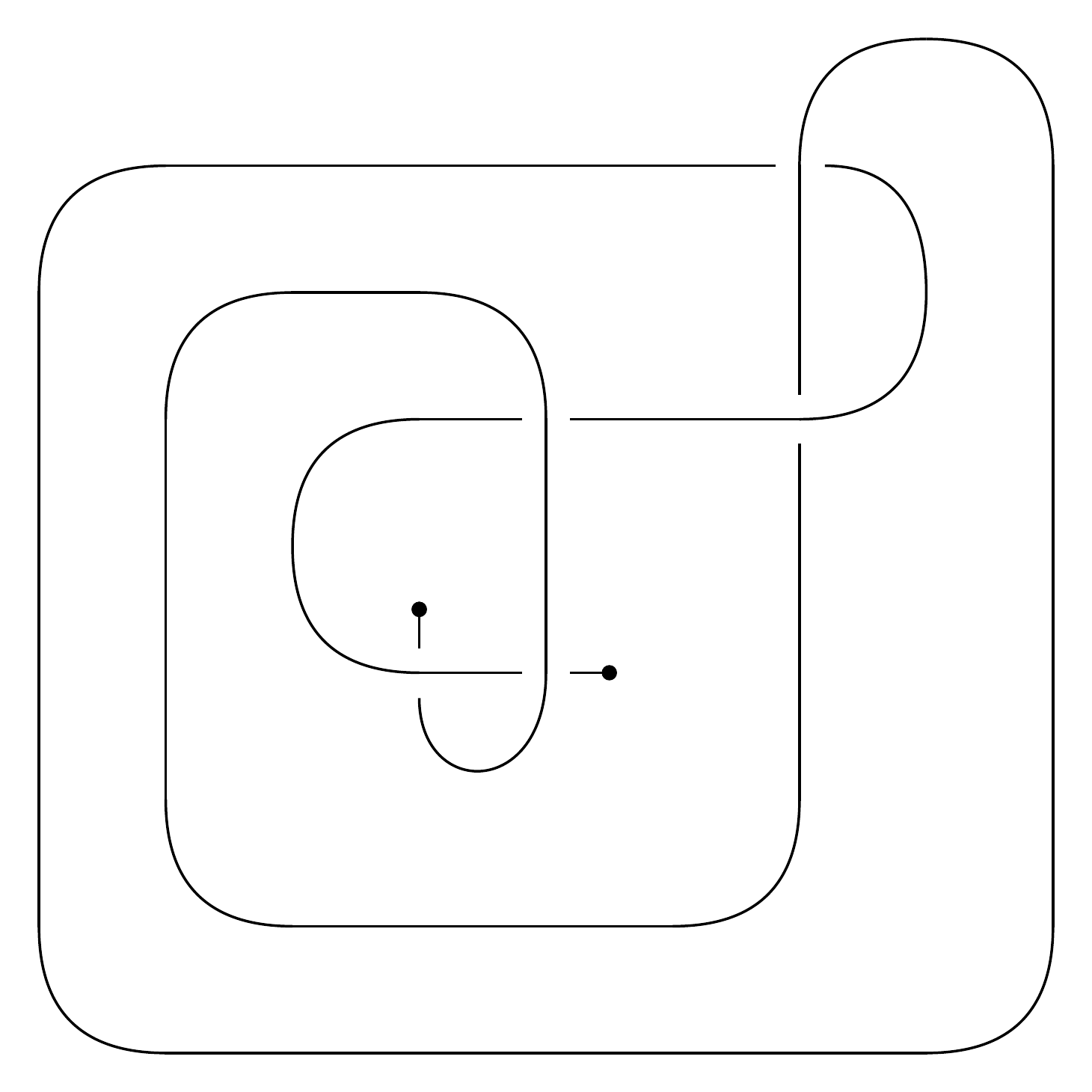}\\
\textcolor{black}{$5_{867}$}
\vspace{1cm}
\end{minipage}
\begin{minipage}[t]{.25\linewidth}
\centering
\includegraphics[width=0.9\textwidth,height=3.5cm,keepaspectratio]{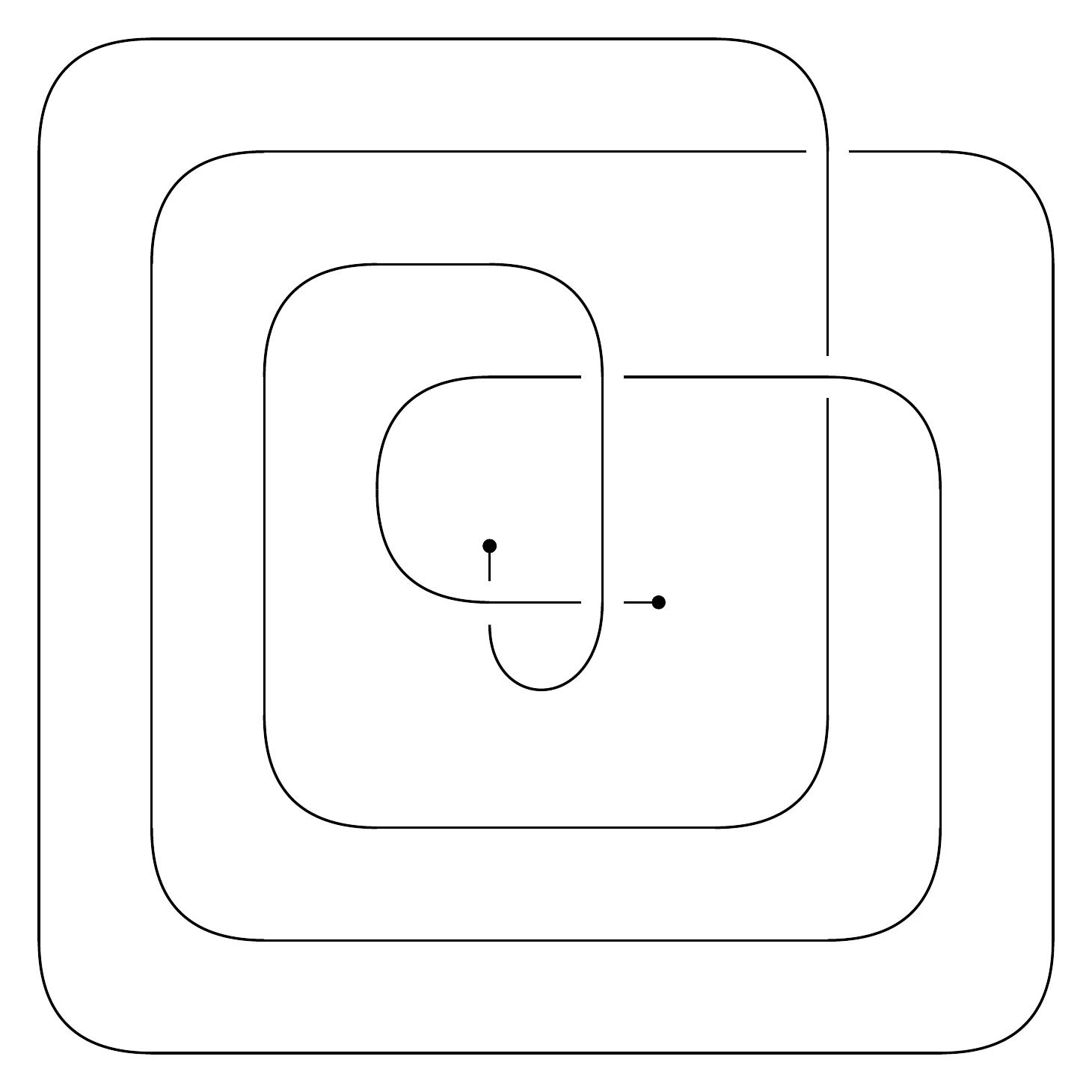}\\
\textcolor{black}{$5_{868}$}
\vspace{1cm}
\end{minipage}
\begin{minipage}[t]{.25\linewidth}
\centering
\includegraphics[width=0.9\textwidth,height=3.5cm,keepaspectratio]{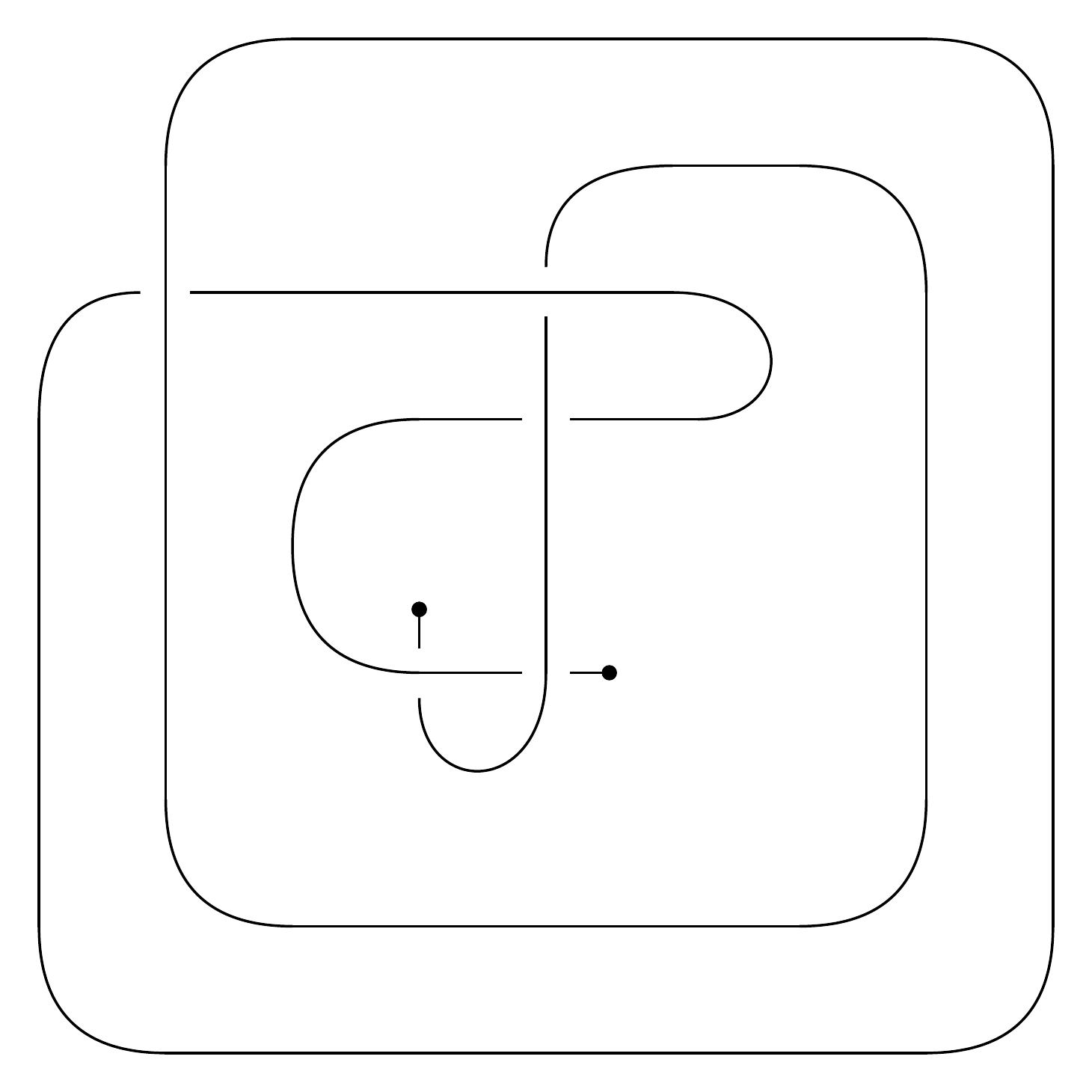}\\
\textcolor{black}{$5_{869}$}
\vspace{1cm}
\end{minipage}
\begin{minipage}[t]{.25\linewidth}
\centering
\includegraphics[width=0.9\textwidth,height=3.5cm,keepaspectratio]{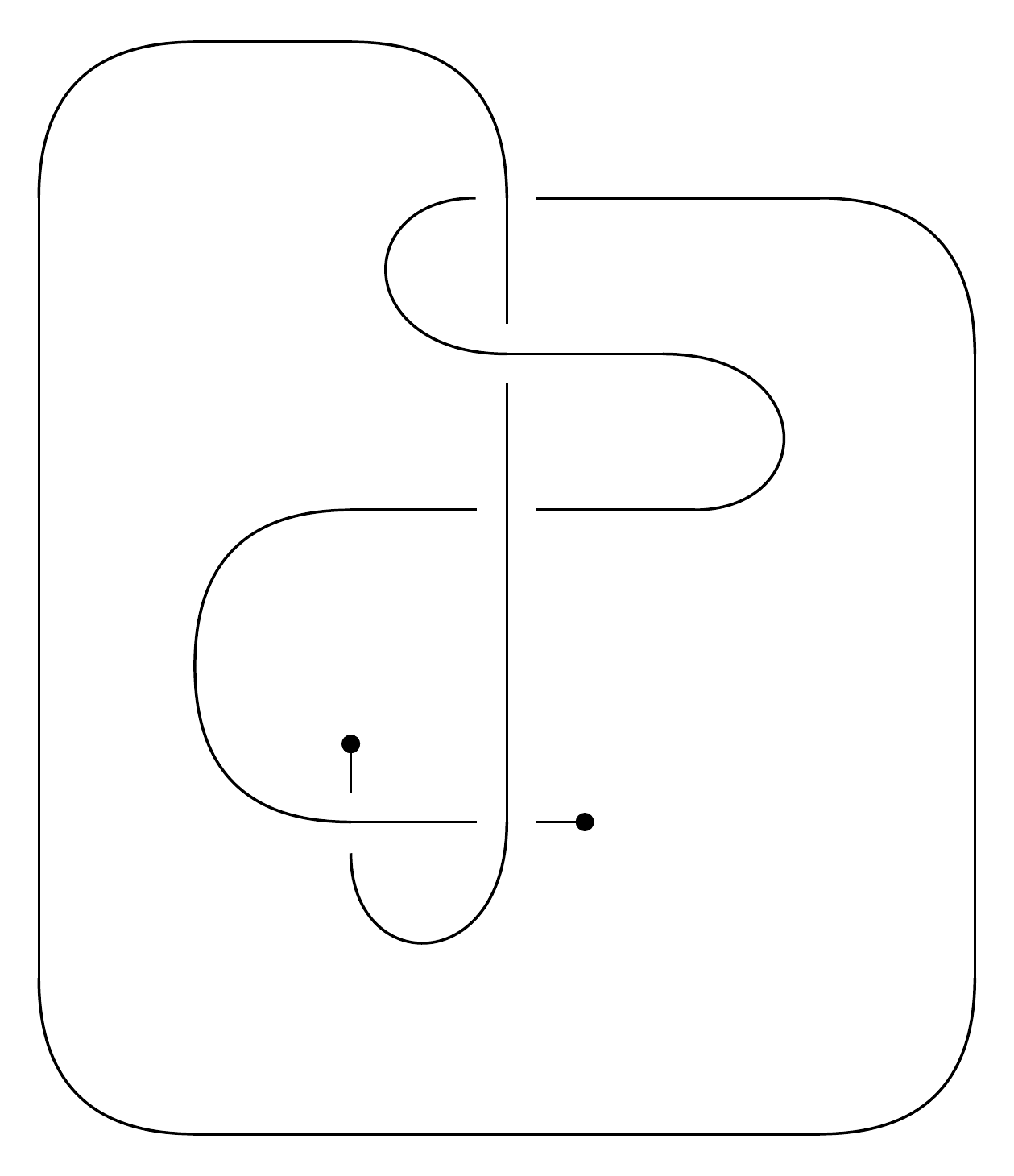}\\
\textcolor{black}{$5_{870}$}
\vspace{1cm}
\end{minipage}
\begin{minipage}[t]{.25\linewidth}
\centering
\includegraphics[width=0.9\textwidth,height=3.5cm,keepaspectratio]{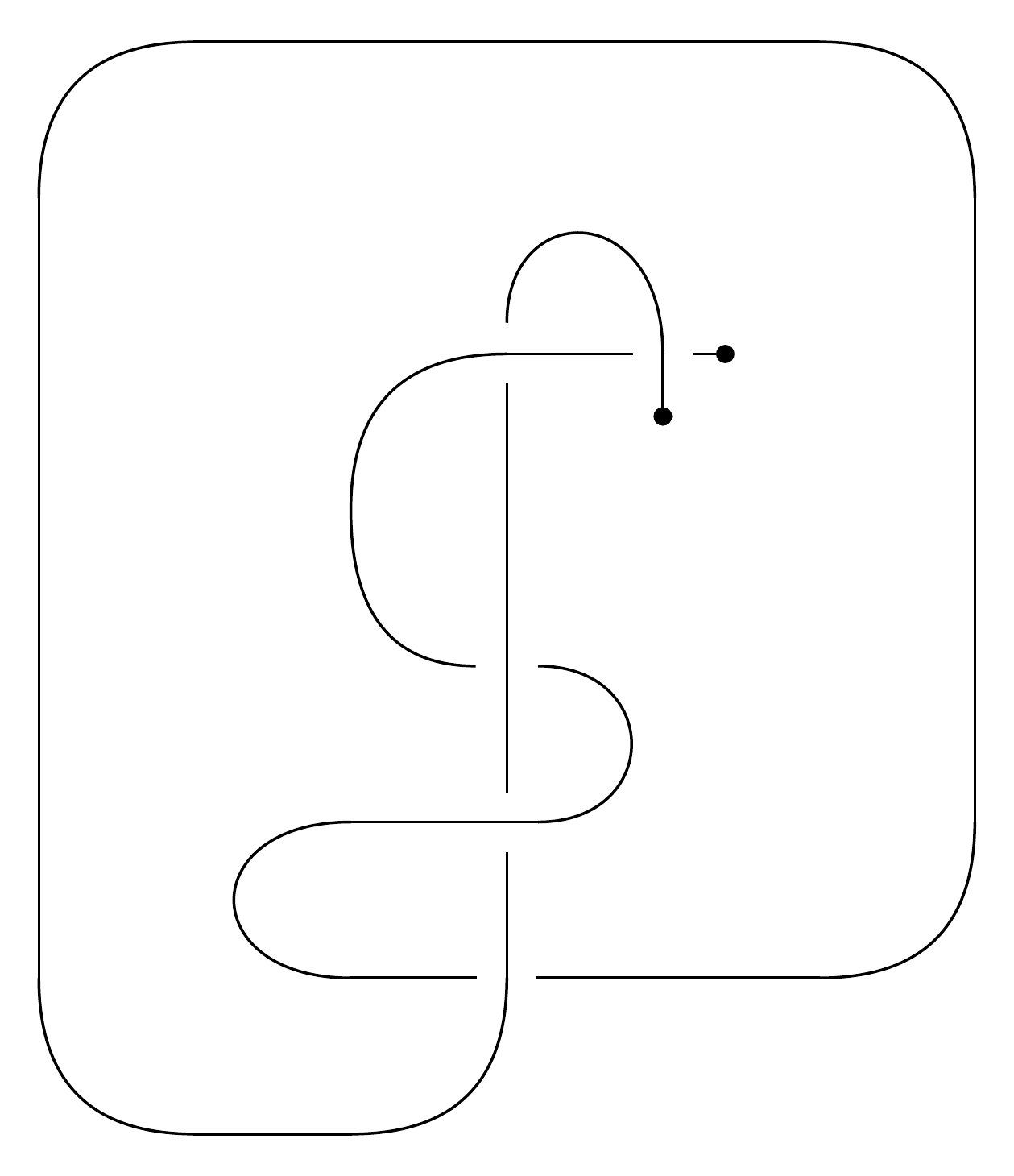}\\
\textcolor{black}{$5_{871}$}
\vspace{1cm}
\end{minipage}
\begin{minipage}[t]{.25\linewidth}
\centering
\includegraphics[width=0.9\textwidth,height=3.5cm,keepaspectratio]{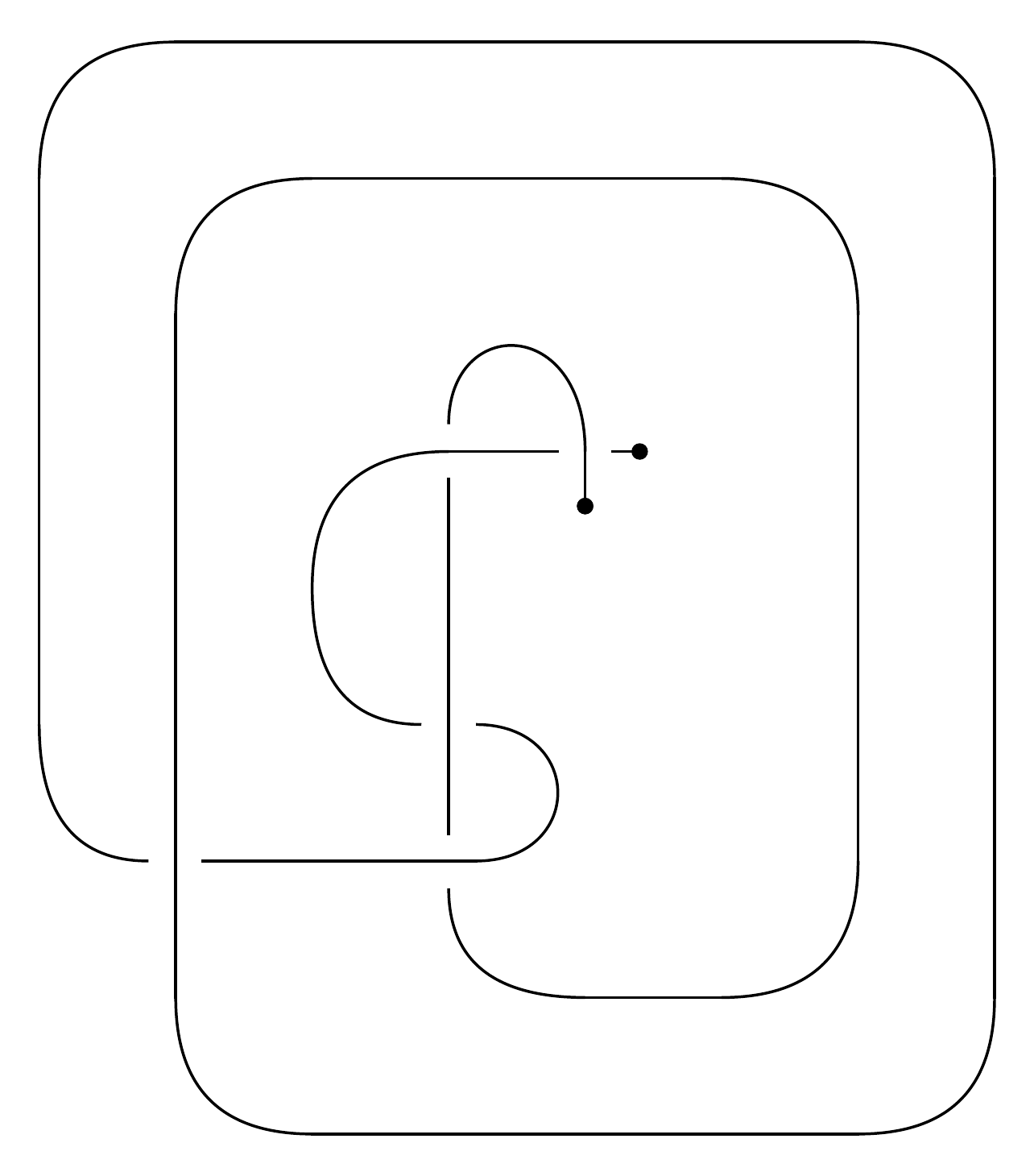}\\
\textcolor{black}{$5_{872}$}
\vspace{1cm}
\end{minipage}
\begin{minipage}[t]{.25\linewidth}
\centering
\includegraphics[width=0.9\textwidth,height=3.5cm,keepaspectratio]{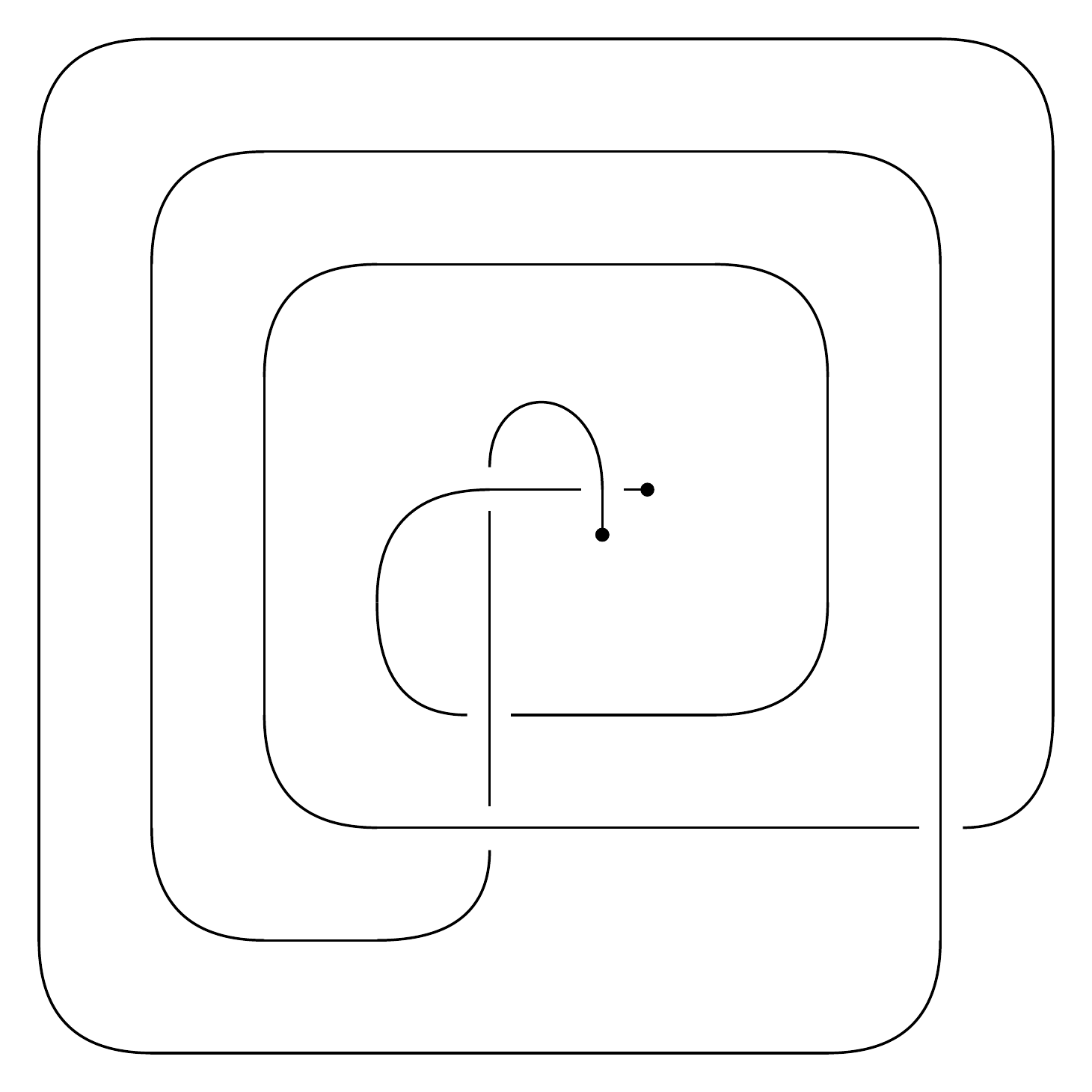}\\
\textcolor{black}{$5_{873}$}
\vspace{1cm}
\end{minipage}
\begin{minipage}[t]{.25\linewidth}
\centering
\includegraphics[width=0.9\textwidth,height=3.5cm,keepaspectratio]{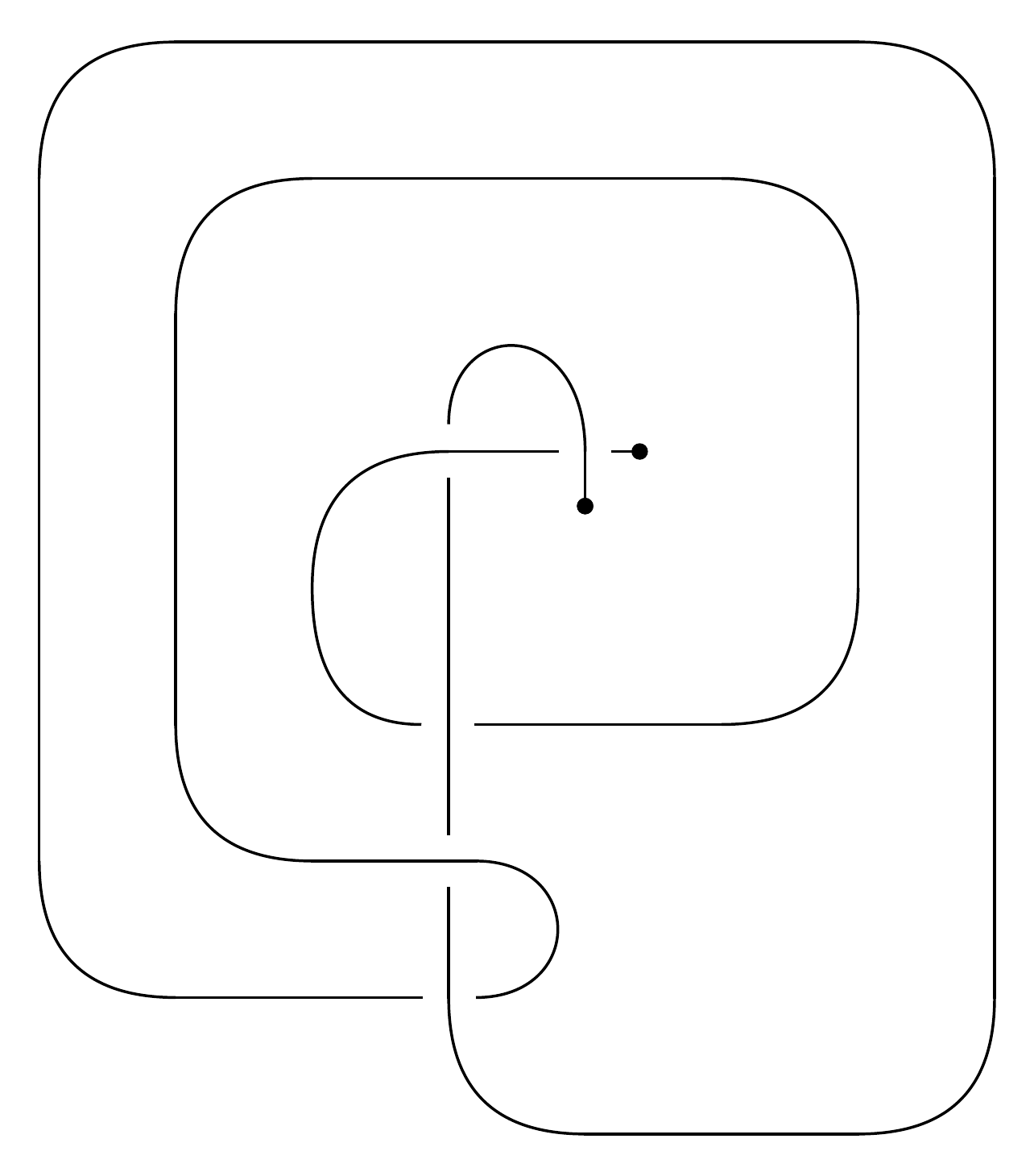}\\
\textcolor{black}{$5_{874}$}
\vspace{1cm}
\end{minipage}
\begin{minipage}[t]{.25\linewidth}
\centering
\includegraphics[width=0.9\textwidth,height=3.5cm,keepaspectratio]{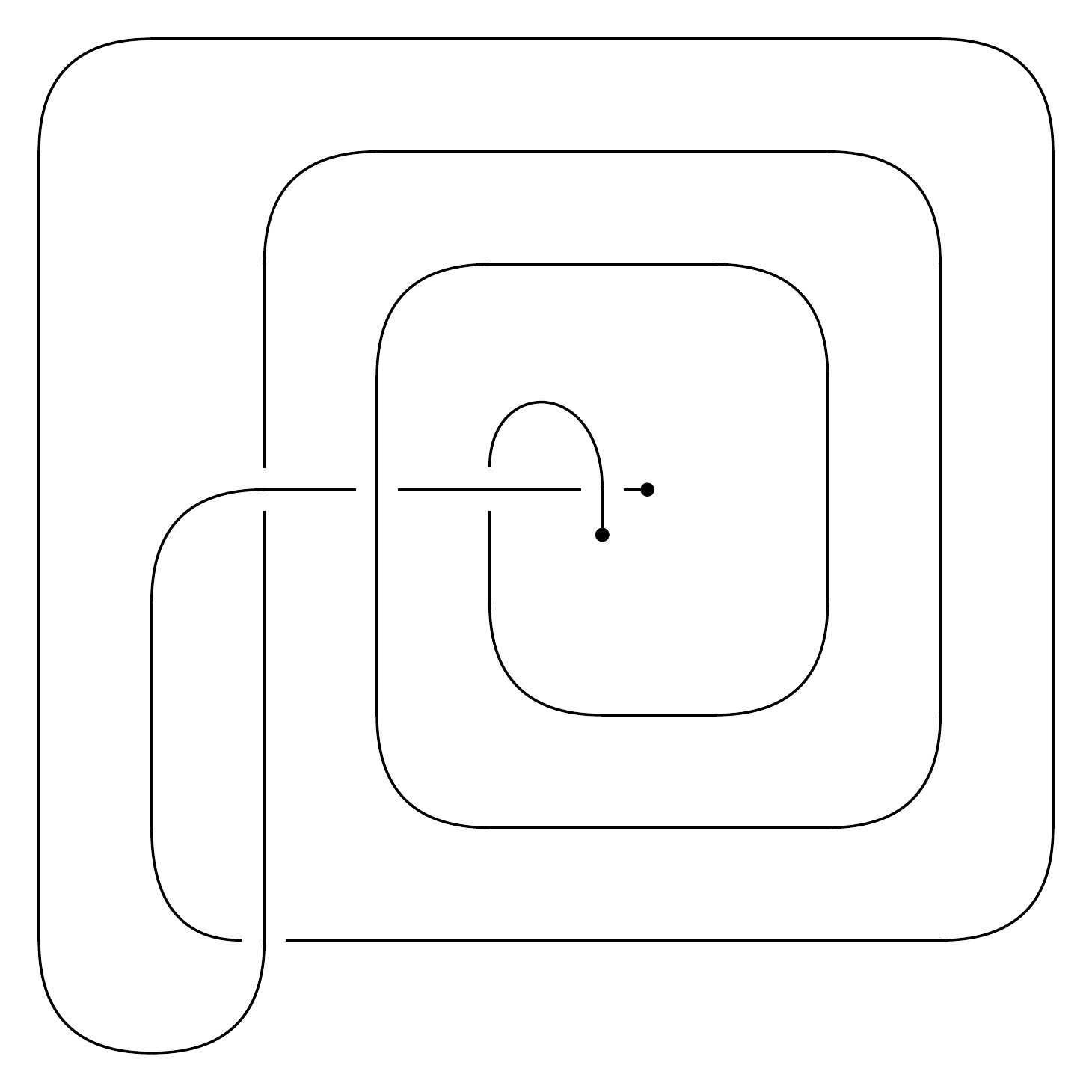}\\
\textcolor{black}{$5_{875}$}
\vspace{1cm}
\end{minipage}
\begin{minipage}[t]{.25\linewidth}
\centering
\includegraphics[width=0.9\textwidth,height=3.5cm,keepaspectratio]{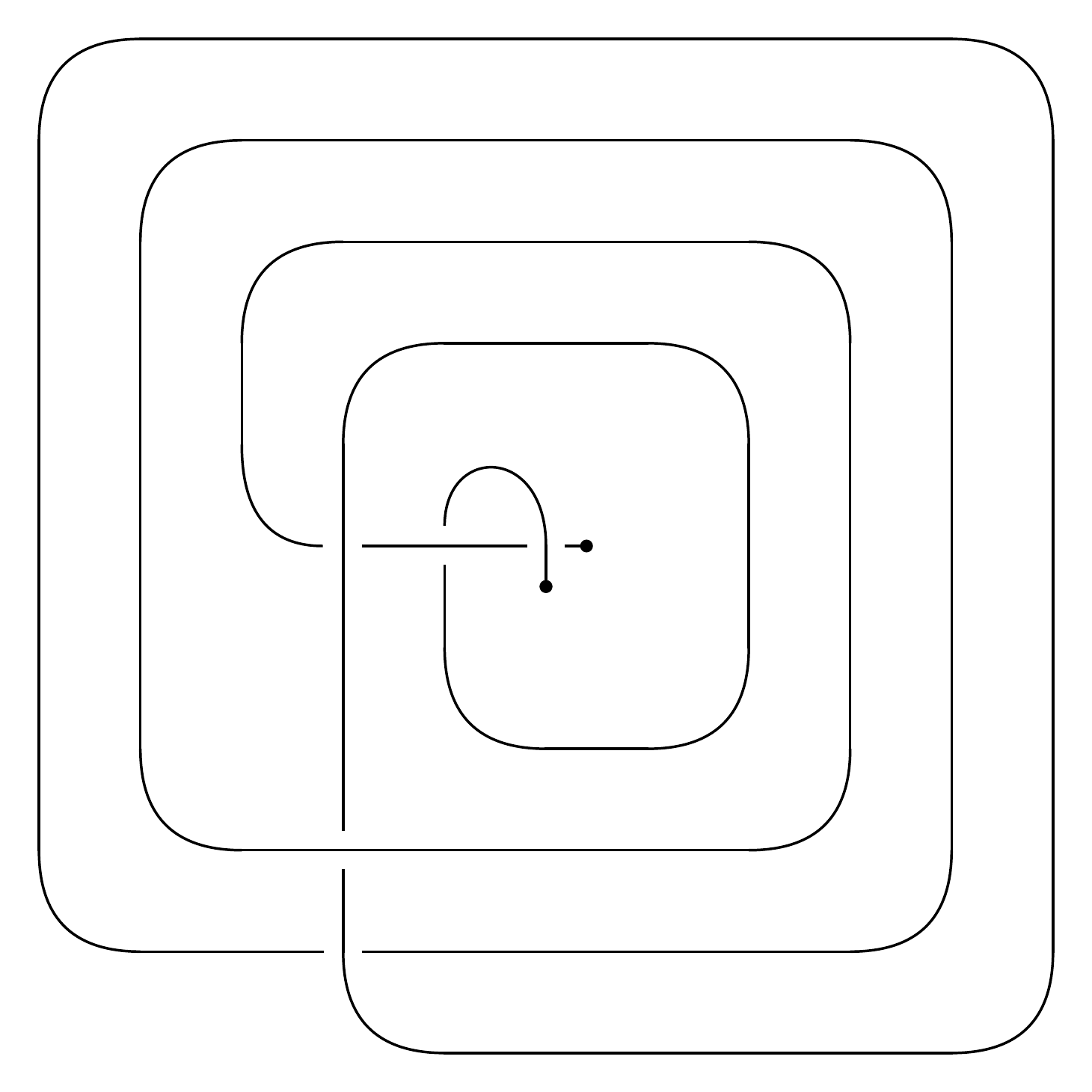}\\
\textcolor{black}{$5_{876}$}
\vspace{1cm}
\end{minipage}
\begin{minipage}[t]{.25\linewidth}
\centering
\includegraphics[width=0.9\textwidth,height=3.5cm,keepaspectratio]{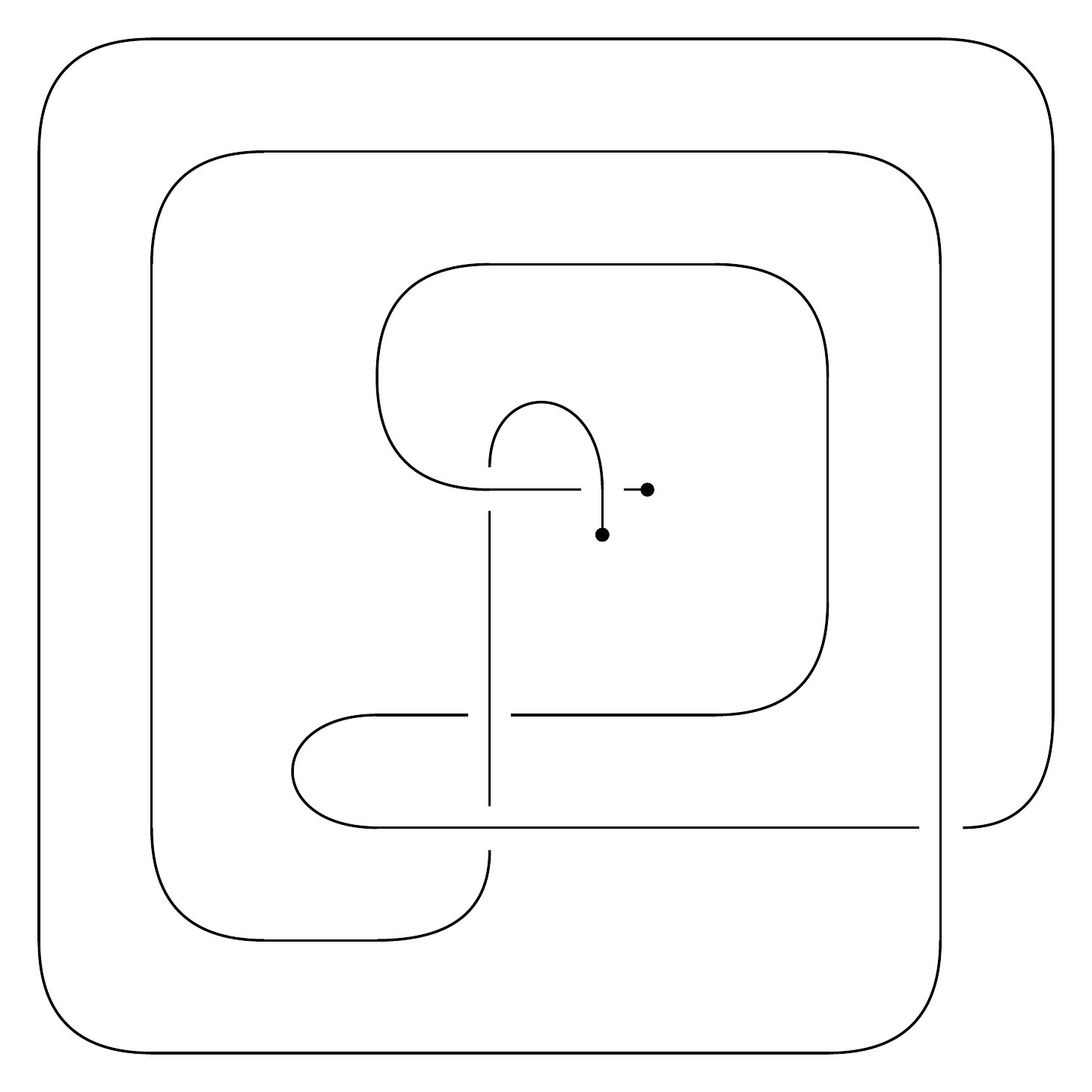}\\
\textcolor{black}{$5_{877}$}
\vspace{1cm}
\end{minipage}
\begin{minipage}[t]{.25\linewidth}
\centering
\includegraphics[width=0.9\textwidth,height=3.5cm,keepaspectratio]{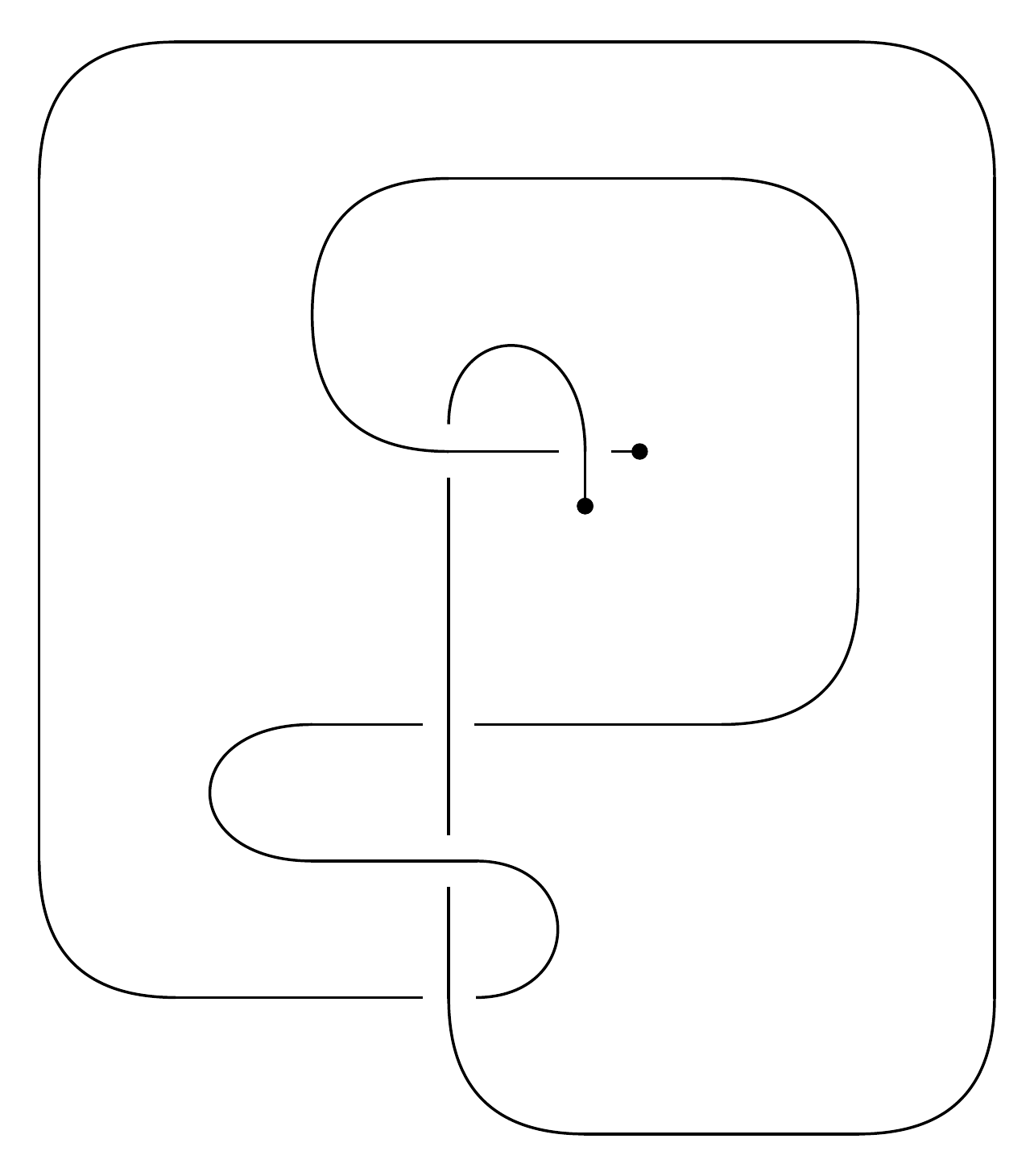}\\
\textcolor{black}{$5_{878}$}
\vspace{1cm}
\end{minipage}
\begin{minipage}[t]{.25\linewidth}
\centering
\includegraphics[width=0.9\textwidth,height=3.5cm,keepaspectratio]{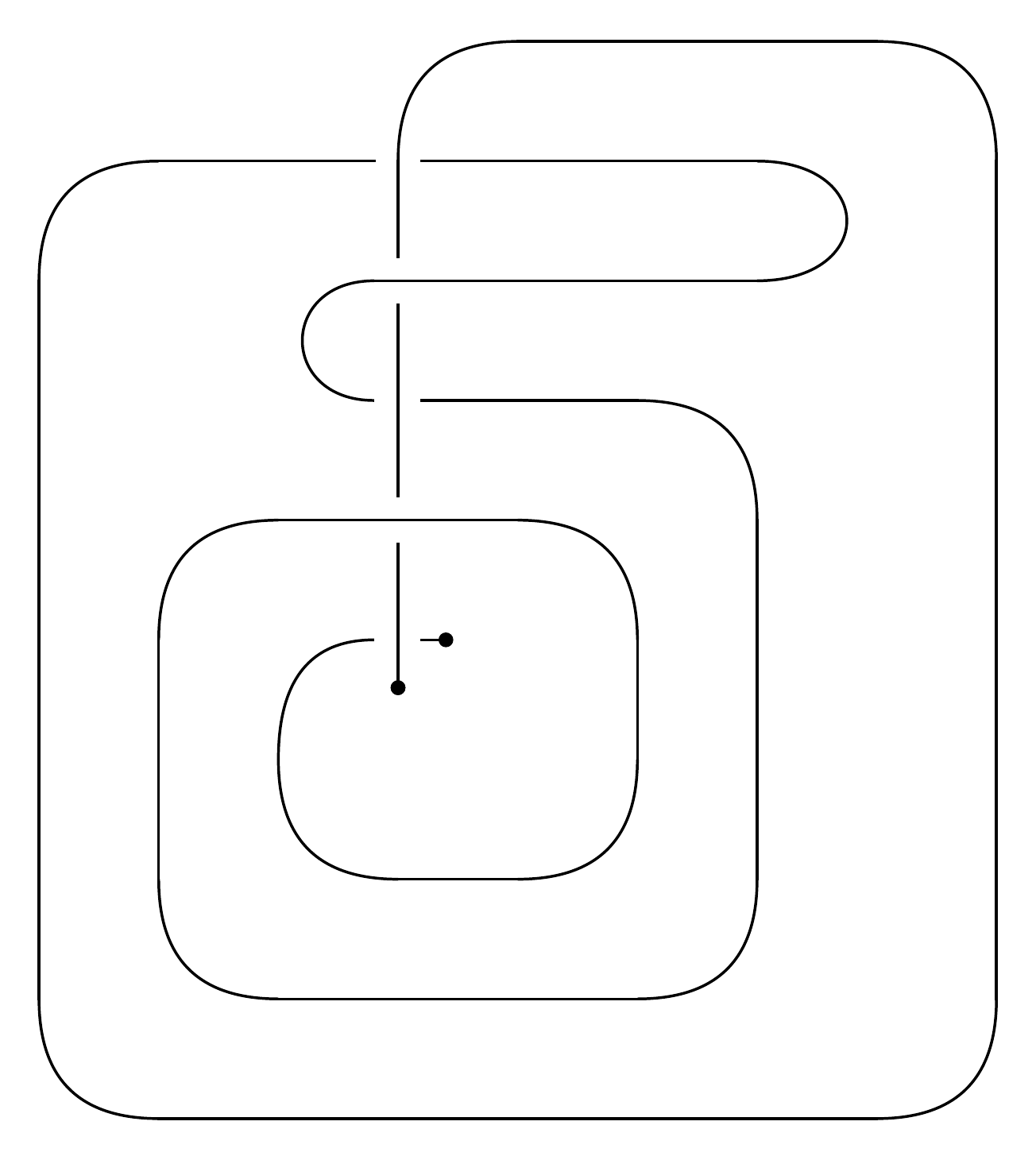}\\
\textcolor{black}{$5_{879}$}
\vspace{1cm}
\end{minipage}
\begin{minipage}[t]{.25\linewidth}
\centering
\includegraphics[width=0.9\textwidth,height=3.5cm,keepaspectratio]{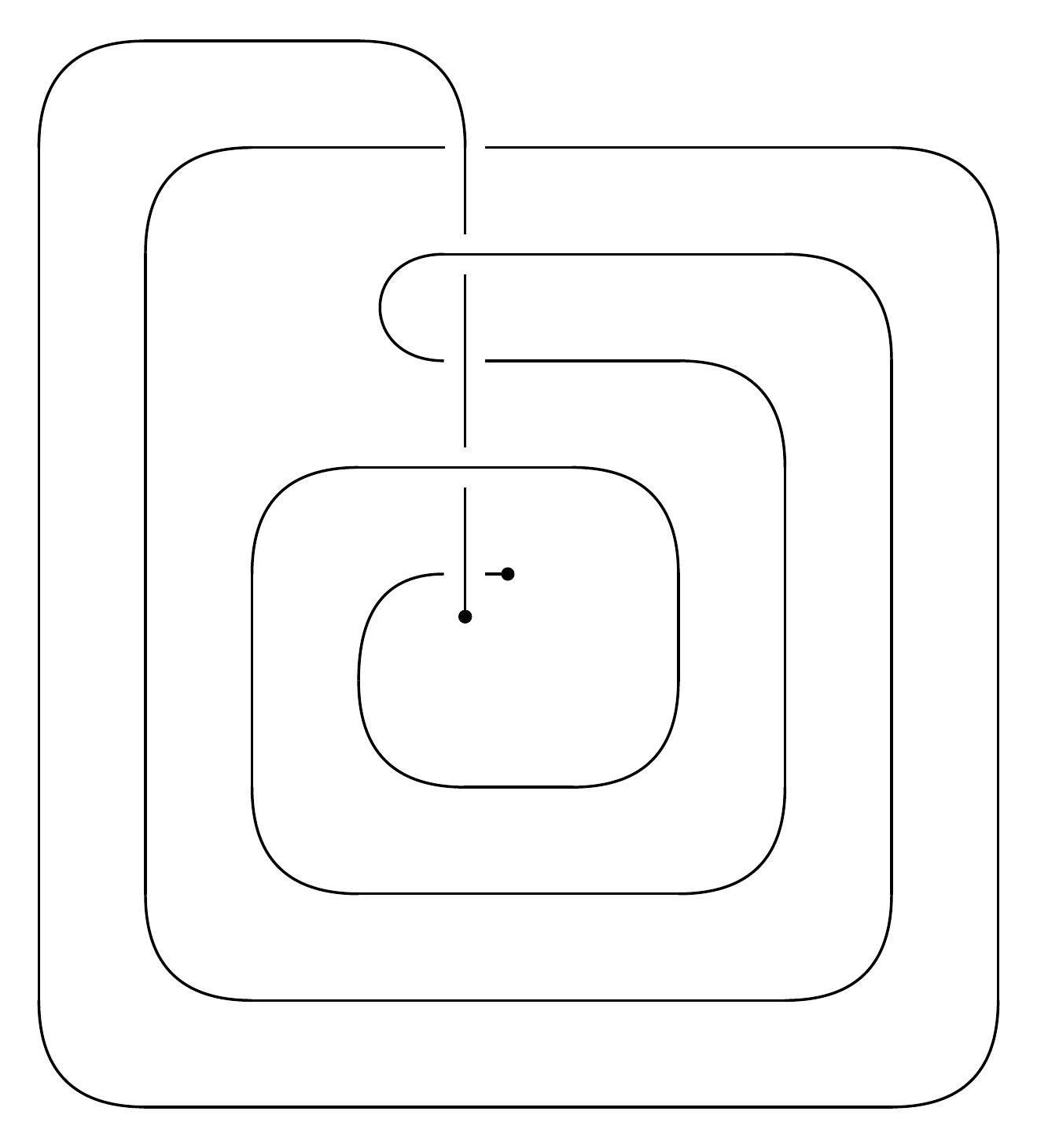}\\
\textcolor{black}{$5_{880}$}
\vspace{1cm}
\end{minipage}
\begin{minipage}[t]{.25\linewidth}
\centering
\includegraphics[width=0.9\textwidth,height=3.5cm,keepaspectratio]{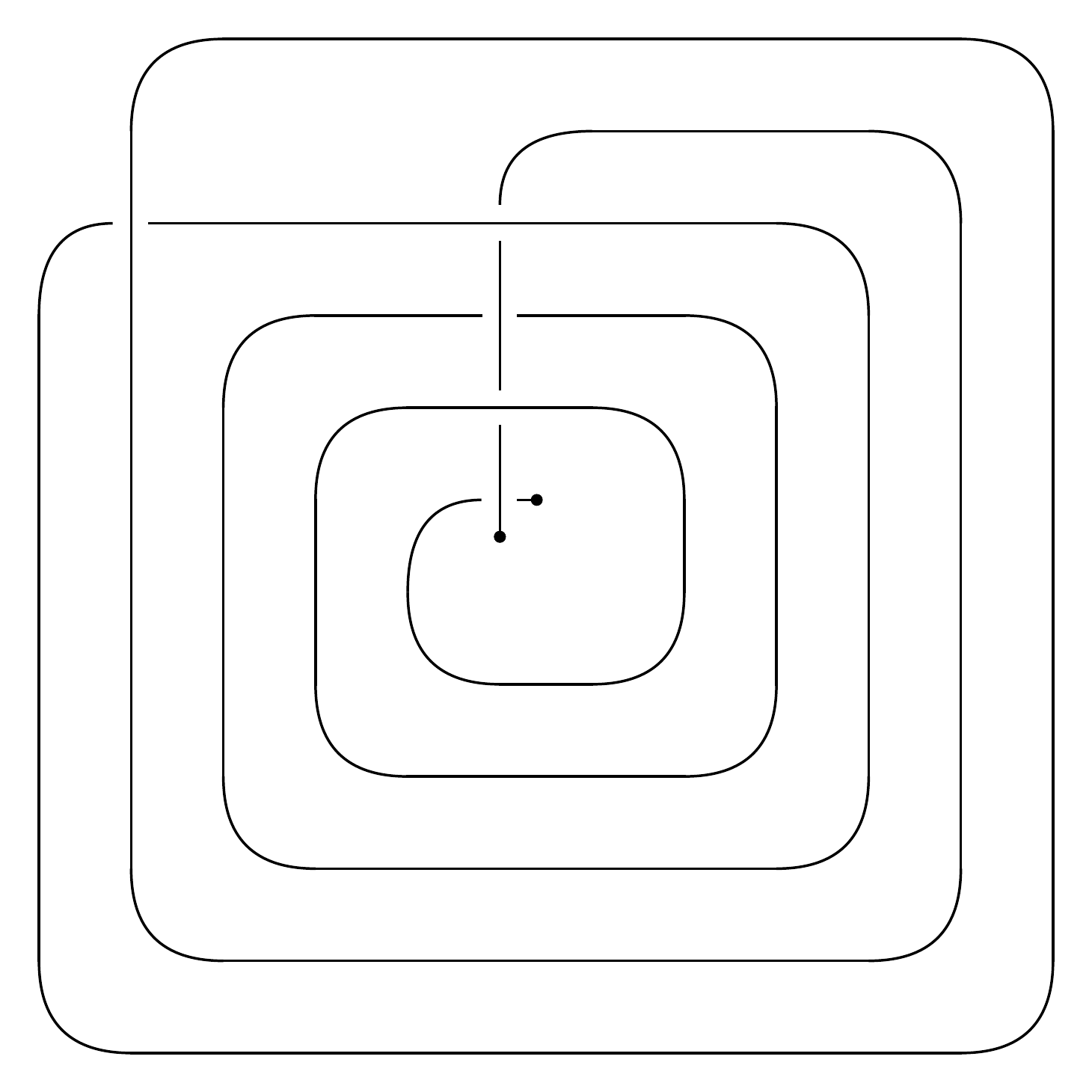}\\
\textcolor{black}{$5_{881}$}
\vspace{1cm}
\end{minipage}
\begin{minipage}[t]{.25\linewidth}
\centering
\includegraphics[width=0.9\textwidth,height=3.5cm,keepaspectratio]{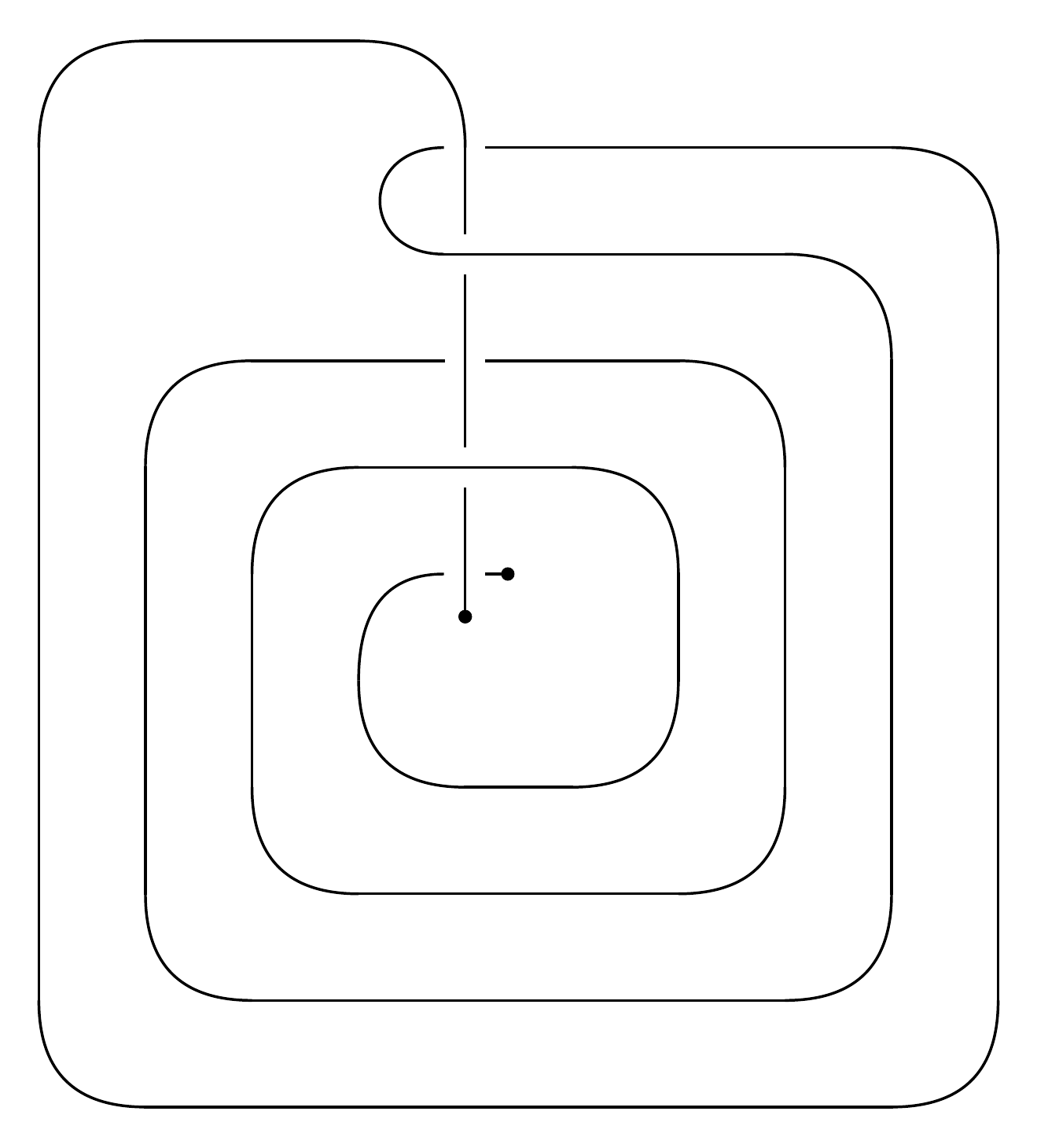}\\
\textcolor{black}{$5_{882}$}
\vspace{1cm}
\end{minipage}
\begin{minipage}[t]{.25\linewidth}
\centering
\includegraphics[width=0.9\textwidth,height=3.5cm,keepaspectratio]{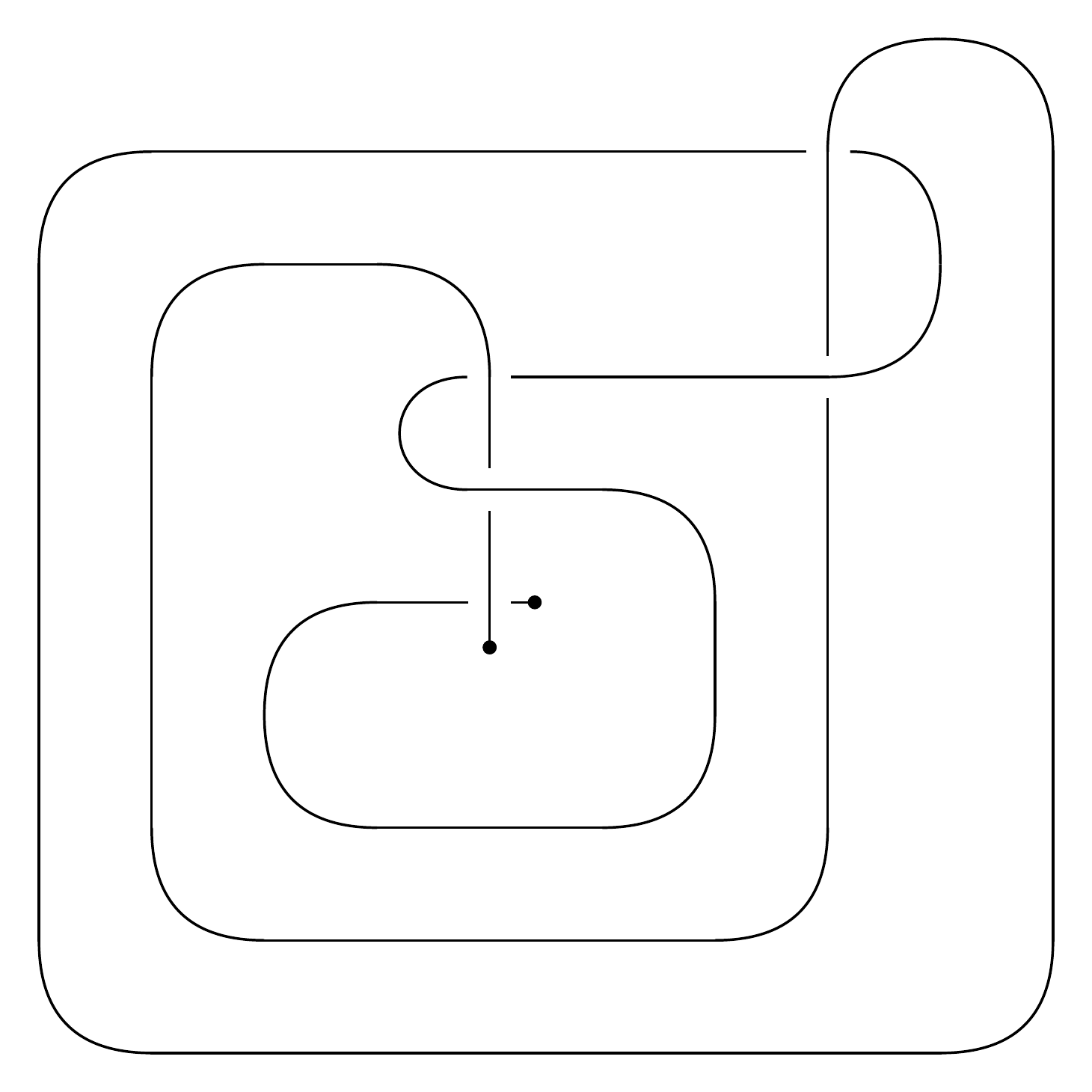}\\
\textcolor{black}{$5_{883}$}
\vspace{1cm}
\end{minipage}
\begin{minipage}[t]{.25\linewidth}
\centering
\includegraphics[width=0.9\textwidth,height=3.5cm,keepaspectratio]{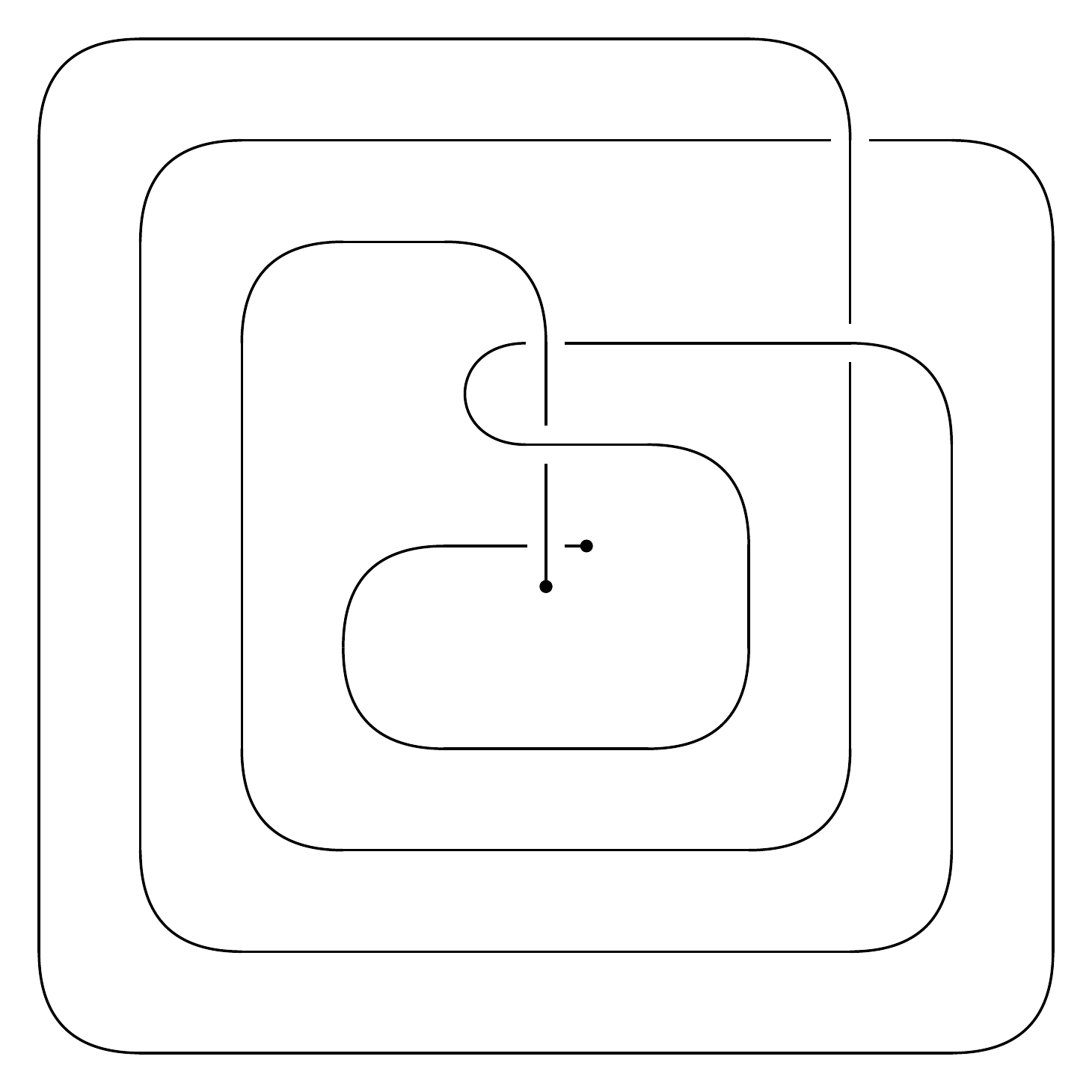}\\
\textcolor{black}{$5_{884}$}
\vspace{1cm}
\end{minipage}
\begin{minipage}[t]{.25\linewidth}
\centering
\includegraphics[width=0.9\textwidth,height=3.5cm,keepaspectratio]{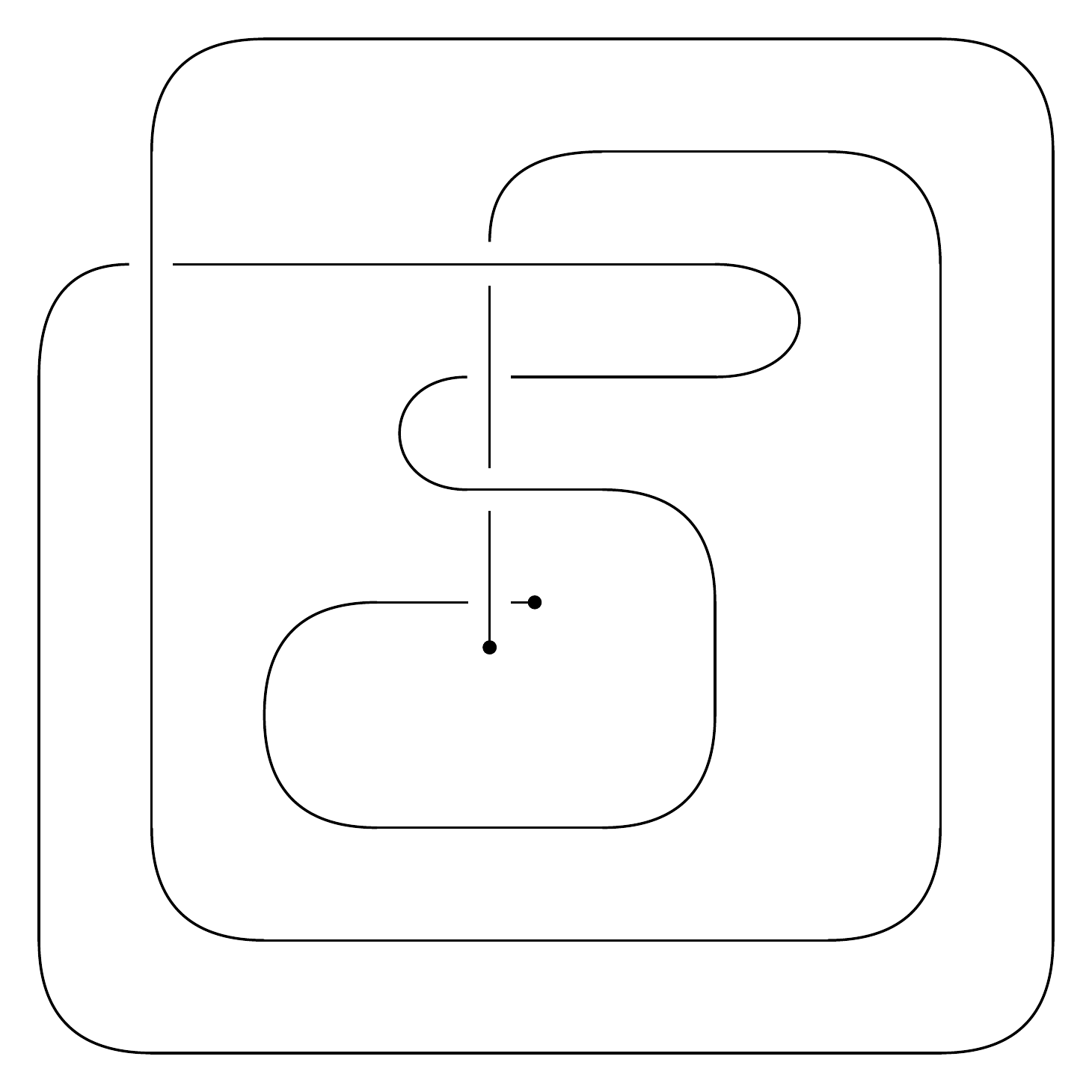}\\
\textcolor{black}{$5_{885}$}
\vspace{1cm}
\end{minipage}
\begin{minipage}[t]{.25\linewidth}
\centering
\includegraphics[width=0.9\textwidth,height=3.5cm,keepaspectratio]{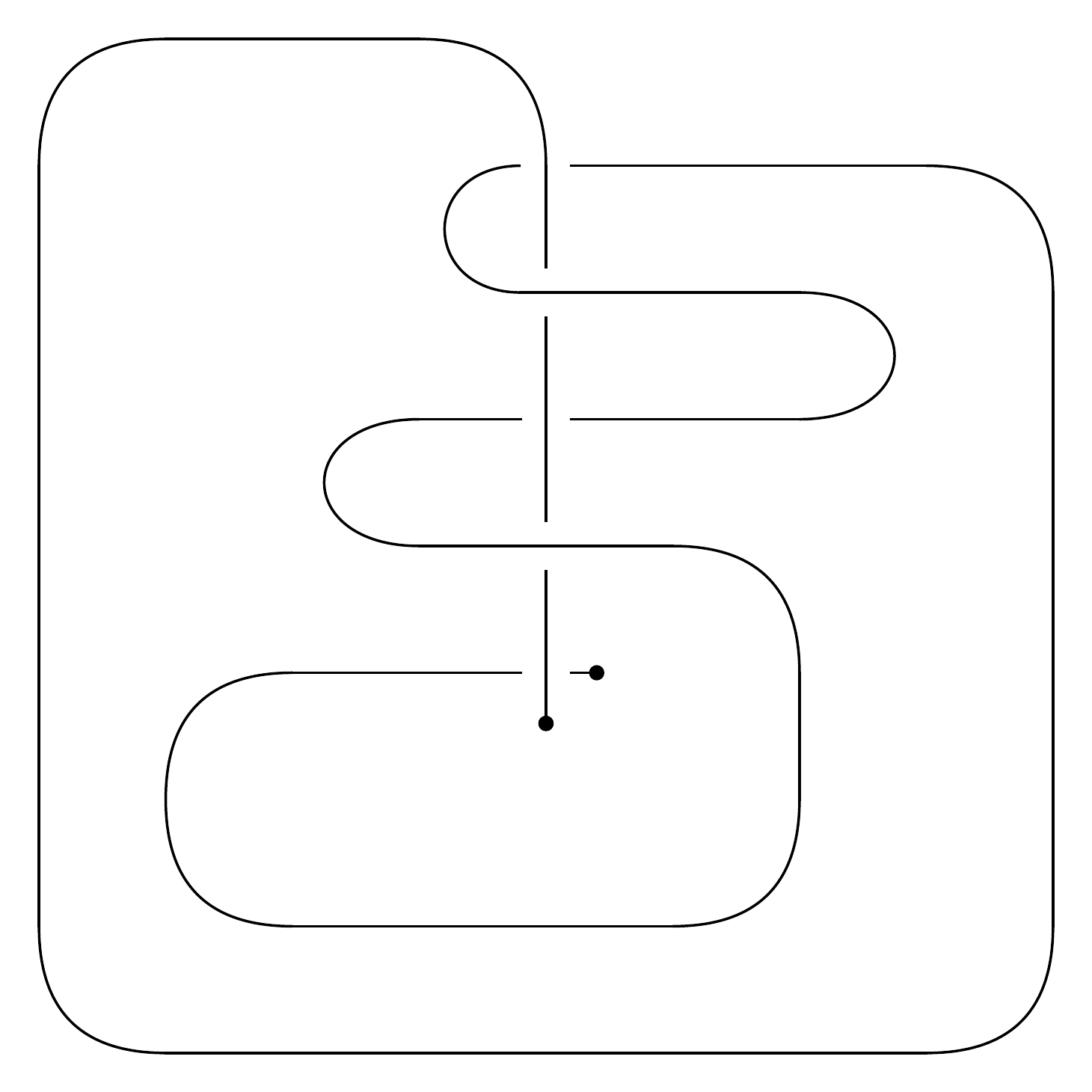}\\
\textcolor{black}{$5_{886}$}
\vspace{1cm}
\end{minipage}
\begin{minipage}[t]{.25\linewidth}
\centering
\includegraphics[width=0.9\textwidth,height=3.5cm,keepaspectratio]{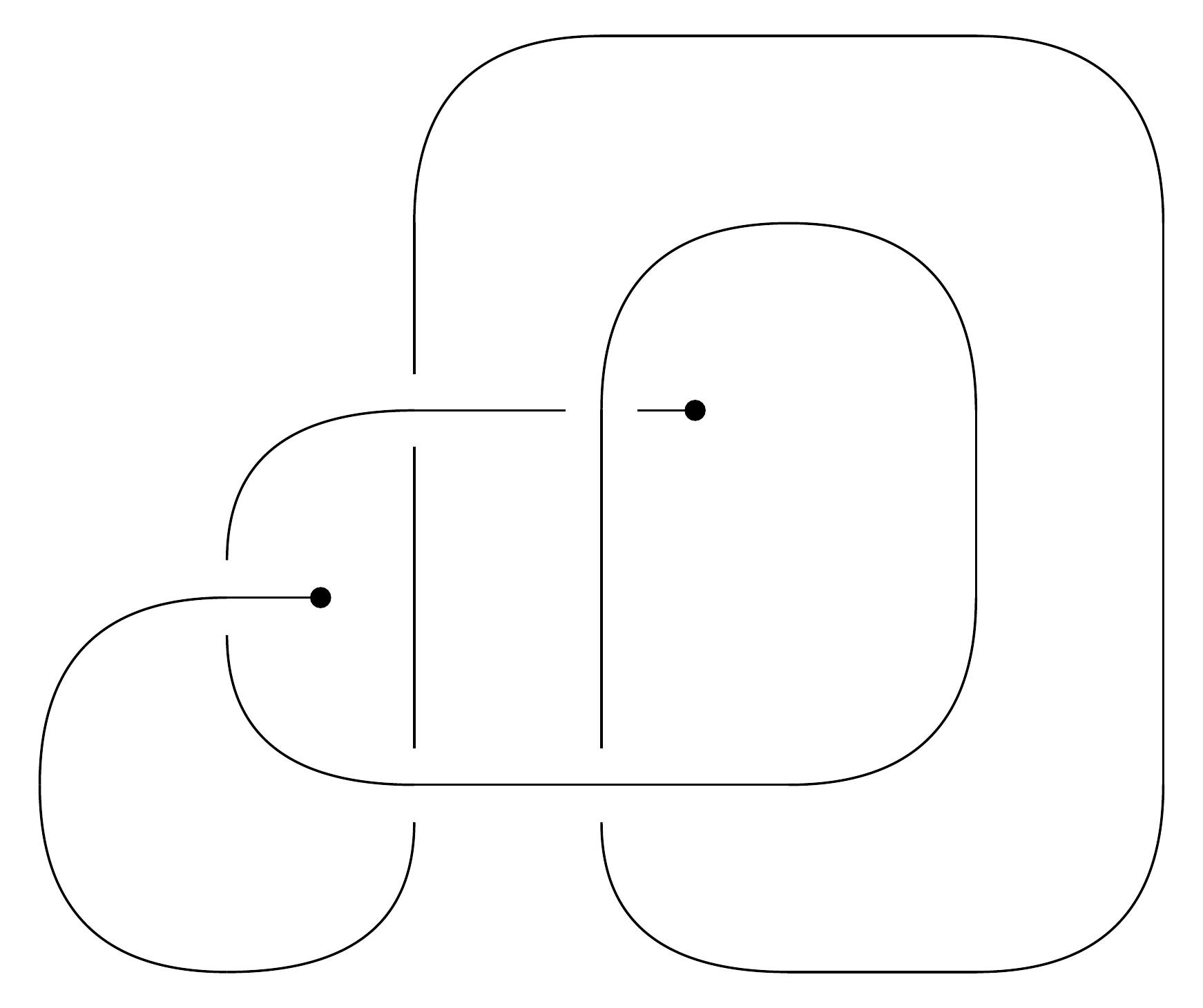}\\
\textcolor{black}{$5_{887}$}
\vspace{1cm}
\end{minipage}
\begin{minipage}[t]{.25\linewidth}
\centering
\includegraphics[width=0.9\textwidth,height=3.5cm,keepaspectratio]{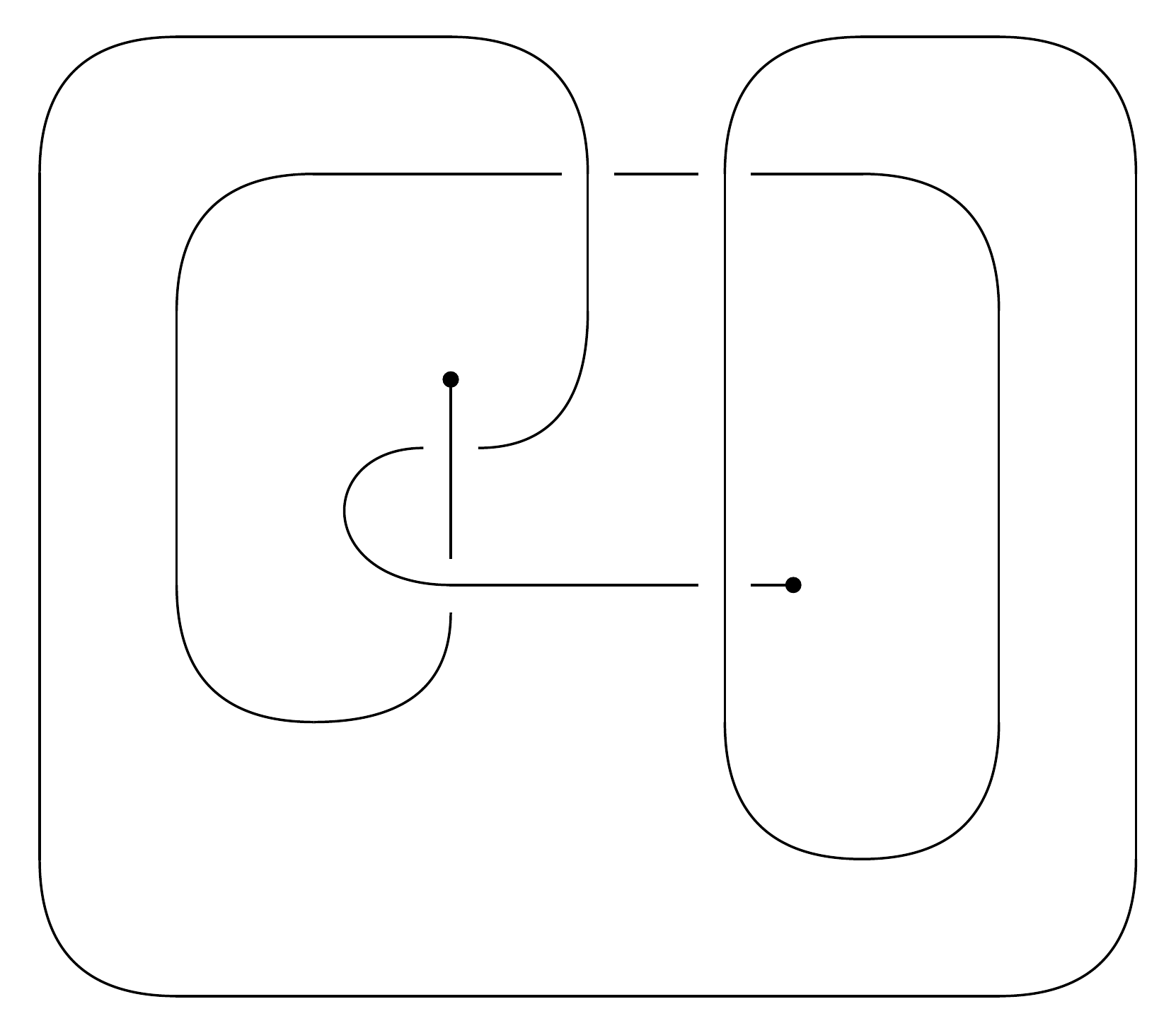}\\
\textcolor{black}{$5_{888}$}
\vspace{1cm}
\end{minipage}
\begin{minipage}[t]{.25\linewidth}
\centering
\includegraphics[width=0.9\textwidth,height=3.5cm,keepaspectratio]{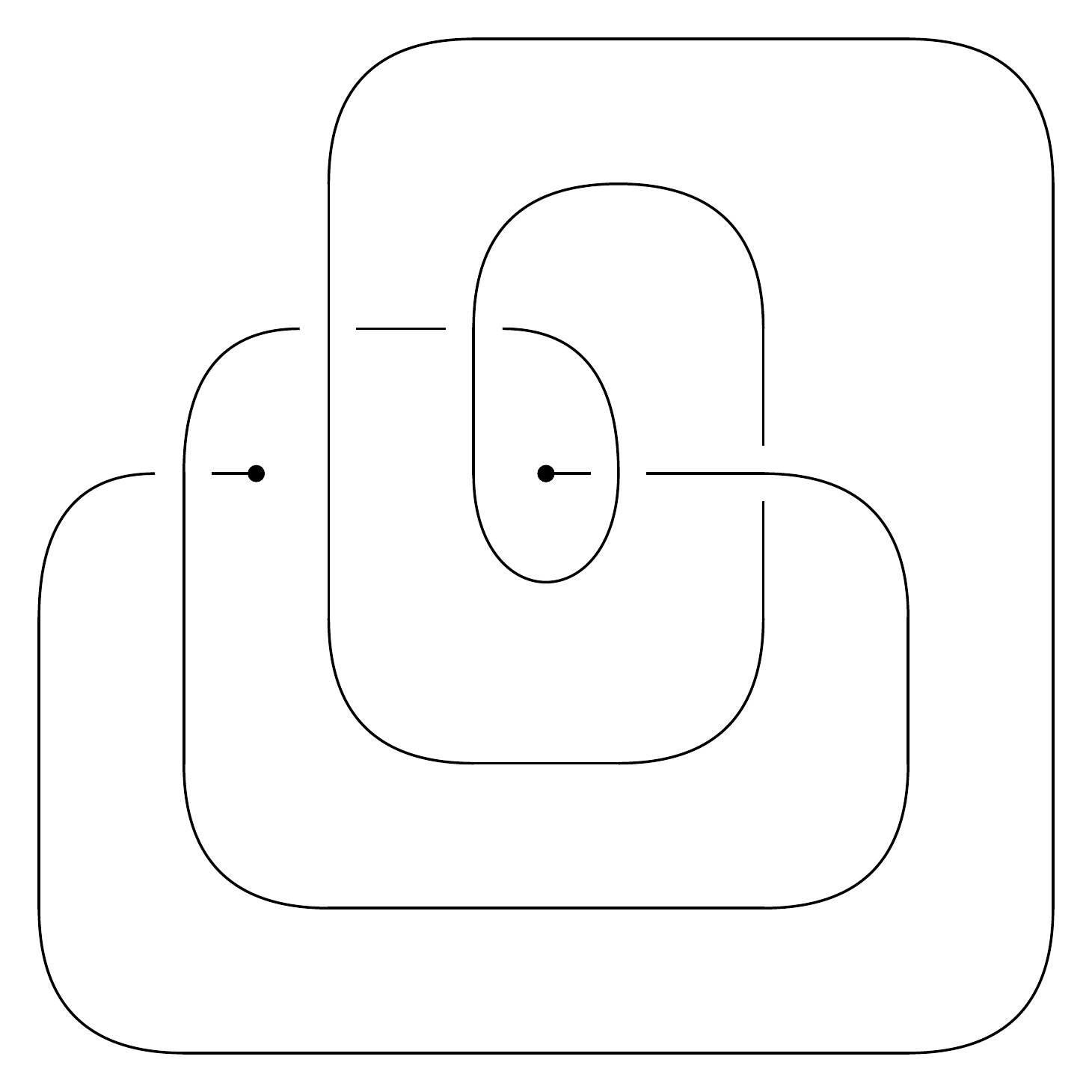}\\
\textcolor{black}{$5_{889}$}
\vspace{1cm}
\end{minipage}
\begin{minipage}[t]{.25\linewidth}
\centering
\includegraphics[width=0.9\textwidth,height=3.5cm,keepaspectratio]{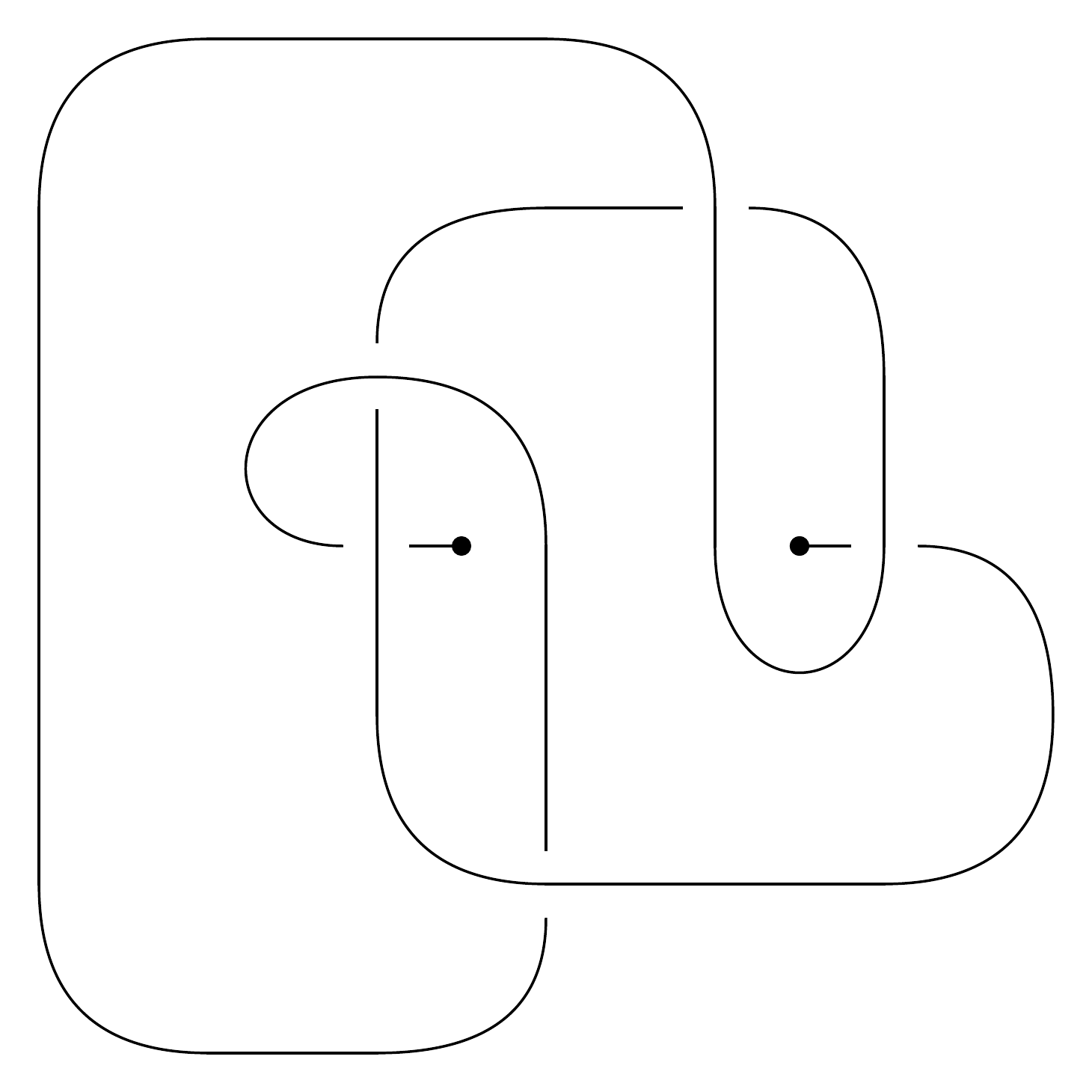}\\
\textcolor{black}{$5_{890}$}
\vspace{1cm}
\end{minipage}
\begin{minipage}[t]{.25\linewidth}
\centering
\includegraphics[width=0.9\textwidth,height=3.5cm,keepaspectratio]{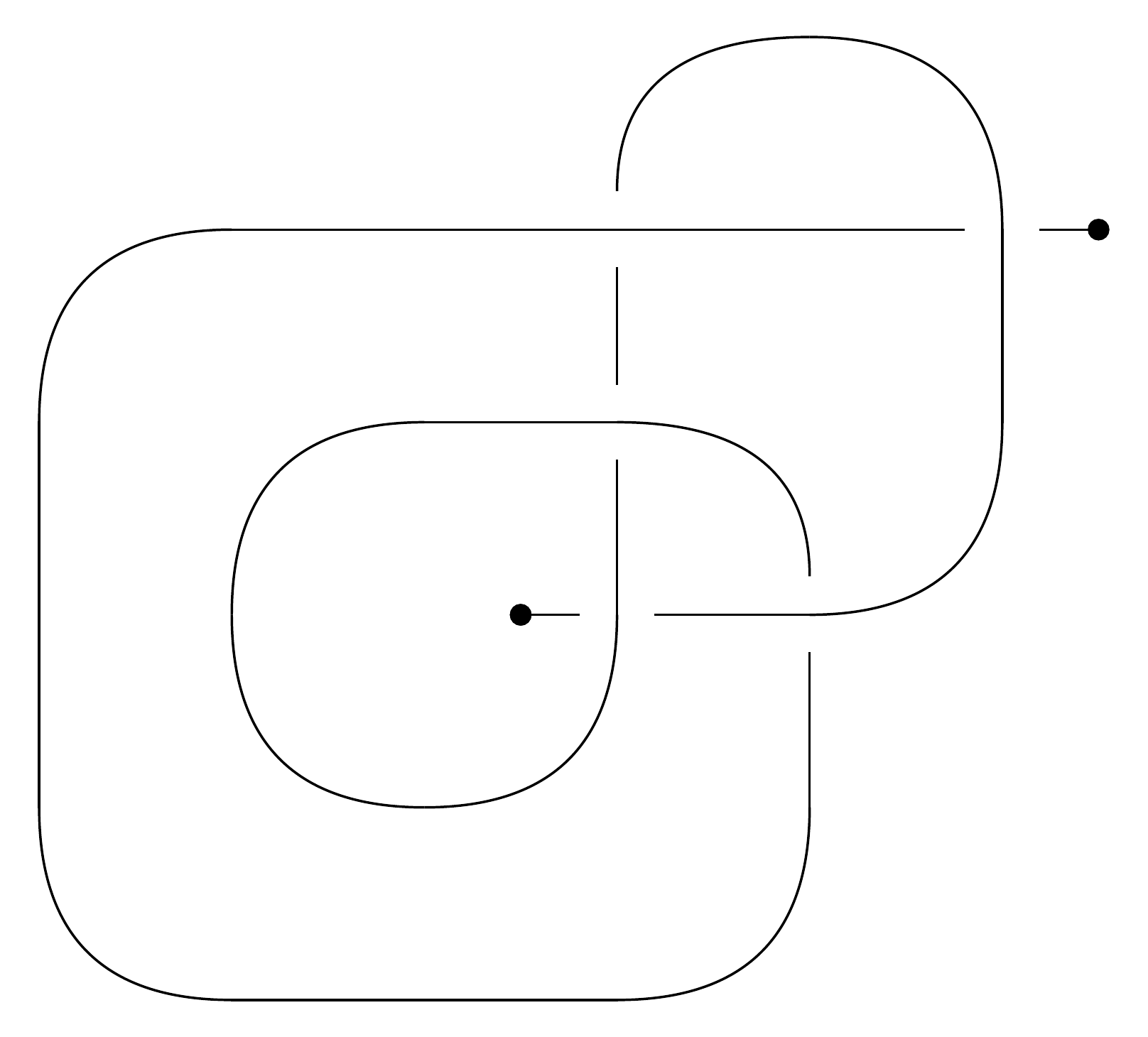}\\
\textcolor{black}{$5_{891}$}
\vspace{1cm}
\end{minipage}
\begin{minipage}[t]{.25\linewidth}
\centering
\includegraphics[width=0.9\textwidth,height=3.5cm,keepaspectratio]{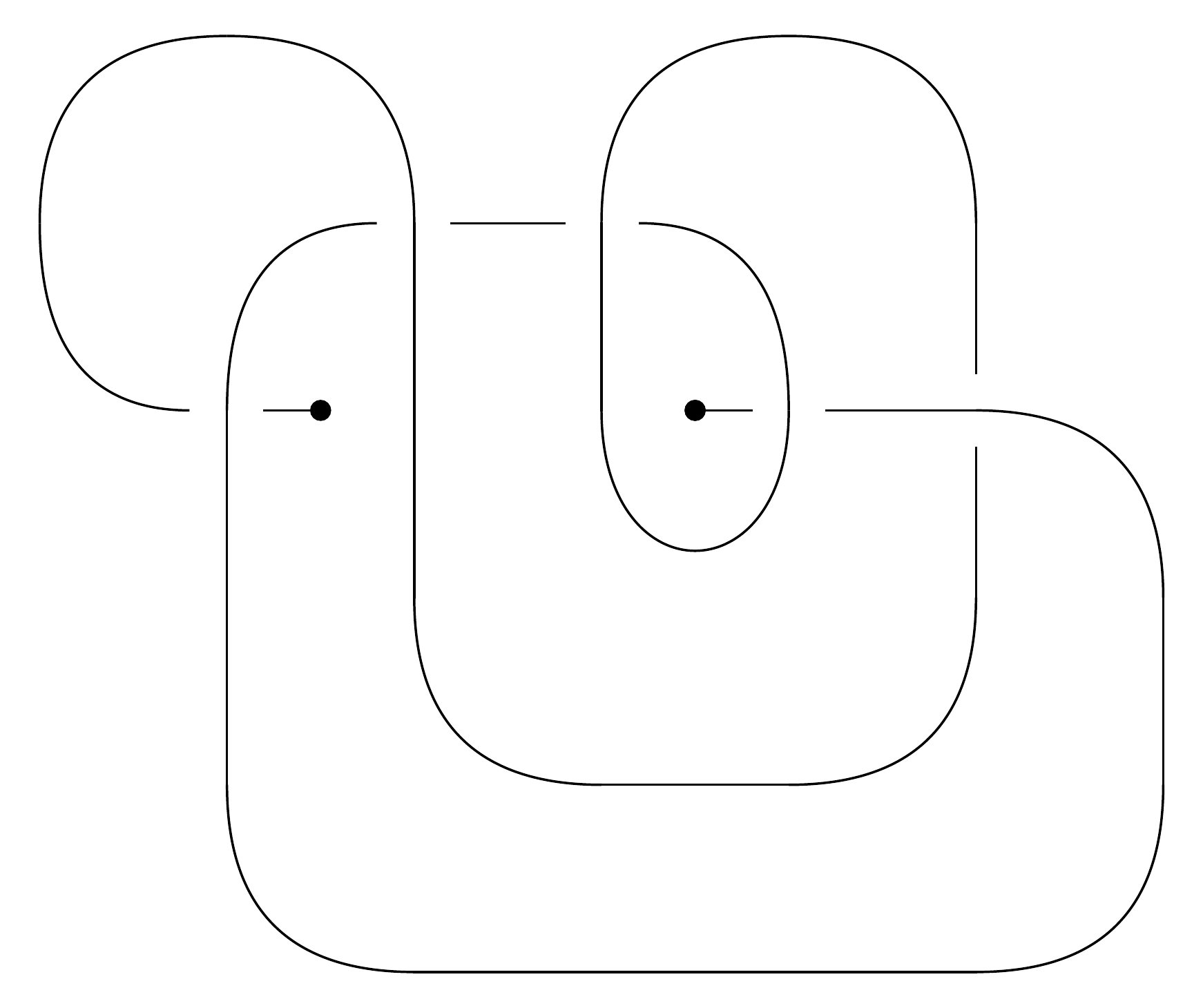}\\
\textcolor{black}{$5_{892}$}
\vspace{1cm}
\end{minipage}
\begin{minipage}[t]{.25\linewidth}
\centering
\includegraphics[width=0.9\textwidth,height=3.5cm,keepaspectratio]{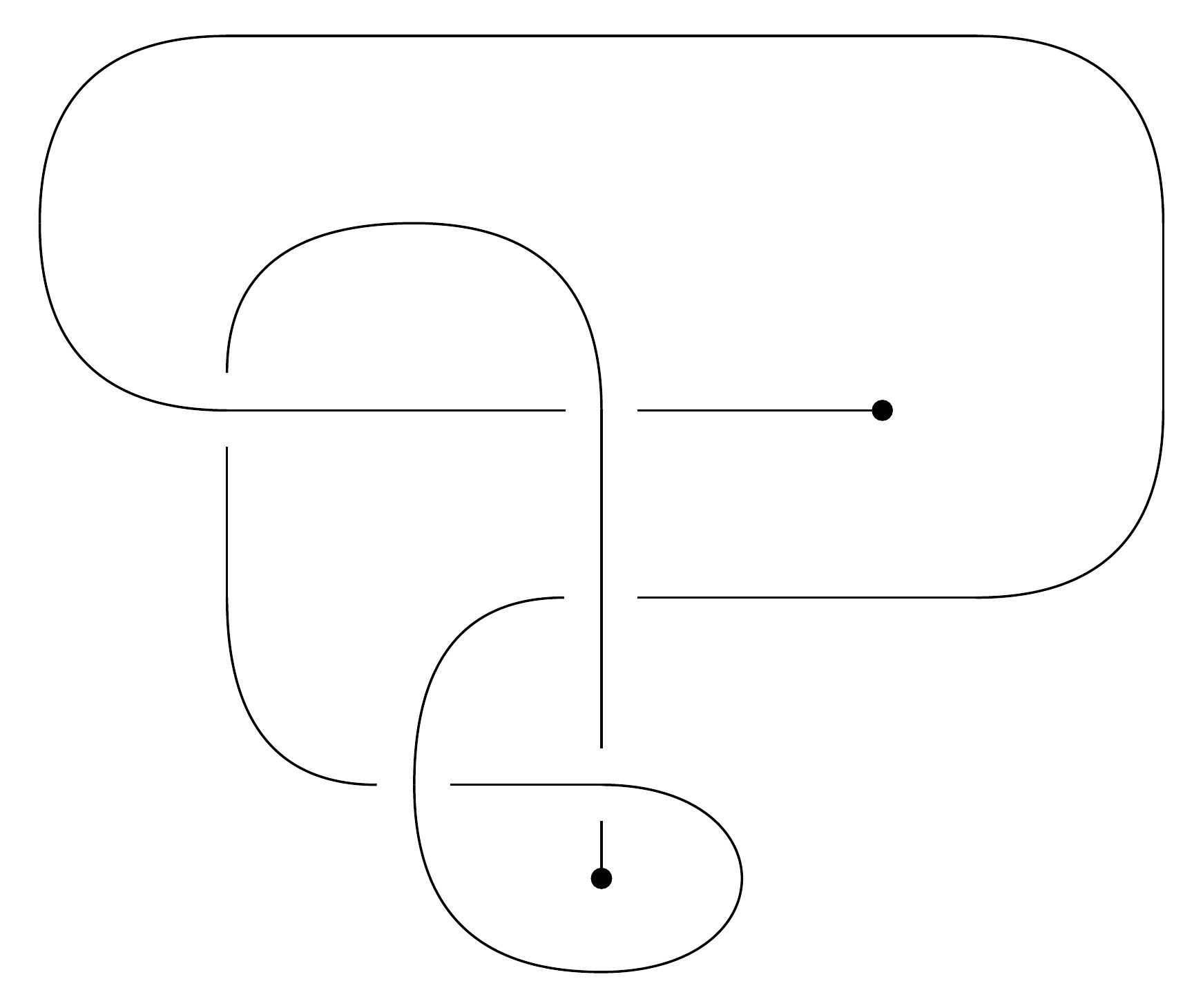}\\
\textcolor{black}{$5_{893}$}
\vspace{1cm}
\end{minipage}
\begin{minipage}[t]{.25\linewidth}
\centering
\includegraphics[width=0.9\textwidth,height=3.5cm,keepaspectratio]{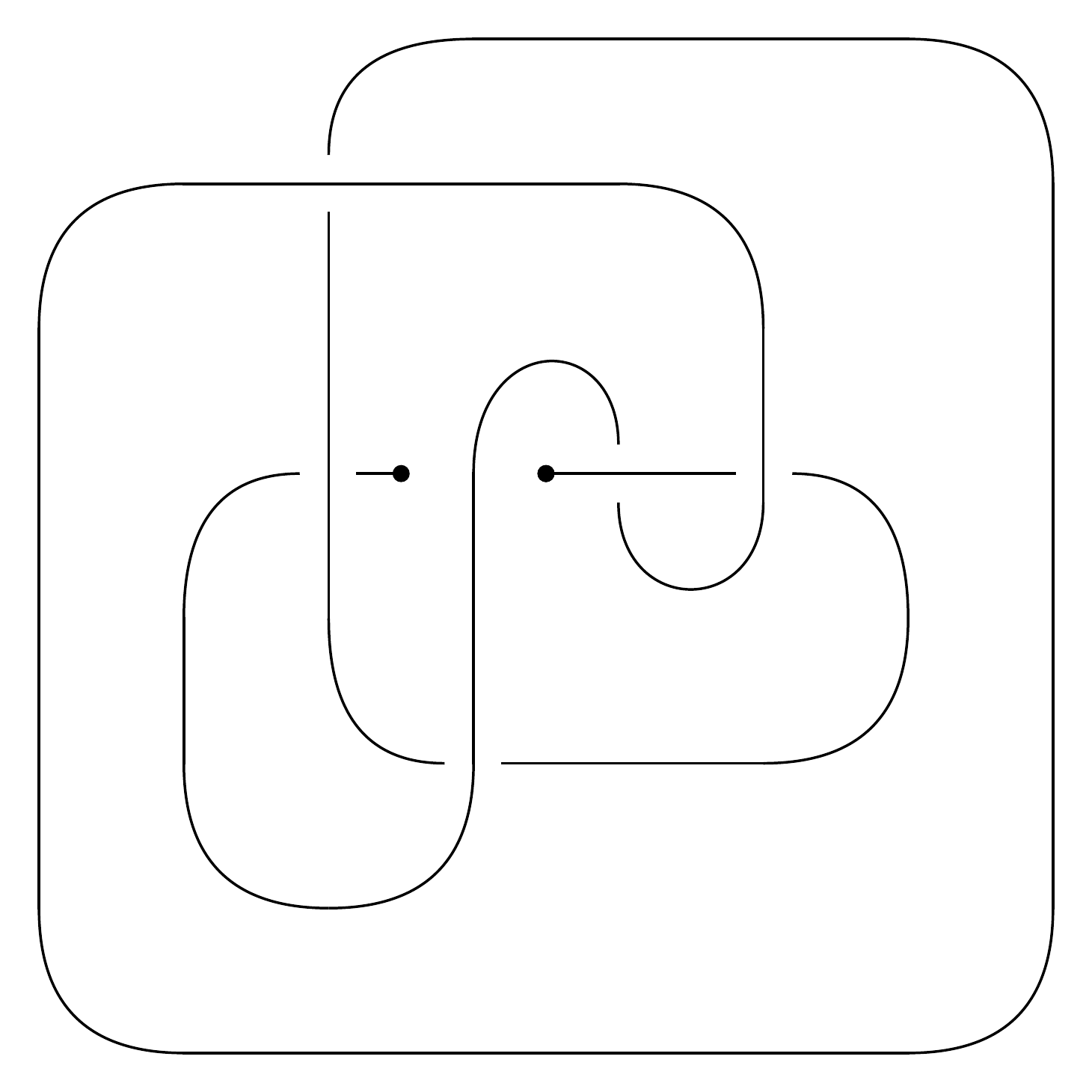}\\
\textcolor{black}{$5_{894}$}
\vspace{1cm}
\end{minipage}
\begin{minipage}[t]{.25\linewidth}
\centering
\includegraphics[width=0.9\textwidth,height=3.5cm,keepaspectratio]{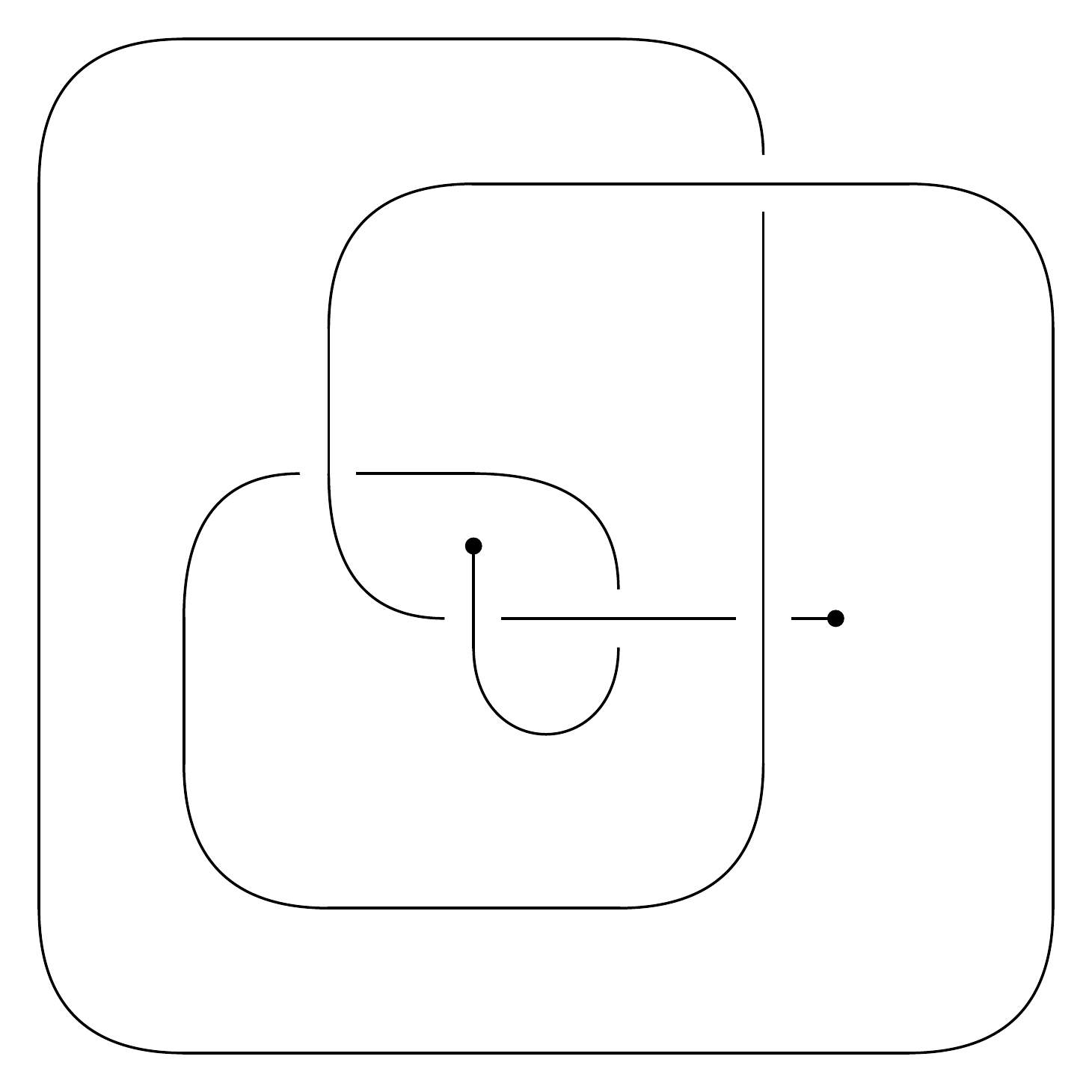}\\
\textcolor{black}{$5_{895}$}
\vspace{1cm}
\end{minipage}
\begin{minipage}[t]{.25\linewidth}
\centering
\includegraphics[width=0.9\textwidth,height=3.5cm,keepaspectratio]{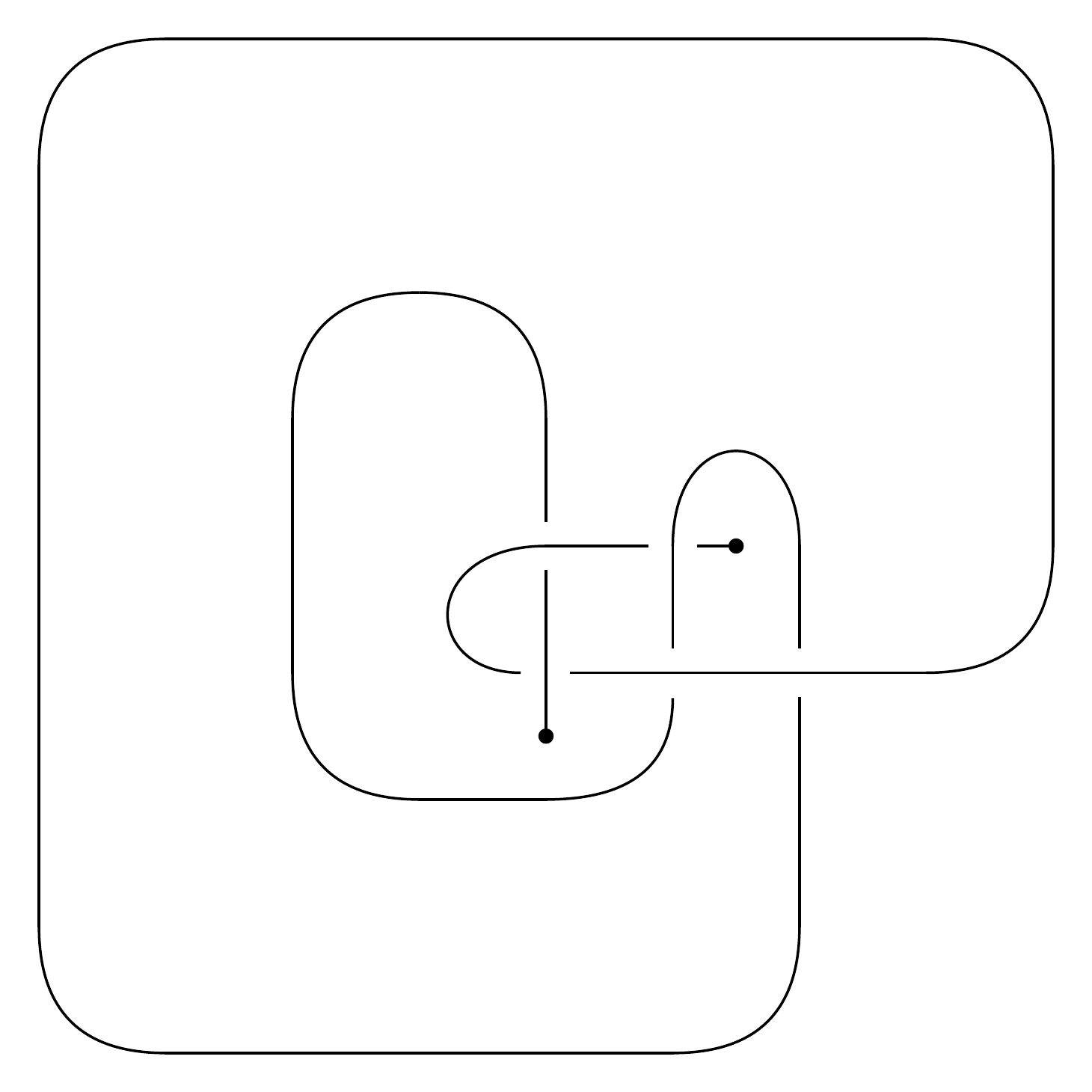}\\
\textcolor{black}{$5_{896}$}
\vspace{1cm}
\end{minipage}
\begin{minipage}[t]{.25\linewidth}
\centering
\includegraphics[width=0.9\textwidth,height=3.5cm,keepaspectratio]{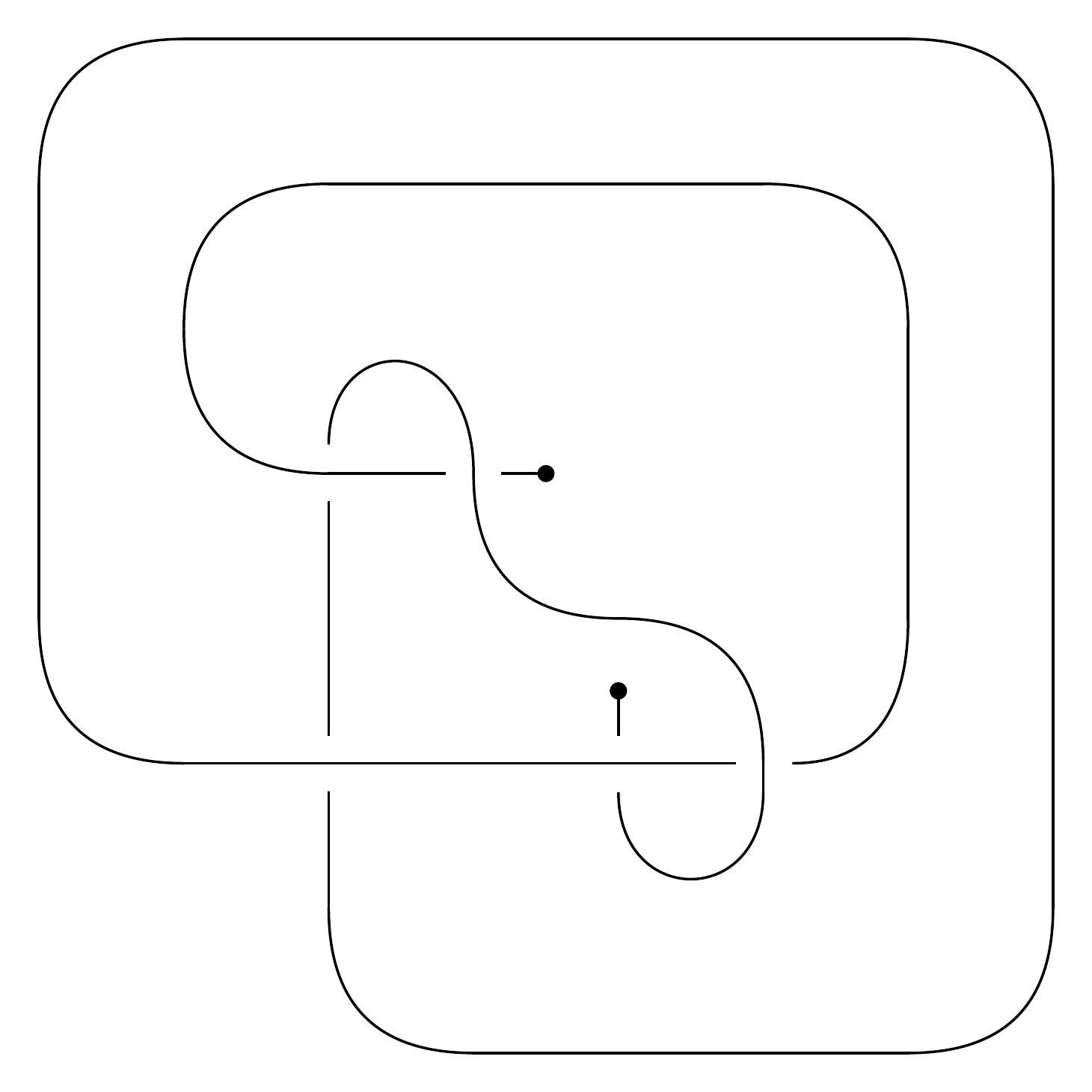}\\
\textcolor{black}{$5_{897}$}
\vspace{1cm}
\end{minipage}
\begin{minipage}[t]{.25\linewidth}
\centering
\includegraphics[width=0.9\textwidth,height=3.5cm,keepaspectratio]{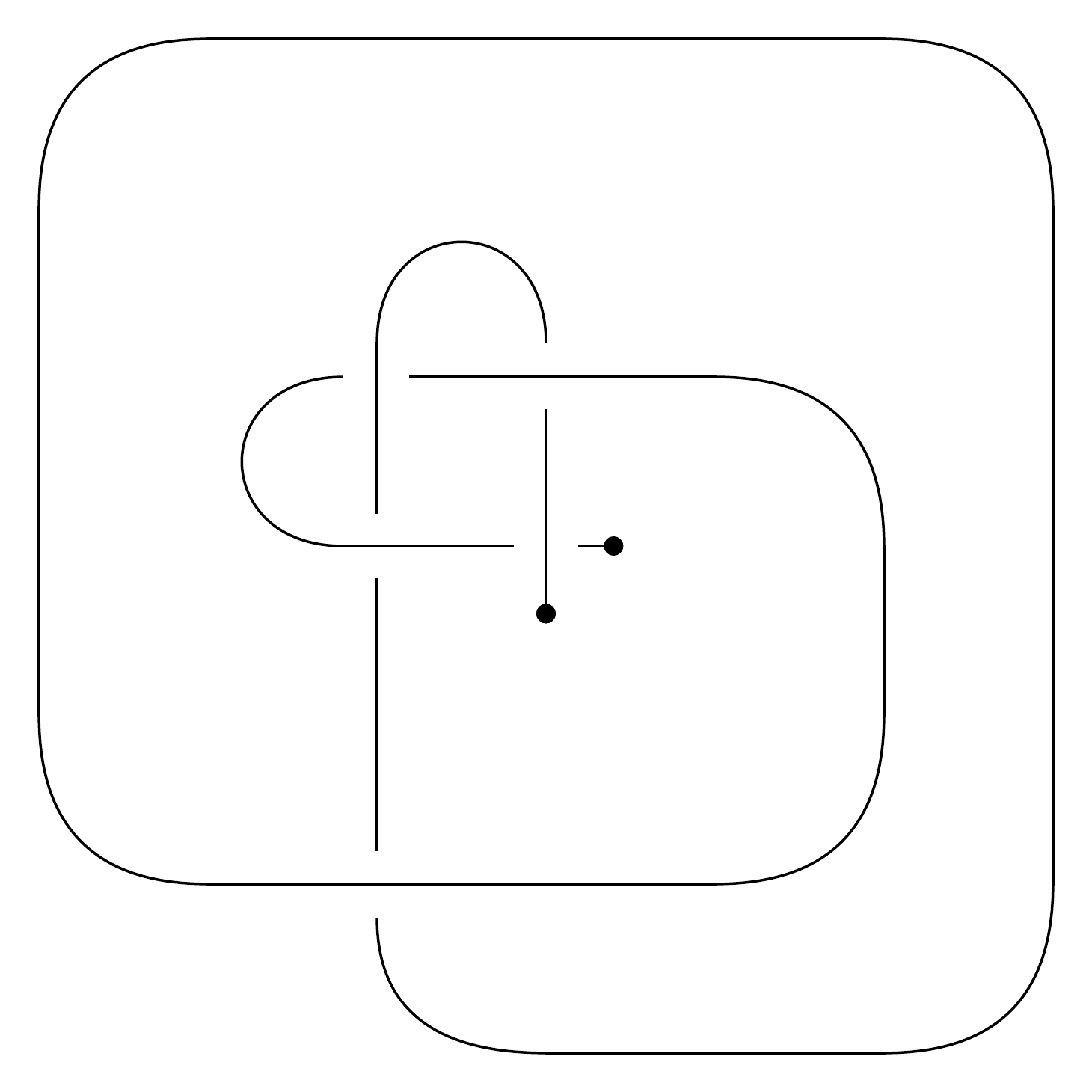}\\
\textcolor{black}{$5_{898}$}
\vspace{1cm}
\end{minipage}
\begin{minipage}[t]{.25\linewidth}
\centering
\includegraphics[width=0.9\textwidth,height=3.5cm,keepaspectratio]{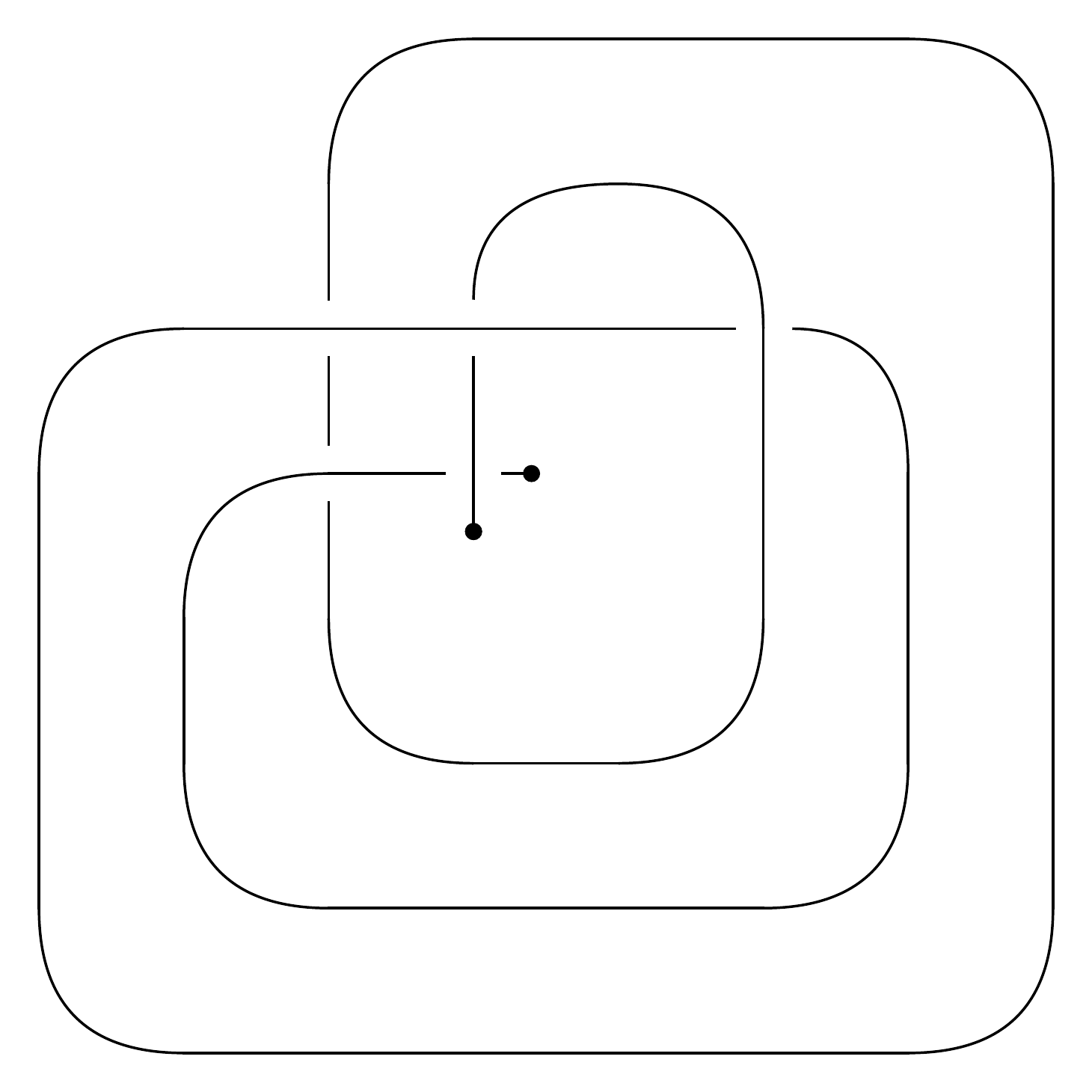}\\
\textcolor{black}{$5_{899}$}
\vspace{1cm}
\end{minipage}
\begin{minipage}[t]{.25\linewidth}
\centering
\includegraphics[width=0.9\textwidth,height=3.5cm,keepaspectratio]{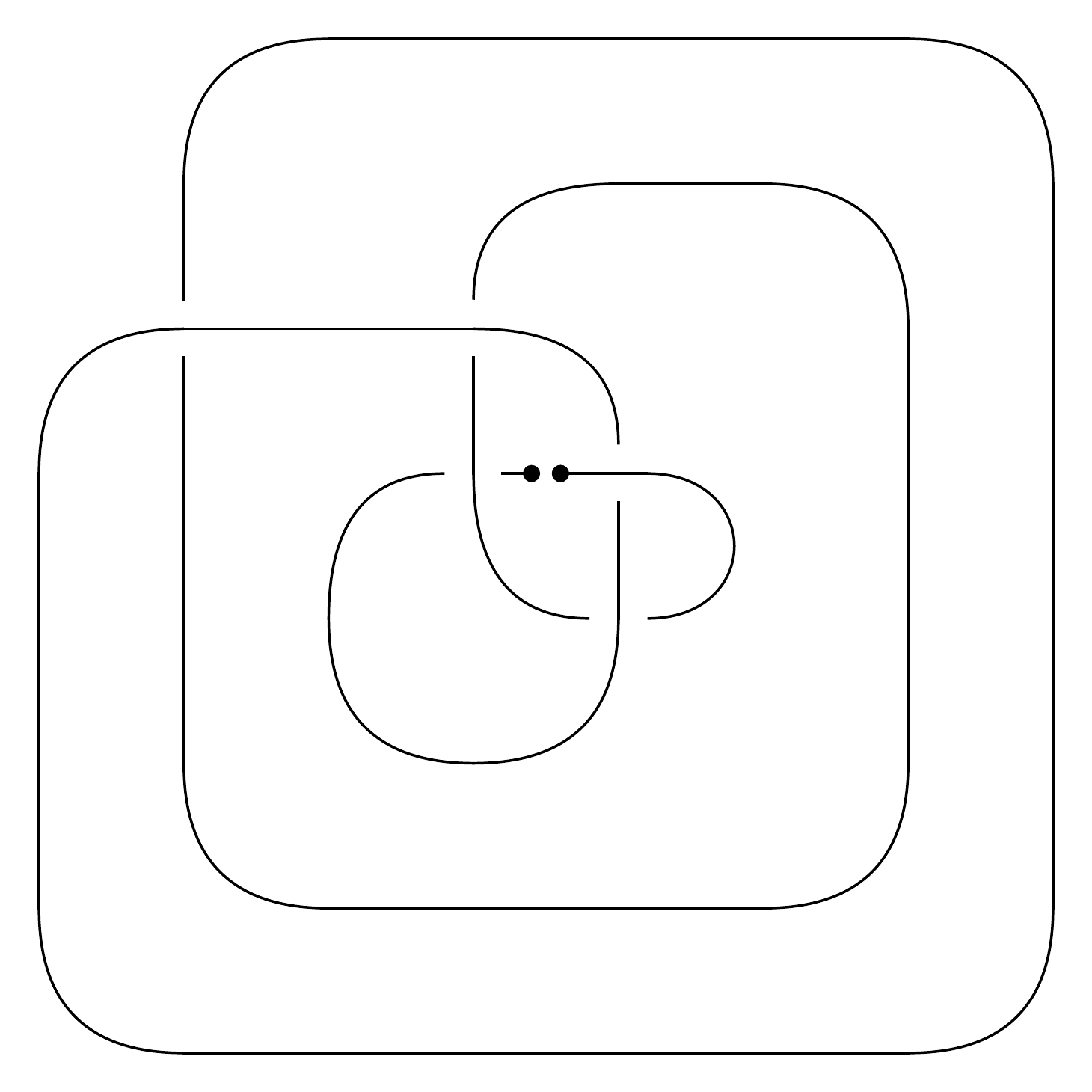}\\
\textcolor{black}{$5_{900}$}
\vspace{1cm}
\end{minipage}
\begin{minipage}[t]{.25\linewidth}
\centering
\includegraphics[width=0.9\textwidth,height=3.5cm,keepaspectratio]{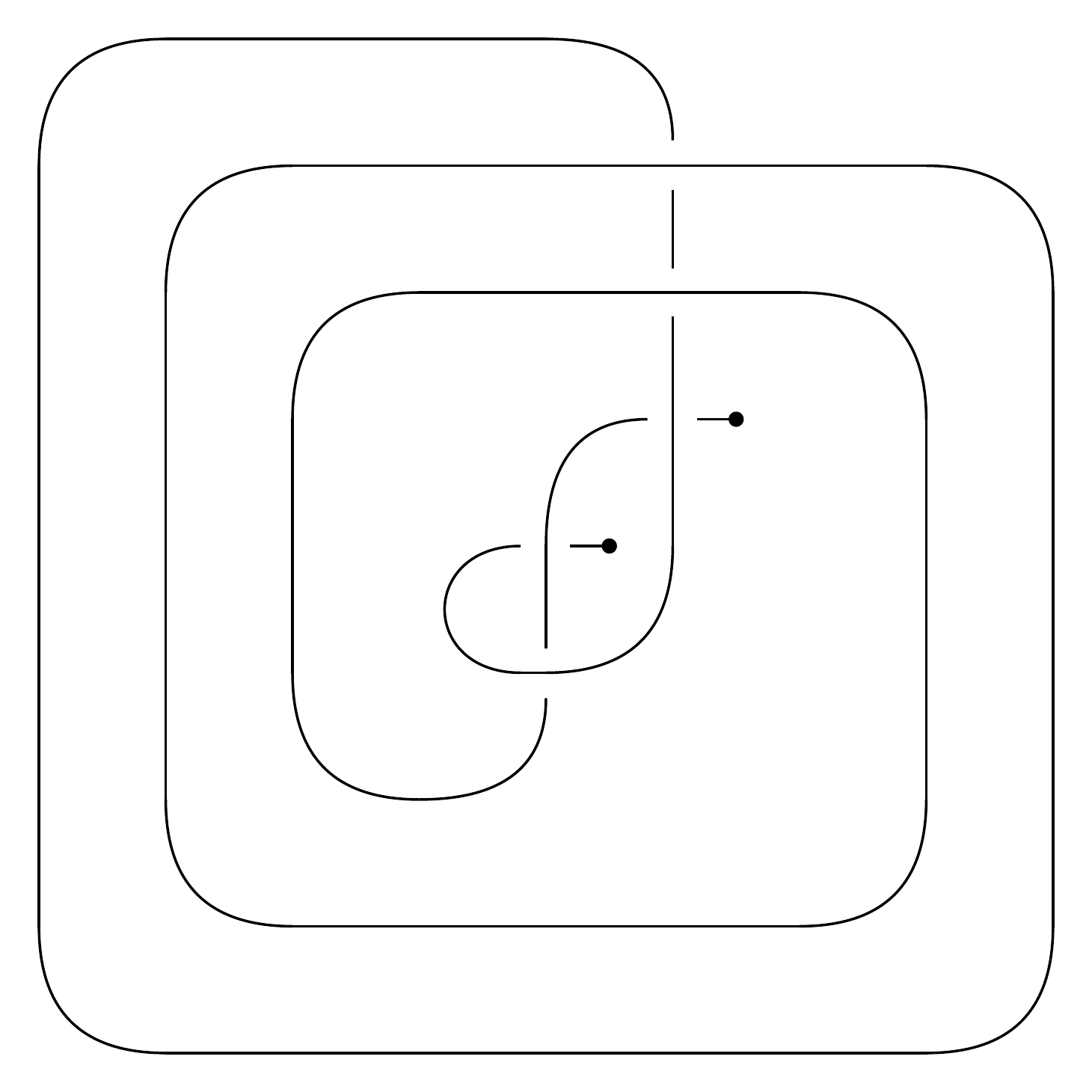}\\
\textcolor{black}{$5_{901}$}
\vspace{1cm}
\end{minipage}
\begin{minipage}[t]{.25\linewidth}
\centering
\includegraphics[width=0.9\textwidth,height=3.5cm,keepaspectratio]{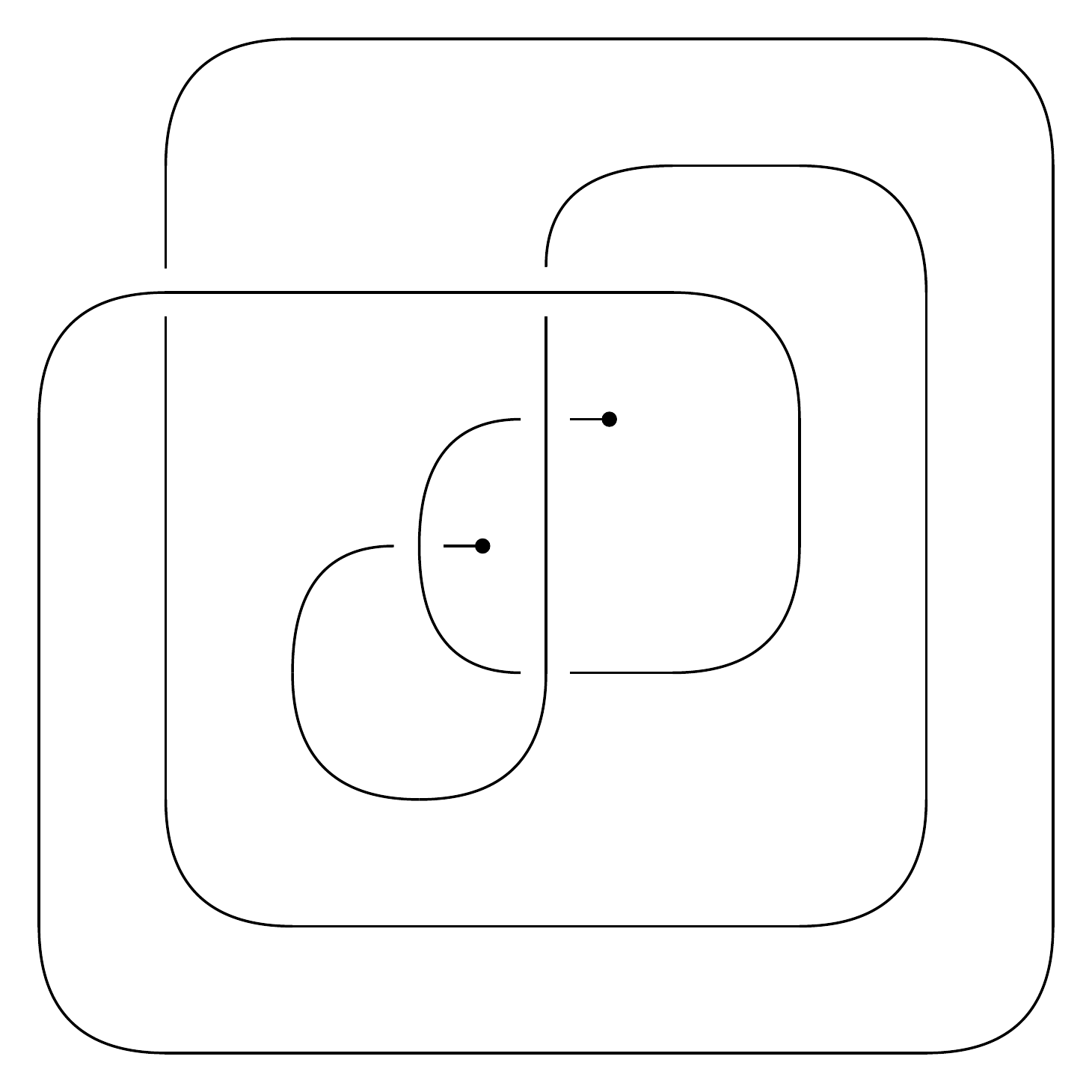}\\
\textcolor{black}{$5_{902}$}
\vspace{1cm}
\end{minipage}
\begin{minipage}[t]{.25\linewidth}
\centering
\includegraphics[width=0.9\textwidth,height=3.5cm,keepaspectratio]{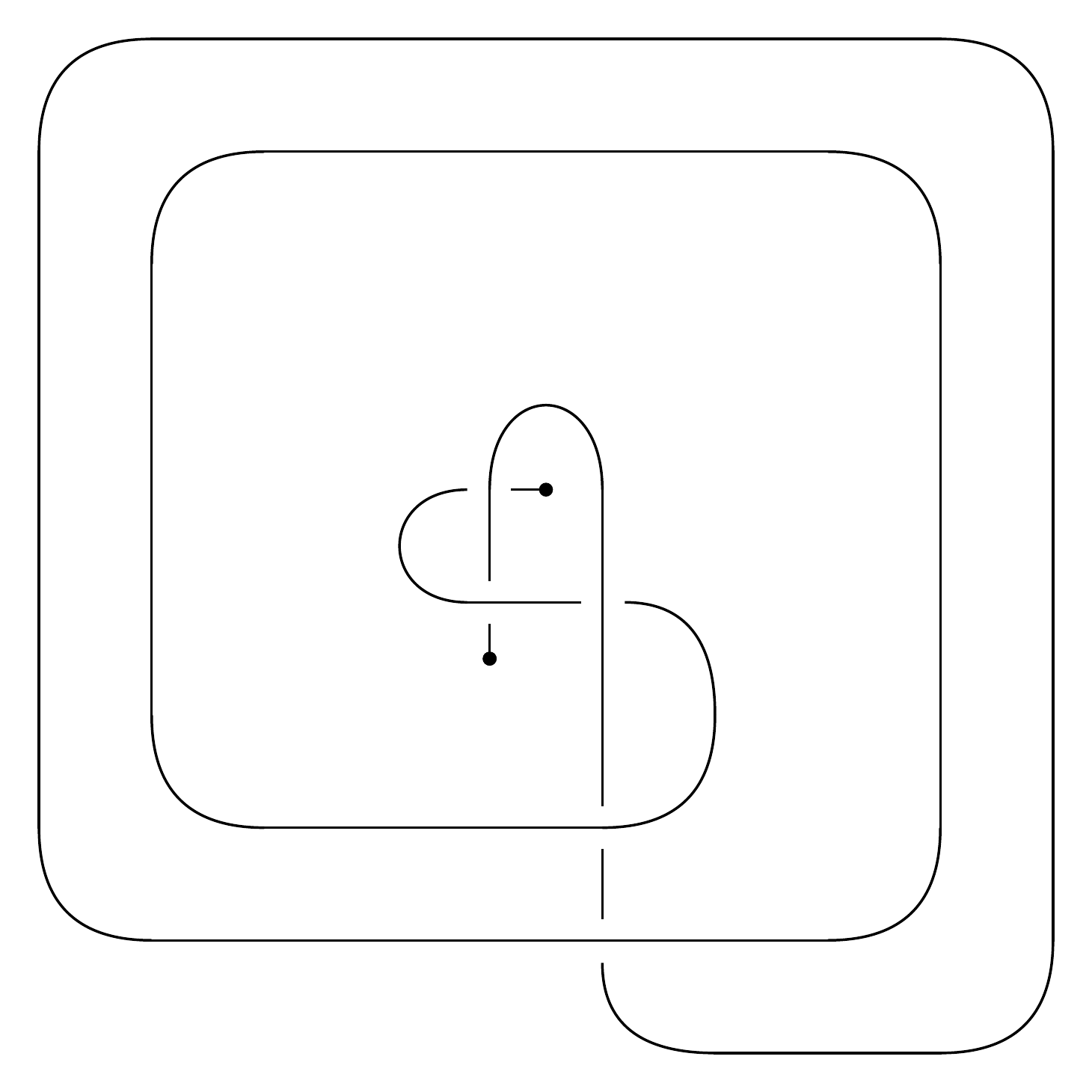}\\
\textcolor{black}{$5_{903}$}
\vspace{1cm}
\end{minipage}
\begin{minipage}[t]{.25\linewidth}
\centering
\includegraphics[width=0.9\textwidth,height=3.5cm,keepaspectratio]{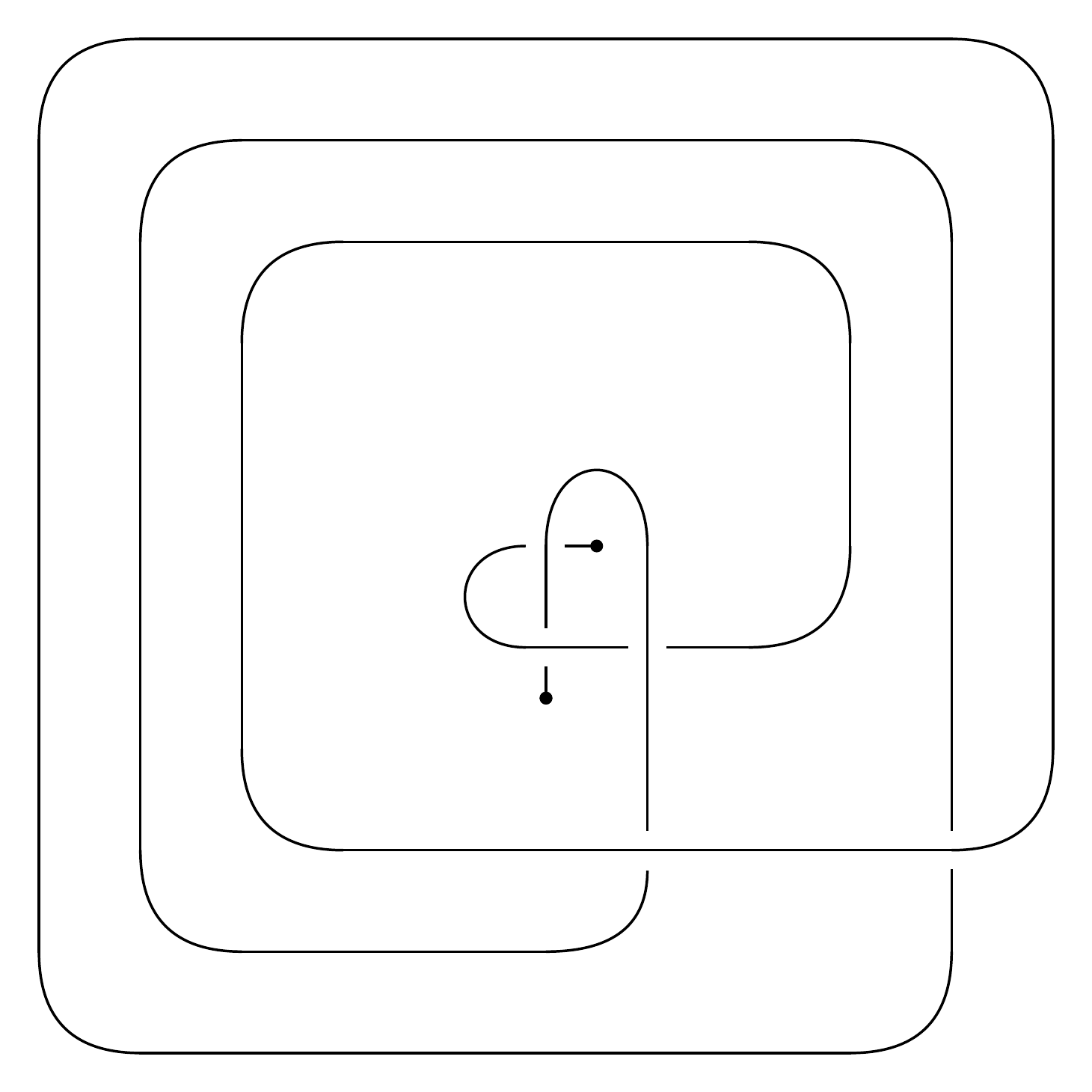}\\
\textcolor{black}{$5_{904}$}
\vspace{1cm}
\end{minipage}
\begin{minipage}[t]{.25\linewidth}
\centering
\includegraphics[width=0.9\textwidth,height=3.5cm,keepaspectratio]{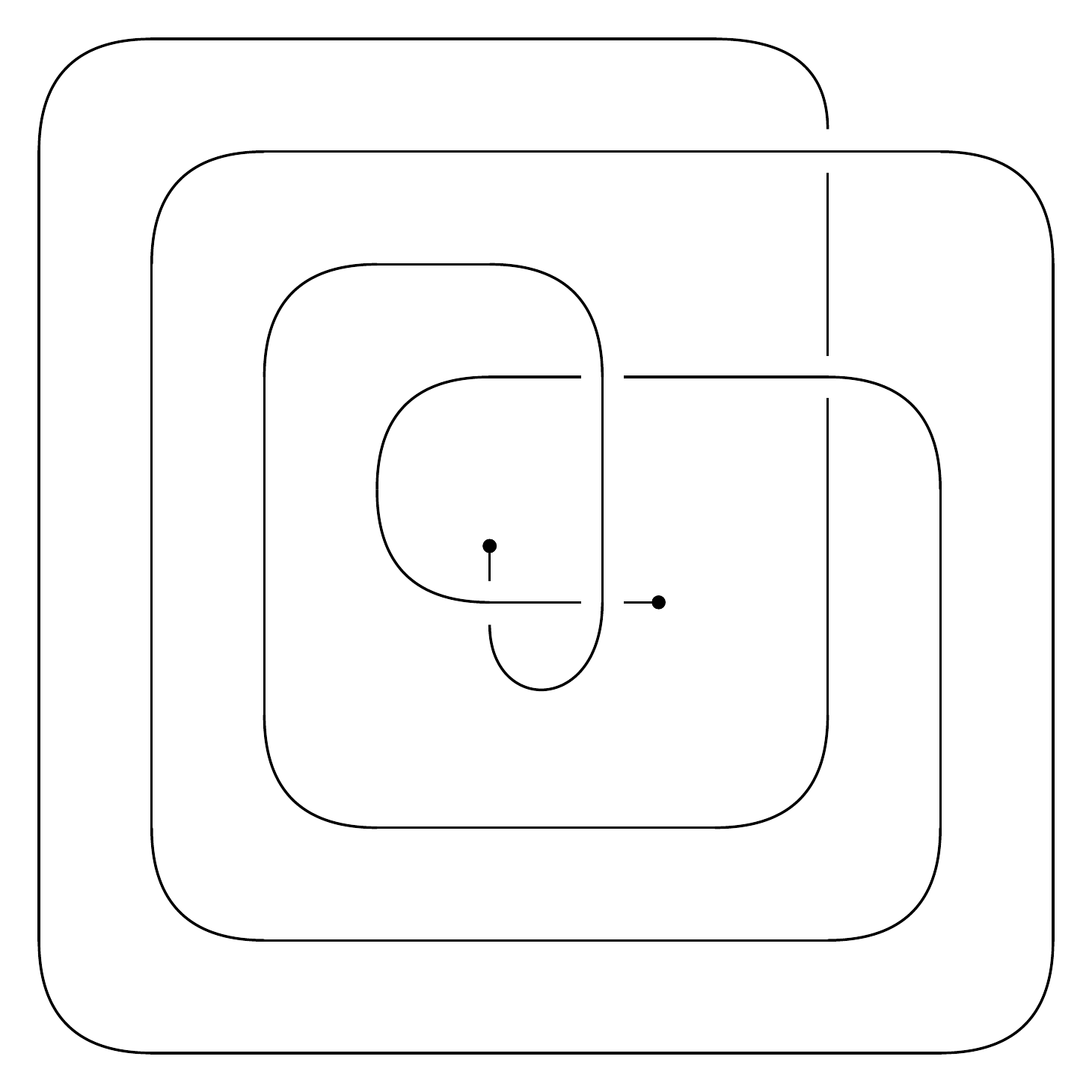}\\
\textcolor{black}{$5_{905}$}
\vspace{1cm}
\end{minipage}
\begin{minipage}[t]{.25\linewidth}
\centering
\includegraphics[width=0.9\textwidth,height=3.5cm,keepaspectratio]{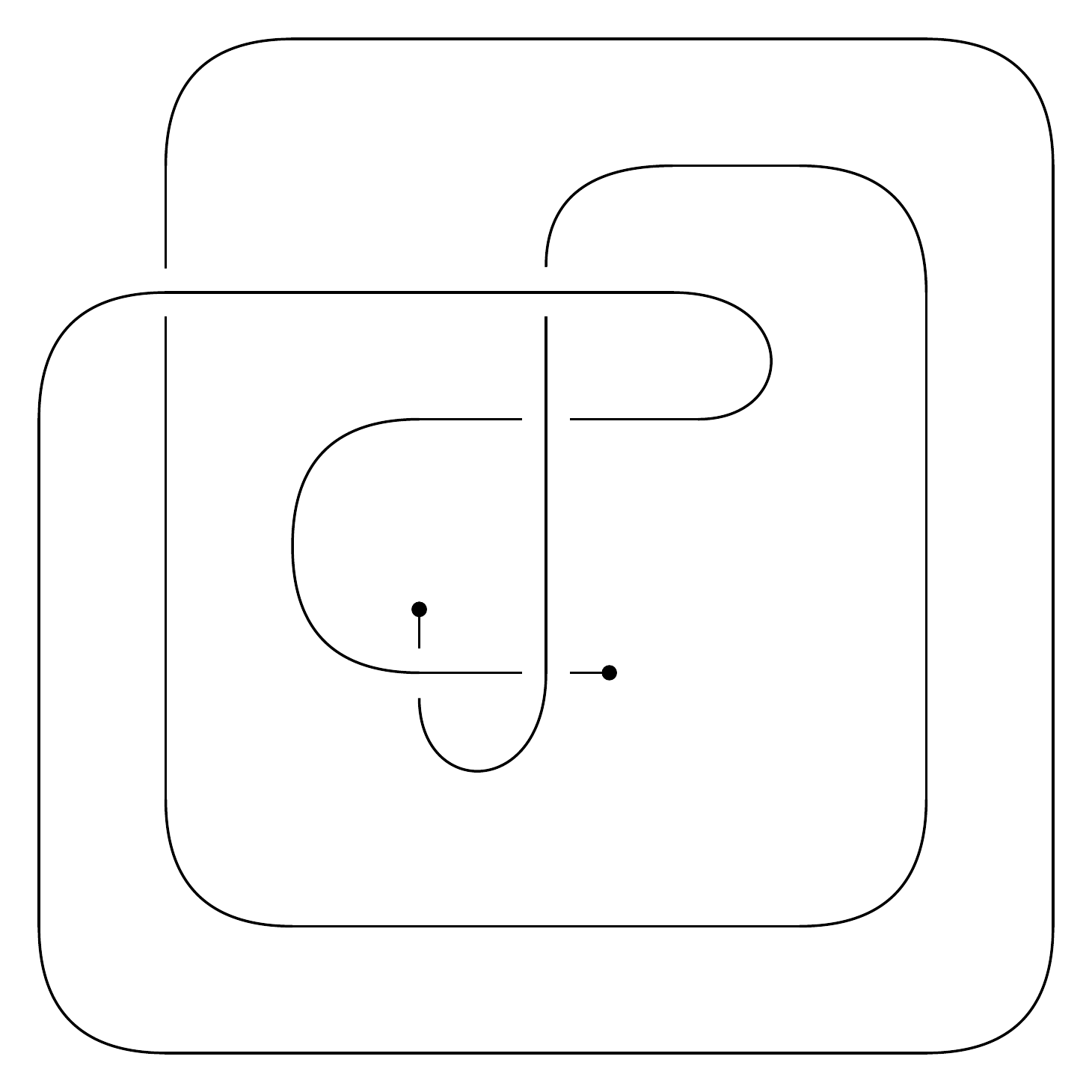}\\
\textcolor{black}{$5_{906}$}
\vspace{1cm}
\end{minipage}
\begin{minipage}[t]{.25\linewidth}
\centering
\includegraphics[width=0.9\textwidth,height=3.5cm,keepaspectratio]{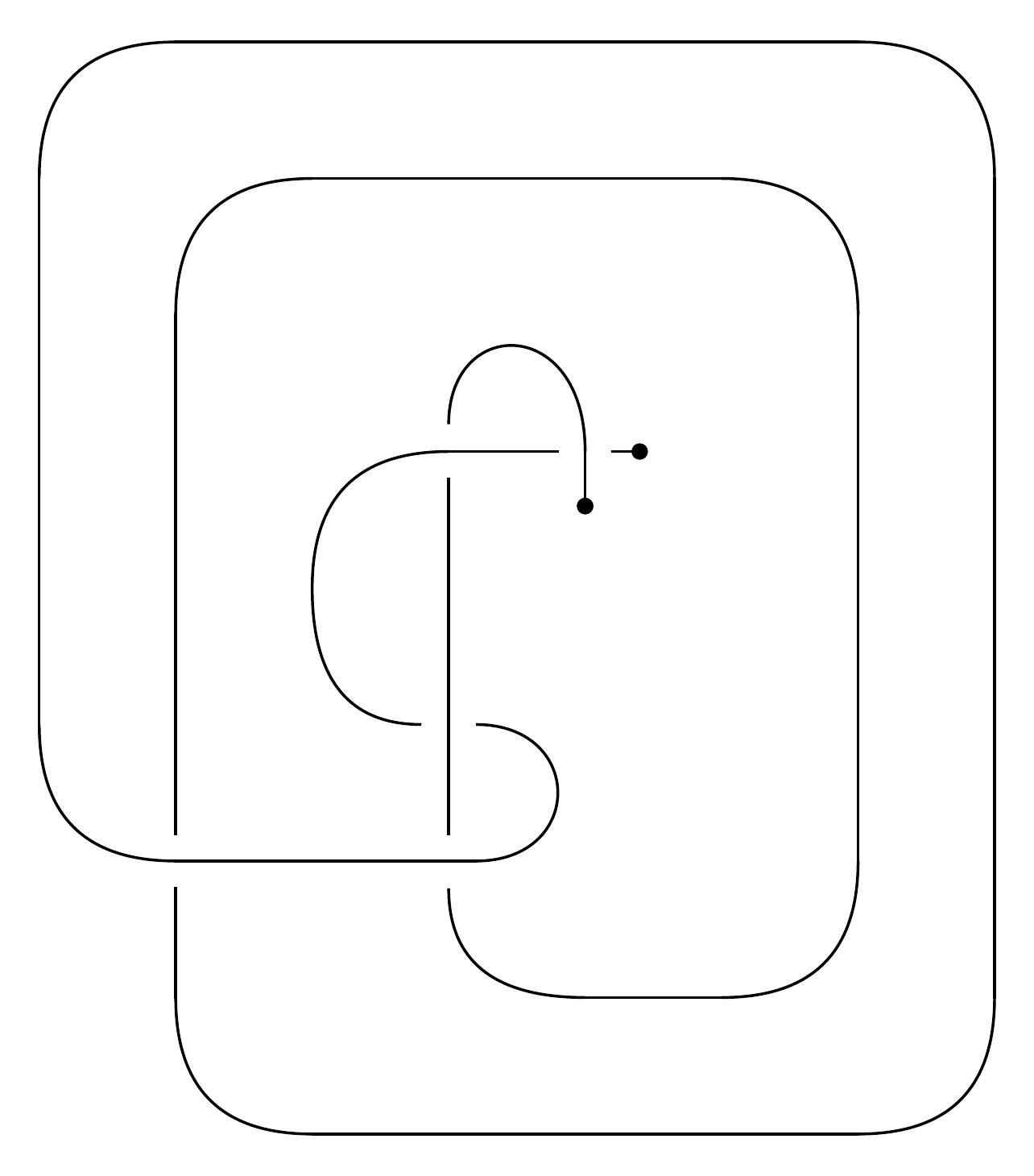}\\
\textcolor{black}{$5_{907}$}
\vspace{1cm}
\end{minipage}
\begin{minipage}[t]{.25\linewidth}
\centering
\includegraphics[width=0.9\textwidth,height=3.5cm,keepaspectratio]{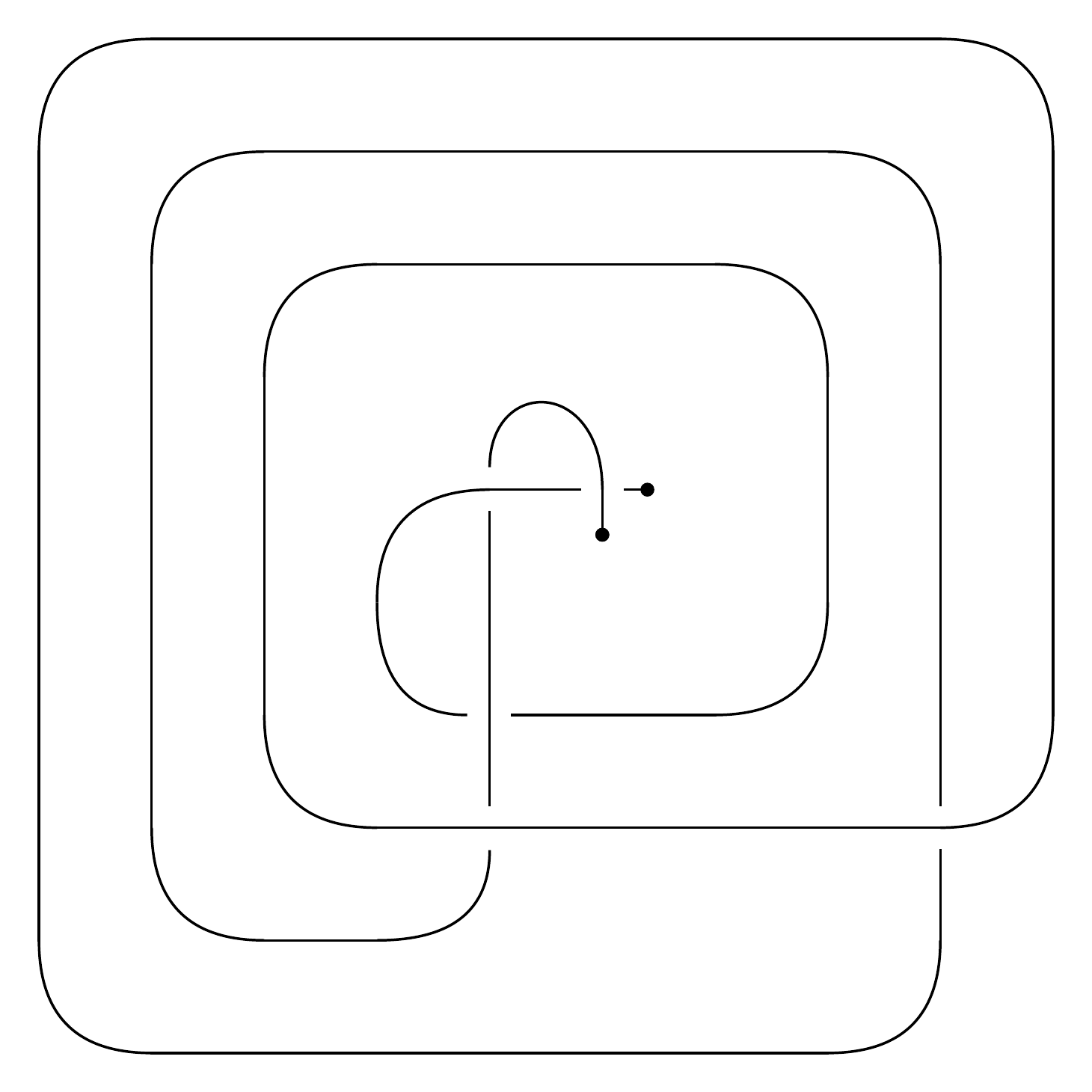}\\
\textcolor{black}{$5_{908}$}
\vspace{1cm}
\end{minipage}
\begin{minipage}[t]{.25\linewidth}
\centering
\includegraphics[width=0.9\textwidth,height=3.5cm,keepaspectratio]{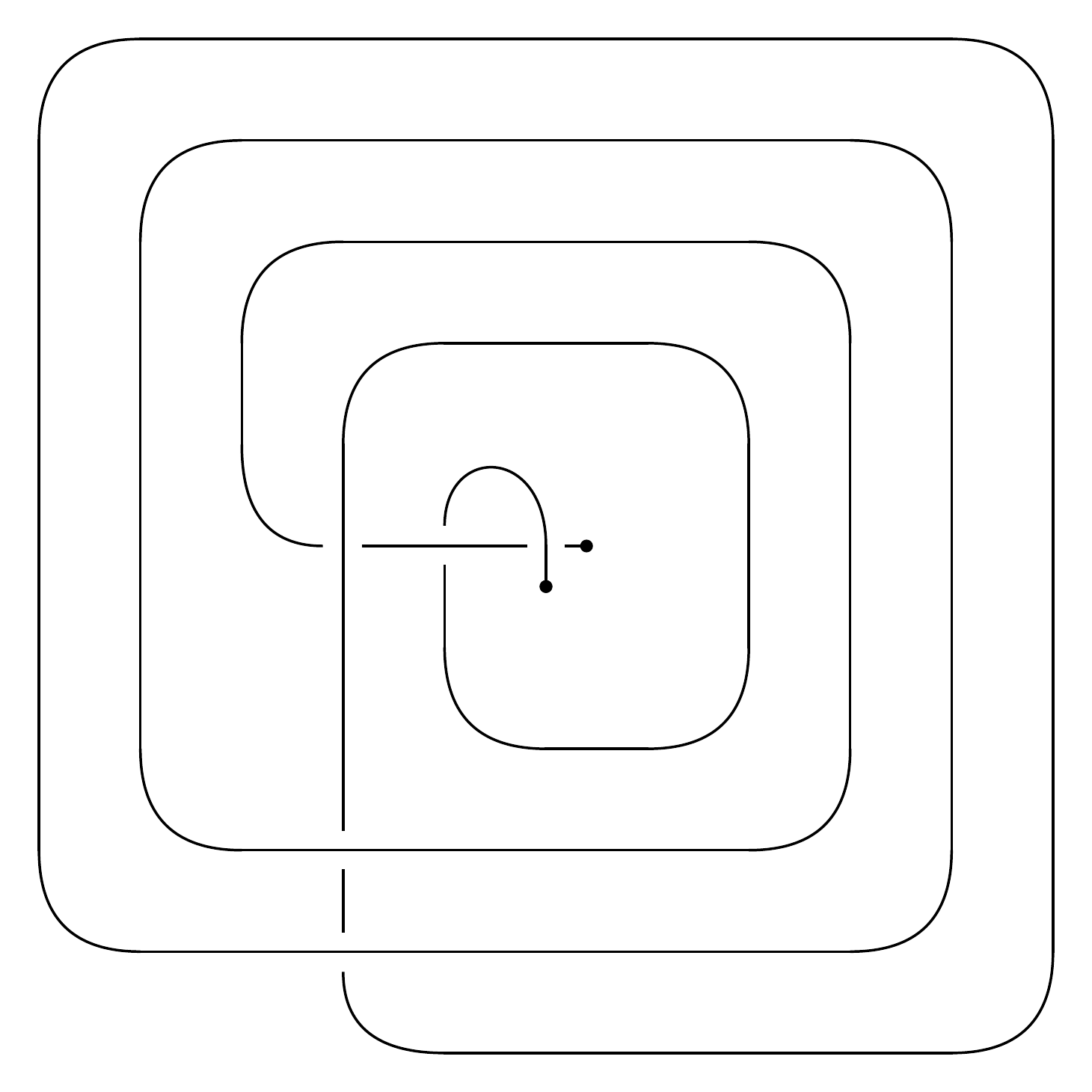}\\
\textcolor{black}{$5_{909}$}
\vspace{1cm}
\end{minipage}
\begin{minipage}[t]{.25\linewidth}
\centering
\includegraphics[width=0.9\textwidth,height=3.5cm,keepaspectratio]{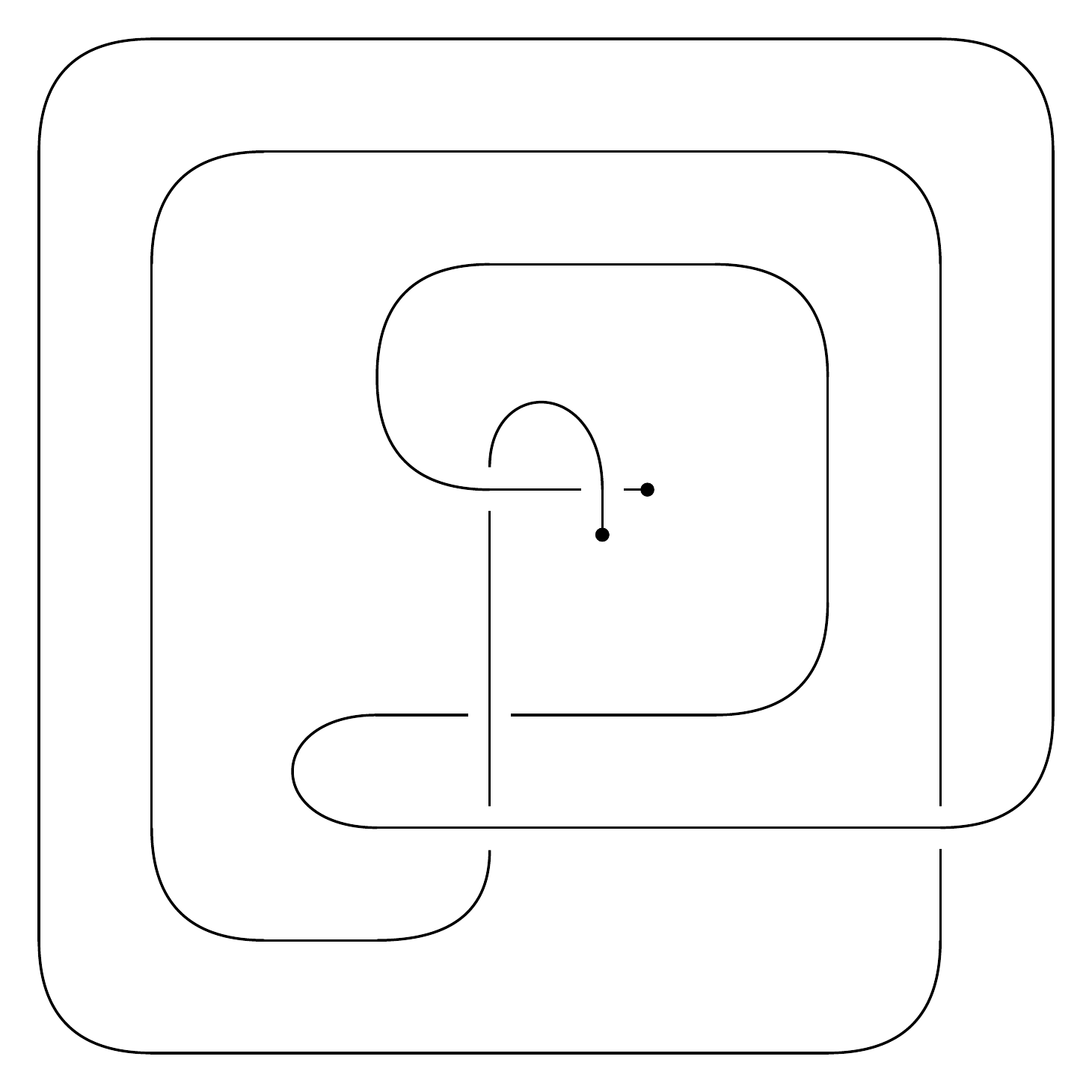}\\
\textcolor{black}{$5_{910}$}
\vspace{1cm}
\end{minipage}
\begin{minipage}[t]{.25\linewidth}
\centering
\includegraphics[width=0.9\textwidth,height=3.5cm,keepaspectratio]{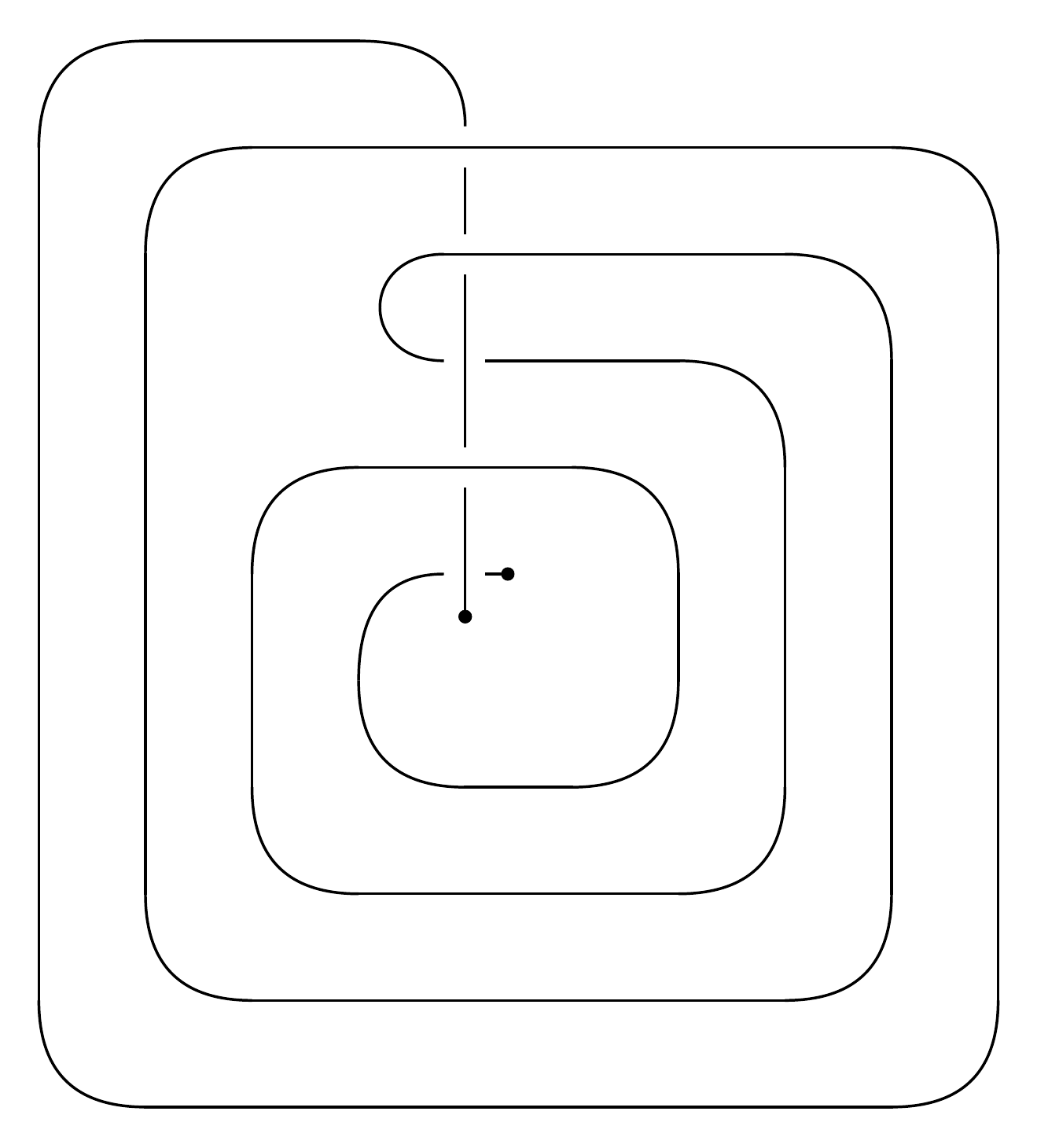}\\
\textcolor{black}{$5_{911}$}
\vspace{1cm}
\end{minipage}
\begin{minipage}[t]{.25\linewidth}
\centering
\includegraphics[width=0.9\textwidth,height=3.5cm,keepaspectratio]{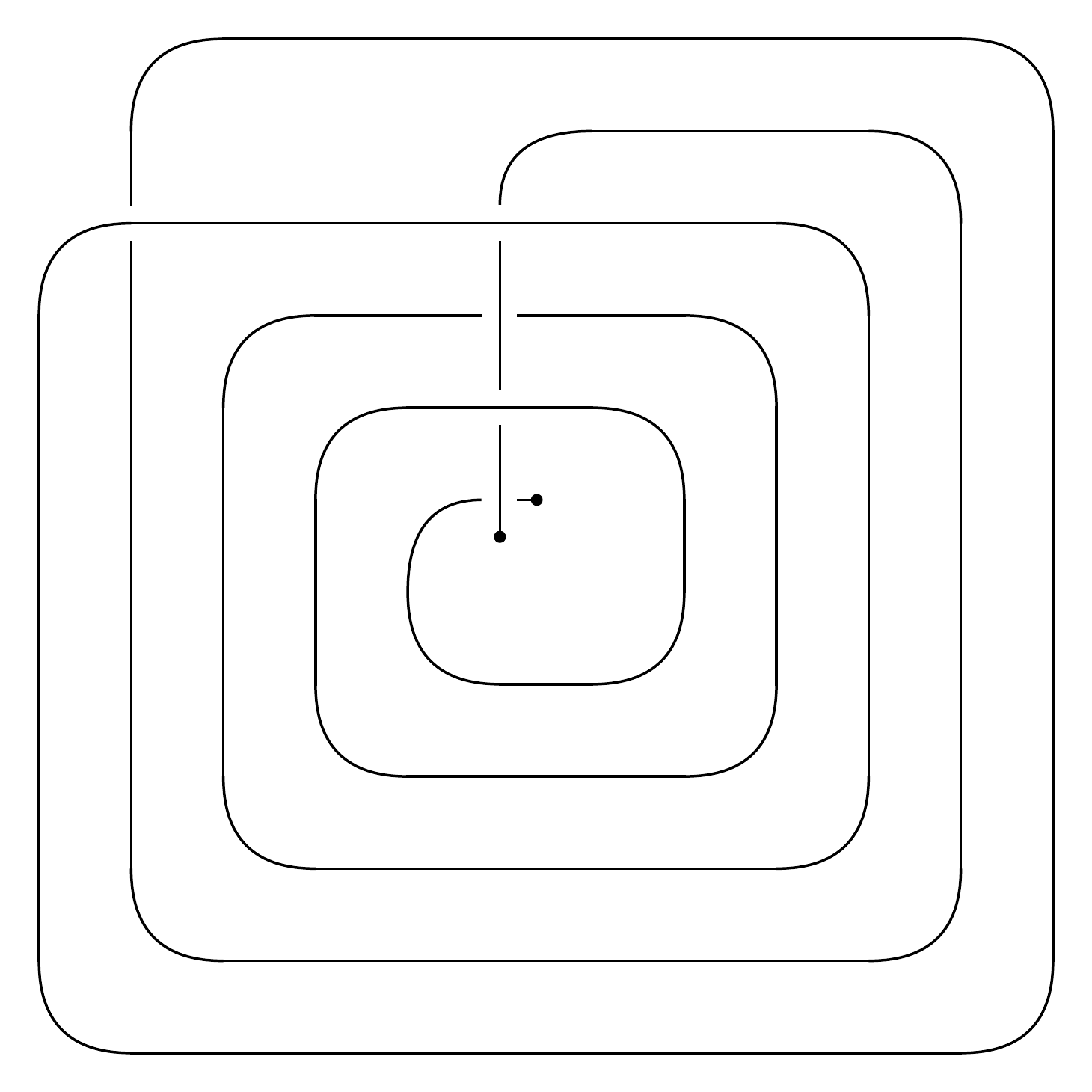}\\
\textcolor{black}{$5_{912}$}
\vspace{1cm}
\end{minipage}
\begin{minipage}[t]{.25\linewidth}
\centering
\includegraphics[width=0.9\textwidth,height=3.5cm,keepaspectratio]{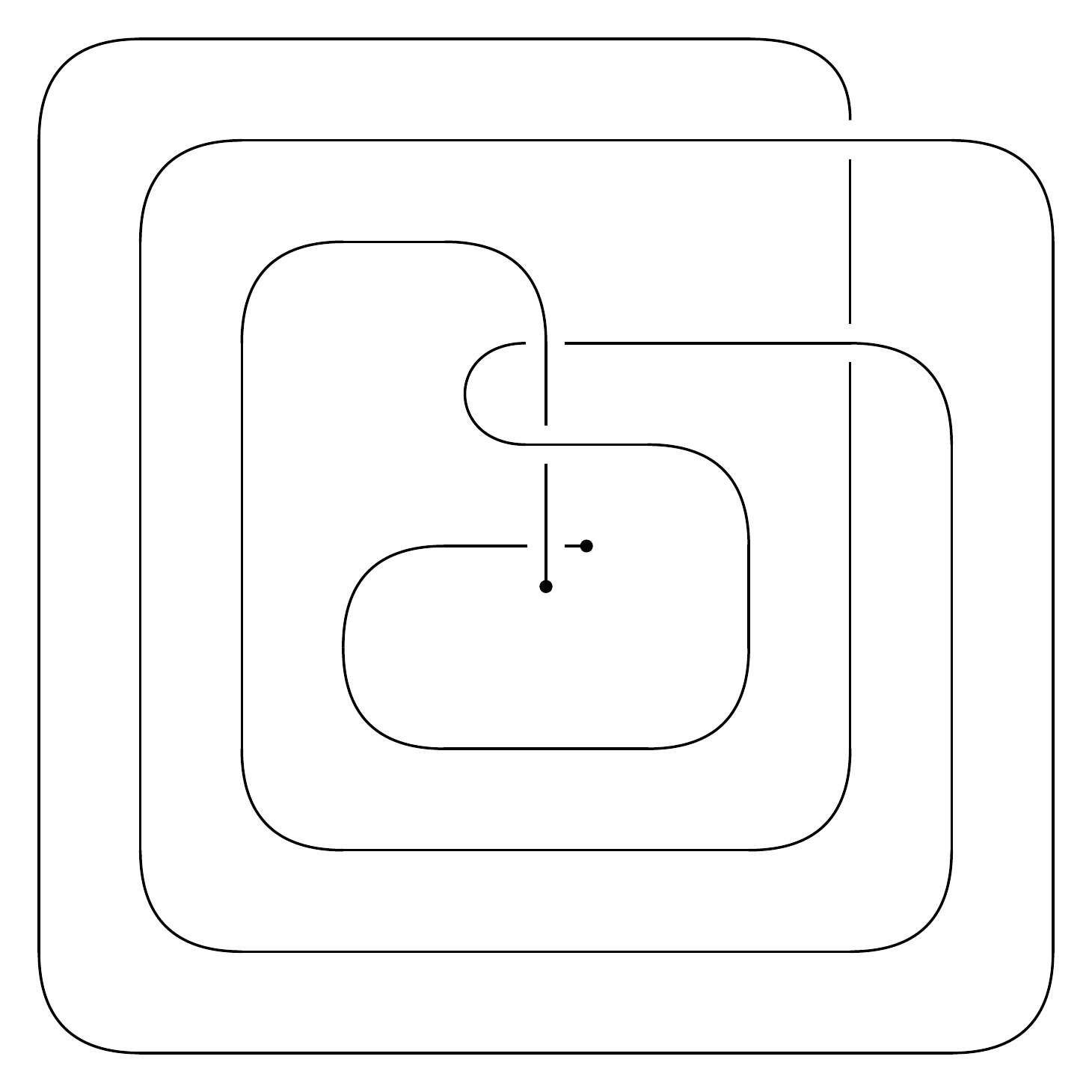}\\
\textcolor{black}{$5_{913}$}
\vspace{1cm}
\end{minipage}
\begin{minipage}[t]{.25\linewidth}
\centering
\includegraphics[width=0.9\textwidth,height=3.5cm,keepaspectratio]{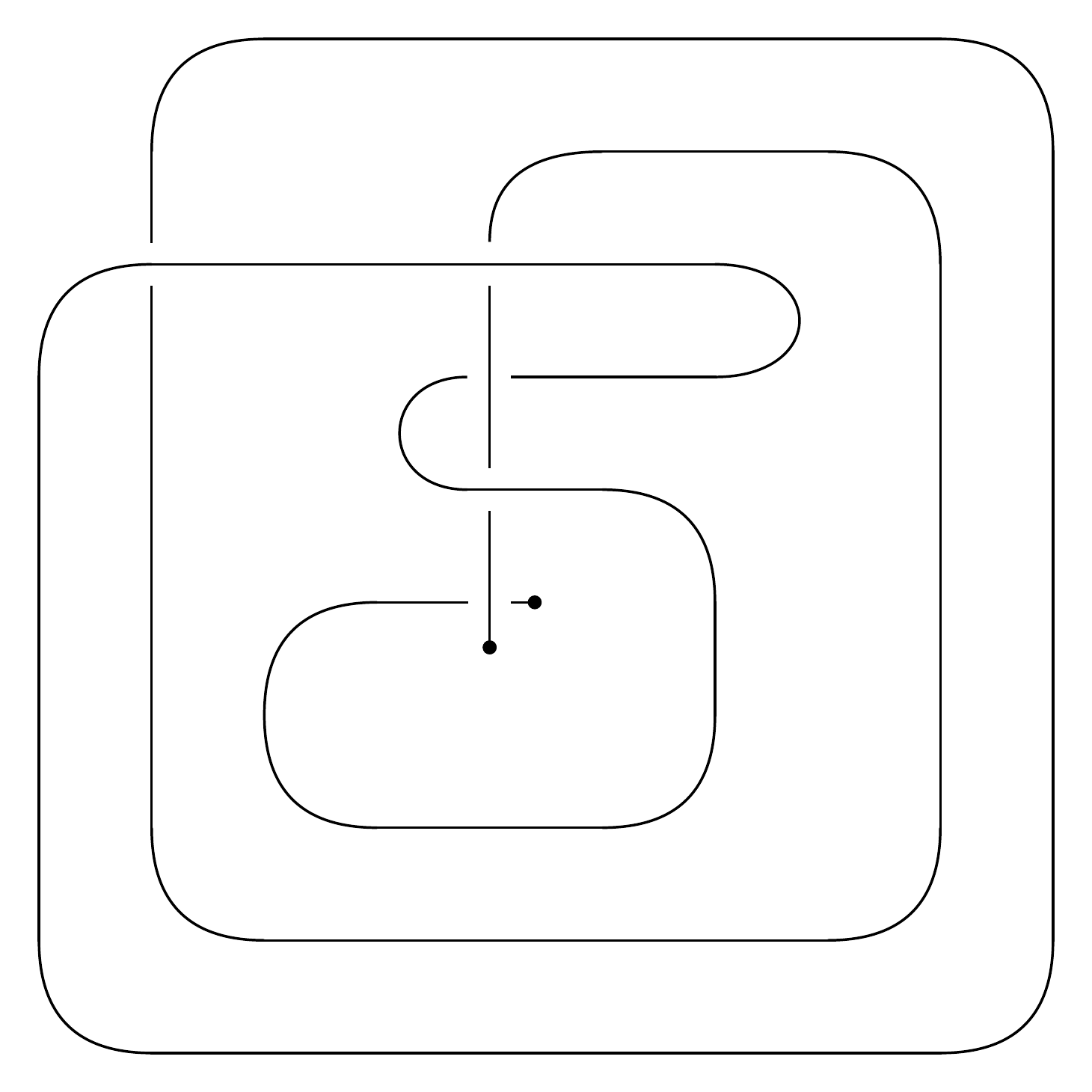}\\
\textcolor{black}{$5_{914}$}
\vspace{1cm}
\end{minipage}
\begin{minipage}[t]{.25\linewidth}
\centering
\includegraphics[width=0.9\textwidth,height=3.5cm,keepaspectratio]{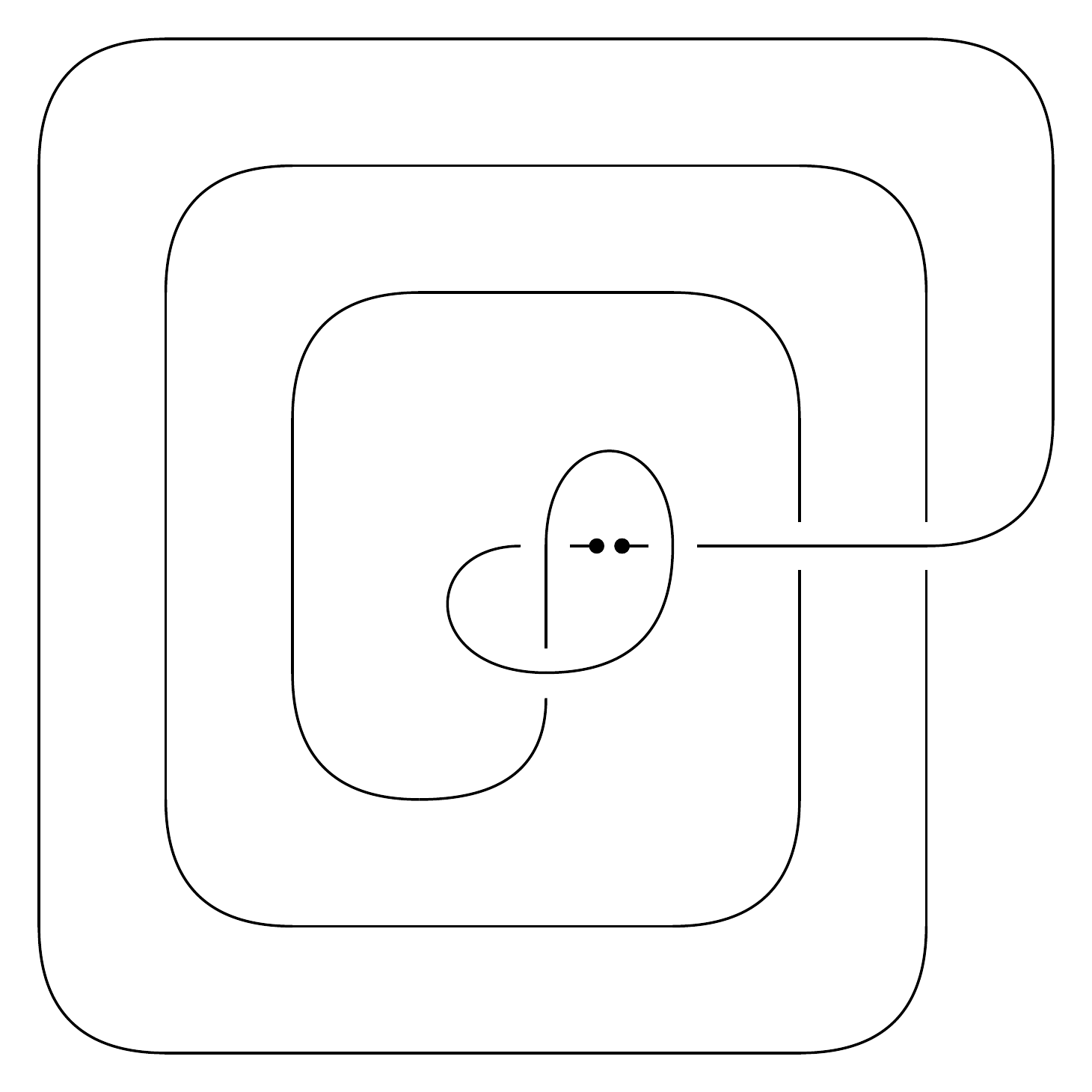}\\
\textcolor{black}{$5_{915}$}
\vspace{1cm}
\end{minipage}
\begin{minipage}[t]{.25\linewidth}
\centering
\includegraphics[width=0.9\textwidth,height=3.5cm,keepaspectratio]{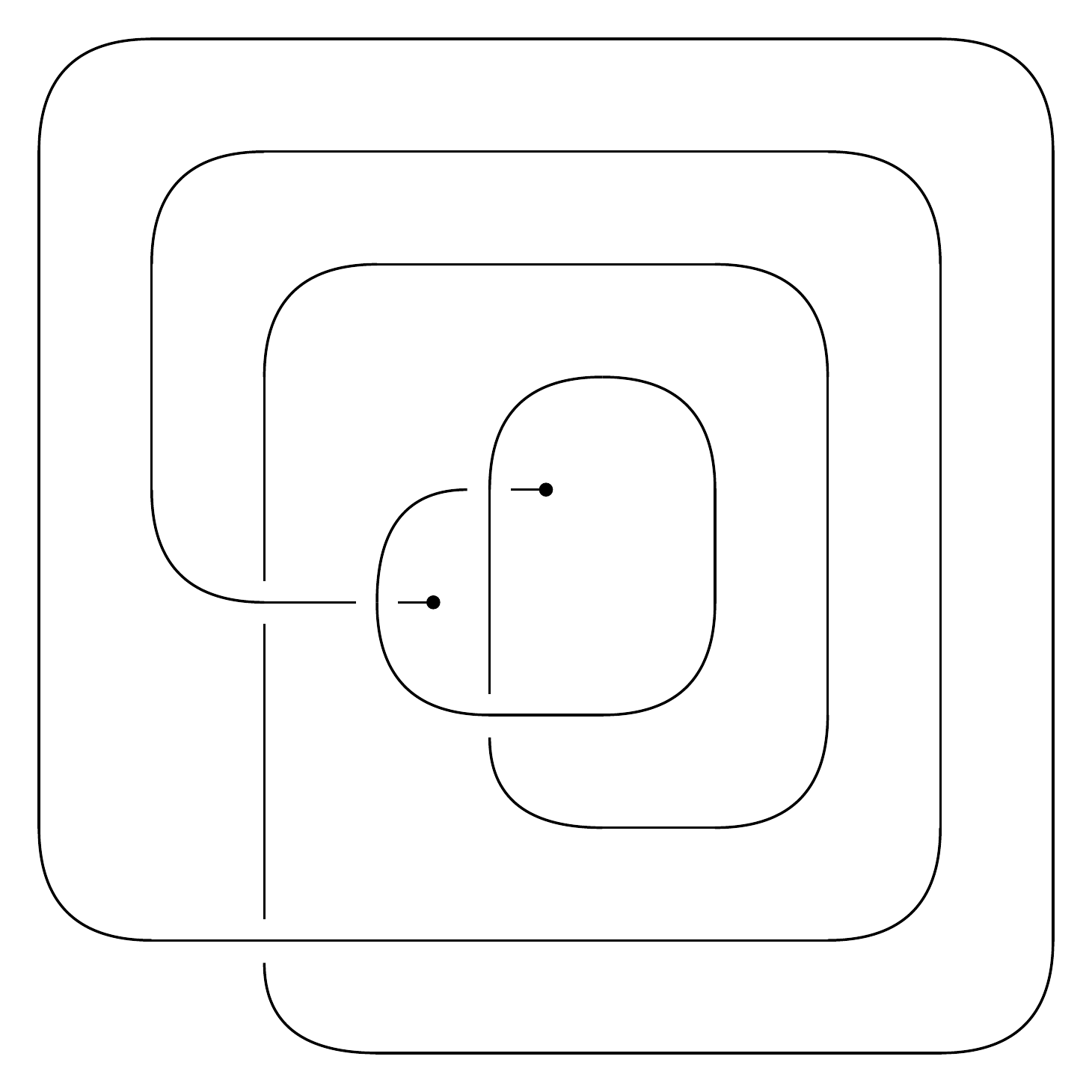}\\
\textcolor{black}{$5_{916}$}
\vspace{1cm}
\end{minipage}
\begin{minipage}[t]{.25\linewidth}
\centering
\includegraphics[width=0.9\textwidth,height=3.5cm,keepaspectratio]{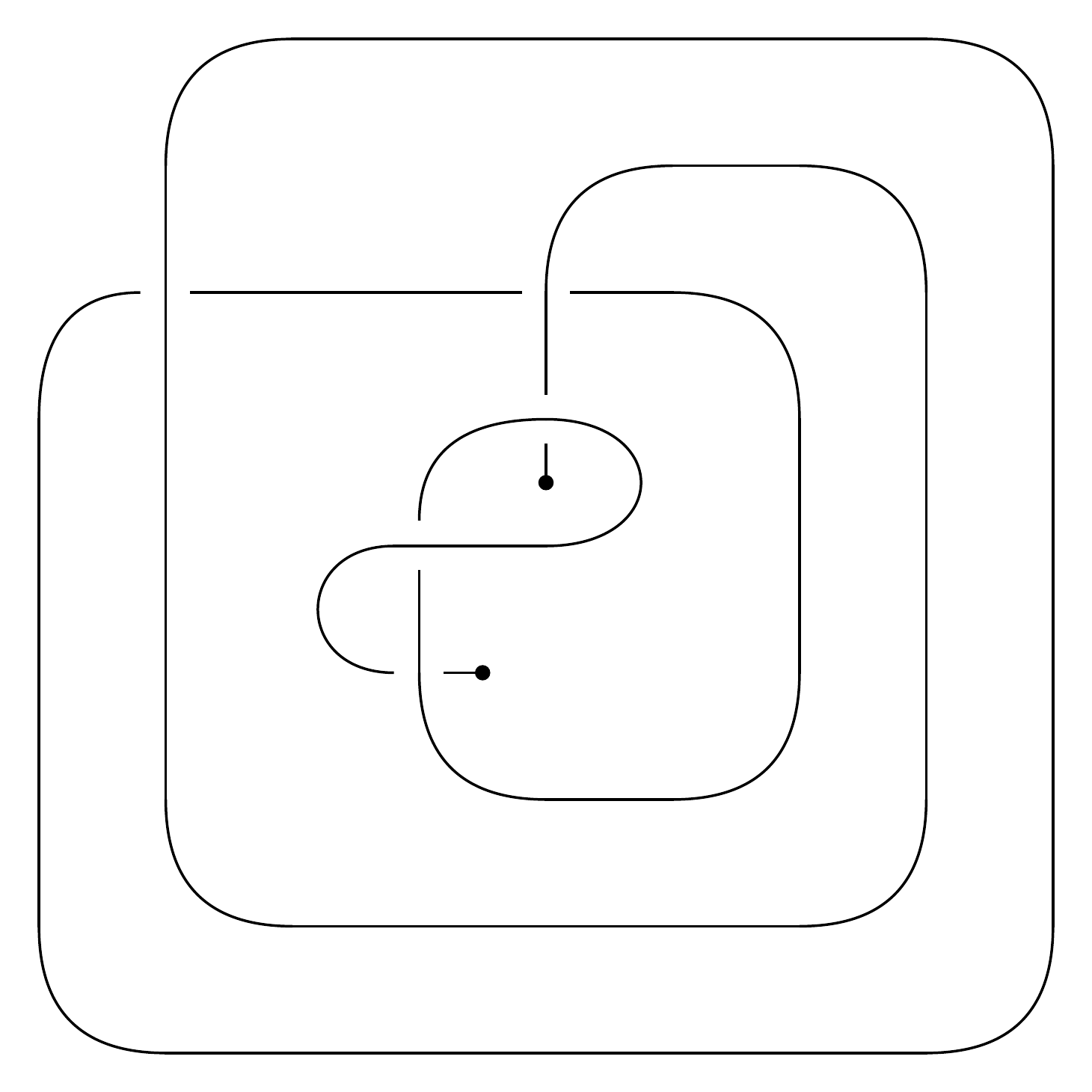}\\
\textcolor{black}{$5_{917}$}
\vspace{1cm}
\end{minipage}
\begin{minipage}[t]{.25\linewidth}
\centering
\includegraphics[width=0.9\textwidth,height=3.5cm,keepaspectratio]{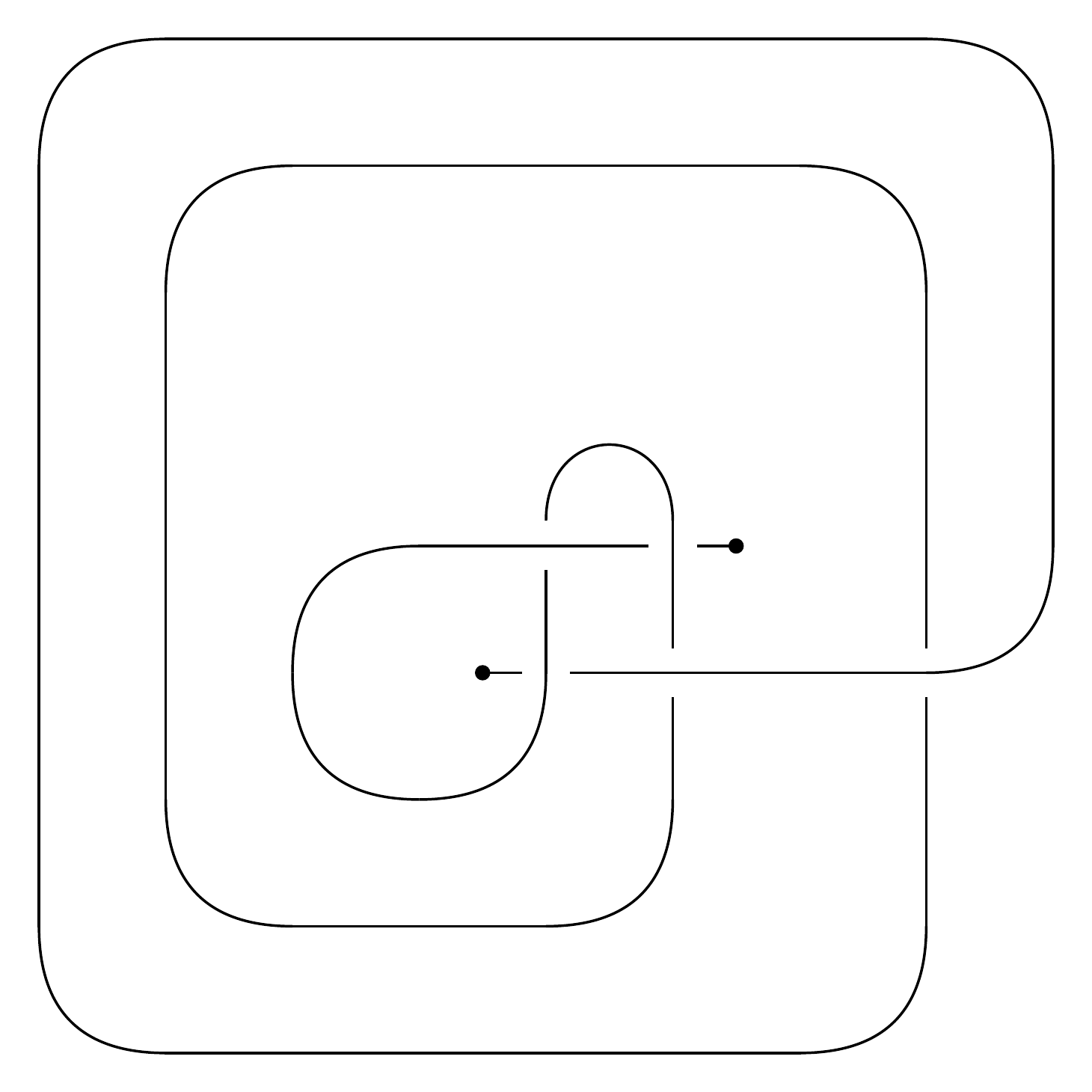}\\
\textcolor{black}{$5_{918}$}
\vspace{1cm}
\end{minipage}
\begin{minipage}[t]{.25\linewidth}
\centering
\includegraphics[width=0.9\textwidth,height=3.5cm,keepaspectratio]{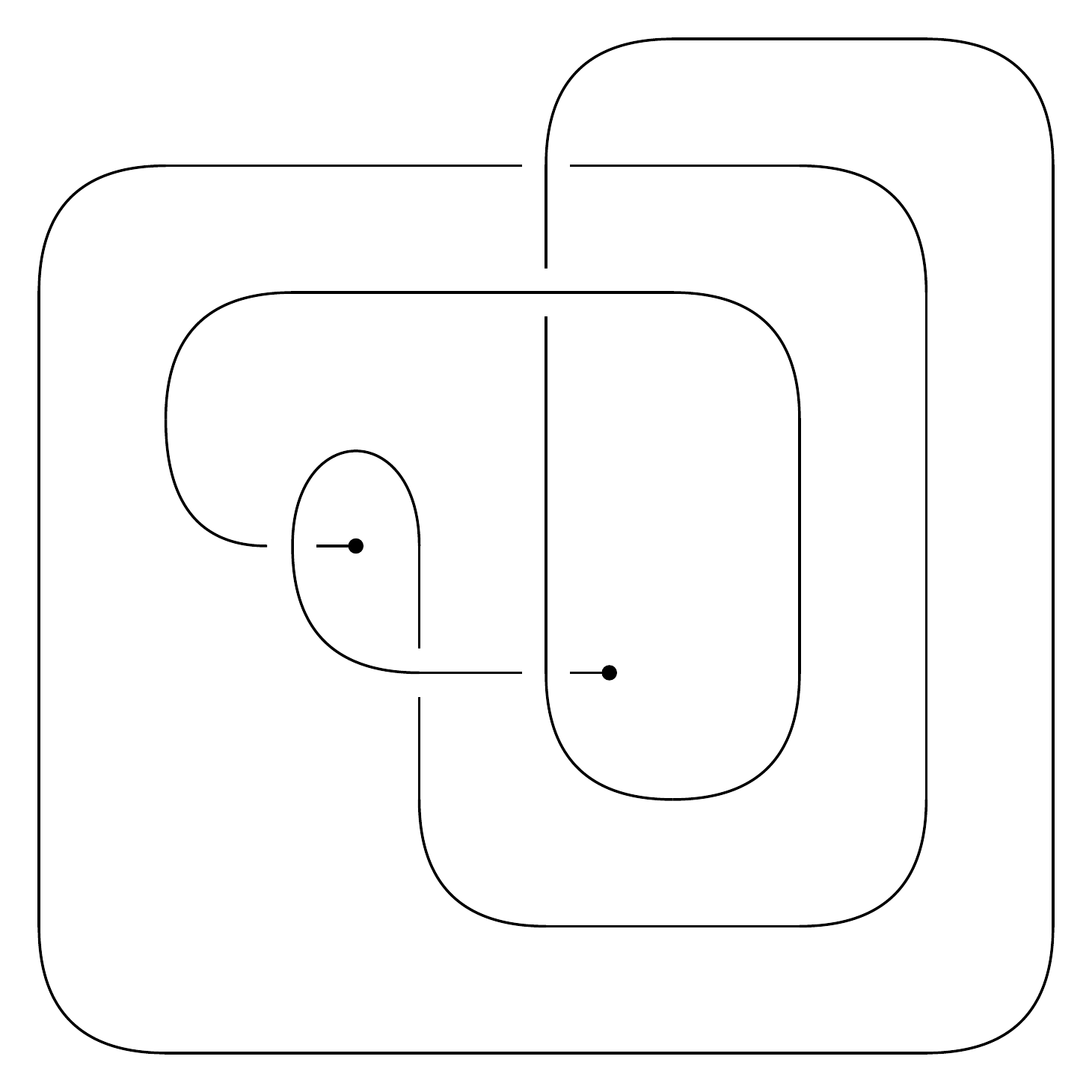}\\
\textcolor{black}{$5_{919}$}
\vspace{1cm}
\end{minipage}
\begin{minipage}[t]{.25\linewidth}
\centering
\includegraphics[width=0.9\textwidth,height=3.5cm,keepaspectratio]{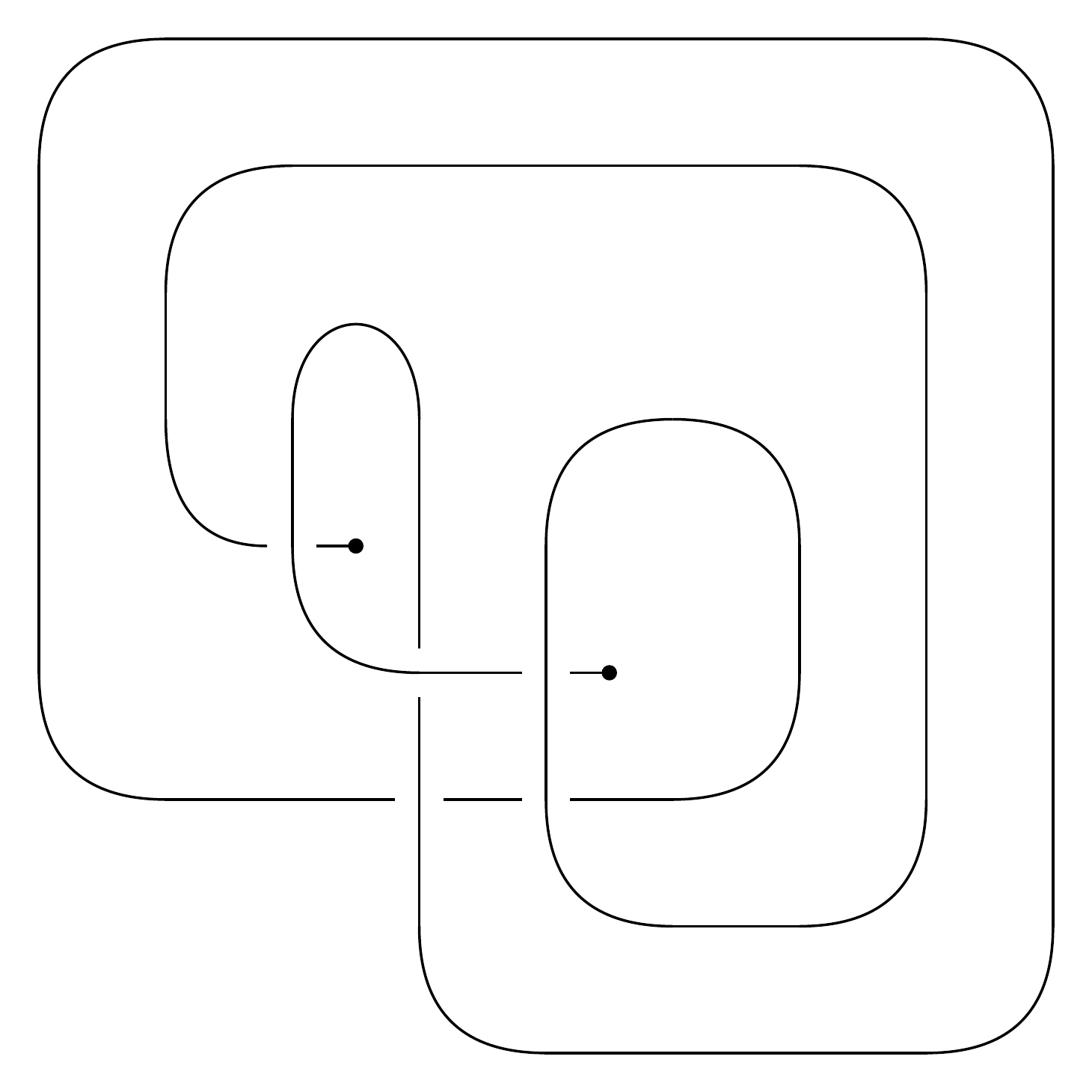}\\
\textcolor{black}{$5_{920}$}
\vspace{1cm}
\end{minipage}
\begin{minipage}[t]{.25\linewidth}
\centering
\includegraphics[width=0.9\textwidth,height=3.5cm,keepaspectratio]{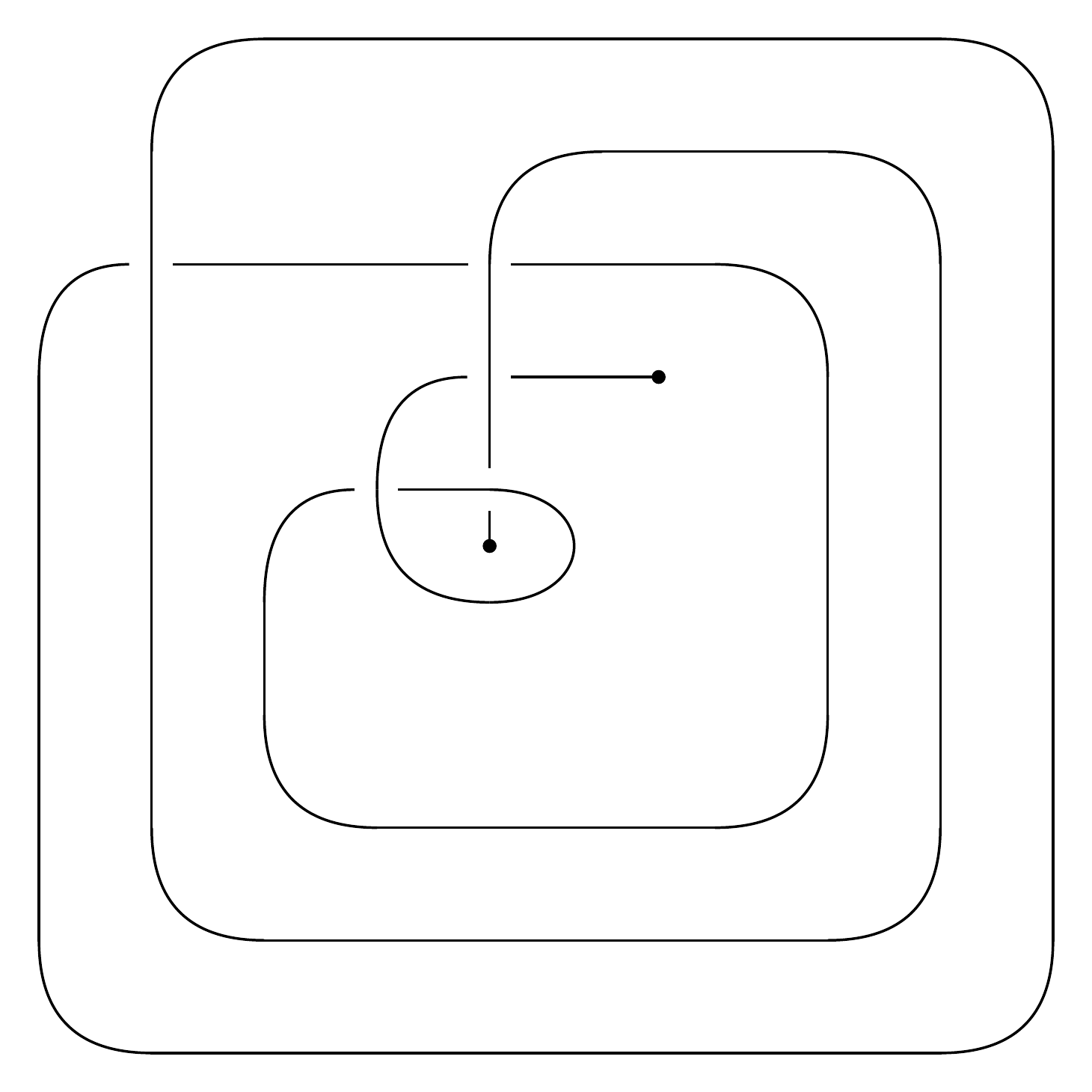}\\
\textcolor{black}{$5_{921}$}
\vspace{1cm}
\end{minipage}
\begin{minipage}[t]{.25\linewidth}
\centering
\includegraphics[width=0.9\textwidth,height=3.5cm,keepaspectratio]{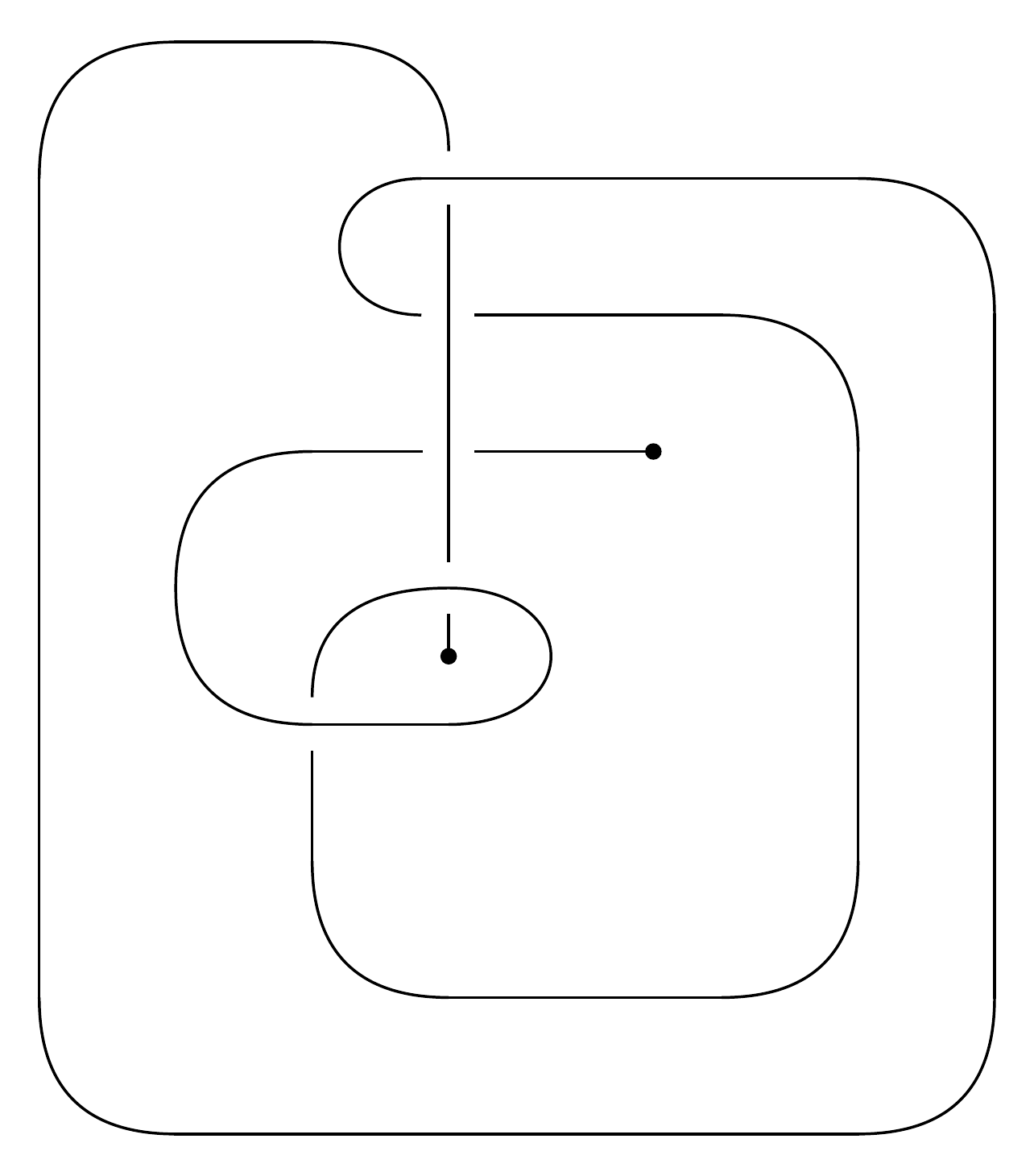}\\
\textcolor{black}{$5_{922}$}
\vspace{1cm}
\end{minipage}
\begin{minipage}[t]{.25\linewidth}
\centering
\includegraphics[width=0.9\textwidth,height=3.5cm,keepaspectratio]{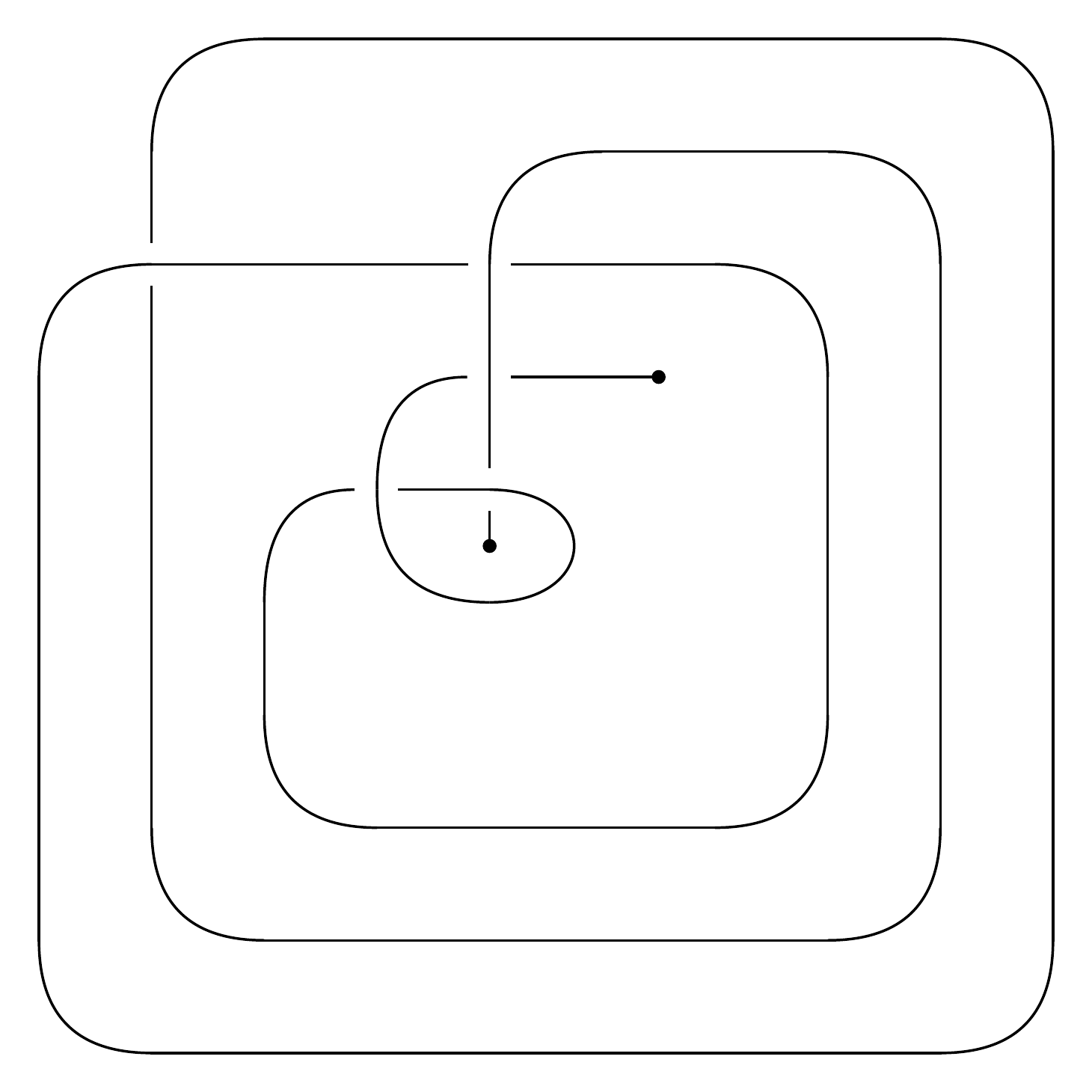}\\
\textcolor{black}{$5_{923}$}
\vspace{1cm}
\end{minipage}
\begin{minipage}[t]{.25\linewidth}
\centering
\includegraphics[width=0.9\textwidth,height=3.5cm,keepaspectratio]{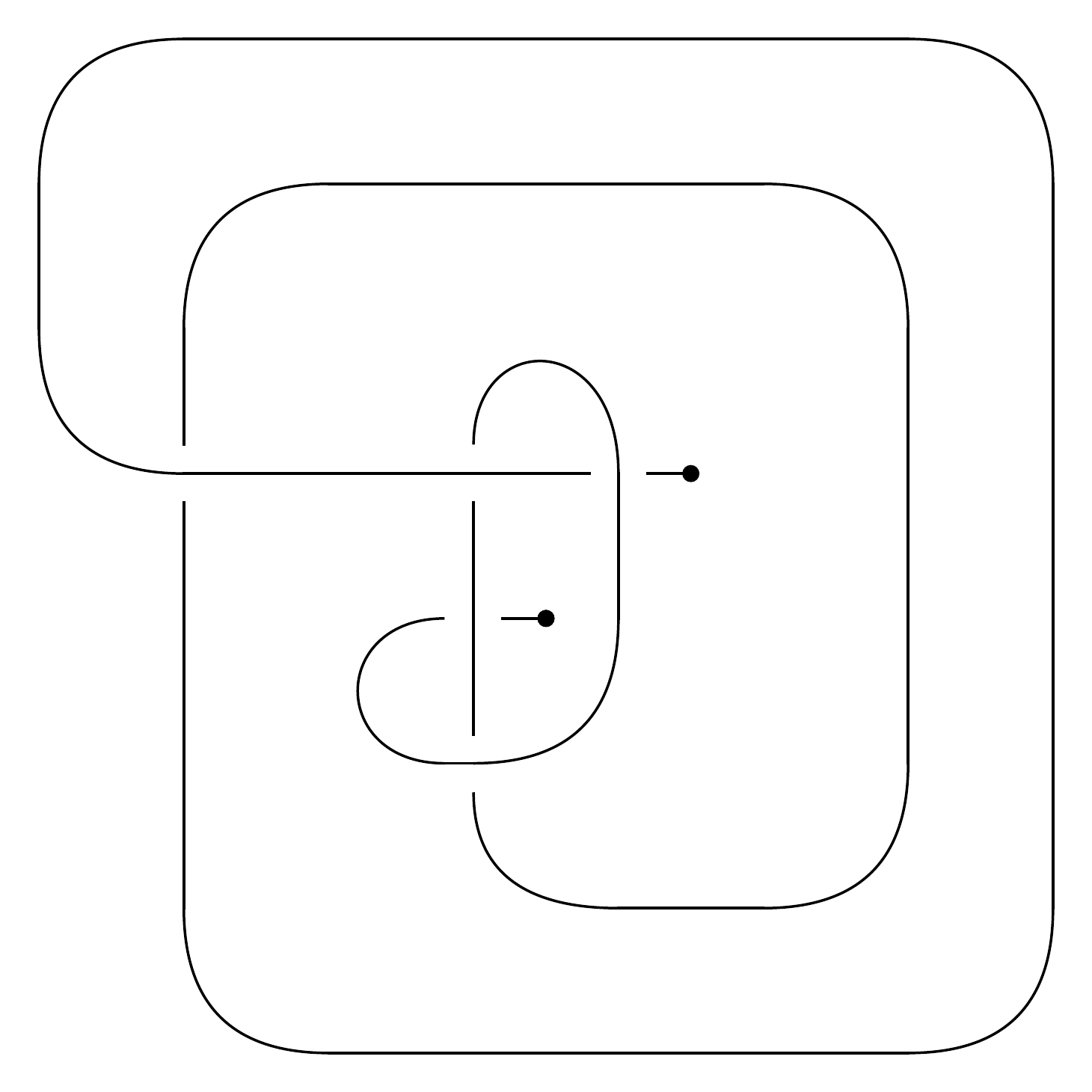}\\
\textcolor{black}{$5_{924}$}
\vspace{1cm}
\end{minipage}
\begin{minipage}[t]{.25\linewidth}
\centering
\includegraphics[width=0.9\textwidth,height=3.5cm,keepaspectratio]{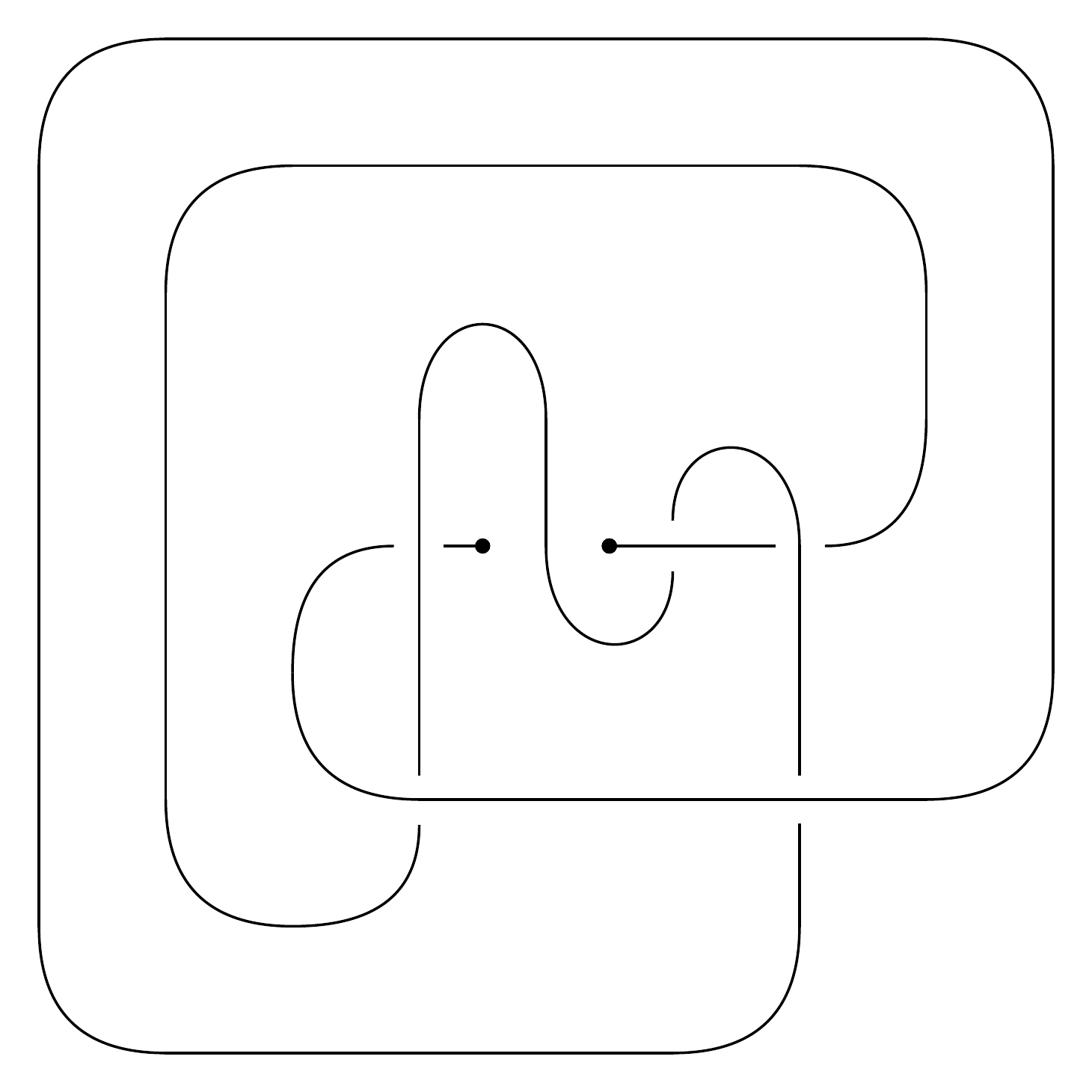}\\
\textcolor{black}{$5_{925}$}
\vspace{1cm}
\end{minipage}
\begin{minipage}[t]{.25\linewidth}
\centering
\includegraphics[width=0.9\textwidth,height=3.5cm,keepaspectratio]{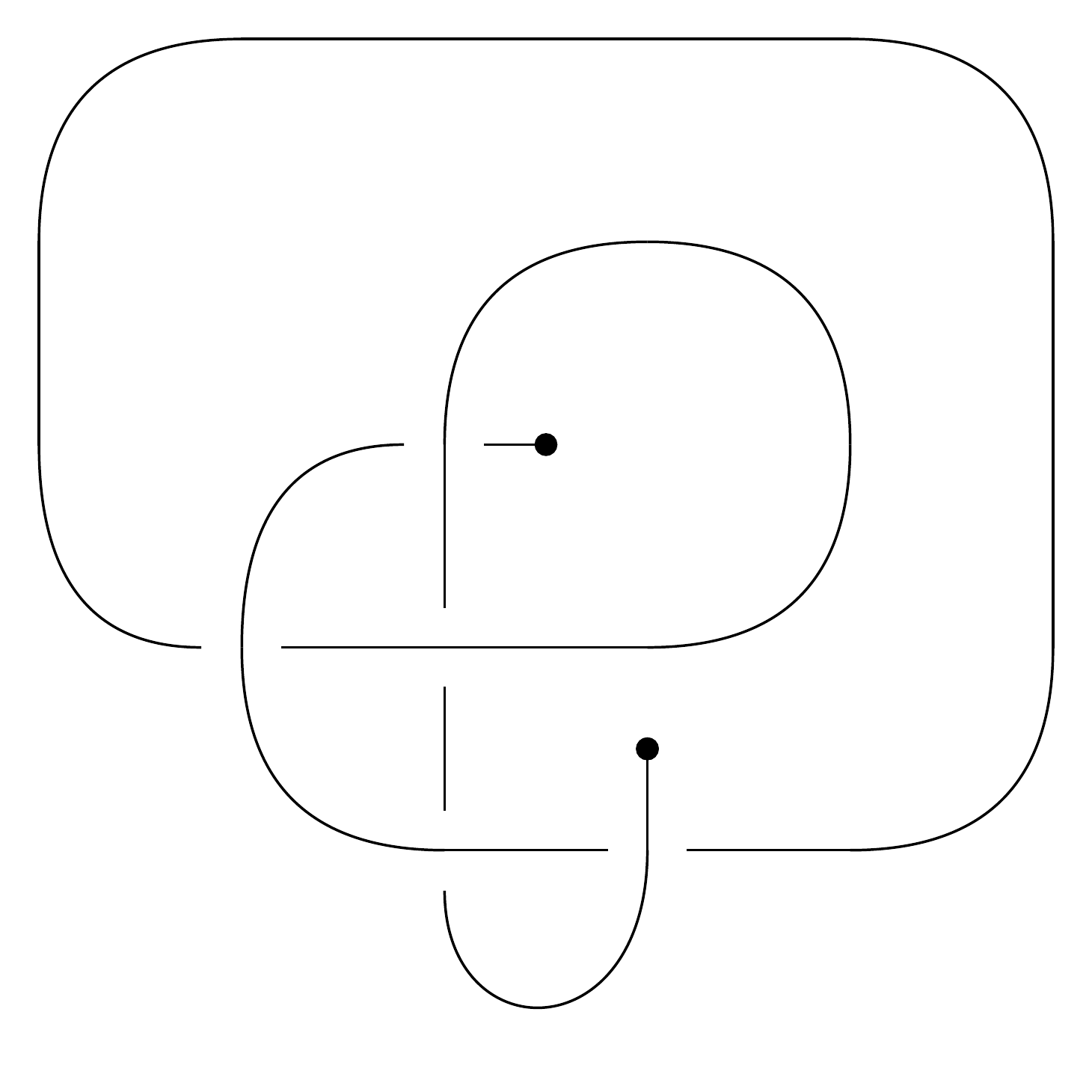}\\
\textcolor{black}{$5_{926}$}
\vspace{1cm}
\end{minipage}
\begin{minipage}[t]{.25\linewidth}
\centering
\includegraphics[width=0.9\textwidth,height=3.5cm,keepaspectratio]{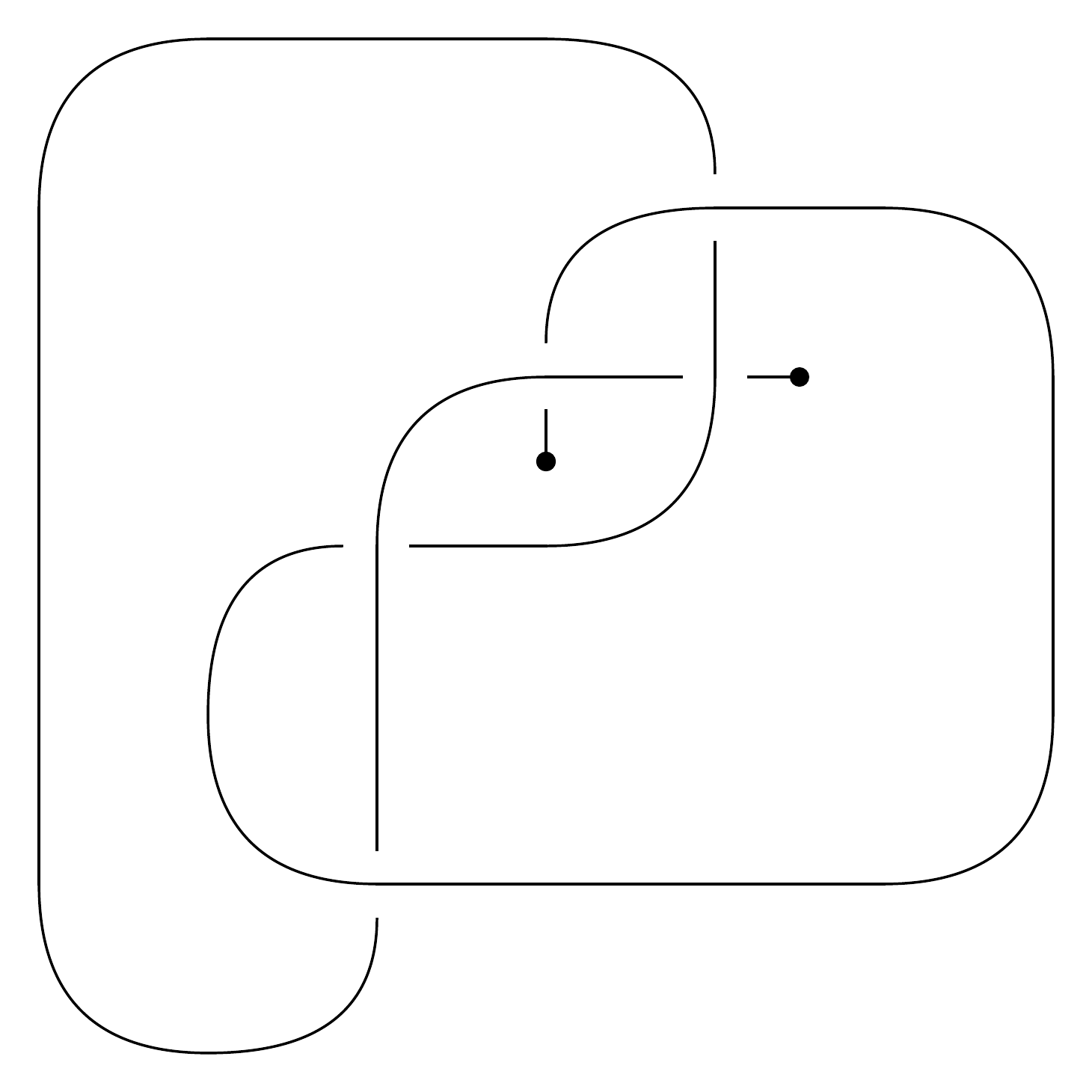}\\
\textcolor{black}{$5_{927}$}
\vspace{1cm}
\end{minipage}
\begin{minipage}[t]{.25\linewidth}
\centering
\includegraphics[width=0.9\textwidth,height=3.5cm,keepaspectratio]{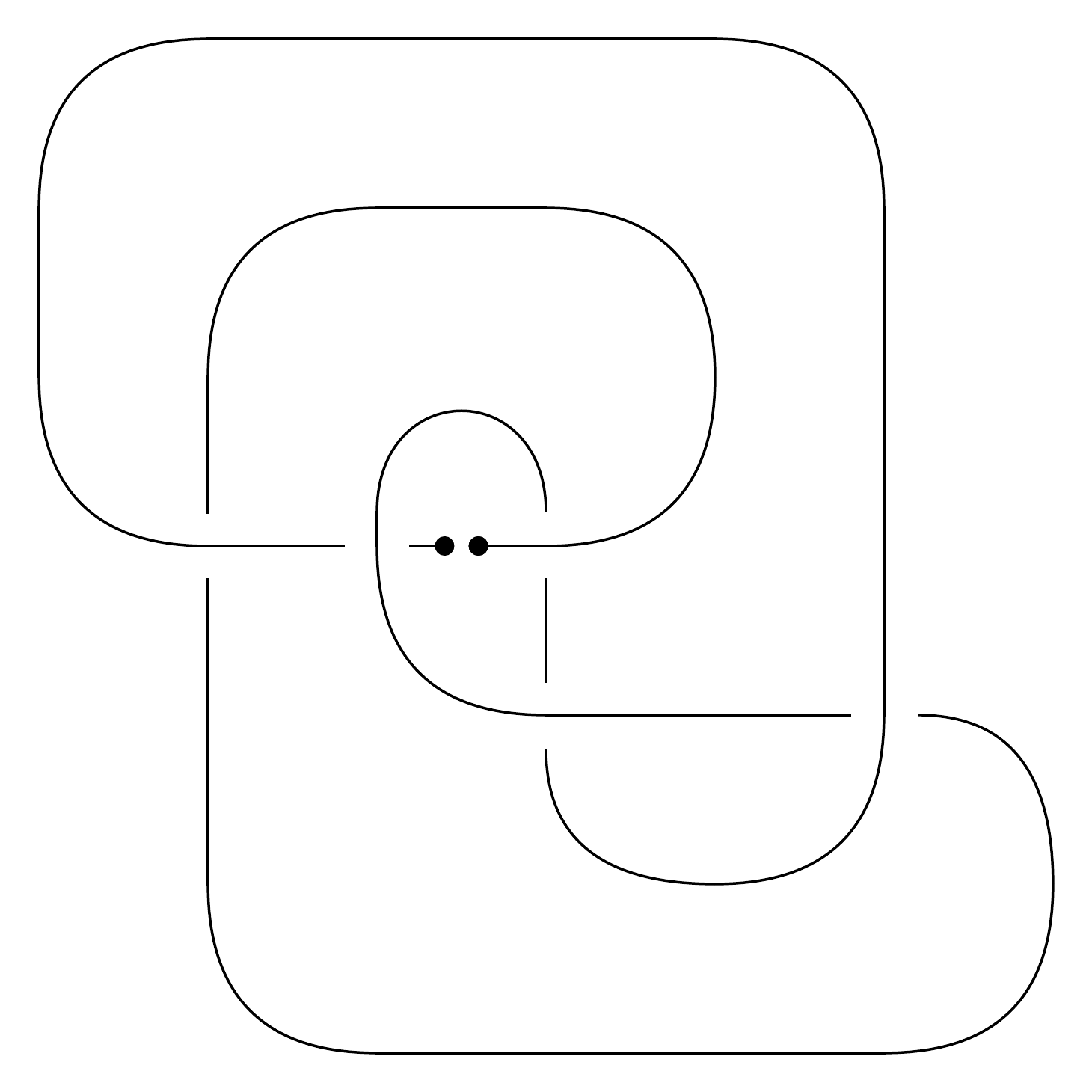}\\
\textcolor{black}{$5_{928}$}
\vspace{1cm}
\end{minipage}
\begin{minipage}[t]{.25\linewidth}
\centering
\includegraphics[width=0.9\textwidth,height=3.5cm,keepaspectratio]{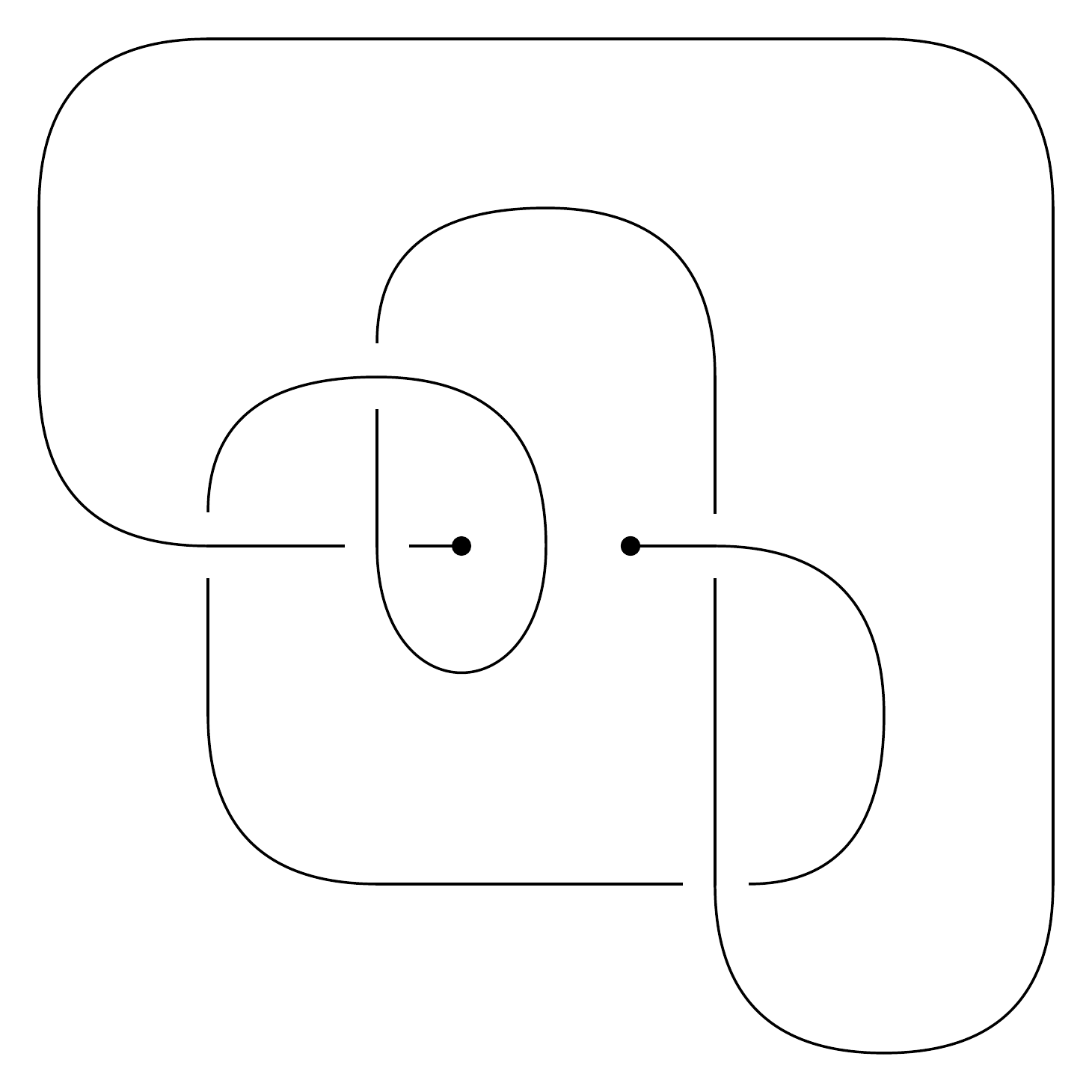}\\
\textcolor{black}{$5_{929}$}
\vspace{1cm}
\end{minipage}
\begin{minipage}[t]{.25\linewidth}
\centering
\includegraphics[width=0.9\textwidth,height=3.5cm,keepaspectratio]{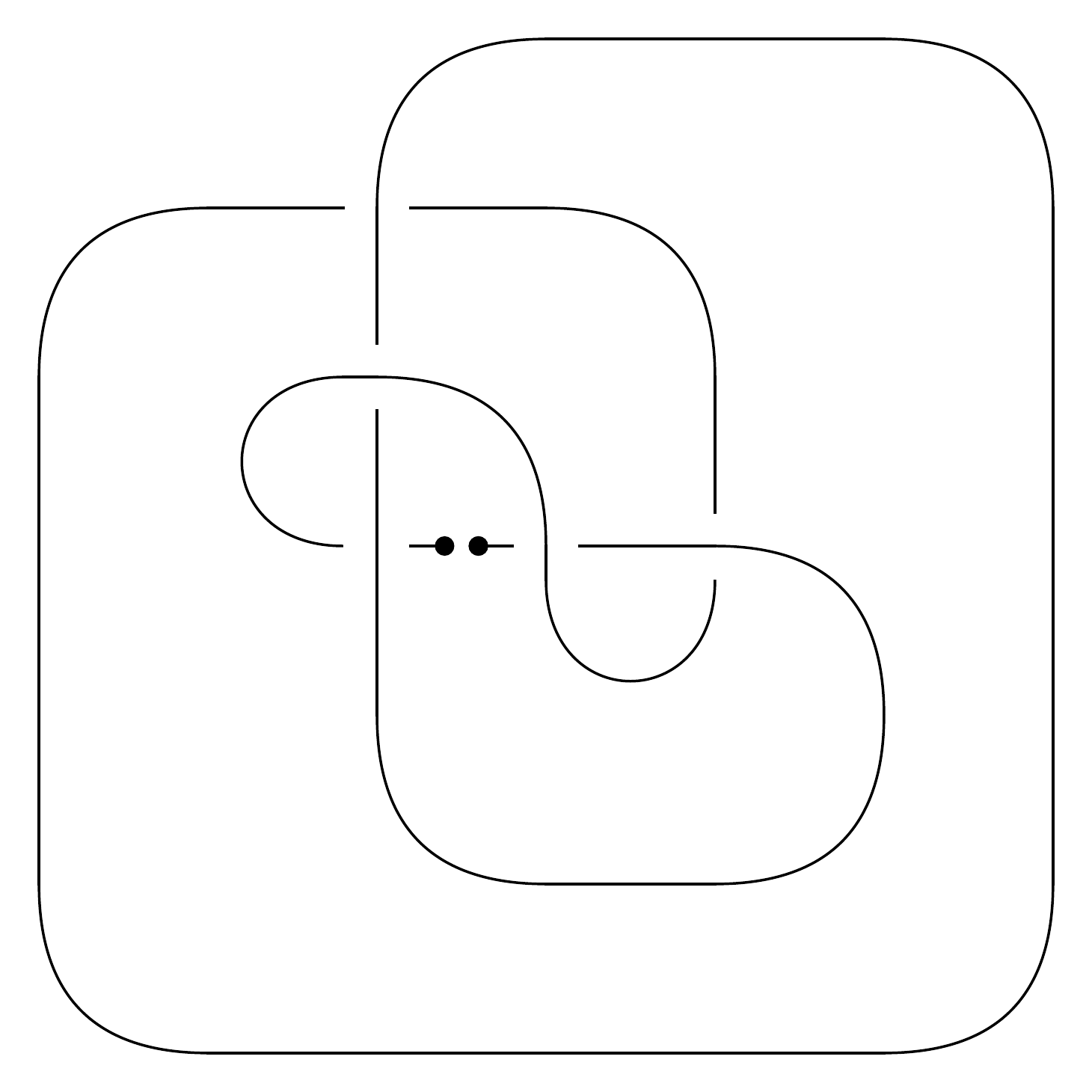}\\
\textcolor{black}{$5_{930}$}
\vspace{1cm}
\end{minipage}
\begin{minipage}[t]{.25\linewidth}
\centering
\includegraphics[width=0.9\textwidth,height=3.5cm,keepaspectratio]{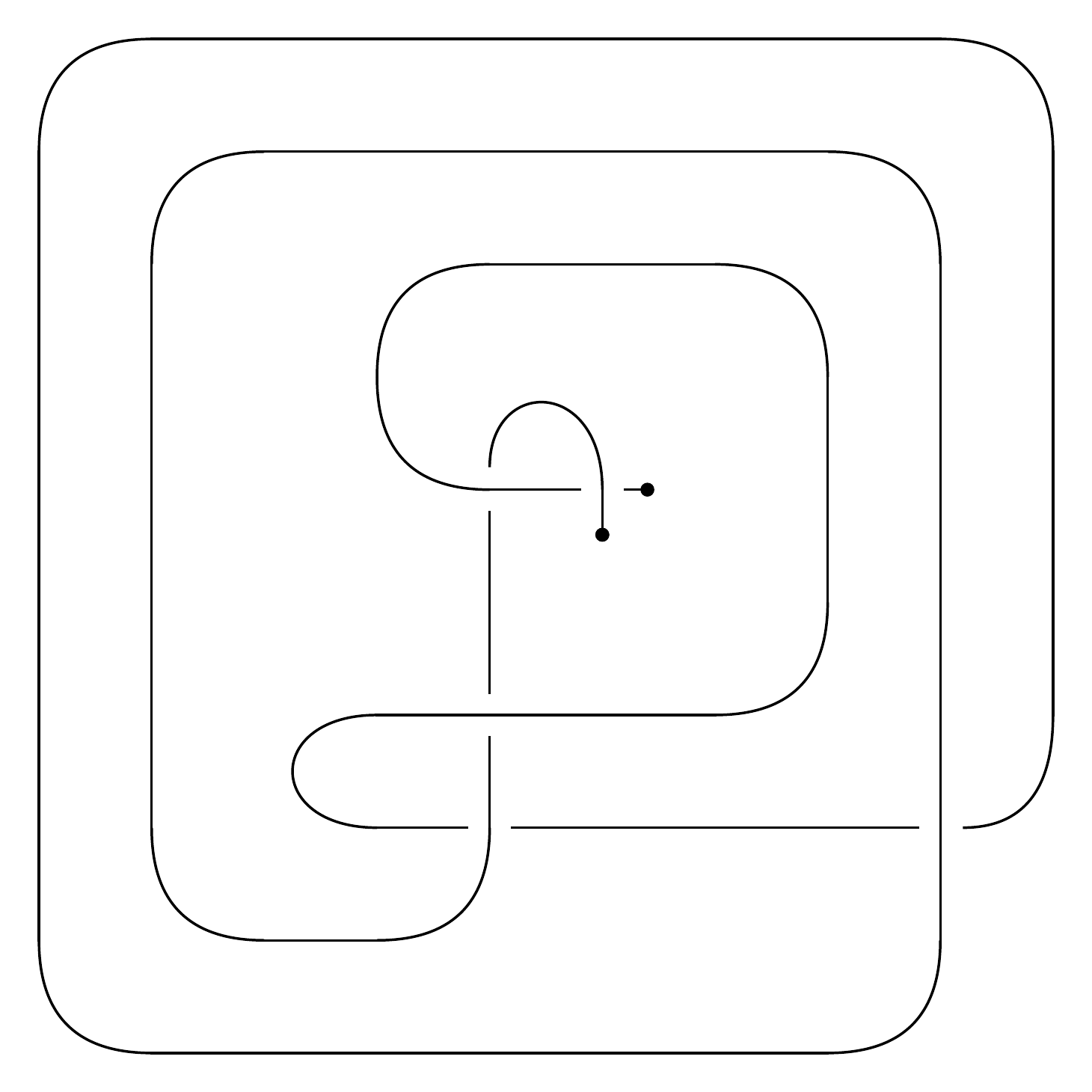}\\
\textcolor{black}{$5_{931}$}
\vspace{1cm}
\end{minipage}
\begin{minipage}[t]{.25\linewidth}
\centering
\includegraphics[width=0.9\textwidth,height=3.5cm,keepaspectratio]{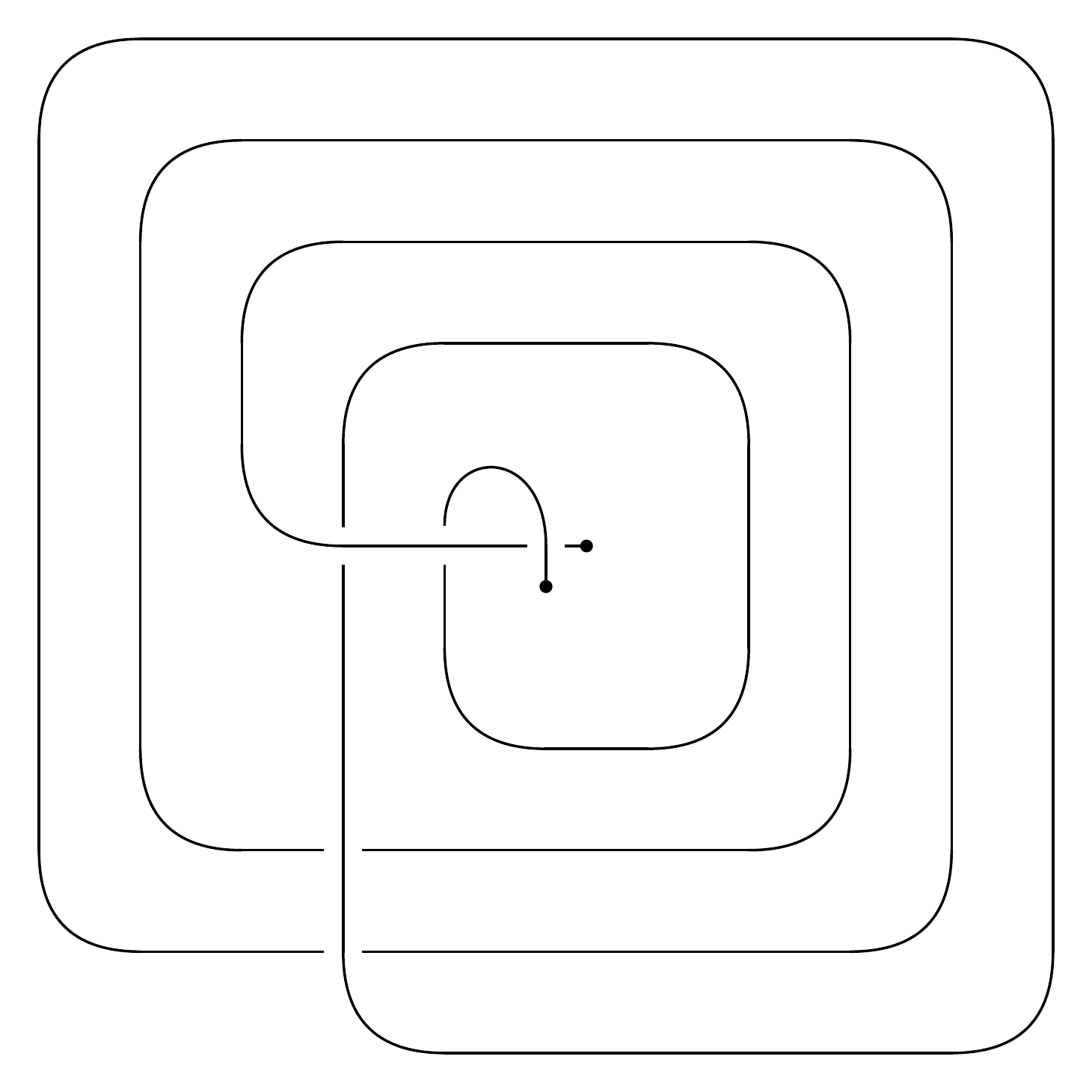}\\
\textcolor{black}{$5_{932}$}
\vspace{1cm}
\end{minipage}
\begin{minipage}[t]{.25\linewidth}
\centering
\includegraphics[width=0.9\textwidth,height=3.5cm,keepaspectratio]{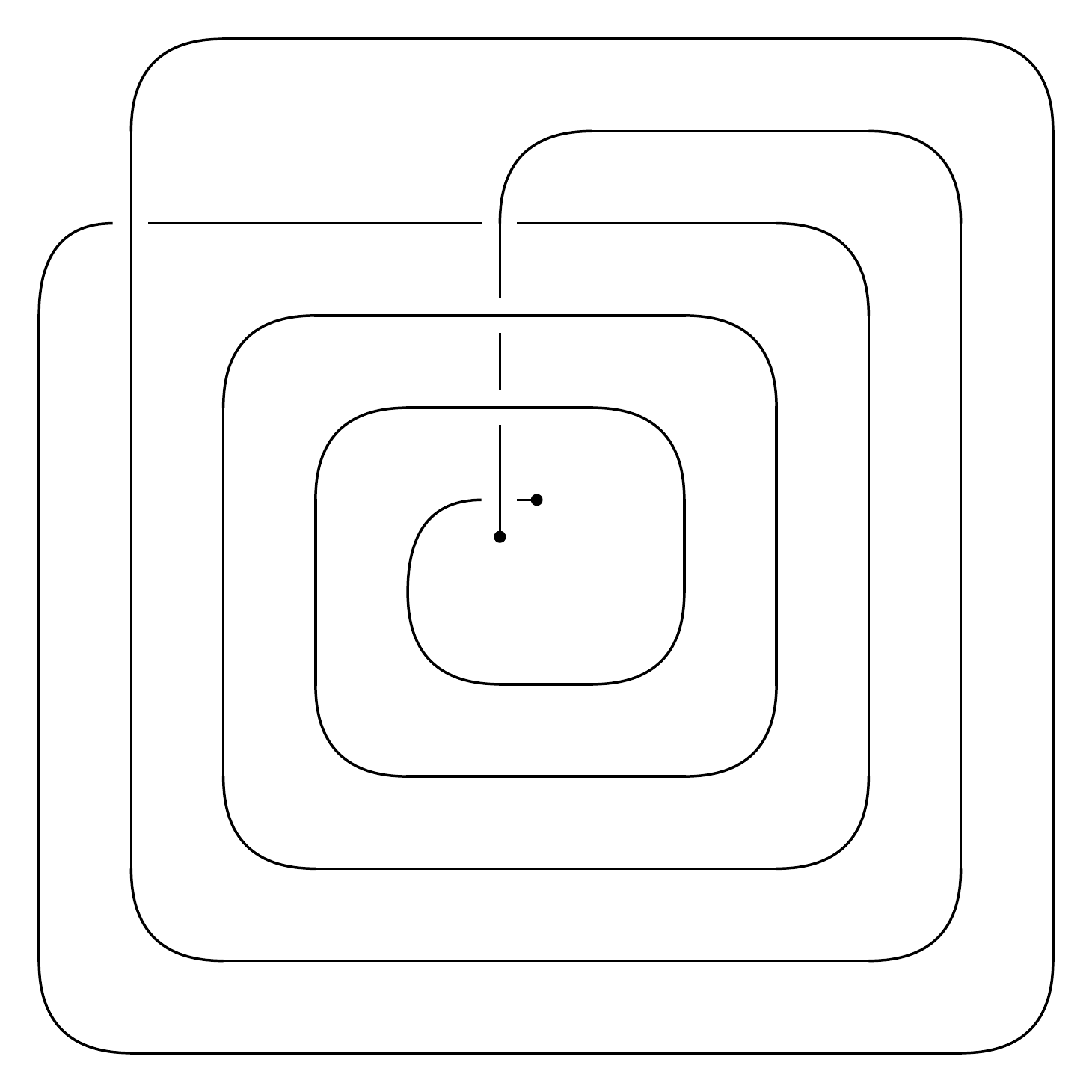}\\
\textcolor{black}{$5_{933}$}
\vspace{1cm}
\end{minipage}
\begin{minipage}[t]{.25\linewidth}
\centering
\includegraphics[width=0.9\textwidth,height=3.5cm,keepaspectratio]{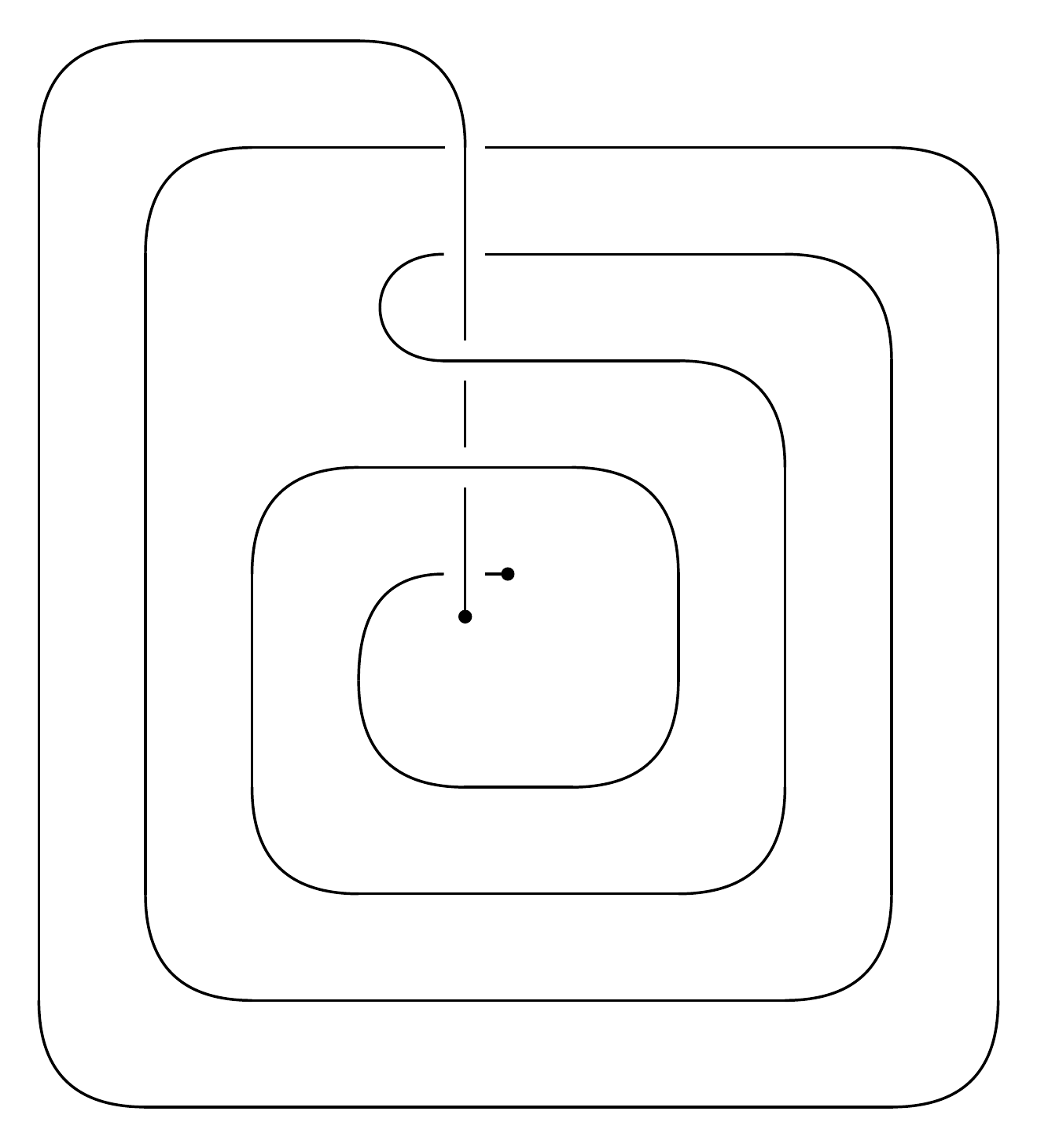}\\
\textcolor{black}{$5_{934}$}
\vspace{1cm}
\end{minipage}
\begin{minipage}[t]{.25\linewidth}
\centering
\includegraphics[width=0.9\textwidth,height=3.5cm,keepaspectratio]{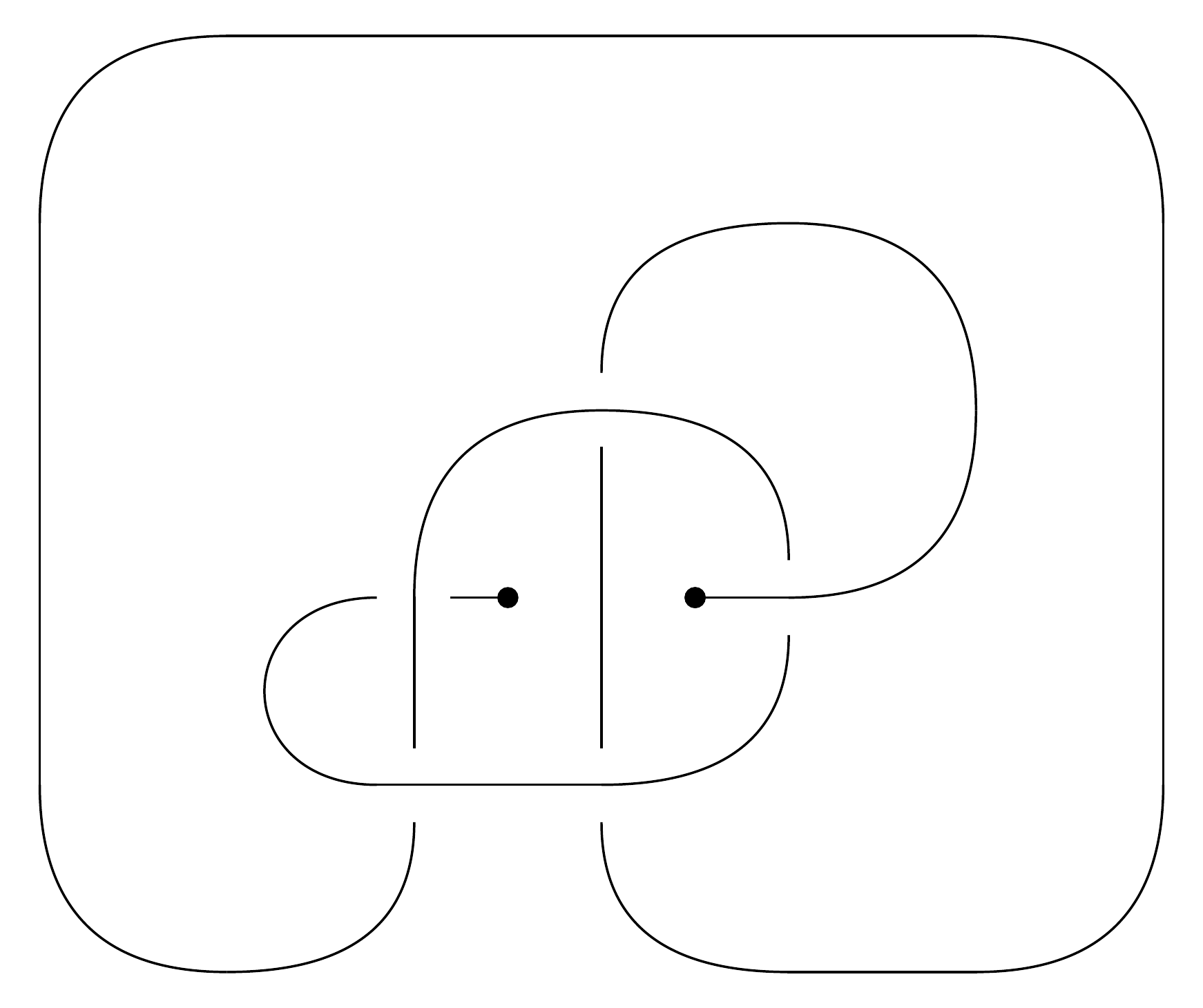}\\
\textcolor{black}{$5_{935}$}
\vspace{1cm}
\end{minipage}
\begin{minipage}[t]{.25\linewidth}
\centering
\includegraphics[width=0.9\textwidth,height=3.5cm,keepaspectratio]{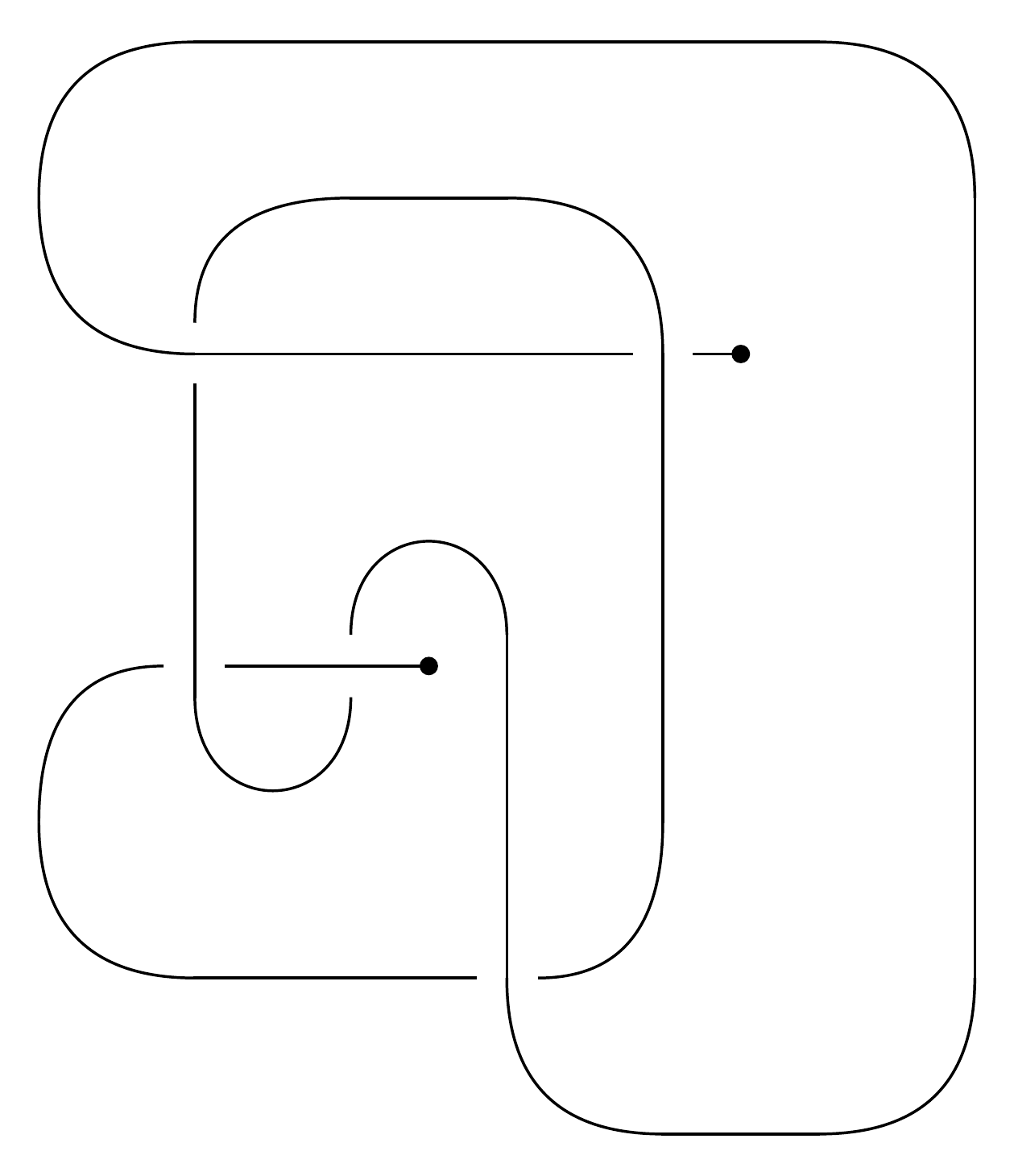}\\
\textcolor{black}{$5_{936}$}
\vspace{1cm}
\end{minipage}
\begin{minipage}[t]{.25\linewidth}
\centering
\includegraphics[width=0.9\textwidth,height=3.5cm,keepaspectratio]{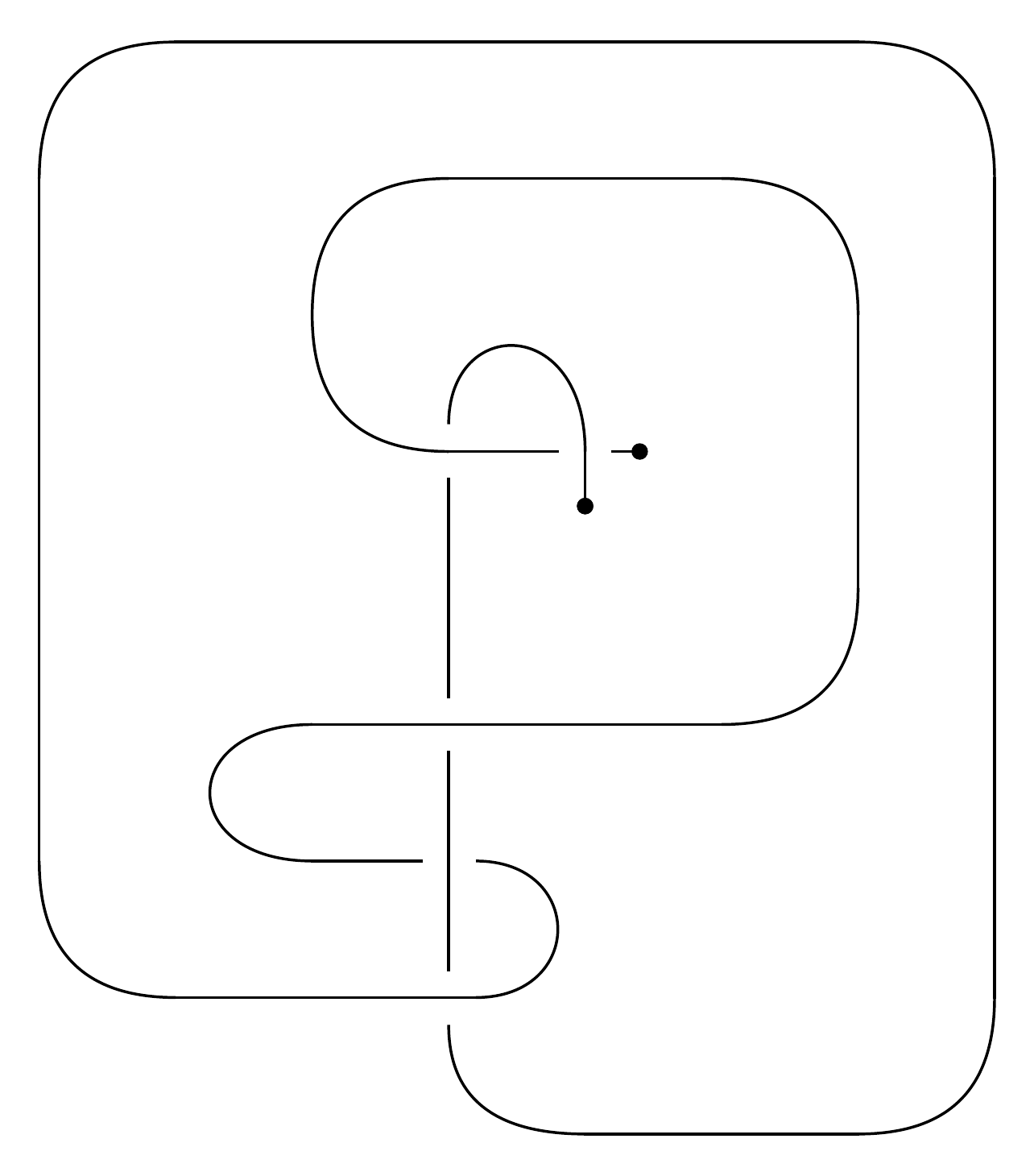}\\
\textcolor{black}{$5_{937}$}
\vspace{1cm}
\end{minipage}
\begin{minipage}[t]{.25\linewidth}
\centering
\includegraphics[width=0.9\textwidth,height=3.5cm,keepaspectratio]{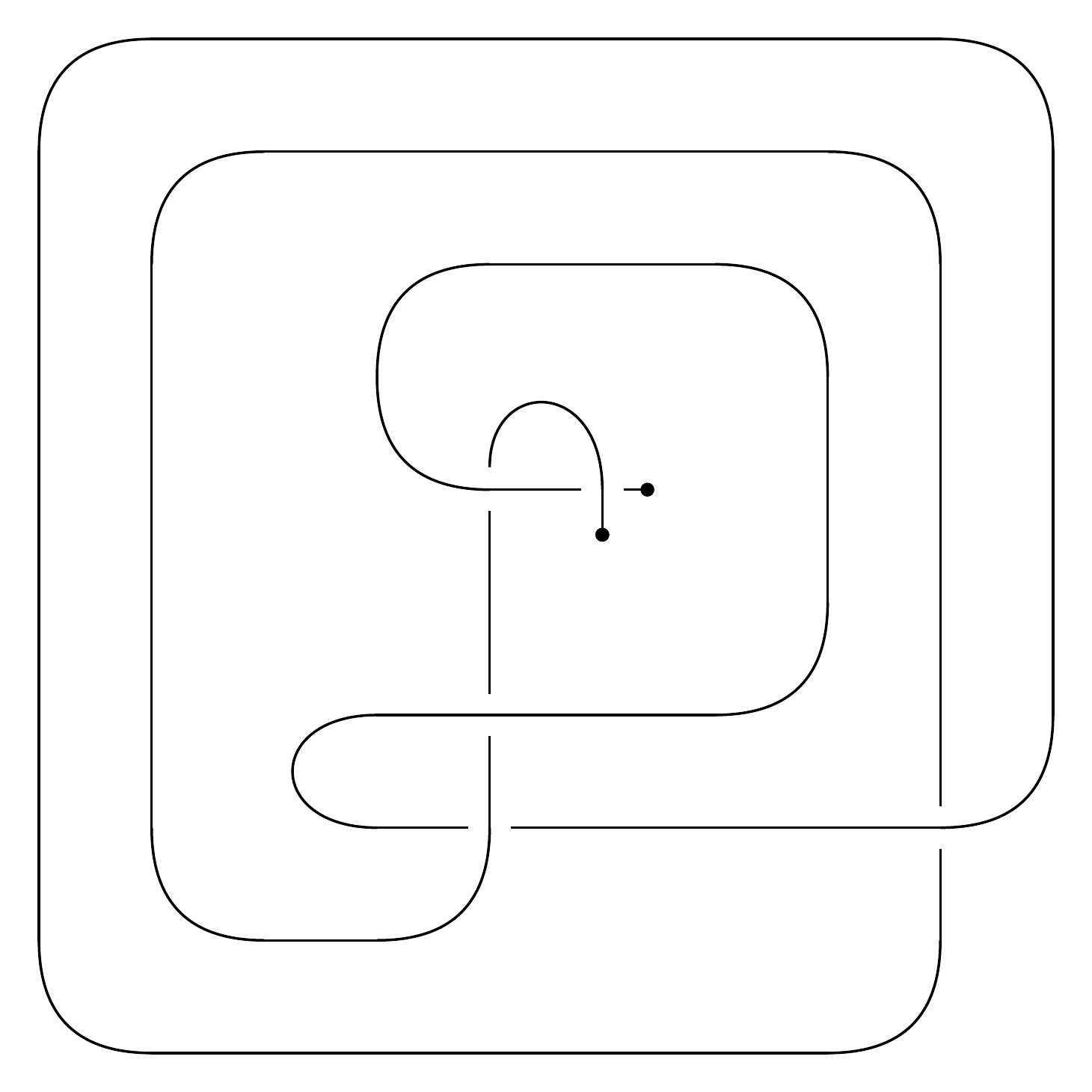}\\
\textcolor{black}{$5_{938}$}
\vspace{1cm}
\end{minipage}
\begin{minipage}[t]{.25\linewidth}
\centering
\includegraphics[width=0.9\textwidth,height=3.5cm,keepaspectratio]{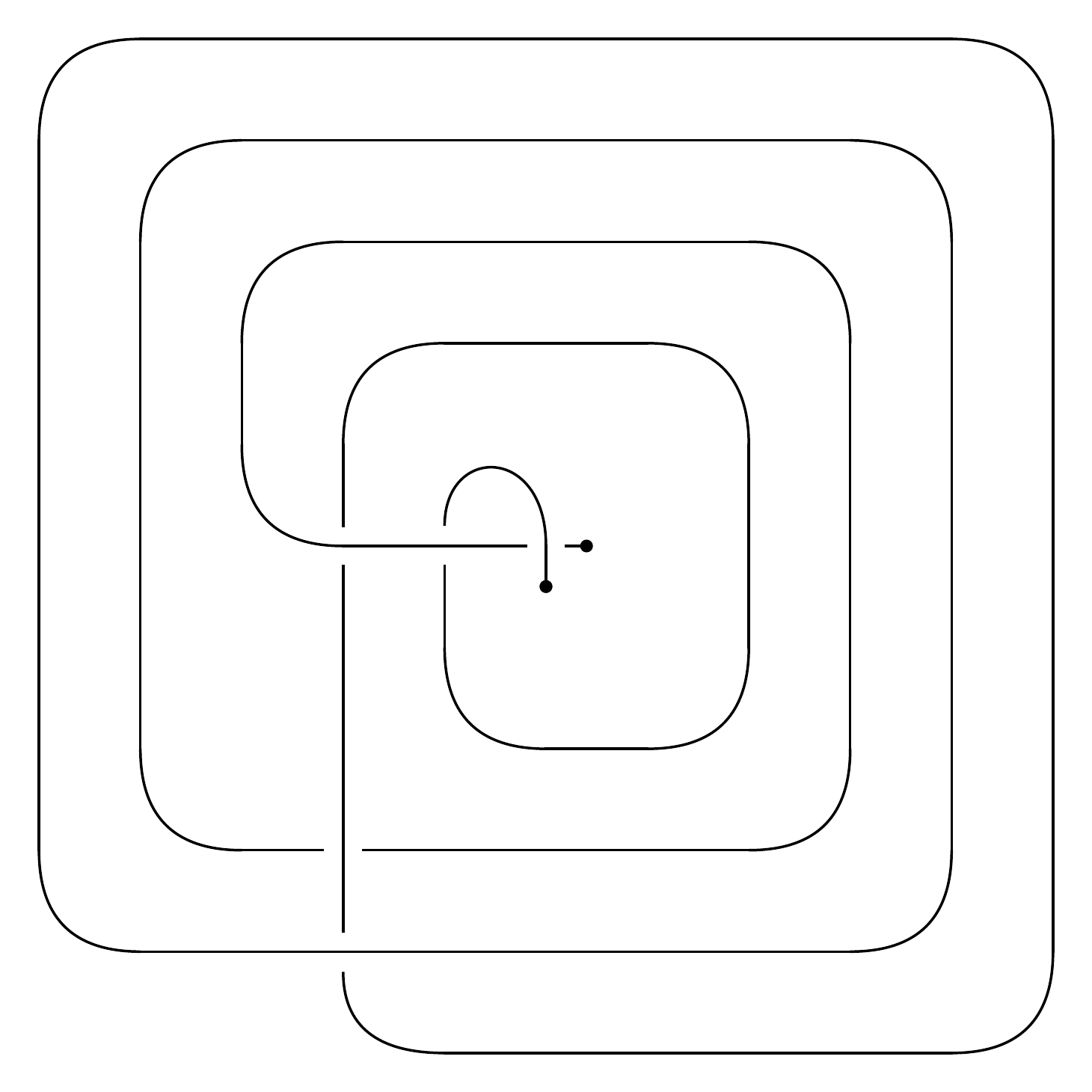}\\
\textcolor{black}{$5_{939}$}
\vspace{1cm}
\end{minipage}
\begin{minipage}[t]{.25\linewidth}
\centering
\includegraphics[width=0.9\textwidth,height=3.5cm,keepaspectratio]{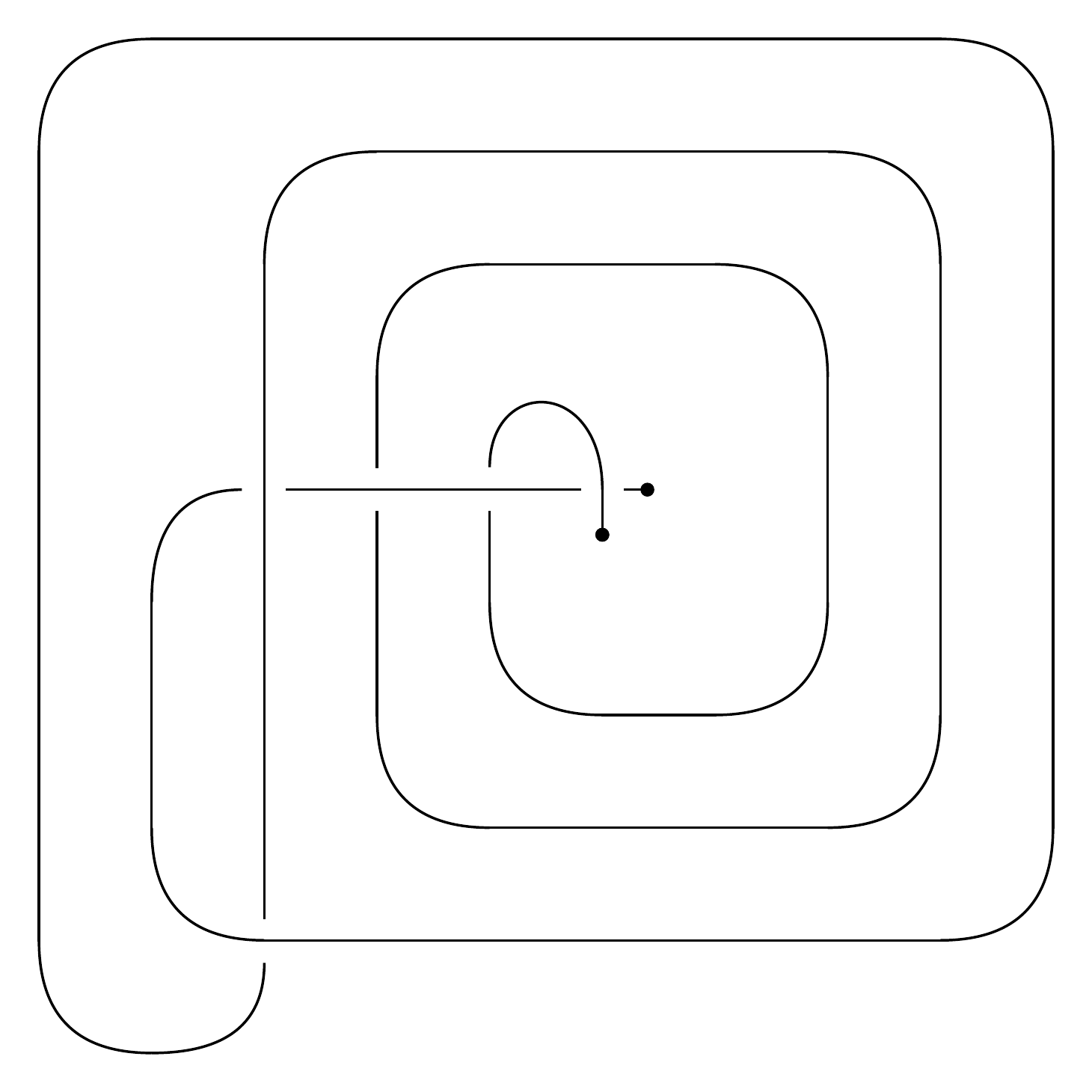}\\
\textcolor{black}{$5_{940}$}
\vspace{1cm}
\end{minipage}
\begin{minipage}[t]{.25\linewidth}
\centering
\includegraphics[width=0.9\textwidth,height=3.5cm,keepaspectratio]{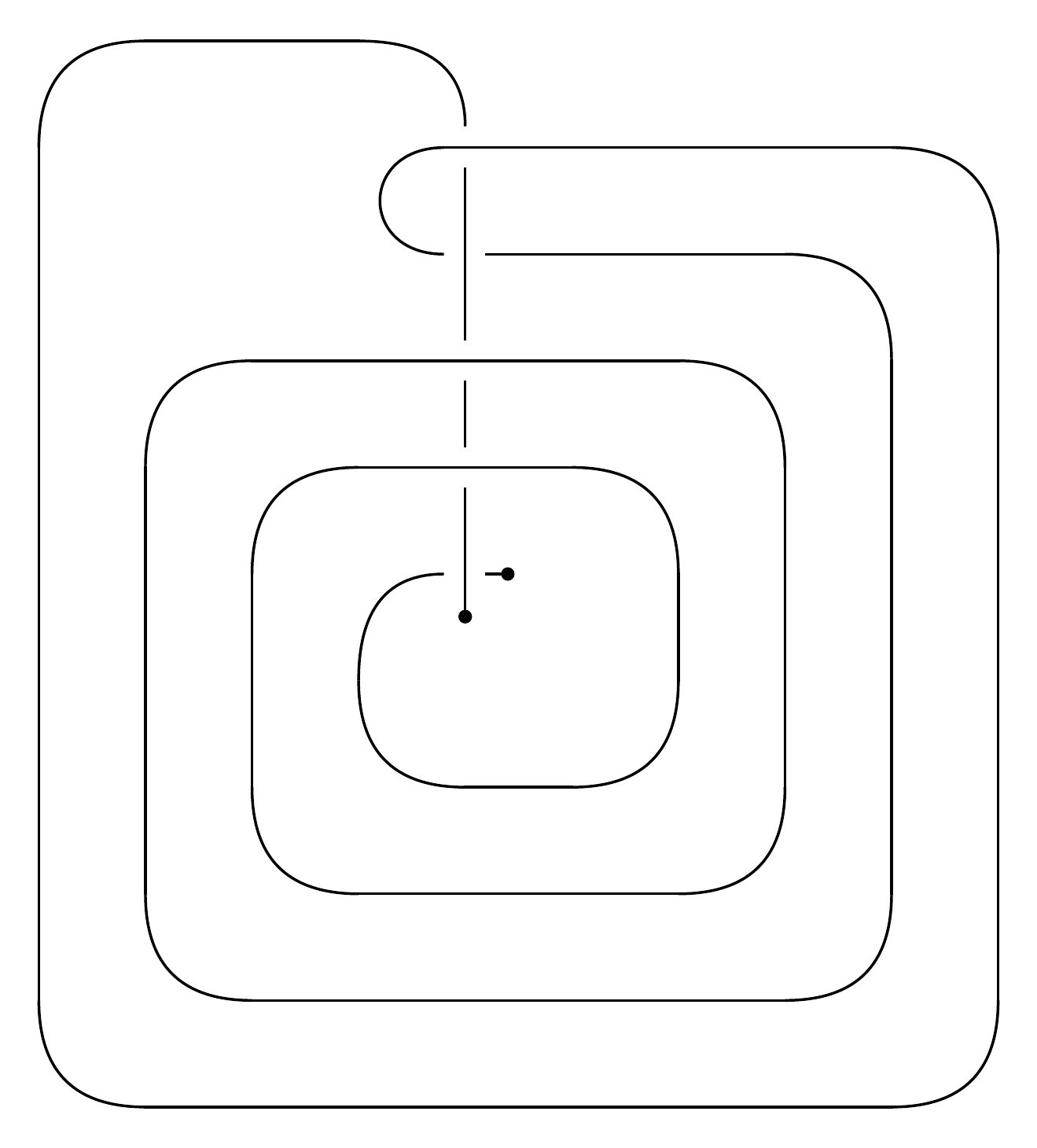}\\
\textcolor{black}{$5_{941}$}
\vspace{1cm}
\end{minipage}
\begin{minipage}[t]{.25\linewidth}
\centering
\includegraphics[width=0.9\textwidth,height=3.5cm,keepaspectratio]{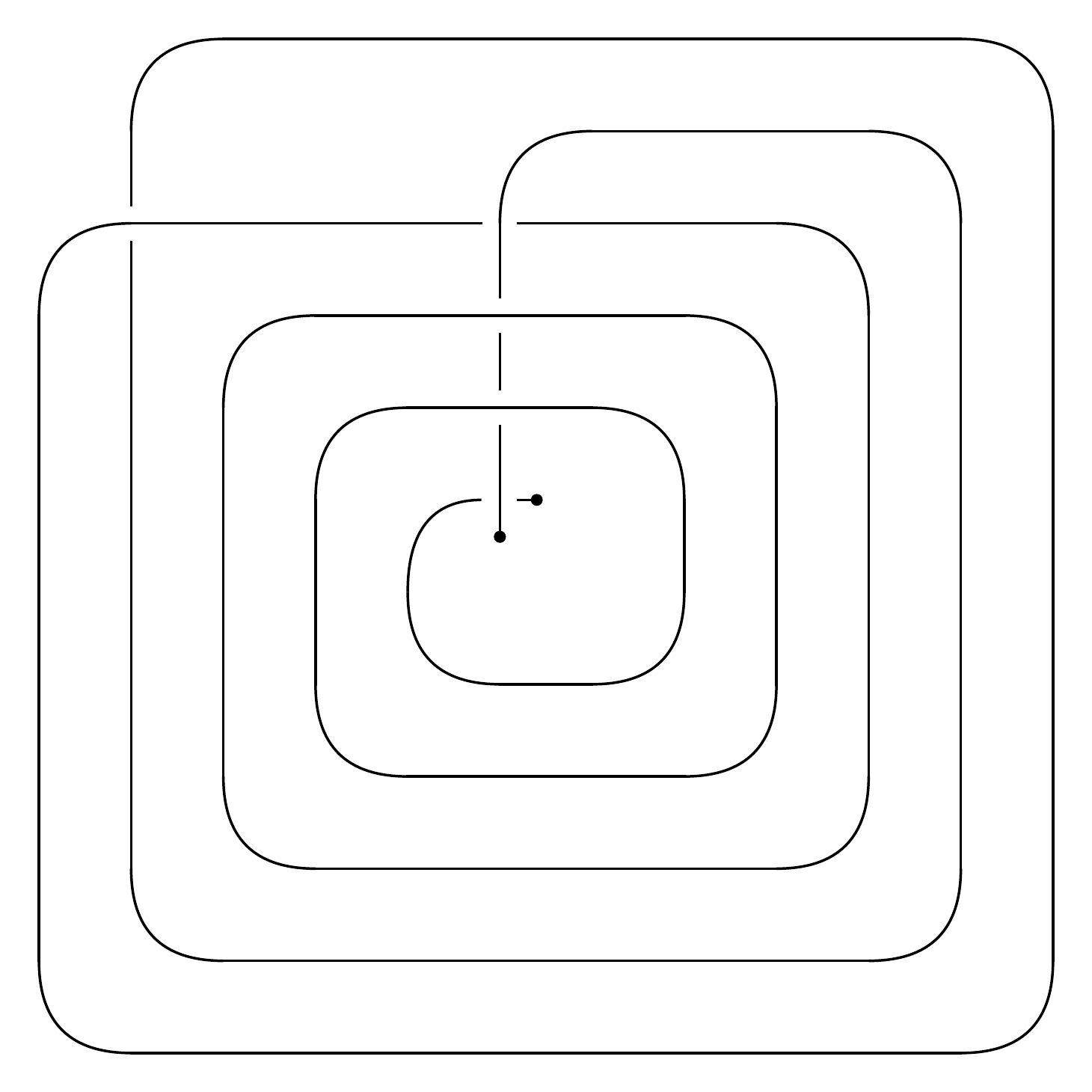}\\
\textcolor{black}{$5_{942}$}
\vspace{1cm}
\end{minipage}
\begin{minipage}[t]{.25\linewidth}
\centering
\includegraphics[width=0.9\textwidth,height=3.5cm,keepaspectratio]{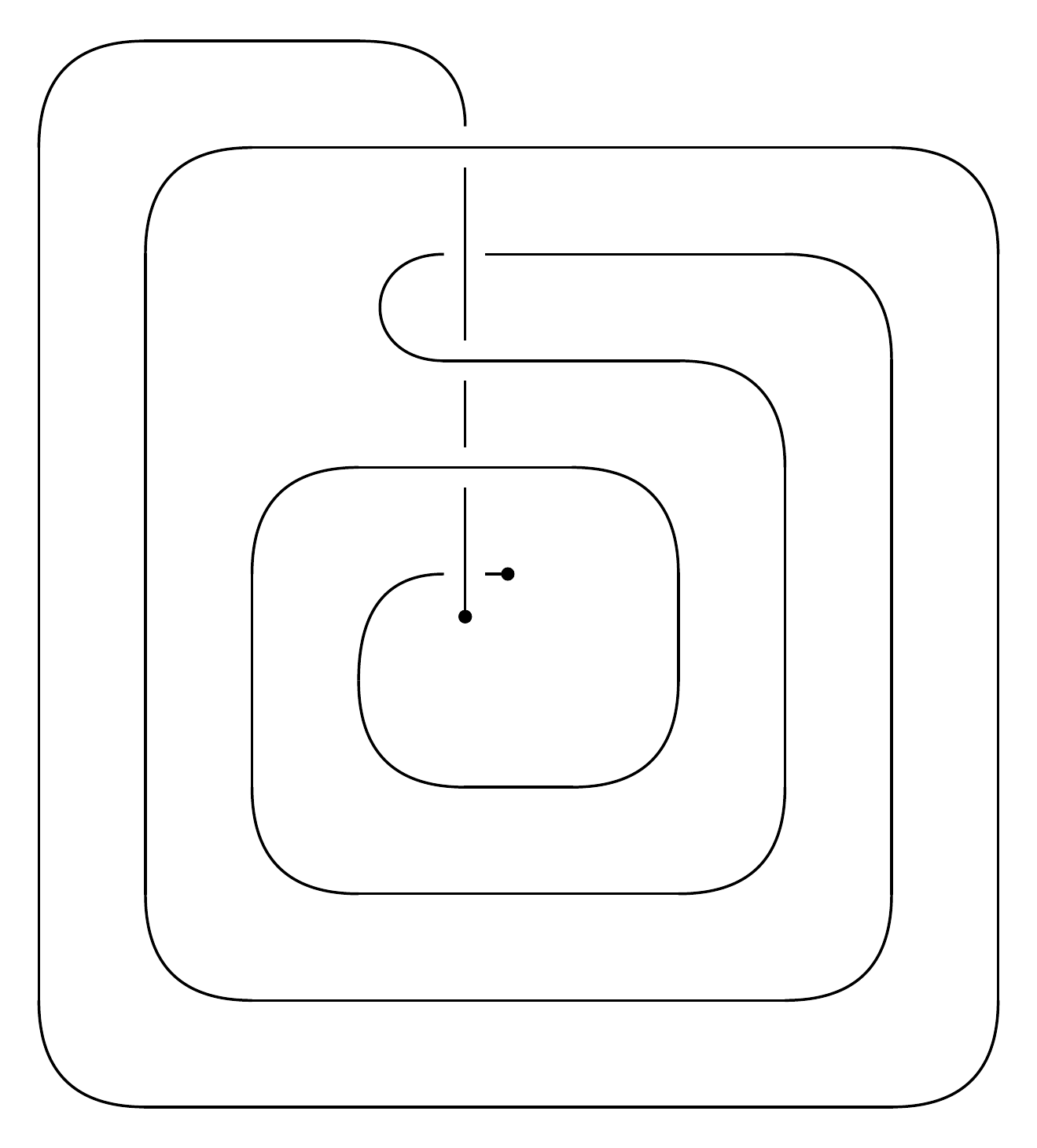}\\
\textcolor{black}{$5_{943}$}
\vspace{1cm}
\end{minipage}
\begin{minipage}[t]{.25\linewidth}
\centering
\includegraphics[width=0.9\textwidth,height=3.5cm,keepaspectratio]{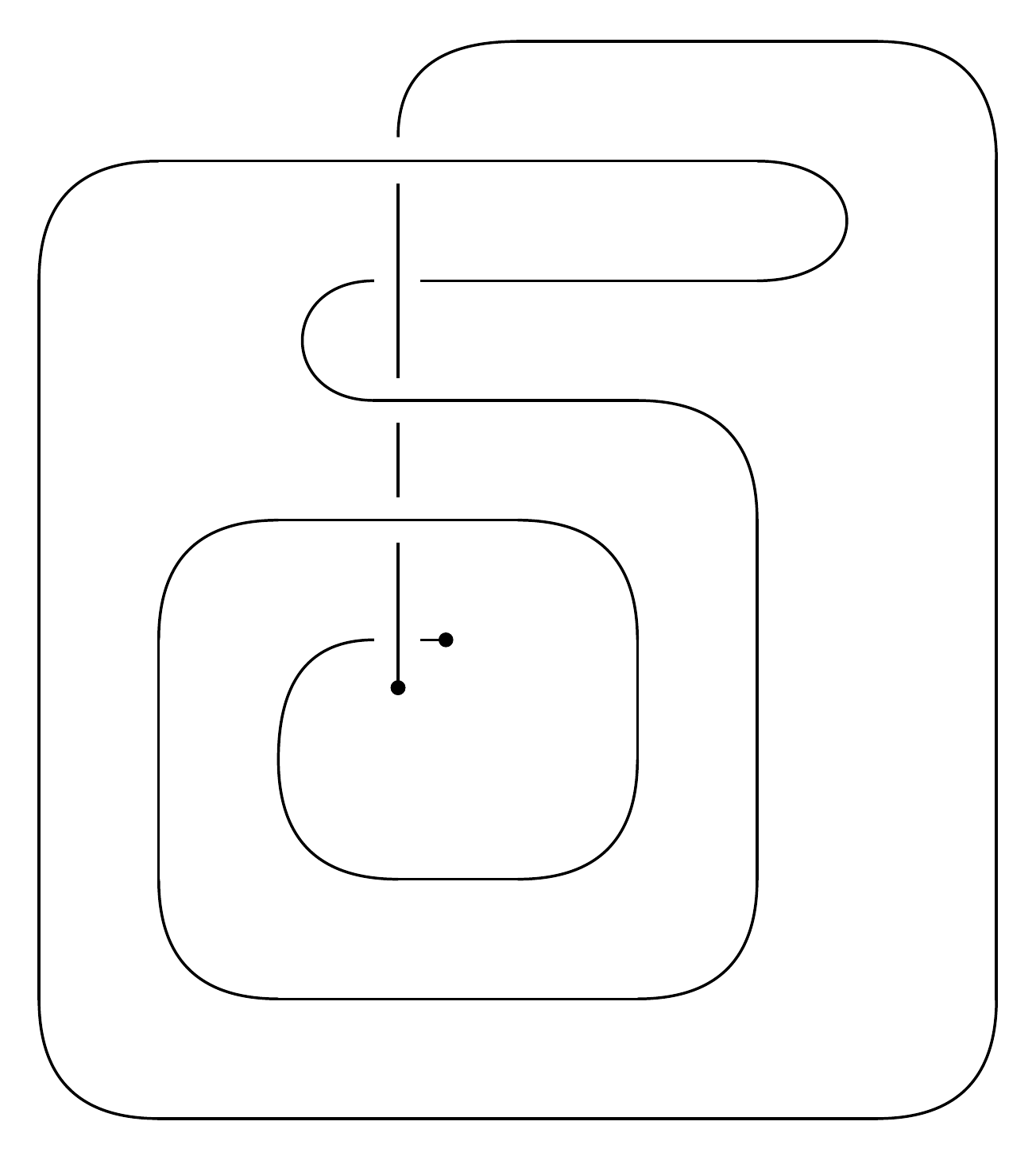}\\
\textcolor{black}{$5_{944}$}
\vspace{1cm}
\end{minipage}
\begin{minipage}[t]{.25\linewidth}
\centering
\includegraphics[width=0.9\textwidth,height=3.5cm,keepaspectratio]{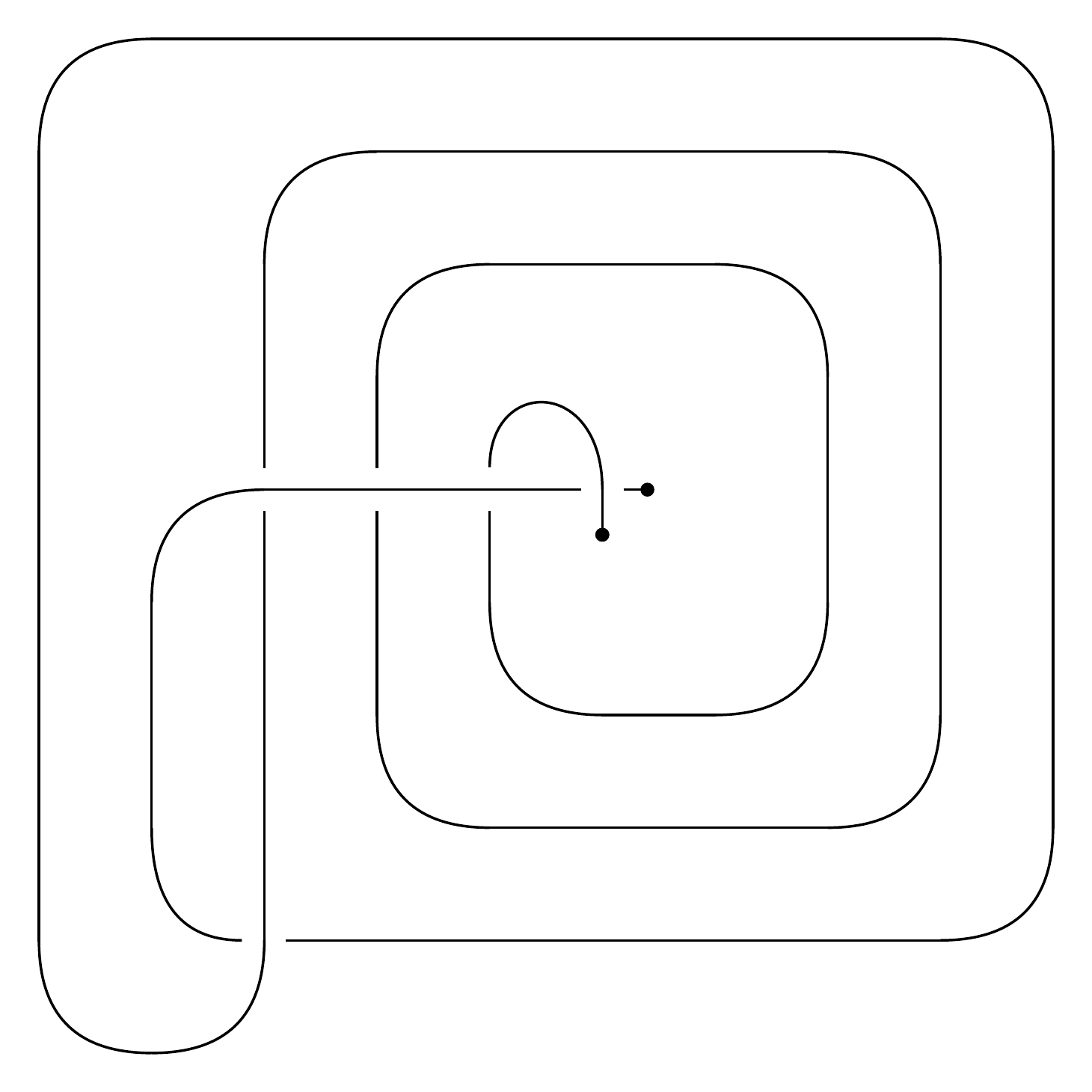}\\
\textcolor{black}{$5_{945}$}
\vspace{1cm}
\end{minipage}
\begin{minipage}[t]{.25\linewidth}
\centering
\includegraphics[width=0.9\textwidth,height=3.5cm,keepaspectratio]{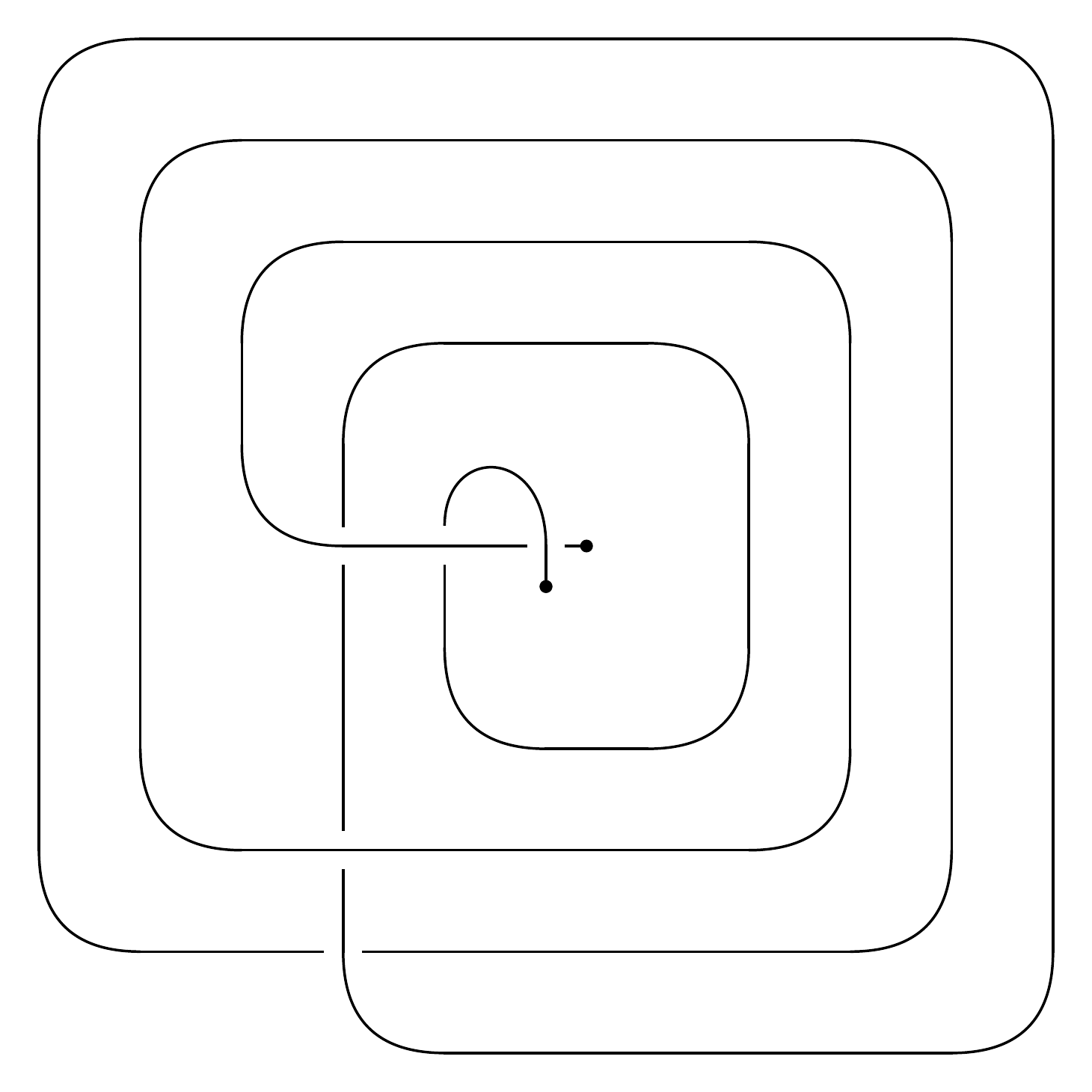}\\
\textcolor{black}{$5_{946}$}
\vspace{1cm}
\end{minipage}
\begin{minipage}[t]{.25\linewidth}
\centering
\includegraphics[width=0.9\textwidth,height=3.5cm,keepaspectratio]{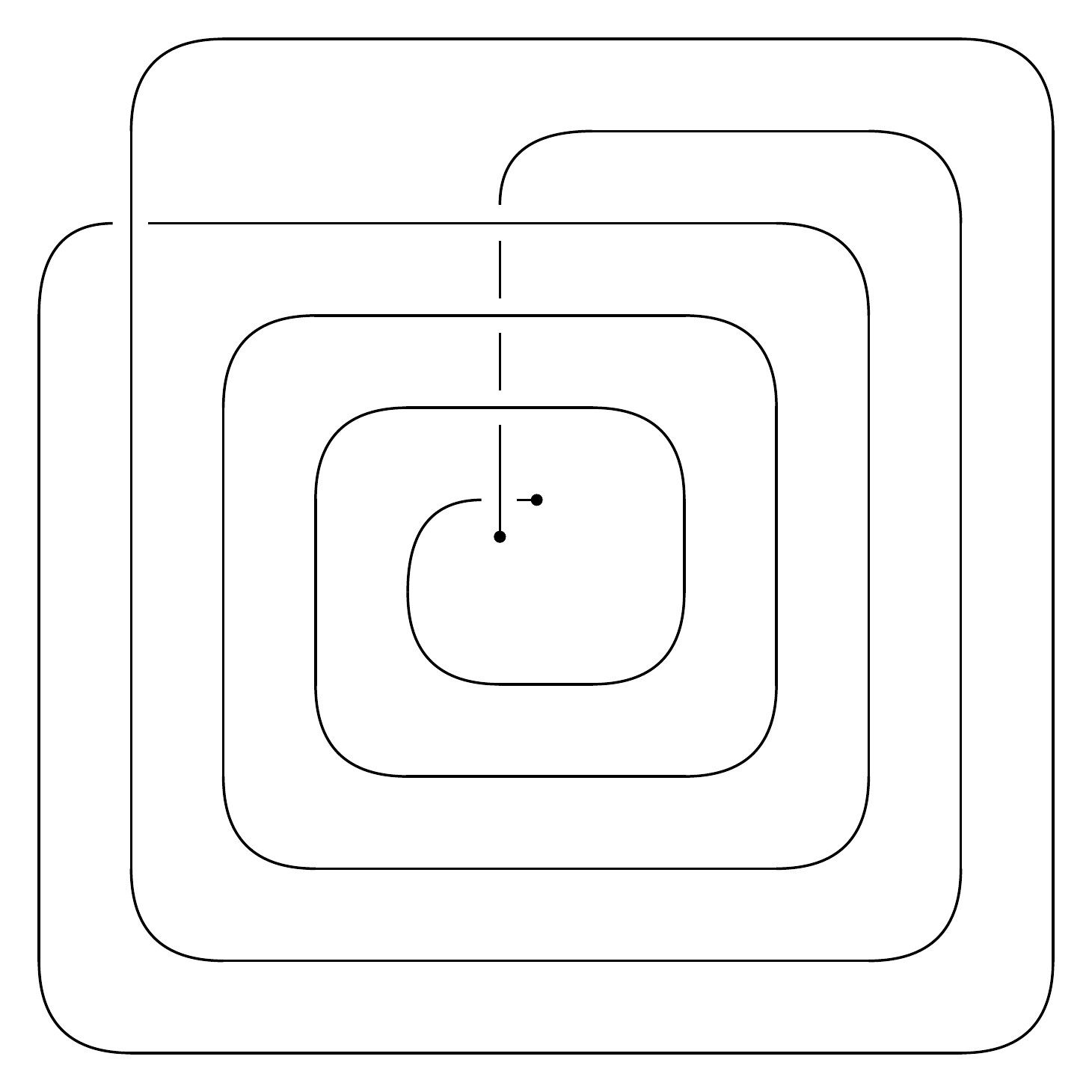}\\
\textcolor{black}{$5_{947}$}
\vspace{1cm}
\end{minipage}
\begin{minipage}[t]{.25\linewidth}
\centering
\includegraphics[width=0.9\textwidth,height=3.5cm,keepaspectratio]{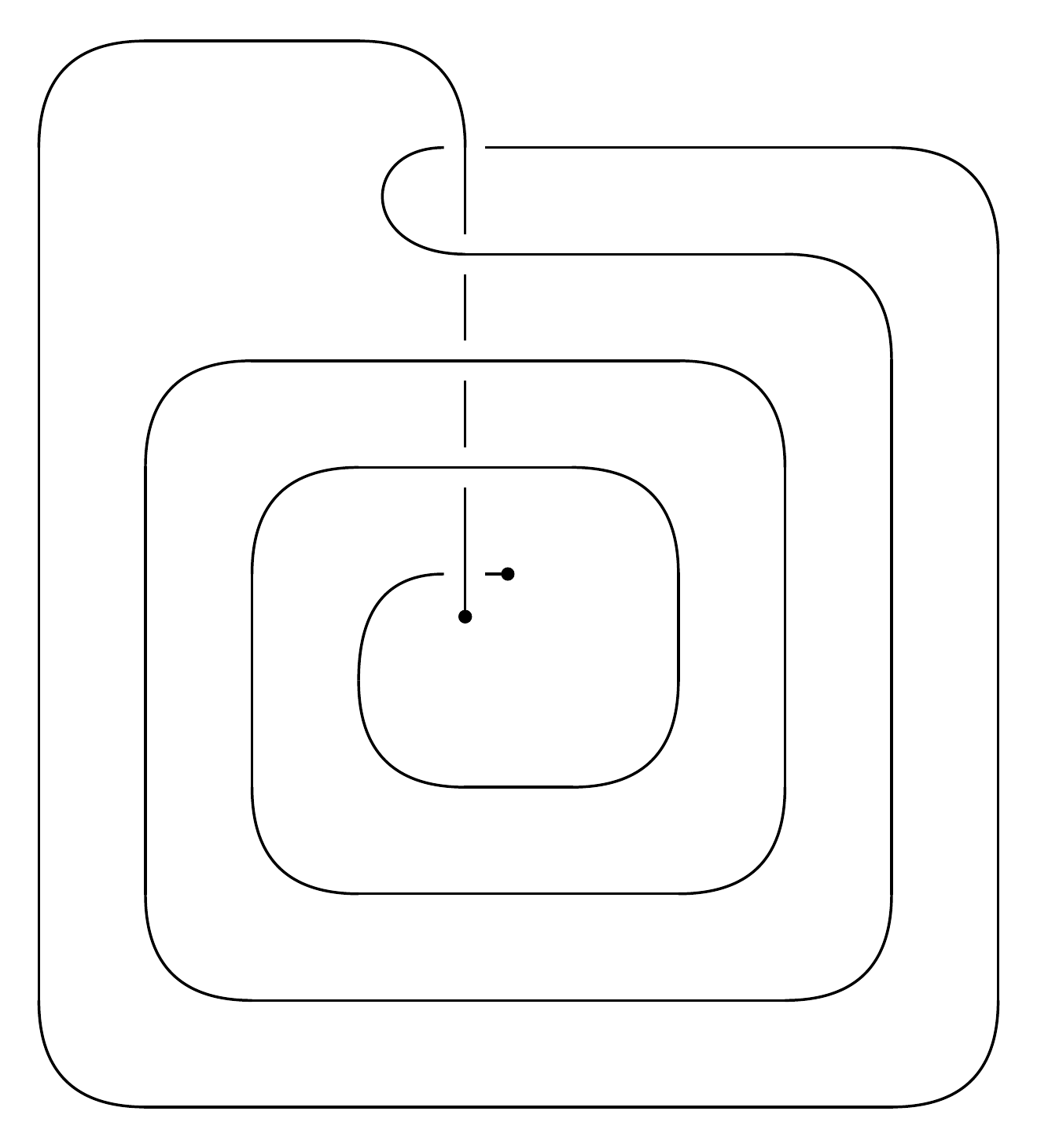}\\
\textcolor{black}{$5_{948}$}
\vspace{1cm}
\end{minipage}
\begin{minipage}[t]{.25\linewidth}
\centering
\includegraphics[width=0.9\textwidth,height=3.5cm,keepaspectratio]{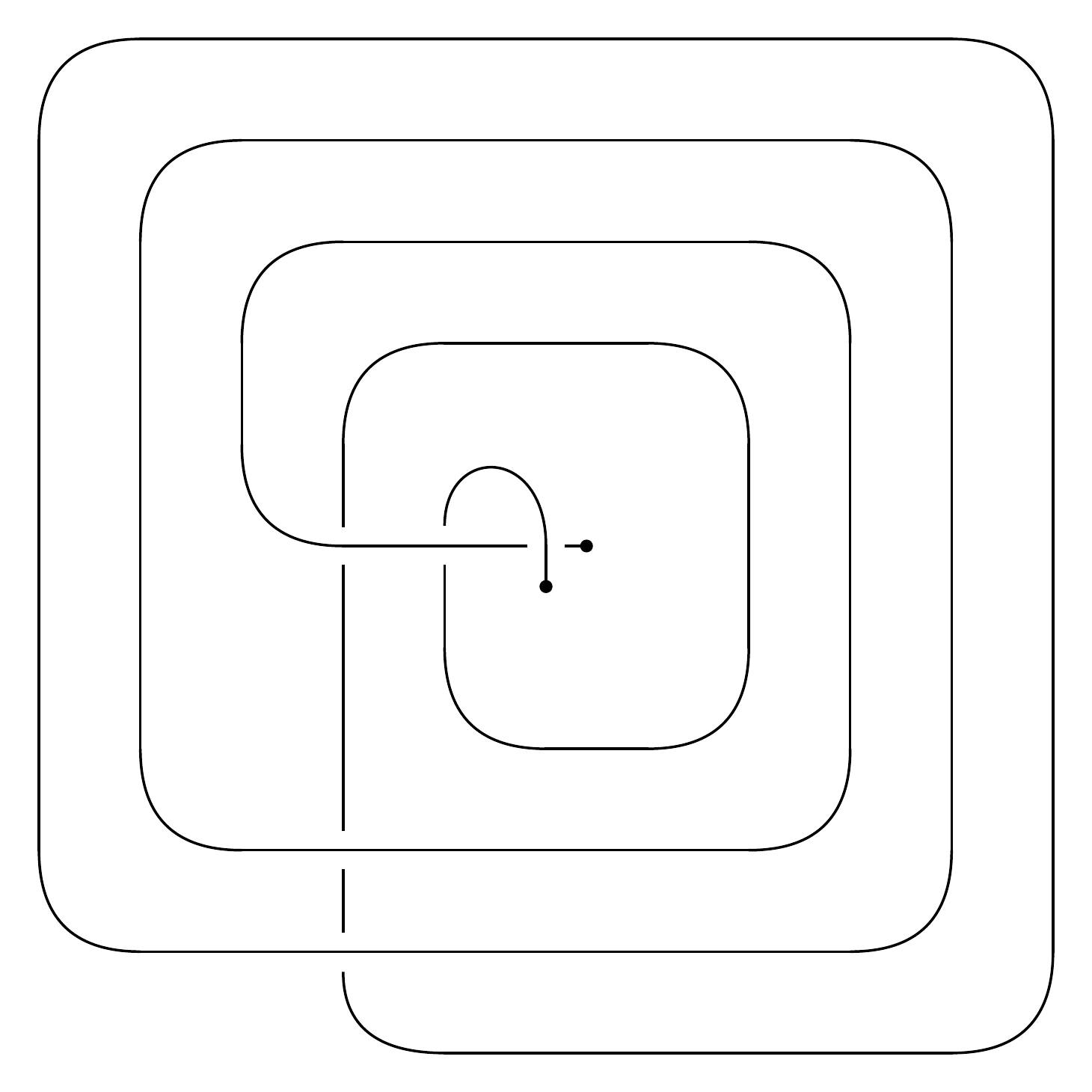}\\
\textcolor{black}{$5_{949}$}
\vspace{1cm}
\end{minipage}
\begin{minipage}[t]{.25\linewidth}
\centering
\includegraphics[width=0.9\textwidth,height=3.5cm,keepaspectratio]{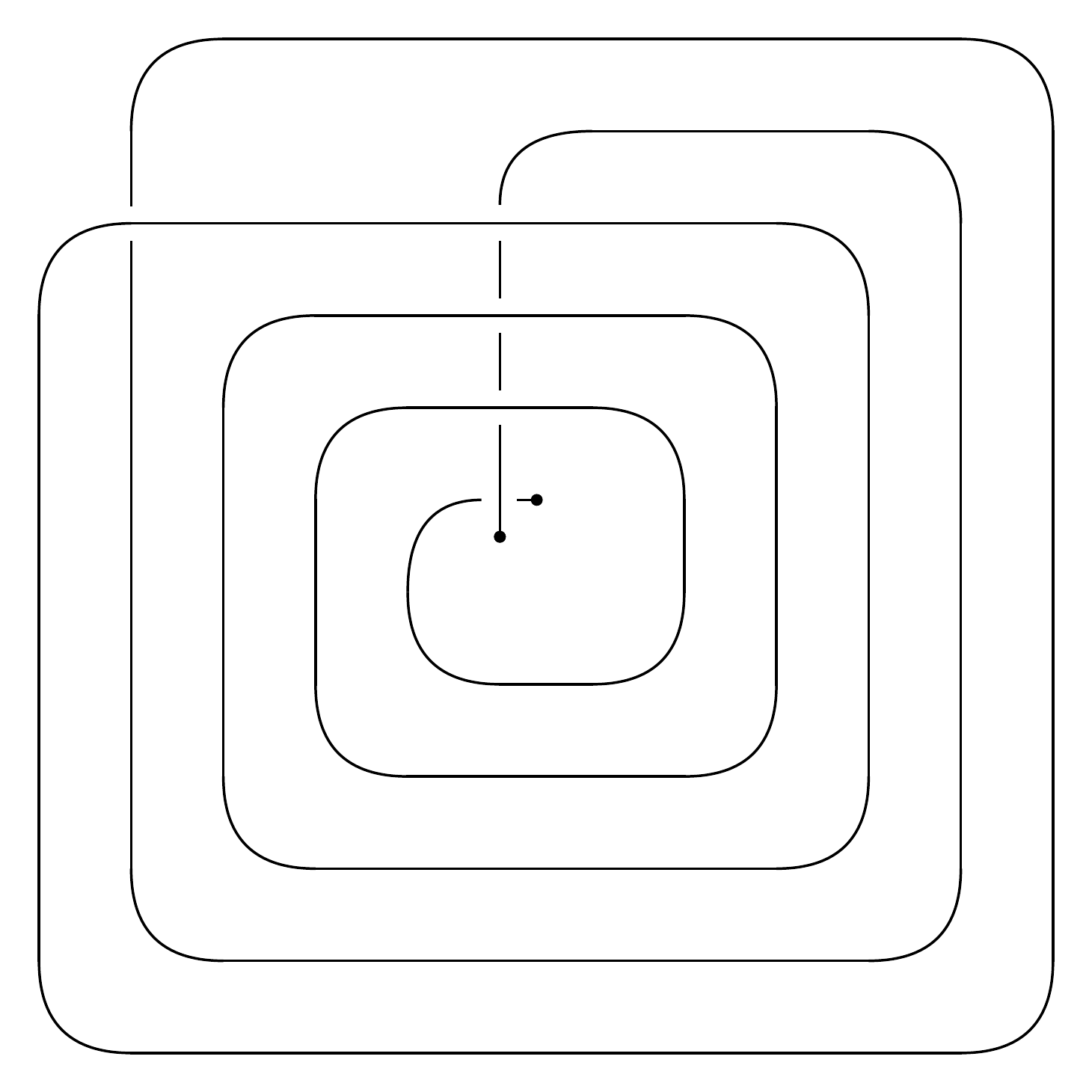}\\
\textcolor{black}{$5_{950}$}
\vspace{1cm}
\end{minipage}

%% file: input_sphere.tex
\begin{minipage}[t]{.25\linewidth}
\centering
\includegraphics[width=0.9\textwidth,height=3.5cm,keepaspectratio]{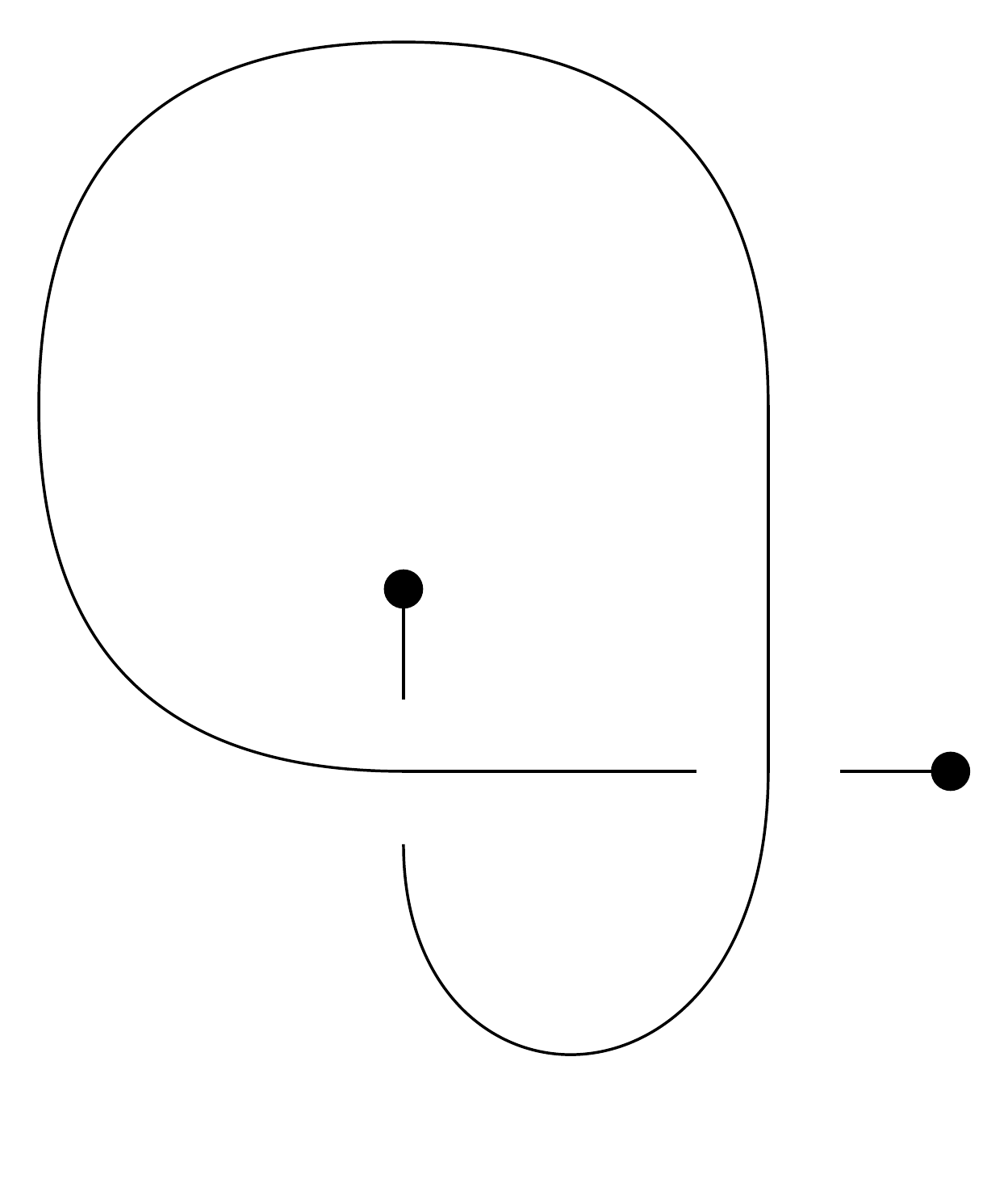}\\
\textcolor{black}{$2_{1}$}
\vspace{1cm}
\end{minipage}
\begin{minipage}[t]{.25\linewidth}
\centering
\includegraphics[width=0.9\textwidth,height=3.5cm,keepaspectratio]{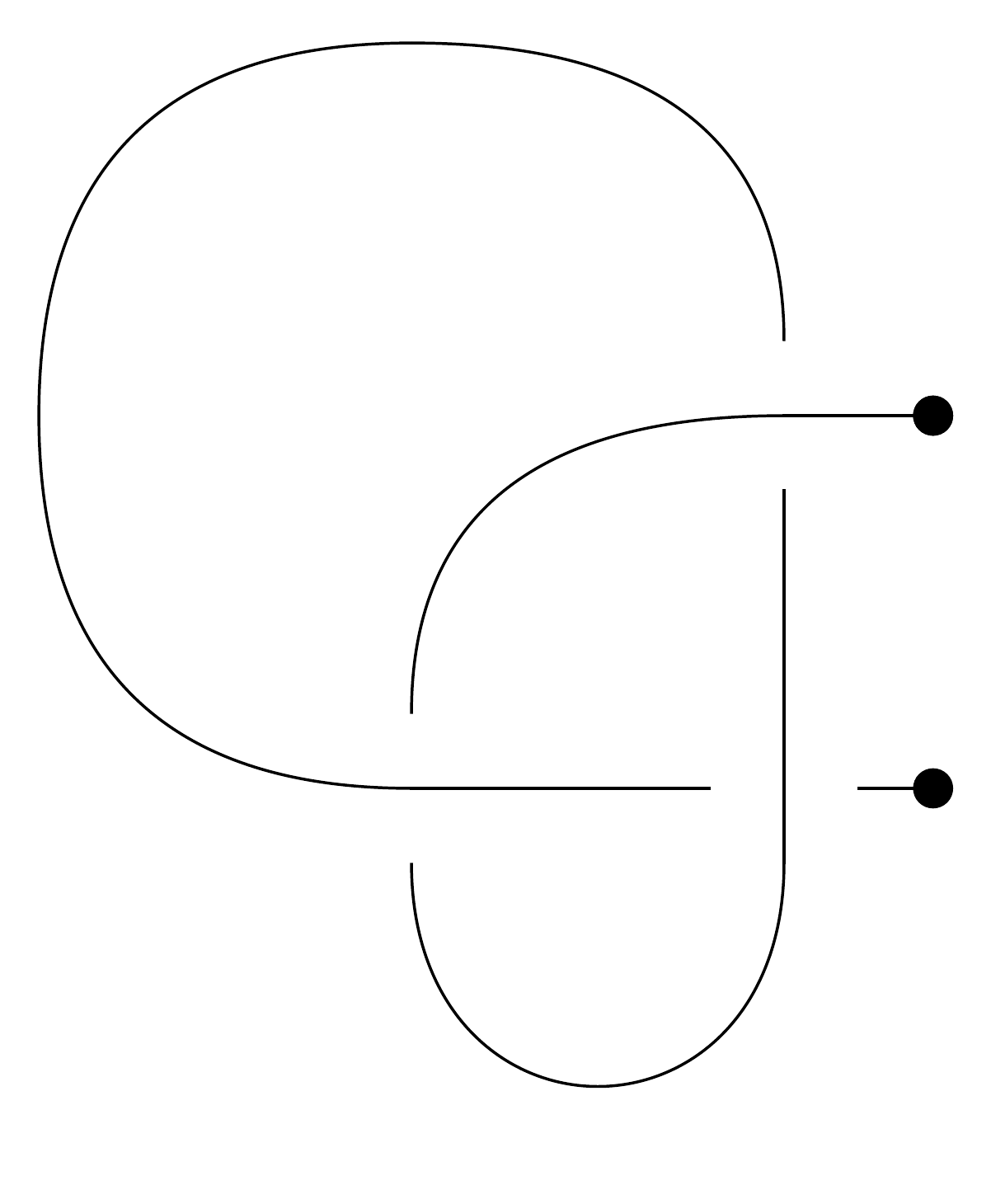}\\
\textcolor{red}{$3_{1}$}
\vspace{1cm}
\end{minipage}
\begin{minipage}[t]{.25\linewidth}
\centering
\includegraphics[width=0.9\textwidth,height=3.5cm,keepaspectratio]{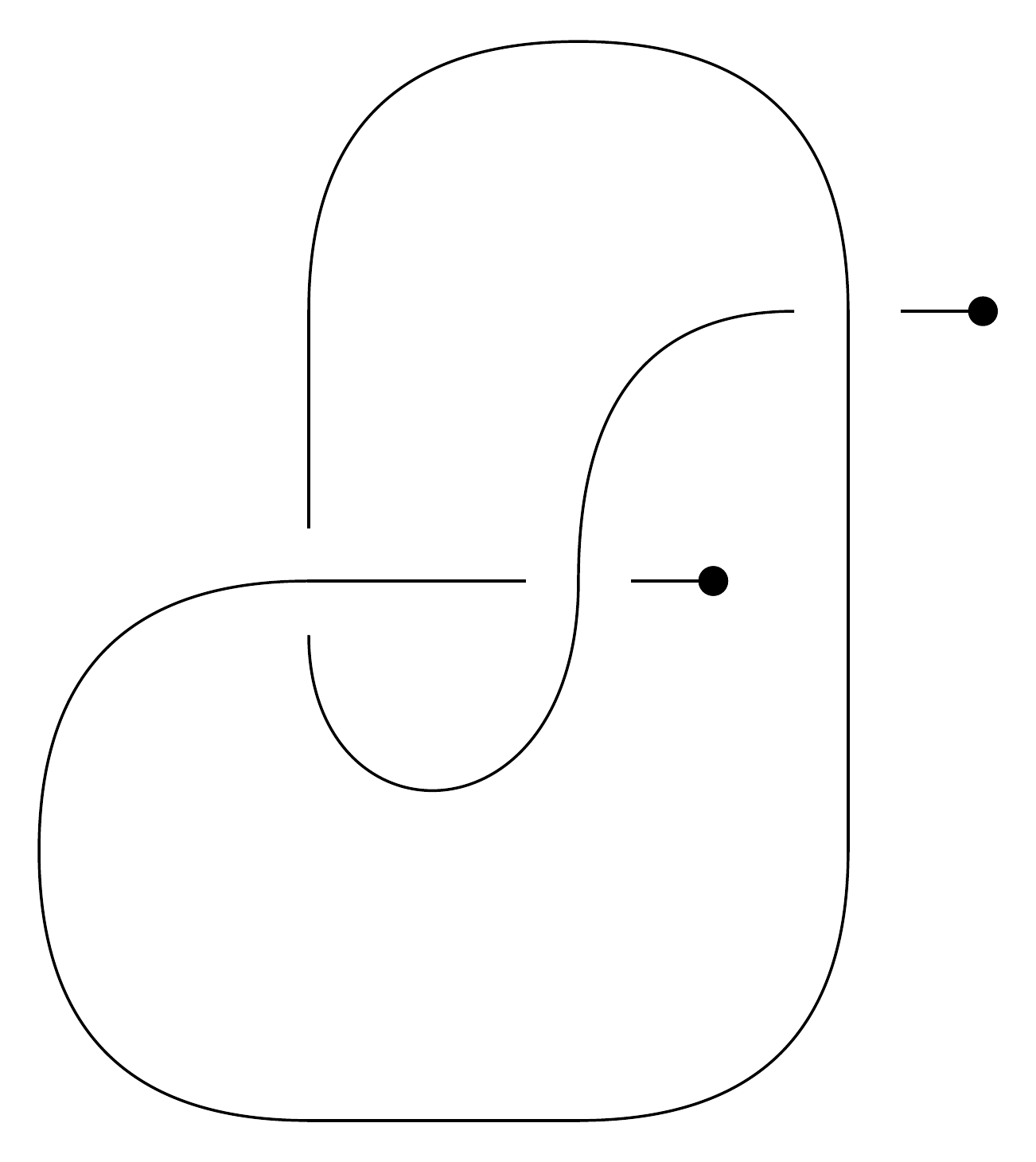}\\
\textcolor{black}{$3_{2}$}
\vspace{1cm}
\end{minipage}
\begin{minipage}[t]{.25\linewidth}
\centering
\includegraphics[width=0.9\textwidth,height=3.5cm,keepaspectratio]{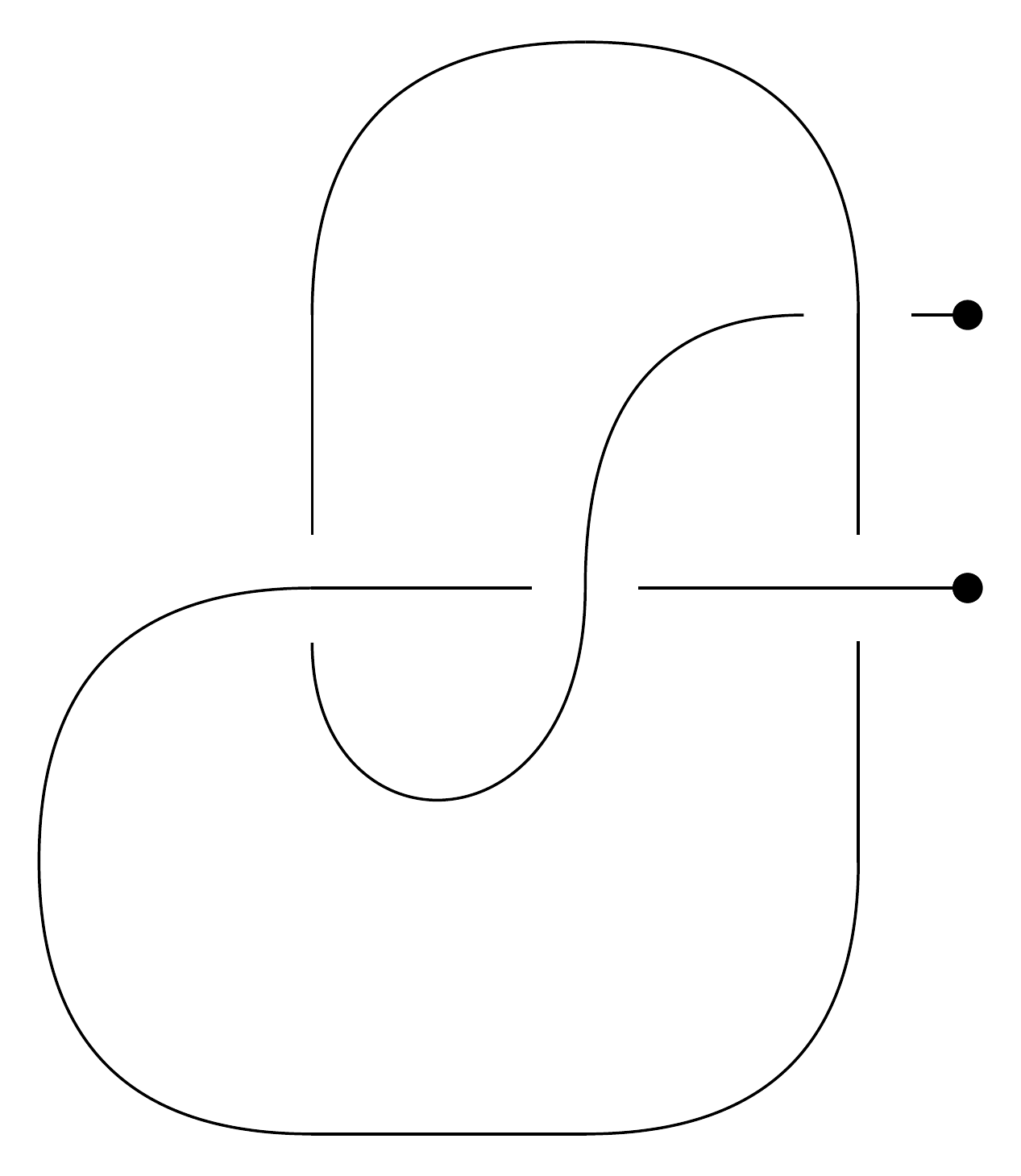}\\
\textcolor{red}{$4_{1}$}
\vspace{1cm}
\end{minipage}
\begin{minipage}[t]{.25\linewidth}
\centering
\includegraphics[width=0.9\textwidth,height=3.5cm,keepaspectratio]{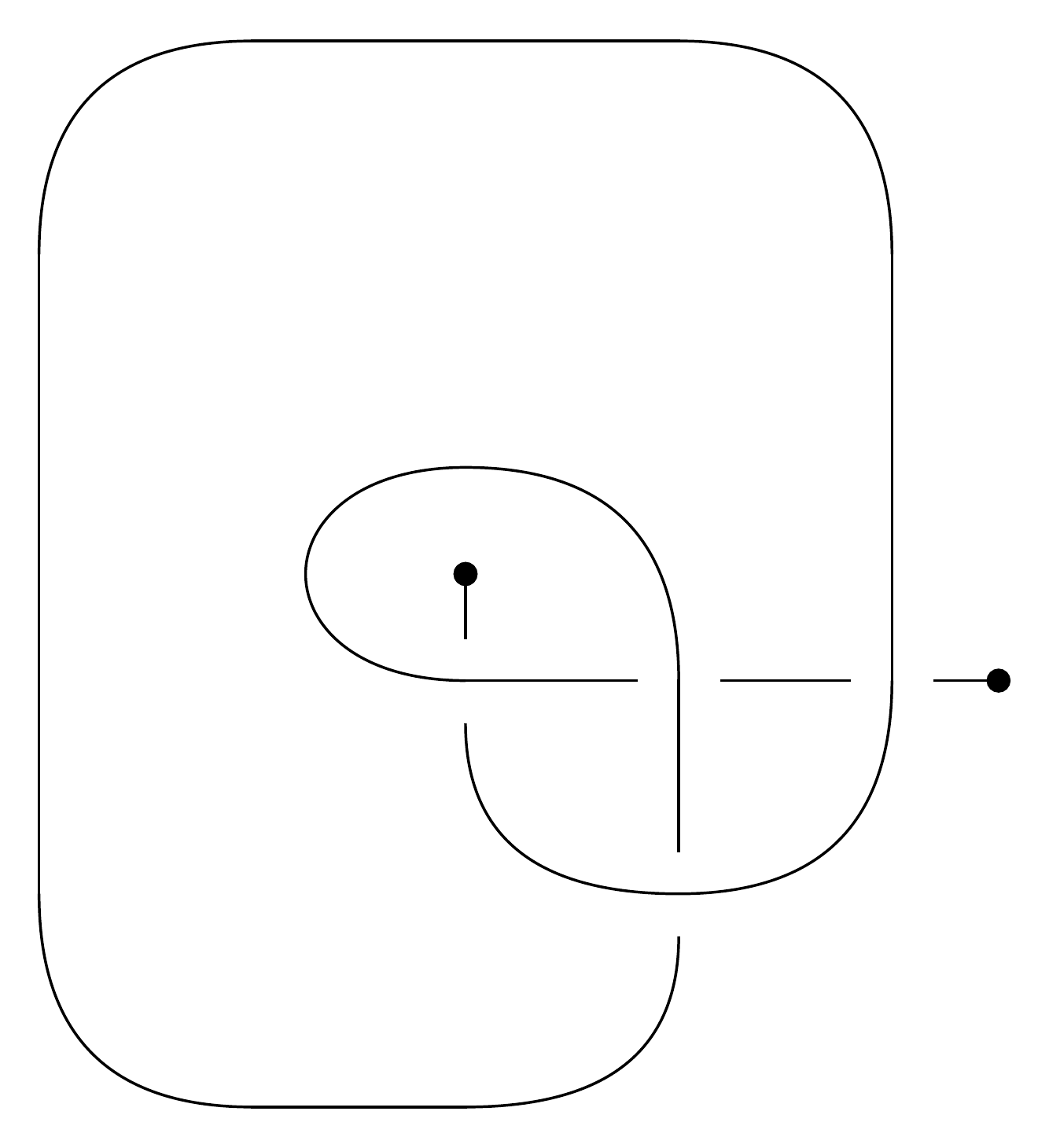}\\
\textcolor{black}{$4_{2}$}
\vspace{1cm}
\end{minipage}
\begin{minipage}[t]{.25\linewidth}
\centering
\includegraphics[width=0.9\textwidth,height=3.5cm,keepaspectratio]{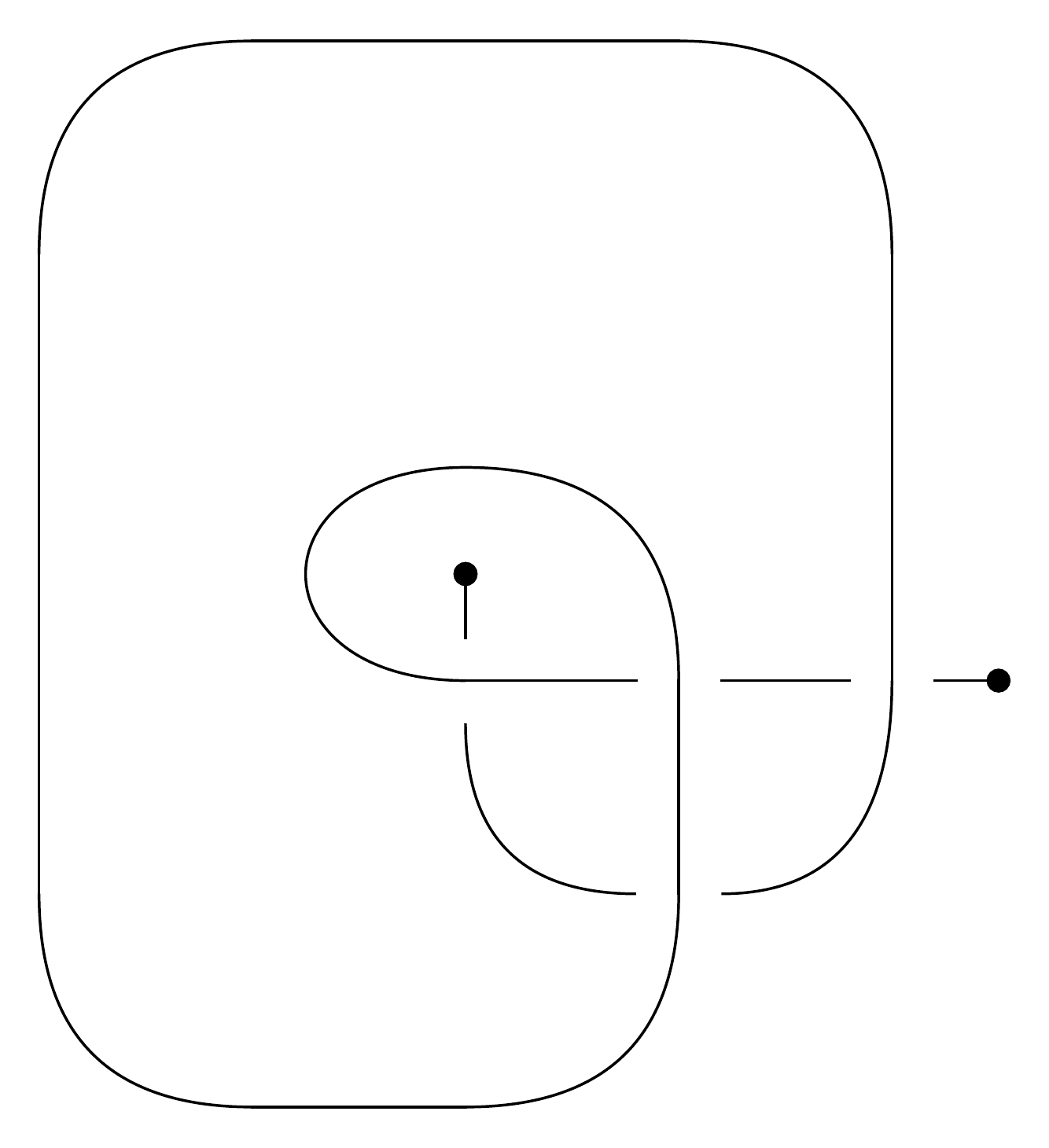}\\
\textcolor{black}{$4_{3}$}
\vspace{1cm}
\end{minipage}
\begin{minipage}[t]{.25\linewidth}
\centering
\includegraphics[width=0.9\textwidth,height=3.5cm,keepaspectratio]{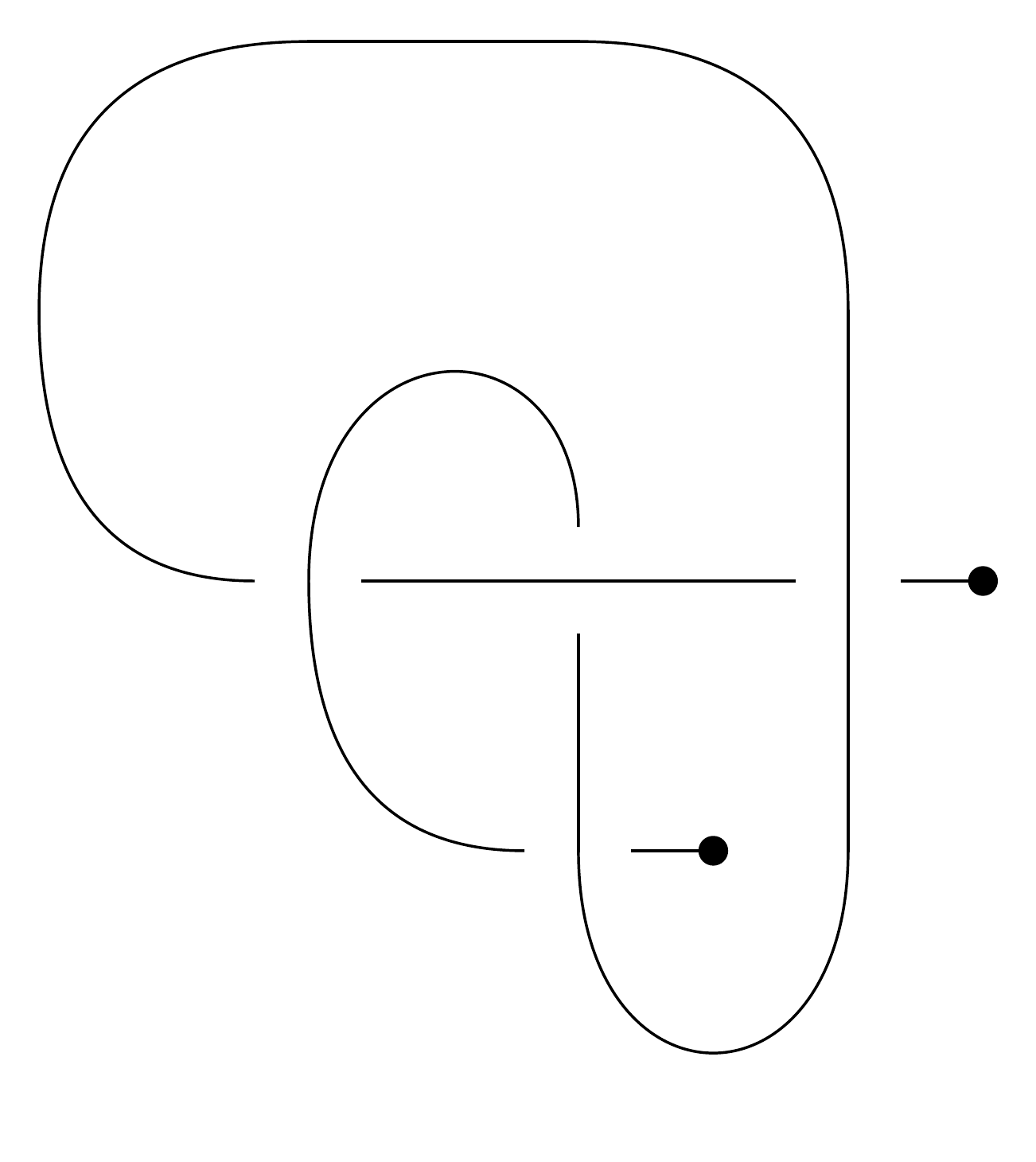}\\
\textcolor{black}{$4_{4}$}
\vspace{1cm}
\end{minipage}
\begin{minipage}[t]{.25\linewidth}
\centering
\includegraphics[width=0.9\textwidth,height=3.5cm,keepaspectratio]{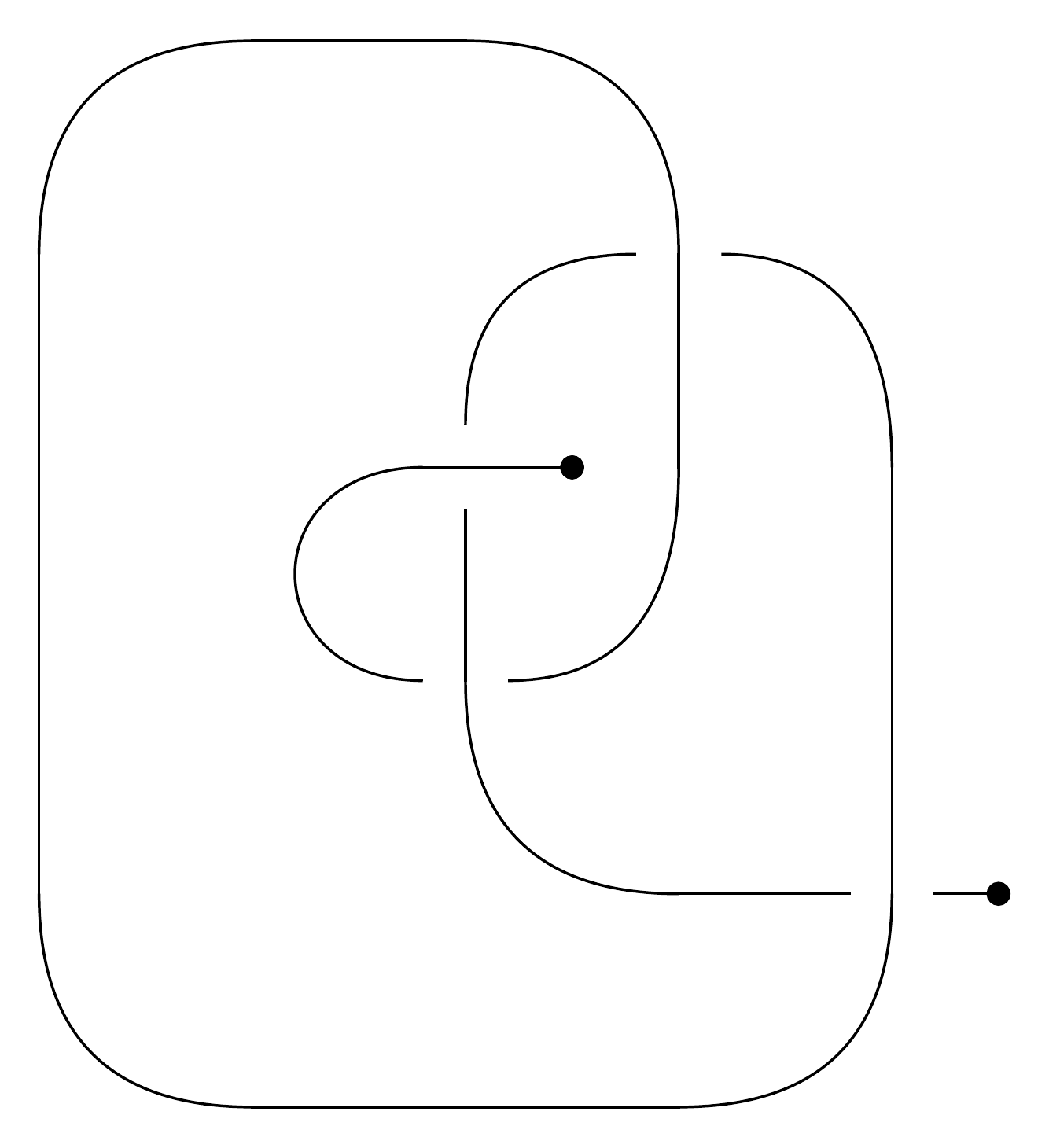}\\
\textcolor{black}{$4_{5}$}
\vspace{1cm}
\end{minipage}
\begin{minipage}[t]{.25\linewidth}
\centering
\includegraphics[width=0.9\textwidth,height=3.5cm,keepaspectratio]{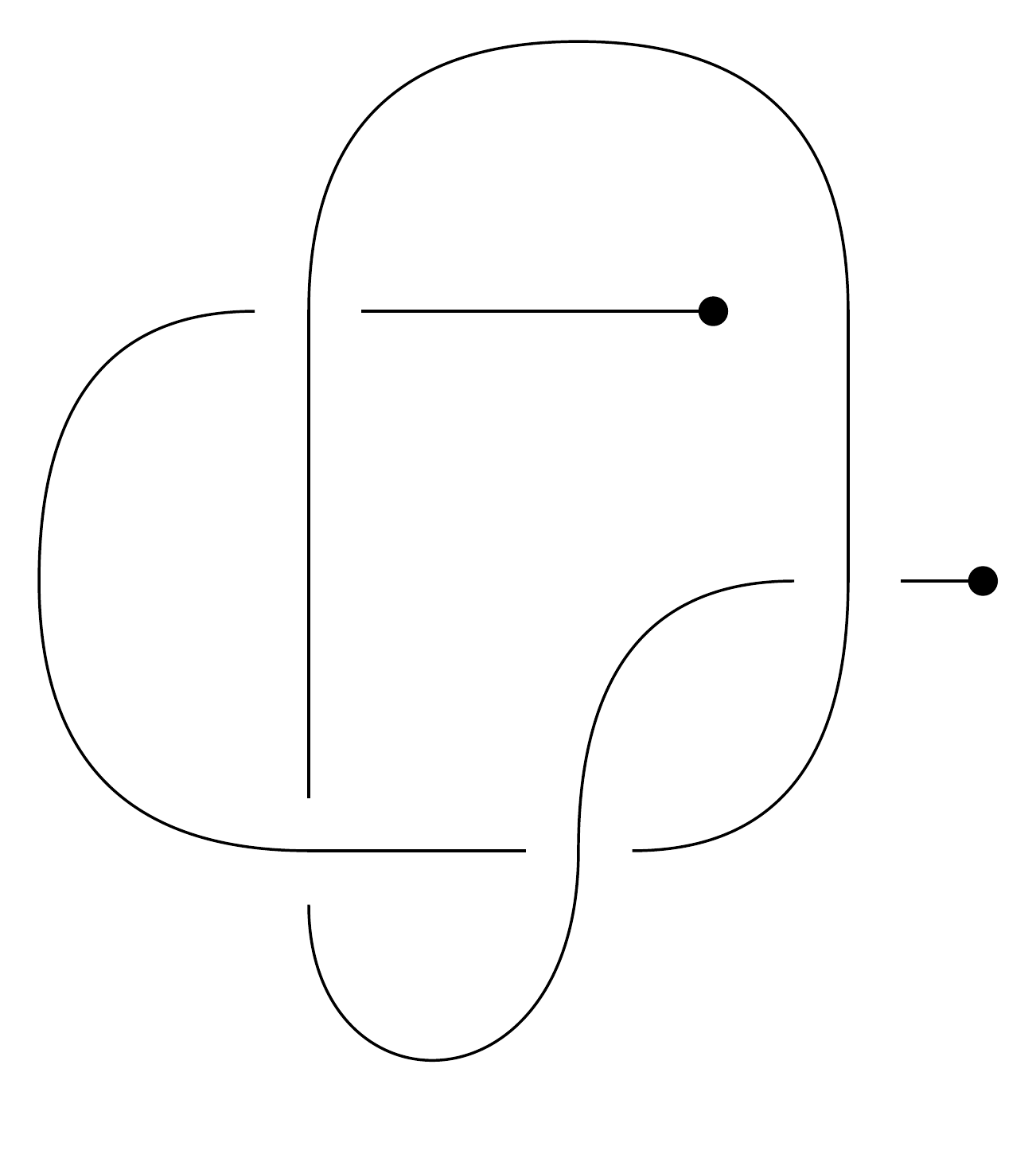}\\
\textcolor{black}{$4_{6}$}
\vspace{1cm}
\end{minipage}
\begin{minipage}[t]{.25\linewidth}
\centering
\includegraphics[width=0.9\textwidth,height=3.5cm,keepaspectratio]{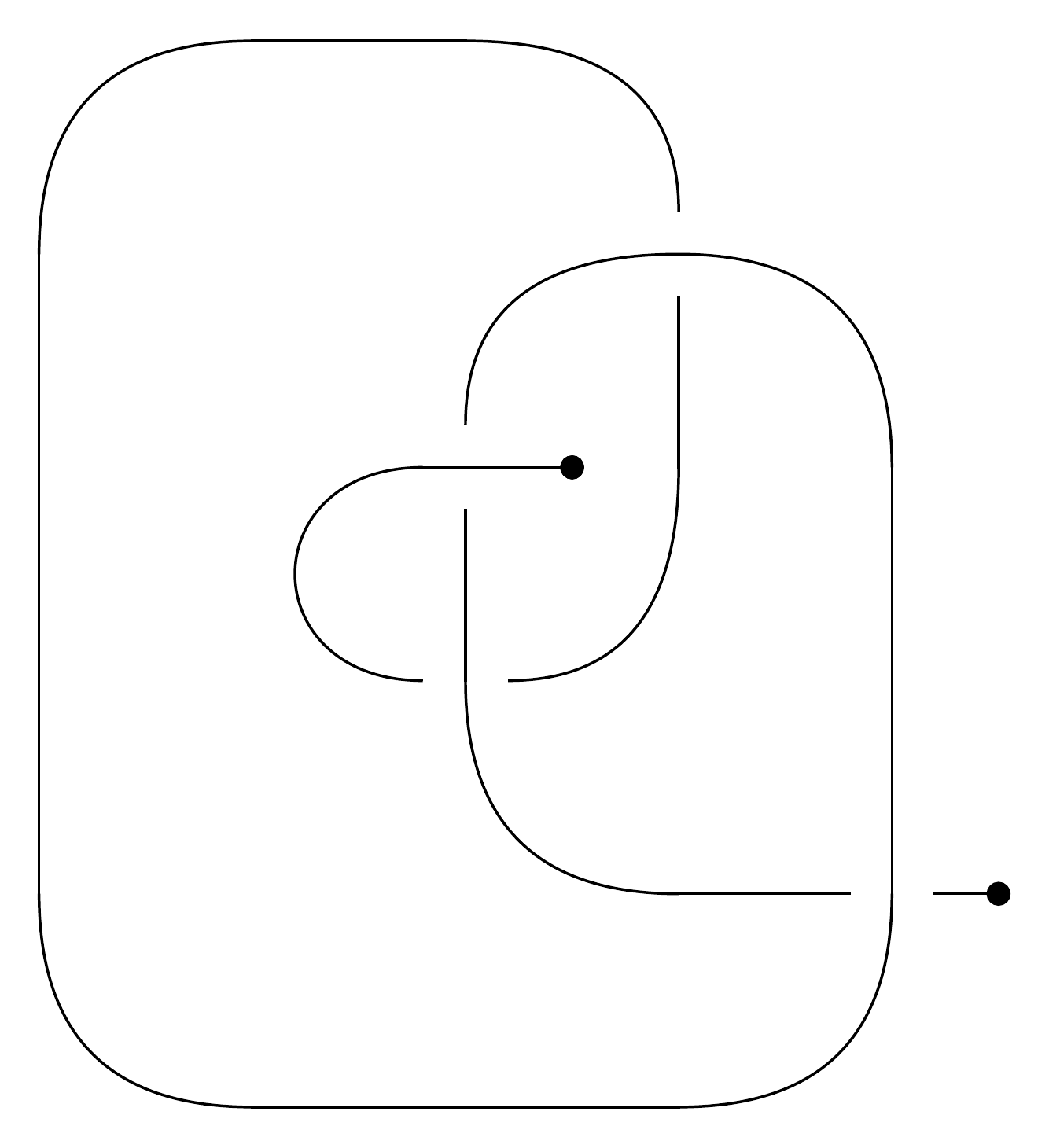}\\
\textcolor{black}{$4_{7}$}
\vspace{1cm}
\end{minipage}
\begin{minipage}[t]{.25\linewidth}
\centering
\includegraphics[width=0.9\textwidth,height=3.5cm,keepaspectratio]{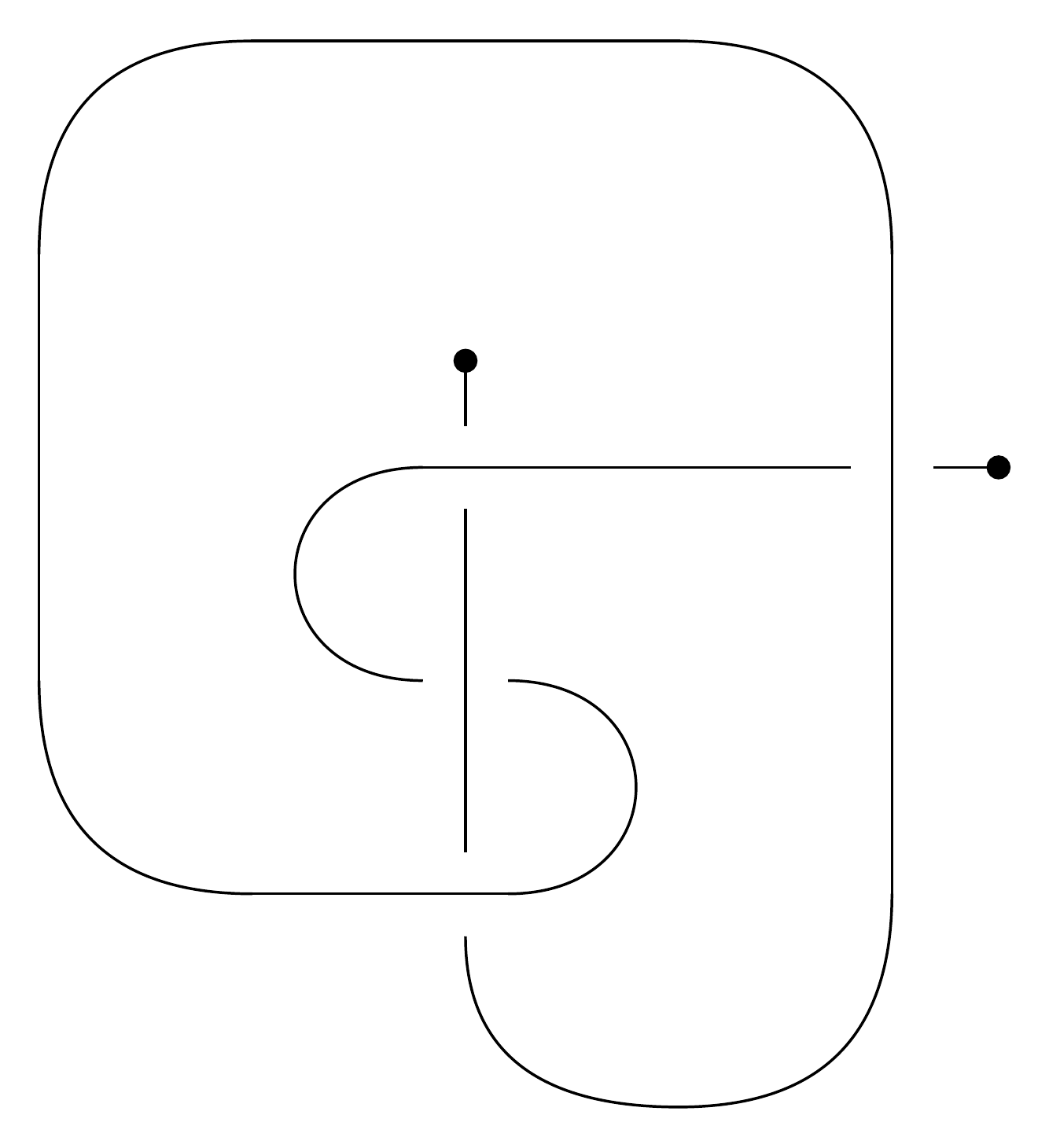}\\
\textcolor{black}{$4_{8}$}
\vspace{1cm}
\end{minipage}
\begin{minipage}[t]{.25\linewidth}
\centering
\includegraphics[width=0.9\textwidth,height=3.5cm,keepaspectratio]{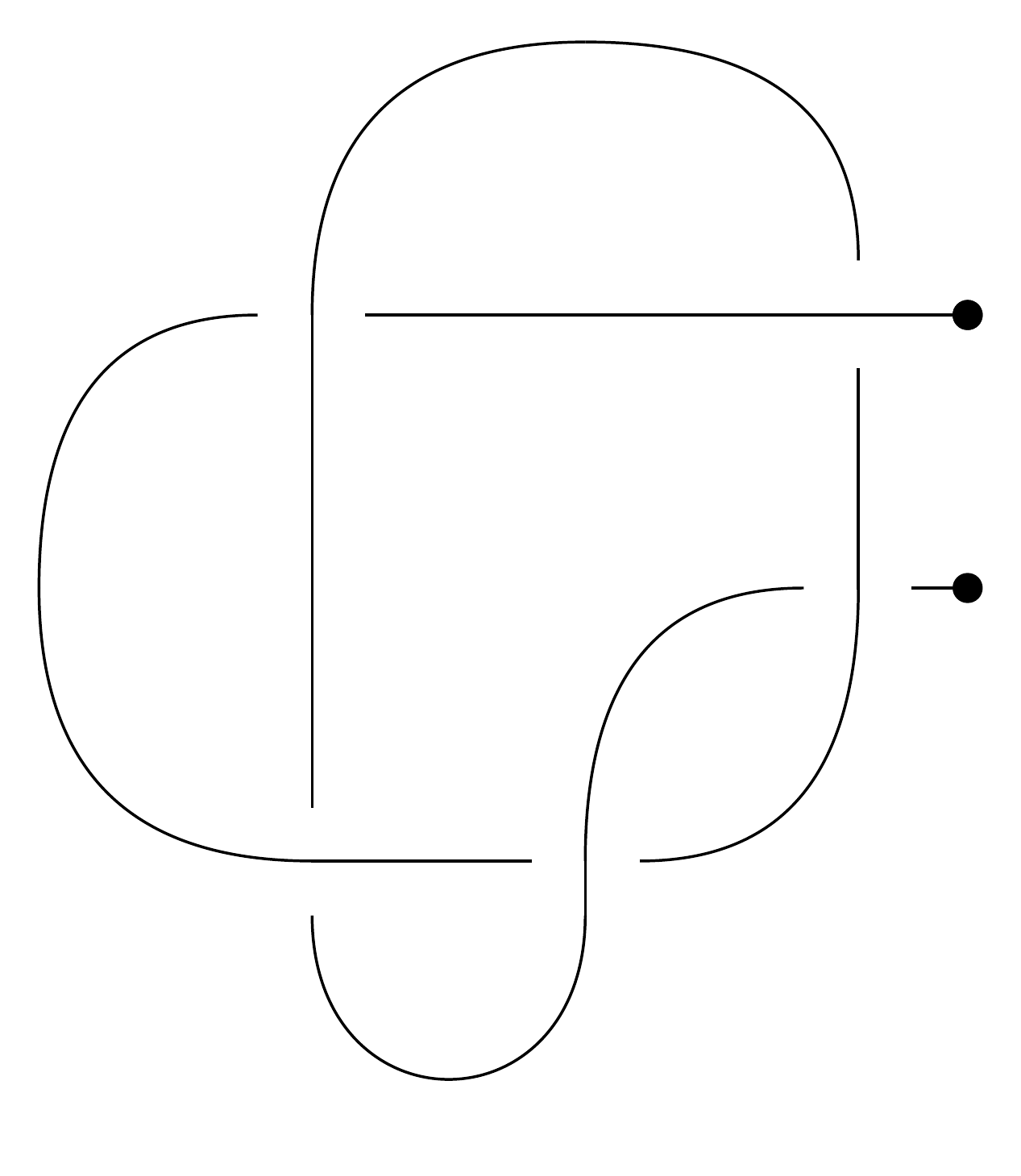}\\
\textcolor{red}{$5_{1}$}
\vspace{1cm}
\end{minipage}
\begin{minipage}[t]{.25\linewidth}
\centering
\includegraphics[width=0.9\textwidth,height=3.5cm,keepaspectratio]{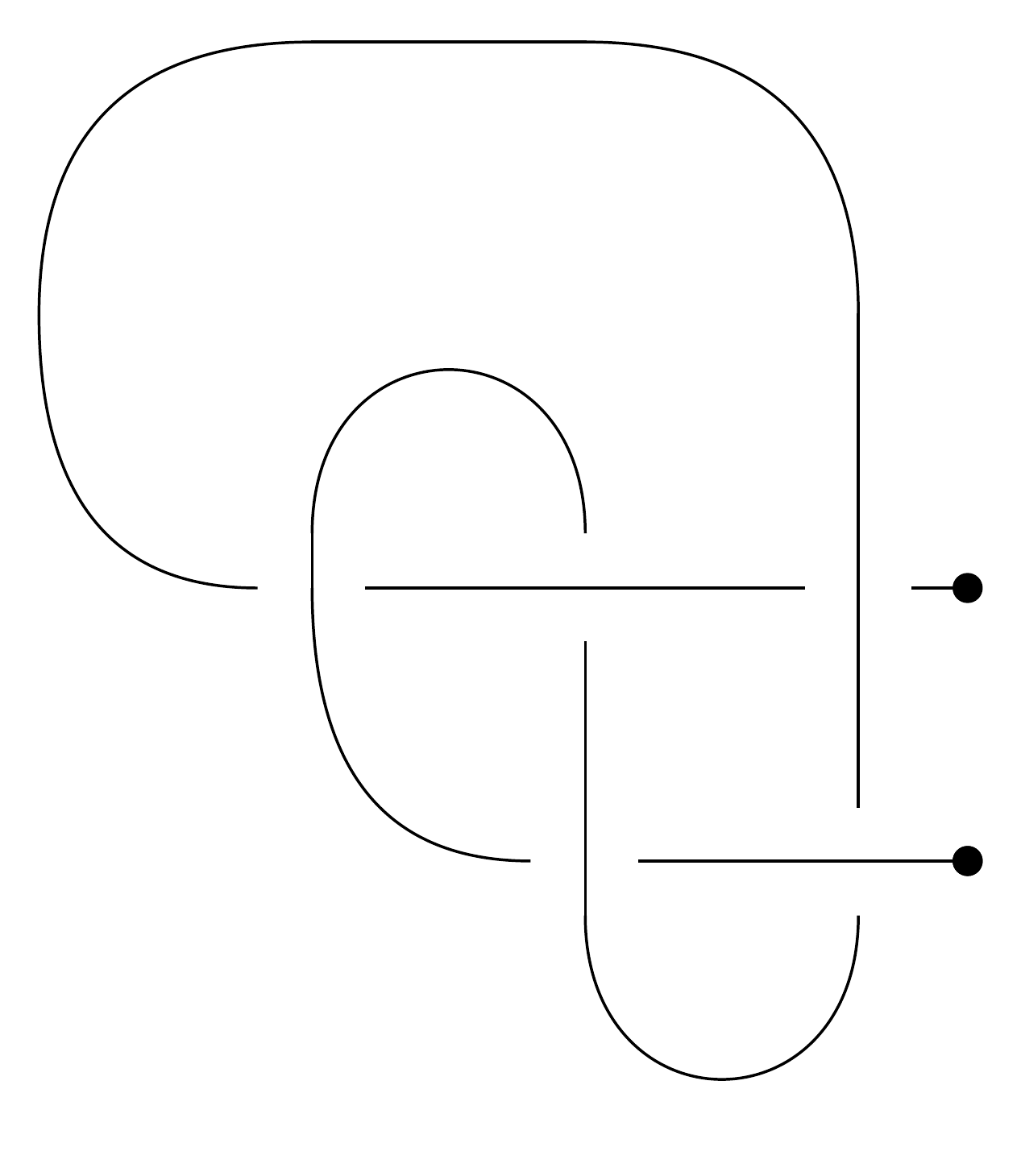}\\
\textcolor{red}{$5_{2}$}
\vspace{1cm}
\end{minipage}
\begin{minipage}[t]{.25\linewidth}
\centering
\includegraphics[width=0.9\textwidth,height=3.5cm,keepaspectratio]{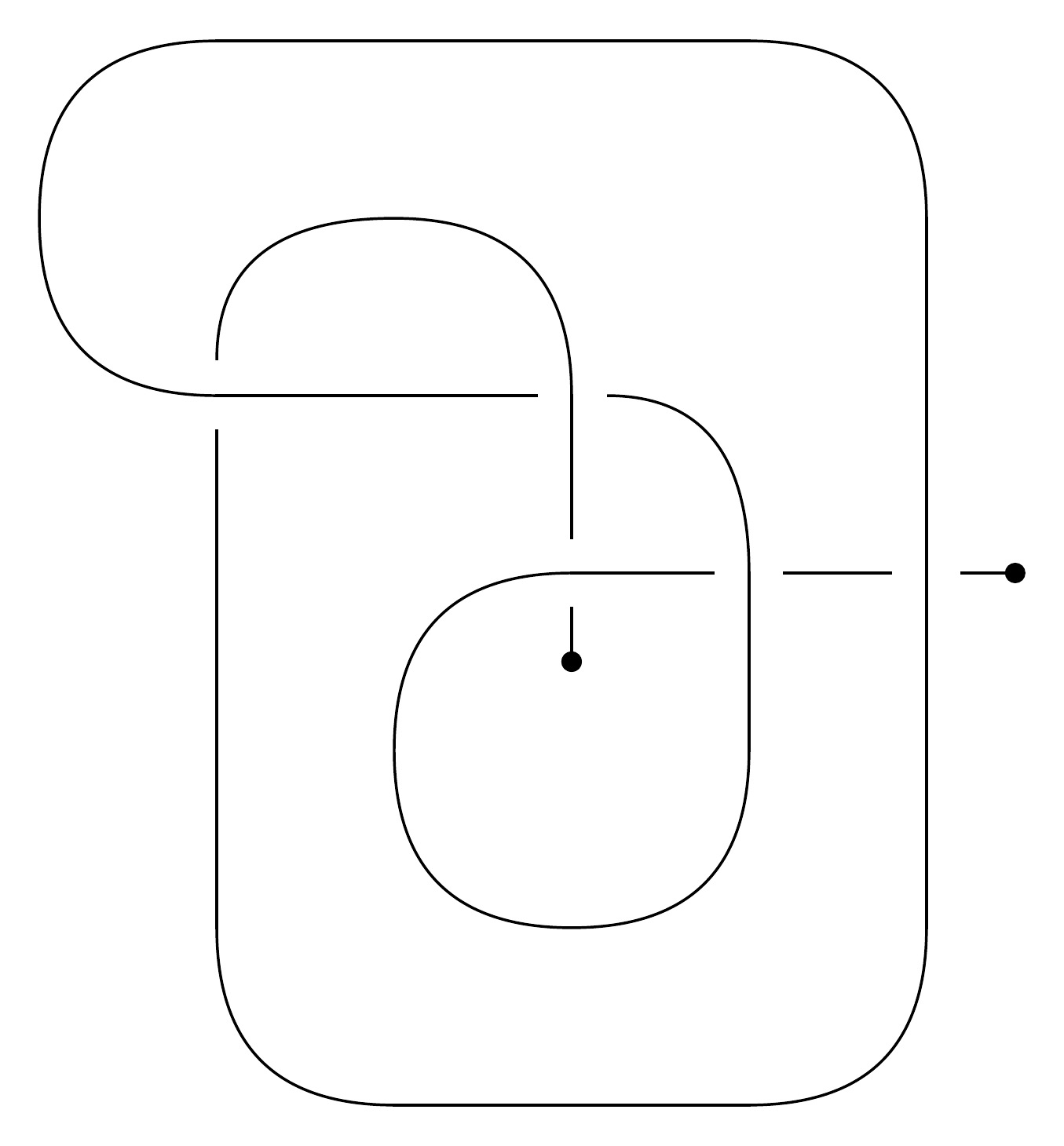}\\
\textcolor{black}{$5_{3}$}
\vspace{1cm}
\end{minipage}
\begin{minipage}[t]{.25\linewidth}
\centering
\includegraphics[width=0.9\textwidth,height=3.5cm,keepaspectratio]{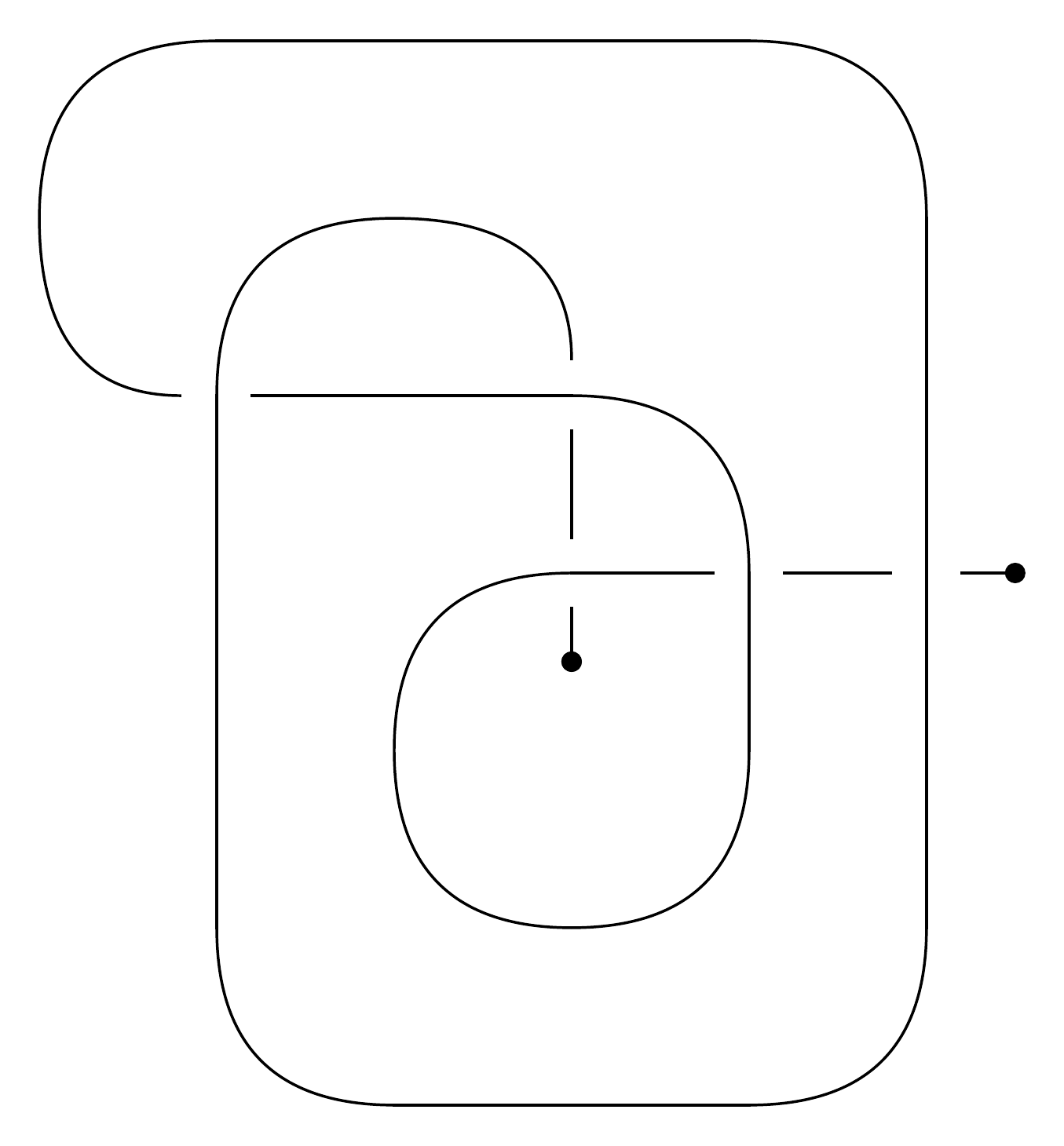}\\
\textcolor{black}{$5_{4}$}
\vspace{1cm}
\end{minipage}
\begin{minipage}[t]{.25\linewidth}
\centering
\includegraphics[width=0.9\textwidth,height=3.5cm,keepaspectratio]{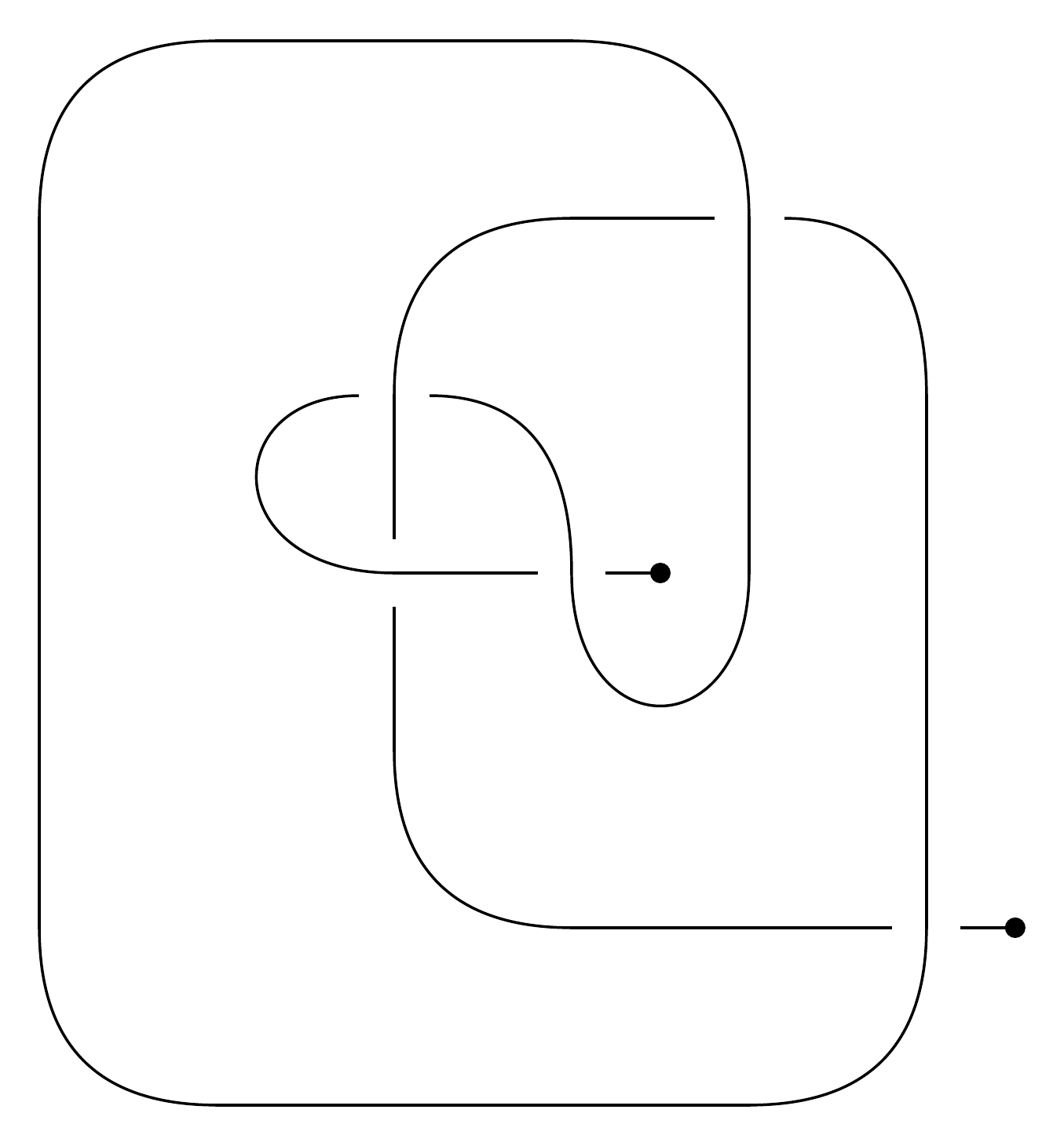}\\
\textcolor{black}{$5_{5}$}
\vspace{1cm}
\end{minipage}
\begin{minipage}[t]{.25\linewidth}
\centering
\includegraphics[width=0.9\textwidth,height=3.5cm,keepaspectratio]{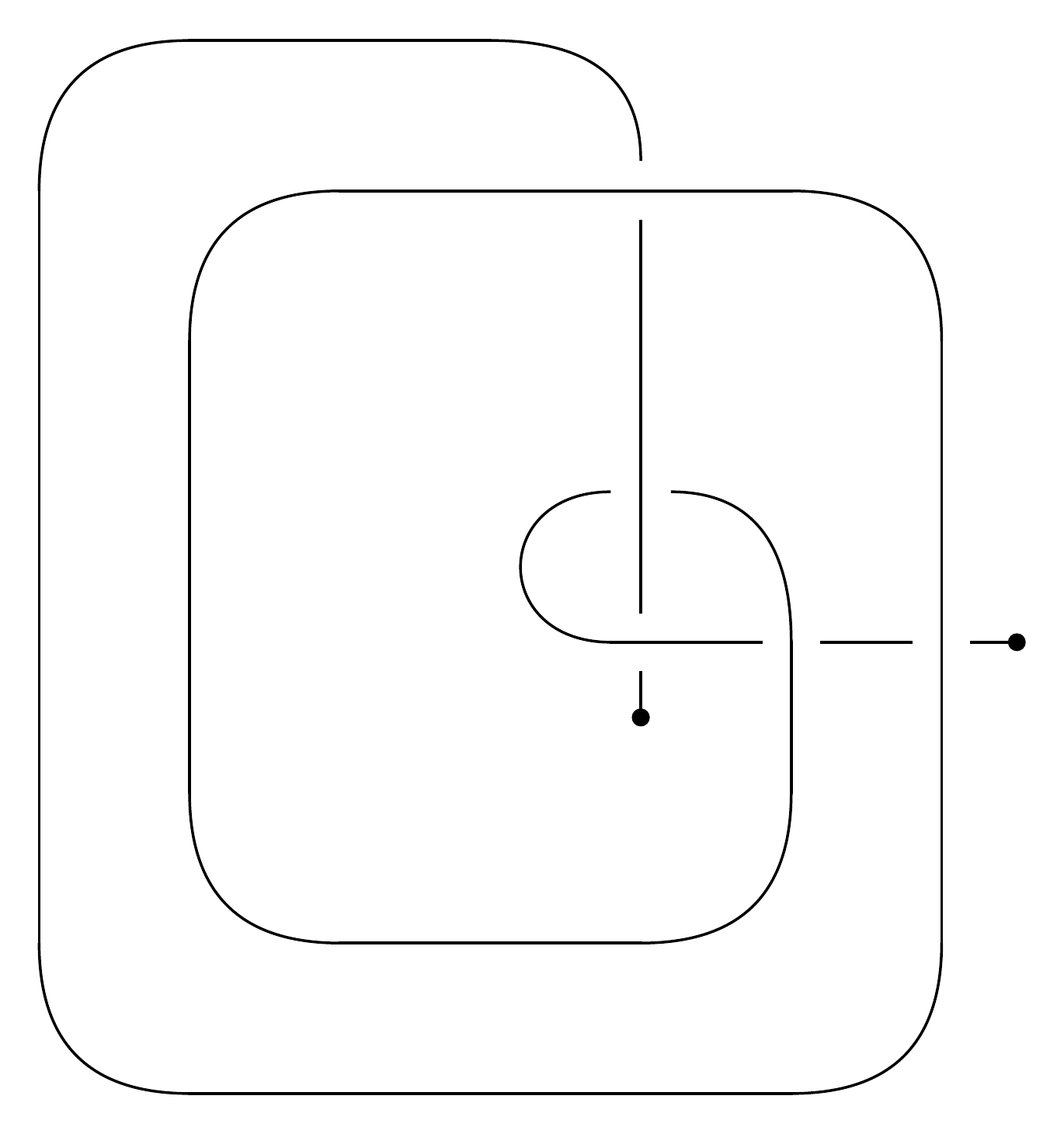}\\
\textcolor{black}{$5_{6}$}
\vspace{1cm}
\end{minipage}
\begin{minipage}[t]{.25\linewidth}
\centering
\includegraphics[width=0.9\textwidth,height=3.5cm,keepaspectratio]{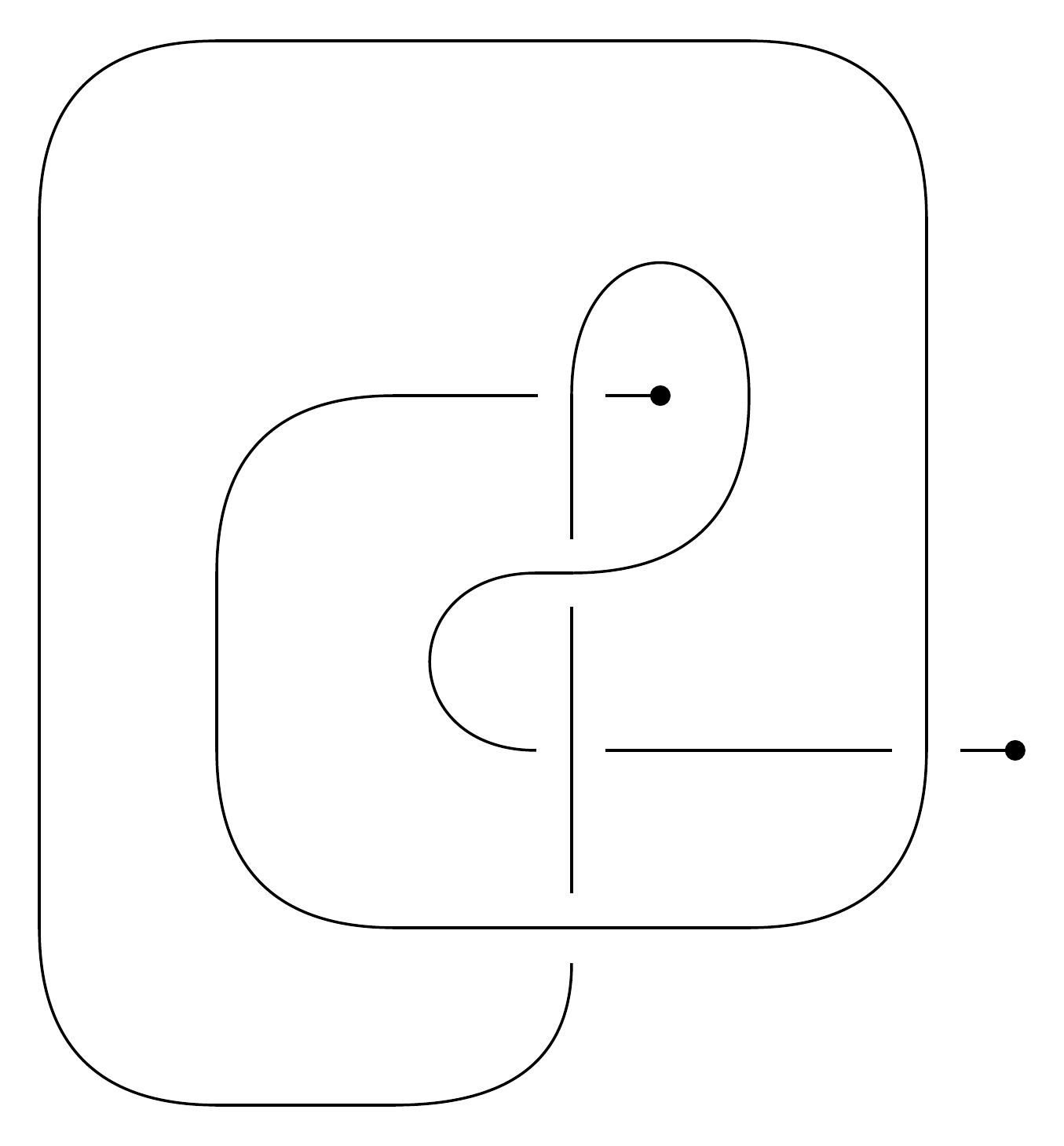}\\
\textcolor{black}{$5_{7}$}
\vspace{1cm}
\end{minipage}
\begin{minipage}[t]{.25\linewidth}
\centering
\includegraphics[width=0.9\textwidth,height=3.5cm,keepaspectratio]{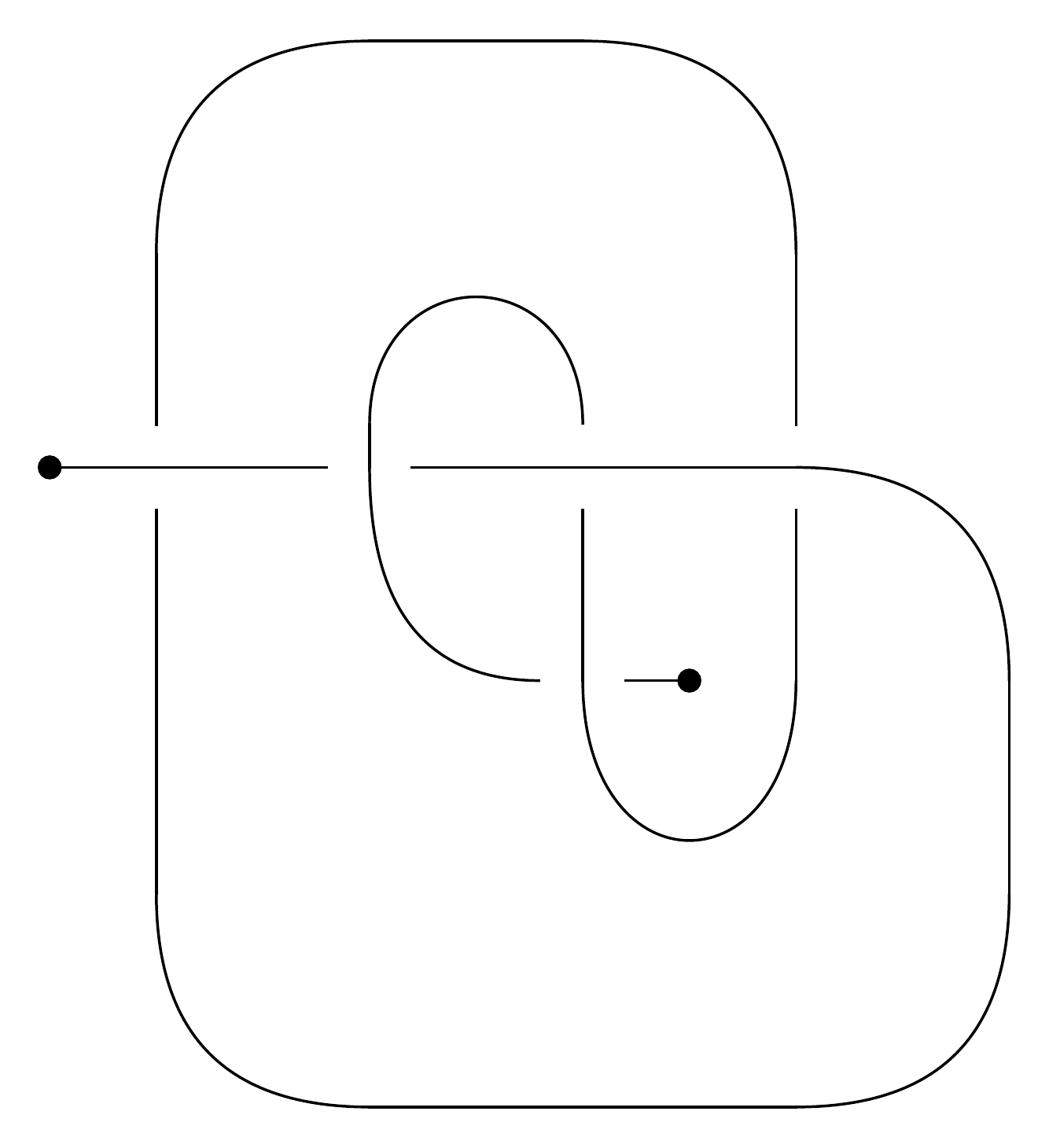}\\
\textcolor{black}{$5_{8}$}
\vspace{1cm}
\end{minipage}
\begin{minipage}[t]{.25\linewidth}
\centering
\includegraphics[width=0.9\textwidth,height=3.5cm,keepaspectratio]{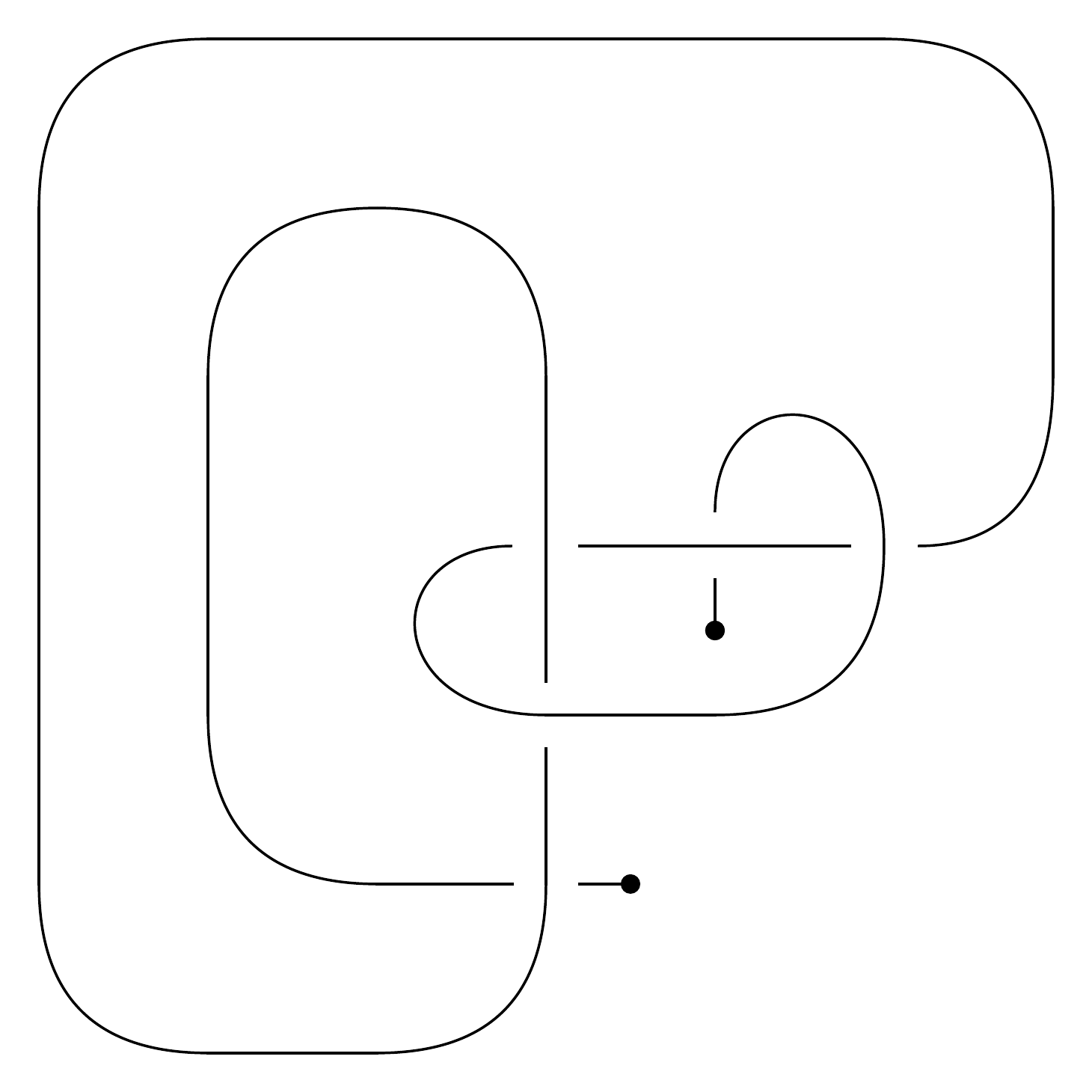}\\
\textcolor{black}{$5_{9}$}
\vspace{1cm}
\end{minipage}
\begin{minipage}[t]{.25\linewidth}
\centering
\includegraphics[width=0.9\textwidth,height=3.5cm,keepaspectratio]{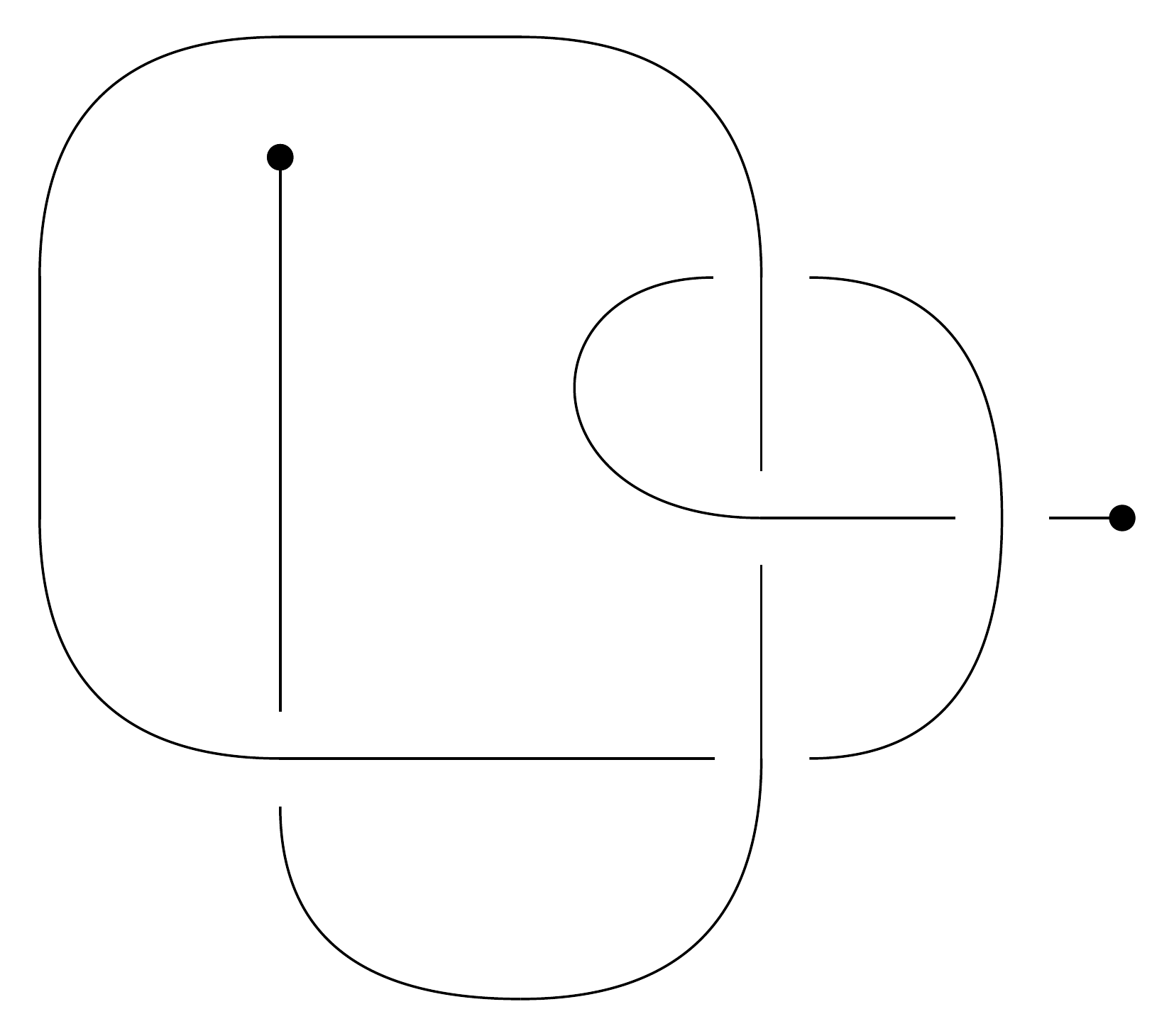}\\
\textcolor{black}{$5_{10}$}
\vspace{1cm}
\end{minipage}
\begin{minipage}[t]{.25\linewidth}
\centering
\includegraphics[width=0.9\textwidth,height=3.5cm,keepaspectratio]{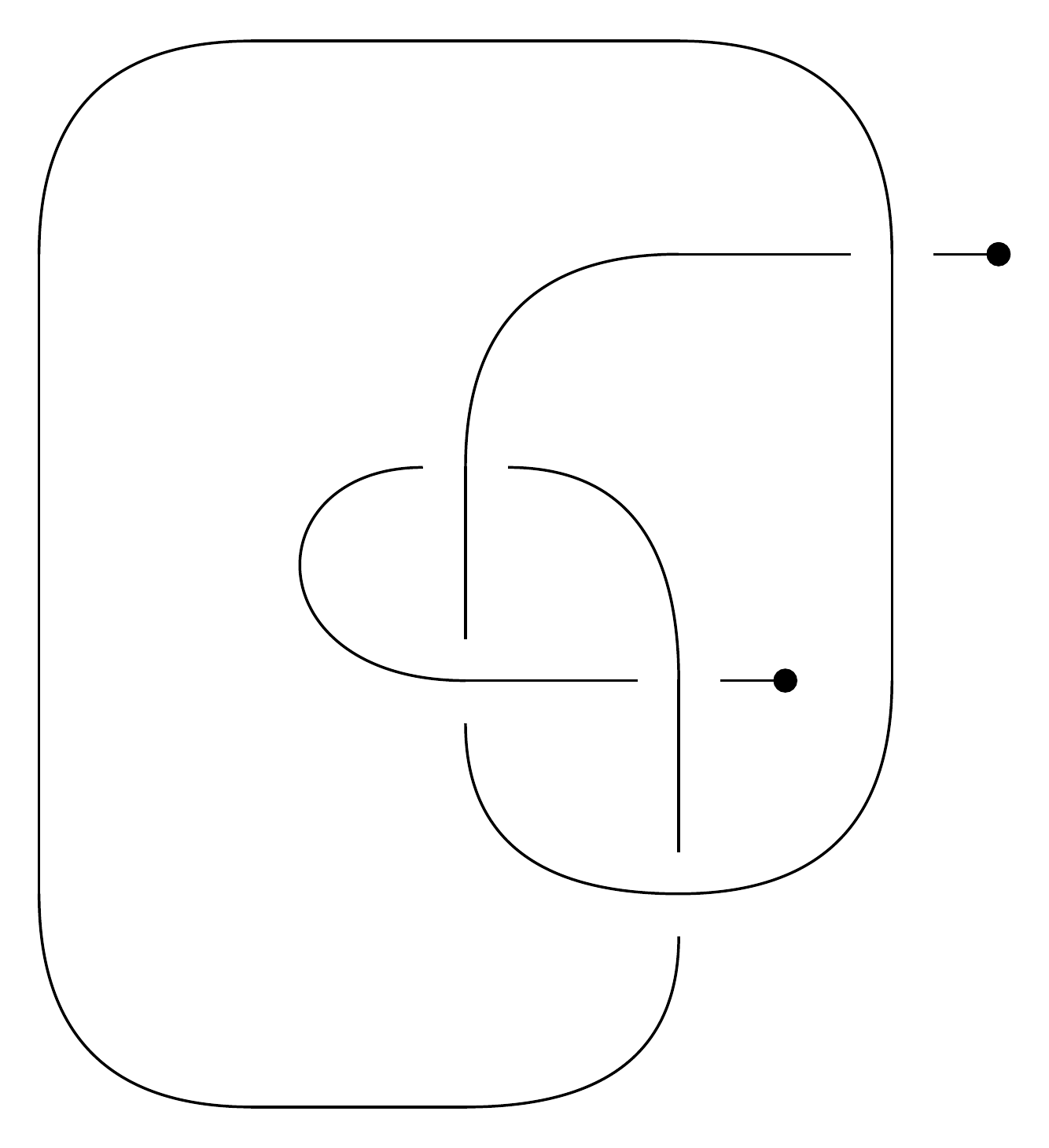}\\
\textcolor{black}{$5_{11}$}
\vspace{1cm}
\end{minipage}
\begin{minipage}[t]{.25\linewidth}
\centering
\includegraphics[width=0.9\textwidth,height=3.5cm,keepaspectratio]{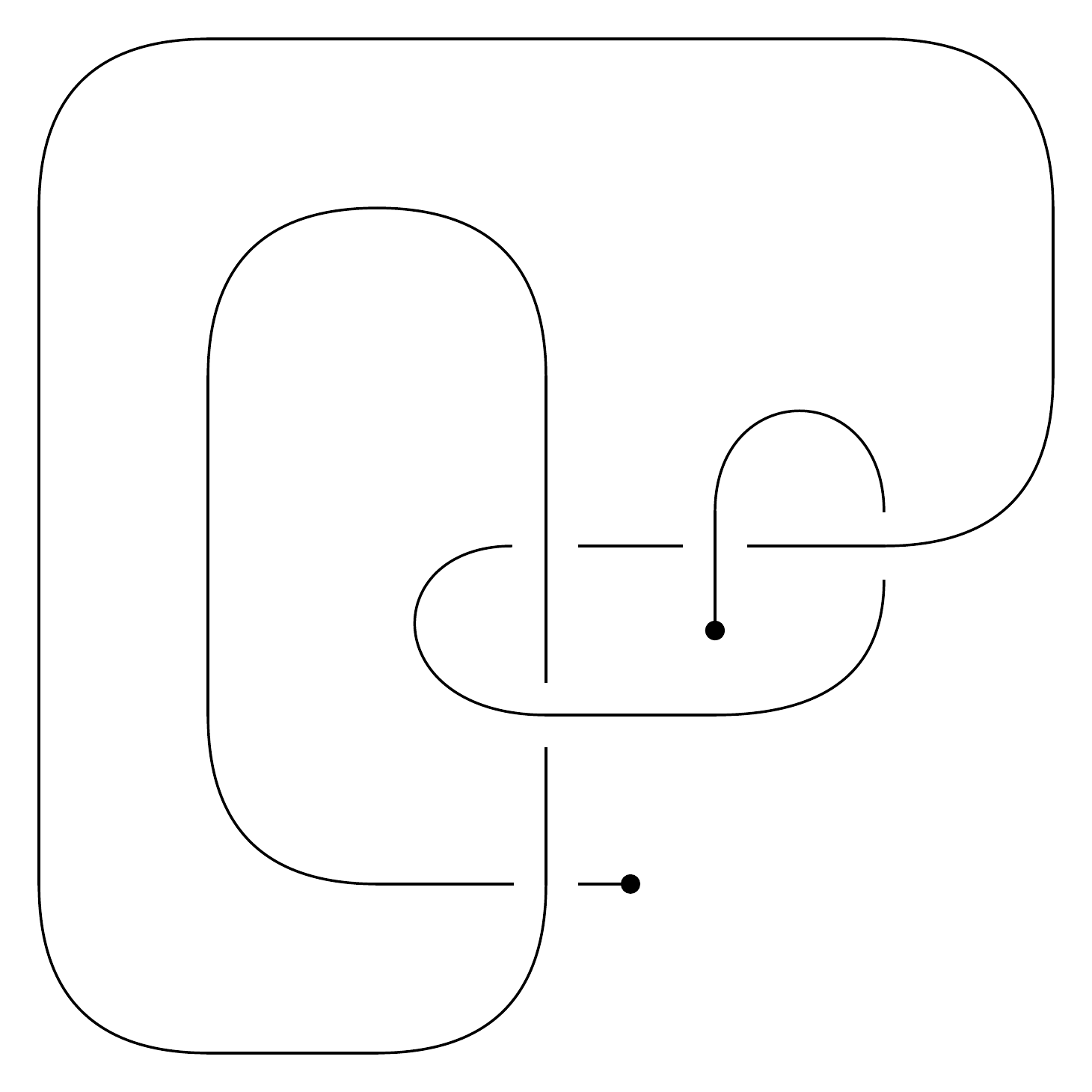}\\
\textcolor{black}{$5_{12}$}
\vspace{1cm}
\end{minipage}
\begin{minipage}[t]{.25\linewidth}
\centering
\includegraphics[width=0.9\textwidth,height=3.5cm,keepaspectratio]{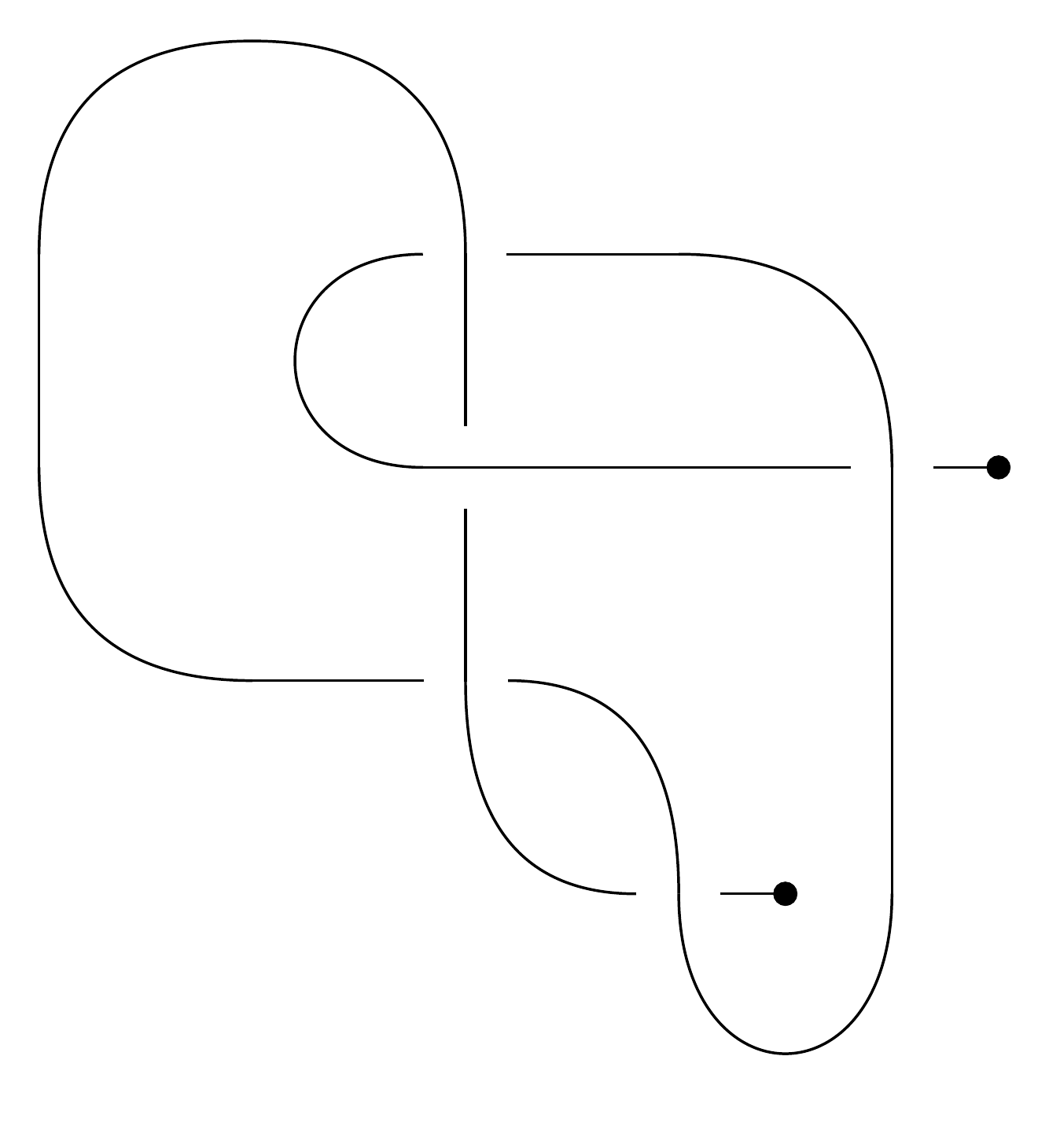}\\
\textcolor{black}{$5_{13}$}
\vspace{1cm}
\end{minipage}
\begin{minipage}[t]{.25\linewidth}
\centering
\includegraphics[width=0.9\textwidth,height=3.5cm,keepaspectratio]{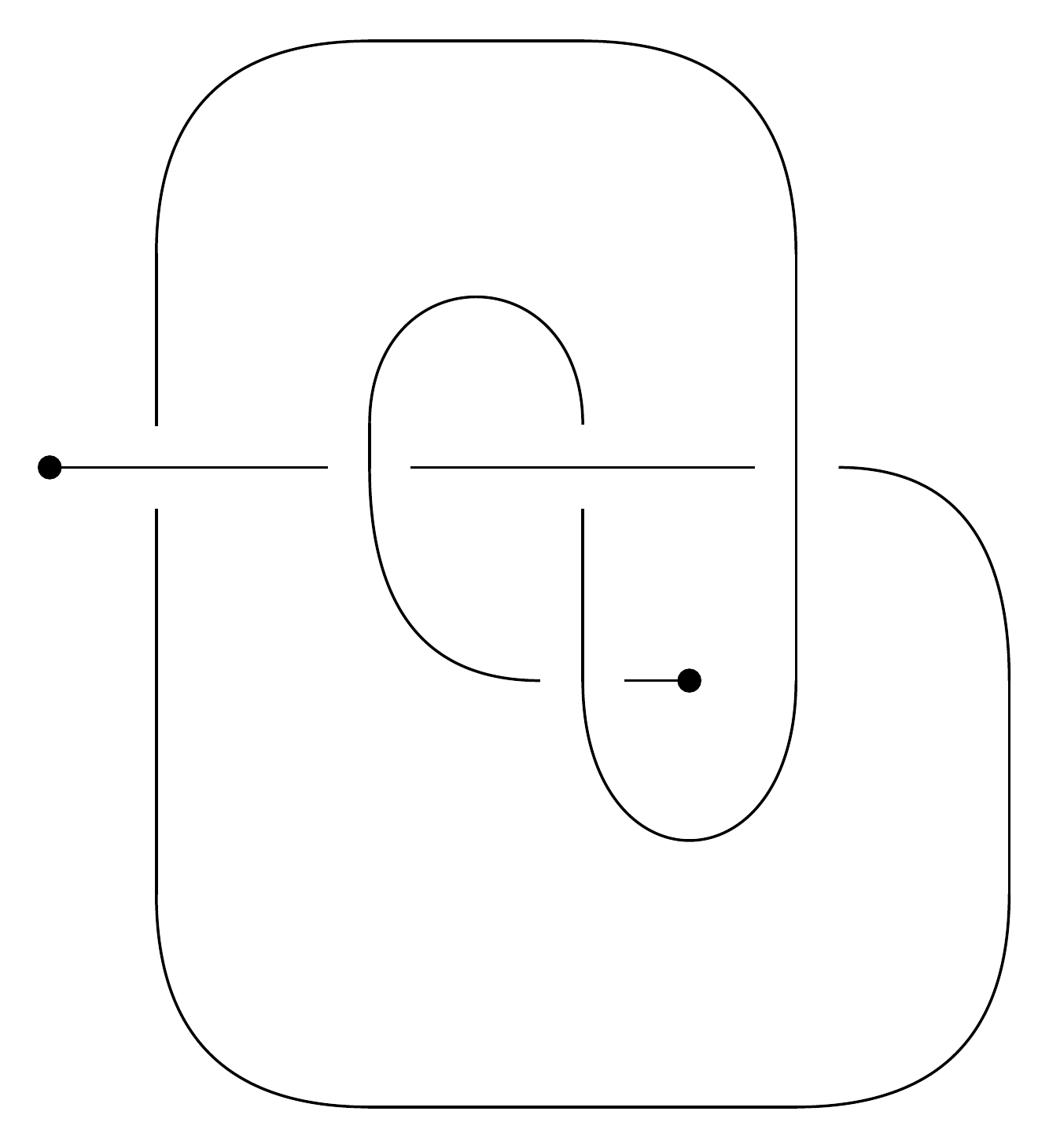}\\
\textcolor{black}{$5_{14}$}
\vspace{1cm}
\end{minipage}
\begin{minipage}[t]{.25\linewidth}
\centering
\includegraphics[width=0.9\textwidth,height=3.5cm,keepaspectratio]{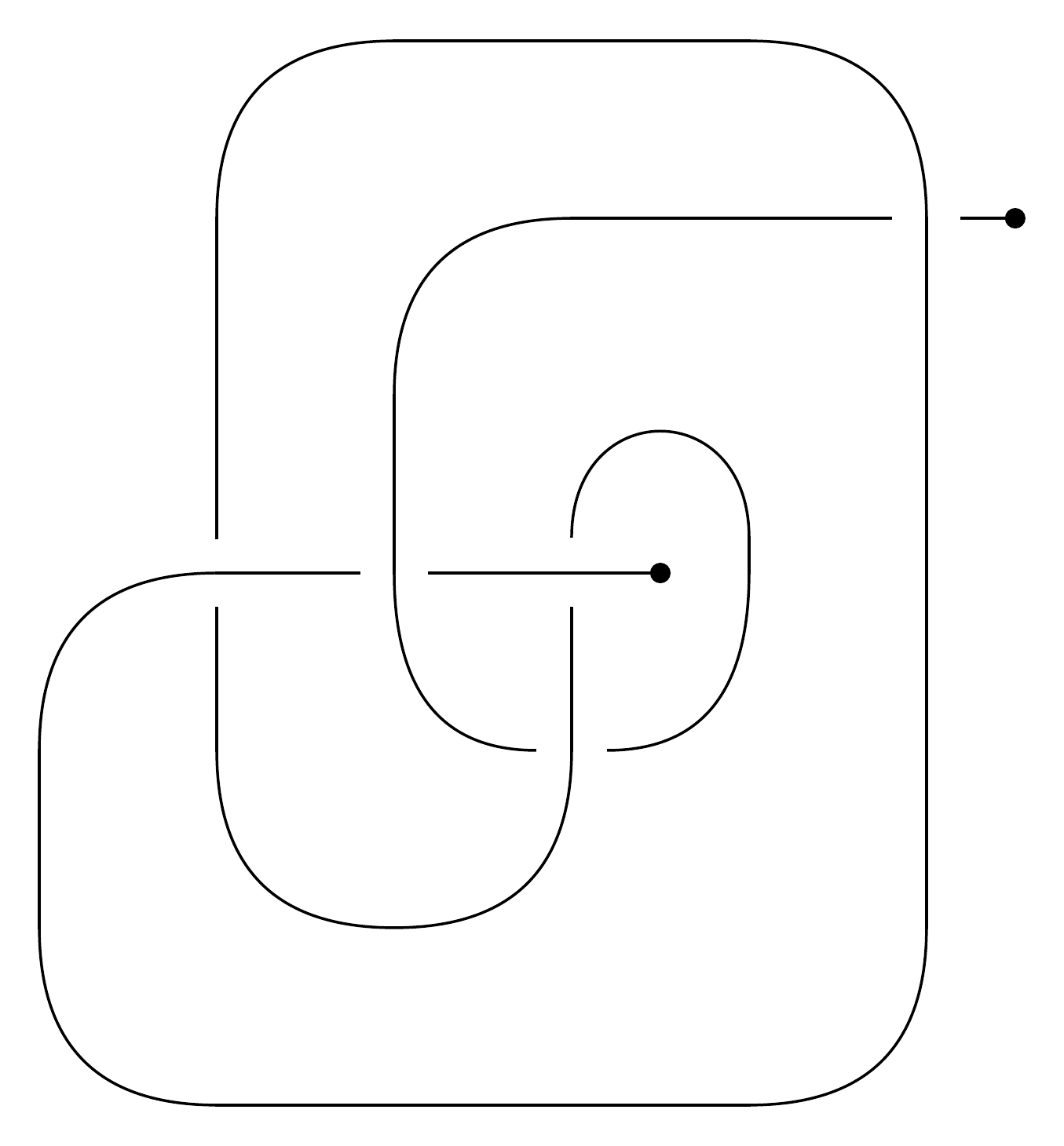}\\
\textcolor{black}{$5_{15}$}
\vspace{1cm}
\end{minipage}
\begin{minipage}[t]{.25\linewidth}
\centering
\includegraphics[width=0.9\textwidth,height=3.5cm,keepaspectratio]{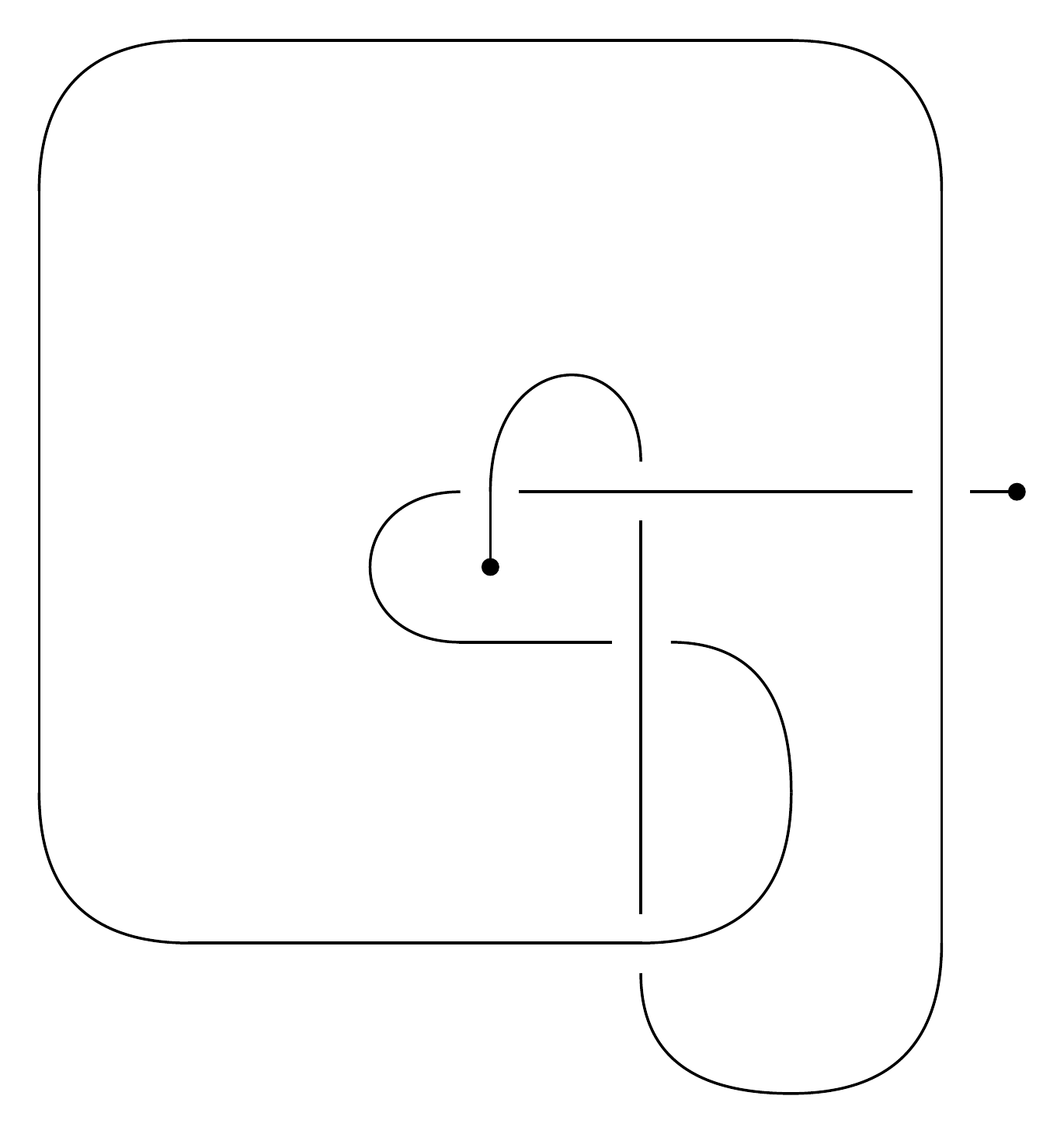}\\
\textcolor{black}{$5_{16}$}
\vspace{1cm}
\end{minipage}
\begin{minipage}[t]{.25\linewidth}
\centering
\includegraphics[width=0.9\textwidth,height=3.5cm,keepaspectratio]{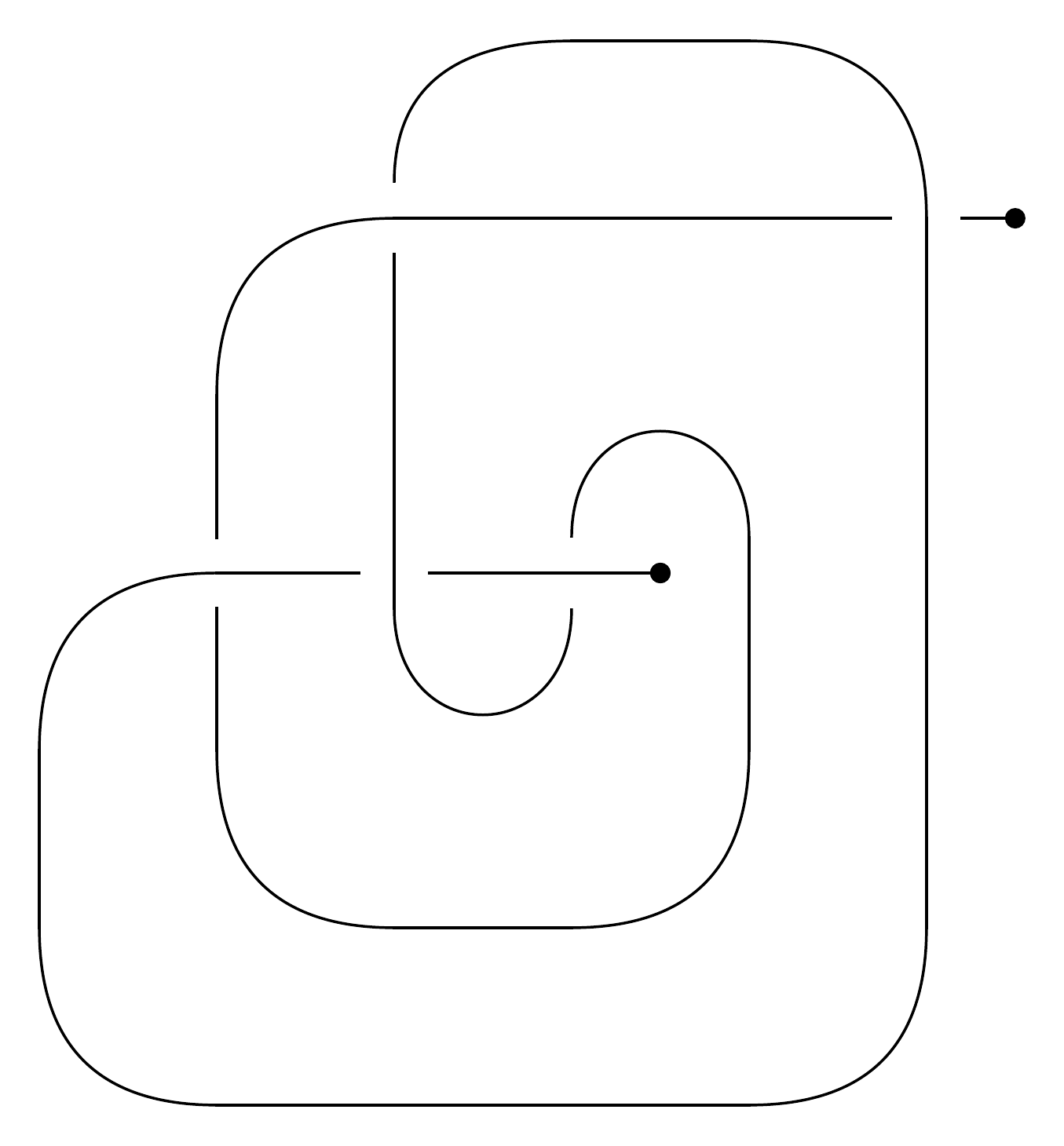}\\
\textcolor{black}{$5_{17}$}
\vspace{1cm}
\end{minipage}
\begin{minipage}[t]{.25\linewidth}
\centering
\includegraphics[width=0.9\textwidth,height=3.5cm,keepaspectratio]{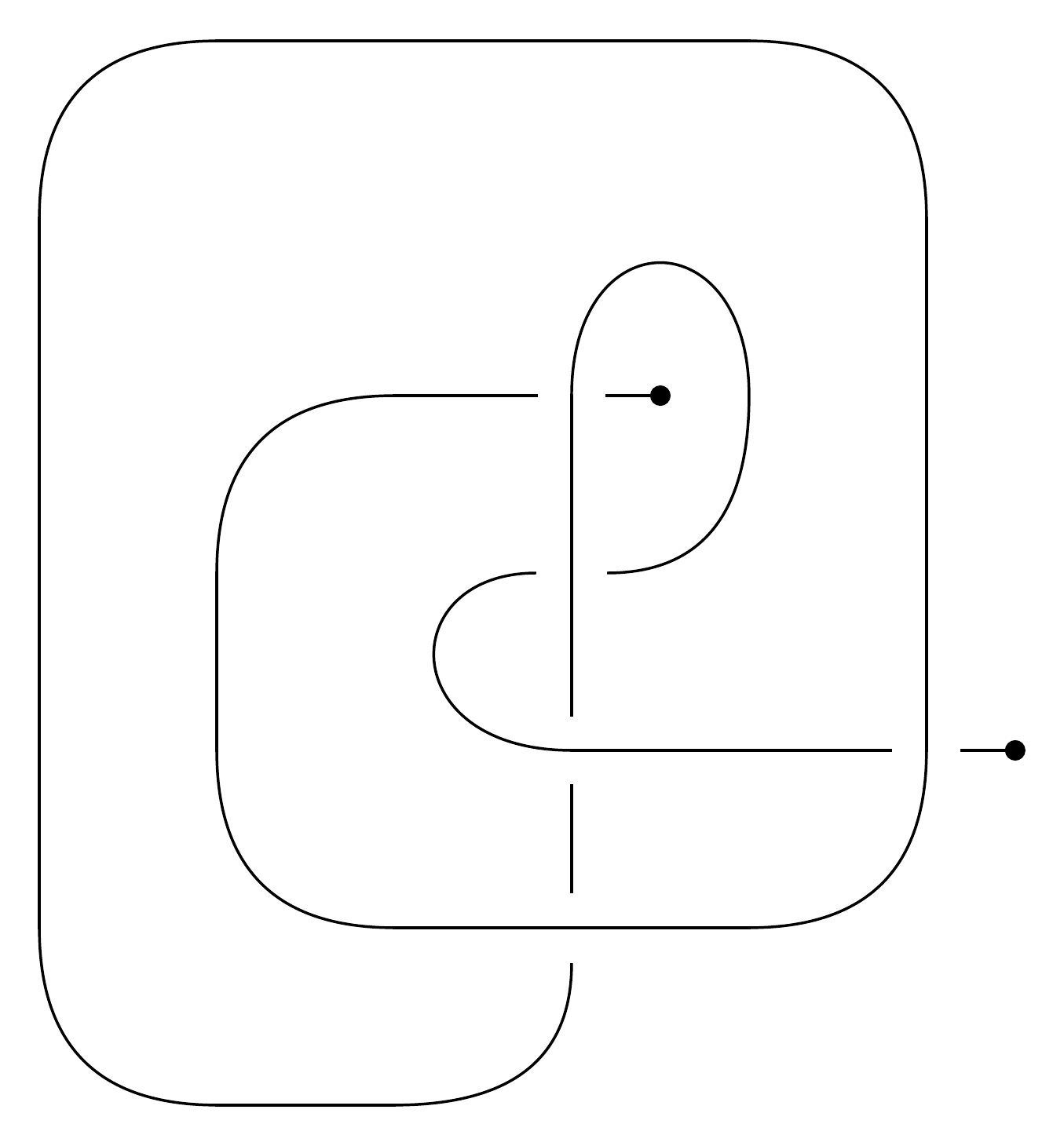}\\
\textcolor{black}{$5_{18}$}
\vspace{1cm}
\end{minipage}
\begin{minipage}[t]{.25\linewidth}
\centering
\includegraphics[width=0.9\textwidth,height=3.5cm,keepaspectratio]{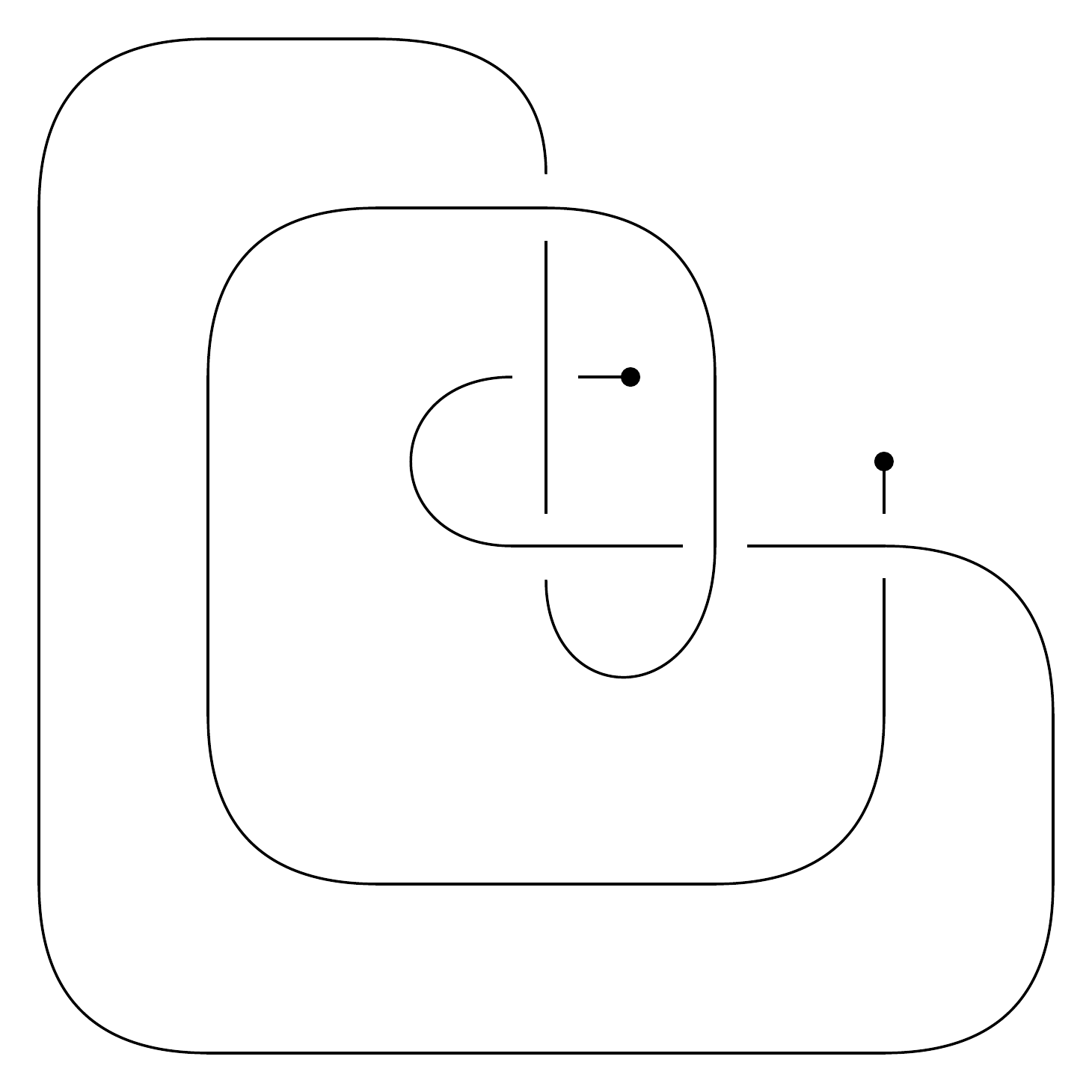}\\
\textcolor{black}{$5_{19}$}
\vspace{1cm}
\end{minipage}
\begin{minipage}[t]{.25\linewidth}
\centering
\includegraphics[width=0.9\textwidth,height=3.5cm,keepaspectratio]{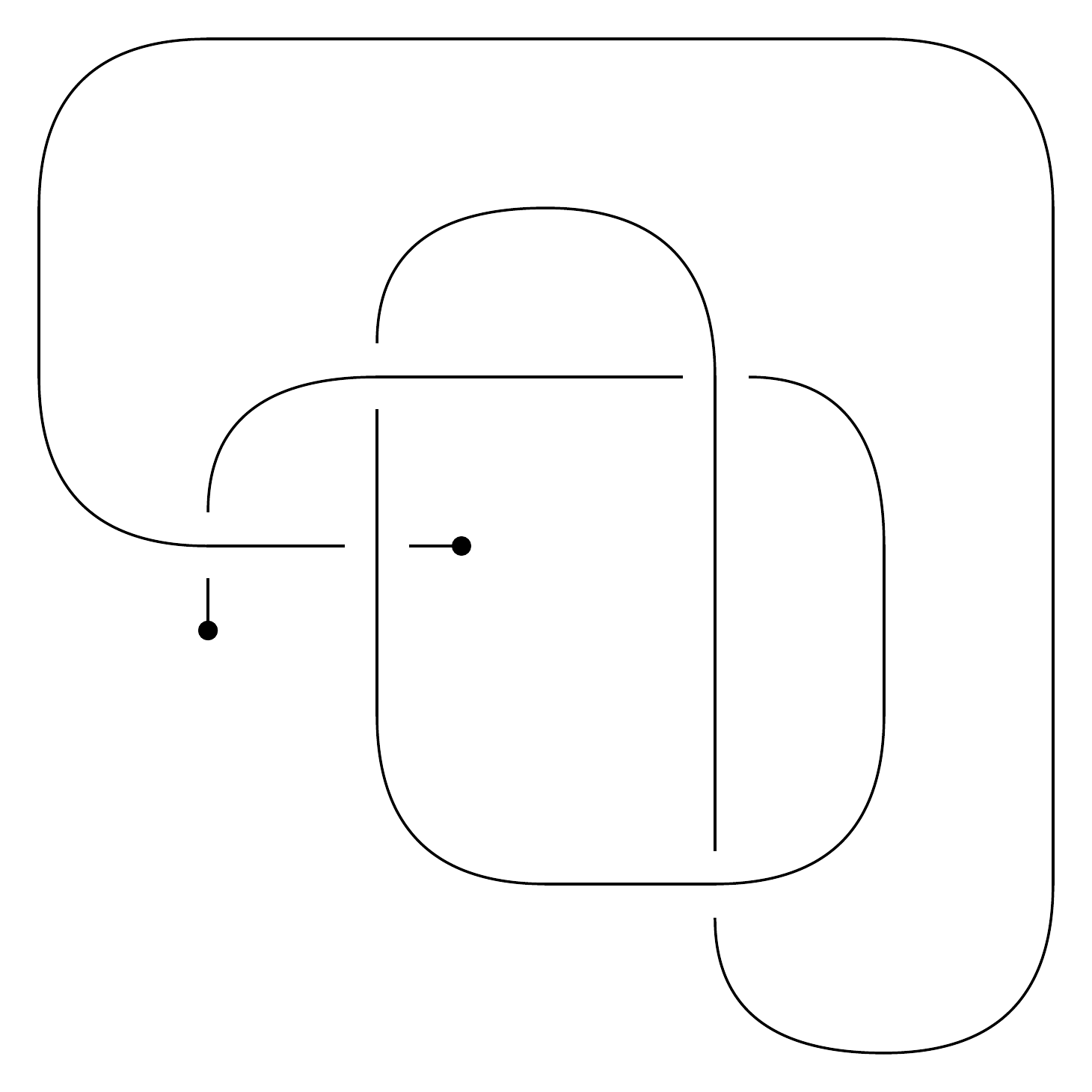}\\
\textcolor{black}{$5_{20}$}
\vspace{1cm}
\end{minipage}
\begin{minipage}[t]{.25\linewidth}
\centering
\includegraphics[width=0.9\textwidth,height=3.5cm,keepaspectratio]{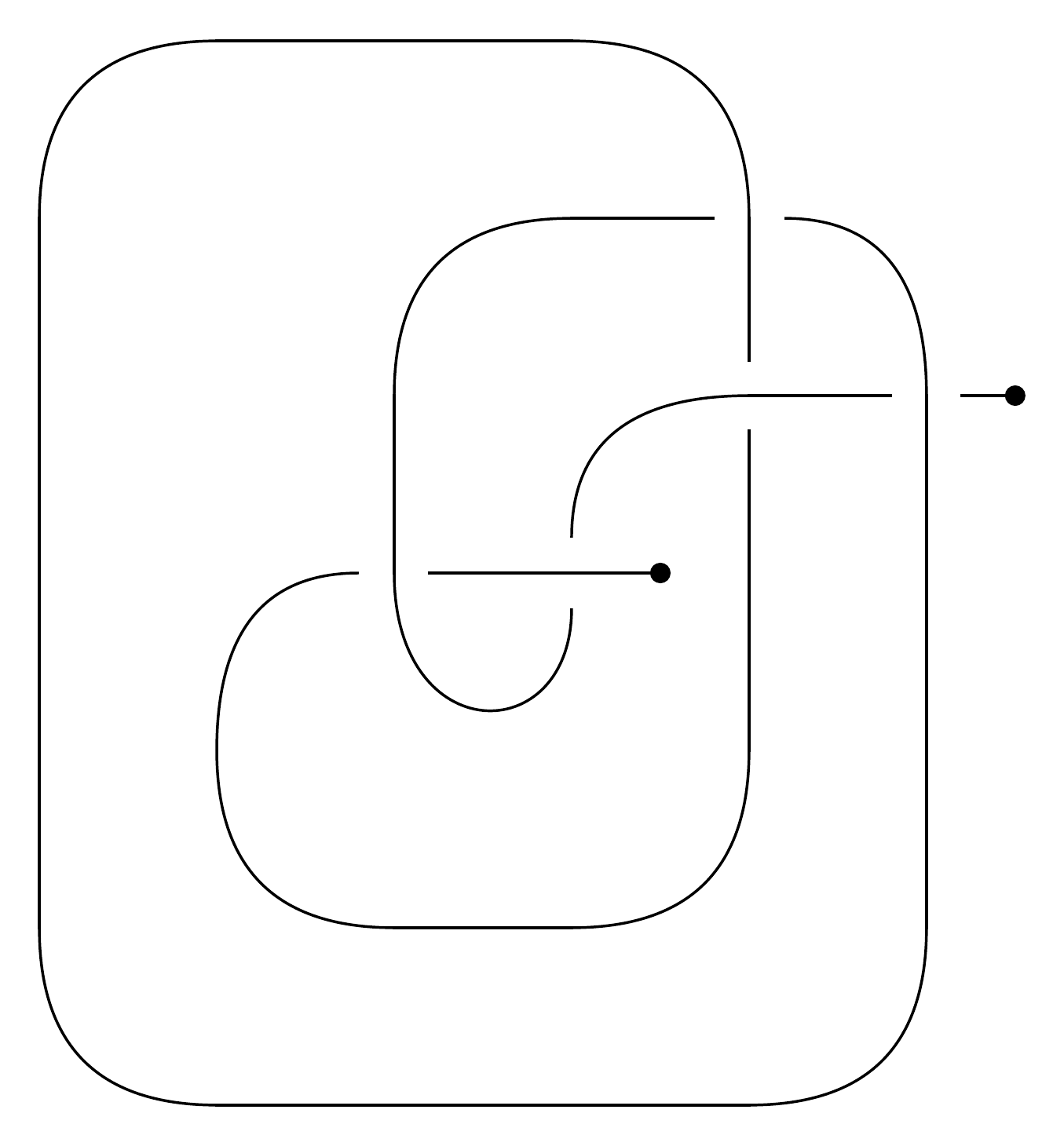}\\
\textcolor{black}{$5_{21}$}
\vspace{1cm}
\end{minipage}
\begin{minipage}[t]{.25\linewidth}
\centering
\includegraphics[width=0.9\textwidth,height=3.5cm,keepaspectratio]{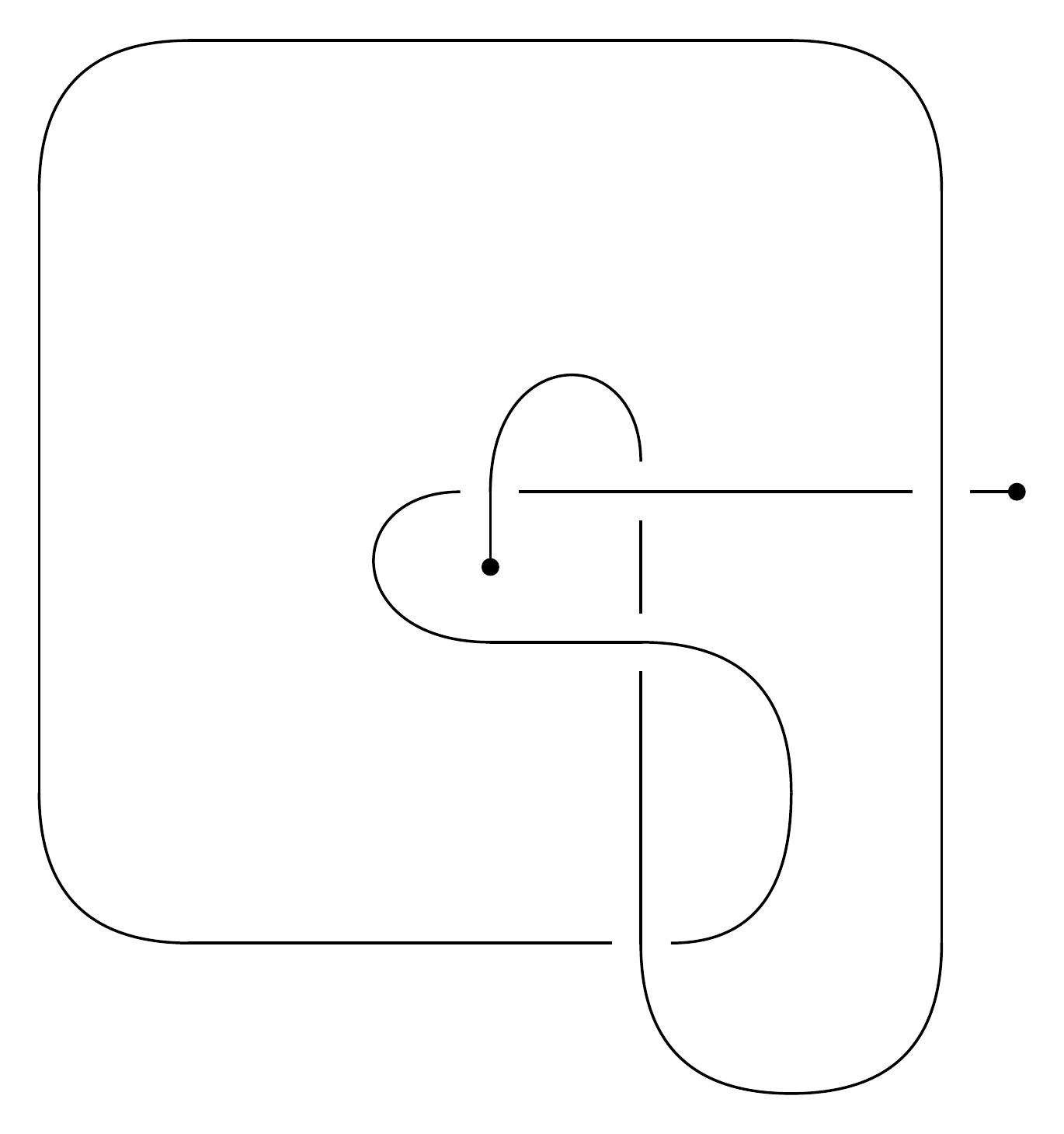}\\
\textcolor{black}{$5_{22}$}
\vspace{1cm}
\end{minipage}
\begin{minipage}[t]{.25\linewidth}
\centering
\includegraphics[width=0.9\textwidth,height=3.5cm,keepaspectratio]{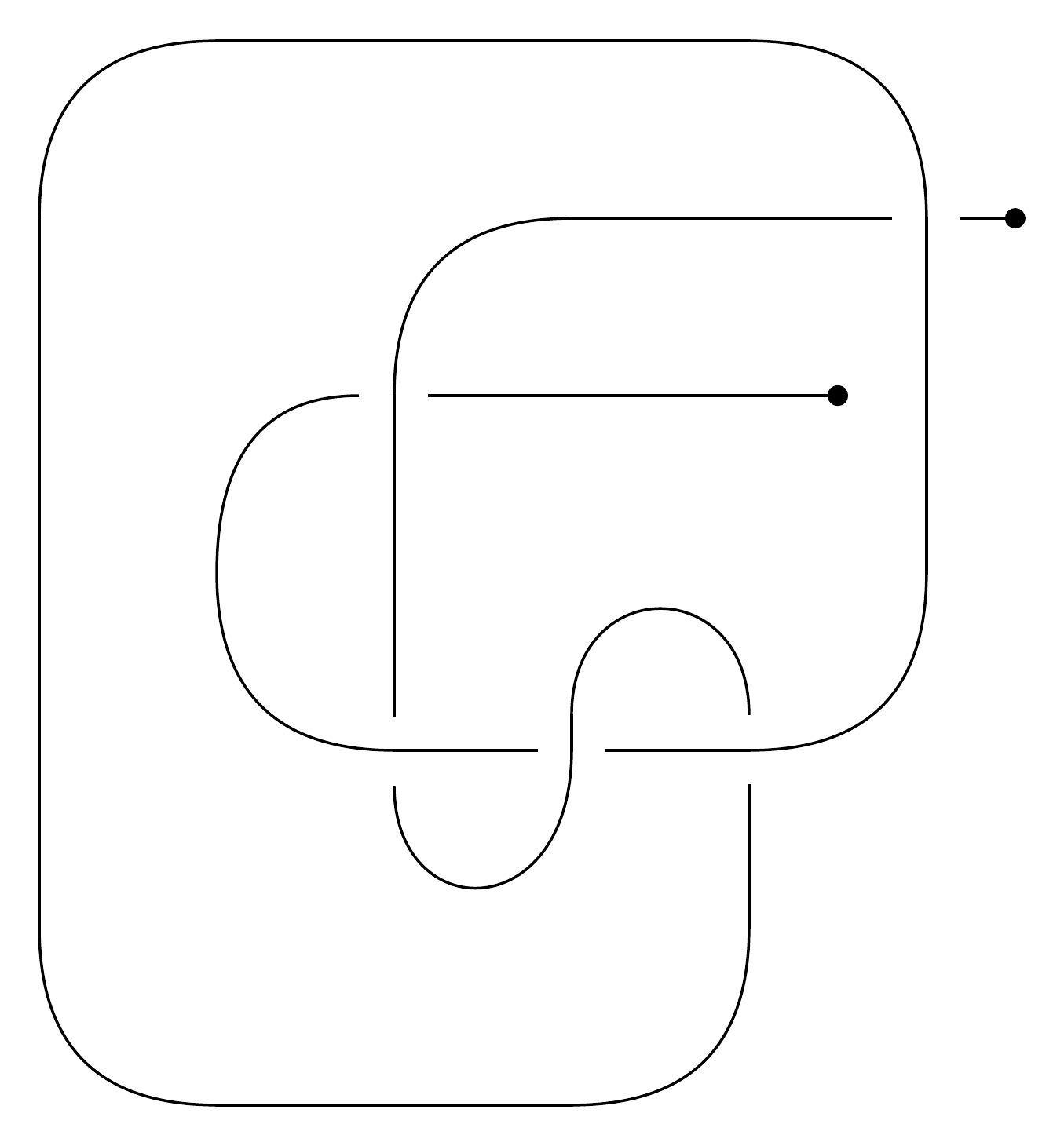}\\
\textcolor{black}{$5_{23}$}
\vspace{1cm}
\end{minipage}
\begin{minipage}[t]{.25\linewidth}
\centering
\includegraphics[width=0.9\textwidth,height=3.5cm,keepaspectratio]{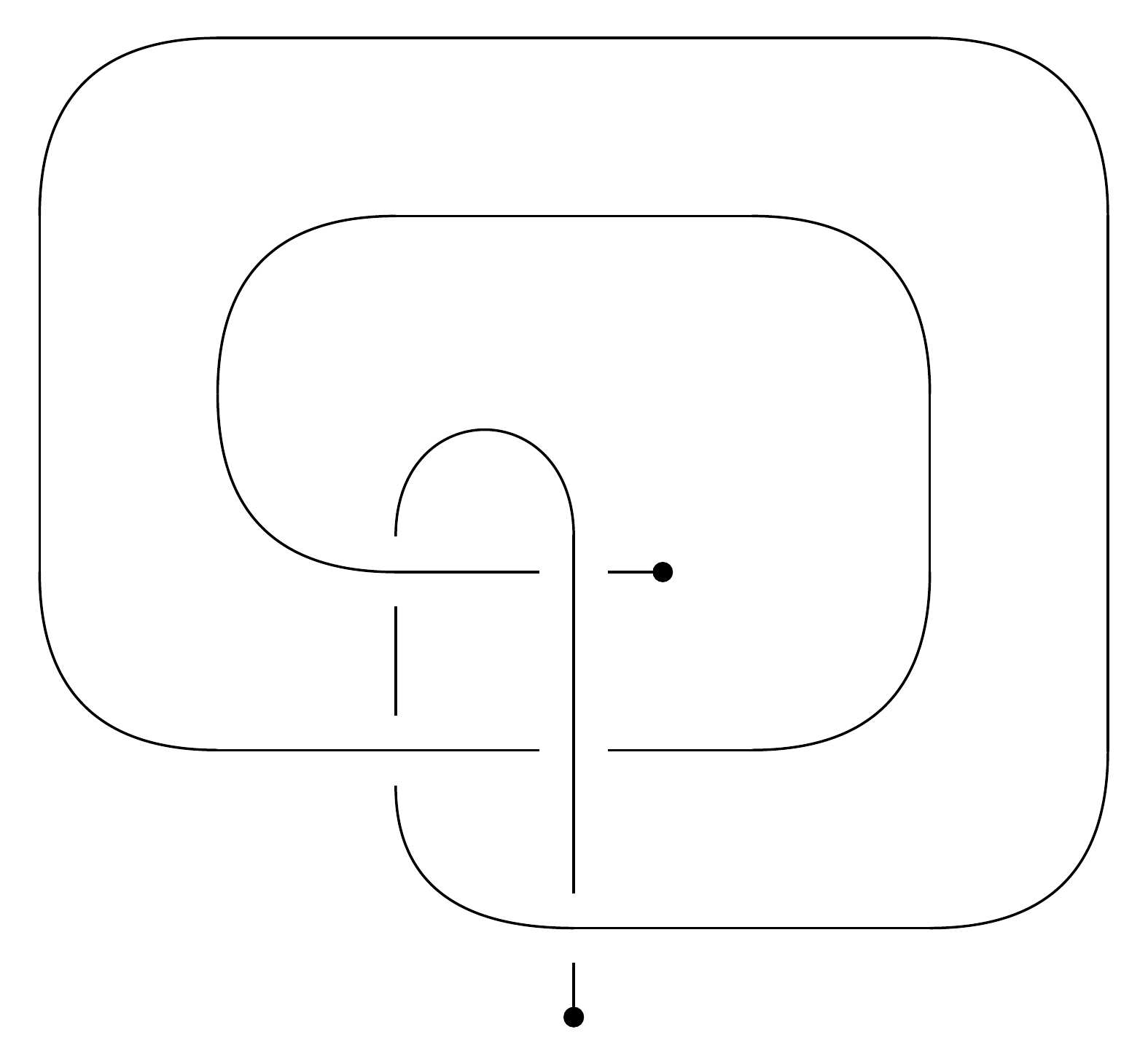}\\
\textcolor{black}{$5_{24}$}
\vspace{1cm}
\end{minipage}
\begin{minipage}[t]{.25\linewidth}
\centering
\includegraphics[width=0.9\textwidth,height=3.5cm,keepaspectratio]{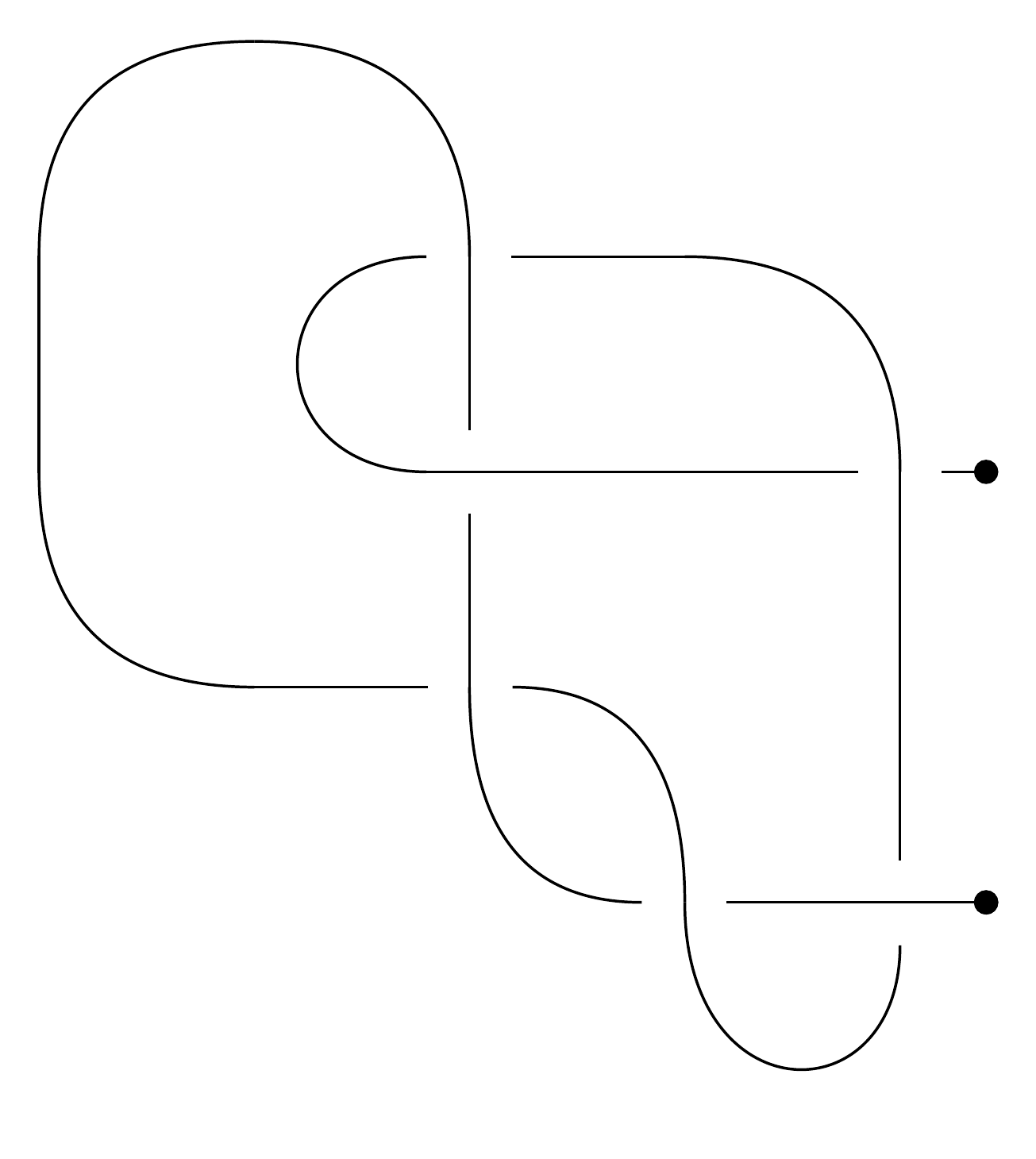}\\
\textcolor{red}{$6_{1}$}
\vspace{1cm}
\end{minipage}
\begin{minipage}[t]{.25\linewidth}
\centering
\includegraphics[width=0.9\textwidth,height=3.5cm,keepaspectratio]{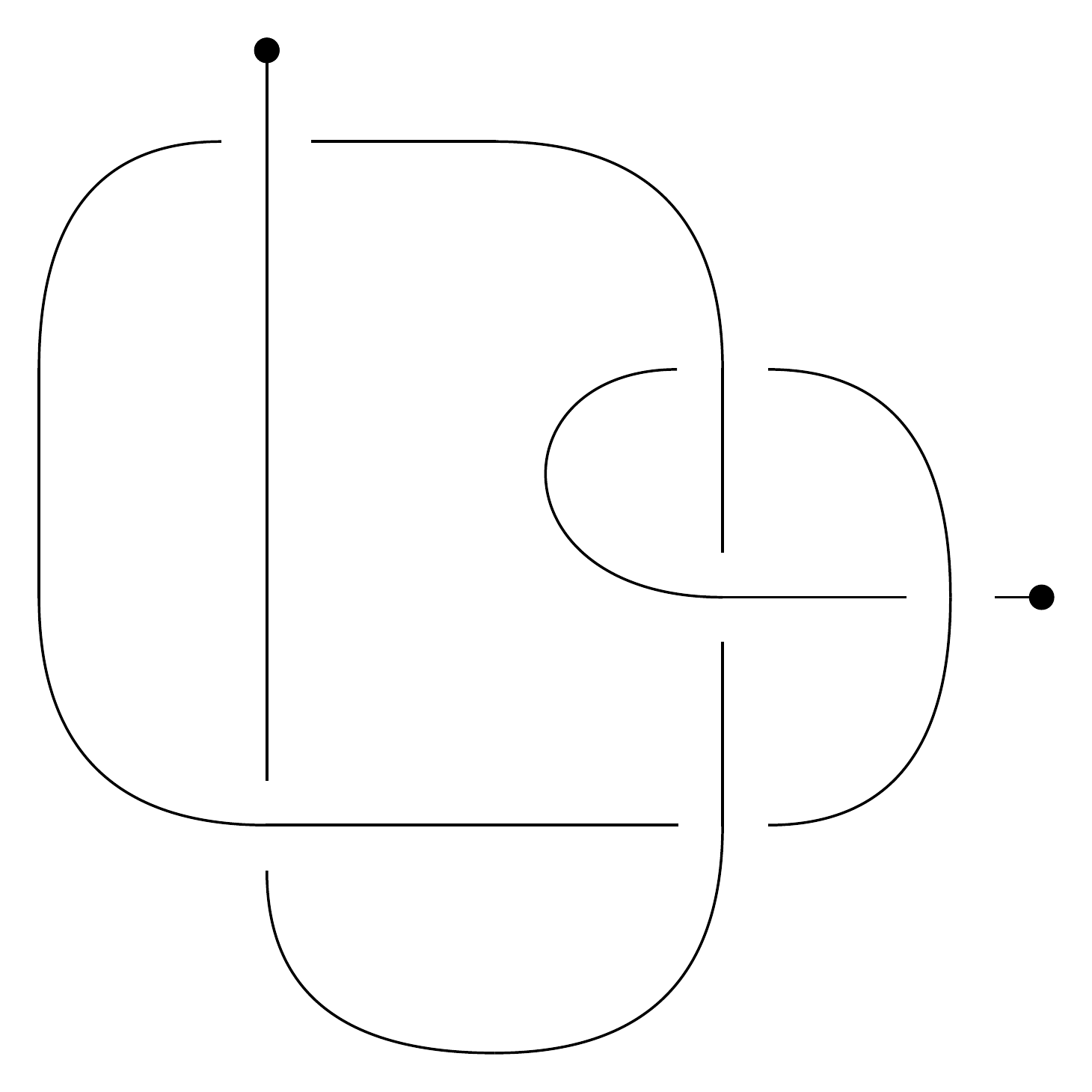}\\
\textcolor{red}{$6_{2}$}
\vspace{1cm}
\end{minipage}
\begin{minipage}[t]{.25\linewidth}
\centering
\includegraphics[width=0.9\textwidth,height=3.5cm,keepaspectratio]{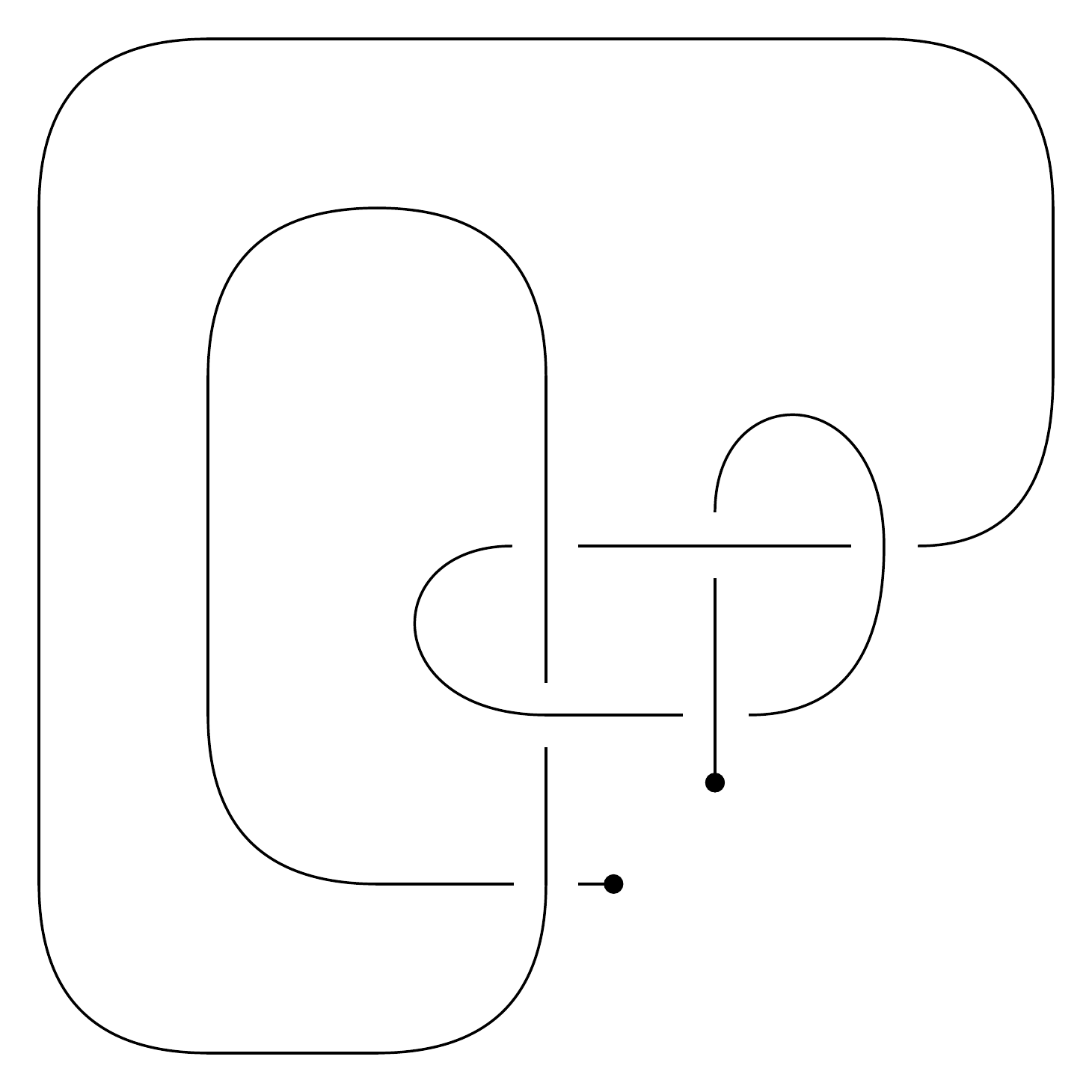}\\
\textcolor{red}{$6_{3}$}
\vspace{1cm}
\end{minipage}
\begin{minipage}[t]{.25\linewidth}
\centering
\includegraphics[width=0.9\textwidth,height=3.5cm,keepaspectratio]{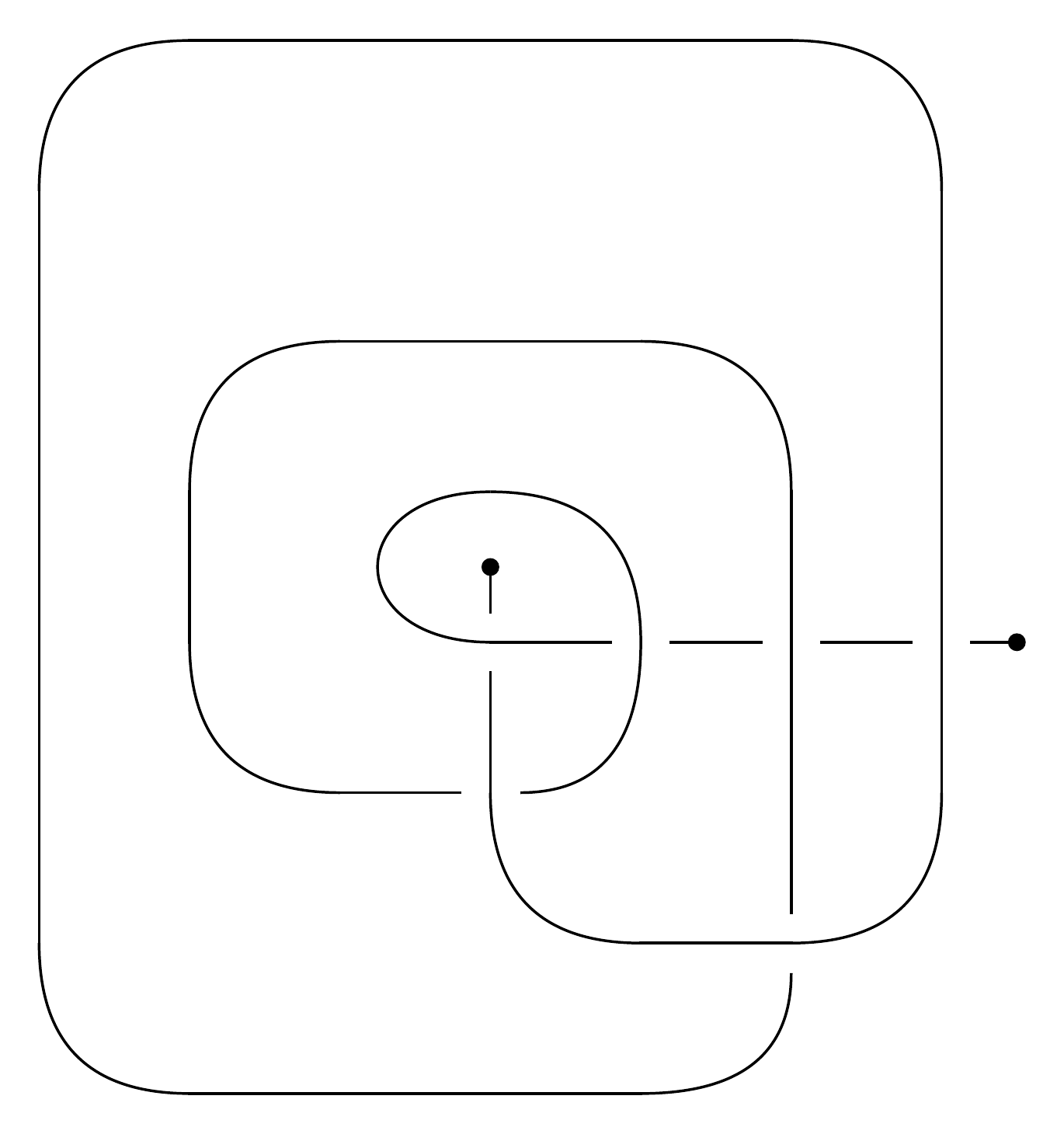}\\
\textcolor{black}{$6_{4}$}
\vspace{1cm}
\end{minipage}
\begin{minipage}[t]{.25\linewidth}
\centering
\includegraphics[width=0.9\textwidth,height=3.5cm,keepaspectratio]{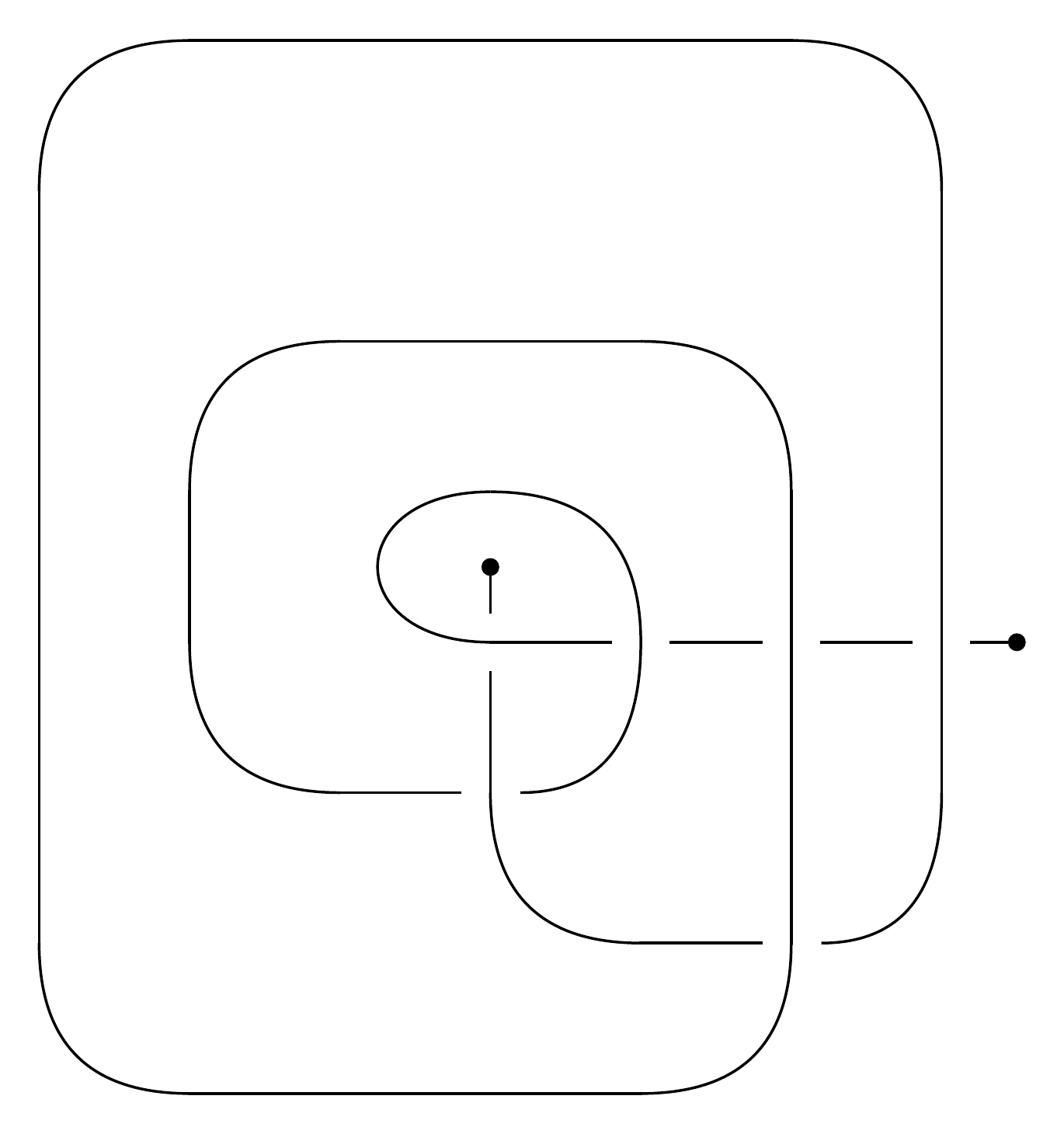}\\
\textcolor{black}{$6_{5}$}
\vspace{1cm}
\end{minipage}
\begin{minipage}[t]{.25\linewidth}
\centering
\includegraphics[width=0.9\textwidth,height=3.5cm,keepaspectratio]{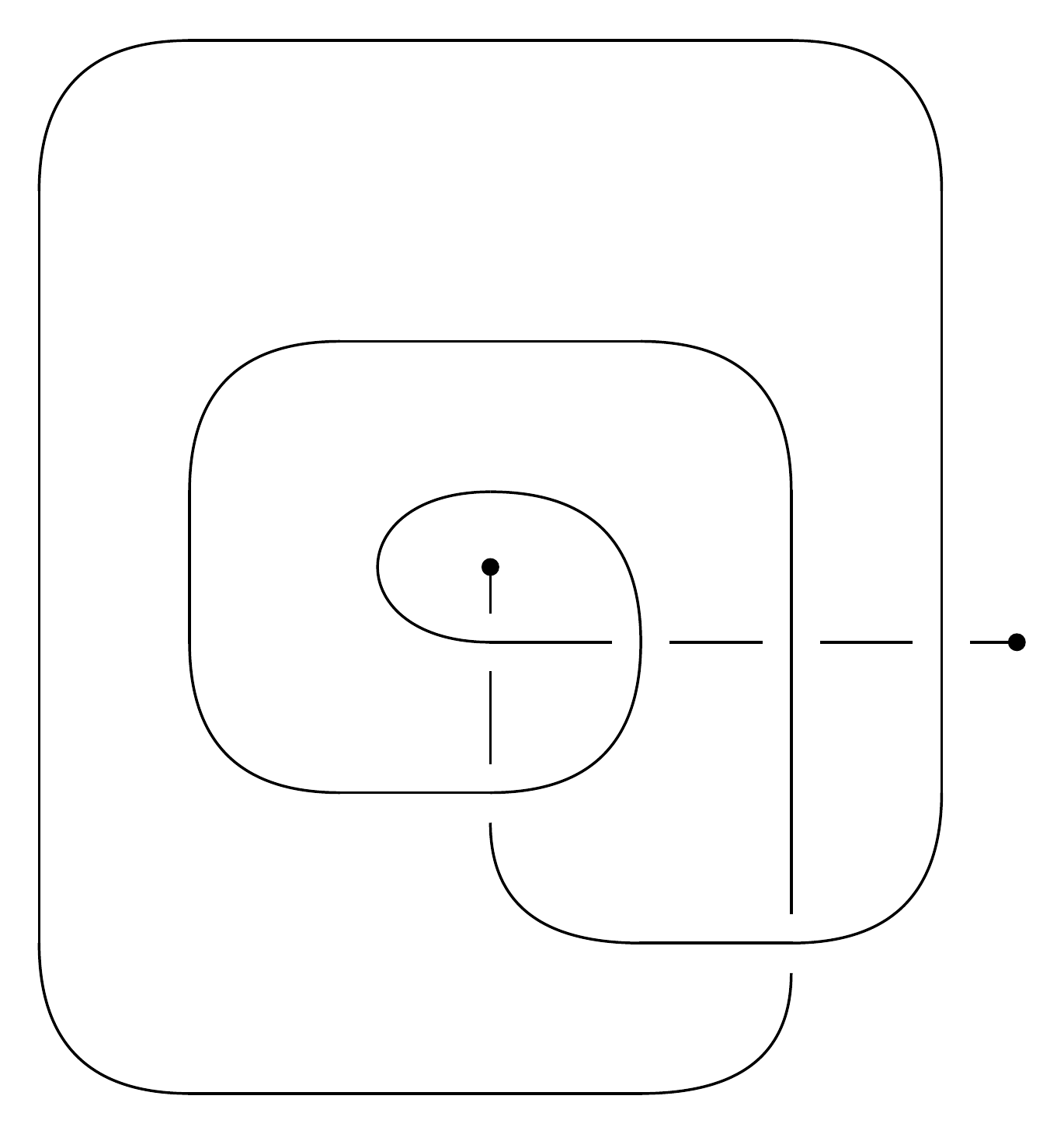}\\
\textcolor{black}{$6_{6}$}
\vspace{1cm}
\end{minipage}
\begin{minipage}[t]{.25\linewidth}
\centering
\includegraphics[width=0.9\textwidth,height=3.5cm,keepaspectratio]{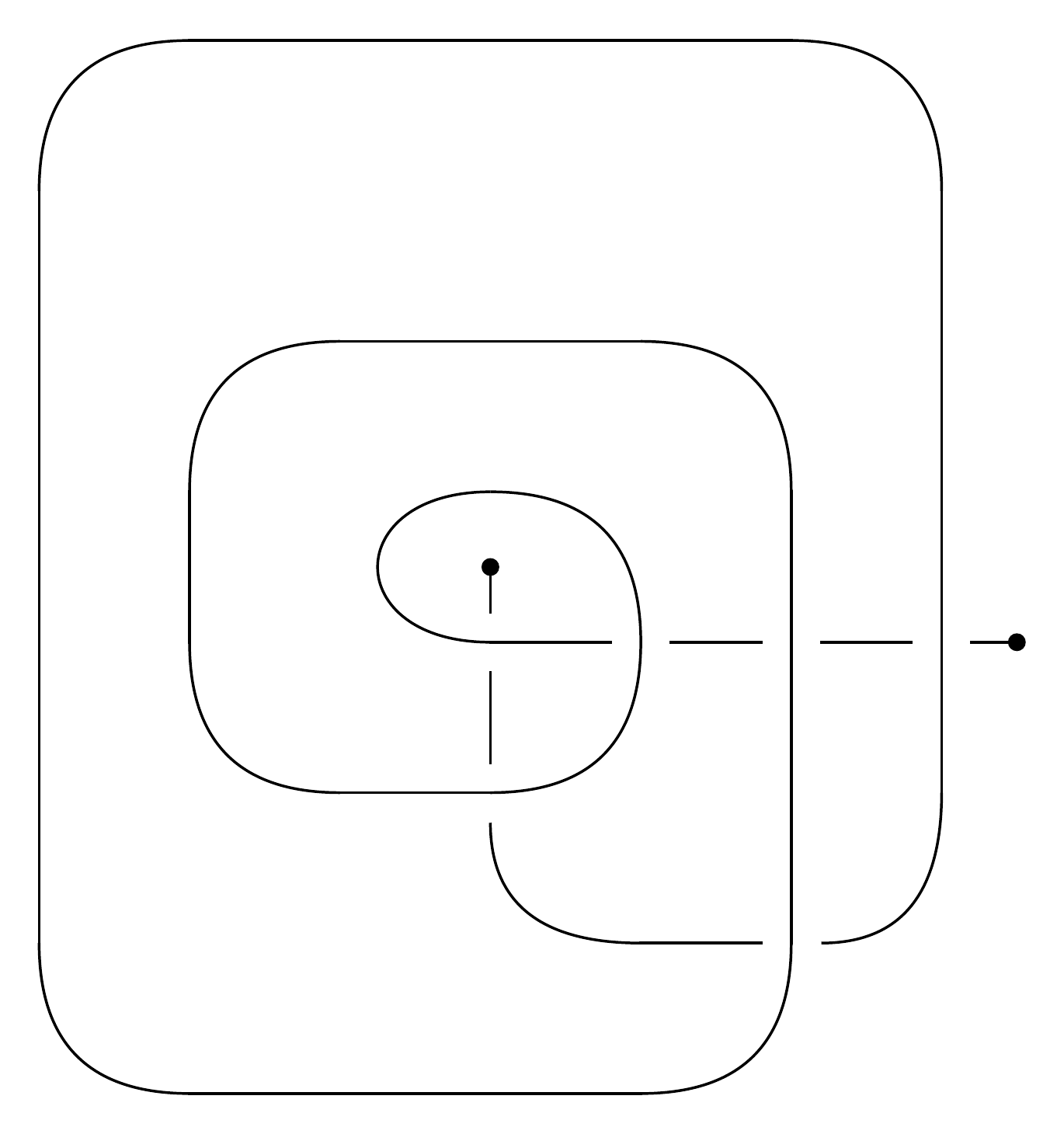}\\
\textcolor{black}{$6_{7}$}
\vspace{1cm}
\end{minipage}
\begin{minipage}[t]{.25\linewidth}
\centering
\includegraphics[width=0.9\textwidth,height=3.5cm,keepaspectratio]{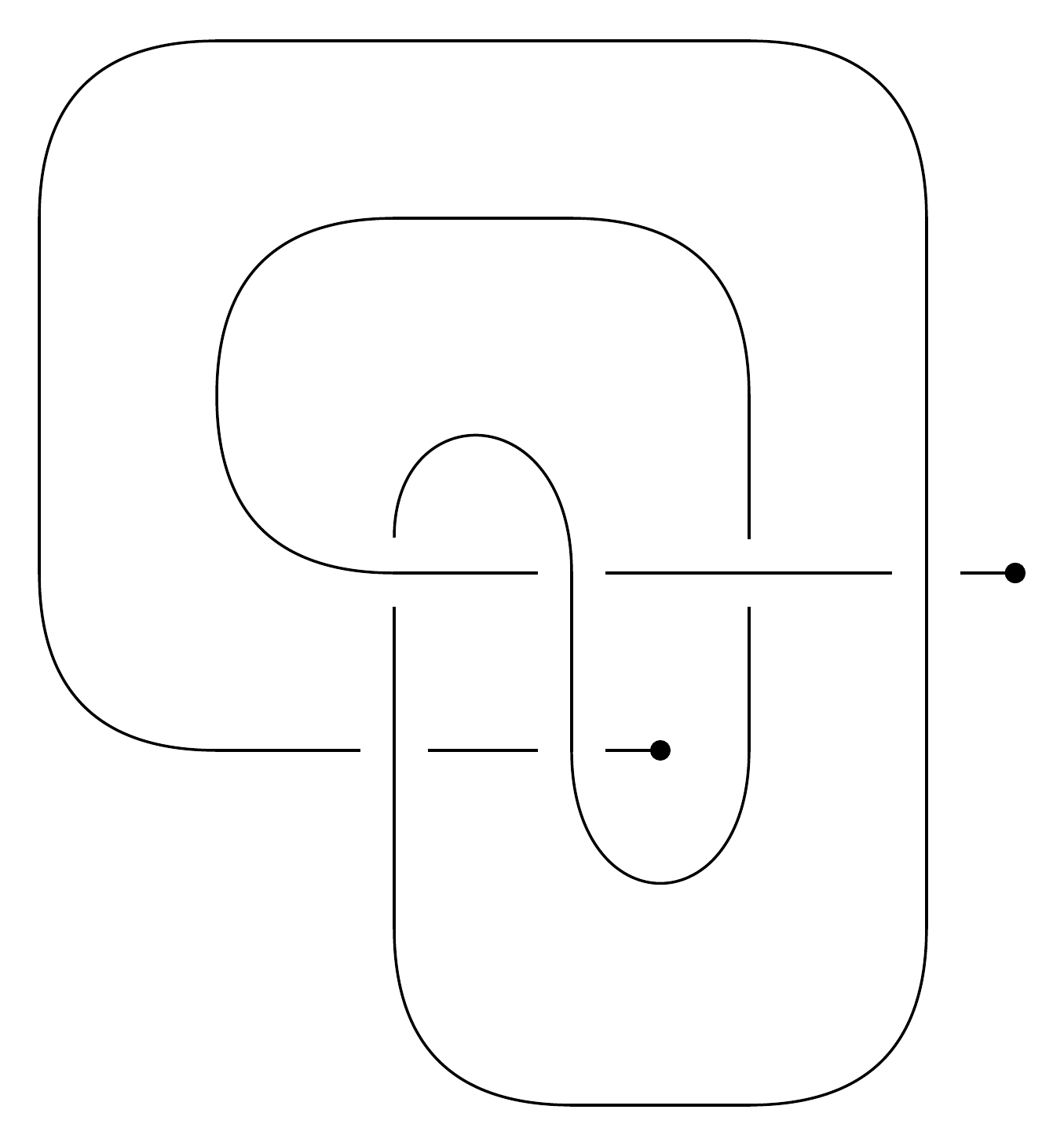}\\
\textcolor{black}{$6_{8}$}
\vspace{1cm}
\end{minipage}
\begin{minipage}[t]{.25\linewidth}
\centering
\includegraphics[width=0.9\textwidth,height=3.5cm,keepaspectratio]{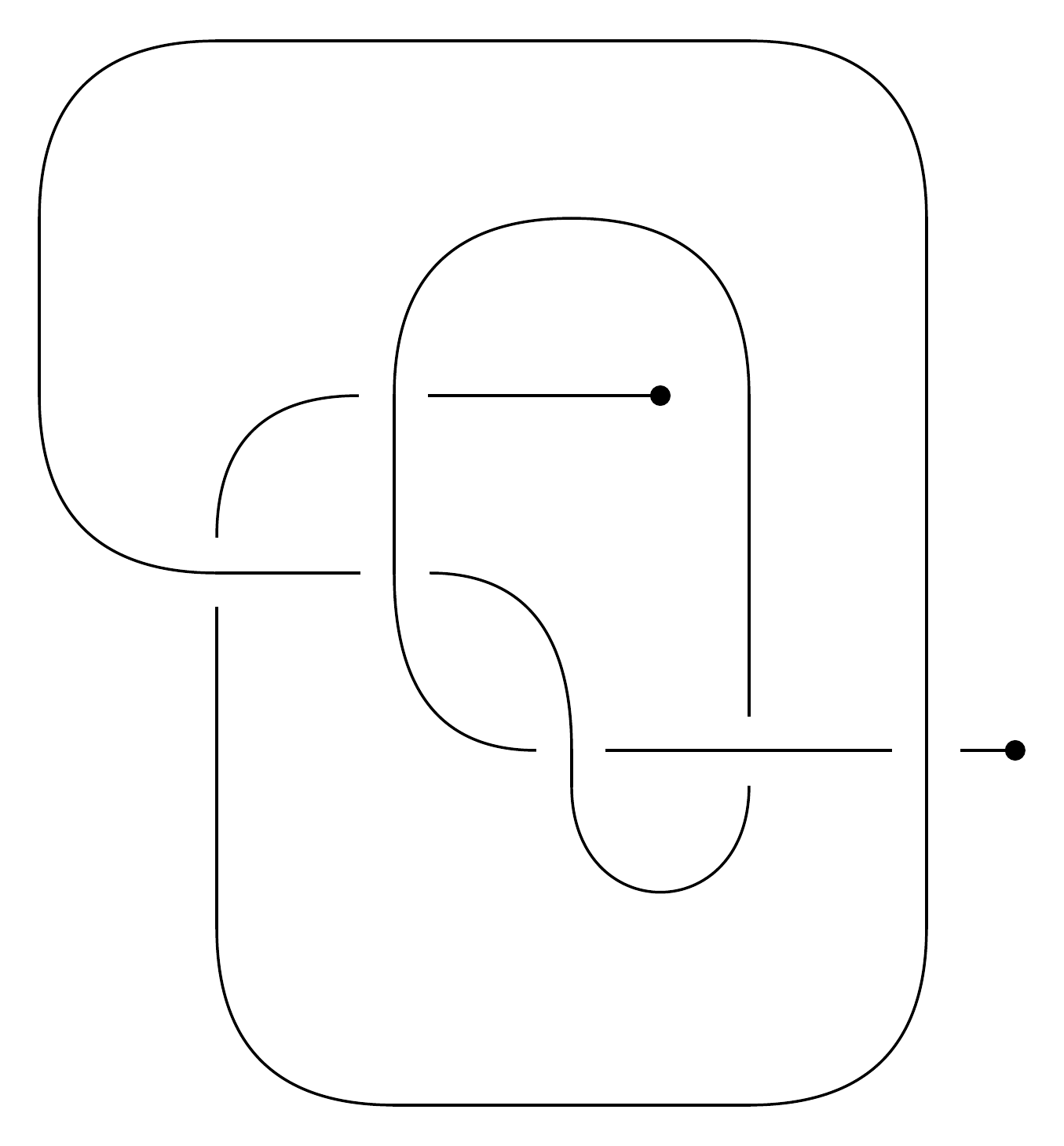}\\
\textcolor{black}{$6_{9}$}
\vspace{1cm}
\end{minipage}
\begin{minipage}[t]{.25\linewidth}
\centering
\includegraphics[width=0.9\textwidth,height=3.5cm,keepaspectratio]{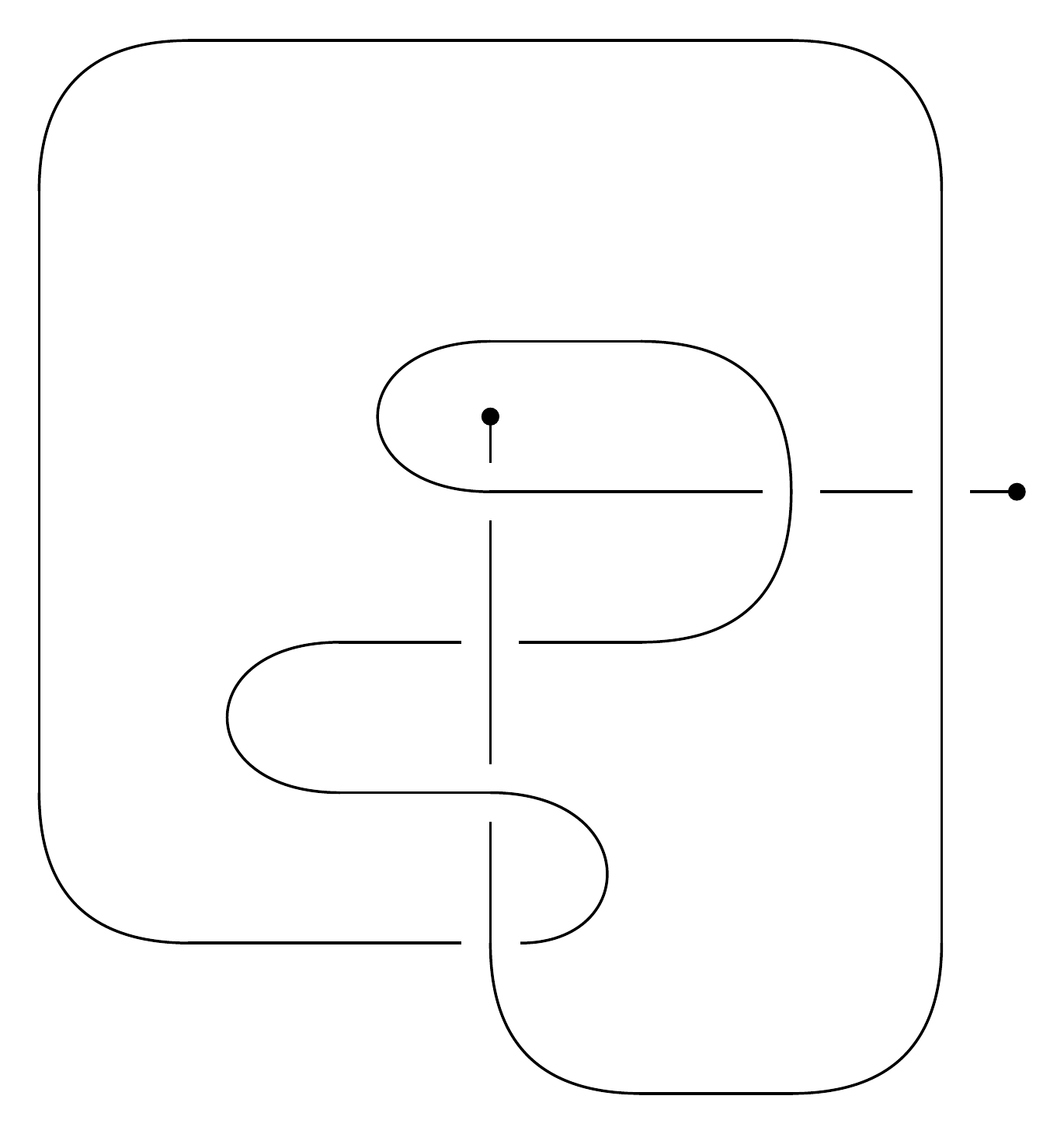}\\
\textcolor{black}{$6_{10}$}
\vspace{1cm}
\end{minipage}
\begin{minipage}[t]{.25\linewidth}
\centering
\includegraphics[width=0.9\textwidth,height=3.5cm,keepaspectratio]{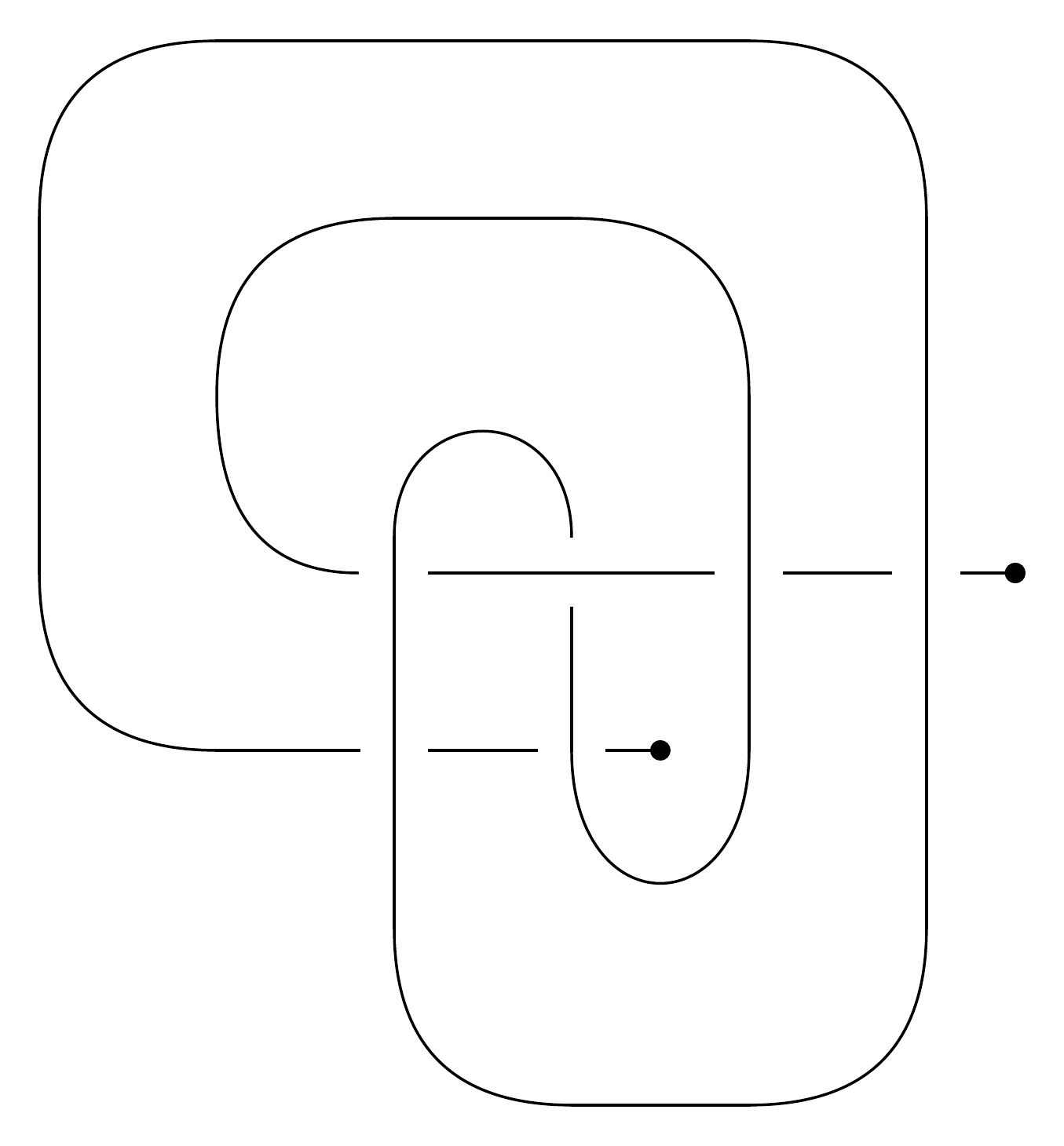}\\
\textcolor{black}{$6_{11}$}
\vspace{1cm}
\end{minipage}
\begin{minipage}[t]{.25\linewidth}
\centering
\includegraphics[width=0.9\textwidth,height=3.5cm,keepaspectratio]{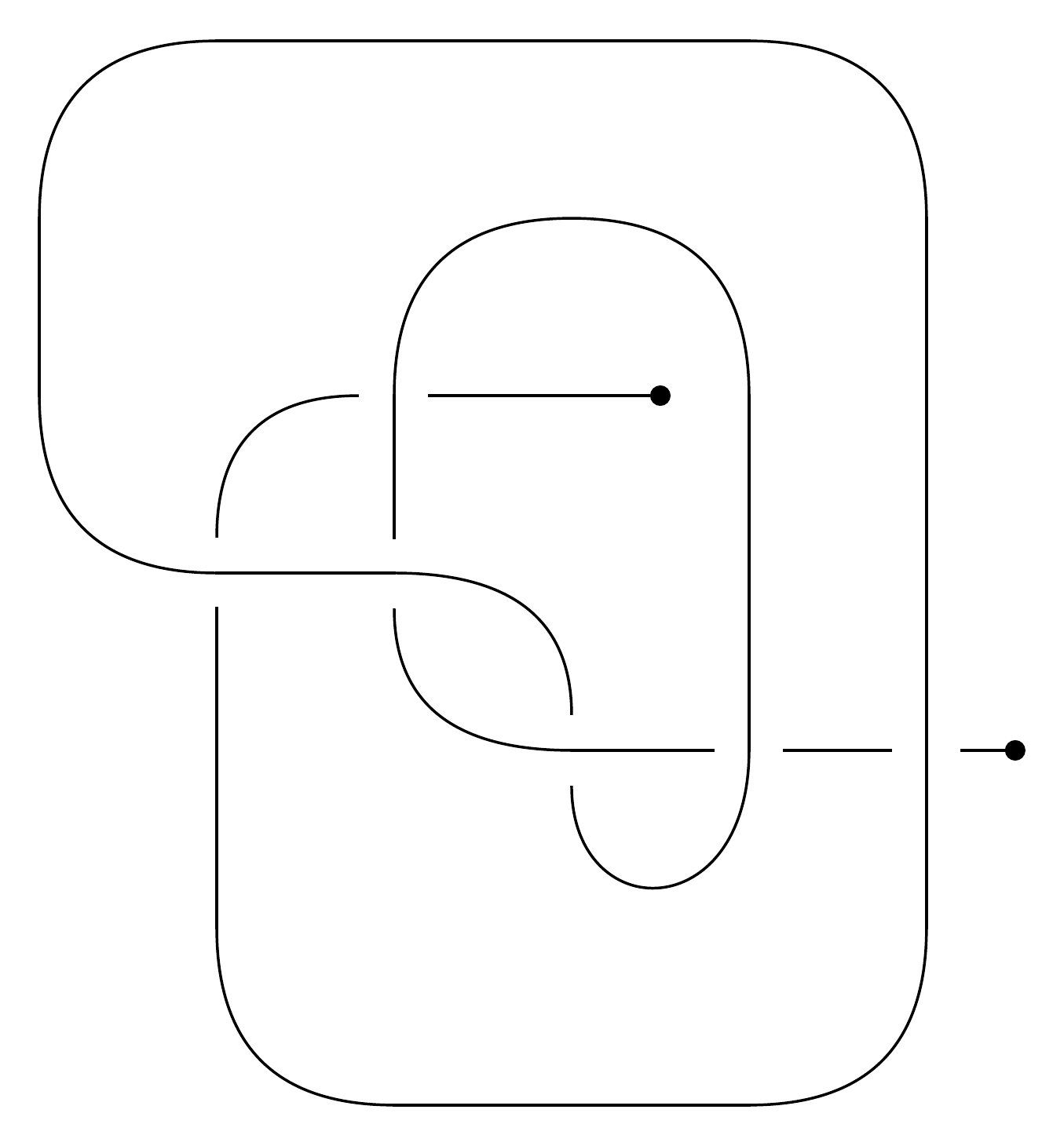}\\
\textcolor{black}{$6_{12}$}
\vspace{1cm}
\end{minipage}
\begin{minipage}[t]{.25\linewidth}
\centering
\includegraphics[width=0.9\textwidth,height=3.5cm,keepaspectratio]{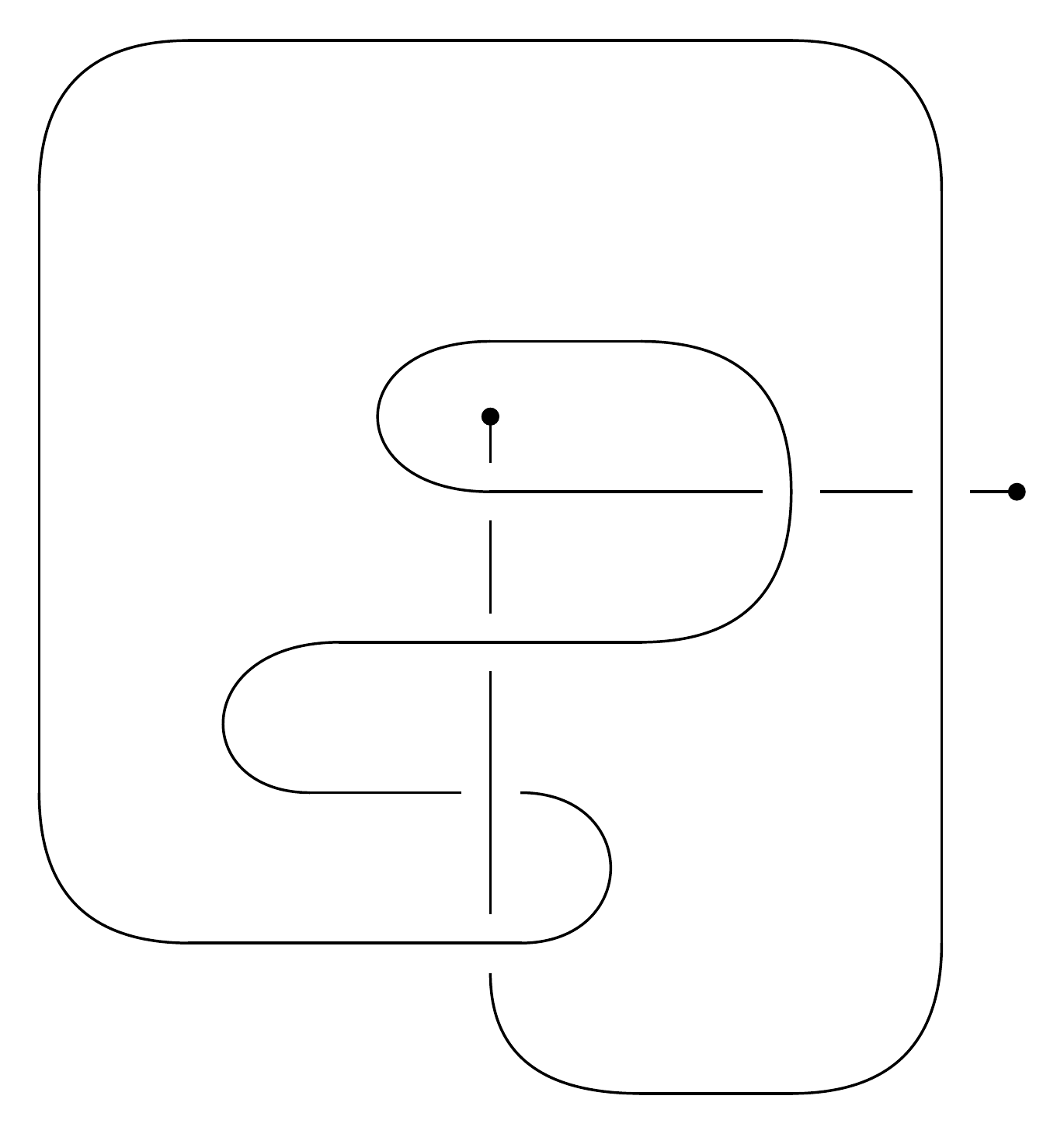}\\
\textcolor{black}{$6_{13}$}
\vspace{1cm}
\end{minipage}
\begin{minipage}[t]{.25\linewidth}
\centering
\includegraphics[width=0.9\textwidth,height=3.5cm,keepaspectratio]{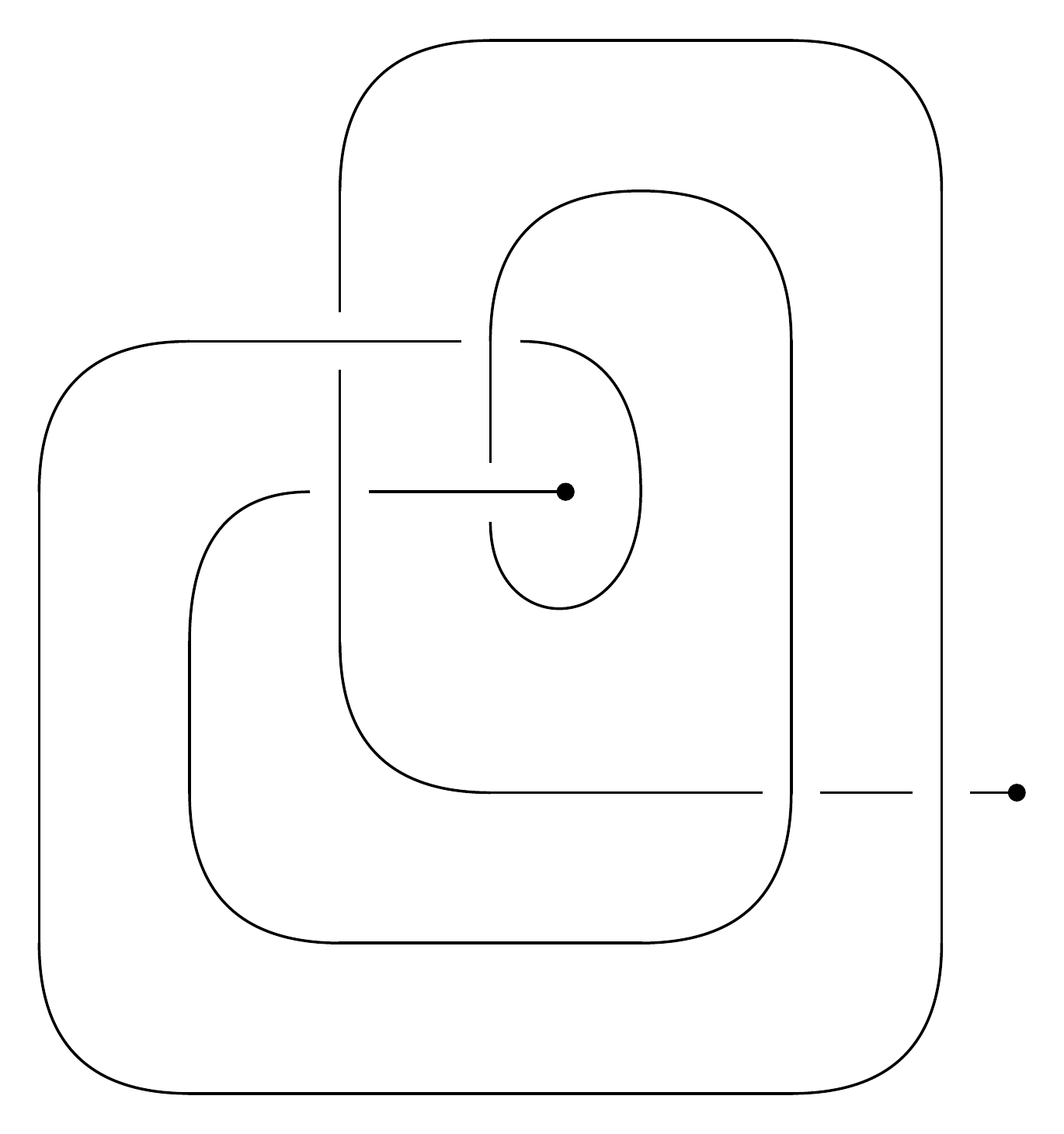}\\
\textcolor{black}{$6_{14}$}
\vspace{1cm}
\end{minipage}
\begin{minipage}[t]{.25\linewidth}
\centering
\includegraphics[width=0.9\textwidth,height=3.5cm,keepaspectratio]{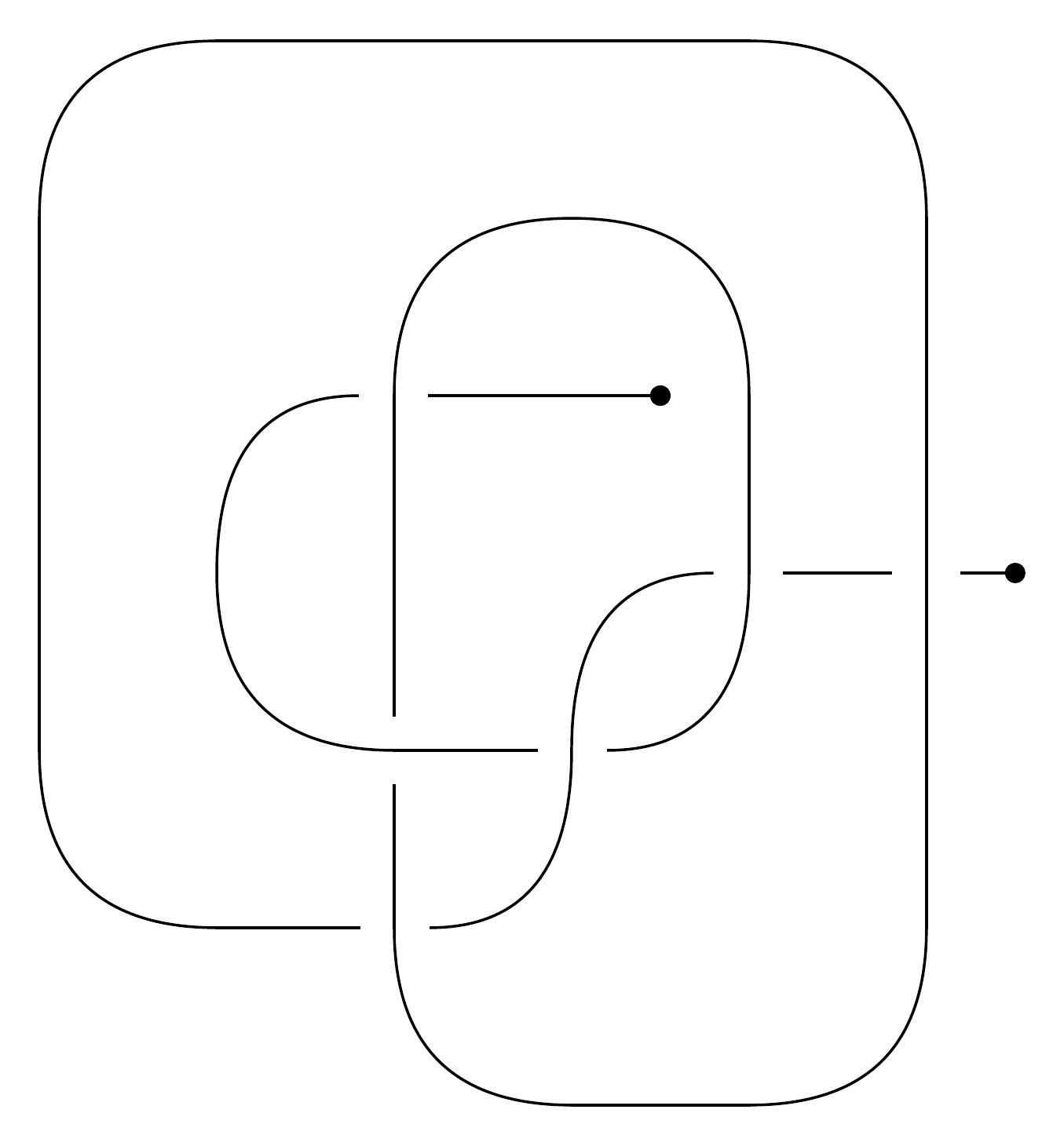}\\
\textcolor{black}{$6_{15}$}
\vspace{1cm}
\end{minipage}
\begin{minipage}[t]{.25\linewidth}
\centering
\includegraphics[width=0.9\textwidth,height=3.5cm,keepaspectratio]{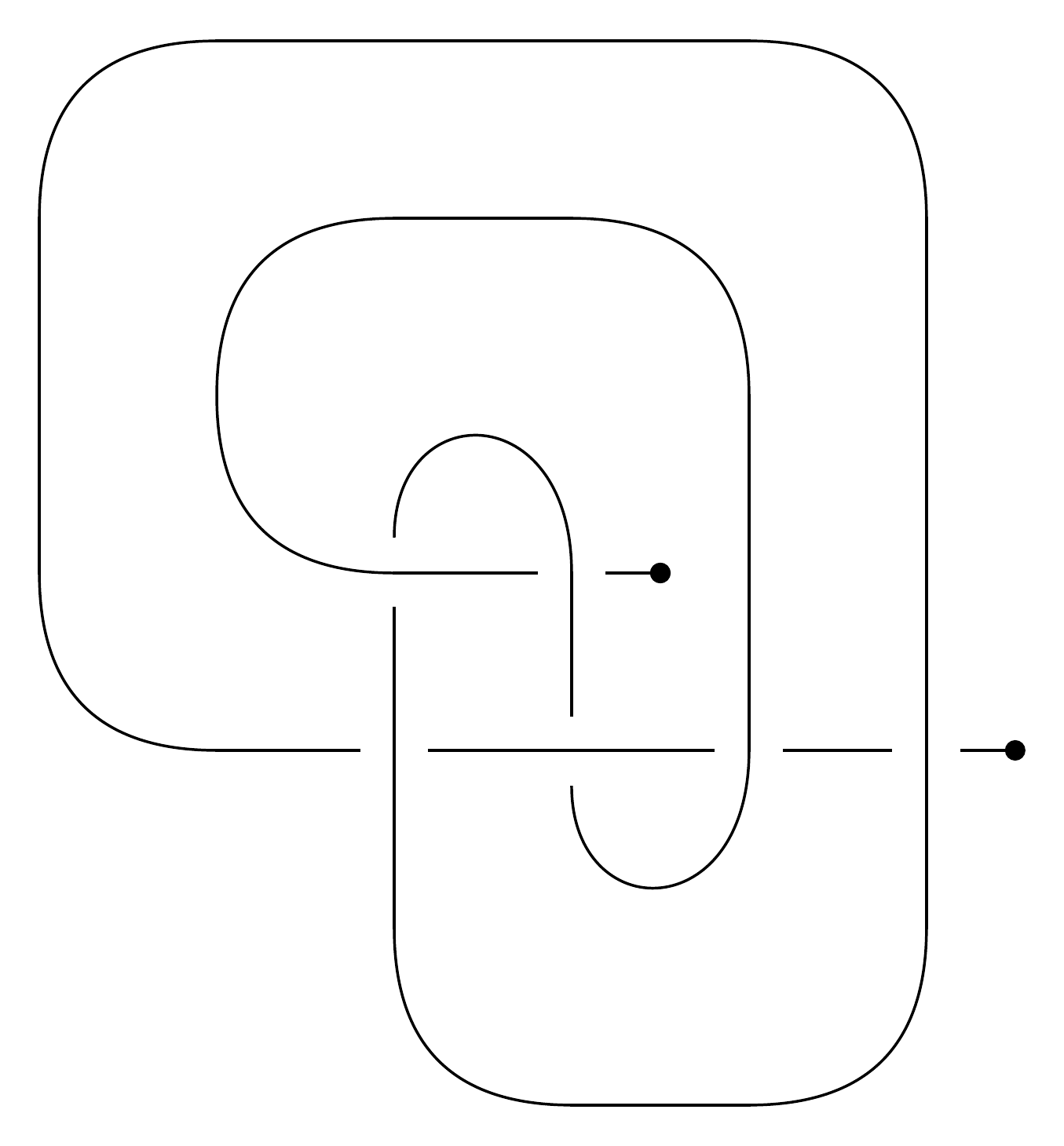}\\
\textcolor{black}{$6_{16}$}
\vspace{1cm}
\end{minipage}
\begin{minipage}[t]{.25\linewidth}
\centering
\includegraphics[width=0.9\textwidth,height=3.5cm,keepaspectratio]{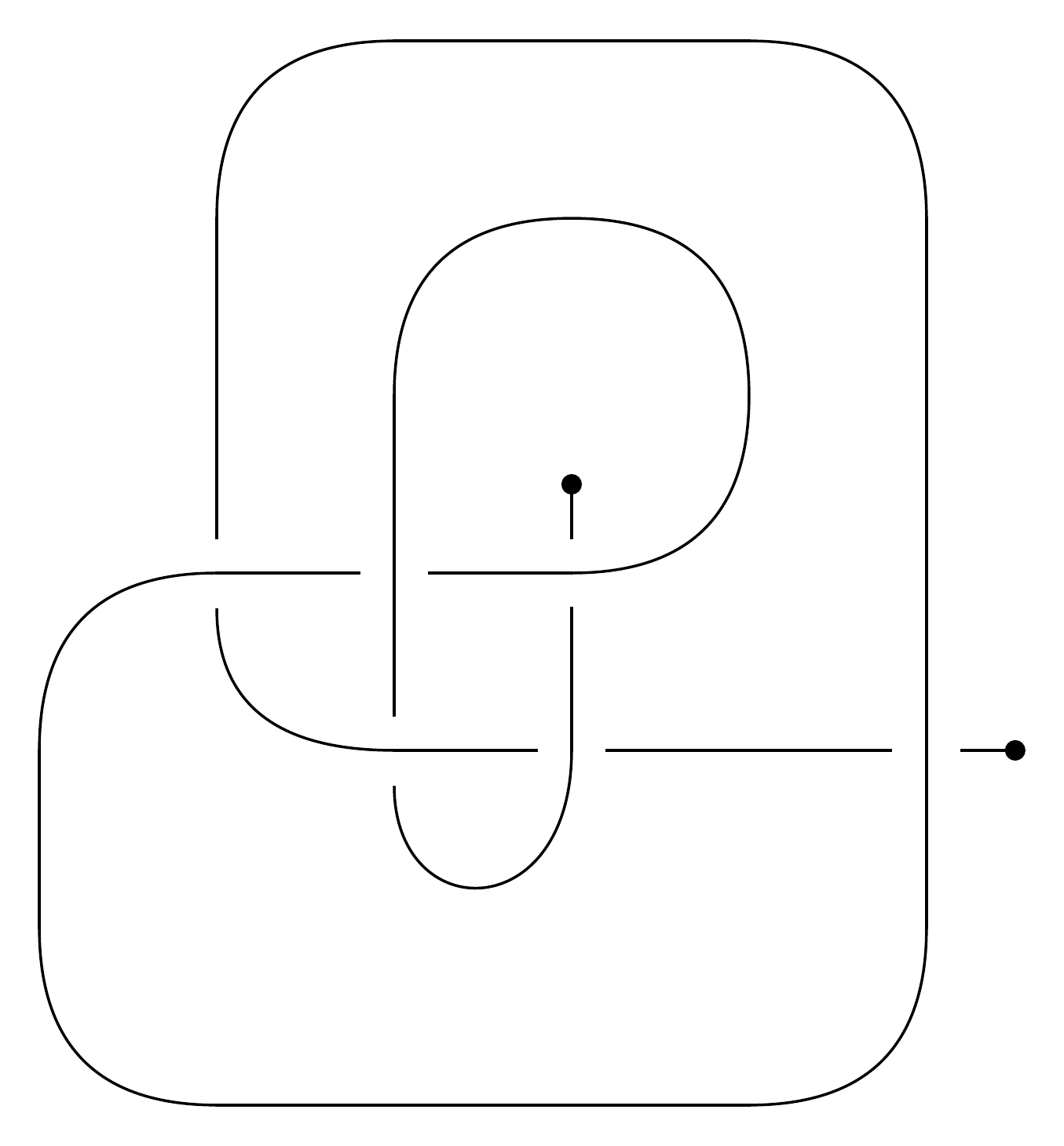}\\
\textcolor{black}{$6_{17}$}
\vspace{1cm}
\end{minipage}
\begin{minipage}[t]{.25\linewidth}
\centering
\includegraphics[width=0.9\textwidth,height=3.5cm,keepaspectratio]{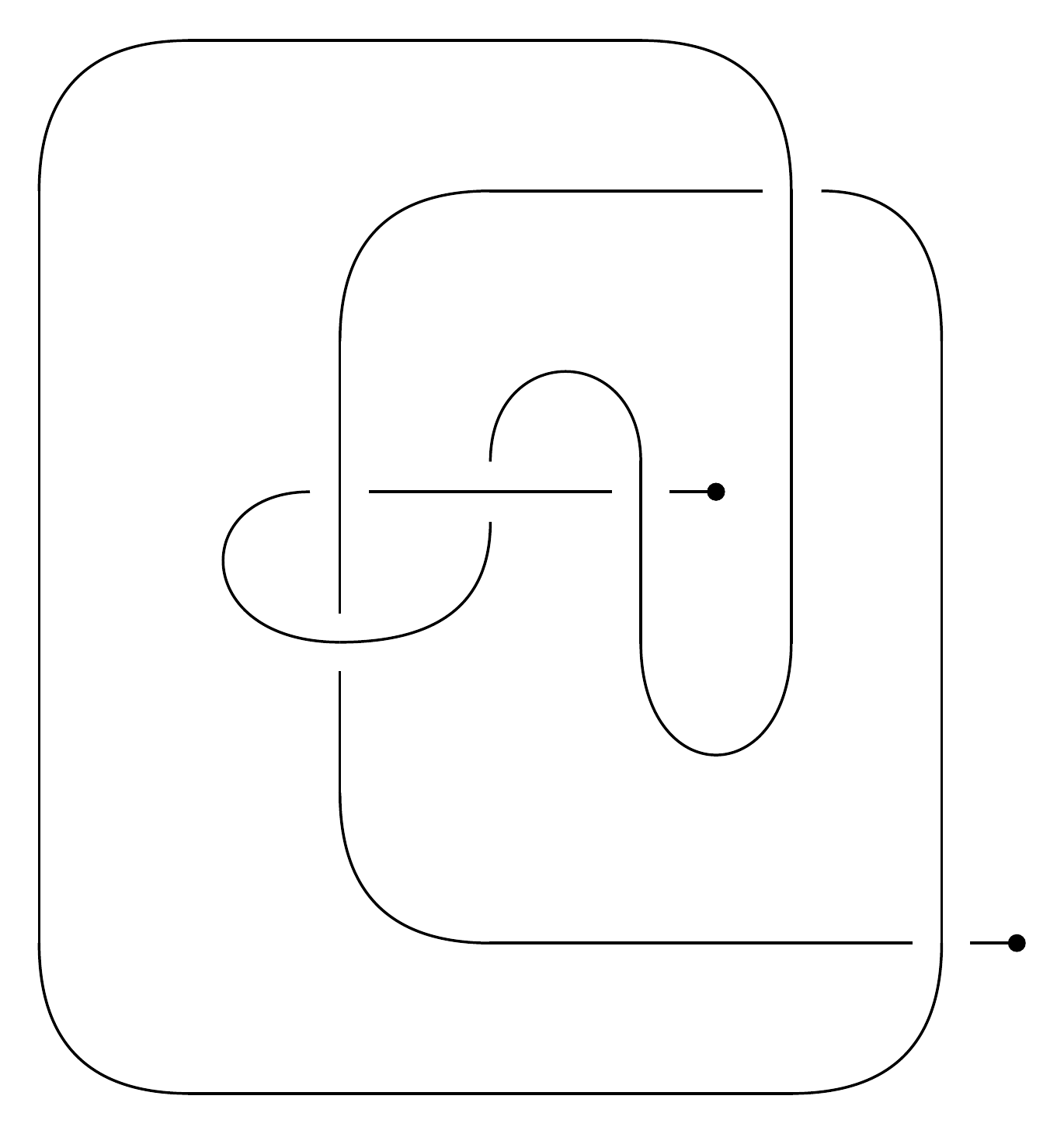}\\
\textcolor{black}{$6_{18}$}
\vspace{1cm}
\end{minipage}
\begin{minipage}[t]{.25\linewidth}
\centering
\includegraphics[width=0.9\textwidth,height=3.5cm,keepaspectratio]{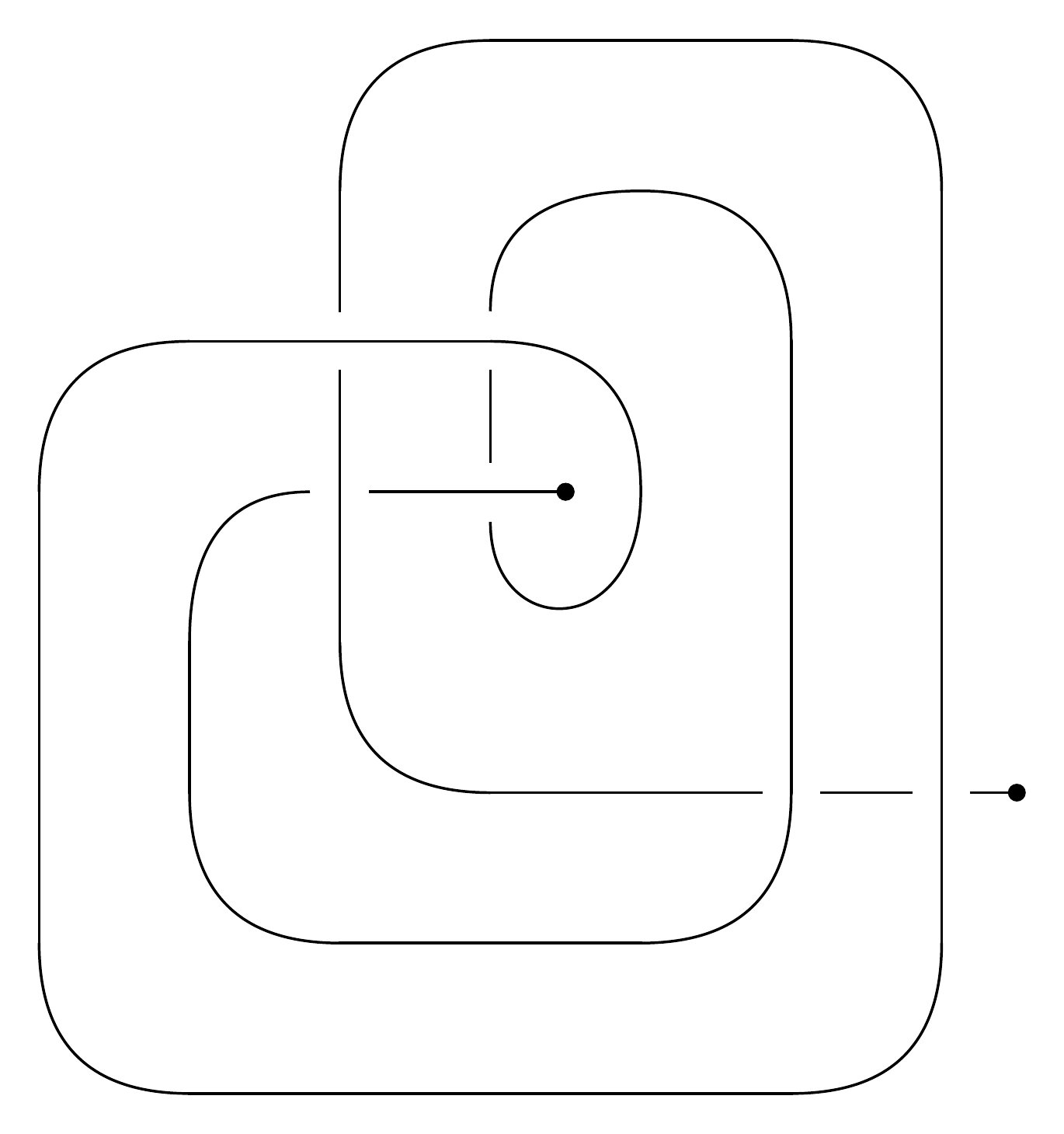}\\
\textcolor{black}{$6_{19}$}
\vspace{1cm}
\end{minipage}
\begin{minipage}[t]{.25\linewidth}
\centering
\includegraphics[width=0.9\textwidth,height=3.5cm,keepaspectratio]{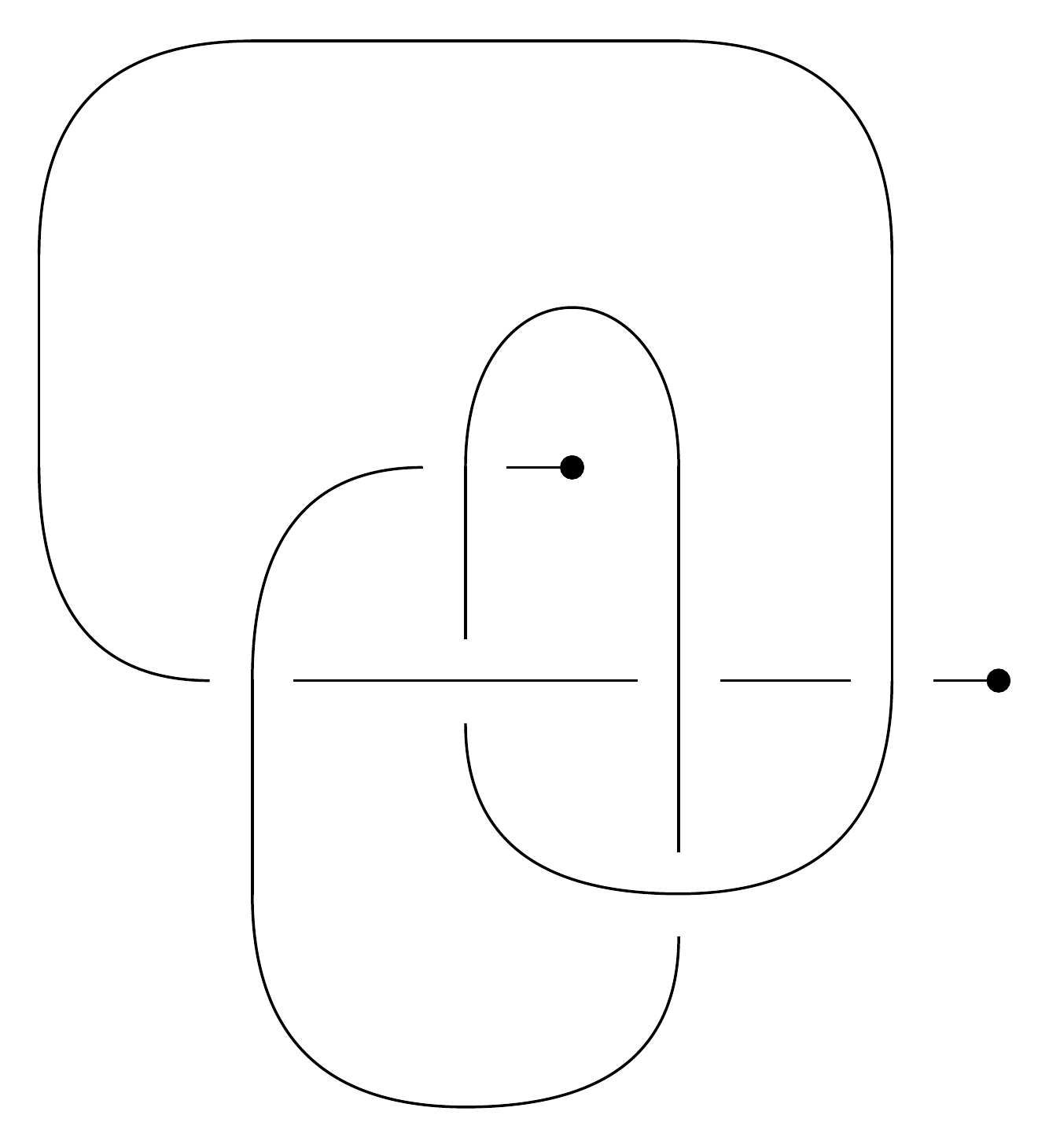}\\
\textcolor{black}{$6_{20}$}
\vspace{1cm}
\end{minipage}
\begin{minipage}[t]{.25\linewidth}
\centering
\includegraphics[width=0.9\textwidth,height=3.5cm,keepaspectratio]{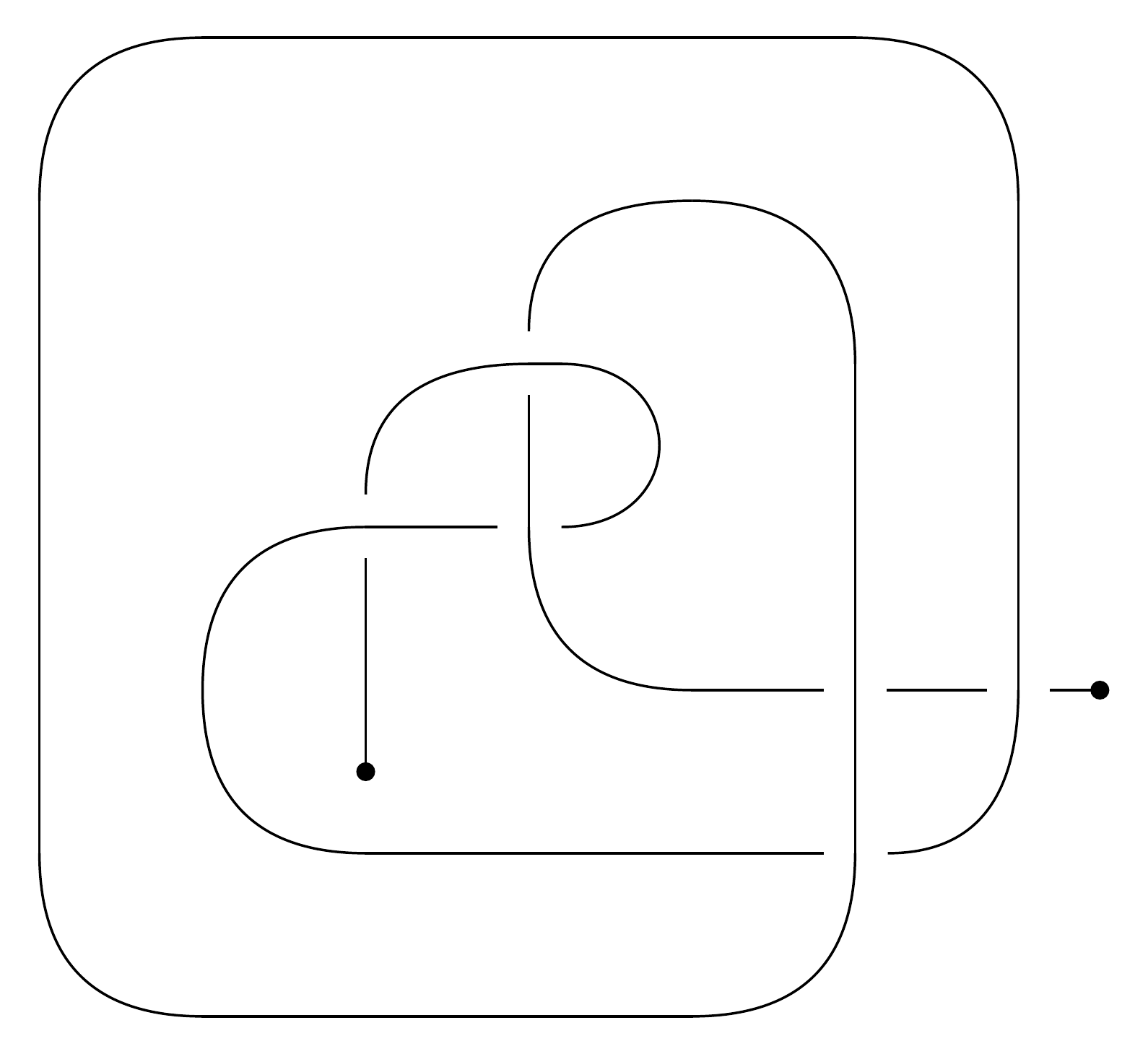}\\
\textcolor{black}{$6_{21}$}
\vspace{1cm}
\end{minipage}
\begin{minipage}[t]{.25\linewidth}
\centering
\includegraphics[width=0.9\textwidth,height=3.5cm,keepaspectratio]{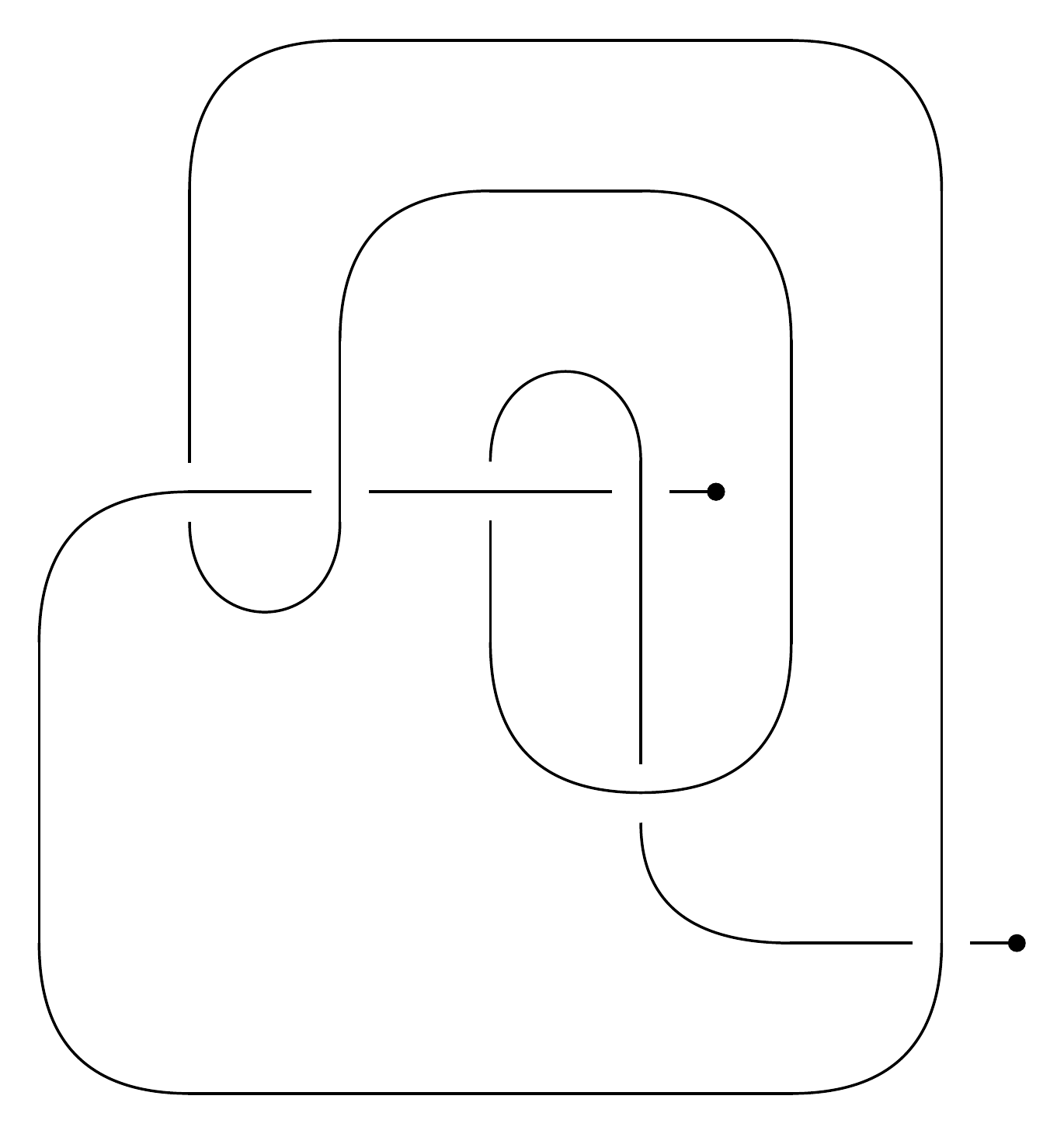}\\
\textcolor{black}{$6_{22}$}
\vspace{1cm}
\end{minipage}
\begin{minipage}[t]{.25\linewidth}
\centering
\includegraphics[width=0.9\textwidth,height=3.5cm,keepaspectratio]{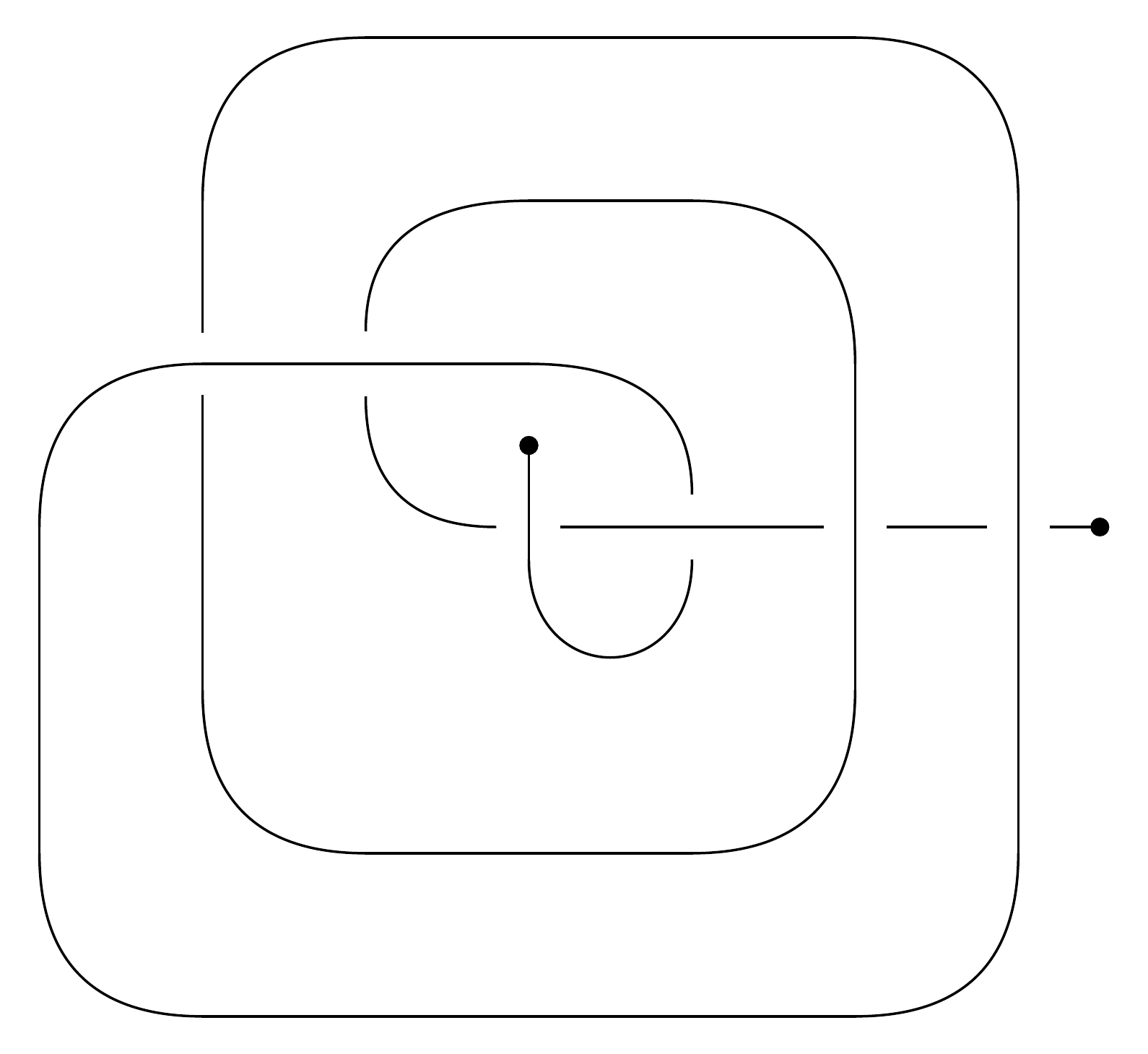}\\
\textcolor{black}{$6_{23}$}
\vspace{1cm}
\end{minipage}
\begin{minipage}[t]{.25\linewidth}
\centering
\includegraphics[width=0.9\textwidth,height=3.5cm,keepaspectratio]{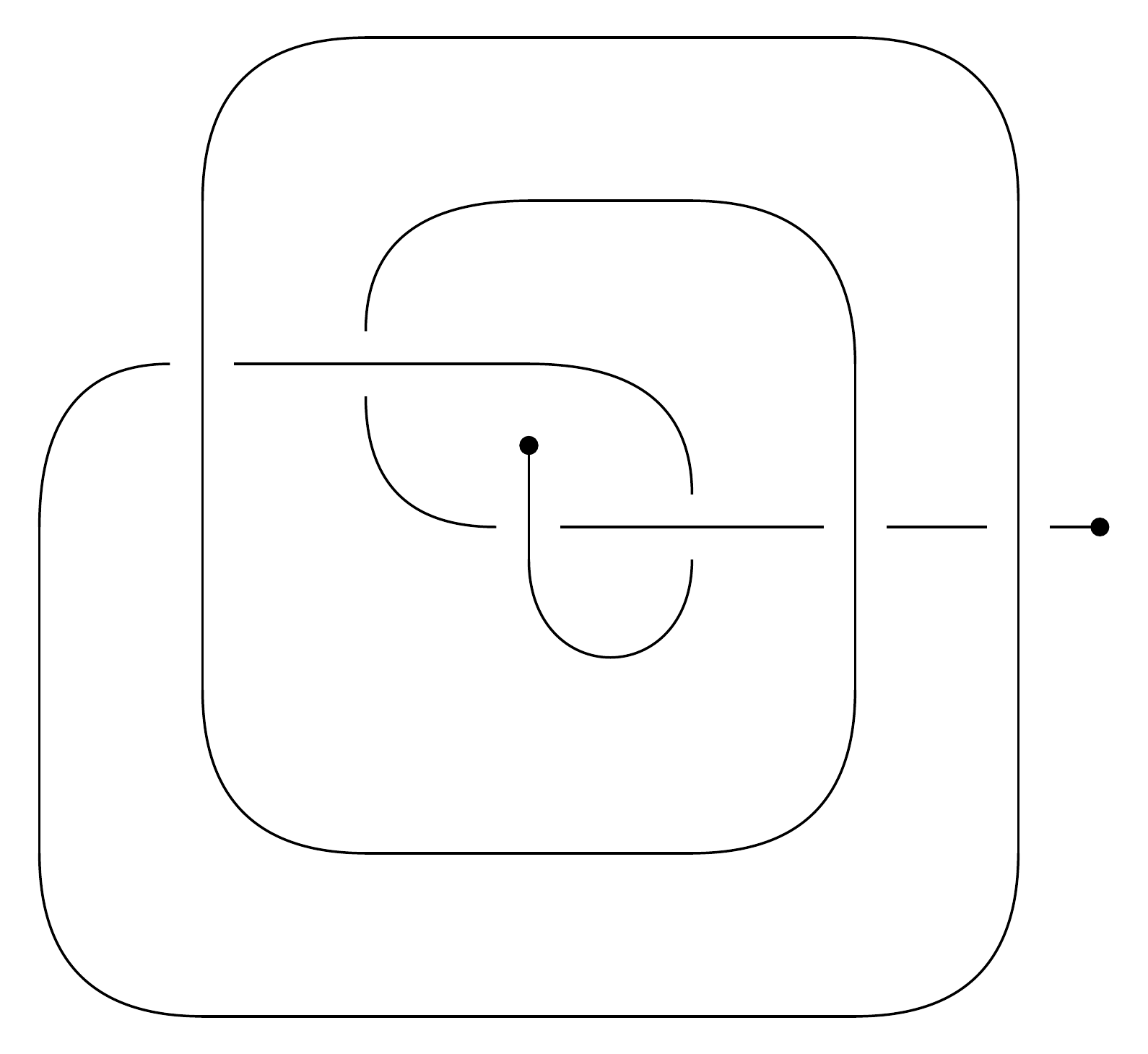}\\
\textcolor{black}{$6_{24}$}
\vspace{1cm}
\end{minipage}
\begin{minipage}[t]{.25\linewidth}
\centering
\includegraphics[width=0.9\textwidth,height=3.5cm,keepaspectratio]{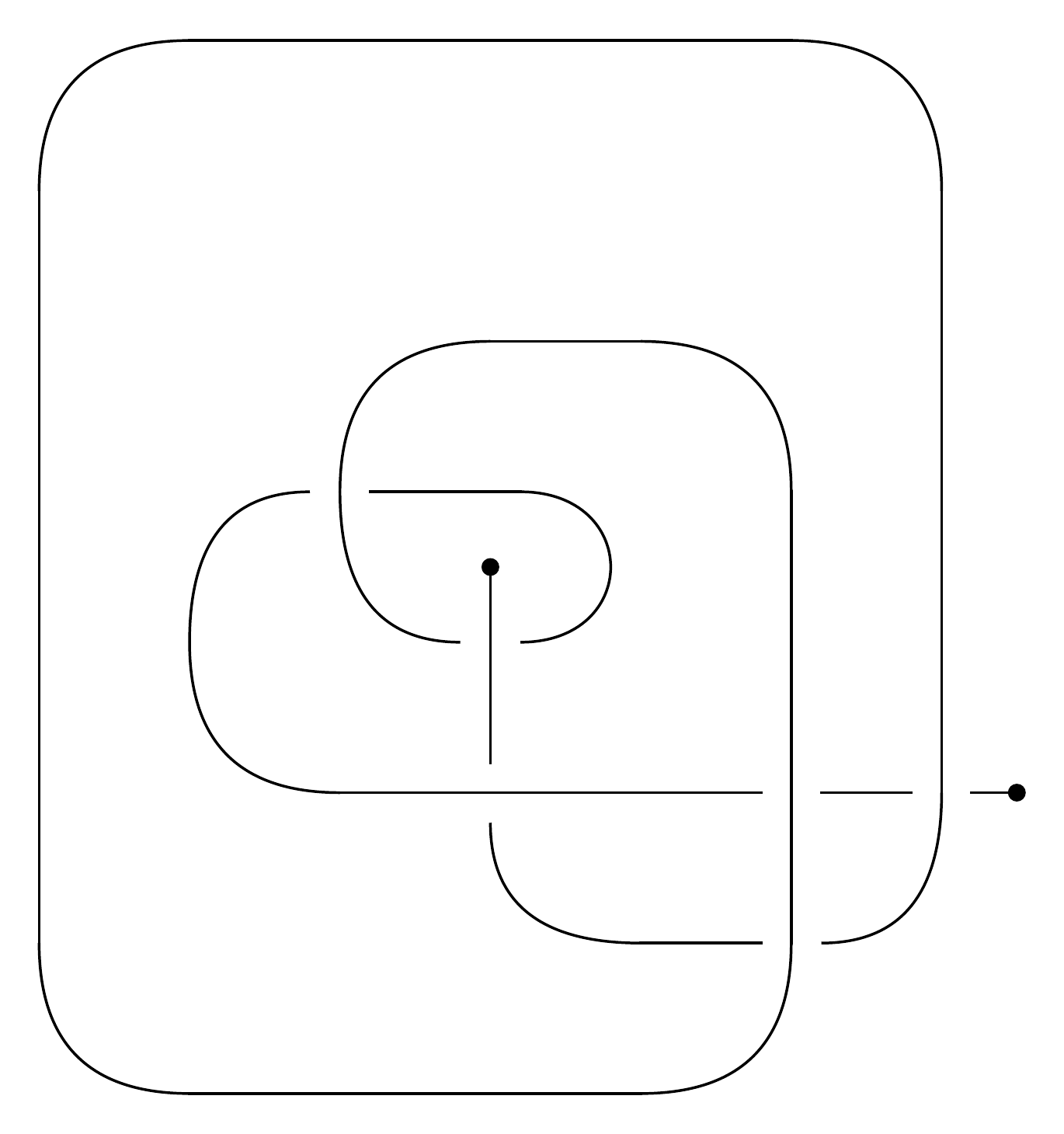}\\
\textcolor{black}{$6_{25}$}
\vspace{1cm}
\end{minipage}
\begin{minipage}[t]{.25\linewidth}
\centering
\includegraphics[width=0.9\textwidth,height=3.5cm,keepaspectratio]{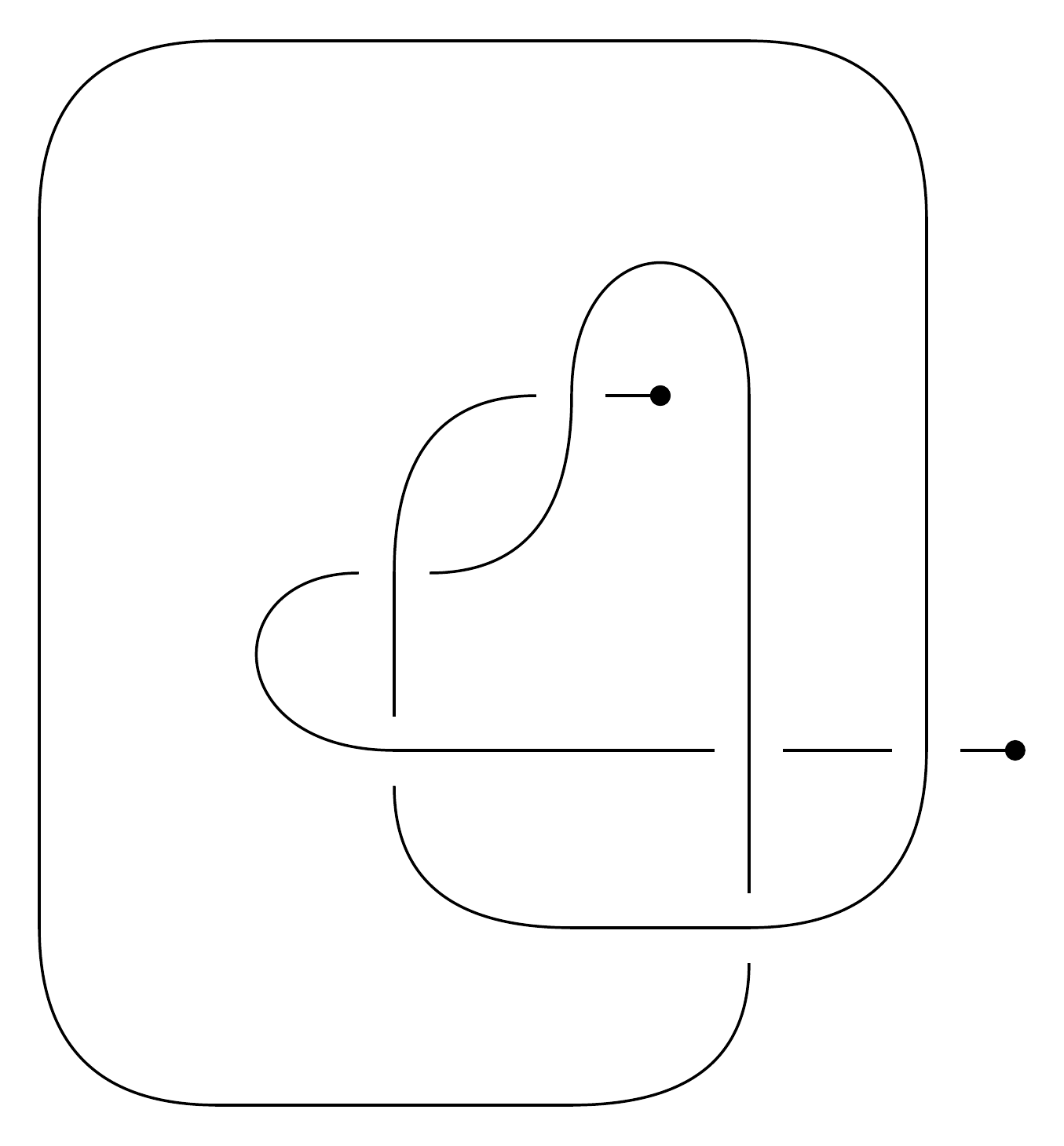}\\
\textcolor{black}{$6_{26}$}
\vspace{1cm}
\end{minipage}
\begin{minipage}[t]{.25\linewidth}
\centering
\includegraphics[width=0.9\textwidth,height=3.5cm,keepaspectratio]{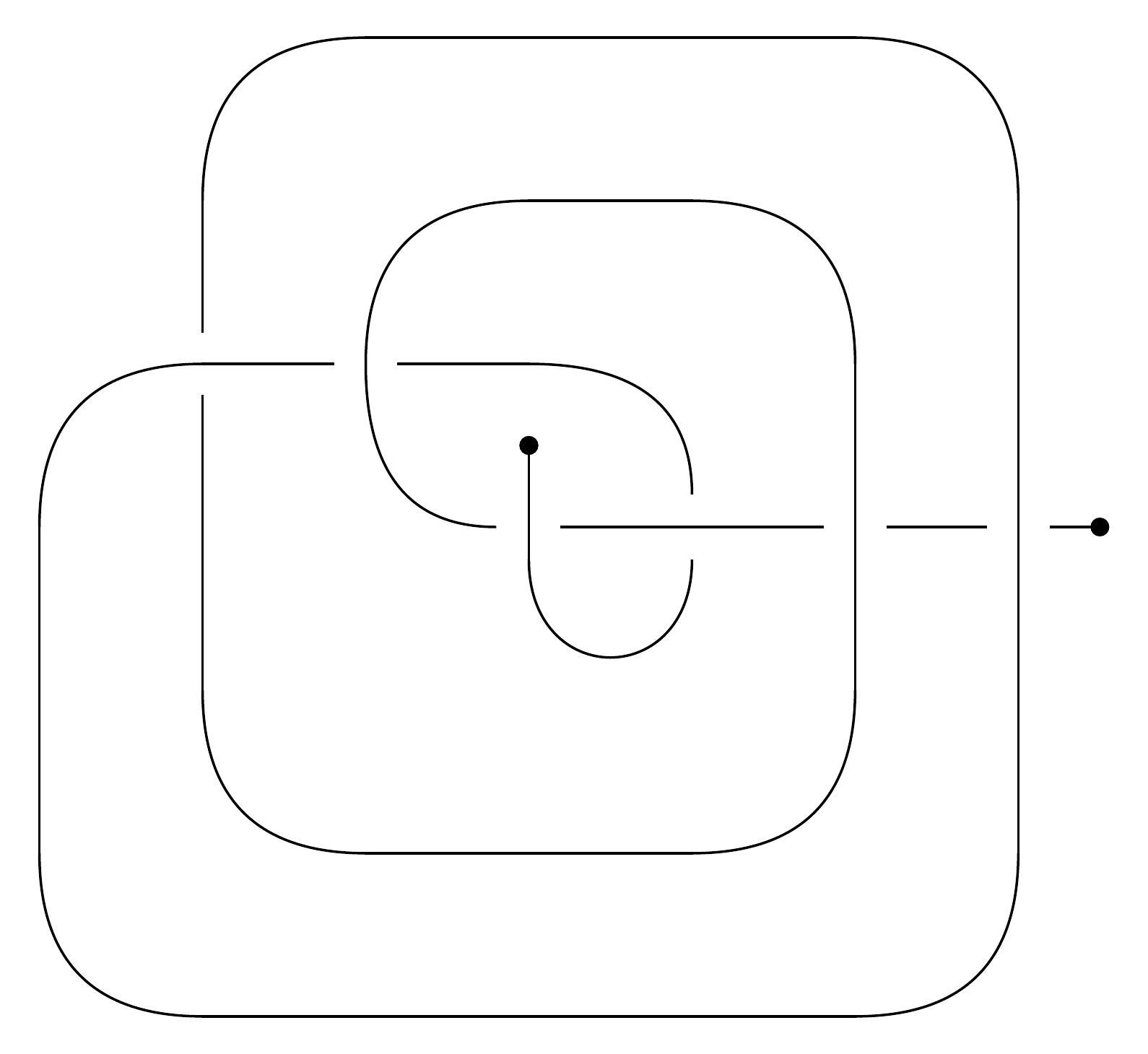}\\
\textcolor{black}{$6_{27}$}
\vspace{1cm}
\end{minipage}
\begin{minipage}[t]{.25\linewidth}
\centering
\includegraphics[width=0.9\textwidth,height=3.5cm,keepaspectratio]{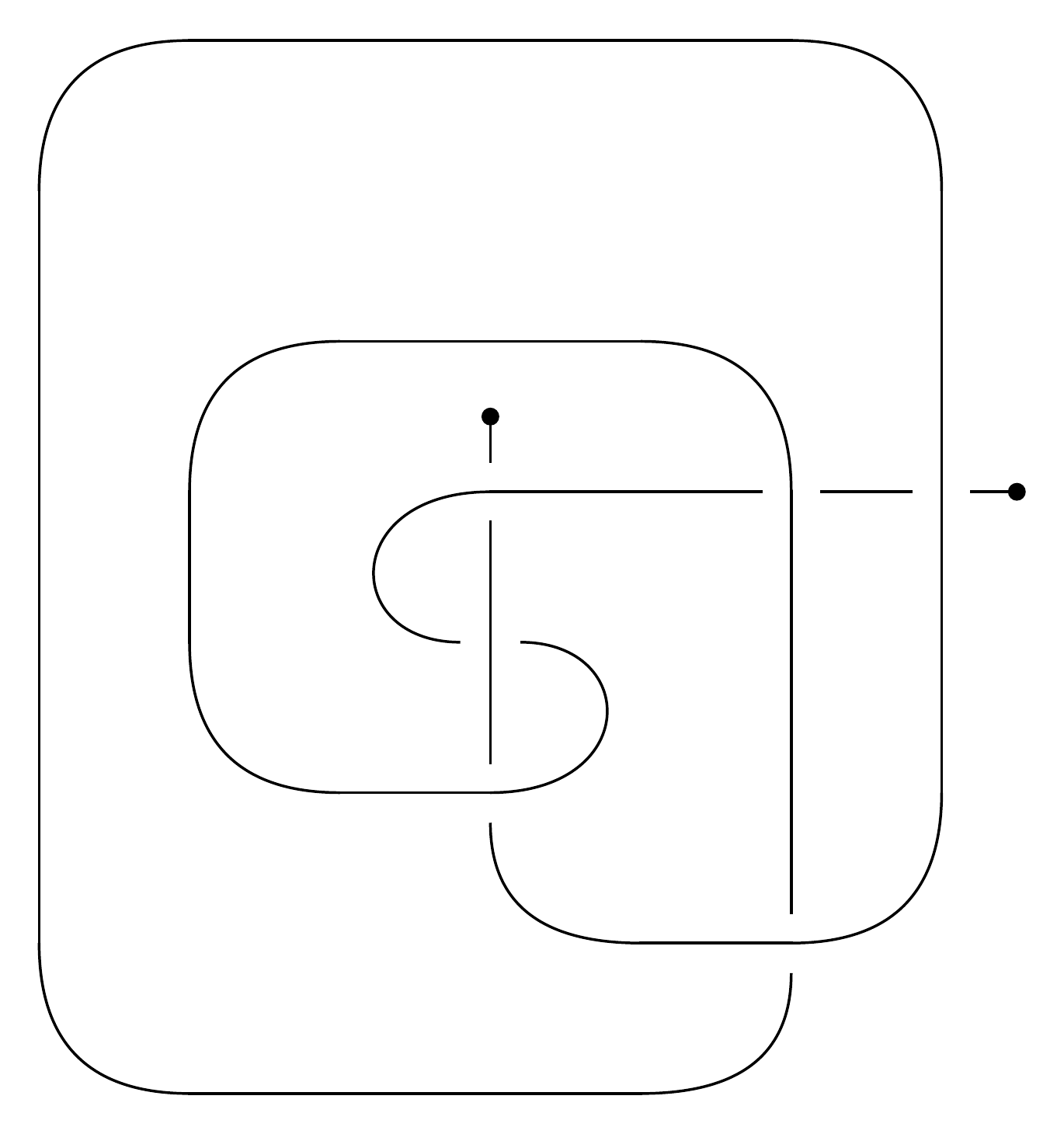}\\
\textcolor{black}{$6_{28}$}
\vspace{1cm}
\end{minipage}
\begin{minipage}[t]{.25\linewidth}
\centering
\includegraphics[width=0.9\textwidth,height=3.5cm,keepaspectratio]{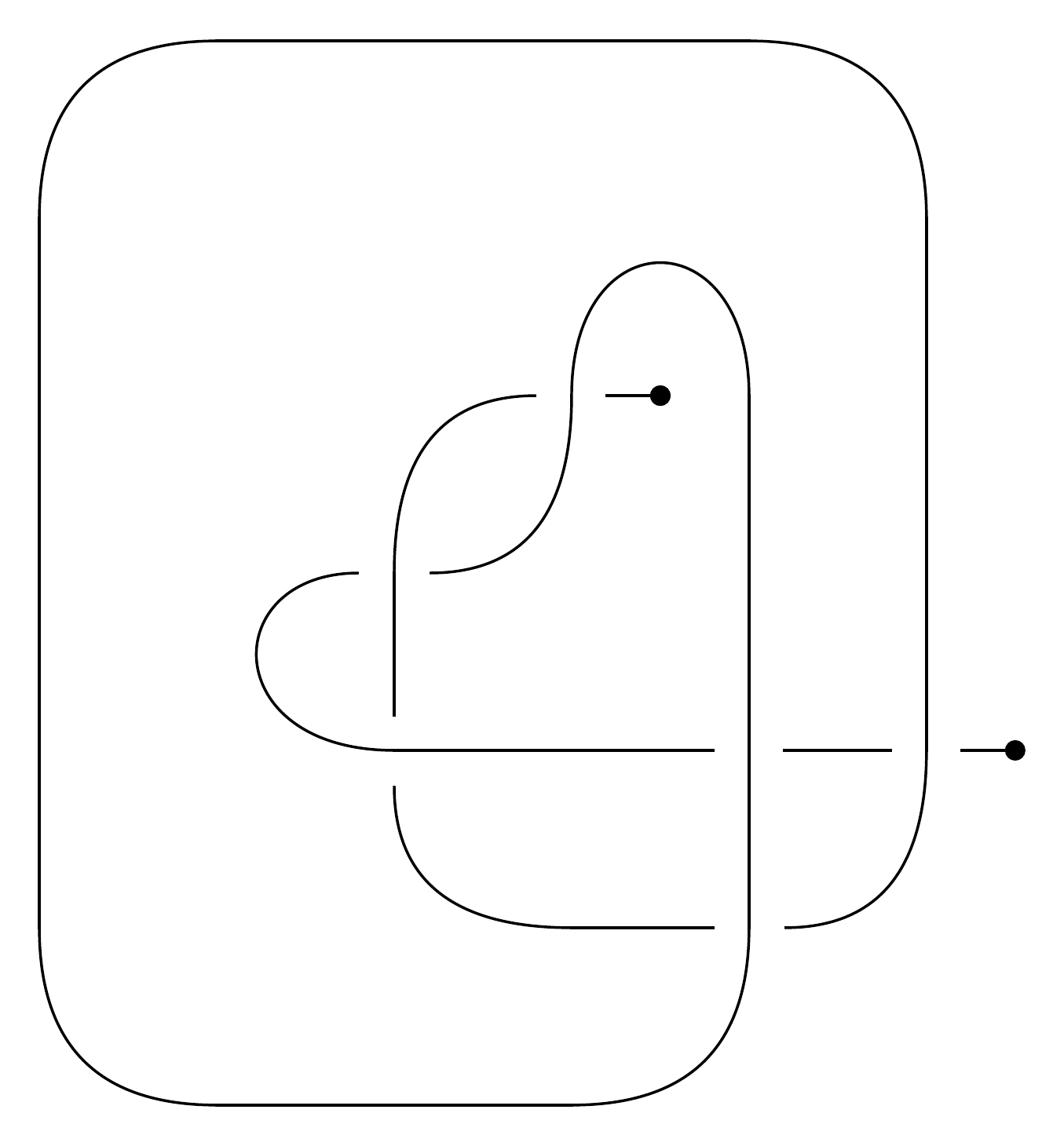}\\
\textcolor{black}{$6_{29}$}
\vspace{1cm}
\end{minipage}
\begin{minipage}[t]{.25\linewidth}
\centering
\includegraphics[width=0.9\textwidth,height=3.5cm,keepaspectratio]{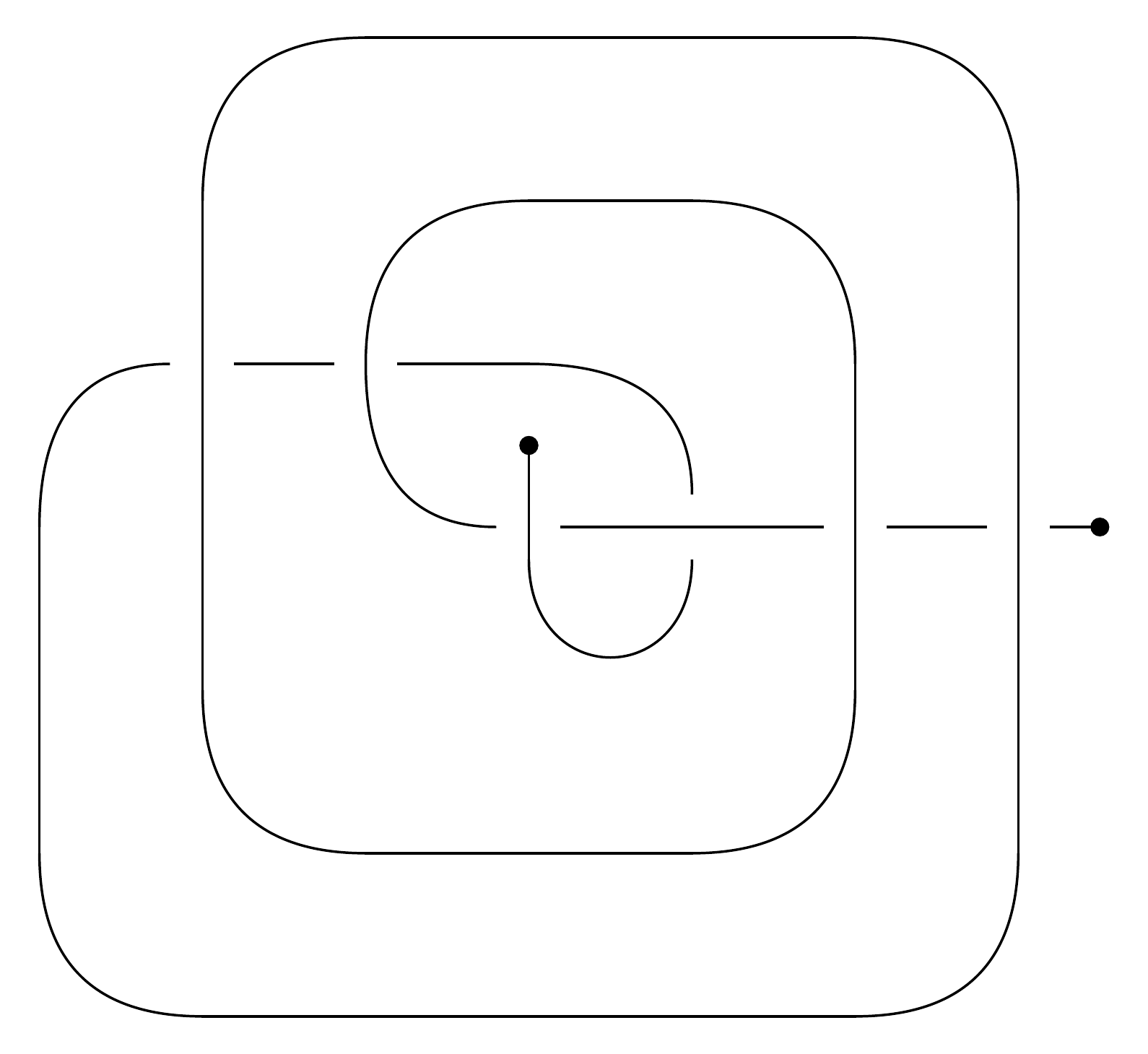}\\
\textcolor{black}{$6_{30}$}
\vspace{1cm}
\end{minipage}
\begin{minipage}[t]{.25\linewidth}
\centering
\includegraphics[width=0.9\textwidth,height=3.5cm,keepaspectratio]{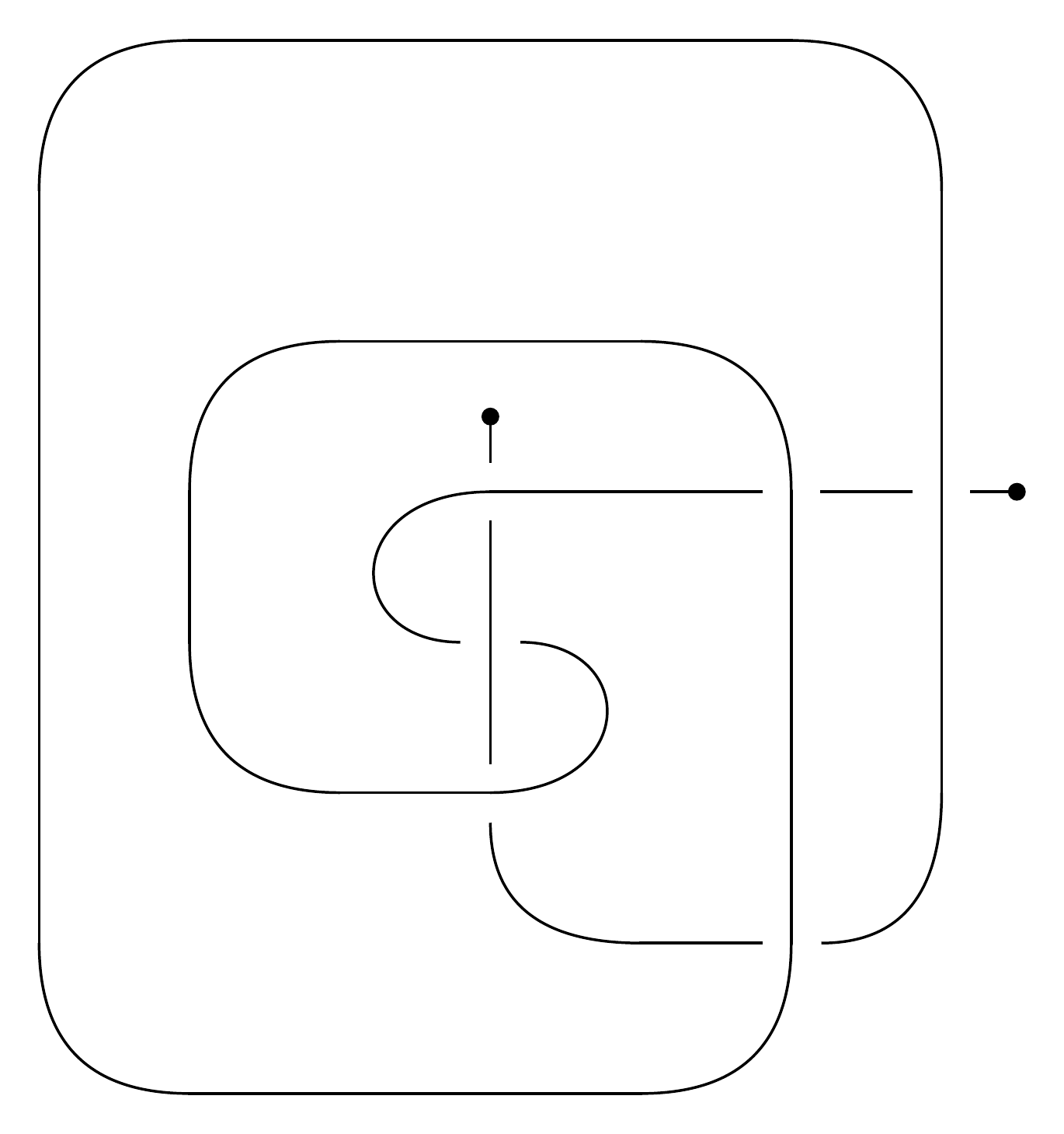}\\
\textcolor{black}{$6_{31}$}
\vspace{1cm}
\end{minipage}
\begin{minipage}[t]{.25\linewidth}
\centering
\includegraphics[width=0.9\textwidth,height=3.5cm,keepaspectratio]{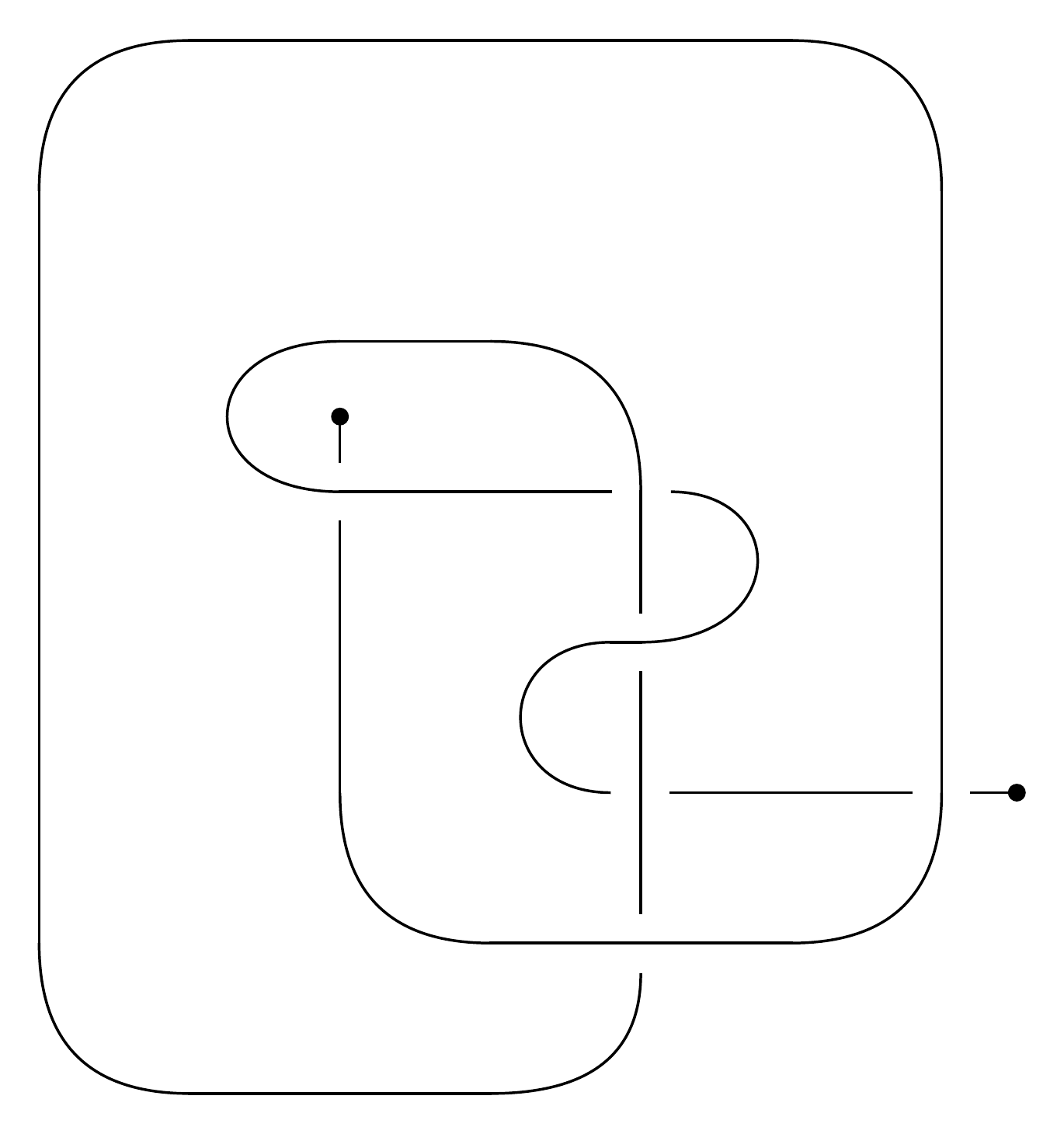}\\
\textcolor{black}{$6_{32}$}
\vspace{1cm}
\end{minipage}
\begin{minipage}[t]{.25\linewidth}
\centering
\includegraphics[width=0.9\textwidth,height=3.5cm,keepaspectratio]{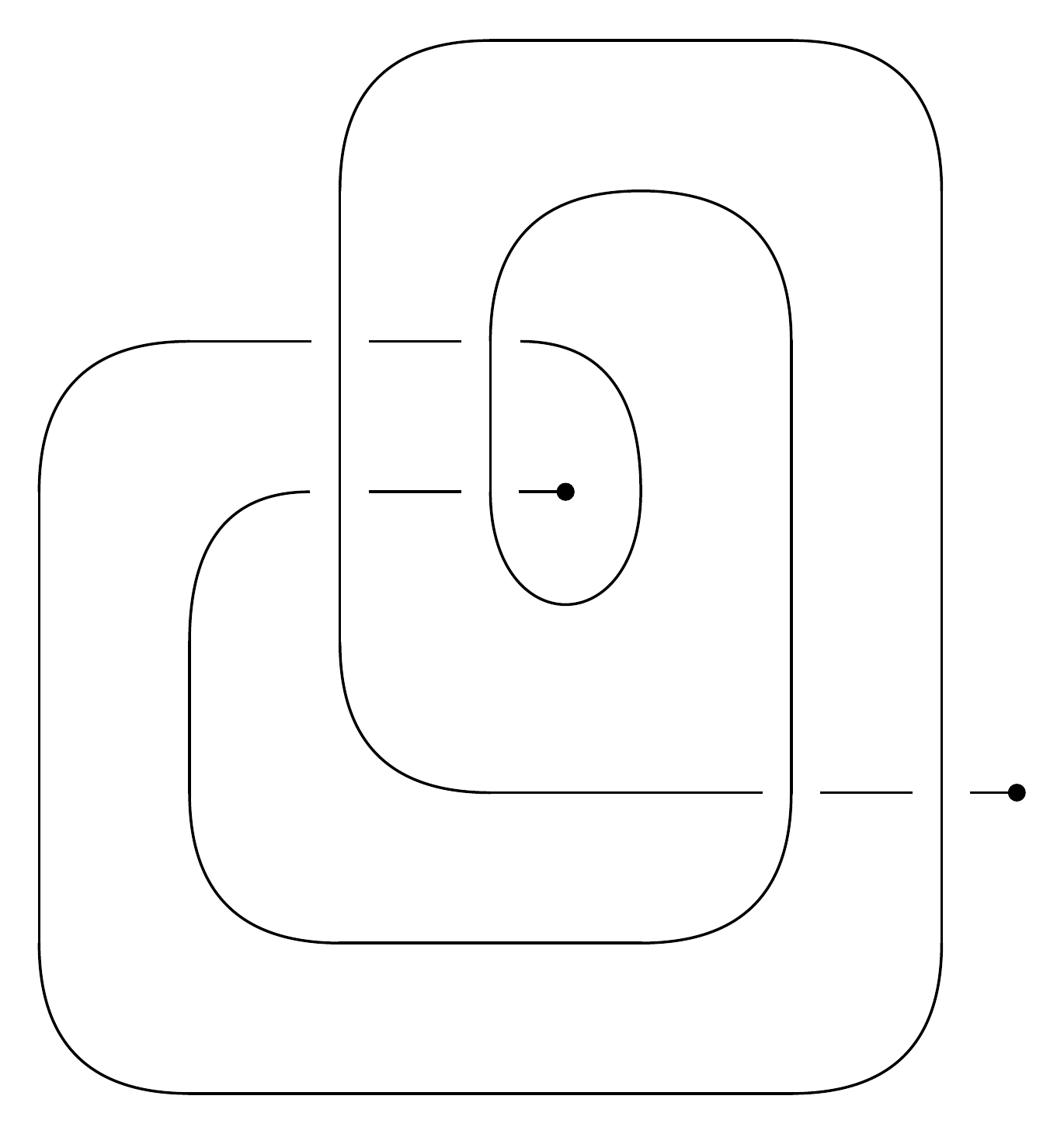}\\
\textcolor{black}{$6_{33}$}
\vspace{1cm}
\end{minipage}
\begin{minipage}[t]{.25\linewidth}
\centering
\includegraphics[width=0.9\textwidth,height=3.5cm,keepaspectratio]{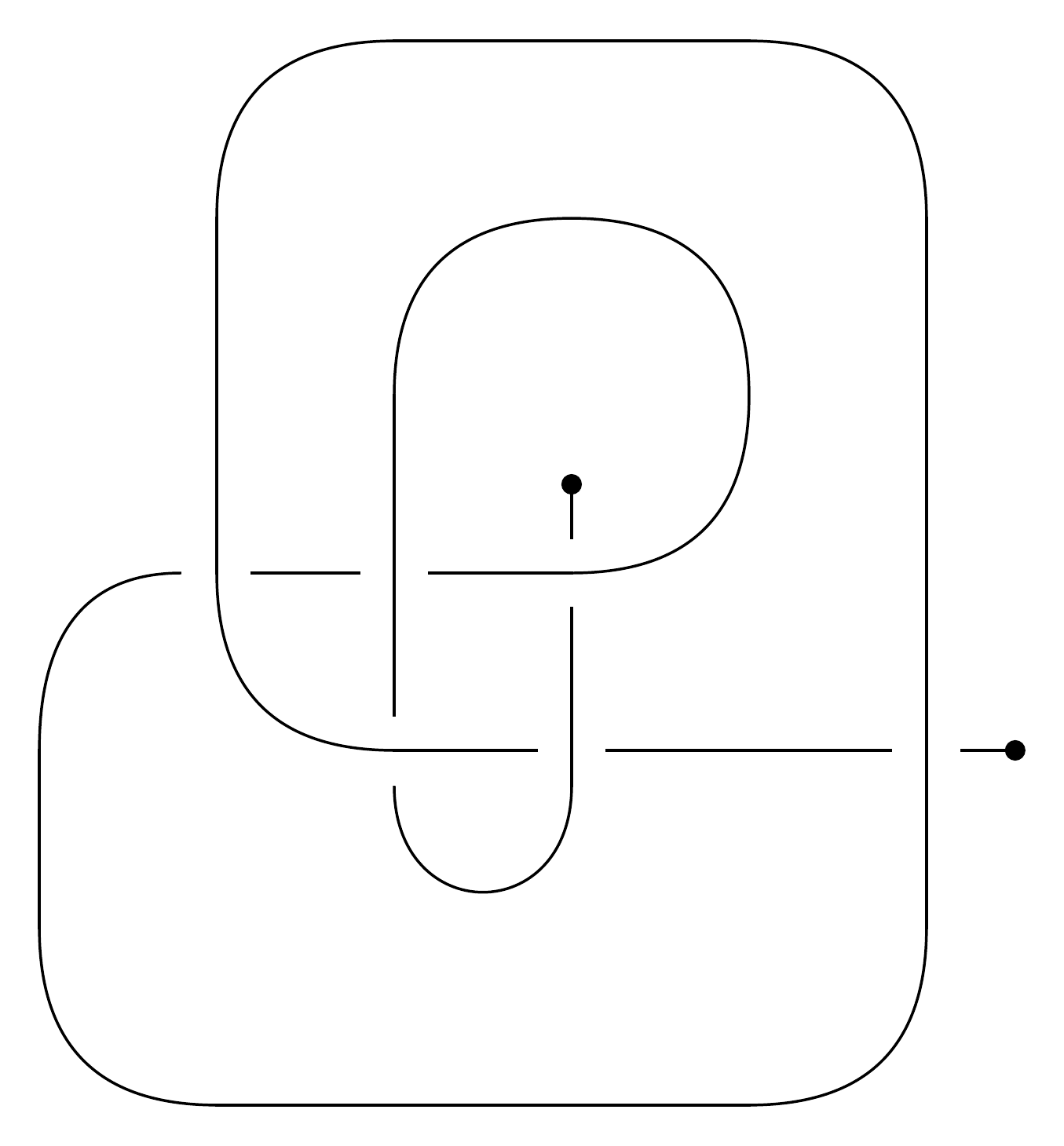}\\
\textcolor{black}{$6_{34}$}
\vspace{1cm}
\end{minipage}
\begin{minipage}[t]{.25\linewidth}
\centering
\includegraphics[width=0.9\textwidth,height=3.5cm,keepaspectratio]{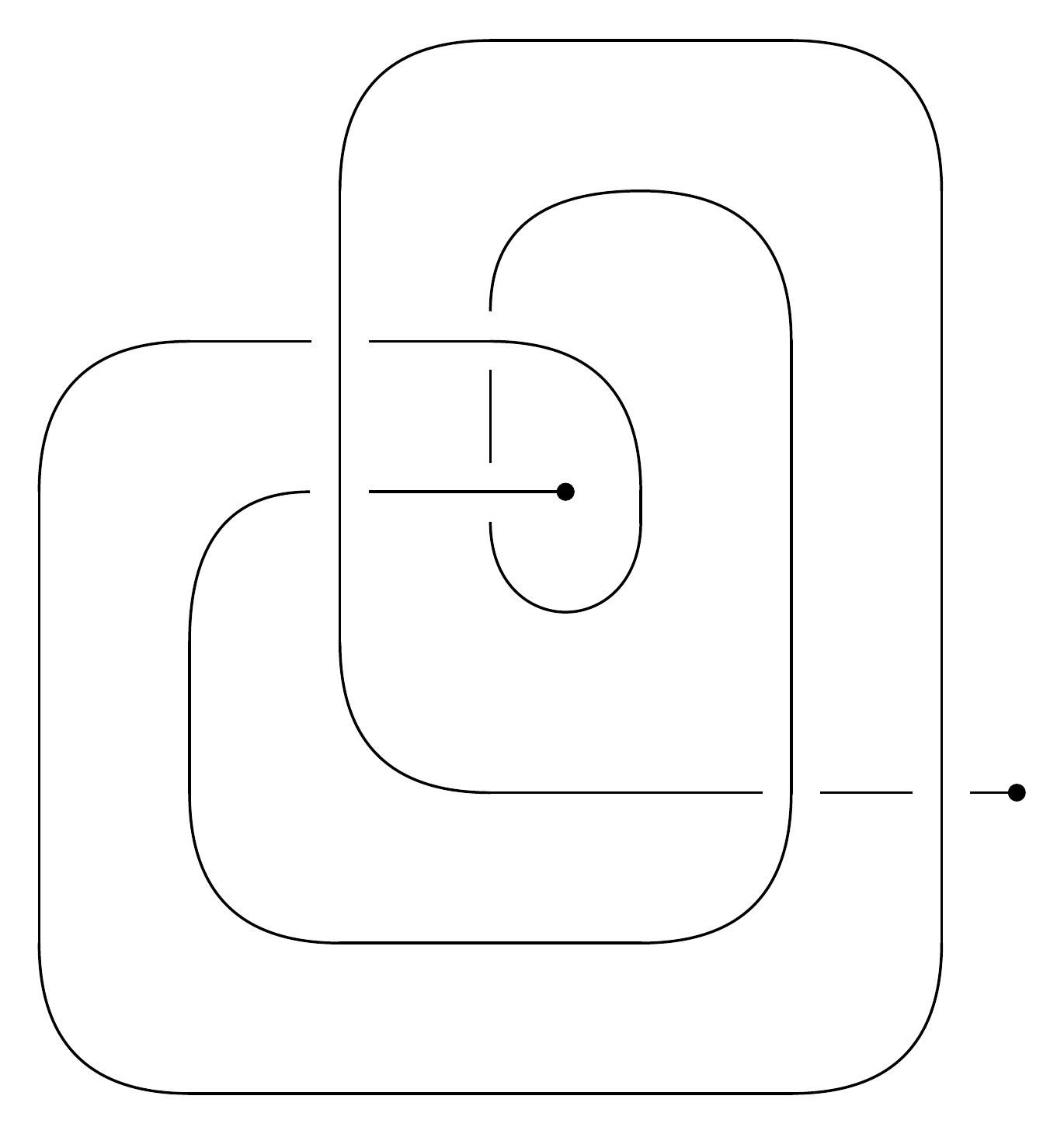}\\
\textcolor{black}{$6_{35}$}
\vspace{1cm}
\end{minipage}
\begin{minipage}[t]{.25\linewidth}
\centering
\includegraphics[width=0.9\textwidth,height=3.5cm,keepaspectratio]{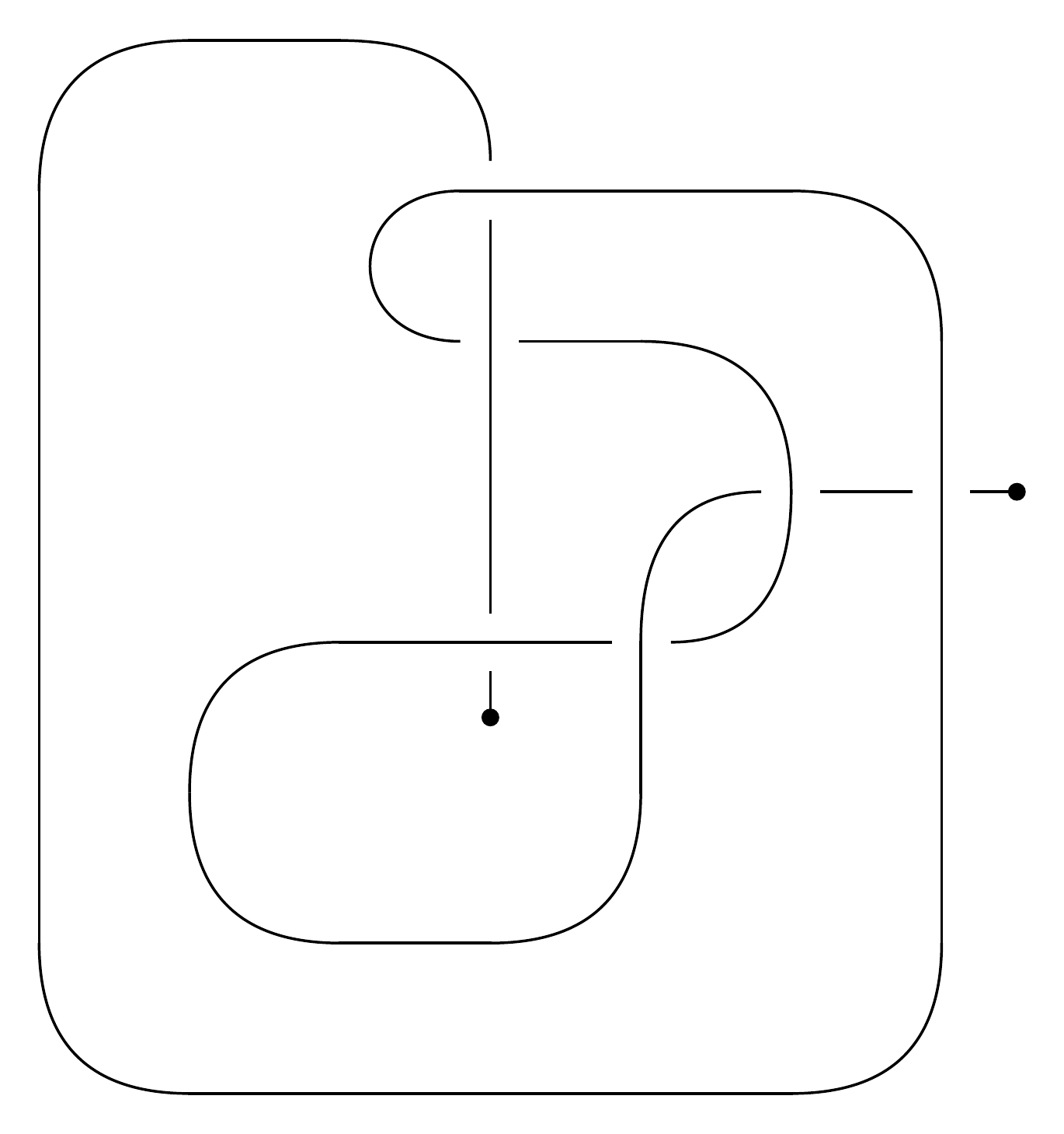}\\
\textcolor{black}{$6_{36}$}
\vspace{1cm}
\end{minipage}
\begin{minipage}[t]{.25\linewidth}
\centering
\includegraphics[width=0.9\textwidth,height=3.5cm,keepaspectratio]{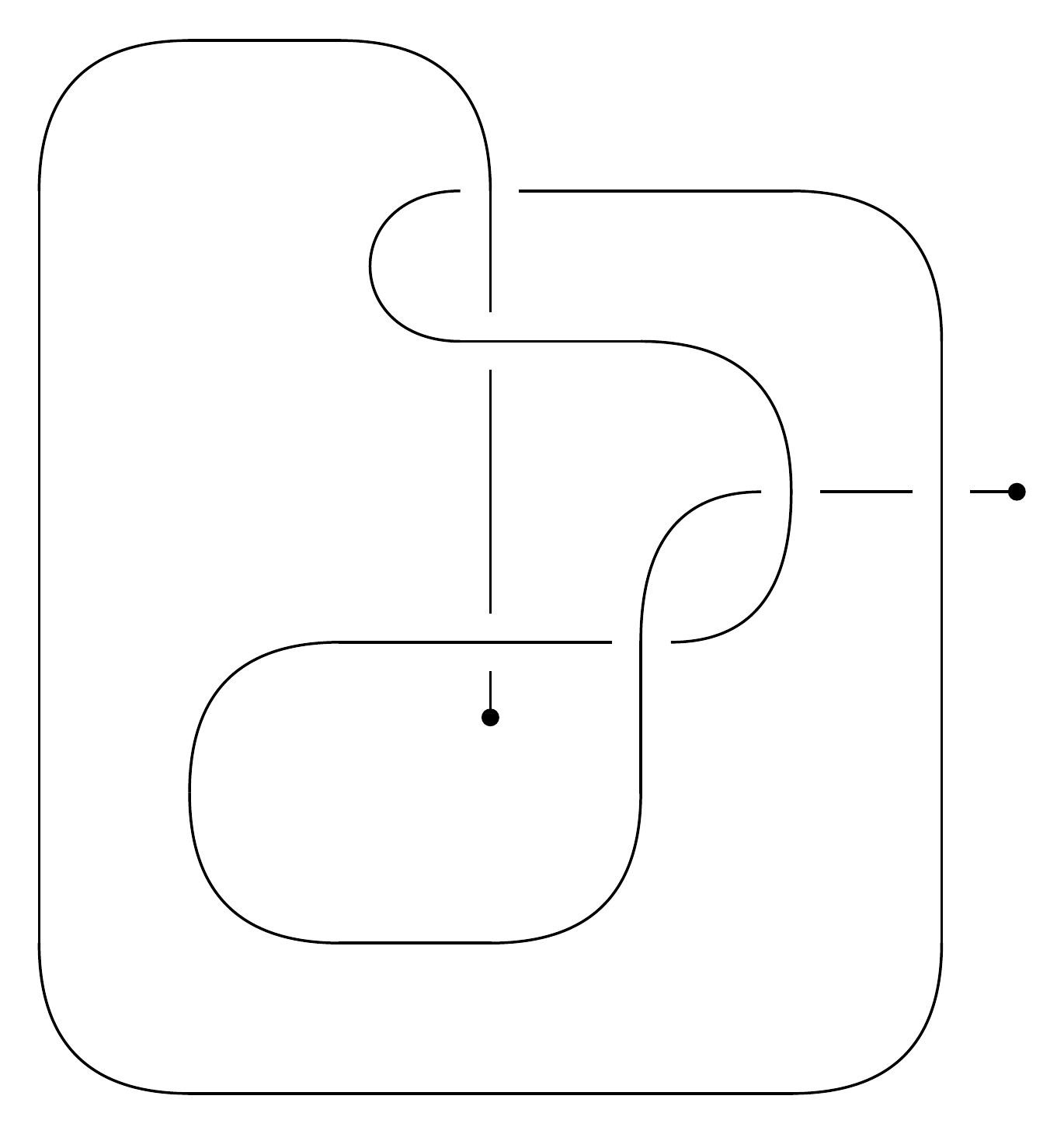}\\
\textcolor{black}{$6_{37}$}
\vspace{1cm}
\end{minipage}
\begin{minipage}[t]{.25\linewidth}
\centering
\includegraphics[width=0.9\textwidth,height=3.5cm,keepaspectratio]{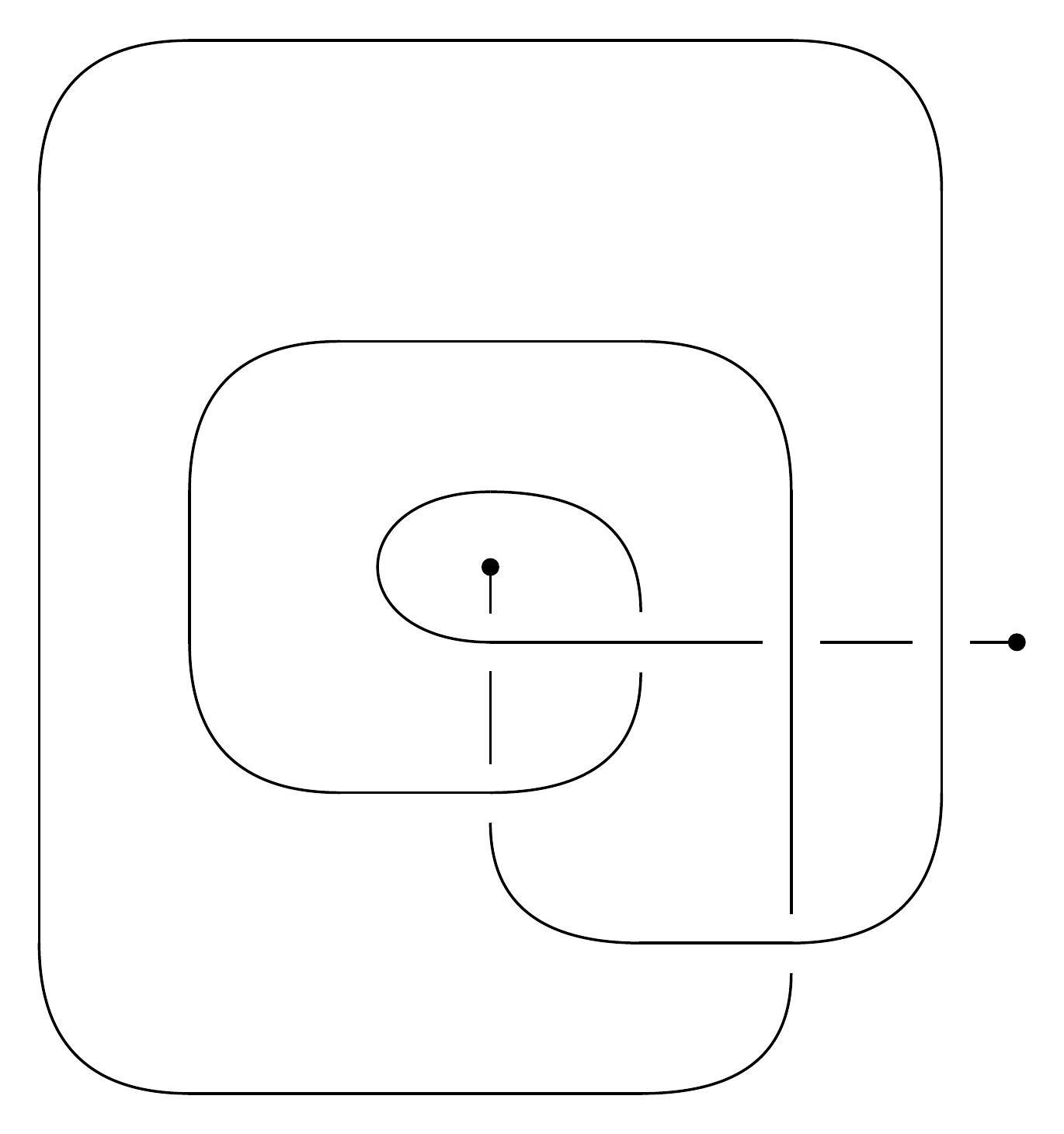}\\
\textcolor{black}{$6_{38}$}
\vspace{1cm}
\end{minipage}
\begin{minipage}[t]{.25\linewidth}
\centering
\includegraphics[width=0.9\textwidth,height=3.5cm,keepaspectratio]{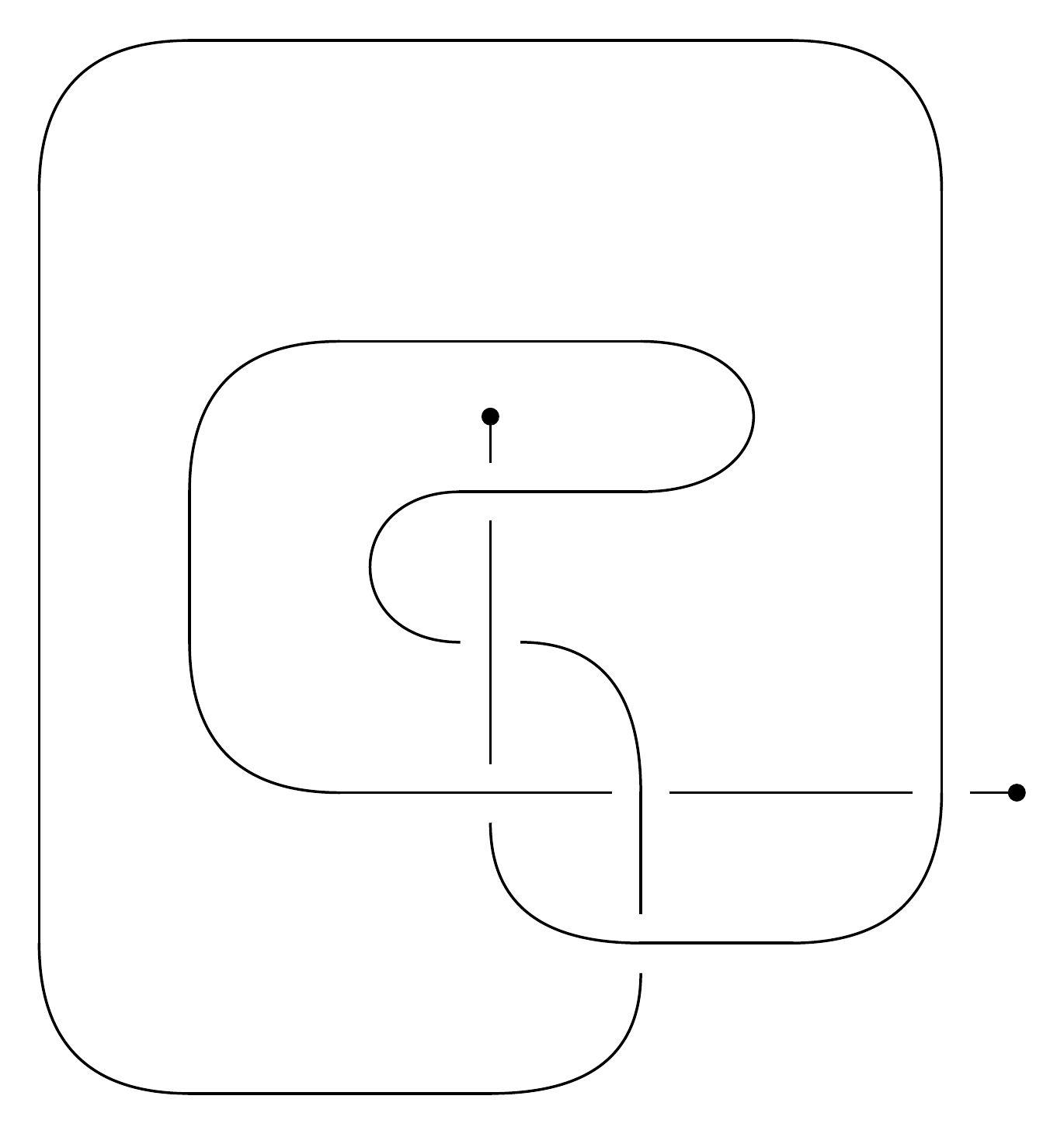}\\
\textcolor{black}{$6_{39}$}
\vspace{1cm}
\end{minipage}
\begin{minipage}[t]{.25\linewidth}
\centering
\includegraphics[width=0.9\textwidth,height=3.5cm,keepaspectratio]{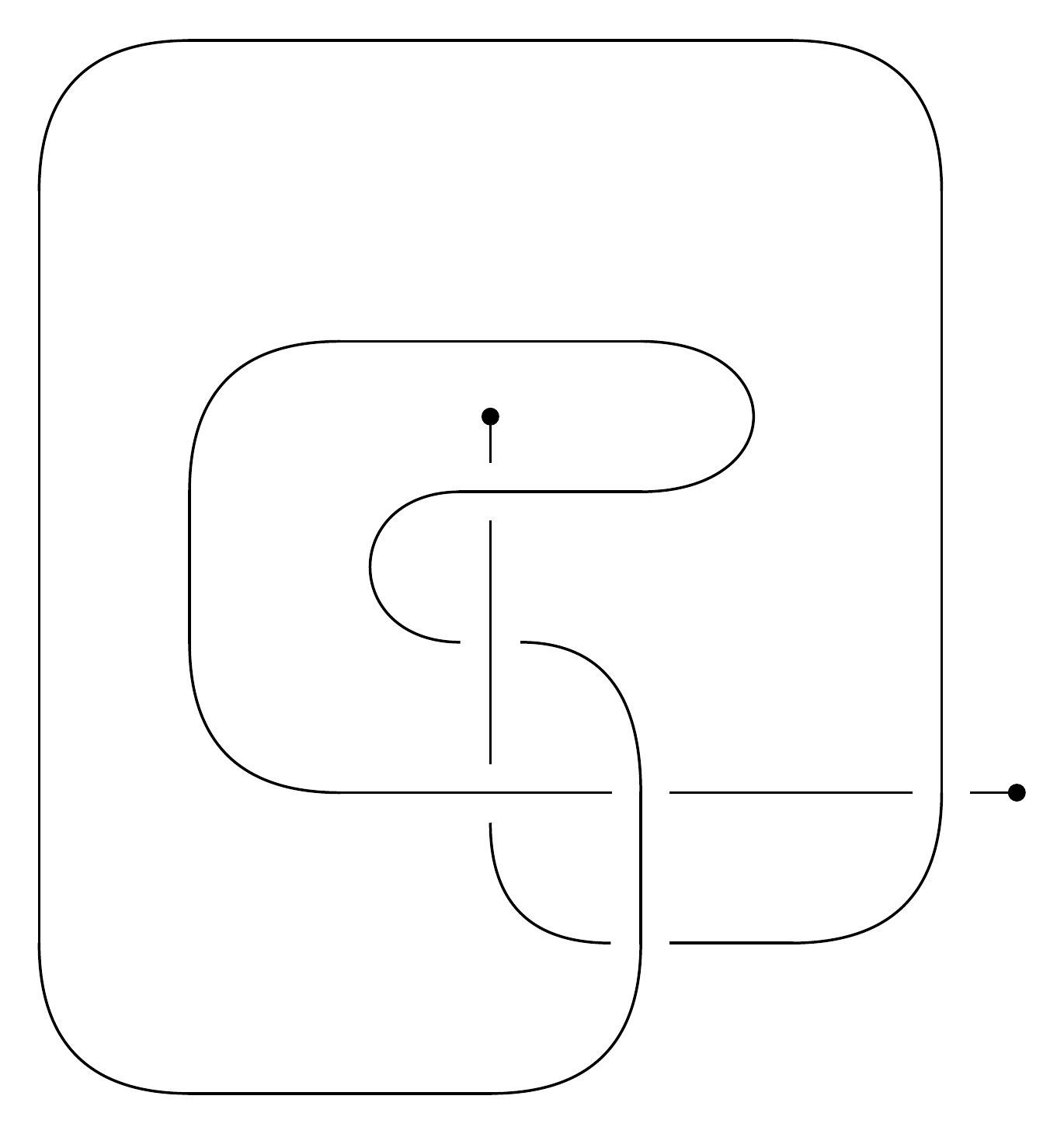}\\
\textcolor{black}{$6_{40}$}
\vspace{1cm}
\end{minipage}
\begin{minipage}[t]{.25\linewidth}
\centering
\includegraphics[width=0.9\textwidth,height=3.5cm,keepaspectratio]{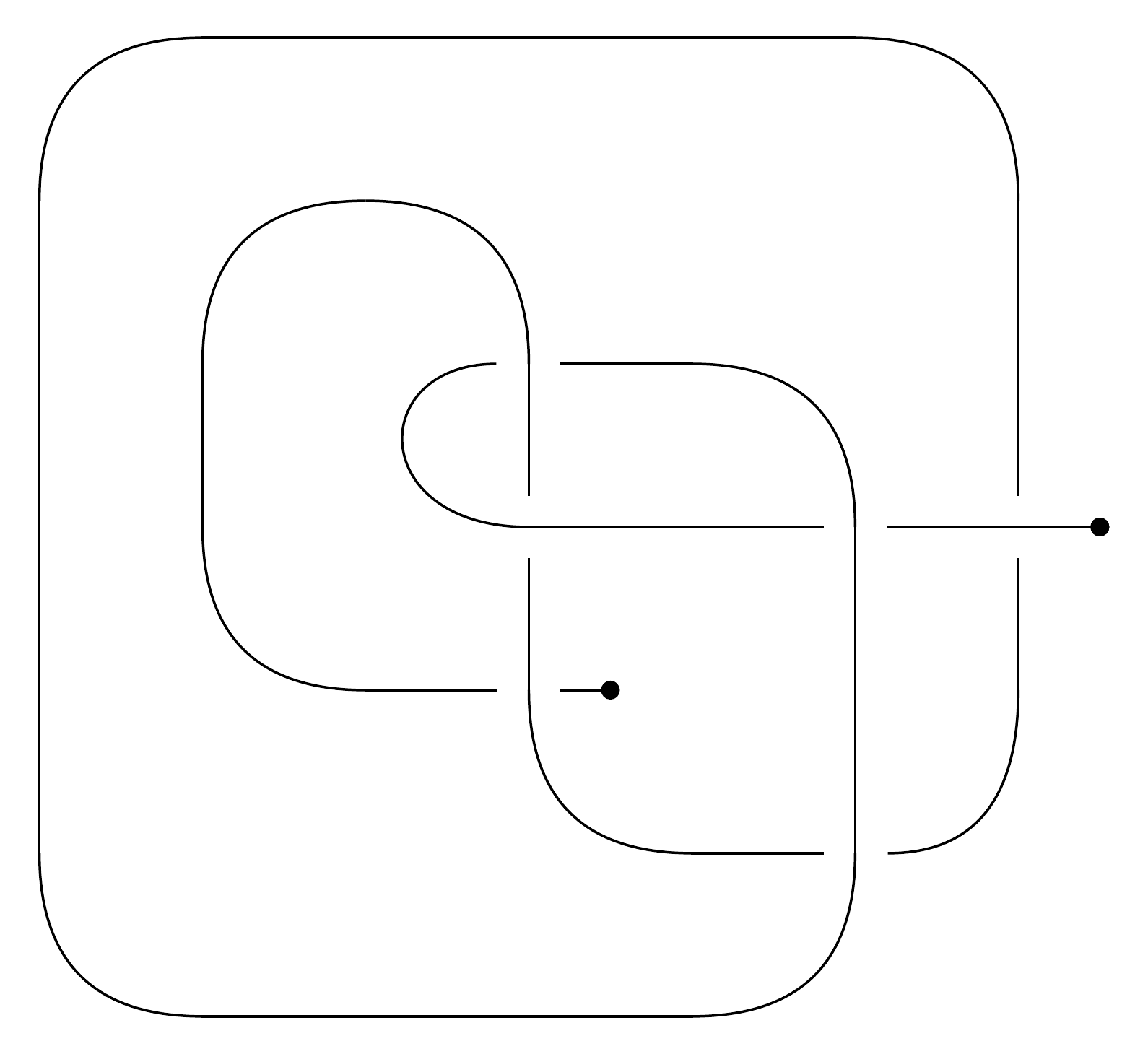}\\
\textcolor{black}{$6_{41}$}
\vspace{1cm}
\end{minipage}
\begin{minipage}[t]{.25\linewidth}
\centering
\includegraphics[width=0.9\textwidth,height=3.5cm,keepaspectratio]{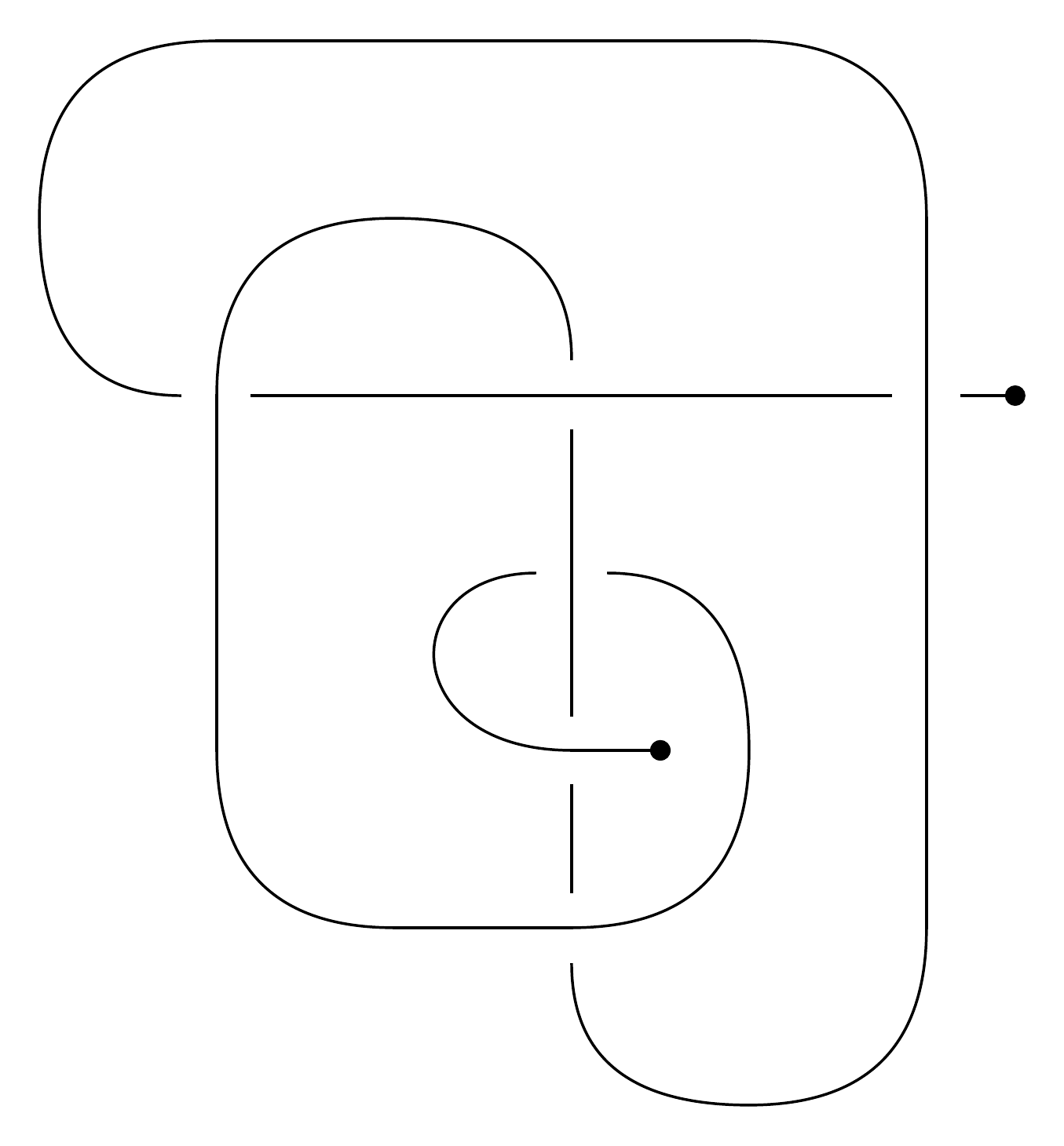}\\
\textcolor{black}{$6_{42}$}
\vspace{1cm}
\end{minipage}
\begin{minipage}[t]{.25\linewidth}
\centering
\includegraphics[width=0.9\textwidth,height=3.5cm,keepaspectratio]{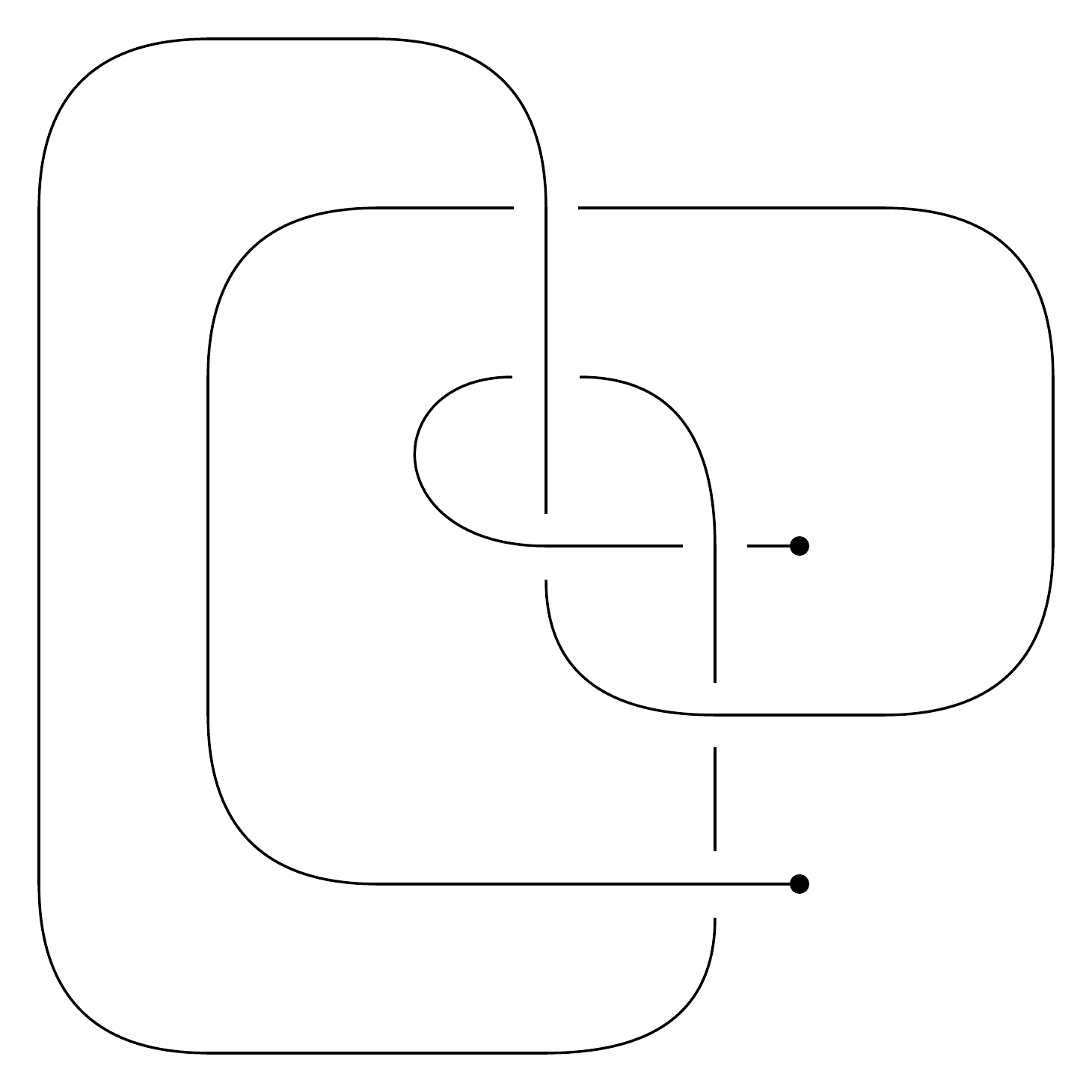}\\
\textcolor{black}{$6_{43}$}
\vspace{1cm}
\end{minipage}
\begin{minipage}[t]{.25\linewidth}
\centering
\includegraphics[width=0.9\textwidth,height=3.5cm,keepaspectratio]{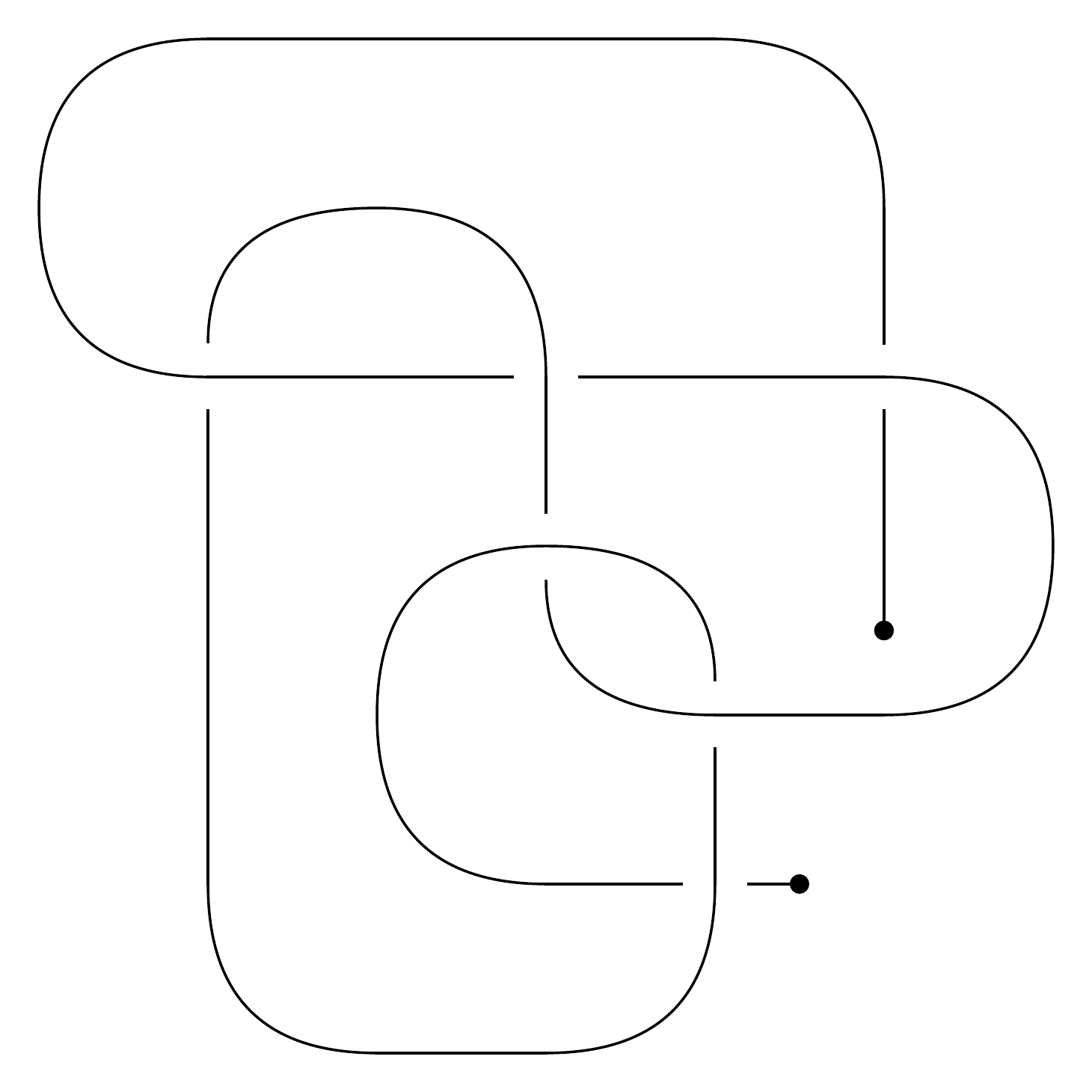}\\
\textcolor{black}{$6_{44}$}
\vspace{1cm}
\end{minipage}
\begin{minipage}[t]{.25\linewidth}
\centering
\includegraphics[width=0.9\textwidth,height=3.5cm,keepaspectratio]{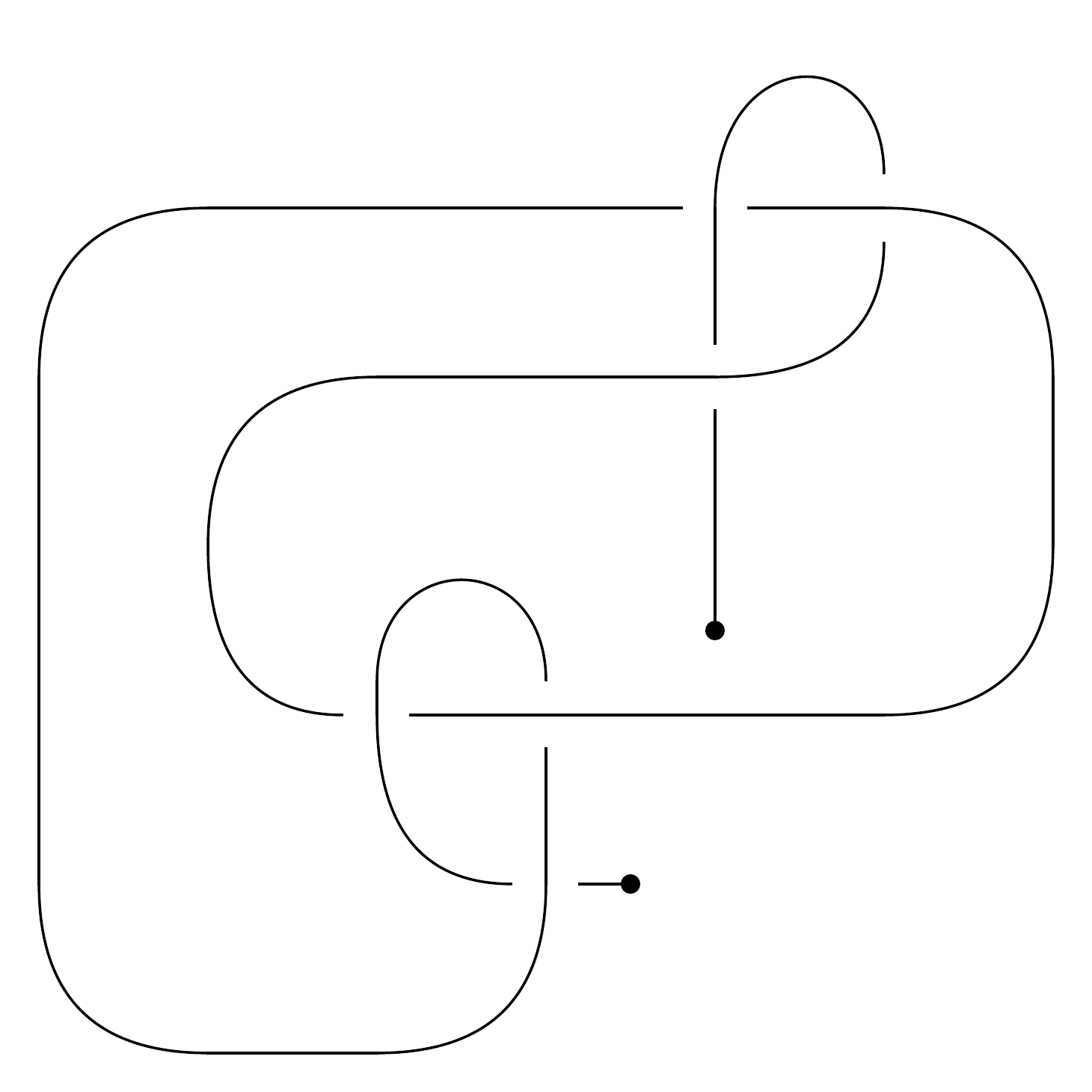}\\
\textcolor{black}{$6_{45}$}
\vspace{1cm}
\end{minipage}
\begin{minipage}[t]{.25\linewidth}
\centering
\includegraphics[width=0.9\textwidth,height=3.5cm,keepaspectratio]{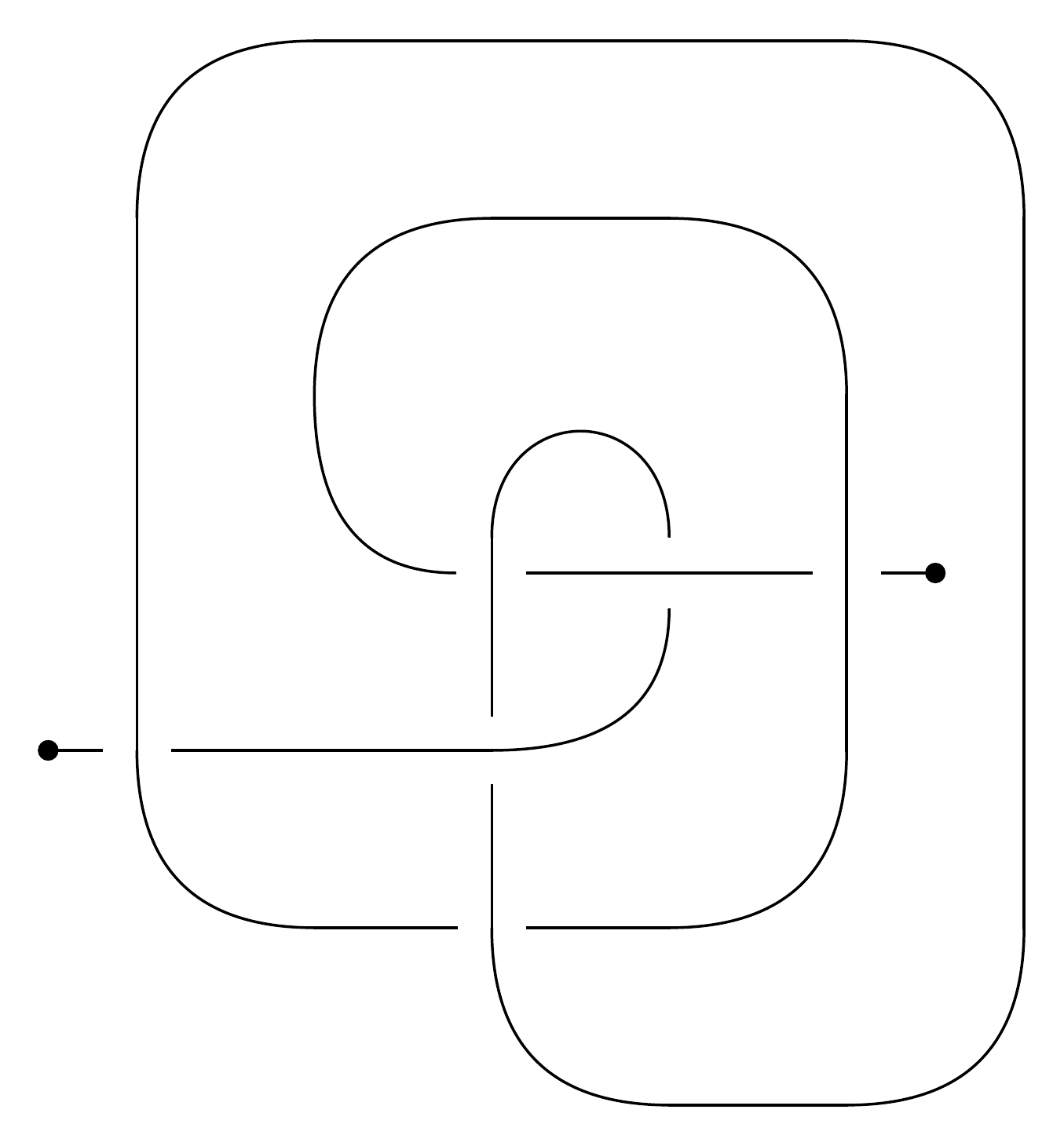}\\
\textcolor{black}{$6_{46}$}
\vspace{1cm}
\end{minipage}
\begin{minipage}[t]{.25\linewidth}
\centering
\includegraphics[width=0.9\textwidth,height=3.5cm,keepaspectratio]{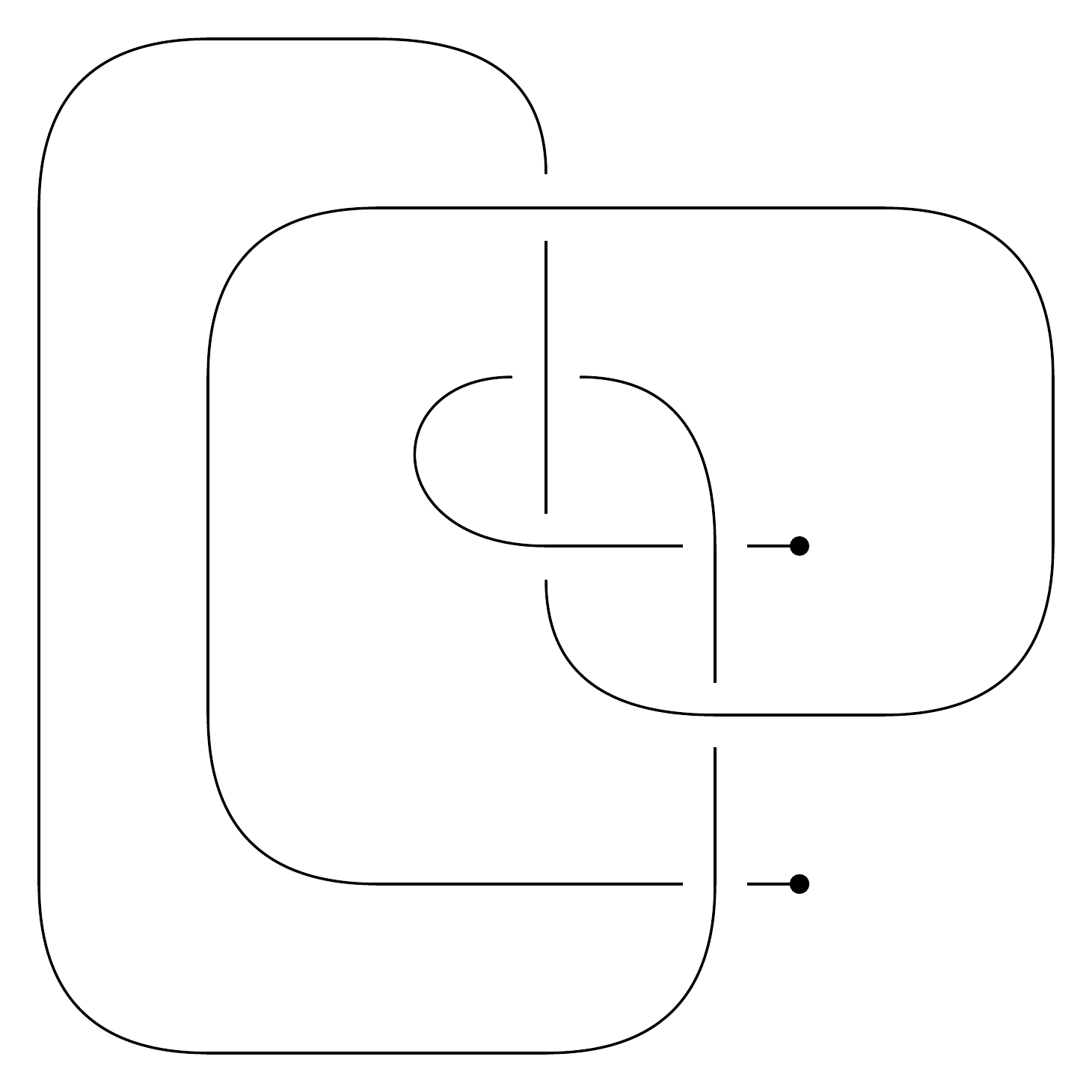}\\
\textcolor{black}{$6_{47}$}
\vspace{1cm}
\end{minipage}
\begin{minipage}[t]{.25\linewidth}
\centering
\includegraphics[width=0.9\textwidth,height=3.5cm,keepaspectratio]{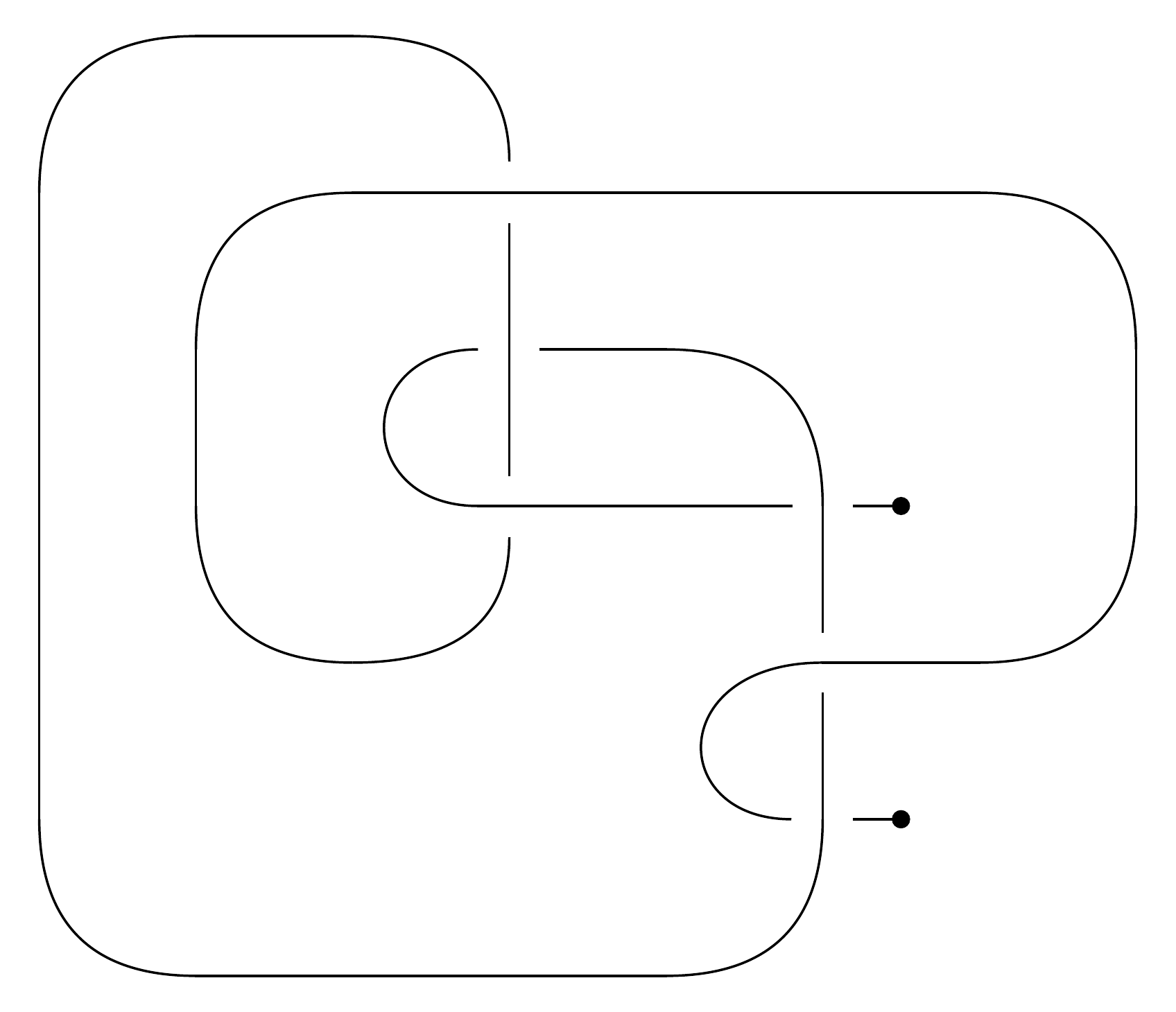}\\
\textcolor{black}{$6_{48}$}
\vspace{1cm}
\end{minipage}
\begin{minipage}[t]{.25\linewidth}
\centering
\includegraphics[width=0.9\textwidth,height=3.5cm,keepaspectratio]{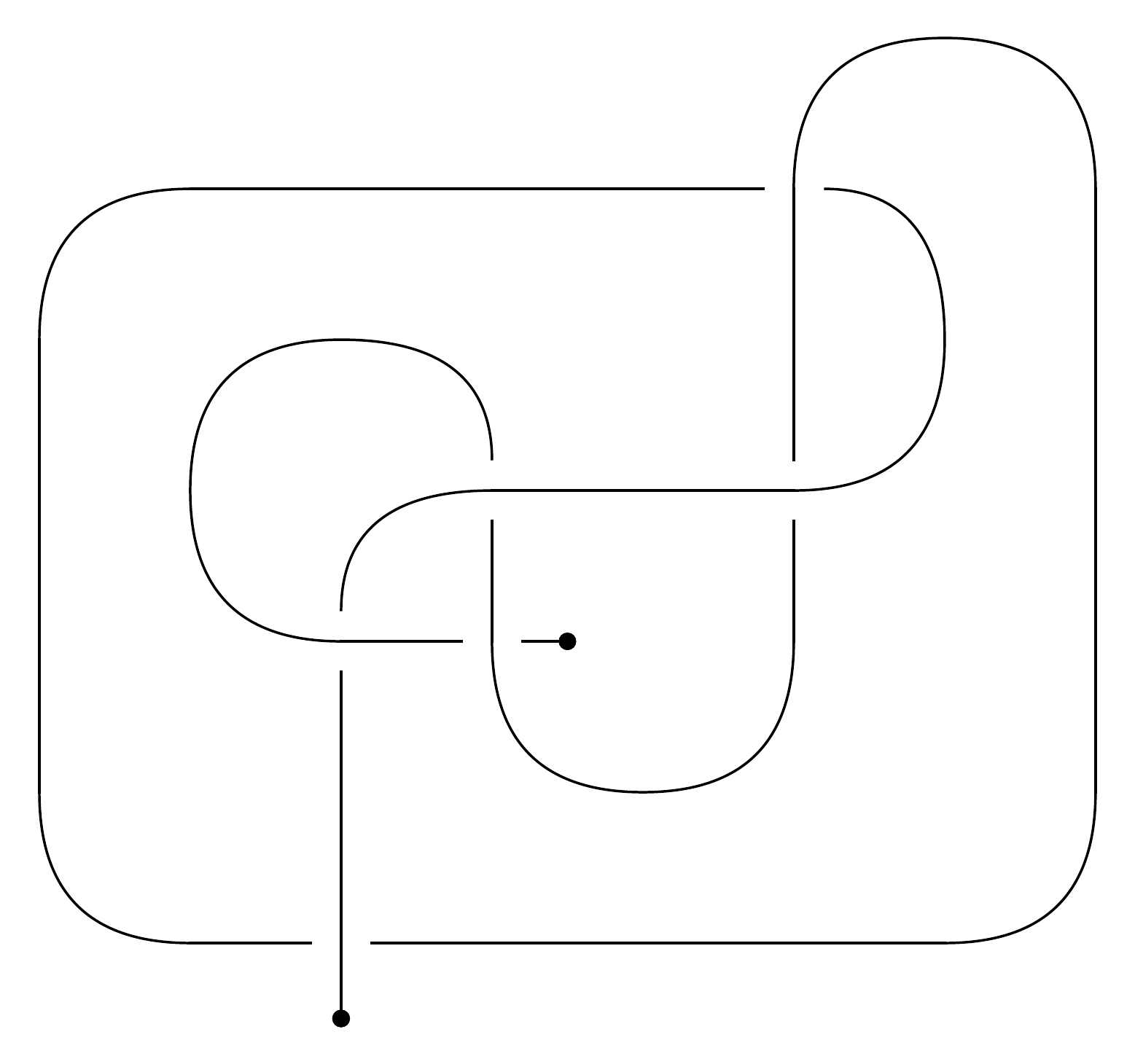}\\
\textcolor{black}{$6_{49}$}
\vspace{1cm}
\end{minipage}
\begin{minipage}[t]{.25\linewidth}
\centering
\includegraphics[width=0.9\textwidth,height=3.5cm,keepaspectratio]{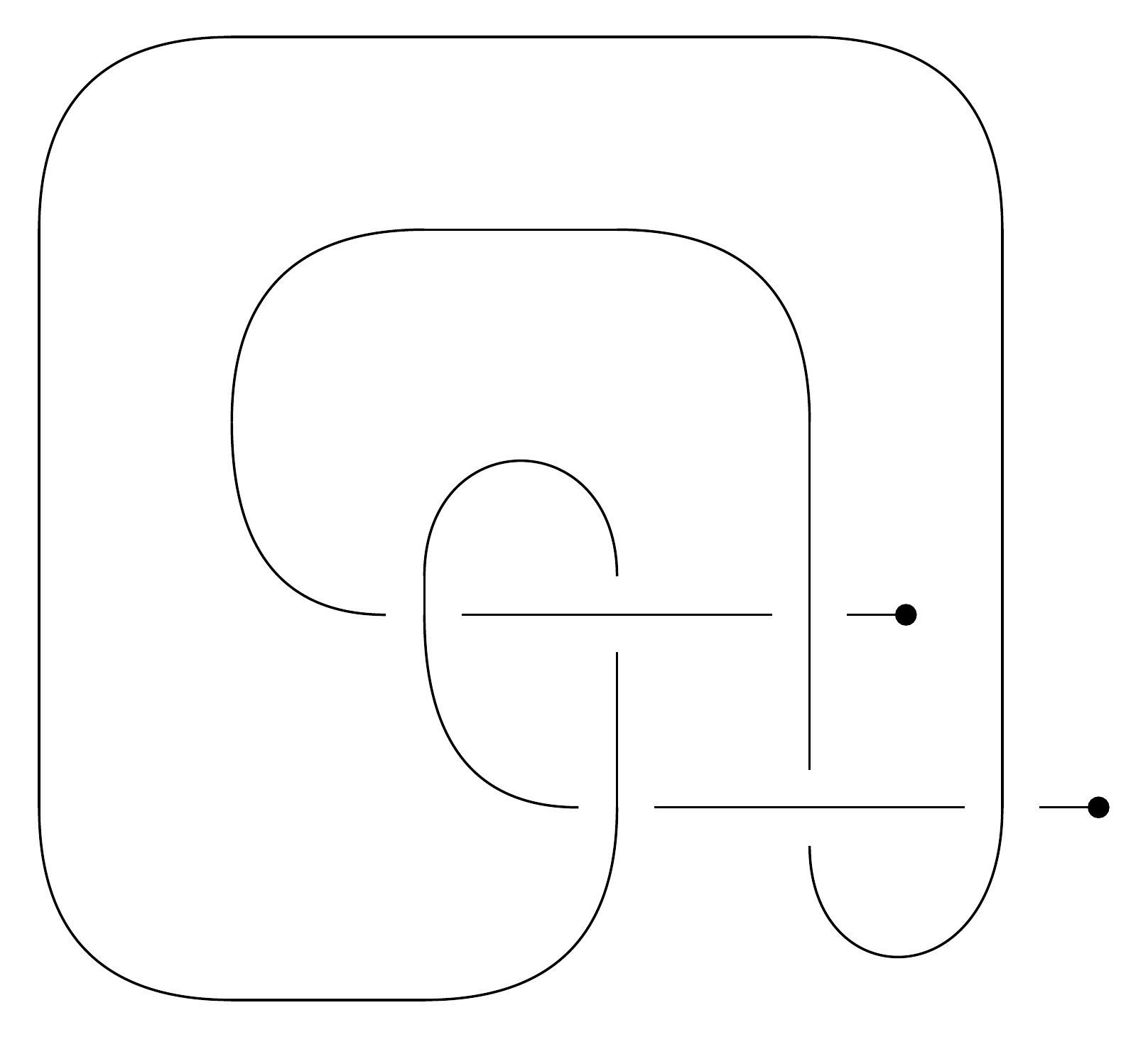}\\
\textcolor{black}{$6_{50}$}
\vspace{1cm}
\end{minipage}
\begin{minipage}[t]{.25\linewidth}
\centering
\includegraphics[width=0.9\textwidth,height=3.5cm,keepaspectratio]{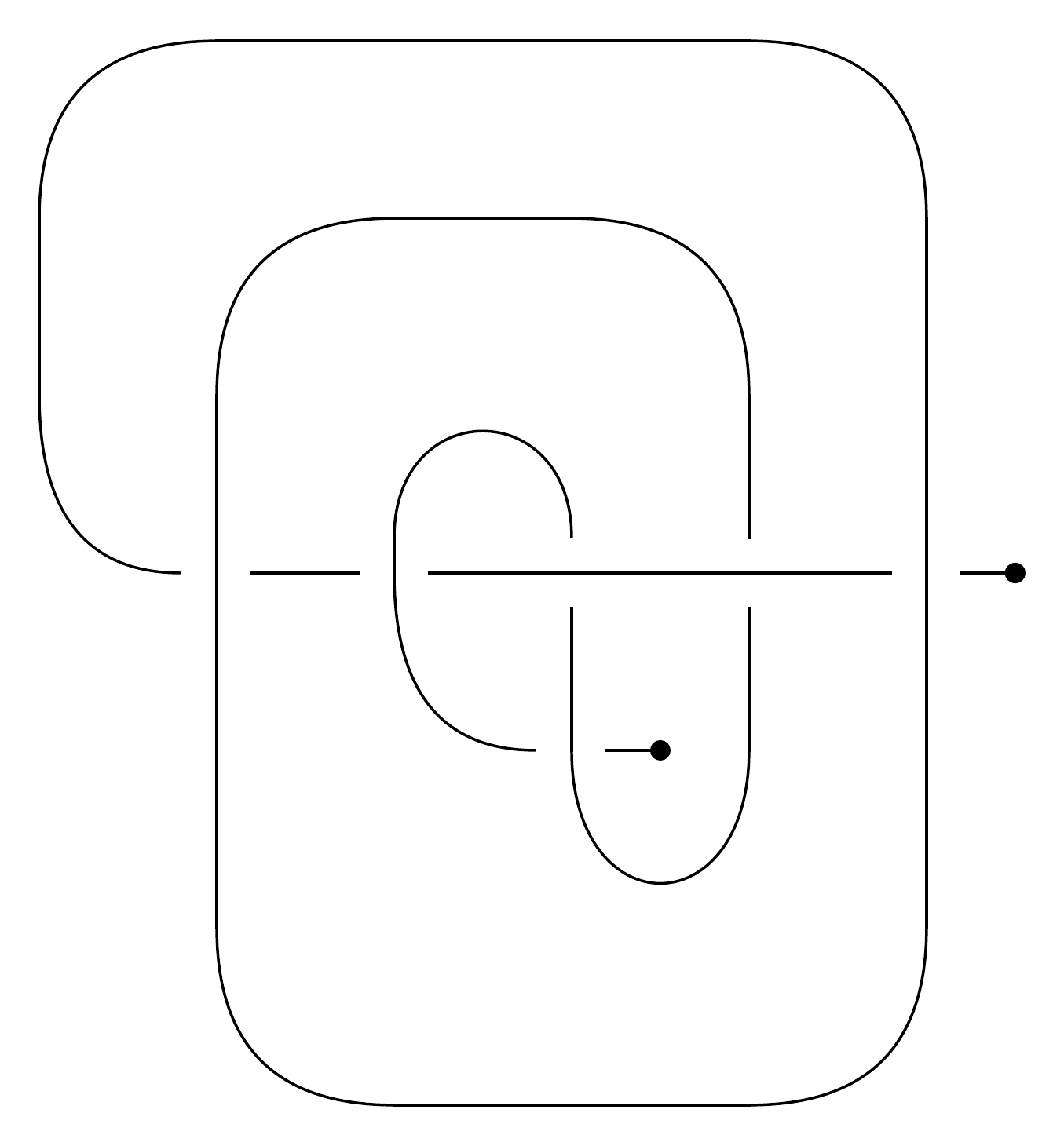}\\
\textcolor{black}{$6_{51}$}
\vspace{1cm}
\end{minipage}
\begin{minipage}[t]{.25\linewidth}
\centering
\includegraphics[width=0.9\textwidth,height=3.5cm,keepaspectratio]{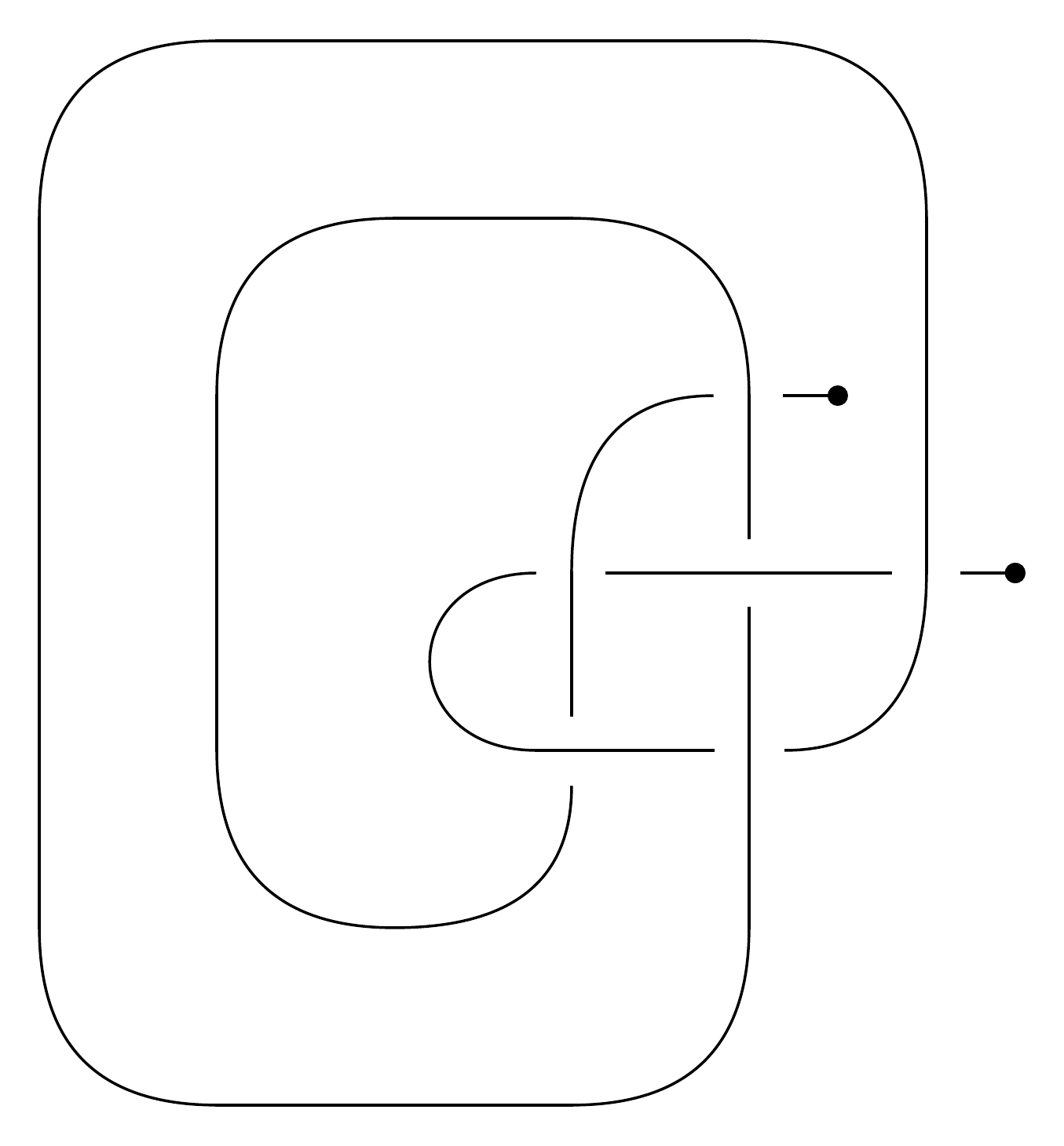}\\
\textcolor{black}{$6_{52}$}
\vspace{1cm}
\end{minipage}
\begin{minipage}[t]{.25\linewidth}
\centering
\includegraphics[width=0.9\textwidth,height=3.5cm,keepaspectratio]{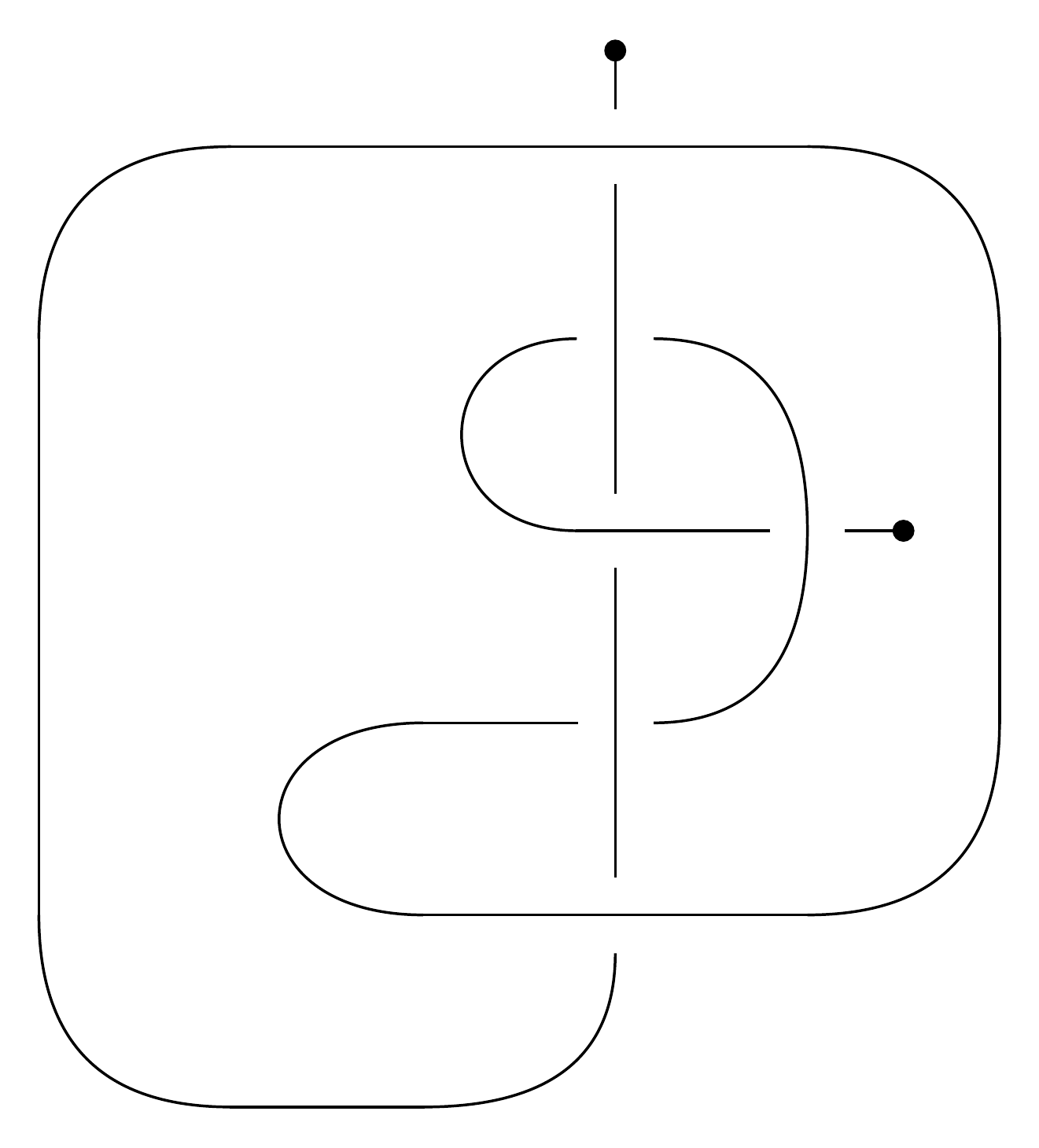}\\
\textcolor{black}{$6_{53}$}
\vspace{1cm}
\end{minipage}
\begin{minipage}[t]{.25\linewidth}
\centering
\includegraphics[width=0.9\textwidth,height=3.5cm,keepaspectratio]{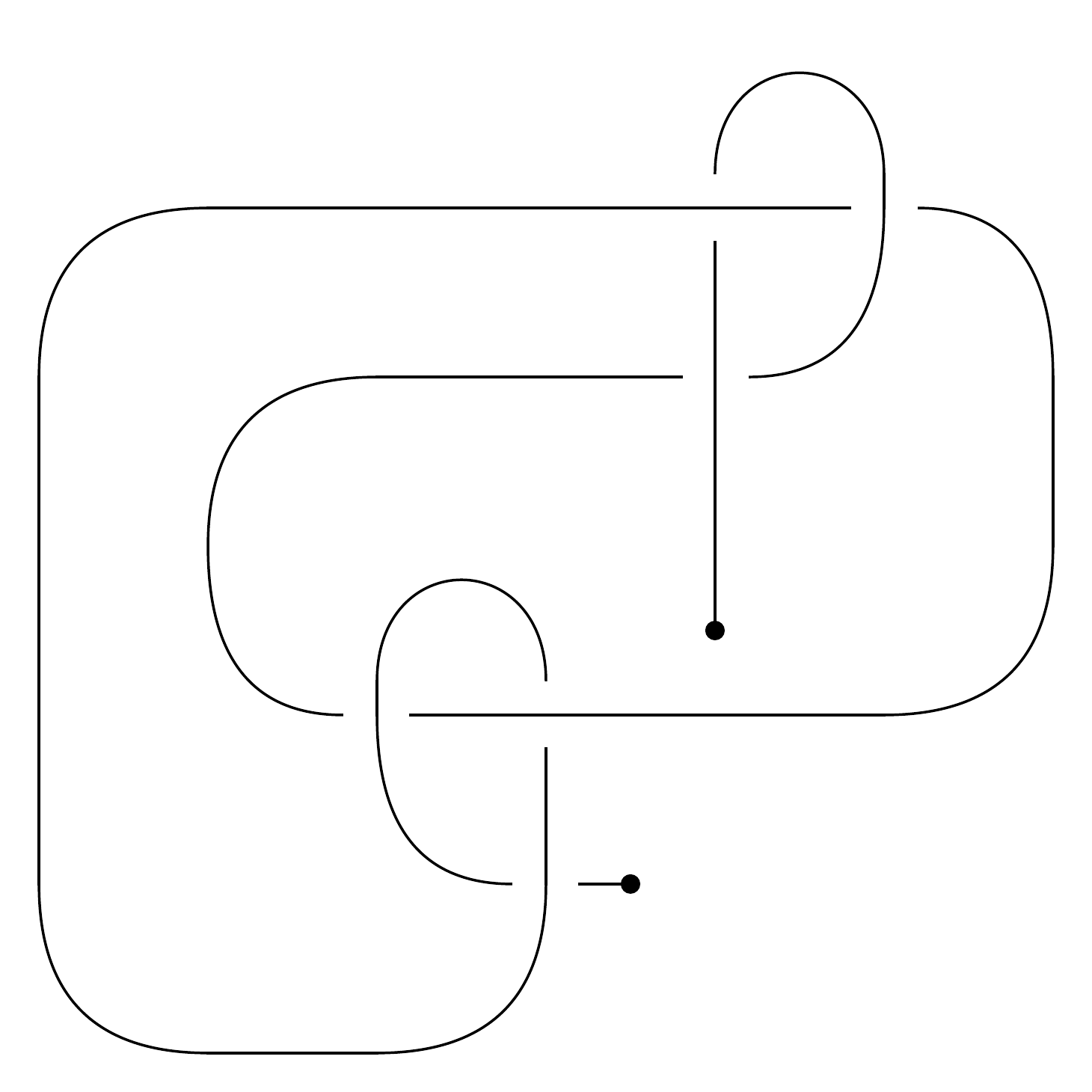}\\
\textcolor{black}{$6_{54}$}
\vspace{1cm}
\end{minipage}
\begin{minipage}[t]{.25\linewidth}
\centering
\includegraphics[width=0.9\textwidth,height=3.5cm,keepaspectratio]{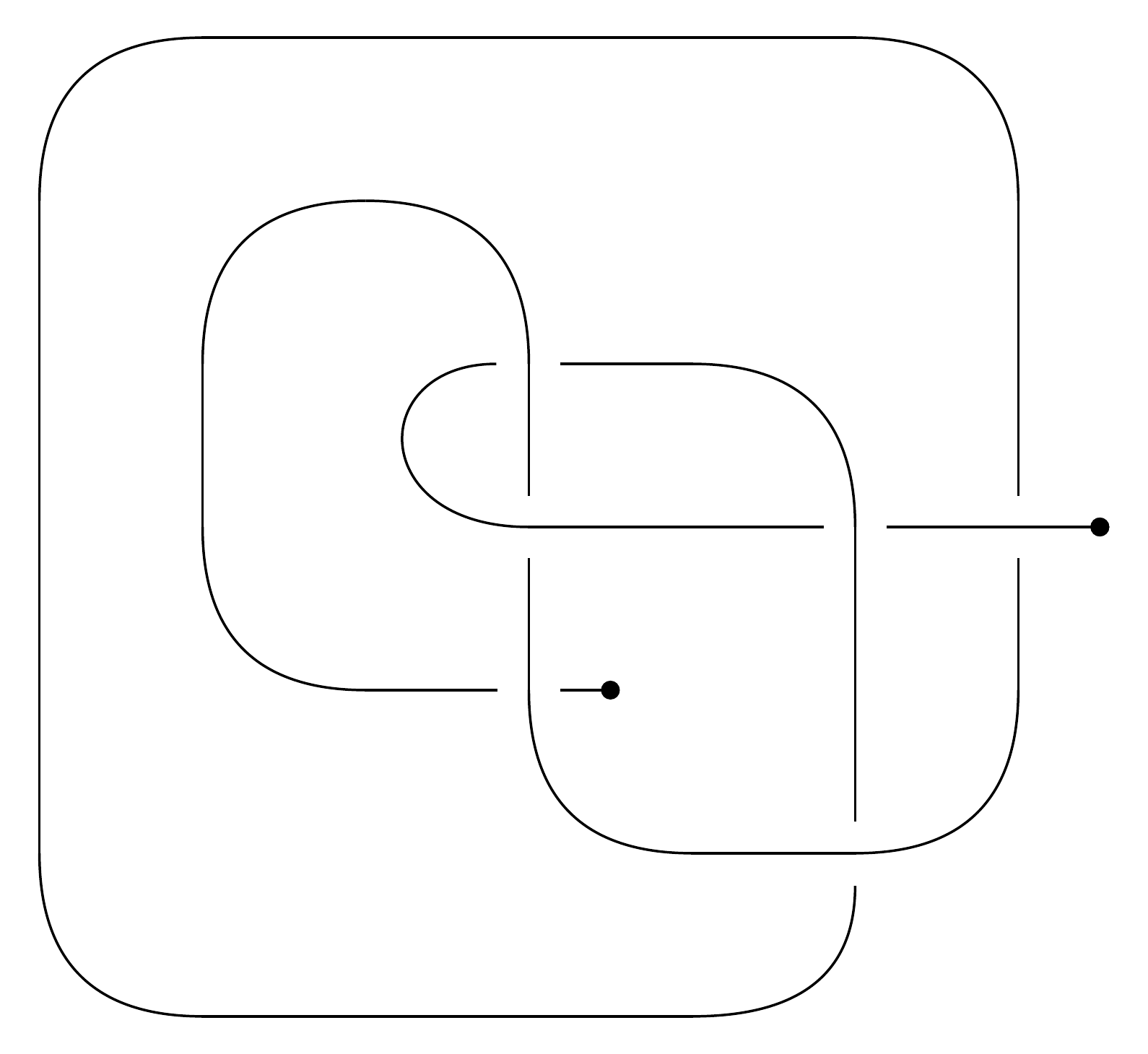}\\
\textcolor{black}{$6_{55}$}
\vspace{1cm}
\end{minipage}
\begin{minipage}[t]{.25\linewidth}
\centering
\includegraphics[width=0.9\textwidth,height=3.5cm,keepaspectratio]{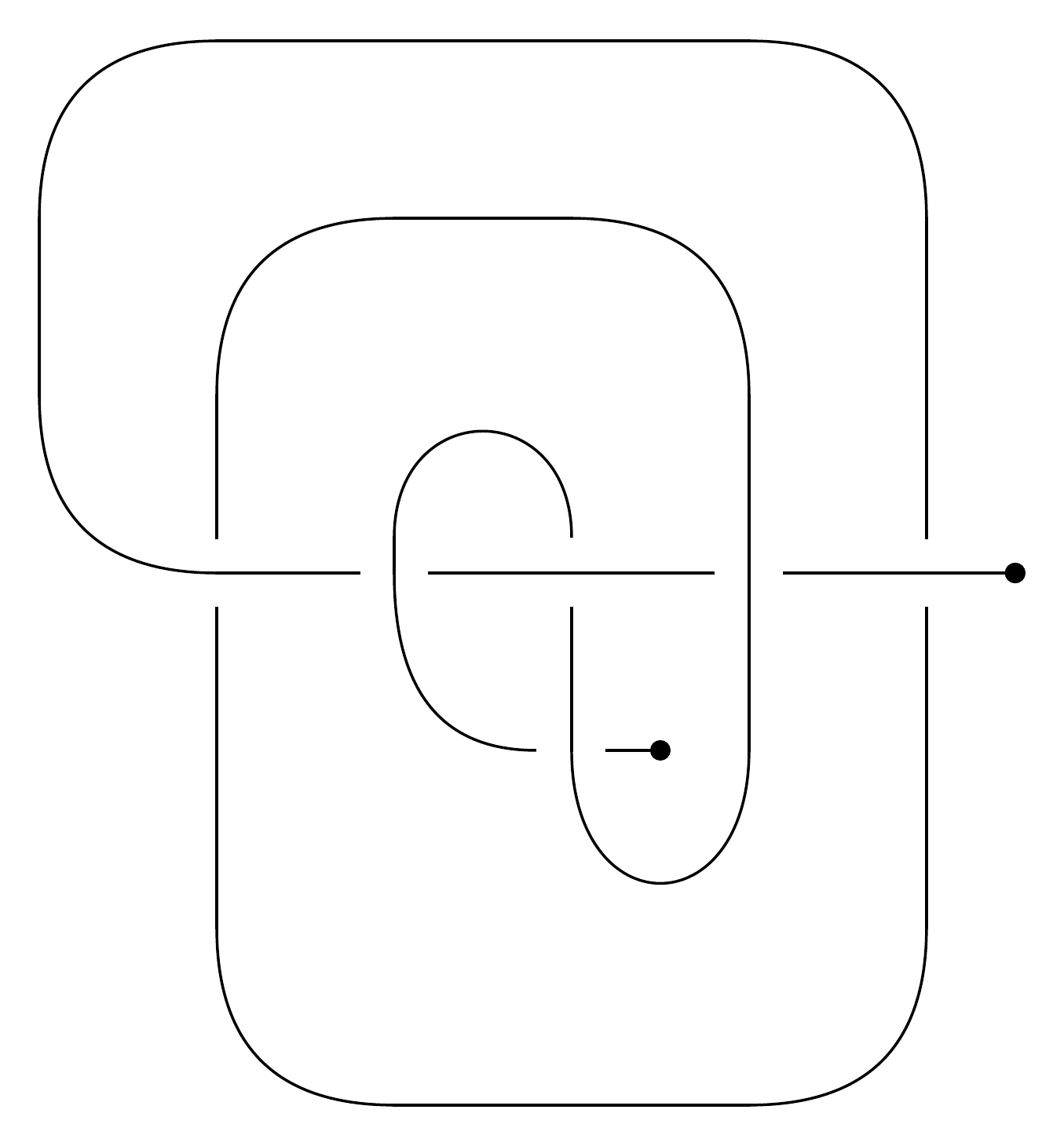}\\
\textcolor{black}{$6_{56}$}
\vspace{1cm}
\end{minipage}
\begin{minipage}[t]{.25\linewidth}
\centering
\includegraphics[width=0.9\textwidth,height=3.5cm,keepaspectratio]{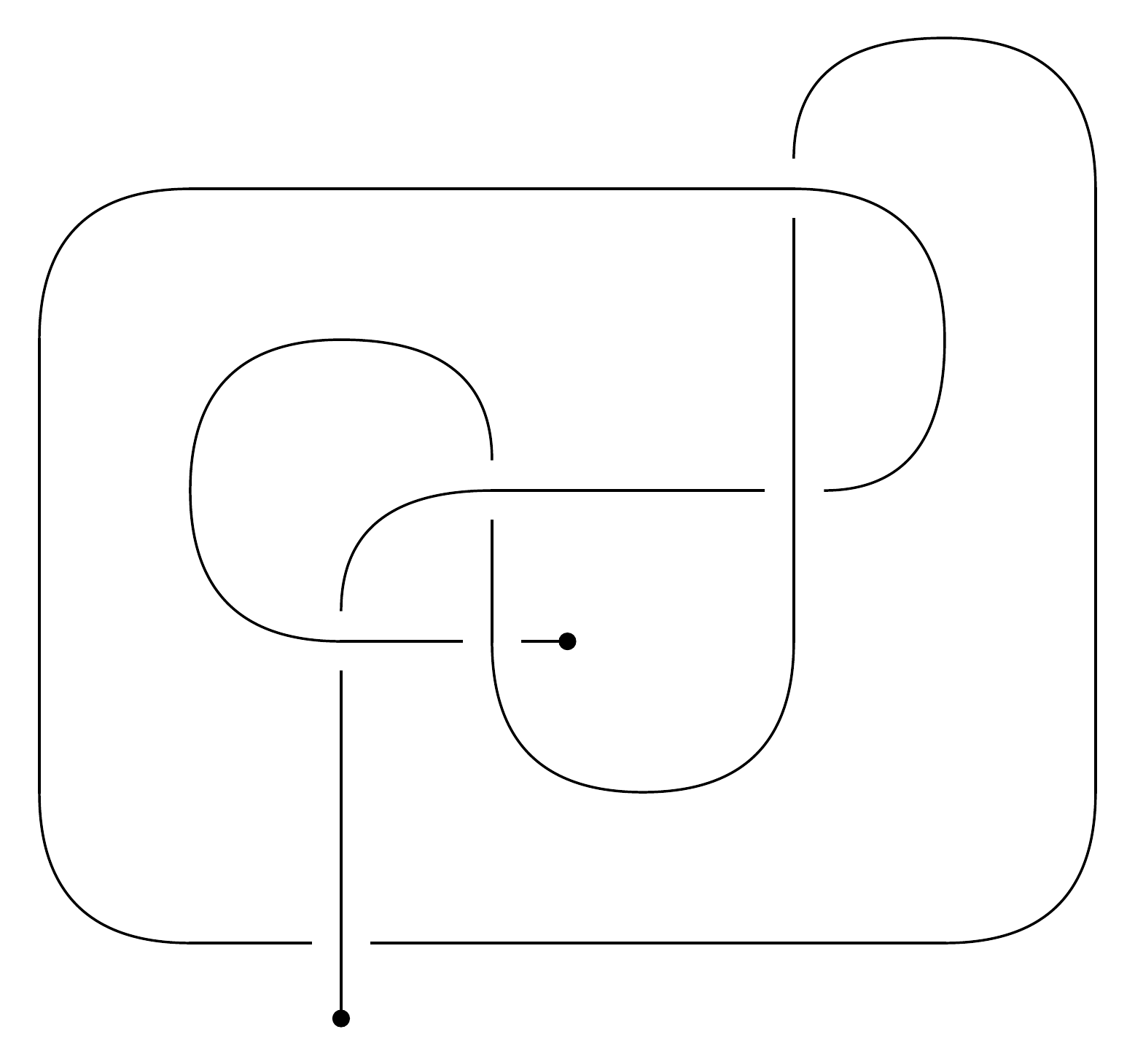}\\
\textcolor{black}{$6_{57}$}
\vspace{1cm}
\end{minipage}
\begin{minipage}[t]{.25\linewidth}
\centering
\includegraphics[width=0.9\textwidth,height=3.5cm,keepaspectratio]{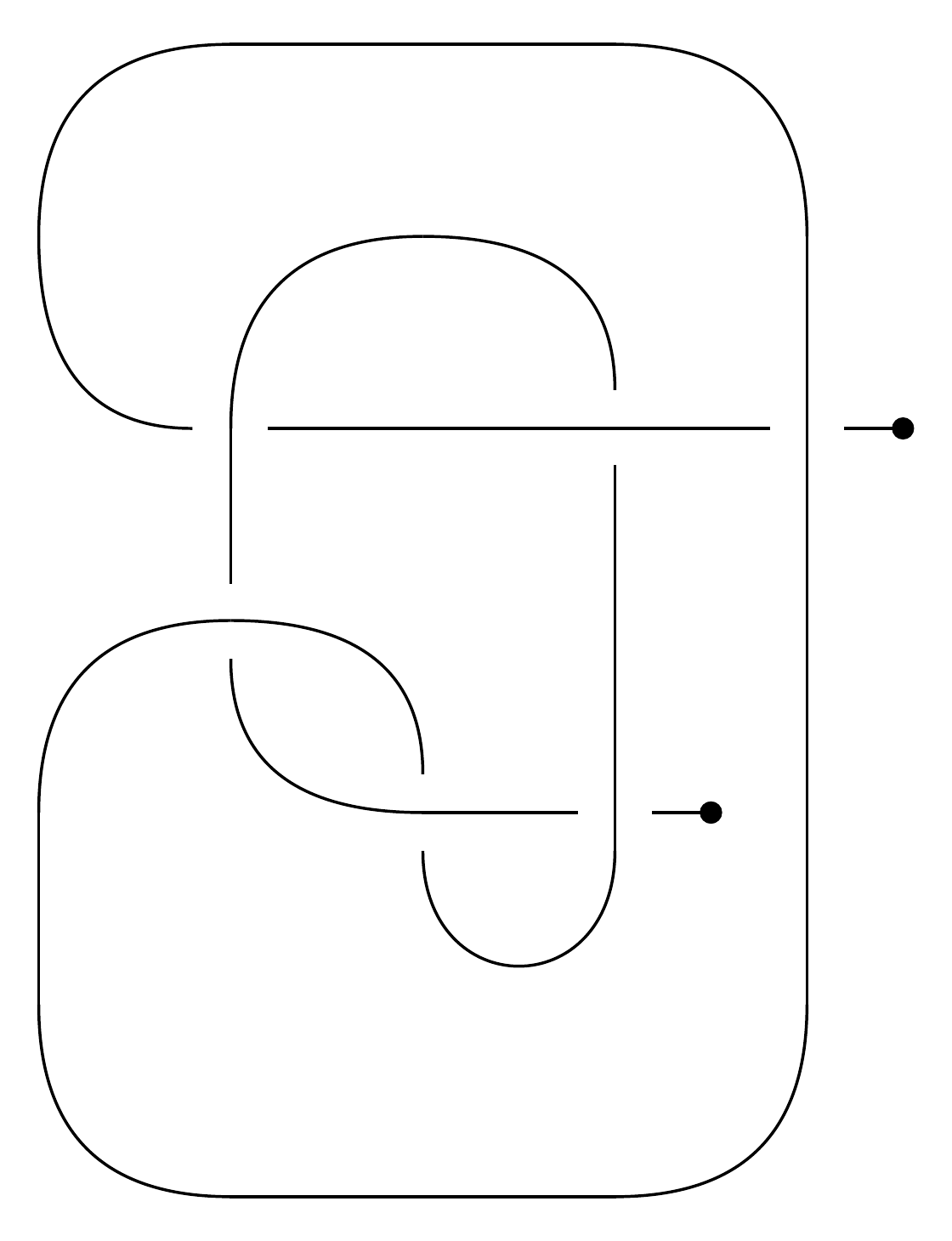}\\
\textcolor{black}{$6_{58}$}
\vspace{1cm}
\end{minipage}
\begin{minipage}[t]{.25\linewidth}
\centering
\includegraphics[width=0.9\textwidth,height=3.5cm,keepaspectratio]{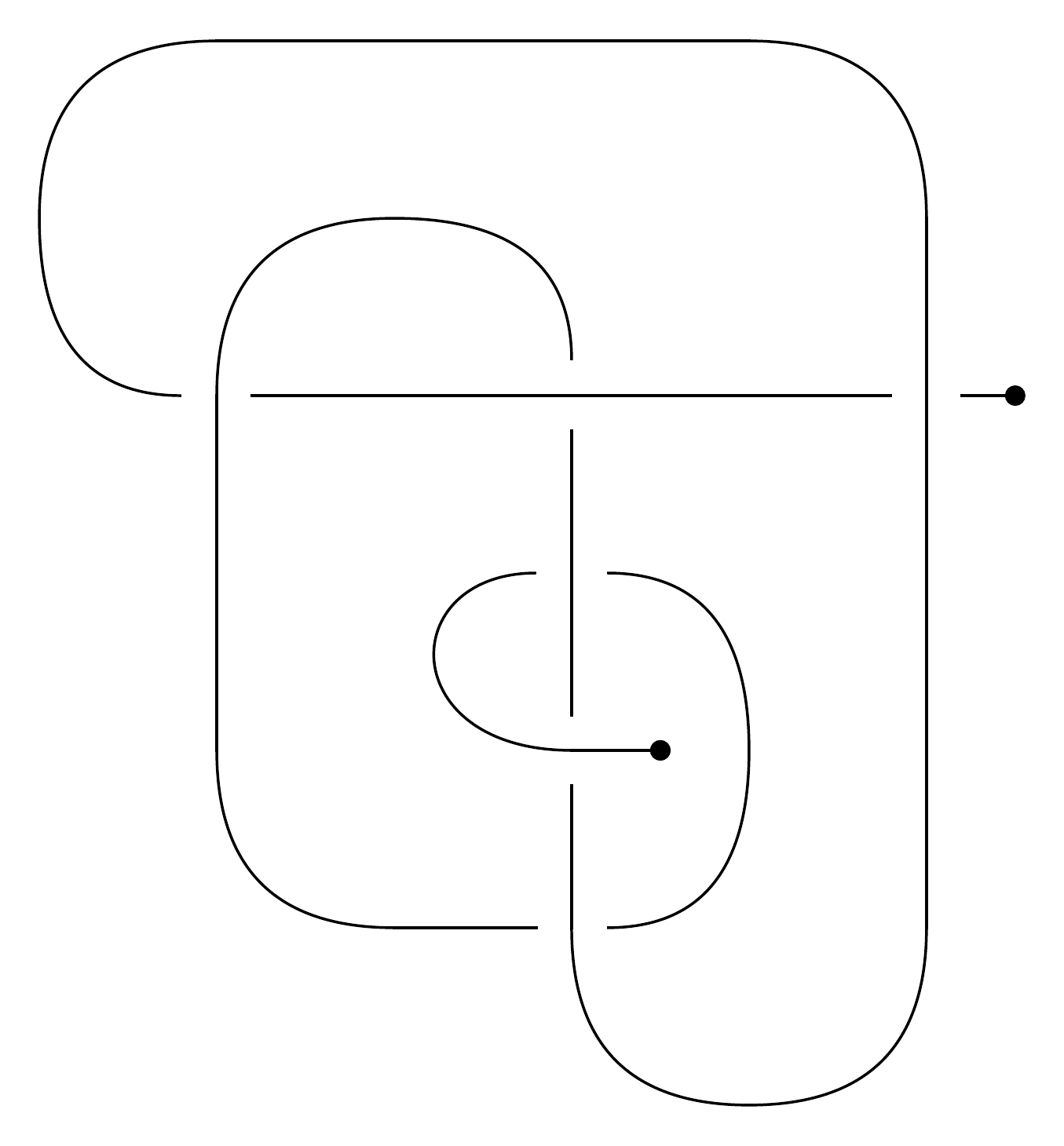}\\
\textcolor{black}{$6_{59}$}
\vspace{1cm}
\end{minipage}
\begin{minipage}[t]{.25\linewidth}
\centering
\includegraphics[width=0.9\textwidth,height=3.5cm,keepaspectratio]{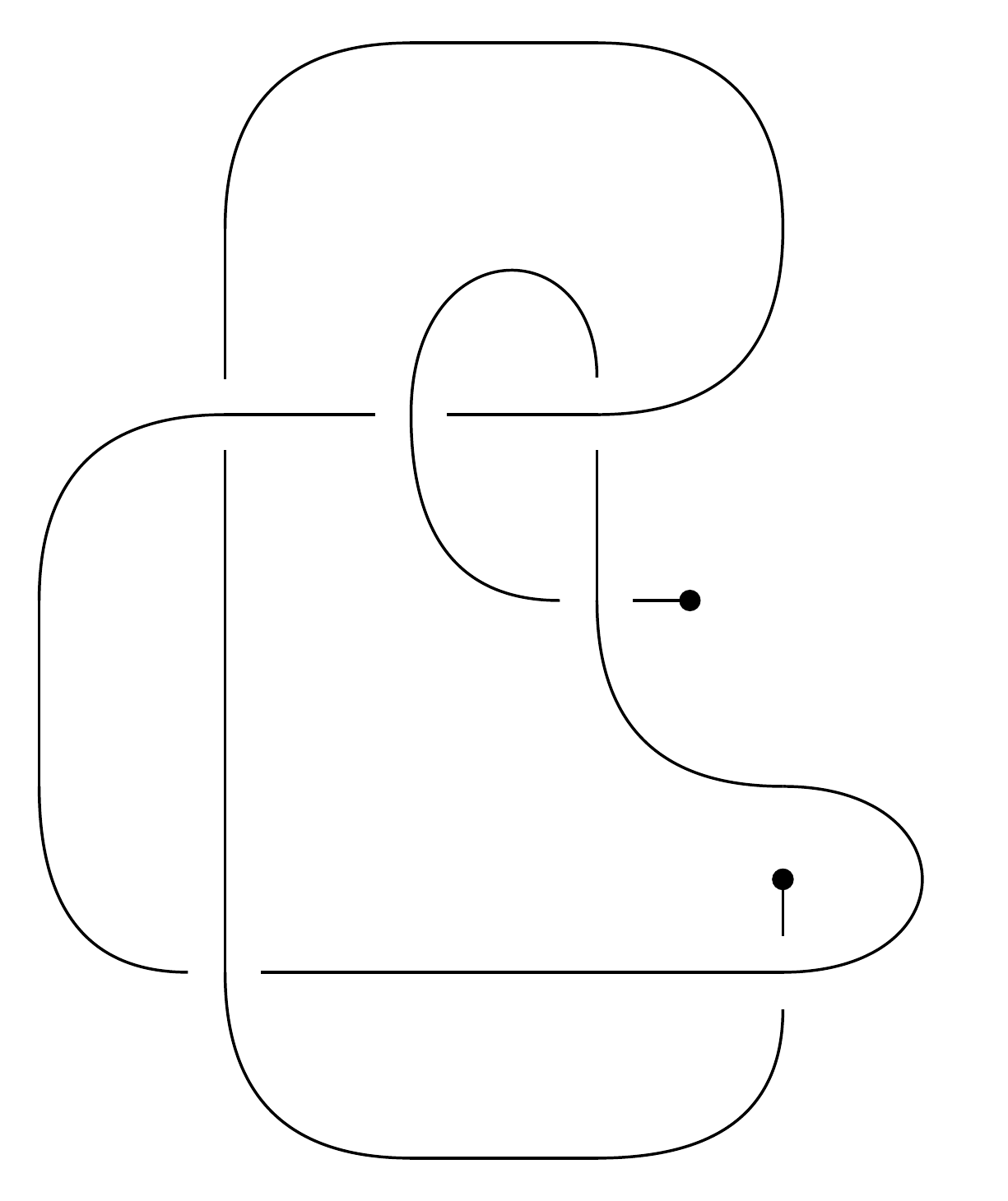}\\
\textcolor{black}{$6_{60}$}
\vspace{1cm}
\end{minipage}
\begin{minipage}[t]{.25\linewidth}
\centering
\includegraphics[width=0.9\textwidth,height=3.5cm,keepaspectratio]{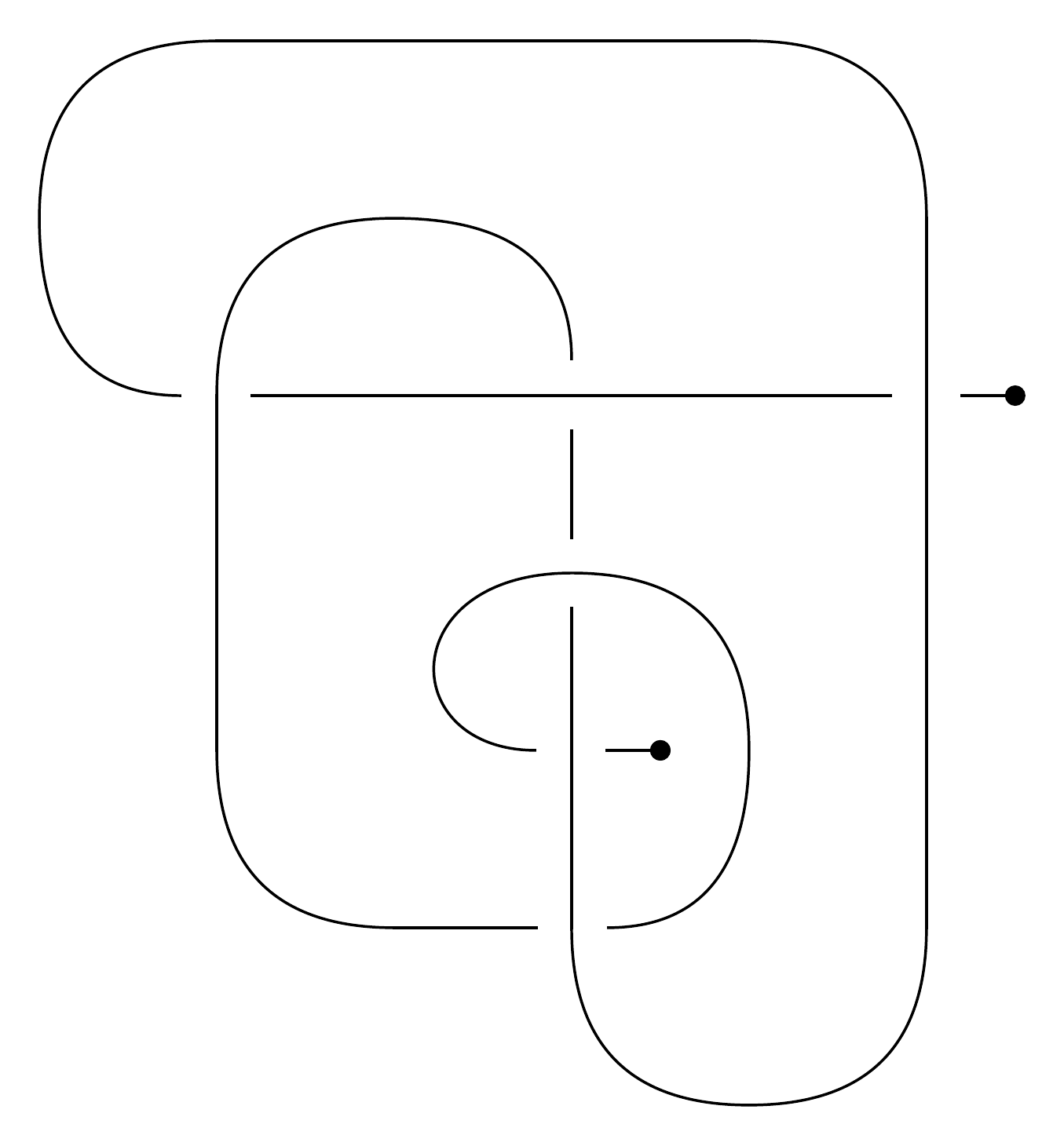}\\
\textcolor{black}{$6_{61}$}
\vspace{1cm}
\end{minipage}
\begin{minipage}[t]{.25\linewidth}
\centering
\includegraphics[width=0.9\textwidth,height=3.5cm,keepaspectratio]{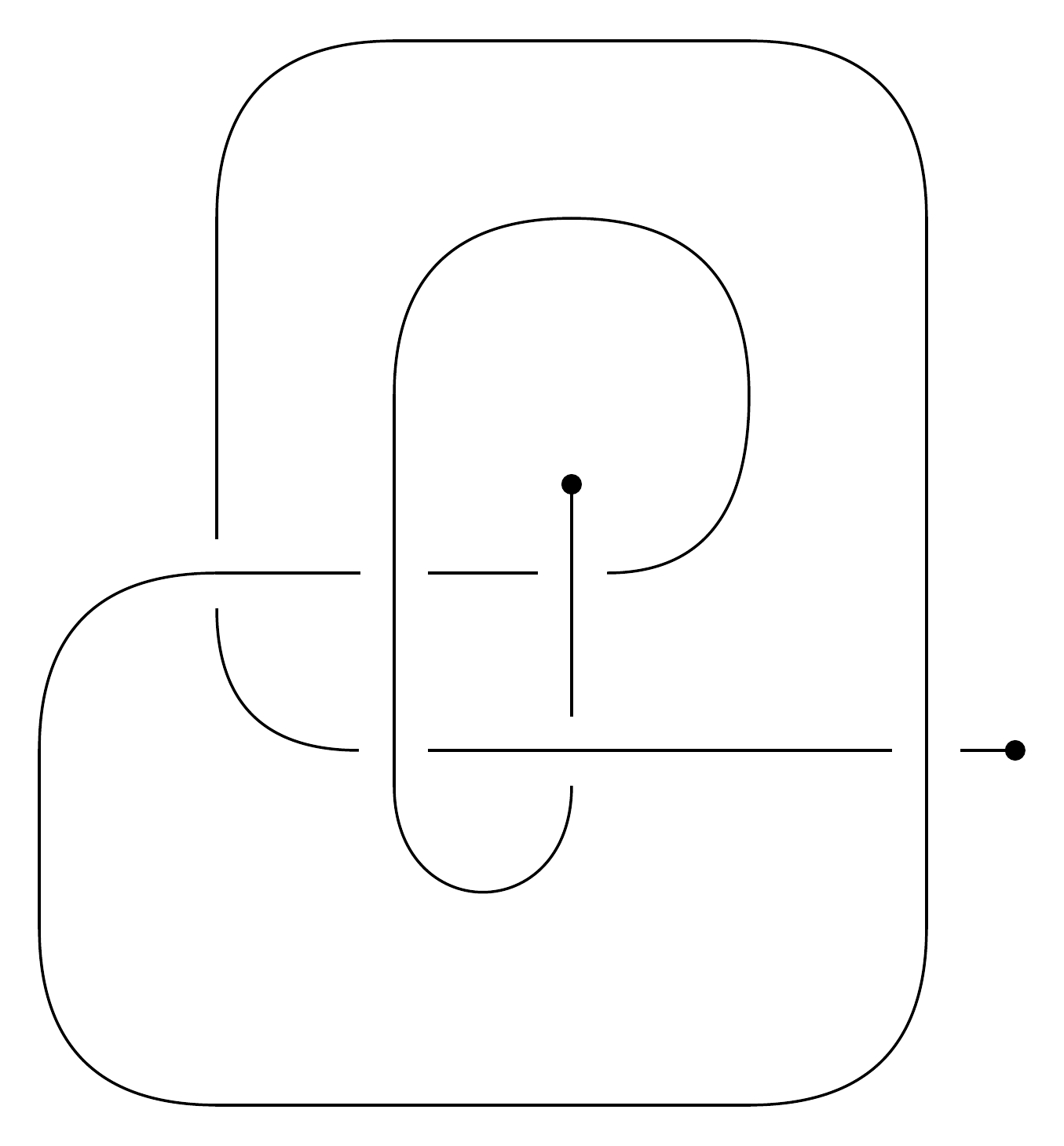}\\
\textcolor{black}{$6_{62}$}
\vspace{1cm}
\end{minipage}
\begin{minipage}[t]{.25\linewidth}
\centering
\includegraphics[width=0.9\textwidth,height=3.5cm,keepaspectratio]{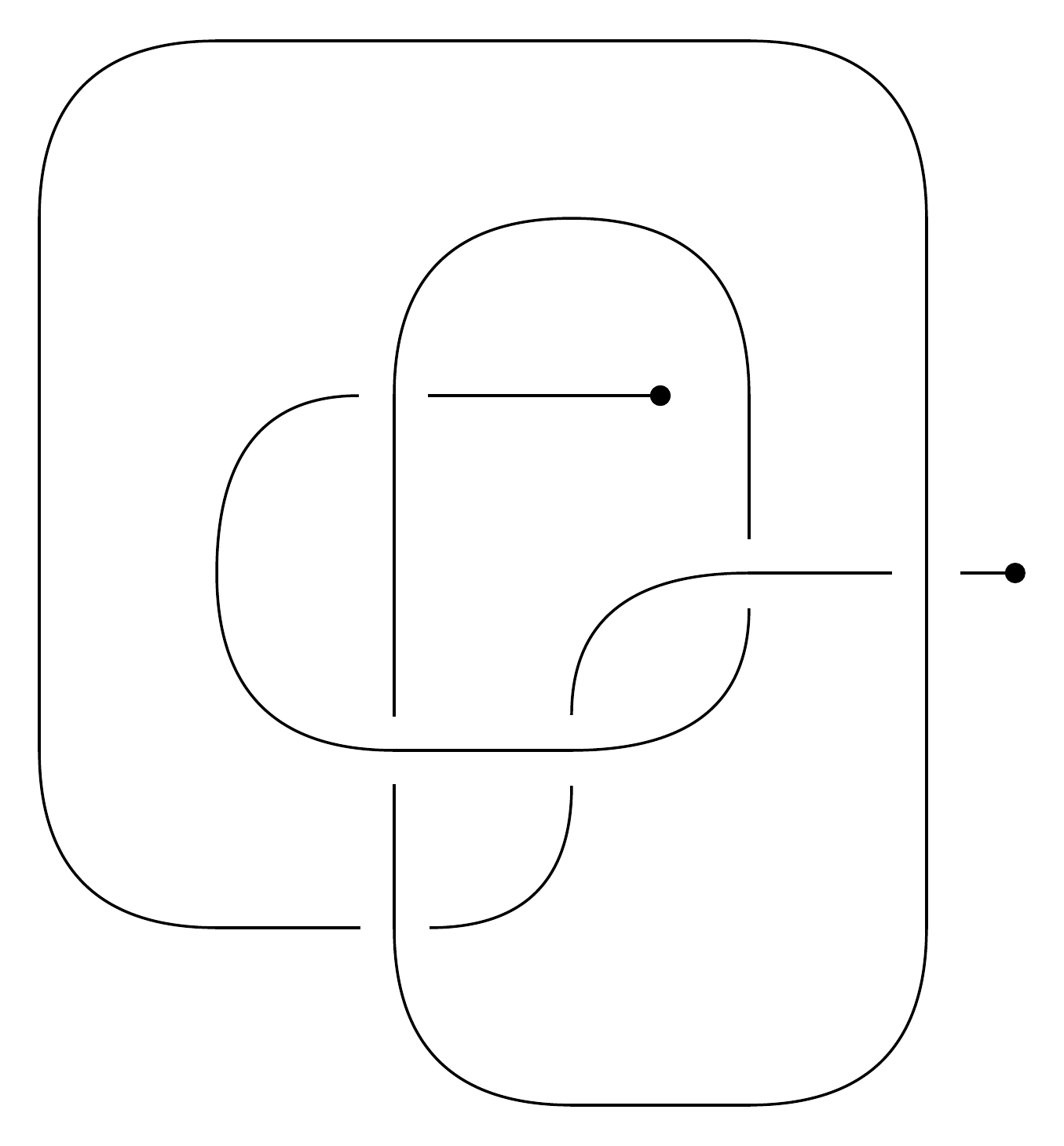}\\
\textcolor{black}{$6_{63}$}
\vspace{1cm}
\end{minipage}
\begin{minipage}[t]{.25\linewidth}
\centering
\includegraphics[width=0.9\textwidth,height=3.5cm,keepaspectratio]{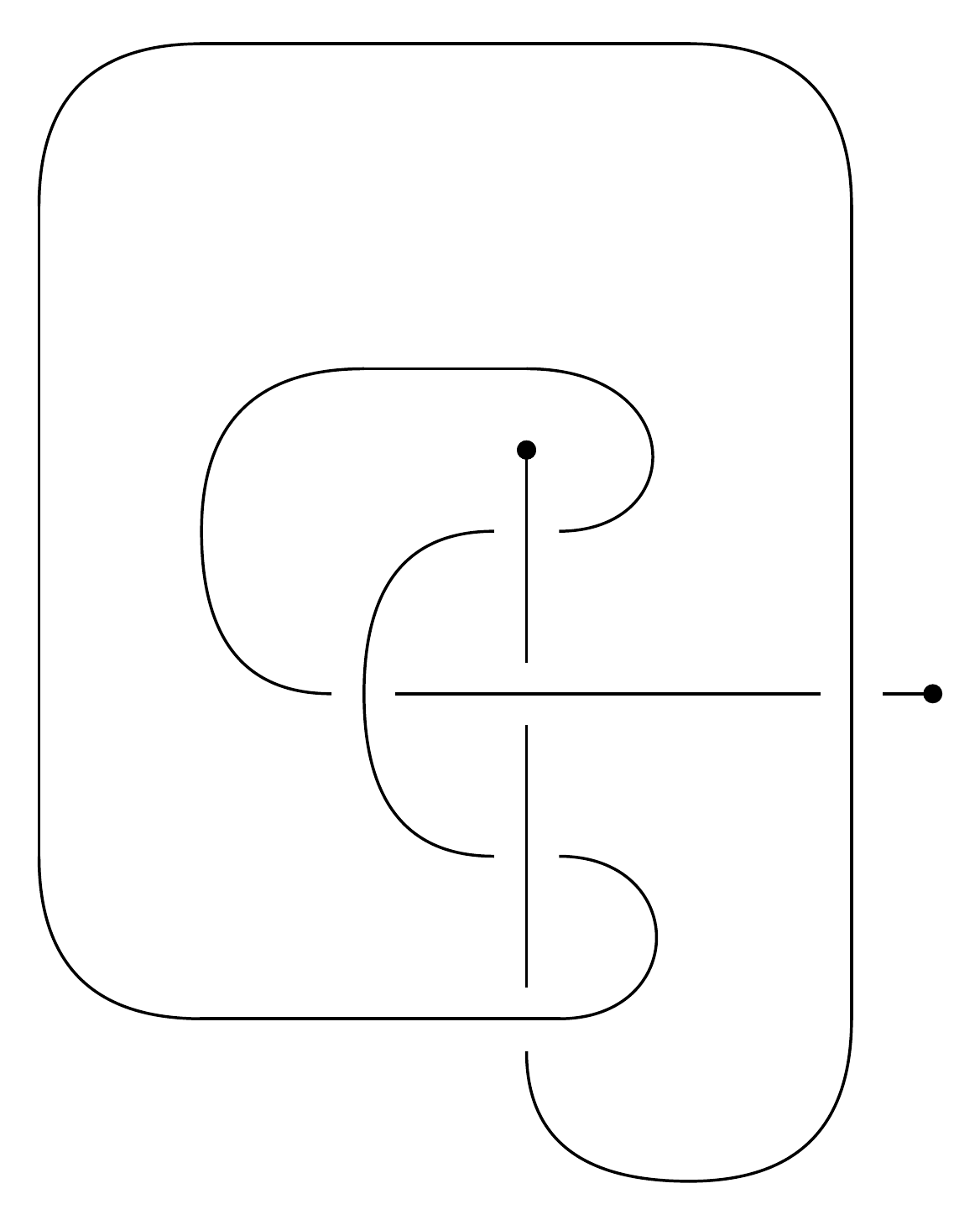}\\
\textcolor{black}{$6_{64}$}
\vspace{1cm}
\end{minipage}
\begin{minipage}[t]{.25\linewidth}
\centering
\includegraphics[width=0.9\textwidth,height=3.5cm,keepaspectratio]{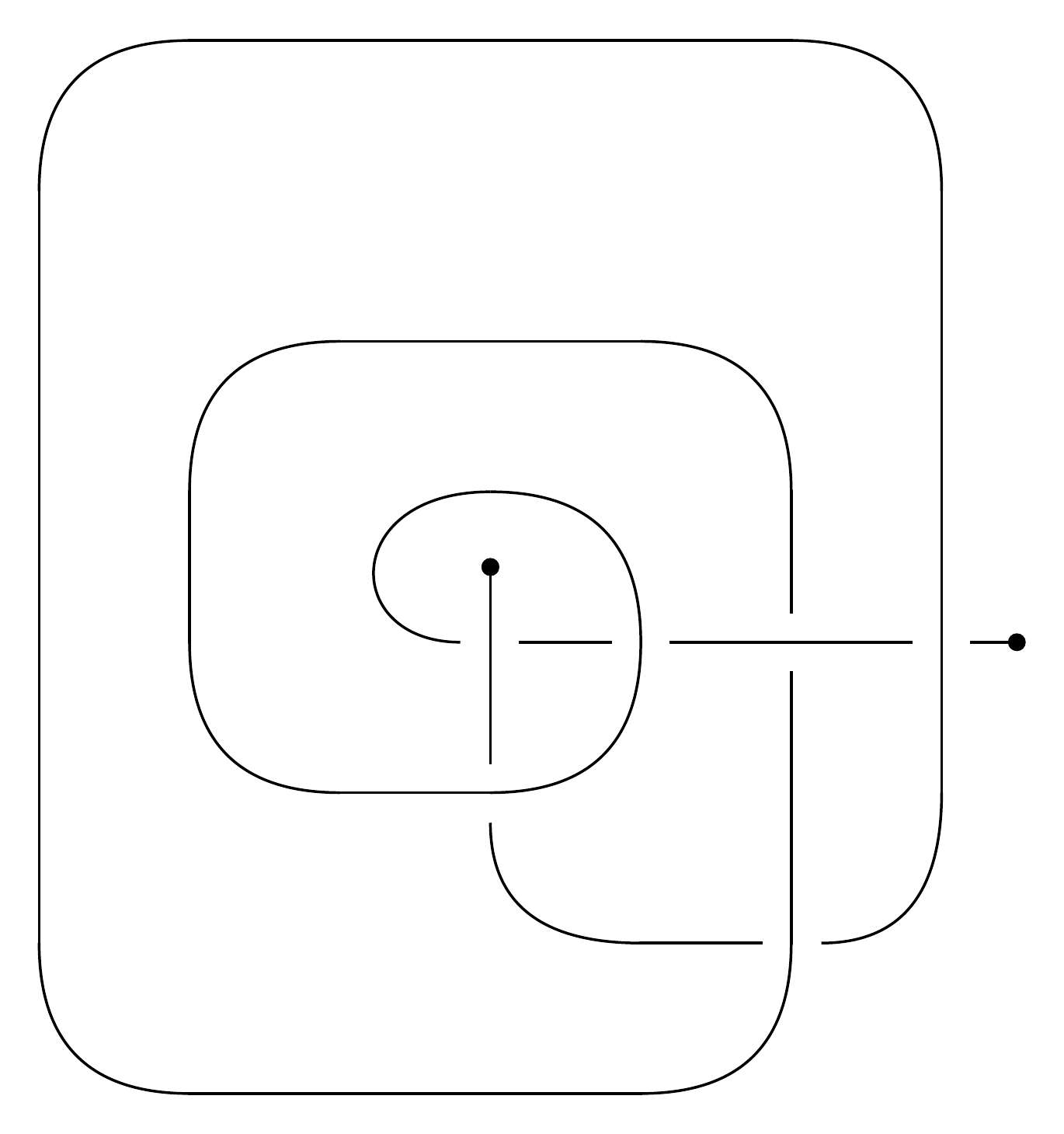}\\
\textcolor{black}{$6_{65}$}
\vspace{1cm}
\end{minipage}
\begin{minipage}[t]{.25\linewidth}
\centering
\includegraphics[width=0.9\textwidth,height=3.5cm,keepaspectratio]{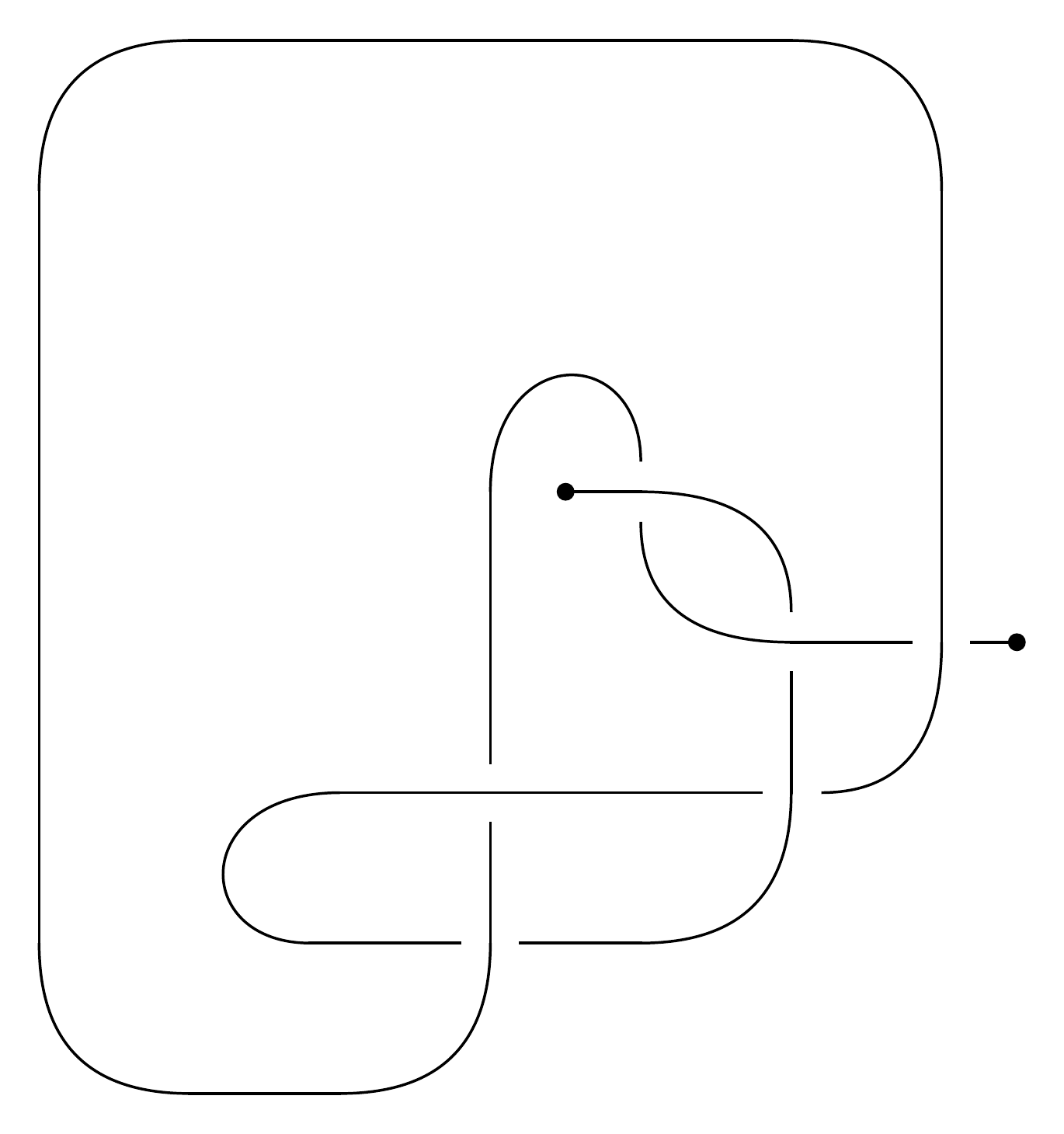}\\
\textcolor{black}{$6_{66}$}
\vspace{1cm}
\end{minipage}
\begin{minipage}[t]{.25\linewidth}
\centering
\includegraphics[width=0.9\textwidth,height=3.5cm,keepaspectratio]{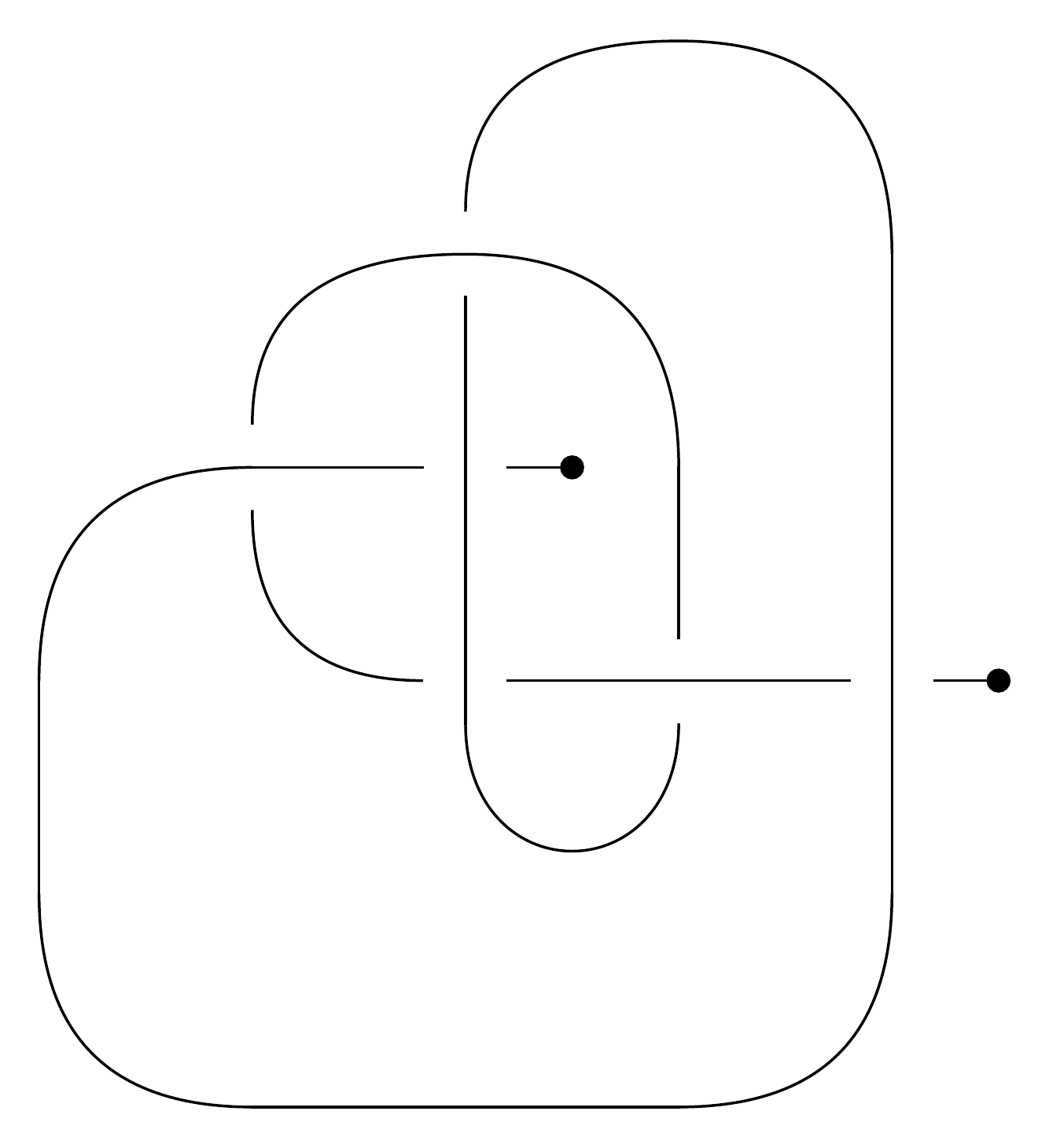}\\
\textcolor{black}{$6_{67}$}
\vspace{1cm}
\end{minipage}
\begin{minipage}[t]{.25\linewidth}
\centering
\includegraphics[width=0.9\textwidth,height=3.5cm,keepaspectratio]{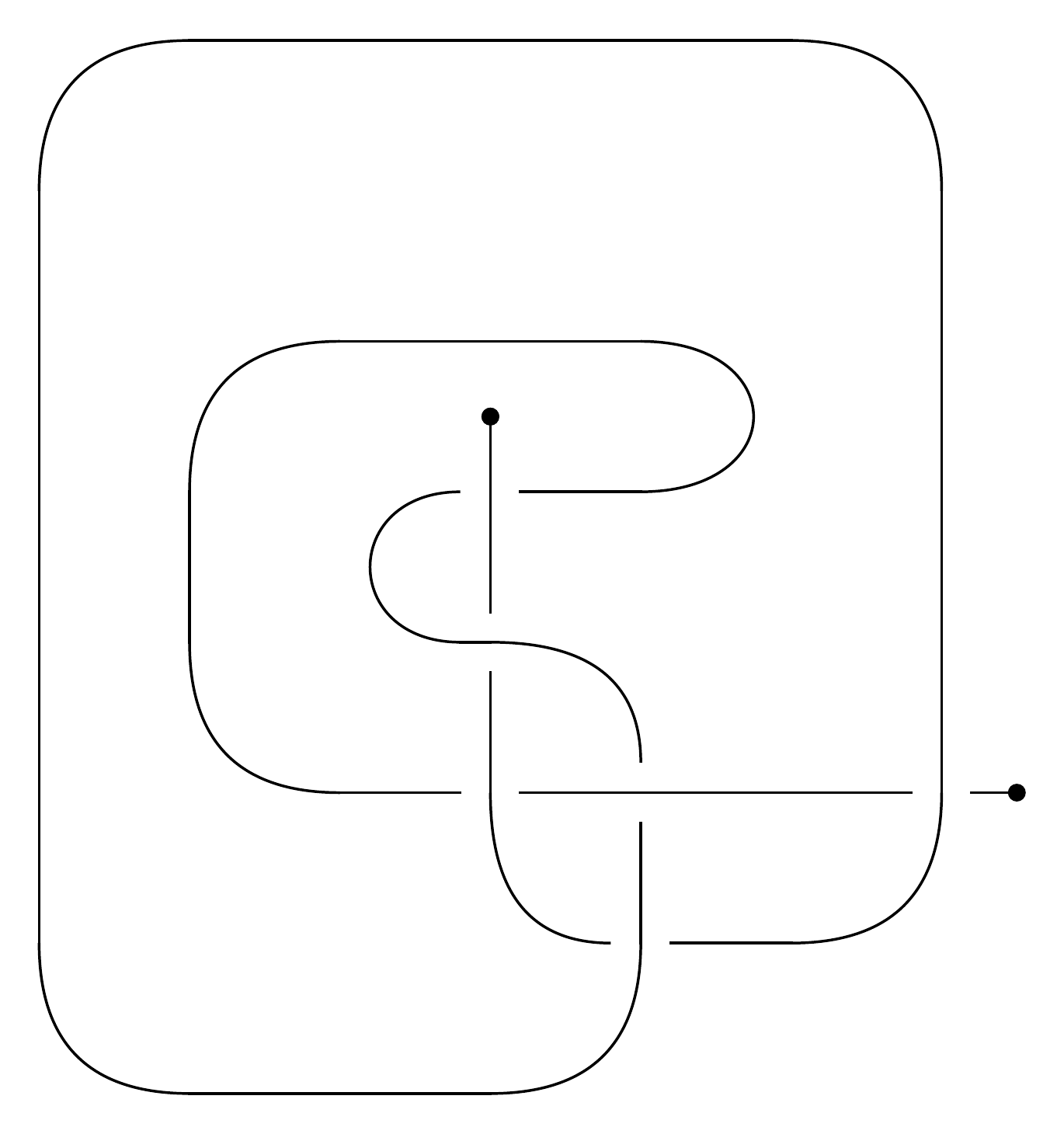}\\
\textcolor{black}{$6_{68}$}
\vspace{1cm}
\end{minipage}
\begin{minipage}[t]{.25\linewidth}
\centering
\includegraphics[width=0.9\textwidth,height=3.5cm,keepaspectratio]{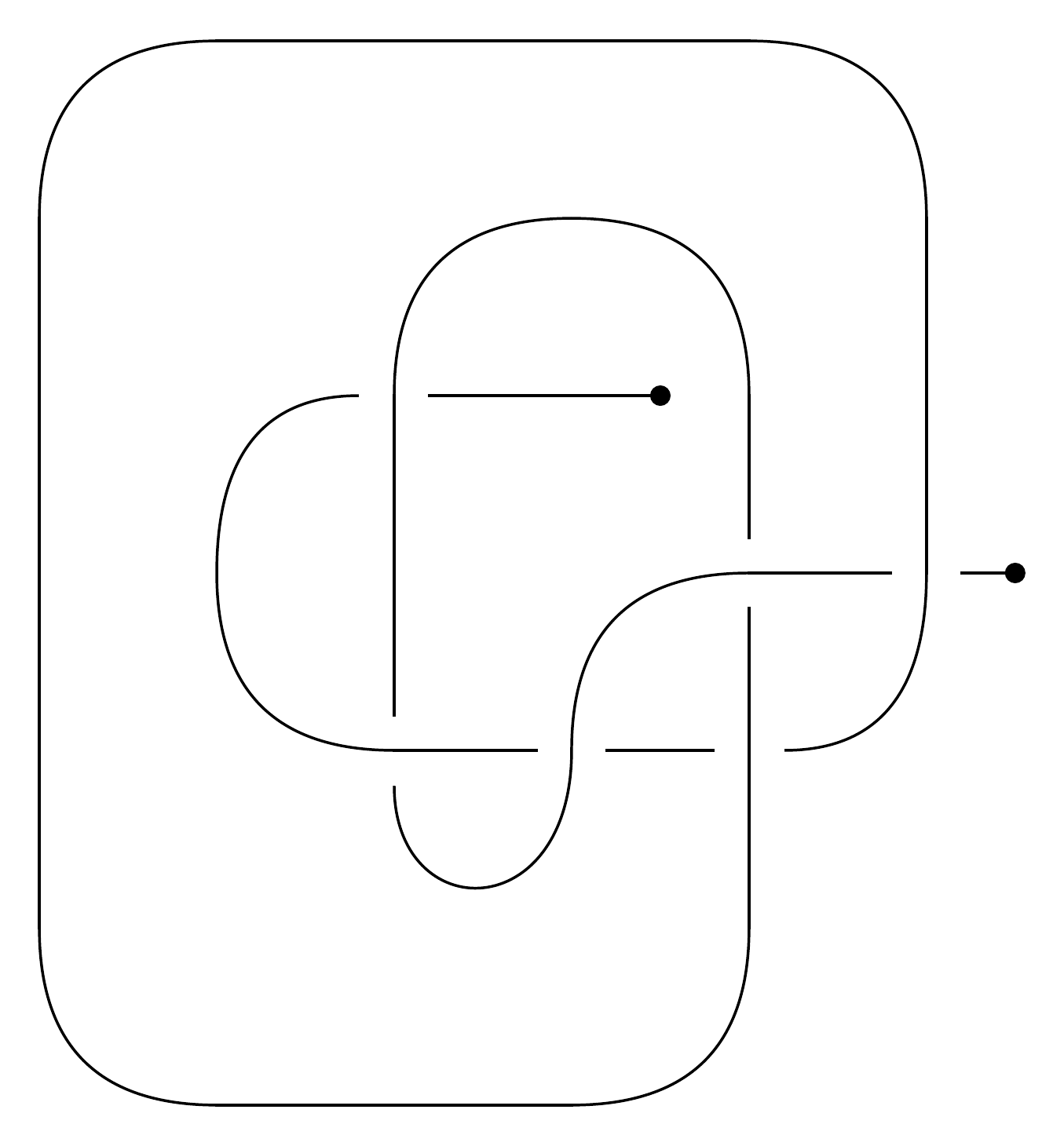}\\
\textcolor{black}{$6_{69}$}
\vspace{1cm}
\end{minipage}
\begin{minipage}[t]{.25\linewidth}
\centering
\includegraphics[width=0.9\textwidth,height=3.5cm,keepaspectratio]{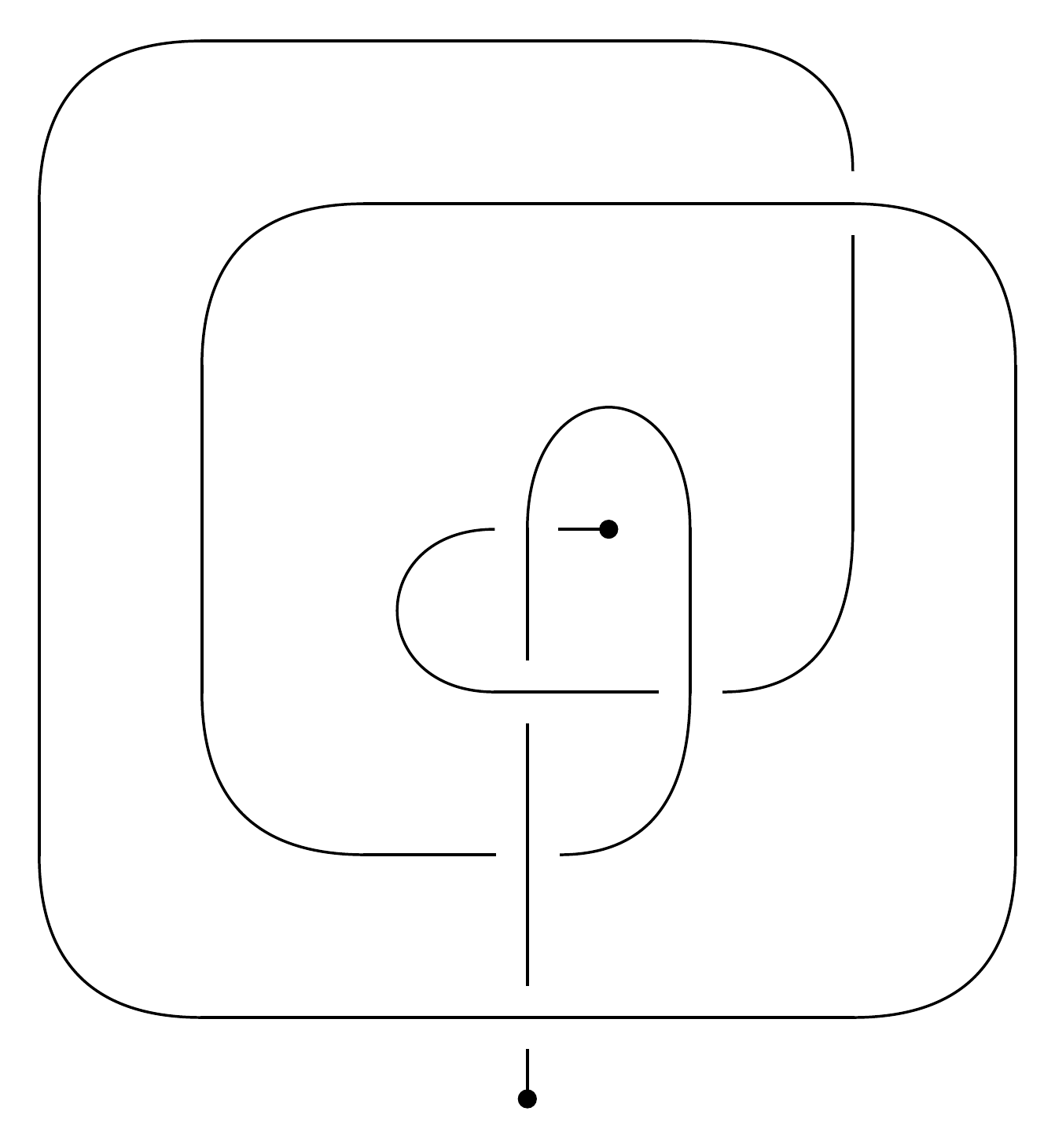}\\
\textcolor{black}{$6_{70}$}
\vspace{1cm}
\end{minipage}
\begin{minipage}[t]{.25\linewidth}
\centering
\includegraphics[width=0.9\textwidth,height=3.5cm,keepaspectratio]{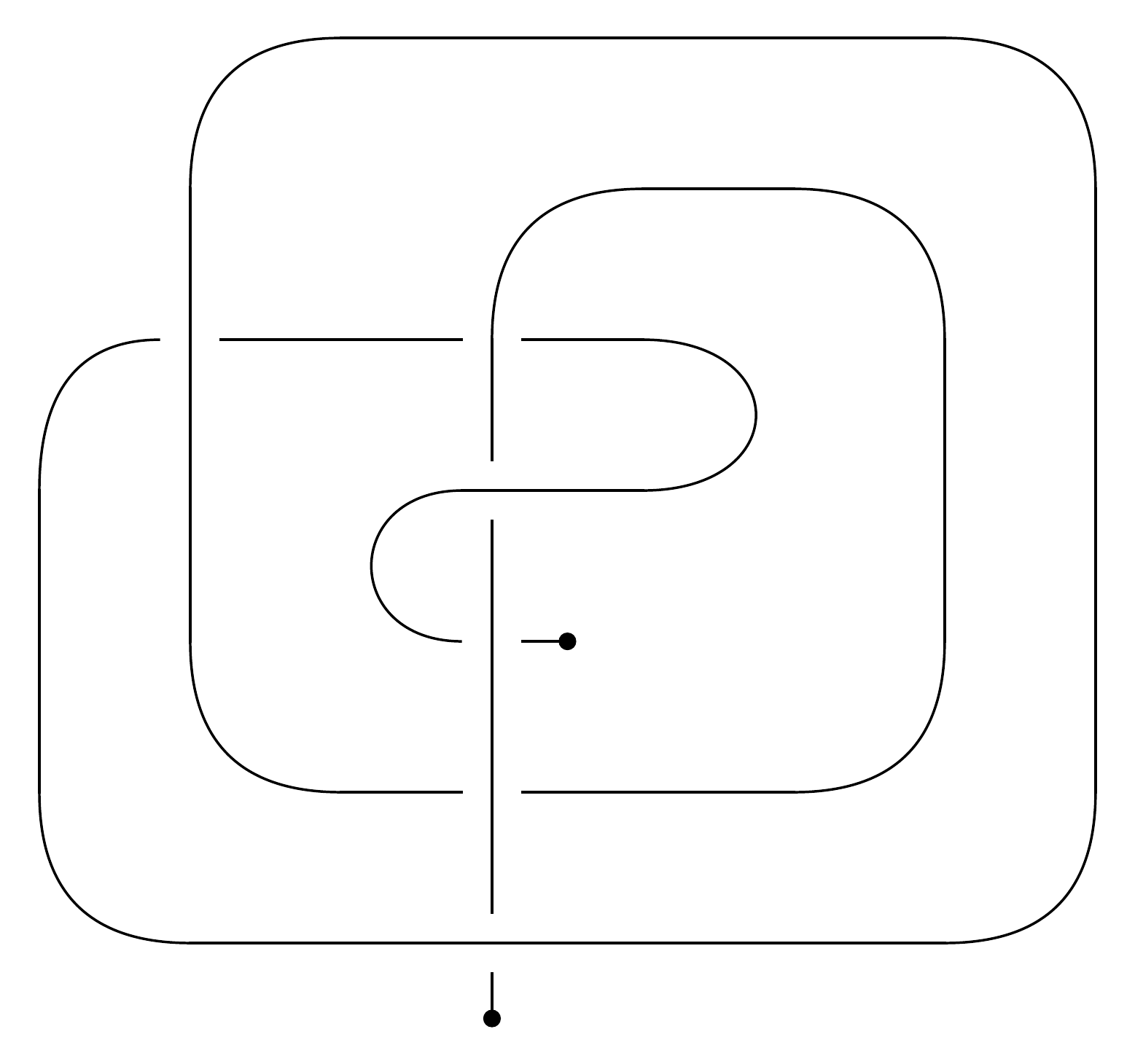}\\
\textcolor{black}{$6_{71}$}
\vspace{1cm}
\end{minipage}
\begin{minipage}[t]{.25\linewidth}
\centering
\includegraphics[width=0.9\textwidth,height=3.5cm,keepaspectratio]{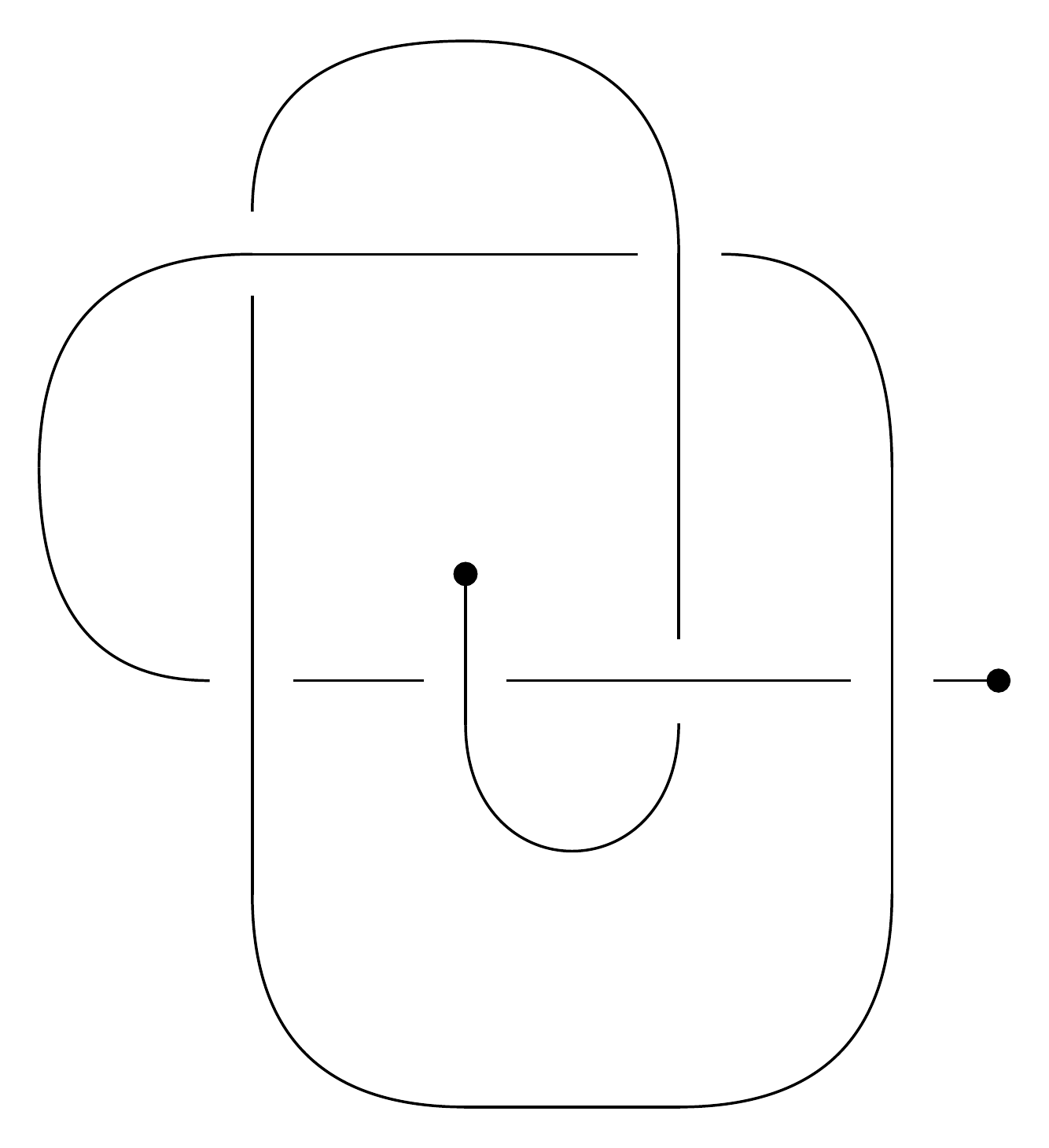}\\
\textcolor{black}{$6_{72}$}
\vspace{1cm}
\end{minipage}
\begin{minipage}[t]{.25\linewidth}
\centering
\includegraphics[width=0.9\textwidth,height=3.5cm,keepaspectratio]{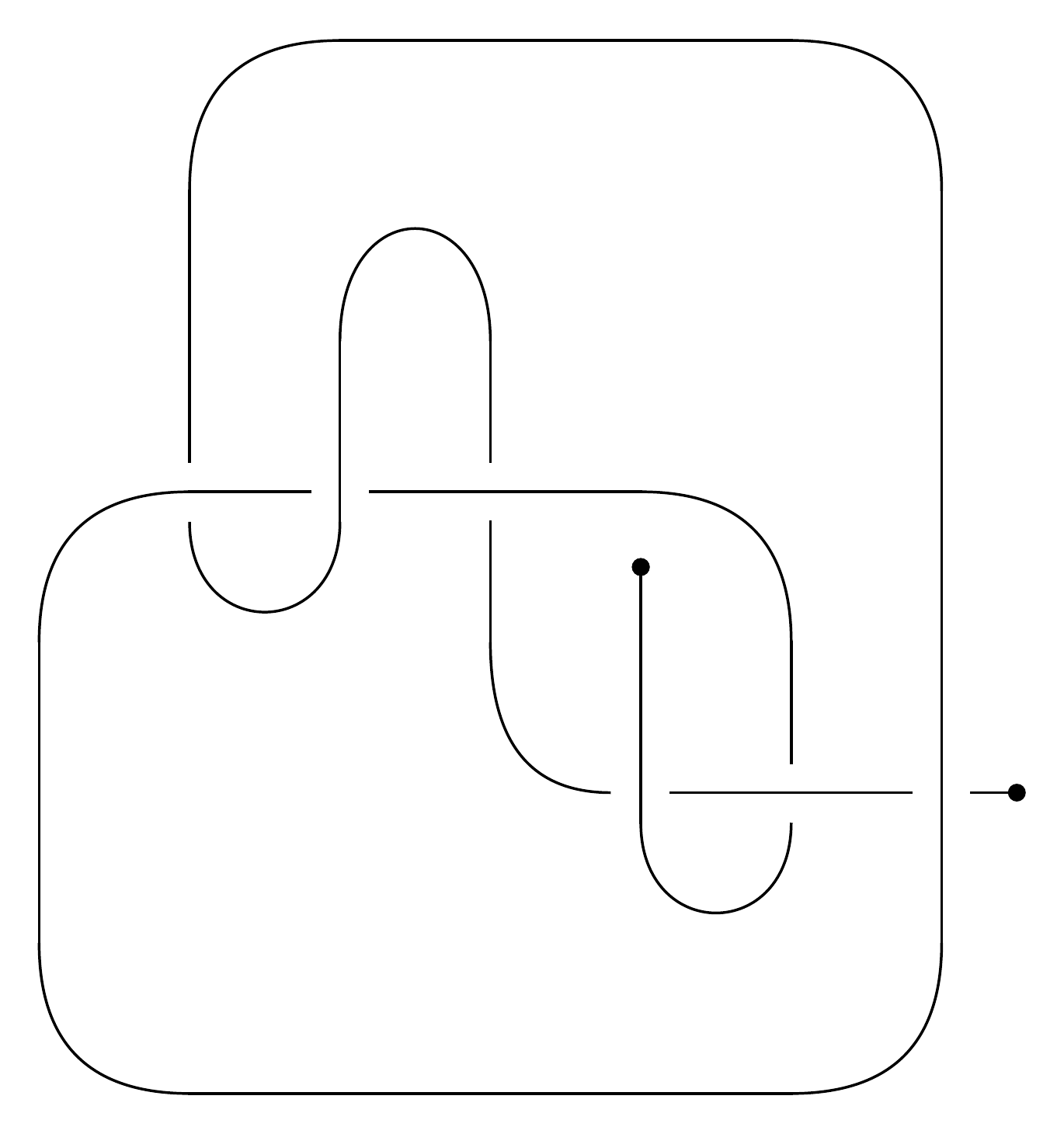}\\
\textcolor{black}{$6_{73}$}
\vspace{1cm}
\end{minipage}
\begin{minipage}[t]{.25\linewidth}
\centering
\includegraphics[width=0.9\textwidth,height=3.5cm,keepaspectratio]{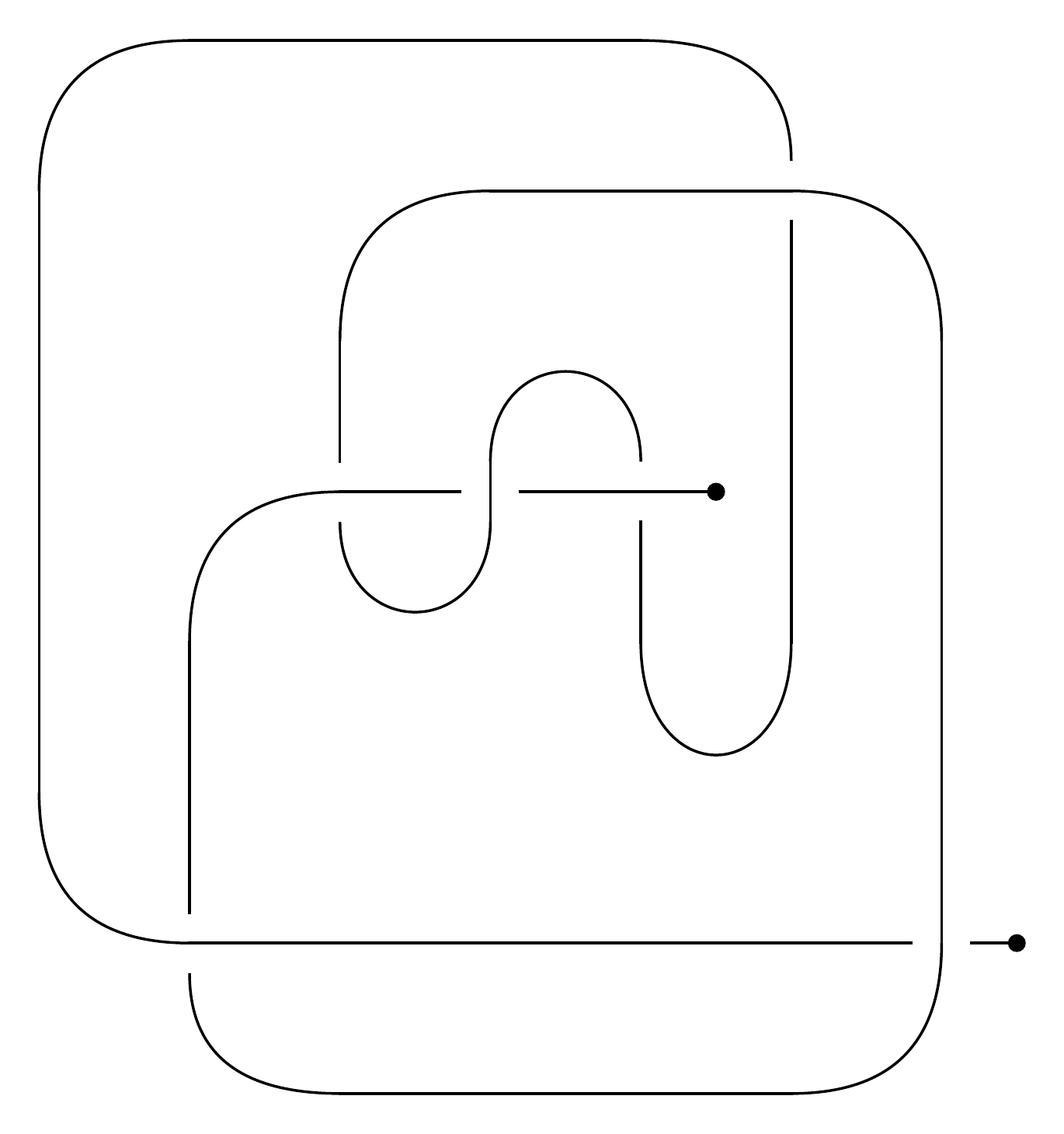}\\
\textcolor{black}{$6_{74}$}
\vspace{1cm}
\end{minipage}
\begin{minipage}[t]{.25\linewidth}
\centering
\includegraphics[width=0.9\textwidth,height=3.5cm,keepaspectratio]{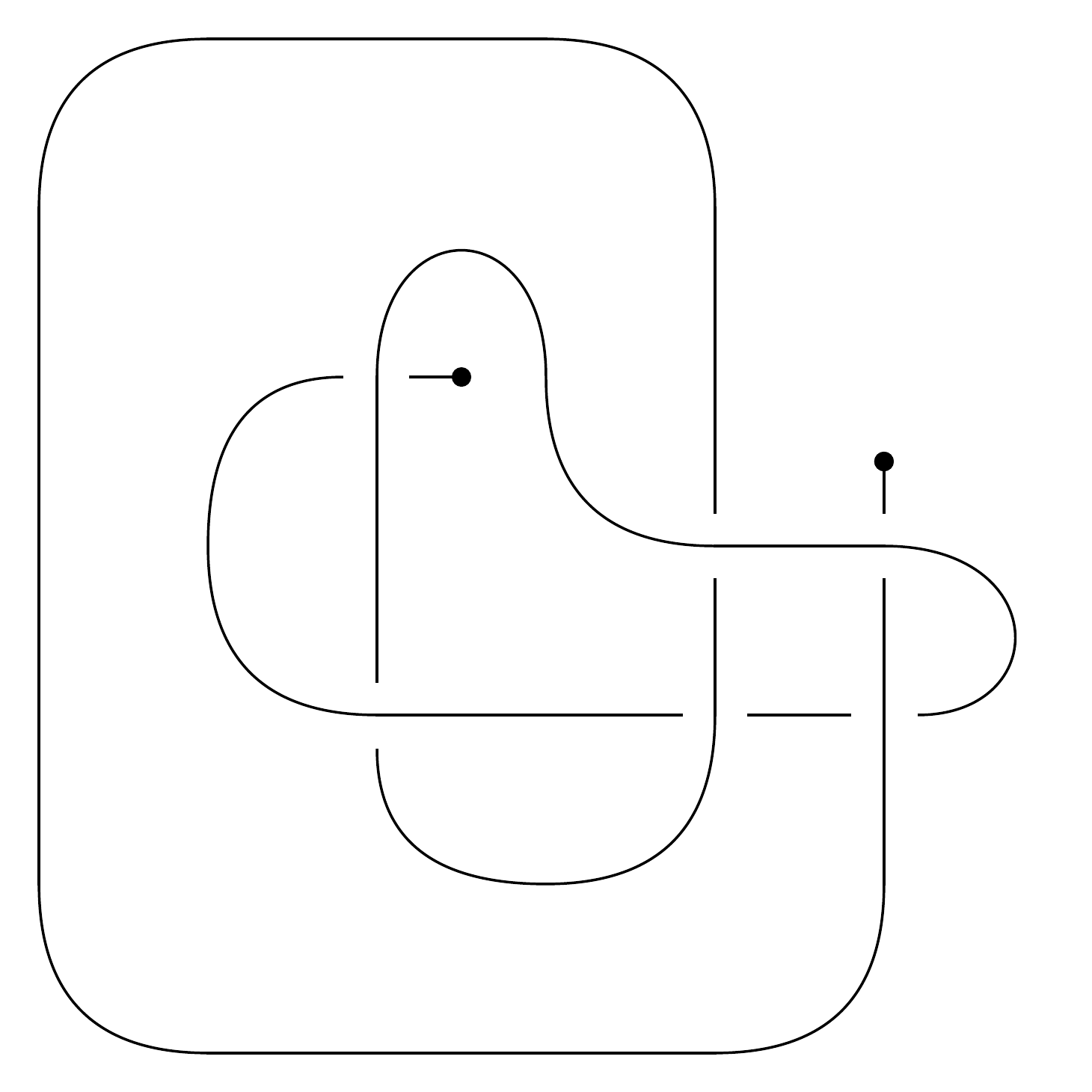}\\
\textcolor{black}{$6_{75}$}
\vspace{1cm}
\end{minipage}
\begin{minipage}[t]{.25\linewidth}
\centering
\includegraphics[width=0.9\textwidth,height=3.5cm,keepaspectratio]{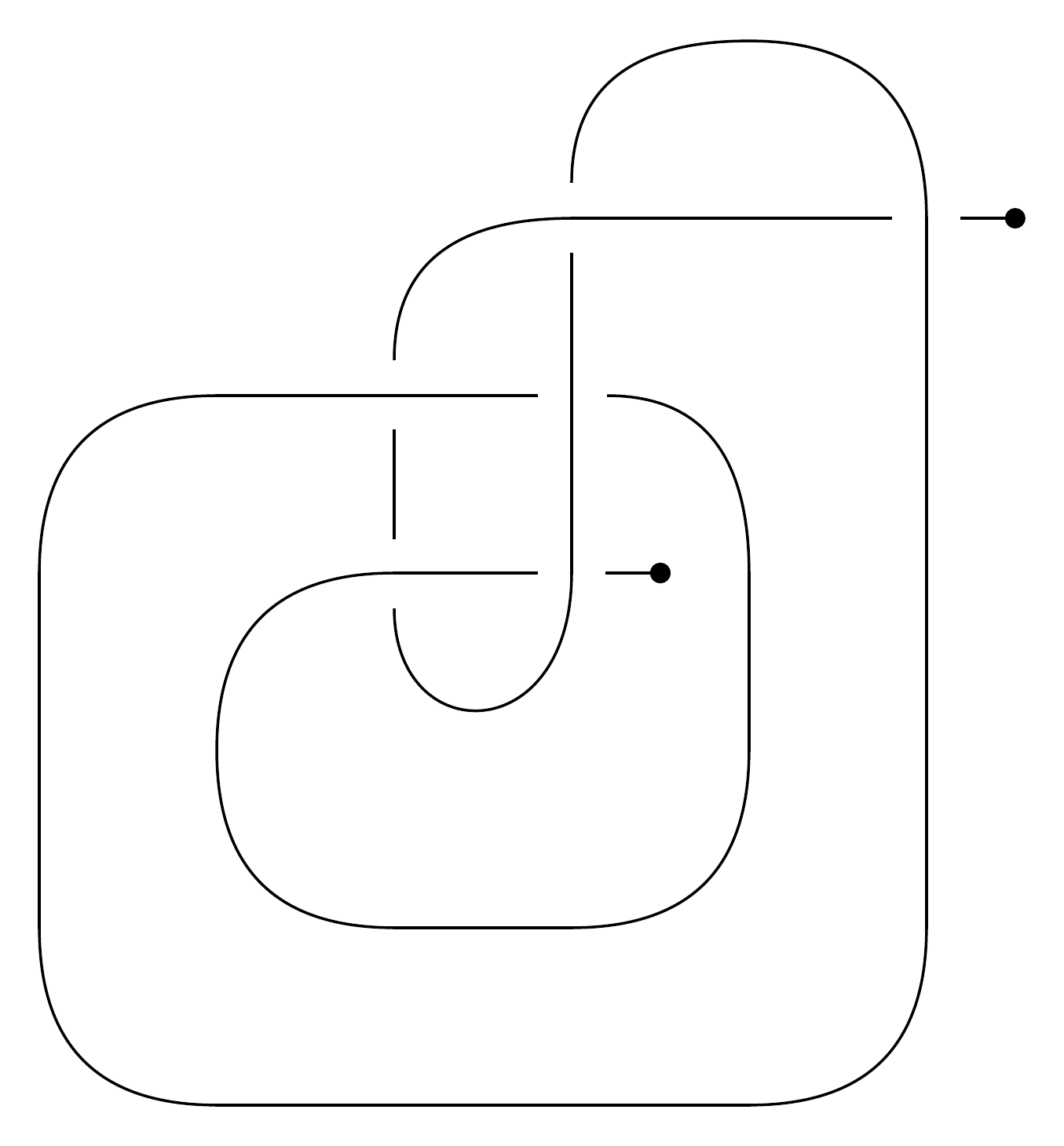}\\
\textcolor{black}{$6_{76}$}
\vspace{1cm}
\end{minipage}
\begin{minipage}[t]{.25\linewidth}
\centering
\includegraphics[width=0.9\textwidth,height=3.5cm,keepaspectratio]{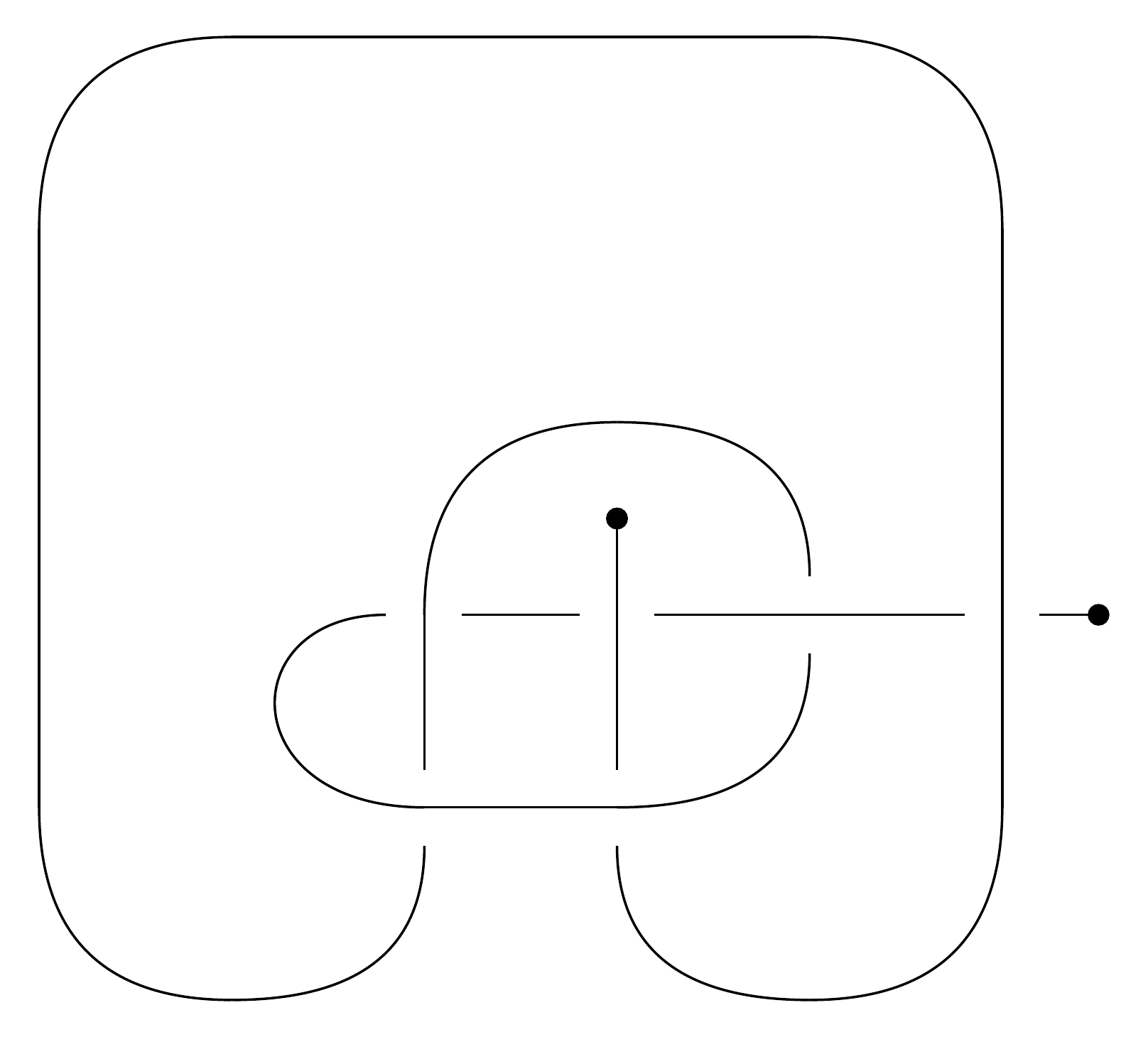}\\
\textcolor{black}{$6_{77}$}
\vspace{1cm}
\end{minipage}
\begin{minipage}[t]{.25\linewidth}
\centering
\includegraphics[width=0.9\textwidth,height=3.5cm,keepaspectratio]{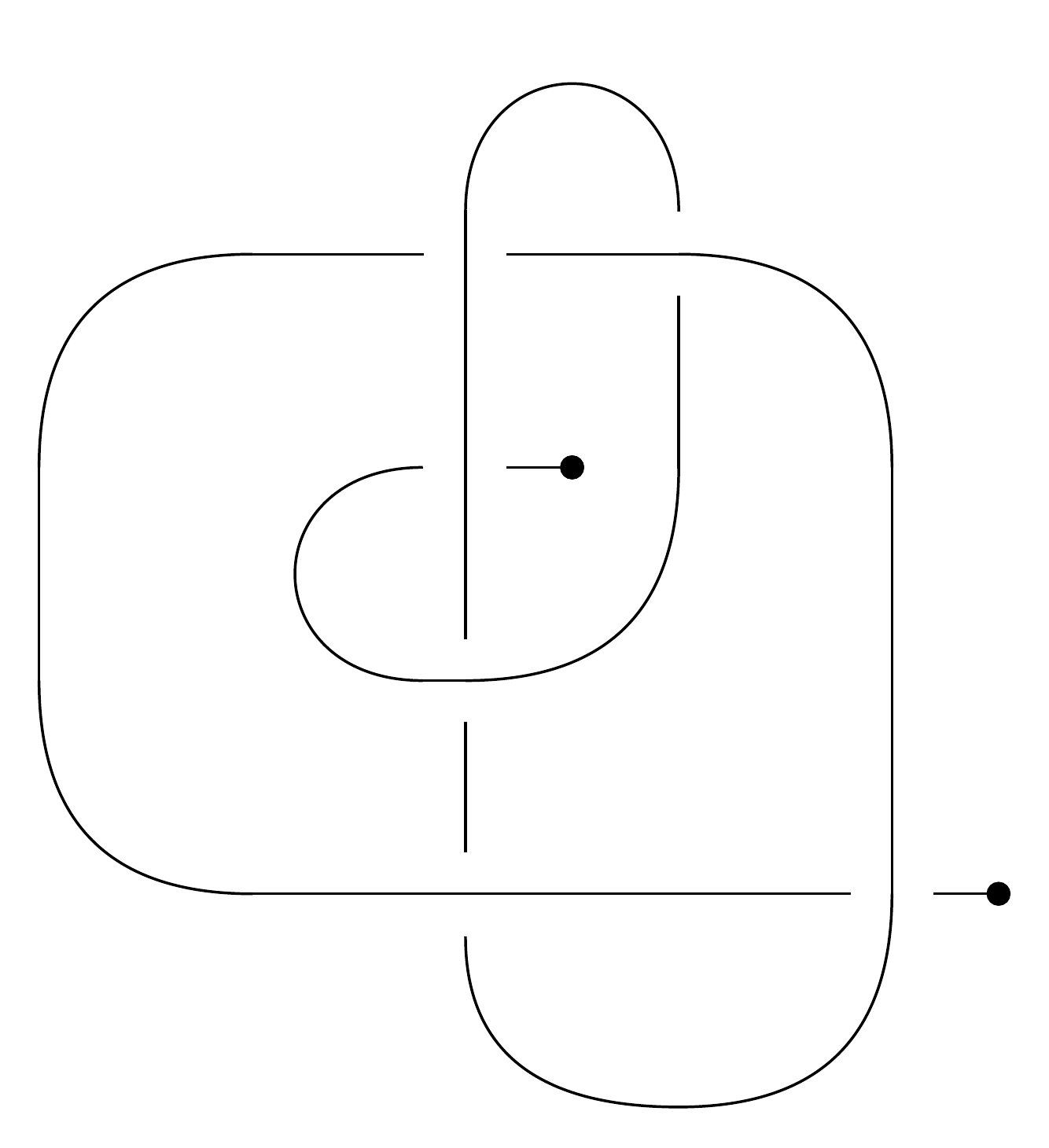}\\
\textcolor{black}{$6_{78}$}
\vspace{1cm}
\end{minipage}
\begin{minipage}[t]{.25\linewidth}
\centering
\includegraphics[width=0.9\textwidth,height=3.5cm,keepaspectratio]{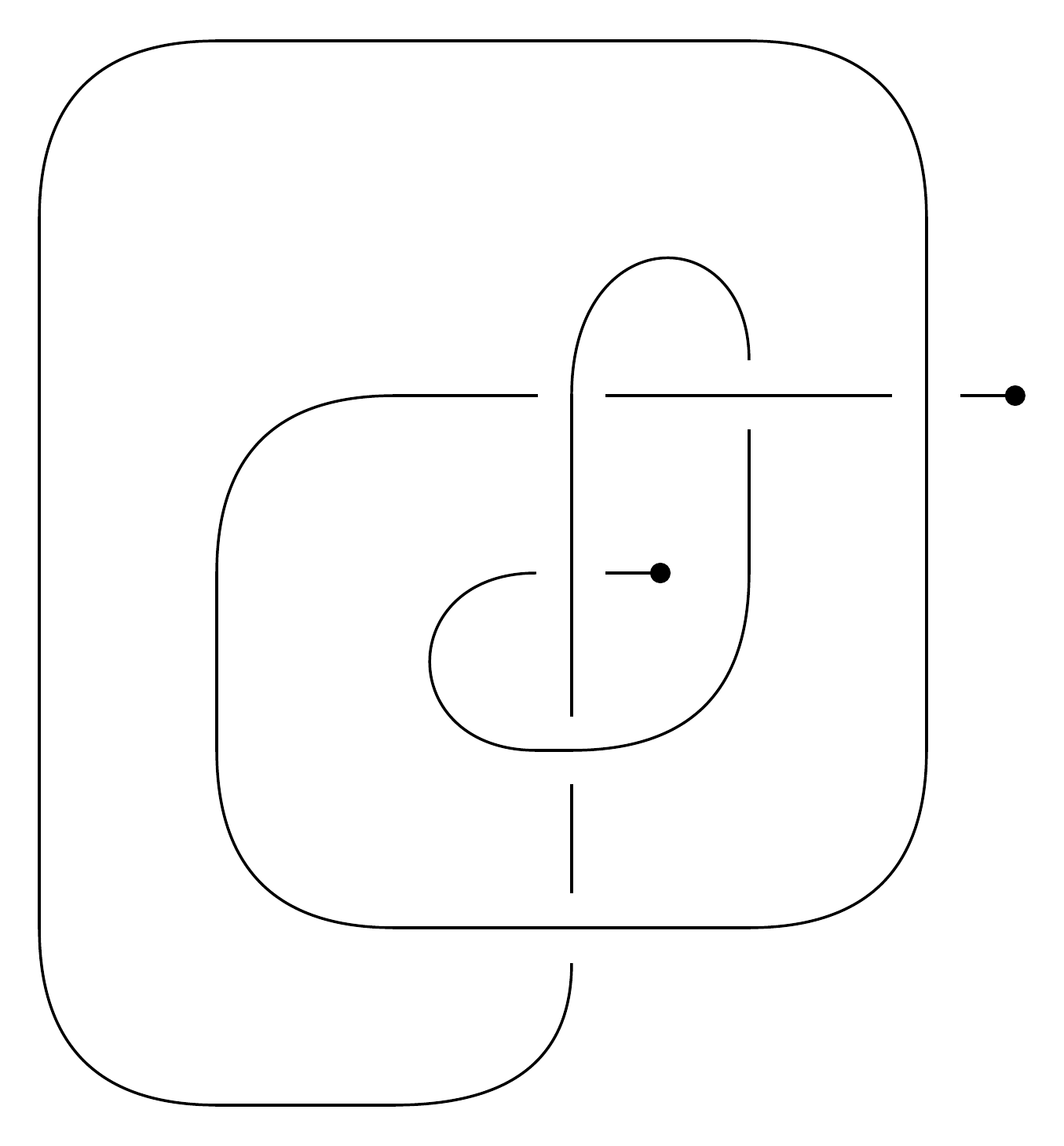}\\
\textcolor{black}{$6_{79}$}
\vspace{1cm}
\end{minipage}
\begin{minipage}[t]{.25\linewidth}
\centering
\includegraphics[width=0.9\textwidth,height=3.5cm,keepaspectratio]{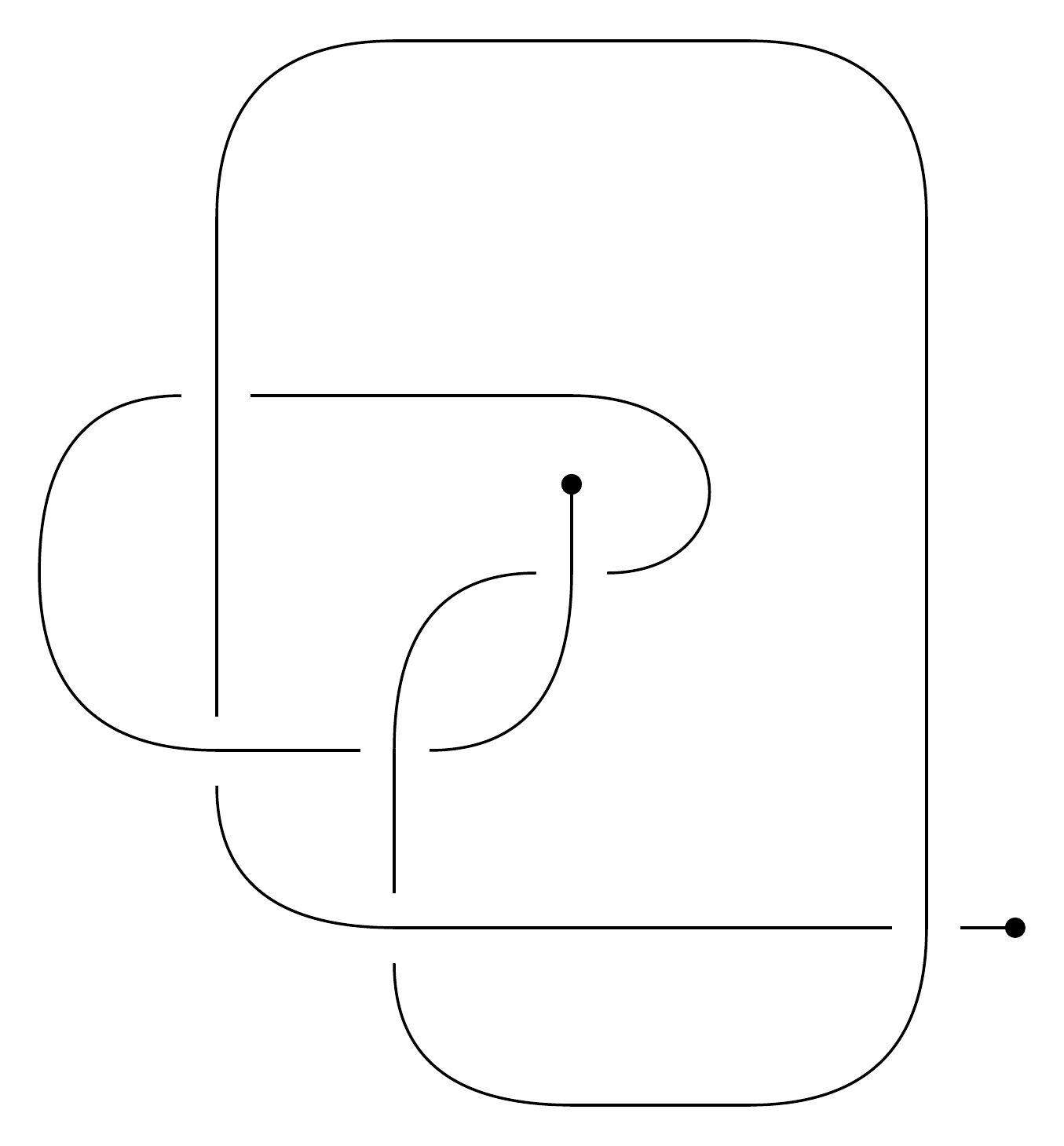}\\
\textcolor{black}{$6_{80}$}
\vspace{1cm}
\end{minipage}
\begin{minipage}[t]{.25\linewidth}
\centering
\includegraphics[width=0.9\textwidth,height=3.5cm,keepaspectratio]{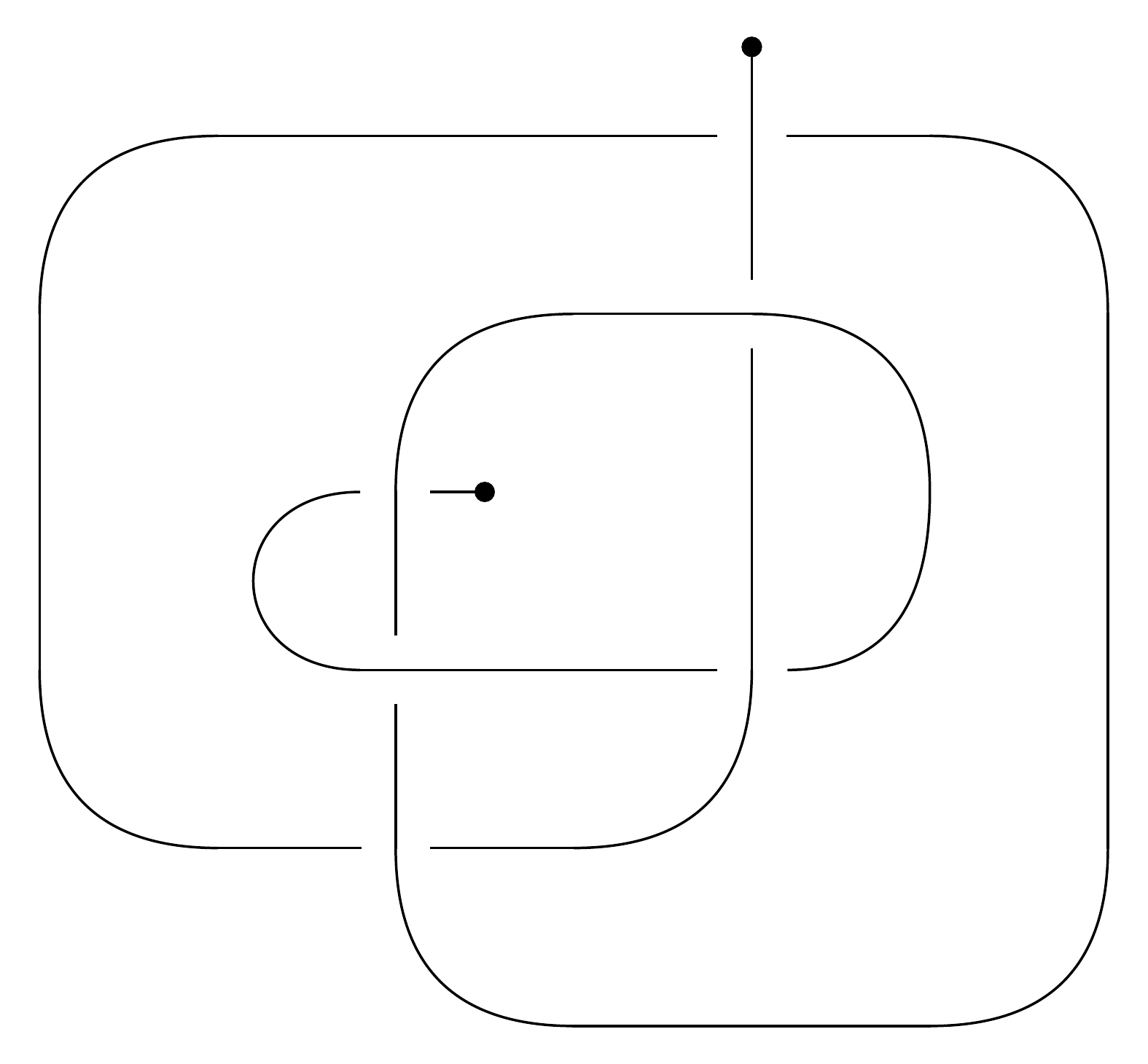}\\
\textcolor{black}{$6_{81}$}
\vspace{1cm}
\end{minipage}
\begin{minipage}[t]{.25\linewidth}
\centering
\includegraphics[width=0.9\textwidth,height=3.5cm,keepaspectratio]{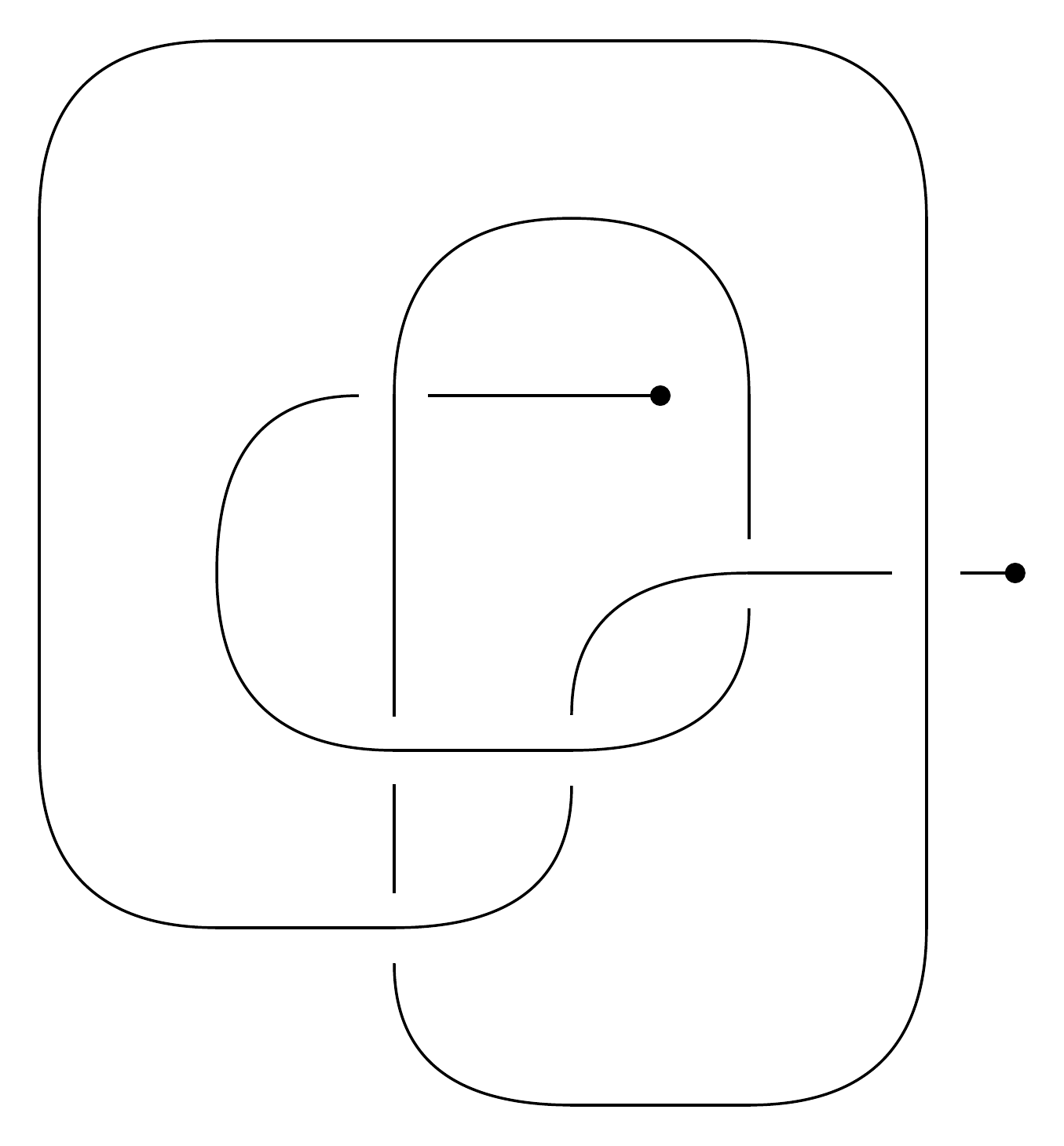}\\
\textcolor{black}{$6_{82}$}
\vspace{1cm}
\end{minipage}
\begin{minipage}[t]{.25\linewidth}
\centering
\includegraphics[width=0.9\textwidth,height=3.5cm,keepaspectratio]{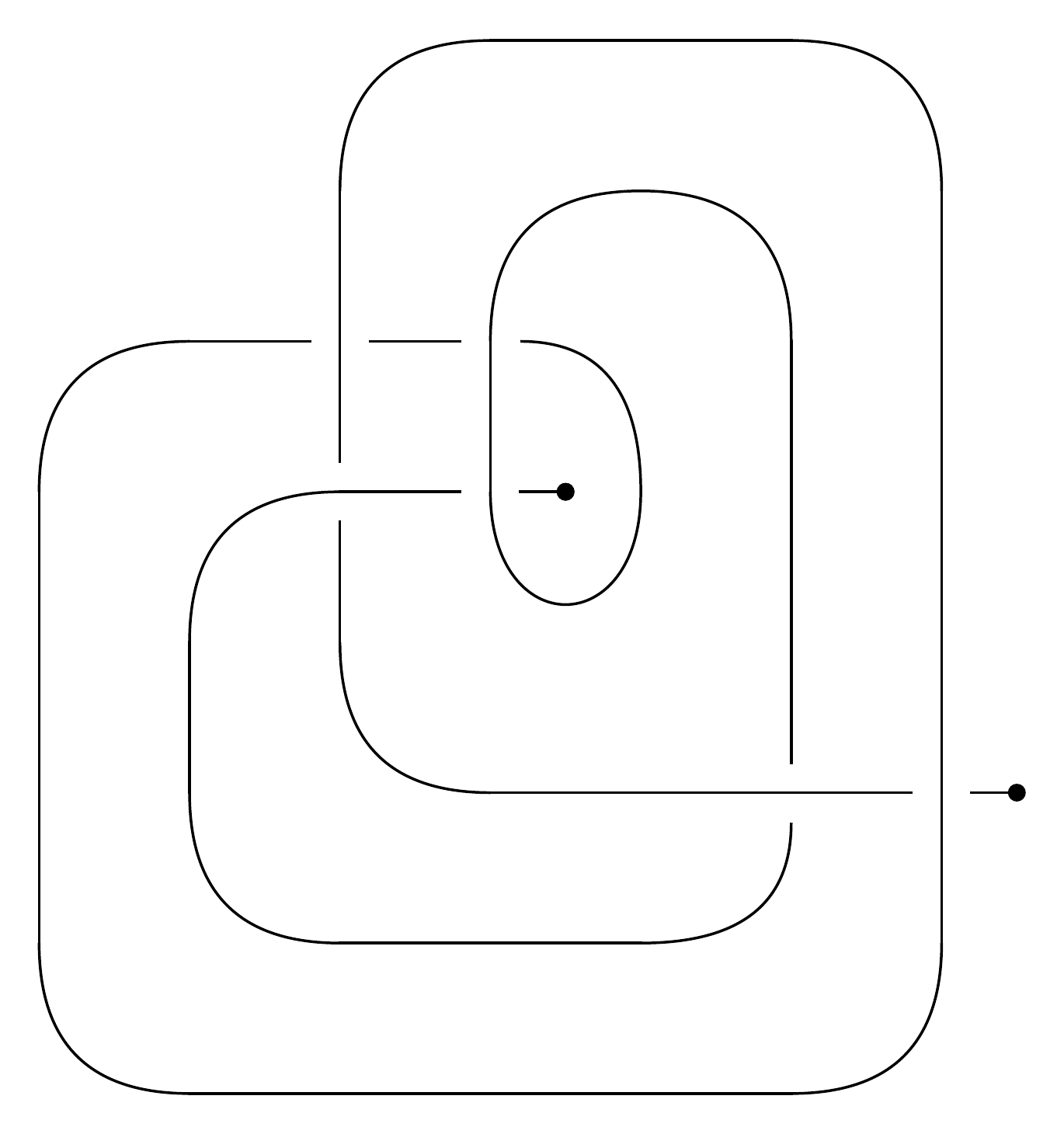}\\
\textcolor{black}{$6_{83}$}
\vspace{1cm}
\end{minipage}
\begin{minipage}[t]{.25\linewidth}
\centering
\includegraphics[width=0.9\textwidth,height=3.5cm,keepaspectratio]{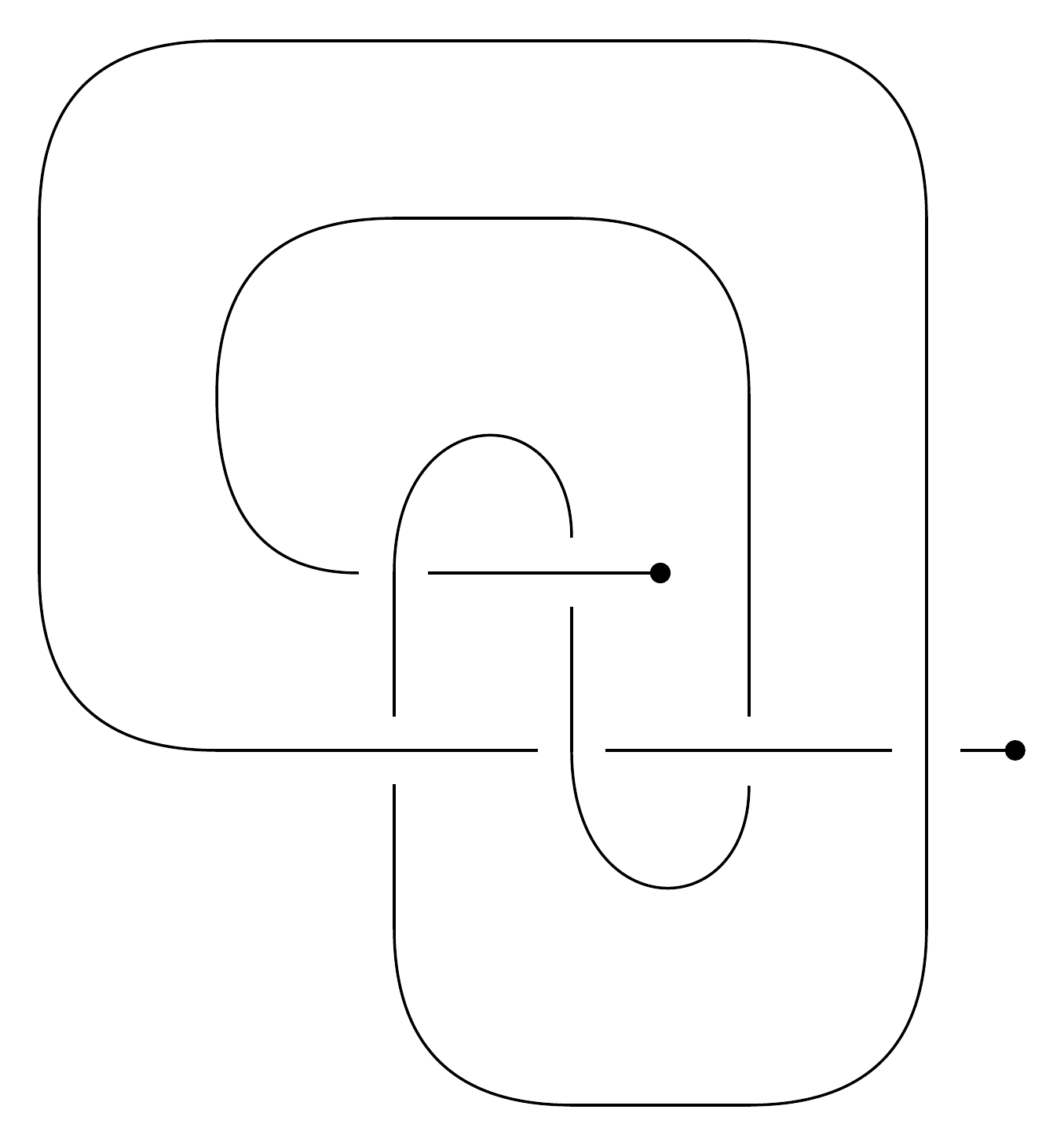}\\
\textcolor{black}{$6_{84}$}
\vspace{1cm}
\end{minipage}
\begin{minipage}[t]{.25\linewidth}
\centering
\includegraphics[width=0.9\textwidth,height=3.5cm,keepaspectratio]{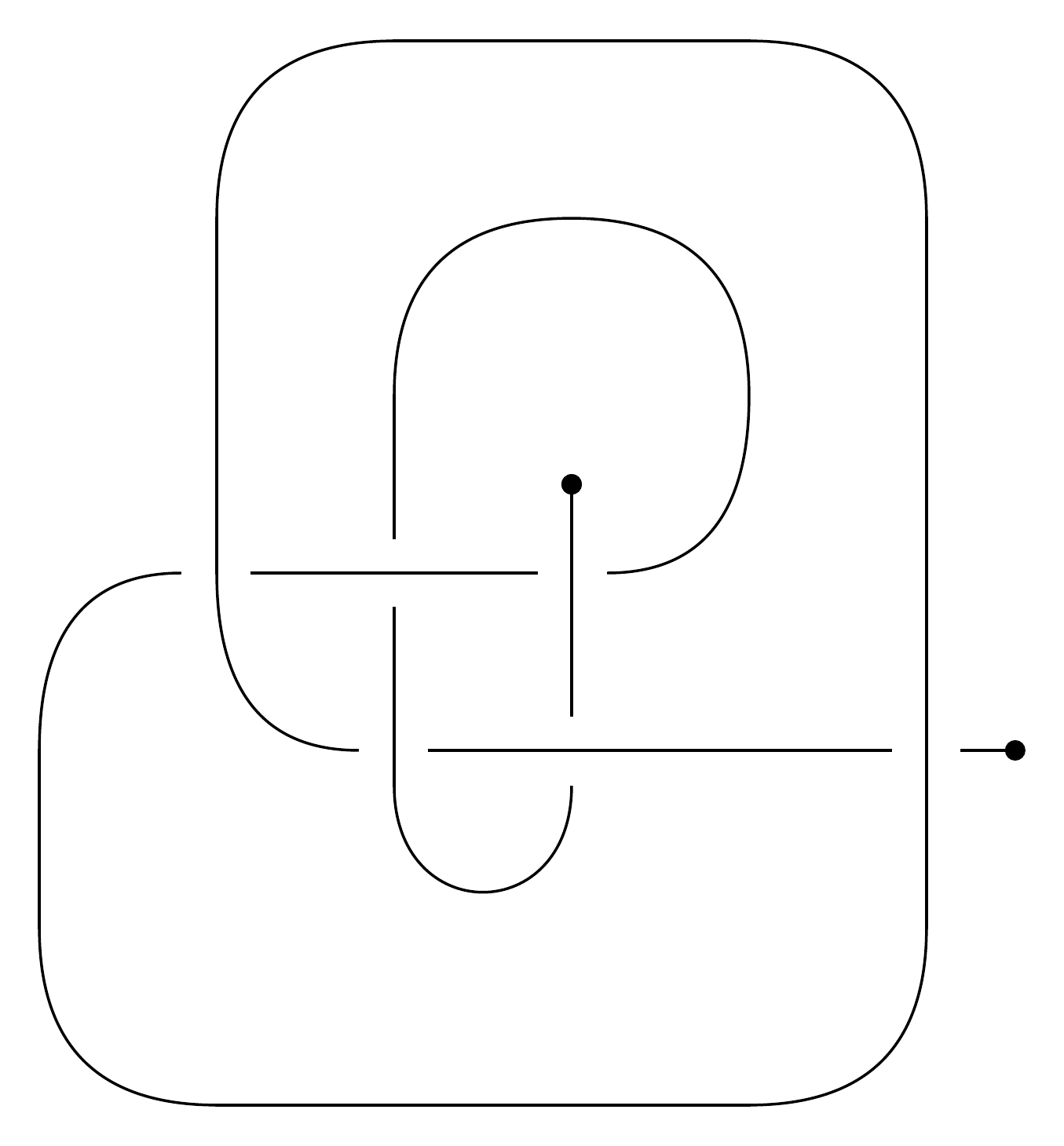}\\
\textcolor{black}{$6_{85}$}
\vspace{1cm}
\end{minipage}
\begin{minipage}[t]{.25\linewidth}
\centering
\includegraphics[width=0.9\textwidth,height=3.5cm,keepaspectratio]{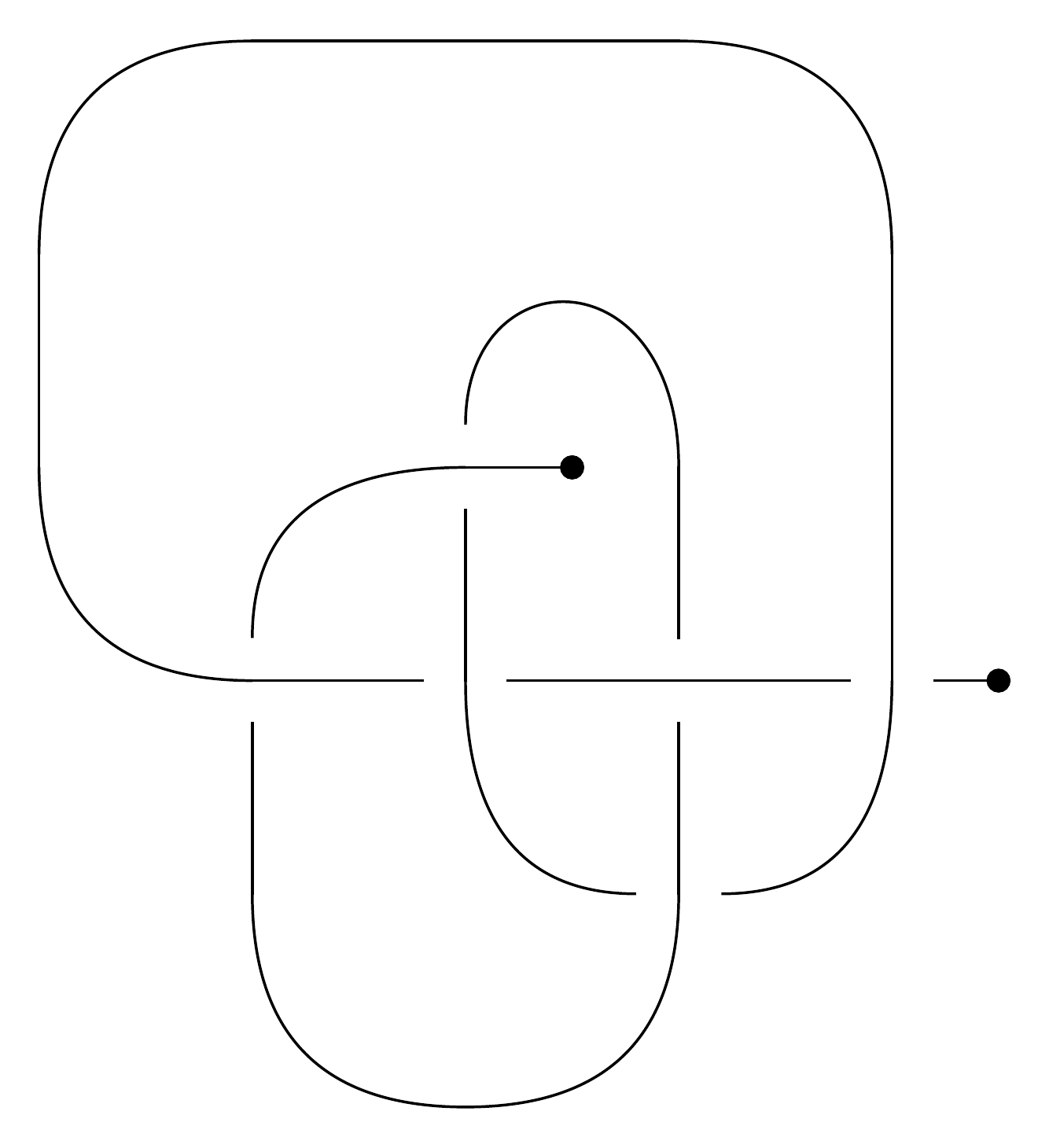}\\
\textcolor{black}{$6_{86}$}
\vspace{1cm}
\end{minipage}
\begin{minipage}[t]{.25\linewidth}
\centering
\includegraphics[width=0.9\textwidth,height=3.5cm,keepaspectratio]{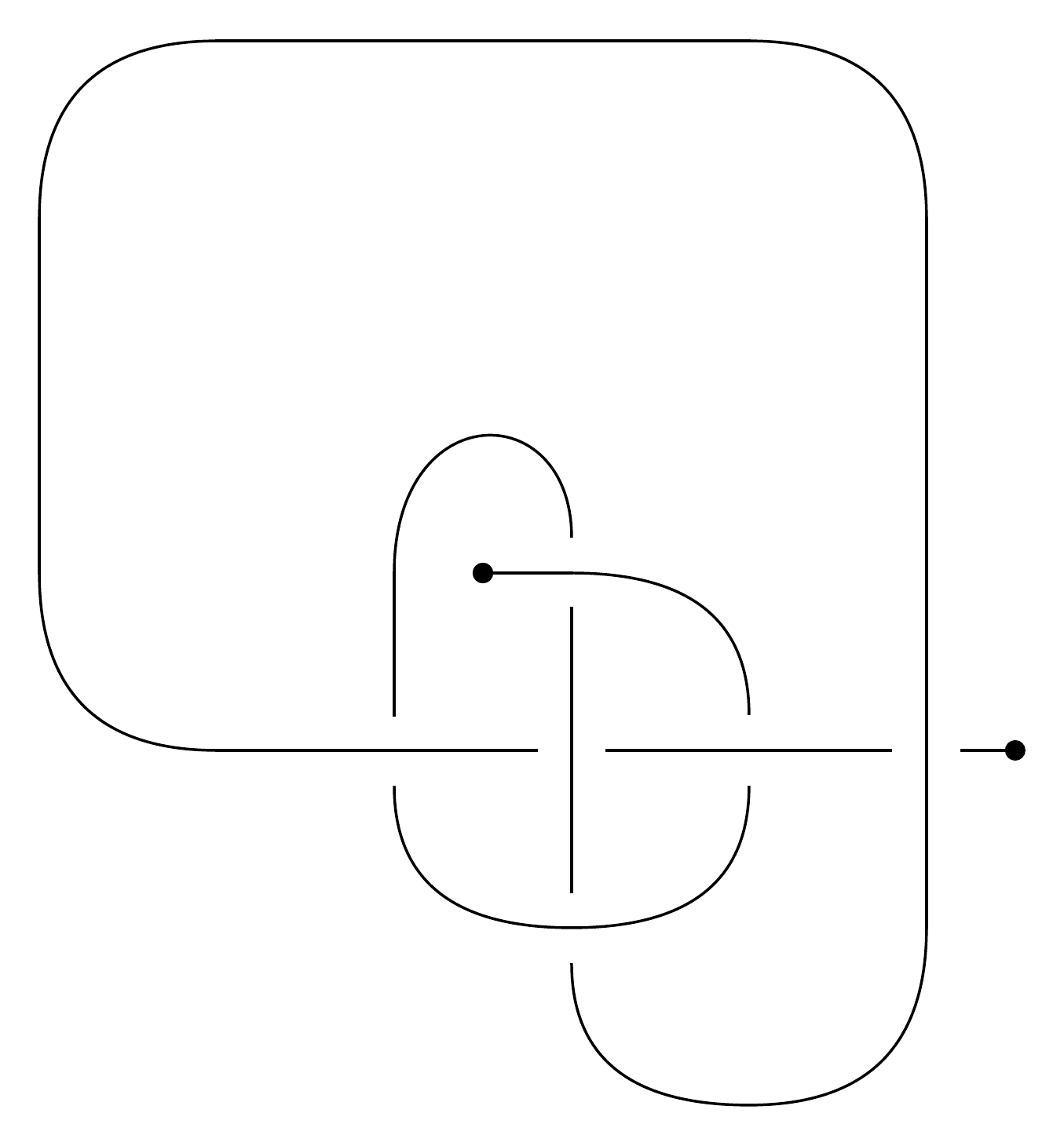}\\
\textcolor{black}{$6_{87}$}
\vspace{1cm}
\end{minipage}
\begin{minipage}[t]{.25\linewidth}
\centering
\includegraphics[width=0.9\textwidth,height=3.5cm,keepaspectratio]{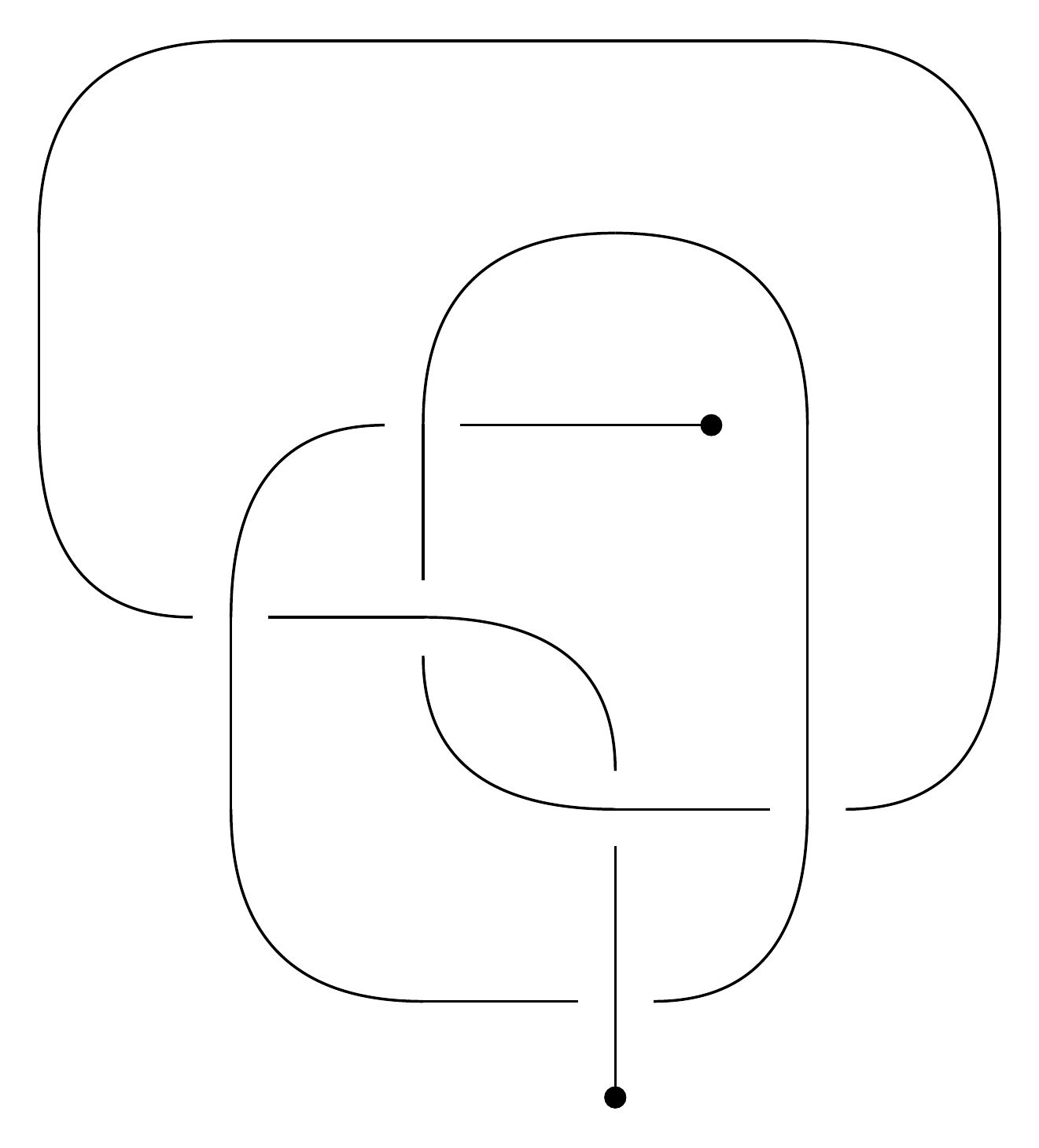}\\
\textcolor{black}{$6_{88}$}
\vspace{1cm}
\end{minipage}
\begin{minipage}[t]{.25\linewidth}
\centering
\includegraphics[width=0.9\textwidth,height=3.5cm,keepaspectratio]{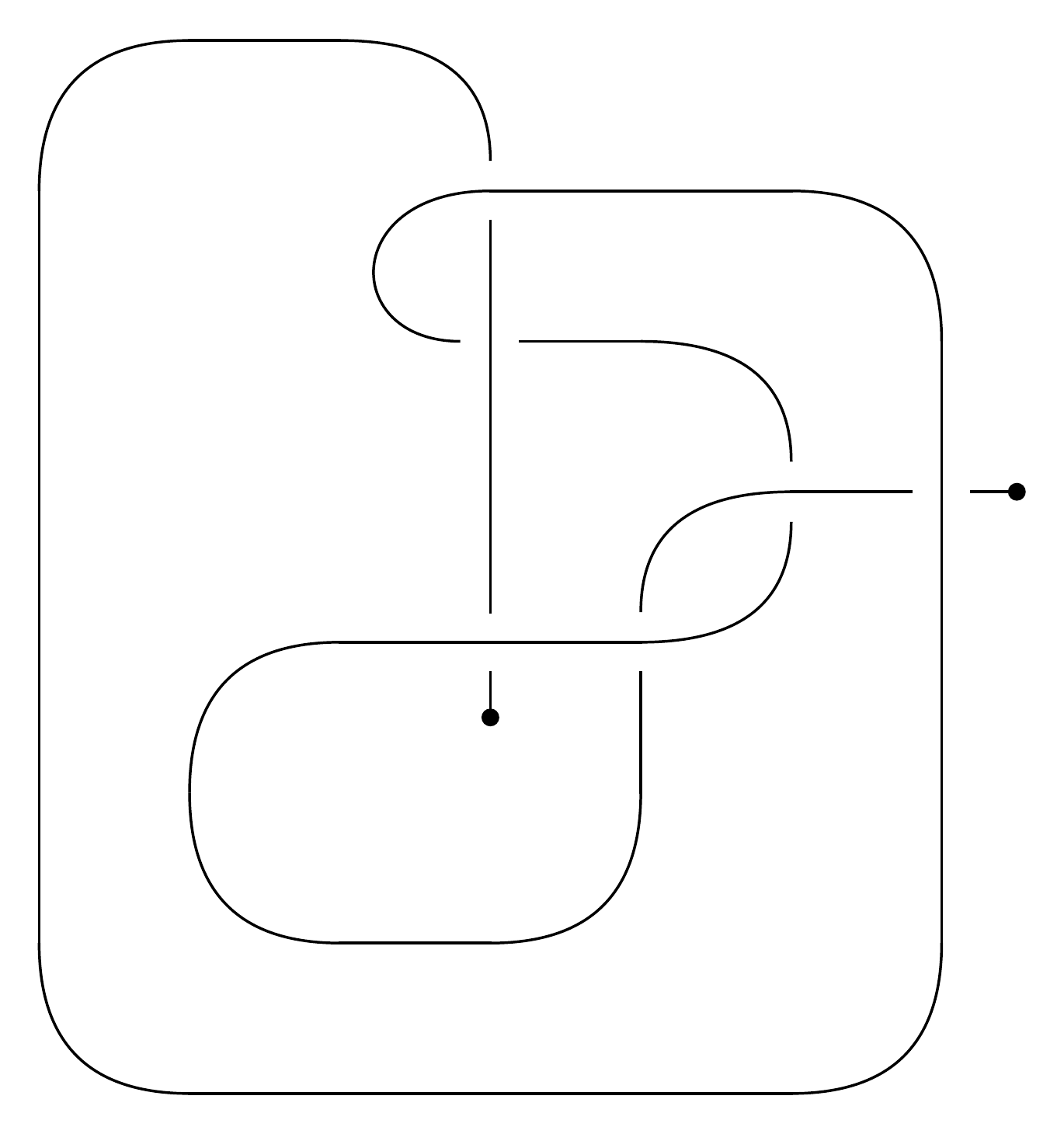}\\
\textcolor{black}{$6_{89}$}
\vspace{1cm}
\end{minipage}
\begin{minipage}[t]{.25\linewidth}
\centering
\includegraphics[width=0.9\textwidth,height=3.5cm,keepaspectratio]{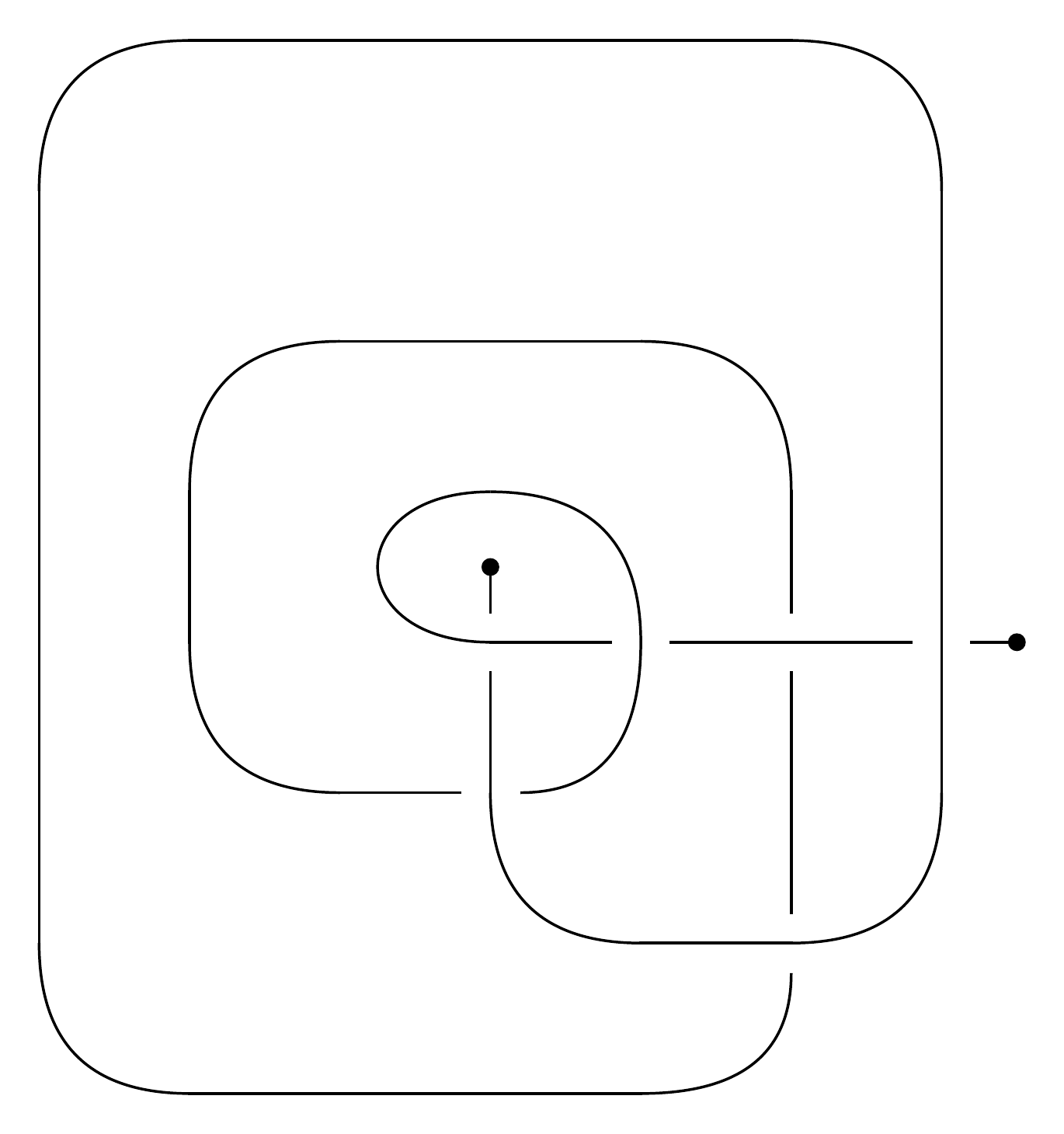}\\
\textcolor{black}{$6_{90}$}
\vspace{1cm}
\end{minipage}
\begin{minipage}[t]{.25\linewidth}
\centering
\includegraphics[width=0.9\textwidth,height=3.5cm,keepaspectratio]{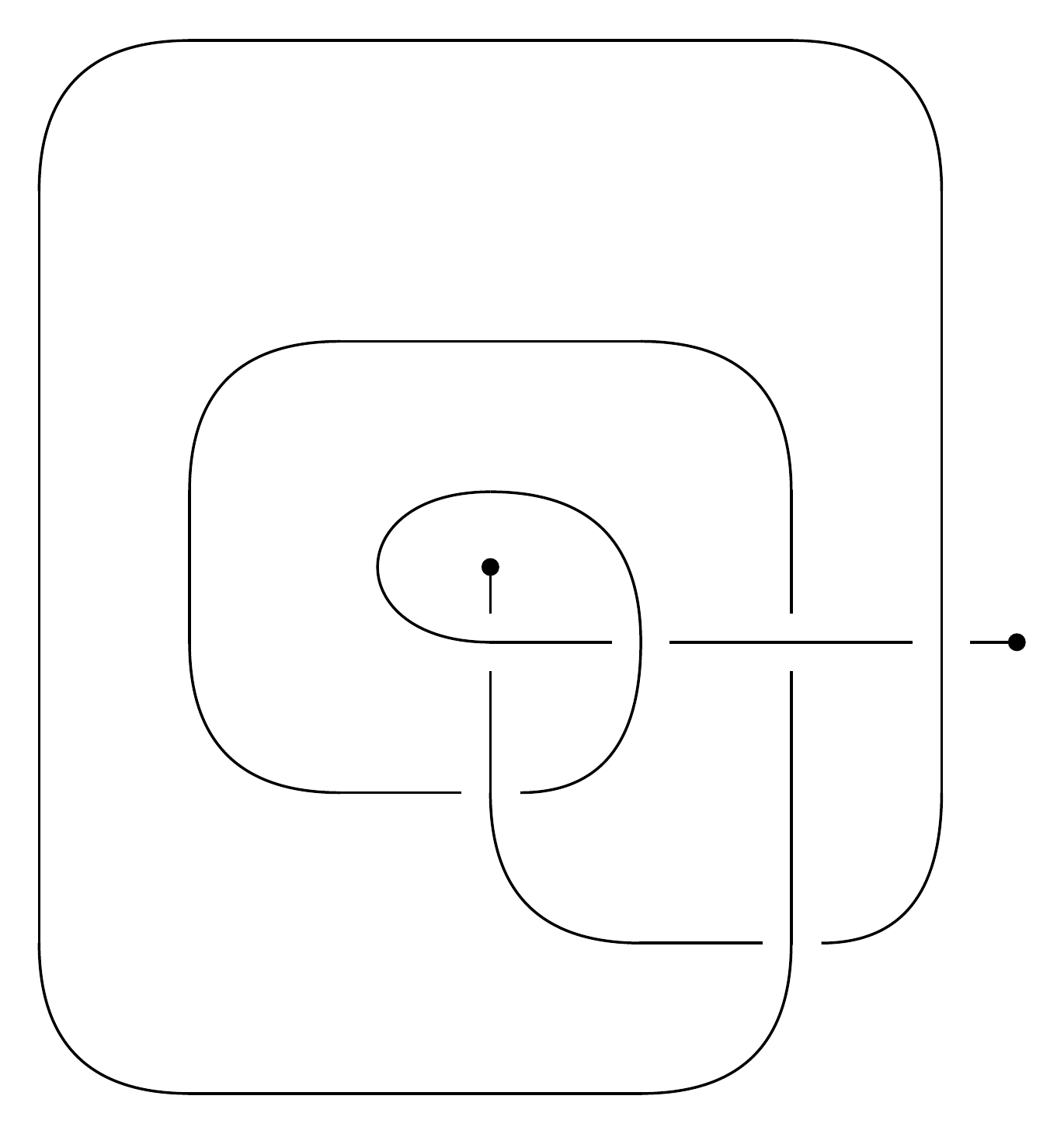}\\
\textcolor{black}{$6_{91}$}
\vspace{1cm}
\end{minipage}
\begin{minipage}[t]{.25\linewidth}
\centering
\includegraphics[width=0.9\textwidth,height=3.5cm,keepaspectratio]{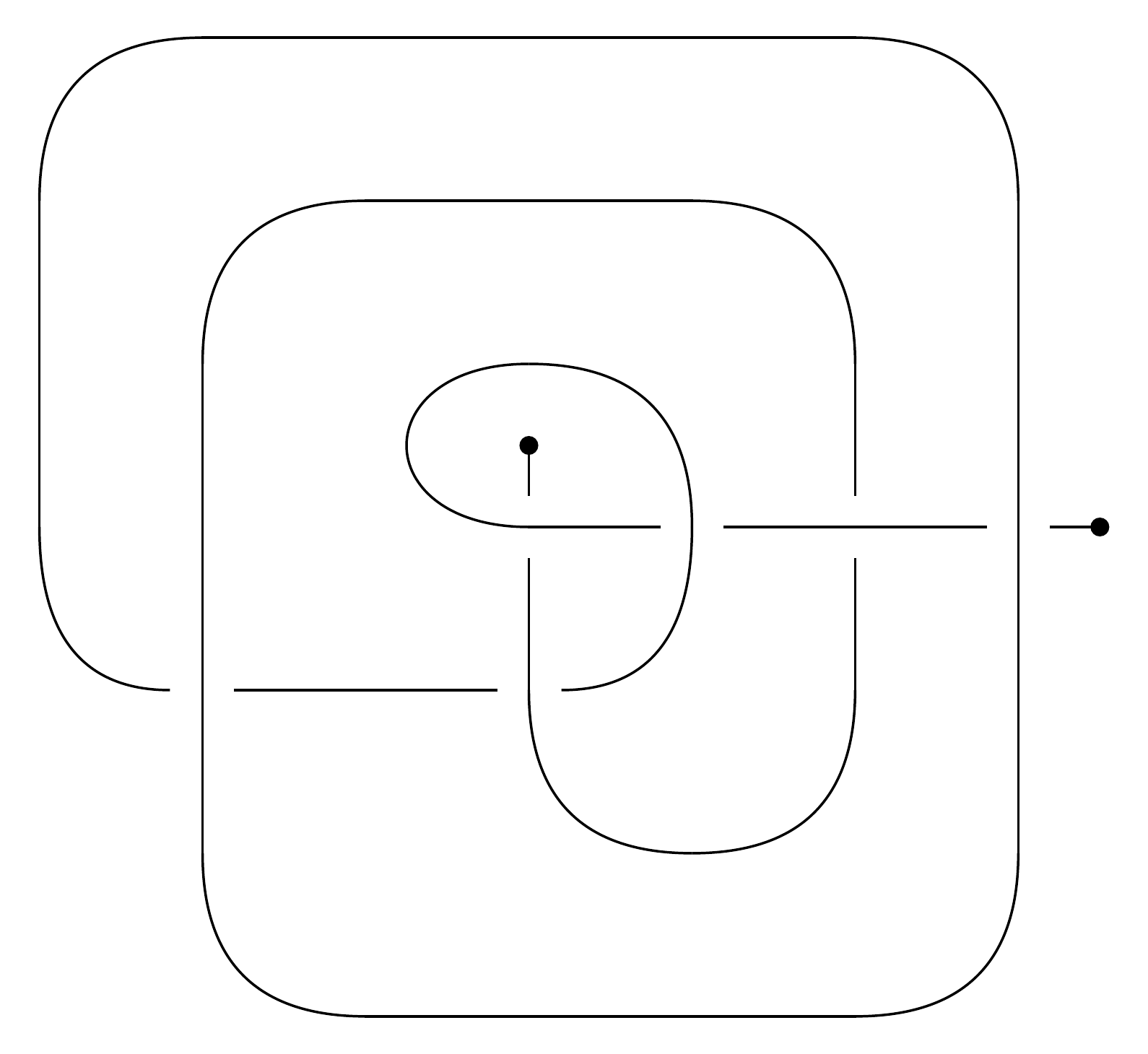}\\
\textcolor{black}{$6_{92}$}
\vspace{1cm}
\end{minipage}
\begin{minipage}[t]{.25\linewidth}
\centering
\includegraphics[width=0.9\textwidth,height=3.5cm,keepaspectratio]{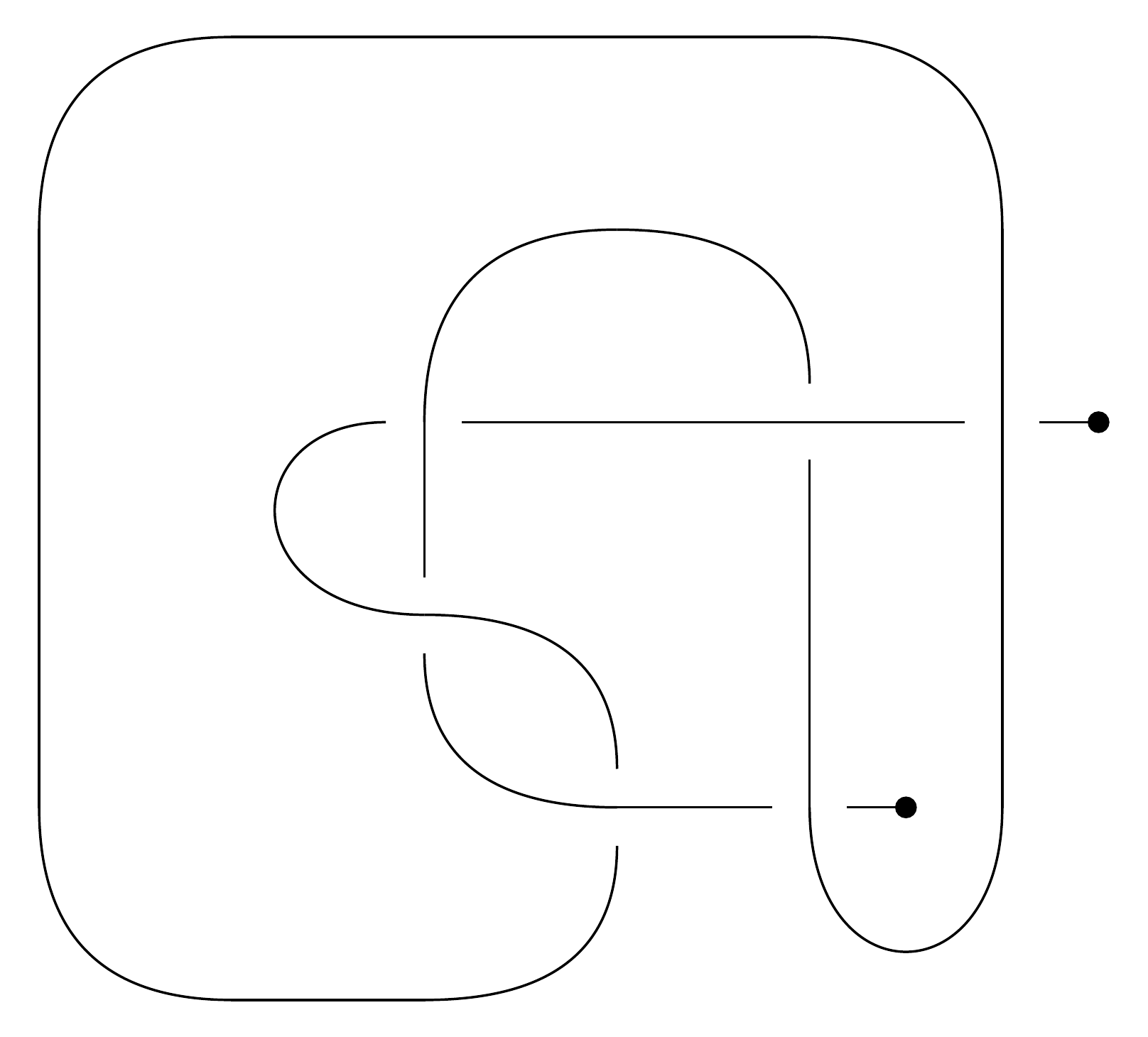}\\
\textcolor{black}{$6_{93}$}
\vspace{1cm}
\end{minipage}
\begin{minipage}[t]{.25\linewidth}
\centering
\includegraphics[width=0.9\textwidth,height=3.5cm,keepaspectratio]{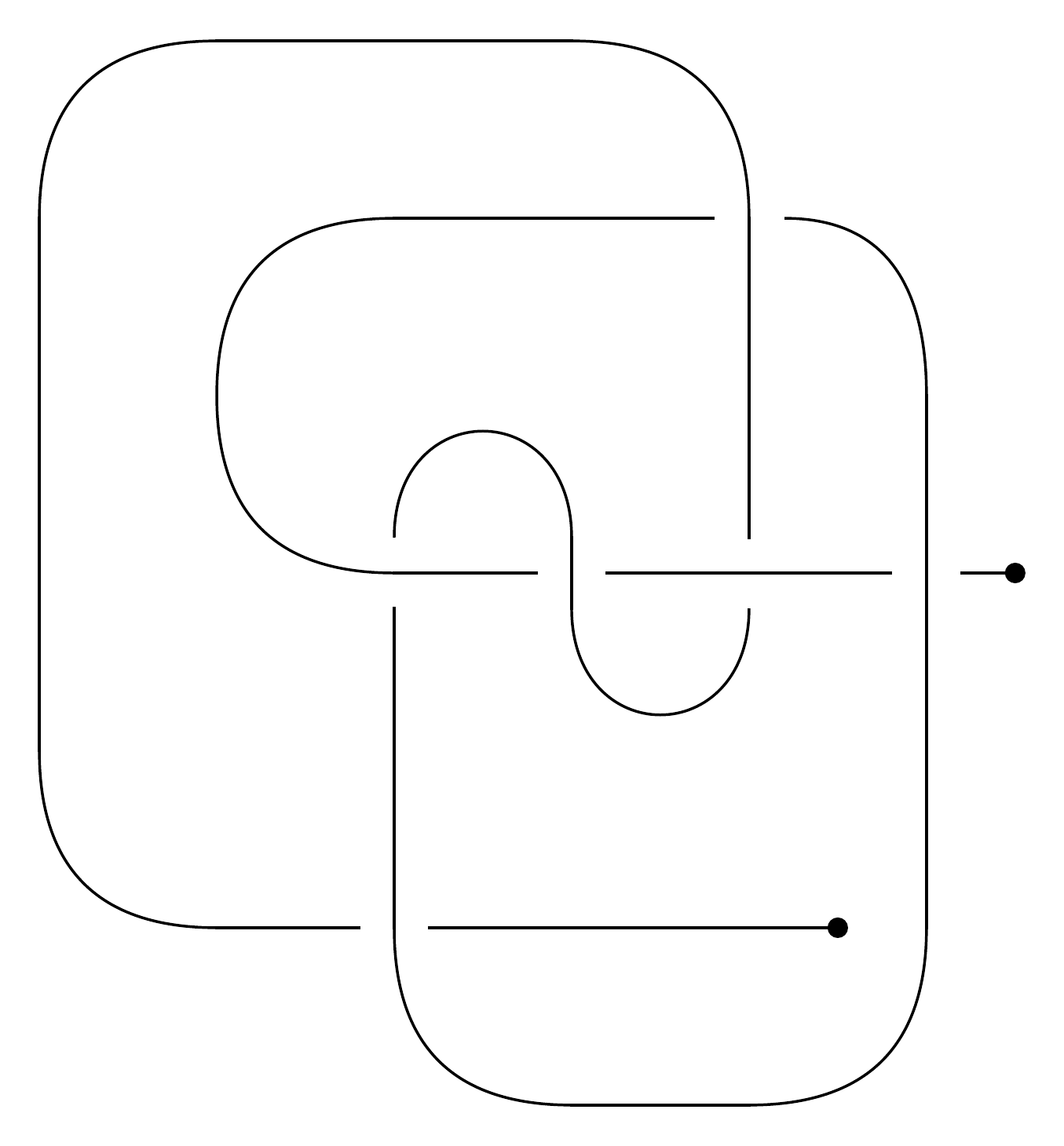}\\
\textcolor{black}{$6_{94}$}
\vspace{1cm}
\end{minipage}
\begin{minipage}[t]{.25\linewidth}
\centering
\includegraphics[width=0.9\textwidth,height=3.5cm,keepaspectratio]{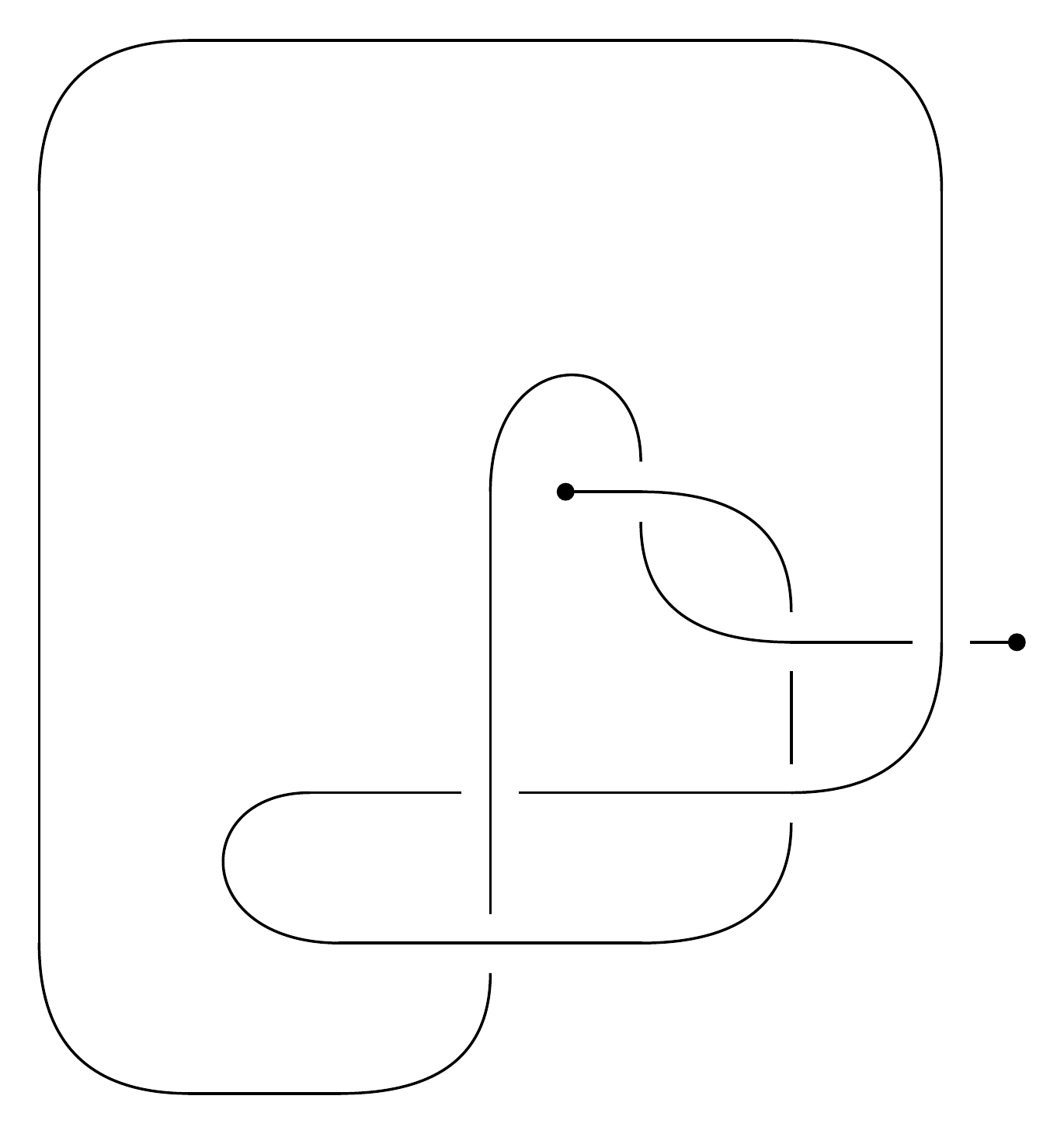}\\
\textcolor{black}{$6_{95}$}
\vspace{1cm}
\end{minipage}
\begin{minipage}[t]{.25\linewidth}
\centering
\includegraphics[width=0.9\textwidth,height=3.5cm,keepaspectratio]{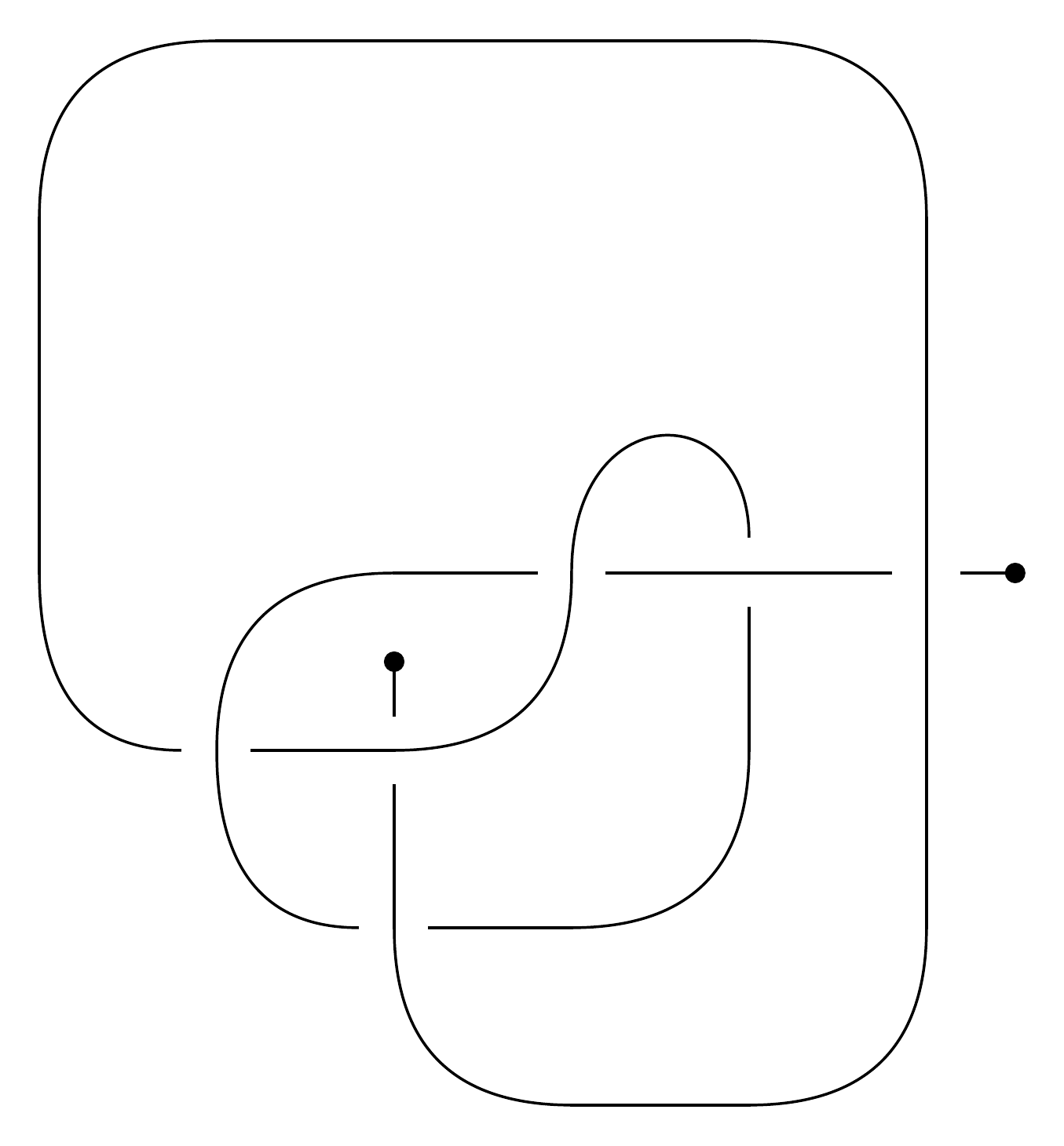}\\
\textcolor{black}{$6_{96}$}
\vspace{1cm}
\end{minipage}
\begin{minipage}[t]{.25\linewidth}
\centering
\includegraphics[width=0.9\textwidth,height=3.5cm,keepaspectratio]{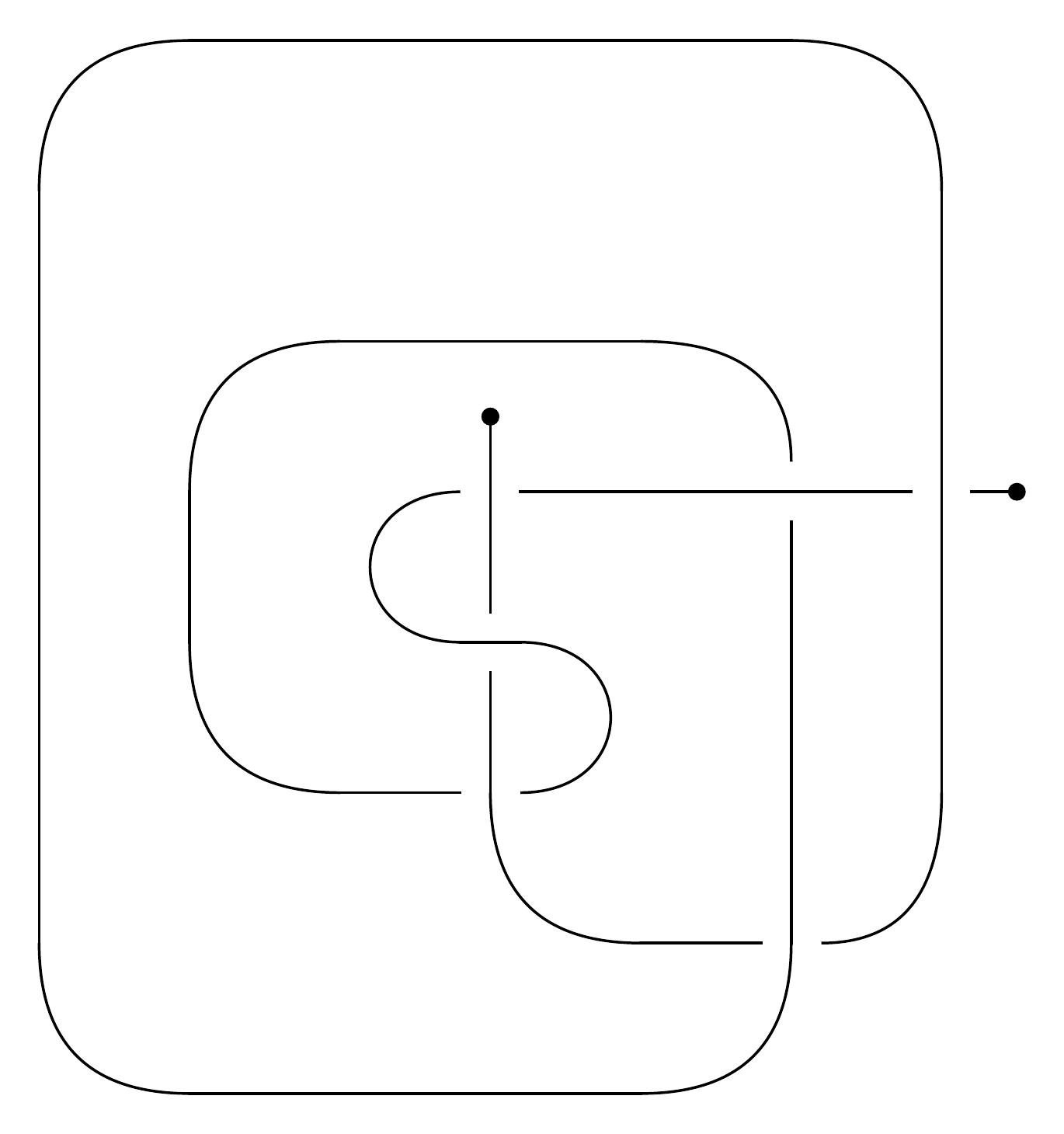}\\
\textcolor{black}{$6_{97}$}
\vspace{1cm}
\end{minipage}
\begin{minipage}[t]{.25\linewidth}
\centering
\includegraphics[width=0.9\textwidth,height=3.5cm,keepaspectratio]{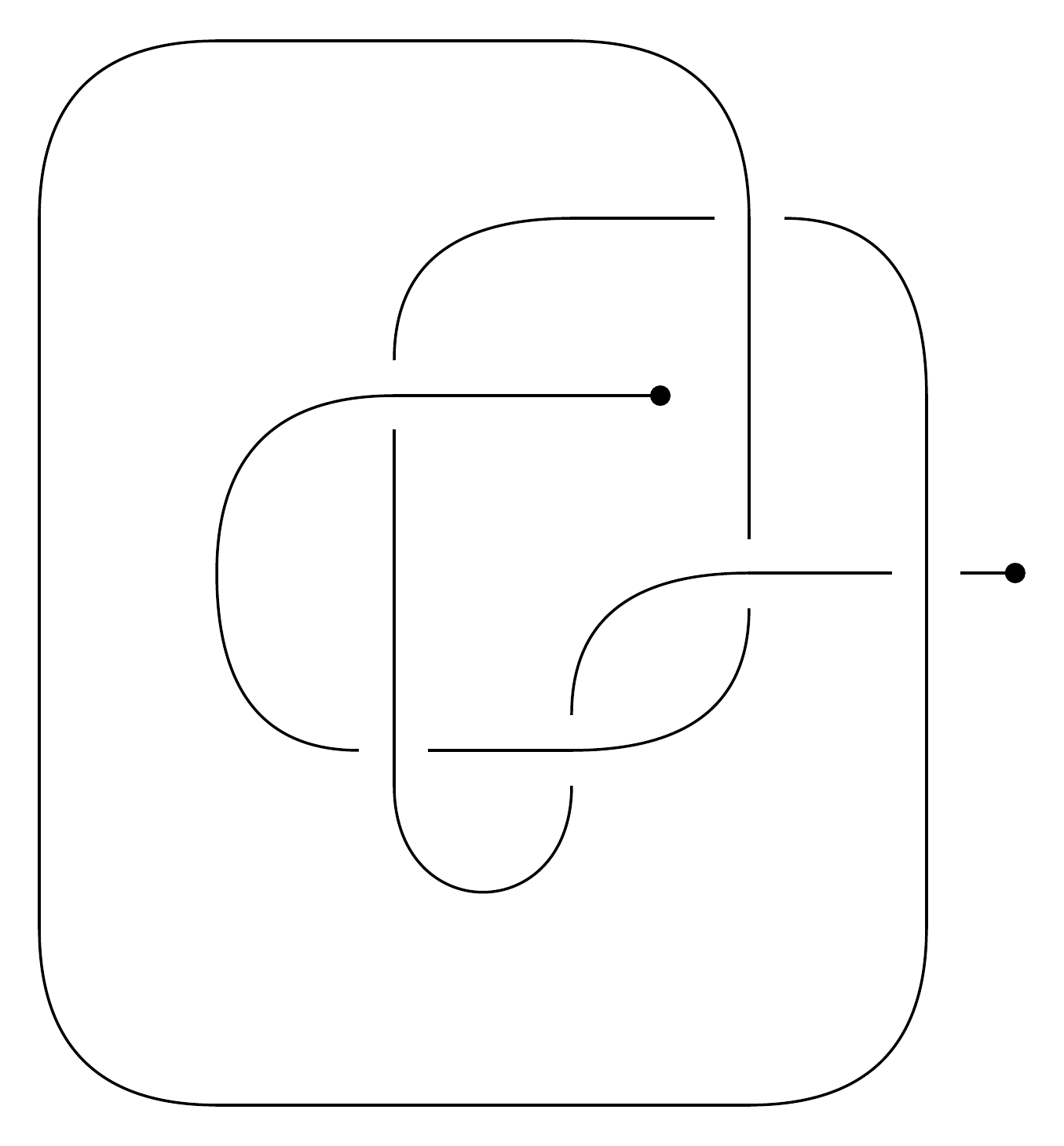}\\
\textcolor{black}{$6_{98}$}
\vspace{1cm}
\end{minipage}
\begin{minipage}[t]{.25\linewidth}
\centering
\includegraphics[width=0.9\textwidth,height=3.5cm,keepaspectratio]{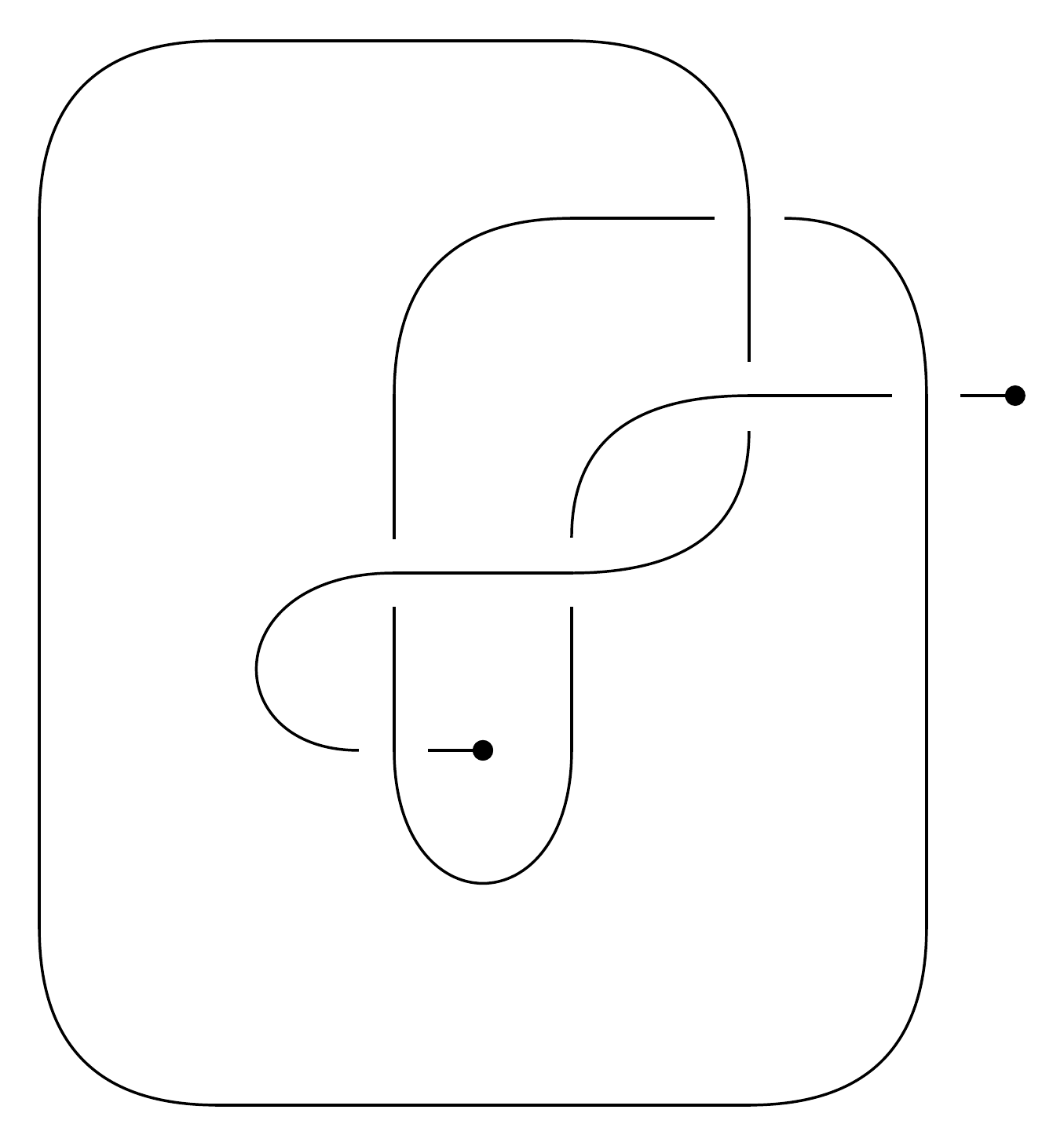}\\
\textcolor{black}{$6_{99}$}
\vspace{1cm}
\end{minipage}
\begin{minipage}[t]{.25\linewidth}
\centering
\includegraphics[width=0.9\textwidth,height=3.5cm,keepaspectratio]{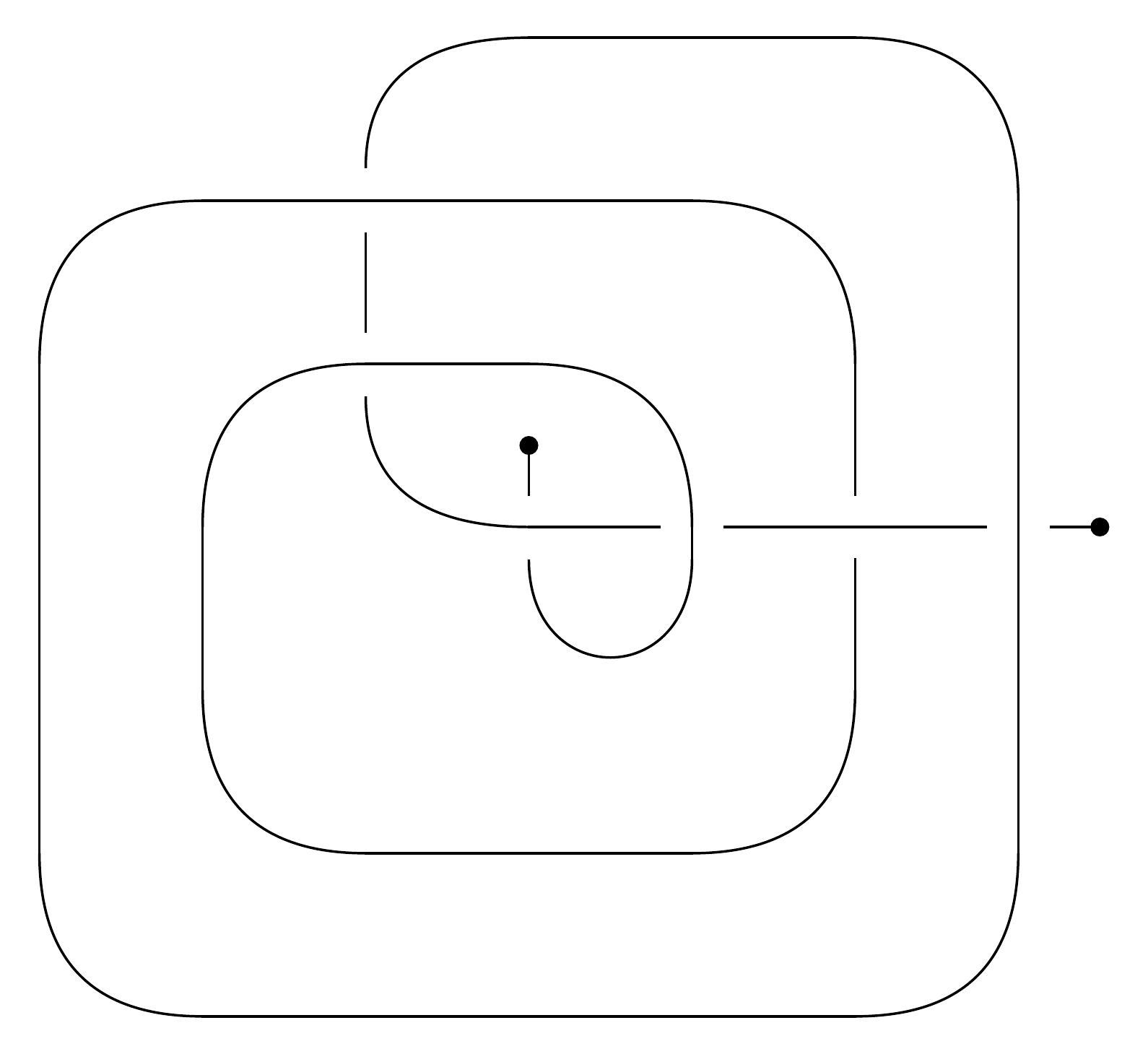}\\
\textcolor{black}{$6_{100}$}
\vspace{1cm}
\end{minipage}
\begin{minipage}[t]{.25\linewidth}
\centering
\includegraphics[width=0.9\textwidth,height=3.5cm,keepaspectratio]{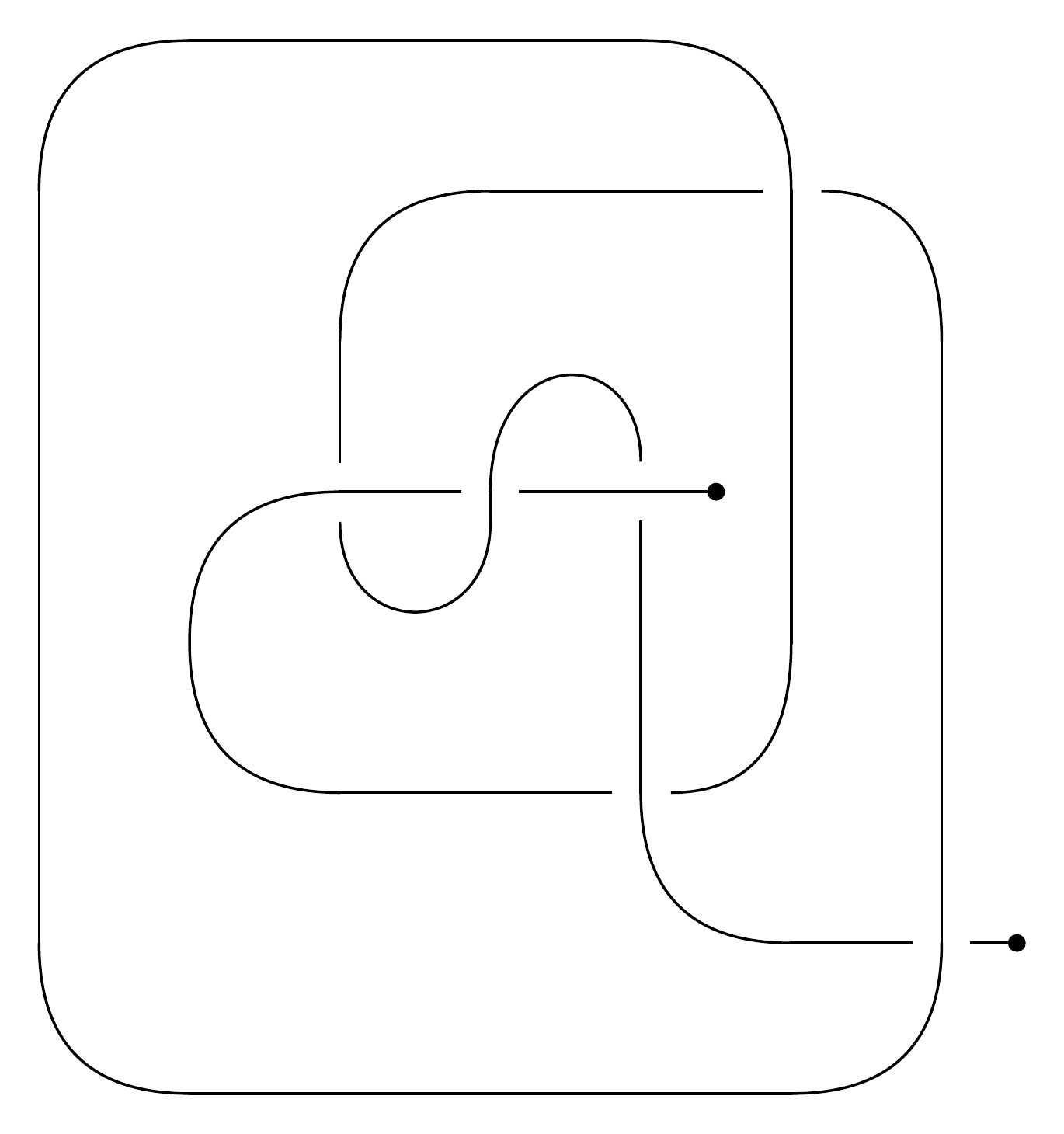}\\
\textcolor{black}{$6_{101}$}
\vspace{1cm}
\end{minipage}
\begin{minipage}[t]{.25\linewidth}
\centering
\includegraphics[width=0.9\textwidth,height=3.5cm,keepaspectratio]{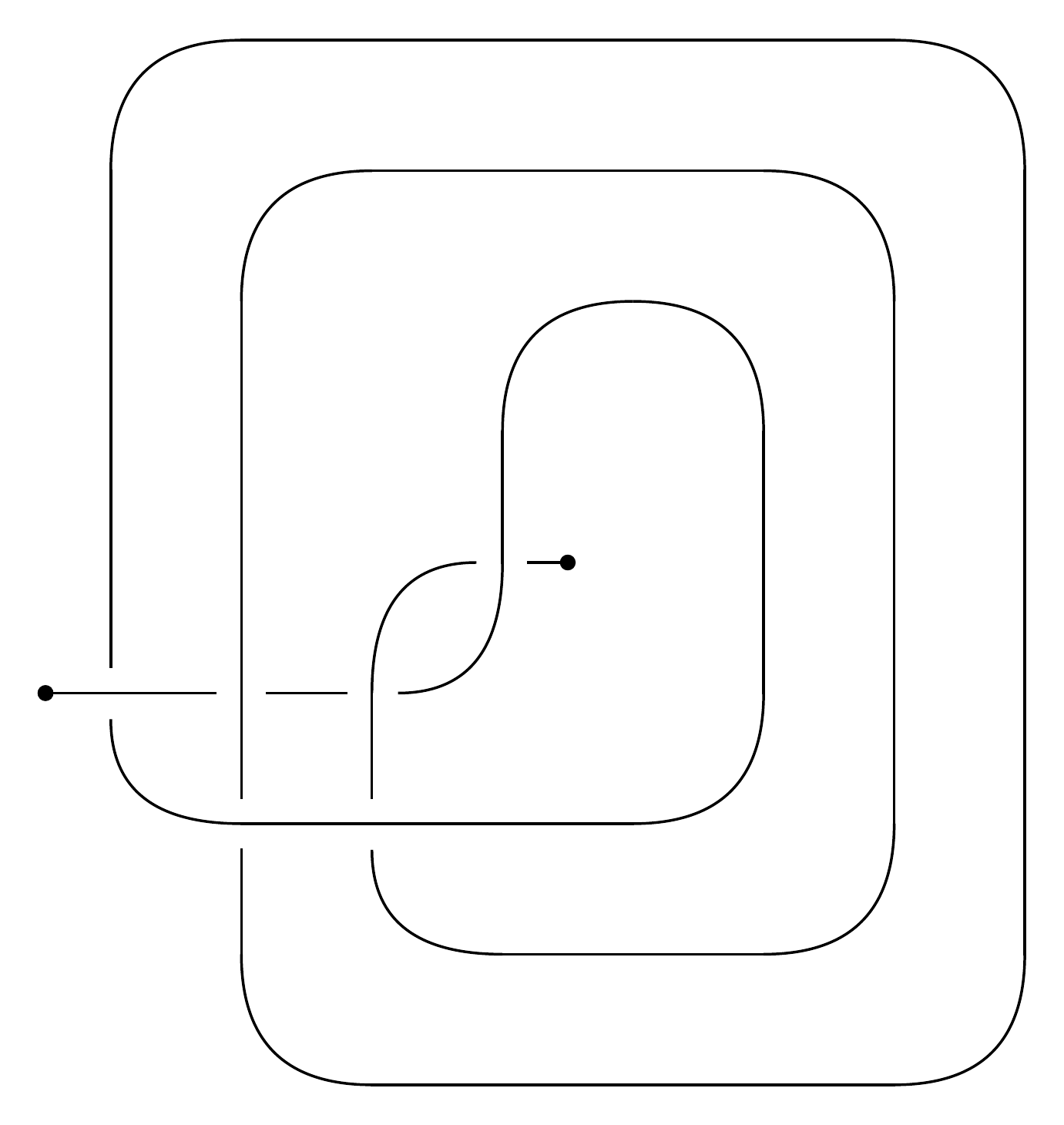}\\
\textcolor{black}{$6_{102}$}
\vspace{1cm}
\end{minipage}
\begin{minipage}[t]{.25\linewidth}
\centering
\includegraphics[width=0.9\textwidth,height=3.5cm,keepaspectratio]{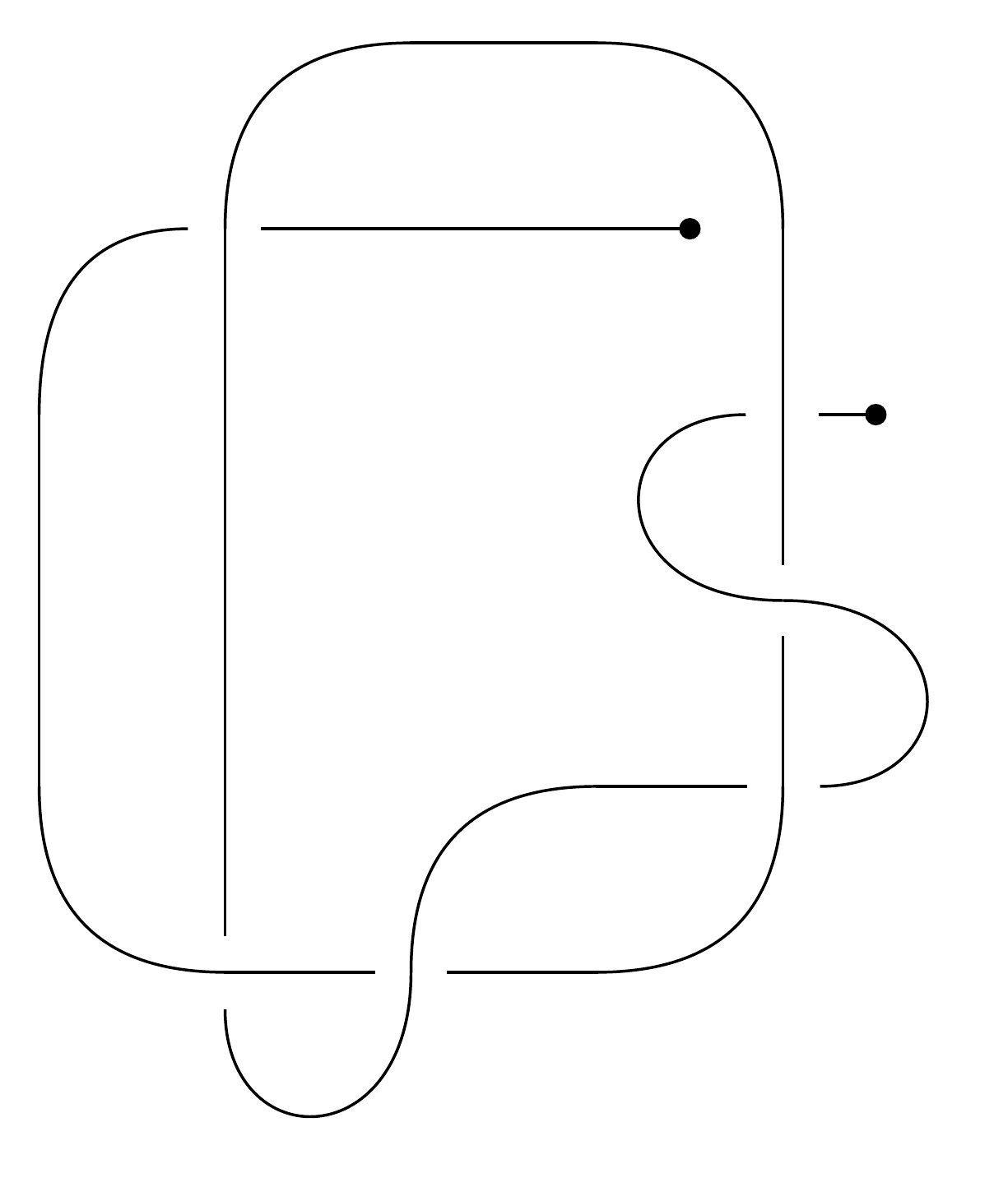}\\
\textcolor{black}{$6_{103}$}
\vspace{1cm}
\end{minipage}
\begin{minipage}[t]{.25\linewidth}
\centering
\includegraphics[width=0.9\textwidth,height=3.5cm,keepaspectratio]{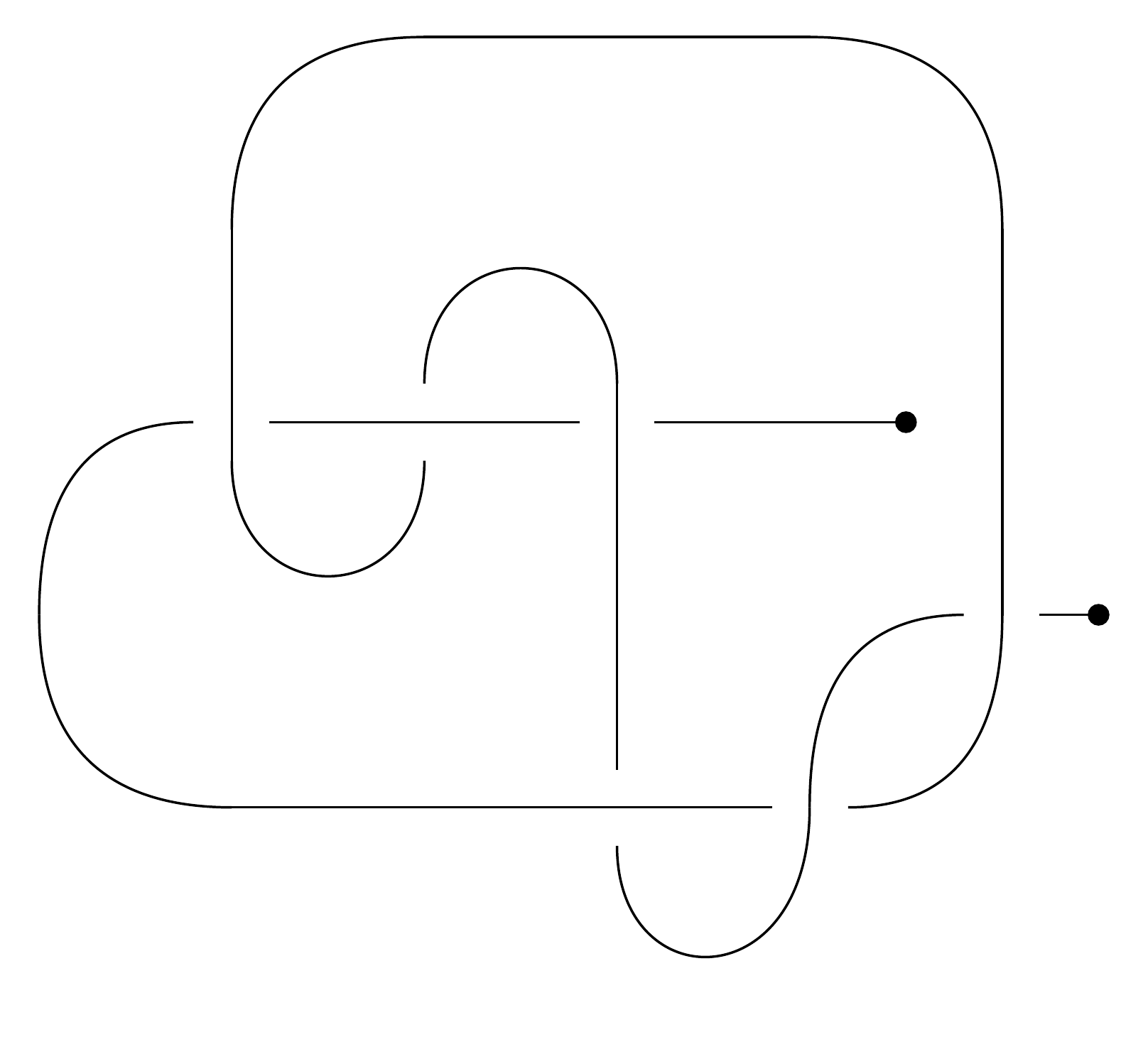}\\
\textcolor{black}{$6_{104}$}
\vspace{1cm}
\end{minipage}
\begin{minipage}[t]{.25\linewidth}
\centering
\includegraphics[width=0.9\textwidth,height=3.5cm,keepaspectratio]{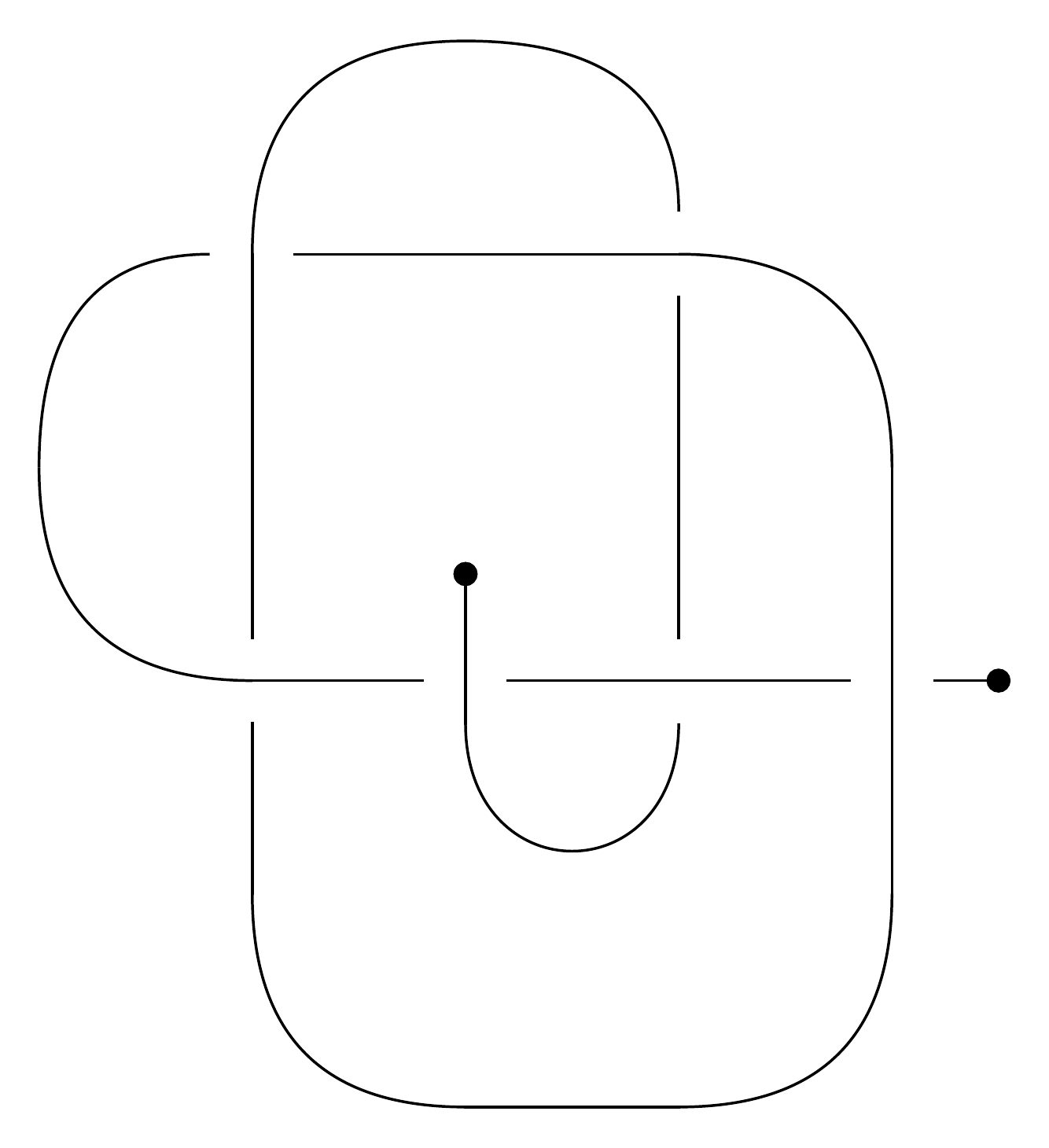}\\
\textcolor{black}{$6_{105}$}
\vspace{1cm}
\end{minipage}
\begin{minipage}[t]{.25\linewidth}
\centering
\includegraphics[width=0.9\textwidth,height=3.5cm,keepaspectratio]{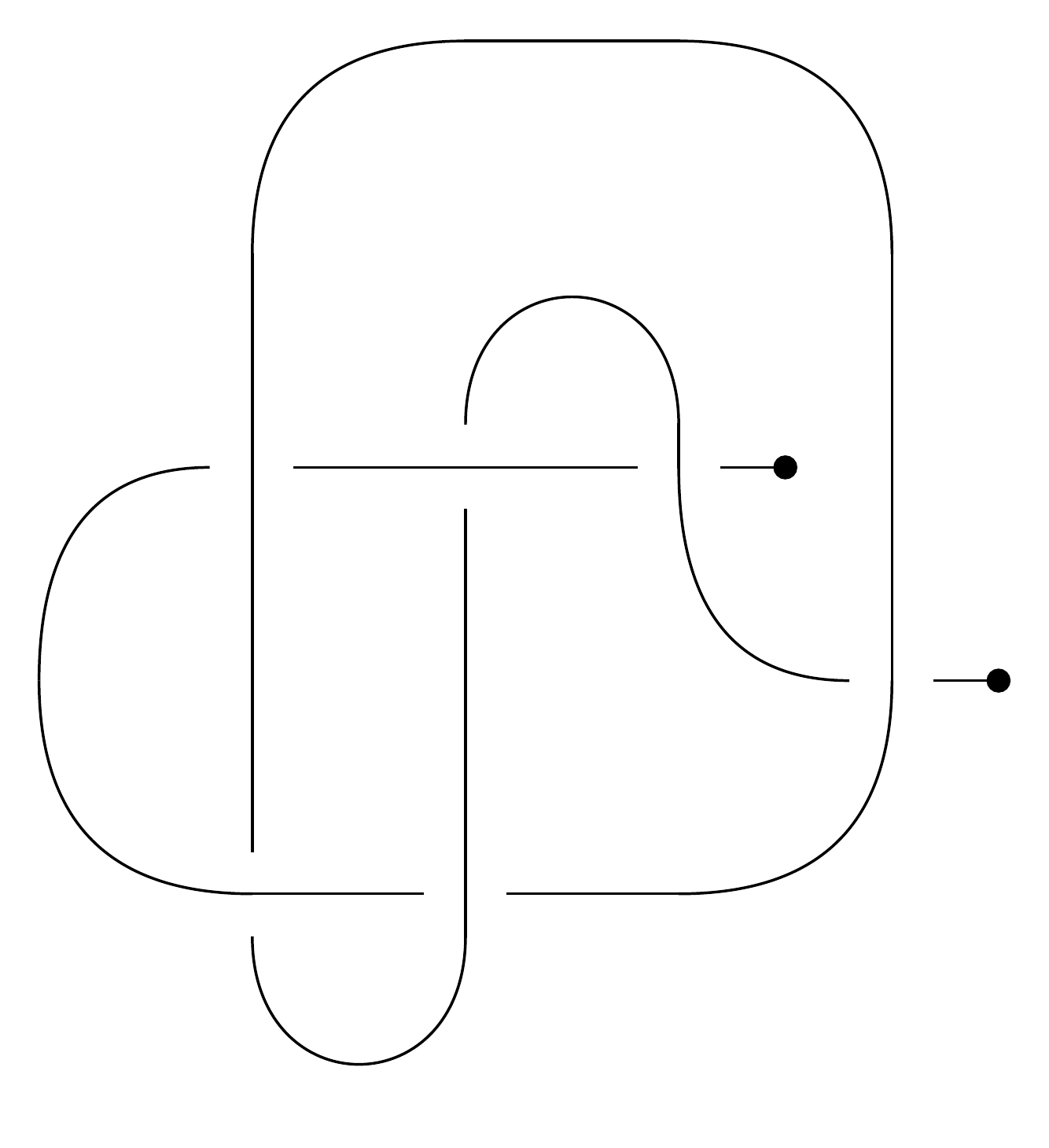}\\
\textcolor{black}{$6_{106}$}
\vspace{1cm}
\end{minipage}
\begin{minipage}[t]{.25\linewidth}
\centering
\includegraphics[width=0.9\textwidth,height=3.5cm,keepaspectratio]{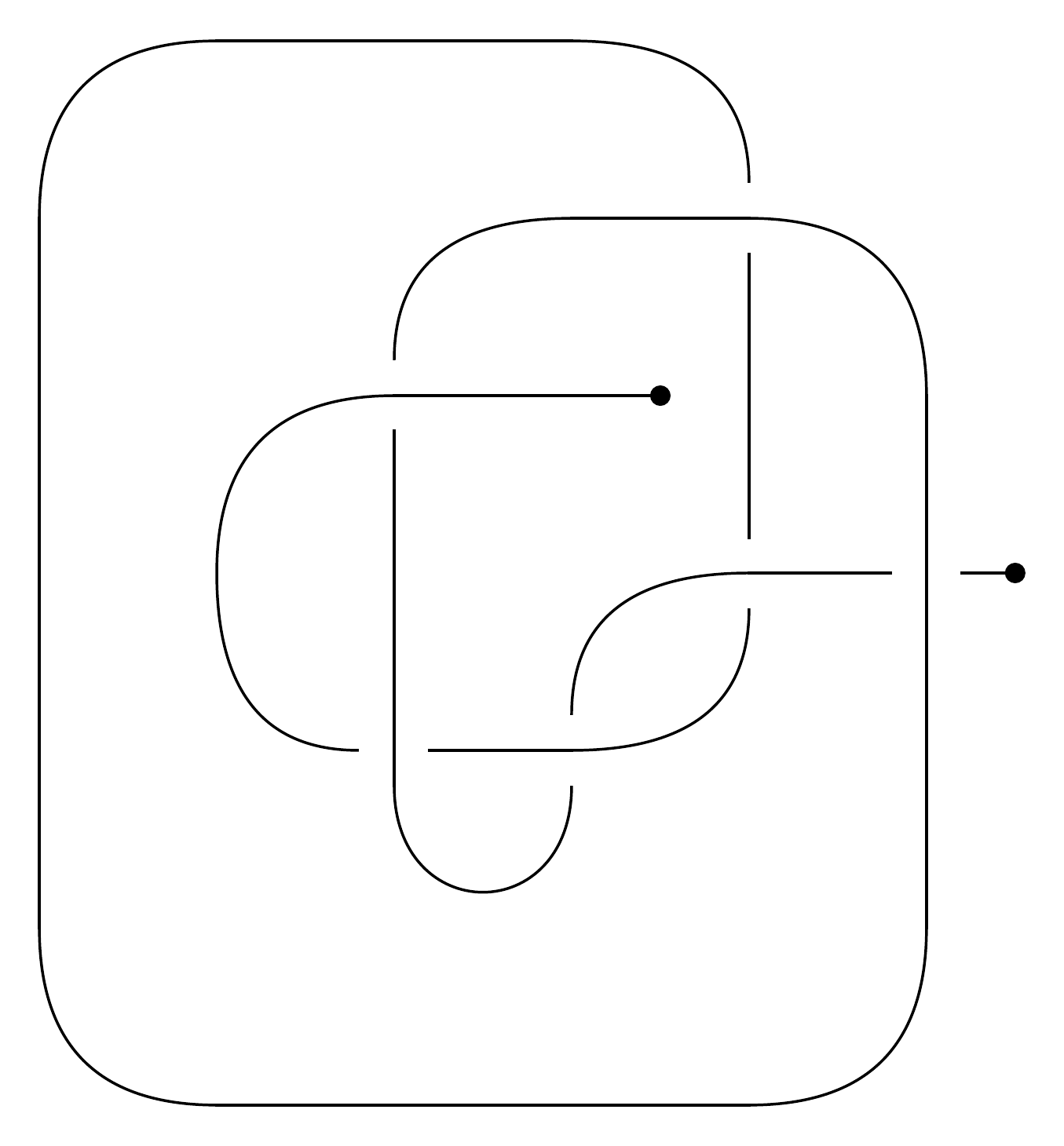}\\
\textcolor{black}{$6_{107}$}
\vspace{1cm}
\end{minipage}
\begin{minipage}[t]{.25\linewidth}
\centering
\includegraphics[width=0.9\textwidth,height=3.5cm,keepaspectratio]{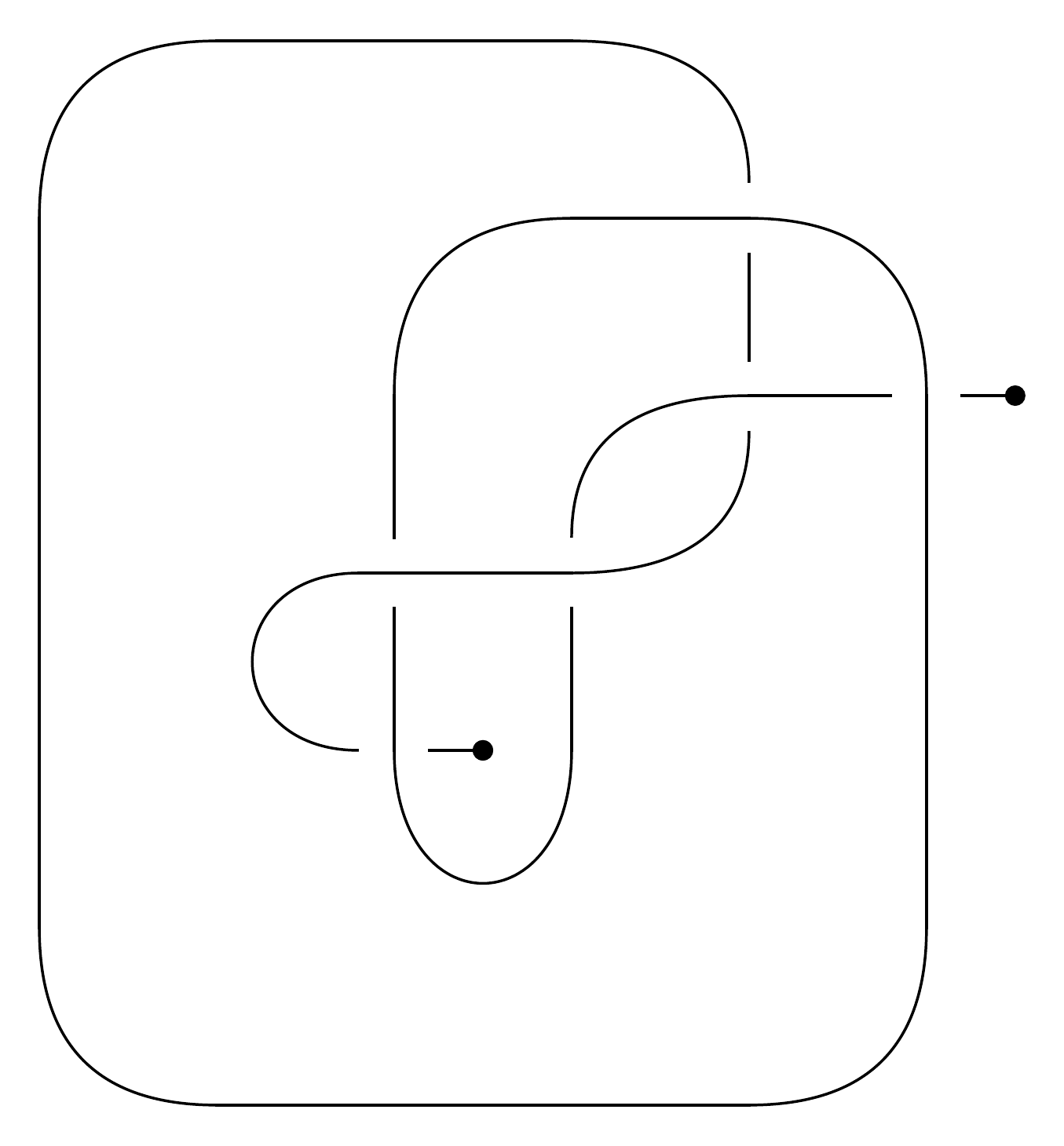}\\
\textcolor{black}{$6_{108}$}
\vspace{1cm}
\end{minipage}
\begin{minipage}[t]{.25\linewidth}
\centering
\includegraphics[width=0.9\textwidth,height=3.5cm,keepaspectratio]{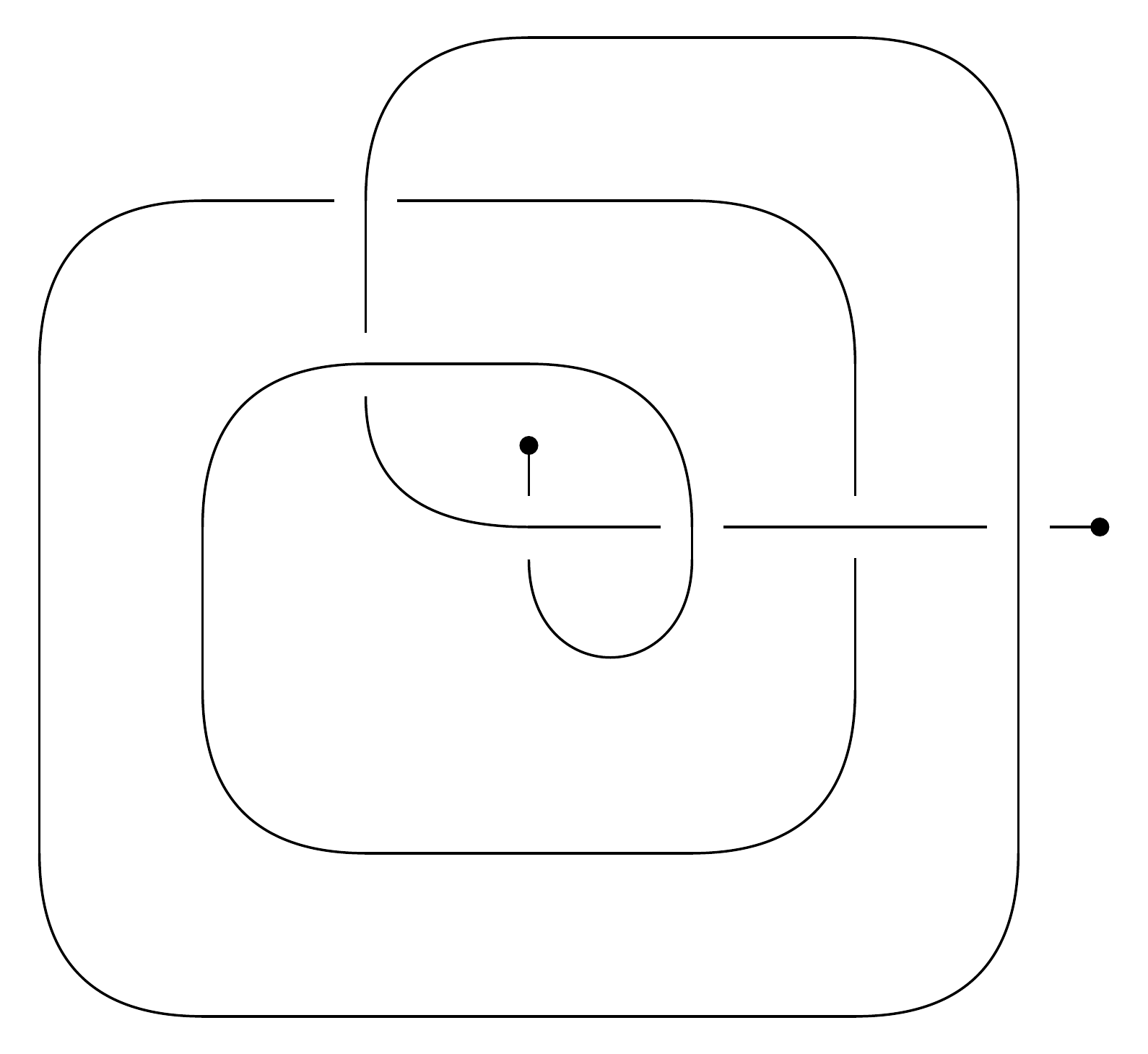}\\
\textcolor{black}{$6_{109}$}
\vspace{1cm}
\end{minipage}
\begin{minipage}[t]{.25\linewidth}
\centering
\includegraphics[width=0.9\textwidth,height=3.5cm,keepaspectratio]{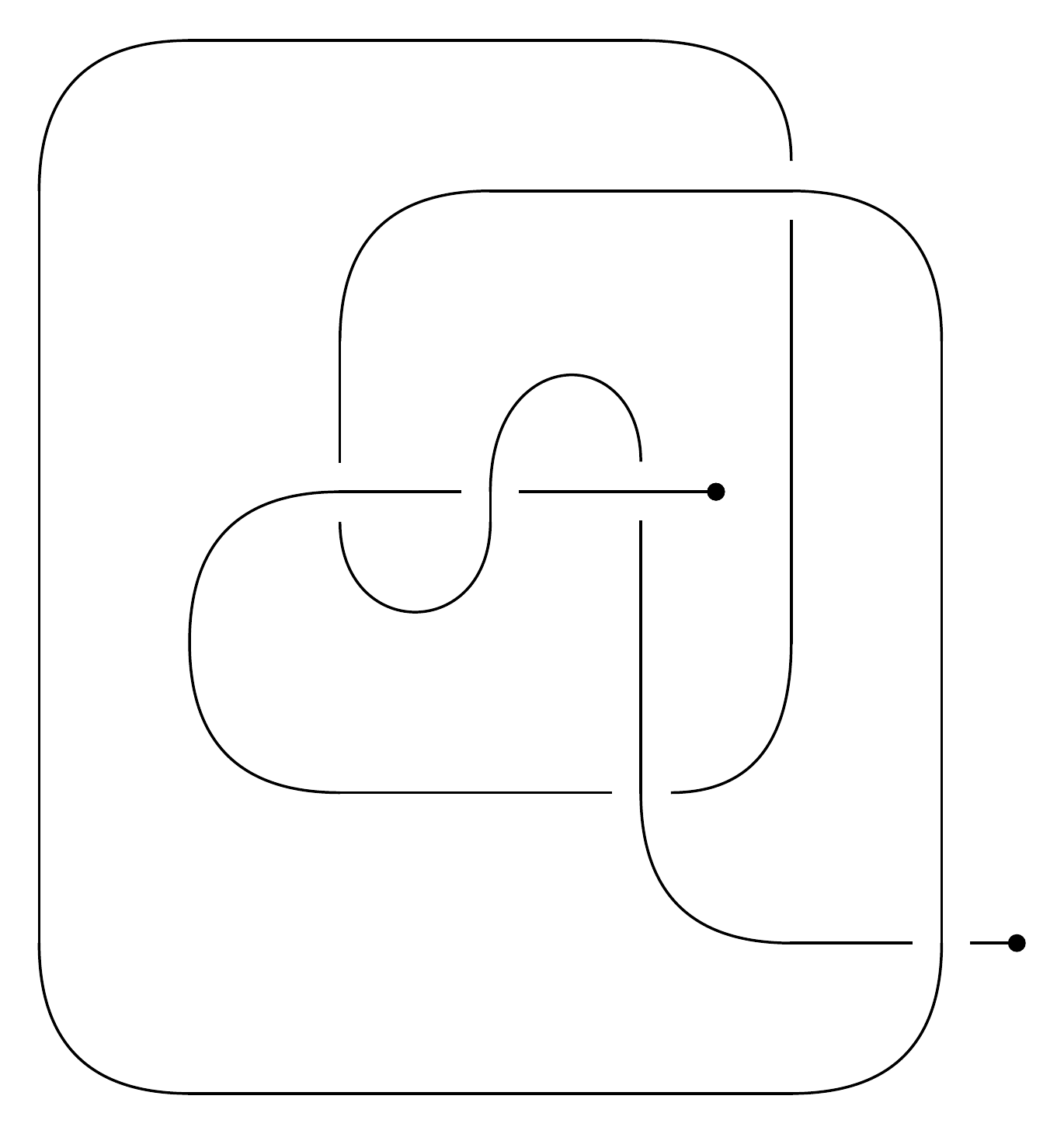}\\
\textcolor{black}{$6_{110}$}
\vspace{1cm}
\end{minipage}
\begin{minipage}[t]{.25\linewidth}
\centering
\includegraphics[width=0.9\textwidth,height=3.5cm,keepaspectratio]{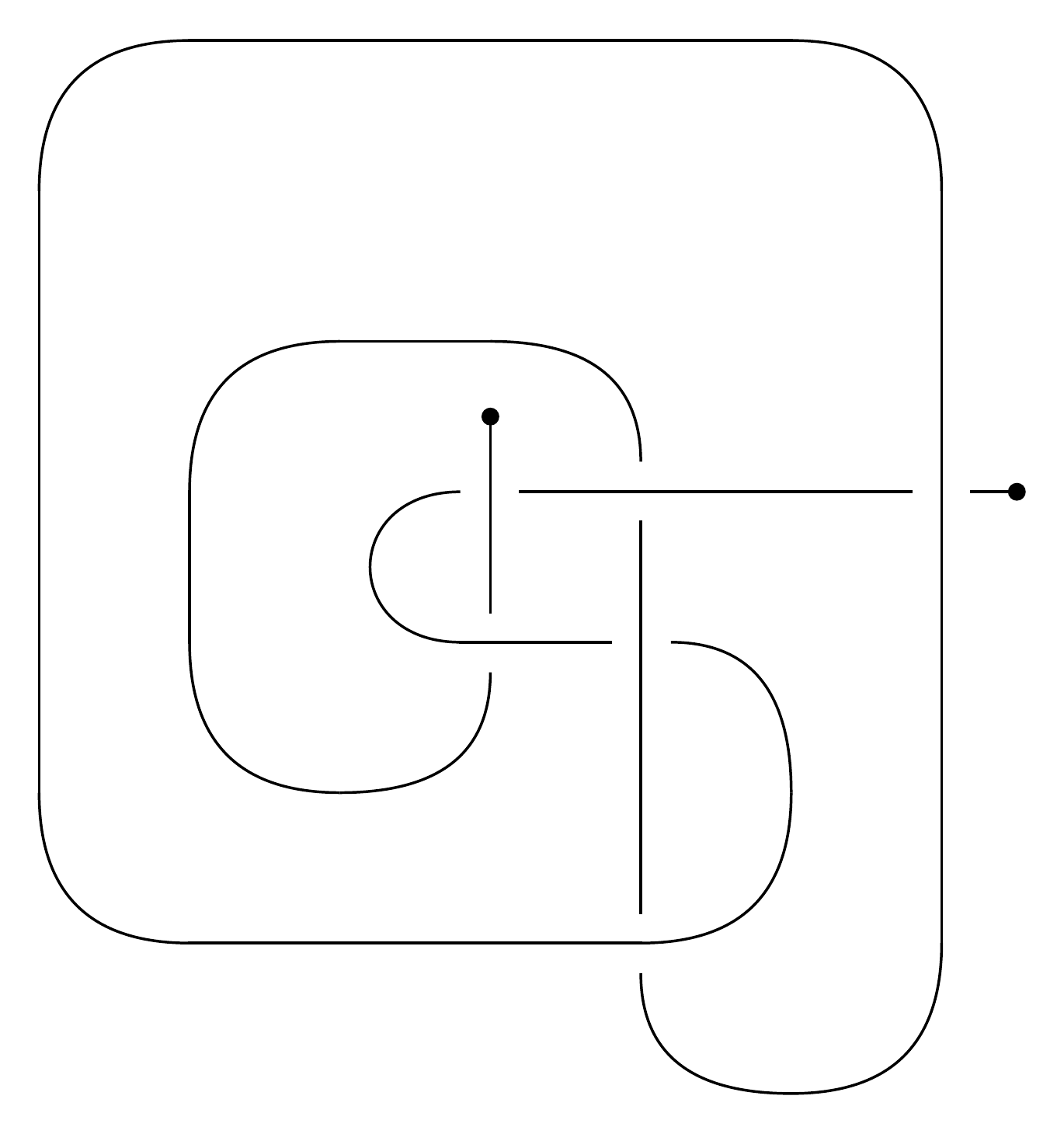}\\
\textcolor{black}{$6_{111}$}
\vspace{1cm}
\end{minipage}
\begin{minipage}[t]{.25\linewidth}
\centering
\includegraphics[width=0.9\textwidth,height=3.5cm,keepaspectratio]{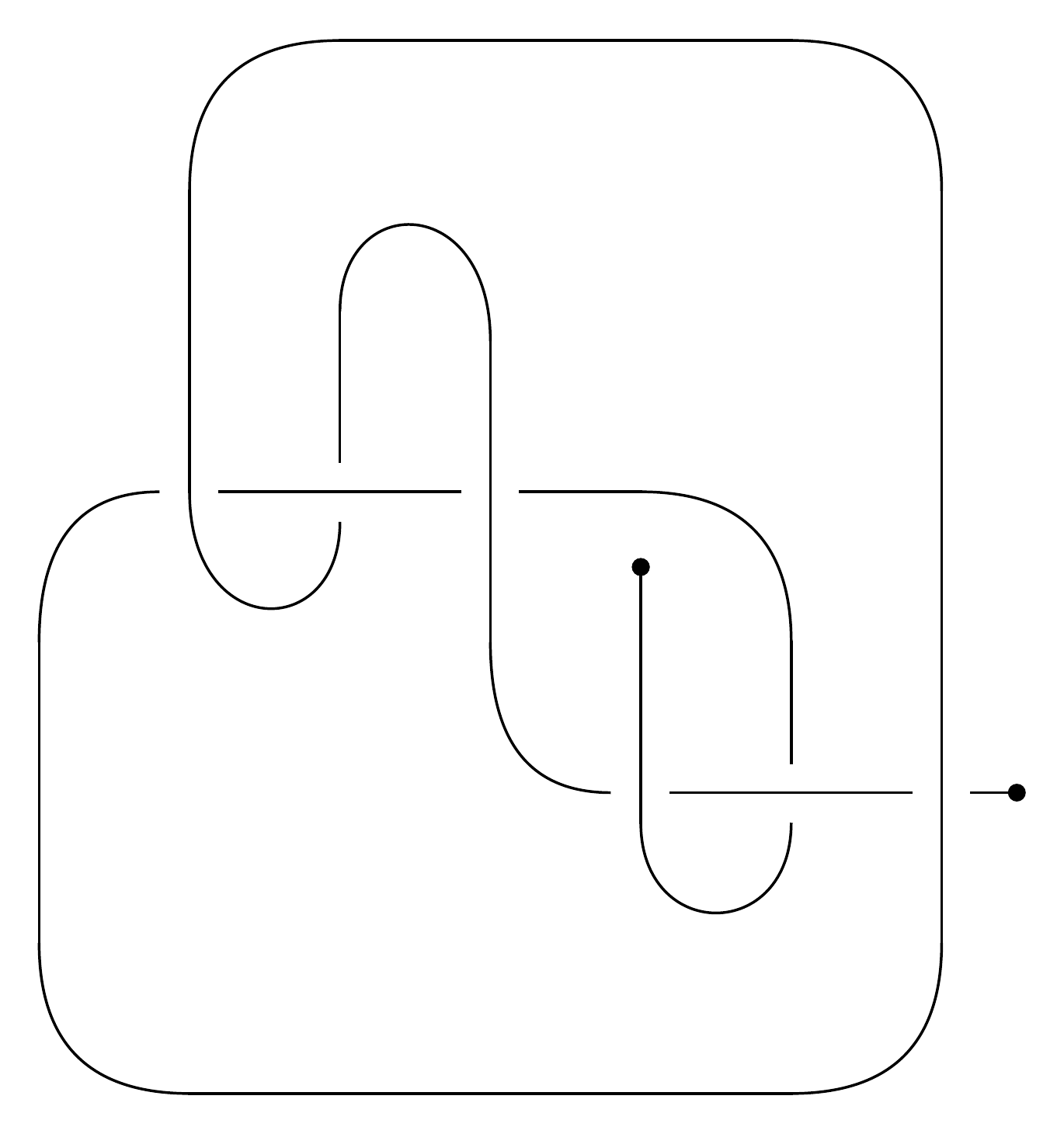}\\
\textcolor{black}{$6_{112}$}
\vspace{1cm}
\end{minipage}
\begin{minipage}[t]{.25\linewidth}
\centering
\includegraphics[width=0.9\textwidth,height=3.5cm,keepaspectratio]{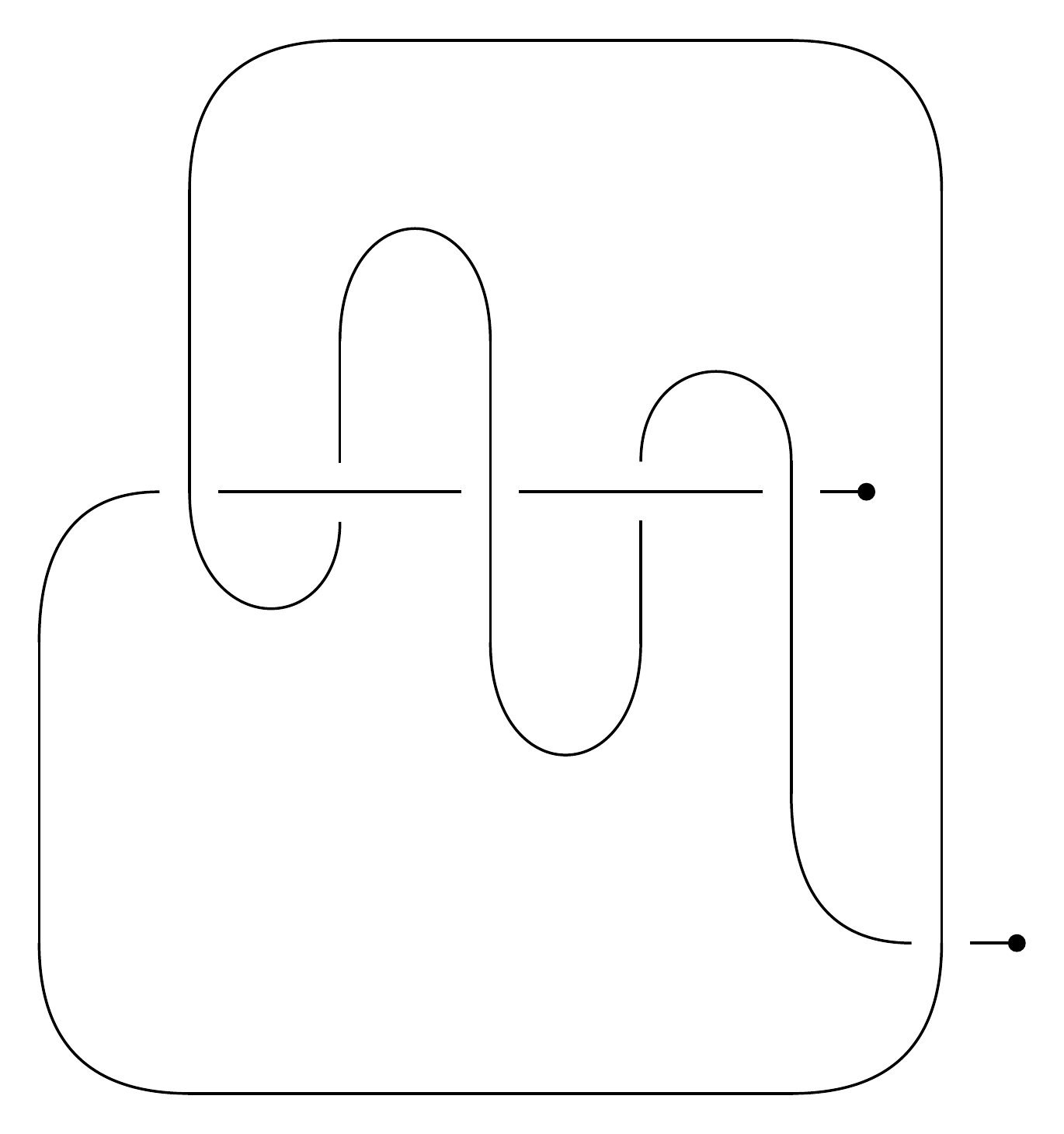}\\
\textcolor{black}{$6_{113}$}
\vspace{1cm}
\end{minipage}
\begin{minipage}[t]{.25\linewidth}
\centering
\includegraphics[width=0.9\textwidth,height=3.5cm,keepaspectratio]{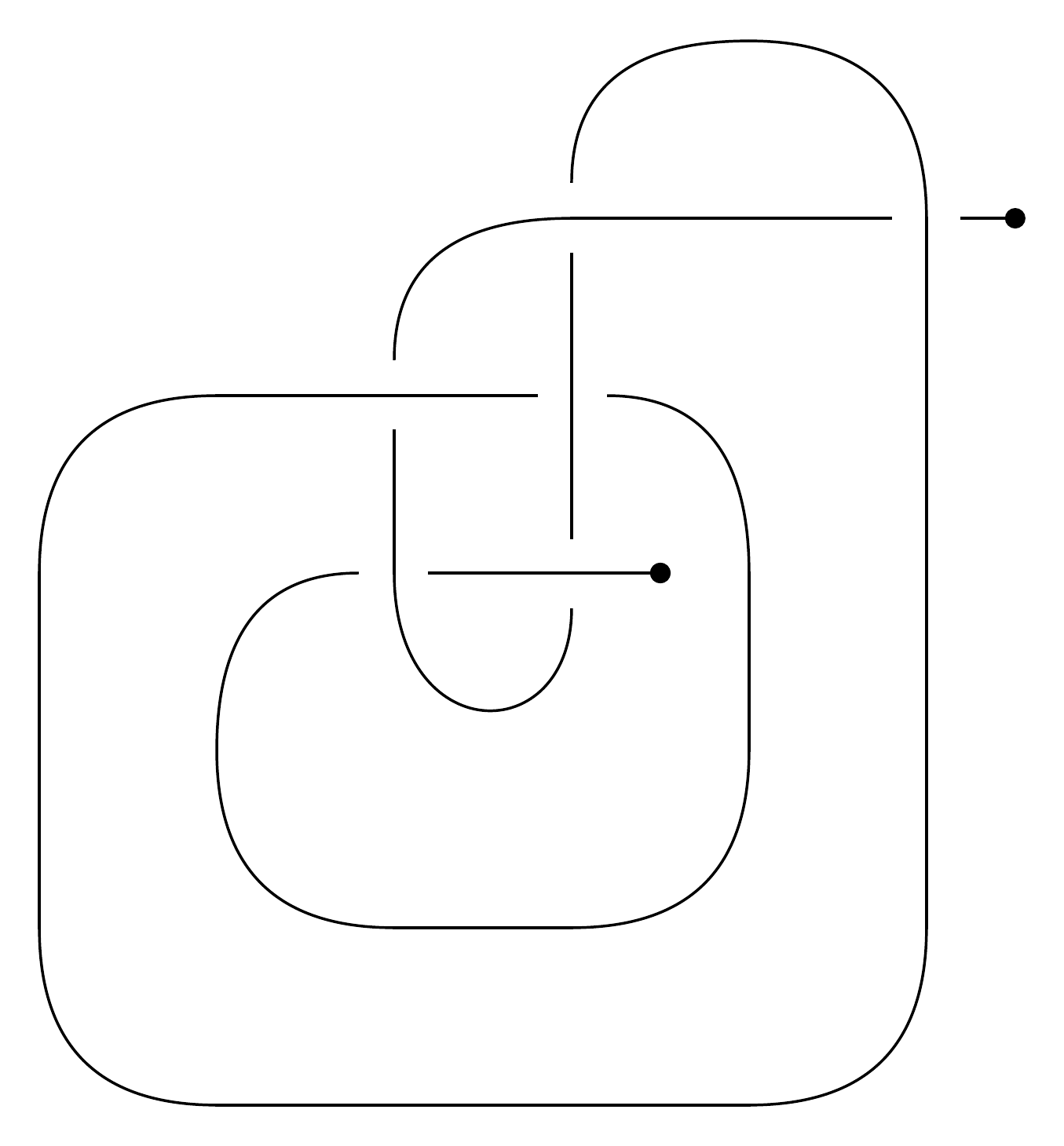}\\
\textcolor{black}{$6_{114}$}
\vspace{1cm}
\end{minipage}
\begin{minipage}[t]{.25\linewidth}
\centering
\includegraphics[width=0.9\textwidth,height=3.5cm,keepaspectratio]{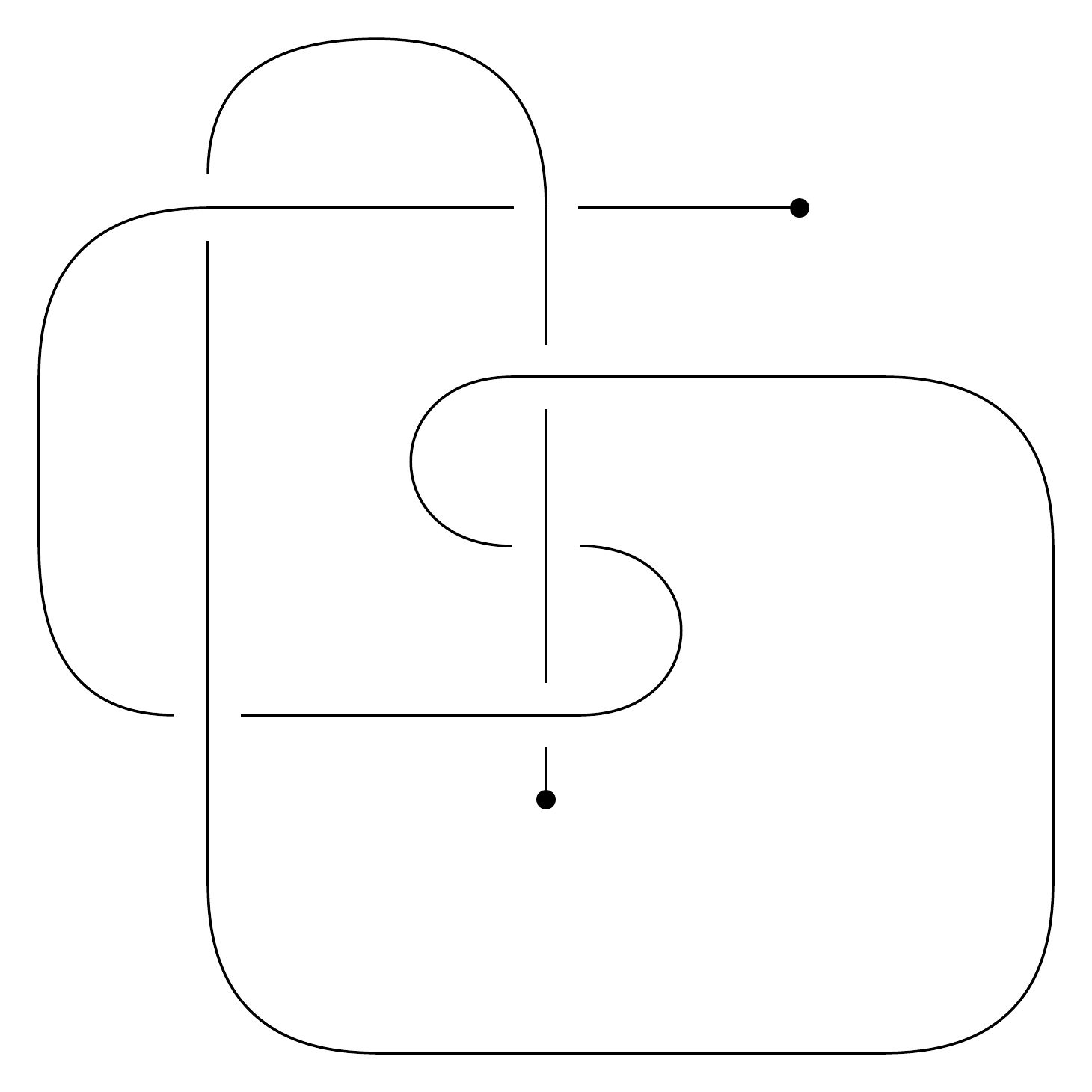}\\
\textcolor{black}{$6_{115}$}
\vspace{1cm}
\end{minipage}
\begin{minipage}[t]{.25\linewidth}
\centering
\includegraphics[width=0.9\textwidth,height=3.5cm,keepaspectratio]{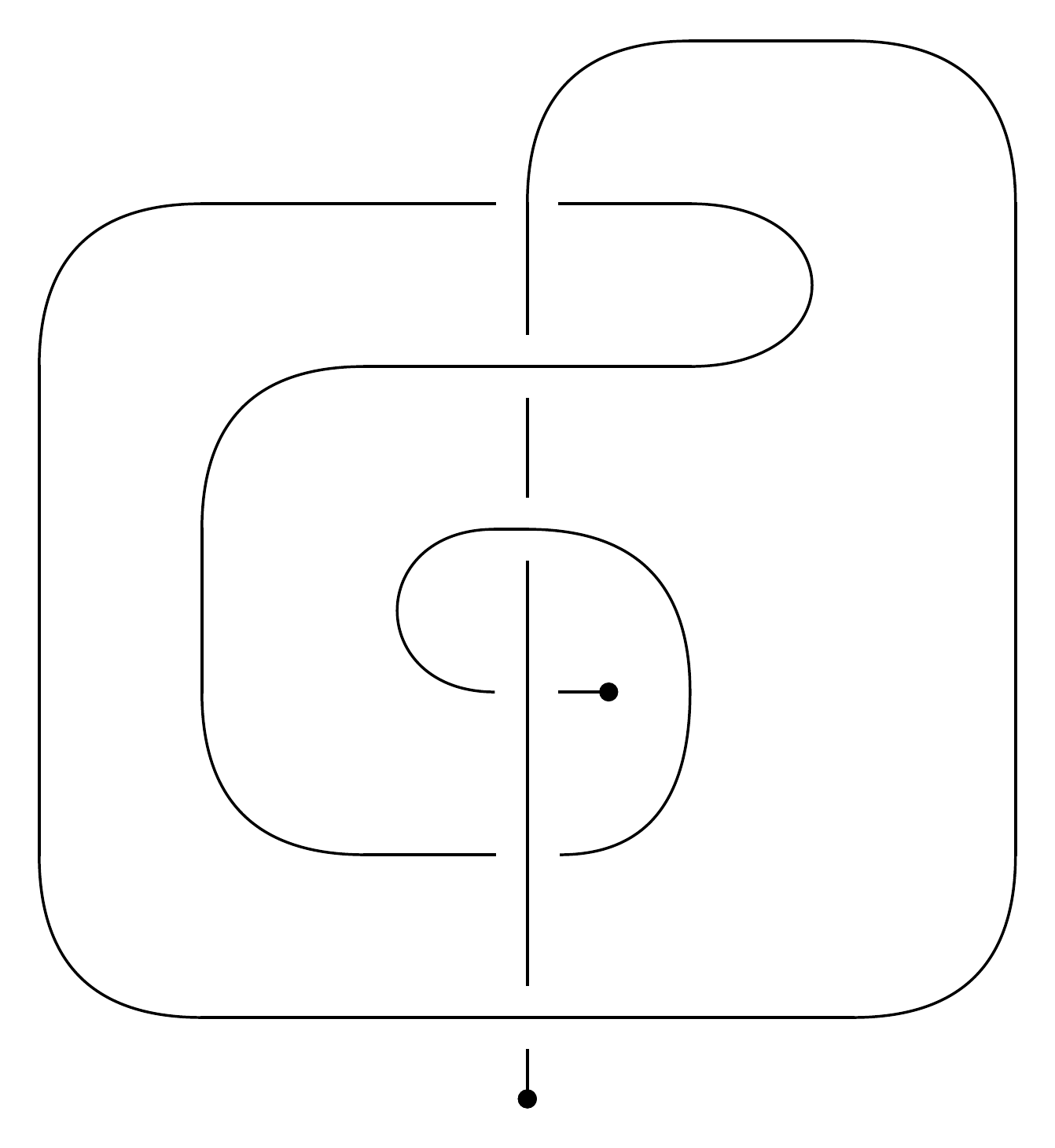}\\
\textcolor{black}{$6_{116}$}
\vspace{1cm}
\end{minipage}
\begin{minipage}[t]{.25\linewidth}
\centering
\includegraphics[width=0.9\textwidth,height=3.5cm,keepaspectratio]{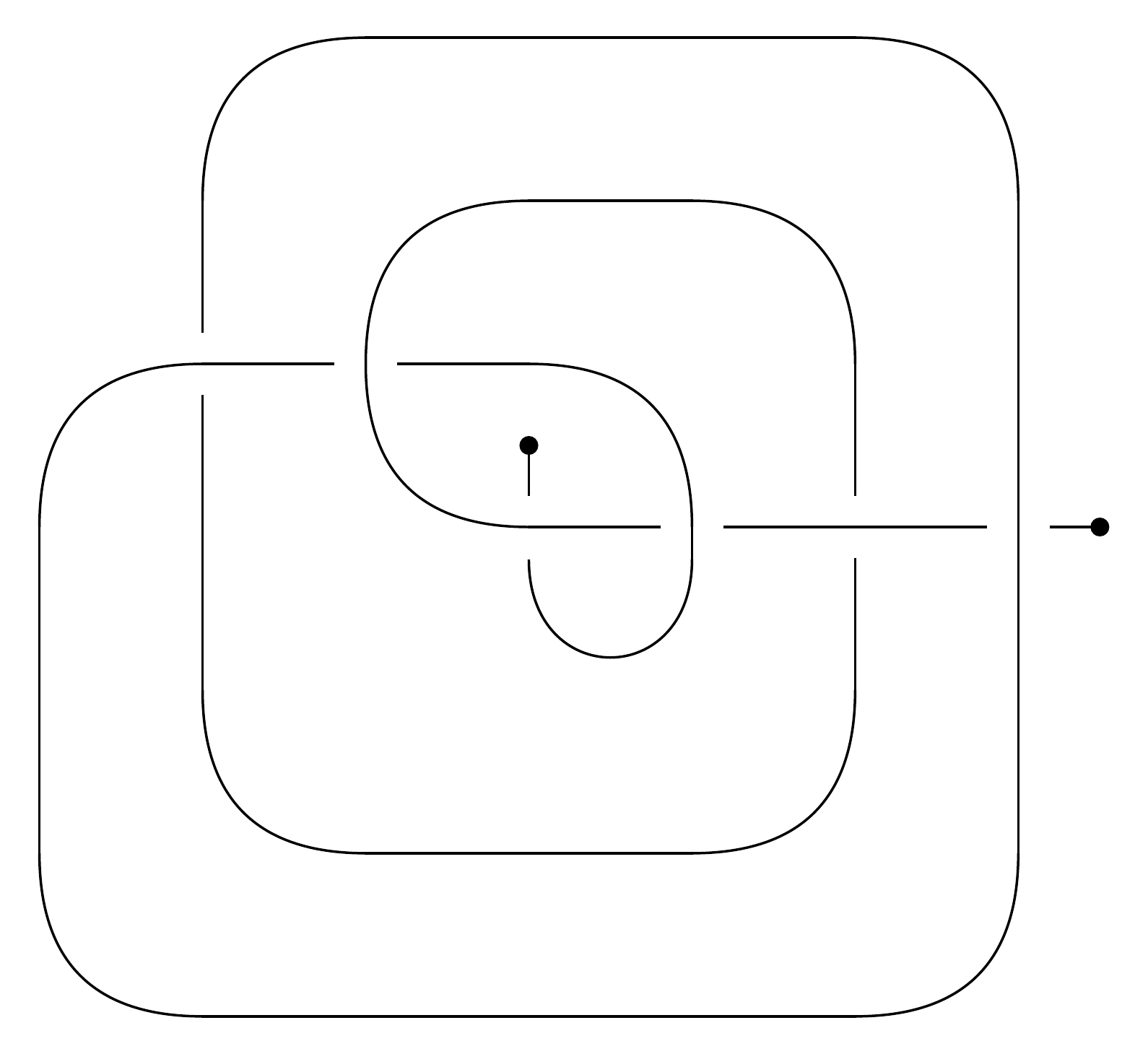}\\
\textcolor{black}{$6_{117}$}
\vspace{1cm}
\end{minipage}
\begin{minipage}[t]{.25\linewidth}
\centering
\includegraphics[width=0.9\textwidth,height=3.5cm,keepaspectratio]{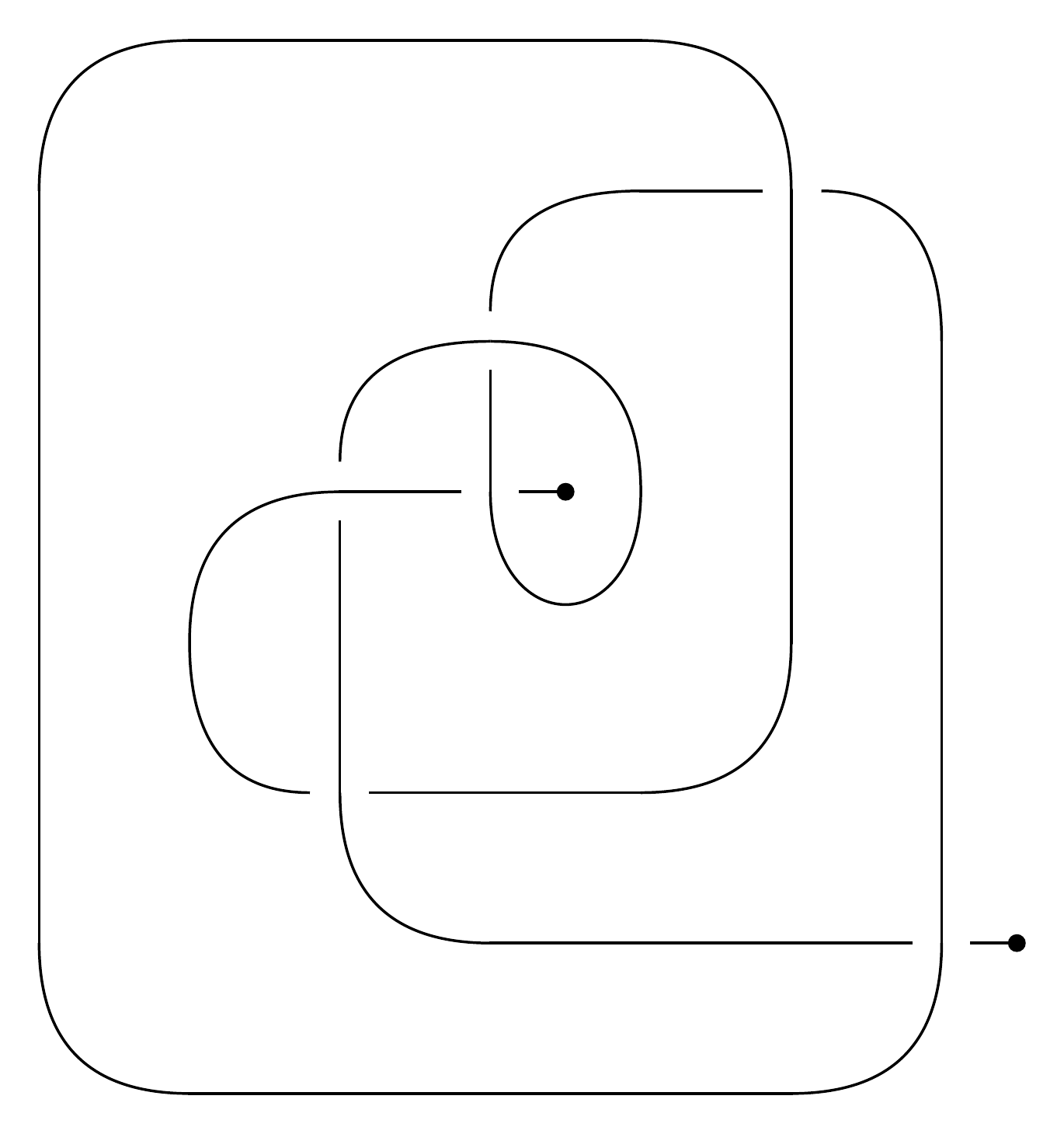}\\
\textcolor{black}{$6_{118}$}
\vspace{1cm}
\end{minipage}
\begin{minipage}[t]{.25\linewidth}
\centering
\includegraphics[width=0.9\textwidth,height=3.5cm,keepaspectratio]{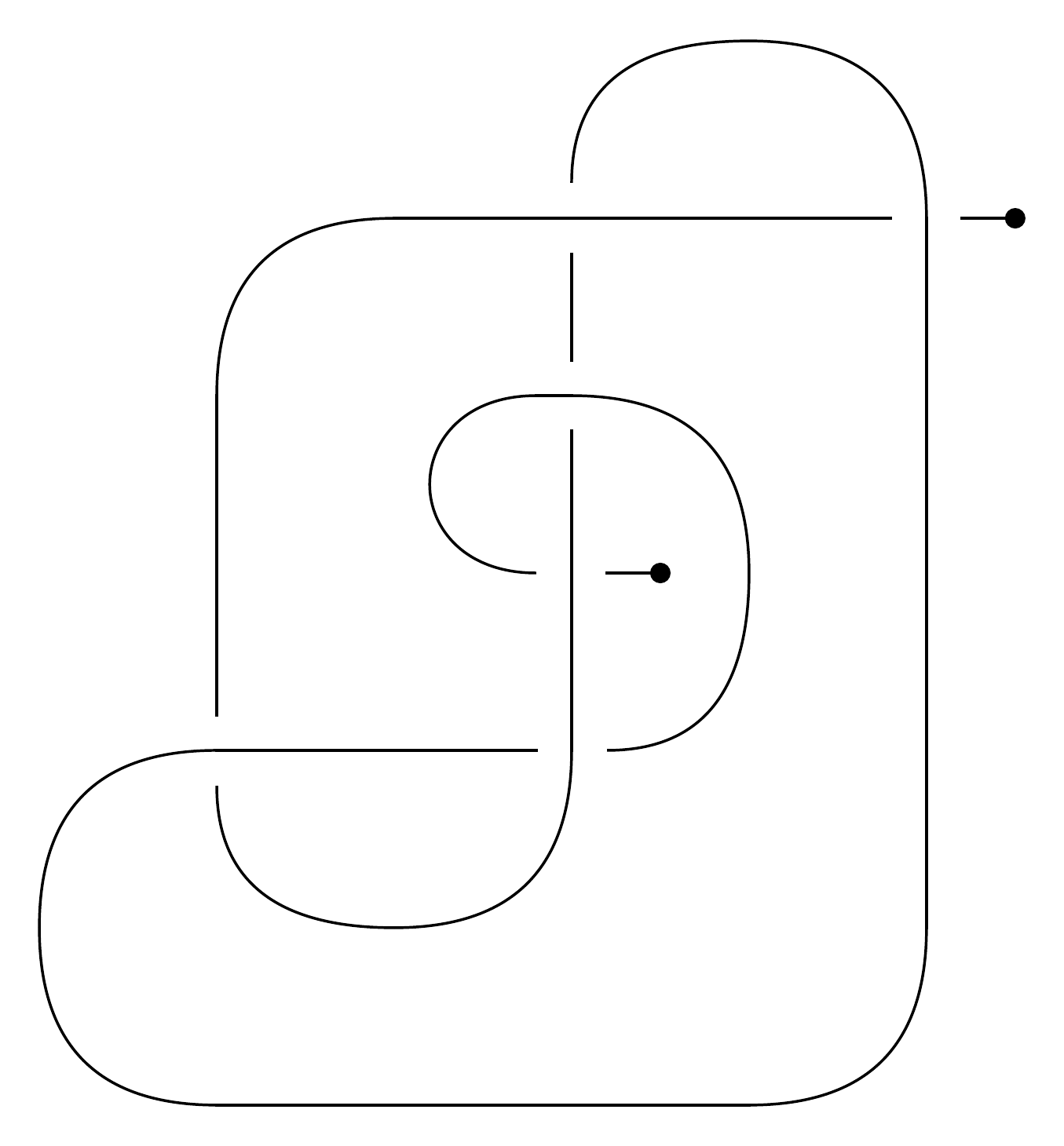}\\
\textcolor{black}{$6_{119}$}
\vspace{1cm}
\end{minipage}
\begin{minipage}[t]{.25\linewidth}
\centering
\includegraphics[width=0.9\textwidth,height=3.5cm,keepaspectratio]{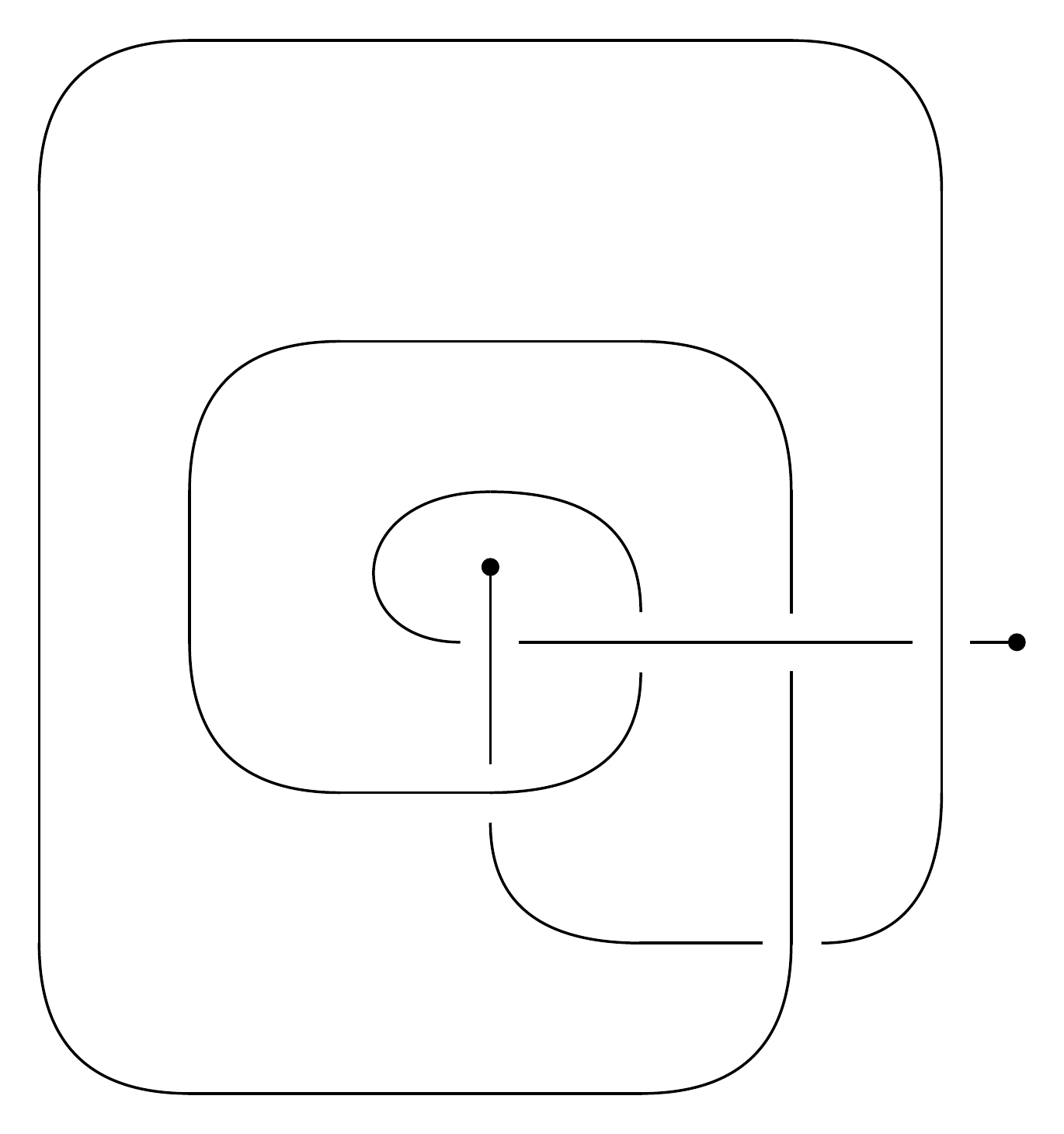}\\
\textcolor{black}{$6_{120}$}
\vspace{1cm}
\end{minipage}
\begin{minipage}[t]{.25\linewidth}
\centering
\includegraphics[width=0.9\textwidth,height=3.5cm,keepaspectratio]{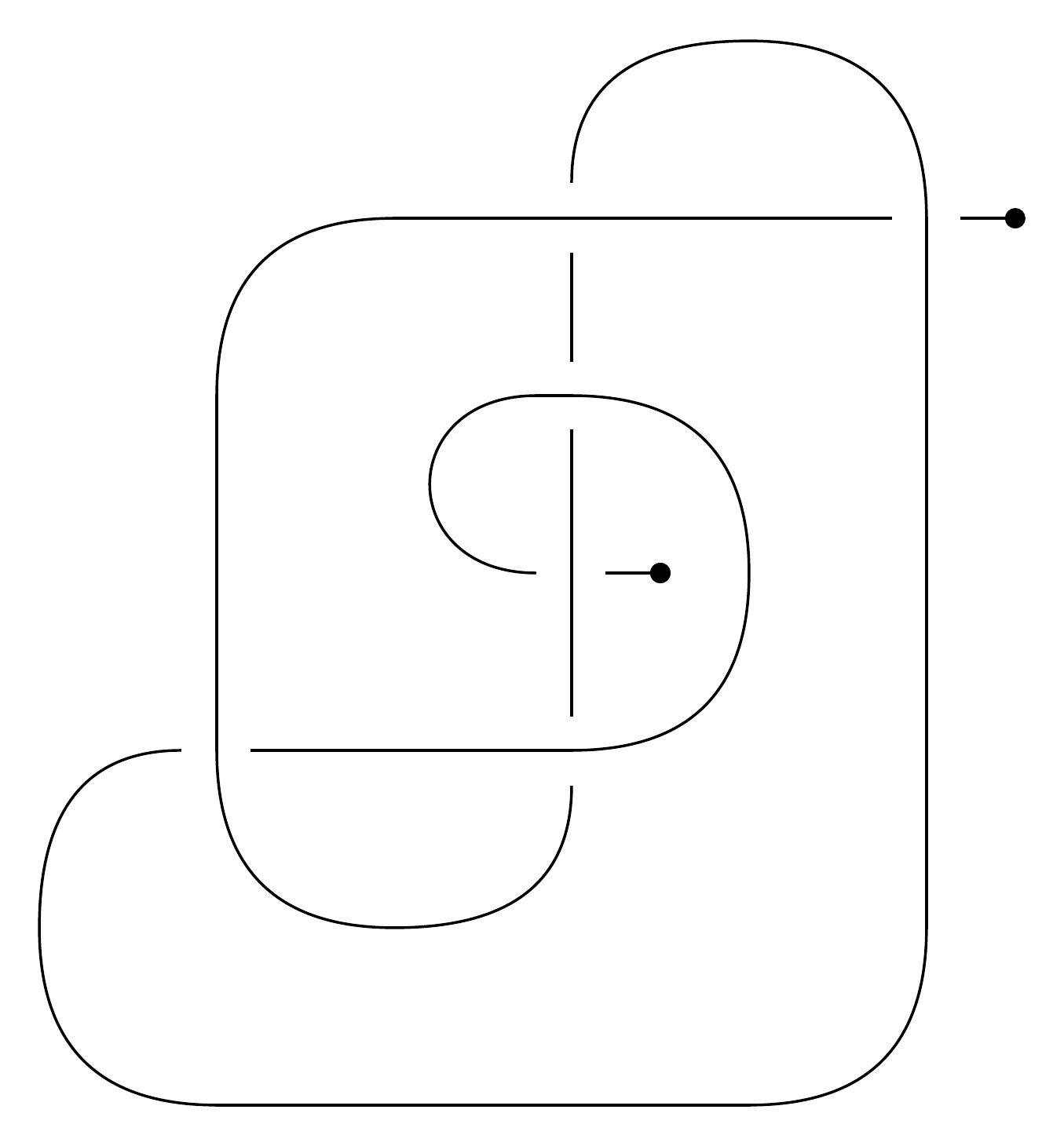}\\
\textcolor{black}{$6_{121}$}
\vspace{1cm}
\end{minipage}

%% file: knotoid_class.bbl
\begin{thebibliography}{99}

\bibitem{Adams}{\sc C. Adams}. \newblock{The Knot Book. Freeman and Company, New York (1994).}

\bibitem{barth} {\sc A. Bartholomew} \newblock \url{http://www.layer8.co.uk/maths/knotoids/index.htm}, Feb. 19 2019.

\bibitem{barbensi} {\sc A. Barbensi, D. Buck, H. A. Harrington and M. Lackenby}, {\em Double branched covers of knotoids}. arxiv:1811.09121 (2018).

\bibitem{chapman} {\sc H. Chapman} {\em Slipknotting in random diagrams}. \newblock arXiv:1803.07114 (2018).

\bibitem{dabrowski2015} {\sc P. Dabrowski-Tumanski, A. I. Jarmolinska and J. I. Sulkowska JI.} {\em Prediction of the optimal set of contacts to fold the smallest knotted protein}. \newblock Journal of Physics-Condensed Matter. 27 no 35 (2015). doi: 10.1088/0953-8984/27/35/354109. 

\bibitem{dabrowski2016} {\sc P. Dabrowski-Tumanski, A. Stasiak, and J.I. Sulkowska}, {\em In Search of Functional Advantages of Knots in Proteins}. \newblock PLoS ONE 11 (2016), e0165986, doi:10.1371/journal.pone.0165986.

\bibitem{knotprot2} {\sc P. Dabrowski-Tumanski, P. Rubach, D. Goundaroulis, J. Dorier, K. C. Millett, E.J. Rawdon, P.Su{\l}\-kowski, A. Stasiak, J.I. Sulkowska.} {\em KnotProt 2.0: a database of proteins with knots and other entangled structures}. \newblock Nucleic Acids Res. 47 D1, (2018) D367-D375. \url{https://doi.org/10.1093/nar/gky1140}.

\bibitem{knoto-id} {\sc J. Dorier, D. Goundaroulis, F. Benedetti and A. Stasiak}, {\em Knoto-ID: a tool to study the entanglement of open protein chains using the concept of knotoids.} Bioinformatics,  34 no 19 (2018),  3402-3404, DOI \url{https://doi.org/10.1093/bioinformatics/bty365}.

\bibitem{kauffmandye2009} {\sc H. A. Dye and L. H. Kauffman}, {\em Virtual crossing number and the arrow polynomial.} \newblock J. Knot. Theor. Ramif. 18 no. 10 (2009), 1335-1357.

\bibitem{gabrov} {\sc B. Gabrov{\v{s}}ek and M. Mroczkowski}, {\em Knots in the solid torus up to 6 crossings} \newblock J Knot Theor Ramif 21 (11) (2012).

\bibitem{gabrov1} {\sc B. Gabrov{\v{s}}ek}, {\em Tabulation of prime knots in lens spaces} \newblock Mediterranean Journal of Mathematics 14 (2) (2017).

\bibitem{goundaroulis_sr2017} {\sc  D. Goundaroulis, J. Dorier, F. Benedetti and A. Stasiak}, {\em Studies of global and local entanglements of individual protein chains using the concept of knotoids}. \newblock Sci. Rep. 7 no 1 (2017), 6309.

\bibitem{goundaroulis_pol2017} {\sc  D. Goundaroulis, N. G{\"u}g{\"u}mc{\"u}, S. Lambropoulou, J. Dorier, A. Stasiak and L. Kauffman}, {\em 
Topological Models for Open-Knotted Protein Chains Using the Concepts of Knotoids and Bonded Knotoids.} \newblock Polymers 9 no 9 (2017), 444.

\bibitem{gugumcu2017} {\sc N. G\"ug\"umc\"u and L. Kauffman}, {\em New invariants of knotoids}.  \newblock European J. Comb. 65 (2017), pp. 186-229. DOI:10.1016/j.ejc.2017.06.004.

\bibitem{gugumcu2017_2} {\sc N. G\"ug\"umc\"u and S. Lambropoulou}, {\em Knotoids, braidoids and applications}. \newblock Symmetry 9 no12 (2017) pp 315.

\bibitem{gugumcu2018} {\sc N. G\"ug\"umc\"u and S. Nelson}, {Biquandle Coloring Invariants of Knotoids}. \newblock arXiv:1803.11308 (2018).

\bibitem{kodokostas} {\sc D. Kodokostas and S. Lambropoulou}, {\em Rail knotoids}. \newblock arXiv:1812.09493 (2018).

\bibitem{korablev2017} {\sc P.G. Korablev, and Y. K. May}, {\em Knotoids and knots in the thickened torus.} Sib. Math. J. 58 no. 5 (2017).

\bibitem{korablev2018} {\sc P. G. Korablev, Y. K. May and V. V. Tarkaev} {\em Classification of low complexity knotoids} Sib \`Elektron. Mat. Izv. 15 (2018), 1237-1244.

\bibitem{mallam2012} {\sc A. L. Mallam and S. E. Jackson},{\em Knot formation in newly translated proteins is spontaneous and accelerated by chaperonins}. \newblock Nat. Chem. Biol. 8 no 2 (2012), 147-53. doi: 10.1038/nchembio.742. 

\bibitem{przybyl} {\sc P. Piera{\'n}ski, P. Sylwester and A. Stasiak}, {\em Tight open knots}. \newblock Eur. Phys. J. E, 6 no 2 (2001), 123-128.

\bibitem{realize} {\sc G. Vijayan and A. Wigderson}, {\em Planarity of edge ordered graphs}, Technical Report 307, Department of Electrical Engineering and Computer Science, Princeton University.

\bibitem{sage}   {\sc The Sage Developers}, {\em {S}ageMath, the {S}age {M}athematics {S}oftware {S}ystem ({V}ersion x.y.z)} (2018) \url{http://www.sagemath.org}   Accessed 6 Dec 2018.

\bibitem{sriramoju2018} {\sc  M. K. Sriramoju, Y. Chen, Y-T. C. Lee, and S-T. D. Hsu} {\em Topologically knotted deubiquitinases exhibit unprecedented mechanostability to withstand the proteolysis by an AAA+ protease} \newblock Sci. Rep. 8 no 1 (2018), 7076.

\bibitem{sulkowska_on2012} {\sc J. I. Sulkowska, J. K. Noel and J. N. Onuchic}, {\em Energy landscape of knotted protein folding}. Proc. Natl. Acad. Sci. U.S.A. 109 no 44 (2012), 17783-8. doi: 10.1073/pnas.1201804109.

\bibitem{sulkowska_2012} {\sc J.I. Sulkowska, E. J. Rawdon, K. C. Millett, J. N. Onuchic and A. Stasiak}  {\em Conservation of complex knotting and slipknotting patterns in proteins}, Proc. Natl. Acad. Sci. U.S.A. 109 (2012), E1715?E1723, doi: 10.1073/pnas.1205918109.

\bibitem{turaev2012} {\sc V. Turaev}, {\em Knotoids}. \newblock Osaka J. Math. 1 no 49 (2012), pp. 195-223.

\bibitem{virnau2006} {\sc P. Virnau, L.A. Mirny and M. Kardar}, {\em Intricate knots in proteins: Function and evolution}. \newblock PLoS Comput. Biol. 2 (2006), 1074-1079.






\end{thebibliography}
